\errorstopmode
\documentclass[10pt,CMpaper]{CMbook2}
\usepackage{ifpdf}
\ifpdf%
\PassOptionsToPackage{pdftex}{graphicx}
\PassOptionsToPackage{pdftex}{color}
\else%
\PassOptionsToPackage{dvips}{graphicx}%
\fi

\usepackage{graphicx,color}
\usepackage{eurosym}
\usepackage{
amsfonts,amssymb}
\usepackage{mathrsfs}
\usepackage{oldstyle}
\usepackage[T1]{fontenc}
\usepackage[utf8]{inputenc} 
\usepackage[USenglish]{babel} %
\usepackage[all]{xy}
\usepackage{pb-diagram}
\usepackage[nohints]{minitoc} 
\usepackage{multicol}
\usepackage{CMmacros2}
\usepackage{makeidx}
\makeindex
\usepackage{lmodern}
\usepackage{microtype}
\usepackage{array}
\usepackage{stmaryrd}
\usepackage{float}
\usepackage{wasysym}
\usepackage{bm}  

\usepackage{enumerate} 
\usepackage[novbox]{pdfsync} 

\DeclareFontFamily{U}{wasy}{}
\DeclareFontShape{U}{wasy}{m}{n}{<5> <6> <7> <8> <9> gen * wasy
      <10> <10.95> <12> <14.4> <17.28> <20.74> <24.88>wasy10  }{}
\DeclareFontShape{U}{wasy}{b}{n}{ <-10> ssub * wasy/m/n
 <10> <10.95> <12> <14.4> <17.28> <20.74> <24.88>wasyb10 }{}
\DeclareFontShape{U}{wasy}{bx}{n}{<5> <6> <7> <8> <9> gen * wasy
 <10> <10.95> <12> <14.4> <17.28> <20.74> <24.88>wasyb10}{}
\DeclareSymbolFont{wasy}{U}{wasy}{m}{n}
\SetSymbolFont{wasy}{bold}{U}{wasy}{bx}{n}
\DeclareFontFamily{U}{lasy}{}
\DeclareFontShape{U}{lasy}{m}{n}{ <5> <6> <7> <8> <9> gen * lasy
      <10> <10.95> <12> <14.4> <17.28> <20.74> <24.88>lasy10  }{}
\DeclareFontShape{U}{lasy}{b}{n}{ <5> <6> <7> <8> <9> gen * lasy
      <10> <10.95> <12> <14.4> <17.28> <20.74> <24.88>lasyb10  }{}
\DeclareFontFamily{U}{stmry}{}
\DeclareFontShape{U}{stmry}{m}{n}
   {  <5> <6> <7> <8> <9> <10> gen * stmary
      <10.95><12><14.4><17.28><20.74><24.88>stmary10%
   }{}
\DeclareFontShape{U}{stmry}{b}{n}
   {  <5> <6> <7> <8> <9> <10> gen * stmary
      <10.95><12><14.4><17.28><20.74><24.88>stmary10%
   }{}

\ifpdf
\makeatletter
\def\toclevel@chapternb{0}
\makeatother
\usepackage[bookmarksopen=false,pdftex=true,breaklinks=true,backref=page,%
    pagebackref=true,plainpages=false,%
    hyperindex=true,pdfstartview=FitH,
    pdfpagelabels=true,colorlinks=false,citecolor=red,colorlinks=false]%
    {hyperref}

\else
\newcommand{\texorpdfstring}[2]{#1}
\fi

\DeclareSymbolFont{lasy}{U}{lasy}{m}{n}
\SetSymbolFont{lasy}{bold}{U}{lasy}{b}{n}
\let\Box\undefined
\DeclareMathSymbol\Box{\mathord}{lasy}{"32}

\title{Commutative algebra: \\[2mm]  Constructive methods \\[6mm]
\LARGE Finite projective modules
}
\author{Henri Lombardi \&\
Claude Quitt\'e}
\date{\today }
\begin{document}
\overfullrule=0cm
\hfuzz14pt

\let\oldref\ref
\renewcommand{\ref}[1]{\hbox{\oldref{#1}}}

\makeatletter
\def\U@@@ref#1\relax#2\relax{\let\showsection\relax#2}
\newcommand{\U@@ref}[5]{\U@@@ref#1}
\newcommand{\U@ref}[1]{\NR@setref{#1}\U@@ref{#1}}
\newcommand{\@iref}{\@ifstar\@refstar\U@ref}
\newcommand{\iref}[1]{\hbox{\@iref{#1}}}
\def\l@chapter#1#2{\ifnum\c@tocdepth>\m@ne\addpenalty{-\@highpenalty}\vskip1.0em\@plus\p@
\setlength\@tempdima{2.8em}\begingroup\parindent\z@\rightskip\@pnumwidth\parfillskip-\@pnumwidth\penalty-2000\leavevmode\bfseries
\advance\leftskip\@tempdima\hskip-\leftskip\hskip0.8em{\boldmath#1}\nobreak\hfil\nobreak\hb@xt@\@pnumwidth{\hss}\par\penalty\@highpenalty\endgroup\fi}
\def\l@chapterbis#1#2{\ifnum\c@tocdepth>\m@ne\addpenalty{-\@highpenalty}\vskip1.0em\@plus\p@
\setlength\@tempdima{2em}\begingroup\parindent\z@\rightskip\@pnumwidth\parfillskip-\@pnumwidth\penalty-2000\leavevmode\bfseries
\advance\leftskip\@tempdima\hskip-\leftskip{\boldmath#1}\nobreak\hfil\nobreak\hb@xt@\@pnumwidth{\hss#2}\par\penalty\@highpenalty\endgroup\fi}
\def\toclevel@chapterbis{0}
\def\l@section{\@dottedtocline{1}{1em}{0em}}
\def\l@subsection{\@dottedtocline{2}{1.75em}{0.25em}}
\def\l@subsubsection{\@dottedtocline{3}{2.5em}{0.25em}}
\def\numberline#1{\hb@xt@\@tempdima{\hss#1\ \hfil}}
\def\contentsline#1#2#3#4{\ifx\\#4\\\csname
l@#1\endcsname{#2}{#3}\else\ifHy@linktocpage\csname
l@#1\endcsname{{#2}}{\hyper@linkstart{link}{#4}{#3}\hyper@linkend }
\else\csname l@#1\endcsname{\hyper@linkstart{link}{#4}{#2}\hyper@linkend
}{#3}\fi\fi}
\makeatother

\setcounter{minitocdepth}{3}
\dominitoc

\ifx\mesMacrosDejaChargees\undefined\let\chargeMesMacros\relax
\else\let\chargeMesMacros \fi
\chargeMesMacros

\gdef\mesMacrosDejaChargees{}


\newcounter{MF}
\newcommand\stMF{\stepcounter{MF}}


\newcommand{\DBxk}{{\Der \gk\gB\xi}}%
\newcommand{\DAxk}{{\Der \gk\gA\xi}}%
\newcommand{\DkXxk}{{\Der \gk\kuX\xi}}%
\newcommand \STR {\MA{\rm SR}}
\newcommand \Str {\MA{\mathsf{Sr}}}
\newcommand \Bdim {\MA{\mathsf{Bdim}}}
\newcommand \Sdim {\MA{\mathsf{Sdim}}}
\newcommand \Cdim {\MA{\mathsf{Cdim}}}
\newcommand \Gdim {\MA{\mathsf{Gdim}}}
\newcommand \DAbul {\rD_{\!\Abul}}

\renewcommand{\labelenumii}{\theenumii.}
\newcommand \perso[1]{}
\newcommand \entrenous[1]{}
\newcommand \hum[1]{}
\newcommand \incertain[1]{}
\newcommand \Oui[1]{}

\newcounter{bidon}
\newcommand{\rdb}{\refstepcounter{bidon}}

\newcommand \ix[1] {\index{#1}\emph{#1}}
\newcommand \ixc[2] {\index{#1!#2}\emph{#1}}
\newcommand \ixx[2] {\index{#1!#2}\emph{#1 #2}}
\newcommand \ixy[2] {\index{#2!#1 ---}\emph{#1 #2}}
\newcommand \ixd[2] {\index{#2!#1}\emph{#1}}
\newcommand \ixe[2] {\index{#2@#1}\emph{#1}}
\newcommand \ixf[3] {\index{#3@#2!#1}\emph{#1}}
\newcommand \ixg[3] {\index{#3@#1!#2}\emph{#1}}

\newcommand \Grandcadre[1]{%
\begin{center}
\begin{tabular}{|c|}
\hline
~\\[-3mm]
#1\\[-3mm]
~\\ 
\hline
\end{tabular}
\end{center}

}

\CMnewtheorem{theorem}{Theorem}{\itshape}   
\CMnewtheorem{thdef}{Theorem and definition}{\itshape}
\CMnewtheorem{plcc}{Concrete local-global principle}{\itshape}
\CMnewtheorem{prcf}{Closed covering principle}{\itshape}
\CMnewtheorem{prmf}{Closed patching principle}{\itshape}
\CMnewtheorem{prvq}{Covering principle by quotients}{\itshape}
\CMnewtheorem{prmq}{Patching principle of quotients}{\itshape}
\CMnewtheorem{plca}{Abstract local-global principle$\mathbf{^*}$}{\itshape}
\CMnewtheorem{plcd}{Dynamical local-global principle}{\itshape}

\CMnewtheorem{proposition}{Proposition}{\itshape}
\CMnewtheorem{propdef}{Proposition and definition}{\itshape}
\CMnewtheorem{lemma}{Lemma}{\itshape}
\CMnewtheorem{corollary}{Corollary}{\itshape}
\CMnewtheorem{fact}{Fact}{\itshape}
\CMnewtheorem{theoremc}{Theorem\etoz}{\itshape}
\CMnewtheorem{lemmac}{Lemma\etoz}{\itshape}
\CMnewtheorem{corollaryc}{Corollary\etoz}{\itshape}
\CMnewtheorem{proprietec}{Property\etoz}{\itshape}
\CMnewtheorem{propositionc}{Proposition\etoz}{\itshape}
\CMnewtheorem{factc}{Fact\etoz}{\itshape}

\CMnewtheorem{remark}{Remark}{}
\CMnewtheorem{remarks}{Remarks}{}
\CMnewtheorem{comment}{Comment}{}
\CMnewtheorem{comments}{Comments}{}
\CMnewtheorem{example}{Example}{}
\CMnewtheorem{examples}{Examples}{}
\CMnewtheorem{definition}{Definition}{}
\CMnewtheorem{definitions}{Definitions}{}
\CMnewtheorem{definota}{Definition and notation}{}
\CMnewtheorem{definotas}{Definitions and notations}{}
\CMnewtheorem{convention}{Convention}{}
\CMnewtheorem{notation}{Notation}{} 
\CMnewtheorem{notations}{Notations}{} 
\CMnewtheorem{question}{Question}{}
\CMnewtheorem{algorithm}{Algorithm}{}

\newcounter{exercise}[chapter]
\newenvironment{exercise}{\ifhmode\par\fi
\vskip-\lastskip\vskip1.5ex\mou\penalty-300 \relax
\everypar{}\noindent
\refstepcounter{exercise}{\bfseries Exercise \theexercise.}\relax
\itshape\ignorespaces}{\par\vskip-\lastskip\vskip1em}

\newcounter{problem}[chapter]
\newenvironment{problem}{\ifhmode\par\fi
\vskip-\lastskip\vskip1em\mou\penalty-300 \relax
\everypar{}\noindent
\refstepcounter{problem}{\bfseries Problem \theproblem.}\relax
\itshape\ignorespaces}{\par\vskip-\lastskip\vskip1em}

\newcommand
\CHAP[1]{ 
\goodbreak\vskip4mm \mou \noindent  {\bf #1}\par\nobreak
\vskip1mm \mou \nobreak
}

\newcommand {\junk}[1]{}
\newcommand {\eop}{\hbox{$\square$}}

\newcommand\DebP{\raisebox{1pt}{\large\bf \Rightcircle}\,}
\newcommand\FinP{\raisebox{1pt}{\large\bf \Leftcircle$\!\!$}}

\def\thefootnote{\arabic{footnote}}

\newenvironment{proof}{\ifhmode\par\fi\vskip-\lastskip\vskip0.5ex\global\insidedemotrue
\everypar{}\noindent{\setbox0=\hbox{$\!\!$\DebP }\global\wdTitreEnvir\wd0\box0}
\ignorespaces}{\enddemobox\par\vskip.5em}
\def\enddemobox{\ifinsidedemo
\ifmmode\hbox{$\square$}\else
\ifhmode\unskip\else\noindent\fi\nobreak\null\nobreak\hfill
\nobreak$\square$\fi\fi\global\insidedemofalse}

\newenvironment{Proof}[1]{\ifhmode\par\fi\vskip-\lastskip\vskip0.5ex\global\insidedemotrue
\everypar{}\noindent{\setbox0=\hbox{\it #1}\global\wdTitreEnvir\wd0\box0}\ignorespaces}{\enddemobox\par\vskip.5em}
\def\enddemobox{\ifinsidedemo
\ifmmode\hbox{$\square$}\else
\ifhmode\unskip\else\noindent\fi\nobreak\null\nobreak\hfill
\nobreak$\square$\fi\fi\global\insidedemofalse}

\newcommand\facile{\begin{proof}
The \dem is left to the reader.
\end{proof}
}

\newcommand\dsp{\displaystyle}

\newcommand\mapright[1]{\smash{\mathop{\longrightarrow}\limits^{#1}}} 
\newcommand\maprightto[1]{\smash{\mathop{\longmapsto}\limits^{#1}}} 
\def\mapdown#1{\downarrow\rlap{$\vcenter{\hbox{$\scriptstyle 
#1$}}$}}


\newcommand\gui[1]{``#1''}

\newcommand\subsec[1]{\goodbreak\rdb\addcontentsline{toc}{subsection}{#1}\subsection*{#1}}

\newcommand\subsect[2]{\goodbreak\rdb\addcontentsline{toc}{subsection}{#2}
\subsection*{#1}}

\newcommand\subsubs[1]{

\goodbreak\rdb\medskip 

{\bf #1}

\smallskip }

\newcommand\subsubsec[1]{\goodbreak\rdb\addcontentsline{toc}{subsubsection}{#1}\subsubsection*{#1}}

\newcommand\subsubsect[2]{\goodbreak\rdb\addcontentsline{toc}{subsubsection}{#2}\subsubsection*{#1}}

\newcommand\Subsubsec[1]{\goodbreak\rdb\addcontentsline{toc}{subsection}{#1}\subsubsection*{#1}}

\newcommand\Subsubsect[2]{\rdb\addcontentsline{toc}{subsection}{#2}\subsubsection*{#1}}

\newcommand\Subsec[1]{\goodbreak\rdb\addcontentsline{toc}{section}{#1}\subsection*{#1}}

\newcommand\Exercices{\rdb\addcontentsline{toc}{section}{Exercises and \pbsz}%
\markright{Exercises and \pbsz}\section*{Exercises and \pbsz}\pagestyle{CMExercicesheadings}\small }

\newcommand\sol{%
\rdb\normalsize\addcontentsline{toc}{subsection}{Solutions of selected exercises}%
\markright{Solutions of selected exercises}\subsection*{Some solutions, or sketches of solutions}%
\pagestyle{CMExercicesheadings}\small}

\newcommand\exer[1]{{\goodbreak\smallskip\noindent\textbf{Exercise \ref{#1}.} }}

\newcommand\prob[1]{{\goodbreak\smallskip\noindent\textbf{Problem \ref{#1}.} }}

\newcommand\Biblio{\goodbreak\rdb\normalsize\addcontentsline{toc}{section}{Bibliographic comments}%
\markright{Bibliographic comments}%
\pagestyle{CMExercicesheadings}\section*{Bibliographic comments} 

}

\newcommand\Intro{\addcontentsline{toc}{section}{Introduction}\markright{Introduction}%
\pagestyle{CMExercicesheadings}\subsection*{Introduction} }

\newcommand\defa[2]{\mni\rdb\textbf{Alternative definition \ref{#1}
\it #2 
}}


\renewcommand\matrix[1]{{\begin{array}{ccccccccccccccccccccccccc} #1 \end{array}}}

\newcommand\MA[1]{\mathop{#1}\nolimits}

\newcommand{\vect}[1]{\mathchoice{\overrightarrow{\strut#1}}%
{\overrightarrow{\textstyle\strut#1}}{\overrightarrow{\scriptstyle#1}}{\overrightarrow{\scriptscriptstyle#1}}}

\newcommand\vab[2]{[\,#1\;#2\,]}

\newcommand\abs[1]{\left|{#1}\right|}
\newcommand\abS[1]{\big|{#1}\big|}
\newcommand\aqo[2]{#1\sur{\gen{#2}}\!}
\newcommand\Aqo[2]{#1\sur{\big\langle{#2}\big\rangle}\!}
\newcommand\bloc[4]{\left[\matrix{#1 & #2 \cr #3 & #4}\right]}

\newcommand\carray[2]{{\left[\begin{array}{#1} #2 \end{array}\right]}}
\newcommand\cmatrix[1]{\left[\matrix{#1}\right]}
\newcommand\clmatrix[1]{{\left[\begin{array}{lllllll} #1 \end{array}\right]}}
\newcommand\crmatrix[1]{{\left[\begin{array}{rrrrrrr} #1 \end{array}\right]}}
\newcommand\dmatrix[1]{\abs{\matrix{#1}}}
\newcommand\Cmatrix[2]{\setlength{\arraycolsep}{#1}\left[\matrix{#2}\right]}

\newcommand\Dmatrix[2]{\setlength{\arraycolsep}{#1}\left|\matrix{#2}\right|}

\newcommand{\Dpp}[2]{{{\partial #1}\over{\partial #2}}}

\newcommand\Dlu[2]{{\rm Dl}_{#1}(#2)}
\newcommand\dessus[2]{{\textstyle {#1} \atop \textstyle {#2}}}
\newcommand\eqdf[1]{\buildrel{#1}\over =}
\newcommand\formule[1]{{\left\{ {\arraycolsep2pt\begin{array}{lll} #1 \end{array}}\right.}}
\newcommand\formul[2]{{\left\{ \begin{array}{#1} #2 \end{array}\right.}}
\newcommand\gen[1]{\left\langle{#1}\right\rangle}
\newcommand\geN[1]{\big\langle{#1}\big\rangle}
\newcommand\impdef[1]{\buildrel{#1}\over \Longrightarrow}

\newcommand\eqdefi{\eqdf{\rm def}}
\newcommand\eqdef{\buildrel{{\rm def}}\over \Longleftrightarrow }

\newcommand{\Kr}[2]{#1\lrb{#2}}

\newcommand \idg[1] {|\,#1\,|}
\newcommand \idG[1] {\big|\,#1\,\big|}

\newcommand \dex[1] {[\,#1\,]}
\newcommand \deX[1] {\big[\,#1\,\big]}

\newcommand \lst[1] {[\,#1\,]}
\newcommand \lsT[1] {\big[\,#1\,\big]}

\newcommand\lra[1]{\langle{#1}\rangle}
\newcommand\lrb[1] {\llbracket #1 \rrbracket}
\newcommand\lrbn {\lrb{1..n}}
\newcommand\lrbzn {\lrb{0..n}}
\newcommand\lrbl {\lrb{1..\ell}}
\newcommand\lrbm {\lrb{1..m}}
\newcommand\lrbk {\lrb{1..k}}
\newcommand\lrbh {\lrb{1..h}}
\newcommand\lrbp {\lrb{1..p}}
\newcommand\lrbq {\lrb{1..q}}
\newcommand\lrbr {\lrb{1..r}}
\newcommand\lrbs {\lrb{1..s}}

\newcommand \fraC[2] {{{#1}\over {#2}}}
\newcommand \meck[2] {\{#1, #2\}}
\newcommand \sat[1] {#1^{\rm sat}}
\newcommand \satu[2] {#1^{\rm sat_{#2}}}
\newcommand \scp[2] {\gen{#1\,|\,#2}\!}
\newcommand\sur[1]{\!\left/#1\right.}
\newcommand\so[1]{\left\{{#1}\right\}}
\newcommand\sO[1]{\big\{{#1}\big\}}
\newcommand\sotq[2]{\so{\,#1\,\vert\,#2\,}}
\newcommand\sotQ[2]{\sO{\,#1\;\big\vert\;#2\,}}
\newcommand\frt[1]{\!\left|_{#1}\right.\!}
\newcommand\Frt[2]{\left.#1\right|_{#2}\!}
\newcommand\sims[1]{\buildrel{#1}\over \sim}
\newcommand\tra[1]{{\,^{\rm t}\!#1}}
\newcommand\Al[1]{\Vi^{\!#1}}
\newcommand\Ae[1]{\gA^{\!#1}}

\newcommand\Snic[1]{$$\nds #1$$}
\newcommand\snic[1]{

{\centering$#1$\par}

}
\newcommand\snac[1]{

{\small\centering$#1$\par}

}
\newcommand\snuc[1]{

{\footnotesize\centering$#1$\par}
}

\newcommand\snucc[1]{

{\footnotesize$$\preskip0pt\postskip0pt\textstyle#1$$}}

\newcommand\snicc[1]{
{$$\preskip0pt\postskip0pt\textstyle#1$$}}

\newcommand\snif[3]{\vspace{#1}\noindent\centerline{$#3$}
\vspace{#2}}

\newcommand\env[2]{{{#2}_{#1}^{\mathrm{e}}}} 
\newcommand\Om[2]{\Omega_{{#2}/{#1}}}
\newcommand\Der[3]{{\rm Der}_{{#1}}({#2},{#3})}

\newcommand \isA[1] {_{#1/\!\gA}}
\newcommand\OmA[1]{\Omega\isA{#1}}

\newcommand \bra[1] {\left[{#1}\right]}
\newcommand \bu[1] {{{#1}\bul}}
\newcommand \ci[1] {{{#1}^\circ}}
\newcommand \wi[1] {\widetilde{#1}}
\newcommand \wh[1]{{\widehat{#1}}}
\newcommand \ov[1] {\overline{#1}}
\newcommand \und[1] {\underline{#1}}

\newcommand{\pref}[1]{\textup{\hbox{\normalfont(\ref{#1})}}}
\newcommand \vref[1] {\ref{#1}}
\newcommand \vpageref[1] {\paref{#1}}
\newcommand \thref[1] {\Thoz~\ref{#1}}
\newcommand \Thref[1] {\Thoz~\ref{#1}}
\newcommand \paref[1] {page~\pageref{#1}}
\newcommand \thrf[1] {\Thoz~\ref{#1}}
\newcommand \thrfs[2] {\Thosz~\ref{#1} and~\ref{#2}}
\renewcommand \rref[1] {\ref{#1} (\paref{#1})}
\newcommand \egrf[1] {\egtz~(\ref{#1})}
\newcommand \Egrf[1] {\Egtz~(\ref{#1})}
\newcommand \egref[1] {\egtz~(\ref{#1})  (\paref{#1})}
\newcommand \Egref[1] {\Egtz~(\ref{#1})  (\paref{#1})}
\newcommand \eqrf[1] {equation~(\ref{#1})}
\newcommand \Eqrf[1] {Equation~(\ref{#1})}
\newcommand \eqvrf[1] {equation (\ref{#1}) (\paref{#1})}
\newcommand \Eqvrf[1] {Equation (\ref{#1}) (\paref{#1})}
\newcommand \prirf[1] {principle~\ref{#1}}
\newcommand \plgrf[1] {\plgz~\ref{#1}}
\newcommand \plgref[1] {\plgz~\ref{#1}}
\newcommand \plgrfs[2] {\plgs \ref{#1} and~\ref{#2}}
\newcommand \Dfref[1] {\Dfnz~\ref{#1}}
\newcommand \Dfsref[2] {\Dfnsz~\ref{#1} and~\ref{#2}}

\newcommand\VRT[1]{\rotatebox{90}{\hbox{$#1$}}}
\newcommand\VRTsubseteq{\VRT{\subseteq}}
\newcommand\VRTsupseteq{\VRT{\supseteq}}
\newcommand\VRTlongrightarrow{\VRT{\longrightarrow}}
\newcommand\VRTlongleftarrow{\VRT{\longleftarrow}}

\newcommand \rC[1]{\MA{{\rm C}_{#1}}}
\newcommand \rF[1]{\MA{{\rm F}_{\!#1}}}
\newcommand \rR[1]{\MA{{\rm R}_{#1}}}
\newcommand \rRs[1]{\MA{\Rs_{#1}}}
\newcommand \ep[1]{^{(#1)}}

\newcommand \bal[1] {^\rK_{#1}}
\newcommand \ul[1] {_\rK^{#1}}
\newcommand \SNw[1] {P_{#1}}
\newcommand \gBtst {\gB[[t]]^{\!\times}}

\newcommand \THo[2]{\rdb
\mni{\bf Theorem {#1}~} {\it #2

}}

\newcommand \THO[2]{\rdb
\mni{\bf Theorem \ref{#1} bis~} {\it #2

}}

\newcommand{\llongrightarrow}{\relbar\joinrel\mkern-1mu\longrightarrow}
\newcommand{\lllongrightarrow}{\relbar\joinrel\mkern-1mu\llongrightarrow}
\newcommand{\llllongrightarrow}{\relbar\joinrel\mkern-1mu\lllongrightarrow}
\newcommand{\lllllongrightarrow}{\relbar\joinrel\mkern-1mu\llllongrightarrow}
\newcommand \lora {\longrightarrow}
\newcommand \llra {\llongrightarrow}
\newcommand \lllra {\lllongrightarrow}
\newcommand \llllra {\llllongrightarrow}
\newcommand \lllllra {\lllllongrightarrow}
\newcommand\simarrow{\vers{_\sim}}
\newcommand\isosim {\simarrow}
\newcommand\vers[1]{\buildrel{#1}\over \lora }
\newcommand\vvers[1]{\buildrel{#1}\over \llra }
\newcommand\vvvers[1]{\buildrel{#1}\over \lllra }
\newcommand\vvvvers[1]{\buildrel{#1}\over \llllra }
\newcommand\vvvvvers[1]{\buildrel{#1}\over \lllllra }

\newcommand \mt {\mapsto}
\newcommand \lmt {\longmapsto}

\renewcommand \le{\leqslant}
\renewcommand \leq{\leqslant}
\renewcommand \preceq{\preccurlyeq}
\renewcommand \ge{\geqslant}
\renewcommand \geq{\geqslant}
\renewcommand \succeq{\succurlyeq}

\newcommand \DeuxCol[2]{%
\sni\mbox{\parbox[t]{.475\textwidth}{#1}%
\hspace{.05\textwidth}%
\parbox[t]{.475\textwidth}{#2}}}

\newcommand \Deuxcol[4]{%
\sni\mbox{\parbox[t]{#1\textwidth}{#3}%
\hspace{.05\textwidth}%
\parbox[t]{#2\textwidth}{#4}}}


\floatstyle{boxed}
\floatname{agc}{Algorithm}
\newfloat{agc}{ht}{lag}[section]

\floatstyle{boxed}
\floatname{agC}{Algorithm}
\newfloat{agC}{H}{lag}[section]

\newenvironment{algor}[1][]%
{\par\smallskip\begin{agc}
\vskip 1mm
\begin{algorithm}{\bfseries#1}
\upshape\sffamily
}
{\end{algorithm}
\end{agc}
}

\newenvironment{algoR}[1][]%
{\par\smallskip\begin{agC}
\vskip 1mm
\begin{algorithm}{\bfseries#1}
\upshape\sffamily
}
{\end{algorithm}
\end{agC}
}

\newcommand\Vrai{\mathsf{True}}
\newcommand\Faux{\mathsf{False}}
\newcommand\ET{\mathsf{ and }}
\newcommand\OU{\mathsf{ or }}

\newcommand\pour[3]{\textbf{for } $#1$ \textbf{ from } $#2$
           \textbf{ to } $#3$ \textbf{ do }}
\newcommand\pur[2]{\textbf{for } $#1$ \textbf{ in } $#2$
            \textbf{ do }}
\newcommand\por[3]{\textbf{for } $#1$ \textbf{ from } $#2$
           \textbf{ to } $#3$  }
\def\sialors#1{\textbf{if } $#1$ \textbf{ then }}
\def\tantque#1{\textbf{while } $#1$ \textbf{ do }}
\newcommand\finpour{\textbf{end for}}
\newcommand\sinon{\textbf{else }}
\newcommand\finsi{\textbf{end if}}
\newcommand\fintantque{\textbf{end while}}
\newcommand\aff{\leftarrow }


\newcommand\Debut{\\[1mm] \textbf{Begin }}
\newcommand\Fin{\textbf{\\ End.}}
\newcommand\Entree{\\[1mm] \textbf{Input: }}
\newcommand\Sortie{\\ \textbf{Output: }}
\newcommand\Varloc{\\ \textbf{Local variables: }}
\newcommand\Repeter{\textbf{Repeat }}
\newcommand\jusqua{\textbf{until }}
\newcommand\hsz{\\ }
\newcommand\hsu{\\ \hspace*{4mm}}
\newcommand\hsd{\\ \hspace*{8mm}}
\newcommand\hst{\\ \hspace*{1,2cm}}
\newcommand\hsq{\\ \hspace*{1,6cm}}
\newcommand\hsc{\\ \hspace*{2cm}}
\newcommand\hsix{\\ \hspace*{2,4cm}}
\newcommand\hsept{\\ \hspace*{2,8cm}}

\newcommand \legendre {\overwithdelims()}
\newcommand \legendr[2] {\Big(\frac {#1}{#2} \Big)}
\newcommand \som {\sum\nolimits}
\newcommand \ds {\displaystyle}
\newcommand \nds {\textstyle}
\newcommand \Som {{\nds\sum\nolimits}}

\newcommand\et{\;\;\hbox{ and }\;\;}


\newcommand \divi {\mid}
\def \nedivi {\not\kern 2.5pt\mid}

\newcommand \vu {\vee} 
\newcommand \vi {\wedge} 
\newcommand \Vu {\bigvee\nolimits}
\newcommand \Vi {\bigwedge\nolimits}

\newcommand \vda {\,\vdash\,}
\newcommand \im {\rightarrow} 
\newcommand \dar[1] {\MA{\downarrow \!#1}}
\newcommand \uar[1] {\MA{\uparrow \!#1}}

\newcommand \Un {\mathbf{1}}
\newcommand \Deux {\mathbf{2}}
\newcommand \Trois {\mathbf{3}}
\newcommand \Quatre {\mathbf{4}}
\newcommand \Cinq {\mathbf{5}}


\newcommand \Pf {{\rm P}_{{\rm f}}}
\newcommand \Pfe {{\rm P}_{{\rm fe}}}

\newcommand \Ex {{\exists}}
\newcommand \Tt {{\forall}}

\newcommand \Lst{\mathsf{Lst}}
\newcommand \Irr{\mathsf{Irr}}
\newcommand \Prim{\mathsf{Prim}}
\newcommand \Rec{\mathsf{Rec}}

\newcommand \tsbf[1]{\textbf{\textsf{#1}}}
\newcommand \FAN{\tsbf{FAN} }
\newcommand \FANz{\tsbf{FAN}}
\newcommand \KL{\FAN}
\newcommand \KLz{\FANz}
\newcommand \KLp{\tsbf{KL}$_2$ }
\newcommand \KLpz{\tsbf{KL}$_2$}
\newcommand \kl{\tsbf{KL}$_1$ }
\newcommand \klz{\tsbf{KL}$_1$}
\newcommand \HC{\tsbf{HC} }
\newcommand \HCz{\tsbf{HC}}
\newcommand \LLPO{\tsbf{LLPO} }
\newcommand \LLPOz{\tsbf{LLPO}}
\newcommand \LPO{\tsbf{LPO} }
\newcommand \LPOz{\tsbf{LPO}}
\newcommand \MP{\tsbf{MP} }
\newcommand \MPz{\tsbf{MP}}
\newcommand \TEM{\tsbf{LEM} }
\newcommand \TEMz{\tsbf{LEM}}
\newcommand \UC{\tsbf{UC} }
\newcommand \UCz{\tsbf{UC}}
\newcommand \UCp{\tsbf{UC}$^+$ }
\newcommand \UCpz{\tsbf{UC}$^+$}
\newcommand \Mini{\tsbf{Min} }
\newcommand \Minip{\tsbf{Min}$^+$ }
\newcommand \Minipz{\tsbf{Min}$^+$}
\newcommand \Minim{\tsbf{Min}$^-$ }
\newcommand \Miniz{\tsbf{Min}}

\newcommand \bul{^\bullet}
\newcommand \eci{^\circ}
\newcommand \esh{^\sharp}
\newcommand \efl{^\flat}
\newcommand \epr{^\perp}
\newcommand \eti{^\times}
\newcommand \etl{^* }
\newcommand \eto{$^*\!$ }
\newcommand \etoz{$^*\!$}
\newcommand \sta{^\star}
\newcommand \ist{_\star}
\newcommand \eo {^{\mathrm{op}}}

\newcommand\Abul {\gA\!\bul}
\newcommand\Ati {\gA^{\!\times}}
\newcommand\Asta {\gA^{\!\star}}
\newcommand\Atl {\gA^{\!*}}
\newcommand\Bti {\gB^{\times}}
\newcommand\Bst {\Bti}
\newcommand\te  {\otimes}

\newcommand \iBA {_{\gB/\!\gA}}
\newcommand \iBk {_{\gB/\gk}}
\newcommand \iBK {_{\gB/\gK}}
\newcommand \iAk {_{\gA/\gk}}
\newcommand \iCk {_{\gC/\gk}}

\newcommand \tgaBG {\gB\{G\}}  
\newcommand \zcoho {Z^1(G, \Bti)}
\newcommand \bcoho {B^1(G, \Bti)}
\newcommand \hcoho {H^1(G, \Bti)}

\newcommand \vep{{\varepsilon}}

\newcommand\equidef{\buildrel{{\rm def}}\over{\;\Longleftrightarrow\;}}

\newcommand \E{\mathaccent19}    
\newcommand \aigu{\mathaccent19}    

\newcommand \noi {\noindent}
\renewcommand \ss {\smallskip}
\newcommand \sni {\ss\noi}
\newcommand \snii {}
\newcommand \ms {\medskip}
\newcommand \mni {\ms\noi}
\newcommand \bs {\bigskip}
\newcommand \bni {\bs\noi}
\newcommand \ce{\centering}
\newcommand \alb{\allowbreak}

\newcommand \ua  {{\underline{a}}}
\newcommand \ub  {{\underline{b}}}
\newcommand \ual {{\underline{\alpha}}}
\newcommand \ube {{\underline{\beta}}}
\newcommand \uc  {{\underline{c}}}
\newcommand \ud  {{\underline{d}}}
\newcommand \udel{{\underline{\delta}}}
\newcommand \ue  {{\underline{e}}}
\newcommand \uf  {{\underline{f}}}
\newcommand \uF  {{\underline{F}}}
\newcommand \ug  {{\underline{g}}}
\newcommand \uh  {{\underline{h}}}
\newcommand \uH  {{\underline{H}}}
\newcommand \uga {{\underline{\gamma}}}
\newcommand \ut  {{\underline{t}}}
\newcommand \uu  {{\underline{u}}}
\newcommand \ux  {{\underline{x}}}
\newcommand \uxi {{\underline{\xi}}}
\newcommand \uy  {{\underline{y}}}
\newcommand \uP  {{\underline{P}}}
\newcommand \uS  {{\underline{S}}}
\newcommand \us  {{\underline{s}}}
\newcommand \uT  {{\underline{T}}}
\newcommand \uU  {{\underline{U}}}
\newcommand \uX  {{\underline{X}}}
\newcommand \uY  {{\underline{Y}}}
\newcommand \uZ  {{\underline{Z}}}
\newcommand \uz  {{\underline{z}}}
\newcommand \uzeta  {{\underline{\zeta}}}
\newcommand \uze {{\underline{0}}}

\newcommand \am {a_1,\ldots,a_m}
\newcommand \an {a_1,\ldots,a_n}
\newcommand \bn {b_1,\ldots,b_n}
\newcommand \bbm {b_1,\ldots,b_m}
\newcommand \br {b_1,\ldots,b_r}
\newcommand \azn {a_0,\ldots,a_n}
\newcommand \bzn {b_0,\ldots,b_n}
\newcommand \czn {c_0,\ldots,c_n}
\newcommand \cq {c_1,\ldots,c_q}
\newcommand \gq {g_1,\ldots,g_q}
\newcommand \un {u_1,\ldots,u_n}
\newcommand \xk {x_1,\ldots,x_k}
\newcommand \Xk {X_1,\ldots,X_k}
\newcommand \xl {x_1,\ldots,x_\ell}
\newcommand \xm {x_1,\ldots,x_m}
\newcommand \xn {x_1,\ldots,x_n}
\newcommand \xp {x_1,\ldots,x_p}
\newcommand \yp {y_1,\ldots,y_p}
\newcommand \xq {x_1,\ldots,x_q}
\newcommand \xzk {x_0,\ldots,x_k}
\newcommand \xzn {x_0,\ldots,x_n}
\newcommand \xhn {x_0:\ldots:x_n}
\newcommand \Xn {X_1,\ldots,X_n}
\newcommand \Xzn {X_0,\ldots,X_n}
\newcommand \Xm {X_1,\ldots,X_m}
\newcommand \Xr {X_1,\ldots,X_r}
\newcommand \xr {x_1,\ldots,x_r}
\newcommand \Yr {Y_1,\ldots,Y_r}
\newcommand \Yn {Y_1,\ldots,Y_n}
\newcommand \ym {y_1,\ldots,y_m}
\newcommand \Ym {Y_1,\ldots,Y_m}
\newcommand \yk {y_1,\ldots,y_k}
\newcommand \yr {y_1,\ldots,y_r}
\newcommand \yn {y_1,\ldots,y_n}
\newcommand \zn {z_1,\ldots,z_n}
\newcommand \Zn {Z_1,\ldots,Z_n}

\newcommand \Sun {$S_1$, $\dots$, $S_n$ }
\newcommand \Sunz {$S_1$, $\dots$, $S_n$}

\newcommand \xpn {x'_1,\ldots,x'_n}
\newcommand \uxp  {{\underline{x'}}}
\newcommand \ypm {y'_1,\ldots,y'_m}
\newcommand \uyp  {{\underline{y'}}}

\newcommand \aln {\alpha_1,\ldots,\alpha_n}
\newcommand \gan {\gamma_1,\ldots,\gamma_n}
\newcommand \xin {\xi_1,\ldots,\xi_n}
\newcommand \xihn {\xi_0:\ldots:\xi_n}

\newcommand \lfs {f_1,\ldots,f_s}

\newcommand \Cin{C^{\infty}}
\newcommand \Ared {\gA\red}
\newcommand \Aqim {\gA\qim}
\newcommand \Amin {\Aqim}
\newcommand \Aqi {\gA_\mathrm{pp}}	

\newcommand \AT {{\gA[T]}}
\newcommand \AX {{\gA[X]}}
\newcommand \AY {{\gA[Y]}}
\newcommand \ArX {{\gA\lra X}}
\newcommand \ArY {{\gA\lra Y}}
\newcommand \BX {{\gB[X]}}
\newcommand \BY {{\gB[Y]}}
\newcommand \kX {{\gk[X]}}
\newcommand \kT {{\gk[T]}}
\newcommand \KT {{\gK[T]}}
\newcommand \KX {{\gK[X]}}
\newcommand \KY {{\gK[Y]}}
\newcommand \QQX {{\QQ[X]}}
\newcommand \VX {{\gV[X]}}
\newcommand \ZZX {{\ZZ[X]}}

\newcommand \AXY {\AuX \lra Y}

\newcommand \AXm {{\gA[\Xm]}}

\newcommand \kGa {{\gk[\Gamma]}}

\newcommand \kXn {{\gk[\Xn]}}
\newcommand \lXn {{\gl[\Xn]}}
\newcommand \AXk {{\gA[\Xk]}}
\newcommand \AXn {{\gA[\Xn]}}
\newcommand \BXn {{\gB[\Xn]}}
\newcommand \CXn {{\gC[\Xn]}}
\newcommand \KXn {{\gK[\Xn]}}
\newcommand \LXn {{\gL[\Xn]}}
\newcommand \RXn {{\gR[\Xn]}}
\newcommand \QQXk {{\QQ[\Xk]}}
\newcommand \QQXn {{\QQ[\Xn]}}
\newcommand \RXzn {{\gR[\Xzn]}}
\newcommand \Rxzn {{\gR[\xzn]}}
\newcommand \Ryn {{\gR[\yn]}}
\newcommand \VXn {{\gV[\Xn]}}
\newcommand \ZZXn {{\ZZ[\Xn]}}
\newcommand \ZZXk {{\ZZ[\Xk]}}

\newcommand \kXr {{\gk[\Xr]}}
\newcommand \lXr {{\gl[\Xr]}}
\newcommand \AXr {{\gA[\Xr]}}
\newcommand \KXr {{\gK[\Xr]}}
\newcommand \LXr {{\gL[\Xr]}}
\newcommand \RXr {{\gR[\Xr]}}
\newcommand \VXr {{\gV[\Xr]}}

\newcommand \kYr {{\gk[\Yr]}}
\newcommand \kyr {{\gk[\yr]}}
\newcommand \kyn {{\gk[\yn]}}
\newcommand \AYr {{\gA[\Yr]}}
\newcommand \KYr {{\gK[\Yr]}}
\newcommand \kXm {{\gk[\Xm]}}
\newcommand \KXm {{\gK[\Xm]}}
\newcommand \KYm {{\gK[\Ym]}}
\newcommand \kYm {{\gk[\Ym]}}
\newcommand \BYm {{\gB[\Ym]}}
\newcommand \lYr {{\gl[\Yr]}}
\newcommand \LYr {{\gL[\Yr]}}
\newcommand \RYr {{\gR[\Yr]}}

\newcommand \Axr {{\gA[\xr]}}
\newcommand \Kxr {{\gK[\xr]}}
\newcommand \kxr {{\gk[\xr]}}
\newcommand \Ayr {{\gA[\yr]}}
\newcommand \Kyr {{\gK[\yr]}}
\newcommand \Ayn {{\gA[\yn]}}
\newcommand \AYn {{\gA[\Yn]}}

\newcommand \kxm {{\gk[\xm]}}
\newcommand \kxn {{\gk[\xn]}}
\newcommand \lxn {{\gl[\xn]}}
\newcommand \Axn {{\gA[\xn]}}
\newcommand \Bxn {{\gB[\xn]}}
\newcommand \Cxn {{\gC[\xn]}}
\newcommand \Kxn {{\gK[\xn]}}
\newcommand \Kyn {{\gK[\yn]}}
\newcommand \Lxn {{\gL[\xn]}}
\newcommand \Rxn {{\gR[\xn]}}

\newcommand \Aux {{\gA[\ux]}}
\newcommand \Auy {{\gA[\uy]}}
\newcommand \Bux {{\gB[\ux]}}
\newcommand \Kuy {{\gK[\uy]}}
\newcommand \Kuu {{\gK[\uu]}}
\newcommand \Kux {{\gK[\ux]}}
\newcommand \kux {{\gk[\ux]}}
\newcommand \kuy {{\gk[\uy]}}
\newcommand \Rux {{\gR[\ux]}}
\newcommand \Ruy {{\gR[\uy]}}

\newcommand \ZZuS {{\ZZ[\uS]}}

\newcommand \AuX {{\gA[\uX]}}
\newcommand \BuX {{\gB[\uX]}}
\newcommand \KuX {{\gK[\uX]}}
\newcommand \kuX {{\gk[\uX]}}
\newcommand \LuX {{\gL[\uX]}}
\newcommand \RuX {{\gR[\uX]}}
\newcommand \ZZuX {{\ZZ[\uX]}}
\newcommand \QQuX {{\QQ[\uX]}}
\newcommand \AuY {{\gA[\uY]}}
\newcommand \BuY {{\gB[\uY]}}
\newcommand \KuY {{\gK[\uY]}}
\newcommand \kuY {{\gk[\uY]}}

\newcommand \Gn  {\gG_n}
\newcommand \Gnk {\gG_{n,k}}
\newcommand \Gnr {\gG_{n,r}}
\newcommand \cGn {{\cG_n}}
\newcommand \cGnk{{\cG_{n,k}}}

\newcommand \GGn {\GG_{n}}
\newcommand \GGnk{{\GG_{n,k}}} 
\newcommand \GGnr{{\GG_{n,r}}}
\newcommand \GA  {\mathbb{AG}}
\newcommand \GAn {\GA_{n}}  
\newcommand \GAq {\GA_{q}}
\newcommand \GAnk{\GA_{n,k}}
\newcommand \GAnr{\GA_{n,r}}

\newcommand \Mm {{\MM_{m}}}
\newcommand \Mn {{\MM_{n}}}
\newcommand \Mk {{\MM_{k}}}
\newcommand \Mq {{\MM_{q}}}
\newcommand \Mr {{\MM_{r}}}
\newcommand \MMn {{\MM_{n}}}

\newcommand \Bo{\BB\mathrm{o}}

\newcommand \GL {\mathbb{GL}}
\newcommand \GLn {{\GL_n}}
\newcommand \Gl {\mathbf{GL}}
\newcommand \Gln {{\Gl_n}}
\newcommand \PGL {\mathbb{PGL}}
\newcommand \SL {\mathbb{SL}}
\newcommand \SLn {{\SL_n}}
\newcommand \EE {\mathbb{E}}
\newcommand \En {\EE_n}
\newcommand \Pn {\PP^n}
\newcommand \An {\AA^{\!n}}
\newcommand \Am {\AA^{\!m}}
\newcommand \Sl {\mathbf{SL}}
\newcommand \Sln {{\Sl_n}}

\newcommand \I  {\mathrm{I}}
\newcommand \G  {\mathrm{G}}

\newcommand \rA {\mathrm{A}}
\newcommand \rD {\mathrm{D}}
\newcommand \rc {\mathrm{c}}
\newcommand \rE {\mathrm{E}}
\newcommand \rG {\mathrm{G}}
\newcommand \rH {\mathrm{H}}
\newcommand \rI {\mathrm{I}}
\newcommand \rJ {\mathrm{J}}
\newcommand \rK {\mathrm{K}}
\newcommand \rL {\mathrm{L}}
\newcommand \rM {\mathrm{M}}
\newcommand \rN {\mathrm{N}}
\newcommand \rP {\mathrm{P}}
\newcommand \rQ {\mathrm{Q}}
\newcommand \rS {\mathrm{S}}
\newcommand \rT {\mathrm{T}}
\newcommand \rZ {\mathrm{Z}}
\newcommand \rd {\mathrm{d}}

\renewcommand \AA{\mathbb{A}}
\newcommand \BB{\mathbb{B}}
\newcommand \CC{\mathbb{C}}
\newcommand \FF{\mathbb{F}}
\newcommand \FFp{\FF_p}
\newcommand \FFq{\FF_q}
\newcommand \GG{\mathbb{G}}
\newcommand \MM{\mathbb{M}}
\newcommand \NN{\mathbb{N}}
\newcommand \PP{\mathbb{P}}
\newcommand \QQ{\mathbb{Q}}
\newcommand \UU{\mathbb{U}}
\newcommand \ZZ{\mathbb{Z}}
\newcommand \ZB{\mathbb{ZB}}
\newcommand \RR{\mathbb{R}}
\newcommand \ASL {\mathbb{ASL}}
\newcommand \AGL {\mathbb{AGL}}

\newcommand \Z{\mathbb{Z}} 

\newcommand \cA {\mathcal{A}}
\newcommand \cB {\mathcal{B}}
\newcommand \cC {\mathcal{C}}
\newcommand \cD {\mathcal{D}}
\newcommand \cE {\mathcal{E}}
\newcommand \Diff {\mathcal{D}}
\newcommand \cG {\mathcal{G}}
\newcommand \cI {\mathcal{I}}
\newcommand \cJ {\mathcal{J}}
\newcommand \cF {\mathcal{F}}
\newcommand \cK {\mathcal{K}}
\newcommand \cL {\mathcal{L}}
\newcommand \cO {\mathcal{O}}
\newcommand \cP {\mathcal{P}}
\newcommand \cR {\mathcal{R}}
\newcommand \cM {\mathcal{M}}
\newcommand \cN {\mathcal{N}}
\newcommand \cS {\mathcal{S}}
\newcommand \cV {\mathcal{V}}
\newcommand \cZ {\mathcal{Z}}

\newcommand \SK {\cS^\rK}
\newcommand \IK {\cI^\rK}
\newcommand \JK {\cJ^\rK}
\newcommand \IH {\cI^\rH}
\newcommand \JH {\cJ^\rH}

\newcommand \ga{\mathbf{a}}
\newcommand \gb{\mathbf{b}}
\newcommand \gc{\mathbf{c}}
\newcommand \gh{\mathbf{h}}
\newcommand \gk{\mathbf{k}}
\newcommand \gl{\mathbf{l}}
\newcommand \gs{\mathbf{s}}
\newcommand \gv{\mathbf{v}}
\newcommand \gw{\mathbf{w}}
\newcommand \gA{\mathbf{A}}
\newcommand \gB{\mathbf{B}}
\newcommand \gC{\mathbf{C}}
\newcommand \gD{\mathbf{D}}
\newcommand \gE{\mathbf{E}}
\newcommand \gF{\mathbf{F}}
\newcommand \gG{\mathbf{G}}
\newcommand \gK{\mathbf{K}}
\newcommand \gL{\mathbf{L}}
\newcommand \gM{\mathbf{M}}
\newcommand \gQ{\mathbf{Q}}
\newcommand \gR{\mathbf{R}}
\newcommand \gS{\mathbf{S}}
\newcommand \gT{\mathbf{T}}
\newcommand \gU{\mathbf{U}}
\newcommand \gV{\mathbf{V}}
\newcommand \gx{\mathbf{x}}
\newcommand \gy{\mathbf{y}}
\newcommand \gX{\mathbf{X}}
\newcommand \gW{\mathbf{W}}
\newcommand \gZ{\mathbf{Z}}

\newcommand\fa{\mathfrak{a}}
\newcommand\fb{\mathfrak{b}}
\newcommand\fc{\mathfrak{c}}
\newcommand\fd{\mathfrak{d}}
\newcommand\fA{\mathfrak{A}}
\newcommand\fB{\mathfrak{B}}
\newcommand\fC{\mathfrak{C}}
\newcommand\fD{\mathfrak{D}}
\newcommand\fI{\mathfrak{i}}
\newcommand\fII{\mathfrak{I}}
\newcommand\fj{\mathfrak{j}}
\newcommand\fJ{\mathfrak{J}}
\newcommand\ff{\mathfrak{f}}
\newcommand\ffg{\mathfrak{g}}
\newcommand\fF{\mathfrak{F}}
\newcommand\fh{\mathfrak{h}}
\newcommand\fl{\mathfrak{l}}
\newcommand\fm{\mathfrak{m}}
\newcommand\fp{\mathfrak{p}}
\newcommand\fq{\mathfrak{q}}
\newcommand\fM{\mathfrak{M}}
\newcommand\fN{\mathfrak{N}}
\newcommand\fP{\mathfrak{P}}
\newcommand\fQ{\mathfrak{Q}}
\newcommand\fR{\mathfrak{R}}
\newcommand\fx{\mathfrak{x}}
\newcommand\fV{\mathfrak{V}}
\newcommand\fZ{\mathfrak{Z}}


\newcommand \Ig {\mathrm{Ig}}
\newcommand \Id {\mathrm{Id}}
\newcommand \In {{\rI_n}}
\newcommand \J {\MA{\mathrm{Jac}}}
\newcommand \JJ {\MA{\mathrm{JAC}}}
\newcommand \Sn {{\mathrm{S}_n}}

\newcommand \DA {\rD_{\!\gA}}
\newcommand \DB {\rD_\gB}
\newcommand \DV {\rD_\gV}
\newcommand \JA {\rJ_\gA}

\newcommand \red {_{\mathrm{red}}}
\newcommand \qim {_{\mathrm{min}}}
\newcommand \rja {\mathrm{Ja}}

\newcommand \ide {\mathrm{e}}

\newcommand \Adj {\MA{\mathrm{Adj}}}
\newcommand \adj {\MA{\mathrm{adj}}}
\newcommand \Adu {\MA{\mathrm{Adu}}}
\newcommand \Ann {\mathrm{Ann}}
\newcommand \Aut {\MA{\mathrm{Aut}}}
\newcommand \BZ {\MA{\mathrm{BZ}}}
\newcommand \car {\MA{\mathrm{char}}}
\newcommand \Cl {\MA{\mathrm{Cl}}}
\newcommand \Coker {\MA{\mathrm{Coker}}}
\renewcommand \det {\MA{\mathrm{det}}}
\renewcommand \deg {\MA{\mathrm{deg}}}
\newcommand \Diag {\MA{\mathrm{Diag}}}
\newcommand \di {\MA{\mathrm{di}}}
\newcommand \disc {\MA{\mathrm{disc}}}
\newcommand \Disc {\MA{\mathrm{Disc}}}
\newcommand \Div {\MA{\mathrm{Div}}}
\newcommand \dv {\MA{\mathrm{div}}}
\newcommand \ev {{\mathrm{ev}}}
\newcommand \End {\MA{\mathrm{End}}}
\newcommand \Fix {\MA{\mathrm{Fix}}}
\newcommand \Frac {\MA{\mathrm{Frac}}}
\newcommand \Gal {\MA{\mathrm{Gal}}}
\newcommand \Gfr {\MA{\mathrm{Gfr}}}
\newcommand \Gram {\MA{\mathrm{Gram}}}
\newcommand \gram {\MA{\mathrm{gram}}}
\newcommand \hauteur {\mathrm{height}}
\newcommand \Ker {\MA{\mathrm{Ker}}}
\renewcommand \ker {\Ker}
\newcommand \Hom {\MA{\mathrm{Hom}}}
\newcommand \Ifr {\MA{\mathrm{Ifr}}}
\newcommand \Icl {\MA{\mathrm{Icl}}}
\renewcommand \Im {\MA{\mathrm{Im}}}
\newcommand \Lin {\mathrm{L}}
\newcommand \LIN {\mathrm{Lin}}
\newcommand \Mat {\MA{\mathrm{Mat}}}
\newcommand \Mip {\mathrm{Min}}
\newcommand \md {\mathrm{lm}} 
\newcommand \mod {\;\mathrm{mod}\;}
\newcommand \Mor {\MA{\mathrm{Mor}}}
\newcommand \poles {\hbox {\rm p\^oles}}
\newcommand \pgcd {\MA{\mathrm{gcd}}}
\newcommand \ppcm {\MA{\mathrm{lcm}}}
\newcommand \Rad {\MA{\mathrm{Rad}}}
\newcommand \Reg {\MA{\mathrm{Reg}}}
\newcommand \rg{\MA{\mathrm{rk}}} 
\newcommand \Res {\mathrm{Res}}
\newcommand \Rs {\MA{\mathrm{Rs}}}
\newcommand \rPr{\MA{\mathrm{Pr}}}
\newcommand \Rv {\mathrm{Rv}}
\newcommand \Syl {\mathrm{Syl}}
\newcommand \Stp {\MA{\mathrm{Stp}}}
\newcommand \St {\mathrm{St}}
\newcommand \Tri {\MA{\mathrm{Tri}}}
\newcommand \Tor {\MA{\mathrm{Tor}}}
\newcommand \tr {\MA{\mathrm{tr}}}
\newcommand \Tr {\MA{\mathrm{Tr}}}
\newcommand \Tsc {\MA{\mathrm{Tsch}}}
\newcommand \Um {\MA{\mathrm{Um}}}
\newcommand \val {\MA{\mathrm{val}}}
\newcommand \Suslin{{\rm Suslin}}

\newcommand \rgl {\rg^\lambda}
\newcommand \rgg {\rg^\gamma}


\newcommand \sfP {\mathsf{P}}
\newcommand \sfC {\mathsf{C}}

\newcommand \GK {\MA{\mathsf{GK}}}
\newcommand \GKO {\MA{\mathsf{GK}_0}}
\newcommand \Gr {\MA{\mathsf{Gr}}}
\newcommand \HO {\MA{\mathsf{H}_0}}
\newcommand \HOp {\MA{\mathsf{H}_0^+}}
\newcommand \Hdim {\MA{\mathsf{Hdim}}}
\newcommand \HeA {{\Heit\gA}}
\newcommand \Heit {\MA{\mathsf{Heit}}}
\newcommand \Hspec {\MA{\mathsf{Hspec}}}
\newcommand \Jdim {\MA{\mathsf{Jdim}}}
\newcommand \jdim {\MA{\mathsf{jdim}}}
\newcommand \Jspec {\MA{\mathsf{Jspec}}}
\newcommand \jspec {\MA{\mathsf{jspec}}}
\newcommand \KO {\MA{\mathsf{K}_0}}
\newcommand \KOp {\MA{\mathsf{K}_0^+}}
\newcommand \KTO {\MA{\wi{\mathsf{K}}_0}}
\newcommand \Kdim {\MA{\mathsf{Kdim}}}
\newcommand \Max {\MA{\mathsf{Max}}}
\newcommand \Min {\MA{\mathsf{Min}}}
\newcommand \OQC {\MA{\mathsf{Oqc}}}
\newcommand \Pic {\MA{\mathsf{Pic}}}
\newcommand \Spec {\MA{\mathsf{Spec}}}
\newcommand \SpecA {\Spec\gA}
\newcommand \SpecT {\Spec\gT}
\newcommand \Vdim {\MA{\mathsf{Vdim}}}
\newcommand \Zar {\MA{\mathsf{Zar}}}
\newcommand \ZF {\MA{\mathsf{ZF}}}
\newcommand \ZarA {{\Zar\gA}}


\newcommand \fRes {\MA{\mathfrak{Res}}}

\newcommand \SPEC {\MA{\mathfrak{spec}}}
\newcommand \SPECK {\SPEC_\gK}
\newcommand \SPECk {\SPEC_\gk}

\newcommand \Ap {{\gA_\fp}}
\newcommand \zg {\ZZ[G]}

\newcommand \num {{no.~}}

\newcommand{\bma}{\bm{a}}
\newcommand{\bmb}{\bm{b}}
\newcommand{\bmc}{\bm{c}}
\newcommand{\bmd}{\bm{d}}
\newcommand{\bme}{\bm{e}}
\newcommand{\bmf}{\bm{f}}
\newcommand{\bmu}{\bm{u}}
\newcommand{\bmv}{\bm{v}}
\newcommand{\bmw}{\bm{w}}
\newcommand{\bmy}{\bm{y}}
\newcommand{\bmx}{\bm{x}}
\newcommand{\bmz}{\bm{z}}

\newcommand\comm{\rdb
\noi{\it Comment. }}

\newcommand\COM[1]{\rdb
\noi{\it Comment #1. }}

\newcommand\comms{\rdb
\noi{\it Comments. }}

\newcommand\rem{\rdb
\noi{\it Remark. }}

\newcommand\REM[1]{\rdb
\noi{\it #1 remark. }}

\newcommand\rems{\rdb
\noi{\it Remarks. }}

\newcommand\exl{\rdb
\noi{\bf Example. }}

\newcommand\EXL[1]{\rdb
\noi{\bf Example: #1. }}

\newcommand\exls{\rdb
\noi{\bf Examples. }}

\newcommand\Pb{\rdb
\noi{\bf Problem. }}

\newcommand\PB[1]{\rdb
\noi{\bf Problem #1. }}

\newcommand \eoe {\hbox{}\nobreak\hfill
\vrule width 1.4mm height 1.4mm depth 0mm \par }

\newcommand\eoq{\hbox{}\nobreak
\vrule width 1.4mm height 1.4mm depth 0mm}

\newcommand \ndsp{\textstyle}

\newcommand \recu {induction }
\newcommand \recuz {induction}
\newcommand \hdr {\recu hypothesis }
\newcommand \hdrz {\recu hypothesis}

\newcommand \cad {i.e.\ }
\newcommand \Cad {I.e.\ }
\newcommand \cadz {i.e.}

\newcommand \ssi {if and only if }
\newcommand \ssiz {if and only if~}

\newcommand \cnes {nec\-es\-sa\-ry and sufficient condition }

\newcommand \spdg {without loss of generality }
\newcommand \spdgz {without loss of generality}
\newcommand \Spdg {Without loss of generality }
\newcommand \Spdgz {Without loss of generality}

\newcommand \Propeq {The following properties are equivalent.}
\newcommand \propeq {the following properties are equivalent.}

\newcommand \Kev {$\gK$-\evc }
\newcommand \Kevs {$\gK$-\evcs }
\newcommand \Kevz {$\gK$-\evcz}
\newcommand \Kevsz {$\gK$-\evcsz}

\newcommand \Ksv {$\gK$-linear subspace }
\newcommand \Ksvs {$\gK$-linear subspaces }
\newcommand \Ksvz {$\gK$-linear subspace}
\newcommand \Ksvsz {$\gK$-linear subspaces}

\newcommand \Lev {$\gL$-\evc }
\newcommand \Levs {$\gL$-\evcs }
\newcommand \Levz {$\gL$-\evcz}
\newcommand \Levsz {$\gL$-\evcsz}

\newcommand \Qev {$\QQ$-\evc }
\newcommand \Qevs {$\QQ$-\evcs }
\newcommand \Qevz {$\QQ$-\evcz}
\newcommand \Qevsz {$\QQ$-\evcsz}

\newcommand \Qsev {$\QQ$-\sevc }
\newcommand \Qsevs {$\QQ$-\sevcs }
\newcommand \Qsevz {$\QQ$-\sevcz}
\newcommand \Qsevsz {$\QQ$-\sevcsz}

\newcommand \kev {$\gk$-\evc }
\newcommand \kevs {$\gk$-\evcs }
\newcommand \kevz {$\gk$-\evcz}
\newcommand \kevsz {$\gk$-\evcsz}

\newcommand \lev {$\gl$-\evc }
\newcommand \levs {$\gl$-\evcs }
\newcommand \levz {$\gl$-\evcz}
\newcommand \levsz {$\gl$-\evcsz}

\newcommand \Alg {$\gA$-\alg}
\newcommand \Algs {$\gA$-\algs}
\newcommand \Algz {$\gA$-\algz}
\newcommand \Algsz {$\gA$-\algsz}

\newcommand \Aslg {$\gA$-sub\alg}
\newcommand \Aslgs {$\gA$-sub\algs}
\newcommand \Aslgz {$\gA$-sub\algz}
\newcommand \Aslgsz {$\gA$-sub\algsz}

\newcommand \Blg {$\gB$-\alg}
\newcommand \Blgs {$\gB$-\algs}
\newcommand \Blgz {$\gB$-\algz}
\newcommand \Blgsz {$\gB$-\algsz}

\newcommand \Clg {$\gC$-\alg}
\newcommand \Clgs {$\gC$-\algs}
\newcommand \Clgz {$\gC$-\algz}
\newcommand \Clgsz {$\gC$-\algsz}

\newcommand \klg {$\gk$-\alg}
\newcommand \klgs {$\gk$-\algs}
\newcommand \klgz {$\gk$-\algz}
\newcommand \klgsz {$\gk$-\algsz}

\newcommand \kslg {$\gk$-sub\alg}
\newcommand \kslgs {$\gk$-sub\algs}
\newcommand \kslgz {$\gk$-sub\algz}
\newcommand \kslgsz {$\gk$-sub\algsz}

\newcommand \llg {$\gl$-\alg}
\newcommand \llgs {$\gl$-\algs}
\newcommand \llgz {$\gl$-\algz}
\newcommand \llgsz {$\gl$-\algsz}

\newcommand \Klg {$\gK$-\alg}
\newcommand \Klgs {$\gK$-\algs}
\newcommand \Klgz {$\gK$-\algz}
\newcommand \Klgsz {$\gK$-\algsz}

\newcommand \Kslg {$\gK$-sub\alg}
\newcommand \Kslgs {$\gK$-sub\algs}
\newcommand \Kslgz {$\gK$-sub\algz}
\newcommand \Kslgsz {$\gK$-sub\algsz}

\newcommand \Llg {$\gL$-\alg}
\newcommand \Llgs {$\gL$-\algs}
\newcommand \Llgz {$\gL$-\algz}
\newcommand \Llgsz {$\gL$-\algsz}

\newcommand \QQlg {$\QQ$-\alg}
\newcommand \QQlgs {$\QQ$-\algs}
\newcommand \QQlgz {$\QQ$-\algz}
\newcommand \QQlgsz {$\QQ$-\algsz}

\newcommand \Rlg {$\gR$-\alg}
\newcommand \Rlgs {$\gR$-\algs}
\newcommand \Rlgz {$\gR$-\algz}
\newcommand \Rlgsz {$\gR$-\algsz}

\newcommand \RRlg {$\RR$-\alg}
\newcommand \RRlgs {$\RR$-\algs}
\newcommand \RRlgz {$\RR$-\algz}
\newcommand \RRlgsz {$\RR$-\algsz}

\newcommand \ZZlg {$\ZZ$-\alg}
\newcommand \ZZlgs {$\ZZ$-\algs}
\newcommand \ZZlgz {$\ZZ$-\algz}
\newcommand \ZZlgsz {$\ZZ$-\algsz}

\newcommand \Amo {$\gA$-module }
\newcommand \Amos {$\gA$-modules }
\newcommand \Amoz {$\gA$-module}
\newcommand \Amosz {$\gA$-modules}

\newcommand \Asub {$\gA$-submodule }
\newcommand \Asubs {$\gA$-submodules }
\newcommand \Asubz {$\gA$-submodule}
\newcommand \Asubsz {$\gA$-submodules}

\newcommand \Amrc {\pro \Amo of constant rank }
\newcommand \Amrcs {\pro \Amos of constant rank }
\newcommand \Amrcz {\pro \Amo of constant rank}
\newcommand \Amrcsz {\pro \Amos of constant rank}

\newcommand \kmrc {\pro \kmo of constant rank }
\newcommand \kmrcs {\pro \kmos of constant rank }
\newcommand \kmrcz {\pro \kmo of constant rank}
\newcommand \kmrcsz {\pro \kmos of constant rank}

\newcommand \Bmo {$\gB$-module }
\newcommand \Bmos {$\gB$-modules }
\newcommand \Bmoz {$\gB$-module}
\newcommand \Bmosz {$\gB$-modules}
\newcommand \Bsmo {$\gB$-submodule }
\newcommand \Bsmos {$\gB$-submodules }
\newcommand \Bsmoz {$\gB$-submodule}
\newcommand \Bsmosz {$\gB$-submodules}

\newcommand \Cmo {$\gC$-module }
\newcommand \Cmos {$\gC$-modules }
\newcommand \Cmoz {$\gC$-module}
\newcommand \Cmosz {$\gC$-modules}

\newcommand \kmo {$\gk$-module }
\newcommand \kmos {$\gk$-modules }
\newcommand \kmoz {$\gk$-module}
\newcommand \kmosz {$\gk$-modules}

\newcommand \ksmo {$\gk$-submodule }
\newcommand \ksmos {$\gk$-submodules }
\newcommand \ksmoz {$\gk$-submodule}
\newcommand \ksmosz {$\gk$-submodules}

\newcommand \Kmo {$\gK$-module }
\newcommand \Kmos {$\gK$-modules }
\newcommand \Kmoz {$\gK$-module}
\newcommand \Kmosz {$\gK$-modules}

\newcommand \Lmo {$\gL$-module }
\newcommand \Lmos {$\gL$-modules }
\newcommand \Lmoz {$\gL$-module}
\newcommand \Lmosz {$\gL$-modules}

\newcommand \Ali {$\gA$-\lin map }
\newcommand \Alis {$\gA$-\lin maps }
\newcommand \Aliz {$\gA$-\lin map}
\newcommand \Alisz {$\gA$-\lin maps}

\newcommand \Bli {$\gB$-\lin map }
\newcommand \Blis {$\gB$-\lin maps }
\newcommand \Bliz {$\gB$-\lin map}
\newcommand \Blisz {$\gB$-\lin maps}

\newcommand \Cli {$\gC$-\lin map }
\newcommand \Clis {$\gC$-\lins maps }
\newcommand \Cliz {$\gC$-\lin map}
\newcommand \Clisz {$\gC$-\lin maps}

\newcommand \kli {$\gk$-\lin map }
\newcommand \klis {$\gk$-\lin maps }
\newcommand \kliz {$\gk$-\lin map}
\newcommand \klisz {$\gk$-\lin maps}

\newcommand \Kli {$\gK$-\lin map }
\newcommand \Klis {$\gK$-\lin maps }
\newcommand \Kliz {$\gK$-\lin map}
\newcommand \Klisz {$\gK$-\lin maps}

\newcommand \ZZmo {$\ZZ$-module }
\newcommand \ZZmos {$\ZZ$-modules }
\newcommand \ZZmoz {$\ZZ$-module}
\newcommand \ZZmosz {$\ZZ$-modules}

\newcommand \ZZsmo {$\ZZ$-submodule }
\newcommand \ZZsmos {$\ZZ$-submodules }
\newcommand \ZZsmoz {$\ZZ$-submodule}
\newcommand \ZZsmosz {$\ZZ$-submodules}

\newcommand \ZZx {{\ZZ[x]}}
\newcommand \Kx {{\gK[x]}}

\newcommand \ac{\agqt closed }
\newcommand \acz{\agqt closed}

\newcommand \acl {\icl \ri}
\newcommand \acls {\icl \ris}
\newcommand \aclsz {\icl \risz}
\newcommand \aclz {\icl \riz}

\newcommand \adk {Dedekind \ri}
\newcommand \adks {Dedekind \ris}
\newcommand \adkz {Dedekind \riz}
\newcommand \adksz {Dedekind \risz}

\newcommand \dDk {Dedekind domain }
\newcommand \dDks {Dedekind domains }
\newcommand \dDkz {Dedekind domain}
\newcommand \dDksz {Dedekind domains}

\newcommand \adp {Pr\"ufer \ri}
\newcommand \adps {Pr\"ufer \ris}
\newcommand \adpsz {Pr\"ufer \risz}
\newcommand \adpz {Pr\"ufer \riz}

\newcommand \adpc {\coh \adp }
\newcommand \adpcs {\cohs \adps }
\newcommand \adpcz {\coh \adpz}
\newcommand \adpcsz {\cohs \adpsz}

\newcommand \adu {\uvl splitting \alg} 
\newcommand \adus {\uvl splitting \algs} 
\newcommand \aduz {\uvl splitting \algz}
\newcommand \adusz {\uvl splitting \algsz}

\newcommand \adv {valuation \ri}
\newcommand \advs {valuation \ris}
\newcommand \advsz {valuation \risz}
\newcommand \advz {valuation \riz}

\newcommand \sdv {valuation sub\ri}
\newcommand \sdvs {valuation sub\ris}
\newcommand \sdvsz {valuation sub\risz}
\newcommand \sdvz {valuation sub\riz}

\newcommand \aG {Galois \alg }
\newcommand \aGs {Galois \algs }
\newcommand \aGz {Galois \algz}
\newcommand \aGsz {Galois \algsz}

\newcommand \agB {Boolean \alg }
\newcommand \agBs {Boolean \algs }
\newcommand \agBz {Boolean \algz}
\newcommand \agBsz {Boolean \algsz}

\newcommand \agsB {Boolean sub\alg }
\newcommand \agsBs {Boolean sub\algs }
\newcommand \agsBz {Boolean sub\algz}
\newcommand \agsBsz {Boolean sub\algsz}

\newcommand \agH {Heyting \alg }
\newcommand \agHs {Heyting \algs }
\newcommand \agHz {Heyting \algz}
\newcommand \agHsz {Heyting \algsz}

\newcommand \agq{algebraic }
\newcommand \agqs{algebraic }
\newcommand \agqz{algebraic}
\newcommand \agqsz{algebraic}

\newcommand \Agq{Algebraic }

\newcommand \agqt{algebraically }
\newcommand \Agqt{Algebraically }
\newcommand \Agqtz{Algebraically}

\newcommand \alg {algebra }
\newcommand \algs {algebras }
\newcommand \algz {algebra}
\newcommand \algsz {algebras}

\newcommand \algo{algo\-rithm }
\newcommand \algos{algo\-rithms }
\newcommand \algoz{algo\-rithm}
\newcommand \algosz{algo\-rithms}
\newcommand \Algo{Algorithm }
\newcommand \Algoz{Algorithm}

\newcommand \algq{algo\-rithmic }
\newcommand \algqs{algo\-rithmic }
\newcommand \algqz{algo\-rithmic}
\newcommand \algqsz{algo\-rithmic}

\newcommand \algqt{algo\-rithmi\-cally }

\newcommand \agsp {\spl \agq}
\newcommand \agsps {\spls \agqs}
\newcommand \agspz {\spl \agqz}
\newcommand \agspsz {\spls \agqsz}

\newcommand \ali {\lin map }
\newcommand \alis {\lin maps }
\newcommand \aliz {\lin map}
\newcommand \alisz {\lin maps}

\newcommand \alo {local \ri}
\newcommand \aloz {local \riz}
\newcommand \alos {local \ris}
\newcommand \alosz {local \risz}
\newcommand \Alo {Local \ri}
\newcommand \Alos {Local \ris}
\newcommand \Aloz {Local \riz}
\newcommand \Alosz {Local \risz}

\newcommand \algb {\lgb \ri}
\newcommand \algbs {\lgb \ris}
\newcommand \algbz {\lgb \riz}
\newcommand \algbsz {\lgb \risz}
\newcommand \Algb {\Lgb \ri}
\newcommand \Algbs {\Lgb \ris}

\newcommand \alrd {\dcd \alo }
\newcommand \alrds {\dcds \alos }
\newcommand \alrdz {\dcd \aloz}
\newcommand \alrdsz {\dcds \alosz}

\newcommand \anar {\ari \ri}
\newcommand \anars {\ari \ris}
\newcommand \anarsz {\ari \risz}
\newcommand \anarz {\ari \riz}
\newcommand \Anars {\Ari \ris}
\newcommand \Anarsz {\Ari \risz}

\newcommand \anor {normal \ri}
\newcommand \anors {normal \ris}
\newcommand \anorsz {normal \risz}
\newcommand \anorz {normal \riz}

\newcommand \apf {\pf \alg}
\newcommand \apfs {\pf \algs}
\newcommand \apfz {\pf \algz}
\newcommand \apfsz {\pf \algsz}

\newcommand \pGa {pre-Galois } 
\newcommand \apG {\pGa\alg} 
\newcommand \apGs {\pGa\algs}
\newcommand \apGz {\pGa\algz}
\newcommand \apGsz {\pGa\algsz}

\newcommand \areg{regular map }
\newcommand \aregs{regular maps }
\newcommand \aregz{regular map}
\newcommand \aregsz{regular maps}

\newcommand \Ari{Arithmetic }
\newcommand \ari{arith\-metic }
\newcommand \ariz{arith\-metic}
\newcommand \aris{arith\-metic }
\newcommand \arisz{arith\-metic}

\newcommand \asfs {\stfes \algs}
\newcommand \asfsz {\stfes \algsz}
\newcommand \asf {\stfe \alg}
\newcommand \asfz {\stfe \algz}

\newcommand \ase {\ste \alg}
\newcommand \ases {\stes \algs}
\newcommand \asez {\ste \algz}
\newcommand \asesz {\stes \algsz}

\newcommand \atf {\tf \alg}
\newcommand \atfs {\tf \algs}
\newcommand \atfz {\tf \algz}
\newcommand \atfsz {\tf \algsz}

\newcommand \auto {auto\-morphism }
\newcommand \autos {auto\-morphisms }
\newcommand \autoz {auto\-morphism}
\newcommand \autosz {auto\-morphisms}


\newcommand \bdp {partial \dcn basis }
\newcommand \bdpz {partial \dcn basis}
\newcommand \bdf {\fap basis }
\newcommand \bdfz {\fap basis}

\newcommand \bb {bien bon } 
\newcommand \bbs {bien bons }
\newcommand \bbz {bien bon}

\newcommand \bdg {Gr\"obner basis }
\newcommand \bdgs {Gr\"obner bases }
\newcommand \bdgz {Gr\"obner basis}
\newcommand \bdgsz {Gr\"obner bases}

\newcommand \cac{\ac field }
\newcommand \cacz{\ac field}
\newcommand \cacs{\ac fields }
\newcommand \cacsz{\ac fields}

\newcommand \calf{Computer Algebra }
\newcommand \calfz{Computer Algebra}

\newcommand \cara{char\-ac\-ter\-is\-tic }
\newcommand \caras{char\-ac\-ter\-is\-tic }
\newcommand \caraz{char\-ac\-ter\-is\-tic}
\newcommand \carasz{char\-ac\-ter\-is\-tic}
\newcommand \Cara{Characteristic }
\newcommand \Caras{Characteristics }

\newcommand \carn{char\-ac\-te\-ri\-za\-tion }
\newcommand \carns{char\-ac\-te\-ri\-za\-tions }
\newcommand \carnz{char\-ac\-te\-ri\-za\-tion}
\newcommand \Carn{Characterization }
\newcommand \Carns{Characterizations }

\newcommand \carar{char\-ac\-ter\-ize }

\newcommand \care{char\-ac\-ter\-ized }
\newcommand \caree{char\-ac\-ter\-ized }
\newcommand \cares{char\-ac\-ter\-ized }
\newcommand \carees{char\-ac\-ter\-ized }

\newcommand \carf{finite char\-ac\-ter }
\newcommand \carfz{finite char\-ac\-ter}

\newcommand \cdi{dis\-crete field }
\newcommand \cdis{dis\-crete fields }
\newcommand \cdiz{dis\-crete field}
\newcommand \cdisz{dis\-crete fields}

\newcommand \cdac{\ac \cdi}
\newcommand \cdacs{\ac \cdis \ac}
\newcommand \cdacz{\ac \cdiz}
\newcommand \cdacsz{\ac \cdisz}

\newcommand \cdr{splitting field }
\newcommand \cdrz{splitting field}
\newcommand \cdrs{splitting fields }
\newcommand \cdrsz{splitting fields}

\newcommand \cdv{change of variables } 
\newcommand \cdvs{changes of variables } 
\newcommand \cdvz{change of variables} 
\newcommand \cdvsz{changes of variables} 
\newcommand \Cdv{Change of variables } 
\newcommand \Cdvs{Changes of variables } 
\newcommand \Cdvz{Change of variables} 
\newcommand \Cdvsz{Changes of variables} 

\newcommand \cli{integral closure }
\newcommand \cliz{integral closure}

\newcommand \coe {co\-ef\-fi\-cient }
\newcommand \coes {co\-ef\-fi\-cients }
\newcommand \coez {co\-ef\-fi\-cient}
\newcommand \coesz {co\-ef\-fi\-cients}

\newcommand \coh {coherent }
\newcommand \cohs {coherent }
\newcommand \cohz {coherent}
\newcommand \cohsz {coherent}
\newcommand \Coh {Coherent }

\newcommand \cohc {coherence }
\newcommand \cohcz {coherence}
\newcommand \Cohc {Coherence }

\newcommand \cori {\coh \ri}
\newcommand \coriz {\coh \riz}
\newcommand \como {\coh module }
\newcommand \comoz {\coh module}
\newcommand \coris {\coh \ris}
\newcommand \corisz {\coh \risz}
\newcommand \comos {\coh modules }
\newcommand \comosz {\coh modules}

\newcommand \coli {\lin combination }
\newcommand \colis {\lin combinations }
\newcommand \coliz {\lin combination}
\newcommand \colisz {\lin combinations}

\newcommand \com {comaximal }
\newcommand \comz {comaximal}
\newcommand \come {comaximal }
\newcommand \comez {comaximal}

\newcommand \comt {comaximality }
\newcommand \comtz {comaximality}
\newcommand \Comt {Comaximality }

\newcommand \cpb {compatible } 
\newcommand \cpbs {compatible } 
\newcommand \cpbz {compatible} 
\newcommand \cpbsz {compatible} 

\newcommand \cpbt {compatibility } 
\newcommand \cpbtz {compatibility}

\newcommand \coo {coordinate }
\newcommand \coos {coordinates }
\newcommand \cooz {coordinate}
\newcommand \coosz {coordinates}

\newcommand \cop {complementary }
\newcommand \cops {complementary }
\newcommand \copz {complementary}
\newcommand \copsz {complementary}

\newcommand \cor {coregular }
\newcommand \core {coregular }
\newcommand \corz {coregular}
\newcommand \corez {coregular}

\newcommand \crc{character }
\newcommand \crcs{characters }
\newcommand \crcz{character}
\newcommand \crcsz{characters}

\newcommand \cro {crossed }
\newcommand \cros {crossed }
\newcommand \croz {crossed}
\newcommand \crosz {crossed}


\newcommand \dcd {\rdt dis\-crete }
\newcommand \dcdz {\rdt dis\-crete}
\newcommand \dcds {\rdt dis\-crete }
\newcommand \dcdsz {\rdt dis\-crete}

\newcommand \dcn {decomposition }
\newcommand \dcns {decompositions }
\newcommand \dcnz {decomposition}
\newcommand \dcnsz {decompositions}

\newcommand \dcnb {bounded \dcn }
\newcommand \dcnbz {bounded \dcnz}
\newcommand \dcnbs {bounded \dcns }
\newcommand \dcnbsz {bounded \dcnsz}

\newcommand \dcnc {complete \dcn}
\newcommand \dcncz {complete \dcnz}
\newcommand \dcncs {complete \dcns}
\newcommand \dcncsz {complete \dcnsz}

\newcommand \dcnp {partial \dcn}
\newcommand \dcnps {partial \dcns}
\newcommand \dcnpz {partial \dcnz}
\newcommand \dcnpsz {partial \dcnsz}

\newcommand \dcp {decomposable }
\newcommand \dcps {decomposable }
\newcommand \dcpz {decomposable}
\newcommand \dcpsz {decomposable}
\newcommand \Dcp {Decomposable }

\newcommand \ddk {Krull dimension }
\newcommand \ddkz {Krull dimension}
\newcommand \ddi {of dimension at most~}
\newcommand \ddiz {of dimension at most}

\newcommand \ddp {Pr\"ufer domain }
\newcommand \ddps {Pr\"ufer domains }
\newcommand \ddpz {Pr\"ufer domain}
\newcommand \ddpsz {Pr\"ufer domains}

\newcommand \ddv {valuation domain }
\newcommand \ddvz {valuation domain}
\newcommand \ddvs {valuation domains }
\newcommand \ddvsz {valuation domains}

\newcommand \dgn {de\-gen\-er\-ate } 
\newcommand \dgne {de\-gen\-er\-ate }
\newcommand \dgns {de\-gen\-er\-ate }
\newcommand \dgnes {de\-gen\-er\-ate }
\newcommand \dgnz {de\-gen\-er\-ate}
\newcommand \dgnez {de\-gen\-er\-ate}
\newcommand \dgnsz {de\-gen\-er\-ate}
\newcommand \dgnesz {de\-gen\-er\-ate}

\newcommand \Demo{Proof }

\newcommand \dem{proof }
\newcommand \demz{proof}
\newcommand \dems{proofs }
\newcommand \demsz{proofs}

\newcommand \denb{countable }
\newcommand \denbs{countable }
\newcommand \denbz{countable}
\newcommand \denbsz{countable}

\newcommand \deno{denominator }
\newcommand \denos{denominators }
\newcommand \denoz{denominator}
\newcommand \denosz{denominators}

\newcommand \deter {de\-ter\-mi\-nant }
\newcommand \deters {de\-ter\-mi\-nants }
\newcommand \deterz {de\-ter\-mi\-nant}
\newcommand \detersz {de\-ter\-mi\-nants}

\newcommand \dig{diagonalizable }
\newcommand \digs{diagonalizable }
\newcommand \digz{diagonalizable}
\newcommand \digsz{diagonalizable}

\newcommand \Dil{Differential }
\newcommand \dil{differential } 
\newcommand \dils{differentials }
\newcommand \dile{differential }
\newcommand \diles{differentials }
\newcommand \dilz{differential}
\newcommand \dilsz{differentials}
\newcommand \dilez{differential}
\newcommand \dilesz{differentials}

\newcommand \din{diagonalization }
\newcommand \dinz{diagonalization}
\newcommand \dins{diagonalizations }
\newcommand \dinsz{diagonalizations}

\newcommand \dit{distributivity } 
\newcommand \ditz{distributivity}

\newcommand \dfn{definition }
\newcommand \dfns{definitions }
\newcommand \dfnz{definition}
\newcommand \dfnsz{definitions}

\newcommand \Dfn{Definition }
\newcommand \Dfns{Definitions }
\newcommand \Dfnz{Definition}
\newcommand \Dfnsz{Definitions}

\newcommand \dlg{enlargement } 
\newcommand \dlgz{enlargement}
\newcommand \Dlg{Enlargement }

\newcommand \Dmo {$\gD$-module }
\newcommand \Dmos {$\gD$-modules }
\newcommand \Dmoz {$\gD$-module}
\newcommand \Dmosz {$\gD$-modules}

\newcommand \discri{discriminant }
\newcommand \discris{discriminants }
\newcommand \discriz{discriminant}
\newcommand \discrisz{discriminants}

\newcommand \dvn {derivation }
\newcommand \dvns {derivations }
\newcommand \dvnz {derivation}
\newcommand \dvnsz {derivations}

\newcommand \dvz {zerodivisor }
\newcommand \dvzs {zerodivisors }
\newcommand \dvzz {zerodivisor}
\newcommand \dvzsz {zerodivisors}

\newcommand \dve {divisibility }
\newcommand \dvez {divisibility}

\newcommand \dvee {with explicit \dve}
\newcommand \dveez {with explicit \dvez}


\newcommand \eco {\com \elts}
\newcommand \ecoz {\com \eltsz}

\newcommand \ecr {\cor \elts}
\newcommand \ecrz {\cor \eltsz}

\newcommand \eds {sca\-lar ex\-ten\-sion }
\newcommand \edsz {sca\-lar ex\-ten\-sion}
\newcommand \Eds {Sca\-lar ex\-ten\-sion }
\newcommand \Edsz {Sca\-lar ex\-ten\-sion}

\newcommand \egmt {also }
\newcommand \egmtz {also}
\newcommand \Egmt {Also }
\newcommand \Egmtz {Also}

\newcommand \egt {equality } 
\newcommand \egts {equalities }
\newcommand \egtz {equality}
\newcommand \egtsz {equalities}

\newcommand \Egt {Equa\-lity }
\newcommand \Egtz {Equa\-lity}
\newcommand \Egts {Equa\-lities }
\newcommand \Egtsz {Equa\-lities}

\newcommand \eli{eli\-mi\-na\-tion }
\newcommand \eliz{eli\-mi\-na\-tion}
\newcommand \Eli{Eli\-mi\-na\-tion }
\newcommand \Eliz{Eli\-mi\-na\-tion}

\newcommand \elgbm{\imlgz\elr \lgb machinery }
\newcommand \elgbmz{\imlgz\elr \lgb machinery}
\newcommand \elgbmd{\imlgdz\elr \lgb machinery }
\newcommand \elgbmdz{\imlgdz\elr \lgb machinery}

\newcommand \Elr{Ele\-men\-tary }
\newcommand \elr{ele\-men\-tary }
\newcommand \elrs{ele\-men\-tary }
\newcommand \elrz{ele\-men\-tary}
\newcommand \elrsz{ele\-men\-tary}

\newcommand \elrt{ele\-men\-ta\-rily }
\newcommand \elrtz{ele\-men\-ta\-rily}

\newcommand \elt{el\-e\-ment }
\newcommand \elts{el\-e\-ments }
\newcommand \eltz{el\-e\-ment}
\newcommand \eltsz{el\-e\-ments}

\def \endo {endo\-mor\-phism }
\def \endos {endo\-mor\-phisms }
\def \endoz {endo\-mor\-phism}
\def \endosz {endo\-mor\-phisms}

\def \entrel {en\-tail\-ment re\-la\-tion }
\def \entrelz {en\-tail\-ment re\-la\-tion}
\def \entrels {en\-tail\-ment re\-la\-tions }
\def \entrelsz {en\-tail\-ment re\-la\-tions}

\def \enum {enumerable }
\def \enums {enumerable }
\def \enumz {enumerable}
\def \enumsz {enumerable}

\newcommand \Erg {$E$-regular }
\newcommand \Erge {$E$-regular }
\newcommand \Ergs {$E$-regular }
\newcommand \Erges {$E$-regular }
\newcommand \Ergez {$E$-regular}
\newcommand \Ergz {$E$-regular}
\newcommand \Ergsz {$E$-regular}
\newcommand \Ergesz {$E$-regular}

\newcommand\evc{vector space }
\newcommand\evcs{vector spaces }
\newcommand\evcz{vector space}
\newcommand\evcsz{vector spaces}

\newcommand\sevc{vector subspace }
\newcommand\sevcs{vector subspaces }
\newcommand\sevcz{vector subspace}
\newcommand\sevcsz{vector subspaces}

\newcommand \eqn{equation }
\newcommand \eqns{equations }
\newcommand \eqnz{equation}
\newcommand \eqnsz{equations}

\newcommand \Eqn{Equation }
\newcommand \Eqns{Equations }
\newcommand \Eqnz{Equation}
\newcommand \Eqnsz{Equations}

\newcommand \eqv {equi\-valent }
\newcommand \eqve {equi\-valent }
\newcommand \eqvs {equi\-valent }
\newcommand \eqves {equi\-valent }
\newcommand \eqvz {equi\-valent}
\newcommand \eqvez {equi\-valent}
\newcommand \eqvsz {equi\-valent}
\newcommand \eqvesz {equi\-valent}

\newcommand \eqvt {equiva\-lently }
\newcommand \eqvtz {equiva\-lently}
\newcommand \Eqvt {Equiva\-lently }
\newcommand \Eqvtz {Equiva\-lently}

\newcommand \eqvc {equi\-va\-lence }
\newcommand \eqvcs {equi\-va\-lences }
\newcommand \eqvcz {equi\-va\-lence}
\newcommand \eqvcsz {equi\-va\-lences}

\newcommand\evn{evaluation }
\newcommand\evnz{evaluation}
\newcommand\evns{evaluations }
\newcommand\evnsz{evaluations}

\newcommand \fab {bounded \fcn}
\newcommand \fabz {bounded \fcnz}
\newcommand \fabs {bounded \fcns}
\newcommand \fabsz {bounded \fcnsz}

\newcommand \fac {total \fcn}
\newcommand \facz {total \fcnz}
\newcommand \facs {total \fcns}
\newcommand \facsz {total \fcnsz}

\newcommand \fap {partial \fcn}
\newcommand \faps {partial \fcns}
\newcommand \fapz {partial \fcnz}
\newcommand \fapsz {partial \fcnsz}

\newcommand \fcn {facto\-ri\-zation }
\newcommand \fcns {facto\-ri\-zations }
\newcommand \fcnz {facto\-ri\-zation}
\newcommand \fcnsz {facto\-ri\-zations}

\newcommand \fima {maximal filter }
\newcommand \fimas {maximal filters }
\newcommand \fimaz {maximal filter}
\newcommand \fimazs {maximal filters}

\newcommand \fit {faithfully }

\newcommand \fip {prime filter }
\newcommand \fips {prime filters }
\newcommand \fipz {prime filter}
\newcommand \fipsz {prime filters}

\newcommand \fipma {maximal \fip }
\newcommand \fipmas {maximal \fips }
\newcommand \fipmaz {maximal \fipz}
\newcommand \fipmazs {maximal \fipsz}

\newcommand \fdi {strongly dis\-crete }
\newcommand \fdis {strongly dis\-crete }
\newcommand \fdisz {strongly dis\-crete}
\newcommand \fdiz {strongly dis\-crete}

\newcommand \fmt {formally }
\newcommand \fmtz {formally}

\newcommand \fnt {\fmtz unramified } 
\newcommand \fnts {\fmt unramified }
\newcommand \fntz {\fmt unramified}
\newcommand \fntsz {\fmt unramified}

\newcommand \fpt {\fit flat }
\newcommand \fpte {\fit flat }
\newcommand \fpts {\fit flat }
\newcommand \fptes {\fit flat }
\newcommand \fptz {\fit flat}
\newcommand \fptez {\fit flat}
\newcommand \fptsz {\fit flat}
\newcommand \fptesz {\fit flat}

\newcommand \frg{regular function }
\newcommand \frgs{regular functions }
\newcommand \frgz{regular function}
\newcommand \frgsz{regular functions}

\newcommand \Frg {$F$-regular }
\newcommand \Frge {$F$-regular }
\newcommand \Frgs {$F$-regular }
\newcommand \Frges {$F$-regular }
\newcommand \Frgz {$F$-regular}
\newcommand \Frgez {$F$-regular}
\newcommand \Frgsz {$F$-regular}
\newcommand \Frgesz {$F$-regular}

\newcommand \ftr {trace form }
\newcommand \ftrz {trace form}

\newcommand\gaq{\agq \gmt}
\newcommand\gaqz{\agq \gmtz}

\newcommand\gmt{geometry }
\newcommand\gmts{geometries }
\newcommand\gmtz{geometry}
\newcommand\gmtsz{geometries}

\newcommand\gmq{geometric }
\newcommand\gmqs{geometric }
\newcommand\gmqz{geometric}
\newcommand\gmqsz{geometric}
\newcommand\Gmq{Geometric }

\newcommand\gmqt{geometrically }
\newcommand\gmqtz{geometrically}
\newcommand \Gmqtz{Geometrically}

\newcommand\gne{generalized }
\newcommand\gnee{generalized }
\newcommand\gnes{generalized }
\newcommand\gnees{generalized }
\newcommand\gnez{generalized}
\newcommand\gneez{generalized}
\newcommand\gnesz{generalized}
\newcommand\gneesz{generalized}

\newcommand\gnl{general }
\newcommand\gnle{general }
\newcommand\gnls{general }
\newcommand\gnles{general }
\newcommand\gnlz{general}
\newcommand\gnlez{general}
\newcommand\gnlsz{general}
\newcommand\gnlesz{general}
\newcommand\Gnl{General }
\newcommand\Gnlz{General}

\newcommand\gnlt{generally }
\newcommand\gnltz{generally}
\newcommand\Gnlt{Generally }
\newcommand\Gnltz{Generally}

\newcommand\gnn{generalization }
\newcommand\gnns{generalizations }
\newcommand\gnnz{generalization}
\newcommand\gnnsz{generalizations}

\newcommand\gnq{generic }
\newcommand\gnqs{generic }
\newcommand\gnqz{generic}
\newcommand\gnqsz{generic}

\newcommand\gns{generalize }
\newcommand\gnss{generalizes }
\newcommand\gnr{generalize } 

\newcommand\gnt{generality }
\newcommand\gnts{generalities }
\newcommand\gntz{generality}
\newcommand\gntsz{generalities}

\newcommand \Grg {$G$-regular }
\newcommand \Grge {$G$-regular }
\newcommand \Grgs {$G$-regular }
\newcommand \Grges {$G$-regular }
\newcommand \Grgz {$G$-regular}
\newcommand \Grgez {$G$-regular}
\newcommand \Grgsz {$G$-regular}
\newcommand \Grgesz {$G$-regular}

\newcommand \grl{$l$-group }
\newcommand \grls{$l$-groups }
\newcommand \grlz{$l$-group}
\newcommand \grlsz{$l$-groups}

\newcommand \sgrl{$l$-subgroup }
\newcommand \sgrls{$l$-subgroups }
\newcommand \sgrlz{$l$-subgroup}
\newcommand \sgrlsz{$l$-subgroups}

\newcommand\gtr{generator }
\newcommand\gtrs{generators }
\newcommand\gtrz{generator}
\newcommand\gtrsz{generators}


\newcommand \homo {homomorphism }
\newcommand \homos {homomorphisms }
\newcommand \homosz {homomorphisms}
\newcommand \homoz {homomorphism}

\newcommand \hmg {homogeneous }
\newcommand \hmgs {homogeneous }
\newcommand \hmgz {homogeneous}
\newcommand \hmgsz {homogeneous}
\newcommand \Hmg {Homogeneous }

\newcommand \Hrg {$H$-regular }
\newcommand \Hrge {$H$-regular }
\newcommand \Hrgs {$H$-regular }
\newcommand \Hrges {$H$-regular }
\newcommand \Hrgz {$H$-regular}
\newcommand \Hrgez {$H$-regular}
\newcommand \Hrgsz {$H$-regular}
\newcommand \Hrgesz {$H$-regular}

\newcommand \icl {integrally closed }
\newcommand \iclz {integrally closed}

\newcommand \id {ideal }
\newcommand \ids {ideals }
\newcommand \idz {ideal}
\newcommand \idsz {ideals}

\newcommand \ida {\agq \idt }
\newcommand \idas {\agq \idts }
\newcommand \idasz {\agq \idtsz}
\newcommand \idaz {\agq \idtz}

\newcommand \idc  {Cramer's \idt }
\newcommand \idcs {Cramer's \idts }
\newcommand \idcsz {Cramer's \idtsz}
\newcommand \idcz {Cramer's \idtz}

\newcommand \idd {de\-ter\-mi\-nantal \id}
\newcommand \idds {de\-ter\-mi\-nantal \ids}
\newcommand \iddz {de\-ter\-mi\-nantal \idz}
\newcommand \iddsz {de\-ter\-mi\-nantal \idsz}

\newcommand \idema {maximal \id}
\newcommand \idemas {maximal \ids}
\newcommand \idemaz {maximal \idz}
\newcommand \idemasz {maximal \idsz}

\newcommand \idep {prime \id}
\newcommand \idepz {prime \idz}
\newcommand \ideps {prime \ids}
\newcommand \idepsz {prime \idsz}

\newcommand \idemi {minimal \idep}
\newcommand \idemis {minimal \ideps}
\newcommand \idemiz {minimal \idepz}
\newcommand \idemisz {minimal \idepsz}

\newcommand \iDKM {\index{Dedekind-Mertens}}
\newcommand \DKM {\iDKM Dedekind-Mertens }
\newcommand \DKMz {\iDKM Dedekind-Mertens}

\newcommand \idf {Fitting \id}
\newcommand \idfs {Fitting \ids}
\newcommand \idfz {Fitting \idz}
\newcommand \idfsz {Fitting \idsz}

\newcommand \idm {idem\-po\-tent }
\newcommand \idms {idem\-po\-tents }
\newcommand \idmz {idem\-po\-tent}
\newcommand \idmsz {idem\-po\-tents}

\newcommand \idme {\idm}
\newcommand \idmes {\idms}
\newcommand \idmez {\idmz}
\newcommand \idmesz {\idmsz}

\newcommand \idst{\spt \idm}
\newcommand \idstz{\spt \idmz}

\newcommand \idp {prin\-ci\-pal \id}
\newcommand \idps {prin\-ci\-pal \ids}
\newcommand \idpsz {prin\-ci\-pal \idsz}
\newcommand \idpz {prin\-ci\-pal \idz}

\newcommand \idt {iden\-ti\-ty }
\newcommand \idts {iden\-ti\-ties }
\newcommand \idtz {iden\-ti\-ty}
\newcommand \idtsz {iden\-ti\-ties}

\newcommand \idtr {in\-de\-ter\-mi\-nate }
\newcommand \idtrs {in\-de\-ter\-mi\-nates }
\newcommand \idtrz {in\-de\-ter\-mi\-nate}
\newcommand \idtrsz {in\-de\-ter\-mi\-nates}

\newcommand \ifr {fractional \id}
\newcommand \ifrs {fractional \ids}
\newcommand \ifrz {fractional \idz}
\newcommand \ifrsz {fractional \idsz}

\newcommand \iJG{\index{Gauss-Joyal Lemma}}

\newcommand \ihi {\index{Hilbert} }
\newcommand \ihiz {\index{Hilbert}}
\newcommand \imlgd {\index{elementary local-global machinery of \qirisz} }
\newcommand \imlgdz {\index{elementary local-global machinery of \qirisz}}
\newcommand \imlg {\index{elementary local-global machinery of \zedr \risz} }
\newcommand \imlgz {\index{elementary local-global machinery of \zedr \risz}}
\newcommand \imlb {\index{basic local-global machinery (with \idepsz)} }
\newcommand \imlbz {\index{basic local-global machinery (with \idepsz)}}
\newcommand \imla {\index{local-global machinery of \anarsz} }
\newcommand \imlaz {\index{local-global machinery of \anarsz}}
\newcommand \imlma {\index{local-global machinery of \idemasz} } 
\newcommand \imlmaz {\index{local-global machinery of \idemasz}}
\newcommand \iplg {\index{basic local-global principle} }
\newcommand \iplgz {\index{basic local-global principle}}

\newcommand \imd {immediate }
\newcommand \imde {immediate }
\newcommand \imds {immediate }
\newcommand \imdes {immediate }
\newcommand \imdz {immediate}
\newcommand \imdez {immediate}
\newcommand \imdsz {immediate}
\newcommand \imdesz {immediate}

\newcommand \imdt {immediately }
\newcommand \imdtz {immediately}

\newcommand \ird {irreducible }
\newcommand \irds {irreducible }
\newcommand \irdz {irreducible}
\newcommand \irdsz {irreducible}

\newcommand \ing {\gne inverse }
\newcommand \ings {\gnes inverses }
\newcommand \ingz {\gne inverse}
\newcommand \ingsz {\gnesz inverses}

\newcommand \inoe {\index{Noether!position}}
\newcommand \iNoe {\inoe\Noe}
\newcommand \iNoez {\inoe\Noez}

\newcommand \iMP {Moore-Penrose inverse }
\newcommand \iMPz {Moore-Penrose inverse}
\newcommand \iMPs {Moore-Penrose inverses }
\newcommand \iMPsz {Moore-Penrose inverses}

\newcommand \ipp {potential \idep}
\newcommand \ipps {potential \ideps}
\newcommand \ippz {potential \idepz}
\newcommand \ippsz {potential \idepsz}

\newcommand \iso {iso\-mor\-phism }
\newcommand \isos {iso\-mor\-phisms }
\newcommand \isosz {iso\-mor\-phisms}
\newcommand \isoz {iso\-mor\-phism}

\newcommand \isoc {iso\-mor\-phic }
\newcommand \isocz {iso\-mor\-phic}

\newcommand \itf {\tf \id}
\newcommand \itfs {\tf \ids}
\newcommand \itfz {\tf \idz}
\newcommand \itfsz {\tf \idsz}

\newcommand \inv {inverse }
\newcommand \invs {inverses }
\newcommand \invz {inverse}
\newcommand \invsz {inverses}

\newcommand \iv {invertible }
\newcommand \ivs {invertible }
\newcommand \ivz {invertible}
\newcommand \ivsz {invertible}

\newcommand \KRA {\index{Kronecker!trick}Kronecker }
\newcommand \KRO {\index{Kronecker!\tho  (1)}Kronecker }
\newcommand \KRN {\index{Kronecker!\tho  (2)}Kronecker }
\newcommand \KRAz {\index{Kronecker!trick}Kronecker}
\newcommand \KROz {\index{Kronecker!\tho  (1)}Kronecker}
\newcommand \KRNz {\index{Kronecker!\tho  (2)}Kronecker}


\newcommand \lgb {local-global }
\newcommand \lgbe {local-global }
\newcommand \lgbs {local-global }
\newcommand \lgbz {local-global}
\newcommand \lgbez {local-global}
\newcommand \lgbsz {local-global}
\newcommand \lgbes {local-global }
\newcommand \Lgb {Local-global }

\newcommand \lin {lin\-e\-ar }
\newcommand \lins {lin\-e\-ar }
\newcommand \linz {lin\-e\-ar}
\newcommand \linsz {lin\-e\-ar}

\newcommand \lint {lin\-e\-arly }

\newcommand \lmo {\lot cyclic }
\newcommand \lmos {\lot cyclic }
\newcommand \lmoz {\lot cyclic}
\newcommand \lmosz {\lot cyclic}
\newcommand \Lcy {\Lot cyclic }
\newcommand \Lcyz {\Lot cyclic}

\newcommand \lnl {\lot \nl}
\newcommand \lnls {\lot \nls}
\newcommand \lnlz {\lot \nlz}
\newcommand \lnlsz {\lot \nlsz}
\newcommand \Lnls {\Lot \nls }

\def \lot {lo\-cal\-ly }
\def \lotz {lo\-cal\-ly}
\def \Lot {Locally }

\def \lon {lo\-cal\-i\-za\-tion }
\def \lons {lo\-cal\-i\-za\-tions }
\def \lonz {lo\-cal\-i\-za\-tion}
\def \lonsz {lo\-cal\-i\-za\-tions}
\def \Lon {Lo\-cal\-i\-za\-tion }

\def \lop {\lot prin\-ci\-pal }
\def \lops {\lot prin\-ci\-pal }
\def \lopsz {\lot prin\-ci\-pal}
\def \lopz {\lot prin\-ci\-pal}
\def \Lop {\Lot prin\-ci\-pal }

\newcommand \losd {\lot \sdz }
\newcommand \losdz {\lot \sdzz}
\newcommand \Losd {\Lot \sdz }

\newcommand \lsd {pf-\ri }
\newcommand \lsds {pf-\ris }
\newcommand \lsdsz {pf-\risz}
\newcommand \lsdz {pf-\riz}
\newcommand \Lsd {Pf-\ri }
\newcommand \Lsds {Pf-\ris }
\newcommand \Lsdz {Pf-\riz}
\newcommand \Lsdsz {Pf-\risz}
\newcommand \mdi {module of \diles }
\newcommand \mdiz {module of \dilesz}
\newcommand \mdis {modules of \diles }
\newcommand \mdisz {modules of \dilesz}

\newcommand \mlm {\lmo module }
\newcommand \mlms {\lmo modules }
\newcommand \mlmz {\lmo module}
\newcommand \mlmsz {\lmo modules}

\newcommand \mlmo {cyclic \lon matrix } 
\newcommand \mlmos {cyclic \lon matrices } 
\newcommand \mlmoz {cyclic \lon matrix} 
\newcommand \mlmosz {cyclic \lon matrices} 

\newcommand \mlp {principal \lon matrix }
\newcommand \mlpz {principal \lon matrix}
\newcommand \mlps {principal \lon matrices }
\newcommand \mlpsz {principal \lon matrices}

\newcommand \mlr {\elr operation }
\newcommand \mlrs {\elr operations }
\newcommand \mlrz {\elr operation}
\newcommand \mlrsz {\elrs operations}

\newcommand \mlrr {\elr row operation }
\newcommand \mlrrs {\elr row operations }
\newcommand \mlrrz {\elr row operation}
\newcommand \mlrrsz {\elrs row operations}

\newcommand \mlrc {\elr column operation }
\newcommand \mlrcs {\elr column operations }
\newcommand \mlrcz {\elr column operation}
\newcommand \mlrcsz {\elrs column operations}

\newcommand \mo {monoid }
\newcommand \mos {monoids }
\newcommand \mosz {monoids}
\newcommand \moz {monoid}

\newcommand \moco {\com\mos}
\newcommand \mocoz {\com\mosz}

\newcommand \molo {\lon morphism }
\newcommand \molos {\lon morphisms }
\newcommand \moloz {\lon morphism}
\newcommand \molosz {\lon morphisms}

\newcommand \mom {monomial }
\newcommand \moms {monomials }
\newcommand \momz {monomial}
\newcommand \momsz {monomials}

\newcommand \mpf {\pf module }
\newcommand \mpfs {\pf modules }
\newcommand \mpfz {\pf module}
\newcommand \mpfsz {\pf modules}

\newcommand \mpl {flat module }
\newcommand \mpls {flat modules }
\newcommand \mplz {flat module}
\newcommand \mplsz {flat modules}

\newcommand \mpn {\pn matrix }
\newcommand \mpns {\pn matrices }
\newcommand \mpnz {\pn matrix}
\newcommand \mpnsz {\pn matrices}
\newcommand \Mpn {Presentation matrix }

\newcommand \mprn {\prn matrix }
\newcommand \mprns {\prn matrices }
\newcommand \mprnz {\prn matrix}
\newcommand \mprnsz {\prn matrices}
\newcommand \Mprn {\Prn matrix }
\newcommand \Mprns {\Prn matrices }

\newcommand \mptf {\ptf module }
\newcommand \mptfs {\ptf modules }
\newcommand \mptfz {\ptf module}
\newcommand \mptfsz {\ptf modules}
\newcommand \Mptf {\Ptf module }
\newcommand \Mptfs {\Ptf modules }
\newcommand \Mptfsz {\Ptf modules}

\newcommand \smptf {\ptf submodule }
\newcommand \smptfs {\ptf submodules }
\newcommand \smptfz {\ptf submodule}
\newcommand \smptfsz {\ptf submodules}

\newcommand \mrc {\pro module of constant rank }
\newcommand \mrcz {\pro module of constant rank}
\newcommand \mrcs {\pro modules of constant rank }
\newcommand \mrcsz {\pro modules of constant rank}

\newcommand \mtf {\tf module }
\newcommand \mtfs {\tf modules }
\newcommand \mtfz {\tf module}
\newcommand \mtfsz {\tf modules}
\newcommand \Mtf {\Tf module }
\newcommand \Mtfs {\Tf modules }
\newcommand \smtf {\tf submodule }
\newcommand \smtfs {\tf submodules }
\newcommand \smtfz {\tf submodule}
\newcommand \smtfsz {\tf submodules}
\newcommand \Smtf {\Tf submodule }
\newcommand \Smtfs {\Tf submodules }


\newcommand \imN {\index{Newton!method}\index{method!Newton's ---}} 
\newcommand \isN {\index{Newton!sums}}

\newcommand \ncr{nec\-es\-sa\-ry }
\newcommand \ncrs{nec\-es\-sa\-ry }
\newcommand \ncrz{nec\-es\-sa\-ry}
\newcommand \ncrsz{nec\-es\-sa\-ry}

\newcommand \ncrt{nec\-es\-sa\-rily }
\newcommand \ncrtz{nec\-es\-sa\-rily}

\newcommand \ndz {regular }
\newcommand \ndzs {regular }
\newcommand \ndzz {regular}
\newcommand \ndzsz {regular}
\newcommand \ndze {regular }
\newcommand \ndzes {regular }
\newcommand \ndzez {regular}
\newcommand \ndzesz {regular}

\newcommand \nl {simple }
\newcommand \nlz {simple}
\newcommand \nls {simple }
\newcommand \nlsz {simple}

\newcommand \noe {Noetherian }
\newcommand \noes {Noetherian }
\newcommand \noee {Noetherian }
\newcommand \noees {Noetherian }
\newcommand \noez {Noetherian}
\newcommand \noesz {Noetherian}
\newcommand \noeez {Noetherian}
\newcommand \noeesz {Noetherian}
\newcommand \Noe {Noether }
\newcommand \Noez {Noether}
\newcommand \noet {Noetherianity }
\newcommand \noetz {Noetherianity}
\newcommand \noep {\Noe position }
\newcommand \noepz {\Noe position}

\newcommand \inst {\index{Nullstellensatz}}
\newcommand \nst {\inst Nullstellensatz }
\newcommand \nstz {\inst Nullstellensatz}
\newcommand \nsts {\inst Nullstellens\"atze }
\newcommand \nstsz {\inst Nullstellens\"atze}


\newcommand \op{operation }
\newcommand \ops{operations }
\newcommand \opz{operation}
\newcommand \opsz{operations}
\newcommand \opari{\ari\op}
\newcommand \oparis{\ari\ops}
\newcommand \opariz{\ari\opz}
\newcommand \oparisz{\ari\opsz}

\newcommand \oqc {compact-open subspace }
\newcommand \oqcs {compact-open subspaces }
\newcommand \oqcz {compact-open subspace}
\newcommand \oqcsz {compact-open subspaces}

\newcommand \ort{orthogonal }
\newcommand \orte{orthogonal }
\newcommand \orts{orthogonal }
\newcommand \ortes{orthogonal }
\newcommand \ortz{orthogonal}
\newcommand \ortez{orthogonal}
\newcommand \ortsz{orthogonal}
\newcommand \ortesz{orthogonal}


\newcommand \pa {saturated pair }
\newcommand \pas {saturated pairs }
\newcommand \paz {saturated pair}
\newcommand \pasz {saturated pairs}

\newcommand \paral{parallel }
\newcommand \parals{parallel }
\newcommand \paralz{parallel}
\newcommand \paralsz{parallel}

\newcommand \paralm{parallel }

\newcommand \pb{problem }
\newcommand \pbs{problems }
\newcommand \pbz{problem}
\newcommand \pbsz{problems}
\newcommand \Pbm{Problem }
\newcommand \Pbms{Problems }
\newcommand \Pbmz{Problem}
\newcommand \Pbmsz{Problems}

\newcommand \pf {fi\-nite\-ly presented }
\newcommand \pfz {fi\-nite\-ly presented}
\newcommand \Fp {Finitely presented }
\newcommand \Fpz {Finitely presented}

\newcommand \plc {\rdt\zed}
\newcommand \plcs {\rdt\zed}
\newcommand \plcz {\rdt\zedz}
\newcommand \plcsz {\rdt\zedz}

\newcommand \plg {local-global principle }
\newcommand \plgs {local-global principles }
\newcommand \plgz {local-global principle}
\newcommand \plgsz {local-global principles}
\newcommand \Plg {Local-global principle }
\newcommand \Plgs {Local-global principles }
\newcommand \Plgz {Local-global principle}
\newcommand \Plgsz {Local-global principles}

\newcommand \plga {abstract \plg}
\newcommand \plgas {abstract \plgs}
\newcommand \plgaz {abstract \plgz}
\newcommand \plgasz {abstract \plgsz}

\newcommand \plgc {concrete \plg}
\newcommand \plgcz {concrete \plgz}
\newcommand \plgcs {concrete \plgs}
\newcommand \plgcsz {concrete \plgsz}

\newcommand \PLCC[2]{\rdb
\mni{\bf Concrete local-global principle \ref{#1} bis~} {\it #2

}}

\newcommand \pn {presentation }
\newcommand \pns {presentations }
\newcommand \pnz {presentation}
\newcommand \pnsz {presentations}

\newcommand \pog {\hmg\pol}
\newcommand \pogs {\hmg\pols}
\newcommand \pogz {\hmg\polz}
\newcommand \pogsz {\hmg\polsz}

\newcommand \Pol {Polynomial }
\newcommand \Pols {Polynomials }

\newcommand \pol {polynomial }
\newcommand \pols {polynomials }
\newcommand \polz {polynomial}
\newcommand \polsz {polynomials}

\newcommand \Polcar {\Cara \pol}
\newcommand \polcar {\cara\pol}
\newcommand \polcarz {\cara\polz}
\newcommand \polcars {\cara\pols}
\newcommand \polcarsz {\cara\polsz}

\newcommand \polfon {fundamental \pol }
\newcommand \polfonz {fundamental \polz}

\newcommand \poll{polynomial }
\newcommand \polls{polynomial }
\newcommand \pollsz{polynomial}
\newcommand \pollz{polynomial}
\newcommand \polle{polynomial }
\newcommand \polles{polynomial }
\newcommand \pollesz{polynomial}
\newcommand \pollez{polynomial}

\newcommand \pollt{polynomially }
\newcommand \polltz{polynomially}

\newcommand \polmin {minimal \pol}
\newcommand \polminz {minimal \polz}
\newcommand \polmins {minimal \pols}
\newcommand \polminsz {minimal \polsz}

\newcommand \polmu {rank \pol}
\newcommand \polmus {rank \pols}
\newcommand \polmuz{rank \polz}

\newcommand \polu {\mon \pol }
\newcommand \polus {\mon \pols }
\newcommand \poluz {\mon \polz }
\newcommand \polusz {\mon \polsz }

\newcommand \ppv {primitive by values } 
\newcommand \ppvs {primitive by values }
\newcommand \ppvz {primitive by values}
\newcommand \ppvsz {primitive by values}

\newcommand \pppv {primitive \pol by values } 
\newcommand \pppvs {primitive \pols by values }
\newcommand \pppvz {primitive \pol by values}
\newcommand \pppvsz {primitive \pols by values}

\newcommand \prc  {\pro of constant rank }
\newcommand \prcs  {\pros of constant rank }
\newcommand \prcz  {\pro of constant rank}
\newcommand \prcsz  {\pros of constant rank}

\newcommand \prca {\rca principle }
\newcommand \prcc {\rcc principle }
\newcommand \prce {\rce principle }

\newcommand \prf {closed covering principle }
\newcommand \prfz {closed covering principle}

\newcommand \Prmt {Precisely }
\newcommand \Prmtz {Precisely}
\newcommand \prmt {precisely }
\newcommand \prmtz {precisely}

\newcommand \prn {pro\-jec\-tion }
\newcommand \prns {pro\-jec\-tions }
\newcommand \prnz {pro\-jec\-tion}
\newcommand \prnsz {pro\-jec\-tions}
\newcommand \Prn {Projection }

\newcommand \pro {pro\-jec\-tive }
\newcommand \pros {pro\-jec\-tive }
\newcommand \proz {pro\-jec\-tive}
\newcommand \prosz {pro\-jec\-tive}
\newcommand \Pro {Projective }

\newcommand \prof {depth }
\newcommand \profz {depth}
\newcommand \profs {depths }
\newcommand \profsz {depths}
\newcommand \Prof {Depth }

\newcommand \prr {pro\-jec\-tor }
\newcommand \prrs {pro\-jec\-tors }
\newcommand \prrz {pro\-jec\-tor}
\newcommand \prrsz {pro\-jec\-tors}

\newcommand \Prt {Property }
\newcommand \Prts {Properties }
\newcommand \prt {property }
\newcommand \prts {properties }
\newcommand \prtz {property}
\newcommand \prtsz {properties}

\newcommand \ptf {\tf \pro }
\newcommand \ptfz {\tf \proz}
\newcommand \ptfs {\tf \pro }
\newcommand \ptfsz {\tf \proz}
\newcommand \Ptf {\Tf \pro }

\newcommand \qc {quasi-compact }
\newcommand \qcs {quasi-compact }
\newcommand \qcz {quasi-compact}
\newcommand \qcsz {quasi-compact}

\newcommand \qf {quasi-free }
\newcommand \qfz {quasi-free}

\newcommand \qiri {pp-ring }
\newcommand \qiris {pp-rings }
\newcommand \qirisz {pp-rings}
\newcommand \qiriz {pp-ring}
\newcommand \Qiri {Quasi-integral }

\newcommand \sqiri {pp-subring }
\newcommand \sqiris {pp-subrings }
\newcommand \sqirisz {pp-subrings}
\newcommand \sqiriz {pp-subring}

\newcommand \qiv {quasi-inverse }
\newcommand \qivs {quasi-inverse }
\newcommand \qivz {quasi-inverse}
\newcommand \qivsz {quasi-inverse}

\newcommand \qnl {quasi-\nl}
\newcommand \qnls {quasi-\nl}
\newcommand \qnlz {quasi-\nlz}
\newcommand \qnlsz {quasi-\nlz}

\newcommand \qreg {quasi-regular }
\newcommand \qregz {quasi-regular}
\newcommand \qregs {quasi-regular }
\newcommand \qregsz {quasi-regular}

\newcommand \qtf {quantifier }
\newcommand \qtfs {quantifiers }
\newcommand \qtfz {quantifier}
\newcommand \qtfsz {quantifiers}

\newcommand \qtn {quantification }
\newcommand \qtnz {quantification}

\newcommand \ri {ring }
\newcommand \riz {ring}
\newcommand \ris {rings }
\newcommand \risz {rings}

\newcommand \Ri {Ring }
\newcommand \Riz {Ring}
\newcommand \Ris {Rings }
\newcommand \Risz {Rings}

\newcommand \RC {R-compatible }
\newcommand \RCs {R-compatible }
\newcommand \RCz {R-compatible}
\newcommand \RCsz {R-compatible}

\newcommand \rcm {patching }
\newcommand \rcmz {patching}
\newcommand \rcms {patchings }

\newcommand \rcc {concrete \rcm }
\newcommand \rccz {concrete \rcmz}
\newcommand \rca {abstract \rcm }
\newcommand \rcaz {abstract \rcmz}
\newcommand \rce {\rcc of equalities }

\newcommand \rde {dependence re\-la\-tion }
\newcommand \rdes {dependence re\-la\-tions }
\newcommand \rdesz {dependence re\-la\-tions}
\newcommand \rdez {dependence re\-la\-tion}
\newcommand \rdi {integral \rde }
\newcommand \rdis {integral \rdes }
\newcommand \rdiz {integral \rdez}
\newcommand \rdisz {integral \rdesz}

\newcommand \rdl {\lin \rde }
\newcommand \rdls {\lin \rdes }
\newcommand \rdlsz {\lin \rdesz}
\newcommand \rdlz {\lin \rdez}
\newcommand \rdt {residually }
\newcommand \Rdt {Residually }

\newcommand \Recol {Patching }
\newcommand \recol {patching }
\newcommand \Recolz {patching}

\newcommand \reg {regular }
\newcommand \regs {regular }
\newcommand \regz {regular}
\newcommand \regsz {regular}

\newcommand \rpf {\pf reduced }
\newcommand \rpfs {\pf reduced }
\newcommand \rpfz {\pf reduced}
\newcommand \rpfsz {\pf reduced}

\newcommand \scf {finitary scheme }
\newcommand \scfs {finitary schemes }
\newcommand \scfz {finitary scheme}
\newcommand \scfsz {finitary schemes}

\newcommand \scl {\elr scheme }
\newcommand \scls {\elr schemes}
\newcommand \sclz {\elr scheme}
\newcommand \sclsz {\elr schemes}

\newcommand \sdo {\ort\sdr}
\newcommand \sdos {\ort\sdrs}
\newcommand \sdoz {\ort\sdrz}
\newcommand \sdosz {\ort\sdrsz}

\newcommand \sdr {direct sum }
\newcommand \sdrs {direct sums }
\newcommand \sdrz {direct sum}
\newcommand \sdrsz {direct sums}

\newcommand \sdz {without \dvzs}
\newcommand \sdzz {without \dvzsz}

\newcommand \seqreg {regular sequence }
\newcommand \seqregz {regular sequence}
\newcommand \seqregs {regular sequences }
\newcommand \seqregsz {regular sequences}

\newcommand \seqqreg {\qreg sequence}
\newcommand \seqqregz {\qreg sequence}
\newcommand \seqqregs {\qreg sequences}
\newcommand \seqqregsz {\qreg sequences}

\newcommand \sErg {$E$-regular sequence }
\newcommand \sErgs {$E$-regular sequences }
\newcommand \sErgz {$E$-regular sequence}
\newcommand \sErgsz {$E$-regular sequences}

\newcommand \sfio {fundamental system of orthogonal \idms }
\newcommand \sfios {fundamental systems of orthogonal \idms }
\newcommand \sfioz {fundamental system of orthogonal \idmsz}
\newcommand \sfiosz {fundamental systems of orthogonal \idmsz}

\newcommand \sgr {\gtr set }
\newcommand \sgrs {\gtr sets }
\newcommand \sgrz {\gtr set}
\newcommand \sgrsz {\gtr sets}

\newcommand \sing {singular }
\newcommand \sings {singular }
\newcommand \singz {singular}
\newcommand \singsz {singular}

\newcommand \sli {\sys of \lin \eqns }
\newcommand \slis {\syss of \lin \eqns }
\newcommand \slisz {\syss of \lin \eqnsz}
\newcommand \sliz {\sys of \lin \eqnsz}

\newcommand \smq {symmetric }
\newcommand \smqs {symmetric }
\newcommand \smqz {symmetric}
\newcommand \smqsz {symmetric}

\newcommand \Smq {Symmetric }

\newcommand \spl {separable }
\newcommand \spls {separable }
\newcommand \splz {separable}
\newcommand \splsz {separable}
\newcommand \sply {separably }
\newcommand \splyz {separably}
\newcommand\spt{separability }
\newcommand\sptz{separability}
\newcommand\Spt{Separability }
\newcommand\Sptz{Separability}

\newcommand \ste {strictly \'etale }
\newcommand \stes {strictly \'etale }
\newcommand \stez {strictly \'etale}
\newcommand \stesz {strictly \'etale}
\newcommand \Ste {Strictly \'etale }

\newcommand \stf {strictly finite }
\newcommand \stfs {strictly finite }
\newcommand \stfz {strictly finite}
\newcommand \stfsz {strictly finite}
\newcommand \stfe {strictly finite }
\newcommand \stfes {strictly finite }
\newcommand \stfez {strictly finite}
\newcommand \stfesz {strictly finite}
\newcommand \Stf {Strictly finite }

\newcommand \stl {stably free }
\newcommand \stls {stably free }
\newcommand \stlz {stably free}
\newcommand \stlsz {stably free}

\newcommand \sul {additional } 
\newcommand \suls {additional }
\newcommand \sulz {additional}
\newcommand \sulsz {additional}

\newcommand \supl {complementary } 
\newcommand \supls {complementary } 
\newcommand \suplz {complementary} 
\newcommand \suplsz {complementary} 

\newcommand \Sus {Suslin\index{Suslin} }
\newcommand \Susz {Suslin\index{Suslin}}

\newcommand \susi {\sing sequence }
\newcommand \susis {\sing sequences }
\newcommand \susiz {\sing sequence}
\newcommand \susisz {\sing sequences}

\newcommand \Sut {Support }
\newcommand \Suts {Supports }
\newcommand \sut {support }
\newcommand \suts {supports }
\newcommand \sutz {support}
\newcommand \sutsz {supports}

\newcommand \syc {coordinate \sys}
\newcommand \sycs {coordinate \syss}
\newcommand \sycz {coordinate \sysz}
\newcommand \sycsz {coordinate \syssz}

\newcommand \sml {strict semi-local }
\newcommand \smls {strict semi-local }
\newcommand \smlz {strict semi-local}
\newcommand \smlss {strict semi-local}

\newcommand \syp {\poll \sys}
\newcommand \syps {\poll \syss}
\newcommand \sypz {\poll \sysz}
\newcommand \sypsz {\poll \syssz}

\newcommand \sys {system }
\newcommand \syss {systems }
\newcommand \sysz {system}
\newcommand \syssz {systems}
\newcommand \Sys {System }
\newcommand \Sysz {System}

\newcommand \slgb {semi-local }
\newcommand \slgbs {semi-local }
\newcommand \slgbz {semi-local}
\newcommand \slgbsz {semi-local}
\newcommand \Slgb {Semi-local }
\newcommand \Slgbs {Semi-local }

\newcommand \spb {separable }  
\newcommand \spbs {separable }
\newcommand \spbz {separable}
\newcommand \spbsz {separable}
\newcommand \Spb {Separable }

\newcommand \spby{separably }
\newcommand \Spby{Separably }

\newcommand \srg {regular sequence }
\newcommand \srgs {regular sequences }
\newcommand \srgz {regular sequence}
\newcommand \srgsz {regular sequences}

\newcommand \SSO {Serre's Splitting Off\index{Serre's Splitting Off \thoz} }
\newcommand \SSOz {Serre's Splitting Off\index{Serre's Splitting Off \thoz}}


\newcommand \tf {fi\-nite\-ly gen\-e\-ra\-ted }
\newcommand \tfz {fi\-nite\-ly gen\-e\-ra\-ted}
\newcommand \Tf {Finitely gen\-e\-ra\-ted }
\newcommand \Tfz {Finitely gen\-e\-ra\-ted}

\newcommand \Tho {Theorem }
\newcommand \Thoz {Theorem}
\newcommand \Thos {Theorems }
\newcommand \Thosz {Theorems}
\newcommand \tho {theorem }
\newcommand \thos {theorems }
\newcommand \thoz {theorem}
\newcommand \thosz {theorems}

\newcommand \thoc {\thoz\etoz~}
\newcommand \Thoc {\Thoz\etoz~}

\newcommand \torf {torsion-free }
\newcommand \torfz {torsion-free}

\newcommand \styc {trace \sys of \coos}
\newcommand \stycz {trace \sys of \coosz}

\newcommand \trdi {distributive lattice }
\newcommand \trdis {distributive lattices }
\newcommand \trdiz {distributive lattice}
\newcommand \trdisz {distributive lattices}

\newcommand \trf {the result follows }
\newcommand \Trf {The result follows }
\newcommand \trfz {the result follows}
\newcommand \Trfz {The result follows}

\newcommand \strdi {distributive sublattice }
\newcommand \strdis {distributive sublattices }
\newcommand \strdiz {distributive sublattice}
\newcommand \strdisz {distributive sublattices}

\newcommand \trel {\elr transformation }
\newcommand \trels {\elr transformations }
\newcommand \trelz {\elr transformation}
\newcommand \trelsz {\elr transformations}

\newcommand \umd {unimodular }
\newcommand \umds {unimodular }
\newcommand \umdz {unimodular}
\newcommand \umdsz {unimodular}
\newcommand \Umd {Unimodular }

\newcommand \uny {unitary } 
\newcommand \unys {unitary } 
\newcommand \unyz {unitary} 
\newcommand \unysz {unitary} 

\newcommand \mon {monic } 
\newcommand \monz {monic} 
\newcommand \Mon {Monic } 
\newcommand \Monz {Monic} 
\newcommand \mons {\mon } 
\newcommand \monsz {\monz} 

\newcommand \uvl {universal }
\newcommand \uvle {universal }
\newcommand \uvls {universal }
\newcommand \uvles {universal }
\newcommand \uvlz {universal}
\newcommand \uvlez {universal}
\newcommand \uvlsz {universal}
\newcommand \uvlesz {universal}
\newcommand \Uvl{Universal }

\newcommand \vfn {verification }
\newcommand \vfns {verifications }
\newcommand \vfnz {verification}
\newcommand \vfnsz {verifications}

\def \vmd {\umd vector }
\def \vmds {\umd vectors }
\def \vmdz {\umd vector}
\def \vmdsz {\umd vectors}

\newcommand \vrt {variety }
\newcommand \vrts {varieties }
\newcommand \vrtz {variety}
\newcommand \vrtsz {varieties}
\newcommand \vgq {\agq \vrt}
\newcommand \vgqs {\agq \vrts}
\newcommand \vgqz {\agq \vrtz}
\newcommand \vgqsz {\agq \vrtsz}

\newcommand \wir {with respect to }

\newcommand \zed {zero-dimensional }
\newcommand \zedz {zero-dimensional}
\newcommand \zede {zero-dimensional }
\newcommand \zedez {zero-dimensional}
\newcommand \zeds {zero-dimensional }
\newcommand \zedsz {zero-dimensional}
\newcommand \zedes {zero-dimensional }
\newcommand \zedesz {zero-dimensional}
\newcommand \Zed {Zero-dimensional }

\newcommand \zedr {reduced \zed }
\newcommand \zedrs {reduced \zeds }
\newcommand \zedrz {reduced \zedz}
\newcommand \zedrsz {reduced \zedz}
\newcommand \Zedr {Reduced \zed }

\newcommand \zmt {Zariski-Grothendieck \tho}
\newcommand \zmtz {Zariski-Grothendieck \thoz}


\newcommand \cof {constructive }
\newcommand \cofs {constructive }
\newcommand \cofz {constructive}
\newcommand \cofsz {constructive}

\newcommand \cov {constructive }
\newcommand \covz {constructive}
\newcommand \covsz {constructive}
\newcommand \covs {constructive }
\newcommand \Cov {Constructive }

\newcommand \coma {\cov \maths }
\newcommand \comaz {\covs \mathsz}
\newcommand \Coma {\Cov \maths }
\newcommand \clama {classical \maths }
\newcommand \clamaz {classical \mathsz}
\newcommand \Clama {Classical \maths }

\renewcommand \cot {constructively }
\newcommand \cotz {constructively}

\newcommand \matn {mathematician }
\newcommand \matnz {mathematician}
\newcommand \matne {mathematician }
\newcommand \matns {mathematicians }
\newcommand \matnsz {mathematicians}
\newcommand \matnes {mathematicians }

\newcommand \maths {mathematics }
\newcommand \mathsz {mathematics}
\newcommand \mathe {mathematical }
\newcommand \mathz {mathematical}

\newcommand \prco {\cof proof }
\newcommand \prcos {\cof proofs }
\newcommand \prcoz {\cof proof}
\newcommand \prcosz {\cof proofs}

\newcommand \wirt {with repect to }
\newcommand \wrt {w.r.t.\ }
\newcommand \resp {resp.\ }

\newcommand\cm{cm}
\makeatletter

\DeclareRobustCommand\\{%
  \let \reserved@e \relax
  \let \reserved@f \relax
  \@ifstar{\let \reserved@e \vadjust \let \reserved@f \nobreak
             \@xnewline}%
          \@xnewline}
\makeatother

\newcommand{\blocs}[8]{%
{\setlength{\unitlength}{.0833\textwidth}
\tabcolsep0pt\renewcommand{\arraystretch}{0}%
\begin{tabular}{|c|c|}
\hline
\parbox[t][#3\cm][c]{#1\cm}{\begin{minipage}[c]{#1\cm}
\centering#5
\end{minipage}}&
\parbox[t][#3\cm][c]{#2\cm}{\begin{minipage}[c]{#2\cm}
\centering#6
\end{minipage}}\\
\hline
\parbox[t][#4\cm][c]{#1\cm}{\begin{minipage}[c]{#1\cm}
\centering#7
\end{minipage}}&
\parbox[t][#4\cm][c]{#2\cm}{\begin{minipage}[c]{#2\cm}
\centering#8
\end{minipage}}\\
\hline
\end{tabular}
}}


\newcommand\tri[7]{
$$\quad\quad\quad\quad
\vcenter{\xymatrix@C=1.5cm
{
#1 \ar[d]_{#2} \ar[dr]^{#3} \\
{#4} \ar[r]_{{#5}}   & {#6} \\
}}
\quad\quad \vcenter{\hbox{\small {#7}}\hbox{~\\[1mm] ~ }}
$$
}


\newcommand\carre[8]{
$$
\xymatrix @C=1.2cm{
#1\,\ar[d]^{#4}\ar[r]^{#2}   & \,#3\ar[d]^{#5}   \\
#6\,\ar[r]    ^{#7}    & \,#8  \\
}
$$
}

\newcommand\pun[7]{
$$\quad\quad\quad\quad
\vcenter{\xymatrix@C=1.5cm
{
#1 \ar[d]_{#2} \ar[dr]^{#3} \\
{#4} \ar@{-->}[r]_{{#5}\,!}   & {#6} \\
}}
\quad\quad \vcenter{\hbox{\small {#7}}\hbox{~\\[1mm] ~ }}
$$
}

\newcommand\puN[8]{
$$\hspace{#8}
\vcenter{\xymatrix@C=1.5cm
{
#1 \ar[d]_{#2} \ar[dr]^{#3} \\
{#4} \ar@{-->}[r]_{{#5}\,!}   & {#6} \\
}}
\quad\quad \vcenter{\hbox{\small {#7}}\hbox{~\\[1mm] ~ }}
$$
}

\newcommand\Pun[8]{
$$\quad\quad\quad\quad
\vcenter{\xymatrix@C=1.5cm
{
#1 \ar[d]_{#2} \ar[dr]^{#3} \\
{#4} \ar@{-->}[r]_{{#5}\,!}   & {#6} \\
}}
\quad\quad
\vcenter{\hbox{\small {#7}}
\hbox{~\\[3.5mm] ~ }
\hbox{\small {#8}}
\hbox{~\\[-3.5mm] ~ }}
$$
}


\newcommand\PUN[9]{
$$\quad\quad
\vcenter{\xymatrix@C=1.5cm
{
#1 \ar[d]_{#2} \ar[dr]^{#3} \\
{#4} \ar@{-->}[r]_{{#5}\,!}   & {#6} \\
}}
\quad\quad
\vcenter{
\hbox{\small {#7}}
\hbox{~\\[-3mm] ~}
\hbox{\small {#8}}
\hbox{~\\[-3mm] ~}
\hbox{\small {#9}}
\hbox{~\\[-3.5mm] ~ }}
$$
}

\newcommand\Pnv[9]{
$$\quad\quad\quad\quad
\vcenter{\xymatrix@C=1.5cm
{
#1 \ar[d]_{#2} \ar[dr]^{#3} \\
{#4} \ar@{-->}[r]_{{#5}\,!}   & {#6} \\
}}
\quad\quad
\vcenter{
\hbox{\small {#7}}
\hbox{~\\[1mm] ~}
\hbox{\small {#8}}
\hbox{~\\[-1mm] ~}
\hbox{\small {#9}}
\hbox{~\\[0mm] ~ }}
$$
}

\newcommand\PNV[9]{
$$\quad\quad\quad\quad
\vcenter{\xymatrix@C=1.5cm
{
#1 \ar[d]_{#2} \ar[dr]^{#3} \\
{#4} \ar@{-->}[r]_{{#5}\,!}   & {#6} \\
}}
\quad\quad
\vcenter{
\vspace{4mm}
\hbox{\small {#7}}
\hbox{~\\[-1.7mm] ~}
\hbox{\small {#8}}
\hbox{~\\[-1.7mm] ~}
\hbox{\small {#9}}
\hbox{~\\[2mm] ~ }
}
$$
}


\newdimen\xyrowsp
\xyrowsp=3pt
\newcommand{\SCO}[6]{
\xymatrix @R = \xyrowsp {
                                  &1 \ar@{-}[dl] \ar@{-}[dr] \\
#3 \ar@{-}[ddr]                   &   & #6 \ar@{-}[ddl] \\
                                  &\bullet\ar@{-}[d] \\
                                  &\bullet   \\
#2 \ar@{-}[ddr] \ar@{-}[uur]      &   & #5 \ar@{-}[ddl] \ar@{-}[uul] \\
                                  &\bullet \ar@{-}[d] \\
                                  &\bullet  \\
#1 \ar@{-}[uur]                   &   & #4 \ar@{-}[uul] \\
                                  & 0 \ar@{-}[ul] \ar@{-}[ur] \\
}
}


\makeatletter
\newif\if@borderstar
\def\bordercmatrix{\@ifnextchar*{%
  \@borderstartrue\@bordercmatrix@i}{\@borderstarfalse\@bordercmatrix@i*}%
}
\def\@bordercmatrix@i*{\@ifnextchar[{%
  \@bordercmatrix@ii}{\@bordercmatrix@ii[()]}
}
\def\@bordercmatrix@ii[#1]#2{%
  \begingroup
    \m@th\@tempdima.875em\setbox\z@\vbox{%
      \def\cr{\crcr\noalign{\kern 2\p@\global\let\cr\endline}}%
      \ialign {$##$\hfil\kern.2em\kern\@tempdima&\thinspace%
      \hfil$##$\hfil&&\quad\hfil$##$\hfil\crcr\omit\strut%
      \hfil\crcr\noalign{\kern-\baselineskip}#2\crcr\omit%
      \strut\cr}}%
    \setbox\tw@\vbox{\unvcopy\z@\global\setbox\@ne\lastbox}%
    \setbox\tw@\hbox{\unhbox\@ne\unskip\global\setbox\@ne\lastbox}%
    \setbox\tw@\hbox{%
      $\kern\wd\@ne\kern-\@tempdima\left\@firstoftwo#1%
        \if@borderstar\kern.2em\else\kern -\wd\@ne\fi%
      \global\setbox\@ne\vbox{\box\@ne\if@borderstar\else\kern.2em\fi}%
      \vcenter{\if@borderstar\else\kern-\ht\@ne\fi%
        \unvbox\z@\kern-\if@borderstar2\fi\baselineskip}%
\if@borderstar\kern-2\@tempdima\kern.4em\else\,\fi\right\@secondoftwo#1 $%
    }\null\;\vbox{\kern\ht\@ne\box\tw@}%
  \endgroup
}
\makeatother

\newcounter{nbtotalexos}
\newcounter{nbtotalprob}
\newcommand*{\incrementeexosetprob}{%
\addtocounter{nbtotalexos}{\value{exercise}}
\addtocounter{nbtotalprob}{\value{problem}}}
\newcommand*{\finincrementeexosetprob}{%
\addtocounter{nbtotalexos}{-1}\refstepcounter{nbtotalexos}\label{nombreexos}%
\addtocounter{nbtotalprob}{-1}\refstepcounter{nbtotalprob}\label{nombreprob}%
}

\romanpagenumbers

\begingroup

\pagestyle{CMcadreseul}

\def\entree#1.{\setbox0=\hbox{#1. --- }\hangindent \wd0\hangafter1\noindent\box0\ignorespaces}

\setpagenumber1

\begin{center}\Large\bf
\CMauthor 
\end{center}
\vskip2cm
\begin{center}\Huge\bf
\CMtitle 
\end{center}
\vskip1cm
\begin{center}\Large{\bf
Course and exercises}

English translation by Tania K. Roblot

{\large Last corrections \today, \rm \normalsize see page \pageref{labcor}
}

\end{center}

\break

\cleardoublepage
\thispagestyle{CMcadreseul}
\begingroup
\parskip0pt

\noindent {\sc Henri Lombardi. } Ma\^{\i}tre de Conf\'erences at the Universit\'e de Franche-Comt\'e.
His research focuses on constructive mathematics, real algebra and algorithmic complexity.
He is one of the founders of the international group M.A.P. (Mathematics, Algorithms, Proofs), created in 2003: see the site \url{https://mapcommunity.github.io/}

\smallskip \hfill{\tt henri.lombardi@univ-fcomte.fr}\\
\null\hfill \url {http://hlombardi.free.fr}
  
\smallskip

\bigskip  {\sc Claude Quitt\'e. } Ma\^{\i}tre de Conf\'erences at the Universit\'e de Poitiers.
His research focuses on effective commutative algebra and computer algebra. 

\smallskip \hfill{\tt claude.quitte@math.univ-poitiers.fr}

\vfill
\begingroup
\noindent 
Mathematics Subject Classification (2010) 

-- Primary: 13 ~ Commutative Algebra.

-- Secondary:

\leftskip4ex
03F ~ Proof theory and constructive mathematics. 

\leftskip4ex
06D ~ Distributive lattices.

\leftskip4ex
14Q ~ Computational aspects of algebraic geometry.

\vfill
\vfill

\endgroup

\newpage
\thispagestyle{CMcadreseul}
\null\vfill\vfill
\hfill{\it to James Brewer}
\vfill
\null
\newpage

\endgroup

\let\oldshowchapter\showchapter
\let\oldshowsection\showsection
\def\showchapter#1{\edef\temp{\thechapter}\def\ttemp{#1}\ifx\ttemp\temp\relax\def\showsection##1{##1}\else\let\showsection\oldshowsection\oldshowchapter{#1}\fi}
\chapter*{Preface of the French edition}
\markboth{Preface}{Preface}

This book is an introductory course to basic commutative \alg with a particular emphasis on \mptfsz, which constitutes the \agq version of the 
vector bundles in \dile \gmtz.

We adopt the \cof point of view, with which all existence \thos have an explicit \algq content.
 In particular, when a \tho affirms the existence of an object -- the solution of a \pb --
a construction \algo of the object can always be extracted from the given \demz.

We revisit with a new and often simplifying eye 
several abstract classical theories.  
In particular, we review theories which did not have any \algq content in their \gnl natural framework, such as Galois theory, the \adksz, 
the \mptfs or the \ddkz. 

\Cov \alg is actually an old discipline,
developed among others by Gauss and Kronecker.
We are in line with the modern \gui{bible} on the subject,  
which is the book by Ray Mines, Fred Richman and Wim Ruitenburg, {\it  A Course in
Constructive Algebra}, published in 1988.
We will cite it in abbreviated form~\cite{MRR}.

This work corresponds to an MSc graduate level, at least up to Chapter~\ref{chapNbGtrs}, but only requires  as prerequisites the basic notions concerning group theory, \lin \alg over fields,
 \detersz,
modules over commutative \risz, as well as the \dfn 
of quotient and localized \risz.
A familiarity with \pol \risz, the \ari \prts of $\ZZ$ and Euclidian \ris 
is \egmt desirable. 

Finally, note that we consider the  exercises and  \pbs (a little over $320$ in total) 
as an essential part of the book.

We will try to publish the maximum amount of missing solutions, as well as additional exercises on the web page of one of the authors:\\ 
\url{http://hlombardi.free.fr/publis/LivresBrochures.html}

\subsection*{Acknowledgements.} We would like to thank all the colleagues who encouraged us in our project, gave us some truly helpful assistance or provided us with valuable information. Especially  MariEmi Alonso, Thierry Coquand, Gema D\'{\i}az-Toca, Lionel Ducos,  M'hammed El Kahoui, Marco Fontana, Sarah Glaz, Laureano Gonz\'alez-Vega, Emmanuel Hallouin,  Herv\'e Perdry, Jean-Claude Raoult, Fred Richman, Marie-Fran\c{c}oise Roy,
  Peter Schuster and Ihsen Yengui. 
Last but not least, a special mention for our \LaTeX~expert, Fran\c{c}ois P\'etiard.

\goodbreak
Finally, we could not forget to mention the Centre International de Recherches Math\'ematiques \`a Luminy and the Mathematisches Forschungsinstitut Oberwolfach, 
who welcomed us for research visits during the preparation of this book, offering us invaluable working conditions. 

\vspace{5mm}
\begin{flushright}
Henri Lombardi, Claude Quitt\'e\\
August 2011
\end{flushright}

\chapter*{Preface of the English edition}

\vspace{3em}
In this edition, we have corrected the errors that we either found ourselves or that were signalled to us. 

We have added some exercise solutions as well as some additional content. Most of that additional content is corrections of exercises, or new exercises or \pbsz.

The additions within the course are the following. A paragraph on the null tensors added as the end of Section \ref{secStabPf}. 
The paragraph on the quotients of flat modules at the end of Section \ref{secPlatDebut} has been fleshed out. 
We have added Sections \iref{secPlgcor} and \iref{secPlgprof2} 
in Chapter \ref{chapPlg} devoted to the \plgsz.  

\medskip 
None of the numbering has changed, except for the \plg XII-7.13 which has become \ref{plcc.ddk}. 

\medskip There are now 297 exercises and 42 \pbsz.

\medskip Any useful precisions are on the site:

\url{http://hlombardi.free.fr/publis/LivresBrochures.html}

\medskip 
\subsection*{Acknowledgements.}
We cannot thank Tania K. Roblot enough for the work achieved translating the book into English. 

\vspace{5mm}
\begin{flushright}
Henri Lombardi, Claude Quitt\'e\\
May 2014
\end{flushright}

\subsection*{This is the web updated version of the book}
\label{labcor}

Except for the corrections indicated below, it is the same text as the one of the printed book. The unique structural modifications concern the table of contents: the general table of contents is shortened, and there is a detailed table of contents at the beginning of each chapter. 

\subsection*{Corrections to the printed book}
Chapter VI: the title of section VI-5 is fixed as \gui{Dualizing \lin forms, \ste \algs} instead of \gui{Dualizing \lin forms, \stf \algs}
 
Solution of Problem 3 in Chapter XII: \paref{allnonzero} replace \gui{all nonzero} by \gui{not all zero}.

Chapter XIII. Exercise 17 item \emph{3}. The solution is changed.

In Section XV-9, the proof of Lemma 9.3 is not correct. It is necessary to give directly a proof of the $(a,b,(ab))$ trick for depth 2, allowing us to prove the \plgcz. So  9.3, 9.4, 9.5 and 9.6 become 9.6, 9.3, 9.4 and 9.5.

More details on \url{http://hlombardi.free.fr/publis/LivresBrochures.html}

\cleardoublepage
\setcounter{tocdepth}{1}
\tableofcontents



~

\vskip1cm
\rdb
{\LARGE \bf Foreword}
\pagestyle{CMExercicesheadings}
\thispagestyle{CMchapitre}
\addcontentsline{toc}{chapterbis}{Foreword}
\markboth{Foreword}{Foreword}

\vskip1cm

\begin{flushright}
{\small
 Quant  \`a moi, je proposerais de s'en tenir
 aux r\`egles
suivantes:


\item 1. Ne jamais envisager que des objets susceptibles
d'\^etre d\'efinis
\\
en un nombre fini de mots;
\\
\item 2. Ne jamais perdre de vue que toute proposition
 sur l'infini doit
\\
\^etre la traduction, l'\'enonc\'e abr\'eg\'e
de propositions sur le
fini;
\\
\item 3. \'Eviter les classifications et les d\'efinitions
non pr\'edicatives.


\medskip \rm Henri Poincar\'e,

dans  {\it La logique de l'infini }
(Revue de M\'etaphysique et de Morale, 1909).
\\
R\'e\'edit\'e dans  {\it
Derni\`eres pens\'ees}, Flammarion.}{\footnote{The official translation by John W. Bolduc (1963) is as follows: \gui{As for me, I would propose that we be guided by the following rules:
\item 1. Never consider any objects but those capable of being defined
in a finite number of words;

\item 2. Never lose sight of the fact that every proposition concerning infinity must
be the translation, the precise statement of propositions concerning the finite;

\item 3. Avoid nonpredicative classifications and definitions.}
}}
\end{flushright}

\vspace{.3cm}

This book is an introductory course to basic commutative \alg with a particular emphasis on \mptfsz, which constitutes the \agq version of the 
vector bundles in \dile \gmtz.

As indicated in the preface, we adopt the \cov method, with which all existence \thos have an explicit \algq content.
\Coma can be seen as the most theoretical branch of \calfz, which handles \maths which \gui{run on a computer.} 
Our course is nevertheless distinguishable from usual \calf courses in two key aspects.

First of all, our \algos are often only implied, 
underlying the \demz, and are in no way optimized for the fastest execution,
as one might expect when aiming for an efficient implementation.

Second, our theoretical approach is entirely \covz, whereas  \calf courses typically have little concern for this issue.
The philosophy here is then not, as is customary \gui{black or white, the good cat is one that catches the mouse}{\footnote{Chinese proverb.}} but rather follows \gui{Truth includes not only the result but also the path to it. The investigation of truth must itself be true; true investigation is developed truth, the dispersed elements of which are brought together in the result.}{\footnote{Karl Marx, Comments on the latest Prussian censorship instruction, 1843 (cited by Georges Perec in \emph{Les Choses}); transcribed here as by Sally Ryan on \url{http://www.marxists.org/archive/marx/works/1842/02/10.htm}.}}

We often speak of two points of view on a given subject: classical and \cofz.
In particular, we have marked with a star the statements 
(theorems, lemmas, \dots) which are true in \clamaz, but for which we do not give a \cov \dem and which often cannot have one.
These \gui{starred} statements will then likely never be implemented on a machine, but are often useful as intuition guides, and to at least link with the usual presentations written in the style of \clamaz.

As for the \dfnsz, we \gnlt first give a \cov variant, 
even if it means showing the equivalence with the usual \dfn in the context of \clamaz.

The reader will notice that in the \gui{starred} \dems we freely use Zorn's lemma and the Law of Excluded Middle (\TEMz)\footnote{The Law of the Excluded Middle 
 states that $P\vee \lnot P$ is true for every proposition $P$. This principle is accepted in \clamaz. See \paref{FirstTEM} for a first explanation regarding the refusal of \TEM in \comaz.}, whereas the other \dems always have a direct translation into an \algoz.

\Cov \alg is actually an old discipline,
developed by Gauss and Kronecker,  among others.
As also specified in the preface, we are in line with the modern \gui{bible} on the subject,  
which is the book by Ray Mines, Fred Richman and Wim Ruitenburg, {\it  A Course in
Constructive Algebra}, published in 1988.
We will cite it in abbreviated form~\cite{MRR}.
Our work is however self-contained and we do not demand \cite{MRR} as a prerequisite.
The books on \coma by Harold M.~Edwards \cite{Ed,Ed2}  and the one of Ihsen Yengui \cite{Yengui}  are also recommended.

\subsection*{The work's content}

We begin with a brief commentary
 on the choices that have been made regarding the covered themes. 

The theory of \mptfs is one of the unifying themes of this work.  
We see this theory in abstract form as an \agq theory of vector bundles, and in concrete form as that of \idm matrices. 
The comparison of the two views is sketched in the introductory chapter.

The theory of \mptfs itself is treated in  
Chapters~\ref{chap ptf0} (first \prtsz),~\ref{chap AlgStricFi} (\algs which are \mptfsz),~\ref{chap ptf1} (rank theory and examples),~\ref{chapNbGtrs} (Serre's Splitting Off theorem) 
and~\ref{ChapMPEtendus} (extended \mptfsz).

Another unifying theme is provided by \plgsz, as in \cite{Kun} for example. It is a highly efficient conceptual framework, even though it is a little vague.
From a \cof point of view, we replace the \lon at an arbitrary \idep with a finite number of \lons at \mocoz. 
The notions which respect the \plg are considered \gui{good notions,} in the sense that they are ready 
for the passage of commutative rings to Grothendieck schemes, which we will unfortunately be unable to address due to the restricted size of this book.

Finally, one last recurrent theme is found in the method, quite common in computer algebra, 
called \emph{the lazy \evnz}, or in its most advanced form, \emph{the dynamic \evnz} method. This method is necessary if one wants to set up an \algq processing of the questions which a priori require the solution to a \fcn \pbz. 
This method has \egmt led to  the development of the \lgb \cov machinery found in Chapters~\ref{chap mpf} 
and~\ref{chapPlg}, as well as the \cov theory of the \ddk (Chapter~\ref{chapKrulldim}), with important applications in the last chapters.

We now proceed to a more detailed description of the contents of the book.

\smallskip In Chapter~\ref{chapMotivation}, 
we explain the close relationship that can be established between the notions of vector bundles in \dile \gmt and of \mptfs in commutative algebra. 
This is part of the \gnl algebraization process in \mathsz, 
a process that can often simplify, abstract and generalize surprisingly well concepts from particular theories.

\smallskip Chapter~\ref{chapSli} is devoted to \slis  over a commutative ring, treated as \elrz. 
It requires almost no theoretical apparatus, 
apart from the question of \lon at a \moz, of which we give a reminder in Section~\ref{secPrelimCh2}.
We then get to our subject matter by putting in place the \plgc for solving \slis
(Section~\ref{secPLGCBasic}),
a simple and effective tool that will be repeated and varied constantly. 
From a \cof point of view, solving \slis \imdt renders as central the concept of \coh rings that we treat in Section~\ref{secAnneauxCoherents}. \Coh rings are those for which we have a minimal grip on the solution of homogeneous \slisz. Very surprisingly, this concept does not appear in the classical commutative \alg treatises. That is because in general the concept is completely obscured by that of a \noe ring. This obscuration does not occur in \coma where \noet does not \ncrt imply \cohcz. We develop in Section~\ref{sec sfio} the question of finite products of \risz, with the notion of a \sfio and the Chinese Remainder \thoz. The long Section~\ref{secCramer} is devoted to many variations on the theme of \detersz. Finally, Section~\ref{secPLGCBasicModules} returns to the basic \plg in a slightly more \gnl version devoted to exact sequences of modules.

\smallskip Chapter~\ref{chapGenerique} develops the method of indeterminate \coesz, first developed by Gauss. 
Numerous \thos of existence in commutative \alg rely on \gui{\idas under conditions} and thus on memberships 
$g\in\gen{\lfs}$ in a \riz~$\ZZ[c_1,\ldots,c_r,\Xn]$, where the $X_i$'s are the variables and the $c_j$'s are the parameters of the \tho under consideration. In this sense, we can consider that commutative \alg is a vast theory of \idasz, which finds its natural framework in the method of indeterminate \coesz, \cad the method 
in which the parameters of the given \pb are taken as \idtrsz.
In that assurance we are, to the extent our powers allow to, systematically \gui{chasing \idasz.}
 This is the case not only in the \gui{purely computational} Chapters~\ref{chapSli} and~\ref{chapGenerique}, but throughout the book. In short, rather than simply assert in the context of an existence theorem \gui{there is an algebraic identity which certifies this existence,} 
we have tried each time to give the algebraic identity itself. 

Chapter~\ref{chapGenerique} can be considered as a basic algebra course with $19^{\textrm{th}}$ century methods. 
Sections~\ref{secAnnPols},~\ref{secLemArtin} and~\ref{secThKro} provide certain \gnts about \polsz,
featuring in particular the algorithm for partial factorization, the \gui{theory of \idasz}
(which explains the method of indeterminate \coesz), the \elrs \smqs \polsz, the \DKM lemma and the \KROz's \thoz.  The last two results are basic tools which give precise information on the \coes of the product of two \polsz; they are often used in the rest of this manuscript. Section~\ref{sec0adu} introduces the \adu of a \mon \pol over an arbitrary commutative \riz, which is an efficient substitute for the field of the roots of a \pol over a field. 
Section~\ref{secDisc} is devoted to the discriminant and explains in what precise sense a \gnq matrix is \digz. With these tools in hand, we can treat the basic Galois theory in Section~\ref{secGaloisElr}. The \elr theory of elimination via the resultant is given in Section~\ref{secRes}. We can then give the basics of \agq number theory with the \tho of unique \dcn into prime factors for a \itf of a number field (Section~\ref{secApTDN}). Section~\ref{secChap3Nst} shows Hilbert's \nst as an application of the resultant. Finally, Section~\ref{secNewton} on Newton's method in algebra closes this chapter.

\smallskip Chapter~\ref{chap mpf} is devoted to the study of the \elrs \prts of \mpfsz. 
These modules play a role for rings similar to that played by finite dimensional vector spaces for fields: the theory of \mpfs  is a more abstract, and often  profitable, way to address the issue of \slisz.
Sections~\ref{sec pf chg} to~\ref{secStabPf} show the basic stability \prts as well as the important example of the \id of a zero 
for a \syp (on an arbitrary commutative ring). 
We then focus on the classification \pb of \mpfs over a given \riz. 
Working towards principal ideal domains (PID), 
for which the classification \pb is completely solved 
(Section~\ref{secBézout}), we will encounter \qiris 
(Section~\ref{subsecAnneauxqi}), which are the \ris where the annihilator of an \elt is always generated by an \idmz. 
This will be the opportunity to develop 
an \emph{\elr \lgbe machinery} 
which conveniently reformulates a \cot established result for integral \ris into 
the analogous result for pp-\risz.\imlgd
This proof-rewriting machinery is \elr as it is founded on the \dcn of a \ri into a finite product of \risz.
The interesting thing 
is that this \dcn is obtained via 
a rereading of the \cov \dem written in 
the integral case; here we see that in \coma the \dem is often even more important than the result.
Similarly, we have an \emph{\elr \lgbe machinery\imlgz} which conveniently reformulates a \cot established result for \cdis into 
the analogous result for \zedrs \ris (Section~\ref{secKrull0dim}). 
The \zeds \ris elementarily defined here constitute 
an important intermediate step to \gnr specific results regarding \cdis to arbitrary commutative \risz: they are a 
key tool of commutative \algz.
Classically, these appear in the literature in their \noee form, \cad that of Artinian \risz. 
Section~\ref{sec Fitt} introduces very important invariants: the \idfs of a \mpfz. 
Finally, Section~\ref{subsecIdealResultant} applies this notion to introduce 
the resultant \id of a \itf over a \pol \ri when the \id contains a monic \polz, and to prove a \tho of \agq \eli over an arbitrary \riz.

\smallskip Chapter~\ref{chap ptf0} is a first approach to the theory of \mptfsz.
Sections~\ref{subsecPropCarPTFS} to~\ref{secPtfCoNat} state basic \prts along with the important example of the \zeds \risz. 
Section~\ref{secMPTFlocLib} states the local structure \thoz : 
a module is \ptf \ssi it becomes free after \lon at suitable \ecoz. Its \cov \dem is a rereading of a result established in Chapter~\ref{chapSli} for \gui{well conditioned} \slis (\thref{theoremIFD}).
Section~\ref{Idpp} develops the example of the \lmos \pro modules. Section~\ref{subsec det ptf} introduces the \deter of an \endo of a \mptfz. 
This renders the \dcn of such a module into a direct sum of its components of constant rank accessible. 
Finally, Section~\ref{secSalutFini}, which we were not too sure where to place in this work, 
hosts some \suls considerations on \emph{\prts of \carfz}, a concept introduced in Chapter~\ref{chapSli} to discuss the connections between concrete \plgs and abstract \plgsz.

\smallskip Chapter~\ref{chap AlgStricFi} is essentially devoted to \algs which are \mptfs over their base \risz. 
We call these \asfsz. 
When applied to commutative \risz, they constitute a natural \gnn of the concept of a finite \alg over a field. 
The icing on the cake being the important case of \aGsz, which generalize Galoisian extensions of \cdis to commutative \risz.

Section~\ref{secEtaleSurCD} treats the case where a base \ri is a \cdiz.
It provides the \covs versions of the structure \tho obtained in \clamaz.
The case of \'etale \algs (when the \discri is \ivz) is particularly enlightening. We discover that the classical \thos always implicitly assume that we know how to factorize the \spls \pols over the base field. 
The \cov \dem of the primitive \elt \tho \ref{thEtalePrimitif} is significant for its deviation from the classical \demz.
Section~\ref{sec2GaloisElr} applies the previous results to complete the basic Galois theory started in Section~\ref{secGaloisElr} by characterizing the Galoisian extensions of \cdis like the \'etale and normal extensions. 
Section~\ref{sec1Apf} is a brief introduction to \apfsz, by focusing on  
integral \algsz 
\footnote{By \gui{integral \algz} we mean \emph{an \alg that is  integral on its base \riz}, not to be confused with an \alg that is an integral \riz.}, with a weak \nst and the Lying Over lemma.
Section~\ref{subsecAlgStfes} introduces \stfes \algs over an arbitrary \riz. 
In Sections~\ref{secAlgSte} and~\ref{secAlgSpb}, the related concepts of a \ste \alg and of a \spb \alg are introduced. These generalize the concept of an \'etale \alg over a \cdiz. 
In Section~\ref{secAGTG}, we \cofz ly present the basics of the theory of \aGs for commutative \risz.
It is in fact an Artin-Galois theory since it adopts the approach Artin had developed for the case of fields, starting directly from a finite group of \autos of a field, the base field appearing only as a byproduct of subsequent constructions.

\smallskip In Chapter~\ref{ChapGalois},
the dynamic method -- a cornerstone of modern methods in \cov \alg -- is implemented to deal with the field of roots of a \pol and the Galois theory in the \spb case. From a \cof point of view, 
we need to use the dynamic method when 
we do not know how to factorize the \pols over the base field. \\
For training purposes, 
Section~\ref{secNstSCA} 
begins by establishing results in a \cov form for the \nst when we do not know how to factorize the \pols over the base field. \Gnl considerations on the dynamic method are developed in Section~\ref{subsecDyna}.
More details on the course of the festivities are given in the introduction of the chapter.

\smallskip Chapter~\ref{chap mod plats} is a brief introduction to flat modules and to flat and \fptes \algsz.
Intuitively speaking, an \Alg $\gB$ is flat when the homogeneous \slis over $\gA$ have \gui{no more} solutions in 
$\gB$ than in $\gA$, and it is faithfully flat if this statement is \egmt true for nonhomogeneous \slisz.
These crucial notions of commutative \alg were introduced by Serre in~\cite[GAGA,1956]{Serregaga}. We will only state the truly fundamental results. This is \egmt when we will introduce the concepts of a 
\lsd (i.e. a \ri whose principal \ids are flat),
of a \torf module (for an arbitrary \riz),
of an \anar and of a Pr\"ufer \riz. As always, we focus on the \plg when it applies.

\smallskip Chapter~\ref{chap Anneaux locaux} discusses local \ris and some \gnnsz. Section~\ref{secAloc1} introduces the \cov terminology for some common classical concepts, including the important concept of a Jacobson radical. A related concept is that of a \plc \ri (a \ri $\gA$ such that $\gA/\Rad\gA$ is \zedz). It is a robust concept, which never uses \idemasz, and most of the \thos in the literature regarding semi-local \ris (in \clama they are the \ris which only have a finite number of \idemasz) apply to \plcs \risz. 
Section~\ref{secAloc2} lists some results which show that on a \alo we reduce the solution of particular \pbs to the case of fields. 
Sections~\ref{secLoc1+fa} and~\ref{secExlocGeoAlg} establish, based on \gmqs examples (\cad regarding the study of \sypsz), a link between the notion of a local study in the topological intuitive sense and the study of certain \lons of \ris 
(in the case of polynomial systems over a discrete field these \lons are \alosz).
In particular we introduce the notions of tangent and cotangent spaces at zero of a \sypz.
Section~\ref{secRelIdm} is a brief study of \dcps \risz, including the particular case from \clama of decomposed \ris (finite products of \alosz), which play an important role in the theory of Henselian \alosz.
Finally, Section~\ref{secAlocglob} treats the notion of a \lgb \riz,
which generalizes both the concept of \alos and that of \zed \risz. These \ris verify very strong \lgb \prtsz; e.g. the \mrcs are always free. Moreover, the class of \lgb \ris is stable under integral extensions.

\smallskip Chapter~\ref{chap ptf1} continues the study of \mptfs started in Chapter~\ref{chap ptf0}.
In Section~\ref{sec ptf loc lib}, 
we return to the question of the \carn of \mptfs as  \lot free modules, 
\cad of the local structure \thoz. 
We give a matrix version of it (\thref{th matproj}), which summarizes and clarifies the different statements of the \thoz. 
Section~\ref{subsecCalRang} is devoted to the \ri of ranks over~$\gA$.
In \clamaz, the rank of a \mptf is defined as a \lot constant function in the Zariski spectrum. 
We give here an \elr theory of the rank which does not call upon \idepsz. 
In Section~\ref{secAppliLocPtf}, we provide some simple applications of the local structure \thoz. 
Section~\ref{secGrassman} introduces Grassmannians.
In Section~\ref{subsecClassifMptfs}, we introduce the \gnl \pb of the complete classification of the \mptfs over a fixed \ri $\gA$. 
This classification is a fundamental and difficult \pbz, which does not have a \gnl \algq solution. 
Section~\ref{secAppliIdenti} presents a nontrivial example for which this classification can be obtained.

\smallskip Chapter~\ref{chapTrdi} is devoted to \trdis and 
lattice ordered groups (l-groups).
The first two sections describe these \agqs structures along with their basic \prtsz. These structures are important in commutative \alg for several reasons.

First, the \dve theory has as its \gui{ideal model} the natural numbers' \dve theory. The structure of the multiplicative \moz~\hbox{$(\NN\etl,\times ,1)$} makes it the positive part of an \grlz. In commutative \algz, this can be generalized in two possible ways. The first \gnn is the theory of integral \ris whose \itfs form a \trdiz, called \ddpsz, which we will study in Chapter~\ref{ChapAdpc}; 
their nonzero \itfs form the positive part of an \grlz.
The second is the theory of gcd \ris that we study in Section~\ref{secAnnPgcd}. 
Let us notify the first appearance of the Krull dimension $\leq 1$ in \Thref{propGCDDim1}: an integral gcd \ri with dimension $\leq 1$ is a Bézout \riz.

Secondly, the \trdis act as the \cov counterpart of various 
spectral spaces which have emerged as powerful tools of abstract \algz.
The relationship between \trdis and spectral spaces will be discussed in Section~\ref{secEspSpectraux}.
In Section~\ref{secZarAcom}, we set up the Zariski lattice of a commutative \ri $\gA$, which is the constructive counterpart of the famous Zariski spectrum. Our goal here is to establish a parallel between the construction of the \zede reduced closure of a \ri (denoted by~$\Abul$) and that of the \agB generated by a \trdi (which is the subject of \Thref{thZedGenEtBoolGen}).
The object $\Abul$ constructed as above essentially contains the same information as the product of the \ris $\Frac(\gA\sur\fp)$ for all \ideps $\fp$ of $\gA$(\footnote{This product is not accessible in \comaz, $\Abul$ is its perfectly effective \cof substitute.}). This result is closely related to the fact that the Zariski lattice of $\Abul$ is the \agB generated by the Zariski lattice of $\gA$.

A third reason to be interested in \trdis is \cov logic (or intuitionistic logic). In this logic, the set of truth values of classical logic, that is the two-element \agB $\so{\Vrai,\Faux}$, 
is replaced by a quite mysterious \trdiz. 
\Cov logic is informally discussed in the Annex. 
In Section~\ref{secEntRelAgH}, we set up the tools that provide a framework for a formal \agq study of \cov logic: 
\entrels and \agHsz.
In addition, \entrels and \agHs have their own use in the \gnl study of \trdisz. For example, the Zariski lattice of a \coh \noe \ri is a \agH (Proposition~\ref{propNoetAgH}).

\smallskip Chapter~\ref{ChapAdpc} deals with \anarsz, \adps and \adksz. \Anars are \ris for which the lattice of \itfs is distributive.
A \adp is a reduced \anar and is characterized by the fact that all of its \ids are flat. A \adpc is the same thing as an arithmetic pp-ring. 
It is characterized by the fact that its \itfs are projective. A \adk is a \noe and \fdi \adpc (in \clamaz,  with \TEMz, every \ri   is \fdi and every \noe \ri is \cohz). These \ris first appeared as 
the rings of integers of number fields.  
The paradigm in the integral case 
is the unique \dcn into prime factors of any nonzero \itfz. The general \aris \prts of \itfs are mostly verified by all the \anarsz.
For the most subtle \prts concerning the \fcns of the \itfsz, and in particular the \dcn into prime factors, a \noee assumption, or at least a dimension~$\leq 1$ assumption, is essential.
In this chapter, we wanted to show the progression of the \prts satisfied by the \ris as we strengthen the assumptions from the \anars to the \fac\adksz. 
We focus on the simple \algq character of the \dfns in the \cof framework. 
Certain \prts only depend on dimension $\leq 1$, and we wanted to do justice to pp-\ris \ddi1. We \egmt carried out  a more progressive and more elegant study of the \pb of the \dcn into prime factors than in the presentations which allow \TEMz.
For example, \Thosz~\ref{thAESTE} and~\ref{lemthAESTE} provide precise \cov versions of the \tho concerning the normal finite extensions of \adksz, with or without the total \fcn\prtz.

The chapter begins with a few epistemological remarks on the intrinsic interest of addressing the \fcn \pbs with the  \fap \tho rather than the total \fcn one. 
To get a good idea of how things unfold, simply refer to the table of contents at the beginning of the chapter on \paref{ChapAdpc} and to the table of \thos on \paref{tdtChapAdpc}.

\smallskip Chapter~\ref{chapKrulldim} is devoted to the \ddk of commutative \risz, of their morphisms and of \trdisz, and to the valuative dimension of commutative \risz. 

Several important notions of dimension in classical commutative \alg are dimensions of spectral spaces. These very peculiar topological spaces have the \prt of being fully described (at least in \clamaz) by their compact-open subspaces, which form a \trdiz. It so happens that the corresponding \trdi \gnlt has a simple interpretation, without any recourse to spectral spaces. In 1974, Joyal showed how to \cot define the \ddk of a \trdiz. Since this auspicious day, the theory of dimension which seemed bathed in ethereal spaces -- that are invisible when you do not trust the axiom of choice -- has become (at least in principle) an \elr theory, 
without any further mysteries. 

Section~\ref{secEspSpectraux} describes the approach of the \ddk in \clamaz.
It also explains how to interpret the \ddk of such a space in terms of the \trdi of its \oqcsz.
Section~\ref{secDefConsDimKrull} states the \cov \dfn of the \ddk of a commutative \riz, denoted by $\Kdim\gA$, and draws some consequences. 
Section~\ref{secKrullElem} states some more advanced \prtsz, in particular the \plg and the \prf for the \ddkz. 
Section~\ref{secDDKExtEn} deals with the \ddk of integral extensions and Section~\ref{secDimGeom} that of \gmqs \ris (corresponding to \sypsz) on the \cdisz.
Section~\ref{secDDKTRDIS} states the \cov \dfn of the \ddk of a \trdi and shows that the \ddk of a commutative \ri and that of its Zariski lattice coincide.
Section~\ref{secKdimMor} is devoted to the dimension of the morphisms between commutative \risz. The \dfn uses the reduced \zed closure of the source \ri of the morphism. To prove the formula which defines 
the upper bound of $\Kdim \gB$ from
$\Kdim \gA$ and $\Kdim\rho$ (when we have a morphism $\rho:\gA\to\gB$),
we must introduce the minimal pp-closure of a commutative \riz. This object is a \cov counterpart of the product of all the $\gA/\fp$, when~$\fp$ ranges over the \idemis of~$\gA$.
Section~\ref{secValdim} introduces the valuative dimension of a commutative \ri and in particular uses this concept to prove the following important result: 
for a nonzero \anar $\gA$, we have $\Kdim\AXn=n+\Kdim\gA$.
Section~\ref{secGoingUp} states \cov versions of the Going up and Going down \Thosz.

\smallskip In Chapter~\ref{chapNbGtrs}, titled \emph{Number of \gtrs of a module}, we establish the \elrz, non-\noee and \cov versions of the \gui{great} \thos of commutative \algz, their original form due to \KRNz, Bass, Serre, Forster and Swan. These results relate to the number of radical \gtrs of a \itfz, the number of \gtrs of a module, the possibility of producing a free submodule as a direct summand 
in a module, and the possibility to simplifying \isosz, in the following way: if $M\oplus N\simeq M'\oplus N$ then $M\simeq M'$.
They involve the \ddk or other, more sophisticated dimensions introduced by R. Heitmann as well as by the authors of this work and T. Coquand.

Section~\ref{secKroBass} is devoted to Kronecker's \Tho and its extensions 
(the most advanced, non-\noeez, is due to R. Heitmann \cite{Hei84}).
\KRNz's \Tho is usually stated in the following form: an \vgq in $\CC^n$ can always be defined by $n+1$ \eqnsz. The form due to Heitmann is that in a \ri of \ddk less than or equal to $n$, for all \itf $\fa$ there exists an \id $\fb$ generated by at most $n+1$ \elts of~$\fa$ such that $\sqrt \fb=\sqrt\fa$.
The \dem also gives Bass' stable range \Thoz.  
The latter \tho was improved by involving \gui{better} dimensions than the \ddkz. This is the subject of Section~\ref{subsecDimHeit} where the \emph{Heitmann dimension} is defined, discovered while carefully reading Heitmann's \dems (Heitmann uses another dimension, a priori a little worse, which we \egmt explain in \cofs terms).
In Section \ref{secSOSFSa}, we explain which matrix \prts of a \ri allow \SSO \thoz, Forster-Swan's \tho 
(controlling the number of \gtrs of a \mtf according to the local number of \gtrsz) and Bass' simplification \thoz.
Section~\ref{secSUPPORTS} introduces the concepts of support 
(a mapping from a \ri to a \trdi satisfying certain axioms) and of $n$-stability.
The latter was defined by Thierry Coquand, after having analyzed one of Bass' \dems which establishes that the \mptfs over a \riz~$\gV[X]$, where~$\gV$ is a valuation \ri of finite \ddkz, are free.
In the final section, we prove that the crucial matrix \prt introduced in Section~\ref{secSOSFSa} is satisfied, on one hand by the $n$-stable \risz, and on the other by the \ris of Heitmann dimension~$<n$.

\smallskip Chapter~\ref{chapPlg} is devoted to the \plg and its variants.
Section~\ref{subsecMoco} introduces the notion of the covering of a \mo by a finite family of \mosz, which \gnss the notion of \mocoz. 
The covering Lemma~\ref{lemRecouvre} will be decisive in Section~\ref{secMachLoGlo}.
Section~\ref{subsec loc glob conc} states some \plgcsz. This is to say that some \prts are globally true as soon as they are locally true.
Here, \gui{locally} is meant in the \cof sense: after \lon at a finite number of \mocoz. Most of the results have been established in the previous chapters. Grouping them shows the very broad scope of these principles. Section~\ref{subsec loc glob abs} restates some of these principles as \plgasz. Here, \gui{locally} is meant in the abstract sense: after \lon at any arbitrary \idepz. 
We are mainly interested in comparing  the abstract principles and the corresponding \plgcsz.
Section~\ref{secColleCiseaux} explains the construction  of \gui{global} objects from objects of the same type defined only locally, as is usual in \dile \gmtz. It is the impossibility of this construction when seeking to glue certain \ris together which is at the root of Grothendieck schemes. 
In this sense, Sections~\ref{subsec loc glob conc} and~\ref{secColleCiseaux}
constitute the basis from which we can develop the theory of schemes in a completely \cof framework.

The following sections are of a different nature. 
Methodologically, they are devoted to the decryption of different variations of the \plg in \clamaz.
For example, the \lon at every \idepz, the passage to the quotient by every \idema or the \lon at every \idemiz, each of which applies in particular situations. 
Such a decryption certainly presents a confusing character insofar as it takes as its starting point a classical \dem that uses \thos in due and proper form, but where the \cof decryption of this \dem is not only given by the use of \cof \thos in due and proper form. 
One must also look at what the classical \dem does with its purely \id objects (e.g. \idemasz) to understand how it gives us the means to construct a finite number of \elts that will be involved in a \cof \tho (e.g. a \plgc) to reach the desired result.
Decrypting such a \dem we use the \gnl dynamic method presented in Chapter~\ref{ChapGalois}. We thus describe \emph{local-global machineries} that are significantly less \elrs than those in Chapter~\ref{chap mpf}: the basic \cov \lgbe machinery \gui{with \idepsz} (Section~\ref{secMachLoGlo}), the \cov \lgbe machinery \gui{with \idemasz}  (Section~\ref{subsecLGIdeMax}) and the \cov \lgbe machinery \gui{with \idemisz} (Section~\ref{subsecLGIdepMin})\imlbz.  
By carrying out \gui{Poincar\'e's program} 
used as an epigraph for this foreword, our \lgb machineries take into account an essential remark made by Lakatos that the most interesting and robust thing in a theorem is always its proof, even if it can be criticized in some respects (see \cite{La}).
  
In Sections \ref{secPlgcor} and \ref{secPlgprof2}, we examine to what extent certain \plgs remain valid when we replace in the statements the lists of \eco by lists of depth $\geq 1$ or of depth $\geq 2$.

\smallskip In Chapter~\ref{ChapMPEtendus}, we treat the question of \mptfs over \ris of \polsz. The decisive question is to establish for which classes of \ris the \ptfs modules over a \pol \ri  are derived by \eds of a \ptf module over the \ri itself (possibly by putting certain restrictions on the considered \mptfs or on the number of variables in the \pol \riz). Some \gnts on the extended modules are given in Section~\ref{sec etendus}.
The case of the \mrcsz~$1$, which is fully clarified by Traverso-Swan-Coquand's \thoz, is dealt with in Section~\ref{sec.Traverso}.
Coquand's \cov \dem uses the \cov \lgb machinery with \idemis in a crucial way.
Section~\ref{subsecQPatch} deals with Quillen and Vaserstein's patching \thosz,
which state that certain objects are obtained by \eds (from the base \ri to the \pol \riz) \ssi this \prt is \lot satisfied. 
We also have a sort of converse to Quillen's patching due to Roitman, in a \cov form. 
Section~\ref{sec.Horrocks} is devoted to Horrocks' \thosz.
The \cov \dem of Horrocks' global \tho is obtained from the \dem of Horrocks' local \tho by using the basic \lgb machinery and concluding with Quillen's \cof patching. 
Section~\ref{sec.QS} gives several \prcos of Quillen-Suslin's \tho (the \mptfs over a \pol \ri on a \cdi are free) founded on different classical \demsz. 
Section~\ref{sec.Etendus.Valuation} establishes Lequain-Simis' \tho (the \mptfs over a \pol \ri on an \anar are extended). 
The \dem uses the dynamic method presented in Chapter~\ref{ChapGalois}. This allows us to establish Yengui's induction \thoz, 
a \cov variation of Lequain-Simis' induction.

\smallskip In Chapter~\ref{ChapSuslinStab}, we prove \gui{Suslin's Stability \Thoz} in the special case of \cdisz. Here also, we use the basic \lgb machinery presented in Chapter~\ref{chapPlg} to obtain a \cov \demz.

\smallskip The Annex describes a Bishop style \cov set theory. 
It can be seen as an introduction to \cov logic. 
In it we explain the Brouwer-Heyting-Kolmogorov semantic for connectives and quantifiers. We discuss certain weak forms of  \TEM along with several problematic principles in \comaz.

\subsection*{Some epistemological remarks}

In this work, we hope to show that classical commutative \alg books such as \cite{Eis}, \cite{Kun}, \cite{Laf}, \cite{Mat},
\cite{Gla}, \cite{Kapl}, \cite{Atiyah}, \cite{Nor}, \cite{Gil}, \cite{Lam06}
(which we highly recommend), or even \cite{Bou}
and the remarkable work available on the web \cite{Stacks},
could be entirely rewritten from a \cof point of view, dissipating the veil of mystery which surrounds the nonexplicit existence \thos of \clamaz.
Naturally, we hope that the readers will take advantage of our work to 
take a fresh look at the classical Computer Algebra 
books like, for instance, \cite{CLS}, \cite{KrRo}, \cite{ElkMo},  \cite{vzGaGe}, \cite{Mora}, \cite{CCS} \linebreak
or \cite{Singular}.

\smallskip 
Since we want an \algq processing of commutative \alg we cannot use all the tricks that arise from the systematic use of Zorn's Lemma and the Law of Excluded Middle in \clamaz.
Undoubtedly, the reader understands that it is difficult to implement Zorn's lemma in Computer Algebra.
\label{FirstTEM}The refusal  of \TEMz, however, must seem
harder to stomach.\index{LEM@\TEM} 
It is simply a practical observation on our part.
If in a classical \dem there is a reasoning that leads to a computation in the form  \gui{if $x$ is invertible, do this, otherwise do that,} then, clearly, it directly translates into an \algo only when there is an invertibility test for the \ri in question.
It is in stressing this difficulty, which we must constantly work around, that we are often led to speak of two points of view on the same subject: classical and \cofz.

We could argue forever about whether \coma is part of \clamaz, the part that deals exclusively with the explicit aspect of things, or conversely whether it is \clama which is a part of \comaz, the part whose \thos are \gui{starred,} \cad which systematically add \TEM and the axiom of choice in their assumptions. One of our objectives is to tip the balance in the second direction, not for  philosophical debate but for practical purposes.

Finally, let us mention two striking traits of this work compared to classical texts on commutative \algz.

The first is that \noet is left on the backburner.
Experience shows that indeed \noet is often too strong an assumption, which hides the true \algq nature of things. 
For example, such a \tho usually stated for \noes \ris and \mtfsz, when its \dem is examined to extract an \algoz, turns out to be a \tho on \coh \ris and \mpfsz.
The usual \tho is but a corollary of the right \thoz, but with two non\cof arguments allowing us to deduce \cohc  and finite presentation from \noet and finite generation in \clamaz. 
A \dem in the more satisfying framework of \cohc and \mpfs is often already published in research articles, although rarely in an entirely \cov form, but \gui{the right statement} is \gnlt missing.{\footnote{This \noee professional bias has produced a linguistic shortcoming in the English literature which consists in taking \gui{local \riz} to mean \gui{\noe local \riz.}}}

The second striking trait of this work is the almost total absence of negation in the \cof statements. For example, instead of stating that for a nontrivial \ri $\gA$, two free modules of respective rank $m$ and $n$ with~$m>n$ cannot be isomorphic, we prefer to say without any assumption about the \ri that if these modules are isomorphic, then the \ri is trivial (Proposition~\ref{propDimMod1}). This nuance may seem quite slight at first, but it has an \algq importance. It will allow us to substitute a \dem from \clama using a \ri $\gA=\gB\sur\fa$, which would conclude that~$1\in\fa$ by contradiction, with a fully \algq \dem that constructs~$1$ as an \elt of the \id $\fa$ from an \iso between $\Ae m$ and~$\Ae n$.

For a \gnl presentation of the ideas which led to the new methods used in \cov \alg in this work, we suggest to the reader the summary article~\cite[Coquand\&Lombardi, 2006]{CL05}.

\vspace{5mm}
\begin{flushright}
Henri Lombardi, Claude Quitt\'e\\
August 2011
\end{flushright}

\newpage
\ifodd\arabic{page}\null\thispagestyle{CMcadreseul}\newpage\else\relax\fi

\null\vfil
\centerline{\includegraphics[width=8cm]{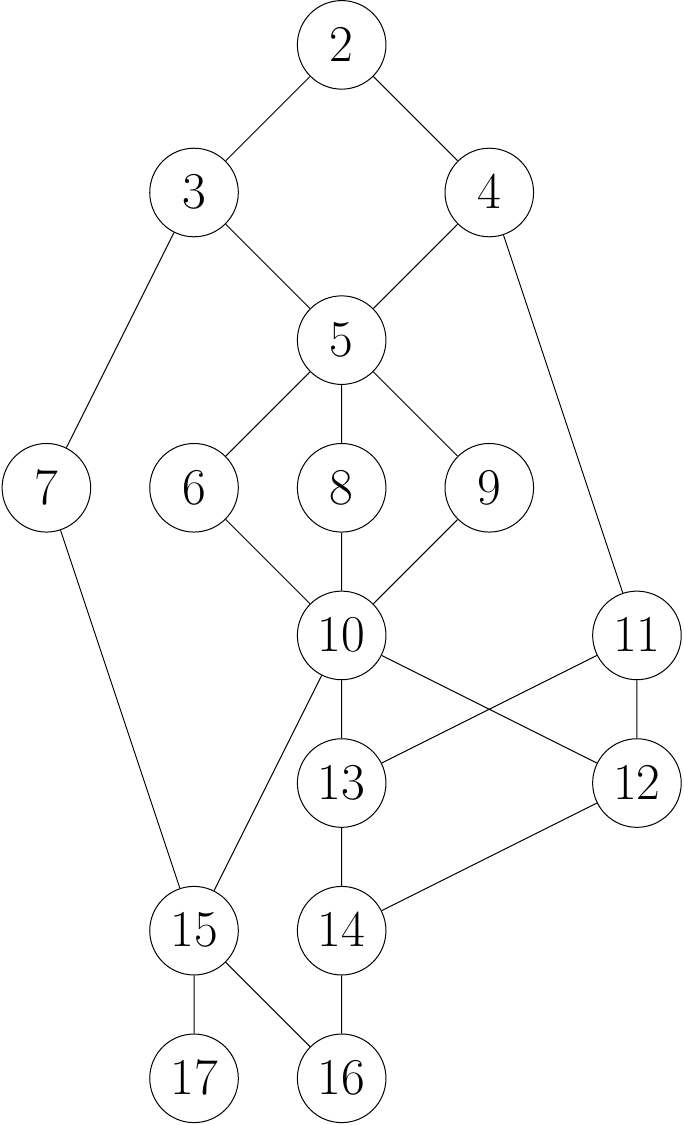}}
\kern3.4truemm

\newpage

The flowchart on the previous page shows the dependence relations between the different chapters

\bigskip 
\begin{enumerate}\setcounter{enumi}{1}
\item \nameref{chapSli}\\
\Coh \ris and modules. A little bit of exterior \algz.
\item \nameref{chapGenerique}\\
Dedekind-Mertens and Kronecker's lemmas. 
Basic Galois theory. Classical \nstz.  
\item \nameref{chap mpf}\\
Category of \mpfsz. Zero-dimensional \risz. 
\Elr \lgb machineries. Fitting \idsz.
\item \nameref{chap ptf0}\\
Local structure \thoz. Determinant. Rank.
\item \nameref{chap AlgStricFi}
\item \nameref{ChapGalois}\\
General \nst (without \agq closure). General Galois theory (without  \fcn \algoz).
\item \nameref{chap mod plats}\\
Flat and \fptes \algsz.
\item \nameref{chap Anneaux locaux}\\
Decomposable \riz. Local-global \riz.
\item \nameref{chap ptf1}
\item \nameref{chapTrdi}\\
GCD \riz. Zariski lattice of a commutative \riz. Entailment relations.
\item \nameref{ChapAdpc}\\
Integral extensions. Dimension $\leq 1$. Factorization of \itfsz.
\item \nameref{chapKrulldim}\\
Krull dimension. Dimension of morphisms. Valuative dimension.
Dimension of integral and \polle extensions.
\item \nameref{chapNbGtrs}\\
Kronecker's, Bass' and Forster-Swan's \thosz. Serre's Splitting Off theorem. Heitmann dimension.
\item \nameref{chapPlg}
\item \nameref{ChapMPEtendus}\\
Traverso-Swan-Coquand's, Quillen-Suslin's, Bass-Lequain-Simiss \thosz.
\item \nameref{ChapSuslinStab}
%
%
\end{enumerate}

\newpage\thispagestyle{CMcadreseul}
\incrementeexosetprob 
\cleardoublepage
\pagestyle{CMheadings}
\arabicpagenumbers
\setpagenumber1


\chapter{Examples} 
\label{chapMotivation}
\minitoc
\subsection*{Introduction} 
\addcontentsline{toc}{section}{Introduction}

Throughout the manuscript, unless explicitly stated otherwise, \ris are commutative and unitary, and a \ri \homo $\varphi:\gA\to\gB$ must satisfy~$\varphi(1_\gA)=1_\gB$.

Let $\gA$ be a \riz.
We say that an \Amo $M$ is  
a \emph{finite rank free module}%
\index{module!free --- of finite rank}
when it is isomorphic to a module $\Ae n$. 
We say that it is a \emph{finitely generated projective} module when there exists an \Amo $N$ such that $M\oplus N$ is a finite rank free module. 
This is equivalent to saying that $M$ is isomorphic to the image of a \emph{projection matrix}%
\index{matrix!projection ---}\index{module!finitely generated projective ---}
 (a matrix $P$ such that $P^2=P$). That is, the \prn matrix 
onto $M$ 
along $N$, 
\prmt defined as follows: 
 $$ 
 M\oplus N \lora M\oplus N, \quad  x+y\longmapsto x \qquad \hbox{for $x\in M$ and $y\in N$}. 
 $$ 
A \mprn is also called a \ix{projector}.

When we have an \iso $M\oplus \gA^\ell\simeq \gA^k$, the \mptf $M$ is called \emph{stably free}.
\index{module!stably free ---}
\index{finitely generated projective!module} 

While over a field or over a PID the \mptfs are free (over a field they are finite dimensional \evcsz),
over a \gnl commutative \ri the classification of the \mptfs is both an important and a difficult \pbz.

Kronecker and Dedekind have proven that a nonzero \itf in the ring of integers of a number field
is always \iv (thus \ptfz), but that it is rarely free (\cad principal). This is a fundamental phenomenon, which is at the root of the modern development of number theory.  

\smallskip 
In this chapter, we try to explain why the notion of a \mptf is important by giving meaningful examples from \dile \gmtz.

The datum of a vector bundle on a smooth compact manifold $V$ is in fact \eqve to the datum of a \mptf over the \ri~\hbox{$\gA= \Cin(V)$} of smooth functions over $V$; to a vector bundle we associate the \Amo of its sections, this  is always a \mptf but it is free only when the bundle is trivial. 

The tangent bundle corresponds to a module built by a purely formal procedure from the \ri $\gA$. 
In the case where the manifold $V$ 
is a sphere, the module of the sections of the tangent bundle is \stlz. An important result about the sphere is that there exist no smooth everywhere nonzero vector fields. 
This is equivalent to the fact that the module of sections of the tangent bundle is not free.    

\medskip 
We try to be as explicit as possible, but in this motivating chapter we freely use the reasonings of \clama without worrying about being completely rigorous from a \cof point of view.

\newpage
\section{Vector bundles on a smooth compact manifold} 
\label{subsecChampVect}

Here, we give some motivations for \mptfs and \lon by explaining the example of vector bundles on a compact smooth manifold. 
Two important particular cases are tangent and cotangent bundles corresponding to $\Cin$ vector fields and to $\Cin$ \dile forms.

We will use the term \gui{smooth} as a synonym for \gui{of class~$\Cin$.}

We will see that the fact that the sphere cannot be combed admits a purely \agq interpretation.

\ss In this section, we consider a smooth real differentiable manifold $V$ 
 and we denote by $\gA= \Cin(V)$  the real \alg of global smooth functions on the manifold.

\Subsubsect{Some \lons of the \alg of continuous functions}{Some \lons} 

Let us first consider an \elt $f\in\gA$ along with the open set (open subset of the manifold $V$ to be precise)

\snic{U=\sotq{x\in V}{f(x)\neq 0}} 

and let us see how we can interpret the \alg $\gA[1/f]$: two \eltsz~$g/f^k$ and $h/f^k$ are equal in $\gA[1/f]$ \ssi for some exponent $\ell$ we have $gf^\ell=hf^\ell$ which means precisely~$g\frt{U}=h\frt{U}.$

It follows that we can interpret $\gA[1/f]$ as a sub-\alg of the \alg of smooth functions on $U$: 
this sub-\alg has as \elts the functions which can be written as $(g\frt{U})/(f\frt{U})^k$ (for a given exponent $k$) with $g\in\gA$, which a priori introduces certain restrictions on the behavior of the function on the border of~$U$.

To avoid having to deal with this difficult problem, we use the following lemma.
\begin{lemma} 
\label{lemlocvar1} 
Let $U'$ be an open subset of $V$ containing the support of a function $f$. Then, the natural map 
(by restriction),

\snic{\hbox{from  }\Cin(V)[1/f]=\gA[1/f]\hbox{  to   }\Cin(U')[1/f\frt{U'}],}

is an \isoz. 
\end{lemma}
\begin{proof}
Recall that the support of a function $f$ is the adherence of the open subset  $U.$ We have a restriction \homo $h\mapsto h\frt{U'}$ from $\Cin(V)$ to $\Cin(U')$ that induces a \homo $\varphi:\Cin(V)[1/f]\to\Cin(U')[1/f\frt{U'}]$.
We want to prove that $\varphi$ is an \isoz.
If $g\in \Cin(U')$, then the function $gf$, which equals zero on $U'\setminus \ov U$, 
can be extended to a smooth function on the whole of $V$ by making it zero outside of $U'$.
We continue to denoted it by $gf$. 
So, the reciprocal \iso of $\varphi$ is given by $g\mapsto gf/f$
and $g/f^m\mapsto gf/f^{m+1}$.
\end{proof}

A \emph{germ of a smooth function at a point $p$} of the manifold $V$ is given by a pair $(U,f)$ where $U$ is an open subset  containing $p$ and $f$ is a smooth \hbox{function $U\rightarrow \RR$}.
Two pairs $(U_1,f_1)$ and $(U_2,f_2)$ define the same germ if there exist an open subset  $U\subseteq U_1\cap U_2$ containing $p$ such that $f_1\frt{U}=f_2\frt{U}$. The germs of smooth functions at a point $p$ form an \RRlg that we denote by~$\gA_{p}.$

We then have the following little \gui{\agq miracle.}

\begin{lemma} 
\label{lemlocvar2} 
The \alg $\gA_{p}$ is naturally isomorphic 
to the localization $\gA_{S_p}$, 
where $S_p$ is the multiplicative part of nonzero functions at a point $p$.
\end{lemma}
\begin{proof}
First, we have a natural map $\gA\to\gA_p$ that associates to a function defined 
on~$V$ its germ at $p$. 
It follows immediately that the image of $S_p$ is made of invertible
\elts of $\gA_p$.
Thus, we have a factorization of the above natural map 
which provides a \homoz~$\gA_{S_p}\to\gA_p$.

Next, we define a \homo  $\gA_p\to\gA_{S_p}$.
If $(U,f)$ defines the germ~$g$ then consider a function~$h\in \gA$ which is equal to $1$ on an open subset ~${U'}$ containing $p$ with $\overline{U'}\subseteq U$ and which equals zero outside of $U$ 
(in a chart we will be able to take $U'$ to be an open ball with center~$p$).
So, each of the three pairs $(U,f)$, $(U',f\frt{U'})$ and $(V,fh)$ define the same germ $g.$ 
Now,~$fh$ defines an \elt of $\gA_{S_p}$. It remains to check that the correspondence  that we have just established does indeed produce a \homo of the \alg $\gA_{p}$ on the \alg $\gA_{S_p}$: 
no matter how the germ is represented as a pair $(U,f)$, the \elt $fh/1$ of  $\gA_{S_p}$ only depends on the germ~$g$.  

Finally, we check that the two \homos of \RRlgs that we have defined are indeed inverse \isos of each other.
\end{proof}

In short, we have algebrized the concept of a germ of a smooth function. Except that the \mo $S_p$ is defined from the manifold $V$, not only from the \algz~$\gA$. 

However, if $V$ is compact, the \mos $S_p$ are precisely the complements of the \idemas of $\gA$. In fact, on the one hand, whether $V$ is compact or not, the set of $f\in \gA$ zero at $p$ always constitutes a \idema $\fm_p$ with a residual field equal to $\RR$. On the other hand, if $\fm$ is a \idema of $\gA$ the intersection of the 
$Z(f)=\sotq{x\in V}{f(x)=0}$ for each $f\in \fm$ is a non-empty compact subset (note that $Z(f)\cap Z(g)=Z(f^2+g^2)$). Since the \id is maximal, this compact subset is \ncrt reduced to one point~$p$ and we then get~$\fm=\fm_p$. 

\Subsubsect{Vector bundles and \mptfsz}{Vector bundles} 

Now recall the notion of a \emph{vector bundle} over $V$. 
\\
A vector bundle is given by a smooth manifold $W$, 
a smooth surjective mapping $\pi:W\rightarrow V$, and a structure of a finite dimensional vector space on every fiber~$\pi^{-1}(p)$. In addition, locally, all this must be diffeomorphic to the following simple situation, called trivial:  
$$
\pi_1: (U\times{\RR}^m) \rightarrow U,\;(p,v)\mapsto p,
$$ 
with $m$ that can depend on $U$ if $V$ is not connected. 
This means that the structure of the (finite dimensional) vector space on the fiber over $p$ must \gui{properly} depend on $p$.

Such an open set (or subset) $U$, which trivializes the bundle, is called a \emph{distinguished open set (or subset)}.

A \emph{section} of the vector bundle  $\pi:W\rightarrow V$ is by \dfn a mapping~\hbox{$\sigma:V\rightarrow W$} such that 
$\pi\circ\sigma=\Id_V$. We will denote by $\Gamma(W)$ the set of smooth sections of this bundle. It is equipped with a natural \Amo  structure. 

Now suppose that the manifold $V$ is compact.  
As the bundle is locally trivial there exists a finite covering of $V$ by distinguished open subsets $U_i$ and a 
partition of the unity 
$(f_i)_{i\in \lrbs}$ \emph{subordinate to the open cover $U_i$}: the support of $f_i$ is a compact set $K_i$ contained in~$U_i$.

We notice from Lemma~\ref{lemlocvar1} that the \algs 
$\gA[1/f_i]=\Cin(V)[1/f_i]$ and $\Cin(U_i)[1/f_i]$ are naturally isomorphic.

If we localize the \ri $\gA$ and the module $M=\Gamma(W)$ by making $f_i$ invertible, we obtain the \ri $\gA_i=\gA[1/f_i]$ and the module $M_i$. Let $W_i=\pi^{-1}(U_i)$. Then, $W_i\to U_i$ is \gui{isomorphic} to $\RR^{m_i}\times U_i\to U_i$. Thus it boils down to taking a section of the bundle $W_i$, or to taking the $m_i$ functions~$U_i\to\RR$ which make a section of the bundle $\RR^{m_i}\times U_i\to U_i$. In other words, the module of the sections of $W_i$ is free and of rank~$m$. 

Since a module that becomes free after \lon in a finite number of \eco is \ptf (\plgref{plcc.cor.pf.ptf}), we then get the direct part (point \emph{1}) of the following \thoz. 

\begin{theorem} 
\label{thFVMPTF} 
Let $V$ be a smooth compact manifold, and \hbox{let $\gA= \Cin(V)$}.
\begin{enumerate}
\item If $W\vvers\pi V$ is a vector bundle on $V$, the \Amo of the smooth sections of $W$ is a \ptf module.
\item Conversely, every \ptf \Amo is isomorphic to the module of the smooth sections of a vector bundle on $V$.
\end{enumerate}
\end{theorem}

Let us consider the converse part of the \thoz: if we take a \ptf \Amo $M$, we can construct a vector bundle $W$ over~$V$ for which the module of sections is isomorphic to $M$. We proceed as follows. 
Consider a \mprn $F=(f_{ij})\in\Mn(\gA)$ such that~\hbox{$\Im F\simeq M$}
and set 
$$
W=\sotq{(x,h)\in V \times \RR^n}{h\in\Im F\frt x},
$$
where $F\frt x$ designates the matrix $(f_{ij}(x))$.
The reader will then be able to show that~$\Im F$ is identified with the module of sections $\Gamma(W)$:
to the \elt $s\in\Im F$ is matched the section $\wi{s}$ defined by $x\mapsto \wi{s}(x)=(x,s\frt x)$.\perso{it would also be pretty good to understand the \ri of functions of the manifold $W$}
In addition, in the case where $F$ is the standard \mprn  
$$
\I_{k,n}=\blocs{.8}{.6}{.8}{.6}{$\I_k$}{$0$}{$0$}{$0_r$}\quad (k+r=n),
$$ 
 then $W$ is clearly trivial; it is equal to $V\times \left(\RR^k\times \so{0}^{r}\right)$.
Finally, a \mptf becomes free after \lon at the appropriate \eco (\thref{prop Fitt ptf 2}, point~\emph{3}, 
or \thref{propPTFDec}, more precise matrix form).
Consequently, the bundle $W$ defined above is locally trivial; it is indeed a vector bundle.

\Subsubsec{Tangent vectors and \dvns} 

A decisive example of a vector bundle is the tangent bundle, for which the \elts are the pairs $(p,v)$ where $p\in V$ and $v$ is a tangent vector at the point $p$. 

When the manifold $V$ is a manifold immersed in a space $\RR^n$, a tangent vector $v$ at the point $p$ can be identified 
with the \dvn at the point $p$ in the direction of $v$. 

When the manifold $V$ is not a manifold immersed in a space $\RR^n,$ a tangent vector $v$ can be \emph{defined} as a \emph{\dvn at the point $p$,} \cad as 
an $\RR$-\lin form 
$v:\gA\rightarrow \RR$ which satisfies Leibniz's rule
\index{derivation!at a point of a manifold}
\begin{equation}\label{Leibniz1}
v(fg)=f(p)v(g)+g(p)v(f).
\end{equation}

We can check with a few computations that the tangent vectors at $V$ indeed form a vector bundle $\mathrm{T}_V$ over~$V$.

\ss To a vector bundle $\pi:W\rightarrow V$ is associated the $\gA$-module $\Gamma (W)$ formed by the smooth sections of the bundle.
In the tangent bundle case, $\Gamma({\mathrm{T}_V})$ is nothing else but the $\gA$-module of the usual (smooth) vector fields.

Just as a tangent vector at the point $p$ is identified with a \dvn at the point $p$, which can be defined in \agq terms (\eqrf{Leibniz1}),  a (smooth) tangent vector field can be identified with an \elt of the~\emph{\hbox{$\gA$-module} of the 
\dvns of the \RRlg $\gA$}, defined as follows. 

A \dvn of an \RRlg 
$\gB$ in a \Bmo $M$ is an~$\RR$-\lin mapping $v:\gB\rightarrow M$ which satisfies Leibniz's rule%
\index{derivation!module of ---}
\begin{equation}\label{Leibniz2}
v(fg)=f\,v(g)+g\,v(f).
\end{equation}
The \Bmo of \dvns of $\gB$ in $M$ is denoted by $\Der{\RR}{\gB}{M}$.%
\index{derivation!of an \alg in a module}
\\
When we \gui{simply} refer to a \dvn of an \RRlg $gB$, 
what we mean 
is a \dvn with values in $\gB$. 
When the context is clear we write $\mathrm{Der}(\gB)$ 
as an abbreviation for~$\Der{\RR}{\gB}{\gB}$.%
\index{derivation!of an \algz}  

The \dvns at a point $p$ are then the \elts of $\Der{\RR}{\gA}{\RR_p}$
where $\RR_p=\RR$ is equipped with the  \Amo structure given by the \homo $f\mapsto f(p)$ of $\gA$ in $\RR$.
Thus $\Der{\RR}{\gA}{\RR_p}$ is an abstract \agq version of the tangent space at the point $p$ at the manifold $V$. 
\index{tangent space}\index{space!tangent ---}

\ss A smooth manifold is called \emph{parallelizable} if it has a (smooth)  field of bases
($n$ smooth sections of the tangent bundle that give a base at every point). This boils down to saying that the tangent bundle is trivial, 
or even that the~\hbox{\Amoz} of sections of this bundle, 
the module $\mathrm{Der}(\gA)$ of \dvns of~$\gA,$ is free.

\Subsubsec{Differentials and cotangent bundle} 

The dual bundle of the tangent bundle, called the cotangent bundle, has the \dil forms on the manifold~$V$ as its sections.

The corresponding \Amoz, called the \mdiz, 
can be defined \emph{by \gtrs and relations} in the following way.

\rdb\label{dilesK}
\Gnltz, if $(f_i)_{i\in I}$ is a family of \elts that generate \hbox{an \RRlgz} $\gB$, the \emph{\Bmo of (K\"ahler) \dilesz}%
\index{module!of (K\"ahler) differentials}%
\index{Kahler@K\"ahler!differential}%
\index{differential!(K\"ahler) ---}
 of~$\gB$, denoted by~$\Om{\RR}{\gB}$, is generated by the (purely formal) $\rd f_i$'s subject to the relations \gui{derived from} the relations that bind the $f_i$'s: if  $P\in \RR[z_1,\ldots ,z_n]$ and \hbox{if $P(f_{i_1},\ldots ,f_{i_n})=0$}, the derived relation is 
$$
\som_{k=1}^n \Dpp P {z_k}(f_{i_1},\ldots , f_{i_n})\rd f_{i_k} =0.
$$
Furthermore, we have the canonical mapping $\rd:\gB\to\Om{\RR}{\gB}$ available, 
defined \hbox{by $\rd f=$ the class of $f$}
(if $f=\som \alpha_if_i$, with $\alpha_i\in\RR$, $\rd f=\som \alpha_i\rd f_i$), which is a \dvnz.%
\footnote{For further details on the subject see \Thosz~\ref{thDerivUniv} and \ref{thDerivUnivPF}.}%

We then prove that, for every \RRlg $\gB$, the \Bmo of \dvns of~$\gB$ is the dual module 
of the \Bmo of K\"ahler \dilsz.

In the case where the \Bmo of \diles of $\gB$ is a \ptf module
(for example when $\gB=\gA$), then it is itself the dual module of the \Bmo of \dvns of~$\gB$.

\Subsubsect{The smooth \agq compact manifolds case}{The smooth \agq case} 

In the case of a smooth compact real \emph{\agqz} manifold $V$, 
the \alg $\gA$ of smooth functions on $V$ has as sub-\alg that of the \pol functions, denoted by~$\RR[V]$.

The modules of vector fields and \dil forms can be defined as above in terms of the \alg $\RR[V]$.

Every \mptf $M$ on $\RR[V]$ corresponds to a vector bundle~\hbox{$W\to V$} that we qualify as \emph{strongly \agqz}. The smooth sections of this vector bundle form an \Amo that is (isomorphic to) the module obtained from $M$ by scalar extension to~$\gA$.

So, the fact that 
the manifold is parallelizable can be tested on an \elr level, that of the module~$M$.  

Indeed the assertion \gui{the $\gA$-module of smooth sections of $W$ is free} concerning the smooth case is equivalent to 
the corresponding assertion in the \agq case \gui{the $\RR[V]$-module 
$M$ is free.} Proof sketch: Weierstrass' approximation \tho allows us to 
approximate a smooth section by a \pol section, and a ``smooth basis'' ($n$ smooth sections of the bundle that at every point give a basis), by a \pol one.

\ms Let us now examine the smooth compact surfaces case. 
Such a surface is parallelizable \ssi it is orientable and has an everywhere nonzero vector field. 
Figuratively the latter condition reads: the surface can be combed. 
The integral curves of the vector field then form a \emph{beautiful curve family}, \cad a locally rectifiable curve family.

Thus for an orientable smooth compact \agq surface $V$ \propeq
\begin{enumerate}
\item There exists an everywhere nonzero vector field.
\item There exists a beautiful curve family.
\item The manifold is parallelizable.
\item The K\"ahler \mdi of $\RR[V]$ is free.
\end{enumerate}

As previously explained, the latter condition stems from pure \alg (see also Section~\ref{subsecFDVL}).
 
Hence the possibility of an \gui{\agqz} \dem 
of the fact that 
the sphere cannot be combed. This has been done by Richard Swan in \cite{Swan93}. He uses some advanced tools out of the scope of our book.

Les outils qu'il utilise sont un peu trop avancés pour que nous puissions rendre compte de sa démonstration dans l'ouvrage présent.


\Subsubsect{The differential module and the module of \dvns of a \pf \algz}{Derivations of a \pf \algz} 

Let $\gR$ be a commutative \riz. 
For a  \pf \Rlg  $$\gA=\aqo{\gR[\Xn]}{\lfs}=\gR[\xn],$$ 
the \dfns of the module of \dvns and 
the module of \dils 
are updated as follows. 

We denote by $\pi:\gR[\Xn]\to\gA,\;g(\uX)\mapsto g(\ux)$ the canonical projection.

We consider the Jacobian matrix of the \sys of equations $\lfs$, 
$$
J(\uX)=
\cmatrix{ 
 \Dpp{f_1}{X_1}(\uX)  & \cdots & 
\Dpp{f_1}{X_n}(\uX)      \cr 
 \vdots   &  & 
\vdots       \cr 
 \Dpp{f_s}{X_1}(\uX)  & \cdots & 
\Dpp{f_s}{X_n}(\uX)       
}.
$$
The matrix $J(\ux)$ defines an \Ali $\gA^n\to\gA^s$.
So, we have two natural \isos $\Om{\gR}{\gA}\simeq 
\Coker\tra{J(\ux)}$
and  $\mathrm{Der}(\gA)\simeq \Ker J(\ux)$. The first \iso results from the \dfn of the \mdiz. The second can be clarified as follows: if $u=(u_1,\ldots ,u_n)\in\Ker J(\ux)$, we associate with it \gui{the partial \dvn in the direction of the tangent vector~$u$} (actually it is rather a vector field) defined by
$$
\delta_u: \gA\to\gA,\;\pi(g)\mapsto \som_{i=1}^n u_i 
\Dpp{g}{X_i}(\ux).
$$
So, $u\mapsto \delta_u$ is the \iso in question.

\begin{exercise} 
\label{exoDerDifPF} 
Prove the statement made about the module of \dvnsz. Then confirm from it that $\mathrm{Der}(\gA)$ is the dual module of $\Om{\gR}{\gA}$: if $\varphi:E\to F$ is a \ali between finite rank free modules, we always have
$\Ker \varphi \simeq (E\sta/\Im \tra{\varphi})\sta$.
\end{exercise}

In the remainder of this chapter 
we are interested in the smooth case, in which the purely \agq concepts coincide with the analogous concepts from \dil \gmtz.

\section[\Dil forms on a smooth affine manifold]{\Dil forms with \pol \coes on a smooth affine manifold} 
\label{subsecFDVL}

\vspace{4pt}
\subsect{The module of \dil forms with \pol \\
\coes on the sphere}{The sphere case} \label{subsecFDSphere}

Let $S=\sotq{(\alpha,\beta,\gamma)\in\RR^3}{
\alpha^2+\beta^2+\gamma^2=1}$. 
The \ri of \pol functions over $S$ is the \RRlg 
%
$$
\gA=\aqo{\RR[X,Y,Z]}{X^2+Y^2+Z^2-1}=\RR[x,y,z].
$$ 
The \Amo of \dil forms with \pol \coes on $S$ is
$$ \Om{\RR}{\gA}= (\gA\ {\rd}x\oplus \gA\ {\rd}y\oplus \gA\ 
{\rd}z)/\gen{x{\rd}x+y{\rd}y+z{\rd}z }\simeq \gA^3/\gA v,
$$
where $v$ is the column vector $\tra{ [\,x\;y\;z\,] }$.\\ 
This vector is \emph{unimodular} (this means that its \coos are \eco of $\gA$) since $[\,x\;y\;z\,]\cdot v=1$. Thus, the matrix 
$$
P=v\cdot[\,x\;y\;z\,]=
\cmatrix{ 
 x^2   &  xy   &   xz   \cr 
  xy  &  y^2   &  yz    \cr
 xz   &  yz   & z^2
}
$$
satisfies $P^2=P$, $P\cdot v=v$, $\Im(P)=\gA v$ such that by posing
$Q=\I_3-P$ we get  
$$
\Im(Q)\simeq \gA^3/\Im(P)\simeq \Om{\RR}{\gA}, 
\;\;\hbox{and}\;\; \Om{\RR}{\gA}\oplus \Im(P)\simeq \Om{\RR}{\gA}
\oplus \gA\simeq \gA^3.
$$
This highlights the fact that $\Om{\RR}{\gA}$ is a \stl \pro \Amo of rank $2$.

The previous considerations continue to hold if we substitute $\RR$ by a field of  \cara $\neq 2$ or even by a commutative \ri $\gR$ where $2$ is invertible.

An interesting \pb that arises is to ask 
for which \ris $\gR$, precisely, 
is the \Amo $\Om{\gR}{\gA}$ free.

\subsect{The module of \dil forms with \pol \\
\coes on a smooth \agq manifold}{The smooth \agq manifold case} 
 
\vspace{2pt}   
\subsubsec{The smooth hypersurface case} 
Let $\gR$ be a commutative \riz, and $f(\Xn)\in \RXn=\RuX$. Consider the $\gR$-\alg 
$$\gA=\aqo{\gR[\Xn]}{f}=\Rxn=\Rux .$$ 
We say that \emph{the hypersurface $S$ defined by $f=0$ is smooth} if for every field~$\gK$ \gui{extension of $\gR$}
({\footnote{In this introductory chapter, when we use the incantatory figurative expression \emph{field~$\gK$ \gui{extension of $\gR$,}} we simply mean that $\gK$ is a field with an \Rlg structure. This boils down to saying that a sub\ri of $\gK$ is isomorphic to a (integral) quotient of $\gR$, and that the \iso is given.
Consequently the \coes of $f$ can be \gui{seen} in $\gK$ and the speech following the incantatory expression does indeed have a precise \agq meaning. In Chapter~\ref{chapGenerique}
we will define a \ri extension as an \emph{injective} \homoz. This \dfn directly conflicts with the figurative expression used here if $\gR$ is not a field. This explains the inverted commas used in the current chapter.}})
 and for every point $\uxi 
=(\xin)\in \gK^n$ satisfying~\hbox{$f(\uxi )=0$} one of the \coos 
$(\partial f/\partial X_i)(\uxi )$ is nonzero.
By the formal \nstz, this is equivalent to the existence of 
$F$, $B_1$, \ldots, $B_n$ in~$\RuX$ satisfying
$$ 
Ff+B_1 \Dpp{f}{X_1}+\cdots +B_n{\partial f\over 
\partial X_n}=1.
$$
Let $b_i=B_i(\ux )$ be the image of $B_i$ in $\gA$ and
$\partial_if=(\partial f/\partial X_i)(\ux )$. We thus have in~$\gA$
$$ 
b_1 \,\partial_1f+\cdots +
b_n\, \partial_nf=1.
$$
The \Amo of \dil forms with \pol \coes on $S$ is
$$ \Om{\gR}{\gA}= \aqo{(\gA\ \rd x_1\oplus \cdots \oplus \gA\ \rd 
x_n)}{\rd f}\simeq \gA^n/\gA v,
$$
where $v$ is the column vector 
$\tra{[\,\partial_1f\;\cdots\;\partial_nf\,]}$. 
This vector is \umd since $[\,b_1\;\cdots\;b_n\,]\cdot v=1$. So, the matrix 
$$P=v\cdot[\,b_1\;\cdots\;b_n\,]=
\cmatrix{ 
 b_1\partial_1f   &  \ldots   &   b_n\partial_1f   \cr 
  \vdots  &           &  \vdots    \cr
 b_1\partial_nf   &  \ldots   &   b_n\partial_nf
}
$$
satisfies $P^2=P$, $P\cdot v=v$, $\Im(P)=\gA v$ such that by posing 
$Q=\I_n-P$ we get  
$$\Im(Q)\simeq \gA^n/\Im(P)\simeq \Om{\gR}{\gA} \;\;\hbox{and}\;\; 
\Om{\gR}{\gA}\oplus \Im(P)\simeq \Om{\gR}{\gA}\oplus \gA\simeq \gA^n.$$
This highlights the fact that $\Om{\gR}{\gA}$ is a \stl \pro \Amo of rank $n-1$.

\subsubsect{The smooth complete intersection case}{The case of a complete intersection} 

We treat the case of using two \eqns to define a smooth complete intersection. The \gnn to an arbitrary number of \eqns is straightforward.

Let $\gR$ be a commutative \riz, and $f(\uX)$, $g(\uX)\in \gR[\Xn]$. Consider the $\gR$-\alg 
$$\gA=\aqo{\gR[\Xn]}{f,g}=\gR[\xn]=\gR[\ux ].$$ 
The Jacobian matrix of the \sys of \eqns $(f,g)$ is
$$J(\uX)=
\cmatrix{ 
 {\partial f\over \partial X_1}(\uX)  & \cdots & 
{\partial f\over \partial X_n}(\uX)      \cr 
 {\partial g\over \partial X_1}(\uX)  & \cdots &
{\partial g\over \partial X_n}(\uX)      
}.
$$ 
We say that \emph{the \agq manifold $S$ defined by $f=g=0$ is smooth and of codimension $2$} if, for every field $\gK$ \gui{extension of $\gR$} and for every point~\hbox{$(\uxi) =(\xin)\in \gK^n$} 
satisfying~\hbox{$f(\uxi )=g(\uxi )=0,$} then one of the~\hbox{$2\times 2$} minors of the Jacobian matrix $J_{k,\ell}(\uxi ),$ where
$$J_{k,\ell}(\uX)=
\left\vert\matrix{ 
 {\partial f\over \partial X_k}(\uX)  &  
{\partial f\over \partial X_\ell}(\uX)      \cr 
 {\partial g\over \partial X_k}(\uX)  & 
{\partial g\over \partial X_\ell}(\uX)      
}\right\vert
$$ 
is nonzero.

\noi By the formal \nstz, this is equivalent to the existence of \pols 
$F,\,G$ and $(B_{k,\ell})_{1\leq k<\ell\leq n}$ in $\RuX$ which satisfy 
$$ Ff+Gg+\som_{ 1\leq k<\ell\leq n} {B_{k,\ell}(\uX) J_{k,\ell}(\uX) }=1.
$$
Let $b_{k,\ell}=B_{k,\ell}(\ux )$ be the image of $B_{k,\ell}$ in $\gA$ 
and
$j_{k,\ell}=J_{k,\ell}(\ux )$. We therefore have in~$\gA$
$$\preskip.0em \postskip.4em \som_{1\leq k<\ell\leq n}{b_{k,\ell}\,j_{k,\ell}}=1.\eqno(*)
$$
The \Amo of \dil forms with \pol \coes on $S$ is
$$ 
\Om{\gR}{\gA}= \aqo{(\gA\ \rd x_1\oplus \cdots \oplus \gA\ \rd x_n)}
{\rd f,\rd g}\simeq \gA^n/\Im\tra J,
$$
where $\tra J$ is the Jacobian matrix transpose (taken in~$\gA$):
$$\tra J=\tra J(\ux )=
\cmatrix{ 
\partial_1f    &  \partial_1g   \cr 
    \vdots     &  \vdots        \cr 
\partial_nf    &  \partial_ng    
}
.$$ 
Equality $(*)$ implies that the Jacobian matrix $J(\ux )$ defines a surjective \aliz, 
 and its transpose defines an injective \aliz: 
more \prmtz, if we let 
$$T_{k,l}(\ux )=
\cmatrix{ 
0 &\cdots&0& \partial_\ell g&0&\cdots &0& -\partial_k g&0&\cdots   &0 \cr 
0 &\cdots&0& -\partial_\ell f&0&\cdots &0& \partial_k f&0&\cdots   &0  
}
$$
and $T=\sum_{1\leq k<\ell\leq n}{b_{k,\ell}T_{k,l}}$, then $T\cdot \tra 
J=\I_2=J\cdot \tra {\,T}$ and the matrix $P=\tra J\cdot T$ satisfies 
$$
P^2=P,\;P\cdot \tra J=\tra J,\;\Im P=\Im\tra J\simeq \gA^2,
$$ 
so that by posing $Q=\I_n-P$ we get  
$$
\Im Q \simeq \gA^n/\Im P\simeq \Om{\gR}{\gA} \;\;\hbox{and}\;\; 
\Om{\gR}{\gA}\oplus \Im P\simeq \Om{\gR}{\gA}\oplus \gA^2\simeq 
\gA^n.
$$
This highlights the fact that $\Om{\gR}{\gA}$ is a \stl \pro \Amo of rank $n-2$.

\subsubsect{The \gnl case}{\Gnl case} 

We treat the case of using $m$ \eqns to define a smooth manifold of codimension~$r$.

Let $\gR$ be a commutative \riz, and $f_i(\uX)\in \RXn$, 
$i=1,\ldots,m$. Consider the $\gR$-\alg 
$$\gA=\aqo{\RXn}{f_1,\ldots ,f_m}=\Rxn=\gR[\ux].$$ 
The Jacobian matrix of the \sys of \eqns $(f_1,\ldots ,f_m)$ is
$$J(\uX)=
\cmatrix{ 
 {\partial f_1\over \partial X_1}(\uX)  & \cdots & 
{\partial f_1\over \partial X_n}(\uX)      \cr 
 \vdots   &  & 
\vdots       \cr 
 {\partial f_m\over \partial X_1}(\uX)  & \cdots & 
{\partial f_m\over \partial X_n}(\uX)       
}.
$$ 

\vspace{-.9em}
\pagebreak	
We say that \emph{the \agq manifold $S$ defined by $f_1=\cdots=f_m=0$ is smooth and of codimension $r$} if the Jacobian matrix taken in $\gA$ is \gui{of rank $r$,} \cad 
\Grandcadre{every minor of order $r+1$ is zero, \\and the minors of order $r$ are \comz.}

This implies that for every field $\gK$ \gui{extension of $\gR$} and at every point \hbox{$(\uxi)\in \gK^n$} of the manifold of the zeros of the $f_i$'s in $\gK^n$, the tangent space is of codimension~$r$. If the \riz~$\gA$ is reduced, this \gui{geometric} condition is in fact sufficient (in \clamaz).  

Let $J_{k_1,\ldots,k_r}^{i_1,\ldots,i_r}(\uX)$ be the $r\times r$ minor extracted from the rows $i_1,\ldots ,i_r$ and from the columns $k_1,\ldots,k_r$ of $J(\uX)$, 
and taken in $\gA$: $j_{k_1,\ldots,k_r}^{i_1,\ldots,i_r}= 
J_{k_1,\ldots,k_r}^{i_1,\ldots,i_r} (\ux )$. 

The condition on $r\times r$ minors indicates the existence of \elts $b_{k_1,\ldots,k_r}^{i_1,\ldots,i_r}$  of~$\gA$
such that
$$ \sum_{1\leq k_1<\cdots<k_r\leq n, 1\leq i_1<\cdots<i_r\leq m} 
{b_{k_1,\ldots,k_r}^{i_1,\ldots,i_r}\,j_{k_1,\ldots,k_r}^{i_1,\ldots,i_r
}}=1.
$$

The \Amo of \dil forms with \pol \coes on $S$ is
$$ 
\Om{\gR}{\gA}= \aqo{(\gA\ \rd x_1\oplus \cdots \oplus \gA\ \rd x_n)}{\rd 
f_1,\ldots ,\rd f_m}\simeq \gA^n/\Im\tra J,
$$
where $\tra J=\tra J(\ux )$ is the Jacobian matrix transpose (seen in~$\gA$).

We will see that $\Im\tra J$ is the image of a projection matrix of \hbox{rank $n-r$}.
This will highlight the fact that $\Om{\gR}{\gA}$ is a \pro \Amo of \hbox{rank 
$n-r$} (but a priori it is not \stlz).

To do so it suffices to compute a matrix $H$ of $\gA^{m\times n}$ 
such that $\tra J\,H\,\tra J=\tra J$, as then the matrix $P=\tra J\,H$ 
is the sought projection matrix.

We are therefore reduced to solve a \sli whose unknowns are the \coes of the matrix $H$. However, the solution of a \sli is essentially  a local matter, and if we localize by rendering a minor of order~$r$ invertible, the solution is not too difficult to find, knowing that every minor of order $r+1$ 
is zero.  

Here is an example of how this can work.

\begin{exercise} 
\label{exorecolle} 
{\rm  In this exercise, we perform a 
patching 
in the most naive way possible. Let $A \in \Ae{n \times m}$ be a matrix of rank~$r$.
We want to construct a matrix $B\in \Ae{m \times n}$ such that $ABA=A$.
Note that if we have a solution for a matrix $A$, we ipso facto have a solution for every \eqv matrix.
\begin{enumerate}
\item [\emph{1.}] Treat the case where $A=\I_{r,n,m}=\blocs{1}{1.6}{1}{.8}{$\I_{r}$}{$0$}{$0$}{$0$}$
\item [\emph{2.}] Treat the case where $PAQ=\I_{r,n,m}$ with $P$ and $Q$ \ivz.
\item [\emph{3.}] Treat the case where $A$ has an \iv minor of order $r$.
\item [\emph{4.}] Treat the \gnl case.
\end{enumerate}
}
\end{exercise}
%
\begin{Proof}{\bf Solution. } \emph{1.} Take $B=\tra A$.

 \emph{2.} Take $B=Q\tra{(PAQ)}P$.

\emph{3.} Suppose \spdg that the \iv minor is in the north-west corner.
Let $s=n-r$, $t=m-r$.
We write $\delta_1=\det R,$
$$  A=\blocs{1.0}{1.2}{1.0}{.8}{$R$}{$-V$}{$-U$}{$W$}\,,\;
L=\blocs{1.0}{.8}{1.0}{.8}{$\I_r$}{$0$}{$U\wi R$}{$\delta_1\I_s$}\,,
\; 
C=\blocs{1.0}{1.2}{1.0}{1.2}{$\I_r$}{$\wi R V$}{$0$}{$\delta_1\I_t$}.
$$
We get
$LA=\blocs{1.0}{1.2}{1.0}{.8}{$R$}{$-V$}{$0$}{$W'$}$ with $W'=-\delta_1 U\wi R V + W$,
then  

\snic{LAC=\blocs{1.0}{1.2}{1.0}{.8}{$R$}{$0$}{$0$}{$\delta_1 W'$}\,.}

Since the minors of order $r+1$ of $A$ are zero, we get $\delta_1^2 W'=0$
. 
 Thus let~$M=\blocs{1.0}{.8}{1.0}{1.2}{$\wi R$}{$0$}{$0$}{$0$}\,$,
hence $(LAC)M(LAC)=\blocs{1.0}{1.2}{1.0}{.8}{$\delta_1R$}{$0$}{$0$}{$0$}=\delta_1 LAC$. 
\\[1mm]
With $B_1=CML$ this gives $$LAB_1AC=(LAC)M(LAC)=\delta_1 LAC,$$ thus by multiplying on the left by $\wi L$ and on the right by $\wi C$ 
$$ \delta_1^{s+t}AB_1A=\delta_1^{s+t+1}A.
$$
Whence the solution $B=B_1/\delta_1$ since we supposed that $\delta_1$ is \ivz.

\emph{4.} The precomputation made with the minor $\delta_1$ did not require that it be \ivz. It can be done with each of the minors $\delta_\ell$ of order $r$ of $A$. 
This results in as many \egts $\delta_\ell^{s+t}AB_\ell A=\delta_\ell^{s+t+1}A$. 
\\
 A \coli $\som_\ell a_\ell\delta_\ell=1$, raised to a sufficient power, 
results in an \egt $\som_\ell b_\ell\delta_\ell^{s+t+1}=1$, hence $ABA=A$
for $B=\som_\ell b_\ell \delta_\ell^{s+t} B_\ell$.
\end{Proof}

\rems ~\\
1) We will return to the \egt $ABA = A$ when using a Cramer-style magical formula, cf. \Tho~\ref{propIGCram}.

2) In the last example, we were directly inspired by the \gui{Rank \Thoz}
which states that if a smooth mapping $\varphi:U\to \RR^k$ has constant rank~$r$
at every point of $V=\sotq{x\in U}{\varphi(x)=0}$, then $V$ is a smooth sub-manifold of codimension $r$ of the open subset $U\subseteq \RR^n$. 
It turns out that the analogue we have developed here does not always work correctly. For example with $\gR=\FF_2$, $
f_1 = X^2 + Y$ and  $f_2 = Y^2$, the manifold~$V$ is reduced to a point, the origin 
(even if we pass to the \agq closure of~$\FF_2$), in which the Jacobian matrix is of rank $1$: $\cmatrix{0&1\cr0&0}$. However, $V$ is not a curve, it is a multiple point.
This means that the Rank \Tho poses some \pbs in nonzero \caraz.
Our \dfn is therefore abusive when $\gR$ is not a \QQlgz.      
\eoe

\newpage \thispagestyle{CMcadreseul}
 
\incrementeexosetprob


\chapter{The basic local-global principle and \slisz}
\label{chapSli}
\minitoc

In this Chapter, as in the entirety of this manuscript unless explicitly stated otherwise, 
\ris are commutative and \unysz, and \homos between \ris preserve the multiplicative identities. In particular, a subring has the same multiplicative identity as the whole \riz.
\perso{compiled \today}

\subsection*{Introduction}
\addcontentsline{toc}{section}{Introduction}
Solving \slis is an omnipresent theme of commutative algebra, in particular in its most developed form for which homological methods are at use.
In this chapter, we recall some classical results on this topic, which we will come back to often throughout this work.

Particular attention is given to the basic \plgz, the notion of a \coh module and some variants of Cramer's formula.\iplg

\section[Some facts concerning \lonsz]{Some facts concerning quotients and \lonsz}
\label{secPrelimCh2}

Let us begin by recalling the following result on quotients. Let $\fa$ be an \id of a ring $\gA$. When needed, the canonical mapping will be denoted by $\pi_{\gA,\fa}:\gA\to\gA\sur{\fa}$.

The quotient ring $(\gA\sur{\fa},\pi_{\gA,\fa})$
is characterized,
\emph{up to unique \isoz},
by the following \uvl \prtz.%
\index{ring!quotient --- by the \id $\fa$}\label{NOTAgAfa}

\pagebreak
\begin{fact}
\label{factUnivQuot} \emph{(Characteristic property of the quotient by the ideal $\fa$)}\\
A \ri \homo  $\psi : \gA\to\gB$ is factorized by $\pi_{\gA,\fa}$
\ssi
$\fa\subseteq \Ker \psi$, meaning  $\psi (\fa)\subseteq \so{0_\gB}$. In this case, the factorization is
unique.
\vspace{-.1em}
\pun{\gA}{\pi_{\gA,\fa}}{\psi}{\gA\sur{\fa}}{\theta}{\gB}{\homos vanishing on $\fa$.}

\vspace{-3mm}
\noindent \emph{Explanation regarding the figure.} In a figure of the type found above,
everything but the morphism $\theta$ corresponding to the dotted arrow is given. 
The exclamation mark signifies that $\theta$ makes the diagram commute and that it is \emph{the unique}  
morphism with this \prtz.
\end{fact}
 
We denote by $M\sur{\fa M}$ the $\gA\sur{\fa}$-module obtained from the quotient of the \Amo $M$ by
the submodule generated by the elements $ax$ for $a\in\fa$ and $x\in M$.
This module can thus be defined through the extension of scalars to $\gA\sur{\fa}$ from the \Amo $M$ (see \paref{pageChgtBase}, and exercise \ref{exoEdsQuot}).

\medskip Let us move on to \lonsz, which are very analogous to quotients (we will return to this analogy in further detail on \paref{secIDEFIL}). 
In this work, when referring to a \ix{monoid} contained within a \ri (i.e. a submonoid of a \riz) we always assume a subset of the ring which contains 1 and is closed under multiplication.

\rdb
For a given \ri $\gA$, we denote by $\Ati$ the multiplicative group of  invertible \eltsz, also called the \ix{group of units}.
\label{NOTAAst}

If $S$ is a \moz, we denote by~$\gA_S$ or $S^{-1}\gA$ the localization of $\gA$
at  $S$.  Every \elt of $\gA_S$ can be written in the form $x/s$ with $x\in \gA$
and $s\in S$.
\\
By \dfn we have \hbox{$x_1/s_1=x_2/s_2$} if there exists an $s\in S$ such that
$ss_2x_1=ss_1x_2$.
When needed, we will denote by~\hbox{$j_{\gA,S}:\gA\to\gA_S$} the canonical mapping $x\mapsto x/1$. 

The localized ring $(\gA_S,j_{\gA,S})$ is characterized,
\emph{up to unique isomorphism},
by the following \uvl \prtz.%
\index{ring!localized --- at a \mo $S$}

\begin{fact}
\label{factUnivLoc} \emph{(Characteristic property of the \lon at $S$)}%
\index{localization!characteristic property of the ---}\\
A \ri \homo  $\psi : \gA\to\gB$ is factorized by $j_{\gA,S}$
\ssi
$\psi (S)\subseteq \gB\eti$. When this is the case, the factorization is
unique.
\vspace{-.5em}
\puN{\gA}{j_{\gA,S}}{\psi}{S^{-1}\gA}{\theta}{\gB}{\homos
which send $S$ into $\gB\eti$.}{8mm}
\end{fact}

\vspace{-1em}
Similarly, we denote by $M_S=S^{-1}M$ the $\gA_S$-module obtained by \lon of the \Amo $M$ at $S$.
Every \elt of $M_S$ is of the form $x/s$ with $x\in M$ and $s\in S$.
By \dfnz, we have $x_1/s_1=x_2/s_2$ if there exists an $s\in S$ such that
$ss_2x_1=ss_1x_2$.
This module $M_S$ can also be defined through an extension of scalars
to $\gA_S$ from the \Amo $M$ (see \paref{pageChgtBase}, and exercise \ref{exoEdsQuot}).
\index{module!\lon at $S$}

\rdb\label{NotaSatmon}
\smallskip The \mo $S$ contained in a \ri $\gA$ is called
\emph{saturated} when
\index{monoid!saturated ---}%
\index{saturated!monoid}
$$\forall s, t \in \gA \;\;( st\in S \;\Rightarrow\; s\in S)
$$
is satisfied. A saturated \mo is also called a \emph{filter}.
A \emph{principal filter} is a filter generated by a single \eltz; that is,
it is just the set of divisors of some arbitrary power of that \eltz.
We denote by $\sat{S}$ the saturation of the \mo $S$; 
it is obtained by adding all elements dividing an element of $S$.
When we  saturate a \moz, the \lon remains unchanged.{\footnote{In fact, depending on the specific construction chosen to define \lonz, we would either have an \egt or a canonical \iso between the two \lonsz.}}
Two \mos $S_1$ and $S_2$ are said to be \emph{\eqvsz} if they have the same saturation. 
We then write $\gA_{S_1}=\gA_{S_2}$.%
\index{filter!of a commutative ring}%
\index{monoid!equivalent ---}%
\index{equivalent!monoids}%
\index{filter!principal --- of a commutative ring}%
\index{principal!filter of a commutative ring}%

\Grandcadre{It is possible to localize by a \mo which contains $0$. \\
The result is then the \emph{trivial} \ri (recall that a \ri is trivial if it is\\ reduced to a single \eltz, \cad if $1=0$).}
\index{trivial!ring}
\index{ring!trivial ---}

\rdb
If $S$ is generated by $s\in \gA$, \cad if
$S=s^\NN\eqdefi\sotq{s^k}{k\in\NN}$, we denote by~$\gA_s$ or $\gA[1/s]$ the localized ring
$S^{-1}\gA$, which is isomorphic to
$\aqo{\gA[T]}{sT-1}$.\label{NOTAA[1/s]}

 \rdb
\medskip
In a ring, the \ixc{conductor}{of an ideal into another}
of an \id $\fa$ into an \id $\fb$
is the \id

\snic{(\fb:\fa)_\gA=\sotq{a\in\gA}{a \fa\subseteq \fb}.}

More \gnltz, if $N$ and $P$  are submodules of an \Amoz~$M$,
we define the \ixc{conductor}{of a submodule into another}
of $N$ into $P$ as the \id\index{ideal!conductor ---}\label{NOTATransp} 

\snic{(P:N)_\gA=\sotq{a\in\gA}{aN\subseteq P}.}

\rdb
Recall also that the \emph{annihilator} of an \elt $x$ from an \Amo $M$ is the \id $\Ann_\gA(x)=(\gen{0_\gA}:\gen{x})=\sotq{a\in\gA}{ax=0}$.

The \emph{annihilator of a module $M$} is the \id $\Ann_\gA(M)=(\gen{0_M}:M)_\gA$.
A module or an ideal is \ixc{faithful}{module} 
if its annihilator is reduced to $0$.%
\index{annihilator!of an element}%
\index{annihilator!of a module}%
\index{module!faithful ---}%
\label{NOTAAnn}%
\index{faithful!ideal}\index{ideal!faithful ---}

\rdb
The following notations are \egmt useful for a submodule $N$ of~$M$.

\snic{
(N:\fa)_M=\sotq{x\in M}{\,x\,
\fa\subseteq N}.\label{NOTAAnn2}
}

\snic{
(N:\fa^\infty)_M=\sotq{x\in M}{\exists n,\,x\,
\fa^n\subseteq N}.
}

The latter submodule is called the \emph{saturation} of $N$ by $\fa$ in $M$. 
\index{saturation!of $N$ by $\fa$ in $M$}

\smallskip  We say that an \elt $x$ of an \Amo $M$ is {\em  \ndzz}
(if $M=\gA$ we also say that $x$ is a \ix{nonzerodivisor})
if the sequence
$$\preskip.2em \postskip.3em 
0\vers{}\gA\vers{.x}M 
$$
is exact; in other words if $\Ann(x)=0$. If $0_\gA$ is a regular in~$\gA$, the \ri is trivial.%
\index{regular!element}

\smallskip
When the context is unambiguous, we omit the $\gA$ or $M$ subscript to simplify the previous notations regarding conductors.

\rdb
The \ixx{total}{\ri of fractions} or \ixx{total}{quotient ring} of $\gA$,
denoted by $\Frac\gA$, is the localized ring $\gA_S$, where $S$ is the \mo of \ndzs \elts of $\gA$, denoted by $\Reg \gA$.\label{NOTATotFrac}
\index{ring!total --- of fractions}
\index{ring!total quotient ---}

\goodbreak
\begin{fact}
\label{factKerAAsMMs}~
\begin{enumerate}
\item The kernel of the natural \homo $j_{\gA,s}:\gA \to \gA_s=\gA[1/s]$ is 
 the \id $(0:s^\infty)_\gA$. It is reduced to $0$ \ssi $s$ is \ndzz.

\item Similarly, the kernel of the natural \homo of $M$ to $M_s=M[1/s]$ is the \Asub $(0:s^\infty)_M$.

\item The natural \homo $\gA\to\Frac\gA$ is injective.
\end{enumerate}
\end{fact}

\begin{fact}
\label{fact.bilocal}
If $S\subseteq S'$ are two \mos of $\gA$ and $M$ is an \Amoz, 
we have two canonical identifications 
$(\gA_S)_{S'}\simeq \gA_{S'}$  and
$(M_S)_{S'}\simeq M_{S'}$.
\end{fact}

\section{The basic local-global principle}
\label{secPLGCBasic}
We will study the \gnl workings
of the \plg
in commutative \alg 
in Chapter~\ref{chapPlg}.
However, we will encounter it at every turn,
under different forms adapted to each situation. 
In this section, an essential instance of this principle is given as it is so simple and efficient that it would be a pity  
to go without it any longer.

The \plg affirms that certain \prts are true
\ssi they are true after ``sufficiently many'' \lonsz. 
In \clama we often invoke \lon at every maximal \idz. It is a lot of work and seems a bit mysterious, especially from an \algq point of view. We will use simpler (and less intimidating) versions in which only a finite number of \lons are used.

\subsec{Comaximal \lons and the local-global principle}

The following \dfn corresponds to the intuitive idea that certain (finite) systems of \lons of a \ri $\gA$ are ``sufficiently numerous'' to capture all the information contained within $\gA$.

\goodbreak
\begin{definition}\label{def.moco0} ~
\begin{enumerate}
\item   Let $s_1$, $\dots$, $s_n$ be \eltsz. if
$\gen{1} = \gen{s_1,\dots,s_n}$ then $s_1$, $\dots$, $s_n$ are said to be \ixc{comaximal}{elements}. 
%
\item   Let $S_1$, $\dots$, $S_n$ be \mosz. 
If for every $s_1\in S_1$, \dots, $s_n\in S_n$, the $s_i$'s are \com then $S_1$, $\dots$, $S_n$ are called \ixc{comaximal}{monoids}.

\end{enumerate}
\end{definition}

\rdb
\noindent \textbf{Two fundamental examples.} \label{explfonda}

\noindent 
1) If $s_1$, $\dots$, $s_n$ are \com then the \mos they generate are \comz.
Indeed, consider every $s_i^{m_i}$ ($m_i\geq 1$) in the \mos $s_i^{\NN}$ and let $a_1, \dots, a_n$ be such that $\sum_{i=1}^na_is_i=1$. By raising the latter \egt to the power of $1-n+\sum_{i=1}^nm_i$ and by conveniently regrouping the terms in the resulting sum, we get an \egt of the form $\sum_{i=1}^nb_is_i^{m_i}=1$, as required.

\noindent 
2)
If $a=a_1\cdots a_n\in\gA$, then
the \mos $a^\NN$, $1+a_1\gA,$ \ldots, $1+a_n\gA$ are \comz.
Indeed, take an \elt $b_i=1-a_ix_i$ in each \mo $1+a_i\gA$ and an \elt $a^{m}$ in the \mo $a^{\NN}$. 
We need to prove that the \id $\fm=\gen{a^{m},b_1,\dots,b_n}$ contains $1$.
However, modulo  $\fm$  we have $1=a_ix_i$, thus $1=a\prod_ix_i=ax$, and we finally obtain~$1=1^{m}=a^{m}x^{m}=0$.
\eoe

\smallskip  Here is a \carn from \clamaz.
\begin{factc}
 \label{factMoco}
 {\rm  Let $S_1$, $\dots$, $S_n$ be \mos in a nontrivial \ri $\gA$ (i.e., $1\neq_\gA0$).
The \mos $S_i$  are \com \ssi for every \idep (resp.\,for every \idemaz) $\fp$ one of the $S_i$ is contained within~$\gA\setminus\fp$.
 } 
\end{factc}
%
\begin{proof}
Let $\fp$ be a \idepz.  If none of the $S_i$'s are contained in $\gA\setminus\fp$ then for each $i$ there exists some $s_i \in S_i
\cap \fp$. Consequently, $s_1$, \ldots, $s_n$ are not \comz.
\\
 Conversely, suppose that for every \idema $\fm$ one of the $S_i$'s is contained within~$\gA\setminus\fm$ and let $s_1\in S_1$, $\ldots$, $s_n\in S_n$ then the \id
$\gen{s_1,\ldots,s_n}$ is not contained in any \idemaz. Thus it contains $1$. 
\end{proof}

\rdb
\medskip
We denote by $\Ae{m\times p}$ or $\MM_{m,p}(\gA)$
the \Amo of  $m$-by-$p$ matrices with \coes in $\gA$, and $\Mn(\gA)$ means $\MM_{n,n}(\gA)$.
The group of \iv matrices is denoted by $\GLn(\gA)$,
the subgroup consisting of 
the matrices of \deterz~$1$
is denoted by $\SLn(\gA)$.
The subset of  $\Mn(\gA)$ consisting of the \mprns
(\cad matrices $F$ such that $F^2=F$) is denoted by~$\GAn(\gA)$.
The acronyms are explained as follows: $\GL$ for \lin group,
$\SL$ for special \lin group and $\GA$ for affine Grassmannian.
\label{NOTAmatrices}

\begin{plcc}\label{plcc.basic}   
\emph{(Basic \lgb principle\iplgz,
concrete gluing
 of solutions of a \sliz)} 
\\
Let $S_1$, $\dots$, $S_n$ be \moco of $\gA$,  $B$
a matrix \hbox{of $\Ae{m\times p}$} and $C$ a column vector of $\Ae{m}$.
Then \propeq
\begin{enumerate}
\item  {The \sli $BX=C$ has a solution in $\gA^{p}$}.
\item  {For $ i\in\lrbn$,
the \sli $BX=C$ has a solution in~$\gA_{S_i}^{p}$}.
\end{enumerate}
This principle \egmt holds for \slis with \coes in an \Amo $M$.
\end{plcc}
\begin{proof}
\emph{1} $\Rightarrow$ \emph{2.} Clearly true.
\\
\emph{2} $\Rightarrow$ \emph{1.}
For each $i$, we have  $Y_i \in \gA^{p}$ and $s_i \in S_i$ such that
$B  (Y_i/s_i) =  C$ in~$\gA_{S_i}^m$. 
This means that we have some $t_i \in S_i$
such that  $t_i\,B  Y_i = s_i t_i\,C$ in~$\Ae{m}$.
Using 
$\sum_i  a_is_i t_i =1$, we get a solution 
in
$\gA$: $X=\sum_i a_i t_i Y_i$.
\end{proof}

\rem As to the merits, this \plgc boils down to the following remark
when speaking of an integral \ri (a \ri is said to be \emph{integral} if every \elt is null or \ndzz\footnote{This notion is discussed in further detail on \paref{subsecAnneauxqi}.}). If every $s_i$ is \ndzs and 
if\index{ring!integral ---}\index{integral!ring}
$$\preskip.4em \postskip.4em
\frac{x_1}{s_1}=\frac{x_2}{s_2}
=\cdots=\frac{x_n}{s_n},
$$
then the common value of these fractions, when $\sum_is_iu_i=1$,
is also equal to
$$\preskip-.2em \postskip.5em
\frac{x_1u_1+\cdots+x_nu_n}{s_1u_1+\cdots+s_nu_n}
=x_1u_1+\cdots+x_nu_n.
$$
This principle could then also be called ``the art of shrewdly getting rid of denominators.''
Arguably, the most remarkable thing is that this holds in full \gntz, even if the \ri is not integral.
Our thanks go to Claude Chevalley for introducing arbitrary \lonsz.
In some scholarly works, we find the following reformulation 
(at the cost of an information loss regarding the concreteness of the result):
\hbox{the \Amoz} $\bigoplus_\fm\!\gA_{1+\fm}$ (where $\fm$ ranges over every \idema of
$\gA$) is faithfully flat.
\eoe

\begin{corollary}
\label{corplcc.basic}
Let $S_1$, $\dots$, $S_n$ be \moco of $\gA$,  $x\in\gA$
 and $\fa,\, \fb$ be two \itfs of $\gA$.
Then, we have the following \eqvcsz.
\begin{enumerate}
\item  $x=0$ in $\gA$ \ssi for $i\in\lrbn,$ 
$x=0$ on~$\gA_{S_i}$.
\item  $x$ is \ndz in $\gA$ \ssi for $i\in \lrbn,$ 
$x$ is \ndz in~$\gA_{S_i}$.
\item
{$\fa=\gen{1}$ in $\gA$} \ssi for $i\in \lrbn,$ 
$\fa=\gen{1}$ in~$\gA_{S_i}$.
\item
{$\fa\subseteq\fb$ in $\gA$} \ssi for $i\in \lrbn,$
$\fa\subseteq\fb$ in~$\gA_{S_i}$.
\end{enumerate}
\end{corollary}

%
\facile

\rem In fact, as we will see in the \plgrf{plcc.basic.modules}, \ids do not need to be \tfz.
\eoe

\subsubsection*{Examples}

Let us give some simple examples of applications of the basic \plgcz.
A typical application of the first example (Fact \ref{factExl1Plg}) 
is where the module $M$ in the statement is a nonzero \id of a Dedekind ring.\iplg
A module $M$ is said to be \emph{\lmoz} if after each
\lon at \moco $S_1,\ldots ,S_n$, it is generated by a single \eltz.%
\index{module!\lmo ---}%
\index{locally!cyclic module}

\begin{fact}\label{factExl1Plg}
Let $M=\gen{a,b}=\gen{c,d}$ be a module with two \sgrsz. 
Suppose this module is faithful and \lmoz.
Then, there exists a matrix $A\in\SL_2(\gA)$ such that
$\vab a b \,A=\vab c d $.
\end{fact}
\begin{proof}
If \smashbot{$A=\cmatrix{x&y\cr z&t}$}, the cotransposed matrix must be equal to

\snic{B=\Adj A=\cmatrix{t&-y\cr -z&x}.}

In particular, we mean to solve the following \sliz:
$$
\vab a b \,A=\vab c d , \quad \vab c d \,B=\vab a b  \eqno (*)
$$
where the unknowns are $x$, $y$, $z$, $t$. Note that $A\,B=\det(A)\;\I_2$.\\
Conversely, if this \sli is solved,
we will have $\vab a b =\vab a b \,A\,B$. So $\big(1-\det(A)\big)\vab a b =\vab 0 0$, and since the module
is faithful, $\det(A)=1$.\\
We have some \moco $S_i$ such that $M_{S_i}$ is generated by $g_i/1$ for some $g_i\in M$.
To solve the \sli it suffices to solve it after localizing at each of the $S_i$'s. 
\\
In the \ri $\gA_{S_i}$,
we have the \egts $a=\alpha_ig_i$,  $b=\beta_ig_i$, $g_i=\mu_ia+\nu_ib$, \hbox{thus
$\big(1-(\alpha_i\mu_i+\beta_i\nu_i)\big)\,g_i=0$}. 
\\
The module $M_{S_i}=\gen{g_i}$
 stays faithful, so $1=\alpha_i\mu_i+\beta_i\nu_i$ in $\gA_{S_i}$.
Therefore:

\snic{\vab a b\, E_i=\vab{g_i}0 \;$  with  $\; E_i=\cmatrix{\mu_i&-\beta_i\cr \nu_i&\alpha_i}$  and  $\det(E_i)=1.}

Similarly we obtain $\vab c d \,C_i=\vab{g_i}0$ for some matrix $C_i$ with
\deterz~1 in $\gA_{S_i}$. By taking $A_i=E_i\,\Adj(C_i)$ 
we get $\vab a b\, A_i=\vab c d$ and $\det(A_i)=1$ in $\gA_{S_i}$. 
Thus the \sli $(*)$ has a solution in~$\gA_{S_i}$.
\end{proof}

\rdb
Our second example is given by the Gauss-Joyal Lemma:
point~\emph{\ref{i1lemPrimitf}}
 in the following lemma is proven by applying the basic \plgz.\iplg
Before stating this result, we first need to recall some \dfnsz.
\rdb

An \elt $a$ of a \ri is said to be \ix{nilpotent}
if $a^n=0$ some integer $n\in\NN$.
The nilpotent \elts of a \ri $\gA$ form an \id called
\emph{the nilradical}, or the \ixy{nilpotent}{radical} of the \riz.
A \ri is \emph{reduced} if its nilradical equals~$0$.
More \gnltz, the nilradical of an \id $\fa$ of $\gA$ is the \id consisting of \elts
$x\in\gA$, such that each $x$ has some power in $\fa$. We denote the nilradical of an \id $\fa$ of $\gA$ by $\sqrt{\fa}$ or by
$\DA(\fa)$. We also use $\DA(x)$ to denote $\DA(\gen{x})$.
An \id $\fa$ is called \emph{a radical \idz} when it is equal to its nilradical.
The \ri $\gA/\DA(0)=\gA\red$ is \emph{the reduced \ri associated with $\gA$}.%
\index{ring!reduced ---}\label{NOTADA}%
\index{reduced!ring}%
\index{nilradical!of a ring}%
\index{nilradical!of an ideal}%
\index{radical!ideal}%
\index{ideal!radical ---}

For some \pol $f$ of $\AXn=\AuX$, we call the
\emph{content} of $f$ and denote by $\rc_{\gA,\uX}(f)$ or  $\rc(f)$
 the \id generated by the \coes of $f$.
The \pol  $f$ is said to be \ixc{primitive}{polynomial} (in $\uX$) when  $\rc_{\gA,\uX}(f)=\gen{1}$.%
\index{polynomial!primitive ---}%
\index{content!of a \polz}

When a \pol $f$ of $\AX$ is given in the form
$f(X)=\sum_{k=0}^{n}a_kX^k$, we say that $n$ is the \ix{formal degree}
of $f$, and $a_n$ is its \emph{formally leading \coez}.
Finally, if $f$ is null, its formal degree is $-1$.%
\index{formally leading!\coe}

\begin{lemma}\label{lemGaussJoyal}\iJG~
\begin{enumerate}
\item \label{i1lemPrimitf} \emph{(Poor man's Gauss-Joyal)}
The product of two primitive \pols is a primitive \polz.
\item \label{i2lemPrimitf}  \emph{(Gauss-Joyal)}
For $f$, $g\in\AuX$, there exists a $p\in\NN$ such that 
$$\preskip.2em \postskip.1em 
\big(\rc(f)\rc(g)\big)^p\subseteq\rc(fg). 
$$%
\item \label{i3lemPrimitf} \emph{(Nilpotent elements in $\AuX$)}
An \elt $f$ of $\AuX$ is nilpotent \ssi all of its \coes are nilpotent.
In other words, we have the following \egtz:  $(\AuX)\red=\Ared[\uX]$.
\item \label{i4lemPrimitf} \emph{(Invertible elements in $\AuX$)}
An \elt $f$ of $\AuX$ is \iv \ssi $f(\uze)$ is \iv and $f-f(\uze)$
is nilpotent. In other words, $\AuX\eti = \Ati+\DA(0)[\uX]$
and in particular $(\Ared[\uX])\eti=(\Ared)\eti$.
\end{enumerate}
\end{lemma}

\begin{proof}
Note that, a priori, we have the following inclusion:
$\rc(fg)\subseteq \rc(f)\rc(g)$.

\noindent 
\emph{\ref{i1lemPrimitf}. For univariate \pols $f,\,g\in\AX$.} We have $\rc(f)=\rc(g)=\gen{1}$.
Consider the quotient \ri $\gB=\gA\sur{\DA\big(\rc(fg)\big)}$. We need to prove that this \ri is trivial. 
It suffices to do so after \lon at \ecoz, for example at the \coes of $f$. That is, we can suppose that some \coe of $f$ is \ivz.
Let us give a proof of a sufficiently \gnl example. Suppose

\snic{f(X)=a+bX+X^2+cX^3+\dots$ and $g(X)=g_0+g_1X+g_2X^2+\dots}

In the \ri 
$\gB$ we have $ag_0=0$, $ag_1+bg_0=0$, $ag_2+bg_1+g_0=0$, \hbox{thus  $bg_0^2=0$}, then $g_0^3=0$, thus $g_0=0$.
We then have $g=Xh$ and $c(fg)=c(fh)$. Moreover, since the formal degree of $h$
is smaller than that of $g$, we can conclude by \recu on the formal degree that $g=0$. As $\rc(g)=\gen{1}$, the \ri is trivial. 

\noindent 
\emph{\ref{i2lemPrimitf}. For univariate \polsz.} Consider a \coe $a$ of $f$ and a \coez~$b$ of $g$. We prove that $ab$ is nilpotent in 
$\gB=\gA\sur{\rc(fg)}$. This boils down to proving that $\gC=\gB[1/(ab)]$ is trivial.
However, in $\gC$, $f$ and $g$ are primitive,
so point \emph{\ref{i1lemPrimitf}}
implies that $\gC$ is trivial.

\emph{\ref{i2lemPrimitf}} and \emph{\ref{i1lemPrimitf}. General case.} 
Point \emph{\ref{i2lemPrimitf}.} is proved by \recu on the number of variables from the univariate case.
Indeed, for $f\in \AX[Y]$ we have the \egt  
$
\rc_{\gA,X,Y}(f)=\gen{\rc_{\gA,X}(h)\mid h\in\rc_{\AX,Y}(f)}. 
$
Then we deduce point \emph{\ref{i1lemPrimitf}} from it.

\emph{\ref{i3lemPrimitf}.} Note that $f^2=0$ implies
$\rc(f)^p=0$ for some $p$ from point~\emph{\ref{i2lemPrimitf}.}

\emph{\ref{i4lemPrimitf}.}
The condition is sufficient: in a \riz, if $x$ is nilpotent, then
$1-x$ is \iv because $(1-x)(1+x+\cdots+x^n)=1-x^{n+1}$. Thus if 
$u$ is \iv and $x$ nilpotent, $u+x$ is \ivz.
To see that the condition is \ncr it suffices to deal with the univariate case (we conclude by \recu on the number of  variables).
Let $fg=1$ with $f=f(0)+XF(X)$ \hbox{and  $g=g(0)+XG(X)$}. We obtain $f(0)g(0)=1$.
Let $n$ be the formal degree of $F$ and $m$ that of $G$.
We must prove that $F$ and $G$ are nilpotent.
\\
If $n=-1$ or $m=-1$, the result is obvious. We reason by \recu on $n+m$
assuming that $n$, $m\geq0$, $F_n$ and $G_m$ being the \fmt leading \coesz. By \hdr 
the result is obtained for the \ris  $(\aqo{\gA}{F_n})[X]$ and $(\aqo{\gA}{G_m})[X]$. Since $F_nG_m=0$, we can conclude with the following lemma. 

\noindent NB: some details are given in exercise~\ref{exoNilIndexInversiblePol}.
\end{proof}
%

\begin{lemma}\label{lemNilpotProd}
Let $a$, $b$, $c\in\gA$. If $c$ is nilpotent modulo $a$ and modulo~$b$, and if $ab=0$,
then $c$ is nilpotent.
\end{lemma}
\begin{proof}
We have $c^n=xa$ and $c^m=yb$ therefore $c^{n+m}=xyab=0$.
\end{proof}

\rem
We can reformulate this lemma in a more structural manner as follows. For two \ids $\fa,\, \fb$ consider the canonical morphism 
$$\preskip.3em \postskip.3em 
\gA \to \gA\sur{\fa} \times \gA\sur{\fb} 
$$ whose kernel is $\fa\cap\fb$. If an \elt of $\gA$ is nilpotent
modulo $\fa$ and modulo~$\fb$, it is also nilpotent modulo $\fa\cap\fb$, thus also
modulo $\fa\fb$, as $(\fa\cap\fb)^2 \subseteq \fa\fb$.
This touches on the ``{closed covering principle},'' see \paref{prcf1}.
\eoe


\subsec{Finite character \prtsz}

The basic \plgc can be reformulated as a ``transfer principle.''

\CMnewtheorem{ptfb}{Basic Transfer principle}{\itshape}
\begin{ptfb}\label{lem}
\label{PrTransfertBasic}\index{transfer principle}\iplg\\
For some \sli in a \ri $\gA$ the \elts $s$ such that
the \sli has a solution in $\gA[1/s]$ form an \id of $\gA$. 
\end{ptfb}

Firstly, we invite the reader to prove that this transfer principle is \eqv to the basic \plgcz.

We now provide a detailed analysis of what is going on. 
The \eqvc
actually relies on the following notion.

\begin{definition}\label{defiPropCarFini}
A \prt $\sfP $ concerning commutative \ris and modules is called a \emph{finite character property} if it is preserved by \lon
and if, when it holds for $S^{-1}\gA$, then it also holds for
$\gA[1/s]$ for some $s\in S$.\index{finite character property}
	\index{property!finite character ---}
\end{definition}
%
\begin{fact}
\label{fact1PropCarFin}
Let  $\sfP $ be a \carf \prtz.
Then the \plgc for $\sfP $ is \eqv to the transfer principle for~$\sfP$.
In other words, the following principles are \eqvsz.
\begin{enumerate}
\item  If the \prt $\sfP$ is true after \lon at every monoid in a family of \mocoz, then it is true.
\item 
The set of \elts $s$ (in a given \riz) such that the \prt $\sfP$ 
is true after \lon at $s$ is an ideal.
\end{enumerate}
\end{fact}
%
\begin{proof}
Let $\gA$ be a \ri which provides the context for the \prt $\sfP$.
Now consider the set
$I=\sotq{s\in\gA}{\sfP  \mathrm{\;is\; true\; for\;} \gA_s}$.

\noindent
\emph{1 $\Rightarrow$ 2}. Suppose \emph{1}.
Let $s,t\in I$,
$a,b\in\gA$ and $u=as+bt$.
The \eltsz~$s$ and $t$ are \com in $\gA_u$.
Since  $\sfP $ is closed under \lonz,
$\sfP $ is true for $(\gA_u)_s=(\gA_s)_u$ and $(\gA_u)_t=(\gA_t)_u$.
By applying \emph{1}, $\sfP $ is true for $\gA_u$, i.e., $u=as+bt\in I$.

\noindent
\emph{2 $\Rightarrow$ 1}. 
Suppose \emph{2} and let $(S_i)$ be the considered family of \mocoz. Since we have a \prt of \carfz, we find in each~$S_i$ an \elt $s_i$ such that $\sfP $ is true after \lon at $s_i$. Since the $S_i$'s are \com the $s_i$'s are \ecoz.
By applying \emph{2}, we get $I=\gen{1}$. Finally, the \lon at $1$ provides the answer.
\end{proof}

Most of the \plgcs which we will consider in this manuscript apply to \carf \prtsz.
One may thus replace any \plgc with its corresponding transfer principle.

\smallskip  
For \carf \prts we have an \eqvc in \clama  between two notions,
one concrete and the other abstract.

\begin{factc}
\label{fact2PropCarFin}
Let  $\sfP $ be a \carf \prtz.
Then, in \clama  \propeq
\begin{enumerate}
\item There exist \moco such that the \prt $\sfP$ is true after \lon
at each \moz.
\item The \prt $\sfP $ is true after \lon at every \idemaz.
\end{enumerate}
\end{factc}
\begin{proof}
\emph{1 $\Rightarrow$ 2}.
Let $(S_i)$ be the family of \moco under consideration. Since it is a \carf \prtz, we find in each $S_i$ some \elt $s_i$ such that $\sfP $ is true
after \lon at $s_i$. Since the $S_i$'s are \com the~$s_i$'s are \ecoz.
Let $\fm$ be a \idemaz. Some~$s_i$ is not in $\fm$. The \lon at $1+\fm$
 is a \lon of the \lon at $s_i$.
Thus $\sfP $ is true after \lon at $1+\fm$.
\\
\emph{2 $\Rightarrow$ 1}.  For each \idema $\fm$ select an
$s_\fm\notin\fm$ such that the \prt $\sfP $ is true after \lon at
$s_\fm$. The set of $s_\fm$ generates an \id which is not contained in any
\idemaz, therefore it is the \idz~$\gen{1}$. A finite family of some of these $s_\fm$ is then a \sys of \ecoz. The family of \mos generated by these \elts is suitable.
\end{proof}

We immediately obtain the following corollary.

\pagebreak
\begin{factc}
\label{factPropCarFin}
Let  $\sfP $ a \carf \prtz.
Then the \plgc for $\sfP $ is \eqv (in \clamaz) to the \plga for~$\sfP $.
In other words, the following principles are \eqvsz.
\begin{enumerate}
\item If the \prt $\sfP $ is true after \lon
at each \mo in a family of \mocoz, then it is true.
\item If the \prt $\sfP $ is true after \lon at every \idemaz, then it is true.
\end{enumerate}
\end{factc}


\rem Let us give a direct \dem of the \eqvc from \clama between the transfer principle 
and the \plga  for the \prt $\sfP $ (which we assume is of \carfz).

\noindent
\emph{Transfer $\Rightarrow$ Abstract}. Suppose the \prt is true after \lon at
every \idemaz.  The ideal given by the transfer principle cannot be strict,\footnote{An ideal $\fa$ is said to be strict when $1\notin \fa$.} otherwise it would be contained in a maximal \id $\fm$, which contradicts the fact that the \prt is true after \lon at some $s\notin\fm$.

\noindent
\emph{Abstract $\Rightarrow$ Transfer}. For each \idema $\fm$ select an
$s_\fm\notin\fm$ such that the \prt $\sfP $ is true after \lon at
$s_\fm$. The set of $s_\fm$ generates an \id not contained in any
\idemaz, thus it is the \id $\gen{1}$. We can then conclude by the transfer principle:
the \prt is true after \lon at $1$!
\eoe

\medskip \comm \label{comabstraitconcret}
The advantage of localizing at a \idep is that the result is a local \riz, which has very nice \prts (see Chapter~\ref{chap Anneaux locaux}). 
The disadvantage is that the proofs which use an \plga instead of its corresponding \plgc are non-\covs 
to the extent that 
the only access we have (in a \gnl situation) to the \ideps is given by Zorn's Lemma. 
Furthermore even 
Fact \ref{factMoco} 
is obtained by contradiction, 
 which removes any \algq trait from the corresponding ``construction.''\\
Some \plgcs do not have a corresponding abstract version as the \prt they are affiliated with is not of \carfz.
This is the case with the \plgcs  for \mtfs and for \cohs \ris (\vref{plcc.tf} and \vref{plcc.coh} respectively).

We will systematically make efficient and \cof use of the basic concrete \plg and its consequences.\iplg 
Often, we will draw inspiration from some \plgaz's \dem found in \clamaz.
In Chapter~\ref{chapPlg} we will develop a \gnl \lgb machinery to fully exploit the classical \lgb proofs in a \cov manner.
\eoe

%
\subsubsection*{Abstract version of the basic \plgz}

Since we are dealing with a \carf \prtz, \clama provides the following abstract version of the basic \plgz.\iplg

\begin{plca}\label{plca.basic} 
\emph{(Abstract basic \lgb principle: \rca of
solutions of a \sliz)}~\\
Let $B$ be a matrix $\in \Ae{m\times p}$ and
$C$ a column vector of $\Ae{m}$.
Then \propeq
\begin{enumerate}
\item  The \sli $BX=C$ has a solution in~$\gA^{p}$.
\item  For every \idema $\fm$
the \sli $BX=C$ has a solution in~$(\gA_{1+\fm})^p$.
\end{enumerate}
\end{plca}

\subsec{Forcing \comtz}

\Lon at an \elt $s\in\gA$ is a fundamental operation in commutative \alg  for forcing the invertibility of $s$. \perso{ce
paragraphe ne semble pas \^etre
au bon endroit, mais where le mettre?
en jan 07 la seule utilisation du r\'esultat sert \`a
monter que la puissance ext\'erieure d'une matrice injective
est une matrice injective.}

Sometimes you may need to make $n$ \elts $a_1,\ldots
,a_n$ of a \ri $\gA$ \comz. To this end we introduce the \ri
$$
\gB={\AXn}\big/\geN{1-\Som_ia_iX_i}=\Axn.
$$

\begin{lemma}
\label{lemKerCom}
The kernel of the natural \homo $\psi :\gA\to\gB$ is \hbox{the \id
$(0:\fa^\infty)$}, where
$\fa=\gen{a_1,\ldots ,a_n}$. In particular, the \homo is injective
\ssi $\Ann\,\fa=0$.
\end{lemma}
\begin{proof}
Let $c$ be an \elt of the kernel. Considering the \iso 

\snic{\aqo{\gB}{(x_j)_{j\neq i}}\simeq
\gA[1/a_i],}

we have $c=_{\gA[1/a_i]}0$. Thus $c\in (0:a_i^\infty)$. From this we deduce that \hbox{$c\in(0:\fa^\infty)$}. Conversely if $c\in(0:\fa^\infty)$, there exists an $r$ such that $ca_i^r=0$ for each $i$, and \hbox{therefore
$\psi(c)=\psi(c)(\sum a_ix_i)^{nr}= 0$}.
\end{proof}

\section{Coherent \ris and modules}
\label{secAnneauxCoherents}

\subsection*{A fundamental notion}
\addcontentsline{toc}{subsection}{A fundamental notion}

A \ri $\gA$ is called \ixc{coherent}{ring}\index{ring!coherent ---}
if every \lin equation

\snic{LX=0 \;\hbox{ with }\; L\in \Ae{1{\times}n}\;\hbox{ and }\;X\in \Ae{n{\times}1}}

has for solutions the
\elts of a \tf \Asub of $\Ae{n{\times}1}$. In other words,
\begin{equation}\label{eqAnCoh}
\formule{\forall n\in\NN,\, \forall L\in \Ae{1{\times}n},\, \exists
m\in\NN,
\, \exists G\in \Ae{n{\times}m},\,\forall X\in \Ae{n{\times}1} \,,\\[2mm]
\quad \quad LX=0\quad \Longleftrightarrow\quad
 \exists Y\in \Ae{m{\times}1},\; X=GY\;.}
\end{equation}
This means that we have some control over the solution space of the homogeneous \sli \hbox{$LX=0$.} 

Clearly, a finite product of \ris is \coh \ssi each factor is \cohz.

More \gnltz, given $V=(v_1,\ldots ,v_n)\in M^n$ where $M$ is an \Amoz, the \Asub of $\Ae{n}$ defined as the kernel of the \ali
$$
\breve{V}:\Ae{n}\lora M,\quad   (\xn)\lmt\Som_ix_iv_i
$$
is called the \ix{syzygy module} \emph{between the $v_i$'s}.
More specifically, we say that it is the \emph{syzygy module of (the vector) $V$}. An \elt $(\xn)$ of this kernel is called a \emph{\rdlz} or a \emph{syzygy} between the $v_i$'s. When $V$ is a \sgr of $M$ the syzygy module between the $v_i$'s is often called the \emph{(first) syzygy module of $M$}.%
\index{syzygy (\rdlz)}\index{linear dep@\rdl (syzygy)}%
\index{syzygy module!of a vector $V\in M^{n}$}%
\index{module!syzygy ---}\index{dependence relation!linear ---}%

By slight abuse of terminology, we indifferently refer to the term \emph{syzygy} to mean 
the equality $\som_ix_iv_i=0$ or the element $(\xn)\in \Ae{n}$.
The \Amo $M$ is said to be \ixc{coherent}{module}\index{module!coherent ---} if for every $V\in M^n$
the syzygy module is \tfz, in other words if we have:
\begin{equation}\label{eqMoCoh}
\formule{
\forall n\in\NN,\,\forall V\in M^{n{\times}1},\, \exists m\in\NN
,\, \exists G\in \Ae{m{\times}n},\,\forall X\in \Ae{1{\times}n}\,, \\[2mm]
\quad \quad
XV=0\quad \Longleftrightarrow\quad
\exists Y\in \Ae{1{\times}m},\; X=YG\;.}
\end{equation}
A \ri $\gA$ is then \coh \ssi it is \coh as an \Amoz.

Notice that we used a transposed notation in equation (\ref{eqMoCoh}) with respect to equation (\ref{eqAnCoh}). 
This was to avoid writing the sum $\som_ix_iv_i$ as $\som_iv_ix_i$ with $v_i\in M$ and $x_i\in\gA$.
For the remainder of this work, we will \gnlt not use this transposition,
as it seems preferable to keep to the usual form $AX=V$ for a \sliz, even when the matrices $A$
and $V$ have their \coes in a module $M$.

\begin{proposition}
\label{propCoh1}  Let $M$  be a \coh \Amoz.\\
Any {\em homogeneous \sli $BX=0$, where $B\in
M^{k{\times}n}$ and $X\in \Ae{n{\times}1}$}, has the \elts
of a \tf \Asub of~$\Ae{n{\times}1}$ as its solution set.
\end{proposition}
\begin{proof}
The \gnl \dem is by \recu on the number of linear \eqns $k$, where the procedure is as follows:
{solve the first \eqnz, then substitute the obtained \gnl solution into the second \eqnz, and so on}. 
So let us for example do the \dem for $k=2$ and take a closer look at this process. 
The matrix $B$ is composed of the rows
 $L$ and $L'$. We then have a matrix $G$ such that
\Snic{LX=0\quad \Longleftrightarrow\quad  \exists Y\in \Ae{m{\times}1},\;
X=GY.}
We now need to solve $L'GY=0$ which is equivalent to the existence of a column vector~$Z$ such that $Y=G'Z$ for a suitable matrix $G'$. Thus $BX=0$ \ssi $X$ can be expressed as $GG'Z$.
\end{proof}

The above proposition is particularly important for \slis on $\gA$ (\cad when~$M=\gA$).

\medskip \comm The notion of a \coh \ri 
is then fundamental 
from an \algq point of view in commutative algebra.
Usually, this notion is hidden behind that of a \emph{\noez} \riz,\footnote{The \cov \dfn of this notion is given after this comment.} and rarely put forward as we have here.
In \clama every \noe \ri $\gA$ is \coh because every submodule of $\Ae{n}$ is \tfz, and every \mtf is \coh for the same reason.
Furthermore, we have the Hilbert\ihi \thoz, which states that \emph{if $\gA$ is \noez, every \tf \Alg is \egmt a \noe \riz,} whereas the same statement does not hold if one replaces ``{\noez}'' with ``{\cohz}.''

From an \algq point of view however, it seems
impossible to find a satisfying \cov formulation 
of \noet which implies \cohc
(see exercise \ref{exo.quo.coh}),
and \cohc is often the most important \prt from an \algq point of view.
Consequently, \cohc cannot be implied (as  is the case in \clamaz)
when we speak of a \noe \ri or module.

The classical \tho stating that in a \noe \ri
every \tf \Amo is \noe is often advantageously replaced by the following \cof         \thoz.\footnote{For the 
non-\noee version see \thrf{propCoh2},
and for the \noee version see \cite[corollary 3.2.8 p.~83]{MRR}.}

 \emph{Over a \coh (resp.\,\noe \cohz) \ri
every \pf \Amo is  \coh (resp.\,\noe \cohz)}.

In fact, as this example shows, \noet is often an unnecessarily strong assumption.
\eoe

\medskip
The following \dfn of a \noe module is \eqv
in \clama  to the usual \dfn 
but it is much better adapted to \cov \alg 
(only the trivial \ri \cot satisfies the usual \dfnz).
\begin{definition}\label{definoetherien}\label{noetherien}\emph{(Richman-Seidenberg theory of Noetherianity, \cite{ric74,sei74b})}\\
An \Amo is called \emph{\noez} if it satisfies the following \emph{ascending chain condition}: any ascending sequence of \tf submodules has two equal consecutive terms. A \ri $\gA$ is called \emph{\noez} if it is \noe as an \Amoz.%
\index{module!Noetherian ---}%
\index{ring!Noetherian ---}%
\index{Noetherian@\noez!module}%
\index{Noetherian@\noez!ring}%
\end{definition}

\smallskip Here is a corollary of proposition \ref{propCoh1}.

\begin{corollary}\emph{(Conductors and \cohcz)}
\label{corpropCoh1} ~\\
Let $\gA$ be a \coh \riz. Then, the conductor of a \itf into another is a \itfz.
More \gnltz, if $N$ and $P$  are two \tf submodules of a \coh \Amoz, then $(P:N)$ is a \itfz.
\end{corollary}

\begin{theorem}\label{propCoh4}
An \Amo $M$ is \coh \ssi the following two conditions hold.
\begin{enumerate}
\item  The intersection of two arbitrary \tf submodules is a \mtfz.
\item  The annihilator of an arbitrary \elt is a \itfz.
\end{enumerate}
\end{theorem}
\begin{proof}
\emph{The first condition is \ncrz.}
Let $g_1$, \ldots, $g_n$ be the \gtrs of the first submodule 
and $g_{n+1}$, \ldots, $g_{m}$ be the \gtrs of the second.
Taking an \elt of the intersection reduces to finding a syzygy
$\som_{i=1}^{m}\alpha_ig_i=0$ between the $g_i$'s. To such a syzygy
 $\alpha =(\alpha_1,\ldots ,\alpha_m)\in\Ae{m}$ corresponds the \elt
$\varphi(\alpha ) = \alpha_1g_1+\cdots +\alpha_ng_n =-(\alpha
_{n+1}g_{n+1}+\cdots+\alpha_mg_{m})$ in the intersection. Thus if $S$ is a \sgr for the syzygies between the $g_i$'s,  $\varphi(S)$
generates the intersection of the two
submodules.

\noindent \emph{The second condition is \ncrz} by \dfnz.

\noindent \emph{The two conditions together are sufficient.} 
Here we give the key idea of the \dem and leave the details to the reader.
Consider the syzygy module  of some $L\in M^n$.
We perform \recu on $n$. For $n=1$ the second condition
applies and gives a \sgr for the syzygies connecting the single \elt of
$L$.
\\
Suppose that the syzygy module of every $L\in M^{n}$
is \tf and consider some $L'\in M^{n+1}$. Let $k\in \lrbn$,
we write
$L'=L_1\bullet L_2$ where $L_1=(a_1,\ldots ,a_{k})$ and
$L_2=(a_{k+1},\ldots ,a_{n+1})$. Let
$M_1=\gen{a_1,\ldots ,a_{k}}$ \hbox{and $M_2=\gen{a_{k+1},\ldots ,a_{n+1}}$}.
Taking a syzygy $\sum_{i=1}^{n+1}\alpha_ia_i=0$ reduces to taking an \elt of the intersection $M_1\cap M_2$ (as above). We thus obtain a \sgr for the syzygies between the $a_i$'s by taking the union of the three following
\syss of syzygies: the system of syzygies between the \elts of $L_1$,
the system of syzygies between the \elts of $L_2$, and that which comes from the \sgr of the intersection $M_1\cap M_2$. 
\end{proof}

In particular, \emph{a \ri is \coh \ssi on the one hand the intersection of the
two \itfs is always a \itfz, and on the other hand the annihilator of an \elt is always a \itfz}.

\medskip \exls
If $\gK$ is a discrete field, every \pf \alg over $\gK$ 
is a \coh \ri (\Thref{thpolcohfd}).
It is also clear that every Bézout domain
(cf. \paref{secBézout}) is a \coh \riz.
\eoe

\pagebreak

\subsection*{Local character of \cohcz}
\addcontentsline{toc}{subsection}{Local character of \cohcz}

\Cohc is a local notion in the following sense.

\begin{plcc}
\label{plcc.coh}
{\em  (\Coh modules)}\\
Consider a \ri $\gA$, let $S_1$, $\ldots$, $ S_n$ be \moco and $M$ 
an~\hbox{\Amoz}.
\begin{enumerate}
\item The module $M$ is \coh \ssi each $M_{S_i}$ is \cohz.
\item The \ri $\gA$  is \coh \ssi each $\gA_{S_i}$ is \cohz.
\end{enumerate} 
\end{plcc}
\begin{proof}
Let  $a=(a_1,\ldots ,a_m)\in M^m$, and $N\subseteq\Ae{m}$ be the syzygy module of~$a$.
We find that for any \mo $S$, $N_S$ is the syzygy module of $a$ in $M_S$. This brings us to prove the following \plgcz.
\end{proof}

\begin{plcc}
\label{plcc.tf}
{\em  (\Tf modules)}\\
Let $S_1$, $\ldots$, $ S_n$ be \moco of $\gA$ and $M$ an \Amoz.
Then,~$M$ is \tf \ssi each $M_{S_i}$ is \tfz.
\end{plcc}
\begin{proof}
Suppose that $M_{S_i}$ is a \tf $\gA_{S_i}$-module for each $i$.
Let us prove that $M$ is \tfz.
Let $g_{i,1}$, \ldots, $g_{i,q_i}$ be \elts of $M$ which generate $M_{S_i}$. 
Let $x\in M$ be arbitrary. For each $i$ we have some $s_i\in S_i$
and some $a_{i,j}\in \gA$ 
such that:
$$\preskip.2em \postskip.4em s_ix=a_{i,1}g_{i,1}+\cdots+a_{i,q_i}g_{i,q_i} \quad {\rm  in} \quad M.
$$
When writing $\sum_{i=1}^{n} b_i s_i =1$, we observe that $x$ is a \coli of the
$g_{i,j}$'s.
\end{proof}

\rdb\label{remplcc.tf}
\rem
Consider the $\ZZ$-submodule $M$ of $\QQ$ generated by the
\eltsz~$1/p$ where $p$ ranges over the set of prime numbers.
We can easily check that $M$ is not \tf but that it becomes \tf after \lon at any \idepz.
This means that the \plgc \ref{plcc.tf} does not have a corresponding ``abstract'' version, in which the \lon at some
\moco would be replaced by the \lon at every \idepz.
Actually, the \prt $\,\sfP\,$ for a module to be \tf is not
a \carf \prtz, as we can see with the module
$M$ above and the \mos $\ZZ\setminus\so0$ or $1+p\ZZ$.
Moreover, the \prt 
satisfies the transfer principle,
but it so happens here that it is of no use.
\eoe

\subsec{About the \egt and the membership tests}
\label{subsecTestDEgalite}
We now introduce several \covs notions
relating to the \egt test and the membership test.

\smallskip  A set $E$ is well defined when we have indicated how to construct its \elts and when we have constructed an \eqvc relation which defines the \egt of two \elts in a set.
We denote by $x=y$  the \egt in $E$, \hbox{or $x=_Ey$} if \ncrz.
The set $E$ is called
\ixc{discrete}{set} when the following axiom holds
\index{set!discrete ---}
$$
\forall x,y\in E \qquad x=y \;\; \hbox{or} \;\;  \lnot (x=y).
$$

Classically, every set is discrete,
as the ``or'' present in the \dfn is understood in an abstract manner.
Constructively, this same ``or'' is understood according to the usual language's meaning: at least one of the two alternatives must occur.
It is thus an ``or'' of an algorithmic nature.
In short, a set is discrete if we have a test for the \egt of two arbitrary \elts of this set.

If we want to be more precise and explain in detail what comprises an \egt test in the set $E$, we will say that it is a construction which,
from two given \elts of $E$, provides a ``yes'' or ``no'' answer to the posed question (are these \elts equal?). However, we could not go into much further detail.
In \coma the notions of integers and of construction are basic concepts. They can be explained and commented on, but not strictly speaking ``defined.'' The \cov meaning of the ``or'' and that of the ``there exists'' 
are as such 
directly dependent of the notion of construction,\footnote{In \clama we may wish to define the notion of construction from the notion of a ``correct program.'' However, what we define in this way is rather the notion of ``{mechanized construction},'' 
and especially 
in the notion of a ``correct program,'' there is the fact that the program must halt after a finite number of steps. This hides a ``{there exists},'' which in \coma refers in an irreducible manner to the notion of construction. On this matter, see Section \ref{AnnexeCalculsMec} of the Annex.} which we do not attempt to define.

\smallskip A \ixx{discrete}{field}
is simply a \ri where the following axiom is satisfied:
\index{field!discrete ---}\index{field}
\begin{equation}\label{eqDefCodi}
\forall x\in \gA \qquad x=0 \;\; {\rm  or} \;\;  x\in\Ati
\end{equation}
The trivial \ri is a discrete field.

\medskip
\rem
The Chinese pivot method (often called  Gaussian elimination)
works \algqz ally 
with discrete fields.
This means that the basic \lin \alg is explicit over discrete fields.\eoe

\medskip
Note that a discrete field $\gA$
is a discrete set \ssi
the test ``$1=_\gA0$?'' is explicit.\footnote{The \gnl notion of a field
in \coma will be defined \paref{corpsdeHeyting}. We will then see that if a field is a discrete set, 
then it is a discrete field.}
Sometimes, however, it is known that a ring constructed during an \algo is a discrete field without knowing whether it is trivial or not.

If $\gA$ is a nontrivial discrete field, the statement ``$M$ is a free
finite dimensional \evcz'' is more precise than the statement ``{$M$ is
a \tf \evcz}'' as in the first case knowing how to extract a basis of the
\sgr is similar to having a test of \lin independence in~$M$.

\rdb
\smallskip 
A subset $P$ of a set $E$ is said to be \ix{detachable}
when the following \prt is satisfied:\label{detachable}
$$
\preskip-.3em \postskip.4em 
\forall x\in E \qquad x\in P \;\; {\rm  or} \;\;
   \lnot(x\in P). 
$$
It amounts to the same to take a detachable part
$P$ of $E$ or to take its
\cara function $\chi_P:E\to\so{0,1}$.

In \comaz, if two sets
$E$ and $F$ are correctly defined, then so is the \ixx{set}
{of functions from $E$ to $F$}, which is denoted by $F^E$. Consequently, the \ixx{set}{of detachable subsets} of a set $E$ is itself correctly defined since it is identified with the set $\so{0,1}^E$ of \cara functions over $E$.

\subsec{Strongly discrete \cohs \ris and modules}

A \ri (resp.\,a module) is said to be \emph{strongly discrete}%
\index{module!strongly discrete ---}%
\index{strongly discrete!ring, module}%
\index{ring!strongly discrete ---}
when the \itfs (resp.\,the \tf submodules) are detachable, \cad if the
quotients by the \itfs (resp.\,by the \tf submodules) are discrete.

This means that we have a test for deciding whether a \lin \eqnz~\hbox{$LX=c$} has a solution or not, and by computing one in the affirmative case.

A key result in \cov \alg and \calf states that
$\ZZ[\Xn]$ is a \fdi \coh \riz. 
\\
More \gnltz, we have the following \cov version 
of the Hilbert \tho\ihi (see \cite{MRR,Lou}).

\emph{If $\gA$ is a \fdi \noe \coh \riz,
 so is any \pf \Algz.}

The following proposition is proven similarly to proposition
\ref{propCoh1}.

\begin{proposition}
\label{propCohfd1}
Over a \fdi \coh module $M$, every \sli  $BX=C$ ($B\in
M^{k{\times}n},\;C\in M^{k{\times}1},\;X\in \Ae{n{\times}1}$) can be tested.
In the affirmative case, a particular solution $X_0$ can be computed.
Furthermore the solutions $X$ are all the \elts of $X_0+N$ where $N$ is a \tf \Asub of $\Ae{n{\times}1}$.
\end{proposition}

\newpage
\section{Fundamental systems of \orts \idmsz} \label{sec sfio}
An \elt $e$ of a \ri is said to be \ix{idempotent} if $e^2=e$.
 In this case, $1-e$ is also an \idmz, called the \emph{complementary idempotent of $e$}, or the \emph{complement of $e$}.
\index{idempotent!complementary ---}\index{complement!of an \idmz}
For two \idms $e_1$ and $e_2$, we have
$$\gen{e_1}\cap\gen{e_2}=\gen{e_1e_2}, \quad \gen{e_1}+\gen{e_2}=\gen{e_1,e_2}=\gen{e_1+e_2-e_1e_2},
$$
where $e_1e_2$ and $e_1+e_2-e_1e_2$ are \idmsz.
Two \idms $e_1$ and $e_2$ are said to be \ixc{orthogonal}{idempotents}
when $e_1e_2=0$. We then have $\gen{e_1}+\gen{e_2}=\gen{e_1+e_2}$.

A \ri is said to be \ixc{connected}{\ri} if every \idm is equal to $0$ or $1$. \index{ring!connected ---}

In the following, we implicitly use the following obvious fact:
for an idempotent $e$ and an \elt $x$, $e$
divides  $x$ \ssi $x=ex$.

The presence of an \idm $\neq 0$, $1$ means that the \ri $\gA$ is isomorphic to a product of two \ris $\gA_1$ and $\gA_2$, and that any computation in $\gA$ can be split into two ``simpler'' computations in $\gA_1$ and $\gA_2$.
We describe the situation as follows.  

\begin{fact}\label{lemCompAnnComm}
For every \iso $\lambda:\gA\to\gA_1\times \gA_2$, there exists a unique
\elt $e\in\gA$
satisfying the following \prtsz.
\begin{enumerate}
\item The \elt $e$ is \idm  (its complement is denoted by $f=1-e$).
\item The \homo 
$\gA\to\gA_1$
identifies $\gA_1$ with $\gA\sur{\gen{e}}$  and with $\gA[1/f]$.
\item The \homo 
$\gA\to\gA_2$
identifies $\gA_2$ with $\gA\sur{\gen{f}}$ and with $\gA[1/e]$.
\end{enumerate}
Conversely, if $e$ is an \idm and $f$ is its complement, the canonical \homo
 $\gA\to\aqo{\gA}{e} \times \aqo{\gA}{f}$ is an \isoz.
\end{fact}
\begin{proof}
The \elt $e$ is defined by $\lambda(e)=(0,1)$.
\end{proof}

Here are some often useful facts.
\begin{fact}\label{fact.loc.idm}\label{fact.loc.idm1}Let $e$ be an \idm of $\gA$, $f=1-e$ and $M$ be an \Amoz.
\begin{enumerate}
\item The \mos $e^\NN=\so{1,e}$
and $1+f\gA$ have the same saturation.

\item As an \Amoz, $\gA$ is the direct sum of $\gen{e}=e\gA$ and $\gen{f}=f\gA$.
The \id $e\gA$ is a \ri 
where $e$ is a neutral \elt of the multiplication. We then have three isomorphic \ris
$$\preskip.2em \postskip.4em
\gA[1/e]=(1+f\gA)^{-1}\gA\;\simeq\;\aqo{\gA}{f}\;\simeq\; e\gA.
$$
These \isos stem from the three canonical mappings
\[\preskip.2em \postskip.4em\arraycolsep2pt
\begin{array}{lllllllllllllllll} 
\gA&\rightarrow&\gA[1/e]&\ {\string:}\ \ & x&\mapsto& x/1,  \\ 
\gA&\rightarrow&\aqo{\gA}{f}&\ {\string:}\ & x&\mapsto& x\,\mod\gen{f},  \\ 
\gA&\rightarrow& e\gA&\ {\string:}\ & x&\mapsto&   e\,x, 
\end{array}
\]  
 which are surjective and have the same kernel.

\item We have three isomorphic \Amos $M[1/e]\simeq M/fM\simeq eM.$
These \isos stem from the three canonical mappings 
\[\preskip.2em \postskip.4em \arraycolsep2pt
\begin{array}{lllllllllllllllll} 
M &\rightarrow&M [1/e]  &\ {\string:}\ \ & x&\mapsto& x/1,  \\ 
M &\rightarrow&{M }/{fM}  &\ {\string:}\ & x&\mapsto& x\,\mod\gen{f},  \\ 
M &\rightarrow& eM   &\ {\string:}\ & x&\mapsto&   e\,x, 
\end{array}
\]  
 which are surjective and have the same kernel.
\end{enumerate}
\end{fact}

In addition, care must be taken that the \id $e\gA$, which is a \ri with $e$ as its neutral \eltz, is not a sub\ri of $\gA$ (unless $e=1$).

\smallskip
In a \ri $\gA$ a \emph{\sfioz}
is a list $(e_1,\ldots ,e_n)$  of \elts of $\gA$ which satisfy the following \egtsz:%
\index{fundamental sys@\sfioz}
$$\preskip.2em \postskip.3em\ndsp 
 e_ie_j=0 \;\hbox{ for }\; i\not= j, \quad
 \hbox{ and }\quad \som_{i=1}^n\, e_i=1  .
$$
This implies that the $e_i$'s are \idmsz. We do not claim that none of them are null.\footnote{This is much nicer to obtain uniform statements. Furthermore this is virtually necessary when we do not have at our disposal an \egt to zero test for \idms in the given \riz.}

\begin{theorem}\label{fact.sfio} \emph{(Fundamental systems of \orts \idmsz)}\\
Let $(e_1,\ldots ,e_n)$ be a \sfio of a \riz~$\gA$, and  $M$ be an \Amoz.
Note that $\gA_i=\aqo{\gA}{1-e_i}\simeq\gA[1/e_i]$. Then:
\[\preskip-.20em \postskip.4em 
\begin{array}{rcl}
\gA  & \simeq  &  \gA_1\times\cdots
\times \gA_n, \\
M  & =  &  {e_1}M\oplus\cdots\oplus {e_n}M .
\end{array}
\]
\end{theorem}

\vspace{1em}
Take note that ${e_1}M$ is an \Amo and an $\gA_1$-module, but that it is not an $\gA_2$-module (unless it is null). 

The following lemma
gives a converse of \Thref{fact.sfio}.

\begin{lemma}
\label{lemfacile}
Let $(\fa_i)_{i\in\lrbn}$  be \ids of $\gA$. We have $\gA=\bigoplus_{i\in\lrbn}\fa_i$ \ssi there exists a \sfio
$(e_i)_{i\in\lrbn}$ such that $\fa_i=\gen{e_i}$ for $i\in\lrbn$.
In this case, the \sfio is uniquely determined.
\end{lemma}
\begin{proof}
Assume that $\gA=\bigoplus_{i\in\lrbn}\fa_i$. We have $e_i\in\fa_i$
such that $\som_ie_i=1$, and since $e_ie_j\in\fa_i\cap\fa_j=\{0\}$ for
$i\neq j$, we indeed obtain a \sfioz. 
Furthermore if $x\in\fa_j$, we have
$x=x\som_ie_i=xe_j$ and thus $\fa_j=\gen{e_j}$. 
The converse is \imdez.
The uniqueness follows from that of writing an \elt as a direct sum.
\end{proof}

Next we give two very useful lemmas.

\begin{lemma}
\label{lem2ide.idem} \emph{(Lemma of the \id generated by an \idmz)}\\
An \id  $\fa$ is generated by an \idm \ssi

\snic{\fa+\Ann\,\fa=\gen{1}.}
\end{lemma}
%
\begin{proof}
First, if $e$ is \idmz, we have $\Ann\,\gen{e}=\gen{1-e}$.
For the reciprocal implication, let $e\in\fa$ such that $1-e\in\Ann\,\fa$. Then $e(1-e)=0$, therefore~$e$ is \idmz, and for every $y\in\fa$, $y=ye$,
thus $\fa\subseteq \gen{e}.$
\end{proof}
%

\begin{lemma}
\label{lem.ide.idem} \emph{(Lemma of the \tf \idm \idz)}\\
If $\fa$ is a \tf \idm \id (i.e., $\fa =\fa^2 $)
in  $\gA$, then $\fa=\gen{e}$ where $e^2=e$ is entirely determined by $\fa$.%
\index{Lemma of the \tf \idm \idz}
\end{lemma}
\begin{proof} We use the \deter trick.
Consider a \sgrz~\hbox{$(a_1,\ldots a_q)$} of $\fa$ and the column vector ${\ua}=\tra{[\,a_1\;\cdots\; a_q\,]}$.
\\
Since $a_j\in\fa^2$ for $j\in\lrbq$, there exists a $C\in\MM_q(\fa)$ such 
\hbox{that ${\ua}=C\,{\ua}$}, \hbox{so
$(\I_q-C)\,{\ua}=\uze$} and $\det(\I_q-C)\,\ua=\uze$. However, $\det(\I_q-C)=1-e$
\hbox{where $e \in \fa$}. Hence $(1-e)\fa=0$, and we apply Lemma~\ref{lem2ide.idem}.\\
Finally, the uniqueness of $e$ follows immediately from Lemma~\ref{lemfacile}.
\end{proof}


Let us finally recall the {Chinese remainder \thoz}, a very efficient tool which hides a \sfioz.
Some \ids $\fb_1$, \ldots, $\fb_\ell$ of a \ri $\gA$ are called \ixc{comaximal}{\ids}
when $\fb_1+\cdots+\fb_\ell=\gen{1}$.

\CMnewtheorem{Thresteschi}{Chinese Remainder Theorem}{\itshape}
\begin{Thresteschi} ~\label{restes chinois}\\
Let $(\fa_i)_{i\in\lrbn}$ be pairwise \com \ids in $\gA$ and $\fa=\bigcap_i \fa_i$.
\\
Then $\fa=\prod_i \fa_i$, and the canonical mapping $\gA/\fa\to\prod_i
\gA/\fa_i$ is an \isoz. 
Now, there exist  $e_1$, $\dots$, $e_n$ in  $\gA$ such that
$\fa_i=\fa+\gen{1-e_i}$ and the $\pi_{\gA,\fa}(e_i)$'s form a \sfio of~$\gA/\fa$. 
\end{Thresteschi}

 As a corollary we obtain the following result.

\begin{lemma}\label{lemDesNoyaux} \emph{(Kernels' Lemma)}\\
Let $P=P_1\cdots P_\ell\in\AX$ and an \Ali $\varphi:M\to M$
 satisfying $P(\varphi)=0$. Assume the $P_i$'s are pairwise \com and let $K_i=\Ker \big(P_i(\varphi)\big)$, $Q_i=\prod_{j\neq i}P_j$.
Then we have\index{Kernels' Lemma}
$$\preskip.2em \postskip.4em\ndsp\mathrigid 2mu
K_i=\Im \big(Q_i(\varphi)\big),\;  M =\bigoplus_{j=1}^\ell K_j \hbox{ and }    \Im \big(P_i(\varphi)\big)=\Ker \big(Q_i(\varphi)\big)=\bigoplus_{j\neq i}K_i.    
$$
\end{lemma}

\smallskip 
\begin{proof}
Consider the \ri $\gB=\aqo{\AX}{P}$. The module $M$ can be seen as a \Bmo by the operation $(Q,y)\mapsto Q\cdot y=Q(\varphi)(y)$. We then apply the Chinese remainder \tho and  \Thref{fact.sfio}.
\\
This \dem summarizes the following computation.
From the \egts $U_{ij}P_i+U_{ji}P_j=1$, we get the \egts $U_iP_i+V_iQ_i=1$ together with an \egt $\sum_iW_iQ_i=1$. Let $p_i=P_i(\varphi)$, $q_i=Q_i(\varphi)$, and so on.\\
 Then, every obtained \endo commutes 
 and we obtain the \egts $p_iq_i=0$, $u_ip_i+v_iq_i=\Id_M$,  $\sum_iw_iq_i=\Id_M$.
The claimed result readily follows.
\end{proof}
%

\section{A little exterior \algz}
\label{secCramer}

\begin{flushright}
{\em That a \hmg \sys of $n$ linear equations with $n$ unknowns \\
admits (over a discrete field) a nontrivial solution\\
 \ssi the \deter of the \sys is zero, \\
here is a fact of utmost importance whose scope we \\
will never finish measuring.}\\
Anonymous \\
\end{flushright}

\begin{flushright}
{\em Eliminate, eliminate, eliminate \\
Eliminate the eliminators of elimination theory!}\\
Mathematical poem (extract)\\
S. Abhyankar \\
\end{flushright}

Some simple examples illustrating these ideas are given in this section.

\subsec{Free submodules as direct summands (Splitting Off)}

Let $k\in\NN$. A \emph{free module of rank $k$} is by \dfn an \Amo isomorphic to $\Ae{k}$. 
\index{module!free --- of rank $k$}
\index{rank!free module of --- $k$}
If $k$ is not specified, we will say \emph{free module of finite rank}. 
\index{module!free --- of finite rank}\index{rank!of a free module}

When $\gA$ is a \cdi we speak of a \emph{finite dimensional \evcz} or a \emph{finite rank \evcz} interchangeably.\index{vector space!finite dimensional ---}%
\index{finite dimensional!vector space}%
\index{dimension!of a vector space}

The modules whose structure is the simplest are the free modules of finite rank. We are thus interested in the possibility of constructing an arbitrary module $M$ in the form $L\oplus N$ where $L$ is a free module of finite rank. A (partial) answer to this question is given by the exterior algebra. 

\begin{proposition}\label{propSplittingOffAlgExt}
\emph{(Splitting Off)} \\
 Let $a_1$, \ldots, $a_k$ be \elts of an \Amo $M$, then \propeq
\begin{enumerate}
\item The submodule 
$L=\gen{a_1,\ldots ,a_k}$ of $M$
is free with basis $(a_1,\ldots ,a_k)$ and is a direct summand of $M$.
\item There exists a $k$-multi\lin alternating form $\varphi:M^k\to\gA$
which satisfies the \egt $\varphi(a_1,\ldots ,a_k)$ $=1$.
\end{enumerate}
\end{proposition}
\begin{proof}
\emph{1} $\Rightarrow$ \emph{2}. If $L\oplus N=M$, if $\pi:M\to L$ is the \prn \paralm to $N$, and if $\theta_j:L\to\gA$ is the $j$-th \coo form for the
basis $(a_1,\ldots ,a_k)$, we define
$$\preskip-.2em \postskip.4em\ndsp 
\varphi(\xk)=\det\Big(\big(\theta_j(\pi(x_i))\big)_{i,j\in \lrbk}\Big). 
$$
\emph{2} $\Rightarrow$ \emph{1}.  We define the \ali $\pi:M\to M$ as
$$\preskip.4em \postskip.4em
\pi(x)=\som_{j=1}^k\,\varphi(\underbrace{a_1,\ldots ,x,\ldots ,a_k}_{(x
\mathrm{\,\,is\,\, in\,\, position\,\,} j)})\,a_j
.
$$
We \imdt have $\pi(a_i)=a_i$ and $\Im\pi\subseteq
L:=\gen{a_1,\ldots ,a_k}$, thus $\pi^2=\pi$ and $\Im\pi=L$.
Finally, if $x=\sum_j\lambda_ja_j=0$, then $\varphi(a_1,\ldots ,x,\ldots
,a_k)=\lambda_j=0$ (with $x$ in position $j$).
\end{proof}

Special case: for $k=1$ we say that the \elt $a_1$ of $M$ is
\ixc{unimodular}{\elt of a module} when there exists a \lin form $\varphi:M\to\gA$ such that $\varphi(a_1)=1$.
The vector $b=(\bn)\in\gA^{n}$ is \umd \ssi the $b_i$'s are \comz.
In this case we also say that the sequence $(\bn)$ is 
\emph{\umdz}.\index{unimodular!sequence (or vector)}\index{sequence!unimodular ---}\index{vector!unimodular ---}

\subsec{The rank of a free module}

As we will see, the rank of a free module is a well-determined integer if the \ri is nontrivial.
In other words, two \Amosz~\hbox{$M\simeq \Ae{m}$} and $P\simeq \gA^{p}$ with $m\neq p$ can only be isomorphic if $1=_\gA0$.\label{Nota1rang}

We will use the notation $\rg_\gA(M)=k$ (or $\rg(M)=k$ if $\gA$ is clear from the context) to indicate that a (supposedly free) module has rank~$k$.

A scholarly \dem consists to say that, if $m>p$, the $m$-th exterior power 
of $P$ is $\{0\}$ whereas that of $M$ is isomorphic to $\gA$
(this is essentially the proof for Corollary~\ref{corprop inj surj det}).

The same \dem can be presented in a more \elr way as follows.
First recall the basic Cramer formula.
If~$B$ is a square matrix of order $n$, we denote by $\wi{B}$ or
$\Adj B$ the \emph{cotransposed} matrix (sometimes called \emph{adjoint}). The \elr form of  
\idcs is then expressed as:%
\label{NOTACotrans}%
\index{cotransposed!matrix}\index{adjoint!matrix}%
\index{matrix!adjoint (cotransposed) ---}%
\begin{equation}\label{eqIDC1}
A\; \Adj( A)=\Adj( A)\; A=\det( A)\;\I_n.
\end{equation}
This formula, in combination with the product formula 

\snic{\det(AB)=\det (A)\det (B),}

\snii
has a couple of implications regarding square matrices. First, that a square matrix $A$ is \iv on one side \ssi $A$ is \iv \ssi its \deter is \ivz. Second, that the inverse of $A$ is equal to $(\det A)^{-1}\Adj A$.
\\
We now consider two \Amos $M\simeq \Ae{m}$ and $P\simeq \gA^{p}$ with $m\geq
p$  and a surjective \ali $\varphi: P\rightarrow M$.
Therefore there exists a \ali $\psi :M\rightarrow P$ such that $\varphi \circ \psi =\Id_M$.
This corresponds to two matrices $A\in\Ae{m\times p}$ and $B\in\Ae{p\times m}$ with $AB=\I_m$. If $m=p$,
the matrix  $A$ is \iv with inverse $B$ and $\varphi$ and $\psi$ are reciprocal \isosz.
If~$m>p$, we have $AB=A_1B_1$ with square $A_1$ and $B_1$ respectively obtained from~$A$ and~$B$ by filling in with zeros ($m-p$ columns for $A_1$, $m-p$ rows for~$B_1$).
$$\preskip-.3em \postskip.2em
\;\;A_1 =\blocs{.5}{1.3}{0}{1.8}{}{}{$0$ \\[1mm] $~\vdots~$ \\[1mm]$0$}{$A$}\,,
\qquad
B_1 =\blocs{0}{1.8}{.5}{1.3}{}{$\begin{array}{ccc}0 & \cdots & 0\end{array}$}{}{$B$}\,,
\qquad
A_1 B_1 = {\rm I}_m.
$$

Thus $1=\det\I_m=\det (AB)=\det(A_1B_1)=\det(A_1)\det(B_1)=0.$

In this \dem we clearly see the commutativity of the \ri appear (which is truly \ncrz). Let us summarize.

\pagebreak

\begin{proposition}
\label{propDimMod1}
Let two \Amos $M\simeq \Ae{m}$ and $P\simeq \gA^{p}$ and a surjective \ali $\varphi: P\rightarrow M$.
\begin{enumerate}\itemsep0pt
\item  If $m=p$, then $\varphi$ is an \isoz. In other words, in a module~$\Ae{m}$ every \sgr of $m$ \elts is a basis.
\item  If $m>p$, then $1=_\gA0$, and if the \ri is nontrivial, $m>p$ is impossible.
\end{enumerate}
\end{proposition}

In the following, this important classification \tho will often appear as a corollary of more subtle \thosz, as for example \Thref{prop unicyc} or \Thoz~\ref{prop quot non iso}.

\subsec{Exterior powers of a module}
\label{subsecPuissExt}

\noi{\bf Terminology.} Recall that any determinant of a square matrix extracted from
$A$ on certain rows and columns is called a
\ix{minor} of  $A$. We speak of a \emph{minor of order $k$} when the extracted square matrix is in $\Mk(\gA)$.
When $A$ is a square matrix, a \emph{principal minor} is a minor corresponding to a matrix extracted on the same set of indices for both the rows and the columns.
For example \hbox{if $A\in\Mn(\gA)$}, the \coe of $X^k$ in the \pol $\det(\In+XA)$ is the sum of the principal minors of order $k$ of $A$.
Finally, a principal minor in the north-west position, i.e. obtained by extracting the matrix on the first lines and first columns, is called a
\emph{dominant principal minor}.\eoe%
\index{minor!of order $k$}%
\index{minor!principal ---}%
\index{minor!dominant principal ---}%
\index{principal!minor}

\smallskip 
Let $M$ be an \Amoz.
A $k$-multi\lin alternating map $\varphi :M^k\to P$ is called a $k$-th \ixc{exterior power}{of a module} of the \Amo $M$ if every multi\lin alternating map $\psi :M^k\to R$ is uniquely expressible in the form~\hbox{$\psi=\theta\circ\varphi$}, where $\theta$ is an \Ali from $P$ to~$R$.

\vspace{-9pt}
\Pnv{M^k}{\varphi}{\psi}{P}{\theta}{R}{}{$k$-multilinear alternating maps}{\alisz.}
\label{PuissExtMod}

\vspace{-7pt}
Clearly $\varphi :M^k\to P$ is unique in the categorical sense,
\cad that for every other exterior power $\varphi' :M^k\to P'$ there is a unique \ali $\theta:P\to P'$ which makes the suitable diagram commutative,
and that $\theta$ is an \isoz.

We then denote $P$ by $\Al{k}\!M$ or $\Vi_\Ae{k} M$  and $\varphi(\xk)$ by $\lambda_{k}(x_{1},\ldots,x_{k})$ or $x_1\vi\cdots\vi x_k$.

The existence of a $k$-th exterior power for every module $M$ results from \gnl considerations analogous to those that we will detail for the tensor product on \paref{ProdTens} in Section~\ref{secStabPf}.

The simplest theory of exterior powers, analogous to the \elr theory of the \deterz, shows that if $M$ is a free module with a basis of $n$ \elts $(a_1,\ldots ,a_n),$ then $\Al{k}M$ is zero if $k>n$, and otherwise it is a free module whose basis is the $n \choose k$ $k$-vectors $a_{i_1}\vi \cdots\vi a_{i_k}$, where $(i_1,\ldots ,i_k)$ ranges over the set of strictly increasing $k$-tuples of \elts of $\lrbn$. \\
In particular, $\Al{n}M$ is free and of rank~1 with $a_{1}\vi \cdots\vi a_{n}$ as its basis.

To every \Ali $\alpha :M\to N$ corresponds a unique \Ali $\Al{k}\alpha :
\Al{k}M\to\Al{k}N$ satisfying the \egt 

\snic{\big(\Al{k}\alpha \big)(x_1\vi\cdots\vi
x_k)=\alpha(x_1)\vi \cdots \vi \alpha(x_k)}

for every $k$-vector $x_1\vi\cdots\vi x_k$ of $\Al{k}M$. The \ali $\Al{k}\alpha$ is called the $k$-th \ixc{exterior power}{of a \aliz} of the \aliz~$\alpha.$

Moreover we have $\big(\Al{k}\alpha\big) \circ
\big(\Al{k}\beta\big) =\Al{k}(\alpha \circ \beta )$ when $\alpha
\circ\beta $ is defined.
In short, each $\Al{k}(\bullet)$ is a functor.

If $M$ and $N$ are free with respective bases $(a_1,\ldots ,a_n)$ and $(b_1,\ldots ,b_m)$, and if~$\alpha$ admits the matrix $H$ on its bases,
then  $\Al{k}\alpha$ admits the matrix denoted by~$\Al{k}H$ on the corresponding bases of $\Al{k}M$ and  $\Al{k}N$. The \coes of this matrix are all the minors of order $k$ of the matrix~$H$.

\subsec{Determinantal ideals}
\label{secIdd}

\begin{definition}\label{defIdDet}
Let $G\in\Ae{n\times m}$ and $ k\in\lrb{1..\min(m,n)}$,
 {\em  the \idd of order $k$ of the matrix $G$}  is the \idz,
denoted by $\cD_{\gA,k}(G)$ or $\cD_k(G)$,  generated by the minors of order $k$ of $G$.
\index{determinantal ideals!of a matrix}
For $k\leq 0$ we set by convention~\hbox{$\cD_k(G)=\gen{1}$}, and for
$k> \min(m,n)$, $\cD_k(G)=\gen{0}$.
\end{definition}

These conventions are natural because they allow us to obtain in full \gnt the following \egtsz.
\begin{itemize}
\item If $H=\blocs{.6}{.9}{.6}{.5}{$\I_r$}{$0$}{$0$}{$G$}$\ , for all $k\in\ZZ$
we have $\cD_k(G)=\cD_{k+r}(H)$.
\item If $H=\blocs{.6}{.9}{.7}{.5}{$0$}{$0$}{$0$}{$G$}$\ , for all $k\in\ZZ$ we have
$\cD_k(H)=\cD_{k}(G)$.
\end{itemize}

\begin{fact}
\label{fact.idd inc}
For every matrix $G$ of type $n\times m$ we have the inclusions
\begin{equation}\label{eqfact.idd inc}
\{0\}=\cD_{1+\min(m,n)}(G) \subseteq\cdots\subseteq\cD_1(G)
\subseteq\cD_0(G)=\gen{1}=\gA
\end{equation}
More \prmt for all $k,r\in\NN$ we have one inclusion
\begin{equation}\label{eqfact.idd inc2}
\cD_{k+r}(G)\subseteq\cD_{k}(G)\,\cD_{r}(G)
\end{equation}
\end{fact}

Indeed, every minor of order $h+1$ is expressed as a \coli of minors of order $h$, and the inclusion (\ref{eqfact.idd inc2}) is obtained via the Laplace expansion of the \deterz. 

\pagebreak

\begin{fact}\label{fact.idd.sousmod}\label{factlnlImage}
Let $G_1\in\Ae{n\times m_1}$,  $G_2\in\Ae{n\times m_2}$ and $H\in \gA^{p\times n}$.
\begin{enumerate}
\item If $\Im G_1\subseteq\Im G_2$, then for any integer $k$ we have $\cD_k(G_1) \subseteq \cD_k(G_2)$.
\item For any integer $k$, we have $\cD_k(HG_1) \subseteq \cD_k(G_1)$.
\item The \idds of a matrix $G\in\Ae{n\times m}$ only depend on the \eqvc class of the submodule image of $G$ (i.e., they only depend on $\Im G$, up to \auto of the module $\Ae{n}$).
\item  In particular, if $\varphi$ is a \ali between free modules of finite rank, the \idds of a matrix of $\varphi$ do not depend on the chosen bases.
We denote them by $\cD_k(\varphi)$ and we call them the
{\em  \idds of the \ali $\varphi$}.
\end{enumerate}
\index{determinantal ideals!of a \ali (free modules)}
\index{determinantal ideals}
\end{fact}
\begin{proof}
\emph{1.} Each column of $G_1$ is a \coli of columns of $G_2$. We conclude with the multilinearity of the \deterz. \\
\emph{2.} Same reasoning by replacing the columns with the rows.\\
Finally, \emph{3} implies \emph{4} and results from the two preceding items.
\end{proof}

\rem A \idd is therefore essentially attached to a \tf submodule $M$ of a free module $L$. However, it is the structure of the inclusion $M\subseteq L$ and not only the structure of $M$ which intervenes to determine the \iddsz. For example $M=3\ZZ\times 5\ZZ$ is a free \hbox{\ZZsmo} of
$L=\ZZ^2$ and its \idds are~\hbox{$\cD_1(M)=\gen{1}$}, $\cD_2(M)=\gen{15}$.
If we replace $3$ and $5$ with $6$ and $10$ for example, we obtain another free submodule, but the structure of the inclusion is different since the \idds are now~$\gen{2}$ and~$\gen{60}$.
\eoe

\begin{fact}
\label{fact.idd prod}
If $G$ and $H$ are matrices such that $GH$ is defined, then, for all $n\geq 0$ we have
\begin{equation}\label{eqfact.idd prod}\preskip-.10em \postskip.3em
\cD_n(GH)\subseteq \cD_n(G)\,\cD_n(H)
\end{equation}
\end{fact}

%
\begin{proof}
The result is clear for $n=1$. 
For $n>1$, we reduce to the \hbox{case $n=1$} by noting that the minors of order~$n$ of~$G$,~$H$ and $GH$
represent the \coes of the matrices \gui{$n$-th exterior power of $G$, $H$ and $GH$}  \big(taking into account the \egt $\Al{n}(\varphi\psi )=\Al{n}\varphi
\circ  \Al{n}\psi$\big).
\end{proof}

The following \egt is immediate.
\begin{equation}\label{eqIDDSDIR}\preskip.3em \postskip.2em
\cD_n(\varphi\oplus\psi)=\som_{k=0}^n\cD_{k}(\varphi)\,\cD_{n-k}(\psi)
\end{equation}

\vspace{-.4em}
\pagebreak	
\subsec{The rank of a matrix}
\label{secRangmat}

\begin{definition}
\label{defRangk} ~\\
A \ali $\varphi$ between free modules of finite rank is said to be 
\begin{itemize}
\item \emph{of rank $\leq k$} if $\cD_{k+1}(\varphi)=0$,
\item \emph{of rank $\geq k$} \hbox{if $\cD_{k}(\varphi)=\gen{1},$} 
\item  \emph{of rank $k$} if it is both of \hbox{rank $\geq k$} and of rank~$\leq k$.
\end{itemize}
\index{rank!of a matrix}%
\index{rank!of a linear map}%
\index{matrix!of rank $\geq k$}%
\index{matrix!of rank $\leq k$}%
\index{matrix!of rank $k$}
\end{definition}

\smallskip 
We will use the notations $\rg(\varphi)\geq k$ and $\rg(\varphi)\leq k$,
in accordance with the preceding definition, without presupposing that $\rg(\varphi)$ is defined.
Only the notation~\hbox{$\rg(\varphi)= k$} will mean that the rank is defined.

We will later generalize this \dfn to the case of \alis between \mptfsz: see the notation \ref{notaRgfi} as well as exercices~\hbox{\ref{exoIDDPTF1},
\ref{exoIDDPTF2}} and \ref{exoLocSimpPtf}.

\medskip \comm
The reader is cautioned that there is no universally accepted \dfn for \gui{matrix of rank $k$} in the literature.
When reading another book, one must first ascertain the \dfn adopted by the author. For example in the case of an integral \ri $\gA$, we often find the rank defined as that of the matrix over the quotient field of $\gA$. 
Nevertheless a matrix of rank $k$ in the sense of \Dfnz~\ref{defRangk} is \gnlt of rank $k$ in the sense of other authors.
\eoe

\medskip
The following \plgc is an immediate consequence of the basic \plgz.

\begin{plcc}
\label{plccRangMat}\emph{(Rank of a matrix)}\\
Let $S_1$, $\dots$, $S_n$ be \moco of $\gA$ and $B$ be a matrix $\in \Ae{m\times p}$.
Then \propeq
\begin{enumerate}
\item  The matrix is of rank $\leq k$ (resp.\,of rank $\geq k$) over~$\gA$.
\item For $i\in\lrbn,$ the matrix is of rank $\leq k$ (resp.\,of rank $\geq k$) over~$\gA_{S_i}$.
\end{enumerate}
\end{plcc}

\subsec{Generalized pivot method}
\label{secPivotdeGauss}

\rdb
\vspace{3pt}

\noi{\bf Terminology.} 

1) Two matrices are said to be
 \emph{\eqvesz} if we can pass from one to the other by left- and right-multiplying by invertible matrices.

\noindent 2)
Two square matrices in $\Mn(\gA)$ are said to be \emph{similar}
when they represent the same \endo of $\Ae n$ over two bases (distinct or not), in other words when they are conjugate with respect to the action $(G,M)\mapsto GMG^{-1}$ of $\GLn(\gA)$ over 
$\Mn(\gA)$.\index{similar!matrices}\index{matrices!similar ---}

\noindent 3)
An \emph{\elr row operation} on a matrix of $n$ rows consists in replacing a row $L_i$ with a row $L_i+\lambda L_j$ where~\hbox{$ i\neq j$}.\\
 We also denote this by $L_i\aff L_i+\lambda L_j$. This corresponds to the left-multiplication by a matrix, said to be
 \emph{\elrz}, denoted by $\rE^{(n)}_{i,j}(\lambda)$
(or, if the context allows it, $\rE_{i,j}(\lambda)$). This matrix is obtained from $\In$ by means of the same \elr row operation.
\\
The right-multiplication by the same matrix $\rE_{i,j}(\lambda)$ corresponds to the \emph{\elr column operation} (for a matrix having $n$ columns) which transforms the matrix $\In$ into $\rE_{i,j}(\lambda)$: $C_j\aff C_j+\lambda C_i$.

\noindent 4)
The subgroup of $\SLn(\gA)$
generated by the \elr matrices is called the \emph{\elr group} and it is denoted by $\En(\gA)$. Two matrices are said to be
 \emph{\elrt \eqvesz} when we can pass from one to the other via \elr row and column operations.%
\index{matrices!equivalent ---}%
\index{equivalent!matrices}%
\index{elementarily equivalent!matrices}%
\index{matrices!elementarily equivalent ---}%
\index{matrix!elementary ---}%
\index{group!elementary ---}%
\index{operation!elementary ---}%
\index{elementary!matrix}%
\index{elementary!group}%
\index{elementary!row operation}%
\index{elementary!column operation}%
\label{NOTAEn}%
\eoe

\CMnewtheorem{lemmininv}{Invertible minor lemma}{\itshape}
\begin{lemmininv}\label{lem.min.inv}\index{Invertible minor lemma}
\emph{(Generalized pivot)}\\
If a matrix $G\in\gA^{q\times m}$ has an invertible minor of order $k\leq \min(m,q)$
, it is \eqve to a matrix
$$\preskip.0em \postskip.4em 
\cmatrix{
    \I_{k}   &0_{k,m-k}      \cr
    0_{q-k,k}&  G_1      }, 
$$
where $\cD_r(G_1)=\cD_{k+r}(G)$ for all $r\in\ZZ$. 
\end{lemmininv}

\begin{proof}
By eventually permuting the rows and the columns we bring the invertible minor to the top left. Next, by right-multiplying (or left-multiplying) by an invertible matrix, we reduce to the form
$$\preskip.4em \postskip.4em 
G'\;=\;
\cmatrix{
\I_k & A     \cr
  B &  C
}, 
$$
then by \elrs row and column operations, we obtain
$$\preskip.4em \postskip.4em
G''\;=\;
\cmatrix{
   \I_{k}   &0_{k,m-k}      \cr
    0_{q-k,k}&     G_1}.
$$
Finally, $\cD_r(G_1)=\cD_{k+r}(G'')=\cD_{k+r}(G)$ for all $r\in\ZZ$.
\end{proof}

As an immediate consequence we obtain the freeness lemma.

\CMnewtheorem{lemli}{Freeness lemma}{\itshape}
\begin{lemli}\label{NOTAIkqm} \label{lem pf libre}\index{Freeness lemma}
Consider a matrix $G\in\gA^{q\times m}$ of {rank $\leq k$ with $1\leq k\leq \min(m,q)$}. If the matrix~$G$ has an invertible minor of order~$k$, then it is \eqve to the matrix
$$\preskip.0em \postskip.4em 
\I_{k,q,m}\;=\;
\cmatrix{
    \I_{k}   &0_{k,m-k}      \cr
    0_{q-k,k}&     0_{q-k,m-k}      }. 
$$
In this case, the image, the kernel and the cokernel of $G$ are free,
respectively of ranks $k$, $m-k$ and $q-k$. Moreover the image and the kernel have free summands.\\ 
If $i_1$, $\ldots$, $i_k$ (resp.\,$j_1$, $\ldots$, $j_k$) are the indexes of rows (resp.\,of columns) of the invertible minor, then the columns $j_1$, $\ldots$, $j_k$ form a basis of the module $\Im G$, and $\Ker G$ is the module of vectors annihilated by the \lin forms corresponding to the rows $i_1$, $\ldots$, $i_k$.
\end{lemli}

\begin{proof}
With the notations of the previous lemma we have {\mathrigid 1mu $\cD_1(G_1)=\cD_{k+1}(G)= 0 $},
\hbox{so $G_1=0$}. The rest is left to the reader.
\end{proof}

The matrix $\I_{k,q,m}$ is called a \emph{standard simple matrix}.
We denote the matrix  $\I_{k,n,n}$ by $\I_{k,n}$ and we call it a \emph{standard \mprnz}.%
\index{matrix!standard simple ---}%
\index{matrix!standard --- of projection}

\begin{definition}%
\label{defnl}%
A \ali between free modules of finite rank 
is said to be
 \emph{simple} if it can be represented by a matrix $\I_{k,q,m}$ over suitable bases. Similarly a matrix is said to be
 \emph{\nlz} when it is \eqve to a matrix $\I_{k,q,m}$.  
\index{simple!matrix}%
\index{simple!\ali}%
\index{matrix!simple ---}%
\end{definition}

\subsec{Generalized Cramer formula}

We study in this subsection some \gnns of the usual Cramer formulas. We will exploit these in the following paragraphs.

For a matrix $A\in\Ae{m{\times}n}$ we denote by $A_{\alpha,\beta}$ the matrix extracted on the rows $\alpha=\{\alpha_1,\ldots ,\alpha_r\}\subseteq\lrbm$ and the columns $\beta=\{\beta_1,\ldots ,\beta_s\}\subseteq\lrbn$.%
\label{NOTAextraite}

Suppose that the matrix $A$ is of rank $\leq k$.
Let $V\in\Ae{m{\times}1}$ be a column vector such that the bordered matrix $[\,A\,|\,V\,]$ is also of rank $\leq k$.
Let us call $A_j$ the $j$-th column of $A$.
Let $\mu_{\alpha,\beta}=\det(A_{\alpha,\beta})$ be the minor of order $k$ of the matrix $A$ extracted on the rows $\alpha=\{\alpha_1,\ldots ,\alpha_k\}$ and the columns $\beta=\{\beta_1,\ldots ,\beta _k\}$.
For $j\in\lrbk$ let $\nu_{\alpha,\beta,j}$ be the \deter of the same extracted matrix, except that the column $j$ has been replaced with the extracted column of $V$ on the rows $\alpha$.
Then, we obtain for each pair $(\alpha,\beta)$ of multi-indices a Cramer \idtz:
\begin{equation}\label{eqMPC1}
\qquad \mu_{\alpha,\beta}\;V=\som_{j=1}^k
\nu_{\alpha,\beta,j}\,A_{\beta_j}\qquad
\end{equation}
due to the fact that the rank of the bordered matrix $[\,A_{1..m,\beta}\,|\,V\,]$ is $\leq k$. This can be read as follows:
\begin{eqnarray}\preskip-.2em \postskip.4em
\qquad \mu_{\alpha,\beta}\;V&=& \left[
\begin{array}{ccc}
    A_{\beta_1} & \ldots  & A_{\beta_k}
\end{array}
\right] \cdot  \left[
\begin{array}{c}
      \nu_{\alpha,\beta,1} \\
    \vdots  \\
     \nu_{\alpha,\beta,k}
\end{array}
\right]  \nonumber \\[.1em]
 &=&\left[
\begin{array}{ccc}
    A_{\beta_1} & \ldots  & A_{\beta_k}
\end{array}
\right] \cdot
\Adj(A_{\alpha,\beta})
 \cdot
 \left[
 \begin{array}{c}
     v_{\alpha_1}  \\
     \vdots   \\
     v_{\alpha_k}
 \end{array} \right] \nonumber\\[.2em]
\label{eqMPC2}
 &=&A \cdot (\I_n)_{1.. n,\beta} \cdot
\Adj(A_{\alpha,\beta}) \cdot (\I_m)_{\alpha,1.. m}\cdot V
\end{eqnarray}

This leads us to introduce the following notation.
\begin{notation}
\label{notaAdjalbe}
{\rm We denote by $\cP_{\ell}$ the set of parts of 
$\lrb{1..\ell}$ and $\cP_{k,\ell}$ the set of parts  of
$\lrb{1..\ell}$ with $k$ \eltsz. 
For $A\in\Ae{m{\times}n}$ and $\alpha\in
\cP_{k,m},\,\beta\in \cP_{k,n}$ we define
$$
\preskip-.20em \postskip.0em
\Adj_{\alpha,\beta}(A):=
(\I_n)_{1.. n,\beta} \cdot
\Adj(A_{\alpha,\beta}) \cdot (\I_m)_{\alpha,1.. m}.
$$
}
\end{notation}

\vspace{-.2em}
For example with the matrix
$$\preskip.2em \postskip.4em
\;\;A=\crmatrix{ 5&-5&7&4\cr
9&-1&2&7\cr13&3&-3&10},$$
and the parts $\alpha=\{1,2\}$ and
$\beta=\{2,3\}$, we obtain
{\small
$$\preskip-.2em \postskip-.2em
A_{\alpha,\beta}=\crmatrix{
-5&7\cr-1&2},\;\Adj(A_{\alpha,\beta})=\crmatrix{
2&-7\cr1&-5}\; \hbox{and}\;
\Adj_{\alpha,\beta}(A)=\crmatrix{ 0&0&0\cr2&-7&0\cr
   1&-5&0\cr0&0&0}.
$$
}

When $\cD_{k+1}([\,A\,|\,V\,])=0$, \egrf{eqMPC2} is written as follows.
\begin{equation}\label{eqGema}
\mu_{\alpha,\beta}\;V\;=\;A \cdot \Adj_{\alpha,\beta}(A) \cdot V
\end{equation}

We thus obtain the following \egtz, 
under the assumption that $A$ is of rank $\leq k$.
\begin{equation}\preskip-.20em \postskip.4em\label{eqIGCram}
\mu_{\alpha,\beta}\;A\;=\;A \cdot \Adj_{\alpha,\beta}(A) \cdot A
\end{equation}

The \idcs (\ref{eqGema}) and (\ref{eqIGCram}) provide the congruences which are not subject to any hypothesis: it suffices for example to read (\ref{eqGema}) in the quotient \ri $\gA/\cD_{k+1}([\,A\,|\,V\,])$ to obtain the congruence~(\ref{eqGema2}).

\begin{lemma}
\label{lemCram}  \emph{(Generalized Cramer formula)} 
Without any assumption on the matrix $A$ or the vector $V$, we have {for
$\alpha\in \cP_{k,m}$} \hbox{and $\beta\in \cP_{k,n}$} the following congruences.
\begin{eqnarray}\label{eqGema2}\preskip-.20em \postskip.2em
\mu_{\alpha,\beta}\;V&\equiv&A \cdot \Adj_{\alpha,\beta}(A) \cdot V
\qquad \mathrm{mod}\quad \cD_{k+1}(\,[\,A\,|\,V\,]\,),
\\[.2em] \label{eqCGCram}
\mu_{\alpha,\beta}\;A&\equiv&A \cdot \Adj_{\alpha,\beta}(A) \cdot A
\qquad \mathrm{mod}\quad \cD_{k+1}(A).
\end{eqnarray}
\end{lemma}

A simple special case is the following where $k=m\leq n$.
\begin{equation}\label{eqGema3}\preskip.3em \postskip.4em
\mu_{1..m,\beta}\;\I_m\;=\;A \cdot \Adj_{1..m,\beta}(A)\quad (\beta\in
\cP_{m,n}).
\end{equation}
This \egt is in fact a direct consequence of the basic \idcz~(\ref{eqIDC1}). Similarly we obtain
\begin{equation}\label{eqGema4}\preskip.2em \postskip.0em
\mu_{\alpha,1..n}\;\I_n\;=\; \Adj_{\alpha,1..n}(A)  \cdot A
\quad (\alpha \in \cP_{n,m},\,n\leq m)
\end{equation}

\subsec{A magic formula}

An immediate consequence of the \idc (\ref{eqIGCram}) is the less usual \idt (\ref{eqIGCram2}) given in the following \thoz. Similarly the \egts (\ref{eqIGCram3}) and (\ref{eqIGCram4}) easily result from (\ref{eqGema3}) and (\ref{eqGema4}).

\begin{theorem}
\label{propIGCram}
Let $A\in\Ae{m\times n}$ be a matrix of rank $k$.
We thus have an \egt
$\som_{\alpha\in \cP_{k,m},\beta\in \cP_{k,n}}
c_{\alpha,\beta}\,\mu_{\alpha,\beta}=1.$
Let 
$$\preskip.3em \postskip.0em \ndsp
B\;=\;\som_{\alpha\in \cP_{k,m},\beta\in
\cP_{k,n}}\,c_{\alpha,\beta}\,\Adj_{\alpha,\beta}(A). 
$$
\begin{enumerate}
\item  We have
\begin{equation}\label{eqIGCram2}\preskip-.2em \postskip.4em
A\cdot B\cdot A=A.
\end{equation}
Consequently $A\, B$ is a \prn matrix of rank $k$ and the submodule 
$\Im A=\Im AB$
is a direct summand in~$\Ae{m}$.

\item  If  $k=m$, then
\begin{equation}\label{eqIGCram3}\preskip-.20em \postskip.4em
A\cdot B=\I_m.
\end{equation}

\item  If $k=n$, then
\begin{equation}\label{eqIGCram4}\preskip-.2em \postskip.4em
B\cdot A=\I_n.
\end{equation}
\end{enumerate}
\end{theorem}

The following \idtz, which we will not use in this work, is even more miraculous.

\begin{proposition}
\label{propIGCram2} \emph{(Prasad and Robinson)}\\
With the assumptions and the notations of \Thref{propIGCram}, 
if we have 

\snic{\forall \alpha,\alpha'\in \cP_{k,m},$
$\forall\beta,\beta '\in \cP_{k,n}$
$\;c_{\alpha,\beta}\,c_{\alpha',\beta'}=
c_{\alpha,\beta'}\,c_{\alpha',\beta},}

then
\begin{equation}\label{eqIGCramPraRo}\preskip-.1em \postskip.0em
B\cdot A\cdot B=B.
\end{equation}
\end{proposition}

\subsec{Generalized inverses and \lnls maps}
 \label{subsecInvGen}

Let $E$ and $F$ be two \Amosz, and $\varphi:E\rightarrow F$ be a \aliz.
We can see this 
as some sort of \gne \sli (a usual \sli corresponds to the free modules of finite rank case). Informally such a \sli is considered to be \gui{well-conditioned} if there is a systematic way to solve the equation $\varphi(x)=y$ for $x$ from a given $y$, when such a solution exists. 
More \prmtz, we ask if there exists a \aliz~\hbox{$\psi:F\rightarrow E$} satisfying  $\varphi\big(\psi(y)\big)=y$ each time there exists a solution~$x$. This amounts to asking
$\varphi\big(\psi\big(\varphi(x)\big)\big)=\varphi(x)$ for all $x\in E$.

This clarifies the importance of the \eqrf{eqIGCram2} and leads to the notion of a \ingz.

The terminology regarding \ings does not seem fully fixed. We adopt that of \cite{Lan}. \\
In the book \cite{Bha}, the author uses the term
\gui{reflexive g-inverse.}

\pagebreak

\begin{definition}
\label{defIng}
Let $E$ and $F$ be two \Amosz, and $\varphi:E\rightarrow F$ be a \aliz.
A \ali $\psi :F\rightarrow E$ is called a \ix{generalized inverse} of $\varphi$ if we have
\begin{equation}\label{eqdefIng}\preskip-.3em \postskip.4em
\varphi \circ\psi \circ\varphi =\varphi \quad \mathrm{and}  \quad \psi
\circ\varphi \circ\psi =\psi.
\end{equation}
A \ali is said to be
 \emph{\lnlz} when it has a \ingz. \index{{locally!simple linear map}}
\end{definition}

The following fact is immediate.
\begin{fact}
\label{factIng0}
When $\psi$ is a \ing of $\varphi$, we have:
\begin{enumerate}
\item [--] $\varphi\, \psi$ and $\psi\, \varphi$ are \prnsz,
\item [--] $\Im\varphi =\Im\varphi\, \psi$,  $\Im\psi  =\Im\psi\,
\varphi$,  $\Ker\varphi =\Ker\psi \,\varphi$, $\Ker\psi
=\Ker\varphi\,\psi$,
\item [--] $E=\Ker\varphi \oplus \Im\psi$ and
 $F=\Ker\psi  \oplus \Im\varphi$,
\item [--] $\Ker\varphi \simeq \Coker\psi$ and
 $\Ker\psi \simeq \Coker\varphi$.
\end{enumerate}
Moreover $\varphi$ and $\psi$ provide by restriction reciprocal \isos $\varphi_1$ and $\psi_1$ between $\Im\psi$ and $\Im\varphi$. In matrix form we obtain:
\vspace{-1mm}
$$
\bordercmatrix [\lbrack\rbrack]{
           & \Im\psi   
                       & \Ker\varphi \cr
\Im\varphi & \varphi_1 
                       &   0    \cr
\Ker\psi   &     0     
                       &0
}~=~\varphi,
\qquad
\bordercmatrix [\lbrack\rbrack]{
           & \Im\varphi      & \Ker\psi \cr
\Im\psi & \psi_1  &   0    \cr
\Ker\varphi   &  0 & 0
}~=~\psi.
$$
\end{fact}

\vspace{1em}
\rems\\
1) If we have a \ali $\psi_0$ satisfying as in \Thref{propIGCram} the \egtz~\hbox{$\varphi\, \psi_0\,\varphi =\varphi$}, we obtain a \ing of $\varphi$ by stating~\hbox{$\psi=\psi_0\,\varphi \,\psi_0$}.
In other words, a \ali $\varphi$ is \lnl \ssi there exists a $\psi$ satisfying
{$\varphi \,\psi \,\varphi =\varphi$}.

\noindent  2) A \nl \ali between free modules of finite rank is \lnl
(immediate \vfnz).

\noindent  3) \Thref{propIGCram} informs us that a \ali which has rank $k$ in the sense of \dfnz~\ref{defRangk} is \lnlz.
\eoe

\begin{fact}
\label{factInvGenCrois} Let $\varphi:\Ae{n}\rightarrow \Ae{m}$ be a \aliz.
\Propeq
\vspace{-2pt}
\begin{enumerate}\itemsep0pt
\item The \ali $\varphi$ is \lnlz.
\item There exists a  $\varphi\bul :\Ae{m}\rightarrow \Ae{n}$ such that

\snic{\Ae{n}=\Ker\varphi\oplus\Im\varphi\bul $ and  $\Ae{m}=\Ker\varphi\bul
\oplus\Im\varphi.}
\item The submodule $\Im\varphi$ is a direct summand in $\Ae{m}$.
\end{enumerate}
\end{fact}
\begin{proof}
\emph{1} $\Rightarrow$ \emph{2.} If $\psi$ is a \ing of $\varphi$, we can take $\varphi\bul=\psi$.
\\
\emph{2} $\Rightarrow$ \emph{3.}  Obvious.
\\
\emph{3} $\Rightarrow$ \emph{1.}
If $\Ae{m}=P\oplus\Im\varphi$, denote by $\pi:\Ae{m}\to\Ae{m}$ the projection
over $\Im\varphi$ parallel to $P$. For each vector $e_i$ of the canonical basis of $\Ae{m}$ there exists an \elt $a_i$ of $\Ae{n}$ such that
$\varphi(a_i)=\pi(e_i)$. We define $\psi :\Ae{m}\rightarrow \Ae{n}$ as
$\psi(e_i)=a_i$. Then, $\varphi \circ  \psi =\pi$ and $\varphi\circ \psi\circ \varphi =\pi\circ \varphi=\varphi$,
and $\psi \circ \varphi\circ  \psi $ is a \ing of $\varphi$.
\end{proof}

The notion of a \lnl \ali is a local notion in the following sense.

\rdb
\begin{plcc}
\label{fact.lnl.loc}
\emph{(Locally simple \alisz)}
Let $S_1$, $\ldots$, $S_n$ be \moco of a \ri $\gA$.
Let $\varphi : \Ae{m}\rightarrow \gA^{q}$ be a \aliz. If every $\varphi _{S_i}:
\gA_{S_i}^m\rightarrow \gA_{S_i}^q$ is \nlz, then $\varphi $ is \lnlz.
More \gnlt  $\varphi $ is \lnl \ssi all the $\varphi _{S_i}$'s are \lnlsz.
\end{plcc}

\begin{proof}
Let us focus on the second statement. To prove that $\varphi $ is \lnl amounts to finding a $\psi$ which satisfies $\varphi \,\psi \,\varphi =\varphi $. This is a \sli in the \coes of the matrix of $\psi$ and we can therefore apply the basic \plgc \ref{plcc.basic}.\iplg
\end{proof}

The terminology of a \lnl \ali is justified by the previous \plg and by the converse given in item \emph{\ref{IFDg}} of \Thref{theoremIFD} (also see the \lnl map lemma in the local \ri case, \paref{lelnllo}).


\subsec{Grassmannians}
\label{subsecGrassmanniennes}

The following \tho serves as an introduction to the grassmannian \vrtsz. It results from Fact~\ref{factInvGenCrois} and \Thref{propIGCram}.

\begin{theorem}\label{propFactDirRangk}
\emph{(\Tf submodules as direct summands of a free module)}
Let $M=\gen{C_1,\ldots,C_m}$ be a \tf submodule of 
$\Ae{n}$ and~\hbox{$C=[\,C_1\;\cdots\;C_m\,]\in\Ae {n\times m}$} be the corresponding matrix.
\begin{enumerate}
\item \Propeq
\begin{enumerate}
\item The matrix $C$ is \lnlz.
\item The module $M$ is a direct summand of $\Ae{n}$.
\item The module $M$ is the image of a matrix $F\in\GAn(\gA)$.
\end{enumerate}
\item \Propeq
\begin{enumerate}
\item The matrix $C$ is of rank $k$.
\item \label{i1propFactDirRangk}The module $M$ is image of a matrix $F\in\GAn(\gA)$ of rank $k$.
\end{enumerate}
\end{enumerate}
\end{theorem}

The \gui{variety} of vector lines in a \Kev of dimension~$n+1$
is, intuitively, of dimension~$n$, as a vector line essentially depends on $n$ parameters (a nonzero vector, up to a multiplicative constant, that makes~$(n+1)-1$
independent parameters).
We call this variety the \pro space of dimension $n$ over $\gK$.
\\
Furthermore, passing from a field $\gK$ to an arbitrary \ri $\gA$, the correct \gnn of a 
\gui{vector line in $\gK^{n+1}$} is
\gui{the image of a \mprn of rank $1$ in $\Ae {n+1}$.}
This leads to the following \dfnsz.

\pagebreak

\begin{definition}\label{defiGrassmanniennes}~
\begin{enumerate}
\item We define the space $\GAnk(\gA)\subseteq\GAn(\gA)$ as the set of \mprns of rank $k$ and $\GGnk(\gA)$ as the set of submodules of $\Ae n$ which are images of matrices of $\GAnk(\gA)$.
\item The space $\GG_{n+1,1}(\gA)$ is again denoted by $\Pn(\gA)$ 
and we call it the \emph{projective space of dimension $n$ over $\gA$}.
\item We denote by $\GGn(\gA)$ the space of all the submodules that are direct summands of~$\Ae n$ 
(i.e., images of a \mprnz).
\end{enumerate}
\index{projective space of dimension $n$ over a \riz}
\end{definition}

The above definition is a little unsatisfactory, insofar as we have not explained how the set
 $\GGnk(\gA)$ is structured. 
Only this structure makes it worthy of the label \gui{space.}

A partial answer is given by the observation that $\GGnk$ is a functor.
More \prmtz,
to every \homo $\varphi:\gA\to\gB$ we associate a natural map $\GGnk(\varphi):\GGnk(\gA)\to\GGnk(\gB)$, so that 
$$\preskip.4em \postskip.3em 
\GGnk(\Id_\gA)=\Id_{\GGnk(\gA)} ,\,\hbox{ and }\, \GGnk(\psi\circ \varphi)=\GGnk(\psi)\circ \GGnk(\varphi), 
$$
when $\psi\circ \varphi$ is defined.

\subsec{Injectivity and surjectivity criteria}
\label{subsecCritInjSurj}

Two famous propositions are contained in the following \thoz.
\begin{theorem}
\label{prop inj surj det}
Let $\varphi:\Ae{n}\rightarrow \Ae{m}$ be a \ali with matrix~$A$.
\begin{enumerate}
\item The map $\varphi$ is surjective \ssi $\varphi$ is of rank $m$,
\cad here $\cD_m(\varphi)=\gen{1}$ (we then say that $A$ is \emph{\umdz}).
\item 
 \emph{(McCoy's theorem)}\index{McCoy!theor@\thoz}
The map $\varphi$ is injective \ssi $\cD_n(\varphi)$ is faithful, \cad if 
the annihilator
of $\cD_n(\varphi)$ is reduced to $\{0\}$.\index{unimodular!matrix}
\end{enumerate}
\end{theorem}

\begin{proof}
\emph{1.} If $\varphi$ is surjective, it admits a right inverse $\psi$, and Fact~\ref{fact.idd prod} gives
$\gen{1}=\cD_m(\I_m)\subseteq \cD_m(\varphi)\cD_m(\psi)$, so $\cD_m(\varphi)=\gen{1}$.
Conversely, if $A$ is of rank $m$, \eqrf{eqIGCram3} shows that $A$ admits a right inverse, and
$\varphi$ is surjective.

\noindent  \emph{2.} Assume that $\cD_n(A)$ is faithful. By \egt
(\ref{eqGema4}), if $AV=0$,  then $\mu_{\alpha,1..n}V=0$ for all the \gtrs $\mu_{\alpha,1..n}$ of  $\cD_n(A)$, and so $V=0$.
\\
For the converse, we will prove by induction on $k$ the following \prtz: \emph{if $k$ column vectors $\xk$ are \lint independent, then the annihilator of the vector $x_1\land\cdots\land x_k$ is reduced to $0$.} For $k=1$ it is trivial. To pass from $k$ to $k+1$ we proceed as follows. Let $z$ be a scalar that annihilates $x_1\land\cdots\land x_{k+1}$. For $\alpha\in \cP_{k,m} $, we denote by $d_\alpha(\yk)$ the minor extracted on the index rows of $\alpha$ for the column vectors~$\yk$ of~$\Ae{m}$. 
Since~$z (x_1\land\cdots\land x_{k+1})=0$, and by the Cramer formulas, 
we have the \egt

\snic{z \,\big(d_\alpha(\xk)x_{k+1}-
d_\alpha(x_1,\ldots, x_{k-1},x_{k+1} )x_k+ \cdots\big) = 0,
}

so $z \,d_\alpha (\xk)=0$. \\
As this is true for any $\alpha$, this gives~\hbox{$z (x_1\land\cdots\land x_{k})=0$}, and by the \hdrz, $z=0$.
\end{proof}
\rem \Thref{prop inj surj det} can also be read in the following way.
\begin{enumerate}\itemsep1pt
\item The \ali $\varphi :\Ae{n}\to\Ae{m}$ is surjective \ssi
the map $\Al{m}\varphi
:\gA^{n\choose m}\to \gA$ is surjective.
\item The \ali $\varphi :\Ae{n}\to\Ae{m}$ is injective \ssi  
the map $\Al{n}\varphi
:\gA  \to\gA^{m\choose n}$ is injective. \eoe
\end{enumerate}

\begin{corollary}
\label{corprop inj surj det}
Let $\varphi:\Ae{n}\rightarrow \Ae{m}$ be an \Aliz.
\begin{enumerate}
\item  If $\varphi$ is surjective and  $n< m$, the \ri is trivial.
\item  If $\varphi$ is injective and  $n> m$, the \ri is trivial.
\end{enumerate}
\end{corollary}

\rem A more positive, \eqvez, but probably even more bewildering formulation of the results of the previous corollary is the following.
\begin{enumerate}
\item If $\varphi$ is surjective, then $X^m$ divides $X^n$ in $\gA[X]$.
\item If $\varphi$ is injective, then $X^n$ divides $X^m$ in $\gA[X]$.
\end{enumerate}
In some way, this is closer to the formulation found in \clamaz:
if the \ri is nontrivial, then  $m\leq n$ in the first case
(resp.\,$n\leq m$ in the second case).\\
The advantage of our formulations is that they work in all cases, without the need to assume that we know how to decide if the \ri is trivial or not. 
\eoe

\begin{corollary}
\label{corInjPuisExt}
If $\varphi :\Ae{n}\rightarrow \Ae{m}$ is injective, the same applies for every exterior power of $\varphi$.
\end{corollary}
\begin{proof}
The annihilator of $\cD_n(\varphi)$ is reduced to $0$ by the previous \thoz. There exists a \ri $\gB\supseteq\gA$ such that the \gtrs of  $\cD_n(\varphi)$ become \com in $\gB$ (Lemma~\ref{lemKerCom}). The \Bli $\varphi_1 :\gB^n\rightarrow \gB^m$ obtained by extending $\varphi$ to $\gB$ is thus of rank $n$ and admits a left inverse $\psi$ (item~\emph{3} of \Thref{propIGCram}), \cad
$\psi \circ\varphi_1 =\Id_{\gB^n}$. Therefore 

\snic{\Al{k} \psi \,\circ\,\Al{k}
\varphi_1 =\Id_{\vi^{\!k} \gB^n}.}

Thus the matrix of $\Al{k}\! \varphi_1$ is injective, and since it is the same matrix as that of $\Al{k}\!\varphi$, the \ali $\Al{k}\! \varphi$  is injective.
\end{proof}

\subsec{Characterization of  \lnl maps} \label{CarLnl}

The following lemma places a bijective correspondence between the \sfios and the non-decreasing sequences of \idms for \dvez.

\pagebreak

\begin{lemma}\label{lem ide-div}
Let $(e_{q+1} =0$, $e_q$, $\ldots$, $e_{1}$, $e_{0}=1)$ be a list of \idms such that $e_{i}$ divides $e_{i+1}$ for $i=0,\ldots ,q$. Then, the \elts 
$r_i:=e_i-e_{i+1}$, for $i\in\lrb{0..q}$,
form a \sfioz. 
Conversely, every \sfioz~\hbox{$(r_{0},\ldots,r_{q})$}
defines such a list of \idms by letting 
$$\preskip.3em \postskip.3em \ndsp
e_j=\sum_{k\geq j}r_{k}\hbox{ for }j\in\lrb{0..q+1}. 
$$
\end{lemma}
\begin{proof} It is clear that $\som_ir_i=1$.
For $0 \leq i <  q$, we have 
$e_{i+1}=e_{i}e_{i+1}$.\\
Hence~$(e_i-e_{i+1})e_{i+1}=0$, \cad
$(r_q+\cdots+ r_{i+1})\cdot r_{i}=0$. We can now easily deduce that $r_ir_j=0$ for $j>i$.
\end{proof}

We denote by $\Diag(a_1,\ldots ,a_n)$ the diagonal matrix of order $n$ whose \coe in position $(i,i)$ is the \elt $a_i$. \label{NOTADiag}

In the following \tho some of the \idms $r_i$ in the \sfio can very well be equal to zero. For example if the \ri is connected and nontrivial, all but one are equal to zero.

\begin{theorem}
\label{theoremIFD}  
\emph{(Locally simple matrix)}%
\index{locally!simple matrix}%
\index{matrix!locally simple ---} \\
Let $G\in \Ae{m\times n}$ be the matrix of 
$\varphi:\Ae{n}\rightarrow \Ae{m}$
 and  $q=\inf(m,n)$.\\
\Propeq
\begin{enumerate}\itemsep=1pt
\item \label{IFDa} 
 The \ali $\varphi$ is \lnlz.

\item \label{IFDb} 
The submodule  $\Im\varphi$ is a direct summand of $\Ae{m}$.

\item \label{IFDc} 
  $\Im\varphi$  is a direct summand of $\Ae{m}$ and  $\Ker\varphi$
  is a direct summand of  $\Ae{n}$.

\item \label{IFDd} 
  There exists a \ali  $\varphi\bul :\Ae{m}\rightarrow \Ae{n}$
  with $\Ae{n}=\Ker\varphi\oplus\Im\varphi\bul $
  and  $\Ae{m}=\Ker\varphi\bul \oplus\Im\varphi$.

\item \label{IFDe} 
  Each \idd $\cD_k(\varphi)$ is \idmz.

\item \label{IFDf1} 
There exists a (unique) \sfio \hbox{$(r_0, r_1, \dots , r_q)$}
 such that on each localized ring $\gA[1/r_k]$ 
 the map $\varphi$ is of rank~$k$.

\item \label{IFDf} 
  Each \idd $\cD_k(\varphi)$ is generated by an \idm $e_k$.
  Then let $r_k=e_k-e_{k+1}$. The $r_k$'s form a \sfioz.
  For every minor $\mu$ of order $k$ of $G$,
  on the localized ring $\gA[1/(r_k\,\mu)]$ the \ali $\varphi$
  becomes \nl of rank~$k$.

\item \label{IFDg} 
  The \ali $\varphi$ becomes \nl after \lon at suitable \ecoz.

\item \label{IFDi} 
Each \idd $\cD_k(\varphi)$ is generated by an \idm $e_k$ and the matrix of $\varphi$ becomes equivalent to the matrix $\Diag(e_1,e_2, \ldots ,e_q)$, eventually filled-in with zeros (for both rows and columns), after \lon at suitable \ecoz.

\item $\!\!\!\sta$ \label{IFDh} 
The \ali $\varphi$ becomes \nl after \lon at any arbitrary maximal \idz.
\end{enumerate}
\end{theorem}

\begin{proof}
The \eqvc of items \emph{\ref{IFDa}}, \emph{\ref{IFDb}}, \emph{\ref{IFDc}}, \emph{\ref{IFDd}}
is already clear (see Facts \ref{factIng0} and~\ref{factInvGenCrois}).
Furthermore, we trivially have \emph{\ref{IFDf}} $\Rightarrow$ \emph{\ref{IFDf1}}
$\Rightarrow$ \emph{\ref{IFDe}} and \emph{\ref{IFDi}} $\Rightarrow$ \emph{\ref{IFDe}.}

\noindent  Since $q=\inf(m,n)$, we have  $\cD_{q+1}(\varphi)=0$.

\noindent  \emph{\ref{IFDa}} $\Rightarrow$ \emph{\ref{IFDe}.} 
We have $GHG=G$ for some matrix $H$ and we apply Fact~\ref{fact.idd prod}.

\noindent  \emph{\ref{IFDe}} $\Rightarrow$ \emph{\ref{IFDf}.}  
The fact that each $\cD_k(\varphi)$ is generated by an \idm $e_k$ results from Fact~\ref{lem.ide.idem}.
The fact that $(r_0,\ldots,r_q)$ is a \sfio results from Lemma~\ref{lem ide-div} (and Fact~\ref{fact.idd inc}). \\
As
$r_ke_{k+1}=0$,  over the \ri $\gA[1/r_k]$, and thus over the \ri $\gA[1/(\mu r_k)]$, where $\mu$ is a minor of order $k$, every minor of order $k+1$ of the matrix~$G$ is null.
Thus,  by the freeness lemma \ref{lem pf libre}, $G$ is \nl of rank $k$.

\noindent  \emph{\ref{IFDf}} $\Rightarrow$ \emph{\ref{IFDi}.} 
Over $\gA[1/r_k]$ and so over $\gA[1/(\mu r_k)]$ ($\mu$ a minor of order $k$), 
we have~$\Diag(e_1,\ldots ,e_q)=\Diag(1,\ldots ,1,0,\ldots ,0)$ with $1$ appearing $k$ times.

\noindent  \emph{\ref{IFDf}} $\Rightarrow$ \emph{\ref{IFDg}.}
Let $t_{k,j}$ be the minors of order $k$ of $G$. 
The \lons are those at $t_{k,j}r_k$.
We must verify that they are \comez.
Each  $e_k$ is in the form $\sum t_{k,j}v_{k,j}$,
so $\sum_{k,j}v_{k,j}(t_{k,j}r_k)= \sum_ke_kr_k =\sum r_k=1$.

\noindent \emph{\ref{IFDg}} $\Rightarrow$ \emph{\ref{IFDa}.} By application of the
\plgrf{fact.lnl.loc} since every \nl map is \lnlz.

\noindent \emph{\ref{IFDg}}   $\Rightarrow$ \emph{\ref{IFDh}.}
(In \clamaz.)
 Because the complement of a \idema always contains at least one \elt in a \sys of \eco (we can assume that the \ri is nontrivial).

\noindent \emph{\ref{IFDh}} $\Rightarrow$ \emph{\ref{IFDg}.}
 (In \clamaz.) For each \idema
$\fm$ we obtain a $s_\fm\notin\fm$ and a matrix $H_\fm$ such that we have $GH_\fm G=G$ in~$\gA[1/s_\fm]$. The \id generated by the $s_\fm$'s is not contained in any \idema and so it is the \id $\gen{1}$. Thus there is a finite number of these $s_\fm$'s which are \comz.

\noindent Let us finish by giving a direct \dem for the implication \emph{\ref{IFDf1}} $\Rightarrow$ \emph{\ref{IFDa}.}\\
 On the \riz~$\gA[1/r_k]$ the matrix $G$ is of rank $k$ so there exists a matrix~$B_k$ satisfying $GB_kG=G$ (\Thref{propIGCram}). This means that on the \riz~$\gA$ we have a matrix~$H_k$ in $\Ae {n\times m}$ satisfying $r_kH_k=H_k$ and $r_kG=GH_kG$. We then take $H=\sum_k H_k$ and obtain $G=GHG$.
\end{proof}

The \eqvc of items \emph{\ref{IFDa}} to \emph{\ref{IFDi}} has been established \cotz, whilst item \emph{\ref{IFDh}} only implies the previous ones in \clamaz.

\subsec{Trace, norm, discriminant, transitivity}
\label{subsecTransDet}

We denote by $\Tr(\varphi)$ and $\rC{\varphi}(X)$ the trace and the \emph{\polcarz} of an \endo $\varphi$ of a free module of finite rank (we take as \polcar of a matrix $F\in \Mn(\gA)$ the \pol $\det(X\I_n-F)$, which has the advantage of being \monz).%
\index{polynomial!characteristic --- of an endomorphism}%
\index{characteristic!polynomial of an endomorphism}%
\label{NOTA1Polcar}

\pagebreak

\begin{notation}\label{notaCTrN}~
\begin{enumerate}
\item [--] If $\gA\subseteq\gB$ and if $\gB$ is a free \Amo of finite rank,
we denote  $\rg_\gA(\gB)$ by $\dex{\gB:\gA}$.
\item [--]
For $a\in\gB$ we then denote by $\Tr_{\gB/\!\gA}(a)$, $\rN_{\gB/\!\gA}(a)$ and $\rC{\gB/\!\gA}(a)(X)$ the trace, the \deter and the \polcar of the multiplication by~$a$,
seen as an \endo of the \Amo $\gB$.
\end{enumerate}%
\index{polynomial!characteristic --- of an element}%
\index{characteristic!polynomial of an element}
\end{notation}

\vspace{.3em}
\begin {lemma}\label{lem1TransDet}
 Assume that
  $\gA\subseteq\gB$ and that $\gB$ is a free \Amo of finite rank $m$.
\begin{enumerate}
\item  Let $E$ be a free \Bmo of finite rank $n$. If $\ue = (e_i)_{i\in \lrbm}$ is a basis of $\gB$ over $\gA$ and $\uf = (f_j)_{j\in \lrbn}$ a basis of $E$ over $\gB$, then
$(e_if_j)_{i,j}$ is a basis of $E$ over $\gA$.
Consequently, $E$ is free over $\gA $ and
$$\preskip.4em \postskip.4em 
\rg_\gA( E) = \rg_\gB  (E) \times \rg_\gA (\gB) . 
$$
\item If $\gB\subseteq\gC$ and if $\gC$ is a free \Bmo of finite rank,
we have
$$\preskip.4em \postskip.4em 
\dex{\gC:\gA}=\dex{\gC:\gB}\,\dex{\gB:\gA}.
$$
\end{enumerate}
\end {lemma}

\vspace{.3em}
\rem
 Let $\gC = \aqo{\gA[Y]}{Y^3}=\gA[y]$, a free \Alg of rank~$3$. Since $y^4 = 0$, $\gB  = \gA
\oplus \gA y^2$ is a sub-\alg of $\gC$, free over $\gA$, whose rank (equal to $2$) does not divide the rank of $\gC$ (equal to $3$). 
The \egtz~\hbox{$\dex{\gC :  \gA} = \dex{\gC :  \gB } \dex{\gB  :\gA}$} does not apply because $\gC$ is not free over $\gB$. 
\eoe


\begin{theorem}\label{Th.transitivity} \emph{(Transitivity formulas for the trace, the \deter and the \polcarz)}
Under the same assumptions, let $u_\gB  :  E \to E$ be a \Bliz. We denote by $u_\gA$ this map when  considered as an \Aliz. We then have the \egtsz:
\begin{enumerate}
\item [ ] $\det(u_\gA) = \rN_{\gB/\!\gA}\big(\det(u_\gB )\big) $,
 $\;\Tr(u_\gA) = \Tr_{\gB /\!\gA  } \big( \Tr(u_\gB )\big)$,
\item [ ] $\rC{u_\gA}(X) = \rN_{\gB [X]/\!{\gA[X]}} \big(
\rC{u_\gB}(X)\big).$
\end{enumerate}
\end{theorem}

\begin {proof} We use the notations of Lemma~\ref{lem1TransDet}.
Let $u_{kj}$ be the \elts of $\gB$ defined by $u(f_j) = \sum_{k=1}^n u_{kj} f_k$. Then the matrix $M$ of $u_\gA$ with respect to the basis $(e_if_j)_{i,j}$ is expressed as 
a block matrix 
$$\preskip.4em \postskip.4em 
 M = \cmatrix {
M_{11} &  \cdots & M_{1n} \cr
\vdots &        &    \vdots \cr
M_{n1} &  \cdots & M_{nn} \cr
}, 
$$
where $M_{kj}$ represents the \Ali $b \mapsto b u_{kj}$ of $\gB $ in $\gB$ to with respect to the basis~$\ue$.
This provides the desired \egt regarding the trace of $u_A$ since 
$$
{ \arraycolsep2pt
\begin{array}{ccccccc}
 \Tr(u_\gA) & =  & \ds \som_{i=1}^n \Tr(M_{ii}) &=&\ds \som_{i=1}^n \Tr_{\gB /\gA}(u_{ii})\\[.5em]
  &  = & \ds \Tr_{\gB /\gA}  \big(\som_{i=1}^n u_{ii}\big) &=& \Tr_{\gB /\gA}\big(\Tr(u_\gB)\big).\
\end{array}
}
$$
As for the \egt for the \deterz, note that the matrices $M_{ij}$ pairwise commute  ($M_{ij}$ is the matrix of the multiplication by $u_{ij}$). We can then apply the following Lemma~\ref{lemdeterblocs}, which gives us:
$$\preskip.4em \postskip.4em 
\det(M) = \det(\Delta) \quad \hbox {with} \quad
\Delta = \som_{\sigma \in \Sn} \varepsilon(\sigma)
M_{1\sigma_1} M_{2\sigma_2} \ldots M_{n\sigma_n}. 
$$
However, $\Delta$ is none other than the matrix of the multiplication by the \elt
$
\som_{\sigma \in \Sn} \varepsilon(\sigma) u_{1\sigma_1} u_{2\sigma_2}
\ldots u_{n\sigma_n}, 
$ 
 i.e., by
$\det(u_\gB )$, thus:
$$\preskip.3em \postskip.4em 
\det(u_\gA) = \det(M) = \rN_{\gB /\gA}\big(\det(u_\gB )\big). 
$$
Finally, the \egt for the \polcar is deduced 
from the one for \deters 
by using the fact that $\rC{u_\gA}(X)$ is the
\deter of the \endo $X\Id_{E[X]} - u_\gA $ 
         of the $\gA[X]$-module $E[X]$ 
whereas $\rC{u_\gB}(X)$ is that of the same map seen as an \endo of the $\gB[X]$-module~$E[X]$.
\end {proof}

In a noncommutative \riz, two \elts $a$ and $b$ are said to be \emph{permutable} or \emph{commuting} if $ab=ba$.

\begin {lemma}\label{lemdeterblocs}
Let $(N_{ij})_{i,j}$ be a family of 
$n^2$ pairwise commuting square matrices, and $N$ the square matrix of order $mn$:
$$
N = \cmatrix {
N_{11} &  \cdots & N_{1n} \cr
\vdots &        &    \vdots \cr
N_{n1} &  \cdots & N_{nn} \cr
}.
$$
Then:
$
\det(N) = \det\big( \som_{\sigma \in \Sn} \varepsilon(\sigma)
\,N_{1\sigma_1} N_{2\sigma_2} \ldots N_{n\sigma_n} \big).
$
\end {lemma}

\begin {proof}
  Let $\Delta$ be the $n \times n$ matrix defined by $\Delta = \sum_{\sigma
\in \Sn} \varepsilon(\sigma) N_{1\sigma_1} N_{2\sigma_2} \ldots
N_{n\sigma_n}$. Thus we must prove that $\det(N) = \det(\Delta)$.

\noindent Let us treat the special cases $n = 2$ then $n = 3$. We replace $\gA$ with $\gA[Y]$ and $N_{ii}$ by $N_{ii} + Y \rI_m$, which has the advantage of making some \deters \ndzs in $\gA[Y]$.
It suffices to establish the \egts with these new matrices, as we finish by making $Y=0$.

\noindent 
The key-element 
of the \dem for $n = 2$ resides in the following \egtz:
$$
\cmatrix {N_{11} & N_{12} \cr N_{21} & N_{22} \cr}
\cmatrix {N_{22} & 0      \cr -N_{21} & \I_m \cr} =
\cmatrix {N_{11}N_{22} - N_{12}N_{21} & N_{12} \cr 0 & N_{22} \cr}.
$$
We then consider the LHS and RHS \deters
$$
\det(N) \det(N_{22}) = \det(N_{11}N_{22} - N_{12}N_{21}) \det(N_{22}),
$$
 next we simplify by $\det(N_{22})$ (which is \ndzz) to obtain the result.

\noindent Case $n = 3$ uses the \egtz:

\snac{
\cmatrix {
N_{11} & N_{12} & N_{13} \cr
N_{21} & N_{22} & N_{23} \cr
N_{31} & N_{32} & N_{33} \cr}
\cmatrix {
N_{22}N_{33} - N_{23}N_{32} & 0   & 0   \cr
N_{31}N_{23} - N_{21}N_{33} & \I_m & 0   \cr
N_{21}N_{32} - N_{22}N_{31} & 0   & \I_m \cr} =
\cmatrix {
\Delta & N_{12} & N_{13} \cr
0      & N_{22} & N_{23} \cr
0      & N_{32} & N_{33} \cr},
}  

which leads to

\snic{\det(N) \det(N_{22}N_{33} - N_{23}N_{32}) = \det(\Delta) \det
\cmatrix {N_{22} & N_{23}\cr N_{32} & N_{33}\cr}.}

Case $n = 2$ provides $\det(N_{22}N_{33} - N_{23}N_{32}) = \det \cmatrix {N_{22} & N_{23}\cr
N_{32} & N_{33}\cr}$. We simplify by this \deter and obtain $\det(N) = \det(\Delta)$.

\noindent   The \gnl case is left as an exercise 
(see Exercise~\ref{exolemdeterblocs}).
\end {proof}

\begin {corollary}
Let $\gA  \subseteq \gB  \subseteq \gC$ be three \ris with $\gC$ free of finite rank over $\gB $ and $\gB $ free of finite rank over $\gA $.  We then have:
\begin{enumerate}
\item [ ] $\rN_{\gC/\!\gA} = \rN_{\gB /\!\gA} \circ  \rN_{\gC/\gB }  $,
 $\quad\Tr_{\gC/\!\gA} = \Tr_{\gB /\!\gA}   \circ \Tr_{\gC/\gB } $,
\item [ ] $\rC{\gC/\!\gA}(c)(X)= \rN_{\gB [X]/\!\gA[X]}\!  \big(\rC{\gC/\gB}(c)( X)\big)\;\; (c\in\gC)$.
\end{enumerate}
\end {corollary}


\subsubsection*{Gram determinants and discriminants}

\begin{definition}\label{defiGram}
Let $M$ be an \Amoz, $\varphi:M\times M\to\gA$ be a \smq bi\lin form and $(\ux)=(\xk)$ be a list of \elts of $M$. We call  the matrix 
$$\preskip.0em \postskip.2em
\Gram_\gA(\varphi,\ux)\eqdefi\big(\varphi(x_i,x_j)\big)_{i,j\in\lrbk}
$$
the \emph{Gram matrix of $(\xk)$ for $\varphi$}.
Its \deter is called the \emph{Gram \deter of $(\xk)$ for $\varphi$}
and is denoted by $\gram_\gA(\varphi,\ux)$.%
\index{matrix!Gram ---}%
\index{determinant!Gram ---}%
\index{Gram!matrix}%
\index{Gram!determinant}%
\end{definition}

If $\gA y_1 + \cdots + \gA y_k \subseteq \gA x_1 + \cdots + \gA x_k$ we have an \egt

\snic{
\gram(\varphi,\yk)=\det(A)^2\gram(\varphi,\xk),}

 where $A$ is a $k\times k$ matrix which expresses the $y_j$'s in terms of the $x_i$'s.

\smallskip We now introduce an important case of a Gram \deterz, the \discriz.
Recall that two \elts $a$, $b$ of a \ri $\gA$ are said to be
 \emph{associated}
if there exists a $u\in\Ati$ such that $a=ub$. In the literature such elements are also referred to as \emph{associates}.%
\index{associated!elements in $\gA$}\index{associates!\emph{see associated elements}}

\begin{propdef}\label{defiDiscTra}
Let $\gC\supseteq\gA$ be an \Alg which is \hbox{a free \Amoz} of finite rank and $x_1$, \dots, $x_k$, $y_1$, \dots, $y_k\in\gC$.
\begin{enumerate}
\item We call 
the \deter of the matrix
$$\preskip.4em \postskip.0em 
\big(\Tr_{\gC/\!\gA}(x_ix_j)\big)_{i,j\in\lrbk} 
$$
the \ix{discriminant} of $(\xk)$. We denote it by $\disc_{\gC/\!\gA}(\xk)$ or $\disc(\xk)$.
\item If $\gA y_1 + \cdots + \gA y_k \subseteq \gA x_1 + \cdots + \gA x_k$ we have
$$\preskip.3em \postskip.3em 
\disc(\yk)=\det(A)^2\disc(\xk), 
$$
 where $A$ is a $k\times k$ matrix which expresses the $y_j$'s 
in terms of the $x_i$'s.
\item In particular, if $(\xn)$ and $(\yn)$ are two bases of the \Alg $\gC$,
 the \elts $\disc(\xn)$ and $\disc(\yn)$ are multiplicatively congruent modulo the squares of $\Ati$.
We call the corresponding \eqvc class the \emph{\discri of the extension} $\gC\sur\gA$. We denote it by $\Disc_{\gC/\!\gA}$.%
\item If $\Disc_{\gC/\!\gA}$ is \ndz and $n=\dex{\gC:\gA}$, a \sys $\un$ in $\gC$
is an $\gA$-basis of $\gC$ \ssi $\disc(\un)$ and $\Disc_{\gC/\!\gA}$
are associated elements.

\end{enumerate}
\end{propdef}

For example when $\gA=\ZZ$ the \discri of the extension
is a well-defined integer, whereas if $\gA=\QQ$, the \discri is \care on the one hand by its sign, and on the other hand by the list of prime numbers contained therein with an odd power.

\begin{proposition}\label{propTransDisc}
Let $\gB$ and $\gC$ be two free \Algs of ranks $m$ and~$n$, respectively
, and consider the product algebra $\gB \times \gC$.

Given a list~\hbox{$(\ux)=(\xm)$} of \elts of $\gB$ and a list $(\uy)=(\yn)$ of \elts of $\gC$, we have:
$$\preskip.2em \postskip.4em
\;\;\disc_{(\gB\times\gC)/\!\gA}(\ux,\uy) =
\disc_{\gB/\!\gA}(\ux) \times \disc_{\gC /\!\gA}(\uy).
$$
In particular,
$\Disc_{(\gB\times \gC)/\!\gA}=\Disc_{\gB/\!\gA}\times\Disc_{\gC /\!\gA}.$
\end {proposition}
\facile

\begin{proposition}\label{prop1TransDisc}
Let $\gB\supseteq\gA$ be a free \Alg of finite rank $p$.  \\
We consider

\vspace{-.1em} 
\begin{itemize}\itemsep=0pt
\item a \Bmo $E$,
\item a \smq $\gB$-bi\lin form
$\varphi_\gB: E \times E \to \gB$,
\item  a basis $(\ub) = (b_i)_{i \in\lrbp}$ of $\gB$ over $\gA$, and
\item  a family $(\ue) = (e_j)_{j \in\lrbn}$ of $n$ \elts of $E$.
\end{itemize}
Let $(\ub  \star  \ue)$ be a family $(b_i e_j)$ of $np$ \elts of $E$ and $\varphi_\gA: E \times E \to \gA$
be the \smq $\gA$-bi\lin form defined by:
$$\preskip.4em \postskip.4em
\varphi_\gA(x, y) = \Tr_{\gB/\!\gA}\big(\varphi_\gB(x, y)\big).
$$
We then have the following transitivity formula:
$$\preskip.4em \postskip.4em
\gram(\varphi_\gA, \ub \star \ue) =
\disc_{\gB/\!\gA}(\ub)^n \times \rN_{\gB/\!\gA} \big(\gram(\varphi_\gB, \ue)\big).
$$
\end{proposition}
\begin {proof}
In the following the indices $i$, $i'$, $k$, $j$, $j'$ satisfy $i$, $i'$, $k\in\lrbp$ \hbox{and $j$, $j'\in\lrbn$}. Let us agree to sort $\ub \star \ue$ in the following order:

\snic{\ub \star \ue = 
b_1 e_1, \ldots, b_p e_1,
b_1 e_2, \ldots, b_p e_2, \ldots,
b_1 e_n, \ldots, b_p e_n .
}

For $x \in \gB$, let $\mu_x:  \gB \to \gB$ be the multiplication by $x$ and $m(x)$ be the matrix of $\mu_x$ with respect to the basis $(b_i)_{i\in\lrbp}$ of $\gB$ over $\gA$. 
Thus we define an \isoz~$m$ of the \ri $\gB$ into a commutative sub\ri of $\MM_p(\gA)$. If we let $m_{ki}(x)$ be the \coes of the matrix $m(x)$,
we then have:

\snic{\dsp\mu_x(b_i) = b_i x = \som_{k=1}^p m_{ki}(x) b_k,}

with $\rN_{\gB/\!\gA}(x) = \det\big(m(x)\big)$.  By letting $\varphi_{jj'} = \varphi_\gB(e_j, e_{j'}) \in \gB$, we have
$$\preskip.3em \postskip.3em 
\varphi_\gA( b_i e_j
b_{i'} e_{j'} ) = \Tr_{\gB/\!\gA}\big(\varphi_\gB( b_i e_j b_{i'} e_{j'} )\big) = \Tr_{\gB/\!\gA} (b_i b_{i'} \varphi_{jj'} ). 
$$
  By using the \egt {\mathrigid 2mu $b_{i'} \varphi_{jj'} =
\som_{k=1}^p m_{ki'}(\varphi_{jj'})\, b_k$, we have with 
$\Tr=\Tr_{\gB/\!\gA}$}:
{\small
$$
\Tr ( b_i b_{i'} \varphi_{jj'} ) =
\Tr \big( \som_{k = 1}^p b_i\, m_{ki'}(\varphi_{jj'})\, b_k \big) =
\som_{k = 1}^p \Tr(b_i b_k)\, m_{ki'}(\varphi_{jj'}).\eqno(*)
$$
}

\vspace{-12pt}\noindent We define $\beta \in \MM_p(\gA)$ by $\beta_{ik} = \Tr_{\gB/\!\gA}(b_i
b_k)$. The right-hand sum in~$(*)$ is none other than the \coe of a product of matrices: $\big(\beta \cdot m(\varphi_{jj'})\big)_{ii'}$.  The Gram \deter of $\ub \star \ue$ for $\varphi_\gA$ is therefore an $np \times np$ matrix comprised      of $n^2$ blocks of $p \times p$ matrices. Here is that matrix if we let~\hbox{$\phi_{jj'}=m(\varphi_{jj'})$} to simplify the expression:
{
\small
$$\preskip.2em \postskip.3em\!\!
\cmatrix {
\beta \phi_{11} &\beta  \phi_{12} & \ldots &\beta  \phi_{1n} \cr
\beta \phi_{21} &\beta  \phi_{22} & \ldots &\beta  \phi_{2n} \cr
      \vdots       &                    &        &       \vdots       \cr
\beta \phi_{n1} &\beta  \phi_{n2} & \ldots &\beta  \phi_{nn} \cr
} =
\cmatrix {
\beta  & 0     & \ldots & 0      \cr
0      & \beta & \ldots & \vdots \cr
\vdots &       & \ddots & 0      \cr
0      &       & \ldots & \beta  \cr
}
\cmatrix {
\phi_{11} & \phi_{12} & \ldots & \phi_{1n} \cr
\phi_{21} & \phi_{22} & \ldots & \phi_{2n} \cr
\vdots       &              &        & \vdots       \cr
\phi_{n1} & \phi_{n2} & \ldots & \phi_{nn} \cr
}\!.
$$
}

\vspace{-10pt}
\noindent 
By taking the \deters we obtain
$$\preskip.0em \postskip.3em
\gram(\varphi_\gA, \ub \star \ue) = \det(\beta)^n \cdot
\det \cmatrix {
\phi_{11} & \phi_{12} & \ldots & \phi_{1n} \cr
\phi_{21} & \phi_{22} & \ldots & \phi_{2n} \cr
\vdots       &              &        & \vdots       \cr
\phi_{n1} & \phi_{n2} & \ldots & \phi_{nn} \cr
}.
$$

By using the fact that the matrices $\phi_{jl}$ pairwise commute, we find that the right-\deter is equal to
$$\preskip.3em \postskip.2em\mathrigid1mu
\det\big( \sum_{\sigma \in \Sn} \varepsilon(\sigma)
\phi_{1\sigma_1} \phi_{2\sigma_2} \ldots
\phi_{n\sigma_n} \big) \,=\;
\det\, m\big(\det(\varphi_{jl})\big) \,=\; \rN_{\gB/\!\gA} \big(\gram(\varphi_\gB, \ue)\big),
$$
as required.
\end {proof}
%

\begin{theorem}\label{thTransDisc} \emph{(Transitivity formula for the \discrisz)}\\
 Let $\gA \subseteq \gB \subseteq \gC$, with $\gB$ free over $\gA$,  $\gC$ free over $\gB$, \hbox{$\dex{\gC: \gB} = n$} {and $\dex{\gB: \gA} = m$}.
Let $(\ub) = (b_i)_{i\in\lrbm}$ be a basis of $\gB$ over~$\gA$,~\hbox{$(\uc) = (c_j)_{j\in\lrbn}$} be a basis of~$\gC$ over $\gB$ and let $(\ub \star \uc)$ be the basis $(b_i c_j)$ of $\gC$ over $\gA$.
Then:
\[\arraycolsep2pt\preskip.3em \postskip.3em
\begin{array}{rcl}
\disc_{\gC/\!\gA}(\ub \star \uc)  &  = & \disc_{\gB/\!\gA}(\ub)^{\dex{\gC: \gB}}\ \rN_{\gB/\!\gA}\big(\disc_{\gC/\gB}(\uc)\big) , \\ [1.3mm]
\hbox{and so}\quad \Disc_{\gC/\!\gA}  &  = & \Disc_{\gB/\!\gA}^{~~~\dex{\gC: \gB}}\
 \rN_{\gB/\!\gA}(\Disc_{\gC/\gB}).  \end{array}
\]
\end {theorem}

\vspace{1em}
\begin {proof}
Direct application of Proposition~\ref{prop1TransDisc}.
\end {proof}
%

\section{Basic local-global principle for modules}
\label{secPLGCBasicModules}

This section's results will not be used before Chapter~\ref{chap ptf0}. 

We are about to give a slightly more \gnl version of the basic \plg \vref{plcc.basic}. This new principle concerns arbitrary \Amos and \alisz, whilst the basic principle can be considered as the special case where the modules are free and of finite rank. The \dem is essentially the same as that of the basic principle.\iplg
\\
Beforehand, we start with a brief review of exact sequences and we establish some \elr \prts of the \lon regarding modules.

\subsec{Complexes and exact sequences}
\index{exact sequence!of \alisz}

When we have successive \alis 

\snic{M\vers{\alpha}N\vers{\beta}P\vers{\gamma}Q\;,}

we say that they form a \ix{complex} if the composition of 
any two successive \alis is null.
We say that the sequence is \emph{exact in $N$} if $\Im\alpha=\Ker\beta$. The entire sequence is said to be
 exact if it is exact in $N$ and $P$. This extends to sequences of arbitrary length.

This \gui{abstract} language has an \imde counterpart in terms of \slis when we are dealing with free modules of finite rank. For example if $N=\Ae n$, $P=\Ae m$ and if we have an exact sequence   

\snic{0\to M\vers{\alpha}N\vers{\beta}P\vers{\gamma}Q\to 0\;,}

The \ali $\beta$ is represented by a matrix associated with a \sys of~$m$ \lin \eqns with $n$ unknowns, the module $M$, isomorphic to $\Ker \beta$, represents the 
defect of injectivity 
of $\beta$ and the module $Q$, isomorphic to $\Coker \beta$, represents its 
defect of surjectivity
of $\beta$.

\smallskip An exact complex of the type
\[
\begin{array}{ccccccccccccccccccccccccccccccc} 
    0&\to &M_m&\vvers{u_m}& M_{m-1}&\lora& \cdots\cdots\cdots& \vvers{u_1}& M_0&\to& 0
 \end{array}
\]
with $m\geq 3$ is called a \emph{long exact sequence (of length $m$)}.

If $m=2$, we say that we have a \emph{short exact sequence}. In this case $M_2$ can be identified with a submodule of $M_1$, and, modulo this identification, $M_0$ can be identified with $M_1/M_2$.

\rdb
An important fact to note is that every long exact sequence of length $m$ \gui{can be decomposed into} $m-1$ short exact sequences according to the following schema. \label{sexaseco}
\[ 
\begin{array}{ccccccccccccccccccccccccccccccc} 
 0 & \to  & E_2 & \vvers{\iota_2} & M_1 & \vvers{u_1} & M_0 &\to & 0\\[1mm] 
 0 & \to  & E_3 & \vvers{\iota_3} & M_2 & \vvers{v_2} & E_2 &\to & 0\\[1mm] 
   & \vdots   &   &   &   &  &   &\vdots \\[1mm] 
 0 & \to  & E_{m-1} & \vvers{\iota_{m-1}} & M_{m-2} & \vvers{v_{m-2}} & E_{m-2} &\to & 0\\[1mm] 
 0 & \to  & M_m & \vvers{u_m} & M_{m-1} & \vvers{v_{m-1}} & E_{m-1} &\to & 0\\
 \end{array}
\]
with $E_i=\Im u_{i+1}\subseteq M_i$ for $i\in\lrb{2..m-1}$, 
the $\iota_k$'s canonical injections, and the $v_k$'s obtained from the $u_k$'s by restricting the range to $\Im u_k$.

\smallskip  An important theme of commutative \alg is provided by the transformations that preserve, or do not preserve, exact sequences.

Here are two basic examples, which use the modules of \alisz.

Let $\Lin_\gA(M,P)$ be the \Amo of  \Alis from $M$ to~$P$ and $\End_\gA(M)$ designate $\Lin_\gA(M,M)$ (with its \ri structure \gnlt noncommutative).
The \emph{dual module} of $M$, $\Lin_\gA(M,\gA)$, will in \gnl be denoted by~$M\sta$.%
\label{NOTAAlis}\index{module!dual ---}

\begin{fact}\label{fact1HomEx}
If $0\to M\vers{\alpha}N\vers{\beta}P$ is an exact sequence of \Amosz, and if~$F$ is an \Amoz, then the sequence 
$$\preskip.4em \postskip.0em 
0\to \Lin_\gA(F,M)\lora \Lin_\gA(F,N)\lora \Lin_\gA(F,P) 
$$
is exact.
\end{fact}
%
\begin{proof} \emph{Exactness in $\Lin_\gA(F,M)$.}
Let $\varphi\in\Lin_\gA(F,M)$ such that $\alpha\circ \varphi=0$. 
Then, since the first sequence is exact in $M$, for all $x\in F$, $\varphi(x)=0$, so~$\varphi=0$.

\noindent \emph{Exactness in $\Lin_\gA(F,N)$.}  
Let $\varphi\in\Lin_\gA(F,N)$ such that $\beta\circ \varphi=0$.
Then, since the first sequence is exact in $N$, for all $x\in F$, $\varphi(x)\in\Im\alpha$. \\
Let $\alpha_1:\Im\alpha\to M$ be the inverse of the bijection $\alpha$ (regarding the codomain of~$\alpha$ as  $\Im\alpha$) and $\psi=\alpha_1\, \varphi$.
\\
 We then obtain the \egts
$\Lin_\gA(F,\alpha)(\psi)=\alpha\, \alpha_1\, \varphi=\varphi$.
\end{proof}
%
\begin{fact}\label{fact2HomEx}
If $N\vers{\beta}P\vers{\gamma}Q\to 0$ is an exact sequence of \Amos and if~$F$ is an \Amoz, then the sequence 
$$\preskip.4em \postskip.0em 
0\to \Lin_\gA(Q,F)\lora \Lin_\gA(P,F)\lora \Lin_\gA(N,F) 
$$
is exact. 
\end{fact}
%
\begin{proof} \emph{Exactness in $\Lin_\gA(Q,F)$.} If $\varphi\in\Lin_\gA(Q,F)$ satisfies $\varphi\circ \gamma=0$, then, since~$\gamma$ is surjective, $\varphi=0$.

\noindent \emph{Exactness in $\Lin_\gA(P,F)$.} If $\varphi:P\to F$ satisfies $\varphi\circ \beta=0$, then $\Im\beta\subseteq\Ker\varphi$ and $\varphi$ is factorized by $P\sur{\,\Im\beta}\simeq Q$,
that is $\varphi=\psi\circ \gamma$ for a \ali $\psi:Q\to F$, \cad $\varphi\in\Im\Lin_\gA(\gamma,F)$.
\end{proof}

\vspace{-.7em}
\pagebreak

\begin{fact} \label{factDualReflexif}
Let $\beta:N\to P$ be a \ali and $\gamma:P\to\Coker\beta$ be the canonical \prnz.  
\begin{enumerate}
\item The canonical map $\tra \gamma:(\Coker\beta)\sta \to P\sta$ induces an \iso of $(\Coker\beta)\sta$ on $\Ker\!\tra{\beta}$.

\item If the canonical \alis $N\to N^{\star\star}$ and $P\to P^{\star\star}$ are \isosz, then the canonical surjection of $N\sta$ in $\Coker\!\tra{\beta}$ provides by duality an \iso of $(\Coker\! \tra{\beta})\sta $ on $\Ker\beta$.
\end{enumerate}

\end{fact}

\begin{proof}
\emph{1.} We apply Fact~\ref{fact2HomEx} with $F=\gA$.

 \emph{2.} 
We apply item \emph{1} to the \ali $\tra\beta$ by identifying $N$ and~$N^{\star\star}$, as well as $P$ and~$P^{\star\star}$, and thus also $\beta$ and $^{\rm t}{(\!\tra{\beta})}$.
\end{proof}

\rem It is possible to slightly weaken the assumption by requiring that the \ali $P\to P^{\star\star}$ be injective.
\eoe

\subsec{\Lon and exact sequences}

\begin{fact}\label{fact.sexloc} Let $S$ be a \mo of a \ri $\gA$.
\begin{enumerate}
\item
 If $M$ is a submodule of $N$, we have the canonical identification of $M_S$ with a submodule of $N_S$ and of $(N/M)_S$ with $N_S/M_S$.\\
In particular,  for every \id $\fa$ of $\gA$, the \Amo $\fa_S$ is canonically identified with the \id $\fa\gA_S$ of $\gA_S$.
\item   If $\varphi:M\rightarrow N$ is an \Aliz, then:
\begin{enumerate}
\item    ${\rm  Im}(\varphi_S)$  is canonically identified with
$\big({\rm  Im}(\varphi)\big)_S$\,,
\item    ${\rm  Ker}(\varphi_S)$  is canonically identified with
$\big({\rm  Ker}(\varphi)\big)_S$\,,
\item    ${\rm  Coker}(\varphi_S)$  is canonically identified with
$\big({\rm  Coker}(\varphi)\big)_S$\,.
\end{enumerate}
%
\item   If we have an exact sequence of \Amos
$$\preskip.2em \postskip.3em
M\vers{\varphi}N\vers{\psi}P\;,
$$
then the sequence of $\gA_S$-modules
$$\preskip.2em \postskip.0em
M_S\vers{\varphi_S}N_S\vers{\psi_S}P_S
$$
is \egmt exact.
\end{enumerate}
\end{fact}

\begin{fact}\label{fact.LocIntersect}
If $M_1$, $\ldots$, $M_r$ are submodules of $N$ and $M=\bigcap_{i=1}^rM_i$,
then by identifying the modules $(M_i)_S$ and $M_S$ with submodules of $N_S$ we obtain $M_S=\bigcap_{i=1}^r(M_i)_S$.
\end{fact}

\vspace{.3em} 
\begin{fact}
\label{fact.transporteur}
Let $M$ and $N$ be two submodules of an \Amo $P$, \emph{with $N$ \tfz}.
Then, the conductor \id $(M_S:N_S)$ is identified with $(M:N)_S$, via the natural maps of $(M:N)$ in  $(M_S:N_S)$ and $(M:N)_S$.
\end{fact}

This is particularly applied to the annihilator of a \itfz.

\subsec{Local-global principle for exact sequences of modules}

\begin{plcc}
\label{plcc.basic.modules}\emph{(For exact sequences)}\\
Let $S_1$, $\ldots$, $ S_n$  be \moco of $\gA$, $M$, $N$, $P$ be \Amos and $\varphi:M\to N$, $\psi:N\to P$ be two \alisz. We write $\gA_i$ \hbox{for $\gA_{S_i}$},
$M_i$ for $M_{S_i}$ etc. \Propeq
\begin{enumerate}\itemsep0pt
\item  The sequence \smash{$M\vers{\varphi}N\vers{\psi}P$} is exact.
\item  For each $ i\in\lrbn,$ the sequence $M_i\vers{\varphi_i}N_i\vers{\psi_i}P_i$ is exact.
\end{enumerate}
As a consequence, $\varphi$ is injective (resp.\,surjective) \ssi for each $ i\in\lrbn,$  $\varphi_i$ is injective (resp.\,surjective)
\end{plcc}
\begin{proof}
We have seen that \emph{1} $\Rightarrow$ \emph{2} in Fact~\ref{fact.sexloc}.\\
Assume \emph{2}. Let $\mu_i:M\to M_i,$  $\nu_i:N\to N_i,$  $\pi_i:P\to
P_i$ be the canonical \homosz. Let $x\in M$ and $z=\psi\big(\varphi(x)\big)$. We thus have
$$\preskip.4em \postskip.4em 
0=\psi_i\big(\varphi_i(\mu_i(x))\big)=\pi_i\big(\psi(\varphi(x))\big)=\pi_i(z), 
$$
for some $s_i\in S_i$, $s_iz=0$ in $P$. We conclude that $z=0$ by using the comaximality of the $S_i$'s:
$\som_i u_is_i=1$. Now let $y\in N$ such \hbox{that $\psi(y)=0$}. For each $i$ there exists some $x_i\in M_i$ such that $\varphi_i(x_i)=\nu_i(y)$. \\
We write $x_i=_{M_i}a_i/s_i$ with $a_i\in M$ and $s_i\in S_i$. The \egt $\varphi_i(x_i)=\nu_i(y)$ means that for some $t_i\in S_i$ we have $t_i\varphi(a_i)=t_is_iy$ in~$N$. If $\som_i v_it_is_i=1$,
we can deduce that $\varphi\big(\som_iv_it_ia_i\big)=y$. Thus $\Ker\psi$ is indeed included in $\Im\varphi$.
\end{proof}

\begin{plca}
\label{plca.basic.modules}\emph{(For exact sequences)}\\
Let $M$, $N$, $P$ be \Amosz, and $\varphi:M\to N$ \hbox{and $\psi:N\to P$} be two
\alisz.  \Propeq
\begin{enumerate}\itemsep0pt
\item  The sequence \smash{$M\vers{\varphi}N\vers{\psi}P$} is exact.
\item  For every \idema $\fm$ the sequence $M_\fm\vers{\varphi_\fm}N_\fm\vers{\psi_\fm}P_\fm$ is exact.
\end{enumerate}
As a consequence, $\varphi$ is injective (resp.\,surjective) \ssi for every \idema $\fm$,  $\varphi_\fm$ is injective (resp.\,surjective).
\end{plca}
%
\begin{proof}
The \prt $x=0$ for an \elt $x$ of a module is a \carf \prtz.
Similarly for the \prt $y\in\Im\varphi$. Thus, even if the \prt \gui{the sequence is exact} is not of \carfz, it is a conjunction of  \carf \prtsz, and we can apply Fact\etoz~\ref{fact2PropCarFin} to deduce the abstract \plg from the concrete \plgz.
\end{proof}

Let us finally mention a \plgc for \mosz.

\pagebreak

\begin{plcc}
\label{plcc.basic.monoides}\emph{(For \mosz)}\\
Let $S_1$, $\ldots$, $ S_n$ be \moco of $\gA$, $V$ be a \moz. \Propeq
\begin{enumerate}\itemsep0pt
\item  The \mo $V$ contains $0$.
\item  For $ i\in\lrbn,$ the \mo $V$ seen in $\gA_{S_i}$ contains $0$.
\end{enumerate}
\end{plcc}
%
\begin{proof}
For each $i$ we have some $v_i\in V$ and some $s_i\in S_i$ such that $s_iv_i=0$.
Let~$v=\prod_i v_i\in V$. Then, $v$ is zero in the $\gA_{S_i}$'s, thus in $\gA$.
\end{proof}
%

\Exercices{

\begin{exercise}
\label{exo2Lecteur}
{\rm  We recommend the reader to do the \dems which are not given, are sketched, are
left to the reader,
etc\ldots
\, In particular, consider the following cases.
\begin{itemize}
\item \label{exofactKerAAsMMs} Check Facts \ref{factUnivLoc} to \ref{fact.bilocal}. 
\item  \label{exocorplcc.basic}
 Prove Corollary \ref{corplcc.basic}.

\item  \label{exolemGaussJoyal}
In Lemma \ref{lemGaussJoyal} compute suitable exponents for the items \emph{\ref{i2lemPrimitf}},
  \emph{\ref{i3lemPrimitf}} and \emph{\ref{i4lemPrimitf}}, 
  by making the \dem completely explicit.
\item  \label{exocorpropCoh1}
  Prove Corollary \ref{corpropCoh1}.
\label{exopropCoh4}
  Give a more detailed \dem of \Thref{propCoh4}.
\label{exoplcc.coh}
  Check the details in the \dem of the \plgrf{plcc.coh}.
\label{exopropCohfd1}
  Prove Proposition \ref{propCohfd1}.
\item \label{exofact.transporteur} 
Check Facts \ref{fact.sexloc} to \ref{fact.transporteur}. For Fact \ref{fact.LocIntersect} we use the exact sequence $0\to M\to N\to \bigoplus_{i=1}^rN/M_i$ which is preserved by \lonz.

\end{itemize}
}
\end{exercise}

\vspace{-1em}
\begin{exercise}\label{exoNilpotentChap2} 
{\rm (Also see exercise \ref{exoNilIndexInversiblePol})
\begin{enumerate}\itemsep0pt
\item  \emph{(Invertible elements in $\gB[T]$, cf. Lemma \ref{lemGaussJoyal})}\\
Let two \pols $f= \sum_{i=0}^n a_i T^i$, $g= \sum_{j=0}^m b_j T^j$
 with  $fg = 1$.
Show that the \coes $a_i$, $i \ge 1$, $b_j$, $j \ge 1$ are nilpotent elements and that $a_n^{m+1} = 0$.

\item  \emph{(\Cara \pol of a nilpotent matrix)}\\
Let $A \in \Mn(\gB)$ be a nilpotent matrix and 
$\rC{A}(T) = T^{n} + \sum_{k=0}^{n-1}a_kT^k
$ be its \polcarz.
\begin{enumerate}
\item Show that the \coes $a_i$ are nilpotent elements.
\item Precisely, if $A^e = 0$, then  $\Tr(A)^{(e-1)n + 1} = 0$ and

\snic{
 a_i^{e_i} = 0  \quad \hbox {where} \quad
e_i = (e-1) {n \choose i} + 1\quad(i=0,\ldots,n-1).}
\end{enumerate}
\end{enumerate}
}
\end{exercise}

\vspace{-1em}
\begin{exercise}
 \label{exolemUMD}
 {\rm  Let $x=(\xn)\in\Ae{n}$ be a vector and $s\in\gA$.
\begin{enumerate}\itemsep0pt
\item If $x$ is \umd in $\aqo{\gA}{s}$ and in $\gA[1/s]$, it is \umd in $\gA$.
\item Let $\fb$ and $\fc$ be two \ids of $\gA$. If $x$ is \umd modulo
$\fb$ and modulo $\fc$, then it is also \umd modulo $\fb\fc$.
\end{enumerate}
} \end{exercise}

\vspace{-1em}
\pagebreak

\begin{exercise} 
\label{exoPlgb1} 
(A typical application of the basic \plgz)\\
{\rm   Let
$x=(\xn)\in\Ae{n}$ be \emph{\umdz}.
For $d\geq 1$, we denote by $\AXn_d$ the \Asub of the \pogs of degree $d$ and 

\snic{I_{d,x}=\sotq{f\in\AuX_d}{f(x)=0}, \hbox{ \Asub of  }\AuX.\label{NOTAAXd}}
\begin{enumerate}\itemsep0pt
\item If  $x_1\in\Ati$,  every $f\in I_{d,x}$ is a \coli \hbox{of the $x_1X_j-x_jX_1$} with \pogs
of degree~$d-1$ for \coesz. 
\item \Gnltz, every $f\in I_{d,x}$ is a \coli of the $(x_kX_j-x_jX_k)$ with \pogs of degree~$d-1$ for \coesz.
\item Let $I_x=\bigoplus_{d\geq 1}I_{d,x}$. Show that $I_x=\sotq{F}{F(tx)=0}$ (where $t$ is a new \idtrz). Show that $I_x$ is \emph{saturated}, i.e., if $X_j^mF\in I_x$ for some $m$ and for each $j$, then $F\in I_x$.
\end{enumerate}
}
\end{exercise}

\vspace{-1em}
\begin{exercise}
\label{exoGaussJoyal} (Variations of the \iJG Gauss-Joyal Lemma \ref{lemGaussJoyal})
\\
{\rm Show that the following statements are \eqves
(each statement is \uvlez, i.e., valid for all \pols and every commutative \ri $\gA$):
\begin{enumerate}\itemsep0pt
\item  $\rc(f)=\rc(g)=\gen{1}\; \Rightarrow\;  \rc(fg)=\gen{1}$,
\item  $(\exists i_0,j_0\;\;f_{i_0}=g_{j_0}=1)\;\Rightarrow\;\rc(fg)=\gen{1}$,
\item  $\exists p\in\NN, \; \; \big(\rc(f)\rc(g)\big)^p\subseteq \rc(fg)$,
\item  \emph{(Gauss-Joyal)} $\; \DA\big({\rc(f)\rc(g)}\big) = \DA\big({\rc(fg)}\big)$.
\end{enumerate}
}
\end{exercise}

\vspace{-1em}
\begin{exercise}\label{exoNormPrimitivePol}
{(Norm of a primitive \pol through the use of a null \riz)}\\
{\rm  
Let $\gB$ be a free \Alg of finite rank, $\uX = (\Xn)$ be \idtrsz, $Q \in \gB[\uX]$ and $P = \rN_{\gB[\uX]/\gA[\uX]}(Q) \in \gA[\uX]$.  
Show that if $Q$ is primitive, then so is $P$.  
\emph{Hint:} check that $\gA\cap \rc_\gB(P) = \rc_\gA(P)$, 
consider the sub\ri $\gA' = \gA/\rc_\gA(P)$ of $\gB' = \gB/\rc_\gB(P)$ and the $\gA'$-\lin map \gui{multiplication by $Q$,}
$m_Q : \gB'[\uX] \to \gB'[\uX],\;R\mapsto QR$.

}

\end {exercise}

\vspace{-1em}
\begin{exercise}
\label{exoCohfd1}
{\rm  Show that a \coh \ri $\gA$
is \fdi \ssi the test
\gui{$1\in\gen{a_1,\ldots ,a_n}$?}
is explicit for every finite sequence $(a_1,\ldots ,a_n)$ in $\gA$.
}
\end{exercise}

\vspace{-1em}
\begin{exercise}
\label{exo.quo.coh} (An example of a \coh \noe \ri with a
\emph{non}-\coh quotient.) \\
{\rm
Consider the \ri $\ZZ$ and an \id $\fa$ generated by an infinite sequence of \eltsz, all zeros besides eventually one, which is then equal to $3$
(for example we place a~$3$ the first time, if it ever occurs, that a zero of the Riemann zeta function\footnote{Here we enumerate the zeros $a_n+ib_n$ with $b_n>0$ by order of magnitude.} has  real part not equal to $1/2$).
If we are able to provide a finite \sys of \gtrs for the annihilator of $3$ in $\ZZ/\fa$, we are able to say whether the infinite sequence is identically zero or not.
This would mean that there exists a sure method to solve conjectures of the Riemann type.

\comm
As every reasonable \cov \dfn of  \noet seems to demand that a \noe \riz's quotient remains \noez, and given the above \gui{counterexample,}
we cannot hope to have a \cov \dem of the \tho
of \clama which states that every \noe \ri is \cohz.
\eoe
}
\end{exercise}

\vspace{-1em}
\begin{exercise}\label{exoIdempAX} (Idempotents of $\AX$)\\
{\rm  
Prove that every \idm of $\AX$ is an \idm of $\gA$.
}
\end {exercise}

\vspace{-1em}
\begin{exercise}
\label{exoIdmsSupInf}
{\rm  Let $u$ and $v$ be two \idms and $x$ be an \elt of $\gA$. \\
The \elt $1-(1-u)(1-v)=u+v-uv$ is denoted by $u\vu v$.
\begin{enumerate}\itemsep0pt
\item Show that $x\in u\gA\,\Leftrightarrow\,ux=x$.
In particular, $u\gA=v\gA \,\Leftrightarrow\, u=v$.
\item  The \elt $uv$
is \und{the} least common multiple of $u$ and $v$ amongst the \idms of $\gA$
(i.e., if $w$ is an \idmz, $w\in u\gA\cap v\gA \,\Leftrightarrow\,
w\in uv\gA$). Actually, we even have $u\gA\cap v\gA=uv\gA$. We write $u\vi v=uv$.
\item Prove the \egt $u\gA+ v\gA=(u\vu v)\gA$.
Infer that $u\vu v$ is \und{the} greatest common divisor of $u$ and $v$ amongst the \idms of $\gA$
(in fact an arbitrary \elt of $\gA$ divides $u$ and $v$ \ssi it divides $u\vu v$).
\item By a sequence of \mlrsz, transform the matrix $\Diag(u,v)$ into the matrix $\Diag(u\vu v,u\vi v)$.\\
 From it, deduce that the two \Amos $u\gA\oplus v\gA$ and ${(u\vu v)}\gA\oplus {(u\vi v)}\gA$
are isomorphic.
\item Show that the two \ris $\aqo{\gA}{u}\times \aqo{\gA}{v}$ and $\aqo{\gA}{u\vu v}\times \aqo{\gA}{u\vi v}$ are isomorphic.
\end{enumerate}
}
\end{exercise}

\vspace{-1em}
\begin{exercise}
\label{exoFracSfio}
{\rm Let $\gA$ be a \ri and $(e_1,\ldots,e_n)$ be a \sfio of $\Frac\gA=\gK$.
We write $e_i=a_i/d$ with $a_i\in\gA$ and $d\in\Reg\gA$.
We then have $a_ia_j=0$ for $i\neq j$ and $\sum_ia_i$ \ndzz.
\\
\emph{1.} 
Establish a converse.  
\\
\emph{2.} Show that $\gK[1/e_i]\simeq \Frac\big(\gA\sur{\Ann_\gA(a_i)}\!\big)$ and $\gK\simeq\prod_i\Frac\big(\gA\sur{\Ann_\gA(a_i)}\!\big)$. 
}
\end{exercise}

\vspace{-1em}
\begin{exercise}
\label{exoFracZed} (Separating the \irds components)
\\
{\rm 
\emph{1.} Let $\gA=\QQ[x,y,z]=\aqo{\QQ[X,Y,Z]}{XY,XZ,YZ}$ and $\gK=\Frac\gA$.
What are the zeros of $\gA$ in $\QQ^3$ (i.e.  $(x,y,z)\in\QQ^3$
such that $xy=yz=zx=0$)?
Give a reduced form of the \elts of $\gA$. Show that $x+y+z\in\Reg\gA$. Show that the \elts $\frac{x}{x+y+z}$, $\frac{y}{x+y+z}$ and $\frac{z}{x+y+z}$ form a \sfio in $\gK$. 
Show that $\gK\simeq\QQ (X) \times\QQ (Y) \times\QQ (Z) $.

\emph{2.} Let $\gB=\QQ[u,v,w]=\aqo{\QQ[U,V,W]}{UVW}$ and $\gL=\Frac\gB$.\\
What are the zeros of $\gB$ in $\QQ^3$?
Give a reduced form of the \elts of $\gB$. 
Show that $\gL\simeq\QQ(U,V)\times\QQ(V,W)\times\QQ(W,U)$.  
 
}
\end{exercise}

\vspace{-1em}
\begin {exercise} \label {exoIdempotentE2}
       (Idempotent and \elr group)\\
{\rm
Let $a \in \gA$ be an \idmz. For $b \in \gA$, 
give
a matrix $A \in \EE_2(\gA)$  and an \elt $d\in\gA$ such that $A \cmatrix {a\cr b\cr} = \cmatrix {d\cr 0\cr}$. In particular, explain why $\gen {a,b} = \gen {d}$.\\
Moreover, 
prove that if $b$ is \ndz (resp.\,\ivz) modulo $a$, then $d$ is \ndz 
(resp.\,\ivz). Finally, if $b$ is \idmz, $d=a\vu b=a+b-ab$.
}

\end {exercise}

\vspace{-1em}
\pagebreak

\begin{exercise}\label{exoSfio} 
{\rm  Let $(r_1,\ldots,r_m)$ be a finite family of \idms in a \ri $\gA$. Let $s_i=1-r_i$ and, for a subset $I$ of $\lrbm$, let $r_I=\prod_{i\in I}r_i\prod_{i\notin I}s_i$.

\noindent \emph{1.} 
Show that the diagonal matrix $D=\Diag(r_1,\ldots,r_m)$ is similar to a matrix  $D'=\Diag(e_1,\ldots,e_m)$ where the $e_i$'s are \idms which satisfy: $e_i$ divides $e_j$ if $j> i$. You can start with the $n=2$ case and use Exercise \ref{exoIdmsSupInf}. Show that  $\gen{e_k}=\cD_k(D)$ for all $k$.

\noindent \emph{2.} 
Show that we can write $D'=PDP^{-1}$ with $P$ a \emph{generalized permutation matrix}, \cad a matrix which can be written as $\sum_jf_jP_j$ where the $f_j$'s form a \sfio and each $P_j$ is a permutation matrix.
\index{matrix}{generalized permutation ---}
Suggestions:
\begin{itemize}\itemsep0pt
\item  The $r_I$'s form a \sfioz.
The diagonal matrix $r_ID$ has the \elt $r_I$ as its \coe in position $(i,i)$ if $i\in I$ and $0$ otherwise.
The matrix $P_I$ then corresponds to a permutation bringing the \coes $r_I$ to the head of the list.
Finally, $P=\sum_Ir_IP_I$. Note that the test \gui{$r_I=0$?} is not \ncrz!
\item We can also treat the $m=2$ case: find $P=e
\cmatrix{1&0\cr0&1}+ f\cmatrix{0&1\cr1&0}$ with $f=r_2s_1$, $e=1-f$,
and $D'=\Diag(r_1\vu r_2,r_1\vi r_2)$.
\\
Next we treat the $m>2$ case step by step.
\end{itemize}

}
\end{exercise}

\vspace{-1em}
\begin{exercise}
\label{exoChinois}
{\rm  Recall the \dem of the Chinese Remainder \Tho (\paref{restes
chinois}) and explicitly give the \idmsz.}
\end{exercise}

\vspace{-1em}
\begin{exercise} \label {exoFacileGrpElem1}
       (\Elr Group: first steps) {\rm $\MM_2(\gA)$ case.\\
 \emph{1.} 
Let $a \in \gA$. 
Determine a matrix $P \in \EE_2(\gA)$
such that $P\cmatrix {a\cr 0\cr} = \cmatrix {0\cr
a\cr}$. Same for \smashtop{$\cmatrix {\varepsilon a\cr 0\cr}
\mapsto \cmatrix {a\cr 0\cr}$} where $\varepsilon \in \Ati$.

\noindent \emph{2.} 
Write the matrices \halfsmashtop{$\crmatrix {0 & -1\cr 1 & 0\cr}$ and $\crmatrix {-1 & 0\cr 0 & -1\cr}$} as \elts of $\EE_2(\gA)$.

\noindent \emph{3.} 
Show that every triangular matrix of $\SL_2(\gA)$ is in $\EE_2(\gA)$.

\noindent \emph{4.} 
Let $u = \cmatrix {x\cr y}$, $v = \cmatrix {y\cr x}$, $w = \crmatrix {-y\cr
x}$ with $x,y\in\gA$.  
Show that \hbox{$v\in \GL_2(\gA)\cdot u$} and $w\in \EE_2(\gA)\cdot u$, but not \ncrt
$v \in \SL_2(\gA)\cdot u$. For example, if $x$, $y$ are two \idtrs over a \ri~$\gk$, $\gA=\gk[x,y]$ and $v = Au$, with $A \in\GL_2(\gA)$, then
$\big(\det(A)\big)(0,0)=-1$. 
Consequently, we have \hbox{$\det (A) \in -1 +\rD_\gk(0) \gen{x,y}$} (Lemma~\ref{lemGaussJoyal}), therefore $\det(A) = -1$ if $\gk$ is reduced. In addition, if $\det(A) =1$, then $2=0$ in~$\gk$.
As a result, $v \in \SL_2(\gA)\cdot u$ \ssi $2=0$ in~$\gk$.
}
\end {exercise}



\vspace{-1em}
\begin{exercise} \label {exoFacileGrpElem2}
       (\Elr group: next steps)\\
{\rm
\emph{1.} Let $A \in \MM_{n,m}(\gA)$  with an \iv \coe and $(n,m)\neq(1,1)$. Determine matrices $P \in \EE_n(\gA)$ \hbox{and $Q \in\EE_m(\gA)$} such that
\smashbot{$P A\, Q = \cmatrix {1 & 0_{1,m-1}\cr 0_{n-1,1}& A'}$}.
Example: with $a\in\Ati$ give $P$ for $P\,\cmatrix {a \cr 0\cr} = \cmatrix {1 \cr 0\cr}$
 (Exercise~\ref{exoFacileGrpElem1} item~\emph{1}).

%
\noindent \emph{2.} 
Let $A \in \MM_2(\gA)$ with an \iv \coez. 
Compute matrices $P$ \hbox{and $Q \in \EE_2(\gA)$} such that:
$P A\, Q = \cmatrix {1 & 0\cr 0 & \delta\cr}$ with $\delta = \det(A)$.\\
Every matrix $A \in \SL_2(\gA)$ with an \iv \coe belongs to $EE_2(\gA)$. Make the following cases explicit:
$$
\cmatrix {a & 0\cr 0 & a^{-1}\cr},\qquad
\cmatrix {0 & a\cr -a^{-1} & 0\cr},\qquad
\hbox {with $a \in \Ati$.}
$$
Write the following matrices (with $a \in
\Ati$) in $\EE_2(\gA)$:
$$
\cmatrix {a & b\cr 0 & a^{-1}\cr},\qquad
\cmatrix {a & 0\cr b & a^{-1}\cr},\qquad
\cmatrix {0 & a\cr -a^{-1} & b\cr},\qquad
\cmatrix {b & a\cr -a^{-1} & 0\cr}.
$$

\noindent \emph{3.} 
Prove that if $A = \Diag(a_1, a_2, \ldots, a_n) \in \SLn(\gA)$,
then $A\in\EE_n(\gA)$.

\noindent \emph{4.} 
Show that every triangular matrix $A \in \SL_n(\gA)$ belongs to $\EE_n(\gA)$.
}
\end {exercise}


\vspace{-1em}
\begin{exercise} \label {exoFacileGrpElem4}
                 (Division matrices $D_q$ of \deter $1$)\\
{\rm
\noindent A \gui {\gnl division} $a = bq - r$ can be expressed with matrices:
$$
\crmatrix {0 & 1\cr -1& q\cr} \cmatrix {a\cr b\cr} = \cmatrix {b\cr r\cr}.
$$
This leads to the introduction of the matrix $D_q = \crmatrix {0 & 1\cr -1& q\cr} \in
\SL_2(\gA)$.

\noindent Show that $\EE_2(\gA)$ is the \mo generated by the $D_q$ matrices.

}
\end {exercise}


\vspace{-1em}
\begin{exercise} \label {exoFacileGrpElem5} 
%
{\rm
 Let $\gA$ be a \ri and $A$, $B \in \MMn(\gA)$. Assume that we have some~$i\in\gA$ with $i^2 = -1$ and that $2\in\Ati$. Show that the matrices of $\MM_{2n}(\gA)$

\snic{M = \cmatrix {A & -B\cr B & A\cr}\; $ and $\;
M' = \cmatrix {A+iB & 0\cr 0 & A-iB\cr}}

are \emph{\elrt similar}, (i.e., $\Ex P \in \EE_{2n}(\gA),\;P M P^{-1} = M'$).\\
\emph{Hint}: first treat the $n=1$ case
.
}

\end {exercise}


\vspace{-1em}
\begin{exercise} \label {exoFacileGrpElem6}
{\rm
\noindent
For $d \in \Ati$ and $\lambda \in \gA$ compute the matrix

\snic{
\Diag(1, \dots, d, \dots, 1)\cdot \rE_{ij}(\lambda)\cdot
\Diag(1, \dots, d^{-1}, \dots, 1).}


Show that the subgroup of diagonal matrices of $\GL_n(\gA)$ normalizes~$\En(\gA)$.
}
\end {exercise}

\vspace{-1em}
\begin{exercise}\label{exoTroisiemeLemmeLiberte} (A freeness lemma, or a Splitting Off: reader's choice) 
 \\
{\rm
Let $F \in \GAn(\gA)$ be a \prr with an \iv principal minor of order $k$. Show that $F$ is similar to a matrix $\cmatrix {\I_k & 0\cr 0 & F'\cr}$ where $F' \in \GA_{n-k}(\gA)$.

\noindent
The \mptf $P\eqdefi\Im F\subseteq \Ae n$ admits a free direct summand with $k$ columns of $F$ for its basis.
}
\end {exercise}


\vspace{-1em}
\begin{exercise}\label{exoABArang1}
{\rm
Let $A \in \Ae {n \times m}$ be of rank $1$.
 Construct $B \in \Ae {m \times n}$ such that $ABA = A$ and verify that $AB$ is a \prr of rank $1$. Compare your solution with that which would result from the \dem of \Thref{propIGCram}.
}

\end {exercise}

\vspace{-1em}
\begin{exercise}
\label{exoABAabstrait}
{\rm
 \emph{This exercise constitutes an abstraction of the computations that led to  \Thref{propIGCram}.} 
Consider an \Amo $E$ \gui {with enough \lins forms},
i.e.  if $x \in E$ satisfies $\mu(x) = 0$ for all $\mu\in E\sta$, then $x = 0$. This means that the canonical map from $E$ to its bidual, $E \to E{\sta}{\sta}$, is injective. This condition is satisfied if $E$ is a \emph{reflexive} module,
 i.e. $E \simeq E{\sta}{\sta}$, e.g. a \mptfz%
, or a free module of finite rank%
.

\noindent For $x_1$, $\ldots$, $x_n \in E$, denote by $\Vi_r(x_1, \ldots, x_n)$
the \id of $\gA$ generated by the evaluations of every $r$-multi\lin  alternating form of~$E$ at every $r$-tuplet of \elts of $\{x_1, \ldots, x_n\}$.

\noindent \emph{Assume that $1 \in \Vi_r(x_1, \ldots, x_n)$ and $\Vi_{r+1}(x_1, \ldots, x_n) = 0$.}

\noindent We want to prove that the submodule $\som \gA x_i$ is a direct summand in~$E$ 
by explicitly giving a \prr $\pi : E \to E$ whose image is this submodule.
\begin{enumerate}\itemsep0pt
\item \emph{(Cramer's formulas)}
Let $f$ be an $r$-multi\lin alternating form over~$E$. Show, for $y_0$, $\ldots$, $y_r \in \som \gA x_i$, that
$$
\som_{i=0}^r
(-1)^i f(y_0, \ldots, y_{i-1}, \widehat {y_i}, y_{i+1}, \ldots, y_r)\, y_i = 0
.$$
Or, for $y$, $y_1$, $\ldots$, $y_r \in \som \gA x_i$, that
$$
f(y_1, \ldots, y_r)\,y =
\som_{i=1}^r f(y_1, \ldots, y_{i-1}, y, y_{i+1}, \ldots, y_r)\, y_i
.$$

\item
Give $n$ \lins forms $\alpha_i \in E\sta$ such that
the \ali 

\snic{\pi : E \to E$, $\; x\mapsto \som_i \alpha_i(x) x_i}

  is a \prr 
onto
$\som \gA x_i$. 
 We define $\psi : \Ae{n} \to E$  by $e_i \mapsto x_i$ and $\varphi : E \to \Ae{n}$  by $\varphi(x) = \big(\alpha_1(x), \ldots, \alpha_n(x)\big)$.
Arrange for $\pi = \psi \circ \varphi$ and $\pi \circ \psi = \psi$,
so that $\psi \circ \varphi \circ \psi = \psi$.

\item \emph{(New \dem of \Thref{propIGCram})}
Let $A \in \Ae {m \times n}$ be a matrix of rank~$r$.
Show that there exists a $B \in \Ae {n \times m}$ such that~$A\,B\,A = A$.
\end{enumerate}
}
\end{exercise}

\vspace{-1em}
\begin{exercise}
\label{exoPrepBinetCauchy}
{\rm  
Let $A\in \Ae {n\times m}$ and $B\in\Ae {m\times n}$.

\noindent  \emph{1.} We have the following commutativity formula:
$$\preskip.1em \postskip.1em\det(\I_m+XBA)=\det(\I_n+XAB).$$

\emph{First \demz}. First treat the case where $m=n$, for example by the method of undetermined \coesz.
If $m\neq n$, $A$ and $B$ can be completed with rows and columns of $0$'s to turn them into square matrices $A_1$ and $B_1$ of size $q=\max(m,n)$ as in the \dem given \paref{eqIDC1}. Then check that $\det(\I_m+XBA)=\det(\I_q+XB_1A_1)$ and $\det(\I_n+XAB)=\det(\I_q+XA_1B_1)$.

\emph{Second \demz}. Consider an undetermined $X$ and the matrices

\snic{B'=\cmatrix{XB&\I_m\cr\I_n&0_{n,m}}\quad \mathrm{and}\quad
A'=\cmatrix{A&\I_n\cr\I_m&-XB}.}

Compute $A'B'$ and $B'A'$.

\noindent  \emph{2.} What can be deduced about the \polcars of $A\,B$ and $B\,A$?
}
\end{exercise}

\vspace{-1em}
\begin{exercise}\label{exoBinetCauchy} (Binet-Cauchy formula)
\index{Binet-Cauchy!formula}\\
 {\rm  We use the notations on \paref{notaAdjalbe}.
    For two matrices $A\in \Ae {n\times m}$ and $B\in\Ae {m\times n}$,
    prove that we have the Binet-Cauchy formula:
$$
\preskip.4em \postskip.0em
\det(BA)=\som_{\alpha \in \cP_{m,n}
}\det(B_{1..m,\alpha})\det(A_{\alpha,1..m}).
$$

\noindent  \emph{First \demz}.  Use the formula $\det(\I_m+XBA)=\det(\I_n+XAB)$  (Exercise \ref{exoPrepBinetCauchy}).
Then consider the \coe of $X^m$ in each of the \pols
$\det(\I_m+XBA)$ and $\det(\I_n+XAB)$.

\noindent  \emph{Second \demz}. The matrices $A$ and $B$ represent \alis $u :\Ae{m}\to\Ae{n}$ and $v :\Ae{n}\to\Ae{m}$. \\
Then consider the matrices of $\Al mu $, $\Al mv $ and $\Al m(v\circ u )$ with respect to the bases naturally associated with the canonical bases \hbox{of $\Ae{n}$} and~$\Ae{m}$.
\\
Conclude by writing $\Al m(v\circ u )=\Al m v \circ\Al  mu $.

\noindent  \emph{Third \demz}. In the product $BA$
insert between $B$ and $A$ a diagonal matrix $D$ having \idtrs $\lambda_i$ for \coesz, and see which is the \coe of $\lambda_{i_1}\cdots\lambda_{i_m}$
in the \pol $\det(BDA)$ (to do this take $\lambda_{i_1}=\cdots=\lambda_{i_m}=1$ and let the other be null). Conclude by letting all the $\lambda_i$'s be equal to $1$.
}
\end{exercise}

\vspace{-1em}
\begin{exercise}
\label {exoDetExtPower}
\noindent {\rm
Let $u \in \End_\gA(\Ae{n})$. For $k \in \lrb{0..n}$, let $u_k =
\Al k(u)$. \\
Show that $\det(u_k) =
\det(u)^{n-1 \choose k-1}$ and that
$$\preskip.4em \postskip-1.1em 
\det(u_k)
\det(u_{n-k}) = \det(u)^{n \choose k}\ . 
$$
}
\end{exercise}

\vspace{-1em}
\begin{exercise}
 \label{exoMatInjLocSimple}
 {\rm  For $A\in \Ae {n\times r}$ prove that \propeq
 \begin{enumerate}
\item The matrix $A$ is injective and \lnlz.
\item There exists a matrix $B\in \Ae {r\times n}$ such that $B\,A=\I_r$.
\item The \idd $\cD_r(A)=\gen{1}$.
\end{enumerate}
\emph{Hint}: See \Thos \vref{propIGCram}, \vref{prop inj surj det} and \vref{theoremIFD}.
 } \end{exercise}

\vspace{-1em}
\begin{exercise}
\label{exolemdeterblocs}
{\rm  Treat the \gnl case in the \dem of Lemma \ref{lemdeterblocs}.
}
\end{exercise}

\vspace{-1em}
\begin{exercise}
\label{exoGram} 
{\rm If $\gram_\gA(\varphi,\xn)$ is \ivz, the submodule $\gA x_1+\cdots+\gA x_n$ is free with $(\xn)$ as its basis.
}
\end{exercise}


\vspace{-1em}
\begin{problem}\label{exoRationaliteLineaire} {(Gauss' pivot, $A\,B\,A=A$, and \lin rationality)}\\
{\rm  
Let $\gK$ be a \cdiz. If $x \in \gK^n$ is a nonzero vector, its \emph {pivot index} $i$ is the least index $i$ such that $x_i \ne 0$. We say that the \coe $x_i$ is the \emph {pivot} of~$x$.  The \emph{height} $h(x)$ of $x$ is the integer $n-i+1$ and it is agreed that $h(0) = 0$. For example, for $n = 4$ and $x = \cmatrix {0\cr 1\cr *\cr *\cr}$, the pivot index of $x$ is $i = 2$, and $h(x) = 3$. 
The following notions of \gui{staggering} are relative to this height $h$.  \\
We say that a matrix $A \in \MM_{n,m}(\gK)$  \emph {has staggered columns} if the nonzero 
columns of $A$ have distinct heights; we say that it is \emph{strictly staggered} 
if, additionally, the rows at the pivot indices are vectors of the canonical basis of $\gK^m$ (these vectors are \ncrt distinct). Here is a strictly staggered matrix
($0$ has been replaced by a dot):
$$\preskip-.2em \postskip.4em 
\cmatrix {
\cdot      &\cdot      &\cdot      &1      &\cdot      &\cdot\cr
\cdot      &\cdot      &\cdot      &a_{24} &\cdot      &\cdot\cr
\cdot      &\cdot      &1      &\cdot      &\cdot      &\cdot\cr
\cdot      &\cdot      &a_{43} &a_{44} &\cdot      &\cdot\cr
1      &\cdot      &\cdot      &\cdot      &\cdot      &\cdot\cr
\cdot      &1      &\cdot      &\cdot      &\cdot      &\cdot\cr
a_{71} &a_{72} &a_{73} &a_{74} &\cdot      &\cdot\cr
\cdot      &\cdot      &\cdot      &\cdot      &1      &\cdot\cr
a_{91} &a_{92} &a_{93} &a_{94} &a_{95} &\cdot\cr
}.
$$
\emph {1.}
Let $A \in \MM_{n,m}(\gK)$ be strictly staggered; 
we define $\ov {A} \in \MM_{n,m}(\gK)$ by annihilating the nonpivot \coes (the $a_{ij}$'s in the above exercise) and $B = \tra {\,\ov A} \in \MM_{m,n}(\gK)$. Check that $ABA = A$. \\
Describe the \prrs $AB$, $BA$ and the \dcn $\gK^n = \Im AB \oplus \Ker AB$.

\emph {2.}
Let $A \in \MM_{n,m}(\gK)$ be an arbitrary matrix. How do you obtain $Q \in
\GL_m(\gK)$ such that $A' = AQ$ is strictly staggered?
How do you compute $B \in \MM_{m,n}(\gK)$ satisfying $ABA = A$?

\emph {3.}
Let $A \in \MM_{n,m}(\gK)$ and $y \in \gK^n$. Assume that the \sli $Ax = y$ admits a solution $x$ on an over\ri of $\gK$. Show that it admits a solution on $\gK$.

\emph {4.}
Let $\gK_0 \subseteq \gK$ be a subfield and $E$, $F$ be 
two \supl $\gK$-linear subspaces 
of $\gK^n$. Assume that $E$ and $F$ are generated by vectors with components in $\gK_0$. Show that $\gK_0^n = (E\cap \gK_0^n) \oplus (F\cap \gK_0^n)$.

Let $E \subseteq \gK^n$ be 
a $\gK$-linear subspace. 
We say that \emph{$E$ is $\gK_0$-rational} if it is generated by vectors with components in $\gK_0$.

\emph {5.}
Let $F$ be a complementary subspace  of $E$ in $\gK^n$   
generated by vectors of the canonical basis of $\gK^n$: $\gK^n = E \oplus F$ and $\pi : \gK^n \twoheadrightarrow E$ be the associated \prnz.
\begin {itemize}
\item [\emph {a.}]
Show that $E$ is $\gK_0$-rational \ssi $\pi(e_j) \in \gK_0^n$ for every vector $e_j$ of the canonical basis.

\item [\emph {b.}]
Deduce the existence of a smaller field of rationality for $E$.

\item [\emph {c.}]
What is the field of rationality of the image in $\gK^n$ of a strictly staggered matrix?
\end {itemize}

}

\end {problem}

\vspace{-1em}
\begin{problem}
\label{exoPlgb2} ~\\
{\rm  
\emph{1.} \emph{Partial \fcn algorithm.}
 Given two integers $a$ and $b$ prove that we can \gui{efficiently} compute a finite family of pairwise coprime positive integers $p_i$ such that $a=\pm\prod_{i=1}^np_i^{\alpha_i}$ and $b=\pm\prod_{i=1}^np_i^{\beta_i}$.

\emph{2.} Consider a \sli $AX=B$ in $\ZZ$ which admits an infinity of solutions in $\QQ^m$.
To know if it admits a solution in $\ZZ^m$ we can try a \lgb method. Start by determining a solution in $\QQ$, which is a vector $X\in\QQ^{m}$. Find an integer $d$ such that $dX\in\ZZ^{m}$, such that $X$ has \coes in $\ZZ[1/d]$.
It then suffices to construct a solution in each localized ring $\ZZ_{1+p\ZZ}$ for the prime $p$'s which divide $d$ and to apply the \plgc \ref{plcc.basic}. To know if there is a solution in $\ZZ_{1+p\ZZ}$ and to construct one, we can use the pivot method, provided we take as pivot an \elt of the matrix (or rather the remaining part of the matrix) which divides every other \coez, \cad a \coe wherein $p$ appears with a minimum exponent.\\
The drawback of this method is that it requires factorizing $d$, which can render it unfeasible. \\
However, we can slightly modify the method in order to avoid having to completely factorize $d$. We will use the partial factorization \algoz. Start as if $d$ were a prime number. More \prmt work with the \ri $\ZZ_{1+d\ZZ}$.
Check whether a \coe of the matrix is comaximal to $d$. If one is found, use it as your pivot.
Otherwise no \coe of the matrix is comaximal to $d$ and (by using if \ncr the partial factorization \algoz) we have one of the following three cases:
\begin{itemize}
\item $d$ divides all the \coes of the matrix, in which case, either it also divides the \coes of $B$ and it is reduced to a simpler \pbz, or it does not divide any \coe of $B$ and the \sli has no solution,
\item  $d$ is written as a product of pairwise comaximal factors $d=d_1\cdots d_k$ with $k\geq 2$, in which case we can then work with the \lons at the \mos  $(1+d_1\ZZ)$, \ldots,   $(1+d_k\ZZ)$,
\item   $d$ is written as a pure power of some $d'$ dividing $d$, which, with~$d'$ in place of $d$, brings us to a similar but simpler \pbz.
\end{itemize}
Check that we can recursively exploit the idea expressed above. Write an \algo and test it.
Examine whether the obtained \algo runs in a reasonable time.
}
\end{problem}

}

\sol{

\exer{exoNilpotentChap2}{
\emph{1.} Assume \spdg $a_0 = b_0 = 1$.
When you write $fg = 1$, you get 

\snic{0 = a_n b_m$, $0 = a_n b_{m-1}
+ a_{n-1} b_m$, $0 = a_n b_{m-2} + a_{n-1} b_{m-1} + a_{n-2} b_{m},}

and so on up to degree~$1$.\\
Then prove by \recu over~$j$ that $\deg(a_n^j g) \leq m-j$. 
\\
In particular, for $j = m+1$, we get $\deg(a_n^{m+1} g) \leq -1$, i.e.
$a_n^{m+1} g = 0$. Whence $a_n^{m+1} = 0$.
Finally, by reasoning modulo $\DB(0)$, we obtain successively nilpotent $a_j$'s for $j=n-1$, $\ldots$, $1$.

\noindent
\emph{2a.} Consider the \pols over the commutative \ri $\gB[A]$:

\snic{
f(T) = \det(\In - TA) \; \hbox{ and } \;
g(T) = \det(\In + TA + T^2A^2 + \cdots + T^{e-1}A^{e-1}).}

We have $f(T) g(T) = \det(\In - T^eA^e)=1$.
The \coe of degree~$n-i$ of $f$ is $\pm a_i$. Apply \emph{1}.
}

\noindent
\emph{2b.} It suffices to prove that $\Tr(A)^{(e-1)n + 1} = 0$, because $a_i = \pm
\Tr\big(\bigwedge^{n-i}(A)\big)$.\\
 Consider the \deter defined with respect to a fixed basis $\cB$ of $\Ae n$. If we take the canonical basis formed by the $e_i$'s, we have an obvious \egt 
$$
\Tr(f)= \det_\cB(f(e_1),e_2,\dots,e_n)+\cdots+\det_\cB\big(e_1,e_2,\dots,f(e_n)\big).
$$
It can be written in the following form:
$$
\Tr(f)\det_\cB(e_1,\dots,e_n)= \det_\cB(f(e_1),e_2,\dots,e_n)+\cdots+
\det_\cB\big(e_1,e_2,\dots,f(e_n)\big).
$$
In this form we can replace the $e_i$'s by any \sys of $n$ vectors of $\Ae{n}$: both sides are  $n$-multi\lin alternating forms (at the $e_i$'s) over $\Ae n$, therefore are equal because they coincide on a basis.
\\
 Thus, multiplying a \deter by $\Tr(f)$ reduces to replacing it by a sum of \deters in which $f$ acts on each vector.
\\
One deduces that the expression $\Tr(f)^{n(e-1)+1}\det_\cB(e_1,\dots,e_n)$ is equal to a sum of which each term is a \deter of the form 
$$
 \det_\cB\big(f^{m_1}(e_1),f^{m_2}(e_2),\dots,,f^{m_n}(e_n)\big),
$$
with $\som_i m_i=n(e-1)+1 $,
therefore at least one of the exponents $m_i$ is $\geq e$.

\noindent \rem This solution for the bound $n(e-1)+1$ is due to Gert Almkvist. See on this matter:
{\sc Zeilberger D.} {\it Gert Almkvist's generalization of a mistake of Bourbaki.}
Contemporary Mathematics {\bf 143} (1993), p.~609--612. \eoe

\exer{exolemUMD}{\emph{1.} Let $\fa=\gen{x_1,\ldots,x_n}$. Obtaining $s^r\in\fa$ (for some $r$), \hbox{and $1-as\in\fa$} (for some $a$). Write $1=a^rs^r+(1-as)(1+as+\cdots) \in \fa$.\\
\emph{2.}
$\fa+\fb=\gen{1}$, $\fa+\fc=\gen{1}$ and $(\fa+\fb)(\fa+\fc)\subseteq \fa+\fb\fc $,
therefore
$\fa+\fb\fc=\gen{1}$.

}

\exer{exoPlgb1}{
\emph{1.}
Since $f$ is \hmgz, we have $f(tx) = 0$ for a new
\idtr $t$. Whence each $U_i \in \gA[X_1,
\ldots, X_n, t]$ such that $f = \sum_{i=1}^n (X_i - tx_i)U_i$.
\\ By making $t := x_1^{-1}X_1$, we obtain each $v_i \in \gA[X_1, \ldots,
X_n]$ such that 

\snic{f = \sum_{i=2}^n (x_1X_i - x_iX_1)v_i.}

Finally, since $f$ is \hmg of degree $d$, we can replace $v_i$ by its \hmg component of degree $d-1$.
\\
\emph{2.}
Consider the \egt $f = \sum_{k,j}(x_k X_j - x_jX_k)u_{kj}$, where the $u_{kj}$'s are \pogs of degree $d-1$. It is a \sli in the \coes of the~$u_{kj}$'s. Since this \sys admits a solution over each localized~$\gA_{x_i}$ and that the~$x_i$'s are \comz, it admits a solution over~$\gA$.
\\
\emph{3.}
If $F = \sum_d F_d$ is the \dcn of $F \in \gA[X_1, \ldots, X_n]$ into \hmg components, we have $F(tx) = 0$ \ssi $F_d(x) = 0$ for all $d$, whence the first item of the question. For the saturation, we prove that if $X_iF \in I_x$ for all $i$, then $F \in I_x$.
But we have $x_iF(tx) = 0$. Therefore, by comaximality of the $x_i$'s, we get $F(tx) = 0$,
i.e. $F \in I_x$.
}

\exer{exoNormPrimitivePol} 
The \pol $Q$, regarded as a \pol with \coes in $\gB'$, remains primitive and therefore \ndz (Gauss-Joyal, item \emph{\iref{i2lemPrimitf}} of Lemma~\ref{lemGaussJoyal}). 
Since $m_Q$ is injective, its \deter $\det(m_Q) = P \in \gA'[\uX]$ is \ndz (\Thref{prop inj surj det}, item~{\emph 2}). But $P$ is \egmt null in $\gA'[\uX]$. Thus $\gA'$ is the null \riz, in other words $1 \in \rc_\gA(P)$.

\exer{exoIdempAX}
Let $f(X)$ be an \idm of $\AX$. Clearly $e=f(0)$ is \idmz.
We want to prove that $f=e$. For this we can reason separately modulo $e$ and modulo $1-e$.
\\
  If $e=0$, then $f=Xg$. We have $(Xg)(1-Xg)=0$, or $1-Xg$ is \ndzz, thus~$g=0$.
\\
 If $e=1$, consider the \idm $1-f$ and we are reduced to the previous case.

\exer{exoIdmsSupInf}
For question \emph{5} first prove the result when $uv=0$.
In the \gnl situation, write $u'=1-u$ and $v'=1-v$. We then have a \sfio $(uv,uv',u'v,u'v')$ and by applying the previous special case we see that the two \ris are isomorphic to the product  
\hbox{$\aqo{\gA}{uv'}\times (\aqo{\gA}{uv})^2\times \aqo{\gA}{u'v}$}.


\exer{exoFracSfio}\emph{2.} We have $\gK[1/e_i]\simeq \gK\sur{\Ann_\gK(e_i)}$ and $\Ann_\gK(e_i)=\Ann_\gA(a_i)\gK$. For an \elt $x$ of~$\gA$, write $dx=\sum_{i\in\lrbn}x_i$ in $\gK$, with $x_i=e_idx=a_ix$. The \dcn is thus entirely in~$\gA$. 
Since $dx\equiv x_i\mod \Ann_\gA(a_i)$  the component $\gK\sur{\Ann_\gK(e_i)}$ of the product, when seen as the \id $e_i\gK$, is formed from the \elts of the form~$a_ix/y$ with $x\in\gA$ and~$y$ \ndz in~$\gA$. 
But~$y$ is \ndz in~$\gA$ \ssi each $y_i=a_iy$ is \ndz modulo $\Ann_\gA(a_i)$, 
so that $\gK\sur{\Ann_\gK(e_i)}$ is identified with $\Frac(\gA\sur{\Ann_\gA(a_i)}\!)$.


\exer{exoFracZed}
\emph{1.} The zeros of $\gA$ are the three \gui{coordinate axes.}
\\
Every \elt of $\gA$ is uniquely written in the form 

\snic{u=a+xf(x)+yg(y)+zh(z),}

with $f$, $g$, $h\in\QQ[T]$. This implies that $x+y+z$ is \ndz
because 
$$\preskip.3em \postskip.3em 
(x+y+z)u=x\big(a+xf(x)\big)+y\big(a+yg(y)\big)+z\big(a+zh(z)\big). 
$$
So the \elts $\frac{x}{x+y+z}$, $\frac{y}{x+y+z}$ and $\frac{z}{x+y+z}$ form a \sfio of $\gK$. 
Conclude with Exercise \ref{exoFracSfio} by noting that $\Ann_\gA(x)=\gen{y,z}$, and thus that
$$\preskip.2em \postskip.3em 
\gA\sur{\Ann_\gA(x)}\simeq \QQX. 
$$
\emph{2.} The zeros of $\gB$ are the three \gui{coordinate planes.} The \sfio in $\gL$ is given by \smashbot{$\frac{uv}{uv+vw+wu}$, $\frac{vw}{uv+vw+wu}$} and $\frac{wu}{uv+vw+wu}$.

\vspace{2pt}
\exer{exoIdempotentE2}
It suffices to solve the question modulo $a$ and modulo $1-a$.\\
Modulo $a$: $\cmatrix{a\cr b}=\cmatrix{0\cr b}\mapsto \cmatrix{b\cr b}\mapsto \cmatrix{b\cr 0}$.\\          
Modulo $1-a$, $\cmatrix{a\cr b}=\cmatrix{1\cr b}\mapsto  \cmatrix{1\cr 0}$.
By patching: $d=(1-a)b+a$ with for example the matrix
$A=A_2A_1,$ where
$$  A_1=(1-a)\bloc{1}{1}{0}{1}+ a \bloc{1}{0}{0}{1}=
\bloc{1}{1-a}{0}{1} , $$
$$
A_2=(1-a)\bloc{\phantom-1}{0}{-1}{1}+ a \bloc{\phantom-1}{0}{-b}{1}=\bloc{1}{0}{a-ab-1}{1}  $$
and
$$\preskip-.4em \postskip-.0em
A=\bloc{1}{1-a}{a-ab-1}{a}.
$$

\exer{exoFacileGrpElem4}{
The matrix $D_0 = \crmatrix {0 & -1\cr 1& 0\cr}$ transforms $\cmatrix {x\cr
y\cr}$ into $\crmatrix {-y\cr x\cr}$, so $D_0^2 = -\I_2$ and $D_0^3 = -D_0
= D_0^{-1}$. \\
We also have $D_0=\rE_{12}(1)\rE_{21}(-1)\rE_{12}(1)$,
$D_0 D_q= -\rE_{12}(q)$ and $ D_q D_0=- \rE_{21}(q)$.
}

\exer{exoTroisiemeLemmeLiberte}
Let $(e_1, \ldots, e_n)$ be the canonical basis of~$\Ae n$ and $(f_1, \ldots, f_n)$ the $n$ columns of~$F$.  We can assume that the \iv principal minor is in the north-west position
such that $(f_1, \ldots, f_k, e_{k+1}, \ldots, e_n)$ is a basis of~$\Ae n$. \\
Since $F(f_j) = f_j$, the matrix of~$F$ with respect to this basis is $G \eqdf {\rm def} \cmatrix {\I_k & *\cr 0 & *\cr}$.

The matrix $G$ is \idme as well as its transposed $G'$.
Apply to the \prr $G'$ the operation that we just subjected to~$F$.\\
Since $G'(e_j) \in \bigoplus
_{i \ge k+1} \gA e_i$ for $j \ge k+1$,
the matrix of $G'$ with respect to the new basis is of the form $H = \cmatrix {\I_k & 0\cr 0 & *\cr}$, whence the result because $F$ is similar \hbox{to $\tra{H}$}.

\exer{exoABArang1}{
We have each $b_{ji} \in \gA$ such that $1 = \sum_{i,j} b_{ji} a_{ij}$.  Let $B \in \Ae {m \times n}$ be defined by $B = (b_{ji})$. Check that $ABA = A$: $(ABA)_{ij} = \sum_{l,k} a_{il} b_{lk} a_{kj}$.  \\
But $\dmatrix {a_{il} & a_{ij}\cr
a_{kl} & a_{kj}\cr} = 0$, so $(ABA)_{ij} = \sum_{l,k} a_{ij} a_{kl}
b_{lk} = a_{ij} \sum_{l,k} a_{kl} b_{lk} = a_{ij}$. Consequently, $AB$ is a \prrz.  \\
Let us prove that $AB$ is of rank~$1$. We have
$\Tr(AB) = \sum_i (AB)_{ii} = \sum_{i,j} a_{ij} b_{ji} = 1$,
\hbox{thus $\cD_1(AB)=1$}.
Furthermore,
  $\cD_2(AB)\subseteq\cD_2(A)=0$.
}

\exer{exoABAabstrait} {\emph{1.}
Fix a \lin form $\mu$. The map $E^{r+1} \to \gA$ defined by

\snic{(y_0, \ldots, y_r) \mapsto \som_{i=0}^r
(-1)^i f(y_0, \ldots, y_{i-1}, \widehat {y_i}, y_{i+1}, \ldots, y_r)\, \mu(y_i),}

where teh  symbol $\widehat {y_i}$ denotes the omission of the \eltz, is an  $(r+1)$-multi\lin alternating form.

  According to the hypothesis $\Vi_{r+1}(x_1,
\ldots, x_n) = 0$ and the injectivity of $E \mapsto E{\sta}{\sta}$,
we obtain

\snic{
\som_{i=0}^r
(-1)^i f(y_0, \ldots, y_{i-1}, \widehat {y_i}, y_{i+1}, \ldots, y_r)\, y_i = 0.}

Write $y$ instead of $y_0$ and execute the following operation: in the expression

\snic{
(-1)^i f(y, \ldots, y_{i-1}, \widehat {y_i}, y_{i+1}, \ldots, y_r),}

bring $y$ between $y_{i-1}$ and $y_i$. The permutation thus executed necessitates a multiplication by $(-1)^{i-1}$. We then obtain the second \egt in which all the signs \gui{have disappeared.} For example with $r = 4$, the expression
$$
\begin {array} {c}
f(\widehat {y}, y_1, y_2, y_3, y_4)y -
f(y, \widehat {y_1}, y_2, y_3, y_4)y_1 +
f(y, y_1, \widehat {y_2}, y_3, y_4)y_2 -
\\[.1em]
f(y, y_1, y_2, \widehat {y_3}, y_4)y_3 +
f(y, y_1, y_2, y_3, \widehat {y_4})y_4 =
\\[.2em]
f(y_1, y_2, y_3, y_4)y -
f(y, y_2, y_3, y_4)y_1 +
f(y, y_1, y_3, y_4)y_2 -
\\[.1em]
f(y, y_1, y_2, y_4)y_3 +
f(y, y_1, y_2, y_3)y_4
\\
\end {array}
$$
is none other than
$$
\begin {array} {c}
f(y_1, y_2, y_3, y_4)y -
f(y, y_2, y_3, y_4)y_1 -
f(y_1, y, y_3, y_4)y_2 -
\\[.1em]
f(y_1, y_2, y, y_4)y_3 -
f(y_1, y_2, y_3, y)y_4.
\end {array}
$$
A faster proof: apply a \lin form $\mu$ to the last expression above, check that the obtained map $(y, y_1, y_2, y_3, y_4) \mapsto \mu(\ldots)$ is 
$5$-multi\lin alternating, and therefore is null by the assumptions.

\noindent  \emph{2.}
Treat the $r = 3$ case. We have an assumption
$$
1 = \som_{ijk} \alpha_{ijk} f_{ijk}(x_i,x_j,x_k), \qquad
\hbox {$f_{ijk}$ $3$-multi\lin alternating over $E$.}
$$
Define $\pi : E \to E$ by:
$$
\pi(x) = \som_{ijk} \alpha_{ijk}
[ f_{ijk}(x,x_j,x_k)x_i +  f_{ijk}(x_i,x,x_k)x_j +  f_{ijk}(x_i,x_j,x) x_k].
$$
Clearly, the image of $p$ is contained in the submodule $\som \gA x_i$.
In addition, \hbox{for $ {x \in \som \gA x_i}$}, we have 
$$
f_{ijk}(x,x_j,x_k)x_i +  f_{ijk}(x_i,x,x_k)x_j +  f_{ijk}(x_i,x_j,x) x_k =
f_{ijk}(x_i,x_j,x_k)x.
$$
Whence $\pi(x) = x$: the endomorphism $\pi : E \to E$ is a \prr onto
$\som \gA x_i$. Notice that $p$ is of the form $\pi(x) = \som_i \alpha_i(x)
x_i$ i.e. $\pi = \psi \circ \varphi$ and that $\pi \circ \psi = \psi$.

\noindent  \emph{3.}
The module $E$ in question is $\Ae{m}$ and the vectors $x_1$, $\ldots$, $x_n$ are the columns of~$A$. We have $\psi = A : \Ae{n} \to \Ae{m}$, and if we let $B \in \Ae {n \times m}$ be the matrix of $\varphi : \Ae{m} \to \Ae{n}$, we indeed have~$ABA = A$.  So, the \ali $AB : \Ae{m} \to \Ae{m}$ is 
a \prr having the same image as $A$.
}

\exer{exoDetExtPower}{
Let us first see the case where $u=\Diag(\lambda_1, \ldots, \lambda_n)$. We have a basis $(e_I)$ of $\Vi^k(\Ae{n})$ indexed by the subsets $I \subseteq \{1,\ldots, n\}$ of cardinality~$k$:

\snic{e_I = e_{i_1} \wedge \cdots \wedge e_{i_k} \qquad
I = \so{i_1 < \cdots < i_k }.
}

Then, $u_k$ is diagonal with respect to the basis $(e_I)$: $u_k(e_I) = \lambda_I e_I$ with $ {\lambda_I = \prod_{i \in I} \lambda_i}$.
It follows that $ {\det(u_k) = \prod_{\#I = k} \prod_{i \in I} \lambda_i}$.  It remains to determine, for some $j$ given in~$\lrbn$, the number of occurrences of $\lambda_j$ in the above product. In other words, how many subsets $I$, of cardinality~$k$, contain~$j$? As many as there are subsets of cardinality $k-1$ contained in $\{1, \cdots, n\} \setminus \{j\}$,
i.e.  ${n-1 \choose k-1}$. The result is proven for a generic matrix. Thus it is true for any matrix. The second point follows from the \egts

\snic{\dsp{n-1 \choose k-1} + {n-1 \choose n-k-1} =
{n-1 \choose k-1} + {n-1 \choose k} = {n \choose k}.
}
}

\exer{exolemdeterblocs}{The \gnl case is treated by \recu on $n$.
Consider the \pol \ri $\ZZ[(x_{ij})]$ with $n^2$ \idtrs and the \uvle matrix
$A=(x_{ij})$ with \coes in this \riz. Let $\Delta_{1k} \in \ZZ[(x_{ij})]$ be the cofactor of $x_{1k}$ in $A$. These cofactors satisfy the identities:
$$
\som_{j = 1}^n x_{1j} \Delta_{1j} =
\det A,
\;\;\;\;
\som_{j = 1}^n x_{ij} \Delta_{1j} = 0 \quad \hbox {for } i > 1  .
$$
  Since the $N_{kl}$'s pairwise commute, the specialization $x_{kl} \mapsto N_{kl}$ is legitimate. Let $N'_{1j} = \Delta_{1j}(x_{kl} \mapsto N_{kl})$, then we have

\snic{
N'_{11} =
\sum_{\sigma \in \rS_{n-1}} \varepsilon(\sigma)
N_{2\sigma_2} N_{3\sigma_3} \ldots N_{n\sigma_n}.}

Let us define $N'$ by:
$$
N' = \cmatrix {
N'_{11}      &  0    & \cdots & 0      \cr
N'_{12}      &  \rI_m  &        & \vdots \cr
\vdots       &  \vdots     & \ddots & 0      \cr
N'_{1n}      &  0    & \cdots       & \rI_m    \cr}\,,
\; \hbox {so that} \;\;
N N' = \cmatrix {
\Delta  & N_{12} & \cdots & N_{1n} \cr
0       & N_{22} & \cdots & N_{2n} \cr
\vdots  &        &        & \vdots \cr
0       & N_{n2} & \cdots & N_{nn} \cr
}.
$$
By taking  \detersz, we get
$$
\det(N) \det(N'_{11}) = \det(\Delta)
\det \cmatrix {
N_{22} & \cdots & N_{2n} \cr
\vdots &        & \vdots \cr
N_{n2} & \cdots & N_{nn}  }.
$$
The \hdr provides the \egts
$$
\det \cmatrix {
N_{22} & \cdots & N_{2n} \cr
\vdots &        & \vdots \cr
N_{n2} & \cdots & N_{nn} \cr} =
\det\Big( \sum_{\sigma \in \rS_{n-1}} \varepsilon(\sigma)
N_{2\sigma_2} N_{3\sigma_3} \cdots N_{n\sigma_n} \Big)
= \det(N'_{11}).
$$
Simplification by the \ndz \elt $\det(N'_{11})$ gives the \egt $\det(N) = \det(\Delta)$.
}

\prob{exoRationaliteLineaire} 
\emph {1.}
If $A_j$ is a nonzero column of $A$, we have $BA_j = e_j$ and therefore $ABA_j = A_j$; thus $AB$ is the identity over $\Im A$, so $ABA = A$. The matrix $AB$ is lower triangular, and its diagonal \coes are $0,1$. The matrix $BA$ is diagonal and its diagonal \coes are $0,1$.

\snic {
B = \cmatrix {
\cdot&   \cdot&   \cdot&   \cdot&   1&  \cdot&  \cdot&  \cdot&  \cdot\cr
\cdot&   \cdot&   \cdot&   \cdot&   \cdot&  1&  \cdot&  \cdot&  \cdot\cr
\cdot&   \cdot&   1&   \cdot&   \cdot&  \cdot&  \cdot&  \cdot&  \cdot\cr
1&   \cdot&   \cdot&   \cdot&   \cdot&  \cdot&  \cdot&  \cdot&  \cdot\cr
\cdot&   \cdot&   \cdot&   \cdot&   \cdot&  \cdot&  \cdot&  1&  \cdot\cr
\cdot&   \cdot&   \cdot&   \cdot&   \cdot&  \cdot&  \cdot&  \cdot&  \cdot\cr
},
\qquad
BA = \cmatrix {
1&  \cdot&  \cdot&  \cdot&  \cdot&  \cdot\cr
\cdot&  1&  \cdot&  \cdot&  \cdot&  \cdot\cr
\cdot&  \cdot&  1&  \cdot&  \cdot&  \cdot\cr
\cdot&  \cdot&  \cdot&  1&  \cdot&  \cdot\cr
\cdot&  \cdot&  \cdot&  \cdot&  1&  \cdot\cr
\cdot&  \cdot&  \cdot&  \cdot&  \cdot&  \cdot\cr
},}

\snic {
AB = \cmatrix {
1&      \cdot&       \cdot&   \cdot&       \cdot&       \cdot&   \cdot&       \cdot&  \cdot\cr
a_{24}& \cdot&       \cdot&   \cdot&       \cdot&       \cdot&   \cdot&       \cdot&  \cdot\cr
\cdot&      \cdot&       1&   \cdot&       \cdot&       \cdot&   \cdot&       \cdot&  \cdot\cr
a_{44}& \cdot&  a_{43}&   \cdot&       \cdot&       \cdot&   \cdot&       \cdot&  \cdot\cr
\cdot&      \cdot&       \cdot&   \cdot&       1&       \cdot&   \cdot&       \cdot&  \cdot\cr
\cdot&      \cdot&       \cdot&   \cdot&       \cdot&       1&   \cdot&       \cdot&  \cdot\cr
a_{74}& \cdot&  a_{73}&   \cdot&  a_{71}&  a_{72}&   \cdot&       \cdot&  \cdot\cr
\cdot&      \cdot&       \cdot&   \cdot&       \cdot&       \cdot&   \cdot&       1&  \cdot\cr
a_{94}& \cdot&  a_{93}&   \cdot&  a_{91}&  a_{92}&   \cdot&  a_{95}&  \cdot\cr
}.}

The \supl subspace $\Ker AB$ of $\Im A = \Im AB$ in $\gK^n$ admits as its basis the~$e_i$'s for the indices $i$ of the rows that do not contain a pivot index.\\
In the example, $(e_2, e_4, e_7, e_9)$ is a basis of $\Ker AB$.

\noindent\emph {2.}
We obtain $(Q, A')$ by Gauss' (classical) pivot method. 
If the matrix $B' \in M_{n,m}(\gK)$ satisfies $A'B'A' = A'$, then $AQB'AQ = AQ$, therefore the matrix $B = QB'$ satisfies $ABA = A$.

\noindent\emph {3.}
Consider a matrix $B \in \MM_{m,n}(\gK)$ such that $ABA = A$. Then, if $y = Ax$ for some $m$-vector with \coes in an over\ri of $\gK$, 
we have $A(By) = y$, whence the existence of a solution on $\gK$, namely $By$.

\noindent\emph {4.}
Let $(u_1, \ldots, u_r)$ be a \sgr of the $\gK$-\evc $E$, constituted of vectors of $\gK_0^n$; similarly for $(v_1, \ldots, v_s)$ and $F$. Let $z \in \gK_0^n$, which we want to express in the form $z = x_1 u_1 + \cdots + x_r u_r + y_1 v_1 + \cdots + y_s v_s$ with each~$x_i, y_j \in \gK_0$. We thus obtain a $\gK_0$-\lin \sys from the unknowns $x_i$'s, $y_j$'s which admits a solution on $\gK$, therefore \egmt on $\gK_0$.

\noindent\emph {5.a.}
If every $\pi(e_j)$ is in $\gK_0^n$, then the subspace $E$, generated by the $\pi(e_j)$'s, is \hbox{$\gK_0$-rational}. 
Conversely, if $E$ is \hbox{$\gK_0$-rational}, since $F$ is also $\gK_0$-rational, by the previous question we have 
$\pi(e_j) \in \gK_0^n$ for all $j$.

\noindent\emph {b.}
Now trivial: $\gK_0$ is the subfield generated by the components of the $\pi(e_j)$~vectors.

\noindent\emph {c.}
The field of rationality of a strictly staggered matrix is the subfield generated by the \coes of the matrix.
For example with $E = \Im A \subset \gK^5$:

\snic{
A = \bordercmatrix [\lbrack\rbrack]{
    & w_1 & w_2 & w_3 \cr
e_1 & 1   & 0   & 0 \cr
e_2 & a   & 0   & 0 \cr
e_3 & 0   & 1   & 0 \cr
e_4 & 0   & 0   & 1 \cr
e_5 & b   & c   & d \cr
},}

we get $E = \gK w_1 \oplus \gK w_2 \oplus \gK w_3$ and we have
$\gK^5 = E \oplus F$ with $F = \gK e_2 \oplus \gK e_5$. 
Since 
$$\preskip-.2em \postskip.4em 
e_1 - w_1 \in F, \quad e_3 - w_2 \in F, \quad e_4 - w_3 \in F,  
$$
we have $\pi(e_1) = w_1$, $\pi(e_3) = w_2$, $\pi(e_4) = w_3$
and $\pi(e_2) = \pi(e_5) = 0$. The field of rationality of $E$ is $\gK_0 = \gk(a,b,c,d)$, where $\gk$ is the prime subfield of $\gK$.

}

\Biblio

The \iJG Gauss-Joyal Lemma is in~\cite{Es}, which gives it its name.
On the \gnl subject of comparison between the \ids $\rc(f)\rc(g)$ and $\rc(fg)$ see~\cite{CDLQ02,Glaz2,Nor2} and, in this work, Sections \ref{secLemArtin} and \ref{secThKro}
and Proposition~\ref{propLG}.

Regarding the \cof treatment of  \noetz, see
\cite{MRR,JL,Per1,Per2,ric74,sei74b,sei84,Ten}.

The whole of Section \ref{secCramer} can be more or less found in \cite{Nor}.
For example the formula (\ref{eqIGCram}) on \paref{eqIGCram} is found in a related form in \Tho 5 on page~10.
Likewise, our Cramer-style magic formula (\ref{eqIGCram2}) on \paref{eqIGCram2} is very similar to \Thoz~6 on page 11: Northcott attaches central importance to the matrix equation $A\,B\,A=A$.
On this subject, see also \cite{RM} and \cite[D\'{\i}az-Toca\&al.]{DiGLQ}.

Proposition \ref{propIGCram2} is in \cite{Bha} \Thoz~5.5.

Concerning \Thref{theoremIFD}: in \cite{Nor} \Tho 18 on page 122 establishes the \eqvc of items  \emph{\ref{IFDa}} and \emph{\ref{IFDe}}
by a method which is not entirely \covz, but \Tho 5 page 10
would allow us to give an explicit formula for the implication \emph{\ref{IFDe}} $\Rightarrow$ \emph{\ref{IFDa}}.

\newpage \thispagestyle{CMcadreseul}
\incrementeexosetprob


\chapter{The method of undetermined coefficients}
\label{chapGenerique}
\perso{compil\'e le \today}
\minitoc

\Intro

\begin{flushright}
{\em Weil Gauss ein echter Prophet der Wissenschaft ist,\\
deshalb reichen die Begriffe,\\
die er aus der Tiefe der Wissencshaft sch\"opft,\\
weit hinaus \"uber den Zweck, \\
zu welchem sie aufgestellt wurden.
}\\
Kronecker\\
Vorlesungen Sommersemester 1891. Le{\c c}on 11 \cite{BoSc} \\[1mm]
Approx. transl.\\
{\em Because Gauss is a true Prophet of Science,\\
the concepts that he draws from the depths of Science\\
go beyond the purpose for which they were established.
}
\end{flushright}

In 1816, Gauss published a fundamental article \cite{Gauss} in which he corrects (without citing) the \dem of the fundamental \tho of \alg given by Laplace a few years beforehand. Laplace's \dem is itself remarkable as it is \gui{purely \agqz:} it claims only two very \elrs \prts for real numbers: the existence of the square root of a non-negative number  and that of a zero for a \pol of odd degree.
 
 Gauss' goal is to treat this \tho without using a (hypothetical) field of imaginary numbers, over which an arbitrary \pol would be decomposed into \lins factors. Laplace's \dem implicitly assumes the existence of such a field~$\gK$ containing~$\CC=\RR[i]$, and shows that the \dcn into products of \lins factors actually takes place in~$\CC[X]$. 
 
Gauss' \dem dispenses with the assumption about the field~$\gK$ and constitutes a tour de force 
that shows that you can handle things in a purely formal way. 
He proves the existence of the gcd of two \pols by using Euclid's \algo as well as the corresponding Bézout relation. He shows that every \smq \pol is uniquely expressed as a \pol of \elr \smq functions (by introducing a lexicographical order on the \momsz). He defines the \discri of a \polu purely formally. He shows (without resorting to roots) that every \pol can be decomposed into a product of \pols with a nonzero \discriz. He shows (without resorting to roots) that a \pol admits a square factor \ssi its \discri is zero (he works in zero \caraz). Finally, Gauss makes Laplace's \dem work in a purely formal way, without resorting to a \cdrz, by only using resultants and \discrisz. 
 
 In short, he establishes a \gui{\gnl method of undetermined \coesz} on a firm basis. This was to be systematically reused, 
in particular by Leopold Kronecker, Richard Dedekind, Jules Drach, Ernest Vessiot\ldots

\medskip In this chapter, we introduce the method of undetermined \coes and we give some of its applications.

We begin with some \gnts about \pol \risz. The \DKM lemma and \KROz's theorem are two basic tools which provide precise information about the \coes of a product of two \polsz. These two results will often be used in the remainder of this work.
 
Here we study the \elr \prts of the \discri and the resultant, and we introduce the fundamental tool that is the \adu of a \mon \polz. The latter allows for a simplification of purely formal \dems such as Gauss' by providing a formal substitute for the \polz 's \gui{\cdrz.} 

All of this is very consistent and works with arbitrary commutative \risz. 
The reader will only notice the apparition of fields from 
Section~\ref{secGaloisElr}.

The applications that we treat relate to basic Galois theory, the first steps in algebraic number theory, and Hilbert's \nstz.
We have \egmt dedicated a section to 
Newton's method in \algz.

\subsec{A few words on finite sets}\label{Deux mots}

  A set~$E$ is said to be
 \emph{finite} when we explicitly have a bijection 
  between~$E$ and an initial segment~$\sotq{x\in\NN}{x<n}$ of~$\NN$. It is said to be
 \emph{finitely enumerable} when we explicitly have a surjection of a finite set~$F$ onto~$E$.%
\index{finitely enumerable!set}%
\index{finite!set}%
\index{set!finite ---}%
\index{set!finitely enumerable ---}

\rdb \label{NOTAPfPfe}
In \gnl the context is sufficient to distinguish between the two notions. Sometimes, it is advantageous to be very precise. 
We will make the distinction when necessary by using the notation~$\Pf$ or~$\Pfe$:
we will denote by~$\Pf(S)$ \emph{the set of finite subsets} of the set~$S$ and~$\Pfe(S)$ \emph{the set of finitely enumerable subsets} of $S$.
In \coma when~$S$ is discrete (resp.\,finite), we have the \egt $\Pf(S)=\Pfe(S)$ and it is a discrete set (resp.\,finite).%
\footnote{In \coma
we \gnlt refrain from considering the \gui{set of all subsets of a set,} even finite, because it is not a \gui{reasonable} set:
it does not seem possible to give a clear \dfn of its \elts
(see the discussion \paref{P(X)}).
When we used the notation~$\cP_{\ell}$ for \gui{the set of subsets of~$\so{1,\ldots ,\ell}$,} on \paref{notaAdjalbe}, it was in fact the set of finite subsets of~$\so{1,\ldots ,\ell}$.
}
When~$S$ \emph{is not} discrete,~$\Pf(S)$ \emph{is not} equal to~$\Pfe(S)$.\index{set!of finite subsets}%
\index{set!of finitely enumerable subsets}

Also note that when~$S$ is a finite set every detachable subset (cf. \paref{subsecTestDEgalite}) is finite: the set of finite subsets is then equal to the set of detachable subsets.

The finitely enumerable subsets are omnipresent in the usual \mathe sense. For example when we speak of a \itf we mean an \id generated by a finitely enumerated subset and not by a finite subset. 
Similarly, when we speak of a \emph{finite family}~$(a_i)_{i\in I}$
in the set~$E$, we mean that~$I$ is a finite set, therefore the subset~\hbox{$\sotq{a_i}{i\in I}\subseteq E$} is finitely enumerated.%
\index{finite!family}

\rdb \label{NOTAEnumNonempty}
Finally, a nonempty set~$X$  is said to be
 \emph{enumerable} if there is a surjective map~$x=(x_n):\NN\to X$.%
\index{enumerable!nonempty --- set}

\pagestyle{CMheadings}
\section{\Pol \risz}
\label{secAnnPols}

\subsec{Partial \fcn \algoz}

We assume the reader to be familiar with the extended Euclid \algo which computes the \mon gcd of two \polus in~$\KX$ when~$\gK$ is a \cdi (see for example \Pbmz~\ref{exoAnneauEuclidien}).
\begin{lemma}\label{lemPartialDec}
If~$\gK$ is a \cdiz, we have a \emph{\fapz} \algo for the finite families of \polus in~$\KX$: a \fap for a finite family~$(g_1,\ldots,g_r)$ is given by a finite pairwise comaximal family~$(\lfs)$ of \polus and by the expression of each~$g_i$ in the form
$$\preskip-.4em \postskip.3em
g_i=\prod\nolimits_{k=1}^sf_k^{m_{k,i}}\; (m_{k,i}\in\NN).
$$
The family~$(\lfs)$ is called a \emph{\bdfz} for the family~$(g_1,\ldots,g_r)$.%
\index{partial factorization algorithm}%
\index{factorization!partial ---}%
\index{partial factorization!basis} \end{lemma}
%

\begin{proof}
If the~$g_i$'s are pairwise comaximal, there is nothing left to prove.
Otherwise, assume for example that~$\pgcd(g_1,g_2)=h_0$, $g_1=h_0h_1$ and~$g_2=h_0h_2$ with~$\deg(h_0)\geq1$. We replace the family~$(g_1,\ldots,g_r)$ with the family~$(h_0,h_1,h_2,g_3,\ldots,g_r)$. We note that the sum of the degrees has decreased. We also note that we can delete from the list the \pols equal to $1$, or any repeats of a \polz.
We finish by \recu on the sum of the degrees. The details are left to the reader.
\end{proof}
\vspace{-1pt}

\subsec{\Uvl property of \pol \risz}

\vspace{1pt}
A \pol \ri $\AXn$ satisfies the \uvl \prt which defines it as \emph{the commutative \ri freely generated by~$\gA$ and~$n$ new \eltsz.} This is the \prt described by means of the \evn \homo in the following terms.

\begin{proposition}\label{propApolLibre}
Given two commutative \risz~$\gA$ and~$\gB$, a \homo~$\rho:\gA\to\gB$ and~$n$ \elts $b_1$, $\ldots$, $b_n\in\gB$ there exists a unique \homo $\varphi:\AXn=\AuX\to\gB$ which extends~$\rho$ and which takes the~$X_i$'s to the~$b_i$'s.
\vspace{-.3em}
\Pun{\gA}{j}{\rho}{\AuX}{\varphi}{\gB}{~}{$\varphi(X_i)=b_i$, $i\in\lrbn$. \qquad~}
\end{proposition}

\vspace{-1.3em}
This \homo $\varphi$ is called \emph{the \evn \homoz} (of every~$X_i$ to~$b_i$).
If~$P\in\AuX$ has as its image~$P^{\rho}$ in~$\BXn$, we obtain the \egt $\varphi(P)=P^{\rho}(b_1,\ldots,b_n)$. The \evn \homo is also called a \emph{specialization}, and we say that~$\varphi(P)$ is obtained by \emph{specializing} each~$X_i$
to~$b_i$.
When~$\gA\subseteq\gB$, the \elts $b_1$, $\ldots$, $b_n\in\gB$ are said to be
\emph{\agqt independent over~$\gA$} if the corresponding \evn \homo is injective.%
\index{algebraically independent!elements over a subring}%
\index{evaluation homomorphism}%
\index{specialization}

By Proposition~\ref{propApolLibre} every computation made in~$\AuX$ is transferred into~$\gB$ by means of the \evn \homoz.

Clearly, $\Sn$ acts as a group of \autos of~$\AuX$ by permutation of the \idtrsz: $(\sigma,Q)\mapsto Q(X_{\sigma1},\ldots, X_{\sigma n})$.

The following corollary results \imdt from Proposition~\ref{propApolLibre}.


\begin{corollary}\label{propZXnLibre}
Given $n$ \elts $b_1$, $\ldots$, $b_n$ in a commutative \ri $\gB$, there exists a unique \homo
$\varphi:\ZZXn\to\gB$ which takes every~$X_i$ to~$b_i$.
\end{corollary}

\vspace{-.5em}
\pagebreak

\subsec{\Agq identities}

An \ida is an \egt between two \elts of $\ZZXn$ defined differently. It gets automatically transferred into every commutative \ri by means of the previous corollary.

Since the \ri $\ZZXn$ has particular \prtsz, it happens that some \idas are easier to prove in  $\ZZXn$ than in \gui{an arbitrary \riz~$\gB$.}
Consequently, if the structure of a \tho reduces to a family of \idasz, which is very frequent in commutative \algz, it is often in our interest to use a \ri of \pols with \coes in~$\ZZ$ by taking as its \idtrs the relevant \elts in the statement of the \thoz.

The \prts of the \risz~$\ZZuX$ which may prove useful are numerous.
The first is that it is an integral \riz. So it is a sub\ri of its quotient field~$\QQ(\Xn)$ which offers all the facilities of \cdisz.

The second is that it is an infinite and integral \riz.
Consequently, \gui{all bothersome but rare cases can be ignored.} A case is rare when it corresponds to the annihilation of a \polz~$Q$ that evaluates
to zero everywhere.
It suffices to check the \egt corresponding to the \ida when it is evaluated at the points of~$\ZZ^n$ which do not annihilate~$Q$.
Indeed, if the \ida we need to prove is~$P=0$, we get that the \polz~$PQ$ defines the function over~$\ZZ^n$ that evaluates to zero everywhere, this implies that~$PQ=0$ and thus~$P=0$ since~$Q\neq0$ and~$\ZZuX$ is integral.
This is sometimes called the \gui{extension principle for \idasz.}%
\index{extension principle for algebraic identities}

Other remarkable \prts of~$\ZZ[\uX]$ could sometimes be used,
like the fact that it is a unique factorization domain (UFD) as well as being a 
\fdi \coh \noe \ri of finite \ddkz.

\subsubsection*{An example of application}

\begin{lemma}\label{lemPrincipeIdentitesAlgebriques}
For $A$, $B \in \Mn(\gA)$, we have the following results.
\begin{enumerate}
\item $\wi {AB} = \wi{B}\wi{A}$.
\item $\rC{AB} = \rC{BA}$. 
\item  $\wi {PAP^{-1}} = P\wi {A}P^{-1}$ for $P \in \GLn(\gA)$.
\item  $\wi{\wi{A}} = \det(A)^{n-2} A$ if $n \ge 2$.
\item  \emph{(The Cayley-Hamilton \thoz)} $\rC A(A) = 0$.
\item
  If $\Gamma_A(X)=(-1)^{n+1}\big(\rC A(X)-\rC A(0)\big)\sur{X}$,
we have $\wi{A} = \Gamma_A(A)$ ($n\geq2$).\\
We also have 
$\Tr \big(\wi{A}\big)=(-1)^{n+1}\Gamma_A(0)$.
\item  \emph{(Sylvester's identities)} Let $r\geq 1$ and $s\geq 2$ such that $n=r+s$. Let $C\in\MM_r(\gA)$, $F\in\MM_s(\gA)$,
$D\in\MM_{r,s}(\gA)$, $E\in\MM_{s,r}(\gA)$ be the matrices extracted form $A$
as below 
$$\preskip.0em \postskip.4em
A=\blocs{.9}{.6}{.9}{.6}{$C$}{$D$}{$E$}{$F$}\;.
$$
Let $\alpha_i=\so{1,\ldots,r,r+i}$ and $\mu_{i,j}=\det(A_{\alpha_i,\alpha_j})$
for~$i,j\in\lrb{1..s}$. Then:
$$\preskip.2em \postskip.2em
\det(C)^{s-1}\det(A)=\det\big((\mu_{i,j})_{i,j\in\lrbs} \big).
$$
\item  If $\det A = 0$, then $\Al{2}\wi A=0$.
\end{enumerate}%
\index{Cayley-Hamilton}%
\index{Sylvester!identities}
\end{lemma}
\begin{proof}
We can take all the matrices with undetermined \coes over~$\ZZ$ and localize
the \ri at~$\det P$.
In this case~$A$, $B$, $C$ and~$P$ are \iv in the quotient field of the \riz~$\gB=\ZZ[(a_{ij}),(b_{ij}),(p_{ij})]$.
Furthermore, the matrix $\wi A$ satisfies the \egt $\wi A A=\det(A)\,\In$, which characterizes it since~$\det A$ is \ivz. This provides item \emph{1} via the \egt $\det(AB)=\det (A)\det (B)$, items \emph{3} and~\emph{4}, and item \emph{6} via item \emph{5} and the \egt $\rC A(0)=(-1)^n\det A$.

For item~\emph{2} we note that~$AB=A(BA)A^{-1}$. 

For Cayley-Hamilton's \thoz, we first treat the case of the \emph{companion matrix of a \poluz} $f=T^n-\sum_{k=1}^na_kT^{n-k}$:  \label{matrice.compagne}
\index{matrix!companion --- of a polynomial}%
\index{companion!matrix of a polynomial}
$$
P =\cmatrix{
0&\cdots&\cdots&\cdots&0&a_n\cr
1&0&&&\vdots&a_{n-1}\cr
0&\ddots&\ddots&&\vdots&\vdots\cr
\vdots&\ddots&\ddots&\ddots&\vdots&\vdots\cr
\vdots& & \ddots&1&0&a_2\cr
0&\cdots&\cdots&0&1&a_1
}.
$$
This is the matrix of the 
\gui{multiplication by $t$,} $\mu_t:y\mapsto ty$  (where $t$ is the class of~$T$) in the quotient \ri $\aqo{\AT}{f(T)}=\gA[t]$, expressed over the basis of the \moms ordered by increasing degrees.
Indeed, on the one hand a direct computation shows that~$\rC P(T)=f(T)$. On the other hand~$f(\mu_t)=\mu_{f(t)}=0$, thus~$f(P)=0$.\\
Moreover, in the case of the \gnq matrix, the \deter of the family $(e_1, Ae_1, \ldots, A^{n-1}e_1)$ is \ncrt nonzero, therefore the \gnq matrix is similar to the companion matrix of its \polcar over the quotient field of~$\ZZ[(a_{ij})]$.

\emph{7.} Since $C$ is \ivz, we can use the \gne Gauss' pivot, by left-multiplication by a matrix \smashbot{$\blocs{.9}{.6}{.9}{.6}{$C^{-1}$}{$0$}{$E'$}{$\I_s$}\,,$}
this reduces to the case where~$C=\I_r$ and~$E=0$. 

\smallskip Finally, item \emph{8} results from Sylvester's identity (item \emph{7}) with~$s=2$.
\end{proof}

\rem Item \emph{3} allows us to define the \emph{cotransposed \endoz} of an \endo of a free module of finite rank, from the cotransposed matrix.%
\index{cotransposed!endomorphism}

\subsubsection*{Weights, \hmgs \polsz}

We say that we have defined a weight on a \pol \alg $\AXk$ when we attribute to each \idtr $X_i$ a weight~$w(X_i)\in\NN$.
We then define the weight of the \mom $\uX^{\underline{m}}=X_1^{m_1}\cdots X_k^{m_k}$ as

\snic{w(\uX^{\underline{m}})=\sum_im_iw(X_i),}

so that~$w(\uX^{\underline{m}+\underline{m'}})=
w(\uX^{\underline{m}})+w(\uX^{\underline{m'}})$. The degree of a \polz~$P$ for this weight, \gnlt denoted by $w(P)$, is the greatest of the weights of the \moms appearing with a nonzero \coez.
This is only well-defined if we have a test of \egt to~$0$ in~$\gA$ at our disposal. In the opposite case we simply define the statement \gui{$w(P)\leq r$.}

A \pol is said to be
 \ixc{homogeneous}{polynomial} (for a weight~$w$) if all of its \moms have the same weight.

When we have an \ida and a weight available, each \hmg component of the \ida provides a particular \idaz.

We can also define weights with values in some \mos with a more complicated order than~$(\NN,0,+,\geq)$. We then ask that this \mo be the positive part of a product of totally ordered Abelian groups, or more \gnlt a \mo with gcd (this notion will be introduced in Chapter~\ref{chapTrdi}).

\subsec{\Smq \polsz}

We fix $n$ and $\gA$ and we let $S_1$, $\ldots$, $S_n$ be the \emph{\elr \smq \pols 
at the~$X_i$'s} in~$\AXn$. They are defined by the \egt
$$\preskip.4em \postskip.4em
T^n+S_1T^{n-1}+S_2T^{n-2}+\cdots+S_n=\prod\nolimits_{i=1}^n(T+X_i).
$$
We have $S_1=\sum_iX_i$, $S_n=\prod_iX_i$,
$S_k=\sum_{J\in\cP_{k,n}}\prod_{i\in J}X_i$.
Recall the following well-known \tho (a \dem is suggested in Exercise~\ref{exothSymEl}).%
\index{polynomial!elementary symmetric ---}

\begin{theorem}\label{thSymEl} \emph{(\Elr \smq \polsz)}
\begin{enumerate}
\item A \polz~$Q\in\AXn=\AuX$, invariant under permutations of the variables, is uniquely expressible as a \pol in the \elr \smq functions $S_1$, \ldots, $S_n$. In other words
\begin{itemize}
  \item the sub\ri of the fixed points of $\AuX$ by the action of the \smq group $\Sn$ is the \ri $\gA[S_1,\ldots,S_n]$ generated by~$\gA$ and the~$S_i$'s, and 
  \item the $S_i$'s are \agqt independent over $\gA$.
\end{itemize}
\item Let us denote by $d(P)$ the total degree of $P\in\AuX$ when each $X_i$ is affected by the weight~$1$, and~$d_1(P)$ its degree in~$X_1$.
Let~$\delta(Q)$ be the total degree of~$Q\in\gA[S_1,\ldots,S_n]$ when each variable~$S_i$ is affected by the weight~$i$ and~$\delta_1(Q)$ its total degree when each variable~$S_i$ is affected by the weight~$1$.
Assume that~$Q(S_1,\ldots,S_n)$ is evaluated in~$P(\uX)$.
\begin{itemize}
\item [a.] $d(P)=\delta(Q)$, and if $Q$ is $\delta$-\hmgz, then $P$ is $d$-\hmgz.
\item [b.] 
$d_1(P)=\delta_1(Q)$.
\end{itemize}

\item $\AXn$ is a free module of rank $n!$ over $\gA[S_1,\ldots,S_n]$ and a basis is formed by the \moms $X_1^{k_1}\cdots X_{n-1}^{k_{n-1}}$ such that~$k_i\in\lrb{0.. n-i}$ for each~$i$.
\end{enumerate}
\end{theorem}

\begin{corollary}\label{corthSymE1}
On a \ri $\gA$ consider the \gnq \pol 
$$\preskip.3em \postskip.2em 
f=T^n+f_1T^{n-1}+f_2T^{n-2}+\cdots+f_n, 
$$
where the~$f_i$'s are the \idtrsz.
We have an injective \homo $j:\gA[f_1,\ldots,f_n]\to\AXn$ such that the~$(-1)^{k}j(s_k)$'s
are the \elr \smq \pols in the~$X_i$'s.
\end{corollary}

In short we can always reduce to the case where~$f(T)=\prod_i(T-X_i)$, where the~$X_i$'s are other \idtrsz.

\begin{corollary}\label{corthSymE2}
On a \ri $\gA$ consider the \gnq \pol 
$$\preskip.3em \postskip.2em 
f=f_0T^n+f_1T^{n-1}+f_2T^{n-2}+\cdots+f_n, 
$$
where the~$f_i$'s are \idtrsz. 
We have an injective \homo  
{
$j:\gA[f_0,\ldots,f_n]\to\gB = \gA[F_0,\Xn]$}, with the following \egt in~$\gB[T]$.
$$\preskip.0em \postskip.4em 
\dsp j(f_0)\,T^n+j(f_1)\,T^{n-1}+\cdots+j(f_n)=F_0\,\prod\nolimits_i\,(T-X_i)\,. 
$$
\end{corollary}

In short, we can always reduce to the case where~$f(T)=f_0\prod_i(T-X_i)$, 
with \idtrs $f_0$, $X_1$, \dots, $X_n$.

\begin{proof}
It suffices to see that if $f_0$, $g_1$, $\ldots$, $g_n\in\gB
$ are \agqt independent over~$\gA$, then the same goes for $f_0$, $f_0g_1$, $\ldots$, $f_0g_n$.
It suffices to verify that $f_0g_1$, $\ldots$, $f_0g_n$ are \agqt independent over~$\gA[f_0]$. This results from $f_0$ being \ndz and from $g_1$, $\ldots$, $g_n$ being \agqt independent over~$\gA[f_0]$.
\end{proof}
%

\section{Dedekind-Mertens lemma}
\label{secLemArtin}

Recall that for a \polz~$f$ of $\AXn=\AuX$, we call the \gui{content of~$f$} and denote by~$\rc_{\gA,\uX}(f)$ or~$\rc(f)$ the \id generated by the \coes of~$f$.

Note that we always have $\rc(f)\rc(g){\supseteq}\rc(fg)$
and thus $\rc(f)^{k+1}\rc(g){\supseteq}\rc(f)^k\rc(fg)$ for all~$k\geq0$.
For~$k$ large enough this inclusion becomes an \egtz.

\pagebreak

\CMnewtheorem{lemDKME}{\DKM lemma}{\itshape}

\begin{lemDKME}\label{lemdArtin}%
 ~\\
For $f$, $g\in\AT$ with $m\geq\deg g$ we have
 \fbox{$\rc(f)^{m+1}\rc(g)=\rc(f)^m\rc(fg)$}. 
\end{lemDKME}
\begin{proof}
First of all, notice that the products $f_ig_j$ are the \coes of the \polz~$f(Y)g(X)$. Similarly, for some \idtrs $Y_0,\ldots ,Y_m$, the content of the \polz~$f(Y_0)\cdots f(Y_m)g(X)$ is equal to~$\rc(f)^{m+1}\rc(g)$.
\\
 Let~$h=fg$.
Imagine that in the \riz~$\gB=\gA[X,Y_0,\ldots ,Y_m]$ we are able to show the membership of the \pol $f(Y_0)\cdots f(Y_m)g(X)$ in the \id
$$\preskip.3em \postskip.3em\ndsp 
\som_{j=0}^m \big(h(Y_j)\,\prod\nolimits_{k,k\neq j} \geN{f(Y_k)}\big). 
$$
We would immediately deduce that $\rc(f)^{m+1}\rc(g)\subseteq \rc(f)^m\rc(h)$.
\\
This is more or less what is going to happen.
We get rid of the \denos in Lagrange's interpolation formula (we need at least~$1+\deg g$ interpolation points):
$$\preskip-.4em \postskip.4em\ndsp 
g(X)=\sum_{j=0}^m\; \; \fraC{\prod\nolimits_{k,k\neq j}(X-Y_k)}
{\prod\nolimits_{k,k\neq j}(Y_j-Y_k)}\;\,g(Y_j)\,. 
$$
In the \ri $\gB$, by letting $\Delta=\prod_{j\neq k}(Y_j-Y_k)$, we get:
$$\preskip.4em \postskip.4em\ndsp 
\Delta \cdot g(X)\;\in \;\som_{j=0}^m \geN{g(Y_j)}. 
$$
Thus by multiplying by $f(Y_0)\cdots f(Y_m)$:
$$\preskip.4em \postskip.4em\ndsp 
\Delta \cdot f(Y_0)\cdots f(Y_m)\cdot g(X)\;\in\; \som_{j=0}^m h(Y_j)\,
\prod\nolimits_{k,k\neq j} \geN{f(Y_k)}. 
$$
If we show that for any $Q\in \gB$ we have
$\rc(Q)=\rc(\Delta \cdot Q)$, the previous membership gives 
$\rc(f)^{m+1}\rc(g)\subseteq \rc(f)^m\rc(h).$
\\
 Note that $\rc(Y_i\,Q)=\rc(Q)$ and especially that
$$ \rc\big(Q(Y_0\pm Y_1,Y_1,\ldots ,Y_m)\big)\subseteq \rc\big(Q(Y_0,Y_1,\ldots ,Y_m)\big).
$$
Therefore, by putting $Y_0=(Y_0\pm Y_1)\mp Y_1$, 
$\rc\big(Q(Y_0\pm Y_1,Y_1,\ldots,Y_m)\big)= \rc(Q).$  
 The following \pols thus all have the same content:

\snac{Q,\;Q(Y_0+ Y_1,Y_1,\ldots ,Y_m),\;Y_0\,Q(Y_0+ Y_1,Y_1,\ldots
,Y_m),\;(Y_0-Y_1)\,Q(Y_0,Y_1,\ldots ,Y_m).}

Whence $\rc(Q)=\rc(\Delta \cdot Q)$.
\end{proof}

We deduce the following corollaries.

\begin {corollary} \label {corLDM}
If $f_1$, $\ldots$, $f_d$ are $d$ \pols (with one \idtrz) of degree~$\le \delta$, with~$e_i = 1 + (d-i)\delta$ we have

\snic {
\rc(f_1)^{e_1} \rc(f_2)^{e_2} \cdots \rc(f_d)^{e_d}  \subseteq \rc(f_1f_2\cdots f_d).
}
\end {corollary}

%
\begin {proof}
Let $f = f_1$ and $g = f_2 \cdots f_d$. We have $\deg g \le (d-1)\delta$ and~$e_1 = 1 + (d-1)\delta$. \DKM lemma thus gives:

\snac {\rc(f)^{e_1} \rc(g) = \rc(f)^{(d-1)\delta} \rc(fg) \subseteq \rc(fg), $ \cad $
\rc(f_1)^{e_1} \rc(f_2 \cdots f_d) \subseteq \rc(f_1f_2\cdots f_d).
}

We finish by \recu on $d$.
\end {proof}

\vspace{-.7em}
\pagebreak

\begin{corollary}\label{corlemdArtin}  Let $f$ and $g\in \AT.$
\begin{enumerate}
\item \label{i1corlemdArtin} If $\Ann_\gA\big(\rc(f)\big)=0$, then $\Ann_\AT(f)=0$ \emph{(McCoy's Lemma)}.\index{McCoy's Lemma}%
\item \label{i2corlemdArtin} If $\gA$ is reduced, then $\Ann_\AT(f)= \Ann_\gA\big(\rc(f)\big)[T]$.
\item \label{i3corlemdArtin} The \polz~$f$ is nilpotent \ssi each of its \coes is nilpotent.
\item If $\rc(f)=1$, then $\rc(fg)=\rc(g)$.
\end{enumerate}
\end{corollary}

\begin{proof}
Let $g\in \Ann_\AT(f)$ and $m \ge \deg(g)$. \DKM lemma implies:
$$\preskip.4em \postskip.3em 
\rc(f)^{1+m}g=0.\eqno(*) 
$$
\emph{1.}
So $\Ann_\gA \rc(f)=0$ implies $g=0$.
\\
\emph{2.}
Since the \ri is reduced, $(*)$ implies $\rc(f)g=0$. Thus every \polz~$g$ annihilated by $f$ is annihilated by~$\rc(f)$.\\
Furthermore, $\Ann_\gA\big(\rc(f)\big)=\gA\cap \Ann_\AT(f)$
and thus the inclusion

\snic{\Ann_\AT(f)\supseteq \Ann_\gA\big(\rc(f)\big)[T]}

is always true (whether $\gA$ is reduced or not).\\
\emph{3.} If $f^2=0$, the \DKM lemma implies $\rc(f)^{2+\deg f}=0$.\\
\emph{4.} Immediately from $\rc(f)^{m+1}\rc(g)=\rc(f)^m\rc(fg)$.
\end{proof}
%

\section{One of Kronecker's \thosz}
\label{secThKro}

\vspace{3pt}
\subsec{\Algs and integral \eltsz}

We first introduce the terminology of \Algsz. The \algs that we consider in this work are associative, commutative and unitary, unless stated otherwise.

\begin{definition}
\label{def0Alg} ~
\begin{enumerate}
\item An \emph{\Algz} is a commutative \ri $\gB$ with a \homo of commutative \ris $\rho:\gA\to\gB$. That makes~$\gB$ an \Amoz. When $\gA\subseteq\gB$, or more \gnlt if~$\rho$ is injective, we say that~$\gB$ is an \ix{extension} of~$\gA$.%
\index{algebra!over a ring}

\item A \emph{morphism} of the \Alg $\gA\vers{\rho}\gB$ to \hbox{the \Alg $\gA\vers{\rho'}\gB'$} is a \homo of \risz~$\gB\vers{\varphi}\gB'$ satisfying $\varphi\circ\rho=\rho'$. The set of morphisms of \Algs of~$\gB$ to~$\gB'$ is denoted by~$\Hom_\gA(\gB,\gB')$.
\vspace{-.5em}
\tri{\gA}{\rho}{\rho'}{\gB}{\varphi}{\gB'}{}
\end{enumerate}
\end{definition}

\vspace{.1em}
\rems ~\\
1) We chose not to reserve the terminology \gui{extension} for the case of fields. This will require us to use in the cases of fields statements such as \gui{$\gL$ is a field extension of $\gK$} or \gui{$\gL$ is a field, extending~$\gK$}
from this point on.

2) Every \ri is uniquely a \ZZlg and every \homo of \ris is a morphism of the corresponding \ZZlgsz. The category of commutative \ris can be regarded as a special case among the categories of \algs defined above.
\eoe \rdb

\mni
{\bf Notation.} If $b\in\gB$ and $M$ is a \Bmoz, we denote by $\mu_{M,b}$ or $\mu_b$ the multiplication by~$b$ in~$M$: $y\mapsto by,\, M\to M$. This can be regarded as a \Bliz, or, if~$\gB$ is an \Algz, as an \Ali for the \Amo structure of~$M$. \label{NOTAmux}

\begin{definition}
\label{defEntierAnn0} Let  $\gA\subseteq\gB$ be \risz.
\begin{enumerate}
\item An \elt $x\in\gB$ is said to be \ixc{integral}{element over a ring}
over~$\gA$ if there exists some integer~$k\geq 1$ such \hbox{that $x^k=a_1x^{k-1}+a_2 x^{k-2}+\cdots +a_k$}
with each~$a_h\in\gA$. If~$\gA$ is a \cdiz, we also say that~$x$ is \ixc{algebraic}{element over a discrete field} over~$\gA$.
\item In this case, the \polu $P=X^k-(a_1X^{k-1}+a_2 X^{k-2}+\cdots +a_k)$ is called an \emph{\rdiz} of $x$ over~$\gA$.
In fact, by abuse of language we also say that the \egt $P(x)=0$ is an \emph{\rdiz}. If~$\gA$ is a \cdiz,~we also speak of an \emph{\agq \rdez}.%
\index{dependence relation!integral ---}%
\index{dependence relation!algebraic ---}
\item The \riz~$\gB$ is said to be
 \emph{integral} over $\gA$ if every \elt of~$\gB$ is integral
over~$\gA$. We will also say that the \Alg $\gB$ is \emph{integral}.
If~$\gA$ and~$\gB$ are \cdisz,~we say that ~$\gB$ is \emph{\agqz} over~$\gA$.%
\index{algebra!integral ---}%
\index{ring!integral --- over a subring}%
\index{integral!ring over a subring}%
\index{algebraic!field over a subfield}
\item
If $\rho:\gC\to\gB$ is a \Clg with~$\rho(\gC)=\gA$, we say that the \algz~$\gB$ is integral over~$\gC$ if it is integral over~$\gA$.
\end{enumerate}
\end{definition}

\subsec{The \tho}

\begin{theorem}\label{thKro} \emph{(\KROz's \Thoz) \cite{Kro1}}\\
In $\gB[T]$, consider the \pols
$$\preskip.3em \postskip.3em \mathrigid 2mu
f=\sum_{i=0}^n(-1)^if_iT^{n-i},\;  g=\sum_{j=0}^m(-1)^j g_j T^{m-j} \;\, \mathit{ and } \;\, h=fg=\sum_{r=0}^{p}(-1)^r h_rT^{p-r},
$$
where $p=m+n$. Let $\gA=\gZ[h_0,\ldots,h_p]$ be the sub\ri generated by the \coes of $h$ ($\gZ$ is the sub\ri of $\gB$ generated by $1_\gB$). 
\begin{enumerate}
\item Each $f_ig_j$ is integral over $\gA$.
\item In the case where we take \idtrs over the \riz~$\ZZ$ for $f_i$ and $g_j$, we find an \rdi over $\gA$ for $z_{i,j}=f_ig_j$ which is \hmg for different \syss of weights attributed to the \momsz:
\begin{enumerate}
\item the respective weights of $z_{k,\ell}$ and $h_r$ are $k+\ell$ and $r$.
\item the respective weights of $z_{k,\ell}$ and $h_r$ are $p-k-\ell$ and $p-r$.
\item the weights of $z_{k,\ell}$ and $h_r$ are $w(z_{k,\ell})=w(h_r)=1$.
\end{enumerate}
Naturally these \rdis are then applicable in every \riz.
\end{enumerate}
\end{theorem}
\begin{proof}
It suffices to treat item \emph{2.}
\\
 Let us first examine an intermediate \gnq case.
We take $f_0=g_0=1$ and \idtrs over $\ZZ$ for the other $f_i$'s and $g_j$'s. The \pols $f$ and $g$ are thus \polus in $\gB[T]$ with $\gB=\ZZ[f_1,\ldots,f_n,g_1,\ldots,g_m]$, and $\gA=\ZZ[h_1,\ldots,h_p]$.\\
Assume \spdg \hbox{that $\gB\subseteq\gC=\ZZ[\xn,\alb\ym]$}, where each $x_i$ and $y_j=x_{n+j}$ are \idtrsz, the $f_i$'s are the \elr \smq \pols in the $x_i$'s, and the $g_j$'s are the \elr \smq \pols in the $y_j$'s (apply Corollary~\ref{corthSymE1} twice).
If we attribute to $x_i$ and $y_j$ a weight of 1, the $z_{k,\ell}$ and $h_r$ are \hmgs and we obtain the weights described in \emph{2a}.
To compute an \rdi for $f_ig_j$ (with eventually $i$ or $j=0$) over~$\gA$, consider the subgroup $H_{i,j}$ \hbox{of $\rS_p$} formed by the $\sigma$'s which satisfy $\sigma(f_ig_j)=f_ig_j$
(this subgroup contains at least all the permutations which stabilises $\lrbn$).
Then consider the \pol
$$\preskip.4em \postskip.3em
P_{i,j}(T)=\prod\nolimits_{\tau\in\rS_p/H_{i,j}}\big(T-\tau(f_ig_j)\big),\eqno(*)
$$
where $\tau\in\rS_p/H_{i,j}$ means that we take exactly one $\tau$ from each left coset of~$H_{i,j}$. Then, $P_{i,j}$ is \hmg for the weights $w_a$ described in~\emph{2a} ($i,j$ being fixed, we denote by $w_a$ the weights \emph{2a}, with $w_a(T)=w_a(z_{i,j})$).
Moreover,~$P_{i,j}$ is \smq in the $x_k$'s ($k\in\lrbp$). It is therefore uniquely expressible as a \polz~$Q_{i,j}(\uh,T)$ in each $h_r$ and $T$, and $Q_{i,j}$ \hbox{is $w_a$-\hmg} (\Thref{thSymEl} items \emph{1} and \emph{2a}).
The degree in $T$ of~$Q_{i,j}$ \hbox{is $d_{i,j}=(\rS_p:H_{i,j})$}.
For $R\in\gC[T]$, denote by $\delta(R)$ the integer $\deg_{x_1}(R)+\deg_T(R)$. We see that~$\delta$ is a weight, and that \hbox{$\delta(f_ig_j)=w(f_ig_j)\leq1$}, $\delta(h_r)=w(h_r)\leq1$ (\hbox{with $w(h_r)=1$} if $i$, $j$, $r\geq1$). 
Moreover, each factor of $P_{i,j}$ in~$(*)$ is of weight~$1$ (but not necessarily \hmg because we can have $\delta(\sigma\big(f_ig_j)\big)=0$).
This gives $\delta(Q_{i,j})=d_{i,j}$ when the \pol is evaluated in $\gC[T]$. Moreover, by  \Thref{thSymEl} item~\emph{2b}, when we write a \smq \pol in $(\xp)$, say $S(\ux)$, as a \polz~$S_1(\uh)$ in the $h_i$'s, \hbox{we have $\delta(S)=w(S_1)$}. Thus $w(Q_{i,j})=d_{i,j}$.
\\
To treat item~\emph{2} itself it suffices to \gui{homogenize.}
Indeed, if we let $\wi f_i=f_i/f_0$ and $\wi g_j=g_j/g_0$, which is legitimized by Corollary~\ref{corthSymE2}, for $\wi f_i$ and $\wi g_j$ we return to the previous situation with regard to the weights \emph{2a}. We obtain a \hmg \rdi for $\wi z_{i,j}=\wi f_i\wi g_j$ over the sub\ri generated by the~$\wi h_{r}$
$$\preskip.0em \postskip.4em
Q_{i,j}(\wi h_{1},\ldots,\wi h_{p},\wi z_{i,j})=0,
$$
with $\wi z_{i,j}=f_ig_j/h_0$ and $\wi h_{r}=h_{r}/h_0$.
\\
We multiply the \ida obtained by $h_0^{d_{i,j}}$ so that we obtain a \polu in $z_{i,j}$.
\\
All the \denos have vanished because $w(Q_{i,j})=d_{i,j}$. We obtain
$$\preskip.4em \postskip.4em
R_{i,j}(h_0,\ldots,h_p,f_ig_j)=0,
$$
where $R_{i,j}(h_0,\ldots,h_p,T)$ is unitary in $T$ and homogeneous for the weights $w_a$ and~$w$.
\\
What remains is the question of the homogeneity for the weights $w_b$ in \emph{2b}:
it suffices to note that we have for all $R\in\AT$ the \egt
%
$w_a(R)+w_b(R)=pw(R).$
\end{proof}
\exl In the case where $m=n=2$, the indicated computation gives the following results.
\\
When $f_0=g_0=1$ the \coe $g_1$ annihilates the \polz
$$\arraycolsep2pt\begin{array}{rcl}
p_{01}(t)&=& t^6 - 3 h_1 t^5 + (3 h_1^2 + 2 h_2) \,t^4 + (-h_1^3 - 4 h_1 h_2) \,t^3 +\\[1mm]
&&(2 h_1^2 h_2 + h_1 h_3 + h_2^2 - 4 h_4) \,t^2 + (-h_1^2 h_3 - h_1 h_2^2 + 4 h_1 h_4) \,t
\\[1mm]
&&-h_1^2 h_4 + h_1 h_2 h_3 - h_3^2\,,
\end{array}$$
so in the \gnl case $f_0g_1$ annihilates the \polz
$$\arraycolsep2pt\begin{array}{rcl}
q_{01}(t)&=& t^6 - 3 h_1 t^5 + (3 h_1^2 + 2h_0 h_2) \,t^4 + (-h_1^3 - 4 h_0 h_1 h_2) \,t^3 +\\[1mm]
&&(2 h_0 h_1^2 h_2 + h_0^2 h_1 h_3 +h_0^2  h_2^2 - 4h_0^3 h_4) \,t^2 +
\\[1mm]
&&(-h_0^2 h_1^2 h_3 - h_0^2 h_1 h_2^2 + 4 h_0^3h_1 h_4) \,t
-h_0^3h_1^2 h_4 + h_0^3h_1 h_2 h_3 - h_0^4h_3^2\,.
\end{array}$$
When $f_0=g_0=1$ the \coe $g_2$ annihilates the \polz
$$\arraycolsep2pt\begin{array}{rcl}
p_{02}(t)      & =& t^6 - h_2 t^5 + (h_1 h_3 - h_4) \,t^4 + (-h_1^2 h_4 + 2 h_2 h_4 - h_3^2) \,t^3 +  \\[1mm]
      &  & (h_1 h_3 h_4 - h_4^2) \,t^2 - h_2 h_4^2 \,t + h_4^3
\,,\end{array}$$
so $f_0g_2$ annihilates the \pol
$$\arraycolsep2pt\begin{array}{rcl}
q_{02}(t)      & =& t^6 - h_2 t^5 + (h_1 h_3 - h_0 h_4) \,t^4 + (-h_1^2 h_4 + 2 h_0 h_2 h_4 - h_0 h_3^2) \,t^3 +  \\[1mm]
      &  & (h_0 h_1 h_3 h_4 - h_0^2 h_4^2) \,t^2 - h_0^2 h_2 h_4^2 \,t + h_0^3 h_4^3\,.
\end{array}$$
When $f_0=g_0=1$ the \coe $f_1  g_1$ annihilates the \pol
$$\arraycolsep2pt\begin{array}{rcl}
p_{11}(t)&=&t^3 - 2 h_2 t^2 + (h_1 h_3 + h_2^2 - 4 h_4) \,t +
      h_1^2 h_4 - h_1 h_2 h_3 +
h_3^2\,.
\end{array}
$$
When $f_0=g_0=1$ the \coe $f_1  g_2$ annihilates the \pol
$$\arraycolsep2pt\begin{array}{rcl}
p_{12}(t)&=&t^6 - 3 h_3 t^5 + (2 h_2 h_4 + 3 h_3^2) \,t^4 + (-4 h_2 h_3 h_4 - h_3^3) \,t^3
+ \\[1mm]
      &  &(h_1 h_3 h_4^2 + h_2^2 h_4^2 + 2 h_2 h_3^2 h_4 - 4 h_4^3) \,t^2 + \\[1mm]
      &  &(-h_1 h_3^2 h_4^2 -
h_2^2 h_3 h_4^2 + 4 h_3 h_4^3) \,t - h_1^2 h_4^4 + h_1 h_2 h_3 h_4^3 - h_3^2 h_4^3\,.
\end{array}$$

\pagebreak

\begin{corollary}\label{corthKro}
 \emph{(Multivariate \KROz's \Thoz) }\\
In $\gB[\Xk]$ consider the \pols

\snic{f=\som_{\alpha}f_\alpha X^{\alpha},\;\;  g=\som_{\beta}b_\beta X^{\beta} \;\;\hbox{ and }\;\; h=fg=\som_{\gamma}h_\gamma X^{\gamma},}

(here, $\alpha$, $\beta$, $\gamma$ are multi-indices, and if $\alpha=(\alpha_1,\ldots,\alpha_k)$,
$X^{\alpha}$ is a notation for $X_1^{\alpha_1}\cdots X_k^{\alpha_k}$).
Let $\gA=\gZ[(h_\gamma)]$ be the sub\ri generated by the \coes of $h$ ($\gZ$ is the sub\ri of $\gB$ generated by $1_\gB$). Then, each $f_\alpha g_\beta$ is integral over $\gA$.
\end{corollary}
\begin{proof}
We apply what is termed \emph{\KRAz's trick}: let~$X_j=T^{n^j}$ with large enough $n$. This transforms $f$, $g$ and $h$ into \pols $F(T)$, $G(T)$, $H(T)$ whose \coes are those of~$f$,~$g$ and~$h$, respectively.
\end{proof}
%


\section[The \adu (1)]{The \adu  for a \polu over a commutative \ri (1)} \label{sec0adu}

\emph{Disclaimer.}
In a context where we manipulate \algsz, it is sometimes preferable to keep to the intuition that we want to have a field as the base ring, 
even if it is only a commutative \riz. In which case we choose to give a name such as $\gk$ to the base \riz.
This is what we are going to do in this section dedicated to the \aduz.
\\
When we are truly dealing with a discrete field, we will use $\gK$ instead.

We now proceed to the inverse operation to that which passes from the \pol \ri to the sub\ri of \smq \polsz.
\\
In the presence of a \polu $f=T^n+\sum_{k=1}^n(-1)^{k} s_kT^{n-k} \in \kT$ over a \ri $\gk$, we want to have at our disposal an extension of $\gk$ where the \pol is decomposed into \lin factors.
Such an extension can be constructed in a purely formal way. The result is called the \aduz.

\begin{definota}\label{definotaAdu}
Let $f=T^n+\sum_{k=1}^n(-1)^{k} s_kT^{n-k} \in \kT$ be a \polu of degree~$n$.
We denote by $\Adu_{\gk,f}$ the \emph{\adu of $f$ over $\gk$} defined as follows%
\index{universal splitting algebra!of $f$ over $\gk$}%
\index{ideal!of the symmetric relators}%
\index{relator!symmetric ---}
$$\preskip.4em \postskip.4em
\Adu_{\gk,f}=\kXn/\cJ(f)=\kxn
,
$$
where $\cJ(f)$ is the \emph{\id of \smq relators} \ncr to 
identify $\prod_{i=1}^n(T-x_i)$ with $f(T)$  in the quotient.
Precisely, if  $S_1$, $ S_2$, \ldots, $S_n$ are the \elr \smq functions of the~$X_i$'s,
 the \id $\cJ(f)$ is given by
$$\preskip.4em \postskip.4em
\cJ(f)=\gen{S_1-s_1,\;S_2-s_2,\ldots ,\;S_n-s_n}.$$
\end{definota}

The \adu $\gA=\Adu_{\gk,f}$ can be \caree by the following  \prtz.

\begin{fact}
\label{factEvident} \emph{(\Uvl \dcn \algz, \cara \prtz)}
\begin{enumerate}
\item Let $\gC$ be a \klg such that $f(T)$ is decomposed into a product of factors $T-z_i$. Then, there exists a unique \homo of \klgs \hbox{of $\gA$} \hbox{to $\gC$} which sends the $x_i$'s to the $z_i$'s.
\item This characterizes the \adu $\gA=\Adu_{\gk,f}$, up to unique \isoz.
\item Moreover, if $\gC$ is generated (as a \klgz) by the $z_i$'s, the \adu is \isoc to a quotient of $\gA$%
.
\end{enumerate}
\end{fact}
%
\begin{proof}
For item \emph{1} we use Proposition~\ref{propApolLibre}, which describes the \algs of \pols
as \algs freely generated by the \idtrsz, and Fact~\ref{factUnivQuot}, which describes the quotient \ris as those which allow us to uniquely factorize certain \homosz. Item \emph{2} results from the ascertainment that an object that solves a \uvl \pb is alway unique up to unique \isoz.
\end{proof}

By taking $\gC=\gA$ we obtain that every permutation of $\so{1,\ldots,n}$ produces a (unique) $\gk$-\auto of $\gA$.
\\
Stated otherwise: the group $\Sn$ of permutations of $\so{\Xn}$ acts \hbox{on $\kXn$} and fixes the \id $\cJ(f)$, thus the action passes to the quotient and this defines $\Sn$ as a group of \autos of the \aduz.

  To study the \adu we introduce \ix{Cauchy modules} which are the following \polsz:
{\small
$$\preskip.4em \postskip.6em
{\arraycolsep2pt\begin{array}{rcl}
f_1(X_1)& =  & f(X_1)   \\
f_2(X_1,X_2)& =  &  \big(f_1(X_1)-f_1(X_2)\big)\big/(X_1-X_2)    \\
 & \vdots  &  \\
f_{k+1}(X_1,\ldots ,X_{k+1})& =  &  {\dsp \frac{f_k(X_1,\ldots ,X_{k-1},X_k)-
f_k(X_1,\ldots ,X_{k-1},X_{k+1})}{X_k-X_{k+1}}}    \\
 & \vdots  &       \\
f_n(X_1,\ldots ,X_n)& =  &{\dsp \frac{f_{n-1}(X_1,\ldots ,X_{n-2},X_{n-1})-f_{n-
1}(X_1,\ldots ,X_{n-2},X_{n})}{X_{n-1}-X_{n}}} .    \\
\end{array}}$$
}

The following fact results from the \cara \prt of the \adusz.

\begin{fact}\label{factAduAdu}
With the previous notations for the Cauchy modules, \hbox{let $\gk_1=\gk[x_1]$} and $g_2(T)=f_2(x_1,T)$.
Then, the canonical $\gk_1$-\lin map~\hbox{$\Adu_{\gk,f}\to\Adu_{\gk_1,g_2}$} (which sends each $x_i$ ($i\geq2$) of $\Adu_{\gk,f}$ to the~$x_i$'s \hbox{of $\Adu_{\gk_1,g_2}$}) is an \isoz.
\end{fact}

\mni \exls (Cauchy modules)\\
 With $n=4$,
{\small
$$\arraycolsep2pt\begin{array}{rcl}
f_1(x)& =  & x^4-s_1x^3+s_2x^2-s_3x+s_4   \\
f_2(x,y)& =  &  (y^{3}+ y^{2}x+ yx^{2} + x^{3})- s_1(y^{2}+yx+x^2)  +
s_2(y+x)  - s_3 \\
&=&y^{3} +y^{2}(x- s_1) + y(x^{2}- s_1x+ s_2) + (x^{3}- s_1x^{2}+ s_2x-
s_3)  \\
f_{3}(x,y,z)& =  &     (z^2+ y^2+ x^2+ zy + zx+ yx) - s_1(z+y+x)   + s_2
\\
&=& z^{2} +z (y+x- s_1)   + \big((y^{2} + yx + x^{2})- s_1(y+ x)  + s_2\big)
\\
f_4(x,y,z,t)& =  &t+z+y+x-s_1 .
\end{array}$$
}
For $f(T)=T^6$,

\snac{
\arraycolsep2pt\begin{array}{rcl}
f_2(x,y)& =  &  y^{5}+ y^{4}x+ y^3x^{2} + y^2x^{3}+yx^4+x^5 \\
f_{3}(x,y,z)& =  &   (z^{4}+ y^{4}+ x^{4})+ ( z^{2}y^{2} + z^{2} x^{2}+
y^{2}x^{2}) +  \\
&&( zy^{3} + zx^{3}+ y z^{3}+  y x^{3}+ x z^{3} +x y^{3})  +\\
&& ( zy x^{2}  + zx  y^{2} + yx  z^{2}    )             \\
f_4(x,y,z,t)& =  & (t^{3}+ z^{3}+ y^{3} +x^{3}) + (tzy+tyx+tzx+zyx)+ \\
&&t^{2}(z+y+x)+ z^{2}(t+y+x)+ \\
&&y^{2}(t+z+x)+ x^{2} (t+z+y)       \\
f_5(x,y,z,t,u)&=&  (u^{2} + t^{2}+ z^{2}+ y^{2} +x^{2}) +\\
&& (x u + x t + x z +
x y + t u + z u + z t + y u + y t + y z)\\
f_6(x,y,z,t,u,v)&=& v+u+t+z+y+x .
\end{array}}

More \gnltz, for $f(T)=T^n$, $f_k(t_1,\ldots ,t_k)$ is the sum of all the \moms of degree $n+1-k$ in $t_1,\ldots ,t_k$.\\
By linearity, this allows us
to obtain an explicit, precise description of the Cauchy modules for an arbitrary \polz.
\eoe

\medskip
By the remark following the last example, the \polz~$f_i$ is \smq in the variables $X_1$, $\ldots$, $ X_i$, \mon in $X_i$, of total degree $n-i+1$.

Fact~\ref{factEvident} implies that the \id $\cJ(f)$ is equal to the \id generated by the Cauchy modules. Indeed, the quotient \ri by the latter \id clearly realizes the same \uvl \prt as the quotient \ri by $\cJ(f)$.

Thus the \adu is a free \kmo of rank~$n!$. More \prmtz, we obtain the following result.

\begin{fact}
\label{factBase}
The \kmo $\gA=\Adu_{\gk,f}$ is free and a basis is formed by the \gui{\momsz}
$x_1^{d_1}\cdots x_{n-1}^{d_{n-1}}$ such that for $k=1,\ldots ,n-1$ we have $d_k\leq n-k$.
\end{fact}

\begin{corollary}\label{corfactBase}
Considering the \adu of the \gnq \polu $f(T)=T^n+\sum_{k=1}^n(-1)^kS_kT^{n-k}$, where the $S_i$'s are \idtrsz, we get an \alg of \pols $\kxn$ with each $S_i$ identifiable with \elr \smq \pols in each~$x_i$.
\end{corollary}

\comm (For those who know Gr\"obner bases)\\
In the case where $\gk$ is a \cdiz, the Cauchy modules can be seen as a Gr\"obner basis of the \id $\cJ(f)$, for the lexicographic \mom order with $X_1<X_2<\cdots< X_n$.
\\
In fact, even if $\gk$ is not a \cdiz, the Cauchy modules still work as a Gr\"obner basis: Every \pol in the $x_i$'s can be re-expressed over the previous \mom basis by successive divisions by the Cauchy modules. 
We first divide by $f_n$ with respect to the variable~$X_n$, which cancels it out.
Next we divide by $f_{n-1}$ with respect to the variable $X_{n-1}$, which brings it to a degree $\leq 1$, and so on.
\eoe

\penalty-2500
\section{Discriminant, \dinz}
\label{secDisc}
\vspace{4pt}
\subsec{Definition of the \discri of a \polu}

We define the \ix{discriminant} of a univariate \polu $f$ over a commutative \ri $\gA$ starting with the case where $f$ is the \gnq \polu of degree $n$:
$$
\mathrigid1mu f(T)\,=\,T^n-S_1T^{n-1}+S_2T^{n-2}+\cdots+(-1)^nS_n\;\in\,\ZZ[S_1,\ldots,S_n][T]\,=\,\ZZuS[T].
$$
We can write $f(T)=\prod_i(T-X_i)$ in $\ZZXn$  
(Corollary~\ref{corthSymE1}), and we set
\begin{equation}\preskip.0em \postskip.4em
\label{eqDiscri}
\disc_T(f)=  (-1)^{n(n-1)/2}\, \prod\nolimits_{i=1}^n f'(X_i)
=\prod\nolimits_{1\leq i< j\leq n}(X_i-X_j)^2.
\end{equation}
As this \pol in the $X_i$'s is clearly variable permutation invariant, there exists a unique \pol in the $S_i$'s, $D_n(S_1,\ldots,S_n)\in\ZZuS$,
which is equal to $\disc_T(f)$. In short, the auxiliary variables $X_i$ can indeed vanish.

Then, for a \gui{concrete} \pol
$$\preskip.3em \postskip.2em
g(T)=T^n-s_1T^{n-1}+s_2T^{n-2}+\cdots+(-1)^ns_n\in\AT,
$$
 we define $\disc_T(g)=D_n(s_1,\ldots,s_n)$.

Naturally, if it happens that $g(T)=\prod_{i=1}^n(T-b_i)$ in a \ri $\gB\supseteq\gA$, we would then obtain $\disc_T(g)=\prod\nolimits_{1\leq i< j\leq n}(b_i-b_j)^2$ by evaluating the formula~(\ref{eqDiscri}).
In particular, by using the \adu we could directly define the \discri by this formula.

A \polu is said to be
 \emph{separable} when its \discri is \ivz.%
\index{separable!monic polynomial}%
\index{polynomial!separable monic ---}

\subsec{Diagonalization of the matrices on a \riz}

Let us first recall that if $f\in\AT$, a \emph{zero of $f$ in an \Alg $\gB$} (given by a \homo $\varphi:\gA\to\gB$) is a $y\in\gB$ which annihilates the \polz~$f^{\varphi}$, the image of $f$ in $\gB[T]$.

In addition, the zero $y$ is said to be
 \emph{simple} if $f'(y)\in\gB\eti$ (we also say that it is a \emph{simple root} of $f$).%
\index{simple!zero}%
\index{simple!root}%
\index{root!simple ---}%
\index{zero!simple --- of a polynomial}%
\index{zero!of a polynomial}

Here, we are interested in the \dins of matrices on an arbitrary commutative \riz, when the \polcar is \emph{\splz}.

First of all, we have the classical \gui{Kernels' Lemma}~\ref{lemDesNoyaux}.\\
Next is a \gnn of the \tho which states (in the \cdi case) that a simple zero of the \polcar defines a proper subspace of dimension $1$.

\begin{lemma}\label{lemValProp}%
Let $n\geq2$, $a\in\gA$ and $A\in\Mn(\gA)$ be a matrix whose \polcar $f(X)=\rC A(X)$ admits $a$ as simple zero. Let   \hbox{$g=f/(X-a)$}, $h=X-a$, $K=\Ker h(A)$ and $I=\Im h(A)$.
\begin{enumerate}
\item We have $K=\Im g(A)$, $I=\Ker g(A)$ and $\Ae n=I\oplus K$.
\item The matrix $g(A)$ is of rank $1$, and $h(A)$ of rank $n-1$.
\item If a \polz~$R(X)$ annihilates $A$,
then $R(a)=0$, \cad $R$ is a multiple of $X-a$.
\item The principal minors of order $n-1$ of $A-a\In$ are \comz.
We  localize by inverting such a minor, the matrix $g(A)$ becomes simple of rank $1$, the modules $I$ and $K$ become free of rank $n-1$ and $1$.
\end{enumerate}
\end{lemma}
\begin{proof}
Suppose \spdg that $a=0$.
\\
Then, $f(X)=Xg(X)$, $h(A)=A$, $g(A)=\pm\wi A$, $\Tr\big(g(A)\big)=g(0)$ (Lemma~\ref{lemPrincipeIdentitesAlgebriques} item \emph{6}), and $g(0)=f'(0)\in\Ati$.

 \emph{1.}
 We write $g(X)=Xk(X)+g(0)$. This shows that the \pols $g(X)$ and $X$ are \comz.
Given the Cayley-Hamilton \Thoz, the Kernels' Lemma 
applies and gives item~\emph{1.}

\emph{2.} Let $\mu_1$, $\ldots$, $\mu_n$ be the principal minors of order $n-1$ of $A$. \\
Since $g(A)=\pm\wi A$, we get $g(0)=\Tr \big(g(A)\big)=\pm \Tr\wi A =\pm\sum_i\mu_i$. This shows that $\rg\big(h(A)\big)=n-1$ and $\rg\big(g(A)\big)\geq1$.
Finally, we know that $\rg(\wi A)\leq1$ by Lemma~\ref{lemPrincipeIdentitesAlgebriques} item~\emph{8.}.

\emph{3.} Suppose $R(A)=0$.
By multiplying by $\wi{A}$, we obtain $R(0)\wi{A} = 0$ (since $\wi{A} A = 0$). By taking the trace, $R(0)\Tr(\wi{A}) = 0$ thus $R(0) = 0$.
\\
Note that item \emph{3} also results from item \emph{4.}

\emph{4.} We have already seen that the $\mu_i$'s are \comz.
After \lon at some $\mu_i$, the matrix $g(A)$ becomes simple of rank $1$ under the freeness lemma~\paref{lem pf libre}. Therefore $I$ and $K$ become free of rank $n-1$ and $1$.
\end{proof}
%

\begin{proposition}\label{propDiag} \emph{(Diagonalization of a matrix whose \polcar is \splz)}  
Let $A\in\Mn(\gA)$ be a matrix whose \polcar $\rC A(X)$ is \splz, and $\gA_1 \supseteq\gA$ be a \ri on which we can write $\rC A(X)=\prod_{i=1}^n(X-x_i)$ (for example, $\gA_1=\Adu_{\gA,f}$). \\
Let $K_i=\Ker(A-x_i\In)\subseteq \gA_1^n$.
\begin{enumerate}
\item $\gA_1^n=\bigoplus_i K_i$.
\item Each $K_i$ is the image of a matrix of rank $1$.
\item Every \polz~$R$ which annihilates $A$ is a multiple of $\rC A$.
\item After \lon at  \eco of $\gA_1$ the matrix is \digz, similar to $\Diag(\xn)$.
\end{enumerate}
\end{proposition}
NB: if $\alpha\in\End_{\gA_1}(\gA_1^n)$ has for matrix $A$, we have $\alpha\frt{K_i}=x_i\,\Id_{K_i}$ for each~$i$.
\begin{proof}
This is an immediate consequence of the Kernels' Lemma and Lemma~\ref{lemValProp}. 
To render the matrix \dig it suffices to invert some product $\nu_1\cdots\nu_n$
where each $\nu_i$ is a principal minor of order $n-1$ of the matrix $A-x_i\In$ (which a priori makes $n^n$ \come \lonsz).
\end{proof}

\rem An analogous result concerning a matrix that annihilates a \spl \polz~$\prod_{i}(X-x_i)$ is given in Exercise~\ref{exo2Diag}. The proof is \elrz.
\eoe

\subsec{The \gnq matrix is \dig}

Consider $n^2$ \idtrs $(a_{i,j})_{i,j\in \lrbn}$ and let $A$ be the corresponding matrix (it has \coes in $\gA=\ZZ[(a_{i,j})]$).

\begin{proposition}\label{propMatGenDiag}
The \gnq matrix $A$ is \dig over a \riz~$\gB$ containing $\ZZ[(a_{i,j})]=\gA$.
\end{proposition}
\begin{proof}
Let
$f(T)=T^n-s_1T^{n-1}+\cdots+(-1)^ns_n$ be the \polcar of $A$.
Then the \coes $s_i$ are \agqt independent over $\ZZ$.  
To realize this, it suffices to specialize $A$ as the companion matrix of a \gnq \poluz.
\\
In particular, the \discri $\Delta=\disc(f)$ is nonzero in the integral \riz~$\gA$.
Then consider the \riz~$\gA_1=\gA[1/\Delta]\supseteq\gA$ and the \adu $\gC=\Adu_{\gA_1,f}$. Let the $x_i$ be the \elts of $\gC$ such that $f(T)=\prod_i(T-x_i)$.
\\
Finally, apply Proposition~\ref{propDiag}. If we want to obtain a \dig matrix, we invert for instance $a=\prod_i\det\big((A-x_i\In)_{1..n-1,1..n-1}\big)$. This is an \elt of $\gA$ and it suffices to convince ourselves that it is nonzero by exhibiting a particular matrix, for example the companion matrix of the \pol $X^n-1$.
\\
Ultimately, consider $\gA_2=\gA[1/(a\Delta)]\supseteq\gA$ and take
${\gB=\Adu_{\gA_2,f}\supseteq\gA_2.}$
\end{proof}

The strength of the previous result, \gui{which makes life considerably easier} is illustrated in the following two subsections.

\subsec{An identity concerning \polcarsz}

\begin{proposition}\label{prop1tschir}
Let $A$ and $B\in \Mn(\gA)$ be two matrices which have the same \polcarz, and let $g\in\AT$. Then the matrices $g(A)$ and $g(B)$ have the same \polcarz.
\end{proposition}

\pagebreak

\begin{corollary}\label{corprop1tschir}\label{lemPolCar}~
\begin{enumerate}
\item If $A$ is a matrix with \polcar $f$, and if we can write $f(T)=\prod_{i=1}^n(T-x_i)$ on a \ri $\gA_1\supseteq \gA$, then the \polcar of $g(A)$ is equal to the product $\prod_{i=1}^n\big(T-g(x_i)\big)$.
\item Let $\gB$ be a free \Alg of finite rank $n$ and $x\in\gB$.
Suppose that in $\gB_1\supseteq \gB$, we have $\rC{\gB/\!\gA}(x)(T) = \prod_{i=1}^n(T-x_i)$. Then, for all $g \in \gA[T]$, we have the following \egtsz.
$$\preskip.0em \postskip.2em 
\gB/\!\gA\big(g(x)\big)(T) = \prod\nolimits_{i=1}^n\big(T-g(x_i)\big), 
$$
$\Tr\iBA \big(g(x)\big) = \som_{i=1}^n g(x_i)$   and 
$\rN\iBA \big(g(x)\big) = \prod_{i=1}^n g(x_i)$.
\end{enumerate}
\end{corollary}
\begin{Proof}{\Demo of the proposition and the corollary.} \\
Item \emph{1} of the corollary. Consider the matrix $\Diag(\xn)$ which has the same \polcar as  $A$ and apply the proposition with the \riz~$\gA_1$.
\\
Conversely, if item \emph{1} of the corollary is proven for $\gA_1=\Adu_{\gA,f}$, it implies Proposition~\ref{prop1tschir} since the \pol $\prod_{i=1}^n\big(T-g(x_i)\big)$ computed \hbox{in $\Adu_{\gA,f}$} can only depend on $f$ and $g$.\\
Now note that the structure of the statement of the corollary, item~\emph{1}, when we take $\gA_1=\Adu_{\gA,f}$, is a family of \idas with the \coes of the matrix~$A$ for \idtrsz.
It thus suffices to prove it for the \gnq matrix. However, it is \dig over some over\ri (Proposition~\ref{propMatGenDiag}), and for some \dig matrix the result is clear.
\\
 Finally, item \emph{2} of the corollary is an \imde consequence of item~\emph{1.}
\end{Proof}
%

\subsec{An identity concerning exterior powers}

The following results, analogous to Proposition~\ref{prop1tschir} and to Corollary~\ref{corprop1tschir}, can be proven by following the exact same proof sketch.

\begin{proposition}\label{propPolCarPuissExt}
If $\varphi$ is an \endo of a free \Amo of finite rank, the \polcar of ${\Al{k}\!\varphi}$ only depends on the integer $k$ and on the \polcar of $\varphi$.
\end{proposition}

\begin{corollary}\label{corpropPolCarPuissExt}
If $A\in\Mn(\gA)$ is a matrix with \polcarz~$f$, and if $f(T)=\prod_{i=1}^n(T-x_i)$ in an over\ri of $\gA$, then the \polcar of $\Al{k}\!A$ is equal to the product 
$\prod_{J\in \cP_{k,n}}(T-x_J)$, where $x_J=\prod_{i\in J}x_i.$
\end{corollary}

\penalty-2500
\subsec{Tschirnhaus transformation}

\begin{definition}\label{defiTschir}
Let $f$ and $g\in\AT$ with $f$ a \mon of degree $p$.
Consider the \Alg $\gB=\aqo{\AT}{f}$, which is a free \Amo of rank~$p$. We define the 
\emph{Tschirnhaus transform of $f$ by $g$}, 
denoted by $\Tsc_{\gA,g}(f)$ or~$\Tsc_g(f)$, by the \egt

\snic{\Tsc_{\gA,g}(f)=\rC{\gB/\!\gA}(\ov g),\quad (\ov g $ is the class of $g$ in $\gB).}%
\index{Tschirnhaus transform}
\end{definition}

Proposition~\ref{prop1tschir} and Corollary~\ref{corprop1tschir} give the following result.

\begin{proposition}\label{prop2tschir}
Let $f$ and $g\in\AT$ with \mon $f$ of degree $p$.
\begin{enumerate}
\item If $A$ is a matrix such that $f(T)=\rC A(T)$, we have
$$\preskip.3em \postskip.4em 
\Tsc_g(f)(T)=\rC{g(A)}(T). 
$$
\item If $f(T)=\prod_i(T-x_i)$ on a \ri which contains $\gA$, we have
$$\preskip.3em \postskip.4em\ndsp 
\Tsc_g(f)(T)=\prod\nolimits_i\big(T-g(x_i)\big), 
$$
in particular, with $\gB=\aqo{\AT}{f}$ we get
$$\preskip.3em \postskip.4em\ndsp 
\rN\iBA (g)=\prod\nolimits_ig(x_i) \et
\Tr\iBA (g)=\som_ig(x_i). 
$$
\end{enumerate}
\end{proposition}

\medskip 
\rem We can also write $\Tsc_{\gA,g}(f)(T)=\rN_{\gB[T]/\!\AT}(T-\ov g)$.
In fact for an entirely unambiguous notation we should write $\Tsc(\gA,f,g,T)$ instead of $\Tsc_{\gA,g}(f)$.
An analogous ambiguity is also found in the notation $\rC{\gB/\!\gA}(g)$.
\eoe

\rdb
\subsubsection*{Computation of the Tschirnhaus transform}
Recall that the matrix $C$ of the \endo $\mu_t$ of multiplication by $t$ (the class of $T$ in $\gB$) is called the companion matrix of $f$ (see \paref{matrice.compagne}).
Then the matrix (over the same basis) of $\mu_{\ov g}=g(\mu_t)$ is the matrix $g(C)$. Thus~$\Tsc_g(f)$ is the \polcarz\footnote{The efficient computation of  \deters and \polcars is of great interest in computer algebra.
You can for example consult \cite{Jou}.
Another formula we can use for the computation of the Tschirnhaus transform is $\Tsc_g(f)=\Res_X\big(f(X),T-g(X)\big)$ (see Lemma~\ref{lemResultant}).} of $g(C)$.

\subsec{New version of the discriminant}

Recall (\Dfn \ref{defiDiscTra}) that when $\gC\supseteq \gA$ is a free \Alg of finite rank and $x_1$, \dots, $x_k\in\gC$, we call the \deter  of the matrix $\big(\Tr_{\gC/\!\gA}(x_ix_j)\big)_{i,j\in \lrbk}$ the discriminant of 
$(\xk)$. We denote it by $\disc_{\gC/\!\gA}(\xk)$.\\ 
Moreover, if $(\xk)$ is an $\gA$-basis of $\gC$, we denote by $\Disc_{\gC/\!\gA}$ the multiplicative class of $\disc_{\gC/\!\gA}(\xk)$ modulo the squares of $\Ati$. We call it the {\discri of the extension} $\gC\sur\gA$.

In this subsection, we make the link between the \discri of free \algs of finite rank and the \discri of \polusz.
\\
Let us emphasize the remarkable \crc of the implication \emph{1a} $\Rightarrow$ \emph{1b} in the following proposition.


\begin{proposition}\label{propdiscTra}
 \emph{(Trace-valued discriminant)}\\
Let $\gB$ be a free \Alg of finite rank $n$,
$x \in \gB$ and $f=\rC{\gB/\!\gA}(x)(T)$. We have
$$\preskip-.2em \postskip.4em
\disc(1, x, \ldots, x^{n-1}) = \disc(f) =
(-1)^{n(n-1) \over 2} \rN\iBA \big(f'(x)\big).
$$
We say that $f'(x)$ is the \emph{different of $x$}.%
\index{different!of an element in a finite free algebra}
The following results ensue.
\begin{enumerate}\itemsep1pt
\item \Propeq
\begin{enumerate}\itemsep1pt
\item $\disc(f)\in\Ati$.
\item  $\Disc\iBA \in\Ati$ and $(1,x,\ldots,x^{n-1})$ is an $\gA$-basis of $\gB$.
\item $\Disc\iBA \in\Ati$ and $\gB=\gA[x]$.
\end{enumerate}
\item If $\Disc\iBA $ is \ndzz, \propeq
\begin{enumerate}\itemsep1pt
\item $\Disc\iBA $ and $\disc(f)$ are associated \eltsz.
\item $(1,x,\ldots,x^{n-1})$ is an $\gA$-basis of $\gB$.
\item $\gB=\gA[x]$.
\end{enumerate}
\item The \discri of a \polu $g\in\AT$ represents
(modulo the squares of $\Ati$) the \discri of the extension $\aqo{\AT}{g}$ of $\gA$. We have $\mathrm{disc}_T(g)\in\mathbf{A}^\times$ if and only if $\langle{g(T),g'(T)}\rangle=\mathbf{A}$.
\end{enumerate}
\end {proposition}
\begin{proof}
In an over\ri $\gB'$ of $\gB$, we can write 
$f(T) = (T - x_1)\cdots (T - x_n)$.
 For some $g \in \gA[T]$, by applying Corollary~\ref{lemPolCar}, we obtain the \egts
$$\preskip.2em \postskip.4em 
\Tr\iBA \big(g(x)\big) = g(x_1) + \cdots +
g(x_n) \hbox{  and  }\rN\iBA \big(g(x)\big) = g(x_1) \cdots g(x_n). 
$$
 Let $M \in \Mn(\gA)$ be the matrix intervening in the computation of the \discri of $(1, x, \ldots, x^{n-1})$:
$$\preskip.2em \postskip.4em 
M=\big((a_{ij})_{i,j\in\lrb{0..n-1}}\big),\qquad a_{ij} = \Tr\iBA (x^{i+j}) = x_1^{i+j} + \cdots + x_n^{i+j}\,. 
$$
Let $V \in \Mn(\gB')$ be the Vandermonde matrix having $[\,x_1^i\; \ldots\; x_n^i\,]$ (where $i\in\lrb{0.. n-1}$) for rows. Then $M = V\tra {V}$.  We deduce
$$\preskip.2em \postskip.3em \ndsp
\det(M) = \det(V)^2 = \prod\nolimits_{i < j} (x_i - x_j)^2 = \disc(f)\,. 
$$
This proves the first \egtz. Since
$\rN\iBA \big(f'(x)\big) = f'(x_1) \cdots f'(x_n)$ and 
$
f'(x_i) =\prod_{j \mid j \ne i} (x_i - x_j)
$,
we get
$$\preskip.0em \postskip.4em\ndsp
\rN\iBA \big(f'(x)\big) = \prod\nolimits_{(i,j) \mid j \ne i} (x_i - x_j) =
(-1)^{n(n-1) \over 2} \prod\nolimits_{i < j} (x_i - x_j)^2.
$$
The \dem of the consequences is left to the reader (use Proposition~\ref{defiDiscTra}).
\end{proof}

\subsect{Discriminant of a \aduz}{Discriminant of a \aduz}

The \egt of the \gui{trace-valued} \discri and the \gui{\pollz} \discriz, together with the transitivity formula (\Thref{thTransDisc}), allows us 
to complete the following computation.

\pagebreak

\begin{fact}\label{factDiscriAdu} \emph{(Discriminant of a \aduz)}\\
Let $f$ be a \polu of degree $n\geq2$ of $\kT$ and $\gA=\Adu_{\gk,f}$.
\\ Then
$\Disc\iAk =\big(\disc_T(f)\big)^{n!/2}$.
\end{fact}
\begin{proof}
We use the notations of Section~\ref{sec0adu}.
We reason by \recu on~$n$, the $n=2$ case being clear.
We have $\gA=\gk_1[x_2,\ldots,x_n]$ with 

\snic{\gk_1=\gk[x_1]\simeq\aqo{\gk[X_1]}{f(X_1)}.}

Moreover, $\gA\simeq\Adu_{\gk_1,g_2}$
where 

\snic{g_2(T)=f_2(x_1,T)
=\big(f(T)-f(x_1)\big)\big/(T-x_1)\in\gk_1[T]\subseteq\gA[T].}

The transitivity formula of the \discris then gives the following \egtsz.
$$\preskip.0em \postskip.4em
{\mathrigid1mu \Disc\iAk =
\Disc_{\gk_1/\gk}^{~\dex{\gA: \gk_1}}\
 \rN_{\gk_1/\gk}(\Disc_{\gA/\gk_1})=
 (\disc f)^{(n-1)!}\ \rN_{\gk_1/\gk}(\Disc_{\gA/\gk_1}).}
$$
By using the \hdr we obtain the \egt
$$\preskip.2em \postskip.4em\ndsp 
\Disc_{\gA/\gk_1}=(\disc g_2)^{(n-1)!/2}=
\big(\prod\nolimits_{2\leq i<j\leq n}(x_i-x_j)^2\big)^{(n-1)!/2}. 
$$
For $i \in \lrb {2..n}$, let $\tau_i$ be the transposition $(1,i)$; for $z \in \gk_1$, by Corollary~\ref{lemPolCar}, $\rN_{\gk_1/\gk}(z) = z \prod_{i=2}^n \tau_i(z)$.  
Applied to $z = \prod\nolimits_{2\leq i<j\leq n}(x_i-x_j)^2$, this gives
$$\preskip.2em \postskip.4em\ndsp
\mathrigid 2mu \rN_{\gk_1/\gk}(z) = (\disc f)^{n-2},
\; \hbox { whence } \,
\rN_{\gk_1/\gk}(\Disc_{\gA/\gk_1}) = (\disc f)^{(n-2)\cdot(n-1)!/2}\,,
$$
then
$$\preskip.0em \postskip.4em
\Disc\iAk = (\disc f)^{(n-1)! + (n-2)\cdot(n-1)!/2} = (\disc f)^{n!/2}.
$$

NB: a detailed examination of the previous computation shows that in fact we have computed the \discri of the \gui{canonical} basis of the \adu described in Fact~\ref{factBase}.
\end{proof}
%


\begin{lemma}
\label{lemPolCarAdu} (Same assumptions as for Fact~\ref{factDiscriAdu})
Let $z\in\gA$.

\snic{\rC{\gA/\gk}(z)(T)=\prod\nolimits_{\sigma\in\Sn}\big(T-\sigma(z)\big).}

In particular, $\Tr\iAk (z)=\som_{\sigma\in\Sn}\!\sigma(z)$ and $\rN\iAk (z)=\prod_{\sigma\in\Sn}\!\sigma(z)$.
\end{lemma}
%
\begin{proof}
It suffices to show the formula for the norm, because we then obtain the one for the \polcar by replacing $\gk$ by $\kT$ (which replaces~$\gA$ by $\AT$).
The formula for the norm is proven by \recu on the number of variables by using Fact~\ref{factAduAdu}, the transitivity formula for the norms and Corollary~\ref{lemPolCar}.
\end{proof}
%

\newpage	
\section{Basic Galois theory (1)}
\label{secGaloisElr}

\Grandcadre{In Section~\ref{secGaloisElr},
$\gK$ designates a nontrivial \cdiz.}

\vspace{-8pt}
\subsec{Factorization and zeros}

Recall that a \ri is integral if every \elt is zero or \ndzz.\footnote{The notion is discussed in more detail on \paref{subsecAnneauxqi}.}
A sub\ri of an integral \ri is integral.
A \cdi is an integral \riz. A \ri $\gA$ is integral \ssi its total 
\ri of fractions $\Frac\gA$ is a \cdiz. We say that $\Frac\gA$ is the \ixc{field of fractions}{of an integral ring}, or the \emph{quotient field} of~$\gA$.%
\index{quotient field!of an integral ring}

\begin{proposition}\label{propZerFactPol}
Let $\gA\subseteq \gB$ be \ris and $f\in\AT$ be some \polu of degree~$n$.
\begin{enumerate}
\item If $z$ is a zero of $f$ in $\gB$, $f(T)$ is divisible by $T-z$ in~$\gB[T]$. 
\item Henceforth assume that $\gB$ is integral and nontrivial.\footnote{We could make do without the negative assumption \gui{nontrivial} by reading the assumption that the $z_i$'s are \gui{distinct} as meaning that each $z_i-z_j$ is \ndzz.} If $z_1$, $\dots$, $z_k$ are the  pairwise distinct zeros of $f$ in $\gB$, the \pol $f(T)$ is divisible by
$\prod_{i=1}^k(T-z_i)$ in~$\gB[T]$.
\item In addition, if $k=n$, then $f(T)=\prod_{i=1}^n(T-z_i)$, and the $z_i$'s are the only zeros of $f$ in $\gB$ and in every integral extension of~$\gB$.
\end{enumerate}
\end{proposition}
%
\begin{proof}
The \dem is \imdez. Certain more precise results are in Exercise~\ref{exoLagrange}, which is dedicated to Lagrange interpolation.
\end{proof}
%

\subsec{Strictly finite \algs over a \cdiz}

\begin{definition}\label{defiSTF}~\\
A \Klg $\gA$ is said to be
 \emph{\stfez} if it is a free \Kev of finite dimension.%
\index{strictly finite!algebra over a discrete field}%
\index{algebra!strictly finite --- over a discrete field}
\end{definition}

In other words, we know of a finite basis of $\gA$ as in a \Kevz.
In this case, for some $x\in\gA$, the trace, the norm, the \polcar of (multiplication by) $x$, as well as the \polmin of $x$ over $\gK$, denoted by~$\Mip_{\gK,x}(T)$ or $\Mip_{x}(T)$, can be computed by standard methods of the \lin \alg over a \cdiz.
Similarly, every finite \Kslg of~$\gA$ is \stfe and the intersection of two \stfes sub\algs is \stfez.

\begin{lemma}\label{lemEntReduitConnexe}
Let $\gB\supseteq\gK$ be a \ri integral over $\gK$. \Propeq
\begin{enumerate}
\item $\gB$ is a \cdiz.
\item $\gB$ is \sdzz: $xy=0\Rightarrow (x=0$ \emph{or} $y=0)$.
\item $\gB$ is connected and reduced.
\end{enumerate}
Consequently, if $\gB$ is a \cdiz, every finite \Kslg of $\gB$ is a \cdiz.
 \end{lemma}
\begin{proof}
 The implications \emph{1} $\Rightarrow$ \emph{2}
$\Rightarrow$ \emph{3} are clear.
\\
\emph{3} $\Rightarrow$ \emph{1}.
Each element $x\in\gB$ annihilates a nonzero \pol of $\KX$ that we can assume is of the form $X^k\big(1-XR(X)\big)$. Then $x\big(1-xR(x)\big)$ is nilpotent thus zero.
The \elt $e=xR(x)$ is \idm and $x=ex$. If $e=0$, then~$x=0$. If $e=1$, then $x$ is \ivz.
\end{proof}
%

\begin{lemma}\label{lemdefiSTF}
Let $\gK\subseteq\gL\subseteq\gA$ with $\gA$ and $\gL$ \stfs over $\gK$.
If~$\gL$ is a \cdiz, then $\gA$ is \stfe over $\gL.$
\end{lemma}
\begin{proof}
\Demo left to the reader (or see Fact~\ref{fact1Etale} item \emph{\iref{i3fact1Etale}}).
\end{proof}

If $g$ is an \ird \pol of $\gK[T]$, the quotient \alg $\aqo{\gK[T]}{g}$ is a \stf \cdi over $\gK$.
In fact, as a corollary of the two previous Lemmas we get that every \stfe extension of \cdis is obtained by iterating this construction.

\begin{fact}\label{factStrucExtFiC} \emph{(Structure of a \stfe extension of \cdisz)}\\
Let $\gL=\gK[\xm]$ be a \stf \cdi over $\gK$.\\ 
For $ k\in\lrb{1.. m+1}$, let $\gK_k=\gK[(x_i)_{i<k}]$ and $f_k=\Mip_{\gK_{k},{x_k}}(T)$, such that $\gK_1=\gK$, and for $k\in \lrbm$, $\gK_{k+1}\simeq\aqo{\gK_k[X_k]}{f_k(X_k)}$. \\
Then, for $k<\ell$ in $\lrb{1.. m+1}$, the inclusion $\gK_k\to\gK_\ell$ is a \stfe extension of \cdisz, with 
$$\preskip.2em \postskip.3em \ndsp
\dex{\gK_\ell:\gK_k} = \prod\nolimits_{k\leq i< \ell}\dex{\gK_{i+1}:\gK_i}
  = \prod\nolimits_{k\leq i< \ell}\deg_T(f_i). 
$$
Moreover, if $F_k\in\gK[X_1,\dots,X_{k}]$ is a \polu in $X_k$ for which we have $F_k\big((x_i)_{i<k},X_k\big)=f_k(X_k)$, we get, by \fcn of the \evn \homoz, an \iso
$$\preskip.2em \postskip.4em 
       \aqo{\gK[\Xm]}{F_1,\dots,F_m} \; \simarrow \; \gL. 
$$
\end{fact}

\smallskip 
\begin{definition}
\label{defCorpsdesRacines}
Let $g\in\gK[T]$ be a \poluz, we call  a \cdi
$\gL$ extension of $\gK$ in which~$g$ can be completely decomposed and which is generated like \Klg by the zeros of~$g$ a \emph{\cdr of $g$ over $\gK$}.%
\index{field of roots!of a polynomial}
\end{definition}

Note that $\gL$ is finite over $\gK$ but that we do not ask that $\gL$ be \stf over $\gK$
(in fact, there is no \cov \dem that such a \cdr must be \stf over $\gK$). This necessitates some subtleties in the following \thoz.

\pagebreak

\begin{theorem} \emph{(Uniqueness of the \cdr in the \stf case)}
\label{propUnicCDR}
Let $f\in\gK[T]$ be a \poluz. Assume that there exists a \cdr $\gL$ of $f$ over $\gK$.
\begin{enumerate}
\item Let $\gM\supseteq\gK$ be a \stf \cdi over $\gK$,  generated by $\gK$ and some zeros of $f$ in $\gM$.
 The field $\gM$ is isomorphic to a subfield of~$\gL$.   
\item  Assume that there exists a \cdr of $f$, \stf over~$\gK$. Then every \cdr of $f$ over $\gK$ is isomorphic to $\gL$ (which is thus \stf over $\gK$).
\item Let $\gK_1$, $\gK_2$ be two nontrivial \cdiz, $\tau : \gK_1 \to \gK_2$ be an \isoz, $f_1 \in \gK_1[T]$ be a \poluz, $f_2 = f_1^\tau \in \gK_2[T]$.  If~$\gL_i$ is a \stf field of roots of $f_i$ over $\gK_i$ ($i = 1, 2$), then~$\tau$  extends to an \iso from $\gL_1$ to $\gL_2$.
\end{enumerate}
\end{theorem}
\begin{proof}
We only prove item \emph{1} in a particular case (sufficiently \gnlz). The rest is left to the reader.\\
We write $f(T)=\prod_{i=1}^n(T-x_i)$ in~$\gL[T]$. 
Also suppose \hbox{that $\gM=\gK[y,z]$} with $y\neq z$ \hbox{and $f(y)=f(z)=0$}. \\
We thus have in $\gM[T]$ the \egt $f(T)=(T-y)f_1(T)=(T-y)(T-z)f_2(T)$ (Proposition~\ref{propZerFactPol}). \\ 
Since $f(y)=0$, the \polmin $g(Y)$ of $y$ over $\gK$ divides $f(Y)$ in~$\gK[Y]$. Therefore $\prod_{i=1}^ng(x_i)=0$ in $\gL$, which is a \cdiz, and one of the $x_i$'s, say $x_1$, annihilates $g$. Here we obtain

\snic{
\gK[y] \simeq \aqo{\gK[Y]}{g(Y)}\simeq\gK[x_1] \subseteq\gL.
}

The \cdi $\gK[y]$ is \stf over $\gK$ and  $\gM$ is \stf over $\gK[y]$ (Lemma~\ref{lemdefiSTF}). Then let $h\in\gK[Y,Z]$ be a \polu in $Z$ such that $h(y,Z)$ is the \polmin of $z$ over $\gK[y]$. \\
Since $f_1(z)=0$, the \polz~$h(y,Z)$ divides $f_1(Z)=f(Z)/(Z-y)$ in~$\gK[y][Z]$, thus its image $h(x_1,Z)$ in~$\gK[x_1][Z]$ is an \ird \pol which divides $f(Z)/(Z-x_1)$. So $h(x_1,Z)$ admits as a zero one of the $x_i$'s for $i\in\lrb{2..n}$, say $x_2$, and $h(x_1,Z)$ is the \polmin of~$x_2$ over $\gK[x_1]$. We thus obtain the \isos
$$\preskip.3em \postskip.3em 
\gK[y,z]\simeq\aqo{\gK[y][Z]}{h(y,Z)}\simeq\aqo{\gK[x_1][Z]}{h(x_1,Z)} \simeq
\gK[x_1,x_2]\subseteq\gL. 
$$
Note that we also have $\gK[y,z]\simeq\aqo{\gK[Y,Z]}{g(Y),h(Y,Z)}$.
\end{proof}

\rem A detailed inspection of the previous \dem leads to the conclusion that if $\gL$ is a \stf \cdr over $\gK$, the group of $\gK$-\autos of $\gL$ is a finite group having at \hbox{most $\dex{\gL:\gK}$} \eltsz. 
If we do not assume that $\gL$ is \stf over $\gK$, we only obtain that it is absurd to assume that this group contains more \hbox{than~$\dex{\gL:\gK}$} \eltsz.
\eoe

\subsec{The \elr case of Galois theory}
\label{subsecGaloisElr}

\begin{definota}\label{NOTAStStp}
We will use the following notations when a group $G$ operates over a set~$E$.
\begin{enumerate}
\item [---]
For $x\in E$, \smash{$\St_G(x)=\St(x)\eqdefi\sotq{\sigma\in G}{\sigma(x)=x}$} designates the  \ix{stabilizer} of $x$.
\item [---] $G.x$ designates the orbit of $x$ under $G$, and we write $G.x=\so{x_1,\ldots,x_k}$ as an abbreviation for: \emph{$(x_1,\ldots,x_k)$ is an enumeration without repetition of $G.x$, with $x_1=x$.}
\item [---]
 For $F\subseteq E$, $\Stp_G(F)$ or $\Stp(F)$ designates the pointwise stabilizer of~$F$.
\item [---]
 If $H$ is a subgroup of $G$,
 \begin{enumerate}
\item [--] we denote by $\idg{G:H}$ the index of $H$ in $G$,%
\index{index!of a subgroup in a group}
\item [--] we denote by $\Fix_E(H)=\Fix(H)=E^H$
the subset of \elts fixed by $H$, $\sotq{x\in E}{\forall \sigma\in H,\;\sigma(x)=x}$,
\item [--]  writing $\sigma\in G/H$ means that we take an \elt $\sigma\in G$ in each  left coset of $H$ in $G$.
\end{enumerate}
 \end{enumerate}
When $G$ is a finite group 
operating over a \ri $\gB$, for $b\in \gB$, we write
\[
\Tr_G(b)   =   \sum_{\sigma\in G}\sigma(b), \;
 \rN_G(b)  =   \prod_{\sigma\in G}\sigma(b), \hbox{ and }
\rC{G}(b)(T)   =    \prod_{\sigma\in G}\big(T-\sigma(b)\big).
\]
If $G.b=\so {b_1,\ldots ,b_k}$, (the $b_i$'s pairwise distinct), we write
$$\preskip.2em \postskip.4em\ndsp 
\Rv_{G,b}(T)=\prod_{i=1}^k(T-b_i). 
$$
This \pol is called the \emph{resolvent} of $b$ (relative to $G$).
It is clear that $\big({\Rv_{G,b}}\big)^r=\rC{G}(b)$ with $r=\idG{G:\St_G(b)}$.%
\index{resolvent}
\end{definota}

\rdb
Given an \Alg $\gB$ we denote by $\Aut_\gA(\gB)$ the group of $\gA$-\autos of $\gB$.
\label{NOTAAutAB}

\begin{definition}\label{defiGalGal}
\label{i4defiCorGal}
If $\gL$ is a \stf extension of $\gK$, and a \cdr for a \spl \polu over $\gK$, we say that~$\gL$ is a \ix{Galois extension} of $\gK$, we then denote $\Aut_\gK(\gL)$ by $\Gal(\gL/\gK)$ and we say that it is the \ix{Galois group} of the extension $\gL/\gK$. 
\end{definition}

Note well that in the \dfn of a Galois extension $\gL/\gK$, the fact that $\gL$ is \stf (and not only finite) over $\gK$ is implied.

\begin{propdef}\label{defiCorGal} \emph{(Galois correspondence)}\\
Let $\gL\supseteq\gK$ be a \stf field over~$\gK$.
\begin{enumerate}
\item  \label{i1defiCorGal}
  The group\perso{I think that there is no \cov \dem that this group is finite. Otherwise we would for instance have an \algo to decide if an extension of degree 3 is Galoisian and in characteristic $0$ this reduces to decide whether a discriminant is a square.}
 $\Aut_\gK(\gL)$ is a detachable subgroup of $\GL_\gK(\gL)$.
  If $H$ is a subgroup of $\Aut_\gK(\gL)$, the subfield $\gL^H$ is called the \emph{fixed field of~$H$}.
\item \label{i2defiCorGal} We call  the two mappings $\Fix$ and $\Stp$ between the two following sets the \ix{Galois correspondence}. On the one hand $\cG=\cG_{\gL/\gK}$ is the set of finite subgroups of $\Aut_\gK(\gL)$. On the other hand $\cK=\cK_{\gL/\gK}$ is the set of  \stf subextensions of~$\gL$.
\item \label{i3defiCorGal} In the Galois correspondence each of the two mappings is decreasing.
In addition, $H\subseteq \Stp(\gL^H)$ for all~$H\in\cG$, $\gM\subseteq \gL^{\Stp(\gM)}$ for all~$\gM\in\cK$, $\Stp\circ\Fix\circ\Stp=\Stp$ and $\Fix\circ\Stp\circ\Fix=\Fix$.
%
%
\end{enumerate}
\end{propdef}
\begin{proof}
In item~\emph{\ref{i1defiCorGal}} we have to prove that the subgroup is detachable and, in item~\emph{\ref{i2defiCorGal}}, that $\Fix$ and $\Stp$ indeed act on the two sets as described. This is based on finite dimensional \lin \alg   over the \cdisz. We leave the details to the reader.
\end{proof}

\rem
Even though we can decide if a given \elt of $\GL_\gK(\gL)$ is in~$\Aut_\gK(\gL)$, and even though it is easy to bound the number of \elts of~$\Aut_\gK(\gL)$, there is no sure \gnl method to compute this number.
\eoe

\medskip 
As a consequence of \Thref{propUnicCDR} we have the following corollary.

\begin{theorem}
\label{thPrIs} \emph{(Isomorphism extension \thoz)}\\
Let  $\gL/\gK$ be a Galois extension and $\gM$ be a finite $\gK$-subextension of~$\gL$. Every $\gK$-\homo $\tau:\gM\to\gL$ extends to an \eltz~$\wi\tau$ of~$\Gal(\gL/\gK)$.
\end{theorem}
%
\begin{proof}
$\gL$ is the \cdr of a \spl \pol $g\in\KT$.
We notice that since $\gL$ is \stf over $\gK$, $\gM$ is \stf over $\gK$ and $\gL$ \stf over $\gM$.
  Let $\gM'$ be the image of $\tau$. It is a \stf field over $\gK$, so $\gL$ is \stf over $\gM'$.
Thus $\gL$ is a field of roots of $g$ \stf over $\gM$ and over $\gM'$. By \Thref{propUnicCDR} (item \emph{3}), we can extend $\tau$ to a $\gK$-\iso $\wi\tau : \gL \to \gL$.
\end{proof}

\rdb \label{NOTAGalKf}
 When a \spl \pol over $\gK$ has a \cdr $\gL$ \stf over $\gK$, the group $\Gal(\gL/\gK)$ can also be denoted by $\Gal_\gK(f)$ insofar as \Thref{propUnicCDR} gives the uniqueness of $\gL$ (up to $\gK$-\autoz).

\medskip
\rem
In \coma we have the following results (trivial in \clamaz).
For some subgroup $H$ of a finite group \propeq
\begin{itemize}\itemsep.5pt
\item  $H$  is finite.
\item  $H$  is \tfz.
\item  $H$  is detachable.
\end{itemize}
Similarly for some \Ksv $M$ of a finite dimensional \Kev  \propeq
\begin{itemize}\itemsep.5pt
\item  $M$  is finite dimensional.
\item  $M$  is \tf (i.e., the image of a matrix).
\item  $M$  is the kernel of a matrix.\eoe
\end{itemize}

\vspace{-.3em}
\pagebreak

\begin{propdef}\label{propGaloiselr}\emph{(\Elr Galois situation)}\\
Let $\gA\subseteq\gB$ be two \risz.
An \emph{\elr Galois situation} is defined as follows.
\begin{enumerate}\itemsep0pt
\item [{i.}] We have a \spl \mon \polz~$Q\in\AT$ of degree $d$ and  \elts $y_1$, $y_2$, $\ldots$, $y_d$ of $\gB$ such that
$$\preskip.4em \postskip.4em \ndsp
Q(T)=\prod_{i=1}^d (T-y _i). 
$$
\item [{ii.}]  Let $y=y_1$. Assume for each $i$ that $\gB=\gA[y _i]$ and that $\gen{Q}$ is the kernel of the \homo of \Algs $\AT\to\gB$ which sends $T$ into $y_i$
(whence $\gB=\gA[y]=\gA[y _i]\simeq\aqo{\AT}{Q}$).
For each $i$ there thus exists a unique $\gA$-\auto $\sigma_i$ of $\gB$ satisfying $\sigma_i(y)=y_i$.

\item [{iii.}]  Assume that these \autos form a group, which we denote by~$G$. In particular, $\abs{G}=d=\dex{\gB:\gA}$.
\end{enumerate}
In the \elr Galois situation we have the following results.
\begin{enumerate}\itemsep0pt
\item
\begin{enumerate}\itemsep0pt
\item \label{i1propGaloiselr}
$\Fix_\gB(G)=\gA$.
\item \label{i2propGaloiselr}
For all $z\in\gB$, $\rC{\gB/\gA}(z)(T)=\rC G(z)(T)$.
\end{enumerate}
\item \label{i3propGaloiselr}
Let $H$ be a detachable subgroup of $G$, $\gA'=\gB^H$ and 
$$\preskip.4em \postskip.4em \ndsp
Q_H(T)=\prod_{\sigma\in H} \big(T-\sigma(y)\big). 
$$
Then, we find the \elr Galois situation with $\gA'$, $\gB$, $Q_H$ and $\big(\sigma(y)\big)_{\sigma\in H}$.
In particular, $\gB=\gA'[y]$ is a free $\gA'$-module of rank $\abs{H}=\dex{\gB:\gA'}$.
In addition, $H$ is equal to $\Stp_G(\gA')$.
\end{enumerate}
\end{propdef}
\begin{proof}
\emph{1a.}
Consider some $x=\sum_{k=0}^{d-1} \xi_ky ^k$ in $\gB$ (with each $\xi_k\in\gA$) invariant under the action of $G=\so{\sigma_1,\ldots,\sigma_{d}}$.
\\
We thus have for all $\sigma\in G$, $x=\sum_{k=0}^{d-1} \xi_k\sigma(y )^k$.
If~$V \in \Mn(\gB)$ is the Vandermonde matrix
 $$\preskip.0em \postskip.0em
 V=\cmatrix {1& y_1  &y_1^2& \cdots&y_1^{d-1}\cr
 \vdots&&&&\vdots\cr
 \vdots&&&&\vdots\cr
 1& y_d&y_d^2& \cdots&y_d^{d-1}
 },$$
  we get
$$
V\ \left[\begin{array}{c} \xi _0 \\ \xi _1 \\ \vdots \\ \xi _{d-1} \end{array}\right]
=\left[\begin{array}{c} x \\ x \\ \vdots \\ x \end{array}\right]=
V\ \left[\begin{array}{c} x \\ 0 \\ \vdots \\ 0 \end{array}\right].
$$
Since $\det(\!\tra{V}V)=\disc_T(Q)\in\Ati$, we get
$[\,\xi_0 \ \xi_1 \ \cdots\ \xi_{d-1}\,]=[\,x \ 0\ \cdots\ 0\,]$, \hbox{and $x =\xi _0\in\gA$}.
\\
\emph{1b.}
Since $\gB\simeq\aqo{\AT}{Q}$, Corollary~\ref{corprop1tschir} gives, for $g\in\AY$ \hbox{and~$z=g(y_1)$}, the \egts
$$\preskip.0em \postskip.4em\ndsp
\rC{\gB/\gA}(z)(T)=\prod\nolimits_i\big(T-g(y_i)\big) =
\prod\nolimits_{\sigma\in G}\big(T-\sigma(g(y_1))\big) = \rC G(z)(T).
$$
\emph{2.} It is clear that $\gB=\gA'[\sigma(y)]$ for each $\sigma\in H$ and that $Q_H$ is a \spl \pol of $\gA'[T]$. It remains to see that every \polz~$P\in\gA'[T]$ which annihilates some $y_i=\sigma_i(y)$ ($\sigma_i\in H$) is a multiple of $Q_H$. For all \hbox{\mathrigid 1mu $\sigma\in H$}, since~$\sigma$ is an
$\gA'$-\auto of $\gB$, we have {\mathrigid 1mu $P\big(\sigma(y_i)\big)=\sigma\big(P(y_i)\big)=0$}. Thus~$P$ is divisible by each~$T-\sigma(y)$, for $\sigma\in H$. As these \pols are pairwise \comz, $P$ is a multiple of their product~$Q_H$. \\
Finally, if $\sigma_j\in G$ is an $\gA'$-\auto of $\gB$, $\sigma_j(y)=y_j$ must be a zero of $Q_H$. However, since $Q$ is \splz, the only $y_i$'s that annihilate $Q_H$ are the $\sigma(y)$'s for $\sigma\in H$. Therefore $\sigma_j\in H$.
\end{proof}

\rems 1) In the \elr Galois situation nothing states that the $y_i$'s are the only zeros of $Q$ in $\gB$, nor that the $\sigma_i$'s are the only \hbox{$\gA$-\autosz} of~$\gB$. Take for example $\gB=\gK^3$, and three distinct \elts $a$, $b$, $c$
  in the \cdi $\gK$. The \polz~$Q=(T-a)(T-b)(T-c)$ admits 27 zeros in $\gB$, including six 
 which have $Q$ as \polminz, which makes six $\gK$-\autos of $\gB$.
\\
In addition, if we take $z_1=(a,b,c)$, $z_2=(b,a,b)$ and $z_3=(c,c,a)$, we see that~$Q=(T-z_1)(T-z_2)(T-z_3)$, which shows that the first condition does not imply the second. However, with $y_1=(a,b,c)$, $y_2=(b,c,a)$ \hbox{and $y_3=(c,a,b)$}, we are in the \elr Galois situation.

2) Concerning condition \emph{iii} in the \dfn of the \elr Galois situation, we can easily see that it is equivalent to the fact that each $\sigma_i$ permutes the $y_j$'s.
This condition is not a consequence of the first two, as the following example proves. Consider the following $5 \times 5$ latin square (in each row and each column, the integers are different), which is not the table of a group
$$\preskip.0em \postskip.4em
\cmatrix {
1 & 2 & 3 & 4 & 5\cr
2 & 4 & 1 & 5 & 3\cr
3 & 5 & 4 & 2 & 1\cr
4 & 1 & 5 & 3 & 2\cr
5 & 3 & 2 & 1 & 4\cr
}.
$$
Each row defines a permutation $\sigma_i \in S_5$; thus $\sigma_1 = \Id$, $\sigma_2 = (12453)$, \dots, $\sigma_5 = (154)(23)$. The $\sigma_i$'s do not form a group (which would be of order $5$) because $\sigma_5$ is of order $6$.
Let $\gB = \gK^5$ where $\gK$ is a field having at least $5$ \elts $a_1$, $\ldots$, $a_5$, $y = (a_1, \ldots, a_5) \in \gB$, $y_i = \sigma_i(y)$ and

\snic{Q(T) = \prod\nolimits_i (T - y_i) = \prod\nolimits_i (T - a_i) \in \gK[T].}

Then, in \ref{propGaloiselr}, the first two conditions \emph{i}, \emph{ii} are satisfied but not condition~\emph{iii.}\\
Luckily things are simpler in the field case. \eoe

\pagebreak

\begin{lemma}\label{lemGaloiselr}
Let $\gL=\gK[y]$ be a \stf \cdi over~$\gK$.
Let~$Q$ be the \polmin of $y$ over $\gK$.
If $Q$ is \spl and can be completely factorized in $\gL[T]$, we find ourselves in the \elr Galois situation and the corresponding group $G$ is the group $\Gal(\gL/\gK)$ of all \hbox{the $\gK$-\autosz} of~$\gL$.
\end{lemma}
\begin{proof}
Let $y=y_1$, $\ldots$, $y_d$ be the zeros of $Q$ (of degree $d$) in $\gL$.
Each $y_i$ annihilates~$Q$ and~$Q$ is \ird in $\KT$, so $Q$ is the \polmin of $y_i$ over~$\gK$ and $\gK[y_i]$ is a \Ksv of $\gL$, free and of same dimension $d$, therefore equal to $\gL$.
Finally, since $\gL$ is integral, the $y_i$'s are the only zeros of~$Q$ in~$\gL$, thus every $\gK$-\auto of $\gL$ is some $\sigma_i$, and the $\sigma_i$'s do indeed form a group: the Galois group $G=\Gal(\gL/\gK)$.
\end{proof}
%
\begin{theorem} \emph{(Galois correspondence, the \elr case)}
\label{thGaloiselr} \\
Let $\gL=\gK[y]$ be a \stf \cdi over $\gK$.
Let $Q$ be the \polmin of $y$ over $\gK$.
Assume that $Q$ is \spl and can be completely factorized in~$\gL[T]$. In particular, $\gL$ is a Galois extension of~$\gK$.
We have the following results.
\begin{enumerate}
\item  The two maps of the Galois correspondence are two reciprocal bijections.
\item  For all $\gM\in \cK_{\gL/\gK}$, $\gL/\gM$ is a Galois extension of the Galois group $\Fix(\gM)$ and $\dex{\gL:\gM}=\abs{\Fix(\gM)}$.
\item If $H_1, H_2\in\cG_{\gL/\gK}$ and $\gM_i=\Fix(H_i)\in \cK_{\gL/\gK}$, then
\begin{itemize}
\item  $H_1\cap H_2$
corresponds to the \Kslg generated by $\gM_1\cup \gM_2$,
\item  $\gM_1\cap \gM_2$
corresponds to the subgroup generated by $H_1\cup H_2$.
\end{itemize}
\item If $H_1\subseteq H_2$ in $\cG_{\gL/\gK}$ and $\gM_i=\Fix(H_i)$, then $\gM_1\supseteq\gM_2$ and we have the \egt $\idg{H_2:H_1}=\dex{\gM_1:\gM_2}$.
\item For all $z\in\gL$, $\rC{\gL/\gK}(z)(T)=\rC{\Gal(\gL/\gK)}(z)(T)$.
\end{enumerate}
\end{theorem}
 
\begin{proof}
It suffices to prove the first item.
By Proposition~\ref{propGaloiselr} we have the \egt $\Stp\circ \Fix=\Id_{\cG_{\gL/\gK}}$.\\
Now let $\gM\in \cK_{\gL/\gK}$. Since $\gL=\gK[y]$, we have $\gL=\gM[y]$.
As~$\gL$ is \stf over~$\gM$, we can compute the \polminz~$P$ of~$y$ over~$\gM$. It divides~$Q$ therefore it is \splz. It can be completely factorized in~$\gL[T]$.
Thus, with $\gM$, $\gL=\gM[y]$ and~$P$, we are in the assumptions of Lemma~\ref{lemGaloiselr},
so in the \elr Galois situation.  
The $\gM$-\autos of~$\gL$ are $\gK$-\autos thus they are exactly the \elts of the stabilizer $H=\Stp_G(\gM)$ (where $G=\Gal(\gL/\gK)$).
In this situation item \emph{1b} of Proposition~\ref{propGaloiselr} states that $\Fix(H)=\gM$.
\end{proof}

We have just established that the Galois correspondence is bijective, \cad the fundamental \tho of  Galois theory, in the \elr case. However, it will later turn out that this case is in fact the \gui{\gnlz} case: each time that we have a Galois extension we can reduce to the \elr situation 
   (Theorems~\ref{thResolUniv} and~\ref{thEtalePrimitif}).

\subsect{Construction of a \cdr by means of a Galois resolvent, basic Galois theory}{Construction of a \cdrz}
\label{secResolUniv}

\Grandcadre{In this subsection $f\in\gK[T]$ is a \spl \polu \\ of degree $n$ and $\gA=\Adu_{\gK,f}$ with $f(T)=\prod_i(T-x_i)$ in $\gA$.}

The aim of the current subsection is to prove the following result: if $\gK$ is infinite, and if we know how to factorize the \spl \polus in $\KT$, then we know how to construct a \cdr for any arbitrary \spl \poluz, and the obtained extension falls within the \elr framework of \Thref{thGaloiselr}.

We construct the \cdr by a \gui{uniform} method.
Since it is \stfz,  \thref{propUnicCDR} says that this  \cdr is \isoc to any other.

\begin{theorem}\label{thResolUniv}
We introduce \idtrs $u_1$, \dots, $u_n$.
For $\sigma\in\Sn$ we define $u_\sigma=\sum_iu_ix_{\sigma i}$.
We write
$$\preskip.4em \postskip.4em\ndsp
R(\uu,T)=\prod\nolimits_{\sigma\in\Sn}(T-u_\sigma)\in\gK[\uu,T],
$$
and $D(\uu)=\disc_T(R)\in\Kuu$.
\begin{enumerate}
\item \label{i1thResolUniv} One of the \coes of~$D$ is equal to $\pm\disc(f)^{(n-2)!(n!-1)}$.
\end{enumerate}
In the following we assume that we specialize the $u_i$'s to \elts $a_i\in\gK$ and that $D(\ua)\neq0$ (this is always possible if $\gK$ is infinite).
\begin{enumerate}\setcounter{enumi}{1}
\item \label{i2thResolUniv} For any arbitrary $\sigma\in\Sn$, the \elt $a_\sigma=\sum_ia_ix_{\sigma i}$ admits the \pol $R(\ua,T)\in\KT$
for \polminz, such that
$$\preskip.4em \postskip.2em
\gA=\gK[a_\sigma]\simeq\aqo{\KT}{R(\ua,T)}.
$$
We write $a=a_\Id=\sum_i a_ix_i$.
\item \label{i3thResolUniv} The only \elts of $\gA$ fixed by $\Sn$ are the \elts of $\gK$.
\item  Assume that we know how to decompose $R(\ua,T)$ into a product of \irds factors in $\KT$:
$R(\ua,T)=\prod_{j=1}^\ell Q_j$.
\begin{enumerate}
\item \label{i4thResolUniv} If $\ell=1$, $\gA$ is a field,
the extension $\gA/\gK$ is a \cdr for the \pol $f$, as well as for $R(\ua,T)$, and the situation pertains to \Thref{thGaloiselr}. In particular, $\Gal(\gA/\gK)\simeq\Sn.$
\item \label{i5thResolUniv} If $\ell>1$, then
 $\gA\simeq\prod_j\gK_j$ where 
 
\snic{\gK_j=\gK[\pi_j(a)]=\aqo{\gA}{Q_j(a)}\simeq\aqo{\KT}{Q_j}.}

($\pi_j:\gA\to\gK_j$ is the canonical \prnz.) \\
Let $H_j$ be the subgroup of $\Sn$ that stabilizes the \id $\gen{Q_j(a)}_{\!\gA}$. Then
\begin{enumerate}
\item [--] $\Sn$ operates transitively over the \ids $\gen{Q_j(a)}_{\!\gA}$, 
so that
the $Q_j$'s all have the same degree, $\abs{H_j}=\deg(Q_j)=\dex{\gK_j:\gK}$, and the $\gK_j$'s are  pairwise isomorphic \cdisz,
\item [--] the extension $\gK_1/\gK$ is a \cdr for $f$, as well as for each $Q_j$, and the situation pertains to \Thref{thGaloiselr}, in particular,
$H_1=\Gal(\gK_1/\gK)$.
\end{enumerate}
\end{enumerate}
\end{enumerate}
\end{theorem}

\vspace{-.4em}
\pagebreak

\begin{proof}
\emph{\ref{i1thResolUniv}.} The discriminant $D$ equals (up to sign) the product of the $u_\sigma-u_\tau$ for $\sigma\neq\tau\in\Sn$.
Each $u_\sigma-u_\tau$ is a sum of \elts $u_i(x_{\sigma i}-x_{\tau i})$:
each $u_i$ has  \coe $0$ or some $x_j-x_k$ ($j\neq k$).
The first \mom for the lexicographical order that appears in the product~$D$ is the \mom
$$\preskip.2em \postskip.4em
u_1^{n!(n!-(n-1)!)}u_2^{n!\big((n-1)!-(n-2)!\big)}\cdots u_{n-1}^{n!(2!-1!)},
$$
with \coe a product of \elts of the type $x_i-x_j$  ($i\neq j$). More \prmt if
$\delta=\disc(f)$, the \coe in question will be, up to sign,

\snic{ \delta^{(n-2)!(n!-1)}.}

\emph{\ref{i2thResolUniv}.} We use Proposition~\ref{propdiscTra} since $R(\ua,T)$ is the \polcar of~$a$ (Lemma~\ref{lemPolCarAdu}).
\\
\emph{\ref{i3thResolUniv}.}
See item \emph{1b} of Proposition~\ref{propGaloiselr}.\\
\emph{4a.} This is obvious.\\
\emph{4b.} The fact that $\gA\simeq\prod_j\gK_j$ results from the Chinese remainder \thoz.\\
The \egt $\prod_jQ_j(T)=\prod_{\sigma}(T-a_{\sigma})$ in $\AT$ remains valid in $\gK_1[T]$. \\
Thus, there exists for all $j$ some $\sigma_j$ such that $Q_j\big(\pi_1(a_{\sigma_j})\big)=0$, in other  
words, $Q_j(a_{\sigma_j})\in \gen{Q_1(a)}_{\!\gA}$. Furthermore, in $\gA$ we have $Q_j(a_{\sigma_j})=\sigma_j\big(Q_j(a)\big)$ because $Q_j\in\KT$. So $\sigma_j\big(\gen{Q_j(a)}_{\!\gA}\big)\subseteq \gen{Q_1(a)}_{\!\gA}$.
\\
This gives us a surjection $\sigma_j : \aqo\gA {Q_j(a)}   \to   \aqo\gA {Q_1(a)}$,
\cad a surjection $\aqo\KT {Q_j} \to   \aqo\KT {Q_1}$.
This results in $\deg Q_1\leq\deg Q_j$, and by symmetry $\deg Q_j=\deg Q_1$, whence $\sigma_j\big(\gen{Q_j(a)}_{\!\gA}\big)= \gen{Q_1(a)}_{\!\gA}$.
\\
Thus $\Sn$ operates transitively over the \ids $\gen{Q_j(a)}_{\!\gA}$ and the $\gK_j$'s are pairwise isomorphic.
\end{proof}

\rem The construction of the \cdr suggested here is in fact more or less impractical as soon as the degree $n$ of $f$ is equal to or greater than 7, as it necessitates a factorization of a \pol of degree $n!$.
We propose in~Chapter~\ref{ChapGalois} a less brutal dynamic method that has the additional advantage of not demanding to know how to factorize the \spl \pols of~$\KT$.
The counterpart of this absence of factorization will be that, despite knowing how to compute in \gui{some} \cdrz, a priori we will never be able to determine  it in its entirety (in the sense of knowing its dimension as a \Kevz). The same lack of precision also applies to the Galois group.
\eoe

\medskip
\exl Consider the \pol {\tt p(T)}$\in\QQ[T]$ below.
We ask {\tt Magma} to randomly take some \coli {\tt z} from the {\tt xi} (the zeros \hbox{of {\tt p(T)}} in the \adu $\gA=\Adu_{\QQ,p}$), to compute $\Mip_{\QQ,z}(T)$, and then to factorize it.
The software efficiently gives the \polminz~{\tt pm} of degree $720$ and decomposes it into a product of $30$ factors of degree $24$ (the totality in one or two minutes).
One of these factors is the \polz~{\tt q}. As {\tt q} is very cumbersome, we ask {\tt Magma} to compute a \bdg of the \id generated by the Cauchy modules on the one hand, and by~{\tt q(z)} on the other, which provides a clearer description of the \cdr $\aqo{\gA}{q(z)}$:
{\tt x6} is annihilated by {\tt p}, {\tt x5} is annihilated by a \pol of degree $4$ with \coes in $\QQ[x_6]$, {\tt x1}, \ldots, {\tt x4} are expressed in terms of {\tt x5} and {\tt x6}. The computation of the \bdg takes several hours. {\tt Magma} can then compute the Galois group, which is given by two \gtrsz.
Here are the results:

{\small \def\baselinestretch{1}\label{exemple1Galois}
\begin{verbatim}
p:=T^6 - 3*T^5 + 6*T^4 - 7*T^3 + 2*T^2 + T - 1;
z:=x1 + 2*x2 + 13*x3 - 24*x4 + 35*x5 - 436*x6;
pm:=T^720 + 147240*T^719 + 10877951340*T^718 + 537614218119000*T^717 +
    19994843992714365210*T^716 + 596880113924932859498208*T^715 +
    14896247531385087685472255280*T^714 + ...
q:= T^24 + 4908*T^23 + 13278966*T^22 + 25122595960*T^21 +
    36160999067785*T^20 + 41348091425849608*T^19 +
    38304456918334801182*T^18 + 28901611463650323108996*T^17 +...
//we annihilate q(z): description of the field of roots;
Affine Algebra of rank 6 over Rational Field
Variables: x1, x2, x3, x4, x5, x6
Quotient relations:
x1 + 18/37*x5^3*x6^5 - 45/37*x5^3*x6^4 + 104/37*x5^3*x6^3 - 3*x5^3*x6^2
    + 36/37*x5^3*x6 - 1/37*x5^3 - 27/37*x5^2*x6^5 + 135/74*x5^2*x6^4 -
    156/37*x5^2*x6^3 + 9/2*x5^2*x6^2 - 54/37*x5^2*x6 + 3/74*x5^2 +
    91/37*x5*x6^5 - 455/74*x5*x6^4 + 460/37*x5*x6^3 - 25/2*x5*x6^2 +
    108/37*x5*x6 + 31/74*x5 - 41/37*x6^5 + 205/74*x6^4 - 204/37*x6^3 +
    11/2*x6^2 - 45/37*x6 - 53/74,
x2 + x6 - 1,
x3 + x5 - 1,
x4 - 18/37*x5^3*x6^5 + 45/37*x5^3*x6^4 - 104/37*x5^3*x6^3 + 3*x5^3*x6^2
    - 36/37*x5^3*x6 + 1/37*x5^3 + 27/37*x5^2*x6^5 - 135/74*x5^2*x6^4 +
    156/37*x5^2*x6^3 - 9/2*x5^2*x6^2 + 54/37*x5^2*x6 - 3/74*x5^2 -
    91/37*x5*x6^5 + 455/74*x5*x6^4 - 460/37*x5*x6^3 + 25/2*x5*x6^2 -
    108/37*x5*x6 - 31/74*x5 + 41/37*x6^5 - 205/74*x6^4 + 204/37*x6^3 -
    11/2*x6^2 + 45/37*x6 - 21/74,
x5^4 - 2*x5^3 + x5^2*x6^2 - x5^2*x6 + 4*x5^2 - x5*x6^2 + x5*x6 - 3*x5 +
     x6^4 - 2*x6^3 + 4*x6^2 - 3*x6 - 1,
x6^6 - 3*x6^5 + 6*x6^4 - 7*x6^3 + 2*x6^2 + x6 - 1
// the Galois group;
Permutation group acting on a set of cardinality 6
Order = 24 = 2^3 * 3
    (1, 4)(2, 5)(3, 6)
    (1, 2, 4, 6)
\end{verbatim}
}

Note that $\disc_T(p)=2^4\times 37^3$, which is not unrelated to the \denos appearing in the \bdgz. We will return to this example on \paref{example2Galois} when discussing the dynamic method. 
\eoe

\medskip
\rem
Here we interrupt our treatment of basic Galois theory. We shall resume the current thread in Sections VI-1 and VI-2, which the reader can refer to directly from here (the results of the intermediate chapters will not be used). In Chapter VII we will address a more sophisticated theory which proves to be necessary when we do not have at our disposal a factorization algorithm for the separable polynomials over the base field.
\eoe

\section{The resultant}
\label{secRes}

The resultant is the basic tool of Elimination theory. This is based on the basic Elimination lemma on \paref{LemElimAffBasic}, which is applied to arbitrary rings, and on its Corollary~\ref{cor2LemElimAffBasic} for the geometric case.

\subsec{\Eli theory}

\Eli theory concerns the \syss of \pol \eqns (or \emph{\sypsz}).\index{polynomial systems}\index{elimination!theory}

Such a \sys $(\lfs)$ in $\kXn=\kuX$, where $\gk$ is a \cdiz, can admit some zeros in $\gk^{n}$, or in $\gL^{n}$, where $\gL$ is an overfield of~$\gk$, or even  an arbitrary \klgz.
The zeros depend only on the \id $\fa=\gen{\lfs}$ of $\kuX$ generated by the $f_i$'s. We also call them \emph{the zeros of the \id $\fa$}.

Let $\pi:\gL^{n}\to\gL^{r}$ be the projection which forgets the last $n-r$ \coosz.
If $V\subseteq \gL^{n}$ is the set of zeros of $\fa$ on $\gL$, we are 
interested in as precise a description as possible of the projection
 $W=\pi(V)$, if possible as zeros of a \syp in the variables $(X_1,\dots,X_{r})$.

Here intervenes in a natural way the \emph{\eli \idz} (\eli of the variables $X_{r+1}$, \dots, $X_n$ for the considered \sypz), which is defined by $\fb=\fa\cap\kXr$.
Indeed every \elt of $W$ is clearly a zero of $\fb$.

The converse is not always true (and in any case not at all obvious), but it is true in some good cases: if $\gL$ is an \cac and if the \id is in a 
\noep
(\Thref{thNstfaibleClass}).

A reassuring fact, and easy to establish by the considerations of \lin \alg over  \cdisz, is that the  \eli  \id $\fb$ \gui{does not depend on} the considered base field $\gk$. More \prmtz, if $\gk_1$ is an overfield of $\gk$, we have the following results.
\begin{itemize}
\item The \id $\gen{\lfs}_{\gk_1[\Xn]}$ only depends on the \id $\fa$: \\
it is the \id $\fa_1$ of~\hbox{$\gk_1[\Xn]$} generated by $\fa$. 
\item The \id of \eli $\fb_1=\fa_1\cap\gk_1[\Xr]$ only depends on $\fb$: 
\\
it is the \id of $\gk_1[\Xr]$ generated by $\fb$. 
\end{itemize}

\medskip 
\Elr \Eli theory faces two obstacles.

The first is the difficulty of computing $\fb$ from $\fa$,
\cad of computing some finite \sgr of $\fb$ from the \syp $(\lfs)$.
This computation is rendered possible by the theory of the \bdgsz, which we do not address in this work. In addition this computation is not uniform, unlike the computations linked to resultant theory.

The second obstacle is that one only obtains truly satisfactory results for \hmg \sypsz.
The basic example that shows this is the \deterz. 
Consider a \gnq \sli $(f_1,\dots,f_n)$ \hbox{of $\gk[\ua][\uX]$},
where the variables $a_{ij}$ in $\ua$ represent the~$n^{2}$ \coes of the~$n$ \lin forms $f_i$, and the $X_j$'s are the unknowns.
Then the \id $\gen{\det(\ua)}$ \hbox{of $\gk[\ua]$} is indeed the \eli \id  of the variables $X_j$ for the \sys $(f_1,\dots,f_n)$, provided we only take into account the zeros of the \sys distinct from $\uze=(0,\dots,0)$.
\\
The simplicity of this result should be contrasted with the discussion, in the non-homogeneous framework, of \syss where the $f_i$'s are affine forms.
\\
Furthermore, even though the zeros of the \id $\gen{\det(\ua)}$ correspond effectively to the \syss that admit some zero $\neq \uze$, this \id is not exactly equal to $\gen{f_1,\dots,f_n}\cap \gk[\ua]$, we first need to \emph{saturate} the \id $\fa=\gen{f_1,\dots,f_n}$ w.r.t.\,the \hmg variables $X_j$; \cad add every $g$ to it such that, for each $j\in\lrbn$, $gX_j^{N}\in\fa$ for some large enough $N$. 
In the current case, this saturated \id is the \id $\fa+\det(\ua)\gk[\ua][\uX]$, each $\det(\ua)X_j$ is in $\fa$, and the intersection of the saturation with $\gk[\ua]$ is indeed $\gen{\det(\ua)}$.    

What will be retained from this little introduction to \Eli theory is a \dfnz: 
let $\gk$ be a commutative \riz, $\fa$ be an \id of $\kXn$ \hbox{and $r\in\lrb{0..n-1}$}, we then define the \emph{\eli \id of the variables $X_{r+1}$, \dots, $X_n$ for the \id $\fa$} as being the \id $\fb=\fa\cap\kXr$.\index{elimination!ideal}\index{ideal!elimination ---}

We will remain wary of the fact that if $\gk$ is an arbitrary \riz, the \id $\fa$ can very well be \tf even if $\fb$ is not \tfz.

\subsec{The Sylvester matrix}

In what follows, we do not assume the \ri $\gA$ to be discrete, so much so that the degree of a \pol of $\AX $ is not \ncrt exactly known. From the point of view of computation, in \gnl we have to take the \pols in~$\AX$ in the form of \emph{formal \polsz}, \cad pairs $(f,p)$ where $f$  is a \pol and $p$ is the upper bound of its degree. This notion is \egmt useful when changing the base \ri because a \pol can for instance have its degree decrease without us knowing how to test it (e.g. upon passage to the quotient \riz).%
\index{polynomial!formal ---}

Recall the \dfn of the Sylvester matrix and of the resultant of two \pols (formal \pols of degrees $p$ and $q\geq0$):
$$\arraycolsep2pt
\begin{array}{lcl} 
f & = & a_pX^p+ a_{p-1}X^{p-1}+\cdots+a_0, \\[1mm] 
g & = &b_{q}X^{q}+b_{q-1}X^{q-1}+\cdots+b_0 .
\end{array}
$$
The \emph{Sylvester matrix} of $f$ and $g$ (in degrees $p$ and $q$) is the following matrix \index{matrix!Sylvester ---}\index{Sylvester!matrix}
$$
\Syl_X(f,p,g,q)=
\underbrace{
\left[
\begin{array}{ccccccccccc}
 a_p&\cdots&\cdots&\cdots&\cdots&a_0&& \\[.5mm]
&\ddots& & & & &\ddots& \\[.5mm]
 &&a_p&\cdots&\cdots&\cdots&\cdots&a_0\\[.8mm]
b_{q}&\cdots& \cdots &b_0&&&& \\[.5mm]
 &\ddots& & &\ddots&&& \\[3.5mm]
 &&&\ddots& & &\ddots& \\[.5mm]
 &&&&b_{q}&\cdots&\cdots& b_0  \\[.5mm]
\end{array}
 \right]
}_{p+q}
\matrix{ \left.\matrix{\cr \cr \cr \cr}\right\}&{q} \cr
 \left.\matrix{\cr \cr \cr \cr \cr\cr}\right\}&{p}
\cr}
$$
This matrix can be regarded as the matrix whose rows are the \coos of the \pols $(X^{q-1}f, \ldots , Xf, f, X^{p-1}g, \ldots , Xg, g)$ over the 
{basis $(X^{p+q-1},X^{p+q-2},\ldots ,X,1)$}.

\rdb
The \emph{resultant of $f$ and $g$ (in degrees $p$ and $q$)}, denoted by $\Res_{X}(f,p,g,q)$, is the \deter of this Sylvester matrix\index{resulant!of two polynomials}
\begin{equation}\label{eqResultant}\preskip.3em \postskip.4em
\Res_{X}(f,p,g,q)\eqdefi\det\big(\Syl_X(f,p,g,q)\big).
\end{equation}
If the context is clear, we also denote it by $\Res_X(f,g)$
or $\Res(f,g)$.
We have
\begin{equation}\preskip.4em \postskip.2em
\label{eqres1.0}
{\Res_{X}(f,p,g,q)=(-1)^{pq}\Res_{X}(g,q,f,p)},
\end{equation}
and also, for $a,b\in\gA$,
\begin{equation}\preskip.2em \postskip.3em
\label{eqres1.1}
\Res_{X}(af,p,bg,q) = a^qb^p\Res_{X}(f,p,g,q).
\end{equation}

If $p=q=0$, we obtain the \deter of an empty matrix, \cad $1$. 

\smallskip 
\begin{fact}\label{fact1Res} 
If $p\geq 1$ or $q\geq1$, then
$\Res_X(f,p,g,q)\in\gen{f,g}_{\AX }\cap\gA$. More \prmtz,
for each $n\in\lrb{0..p+q-1}$, there exist $u_n$ and $v_n\in\AX$ such that $\deg u_n<q$, $\deg v_n<p$ and 
\begin{equation}\label{eqResultant2}\preskip.3em \postskip.4em
X^n\,\Res_X(f,g)=u_n(X)f(X)+v_n(X)g(X).
\end{equation}
\end{fact}
\begin{proof}
Let $S$ be the transpose of $\Syl_X(f,p,g,q)$. The columns of $S$ express \pols $X^kf$ or $X^\ell g$ over the basis of the \moms of degree $<p+q$. By using Cramer's formula

\snic{S \,\wi S=\det S\cdot\rI_{p+q}\;,}

we see that each $X^n\Res(f,g)$ (which corresponds to one of the columns of the right-hand side matrix) is a \coli of the columns of~$S$.
\end{proof}

\rem We can also view \Egrf{eqResultant2} in the $n=0$ case as expressing the \deter of the matrix below developed according to the last column (this is in fact the Sylvester matrix in which we have replaced each \coe in the last column by the \gui{name} of its row):
$$
\left[
\begin{array}{ccccccccccc}
 a_p&\cdots&\cdots&\cdots&\cdots&a_0&& X^{q-1}f\\[.5mm]
&\ddots& & & & &\ddots& \\[.5mm]
 &&a_p&\cdots&\cdots&\cdots&\cdots&f\\[2mm]
b_q&\cdots& \cdots &b_0&&&&X^{p-1}g \\[.5mm]
 &\ddots& &&\ddots&&& \\[5.5mm]
 &&&b_q&\cdots&\cdots&b_0& Xg\\[1.5mm]
 &&&&b_{q}&\cdots&\cdots& g  
\end{array}
 \right]. $$
\eoe

\begin{corollary}\label{lem0Resultant}
Let $f$, $g\in\AX$ and $a\in\gB\supseteq\gA$, with $f(a)=g(a)=0$, and~$p\geq1$ or~$q\geq1$, then $\Res_X(f,p,g,q)=0$.
\end{corollary}

Note that if the two degrees are over-evaluated the resultant is annihilated, and the intuitive interpretation is that the two \pols have a common zero \gui{at infinity.}
Whilst if $a_p=1$, the resultant (for $f$ in degree $p$) is the same regardless of the formal degree chosen for $g$. This then allows for an unambiguous switch to the notation $\Res(f,g)$, as in the following lemma.

\begin{lemma}
\label{lemResultant}
Let $f$ and $g\in\AX$ with $f$ \mon of degree $p$.
\begin{enumerate}
\item We write $\gB=\aqo{\AX}{f}$ and denote by $\mu_g$  multiplication by (the class of) $g$ in~$\gB$, which is a free \Amo of rank $p$.
Then
\begin{equation}\preskip.4em \postskip.0em
\label{eq0lemResultant}
\rN\iBA (g)=\det\mu_g=\Res(f,g).
\end{equation}
\item Therefore
\vspace{-.7em}
\begin{eqnarray}\label{eqResultant4}\postskip.4em
\Res(f,gh)&=&\Res(f,g)\,\Res(f,h), \\ 
\Res(f,g+fh)&=&\Res(f,g).\label{eqResultant5}
\end{eqnarray}
%
\item For every square matrix $A\in\MM_p(\gA)$ for which the \polcar is equal to $f$, we have
\begin{equation}\preskip.2em \postskip.4em
\label{eqlemResultant}
\Res(f,g)=\det\big(g(A)\big).
\end{equation}
\item If we write $f=\prod_{i=1}^p(X-x_i)$ in an extension of $\gA$, we obtain
\begin{equation}\label{eqResultant3}\preskip.4em 
\Res(f,g)= \prod\nolimits_{i=1}^pg(x_i).
\end{equation}
%
\end{enumerate}
\end{lemma}

\vspace{-.7em}
\pagebreak

\begin{proof}
\emph{1.} By \elr manipulations of rows, the Sylvester matrix
$$
\Syl_X(f,p,g,q)=
\left[
\begin{array}{ccccccccccc}
 1&a_{p-1}&\cdots&\cdots&\cdots&a_0&& \\[.5mm]
&\ddots&\ddots& & & &\ddots& \\[.5mm]
 &&1&a_{p-1}&\cdots&\cdots&\cdots&a_0\\[2mm]
b_{q}&\cdots& \cdots &b_0&&&& \\[.5mm]
 &\ddots& & &\ddots&&& \\[3.5mm]
 &&&\ddots& & &\ddots& \\[.5mm]
 &&&&b_{q}&\cdots&\cdots& b_0
\end{array}
 \right]
$$
is transformed into the matrix visualized below, in which the rows {$q+1$, $\ldots$, $q+p$} now contain the remainders of the division by $f$ of the \pols $X^{p-1}g$, $\ldots$, $Xg$, $g$.
Thus the $p\times p$ matrix in the south-east corner is exactly the transpose of the matrix of the \endo $\mu_g$ of $\gB$ over the basis of the \moms
and its \deter is equal to that of the Sylvester matrix.
$$
\left[
\begin{array}{ccccccccccc}
 1&a_{p-1}&\cdots&\cdots&\cdots&a_0&& \\[.5mm]
&\ddots&\ddots& & & &\ddots& \\[.5mm]
 &&1&a_{p-1}&\cdots&\cdots&\cdots&a_0\\[2mm]
0&\cdots& 0 &\times &\cdots&\cdots&\cdots& \times \\[.5mm]
\vdots &&\vdots&\vdots& &&&\vdots \\[3.5mm]
\vdots &&\vdots&\vdots& &&&\vdots \\[.5mm]
0&\cdots& 0 &\times &\cdots&\cdots&\cdots& \times
\end{array}
 \right]
$$

\emph{2.} Results from item \emph{1.}
 
\emph{3} and \emph{4.} Result from Proposition~\ref{prop2tschir} via item~\emph{1.}

We can also give the following direct \demsz.

\emph{4.} First of all, from \Eqrf{eqResultant4} we deduce the \smqz al formula
\[
 \Res(f_1f_2,g)=\Res(f_1,g)\;\Res(f_2,g)
\] 
for $f_1$ and $f_2$ \mons (use the \eqnsz~(\ref{eqres1.0}) and (\ref{eqres1.1}) and the fact that in the case where the \coes of $g$ are \idtrs we can assume $g=b_qg_1$ with $g_1$ \monz). Next, a direct computation gives $\Res(X-a,g)=g(a)$.

\emph{3.} We must prove $\Res(\rC A,g)=\det\big(g(A)\big)$ for some \polz~$g$ and an arbitrary matrix $A$. This is an \ida concerning the \coes of~$A$ and of~$g$. We can thus restrict ourselves to the case where the matrix $A$ is the \gnq matrix. Then, it is diagonalized in an over\ri and we conclude by applying item~\emph{4.}
\end{proof}

\rem \label{remreciproquenonnegligeable} Item~\emph{4} offers a non-negligible converse to Corollary~\ref{lem0Resultant}:
if~$\gA$ is integral and if $f$ and $g$ are two \polus of~$\AT$ that are completely factorized in an integral \ri containing~$\gA$, they have a common zero \ssi their resultant is null. \eoe

\medskip 
In the case of a nontrivial \cdi $\gK$ we can do a little better. 

\begin{fact}\label{factResNulPGCD}
Let $f$ and $g\in\KX$ of degrees $p$ and $q\geq1$, with $\Res(f,g)=0$. Then, $f$ and $g$ have a gcd of degree $\geq1$.
\end{fact}
\begin{proof}
The \Kli $(u,v)\mapsto uf+vg$ where $\deg u<q$ and $\deg v<p$ admits as matrix over the bases of \moms the transpose of the Sylvester matrix.
Thus let $(u,v)\neq(0,0)$ in the kernel. The \polz~$uf=-vg$ is of degree $<p+q$. So $\deg \big(\ppcm(f,g)\big)<p+q$, which implies $\deg \big(\pgcd(f,g)\big)>0$.
\end{proof}

\comm The above \dem assumes that we know the \elr theory of \dve (via Euclid's \algoz) in the \ris of type $\KX$.
This theory shows the existence of a gcd and of a lcm with the relation
$$\preskip.4em \postskip.4em 
\qquad\qquad\ppcm(f,g)\, \pgcd(f,g)=\alpha fg, \qquad (\alpha\in\gK\eti). 
$$
Another \dem would consist in saying that in a \cdi $\gL$, which is an extension of $\gK$, the \pols $f$ and $g$ are split (i.e., are decomposed into factors of degree $1$) which implies, given the previous remark, that $f$ and $g$ have a common zero and thus a common factor of degree $>0$.
One must then finish by stating that the gcd is computed by Euclid's \algo and thus does not depend on the chosen base field (which must only contain the \coes of $f$ and $g$). Nevertheless this second \demz, which somewhat gives \gui{the true motivation for the \thoz,}
assumes the existence of $\gL$ (which is not guaranteed from a \cof point of view) and does not avoid the theory of \dve in~$\KX$ via Euclid's \algoz.
\eoe

\CMnewtheorem{lemelibas}{Basic \eli lemma}{\itshape}
\begin{lemelibas}\label{LemElimAffBasic}%
\index{Basic elimination lemma}%
\index{elimination!basic --- lemma}%
\index{elimination!ideal}\index{ideal!elimination ---}\\
Let $f$ and $g\in\AX $ with $f$ \mon of degree $p$. Then,
$R=\Res_X(f,g)$ is well defined and the \eli \id $\fa=\gen{f,g}_{\AX}\cap\gA$ satisfies
$$\preskip.3em \postskip.0em 
\fa^p\; \subseteq\; \Res_X(f,g)\gA\;\subseteq\; \fa. 
$$
In particular
\begin{enumerate}\itemsep2pt
\item $R$ is \iv \ssi $1\in\gen{f,g}$,
\item $R$ is \ndz \ssi $\fa$ is faithful, and 
\item $R$ is nilpotent
\ssi $\fa$ is nilpotent.
\end{enumerate}
 
\end{lemelibas}

\begin{proof}
We already know that $\Res_X(f,g)\in\gen{f,g}_{\AX}$.\\
We use the notations of Lemma~\ref{lemResultant}, item~\emph{1.}
We denote by $x$ the class of~$X$ in~$\gB=\aqo{\AX}{f}$.
A basis of $\gB$ over $\gA$ is $(1,x,\dots,x^{p-1})$.
Let $(\gamma_i)_{i\in\lrbp} $ be \elts of $\fa$. The \elts $\gamma_1$, $\gamma_2x$, $\ldots$, $\gamma_px^{p-1}$ are in~$\Im\mu_g$,
 so the matrix $D=\Diag(\gamma_1,\ldots,\gamma_p)$ can be written in the form $GB$, where~$G$ is the matrix of $\mu_g$ over the basis of the \momsz.
 It follows that
 $$\preskip.4em \postskip.4em\ndsp
 \prod\nolimits_{k=1}^p\gamma_k\,=\,\det D\,=\,\det G\,\det B\,=\,\Res(f,g)\,\det B.
 $$
 Thus the \elt $\prod\nolimits_{k=1}^p\gamma_k$
 of $\fa^p$ belongs to $\gen{\Res(f,g)}_\gA$.
\end{proof}

The basic \eli lemma will be \gne later (Lemmas~\ref{lemElimPlusieurs}
and \ref{LemElimAff}).
The term \gui{\eli \idz} comes from the following facts which result from the previous lemma and from Lemma \ref{lemResultant}.
\begin{corollary}\label{corLemElimAffBasic}
Let $\gA$ be an integral \ri and $f$, $g\in\AX$. If $f$ is \mon and can be completely factorized, \propeq
\begin{enumerate}\itemsep2pt
\item The \eli \id \smash{$\gen{f,g}_{\AX}\cap\gA$} is null.
\item The resultant $\Res_X(f,g)=0$.
\item The \pols $f$ and $g$ have a common root.
\end{enumerate}

\end{corollary}

 A \cdi $\gK$ is said to be \emph{\agqt closed} if every \polu of~$\KX$ can be decomposed into a product of factors $X-x_i$ ($x_i\in\gK$).%
\index{field!algebraically closed discrete ---}%
\index{algebraically closed!discrete field}

\begin{corollary}\label{cor2LemElimAffBasic}
Let $\gK$ be a \cdacz.
\\  
Write $\gA=\KYm$. Let $f$ and $g\in\AX$ with $f$ \mon in~$X$. 
For some arbitrary \elt $\uzeta=(\zeta_1,\ldots,\zeta_m)$ of $\gK^m$, \propeq
\begin{enumerate}
\item $\uzeta$ annihilates all the \pols of the \eli \id $\gen{f,g}\cap{\gA}$.%
\index{elimination!ideal}\index{ideal!elimination ---}
\item $\Res_X\big(f(\uzeta,X),g(\uzeta,X)\big)=0$.
\item $f(\uzeta,X)$ and $g(\uzeta,X)$ have a common root.
\end{enumerate}
Consequently if $V$ is the set of zeros common to $f$ and $g$ in~$\gK^{m+1}$, and if $\pi:\gK^{m+1}\to\gK^{m}$ is the projection that forgets the last coordinate, then $\pi(V)$ is the set of zeros of $\Res_X(f ,g)\in\KYm.$
\end{corollary}

\subsec{Revisiting the \discri}

When $g=\prod_{i=1}^n(X-y_i)$, Lemma~\ref{lemResultant} gives $\Res_X(g,g')=\prod_{i=1}^ng'(y_i)$ and thus
\begin{equation}\preskip-.2em \postskip.5em
\label{eqDiscriRes}
\disc(g)=(-1)^{n(n-1)/2}\Res_X(g,g').
\end{equation}
Since the \egt $g(X)=\prod_{i=1}^n(X-y_i)$ can always be performed in the \adu if $g$ is \monz, we obtain that \Egrf{eqDiscriRes} is valid for every \poluz, over every commutative \riz.

The following fact results therefore from the basic \eli lemma.

\begin{fact}
\label{factDiscUnit}
Consider some \polu $g\in\AX$.
\vspace{-1pt}
\begin{enumerate}
\item [--] $\gen{g(X),g'(X)}=\gen{1}$  \ssi  $\disc g$ is \ivz.
\item [--]  The \id $\gen{g(X),g'(X)}\cap\gA$ is faithful \ssi $\disc g$ is a
 \ndz\elt of $\gA$.
\end{enumerate}
\end{fact}

\vspace{-.5em}
\pagebreak

\begin{fact}\label{factDiscProd}
If $f=gh\in\AX$ with $g$, $h$ \monsz, we have the following \egt
\begin{equation}\preskip-.2em \postskip.4em
\label{eqfactDiscProd}
\disc(f)=\disc(g)\disc(h)\Res(g,h)^2
\end{equation}
\end{fact}

\smallskip 
\begin{proof}
This \imdt results from \Eqns (\ref{eqResultant4}), (\ref{eqResultant5}) \paref{eqResultant5} and (\ref{eqDiscriRes}).
\end{proof}

\begin{corollary}\label{corfactDiscProd} 
Let $f\in\AX$ be \mon and $\gB=\gA[x]=\aqo{\AX}{f}$.
\begin{enumerate}
\item If $f$ possesses a square factor, $\disc f =0$. Conversely, 
if $\disc f =0$ and if $f(X)=\prod(X-x_i)$ in some integral \ri containing $\gA$, two of the zeros $x_i$ are equal.
\item \label{icorfactDiscProd} Assume $f$ is \spl and $f=gh$ ($g$ and $h$ \monsz). 
\begin{enumerate}
\item The \pols $g$ and $h$ are \spls and \comz.
\item
There exists some \idm $e$ of $\gB$ such that $\gen{e}=\gen{\pi(g)}$.\\
We have $\gB\simeq\aqo{\gB}{g}\times \aqo{\gB}{h}$.
\end{enumerate}
\item Assume $\disc f$ is \ndz and $f=gh$ ($g$ and $h$ \monsz). \\
Then, the \elts $\disc g$, $\disc h$ and $\Res(g,h)$ are \ndzsz.
\end{enumerate}
\end{corollary}
\begin{proof}
All this results from Fact~\ref{factDiscProd}, except maybe the \idm $e$ in item~\emph{\ref{icorfactDiscProd}}. If~$gu+hv=1$, then $e=\ov{gu}$ is required.
\end{proof}
%

\begin{corollary}\label{corcorfactDiscProd}
Let $\gK$ be a \cdiz, $f\in\KX$ a \spl \polu and $\gB=\aqo{\KX}{f}$.
In item~\ref{icorfactDiscProd} of the previous corollary, we associate with every divisor $g$ of $f$ the \idm $e$ such that $\gen{\ov g}=\gen{e}$. This establishes a bijection between the \mons divisors of $f$ and the \idms of~$\gB$. This bijection respects \dvez.
\end{corollary}
\begin{proof}
The reciprocal bijection is given by $e=\ov v\mapsto\pgcd(v,f)$.
\end{proof}

We now introduce the notions of prime subfields and of \cara of a \cdiz. \rdb

More \gnltz, if $\gA$ is an arbitrary \riz, we denote by $\Z_\gA$ the \emph{prime sub\ri of} $\gA$ defined as follows:
$$
\Z_\gA=\sotq{n \cdot (m \cdot 1_A)^{-1}}{n,m\in \ZZ,\;m\cdot1_\gA\in\Ati}
.\label{NOTAZA}
$$
If $\rho:\ZZ\to\gA$ is the unique \homo of \ris of $\ZZ$ in $\gA$, the prime sub\ri is therefore \isoc to $S^{-1}\ZZ\sur{\Ker\rho}$, where $S=\rho^{-1}(\Ati)$. A \ri can be called prime if it is equal to its prime sub\riz.
Actually the terminology is only common in the case of fields.\rdb

When $\gK$ is a \cdiz, the prime sub\ri is a subfield, called the \emph{prime subfield of~$\gK$}.
For some $m>0$ we will say that \emph{$\gK$ is of \cara $>m$}, and we write \gui{$\car(\gK)>m$}
if for every  
$n\in\lrbm$, the \elt $n\cdot1_\gK$ is \ivz.\label{NOTACarK}%
\index{prime!subring of a ring}%
\index{characteristic!of a field}%
\index{prime!subfield of a field}%
\index{field!prime ---}

When $\gK$ is nontrivial, if there exists some $m>0$ such that $m\cdot1_\gK=0$, 
then there is a minimum number of them,
which is a prime number $p$, and we say that the field \emph{is of \cara $p$}.
When the prime subfield of $\gK$ is \isoc to $\QQ$, the convention is to speak of a \emph{null \caraz}, but we will also use the terminology \emph{infinite \caraz} in the contexts where it is useful to remain consistent with the previous notation, for instance in Fact~\ref{factPolSepFC}.

We can conceive\footnote{We can also be presented with such cases resulting from a complicated construction in a subtle \demz.} some nontrivial \cdis whose \cara is not well defined from a \cof point of view. However, for a \cdi the statement \gui{$\car(\gK)>m$} is always decidable.

\begin{fact}\label{factPolSepFC}
Let $\gK$ be a \cdi and $f\in\KX$ be a \poluz.
If~$\disc f=0$ and $\car(\gK)>\deg f$, $f$ possesses a square factor 
of degree~\hbox{$\geq1$}.
\end{fact}

\begin{proof}
Let $n=\deg f$. The \polz~$f'$ is of degree $n-1$.
Let $g=\pgcd(f,f')$. We have $\deg g \in\lrb{1..n-1}$ (Fact~\ref{factResNulPGCD}).
We write $f=gh$ therefore
$$\preskip.4em \postskip.4em 
\disc(f)=\Res(g,h)^2\disc(g)\disc(h). 
$$
Thus, $\Res(g,h)=0$, or $\disc(g)=0$, or $\disc(h)=0$. In the first case, the \pols $g$ and $h$ have a gcd $k$ of degree $\geq1$ and $k^2$ divides $f$.
In the two other cases, since $\deg g<\deg f$ and $\deg h<\deg f$, we can finish by \recu on the degree, by noting that if $\deg f=1$, then $\disc f\neq0$, which assures the initialization.
\end{proof}
%

\section{\Agq number theory, first steps}
\label{secApTDN}

Here we give some \gnl applications, in \elr number theory, of the results previously obtained in this chapter. For a glimpse of the many fascinating facets of number theory, the reader should consult the wonderful book \cite{IR}.

\subsec{Integral \algsz} 

We give a few precisions relating to \Dfnz~\ref{defEntierAnn0}.

\begin{definition}
\label{def05Alg}~
\begin{enumerate}
\item
An \Alg $\gB$ is said to be
 \emph{finite}
if $\gB$ is a \tf \Amoz. We also say that \emph{$\gB$ is finite over $\gA$}. 
In the case of an extension, we speak of a \emph{finite extension} of $\gA$.%
\index{algebra!finite ---}
\item  Assume $\gA\subseteq\gB$. The \riz~$\gA$ is said to be
 \ix{integrally closed} in $\gB$ if every \elt of $\gB$ integral over $\gA$ is in~$\gA$.\index{ring!integrally closed --- in \dots}
\end{enumerate}
\end{definition}

\begin{fact}
\label{factEntiersAnn}
Let $\gA\subseteq\gB$ and $x\in\gB.$ \Propeq
\begin{enumerate}
\item The \elt $x$ is integral over $\gA.$
\item The sub\algz~$\gA[x]$ of $\gB$ is finite.
\item There exists a faithful and \tf \Amo 
$M\subseteq\gB$ such that $xM\subseteq M$.
\end{enumerate}
\end{fact}
\begin{proof}
\emph{3} $\Rightarrow$~\emph{1} (a fortiori \emph{2} $\Rightarrow$~\emph{1}.)
Consider a matrix $A$ with \coes in~$\gA$ which represents $\mu_{x,M}$ (the multiplication by $x$ in $M$) on a finite \sgr of $M$. If $f$ is the \polcar of $A$, we have by the Cayley-Hamilton \tho  $0=f(\mu_{x,M})=\mu_{f(x),M}$ and since the module is faithful, $f(x)=0$. \\
The rest is left to the reader.
\end{proof}

We also easily obtain the following fact.
\begin{fact}
\label{factEntiersAnn2}
Let $\gB$ be an \Alg and $\gC$ be a \Blgz.
\begin{enumerate}
\item If $\gC$ is finite over $\gB$ and $\gB$ finite over $\gA$, then $\gC$ is finite over $\gA$.
\item An \Alg generated by a finite number of integral \elts over $\gA$ is finite.
\item The \elts of $\gB$ that are integral over $\gA$ form a \ri  \icl in~$\gB$. 
We call it the \emph{integral closure of $\gA$ in $\gB$}.%
\index{integral closure!of $\gA$ in $\gB\supseteq\gA$}
\end{enumerate}
\end{fact}

\begin{lemma}\label{lemPolEnt}
Let $\gA \subseteq \gB$ and $f \in \BuX$. The \pol $f$ is integral over~$\AuX$ \ssi each \coe of $f$ is integral over $\gA$.
\end{lemma}
\begin{proof} The condition is sufficient, by item~\emph{3}
of the previous lemma. In the other direction consider an \rdi $P(f)=0$ for $f$
(with $P\in\AuX[T]$, \monz). We have in $\gB[\uX,T]$ an \egt
$$
P(\uX,T)=\,\bigr(T - f(\uX) \bigl)\,
\bigr(T^n + u_{n-1}(\uX) T^{n-1} + \cdots + u_{0}(\uX) \bigl)
.$$
Since the \coe of $T^n$ in the second factor is $1$, the multivariate \KROz's \tho 
 implies that each \coe of $f$ is integral over~$\gA$.
\end{proof}
%

\begin{lemma}\label{lemPolCarInt}
Let $\gA \subseteq \gB$, $L$ be a free \Bmo of finite rank and $u \in \End_\gB(L)$ be integral over~$\gA$. Then, the \coes of the \polcar of~$u$ are integral over~$\gA$. In particular, $\det( u)$ and $\Tr(u)$ are integral over~$\gA$.
\end{lemma}
\begin{proof}
Let us first prove that $\det(u)$ is integral over~$\gA$.
Let $\cE=(e_1,\ldots,e_n)$ be a fixed basis of $L$.
The \Amo $\gA[u]$ is a \tf \Amoz, and so the module
$$\preskip.2em \postskip.4em\ndsp
E = \som_{i\in\lrbn, k\geq 0} \gA u^k(e_i) \subseteq L
$$
is a \tf \Amoz, with $u(E) \subseteq E$. Let us introduce the module
$$\preskip-.2em \postskip.4em\ndsp
D = \som_{\ux \in E^n}
\gA\det_{\cE} (\ux) \subseteq \gB.
$$
Since $E$ is a \tf \Amoz, $D$ is a \tf \Amoz, and it is faithful, we have $1 \in D$ as $\det_{\cE} (\cE)=1$.
Finally, the \egtz
$$\preskip.4em \postskip.4em
\det(u) \det_{\cE} (x_1, \ldots, x_n) =
\det_{\cE} \big(u(x_1), \ldots, u(x_n)\big)
$$
and the fact that $u(E) \subseteq E$ show that $\det(u) D \subseteq D$.\\
Next consider $\AX\subseteq\BX$ and the $\BX$-module $L[X]$.\\ We have $X\Id_{L[X]} - u \in \End_{\gB[X]}(L[X])$.
If $u$ is integral over~$\gA$, $X\Id_{L[X]} - u$ is integral over~$\gA[X]$ therefore $\rC{u}(X) =
\det(X\Id_{L[X]} - u)$ is integral over $\gA[X]$. We conclude with Lemma~\ref{lemPolEnt}.
\end{proof}

\begin{corollary}\label{corlemPolcarEntier}
Let $\gA\subseteq\gB\subseteq\gC$ 
          where $\gC$ is a finite free \Bmoz.
Let $x\in\gC$ be integral over $\gA$. Then, $\Tr_{\gC/\gB}(x)$, $\rN_{\gC/\gB}(x)$ and all the \coes of $\rC{\gC/\gB}(x)$ are integral over $\gA$.
If in addition $\gB$ is a \cdiz, the \coes of the \polmin $\Mip_{\gB,x}$ are integral over~$\gA$.
\end{corollary}
\begin{proof}
We apply the previous lemma with $L=\gC$ and $u=\mu_x$.
For the final statement, we use \KROz's \tho and the fact that the \polmin divides the \polcarz.
\end{proof}
%

%
\subsubsection*{Integrally closed \risz}

\begin{definition}\label{defiIntClos}
An integral \ri $\gA$ is said to be \ix{integrally closed} if it is \icl in its quotient field.%
\index{ring!integrally closed ---}
\end{definition}

\begin{fact}
\label{fact.loc.entier}
Let $\gA\subseteq\gB$, $S$ be a \mo of $\gA$, $x\in\gB$ and $s\in S$.
\begin{enumerate}
\item  The \elt $x/s\in\gB_S$ is integral over $\gA_S$ \ssi there exists a~$u\in S$ such that $xu$ is integral over $\gA$.
\item  If $\gC$ is the integral closure of $\gA$ in $\gB$, then $\gC_{S}$ is the integral closure of $\gA_S$ in~$\gB_{S}$.
\item  If $\gA$ is \iclz, then so is $\gA_S$.
\end{enumerate}
\end{fact}
%
\begin{proof}
It suffices to prove item \emph{1.}
 First assume $x/s$ integral over $\gA_S$.
We have for example an \egt in $\gB$
$$\preskip.2em \postskip.2em 
u(x^3+a_2sx^2+a_1s^2x+a_0s^3)=0, 
$$
with $u\in S$ and each $a_i\in\gA$. By multiplying by $u^2$ we obtain

\snic{(ux)^3+a_2us(ux)^2+a_1u^2s^2(ux)+a_0u^3s^3=0}

 in $\gB$. Conversely suppose $xu$ is integral over $\gA$ with $u\in S$.
We have for example an \egt
$$\preskip.4em \postskip.4em 
(ux)^3+a_2(ux)^2+a_1(ux)+a_0=0 
$$
in $\gB$, therefore in $\gB_S$ we have

\snic{x^3+(a_2/u)x^2+(a_1/u^2)x+(a_0/u^3)=0.}

\vspace{-1em}
\end{proof}
%

\begin{plcc}
\label{plcc.entier} \emph{(Integral elements)}\\
Let $S_1$, $\ldots$, $S_n$ be \moco of a \ri $\gA\subseteq\gB$ and $x\in\gB$.
We have the following \eqvcsz.
\begin{enumerate}
\item  The \elt $x$ is integral over $\gA$ \ssi it is integral over each~$\gA_{S_i}$.
\item  Assume $\gA$ is integral, then $\gA$ is \icl \ssi each~$\gA_{S_i}$ is \iclz.
\end{enumerate}
\end{plcc}
\begin{proof}
In item \emph{1} we need to prove that if the condition is locally achieved, then it is globally achieved. Consider then some $x\in\gB$ which satisfies for each $i$ a relation  $(s_ix)^k=a_{i,1}(s_ix)^{k-1}+ a_{i,2} (s_ix)^{k-2}+ \cdots +a_{i,k}$ with $a_{i,h}\in\gA$ and $s_i\in S_i$ (we can assume \spdg that the degrees are the same).
We then use a relation $\sum s_i^ku_i=1$ to obtain an \rdi of $x$ over $\gA$.
\end{proof}

\KROz's \tho easily implies the following lemma.

\begin{lemma}\label{lem0IntClos} 
\emph{(\KROz's \thoz, case of an integral \riz)}\\
Let $\gA$ be \iclz, and $\gK$ be its quotient field.
If we have~$f=gh$ in~$\KT$ with $g$, $h$ \mons and $f\in\AT$, then $g$ and $h$ are also in~$\AT$.
\end{lemma}

\begin{lemma}\label{lemZintClos}
The \ri $\ZZ$ as well as the \ri $\KX$ when $\gK$ is a \cdi are \iclz.
\end{lemma}
\begin{proof}
In fact this holds for every \ri with an integral gcd $\gA$ (see Section~\ref{secGpReticules}).
Let $f(T)=T^n-\sum_{k=0}^{n-1}f_kT^k$ and $a/b$ be a reduced fraction in the quotient field 
 of $\gA$ with $f(a/b)=0$.
By multiplying by $b^n$ we obtain  

\snic{a^n=b\;\sum_{k=0}^{n-1}f_ka^kb^{n-1-k}.}

Since $\pgcd(a,b)=1$,  $\pgcd(a^n,b)=1$. But $b$ divides $a^n$, therefore $b$ is \ivz, and $a/b\in\gA$.
\end{proof}
%

\begin{theorem}\label{thIntClosStab}
If $\gA$ is \iclz, the same goes for~$\AX$.
\end{theorem}
\begin{proof}
Let $\gK=\Frac\gA$. If some \elt $f$ of $\gK(X)$ is integral over $\AX$, it is integral over $\KX$, therefore in $\KX$ because $\KX$ is \iclz.
The result follows by Lemma~\ref{lemPolEnt}; all the \coes of the \polz~$f$ are integral over $\gA$, therefore in $\gA$.
\end{proof}

An interesting corollary of \KROz's \tho is the following \prt (with the same notation as in \Thref{thKro}).

\pagebreak

\begin{proposition}
\label{propArm} Let $f,g\in\AuX$.
Assume that $\gA$ is \iclz, and that $a\in \gA$ divides all the \coes of $h=fg$, then $a$ divides all the $f_{\!\alpha\,} g_\beta$.
In other words
$$ \rc(fg)\equiv 0\;\mod\;a\;\;\iff \;\;
\rc(f)\rc(g)\equiv 0\;\mod\;a.
$$
\end{proposition}
\begin{proof}
Indeed, when considering the \pols $f/a$ and $g$ with \coes in the quotient field of $\gA$,  \KROz's \tho implies that $f_{\!\alpha\,} g_\beta/a$ is integral over $\gA$ because every $h_\gamma/a$ is in $\gA$.
\end{proof}

\subsubsection*{Decomposition of \pols into products of \irds factors}

\begin{lemma}\label{lemKXfactor}
Let $\gK$ be a \cdiz. The \pols of $\KX$ can be decomposed into products of \irds factors \ssi we have an \algo to compute the zeros in $\gK$ of an arbitrary \pol of~$\KX$.
\end{lemma}
\begin{proof}
The second condition is a priori weaker since it amounts to determining the factors of degree $1$ for some \pol of $\KX$. Assume this condition is satisfied. To know whether there exists a \dcn $f=gh$ with~$g$ and~$h$ \mon of fixed degrees $>0$, we apply \KROz's \thoz. We obtain for each \coe of $g$ and $h$ a finite number of possibilities (they are the zeros of \polus that we can explicitly express according to the \coes of $f$).
\end{proof}
%

\begin{proposition}\label{propZXfactor}
In $\ZZ[X]$ and $\QQ[X]$ the \pols can be decomposed into products of \ird factors. 
A nonconstant \pol of $\ZZ[X]$ is \ird in $\ZZ[X]$ \ssi it is primitive and \ird in~$\QQ[X]$.
\end{proposition}
\begin{proof}
For $\QQ[X]$ we apply Lemma~\ref{lemKXfactor}. We must therefore show that we know how to determine the rational zeros of a \polu $f$ with rational \coesz. We can even assume that the \coes of $f$ are integral.
The \elr theory of \dve in $\ZZ$ shows then that if $a/b$ is a zero of~$f$,
$a$ must divide the leading \coe and $b$ the constant \coe of $f$;
there is therefore only a finite number of tests to execute.\\
For $\ZZ[X]$, a primitive \pol $f$ being given, we want to know if there exists a \dcn $f=gh$ with  $g$ and $h$ of fixed degrees $>0$. We can assume $f(0)\neq0$.
We apply \KROz's \thoz. A product $g_0h_j$ for instance must be a zero in $\ZZ$ of a \polu $q_{0,j}$ of $\ZZ[T]$ that we can compute.
In particular, $g_0h_j$ must divide $q_{0,j}(0)$, which only leaves a finite number of possibilities for $h_j$.\\
Finally, for the last item, if some primitive \polz~$f$ in $\ZZ[X]$ can be decomposed in the form $f=gh$ in $\QQ[X]$ we can assume that $g$ is primitive in~$\ZZ[X]$. Let $a$ be a \coe of $h$, then every $ag_j$ is in $\ZZ$ (\KROz's \thoz), and the Bézout relation $\sum_jg_ju_j=1$ shows \hbox{that $a\in\ZZ$}.
\end{proof}

\vspace{-.8em}
\pagebreak

\subsec{Number fields}

We call  a \cdi  $\gK$ a \emph{number field} if it is \stf over $\QQ$.

\vspace{-1pt}
\subsubsection*{Galois closure}

\begin{theorem}\label{thClotAlgQ} \emph{(Splitting field, primitive \elt \thoz)}
\begin{enumerate}
\item
If $f$ is a \spl \polu of $\QQ[X]$
there exists a number field $\gL$ over which we can write $f(X)=\prod_i(X-x_i)$.
In addition, with some $\alpha\in\gL$ we have
$$\preskip.0em \postskip.4em
\;\;\;\gL=\QQ[\xn]=\QQ[\alpha]\simeq\aqo{\QQ[T]}{Q},
$$
where $Q(\alpha)=0$ and the \polu $Q$ is \ird in $\QQ[T]$ and is completely decomposable in $\gL[T]$.
\\
In particular, the extension $\gL/\QQ$ is Galoisian and \Thref{thGaloiselr} applies.
\item Every number field $\gK$ is contained in a Galois extension of the above type. In addition, there exists some $x\in\gK$ such that $\gK=\QQ[x]$.
\end{enumerate}

\end{theorem}
\begin{proof}
\emph{1.} This results from \Thref{thResolUniv} and from Proposition~\ref{propZXfactor}.\\
\emph{2.} A number field is generated by a finite number of \elts that are \agqs over $\QQ$.
Each of these \elts admits a \polmin that is \ird over $\QQ$ and therefore \spl (Fact~\ref{factPolSepFC}).
By taking the lcm~$f$ of these \pols we obtain a \spl \polz. 
By applying item~\emph{1} to~$f$
and by using \Thref{propUnicCDR}, we see that $\gK$ is \isoc to a subfield 
of~$\gL$. Finally, as the Galois correspondence is bijective and as the Galois group $\Gal(\gL/\QQ)$ is finite, the field $\gK$ only contains an explicit finite number of subfields $\gK_i$ \stfs over $\QQ$. If   we choose $x\in\gK$ outside of the union of these subfields
(which are strict \Qsevsz), we \ncrt have $\QQ[x]=\gK$; it is a subfield of $\gK$ \stf over $\QQ$ and distinct from all the $\gK_i$'s.
\end{proof}

\subsubsection*{Cotransposed \eltz} \rdb

If $\gB$ is a free \Alg of finite rank, we can identify $\gB$ with a commutative sub\alg of $\End_\gA(B)$, where $B$ designates the \Amo $\gB$ deprived of its multiplicative structure, by means of the \homo $x\mapsto\mu_{\gB,x}$, \hbox{where $\mu_{\gB,x}=\mu_x$} is the multiplication by $x$ in $\gB$.
Then, since $\wi \mu_x=G(\mu_x)$ for some \polz~$G$ of $\AT$ (Lemma~\ref{lemPrincipeIdentitesAlgebriques} item~\emph{6}), we can define $\wi x$ by the \egt $\wi x=G(x)$, or equivalently $\wi {\mu_x}=\mu_{\wi x}$.
If more precision is \ncrz, we will use the notation $\Adj\iBA (x)$.
This \elt $\wi x$ is called \emph{the cotransposed \elt of $x$}.
We then have the important \egt%
\index{cotransposed!element (in a free algebra)}
\begin{equation}
\label{eqelt0cotransp}
x\ \wi x=x\ \Adj\iBA(x)=\rN\iBA(x).
\end{equation}

\rdb
\rem
\label{factNormeRationnelle}
Let us also note that the applications \gui{norm of} and \gui{cotransposed \elt of} enjoy some \prts of \gui{$\gA$-rationality,} which directly result from their \dfnsz:
if $P\in\gB[\Xk]$, then by taking the $x_i$'s in $\gA$, $\rN\iBA \big(P(\xk)\big)$ and $\Adj\iBA \big(P(\xk)\big)$ are given by \pols of $\gA[\Xk]$.
\\
In fact $\BuX$ is free over $\AuX$ with the same basis as that of $\gB$ over~$\gA$ and $\rN\iBA \big(P(\ux)\big)$ is given by the evaluation at $\ux$ of $\rN_{\BuX/\!\AuX}\big(P(\uX)\big)$ (likewise for the cotransposed \eltz). We will use by abuse of notation $\rN\iBA \big(P(\uX)\big)$.
\\
   Furthermore, if $\dex{\gB:\gA}=n$ and if $P$ is \hmg of degree $d$, then $\rN\iBA \big(P(\uX)\big)$ is \hmg of degree $nd$ and $\Adj\iBA \!\big(P(\uX)\big)$ is \hmg of degree $(n-1)\,d$.
\eoe

\subsec{\Ri of integers of a number field}

If $\gK$ is a number field its \ixc{ring of integers}{of a number field} is the integral closure of  $\ZZ$ in $\gK$.

\begin{propdef}\label{propAECDN} \emph{(Discriminant of a number field)}
 \\
Let $\gK$ be a number field and $\gZ$ its \ri of integers.
\begin{enumerate}
\item  An \elt $y$ of $\gK$ is in $\gZ$ \ssi $\Mip_{\QQ,y}(X)\in\ZZX$. 
\item  We have $\gK=(\NN\etl)^{-1}\gZ$.  
\item  
Assume that $\gK=\QQ[x]$ with $x\in\gZ$. Let $f(X)=\Mip_{\QQ,x}(X)$ be in $\ZZ[X]$ and $\Delta^2$ be the greatest square factor of $\disc_X f$. 
\\
Then, $\ZZ[x]\subseteq\gZ\subseteq  {1 \over \Delta} \ZZ[x]$.

\item The \ri $\gZ$ is a free \ZZmo of rank $\dex{\gK:\gQ}$.
\item  The integer $\Disc_{\gZ/\ZZ}$ is well-defined. We call it the \emph{discriminant of the number field $\gK$}.%
\index{discriminant!of a number field} 
\end{enumerate}
 
\end{propdef}
%
\begin{proof}
\emph{1.} Results from Lemma~\ref{lem0IntClos} (\KROz's \thoz). 
 
\emph{2.} Let $y\in\gK$ and $g(X)\in\ZZX$ be a nonzero \pol that annihilates $y$. 
If $a$ is the leading \coe of $g$, $ay$ is integral over $\ZZ$.
 
\emph{3.}  
Let $\gA=\ZZ[x]$ and $n=\dex{\gK:\gQ}$.  
Let $z\in\gZ$, which we as  $h(x)/\delta$ with~$\delta\in\NN\etl$, $\gen{\delta}+\rc(h)=\gen{1}$ and $\deg h < n$. We have $\gA+\ZZ z\subseteq{1 \over \delta} \gA$ and it thus suffices to prove that $\delta^2$ divides~$\disc_X(f)$.
%
%
 The \ri $\gA$ is a free \ZZmo of rank $n$, with the basis $\cB_0=(1,x,\ldots,x^{n-1})$. Proposition~\ref{propdiscTra} gives
$$\preskip-.2em \postskip.4em 
\Disc_{\gA/\ZZ}=\disc_{\gA/\ZZ}(\cB_0)=\disc_{\gK/\QQ}(\cB_0)=\disc_X f . 
$$
The \ZZmo $M =\gA+\ZZ z$ is \egmt free, of rank $n$ with a basis $\cB_1$, and we obtain the \egts
 
\snic{\disc_X f=\disc_{\gK/\QQ}(\cB_0)= \disc_{\gK/\QQ}(\cB_1)\times d^2,}

 where $d$ is the \deter of the matrix of $\cB_0$ over $\cB_1$ (Proposition~\ref{defiDiscTra}~\emph{2}).
\\
Finally, $d=\pm\delta$ by the following Lemma~\ref{lemSousLibre}, as required.
 
\emph{4.} \Spdg we use the setup of item \emph{3.}
There is only a finite number of \tf \ZZmos between $\ZZ[x]$ and ${1\over\Delta}\ZZ[x]$, 
and for each of them we can test whether it is contained in $\gZ$.
The largest is \ncrt equal to $\gZ$.
\end{proof}
\rems ~\\
1) As a corollary, we see that in the context of item~\emph{3}, if $\disc_X(f)$ has no square factors, then $\gZ=\ZZ[x]$.

2) The \dem of item~\emph {4} does not provide the {practical} means to compute a $\ZZ$-basis of $\gZ$. For some more precise information see \Pbmz~\ref{exoLemmeFourchette}. Actually we do not know of a \gnl \emph{\poll time} \algo to compute a $\ZZ$-basis of $\gZ$.
\eoe

\medskip 
One says that an ideal $\fa$ of a \ri $\gA$ is \emph{principal} when
it is generated by a single element.%
\index{principal!ideal}\index{ideal!principal ---}

\begin{lemma}\label{lemSousLibre}
Let $N \subseteq M$ be two free \Amos of the same {rank $n$}
{with $M = N + \gA z$}. Assume that for some \ndz \elt $\delta \in \gA$, we 
\hbox{have~$\delta z \in N$} and $\delta z = a_1 e_1 + \cdots + a_n e_n$,
where $(e_1, \ldots, e_n)$ is a basis of $N$. Then, the \deter $d$ of a matrix of a basis of $N$ over a basis $M$ satisfies
\begin{equation}\preskip.4em \postskip.4em
\label{eqlemSousLibre}
d \,\gen {\delta, \an} = \gen {\delta}
\end{equation}
In particular, $\gen {\delta, \an}$ is a \idpz, and if $\delta$, $a_1$, \dots, $a_n$ are \comz, then $\gen {d} = \gen {\delta}$. Moreover, $M/N \simeq \aqo {\gA}{d}$.  
\end{lemma}
%
\begin{proof} \Egrf{eqlemSousLibre}  is left to the reader (see 
Exercise~\ref{exolemSousLibre}).
\\
It remains to prove that $M/N \simeq \aqo {\gA}{d}$.
By letting $\ov z$ be the class of~$z$ in~$M/N$, since $M/N\simeq\gA \ov z$, we must prove that $\Ann_\gA(\ov z) = \gen {d}$,
\cad \hbox{that $bz \in N \Leftrightarrow b \in \gen {d}$}. It is clear that $dz \in N$.\\
If $bz \in N$, then $b\delta z \in \delta N$, therefore by writing
$\delta z = a_1e_1 + \cdots + a_ne_n$, we \hbox{get $ba_i \in \gen{\delta}$}, then
$b\gen {\delta,\an} \subseteq \gen {\delta}$. By multiplying by $d$ and by simplifying by~$\delta$, we obtain $b \in \gen {d}$.
\end{proof}
%

\subsubsection*{The multiplicative theory of the \ids of a number field}

\begin{definition}
\label{defiiv}
An \id $\fa$ of a \ri $\gA$ is said to be
 \emph{\ivz} if there exist an \id $\fb$ and a \ndz \elt $a$ such that $\fa\,\fb=\gen{a}$.%
\index{invertible!ideal}%
\index{ideal!invertible ---}
\end{definition}

\begin{fact}
\label{factdefiiv}
Let $\fa$ be an \iv \id of a \ri $\gA$.
\begin{enumerate}
\item  The \id $\fa$ is \tfz.
\item If $\fa$ is generated by $k$ \elts and if $\fa\,\fb=\gen{a}$ with $a$ \ndzz, then $\fb$ is generated by $k$ \eltsz. Furthermore $\fb=(\gen{a}:\fa)$.
\item  We have the rule $\fa\,\fc \subseteq \fa\,\fd \;\Rightarrow\;\fc \subseteq \fd $
for all \idsz~$\fc $ and~$\fd $.
\item If $\fc \subseteq \fa$ there exists a unique $\fd $ such that $\fd \,\fa=\fc $, namely $\fd=(\fc:\fa)$, 
and if $\fc $ is \tfz, so is $\fd$. 
\end{enumerate}
\end{fact}
\begin{proof}
 \emph{3.} If $\fa\,\fc \subseteq \fa\,\fd $ by multiplying by $\fb$ we obtain $a\,\fc \subseteq a\,\fd $,
and since $a$ is \ndzz, this implies $\fc \subseteq \fd$.
\\
\emph{1.}
If $\fa\,\fb=\gen{a}$, we find two \itfs $\fa_1\subseteq\fa$ and $\fb_1\subseteq\fb$
such \hbox{that $a\in\fa_1\,\fb_1$} and thus $\fa\,\fb=\gen{a}\subseteq\fa_1\,\fb_1\subseteq\fa\,\fb_1\subseteq \fa\,\fb$.
From the above, we deduce the \egts $\fa_1\,\fb_1=\fa\,\fb_1=\fa\,\fb$.
Whence $\fb=\fb_1$ by item~\emph{3}. Similarly, $\fa=\fa_1$.
\\
\emph{2.} If $\fa=\gen{a_1,\dots,a_k}$, we find $b_1$, \dots, $b_k\in\fb$ such that $\som_i a_ib_i=a$. \\
By reasoning as in item \emph{1} with $\fa_1=\fa$ and $\fb_1=\gen{b_1,\dots,b_k}$ we obtain the \egt $\fb=\gen{b_1,\dots,b_k}$. 
Since $\fa\,\fb=\gen{a}$, we have $\fb\subseteq (\gen{a}:\fa)$. Conversely, if $x\fa\subseteq \gen{a}$, then 
%
${x\gen{a}=x\,\fa\,\fb\subseteq a\,\fb},$
%
thus $ax=ab$ for some $b\in\fb$ and $x\in\fb$ because $a$ is \ndzz. 
\\
\emph{4.} From $\fa \,\fb=\gen{a}$ we deduce $\fc \,\fb\subseteq \gen{a}$.
All the \elts of $\fc \,\fb$ being multiples of $a$, by dividing them by $a$ we get an \id $\fd $, that we denote by $\fraC 1 a \,\fc \,\fb$, and with which we obtain the \egt $\fa\,\fd=\fraC 1 a \,\fc \,\fb\,\fa =\fraC 1 a \,\fc \,\gen{a}=\fc $ because $a$ is \ndzz.
\\
If $\fc $ is \tfz, $\fd $ is generated by the \elts obtained by dividing each \gtr of $\fc \,\fb$ by $a$.
\\
The uniqueness of $\fd$ results from item \emph{3.}
\\
All is left to prove is that $\fd=(\fc:\fa)$. The inclusion $\fd\subseteq (\fc:\fa)$ is \imdez.   
Conversely, if $x\fa\subseteq \fc$, then $x\gen{a}\subseteq \fc\,\fb$, therefore $x\in\fraC 1 a \,\fc \,\fb=\fd$.
\end{proof}

The following \tho is the key \tho in the multiplicative theory of the \ids of number fields. We provide two \demsz.
Beforehand we invite the readers to acquaint themselves with \Pbmz~\ref{exoPetitKummer} 
which gives Kummer's little \thoz, which solves with minimal costs the question for \gui{almost all} the \itfs of the number fields.%
\index{Kummer!little theorem} 
\Pbmz~\ref{exoCyclotomicRing} is \egmt instructive as it gives a direct \dem of the invertibility of all the nonzero \itfs as well as of their unique \dcn into a product of \gui{prime factors} for the \ri $\ZZ[\root n \of 1\,]$. 

\begin{theorem} \emph{(Invertibility of the \ids of a number field)}\label{th1IdZalpha}\\
Every nonzero \itf of the \ri of integers $\gZ$ of a number field~$\gK$ is \ivz.
\end{theorem}
\begin{proof}
\emph{First \dem  (\`a la Kronecker.\footnote{Actually Kronecker does not use the \emph{cotransposed \eltz} of $\alpha+\beta X+\gamma X^2$ (as stated in the \dfn we have given), but the product of all the conjugates of $\alpha X+\beta Y+\gamma Z$ in a Galois extension.
This introduces a slight variation in the \demz.})}
\\
Take for example $\fa=\gen{\alpha,\beta,\gamma}$.
Let $\gA=\QQ[X]$ and $\gB=\KX$.
The \alg $\gB$ is free over $\gA$ with the same basis as that of $\gK$ over $\QQ$.
Consider the \polz~$g=\alpha+\beta X+\gamma X^2$ which satisfies $\rc_\gZ(g)=\fa$.
Since $\alpha$, $\beta$, $\gamma$ are integral over $\ZZ$, $g$ is integral over $\ZZ[X]$.
Let $h(X)=\Adj\iBA (g)$ be the cotransposed \elt of $g$.
We know that $h$ is expressed as a \pol in~$g$ and in the \coes of the \polcar of $g$. By applying Corollary~\ref{corlemPolcarEntier} we deduce that $h$ has \coes in $\gZ$.
 Let $\fb$ be the \itf of $\gZ$ generated by the \coes of $h$.
 \\
We have $gh=\rN\iBA (g)\in\gZ[X]\cap\QQ[X]=\ZZ[X]$.
Let $d$ be the gcd of the \coes of $gh$. Proposition~\ref{propArm} tells us that an arbitrary \elt of~$\gZ$ divides $d$ \ssi it divides all the \elts \hbox{of $\fa\,\fb$}. 
\\
In particular, $d\gZ\supseteq\fa\,\fb$. Given the Bézout relation that expresses $d$ according to the \coes of $gh$ we also have $d\in\fa\,\fb$. Therefore $d\gZ=\fa\,\fb$.

\emph{Second \dem (\`a la Dedekind.)}\index{Dedekind!inversion of an ideal \`a la ---}
\\
 First of all we notice that it suffices to know how to invert the \ids with two \gtrs by virtue of the following remark. For three arbitrary \ids~$\fa$,~$\fb$,~$\fc$ in a \ri we always have the \egt
$$\preskip.4em \postskip.4em 
(\fa+\fb)(\fb+\fc)(\fc+\fa)=(\fa+\fb+\fc)(\fa\fb+\fb\fc+\fa\fc), 
$$
therefore, if we know how to invert the \ids with $k$ \gtrs ($k\geq2$), we \egmt know how to invert the \ids with $k+1$ \gtrsz.
\\
We thus consider an \id $\gen{\alpha,\beta}$ with  $\alpha\neq0$.
As $\alpha$ is integral over $\ZZ$, we can find $\ov{\alpha} \in \gZ$ such that $\ov{\alpha}\alpha \in \ZZ \setminus \so{0}$.
Thus, even if it means replacing $(\alpha,\beta )$ with $(\ov{\alpha}\alpha,\ov{\alpha}\beta )$,
we restrict ourselves to the study of an \id $\gen{a,\beta}$ with $(a,\beta)\in\ZZ\times\gZ$.
\\
 Let $f \in \ZZ[X]$ be a \polu which is annihilated in $\beta$. We write

\snic{f(X)=(X-\beta)h(X)$,
where $h \in \gZ[X]\,.}

We thus have $f(a X)=(a X-\beta )h(a X)$, which we rewrite as $f_1=g_1h_1$.
Let then $d$ be the gcd of the \coes of $f_1$ in~$\ZZ$. With $\fb=\rc_\gZ(h_1)$ and $\fa=\rc_\gZ(g_1)=\gen{a,\beta}$, we clearly have $d \in \fa \fb$.
Moreover, Proposition~\ref{propArm} tells us that an arbitrary \elt of~$\gZ$ divides all the \elts of $\rc_\gZ(f_1)=\gen{d}$ \ssi it divides all the \elts of the \id $\fa\,\fb$. In particular,~$d$
divides all the \elts of $\fa\,\fb$.
Thus $\fa \fb=\gen{d}$.
\end{proof}

The following \tho shows that the \itfs of a number field with regard to the \elr operations (sum, intersection, product, exact division) behave essentially equivalently to the principal \ids of $\ZZ$. The latter translate the theory of \dve for the natural numbers very precisely.
\\
Recall that in the bijection $n \mapsto n\ZZ$ ($n \in \NN$, $n\ZZ$ a \itf of~$\ZZ$): the product corresponds to the product, \dve corresponds to inclusion; the gcd to the sum; the lcm to the intersection; and the exact division to the conductor.

\begin{theorem} \emph{(The \itfs of a number field)}\label{th2IdZalpha}\\
Let $\gK$ be a number field and $\gZ$ its \ri of integers.
\begin{enumerate}
\item \label{i2th2IdZalpha}
If $\fb$ and $\fc$ are two arbitrary \idsz, and if $\fa$ is some nonzero \itf of $\gZ$, we have the implication
$$\preskip.2em \postskip.4em
\fa\, \fb\subseteq \fa\, \fc\ \ \Rightarrow\ \ \fb\subseteq  \fc\,.
$$
\item \label{i3th2IdZalpha}
 If $\fb\subseteq\fc$ are two \itfsz, there exists some \itf $\fa$ such that $\fa\, \fc=\fb$.
\item \label{i4th2IdZalpha}
 The set of \itfs of $\gZ$ is stable by finite intersections and we have the following \egts (where $\fa$, $\fb$, $\fc$ designates \itfs of~$\gZ$):
\[\preskip.2em \postskip.4em
\arraycolsep2pt\def\mathit#1{\llap{{#1}\hskip1.5cm}}
\begin{array}{lrclc}
\mathit{a.} &(\fa\cap\fb)(\fa+\fb)  &  = & \fa\fb \,,&\qquad\qquad \qquad \\[1mm]
\mathit{b.} &\fa\cap(\fb+\fc)  &  = & (\fa\cap\fb)+(\fa\cap\fc)  \,, \\[1mm]
\mathit{c.} &\fa+(\fb\cap\fc)  &  = & (\fa+\fb)\cap(\fa+\fc)   \,,\\[1mm]
\mathit{d.} &\fa  (\fb\cap\fc)  &  = & (\fa \fb)\cap(\fa \fc)  \,,\\[1mm]
\mathit{e.} &(\fa+\fb)^n  &  = & \fa^n+\fb^n\quad(n\in\NN) \,.
\end{array}
\]
\item \label{i5th2IdZalpha} If $\fa$ is some nonzero \itf of $\gZ$ the \ri $\gZ\sur\fa$ is finite. 
\\
In particular,
we have tests to decide:
\begin{itemize}
\item if some $x\in\gZ$ is in $\fa$,
\item if some $x\in\gZ$ is \iv modulo $\fa$, 
\item if $\fa$ is contained in another \itf $\fb$, 
\item  if $\gZ\sur\fa$ is a \cdi (we then say that $\fa$ is a detachable \idemaz).
\end{itemize}

\item \label{i6th2IdZalpha} Every distinct \itf of $\gen{0}$ and $\gen{1}$ is equal to a product of detachable \iv \idemasz, and this \dcn is unique up to order of the factors.
\end{enumerate}
\end{theorem}
\begin{proof}
\emph{\ref{i2th2IdZalpha}} and \emph{\ref{i3th2IdZalpha}.}
By Fact~\ref{factdefiiv}.
 
\emph{\ref{i4th2IdZalpha}.}
If one of the \itfs is zero everything is clear. We assume they are nonzero in the remainder of the \demz.
 
\emph{\ref{i4th2IdZalpha}a.}
Let $\fc$ such that $\fc(\fa+\fb)=\fa\fb$.
Since $(\fa\cap\fb)(\fa+\fb)\subseteq\fa\fb$, we obtain the inclusion $\fa\cap\fb\subseteq\fc$ (simplification by $\fa+\fb$). Conversely, $\fc\fa \subseteq \fa\fb$, \hbox{thus $\fc \subseteq \fb$}
   (simplification by $\fa$). Similarly $\fc \subseteq \fa$.
 
\emph{\ref{i4th2IdZalpha}c.}
We multiply both sides by $\fa+\fb+\fc=(\fa+\fb)+(\fa+\fc)$.\\
The right-hand side gives $(\fa+\fb)(\fa+\fc)$.\\
The left-hand side gives $\fa(\fa+\fb+\fc)+ \fa(\fb\cap\fc)+(\fb+\fc)(\fb\cap\fc)$. \\
Both cases result in $\fa(\fa+\fb+\fc)+\fb\fc$.
 
\emph{\ref{i4th2IdZalpha}b.}
For the inclusion, the \itfs form a lattice (the supremum is the sum and the infimum is the intersection). We come to see that one of the laws is distributive with respect to the other. Classically, in a lattice this implies the other \dit (see \paref{DistriTrdi}).
 
\emph{\ref{i4th2IdZalpha}d.}
The map $\fx\mapsto\fa\,\fx$ (of the set of \itfs to the set of \itfs which are multiples of $\fa$) is an \iso of the order structure by item~\emph{\ref{i2th2IdZalpha}.} This implies that the map transforms $\fb\cap\fc$ into the infimum of $\fa\fb$ and $\fa\fc$ inside the set of \itfs that are multiples of $\fa$. 
 It thus suffices to establish that $\fa\fb\cap\fa\fc$ is a multiple of $\fa$. This results from item~\emph{\ref{i3th2IdZalpha}.}
 
\emph{\ref{i4th2IdZalpha}e.}
For example with $n=3$, $(\fa+\fb)^3=\fa^3+\fa^2\fb+\fa\fb^2+\fb^3$.
\\
By multiplying $(\fa+\fb)^3$ and $\fa^3+\fb^3$ by $(\fa+\fb)^2$
we find in both cases 
$$\preskip.4em \postskip.4em 
\fa^5+\fa^4\fb+\cdots+\fa\fb^4+\fb^5. 
$$

\emph{\ref{i5th2IdZalpha}.} 
View $\gZ$ as a free \ZZmo of rank $n=\dex{\gK:\QQ}$. It is obvious that a \itf $\fa$ containing the integer $m\neq0$ can be explicitly expressed as a \tf \ZZsmo of $\ZZ^n$ containing $m\ZZ^n$. 
 
\emph{\ref{i6th2IdZalpha}.}
Let $\fa$ be a \itf $\neq\gen{0},\gen{1}$. 
The \tf \idemas of $\gZ$ containing $\fa$ are obtained by determining the \tf \idemas of $\gZ\sur\fa$ (which is possible because the \ri $\gZ\sur\fa$ is finite). 
If $\fp$ is a \tf \idema containing $\fa$, we can write $\fa=\fb\,\fp$.
Furthermore, 
we have the \egt $\idg{\gZ:\fa}=\idg{\gZ:\fb} \,\idg{\fb:\fa}$. We then obtain the \dcn into products of \tf \idemas by \recu on $\idg{\gZ:\fa}$. The uniqueness results from the fact that if a \tf \idema $\fp$ contains a product of \tf \idemasz, it is \ncrt equal to one of them, otherwise it would be comaximal with the product. 
\end{proof}

We end this section with a few generalities concerning \emph{the \ids that avoid the conductor}. The situation in number theory is the following.
We have a number field $\gK=\QQ[\alpha]$ with  $\alpha$ integral over $\ZZ$.
We denote by $\gZ$ the \ri of integers of $\gK$, \cad the integral closure of $\ZZ$ in $\gK$. Even though it is possible in principle, it is not easy to obtain a basis of $\gZ$ as a \ZZmoz, nor is it easy to study the structure of the \mo of the \itfs of~$\gZ$.

Assume that we have a \ri $\gZ'$ which constitutes an approximation of $\gZ$ in the sense that $\ZZ[\alpha]\subseteq\gZ'\subseteq\gZ$. For example \hbox{let $\gZ'=\ZZ[\alpha]$} initially. We are interested in the multiplicative structure of the group of \ifrs of~$\gZ$,%
\footnote{A \ifr of $\gZ$ is a $\gZ$-submodule of $\gK$ equal to $\fraC{1}m\,\fa$ for \hbox{some $m\in\ZZ\sta$} and a \itf $\fa$ of $\gZ$, 
cf.~\paref{NOTAIfr}.} and we want to rely on that of $\gZ'$ to study it in detail.

The following \tho states that \gui{this works very well for most \idsz, \cad for every one that avoids the conductor of $\gZ$ into $\gZ'$.}  

\begin{definition}\label{defiConducteur}~
\\
Let $\gA$, $\gB$ be two \ris such that $\gA \subseteq \gB$, and let $\fa$ and $\fb$ be respective \ids of $\gA$ and $\gB$. 
\begin{enumerate}
\item The \emph{conductor of $\gB$ into
$\gA$} is $(\gA:\gB) = \sotq {x \in \gB} {x\gB \subseteq \gA}$.%
\index{conductor!of a ring into a subring}
\item The \emph{extension of $\fa$} is the \id $\fa\gB$ of $\gB$.%
\index{extension!of an ideal (in an overring)}
\item The \emph{contraction of $\fb$} is the \id $\gA \cap \fb$ of $\gA$.%
\index{contraction!of an ideal (in a subring)}
\end{enumerate} 
\end{definition}

\begin{theorem}\label{propEvitementConducteur}\index{Dedekind!ideals that avoid the conductor} 
\emph{(Dedekind's \thoz, \ids that avoid the conductor)}\\
Let $\gA$, $\gB$ be two \ris such that $\gA \subseteq \gB$ and $\ff$ be the conductor of $\gB$ into~$\gA$. 
\begin{enumerate}
\item The \id $\ff$ is the annihilator of the \Amo $\gB\sur\gA$. It is simultaneously an \id of $\gA$ and an \id of $\gB$, and it is the greatest \id with this \prtz.
\end{enumerate}
 We denote by $\cA$ (resp.\,$\cB$) the class of \ids of $\gA$ (resp.\,of $\gB$) \com to~$\ff$.
\begin{enumerate}\setcounter{enumi}{1}
\item For $\fa \in \cA$, we have $\gA\sur\fa \simeq \gB\sur{\fa\gB}$ and for $\fb \in \cB$, we have $\gB\sur\fb \simeq \gA\sur{\gA\cap\fb}$.

\item $\cA$ is stable under multiplication, sum, intersection and satisfies 

\snic{\fa \in \cA, \,\fa'
\supseteq \fa \;\;\;\Longrightarrow\;\;\; \fa' \in \cA.}

In particular, $\fa_1\fa_2 \in \cA$ \ssi $\fa_1$ and $\fa_2 \in \cA$. 
The same \prts are valid for~$\cB$.
\item The extension and the contraction, restricted respectively to $\cA$ and $\cB$, are 
inverses of each other. They preserve multiplication, inclusion, intersection and the \tf \crcz.
\item Assume that $\gB$ is integral. Then, an \id $\fa\in\cA$ is \iv in $\gA$ \ssi $\fa\gB$ is \iv in $\gB$. Similarly, an \id $\fb\in\cB$ is \iv in $\gB$ \ssi $\gA \cap \fb$ is \iv in~$\gA$.

\end{enumerate}
\end{theorem}
 
%
\begin{proof}
We only prove a few \prtsz. 
Notice that we always have the inclusions $\fa \subseteq \gA\cap \fa\gB$ and $(\gA\cap\fb)\gB \subseteq \fb$.

Let $\fa \in \cA$, so $1 = a + f$ with $a \in \fa$ and $f \in \ff$;
a fortiori, $1 \in \fa\gB + \ff$. Let us prove that $\gA\cap\fa\gB = \fa$.
We take $x \in \gA\cap\fa\gB$ and we write 

\snic{x = xf + xa \in \fa\gB\ff + \fa \subseteq
\fa\gA + \fa = \fa\,.}

Hence the result. We also see that $\gB = \gA + \fa\gB$,
so the composed morphism $\gA \to \gB\sur{\fa\gB}$ is surjective with kernel $\fa$, which gives an \iso $\gA\sur\fa \simeq \gB\sur{\fa\gB}$.

Let $\fb \in \cB$, so $1 = b + f$ with $b \in \fb$, $f \in \ff$.
Since $\ff \subseteq \gA$, we have $b \in \gA \cap \fb$ therefore $1 \in \gA\cap \fb + \ff$. Let us prove that $(\gA \cap\fb)\gB = \fb$.
If  $x \in \fb$, then

\snac {
x = (b+f)x = bx + xf \in (\gA \cap\fb)\gB + \fb\ff \subseteq 
(\gA \cap\fb)\gB + \gA\cap\fb \subseteq (\gA \cap\fb)\gB.
}

Thus $\fb \subseteq (\gA \cap\fb)\gB$ then $\fb = (\gA \cap\fb)\gB$.
In addition, since $\gB = \fb + \ff \subseteq \fb + \gA$, the composed morphism $\gA \to \gB\sur\fb$  is surjective, with kernel $\gA\cap \fb$, which gives an \iso $\gA\sur{\gA\cap
\fb} \simeq \gB\sur\fb$.

The extension is multiplicative, so the contraction (restricted to $\cB$) which is its inverse, is \egmt multiplicative. The contraction is compatible with the intersection, so the extension (restricted to $\cA$) which is its inverse, is \egmt compatible with the intersection.

Let $\fb=\gen{b_1, \ldots, b_n}_\gB \in \cB$. Let us prove that $\gA\cap\fb$ is \tfz.  \\
We write $1 = a + f^2$ with  $a \in \fb$, $f \in \ff$. Since $f \in \gA$, we have $a \in \gA\cap\fb$. We prove that $(a, fb_1, \ldots, fb_n)$ is a \sgr of $\gA\cap\fb$.  \\
Let $x \in \gA\cap\fb$ which we write as $x = \sum_i y_ib_i$ with  $y_i \in \gB$, then
$$\preskip.4em \postskip.4em\ndsp 
x = \sum_i (y_i a + y_i f^2) b_i =  
xa + \sum_i (y_if)fb_i\in \gen {a,fb_1, \ldots, fb_n}_\gA. 
$$
For an \id $\fb\in\cB$ (not \ncrt \tfz), we have in fact proved the following result: if $1 = a + f^2$ with  $a \in \fb$ and $f \in \ff$, then $\gA\cap\fb = \gA a + f(f\fb)$ (and $f\fb$ is an \id of $\gA$).

Let $\fb \in \cB$ be an \iv \idz, let us prove that $\fa = \gA\cap\fb$ is an \iv \idz. We write $1 = a + f$ with  $a \in \fb$ and $f \in \ff$, such that $a \in \fa$. \\
If $a = 0$, then $1 = f \in \ff$, so $\gA = \gB$ and there is nothing left to prove. Otherwise, $a$ is \ndz and there exists an \id $\fb'$ of $\gB$ such that $\fb\fb' = a\gB$.  \\
Since the \ids $a\gB$, $\fb$ and $\fb'$ are \com to $\ff$, we can apply the multiplicative \crc of the contraction to the \egt $\fb\fb' = a\gB$ to obtain the \egt $\fa\fa' = a\gA$ with $\fa' = \gA\cap\fb'$.
\end{proof}

\entrenous{
\rem
1) Donner des exemples for lesquels $\gA \cap \fa\gB \ne \fa$ et
$(\gA \cap \fb)\gB \ne \fb$.

2) Il faudrait aussi peut-\^etre faire le lien with  certains exos.
En fait le remarquable \tho pr\'ec\'edent is insuffisamment exploit\'e.
\eoe

3) tr\`es improbable
\begin{proposition}\label{prop0ITFSCDN}
\emph{(\tho un and demi)}
\hum{semble difficile \`a faire ici, voir le corollaire \ref{corpropZerdimLib}}
\end{proposition}

}

\section{Hilbert's \nstz}\label{secChap3Nst}\ihi   

In this section we illustrate the importance of the resultant by showing how Hilbert's \nst can be deducted from it. We will use a \gnn of the basic \eli lemma \ref{LemElimAffBasic}.

\subsec{The \agq closure of $\QQ$ and of finite fields}

Let $\gK\subseteq\gL$ be \cdisz. We say that \emph{$\gL$ is an \agq closure of~$\gK$} if~$\gL$ is  \agq over $\gK$ and \agqt closed.%
\index{closure!algebraic ---}

The reader will concede that $\QQ$ and the fields $\FFp$ possess an \agq closure.
This will be discussed in further detail in Section~\ref{secEtaleSurCD},
especially with \Thref{thClsep}.

\subsec{The classical \nst (\agqt closed case)}

The \nst is a \tho which concerns the \syss of \pol \eqns over a \cdiz.
Very informally, its meaning can be described as follows: a \gmq statement \ncrt possesses an \agq certificate. Or even: a \dem in commutative \alg can (almost) always be summarized by simple \idas if it is sufficiently \gnlez.

If we have \cdis $\gK\subseteq\gL$, and if $(\uf)=(\lfs)$ is a \sys of \pols in $\KXn=\KuX$,
we say that $(\xin)=(\uxi)$ is a \emph{zero of $(\uf)$ in~$\gL^n$}, or a \emph{zero of $(\uf)$ with \coos in $\gL$}, if the \eqns $f_i(\uxi)=0$ are satisfied.
Let $\ff=\gen{\lfs}_\KuX$. Then, all the \pols $g\in\ff$ are annihilated in such a $(\uxi).$
We therefore equally refer to~$(\uxi)$ as a \emph{zero of the \id $\ff$ in $\gL^n$} or as \emph{having \coos in $\gL$}.

 We begin with an almost obvious fact. 

\pagebreak

\begin{fact}\label{factGCDDeg} Let $\gk$ be a commutative \ri and $h\in\kX$ a \polu of degree $\geq1$. 
 \begin{itemize}
\item If some multiple of $h$ is in $\gk$, this multiple is null.
\item Let $f$ and $g\in\kX$ of respective formal degrees $p$ and $q$. If $h$ divides~$f$ and~$g$, then~$\Res_X(f,p,g,q)=0$.
\end{itemize}
\end{fact}

We now present a \gnn of the basic \eli lemma \ref{LemElimAffBasic}.
 
\begin{lemma}\label{lemElimPlusieurs} \emph{(\Eli of a variable between several \polsz)}
\label{lemElimParametre}%
\index{elimination!of a variable}
\\
Let $f$,  $g_1$, $\ldots$, $g_r$  $\in\kX$ ($r\geq1$), with  $f$  \mon of degree $d$. \\
Let $\ff=\gen{f,g_1,\ldots,g_r}$ and $\fa=\ff\cap\gk$ (this is the \eli \id of the variable $X$ in $\ff$).
Also let
\[\preskip.2em \postskip.4em 
\begin{array}{c} 
g(T,X)=g_1+Tg_2+\cdots+T^{r-1}g_r\in\gk[T,X],    \\[.3em] 
R(T)=R(f,g_1,\ldots,g_r)(T)=\Res_X\big(f,g(T,X)\big)\in\gk[T],    \\[.25em] 
\fb=\fR(f,g_1,\ldots,g_r)\eqdefi\rc_{\gk,T}\big(R(f,g_1,\ldots,g_r)(T)\big)\subseteq\gk.  
 \end{array}
\]
\begin{enumerate}
\item The \id $\fb$ is generated by $d(r-1) + 1$ \elts and we have the inclusions 
\begin{equation}\preskip-.4em \postskip.3em
\label{eq0lemElimPlusieurs} 
\fb\subseteq\fa\subseteq\sqrt\fb=\sqrt\fa\,.
\end{equation}
More \prmtz, let $e_i = 1 + (d-i)(r-1)$, $i \in \lrb{1..d}$, then for arbitrary \elts $a_1$, $\ldots$, $a_d \in \fa$, we have
$$\preskip.3em \postskip.3em
a_1^{e_1} a_2^{e_2} \cdots a_d^{e_d}\, \in \,\fR(f, g_1,\ldots, g_r)\,.
$$
In particular, we have the following \eqvcs 
\begin{equation}\preskip.2em \postskip.2em
\label{eq01lemElimPlusieurs} 1\in \fb\; \iff\; 1\in\fa \;\iff\; 1\in\ff\,.~\;\; 
\end{equation}

\item If  $\gk$ is a \cdi contained in a discrete \cacz~$\gL$, let $h$ be the \mon gcd of $f$, $g_1$, $\ldots$, $g_r$ and $V$ be the set of zeros of~$\ff$ in~$\gL^n$. Then, we have the following \eqvcs
\begin{equation}\preskip.3em \postskip.4em
\label{eqlemElimPlusieurs} 1\in \fb\; \iff\; 1\in\fa \;\iff\; 1\in\ff \;\iff\;h=1\;\iff\;
V=\emptyset~\;\;~\;\;
\end{equation}

%
%
\end{enumerate}
\end{lemma}
%
\begin{proof}
\emph{1.} We know that $R(T)$ is of the form 

\snic{u(T,X)f(X)+v(T,X)g(T,X),}

so each \coe of $R(T)$ is a \coli of $f$ and the~$g_i$'s in~$\kX$. This gives the inclusion $\fb\subseteq\fa$. The in\egt $\deg_T(R)\leq d(r-1)$ gives the
majoration $d(r-1) + 1$ for the number of \gtrs of $\fb$.
\\
If  $f_1$, $\ldots$, $f_d$ are $d$ \pols (with one \idtrz) of degree $< r$, we deduce from the \DKM lemma (see Corollary~\ref{corLDM}) the following inclusion.
$$\preskip-.1em \postskip.2em 
 \rc(f_1)^{e_1} \rc(f_2)^{e_2} \cdots \rc(f_d)^{e_d}  \subseteq \rc(f_1f_2\cdots f_d).\eqno(\star)
$$

Assume $f(X) = (X-x_1) \cdots (X-x_d)$. Then let for $i \in \lrb{1..d}$
$$\preskip.4em \postskip.4em 
f_i(T) = g_1(x_i) + g_2(x_i)T + \cdots + g_r(x_i)T^{r-1}, 
$$
such that $f_1 f_2 \cdots f_d = \Res_X(f, g_1 + g_2T + \cdots +
g_rT^{r-1})$.  
\\
Thus, for $a_j \in \fa= \gen {f, g_1, \ldots, g_r}_{\kX} \cap \gk$, by evaluating at $x_i$, we obtain~\hbox{$a_j \in \gen {g_1(x_i), \ldots, g_r(x_i)} =
\rc(f_i)$}. By applying the inclusion $(\star)$ we obtain the membership $a_1^{e_1} a_2^{e_2} \cdots a_d^{e_d}\in\fb$.
\\
Let us move on to the \gnl case. Consider the \aduz~\hbox{$\gk'=\Adu_{\gk,f}$}.
The previous computation is valid for $\gk'$. Since $\gk'=\gk\oplus E$
as a \kmoz, we have the \egt $(\fb\gk')\cap \gk=\fb$. For some $a_j\in\fa$, this allows us to conclude that $a_1^{e_1} a_2^{e_2} \cdots a_d^{e_d}\in\fb$, because the product is \hbox{in $(\fb\gk')\cap \gk$}.

\emph{2.} 
By \dfn of the gcd, we have $\ff = \gen {h}$. Moreover, $h=1\Leftrightarrow V=\emptyset$. So the rest clearly follows by item \emph{1.}\\
\emph{Here is however a more direct \dem for this particular case, which gives the point of origin of the magical \dem of {1.}}
\\
Assume that $h$ is equal to $1$; then in this case $1\in\ff$ and $1\in\fa$. Assume next that $h$ is of degree $\geq 1$; then $\fa=\gen{0}$. We therefore have obtained the \eqvcs $\;1\in\fa \iff 1\in\ff \iff \deg(h)=0\;$ and $\;\fa=\gen{0} \iff \deg(h) \ge 1$.
\\
Let us now prove the \eqvc $\;\deg(h) \ge 1 \iff \fb = \gen{0}$. \\
 If $\deg(h)\geq1$, then $h(X)$ divides $g(T,X)$, so $R(f,g_1,\ldots,g_r)(T)=0$ (Fact~\ref{factGCDDeg}), \cad $\fb=\gen{0}$.
\\ Conversely, assume $\fb=\gen{0}$. Then, for all values of the parameter~$t\in\gL$, the \pols $f(X)$ and $g(t,X)$ have a common zero in $\gL$ ($f$ is \mon and the resultant of both \pols is null).
\\
Consider the zeros $\xi_1$, $\ldots$, $\xi_d \in \gL$ of $f$. By taking $d(r-1)+1$ distinct values of $t$, we find some $\xi_\ell$ such that $g(t, \xi_\ell) = 0$ for at least $r$ values of $t$.  This implies that $g(T,\xi_\ell)$ is zero everywhere, \cad that $\xi_\ell$ annihilates all the $g_i$'s, and that $h$ is a multiple of $X-\xi_\ell$, therefore $\deg(h) \ge 1$.
\end{proof}

Item \emph{2} of Lemma~\ref{lemElimPlusieurs} gives the following corollary.
\begin{corollary}\label{corlemElimPlusieurs}
Let $\gK$ be a nontrivial \cdi contained in an \cac $\gL$. Given the hypotheses of Lemma~\ref{lemElimPlusieurs}, with the \ri $\gk=\gK[X_1,\ldots,X_{n-1}]$, then, for  $\alpha=(\alpha_1,\ldots,\alpha_{n-1})\in\gL^{n-1}$ 
\propeq
\begin{enumerate}
\item There exists a $\xi\in\gL$ such that $(\alpha,\xi)$ annihilates $(f,g_1,\ldots,g_r)$.
\item $\alpha$ is a zero of the \id $\fb=\fR(f,g_1,\ldots,g_r)\subseteq\gk$.
\end{enumerate}
Note: if the total degree of the \gtrs of $\,\ff$ is bounded
above by~$d$,
we obtain as \gtrs of $\,\fb$, $d(r-1)+1$ \pols of total degree bounded by~$2d^2$.  
\end{corollary}

\rem The above corollary has the desired structure to step through
an \recu which allows for a description of the zeros of $\ff$ in $\gL^n$. 
\\
Indeed, by starting from the \itf $\ff\subseteq\KXn$ we produce a \itf $\fb\subseteq\gk$ with the following \prtz: \emph{the zeros of $\ff$ in $\gL^n$ are exactly projected onto the zeros of $\fb$ 
in $\gL^{n-1}$}. More \prmtz, above each zero of $\fb$ in $\gL^{n-1}$ there is a finite, nonzero number of zeros of $\ff$ in $\gL^n$, bounded by $\deg_{X_n}(f)$.
\\
 So either all the \gtrs of $\fb$ are zero and the process describing the zeros of $\ff$ is complete, or one of the \gtrs of $\fb$ is nonzero and we are ready to do to $\fb\subseteq\gK[X_1,\ldots,\alb X_{n-1}]$ what we did to $\ff\subseteq\KXn$ \emph{on the condition however that we find} a \polu in $X_{n-1}$ in the \idz~$\fb$.
\\
 This final question is resolved by the following \cdv lemma. 
\eoe

\begin{lemma}
\label{lemCDV} \emph{(\Cdv lemma)}\\
Let $\gK$ be an infinite \cdi and  $g\neq0$ in $\KuX=\KXn$ of degree $d$.
There exists $(a_1, \ldots, a_{n-1}) \in \gK^{n-1}$ such that the \pol

\snic {
g(X_1 + a_1X_n, \ldots, X_{n-1} + a_{n-1} X_n, X_n) 
}

is of the form $aX_n^d + h$ with  $a \in \gK^\times$ and $\deg_{X_n} h < d$.
\end{lemma}
\begin{proof}
Let $g_d$ be the \hmg components of degree $d$ of $g$. Then

\snic {
g(X_1 + a_1X_n, \ldots, X_{n-1} + a_{n-1} X_n, X_n) =
g_d(a_1, \ldots, a_{n-1}, 1) X_n^d + h, 
}

with  $\deg_{X_n} h < d$. Since $g_d(X_1, \ldots, X_n)$ is nonzero and \hmgz,
the \polz~$g_d(X_1, \ldots, X_{n-1}, 1)$ is nonzero. \\
There thus exists $(a_1, \ldots, a_{n-1}) \in \gK^{n-1}$ such that $g_d(a_1, \ldots, a_{n-1}, 1) \ne 0$.
\end{proof}

We now obtain a \gui{weak \nstz} (\cad the \eqvc between~$\,V=\emptyset\,$  and~$\,\gen{\lfs}=\gen{1}\,$ in the \thoz)
and a \gui{\iNoe position} which gives a description of $V$ in the nonempty case.

\begin{theorem}\label{thNstfaibleClass}\emph{(Weak \nst and \iNoe position)}\\
Let $\gK$ be an infinite \cdi contained in an \cac $\gL$ and $(\lfs)$ a \syp in $\KXn$.%
\index{polynomial system}
\\
 Let $\ff=\gen{\lfs}_\KuX$ and $V$ be the \vrt of the zeros of $(\lfs)$ in~$\gL^n$.
\begin{enumerate}
\item Either $\gen{\lfs}=\gen{1},$ and $V=\emptyset$.
\item Or $V\neq\emptyset$. Then there exist an integer $r\in\lrb{0..n}$, a $\gK$-\lin \cdv (the new variables are denoted by $\Yn$), and  \itfs $\ff_j\subseteq\gK[Y_1,\ldots,Y_j]$
$(j\in\lrb{r..n})$, which satisfy the following \prtsz.
\begin{itemize}
\item [$\bullet$] We have $\ff\,\cap\,\gK[\Yr]=0$. In other words, the \ri  $\gK[\Yr]$ is identified with a sub\ri of the quotient \ri $\KuX\sur\ff$.
\item [$\bullet$] Each $Y_j$ $(j\in\lrb{r+1..n})$ is integral over $\gK[\Yr]$ modulo~$\ff$.
In other words the \ri $\KuX\sur\ff$ is integral over the sub\ri $\gK[\Yr]$. 
\item [$\bullet$] We have the inclusions $\gen{0}=\ff_r\subseteq\ff_{r+1} \subseteq \ldots \subseteq\ff_{n-1}\subseteq\ff$ and for each~\hbox{$j\in\lrb{r..n}$}
we have the \egt $\sqrt {\ff\,} \cap \gK[Y_1,\ldots,Y_j]=\sqrt {\ff_j}$.
\item [$\bullet$] For the new \coos corresponding to the $Y_i$'s,
let $\pi_j$ be the \prn $\gL^n\to\gL^j$
which forgets the last \coos ($j\in\lrbn$).
For each $j\in\lrb{r..n-1}$ the \prn of the \vrt $V\subseteq\gL^n$ over~$\gL^j$ 
is exactly the \vrt $V_j$ of the zeros of $\ff_j$. In addition, for each \elt~$\alpha$ of~$V_j$, the fiber $\pi_j^{-1}(\alpha)$ is finite, nonempty, with a uniformly bounded number of \eltsz.
\end{itemize}
\end{enumerate}
In particular
\begin{itemize}
\item Either $V$ is empty (and we can concede that~$r = -1$).
\item Or $V$ is finite and nonempty, $r=0$ and the \coos of the points of $V$ are \agqs over~$\gK$. 
\item Or $r\geq1$ and the \prn $\pi_r$ surjectively sends $V$ onto $\gL^r$ (so $V$ is infinite).
In this case, if $\alpha\in\gK^r$, the \coos of the points of~$\pi_r^{-1}(\alpha)$ are \agqs over~$\gK$.
\end{itemize}
\end{theorem}
%
\begin{proof}
We reason as stated in the remark preceding the \cdv lemma. Note that the first step of the process only takes place if the initial \syp is nonzero, in which case the first operation consists in a \lin \cdv which makes one of the $f_i$'s \mon in~$Y_n$.
\end{proof}

\rems ~\\
 1) The number $r$ above corresponds to the maximum number of \idtrs for a \pol \ri $\gK[Z_1,\ldots,Z_r]$ which is \isoc to a \Kslg of $\aqo\KuX{\lfs}$. This is related to \ddk theory which will be presented in Chapter~\ref{chapKrulldim} (see especially \Thref{thDKAG}).

2) Assume that the degrees of the $f_j$'s are bounded above by $d$.
\\
 By basing ourselves on the result stated at the end of Corollary~\ref{corlemElimPlusieurs}, we can give some bounds in the previous \tho by computing a priori, solely according to the integers $n$, $s$, $j$ and $d$,
\begin{itemize}
\item on the one hand an upper bound for the number of \gtrs for each \idz~$\ff_j$,
\item on the other hand an upper bound for the degrees of these \gtrsz.
\end{itemize}

3) The computation of the \ids $\ff_j$ as well as all the statements of the \tho which do not concern the \vrt $V$ are valid even when we do not know of some \cac $\gL$ containing $\gK$.
To do this, we only use Lemmas~\ref{lemElimPlusieurs} and~\ref{lemCDV}.
We will look at this in more detail in \Thosz~\ref{thNstfaibleClassSCA} and~\ref{thNstNoe}.
\eoe

\medskip
The restriction introduced by the hypothesis \gui{$\gK$ is infinite} will vanish in the classical \nst because of the following fact.

\begin{fact}\label{factNstRationnel}
Let $\gK\subseteq\gL$ be \cdis and $h$, $f_1$, \dots, $f_s\in\KXn$, then 
$h\in\gen{\lfs}_\KXn\iff h\in\gen{\lfs}_\LXn$.  
\end{fact}
%
\begin{proof}
Indeed, an \egt $h=\sum_ia_if_i$, once the degrees of the $a_i$'s are fixed, can be seen as a \sli whose unknowns are the \coes of the $a_i$'s.
The fact that a \sli admits a solution does not depend on the field in which we look for the solution, so long as it contains the \coes of the \sliz; the pivot method is a completely rational process.
\end{proof}

As a corollary of the weak \nst and from the previous fact we obtain the classical \nstz. 

\pagebreak

\begin{theorem}\label{thNstClass} \emph{(Classical \nstz)}
\\
Let $\gK$ be a \cdi contained in an \cac $\gL$ and $g$, $f_1$, \dots, $f_s$ 
be  some \pols in $\KXn$. 
Let $V$ be the \vrt of the zeros of $(\lfs)$ in $\gL^n$. 
Then either 1.\ there exists a point~$\xi$ of $V$ such that $g(\xi)\neq0$,
or 2.\ there exists an integer $N$ such that $g^N\in\gen{\lfs}_\KuX$.
\end{theorem}
%
\begin{proof} The $g=0$ case is clear, so we suppose $g\neq 0$.
We apply the \emph{Rabinovitch trick}\index{Rabinovitch!trick}, \cad we introduce an additional \idtr $T$ and we notice that~$g$ is annihilated at the zeros of $(\lfs)$ \ssi the \sys \hbox{$(1-gT,\lfs)$} admits no solution. Then we apply the weak \nst to this new \sypz, with  $\gL$ (which is infinite) instead of $\gK$.
We obtain \hbox{in $\KuX[T]$} (thanks to Fact~\ref{factNstRationnel}) an \egt
$$
\big(1-g(\uX)T\big)a(\uX,T) + f_1(\uX)b_1(\uX,T)+\cdots+f_s(\uX)b_s(\uX,T)=1.
$$
In the localized \ri $\KuX[1/g]$, we perform the substitution $T = 1/g$.
More \prmtz, by remaining in $\gK[\uX,T]$, if $N$ is the greatest of the degrees in~$T$ of the~$b_i$'s, we multiply the previous \egt by $g^N$ and we replace \hbox{in $g^Nb_i(\uX,T)$} \hbox{each $g^NT^k$} \hbox{by $g^{N-k}$} modulo~\hbox{$(1-gT)$}. We then obtain an \egt
$$\preskip.2em \postskip.4em
~~~\big(1-g(\uX)T\big)a_1(\uX,T) + f_1(\uX)c_1(\uX)+\cdots+f_s(\uX)c_s(\uX)=g^N,
$$
in which $a_1=0$ \ncrtz, since if we look at $a_1$ in~$\KuX[T]$, its \fmt leading \coe in $T$ is zero.
\end{proof}

\rem Note that the separation of the different cases in \Thosz~\ref{thNstfaibleClass} and~\ref{thNstClass} is explicit.
\eoe

\begin{corollary}\label{corthNstClass}
Let $\gK$ be a \cdi contained in an \cac $\gL$ and $\fa=\gen{\lfs}$, $\fb$ be two \itfs \hbox{of $\KXn$}.
Let $\gK_0$ be the subfield of $\gK$ generated by the \coes of \hbox{the $f_i$'s}.
\\
\Propeq
\begin{enumerate}
\item $\fb\subseteq\rD_{\KuX}(\fa)$.
\item $\fb \subseteq\rD_{\LuX}(\fa)$.
\item  Every zero of $\fa$ in $\gL^n$ is a zero of $\fb$.
\item  For every subfield $\gK_1$ of $\gL$
finite over $\gK_0$, every zero of $\fa$  in $\gK_1^n$ is a zero of~$\fb$.
\end{enumerate}
In particular, $\rD_{\KuX}(\fa)=\rD_{\KuX}(\fb)$ \ssi $\fa$ and $\fb$ have the same zeros in $\gL^n$.
\end{corollary}
%
\begin{proof}
Immediate consequence of the \nstz. 
\end{proof}

\vspace{-.7em}
\pagebreak

\subsec{The formal \nstz}
We now move onto a \emph{formal \nstz}, formal in the sense that it applies (in \clamaz) to an arbitrary \id over an arbitrary \riz. 
Nevertheless to have a \cof statement we will be content with a \pol \ri $\ZZuX$ for our arbitrary \ri and  a \itf for our arbitrary \idz. 

Although this may seem very restrictive, practice shows that this is not the case because we can (almost) always apply the method of undetermined coefficients to a commutative algebra problem;
 a method which reduces the \pb to a \poll \pb over~$\ZZ$. An illustration of this will be given next.

Note that to read the statement, when we speak of a zero of some $f_i\in\ZZuX$ over a \ri $\gA$, one must first consider $f_i$ modulo $\Ker\varphi$, where $\varphi$ is the unique \homo $\ZZ\to\gA$, with $\gA_1\simeq\ZZ\sur{\Ker\varphi}$ as its image. 
This thus reduces to a \polz~$\ov {f_i}$ of $\gA_1[\uX]\subseteq\gA[\uX]$.

\begin{theorem}\label{thNSTsurZ}
\emph{(\nst over $\ZZ$, formal \nstz)} 
\\
Let $\ZZuX=\ZZXn$.
Consider $g$, $f_1$, \dots, $f_s$ in $\ZZuX$
\begin{enumerate}
\item For the \sys $(\lfs)$ \propeq
\begin{enumerate}
\item $1\in\gen{\lfs}$.
%
%
\item The \sys does not admit a zero on any nontrivial \cdiz. 
\item The \sys does not admit a zero on any finite field or on any finite extension of $\QQ$.
\item The \sys does not admit a zero on any finite field.
\end{enumerate}
\item \Propeq
\begin{enumerate}
\item $\exists N \in \NN,\;g^N\in\gen{\lfs}$.
%
%
\item The \polz~$g$ is annihilated at the zeros of the \sys $(\lfs)$ on any \cdiz.
\item The \polz~$g$ is annihilated at the zeros of the \sys $(\lfs)$ on every finite field and on every finite extension of $\QQ$.
\item The \polz~$g$ is annihilated at the zeros of the \sys $(\lfs)$ on every finite field.
\end{enumerate}
\end{enumerate}
\end{theorem}

%
\begin{proof}
It suffices to prove the weak version \emph{1}, as we can then get the \gnl version \emph{2} by applying the Rabinovitch trick.
Regarding the weak version, the difficult task is the implication \emph{d} $\Rightarrow$ \emph{a}.

Let us first deal with \emph{c} $\Rightarrow$ \emph{a}. Apply the weak \nst by considering~\hbox{$\ZZ\subseteq\QQ$}. 
This gives the membership
$$
m \in \gen {f_1, \ldots, f_s}_{\ZZuX}
\quad\hbox{ with } m\in\ZZ\setminus\so0\eqno(\star_\QQ)
.$$ 
By applying the weak \nst with an \agq closure $\gL_p$ of $\FFp$
we also obtain for each prime number $p \divi m$ a membership
$$
1 \in \gen {f_1, \ldots, f_s}_{\ZZuX} +p\ZZuX \eqno(\star_{\FFp})
.$$ 
However, in any \riz, for three arbitrary \ids $\fa, \fb, \fc$,
we have the inclusion $(\fa + \fb)(\fa + \fc) \subseteq \fa + \fb\fc$.
By expressing the above $m$ in $(\star_\QQ)$ in the form $\prod_jp_j^{k_j}$
with prime $p_j$'s, we therefore obtain
$$
1 \in \gen {f_1, \ldots, f_s}_{\ZZuX} + m\ZZuX.
$$
This membership, joint with $(\star_\QQ)$, provides 
$1 \in \gen {f_1, \ldots, f_s}_{\ZZuX}$.

\emph{d} $\Rightarrow$ \emph{c.} We show that a zero $(\uxi)$ of the \sys $(\lfs)$ in a finite extension of $\QQ$ leads to a zero of $(\lfs)$ in a finite extension of $\FFp$ for all the prime numbers, except for a finite number of them.
\\
 Indeed, let $\gQ=\QQ[\alpha]\simeq \aqo{\QQX}{h(X)}$ (with $h$ \ird and \mon in~$\ZZX$) be a finite extension of $\QQ$ and $(\uxi)\in\gQ^n$ be a zero of $(\lfs)$.
If~\hbox{$\xi_j=q_j(\alpha)$} with $q_j\in\QQX$ for $j\in\lrbn$, this means that 
$$
f_i(q_1,\ldots,q_n)\equiv 0\mod h \quad\hbox{in }\QQX, \;i\in\lrbs.
$$ 
This remains true in $\FFp[X]$ as soon as none of the \denos appearing in the $q_j$'s is a multiple of $p$, provided one takes the fractions from $\FFp$
$$
\ov{f_i}(\ov{q_1},\ldots,\ov{q_n})\equiv 0\mod \ov{h} \quad\hbox{in }\FFp[X], \;i\in\lrbs.
$$ 
For such a $p$, we take an \ird \mon divisor $h_p(X)$ of $\ov{h}(X)$ in~$\FFp[X]$ and consider the finite field $\gF=\aqo{\FFp[X]}{h_p(X)}$ with $\alpha_p$ the class of~$X$. Then, $\big(q_1(\alpha_p),\ldots,q_n(\alpha_p)\big)$ is a zero of $(\lfs)$ in~$\gF^n$.
\end{proof}

We have the following \imd corollary, with \itfsz. 

\begin{corollary}\label{corthNSTsurZ}
\emph{(\nst over $\ZZ$, formal \nstz, 2)} 
\\
Write $\ZZuX=\ZZXn$.
For two \itfs $\fa$, $\fb$ of $\ZZuX$ \propeq
\begin{enumerate}
\item $\rD_{\ZZuX}(\fa)\subseteq\rD_{\ZZuX}(\fb)$.
\item  
 $\rD_{\gK}\big(\varphi(\fa)\big)\subseteq\rD_{\gK}\big(\varphi(\fb)\big)$ for every \cdi $\gK$ and every \homo $\varphi:\ZZuX\to\gK$.
%
%
\item Idem but restricted to \agq extensions of $\QQ$ and to finite fields.
\item Idem but restricted to finite fields.
\end{enumerate}
\end{corollary}

\mni\textbf{An application example}

\mni Consider 
the following result, already proven in Lemma~\ref{lemGaussJoyal}:
\emph{An \eltz~$f$ of $\AuX$ is \iv \ssi $f(\uze)$ is \iv and $f-f(\uze)$
is nilpotent. In other words $\AuX\eti = \Ati+\DA(0)[\uX]$.}
\\
We can assume that $fg=1$ with $f=1+Xf_1$ and $g=1+Xg_1$. Consider the \coes of $f_1$ and $g_1$ as being \idtrsz. We are brought to prove the following result.
\\
\emph{An \egt $\;f_1+g_1+Xf_1g_1=0\;\;(*)\;$ implies that the \coes of~$f_1$ are nilpotent.}
\\
However, since the \idtrs are evaluated in a field, the \coes of~$f_1$ are annihilated at the zeros of the \syp in the \idtrs given by the \egtz~$(*)$.
We conclude with the formal \nstz.

 When compared with the \dem given for item \emph{4} of Lemma~\ref{lemGaussJoyal}, we can assert that the one given here is both simpler (no need to find a more subtle computation) and cleverer (usage of the formal \nstz).

\emph{Note.} Another example is given in the solution to \Pbmz~\ref{exoChasserIdeauxPremiers1}.
\eoe

\section{Newton's method in \algz}  \label{secNewton}

Let $\gk$ be a \ri and $f_1$, $\ldots$, $f_s\in\kuX= \gk[\Xn]$.
The \emph{Jacobian matrix} of the \sys is the matrix%
\index{matrix!Jacobian ---}\imN
$$\preskip.4em \postskip.4em \JJ_{\Xn}(f_1,\ldots,f_s) =
\Big( \frac{\partial f_i}{\partial X_j}
\Big)_{i\in \lrbs,j\in \lrbn }
\in \kuX^{s\times n}.
$$
It is also denoted by $\JJ_{\uX}(\uf)$ or $\JJ(\uf)$.
It is visualized as follows
$$\preskip.3em \postskip.4em
\bordercmatrix [\lbrack\rbrack]{
    & X_1                     & X_2                     &\cdots  & X_n \cr
f_1 & \Dpp {f_1}{X_1} &\Dpp {f_1}{X_2}  &\cdots  &\Dpp {f_1}{X_n} \cr
f_2 & \Dpp {f_2}{X_1} &\Dpp {f_2}{X_2}  &\cdots  &\Dpp {f_2}{X_n} \cr
f_i & \vdots                  &                         &        & \vdots              \cr
 & \vdots                  &                         &        & \vdots              \cr
f_s & \Dpp {f_s}{X_1} &\Dpp {f_s}{X_2}  &\cdots  &\Dpp {f_s}{X_n} \cr
}
.$$
If  $s=n$, we denote by $\J_{\uX}(\uf)$
or $\J_{\Xn}(f_1,\ldots,f_n)$ or $\J(\uf)$ 
the \ixc{Jacobian}{of a polynomial system} of the \sys $(\uf)$, 
\cad the \deter of the Jacobian matrix.

In analysis Newton's method to approximate a root of a differentiable function $f:\RR\to\RR$ is the following. Starting from a point $x_0$ \gui{near a root,} at which the derivative is \gui{far from $0$}, we construct a series $(x_m)_{m\in\NN}$ by \recu by letting
$$\preskip.0em \postskip.4em
\;\;x_{m+1}=x_m-\frac{f(x_m)}{f'(x_m)}.
$$
The method can be generalized for a \sys of $p$ \eqns with $p$ unknowns. A solution of such a \sys is a zero of a function~\hbox{$f:\RR^p\to\RR^p$}. We apply \gui{the same formula} as above
$$\preskip.3em \postskip.4em
x_{m+1}=x_m-f'(x_m)^{-1}\cdot f(x_m),
$$
where $f'(x)$ is the \dile (the Jacobian matrix) of $f$ at the point $x\in\RR^p$, which must be \iv in a neighborhood of $x_0$.

This method, and other methods of the infinitesimal calculus, 
can also be applied in certain cases in \algz, by replacing the Leibnizian infinitesimals by the nilpotent \eltsz.

If for instance $\gA$ is a \QQlg and $x\in\gA$ is nilpotent, the formal series

\snic{1+x+x^2/2+x^3/6+\ldots }

which defines $\exp(x)$ only has a finite number of nonzero terms in
$\gA$ and therefore defines an \elt $1+y$ with $y$ nilpotent. 
Since the \egt

\snic{\exp(x+x')=\exp(x)\exp(x'),}

holds in analysis,  it is also valid with regard to formal series over $\QQ$. So when $x$ and $x'$ are nilpotents in $\gA$ we will obtain the same \egt in $\gA$. Similarly the formal series
$$\preskip.2em \postskip.2em 
y-y^2/2+y^3/3-\ldots  
$$
 which defines $\log(1+y)$, only has a finite number of terms in $\gA$ when $y$ is nilpotent and allows for a definition of $\log(1+y)$ as a nilpotent \elt of~$\gA$. 
 Furthermore, for nilpotent $x$ and $y$, we obtain the \egts 
 
 \snic{\log\big(\exp(x)\big)=x\hbox{ and }\exp\big(\log(1+y)\big)=1+y}
 
as consequences of the corresponding \egts for the formal series.

In a similar style we easily obtain, by using the inverse formal series of $1-x$, the following result.
\begin{lemma}
\label{lemMatInvNil} \emph{(\Rdt \iv \elts lemma)}
\begin{enumerate}
\item If $ef\equiv1$ modulo the nilradical, then $e$ is \iv and
$$\preskip.2em \postskip.4em \ndsp
e^{-1}=f\,\sum_{k\geq 0}(1-ef)^k . 
$$
\item A square matrix $E\in\Mn(\gA)$ \iv modulo the nilradical is \ivz. Assume that $d\det(E)\equiv1$ modulo the nilradical.\\
Let $F=d\wi{E}$ (where $\wi{E}$ is the cotransposed matrix of $E$).
 Then,  $E^{-1}$  is in the sub\ri of $\Mn(\gA)$ generated by the \coes of the \polcar of $E$, $d$  and $E$.\\ 
 More \prmtz,
 the matrix $\In-EF = \big(1 - d\det(E)\big)\In$ is nilpotent~and

 \snic{E^{-1}=F\,\sum_{k\geq 0}\big(1 - d\det(E)\big)^k.}
\end{enumerate}
\end{lemma}

Let us move on to Newton's method. 

\begin{theorem}
\label{thNewtonLin} \emph{(Newton's \lin method\imN)}\\
Let $\fN$ be an \id of $\gA$,
 $\uf=\tra{[\,f_1\;\cdots\;f_n\,]}$ be a vector whose \coos are \pols in $\AXn$, and $\ua=\tra{(a_1,\ldots ,a_n)}$ in~$\Ae n$ be a \emph{approximated simple zero} of the \sys in the following sense.
\begin{enumerate}
\item [--] The Jacobian matrix $J({\ua})$ of $\uf$ at point $\ua$ is
\iv modulo~$\fN$; let $U\in\Mn(\gA)$ be such an inverse.
\item [--] The vector $\uf(\ua)$ is null modulo~$\fN$.
\end{enumerate}
Consider the sequence $(\ua^{(m)})_{m\geq 1}\in\gA^n$  defined by Newton's \lin iteration
$$\preskip.2em \postskip.4em
\ua^{(1)}=\ua, \,\,\,\,\ua^{(m+1)}=\ua^{(m)}- U \cdot \uf(\ua^{(m)})
.$$
\begin{enumerate}
\item [a.]
This sequence satisfies the following $\fN$-adic requirements:
$$\preskip.2em \postskip.4em
\ua^{(1)}\equiv\ua\,\,\,\mod\,\fN, \,\mathrm{\,and\,} \Tt m,\,\,\,
\ua^{(m+1)}\equiv\ua^{(m)} \,\, \mathrm{and}\,\,
\uf(\ua^{(m)})\equiv 0 \,\,\,\,\mod\,\fN^m.
$$
\item [b.] This sequence is unique in the following sense, if $\ub^{(m)}$ is another sequence satisfying the requirements of a., then for all $m$,
$\ua^{(m)}\equiv\ub^{(m)}\,\,\,\mod\,\fN^m.$
\item [c.] Let $\gA_1$ be the sub\ri generated by the \coes of the $f_i$'s, by those of $U$ and by the \coos of $\ua$.
In this \ri let $\fN_1$ be the \id generated by the \coes of $\In-UJ(\ua)$ and the \coos of $\ua$.
If the \gtrs of $\fN_1$ are nilpotent, the sequence converges in a finite number of steps towards a zero of the \sys $\uf$, and it is the unique zero of the \sys congruent to $\ua$ modulo~$\fN_1$.
\end{enumerate}
\end{theorem}

Under the same assumptions, we have the following quadratic method.

\begin{theorem}
\label{thNewtonQuad} \emph{(Newton's quadratic method\imN)}\\
Let us define the sequences $(\ua^{(m)})_{m\geq
0}$ in $\gA^n$ and $(U^{(m)})_{m\geq 0}$ in $\Mn(\gA)$ by the following Newton quadratic iteration 
$$\begin{array}{lcl}
\ua^{(0)}=\ua,&  \quad \quad  & \ua^{(m+1)}=\ua^{(m)}- U^{(m)} \cdot
\uf(\ua^{(m)}),    \\[1mm]
U^{(0)}=U,&   & U^{(m+1)}=U^{(m)}\,\left(2\I_n-J(\ua^{(m+1)})U^{(m)}\right).
\end{array}$$
Then, we obtain for all $m$ the following congruences:
$$
\begin{array}{lcll}
 \ua^{(m+1)}\equiv\ua^{(m)} & \,\,\mathrm{and}\,\, &
  U^{(m+1)}\equiv U^{(m)} &\,\mod \,\fN^{2^m}    \\
  \uf(\ua^{(m)})\equiv 0 &  \,\,\mathrm{and}\,\,  &
   U^{(m)}\,J(\ua^{(m)})\equiv \In  &\,\mod \,\fN^{2^m}.
\end{array}
$$
\end{theorem}

The \dems are left to the reader (cf. \cite{ValuHallouin}) by observing that the iteration concerning the inverse of the Jacobian matrix can be justified by Newton's \lin method or by the following computation in a not \ncrt commutative \ri
$$\preskip-.2em \postskip.2em
(1-ab)^2=1-ab'\quad \mathrm{with}\quad b'=b(2-ab).
$$

\pagebreak

\begin{corollary}
\label{corIdmNewton} \emph{(Residual \idms lemma)}
\begin{enumerate}
\item For every commutative \ri $\gA$:
\begin{enumerate}
\item two equal \idms modulo $\DA(0)=\sqrt{\gen{0}}$ are equal;
\item every \idm $e$ modulo an \id $\fN$ is uniquely lifted to an \idm $e'$ modulo $\fN^2$;
Newton's quadratic iteration is given by $e\mapsto 3e^2-2e^3$.
\end{enumerate}
\item Similarly every matrix $E\in\Mn(\gA)$ \idme modulo~$\fN$ is lifted to a matrix $F$ \idme modulo $\fN^2$. The \gui{lifting} $F$ is unique provided that $F\in\gA[E]$. Newton's quadratic iteration is given by $E\mapsto 3E^2-2E^3$.
\end{enumerate}
\end{corollary}
\begin{proof}
\emph{1a.} Left to the reader. A stronger version is proven in Lemma~\ref{lemIdmIsoles}. 
\\
\emph{1b.} Consider the \hbox{\polz~$T^2-T$}, and note that $2e-1$ is \iv modulo~$\fN$ since $(2e-1)^2=1$ modulo $\fN$.
\\
\emph{2.} Apply item \emph{1} with the commutative \ri $\gA[E]\subseteq\End(\gA^n)$.
\end{proof}

\Exercices

\begin{exercise}
\label{exoLagrange} (Lagrange interpolation)~
{\rm  Let $\gA$ be a commutative \riz.\index{Lagrange interpolation}%
\index{interpolation!Lagrange ---}
Prove the following statements.
\begin{enumerate}\itemsep=0pt
\item Let $f$, $g\in\AX$ and $a_1$, $\ldots$, $a_k$ be \elts of $\gA$ such that $a_i-a_j\in\Reg\gA$ for $i\neq j$.
\begin{enumerate}
\item If the $a_i$'s are zeros of $f$, $f$ is a multiple of $(X-a_1)\cdots(X-a_k)$.
\item If $f(a_i)=g(a_i)$ for $i\in\lrbk$ and if $\deg(f-g)<k$, then $f=g$.
\end{enumerate}  

\item If $\gA$ is integral and infinite, the \elt $f$ of $\AX$ is \care by the \pol function that $f$ defines over $\gA$. 
\item \emph{(Lagrange interpolation \polz)}
Let $(\xzn)$ be in $\gA$ such that each $x_i-x_j\in\Ati$
(for $i\neq j$).
Then, for $(y_0, \ldots, y_n)$ in $\gA$ there exists exactly one \polz~$f$ of degree $\leq n$ such that for each $j\in\lrb{0..n}$ we have $f(x_j)=y_j$.
\\
More \prmtz, the \polz~$f_i$ of degree $\leq n$ such that $f_i(x_i)=1$
and $f_i(x_j)=0$ for $j\neq i$ is equal to

\snic{f_i= \frac{\prod\nolimits_{j\in\lrb{0..n},j\neq i}(X-x_j)}{\prod\nolimits_{j\in\lrb{0..n},j\neq i}(x_i-x_j)},}

and the interpolation \pol $f$ above is equal to $\som_{i\in\lrbzn}y_if_i$.  

\item With the same assumptions, letting $h=(X-x_0)\cdots(X-x_n)$,
we obtain an \iso of \Algsz: $\aqo\AX h \to \gA^{n+1},\;\ov g\mapsto\big(g(x_0),\ldots,g(x_n)\big)$.
\item Interpret the previous results with  \lin \alg (Vandermonde matrix and \deterz) and with  the Chinese remainder \tho (use the pairwise \com \ids $\gen{X-x_i}$).
\end{enumerate}
 
}
\end{exercise}

\vspace{-.6em}
\pagebreak

\begin{exercise}
\label{exoGensIdealEnsFini} 
       {(Generators of the \id of a finite set)} {\rm See also Exercise~\ref{exoGensPolIdeal}.\\ 
Let $\gK$ be a discrete field and $V \subset \gK^n$ be a finite set. 
Following the steps below show that the \idz~$\fa(V) = \sotq{ f \in \Kux}{ \forall\ w \in V,\ f(w) = 0}$ is generated by $n$ \elts (note that this bound does not depend on $\#V$ and that the result is clear for $n=1$).
We denote by~$\pi_n : \gK^n \to \gK$ the $n^{\rm th}$ \prn and for each~$\xi \in \pi_n(V)$, 

\snic{V_\xi = \sotq {(\xi_1, \ldots, \xi_{n-1}) \in \gK^{n-1}} 
         {(\xi_1, \ldots, \xi_{n-1}, \xi) \in V}.}

\begin{enumerate}
\item
Let $U \subset \gK$ be a finite subset and to each $\xi \in U$, associate a \pol

\snic{Q_\xi
\in \gK[x_1, \ldots, x_{n-1}].}

Find a \polz~$Q \in \Kux$ satisfying $Q(x_1, \ldots, x_{n-1},\xi) = Q_\xi$ 
for all~$\xi \in U$.
\item
Let $V \subset \gK^n$ be a set such that $\pi_n(V)$ is finite. Suppose that for each $\xi \in \pi_n(V)$, the \id $\fa(V_\xi)$ is generated by $m$ \polsz. Show that~$\fa(V)$ is generated by $m+1$ \polsz.
Conclude the result.
\end{enumerate}
}
\end{exercise}

\vspace{-1em}
\begin{exercise}
\label{exothSymEl}
(Detailed \dem of \Thref{thSymEl})\\
{\rm Consider the \ri $\AXn=\AuX$ and let $S_1$, $\ldots$, $S_n$ be the \elr \smq functions of $\uX$. All the considered \pols are formal \polsz, because we do not assume that $\gA$ is discrete. 
We introduce another system of \idtrsz,
$(\underline{s}) = (s_1, \ldots, s_n)$, and on the \ri $\gA[\underline{s}]$ we define the weight $\delta$ by~$\delta(s_i) = i$ (a \fmt nonzero \pol has a well-defined formal weight). 
\\
Denote by $\varphi:\gA[\underline{s}]\to\AuX$ the \evn \homo defined by $\varphi(s_i)=S_i$.
\\
Consider on the \moms of $\AuX = \gA[X_1,\ldots, X_n]$ the {\tt deglex} order for which two \moms are first compared according to their total degree, then according to the lexicographical order with  $X_1 > \cdots > X_n$. This provides for some $f \in \AuX$ (\fmt nonzero) a notion of a \emph{\fmt leading \momz} that we denote by $\md(f)$. This \gui{monomial order} is clearly \isoc to $(\NN,\leq)$.

\smallskip 
\emph{0.} Check that every \smq \pol (\cad invariant under the action of $\Sn$) of~$\AuX$ is equal to some \fmt \smq \polz, \cad invariant under the action of~$\Sn$ as a formal \polz.

\emph{1. (Injectivity of $\varphi$)}\,\\
Let $\alpha = (\alpha_1, \ldots, \alpha_n)$ be a decreasing 
exponent sequence ($\alpha_1 \ge \cdots \ge \alpha_n)$. 
\\
Let $\beta_i=\alpha_i-\alpha_{i+1}$
($i\in\lrb{1..n-1}$). Show that

\snic{\md(S_1^{\beta_1} S_2^{\beta_2} \cdots
S_{n-1}^{\beta_{n-1}} S_n^{\alpha_n}) =
X_1^{\alpha_1} X_2^{\alpha_2} \cdots X_n^{\alpha_n}.}

Deduce that $\varphi$ is injective.

\emph{2. (End of the \dem of items 1  and 2  of \Thref{thSymEl})}\,
Let $f \in \AuX$ be a \fmt \smqz, \fmt nonzero \polz, and $\uX^\alpha = \md(f)$. 

\begin{itemize}
\item Show that $\alpha$ is decreasing.
 Deduce an \algo to express every \smq \pol of $\AuX$ as a \pol in $(S_1, \ldots, S_n)$ with \coes
in $\gA$, \cad in the image of $\varphi$. The halting of the \algo can be proven by \recu on the monomial order, \isoc to $\NN$.
\item As an example, write the symmetrized \pol of the \mom  $X_1^4X_2^2X_3$ 
in $\gA[X_1,\ldots,X_4]$ as a \pol in the $S_i$'s.

\end{itemize}

\goodbreak
\emph{3.} 
\emph{(\Demo of item 3 of the \thoz)} 

\begin{itemize}
\item
Let $g(T) \in \gB[T]$ be a \polu of degree $n \ge 1$. 
Show that $\gB[T]$ is a free $\gB[g]$-module with basis $(1, T, \ldots, T^{n-1})$.\\
Deduce that $\gA[S_1,\ldots,S_{n-1}][X_n]$ is a free module \hbox{over $\gA[S_1,\ldots,S_{n-1}][S_n]$}, with basis $(1, X_n, \ldots, X_n^{n-1})$.
\item
Denote by $\uS' = (S'_1, \ldots, S'_{n-1})$ the \elr \smq functions of the variables~$(X_1, \ldots,\alb X_{n-1})$. 
Show that $\gA[\uS', X_n] = \gA[S_1, \ldots, S_{n-1}, X_n]$.
\item
Deduce from the two previous items that $\gA[\uS',X_n]$ is a free $\gA[\uS]$-module with basis $(1, X_n, \ldots, X_n^{n-1})$.
\item Conclude by \recu on $n$ that the family

\snic{\sotq{X^\alpha}{\alpha=(\aln)\in\NN^n,\,\forall k\in\lrbn, \,\alpha_k<k}}

forms a basis of $\AuX$ over $\gA[\uS]$.
\end{itemize}

\emph{4.} 
\emph{(Another \dem of item 3 of the \thoz, and even more, after reading Section~\ref{sec0adu})}\,\,  
Prove that $\AuX$ is canonically \isoc to the \adu of the \polz~$t^n+\sum_{k=1}^{n}
(-1)^ks_kt^{n-k}$ over the \riz~$\gA[s_1,\ldots,s_n]$.
}
\end{exercise}

\vspace{-1em}
\begin{exercise}
\label{exoPolSym1}
{\rm
Let $S_1$, \ldots, $S_n \in \AuX = \gA[X_1, \ldots, X_n]$ be the $n$ \elr \smq functions.
\begin{enumerate}\itemsep=0pt
\item
For $n = 3$, check that $X_1^3 + X_2^3 + X_3^3 = S_1^3 - 3S_1S_2 + 3S_3$.
Deduce that for all $n$, $\sum_{i = 1}^n X_i^3 = S_1^3 - 3S_1S_2 + 3S_3$.
\item
By using a method analogous to the previous question, express the \pols $\sum_{i \ne j} X_i^2 X_j$, $\sum_{i \ne j} X_i^3 X_j$, $\sum_{i < j} X_i^2 X_j^2$ in terms of the \elr \smq functions.

\item
State a \gnl result.
\end {enumerate}
}
\end{exercise}

\vspace{-1em}
\begin{exercise}\label{exoNewtonSum}
 {(The Newton sums and the complete \smq functions)}\\
{\rm
Let $S_i \in \AuX = \gA[X_1, \ldots, X_n]$ be the \elr \smq functions by agreeing to take $S_i = 0$ for $i > n$ and $S_0 = 1$.
\\
For $r\ge 1$, define the \ix{Newton sums} by $\SNw{r} = X_1^r + \cdots +
X_n^r$. Work in the \ri of formal series~$\gA[\uX][[t]]$ and introduce the series\isN
$$\preskip.3em \postskip-.1em\ndsp 
P(t) = \sum_{r \ge 1} \SNw{r}\, t^r  \quad \hbox{ and } \quad
E(t) = \sum_{r \ge 0} S_r \, t^r. 
$$
\begin {enumerate}\itemsep=0pt
\item
Check the \egt
$P(t) = \sum_{i = 1}^n {X_i \over 1 - X_i t}$.
\item
When $u \in \gB[[t]]$ is invertible, considering the logarithmic derivative

\snic{D_{\rm log}(u) = u' u^{-1},}

show that we get a morphism of groups $D_{\rm log} :
(\gBtst, \times) \to (\gB[[t]], +)$.
\item
By using the logarithmic derivation, prove
{\it Newton's relation}

\snic {
P(-t)= {E'(t) \over E(t)}, \;\;\hbox { or }  P(-t)E(t) = E'(t).
}

\item
For $d \ge 1$, deduce Newton's formula

\snic {
\sum_{r = 1}^d (-1)^{r-1} \SNw{r}\, S_{d-r} = d\,S_d.
}
\end {enumerate}
For $r \ge 0$, we define the \ix{complete symmetric function of degree~$r$} by

\snic{H_r = \som_{|\alpha| = r} \uX^\alpha.}

Thus $H_1 = S_1$, $H_2 = \sum_{i \le
j} X_iX_j$, $H_3 = \sum_{i \le j \le k} X_iX_jX_k$. We define the series

\snic {
H(t) = \sum_{r \ge 1} H_r\, t^r.
}
\begin {enumerate} \setcounter {enumi}{4}
\item
Show the \egt
%
$H(t) = \sum_{i = 1}^n {1 \over 1 - X_i \,t}$.
\item
Deduce the \egt $H(t)\,E(-t) = 1$, then for $d \in \lrb {1..n}$,

\snic {
\sum_{r = 0}^{d} (-1)^r S_r\, H_{d-r} = 0, \quad
H_d \in \gA[S_1, \ldots, S_d], \quad
S_d \in \gA[H_1, \ldots, H_d].
}

\item
Consider the \homo  $\varphi : \gA[S_1,\ldots, S_n]\to\gA[S_1,\ldots, S_n]$ defined \hbox{by $\varphi(S_i) = H_i$}. Show that $\varphi(H_d) = S_d$ for $d \in \lrb {1..n}$. 
Thus 
\begin{itemize}
\item $\varphi \circ \varphi =
\I_{\gA[\uS]}$,
\item  $ H_1$, \ldots, $H_n$ are \agqt independent over $\gA$,
\item  $\gA[\uS] = \gA[\uH]$, and expressing~$S_d$ in terms of $H_1$, $\ldots$, $H_d$
is the same as expressing~$H_d$ in terms of $S_1$,~$\ldots$,~$S_d$.
\end{itemize}
\end{enumerate}
}
\end{exercise}


\vspace{-1em}
\begin{exercise}
\label{exolemArtin} (Equivalent forms of the \DKM lemma)\\
{\rm Prove that the following assertions are \eqves
(each of the assertions is \uvlez, \cad valid for all \pols and all commutative \risz):

\emph{1.} $\rc(f)=\gen{1}\; \Longrightarrow\;  \rc(g)=\rc(fg)$.

\emph{2.} $\exists p\in\NN \; \; \rc(f)^{p}\rc(g)\subseteq \rc(fg)$.

\emph{3.} \emph{(Dedekind-Mertens, weak form)} 
$\;\;\;\;\;\exists p\in\NN \; \;
\rc(f)^{p+1}\rc(g)=\rc(f)^p\rc(fg)$.

\emph{4.} $\Ann\big(\rc(f)\big)=0\; \; \Longrightarrow\; \;
\Ann\big(\rc(fg)\big)=\Ann\big(\rc(g)\big)$.

\emph{5.} \emph{(McCoy)} $\;\;\;\;\;\big(\Ann(\rc(f)\big)=0, \; fg=0)\; \Longrightarrow\;
g=0$.

\emph{6.} $(\rc(f)=\gen{1}, \; fg=0)\; \Longrightarrow\;  g=0$.
}
\end{exercise}

\vspace{-1em}
\begin{exercise}
\label{exoMcCoy}
{\rm
Let $\fc = \rc(f)$ be the content of $f \in \gA[T]$. \DKM lemma gives 
$\Ann_\gA(\fc)[T] \subseteq \Ann_{\gA[T]}(f) \subseteq \DA(\Ann_\gA\big(\fc)\big)[T]$.
Give an example for which there is no \egtz.
}
\end{exercise}

\vspace{-1em}
\begin{exercise}
\label{exothKro}
{\rm Deduce \KROz's \tho (\paref{thKro}) from the \DKM lemma.
 }
\end{exercise}

\vspace{-1em}
\begin{exercise}
\label{exoModCauBase} (Cauchy modules) 
{\rm We can give a very precise explanation for the fact that the \id $\cJ(f)$  (\Dfn~\ref{definotaAdu}) is equal to the \id generated by the Cauchy modules. This works with a beautiful formula. Let us introduce a new variable $T$. Prove the following results.
\begin{enumerate}\itemsep=0pt
\item In $\gA[X_1,\ldots ,X_n,T]=\gA[\uX,T]$, we have
\begin{equation}\label{eqModCau}
\arraycolsep2pt\begin{array}{rcl}
f(T)&   =&f_1(X_1)+(T-X_1)f_2(X_1,X_2)+
\\[1mm]
&&\quad (T-X_1)(T-X_2)f_3(X_1,X_2,X_3)+\cdots+
\\[1mm]
&   &\quad\quad
(T-X_1)\cdots (T-X_{n-1})f_n(X_1,\ldots ,X_n)+\\[1mm]
&&\quad\quad\quad (T-X_1)\cdots (T-X_{n})
\end{array}
\end{equation}
\item In the $\gA[\uX]$-submodule of $\gA[\uX,T]$ formed by the \pols of degree $\leq n$ in $T$,
the \polz~$f(T)-(T-X_1)\cdots (T-X_{n})$ possesses two different expressions.
\begin{itemize}
\item On the one hand, over the basis $(1,T,T^2,\ldots ,T^n)$,
its \coos are

\centerline{$\big((-1)^n(s_n-S_n),\ldots ,(s_2-S_2),
-(s_1-S_1),0\big)$.}
\item On the other hand, over the basis 

\centerline{$\big(1,(T-X_1),(T-X_1)(T-X_2),\ldots ,(T-
X_1)\cdots (T-X_{n})\big)$,}

its \coos are $(f_1,f_2,\ldots ,f_{n},0)$.
\end{itemize}
Consequently over the \ri $\AXn$,
each of the two vectors

\snic{\big((-1)^n(s_n-S_n),\ldots ,(s_2-S_2),-(s_1-S_1)\big)
\quad \hbox{ and }\quad (f_1,\ldots ,f_{n-1},f_n)\qquad}

are expressed in terms of the other by means of an unipotent matrix (triangular with $1$'s along the diagonal).
\end{enumerate}

}
\end{exercise}

\vspace{-1em}
\begin{exercise}\label{exoPrimePowerRoot} {(The \polz~$X^p - a$)}
{\rm
Let $a \in \Ati$ and $p$ be a prime number. Suppose that the \polz~$X^p - a$ has in $\gA[X]$ a nontrivial \mon divisor. Show that $a$ is a $p^{\rm th}$ power in $\gA$.
}
\end{exercise}

\vspace{-1em}
\begin{exercise}
\label{exoPrincipeIdentitesAlgebriques} 
(With the extension principle of 
\idasz)\\
{\rm
Let $S_n(\gA)$ be the submodule of $\Mn(\gA)$ consisting of the \smq matrices.  \\
For $A \in S_n(\gA)$, let $\varphi_{\!A}$ be the \endo of $S_n(\gA)$ defined by $S \mapsto \tra {A} S A$.  Compute $\det(\varphi_{\!A})$ in terms of $\det(A)$. Show that $\rC{\varphi_{\!A}}$ only depends on~$\rC {\!A}$.
}
\end{exercise}

\vspace{-1em}
\begin{exercise}
\label{exoFreeFracTransfert}
{\rm
Let $\gB\supseteq\gA$ be an integral \Alg which is a free \Amo of rank $n$, $\gK = \Frac(\gA)$ and~$\gL = \Frac(\gB)$. Show that every basis of $\gB/\gA$ is a basis of~$\gL/\gK$.
}
\end{exercise}

\vspace{-1em}
\begin{exercise}
\label{exoResultant}
{\rm
 Let $f\in\AX$, $g\in\AY$, $h\in\gA[X,Y]$. Prove that\perso{in le cours???}

\snic{\Res_Y\big(g,\Res_X(f,h)\big)=\Res_X\big(f,\Res_Y(g,h)\big).}

}
\end{exercise}

\vspace{-1em}
\begin{exercise}\label{exoSommesNewton}
{(\isN Newton sums and $\Tr(A^k)$)}
{\rm
Let $A \in \Mn(\gB)$ be a matrix. 
\\
Let $\rC{A}(X)=X^n+\sum_{j=1}^n (-1)^j s_j X^{n-j}$, $s_0=1$ and $p_k=\Tr(A^k)$.
 
\emph{1.}
Show that the $p_k$'s and $s_j$'s are linked by Newton's formulas for the sums of the $k^{\rm th}$ powers (Exercise~\ref{exoNewtonSum}): $\sum_{r = 1}^d (-1)^{r-1} p_{r} s_{d-r} = ds_d$
($d\in \lrbn$).
 
\emph{2.}
If $\Tr(A^k) = 0$ for $k \in \lrb {1..n}$, and if $n!$ is \ndz in
$\gB$, then $\rC{A}(X) = X^n$.

NB: this exercise can be considered as a variation on the theme of Proposition~\ref{prop2tschir}.
}
\end {exercise}

\vspace{-1em}
\begin{exercise}
\label{exoCorpsFiniEltPrimitif}
{\rm
Let $\gK \subseteq \gL$ be two finite fields, $q = \#\gK$ and $n = \dex{\gL:\gK}$. The sub\ri of $\gK$ generated by $1$ is a field $\FF_p$ where $p$ is a prime number, \hbox{and $q=p^r$} for an integer $r>0$.
\ix{Frobenius' automorphism} of (the $\gK$-extension)~$\gL$ is given by~$\sigma : \gL \to \gL,\,\sigma(x) = x^q$.
 
\emph{1.}
Let $R$ be the union of the roots in $\gL$ 
                      of the polynomials $X^{q^d} - X$ 
with $1 \le d < n$. \\
Show that $\#R < q^n$ and that for $x \in \gL \setminus R$, $\gL = \gK[x]$.

\emph{2.}
Here $\gK = \FF_2$ and $\gL = \FF_2[X]/\gen {\Phi_5(X)} = \FF_2[x]$ where
$\Phi_5(X)$ is the cyclotomic \pol $X^4 + X^3 + X^2 + X + 1$.
Check that $\gL$ is indeed a field; $x$ is a primitive \elt of $\gL$ over $\gK$ but it is not a \gtr of the multiplicative group~$\gL\eti$.

\emph{3.}
For $x \in \gL\eti$, let $o(x)$ be its order in the multiplicative group $\gL\eti$. \\
Show that $\gL = \gK[x]$ \ssi the order of $q$ in the group $\left(\aqo\ZZ{o(x)}\right)\eti$ is~$n$.
}
\end{exercise}

\vspace{-1em}
\begin{exercise}
\label{exoRacinesUniteCyclique}
{\rm  The aim of the exercise is to prove that in a \cdi the group of  $n^{\rm th}$ roots of unity is cyclic. 
Consequently the multiplicative group of a finite field is cyclic. We prove a result that is barely more \gnlz.
\\
 Show that in a nontrivial commutative \ri $\gA$, if \elts $(x_i)_{i\in\lrbn}$ form a group $G$ for the multiplication, and if $x_i-x_j$ is \ndz for every pair $i,j$ ($i\neq j$), then $G$ is cyclic.
\\
\emph{Hint}: by the structure \tho of  finite Abelian groups, a finite Abelian group, additively denoted, in which every \eqn $dx=0$ admits at most~$d$ solutions is cyclic.
Also use Exercise~\ref{exoLagrange}.
}
\end{exercise}

\vspace{-1em}
\begin{exercise}
\label{exoCorpsFinis} (Structure of finite fields, Frobenius' automorphism){\rm
  
\emph{1.}
Prove that two finite fields which have the same order are \isocz.

\emph{2.} If $\gF\supseteq\FF_p$ is a finite field of order $p^r$, prove that $\tau:x\mapsto x^p$ defines an \auto of $\gF$. This is called \emph{Frobenius' automorphism}. Show that the group of \autos of $\gF$ is a cyclic group of order $r$ generated by $\tau$.  

\emph{3.} In the previous case, $\gF$ is a Galois extension of $\FF_p$. Describe the Galois correspondence.
 
NB: We often denote by $\FF_q$ a finite field of order $q$, 
knowing that it is a slightly ambiguous notation if $q$ is not prime.

}
\end{exercise}

\vspace{-1em}
\begin{exercise}
\label{exo2CorpsFinis} (\Agq closure of $\FF_p$) 
{\rm

\emph{1.}  For each integer $r>0$ construct a field $\FF_{p^{r!}}$ of order $p^{r!}$. By proceeding by \recu we have an inclusion $\imath_r:\FF_{p^{r!}}\hookrightarrow \FF_{p^{(r+1)!}}$.
 
\emph{2.} 
Construct a field $\FF_{p^{\infty}}$ by taking the union of the $\FF_{p^{r!}}$ via the inclusions $\imath_r$.
Show that $\FF_{p^{\infty}}$ is an \agqt closed field that contains a (unique) copy of each finite field of \cara $p$.
}
\end {exercise}

\entrenous{ A caser ?

If  $\gF$ is un corps finite  
$\gF(T_1,\ldots,T_n)$ is pleinement factoriel au sens of Richman
}

\vspace{-1em}
\begin{exercise}\label{exoPpcmPolsSeparables} {(Lcm of separable \polsz)}
{\rm
  
\emph{1.} Let $x$, $x'$, $y$, $y'\in \gB$. Show that $\gen {x,x'} \gen {y,y'} \gen {x,y}^2 \subseteq \gen {xy, x'y + y'x}.$\\ 
Deduce that the product of two \spl and \com \polus in~$\gA[T]$ is a \spl \polz.

\emph{2.} 
If  $\gA$ is a \cdiz, the lcm of several \spl \pols is \splz.
}

\end{exercise}

\vspace{-1em}
\begin{exercise}
\label{exolemSousLibre} (Index of a \tf submodule in a free module)
{\rm  

\emph{1.}
Let $A\in\Ae{m\times n}$ and $E=\Im(A)\subseteq \Ae m$.
Show that $\cD_m(A)$ only depends on $E$. 
\\
We call this \id the \ixc{index}{of a submodule in a free module} \emph{of $E$ in $L=\Ae m$}, and we denote it by $\idg{L:E}_\gA$ \hbox{(or $\idg{L:E}$)}. Note that this index is null as soon as $E$ is not sufficiently close to $L$, for example if $n<m$.
\\
Check that in the case where $\gA=\ZZ$ we find the usual index of the subgroup of a group for two free Abelian groups of the same rank. 
 
 \emph{2.} If $E\subseteq F$ are \tf submodules of $L\simeq \Ae m$,
we have $\idg{L:E}\subseteq \box0{\idg{L:F}}$.
 
 \emph{3.} In addition, if $F$ is free and of rank $m$, we have the transitivity formula
$$\preskip.4em \postskip.4em 
\idg{L:E} = \idg{L:F}\,\idg{F:E}. 
$$

\emph{4.} If $\delta$ is a \ndz \elt of $\gA$, we have $\idg{\delta L:\delta E}=\idg{L:E}$. 
Deduce the \egt (\ref{eqlemSousLibre}) (\paref{eqlemSousLibre}) stated in Lemma~\ref{lemSousLibre}.
}
\end{exercise}

\vspace{-1em}
\begin{exercise}
\label{exoNbGtrsInverse} (Remark on Fact~\ref{factdefiiv})
{\rm   Let $\fa$ and~$\fb$ be two \ids in a \ri $\gA$ such that $\fa\,\fb=\gen{a}$
with $a$ \ndzz. Show that if $\fa$ is generated by $k$ \eltsz, we can find in $\fb$ a \sgr of $k$ \eltsz.
 
}
\end{exercise}

\vspace{-1em}
\begin{exercise}
\label{exoDecompIdeal} (Decomposition of an \id into a product of \iv \idemasz)
{\rm Consider a nontrivial integral \ri \emph{with explicit \dvez}\footnote{We say that an arbitrary ring is with explicit \dve if
 we have an \algo that tests, for $a$ and $b\in\gA$, if $\exists x,\, a=bx$, and in case of a positive answer, gives a suitable $x$.}%
\index{explicit divisibility!ring with ---}\index{ring!with explicit divisibility} $\gA$.  
\vspace{-4pt}
\begin{enumerate}
\item If $\fa$ is an \iv \id and if $\fb$ is a \itfz, prove that there is a test for~$\fb \subseteq \fa$.
\end{enumerate}
\vspace{-4pt}
Let $\fq_1$, $\ldots$, $\fq_n$ be \idemas (in the sense that the quotient \ris $\gA\sur{\fq_k}$ are nontrivial \cdisz), $\fb$ be a \itf and $a$ be a \ndz \elt of $\gA$ satisfying $a\gA=\fq_1 \cdots \fq_n \subseteq \fb$.
\vspace{-4pt}
\begin{enumerate} \setcounter{enumi}{1}
\item  Show that the $\fq_i$'s are \iv and $\fb$ is the product of some of the $\fq_i$'s (and thus it is \ivz). Furthermore, this \dcn of $\fb$ into a product of \tf \idemas is unique up to order of the factors.
\end{enumerate}

}
\end{exercise}
\vspace{-1em}
\begin{exercise}\label{exoSymboleLegendre}
{(Legendre symbol)}\index{Legendre symbol}\\
{\rm  
Let $\gk$ be a finite field of odd cardinality $q$; we define the \emph{Legendre symbol}

\snic{\dsp
{ \legendr \bullet \gk} :
\gk\eti \lora \{\pm 1\}, \;
x \;\longmapsto \;\formule{
\phantom-1 \hbox{ if $x$ is a square in } \gk\eti ,
\\  
-1 \hbox{ otherwise.}
}
}

Show that ${. \legendre \gk}$ is a group morphism and that ${x \legendre \gk} = x^{q-1 \over 2}$.\\
In particular,~$-1$ is a square in $\gk\eti$
\ssi $q \equiv 1 \bmod 4$.

NB: if $p$ is an odd prime number and $x$ is an integer \com to $p$ we find 
Legendre's symbol ${x \legendre p} $ in the form ${{\ov x} \legendre \FFp}$.
}
\end {exercise}

\vspace{-1em}
\begin{exercise}\label{exoRabinovitchTrick}
 {(Rabinovitch's trick)}\\
{\rm  
Let $\fa\subseteq \gA$ be an \id and $x \in \gA$.
Consider the following \id of $\AT$:

\snic {
\fb = \gen {\fa, 1-xT} = \fa[T] + \gen {1-xT}_{\gA[T]}.
}

Show the \eqvc $x \in \sqrt\fa \iff 1 \in \fb$.

}
\end {exercise}

\vspace{-1em}
\begin{exercise}
\label{exoDunford} (Jordan-Chevalley-Dunford decomposition)\index{Dunford!Jordan-Chevalley- --- decomposition}
\\
 {\rm  
Let $M\in\Mn(\gA)$. Suppose that the \polcar of $M$ divides a power of a \spl \pol $f$. 

\emph{1.} Show that there exist $D$, $N\in\Mn(\gA)$ such that:
\begin{itemize}
\item  $D$ and $N$ are \pols in $M$ (with \coes in $\gA$).
\item  $M=D+N$.
\item  $f(D)=0$.
\item  $N$ is nilpotent.
\end{itemize}
%
\emph{2.} 
Prove the uniqueness of the above decomposition, including by weakening the first constraint, by only requiring that $DN=ND$.

}
\end{exercise}

\vspace{-1em}
\begin{exercise}\label{exoDunfordBis}
{(\Spby integral \eltsz)}
\\{\rm  
Let $\gA \subseteq \gB$. We say that $z \in \gB$ is \emph{\spby integral} over $\gA$ if $z$ is a root of a \spl \polu of $\AT$.
Here we are looking for an example for which the sum of two \spby integral \elts is a nonzero nilpotent and non\spby integral \eltz.
\\
Let $\gB = \gA[x] = \aqo{\gA[X]}{X^2 + bX + c}$. Suppose that $\Delta = b^2 - 4c$ is a unit of $\gA$. 
For $a\in\gA$, compute the \polcar of $ax$ over $\gA$ and its \discriz.  
Deduce an example as stated when $\DA(0)\neq 0$.

}

\end{exercise}

\vspace{-1em}
\begin{problem}\label{exoDiscriminantsUtiles} {(Some useful resultants and discriminants)}\\
{\rm
\emph {1.}
Show that $\disc(X^n + c) = (-1)^{n(n-1) \over 2} n^n c^{n-1}$.
More \gnltz, prove for~$n \ge 2$ the \egt

\snic {
\disc(X^n + bX + c) = (-1)^{n(n-1) \over 2}
\bigl(n^n c^{n-1} + (1 - n)^{n-1} b^n \bigr).
}

\emph {2.}
For $n, m \in \NN^*$, by letting $d = \pgcd(n,m)$, $n_1=\frac n d$ and $m_1=\frac m d$ prove the \egt
$$\preskip-.2em \postskip.0em 
\Res(X^n - a, X^m - b) = (-1)^n (b^{n_1} - a^{m_1})^d.
$$
More \gnlt
$$\preskip.1em \postskip.2em 
\Res(\alpha X^n - a,n, \beta X^m - b,m) =
(-1)^n (\alpha^{m_1}b^{n_1} - \beta^{n_1}a^{m_1})^d. 
$$
\emph {3.}
Notations as in item \emph{2}, with $1 \le m \le n-1$.  Then 
prove

\snic {\mathrigid 2mu
\disc(X^n + bX^m + c) = (-1)^{n(n-1) \over 2} c^{m-1}
\bigl(n^{n_1} c^{n_1 - m_1} - 
(n - m)^{n_1 - m_1} m^{m_1} (-b)^{n_1}\bigr)^d.
}

\emph {4.}
For $n \in \NN^*$, denote by $\Phi_n$ the cyclotomic \pol of level $n$
(see \Pbmz~\ref{exoPolCyclotomique}).
Then, for prime $p$ $\ge 3$\index{polynomial!cyclotomic ---}
prove
$$\preskip.0em \postskip.2em 
\disc(\Phi_p) = (-1)^{p-1 \over 2} p^{p-2}. 
$$
\emph {5.}
Let $p$ be prime and $k\ge 1$. Then 
prove that $\Phi_{p^k}(X) = \Phi_p(X^{p^{k-1}})$ and
$$\preskip.1em \postskip.2em 
\disc(\Phi_{p^k}) = (-1)^{\varphi(p^k) \over 2} p^{(k(p-1)-1) p^{k-1}}
\qquad (p,k) \ne (2,1), 
$$
with for $p \ne 2$, $(-1)^{\varphi(p^k)\over 2} = (-1)^{p-1\over 2}$. 
For $p=2$, prove that we obtain $\disc(\Phi_4) = -4$ 
and~$\disc(\Phi_{2^k}) = 2^{(k-1) 2^{k-1}}$ for $k \ge 3$.
In addition, prove that $\disc(\Phi_2) = 1$.

\emph {6.}
Let $n \ge 1$ and $\zeta_n$ be an $n^{\rm th}$ primitive root of the unit.
\\
If $n$ is not the power of a prime number, then prove that $\Phi_n(1) = 1$,
and $1 - \zeta_n$ is \iv in $\ZZ[\zeta_n]$. 
\\
If $n = p^k$ with $p$ prime, $k \ge 1$, then prove that $\Phi_n(1) = p$. Finally, prove that $\Phi_1(1) = 0$.

\emph {7.}
Let $\Delta_n = \disc(\Phi_n)$.  For coprime $n$, $m$, prove that we have the multiplicativity formula $\Delta_{nm} = \Delta_n^{\varphi(m)} \Delta_m^{\varphi(n)}$ and the \egt
$$\preskip.4em \postskip.4em\ndsp 
\Delta_n = (-1)^{\varphi(n) \over 2}
\frac{n^{\varphi(n)}} { \prod_{p \mid n} p^{\varphi(n) \over p-1}} \qquad
\hbox {for }n \ge 3. 
$$
}
\end{problem}

\vspace{-2em}
\begin{problem}
\label{exoAnneauEuclidien} (Euclidean \risz, the $\ZZ[i]$ example)
\\
{\rm  A \emph{Euclidean stathm} is a map $\varphi:\gA\to\NN$ that satisfies the following \prtsz\footnote{In the literature we sometimes find a \gui{Euclidean stathm} defined as a map $\varphi:\gA\to\NN\cup\so{-\infty}$, or $\varphi:\gA\to\NN\cup\so{-1}$
(the minimum value being always equal to $\varphi(0)$).} 
(roughly speaking, we copy the Euclidean division in $\NN$)
\begin{itemize}
\item $\varphi(a)=0\iff a=0$.
\item $\forall a,b\neq0,\;\exists q,r,\;\;\;a=bq+r \hbox{ and } \varphi(r)<\varphi(b)$. 
%
\end{itemize}
 A \emph{Euclidean \riz} is a nontrivial integral \ri given with a Euclidean stathm. Note that the \ri is discrete. We can then do with the \gui{division} given by the stathm the same thing we do in $\ZZ$ with  Euclidean division.\index{ring!Euclidean ---}%
\index{Euclidean!ring}
\\
The most renowed examples are the following.
\begin{itemize}
\item $\ZZ$, with $\varphi(x)=\abs{x}$,
\item  $\KX$ ($\gK$ a \cdiz), with  $\varphi(P)=1+\deg(P)$ for $P\neq0$,
\item $\ZZ[i]\simeq\aqo\ZZX{X^2+1}$, with  $\varphi(m+in)=m^2+n^2$, 
\item $\ZZ[i\sqrt 2]\simeq\aqo\ZZX{X^2+2}$, with  $\varphi(m+i\sqrt 2 n)=m^2+2n^2$. 
\end{itemize}
In addition, in these examples we have the \eqvc $\;x\in\Ati \iff \varphi(x)=1$.

\begin{enumerate}
\item \label{i1exoAnneauEuclidien} \emph{(Extended Euclidean algorithm)}
For all $a,b$, there exist $u,v,a_1,b_1,g$ such that 
$$ \cmatrix{g\cr 0}=\bloc{u}{v}{-b_1}{a_1} \cmatrix{a\cr b} \;\hbox{ and } \; ua_1+vb_1=1
.$$
In particular, $\gen{a,b}=\gen{g}$ and $g$ is a gcd of $a$ and $b$. 
If $(a,b)\neq(0,0)$, $\frac{ab}g$ is a lcm of $a$ and $b$.
\item \label{i2exoAnneauEuclidien}
\begin{enumerate}
\item Show that the \ri $\gA$ is principal.
\item  Let us make the following assumptions.
 \begin{itemize}
\item $\Ati$ is a detachable subset of $\gA$.
\item We have a primality test at our disposal for the \elts of $\gA\setminus\Ati$ in the following sense: given $a\in \gA\setminus\Ati$ we know how to decide if $a$ is \irdz, and in case of a negative response, write $a$ in the form $bc$ \hbox{with $b$, $c\in \gA\setminus\Ati$}. 
\end{itemize}
Show then that $\gA$ satisfies the \gui{fundamental \tho of \ariz}
(unique \dcn into prime factors, up to association).
\end{enumerate} 
\end{enumerate}

\vspace{-.5em}
\emph{The $\ZZ[i]$ example.}
 Recall that $z = m+in \mapsto \ov z = m-in$ is an \auto \hbox{of $\ZZ[i]$} and that the norm $\rN=\rN_{\ZZ[i]/\ZZ}$ ($\rN(z)=z\ov z$) is a Euclidean stathm. Take an \elt of $\ZZ[i]$ close to $a/b\in\QQ[i]$ for the above $q$ and check \hbox{that $\rN(r)\leq\rN(b)/2$}.\\
To know which are the \irds \elts of $\ZZ[i]$, it suffices 
to know how to decompose in $\ZZ[i]$ each prime number $p$ of $\NN$.
\\  
This amounts to determining the \ids containing $p\ZZ[i]$, \cad the \ids of $\gZ_p:=\aqo{\ZZ[i]}p$. But $\gZ_p\simeq \aqo{\FFp[X]}{X^2+1}$. We are thus reduced to finding the divisors of $X^2+1$, therefore to factorizing $X^2+1$, \hbox{in $\FFp[X]$}.  
\begin{enumerate}\setcounter{enumi}{2}
\item \label{i3exoAnneauEuclidien} 
Show that a priori three cases can arise.
\begin{itemize}
\item $X^2+1$ is \ird in $\FFp[X]$, and $p$ is \ird in $\ZZ[i]$.
\item $X^2+1=(X+u)(X-u)$ in $\FFp[X]$ with $u\neq -u$, and then 

\snic{ \gen{p}=\gen{i+u,p}\gen{i-u,p}=\gen{m+in}\gen{m-in}$
 and $p=m^2+n^2.}

\vspace{2pt}
\item \label{i4exoAnneauEuclidien}
$X^2+1=(X+u)^2$ in $\FFp[X]$, and then $\gen{p}=\gen{i+u}^2$. This only happens for $p=2$, with   $2=(-i)(1+i)^2$ (where $-i\in\ZZ[i]\eti$). 
\end{itemize}

\item \label{i5exoAnneauEuclidien}
If $p\equiv 3\mod 4$, then $-1$ is not a square in $\FFp$.
If $p\equiv 1\mod 4$, then~$-1$ is a square in~$\FFp$. In this case give an efficient \algo
to write $p$ in the form $m^2+n^2$ in $\NN$.
\item \label{i6exoAnneauEuclidien} Let $z\in\ZZ[i]$. We can write $z=m(n+qi)$ with $m,n,q\in\NN$ $\pgcd(n,q)=1$. Give an efficient \algo to decompose $z$ into prime factors in $\ZZ[i]$ knowing a \dcn into prime factors of $\rN(z)=m^2(n^2+q^2)$ in $\NN$. 
\\
Given a \dcn into prime factors of $s\in\NN$, describe under which condition $s$ is a sum of two squares, as well as the number of expressions~\hbox{$s=a^2+b^2$} with $0<a\leq b$ in $\NN$.

\item \label{i7exoAnneauEuclidien} Say in which (relatively rare) cases we can generalize the previous procedure to decompose into prime factors the \itfs of a \ri $\ZZ[\alpha]$, when $\alpha$ is an \agq integer.
\end{enumerate}
 
}
\end{problem}

\vspace{-1em}
\begin{problem}
\label{exoPetitKummer} (Kummer's little \thoz)
\\
{\rm \Pbmz~\ref{exoAnneauEuclidien} can be generalized for \ris of principal integers of the form $\ZZ[\alpha]$, but this case is relatively rare. On the contrary, Kummer's little \tho gives the \dcn of a prime number (in $\NN$) into products of $2$-generated \idemas for almost all the prime numbers, in all the \ris of integers. This shows the intrinsic superiority of the \gui{\id numbers} introduced by Kummer. Furthermore, the argument is extremely simple and only requires the Chinese remainder  \thoz.
However, the prime numbers that do not fall under the scope of Kummer's little \tho constitute in fact the heart of \agq number theory. Those are the ones that required a fine tuning of the theory (according to two distinct methods due to Kronecker and Dedekind), without which all decisive progress would not have been possible.

 Consider a zero $\alpha$ of an \ird \poluz~\hbox{$f(T)\in\ZZ[T]$},
such {that $\ZZ[\alpha]\simeq\aqo{\ZZ[T]}{f(T)}$}. Let $\Delta=\disc(f)$.
\begin{enumerate}
\item Let $p$ be a prime number which does not divide $\Delta$.
\begin{itemize}
\item  Show that $f(T)$ is \spl in $\FFp[T]$. 
\item  Decompose $f(T)$ in $\FFp[T]$ in the form $\prod_{k=1}^{\ell}Q_k(T)$ with  distinct \mon \ird $Q_k$'s. Let $q_k=Q_k(\alpha)$ (in fact it is only defined modulo $p$, but we can lift $Q_k$ in $\ZZ[T]$).
Show that in $\ZZ[\alpha]$ we have $\gen{p}=\prod_{k=1}^{\ell}\gen{p,q_k}$ and that the \ids $\gen{p,q_k}$ are maximal, distinct and invertible.
In particular, if $\ell=1$, $\gen{p}$ is maximal. 
\item  Show that this \dcn remains valid in every \ri $\gA$ such that $\ZZ[\alpha]\subseteq\gA\subseteq\gZ$, where $\gZ$ is the \ri of integers of $\QQ[\alpha]$.
\end{itemize}
\item Let $a\in\ZZ[\alpha]$ such that $A=\rN_{\ZZ[\alpha]/\ZZ}(a)$ is \com to $\Delta$. Let $\fa=\gen{b_1,\ldots,b_r}$ be a \itf of $\ZZ[\alpha]$ containing $a$. Show that in $\ZZ[\alpha]$ the \idz~$\fa$ is \iv and can be decomposed into products of \idemas that divide the prime factors of $A$. 
Finally, this \dcn is unique up to order of the factors and all of this remains valid in every \ri $\gA$ as above.
\end{enumerate} 
}
\end{problem}

\vspace{-1em}
\begin{problem}\label{exoPolCyclotomique}  {(The cyclotomic \pol $\Phi_n$)}\index{polynomial!cyclotomic ---}\\
{\rm
In $\gA[X]$, the \polz~$X^n - 1$ is \spl \ssi $n \in \Ati$
.\\
Let $\gQ_n$ be a \cdr over $\QQ$ for this \polz.
Let $\UU_n$ be the group of $n^{\rm th}$ roots of the unit in $\gQ_n$.
It is a cyclic group of order~$n$, which therefore has $\varphi(n)$ \gtrs ($n^{\rm th}$ primitive roots of the unit). We define $\Phi_n(X) \in \gQ_n[X]$ by $\Phi_n(X) = \prod_{o(\xi) = n} (X - \xi)$. It is a \polu of degree $\varphi(n)$. We have the fundamental \egt

\snic {
X^n - 1 = \prod_{d \divi n} \Phi_d(X),
}

which allows us to prove by \recu on $n$ that $\Phi_n(X) \in \ZZ[X]$.

\begin {enumerate}
\item
Following the steps below, prove that
$\Phi_n(X)$ is \ird in $\ZZ[X]$ (therefore in $\QQ[X]$, Proposition~\ref{propZXfactor}).
Let $f$, $g$ be two \polus of $\ZZ[X]$ with $\Phi_n = fg$ and $\deg f \ge 1$; you must prove that $g = 1$.

\begin {itemize}
\item [a.]
It suffices to prove that $f(\xi^p) = 0$ for every prime $p \nedivi n$ and for every zero $\xi$ of $f$ in $\gQ_n$.
\item [b.]
Suppose that $g(\xi^p) = 0$ for some zero $\xi$ of $f$ in $\gQ_n$.
Examine what happens in $\FF_p[X]$ and conclude the result.
\end {itemize}

\item
Let us fix a root $\xi_n$ of $\Phi_n$ in $\gQ_n$.
\\
Show that $\gQ_n = \QQ(\xi_n)$ and that with $(\QQ,\gQ_n,\Phi_n)$, we are in the \elr Galois situation of Lemma~\ref{lemGaloiselr}.
\\
 Describe the explicit \isos of the groups

\snic {
\Aut(\UU_n)  \simeq (\ZZ/n\ZZ)^{\!\times} \simeq \Gal(\gQ_n/\QQ).
}

\item
Let $\gK$ be a field of \cara $0$. What can be said of a \cdrz~$\gL$ of $X^n-1$ over $\gK$?

\end {enumerate}
}
\end {problem}

\vspace{-1.1em}
\pagebreak

\begin{problem}\label{exoCyclotomicRing}
(The \ri $\ZZ[\!\root n \of 1]$: \ddpz, \fcn of \idsz)
\\
 {\rm  
Let $\Phi_n(X) \in \ZZ[X]$ be the cyclotomic \pol of order $n$, \ird over $\QQ$. 
Let $\gQ_n  = \QQ(\zeta_n)\simeq \aqo {\QQ[X]} {\Phi_n}$. The multiplicative group $\UU_n$ generated by $\zeta_n$ ($n^{\rm th}$ primitive root of the unit) is cyclic of order~$n$.
\\
Among other things we will prove 
that the \ri $\gA = \ZZ[\UU_n] = \ZZ[\zeta_n]\simeq\aqo{\ZZ[X]}{\Phi_n}$ is a \emph{Pr\"ufer domain}\index{Pru@Pr\"ufer!domain}\index{domain!Pru@Pr\"ufer ---}: an integral \ri whose nonzero \itfs are \iv (cf. Section~\ref{secIplatTf} and Chapiter~\ref{ChapAdpc}).

 \emph{1.}
Let $p \in \NN$ be a prime number. The steps below show that $\sqrt {p\gA}$ is a
\idp and express it as a finite product of 2-generated \iv \idemasz. Consider the distinct \ird factors of $\Phi_n$ modulo $p$ that we lift to the \polus $f_1, \ldots, f_k \in \ZZ[X]$. Let $g = f_1 \cdots f_k$ (such that $\ov g$ is the subset without a square factor of $\Phi_n$ modulo $p$) and $\fp_i = \gen {p, f_i(\zeta_n)}$ for $i \in \lrbk$.
\begin{enumerate}\itemsep0pt
\item [\emph{a.}]
Show that $\fp_i$ is a \idema and that

\snic {
\sqrt {p\gA} = \gen {p, g(\zeta_n)} = \fp_1 \ldots \fp_k
}

\item [\emph{b.}]
If $p$ does not divide $n$, prove that~$\ov g = \ov {\Phi_n}$, thus $\sqrt
{p\gA} = \gen {p}$ is a \idpz.

\item [\emph{c.}]
Suppose that $p$ divides $n$ and write $n = mp^k$ with $k \ge 1$,
$\pgcd(m, p) = 1$.\\  
By studying the \fcn of $\Phi_n$ modulo $p$, prove that~$\ov g = \ov {\Phi_m}$. Deduce that~$\sqrt {p\gA} = \gen {p, \Phi_m(\zeta_n)}$. Then prove that $p \in \gen {\Phi_m(\zeta_n)}$, and therefore that $\sqrt {p\gA} = \gen {\Phi_m(\zeta_n)}$ is a \idpz.

\item [\emph{d.}]
Deduce that $p\gA$ is a product of the form $\fp_1^{e_1} \ldots \fp_k^{e_k}$.
\end{enumerate}
 \emph{2.}
Let $a \in \ZZ \setminus \{0\}$; prove that $a\gA$ is a product of \iv \idemas with two \gtrsz. Deduce that in $\gA$ every nonzero \itf can be decomposed into a product of $2$-generated \iv \idemas and that the \dcn is unique up to factor order.
}

\end{problem}

\vspace{-1em}
\begin{problem}\label{exoSommeGauss}
{(An \elr \prt of Gauss sums)}
\\ 
{\rm  
Let $\gk$ be a finite field of cardinality $q$ and $\gA$ be an integral \riz. Consider
\begin{itemize}
\item a \gui{multiplicative character} $\chi : \gk\eti \to \Ati$,  \cad a morphism of multiplicative groups,
\item an \gui{additive character}  $\psi : \gk \to \Ati$,  \cad a morphism of groups 
$$\preskip.3em \postskip.3em 
\psi : (\gk, +) \to (\Ati, \times). 
$$
\end{itemize}
  Suppose that neither $\chi$ nor $\psi$ are trivial and that $\chi$ is extended to the whole of $\gk$ \hbox{via $\chi(0) = 0$}. Finally, the Gauss sum of $\chi$ is defined, with respect to $\psi$, by

\snic {
G_\psi(\chi) = \sum_{x \in \gk} \chi(x) \psi(x) = 
\sum_{x \in \gk\eti} \chi(x) \psi(x).
}

We aim to prove that

\snic {
G_\psi(\chi)G_\psi(\chi^{-1}) = q\chi(-1) 
,}


and give \aris applications of this result (Question \emph{4}).

\emph{1.}
Let $G$ be a finite group and $\varphi : G \to \Ati$ be a nontrivial \homoz.  Show that $\sum_{x \in G} \varphi(x) = 0$. 

\emph{2.}
Show that
$$\preskip.2em \postskip.4em\ndsp 
\sum_{x + y = z} \chi(x) \chi^{-1}(y) = \cases {
-\chi(-1) &if $z \ne 0$,\cr
(q-1) \chi(-1) &otherwise.\cr} 
$$
\emph{3.}
Deduce that $G_\psi(\chi)G_\psi(\chi^{-1}) = q\chi(-1)$.

\emph{4.}
Consider $\gk = \FF_p$ where $p$ is an odd prime number, $\gA =\QQ(\!\root p \of 1)$, and $\zeta$  a $p^{\rm th}$ primitive root of the unit in $\gA$. The characters $\psi$ and $\chi$ are defined by

\snic {
\psi(i \bmod p) = \zeta^i, \qquad  \chi(i \bmod p) = {i \legendre p}
\quad \hbox {(Legendre symbol).}
}

\begin{enumerate}
\item [\emph{a.}]
Then, $\chi = \chi^{-1}$, the Gauss sums $G_\psi(\chi)$, 
$G_\psi(\chi^{-1})$ are equal to
$$
\preskip.3em \postskip.1em\ndsp 
\tau \eqdefi \sum_{i \in \FF_p^*} {i \legendre p} \zeta^i, 
$$
and by letting $p^* = (-1)^{p-1 \over 2} p$ (such that $p\etl\equiv1\mod 4$), we obtain

\snic {
\tau^2 = p^*,  \quad \hbox {in particular,} \quad
\QQ(\sqrt {p^*}) \subseteq \QQ(\!\root p \of 1).
}

\item [\emph {b.}]
Define $\tau_0 = \sum_{i \in \FF_p^{\times 2}} \zeta^i, \, \tau_1 = \sum_{i \in \FF_p\eti \setminus \FF_p^{\times 2}} \zeta^i$ such that $\tau  = \tau_0 - \tau_1.$
Show that $\tau_0$ and $\tau_1$ are the roots of $X^2 + X + {1 - p^* \over 4}$ and that the \riz~\hbox{$\ZZ[\tau_0] = \ZZ[\tau_1]$} is the \ri of integers of $\QQ(\sqrt {p^*})$.

\end{enumerate}
}
\end {problem}

\vspace{-1em}
\begin{problem} 
\label{exoDedekindPolynomial}\index{Dedekind!polynomial} 
{(The Dedekind \pol $f(X) = X^3 + X^2 -2X + 8$)}
\\ 
{\rm
The aim of this problem is to provide an example of a \ri $\gA$ of integers of a number field which is not a monogenic $\ZZ$-\algz.\footnote{An \Alg $\gB$ is said to be \emph{monogenic} when it is generated, as an \Algz, by a unique element $x$. So $\gB=\gA_1[x]$ where $\gA_1$ is the image of $\gA$ in $\gB$.}

\begin{enumerate}\itemsep0pt\mou
\item
Show that $f$ is \ird in $\ZZ[X]$ and that $\disc(f) = -2\,012 = -2^2
\times 503$.

\item
Let $\alpha$ be a root of $f(X)$. Show that $\beta = 4\alpha^{-1}$ is integral over~$\ZZ$, that  
$$
\preskip.3em \postskip.3em 
\gA = \ZZ \oplus \ZZ\alpha \oplus \ZZ\beta 
$$
is the \ri of integers of $\QQ(\alpha)$ and that $\Disc_{\gA/\ZZ} = -503$.

\item
Show that the prime number $p = 2$ is completely decomposed in $\gA$, in other words that $\gA/2\gA \simeq \FF_2 \times \FF_2 \times \FF_2$.
Deduce that $\gA$ is not a monogenic \ZZlgz.

\item
\emph{(Avoiding the conductor, Dedekind)} Let $\gB \subseteq \gB'$ be two \risz, $\ff$ be an \id of $\gB$ satisfying $\ff\gB' \subseteq \gB$; a fortiori $\ff\gB' \subseteq\gB'$ and $\ff$ is also an \id of $\gB'$. Then, for every \id $\fb$ of $\gB$ such that $1 \in \fb+\ff$, by letting $\fb' = \fb\gB'$, the canonical morphism $\gB/\fb \to \gB'/\fb'$ is an \isoz.

\item
Deduce that $2$ is an \emph{essential divisor} of $\gA$; by that we mean that $2$ divides  the index $\idg{\gA:\ZZ[x]}$ for any primitive \eltz~$x$ of~$\QQ(\alpha)/\QQ$ integral over~$\ZZ$.
\end{enumerate}
}
\end{problem}

\vspace{-1em}
\begin{problem}\label{exoGaloisNormIdeal} {(Norm of an \id in quasi-Galoisian context)}
 \\
{\rm
Let $(\gB, \gA, G)$ where $G \subseteq \Aut(\gB)$ is a finite group, and $\gA = \gB^G=\Fix_\gB(G)$. If~$\fb$ is an \id of $\gB$, let $\rN'_G(\fb) = \prod_{\sigma \in G} \sigma(\fb)$ (\id of $\gB$) and $\rN_G(\fb) = \gA \cap \rN'_G(\fb)$ (\id of $\gA$).
\begin{enumerate}\itemsep0pt\mou
\item 
Show that $\gB$ is integral over $\gA$.
\item 
Let $\gB = \ZZ[\sqrt d]$ where $d \in \ZZ$ is not a square, $\tau$ be the \auto (also denoted by~\hbox{$z \mapsto \ov z$}) defined by $\sqrt d \mapsto -\sqrt d$, and $G = \gen {\tau}$. Therefore $\gA = \ZZ$. Suppose that~\hbox{$d \equiv 1 \bmod 4$} and let $\fm = \gen {1+\sqrt d, 1-\sqrt d}$.
\begin{enumerate}
\item [\emph{a.}]
We have $\fm = \ov\fm$, $\rN'_G(\fm) = \fm^2 = 2\fm$ and $\rN_G(\fm) = 2\ZZ$.
Deduce that $\fm$ is not \iv and that we do not have $\rN'_G(\fm) = \rN_G(\fm)\gB$.
\item [\emph{b.}]
Show that $\ZZ[\sqrt d]/\fm \simeq \FF_2$; thus $\fm$ is of index $2$ in
$\ZZ[\sqrt d]$ but $2$ is not the gcd of the $\rN_G(z)$'s, $z \in \fm$.
Also check that  
$\fb \mapsto \idg{\gB:\fb}$
is not multiplicative over the nonzero \ids of $\gB$.
\end{enumerate}

\item 
Suppose that $\gB$ is \icl and that $\gA$ is a Bézout domain. 
\\
Let $\fb \subseteq \gB$ be a \itfz.
\begin{enumerate}
\item [\emph{a.}]
Give a $d \in \gA$ such that $\rN'_G(\fb) = d\gB$. In particular, if $\fb$ is nonzero, it is \ivz.
Thus, $\gB$ is a \ddpz.
\item [\emph{b.}]
Show that $\rN_G(\fb) = d\gA$, therefore $\rN'_G(\fb) = \rN_G(\fb)\gB$.
\item [\emph{c.}]
Suppose that
$\gB/\fb$ is \isoc as an \Amo to $\aqo{\gA}{a_1} \times \cdots \times \aqo{\gA}{a_k}$.
\\
Show that $\rN_G(\fb) = \gen {a_1 \cdots a_k}_\gA$.

\item [\emph{d.}]
Suppose $\#G = 2$.  Express, in terms of a finite \sgr of $\fb$, \elts $z_1$, \dots, $z_m \in \fb$ such that $\rN_G(\fb) = \gen {\rN(z_i), i \in \lrb {1..m}}_\gA$.
\end{enumerate}
\end{enumerate}

}
\end{problem}

\vspace{-1em}
\begin{problem}\label{exoLemmeFourchette} 
 {(Forking lemma)}\index{Forking lemma}
{\rm
\begin{enumerate}
\item Let $\gA$ be an \acl 
with quotient field $\gk$, 
$\gK$ be a \spl finite extension of $\gk$ of degree $n$, $\gB$ be the integral closure of $\gA$ in $\gK$.
Show that there exists a basis $(\ue) = (e_1, \ldots, e_n)$ of $\gL/\gK$ contained in $\gB$. Let~\hbox{$\Delta = \disc(\ue)$} and $(\ue') = (e'_1, \ldots, e'_n)$ be the trace-dual basis of $(\ue)$.
Show the inclusions
$$
\preskip.3em \postskip.0em\ndsp 
\bigoplus_{i=1}^n \gA e_i \;\subseteq\; \gB \;\subseteq\;
\bigoplus_{i=1}^n \gA e'_i \;\subseteq\;
\Delta^{-1} \bigoplus_{i=1}^n \gA e_i. 
$$
\end{enumerate}

In the following $\gA=\ZZ$ and $\gk=\QQ$; $\gK$ is thus a number field
and $\gB=\gZ$ is its \ri of integers.   
Consider some $x\in\gZ$ such that $\gK=\QQ[x]$. 
\\
Let $f(X)=\Mip_{\QQ,x}(X)\in\ZZ[X]$ and $\delta^2$ be the greatest square factor of $\disc_X( f)$. 
\\
By Proposition~\ref{propAECDN}, $\gZ$ is a free $\ZZ$-module of rank $n = \dex{\gL : \QQ}$, and we \hbox{have $\ZZx\subseteq\gZ\subseteq \fraC 1 \delta \ZZx$}.
  This is slightly more precise than the result from item~\emph{1.}
\\
Consider a \tf \ZZlg $\gB$ intermediate between $\ZZx$ and $\gZ$.
As it is a \tf \ZZmoz, $\gB$ is \egmt a free \ZZmo of rank $n$.
The most important case is that where $\gB=\gZ$.
\\ 
The aim of the \pb is to find a $\ZZ$-basis of $\gB$ of the form
$$\preskip.2em \postskip.2em\ndsp
\cB=\big(\fraC{g_0}{d_0},\fraC{g_1(x)}{d_1},\fraC{g_2(x)}{d_2},\dots,\fraC{g_{n-1}(x)}{d_{n-1}}\big)  
$$   
with $g_k\in\ZZ[X]$  of degree  $k$  for all $k$, and each $d_k>0$ ass mall as  possible. 
Establish this result 
with \polus $g_k$ and $1=d_0\mid d_1\mid d_2\mid\cdots\mid d_{n-1}$.
\\
The field $\gK$ is a \Qev with basis $(1,x,\dots,x^{n-1})$ and {for $ k \in \lrb{0..n-1}$}, let $\pi_k : \gK \to \QQ$ be the linear component form \hbox{over $x^k$} and
$$\preskip.2em \postskip.2em\ndsp
Q_k = \bigoplus\nolimits_{i=0}^k \QQ\, x^i,\; Z_k = \frac1\delta \bigoplus\nolimits_{i=0}^k \ZZ\, x^i, \quad \hbox{and }{  F_k = Q_k \cap \gB=  Z_k \cap \gB}.
$$
It is clear that $Q_0 = \QQ$, $Q_{n-1} = \gK$, $F_0 = \ZZ$ and $F_{n-1} = \gB$.
\begin{enumerate}\setcounter{enumi}{1}\itemsep0pt
\item Show that the \ZZmo $F_k$ is free and of rank $k+1$.  
\\ The \ZZmo $\pi_k(F_k)$ is a \tf \ZZsmo of $\fraC1\delta\ZZ$.
Show that it is of the form $\fraC1{d_k}\ZZ$ for some $d_k$ that divides $\delta$. NB: $d_0=1$. 
\item Let $y_k$ be an \elt of $F_k$ such that $\pi_k(y_k)=\fraC1{d_k}$.
\\
Write $y_k$ in the form $f_k(x)/d_k$, with $f_k\in\QQ[X]$ \mon and of degree~$k$. Clearly $y_0=1$. However, the other $y_i$'s are not uniquely determined. 
Show that $(1,y_1,\dots,y_k)$ is a $\ZZ$-basis of $F_k$.
\item Show that if $i+j\leq n-1$, we have $d_id_j\mid d_{i+j}$. 
In particular $d_i$ divides $d_k$ \hbox{if $1\leq i<k\leq n-1$}.
Also deduce \hbox{that $d_1^{n(n-1)/2}$} divides~$\delta$. 
\item Show that $d_ky_k\in\ZZ[x]$ for each $k \in \lrb{0..n-1}$.  
Deduce \hbox{that $f_k\in\ZZ[X]$} and that $\big( 1,  f_1(x),\dots,  f_{n-1}(x)\big)$
is a $\ZZ$-basis of $\ZZx$.  
\item Show that   
$\cB={\big( 1, \fraC 1 {d_1} f_1(x),\dots, \fraC 1 {d_{n-1}} f_{n-1}(x)\big)}
$ 
is a $\ZZ$-basis of $\gB$ adapted to the inclusion $\ZZx\subseteq \gB$. 
The $d_i$'s are therefore the invariant factors of this inclusion, and $\prod_{i=1}^{n-1}d_i$ is equal to  
the index $\idG{\gB:\ZZx}$ that divides $\delta$. 
\end{enumerate}

}
\end{problem}

\vspace{-1em}
\begin{problem}
\label{exoPolynomialAutomorphism} 
  {(Changing variables, \pol automorphisms and \imN Newton's method)}
    \\
{\rm
Let $F = (F_1, \ldots, F_n)$ with $F_i \in \gA[X] = \gA[X_1, \ldots,
X_n]$ and $\theta_F : \gA[X] \to \gA[X]$ be the morphism of $\gA$-\algs performing $X_i \mapsto F_i$; we therefore have $\theta_F(g) = g(F)$. Assume that $\gA[X] = \gA[F]$: there thus exists a $G_i \in \gA[X]$ satisfying $X_i = G_i(F)$, which is classically written (with some slight abuses) as $X = G(F)$ and at times $X = G \circ F$ 
(in the sense of maps of $\gA[X]^n$ to $\gA[X]^n$). 
\\
Note the converse as $\theta_F \circ \theta_G = \I_{\gA[X]}$.
\\
Here we will prove that $\theta_G \circ \theta_F = \I_{\gA[X]}$, or $X = F(G)$.
\\
Consequently (cf. Question 1) $G$ is uniquely determined, $\theta_F$ is an \auto of $\gA[X]$ and $F_1, \ldots, F_n$ are \agqt independent over $\gA$.
\\
The idea consists in using the \ri of formal series $\gA[[X]]$ or at least the quotient \ris $\gA[X]/\fm^d$ where $\fm = \gen {X_1, \ldots, X_n}$. Let $F = (F_1, \ldots, F_n) \in \gA[[X]]^n$. Study for which condition there exists a $G = (G_1, \ldots, G_n)$, $G_i \in \gA[[X]]$ without a constant term, satisfying $F(G) = X$. We then have $F(0) = 0$, and by letting $J_0 = \JJ(F)(0)$, we obtain $J_0 \in \GLn(\gA)$ (since $\JJ(F)(0) \circ \JJ(G)(0) = \I_{\Ae n}$). \\
We will prove the converse: in the case where $F(0) = 0$ and $J_0 \in \GLn(\gA)$, there \hbox{exists a $G = (G_1, \ldots, G_n)$}, with $G_i \in \gA[[X]]$, $G_i(0)=0$, and $F(G) = X$.
\begin{enumerate}\itemsep0pt
\item
By assuming this converse, prove that $G$ is unique and that
$G(F) = X$.
\item
Let $\gS \subset \gA[[X]]$ be the set of formal series without a constant term; $\gS^n$ is, with respect to the composition law, a \mo whose neutral \elt is $X$. Recall  Newton's method for solving an \eqn $P(z) = 0$ in $z$: introduce the iterator $\Phi : z \mapsto z - P'(z)^{-1}P(z)$ and the sequence $z_{d+1} = \Phi(z_d)$ with an adequate $z_0$; or a variant $\Phi : z \mapsto z - P'(z_0)^{-1}P(z)$. To solve $F(G) - X = 0$ in $G$, check that this leads to the iterator over $\gS^n$

\snic{
\Phi : G \mapsto G - J_0^{-1} \cdot (F(G) - X)
}

\item 
Introduce $\val : \gA[[X]] \to \NN \cup \{\infty\}$: $\val(g) = d$ means that $d$ is the (total) minimum degree of the \moms of $g$, agreeing that $\val(0) = +\infty$. We therefore have $\val(g) \ge d$ \ssi $g \in \fm^d$. For $g$, let $h \in \gA[[X]]$ \hbox{and $G$, $H \in \gA[[X]]^n$}
$$\preskip.3em \postskip.3em\ndsp 
d(f,g) = {1 \over 2^{\val(f-g)}}, \qquad
d(F,G) = \max_i d(F_i, G_i). 
$$
Show that $\Phi$ is a contracting map: $d\big(\Phi(G), \Phi(H)\big) \le d(G, H)/2$. 
Deduce that~$\Phi$ admits a unique fixed point $G \in \gS^n$, the unique
solution of $F(G) = X$.
\item
Solve the initial \pb with respect to \polsz.
\item
Check that the following \syss are \cdvs and make their inverses explicit (in $\ZZ[X,Y,Z]$ then in $\ZZ[X_1, X_2, X_3, X_4, X_5]$):
\[\preskip.2em \postskip.4em 
\begin{array}{ccc} 
 (X - 2fY - f^2 Z,\ Y + fZ,\ Z) \;\hbox { with  }\, f = XZ + Y^2 ,  \\[.8mm] 
 (X_1 + 3X_2X_4^2 - 2X_3X_4X_5,\  X_2 + X_4^2  X_5,\  X_3 + X_4^3,\
X_4 + X_5^3,\   X_5). 
\end{array}
\]
\end {enumerate}
}
\end {problem}


\vspace{-1em}
\sol


\exer{exoGensIdealEnsFini}
\emph{1.}
Lagrange interpolation:
${
Q = \sum_{\xi \in U} \left(\prod_{\zeta \in U \setminus \{\xi\}} 
{x_n-\zeta \over \xi-\zeta}\right) Q_\xi
}.$
\\  \emph{2.}
Assume that each $\fa(V_\xi) \subset \gK[x_1, \ldots, x_{n-1}]$ (for $\xi \in \pi_n(V)$) is generated by~$m$ \pols

\snic {
\fa(V_\xi) =  \gen{f^\xi_j, j \in \lrbm}, \quad 
f^\xi_j \in \gK[x_1, \ldots, x_{n-1}].
}

By item \emph{1}, there exists an $f_j \in \Kux$ satisfying $f_j(x_1, \ldots, x_{n-1},\xi) = f_j^\xi$ for all~$\xi \in \pi_n(V)$. Then prove, based on item~\emph{1}, that

\snic {
\fa(V) = \langle P, f_1, \ldots, f_m \rangle \qquad
\hbox {with} \  P = \prod_{\xi \in \pi_n(V)} (x_n - \xi).
}

Conclude by \recu on~$n$.


\exer{exothSymEl}
 \emph{4.} 
Consider the \pol \ri $\gB=\gA[s_1,\ldots,s_n]$ where the $s_i$'s are \idtrsz, then the \pol 
%
$
f(t)=t^n+\sum_{k=1}^{n}(-1)^ks_kt^{n-k}\in\gB[t].$
\\ 
Consider also the \adu

\snic{\gC=\Adu_{\gB,f}=\Bxn=\Axn,}

with, in $\gC[t]$, the \egt $f(t)=\prod_{i=1}^n(t-x_i)$.
\\
 Let $\rho:\AXn\to\Axn$ and $\varphi:\gA[s_1,\ldots,s_n]\to\gA[S_1,\ldots,S_n]$ be the \evn\homos $X_i\mapsto x_i$ and $s_i\mapsto S_i$. 

\vspace{-7pt}
\Deuxcol{.7}{.15}
{
We clearly have $\rho(S_i)=s_i$. Therefore, by letting $\rho_1$ be the restriction of $\rho$ to $\gA[\uS]$ and $\gA[\und{s}]
$, we have $\varphi\circ \rho_1=\Id_{\gA[\uS]}$ \hbox{and $\rho_1\circ \varphi=\Id_{\gA[\und{s}]}$}.
This shows that the $S_i$'s are \agqt independent over $\gA$ and we can identify $\gA[\uS]$ and~$\gA[\und{s}]=\gB$.
}
{
\xymatrix @R = 10pt{
\gA[\uX]\ar@<0.7ex>[r]^\rho &
\gA[\ux]\ar@<0.7ex>[l]^\psi 
\\ 
\gA[\uS]\ar@{^{(}->}[u]\ar@<0.7ex>[r]^{\rho_1} &
   \ar@<0.7ex>[l]^\varphi \gA[\und s]\ar@{^{(}->}[u]
\\ 
}
}

By the \uvl \prt of the \aduz, there exists a (unique) $\gB$-\homo $\psi:\gC\to\AuX$
which sends $x_i$ onto $X_i$. It follows that $\rho$ and~$\psi$ are two mutually reciprocal \isosz. Thus the $x_i$'s are \agqt independent over $\gA$ and $\AuX$ is free and of rank $n!$ over $\gA[\uS]=\gB$, with the prescribed basis.
\\
 NB: this \dem does not seem to simply give the fact that the \smq \pols of $\AuX$ are in $\gA[\uS]$.

\exer{exoPolSym1}
\emph{1.}
Let $f = (X_1^3 + X_2^3 + \cdots + X_n^3) - (S_1^3 - 3S_2S_2 + 3S_3)$.
It is a \hmg \smq \polz, therefore $f = g(S_1, \ldots, S_n)$ where $g = g(Y_1, \ldots, Y_n)$ is \hmg in weight, of weight $3$ with respect to the weight $\alpha_1 + 2\alpha_2 + \cdots + n\alpha_n$.\\
The \egt $\alpha_1 + 2\alpha_2 + \cdots + n\alpha_n = 3$ implies $\alpha_i = 0$ for~$i > 3$, so $g$ only depends on $Y_1$, $Y_2$, $Y_3$, say $g = g(Y_1, Y_2, Y_3)$. In the \egt 

\snic{(X_1^3 + X_2^3 + \cdots + X_n^3) - (S_1^3 - 3S_2S_2 + 3S_3) = g(S_1, S_2, S_3),}

put $X_i := 0$ for~$i > 3$; we obtain $g(S'_1, S'_2, S'_3) = 0$ where $S'_1, S'_2, S'_3$ are the \elr \smq functions of $X_1, X_2, X_3$. Deduce that $g = 0$ then~$f = 0$.
 
\emph{2.}
For the first, we can assume $n = 3$; we find $S_1S_2 - 3S_3$. For the other two which are \hmgz, \smqz, of degree $4$ we work with $4$ \idtrs and we obtain $S_1^2 S_2 - 2S_2^2 - S_1S_3 + 4S_4$ and
$S_2^2 - 2S_1S_3 + 2S_4$.

 \emph{3.}
Let $n > d$ and $f(X_1, \ldots, X_{n})$ be a \hmg \smq \pol of degree $d$. Let $h\in\gA[X_1,\ldots,X_d]=f(X_1,\ldots,X_d,0,\ldots,0)$.
If $h=0$, then $f = 0$.
\\
We can translate this result by saying that we have \isos of $\gA$-modules at the level of the \hmg \smq components of degree $d$:

\snuc{\cdots \to
\gA[X_1, \ldots, X_{d+2}]^{\rm sym.}_d
\vvvers{X_{d+2} := 0}
\gA[X_1, \ldots, X_{d+1}]^{\rm sym.}_d
\vvvers{X_{d+1} := 0}
\gA[X_1, \ldots, X_{d}]^{\rm sym.}_d
\;.}

\exer{exoMcCoy}
Let $\gA = \ZZ[U,V]/\gen {U^2, V^2} = \ZZ[u,v] = \ZZ \oplus \ZZ u \oplus
\ZZ v \oplus \ZZ uv$. \\
\emph{a.} We take $f = uT + v$ so $\fc = \gen {u,v}$.  We then have

\snuc{\Ann(u) = \gA u$, $\Ann(v) = \gA v$, $\Ann(\fc) = \Ann(u) \cap \Ann(v) = \gA uv$ and $\rD\big(\Ann(\fc)\big) = \fc.}

\sni
\emph{b.} Let $g = uT - v$.  We have $fg = 0$ but $g \notin \Ann(\fc)[T]$; we have $u \in \rD\big(\Ann(\fc)\big)$ but $u \notin \Ann_{\gA[T]}(f)$ (idem for $v$).

\exer{exoModCauBase}
{It suffices to prove item \emph{1.}
We have 
$$\preskip.3em \postskip.3em
f(T)=f(X_1)+(T-X_1)f_2(X_1,T)$$ by \dfn of $f_1=f$ and $f_2$.
Similarly 
$$
\preskip.3em \postskip.3em
f_2(X_1,T)=f_2(X_1,X_2)+(T-X_2)f_3(X_1,X_2,T)
$$ 
by \dfn of $f_3$. So 
$$\preskip.3em \postskip.3em
f(T)=
f(X_1)+(T-X_1)f_2(X_1,X_2)+(T-X_1)(T-X_2)f_3(X_1,X_2,T).
$$
Continue until
$$\arraycolsep2pt\preskip.4em \postskip-.2em
\begin{array}{rcl}
f_{n-1}(X_1,\ldots ,X_{n-2},T)&=&f_{n-1}(X_1,\ldots ,X_{n-2},X_{n-1})
\;+\\[1mm]
&   &
(T-X_{n-1})f_n(X_1,\ldots ,X_{n-1},T),
\end{array}
$$
which gives
$$\preskip.3em \postskip.3em
\arraycolsep2pt\begin{array}{rcl}
f(T)&   =&f_1(X_1)+(T-X_1)f_2(X_1,X_2)+(T-X_1)(T-X_2)f_3(X_1,X_2,X_3)
\\[1mm]
&   &
+\,\cdots\,+\,(T-X_1)\cdots (T-X_{n-1})f_n(X_1,\ldots ,X_{n-1},T).
\end{array}
$$
Finally, $f_n(X_1,\ldots ,X_{n-1},T)$ is \mon of degree $1$ in $T$
so
$$\preskip.3em \postskip.3em 
f_n(X_1,\ldots ,X_{n-1},T)=f_n(X_1,\ldots ,X_{n-1},X_n)+(T-X_n). 
$$
Note that this proves in particular that $f_n=S_1-s_1$.
}

\exer{exoPrimePowerRoot}
Let $f \in \gA[X]$ be \mon of degree $d$, with $f \divi X^p - a$ and $1 \le d \le p-1$. In a \ri $\gB \supseteq \gA$, we write $f(X) = \prod_{i=1}^d (X - \alpha_i)$, therefore $\alpha_i^p = a$ and~$\prod_i \alpha_i = b$ with~$b = (-1)^d f(0) \in \gA$. By lifting to the power $p$, $a^d = b^p$. However, $\pgcd(d,p) = 1$, so~$1 = ud + vp$, then $a = a^{ud} a^{vp} = (b^u a^v)^p$.

\exer{exoPrincipeIdentitesAlgebriques}
Let $e_{ij}$ be the matrix of $\Mn(\gA)$ having a single nonzero \coez, the \coe in position $(i,j)$, equal to $1$.
The module $S_n(\gA)$ is free and a basis is formed by the $e_{ii}$'s for $i \in \lrbn$ and the $e_{ij} + e_{ji}$ for $1 \leq i < j < n$.  It suffices to treat the case where 
$A = \Diag(\lambda_1, \ldots, \lambda_n)$.
Then, $\varphi_A= \Diag(\lambda_1^{2}, \ldots, \lambda_n^{2})$, \hbox{and $\varphi_A(e_{ij}+e_{ji})=\lambda_i\lambda_j (e_{ij}+e_{ji})$}.
Whence $\det(\varphi_A)=(\det A)^{n+1}$.

\exer{exoFreeFracTransfert}
Let $\ue = (e_1, \ldots, e_n)$ be a basis of $\gB/\gA$. Clearly $\ue$ is a $\gK$-free family. Let $x = b/b' \in \gL$ with $b \in \gB$, $b' \in \gB \setminus \{0\}$; we write

\snic{
x = (b\wi{b'}) / (b'\wi{b'}) =
b\wi{b'} / \rN\iBA(b') \in
\gK e_1 + \cdots + \gK e_n.}


\exer{exoSommesNewton}
\emph{1.} It suffices to prove it for the \gnq matrix $(a_{ij})_{i,j\in\lrbn}$ with \coes
in $\gA=\ZZ[(a_{ij})_{i,j\in\lrbn}]$. This matrix is \dig in an over\ri of~$\gA$.

\emph{2.} Follows immediately from \emph{1}.

\exer{exoCorpsFiniEltPrimitif}
\emph{1.} We have $\#R \le \sum_{d=1}^{n-1} q^d < 1 + q + \cdots + q^{n-1} 
= {q^n - 1\over q-1}$.
\\
 A fortiori, $\#R < q^n-1 < q^n$. Let $x \in \gL\setminus R$ and $d = \dex{\gK[x] : \gK}$. We have $x^{q^d} = x$, and since $x \notin R$, then $d = n$.

 \emph{2.}
The cyclotomic \pol is \ird in $\FF_2[X]$. Indeed, the only \ird \pol of degree $2$ of $\FF_2[X]$ is $X^2 + X + 1$, $\Phi_5$ has no root in~$\FF_2$, and $\Phi_5 \ne (X^2 + X + 1)^2$. We have $\#\gL = 2^4 = 16$, $\#\gL\eti = 15$, but $x^5 = 1$.

 \emph{3.}
Let $\sigma : \gL \to \gL$ be the Frobenius \auto of $\gL/\gK$, \cad $\sigma(x) = x^q$. We can easily check that $\gL = \gK[x]$ \ssi the $\sigma^i(x)$'s, $i \in \lrb{0..n-1}$, are pairwise distinct.
This condition is equivalent to $\,\sigma^k(x) = x \;\Rightarrow\; k \equiv 0 \bmod n\,$, \cad $\,x^{q^k} = x \;\Rightarrow\; k \equiv 0 \bmod n\,$. But

\snic {
x^{q^k} = x \iff x^{q^k-1} = 1 \iff o(x) \divi q^k - 1 \iff
q^k \equiv 1 \bmod o(x).
}

We then deduce, for $x \in \gL\eti$, that $\gL=\gK[x]$ \ssi the order of $q$ in the group of  \ivs \elts modulo $o(x)$ is exactly $n$.


\exer{exoPpcmPolsSeparables}
 \emph{1.}
We have $\gen {g,g'} \gen {g,h} \subseteq \gen {g, g'h} = \gen {g, g'h + gh'}$.
\\
Similarly, $\gen {h,h'} \gen {g,h} \subseteq \gen {h, g'h + gh'}$.
By evaluating the product we get

\snic {
\gen {g,g'} \gen {h,h'} \gen {g,h}^2 \subseteq 
\gen {g, g'h + h'g} \gen {h, g'h + h'g} \subseteq \gen {gh, g'h + h'g}.
} 

For the second item of the question we apply the result established above and Fact~\ref{factDiscUnit}. NB: this also results from \Eqnz~(\ref{eqfactDiscProd}),  Fact~\ref{factDiscProd}.

 \emph{2.} It suffices to treat the case of two \spl \pols $f$, $g \in \gA[T]$. 
\\
Let $h = \pgcd(f,g)$. We have $f = hf_1$, $g = hg_1$, with $\pgcd(f_1, g_1) = 1$.
\\
Since $g$ is \splz, $\pgcd(h, g_1) = 1$, so $\pgcd(hf_1, g_1) = 1=\pgcd(f, g_1)$.  \\
The \polsz~$f$,~$g_1$ are \splsz, \comz, therefore their product $\ppcm(f, g)$ is \splz.

\exer{exolemSousLibre} 
 \emph{1} and \emph{2.} These are special cases of what is stated in 
 Fact~\ref{fact.idd.sousmod}.

 \emph{3.} Suppose $L=\Ae m$. If $A\in\Mm(\gA)$ is a matrix whose columns form a basis of $F$, it is injective and its \deter is \ndzz. If $B$ is a matrix corresponding to the inclusion $F\subseteq E$, we have 

\snic{\idg{L:F}=\gen{\det A}$, $\;\idg{F:E}=\cD_m(B)\;$ and $\;\idg{L:E}=\cD_m(AB),}

hence the desired \egtz.

 \emph{4.} 
We have $\idg{N:\delta N}  = \gen {\delta^n}$.
We also have $\idg{N:\delta M} = \delta^{n-1} \gen {\delta, \an}$: take for the \sgr of $\delta M$ the family $\delta e_1, \dots, \delta e_n, \delta z$ where $e_1, \ldots, e_n$ is a basis of $N$ (we use $M = N + \gA z$), and compute the \idd of order $n$ of a matrix of the following type (for $n = 3$) 
$\cmatrix {
\delta & 0     & 0      & a_1 \cr
0      &\delta & 0      & a_2 \cr
0      &0      & \delta & a_3 \cr
}.$
\\
Then
$$\preskip-.4em \postskip.4em 
idg{N:\delta N} \,=\, \idg{N:\delta M}\, \idg{\delta M: \delta N} \, =\,  \idg{N:\delta M}\,\idg{M:N}, 
$$
\cad $\gen {\delta^n} = \idg{M:N}\; \delta^{n-1} \gen {\delta, \an}$.
\\
By simplifying by $\delta^{n-1}$ we obtain the \egt $\gen{\delta}=d\,\gen {\delta, \an}$.


\exer{exoDecompIdeal} \emph{1.} If $\fa\fa'=a\gA$ with $a$ \ndzz, then $\fb\subseteq\fa$ is equivalent to $\fb\fa'\subseteq a\gA$. Note that the test provides a \itf $\fc=\fb\fa'/a$ such that $\fa\fc=\fb$ in case of a positive response, and an \elt $b\notin\fa$ among the \gtrs of $\fb$ in case of a negative response.

 \emph{2.}
It is clear that the $\fq_i$'s are \iv (and thus \tfz).
\\
Perform the tests $\fb\subseteq\fq_i$. If a response is positive, for instance $\fb\subseteq\fq_1$, write~$\fc\fq_1=\fb$, whence $\fq_2\cdots\fq_n\subseteq\fc$, and finish by \recuz.
\\
 If all the tests are negative, we have some $x_i\in\fb$ and $y_i\in\gA$
such that $1-x_iy_i\in\fq_i$ (here suppose that the quotient \risz~$\gA\sur{\fq_i}$ are \cdisz), whence, by evaluating the product, $1-b\in \fq_1\cdots\fq_n\subseteq\fb$ with $b\in\fb$, so $1\in\fb$.
\\
Finally, we address the uniqueness question. Assume that $\fb=\fq_1\cdots\fq_k$.
\\
It suffices to prove that if a \tf \idema $\fq$ contains $\fb$, it is equal to one of the $\fq_i$'s ($i\in\lrbk$). 
\\
 Since we can test $\fq\subseteq\fq_i$, if each of the tests are negative we explicitly have $1\in\fq+\fq_i$ for each $i$ and so $1\in\fq+\fb$.
\\
NB: if we do not assume that $\fb$ is \tf and $\gA$ has explicit \dvez, the \dem of Kummer's little \tho would require that we at least know how to test~$\fq\subseteq\fb$ for every \gui{subproduct} $\fq$ of $\fq_1\cdots\fq_n$.


\exer{exoRabinovitchTrick} 
Assume $x \in \sqrt\fa$; as $\fa\subseteq\fb$, in $\gA[T]\sur\fb$, $\ov
x$ is nilpotent and invertible (since $\ov x \ov T = 1$), therefore $\gA[T]\sur\fb$ is the null \riz, \cad $1 \in \fb$.
\\
Conversely, suppose $1 \in \fb$ and reason in the \ri $\gA[T]/\fa[T] = (\gA\sur\fa)[T]$. Since $1 \in \fb$, $1-xT$ is \iv in this \riz, therefore $x$ is nilpotent in~$\gA\sur\fa$, \cad $x \in \sqrt\fa$.


\exer{exoDunford} \emph{(Decomposition of Jordan-Chevalley-Dunford)} \\ 
 \emph{Existence.}
Look for a zero $D$ of $f$,  a \gui{neighbor of $M$,} (i.e., with $M-D$ nilpotent), in the commutative \ri $\gK[M]$.
We have by hypothesis $f(M)^{k}=0$ for some $k\leq n$, and if $uf^{k}+vf'=1$, we obtain ${v(M)f'(M)=\In.}$

Consequently, the Newton method\imN, starting \hbox{with $x_0=M$}, 
gives the solution in $\gK[M]$ in $\lceil\log_2(k)\rceil$ iterations.

\emph{Uniqueness.}
The solution is unique, under the condition $f(D)=0$, in every commutative \ri containing~$\gK[M]$, for example in $\gK[M,N]$ if the pair $(D,N)$ solves the given \pbz.\\
When we only assume that the \polmin of $D$ is \splz, the uniqueness is more delicate.\\
A solution would be to directly prove that the \polcar of~$D$ is \ncrt equal to that of $M$, but it is not that simple.\footnote{In zero characteristic, 
one trick consists in retrieving the \polcar of a matrix $A$ from the $\Tr(A^{k})$ by following Le Verrier's method.}\\
Let us call $(D_1,N_1)$ the solution in $\gK[M]$ given by Newton's method.
Since $D$ and $N$ commute, they commute with $M=D+N$ and so with~$D_1$ and~$N_1$ because they belong to $\gK[M]$. From this we deduce that $D-D_1$ is nilpotent because it is equal to $N_1-N$ with~$N$ and~$N_1$ being nilpotents that commute. But the \algz~$\gK[D,D_1]$ is \'etale
by \Thref{corlemEtaleEtage}, so it is reduced, and $D=D_1$.


\exer{exoDunfordBis} 
We have $\gB = \gA[x] = \gA \oplus \gA x$ with $x$ \sply integral over $\gA$.
Let $\,z \mapsto \wi z\;$ be the \auto of the \Alg $\gB$ which swaps $x$ and $-b-x$.
\\
For~$z \in \gB$, we have ${\rm C}\iBA(z)(T) = (T - z)(T - \wi z)$.\\ 
Thus ${\rm C}\iBA(ax)(T) = T^2 + abT + a^2c$, and its \discri is equal to $a^2\Delta$.
\\
Let $\vep\in \gA$ be nonzero nilpotent and let $y = (\vep-1)x$. Then, $y$ is \spby integral over $\gA$ because $(\vep-1)^2\Delta$ is \ivz.
Furthermore, the \elt $z=x + y = \vep x$ is nonzero nilpotent. Assume that $\vep^2=0$ and let $g\in\AX$ be a \polu that annihilates $z$, we will prove that $g$ is not \splz.\\
Let us write $g(X)=u+vX+X^2h(X)$, then $z^2=0$, so $u+vz=0$.\\
Since $\gB=\gA\oplus\gA x$, we obtain $u=v\vep=0$, then \hbox{$g(X)=X\ell(X)$} \hbox{with  $\ell(0)=v$} non-\iv (otherwise, $\vep=0$). Finally, $\disc(g)=\disc(\ell)\ \Res(X,\ell)^2=\disc(\ell)\ v^2$ is non-\ivz.


\prob{exoDiscriminantsUtiles}~\\
\emph {1.}
Let $f(X) = X^n + c = (X-x_1) \cdots (X-x_n)$. Then, $f' = nX^{n-1}$ and

\snic {
\Res(f,f') =  f'(x_1) \cdots f'(x_n) =
n^n (x_1 \cdots x_n)^{n-1} = n^n \bigl((-1)^n c\bigr)^{n-1} =
n^n c^{n-1}.
}


Variant:
$$\preskip-.4em \postskip.4em\ndsp 
\Res(f', f) = n^n \Res(X^{n-1}, f) = n^n \prod_{i=1}^{n-1} f(0) 
= n^n c^{n-1}. 
$$
\emph {2.}
Let $f(X) = X^n + bX + c = (X - x_1) \cdots (X - x_n)$;

\snic {
\disc(f) = (-1)^{n(n-1)\over 2} \prod_{i=1}^n y_i 
\quad \hbox{with}\quad  y_i = f'(x_i) = n x_i^{n-1} + b\,. 
}

To compute the product of the $y_i$'s, we compute the product $P$ of the $x_i y_i$'s (that of the $x_i$'s is equal to~$(-1)^n c$). 
We have $x_iy_i = nx_i^n + bx_i = ux_i + v$,
with $u = (1-n)b$, $v = -nc$. We use the \elr \smq functions $S_j(x_1, \ldots, x_n)$ (almost all null)
$$
\preskip-.1em \postskip.1em\ndsp 
\;\;\;\prod_{i=1}^n (ux_i + v) = \sum_{j=0}^n u^j S_j(x_1, \ldots, x_n) v^{n-j}. 
$$
We get

\snic {
P = v^n + u^n S_n + u^{n-1} S_{n-1} v =
v^n + u^n (-1)^n c + u^{n-1} (-1)^{n-1} b v\,,
}

\cadz, by replacing $u$ and $v$ by their values
\[ \preskip.4em \postskip.4em
\begin {array}{rcl}
P &=&
 (-1)^n n^n c^n + (n-1)^n b^n c - n (n-1)^{n-1} b^n c \\
&=&
 (-1)^n n^n c^n + b^n c \bigl( (n-1)^n - n (n-1)^{n-1} \bigr) \\
&=& (-1)^n n^n c^n - b^n c  (n-1)^{n-1}.
\\
 \end{array}
\]
By dividing by $(-1)^n c$, we obtain the product of the~$y_i$'s then the stated formula.
 
\emph {3.}
Left to the sagacity of the reader who can consult \cite{Swan62}.
 
\emph {4.}
By letting $\Delta_p = \disc(\Phi_p)$, we have the \egt
$$\preskip.4em \postskip.4em 
 \disc(X^p - 1) = \Res(X-1,\Phi_p)^2 \disc(X-1)\Delta_p =
\Phi_p(1)^2 \Delta_p  = p^2 \Delta_p.
$$
By using  $\disc(X^n - 1) = (-1)^{n(n-1) \over 2} n^n (-1)^{n-1}$, 
we obtain
$$\preskip.2em \postskip.4em
\Delta_2 = 1, \qquad  \Delta_p =(-1)^{p-1\over 2} p^{p-2}
\quad\hbox {for $p \ge 3$}.
$$
\emph {5.}
Let $q = p^{k-1}$; let us first prove that $r := \Res(X^q - 1, \Phi_{p^k}) = p^q$.\\
With $X^q - 1 = \prod_{i=1}^{q} (X-\zeta_i)$, we have $r = \prod_{i=1}^{q}
\Phi_{p^k}(\zeta_i)$. In addition

\snic {
\Phi_{p^k}(X) = {Y^p - 1 \over Y-1} = Y^{p-1} + \cdots + Y + 1
\quad \hbox{with} \quad Y = X^q.
}

By making $X := \zeta_i$, we must make $Y := 1$, we obtain $\Phi_{p^k}(\zeta_i) = p$,
then $r = p^q$.
\\
Let $D_k = \disc(X^{p^k} - 1)$. Since $X^{p^k} - 1 = (X^q - 1)
\Phi_{p^k}(X)$, we have

\snic {
D_k = \Res(X^q-1, \Phi_{p^k})^2 D_{k-1} \disc(\Phi_{p^k}) =
p^{2q} D_{k-1} \disc(\Phi_{p^k}).
}

We use $\disc(X^n - 1) = (-1)^{n(n-1) \over 2} n^n (-1)^{n-1}$ for $n = p^k$ and $q$ 

\snic {
D_k/D_{k-1} = \varepsilon\, p^N, \quad
\varepsilon = \pm 1, \quad N = kp^k - (k-1)q = \big(k(p-1)+1\big)\,q\,.
}


For $\disc(\Phi_{p^k})$ to be obtained, $D_k/D_{k-1}$ must be divided by $p^{2q}$, which replaces the exponent~$N$ with $N - 2q = (k(p-1)-1) q$. As for the sign $\varepsilon$, for odd $p$,

\snic {
\varepsilon = (-1)^{{p^k - 1} \over 2} (-1)^{{q - 1} \over 2} =
(-1)^{{p^k - q} \over 2} = (-1)^{{p - 1} \over 2} 
.}

For $p = 2$, $\varepsilon = 1$ for $k \ge 3$ or $k = 1$ and $\varepsilon 
= -1$ for $k = 2$.
 
\emph{6.}
If $n$ is not the power of a prime, we can write $n = mp^k$ with $p$ prime, $\pgcd(m,p) = 1$, $k \ge 1$ and $m \ge 2$.
Then, $\Phi_n(X) = \Phi_m(X^{p^k}) /\Phi_m(X^{p^{k-1}})$, an \egt in which we put $X = 1$ to obtain $\Phi_n(1) = 1$.
The other items are easy.
 
\emph{7.}
Let $f$, $g$ be two \polusz, with $d = \deg f$, $e = \deg g$ and $d, e \ge 1$. Let $\gA[x] = \aqo{\gA[X]}{f(X)}$, $\gA[y] = \aqo{\gA[Y]}{g(Y)}$.
Let $f \otimes g$ be the \polcar of $x \otimes y$ in $\gA[x]\te_\gA\gA[y] = \aqo {\gA[X,Y]} {f(X),g(Y)}$. It is a \polu of degree $d\,e$. Since $f(X) = \prod_i (X - x_i)$, $g(Y) = \prod_j(Y - y_j)$, we obtain $(f\otimes g)(T) = \prod_{i,j} (T - x_iy_j)$. We easily see that

\snic {
\disc(f \otimes g) = \prod_{(i,j) < (i',j')} (x_iy_j - x_{i'}y_{j'})^2 =
\disc(f)^e \disc(g)^d f(0)^e g(0)^d \pi\,,
}

where $\pi \in \gA$ is the product $\prod_{i\ne i' ,\, j\ne j'} (x_iy_j - x_{i'}y_{j'})$.
\\
Let $n$, $m \ge 2$ with $\pgcd(n,m) = 1$ and $\zeta_n$, $\zeta_m$, $\zeta_{nm}$
be the roots of the unit of respective orders $n$, $m$, $nm$. By the Chinese remainder \thoz, 
we obtain the \egt $\Phi_{nm} = \Phi_n \otimes \Phi_m$. 
As $\Phi_n(0) = \Phi_m(0) = 1$ (since $n, m \ge 2$), we have the \egt
$$
\preskip.0em \postskip.3em 
\Delta_{nm} = \Delta_n^{\varphi(m)}\Delta_m^{\varphi(n)} \,\pi\,, 
$$
where $\pi \in \ZZ$ is the following product.
$$\preskip.2em \postskip.4em \ndsp
 \prod_{i\ne i' ,\, j\ne j'}
(\zeta_n^i \zeta_m^j - \zeta_n^{i'}\zeta_m^{j'}) \hbox{, for  }i,i'\in(\ZZ\sur{n\ZZ})\eti\hbox{ and }j, j'\in(\ZZ\sur{m\ZZ})\eti.
$$
Let $C \subset (\ZZ\sur{nm\ZZ})\eti \times (\ZZ\sur{nm\ZZ})\eti$ be the set of pairs $(a,b)$ with  $a, b$ \iv modulo $nm$, $a \not\equiv b \bmod n$, $a \not\equiv b \bmod m$. The Chinese remainder \tho gives~us

\snic {
\pi = \prod_{(a, b) \in C}(\zeta_{nm}^a - \zeta_{nm}^b). 
}


Let $z \mapsto \ov z$ be  complex conjugation. Then, $\pi$ is of the form $z\ov z$, therefore  $\pi\in\NN\etl$.
\\
Indeed, $(a,b) \in C \Rightarrow (-a,-b) \in C$ with $(a,b) \ne (-a,-b)$.\\ 
Furthermore, for $c \in \ZZ$ a non-multiple of $n$ or $m$, consider the \elt $\zeta_{nm}^c$ which is of order $nm/\pgcd(c,nm) = n'm'$ with 
$n' = n/\pgcd(c,n) > 1$, $m' > 1$ and $\pgcd(n',m') = 1$. 
Therefore $n'm'$ is not a power of a prime number, and, by the previous question, $1 - \zeta_{nm}^c$ is \iv in~$\ZZ[\zeta_{nm}^c]$, a fortiori in~$\ZZ[\zeta_{nm}]$. We deduce that $\pi$ is invertible in~$\ZZ[\zeta_{nm}]$, therefore in $\ZZ$. 
\\
Recap: 
$\pi = 1$, and $\Delta_{nm} = \Delta_n^{\varphi(m)}\Delta_m^{\varphi(n)}$.
\\
Finally, if the formula that gives the cyclotomic \discri is satisfied for two pairwise \com integers $n$ and $m$, it is satisfied for the product $nm$ (use the first item). 
However, it is true for integers
which are powers of a prime by Question~\emph{5}, therefore it is true for every integer $\ge 3$.


\prob{exoAnneauEuclidien}
\emph{\ref{i5exoAnneauEuclidien}.} Consider $p\equiv1\mod 4$. The \polz~$Y^{p-1 \over 2} - 1 \in \FF_p[Y]$ is of degree $< \#\FF_p\eti$.
There thus exists a  non-root $y \in \FF_p\eti$ of this \polz; let $x = y^{p-1 \over 4}$ so that $x^2 = y^{p-1 \over 2} \ne 1$; but $x^4 = 1$ thus $x^2 =
-1$. Actually, for half of the $y \in \FF_p\eti$, we have $y^{p-1 \over 2} = 1$ (the squares), and for the other half (the non-squares), we have~$y^{p-1 \over 2} = -1$.
\\
Let us address the question of the efficient \algoz. What we mean by this is that the execution time has a small power of the number of digits of~$p$ as its order of magnitude.
\\
We first determine some $x \in \FF_p$ such that $x^2 = -1$. For that we randomly draw integers $y$ over $\lrb{2..(p-1)/2}$ and we compute $y^{p-1 \over 4}$ in $\FFp$ (for that we use an efficient \algo of exponentiation modulo $p$). The probability of failure (when the result is $\pm1$) is of $1/2$ at each draw.  
\\
Once such an $x$ is found, it remains to compute $\pgcd(x+i,p)$ with the Euclidean \algoz. As the norm is divided by at least $2$ at each step, the \algo is efficient.
\\
NB: the brute force method which would consist in saying \gui{since $p\equiv 1 \mod 4$, it possesses a factor of the form $m+in$, and all that is left to do is try out every $m<p$} quickly proves to be impractical as soon as $p$ is large enough. 

\emph{\ref{i6exoAnneauEuclidien}.} The \dcn of the prime divisors of $m$ is treated in the previous item. It remains to decompose $n+qi$. 
\\
 Regarding the \dcn of $n^2+q^2$, we already know that the only prime numbers therein are $2$ (with the exponent $1$) or some $p\equiv 1\mod 4$. \\
 If $u+vi$ is the factor of some $p$ that divides $n^2+q^2$, then $u+vi$ or $u-vi$ \hbox{divides $n+qi$}. If $p$ appears with the exponent $k$ in $n^2+q^2$, and if $u+vi$ divides $n+qi$, then $u+vi$ appears with the exponent $k$ in~$n+qi$.
\\
If $s=2^k\prod_ip_i^{m_i}\prod_jq_j^{n_j}$ with every $p_i\equiv 3\mod 4$
 and every $q_j\equiv 1 \mod 4$, then the condition insuring that $s$ is the sum of two squares is that every $m_i$ be even.
\\
 Note that an expression $s=a^2+b^2$ with $0<a\leq b$ corresponds to two conjugated \elts $a\pm i b$ defined up to association  (for example multiplying by $i$ comes down to permuting $a$ and $b$). 
It follows that in the case where $s$ is the sum of two squares, the number of expressions of $s$ as a sum of the squares is equal to $(1/2)\prod_j(1+n_j)$ unless the $n_j$'s are all even, in which case we add or subtract $1/2$ depending on whether we consider that an expression $a^2+0^2$ is or is not legitimate as a sum of two squares. 
\\
 For example with  $5=\rN(a)$, $a=2+i$ and $13=\rN(b)$, $b=3+2i$ we obtain
{\footnotesize 
$$\!\!{\begin{array}{rcl} 
5=\rN(a)  &\hbox{gives}& 5=2^2+1^2,  \\[1mm] 
10=\rN\big(a(1+i)\big)=\rN(1+3i)  &\hbox{gives}&  10=1^2+3^2, \\[1mm] 
5^3=\rN(a^3)=\rN(5a)  &\hbox{gives}& 125=2^2+11^2=10^2+5^2,  \\[1mm] 
5^4=\rN(a^4)=\rN(5a^2)=\rN(25)  &\hbox{gives}& 625=7^2+24^2=15^2+20^2=25^2+0,  \\[1mm] 
5^2\times 13=\rN(a^2 b)=\rN(a^2 \ov b)=\rN(5 b)  &\hbox{gives}& 325=18^2+1=17^2+6^2=15^2+10^2.   
\end{array}}$$}

\vspace{-8pt}
Similarly 
$5^3\times 13=\rN(a^3b)=\rN(a^3\ov b)=\rN(5ab) = \rN(5a\ov b)$
gives

\snic{1625=16^2+37^2=28^2+29^2= 20^2+35^2= 40^2+5^2.}

An analogous computation gives 

\snic{1105=5\times 13\times 17=9^2+32^2=33^2+4^2=23^2+24^2= 31^2+12^2.}


\prob{exoPetitKummer} \emph{1.} The \discri can be specialized and $\Delta$ 
is \iv modulo~$p$.
\\
 Next note that $\aqo{\ZZ[\alpha]}p\simeq\FFp[t]:=\aqo{\FFp[T]}{f(T)}$. 
This already implies that the \ids $\gen{q_k,p}$ are maximal in $\ZZ[\alpha]$.
For $j\neq k$, $\gen{Q_j(t)}+\gen{Q_k(t)}=\gen{1}$
\hbox{in $\FFp[t]$}, \hbox{so $\gen{q_j}+\gen{q_k}+\gen{p}=\gen{1}$} 
\hbox{in $\ZZ[\alpha]$}. Whence $\gen{q_j,p}+\gen{q_k,p}=\gen{1}$. 
\\
By the Chinese remainder \thoz, 
the product of the $\gen{q_k,p}$ is therefore equal to their intersection, which is equal to $\gen{p}$ because the intersection of the $\gen{Q_j(t)}$ in $\FFp[t]$
is equal to their product, which is null. 
\\
 Note that the \egt $\gen{p}=\prod_{k=1}^{\ell}\gen{p,Q_k(\alpha)}$ is maintained in every \ri containing $\ZZ[\alpha]$. Similarly for the \com character of the \idsz.
\\
 If we move from $\ZZ[\alpha]$ to $\gA$, then the only thing left to check is that \hbox{the $\gen{p,q_k}$'s} remain as \idemasz. This is indeed the case and the quotient fields are \isocz. Indeed, every \elt of $\gA$ is of the form $a/m$ where $a\in\ZZ[\alpha]$ and~$m^2$ divides $\Delta$ 
 (Proposition~\ref{propAECDN}). Since $m$ is \com to $p$ the natural \homo $\aqo{\ZZ[\alpha]}{p,q_k}\to\aqo\gA{p,q_k}$ is an \isoz. 

\emph{2.} Apply Exercise~\ref{exoDecompIdeal}.


\prob{exoPolCyclotomique}
\emph{1a.}
For primes $p_1$, $p_2,$ \ldots\, that do not divide $n$, we deduce that $f(\xi^{p_1p_2\ldots}) = 0$, \cad $f(\xi^m) = 0$ for every $m$ such that $\pgcd(n,m) = 1$, or even that~$f(\xi') = 0$ for every $\xi'$, $n^{\rm th}$ primitive root of the unit. So $f = \Phi_n$.
\\
 \emph{1b.}
Let $h(X) = \pgcd_{\QQ[X]}\big(f(X),g(X^p)\big)$. 
By \KROz's \tho $h\in\ZZ[X]$. We have $h(\xi) = 0$, therefore $\deg h \ge 1$. Let us reason modulo $p$. We have $g(X^p) = g(X)^p$, \hbox{so $\ov h \divi \ov f$} \hbox{and $\ov h \divi \ov g^p$}. If  $\pi$ is an \ird factor of $\ov h$, $\pi^2$ is a square factor of $X^n -\ov 1$, but $X^n - \ov 1$ is \spl in~$\FF_p[X]$.
\\
Note: the \discri of the \polz~$X^n + c$ is $(-1)^{\frac{n(n-1)}2}n^n c^{n-1}$, in particular that of $X^n-1$ is $(-1)^{\frac{(n+2)(n+3)}2}n^n$.
 
 \emph{2.}
If $G$ a cyclic group of order~$n$, 
we have the classical \isos 

\snuc{\End(G)\simeq  \ZZ/n\ZZ\hbox{ (as \risz) and }\Aut(G)\simeq\big((\ZZ/n\ZZ)^{\!\times},\times \big) \hbox{ (as groups)}.}

Whence canonical \isos  $\Aut(\UU_n)\simeq(\ZZ/n\ZZ)^{\!\times}\simeq\Gal(\gQ_n/\QQ)$.\\ 
If $m \in (\ZZ/n\ZZ)^{\!\times}$, we obtain the \auto $\sigma_m$ of $\gQ_n$  defined by $\sigma_m(\zeta)=\zeta^{m}$ for $\zeta\in \UU_n$.
 
\emph{3.}
Assume know a field of roots $\gL$ as a \stfe extension of $\gK$. The map $\sigma \mapsto \sigma\frt{\UU_n}$ is an injective morphism of $\Aut_\gK(\gL)$ into $\Aut(\UU_n)$. In particular, $\Aut_\gK(\gL)$ is \isoc to a subgroup of $(\ZZ/n\ZZ)^{\!\times}$. Moreover, for every $n^{\rm th}$ primitive root of  the unit $\xi$ in $\gL$, we have $\gL = \gK(\xi)$. So, every \ird factor of $\Phi_n(X)$ in $\gK[X]$ has the same degree $\dex{\gL : \gK}$. However, it is not a priori obvious to determine what type of operation on $\gK$ is \ncr to factorize $\Phi_n(X)$ in $\gK[X]$.
We now give an example where we can determine with certainty $\dex{\gL : \gK}$: let  $p\ge 3$ be a prime, $p^* = (-1)^{p-1 \over 2}p$ and $\gK = \QQ(\sqrt {p^*})$. Then $\gK \subseteq \gQ_p$ (Gauss), the only $p^{\rm th}$ root of the unit contained in $\gK$ is $1$ and $\Phi_p(X)$ can be factorized in $\gK[X]$ as a product of two \ird
\pols of the same degree~$p-1 \over 2$.


\prob{exoCyclotomicRing} 
\emph{1a.}
On the one hand we have $\gA\sur\fp_i \simeq \aqo {\FF_p[X]}{\ov {f_i}}$ so $\fp_i$
is maximal. On the other hand, let $\ov\gA = \gA\sur{p\gA} \simeq \aqo {\FF_p[X]} {\ov {\Phi_n}}$ and $\pi : \gA \twoheadrightarrow \ov\gA$ be the canonical surjection; then $\sqrt {p\gA} = \pi^{-1} \big(\rD_{\ov\gA}(0)\big)$ and

\snic {
\rD_{\ov\gA}(0) = \aqo {\gen {\ov g}} {\ov {\Phi_n}} \simeq
\aqo {\gen {\ov {f_1}}} {\ov {\Phi_n}} \times \cdots \times
\aqo {\gen {\ov {f_k}}} {\ov {\Phi_n}}
,}

hence the result.

\emph {1b.}
Results from the fact that $\Phi_n$ is \spb modulo $p$.

\emph {1c.}
We easily check the following \egts in $\ZZ[X]$
$$\preskip.2em \postskip.4em
\Phi_n(X) = \Phi_{mp}(X^{p^{k-1}}) =
{\Phi_{m}(X^{p^{k}}) \over \Phi_{m}(X^{p^{k-1}}) }
,$$
and thus in $\FF_p[X]$, by letting $\varphi$ be the Euler's indicator function
$$\preskip.4em \postskip.4em
\Phi_n(X) = {\Phi_{m}(X)^{p^{k}} \over \Phi_{m}(X)^{p^{k-1}}} =
\Phi_{m}(X)^{\varphi(p^{k})} \qquad \bmod p
.$$
The \polz~$\Phi_m$ is \spb modulo $p$ so the subset without a square factor of $\Phi_n$ modulo $p$ is $\ov g = \ov {\Phi_m}$; whence
$\sqrt {p\gA} = \gen {p, \Phi_m(\zeta_n)}$.
\\
Let us prove that $p \in \gen {\Phi_m(\zeta_n)}$. If $\zeta_p \in \UU_n$ is a $p^{\rm th}$ primitive root of the unit, we have the \egt

\snic {
\Phi_p(X) = \sum_{i = 0}^{p-1} X^i = \prod_{j=1}^{p-1} (X - \zeta_p^j),
}

hence, by making $X := 1$

\snic {
p = \prod_{j=1}^{p-1} (1 - \zeta_p^j) \in \gen {1 - \zeta_p}.
}


By applying this to $\zeta_p = \zeta_n^{mp^{k-1}}$, we obtain $p \in 
\big\langle{1 - \zeta_n^{mp^{k-1}}}\big\rangle$. \\
However, $X^{mp^{k-1}} - 1$ is a multiple of $\Phi_m$ in $\ZZ[X]$, therefore $\zeta_n^{mp^{k-1}} - 1$ is a multiple of $\Phi_m(\zeta_n)$ in $\gA$, whence $p \in \gen {\Phi_m(\zeta_n)}$.

\emph{1d.}
As $\sqrt {p\gA} = \fp_1 \cdots \fp_k = \gen {\Phi_m(\zeta_n)}$ is finitely generated, there is an exponent $e$ such that $(\fp_1 \cdots \fp_k)^e \subseteq p\gA$ and we apply Exercise~\ref{exoDecompIdeal}.  \\
Note: we can take $e = \varphi(p^k) = p^k - p^{k-1}$.

\emph{2.}
The first item is immediate. Next, if $\fa$ is a nonzero \itf of $\gA$, it contains a nonzero \elt $z$. Then, $a = \rN_{\gQ_n\sur\QQ} (z) = z\wi z$ is a nonzero integer belonging to $\fa$. We write $a\gA \subseteq \fa$ as a product of \iv \idemas and we again apply Exercise~\ref{exoDecompIdeal} to the \id $\fa$.

\prob{exoSommeGauss}
\emph{1.}
Let $x_0 \in G$ such that $\varphi(x_0) \ne 1$.  \\
We write $\sum_{x \in G} \varphi(x) = \sum_{x \in G} \varphi(xx_0)$, therefore $S \varphi(x_0) = S$ with $S = \sum_{x \in G} \varphi(x)$, \cad $\big(1 - \varphi(x_0)\big) S = 0$, whence $S = 0$.
\\
\emph{2.}
First note that $\chi^{-1}(-1) = \chi(-1)$ since $\chi(-1)^2 = \chi\big((-1)^2\big) = 1$. We write

\snic {
\sum_{x+y = z} \chi(x) \chi^{-1}(y) = 
\sum_{x \ne 0,z} \chi\left( {x \over z-x} \right)
.}

If  $z \ne 0$, the map $x \mapsto {x \over z-x}$ is a bijection of $\gk \cup
\{\infty\}$ onto $\gk \cup \{\infty\}$ which transforms $z$ into $\infty$, $\infty$
into $-1$, $0$ into $0$, which gives a bijection of $\gk\eti \setminus \{z\}$
onto~$\gk\eti \setminus \{-1\}$. We can therefore write

\snic {
\sum_{x+y = z} \chi(x) \chi^{-1}(y) = 
\sum_{v \in \gk\eti  \setminus \{-1\}} \chi(v) =
\sum_{v \in \gk\eti} \chi(v) - \chi(-1) = 0 -\chi(-1).
}

If  $z = 0$ we have the \egt

\snic {
\sum_{x+y=z} \chi(x) \chi^{-1}(y) = \sum_{x \ne 0} \chi(-1) = (q-1)\chi(-1).
}

\emph{3.}
We write

\snic {
G_\psi(\chi)G_\psi(\chi^{-1}) = \sum_{x,y} \chi(x) \chi^{-1}(y) \psi(x+y) =
\sum_{z \in \gk} S(z) \psi(z), 
}


with  $S(z) = \sum_{x+y = z} \chi(x) \chi^{-1}(y)$.  Whence

\snic {\arraycolsep2pt
\begin {array}  {rcl}
G_\psi(\chi)G_\psi(\chi^{-1}) &=& 
(q-1) \chi(-1) - \chi(-1) \sum_{z \ne 0} \psi(z) 
\\[1mm]
&=&
q\chi(-1) - \chi(-1) \sum_{z \in \gk} \psi(z) = q\chi(-1).
\\
\end {array}
}

\emph{4.}
The first item is immediate.
We easily have $\tau_0 \tau_1 = {1 - p^* \over 4}$. The rest follows.

\prob{exoDedekindPolynomial}
\emph{1.}  
If  $g(x) = 0$, with  $x \in \ZZ$ and $g(X) \in \ZZ[X]$ \monz, then $x \divi
g(0)$. Here $\pm1, \pm2, \pm4, \pm8$ are not roots of $f(X)$, therefore this \pol is \irdz. The \discri of the \polz~$X^3 + aX^2 + bX + c$ is

\snic {
18abc - 4a^3c + a^2b^2 - 4b^3 - 27c^2 ,\quad
\hbox{hence the result for $a = 1$, $b = -2$, $c = 8$.}
}

\emph{2.} 
The \elt $\beta = 4\alpha^{-1} \in \QQ(\alpha)$ is integral over~$\ZZ$
since 

\snuc {
\alpha^3 + \alpha^2 - 2\alpha + 8 = 0 \buildrel /\alpha^3 \over \Longrightarrow
    1 + \alpha^{-1} - 2\alpha^{-2} + 8\alpha^{-3} = 0
\buildrel \times 8 \over \Longrightarrow
8 + 2\beta - \beta^2 + \beta^3 = 0.
}

To check that $\gA = \ZZ \oplus \ZZ\alpha \oplus \ZZ\beta$ is a \riz, it suffices to see that $\alpha^2, \alpha\beta, \beta^2 \in
\gA$.  It is clear for $\alpha\beta = 4$. We have $\alpha^2 + \alpha - 2 +
2\beta = 0$, so $\alpha^2 = 2 - \alpha - 2\beta$,
and since~$\beta^3 - \beta^2 + 2\beta + 8 = 0$,  $\beta^2 = \beta - 2 - 8\beta^{-1} = \beta - 2 - 2\alpha$. 
\\
The expression of~$(1, \alpha, \alpha^2)$ over the basis $(1, \alpha, \beta)$ is provided by the matrix
$$\preskip.4em \postskip.4em 
\matrix {
 \phantom{\J\limits_{i}}                 &
\matrix {1 & \alpha & \alpha^2}
                                   \cr
~~\matrix{1\cr\alpha\cr\beta}         &
\crmatrix{   1 & 0      &  2 \cr
            0 & 1      & -1 \cr
            0 & 0      & -2 \cr }
}. 
$$
The \ri $\ZZ[\alpha]$ is therefore of index~$2$ in~$\gA$; but
$$\preskip.2em \postskip.4em 
\Disc_{\ZZ[\alpha]/\ZZ} =  \idg{\gA : \ZZ[\alpha]} ^2 \cdot \Disc_{\gA/\ZZ}
\quad \hbox{so} \quad \Disc_{\gA/\ZZ} = -503. 
$$
Since the \discri of~$\gA$ is squarefree, $\gA$ is the \ri of integers of~$\QQ(\alpha)$.

\emph{3.} 
Let us prove that $\alpha$, $\beta$ and $\gamma := 1 + \alpha + \beta$ form, modulo~$2$, a \sfio
$$\preskip.3em \postskip.2em
\alpha + \alpha^2 = 2 - 2\beta,\quad \beta^2 - \beta = 2 - 2\alpha,\quad
\alpha\beta = 4,
$$
hence modulo~$2$

\snic {
\alpha \equiv \alpha^2,\quad\ \beta \equiv \beta^2,\quad
\gamma^2 \equiv \gamma,\quad \alpha + \beta + \gamma \equiv 1,\quad
\alpha \beta \equiv 0, \quad \alpha \gamma \equiv 0,\quad
\beta \gamma \equiv 0.
}

We therefore have $\gA/2\gA = \FF_2 \ov \alpha \oplus \FF_2 \ov \beta \oplus \FF_2 \ov \gamma$. If we want to compute the \fcn of~$2$ in~$\gA$, we notice that $(\alpha, \beta, \gamma)$ is a $\ZZ$-basis of~$\gA$ and by denoting by~$\pi$ the morphism of reduction modulo~$2$, $\pi :\gA \to \gA/2\gA$, the \ideps of~$\gA$ over $2$ are the inverse images of the \ideps of $\gA/2\gA$. For example~$\fa = \pi^{-1}(\{0\} \,\oplus\, \FF_2 \ov \beta \,\oplus\, \FF_2 \ov \gamma) = \gen {2\alpha, \beta, \gamma}$.  
Thus by letting~$\fb = \gen {\alpha, 2\beta, \gamma}$ and~$\fc = \gen {\alpha, \beta, 2\gamma}$, we have~$\gA/\fa \simeq \gA/\fb \simeq \gA/\fc \simeq \FF_2$ and~$2\gA = \fa\fb\fc = \fa \cap \fb \cap \fc$.
\\
In \gnlz, let $\gK$ be a number field satisfying $\dex{\gK : \QQ} \ge 3$ and~$2$ be completely decomposed in the \ri of integers $\gZ_\gK$. Then,~$\gZ_\gK$  is not monogenic, \cad there exists no $x \in \gZ_\gK$ such that $\gZ_\gK = \ZZ[x]$.  Indeed, $\gZ_\gK/2\gZ_\gK \simeq \FF_2^n$ and $\FF_2^n$ does not admit any  primitive \elt over $\FF_2$ if $n \ge 3$.

\emph{4.} 
By multiplying $1 \in \ff + \fb$ by $\gB'$, we obtain $\gB' \subseteq \ff\gB' + \fb' \subseteq \gB + \fb'$, which shows that $\gB \to \gB'/\fb'$ is surjective. Let us prove that $\gB \to \gB'/\fb'$ is injective, i.e.~$\fb' \cap \gB = \fb$. By multiplying $1 \in \ff+\fb$ by $\fb'\cap\gB$ we obtain the inclusions

\snic{\fb' \cap \gB \subseteq (\fb' \cap \gB)\ff + (\fb' \cap \gB)\fb
\subseteq \fb\gB'\ff + \fb \subseteq \fb\gB + \fb \subseteq \fb.}

\emph{5.} 
In the previous context, let $x \in \gZ_\gK$ be of degree $n = \dex{\gK : \QQ}$. \\
Let 
$d = \idg{\gZ_\gK : \ZZ[x]}$. We have $d\gZ_\gK \subseteq \ZZ[x]$ and $d$ can serve as conductor of $\gZ_\gK$ into $\ZZ[x]$. 
If $2 \nedivi d$, by the Dedekind avoidance, $\gZ_\gK/2\gZ_\gK \simeq \ZZ[x]/2\ZZ[x] = \FF_2[\ov x]$. 
But $\gZ_\gK/2\gZ_\gK \simeq \FF_2^n$ does not admit a primitive \elt over $\FF_2$ for $n \ge 3$.

\prob{exoGaloisNormIdeal}
\emph{1.} $z \in \gB$ is a root of $\prod_{\sigma \in G}(T - z)$, a \polu with \coes in $\gA$.

\emph{2.} 
$\ov\fm = \fm$ is clear. Let us compute $\fm^2$ by letting $d = 4q+1$,
so $1+d = 2(2q+1)$:

\snic {\arraycolsep2pt
\begin{array} {rcl}
\fm^2 &=& \gen {1 + 2\sqrt d + d, 1-d, 1 - 2\sqrt d + d} \\[1mm]
&=& 2 \gen {2q+1 + \sqrt d, 2q, 2q+1 - \sqrt d} =
2 \gen {1 + \sqrt d, 1 - \sqrt d} = 2\fm. 
\end{array}
}

In addition, as a $\ZZ$-module, $\fm = \ZZ(1+\sqrt d) \oplus \ZZ(1-\sqrt d)
= 2\ZZ \oplus \ZZ(1\pm\sqrt d)$. We cannot simplify $\fm^2 = 2\fm$ by $\fm$ (because $\fm \neq 2\gB$ seeing that $1 \pm \sqrt d \notin 2\gB$), therefore $\fm$ is not \ivz. We have $\rN_G(\fm) = 2\ZZ$ therefore $\rN_G(\fm)\gB = 2\gB \ne
\rN'_G(\fm)$.
\\
The canonical map $\ZZ \to \gB/\fm$ is surjective (since $x+y\sqrt d \equiv x+y \bmod\fm$) with kernel $2\ZZ$, so $\FF_2 \simeq \gB/\fm$, and $x + y\sqrt d \mapsto (x+y) \bmod 2$ defines a surjective morphism \emph{of \risz} $\gB \twoheadrightarrow \FF_2$, with kernel $\fm$.
\\
Let $\rN(\fb) = \#(\gB/\fb)$ for nonzero $\fb$. If $z = x(1+\sqrt d) + y(1-\sqrt d) \in\fm$ with~$x$,~$y \in \ZZ$, then $\rN_G(z) = (x+y)^2 - d(x-y)^2 \equiv 4xy \bmod 4$. \\
So $\rN_G(z) \in 4\ZZ$ for $z \in \fm$, but $\rN(\fm) = 2$. We have $\rN(\fm^2) = \rN(2\fm) = 4\rN(\fm) = 8$, but~$\rN(\fm)^2 = 4$.

\emph{3.}
Let $\fb = \gen {b_1, \ldots, b_n}$ and let $\uX = (X_1, \ldots, X_n)$ be $n$ \idtrsz.  Let us introduce the normic \pol $h(\uX)$

\snic {
h(\uX) = \prod\nolimits_{\sigma \in G} h_\sigma(\uX)
\quad \hbox{with} \quad
h_\sigma(\uX) = \sigma(b_1) X_1 + \cdots +  \sigma(b_n) X_n.
}

We have $h(\uX) \in \gA[\uX]$. Let $d$ be a \gtr of $\rc(h)_\gA$.
As $\gB$ is \icl and $\rc(h)_\gB = d\gB$ is principal, we can apply Proposition~\ref{propArm}: we then have~$\prod_\sigma \rc(h_\sigma)_\gB = \rc(h)_\gB = d\gB$, \cad  $\rN'_G(\fb) = d\gB$.
\\
Since  $\gA$ is Bézout, it is \iclz. 
Let $a \in \gA \cap d\gB$.
Then the \elt $a/d \in \Frac(\gA)$ is integral over $\gA$ (because $a/d \in \gB$) so $a/d \in \gA$, \cad $a \in d\gA$. \\
Recap: $\gA \cap d\gB = d\gA$ \cad $\rN_G(\fb) =
d\gA$.
\\
By \dfnz, the \evns of the normic \pol $h$ over $\gB^n$ are the norms of \elts of the \id $\fb$; they belong to the \id of $\gA$ generated by the \coes of the normic \polz, this \id of $\gA$ being $\rN_G(\fb)$.
\\
If $\#G = 2$, the \coe of $X_1X_2$ in $h$ is

\snic {
h(1,1, \ldots, 0) - h(1, 0, \ldots,0) - h(0, 1, \ldots,0) =
\rN_{G}(b_1+b_2) - \rN_{G}(b_1) - \rN_{G}(b_2).
}

This in fact reduces to writing $b_1\ov{b_2} + b_2\ov{b_1} = \rN_{G}(b_1+b_2) - \rN_{G}(b_1) - \rN_{G}(b_2)$. Similarly, the \coe of $X_iX_j$ in $h$ is, for $i \ne j$, $\rN_{G}(b_i+b_j) - \rN_{G}(b_i) - \rN_{G}(b_j)$. Consequently, the \id of $\gA$ generated by the norms $\rN_G(b_i)$ and $\rN_G(b_i + b_j)$ contains all the \coes of $h(\uX)$. It is therefore the \id $\rN_G(\fb)$.


\prob{exoLemmeFourchette} 
 \emph{(Forking lemma)}\index{Forking lemma}
\\
\emph{1.}
For $x \in \gL$, we have $x = \sum_j \Tr_{\gL/\gK}(xe_j)e'_j$.
\\
If $x \in \gB$, then $\Tr_{\gL/\gK}(xe_j)$ is an \elt of $\gK$ integral over $\gA$ so in $\gA$. This proves the middle inclusion.
\\
By writing $e_i = \sum_j \Tr_{\gL/\gK}(e_ie_j)e'_j$, we obtain 

\snic{\tra {\ue} = A \tra {\ue'}$ where $A = \big(\Tr_{\gL/\gK}(e_ie_j)\big)
\in \Mn(\gA)$, with  $\det(A)=\Delta,}

the right-hand side inclusion.

\emph{2.} The \ZZmo $F_k$ is the intersection of $\gB$ and $Z_k$, which are two sub\mtfs of $Z_{n-1}$, free, of rank $n$. It is therefore a free \ZZmo of finite rank,
and the two inclusions $\delta Z_k\subseteq F_k\subseteq Z_k$ show that $F_k$ is of rank $k+1$.   
\\
 The \ZZmo $\pi_k(F_k)$ is a \tf sub\ZZmo of $\fraC1\delta\ZZ$.
Therefore it is generated by $a_k/\delta$ (where $a_k$ is the gcd of the numerators of the \gtrsz). 
\\
Finally, as $1=\pi_k(x^{k})$, $a_k$ must divide $\delta$ and we write $\fraC {a_k}\delta=\fraC1{d_k}$.

\emph{3.} Let $k\geq 1$ and $z\in F_k$. If $\pi_k(z)=a/d_k$ (with $a\in\ZZ$) we have $\pi_k(z-ay_k)=0$.
\\ 
So
$z-ay_k\in F_{k-1}$. Thus $F_k=\ZZ y_k\oplus F_{k-1}$ and we conclude by \recu on $k$ that $z\in\bigoplus_{i=0}^{k}\ZZ y_k$.

\emph{4.} We have $y_iy_j\in F_{i+j}$ so $\fraC 1 {d_id_j}=\pi_{i+j}(y_iy_j)\in \fraC 1 {d_{i+j}} \ZZ$. In other words $d_{i+j}$ is a multiple of $d_id_j$.

\emph{5} and \emph{6.} Let us first prove that $d_kF_k\subseteq \ZZ[x]$ by \recu on $k$. 
The base case~\hbox{$k=0$} is clear.
We then use the fact that $xy_{k-1} \in F_k$ \hbox{and $\pi_k(xy_{k-1})=\fraC1{d_{k-1}}$}, therefore  
$$\preskip.0em \postskip.4em\ndsp 
xy_{k-1} = \fraC  {d_k} {d_{k-1}} y_k + w_{k-1} \quad \hbox{with  } w_{k-1} \in F_{k-1}. 
$$
We get $d_ky_k = xd_{k-1}y_{k-1} - d_{k-1}w_{k-1}$ and the right-hand side is in $\ZZ[x]$, by the \hdrz. Therefore $d_ky_k\in\ZZx$ and 
$$\preskip.4em \postskip.4em 
d_kF_k=d_k(\ZZ y_k \oplus F_{k-1})=\ZZ d_ky_k \oplus d_k F_{k-1}\subseteq \ZZx + d_{k-1} F_{k-1}\subseteq \ZZx. 
$$
We have defined $f_k(X)$ monic, of degree $k$ in $\QQ[X]$, by the \egt $f_k(x)=d_ky_k$. 
\\
Since $(1,\dots,x^{n-1})$ is as much a $\ZZ$-basis of $\ZZx$ as a $\QQ$-basis of $\QQ[x]$, and since $d_ky_k\in\ZZ[X]$, we obtain {$f_k\in\ZZ[X]$}. 

The rest follows easily.


\prob{exoPolynomialAutomorphism}
\emph{1.} If $F(G) = X$, we have $\JJ(F)(0) \circ \JJ(G)(0) = \I_{\Ae n}$. \\
As $\JJ(G)(0)$ is invertible, we apply the result to $G$.
We have $H \in \gS^n$ \hbox{with  $G(H) = X$}. Then $F=F\circ G\circ H = H$. Therefore $F$, $G$ are inverses of each other (as transformations of $\gS^n$).

\emph{2.}
Immediate. We can a posteriori verify $\Phi(\gS^n) \subseteq \gS^n$ as well as the \eqvcz
$$\preskip.3em \postskip.3em 
\Phi(G) = G \iff F(G) = X. 
$$

\emph{3.}
We write $F(X) = J_0 \cdot X + F_2(X)$, where the vector $F_2(X)$ is of degree $\ge 2$ in~$X$. Then, $J_0^{-1} \cdot \big(F(G) - F(H)\big) = G - H + J_0^{-1} \cdot \big(F_2(G) - F_2(H)\big)$.
\\
Then $\Phi(G) - \Phi(H) = -J_0^{-1} \cdot \big(F_2(G) - F_2(H)\big)$. Assume $G_i - H_i \in \fm^d$ ($d \ge 1$), and let us prove that each component of $\Phi(G) - \Phi(H)$ belongs to $\fm^{d+1}$. The result will be the desired in\egtz. Such a component is an $\gA$-\lin combination of $G^\alpha - H^\alpha$ with $\alpha \in \NN^n$ and $|\alpha| \ge 2$. To simplify the notation, let $n = 3$ and write

\snuc {
G^{\alpha} - H^{\alpha} =
(G_1^{\alpha_1} - H_1^{\alpha_1}) G_2^{\alpha_2} G_3^{\alpha_3} +
(G_2^{\alpha_2} - H_2^{\alpha_2}) H_1^{\alpha_1} G_3^{\alpha_3} +
(G_3^{\alpha_3} - H_3^{\alpha_3}) H_1^{\alpha_1} H_2^{\alpha_2}.
}

Since the $H_i$'s, $G_i$'s are constant-free, we have  $G^\alpha - H^\alpha \in \fm^{d+1}$, except perhaps for
$(\alpha_2, \alpha_3) = (0,0)$ or $(\alpha_1, \alpha_3) = (0,0)$ or $(\alpha_1, \alpha_3) = (0,0)$.  It remains to look at the special cases, for example $\alpha_2 = \alpha_3 = 0$. In this case, since $\alpha_1 - 1 \ge 1$,

\snic {
G^\alpha - H^\alpha =
G_1^{\alpha_1} - H_1^{\alpha_1} =
(G_1 - H_1) \sum_{i+j = \alpha_1-1} G_1^i H_1^j \in \fm^{d+1}.
}

We have therefore established $d\big(\Phi(G), \Phi(H)\big) \le d(G, H)/2$. This guarantees in particular that  there exists at most one fixed point of $\Phi$. 
Let $G^{(0)}\in \gS^n$, for example $G^{(0)} = 0$, and the sequence $G^{(d)}$ defined by \recu by means \hbox{of $G^{(d+1)} = \Phi(G^{(d)})$}. \\
For $d \ge 1$, each component of $G^{(d)} - G^{(d-1)}$ is in $\fm^d$, which allows us to define $G \in \gS^n$ by $G = \sum_{d \ge 1} \big(G^{(d)} - G^{(d-1)}\big).$
\\
Then, $G$ is the limit of the $G^{(d)}$ for $d \mapsto \infty$, it is a fixed  point of $\Phi$, \cad $F(G) = X$.

\emph{4.}
Assume $G(F) = X$, so $G\big(F(0)\big) = 0$. \\
Let $\wi F = F - F(0)$, $\wi G = G\big(X + F(0)\big)$. Then, $\wi F(0) = \wi G(0) = 0$ \hbox{and $\wi G(\wi F) = X$}. Hence $\wi F(\wi G) = X$, then $F(G) = X$.

\emph{5.}
Check in both cases that $\J(F) = 1$. For the first, we obtain $G$
(of same maximum degree as $F$) by iterating $\Phi$ four times:

\snic {
G = (-X^2Z^3 - 2XY^2Z^2 + 2XYZ + X - Y^4Z + 2Y^3,\
    -XZ^2 - Y^2Z + Y,\  Z).
}

For the second, we obtain $G = (G_1, \ldots, G_5)$ by iterating $\Phi$ four times:

\snic {\arraycolsep2pt
\begin{array} {rclrcl}
G_1 &=& X_1 - 3X_2X_4^2 + 6X_2X_4X_5^3 - 3X_2X_5^6 + 2X_3X_4X_5 - 2X_3X_5^4 +
\\[1mm]
&&X_4^4X_5 - 4X_4^3X_5^4 + 6X_4^2X_5^7 - 4X_4X_5^{10} + X_5^{13},
\\[1mm]
G_2&=&X_2 - X_4^2X_5 + 2X_4X_5^4 - X_5^7,
\\[1mm]
G_3&=&X_3 - X_4^3 + 3X_4^2X_5^3 - 3X_4X_5^6 + X_5^9,
\\[1mm]
G_4&=&X_4 - X_5^3,
\quad \quad 
G_5\;=\;X_4.
\end{array}
}

Note that the maximum degree of $G$ is $13$ whereas that of $F$ is $3$.


\Biblio

\vspace{8pt}
The proof of the \DKM lemma~\ref{lemdArtin} on \paref{lemdArtin} is taken from Northcott~\cite{Nor2} (he attributes it to Artin).

\KROz's \thref{thKro} on \paref{thKro} is found in \cite[Kronecker]{Kro1}. It is \egmt proven by Dedekind~\cite{Ddk1} and Mertens~\cite{Mer}.

Concerning the resultants and subresultants in one variable, a reference work is \cite{AJ}.
However, we regret the lack of a bibliography. Even if the results are either very old or completely new, we do not see the use of hiding the exact sources.
Another important book for  \algq questions on the subject is \cite{BPR}.

The construction of an abtract \cdr for a \spl \pol given in \Thref{thResolUniv} is (almost exactly) that described by Jules Drach in \cite{Drach}, 
which also seems to be where the universal splitting algebra as a fundamental tool for studying algebraic extensions of fields was introduced.


The telegraphical \dem of \Thref{thIntClosStab} was suggested to us by Thierry Coquand.

The Kronecker approach regarding the theory of \ids of number fields is the subject of a historical survey in~\cite[Fontana\&Loper]{FL}.

The \dem of the \nst given in Section \ref{secChap3Nst} is inspired by the one in \cite{BPR}, itself inspired by a van der Waerden \demz.

\newpage \thispagestyle{CMcadreseul}
\incrementeexosetprob


\chapter{\Fp modules}
\label{chap mpf} 
\perso{compil\'e le \today}
\minitoc

\subsection*{Introduction}
\addcontentsline{toc}{section}{Introduction}

Over a \ri the \mpfs play a similar role as that of the finite dimensional \evcs  over a field: the theory of \mpfs is a slightly more abstract, and at times more profitable, way to approach the subject of \slisz.

\smallskip In the first sections of the chapter, 
we provide the basics of the theory of  \mpfsz.

In Section~\ref{secBézout}, we treat the example of  \mpfs over PIDs, and in Section~\ref{secKrull0dim} that of  \mpfs over \zeds \risz.

Finally, Section~\ref{sec Fitt} is dedicated to important invariants that are  \idfsz, and Section~\ref{subsecIdealResultant} introduces the resultant \id as a direct application of the \idfsz.

\newpage\section{Definition, changing \sgrz} \label{sec pf chg} 

A \ixy{finitely presented}{module} is an \Amo $M$ given by a finite number of \gtrs and  relations. Therefore it is a module with a finite \sgr having a \tf syzygy module.
Equivalently, it is a module~$M$ \isoc to the cokernel of a \ali
$$
\gamma:  \Ae m\longrightarrow  \gA^q.
$$
The matrix $G\in \gA^{q\times m}$ of $\gamma$ has as its columns a \sgr of the syzygy  module  between the \gtrs $g_i$ which are the images of the canonical base of $\gA^q$ by the surjection $\pi: \gA^q\rightarrow M$.
Such a matrix is called a
\index{matrix!presentation ---}
\emph{\mpn of the module $M$ for the \sgr $(\gq)$}.
This translates into
\begin{itemize}
\item  $[\,g_1\;\cdots\;g_q\,]\,G=0$, and
\item  every syzygy between the $g_i$'s is a \coli of the columns of~$G$, i.e.: if $[\,g_1\;\cdots\;g_q\,]\,C=0$ with $C\in \gA^{q\times 1}$, there exists a $C'\in \Ae {m\times 1}$  such that~$C=G\,C'$.
\end{itemize}

\smallskip \exls \rdb 
1) A free module of rank $k$ is a \mpf presented by a matrix column formed of $k$ zeros.\footnote{If we consider that a matrix is given by two integers $q,m\geq 0$ and a family of \elts of the \ri indexed by the pairs $(i,j)$ with $i\in\lrbq,\; j \in \lrbm$,
we can accept an empty matrix of type $k\times 0$, which would be the canonical matrix to present a free module of rank $k$.}
More \gnlt every simple matrix is the \mpn of a free module of finite rank.

 2) \label{exl1pf}
 Recall that a \mptf is a module $P$ \isoc to the image of a \mprn $F\in\Mn(\gA)$ for a specific integer~$n$. Since $\Ae n=\Im (F)\oplus \Im(\In-F)$, we obtain $P\simeq \Coker(\In-F)$. This shows that every \mptf is \pfz.

 3) Let $\varphi :V\to V$ be an \endo of a 
finite-dimensional \evc 
 over a discrete field~$\gK$.
Consider $V$ as a $\KX$-module with the following external law
$$
\formule{\gK[X]\times V&\to& V\\[1mm] (P,u)&\mapsto& P\cdot u:=P(\varphi)(u).}
$$
Let $(\un)$ be a basis of $V$ as a \Kev and $A$ be the matrix of $\varphi$ with respect to this basis. Then we can show that a \mpn of $V$ as a $\KX$-module for the \sgr $(\un)$ is the matrix $X\,\In-A$ (see Exercise~\ref{exoAXmodule}).  \eoe

\setcounter{subsection}{-1}
\begin{lemma} \label{factchangesgrmpf}
When we change a finite \sgr for a given \mpfz, the syzygies between the new \gtrs form a \tf module again.
\end{lemma}
\begin{proof}
Suppose that indeed, with $M\simeq \Coker G$,  another \sgr of the \Amo $M$ is $(h_1,\ldots,h_r)$.
We therefore have matrices~$H_1\in\gA^{q\times r}$ and~$H_2\in\Ae {r\times q}$ 
such that

\snic{[\,g_1\;\cdots\;g_q\,]\,H_1=[\,h_1\;\cdots\;h_r\,]$ and
$[\,h_1\;\cdots\;h_r\,]\,H_2=[\,g_1\;\cdots\;g_q\,].}

Then, the syzygy module between the $h_j$'s is generated by the columns of~$H_2G$ and that  of $\I_r-H_2H_1$.
Indeed on the one hand we clearly have 

\snic{[\,h_1\;\cdots\;h_r\,]\,H_2\,G=0\;$ and
 $\;[\,h_1\;\cdots\;h_r\,]\,(\I_r-H_2H_1)=0.}

On the other hand, if we have a syzygy $[\,h_1\;\cdots\;h_r\,]\,C=0$, we deduce $[\,g_1\;\cdots\;g_q\,]\,H_1C=0$, so $H_1C=GC'$ for some column vector $C'$ and
$$
C=\big((\I_r-H_2H_1)+H_2H_1\big)C=(\I_r-H_2H_1)C+H_2GC'=HC'',
$$
where $H=\lst{\I_r-H_2H_1\mid H_2G}$ and $C''=\Cmatrix{1pt}{C\cr C'}$.
\end{proof}

\rdb
This possibility of replacing a \sgr by another while preserving a finite number of relations is an extremely \gnl phenomenon. It applies to every form of \agq structure which can be defined by \gtrs and relations.
For example, 
it applies to those structures for which every axiom is a universal equality.
 Here is how this works (it suffices to verify that the reasoning applies in each case).\label{nouveausgr}

Assume that we have \gtrs $g_1, $\ldots$, g_n$ and relations

\snic{R_1(g_1,\alb\ldots,g_n)$, \ldots,  $R_s (g_1,\ldots,g_n),}

which \gui{present} a structure $M$.
\\
If we have other \gtrs $h_1$, $\ldots$, $h_m$,
we express them in terms of the $g_j$'s in the form $h_i = H_i(g_1,\ldots,g_n)$. Let $S_i(h_i,g_1,\ldots,g_n)$ be this relation.
\\
We similarly express the $g_j$'s in terms of the $h_i$'s $g_j = G_j(h_1,\ldots, h_m)$.
Let~$ T_j(g_j, h_1,\ldots, h_m)$ this relation.

The structure does not change if we replace the \pn
$$\preskip.3em \postskip-.2em 
(g_1,\ldots,g_n\; ; \alb\; R_1, \ldots, R_s) 
$$
 with
$$\preskip-.2em \postskip.3em 
(g_1,\ldots,g_n, h_1,\ldots, h_m\; ; \; R_1,\alb \ldots,\alb R_s,
S_1,\ldots,S_m). 
$$
As the relations $T_j$ are satisfied, they are consequences of the relations $R_1$, $\ldots$, $R_s,$ $S_1$, $\ldots$, $S_m$,
therefore the structure is always the same with the following \pn

\snic{(g_1,\ldots,g_n, h_1,\ldots, h_m\; ; \; R_1, \ldots, R_s, S_1,\ldots,S_m,
T_1,\ldots,T_n).}

Now in each of the relations $R_k$ and $S_\ell,$ we can replace each~$g_j$ with its expression in terms of the $h_i$'s (which is given in~$T_j$) and this still does not change the presented structure.
We obtain
$$\preskip.3em \postskip.4em 
(g_1,\ldots,g_n, h_1,\ldots, h_m\; ; \; R'_1, \ldots, R'_s, S'_1,\ldots,S'_m,
T_1,\ldots,T_n). 
$$
Finally, if we subtract the pairs $(g_j;T_j)$ one-by-one, it is clear that the structure will still remain unchanged, so we obtain the finite \pn

\snic{ (h_1,\ldots, h_m\; ; \; R'_1, \ldots, R'_s, S'_1,\ldots,S'_m).}

\medskip 
In the case of finitely presented modules this reasoning can be expressed in matrix form.

First of all we note that we do not change the structure of $M$ when we subject the \mpn $G$ to one of the following transformations.\rdb

\begin{enumerate}\label{ManipMpns}
\item Adding a null column (this does not change the syzygy module between fixed \gtrsz).
\item Deleting a null column, except to obtain an empty matrix.
\item Replacing $G$, of type $q\times m$, with $G'$ of type $(q+1)\times (m+1)$ obtained from $G$ by adding a null row on the bottom then a column to the right with $1$ in the position $(q+1,m+1)$, (this reduces to adding a vector among the \gtrsz, by indicating its dependence with respect to the previous \gtrsz)

\snic{G\;\mapsto \;G'\;=\;
\cmatrix{
    G      &C           \cr
    0_{1,m}&1
}.}
\item 
The inverse of the previous operation, except in the case of an empty matrix.
\item Adding to a column a \coli of the other columns (this does not change the syzygy module between fixed \gtrsz).
\item Adding to a row a \coli of the other rows,
(for example if we let $L_i$ be the $i^{\rm th}$ row, replacing $L_1$ with $L_1+\gamma L_2$ reduces to replacing the \gtr $g_2$ with $g_2-\gamma g_1$).
\item Permuting columns or rows.
\end{enumerate}

We then see that if $G$ and $H$ are two \pn matrices of the same module~$M$, we can pass from one to the other by means of the transformations described above. Slightly better: we see that for every finite \sgr of~$M$, we can construct from $G$, by using these transformations, a \mpn of~$M$ for the new \sgrz. Note that consequently, a change of basis of $\gA^q$ or $\Ae m,$ which corresponds to the multiplication of $G$ (either on the left or right) by an \iv matrix,  can be realized by the operations previously described.
\\
More precisely, we obtain the following result.

\pagebreak

\begin{lemma}
\label{lem pres equiv} 
Let $G\in \gA^{q\times m}$ and $H\in \Ae {r\times n}$ be two matrices.
Then \propeq
\begin{enumerate}
\item  The matrices $G$ and $H$ present \gui{the same} module, \cad their cokernels are \isocz.
\item  The two matrices of the figure below are
\elrt \eqvesz.\perso{Ne devrait-on pas mettre en
corollaire une bonne variante du lemme of Schanuel?}
\item  The two matrices of the figure below are \eqvesz.
\end{enumerate}
\vspace{-1em}
\begin{figure}[htbp]
{\small
\begin{center}
\[
\begin{array}{c|p{35pt}|p{15pt}|p{25pt}|p{20pt}|}
\multicolumn{1}{c}{} & \multicolumn{1}{c}{m} & \multicolumn{1}{c}{r} &
\multicolumn{1}{c}{q} & \multicolumn{1}{c}{n} \\
\cline{2-5}
\vrule height20pt depth13pt width0pt q\; & \hfil G \hfil &\hfil
0\hfil &\hfil 0\hfil &\hfil 0\hfil \\
\cline{2-5}
\vrule height15pt depth8pt width0pt r\; & \hfil 0  &\hfil
${\rm I}_{r}$\hfil & \hfil 0\hfil &\hfil0\hfil \\
\cline{2-5}
\end{array}
\]
\vspace{6pt}
\[
\begin{array}{c|p{35pt}|p{15pt}|p{25pt}|p{20pt}|}
\cline{2-5}
\vrule height20pt depth13pt width0pt q\; & \hfil 0 \hfil&\hfil
0\hfil &\hfil $\I_{q}$\hfil &\hfil 0\hfil \\
\cline{2-5}
\vrule height15pt depth8pt width0pt r\; &\hfil 0 \hfil&\hfil
0\hfil & \hfil 0\hfil & \hfil H\hfil \\
\cline{2-5}
\end{array}\]

\end{center}}
\caption{\emph{The two matrices}}
\label{fig}
\end{figure}
\end{lemma}

\vspace{-4pt}
As a first consequence of Lemma \ref{factchangesgrmpf} we obtain a more abstract reformulation of  \cohc as follows.

\begin{fact}
\label{factCohFA}
A \ri is \coh \ssi every \itf is \pf (as \Amoz).
An \Amo is \coh \ssi every \tf submodule is \pfz.
\end{fact}

\subsec{A digression on the \agq computation} \label{DigCalcAlg}
Besides their direct relationship to solving systems of linear equations another reason for the importance of  \mpfs is the following.

Each time an \agq computation reaches an \gui{interesting result} in an \Amo $M$ this computation has only involved a finite number of \elts $x_1$, \dots, $x_n$ of $M$ and a finite number of syzygies between \hbox{the $x_j$'s}, so that there exist a \mpf $P=\Ae{n}\sur{R}$ and a surjective \ali $\theta:P\to x_1\gA+\cdots+x_n\gA\subseteq M$ which sends \hbox{the $e_j$'s} onto the $x_j$'s. Note that $e_j$ designates the class modulo $R$ of \hbox{the $j^{\rm th}$} vector of the canonical basis of~$\Ae{n}$. It must also be true of the above that the \gui{interesting result}  had already been held in $P$ for the~$e_j$'s.

In a more scholarly language we express this idea as follows.\\
\emph{Every \Amo is a filtering colimit  (or filtering inductive limit) of \pf \Amosz.\label{factLimIndFiltPf}
}
\\
However, this statement requires a more subtle treatment in \comaz, and we therefore only indicate its existence.

\section{\Fp \idsz} \label{secIdPf}

\vspace{3pt}
Consider a \ri $\gA$ and a \sgr $(\an)=(\ua)$ for a \itf $\fa$ of $\gA$.
We are interested in the \Amo structure of $\fa$. 

\subsec{Trivial syzygies 
} \label{secRelTrivSeqReg}

Among the syzygies between the $a_i$'s 
there are
what we call the \emph{trivial syzygies} (or \emph{trivial relators} if we see them as \agq \rdes over $\gk$ when $\gA$ is a \klgz): 

\snic{a_ia_j-a_ja_i=0\;$ for $\;i\neq j.}

If $\fa$ is \pfz, we can always take a \mpn of $\fa$ for the \sgr $(\ua)$ in the form

\snic{W=[\,R_\ua\mid U\,],}

where $R_\ua$ is \gui{the} ${n\times  n(n-1)/2}$ \emph{matrix of trivial syzygies} (the order of the columns is without importance). For example, for $n = 4$%
\index{syzygy!trivial ---}%
\index{matrix!of trivial syzygies}
$$\preskip.4em \postskip.4em 
R_{\ua} = \cmatrix {
a_2 & a_3  &  0  & a_4  &  0  &   0\cr
-a_1&   0 &  a_3 &   0 &  a_4 &   0\cr
 0 & -a_1 & -a_2 &   0 &   0 &  a_4\cr
 0 &   0 &   0 & -a_1 & -a_2 & -a_3}. 
$$
 
\begin{lemma}\label{lemDnRz} \emph{(Determinantal \ids of the matrix of trivial syzygies)}
Using the above notations, we have the following results.
\begin{enumerate}
\item
$\cD_n(R_\ua) = \{0\}$.
\item If $1\leq r <n$, then
$\cD_r(R_\ua) = \fa^r$ and 
$$\preskip.3em \postskip.2em 
\fa^r+\cD_r(U)\;\subseteq\; \cD_r(W)  \;\subseteq\;\fa+ \cD_r(U). 
$$
In particular, we have the \eqvc
$$\preskip.3em \postskip.2em 
1\in\cD_{\gA,r}(W)\iff 1\in\cD_{\gA/\!\fa,r}(\ov U) \quad\hbox{ where }\;\ov U=U\mod\fa.
$$
\item
$\cD_n(W) = \cD_n(U)$.
\end{enumerate}
 
\end{lemma}
%
\begin{proof}
\emph {1.} These are \idas and we can take for $a_1$, \dots,~$a_n$ \idtrs over $\ZZ$. 
Since $[\,a_1\;\cdots\;a_n\,] \cdot R_\ua = 0$, we obtain the \egtz~\hbox{$\cD_n(R_\ua) \,[\,a_1\;\cdots\;a_n\,]  = 0$}. The result follows since $a_1$ is \ndzz.

\emph {2.} The inclusion $\cD_r(R_\ua)\subseteq\fa^r$ is obvious for all $r\geq0$. 
For the reverse inclusion, let us take for example $r=4$ and $n\geq5$ and show that
$$\preskip.3em \postskip.4em
\so{a_1^4,\, a_1^3a_2,\, a_1^2a_2^2,\, a_1^2a_2a_3,\, a_1a_2a_3a_4}\subseteq\cD_4(R_\ua).
$$
It suffices to consider the matrices 
 below (we have deleted the $0$'s and replaced~$\pm a_i$ with~$i$ to clarify the structure) extracted from $R_\ua$, and the minors extracted on the last $4$ rows.
\[\preskip.4em \postskip.4em 
\begin{array}{ccccc} 
\cmatrix{ 2 &  3  &  4  &  5  \cr
         1 &  &  &  \cr  
        &  1 &   &  \cr 
        &  &  1 &   \cr 
        &  &  &  1},\,  
&   
\cmatrix{ 2 &  3  & 4   &   \cr
         1 &  &  &  5 \cr  
        &  1 &   &  \cr 
        &  &  1 &   \cr 
        &  &  &  2},\,  
& 
\cmatrix{ 2 &  3  &   &   \cr
         1 &  &  4 &  5 \cr  
             &  1 &   &  \cr 
             &      &  2 &   \cr 
             &      &      &  2}, \, 
\\[10mm]   
\cmatrix{ 2 &  3  &   &   \cr
         1 &  &  4 &  \cr  
             &  1 &   &  5 \cr 
             &      &  2 &   \cr 
             &      &      &  3},  
&   
\cmatrix{ 2 &   &   &   \cr
         1 &  3 &  &  \cr  
        &  2 &  4  &  \cr 
        &  &  3 &  5  \cr 
        &  &  &  4}.  
\end{array}
\]

The inclusion $\fa^r+\cD_r(U)\subseteq \cD_r(W)$ results from $\cD_r(R_\ua)+\cD_r(U)\subseteq \cD_r(W)$ and from the \egt $\cD_r(R_\ua) = \fa^r$. The inclusion $\cD_r(W)  \subseteq\fa+ \cD_r(U)$ is \imdez. 
Finally, the final \eqvc results from the previous inclusions and from the \egt 
$$\preskip.3em \postskip.3em
\cD_{\gA/\!\fa,r}(\ov U)=\pi_{\gA,\fa}^{-1}\big(\fa+\cD_r(U)\big).
$$

\emph {3.} We must show that if a matrix $A\in\Mn(\gA)$ extracted from $W$ contains a column in $R_\ua$, then $\det A=0$.
Take for example the first column of $A$ equal to the first column of $R_\ua$, $\tra{[\,a_2\;-a_1\;0\;\cdots\;0\,]}$. When $z_i=a_i$, Lemma~\ref{lem acheval} below implies 
$\det A=0$, because the $s_j$'s are null.
\end{proof}

\rdb\label{NOTAProdScal}
Recall that $A_{\alpha,\, \beta}$ is the submatrix of $A$ extracted on the rows $\alpha$ and the columns $\beta$.
Let us also introduce the notation for a \gui{scalar product}

\snic{\scp{x}{y}\eqdefi\sum_{i=1}^n x_iy_i}

for two column vectors $x$ and $y$.

\begin{lemma} \label{lem acheval} 
Let $A 
\in \Mn(\gA)$,  $A_j=A_{1..n,j}$, and $z = \tra{[\,z_1\; \cdots\; z_n\,]}\in \Ae{n\times 1}$
with $A_1 = \tra{[\,z_2\; -z_1\; 0\;\cdots\;0\,]}$.
By letting $s_j = \scp {z} {A_j}$ for $j\in \lrb {2..n}$, we have
$$\preskip-.4em \postskip.4em 
\det A  = \som_{j=2}^n (-1)^j \,s_j \,\det(A_{3..n,\, 2..n \setminus \{j\}}). 
$$
In particular, $\det A  \in \gen {s_2, \ldots, s_n}$.

\end{lemma}
\begin{proof}
Let $B=A_{3..n,2..n}$, $B_j=A_{3..n,j}$ and $B_{\hat \jmath}=A_{3..n,\, 2..n \setminus \{j\}}$.  The Laplace expansion of the \deter of $A$ according to the two first rows gives the \egtz:
$$\preskip.2em \postskip.2em \ndsp
\mathrigid1mu 
\det A  = \sum\limits_{j=2}^n (-1)^j \,\dmatrix{\phantom-z_2&a_{1j}\cr -z_1&a_{2j}} \,\det(B_{\hat \jmath}) 
= \sum\limits_{j=2}^n (-1)^j \,(z_1a_{1j}+z_2a_{2j}) \,\det(B_{\hat \jmath}). 
$$
 The gap between this \egt and the desired \egt is

\snic{
\qquad \qquad 
\som_{j=2}^n (-1)^j \,(z_3a_{3j}+\cdots+z_na_{nj}) \,\det(B_{\hat \jmath}).\qquad \qquad(*)
}

Cramer's syzygies between the columns of a matrix with $m=n_2$ gives for~$B$ the \egts  
$$\preskip.2em \postskip.2em\ndsp 
\som_{j=2}^n (-1)^j \det(B_{\hat \jmath})\, B_j=0, \hbox{ a fortiori } \som_{j=2}^n (-1)^j \scp{y}{B_j}\,\det(B_{\hat \jmath})=0,
 $$ for any vector $y\in\Ae{(n-2)\times 1}$.  
By taking $y=\tra{[\,z_3\; \cdots\; z_n\,]}$, we see that the gap $(*)$ is null.
\end{proof}
%

\subsec{Regular sequences}

\begin{definition}
\label{defSeqReg}  
A sequence $(a_1, \ldots, a_k)$ in a \ri $\gA$ is 
\ixc{regular}{sequence}
if each~$a_i$ is \ndz in the \riz~
$\aqo{\gA}{a_j\,;\,j<i}$.%
\index{sequence!regular ---}
\end{definition}

\rem Here we have kept Bourbaki's \dfnz. Most authors also require that the \id $\gen{a_1,\ldots ,a_k}$ does not contain $1$. 
\eoe

\smallskip As a first example, for every \ri $\gk$, the sequence $(X_1, \ldots , X_k)$  is regular in $\gk[X_1,\ldots ,X_k]$.

 Our goal is to show that an \id generated by a regular sequence is a \pf module.

 We first establish a small lemma and a proposition.
\\
Recall that a matrix $M = (m_{ij}) \in \Mn(\gA)$ is said to be
 \emph{alternating} if it is the matrix of an alternating bi\lin form,
i.e. $m_{ii} = 0$ and $m_{ij} + m_{ji} = 0$ for~$i$,~$j\in\lrbn$.%
\index{matrix!alternating ---}\index{alternating!matrix}

The \Amo of  alternating matrices 
is free and of rank  ${n(n-1)}\over 2$ and admits a natural basis.
For example, for $n = 3$,
$$
\mathrigid1mu 
\crmatrix {0 & a & b \cr -a & 0 & c \cr -b & -c & 0 \cr} =
a\crmatrix {0 & 1 &  0 \cr -1 & 0 & 0 \cr 0 & 0 & 0 \cr} +
b\crmatrix {0 & 0 & 1 \cr 0 & 0 & 0 \cr -1 & 0 & 0 \cr} +
c\crmatrix {0 & 0 & 0 \cr 0 & 0 & 1 \cr 0 & -1 & 0 \cr}.
$$

\medskip 
\begin{lemma}\label{PetitLemmeAlterne}
Let $a = \tra{[\,\ua\,]} = \tra{[\,a_1\;\cdots\;a_n\,]} \in \Ae {n\times 1}$.
\begin {enumerate}
\item
Let $M  \in \Mn(\gA)$ be an alternating matrix;
we have $\scp {Ma}{a} = 0$.

\item
A $u \in \Ae {n\times 1}$ is in $\Im R_\ua$ \ssi there exists an alternating matrix  $M\in \Mn(\gA)$ such that $u = Ma$.
\end{enumerate}
\end{lemma}
\begin{proof}\emph{1.} Indeed, $\scp {Ma}{a} = \varphi(a,a)$, where $\varphi$
is an alternating bi\lin form.

\emph{2.} For example, for the first column of $R_\ua$ with $n=4$, we have
$$\preskip.4em \postskip.4em 
\crmatrix {0 & 1 & 0 & 0\cr -1 & 0 & 0 & 0\cr 0 & 0 & 0 & 0\cr 0 & 0 & 0 & 0\cr}
\crmatrix {a_1\cr a_2\cr a_3\cr a_4\cr} =
\cmatrix {a_2\cr -a_1\cr 0\cr 0\cr}, 
$$
and the ${n(n-1)}\over 2$ columns of $R_\ua$ thus correspond to ${n(n-1)}\over 2$ alternating matrices forming the natural basis of the \Amo of  alternating matrices of~$\Mn(\gA)$.
\end {proof}
\pagebreak

\begin{proposition}\label{prop1sregpf}
Let $(z_1, \ldots, z_n)=(\uz)$ be a regular sequence of \elts of~$\gA$ and $z = \tra{ [\,z_1 \;\cdots\; z_n\,]} \in \Ae{n\times 1}$. If $\scp {u}{z} = 0$, there exists an alternating matrix~$M \in \Mn(\gA)$ such that $u = Mz$, and therefore $u\in\Im R_\uz$. 
\end{proposition}
\begin{proof}
We reason by \recu on $n$. 
For $n = 2$, we start from  $u_1z_1 + u_2z_2 =0$. Therefore $u_2z_2 = 0$ in $\aqo{\gA}{z_1}$, and since $z_2$ is \ndz modulo $z_1$, we~have~$u_2 = 0$ in $\aqo{\gA}{z_1}$, say $u_2 = -az_1$ in $\gA$. We get $u_1z_1 -az_2z_1 = 0$, and as~$z_1$ is \ndzz, $u_1 = az_2$, which is written as $\cmatrix {u_1\cr u_2\cr} = \crmatrix {0 & a \cr -a & 0\cr} \cmatrix {z_1\cr z_2\cr}$.

For $n+1$ ($n\geq2$), we start from $u_1z_1 + \cdots + u_{n+1}z_{n+1} = 0$. By using the fact that $z_{n+1}$ is \ndz modulo $\gen{\zn}$,
we obtain $u_{n+1} \in \gen{\zn}$, which we write as $a_1z_1 + \cdots + a_{n}z_{n} + u_{n+1} = 0$. Whence

\snic{
(u_1 - a_1z_{n+1}) z_1 + \cdots + (u_{n} - a_{n}z_{n+1}) z_{n} = 0.}

By \hdrz, we know how to construct an alternating matrix $M \in \MM_{n}(\gA)$  with
$$\preskip.2em \postskip.4em
\mathrigid1mu 
\cmatrix {u_1 - a_1z_{n+1}\cr \vdots\cr u_{n} - a_{n} z_{n+1}\cr} =
M \cmatrix {z_1\cr \vdots\cr z_{n}\cr},
\hbox { i.e. }
\cmatrix {u_1\cr \vdots\cr u_{n}\cr} =
M \cmatrix {z_1\cr \vdots\cr z_{n}\cr} +
z_{n+1}\cmatrix {a_1\cr \vdots\cr a_{n}\cr},
$$
and we obtain the desired result
$$\preskip.4em \postskip.2em
\cmatrix {u_1\cr \vdots\cr u_{n}\cr u_{n+1}\cr} =
\cmatrix {
     &       &          & a_1  \cr
     & M     &          & \vdots \cr
     &       &          & a_{n} \cr
-a_1 & \dots & -a_{n} & 0 \cr}
\cmatrix {z_1\cr \vdots\cr z_{n}\cr z_{n+1}\cr}.
$$

\vspace{-10pt}
\end{proof}

\begin{theorem}\label{propsregpf}
If $(z_1, \ldots, z_n)$ is a regular sequence of \elts of $\gA$, the \id $\gen{\zn}$ is a \pf \Amoz.
More \prmtz, we have the exact sequence

\snic{\Ae {n(n-1)/2} \vvvers {R_\uz}  \Ae n
\vvvvvers {(z_1, \ldots, z_n)}  \gen{\zn} \lora 0.}
\end{theorem}

\rem
The objects defined above constitute an introduction to the first degree of the \emph{Koszul complex} of $(\zn)$.
\eoe

\begin{proof}
This results from Proposition~\ref{prop1sregpf} and from Lemma~\ref{PetitLemmeAlterne}.
\end{proof}
%

\subsec{A \gmt example} \label{secExempleGeo}

Let us begin with a most useful and obvious fact. 
\begin{propdef}\label{prdfCaracAlg}
\emph{(Characters of an \algz)} 
\\
 Let $\imath:\gk\to\gA$ be an \algz.
\begin{itemize}
\item A \homo of \klgs $\varphi:\gA\to\gk$ is called a 
\ixc{character}{of an \algz}.
\item  If $\gA$ has a character $\varphi$, then $\varphi\circ \imath=\Id_\gk$, $\imath\circ \varphi$ is a \prr and $\gA=\gk.1_\gA\oplus\Ker\varphi$. 
In particular,  $\gk$ may be identified with $\gk.1_\gA$. 
\end{itemize}
\end{propdef}
\facile

Now let $(\uf)=(f_1, \ldots, f_s)$ be a \syp over a \riz~$\gk$,
with each $f_i \in \kuX
= \kXn$.  We let 
$$\gA = \kxn = \aqo{\kuX}{\uf}.
$$

In this subsection, we will formally say that $\gA$ is the \emph{\ri of 
the affine \vrt $\uf = \uze$}. 

For the \alg $\gA$, the characters $\varphi:\gA\to\gk$ are given by the zeros in~$\gk^n$ of the \syp $(\lfs)$ 
$$\preskip.2em \postskip.4em
(\uxi)=(\xin)=\big(\varphi(x_1),\ldots,\varphi(x_n)\big),\quad\uf(\uxi)=\uze.
$$
In this case, we say that $(\uxi)\in\gk^ n$ is a point of the \vrt $\uf = \uze$.

The \id 
$$\preskip-.20em \postskip.4em
\fm_\uxi\eqdefi\gen{x_1-\xi_1,\ldots,x_n-\xi_n}_\gA
$$ 
is called the \emph{\id of the point $(\uxi)$ in the \vrtz}. 
We then have as a special case of Proposition~\ref{prdfCaracAlg}: $\gA=\gk\oplus\fm_\uxi$, with $\fm_\uxi=\Ker\varphi$.
\index{ideal!of a point}

{In this subsection we show that the \id $\fm_\uxi$ is a \pf \Amo by making a \mpn for the \sgr $(x_1-\xi_1, \ldots, x_n-\xi_n)$ explicit.}

 By translation, it suffices to treat the case where $\uxi=\uze$, which we assume henceforth.
\\
 The simplest case, that for which there is no \eqnz, has already been treated in  \Thref{propsregpf}.
\\
Let us observe that every $f \in \kuX$ such that $f(\uze) = 0$ is written, in many ways, in the form

\snic{
f = X_1u_1 + \cdots + X_n u_n, \qquad  u_i \in \kuX.}

If $X_1v_1 + \cdots + X_nv_n$ is another expression of $f$, we obtain by subtraction a syzygy between the $X_i$'s in $\kuX$, and so

\snic{
\tra{ [\, v_1\;\cdots\; v_n\,]} - \tra{[\, u_1\;\cdots\; u_n\,]} \in \Im {R_{\uX}}.}


For the \syp $(f_1, \ldots, f_s)$, we thus define (in a non-unique manner) a family of \pols $(u_{ij})_{i \in\lrbn, j \in\lrbs}$, with $f_j = \sum_{i=1}^n X_iu_{ij}$. This gives a matrix $U(\uX) = (u_{ij})$ and its image $U(\ux) = \big(u_{ij}(\ux)\big) \in \Ae{n\times s}$.

\begin{theorem} \label{thidptva}
For a \syp over a \ri $\gk$ and a zero~$(\uxi)\in\gk^n$, the \id $\fm_\uxi$ of the point $(\uxi)$ is a \pf \Amoz.%
\index{polynomial system}
\\
More \prmtz, with the previous notations, for the $\uxi=\uze$ case the matrix $W = [\,R_\ux\,|\,U(\ux)\,]
$ is a \mpn of the \id $\fm_\uze$ for the \sgr $(\xn)$. In other words we have an exact sequence
$$\preskip.2em \postskip.4em \ndsp
\Ae{m} \vvvvvers {[\,R_\ux\,|\,U\,]}
 \Ae n
\vvvvvers {(x_1, \ldots, x_n)}  \fm_\uze\lora 0 \quad\quad(m = {n(n-1) \over 2} + s). 
$$
\end{theorem}
\begin{proof}
 Take for example $n = 3$, $s = 4$, $X = \tra {[\,X_1 \;X_2\;X_3 \,]}$
and to save on indices let us write $f_1 = X_1a_1 + X_2a_2 + X_3a_3$, and $f_2$, $f_3$, $f_4$ by using the letters $b$, $c$, $d$. We claim to have the following \mpn for the \sgr $(x_1,x_2,x_3)$ of $\fm_\uze$
$$\preskip.3em \postskip.3em\ndsp
\cmatrix {
x_2  & x_3  & 0    & a_1(\ux) & b_1(\ux) & c_1(\ux) & d_1(\ux)\cr
-x_1 & 0    & x_3  & a_2(\ux) & b_2(\ux) & c_2(\ux) & d_2(\ux)\cr
0    & -x_1 & -x_2 & a_3(\ux) & b_3(\ux) & c_3(\ux) & d_3(\ux)\cr}.
$$
We define $A = \tra {[\,a_1 \;a_2\;a_3 \,]}$ in $\kuX^3$
(as well as $B,C,D$) so that

\snic{f_1 =\scp {{A}}{X},\; f_2 = \scp {{B}}{X} \;\;\ldots.}

Let $v_1(\ux)x_1 + v_2(\ux)x_2 + v_3(\ux)x_3 =0$ be a syzygy in~$\gA$. 
We lift it in $\kuX$
$$\preskip.3em \postskip.0em 
v_1X_1 + v_2X_2 + v_3X_3 \equiv 0 \mod \gen{\uf}, 
$$
which we write

 \snic{
v_1X_1 + v_2X_2 + v_3X_3 = \alpha f_1 + \beta f_2 + \gamma f_3 + \delta f_4,
\qquad \alpha,\, \beta,\, \gamma,\, \delta \in \kuX.
}

Therefore, with $V = \tra {[\,v_1 \;v_2\;v_3 \,]}$,
$V - (\alpha{A} + \beta{B} + \gamma{C} + \delta{D})$
is a syzygy for $(X_1, X_2, X_3)$, which implies by Proposition~\ref{prop1sregpf}

 \snic{
V - (\alpha{A} + \beta{B} + \gamma{C} + \delta{D}) \in
\Im {R_{\uX}}.
}

Thus, $V \in \Im \,[\,R_{\uX}\,|\,U(\uX)\,]$, and 
$\tra {[\, v_1(\ux)\;v_2(\ux)\;v_3(\ux)\,] } \in \Im \,[\,R_\ux\,|\,U(\ux)\,]$.
\end{proof}

\section{The category of \mpfs} \label{secCatMpf}

The category of  \mpfs over $\gA$ can be constructed from the category of  free modules of finite rank over $\gA$ by a purely categorical procedure.
\begin{enumerate}\itemsep0pt
\item  A \mpf $M$ is described by a triplet 

\snic{(\rK_M,\alb\rG_M,\rA_M),}

where
$\rA_M$ is a \ali between the free modules of finite ranks~$\rK_M$ and~$\rG_M$.
We have~$M\simeq \Coker\rA_M$ and~$\pi_M:\rG_M\rightarrow M$ is the surjective \ali with kernel $\Im\rA_M$.
The matrix of the \aliz~$\rA_M$ is a \mpn of~$M$.

\item  A \ali $\varphi$ of the module $M$ (described by $(\rK_M,\rG_M,\rA_M)$) to the module $N$ (described by $(\rK_N,\rG_N,\rA_N)$) is described by two \alis $\rK_\varphi:\rK_M\rightarrow \rK_N$ and $\rG_\varphi:\rG_M\rightarrow \rG_N$
subject to the commutation relation $\rG_\varphi\circ\rA_M=\rA_N\circ\rK_\varphi$.
$$\preskip.1em \postskip-.1em 
\xymatrix {
\rK_M \ar[r]^{\rA_M} \ar[d]_{\rK_\varphi} & \rG_M \ar[d]^{\rG_\varphi}
\ar@{->>}[r]^{\pi_M}
                   & M \ar[d]^{\varphi} \\
\rK_N \ar[r]_{\rA_N}                    & \rG_N \ar@{->>}[r]_{\pi_N}
                   & N \\} 
$$
\item  The sum of two \alis $\varphi$ and $\psi$ of $M$ to $N$ represented by $(\rK_\varphi,\rG_\varphi)$ and $(\rK_\psi,\rG_\psi)$ is represented by $(\rK_\varphi+\rK_\psi,\alb \rG_\varphi+\rG_\psi)$.
\\
The \ali $a\varphi$ is represented by $(a\rK_\varphi,\alb a\rG_\varphi)$.

\item  To represent the composite of two \alisz, we compose their representations.

\item  Finally, the \ali $\varphi$ of $M$ to $N$ represented by $(\rK_\varphi,\rG_\varphi)$ is null \ssi there exists a $Z_\varphi:\rG_M\rightarrow \rK_N$ satisfying~$\rA_N\circ Z_\varphi=\rG_\varphi$.
\end{enumerate}

This shows that the \pbs concerning  \mpfs can always be interpreted as \pbs regarding matrices, and are often reduced to
\label{exMpf} \pbs concerning the solution of \slis over $\gA$.
\\
For example, given $M$, $N$ and $\varphi$, if we look for a \ali $\sigma:N\rightarrow M$ satisfying $\varphi\circ\sigma=\Id_N$,
we must find \alis $\rK_\sigma:\rK_N\rightarrow \rK_M$,
$\rG_\sigma:\rG_N\rightarrow \rG_M$  and $Z:\rG_N\rightarrow \rK_N$  
satisfying
$$ 
\rG_\sigma\circ\rA_N=\rA_M\circ\rK_\sigma \quad {\rm and}\quad
\rA_N\circ Z=\rG_\varphi\circ\rG_\sigma-\Id_{\rG_N}.
$$
This is none other than a \sli having as unknowns the \coes of the matrices
of the \alis $\rG_\sigma$, $\rK_\sigma$ and~$Z$.

Analogously,
if we have $\sigma:N\to M$ and if we want to know whether there exists a $\varphi :M\to N$ satisfying $\varphi\circ\sigma=\Id_N$,
We will have to solve a \sli whose unknowns are the \coes of the matrices of the \alis $\rG_\varphi$, $\rK_\varphi$ and~$Z$.

Similarly, if we have $\varphi:M\to N$ and if we want to know whether $\varphi$ is \lnlz, 
we must determine whether there exists a $\sigma:N\to M$ 
satisfying~$\varphi \circ\sigma\circ\varphi =\varphi $, and we obtain a \sli having as its unknowns the \coes of the matrices of $\rG_\sigma$, $\rK_\sigma$ and~$Z$.

We deduce the corresponding \plgsz.

\begin{plcc}
\label{plcc.scinde}  
\emph{(For certain \prts of the \alis between \mpfsz)}\\
Let $S_1$, $\dots$, $S_n$ be \moco of $\gA$ and $\varphi :M\to N$ be a \ali between \mpfsz.
Then \propeq
\begin{enumerate}
\item  The \Ali $\varphi$ admits a left-inverse (resp.\,admits a right-inverse, resp.\,is \lnlz).
\item For $i\in\lrbn$,  
the $\gA_{S_i}$-\ali $\varphi_{S_i}:M_{S_i}\to N_{S_i}$ admits a left-inverse (resp.\,admits a right-inverse, resp.\,is \lnlz).
\end{enumerate}
\end{plcc}

\vspace{-.5em}
\pagebreak	
\section{Stability \prts} \label{secStabPf}

\begin{proposition}
\label{propPfInter}
Let $N_1$ and $N_2$ be two \tf \Asubs  
of an \Amo $M$.
If $N_1+N_2$ is \pfz, then $N_1\cap N_2$ is \tfz.
\end{proposition}
%
\begin{proof}
We can follow almost word for word the \dem of item \emph{1} of \Thref{propCoh4}
(\ncr condition).
\end{proof}
%

\begin{proposition}
\label{propPfSex}
Let $N$ be an \Asub of $M$ and $P=M/N$.
\begin{enumerate}
\item If $M$ is \pf and $N$ \tfz, then $P$ is \pfz.
\item If $M$ is \tf and $P$ \pfz, then $N$ is \tfz.
\item If $P$ and $N$ are \pfz, then $M$ is \pfz. More \prmtz, if $A$ and $B$ are \mpns for $N$ and $P$, we have a \mpn 
$D=\blocs{1}{.8}{.7}{.6}{$A$}{$C$}{$0$}{$B$}$
for $M$.
\end{enumerate}
\end{proposition}

%
\begin{proof}
\emph{1.} We can suppose that $M=\gA^p/F$ with $F$ \tfz. 
If $N$ is \tfz, it is of the form  $N=(F'+F)/F$ where $F'$ is \tfz, so $P\simeq \gA^p/(F+F')$.
\\
\emph{2.} We write $M=\gA^p/F$ and $N=(F'+F)/F$. We have $P\simeq \gA^p/(F'+F)$,  
so $F'+F$  (and also $N$) is \tf
(Section~\ref{sec pf chg}).
\\
\emph{3.} Let $x_1$, \ldots, $x_m$ be \gtrs of $N$ and $x_{m+1}$, \ldots, $x_n$ be \elts of~$M$ whose classes modulo $N$ generate $P$. Every syzygy on $(\ov{x_{m+1}}, \ldots, \ov{x_n})$ in~$P$ gives a syzygy on $(\xn)$ in $M$. Similarly, every syzygy on~$(\xn)$ in $M$ gives a syzygy on $(\ov{x_{m+1}}, \ldots, \ov{x_n})$ in $P$.
 
If $A$ is a \mpn of $N$ for $(\xm)$ and if $B$ is a \mpn of $P$ for  $(\ov{x_{m+1}},\ldots ,\ov{x_n})$, we obtain a \mpn $D$ of $M$ for $(\xn)$ in the desired format.
\end{proof}
Note that in the \dem of item \emph{2} the submodules $F$ and $F'$ are not \ncrt \tfz.

\subsection*{Coherence and finite \pnz}
\addcontentsline{toc}{subsection}{Coherence and finite \pnz}

Propositions~\ref{propCoh1} and~\ref{propCohfd1}
(where we take $\gA$ as the \Amoz~$M$)
can be reread in the form of the following \thoz.

\begin{theorem}
\label{propCoh2}
On a \coh \ri every \mpf is \cohz.
On a \fdi \coh \ri every \mpf is \fdi and \cohz.
\end{theorem}

\vspace{-.5em}
\pagebreak
\begin{proposition}\label{propCohpfKer}
Let $\gA$ be a \coh \ri and $\varphi:M\to N$ be a \ali between \pf \Amosz,
then $\Ker\varphi$, $\Im\varphi$ and~$\Coker\varphi$ are \mpfsz.
\end{proposition}

\begin{proposition}
\label{propCohSex}
Let $N$ be a \tf \Asub of $M$.
\begin{enumerate}
\item If $M$ is \cohz, $M/N$ is \cohz.
\item If $M/N$ and $N$ are \cohsz, $M$ is \cohz.
\end{enumerate}
\end{proposition}
\begin{proof}
\emph{1.} Consider a \tf submodule $P=\gen{\ov{x_1},\ldots
,\ov{x_\ell}}$ of $M/N$. Then $P\simeq (\gen{\xl}+N)/N$. We conclude by Proposition~\ref{propPfSex} that it is \pfz.

\emph{2.} Let $Q$ be a \tf submodule of $M$.
The module $(Q+N)/N$ is \tf in $M/N$ therefore \pfz.
Since $(Q+N)/N$ and $N$ are \pfz, so is $Q+N$  (Proposition~\ref{propPfSex}).
Therefore $Q\cap N$ is \tf (Proposition~\ref{propPfInter}). Since $N$ is \cohz, $Q\cap N$ is \pfz. Since $Q/(Q\cap N)\simeq  (Q+N)/N$ and $Q\cap N$ are \pfz, $Q$ is \pf (Proposition~\ref{propPfSex}).
\end{proof}

\subsect{Tensor product, exterior powers, symmetrical powers}{Tensor product, exterior powers, symmetrical powers}
 \label{ProdTens}

Let $M$ and $N$ be two \Amosz.
A bi\lin map $\varphi :M\times N\to P$ is called a \ix{tensor product} of the \Amos $M$ and $N$ if every bi\lin map $\psi :M\times N\to R$ is uniquely expressible in the form $\psi=\theta\circ\varphi$, where $\theta$ is an \Ali from $P$ to $R$.
\vspace{-.5em}
\Pnv{M\times N}{\varphi}{\psi}{P}{\theta}{R}{}{bi\lin maps}{\alisz.}

\vspace{-1em}
It is then clear that $\varphi :M\times N\to P$ is unique in the categorical sense,
\cad that for every other tensor product $\varphi' :M\times N\to P'$ there is a unique \ali $\theta:P\to P'$ which renders the suitable diagram commutative, and that $\theta$ is an \isoz.

If $(\ug)$ is a \sgr of $M$ and $(\uh)$ a \sgr of~$N$, a bi\lin map $\lambda : M\times N\rightarrow P$ is known from its values over the \elts of $\ug\times \uh$.
Furthermore, the values $\lambda(x,y)$ are linked by certain constraints, which are derived from syzygies between \elts of~$\ug$ in $M$ and from syzygies between \elts of $\uh$ in $N$. \\
For example, if we have a syzygy $a_1x_1+a_2x_2+a_3x_3=_M 0$ between 
\eltsz~$x_i$ of $\ug$, with the~$a_i$'s in~$\gA$, this provides for each $y\in \uh$ the following syzygy in $P$:
$a_1\lambda(x_1,y)+a_2\lambda(x_2,y)+a_3\lambda(x_3,y)=0$.

Actually, \gui{those are the only essential constraints, and that shows that a tensor product can be constructed.}

More \prmtz, let $x \te y$ instead of $(x,y)$ be an arbitrary \elt of~$\ug\times \uh$. 
Consider then the \Amo $P$ generated by $x \te y$ \elts linked by the syzygies described above ($a_1(x_1\te y)+a_2(x_2\te y)+a_3(x_3\te y)=_P 0$ for the given example).

\begin{proposition}
\label{propPftens} 
(With the above notations)
\begin{enumerate}
\item There exists a unique bi\lin map $\varphi : M\times N\to P$
such that for all $(x,y)\in\ug\times \uh$, we have $\varphi(x,y)=x\te y$.
\item 
With this bi\lin map  $P$ is a tensor product 
of the modules~$M$ and~$N$. 
In particular, if $M$ and $N$ are free with bases $(\ug)$ and $(\uh)$, the module~$P$ is
free with basis $(\ug\te\uh):=({x\te y})_{x\in \ug,\,y\in\uh}$.

\end{enumerate}
\end{proposition}
\facile

Thus, the tensor product of two \Amos exists and can always be defined from \pns of these modules. It is denoted by $M\te_\gA N$.

The fact that follows is more or less a paraphrase of the previous proposition, but it can only be stated once we know that tensor products exist.

\begin{fact}\label{factpropPftens}
\begin{enumerate}
\item  If two modules are \tf (resp.\,\pfz) then so is their tensor product.

\item If $M$ is free with basis $(g_i)_{i\in I}$ and $N$ is free with basis $(h_j)_{j\in J}$, then $M\otimes N$ is free with basis $(g_i\otimes h_j)_{(i,j)\in I\times J}$.

\item If $M\simeq \Coker \alpha$ and $N\simeq \Coker \beta $, with $\alpha  :L_1\rightarrow L_2$ and $\beta :L_3\rightarrow L_4$, the modules $L_i$ being free, then the \Ali

\snic{(\alpha \te\Id_{L_4})\oplus (\Id_{L_2}\te\beta ):
(L_1\otimes L_4)\oplus ( L_2\otimes L_3)\to L_2\otimes L_4}

has as its cokernel a tensor product of $M$ and $N$.
\end{enumerate}

\end{fact}

\comms~ \\
1) The theory of \emph{universal algebra} provides profound reasons why the construction of the tensor product \emph{cannot fail to work}.
But this \gnl theory is a little too heavy to be presented in this work, and it is best to soak up these kinds of things by impregnating examples.

\rdb%
2) The reader accustomed to \clama would not have read without apprehension our \gui{presentation} of the tensor product of $M$ and $N$, which is a module constructed from \pns of $M$ and $N$.\label{CommUnivPasUniv}
If they have read Bourbaki, they will have noticed that our construction is the same as that of the illustrious multi-headed \matnz, except that Bourbaki limits himself to one \gui{natural and universal} \pnz: every module is generated by \underline{all} its \elts linked by \underline{all} their syzygies. If Bourbaki's \gui{\pnz} has the merit of universality, it has the inconvenience of the weight of the hippopotamus.

In fact, in \comaz, we do not have the same underlying \gui{set theory} as in \clamaz. Once we have \emph{given} a module $M$ by means of a \pn $\alpha :L_1\to L_2$, we do not rush to forget $\alpha$ as we pretend to do in \clamaz.\footnote{A detailed inspection of the object $M$ constructed according to the set theory of classical mathematics reveals that the latter does not forget it either.}
On the contrary, from a \cof point of view, the module $M$ is nothing other than \gui{an encoding of the \ali $\alpha$} (for example in the form of a matrix if the \pn is finite), with the \sul information that this is the \pn of a module.
Furthermore, a \gui{quotient set} is not seen as a set of \eqvc classes, but as \gui{the same preset equipped with a coarser \egt relation;} the quotient set of $(E,=_E)$ by the \eqvc relation $\sim$ is simply the set $(E,\sim)$.
Consequently, our construction of the tensor product, consistent with its implementation on a machine, is entirely \gui{natural and universal} in the framework of  constructive set theory (the reader can consult the simple and brilliant Chapter 3 of \cite{B67}, or one of the other classic works of reference on \coma~\cite{Be,BB85,BR,MRR}).

3) To construct the tensor product of nondiscrete modules, we a priori need  the notion
of a free module over a nondiscrete set. 
For the constructive \dfn of this kind of free module,
 see Exercise \ref{propfreeplat}. 
We can avoid  this kind of free module in the following way. 
We do not use \sgrs of $M$ and $N$ in the construction. The \elts of the tensor product $M\otimes_\gA N$ are given as formal sums $\sum_{i=1}^{n}x_i\otimes y_i$ for finitely enumerated families in~$M$ and $N$. 
Now the \pb is 
to give a correct \dfn of the \eqvc relation which gives as quotient set the set underlying the module $M\otimes_\gA N$.
The details are left to the reader.
\eoe

By its \dfnz, the tensor product is \gui{functorial,} i.e. if we have two \Alis $f:M\to M'$ and $g:N\to N'$, then there exists a unique \ali $h:M\te_\gA N \to M'\te_\gA N'$ satisfying the \egts $h(x\te y)=f(x)\te g(y)$ for $x\in M$ and $y\in N$.
This \ali is naturally denoted by $h=f\te g$.

We also have the canonical \isos 

\snic{M\te_\gA N\simarrow N\te_\gA M\;\hbox{ and }\;
      M\te_\gA (N\te_\gA P)\simarrow (M\te_\gA N)\te_\gA P,}

which we express by saying that the tensor product is commutative and associative.

The following fact \imdt results from the description of the tensor product by \gtrs and relations.

\begin{fact}\label{factSuitExTens}
For every exact sequence of \Amos $M\vers{f} N\vers{g} P\to 0$ and for every \Amo $Q$
the sequence

\snic{M\te_\gA Q\vvvvers{f\te \Id_Q} N\te_\gA Q\vvvvers{g\te \Id_Q} P\te_\gA Q\to 0}

is exact.
\end{fact}

We express this fact by saying that \gui{the functor $\bullet\otimes Q$ is right exact.}

We will not recall in detail the statement of the universal \pbs that solve the exterior powers (already given \paref{PuissExtMod}), 
\emph{\smq powers}\index{symmetric powers!of a module} and the \emph{exterior \algz}\index{algebra!exterior --- of a module}%
\index{exterior!algebra of a module} of an \Amoz.

Here are however the corresponding \gui{small diagrams} for the last two.

\vspace{-1.2em}
\Pnv{M^{k}}{\gs_\gA^k}{\psi}{\gS_\gA^k M}{\theta}{N}{}{\smq multi\lin maps~~}{\alisz.}

\vspace{-2.6em}
\PNV{M}{ \lambda_\gA}{\psi}{\Vi_\gA\! M}{\theta}{\gB}%
{\Amos}{$\psi(x)\times \psi(x)=0$ for all $x\in M$}%
{associative \Algsz.}

\vspace{-1.8em}
As a corollary of Proposition~\ref{propPftens} we obtain the following proposition.

\begin{proposition}
\label{propPfPex}
If $M$ is a \pf \Amoz, then
the same goes for $\Al k_\gA\! M$
and for the \smq powers  $\gS_\gA^k M$ ($k\in\NN$).\\
More \prmtz, if $M$ is generated by the \sys $(x_1, \dots,x_n)$ subjected to syzygies~\hbox{$r_j\in\Ae n$}, we obtain the following results.
\begin{enumerate}
\item The module $\Al k_\gA \!M$ is generated by the $k$-vectors 

\snic{x_{i_1}\vi\cdots\vi
x_{i_k}$ for $1\leq i_1<\cdots <i_k\leq n,}

subjected to the syzygies obtained by making the exterior product of the $r_j$ syzygies by the $(k-1)$-vectors $x_{i_1}\vi\cdots\vi x_{i_{k-1}}$.
\item The module $\gS_\gA^k M$ is generated by the $k$-\smq tensors

\snic{\gs(x_{i_1},\dots, x_{i_k})$ for $1\leq i_1\leq \cdots \leq i_k\leq n,}

subjected to the syzygies obtained by making the product of the $r_j$ syzygies by the $(k-1)$-\smq tensors 
$\gs(x_{i_1},\dots, x_{i_{k-1}})$.
\end{enumerate}
\end{proposition}

For example, with $n=4$ and $k=2$ a syzygy $a_1x_1+\cdots +a_4x_4=0$
in~$M$ leads to $4$ syzygies in $\Al2_\gA M$
$$\preskip.2em \postskip.4em
\arraycolsep2pt
\begin{array}{rcl}
a_2 \,(x_1\vi x_2) \,+\,a_3 \,(x_1\vi x_3) \,+\,a_4 \,(x_1\vi x_4)&  = &  0    \\
a_1 \,(x_1\vi x_2)  \,-\,a_3 \,(x_2\vi x_3) \,-\,a_4 \,( x_2\vi x_4)&  = &  0 \\
a_1 \,(x_1\vi x_3) \,+\,a_2 \,(x_2\vi x_3) \,-\,a_4 \,(x_3\vi x_4)&  = &  0    \\
a_1 \,(x_1\vi x_4) \,+\,a_2 \,(x_2\vi x_4) \,+\,a_3 \,(x_3\vi x_4)&  = &  0
\end{array}
$$
and to $4$ syzygies in $\gS_\gA^2 M$
$$
\arraycolsep2pt
\begin{array}{rcl}
a_1\,\gs(x_1,x_1) \,+\, a_2\,\gs(x_1,x_2) \,+\,a_3\,\gs(x_1,x_3)
\,+\,a_4\,\gs(x_1,x_4)&  = &  0    \\
a_1\,\gs(x_1,x_2) \,+\, a_2\,\gs(x_2,x_2) \,+\,a_3\,\gs(x_2,x_3)
\,+\,a_4\,\gs(x_2,x_4)&  = &  0    \\
a_1\,\gs(x_1,x_3) \,+\, a_2\,\gs(x_2,x_3) \,+\,a_3\,\gs(x_3,x_3)
\,+\,a_4\,\gs(x_3,x_4)&  = &  0    \\
a_1\,\gs(x_1,x_4) \,+\, a_2\,\gs(x_2,x_4) \,+\,a_3\,\gs(x_3,x_4)
\,+\,a_4\,\gs(x_4,x_4)&  = &  0
\end{array}
$$

\rem
More \gnltz, for every exact sequence

\snic {
K \vers {u} G \vers {p} M \to 0
}

we have an exact sequence
$$
\preskip.0em \postskip.2em 
K \otimes \Al{k-1} G \vers {u'} \Al{k} G \vvvers {\Al{k} p} \Al{k} M \to 0 
$$
with $u'(z \otimes y) = u(z) \vi y$ for $z \in K$, $y \in \Al{k-1} G$. \\ 
On the right-hand side, the surjectivity is immediate and it is clear that $(\Al{k} p) \circ u' = 0$, which allows us to define $p' : \Coker u' \to \Al{k} M$ by passage to the quotient.
It remains to prove that $p'$ is an \isoz. For that, it suffices to construct a \aliz~\hbox{$q' : \Al{k} M \to \Coker u'$} that is the inverse of~$p'$. We do not have a choice: for $x_1$, \ldots, $x_k \in M$ with preimages $y_1$, \ldots, $y_k \in G$ by $p$

\snic {
q'(x_1 \vi \cdots \vi x_k) = y_1 \vi \cdots \vi y_k \bmod \Im u'.
}

We leave it up to the reader to verify that $q'$ is indeed defined and suitable.
\\
The analogous result is valid for the \smq powers.
\eoe

\mni\exl \label{belexemple}
Let $\gB$ be the \ri of \pols $\gA[x,y]$ in the \idtrs $x$ and $y$ over a nontrivial \ri $\gA$. Consider the \id $\fb=\gen{x,y}$ of $\gB$, and look at it as a \Bmo that we denote by $M$.
Then, $M$ admits the \sgr $(x,y)$ for which a \mpn is equal to $\cmatrix{\phantom-y\cr-x}$. Deduce that $M\te_\gB M$ admits $(x\te x,x\te y,y\te x,y\te y)$ as a \sgrz, with a \mpn equal to
$$
\crmatrix{y&0&0&y\cr-x&0&y&0\cr0&y&0&-x\cr0&-x&-x&0}
\leqno{ 
\begin{array}{c} 
x\te x\\ x\te y\\ y\te x\\ y\te y \\[-8.5pt] ~
\end{array}}
$$
We deduce the following annihilators
$$
\preskip.4em \postskip.4em
\begin{array}{c}
\Ann_\gB(x\te y-y\te x)=\fb, \quad \Ann_\gB(x\te y+y\te x)=\Ann_\gA(2)\,\fb,   
\\[1mm]
\Ann_\gB(x\te x)=\Ann_\gB(x\te y)=\Ann_\gB(y\te x)=\Ann_\gB(y\te y)=0.
\end{array}
$$
The dual $M\sta=\Lin_\gB(M,\gB)$ of $M$ is free of rank $1$, generated by the form
\label{NOTAdual}
$$
\preskip.4em \postskip.4em
\alpha :M\longrightarrow \gB,\quad z \longmapsto z,
$$
which only gives partial information on the structure of $M$.
For example, for every \lin form $\beta :M\to\gB$ we have $\beta (M)\subseteq\fb$
and therefore $M$ does not have any free direct summands 
of rank $1$ (cf. Proposition~\ref{propSplittingOffAlgExt}).\\
Similarly, the dual $(M\te_\gB M)\sta$ of $M\te_\gB M$ is free of rank $1$,
generated by the form
$$
\preskip-.4em \postskip.4em
\varphi  :M\te_\gB M\longrightarrow \gB,\quad z\te z' \longmapsto zz',
$$
and $M\te_\gB M$ does not possess a free direct summand of rank $1$.\\
Concerning $\gS_\gB^2 M$, we find that it admits a \sgr equal to~$\big(\gs(x,x),\alb\gs(x,y),\gs(y,y)\big)$,
with the \mpn
$$
\preskip.4em \postskip.4em
\crmatrix{y&0\cr-x&y\cr0&-x}.
\leqno{
\begin{array}{c} 
\gs(x,x)\\ \gs(x,y)\\ \gs(y,y)\\[-8.5pt] ~
\end{array}} 
$$
Concerning $\Al2_\gB M$, we find that it is generated by $x\vi y$
with the \mpn $\lst{x\;y}$ which gives
$$\preskip-.3em \postskip.1em
\;\;\Al2_\gB M\,\simeq\,\gB/\fb\,\simeq\,\gA.
$$
But be careful of the fact that $\gA$ as a \Bmo is a quotient and not a submodule of $\gB$.
\eoe

\rdb
\subsection*{Changing the base \riz}
\addcontentsline{toc}{subsection}{Changing the base \riz}
\label{pageChgtBase}

Let $\rho:\gA\rightarrow \gB$ be an \algz.
Every \Bmo $P$ can be equipped with  an \Amo  structure via $\rho$ by letting $a.x\eqdefi \rho(a)x$.

\begin{definition}
\label{defAliAliExtScal} Let $\gA\vers{\rho}\gB$ be an \Algz.
\begin{enumerate}
\item Let $M$ be an \Amoz. An \Ali map $\varphi :M\to P$, where~$P$ is a \Bmoz, is called a \emph{morphism of scalar extension} 
(from~$\gA$ to~$\gB$ for~$M$),
or a \ix{change of the base ring} 
(from $\gA$ to $\gB$ for~$M$), 
if the following \uvl \prt is satisfied.

\vspace{-15pt}

\PNV{M}{\varphi}{\psi}{P}{\theta}{R}{\Amos}{\Alis}{\Bmosz, \Blis}

\vspace{-20pt}
\emph{For every \Bmo $R$,  
every \Ali $\psi :M\to R$ is uniquely 	
 expressible in the form $\psi=\theta\circ\varphi$, where $\theta\in\Lin_\gB(P,R)$.}%

\item A \Bmo $P$  such that there exist an \Amo $M$ and a morphism of \eds $\varphi :M\to P$ is said to be
\ixd{extended}{module} \emph{from~$\gA$}.
We will also say that $P$ \emph{stems from the \Amo $M$ by \edsz}.
\end{enumerate}%
\index{morphism!of \edsz}%
\index{scalar extension}
\end{definition}

It is clear that a morphism of \eds
$\varphi :M\to P$ is unique in the categorical sense, \cad that for every other morphism of \eds
$\varphi' :M\to P'$, there is a unique $\theta\in\Lin_\gB(P,P')$ which renders the suitable diagram commutative, and that $\theta$ is an \isoz.

If $(\ug)$ is a \sgr of $M$ and $P$ an arbitrary \Bmoz, an \Ali $\lambda : M\rightarrow P$ is known from its values over the \elts $x$ of $\ug$.
In addition, the values $\lambda(x)$ are linked by certain constraints, which are derived from syzygies between \elts of $\ug$ in $M$. For example, if we have a syzygy $a_1x_1+a_2x_2+a_3x_3=_M 0$ between \eltsz~$x_i$ of~$\ug$, with the~$a_i$'s in $\gA$, this provides the following syzygy between the $\lambda(x_i)$'s in~$P$:
$\rho(a_1)\lambda(x_1)+\rho(a_2)\lambda(x_2)+\rho(a_3)\lambda(x_3)=0$.

Actually \gui{those are the only essential constraints, and that shows that a \eds can be constructed.}

More \prmtz, let $\rho\ist(x)$ replace $x$ (an arbitrary \elt of~$\ug$). Consider then the \Bmo $M_{1}$ generated by the $\rho\ist(x)$'s, linked by the syzygies described above
(
$\rho(a_1)\rho\ist(x_1)+\rho(a_2)\rho\ist(x_2)+\rho(a_3)\rho\ist(x_3)=_P 0$
for the given example).

\begin{proposition}
\label{propPfExt}
 (With the above notations)
\begin{enumerate}
\item \label{i1propPfExt}
\begin{enumerate}
\item There exists a unique \Ali $\varphi :M\to M_{1}$
such that for all $x\in\ug$, we have $\varphi(x)=\rho\ist(x)$.

\item This \Ali makes  $M_{1}$ a \eds from $\gA$ to $\gB$ for $M$. We will denote it by $M_{1}=\rho\ist(M)$.\perso{the inconv\'enient is que $M_{1}$ n'est pas en \gnl equal to the image of $M$ par $\varphi=\rho\ist$.}
\label{NOTArhosta}

\item In the case of a \mpfz, if $M$ is (\isoc to the) cokernel of a matrix
$F=(f_{i,j})\in\gA^{q\times m}$, then $M_{1}$
is (\isoc to the) cokernel of
{\em the same matrix seen in $\gB$}, \cad the matrix $F^\rho=\big(\rho(f_{i,j})\big)$.
In particular, if $M$ is free with basis $(\ug)$, then $M_{1}$ is free with basis
$\rho\ist(\ug)$.
\end{enumerate}
%
\item \label{i2propPfExt}
Consequently the \eds from $\gA$ to $\gB$ for an arbitrary \Amo exists and can always be defined from a \pn of this module. If the module is \tf (resp.\,\pfz) the \eds is as well.
\item \label{i3propPfExt}
 Knowing that the \edsz s exist, we can describe the previous construction (in a noncyclic manner) as follows:\\
if~$M\simeq \Coker \alpha$ with $\alpha :L_1\rightarrow L_2$, the modules $L_i$ being free, then the module
$M_{1}= \Coker \big(\rho\ist(\alpha)\big)$
is a \eds from $\gA$ to $\gB$ for the module $M$.
\item \label{i4propPfExt}
The \eds is transitive. If \smash{$\gA\vers\rho\gB\vers{\rho'}\gC$} are two \gui{successive} \algs and if $\rho''=\rho'\circ \rho$ define the \gui{composite} \alg, the canonical \Cli  $\rho''\ist(M)\to \rho'\ist\big(\rho\ist(M)\big)$
is an \isoz.
\item \label{i5propPfExt}
The \eds and the tensor product commute. If~$M$,~$N$ are \Amos and $\rho:\gA\to\gB$ is a \homo of \risz, then the natural \Bli $\rho\ist(M\te_\gA N)\to\rho\ist(M)\te_\gB\rho\ist(N)$ is an \isoz.
\item  \label{i6propPfExt}
Similarly the \eds commutes with the construction of the exterior powers, of the \smq powers and of the exterior \algz.
\item  \label{i7propPfExt}
Seen as an \Amoz, $\rho\ist(M)$ is (uniquely) \isoc to the tensor product $\gB\te_\gA M$ (here $\gB$ is equipped with its \Amo structure via~$\rho$). In addition, the \gui{external law} $\gB\times \rho\ist(M)\to\rho\ist(M),$ which defines the \Bmo structure   of $\rho\ist(M)$, is interpreted via the previous \iso like the \Ali 

\snic{\pi\te_\gA \Id_M:\gB\te_\gA\gB\te_\gA M\lora  \gB\te_\gA M,}

obtained from the \Ali $\pi:\gB\te_\gA\gB\to\gB$ \gui{product in $\gB$} ($\pi(b\te c)=bc$).
\item \label{i8propPfExt}
For every exact sequence of \Amos \smash{$M\vers{f} N\vers{g} P\to 0$} the sequence
$$\preskip-.5em \postskip-.2em 
\rho\ist(M)\vvvvers{\rho\ist(f)} \rho\ist(N) \vvvvers{\rho\ist(g)} \rho\ist(P)\to 0 
$$
is exact.
\end{enumerate}
\end{proposition}
\facile

Thus, a \Bmo $P$ is extended from $\gA$ \ssi it is \isoc to a module $\rho\ist(M)$.
Care must be taken, however, to the fact that an extended \Bmo can be derived from several non-\isoc \Amosz. For example when we extend a $\ZZ$-module to $\QQ$, \gui{we kill the torsion,} and $\ZZ$ and $\ZZ\oplus \aqo{\ZZ}{3}$ both give by \eds a \Qev of dimension~$1$.

\smallskip
\rem With the tensorial notation of item \emph{\ref{i7propPfExt}} the canonical \iso given at item \emph{\ref{i5propPfExt}} is written as

\snic{
\gC\otimes_\gA M \vvers{\varphi} \gC\otimes_\gB(\gB \otimes_\gA M)\simeq
(\gC\otimes_\gB\gB) \otimes_\gA M,}

with $\varphi(c\te x)=c\te(1_\gB\te x)$.
We will come back to this type of \gui{associativity} in the remark that follows Corollary~\ref{corPlatEds}.
\eoe

\subsec{Modules of \alis}

\begin{proposition}
\label{propAliCoh}
  If $M$ and $N$ are \mpfs over a \coh \ri $\gA$, then $\Lin_\gA(M,N)$ is \pfz.
\end{proposition}
\begin{proof}
We use the notations of Section~\ref{secCatMpf}.\\ 
Giving an \elt $\varphi$ of $\Lin_\gA(M,N)$ reduces to giving the matrices of $\rG_\varphi$ and~$\rK_\varphi$ that satisfy the condition $\rG_\varphi\,\rA_M=\rA_N\,\rK_\varphi$.
\\
Since the \ri is \cohz,
the solutions of the \sli form  
a \tf \Amoz, generated for example by the solutions corresponding to \alis $\varphi_1$, \ldots, $\varphi_\ell$ given by pairs of matrices $(\rG_{\varphi_1},\rK_{\varphi_1})$, \ldots, $(\rG_{\varphi_\ell},\rK_{\varphi_\ell})$.
Therefore $\Lin_\gA(M,N)=\gen{\varphi_1,\ldots ,\varphi_\ell}$.\\
Furthermore, a syzygy $\som_i a_i\varphi_i=0$ is satisfied \ssi 
we 
have a \ali $Z_\varphi:\rG_M\rightarrow \rK_N$ satisfying $\rA_N\, Z_\varphi=\som_ia_i\rG_{\varphi_i}$. By taking the corresponding \sliz, whose unknowns are the $a_i$'s on the one hand and the \coes of the matrix of $Z_\varphi$ on the other, we note that the syzygy module for the \sgr $(\varphi_1,\ldots ,\varphi_\ell)$ is indeed \tfz.
\end{proof}

\subsection*{The local character  of the \mpfsz}
\addcontentsline{toc}{subsection}{The local character  of the \mpfsz}

The fact that an \Amo is \pf is a local notion, in the following sense.

\begin{plcc}
\label{plcc.pf}  \emph{(\Fp modules)}\\
Let $S_1$, $\ldots$, $S_n$ be \moco of a \ri $\gA$, and $M$ 
 an \Amoz.
Then, $M$ is \pf \ssi each of the~$M_{S_i}$'s is a \pf $\gA_{S_i}$-module.
\end{plcc}
\begin{proof}
Assume that $M_{S_i}$ is a \pf $\gA_{S_i}$-module for each~$i$.
Let us show that $M$ is \pfz.
\\
By the \plgrf{plcc.tf}, $M$ is \tfz. Let $(g_1,\ldots ,g_q)$ be a \sgr of $M$.
\\
Let $(a_{i,h,1}, $\ldots$, a_{i,h,q})\in \gA_{S_i}^q$ be syzygies between the
$g_j/1\in M_{S_i}$ (in other words, $\sum_j \,a_{i,h,j}g_j=0$ in $M_{S_i}$) for $h=1,$ \ldots$, k_i$, which generate  
the $\gA_{S_i}$-syzygy module between the $g_j/1$. 
\\
Suppose \spdg that the $a_{i,h,j}$'s are of the form $a'_{i,h,j}/1$, with $a'_{i,h,j}\in \gA$. Then there exists some suitable $s_i\in S_i$ such that the vectors
$$\preskip-.3em \postskip.3em 
s_i\,(a'_{i,h,1},\ldots,a'_{i,h,q}) = (b_{i,h,1},\ldots,b_{i,h,q}) 
$$
are $\gA$-syzygies between the~$g_j\in M$.
\\
Let us show that the syzygies thus constructed between the~$g_j$'s generate all the syzygies. With this in mind, consider an arbitrary syzygy $(\cq)$ between the~$g_j$'s. Let us view it as a syzygy between the $g_j/1\in M_{S_i}$, and let us write it as an $\gA_{S_i}$-\lin combination of the vectors $(b_{i,h,1},\ldots,b_{i,h,q})$ in~$\gA_{S_i}^q$. After multiplication by some suitable $s'_i\in S_i$ we obtain an \egt in~$\gA^q$
$$\preskip.1em \postskip.4em
s'_i(\cq)=e_{i,1}(b_{i,1,1},\ldots,b_{i,1,q})+
\cdots+e_{i,k_i}(b_{i,k_i,1},\ldots,b_{i,k_i,q}).
$$
We write $\sum_{i=1}^{n} u_i s'_i =1$. We see that $(\cq)$ is an $\gA$-\lin combination of the $(b_{i,h,1},\ldots,b_{i,h,q})$.
\end{proof}

\subsec{Null tensors}

Let $M$ and $N$ be two arbitrary \Amosz, and $t=\sum_{i\in \lrbn}x_i\te y_i\in M\te N$.
The \egt $\sum_{i}x_i\te y_i=0$
does not a priori solely depend on knowing the submodules $\sum_i\gA x_i\subseteq M$ and $\sum_i\gA y_i\subseteq N$. 
\\
Consequently the notation $\sum_{i}x_i\te y_i$ is \gnlt burdened with ambiguity, and is dangerous. We should use the following more precise notation: $\sum_{i}x_i\otimes_{\gA,M,N} y_i$, or at least write the \egts in the form 

\snic{\sum_{i}x_i\te y_i=_{M\te_\gA N} \dots}

This precaution is not needed in the case where the two modules~$M$ and $N$ are flat (see Chapter~\ref{chap mod plats}), for instance when the \riz~$\gA$ is a \cdiz.

\pagebreak

\CMnewtheorem{lemtenul}{Null tensor lemma}{\itshape}
\begin{lemtenul} \label{lem-tenul}\index{Null tensor lemma}
Let $M=\gA x_1+\cdots+\gA x_n$ be a \mtfz, $N$ be another module and $t=\sum_{i\in \lrbn}x_i\te y_i\in M\te_\gA N$. \\
With $X=[\,x_1\,\cdots\,x_n\,]\in M^{1\times n}$ and $Y=\tra[\,y_1\,\cdots\,y_n\,]\in N^{n\times 1}$,
we use the notation $t=X\odot Y$.
\Propeq%
\begin{enumerate}
\item $t=_{M\te_\gA N}0$.
\item We have a $Z\in N^{m\times 1}$ and a matrix $G\in \Ae {n\times m}$ which satisfy
\begin{equation}\label{eqdeff2.plat}\preskip.3em 
XG=_{M^ m} 0\quad {\rm and } \quad GZ=_{N^{n}}Y\,.
\end{equation}
%
\end{enumerate}
\end{lemtenul}
%

\begin{proof}
\emph{2} $\Rightarrow$ \emph{1.} \Gnlt the \egt $X\odot G Z=XG\odot Z$ is guaranteed for every matrix $G$ with \coes in $\gA$ because $x\otimes \alpha z=\alpha x\otimes z$ when~$x\in M$, $z\in N$ and~$\alpha\in\gA$.

\emph{1} $\Rightarrow$ \emph{2.} 
The \egt $t=_{M\te N}0$ comes from a finite number of syzygies within the modules $M$ and $N$. Therefore there exists a submodule $N'$ such that
$$\preskip-.30em \postskip.4em
\gA y_1+\cdots+\gA y_n\subseteq N'=\gA z_1+\cdots+\gA z_m\subseteq N, 
$$
and $X\odot Y=_{M\otimes N'} 0$. We write $Z=\tra[\,z_1\,\cdots\,z_m\,]$.
We then have an exact sequence
$$\preskip-.30em \postskip.4em 
K \vers{a}L\vers\pi N'\to 0 
$$
where $L$ is free, with basis $(\ell_1,\dots,\ell_m)$ and $\pi(\ell_j)=z_j$, which gives an exact sequence
$$\preskip-.30em \postskip.4em M\te K \vvvers{\rI\te a}M\te L\vers{\rI\te \pi} M\te N'\to 0.
$$
If $U\in M^{1\times m}$ satisfies $U\odot Z=_{M\te N'}0$, this means that $U$ seen as an \elt of $M\te L\simeq M^{n}$ (\cad seen as $\sum_j u_j\te_{M\otimes L} \ell_j$) is in the submodule $\Ker(\rI\te \pi)=\Im(\rI\te a)$, in other words

\snic{\sum_j u_j\te_{M\otimes L} \ell_j=\sum_ix_i\otimes \sum_{ij}a_{ij}\ell_j=\sum_j \big(\sum_i a_{ij}x_i\big)\te \ell_j}

for $a_{ij}\in\gA$ that satisfy $\sum_{j}a_{ij}z_j=0$.
In other words $U=XA$ for a matrix~$A$ satisfying $AZ=0$.
\\ 
If we write $Y=HZ$ with $H\in\Ae{n\times m}$, we have $XH\odot Z=0$, which gives an \egt $XH=XA$ with a matrix~$A$ satisfying $AZ=0$. \\
We then let $G=H-A$ and we have $XG=0$ and $GZ=HZ=Y$. 
\end{proof}

\vspace{-.5em}
\section[Classification \pbsz]{Classification \pbs for  \mpfsz}

The first classification \tho concerns  free \Amos of finite rank: two \Amos $M\simeq \Ae m$ and $P\simeq \gA^p$ 
with $m\neq p$ can only be \isoc 
if $1=_\gA0$ (Proposition~\ref{propDimMod1}).

\smallskip 
\rem
Note 
that we use the expression \gui{$M$ is a free module of rank $k$} to mean that $M$ is \isoc to $\Ae k$, even in the case where we ignore whether the \ri $\gA$ is trivial or not. This therefore does not always imply that a priori this integer $k$ is well-determined. 
\eoe

\smallskip 
Rare are the rings for which we can give a \gui{satisfactory} complete classification of the \mpfsz. The case of \cdis is well-known: every \mpf is free (this results from the Chinese pivot or from the freeness lemma).
In this work we treat a few \gnns of this \elr case: the valuation \risz, the PIDs and the reduced \zeds \ris (Sections~\ref{secBézout} and~\ref{secKrull0dim}), and certain Pr\"ufer \ris (Proposition~\ref{propAriCohZed} and \Thref{thMpfPruCohDim}).

Concerning the classification of the \mtfsz, we note the following two important uniqueness results.

\vspace{-.5em}
\subsec{Two results concerning  \mtfsz}

\begin{theorem}
\label{prop unicyc} 
Let $\fa_1\subseteq \cdots\subseteq\fa_n$ and $\fb_1\subseteq \cdots\subseteq\fb_m$ be \ids of $\gA$ with $n \le m$. 
If an \Amo $M$ is \isoc to \hbox{$\gA/\fa_1\oplus\cdots\oplus \gA/\fa_n$} and \hbox{to $\gA/\fb_1\oplus\cdots\oplus \gA/\fb_m$}, then
\begin{enumerate}
\item  we have $\fb_k=\gA$ for $n < k \le m$,
\item  and $\fb_k=\fa_k$ for $1 \le k \le n$.
\end{enumerate}
We say that $(\fa_1, \dots,\fa_n)$ is the list of  \emph{invariant factors\footnote{Note that the list given here can be shortened or extended with terms $\fa_j=\gen{1}$ when we do not have a test for the \egt in question. This is comparable to the list of \coes of a \pol that can be shortened or extended with $0$'s when the \ri is not discrete.}} of the module~$M$.%
\index{invariant factors}\index{factors!invariant --- (of a module)}
\end{theorem}

\begin{proof} 
\emph {1.} It suffices to show that if $n < m$, then $\fb_m=\gA$, in other words that the \ri $\gB := \gA/\fb_m$ is null. By letting $M=\gA/\fa_1\oplus\cdots\oplus \gA/\fa_n$, we have 

\snic{\gB^m=\bigoplus_{j=1}^m\gA/(\fb_j+\fb_m)\simeq M/\fb_mM\simeq
\bigoplus_{i=1}^n\gA/(\fa_i+\fb_m).}

But each $\gA/(\fa_i+\fb_m)$ is a quotient \ri of $\gB$, 
so there exists a surjective \ali from $\gB^n$ onto $\gB^m$ and therefore $\gB$ is null (Proposition~\ref{propDimMod1}).
We assume henceforth \spdg that $m=n$.

\emph {2.} It suffices to show that $\fb_k\subseteq\fa_k$ for $k\in\lrbn$. Notice that for an \id $\fa$ and an \elt $x$ of $\gA$, the kernel of the \aliz~$y\mapsto yx \mod\fa$, from~$\gA$ to $x(\gA/\fa)$ is the \id $(\fa:x)$, and thus that 
$x(\gA/\fa)\simeq \gA/(\fa:x)$.
Now let $x \in \fb_k$. For $j\in\lrb{k..n}$, we have $(\fb_j : x) = \gA$, and therefore

\snic{
xM \simeq \bigoplus_{j=1}^n\gA/(\fb_j:x) =
\bigoplus_{j=1}^{k-1}\gA/(\fb_j:x),
\; \hbox { and } \;
xM \simeq \bigoplus_{i=1}^n\gA/(\fa_i:x).
}

By applying item~\emph {1} to the module $xM$ with the integers $k-1$ and $n$, we
obtain $(\fa_k:x)=\gA$, i.e. $x\in\fa_k$.
\end{proof}

Note that in the previous \thoz, we have not assumed anything regarding the \ids (it is not \ncr that they be \tf nor detachable for the result to be \cot valid).
\begin{theorem}
\label{prop quot non iso} 
Let $M$ be a \tf \Amo and $\varphi\,:\,M\rightarrow M$ be a surjective \aliz.
Then, $\varphi$ is an \iso and its inverse is a \pol in $\varphi$. If a quotient module $M/N$ of $M$ is \isoc to $M$,  then~$N=0$.
\end{theorem}
\begin{corollary}
\label{corInvDInv} If $M$ is a \mtfz, every \elt $\varphi$ right-\iv in $\End_\gA(M)$ is \ivz, and its inverse is a \pol in~$\varphi$.
\end{corollary}
\begin{Proof}{Proof of \Thref{prop quot non iso}. }\\
Let $(\xn)$ be a \sgr of~$M$, $\gB=\gA[\varphi]\subseteq\End_\gA(M)$, and $\fa=\gen{\varphi}$ be the \id of $\gB$ generated by $\varphi$. The \riz~$\gB$ is commutative and we consider $M$ as a \Bmoz. Since the \ali $\varphi$ is surjective, there exists a $P\in \Mn(\fa)$ with $P\,\tra{[\,x_1\;\cdots\;x_n\,]}=\tra{[\,x_1\;\cdots\;x_n\,]}$, \cad

\snic{(\In-P)\,\tra{[\,x_1\;\cdots\;x_n\,]}= \tra{[\,0\;\cdots\;0\,]}}

(where $\In=(\In)_\gB$ is the identity matrix of $\Mn(\gB)$),
and so

\snic{\det(\In-P)\,\tra{[\,x_1\;\cdots\;x_n\,]}
=   
\wi{(\In-P)}\,(\In-P)\,\tra{[\,x_1\;\cdots\;x_n\,]}      
= \tra{[\,0\;\cdots\;0\,]}.
}

Therefore $\det(\In-P)=0_\gB$, but $\det(\In-P)=1_\gB-\varphi\, \psi$ with $\psi\in \gB$ (since~$P$ has \coes in $\fa=\varphi\, \gB$). Thus, $\varphi\, \psi= \psi\,\varphi=\alb 1_\gB=\alb \Id_M$: $\varphi$ is \iv in $\gB$.
\end{Proof}

\section{Quasi-integral \risz}
\label{subsecAnneauxqi}

In the following \dfnz, we infinitesimally modify the notion of an integral \ri usually given in \comaz, not for pleasure, but because our \dfn better corresponds to \algos implementing integral \risz.

\begin{definition}
\label{defqi}
A \ri is said to be
 \emph{integral} if every \elt is null or \ndzz.\footnote{An integral ring is also called a \emph{domain} in the classical literature.
But we prefer to keep \gui{integral ring} in order to distinguish them from rings \gui{witout zerodivisors}. See the \dfn on \paref{eqSDZ}.}
A \ri $\gA$ is said to be
 \emph{quasi-integral} when every \elt admits as its annihilator an (\id generated by an) \idmz.
In the literature, a quasi-integral \ri is sometimes called a 
\emph{pp-ring} (principal \ids are projective, cf. Section~\ref{Idpp}).%
\index{ring!integral ---}%
\index{pp-ring (or quasi-integral ring)}%
\index{ring!quasi-integral --- (or pp-ring)}%
\index{integral!ring}%
\index{integral!quasi --- ring (or pp-ring)}%
\index{quasi-integral!ring (or pp-ring)}%
\end{definition}

As usual, the \gui{or} in the previous \dfn must be read as an explicit or.
An integral \ri is therefore a discrete set \ssi furthermore it is trivial or nontrivial. So, our nontrivial integral rings are precisely the \gui{discrete domains} of~\cite{MRR}.

In this work, sometimes we speak of a \gui{nonzero \eltz} in an integral \riz, but we should actually say \gui{\ndz \eltz} in order not to exclude the trivial \ri case.

\begin{fact}\label{factQIRed}
A \qiri  is reduced.
\end{fact}
%
\begin{proof}
If $e$ is the \idm annihilator of $x$ and if $x^2=0$, then $x\in\gen{e}$,
 therefore $x=ex=0$.
\end{proof}

A \cdi is an integral \riz. A \ri $\gA$ is integral \ssi its total \ri of fractions $\Frac\gA$ is a \cdiz.
A finite product of \qiris  is a \qiriz.
A \ri is integral \ssi it is a connected \qiriz.

\subsec{Equational \dfn of  \qirisz}
In a \qiriz, for $a\in\gA$, let $e_a$ be the unique \idm such that
$\Ann(a)=\gen{1-e_a}$. We have 
$\gA\simeq \gA[1/e_a]\times \aqo{\gA}{e_a}$.\\
In the \ri $\gA[1/e_a]$, the \elt $a$ is \ndzz, and in $\aqo{\gA}{e_a}$, $a$ is null.
\\
We then have $e_{ab}=e_a e_b$, $e_aa=a$ and $e_0=0$.
\\
Conversely, suppose that a commutative \ri is equipped with a unary  law $a\mapsto \ci{a}$ which satisfies the following three axioms
\begin{equation}\label{eqaqis}\preskip.3em \postskip.3em
\ci{a}\,a=a,\quad
\ci{(ab)}=\ci{a}\,\ci{b},\quad
\ci{0}=0.
\end{equation}
Then, for all $a\in\gA$, we have $\Ann(a)=\gen{1-\ci{a}}$, and $\ci{a}$ is \idmz, such that the \ri is a \qiriz.

Indeed, first of all $(1-\ci{a})a = 0$, and if $ax=0$, then
$$\preskip.3em \postskip.3em 
\ci{a}\,x=\ci{a}\,\ci{x}\,x=\ci{(ax)}\,x=\ci{0}\,x=0, 
$$
so $x=(1-\ci{a})x$. Hence $\Ann(a)=\gen{1-\ci{a}}$.
Next let us show that $\ci{a}$ is \idmz. Apply the previous result to $x = 1-\ci{a}$ which satisfies $ax = 0$ (by the first axiom); the \egt $x=(1-\ci{a})x$ %
 gives $x=x^2$, i.e. the \elt $1 - \ci{a}$ is \idmz.

The following splitting lemma is almost immediate.

\CMnewtheorem{lemScindageQi}{Quasi integral splitting lemma}{\itshape}

\begin{lemScindageQi}\label{thScindageQi}
Let $x_1$, \dots, $x_n$ be  $n$ \elts in a \qiriz~$\gA$.
There exists a \sfio $(e_j)$ of cardinality~$2^n$ such that in each of the components~\hbox{$\gA[1/e_j]$}, each~$x_i$ is null or \ndzz.
\end{lemScindageQi}
\begin{proof}
Let $r_i$ be the \idm such that $\gen{r_i}=\Ann(x_i)$, and $s_i=1-r_i$.
By expanding the product $1=\prod_{i=1}^n(r_i+s_i)$ we obtain the \sfio indexed by $\cP_n$: $e_J = \prod_{j\in J}r_j \prod_{k\notin J}s_k$.
We can delete certain \elts of this \sys when we say that they are null.
\end{proof}

%
\subsect{From integral \ris to \qirisz}{Elementary \lgb machinery \num1: from integral \ris to \qirisz}

Knowing how to systematically split  a \qiri into two components  leads to the following \gnl method.
The essential difference with the previous splitting lemma is that we a priori do not know the finite family of \elts which will provoke the splitting.

\rdb
\mni {\bf Elementary \lgb machinery \num1.}\label{MethodeQI}\imlgd
{\it Most \algos that work with  nontrivial integral \ris can be modified in order to work with  \qirisz, by splitting the \ri into two components each time that the \algo written for the integral \ris uses the \gui{is this \elt null or \ndzz?} test. In the first component the \elt in question is null, in the second it is \ndzz.}

\medskip A first example of an application of this \lgb machinery will be given on \paref{exlBezQi}.
However, Corollary~\ref{corlemQI} below could already be obtained from the integral case, where it is obvious, by applying this \lgb machinery.
\\
Let us explain why we speak of \elr \lgb machinery here. Generally a \plg says that a \prtz~$\sfP$ is true \ssi it is true \gui{after \lon at \mocoz.} In the current case, the \moco are generated by \elts $1-r_i$ where the $r_i$'s form a \sfioz. Consequently the \ri is simply \isoc to the product of the localized \risz, and the situation is therefore perfectly simple, \elrz.

\smallskip
\rem \label{remMLGE1}
The reader will have noticed the very informal formulation that we have given for this \lgb machinery: \gui{Most \algos \ldots} This is because it seemed quite difficult to give very precise requirements in advance for the indicated method to work. We could imagine an \algo which works for every integral \riz, but in a completely non-uniform manner, which would make the corresponding tree that we construct in the \qiri case not finite. For example, in the integral case, a given starting configuration would require three tests (to end the computation) if the answers are $0,0,0$, but four tests if the answers are $0,0,1,0$, then five tests if the answers are $0,0,1,1,0$, then six tests if they are $0,0,1,1,1,1$, then seven tests if they are $0,0,1,1,1,0,1$,  etc.  Naturally, we can doubt that such an \algo could exist without the existence of an integral \ri that would fault it at the same time.
In other words, an \algo that is not sufficiently uniform is likely not an \algoz. But we do not assume anything.\\
Even if we have not so far encountered any example of the above type where the \elr \lgb machinery would not apply, we cannot a priori exclude such a possibility.  \eoe

\subsect{Annihilators of the \itfs in  \qirisz}
{Annihilators of the \itfsz}

The following lemma can be considered as an economical variant of the splitting lemma \ref{thScindageQi}.

\pagebreak

\begin{lemma}\label{lemQI}
Let $x_1$, \dots, $x_n$ be \elts of an \Amoz.
\\
 If we have $\Ann(x_i) = \gen{r_i}$ where $r_i$ is an \idm  ($i \in\lrbn$), let 

\snic{s_i=1-r_i$,
$t_1=s_1$, $t_2=r_1s_2$, $t_3=r_1r_2s_3$, $\dots$,
$t_{n+1}=r_1r_2\cdots r_n.}

Then $(t_1, \dots, t_{n+1})$ is a \sfio
and the \elt $x=x_1+t_2x_2+\cdots +t_nx_n$ satisfies

\snic{\Ann(x_1,\dots,x_n) = \Ann(x) = \gen{t_{n+1}}.}

NB: in the component $t_k = 1$ ($k\in\lrbn$), we have 
   $x_k$ \ndz and $x_j = 0$ for $j < k$, and in the component $t_{n+1}=1$,
   we have $x_1=\cdots=x_n=0$.
\end{lemma}

\begin{corollary}
\label{corlemQI}
Over a \qiri  $\gA$ every \tf submodule~$M$ of a free module has as its annihilator an \id $\gen{r}$ with $r$ \idmz, and~$M$ contains an \elt $x$ having the same annihilator.
This applies in particular to a \itf of $\gA$.
\end{corollary}

\begin{Proof}{\Demo of Lemma~\ref{lemQI}.}
We have $t_1 x_1 = x_1$ and

\snic{\arraycolsep2pt
\begin{array}{rl}
1~ = & s_1 + r_1 = s_1 + r_1(s_2 + r_2) ~=~  s_1 + r_1s_2 +
r_1r_2(s_3+r_3)
~=~ \cdots
\\
= & s_1 + r_1s_2 + r_1r_2s_3 +\cdots +r_1r_2\cdots r_{n-1}s_n +
   r_1r_2\cdots r_n
\end{array}}

so $t_1,\dots,t_{n+1}$ is a \sfio and $x = t_1 x_1 + t_2 x_2 + \cdots + t_nx_n$.
It is clear that

\snic{\gen{t_{n+1}} \subseteq \Ann(x_1,\dots,x_n) \subseteq \Ann(x).}

Conversely, let $z\in \Ann(x)$. Then $zx=0$, so $z t_i x_i = z t_i x=0$ 
 for $i \in\lrbn$. Thus, $z t_i \in \Ann(x_i) = \gen{r_i}$ and $z t_i = z t_i r_i = 0$.
Finally, since $z = \sum_{i=1}^{n+1} z t_i$, we have $z = z t_{n+1} \in \gen{t_{n+1}}$.
\end{Proof}


\subsect{Concrete \lgb principle for the \qirisz}{Local-global principle}

The \prt of being a \qiri is local in the following sense.

\begin{plcc}
\label{plcc.aqi}  
\emph{(\qirisz)}\\
Let $S_1$, $\dots$, $S_n$ be \moco of $\gA$.
\Propeq
\begin{enumerate}
\item  The \ri $\gA$ is a \qiriz.
\item For $i=1$, $\dots$, $n,$
each \ri $\gA_{S_i}$ is a \qiriz.
\end{enumerate}
\end{plcc}
%
\begin{proof}
Let $a\in\gA$. For every \mo $S$ of $\gA$ 
we have $\Ann_{\gA_S}(a)=\big(\Ann_{\gA}(a)\big)_S$.
Therefore the annihilator $\fa$ of $a$ is \tf \ssi it is \tf after \lon at the $S_i$'s (\plgref{plcc.tf}).
Next the inclusion $\fa\subseteq\fa^2$ is a matter of the basic \plgc\ref{plcc.basic}.\iplgz
\end{proof}

\vspace{-.7em}
\pagebreak

\section{Bézout \risz}
\label{secBézout}

A \ri $\gA$ is called a \emph{Bézout \riz} when every \itf is principal.
This is the same as saying that every \id with two \gtrs is principal.%
\index{Bézout!ring}\index{ring!Bézout ---}
\begin{equation}\label{eqBézout}\preskip.1em \postskip.3em
\;\;\forall a,\,b\;\;\exists u,\,v,\,g,\,a_1,\,b_1\;\;
\left(au+bv=g,\,a=ga_1,\,b=gb_1\right).
\end{equation}

A Bézout \ri is \fdi \ssi the divisibility relation is explicit.
An integral Bézout \ri is called a \emph{Bézout domain}.%
\index{Bézout!domain}\index{domain!Bézout ---}

A \ixx{local}{ring} is a \ri $\gA$ where is satisfied the following axiom\index{ring!local ---}
\begin{equation}\label{eqAloc}\preskip.2em \postskip.1em
\forall x,\,y\in \gA \qquad x+y \in\Ati\;\Longrightarrow \; ( x\in\Ati\;\mathrm{or}\;
y\in\Ati )\,.
\end{equation}
This is the same as asking
$$
\preskip.1em \postskip.3em
\forall x\in \gA \qquad x \in\Ati\;\; {\rm or} \;\;  1-x \in\Ati\,.
$$
Note that according to this \dfn the trivial \ri is local.
Moreover, the \gui{or} must be understood in the \cof sense:
the alternative must be explicit.
Most of the local \ris with which we usually work in \clama actually satisfy the previous \dfn if we look at it from a \cof point of view.

Every quotient \ri of a \alo is local.
A \cdi is a \aloz.

\begin{lemma}
\label{lemBezloc} \emph{(Bézout always trivial for a \aloz)}\\
A \ri is a local Bézout \ri \ssi it satisfies the following \prtz: $\forall a,b\in\gA$, $a$ divides $b$ or $b$ divides~$a$.
\end{lemma}
\begin{proof}
\emph{The condition is sufficient}. First it is clear that the ring is Bézout. Now assume  $x+y$ to be invertible. If $x$ divides $y$, it divides~\hbox{$x+y$} which divides~$1$, so~$x$ is \ivz. Symmetrically, if $y$ divides $x$, it is \ivz.
So the ring is local. 

\emph{The condition is necessary}. Assume $\gA$ is Bézout and local.
We have $g(1-ua_1-vb_1)=0$. Since $1=ua_1+vb_1+(1-ua_1-vb_1)$, one of the three terms
in the sum is \ivz. If $1-ua_1-vb_1$ is \ivz,  
then $g=a=b=0$. If $ua_1$ is \ivz, then so is $a_1$, and $a$ divides~$g$  which divides $b$. If $vb_1$ is \ivz, then so is $b_1$, and $b$ divides~$g$  which divides~$a$.
\end{proof}

\rdb
Local Bézout \ris are therefore  \gui{valuation \risz} in the Kaplansky sense. We prefer the now usual \dfnz%
\index{ring!valuation ---}\index{valuation!ring}:
a \emph{valuation ring} is a reduced local Bézout \riz.

\subsect{\Fp modules over  valuation \risz}{\Fp modules over the valuation \risz} \label{secpfval}

A matrix $B=(b_{i,j})\in\Ae {m\times n}$ is said to be
 \ixd{in Smith form}{matrix} if every \coe out of the principal diagonal is null, and if for $1\leq i< \inf(m,n)$, the diagonal \coe $b_{i,i}$ divides the following $b_{i+1,i+1}$.%
\index{Smith!matrix in --- form}

\begin{proposition}
\label{propPfVal}
Let $\gA$ be a local Bézout \riz.
\begin{enumerate}
\item Every matrix of $\Ae{m\times n}$ is \elrt \eqve to a matrix
in Smith form.
\item  Every \pf \Amo $M$ is \isoc to a direct sum of modules $\aqo{\gA}{a_i}$:  $M\simeq\bigoplus_{i=1}^p \aqo{\gA}{a_i}$, with in addition, for each $i< p$, $a_{i+1}$ divides $a_{i}$.
\end{enumerate}
\end{proposition}
\begin{proof}
\emph{1.} We use the Gauss pivot method by choosing for first pivot a \coe of the matrix which divides all the others. We finish by \recuz. \\
\emph{2.} Direct consequence of item~\emph{1.}
\end{proof}

\rem
This result is completed by the uniqueness \tho (\Thref{prop unicyc}) as follows.
\begin{itemize}
\item [\emph{1.}] In the reduced matrix in Smith form the \ids $\gen{b_{i,i}}$ are uniquely determined.
\item [\emph{2.}] In the decomposition $\bigoplus_{i=1}^p \aqo{\gA}{a_i}$, the \ids
$\gen{a_{i}}$ are uniquely determined, except that \ids in excessive numbers can be equal to $\gen{1}$: we can delete the corresponding terms, but this only happens without fail when we have an invertibility test in the \riz.
\eoe\end{itemize}

\medskip A \ri $\gA$ is called a \emph{strict Bézout \riz}\index{ring!strict Bézout ---} when every  vector $\vab u v\in\gA^2$ can be transformed into a vector $\vab h 0$ by multiplication by a $2\times 2$ \iv matrix.%
\index{Bézout!strict --- ring}

Now we give an example of how the elementary local-global machinery \num1 (described on \paref{MethodeQI}) is used.

\ms\exl \label{exlBezQi}
We will show that \emph{every  Bézout \qiri is a strict Bézout \riz}.
\\
Let us start with the integral case. Let $u$, $v\in\gA$,

\snic{\exists\,h,a,b,u_1,v_1\quad (h=au+bv,\,u=hu_1,\,v=hv_1).}

If $\Ann(v)=1$, then $v=0$ and $\vab u 0 =\vab u v\;\I_2$. \\
If $\Ann(v)=0$, then  $\Ann(h)=0$, $h(au_1+bv_1)=h$, then $au_1+bv_1=1$.
Finally, {$\vab h 0=\vab u v\crmatrix{a & -v_1 \cr b  &u_1}$}
and the matrix has \deterz~1.
\\
Let us now apply the \elr \lgb machinery \num1 explained on \paref{MethodeQI}. Consider the \idm $e$ such that 

\snic{\Ann(v)=\gen{e}$ and $f=1-e.}

In $\gA[1/e]$, we have  $\vab u 0=\vab u v\;\I_2$. \\
In $\gA[1/f]$, we have  $\vab h 0=\vab u v\crmatrix{a & -v_1 \cr b  &u_1}$.\\
Therefore in $\gA$, we have  \smashbot{$\vab{ue+hf} 0 =\vab u v\cmatrix{fa+e & -fv_1 \cr fb  &fu_1+e}$},
and the matrix has \deterz~1.
\eoe

\subsec{\Fp modules over PIDs}

Assume that $\gA$ is a strict Bézout \riz. If $a$ and $b$ are two \elts on the same row (resp.\,column) in a matrix $M$ with \coes 
 in $\gA$, we can postmultiply (resp.\,premultiply) $M$ by an \iv matrix, which will modify the columns (resp.\,the rows) where the \coesz~$a$ and~$b$ are, which are replaced by~$c$ and~$0$.
When describing this transformation of matrices, we will speak of \emph{Bézout manipulations}.
The \elr manipulations can be seen as special cases of Bézout manipulations.%
\index{manipulation!Bézout ---}

An integral \ri is said to be a
 \ix{principal ideal domain (PID)}\index{domain!principal ideal ---}\index{PID (principal ideal domain)}
when it is Bézout and when every ascending sequence of principal \ids admits two equal consecutive terms (cf. \cite{MRR}). 
In other words a PID is a \noe Bézout domain (see \dfnz~\ref{definoetherien}).
Examples include $\ZZ$ and the \poll \ri $\gK[X]$ when $\gK$ is a \cdiz.

\begin{proposition}
\label{propPfPID}
Let $\gA$ be a PID.
%
\begin{enumerate}
\item  Every matrix $A\in\Ae {m\times n}$ is \eqve to a matrix in Smith form. By letting $b_i$ be the diagonal \coes of the reduced matrix, the \idps $\gen{b_1}\supseteq\cdots\supseteq\gen{b_q}$  ($q=\inf(m,n)$) are invariants of the {matrix} $A$ up to \eqvcz. 
A basis~\hbox{$(e_1,\dots,e_m)$} of~$\Ae{m}$ such \hbox{that $\Im(A)=\sum_{i=1}^{m}\gen{b_i} e_i$}
 is called a \emph{basis adapted to the submodule~\hbox{$\Im(A)$}.}%
\index{basis adapted!to an inclusion} 

\item  For every \pf \Amo $M$, there exist $r$, $p\in\NN$ and \ndz \elts $a_{1},$ \ldots, $a_{p}$, with $a_{i}$ dividing $a_{i+1}$ for $i < p$, such that~$M$ is \isoc to the direct sum $\bigl(\bigoplus_{i=1}^p \aqo{\gA}{a_i}\bigr)\oplus\Ae r$.
\end{enumerate}
If furthermore $\gA$ is  nontrivial and \fdiz, we can ask in item~{2} that no $\gen{a_i}$ be equal to $\gen{1}$.
In this case, we call \emph{invariant factors of the module $M$} the \elts of the list 
$\big(a_1,\ldots,a_p,\MA{\underline{0,\ldots,0}}\limits_{r\;{\it times}}\big)$, 
and the list of invariant factors of $M$ is well-defined\footnote{We find the given \dfn in \Thref{prop unicyc}. We will however note that the order is reversed and that here we have replaced the \idps by their \gtrsz, all of this to conform to the most common terminology.} \gui{up to association.}%
\index{invariant factors}\index{factors!invariant --- (of a module)}
\end{proposition}
\begin{Proof}{Proof idea. }
By the Bézout manipulations on the columns,
we replace the first row with a vector $(g_1,0,\ldots,0)$.
By the Bézout manipulations on the rows,
we replace the first column with a vector $(g_2,0,\ldots,0)$.
We continue the process until we have $g_{k}\gA=g_{k+1}\gA$ for an index $k$.
For example, with odd $k$ this means that the last row operations by means of Bézout manipulations have been mistakenly applied,
since $g_{k}$ divided the first column.
We backtrack by a step,
 and use $g_{k}$ as a Gauss pivot.
We thus obtain a matrix of the form
$$
\blocs{.5}{1.9}{.5}{1.9}{$g$}{$\begin{array}{ccc}0 & \cdots & 0\end{array}$}{$0$ \\[1mm] $~\vdots~$ \\[1mm]$0$}{$B$}
$$
By \recu we obtain a \gui{diagonal} reduced matrix.
We finally verify that we can pass, by Bézout manipulations and \elr manipulations, from a matrix
\halfsmashtop{$\left[\begin{array}{cc}a & 0 \\0 & b\end{array}\right]$
to a matrix $\left[\begin{array}{cc}c & 0 \\0 & d\end{array}\right]$}
where $c$ divides $d$.
\\
Item~\emph{2} is a direct consequence of item~\emph{1}.
\end{Proof}

\rems~\\
1) A simpler \algo can be devised if $\gA$ is \fdiz.
 
2)
We still do not know (in 2014) 
if the conclusion of the previous proposition is true under the sole assumption that~$\gA$ is a Bézout domain. We have neither a \demz, nor a counterexample.\\
However, we do know that the result is true for  Bézout domains of dimension $\leq 1$; see the remark that follows \Thref{thMpfPruCohDim}.
\eoe

\section{Zero-dimensional \risz}
\label{secKrull0dim}

We will say that a \ri is \emph{\zedz} when  the following axiom is satisfied.%
\index{zero-dimensional@\zedz!\riz}\index{ring!zero-dimensional}
\begin{equation}\label{eqZed}\preskip.0em \postskip.4em
\forall x\in \gA~\exists a\in\gA~\exists k\in
\NN\quad \quad x^{k}=ax^{k+1}
\end{equation}

A \ri is said to be
 \ixc{Artinian}{\riz} if it is \zedz, \coh and \noez.%
\index{ring!Artinian ---}

\subsec{Basic \prtsz}
\label{subsecBasicsZerdim}

\begin{fact}
\label{factExZed} ~
\begin{enumerate}
\item [--] Every finite \ri and every \cdi is \zedz.
\item [--] Every quotient \ri and every localized \ri of a \zed \ri is \zedz.
\item [--] Every finite product of \zed \ris is a \zed \riz.
\item [--] A Boolean \alg (cf. Section~\ref{secBoole}) is a \zed \riz.
\end{enumerate}
\end{fact}

\begin{lemma}
\label{lemme:idempotentDimension0} 
\Propeq
\begin{enumerate}
\item  \label{LID001} $\gA$ is \zedz.
\item  \label{LID002} $\forall x\in \gA~\exists s\in\gA~\exists d \in \NN^*$ such that $\gen{x^d} = \gen{s}$ and $s$ \idmz.
\item  \label{LID003}
For every \itf $\fa$ of $\gA$, there exists a $d \in \NN^*$ such that
$\fa^d = \gen{s}$ where $s$ is an \idmz, and in particular,
$\Ann(\fa^d)=\gen{1-s}$ and $\fa^e=\fa^d$ for $e \geq d$.
\end{enumerate}
\end{lemma}
\begin{proof}
\emph{\ref{LID001}} $\Rightarrow$ \emph{\ref{LID002}.} 
For all $x \in \gA$, there exist $a\in \gA$ and $k \in \NN$ such that $x^k = ax^{k+1}$.\\
If $k=0$ we have $\gen{x}=\gen{1}$, we take $s=1$ and $d=1$.
\\
If $k\geq1$, we take $d=k$; by multiplying $k$ times by $ax$, we obtain the \egts $x^k = a x^{k+1} = a^2 x^{k+2} = \cdots = a^k x^{2k}$.
Therefore the \elt $s=a^k x^k$ is an \idmz,~$x^k = sx^k$, and~$\gen{x^k} = \gen{s}$.

\emph{\ref{LID002}} $\Rightarrow$ \emph{\ref{LID001}.}
We have $s=bx^d$ and $x^d s=x^d$. Therefore, by letting $a=bx^{d-1}$,
we obtain the \egts $x^d = bx^{2d} = ax^{d+1}$.

\emph{\ref{LID002}} $\Rightarrow$ \emph{\ref{LID003}.}
If $\fa = x_1\gA + \cdots + x_n\gA$, there exist \idms $s_1,\dots,s_n\in \gA$ and integers $d_1,\dots,d_n \geq 1$ such that $x_i^{d_i}\gA = s_i\gA$.
Let 
$$
\preskip.3em \postskip.3em 
s = 1-(1-s_1)\cdots(1-s_n),
$$
such that $s\gA = s_1\gA + \cdots + s_n\gA$.
It is clear that the \idm $s$ belongs to~$\fa$,
and so to all the powers of~$\fa$.
Moreover, if~$d \geq d_1 + \cdots + d_n - (n-1)$ we have
$$
\preskip.0em \postskip.3em 
\fa^d \subseteq x_1^{d_1}\gA + \cdots + x_n^{d_n}\gA = s_1\gA + \cdots + s_n\gA = s
\gA. 
$$
The result follows since $\fa^d = s\gA$.

Finally, \emph{\ref{LID003}} clearly implies \emph{\ref{LID002}}.
\end{proof}

\begin{corollary}
\label{corZedReg}
If $\fa$ is a faithful \itf of a \zed \riz, then $\fa=\gen{1}$.
In particular, in a \zed \riz, every \ndz \elt is \ivz.
\end{corollary}
\begin{proof}
For $d$ large enough the \id $\fa^d$ is generated by an \idm $s$. This \id is \ndzz, therefore the \idm $s$ is equal to $1$.
\end{proof}

\begin{lemma}\label{lemZerloc}
\emph{(Local \zed \risz)}\\
\Propeq
\begin{enumerate}
\item   $\gA$ is local and \zedz.
\item   Every \elt of $\gA$ is \iv or nilpotent.
\item   $\gA$ is \zed and connected.
\end{enumerate}
\end{lemma}

Consequently a \cdi can also be defined as a reduced local \zed \riz.

\subsec{Reduced \zeds \risz}
\label{subsecAzedred}

\vspace{2pt}
\subsubsec{\Cara \prtsz}

The \eqvcs of the following lemma are easy (see the \dem of the analogous lemma, Lemma~\ref{lemme:idempotentDimension0}).
\begin{lemma}\label{lemZerRed}
\emph{(Reduced \zeds \risz)}\\
\Propeq
\begin{enumerate}
\item  \label{i1lemZerRed}
The \ri $\gA$ is reduced and \zedz. 
\item  \label{i2lemZerRed}
Every \idp is \idm
(i.e., $\forall a\in \gA,\;a\in\gen{a^{2}}$).
\item  \label{i3lemZerRed}
   Every \idp is generated by an \idmz.
\item  \label{i4lemZerRed}
   Every \itf  is generated by an \idmz.
\item  \label{i5lemZerRed}
 For every finite list  $(a_1,\ldots ,a_k)$ of \elts of $\gA$, there exist \orts \idms $(e_1,\ldots ,e_{k})$ such that for $j\in\lrbk$

\snic{\gen{a_1,\ldots ,a_j}=\gen{a_1e_1+\cdots +a_je_j}=\gen{e_1+\cdots+e_j}.}

\item  \label{i6lemZerRed}
   Every \id is \idmz.
\item  \label{i6bislemZerRed}
   The product of two \ids is always equal to their intersection.
\item  \label{i7lemZerRed}
   The \ri $\gA$ is a \qiri and every \ndz \elt is \ivz.
\end{enumerate}
\end{lemma}

\begin{fact}\label{factQoQiZed} 
\begin{enumerate}
\item Let $\gA$ be an arbitrary \riz. If $\Ann(a) = \gen{\vep}$  with $\vep$ \idmz,
then the \elt $b = a+\vep$ is \ndz and $ab = a^2$.
\item If $\gA$ is a \qiriz, $\Frac\gA$ is \zedr and every \idm of $\Frac\gA$ is in $\gA$.
\end{enumerate}
\end{fact}
\begin{proof}
\emph {1.}
Work modulo $\vep$ and modulo $1-\vep$.

\emph {2.}
For some $a\in\gA$, we must find $x\in\Frac\gA$ such that $a^2x=a$.\\
Let $b=a+(1-e_a)\in\Reg\gA$, then $ab=a^{2}$, and we take $x=b^{-1}$.
\\
Now let $a/b$ be an \idm of $\Frac\gA$. We have $a^2=ab$.\\
\hspace*{.1em} --- Modulo ${1-e_a}$, we have $b=a$ and $a/b=1=e_a$ (because $a$ is \ndzz).\\ 
\hspace*{.1em} --- Modulo ${e_a}$, we have $a/b=0=e_a$ (because $a=0$). In short, $a/b=e_a$. 
\end{proof}

\begin{fact}
\label{factZerRedCoh}
A \zedr \ri is \cohz.
It is \fdi \ssi there is an \egt to zero test for the \idmsz.
\end{fact}

We also easily  obtain the following \eqvcsz.

\begin{fact}
\label{factZerRedConnexe}
For a \zed \ri $\gA$ \propeq
\begin{enumerate}
  \item $\gA$ is connected (resp.\,$\gA$ is connected and reduced).
  \item $\gA$ is local (resp.\,$\gA$ is local and reduced).
  \item $\Ared$ is integral (resp.\,$\gA$ is integral).
  \item $\Ared$ is a \cdi (resp.\,$\gA$ is a \cdiz).
\end{enumerate} 
\end{fact}

\subsubsect{Equational definition of  reduced \zed \risz}
{Equational \dfnz}

\vspace{2pt}
A not \ncrt commutative \ri satisfying
$$
\forall x\,\exists a\quad xax=x
$$
is often qualified as \ix{Von Neumann regular}. In the commutative case they are the reduced \zed \risz. We also call them \emph{absolutely flat \risz}, because they are \egmt \cares by the following \prtz: every \Amo is flat
(see Proposition~\ref{propEvcPlat}).\index{ring!absolutely flat ---}

In a commutative \riz, two \elts $a$ and $b$ are said to be \emph{quasi-inverse}
if we have
\begin{equation}\label{eqQuasiInv}\preskip-.2em \postskip.4em
a^2b=a,\quad \quad  b^2a=b
\end{equation}

We also say that $b$ is \und{the} \ix{quasi-inverse} of $a$. Indeed, we check that it is unique. That is, if $a^2b=a=a^2c$,
 $b^2a=b$  and  $c^2a=c$, then
$$\preskip.4em \postskip.4em
c-b=a(c^2-b^2)=a(c-b)(c+b)=a^2(c-b)(c^2+b^2)=0,
$$
since $ab=a^2b^2$, $ac=a^2c^2$ and $a^2(c-b)=a-a=0$.

Moreover, if $x^2y=x$, we check that $xy^2$ is the \qiv of $x$.

Thus \emph{a \ri is \zedr \ssi every \elt admits a \qivz.}

As the \qiv is unique, a \zedr \ri can be regarded as a \ri fitted with an additional unary law $a\mapsto a\bul$ subject to  axiom (\ref{eqQuasiInv}) with $a\bul$ instead of $b$.\\
Note that $(a\bul)\bul=a$ and $(a_1a_2)\bul=a_1\bul a_2\bul$.

\begin{fact}
\label{factQuasiInv}
A \zedr \ri $\gA$ is a \qiriz, with the \idm $e_a=aa\bul$: $\Ann(a)=\gen{1-e_a}$.
We have $\gA\simeq \gA[1/e_a]\times \aqo{\gA}{e_a}$.
In $\gA[1/e_a]$, $a$ is invertible, and in $\aqo{\gA}{e_a}$, $a$ is null.
\end{fact}

\subsubsection*{\Zed splitting lemma}

The following splitting lemma is almost \imdz. The proof resembles that of the quasi-integral splitting lemma \ref{thScindageQi}.

\begin{lemma}\label{thScindageZed}
Let $(x_i)_{i\in I}$ be a finite family of \elts in a \zed \ri $\gA$.
There exists a \sfio $(e_1, \ldots, e_n)$ such that in each component $\gA[1/e_j]$, each $x_i$ is nilpotent or \ivz.
\end{lemma}
%

\subsubsect{From \cdis to \zedr \risz}{\Elr \lgb machinery \num2:
from \cdis to \zedr \risz}

\Zedr \ris look a lot like finite products of \cdisz, and this manifests itself \prmt as follows.

\medskip 
{\bf \Elr \lgb machinery \num2.
}\label{MethodeZedRed}\imlg\\
{\it Most \algos that work with nontrivial discrete fields can be modified in order to work with reduced \zed \risz, by splitting the \ri into two components each time that the \algo written for  \cdis uses the \gui{is this \elt null or \ivz?} test.
In the first component the \elt in question is null, in the second it is \ivz.}

\medskip
\rems \label{remzedred1}
1) We used the term \gui{most} rather than \gui{all} since the statement of the result of the \algo for the \cdis must be written in a form that does not specify that a \cdi is connected.\\
2) Moreover, the same remark as the one we made on \paref{remMLGE1} concerning the \elr \lgb machinery \num1 applies here. The \algo given in the \cdi case must be sufficiently uniform in order to avoid leading to an infinite tree when we want to transform it into an \algo for the \zedr \risz.
\eoe

\smallskip We immediately give an application example of this machinery.

\pagebreak	

\begin{proposition}\label{exlZedBez}
\label{propZedBez} For a \ri $\gA$ \propeq
\begin{enumerate}
\item $\gA$ is a \zedr \riz.
\item $\gA[X]$ is a strict Bézout \qiriz.
\item $\gA[X]$ is a Bézout \qiriz.
\end{enumerate}
\end{proposition}
\begin{proof}
\emph{1 $\Rightarrow$ 2.}
For \cdis this is a classical fact: we use Euclid's extended \algo to compute in the form $g(X)=a(X)u(X)+b(X)v(X)$ a gcd of~$a(X)$ and $b(X)$. In addition, we obtain a matrix $\crmatrix{u& -b_1\cr v&a_1}$ with \deterz~$1$ which transforms $[\,a \;\;b\,]$ into $[\,g \;\;0\,]$.
This matrix is the product of matrices $\cmatrix{0&-1\cr 1&-q_i}$ where the $q_i$'s are the successive quotients.\\
Let us move on to the \zedr \ri case (so  a \qiriz).
First of all $\gA[X]$ is a \qiri as the annihilator of a \pol is the intersection of the annihilators of its \coes 
(see Corollary~\ref{corlemdArtin}~\emph{\iref{i2corlemdArtin}}), hence generated by the product of the corresponding \idmsz.
The \gui{strict Bézout} character of the \algo which has just been explained for \cdis a priori stumbles upon the obstacle of the non-invertibility of the leading \coes in the successive divisions. Nonetheless this obstacle is avoided each time by considering a suitable \idmz~$e_i$, the annihilator of the \coe to be inverted. In $\gA_i[1/e_i]$, (\hbox{where $\gA_i=\gA[1/u_i]$} is the \gui{current} \ri with a certain \idm $u_i$) the divisor \pol has a smaller degree than expected and we start again with the following \coez. In $\gA_i[1/f_i]$, ($f_i=1-e_i$ in $\gA_i$), the leading \coe of the divisor is invertible and the division can be executed.
In this way we obtain a computation tree whose leaves have the desired result. At each leaf the result is obtained in a localized \ri 
$\gA[1/h]$ for a certain \idm $h$, and the $h$'s at the leaves of the tree form a \sfioz. This allows us to glue together all the \egtsz.\footnote{For a more direct \demz, see Exercise~\ref{exoZerRedBez}.}

 \emph{3 $\Rightarrow$ 1.}
This results from the following lemma.

 \textit{\textbf{Lemma}}. \emph{For an arbitrary \ri $\gA$, if the \id $\gen{a,X}$ is a \idp of $\gA[X]$, then $\gen{a}=\gen{e}$ for a certain \idm $e$.}

Suppose that $\gen{a,X}=\gen{p(X)}$ with $p(X)q(X)=X$.
We therefore have  

\snic{\gen{a}=\gen{p(0)}$,  $\;p(0)q(0)=0\;$
and $\;1=p(0)q'(0)+p'(0)q(0),}

hence $p(0)=p(0)^2q'(0)$. Thus, $e=p(0)q'(0)$ is \idm
and $\gen{a}=\gen{e}$.
\end{proof}

\rem The notion of a \zedr \ri can be seen as the non-\noe analogue of the notion of a \cdiz, since if the \agB of the \idms is infinite, the \noet is lost. Let us illustrate this with the example of the \nstz, for which it is not a priori clear if the \noet is an essential ingredient or a simple accident.
A precise \cof statement of Hilbert's \nstz\ihi (weak form) 
is formulated as follows.

\emph{Let $\gk$ be a nontrivial \cdiz, 
$(\lfs)$ be a list of \elts of $\kuX$
,
and $\gA=\aqo\kuX\uf$ be the quotient \algz.
  Then, \hbox{either
 $1\in\gen{\lfs}$}, or there exists a quotient of $\gA$ that  
 is a nonzero finite dimensional \kevz.}

As the \dem is given by a uniform \algo (for further details see \Thref{thNstNoe} and Exercise \ref{exothNst1-zed}), we obtain by applying the \elr \lgb machinery \num2 
the following result, without disjunction, which implies the previous \nst for a nontrivial \cdi (this example also illustrates the first remark on \paref{remzedred1}). 
An \Amo $M$ is said to be \emph{quasi-free}%
\index{quasi-free!module}\index{module!quasi-free ---} if it is \isoc to a finite direct sum of \ids $\gen{e_i}$ with the $e_i$'s \idmz.
We can then in addition require that $e_ie_j=e_j$ if $j>i$, since for two \idms $e$ and $f$, we have

\snic{\gen{e}\oplus\gen{f}\simeq\gen{e\vu f}\oplus\gen{e\vi f}$, where
$e\vi f=ef$ and $e\vu f=e+f-ef.}

\emph{Let $\gk$ be a \zedr \riz, $(\lfs)$ be a list of \elts of $\kXn$ and $\gA$ be the quotient \algz. Then the \id $\gen{\lfs}\,\cap\, \gk$ is generated by an idempotent ${e}$, and by letting $\gk_1=\aqo{\gk}{e}$, there exists a quotient~$\gB$ of $\gA$ which is a quasi-free $\gk_1$-module,  the  natural \homo $\gk_1\to\gB$} being injective.
 \eoe

\subsubsect{\Fp modules over  \zedr \risz}{\Fp modules}

\begin{theorem}\label{thZerDimRedLib}
\emph{(The \zedr \ri paradise)}\\
Let $\gA$ be a \zedr \riz.
\begin{enumerate}
\item Every matrix is \eqve to a matrix in Smith form with \idms on the principal diagonal.
\item Every \mpf is quasi-free.
\item
Every \tf submodule  of a \mpf is a direct summand.
\end{enumerate}
\end{theorem}
\begin{proof}
The results are classical ones for the \cdi case (a \prco can be based on the pivot method).
The \elgbm \num2 then provides (for each of the three items)
the result separately in each $\gA[1/e_j]$,
after splitting the \ri into a product of localized \ris $\gA[1/e_j]$
 for a \sfio $(e_1, $\ldots$, e_k)$.
But the result is in fact formulated in such a way that it is globally true as soon as it is true in each of the components.
\end{proof}

\vspace{-.7em}
\pagebreak


\subsec{\Zed \sypsz}
\label{subsecSypZerdim}

In this subsection we study a particularly important example of a \zed \riz, provided by the quotient \algs associated with \zed \syps over  \cdisz.

Recall the context studied in Section~\ref{secChap3Nst} dedicated to Hilbert's \nstz.
If $\gK\subseteq\gL$ are \cdisz, and if $(\lfs)$ is a \syp in $\KXn=\KuX$,
we say that $(\xin)=(\uxi)$ is a {zero of $\uf$ in $\gL^n$} if the \eqns $f_i(\uxi)=0$ are satisfied.
\\
The study of the \vrt of the zeros of the \sys is closely related to that of the \emph{quotient \alg associated with the \sypz}, namely 
$$\preskip.4em \postskip.4em 
\gA=\aqo\KuX\uf=\Kux \quad (x_i \hbox{ is the class of } X_i \hbox{ in }\gA). 
$$
Indeed, it amounts to the same to take a zero $(\uxi)$ of the \syp in $\gL^n$ or to take a \homo of \Klgs  $\psi:\gA\to\gL$ ($\psi$~is defined by~\hbox{$\psi(x_i)=\xi_i$} for $i\in\lrbn$).
For $h\in\gA$, we write $h(\uxi)=\psi(h)$ for the \evn of $h$ at $\uxi$.

When $\gK$ is infinite, \Thref{thNstfaibleClass} gives us a \noep \inoe by a \lin \cdvz, and an integer $r\in\lrb{-1..n}$ satisfying the following \prts
(we do not change the name of the variables, which is a slight abuse). 
\begin{enumerate}
\item If $r=-1$, then $\gA=0$, \cad $1\in\gen{\uf}$.
\item If $r=0$, each $x_i$ is integral over $\gK$, and $\gA\neq 0$.
\item If $0<r<n$, then $\KXr\,\cap\,\gen{\uf}=0$ 
and the $x_{i}$ for ${i\in\lrb{r+1..n}}$ are integral over $\Kxr$ (which is \isoc to $\KXr$).
\item If $r=n$, $\gen{\uf}=0$ and $\gA=\KuX$
\end{enumerate}

\begin{lemma}\label{lemSypZD1}
\emph{(Precisions on \Thref{thNstfaibleClass})}%
\begin{enumerate}
\item In the case where $r=0$, the quotient \alg $\gA$ is finite over $\gK$.
\item If the quotient \alg $\gA$ is finite over $\gK$, it is \stfe over~$\gK$,
and it is a \zed \riz.
We then say that \emph{the \syp is \zedz}.

\item If the \ri $\gA$ is \zedz, then $r\leq 0$.
\item Every \stfe \alg over the \cdi $\gK$ can be seen as (is \isoc to) the quotient \alg of a \zed \syp over $\gK$.
\end{enumerate}
\end{lemma}
%
\begin{proof}
\emph{1.} 
Indeed, if $p_i(x_i)=0$ for $i\in\lrbn$, 
the \algz~$\gA$ is a quotient of 
$$\preskip.3em \postskip.3em\ndsp 
\gB=\KuX\big/\geN{\big(p_i(X_i)\big)\phantom{\!\!)}_{i\in\lrbn}}, 
$$
which is a finite dimensional \Kevz.

\emph{2.} We start as we did in item \emph{1}. To obtain the \alg $\gA$, it suffices to take the quotient of $\gB$ by the \id $\gen{f_1(\uz),\dots,f_s(\uz)}$ (where the $z_i$'s are the classes of $X_i$'s in $\gB$). We easily see that this \id is a \tf \lin subspace of $\gB$, so the quotient is again a finite dimensional \Kevz. Thus, $\gA$ is \stfe over $\gK$.\\
Let us show that $\gA$ is \zedz. Every $x\in\gA$ annihilates its \polminz, say $f(T)$, so that we have an \egt $x^k\big(1+xg(x)\big)=0$
(multiply $f$ by the inverse of the nonzero \coe of lower degree).

\emph{4.} The \alg $\gA$ is generated by a finite number of \elts $x_i$  
(for example a basis as a \Kevz),  each of 
which annihilate their \polminz, say $p_i(T)$. Thus $\gA$ is a quotient of an \alg 
$$\preskip.4em \postskip.4em 
\gB=\Aqo\KuX{\big(p_i(X_i)\big)_{i\in\lrbn}}=\gA[\zn]. 
$$
The corresponding surjective morphism, 
from $\gB$ over $\gA$, is a \ali whose kernel can be computed (since $\gA$ and $\gB$ are finite dimensional \evcsz), for example by specifying a \sgr $\big(g_1(\uz),\dots,g_\ell(\uz)\big)$.
\\
In conclusion,
 the \alg $\gA$ is \isoc to the quotient \alg associated with the \syp
$\big(p_1(X_1),\dots,p_n(X_n),g_1(\uX),\dots,g_\ell(\uX)\big)$.

\emph{3.}  This point results from the following two lemmas.
\end{proof}

\rem Traditionally, we reserve the term \zed \syp to the $r=0$ case, but the quotient \alg is also \zed when $r=-1$.
\eoe

\begin{lemma}\label{lemnonZR}
If the \ri $\gC[\Xr]$ is \zed with $r>0$, then the \ri $\gC$ is trivial.  
\end{lemma}
%
\begin{proof}
We write $X_1^m\big(1-X_1P(\Xr)\big)=0$. The \coe of $X_1^m$ is both equal to $0$ and to $1$.
\end{proof}
%

\begin{lemma}
\label{lemZrZr2}
Let~$\gk\subseteq \gA$ and~$\gA$ be integral over~$\gk$. If
$\gA$ is a \zed \riz,
$\gk$ is a \zed \riz.
\end{lemma}
\begin{proof}
Let $x\in\gk$, then we have a $y\in\gA$ such that $x^{k}=yx^{k+1}$.
Suppose for example that $y^3+b_2y^2+b_1y+b_0=0$ with $b_i\in\gk$.
\\
Then, $x^{k}=yx^{k+1}=y^2x^{k+2}=y^3x^{k+3}$, and so

\snic{\arraycolsep2pt\begin{array}{rclcl}
0& =  &  (y^3+b_2y^2+b_1y+b_0)x^{k+3} &   &   \\[1mm]
& =   & x^{k}+b_2x^{k+1}+b_1x^{k+2}+b_0x^{k+3}  & =  & x^{k}\big(1+x(b_2+b_1x+b_0x^2)\big).
\end{array}}

\vspace{-.5em}
\end{proof}

\begin{theorem}
\label{thSPolZed} \emph{(\Zed \sys over a \cdiz)}\\
Let~$\gK$ be a \cdi and $(\lfs)$ in $\KXn=\KuX$.\\
Let~$\gA=\aqo\KuX\uf$ be the quotient \alg associated with this \sypz.\\
\Propeq
\begin{enumerate}
\item $\gA$ is finite over~$\gK$.
\item $\gA$ is \stfe over~$\gK$.
\item $\gA$ is a \zed \riz.
\end{enumerate}
If~$\gK$ is contained in a \cdacz~$\gL$, these \prts are also \eqves to the following.
\begin{enumerate} \setcounter{enumi}{3}
\item The \syp has a finite number of zeros in~$\gL^n$.
\item The \syp has a bounded number of zeros in~$\gL^n$.
\end{enumerate}%
\index{zero-dimensional@\zedz!\sypz}\index{polynomial system!zero-dimensional ---}
\end{theorem}
\begin{proof} When $\gK$ is infinite, we obtain the \eqvcs by applying Lemma~\ref{lemSypZD1} and \Thref{thNstfaibleClass}.
\\
In the \gnl case, we can \egmt obtain a \noep \inoe by using a 
(not \ncrt \linz) 
\gnl \cdv as described in Lemma~\ref{lemNoether}
(see \Thref{thNstNoe}).  
\end{proof}

A variation on the previous \tho is given in \Thref{thNst0}.

\medskip 
\rem Rather than using a non-\lin \cdv as proposed in the previous \demz, we can resort to using the technique of \gui{changing the base field.} This works as follows. Consider an infinite field $\gK_1\supseteq \gK$, for example $\gK_1=\gK(t)$, or an \cac $\gK_1$ containing $\gK$ if we know how to construct one. 
Then the \eqvc of items~\emph{1}, \emph{2} and \emph{3} is assured for the \alg $\gA_1$ for the same \syp seen on $\gK_1$. The \alg $\gA_1$ is obtained from $\gA$ by \eds from $\gK$ to $\gK_1$.
It remains to see that each of the three items is satisfied for $\gA$ \ssi it is satisfied for $\gA_1$. A task we leave to the reader.\footnote{See on this subject \Thosz~\ref{thSurZedFidPlat}, \ref{propFidPlatTf} and~\ref{propFidPlatPrAlg}.} 
\eoe


\begin{theorem}\label{thStickelberger}  \emph{(Stickelberger's \thoz)}\\
Same context as in \Thref{thSPolZed}, now with $\gK$ being an \cacz.%
\index{Stickelberger!theorem}
\begin{enumerate}
\item The \syp admits a finite number zeros over $\gK$.  
\\
We write them as
$\uxi_1$, \ldots, $\uxi_\ell$. 
\item For each $\uxi_k$ there exists an \idm $e_k\in\gA$
satisfying $e_k(\uxi_j)=\delta_{j,k}$ (Kronecker symbol)
for all $j\in\lrb{1..\ell}$. 
\item The \idms $(e_1,\dots,e_\ell)$ form a \sfioz.
\item Each \alg $\gA[1/e_k]$ is a \zed \alo (every \elt is \iv or nilpotent).
\item
Let $m_k$ be the dimension of the \Kev $\gA[1/e_k]$. \\
We have
$[\gA:\gk]=\som_{k=1}^\ell m_k$ and
for all $h\in\gA$ we have

\snic{\rC{\gA\sur\gk}(h)(T)=\prod\nolimits_{k=1}^\ell\big(T-h(\uxi_k)\big)^{m_k}.
}

In particular, $\Tr\iAk(h)=\som_{k=1}^\ell m_kh(\uxi_k)$ and 
$\rN\iAk(h)=\prod_{k=1}^\ell h(\uxi_k)^{m_k}$.
\item Let $\pi_k:\gA\to\gK,\;h\mapsto h(\uxi_k)$ be the \evn at $\uxi_k$,
and $\fm_k=\Ker\pi_k$. Then $\gen{e_k-1}=\fm_k^{m_k}$ and $\fm_k=\sqrt{\gen{e_k-1}}$.

\end{enumerate}
\end{theorem}
%
\begin{proof} Let $V=\big\{\uxi_1,\dots,\uxi_\ell\big\}$ be the \vrt of zeros of the \sys in $\gK^n$. 

\emph{2} and \emph{3.}
We have multivariate Lagrange interpolating \pols 
 $L_k\in\KuX$ which satisfy $L_k(\uxi_j)=\delta_{j,k}$. Consider the~$L_k$'s as \elts of $\gA$.\\ 
Since $\gA$ is \zedz, there exist an integer $d$ and an \idm $e_k$ with $\gen{e_k}=\gen{L_k}^d$, therefore $e_kL_k^d=L_k^d$ and $L_k^db_k=e_k$ for a certain $b_k$.
This implies that $e_k(\uxi_j)=\delta_{j,k}$. \\
For $j\neq k$, $e_je_k$ is null over $V$, so by the \nstz, $e_je_k$ is nilpotent in~$\gA$. As it is an \idmz, $e_je_k=0$.\\
The sum of the $e_j$'s is therefore an \idm $e$. This \elt vanishes nowhere, \cad it has the same zeros as $1$. By the \nstz, we obtain $1\in\sqrt{\gen{e}}$. Thus $e=1$ because it is an \iv \idm of $\gA$.

\emph{4.} The \Klg $\gA_k=\gA[1/e_k]=\aqo\gA{1-e_k}$ is the quotient \alg associated with the \syp $(\lfs,1-e_k)$ which admits $\uxi_k$ as its only zero.
Consider an arbitrary \elt $h\in\gA_k$. By reasoning as in the previous item, we obtain by the \nst that if~$h(\uxi_k)=0$, then~$h$ is nilpotent, and if~$h(\uxi_k)\neq 0$, 
then~$h$ is \ivz.

\emph{5.} Since $\gA\simeq\prod_{k=1}^\ell\gA_k$, it suffices to prove that for $h\in\gA_k$, we have the \egt $\rC{\gA_k\sur\gk}(h)(T)=\big(T-h(\uxi_k)\big)^{m_k}$. We identify $\gK$ with its image in~$\gA_k$. The \elt $h_k=h-h(\uxi_k)$ vanishes in $\uxi_k$, 
so it is nilpotent. If~$\mu$ designates  multiplication by $h_k$ in $\gA_k$, $\mu$ is a nilpotent \endoz.
With respect to a suitable basis, its matrix is 
strictly lower triangular
and that of the multiplication by $h$ is triangular with $h(\uxi_k)$'s on the diagonal, therefore its \polcar is $\big(T-h(\uxi_k)\big)^{m_k}$.

\emph{6.} We clearly have $e_k-1\in\fm_k$. If $h\in\fm_k$, the \elt $e_kh$ is null everywhere over $V$, so nilpotent. Therefore $h^Ne_k=0$ for a certain $N$ and $h\in\sqrt{\gen{e_k-1}}$. To show that $\fm_k^{m_k}=\gen{e_k-1}$, we can locate ourselves in~$\gA_k$, where $\gen{e_k-1}=0$. In this \riz, the \id $\fm_k$ is a \Kev of dimension~$m_k-1$. The successive powers of $\fm_k$ then form a decreasing sequence of finite dimensional \Ksvsz, which stabilizes as soon as two consecutive terms are equal.
Thus $\fm_k^{m_k}$ is a \tf strict \idm \idz, therefore null. 
\end{proof}
\rems~\\ 1) The fact that the \syp is \zed results from a rational computation in the field of \coes (in a \noepz ing \inoe or computation of a \bdgz).

\smallskip 
2) Item~\emph{5} 
of Stickelberger's \tho allows us to compute all the useful information on the zeros of the \sys by basing ourselves on the only trace form. In addition, the trace form can be computed in the field of \coes of the \pols of the \sysz.
This has important applications in computer algebra (see for example~\cite{BPR}).
\eoe

\medskip  For examples, consult Exercise~\ref{exoFreeAlgebraPresentation}
and \Pbmz~\ref{exoDimZeroXcYbZa}. For a purely local study of the isolated zeros, see Section~\ref{secExlocGeoAlg}.

\section{Fitting \ids}
\label{sec Fitt} 

The theory of the \idfs of \mpfs is an extremely efficient computing machinery from a theoretical \cof point of view. It has an \gui{\eli theory} side in so far as it is entirely based on computations of \detersz, and it more or less disappeared for a while from the literature under the influence of the idea that we had to \gui{eliminate the \eliz} 
to escape the quagmire 
of computations whose meaning seemed unclear.

The \idfs are becoming fashionable once again and it is for the best.
For more details, please consult~\cite{Nor}.

\penalty-2500
\subsec{Fitting \ids of a \mpfz}
\label{sec Fitt pres fin}

\begin{definition}\label{def ide fit} \\
If $G\in\gA^{q\times m}$ is a \mpn of an \Amo $M$ given by~$q$ \gtrsz, the \emph{\idfs of} $M$ are the \ids

\snic{\cF_{\gA,n}(M)=\cF_{n}(M):= \cD_{\gA,q-n}(G)}

where $n$ is an arbitrary integer.\index{ideal!Fitting ---}%
\index{Fitting!ideal}
\end{definition}

This \dfn is legitimized by the following easy but fundamental lemma.

\begin{lemma}
\label{lemFitInv}
The \idfs of the \mpf $M$ are well-defined, in other words these \ids do not depend on the chosen \pn $G$ for $M$.
\end{lemma}
\begin{proof}
To prove this lemma we must essentially show that the \ids $\cD_{q-n}(G)$ do not change,
\begin{enumerate}
\item on the one hand, when we add a new syzygy, a \coli of the already 
present syzygies,
\item on the other hand, when we add a new \elt to a \sgrz, with a syzygy that expresses this new \elt in relation to
the previous \gtrsz. 
\end{enumerate}
The details are left to the reader.
\end{proof}

We immediately have the following facts.

\pagebreak	

\begin{fact}\label{fact.idf inc} 
For every \mpf $M$ with $q$ \gtrsz, we have the inclusions

\snic{\gen{0} = \cF_{-1}(M) \subseteq \cF_{0}(M) \subseteq \cdots \subseteq
\cF_q(M)= \gen{1}.}

If $N$ is a \pf quotient module of $M,$ we have the %
inclusions $\cF_{k}(M) \subseteq
\cF_{k}(N)$ for all $k\geq 0$.
\end{fact}

\rem In particular, if $\cF_r(M)\neq \gen{1}$ the module $M$ cannot be generated by $r$ \eltsz. We will see (lemma of the number of local \gtrs \paref{lemnbgtrlo}) that the meaning of the \egt $\cF_r(M)= \gen{1}$ is that the module is \emph{locally} generated by $r$ \eltsz.
\eoe

\begin{fact}\label{fact.idf libre} 
Let $M$ be a rank $k$ free \Amoz. Then,

\snic{\cF_{0}(M) = \cdots = \cF_{k-1}(M)  =  \gen{0} \subseteq
\cF_k(M)=\gen{1}.}

More \gnltz, if $M$ is quasi-free \isoc to $\bigoplus_{1\leq i\leq k}\gen{f_i}$, where the~$f_i$'s are \idms such that $f_if_j=f_j$ if $j>i$,
then $\cF_k(M)=\gen{1}$ and~$\cF_i(M)=\gen{1-f_{i+1}}$ for $0\leq i< k$.
\end{fact}

Note that this provides a clever \dem 
that if a module is free with two distinct ranks, the \ri is trivial.

\ms\exls ~
\\ 1. For a finite Abelian group $H$ considered as a $\ZZ$-module, the \id $\cF_0(H)$ is generated by the order of the group whilst the annihilator is generated by its exponent. 
In addition, the structure of the group is entirely \caree by its \idfsz. A \gnn is given in Exercise~\ref{exoFitt0}.

\noi 2. Let us reuse the \Bmo $M$ of Example on~\paref{belexemple}.
The computation gives the following results.
\begin{itemize}
\item For $M$: $\cF_0(M)=0$, $\cF_1(M)=\fb$ and $\cF_2(M)=\gen{1}$,
\item for $M'=M\te M$: $\cF_0(M')=0$, $\cF_1=\fb^3$, $\cF_2=\fb^2$, $\cF_3=\fb$ and $\cF_4=\gen{1}$,
\item for $M''=\gS^2(M)$: $\cF_0 (M'')=0$, $\cF_1=\fb^2$, $\cF_2=\fb$ and $\cF_3=\gen{1}$,
\item for $\Al2M$: $\cF_0(\Al2M)=\fb$ and $\cF_1(\Al2M)=\gen{1}$.
 \eoe
\end{itemize}

\begin{fact}\label{fact.idf.change}\label{fact.idf loc} 
\emph{(Changing the base \riz)}\\
Let $M$ be a \pf \Amoz, $\rho:\gA\rightarrow \gB$ be a \homo of \risz, and $\rho\ist(M)$ be the $\gB$-module obtained by \eds to $\gB$.
We have for every integer $n\geq 0$ the \egt 
$\gen{\rho\big(\cF_{n}(M)\big)} = \cF_{n}\big(\rho\ist(M)\big).$
\\
In particular, if $S$ is a \moz, we have $\cF_{n}(M_S) =  \big(\cF_{n}(M)\big)_S.$
\end{fact}

The two following facts are less obvious.

\pagebreak	

\begin{lemma} \emph{(Annihilator and first \idfz)}\label{fact.idf.ann} 
\\
Let $M$ be a \pf \Amo generated by $q$ \eltsz, we have

\snic{\Ann(M)^q\subseteq \cF_{0}(M) \subseteq \Ann(M).}
\end{lemma}
\begin{proof}
Let $(\xq)$ be a \sgr of $M$, $X={[\,x_1\;\cdots\;x_q\,]}$ and $G$ a \mpn for $X$.
Let $a_1$, \ldots, $a_q\in\Ann(M)$. Then, the diagonal matrix $\Diag(a_1,\ldots ,a_q)$ has as its columns \colis of the columns of $G$, so its \deter $a_1\cdots a_q$ belongs to~$\cF_{0}(M)$. This proves the first inclusion. \\
Let $\delta$ be a minor of order $q$ extracted from $G$. We will show that $\delta \in\Ann(M)$, hence the second inclusion. If $\delta$ corresponds to a submatrix $H$ of $G$ we~have~$X\,H=0$, therefore $\delta X=0$, and this indeed means that $\delta \in\Ann(M)$.
\end{proof}

\begin{fact}\label{fact.idf.sex}  
\emph{(Fitting \ids and exact sequences)}\\
Let $0\rightarrow N\rightarrow M\rightarrow P\rightarrow 0$ be an exact sequence of \mpfsz. For all $p\geq 0$ we~have

\snic{\cF_p(M)\;\supseteq\;\sum_{r\geq0,s\geq0,r+s= p}\cF_r(N)\cF_{s}(P) ,}

and if $M\simeq N \oplus P$, the inclusion is an \egtz.
\end{fact}
\begin{proof}
We can consider that $N\subseteq M$ and $P=M/N$. We use the notations of item~\emph{3} of Proposition~\ref{propPfSex}. We have a \mpn $D$ of~$M$ which is written \gui{in a triangular form}
$D=\cmatrix{A&C\cr 0&B}.$
Then every product of a minor of order $k$ of $A$ and of a minor of order $\ell$ of $B$ is equal to a minor of order $k+\ell$ of $D$.
This implies the stated result for \idfsz. 
\\
The second case is clear, with $C=0$.
\end{proof}

\exl On the \poll \ri $\gA=\ZZ[a,b,c,d]$, let us consider the module
$M=\gA g_1+ \gA g_2=\Coker F$ where $F=\cmatrix{a&b\cr c&d}$. Here $g_1$ and $g_2$ are images of the natural basis $(e_1,e_2)$ of $\Ae2$. Let $\delta=\det(F)$. \\
It is easily seen that $\delta\, e_1$ is a basis of the submodule  $\Im F\cap e_1\gA$
of $\Ae2$. 
\\
Let $N=\gA g_1$ and $P=M/N$.
Then the module $N$ admits the \mpn $[\,\delta\,]$ for the \sgr $(g_1)$ and $P$ admits the \mpn $[\,c\;d\,]$ for the \sgr $(\ov{g_2})$. 
Consequently, we get $\cF_0(M)=\cF_0(N)=\gen{\delta}$ and $\cF_0(P)=\gen{c,d}$. So the inclusion
$\cF_0(N)\cF_0(P)\subseteq \cF_0(M)$ is strict.   
\eoe

\subsec{Fitting \ids of a \mtfz}
\label{sec Fitt tf} 

We can generalize the \dfn of the \idfs to an arbitrary \mtf $M$ as follows.
If $(\xq)$
is a \sgr of $M$ and if $X=\tra{[\,x_1\;\cdots\;x_q\,]}$,
we define $\cF_{q-k}(M)$ as
the \id generated by all the minors of order $k$ of every matrix $G\in\Ae {k\times q}$ satisfying~$GX=0$.
An alternative \dfn is that each $\cF_j(M)$ is the sum of all the $\cF_j(N)$'s
where~$N$ ranges over the \mpfs that are surjectively sent onto~$M$.

This shows that the \ids defined thus do not depend on the considered \sgrz.

\ss The following remark is often useful.

\begin{fact}
\label{facttfpf}
Let $M$ be a \tf \Amoz.
\begin{enumerate}
\item If $\cF_k(M)$ is a \itfz, then $M$ is the quotient of a \mpf $M'$
for which  $\cF_k(M')=\cF_k(M)$.
\item If all the \idfs are \tfz,  then $M$ is the quotient of a \mpf $M'$
having the same \idfs as~$M$.
\end{enumerate}
\end{fact}

\section{Resultant \id}
\label{subsecIdealResultant}

In what follows, we consider a \ri $\gk$ that we do not assume to be discrete.
The resultant of two \pols is at the heart of \eli theory.
If $f$, $g\in\kX$ with $f$ monic,
the basic \eli lemma \paref{LemElimAffBasic} 
can be read in the \alg $\gB=
\aqo{\kX}{f}$ by writing
$$ 
\rD_\gB(\ov g)\,\cap\, \gk=\rD_\gk\big(\Res_X(f,g)\big).
$$
It can then be \gnee with the following result, which can be regarded as a very precise formulation of the Lying Over (see Lemma~\ref{lemLingOver}).

\CMnewtheorem{lemeligen}{\Gnl \eli lemma}{\itshape}
\begin{lemeligen}\label{LemElimAff}%
\index{general elimination lemma}%
\index{elimination!general --- lemma}%
\index{elimination!ideal}\index{ideal!elimination ---}~  
\begin{enumerate}
\item Let $\gk\vers{\rho}\gC$ be an \alg which is a  \kmo generated by $m$ \eltsz, $\fa=\cF_{\gk,0}(\gC)$ its first \idf and $\fc=\Ker\rho$.
Then
\begin{enumerate}
\item $\fc=\Ann_\gk(\gC)$,
\item \fbox{$\fc^m\subseteq\fa\subseteq\fc$} and so \fbox{$\rD_\gk(\fc)=\rD_\gk(\fa)$},
\item if by some \eds $\varphi:\gk\to\gk'$ we obtain the \alg $\rho':\gk'\to\gC'$,
then the \id $\fa':=\cF_0(\gC')$ is equal to $\varphi(\fa)\gk'$ and as a $\gk'$-module, it is \isoc to $\gk'\otimes_\gk\fa\simeq\varphi\ist(\fa)$.
\end{enumerate}
\item Let $\gB\supseteq\gk$ be a \klg which is a free \kmo of rank $m$, and $\fb$ be a \itf of~$\gB$.
\begin{enumerate}
\item The \eli \id $\fb\,\cap\,\gk$ is the kernel of the canonical \homo $\rho:\gk\to\gB\sur\fb$, \cad the annihilator of the \kmoz~$\gB\sur\fb$.

\item The \kmo $\gB\sur\fb $ is \pf and we have

\snic{
\fbox{$(\fb\,\cap\,\gk)^m\subseteq \cF_0 (\gB \sur \fb )\subseteq\fb\,\cap\,\gk$}  
\hbox{ and }  
\fbox{$\rD_\gB(\fb)\,\cap\, \gk=\rD_\gk \big( \cF_0 (\gB \sur \fb )\big)$}
.
}

We denote by $\fRes(\fb):=\cF_{\gk,0}(\gB\sur\fb )$ what we call the \emph{resultant \id of $\fb$}.
\end{enumerate}
\end{enumerate}
\end{lemeligen}

\begin{proof}
\emph{1a} and \emph{1b.} Indeed, $a\in\gk$ annihilates $\gC$ \ssi
it annihilates $1_{\gC}$, \ssi $\rho(a)=0$. 
The desired double inclusion is therefore given by Lemma~\ref{fact.idf.ann} (also valid for \mtfsz).

\emph{1c.} The \idfs are well-behaved under \edsz.

\emph{2.} Apply item~\emph{1} with $\gC=\gB\sur\fb$.
\end{proof}

\rems  
 1) The resultant \id in item \emph{2} can be \prmt described as follows. If $\fb=\gen{b_1,\ldots ,b_s}$ we consider the \emph{\gnee Sylvester mapping}

\snic{\psi:\gB^s\to\gB,\quad (y_1,\ldots ,y_s)\mapsto\psi(\uy)=\som_iy_ib_i.}

It is a \kli between free \kmos of ranks $ms$ and $m$.
Then, we have $\fRes(\fb)=\cD_m(\psi)$. \rdb

2) There are many \gtrs for the \id $\fRes(\fb)$.
Actually, there exist diverse techniques to decrease the number of \gtrs by replacing~$\fRes(\fb)$ by a \itf having considerably fewer \gtrs but having the same nilradical.
On this subject see the work given in Section~\ref{secChap3Nst}, especially Lemma~\ref{lemElimParametre},  
the results of Chapter~\ref{chapKrulldim} on the number of radical \gtrs of a radically \tf \id (\Thref{thKroH}),
and the paper~\cite{DiGLQ}.
\eoe%
\index{Sylvester!generalized --- mapping}%

\smallskip 
Now here is a special case of the \gnl \eli lemma.
This \tho completes Lemma~\ref{lemElimPlusieurs}.

\begin{theorem}
\label{thElimAff} \label{corLemmeElim} \hspace*{-.5em}\emph{(\Agq \eli \thoz: the resultant \idz)}%
\index{ideal!resultant ---}%
\index{resultant!ideal}%
\index{elimination!\agq --- theorem}%
\\
Let $(f,g_1,\ldots,g_r)$ be \pols of $\kX$ with $f$ \mon of degree $m$.
We let 
$$
\ff=\gen{f,g_1,\ldots,g_r}\subseteq\kX \;\hbox{and}\;\gB=\aqo{\kX}{f}.
$$
Let $\psi:\gB^r\to\gB$ be the \emph{\gnee Sylvester mapping} defined by
$$\preskip.4em \postskip.3em\ndsp
(y_1,\ldots ,y_r)\mapsto \psi(\uy)=\som_iy_i\ov{g_i}.
$$
It is a \kli between free \kmos of respective ranks $mr$ and $m$.
Let $\fa$ be the \idd $\cD_m(\psi)$. 
\begin{enumerate}
\item We have $\fa=\cF_{\gk,0}(\kX\sur\ff),$  and 
$$\preskip.4em \postskip.3em
(\ff\,\cap\,\gk)^m\subseteq\fa\subseteq\ff\,\cap\,\gk,\quad \mathit{and\;so}\quad
\rD_\kX(\ff)\,\cap\, \gk=\rD_\gk(\fa).
$$
\item Assume that $\gk=\gA[Y_1,\ldots,Y_q]$ and that $f$ and the $g_i$'s are of total degree~$\leq d$ in $\gA[\uY,X]$. Then the \gtrs of $\cD_m(\psi)$ are of total degree~$\leq d^2$ in $\gA[\uY]$. 
\item The \id $\fa$ does not depend on $\ff$ (under the sole assumption that $\ff$ contains a \poluz).  
We call it \emph{the resultant \id of $\ff$ \wrt the \idtr $X$} and we denote it by $\fRes_X(f,g_1,\ldots,g_r)$ or $\fRes_{X}(\ff)$, or~$\fRes(\ff)$. 
\item If by some \eds $\theta:\gk\to\gk'$ we obtain the \id $\ff'$ of~$\gk'[X]$,
then the \id $\fRes_{X}(\ff')\subseteq\gk'$ is equal to $\theta\big(\fRes_{X}(\ff)\big)\gk'$,
and as a module it is \isoc to $\gk'\otimes_\gk\fRes_{X}(\ff)\simeq\theta\ist\big(\fRes_{X}(\ff)\big)$.

\end{enumerate}
\end{theorem}
NB. Consider the basis $\cE=(1,\ldots ,X^{m-1})$ of $\gB$ over $\gk$.
Let $F\in\gk^ {m\times mr}$ be the matrix of $\psi$ 
for the deduced bases of $\cE$. 
Its columns are the $X^jg_k \mod f$ for $j\in\lrb{0..m-1}$, $k\in\lrbr$ written over the basis $\cE$.
We say that $F$ is a \emph{\gnee Sylvester matrix}. By \dfn we have $\fRes_{X}(\ff)=\cD_m(F)$.%
\index{matrix!generalized Sylvester ---}
\eoe

\begin{proof} Let $\fb=\ff \mod f=\gen{\ov{g_1},\ldots,\ov{g_r}}\subseteq\gB$.
We apply items~\emph{2} and~\emph{1c} of the \gnl \eli lemma by noticing that $\kX\sur\ff\simeq\gB\sur\fb$, with~$\ff\cap\gk=\fb\cap\gk$.
\end{proof}

\rem 
Thus \Thref{thElimAff} establishes a very narrow link between the \eli \id and the resultant \idz. 
The advantages introduced by the resultant \id over the \eli \id are the following\index{ideal!elimination ---}
\begin{itemize}
\item the resultant \id is \tfz,
\item its computation is \emph{uniform},
\item it is well-behaved under \edsz. 
\end{itemize}
Note that in the case where $\gk=\gK[Y_1,\ldots,Y_q]$, for $\gK$  a \cdiz, the \eli \id is also \tf but its computation, for instance via  \bdgsz, is not uniform.
\\
However, the resultant \id is only defined when $\ff$ contains a \polu and this limits the scope of the \thoz.
\eoe

\Exercices

\begin{exercise}
\label{exo3Lecteur}
{\rm  We recommend that the \dems which are not given, or are sketched, or
left to the reader,
etc, be done.
But in particular, we will cover the following cases.
\begin{itemize}\itemsep0pt
\item \label{exo lem pres equiv} 
Give a detailed \dem of Lemma~\ref{lem pres equiv}.
\item \label{exopropCoh1.2} 
Explain why Propositions~\ref{propCoh1} and~\ref{propCohfd1}
(when we take $\gA$ as an \Amo $M$)
can be read in the form of \Thref{propCoh2}.

\item \label{exopropCohpfKer}
Prove Propositions~\ref{propPfInter} and~\ref{propCohpfKer}.
Give a detailed \dem of Propositions~\ref{propPftens} and~\ref{propPfExt}. Show that
$\gA\sur{\fa}\te_\gA \gA\sur{\fb}\simeq \gA\sur{(\fa+\fb)}.$

\item \label{exobelexample}
Justify the statements contained in the Example on \paref{belexemple}.

\item \label{exoZerRed}
Prove Lemmas or Facts~\ref{lemZerloc}, \ref{lemZerRed},
\ref{factZerRedCoh} and~\ref{factZerRedConnexe}.

\item \label{exothZerDimRedLib}
Give \algos for the three items of \Thref{thZerDimRedLib}.

\item \label{exoFittsex}  \label{exoFitt5}
Prove Fact~\ref{facttfpf}.
\end{itemize}}
\end{exercise}

\vspace{-1em}
\begin{exercise}
\label{exoptfpf}
{\rm  Let $M\subseteq N$ be \Amos with $M$ as  direct factor in $N$.
Prove that if~$N$ is \tf (resp.\,\pfz), then so is $M$.
 }
\end{exercise}

\vspace{-1em}
\begin{exercise}\label{exoAXmodule} {(Structure of an $\AX$-module over $\Ae n$ associated with
$A\in\Mn(\gA)$)}\\
{\rm  
Let $\gA$ be a commutative ring and $A\in\Mn(\gA)$. We give $\Ae n$ the structure of an $\AX$-module by letting 

\snic{Q \cdot x = Q(A) \cdot x\hbox{ for }Q \in \AX\hbox{ and }x \in \Ae n.}

We aim to give a \mpn for this $\AX$-module. This generalizes Example~3)
on \paref{exl1pf} given at the beginning of Section~\ref{sec pf chg}, where $\gA$ is a \cdiz. 

Let $\theta_A : \AX^n \twoheadrightarrow \Ae n$ be the unique $\AX$-morphism which transforms the canonical basis of $\AX^n$ into that of $\Ae n$. 
By labeling these two canonical bases by the same name $(e_1, \ldots, e_n)$, 
$\theta_A$ therefore transforms $Q_1 e_1 + \cdots + Q_n e_n$ into $Q_1(A)\cdot e_1 + 
\cdots + Q_n(A)\cdot e_n$. We will show that the sequence below is exact

\snic {
\AX^n\vvvvvers{X\In-A} {\AX^n} \vvers{\theta_A} \Ae n \to 0
}

In other words $\Ae n$ is a \pf $\AX$-module and $X\In-A$ is a \mpn for the \sgr $(e_1, \ldots, e_n)$.

\emph {1.}
Show that we have a direct sum of $\gA$-modules $\AX^n = \Im(X\In-A) \oplus \Ae n$.

\emph {2.}
Conclude the result.}
\end{exercise}

\vspace{-1em}
\begin{exercise}\label{exoTenseursNuls}
 (Description of the null tensors)\\
{\rm  
Let $M$ and $N$ be two \Amos and $z=\sum_{i\in \lrbn}x_i\te y_i\in M\te N$.

\emph{1.} Show that $z=0$ \ssi there exists a \tf submodule~$M_1$ of $M$ 
such that we have $\sum_{i\in\lrbn}x_i\te y_i=_{M_1\te N}0$.

\emph{2.} We write $M_1=\gA x_1+\cdots+\gA x_p$ where $p\geq n$.
Let $y_k=_N0$ for $n<k\leq p$. Use the null tensor lemma with the \egt $\sum_{i\in\lrbp}x_i\te y_i=_{M_1\te N}0$ to give a \carn of the null tensors in the \gnl setting.}
\end{exercise}

\vspace{-1em}
\begin{exercise}
\label{exoEdsQuot}
{\rm  Let $M$ be an \Amoz, $\fa$ be an \id and $S$ be a \mo of $\gA$.

\emph{1.} Show that the canonical \ali $M\to M\sur{\fa M}$ solves the \uvl \pb of the \eds for the \homo $\gA\to\gA\sur\fa$ (i.e. according to \Dfn~\ref{defAliAliExtScal}, this \ali is a morphism of \eds from $\gA$ to $\gA\sur\fa$ for $M$). 
Deduce that the natural \ali $\gA\sur\fa\otimes _\gA M\to M\sur{\fa M}$ is an \isoz.

\emph{2.} Show that the canonical \ali $M\to M_S$ solves the \uvl \pb of the \eds for the \homo $\gA\to\gA_S$. Deduce that the natural \ali $\gA_S\otimes _\gA M\to M_S$ is an \isoz.
 }
\end{exercise}

\vspace{-1em}
\begin{exercise}
\label{exoBézoutstrict}
{\rm
Prove that every matrix over a  Bézout domain is \eqve to a matrix of the form $\cmatrix{T&0\cr0&0}$, where $T$ is triangular and the \elts on the diagonal of $T$ are nonzero (naturally, the rows or columns indicated as zero can be absent).
This \eqvc can be obtained by Bézout manipulations.\\
Generalize to \qiris by using the \gnl method explained on \paref{MethodeQI}. }
\end{exercise}

\vspace{-1em}
\begin{exercise}\label{exoAnneauBézoutStrict}
{(Strict Bézout \risz)}\\
{\rm 
\emph {1.}
For a \ri $\gA$, show that \propeq
\vspace{-2pt}
\begin{itemize}\itemsep=0pt
\item [\emph {a.}]
If $A \in \Ae{n \times m}$, there exists a $Q \in \GL_m(\gA)$
such that $AQ$ is a lower triangular matrix.

\item [\emph {b.}]
Same as item \emph {a} with $(n,m) = (1,2)$, i.e. $\gA$ is a strict Bézout \riz.

\item [\emph {c.}]
For $a$, $b \in \gA$, there exist \com $x$, $y \in \gA$  such that $ax + by = 0$.

\item [\emph {d.}]
For $(\ua) = (a_1, \ldots, a_n)$ in $\gA$, there exist a $d\in \gA$ and a \vmd $(\underline {a'}) = (a'_1, \ldots, a'_n)$  satisfying $(\ua) = d(\underline {a'})$; we then have $\gen {\ua} = \gen {d}$.

\item [\emph {e.}]
Same as item \emph {d} with $n=2$.
\end{itemize}

\emph{2.}
Show that the class of strict Bézout \ris is stable under finite products, quotients and  \lonz.

In the following, we assume that $\gA$ is a strict Bézout \riz.

\emph {3.}
Let $a$, $b$, $d_2 \in \gA$ such that $\gen {a,b} = \gen {d_2}$. Show that there exist \com $a_2$, $b_2 \in \gA$ such that $(a,b) = d_2(a_2,b_2)$. We can consider $d_1$, $a_1$, $b_1$, $u_1$, $v_1$ where $(a,b) = d_1(a_1,b_1)$, $1 = u_1a_1 + v_1b_1$ and introduce

\snuc {
(\star) \quad
\cmatrix {a_2\cr b_2} = \crmatrix {v_1 & a_1\cr -u_1 & b_1}
\cmatrix {\vep \cr k_{12}} \; \hbox {where} \;
d_1 = k_{12}d_2,\ d_2 = k_{21}d_1,\ \vep = k_{12}k_{21} - 1. 
}

%
\emph {4.}
Same as in the previous item but with an arbitrary number of \eltsz; \cad for given $(\ua) = (a_1, \ldots, a_n)$ in $\gA$ and $d$ satisfying $\gen {\ua} = \gen {d}$, there exists $(\underline {a'}) = (a'_1, \ldots, a'_n)$, \comz,
such that $\ua = d\underline {a'}$.

\emph {5.}
Show that every diagonal matrix $\Diag(\an)$ is $\SL_n$-\eqve to a diagonal matrix $\Diag(\bn)$ where $b_1 \divi b_2 \divi \cdots \divi b_n$.\\
 Moreover, if we let $\fa_i = \gen {a_i}$, $\fb_i = \gen {b_i}$, we have $\fb_i = S_i(\fa_1, \ldots, \fa_n)$ where $S_i$ is the \gui {$i^{\rm th}$ \elr \smq function of $\fa_1, \ldots, \fa_n$} obtained by replacing each product with an intersection. For example,

\snic {
S_2(\fa_1, \fa_2, \fa_3) = (\fa_1\cap\fa_2) + (\fa_1\cap\fa_3) + (\fa_2\cap\fa_3).
}

In particular, $\fb_1 = \sum_i \fa_i$, $\fb_n = \bigcap_i \fa_i$.
Moreover $\prod_i \gA\sur{\fa_i} \simeq \prod_i \gA\sur{\fb_i}$.
\\
This last result will be \gne to \anars (Corollary~\ref{corthAnar}). 
\\
Other \gui{true} \elr \smq functions of \ids intervene in Exercise~\ref{exoFitt0}.}
\end{exercise}

\vspace{-1em}
\begin{exercise}\label{exoSmith}\index{Smith!ring}\index{ring!Smith ---}
{(Smith \risz, or elementary divisor rings)}\\
{\rm Define a \emph{Smith ring}
as a \ri over which every matrix admits a reduced Smith form
(cf. Section \ref{secBézout}, \paref{secpfval}).
Such a \ri is a strict Bézout \ri
(cf. Exercise~\ref{exoAnneauBézoutStrict}).
Since over a strict Bézout \riz, every square diagonal matrix is \eqv to a Smith matrix
(Exercise~\ref{exoAnneauBézoutStrict}, question~\emph {5}),
a \ri is a Smith \ri \ssi every matrix is \eqv to a \gui{diagonal} matrix, without the condition of \dve over the \coesz.
These \ris have been  studied in particular by Kaplansky in \cite{Kap}, including the noncommutative case, then by Gillman \& Henriksen in \cite{GillmanHenriksen}.
Here we will limit ourselves to the commutative case.
\\
Show that \propeq
\begin{itemize}\itemsep1pt
  \item [\emph{1.}] $\gA$ is a Smith \riz.
  \item [\emph{2.}] $\gA$ is a strict Bézout \ri and every triangular matrix in $\MM_2(\gA)$
is \eqve to a diagonal matrix.
  \item [\emph{3.}] $\gA$ is a strict Bézout \riz, and if $1 \in \gen{a,b,c}$, 
  then there exist  $(p,q)$, $(p',q')$ such that $1=pp'a+qp'b+qq'c$.
  \item [\emph{4.}] $\gA$ is a strict Bézout \riz, and if $\gen{a,b,c} = \gen{g}$, 
  then there exist  $(p,q)$, $(p',q')$ such that $g=pp'a+qp'b+qq'c$.
\end{itemize}
This gives a nice structure \tho for \mpfsz, by taking into account  the uniqueness of \Thref{prop unicyc}.
Also note that this \tho implies the uniqueness of the Smith reduced matrix of a matrix $A$ (by considering the cokernel module) in the following sense. By denoting by $b_i$ the diagonal \coes of the reduced matrix, the \idps $\gen{b_1}\supseteq\cdots\supseteq\gen{b_q}$  ($q=\inf(m,n)$) are invariants of the {matrix} $A$ up to \eqvcz. 
\\
In terms of modules, these \idps characterize, up to \auto of~$\Ae{m}$, the  inclusion morphism $P=\Im(A)\to \Ae{m}$. 
\\
A basis~\hbox{$(e_1,\dots,e_m)$} of $\Ae{m}$ such that $P=b_1\gA\, e_1+\cdots + b_m\gA\, e_m$ is called a \emph{basis of~$\Ae{m}$ adapted to the submodule $P$.}%
\index{basis adapted!to an inclusion} 
\\
Let $b_r=0$ if $m\geq r>n$, then we have $\gen{b_1}\supseteq\cdots\supseteq\gen{b_r}$.
The \idpsz~\hbox{$\neq \gen{1}$} 
of this list are the {invariant factors} of the module $M=\Coker(A)$.
\Thref{prop unicyc} tells us that this list characterizes the structure of the module~$M$.   
\\
Finally, note that the Smith \ris are stable under finite products, \lon and passage to the quotient.
}

\end{exercise}

\vspace{-1em}
\begin{exercise}\label{exoUnitsOfSomeRings}
{(\Elr example of determination of the group of units)}

{\rm 
\emph {1.}
Let $\gk$ be a reduced \ri and $\gA = \aqo{\gk[Y,Z]}{YZ} =\gk[y,z]$ with $yz = 0$. Show, by using a \noepz ing \inoe of $\gA$ over $\gk$, that $\Ati = \gk\eti$.

\emph {2.}
Let $\gA = \aqo{\ZZ[a,b,X,Y]}{X-aY,Y-bX} = \ZZ[\alpha,\beta,x,y]$
with $x = \alpha y$ and $y = \beta x$.
Show that $\Ati = \{\pm 1\}$;
we therefore have $\gA x = \gA y$ but $y \notin \Ati x$.}
\end{exercise}

\vspace{-1em}
\begin{exercise}\label{exoUAtoUB}  
{(Sufficient conditions for the surjectivity of $\Ati \to (\gA/\fa)\eti$)}
\\
{\rm 
Also see Exercise~\ref{exoLgb2}.\\  
For an \id $\fa$ of a \ri $\gA$, we consider the \prt $(\star)$

\snic {
(\star)\qquad\qquad\qquad
\Ati \to (\gA/\fa)\eti \hbox { is surjective}
,}

i.e. for $x \in \gA$ \iv modulo $\fa$, there exists a
$y \in \Ati$ such that $y \equiv x \bmod \fa$, or if 
$\gA x + \fa$ meets $\Ati$, then $x + \fa$ meets $\Ati$.

\emph{1.}
 Show that $(\star)$ is satisfied when $\gA$ is \zedz.

\emph{2.}
 If $(\star)$ is satisfied for all principal \ids $\fa$, then it also is for all \idsz~$\fa$.

\emph{3.}
 Assume $(\star)$ is satisfied. Let $x$, $y$ be two \elts of an \Amo such  
that $\gA x = \gA y$; show that $y = ux$ for some $u \in \Ati$.\\
NB: Exercise~\ref{exoUnitsOfSomeRings}
provides an example of a \ri $\gA$ with $x$, $y \in \gA$
and $\gA x = \gA y$,  but $y \notin \Ati x$.

\emph{4.}
 Let $\gA' = \gA\sur{\Rad \gA}$, $\pi : \gA \twoheadrightarrow \gA'$
be the canonical surjection and $\fa' = \pi(\fa)$. Show that if $(\star)$ is satisfied for $(\gA',\fa')$, then it is satisfied for $(\gA, \fa)$.
}

\end{exercise}

\vspace{-1.2em}
\pagebreak	

\begin{exercise}
\label{exoCalculT(M)} (Computation of a torsion submodule) 
\\
 {\rm  
Let $\gA$ be an integral \coh \ri and $M$ be a \pf \Amoz. 
 Then the torsion submodule of $M$ is a \mpfz.
\\
 More \prmtz, if we have a \mpn $E$ for $M$ with an exact sequence 
$$\preskip.3em \postskip.3em 
\Ae n \vvers{E} \Ae\ell \vvers\pi M\to 0
$$
and if $F$ is a matrix such that we have an exact sequence
$$\preskip.2em \postskip.3em 
\Ae m \vvers{F} \Ae\ell \vvers {\tra E\,} \Ae n 
$$
(the existence of the matrix $F$ results from the fact that $\gA$ is \cohz) then the torsion submodule $\rT(M)$ of $M$ is equal to $\pi(\Ker{\tra F})$ and \isoc to $\Ker{\tra F}/\Im E$.

 Also show that the result can be \gnee to the case where $\gA$ is a \coh  \qiriz. }
\end{exercise}

\vspace{-1em}
\begin{exercise}
\label{exoZerRedBez} (Euclid's \algo in the \zedr case)\\
 {\rm Here we give a more uniform version of the \dem of Proposition~\ref{propZedBez} and we generalize it. Consider a \zedr \ri $\gA$.

\emph{1.} 
Let $\gB$ be an arbitrary \ri and $b\in\gB$ such that $\gen{b}$ is generated by an \idmz.  For $a\in\gB$, find a matrix $M\in\EE_2(\gB)$ and $d\in \gB$ satisfying the \egt $M\cmatrix{a\cr b}=\cmatrix{d\cr 0}$. In particular $\gen{a,b}=\gen{d}$.

\emph{2.}  Give a \gui{uniform} Euclidean \algo for two \pols of $\gA[X]$.

\emph{3.}  
The \ri $\AX$ is a Smith \riz: give an \algo which reduces every matrix over $\AX$ to a Smith form by means of \elr manipulations of rows and of columns.  
}

\end{exercise}

\vspace{-1em}
\begin{exercise}
\label{exoZDpiv} (Syzygies in dimension $0$)\\
{\rm Here we give the \gnn of the \tho according to which $n+1$ vectors of $\gK^n$ are \lint dependent, from the \cdis case to that of the \zedr \risz.
Note that the syzygy, to be worthy of the name, must have \com \coesz. 

Let $\gK$ be a \zedr \riz, and $y_1$, \ldots, $y_{n+1}\in\gK^{n}$.

\emph{1}. Construct a \sfio $(e_j)_{j\in \lrb{1..n+1}}$ such that, in each component $\gK[1/e_j]$, the vector $y_j$ is a \coli of the $y_i$'s that precede it.

\emph{2}. Deduce that there exists a \sys of \com \elts $(a_1, \ldots, a_{n+1})$ in~$\gK$ 
such that $\sum_i a_iy_i = 0$.

\sni
\rems 1) Recall the convention according to which we accept that certain \elts of a \sfio are null. We see in this example that the statement of the desired \prt is greatly facilitated by it.

 2)
We can either give an adequate working of the matrix of the $y_i$'s by \elrs manipulations by basing ourselves on Lemma~\ref{lemQI}, or treat the \cdis case then use the \elgbm \num2 (\paref{MethodeZedRed}).\eoe
}
\end{exercise}

\vspace{-1.2em}
\pagebreak	

\begin{exercise}
\label{exoZedLG}
{\rm Let $S_1$, $\dots$, $S_n$ be \moco of $\gA$. Show that $\gA$ is \zed \ssi
each of the $\gA_{S_i}$'s is \zedz.
}
\end{exercise}

\vspace{-1em}
\begin{exercise}\label{exoFreeAlgebraPresentation}
(Presentation of an \alg which is free and finite as a module)
\\
{\rm
Let~$\gB$ be a free \Alg of rank $n$ with basis $\ue = (e_1, \ldots, e_n)$. We let

\snic{\varphi : \AuX = \gA[X_1, \ldots, X_n]\twoheadrightarrow \gB}

be the (surjective) \homo of \Algs which performs $X_i \mapsto e_i$. Let $c_{ij}^k$ be the structure constants defined by $e_ie_j = \sum_k c_{ij}^k e_k$.
Consider $a_1$, \ldots, $a_n \in \gA$ defined by $1 = \sum_k a_k e_k$ and let

\snic{
R_0 = 1 - \som_k a_k X_k, \qquad
R_{ij} = X_iX_j - \sum c_{ij}^k X_k.}

Let $\fa = \gen {R_0, R_{ij}, i \leq j}$.
Show that every $f \in \AuX$ is congruent modulo $\fa$ to a \hmg \pol of degree $1$. Deduce that $\Ker\varphi = \fa$.}
\end{exercise}

\vspace{-1em}
\begin{exercise}\label{exoFitt0}
(Some computations of \idfsz)\\
{\rm
\emph{1.} 
Determine the \idfs of an \Amo presented by a matrix in Smith form.

\emph{2.}   
Determine the \idfs of $\gA\sur\fa$.

\emph{3.}   
Let $E$ be a \tf \Amo and $\fa$ be an \idz. Show that
$$\preskip.2em \postskip.4em 
\cF_k(E\oplus \gA\sur{\fa})= \cF_{k-1}(E)+\cF_{k}(E)\,\fa. 
$$
\emph{4.}   
Determine the \idfs of the \Amo $M=\gA\sur{\fa_1}\oplus\cdots\oplus \gA\sur{\fa_n}$ in the case \hbox{where $\fa_1\subseteq\fa_2\subseteq \cdots\subseteq\fa_n$}. 

\emph{5.} 
 Determine the \idfs of the \Amo $M=\gA\sur{\fa_1}\oplus\cdots\oplus$ $\gA\sur{\fa_n}$
without making any inclusion assumptions for the \ids $\fa_k$.
\\
Compare $\cF_0(M)$ and $\Ann(M)$.}
\end{exercise}

\vspace{-1em}
\begin{exercise}\label{exoFitt4}
(The \idfs of a \tf \Amoz)\\
{\rm
Show that Facts~\ref{fact.idf inc}, \ref{fact.idf.change}, \ref{fact.idf.sex} and Lemma~\ref{fact.idf.ann} remain valid for \mtfsz.}
\end{exercise}

\vspace{-1em}
\begin{exercise}
\label{exoFitt6}
{\rm  One of the \cara \prts of  \emph{\adpsz} (which will be studied in Chapter~\ref{ChapAdpc}) is the following: if $A\in\Ae {n\times m}$, $B\in\Ae {n\times 1}$, and if the \idds of $A$ and $[\,A\,|\,B\,]$ are the same, then the \sliz~$AX=B$ admits a solution.

\emph{1.} Let $M$ be a \mtf over a \adp and $N$ be a quotient of~$M$. Show that
if~$M$ and $N$ have the same \idfsz, then $M=N$.

\emph{2.} Show that if a \mtf $M$ over a \adp has \tf \idfsz, then it is a \mpfz.}
\end{exercise}

\vspace{-1em}
\begin{exercise}\label{exoAutresIdF}\index{ideal!Kaplansky ---}\index{Kaplansky!ideal}
(Kaplansky \idsz) \\
{\rm For an \Amo $M$ and an integer $r$ we denote by $\cK_r(M)$ the \id which is a sum of all the conductors $\big(\!\gen{m_1,\ldots m_r}:M\big)$ for all the \syss $(m_1,\ldots m_r)$ in $M$.
We call it \emph{the Kaplansky \id of order $r$ of the module $M$}.
Thus,  $\cK_0(M)=\Ann (M)$, and if $M$ is generated by $q$ \eltsz,
 we have $\cK_q(M)=\gen{1}$.
\begin{itemize}\itemsep0pt
\item Show that if $\cK_q(M)=\gen{1}$, $M$ is \tfz.
\item Show that if $M$ is \tfz, then for every integer $r$ we have the inclusions
$$\preskip-.2em \postskip.0em 
\cF_r(M) \subseteq  \cK_r(M) \subseteq \sqrt{\cF_r(M)}= \sqrt{\cK_r(M)}. 
$$
\end{itemize}
NB: see also 
Exercise~\ref{exoVariationLocGenerated}.}
\end{exercise}

\vspace{-1em}
\begin{exercise}\label{exoPetitsExemplesElim} 
{(An \elr example of resultant \idsz)} \\
{\rm  
Let $f$, $g_1$, \ldots, $g_r \in \AX$, $f$ be \mon of degree $d \ge 1$ and $\ff = \gen {f, g_1, \ldots, g_r} \subseteq \AX$.  We will compare the \id 
$$\preskip.2em \postskip.4em 
\fa = \fR(f, g_1,\ldots,g_r) = \rc_T\big(\Res(f, g_1 + g_2T + \cdots + g_rT^{r-1})\big)
$$
(Section~\ref{secChap3Nst}), and the resultant \id 
$\fb = \fRes(\ff) = {\cF_{\gA,0}(\AX\sur\ff)}$ (see the \gnl \eli lemma of Section~\ref{subsecIdealResultant}).

\emph {1.}
Let $\fa' = \rc_\uT\big(\Res(f, g_1T_1 + g_2T_2 + \cdots +
g_rT_r)\big)$. Show the inclusions
$$\preskip.3em \postskip.3em 
\fa \subseteq \fa' \subseteq \fb \subseteq \ff \cap \gA
. 
$$
\emph {2.}
Let $\gA = \ZZ[a,b,c]$ where $a$, $b$, $c$ are three \idtrsz,
$f = X^d$, $g_1=a$, $g_2=b$ and $g_3=c$. Determine the \ids $\ff \cap \gA$, $\fa$, $\fa'$, $\fb$ and check that they are distinct.  
\Egmt check that $\fR(f, g_1, g_2, g_3)$ depends on the order of the $g_i$'s. \\
Do we have $(\ff\cap\gA)^d \subseteq \fa$?}
\end{exercise}

\vspace{-1em}
\begin{exercise}\label{exoRelateursViaElimIdeal} 
{(Relators and \eli \idz)}\index{elimination!ideal}\index{ideal!elimination ---}
\\
 {\rm  
Let $f_1(\uX), \ldots, f_s(\uX) \in \kuX = \gk[\Xn]$ ($\gk$ is a commutative \riz). 
\\
Let $\fa \subseteq \kuY = \gk[Y_1, \ldots, Y_s]$ be the \id of the relators over $\gk$ of $(f_1, \ldots, f_s)$, \hbox{\cad $\fa = \ker\varphi$}, where $\varphi : \kuY \to \kuX$ is the \evn morphism $Y_i \mapsto f_i$. 
\\
Let $g_i = f_i(\uX) - Y_i \in
\gk[\uY, \uX]$ and $\ff = \gen {g_1, \ldots, g_s}$. \\
Prove that $\fa = \ff \cap \kuY$.
Thus, $\fa$ is the \eli \id of the variables $X_j$ in the \syp of the $g_i$'s.
}
\end{exercise}


\vspace{-1em}
\begin{problem}\label{exoDimZeroXcYbZa}
 {(An example of a \zed \sysz)}\\
{\rm  
Let $\gk$ be a \ri and $a$, $b$, $c \in \NN^*$ with $a \le b \le c$ and at least one strict in\egtz. We define three \pols $f_i \in \gk[X,Y,Z]$

\snic {
f_1 = X^c+Y^b+Z^a,\quad  f_2 = X^a+Y^c+Z^b,\quad  f_3 = X^b+Y^a+Z^c
.}


This is a matter of studying the \sys defined by these three \polsz. 
We denote by $\gA = \gk[x,y,z]$ \hbox{the \klg} $\aqo{\gk[X,Y,Z]}{f_1,f_2,f_3}$.

\emph {1.}
For an arbitrary \ri $\gk$, is $\gA$ free and finite over $\gk$? If so, compute a basis and give the dimension.

\emph {2.}
Give a detailed study of the \sys for $\gk = \QQ$ and $(a,b,c) = (2,2,3)$. That is, determine all the zeros of the \sys in a certain finite extension of $\QQ$ (to be specified), their number and their multiplicities.

\emph {3.}
Is the localized \alg $\gA_{1+\gen{x,y,z}}$ free over $\gk$?  If so, give a basis.
}

\end{problem}

\vspace{-1em}
\begin{problem}
\label{exoIdealResultantGenerique} (The generic resultant \idz)\\
{\rm 
Let $d$, $r$ be two fixed integers with $d\ge 1$. In this exercise we study the generic resultant \id $\fb = \fRes(f, g_1, \ldots, g_r)$ where $f$ is \mon of degree $d$, \hbox{and $g_1$, \ldots, $g_r$} are of degree $d-1$, the \coes of these \pols being \idtrs over~$\ZZ$. The base \ri is therefore $\gk =\ZZ[(a_i)_{i\in \lrb{1..d}}, (b_{ji})_{j\in \lrb{1..r}, i\in
\lrb{1..d}}]$ with

\snic{
f = X^d + \sum_{i=1}^{d} a_iX^{d-i} \quad\hbox{and}\quad  g_j = \sum_{i=1}^{d} b_{ji}X^{d-i}.
}

\emph{1.}  
Put weights on the $a_i$'s and $b_{ij}$'s such that $\fb$ is a \hmg \idz.

\emph{2.}
If $S$ is the \gnee Sylvester matrix of $(f, g_1, \ldots, g_r)$,
specify the weight of the \coes of $S$ and those of its minors of order $d$.

\emph{3.}
Using a Computer Algebra \sysz, study the minimal number of \gtrs of~$\fb$. We could replace $\ZZ$ with $\QQ$, 
introduce the \id $\fm$ of $\gk$ generated by all the \idtrs and consider $E = \fb\sur{\fm\fb}$ which is a finite dimensional \evc over $\gk\sur{\fm} = \QQ$.

}
\end{problem}

\vspace{-1em}
\begin{problem}\label{exoNakayamaHomogeneRegularSequence}
{(Homogeneous Nakayama lemma and \ndz sequences)}
\\
{\rm  
\noindent\emph{1.} \emph{(Regular sequence and \agq independence)}
Let $(\an)$ be a \seqreg of a \ri $\gA$ and $\gk \subseteq \gA$ be a sub\ri such that $\gk \cap \gen {\an} = \{0\}$. Show that $a_1$, \dots, $a_n$ are \agqt independent over~$\gk$.

\emph{2.} {\it(Homogeneous Nakayama lemma)}
Let $\gA = \gA_0 \oplus \gA_1 \oplus \gA_2 \oplus \dots$ be a graded \ri \hbox{and $E = E_0 \oplus E_1 \oplus E_2 \oplus \dots$} be a graded \Amoz. \\
We denote by $\gA_+$ the \id $\gA_1 \oplus \gA_2 \oplus \dots$, so that $\gA/\gA_+ \simeq \gA_0$.
\vspace{-2pt}
\begin{itemize}\itemsep=0pt
\item [\emph {a.}]
Show that if $\gA_+ E = E$, then $E = 0$.

\item [\emph {b.}]
Let $(e_i)_{i \in I}$ be a family of \hmg \elts of $E$. 
Show that if the $e_i$'s generate the $\gA_0$-module $E/\gA_+E$, then they generate the \Amoz~$E$.
\end{itemize}
Note that we do not assume that $E$ is \tfz.

\emph{3.}
Let $\gB = \gB_0 \oplus \gB_1 \oplus \gB_2 \oplus \dots$ be a graded \ri and $h_1$, \ldots, $h_d$ be \hmg \elts of the \id $\gB_+$. Let $\fb = \gen {h_1, \ldots, h_d}$ and $\gA = \gB_0[h_1, \ldots, h_d]$. We therefore have $\gB_0 \cap \fb = \{0\}$, and $\gA$ is a graded sub\ri of $\gB$. Finally, let $(e_i)_{i \in I}$ be a family of \hmg \elts of $\gB$ that generate the $\gB_0$-module $\gB/\fb$.
%
\begin{itemize}
\item [\emph {a.}]
Verify that $\gA_0 = \gB_0$ and that $\fb = \gA_+ \gB$ then show that the $e_i$'s form a \sgr of the \Amo $\gB$.
\item [\emph {b.}]
Suppose that $(h_1, \ldots, h_d)$ is a \seqreg and that the $e_i$'s form a basis of the $\gB_0$-module $\gB/\fb$. Show that $ h_1$, \ldots, $h_d$ are \agqt independent over~$\gB_0$ and that the $e_i$'s form a basis of the \Amo $\gB$.
\end{itemize}
\emph{Recap}:
Let $\gB = \gB_0 \oplus \gB_1 \oplus \gB_2 \oplus \dots$ be a graded \ri and $(h_1, \ldots, h_d)$ be a \hmg \seqreg of the \id $\gB_+$. If $\aqo{\gB}{h_1, \ldots,h_d}$ is a free $\gB_0$-module, then $\gB$ is a free $\gB_0[h_1, \ldots, h_d]$-module and $\gB_0[h_1, \ldots, h_d]$ is a \ri of \pols in $(h_1, \ldots, h_d)$.

\emph {4.}
As a converse.
Let $\gB = \gB_0 \oplus \gB_1 \oplus \gB_2 \oplus \dots$ be a graded \ri \hbox{and $h_1$, \ldots, $h_d$} be \hmg \elts of the \id $\gB_+$, \agqt independent over $\gB_0$. If~$\gB$ is a free $\gB_0[h_1, \ldots, h_d]$-module,
then the sequence $(h_1, \ldots, h_d)$ is \ndzez.

}
\end{problem}


\pagebreak	
\sol

\exer{exoptfpf} {It suffices to apply Proposition~\ref{propPfSex}. Directly:
we consider a \prr $\pi:N\to N$ having $M$ as its image. 
If $X$ is a \sgr of~$N$, then~$\pi(X)$ is a \sgr of $M$. If $N$ is \pfz, the syzygy module  for $\pi(X)$ is obtained by taking the syzygies for $X$ in~$N$ and the syzygies
$\pi(x)=x$ for each \elt $x$ of~$X$.
}


\exer{exoAXmodule} We start by noting that $\theta_A \circ (X\In - A) = 0$
and that $\theta_A$ is the identity over $\Ae n$.

\emph {1.} 
Let us show that $\Im(X\In-A) \cap \Ae n = 0$.  Let $x \in \Im(X\In-A) \cap \Ae n$,
the preliminary computations give $\theta_A(x)=x$ and $\theta_A(x) = 0$.
Let us show that $\AX^n = \Im(X\In-A) + \Ae n$.
It suffices to show that $X^k e_i \in \Im(X\In-A) + \Ae n$ for $k\geq 0$ and $i\in\lrbn$. If $k=0$ it is clear. For $k>0$ we write

\snic {
X^k\In-A^k = (X\In-A) \, \sum_{j+\ell = k-1} X^j A^\ell.
}

By applying this \egt to $e_i$, we obtain $X^k e_i - A^k e_i \in \Im(X\In-A)$,
{so $X^k e_i$ belongs to $\Im(X\In-A) + A^k e_i \subseteq \Im(X\In-A) + \Ae n$}.

\emph {2.}
Let $y \in \Ker\theta_A$.  Let $y = z + w$ with $z \in \Im(X\In-A)$ and $w \in \Ae n$. \\
Therefore $0 =\theta_A(y) = \theta_A(z) + \theta_A(w) = 0 + w$ and $y = z \in \Im(X\In-A)$.


\exer{exoTenseursNuls}
 {(Description of the null tensors, \gnl situation)}\\
\emph{1.} This results from the \dfn of the tensor product and from the fact that in algebra, computations are finite.

\emph{2.} Let $X=[\,x_1\,\cdots\,x_p\,]\in M_1^{1\times p}$, $Y=\tra[\,y_1\,\cdots\,y_p\,]\in N^{p\times 1}$.\\  
We have $M_1=\gA x_1+\cdots+\gA x_p$ and $\sum_{i\in\lrbp}x_i\te y_i=_{M_1\te N}0$, and by the null tensor lemma,  this \egt holds \ssi there exist $q\in\NN$, $G\in \gA^{p\times q}$ and $Z=\tra[\,z_1\,\cdots\,z_q\,]\in N^{q\times 1}$ that satisfy
$$\preskip.2em \postskip.0em 
XG=_{M^q} 0\quad {\rm and } \quad GZ=_{N^{p}}Y\,.
$$


\exer{exoAnneauBézoutStrict}
\emph{1} and \emph{2} are left to the reader.

\emph{3.}
By construction, $\vep$ annihilates $d_2$ (\cad annihilates $a,b$). We therefore have the \egts

\snuc{
d_2\cmatrix {a_2\cr b_2} = 
\crmatrix {v_1 & a_1\cr -u_1 & b_1} \cmatrix {d_2\vep \cr d_2k_{12}} =
\crmatrix {v_1 & a_1\cr -u_1 & b_1} \cmatrix {0 \cr d_1} =
d_1\cmatrix {a_1\cr b_1} = \cmatrix {a\cr b} 
}

It remains to see that $1 \in \gen {a_2, b_2}$. By inverting the $2\times 2$ matrix in $(\star)$ (of \deterz~1), we see that the \id $\gen {a_2, b_2}$ contains $\vep$ and $k_{12}$, so it contains $1 = k_{12}k_{21} - \vep$.

\emph{4.}
By \recu on $n$, $n=2$ being the previous question.
Suppose $n \ge 3$. By \recuz, there exist  $d$ and \com $b_1$, \ldots, $b_{n-1}$ such that 

\snic{(a_1, \ldots, a_{n-1}) = d(b_1, \ldots, b_{n-1}),\hbox{ so }\gen{\ua} = \gen {d, a_n}.}

Item \emph {3}
gives  \com $u$,~$v$ and~$\delta$ such that $(d,a_n) = \delta(u,v)$.
\\
Then
$(\an) = (db_1, \ldots, db_{n-1}, \delta v) =
\delta (ub_1, \ldots, ub_{n-1}, v),$ 
and $\gen {1}=\gen {u,v} =\gen {ub_1, \ldots, ub_{n-1}, v}$.

\emph{5.}
First for $n=2$ with $(a,b)$. There is a $d$ with $(a,b) = d(a',b')$ and $1 = ua' + vb'$. Let $m = da'b' = ab' = ba' \in \gen {a}\cap\gen{b}$; we have $\gen {a} \cap \gen{b} = \gen{m}$ because if $x \in \gen {a} \cap \gen{b}$,
then $x=x(ua' + vb') \in \gen {ba'} + \gen {ab'} = \gen{m}$.
The $\SL_2(\gA)$-\eqvc is provided by the \egt below

\snic {
\cmatrix {1 &-1 \cr vb' & ua'} \cmatrix {a &0 \cr 0 & b} =
\cmatrix {d & 0 \cr 0 & m} \cmatrix {a' &-b' \cr v & u}
.}

For $n \ge 3$. By using the $n=2$ case for the positions $(1,2)$, $(1,3)$, \dots, $(1,n)$, we obtain $\Diag(a_1, a_2, \ldots, a_n) \sim \Diag(a'_1, a'_2, \ldots, a'_n)$ with $a'_1 \divi a'_i$ for $i \ge 2$.
\\
By \recuz, $\Diag(a'_2, \ldots, a'_n) \sim \Diag(b_2, \ldots, b_n)$
where $b_2 \divi b_3 \cdots \divi b_n$. We then check that $a'_1 \divi b_2$ and we let $b_1 = a'_1$.
The scrupulous reader will check the \prt regarding \elr \smq functions.


\exer{exoSmith} \emph{(Smith \risz, or elementary divisor rings)}
\\
Preliminary computation with $A = \cmatrix {a& b\cr 0 &c}$
and $B$ of the form

\snic {
B = \cmatrix {p' & q'\cr *  &* } A \cmatrix {p & * \cr q &* }
.}

The \coe $b_{11}$ of $B$ is equal to $b_{11} = p'(pa + qb) + q'qc$.

$\emph {2} \Rightarrow \emph {3.}$ 
The matrix $A$ is \eqve to a diagonal matrix $\Diag(g,h)$, which gives $(p,q)$ and $(p',q')$ \com with $g = p'(pa + qb) + q'qc$ (preliminary computation), and we have $\gen {a,b,c} = \gen{g,h}$. As $\gA$ is a strict Bézout \riz, we can suppose that $g \divi h$ and since $1 \in \gen {a,b,c}$, $g$ is \iv and so $1$ is expressed as required.

$\emph {3} \Rightarrow \emph {4.}$ 
Beware, here $g$ is imposed. But by question \emph{4} of Exercise~\ref{exoAnneauBézoutStrict},
we can write $(a,b,c) = g\,(a',b',c')$ with $(a',b',c')$ \comz. We apply item~\emph {3} to $(a',b',c')$ and multiply the obtained result by $g$.

$\emph {4} \Rightarrow \emph {2.}$ 
Let $A \in \MM_2(\gA)$ be triangular, $A = \cmatrix {a& b\cr 0 &c}$.  With the parameters of item~\emph {4}, we construct (preliminary computation) a matrix $B$ \eqve to $A$ with \coe $b_{11} = g$. As $g$ divides all the \coes of $B$, we have $B \,\sims{\EE_2(\gA)}\, \Diag(g,h)$.

\emph {1} $\Leftrightarrow$ \emph {2.} Left to the reader (who can consult Kaplansky's paper).

\exer{exoUnitsOfSomeRings} 
\emph {1.}
Let $s = y+z$. Then $\gk[s]$ is a \pol \ri in $s$, and~$y$,~$z$ 
are integral over $\gk[s]$, because they are zeros of $(T-y)(T-z) = T(T-s) \in \gk[s][T]$. 
We easily check that $\gA$ is free over $\gk[s]$ with $(1,y)$ as its basis.  For $u$, $v \in \gk[s]$, the norm over $\gk[s]$ of $u + vy$ is

\snic {
\rN_{\gA/\gk[s]}(u+vy) = (u+vy)(u+vz) = u^2 + suv = u(u + sv)
.}

The \elt $u + vy$ is \iv in $\gA$ \ssi $u(u + sv)$ is \iv in $\gk[s]$. As $\gk$ is reduced, $(\gk[s])\eti = \gk\eti$. Therefore $u \in \gk\eti$ and $v = 0$.

\emph {2.}
We have $\gA = \ZZ[\alpha,\beta,y] = \aqo{\ZZ[a,b,Y]}{(ab-1)Y}$ with
$y(\alpha\beta-1) = 0$. Let $t$ be an \idtr over $\ZZ$ and $\gk = \ZZ[t,t^{-1}]$. Consider the \klg $\gk[y,z]$ with the sole syzygy $yz= 0$.
We have a morphism $\gA \to \gk[y,z]$ which performs

\snic{\alpha \mapsto t(z+1)$, $\beta \mapsto t^{-1}$, $y \mapsto y,}

and we check that it is an injection.\\
 Then an \elt $w \in \Ati$ is also in~$\gk[y,z]^{\times}$, and as~$\gk$ is reduced, $w \in \gk\eti$.
Finally, the units of $\gk = \ZZ[t,t^{-1}]$ are the $\pm t^k$ with $k \in \ZZ$, so $w = \pm 1$.


\exer{exoUAtoUB} \emph {1.}
We know that $\gen {x^n} = \gen {e}$. We look for $y \in \Ati$ such that $y \equiv x \bmod \fa$ over the components $\gA_e$ and $\gA_{1-e}$. First, we have $x^n (1-ax) = 0$ and $x$ \iv modulo $\fa$, so $ax \equiv 1 \bmod \fa$ then $e \equiv 1 \bmod\fa$, \cad $1-e \in \fa$.  \\
In the component $\gA_e$, $x$ is \ivz, so we can take $y=x$. In the component $\gA_{1-e}$, $1\in \fa$, so we can take $y=1$.  Globally, we therefore propose that $y = ex + 1-e$ which is indeed \iv (with inverse $ea^{n}x^{n-1} + 1-e$) and which satisfies $y \equiv x \bmod
\fa$.  Remark: $y = ex + (1-e)u$ with $u \in \Ati$ is also suitable.

\emph {2.}
Let $x$ be \iv modulo $\fa$ so $1-ax \in \fa$ for some $a \in \gA$. \\
Then, $x$ is \iv modulo the principal \id $\gen {1-ax}$,
therefore there exists a $y \in \Ati$ such that $y \equiv x \bmod \gen {1-ax}$,
a fortiori $y \equiv x \bmod \fa$.

\emph {3.}
We write $y=bx$, $x=ay$ so $(1-ab)x = 0$; $b$ is \iv modulo 
$\gen {1-ab}$ so there exists a $u \in \Ati$ such that $u \equiv b \bmod
\gen{1-ab}$ whence $ux = bx = y$.

\emph {4.}
Let $x$ be \iv modulo $\fa$. Then $\pi(x)$ is \iv modulo $\fa'$, whence $y \in \gA$ such that~$\pi(y)$ is \iv in $\gA'$ and $\pi(y) \equiv \pi(x) \bmod \fa'$. Then, $y$ is \iv in $\gA$ and $y-x \in \fa + \Rad\gA$, \cad $y = x + a + z$ with $a \in \fa$ and $z \in \Rad\gA$.
Thus, the \elt $y-z$ is \iv in $\gA$, and $y-z \equiv x \bmod \fa$.

\exer{exoCalculT(M)}  
Let us call $\gA_1$ the quotient field of $\gA$ and let us put an index $1$ to indicate that we are performing a \eds from $\gA$ to $\gA_1$. Thus $M_1$ is the $\gA_1$-\evc corresponding to the exact sequence

\snic{\gA^n_1 \vers{E_1} \gA^\ell_1 \vers{\pi_1} M_1\to 0}

and the submodule $\rT(M)$ of $M$ is the kernel of the natural \Ali from $M$ to $M_1$,
\cad the module $\pi(\Ae\ell \cap \Ker \pi_1)$, or the module $\pi(\Ae\ell \cap \Im E_1)$
 (by regarding $\Ae\ell$ as a submodule of $\gA^\ell_1$).
\\
 The exact sequence
$\Ae m \vers{F} \Ae\ell \vers {\tra E} \Ae n$
gives by \lon the exact sequence
$$\preskip.1em \postskip.4em 
\gA^ m_1 \vers{F_1} \gA^\ell_1 \vers {\tra {E_1}} \gA^n_1 
$$ 
and since $\gA_1$ is a \cdi this gives by duality the exact sequence
$$\preskip.1em \postskip.4em 
\gA^n_1 \vers{E_1} \gA^\ell_1 \vers {\tra {F_1}} \gA^ m_1. 
$$
Thus $\Im E_1=\Ker \!{\tra {F_1}}$, so $\Ae\ell \cap \Im E_1=\Ae\ell \cap \Ker\! {\tra {F_1}}$. Finally, we have the \egtz~\hbox{$\Ae\ell \cap \Ker\! {\tra {F_1}}=\Ker\! {\tra F}$}
because the natural morphism $\gA\to \gA_1$ is injective.
\\
Conclusion: $\rT(M)$ is equal to $\pi(\Ker \!{\tra F})$, \isoc to
$\Ker\!{\tra F}/\Im E$, and therefore is \pf (because $\gA$ is \cohz).

If $\gA$ is a  \coh \qiriz, the total \ri of fractions $\gA_1=\Frac\gA$ is \zedrz, and all the arguments given in the integral case work similarly.


\exer{exoZerRedBez} 
All the results can be obtained from the \cdi case,   for which the \algos are classical, by using the \elgbm of  \zed \risz. Here we will clarify this  very \gnl affirmation.
\\
Let us put two preliminary remarks for an arbitrary \ri $\gA$.
\\
First, let $e$ be \idm and $E$ be an \elr matrix modulo~${1-e}$. If we lift~$E$ to a matrix $F\in \Mn(\gA)$, then the matrix $(1-e)\In+eF\in\En(\gA)$ is \elrz, it acts like $E$ in the component $\aqo\gA {1-e}$, and it does nothing in the component $\aqo\gA {e}$. This allows us to understand how we can retrieve the desired results over $\gA$ by using analogous results modulo the \idms $1-e_i$ when we have a \sfio $(e_1,\ldots,e_k)$ (provided by the \algo that we build).\\ 
Second, if $g \in \AX$ is \mon of degree $m \ge 0$, for all $f \in \AX$, we can divide $f$ by $g$: $f = gq + r$ with $r$ of formal degree $m-1$.

 \emph {1.}
Let $e$ be the \idm such that $\gen {e} = \gen {b}$. It suffices to solve the question modulo $e$ and $1-e$. In the branch $e = 1$, $b$ is \ivz, $\gen{a,b} = \gen {1}$ and the \pb is solved (Gauss pivot). 
In the branch $e = 0$, $b$ is null and the \pb is solved.
If $e = bx$, we find $d = e + (1-e)a$ and

\snic {
M = \rE_{21}(-be) \rE_{12}\big(ex(1-a)\big) =
\cmatrix {1 & ex(1-a)\cr -eb & ae + (1-e)\cr}
.}

\emph {2.} 
We start from two \pols $f$ and $g$. We will build a \pol $h$ and a matrix $M\in\EE_2(\AX)$ such that $M\cmatrix{f\cr g}=\cmatrix{h\cr 0}$. A fortiori $\gen{f,g}=\gen{h}$.
 \\
We proceed by \recu on $m$, the formal degree of $g$, with formally leading \coez~$b$. 
If we initiate the \recu at $m=-1$, $g = 0$ and $\I_2\cmatrix{f\cr g}=\cmatrix{f\cr 0}$, we can treat $m=0$, with $g \in \gA$ and use item \emph {1} ($\gB = \AX$, $a=f$, $b=g$). But it is pointless to treat this case separately (and so we no longer use item~\emph {1}).  Indeed, if~$e$ is the \idm such that $\gen {e} = \gen {b}$, it suffices 
to solve the question modulo~$e$ and $1-e$ and what follows 
holds for all $m \ge 0$.
\\
In the branch $e = 1$, $b$ is \ivz, and since $m \ge 0$, we can perform a classical Euclidean division of $f$ by $g$: $f = qg -r$ with the formal degree of $r$ equal to $m-1$.  We get a matrix
$N\in\EE_2(\AX)$ such that {$N\cmatrix{f\cr g}=\cmatrix{g\cr r}$}, namely \smashtop{$N = \crmatrix {0 &1\cr -1 &q\cr}$}. We can then apply the \hdrz.
\\
In the branch $e = 0$, $g$ is of formal degree $m-1$ and the \hdr applies.

In the following, we use item~\emph{2} by saying that we pass from~$\tra{\vab f  g}$ to~$\tra{\vab  h 0}$ by means of \gui{Bézout manipulations.}

 \emph{3.} By relying on the result of item~\emph{2} we are inspired by the \dem of Proposition~\ref{propPfPID} (a PID is a Smith \riz).  If we were in a nontrivial \cdiz, the \algo would terminate in a finite number of steps which can be directly bounded in terms of $(D,m,n)$, where $D$ is the maximum degree of the \coes of the matrix $M\in \AX ^{m\times n}$ that we want to reduce to the Smith form. It follows that when $\gA$ is \zedr the number of splittings produced by the gcd computations (as in item~\emph{2}) is also bounded in terms of $(D,m,n)$, where~$D$ is now the maximum formal degree of the entries of the matrix. 
 This shows that the complete \algoz, given the preliminary remark, also terminates in a number of steps bounded in terms of $(D,m,n)$.

\rem  The \algos do not require that $\gA$ be discrete.
\eoe

\exer{exoFreeAlgebraPresentation}
It is clear that $\fa\subseteq\Ker\varphi$.
Let $\cE \subseteq \AuX$ be the set of \polsz~$f$ congruent modulo $\fa$ to a \hmg
\pol of degree $1$.\\ We have $1 \in \cE$ and $f \in \cE \Rightarrow X_if \in \cE$ because if $f \equiv \sum_j \alpha_j X_j \mod {\fa}$, then

\snic{
X_if \equiv \som_j \alpha_j X_iX_j \equiv \som_{j,k}
\alpha_j c_{ij}^k X_k \mod {\fa}.}

Therefore $\cE = \AuX$. Let $f \in \Ker\varphi$. We write $f \equiv \sum_k \alpha_k X_k \mod {\fa}$.
\\
Then $\varphi(f)=0 = \sum_k \alpha_k e_k$,
so $\alpha_k = 0$, then $f \in \fa$.

\exer{exoFitt0}~\\
\emph{2.} If $\fa$ is \tf a \mpn of the module $M=\gA\sur\fa$ is a matrix row $L$ having for \coes \gtrs of the \idz. We deduce that $\cD_1(L)=\fa$. Therefore $\cF_{-1}(M)=0\subseteq\cF_0(M)=\fa\subseteq\cF_1(M)=\gen{1}$. 
The result can be \gnee to an arbitrary \id $\fa$.

\emph{3.} Results from \emph{2} and Fact~\ref{fact.idf.sex}.

\emph{4} and \emph{5.} In the \gnl case by applying \emph{2} and \emph{3} we find 

\snic{\cF_0(M)=\prod_{i=1}^n \fa_i$, $\cF_{n-1}(M)=\som_{i=1}^n \fa_i,}

and for the intermediate \ids the \gui{\smq functions} 

\snic{\cF_{n-k}(M)=\som_{1\leq i_1<\ldots< i_k\leq n} \prod_{\ell=1}^k \fa_{i_\ell}.}

In addition, $\Ann (M)=\fa_1\cap\cdots\cap\fa_n$.
\\
When $\fa_1\subseteq\fa_2\subseteq \cdots\subseteq\fa_n$ the result is a little simpler 

\snic{\cF_{n-1}(M)=\fa_n$, $\cF_{n-2}(M)=\fa_n\,\fa_{n-1}$, \ldots\, $\cF_{n-k}(M)=\fa_n\cdots \fa_{n-k+1}.}

We then find for item \emph{1} the result of the direct computation given by the \idds of a matrix in Smith form.

\exer{exoFitt6}{
Let us prove item \emph{1} (afterwards we can apply Fact~\ref{facttfpf}).\\
Take $M=\gen{g_1,\ldots ,g_q}$. Consider a syzygy $\som_i\alpha_ig_i=_N0$.
The aim is to show that the column vector $V=(\alpha_1,\ldots ,\alpha_q)$ is a syzygy in $M$.

\emph{First case, $M$ is \pfz.} \\
Adding the $V$ column to a \mpn $F$ of $M$ for $(g_1,\ldots,g_q)$ does not change the \idds of this matrix, so $V$ is a \coli of the columns of $F$. 

\emph{Second case, $M$ is \tfz.} \\
Since $\cD_1(V)\subseteq \cF_{q-1}(M)$, there exists a matrix $F_1$ of syzygies for $(g_1,\ldots ,g_q)$ in~$M$ with $\cD_1(V)\subseteq \cD_{1}(F_1)$. Since  $\cD_2(V|F_1)\subseteq \cF_{q-2}(M)$, there exists a matrix~$F_2$ of syzygies for $(g_1,\ldots ,g_q)$ in $M$ with $\cD_2(V|F_1)\subseteq \cD_{2}(F_1|F_2)$, but also of course $\cD_1(V|F_1)\subseteq \cD_{1}(F_1|F_2)$, and so on until there exists a matrix $F=[\,F_1\mid\cdots \mid F_q\,]$ of syzygies for $(g_1,\ldots ,g_q)$ in $M$ such that the \idds of $[V|F]$ are contained in those of $F$.
Therefore $V$ is a \coli of the columns of~$F$.
}

\exer{exoAutresIdF}
If a Kaplansky \id is equal to $1$, 
then the module is finitely generated, because the module is finitely generated in the localized \ris $\gA[1/a_i]$'s with the $a_i$'s being \comz.
\\
Key idea: the Kaplansky \ids are a little more \gnlz, but apparently useless in the case where the module is not \tfz. The Kaplansky \ids present the advantage over the \idfs of allowing a \carn of the \mtfsz.

For the second item, here is what happens.
\\
If $a$ is a typical \gtr of $\cK_r(M)$ and if $M$ is generated by $(g_1,...,g_q)$,
we know that there exists $(h_1,...,h_r)$ in $M$
such that $a M $ is contained in $\gen{h_1,...,h_r}$.
\\
A matrix of syzygies for the \sgr $(g_1,...,g_q,h_1,...,h_r)$ is then of the following form
$\cmatrix{   a \I_q \cr    B}$
with $B$ of size $r\times q$.
This simply means that we can express $a g_j$ in terms of the $h_i$'s.
Therefore in the \idf of order $r$ of the module
there is a typical \gtr which is the \deter of $a \I_q$  \cadz~$a^q$.
Thus, every typical \gtr of the Kaplansky \id 
is in the nilradical of the corresponding Fitting \idz. 
Note that the exponent that intervenes here is simply the number of \gtrs of the module.
\\
Now if $a$ is a typical \gtr of  $\cF_r(M)$ we obtain  $a$ as a minor of order $q-r$ for a matrix of syzygies between $q$ \gtrs $(g_1, ..., g_q)$.
Even if it involves renumbering the \gtrsz, this matrix can be expressed as
$\cmatrix{N\cr D}$
where $D$ is a square matrix of order $q-r$, $N$ is an $r\times (q-r)$ matrix, and $\det D = a$.
\\
By \colis of the columns
(\prmt by right-multiplying by the cotransposed matrix of $D$) 
we obtain other syzygies for the same \gtrs in the form 
$\cmatrix{N'
\cr a \I_{q-r}}$
and this implies  that the last $q-r$ \gtrs multiplied by $a$ fall in the module generated by the first $r$ \gtrsz.
In short every typical \gtr of the Fitting \id is also a typical \gtr of the corresponding Kaplansky \idz.

\exer{exoPetitsExemplesElim}

\emph {2.}
We have $\ff \cap\gA = \gen {a,b,c}$ (if $x \in \gA$ satisfies $x \in \gen {X^d, a, b, c}_\AX$, make $X := 0$),
and also~$\fb = \gen {a,b,c}^d$. The \id $\fa$ is the content in $T$ of the \pol $(a + bT + cT^2)^d$ whereas $\fa'$ is the content in $\uT$ of the \pol $(aT_1 + bT_2 + cT_3)^d$. For example for $d = 2$:

\snic {
\fa = \gen {a^2, ab^2, 2ab, 2ac+b^2, b^3, b^2c, 2bc, c^2}, \quad
\fa' = \gen{  a^2, 2ab, 2ac, b^2, 2bc, c^2}.
}

We have $\fa \subsetneq \fa' \subsetneq \fb \subsetneq \ff \cap \gA$ and $\fb = (\ff\cap\gA)^d$. We also see that $\fa$ is not symmetrical in $a$, $b$, $c$. Still for $d = 2$, we have $(\ff\cap \gA)^4 \subsetneq \fa$ and $(\ff\cap \gA)^3 \not\subset \fa'$.  For  arbitrary $d$, it seems that $(\ff\cap \gA)^{3d-2} \subseteq \fa$.

\exer{exoRelateursViaElimIdeal} 
Let $\wi\varphi : \gk[\uX,\uY] \to \kuX$ be the \evn morphism $Y_i \mt f_i$, the base \ri being $\kuX$.  We~have

\snic {
\ker\wi\varphi = \gen{Y_1 - f_1, \ldots, Y_s - f_s} =
\gen{g_1, \ldots, g_s},
}

and since $\wi\varphi$ extends $\varphi$, $\ker\varphi = \kuY \cap \ker\wi\varphi$, as required.


\prob{exoDimZeroXcYbZa} 
First, the cycle $\sigma = (1,2,3)$ performs $\sigma(f_1) = f_2$, $\sigma(f_2) = f_3$ 
and~\hbox{$\sigma(f_3) = f_1$}. Therefore $C_3 = \gen{\sigma}$ operates on $\gA = \gk[x,y,z]$.
If in addition $a=b$ or $b=c$, then $\{f_1, f_2, f_3\}$ is invariant under $\rS_3$.
Finally, note that the origin is a zero of the \sysz, but also that solutions with $x = y = z \neq 0$ (in an extension of $\gk$) exist.

\emph {1.}
There are two cases: the $a \le b < c$ case, the easier one to study (case~I), and the $a < b=c$ case (case~II).

$\bullet$ case~I ($b < c$).
\\
Consider on the \moms of $\gk[X,Y,Z]$ the  order {\tt deglex}  (see Exercise~\ref{exothSymEl}). 
Let us show that $\gA = \sum_{p,q,r:\max(p,q,r) < c} \gk\, x^py^qz^r$. \\
Let $m = x^iy^j z^k$ with $\max(i,j,k) \ge c$.
If $i \ge c$, we replace in $m$, $x^c$ \hbox{with $x^{i-c}x^c = -x^{i-c}(y^b+z^a)$}. Similarly if $j \ge c$ or if $k \ge c$. We then get
$$\preskip.4em \postskip.4em 
m = -(m_1 + m_2) \hbox { with }  m_1, m_2 = \cases {
x^{i-c}y^{b+j}z^k,\ x^{i-c}y^jz^{a+k} & if $i \ge c$,\cr
x^{a+i}y^{j-c}z^k,\ x^iy^{j-c}z^{b+k} & if $j \ge c$,\cr
x^{b+i}y^jz^{k-c},\ x^iy^{a+j}z^{k-c} & if $k \ge c$.
} 
$$
We then see that $m_1 < m$ and $m_2 < m$; we finish by \recuz.
The reader will check that the $x^py^qz^r$ with $p$, $q$, $r < c$ form a $\gk$-basis of $\gA$.  For those familiar with the material: when $\gk$ is a \cdiz, $(f_1, f_2, f_3)$ is a Gr\"obner basis for the monomial order {\tt deglex}.
Recap: $\dim_\gk\gA = c^3$.

$\bullet$ case~II ($a < b=c$). This case is more difficult. 
\\
First suppose that $2$ is \iv in~$\gk$. We introduce

\snic {\arraycolsep2pt
\begin {array} {rcl}
g_1 &=& -f_1 + f_2 + f_3 = 2Z^c + X^a + Y^a - Z^a, \cr
g_2 &=& f_1 - f_2 + f_3 = 2X^c - X^a + Y^a + Z^a, \cr
g_3 &=& f_1 + f_2 - f_3  = 2Y^c + X^a - Y^a + Z^a .
\end {array}
}

We then have

\snic {
2f_1 = g_2 + g_3, \quad  2f_2 = g_1 + g_3,\quad 2f_3 = g_1 + g_2,
}


such that $\gen {f_1, f_2, f_3} = \gen {g_1, g_2, g_3}$. Then we can operate with the $g_j$'s as we did with the $f_i$'s in case~I. If $\gk$ is a \cdiz,
$(g_1, g_2, g_3)$ is a Gr\"obner basis for the  graded lexicographic order {\tt deglex}.  \\
Recap: $\dim_\gk\gA = c^3$ and the $x^py^qz^r$'s with
$p$, $q$, $r < c$ form a $\gk$-basis of $\gA$.  

 $\bullet$
Case II with a \cdi $\gk$ of \cara $2$ is left to the sagacity of the reader. The \ri $\gA$ is not always \zedz! 
This happens for example when $\gk = \FF_2$ and $(a,b) = (1,3)$, $(1,7)$, $(2,6)$, $(3,9)$.
When it is \zedz, it seems that $\dim_\gk\gA < c^3$.

\emph {2.}
For $(a,b,c) = (2,2,3)$, we know that $\dim_\gk \gk[x,y,z] = 3^3 = 27$. We use Stickelberger's \thoz~\ref{thStickelberger}, except that we do not know the zeros of the \sysz. We check, with the help of a Computer Algebra \sysz, that the \polcar of $x$ over~$\gk$ can be factorized into \irds \pols $(\gk = \QQ$)

\snic {
\rC{x} = t^8 (t+2) (t^3-t^2+1)^2  
       (t^4 - 2t^3 + 4t^2 - 6t + 4)  (t^4 + t^3 + t^2 - t + 2)^2
,}


but the \fcn of $\rC{x+2y}$ is of the type $1^8 \cdot 1^1 \cdot 4^1
\cdot 4^1 \cdot 4^1 \cdot 6^1$. Consequently, the \prn $(x,y,z) \mapsto
x$ does not separate the zeros of the \sysz, whereas the \prn $(x,y,z) \mapsto
x+2y$ does. Moreover, we see that the origin is the only zero with multiplicity (equal to $8$). Thanks to the \fcn of $\rC{x}$ and by performing a few additional small computations,
we obtain

\begin{itemize}
\item
Another zero defined over $\gk$, $(x,y,z) = (-2,-2,-2)$, which is simple. 
\item
If $\alpha$, $\beta$, $\gamma$ are the three distinct roots of $t^3 - t^2 +1$, we obtain $6$ simple zeros by making the group $\rS_3$ act on the zero $(\alpha, \beta, \gamma)$. If $s_1$, $s_2$, $s_3$ are the \elr \smq functions of $(X,Y,Z)$, then, over $\QQ$, we have the \egt of \ids $\gen {f_1, f_2, f_3, s_1-1} = \gen {s_1-1, s_2,s_3+1}$, \cad the \alg of these $6$ zeros is the \adu of the \pol $t^3 - t^2 + 1$.
\item
Let $\delta_i$ be a root of $t^4 + t^3 + t^2 - t + 2$ ($i\in\lrb{1..4}$).\\
By letting $y = x=\delta_i$ and $z = 2/(x+1) = -(x^3 + x - 2)/2$, we obtain a zero of the \sysz. The \polmin of $z$ over $\QQ$ is the one we see in the \fcn of $\rC{x}$: $t^4 - 2t^3 + 4t^2 - 6t + 4$.
We thus obtain four simple zeros of the \sysz.
\item
We can make $\rA_3$ act on the four previous zeros.
\end{itemize}
We have therefore obtained $1 + 6 + 3\times 4 = 19$ simple zeros and a zero of multiplicity~$8$. This adds up as required.
\\
Remark: whereas $\dim_\gk \gk[x,y,z] = 27$, we~have

\snic{
\begin {array} {c}
\dim_\gk \gk[x] = \dim_\gk \gk[y] = \dim_\gk \gk[z] = 14,
\\[1mm]
\dim_\gk \gk[x,y] = \dim_\gk \gk[x,z] = \dim_\gk \gk[y,z] = 23.
\end {array}
}

Thus, neither $\gk[x,y]$ nor $\gk[x,y,z]$ are free over $\gk[x]$,
and $\gk[x,y,z]$ is not free over~$\gk[x,y]$.

\emph{3.}
If $\gk$ is a \cdiz, in case~I in \cara $\ne 2$, we find, experimentally, that the local \alg of the origin is $\aqo{\gk[X,Y,Z]}{X^a, Y^a, Z^a}$ and so the multiplicity of the origin would be $a^3$.  As for case~II, this seems quite mysterious.


\prob{exoIdealResultantGenerique}
\emph{1.}
We put the following weights on $\gk[X]$: $X$ is of weight $1$, and the weight of $a_i$ and $b_{ji}$ is $i$. Thus $f$ and $g_j$ are \hmgs of weight $d$. We easily check for all $k \ge 0$ that $(X^k g_j) \bmod f$ is \hmg of weight $d+k$.

\emph{2.}
We index the $d$ rows of $S$ by $1, \ldots, d$, the $i^{\rm th}$ row corresponding to the weight~$i$ via $i \leftrightarrow X^{d-i} \leftrightarrow a_i$. The matrix $S$ is the horizontal concatenation of $r$ square matrices of order $d$,
the $j^{\rm th}$ square matrix being that of the
 multiplication by $g_j$
modulo $f$ in the basis $(X^{d-1}, \ldots, X, 1)$.  
If we number the columns of the first square submatrix of order $d$ of $S$ (corresponding to $g_1$) by $(0, 1, \ldots, d-1)$, then the \coe of index $(i,j)$ is \hmg of weight $i+j$. Similarly for the other \coes with analogous conventions.\\
For example, for $d = 3$, if $f = X^3 + a_1X^2 + a_2X + a_3$, $g = b_1X^2 + b_2X + b_3$, the matrix of the multiplication by $g \bmod f$ is

\snuc{
\bordercmatrix [\lbrack\rbrack]{
                          &g    &Xg \bmod f      &X^2g \bmod f  \cr
X^{d-1}\leftrightarrow 1  &b_1  & -a_1b_1 + b_2  &a_1^2b_1 - a_1b_2 - a_2b_1 + b_3 \cr
X^{d-2}\leftrightarrow 2  &b_2  & -a_2b_1 + b_3  &a_1a_2b_1 - a_2b_2 - a_3b_1 \cr
X^{d-3}\leftrightarrow 3  &b_3  & -a_3b_1        &a_1a_3b_1 - a_3b_2 \cr
}
\quad \hbox {of weights} \quad
\cmatrix {1 & 2 & 3\cr 2 & 3 & 4\cr 3 &4 &5 \cr}.
}

Let $M$ be a submatrix of order $d$ of $S$, $(k_1, \ldots, k_d)$ the exponents of $X$ corresponding to its columns ($k_i \in\lrb{0.. d-1}$, and the columns are $X^{k_i} g_j \bmod f$).\\
Then, $\det(M)$ is \hmgz, and its weight is the sum of the weights of the diagonal \coesz, \cad

\snic {
(1 + k_1) + (2 + k_2) + \cdots + (d + k_d) = d(d+1)/2 + \sum_{i=1}^d k_i.
}

For example, the weight of the first minor of order $d$ of $S$ (corresponding to  multiplication by $g_1$) is $d(d+1)/2 + \sum_{k=0}^{d-1}k = d^2$.  
\\
The weight of each of the $rd\choose d$ minors is bounded below by $d(d+1)/2$ (bound
obtained for $k_i=0$) and bounded above by~$ d(3d-1)/2$ (bound obtained for $k_i = d-1$). These bounds are reached if $r \ge d$.

\emph{3.}
The number $\dim_\QQ E$ is the lower bound of the cardinality of any arbitrary \sgr of~$\fb$. We experimentally find, for small values of $r$ and $d$, that~$\dim_\QQ E = r^d$. But we can do better. Indeed, the consideration of graded objects allows us to assert the following result (\hmg Nakayama lemma, \pb~\ref{exoNakayamaHomogeneRegularSequence}): every graded family of $\fb$ whose image in $E$ is a \hmg \sgr of the graded $\QQ$-\evc $E$ is a (\hmgz)
\sgr of $\fb$. In~particular, there exists a \hmg \sgr of $\fb$ of cardinality $\dim_\QQ E$, conjecturally, $r^d$. We can go further by examining the weights of the minimal \hmg \sgrs of $\fb$. Those are unique and provided by the (finite) series of the graded $\QQ$-\evc $E$. For example, for $d=5$, $r=2$, this series is

\snic {
6t^{25} + 4t^{24} + 6t^{23} + 6t^{22} + 6t^{21} + 2t^{20} + 2t^{19},
}

which means that in any minimal \hmg \sgr of~$\fb$,
there are  $6$ \pols of weight~$25$, 
$4$ \pols of weight~$24$, \dots, $2$ \pols of weight $19$ (with $6 + 4 + \cdots + 2 = 32 = 2^5 = r^d$).
In this example, the number ${rd \choose d}$ of minors of order $d$ of $S$ is $252$.

Conjecturally, it would seem that $\fb$ is generated by \hmg \pols of weight $\le d^2$, with $d+r-1 \choose r-1$ \pols of weight $d^2$ exactly.

\prob{exoNakayamaHomogeneRegularSequence} ~\\
\emph {1.}
We perform a \dem by \recu on $n$. 
\\ 
\emph{Case $n=0$}: trivial result.
\\
\emph{For $n \ge 1$}, we consider $\gA' = \aqo{\gA}{a_1}$. We have $\gk \hookrightarrow \gA'$ because $\gk \cap \gen{a_1} = \{0\}$. The sequence $(\overline {a_2}, \cdots, \overline {a_n})$ in $\gA'$ satisfies the right assumptions for the \recu on $n$. Suppose $f(a_1, \ldots, a_n) = 0$ with $f\in\kXn$ and $\deg_{X_1}(f)\le d$. \\
We write $f = X_1q(X_1, \ldots, X_n) + r(X_2, \ldots, X_n)$ with
$q$, $r$ with \coes in $\gk$ and~$q$ of degree $\le d-1$ in $X_1$. 
In $\gA'$, we have $r(\overline {a_2}, \ldots, \overline {a_n}) = 0$. 
By \recu on~$n$, we have $r=0$.
Since $a_1$ is \ndzz, $q(a_1, \ldots, a_n)= 0$. By \recu on $d$, we obtain $q=0$, so $f=0$.

\emph {2a.}
By \dfnz, $\gA_+E \subseteq E_1 \oplus E_2 \oplus \dots$; and since $\gA_+E = E$,
we get $E_0 = 0$. Then $\gA_+E \subseteq E_2 \oplus E_3 \oplus \dots$, and by using $\gA_+E = E$ again, we get $E_1 = 0$, and so on. So $E_n = 0$ for all $n$, therefore $E = 0$.

\emph {2b.}
Let $F$ be the \Asub of $E$ generated by the $e_i$'s. It is a graded submodule because the $e_i$'s are \hmgsz. The hypothesis is equivalent to $F + \gA_+ E = E$ or $\gA_+(E/F) = E/F$. By question~\emph{2a}, we have $E/F = 0$ \cad $E = F$; the~$e_i$'s generate the \Amo~$E$.

\emph {3a.}
It is clear that $\gA_0 = \gB_0$ and  $\fb = \gA_+ \gB$. By applying the previous question \hbox{to the graded \Amoz} $\gB$ and to $e_i$,
we obtain that the $e_i$'s form a \sgr of the \Amo $\gB$.

\emph {3b.}
Let $S = \sum_i \gB_0 e_i$ (actually, it is a direct sum).\\ 
Let us show that $\gen {h_1, \ldots, h_d} \cap S = \{0\}$. If $s = \sum_i \lambda_i e_i \in \gen {h_1, \ldots, h_d}$ with $\lambda_i \in \gB_0$, then by reducing modulo $\gen {h_1, \ldots, h_d}$, we get $\sum_i \lambda_i \overline {e_i} = 0$, therefore $\lambda_i = 0$ for all $i$ and $s = 0$.
\\
For $\alpha = (\alpha_1, \ldots, \alpha_d) \in \NN^d$, let $h^\alpha = h_1^{\alpha_1} \cdots h_d^{\alpha_d}$. Let us show that
$$
\som_\alpha s_\alpha h^\alpha = 0 \hbox { with } s_\alpha \in S
\;\Longrightarrow\; s_\alpha = 0 \;\; \hbox { for all }\alpha.
\leqno (\star)
$$
For this, we will prove by (decreasing) \recu on $i$, that
$$
\preskip.4em \postskip.4em 
\bigl( f \in S[X_i, \ldots, X_d] \hbox { and }
f(h_i, \ldots, h_d) \equiv 0 \bmod \gen {h_1, \ldots, h_{i-1}} \bigr)
\;\Longrightarrow\; f = 0. 
$$

\emph{First for $i = d$}. The hypothesis is $s_m h_d^m + \cdots + s_1 h_d + s_0 \equiv 0 \bmod \gen {h_1, \ldots, h_{d-1}}$ and we want $s_k = 0$ for all $k$.  We have $s_0 \in S \cap \gen {h_1, \ldots, h_d} = \{0\}$. We can simplify the congruence by $h_d$ (which is \ndz modulo $\gen {h_1, \ldots, h_{d-1}}$) to obtain $s_m h_d^{m-1} + \cdots + s_1 \equiv 0 \bmod \gen {h_1, \ldots, h_{d-1}}$. By iterating the process, we obtain that all the $s_k$'s are null.
\\
\emph{Passing from $i+1$ to $i$}. \\
Let $f \in S[X_i, \ldots, X_d]$ of degree $\le m$ with $f(h_i, \ldots, h_d) \equiv 0 \bmod \gen {h_1, \ldots, h_{i-1}}$. We write $f = X_i q(X_i, \ldots X_d) + r(X_{i+1}, \ldots, X_d)$ with $q$, $r$ with \coes in~$S$ and~$q$ of degree $\le m-1$.  We therefore have $r(h_{i+1}, \ldots, h_d) \equiv 0 \bmod \gen {h_1, \ldots, h_i}$, hence by
\recu on $i$, $r = 0$. We can simplify the congruence by $h_i$ (which is \ndz modulo $\gen {h_1, \ldots, h_{i-1}}$) to obtain $q(h_i, \ldots, h_d) \equiv 0 \bmod \gen {h_1, \ldots, h_{i-1}}$. \hbox{Therefore $q = 0$} by \recu on $m$, then $f = 0$.
\\
Recap: we therefore have the result for $i = 1$ and this result is none other than~$(\star)$.
\\
Once $(\star)$ is proved, we can show that the $e_i$'s are \lint independent over~$\gA$. Let $\sum_i a_i e_i = 0$ with $a_i \in \gA$; we write $a_i = \sum_{\alpha} \lambda_{\alpha,i} h^\alpha$ and
$$
\som_i a_ie_i = \som_{i,\alpha} \lambda_{i,\alpha} h^\alpha e_i =
\som_{\alpha} s_\alpha h^\alpha 
\quad \hbox {with} \quad s_\alpha = \som_i \lambda_{i,\alpha} e_i \in S.
$$
Therefore $s_\alpha = 0$ for all $\alpha$, then $\lambda_{i,\alpha} = 0$ for all $i$, and  $a_i = 0$.

\emph {4.}
\Gnlt, if $(a_1, \ldots, a_d)$ is a \seqreg of a \ri $\gA$, it is $L$-\ndze for all free $\gA$-modules $L$ (left to the reader). We apply this to the \ri $\gA = \gB_0[h_1, \ldots, h_d]$, to the sequence $(h_1, \ldots, h_d)$ (which is indeed a \seqreg of~$\gA$) and to $L = \gB$ (which is a free \Amo by the hypothesis).




\Biblio
Bourbaki (Algebra, Chapter X, or Commutative algebra, Chapter I) 
calls what we have called a coherent module (in accordance with the common usage, especially in the English literature) a \emph{pseudo coherent module}, and what we call a finitely presented coherent module Bourbaki calls a coherent module.
 This is naturally to be related to the \gui{Faisceaux Alg\'ebriques Coh\'erents} 
\hbox{by J.-P. Serre} (precursors of sheaves of modules on a Grothendieck scheme) which are locally given by \pf \coh modules.
It should also be noted that \cite{Stacks} adopts the Bourbaki's \dfn for \coh modules.

\smallskip 
\Thref{prop unicyc} is taken from \cite{MRR} Chap.\ V, Th.\ 2.4.
\Thref{prop quot non iso} is taken from \cite{MRR}  Chap.\ III, Exercise 9 p.\ 80.

The standard reference for \idfs is \cite{Nor}.

As for purely equational algebraic structures and universal algebra
 one can consult~\cite{BuSa}. 

A first introduction to categories is found in~\cite{Cohn}. 
\\ Dedicated books on the subject that we can recommend are \cite{MACL} and~\cite{LaRo}.

The Kaplansky \ids of a module $M$ studied in Exercise~\ref{exoAutresIdF} are used in~\cite[Chap.~IV]{Kun} and~\cite[Chap. 9]{IRa}.

The strict Bézout \ris (Exercise~\ref{exoAnneauBézoutStrict}) and the Smith \ris 
(Exercise~\ref{exoSmith}) have been studied by Kaplansky in \cite{Kap} in a more \gnl framework of not \ncrt commutative \risz. He respectively calls them \gui{Hermite rings} and \gui{elementary divisor rings.} But this terminology is not fixed.
In \cite{Lam06}, where Exercise~\ref{exoAnneauBézoutStrict} finds its source, Lam uses~\hbox{\emph{K-Hermite ring}} for strict Bézout \riz. That is to be distinguished from \emph{Hermite ring}: today a \ri $\gA$ is called a 
\emph{Hermite \riz} if every stably free \Amo is free, \cad if every \vmd 
is completable (see Chapter~\ref{chap ptf0}, Section~\ref{ModStabLibre}).
As for the \gui{\elr divisors,} they are now often used in a more particular sense. For example, 
in the literature it is often said
 that the \ZZmo $$\aqo\ZZ{900}\oplus\aqo\ZZ{10}\simeq \aqo\ZZ{25}\oplus\aqo\ZZ{5}\oplus\aqo\ZZ{4}\oplus\aqo\ZZ{2}\oplus\aqo\ZZ{9}$$ admits for invariant factors the list $(10,900)$ and for \elr divisors the unordered list $(25,5,4,2,9)$.%
\index{ring!Hermite ---}

Exercise~\ref{exoCalculT(M)} was provided to us by Thierry Coquand.

\newpage \thispagestyle{CMcadreseul}
\incrementeexosetprob


\chapter{\Mptfsz, 1}
\label{chap ptf0}
\perso{compil\'e le \today}
\minitoc

\section{Introduction}

Recall that a \mptf is a module \isoc to a direct summand in a free \Amo of finite rank. 
This notion happens to be the natural \gnnz, for  modules over a commutative \riz, of the notion of a finite dimensional \evc over a \cdiz.
This chapter develops the basic theory of these modules.

\smallskip  
One of the initial motivations of this book was to understand \emph{in concrete terms} the following \thos concerning \mptfsz.

\begin{theorem}\label{th.ptf.loc} {\em  (Local structure \tho for \mptfsz)}
An \Amo $P$ is \ptf \ssi it is {\em locally free} in the following sense. There exist \ecoz~$s_1$,~$\ldots$,~$s_{\ell}$ in~$\gA$ such that the modules~$P_{s_i}$ obtained from $P$ by \eds to the \risz~$\gA_{s_i}=\gA[1/s_i]$ are free.
\end{theorem}

\begin{theorem}\label{th.ptf.Fitting} {\em  (\Carn of  \mptfs by their  \idfsz)} A \pf \Amo is \pro \ssi its \idfs are (\idps generated by) idempotents.
\end{theorem}


\begin{theorem}\label{th.ptf.idpt} {\em  (Decomposition of a \mptf into a direct sum of modules of constant rank)} If $P$ is a \ptf \Amo generated by $n$ \eltsz, there exists a \sfio $(r_0, r_1, \ldots, r_n)$ (some eventually null) 
such that each  $r_kP$ is a projective module of rank~$k$ over the \riz~$\aqo{\gA}{1-r_k}$. Then $P=\bigoplus_{k>0}r_kP$ and $\Ann(P)=\gen{r_0}$.
\end{theorem}

In this direct sum we can naturally limit ourselves to the indices $k>0$ such that $r_k\neq 0$.

\pagebreak	

\begin{theorem}\label{th.ptf.plat} {\em  (\Carn of  \mptfs by their flatness)} 
	A \pf \Amo is projective \ssi it is flat.
\end{theorem}

In this chapter we will prove the first three of these \thosz.
They will be taken up again with new \dems in Chapter~\ref{chap ptf1}.
The fourth one will be proven in Chapter~\ref{chap mod plats}, which is dedicated to flat modules.

Other important \thos regarding \mptfs will be proven in Chapters~\ref{chap ptf1}, \ref{chapNbGtrs} and~\ref{ChapMPEtendus}.
The theory of  \algs which are \mptfs (we will call them \asfsz) is developed in Chapter~\ref{chap AlgStricFi}.

\section{Generalities}
\label{subsecPropCarPTFS}

Recall that a \mptf is \pf (Example 2, \paref{exl1pf}).

\subsec{Characteristic \prtsz}

When $M$ and $N$ are two \Amosz, we have a natural \Ali 
$\theta_{M,N}:M\sta\te N \to \Lin_\gA(M,N)$ given by
\begin{equation}\preskip.4em \postskip.35em
\label{NOTAthetaMN}
\theta_{M,N}(\alpha \te y)= \big(x\mapsto \alpha (x)y\big).
\end{equation}
We also write $\theta_M$ for $\theta_{M,M}$.\label{NOTAthetaM}

\medskip \rem We sometimes write $\alpha \te y$ for $\theta_{M,N}(\alpha \te y)$
but it is certainly not recommended when $\theta_{M,N}$ is not injective.
\eoe

\medskip The following \tho gives some immediately equivalent \prtsz.
\begin{theorem} {\em  (\Ptf modules)}\label{propdef ptf}\\
For an \Amo  $P$, \propeq
\begin{enumerate}
\item  [$(a)$] $P$ is a \emph{finitely generated projective module}, \cad there exist an integer~$n$, an \Amo $N$ and an \iso of $P\oplus N$ over~$\Ae n$.
\index{module!finitely generated projective ---}%
\index{finitely generated projective!module}%
\index{projective!finitely generated --- module}

\item  [$(b1)$] There exist an integer $n$, \elts $(g_i)_{i\in\lrbn}$ of $P$ and \lins forms $(\alpha_i)_{i\in\lrbn}$ over $P$ such that for all $x\in P,   \;    x = \sum_i \alpha_i (x)\,  g_i$.

\item  [$(b2)$] The module $P$ is finitely generated, and for every finite \sys of \gtrs
$(h_i)_{i\in\lrbm}$ of $P$ there exist \lin forms
$(\beta_i)_{i\in\lrbm}$ over $P$ such that for all~\hbox{$x\in P$},       $x = \sum_i  \beta_i (x) \, h_i$.

\item  [$(b3)$] The image of $P\sta\te_\gA P$ in $\Lin_\gA(P,P)$ under the canonical \homoz~$\theta_P$ contains $\Id_P$.

\item  [$(c1)$] There exist an integer $n$ and two \alis   
$\varphi : P\rightarrow\Ae n$  
 and~$\psi : \Ae n\rightarrow P$, such that $\psi\circ \varphi= \Id_P$. We then have~$\Ae n=\Im(\varphi )\oplus\Ker(\psi)$ and~$P\simeq\Im(\varphi\circ\psi )$.

\item  [$(c2)$] The module $P$ is \tfz, and for every surjective \ali 
$\psi :\Ae m\rightarrow P$, there exists a \aliz~$\varphi : P\rightarrow \Ae m$  such that~$\psi\circ \varphi= \Id_P$.  
We then have~$\Ae m=\Im(\varphi )\oplus\Ker(\psi)$ and~$P\simeq\Im(\varphi\circ\psi )$.

\item  [$(c3)$] Like $(c2)$ but by replacing $\Ae m$ by an arbitrary \Amo $M$: the module $P$ is \tfz, and for every surjective \ali $\psi :M\rightarrow P$,
$\Ker(\psi)$ is a direct summand.

\item  [$(c4)$]
The module $P$ is \tf and the functor $\Lin_\gA(P,\bullet)$ transforms the surjective \alis
into surjective maps. \\
In other words, for all \Amos $M$, $N$, for every surjective \ali $\psi : M\rightarrow N$ and every \ali $\Phi:P\rightarrow N$, there exists a \ali $\varphi : P\rightarrow M$ such that $\psi\circ \varphi= \Phi$.

\snic{
\xymatrix {
                                         & M\ar@{>>}[d]^{\psi} \\
P\ar@{-->}[ur]^{\varphi} \ar[r]_{\Phi} & N \\
}}

\end{enumerate}
\end{theorem}

\begin{proof}
Item $(b1)$ (resp.\,$(b2)$) is simply a reformulation of $(c1)$ (resp.\,$(c2)$).
\\
Item $(b3)$ is simply a reformulation of~$(b1)$.
\\
We trivially have $(c3) \Rightarrow (c2)\Rightarrow  (c1)$.

$(a) \Rightarrow(c1)$ Consider the canonical maps 
$$\preskip.4em \postskip.4em 
P\rightarrow P\oplus N \hbox{ and } P\oplus N\rightarrow P. 
$$ 

 $(c1) \Rightarrow(a)$ 
Consider $\pi=\varphi\circ\psi$. We have $\pi^2=\pi$. 
This defines a \prn of $\Ae n$ over $\Im\pi=\Im \varphi\simeq P$ \paralm to~$N=\Ker\pi=\Ker\psi$.

  $(b1)\Rightarrow (c4)$ If $\Phi(g_i)=\psi(y_i)$ ($ i\in\lrbn$), we let $\varphi(x)=\sum \alpha_i(x)\,y_i$. We then have for all $x\in P$,
$$\preskip.4em \postskip.4em\ndsp
\Phi(x) =  \Phi\big(\sum \alpha_i(x)\,g_i\big) = \sum \alpha_i(x)\,
\psi(y_i) =\psi \big(\sum \alpha_i(x)\,y_i\big) = \psi \big(\varphi(x)\big).
$$

 $(c4)\Rightarrow (c3)$ We take $N=P$ and $\Phi=\Id_P$.
\end{proof}

We also directly have $(b1)\Rightarrow (b2)$ as follows: by expressing the $g_i$'s as \colis of the $h_j$'s we obtain the $\beta_j$'s from the~$\alpha_i$'s.

In practice, according to the original \dfnz, we consider a \mptf as (an isomorphic copy of) the image of a \mprnz~$F$. Such a matrix, or the \ali that it represents, is again called a \emph{projector}. More \gnltz, every idempotent \endo of a module $M$ is called a \ix{projector}.

When we see a \mptf according to the \dfn $(c1)$, the \prn matrix is that of the \ali $\varphi\circ\psi$. Similarly, if we use the \dfn $(b1)$, the \mprn is the one which has for \coes every $\alpha_i(g_j)$ in position $(i,j)$.

A \sys $\big((g_1,\ldots,g_n),(\aln)\big)$ that satisfies $(b1)$ is called a \ixe{\sycz}{coordinate system} for the \pro module $P$.
Some authors speak of a \emph{basis} of the \mptfz, but we will not be following their lead on this matter.

\begin{fact}
\label{factDualPTF}\label{lemDual} \emph{(Dual of a \mptfz, 1)}\\
Let $\big((g_1,\ldots,g_n),(\aln)\big)$ be a \syc for a \mptf $P$. Then
\begin{enumerate}
\item  [--]   the $g_i$'s generate $P$,
\item  [--]   the $\alpha_j$'s generate $\Lin(P,\gA)=P\sta$,
\item  [--]   the module $P\sta$ is \ptfz,
\item  [--]   the module $(P\sta)\sta$ is canonically \isoc to $P$,
\item  [--]   via this canonical identification, $\big((\aln),(g_1,\ldots,g_n)\big)$ is a \syc for $P\sta$.
\end{enumerate}
In particular, if $P$ is (\isoc to) the image of a \prn matrix~$F$,
the dual module $P\sta$ is (\isoc to) the image of the \prn matrix~$\tra{F}$.
\end{fact}
\begin{proof}
The first item is clear. All the rest is clear from the moment where we show that $\lambda  =\sum \lambda(g_i)\, \alpha_i $ for all $\lambda \in P\sta$, and this \egt is proven by evaluating both sides at an arbitrary \elt $x$ of~$P$:
$$\preskip.3em \postskip-1em\ndsp 
\lambda(x)= \lambda\bigl(\sum\, \alpha_i (x)\,  g_i\bigr)=
\sum\, \alpha_i (x)\,  \lambda(g_i)=
\bigl(\sum\, \lambda(g_i) \, \alpha_i \bigr)(x). 
$$
\end{proof}

\begin{theorem}
\label{prop pf ptf}
Let $\Ae m\vers{\psi} \gA^q\vers{\pi}P\to 0$ be a \pn of a module~$P$.
Then, $P$ is \ptf \ssi $\psi$ is \lnlz.
\end{theorem}
Recall that \gui{$\psi$ is \lnlz} means that there exists a $\varphi :\gA^q\rightarrow \Ae m$ satisfying $\psi\,\varphi\,\psi=\psi$. Moreover, by \thref{propIGCram} every \ali which has a rank in the sense of \Dfnz~\ref{defRangk} is \lnlz.

\begin{proof} If $\psi$ is \lnlz, Fact~\ref{factInvGenCrois} tells us that $\Im\psi$ is a direct summand, and $\Coker\psi$ is \isoc to a 
complementary submodule of $\Im\psi$.
Conversely, if the module $P:=\Coker\psi$ is projective, we apply the \prtz~$(c2)$ of \thref{propdef ptf} to the \prnz~$\pi:\gA^q\rightarrow P$.
We obtain~$\tau:P\rightarrow \gA^q$ with~$\pi\circ\tau=\Id_P$, such that $\gA^q=\Im \tau \oplus \Im\psi$. Therefore $\Im\psi$ is \ptf and we can apply the \prtz~$(c2)$ to $\psi :\Ae m\rightarrow \Im\psi$, which gives us $\varphi$ over the component $\Im\psi$ (and we take for example $0$ over~$\Im \tau$).
\end{proof}

\subsec{Local-global principle}

The fact that an \Amo is \ptf is a local notion in the following sense.

\pagebreak	

\begin{plcc}
\label{plcc.cor.pf.ptf}
\emph{(\Ptf modules)}
Let $S_1$, $\ldots$, $S_n$ be \moco of $\gA$ and $P$ be an \Amoz. 
If the~$P_{S_i}\!$'s are free,~$P$ is \tf and \proz. 
 \\
More \gnltz, the module $P$ is \tf and \pro \ssi the $P_{S_i}\!$'s are \ptf $\gA_{S_i}$-modules.
\end{plcc}
\begin{proof}
This results from \thref{prop pf ptf}, from the \plgref{plcc.pf} for \mpfs and from the \plgref{fact.lnl.loc} for \lnl \alisz.
\end{proof}

The \plgref{plcc.cor.pf.ptf} establishes the implication \gui{if} in \thref{th.ptf.loc}. The converse \gui{only if} has been proven in \thref{theoremIFD} which will give us \thref{prop Fitt ptf 2}.
We will give for this converse a more precise statement and a more conceptual \dem with \thref{th ptf loc free}.

\subsec{\Pro modules and Schanuel's lemma}

The notion of a \pro module can be defined for modules which are not \tfz.
In the following we will rarely use such modules, but it is however useful to give some precisions on this subject.

\begin{definition}\label{defiMPRO}
An \Amo $P$ (not \ncrt \tfz) is said to be \ixc{projective}{module} if it satisfies the following \prtz.
\\
For all \Amos $M,\,N$, for every surjective \ali $\psi : M\rightarrow N$ and every \ali $\Phi:P\rightarrow N$, there exists a \ali $\varphi : P\rightarrow M$ such that
$\psi\circ \varphi= \Phi$.
$$\preskip-.6em \postskip.0em 
\xymatrix {
                                       & M\ar@{>>}[d]^{\psi} \\
P\ar@{-->}[ur]^{\varphi} \ar[r]_{\Phi} & N \\
} 
$$
\end{definition}
%

Thus, given the \carn $(c4)$ in \thref{propdef ptf}, an \Amo is \tf \pro \ssi it is \pro and \tfz.
\\
In the following fact, the last \prt resembles the implication \hbox{$(c4)$ $\Rightarrow$  $(c3)$} in this \thoz.
\\
\rdb
A \ali $\varphi:E\to F$ is called a \emph{split surjection} if there exists a $\psi:F\to E$ with $\varphi\circ \psi=\Id_F$.
In this case we say that $\psi$ is a \emph{section} of $\varphi$, and we have~$E=\Ker\varphi\oplus \psi(F)\simeq \Ker\varphi\oplus F.$
\\
A short exact sequence is said to be \emph{split} if its surjection is split. 
\index{split!surjection}\index{surjection!split ---}%
\index{split!short exact sequence}%
\index{short exact sequence!split ---}%
\index{section!of a split surjection}
\begin{fact}\label{factdefiMPRO} ~
\begin{enumerate}
\item A free module whose basis is a set in bijection with $\NN$ is \proz.
For example the \ri of \pols $\AX$ is a \pro \Amoz.
\item Every module that is a direct summand in a \pro module is \proz.
\item 
If $P$ is \proz, every short exact sequence $0\to N\to M\to P\to 0$ splits.

\end{enumerate}
\end{fact}

\comm In \coma the free modules \emph{are not always} \proz.
Furthermore, it seems impossible to represent every module as a quotient of a free and \pro module. 
Similarly it seems impossible to place every \pro module as a direct summand in a free and \pro module.
For more details on this matter consult Exercise~\ref{propfreeplat}
and~\cite{MRR}.
\eoe
\perso{peut-\^etre give un exo? cela risque d'\^etre d\'emoralisant.}

\begin{lemma}\label{lemScha}
Consider two surjective \Alis with the same image
$ P_1  \vers{\varphi_1}   M  \rightarrow  0  $, $
   P_2  \vers{\varphi_2}   M  \rightarrow  0$
with the modules $P_1$ and $P_2$ being \proz.
\begin{enumerate}
\item  There exist reciprocal \isos $\alpha$, $\beta:P_1\oplus P_2\to P_1\oplus P_2$ such that $(\varphi_1\oplus 0_{P_2})\circ \alpha=0_{P_1} \oplus \varphi_2$ and  $\varphi_1\oplus 0_{P_2}=(0_{P_1} \oplus \varphi_2)\circ \beta$.
\item If we let $K_1=\Ker\varphi_1$ and $K_2=\Ker\varphi_2$, we obtain by restriction of~$\alpha$ and~$\beta$ reciprocal \isos between $K_1\oplus P_2$ and $P_1\oplus K_2$.
\end{enumerate}
\end{lemma}
\begin{proof}
There exists a $u:P_1\to P_2$ such that $\varphi_2\circ u=\varphi_1$ and $v:P_2\to P_1$ such that~$\varphi_1\circ v=\varphi_2$.
 $$
\xymatrix @R = 0.3cm{
P_1 \ar@{-->}[dd]_{u} \ar[rd]^{\varphi_1} \\
                                        & M \\
P_2 \ar@{>>}[ru]_{\varphi_2} \\
}
\qquad\qquad
\xymatrix @R = 0.3cm{
P_1 \ar@{>>}[rd]^{\varphi_1} \\
                                        & M \\
P_2 \ar@{-->}[uu]_{v}\ar[ru]_{\varphi_2} \\
}
\qquad\qquad
\xymatrix @R = 0.4cm{
P_1 \oplus P_2 \ar@/^/@{-->}[dd]^{\beta}
                 \ar[rd]^{~\varphi_1 \oplus 0_{P_2}} \\
                                        & M \\
P_1 \oplus P_2 \ar@/^/@{-->}[uu]^{\alpha}
    \ar[ru]_{~0_{P_1} \oplus \varphi_2} \\
}
$$
We verify that $\alpha$ and $\beta$ defined by the matrices below are suitable.
 $$
 \alpha=\Cmatrix{2pt}{\Id_{P_1}-vu&  v\cr -u& \Id_{P_2} }
\quad\quad
 \beta=\Cmatrix{2pt}{\Id_{P_1}& -v\cr u& \Id_{P_2}-uv }
  .$$
NB: the matrix $\beta$ is a sophisticated variant of what would be the cotransposed matrix of $\alpha$ if $\Id_{P_1}$, $\Id_{P_2}$, $u$ and $v$ were scalars.
\end{proof}

%

\begin{corollary}\label{corlemScha} \emph{(Schanuel's lemma)}
Consider two exact sequences
\vspace{-1mm}
\[
\begin{array}{ccccccccc}
0 &\rightarrow& K_1& \vers{j_1} & P_1& \vers{\varphi_1} & M& \rightarrow& 0   \\
0 &\rightarrow& K_2& \vers{j_2} & P_2& \vers{\varphi_2} & M& \rightarrow& 0
 \end{array}
\]
with the modules $P_1$ and $P_2$ being \proz. Then, $K_1\oplus P_2\simeq K_2\oplus P_1.$
\end{corollary}

\subsec{The category of  \mptfsz}
\label{secCatMptf}
\vspace{3pt}
\subsubsection*{A purely categorical construction}

The category of  \mptfs over $\gA$ can be constructed from the category of  free modules of finite rank over $\gA$ by a purely categorical procedure.
\begin{enumerate}
\item  A \mptf $P$ is described by a pair $(\rL_P,\rPr_P)$ where~$\rL_P$ is a free module of finite rank and $\rPr_P\in\End(\rL_P)$ is a \prrz.
We have $P\simeq \Im\rPr_P\simeq \Coker(\Id_{\rL_P}-\rPr_P)$.
\item  A \ali $\varphi$ from the module $P$ (described by $(\rL_P,\rPr_P)$) to the module $Q$ (described by $(\rL_Q,\rPr_Q)$) is described by a \ali $\rL_\varphi:\rL_P\rightarrow \rL_Q$  subjected to commutation relations

\snic{\rPr_Q\circ\,\rL_\varphi=\rL_\varphi=\rL_\varphi\circ\rPr_P.}

In other words $\rL_\varphi$ is null over $\Ker(\rPr_P)$ and its image is contained in $\Im(\rPr_Q)$.
\item  The identity of $P$ is represented by $\rL_{\Id_P}=\rPr_P$.
\item  The sum of two \alis $\varphi$ and $\psi$ from $P$ to $Q$ represented by $\rL_\varphi$ and $\rL_\psi$ is represented by $\rL_\varphi+\rL_\psi$.
The \ali $a\varphi$ is represented by~$a\rL_\varphi$.
\item  To represent the composition of two \alisz, we compose their representations.
\item  Finally, a \ali $\varphi$ from $P$ to $Q$ represented by $\rL_\varphi$ is null \ssi $\rL_\varphi=0$.
\end{enumerate}

This shows that the \pbs relating to the \mptfs can always be interpreted as \pbs regarding \prn matrices, and often come down \label{exMptf} to \pbs about solving \slis over $\gA$.

\smallskip  An \eqve category, better adapted to computations, is the category whose objects are the \mprns with \coes in $\gA$, a morphism from $F$ to $G$ being a matrix $H$ of a suitable format satisfying the \egts 

\snic{GH=H=HF.}

\subsubsection*{Using \sycsz}

The following fact uses the assertions of the previous paragraph 
while taking the coordinate system point of view.

\begin{fact}\label{factMatriceAlin}
Let $P$ and $Q$ be two \mptfs with \sycs 

\snic{\big((\xn),(\aln)\big)$ and $\big((\ym),(\beta_1,\ldots,\beta_m)\big),}

and let $\varphi:P\to Q$ be an \Aliz.\\
Then, we can encode $P$ and $Q$ by the matrices 

\snic{F\eqdefi \big(\alpha_i(x_j)\big)_{i,j\in\lrbn}\quad\et\quad G\eqdefi \big(\beta_i(y_j)\big)_{i,j\in\lrbm}.}

More precisely, we have the \isos

\snic{
\begin{array}{rcl} 
\pi_1: P\to\Im F\,,  &   & x\mapsto \tra[\,\alpha_1(x) \;\cdots\; \alpha_n(x)\,],  
 \\[1mm] 
\pi_2: Q\to\Im G\,,  &   & y\mapsto \tra[\,\beta_1(y) \;\cdots\; \beta_m(y)\,].   \end{array}
}

As for the \ali $\varphi$, it is encoded by the matrix 

\snic{H\eqdefi \big(\beta_i(\varphi(x_j))\big)_{i\in\lrbn,j\in\lrbm}}

which satisfies $GH=H=HF$. The matrix $H$ is that of the \ali 

\snic{\Ae n\to \Ae m, \quad\pi_1(x)+z\;\mapsto\; \pi_2 \big(\varphi(x)\big)\quad$ if $ x\in P $ and $z\in\Ker F.}

 We say that \emph{the matrix $H$ represents the \ali $\varphi$ in the \sycs $\big((\ux),(\ual)\big)$
and $\big((\uy),(\und{\beta})\big)$}.%
\index{matrix!of a \ali in \sycsz}
\end{fact}

\subsubsection*{Application: the \isos between \tf \pro modules}

Lemma \ref{propIsoIm} says that, for $F\in\GA_m(\gA)$ and $G\in\GA_n(\gA)$, if 
$\Im F$ and~$\Im G$ are \isocz, even if it means \gui{enlarging} 
the matrices $F$ and $G$, they can be assumed to be similar.

In the sequel, we use the notation $\Diag(M_1,\dots,M_k)$ more freely than we have until now. 
 Instead of a list of \elts of the \riz, we consider for $(M_1,\dots,M_k)$ a list of square matrices. The matrix represented as such is usually called a \emph{block diagonal matrix}.%
\label{Notadiagblocs}%
\index{matrix!block diagonal ---}%
\index{block diagonal!matrix}

\begin{lemma}
\label{propIsoIm} \emph{(Enlargement lemma)}\\
Consider the matrix encoding of the category of \mptfsz.
If an \iso $\varphi$ of $\Im F$ over $\Im G$ is encoded by $U$ and its inverse encoded by $U',$ we obtain a matrix $A\in\EE_{n+m}(\gA)$
$$A=\Cmatrix{2pt}{\I_m-F&-U'\cr U&\In-G}=\Cmatrix{2pt}{\I_m&0 \cr U &\In }\,
\Cmatrix{2pt}{ \I_m&-U' \cr 0& \In}\,\Cmatrix{2pt}{\I_m &0 \cr U&\In },$$
with
\begin{equation}\label{eqpropIsoIm}\preskip-.2em 
\bloc{0_m}{0}{0}{G}=A\,\bloc{F}{0}{0}{0_n}\,A^{-1}.
\end{equation}
Conversely, a conjugation between $\Diag(0_m,G)$ and $\Diag(F,0_n)$ provides an \iso between
$\Im F$ and $\Im G$.
\end{lemma}
\begin{proof} The following matrix
$$
\preskip.2em 
\bordercmatrix [\lbrack\rbrack]{
       & \Im F & \Ker F   & \Im G & \Ker G \cr
\Im F &    0    &  0 &       -\varphi^{-1}   &  0\cr
\Ker F &    0    &    \Id   &    0  &  0\cr
\Im G &    \varphi     &    0   &    0  &  0\cr
\Ker G &    0    &    0   &    0  &  \Id
}  ,$$
once $\Im F \oplus \Ker F$ is replaced by $\Ae m$ and $\Im G \oplus \Ker G$ by $\Ae n$, gives the matrix $A$. The presence of the $-$ sign is due to the classical \dcn into a product of \elr matrices
$$
\preskip.4em \postskip-.4em
\Cmatrix{2pt}{0&-a^{-1}\cr a&0} =\Cmatrix{2pt}{1&0\cr a&1} \,\Cmatrix{2pt}{1&-a^{-1}\cr 0&1}
\,\Cmatrix{2pt}{1&0\cr a&1}.
$$
\vspace{-1.5em}
\end{proof}

\vspace{.1em}
\pagebreak

\subsubsection*{When the image of a \prn matrix is free}

If a \prr $P\in\GAn(\gA)$ has as its image a free module of rank $r$, its kernel is not systematically free, and the matrix is therefore not necessarily similar to the standard matrix~$\I_{r,n}$.

It is interesting to find a simple characterization of the fact that the image is~free.

\begin{proposition}
\label{propImProjLib} \emph{(\Mprns whose image is free)}
\\
Let $P\in\Mn(\gA)$. The matrix $P$ is \idme and its image is free of rank~$r$  \ssi there exist two matrices $X\in\Ae {n\times r}$ and $Y\in\Ae {r\times n}$ such that~$YX=\I_r$ and~$P=XY$. In addition we have the following.
\begin{enumerate}
\item  $\Ker P=\Ker Y$, $\Im P=\Im X\simeq \Im Y$, and the columns of $X$ form a basis of~$\Im P$.
\item  For every matrices $X'$, $Y'$ of the same respective formats as $X$ and $Y$, and such that~$P=X'Y'$, there exists a unique matrix~$U\in\GL_r(\gA)$ such that 
$$\preskip-.2em \postskip.3em 
X'=X\,U\,\hbox{ and }\,Y=U\,Y'. 
$$
In fact,
$U=YX'$, $U^{-1}=Y'X$ and $Y'X'=\I_r$.
\end{enumerate}
\end{proposition}
\begin{proof}
Suppose that $P$ is \idme with a free image of rank $r$.
For columns of $X$ we take a basis of $\Im P$. 
Then, there exists a unique matrix $Y$ such that $P=XY$. 
Since $PX=X$ (because $\Im X\subseteq\Im P$ and $P^2=P$), we obtain~$XYX=X$.
Since the columns of $X$ are independent and $X(\I_r-YX)=0$, we obtain~$\I_r=YX$.
\\
Conversely, suppose $YX=\I_r$ and $P=XY$. 
Then

\snic{P^2=XYXY=X\I_rY=XY=P\;$ and $\;PX=XYX=X.}

Therefore $\Im P=\Im X$. In addition, the columns of $X$ are independent because~$XZ=0$ implies $Z=YXZ=0$.
 
\emph{1.} The sequence $\Ae n\vvvers{\In-P}\Ae n\vers{Y}\Ae r$ is exact. Indeed, $Y(\In-P)=0$, and if~$YZ=0$, then $PZ=0$, so $Z=(\In-P)Z$. Thus
\[\preskip.3em \postskip.4em 
\begin{array}{ccc} 
\Ker Y= \Im(\In-P)=\Ker P,\hbox{ and } \\[.3em] 
\Im Y\simeq
\Ae n\sur{\Ker Y}=\Ae n\sur{\Ker P\simeq \Im P}.  
 \end{array}
\]

\emph{2.} Now if $X'$, $Y'$ are of the same respective formats as $X$, $Y$, and if $P=X'Y'$,
we let $U=YX'$ and $V=Y'X$.
Then 
\begin{itemize}
\item $UV=YX'Y'X=YPX=YX=\I_r$,
\item  $X'V=X'Y'X=PX=X$, therefore $X'=XU$,
\item  $UY'=YX'Y'=YP=Y$, therefore $Y'=VY$.
\end{itemize}
Finally, $Y'X'=VYXU=VU=\I_r$.
\end{proof}

\vspace{-.7em}
\pagebreak	

\section[On \zed \risz]{\Mptfs over \zed \risz}
\label{secMPTFzed}

The following \tho \gnss \thref{thZerDimRedLib}.
\perso{les \algbs v\'erifient les points 2,3,4, mais les
\alos ne v\'erifient pas en \gnl les points 5 and 6}

\begin{theorem}
\label{propZerdimLib}
Let $\gA$ be a \zed \riz.
\begin{enumerate}
\item \label{ite1propZerdimLib} If $\gA$ is reduced every \mpf $M$ is \qfz, and every \smtf of $M$ is a direct summand

\item  \label{ite2propZerdimLib} \emph{(\Zed freeness lemma)}\\
Every \ptf \Amo is \qfz.

\item    \label{ite3propZerdimLib}
Every matrix $G\in\gA^{q\times m}$ of rank $\geq k$ is \eqve to a matrix
$$\preskip.2em \postskip.2em 
\Cmatrix{2pt}{
    \I_{k}   &0_{k,m-k}      \cr
    0_{q-k,k}&  G_1      } 
$$
with $\cD_r(G_1)=\cD_{k+r}(G)$ for all $r\geq 0$. In particular, every matrix of rank $k$ is \nlz.

\item    \label{ite4propZerdimLib}
Every \mpf $M$ such that $\cF_r(M)=\gen{1}$
(\cad \lot generated by $r$ \eltsz, cf.\ \Dfnz~\ref{deflocgenk})
is generated by $r$ \eltsz.

\item \label{ite5propZerdimLib}
\emph{(Incomplete basis \thoz)}\\
If a submodule $P$ of a \mptf $Q$ is \ptfz,
it has a \supl submodule.
If $Q$ is free of rank $q$ and~$P$ free of rank $p$, every \supl subspace is free of rank $q-p$.

\item  \label{ite6propZerdimLib}
Let $Q$ be a \ptf \Amo and $\varphi :Q\to Q$ an \endoz. \Propeq
\begin{enumerate}\itemsep0pt
\item $\varphi$ is injective.
\item $\varphi$ is surjective.
\item $\varphi$ is an \isoz.
\end{enumerate}
\end{enumerate}
\end{theorem}

\begin{proof}
Item \emph{\ref{ite1propZerdimLib}} is a reminder of \thref{thZerDimRedLib}.

\emph{\ref{ite2propZerdimLib}.}
We consider a \mpn $A$ of the module and we start by noting that since the module is \proz, $\cD_1(A)=\gen{e}$
with $e$ \idmz. 
We may assume that the first step of the computation is performed
 at the level of the \ri $\gA_{e}=\gA[1/e]$.
We are reduced to the case where $\cD_1(A)=e=1$, which we assume henceforth.
We apply item~\emph{\ref{ite3propZerdimLib}} with $k=1$ and conclude by \recuz.

Item~\emph{\ref{ite3propZerdimLib}} resembles an  \iv minor lemma (\ref{lem.min.inv}) without an \iv minor in the hypothesis.  
We apply with the \riz~$\Ared$ item~\emph{1} of \thref{thZerDimRedLib}.
We then obtain the desired matrix, but only modulo $\DA(0)$.
\\
We notice that the matrix $\rI_k+R$ with $R\in \Mk \big(\DA(0)\big)$ has an \iv \deterz, which allows us to apply the \iv minor lemma.

\emph{\ref{ite4propZerdimLib}.} Results from item~\emph{\ref{ite3propZerdimLib}} applied to a \mpn of the module.

\emph{\ref{ite5propZerdimLib}.}
Let us first look at the second case. Consider the matrix $G$ whose column vectors form a basis for the submodule $P$. Since $G$ is the matrix of an injective \aliz, its \idd of order~$p$ is \ndzz, therefore equal to $\gen{1}$ (Corollary~\ref{corZedReg}).
It remains to apply item~\emph{\ref{ite3propZerdimLib}}.
In the \gnl case, if $P$ is generated by $p$ \eltsz, let us consider a~$P'$ such that $P\oplus P'\simeq\gA^p$.
The module $Q\oplus P'$ is \ptfz, therefore is a direct summand in a module  $L\simeq \Ae n$. 
Then,
by the second case,~$P\oplus P'$ is a direct summand in $L$.
We deduce that $P$ is the image of a \prn {\mathrigid 2mu $\pi:L\to L$}. Finally, the restriction of $\pi$ to $Q$ is a \prn whose image is~$P$.

\emph{\ref{ite6propZerdimLib}}. 
We already know that \emph{b} and \emph{c} are \eqvs because $Q$ is \tf (\thref{prop quot non iso}). 
To prove that \emph{a} implies \emph{b}, we can assume that $Q$ is free (even if that means considering $Q'$ such that $Q\oplus Q'$ is free). Then,~$\varphi$ is represented by a matrix whose \deter is \ndz therefore \ivz.
\end{proof}

The previous \tho admits an important corollary in number theory.

\begin{corollary} \label{corpropZerdimLib} \emph{(One and a half \thoz)}
\begin{enumerate}
\item Let $\fa$ be an \id of $\gA$.
Assume that it is a \pf \Amo with $\cF_1(\fa)=\gen{1}$ and that there exists an $a\in\fa$ such that the \ri $\gB=\aqo{\gA}{a}$ is \zedz. Then, there exists a~\hbox{$c\in\fa$} such that~\hbox{$\fa=\gen{a,c}=\gen{a^m,c}$} for all $m\geq1$.
\item  Let $\gZ$ be the \ri of integers of a number field~$\gK$ and~$\fa$ be a nonzero \itf of $\gZ$. For all $a\neq0$ in $\fa$ there exists a~$c\in\fa$ such \hbox{that $\fa=\gen{a,c}=\gen{a^m,c}$} for all $m\geq1$.%
\index{one and a half!\Tho ---}
\end{enumerate}
\end{corollary}
\begin{proof}
\emph{1.} 
The \Bmo $\fa/a\fa$ is obtained from the \Amo $\fa$ by \eds from $\gA$ to $\gB$,
so its first \idf remains equal to $\gen{1}$. We apply item~\emph{\ref{ite4propZerdimLib}} of \thref{propZerdimLib}: there exists some $c\in\fa$ such that $\fa/a\fa=\gen{c}$ as a \Bmoz. This means that $\fa=c\gA+a\fa$ and gives the desired result.

 \emph{2.} If $\fa=\gen{\xn}$ is a \itf of $\gZ$, there exists a \itf $\fb$ such that $\fa\fb=\gen{a}$ (\thref{th1IdZalpha}).
\\
Let $\ux=[\,x_1 \;\cdots\;x_n \,]$. Therefore there exist $y_1$, \dots, $y_n$ in $\fb$ such that $\ux\tra{\uy}=\sum_i x_iy_i=a$. If $y_ix_j=\alpha_{ij}a$, we~have~$\alpha_{ii}x_k=\alpha_{ki}x_i$. Therefore, the \id $\fa$ becomes principal in $\gZ[1/\alpha_{ii}]$, equal to~$\gen{x_i}$, 
which is free of rank~$1$ (we can assume that the $x_i$'s are nonzero). \\
Since~$\sum_i \alpha_{ii}=1$, the $\alpha_{ii}$'s are \comz, therefore $\fa$ is \ptf and~$\cF_1(\fa)=\gen{1}$ (this is true \lot and therefore globally).\\
To apply item~\emph{1} it remains to verify that $\aqo{\gZ}{a}$ is \zedz.
The \elt $a$ annihilates a \polu $P\in\ZZ[X]$ of nonzero constant \coez, which we write as $aQ(a)=r\neq 0$.
Therefore, $\aqo{\gZ}{a}$ is a quotient \ri 
of $\gC=\aqo{\gZ}{r}$. It suffices to show that $\gC$ is \zedz. 
Let~\hbox{$\gA=\aqo{\ZZ}{r}$}.
Let $\ov u \in\gC$. Since $u$ annihilates a \polu $R\in\ZZ[T]$ of degree $n$, the \ri $\gA[\ov u]$ is a quotient \ri of the \ri $\aqo{\gA[T]} {\ov R(T)}$, 
which is a free \Amo of rank $n$, and so is finite. Therefore we can explicitly find $k\geq 0$ and $\ell\geq 1$ such that ${\ov u}^{k}(1-{\ov u}^{\ell})=0$.  
\end{proof}

\rem The matrix $A=(\alpha_{ij})$ satisfies the following \egts
$$\preskip.4em \postskip.4em
\tra{\uy}\,\ux=aA,\;A^2=A,\;\cD_2(A)=0,\;\Tr(A)=1,\;\ux A=\ux.
$$
We deduce that $A$ is a \mprn of rank $1$.\\ 
Moreover, we~have~$\ux(\In-A)=0$, and if $\ux \tra\uz=0$, then $\tra{\uy}\,\ux \tra\uz=0=aA\tra\uz$, so $A\tra\uz=0$ and~$\tra\uz=(\In-A)\tra\uz$. This shows that $\In-A$ is a \mpn of~$\fa$ (over the \sgrz~$(\xn)$). Therefore,~$\fa$ is \isoc as a~$\gZ$-module to~$\Im A\simeq\Coker(\In-A)$.
\eoe

\section{Stably free modules} 
\label{ModStabLibre}
\index{module!stably free ---}
\index{stably free!module}

Recall that a module $M$ is said to be \emph{\stlz} if it 
\emph{is a 
direct complement of a free module in a free module}, in other words if there exists an \iso between~$\Ae n$ and~$M\oplus \Ae r$ for two integers $r$ and $n$.\\ 
We will then say that $M$ is \emph{of rank} $s=n-r$.\footnote{This notion of rank will be \gneez, \Dfnsz~\ref{def ptf rank constant} and~\ref{defiRang}, and the reader will be able to note that those are indeed \gnnsz.}
The rank of a \stl module over a nontrivial \ri is well-defined. Indeed, if $M\oplus \Ae r\simeq\Ae {n}$ and~$M\oplus \Ae {r'}\simeq\Ae {n'}$, then we have $\Ae r\oplus \Ae {n'} \simeq \Ae {r'}\oplus \Ae n$ by Shanuel's lemma~\ref{corlemScha}.
From an \iso $M\oplus \Ae r\to\Ae n$, we obtain the \prn $\pi:\Ae n\to\Ae n$ over $\Ae r$ \paralm to $M$. This also gives a surjective \Ali $\varphi:\Ae n\to\Ae r$ with $\Ker\pi=\Ker\varphi\simeq M$: it suffices to let $\varphi(x)=\pi(x)$ for all $x\in \Ae n$.

Conversely, if we have a surjective \ali $\varphi : \Ae {n}\to\Ae r$, there exists a $\psi : \Ae r\to\Ae n$ such that $\varphi \circ \psi =\Id_{\Ae r}$.
Then $\pi=\psi\circ \varphi:\Ae {n}\to\Ae n$ is a \prnz, with $\Ker\pi=\Ker\varphi$, $\Im\pi=\Im\psi$ and $\Ker\pi\oplus \Im\pi =\Ae n$, and since $\Im\pi\simeq\Im\varphi=\Ae r$, the module

\snic{M=\Ker\varphi=\Ker\pi\simeq\Coker\pi=\Coker\psi}

is \stlz, and \isoc to $\Im(\Id_{\Ae n}-\pi)$.
Recall that by  
\thref{prop inj surj det},
saying that $\varphi:\Ae n \to \Ae r$ is surjective amounts to saying that~$\varphi$ is of rank $r$,
\cad that $\cD_r(\varphi)=\gen{1}$ in this case.

Finally, if we start from an injective \ali $\psi : \Ae r\to\Ae n$, saying that there exists a  $\varphi : \Ae {n}\to\Ae r$ such that
$\varphi \circ \psi =\Id_{\Ae r}$ amounts to saying that $\cD_r(\psi)=\gen{1}$ (\thref{propIGCram}).
Let us summarize the previous discussion.

\begin{fact}\label{factStablib} For a module $M$ \propeq
\begin{enumerate}
\item $M$ is \stlz.
\item $M$ is \isoc to the kernel of a surjective matrix.
\item $M$ is \isoc the cokernel of an injective matrix of maximum rank.
\perso{j'ai \'ecrit un truc plus lisible que \gui{injective \lnlz}
voir exo \ref{exoMatInjLocSimple}}
\end{enumerate}
\end{fact}

This result can allow us to define a new encoding, specific to \stl modules. Such a module will be encoded by the matrices of the \alis $\varphi$ and $\psi$. As for the dual of $M$ it will be encoded by the transposed matrices, as indicated in the following fact.

\begin{fact}\label{factStablibDual}
Using the previous notations, $M\sta$ is \stlz, canonically \isoc to $\Coker\!\tra{\varphi} $ and to $\Ker\!\tra{\psi}.$
\end{fact}

This is a special case of the following more \gnl result (see also Fact~\ref{factDualReflexif}).

\begin{proposition}\label{propDualSurjScind}
Let $\varphi:E\to F$ be a split surjection and $\psi:F\to E$ be a section of $\varphi$.
Let $\pi:E\to E $ be the \prn $\psi\circ \varphi$, and $j:\Ker\varphi\to E$ be the canonical injection.
\begin{enumerate}
\item $E=\Im\psi\oplus\Ker\varphi$, $\Ker\varphi=\Ker\pi\simeq\Coker\pi=\Coker\psi$.
\item $\Ker\!\tra{j}=\Im\!\tra{\varphi}$ and $\tra{j}$ is surjective, which by \fcn gives a canonical \iso $\Coker \tra{\varphi} \simarrow (\Ker\varphi)\sta$.
\end{enumerate}
\end{proposition}
\begin{proof}
The \ali $\psi$ is a \ing of $\varphi$ (\Dfn~\ref{defIng}).
We therefore have $E=\Im\psi\oplus\Ker\varphi$, and $\psi$ and $\varphi$ define reciprocal \isos between $F$ and $\Im \psi$. The proposition easily follows (see Fact~\ref{factIng0}).
\end{proof}
%

\subsec{When is a \stl module free?}

We then obtain the following results, formulated in terms of the kernel of a surjective matrix.

\begin{proposition}\label{propStabliblib}\emph{(When a \stl module is free, 1)}\\
Let $n=r+s$ and $R\in\Ae {r\times n}$. \Propeq

\begin{enumerate}
\item $R$ is surjective and the kernel of $R$ is free.
\item There exists a matrix $S\in\Ae {s\times n}$ such that the matrix
\smash{$\Cmatrix{2pt}{S\cr R}$} is \ivz.
\end{enumerate}
In particular, every \stl  module of rank $1$ is free.

\end{proposition}
\begin{proof}
\emph{1 $\Rightarrow$ 2.} If $R$ is surjective, there exists an $R'\in\Ae {n\times r}$ with $RR'=\I_r$.
\\
The matrices $R$ and $R'$ correspond to the \alis $\varphi$ and $\psi$ in the preliminary discussion. In particular, we have $\Ae n=\Ker R\oplus \Im R'$. Consider a matrix $S'$ whose column vectors constitute a basis of the kernel of $R$. Since $\Ae n=\Ker R\oplus \Im R'$, the matrix $A'=\lst{S'\mid  R'}$ has as its columns a basis of $\Ae n$.
It is \iv and its \inv is of the form \smashbot{$\Cmatrix{2pt}{S\cr R}$} because $R$ is the only matrix that satisfies $R\,A'=\lst{0_{r,n-r}\mid \I_r}$.
 
\emph{2 $\Rightarrow$ 1.} Let {$A=\Cmatrix{2pt}{S\cr R}$} and let $A'=A^{-1}$, which we write in the form $\lst{S'\,|\, R'}$. We have $RS'=0_{r,n-r}$, therefore 
$$\preskip.4em \postskip.2em 
\Im S'\subseteq \Ker R\eqno(\alpha), 
$$
and $RR'=\I_r$.
Therefore 
$$\preskip.2em \postskip.4em 
\qquad\Ker R\oplus \Im R'=\Ae n=\Im S'\oplus \Im R'\eqno (\beta). 
$$
Finally, $(\alpha)$ and $(\beta)$ imply $\Im S'= \Ker R$.

If $M$ is a \stl module of rank $1$, it is the kernel of a surjective matrix $R\in\Ae{(n-1)\times n}$. Since the matrix is surjective,
we \hbox{obtain $1\in\cD_{n-1}(R)$}, and this gives the row $S$ to complete $R$ as an \iv matrix (develop the \deter according to the first row).
\end{proof}

\begin{corollary}\label{prop2Stabliblib}\emph{(When a \stl module is free, 2)}\\
Consider $R\in\Ae {r\times n}$ and $R'\in\Ae {n\times r}$ with $RR'=\I_r$, $s:=n-r$. Then, the modules $\Ker R$ and $\Coker R'$ are \isoc and \propeq
\begin{enumerate}
\item The kernel of $R$ is free.
\item There exists a matrix $S'\in\Ae {s\times n}$ such that $\lst{S'\mid R'}$ is \ivz.
\item There exist a matrix $S'\in\Ae {s\times n}$ and a matrix $S\in\Ae {s\times n}$ such that
$$
\blocs{1.4}{0}{.6}{.8}{$S$}{}{$R$}{}\;\;
\blocs{.6}{.8}{1.4}{0}{$S'$}{$R'$}{}{}
=\In.
$$
\end{enumerate}
\end{corollary}

\vspace{.2em} 
Recall that a vector $x\in\gA^q$ is said to be \emph{\umdz} when its \coos are \ecoz. It is said to be \emph{completable}
if it is the first vector (row or column) of an \iv matrix.
\index{unimodular!vector}
\index{completable!vector}

\begin{proposition}\label{corpropStabliblib}
\Propeq \perso{peut on l'am\'eliorer with $=m$ and $q=m+1$?}
\begin{enumerate}
\item Every \stl \Amo of rank $\geq m$ is free.
\item Every \umd vector in $\gA^{{q}\times 1}$ with $q>m$ is completable.
\item Every \umd vector in $\gA^q$ with $q>m$ generates the direct complement of a free module in~$\gA^q$.
\end{enumerate}
\end{proposition}
\begin{proof}
Items \emph{2} and \emph{3} are clearly \eqvsz.

\emph{1} $\Rightarrow$ \emph{3.} Let $x \in \gA^q$ be a \umd vector with $q > m$.
Then, we can write $\gA^q = M \oplus \gA x$, and $M$ is \stl of rank $q-1 \ge m$, therefore free.

\emph{3} $\Rightarrow$ \emph{1.} Let $M$ be a \stl \Amo of rank $n \ge m$. We can write $L=M\oplus \gA x_1\oplus\cdots\oplus\gA x_r$, where $L\simeq \Ae{n+r}$.
If $r=0$, there is nothing left to do.
Otherwise, $x_r$ is a \umd vector in $L$, therefore by hypothesis~$\gA x_r$ admits a free \supl subspace in $L$. Thus, $L\sur{\gA x_r}\simeq \Ae{n+r-1}$, and similarly with $M\oplus \gA x_1\oplus\cdots\oplus\gA x_{r-1}$, which is \isoc to $L\sur{\gA x_r}$.
We can therefore conclude by \recu on $r$ that $M$ is free.
\end{proof}

\vspace{-.7em}
\pagebreak	

\subsec{Bass' stable range}

The notion of a stable range is linked to the \elr manipulations (of rows or columns) and allows us to some extent to control the \stl modules.

\begin{definition}\label{defiStableRange}
Let $n\geq 0$. A \ri $\gA$ is said to be of \emph{Bass'} \ix{stable range} \emph{less than or equal to} $n$ when we can \gui{shorten} the \vmds of length $n+1$ in the following sense

\snic{
\hbox { $1 \in \gen {a,\an}\;\Longrightarrow\;\exists \,\xn,\;1 \in \gen {a_1 + x_1a, \ldots, a_n + x_na}$.}
}

In this case we write \gui{$\Bdim \gA < n$.}  
\end{definition}

In the acronym $\Bdim$, $\mathsf{B}$ alludes to \gui{Bass.}

The notation $\Bdim \gA < n$ is legitimized on the one hand by item~\emph{1} in the following fact, and on the other hand by results to come which compare  $\Bdim$ to natural dimensions in commutative \algz.\footnote{See for example the results in Chapter~\ref{chapNbGtrs} which establish a comparison with the Krull and Heitmann dimensions.}
 
Item \emph{3.} uses the \id $\Rad \gA $ which will be defined in Chapter~\ref{chap Anneaux locaux}. The thing to note is that an \elt of $\gA$ is \iv \ssi it is \iv modulo $\Rad\gA$.

\begin{fact}\label{factStableRange} Let $\gA$ be a \ri and $\fa$ be an \idz.
\begin{enumerate}
\item If $\Bdim\gA< n$ and $n<m$ then $\Bdim\gA< m$.  
\item For all $n\geq 0$, we have $\Bdim\gA< n\Rightarrow\Bdim\gA/\fa < n$.
Abbreviated, we write this implication in the form: $\Bdim\gA/\fa \le \Bdim\gA$.
\item We have $\Bdim(\gA/\!\Rad \gA) =\Bdim\gA$ (by using the same abbreviation).
\end{enumerate}
\end{fact}

\begin{proof} \emph{1.} We take $m=n+1$.
Let $(a, a_0,\ldots, a_n)$ with $1 \in \gen {a, a_0, \ldots, a_n}$.  \\
We have $1 = ua + va_0 + \dots $, so $1 \in \gen {a', a_1,\ldots, a_n}$ with $a' = ua + va_0$. \\
Therefore we have $x_1$, \dots, $x_n$ in $\gA$  with $1 \in \gen {a_1 + x_1a', \ldots, a_n + x_na'}$, and consequently ${1 \in \gen {a_0 + y_0a, \dots, a_n + y_na}}$ with $y_0 = 0$ and $y_i = x_iu$ \hbox{for $i \ge 1$}.

\emph{2} and \emph{3.} Left to the reader.
\end{proof}
%

\begin{fact}
\label{corBass} \emph{(Unimodular vectors and \trelsz)}\\
Let $n\geq 0$. If $\Bdim \gA <  n$ and $V\in\Ae{n+1}$ is \umdz, it can be transformed into the vector $(1,0\ldots ,0)$ by \mlrsz.
\end{fact}
\begin{proof}
Let $V=(v_0,v_1,\ldots ,v_{n})$, with $1\in\gen{v_0,v_1,\ldots ,v_{n}}$.  Applying the \dfn with $a=v_0$, we obtain $x_1$, \dots, $x_n$ such that
$$1\in  \gen{v_1+x_1v_0,\alb\ldots ,\alb v_n+x_nv_0}.$$
The vector $V$ can be transformed by \mlrs into the vector $V'=(v_0,v_1+x_1v_0,\ldots ,v_n+x_nv_0)=(v_0,\alb v'_1,\alb \ldots ,\alb v'_{n})$, and we have $y_i$'s such that $\sum_{i=1}^ny_iv'_i=1$.
By \mlrsz, we can transform~$V'$ into $(1,v'_1,\ldots ,v'_{n})$, and then into $(1,0,\ldots ,0)$.
\end{proof}

Proposition~\ref{corpropStabliblib} and Fact~\ref{corBass} give the following \gui{Bass' \thoz.}
Actually, the real Bass' \tho is rather the conjunction of the following \tho with a \tho that provides a sufficient condition to have $\Bdim\gA<n$.
We will present several different variants in \thrfs{Bass0}{Bass} and Fact~\ref{factGdimBdim}.

\begin{theorem}
\label{corBass2} \emph{(Bass' \thoz, \stl modules)}\\
If $\,\Bdim\gA<n$, every \stl \Amo of rank $\geq n$ is free.
\end{theorem}

\section{Natural constructions}\label{secPtfCoNat}

\begin{proposition}
\label{propPtfExt} \emph{(Changing the base \riz)}\\
If $P$ is a \ptf \Amo and if $\rho:\gA\rightarrow \gB$ is a \ri \homoz, then the $\gB$-module $\rho\ist(P)$ obtained by \eds to $\gB$
is \ptfz. If $P$ is \isoc to the image of a \mprn $F=(f_{i,j})$, $\rho\ist(P)$ is \isoc to the image of {\em the same matrix seen in $\gB$}, \cad the \mprn $F^\rho= \big(\rho(f_{i,j})\big)$.
\end{proposition}
\begin{proof}
Changing the base \ri preserves the direct sums and the \prnsz.
\end{proof}

In the following proposition, we can a priori take as the sets of indices $I=\lrbm$ and $J=\lrbn$, but the set $I\times J$, which serves as a set of indices for the square matrix that defines the \ix{Kronecker product} of the two matrices $F$ and $G$ is not equal to $\lrb{1.. mn}$. This is an important argument in favor of the \dfn of  matrices \`a la Bourbaki, 
\cad  with finite row and column index sets which are not \ncrt of the type $\lrbm$.

\begin{proposition}
\label{propTensptf} \emph{(Tensor product)}\\
If $P$ and $Q$ are \pro modules represented by the \mprns $F=(p_{i,j})_{i,j \in I} \in \gA^{I\times I}$ and $G=(q_{k,\ell})_{k,\ell\in J}\in \gA^{J\times J}$, then the tensor product $P\te Q$ is a \mptf represented by the Kronecker product
$$\preskip.0em \postskip.3em 
F\te G=(r_{(i,k),(j,\ell)})_{(i,k),(j,\ell)\in I\times J}, 
$$
where $r_{(i,k),(j,\ell)}=p_{i,j}q_{k,\ell}$.
\end{proposition}
\begin{proof}
Suppose $P\oplus P'=\Ae m$ and $Q\oplus Q'=\Ae n$. The matrix $F$ (resp.\,$G$) represents the \prn over $P$ (resp.\,$Q$) \paralm to $P'$ (resp.\,$Q'$).
Then, the Kronecker product matrix $F\otimes G$ represents the \prn of $\Ae m\otimes \Ae n$ over $P\otimes Q$, \paralm to the subspace $(P'\otimes Q)\oplus (P\otimes Q')\oplus (P'\otimes Q')$.
\end{proof}

\vspace{-.7em}
\pagebreak	

\begin{proposition}
\label{propDual} \emph{(Dual of a \mptfz, 2)}\\
If $P$ is represented by the \mprn $F=(p_{i,j})_{i,j\in I}\in \gA^{I\times I}$, then the dual of $P$ is a \mptf represented by the transposed matrix of~$F$. If $x$ is a column vector in $\Im F$ and $\alpha$ a column vector in the image of $\tra{F}$, the scalar $\alpha(x)$ is the unique \coe of the matrix $\tra{\alpha}\,x$.
\end{proposition}
\begin{proof}
This results from Fact~\ref{factDualPTF}.
\end{proof}

\begin{proposition}
\label{propAliPtfs} \emph{(Modules of \alisz)}
\begin{enumerate}
\item If $P$ or $Q$ is \ptfz, the natural \homo (\paref{NOTAthetaMN})  
$$\preskip-.2em \postskip.2em 
\theta_{P,Q}:P\sta\te Q\to\Lin_\gA(P,Q) 
$$
is an \isoz.
\item If $P$ and $Q$ are \ptfsz, the module $\Lin_\gA(P,Q)$ is a \mptf canonically \isoc to $P\sta\te Q$, represented by the matrix $\tra{F}\te G$.
\item An \Amo $P$ is \ptf \ssi the natural \homo $\theta_{P}$ is an \isoz.
\end{enumerate}
\end{proposition}
\begin{proof}
\emph{1.}
Suppose $P\oplus P'=\Ae m$.
We have \isos
$$\arraycolsep2pt\begin{array}{rcl}
\Lin_\gA(\Ae m,Q)& \simeq& \Lin_\gA(P,Q)\oplus \Lin_\gA(P',Q),     \\[1mm]
(\Ae m)\sta\te Q&    \simeq&  (P\oplus P')\sta \te Q \\
&    \simeq&  (P\sta\oplus (P')\sta) \te Q \\
&    \simeq&  \big(P\sta\te Q \big)\oplus \big((P')\sta\te Q \big).
\end{array}$$
These \isos are compatible with the natural \homos 
$$\arraycolsep2pt\begin{array}{rcl}
Q^m\simeq (\Ae m)\sta\te Q&  \longrightarrow &
\Lin_\gA(\Ae m,Q)\simeq Q^m,   \\
P\sta\te Q&  \longrightarrow &\Lin_\gA(P,Q),     \\
(P')\sta\te Q&  \longrightarrow &\Lin_\gA(P',Q).
\end{array}$$
As the first is an \isoz, so are the others.\\
The case where $Q$ is \ptf is treated analogously.

\emph{2.} Special case of item~\emph{1}. 

\emph{3.} Results from item~\emph{1} and from the fact that $P$ is \ptf if the image of $\theta_P$ contains $\Id_P$ (\thref{propdef ptf}~$(b3)$).
\end{proof}

By using the commutation of the \eds with the tensor product we then obtain the following corollary.
\begin{corollary}\label{corpropAliPtfs}
If $P$ or $Q$ is \ptf (over $\gA$), and if $\gA\vers{\rho}\gB$ is an \algz, the natural \homo  
 
\snic{\rho\ist \big(\Lin_\gA(P,Q)\big)\to
\Lin_\gB \big(\rho\ist(P),\rho\ist(Q)\big)}

is an \isoz.
\end{corollary}
\vspace{-.5em}
\pagebreak

\section{Local structure \thoz}
\label{secMPTFlocLib}

In this work, we give several \dems of the local structure \tho for \mptfsz. 
The shortest path to the solution of this question is that provided by Fitting ideals.
This is the object of this section.

There is a lightning method 
based a kind of magic formula given in Exercise~\ref{exo7.1}. This miracle solution is actually directly inspired by another approach to the \pbz, based on a \gui{dynamic reread} 
of the local freeness lemma (\paref{lelilo}).
This dynamic reread is explained on \paref{quasiglobaldynamique} in Section~\ref{secMachLoGlo}.

However, we consider a more enlightening approach is that based entirely on projection matrices and on the more structural explanations involving
 the systematic use of the \deter of the \endos of \mptfsz.
This will be done in Chapter~\ref{chap ptf1}.

\begin{theorem} 
\label{prop Fitt ptf 2} \label{prop Fitt ptf 1}
\emph{(Local structure and \idfs of a \mptfz, 1)}
\begin{enumerate}
\item A \pf \Amo $P$ is \ptf \ssi its \idfs are (generated by) idempotents.
\item More \prmt for the converse, suppose that a \pf \Amo $P$ has \idms \idfsz, and \hbox{that~$G\in\gA^{q\times n}$} is a \mpn of $P$, corresponding to a \sys of~$q$ \gtrsz. 
\\
Let $f_h$ be the \idm that generates $\cF_h(P)$, \hbox{and~$r_h:=f_h-f_{h-1}$}. 
\begin{enumerate}
\item $(r_0,\ldots,r_q)$ is a \sfioz.
\item Let $t_{h,j}$ be a minor of order $q-h$ of $G$, and $s_{h,j}:=t_{h,j}r_h$. Then, the~$\gA[1/{s_{h,j}}]$-module $P[1/{s_{h,j}}]$ is free of rank~$h$.
\item The \elts $s_{h,j}$ are \comz.
\item We have $r_k=1$ \ssi the matrix $G$ is of rank $q-k$.
\item The module $P$ is \ptfz.
\end{enumerate}
\item  In particular, a \mptf becomes free after \lon at a finite number of \ecoz.
\end{enumerate}
\end{theorem}
%
\begin{proof}
\Thref{prop pf ptf} tells us that the module $P$ presented by the matrix $G$ is \pro \ssi the matrix $G$ is \lnlz. We then apply the \carn of  \lnl matrices by their \idds given in \thref{theoremIFD}, as well as the precise description of the structure of the \lnl matrices given in this \tho (items~\emph{5} and~\emph{7} of the \thoz).

Note: item~\emph{3} can be obtained more directly by applying \thref{theoremIFD} to an \idme matrix (therefore \lnlz) whose image is \isoc to the module $P$. 
\end{proof}
Thus, the \mptfs are locally free, in the strong sense given in \thref{th.ptf.loc}.

\smallskip In Section~\ref{sec ptf loc lib} we give an alternative \dem of \thref{th.ptf.loc}, more intuitive and more enlightening than the one we just gave. In addition, the \eco that provide free \lons are fewer.

\medskip \rem  
\label{rem.test.projectif}
We can therefore test if a \mpf is \pro or not when we know how to test whether its \idfs are idempotents or not. This is possible if we know how to test the membership $x\in \gen{a_1,\ldots ,a_h}$ for every \sys $(x,a_1,\ldots,a_h)$ of \elts of $\gA$, \cad if the \ri is strongly discrete.
One can now compare with \cite{MRR} Chap. III Exercise 4 p. 96.
\eoe

\rdb
\subsection*{Annihilator of a \mptfz}\label{subsubsecAnnMptf}

\begin{lemma}\label{lemAnnMptf}
The annihilator of a \mptf $P$ is equal to its first \idf $\cF_0(P)$, generated by an \idmz.
\end{lemma}
\begin{proof}
We know that the \idfs are generated by \idmsz.
We also know that $\cF_{0}(P)\subseteq\Ann(P)$ (Lemma~\ref{fact.idf.ann}).\\
Let us look at the opposite inclusion. Fact~\ref{fact.transporteur} implies that the annihilator of a \mtf is well-behaved under \lonz, so for every \mo $S$, we have $\Ann_{\gA_S}(P_S)= \big(\Ann_\gA(P)\big)_S$. We know that the same holds for the \idfs of a \mpfz. 
Moreover, to prove an inclusion of \idsz, we can localize at some \ecoz. We therefore choose \eco that render the module $P$ free, in which case the result is obvious.
\end{proof}

The previous \dem 
is an example of the strength of the local structure \tho (item~\emph{3} of \thref{prop Fitt ptf 2}).
The following section describes another such example.

\section[\Lcy \pro modules]{\Lcy \pro modules and \ptf \idsz}
\label{secIdProj}\label{Idpp}

\vspace{4pt}
\subsec{\Lcy modules}
\label{sec mlm}

An \Amo $M$ is said to be \ixc{cyclic}{module} 
if it is generated by a single \eltz: $M=\gA{a}$. In other words, if it is \isoc to a quotient~$\gA/\fa$.


\ms In \clama a module is said to be {\em \lmoz} if it becomes cyclic
after \lon at any arbitrary \idepz.
It seems difficult to provide an \eqv statement that makes sense in \comaz.
Recall also that the remark \paref{remplcc.tf} shows that the notion does not seem pertinent when the module is not assumed to be \tfz. Nevertheless when we restrict ourselves to the \mtfs there is no issue.
The following \dfn has already been given before Fact~\ref{factExl1Plg}.

\pagebreak	

\begin{definition}
\label{defmlm} A \emph{\tfz} \Amo $M$ is said to be {\em \lmoz} if there exist \moco $S_1$, \ldots, $S_n$ of $\gA$ such that each
$M_{S_j}$ is cyclic 
as an $\gA_{S_j}$-module. In the case of an \id we speak of a \emph{\lop \idz}.%
\index{locally!cyclic module}\index{module!locally cyclic ---}%
\index{locally!principal ideal}\index{ideal!locally principal ---}
\end{definition}

Note that the \prt \gui{concrete local} in the previous \dfnz, without the hypothesis that $M$ is \tfz, implies that $M$ is \tf (\plgrf{plcc.tf}).

\ss We will need the following remark.
\begin{fact} \emph{(Successive \lons lemma, 1)}
\label{factLocCas}\index{Successive \lons lemma, 1}
 \\
If $s_1$, \ldots, $s_n$ are \eco of $\gA$ and if for each $i$,
we have \elts
$
s_{i,1},\; \ldots ,\;
s_{i,k_i},
$
 \com in $\gA[1/s_i]$,
then the \elts $s_{i}s_{i,j}$ are \com in~$\gA$.
\end{fact}

Item \emph{3} of the following theorem presents an efficient computational machinery for locally cyclic modules.

\begin{theorem}
\label{propmlm} {\em (\Lcy \tf modules)}\\
Let $M=\gA x_1+\cdots+\gA x_n$ be a \mtfz.
\Propeq
\begin{enumerate}
\item [1.\phantom{*}] The module $M$ is \lmoz.
\item [2.\phantom{*}] There exist $n$ \eco $s_i$ of $\gA$ such that for each $i$ we have $M=_{\gA_{s_i}}\gen{x_i}$.
\item [3.\phantom{*}] There exists a matrix $A = (a_{ij})\in\Mn(\gA)$ that satisfies
\begin{equation}\label{eqmlm}\preskip.3em \postskip.3em
\left\{\arraycolsep2pt
\begin{array}{rcl}
\sum a_{ii}      &=&  1\\[.2em]
a_{\ell j}x_{i} & =& a_{\ell i}x_{j} \qquad \forall i, j,
\ell \in \lrbn
\end{array}
\right.
\end{equation}
in other words, for each row $\ell$, the following matrix is formally of 
\hbox{rank $\leq 1$} (its minors of order 2 are null)
$$\preskip.3em \postskip-.4em
\Cmatrix{2pt}{
a_{\ell 1}&\cdots &a_{\ell n}\cr
x_1&\cdots &x_n
}
.$$
\item [4.\phantom{*}]$\Vi_\Ae 2(M)=0$.
\item [5.\phantom{*}] $\cF_1(M)=\gen{1}$.
\item [6*.] After \lon at any \idepz, $M$ is cyclic.
\item [7*.] After \lon at any \idemaz, $M$ is cyclic.
\end{enumerate}
\end{theorem}
\begin{proof}
\emph{3} $\Rightarrow$ \emph{2} $\Rightarrow$ \emph{1.} Clear, with $s_i=a_{ii}$ in item~\emph{2}.

Let us show that a cyclic module 
satisfies condition~\emph{3}.
\\
If $M = \gen{g}$, we~have~$g = \sum_{i=1}^{n} u_i x_i$ and  
$x_i = gy_i$. 
Let $b_{ij} = u_iy_j$.\\ 
Then, for all $i$, $j$, $\ell \in \lrbn$, we have $b_{\ell j} x_{i} = u_\ell  y_i y_j g= b_{\ell i} x_{j}$. 
In addition
$$\preskip.4em \postskip.4em\ndsp 
g = \som_{i=1}^{n} u_i x_i = \som_{i=1}^{n} u_i y_i g
= \big( \som_{i=1}^{n} b_{ii} \big) g. 
$$

Let $s=1-\som_{i=1}^{n} b_{ii}$.
We have $sg=0$, and so $sx_k=0$ for all $k$.
\\
Take $a_{ij} = b_{ij}$ for $(i, j) \neq (n,n)$ and $a_{nn} = b_{nn} + s $.
Then, the matrix~$(a_{ij})$ indeed satisfies Equations~(\ref{eqmlm}).

\emph{1} $\Rightarrow$  \emph{3.}
The \prt \emph{3} can be seen as the existence of a solution of a \sli whose \coes are expressed in terms of the \gtrsz~$x_i$.
However, a cyclic module  
satisfies the \prtz~\emph{3.}
We can therefore apply the basic \plgz.\iplg

Thus, \emph{1} $ \Leftrightarrow$ \emph{2} $\Leftrightarrow$ \emph{3.}

\emph{1} $\Rightarrow$  \emph{4} and \emph{1} $\Rightarrow$ 
\emph{5.} Because the functors 
$\Vi_\Ae 2\bullet$ and $\cF_1(\bullet)$ are well-behaved under \lonz.

\emph{5} $\Rightarrow$ \emph{1}. $M$ is the quotient module
of a \mpf $M'$ such that $\cF_1(M')=\gen{1}$. We can therefore suppose \spdg that~$M$ is \pf with a \mpn $B\in\Ae{n\times m}$.
By hypothesis, the minors of order $n-1$ of the matrix $B$ are \comz. 
When we invert one of these minors, by the \iv minor lemma (\paref{lem.min.inv}), the matrix $B$ is \eqve to a matrix

\snic{
\Cmatrix{2pt}{
   \I_{n-1}   &0_{n-1,m-n+1}      \cr
    0_{1,n-1}&     B_1},}

and the matrix $B_1\in\Ae {1\times (m-n+1)}$ is also a \mpn of~$M$.
\\
Assume \emph{4} and $n\geq 2$, and let us show that $M$ is, after \lon at suitable \ecoz, generated by $n-1$ \eltsz. This will be sufficient to show (by using an \recu on $n$) that \emph{4} implies \emph{1}, by using Fact~\ref{factLocCas}. The module $\Al2_\gA(M)$ is generated by the \elts $v_{j,k}=x_j\vi x_k$ ($1\leq j<k\leq n$) and the syzygies between the $v_{j,k}$'s are all obtained from the syzygies between the $x_i$'s. Therefore if $\Al2_\gA(M)=0$, $M$ is the quotient of a \mpf $M'$ such that $\Al2_\gA(M')=0$. We then suppose \spdg that $M$ is \pf with a \mpn $A=(a_{ij})$. A \mpn $B$ for $\Al2_\gA(M)$ with the \gtrs $v_{j,k}$ is obtained as indicated in Proposition~\ref{propPfPex}. It is a matrix of format ${{n(n-1)}\over{2}}\times m$ (for some suitable $m$), and each \coe of $B$ is null or equal to some $a_{ij}$. This matrix is surjective, therefore $\cD_{n(n-1)/2}(B)=\gen{1}$ and the $a_{ij}$'s are \comz. However, when we pass from $\gA$ to $\gA[1/a_{ij}]$, $x_i$ becomes \coli of the $x_k$'s ($k\neq i$) and $M$ is generated by $n-1$ \eltsz.

\emph{1} $\Rightarrow$ \emph{6*} $\Rightarrow$ \emph{7*}. Obvious.

The \dem that \emph{7*} implies~\emph{3} is non\covz: in the \dem that \emph{1} implies~\emph{3}, we replace the existence of a solution of a \sli under the basic \plg by the existence of a solution under the corresponding abstract \plgz.\iplg
\end{proof}


A matrix $(a_{ij})$ which satisfies Equations~(\ref{eqmlm}) is called a \emph{cyclic \lon matrix for the $n$-tuple $(\xn)$}. If the $x_i$'s are \elts of $\gA$, they generate a \lop \id and we speak of a \emph{principal localization matrix}.%
\index{matrix!principal \lon ---}%
\index{matrix!cyclic \lon --- for the $n$-tuple $(\xn)$}%
\index{localization!cyclic --- matrix}%
\index{localization!principal --- matrix}%

\medskip
\rem In the case of a module generated by $2$ \elts $M=\gA x+ \gA y$, Equations~(\ref{eqmlm}) are very simple and a \mlmo for $(x,y)$ is a matrix $\cmatrix{1-u&-b\cr-a&u}$ which satisfies
\begin{equation}\label{eqmlm2gen}\preskip.3em \postskip.4em
\dmatrix{1-u&-b\cr x&y}= \dmatrix{-a&u\cr x&y}=0,\,\mathrm{i.e.}\;(1-u)y=bx \;\hbox{ and }\;
ux=ay.
\end{equation}
\vspace{-6pt}
\eoe

\vspace{3pt}
\begin{proposition}\label{pmlm}
Let $M = \gA x_1+\cdots+\gA x_n $ be a \tf \Amoz.
\begin{enumerate}
\item 
If $M$ is \lmo and if $A= (a_{ij})$ is a \mlmo for $(x_1,\dots,x_n)$, we have the following results.
\begin{enumerate}
\item  
$\Cmatrix{2pt}{x_1 \;\cdots \; x_n}\;
A =
\Cmatrix{2pt}{x_1 \;\cdots \; x_n}
$.
\item  
    The \ids  $\cD_2(A)$ and $\cD_1(A^2 - A)$ annihilate $M$.
\item   
One has {\mathrigid 2mu $a_{ii}M\subseteq \gA x_i$, and over the \ri $\gA_i = \gA[\fraC1{a_{ii}}]$,  $M =_{\gA_i} 
\gA_i {x_i}$}.
\item   
$ \gen{a_{1j},\dots,a_{nj}} M = \gA{x_j}$.
\item   
More \gnltz, for any arbitrary \elt $y = \sum\alpha_{i}x_i$ of $M$, if we let $\alpha ={^t[\alpha_1\;\cdots\;\alpha_n}]$ and $\beta  = A\, \alpha$, then $y=\sum_i\beta_ix_i$ and we obtain an \egt of square matrices with \coes in~$M$:
\begin{equation}\label{eqpmlm}\preskip.4em \postskip.4em
\beta x = \Cmatrix{2pt}{\beta_1\cr \vdots\cr \beta_n}\Cmatrix{2pt}{x_1 \;\cdots \; x_n}=
A\,y,\;\;\mathit{i.e.}\quad \forall i,j\;\beta_ix_j=a_{ij}y
\end{equation}
In particular,
${\gen{\beta _{1},\dots,\beta _{n}} M = \gA y.}$
\end{enumerate}
\item 
\Propeq
\begin{enumerate}
\item [--] $M$  is \isoc to the image of a \mprn of rank $1$.
\item [--] $M$ is faithful (\cad $\Ann(M)=0$) and \lmoz.
\end{enumerate} 
In this case, let $A$ be a \mlmo for $(\xn)$. We obtain
\begin{enumerate}
\item [--] $A$ is a \mprn of rank $1$,
\item [--] the following sequence is exact
\smash{$
\Ae n\vvvvers{\In-A} {\Ae n}\vvvvvers{[\,x_1\;\cdots\;x_n\,]}M\to 0,
$}
\item [--]  $M\simeq\Im A$.
\end{enumerate}
\end{enumerate}
\end{proposition}

\begin{proof}
\emph{1.} Item~\emph{1c} is clear, and \emph{1d} is a special case of \emph{1e}.

\emph{1a.} The $j^{\rm th}$ \coo of $\Cmatrix{2pt}{x_1 \;\cdots \; x_n}\, A$ is written as
$${
\som_{i=1}^{n}a_{ij}x_i = \som_{i=1}^{n}a_{ii}x_j = x_j}.$$

\emph{1b.} Let us show that every minor of order $2$ of $A$ annihilates $x_i$. Consider the following matrix
$$
\Cmatrix{2pt}{
a_{ji}  &  a_{j\ell}  &  a_{jh}\cr
a_{ki}  &  a_{k\ell}  &  a_{kh}\cr
x_i     &  x_\ell     &  x_h
}.
$$
Its \deter is null (by expanding with respect to the first row) 
and the expansion with respect to the first column provides
$$
  (a_{j\ell} a_{kh} - a_{jh} a_{k\ell})x_i =0.
$$
Let us show that $A^2 = A$ modulo $\Ann(M)$. What follows is written modulo this annihilator. We come to show that the minors of order $2$ of $A$ are null.
Thus $A$ is a \mlmo for each of its \hbox{rows $L_i$}.
By item~\emph{1a} applied to $L_i$, we have $L_i\, A=L_i$, and \hbox{so $A^2=A$}.

\emph{1e.} Let $\ux=\Cmatrix{2pt}{x_1 \;\cdots \; x_n}$.
On the one hand 

\snic{\dsp\som_i\beta_ix_i=\ux\beta=\ux A\alpha=\ux\alpha=\som_i\alpha_ix_i.}

On the other hand,

\snic
{\dsp
\beta_i x_{j} = \som_{k}\alpha_{k}\, a_{ik} \,x_j=
\som_{k}\alpha_{k}\, a_{ij} \,x_k =
a_{ij}\,\big(\som_{k}\alpha_{k} x_k\big) = a_{ij}\, y.
}

This shows \Egrf{eqpmlm} and we deduce that $\gen{\beta _{1},\dots,\beta _{n}} M =  \gA y$.

\sni \emph{2.}
First assume that $M$ is \isoc to the image of a \mprn $A$ of rank $1$. Let $x_i$ be the $i^{\rm th}$ column of $A$. As~$\cD_2(A)=0$, we have the \egts $a_{\ell j}x_{i}  = a_{\ell i}x_{j}$ for $i, j, \ell \in \lrbn$. This implies that over the \ri $\gA[1/a_{\ell j}]$, $M$ is generated by $x_j$, and since $\cD_1(A)=\gen{1}$, the module is \lmoz. Finally, let~$b\in\Ann(M)$, \hbox{then~$bA=0$}, and~$\cD_1(A)=\gen{1}$ implies $b=0$; the module is faithful.
\\
Now assume that $M$ is \lmoz, and that $A$ is a \mlmo for a \sgr $(\xn)$.
If $M$ is faithful, given \emph{1b}, we have $\cD_2(A)=0$ and $A^2=A$ so $A$ is a \mprn of rank $\leq1$. Since $\Tr(A)=1$, $A$ is of rank $1$. Given~\emph{1a}, the matrix $\I_n-A$ is a matrix of syzygies 
for $(\xn)$. Now let $\sum_{i=1}^n\alpha_ix_i=0$ be an arbitrary syzygy of the $x_i$'s.
As in \emph{1e}, let

\snic{\beta=\tra{[\,\beta_1\;\cdots\;\beta_n\,]} = A\;
\tra{[\,\alpha_1\;\cdots\;\alpha_n\,]}, }

we obtain $\gen{\beta_{1},\dots,\beta_{n}} M =  0$ and, since $M$ is faithful, $\beta=0$.\\
Thus,
$\,A\,\tra{[\,\alpha_1\;\cdots\;\alpha_n\,]}=0$ and $(\In-A)\,\tra{[\,\alpha_1\;\cdots\;\alpha_n\,]}=\tra{[\,\alpha_1\;\cdots\;\alpha_n\,]}$:
every syzygy 
for $(\xn)$ is a \coli of the columns of $\I_n-A$. This shows that $\I_n-A$ is a \mpn of $M$ for the \sgr $(\xn)$.
Since $A^2=A$, we have $M\simeq\Coker(\I_n-A)\simeq\Im A$.
\end{proof}

\subsec{Cyclic \pro modules} 

The following description applies in particular to  \pro \idpsz.
\begin{lemma}
\label{lemIdpPtf}
For a cyclic module $M$, \propeq
\begin{enumerate}
\item $M$ is a \ptf \Amoz.
\item $\Ann(M)=\gen{s}$ with $s$ \idmz.
\item $M\simeq \gen{r}$ with $r$ \idmz.
\end{enumerate}
\end{lemma}
\begin{proof}
The implications \emph{2 $\Rightarrow$ 3
$\Rightarrow$ 1} are obvious, and the implication \emph{1 $\Rightarrow$ 2} is given in Lemma~\ref{lemAnnMptf}.
\end{proof}

We deduce that a \ri $\gA$ is \ix{quasi-integral} 
\ssi every \idp is \proz,
which justifies the English terminology \emph{pp-ring} (principal ideals are projective).%
\index{ring!quasi-integral --- (or pp-ring)}%
\index{pp-ring (or quasi-integral ring)}%

\subsec{\Lcy \pro modules}\label{subsubseclmoproj}

The following lemma \gnss the \eqvc given in Proposition~\ref{pmlm} 
between faithful locally cyclic modules and images of projection matrices of rank 1.

\begin{lemma}\label{lemLmoProj} 
\Propeq
\begin{enumerate}
\item $M$ is \lmo and $\Ann (M)$ is generated by an \idmz.
\item $M$ is \ptf and \lmoz.
\item $M$ is \isoc to the image of a \mprn of rank~$\leq1$.
\end{enumerate}
\end{lemma}
\begin{proof}
\emph{1} $\Rightarrow$ \emph{2.} We localize at \eco that render the module cyclic 
and we apply Lemma~\ref{lemIdpPtf}.
 
In \emph{2} and \emph{3} we let $F$ be a square \mprn of order $n$, with~$M$ as its image. After \lon at \eco it becomes similar to a standard \mprn $\I_{k,n}$, $k$ depending on the \lonz.
 
\emph{2} $\Rightarrow$ \emph{3.} If $k>1$,
we obtain at the corresponding \lon $\cF_1(M)=\gen{0}$.
As we have already $\cF_1(M)=\gen{1}$, the \lon is trivial. The rank of $F$ is therefore~$\leq1$ at all the \lonsz.
 
\emph{3} $\Rightarrow$ \emph{1.}  After \lonz, as the matrix is of rank $\leq1$, we have $k\leq1$. The module  therefore becomes cyclic. 
Moreover, by Lemma~\ref{lemAnnMptf}, $\Ann (M)$ is generated by an \idmz.
\end{proof}

\subsec{\Ptf \idsz}\label{subsubsecIdproj}
Recall that an \id $\fa$ is said to be \emph{faithful} if it is faithful as an \Amoz.%
\index{faithful!ideal}\index{ideal!faithful ---}

\smallskip 
\rem In the most common terminology, an \id is called \ndz if it contains a \ndz \eltz. A fortiori this is a faithful \idz. We will not use this terminology as we find it ambiguous. \eoe

\begin{lemma}\label{lemIdproj}~
\begin{enumerate}
\item \label{i1lemIdproj} If $\fa\subseteq\fb$ with $\fa$ \tf and $\fb$ \lopz, there exists a \itf $\fc$ such that $\fb\fc=\fa$.

\item \label{i2lemIdproj} An \id $\fa$ is \ptf \ssi it is \lop and its annihilator is generated by an \idmz.

\item \label{i20lemIdproj} An \id $\fa$ is \qf \ssi it is principal and its annihilator is generated by an \idmz.

\item \label{i3lemIdproj} Let $\fa_1$ and $\fa_2$ be \ids and $\fb$ be a faithful \ptf \idz.
\hbox{If $\fb\fa_1=\fb\fa_2$,} then $\fa_1=\fa_2$.

\item \label{i6lemIdproj}
An \id is \iv \ssi it is \lop and it contains a \ndz \eltz.
\end{enumerate}
\index{invertible!ideal}
\index{ideal!invertible ---}
\end{lemma}

\smallskip %
\begin{proof}
\emph{\ref{i1lemIdproj}.} It is enough to show that for an arbitrary $a\in\fb$ there exists a \itf $\fc$ such that $\fb\,\fc=\gen{a}$. 
This is given by item~\emph{1e} of Proposition~\ref{pmlm} when
$M=\fb$.

\emph{\ref{i2lemIdproj}.} The direct implication uses Corollary~\ref{corprop inj surj det}: if a \ali $\Ae k\to\gA$ is injective with $k>1$, the \ri is trivial. Therefore at each \lonz, the \id $\fa$ is not only free but principal.
The converse implication is in Lemma~\ref{lemLmoProj}.

\emph{\ref{i20lemIdproj}.} For the direct implication, we write $\fa\simeq \bigoplus_{i\in \lrbn} \gen{e_i}$, where the $e_i$'s are \idms with $e_{i+1}$ being a multiple of $e_i$. We want to show that if $n> 1$, $e_2=0$.
We localize the injection $\fa \to \gA$ at $e_2$ and we obtain an injection
$$
\preskip.0em \postskip.4em \ndsp
~~~~\gA_{e_2}\oplus \gA_{e_2}\simeq e_1\gA_{e_2}\oplus e_2\gA_{e_2}\hookrightarrow \bigoplus e_i\gA_{e_2}\simeq \fa\gA_{e_2} \hookrightarrow \gA_{e_2}, 
$$
so $\gA_{e_2}$ is null (Corollary~\ref{corprop inj surj det}).

\emph{\ref{i3lemIdproj}.} The \id $\fb$ becomes free (after \lonz), and cyclic
by item~\emph{\ref{i2lemIdproj}.} 
If in addition it is faithful, its annihilator is null, and the \gtr is a \ndz \eltz.

\emph{\ref{i6lemIdproj}.} Item \emph{\ref{i1lemIdproj}} implies that a \lop \id that contains a \ndz \elt is \ivz. Conversely, let $\fa=\gen{a_1,\ldots,a_n}$ be an \iv \idz. There exists a $c$ \ndz in $\fa$ and an \id $\fb$ such that $\fa\,\fb=\gen{c}$. Let $b_1$, \ldots, $b_n\in \fb$ with $\sum_ia_ib_i=c$. We have for each $i$, $j\in\lrbn$ \hbox{some $c_{ij}\in\gA$} such that $b_ia_j=c\,c_{ij}$. By using the fact that $c$ is \ndz we verify without difficulty that the matrix $(c_{ij})_{1\leq i,j\leq n}$ is a \mlp for $(a_1,\ldots,a_n)$.
\end{proof}
%


\section[Determinant, \polfon and \polmuz]{Determinant, \polcarz, 
\polfon and \polmuz}
\label{subsec det ptf}

\vspace{4pt}
If $M$ is an \Amoz, we denote by $M[X]$ the $\AX$-module obtained by \edsz.

When $\gA$ is an integral \riz, if $P$ is a \ptf module, \isoc to the image of a \prr $F\in\GAn(\gA)$, by \eds to the quotient field  we obtain a 
 vector space 
$P'$ of finite dimension, say $k$.
We deduce that the \polcar of the matrix $F$ is equal to $(X-1)^kX^{n-k}$. Even simpler, the \deter of the multiplication by $X$ in $P'[X]$ is equal to $X^k$, \cad

\snic{\det \big((\In -F)+XF \big)=X^k.}

\smallskip
When $\gA$ is an arbitrary \riz, we will see that we can define the analogue of the above \pol $X^k$.
First of all, we introduce the \deter of an \endo of a \mptfz.

\subsec{The \deterz, the \polcar and the cotransposed \endoz}

\begin{thdef} 
\label{propdef det ptf} Let $P$ be a \mptfz.
\begin{enumerate}
\item  Let $\varphi\in\End(P)$. Suppose that $P\oplus Q_1$ is \isoc to a free module and let $\varphi_1=\varphi\oplus \Id_{Q_1}$. 
\begin{enumerate}
\item  The \deter of $\varphi_1$ only depends on $\varphi$. 
The scalar defined as such is called the {\em \deter of the \endoz} $\varphi$%
\index{determinant!of an endomorphism}.
We denote it by $\det(\varphi)$ or $\det \varphi$.\label{NOTAdet}
\item  The \deter of the \endo $X\Id_{P[X]}-\varphi$ of $P[X]$ is called the {\em \polcar of the \endoz}~$\varphi$. \\
We denote it by $\rC\varphi(X)$; we have $\rC{-\varphi}(0)=\det\varphi$.\label{NOTAPolcar}%
\index{polynomial!characteristic --- of an endomorphism}
\item  Consider the cotransposed \endo $\Adj(\varphi_1)=\wi{\varphi_1}$ of $\varphi_1$. It operates on~$P$ and the \endo of $P$ defined as such only depends on $\varphi$. We call it {\em the cotransposed \endo of} $\varphi$%
\index{cotransposed!endomorphism}
 and we denote it by $\wi{\varphi}$ or $\Adj(\varphi)$.\label{NOTACotransp}
\item Let $\rho:\gA\to\gB$ be a morphism. By \eds from~$\gA$ to~$\gB$,
we get a \mptf $\rho\ist(P)$ with an \endo $\rho\ist(\varphi)$.
Then we have the following good \gui{functorial} \prts
 
\vspace{-.3em}
 \[ 
 \begin{array}{ccccccc} 
 \det\big(\rho\ist(\varphi)\big)    =  \rho\big(\det(\varphi)\big),~&&\rC{\rho\ist(\varphi)}(X)    =  \rho \big(\rC\varphi(X)\big),  \\[.3em] 
  \Adj\big(\rho\ist(\varphi)\big)   =  \rho\ist\big(\Adj(\varphi)\big). 
  \end{array}
 \] 
\end{enumerate}
%
\item If $\psi : P\rightarrow P$ is another \endo of $P$, we~have
$$\preskip.4em \postskip.4em 
\det(\varphi\circ\psi)=\det(\varphi)\det(\psi). 
$$
\item  If $P'$ is another \mptf and if \smash{$\psi=\blocs{.7}{.8}{.7}{.8}{$\varphi$}{$\gamma$}{$0$}{$\varphi'$}$} is an \endo of $P\oplus P'$ \gui{block-triangular,}
we have \\
$
\det(\psi)=d d' \;\hbox{ and }\;
\wi{\psi}=\blocs{.7}{.8}{.7}{.8}{$d'\wi{\varphi}$}{$\eta$}{$0$}{$d\wi{\varphi'}$}
\,,\hbox{ where } d=\det(\varphi),\,d'=\det(\varphi') . 
$
\item  If $\varphi : P\rightarrow P$ and $\varphi' : P'\rightarrow P'$ are \endos of \mptfsz, and if $\alpha \circ \varphi =\varphi'\circ \alpha  $ for an \iso  $\alpha : P\rightarrow P'$, then $\det(\varphi )=\det(\varphi' )$.
\item  The \ali $\varphi$ is an \iso (resp.\,is injective) \ssi $\det(\varphi)$ is \iv (resp.\,is \ndzz).
\item  We have the \gui{classical} \egt
$$\preskip.2em \postskip.2em 
\wi{\varphi}\circ \varphi \,=\, \varphi\circ\wi{\varphi}
\,=\,\det(\varphi)\,\Id_P. 
$$
\item  The Cayley-Hamilton \tho applies: $\rC\varphi(\varphi)=0$.

\item 
Let
$$\preskip-.0em \postskip.4em 
\begin{array}{rcccl}
\Gamma_{\varphi}(X)& :=  &  -\frac{\rC{-\varphi}(-X)-\rC{-\varphi}(0)}{X} & = &
\frac{-\rC{-\varphi}(-X)+\det(\varphi)}{X}\,,
\end{array} 
$$
such that
$ \rC{-\varphi}(-X) =
-X\Gamma_{\varphi}(X)+\det(\varphi)$.
Then  $\widetilde{\varphi}=\Gamma_{\varphi}(\varphi)$.
\end{enumerate}
\end{thdef}

\begin{proof}  We remark that the \dfns given in item \emph{1} indeed reproduces the usual objects of the same name in the case where the module is free. Similarly the formula in \emph{8} gives, when $\varphi$ is an \endo of a free module,  the same~$\Gamma_\varphi$
as the formula of Lemma~\ref{lemPrincipeIdentitesAlgebriques}. So there is no conflict of notation. 

\emph{1a.} Assume that $\Ae m\simeq P\oplus Q_1$ and $\Ae n\simeq P\oplus Q_2$,  and consider the direct sum
$$\preskip-.5em \postskip.2em 
\Ae {m+n}\simeq P\,\oplus\, Q_1 \oplus P\,\oplus\, Q_2. \eqno(*)
$$
Let 
$\varphi_1=\varphi\,\oplus\, \Id_{Q_1}$ \hbox{and $\varphi_2=\varphi\,\oplus \,\Id_{Q_2}$}. One has to show the \egt  $\det\varphi_1=\det\varphi_2$.
Consider the \endo of $\Ae {m+n}$
$$\preskip.3em \postskip.3em 
\phi=\varphi\,\oplus\,\Id_{Q_1}\oplus\,\Id_{P}\oplus\,\Id_{Q_2}, 
$$ 
such that $\phi$ is conjugated of $\varphi_1\oplus \Id_{\Ae n}$ and of $\varphi_2\oplus \Id_{\Ae m}$. Hence~\hbox{$\det\phi=\det\varphi_1$} \hbox{and $\det\phi=\det\varphi_2$}.

\emph{1c.} We proceed similarly. 
The cotransposition of the \endos satisfies item~\emph{3} in the case of free modules, so $\wi\phi$ operates on $P\oplus Q_1$ and is restricted at $\wi{\varphi_1}$.
Moreover, since $\wi\phi=\Gamma_\phi(\phi)$, $\wi\phi$ operates on each component in the direct sum $(*)$. 
Similarly $\wi\phi$ operates on $P\oplus Q_2$ and is restricted at $\wi{\varphi_2}$.
Therefore $\wi\varphi_1$ and $\varphi_2$ both operate on $P$ in the same way that $\wi\psi$ does. Note that $\wi\varphi=\Gamma_\phi(\varphi)$.
 
\emph{1d.} 
This is a direct consequence of the \dfnsz.

\pagebreak	

All the remaining items of the \tho are consequences of the free case (where the results are clear), of the local structure \tho and of item \emph{1d}. Indeed the statements can be certified by verifying them after \lon at \ecoz,
and the \mptfs we consider become simultaneously free after \lon at a suitable \sys of \ecoz. 
Nevertheless we give more direct proofs.

The assertions \emph{2}, \emph{3}, \emph{4} and \emph{5} easily result from the \dfnsz, knowing that the results are true in the free case.

\emph{6.} We have defined $\wi{\varphi}$  as the restriction of
 $\wi{\varphi_1}$ at $P$. Since $\varphi_1$ is  an \endo of a free module, we get
$$\preskip.2em \postskip.3em 
\wi{\varphi_1}\circ \varphi_1 =  \det(\varphi_1)\,\Id_{P\oplus Q_1}, 
$$
which gives by restriction at $P$ the desired \egt $\wi{\varphi}\circ \varphi =  \det(\varphi)\,\Id_{P}$, since $\det\varphi=\det\varphi_1$.

\emph{7.} We can reproduce the following \demz, classical in the case of free modules. Consider the \endo 
$$\preskip.4em \postskip.4em 
\psi=X\Id_{P[X]}-\varphi \in\End_{\AX}(P[X]).
$$ 
 By item~\emph{6} we have
$$\preskip.0em \postskip.4em 
~\qquad\qquad\wi{\psi}\psi=\psi\wi{\psi}=\rC \varphi (X)\,\Id_{P[X]}.\eqno (+) 
$$
Moreover, $\wi{\psi}$ is a \pol in $X$ with \coes in $\gA[\varphi]$. Therefore we can write $\wi{\psi}=\som_{k\geq 0}\phi_kX^k$, where each $\phi_k:P\rightarrow P$ is a \pol in~$\varphi$. 
By letting $\rC \varphi(X) = \som_{k\geq 0}a_kX^k$ and by identifying both sides of the \egtz~$(+)$ we obtain (by agreeing to $\phi_{-1} = 0$) 
$$\preskip.4em \postskip.4em 
\phi_{k-1}-\phi_k\varphi =a_k\Id_P\hbox{ for all }k \ge 0. 
$$
Then, $\rC \varphi (\varphi ) =\som_{k\ge 0}(\phi_{k-1}-
\phi_k\varphi)\varphi^k=0$.

\emph{8.} The \pol $\Gamma_\varphi$ has been defined in order to satisfy
$$\preskip.4em \postskip.4em 
\rC{-\varphi}(-X) = -X\Gamma_{\varphi}(X)+\det(\varphi). 
$$
By evaluating $X := \varphi$, 
we obtain $\varphi\Gamma_\varphi(\varphi) =\det(\varphi)\Id_P$ (Cayley-Hamilton \thoz), so $\varphi\Gamma_\varphi(\varphi) =\varphi\wi\varphi$.
By replacing
$\varphi$ by  $\theta := T\Id_{P[T]} + \varphi$, we obtain $\theta\Gamma_\theta(\theta) = \theta\wi\theta$, then $\Gamma_\theta(\theta)=\wi\theta$,
because $\theta$ is a \ndz \elt of $\gA[T,\varphi] = \gA[\varphi][T]$.  We finish the proof by putting $T := 0$.
\end{proof}
\rem
\label{rem det}
The \deter of the identity map  
of every \mptfz, including the module reduced to $\so{0}$, is equal to $1$ (by following the above \dfnz).
\eoe

\begin{corollary}
\label{corDetAnn}
Let $\varphi:P\rightarrow P$ be an \endo of a \mptfz, and $x\in P$ satisfying $\varphi(x)=0$, then $\det(\varphi)x=0$.
\end{corollary}
\begin{proof} Results from $\wi{\varphi}\circ \varphi=\det(\varphi)\Id_P$.
\end{proof}

\vspace{-.7em}
\pagebreak	

%

\subsec{The \polfon and the \polmuz}

We are interested in the \polcar of the identity of a \mptfz. 
It is however simpler to introduce another \pol which is directly related to it and which is the analogue of the \pol $X^k$, which we spoke of at the beginning of Section~\ref{subsec det ptf}.

\begin{definotas}
\label{nota RM}
Let $P$ be a \ptf \Amo and $\varphi$ an \endo of $P$. Consider the $\AX$-module $P[X]$ and define the \pols $\rF{\gA,\varphi}(X)$ and $\rR {\gA,P}(X)$ (or $\rF{\varphi}(X)$ and $\rR {P}(X)$ if the context is clear)
by the following \egtsz

\snic{\rF{\varphi}(X)= \det(\Id_{P[X]}+X\varphi) \quad {\rm  and} \quad
\rR{P}(X)=\det(X\Id_{P[X]}).
}

Therefore $\rR{P}(1+X)=\rF{\,\Id_{P}}(X)$.
\begin{itemize}\itemsep1pt
\item The \pol $\rF{\varphi}(X)$ is called the {\em \polfonz} of the \endo $\varphi$.
\item The \coe of $X$ in the \polfon is called the \ixc{trace}{of an \endo of a \mptfz} of~$\varphi$ and 
is denoted by $\Tr_P(\varphi)$.
\item The \pol $\rR{P}(X)$ is called the {\em \polmuz}\footnote{This terminology is justified by the fact that for a free module of rank $k$ the \polmu is equal to~$X^k$, as well as by \thref{th ptf sfio}.} of the module $P$.%
\label{NOTAPolfon}
\end{itemize}%
\index{rank!\pol of a \mptfz}%
\label{NOTAPolmu}%
\index{polynomial!fundamental ---}%
\index{polynomial!rank ---}%
\end{definotas}

Note that

\snic{
\rF{\varphi}(0)=1=\rR{P}(1),\;\;\rC{\varphi}(0)=\det(-\varphi),\;\;\hbox{and}\;\;
\rF{a\varphi}(X)=\rF{\varphi}(aX),}

but $\rC{\varphi}(X)$ is not always monic 
(see the Example on \paref{ex ideal idempotent}).

Also note that for all $a\in \gA$ we get
\begin{equation}\preskip.4em \postskip.4em
\label{eq2polmu}
\det(a\varphi)=\det(a\,\Id_P)\det(\varphi)=
\rR{P}(a)\det(\varphi).
\end{equation}
We deduce the following \egts
$$\preskip.3em \postskip.2em 
\begin{array}{c}
 \rR{P}(0)  =    \det(0_{\End_\gA(P)}),\,\,\,\,\,\\[1.2mm]  
\mathrigid1mu \rC{-\varphi}(-X)\,=\,\det(\varphi -X\Id_{P[X]})  \, = \,   \det \big(-(X\Id_{P[X]}-\varphi)\big)\, =\, \rR{P}(-1)\rC{\varphi}(X),\\[1.2mm]
\det(\varphi)=\rR{P}(-1)\,\rC{\varphi}(0).\,\,\,\,
\end{array}
$$
The last \egt replaces the \egt $\det(\varphi)=(-1)^k\rC{\varphi}(0)$ valid for the free modules of rank $k$.

\smallskip 


We will say that a \pol $R(X)$ is \ixc{multiplicative}{polynomial} when $R(1)=1$ and $R(XY)=R(X)R(Y)$.\index{polynomial!multiplicative ---}

\begin{theorem}
\label{th ptf sfio} \emph{(The \sfio associated with a \mptfz)}
\index{fundamental sys@\sfioz!associated to a \mptfz}
\begin{enumerate}
\item If $P$ is a \mptf over a \ri $\gA$ the \polmu
$\rR{P}(X)$ is multiplicative.

\item \label{NOTAide} In other words, the \coes of $\rR{P}(X)$ form a \sfioz.
If $\rR{P}(X)=r_0+r_1X+\cdots+r_nX^n$, we denote $r_h$ by $\ide_h(P)$: it is called {\em
the \idm associated with the integer $h$ and with the module $P$}  (if $h>n$ we let $\ide_h(P):=0$).

\item \label{th ptf sfio item reg} Every \polmu $\rR{P}(X)$ is a \ndz \elt of $\AX$.

\item \label{remRang} A \gnn of the \egt $\,\rg(P\oplus Q)= \rg(P)+\rg(Q)\,$ regarding the ranks of the free modules is given for the \mptfs by
$$\preskip-.4em \postskip.3em
\rR{P\oplus Q}(X)=\rR{P}(X)\,\rR{Q}(X)\,.
$$
\item \label{remRang2} If $P\oplus Q\simeq\Ae n$ and $\rR{P}(X)=\sum _{k=0}^nr_kX^k$, then $\rR{Q}(X)=\sum _{k=0}^nr_kX^{n-k}$.

\item \label{remRang3} The \egt $\rR{P}(X)=1$ characterizes, among the \mptfsz, the module $P=\so{0}$. It is also equivalent to $\ide_0(P)=\rR{P}(0)=1$.
\end{enumerate}
\end{theorem}

\begin{proof}
\emph{1}  and \emph{2.} If $\mu_a$ designates  multiplication by $a$ in $P[X,Y]$, we clearly have the \egt $\mu_X\mu_Y=\mu_{XY}$, so $\rR{P}(X)\,\rR{P}(Y)=\rR{P}(XY)$ (\thref{propdef det ptf}.\emph{2}). Since $\rR{P}(1)=\det(\Id_P)=1$, we deduce that the \coes of~$\rR{P}(X)$ form a \sfioz.

\emph{3.} Results from McCoy's lemma (Corollary~\ref{corlemdArtin}). We could also prove it using the basic \plg (by localizing at the~$r_i$'s).\iplg

\emph{4.} Results from item~\emph{3} in \thref{propdef det ptf}.

\emph{5.} Results from items~\emph{3} and~\emph{4} since $\big(\som_{k=0}^nr_kX^k\big)\big(\som_{k=0}^nr_{n-k}X^k\big)=X^n$.

\emph{6.} We have $r_0=\det(0_{\End(P)})$. Since the $r_i$'s form a \sfioz, the \egts $\rR P=1$ and $r_0=1$ are \eqvesz.\\
If $P=\so{0}$, then $0_{\End(P)}=\Id_P$, so $r_0=\det(\Id_P)=1$.\\
If $r_0=1$, \hbox{then $0_{\End(P)}$} is \ivz, therefore $P=\so{0}$.
\end{proof}

If $P$ is a free \Amo of rank $k$, we have $\rR{P}(X)=X^k$, the following \dfn is therefore a legitimate extension from free modules to \mptfsz.

\begin{definition}
\label{def ptf rank constant}
A \mptf $P$ is said to be {\em of rank equal to $k$} if $\rR{P}(X)=X^k$. If we do not specify the value of the rank, we simply say that the module is {\em of constant rank}.
We will use the notation $\rg(M)=k$ to indicate that a module 
(assumed to be \pro of constant rank) is of rank $k$. \label{Nota2rang}
\index{rank!module of constant ---}
\end{definition}

Note that by Proposition~\ref{propAnnul}, every \pro module of rank $k>0$ is faithful.

\pagebreak	

\begin{fact}\label{factPolCarRangConstant}
The \polcar of an \endo of a \mrc $k$ is \mon of degree $k$.
\end{fact}
\begin{proof} We can give an elegant direct \dem (see Exercise~\ref{exoPrecisionsDet3}). We could also avoid all effort and use a \lon argument, by relying on the local structure \tho and on Fact~\ref{fact.det loc}, which asserts that everything goes well for the \polcar by \lonz.
\end{proof}

The convention in the following remark allows  for a more uniform formulation of the \thos and the proofs hereinafter.


\medskip \rem
\label{conven rgc}
When the \ri $\gA$ is reduced to $\so{0}$, all the \Amos are trivial. Nevertheless, in accordance with the above \dfnz, the null module over the null \ri is a \mrc equal to $k$, for any value of the integer $k\geq 0$.
Moreover, it is immediate that if a \mptf $P$ has two distinct constant ranks, then the \ri is trivial. We have $\rR{P}(X)=1_\gA X^h=1_\gA X^k$ with $h\neq k$ therefore the
\coe of~$X^h$ is equal to both $1_\gA$ and $0_\gA$. 
\eoe


\subsec{Some explicit computations}

\vspace{3pt}
The \polfon of an \endo $\varphi$ is easier to use than the \polcarz.
This comes from the fact that the \polfon is invariant when we add \gui{as a direct sum} a null \endo to $\varphi$. This allows us to systematically and easily reduce the computation of a \polfon to the case where the \pro module is free. Precisely, we are able to compute the previously defined \pols by following the lemma stated below.

\begin{lemma}
\label{lem calculs} {\rm (Explicit computation of the \deterz, of the \polfonz, of the \polcarz, of the \polmu and of the cotransposed \endoz)} \\
Let $P\simeq \Im F$ be an \Amo with $F\in\GAn(\gA)$. Let $Q=\Ker(F)$, such that $P\oplus Q\simeq \Ae n$, \hbox{and $\I_n-F$} is the matrix of the \prnz~$\pi_Q$ over~$Q$ \paralm to~$P$.
An \endo $\varphi$ of $P$ is \care by the matrix~$H$ of the \endo $\varphi_0= \varphi\oplus 0_Q$ of $\Ae n$. Such a matrix $H$ is subjected to the unique restriction $F\cdot H\cdot F=H$. Let $G=\I_n-F+H$.

\begin{enumerate}\itemsep0pt
\item  Computation of the \deterz:
$$\preskip.4em \postskip.2em
\det(\varphi)=\det(\varphi\oplus\Id_Q)=\det(G).
$$
\item  Therefore also
$$\preskip.4em 
\arraycolsep2pt\begin{array}{rcl}
\det(X\Id_{P[X,Y]}+Y\varphi)&=&  \det \big((X\Id_{P[X,Y]}+Y\varphi)\oplus\Id_Q \big)\,=\\[1mm]
\det(\I_n-F+XF+YH)&  = &  \det(\I_n+(X-1)F+YH).
\end{array}$$

\vspace{-1em}
\item  Computation of the \polmu of $P$:
$$\preskip.4em \postskip.4em
\rR{P}(1+X)=\det \big((1+X)\Id_{P[X]} \big)=\det(\I_n+XF),
$$ 
in particular,

\snic{ \rR{P}(0)=\det(\I_n-F),}

\smallskip
and $\rR{P}(1+X)=1+u_1X+\cdots+u_nX^n$, where $u_h$ is the sum of the principal minors of order $h$ of the matrix $F$.

\item  Computation of the \polfon of $\varphi$:
$$\preskip.4em \postskip.4em
\rF{\varphi}(Y)=\det(\Id_{P[Y]}+Y\varphi) =\det(\I_n+YH)=
1+\som_{k=1}^n v_k Y^k,
$$
where $v_k$ is the sum of the principal minors of order $k$ of the matrix~$H$.
In particular, $\Tr_P(\varphi)=\Tr(H)$.

\item  Computation of the \polcar of $\varphi$:
$$\preskip.4em \postskip.4em
\rC{\varphi}(X)=\det(X\Id_{P[X]}-\varphi) =  \det(\I_n-H+(X-1)F).
$$
\item  Computation of the cotransposed \endo $\wi{\varphi}$ of $\varphi$: it is defined by the matrix
$$\preskip.0em \postskip.4em
\wi{G}\cdot F =  F \cdot \wi{G}=  \wi{G}-\det(\varphi)(\In-F).
$$
\end{enumerate}
\end{lemma}

For the last item we apply item~\emph{3} of \thref{propdef det ptf} with $\varphi$ and $\Id_Q$ by remarking that $G$ is the matrix of $\psi=\varphi\oplus \Id_Q=\varphi_0+\pi_Q$.

Note that the \polcar of $\Id_P$ is equal to $\rR{P}(X-1)$.

The following fact is an immediate consequence of Proposition~\ref{propPtfExt} and of the previous lemma.

\begin{fact}
\label{fact.det loc}
The \deterz, the cotransposed \endoz, the \polcarz, the \polfon and the \polmu are well-behaved under \eds via a \homo $\gA\rightarrow \gB$. \\
In particular, if $\varphi :P\rightarrow P$ is an \endo of a \ptf $\gA$-module and $S$ a \mo of $\gA$, then $\det(\varphi)_S=\det(\varphi_S)$ (or, if we prefer, $\det(\varphi)/1=_{\gA_S}\det(\varphi_S)$). The same thing holds for the cotransposed \endoz, the \polfonz, the \polcar and the \polmuz.
\end{fact}

\exl
\label{ex ideal idempotent}
Let $e$ be an \idm of $\gA$ and $f=1-e$. The module $\gA$ is a direct sum of the submodules $e\gA$ and $f\gA$ which are therefore \ptfsz.
The $1\times 1$ matrix having for unique \coe $e$ is a matrix $F$ whose image is $P=e\gA$. For $a\in \gA$ consider $\mu_a=\mu_{P,a}\in\End_\gA(P)$. The matrix $H$ 
has for unique \coe $ea$.
We then have, by applying the previous formulas,
$$\begin{array}{c}
\det(0_{e\gA})=f,\; \rR{e\gA}(X)=f+eX,\;\rC{\Id_{e\gA}}(X)=f-e+eX,\\[1mm]
\; \det(\mu_a)=f+ea,\\[1mm]
\rF{\mu_a}(X)=1+eaX,\;\rC{\mu_a}(X)=1-ea+e(X-1)=f-ea+eX.
\end{array}$$
Note that the \polcar of $\mu_a$ is not monic \hbox{if $e\neq 1,0$},
and we indeed have the Cayley-Hamilton \tho
$$\preskip.4em \postskip.4em 
\rC{\mu_a}(\mu_a)=(f-ea)\Id_{e\gA}+e\mu_a= (f-ea+ea) \Id_{e\gA} = f\Id_{e\gA} =
0_{e\gA}. 
$$
\vskip-2.7em\eoe

\subsubsection*{With a \sycz}
When we use a \syc Lemma~\ref{lem calculs} leads to the following result.

\begin{fact}\label{factMatriceEndo}
Let $P$ be a \mptf with a \syc $\big((\xn),(\aln)\big)$ and $\varphi$ be an \endo of $P$.\\
Recall (Fact~\ref{factMatriceAlin}) that we can encode $P$ by the matrix 

\snic{F\eqdefi \big(\alpha_i(x_j)\big)_{i,j\in\lrbn}}

($P$ is \isoc to $\Im F\subseteq\Ae n$ by means of $x\mapsto \pi( x) =\tra[\,\alpha_1(x) \;\cdots\; \alpha_n(x)\,]$). In addition the \endo $\varphi$ is represented by the matrix 

\snic{H\eqdefi \big(\alpha_i(\varphi(x_j))\big)_{i,j\in\lrbn}}

which satisfies $H=HF=FH$. 
\begin{enumerate}
%
%
\item We have $\rF\varphi(X)=\det (\In+XH)$ and $\Tr(\varphi)=\Tr(H)=\sum_i\alpha_i \big(\varphi(x_i)\big)$.
\item 
For $\nu \in P\sta$ and $x$, $y \in P$, recall that $\theta_P(\nu \te x)(y)=\nu(y)x$.
The trace of this \endo is given by $\Tr_P \big(\theta_P(\nu \te x)\big)=\nu(x)$.
\end{enumerate}
\end{fact}
\begin{proof}
The matrix $H$ is also that of the \Ali $\varphi_0$ introduced in Lemma~\ref{lem calculs} 

\snic{\pi( x) + y\mapsto \pi \big(\varphi(x)\big)$ with $\pi( x)\in \Im F$ and $y\in \Ker F.}

Item~\emph{2} therefore results from Lemma~\ref{lem calculs}.

\emph{3.} By item~\emph{2}, we have

\snic {\mathrigid1mu
\Tr \big(\theta_P(\nu \te x)\big) = 
\sum_i \alpha_i(\nu(x_i)x) = \sum_i \nu(x_i) \alpha_i(x) = 
\nu\big(\sum_i \alpha_i(x)x_i\big)= \nu(x). 
}
\end{proof}
%

\begin{lemma}\label{lemTraceProT}
Let $M$, $N$ be two  \ptfs \kmos and let $\varphi\in\alb\End_\gk(M)$ and $\psi\in\End_\gk(N)$ be \endosz. \\
Then, $\Tr_{M\te N}(\varphi\otimes \psi)=\Tr_M(\varphi)\Tr_N(\psi)$.
\end{lemma}
\begin{proof}
Consider \sycs for $M$ and $N$ and apply the formula for the trace of the \endos (Fact~\ref{factMatriceEndo}).
\end{proof}

\vspace{-.7em}
\pagebreak	

\subsec{The annihilator of a \mptfz}
\label{subsecAnnul}

We have already established certain results regarding this annihilator by relying on the local structure \tho for  \mptfsz, proven by using the \idfs (see Lemma~\ref{lemAnnMptf}).

Here we give some additional results by using a \dem that does not rely on the local structure \thoz.
\begin{proposition}
\label{propAnnul}
Let $P$ be a \ptf \Amoz. Consider the \id $J_P= \gen{ \alpha (x) \,|\,\alpha \in P\sta,\; x\in P  }$. Let $r_0=\rR{P}(0)=\ide_0(P)$.
\begin{enumerate}
\item $\gen{r_0}=\Ann(P)=\Ann(J_P)$.
\item $J_P=\gen{s_0}$, where $s_0$ is the \idm $1-r_0$.
\end{enumerate}
\end{proposition}
\begin{proof} We obviously have $\Ann(P)\subseteq \Ann(J_P)$.
Let $\big((x_i)_{i\in\lrbn},(\alpha_i)_{i\in\lrbn} \big)$ be a \syc over $P$. Then

\snic{J_P=
\gen{ \alpha_i(x_j)\,;\,{i,j\in\lrbn}},
}

and the \mprn $F= \big(\alpha_i(x_j)\big)_{i,j\in\lrbn}$ has an image \isoc 
to~$P$. By \dfnz, $r_0$ is the \idm $r_0=\det(\I_n-F)$. 
\\ Since $(\I_n - F)F = 0$, we have $r_0F = 0$, \cad $r_0P=0$. \\
Therefore $\gen{r_0}\subseteq \Ann(P)\subseteq \Ann(J_P)$ and $J_P\subseteq \Ann(r_0)$.\\
Moreover, we have $\I_n-F\equiv \I_n$ modulo $J_P$, so by taking the \detersz, we~have~$r_0\equiv 1$ modulo $J_P$, \cad $s_0\in J_P$, then $\Ann(J_P)\subseteq\Ann(s_0)$.\\
We can therefore conclude
$$\preskip.4em \postskip.4em\mathrigid 2mu 
\gen{r_0}\subseteq \Ann(P)\subseteq\Ann(J_P)\subseteq\Ann(s_0)=\gen{r_0}
\hbox{ and }\gen{s_0}\subseteq J_P\subseteq \Ann(r_0) = \gen{s_0}. 
$$
\end{proof}

\subsec{Canonical decomposition of a \pro module}
\label{subsec decomp ptf}

\begin{definition}\label{defi comp ptf}
Let $P$ be a \ptf \Amo and $h\in\NN$.
If $r_h=\ide_h(P)$,
we denote by $P\ep{h}$ the \Asub $r_hP$. 
It is called the {\em component of the module $P$ in rank~$h$}.
\label{NOTAep}
\end{definition}

Recall that, for an \idm $e$ and an \Amo $M$, the module obtained by \eds to $\gA[1/e]\simeq\aqo{\gA}{1-e}$ can be identified with the submodule~$eM$, itself \isoc to the quotient module ${M}\sur{(1-e)M}$.

\begin{theorem}
\label{propdef comp ptf}\label{th decomp ptf}
Let $P$ be a \ptf \Amoz.
\begin{enumerate}
\item  The module $r_hP=P\ep{h}$ is a \pro $\gA[1/r_h]$-module of rank $h$.
\item The module $P$ is the direct sum of the $P\ep{h}$'s.
%
\item The \id $\gen{r_0}$ is the annihilator of the \Amo $P$.
\item For $h>0$, $P\ep{h}= \so{0}$ implies $r_h= 0.$
\end{enumerate}
\end{theorem}
\begin{proof}
\emph{1.} Localize at $r_h$: we obtain $\; \rR{P\ep{h}}(X)=_{\gA[1/r_h]}\rR{P}(X) =_{\gA[1/r_h]}X^h$.

\emph{2.} Because the $r_h$'s form a \sfioz.

\emph{3.} Already proved (Proposition~\ref{propAnnul}).

\emph{4.} Results \imdt from item~\emph{3.}
\end{proof}

Note that, except if $r_h=1$ or $h=0$, the module $r_hP$ is not of constant rank  when considered as an \Amoz..

The previous \tho gives a \gui{structural} \dem of \thref{th.ptf.idpt}.

\medskip
\rem
\label{rem th decomp ptf}
If $P$ is (\isoc to) the image of a \mprn $F$ the \idms $r_k=\ide_k(P)$ attached to the module $P$ can be linked to the \polcar of the matrix $F$ as follows

\snic{\det(X\In-F)=\som_{k=0}^n r_kX^{n-k}(X-1)^k.
}

(Note that the $X^{n-k}(X-1)^k$ form a basis of the module of  \pols of degree~$\leq n$,
triangular with respect to the usual basis.) \eoe

\subsec{Rank \pol and \idfsz}
\label{subsecRangRang}


The \dem of  \thref{corth.ptf.sfio} that follows relies on \thref{prop Fitt ptf 2}, which asserts that a \mptf becomes free after \lon at  \ecoz.

We have placed the \tho here because it answers to the questions that we naturally ask ourselves after \thref{th ptf sfio}.
First, check that a \mprn is of rank $k$ \ssi its image is a \mrc $k$. 
More \gnltz, characterize the \sfio that occurs in the \polmu in terms of the \idfs of the module.

Actually, we can give an alternative \dem of \thref{corth.ptf.sfio} without taking the route of a localization argument, by making use of exterior powers (see Proposition~\ref{prop puissance ext}).

Let us point out that for a \mpf $M$ the \egt $\cF_h(M)=\gen{1}$ means that $M$ is locally generated by $h$ \elts (we have seen this in the case~$h=1$ in \thref{propmlm}, in the \gnl case, see the local number of \gtrs lemma on \paref{lemnbgtrlo} and \Dfnz~\ref{deflocgenk}).

\begin{theorem}
\label{corth.ptf.sfio} 
\emph{(Local structure and \idfs of a \mptfz, 2)}\\
Let $F\in\GAq(\gA)$, $P\simeq\Im  F$ and $\rR{P}(X)=\sum_{i=0}^qr_iX^i$.
\begin{enumerate}
\item  Let $S(X)=\rR{P}(1+X)=1+u_1X+\cdots+u_qX^q$ ($u_h$ is the sum of the principal minors of order $h$ of the matrix $F$).\\
We have, for all $h\in\lrb{0..q}$,
$$\preskip.2em \postskip.2em 
\formule{\cD_h(F)=\gen{r_h+\cdots+r_q}=\gen{r_h,\ldots ,r_q}= \gen{u_h,\ldots
,u_q}\\[1mm]
\cF_h(P)= \gen{r_0+\cdots +r_h}= \gen{r_0,\ldots ,r_h}} 
$$
\vspace{-.7em}
\pagebreak	

\item  In particular%
\begin{enumerate}
\item  $\rg(F)=h\iff \rg(P)=h$, 
\item  $\rg(F)\leq h\iff \deg \rR{P}\leq h,$
\item  $\rg(F)> h\iff r_0=\cdots=r_h=0 \iff \cF_{h}(P)=0.$
\end{enumerate}
\end{enumerate}
\end{theorem}
\begin{proof}
The \egt $ \gen{u_h,\ldots,u_q}=\gen{r_h,\ldots ,r_q}$ results from the \egts

\snic{S(X)=\rR{P}(1+X)\;\hbox{ and }\;\rR{P}(X)=S(X-1).}

To check the \egts 
$\cD_h(F)=\gen{r_h+\cdots+r_q}=\gen{r_h,\ldots ,r_q}$ and 

\snic{\cD_{q-h}(\I_q-F)= \gen{r_0+\cdots +r_h}= \gen{r_0,\ldots ,r_h},}

it suffices to do it after \lon at \ecoz. However, the kernel and the image of $F$ become free after \lon at \eco (\thref{theoremIFD} or \thref{prop Fitt ptf 2}), and the matrix therefore becomes similar to a standard \mprnz.  
\end{proof}

\section{Properties of finite character}
\label{secSalutFini}

The purpose of this section is to illustrate the idea that the good concepts in \alg are those that are controllable by finite procedures. 

We have in mind to highlight \gui{good \prtsz.} There are naturally those that  submit to the \lgb principle: for the \prt to be true it is sufficient and \ncr that it be true after \lon at \mocoz. 
It is a phenomenon that we have frequently encountered, and will continue to encounter hereafter.

Recall that a \prt is said to be \gui{of \carfz} if it is preserved by \lon (by passing from $\gA$ to $S^{-1}\gA$) and if, when it is satisfied after \lon at $S$, then it is satisfied after \lon at $s$ for some $s\in S$.

In Fact\eto \ref{factPropCarFin} we proved in \clama that for the \carf \prtsz, the \plgc (\lon at \mocoz) is \eqv to the \plga (\lon at all the \idemasz).
However, a \prco of the \plgc a priori contains more precise information than a classical \dem of the \plgaz.

\begin{proposition}\label{propFiniBon1} Let $S$ be a \mo of $\gA$.
\begin{enumerate}
  \item Let $AX=B$ be a \sli over $\gA$. Then, if it admits a solution in $\gA_{S}$, there exists an $s\in S$ such that it admits a solution in~$\gA_{s}$.
  \item Let $M$ and $N$ be two \Asubs of a same module, \emph{with $M$ \tfz}. Then, if $M_{S}\subseteq N_{S}$, there exists an $s\in S$ such that $M_{s}\subseteq N_{s}$.
  \item \label{propFiniBon1-3}
  Let $\gA$ be a \emph{\cohz} \riz, $M$, $N$, $P$ be \emph{\pf} \Amosz, and $\varphi:M\to N$, $\psi:N\to P$ be two \alisz. \\
If the sequence $M\vers{\varphi}N\vers{\psi}P$ becomes exact after \lon at $S$ there exists an $s\in S$ such that the sequence becomes exact after \lon at $s$.
  \item \label{propFiniBon1-4} Let $M$ and $N$ be two \emph{\pf} \Amosz. Then, \hbox{if $M_{S}\simeq N_{S}$,} there exists an $s\in S$ such that $M_{s}\simeq N_{s}$.
  \item Let $M$ be a \pf \Amoz. If $M_S$ is free, there exists some $s \in S$ such that $M_s$ is free. Similarly, if $M_S$ is \stlz, there exists some $s \in S$ such that $M_s$ is \stlz.
  \item If a \mpf becomes \pro after \lon at $S$, it becomes \pro after \lon at an \elt $s$ of $S$.
%
\end{enumerate}
\end{proposition}
\begin{proof}
Let us prove item~\emph{\ref{propFiniBon1-3}}.
We first find some $u\in S$ such that $u\,\psi \big(\varphi(x_{j})\big)=0$ for \gtrs $x_{j}$'s of $N$. We deduce that $\psi\circ\varphi$ becomes null after \lon at $u$. Moreover, the hypotheses assure us that $\Ker\psi$ is \tfz. Let $\yn$ be \gtrs of $\Ker\psi$.
For each of them we find a $z_{j}$ in $N$ and an $s_{j}\in S$ such that $s_{j}(\varphi(z_{j})-y_{j})=0$. We take for $s$ the product of $u$ and the $s_{j}$'s.\\
Let us prove item~\emph{\ref{propFiniBon1-4}}. 
Let $G$ and $H$ be \mpns for $M$ and $N$.
Let $G_{1}$ and $H_{1}$ be the two matrices given in Lemma~\ref{lem pres equiv}. By hypothesis there exist two square matrices $Q$ and $R$ with \coes in $\gA$ such that $v=\det(Q)\det(R)\in S$ and 
$Q\,G_{1}=_{\gA_{S}} H_{1}\,R$.
This means that we have over~$\gA$ an \egt
$$\preskip.0em \postskip.4em
w\,(Q\,G_{1}- H_{1}\,R)=0,\quad w \in S
.$$
It therefore suffices to take $s=vw$.
\end{proof}

We have seen that the \eds is well-behaved with respect to tensor products, exterior powers and symmetrical powers.
For the functor $\Lin_\gA$ things do not always go so well. The following are important results for the remainder of this work.

\begin{proposition}
\label{fact.hom egaux}
Let $f:M\rightarrow N$ and $g:M\rightarrow N$ be two \alis between \Amosz, {\em with $M$ finitely generated}.  Then, $f_S=g_S$ \ssi there exists an $s\in S$ such that $sf=sg$.
In other words, the canonical map $\big(\Lin_\gA(M,N)\big)_S\rightarrow \Lin_{\gA_S}(M_S,N_S)$ is injective.
\end{proposition}

\smallskip 
\begin{proposition}\label{fact.homom loc pf}
Let $M$ and $N$ be two \Amos and $\varphi:M_S\rightarrow N_S$ be an \Aliz. We assume that $M$ is \pfz, or that~$\gA$ is integral, $M$ \tf and $N$ torsion-free (\cad $a\in\gA$, $x\in N$, $ax=0$ implies $a=0$ or $x=0$).\\ 
Then, there exists an \Ali $\phi:M\rightarrow N$ and some $s\in S$ such that
$$\forall x\in M\quad  \varphi(x/1)  = \phi(x)/s ,$$
and the canonical map $\big(\Lin_\gA(M,N)\big)_S\rightarrow \Lin_{\gA_S}(M_S,N_S)$ is bijective.
\end{proposition}
\begin{proof}
The second case, which is easy, is left to the reader. To follow the \dem of the first case one must look at the following figure.
Suppose that $M$ is the cokernel of the \ali $g:\Ae m\rightarrow \gA^q$ with a matrix $G=(g_{i,j})$ with respect to the canonical bases, then by Fact~\ref{fact.sexloc} the module~$M_S$ is the cokernel of the \ali $g_S:\gA_S^m\rightarrow \gA_S^q$, represented by the matrix $G_S=(g_{i,j}/1)$ over the canonical bases.
Let 
$$\preskip.4em \postskip.4em \mathrigid 2mu
\Ae m\vers{j_m} \gA_S^m,\, \gA^q\vers{j_q} \gA_S^q,\
M\vers{j_M} M_S,\ N\vers{j_N} N_S,\ \gA^q\vers{\pi} M,\
\gA_S^q\vers{\pi_S} M_S, 
$$
be the canonical maps. Let $\psi:=\varphi \circ \pi_S$, so that $\psi\circ g_S=0$. 
Therefore $\psi\circ g_S\circ j_m=0=\psi \circ j_q \circ g$. 
There exists some $s\in S$, a common denominator for the images under $\psi$ of the vectors of the canonical basis. Hence a \ali $\Psi:\gA^q\rightarrow N$ with $(s\psi)\circ j_q=j_N\circ \Psi$.

\vspace{-.6em}
\begin{figure}[htbp]
\centerline{
\xymatrix{
\Ae m \ar[d]_{\displaystyle g} \ar[rrr]^{\displaystyle j_m} & &
&\gA_S^m\ar[d]^{\displaystyle {g_S}} &
\\
\Ae q \ar[d]_{\displaystyle \pi}
   \ar@{-->}[rdd]^(.4){\displaystyle \Psi}\ar[rrr]^{\displaystyle j_q} & &
&\gA_S^q\ar[d]_{\displaystyle\pi_S}\ar[rdd]^{\displaystyle \psi} &
\\
M\ar[rrr]
^{\displaystyle j_M}\ar@{-->}[rd]_{\displaystyle \phi}
& & &M_S\ar[rd]_{\displaystyle \varphi} &
\\
 & N\ar[rrr]^{\displaystyle j_N} & & &N_S}
}
\caption{Localization of the \homosz}
\label{fig1}
\end{figure}
 
\vspace{-.6em}
Thus,~$j_N\circ\Psi\circ  g= s(j_m\circ g_S\circ \psi)=0$.
By Proposition~\ref{fact.hom egaux} applied to~$\Psi\circ g$, the \egt $j_N\circ(\Psi\circ  g)= 0$ in $N_S$ implies that there exists an $s'\in S$ such that $s'(\Psi\circ  g)=0$. Therefore $s' \Psi$ can be factorized in the form $\phi\circ \pi$.
We then obtain 
$$\preskip.2em \postskip.3em 
(ss'\varphi)\circ j_M\circ \pi=ss'(\varphi \circ\pi_S\circ j_q)
= ss'\psi\circ j_q=s'j_N\circ\Psi =
j_N\circ \phi\circ \pi, 
$$
and since $\pi$ is surjective, $ss'\varphi \circ j_M=j_N\circ\phi$. Thus, for all $x\in M$, we have~$\varphi(x/1)=\phi(x)/ss'$.
\end{proof}

\begin{corollary}\label{corfinchapMPF}
Suppose that $M$ and $N$ are \pfz, or that they are \tfz, torsion-free and that $\gA$ is integral. If $\varphi : M_S \to N_S$ is an \isoz, there exist
an $s \in S$ and an \iso $\psi : M_s \to N_s$ such that $\psi_S = \varphi$.
\end{corollary}
\begin{proof}
Let $\varphi' : N_S \to M_S$ be the inverse of $\varphi$. By the previous proposition, there exist $\phi : M \to N$, $\phi' : N \to M$, $s \in S$, $s' \in S$ such that $\varphi = \phi_S/s$, $\varphi' = \phi'_S/s'$. Let $t = ss'$ and define $\psi = \phi_t/s : M_t \to N_t$, $\psi' = \phi'_t/s' : N_t \to M_t$. Then, $(\psi' \circ \psi)_S$ is the identity over $M_S$, and $(\psi \circ\psi')_S$ is the identity over $N_S$. We deduce the existence of a $u \in S$ such that $(\psi' \circ \psi)_{tu}$ is the identity over $M_{tu}$, and~$(\psi \circ\psi')_{tu}$ is the identity over $N_{tu}$. Consequently, $\psi_{tu} : M_{tu} \to N_{tu}$ is an \iso such that $(\psi_{tu})_S = \varphi$.
\end{proof}

\Exercices

\rdb\begin{exercise}
\label{exoptf0Lecteur}
{\rm  We recommend that the \dems which are not given, or are sketched, or
left to the reader,
etc, be done.
But in particular, we will cover the following cases.
\begin{itemize}\itemsep1pt
\item Show Facts~\ref{factdefiMPRO} and~\ref{factMatriceAlin}.
\item Check the details of Lemma~\ref{lem calculs}.
\item \label{exofact.homegaux}
Show Fact~\ref{fact.hom egaux} as well as the second case in Proposition~\ref{fact.homom loc pf}.
\end{itemize}
}
\end{exercise}

\rdb\vspace{-1em}
\begin{exercise}
 \label{exoProjmemeImage} (Projectors having the same image)\\
 {\rm  Let $a$, $c$ be in a not \ncrt commutative \ri $\gB$.
 \Propeq
\begin{itemize}\itemsep1pt
\item  $ac=c$ and $ca=a$.
\item $a^2=a$, $c^2=c$ and $a\gB=c\gB$.
\end{itemize}
In such a case let $h=c-a$ and $x=1+h$. Show the following results. 

\snuc{ha=hc=0$, $ah=ch=h$, $h^2=0$, $x\in\gB\eti$, $ax=c$,
$xa=x^{-1}a=a$ and $\fbox{$x^{-1}ax=c$}.}

It should be noted in passing that the \egt $ax = c$
returns the \egt $a\gB = c\gB$.
\\
Special case. $\gA$ is a commutative \riz, $M$ is an \Amoz, and $\gB=\End_\gA(M)$: two  \prrs that have the same image are similar.

} \end{exercise}

\rdb\vspace{-1em}
\begin{exercise}
\label{exo2.4.1} (Two \eqv \prrs are similar) \\
{\rm In a (not \ncrtz) commutative \ri $\gB$, consider two \emph{\eqvsz} \idms ($a^2= a$, $b^2=b$, $\exists p,\,q\in\gB\eti,\, b=paq$). We will show that they are \emph{conjugate} ($\exists d\in\gB\eti, \;dad^{-1}=b$).
\begin {itemize}\itemsep1pt
\item
In this question, $a$, $b \in \gB$ are \eqv ($b=paq$), but are not assumed to be \idmsz. 
Show that the \elt $c=p^{-1}bp$ satisfies $a\gB = c\gB$.

\item In particular, if $b$ is \idmz, $c$ is a conjugate \idm of $b$ which satisfies $a\gB = c\gB$. Conclude by using the previous exercise.
\end{itemize}
Special case. $\gA$ is a commutative \riz, $M$ is an \Amoz, and $\gB=\End_\gA(M)$: two \eqv \prrs of $M$  
 are similar.
}
\end{exercise}

\rdb\vspace{-1em}

\pagebreak	

\begin{exercise}\label{exoSchanuelVariation} {(An important consequence of Schanuel's lemma~\ref{corlemScha})}
\\
{\rm
\emph{1.} We consider two exact sequences
\vspace{-1mm}
\[\arraycolsep2pt
\begin{array}{ccccccccccccccc}
0 &\rightarrow& K& \to & P_{n-1}& \to & \cdots & \to & P_1 & \vers{u} & P_0
& \to & M& \rightarrow& 0
\\
0 &\rightarrow& K'& \to & P'_{n-1}& \to & \cdots & \to & P'_1 & \vers{u'} & P'_0
& \to & M&  \rightarrow& 0
\end{array}
\]
with the \pro modules $P_i$ and $P'_i$. Then, we obtain an \iso
$$
\preskip.4em \postskip.2em\ndsp 
K \oplus \bigoplus\limits_{i \equiv n-1 \mod 2} P'_i \oplus
\bigoplus\limits_{j \equiv n \mod 2} P_j \quad \simeq\quad 
K' \oplus \bigoplus\limits_{k \equiv n-1 \mod 2} P_k \oplus
\bigoplus\limits_{\ell \equiv n \mod 2} P'_\ell. 
$$
\emph{2.} Deduce that if we have an exact sequence where the $P_i$'s, $i\in \lrbn$, are \pros
$$\preskip-.4em \postskip.4em 
0 \rightarrow P_n \to  P_{n-1} \to  \cdots  \to  P_1 \to P_0 \to M
\rightarrow 0, 
$$
then, for every exact sequence
$$\preskip.4em \postskip.4em 
0 \to K' \to  P'_{n-1} \to  \cdots  \to  P'_1  \to  P'_0 \to M
\to 0,
$$
 where the $P'_i$'s are \prosz, the module $K'$ is \egmt \proz.
}
\end{exercise}

\rdb\vspace{-1em}
\begin{exercise}
\label{exoSuitExPtfs}
{\rm
Consider an exact sequence 
composed of 
\mptfs
$$\preskip.0em \postskip.4em\;\; 
0  \lora
P_n \vvers {u_n} 
P_{n-1} \vvvers {u_{n-1}}  P_{n-2}
\lora \cdots \lora
P_2 \vvers {u_2} 
P_{1} \lora   0\,. 
$$
Show that
${  {\bigoplus\limits_{i\rm\ odd} P_i \;\simeq\;
\bigoplus\limits_{j\rm\ even} P_j}}.$

Deduce that if the $P_i$'s for $i \ge 2$ are \stlsz, similarly for~$P_1$.
}
\end{exercise}

\rdb\vspace{-1em}
\begin{exercise}
\label{exoZerDiR}
{\rm  Show that 
\propeq
\begin{itemize}\itemsep1pt
\item The \ri $\gA$ is \zedrz.
\item The \pf \Amos are always \ptfsz.
\item Every module $\aqo{\gA}{a}$ is \ptfz.
\end{itemize}
(In other words, show the converse for item~\emph{1} in \thref{propZerdimLib}.)
}
\end{exercise}

\rdb\vspace{-1em}
\begin{exercise}\label{exoProjRang1} {(Projectors of rank $1$, see Proposition~\ref{pmlm})}
\\
{\rm
Let $A = (a_{ij}) \in \Mn(\gA)$. We examine \syps in the  $a_{ij}$'s whose zeros define the sub\vrt $\GA_{n,1}(\gA)$ of~$\Mn(\gA)$.
We denote by $\cD'_2(A)$ the \id generated by the minors having one of the \gui{four corners} on the diagonal (not to be mistaken with the principal minors, except when $n=2$).

\emph {1.}
If $A$ is a \prr of rank $\le 1$, then $\Ann\, A$ is generated by $1-\Tr A$ (\idmz). In particular, a \prr of rank $1$ is of trace~$1$.

\emph {2.}
The \egts $\Tr A = 1$ and $\cD'_2(A) = 0$ imply $A^2 = A$ and $\cD_2(A) = 0$. In this case, $A$ is a \prr of rank $1$ 
(but we can have $\Tr A = 1$ and $A^2 = A$ without having $\cD_2(A) = 0$, e.g.\ for a \prr of rank $3$ over a \ri in which $2=0$.) 
Consequently, for an arbitrary matrix $A$ we have
$$\preskip.3em \postskip.25em 
\gen {1-\Tr A} + \cD_1(A^2 - A) \subseteq
\gen {1-\Tr A} + \cD'_2(A) = \gen {1-\Tr A} + \cD_2(A) 
$$
without \ncrt having the left-\egtz.

\emph {3.}
We consider the \pol  $\det\!\big
(\In+(X-1)A\big)$ (if $A\in\GA_n(\gA)$, it is the rank \pol of the module $P=\Im A$) and we denote by $r_1(A)$ its
\coe in $X$. We therefore have the \egt of the three following \idsz, defining 
the sub\vrt $\GA_{n,1}(\gA)$ of~$\Mn(\gA)$:
$$\preskip.2em \postskip.3em 
\gen {1-\Tr A} + \cD'_2(A) = \gen {1-\Tr A} + \cD_2(A) =
\geN {1-r_1(A)} + \cD_1(A^2 - A). 
$$
Specify the cardinality of each \sgrz. 

}
\end{exercise}

\rdb\vspace{-1em}
\begin{exercise}\label{exoMatProjRang1CoeffReg}
{(Projector of rank $1$ having a \ndz \coez)}\\
{\rm 
Let $A = (a_{ij}) \in \GA_n(\gA)$ be a \prr of rank $1$, $L_i$ its row $i$, $C_j$ its column~$j$.

\emph {1.}
Provide a direct proof of the matrix \egt $C_j \cdot L_i = a_{ij}A$. By noticing that $L_i\cdot C_j = a_{ij}$, deduce the \egt of \ids $\gen {L_i}\gen{C_j} = \gen {a_{ij}}$.

\emph {2.}
Suppose $a_{ij}$ is \ndzz; so $\gen {L_i}$ and $\gen{C_j}$ are \iv \idsz, inverses of each other. Provide a direct \dem of the exactitude in the middle of the sequence

\snic {
\Ae n\vvvers{\In-A} {\Ae n}\vvers{L_i} \gen {L_i}\to 0
}

and therefore conclude that $\gen {L_i} \simeq \Im A$.

\emph {3.} Prove that the matrix $A$ is entirely determined by $L_i$ and $C_j$. More \prmtz, if $\gA$ is a \ri with explicit \dvez,
\begin{itemize}
\item compute the matrix $A$,
\item deduce the condition for which the row $L$ and the column $C$ can be the row~$i$ and the column $j$ of a \mprn of rank $1$ (we suppose that the common \coe in position $(i,j)$ is \ndzz). 
%
\end{itemize}
 
\emph {4.}
Let $C \in \Im A$, $\tra L \in \Im \tra A$ and $a = L\cdot C$.
Show the matrix \egt $C \cdot L = aA$ and deduce the \egt of \ids $\gen {L}\lra{\tra C\,} = \gen {a}$. If $a$ is \ndzz, the \ids $\gen {L}$ and~$\lra {\tra C\,}$ are \ivsz, \invs of each other, $\gen {L} \simeq \Im A$ and $\lra {\tra C\,} \simeq \Im \tra A$.

}
\end{exercise}

\rdb\vspace{-1em}
\begin{exercise}
\label{exoFittTfPtf}
{\rm If a \emph{\tfz} \Amo has its \idfs generated by \idmsz, it is \ptfz.
}
\end{exercise}

\rdb\vspace{-1em}
\begin{exercise}
\label{exoRelateursCourts}  (Short syzygies)\\
{\rm
\emph{Notations, terminology.}
Let $(e_1,\ldots,e_n)$ be the canonical bases of $\Ae n$.\\
Let $x_1$, \dots, $x_n$ be \elts of an \Amoz. Let $x=\tra[\,x_1\; \cdots\; x_n\,]$ {and $x^\perp := \Ker (\tra {x})\subseteq \Ae n$} the syzygy module 
between the~$x_i$'s.\\
We will say of a syzygy $z \in x^\perp$ that it is \gui{short} if it possesses at most two nonzero \coosz, \cad if $z \in \gA e_i \oplus \gA e_j$ $(1\leq i\neq j\leq n)$.
\begin{enumerate}\itemsep=0pt
\item 
Let $z \in x^\perp$. Show that the condition \gui{$z$ is a sum of short syzygies} 
is a \lin condition. Consequently, if $z$ is \gui{locally} a sum of short syzygies, 
it is also globally a sum of short syzygies.

\item 
Deduce that if $M=\sum \gA x_i$ is a \mlmz, then every \elt of $x^{\perp}$ is  a sum of short syzygies. 

\item 
If every syzygy 
between three \elts of $\gA$ is a sum of short syzygies, 
then $\gA$ is an \emph{\anarz},
\cad every \id $\gen{x,y}$ is \lopz.

\item 
In question \emph{2} give a global solution by using a \mlmo $A = (a_{ij}) \in \Mn(\gA)$ for $x$.
\end{enumerate}
}
\end{exercise}

\rdb\vspace{-1.2em}
\pagebreak	

\begin{exercise}
\label{exoPetitsRelateurs} (Trivial syzygies)\\
{\rm We use the notations of Exercise~\ref{exoRelateursCourts}. Now $x_1$, \ldots, $x_n \in \gA$. \\
For $z \in \Ae n$ let 
$\scp{z}{x} = \sum z_i x_i
$. The module of syzygies 
$x^\perp$ contains the \gui{trivial syzygies} 
$x_j e_i - x_i e_j$ (which are a special case of short syzygies).\\
In the two first questions, we show that if $x$ is \umdz, then $x^\perp$ is generated by these trivial syzygies. 
We fix $y \in \Ae n$ such that $\scp {x}{y} = 1$.
\begin{enumerate}\itemsep=0pt
\item 
Recall why $\Ae n = \gA y \oplus x^\perp$.

\item 
For $1 \leq i < j \leq n$, we define $\pi_{ij} \,: \,\Ae n \to \Ae n$ by

\snic{\pi_{ij}(z) = (z_i y_j - z_j y_i) (x_j e_i - x_i e_j),}

so  that $\Im \pi_{ij} \subseteq x^\perp \cap (\gA e_i \oplus \gA e_j)$. Show that $\pi = \sum_{i < j} \pi_{ij}$  is the projection over $x^\perp$ \paralm to $\gA y$. Deduce the result on the trivial syzygies. 
See also Exercise~\ref{exoPlgb1}.

\end{enumerate}
We no longer suppose that $x$ is \umdz. Let $M  
\in \Mn(\gA)$ be an alternating matrix.
\begin{enumerate} \setcounter{enumi}{2}\itemsep=0pt
\item 
Show that by letting $z = Mx$, we have $\scp {x}{z} = 0$. 
\item 
In which way is an alternating matrix a \gui{sum of small alternating matrices}? Make the link with the \dfn of $\pi_{ij}$ in question \emph{2}.

\end{enumerate}
}
\end{exercise}

\rdb\vspace{-1em}
\begin{exercise}
 \label{exoProjImLibre} (\Mprns which have a free image)\\
 {\rm  Let $P \in \GAn(\gA)$ be a projector whose image is free of rank $r$; by Proposition~\ref{propImProjLib} there exist $X \in \Ae {n \times r}$, $Y \in \Ae {r \times n}$ satisfying $YX = \I_r$ and $P = XY$.

\emph{1.}
Clarify  
the \dlg lemma (Lemma~\ref{propIsoIm}), in other words compute $A \in \SL_{n+r}(\gA)$ (and its inverse) such that
$$A^{-1} \Diag(0_r, P) A = \I_{r,n+r}.\eqno(*)$$

\emph{2.}
Suppose that $X=\tra Y$ (so $P$ is symmetrical). \\
Verify that we can impose upon $A$ to be \gui{orthonormal}
\cad $\tra A = A^{-1}$. \\
Conversely, if $A \in \SL_{n+r}(\gA)$ is orthonormal and satisfies~$(*)$, then we can write $P = X\tra{X}$ with $X \in \Ae {n \times r}$ and~$\tra{X} X = \I_r$ (the matrix $P$ is therefore symmetrical).
} \end{exercise}

\rdb\vspace{-1em}
\begin{exercise}
 \label{exoStabLibRang1} (Stably free modules of rank $1$)\\
 {\rm Give direct \dem that every \stl module of rank $1$ is free (Proposition~\ref{propStabliblib}), by using the Binet-Cauchy formula (Exercise~\ref{exoBinetCauchy}).\\
Consider two matrices $R\in\Ae {(n-1)\times n}$ and $R'\in\Ae {n\times (n-1)}$ with $RR'=\I_{n-1}$. Show that $\Ker R$ is a free module.
Conclude the result.
}
\end{exercise}

\rdb\vspace{-1em}
\begin{exercise}
 \label{exoStablib} (\Umd vectors, modules $M$ satisfying $M\oplus\gA\simeq\Ae n$)\\ 
 {\rm 
Let $x$, $y \in \Ae n$ be two vectors and $A \in \Mn(\gA)$ a matrix with first column $x$. Construct a matrix $B \in \Mn(\gA)$ as follows: its first row is $\tra {y}$ and its $n-1$ last rows are the $n-1$ last rows of $\wi {A}$, the cotransposed matrix of $A$.
\begin{enumerate}\itemsep=0pt
\item Show that $\det(B) = \det(A)^{n-2} \scp {x}{y}$ and that the $n-1$ last rows of $B$  belong to $x^\perp := \ker \tra {x}$.
\end{enumerate}
From now on assume that $\scp {x}{y} = 1$. We then know that the two \stl modules $x^\perp$ and $y^\perp$ are duals of one another 
(Facts~\ref{factStablib} and~\ref{factStablibDual}); we detail this \prt in a matrix fashion in the case where $x^\perp$ is free.
\begin{enumerate}\itemsep=0pt
\setcounter{enumi}{1}
\item Recall why $\Ae n = \gA x \oplus y^\perp$ and $\Ae n = \gA y \oplus x^\perp$.

\item Suppose that $\gA x$ possesses a free direct complement in $\Ae n$. Show in a matrix fashion that the same holds for $\gA y$ by constructing an $n \times n$ \iv matrix \gui{adapted} to the \dcn $\Ae n = \gA y \oplus x^\perp$.
\end{enumerate}
} \end{exercise}

\rdb\vspace{-1em}
\begin{exercise}\label{exoMatLocSym}  {(Symmetric \mlpz)}
\\
{\rm
Let $(x_1, \ldots, x_n) \in \Ae n$ possess a \smq \mlp $A \in \Mn(\gA)$. 
Let $\fa = \gen {x_1, \ldots,x_n}$. By using \Egrf{eqpmlm} of Proposition~\ref{pmlm}, show that~$\fa^2$ is principal and specifically that:
$\fa^2 = \gen {x_1^2, \cdots, x_n^2} = \gen {x_1^2 + \cdots + x_n^2}$.
}
\end{exercise}

\vspace{-1em}
\begin {exercise}\label{exoPgcdPpcm} (Regarding $\gA/\fa \oplus \gA/\fb \simeq \gA/(\fa \cap \fb) \oplus \gA/(\fa + \fb)$)\\
{\rm See also Exercise~\ref{exoSECSci} and Corollary~\ref{corthAnar}.

 \emph{1.}
 Let $\fa$, $\fb$ be two \ids of $\gA$ satisfying $1 \in (\fa : \fb) + (\fb : \fa)$. Construct a $\theta \in \GL_2(\gA)$ which satisfies
$\theta(\fa \oplus \fb) = (\fa \cap \fb) \oplus (\fa + \fb)$.
Deduce that $\gA/\fa \oplus \gA/\fb$ is \isoc to~$\gA/(\fa \cap \fb) \oplus \gA/(\fa + \fb)$.

 \emph{2.}
Let $a$, $b \in \gA$, $\fa = \gen {a}$, $\fb = \gen {b}$. Suppose that there exists an $A \in \GL_2(\gA)$ such that $A\Cmatrix{2pt} {a\cr b\cr} = \Cmatrix{2pt} {*\cr0\cr}$. Show that $1 \in (\fb : \fa) + (\fa : \fb)$. Find explicit $d$ and $m$  such that~$\fa \cap \fb = \gen {m}$, $\fa + \fb = \gen {d}$, as well as a matrix \eqvc between~$\Diag (a,b)$ and~$\Diag (m,d)$.

 \emph{3.}
Let $a$, $b \in \gA$ with $a\in \gen{a^2}$. Show that $a$, $b$ satisfy the conditions of question~\emph{2.}

 \emph{4.}
Let $\fa$, $\fb$ be two \itfs such that $\fa + \fb$ is \lopz.  Show $1 \in (\fa : \fb) + (\fb : \fa)$, $\fa \cap \fb$ is \tf and $\fa\fb = (\fa \cap \fb)(\fa + \fb)$.

}

\end{exercise}

\vspace{-.2em} The following exercises bring forth some results on the \deterz, the \polcar and the \polfonz.

\rdb
\begin{exercise}
\label{exoPrecisionsDet1}
{\rm
Let $M$ be a \ptf \Amoz,~$e$ be an \idm of $\gA$, $f=1-e$ and $\varphi$ be an \endo of $M$.
It is clear that $M=eM\oplus fM$, so $eM$ and $fM$ are \ptfsz.
We also have $\varphi(eM)\subseteq eM$, and by letting $\varphi_e:eM\rightarrow eM$ be the \endo defined as such, prove that we have

\snic{\begin{array}{c}
\det(\varphi_e)=f+e\det(\varphi)   \quad {\rm and } \quad
\det(e\varphi)=r_0f+e\det(\varphi)     \\[1mm]
\rF{e\varphi}(X)=\rF{\varphi}(eX)=\rF{\varphi_e}(X)=
f+e\rF{\varphi}(X)     \\[1mm]
\rC{\varphi_e}(X)=f+e\rC{\varphi}(X)     \\[1mm]
\rR{eM}(X)=f+e\rR{M}(X)     \\[1mm]
\end{array}}

Furthermore, show that $e\,\det(\varphi)$ is the \deter of $\varphi_e$ as the \endo of the $\gA[1/e]$-module $eM$.
}
\end{exercise}

\rdb\vspace{-1em}

\pagebreak	

\begin{exercise}
 \label{exoPeticalculPolrang}
 {\rm Consider the \qf module $M=\bigoplus _{k\in\lrbn}(r_k\gA)^{k}$, where the $r_k$'s are \ort \idmsz. 
We have $M\simeq e_1\gA \oplus \cdots \oplus e_n\gA$ with $e_k=\sum_{j=k}^nr_j$, and $e_k \divi e_{k+1}$ for $k\in\lrb{1..n-1}$ (cf.\ Lemma~\ref{lem ide-div}, and Exercises~\ref{exoIdmsSupInf} and~\ref{exoSfio}).
\\
Let $r_0=1-\sum_{i=1}^nr_i$ and $s_k=1-r_k$.
\begin{itemize}
  \item [--] Recall why $\rR{r_k\gA}(X)=s_k+r_kX$.
  \item [--] Show that $\rR{M}(X)=r_0+r_1X+\cdots +r_nX^n = \prod\nolimits_{k=1}^n (s_k+r_kX)^k$.

  \item [--] Verify this \egt using a direct computation.
\end{itemize}
}
\end{exercise}

\rdb\vspace{-1em}
\begin{exercise}
\label{exoPrecisionsDet2} (The \deterz, component by component)\\
{\rm
Let $\varphi$ be an \endo of a \mptf $M$ having~$n$ \gtrsz. Let $r_h=\ide_h(M)$ (for $h\in\lrb{0..n}$) and $d=\det(\varphi)$. Denote by~$\varphi\ep{h}$ the \endo of the \Amo $M\ep{h}$ induced by $\varphi$, $d_h=r_hd$, $\delta_h=\det(\varphi\ep{h})$ and $s_h=1-r_h$.

\emph{1.} Show that we have the following \egts 

\snic{d_0=r_0 ,\; \delta_0=1,\; \delta_h=s_h+d_h \;  \hbox{ and }  \; d=d_0+d_1+\cdots+d_n=\delta_1 \times
\cdots \times \delta_n.}

\emph{2.} Furthermore, show that 
$d_h$ is the \deter of $\varphi\ep{h}$ in $\gA[1/{r_h}]$ when we regard~$\varphi\ep{h}$ as an \endo of the~$\gA[1/{r_h}]$-module $M\ep{h}$.

\emph{3.} Similarly, show that we have

\snic{\rF{\varphi\ep{h}}(X)=s_h+r_h\rF{\varphi}(X) \quad {\rm and }\quad
\rC{\varphi\ep{h}}(X)=s_h+
r_h\rC{\varphi}(X).}

}
\end{exercise}

\rdb\vspace{-1em}

\begin{exercise}
\label{exoPrecisionsDet3} (\Polcar and \polfon in the case of constant rank)
{\rm
Let $\varphi$ be an \endo of a module $M$ of constant rank $h$. 
Prove the following facts.\\
The \polcar of $\varphi$ is monic of degree $h$ and the \polfon of $\varphi$ is of degree $\leq h$. The homogenized \pols at degree $h$ of $\rC{\varphi}(X)$ and $\rF{\varphi}(X)$ are respectively equal to $\det(X\Id_M-Y\varphi)$ and $\det(Y\Id_M+X\varphi)$. In other words
we have the \egts

\snic{ \rC{\varphi}(X)=X^h\rF{\varphi}(-1/X)\quad {\rm and }\quad
\rF{\varphi}(X)=(-X)^h\rC{\varphi}(-1/X).}

Furthermore,  $\det(\varphi)=(-1)^h\rC{\varphi}(0)$ is equal to the \coe of $X^h$ in $\rF{\varphi}(X)$.
}
\end{exercise}

\rdb\vspace{-1em}
\begin{exercise}
\label{exoPrecisionsDet4} (\Polcar and \polfonz, \gnl case)\\
{\rm
Let $\varphi$ be an \endo of a \mptf $M$. Let

\snic{\rF{\varphi}(X)=1+v_1X+\cdots+v_nX^n$ and
$\rR{M}(X)=r_0+r_1X+\cdots+r_nX^n.}

Then, show that we have the following \egtsz.

\snic{\arraycolsep2pt\begin{array}{rcl}
  r_hv_k&=&0\;\hbox{ for }0\leq h< k\leq n,\\[1mm]
   \rC{\varphi}(X)   &=& r_0\,+\,\sum\nolimits_{1\leq h\leq n}
r_h\,X^h\,\rF{\varphi}(-1/X),   \\[1mm]
   \rF{\varphi}(-X)   &=&  r_0\,+\,\sum\nolimits_{1\leq h\leq n}
r_h\,X^h\,\rC{\varphi}(1/X),  \\[1mm]
    \det(\varphi -X\Id_M)  &=&   \rR{M}(-1)\, \rC{\varphi}(X), \\[1mm]
   \det(\varphi)   &=&   r_0+r_1v_1+\cdots+r_nv_n=\rR{M}(-1)\,\rC{\varphi}(0).
\end{array}}
}
\end{exercise}


\rdb\vspace{-1em}

\pagebreak	

\begin{problem}\label{exoSuslinCompletableLemma}
 {(Completion of \umd vectors: a result due to \Susz)}\\
{\rm  
A vector of $\Ae n$ is said to be \emph {completable} if it is equal to the first column of a matrix of~$\GLn(\gA)$. It is then \umdz. We want to show the following result.\index{completable!unimodular vector} 
 
\emph{Let $b\in\gA$ and $(\an)\in\Ae n$ such that $(\ov{a_1}, \ldots, \ov{a_n})$ is completable over $\gA\sur{b\gA}$, then $(\an, b^n)$ is completable (over $\gA$)}. 
 
By hypothesis, we have  $A$, $D \in \Mn(\gA)$ satisfying $A\, D \equiv \In \mod b$, with $\vab{a_1}{\cdots\;a_n}$ as the first row of $A$.
We want to find a matrix of~$\GL_{n+1}(\gA)$ whose first row is $\vab{a_1}{\cdots\;a_n\;b^n}$.  Let~$a = \det(A)$.

\emph {1.}
Show that there exists a $C\in \Mn(\gA)$ such that $\Cmatrix{2pt} {A &b\,\In\cr C &D\cr} \in \GL_{2n}(\gA)$. 

Now it is a matter of transforming the top-right corner $b\,\In$ of the above matrix into $B' := \Diag(b^n, 1,\dots,1)$.

\emph {2.}
Show that we can write $B' = bE + aF$ with $E \in \En(\gA)$ and $F \in \Mn(\gA)$.

\emph {3.}
Verify that $\Cmatrix{2pt} {A &b\,\In\cr C &D\cr} \Cmatrix{2pt} {\In &\wi AF\cr 0 &E\cr} =
\Cmatrix{2pt} {A &B'\cr C &D'\cr}$ with $D' \in \Mn(\gA)$.

\emph {4.}
Show that $\Cmatrix{2pt} {A &B'\cr C &D'\cr}$ is \eqve to a matrix $\Cmatrix{2pt} {A &B'\cr C &D''\cr}$ where $D''$ has \hbox{its $n-1$} last columns null. Deduce the existence of an
\iv matrix whose first row is $\vab{a_1}{\cdots\;a_n\;b^n}$.  

\emph {5.}
Example (Krusemeyer). If $(x,y,z) \in \Ae 3$ is \umdz, $(x,y,z^2)$ is completable. More \prmtz, if $ux + vy + wz = 1$, the matrix below is suitable.

\snic {
\Cmatrix{2pt} {x&y&z^2 \cr v^2 &w-uv &-x-2vz\cr -w-uv &u^2&-y+2uz}.
}

What is its \deter (independently from the fact that $ux+vy+wz=1$)?

\emph {6.}
More \gnltz, we have the following result (Suslin): if $(a_0, \an)$ is \umdz, then $(a_0, a_1, a_2^2, \ldots, a_n^n)$ is completable.

\emph {7.}
Show the following result (\Susz's $n!$ theorem): if $(a_0, \an)$ is \umdz, then for exponents $e_0$, $e_1$, \ldots, $e_n$ such that $n!$ divides $e_0\cdot e_1\cdots e_n$, the vector $(a_0^{e_0}, a_1^{e_1}, \ldots, a_n^{e_n})$ is completable.

}

\end{problem}

\rdb\vspace{-1em}
\begin{problem}\label{exoSphereCompletableYengui} {(The $n$-sphere when $-1$ is a sum of $n$ squares, with I.~Yengui)}\\
{\rm  
\emph {1.}
Let $\gA$ be a \ri in which $-1$ is a sum of $2$ squares and $x_0$, $x_1$, $x_2 \in \gA$ satisfying $x_0^2 + x_1^2 + x_2^2 = 1$.
\begin{enumerate}\itemsep=0pt
\item [\emph {a.}]
Show that the vector $(x_0,x_1,x_2)$ is completable by considering a matrix $M = \Cmatrix{2pt} {x_0 & u & a\cr x_1 & v & b\cr x_2 & 0 & c\cr}$ where $u$, $v$ are \lin forms in $x_0$, $x_1$, $x_2$ and $a$, $b$, $c$ are constants.

\item [\emph {b.}]
Give examples of \ris $\gA$ in which $-1$ is a sum of $2$ squares.
\end{enumerate}

\emph {2.}
Suppose that $-1$ is a sum of $n$ squares in the \ri $\gA$.
\begin{enumerate}\itemsep=0pt
\item [\emph {a.}] We use the notation \smash{$A\sims{\cG }B$} from \paref{NOTAfGg}.
Let $x_0$, $x_1$, \ldots, $x_n$ with $x_0^2 + \cdots + x_n^2 = 1$. Show that
$$\preskip-.4em \postskip.4em 
\tra {[\,x_0\;x_1\; \cdots\; x_n\,]} \sims{\EE_{n+1}} \tra {[\,1\; 0\; \cdots\; 0\,]}. 
$$
In particular,  $\tra {[\,x_0\;x_1\; \cdots\; x_n\,]}$ is completable.

\item [\emph {b.}]
Let $m \ge n$, $x_0, x_1, \ldots, x_m$ and $y_{n+1}, \ldots, y_m$ satisfy
\smashbot{$\sum\limits_{i=0}^n x_i^2 + \!\!\!\sum\limits_{j=n+1}^m \!\!y_jx_j = 1$}. 
Show that $\tra {[\,x_0\;x_1\; \cdots\; x_m\,]} \sims{\EE_{m+1}} \tra {[\,1\; 0\; \cdots\; 0\,]}$. 
\end{enumerate}

\emph {3.}
Suppose that there exists an $a \in \gA$ such that $1 + a^2$ is nilpotent. This is the case if $-1$ is a square in $\gA$, or if $2$ is nilpotent.
\begin{enumerate}\itemsep=0pt
\item [\emph {a.}]
Let $x_0$, $x_1 \in \gA$ with $x_0^2 + x_1^2 = 1$. Show that \smashtop{$\crmatrix {x_0 & -x_1\cr x_1 & x_0} \in \EE_2(\gA)$}.

\item [\emph {b.}]
Let $x_0$, $x_1$, \dots,$x_n$ and $y_2$, \ldots, $y_n$ in $\gA$ such that $x_0^2 + x_1^2 + \sum_{i=2}^n x_i y_i = 1$. Show that $\tra {[\,x_0\;x_1\; \cdots\; x_n\,]}
\sims{\EE_{n+1}} \tra {[\,1\; 0\; \cdots\; 0\,]}$.

\item [\emph {c.}]
Let $\gk$ be a \riz, $\gk[\uX,\uY] = \gk[X_0, \Xn, Y_2, \ldots, Y_n]$ and 

\snic{f=1 - \big(X_0^2 + X_1^2 + \sum_{i=2}^n X_i Y_i\big).}

Let $\gA_n = \gk[x_0, \xn, y_2, \ldots, y_n] = \aqo {\gk[\uX,\uY]}{f}$. Give examples for which, for all $n$, $\tra {[\,x_0\;x_1\; \cdots\; x_n\,]}$ is completable without $-1$ being a square in~$\gA_n$.
\end{enumerate}
}

\end{problem}


\sol

\exer{exoSchanuelVariation}
\emph{1.} By \recu on $n$, the $n=1$ case being exactly Schanuel's lemma (Corollary~\ref{corlemScha}). From each exact sequence, we construct another of length minus one
\snic {\arraycolsep4pt
\begin{array}{ccccccccccccc}
0 &\rightarrow& K& \to & P_{n-1}& \to & \cdots & \to &
P_1 \oplus P'_0 &   \smash{\vvvers{u \oplus \I_{{P_0}'}}}
& \Im u \oplus P'_0 & \to & 0
\\
0 &\rightarrow& K'& \to & P'_{n-1}& \to & \cdots & \to &
P_1 \oplus P'_0 &  \vvvers{u' \oplus \I_{P_0}}
& \Im u' \oplus P_0 & \to & 0
\end{array}
}

But we have $\Im u \oplus P'_0 \simeq \Im u' \oplus P_0$ by Schanuel's lemma applied to the two short exact sequences,

\snic {\arraycolsep4pt
\begin{array}{ccccccccc}
0 &\to& \Im u & \to & P_0 & \to & M & \to & 0 \\
0 &\to& \Im u' & \to & P'_0 & \to & M & \to & 0 \\
\end{array}
}

We can therefore apply the \recu (to the two long exact sequences of length $-1$), which gives the desired result. 

 \emph{2.} Immediate consequence of \emph{1.}

\exer{exoSuitExPtfs}
{
Let us show by \recu on $i$ that $\Im u_i$ is a \mptfz. This is true for $i = 1$. Suppose it is true for $i\geq 1$; we have therefore a surjective \ali $P_i \vers{u_i}  \Im u_i$ where $\Im u_i$ is \ptf and thus $P_i \simeq \Ker  u_i \oplus \Im u_i$. But $\Ker u_i = \Im u_{i+1}$ therefore $\Im u_{i+1}$ is \ptfz.
In addition $P_i \simeq \Im u_{i} \oplus \Im u_{i+1}$. Then
$$\preskip.3em \postskip.0em \arraycolsep4pt
\begin {array} {rcl}
P_1 \oplus P_3 \oplus P_5
\oplus \cdots&\simeq
& 
(\Im u_1 \oplus \Im u_2) \oplus (\Im u_3 \oplus \Im u_4) \oplus
 \cdots
\\
&\simeq & 
\! \Im u_1 \oplus (\Im u_2 \oplus \Im u_3) \oplus
(\Im u_4  \oplus \Im u_5) \oplus  \cdots
\\
&\simeq & 
\! P_2 \oplus P_4 \oplus P_6 \oplus \cdots
\end {array}
$$
}

\vspace{-.5em}
\pagebreak	


\exer{exoProjRang1}
Let $A_1$, \ldots, $A_n$ be the columns of $A$ and $t = \Tr A = \sum_i a_{ii}$.

\emph {1.}
Let us first check that $tA_j = A_j$:

\snic {
\hbox{by using }
\left| \matrix {a_{ii} & a_{ij} \cr a_{ki} & a_{kj} \cr} \right| = 0
\hbox { and } A^2 = A, \;
ta_{kj} = \sum_i a_{ii} a_{kj} =  \sum_i a_{ki} a_{ij} =  a_{kj}.
}

Therefore $(1-t)A = 0$, then $(1-t)t= 0$, \cad $t$ \idmz. In addition, if $a A=0$, \hbox{then $a t = 0$}, \cad $a = a(1-t)$.

\emph {2.}
On the localized \ri at $a_{ii}$, 
two arbitrary columns $A_j$, $A_k$ are multiples of~$A_i$ so $A_j \wedge A_k = 0$. Hence globally $A_j \wedge A_k = 0$, and so $\cD_2(A) = 0$. Moreover, by using $\left| \matrix {a_{ik} & a_{ij}\cr a_{kk} & a_{kj}\cr}\right| = 0$, we have $\sum_k a_{ik}a_{kj} = \sum_k a_{ij}a_{kk} = a_{ij} \Tr A = a_{ij}$, \cad $A^2 = A$. 

\emph {3.}
The \sys on the right-hand side is of cardinality $1 + n^2$, the one in the middle of cardinality $1 + {n \choose 2}^2$. To obtain the left-hand side one, we must count the minors that do not have a corner on the diagonal. 
Suppose $n \ge 3$, then there are ${n \choose 2} {n-2 \choose 2}$ minors,
 and ${n \choose 2}^2 - {n \choose 2} {n-2 \choose 2} = (2n-3){n \choose 2}$ remain, 
 hence the cardinality $1 + (2n-3){n \choose 2}$.
For $n = 3$, each \sys is of cardinality $10$.
For $n>3$,    $1 + n^2$ is strictly less than the other two.


\exer{exoMatProjRang1CoeffReg} 
\emph {1.}
We have $\dmatrix {a_{i\ell} & a_{ij}\cr a_{k\ell} & a_{kj}} = 0$, \cad $a_{kj} a_{i\ell} = a_{ij}a_{k\ell}$.\\
This is the \egt $C_j \cdot L_i = a_{ij} A$. As for $L_i \cdot C_j$, this is the \coe in
position $(i,j)$ of $A^2 = A$, \cad $a_{ij}$.

\emph {2.}
We have $L_i\cdot A = L_i$ so $L_i\cdot(\In-A) = 0$.  Conversely, for $u \in \gA^n$ such that $\scp {L_i}{u} = 0$, it must be shown that $u = (\In-A)(u)$, \cad $Au = 0$, \cad $\scp {L_k}{u} = 0$. But $a_{ij}L_k = a_{kj}L_i$ and as $a_{ij}$ is \ndzz, this is immediate.

\emph {3.} The \egt $a_{kj}a_{i\ell}=a_{ij}a_{k\ell}$ shows that $C\cdot L=a_{ij}A$. Moreover, if $\gA$ is with explicit \dvez, we can compute $A$ from $L$ and $C$. 
\\
If we take a row $L$ whose \coes are called $a_{i\ell}$ ($\ell\in\lrbn$) and a column $C$ whose \coes are called $a_{kj}$ ($k\in\lrbn$), with the common \elt $a_{ij}$ being \ndzz, the conditions are the following:
\begin{itemize}
\item each \coe of $C\cdot L$ must be divisible by $a_{ij}$, hence $A=\frac1{a_{ij}}C\cdot L$,
\item we must have $\Tr(A)=a_{ij}$, \cad $L\cdot C =a_{ij}$.
\end{itemize}
Naturally, these conditions are directly related to the invertibility of the \id generated by the \coes of $L$.

\emph {4.}
In the matrix \egt $C \cdot L = (L\cdot C)\,A$ to be proven, each side is bi\lin in $(L, C)$. However, the \egt is true if $\tra L$ is a column of $\tra A$ and $C$ a column of $A$, therefore it remains true for $\tra L \in \Im\tra A$ and $C \in \Im A$. The rest is easy.

\exer{exoFittTfPtf}
$M$ is the quotient module 
of a \mptf $P$ which share the same \idfsz.
If $P\oplus N=\Ae n$, $M\oplus N$ is a quotient of~$\Ae n$ with the same \idfsz. Therefore there is no nonzero syzygy 
between the \gtrs of $\Ae n$ in the quotient $M\oplus N$. 
Therefore 

\snic{M\oplus N=\Ae n \;\hbox{ and }\;P/M\simeq (P\oplus N)/(M\oplus N)=0.}

\exer{exoRelateursCourts}{
\emph{1.} 
A syzygy
$z=\sum z_k e_k$ is a sum of short syzygies 
\ssi there exist syzygies 
$z_{ij} \in \gA e_i \oplus \gA e_j$ such that $z = \sum_{i < j} z_{ij}$. This is interpreted as follows

\snic{ \exists\alpha_{ij}, \beta_{ij} \in
\gA,\; z_{ij} = \alpha_{ij} e_i + \beta_{ij} e_j, \;\scp {z_{ij}}
{x} = 0 \et z = \sum_{i < j} z_{ij}.}

This is equivalent to $z_k = \sum_{k < j} \alpha_{kj} + \sum_{i < k} \beta_{ik}$ ($k\in\lrbn$) and $\alpha_{ij} x_i + \beta_{ij} x_j = 0$ \hbox{(for $i<j$)}. This is indeed a \sli with the \gui{unknowns} $\alpha_{ij}$,~$\beta_{ij}$.

\emph{2.} 
By reasoning locally, we can assume that the $x_i$'s are multiples of $x_1$, which we write as $b_i x_1 + x_i = 0$. Hence the syzygies 
$r_i = b_i e_1 + e_i$ for $i\in\lrb{2..n}$.  
\\
Let $z \in x^\perp$. Let $y = z - (z_2r_2 + \cdots + z_nr_n)$, we have $y_i = 0$ for $i \ge 2$, and so~$y$ is a (very) short syzygy. 
Thus, $z = y + \sum_{i=2} z_i r_i$ is a sum of short syzygies.

\emph{3.} 
Let $x$, $y \in \gA$. We are looking for $s$, $t$ with $s+t = 1$,  $sx \in \gA y$ and $ty \in \gA x$. We write the syzygy 
$(-1,-1,1)$ between $(x, y, x+y)$ as a sum of short syzygies

\snic{
(-1,-1,1) = (0,a,a') + (b,0,b') + (c,c',0).
}

In particular, $a'+b'=1$,
 and the result follows.

\emph{4.} 
By \dfn $\sum_i a_{ii} = 1$ and 
$\left| 
\matrix {a_{ij} & a_{ik} \cr
x_{j} & x_{k}} \right| = 0$. This provides several short syzygies 
$a_{ij} e_k - a_{ik} e_j$.  We keep the $r_{ik} = a_{ii} e_k - a_{ik} e_i$, \cad those that correspond to a \gui{diagonal minor} $\left| \matrix {a_{ii} & a_{ik} \cr x_{i} & x_{k}}
\right|$. For $z \in \Ae n$, let

\snic{y = Az \et  z' = \sum_{i,k} z_k r_{ik} = \sum_{i,k} z_k (a_{ii} e_k - a_{ik} e_i).
}


Then $z=z'+y$: indeed, the \coe of $e_j$ in $z'$ is

\snic{
\big(\sum_i a_{ii}\big) z_j - \sum_{k} a_{jk} z_k = z_j - (Az)_j.
}

Since $A^2-A$ annihilates $M$, $z-y\in x^{\perp}$, so $z \in x^\perp \Rightarrow y  \in x^\perp$.
Each $y_ie_i$ is a (very) short syzygy 
since $y_ix_i = 0$.
Therefore  $z = z' + y=z'+\sum y_i e_i$ is a sum of short syzygies.
}

\exer{exoPetitsRelateurs}
\emph{1.}
We write $z \in \Ae n$ in the form

\snic{z = \scp {x}{z}.y +
(z - \scp {x}{z}.y),}

which provides the \dcn
$\Ae n = \gA y \oplus x^\perp$.

\emph{2.}
For $i \leq j$, define $z_{ij} \in \gA$ by $z_{ij} = z_i y_j - z_j
y_i$ and let
$$\ndsp
z' = \som_{i < j} z_{ij} (x_j e_i - x_i e_j) =
\som_{i \leq j} z_{ij} (x_j e_i - x_i e_j)
.$$
For fixed $k$, the \coe of $e_k$ in the right-hand sum is
\[
\begin{array}{r}
\sum_{j \ge k} z_{kj} x_j - \sum_{i < k} z_{ik} x_i  \;  =  \;
\sum_{j \ge k} (z_ky_j - z_jy_k) x_j - \sum_{i < k} (z_iy_k - z_ky_i) x_i
 \\[2mm]
    =  \; 
  z_k \sum_{j = 1}^n y_jx_j - y_k \sum_{j = 1}^n z_jx_j\;=\;
   z_k - \scp {z}{x} y_k\,,
  \end{array}
\]
which means that $z' = z - \scp {z}{x} y$ and proves the required result.

\emph{3.} If $\psi$ is the alternating bi\lin form associated with the matrix, the \egt $\scp {x}{z} = 0$ simply means that $\psi(x,x)=0$.

\emph{4.}   
We can express an alternating $n \times n$ matrix as a sum of $n(n-1) \over 2$ small alternating matrices. 
For $n = 3$, here is the alternating matrix allowing us to make the connection with question~\emph{2} ($y$ is fixed and it is $z$ that varies):
$$
M_z = \Cmatrix{2pt} {
0                & z_1y_2 - z_2y_1  & z_1y_3-z_3y_1 \cr
-z_1y_2 + z_2y_1 & 0                & z_2y_3-z_3y_2 \cr
-z_1y_3 + z_3y_1 & -z_2y_3 + z_3y_2 & 0 }.
$$
The \dcn of $M_z$ into small alternating matrices provides the $\pi_{ij}$'s. It must be noted that $z\mapsto M_z,\;\Ae n \to \Mn(\gA)$ is a \ali and that $\pi(z)= M_z x$.

\exer{exoProjImLibre} \emph{1.} We follow the method of this course. It leads to letting
$$
A = \cmatrix {0_r & -Y\cr X & \I_n - P\cr},
\quad 
A' = \cmatrix {0_r & Y\cr -X & \I_n - P\cr}.
$$
These matrices satisfy
$$
A\Cmatrix{2pt} {\I_r & 0\cr 0 & 0_n\cr} = \Cmatrix{2pt} {0_r & 0\cr 0 & P\cr} A,
\
AA' = \I_{n+r}.
$$
\emph{2.} Immediate since we have the formulas right in front of us.

\exer{exoStabLibRang1}{
The Binet-Cauchy formula gives $1=\det(RR')=\sum \delta_i \delta'_i$
with 

\snic{\delta_i=\det(R_{1..n-1,1..n\setminus i})$
and $\delta'_i=\det(R'_{1..n\setminus i,1..n-1}).}

Let $S=[\,\delta'_1\;-\delta'_2\;\cdots\;(-1)^{n-1}\delta'_n\,]$.
Check that the square matrix \smash{$A=\Cmatrix{2pt}{S\cr R}$} has \deter $1$.
This shows that $\Ker R$ is free (Proposition~\ref{propStabliblib}).
\\
Actually let $S'=\tra{[\,\delta_1\;-\delta_2\;\cdots\;(-1)^{n-1}\delta_n\,]}$ and $A'=[S'\,R']$. Then $AA'=\In$, and this shows that $S'\in\Ae n$
is a basis of $\Ker R$.
}


\exer{exoStablib}{
\emph{1.} Consider the matrix $B\,A$. By \dfn of $B$, $B\,A$ is upper triangular, with diagonal $(\scp {x}{y}, \delta, \ldots, \delta)$ where $\delta = \det(A)$. By taking the determinant, we obtain $\det(B) \det(A) = \scp {x}{y} \,\delta^{n-1}$. The announced \idas  are therefore true when $\delta$ is \ivz. Since we are dealing with \idasz, they are always true.
The second item of the question is immediate.

 \emph{2.} 
Write $z \in \Ae n$ in the form $z = \scp {y}{z}\,x + (z - \scp {y}{z}\,x)$, which gives us the \dcn $\Ae n = \gA x \oplus y^\perp$.

 \emph{3.} 
The hypothesis boils down to saying that $x$ is the first column of an \iv matrix~$A$. Therefore $y$ is the first row of the \iv matrix $B$ above. The matrix~$\tra {B}$ is adapted to the \dcn $\Ae n = \gA y \oplus x^\perp$.
}

\exer{exoMatLocSym}
Let $x= [\,x_1 \;\cdots\;x_n \,]$. For $\alpha = \tra {[\alpha_1, \ldots, \alpha_n]}$ and $\beta = A\alpha$,
the \egt in question is

\snic {
\beta \,x=\Cmatrix{2pt}{\beta_1\cr \vdots\cr \beta_n}\Cmatrix{2pt}{x_1 \;\cdots \; x_n}=
{x}\,{\alpha}\,A\;\;$ with $\;\;{x}\,{\beta}=x \,\alpha.
}

Take $\alpha_i = x_i$. Since $x A = x$ and $A$ is \smqz, we obtain $A \tra {x} = \tra {x}$, \hbox{\cad $\beta = \tra {x}$}. Hence $\tra {x}\, x = {x}\tra{x}\,A = (x_1^2 + \cdots + x_n^2)A$. \\
Finally, $x_ix_j \in \gen {x_1^2 + \cdots + x_n^2}$ ($i,j\in\lrbn$).

\exer{exoPgcdPpcm}{
\emph{1.}
Let $\alpha \in (\fb : \fa)$ and $\beta \in (\fa : \fb)$ satisfy $1 = \alpha + \beta$. 
Then, the matrix $\theta = \crmatrix {\alpha & \beta\cr -1 & 1\cr}$, with \deter $1$ and with inverse $\theta^{-1} = \crmatrix {1 & -\beta\cr 1 & \alpha\cr}$, is suitable. Indeed

\snic{\crmatrix {\alpha & \beta\cr -1 & 1\cr} \Cmatrix{2pt} {\fa\cr \fb\cr} \subseteq
\Cmatrix{2pt} {\fa\cap \fb\cr \fa+\fb\cr}
\quad \hbox{and} \quad
\crmatrix {1 & -\beta\cr 1 & \alpha\cr} \Cmatrix{2pt} {\fa\cap \fb\cr \fa+\fb\cr}
\subseteq \Cmatrix{2pt} {\fa\cr \fb\cr}.}

On the left-hand side, the upper inclusion comes from the fact that $\alpha \fa + \beta \fb \subseteq \fa \cap \fb$, and the lower one is trivial. 
On the right-hand side, the upper inclusion comes from the fact that \hbox{$\fa \cap \fb + \beta(\fa + \fb) \subseteq \fa$}, and the lower one comes from the fact that $\fa \cap \fb + \alpha (\fa + \fb) \subseteq \fb$. Recap: we have \hbox{$\theta(\fa \oplus \fb) = (\fa \cap \fb) \oplus (\fa + \fb)$} with \halfsmashbot{$\theta = \crmatrix {1 & \beta\cr 0 & 1\cr} \crmatrix {1 & 0\cr -1 & 1\cr} \in \EE_2(\gA)$.}

 \emph{2.}
We can take $A$ of the form $A = \Cmatrix{2pt} {u & v\cr -b' & a'}$ with $ua' + vb' = 1$ and $a'b = b'a$. Let $m = a'b = b'a$, $d = ua + vb$; by inverting $A\Cmatrix{2pt} {a\cr b} = \Cmatrix{2pt} {d\cr 0}$, we obtain $a = da'$ and $b = db'$. It is clear that $\fa \cap \fb = \gen {m}$ and $\fa + \fb = \gen {d}$. We have $a' \in (\fa : \fb)$ and $b' \in (\fb : \fa)$.  Therefore $1 = \alpha + \beta$ with $\alpha = vb' \in (\fb : \fa)$, $\beta = ua' \in (\fa : \fb)$. To explicit a matrix \eqvcz, it suffices to use a matrix $\theta$ from the previous question:
$$\preskip.4em \postskip.4em 
\theta \Cmatrix{2pt} {a\cr 0} = \Cmatrix{2pt} {vm\cr -a} =
v\Cmatrix{2pt} {m\cr 0} - a'\Cmatrix{2pt} {0\cr d},
\;\;
\theta \Cmatrix{2pt} {0\cr b} = \Cmatrix{2pt} {um\cr b} =
u\Cmatrix{2pt} {m\cr 0} + b'\Cmatrix{2pt} {0\cr d}
. 
$$
%
%
Hence the matrix \eqvcz: $\theta \Cmatrix{2pt} {a & 0\cr 0 &b\cr} = \Cmatrix{2pt} {m & 0\cr 0 &d\cr}\Cmatrix{2pt} {v & u\cr -a' &b'\cr}$.

 \emph{3.}
The hypothesis is $a = a^2x$ for some given $x$. Then, the \elt $e = ax$ is \idm and $\gen{a}=\gen{e}$. We must solve $a'b = b'a$, $1 = ua' + vb'$, which is a \sli in $a', b', u, v$. \\
Modulo $1 - e$, we have $ax = 1$. We take $a' = a$, $b' = b$, $u = x$, $v = 0$. \\
Modulo $e$, we have $a = 0$. We take $a' = a$, $b' = 1$, $u = 0$, $v = 1$.\\ 
Therefore globally

\snic{a' = a, \quad b' = axb + (1-ax)1 = 1 - ax + axb, \quad
u = ax^2, \quad v = 1-ax.}

\sni \emph{4.}
Let $\fa = \langle x_1, \ldots, x_n \rangle$ and $\fb = \langle y_1, \ldots, y_m \rangle$.\\
We write $\fa + \fb = \langle z_1, \ldots, z_{n+m}\rangle$ with $z_1 = x_1$, \dots, $z_{n+m} = y_m$.  Let  $s_1, \ldots, s_{n+m}$ be \com such that over $\gA_{s_i}$ we have $\fa +\fb= \gen{z_i}$. \\
In each localized \ri we have 
 $\fa\subseteq \fb$ or $\fb\subseteq \fa$, hence $\so{\fa+\fb, \fa\cap\fb} =\so{\fa,\fb}$, and 

\snic{1 \in (\fa : \fb) + (\fb : \fa)$,  $\,\fa\fb = (\fa\cap\fb)(\fa \oplus\fb)\,$ and  $\,\fa \cap \fb$ is \tfz.$ }

}

\exer{exoPrecisionsDet1} {
We use the notations of Lemma~\ref{lem calculs}.\\
Let us take a look at the \deters of $e\varphi$ and $\varphi_e$.
We have 
%
$$\preskip.4em \postskip.2em\mathrigid1mu
\det(\varphi)=\det(\I_n-F+H),\, \det(e\varphi)=\det(\I_n-F+eH)
\,\,\hbox{and}\,\,\det(\varphi_e)=\det(\I_n-eF+eH).
$$
%
We deduce that 
\[\preskip.0em \postskip.4em 
\begin{array}{ccc} 
e\det(\varphi_e)=\det(e\I_n-eF+eH)=e\det(\varphi) \hbox{ and } \\[.3em] 
f\det(\varphi_e)=\det(f\I_n-feF+feH)= \det(f\I_n)=f. \end{array}
\] 
Therefore  $\det(\varphi_e)=f\det(\varphi_e)+ e\det(\varphi_e)
=f+e\det(\varphi)$.
\\
Similarly $e\det(e\varphi)=\det(e\I_n-eF+eH)=e\det(\varphi)$ and 
$$\preskip.3em \postskip.3em 
f\det(e\varphi)=\det(f\I_n-fF+feH)=f\det(\I_n-F)=f\rR{M}(0)=f\ide_0(M). 
$$
By applying $\det(\varphi_e)=f+e\det(\varphi)$ to the \endos $\Id+X\varphi$, $X\Id-\varphi$ and $X\Id$ of the $\AX$-module $M[X]$ we obtain
$$\preskip.3em \postskip.3em 
 \rF{\varphi_e}(X)= f+e\rF{\varphi}(X),\,\rC{\varphi_e}(X)=f+e\rC{\varphi}(X) \hbox{ and }\rR{eM}(X)=f+e\rR{M}(X).
$$Moreover, the matrix $eH$ simultaneously represents the \endo $e\varphi$ of $M$ and the \endo
$\varphi_e$ of $eM$. 
We therefore have $\rF{\varphi_e}(X)=\rF{e\varphi}(X)=\det(\I_n+eXH) = \rF{\varphi}(eX)$. 

As far as the last assertion is concerned: we must look at $\det(\varphi_e)$ in $\aqo{\gA}{f}$, we obtain $e\det(\varphi)$ modulo $f\gA$, and this corresponds to the \elt $e\det(\varphi)$ of $e\gA$.
}

\exer{exoPrecisionsDet2}
{\emph{1.} We have $\varphi\ep{h}=\varphi_{r_h}$ by applying the notation of Exercise~\ref{exoPrecisionsDet1}. \\
Therefore $\delta_h=s_h+d_h$. We have $\delta_0=1$ because $M\ep{0}=\so{0}$, and since $\delta_0=s_0+d_0$, this gives $d_0=r_0$.
\\
The \egt $d=d_0+d_1+\cdots+d_n$ is trivial.
\\
The \egt $d=\delta_1 \times \cdots \times \delta_n$ results from item~\emph{3} of \thref{propdef det ptf}. We can also prove $d_0+d_1+\cdots+d_n=\delta_1 \times \cdots \times
\delta_n$ by a direct computation.

 \emph{2} and \emph{3.} 
Already seen in Exercise~\ref{exoPeticalculPolrang}.
}

\exer{exoPrecisionsDet3}{
Recall: for $a\in \gA$ we have $\det(a\varphi)=\rR{M}(a)\det(\varphi)=a^h\det(\varphi)$.\\
We then place ourselves on the \ri $\gA[X,1/X]$ and consider the module $M[X,1/X]$. We obtain
$$X^h\rF{\varphi}(-1/X)= \det \big(X(\Id_M-(1/X)\varphi)\big)=\det(X\Id_M-\varphi)= \rC{\varphi}(X).$$
By replacing $X$ by $-1/X$ in $\rC{\varphi}(X)=X^h\rF{\varphi}(-1/X)$ we obtain the other \egtz. The two \pols are therefore of degrees $\leq h$. As the constant \coe of $\rF{\varphi}$ is equal to $1$, we also obtain that $\rC{\varphi}$ is monic.
\\
For the homogenized
 \polsz, 
the same computation works.\\
For the \deter we notice that $\det(-\varphi)=\rC{\varphi}(0)$.
}

\exer{exoPrecisionsDet4}{
We work over the \ri $\gA_{r_h}$ and we consider the module $r_hM$ and the \endo $\varphi\ep{h}$. We obtain a module of constant rank $h$. Therefore $r_h\rF{\varphi}(X)$ \hbox{and $r_h\rC{\varphi}(X)$} are of degrees $\leq h$, and $r_h \big(X^h\rF{\varphi}(-1/X)\big)=r_h\rC{\varphi}(X)$. It remains to sum up the \egts obtained in this way for ${ h\in\lrbn}$.\\
We perform the same computation for the second \egtz. The last two \egts were already known, except for $\det(\varphi)= r_0+r_1v_1+\cdots+r_nv_n$ which can be proven in the same way as the first.
}


\prob{exoSuslinCompletableLemma}
\emph {1.}
Let $C,\, U \in \Mn(\gA)$ such that $AD = \In + bU$, $DA = \In + bC$.
Then 
$$\preskip-.2em \postskip.4em 
 \Cmatrix{2pt} {A &b\In\cr C &D\cr} \Cmatrix{2pt} {D &-b\In\cr -U &A\cr} =
\Cmatrix{2pt} {\In &0\cr * &\In\cr} \in \GL_{2n}(\gA).
$$

\emph {2.}
We work modulo $a$ by noticing that $b$ is \iv modulo $a$. We can therefore, over $\gA\sur{a\gA}$, consider $b^{-1}B'$: it is a diagonal matrix of \deterz~$1$, therefore it belongs to $\En(\gA\sur{a\gA})$ (cf.\ Exercise~\ref{exoFacileGrpElem2}), we lift it to 
a matrix $E \in \En(\gA)$ and we obtain $B' \equiv bE \bmod a$.

\emph {3.}
Immediate.

\emph {4.}
It suffices to use the submatrix $\I_{n-1}$ which occurs in $B'$ to kill the \coes of the last $n-1$ columns of $D'$. The square submatrix of order $n+1$ obtained from $\Cmatrix{2pt} {A &B'\cr C &D''\cr}$ by deleting the rows $2$ to $n$ and the last $n-1$ columns is \iv with its first row being $[\,a_1\;\cdots\;a_n\; b^n\,]$.

\emph {5.}
Modulo $z$, the vector $\vab x y$ is completable in \smashtop{$A := \crmatrix {x &y\cr -v & u\cr}$}. \\
We have $\det(A)= a := ux + vy\equiv1 \mod z$ and we can take $D=\wi A$. \\
We write $DA = a\,\I_2 = \I_2 - wz\I_2$, so $C = -w\I_2$. 
The matrix \smashbot{$\Cmatrix{2pt} {A &z\I_2\cr C &D\cr}$} has \deter $(a + wz)^2$. To find $E$, we use the \egt 

\snic {
\Cmatrix{2pt} {z &0\cr 0&z^{-1}} = \rE_{21}(-1) \rE_{12}(1-z^{-1}) 
\rE_{21}(z) \rE_{12} \big(z^{-1}(z^{-1}-1)\big)
}

and the fact that modulo $a$, $zw\equiv 1$.
The author of the exercise has obtained a matrix $G$ that is more complicated than that of Krusemeyer. With $p = (y+u)w - u$, $q = (x-v)w + v$,

\snic {
G = \Cmatrix{2pt} {x&y&z^2\cr p(w-1)v-w & -p(w-1)u & y+u(z+1)\cr
            -q(w-1)v & q(w-1)u-w & -x+v(z+1)\cr}.
}

We have $\det(G)=1 + (xu+yv+zw-1) (wz+1) (yq-xp+1)$ whereas Krusemeyer's matrix has \deter $(ux + vy + wz)^2$!

\emph {6.}
Immediate by \recuz.


\prob{exoSphereCompletableYengui} 
\emph {1a.}
We have $\det(M) = -(cx_1-bx_2)u + (cx_0-ax_2)v$.\\
With 
$u = -(cx_1+bx_2)$, $v = cx_0+ax_2$, we obtain $\det(M) = cx_0^2 + cx_1^2 - (a^2 + b^2)x_2^2$. It suffices to take $c = 1$ and $a$, $b \in \gA$ such that $-1 = a^2 + b^2$.

\emph {1b.}
Let us show that $-1$ is a sum of two squares if $\gA$ contains a finite field. We can assume that $\gA$ is a field of odd cardinality $q$. 
\\
Consider the sets $A = \sotQ {a^2}{a \in \gA}$ and $B = \sotQ {-1-b^2}{b \in \gA}$. 
\\
They have $(q+1)/2$ \eltsz, so $A \cap B \ne \emptyset$, which gives the result.

Now here is a more \gnl result. If $n \not\equiv 0 \mod 4$, then $-1$ is a sum of two squares in $\ZZ/n\ZZ$. The hypothesis can be written as $\gcd(n, 4) = 1,2$ so $2 \in n\ZZ + 4\ZZ$, $2 = nu + 4v$. Let $m = -1 + nu = -4v + 1$; since $\gcd(4n, m) = 1$, the arithmetic progression $4n\NN + m$ contains a prime number~$p$ (Dirichlet), which satisfies 
$p \equiv m \equiv -1 \mod n$ and $p \equiv m \equiv 1 \mod 4$. \\
By this last congruence, $p$ is a sum of two squares, $p = a^2 + b^2$, so $-1 = a^2 + b^2$ in $\ZZ/n\ZZ$.

We deduce that if $n.1_\gA = 0$ with $n \not\equiv 0 \mod 4$ (this is the case if $n$ is a prime number), then $-1$ is a sum of two squares in $\gA$.

\emph {2a.}
Let $\an$ such that $-1 = \sum_{i=1}^n a_i^2$. We will use

\snic {
\sum_{i=1}^n (x_i-a_ix_0) (x_i + a_ix_0) = 
\sum_{i=1}^n x_i^2 -x_0^2 \sum_{i=1}^n a_i^2 = 1.
}

\vspace{-.1em}
\pagebreak

We have
$$\preskip.2em \postskip.4em 
\Cmatrix{2pt} {x_0\cr x_1\cr \vdots\cr x_n\cr}  \sims{\EE_{n+1}}
\Cmatrix{2pt} {x_0\cr x_1 + a_1x_0\cr \vdots\cr x_n + a_nx_0\cr}  \sims{\EE_{n+1}}
\Cmatrix{2pt} {x_0 + h\cr x_1 + a_1x_0\cr \vdots\cr x_n + a_nx_0\cr}
\quad
\begin{array}{c}
\hbox{with} \\[1em]
h = \sum_{i=1}^n \lambda_i (x_i+a_ix_0)\\
\end{array} 
$$
By taking $\lambda_i = (1-x_0)(x_i - a_ix_0)$, we obtain $h = 1-x_0$ so
$x_0 + h = 1$. 
\\
It is then clear that 

\snic{\tra {[1, x_1+a_1x_0, \ldots,
x_n+a_nx_0]} \sims{\EE_{n+1}} \tra {[1, 0, \ldots, 0]}.}

 Explicitly, by numbering the $n+1$ rows from $0$ to $n$ (instead of $1$ to $n+1$) and by letting

\snuc {
\begin{array} {c}
N = \prod_{i=1}^n \rE_{i,0} \big(-(x_i + a_ix_0)\big)\
\prod_{i=1}^n \rE_{0,i} \big((1-x_0)(x_i - a_ix_0)\big)\  \prod_{i=1}^n \rE_{i,0}(a_i) 
\\[1mm]
M = N^{-1} = \prod_{i=1}^n \rE_{i,0}(-a_i) \
\prod_{i=1}^n \rE_{0,i} \big((x_0-1)(x_i - a_ix_0)\big)\ \prod_{i=1}^n \rE_{i,0}(x_i + a_ix_0),
\end{array}
}

we obtain a matrix $M \in \EE_{n+1}(\gA)$ with first column $\tra {[\,x_0\; \cdots\;x_n\,]}$.

\emph {2b.}
We use $\gB = \aqo {\gA}{x_{n+1}, \ldots, x_m}$. 
The morphisms $\EE_r(\gA) \twoheadrightarrow \EE_r(\gB)$ are surjective. We first obtain
$$\preskip.0em \postskip.4em 
\tra {[x_0, \ldots, x_n]} \sims{\EE_{n+1}(\gB)}  \tra {[1, 0, \ldots, 0]}, 
$$
so some $x'_0, \ldots, x'_n \in \gA$ with in particular $x'_0 \equiv 1 \mod \gen {x_{n+1}, \ldots, x_m}$ such that
$$\preskip.2em \postskip.4em 
\tra {[x_0, \ldots, x_n, x_{n+1}, \ldots, x_m]} \sims{\EE_{m+1}(\gA)} 
\tra {[x'_0, \ldots, x'_n, x_{n+1}, \ldots, x_m].} 
$$
We easily deduce that

\vspace{-.1em}
\snuc {
\tra {[x'_0, \ldots, x'_n, x_{n+1}, \ldots, x_m]} \sims{\EE_{m+1}(\gA)} 
\tra {[1, \ldots, x'_n, x_{n+1}, \ldots, x_m]} \sims{\EE_{m+1}(\gA)} 
\tra {[1, 0, \ldots, 0].}
}

\sni
\emph {3a.}
We have $\crmatrix {x_0 & -x_1\cr x_1 & x_0\cr} \sims{\EE_{2}(\gA)} B =  \Cmatrix{2pt} {x_0 + ax_1& -x_1 + ax_0\cr x_1 & x_0\cr}$.  By using the fact that $1+a^2$ is nilpotent, we see that $x_0 + ax_1$ is \iv because

\snic {
(x_0 + ax_1)(x_0 - ax_1) = x_0^2 + x_1^2 - (1+a^2)x_1^2 = 1 -(1+a^2)x_1^2.
}

The matrix $B \in \SL_2(\gA)$ has an \iv \coe so it is in $\EE_2(\gA)$.

\emph {3b.}
First reason modulo $\gen {x_2, \ldots, x_n}$, then as in question \emph {2b}.

\emph {3c.}
We can take $\gk = \ZZ/2^e\ZZ$ with $e \ge 2$: $-1$ is not a square in $\gk$. \hbox{So $-1$} is not a square in $\gA_n$ either since there are morphisms $\gA_n \to \gk$, for example the \evn morphism at $x_0 = 1$, $x_i = 0$ for $i \ge 1$, $y_j = 0$ \hbox{for $j \ge 2$}.



\Biblio

Regarding \thref{prop Fitt ptf 2} and the \carn of  \mptfs by their \idfs see \cite{Nor} \Tho 18 p.~122 and Exercise 7 p.~49.
\perso{l'exercice 7 page 49 \'etablit l'\eqvc entre conoyau \pro et
application \lnlz.}
Note however that the \dem given by Northcott is not entirely \covz, since he requires an \prca of the \mptfsz.

We have defined the \deter of an \endo of a \mptf as in \cite[Goldman]{Gold}. The difference resides in the fact that our \dems are \covsz.

A study on the feasability of the local structure \tho for \mptfs can be found in \cite[D\'{\i}az-Toca\&Lombardi]{DiL09}.

Proposition~\ref{fact.homom loc pf} regarding $\big(\Lin_\gA(M,N)\big)_S$ is a crucial result that can be found for instance in \cite{Nor}, Exercise~9 p.~50, and in \cite{Kun} (Chapter~IV, Proposition~1.10). This result will be \gne in Proposition~\ref{propPlateHom}.

\Pbmz~\ref{exoSuslinCompletableLemma} is due to \Sus \cite{Sus77b}.

\newpage 
\thispagestyle{CMcadreseul}
\incrementeexosetprob

\chapter{\Stf \algs and \aGsz}
\label{chap AlgStricFi}\relax
\label{chap AlgEtales}\relax
\perso{compil\'e le \today}
\minitoc

\Intro

The chapter is devoted to a natural \gnn for commutative \ris of the notion of a finite \alg over a field. 
In \comaz, to obtain the conclusions for the field case, it is often \ncr to not only assume that the \alg is a \tf \evcz, but more \prmt that the field is discrete and that we know a basis of the \evcz. This is what brought us to introduce the notion of a \stf \alg over a \cdiz.

The pertinent \gnn of this notion to  commutative \ris is given by the \algs which are \mptfs over the base \riz. So we call them \stf \algsz.

Sections~\ref{secEtaleSurCD} and~\ref{sec2GaloisElr} which only concern \algs over \cdis can be read directly after Section~\ref{secGaloisElr}. Similarly for Section~\ref{secAGTG} if we start from a \cdi (certain \dems are then simplified).

Section~\ref{sec1Apf} is a brief introduction to \apfsz, by insisting on the  case of \algs which are integral over the base ring. 

The rest of the chapter is devoted to the \stf \algs themselves.

In Sections~\ref{secAlgSte} and \ref{secAlgSpb}, we introduce the neighboring notions of \ste \alg and of \spb \algz, which \gns the notion of an \'etale \alg over a \cdiz.

In Section~\ref{secAGTG} we give a \cof presentation of the bases of the \aG theory for commutative \risz. This is in fact an Artin-Galois theory, since it uses the approach developed by Artin for the field case by starting directly from a finite group of \autos of a field, the base field only appearing as a subproduct of the constructions that ensue.

\section{\'Etale \algs over a \cdiz}
\label{secEtaleSurCD}

\Grandcadre{In Sections~\ref{secEtaleSurCD} and \ref{sec2GaloisElr},~$\gK$ designates a nontrivial \cdiz}

Recall that a \Klgz~$\gB$ is said to be finite (resp.\,\stfez) if it is \tf as a \Kev (resp.\,if~$\gB$ is a finite dimensional \Kevz).
If~$\gB$ is a finite \Klgz, this does not imply that we know how to determine a basis of~$\gB$ as a \Kevz, nor that $\gB$ is discrete. If it is \stfez, however, then we know of  
 a finite basis of~$\gB$ as a \Kevz. In this case, for some~$x\in\gB$, the trace, the norm, the \polcar of (multiplication by) $x$, as well as the \polmin of $x$ over~$\gK$ can each be computed 
         using standard methods of \lin \alg 
over a \cdiz. Similarly every finite \Kslg of~$\gB$ is \stfe and the intersection of two \stfes sub\algs is \stfz.

\begin{definition}\label{defi1Etale}
\index{etale@\'etale!algebra over a discrete field}%
\index{algebra!etale@\'etale --- over a discrete field}%
\index{separable algebraic!algebra over a discrete field}%
\index{separable algebraic!element over a discrete field}%
\index{algebra!separable algebraic over a discrete field}%
Let~$\gL$ be a \cdi and~$\gA$ an \Llgz.
\begin{enumerate}\itemsep0pt
\item The \algz~$\gA$ is said to be \emph{\'etale (over~$\gL$)} if it is \stfe and if the \discri $\Disc_{\gA/\gL}$ is \ivz.
\item An \elt of~$\gA$ is said to be \emph{\agsp (over~$\gL$)} if it annihilates a \spl \polz.
\item The \algz~$\gA$ is said to be \emph{\agsp (over~$\gL$)} if every \elt of~$\gA$ is \agsp over~$\gL$.
\end{enumerate}
\end{definition}

When~$f$ is a \polu of~$\gL[X]$, the quotient \algz~$\aqo{\gL[X]}{f}$ is \'etale \ssi $f$ is \spl (Proposition \ref{propdiscTra}).

\subsec{Structure \tho for \'etale \algsz}

\vspace{3pt}
Proposition~\ref{propdiscTra} gives the following lemma.

\pagebreak	

\begin{lemma}\label{lem1EtaleCD}
Let~$\gA$ be a \stf \Klg and $a\in\gA$. If the \polcar $\rC{\gA/\gK}(a)(T)$ is \splz, then the \alg is \'etale and~$\gA=\gK[a]$.
\end{lemma}

In Fact~\ref{fact1Etale}, items~\emph{\ref{i1fact1Etale}} and~\emph{\ref{i2fact1Etale}} give more precise statements for certain items of Lemma~\ref{lemZerRed} and of Fact~\ref{factZerRedConnexe} (concerning \gnl reduced \zed \risz), in the case of a reduced \stf \Klgz.

\Gnl results on integral extensions of \zed \ris are given in Section~\ref{sec1Apf} from \paref{EntSurZdim} onwards. 

\begin{fact}\label{fact1Etale} Let~$\gB\supseteq\gK$ be a \stf \algz.
\begin{enumerate}
\item \label{i1fact1Etale} The \algz~$\gB$ is \zedez. If it is reduced then for every $a\in\gB$ there exists a unique \idm $e\in\gK[a]$ such that~$\gen{a}=\gen{e}$. Furthermore, when $e=1$, \cad when $a$ is \ivz, $a^{-1}\in\gK[a]$.
\item \label{i2fact1Etale} \Propeq
\begin{enumerate}
\item $\gB$ is a \cdiz.
\item $\gB$ is without zerodivisors: $xy=0\Rightarrow (x=0$ \emph{or} $y=0)$.
\item $\gB$ is connected and reduced.
\item The \polmin over~$\gK$ of any arbitrary \elt of~$\gB$ is \irdz.
\end{enumerate}
\item \label{i3fact1Etale} If~$\gK\subseteq\gL\subseteq\gB$ and~$\gL$ is a \stf \cdi over~$\gK$, then~$\gB$ is \stfe over~$\gL$. In addition,~$\gB$ is \'etale over~$\gK$ \ssi it is \'etale over~$\gL$ and~$\gL$ is \'etale over~$\gK$.
\item \label{i4fact1Etale} If $(e_1,\ldots,e_r)$ is a \sfio of~$\gB$, $\gB$ is \'etale over $\gK$ \ssi each of the components~$\gB [1/e_i]$ is \'etale over~$\gK$.
\item \label{i5fact1Etale} If~$\gB$ is \'etale, it is reduced.
\item \label{i6fact1Etale} If $\car(\gK) >\dex{\gB:\gK}$ and if~$\gB$ is reduced, it is \'etale.
\end{enumerate}
\end{fact}
\begin{proof}
\emph{\ref{i1fact1Etale}.}
The \elt $a$ of~$\gB$ is annihilated by a \polu of $\KT$ that we express 
in the form $uT^k\big(1 - T\,h(T)\big)$ with $u\in\gK\eti$, $k\geq0$. So~$\gB$ is \zedez. If it is reduced, $a\big(1 - ah(a)\big) = 0$.
Then, $e=ah(a)$ satisfies $a(1-e)=0$ and a fortiori $e(1-e) = 0$. 
Which allows us to conclude.

\emph{\ref{i2fact1Etale}.} The \eqvc of \emph{a}, \emph{b} and \emph{c} is a special case of Lemma~\ref{lemEntReduitConnexe}. 
The implication \emph{d} $\Rightarrow$ \emph{c} is clear.
Let us take a look at \emph{b} $\Rightarrow$ \emph{d.}
Let $x$ be in~$\gB$ and $f(X)$ be its \polmin over~$\gK$. If $f=gh$, with $g$, $h$ \monsz, then $g(x)h(x)=0$, so $g(x)=0$ or $h(x)=0$.
For example $g(x)=0$, and since~$f$ is the \polminz,~$f$ divides $g$, and $h=1$.
 
\emph{\ref{i3fact1Etale}.} Let $(f_1,\ldots,f_s)$ be a~$\gK$-basis of~$\gL$. We can compute an~$\gL$-basis of~$\gB$ as follows. The basis starts with $e_1=1$. Assume we have computed \elts $e_1$, \ldots, $e_r$ of~$\gB$ \lint independent over~$\gL$. The 
$\gL e_i$'s form a direct sum in~$\gB$ and we have a~$\gK$-basis~$(e_if_1,\ldots,e_if_s)$ for each~$\gL e_i$. If $rs=\dex{\gB:\gK}$, we have finished.  
In the opposite case, we can find~$e_{r+1}\in\gB$ which is not in 
$F_r=\gL e_1\oplus \cdots\oplus \gL e_r$. \\
Then,~$\gL e_{r+1}\cap F_r=\so{0}$ (otherwise, we would express $e_{r+1}$ as an~$\gL$-\coli of $(e_1,\ldots,e_r)$),
and we iterate the process by replacing $(e_1, \ldots, e_r)$ with~$(e_1, \ldots, e_{r+1})$.
\\
Once we have a basis of~$\gB$ as an \Levz, it remains to use the transitivity formula of the \discris (\thref{thTransDisc}).
 
\emph{\ref{i4fact1Etale}.} We use the structure \thref{fact.sfio} (\paref{fact.sfio})  for \sfios and the formula for the discriminant of a direct product  of \algs (Proposition~\ref{propTransDisc}).
 
\emph{\ref{i5fact1Etale}.} Let $b$ be a nilpotent \elt of~$\gB$. For all $x\in\gB$  multiplication by~$bx$ is a nilpotent \endo $\mu_{bx}$ of~$\gB$. We can then find a~\hbox{$\gK$-basis} of~$\gB$ in which the matrix of $\mu_{bx}$ is strictly triangular, so $\Tr(\mu_{bx})=\Tr\iBK (bx)=0$. \\
Thus $b$ is in the kernel of the \Kli 

\snic{tr:\gB\to\Lin_\gK(\gB,\gK),\;\;b\mapsto(x\mapsto \Tr\iBK (bx).}

Finally, $tr$ is an \iso since $\Disc\iBK $ is \ivz, so~$b=0$.
 
\emph{\ref{i6fact1Etale}.} With the previous notation, assume~$\gB$ is reduced and we want to show that the \Kli $tr$ is an \isoz.
\\
It suffices to show that $\Ker tr=0$. Suppose $tr(b)=0$, then $\Tr\iBK (bx)=0$ for every $x$ and in particular $\Tr\iBK (b^n)=0$ for all $n>0$. Therefore the \endo $\mu_b$ of multiplication by $b$ satisfies $\Tr(\mu_b^n)=0$ for every $n>0$. The formulas that link the \isN Newton sums to the \elr \smq functions then show that the \polcar of $\mu_b$ is equal to $T^{\dex{\gB:\gK}}$  (cf. Exercise~\ref{exoSommesNewton}).
The Cayley-Hamilton \tho and the fact that~$\gB$ is reduced allow us to conclude that $b=0$.
\end{proof}
%

\begin{theorem} \emph{(Structure \tho for \'etale \Klgsz, 1)}
\label{th1Etale}
\\
Let 
$\gB$ be an \'etale \Klgz.
\begin{enumerate}
\item \label{i1th1Etale} Every \idz~$\gen{b_1,\ldots,b_r}_\gB$ is generated by an \idm $e$ which is a member of~$\gen{b_1,\ldots,b_r}_{\gK[b_1,\ldots,b_r]}$,
and the quotient \alg is \'etale over~$\gK$.
\item \label{i2th1Etale} Let~$\gA$ be a \tf \Kslg of~$\gB$.
\begin{enumerate}
\item $\gA$ is an \'etale \Klgz.
\item There exist an integer $r\geq1$ and a \sfio $(e_1,\ldots,e_r)$ of~$\gA$ such that, for each $i\in \lrbr$,~$\gB[1/e_i]$ is a free module of finite rank over~$\gA[1/e_i]$. In other words,~$\gB$ is a \qf module over~$\gA$.
\end{enumerate}
\item \label{i3th1Etale} $\gB$ is \agsp over~$\gK$.
\item \label{i4th1Etale}
For all $b\in\gB$, $\rC{\gB/\gK}(b)$ is a product of \spl \polsz.
\end{enumerate}
\end{theorem}
\begin{proof}
\emph{\ref{i1th1Etale}.} If the \id is principal this results from Fact~\ref{fact1Etale} item~\emph{\ref{i1fact1Etale}.}
Moreover, for two \idms $e_1$, $e_2$, we~have~$\gen{e_1,e_2}=\gen{e_1+e_2-e_1e_2}$.
Finally, the quotient \alg is itself \'etale over~$\gK$ by the formula for the discriminant of a direct product \algz.

\emph{\ref{i2th1Etale}.} It suffices to prove item \emph{b}, because the result then follows using the transitivity formula for the \discris for each $\gK\subseteq\gA[1/e_i]\subseteq\gB[1/e_i]$ and the formula for the \discri of a direct product \algz.
\\
To prove item \emph{b}, we try to compute a basis of~$\gB$ over $\gA$ by using the indicated method in the case where~$\gA$ is a \cdi for which we know of a \hbox{$\gK$-basis}, given in Fact~\ref{fact1Etale}~\emph{\ref{i3fact1Etale}}.
The \algo is in danger of struggling 
when $e_{r+1}\gA\cap F_r$ is not reduced to $\so{0}$. We then have an \egt $\alpha_{r+1}e_{r+1}=\sum_{i=1}^r\alpha_ie_i$
with all the $\alpha_i$'s in~$\gA$, and $\alpha_{r+1}\neq0$ but not \iv in~$\gA$.
This implies (item~\emph{\ref{i1th1Etale}}) that we find an \idm $e\neq0,1$ in~$\gK[\alpha_{r+1}]\subseteq\gA$.
We then continue with the two \lons at $e$ and $1-e$. Finally, we notice that the number of splits operated thus is a priori bounded by~$\dex{\gB:\gK}$.

\emph{\ref{i3th1Etale}} and \emph{\ref{i4th1Etale}.} Easily result from \emph{\ref{i2th1Etale}.}
\end{proof}

\rem A \gnn of item~\emph{\ref{i1th1Etale}} of the previous \tho is found in Lemmas~\ref{lemKlgEntiere} and~\ref{lemZrZr1}.
\eoe

\medskip
We can construct step by step  
\'etale \Klgs in virtue of the following lemma, which extends Lemma~\ref{lem1EtaleCD}.

\begin{lemma}\label{lemEtaleEtage}
Let~$\gA$ be an \'etale \Klg and $f\in\AT$ be a \spl \poluz.
Then, $\aqo{\AT}{f}$ is an \'etale \Klgz.
\end{lemma}
\begin{proof}
First consider $\aqo{\AT}{f}$ as a free \Alg of rank $\deg f$.
We have $\Disc\iBA =\disc(f)$ (Proposition~\ref{propdiscTra} item \emph{3}). We conclude with the transitivity formula for the \discrisz.
\end{proof}

The two \thos that follow are corollaries.

\begin{theorem}\label{cor1lemEtaleEtage}
Let~$\gB$ be a \Klgz. The \elts of~$\gB$ which are \agsps over~$\gK$ form a sub\algz~$\gA$. In addition, every \elt of~$\gB$ that annihilates a \spl \polu of~$\gA[T]$ is in~$\gA$.
\end{theorem}
\begin{proof}
Let us first show that if $x$ is \agsp over~$\gK$ and $y$ annihilates a \spl \polu $g$ of~$\gK[x][Y]$, then every \elt of~$\gK[x,y]$ is \agsp over~$\gK$.
If $f\in\KX$ is \spl and annihilates $x$, then the sub\algz~$\gK[x,y]$ is a quotient $\aqo{\gK[X,Y]}{f(X),g(X,Y)}$. This \Klg is \'etale by Lemma~\ref{lemEtaleEtage}.\\
Reasoning by \recuz, we can iterate the previous construction.
We obtain the desired result by noting that an \'etale \Klg is \agsp over~$\gK$, and that every quotient of such an \alg is also \agsp over~$\gK$.
\end{proof}

Here is a \gui{\stfez} variant. We give the \dem again because the variations, although simple, point out the precautions we must take in the \stf case.

\pagebreak	

\begin{theorem} \emph{(\Carn of \'etale \Klgsz)}\label{corlemEtaleEtage}
\\
Let~$\gB$ be a \stf \Klg given in the form $\Kxn$.
\Propeq
\begin{enumerate}
\item $\gB$ is \'etale over~$\gK$.
\item  The \polmin over~$\gK$ of each of the $x_i$'s is \splz.
\item $\gB$ is \agsp over~$\gK$.
\end{enumerate}
In particular, a field~$\gL$ that is a Galois extension of~$\gK$ is \'etale over~$\gK$.
\end{theorem}
\begin{proof}
\emph{1} $\Rightarrow$ \emph{3.} By \thref{th1Etale}.

\emph{2} $\Rightarrow$ \emph{1.} Let us first treat the case of a \stf \Klgz~$\gA[x]$ where~$\gA$ is \'etale over~$\gK$ and where the \polminz~$f$ of $x$ over~$\gK$ is \splz. We then have a surjective \homo of the \stf \Klgz~$\aqo{\gA[T]}{f}$ over~$\gA[x]$ and the kernel of this \homo (which is computed as the kernel of a \ali between finite dimensional \Kevsz) is \tfz, therefore generated by an \idmz~$e$. The \Klgz~$\gC=\aqo{\gA[T]}{f}$ is \'etale by Lemma~\ref{lemEtaleEtage}. We deduce that~$\gA[x]\simeq\aqo{\gC}{e}$ is \'etale over~$\gK$.
\\
We can then conclude by \recu on $n$.
\end{proof}
%

\begin{corollary}\label{corcorlemEtaleEtage}
Let $f\in\KT$ be a \poluz. The \adu $\Adu_{\gK,f}$ is \'etale \ssi $f$ is \splz.
\end{corollary}

\rem We have already obtained this result by direct computation of the \discri of the \adu (Fact~\ref{factDiscriAdu}).\eoe

\begin{theorem} \emph{(Primitive \elt \thoz)}
\label{thEtalePrimitif}
\\
Let~$\gB$ be an \'etale \Klgz.
\begin{enumerate}
\item  If~$\gK$ is infinite or if~$\gB$ is a \cdiz,~$\gB$ is a 
monogenic \algz, 
\prmt of the form~$\gK[b]\simeq\aqo{\KT}{f}$ for some $b\in \gB$ and some \spl $f\in\KT$. 
\\
This applies in particular to a field~$\gL$ which is a Galois extension of~$\gK$, such that the extension~$\gL/\gK$ stems from the \elr case studied in \thref{thGaloiselr}.
\item $\gB$ is a finite product of monogenic \'etale \Klgsz.
\end{enumerate}
\end{theorem}
\begin{proof}
\emph{1.} It suffices to treat the case of an \alg with two \gtrsz~$\gB=\gK[x,z]$. We will look for a \gtr of~$\gB$ of the form $\alpha x+\beta z$ with $\alpha$, $\beta\in\gK$.
Let~$f$ and $g$ be the \polmins of $x$ and $z$ over~$\gK$.
We know that they are \splsz.
Let~$\gC=\aqo{\gK[X,Z]}{f(X),g(Z)}=\gK[\xi,\zeta]$.
It suffices to find $\alpha$, $\beta\in\gK$ such that~$\gC=\gK[\alpha\xi+\beta\zeta]$.
To obtain this result, it suffices that the \polcar of $\alpha\xi+\beta\zeta$ be \splz, as we can apply Lemma~\ref{lem1EtaleCD}.
We introduce two \idtrs $a$ and~$b$, and we denote by~$h_{a,b}(T)$ the \polcar of the multiplication by~$a\xi+b\zeta$ in~$\gC[a,b]$ seen as a free~$\gK[a,b]$-\alg of finite rank.
Actually 
$$\preskip.4em \postskip.4em 
\gC[a,b]\simeq \aqo{\gK[a,b][X,Z]}{f(X),g(Z)}. 
$$
Let $d(a,b)=\disc_T(h_{a,b})$.
We make a computation in a \gui{double \aduz} over~$\gC[a,b]$, in which we separately factorize~$f$ and $g$: 

\snic{f(X)=\prod_{i\in\lrbn}(X-x_i) \;\hbox{ and }\;g(Z)=\prod_{i\in\lrbk}(Z-z_j).}

We obtain

\snac{
\pm d(a,b)=\!\!\!\!\!\!\prod\limits_{(i,j)\neq(k,\ell)}\!\!\!\!(a(x_i-x_k)+b\big(z_j-z_\ell)\big)=
(a^{n^2-n}\! \disc f)^{p^2}  (b^{p^2-p} \!\disc g)^n + \ldots
}


In the right-most side of the \egts above we have indicated the term of highest degree when we order the \moms in $a$, $b$ according to a lexicographic order. Thus the \pol $d(a,b)$ has at least an \iv \coez.
It suffices to choose  $\alpha$, $\beta$ such that $d(\alpha,\beta)\in\gK\eti$ to obtain an \elt $\alpha\xi+\beta\zeta$ of~$\gC$ whose \polcar is \splz.
This completes the \dem for the case where~$\gK$ is infinite.
\\
In the case where~$\gB$ is a \cdi we enumerate the integers of~$\gK$ until we obtain $\alpha$, $\beta$ in~$\gK$ with $d(\alpha,\beta)\in\gK\eti$, or until we conclude that the \cara is equal to a prime number $p$. We then enumerate the powers of the \coes of~$f$ and of $g$ until we obtain enough \elts in~$\gK$, or until we conclude that the field~$\gK_0$ generated by the \coes of~$f$ and $g$ is a finite field. In this case,~$\gK_0[x,z]$ is itself a finite field and it is generated by a \gtrz~$\gamma$ of its multiplicative group, so~$\gK[x,z]=\gK[\gamma]$.

\emph{2.} We use the \dem that has just been given for the case where $\gB$ is a \cdiz. If we do not reach the conclusion, it means that the \dem stumbled at a specific place, which reveals that~$\gB$ is not a \cdiz. Since we have a \stf \Klgz, this provides us with an \idm $e\neq0,1$ in~$\gB$.\footnote{For more details see the solution of Exercise~\ref{exothEtalePrimitif}.}
Thus~$\gB\simeq\gB[1/e]\times \gB[1/(1-e)]$. We can then conclude by \recu on $\dex{\gB:\gK}$.
\end{proof}
%

\subsec{\'Etale \algs over a \sply factorial field}

\vspace{3pt}
When every \spl \pol over~$\gK$ can be decomposed into a product of \ird factors, the field~$\gK$ is said to be \emph{separably factorial}.
\index{separably factorial!discrete field}
\index{field!separably factorial}

\begin{lemma}\label{lemsepfactcorps}
A field~$\gK$ is \sply factorial \ssi we have a test for the existence of a zero in~$\gK$ for an arbitrary \spl \pol of $\KT$.
\end{lemma}
\begin{proof}
The second condition is a priori weaker since it amounts to determining the factors of degree $1$ for a \spl \pol of $\KT$. Suppose this condition is satisfied.
 The proof is just about the same as for Lemma~\ref{lemKXfactor}, but asks for a few \sul details. 
\\
Let $f(T)=T^{n}+\sum_{j=0}^{n-1}a_jT^{j}$. We fix an integer $k\in\lrb{2..n-2}$ and we look for the \pols {$g=T^{k}+\sum_{j=0}^{k-1}b_jT^{j}$} that divide $f$. We will show that there is only a finite number of (explicit) possibilities for each of the $b_j$'s.
The \dem of \KROz's \tho uses \uvl \pols  $Q_{n,k,r}(a_0,\dots,a_{n-1},X)\in\ZZ[\ua,X]$, monic
 in $X$, such that $Q_{n,k,r}(\ua,b_r)=0$. 
\\
These \pols can be computed in the \adu $\gA=\Adu_{\gK,f}$ as follows. Let 
$$
\preskip.4em \postskip.4em \ndsp
G(T)=\prod_{i=1}^{k}(T-x_i)=T^{k}+\sum_{j=0}^{k-1}g_jT^{j}. 
$$
We consider the orbit $(g_{r,1},\dots,g_{r,\ell})$ of $g_r$ under the action of $\Sn$, and we obtain  
$$\preskip-.4em \postskip.1em \ndsp
Q_{n,k,r}(\ua,X)=\prod_{i=1}^{\ell}(X-g_{r,i}). 
$$ 
We deduce that
$$\preskip-.2em \postskip.3em \ndsp
\prod\nolimits_{\sigma\in\Sn}\big(X-\sigma(g_r)\big) =Q_{n,k,r}^{n!/\ell}. 
$$ 
Therefore, by Lemma \ref{lemPolCarAdu}, 
$\rC{\gA/\gk}(g_r)(X)=Q_{n,k,r}^{n!/\ell}(X)$.   
Finally, as~$\gA$ is \'etale over $\gK$ (Corollary~\ref{corcorlemEtaleEtage}), the \polcar of $g_r$  is a product of \spl \pols of $\KT$ by  \thref{th1Etale}~\emph{\ref{i4th1Etale}}.
\\
Thus, $b_r$ must be looked for among the zeros of a finite number of \spl \polsz: there is a finite number of possibilities, all of which are explicit.
\end{proof}
%

\begin{theorem} \emph{(Structure \tho for \'etale \Klgsz, 2)}
\label{th2Etale}\\
Suppose~$\gK$ is \sply factorial.
A \Klgz~$\gB$ is \'etale \ssi it is \isoc to a finite product of \'etale fields over~$\gK$.
\end{theorem}
\begin{proof}
Consequence of the primitive \elt \tho (\thref{thEtalePrimitif}).
\end{proof}
%

\begin{corollary} \label{corth1Etale}
If~$\gL$ is an \'etale field over~$\gK$ and if~$\gK$ is \sply factorial, the same goes for~$\gL$.
\end{corollary}
\begin{proof}
Let $f\in\gL[T]$ be a \spl \poluz. The \Llgz~$\gB=\aqo{\gL[T]}{f}$ is \'etale, therefore it is also an \'etale \Klgz. We can therefore find a \sfio  $(e_1,\dots,e_n)$
such that each $\gB[\fraC1{e_i}]$ is connected. 
This is equivalent to factoring~$f$ into a product of \ird factors.
\end{proof}
%

\begin{corollary}
\label{propIdemMini}
\Propeq
\begin{enumerate}
\item  Every \'etale \Klg is \isoc to a product of \'etale fields over~$\gK$.
\item  The field~$\gK$ is \sply factorial.
\item  Every \spl \pol possesses a field of roots which is a \stf extension (thus Galoisian) of~$\gK$.
\item  Every \spl \pol possesses a field of roots which is \'etale over~$\gK$.
\end{enumerate}
\end{corollary}
%
\begin{proof}
For \emph{2} $\Rightarrow$ \emph{4} we use the fact that the \adu for a \spl \pol is \'etale (Corollary~\ref{corcorlemEtaleEtage}) and we apply \thref{th2Etale}.
\end{proof}
%

\begin{corollary}\label{corth3Etale}
If~$\gK$ is \sply factorial and if $(\gK_i)$ is a finite family of \'etale fields over~$\gK$, there exists a Galois extension~$\gL$ of~$\gK$ which contains a copy of each of the~$\gK_i$.
\end{corollary}
%
\begin{proof}
Each~$\gK_i$ is \isoc to a $\aqo\KT {f_i}$
with $f_i$ \spl \irdz. We consider the lcm~$f$ of the $f_i$'s then a \cdr of~$f$.
\end{proof}
%

\subsec{Perfect fields, \spl closure and \agq closure}

\vspace{3pt}
For a field~$\gK$ of finite \cara $p$ the map $x\mapsto x^p$ is an injective \ri \homoz.

In \clama a field~$\gK$ is said to be \ixc{perfect}{field} if it is of infinite \caraz, or if, being of finite \cara $p$, the morphism~$x\mapsto x^p$
is an \isoz.

In \comaz, to avoid the disjunction on the \cara in the \gui{or} above (which cannot be made explicit), we formulate it as follows: 
\emph{if $p$ is a prime number such that $p.1_\gK=0_\gK$, then the \homoz~$\gK\to\gK,\;x\mapsto x^p$ is surjective.}
 
The field of rationals $\QQ$ and the finite fields (including the trivial field) are perfect.

Let~$\gK$ be a field of finite \cara $p$.
An overfield~$\gL\supseteq\gK$ is called a \emph{perfect closure} of~$\gK$ if it is a perfect field and if every \elt of~$\gL$, raised to a certain power $p^k$, is an \elt of~$\gK$.

\index{closure!perfect ---}
 
\begin{lemma}\label{lemClotParf}
A \cdiz~$\gK$ of finite \cara $p$ has a perfect closure~$\gL$,
unique up to unique \isoz. \\
Furthermore,~$\gK$ is a detachable subset of~$\gL$ \ssi there exists a test for \gui{$\exists x\in\gK,\; y=x^p$?} (with extraction of the 
$p$-th root
of $y$ when it  exists). 

\end{lemma}
%
\begin{Proof}{Proof idea. } An \elt of~$\gL$ is encoded by a pair $(x,k)$, where $k \in \NN$ and $x \in \gK$. This encoding represents the 
$p^k$-th root of $x$.
\\
The \egt in~$\gL$, $(x,k)=_\gL(y,\ell)$, is defined by $x^{p^{\ell}}=y^{p^{k}}$ (in~$\gK$), such \hbox{that $(x^p,k+1)=_\gL(x,k)$}. 
\end{Proof}
%

\begin{lemma}\label{lemSqfDec} \emph{(\Algo for squarefree \fcnz)}
\index{algorithm for squarefree factorization}%
\index{factorization!squarefree ---} 
\\
If~$\gK$ is a perfect \cdiz, we have at our disposal an \algo for \emph{squarefree \fcnz} of the lists of \pols of $\KX$ in the following sense. A squarefree \fcn of a family  $(g_1,\ldots,g_r)$ is given by

\vspace{-.5em}
\begin{itemize}
\item a family  $(\lfs)$ of pairwise comaximal \spl \polsz,
\item the expression of each $g_i$ in the form
$$\preskip.2em \postskip.2em 
g_i=\prod\nolimits_{k=1}^sf_k^{m_{k,i}}\; (m_{k,i}\in\NN) . 
$$
\end{itemize}%
\end{lemma}
%

\begin{Proof} {Proof idea. } We start by computing a \bdf for the family $(g_i)_{i\in\lrbr}$ (see Lemma~\ref{lemPartialDec}). 
If some of the \pols in the basis are of the form $h(X^p)$, we know how to express them as $g(X)^p$, and then we replace $h$ by $g$.
We iterate this procedure until all the \pols of the family have a nonzero derivative. Then we introduce the derivatives of the \pols of the family. For this new family we compute a new \bdfz. \\
We iterate the entire procedure until the original goal is reached. The details are left to the reader. 
\end{Proof}

A \cdiz~$\gK$ is said to be \emph{\sply closed} if every \spl \polu of $\KX$ can be decomposed into a product of factors $X-x_i$ ($x_i\in\gK$). 
 
Let~$\gK\subseteq\gL$ be \cdisz. We say that \emph{$\gL$ is a \spl closure of~$\gK$} if~$\gL$ is \sply closed and \spl \agq over~$\gK$.%
\index{field!separably closed ---}%
\index{closure!separable ---}%
\index{separably closed!discrete field}
 
\begin{lemma}\label{lemSepParfait}~
\begin{enumerate}
\item A \cdi is \agqt closed \ssi it is perfect and \sply closed.
\item If a \cdiz~$\gK$ is perfect, every \'etale field over~$\gK$ is perfect.
\item If a perfect \cdi has a \spb closure, it is also an \agq closure.
\end{enumerate}
 \end{lemma}
%
\begin{proof}
\emph{1.} Results from Lemma~\ref{lemSqfDec} and \emph{3} results from \emph{1} and \emph{2.}

\emph{2.} We consider~$\gL$ \'etale over~$\gK$. Let $\sigma:\gL\to\gL:z\mapsto z^p$. \\
We know that~$\gL=\gK[x]\simeq\aqo\KX f$ where~$f$ is the \polmin of $x$ over~$\gK$. The \eltz~$y=x^p$ is a zero of the \pol $f^{\sigma}$, which is \spl and \ird over~$\gK$ because $\sigma$ is an \auto of~$\gK$. We therefore obtain an \iso $\aqo\KX{f^\sigma}\to \gK[y]\subseteq\gL$. Thus~$\gK[y]$ and~$\gL$ are \Kevs of same dimension, so~$\gK[y]=\gL$ and $\sigma$ is surjective. 
\end{proof}
%
 
\begin{theorem}\label{thClsep}
Let~$\gK$ be a \sply factorial and countable \cdiz. 
\begin{enumerate}
\item $\gK$ has a \spl closure~$\gL$, and every \spl closure of~$\gK$ is {$\gK$-\isocz} to~$\gL$.

\item This applies to~$\gK=\QQ$, $\QQ(\Xn)$, $\FF_p$ or $\FF_p(\Xn)$.

\item In addition if~$\gK$ is perfect, then~$\gL$ is an \agq closure of~$\gK$ and every \agq closure of~$\gK$ is~$\gK$-\isoc to~$\gL$.
\end{enumerate}
\end{theorem}
%
\begin{proof}
We only give a sketch of the \dem of item~\emph{1.}
\\
Recall first of all item \emph{2} of \thref{propUnicCDR}: if a \cdr for $f\in\KX$ exists and is \stf over~$\gK$, then every other \cdr for~$f$ over~$\gK$ is \isoc to the first.
\\
Suppose for a moment that we know how to construct a \stf \cdr for every \spl \pol over~$\gK$.
We enumerate all the \spl \polus of $\KX$ in an infinite sequence $(p_n)_{n\in\NN}$.
We call $f_n$ the lcm of the \pols $p_0$, \ldots, $p_n$.
We construct successive \cdrs~$\gK_0$, \ldots,~$\gK_i$, $\ldots$ for these $f_i$'s.
\\
Because of the previously mentioned result, we know how to construct injective \homos of \Klgsz, 
$$\preskip.4em \postskip.4em
\gK_0\vers{\jmath_1}\gK_1\vers{\jmath_2}\cdots\cdots\vers{\jmath_n}\gK_n\vvers{\jmath_{n+1}}\cdots$$
The \spl closure of~$\gK$ is then the colimit of the \sys constructed thus.
It remains to see why we know how to construct a \stf \cdr for every  \spl \polz~$f$ over~$\gK$. If the field is infinite, the fact is given by \thref{thResolUniv}. In the case of a finite field, the study of finite fields directly shows how to construct a \cdrz.
In the most \gnl case, we can construct a \cdr anyway by brute force, by adding the roots one after the other; we consider an \ird factor $h$ of~$f$ and the field~$\gK[\xi_1]=\aqo{\KX}{h}$. 
Over the new field~$\gK[\xi_1]$, we consider an \ird factor~$h_1(X)$ of $f_1(X)={f(X)\over X-\xi_1} $ which allows us to construct~$\gK[\xi_1,\xi_2]$ etc \ldots\ 
This procedure is possible in virtue of Corollary~\ref{corth1Etale} because the successive fields~$\gK[\xi_1]$,~$\gK[\xi_1,\xi_2]$ \ldots\ remain \sply factorial.   
\end{proof}

\rem There exist several ways to construct an \agq closure of~$\QQ$.
The one proposed 
in the previous \tho depends on the chosen enumeration of the \spl \polus of $\QQX$ and it lacks \gmq pertinence. 
From this point of view, the colimit that we construct is actually of significantly less interest than the special \cdrs that we can construct each time we need to. 
There exist other constructions, of a \emph{\gmq} nature, of \agq closures of $\QQ$ which are interesting however as global objects. 
The most renowned is the one based on 
the \agq real number field to which we add an \elt $i=\sqrt{-1}$. 
For each prime number $p$, another very pertinent \agq closure of $\QQ$ is obtained via the intermediate field formed by the $p$-adic \agq numbers.
\eoe

\section{Basic Galois theory (2)} \label{sec2GaloisElr}

This section  complements Section~\ref{secGaloisElr} (see also \Thosz~\ref{corlemEtaleEtage} and~\ref{thEtalePrimitif}).

\smallskip 
\emph{Some remainders.} A Galois extension of~$\gK$ is defined as a \stf field over~$\gK$ which is a \cdr for a \spl \pol of $\KT$.
\thref{thEtalePrimitif} implies that a Galois extension
always  stems from the \elr case studied in \thref{thGaloiselr}.
Finally, \thref{propUnicCDR} says that such a  \cdr is unique up to  \isoz. \eoe

\begin{definition}
\label{defNormale}
An overfield~$\gL$ of~$\gK$ is said to be \emph{normal} (over~$\gK$) if every $x\in\gL$ annihilates a \polu of~$\gK[T]$ which can be decomposed into a product of \lin factors in~$\gL[T]$.
\end{definition}
\index{normal!overfield}
\perso{extension normale: la \dfn is celle of Richman \cite{MRR}.}

\rem Note that if~$\gL$ is a \stf extension of~$\gK$ or more \gnlt if~$\gL$ has a discrete basis as a \Kevz, then the \polmin of an arbitrary \elt of~$\gK$ exists. If the condition of the above \dfn is satisfied, the \polmin itself can be decomposed into \lin factors in~$\gL[T]$.
\eoe

\begin{fact}\label{factNormal1}
Let $f(T)\in\gK[T]$ be a \polu and~$\gL\supseteq\gK $ be a field of roots for~$f$.
Then,~$\gL$ is a normal extension of~$\gK$.
\end{fact}
\begin{proof}
We have~$\gL=\Kxn$ where $f(T)=\prod_{i=1}^n(T-x_{i})$. 
\\
Let $y=h(\xn)$ be an arbitrary \elt of~$\gL$.
Let

\snic{g(\Xn,T)=\prod_{\sigma\in\Sn}\big(T-h^{\sigma}(\uX)\big).}

We clearly have $g(\ux,y)=0$. Moreover, $g(\ux,T)\in\gK[T]$, because each of the \coes of $g(\uX)(T)$ in $\KuX$ is a symmetric \pol in the $X_{i}$'s, hence a \pol in the \elr \smq functions, which are specialized in \elts of~$\gK$ (the \coes of~$f$) by the~$\gK$-\homo $\uX\mapsto\ux$.
\end{proof}

\begin{theorem} \emph{(\Carn of  Galois extensions)}\\
Let~$\gL$ be a \stf field over~$\gK$. 
\Propeq\label{thNormalEtaleGalois}
\begin{enumerate}
\item \label{i1thNormalEtaleGalois}
$\gL$ is a Galois extension of~$\gK$.
\item \label{i2thNormalEtaleGalois} $\gL$ is \'etale and normal over~$\gK$.
\item \label{i3thNormalEtaleGalois} $\Aut_\gK(\gL)$ is finite and the Galois correspondence is bijective.
\item \label{i4thNormalEtaleGalois} There exists a finite group $G\subseteq\Aut_\gK(\gL)$ whose fixed field is~$\gK$.
\end{enumerate}
In this case, in item \ref{i4thNormalEtaleGalois}, we \ncrt have $G=\Gal(\gL/\gK)$.
\end{theorem}
\begin{proof}
\emph{\ref{i1thNormalEtaleGalois}} $\Rightarrow$  \emph{\ref{i2thNormalEtaleGalois}.} This is Fact~\ref{factNormal1}.
 
\emph{\ref{i2thNormalEtaleGalois}} $\Rightarrow$  \emph{\ref{i1thNormalEtaleGalois}} and \emph{\ref{i3thNormalEtaleGalois}.} 
By  the primitive \elt \thoz,~$\gL=\gK[y]$ for some $y$ in~$\gL$.
The \polminz~$f$ of~$y$ over~$\gK$ is \splz, and~$f$ can be completely factorized  in~$\gL[T]$ because~$\gL$ is normal over~$\gK$.
So~$\gL$ is a \cdr for~$f$. Moreover, \thref{thGaloiselr} applies.

\emph{\ref{i4thNormalEtaleGalois}} $\Rightarrow$  \emph{\ref{i2thNormalEtaleGalois}.} It suffices to show that every $x\in\gL$ annihilates a \spl \pol of $\KT$ which can be completely factorized in~$\gL[T]$, because then the extension is normal (by \dfnz) and \'etale (\thref{corlemEtaleEtage}).
Let

\snic{P(T)=\Rv_{G/H,x}(T)=\prod_{\sigma\in G/H}\big(T-\sigma(x)\big)\;$ where $\;H=\St(x).}

The subscript $\sigma\in G/H$ 
means that we take a $\sigma$ in each left coset of $H$ in~$G$.
Hence any two left cosets have the same cardinality.
The \pol $P$ is fixed by $G$, so $P\in\KT$.
Moreover, $\disc (P)=\prod_{i,j\in\lrbk,i< j}(x_i-x_j)^2$ is \ivz.

Finally, since the Galois correspondence is bijective, and since the fixed field of $G$ is~$\gK$, in item \emph{\ref{i4thNormalEtaleGalois}},
we \ncrt have $G=\Gal(\gL/\gK)$.
\end{proof}
%

%
\begin{theorem}
\label{thSG} \emph{(Galois correspondence, complement)}
\\
Let~$\gL/\gK$ be a Galois extension with Galois group $G=\Gal(\gL/\gK)$.
Let $H$ be a detachable subgroup of~$G$, $\sigma$ be an \elt of~$G$,~%
$H_\sigma=\sigma H\sigma^{-1}$.
\begin{enumerate}
\item The field $\sigma(\gL^H)$ is equal to~$\gL^{H_\sigma}$.
\item $\gL^H$ is a Galois extension of~$\gK$ \ssi $H$ is normal in $G$. In this case the Galois group $\Gal(\gL^H/\gK)$ is canonically \isoc to $G/H$.
\end{enumerate}
\end{theorem}
\begin{proof}
\emph{1.} Immediate computation.

\emph{2.} Let~$\gM=\gL^H$. 
By the  primitive \elt \tho write~$\gM=\gK[y]$, such that $H=\St(y)$. The field~$\gM$ is normal over~$\gK$ \ssi for each $\tau\in G$, we have $\tau(y)\in\gM$,
\cad $\tau(\gM)=\gM$. By item~\emph{1} this means that $\tau H \tau^{-1}=H$.
\end{proof}

Now we add some details to \thref{thGaloiselr}.

\begin{theorem}
\label{thCGSynthese} \emph{(Galois correspondence, synthesis)}
\\
Let~$\gL/\gK$ be a Galois extension.
The Galois correspondence works as follows.
\begin{enumerate}
\item For all~$\gM\in \cK_{\gL/\gK}$,~$\gL/\gM$ is a Galois extension with Galois group $\Fix(\gM)$ and $\dex{\gL:\gM}=\#{\Fix(\gM)}$.
\item If $H_1, H_2\in\cG_{\gL/\gK}$ and~$\gM_i=\Fix(H_i)\in \cK_{\gL/\gK}$, then
\begin{itemize}
\item  $H_1\cap H_2$
corresponds to the \Kslg generated by~$\gM_1$ and~$\gM_2$,
\item  $\gM_1\cap \gM_2$ corresponds to the subgroup generated by $H_1$ and $H_2$.
\end{itemize}
\item If $H_1\subseteq H_2$, 
then
\begin{itemize}
\item $\gM_1\supseteq\gM_2$ and $(H_2:H_1)=\dex{\gM_1:\gM_2}$,
\item $\gM_1/\gM_2$ is a Galois extension \ssi $H_1$ is normal in $H_2$.
In this case the group $\Gal(\gM_1/\gM_2)$ is naturally \isoc to $H_2/H_1$.
%
\end{itemize}
\end{enumerate}
\end{theorem}
%

\penalty-2500
\section{\Fp \algsz} \label{sec1Apf}

\vspace{4pt}
\subsec{Generalities}

\Fp \algs are to \syss of \pol \eqns (or \emph{\sypsz}) what \mpfs are to \slisz.\index{polynomial system}

Here we introduce a few basic \gnl facts regarding these \algsz.

The \algs that we consider in this section are associative, commutative and unitary.

\pagebreak

\begin{definition}
\label{defAlg}\label{defi2STF} Let~$\gA$ be a \klgz.
\begin{enumerate}

\item The \algz~$\gA$ is said to be \emph{\tfz} if it is generated by a finite family as a~\hbox{\klgz}. This boils down to saying that it is \isoc to a quotient \alg $\kXn\sur{\fa} $. We then denote it by~$\gA=\kxn$, where $x_i$ is the image of $X_i$ in $\gA$.
This notation does not imply that~$\gA$ is an extension of~$\gk$.%
\index{algebra!finitely generated ---}

\item The \algz~$\gA$ is said to be \emph{\pfz} if it is \pf as a \klgz. This boils down to saying that it is \isoc to an \alg $\kXn\sur{\fa}$, with a \itf $\fa=\gen{\lfs}$.%
\index{algebra!finitely presented ---}

\item The \algz~$\gA$ is said to be \emph{\rpfz}%
\footnote{\gui{\rpfz} expresses a single, well-defined property. Thus it is to be used and considered as a whole (like a single word) and not to be mistakenly subdivided between \gui{\pfz} and \gui{reduced.}}
 if it is \pf as a reduced \klgz. In other words if it is \isoc to a quotient \alg  $\kXn/\sqrt{\fa}$ with a \itfz~$\fa$.%
\index{algebra!finitely presented reduced ---}

\item The \alg $\gA$ is said to be \emph{strictly finite} if~$\gA$ is a \ptf \kmoz. We also say that \emph{$\gA$ is \stf over~$\gk$}.
In the case of an extension, we speak of a \emph{\stf extension} of~$\gk$.%
\index{algebra!strictly finite ---}%
\index{strictly finite!algebra}

\item If~$\gA$ is \stf we denote by

\snic{\Tr\iAk (x),\;\; \rN\iAk (x), \;\;
\rF{\gA/\gk}(x)(T) \;\;\hbox{and}\;\; \rC{\gA/\gk}(x)(T),}

the trace, the \deterz, the \polfon and the \polcar of the \kli $\mu_{\gA,x}\in\End_\gk(\gA)$.
Moreover, by letting $g(T)=\rC{\gA/\gk}(x)(T)$, the \elt $g'(x)$ is called \emph{the different of~$x$}.%
\index{different!of an \elt in a \stf \algz}
\end{enumerate}
\end{definition}

Note that in the case where~$\gk$ is a \cdiz, we indeed find the notion of a \stf \alg given in \Dfnz~\ref{defiSTF}.

\begin{fact}
\label{factPropUnivAPF} \emph{(Universal \prt  of a \apfz)}
\\
The \pf \alg $\aqo{\kXn}{\lfs}=\kxn$ is \caree by the following \prtz: if a \klgz~$\gk\vers{\varphi}\gA$ contains \elts $\yn$ such that the $f_i^\varphi(\yn)$'s are null, there exists a unique \homo of \klgs $\kxn\to\gA$ which sends the $x_i$'s to the~$y_i$'s.
\end{fact}

\subsubsection*{Changing the \sgrz}

\vspace{-.3em}

\begin{fact}
\label{factChscalg}
When we change the \sgr 
for a \pf \algz~$\gA$ the 
ideal of relations between the new \gtrs is again \tfz.
\end{fact}

Refer to Section~\ref{sec pf chg} to verify that what has been explained slightly informally on \paref{nouveausgr} works well in the current case.

\penalty-2500
\subsubsection*{Transitivity (\apfsz)}

\vspace{-.3em}

\begin{fact}
\label{factTransAPF}
If~$\gk\vers \lambda \gA$ and~$\gA\vers\rho \gC$ are two \pf \algsz, then $\gC$ is a \pf \klgz.
\end{fact}
\begin{proof}
Let $\gA=\kuy\simeq\aqo{\kuY}{g_1,\ldots ,g_t}$ and~$\gC=\Aux\simeq\aqo{\AuX}{\lfs}$.
Let $F_{1},\ldots,F_{s}\in\gk[\uY,\uX]$ be \pols such that $F_{i}(\uy,\uX)=f_{i}(\uX)$.
\\
Then,~$\gC=\gk[\rho(\uy),\ux]\simeq\aqo{\gk[\uY,\uX]}{g_1,\ldots ,g_t,F_{1},\ldots,F_{s}}$.
\end{proof}

\subsubsection*{Sub\algsz}

\vspace{-.3em}

\begin{fact}
\label{factSousAlg}
Let~$\gA\subseteq\gC$ be two \tf \klgsz.
If~$\gC$ is a \pf \klg it is also a \pf \Alg (with \gui{the same} \pnz, read in~$\gA$).
\end{fact}
\begin{proof}
Let \spdg $\gC=\kxn\simeq\aqo{\kuX}{\uf}$ and~$\gA=\kxr$.
We have~$\gA\simeq \kXr/\ff$ with $$\ff=\gen{\lfs}\cap\kXr.$$
Let us denote by $\pi:\kXr\to\gA$ the passage to the quotient 
and for $h\in\kXn$, denote by $h^\pi\in\gA[X_{r+1},\ldots ,X_n]$ its image. So
 $$h^\pi=h(\xr,X_{r+1},\ldots ,X_n).$$
Consider the \homo
$$
\gamma:
\begin{array}{l}
 \aqo{\gA[X_{r+1},\ldots ,X_n]}{f_1^\pi,\ldots ,f_s^\pi} \simeq \\[1mm]
 \aqo{\AXn}{X_1-x_1,\ldots ,X_r-x_r,f_1^\pi,\ldots ,f_s^\pi}
\end{array}
\to\gC.
$$
This is the \homo which fixes~$\gA$ and sends $X_k$ to $x_k$ for $k\in\lrb{r+1..n}$. It suffices to show that~$\gamma$ is injective. Every \elt $g$ of~$\gA[X_{r+1},\ldots ,X_n]$ can be written in the form $g=G^\pi$ with $G\in\kXn$.
\\
Suppose that $g$ modulo~$\gen{f_1^\pi,\ldots ,f_s^\pi}$ is in $\Ker\gamma$.
We then have

\snic{g(x_{r+1},\ldots ,x_n)=G(\xn)=0.}

Therefore $G\in\gen{\lfs}$, which gives us $g\in\gen{f_1^\pi,\ldots ,f_s^\pi}$ (after transformation by $\pi$). As required.
\end{proof}

\rem The condition~$\gA\subseteq\gC$ is essential for the \dem to work properly. Moreover, if must be noted that the \id $\ff$ is not \ncrt \tfz.
\eoe

\subsec{The zeros of a \sypz}  \label{secZerosSyp}
Consider a \syp $(\uf)=(\lfs)$ in $\kXn,$ and a \klgz~$\rho:\gk\to\gB$. 
 
\begin{definition}\label{defiZeroSyp}
A \emph{zero of the \sys $(\uf)$ over~$\gB$} is an $n$-tuple 

\snic{(\uxi)=(\xin)\in\gB^n}

satisfying $f_i^\rho(\uxi)=0$ for each $i$.
The set of zeros of $(\uf)$ over~$\gB$ is often symbolically called the
\ixc{variety of the zeros}{of a \sypz} over~$\gB$ of the \sypz,
and thus, we denote it by $\cZ_\gk(\uf,\gB)$ or $\cZ(\uf,\gB)$.
\index{zero!of a \sypz, over an \algz}
\end{definition}

Some zeros are more interesting than others: the closer the \algz~$\gB$ is to~$\gk$, the more interesting is the zero. We pay particular attention to the zeros over~$\gk$, or by default over finite \klgsz.

Two zeros are a priori particularly disappointing. The one provided by the trivial \algz, and the zero $(\xn)$ over the \ixc{quotient algebra}{for a \sypz} associated with the \sypz, \cad
\Grandcadre{$\gA=\kxn = \aqo{\kXn}{\lfs}.$}

Nevertheless this last \alg plays a central role for our \pb because of two findings. The first is the following.

\begin{fact}
\label{factZeros}
For every \klgz~$\gB$ the set of zeros of $(\uf)$ over~$\gB$ is naturally identified with the set of \homos of \klgs from~$\gA$ to~$\gB$.
In particular, the zeros over~$\gk$ are identified with the characters of the \algz~$\gA$.
\end{fact}

%
\begin{Proof}{\Demo on an example. }
Let $\QQ[x,y]=\aqo{\QQ[X,Y]}{X^2+Y^2-1}$. 
Taking a real point $(\alpha,\beta)$ of the circle $X^2+Y^2=1$ amounts to the same thing as taking a morphism $\rho:\QQ[x,y]\vers{}\RR$
(the one which sends $x$ and $y$ to $\alpha$ and~$\beta$).
\end{Proof}

\rdb 
We therefore have a crucial identification, which we write as an \egt
\Grandcadre{$\Hom_{\gk} (\gA,\gB)=\cZ_{\gk}(\uf,\gB)\subseteq \gB^n.$
\label{ZerosCrucial}}

In short the quotient \algz~$\gA$ \emph{intrinsically} summarizes the pertinent information contained in the \syp $(\uf)$.
Which is also why we say that $\cZ_\gk(\uf,\gB)$ is the 
\ixc{variety of the zeros}{of an \alg over another} of~$\gA$ over~$\gB$.

\smallskip 
The second finding (closely related to the previous one by the way) is the following.

From a \gmq point of view \emph{two \syps $(\uf)$ and $(\ug)$ in~$\kuX$ which have the same zeros}, over any arbitrary \klgz, must be considered as \emph{\eqvsz.}
If that is the case, let~$\gA_{1}=\kux$ and~$\gA_{2}=\kuy$ be the two quotient \algs (we do not give the same name to the classes of~$X_i$'s in the two quotients).
Consider the canonical zero~$(\xn)$ of~$(\uf)$ in~$\gA_1$.
 Since $\cZ(\uf,\gA_1)=\cZ(\ug,\gA_1)$,
 we must have $g_j(\ux)=0$ for each~$j.$ This means that $g_j(\uX)$ is null modulo~$\gen{\uf}$.
Similarly, each $f_i$ must be in~$\gen{\ug}$.

Let us summarize this second finding.

\begin{fact}
\label{fact1Zeros}
Two \syps $(\uf)$ and $(\ug)$ in $\kuX$
admit the same zeros, over any \klgz, \ssi they define the same quotient \algz.
\end{fact}

\medskip
\exl
The circles
$x^2+y^2-3=0$ and $x^2+y^2-7=0$
cannot be distinguished by their rational points -- they do not have any (since over~$\ZZ$, the congruence $a^2 + b^2 \equiv 3c^2 \bmod 4$ leads to $a, b, c$ being even), but the quotient \QQlgs are non-\isocz, and we can observe over $\QQ[\sqrt{3},\sqrt{7}]$ that they have distinct zeros \gui{somewhere.}
\eoe

\smallskip
When~$\gk$ is reduced and if we focus on the zeros over the reduced \klgsz, the \algz~$\gA=\aqo{\kuX}{\uf}$ must be replaced by its reduced variant, which is a \rpf \algz
$$\preskip.4em \postskip.4em 
\gA\sur{\DA(0)}=\kXn/\sqrt{\gen{\lfs}}\,. 
$$
We continue this discussion on \paref{subsecNstMorphismes} in the subsection entitled \gui{\nst and equivalence of two categories.}

\subsubsection*{A digression on  \agq computation}
Besides their direct link to the solution of \syps another reason for the importance of \pf \algs is the following.
Each time that an \agq computation reaches an \gui{interesting result} in a \klgz~$\gB$, this computation has only involved a finite number of \elts $y_1$, \dots, $y_n$ of~$\gB$ and a finite number of relations between the $y_i$'s, such that there exist a \pf \klgz~$\gC=\kxn$ and a surjective morphism~$\theta:\gC\to\kyn\subseteq\gB$ which sends the $x_i$'s to the $y_i$'s and such that the \gui{interesting result} has already occurred in~$\gC$ for the $x_i$'s.
In a more scholarly language: every \klg is a filtering colimit  of \pf \klgsz.\footnote{The reader will notice that this subsection is directly copied from the analogous subsection for \mpfsz, \paref{DigCalcAlg}.}

\subsec{The tensor product of two \klgsz}

\vspace{3pt}
The \emph{direct sum} of two \klgsz~$\gA$ and~$\gB$ in the category of \klgs is given by the solution of the following \uvl \pb (\gui{morphism} here means \gui{\homo of \klgsz}).%
\index{direct sum!in a category}
\\
\emph{Find a \klgz~$\gC$ and two morphisms $\alpha:\gA\to\gC$ and $\lambda:\gB\to\gC$ such that, for every \klgz~$\gD$ and for every pair of morphisms $\varphi:\gA\to\gD$ and $\psi:\gB\to\gD$, there exists a unique morphism~$\gamma:\gC\to\gD$ such that
$\varphi=\gamma\circ\alpha$ and $\psi=\gamma\circ\lambda$.}

\vspace{-1.5em}
$$
\xymatrix @R=15pt{
 & \gA\ar[dr]_{\alpha} \ar[drrr]^{\varphi}
\\
\gk\ar@{->}[ru]^{\beta} \ar@{->}[rd]_{\rho} &&\gC\ar@{-->}[rr]^{\gamma !~~~~~~} && \gD
\\
 & \gB\ar[ur]^{\lambda} \ar[urrr]_{\psi}
\\
}
$$

Note that in the category of commutative \risz, the above \uvl \prt means that~$\gC$, with the two morphisms $\alpha$ and $\lambda$, 
is the \ixc{amalgamated sum}{of two arrows of same source in a category}
or the \ixc{push out}{of two arrows of same source in a category}
of the two arrows $\beta:\gk\to\gA$ and $\rho:\gk\to\gB$.
In French, however, the term \emph{carr\'e cocart\'esien} is used, 
formed with the four arrows $\beta$, $\rho$,~$\alpha$ and~$\lambda$ %
 reflecting the above sketch.\footnote{The term \gui{carr\'e cocart\'esien} could be translated as \gui{cocartesian square.}}

\begin{theorem}\label{factSDIRKlg}
Consider two \klgsz~$\gk\vers{\rho}\gB$ and~$\gk\vers{\beta}\gA$.
\begin{enumerate}
\item [A.] \emph{(Direct sum in the category of \klgsz)}\\
The \algs $\gA$ and $\gB$
admit a direct sum~$\gC$ in the category \hbox{of \klgsz}.
Here are different possible descriptions:
\begin{enumerate}
 \item [1.]  If~$\gA=\aqo{\kXn}{\lfs}$,~$\gC=\aqo{\gB[\Xn]}{f_1^\rho,\ldots,f_s^\rho}$
 with the two natural \homosz~$\gA\to\gC$ and~$\gB\to\gC$.
 \item [2.]  If in addition~$\gB=\kyr\simeq\aqo{\kYr}{g_1,\ldots ,g_t}$ is itself a \pf \klgz, we obtain
$$\gC\simeq\aqo{\gk[\Xn,\Yr]}{\lfs,g_1,\ldots ,g_t}.$$
 \item [3.]  \Gnltz, we can consider the \kmoz~$\gC=\gB\otimes_{\gk}\gA$.
 It has  a commutative \ri structure when defining the product as

\snic{(x\otimes a)\cdot(y\otimes b)=xy\otimes ab.}

\snii
We obtain a  \klg structure and we have two natural \homosz~$\gB\to\gC,\,x\mapsto x\otimes 1$ and~$\gA\to\gC,\,a\mapsto 1\otimes a$. This makes~$\gC$ the direct sum of~$\gB$ and~$\gA$.
\item [4.]   If~$\gB=\gk\sur{\fa}$, we obtain~$\gC\simeq\gA/\fb$ where $\fb=\beta(\fa)\gA$.
\item [5.]   If~$\gB=S^{-1}\gk$, we obtain~$\gC\simeq U^{-1}\gA$ where $U=\beta(S)$.
\end{enumerate}
\item [B.]   \emph{(\Edsz)}\\
We can regard~$\gC$ as a \Blgz. We then say that~$\gC$ is the \Blg obtained from~$\gA$ by
\emph{changing the base ring}, or by \ix{scalar extension}. It is then logical to denote it by $\rho\ist(\gA)$.\index{change of the base ring}
\end{enumerate}
\end{theorem}
\facile

We will be mindful of the fact that~$\gk\subseteq \gA$ does not \gnlt imply
$\gB\subseteq \gC$, in particular in case~\emph{4}.

Also note that the tradition is to speak of a \emph{tensor product of \klgsz} rather than of a direct sum. \index{tensor product!of algebras}

\begin{fact}\label{factBimodule}
If~$\gA$ and~$\gB$ are two \klgsz, and $(M,+)$ is an additive group,
 taking an~$\gA\te_\gk\gB$-module  structure over $M$ is the same as taking an external law of the \Amoz~$\,\gA\times M\to M\,$ and an external law of the~\Bmoz~$\,\gB\times M\to M\,$ which both commute and \gui{coincide over~$\gk$.} We also say that $M$ has a    $(\gA,\gB)$-\ix{bimodule} \emph{structure}.
\end{fact}
\begin{proof}
The explanation is the following with~$\gk\vers{\rho}\gB$,~$\gk\vers{\alpha}\gA$.
\\
If we have an~$\gA\te_\gk\gB$-module structure over $M$, we have the two external laws 
\[\preskip-.4em \postskip.4em 
\;\;\begin{array}{ccc} 
\gB\times M\to M,\, (c,m)\mapsto c\cdot m=(1\te c)m, \hbox{ and }  \\[.3em] 
 \gA\times M\to M,\, (b,m)\mapsto b\star m=(b\te 1)m. \end{array}
\]
Since $b\te c=(b\te 1)(1\te c)=(1\te c)(b\te 1)$, we must have $b\star(c\cdot m)=c\cdot(b\star m)$. If $a\in\gk$, $a(1\te 1)=\alpha(a)\te1=1\te \rho(a)$ so we must have $\rho(a)\cdot m=\alpha(a)\star m$. Thus the two laws commute and coincide over~$\gk$. 
\\
 Conversely, from two external laws that commute and coincide over~$\gk$, we can define $(b\te c)m$ by $b\star(c\cdot m)$.
\end{proof}

Here is an important and easy fact regarding the \edsz.
\begin{fact}\label{factEdsAlg}
Consider two \klgsz~$\gk\vers{\rho}\gk'$ and~$\gk\vers{\alpha}\gA$ and let~$\gA'=\rho\ist(\gA)$.
If the \klgz~$\gA$ is \tf (resp.\,\pfz, finite, integral, \stfez)
the same holds for the~$\gk'$-\algz~$\gA'$.
\end{fact}
\facile

\penalty-2500
\subsec{Integral \algsz}
\vspace{4pt}
\subsubsec{The Lying Over lemma}

In this and the following subsection we complete what has already been said on integral \algs in Section~\ref{secApTDN}. 

\smallskip  The following lemma expresses the \cof content of the classical Lying Over lemma of classical mathematics, which asserts that if~$\gB$ is a \ri integral over a sub\riz~$\gA$
there is always a \idep of~$\gB$ above a given \idep of~$\gA$. \index{Lying
Over}

Recall that we denote by $\DA(\fa)$ the nilradical of the \id $\fa$ of~$\gA$.
\begin{lemma}\label{lemLingOver} \emph{(Lying Over)}\\
 Let $\gA\subseteq\gB$ with~$\gB$ integral over~$\gA$ and $\fa$ be an \id of~$\gA$, then $\fa\gB\,\cap\,\gA\subseteq\DA(\fa)$, or (which amounts to the same thing)

\snic{\rD_\gB(\fa\gB)\,\cap\,\gA=\DA(\fa).}

In particular, $1\in\fa\;\Leftrightarrow\; 1\in\fa\gB$.
\end{lemma}
\begin{proof}
If $x\in\fa\gB$ we have $x=\sum a_ib_i$, with $a_i\in \fa,\;b_i\in \gB$.
The $b_i$'s generate a finite \Aslgz~$\gB'$.
Let $G$ be a finite \sgr (with $\ell$ \eltsz) of the \Amoz~$\gB'$.
Let $B_i\in\MM_{\ell}(\gA)$ be a matrix that expresses the multiplication by $b_i$ over $G$. The multiplication by $x$ is expressed by the matrix $\sum a_iB_i$, which has \coes in $\fa$. The \polcar of this matrix, which annihilates $x$ (because~$\gB'$ is a faithful \Amoz), therefore has all of its \coes (except for the leading \coez) in~$\fa$. 
When~$x\in\gA$, this implies $x^{\ell}\in\fa$.
\perso{En fait tout $x\in\fa\gB$
est entier sur l'\id $\fa$. Mais la terminologie n'est pas encore introduite.
Voir le lemme~\ref{lemLingOver2}}
\end{proof}

\rem \label{remlyingover}
Let us indicate how we can deduce the classical Lying Over lemma in \clamaz.
Consider the case where $\fa$ is a \idep and let $S=\gA\setminus\fa$.
Then, $\fa\gB\,\cap\, S = (\fa\gB\,\cap\,\gA)\,\cap\, S$ is empty by Lemma~\ref{lemLingOver}. By \emph{Krull's lemma}, there exists therefore a \idepz~$\fp$ of~$\gB$ such that $\fa\gB\subseteq\fp$ and $\fp\,\cap\, S=\emptyset$, which implies $\fp\,\cap\,\gA=\fa$. It would \egmt be easy to deduce, in \clamaz, Lemma~\ref{lemLingOver} from the Lying Over lemma.\eoe
\index{Krull's lemma}

\medskip \exl
Here we show that the condition \gui{$\gB$ integral over~$\gA$} is crucial 
in the Lying Over lemma. 
Consider~$\gA=\ZZ$,~$\gB=\ZZ[1/3]$ and $\fa=3\ZZ$. Then, we obtain $\fa\gB=\gen{1}$, but $\fa\neq \gen{1}$.\eoe

\subsubsec{Algebras integral over \zed \risz}\label{EntSurZdim}

Here we examine the special case 
of algebras over \zed \risz.

  Algebras  integral over  \cdis   
are an important example of \zed \risz.
In this situation, we give a more precise version on item~\emph{\iref{LID003}} of Lemma~\ref{lemme:idempotentDimension0} as follows (see also \thref{th1Etale}).

\begin{lemma}
\label{lemKlgEntiere}
  An \algz~$\gA$  integral over a \cdi~$\gK$ is \zedz.  
More \prmtz, let \hbox{$\fa=\gen{a_1,\ldots ,a_n}=\gen{\ua}$} be a \itfz. There exist an integer $d$ and an \idm $s\in a_1\gK[\ua]+\cdots +a_n\gK[\ua]$ such \hbox{that $\fa^d=\gen{s}$}.
\end{lemma}
\begin{proof}
An \elt $x$ of~$\gA$ is annihilated by a \polu of $\KX$ that we express as $uX^k\big(1 - X\,h(X)\big)$ where $u\in\gK\eti$, $k\geq0$ and so $x^k\big(1-xh(x)\big)=0$. 

The \idm $e_x$ such that $\gen{e_x}=\gen{x}^d$ for large enough $d$ is then equal to~$\big(xh(x)\big)^k$, and~$d$ is \gui{large enough} as soon as $d\geq k$.

In the case of the \itf $\fa=\gen{a_1,\ldots ,a_n}$, each \idmz~$e_{a_i}$ is an \elt of $a_i\gK[a_i]$. Therefore their gcd, which is the \idm $s$ in the statement of the lemma, is in $a_1\gK[\ua]+\cdots +a_n\gK[\ua]$ (because the gcd of two \idms $e$ and $f$ is $e\vu f=e+f-ef$).
\end{proof}

\begin{lemma}
\label{lemZrZr1}
Let~$\gk$ be a \zed \ri and~$\gA$ an \alg integral over~$\gk$.
\begin{enumerate}\itemsep0pt
\item The \riz~$\gA$ is \zedz.
\item More \prmtz, if  $\fa=\gen{a_1,\ldots ,a_n}
$, there exist an integer $d$ and an \idm $s\in a_1\gk[\ua]+\cdots +a_n\gk[\ua]$ such that~\hbox{$\fa^d=\gen{s}$}.
\item In particular, we obtain for each $a\in\gA$ an \egt
$$\preskip.2em \postskip.0em
a^d\big(1-af(a)\big)=0,
$$
 with some $f(X)\in\kX$ (so, $\big(af(a) \big)^d$ is \idmz).
\end{enumerate}
\end{lemma}
NB: We do not assume that $\rho:\gk\to\gA$ is injective.
\begin{proof}
It suffices to prove item~\emph{2.}\\
By applying the \elgbm from~\paref{MethodeZedRed}, we extend the result of Lemma~\ref{lemKlgEntiere} to the case where~$\gk$ is reduced \zedz. Then we extend the \zed case to the reduced \zed case by passing to the quotient via the nilradical and by using \gui{\imN Newton's method in \algz}
(Section~\ref{secNewton}). More \prmtz, let $\fN=\DA(0)$. By the reduced \zed case, there exist $x_1$, \dots, $x_n$ $\in\gk[\ua]$ such that
$$\preskip.2em \postskip.4em
s=a_1x_1+\cdots +a_nx_n,\hbox{ with }s^2\equiv s  \mod  \fN\hbox{ and }sa_i\equiv a_i
\mod \fN.
$$
The \elt $s$ is congruent modulo $\fN$ to a unique \idm $s_1$, which is written as~$sp(s)$ with $p(T)\in\ZZ[T]$
(Corollary~\ref{corIdmNewton}). Since~$s\in\gk[\ua]$, this gives an \egt $s_1=a_1y_1+\cdots +a_ny_n$ with $y_1$, \dots, $y_n\in\gk[\ua]$. 
In addition,~$s_1a_i\equiv sa_i\equiv a_i$ modulo $\fN$ for each $i$. 
Since $(1-s_1)a_i\in\fN$, there exists a $k_i$ 
such that $(1-s_1)a_i^{k_i}=0$ for each $i$.
Finally, with $k=k_1+\cdots +k_n$, we obtain $\fa^k=\gen{s_1}$.
\end{proof}

Recall that Lemma~\ref{lemZrZr2} establishes the following reciprocal.

\emph{Let~$\gk\subseteq \gA$, with $\gA$ integral over~$\gk$. If $\gA$ is a \zed \riz, then $\gk$ is a \zed \riz.
}

\subsubsec{A weak \nstz}

The following \thoz, for the implication \emph{2} $\Rightarrow$ \emph{3} limited to the case where~$\gA$ is a \cdiz, is often called the \gui{weak \nstz} in the literature, because it can serve as a preliminary to the \nst (in \clamaz). It is to be distinguished from the other weak \nsts already considered in this work.
  
\begin{theorem}
\label{thNst0} \emph{(A weak \nstz)}
\\
Let~$\gK$ be a reduced \zed \ri and~$\gA$ be a \tf \Klgz. For the following \prtsz, we have 1 $\Rightarrow$ 2 $\Leftrightarrow$ 3.
\begin{enumerate}
\item $\gA$ is a local \riz.
\item $\gA$ is \zedz.
\item $\gA$ is finite over~$\gK$.
\end{enumerate}
\end{theorem}
NB: We do not assume that $\rho:\gK\to\gA$ is injective.
\begin{proof} We already know that \emph{3} implies \emph{2.} Let us see that \emph{1} or \emph{2} implies \emph{3}.

We can replace~$\gK$ with $\rho(\gK)$ which is \egmt reduced \zedz. We then have~$\gK\subseteq\gA=\Kxn=\gK[\ux]$. 
Our proof is by \recu on $n$. 
The {\mathrigid 2mu $n=0$} case is trivial. Let us do the inductive step from $n-1$ to~$n$. \\
If~$\gA$ is \zedz, there exist a \pol $R\in\KXn$ and an integer~$\ell$ such that $x_n^\ell\big(1-x_nR(\ux)\big)=0$.
The \pol $X_n^\ell\big(1-X_nR(\uX)\big)$ has one of its \coes equal to $1$ and is therefore primitive.\\
If $\gA$ is local, $x_n$ or $1+x_n$ is \ivz. \Spdg we assume that $x_n$ is \ivz. There exists a \pol $R\in\KXn$ such that $1+x_nR(\ux)=0$.
The \pol $1+X_nR(\uX)$ has one of its \coes equal to $1$ and is therefore primitive.
\\
In both cases, we can perform a \cdv as in Lemma~\ref{lemCDV} (infinite \cdi case) or~\ref{lemNoether} (\gnl case). We then have~$\gA=\Kyn$, and~$\gA$ is finite over~$\gA_1=\gK[y_1,\ldots ,y_{n-1}]\subseteq\gA$. \\
If $\gA$ is \zedz, Lemma~\ref{lemZrZr2} implies that~$\gA_1$ is \zed and we can therefore apply the \hdrz.\\
If $\gA$ is local, item~\emph{3} of \thref{thJacplc} implies that~$\gA_1$ is local and we can therefore apply the \hdrz.
\end{proof}

%
\rem What is new  
for the implication \emph{2} $\Rightarrow$ \emph{3} in \thref{thNst0}, compared to \thref{thSPolZed}
which uses  \iNoe positioning, is therefore the fact that we only assume that the \alg is \tf instead of \pfz.
The two \dems are ultimately based on Lemma~\ref{lemZrZr2} and on a \cdv lemma.
\eoe

\subsubsec{Integral \algs over a \qiriz}

\rdb \label{NOTAReg}
We denote by $\Reg\gA$ the filter of the \ndz \elts of the \riz~$\gA$, such that the total \ri of fractions $\Frac\gA$ is equal to $(\Reg\gA)^{-1}\!\gA$.

\begin{fact}\label{factReduitEntierQi} ~\\
Let~$\gA$ be a \qiriz, $\gK=\Frac \gA$,~$\gL\supseteq\gK$ be a reduced integral \Klg and~$\gB$ be the integral closure of~$\gA$ in~$\gL$.\\
Then,~$\gB$ is a \qiri and $\Frac\gB=\gL=(\Reg\gA)^{-1}\gB$. 
\end{fact}
%
\begin{proof} $\gK$ is \zedr because~$\gA$ is a \qiri (Fact~\ref{factQoQiZed}).
The \riz~$\gL$ is \zed because it is integral over~$\gK$. As it is reduced, it is a \qiriz. As~$\gB$ is \icl in $\gL$, every \idm of~$\gL$ is in~$\gB$, so~$\gB$ is a \qiriz.
\\
Consider some $x\in\gL$ and some \polu $f\in\KX$ which annihilates $x$. By getting rid of the \denos we obtain a \pol 

\snic{g(X)=a_mX^m+a_{m-1}X^{m-1}+\cdots+a_0\in\AX}

which annihilates $x$, with $a_m\in\Reg\gA$. Then, $y=a_mx$, integral over~$\gA$, is in~$\gB$ and $x\in(\Reg\gA)^{-1}\gB$.
\end{proof}

\vspace{-.7em}
\pagebreak	

\subsubsec{Algebras that are \mpfsz}

\vspace{.1em}
\begin{theorem}
\label{propAlgFinPresfin} \emph{(When a \klg is a \pf \kmoz)}
\begin{enumerate}
\item For a \klgz~$\gA$ \propeq
\begin{enumerate}
\item $\gA$ is a \pf \kmoz.
\item $\gA$ is finite and is a \pf \klgz.
\item $\gA$ is \pf and integral over~$\gk$.
\end{enumerate}
\item If these conditions are satisfied and~$\gk$ is \coh (resp.\
\fdi \cohz), then~$\gA$ is \coh (resp.\,\fdi \cohz).
\end{enumerate}
\end{theorem}
\begin{proof}
\emph{1a} $\Rightarrow$ \emph{1b.}
Let $\gA=\som_{i=1}^{m}b_i\gk$ be a \pf \kmoz.
We must give a finite presentation of~$\gA$ as a~\klgz. 
Consider the \sgr $(b_1,\ldots ,b_m)$. On the one hand, we take the 
$\gk$-syzygies given by the \pn of~$\gA$ as a \kmoz. On the other hand we express each $b_ib_j$ as a~$\gk$-\lin combination of the $b_k$'s.
Modulo these last relations, every \pol in the $b_i$'s with \coes in~$\gk$ can be rewritten as a~$\gk$-\lin combination of the $b_i$'s. Therefore it evaluates to~$0$ in~$\gA$ \ssi (as a \polz) it is in the \id generated by all the relations  we have given.
\\
\emph{1b} $\Leftrightarrow$ \emph{1c.} Clear.
\\
\emph{1b} $\Rightarrow$ \emph{1a.}
Suppose that~$\gA$ is finite over~$\gk$ with

\snic{\gA=\kxn
=\aqo{\kuX}{\uf}.}

For each $i$, let $t_i(X_i)\in\gk[X_i]$ be a \polu such that $t_i(x_i)=0$, and $\delta_i=\deg\,t_i$. We have 

\snic{\gA=\aqo{\kuX}{t_1,\ldots ,t_n,h_1,\ldots ,h_s},}

where the $h_j$'s are the reduced $f_j$'s modulo~$\gen{t_1,\ldots,t_n }$.\\
The \gui{\momsz} $\ux^\ud=x_1^{d_1}\cdots x_n^{d_n}$ where $d_1<\delta_1,\ldots,\, d_n<\delta_n$ (which we denote by $\ud<\udel$) form a basis for the \alg $\aqo{\gk[\uX]}{\ut}$ and a \sgr $G$ of the \kmoz~$\gA$. An arbitrary $\gk$-syzygy between these \gtrs is obtained when we write $\som_{j=1}^sg_j(\ux)h_j(\ux)=0$, on the condition that we express it as a~$\gk$-\lin combination of \elts of $G$. We can naturally limit ourselves to the $g_j$'s that are of degree $<\delta_i$ in each variable $X_i$. 
If we fix an index $j\in\lrb{1..s}$ and a \mom $\ux^\ud$ with $\ud<\udel$, we obtain a $\gk$-syzygy between the \elts of $G$ by rewriting $\uX^\ud\, h_j(\uX)$ modulo $\gen{t_1,\ldots,t_n }$ and by saying that the \coli of the \elts of~$G$ obtained as such is null.
These syzygies generate the $\gk$-syzygy module between the \elts of~$G$.

\emph{2}. If~$\gk$ is \coh (resp.\,\fdi \cohz), then we know that~$\gA$ is \coh (resp.\,\fdi \cohz) as a \kmo (since it is \pfz).
Let $(b_i)_{i=1}^m$ be a \sgr of~$\gA$ as a~{\kmo} and $v=(v_1,\ldots ,v_n)\in\Ae n$.
The \idz~$\gen{v_1,\ldots ,v_n}$ is the~\hbox{\kmoz} \tf by 
the $v_ib_j$'s, so it is detachable if~$\gk$ is \fdiz.\\
Moreover, an~$\gA$-syzygy for $v$ can be rewritten as a~$\gk$-syzygy between the $v_ib_j$'s. Therefore a \sgr of the $\gk$-syzygy module between the $v_ib_j$'s gives on a reread a \sgr  of the $\gA$-syzygy module between the~$v_i$'s.
\end{proof}

\subsubsec{Integral \alg over an \aclz}

Here we generalize Proposition \ref{propAECDN}.
\begin{theorem}
\label{Thextent} 
Let $\gA$ be an \aclz, $\gK$ be its quotient field, $\gL$ be a \stf overfield of $\gK$ and $\gB$ be the integral closure of~$\gA$ in~$\gL$.
For $z\in\gL$, let $\mu_{\gL,z}\in \End_\gK(\gL)$ be  multiplication by $z$, and $\nu_z(X)$ and $\chi_z(X)$ be the \polmin and the \polcar of~$\mu_{\gL,z}$ (they are \elts of $\gK[X]$). %
\begin{enumerate}
\item \label{i1Thextent} For  $z\in\gL$, we have $z\in\gB\Longleftrightarrow\nu_z\in \AX\Longleftrightarrow\chi_z\in\AX$. In particular, for $z\in\gB$,
$\rN_{\gL/\gK}(z)$ and $\Tr_{\gL/\gK}(z)\in\gA.$  
\end{enumerate} 
We now suppose that $\gL$ is \'etale over $\gK$, \cad that $\Disc_{\gL/\gK}\in\gK\eti$.
\begin{enumerate} \setcounter{enumi}{1}
\item \label{i2Thextent}
Let $x$ be an \elt of $\gB$ such that $\Kx=\gL$. Let $\Delta_x=\disc(\chi_x)$.
\begin{enumerate}
\item $\gA[x]\simeq \aqo\AX{\chi_x}$, free \Amo of rank $\dex{\gL:\gK}$.
\item We have \hbox{$\gA[x][1/\Delta_x]=\gB[1/\Delta_x]$}, \aclz. 
\item If $\gA$ is a          gcd domain,
if  $\Delta_x=d^{2}b$  and $b$
is squarefree 
then  $\gA[x][1/d]=\gB[1/d]$ and it is an \aclz. 
\end{enumerate} 
\item \label{i3Thextent}
Let $\cB$ be a basis of $\gL$ over $\gK$ contained in $\gB$ and $M\subseteq \gB$ be the \Amo with basis~$\cB$.
\begin{enumerate}
\item The \elt $\Delta=\disc_{\gL/\gK}(\cB)$ is in $\gA$.
\item For all $x\in\gB$,  $\Delta x\in M$, in other words \hbox{$M\subseteq \gB\subseteq \fraC 1\Delta M$}. 
\item If $\gA$ is a gcd domain, for all $x\in\gB$, there exists a $\delta\in\gA$ such that~$\delta^{2}$ divides $\Delta$ \hbox{and  $\delta x\in M$}.
\\
 If in addition $\Delta=d^{2}b$ with $b$ being squarefree,
 $M\subseteq \gB\subseteq \fraC 1 d M$. 
\end{enumerate}
\end{enumerate}
\end{theorem}
%
\begin{proof}
\emph{\ref{i1Thextent}.} 
If $z\in\gB$, it annihilates a \polu
 $h(X)\in\AX$, and the \polz~$\nu_z$ divides~$h$ in $\KX$.
 As $\nu_z$ is monic
and $\gA$ is \iclz, we obtain $\nu_z\in\AX$ by Lemma~\ref{lem0IntClos}. 

Moreover in $\KX$, $\nu_z$ divides  $\chi_z$ and $\chi_z$ divides a power of $\nu_z$, so, still by Lemma~\ref{lem0IntClos},  $\nu_z\in\AX$ is equivalent to $\chi_z\in\AX$. 

\emph{\ref{i2Thextent}a.} Clear: $(1,x,\dots,x^{\dex{\gL:\gK}-1})$ is both a basis of $\gA[x]$ over $\gA$ and of $\gL$ over~$\gK$.
Note that by hypothesis $\chi_x=\nu_x$.

\emph{\ref{i2Thextent}b.}  Consider the special case of \emph{\ref{i3Thextent}b.} where $M=\gA[x]$. 
\\
 We obtain $\gB[1/\Delta_x]=\gA[x][1/\Delta_x]$, and since $\gB$ is \iclz, the same goes for $\gB[1/\Delta_x]$. 

\emph{\ref{i2Thextent}c.}  Special case of \emph{\ref{i3Thextent}c.} with $M=\gA[x]$, reasoning as in~\emph{3c.} 

\emph{\ref{i3Thextent}a.} Immediate consequence of \emph{\ref{i1Thextent}.}

\emph{\ref{i3Thextent}b.} Let $\cB=(\bn)$ and $x=\sum_{i}x_ib_i$ with $x_i\in\gK$.
Consider for example the \coe $x_1$, assumed nonzero.
The $n$-tuple  $\cB'=(x,b_2,\dots,b_n)$ is \hbox{a $\gK$-basis} of $\gL$ contained in $\gB$.
The matrix of $\cB'$ over $\cB$ has as its \deter $x_1$. 
Therefore $x_1^{2}\Delta=x_1^{2}\disc(\cB)=\disc(\cB')\in\gA$. A fortiori $(x_1\Delta)^{2}\in\gA$, and since $\gA$ is \iclz, $x_1\Delta\in\gA$.
Thus all the \coos over $\cB$ of $\Delta x$ are in $\gA$.

\emph{\ref{i3Thextent}c.} When $\gA$ is a gcd domain, we express the \elt $x_1$ as a reduced fraction $x_1=a_1/\delta_1$. Then, since $x_1^{2}\Delta\in\gA$, $\delta_1^{2}$ divides $a_1^{2}\Delta$, and since $\pgcd(a_1,\delta_1)=1$, the \elt $\delta_1^{2}$ divides $\Delta$. We proceed in the same way for each $x_i=a_i/\delta_i$. If $\delta$ is the lcm of the $\delta_i$'s, $\delta^{2}$ is the lcm of the $\delta_i^2$'s, so it divides $\Delta$, and $\delta x\in M$. 
\end{proof}

\section{\Stf \algsz}
\label{subsecAlgStfes}

\subsec{The dual module and the trace}

If $P$ and $Q$ are \ptf \kmosz, we have a canonical \iso $\theta_{P,Q}:P\sta\te_\gk Q\to \Lin_\gk(P,Q)$.

When the context is clear we can identify $\alpha\te x\in P\sta\te_\gk Q$ with the corresponding \kli $y\mapsto \alpha(y)x$. 

In particular, a \syc of $P$, $\big((\xn),(\aln) \big)$, is \cares by the \egtz
\begin{equation}
\label{eq1syc}
\som_{i=1}^n\alpha_i\te x_i=\Id_P.
\end{equation}
Dually we have, modulo the identification of $P$ with $(P\sta)\sta$,
\begin{equation}
\label{eq2syc}
\som_{i=1}^nx_i\te \alpha_i=\Id_{P\sta}.
\end{equation}
This \eqn means that for every~$\gamma\in P\sta$ we have~$\gamma=\sum_{i=1}^n\gamma(x_i)\alpha_i$.

\begin{definota}\label{definotaAsta}
Let~$\gA$ be a \klgz. 
\\
The dual $\Asta$ of the \kmoz~$\gA$ has  an \Amo structure via the external law $(a,\alpha)\mapsto a\centerdot\alpha\eqdefi \alpha\circ \mu_a$, \cad $(a\centerdot\alpha)(x) = \alpha(ax)$.  
\end{definota}

Facts~\ref{factMatriceAlin} and/or~\ref{factMatriceEndo} give the following result.

\begin{fact}\label{factTraceSyc}
Let $\big((\xn),(\aln) \big)$ be a \syc for the \stf \klgz~$\gA$, then the \kli $\mu_{\gA,a}$ is represented in this \sys by the matrix ${\big(\alpha_i(ax_j)\big)_{i,j\in\lrbn}}$ and we have
\begin{equation}\preskip.0em \postskip.0em
\label{eqfactTraceSyc}
\Tr\iAk=\som_{i=1}^nx_i\centerdot\alpha_i,\quad \big(\mathit{\cad}\, \forall a\in\gA,\;\Tr\iAk(a)=\som_{i=1}^n\alpha_i(ax_i)\big).
\end{equation}
\end{fact}

\vspace{-.5em}
\pagebreak

\subsec{Norm and cotransposed \eltz}

We introduce the notion of a \emph{cotransposed \eltz} in a \stf \algz. It suffices to build upon what was said in the case of a free \alg of finite rank on \paref{eqelt0cotransp}.
If~$\gA$ is \stfe over~$\gk$ we can indentify~$\gA$ with a commutative \kslg of $\End_\gk(A)$, where $A$ designates the \kmoz~$\gA$ deprived of its multiplicative structure, by means of the multiplication \homo $x\mapsto\mu_{A,x}=\mu_x$.
Then, since $\wi \mu_x=G(\mu_x)$ for a \polz~$G\in\kT$ (item \emph{6} of \thref{propdef det ptf}), we can define $\wi x$ by the \egt $\wi x=G(x)$, or (what amounts to the same thing) $\wi {\mu_x}=\mu_{\wi x}$. If more precision is \ncr we will use the notation $\Adj\iAk (x)$.
This \elt $\wi x$ is called \emph{the cotransposed \elt of $x$}.
The \egt $\wi {\mu_x}\,\mu_x=\det(\mu_x)\Id_\gA$ then gives
\index{cotransposed!element (in a \stf \algz)}
\begin{equation}
\label{eqeltcotransp}
x\ \Adj\iAk (x)=\rN\iAk (x).
\end{equation}

\begin{lemma}
\label{lemIRAdu} Let~$\gk\vers{\rho}\gA$ be a \stf \algz, $x\in\gA$ and $y\in\gk$.
\begin{enumerate}
\item We have $x\in\Ati$ \ssi $\rN\iAk (x)\in\Ati$.
\\
In this case $x^{-1}=\wi x/\rN\iAk (x)$.
\item $x$ is \ndz in~$\gA$ \ssi $\rN\iAk (x)$ is \ndz in $\gk$. In this case $\wi x$ is \egmt \ndzz.
\item $\rho(\gk)$ is a direct summand in~$\gA$.
\end{enumerate}
Let $e=\ide_0(\gA)$ (such that~$\gen{e}_\gk=\Ann_\gk(\gA)$).
\begin{enumerate}\setcounter{enumi}{3}
\item We have $\rho(y)\in\Ati$ \ssi $y\in(\aqo{\gk}{e})\eti$.
\item $\rho(y)$ is \ndz in~$\gA$ \ssi $y$ is \ndz in $\aqo{\gk}{e}$.
\end{enumerate}
\end{lemma}
NB. If~$\gA$ is a faithful \kmoz, \cad if $\rho$ is injective,
we identify~$\gk$ with $\rho(\gk)$.
Then,~$\gk$ is a direct summand  in~$\gA$, and an \elt $y$ of~$\gk$ is \iv (resp.\,\ndzz) in~$\gk$ \ssi it is \iv (resp.\,\ndzz) in~$\gA$.
\begin{proof}
\emph{1.} In a \mptf an \endo (here $\mu_x$) is a bijection \ssi its \deter is \ivz.

\emph{2.} In a \mptf an \endo is injective \ssi its \deter is \ndzz.

Items~\emph{3}, \emph{4} and~\emph{5} can be proven after \lon at \eco of~$\gk$. By the local structure \tho we are reduced to the case where~$\gA$ is free of finite rank, say $k$.
If $k=0$, then~$e=1$, so~$\gA$ and~$\aqo{\gk}{e}$ are trivial and everything is clear (even if it is a little unsettling). Let us examine the case where $k\geq1$, hence $e=0$, and let us identify~$\gk$ with $\rho(\gk)$.
\\
Items~\emph{4} and \emph{5} then result from items~\emph{1} and \emph{2} because $\rN\iAk (y)=y^{k}$.
\\
For item~\emph{3}, we consider a basis $(b_1,\ldots,b_k)$ of~$\gA$ over~$\gk$ and \elts $a_1$, \ldots, $a_k\in\gk$ such that $\sum_ia_ib_i=1$.
We have $\rN\iAk (\sum_ia_ib_i)=1$. 
Moreover,
for~$y_1$, \dots, $y_k\in\gk$, $\rN\iAk (\sum_i y_ib_i)$ is expressed as a \pog of degree $k$ in $\kuy$
(see the remark on~\paref{factNormeRationnelle}), and so

\snic{\rN\iAk (\sum_ia_ib_i)=\sum_ia_i\beta_i=1}

for suitable $\beta_i\in\gk$.
\\
Let us consider the \elt $\beta\in\End_\gk(\gA)$ defined by $\beta(\sum_ix_i{b_i})=\sum_ix_i\beta_i$. \\
Then, $\beta(1)=1$, so $\beta(z)=z$ for $z\in\gk$,
$\Im \beta=\gk$ and $\beta\circ \beta=\beta$.
\end{proof}

\subsec{Transitivity and rank}\label{subsecTransRang}
\vspace{3pt}
When~$\gA$ is of constant rank $n$, we write $\dex{\gA:\gk}=n$. This generalizes the notation already defined in the free \alg case, and this will be \gne further in Chapter~\ref{chap ptf1} (notation~\ref{notaTraceDetCarAlg}).
In this subsection, $m$ and $n$ are integers.

\begin{fact}\label{factTransptf}
Let~$\gA$ be a \stf \klgz, $M$ be a \ptf \Amo and $\gB$ be a \stf \Algz.
\begin{enumerate}
\item  $M$ is also a \ptf \kmoz.
\item  Suppose $\rg_\gA M=m$ and let $f(T)=\rR\gk(\gA)\in\BB(\gk)[T]$ be the \polmu of~$\gA$ as a \kmoz,
then $\rR\gk(M)=f^m(T)=f(T^m)$.
\item $\gB$ is \stfe over~$\gk$ and $\Tr\iBk=\Tr\iAk\circ \Tr\iBA$.
\end{enumerate}
\end{fact}
\begin{proof} \emph{1.} Assume that~$\gA\oplus E\simeq \gk^r$ (\kmosz) and $M\oplus N\simeq \Ae s$ (\Amosz).
Then~$M\oplus N\oplus E^s\simeq \gk^{rs}$ (\kmosz). We can state this again with \sycs in the following form: if $\big((\xn),(\aln) \big)$ is a \syc for the \kmoz~$\gA$ and $\big((\ym),(\beta_1,\ldots,\beta_m) \big)$ is a \syc for the \Amo $M$, then $\big((x_iy_j),(\alpha_i\circ \beta_j) \big)$ is a \syc for the \kmo $M$.\\
\emph{2.} Left to the reader (who can rely on the previous description of the \sycz, or consult the \dem of Lemma~\ref{lem1TransPtf}).\\
\emph{3.} We work with \sycs as in item~\emph{1} and we apply Fact~\ref{factTraceSyc} regarding the trace.
\end{proof}
%

\begin{theorem}\label{propTraptf}
Let~$\gk\subseteq\gA\subseteq\gB$ be \risz.
Suppose that~$\gB$ is \stf over~$\gA$.
Then
\begin{enumerate}
\item the \ri $\gB$ is \stf over~$\gk$ \ssiz$\gA$ is \stf over~$\gk$,
\item if $\dex{\gA:\gk}=n$ and $\dex{\gB:\gA}=m$, then $\dex{\gB:\gk}=mn$,
\item  if $\dex{\gB:\gk}=mn$ and $\dex{\gB:\gA}=m$, then $\dex{\gA:\gk}=n$.
\end{enumerate}
\end{theorem}
\begin{proof}
\emph{1.} If~$\gB$ is \stf over~$\gk$, then $\gA$ is \stf over~$\gk$; this results from~$\gA$ being a direct summand 
in $\gB$ (Lemma~\ref{lemIRAdu} item~\emph{3}), which is a \ptfz~$\gk$-module.
\\
The converse implication is in Lemma~\ref{factTransptf}.
 
\emph{2} and \emph{3.} Result from item \emph{2} of Fact~\ref{factTransptf}: if $f=T^{n}$ then $f^m(T)=T^{mn}$; if~$f=\sum_kr_kT^{k}$ is a multiplicative \pol  such that $f^m(T)=T^{mn}$ then $f=T^n$, since $f^m(T)=f(T^{m})=\sum_kr_kT^{km}$.
\end{proof}

\rem More \gnl transitivity formulas (in the case of nonconstant rank) are given in Section~\ref{secAppliLocPtf} in the subsection entitled \gui{\nameref{subsecTransPtf}} on \paref{subsecTransPtf} (in particular, see Corollary \ref{lem2TransPtf} and \thref{corthTransPtf}).
\eoe

\section{Dualizing \lin forms, \ste \algs}\label{secAlgSte}

\begin{definition}\label{defdualisante}
\emph{(Non-\dgne \smq bi\lin form, dualizing \lin form, \asez)}
\\
Let $M$ be a \kmo and~$\gA$ be a \klgz.

\begin{enumerate}

\item
If $\phi : M \times M \to \gk$ is a \smq bi\lin form, it is associated with the \kli $\varphi : M \to M\sta$ defined by $\varphi(x) = \phi(x, \bullet) = \phi(\bullet, x)$. \\
We say that $\phi$ is \emph{non-\dgnez} if $\varphi$ is an \isoz.
\index{non-degenerate!bi\lin form}
\index{bilinear form!non-\dgne ---}

\item
If $\lambda \in \Lin_\gk(\gA,\gk)=\Asta$, it is associated with the \smqz~$\gk$-bi\lin form over~$\gA$, denoted by $\Phi_{\gA/\gk,\lambda}=\Phi_\lambda$ and defined by $\Phi_\lambda(x,y) = \lambda(xy)$.
\\
We say that the \lin form $\lambda$ is \ix{dualizing} if $\Phi_\lambda$ is non-\dgnez. 
\\
We call a \emph{Frobenius \algz} an \alg for which there exists a dualizing \lin form.%
\index{linear form!dualizing ---}%
\index{algebra!Frobenius ---}%
\index{Frobenius!algebra}

\item
If~$\gA$ is \stfe over~$\gk$ the form $\Phi_{\Tr\iAk}$ is called the \ix{trace form}.

\item The \alg $\gA$ is said to be \emph{\stez} over~$\gk$ if it is \stfe and if the trace is dualizing, \cad the \ftr is non-\dgnez.
\index{strictly etale@strictly \'etale!algebra}
\index{algebra!strictly \'etale ---}

\end{enumerate}
\end{definition}

\rem If~$\gA$ is free with basis $(\ue)=(e_1,\ldots,e_n)$ over~$\gk$,
the matrix of $\phi$ and that of $\varphi$ coincide (for the suitable bases). Moreover, $\phi$ is non-\dgne \ssi $\Disc\iAk =\disc\iAk (\ue)$ is \ivz.
Note that when~$\gk$ is a \cdi we once again find 
\Dfnz~\ref{defi1Etale} 
for an \'etale \algz.\footnote{We have not given the \gnl \dfn of an \'etale \algz. It so happens that the \'etale \algs over \cdis are always \stes (at least in \clamaz, this is in relation to~\thref{thSepProjFi}), but that it is no longer the case for an arbitrary commutative \riz, hence the necessity to introduce the terminology \gui{\stez} here.}
\eoe

\subsec{Dualizing forms}

\begin{theorem}\label{factCarDua} \emph{(\Carn of the dualizing forms in the \stf case)} \\
Let~$\gA$ be a \klg and $\lambda\in\Asta$. For $x\in\gA$, let $x\sta=x\centerdot\lambda \in\Asta.$
\begin{enumerate}
\item 
If~$\gA$ is \stf and if $\lambda$ is dualizing, then for every \sgr $(x_i)_{i\in\lrbn}$, there exists a \sys $(y_i)_{i\in\lrbn}$ such that we have
\begin{equation}\preskip.4em \postskip.4em
\label{eqDua}
\som_{i=1}^ny_i\sta\te x_i=\Id_\gA, \quad \cad \quad\forall x\in\gA,\;
x = \som_{i=1}^n \lambda(xy_i)x_i.
\end{equation}
Moreover, if~$\gA$ is faithful, $\lambda$ is surjective.

\item
Conversely, if there exist two \syss $(x_i)_{i\in\lrbn}$, $(y_i)_{i\in\lrbn}$ such that  $\sum_iy_i\sta\te x_i=\Id_\gA$, then
\begin{itemize}
\item $\gA$ is \stfez, 
\item the form $\lambda$ is dualizing,  
\item   and we have the \egt $\sum_ix_i\sta\te y_i=\Id_\gA$.
\end{itemize}

\item If~$\gA$ is \stfez, \propeq
\begin{enumerate}
\item  $\lambda$ is dualizing.
\item  $\lambda$ is a basis of the \Amo $\Asta$ 
(which is therefore free of rank $1$). 
\item  $\lambda$ generates the \Amo $\Asta$, i.e.~$\gA\centerdot\lambda=\Asta$. 
\end{enumerate}

\end{enumerate}
\end{theorem}
\begin{proof}
\emph{1.}
On the one hand $y \mapsto y\sta$ is an \iso of~$\gA$ over $\Asta$, and on the other hand every \sgr is the first component of a \sycz.
Let us take a look at the surjectivity. As $\gA$ is faithful we can assume that $\gk\subseteq \gA$. Let $\fa$ be the \id of~$\gk$ generated by the $\lambda(y_i)$'s. \Egrf{eqDua} gives the membership $1=\sum_i\lambda(y_i)x_i\in\fa\gA$. As~$\gA$ is integral over~$\gk$, the Lying Over (Lemma~\ref{lemLingOver}) shows that $1\in\fa$. 

\emph{2.} \Egrf{eqDua} gives  
 $\alpha = \sum_i \alpha(x_i) y_i\sta$ for $\alpha\in\Asta$. This proves that $y \mapsto y\sta$ is surjective. Moreover, if $x\sta=0$, then \hbox{we have $\lambda(xy_i) = 0$}, \hbox{then $x = 0$}. Thus $\lambda$ is dualizing. 
\\
Finally, the \egt $\alpha = \sum_i \alpha(x_i) y_i\sta$ with $\alpha = x\sta$ gives $x\sta = \sum_i \lambda(x_i x) y_i\sta$,  
and since $z\mapsto z\sta$ is a~$\gk$-\isoz, $x= \sum_i \lambda(x_i x)y_i$.

\emph{3.} 
\emph{a} $\Leftrightarrow$ \emph{b.} \gui{$\lambda$ is dualizing} means that $x\mapsto x\sta$ is an \isoz, \cad that~$\lambda$ is an~$\gA$-basis of $\Asta$. 
The implication \emph{c} $\Rightarrow$ \emph{a} results from item~\emph{2} because a \syc is given by $\big((x_i),(y_i\sta)\big)$.
\end{proof}

\exls See Exercises~\ref{exoGroupAlgebra} to \ref{exoFrobeniusAlgExemples} and \Pbmz~\ref{exoBuildingFrobAlgebra}.

1) If $f\in\kX$ is \monz, the \algz~$\gk[x]=\aqo{\gk[X]}{f(X)}$ is a Frobenius \alg (Exercise~\ref{exo1Frobenius}).

2) The \algz~$\gk[x,y]=\aqo{\gk[X,Y]}{X^2,Y^2,XY}$ is not a Frobenius \alg (Exercise~\ref{exoFrobeniusAlgExemples}).
\eoe

\subsubsection*{\Edsz}

\begin{fact}\label{factEdsDualisante}\emph{(Stability of the dualizing forms by \edsz)}\\
Consider two \klgsz~$\gk'$ and~$\gA$
and let~$\gA'=\gk'\te_\gk\gA$.\\ 
If the form $\alpha\in\Lin_\gk(\gA,\gk)$ is dualizing, so is the form $\alpha'\in\Lin_{\gk'}(\gA',\gk')$ obtained by \edsz.
\\
Consequently, \eds preserves the Frobenius \prt of an \algz. 
\end{fact}

\vspace{-.1em}
\pagebreak	

\subsubsection*{Transitivity for dualizing forms}
\vspace{-2pt}

\begin{fact}\label{factDuaTrans}
Let~$\gA$ be a \stf \klgz, $\gB$ be a \stf \Algz, $\beta\in\Lin_\gA(\gB,\gA)$ and $\alpha\in\Lin_\gk(\gA,\gk)$.
\begin{enumerate}
\item If $\alpha$ and $\beta$ are dualizing, so is $\alpha\circ \beta$.
\item If $\alpha\circ \beta$ is dualizing and $\beta$ is surjective (for instance $\gB$ is faithful and $\beta$ is dualizing), then $\alpha$ is dualizing.
%
%
\end{enumerate}
\end{fact}
\begin{proof}  If $\big((a_i),(\alpha_i) \big)$ is a \syc of~$\gA\sur\gk$ and $\big((b_j),(\beta_j) \big)$ is a \syc of~$\gB\sur\gA$, then $\big((a_ib_j),(\alpha_i\circ\beta_j) \big)$ is a \syc of~$\gB\sur\gk$.

 \emph{1.}
For $a \in \gA$, $b \in \gB$, $\eta\in\Lin_\gk(\gA,\gk)$ and $\epsilon\in\Lin_\gA(\gB,\gA)$
 we can easily verify that $ab\centerdot(\eta \circ \epsilon) = (a\centerdot\eta) \circ (b\centerdot\epsilon)$. 
\\
Since $\alpha$ is dualizing, we have $u_i\in\gA$ such that $u_i\centerdot\alpha=\alpha_i$ for $i\in\lrbn$. 
 Since $\beta$ is dualizing, we have $v_j\in\gB$ such that $v_j\centerdot\beta=\beta_j$ for $j\in\lrbm$. 
Then, $u_iv_j\centerdot (\alpha \circ\beta)=\alpha_i\circ\beta_j$, and this shows that $\alpha \circ\beta$ is dualizing.

 \emph{2.} Let $\alpha'\in\Lin_\gk(\gA,\gk)$, which we aim to express in the form $a\centerdot\alpha$. Note that for every $b_0\in\gB$, we have $\big(b_0\centerdot(\alpha'\circ\beta)\big)\frt\gA = \beta(b_0)\centerdot\alpha'$; in particular, if $\beta(b_0)=1$, then $\big(b_0\centerdot(\alpha'\circ\beta)\big)\frt\gA =\alpha'$. Since $\alpha \circ \beta$ is dualizing, there exists a $b\in\gB$ such that $\alpha'\circ \beta = b\centerdot(\alpha\circ \beta)$.  By multiplying this \egt by $b_0 \centerdot$, we obtain, by restricting to~$\gA$, $\alpha' = \big((b_0b) \centerdot (\alpha\circ\beta) \big)\frt\gA = \beta(b_0b) \centerdot \alpha$.
\end{proof}
%

\subsec{\Ste \algsz}

The following \tho is an \imd corollary of \thref{factCarDua}.
\begin{theorem}\label{factCarAste} \emph{(\Carn of  \asesz)} 
Let~$\gA$ be a \stf \klgz. For $x\in\gA$, let $x\sta=x\centerdot\Tr\iAk\in\Asta$. 
\begin{enumerate}
\item 
If~$\gA$ is \stez, then for every \sgr  $(x_i)_{i\in\lrbn}$, there exists a \sys $(y_i)_{i\in\lrbn}$ such that we have 
\begin{equation}\preskip.3em \postskip.3em
\label{eqaste}
\som_{i=1}^ny_i\sta\te x_i=\Id_\gA, \quad \cad \quad\forall x\in\gA,\;
x = \som_{i=1}^n \Tr\iAk(xy_i)x_i.
\end{equation}
Such a pair $\big((x_i),(y_i) \big)$  is called a \emph{\stycz}.\\
In addition, if~$\gA$ is faithful, $\Tr\iAk$ is surjective.
\index{trace system of coordinates}

\item
Conversely, if we have a pair $\big((x_i)_{i\in\lrbn}, (y_i)_{i\in\lrbn}\big)$
that satisfies~\pref{eqaste}, then~$\gA$ is \stez, and we have $\sum_ix_i\sta\te y_i=\Id_\gA$.

\item \Propeq
\begin{enumerate}
\item  $\Tr\iAk$ is dualizing (\cad $\gA$ is \stez).
\item  $\Tr\iAk$ is a basis of the \Amo $\Asta$ 
(which is therefore free of rank~$1$). 
\item  $\Tr\iAk$ generates the \Amo $\Asta$. 
\end{enumerate}
\end{enumerate}
\end{theorem}

\subsubsection*{\Edsz}

The following fact extends Facts~\ref{factEdsAlg} and~\ref{factEdsDualisante}.
\begin{fact}\label{factEdsEtale}
Consider two \klgsz~$\gk'$ and~$\gA$ and let~$\gA'=\gk'\te_\gk\gA$. %
\begin{enumerate}
\item If~$\gA$ is \ste over~$\gk$,
then~$\gA'$ is \ste over~$\gk'$.
\item If~$\gk'$ is \stfe and contains~$\gk$, and if~$\gA'$ is \ste over~$\gk'$, then~$\gA$ is \ste over~$\gk$.
\end{enumerate}
\end{fact}
\begin{proof}
\emph{1.} Left to the reader.
\\
\emph{2.} First assume that~$\gA$ is free over~$\gk$. Let $\Delta=\Disc\iAk =\disc\iAk (\ue)\in\gk$ for a basis $\ue$ of~$\gA$ over~$\gk$. By \eds we obtain the \egtz~$\Delta=\Disc_{\gA'/\gk'}\in\gk'$. If $\Delta$ is \iv in~$\gk'$ it is \iv in~$\gk$ by Lemma~\ref{lemIRAdu}.
In the \gnl case we reduce it back to the previous case by \lon at \eco of~$\gk$.\perso{la \dem fonctionne aussi bien with~$\gk'$ \fpte sur~$\gk$}
\end{proof}
%

\subsubsection*{Transitivity for \asesz}

\begin{fact}\label{propTrEta} 
Let~$\gA$ be a \stf \klg and~$\gB$ be a \stf \Algz. 
\begin{enumerate}
\item If~$\gA$ is \ste over~$\gk$, then~$\gB$ is \ste over~$\gk$.
\item If~$\gB$ is \ste over~$\gk$ and faithful over~$\gA$, then
$\gA$ is \ste over~$\gk$. 
\end{enumerate}
\end{fact}

\begin{proof}
Results from Facts~\ref{factDuaTrans} and \ref{factTransptf}.
\end{proof}

\subsubsection*{Separability and nilpotency}

\begin{theorem}\label{prop2EtaleReduit} \label{propEtaleReduit}
Let~$\gA$ be a \ste \klgz.
\begin{enumerate}
\item If~$\gk$ is reduced, then so is~$\gA$.
\item The \id $\DA(0)$ is generated by the image of ${\rD_\gk(0)}$ in~$\gA$.
\item If~$\gk'$ is a reduced \klgz, ${\gA'} =\gk'\otimes_\gk \gA$ is reduced.
\perso{ce serait super d'avoir la r\'eciproque, if~$\gA$ is \stfe and if for every~$\gk'$ r\'eduit \ldots\ldots then \ldots \stez}

\end{enumerate}
\end{theorem}

\begin{proof}
\emph{1.} We reason more or less as for the case where~$\gk$ is a \cdi (Fact~\ref{fact1Etale}). First suppose that~$\gA$ is free over~$\gk$. Let $a\in\DA(0)$. \\
For all $x\in\gA$  multiplication by $ax$ is a nilpotent \endo  $\mu_{ax}$ of~$\gA$. Its matrix is nilpotent so the \coes of its \polcar are nilpotent (see for example Exercise~\ref{exoNilpotentChap2}), therefore null since~$\gk$ is reduced.
In particular, $\Tr\iAk(ax)=0$. Thus $a$ is in the kernel of the \kli $tr:a\mapsto\big(x\mapsto \Tr\iAk(ax)\big)$.
However, $tr$ is an \iso by hypothesis \hbox{so $a=0$}.
\\
In the \gnl case we reduce it to the case where~$\gA$ is free over~$\gk$ by the local structure \tho for \mptfs (taking into account Fact~\ref{factEdsEtale}~\emph{1}).

 Item \emph{3} results from \emph{1} and from Fact~\ref{factEdsEtale}~\emph{1.} Item \emph{2} results from  \emph{3}, when we consider~$\gk'=\gk\red$.
\end{proof}

The same technique proves the following lemma.

\pagebreak	

\begin{lemma}\label{lemNilpNilp}
If~$\gA$ is \stfe over~$\gk$ and if $a\in\gA$ is nilpotent, the \coes of $\rF{\gA/\gk}(a)(T)$ are nilpotent (except the constant \coez).
\end{lemma}

\subsec{Tensor products}

If $\phi$ and $\phi'$ are two \smq bi\lin forms over $M$ and $M'$, we define a \smq bi\lin form over $M\otimes_\gk M'$, denoted $\phi\otimes\phi'$, by 
$$\preskip.3em \postskip.2em
(\phi\otimes\phi')(x\otimes x', y\otimes y') = \phi(x,y)\phi'(x',y')
\label{NOTAptfb}
.$$

\begin{proposition}
\label{propProTNdg}
 \emph{(Tensor product of two non-\dgnes forms)}
\\
Let $M$, $M'$ be two \ptf \kmos and~$\gA$,~${\gA'}$ two \stf \klgsz.
\begin{enumerate}
\item
If $\phi$ over $M$ and $\phi'$ over $M'$ are two non-\dgne \smq bi\lin forms, so is $\phi \otimes
\phi'$.
\item
If $\lambda \in \Asta$ and $\lambda' \in {\gA'}\sta$ are two dualizing $\gk$-\lin forms, so is $\lambda\otimes\lambda' \in (\gA\otimes_\gk\gA')\sta$.
\end{enumerate}
\end{proposition}
\begin{proof}
\emph{1.}
The canonical \kli $M\sta \otimes_\gk {M'}\sta \to (M\otimes_\gk M')\sta$ is an \iso since $M$, $M'$ are \ptfsz.
Let $\varphi : M \to M\sta$ be the \iso associated with $\phi$, and $\varphi' : M' \to {M'}\sta$ be the one associated with $\phi'$. 
The morphism associated with $\phi\otimes\phi'$ is composed of two \isosz, so it is an \iso
$$
\xymatrix {
M\otimes_\gk M' \ar@{->}[rr]\ar@{->}[dr]_{\varphi\otimes\varphi'}
 && (M \otimes_\gk M')\sta
\\
 & M\sta \otimes_\gk {M'}\sta \ar@{->}[ur]_{\rm can.\ iso.}
\\
}
$$
\emph {2.}
Results from $\Phi_{\lambda\otimes\lambda'} = \Phi_{\lambda}\otimes\Phi_{\lambda'}$.
\end{proof}

The previous proposition and Lemma~\ref{lemTraceProT} give the following result.
\begin{corollary}\label{corlemTraceProT}
Let~$\gA$ and~$\gC$ be two \stf \klgsz.
Then

\snic {\Phi_{\Tr_{(\gA\otimes_\gk\gC)/\gk}}=\Phi_{\Tr\iAk}\otimes \Phi_{\Tr_{\gC/\!\gk}}.}

In particular,~$\gA\otimes_\gk\gC$ is \ste if~$\gA$ and~$\gC$ are \stesz. 
(For the precise computation of the \discriz, see Exercise~\ref{exoTensorielDiscriminant}.)
\end{corollary}

\subsec {Integral \eltsz, \idmsz, \dinz} \label{secEtaleIdmDin}

The following \tho is a subtle consequence of the remarkable Lemma~\ref{lemPolCarInt}. It will be useful in the context of Galois theory for \thref{thZsuffit}. 

\pagebreak

\begin{theorem}\label{thIdmEtale}
Let $\rho:\gk\to\gk'$ be an injective \ri \homo with~$\gk$ \icl in~$\gk'$, and~$\gA$ be a \ste \klgz. By \eds we obtain ${\gA'}=\rho\ist(\gA)\simeq\gk'\otimes_\gk \gA$  \ste over~$\gk'$. 
\begin{enumerate}
\item The \homoz~$\gA\to\gA'$ is injective.
\item The \riz~$\gA$ is \icl in~$\gA'$.
\item Every \idm of~$\gA'$ is in~$\gA$.
\end{enumerate}
\end{theorem}
\begin{proof}
Item \emph{3} is a special case of item \emph{2.}

\emph{1.} Apply the local structure \tho for \mptfs and the \plgrf{plcc.basic.modules} for exact sequences.

\emph{2.}
We can identify~$\gk$ with a sub\ri of~$\gk'$ and~$\gA$ with a sub\ri of~$\gA'$.
Recall that~$\gA$ is finite, therefore integral over~$\gk$. 
It suffices to treat the case where~$\gA$ is free over~$\gk$  (local structure \tho for \mptfs and \plgrf{plcc.entier} for integral \eltsz).\\
Let $(\ue)=(e_1,\ldots,e_n)$ be a basis of~$\gA$ over~$\gk$ and $(\uh)$ the dual basis with respect to the \ftrz. 
If $n=0$ or $n=1$ the result is obvious. 
\\ Suppose $n\geq2$.
Note that $(\ue)$ is also a basis of~$\gA'$ over~$\gk'$.
 In addition, since, for~$a\in\gA$, the \endos $\mu_{\gA,a}$ and $\mu_{\gA',a}$ have the same matrix over $(\ue)$, the \ftr over~$\gA'$ is an extension of the \ftr over~$\gA$ and $(\uh)$ remains the dual basis relative to the \ftr
in~$\gA'$. Let $x=\sum_ix_ie_i$ be an integral \elt of~$\gA'$ over~$\gA$ ($x_i\in\gk'$).
We must prove that the $x_i$'s are in~$\gk$, or (which amounts to the same thing)
integral over~$\gk$. However, $xh_i$ is integral over~$\gk$. The matrix of~$\mu_{\gA',xh_i}$ is therefore an integral \elt of $\Mn(\gk')$ over~$\gk$.
Therefore the \coes of its \polcar are integral over~$\gk$ (Lemma~\ref{lemPolCarInt}),
so in~$\gk$, and in particular $x_i=\Tr_{\gA'/\gk'}(xh_i)\in\gk$.
\end{proof}
%

\begin{lemma}\label{lemDiag}
The cartesian product  $\gk^n$ is a \ste \klgz.
The \discri of the canonical basis is equal to $1$.
If~$\gk$ is a nontrivial connected \riz, this \klg has exactly $n$ \crcs and $n!$ \autos (those that we spot at first sight).
\end{lemma}
%
\begin{proof} The assertion regarding the \discri is clear (Proposition~\ref{propTransDisc}).\\
We obviously have as the \crcs the $n$ natural \prns $\pi_i:\gk^n\to\gk$ over each of the factors, and as the~$\gk$-\autos the $n!$ \autos obtained by permuting the \coosz. 
Let $e_i$ be the \idm defined by $\Ker\pi_i=\gen{1-e_i}$. If $\pi:\gk^n\to\gk$ is a \crcz, the $\pi(e_i)$'s form a \sfio of~$\gk$. Since~$\gk$ is nontrivial and connected, all but one are null, $\pi(e_j)=1$ for example. Then, $\pi=\pi_j$, because they are \klis that coincide over the~$e_i$'s.
Finally, as a consequence every~$\gk$-\auto of~$\gk^n$ permutes the~$e_i$'s. 
\end{proof}

\vspace{-.7em}
\pagebreak

\begin{definition} \emph{(Diagonal \algsz)}  \label{defiDiagAlg} 
\begin{enumerate}
\item A \klg is said to be \emph{diagonal} if it is \isoc to a product \algz~$\gk^n$ for some $n\in\NN$. In particular, it is \stez.

\item Let~$\gA$ be a \ste \klg and~$\gL$ be a \klgz. \\
We say that
\emph{$\gL$ diagonalizes~$\gA$} if~$\gL\otimes _\gk\gA$ is a diagonal \Llgz.
\end{enumerate}
\index{diagonalize}
\end{definition}

\begin{fact}\label{factDiagAlg} \emph{(Monogenic diagonal \algsz)}\\
Let $f\in\kX$ be a \polu of degree $n$ and~$\gA=\aqo\kX{f}$.
\begin{enumerate}
\item The \klgz~$\gA$ is diagonal \ssiz$f$ is \spl and can be decomposed into \lin factors in $\kX$.
\item  In this case, if~$\gk$ is nontrivial connected,~$f$ admits exactly $n$ zeros in~$\gk$, and the \dcn of~$f$ is unique up to the order of the factors. 
\item A \klgz~$\gL$ diagonalizes~$\gA$ \ssiz$\disc(f)$ is \iv in $\gL$ and~$f$ can be decomposed into \lin factors in~$\gL[X]$. 
\end{enumerate}

\end{fact}
\begin{proof} \emph{1.}
If~$f$ is \spl and can be completely factorized, we have an \isoz~$\gA\simeq\gk^n$ by the Lagrange interpolation \tho (Exercise~\ref{exoLagrange}).
\\
Let us show the converse. 
Every \crc $\kX\to\gk$ is an \evn \homoz, so every \crcz~$\gA\to\gk$ is the \evn at a zero of~$f$ in~$\gk$.
Thus the \iso given in the hypothesis is of the form

\snic{\ov g\mapsto \big(g(x_1),\ldots,g(x_n)\big)\quad  (x_i\in\gk\hbox{ and }f(x_i) = 0).}


Then let $g_i$ satisfy $g_i(x_i) = 1$ and, for $j \ne i$, $g_i(x_j) = 0$. For $j \ne i$, the \elt $x_i - x_j$ divides $g_i(x_i) - g_i(x_j) = 1$, so $x_i - x_j$ is \ivz. This implies that $f=\prod_{i=1}^n(X-x_i)$ (again by Lagrange).

\emph{2.} With the previous notations we must show that the only zeros of~$f$ in~$\gk$ are the $x_i$'s. A zero of~$f$ corresponds to a \crc $\pi:\gA\to\gk$. We therefore must prove that~$\gk^n$ does not admit any other \crc than the \prns over each factor. However, this has been proven in Lemma~\ref{lemDiag}.  

\emph{3.} Apply item~\emph{1} to the \Llgz~$\gL\te_\gk\gA\simeq\aqo{\gL[X]} f$.
\end{proof}

\rems~\\
 1) Item~\emph{2} requires~$\gk$ to be connected.

 2) (Exercise left to the reader) If~$\gk$ is a \cdi and if $A$ is a matrix of $\Mn(\gk)$, saying that~$\gL$ diagonalizes $\gk[A]$ means that this matrix is \gui{\digz} in $\Mn(\gL)$, in the (weak) sense that~$\gL^n$
is a direct sum of the eigen-subspaces of $A$.

\rdb%
3) The \dcn of a \riz~$\gA$ into a finite product of nonzero connected \risz, when possible, is unique up to the order of the factors. Each connected factor, \isoc to a localized \riz~$\gA[1/e]$, corresponds in fact to an \emph{indecomposable} \idm $e$.%
\footnote{The \idm $e$ is said to be indecomposable if the \egt $e=e_1+e_2$ with $e_1$, $ e_2$ being \idms and $e_1e_2=0$ implies $e_1=0$ or $e_2=0$.} This can be understood to be a consequence of the structure \tho for finite \agBs (see Theorem~\ref{factagb}). We can also obtain the result by reasoning with a \sfio as in the \dem of Lemma~\ref{lemDiag}.%
\index{indecomposable!idempotent}

4)  In item \emph{2}, the \gui{nontrivial} hypothesis gives a more common statement. Without this hypothesis we would have said in the first part of the sentence: every zero of~$f$ is given by one of the $x_i$'s corresponding to the assumed \dcn of~$f$ into \lin factors.

5) For the most part the previous fact is a more abstract reformulation of the Lagrange interpolation \thoz. 
\eoe

\begin{proposition}\label{propEtaleCdi}
Let~$\gK$ be a \sply factorial \cdi and~$\gB$ be a \stfe \Klgz.
Then,~$\gB$ is \'etale \ssi it is diagonalized by an overfield of~$\gK$ \'etale over~$\gK$.
\end{proposition}
\begin{proof}
Suppose~$\gB$ is \'etale. It is \isoc to a product of fields~$\gK_i$ \'etale over~$\gK$ (\thref{th2Etale}) and there exists a field~$\gL$ \'etale over~$\gK$, which is a Galois extension that contains a copy of each~$\gK_i$ 
(Corollary~\ref{corth3Etale}).
We easily see that~$\gL$ diagonalizes~$\gB$.
\\
Suppose that a field~$\gL$ \'etale over~$\gK$ diagonalizes~$\gB$.
Then, $\Disc_{\gB/\gK}$ is \iv in~$\gL$ therefore in~$\gK$, so $\gB$ is \'etale.
\end{proof}

\entrenous{
In le cas d'un \cdi non \spbz ment factoriel, on doit pouvoir diagonaliser une
\Klg \'etale au moyen d'une autre \Klg \'etale (qui lui ressemble).

En fait on serait int\'eress\'e par un \tho beaucoup plus \gnl
du type suivant.

NB: Le \tho suivant is un peu pr\'ematur\'e if l'on veut vraiment
$\gB$ galoisienne.

\label{thEtaleDin}
\mni{\bf \Tho} \emph{Toute \klg \stez~$\gA$ of rank constant peut \^etre diagonalis\'ee par une \klg
fid\`ele and \stez~$\gB$, and m\^eme par une \aG $(\gk,\gB,G)$.}

S'il y avait une solution \gnq \`a ce \pb cela serait une jolie \gnn of l'\aduz. Lorsque~$\gA$ is une \klg monog\`ene
$\aqo{\kT}{f}$
(with~$f$ \splz) le candidat le plus naturel pour~$\gB$ is $\Adu_{\gk,f}$.

And il semble bien qu'apr\`es \eds fid\`ele, on puisse rendre toute \alg
\ste produit d'\algs \stes monog\`enes (voir Ferrand and sa notion of ``\lotz'')

Tant qu'\`a faire on pourrait
avoir une \aG \gui{assez universelle} modulo la donn\'ee d'\homos
s\'eparants~$\gA\to\gB$ (correspondant aux \homos $t\mapsto x_i$
dans le cas monog\`ene).

If cela marche cela doit s\^urement exister in la litt\'erature,
mais where???
}
%

\section[\Spb \algsz]{\Spb \algsz, \idstz} 
\label{secAlgSpb}

The results in this section will be used in Section~\ref{secAGTG} devoted to \aGsz, but only for \thref{thCorGalGen} which establishes the Galois correspondence in the connected case.
Moreover, they are also very useful when studying \mdisz. Here we will limit ourselves to speaking of \dvnsz.

\begin{definotas}
\label{NotaAGenv}
Let~$\gA$ be a \klgz.
\begin{enumerate}
\item  The \algz~$\gA\otimes_\gk\gA$, called the \emph{enveloping \algz} of~$\gA$, is denoted by~$\env{\gk}{\gA}$.
\item  
This \klg possesses two natural  \Alg structures, respectively given by the \homos $g\iAk:a\mapsto a\te 1$
(left-structure) and $d\iAk:a\mapsto 1\te a$ (right-structure).
\index{algebra!enveloping ---}
 We will use the following abbreviated notation for the two corresponding \Amo structures. For $a\in\gA$ and~$\gamma\in\env\gk\gA$,

\snic{
a\cdot \gamma\,=\, g\iAk(a)\gamma\,=\,(a\te 1) \gamma \,\et\,
 \gamma\cdot a\,=\, d\iAk(a)\gamma\,=\, \gamma(1\te a).
}

\item We will denote by $\rJ\iAk$  (or $\rJ$ if the context is clear) the \id of $\env\gk\gA$ generated by the \elts of the form
$a\te1-1\te a=a\cdot 1_{\env\gk\gA}-1_{\env\gk\gA}\cdot a$. 
\item We also introduce the following \klis
\begin{eqnarray}\preskip.0em \postskip.2em
\label{eqDiAk}
\Delta\iAk:\gA\to\rJ\iAk,\quad a\mapsto a\te1-1\te a.\qquad\\
\label{eqmuAk}
\mu\iAk \,:\,\env\gk\gA\rightarrow \gA,\,\,a\otimes b \mapsto ab
\quad\mathrm{(multiplication)}
\end{eqnarray}
\item In the case where $\gA$ is a \tf $\gk$-\algz, $\gA=\gk[\xn]$, the same holds for $\env{\gk}{\gA}$ and we have the following possible description of the previous objects.
\begin{itemize}
\item $\env{\gk}{\gA}=\gk[\yn,\zn]=\gk[\uy,\uz]$ with
$y_i = x_i\te 1$, $z_i = 1\te x_i$.
\item For $a=a(\ux)\in\gA$, and  $h(\uy,\uz)\in\gk[\uy,\uz]$, we have
\begin{itemize}
\item 
$g\iAk(a)=a(\uy)$,  $d\iAk(a)=a(\uz)$,
\item 
$a\cdot h=a(\uy)h(\uy,\uz)$, $h\cdot a= a(\uz)h(\uy,\uz)$,
\item $\Delta\iAk(a)=a(\uy)-a(\uz)$,
\item and $\mu\iAk(h)=h(\ux,\ux)$.
\end{itemize} 
\item  $\rJ\iAk$ is the \id of $\gk[\uy,\uz]$ generated by the $(y_i-z_i)$'s.
\end{itemize}
\item Finally, in the case where $\gA=\aqo{\kXn}{\lfs}=\gk[\ux]$, in other words when $\gA$ is a \pf \klgz, the same holds for $\env{\gk}{\gA}$ (see \thref{factSDIRKlg}).

\snic {
\env{\gk}{\gA}=\aqo{\gk[\Yn,\Zn]}{\uf(\uY),\uf(\uZ)}=\gk[\uy,\uz].
}

\end{enumerate}
\end{definotas}
Note that
$\mu\iAk(a\cdot\gamma)=a\mu\iAk(\gamma)=\mu\iAk(\gamma\cdot a)$ for~$\gamma\in\env\gk\gA$ and $a\in\gA$.


\subsec{Towards the \spt \idmz}

\begin{fact}\label{fact2OmAbsrait}~

\begin{enumerate}
\item
The map $\mu\iAk$ is a \crc of \Algs
(for the two structures).
\item
We have $\rJ\iAk =\Ker(\mu\iAk )$. So~$\gA\simeq\env\gk\gA\sur{\rJ\iAk}$ and
$$\preskip.3em \postskip.3em
\;\env\gk\gA = (\gA \otimes 1) \oplus \rJ\iAk =
(1 \otimes \gA) \oplus \rJ\iAk 
,$$
and $\rJ\iAk$ is the left- (or right-) \Amo generated by $\Im\Delta\iAk$.
\item In the case where $\gA=\aqo{\kXn}{\lfs}=\gk[\ux]$ we obtain
$$\preskip.2em \postskip.4em
\gk[\uy,\uz]=\gk[\uy]\oplus \gen{y_1-z_1,\dots,y_n-z_n}=\gk[\uz]\oplus \gen{y_1-z_1,\dots,y_n-z_n}.
$$
\end{enumerate}
\end{fact}
\begin{proof}
The inclusion $\rJ\iAk\subseteq\Ker(\mu\iAk)$ is clear. Denoting $\Delta\iAk$ by $\Delta$, we have

\snac
{
\som_ia_i\te b_i = \big(\som_i a_ib_i\big)\te 1 - \som_i a_i \cdot \Delta(b_i) =
1\te \big(\som_i a_ib_i\big) - \som_i \Delta(a_i) \cdot b_i
.}

We deduce that $\Ker(\mu\iAk)$ is the (left- or right-) \Amo generated by~$\Im\Delta$ and therefore that it is contained in $\rJ\iAk$. 
\\
The result follows by~\ref{prdfCaracAlg}.
\end{proof}
%

\exl For~$\gA=\kX$, we have $\env\gk\gA\simeq \gk[Y,Z]$ with the \homos 
\[\preskip.4em \postskip.4em 
\begin{array}{ccl} 
h(X)\mapsto h(Y)  &   & \hbox{  (on the left-hand side) and }  \\[.3em] 
h(X)\mapsto h(Z)  &   &  \hbox{  (on the right-hand side)},   
 \end{array}
 \]
so $h\cdot g=h(Y)g$ and $g\cdot h=h(Z)g$. We also have
$$\preskip.2em \postskip.4em 
\Delta\iAk(h)=h(Y)-h(Z),\, \mu\iAk\big(g(Y,Z)\big)=g(X,X) \hbox{ and }\rJ\iAk=\gen{Y-Z}.
$$
We see that $\rJ\iAk$ is free with $Y-Z$ as its basis over $\env\gk\gA$,
and as a left-\Amoz, it is free with basis $\big((Y-Z)Z^n \big)_{n\in\NN}$. 
\eoe

\begin{fact}\label{factOmAbsrait} We write $\Delta$ to denote $\Delta\iAk$.
\begin{enumerate}
\item For $a,b\in\gA$ we have $\Delta(ab)=\Delta(a)\cdot b+a\cdot\Delta(b)$. More \gnltz,
\vspace{.5mm}

\snic{
\arraycolsep2pt
\begin{array}{rcl}
\Delta(a_1\cdots a_n)  & =  & \Delta(a_1)\cdot a_2\cdots a_{n}+a_1\cdot \Delta(a_2)\cdot a_3\cdots a_n +
\cdots \\[.8mm]
  & &
\quad{}+
a_1\cdots a_{n-2}\cdot\Delta(a_{n-1})\cdot a_n
+
a_1\cdots a_{n-1}\cdot\Delta(a_n).
\end{array}
}
\item If~$\gA$ is a \tf \klgz, generated by $(\xr)$, $\rJ\iAk$ is a \itf of $\env\gk\gA$, generated by $(\Delta(x_1),\ldots,\Delta(x_r))$.
\item  Over the \id $\Ann(\rJ\iAk)$, the two structures of \Amosz, on the left- and right-hand sides, coincide.
In addition, for $\alpha \in \Ann(\rJ\iAk)$ and~$\gamma \in \env\gk\gA$,
we have
\begin{equation}\preskip.2em 
\label{eqa2struct}
\gamma\alpha =  \mu\iAk(\gamma) \cdot \alpha =
\alpha \cdot\mu\iAk(\gamma).
\end{equation}

\end{enumerate}

\end{fact}
%
\begin{proof} \emph{1.} Immediate computation. Item \emph{2} results from it since $\rJ\iAk$ is the \id generated by the image of $\Delta$, and since for every \gui{\momz} in the \gtrsz, for example $x^3y^4z^2$, $\Delta(x^3y^4z^2)$, is equal to a \coli (with \coes in $\env\gk\gA$) of the images of the \gtrs $\Delta(x)$, $\Delta(y)$ and $\Delta(z)$.

 \emph{3.} The \id $\fa=\Ann(\rJ\iAk)$ is an $\env\gk\gA$-module, so it is stable for the two  \Amo laws. 
Let us show that these two structures coincide.
If~$\alpha\in\fa$, for every $a\in\gA$ we have $0=\alpha(a\cdot 1-1\cdot a)=a\cdot \alpha-\alpha\cdot a$.
\\
\Egtz~(\ref{eqa2struct}) stems from the fact that
$\gamma - \mu\iAk(\gamma) \cdot 1 $ and $\gamma - 1 \cdot \mu\iAk(\gamma)$
 are \hbox{in $\Ker\mu\iAk=\rJ\iAk$}.
\end{proof}
%

\begin{lemma}\label{lemmeIdempotentAnnJ} 
The \id $\rJ\iAk$ is generated by an \idm if and only if

\snic{1 \in
\mu\iAk\big(\Ann(\rJ\iAk)\big).}

Moreover, if $1 = \mu\iAk(\vep)$ with $\vep
\in \Ann(\rJ\iAk)$, then $\vep$ is an \idmz, and we have
$$\preskip-.4em \postskip.3em 
\Ann(\rJ\iAk) = \gen {\vep} \hbox{ and } \rJ\iAk = \gen {1 - \vep}, 
$$
such that $\vep$ is uniquely determined.
\end{lemma}

\begin{proof} We omit the~$\gA\sur\gk$ subscript.
If $\rJ = \gen {\vep}$ with an \idm $\vep$, we obtain the \egts $\Ann(\rJ) = \gen {1 - \vep}$ and
$\mu(1 - \vep) = 1$. \\
Conversely, suppose that $1 = \mu(\vep)$ with $\vep \in \Ann(\rJ)$. Then $\mu(1 - \vep) = 0$, so $1 - \vep \in \rJ$, then $(1-\vep)\vep = 0$, \cad $\vep$ is \idmz. \\
The \egt $1 = (1-\vep) + \vep$ implies that $\Ann(\rJ) = \gen {\vep}$ and $\rJ = \gen {1 - \vep}$.
\end{proof}

\vspace{-.7em}
\pagebreak	

\subsubsection*{Bézout matrix of a \syp} 
Let $f_1$, \ldots, $f_s \in \gk[\Xn]=\kuX$.\\ 
We define the \ixy{Bézout}{matrix} of the \sys $(\uf)=(f_1, \ldots, f_s)$ in the variables~$(\Yn,\Zn)$
by%
\index{Bézout!matrix}
\[\preskip.4em \postskip.4em 
\begin{array}{ccc} 
  \BZ_{\uY,\uZ}(\uf)=(b_{ij})_{i\in\lrbs, j \in \lrbn},\;\hbox{ where }   \\[2mm] 
b_{ij} = \frac
{f_i(Z_{1..j-1}, Y_j, Y_{j+1..n}) - f_i(Z_{1..j-1}, Z_j, Y_{j+1..n}) }{
Y_j - Z_j} \,.\end{array}
\]

Thus for $n = 2$, $s=3$:
$$\preskip.2em \postskip.3em 
\BZ_{\uY, \uZ}(f_1, f_2, f_3) = 
\cmatrix {
{f_1(Y_1,Y_2) - f_1(Z_1,Y_2) \over Y_1-Z_1} & 
{f_1(Z_1,Y_2) - f_1(Z_1,Z_2) \over Y_2-Z_2} \cr
\noalign {\smallskip}
{f_2(Y_1,Y_2) - f_2(Z_1,Y_2) \over Y_1-Z_1} & 
{f_2(Z_1,Y_2) - f_2(Z_1,Z_2) \over Y_2-Z_2} \cr
\noalign {\smallskip}
{f_3(Y_1,Y_2) - f_3(Z_1,Y_2) \over Y_1-Z_1} & 
{f_3(Z_1,Y_2) - f_3(Z_1,Z_2) \over Y_2-Z_2} \cr
}. 
$$
For $n=3$, the $i^{\rm th}$ row of the Bézout matrix is
$$\preskip.4em \postskip.3em\ndsp 
\big[{\,
{f_i (Y_1,Y_2,Y_3) - f_i (Z_1,Y_2,Y_3) \over Y_1-Z_1} \;
{f_i (Z_1,Y_2,Y_3) - f_i (Z_1,Z_2,Y_3) \over Y_2-Z_2} \;
{f_i (Z_1,Z_2,Y_3) - f_i (Z_1,Z_2,Z_3) \over Y_3-Z_3}}\,\big]. 
$$
We have the \egt
$$\preskip.2em \postskip.4em
\BZ_{\uY,\uZ}(\uf) \cdot \cmatrix {Y_1 - Z_1\cr \vdots\cr Y_n - Z_n\cr}
= \cmatrix {f_1(\uY) - f_1(\uZ)\cr \vdots\cr f_s(\uY) - f_s(\uZ)\cr}
\eqno(\star)
$$
In addition $\BZ_{\uX,\uX}(\uf)=\JJ_\uX(\uf)$, the Jacobian matrix of $(f_1, \ldots, f_s)$.

Now consider a \tf \klg

\snic{\gA =
\gk[\xn] = \gk[\ux],}

with \pols $f_i$ satisfying $f_i(\ux) = 0$ for every $i$. 
Its enveloping \alg is $\env\gk\gA = \gk[\yn,\zn]$ (using the notation from the beginning of the section).

Then the matrix $\BZ_{\uy,\uz}(\uf)\in\MM_{s,n}(\env{\gk}{\gA})$ has as its image under $\mu\iAk$ the Jacobian matrix $\JJ_\ux(f_1, \ldots, f_s)\in\MM_{s,n}(\gA)$.

For a minor $D$ of order $n$ of $\BZ_{\uy,\uz}(\uf)$,  
\Egtz~$(\star)$ shows that $D \,(y_j - z_j) = 0$ for $j\in\lrbn$. In other words $D  \in \Ann(\rJ\iAk)$. The Bézout matrix therefore allows us to construct \elts of the \id $\Ann(\rJ\iAk)$.  
\\
In addition, $\delta := \mu\iAk(D )$ is the corresponding minor in $\JJ_\ux(\uf)$.

Let us give an application of this theory
to the special case 
when the transposed matrix $\tra{\JJ_{\ux}(\uf)} :\gA^s \to \gA^n$ is surjective, \cad $1\in\cD_n(\JJ_{\ux}(\uf))$. We therefore have an \egt $1 = \sum_{I\in\cP_{n,s}} u_I\delta_I$ in~$\gA$, where $\delta_I$ is the minor of the extracted matrix of $\JJ_{\ux}(\uf)$ on the rows $i\in I$. By letting $\vep = \sum_{I\in\cP_{n,s}} u_ID _I \in \env{\gk}{\gA}$, we obtain $\mu\iAk(\vep) = 1$ with $\vep \in \Ann(\rJ\iAk)$.
\\
Recap: $\vep$ is what we call the \idst of $\gA$, and $\gA$ is a \spb \algz, notions which will be defined later (\Dfnz~\ref{defiSpb}).\\
Therefore, if~$\gA$ is a \pf \klg $\aqo{\gk[\uX]}{\uf}$ and if the \ali $\tra{\JJ_{\ux}(\uf)} :\gA^s \to \gA^n$ is surjective, then $\gA$ is \spbz.
\\
More \gnltz, for a \apf $\gA = \aqo{\gk[\uX]}{\uf}$, we will see that $\Coker\big(\tra{\JJ_{\ux}(\uf)}\big)$ and $\rJ\iAk \big/ \rJ\iAk^2$ are \isoc \Amos (\thref{thDerivUnivPF}).

\subsec{Derivations}

\begin{definition}
Let $\gA$ be a \klg and $M$ be an \Amoz.\label{defiDeriv} \\
We call a \emph{$\gk$-derivation of $\gA$ in $M$}, a \kli $\delta:\gA\to M$ which satisfies the Leibniz \egt
$$
\delta(ab)=a\delta(b)+b\delta(a).
$$%
\index{derivation!of an \alg in a module}%
We denote by $\Der\gk\gA M$ the \Amo of the $\gk$-\dvns of $\gA$ in $M$.%
\index{derivation!module of ---s}
A \dvn with values in $\gA$ is \gui{simply} called a \dvn of $\gA$. When the context is clear, $\mathrm{Der}(\gA)$ 
is an abbreviation for~$\Der{\gk}{\gA}{\gA}$.%
\index{derivation!of an \algz}   
\end{definition}

Note that $\delta(1)=0$ because $1^{2}=1$, and so $\delta\frt\gk=0.$ 

\begin{thdef}\label{thDerivUniv} \emph{(\Uvl derivation)}%
\index{derivation!universal ---}
\\
The context is that of \Dfn~\ref{NotaAGenv}.
\begin{enumerate}
\item 
Over $\rJ/\rJ^{2}$ the two  \Amo structures (on the left- and right-hand sides) coincide.
\item The composite map $\rd:\gA\to\rJ/\rJ^{2}$, defined by $\rd(a)=\ov{\Delta(a)}$, is a $\gk$-\dvnz.
\item It is a \uvl $\gk$-\dvn in the following sense. %
\\
For every \Amo $M$ and every $\gk$-\dvn $\delta:\gA\to M$, there exists a unique \Ali $\theta:\rJ/\rJ^{2}\to M$ such that $\theta\circ \rd=\delta$.

\vspace{-1em}
\Pnv{\gA}{\rd}{\delta}{\rJ/\rJ^{2}}{\theta}{M}{}{$\gk$-\dvnsz}{\Alisz.}

\vspace{-1.2em}
\end{enumerate}
The \Amo $\rJ/\rJ^{2}$, denoted by $\Om{\gk}{\gA}$, is called the \emph{module of (K\"ahler) \diles of~$\gA$.}
\end{thdef}%
\index{module!of (K\"ahler) differentials}%
\index{Kahler@K\"ahler!differential}%
\index{differential!(K\"ahler) ---}
%
\begin{proof}
Items \emph{1} and \emph{2} are left to the reader.\\
\emph{3.} The uniqueness is clear, let us show the existence. 

\pagebreak	
We define the \kli $\tau:\gA\te_\gk\gA\to M$

\vspace{-.7em}
\tri{\gA}{\Delta}{\delta}{\gA\te_\gk\gA}{\tau}{M}
{$\tau(a\te b)=-a\,\delta(b).$}

\vspace{-1em}
The diagram commutes and $\tau$ is 
$\gA$-\lin on the left-hand side.\\
It remains to see that $\tau(\rJ^{2})=0$, because $\theta$ is then defined by restriction and passage to the quotient of $\tau$. 
 We verify that $\tau(\Delta(a)\Delta(b))=b\delta(a)+a\delta(b)-\delta(ab)=0$.
\end{proof}

We now consider the case of a \apf 

\snic{\gA=\aqo{\kXn}{\lfs}=\gk[\ux].}

We use the notations in \ref{NotaAGenv}. 
Recall that the Jacobian matrix of the \syp is defined as
$$\preskip.0em \postskip.6em
 \JJ_{\uX}(\uf)\; =\;
\bordercmatrix [\lbrack\rbrack]{
    & X_1                     & X_2                     &\cdots  & X_n \cr
f_1 & \Dpp {f_1}{X_1} &\Dpp {f_1}{X_2}  &\cdots  &\Dpp {f_1}{X_n} \cr
f_2 & \Dpp {f_2}{X_1} &\Dpp {f_2}{X_2}  &\cdots  &\Dpp {f_2}{X_n} \cr
f_i & \vdots                  &                         &        & \vdots              \cr
 & \vdots                  &                         &        & \vdots              \cr
f_s & \Dpp {f_s}{X_1} &\Dpp {f_s}{X_2}  &\cdots  &\Dpp {f_s}{X_n} \cr
}
.$$
In the following \thoz, we denote by $\rja=\tra{\JJ_{\uX}(\uf)}:\Ae s\to \Ae n$ the \ali defined by the transposed matrix, and by $(e_1,\dots,e_n)$ the canonical basis of~$\Ae n$.
We define 
\[ 
\begin{array}{rcl} 
\delta:\gA\to\Coker(\rja)  &:   & g(\ux)\mapsto   \som_{i=1}^{n} \Dpp {g}{X_i}(\ux)\,\ov{e_i},\\[1.5mm] 
\lambda:\Ae n \to \rJ/\rJ^{2}  &:   &  e_i\mapsto\rd(x_i)=\ov{y_i-z_i}. 
 \end{array}
\] 
 
\medskip 
\begin{theorem}\label{thDerivUnivPF}  \emph{(\Uvl derivation via the Jacobian)}
\begin{enumerate}
\item The map $\delta$ is a $\gk$-\dvn with $\delta(x_i)=\ov{e_i}$.
\item The \Ali $\lambda$ induces by passage to the quotient an \iso $\ov \lambda:\Coker(\rja)\to\rJ/\rJ^{2}$.
%
\end{enumerate}
Consequently, $\delta$ is \egmt a \uvl \dvnz.
$$
\xymatrix @C=.75cm@R=.4cm{ 
         &&\Coker(\rja) \ar[dd]^{\ov\lambda}&\delta(x_i)\ar@{<->}[dd]\\
\gA\ar[urr]^{\delta}\ar[drr]_\rd\\
         &&\rJ/\rJ^{2}& \rd(x_i)\\ 
}
$$ 
\end{theorem}

\vspace{5em}
\begin{proof}
\emph{1.} Left to the reader.

\emph{2.} 
We start by showing the inclusion $\Im(\rja)\subseteq \Ker\lambda$, \cad for each~$k$,

\snic{\lambda\big(\som_{i=1}^{n}\Dpp{f_k}{X_i}(\ux)\,e_i\big)= 0.}

For $g\in\kuX$ we use Taylor's theorem at order $1$:
\snic{g(\uy)\equiv g(\uz)+\som_{i=1}^{n}\Dpp{g}{X_i}(\uz)\,(y_i-z_i) \mod \rJ^{2}.
}

For $g=f_k$ we have $f_k(\uy)=f_k(\uz)=0$, so $\som_{i=1}^{n}\Dpp{f_k}{X_i}(\uz)\,(y_i-z_i)\in\rJ^{2}$. This proves the above \egt by taking into account the \Amo law \hbox{over $\rJ/\rJ^{2}$}.
This shows that $\lambda$ passes to the quotient, with
$$\preskip.4em \postskip.4em 
 \ov\lambda:\delta(x_i)=\ov{e_i}\mapsto\rd(x_i)=\ov{y_i-z_i}.
$$
Moreover, since $\delta$ is a $\gk$-\dvnz, the \uvl \prt of the \dvn $\rd:\gA\to\rJ/\rJ^{2}$ gives us an $\gA$-\lin \fcn

\snic{\rJ/\rJ^{2} \to \Coker(\rja)\;:\;\rd(x_i)\mapsto \delta(x_i).}

It is clear that the two mappings are inverses of each other.
\end{proof}
%

\subsect{\Spt idempotent of a \asez}{\Spt idempotent of a \asez} 

Let~$\gA$ be a \stfe \klgz. For $a\in\gA$, let $a\sta =
a\centerdot\Tr\iAk$.
We have a canonical \kli $\env{\gk}{\gA} \to \End_\gk(\gA)$, composed of the \ali $\env{\gk}{\gA} \to \Asta\te_\gk\gA$, $a \otimes b \mapsto a\sta \otimes b$, and of the natural \iso $\Asta\te_\gk\gA \to \End_\gk(\gA)$.

If~$\gA$ is \ste these \alis are all \isosz.
Then, if $\big((x_i), (y_i) \big)$ is a \stycz, the \elt $\sum_i x_i\otimes y_i$ is independent of the choice of the \sys because its image in $\End_\gk(\gA)$ is~$\Id_\gA$. In particular, $\sum_i x_i\otimes y_i = \sum_i y_i\otimes x_i$. 

The following \tho identifies the \caras \prts of this \elt $\sum_i x_i\otimes y_i$.
These \prts lead to the notion of a \spb \algz.

\begin{theorem}\label{thAlgSteIdmSpb} \emph{(\Spt idempotent of a \asez)}\\
Let~$\gA$ be a \ste \klg and $\big((x_i), (y_i) \big)$ be a \styc of~$\gA$.  Then, the \elt $\vep = \sum_i x_i \te y_i \in \env\gk\gA$ satisfies the conditions of Lemma~\ref{lemmeIdempotentAnnJ}. In particular, $\vep$ is \idm and we have
$$\preskip.3em 
\som_i x_iy_i = 1, \quad  
a\cdot\vep = \vep\cdot a \quad \forall a \in \gA.
$$
\end{theorem}

\smallskip 
NB: We prove the converse (for \asfsz) a little later (\thref{thAlgStfSpbSte}).
\begin{Proof} {Proof in the Galoisian case (to be read after \thref{thA}). }\\
Let $(\gk,\gA,G)$ be a \aGz. Since the result to be proven is independent of the \stycz, we can suppose that the families $(x_i)$ and $(y_i)$ are two \syss of \elts of~$\gA$ satisfying the conditions of item \emph{2} of Artin's \thoz~\ref{thA}.

Saying that $\mu(\vep) = 1$ consists in saying that $\sum_i x_i y_i = 1$, which is what $\big((x_i), (y_i) \big)$ satisfies. To show that $\sum_i ax_i \otimes y_i = \sum_i x_i \otimes ay_i$, it suffices to apply $\psi_G$; we let $(g_\sigma)_\sigma$ be the image of the left-hand side, and $(d_\sigma)_\sigma$ be the image of the right-hand side. We obtain, by letting $\delta$ be the Kronecker symbol,
$$\preskip-.2em \postskip.2em\ndsp 
g_\sigma = \sum_i ax_i\sigma(y_i) = a \delta_{\sigma, \Id}, 
\qquad
d_\sigma = \sum_i x_i\sigma(ay_i) = \sigma(a) \delta_{\sigma, \Id}
. 
$$
We indeed have the \egt since the components of the two families $(d_\sigma)$ and $(g_\sigma)$ are null except at the index $\sigma = \Id$, at which their (common) \hbox{value is $a$}.
\\
Note that $\vep$ is equal to the \elt $\vep_{\Id}$ introduced in Lemma~\ref{lemArtin}. Its image under $\varphi_G$ is the \idm $e_\Id$, which confirms that $\vep$ is \idmz.
\end{Proof}

\vspace{-.3em}
\begin{Proof} {(\Gnlz) Proof in the \ste case. }\\
We write $\Tr$ for $\Tr_{\gA\sur\gk}$ and let $m_\vep : \env\gk\gA \to
\env\gk\gA$ be  multiplication by $\vep$. We have
$$\preskip-.3em \postskip.3em \ndsp
\Tr(ab) = \som_i \Tr(ay_i) \Tr(bx_i), \qquad a, b \in \gA.\eqno(\star)
$$
 Indeed, this easily results from the \egt $a = \sum_i \Tr(ay_i)x_i$.

We rewrite $(\star)$ as the \egt of two $\gk$-\lins forms,
$\env\gk\gA \to \gk$:
$$\preskip.3em \postskip.4em
\Tr_{\gA\sur\gk} \circ\, \mu\iAk = 
\Tr_{\env\gk\gA\sur\gk} \circ\, m_\vep.\eqno(*)
$$
Let us show that $\vep \in \Ann(\rJ)$. Let $z \in \env\gk\gA$,
$z' \in \rJ$. 
By evaluating the \egt $(*)$ at $zz'$, we obtain
$$
\preskip.3em \postskip.3em
\Tr_{\gA\sur\gk}\big(\mu\iAk(zz')\big) = 
\Tr_{\env\gk\gA\sur\gk}(\vep zz').
$$
But $\mu\iAk(zz') = \mu\iAk(z)\mu\iAk(z') = 0$ because $z' \in \rJ = \Ker\mu\iAk$. 
We deduce that $\Tr_{\env\gk\gA\sur\gk}(\vep zz') =
0$ for every $z \in \env\gk\gA$. As $\Tr_{\env\gk\gA\sur\gk}$ is non-\dgne
we obtain $\vep z' = 0$. 
Thus $\vep \in \Ann(\rJ)$.

It remains to show that $\mu\iAk(\vep) = 1$, \cad $s=\sum_i x_iy_i = 1$. 
\\
The \egt 
$\Tr(x) = \sum_i \Tr(xx_iy_i)$ (Fact~\ref{factMatriceEndo})
 says that $\Tr \big((1-s)x \big) = 0$
for all $x \in \gA$, thus $s = 1$.
\end{Proof}

\subsec{\Spb \algsz} \label{subsecAlgSpb}

\begin{theorem}\label{thSepIversen}
For a \klgz~$\gA$ \propeq
\begin{enumerate}
\item \label{i1thSepIversen}
$\gA$ is projective as an $\env\gk\gA$-module.
\item \label{i2thSepIversen}
$\rJ\iAk$ is generated by an \idm of $\env\gk\gA$.
\item \label{i4thSepIversen}
$\rJ\iAk$ is \tf and \idmz.
\item \label{i5thSepIversen}
$1 \in\mu\iAk\big(\Ann(\rJ\iAk)\big)$.
\item \label{i6thSepIversen}
There exist an $n\in\NN$ and $x_1$, \dots,  $x_n$,  $y_1$, \dots,  $y_n\in\gA$ such that
$\sum_ix_iy_i=1$ and for every $a\in\gA$, $\sum_iax_i\otimes y_i=\sum_ix_i\otimes ay_i$. 
\end{enumerate}
In this case we denote by $\vep\iAk$ the unique \idm which generates the \id $\Ann(\rJ\iAk).$
\\
When $\gA$ is a  \tf \klgz, another equivalent \prt is
\begin{enumerate} \setcounter{enumi}{5}
\item \label{i7thSepIversen}   
$\Om{\gk}{\gA}=0$.\footnote{By \Thref{thDerivUnivPF}, if $\gA=\aqo{\kXn}{\lfs}=\gk[\ux]$, $\Om{\gk}{\gA}=0$ means that
the matrix $\rja(\ux)$ (the transposed of the Jacobian matrix) is surjective, \cad $1\in\cD_n(\JJ_{\uX}(\uf)(\ux))$.}
\end{enumerate}
\end{theorem}
\begin{proof}
Since~$\gA\simeq\env\gk\gA\sur{\rJ\iAk}$, items~\emph{\ref{i1thSepIversen}} and~\emph{\ref{i2thSepIversen}} are \eqvs under Lemma~\ref{lemIdpPtf} regarding the cyclic \pro modules. Items~\emph{\ref{i2thSepIversen}} and~\emph{\ref{i4thSepIversen}} are \eqvs under Lemma~\ref{lem.ide.idem} on \tf \idm \idsz.
Lemma~\ref{lemmeIdempotentAnnJ} gives the \eqvc of \emph{\ref{i2thSepIversen}} and \emph{\ref{i5thSepIversen}}.\\
\emph{\ref{i4thSepIversen}} $\Leftrightarrow$ \emph{\ref{i7thSepIversen}.} If~$\gA$ is a \tf \klgz, then $\rJ\iAk$ is a \itf of~$\env\gk\gA$, therefore condition~\emph{\ref{i4thSepIversen}} 
can be reduced to: $\rJ\iAk$ is \idmz, \cad $\Om{\gk}{\gA}=0$.\\ 
Finally, \emph{\ref{i6thSepIversen}} is the concrete form of \emph{\ref{i5thSepIversen}.}
\end{proof}
%

\begin{definition}\label{defiSpb}
We call  an algebra that satisfies the equivalent properties stated in
\thref{thSepIversen}
a \emph{\spb \algz}. The \idm $\vep\iAk\in\env\gk\gA$ is called \emph{the \idstz} of~$\gA$.
\index{idempotent!separability ---}
\index{algebra!separable ---}
\index{separable!algebra}
\end{definition}

\comm
It should be noted that Bourbaki 
uses a notion of separable extension for fields that is quite different to the above definition. 
In \clamaz, \algs over a field~$\gK$  \gui{\spbs in the sense of \Dfnz~\ref{defiSpb}} are the \algs that are \gui{finite and \spbs in the Bourbaki sense} (see \thref{thSepProjFi}). Many authors follow Bourbaki at least for the \agq extensions of fields, whether they are finite or not.
In the case of an \agq \Klg over a discrete field~$\gK$, the \dfn \`a la Bourbaki means that every \elt of the \alg is a zero of a \spl \polu of~$\KT$.
\eoe

\begin{fact}\label{factSpbEds} \emph{(Stability of \spb \algs by \edsz)}\\
Let $\imath:\gk\to\gA$ and $\rho:\gk\to\gk'$ be two \klgs and~$\gA'=\rho\ist(\gA)$. 
We have a canonical \iso $\rho\ist(\env\gk\gA)\to 
{\gA'_{\,\gk'}}^{\!\!\mathrm{e}}$ and the diagram below commutes 

\vspace{-1.3em}
\carre{\env\gk\gA}{\env\gk\rho}{{\gA'_{\,\gk'}}^{\!\!\mathrm{e}}}
{\mu\iAk}{\mu_{\gA'\!/\gk'}}{\gA}{\rho}{\gA'}

\vspace{-1.1em}
In particular, a \spb \alg remains \spb by \edsz.
\end{fact}

\smallskip 
\facile

\medskip 
Now we prove the converse of \thref{thAlgSteIdmSpb},
which requires a preliminary lemma.

\begin {lemma}\label{lemTracAlgEnv}
Let~$\gA$ be a \stfe \klg and $\env\gk\gA $ its enveloping \algz.
\begin{enumerate}
\item
$\env\gk\gA$ is a \stfe left \Alg 
whose trace is given by~$\gamma_{\rm l} \circ (\Id_\gA \te \Tr\iAk)$
(where~$\gamma_{\rm l}:\gA\te_\gk\gk\to\gA$ is the canonical \isoz), \cad for $\alpha = \sum_i a_i\te b_i$:
$$\preskip.0em \postskip.4em\ndsp 
 \Tr_{(\env\gk\gA\sur\gA)_{\rm l}}(\alpha) = 
\sum_i a_i\Tr\iAk(b_i).
$$
Similarly, 
$\env\gk\gA$ is a \stfe right \Alg whose trace is given by~$\gamma_{\rm r}\circ (\Tr\iAk\te\Id_\gA)$, \cad
$
\Tr_{(\env\gk\gA\sur\gA)_{\rm r}}(\alpha) = 
\sum_i \Tr\iAk(a_i) b_i.
$

\item
Over $\Ann(\rJ\iAk)$, the~$\gA$-\lins forms $\Tr_{(\env\gk\gA\sur\gA)_{\rm l}}$, $\Tr_{(\env\gk\gA\sur\gA)_{\rm r}}$ and $\mu\iAk$ coincide, \cad if $\alpha = \sum_i a_i\te b_i \in
\Ann(\rJ\iAk)$,
$$\preskip.4em \postskip.4em\ndsp 
 \sum_i a_ib_i = \sum_i a_i\Tr\iAk(b_i) = \sum_i \Tr\iAk(a_i) b_i.
$$
\end{enumerate}

\end{lemma}

\smallskip 
\begin{proof}
\emph{1.}
This is a \gnl structural result: the trace is preserved by \eds (see Fact~\ref{fact.det loc}).
In other words if~$\gk'$ is a \klgz,~$\gk' \te_\gk \gA$ is a \stfe $\gk'$-\alg whose trace is~$\gamma\circ (\Id_{\gk'} \otimes \Tr\iAk)$ where~$\gamma:\gk'\te_\gk \gk\to\gk'$ is the canonical \isoz.

\emph{2.}
\Gnltz, under the hypotheses that $E$ is a \ptf \Amoz, $x \in E$, $\nu \in E\sta$ and $u=\theta_{E}(\nu\te x) \in \End_\gA(E)$, we obtain the \egt $\Tr_E(u)= \nu(x)$
(see Fact~\ref{factMatriceEndo}). 
\\
We apply this to $E = \env\gk\gA$, $x = \alpha\in E$ and $\nu = \mu\iAk\in E\sta$, by noting then that $u = \theta_{E}(\nu\te \alpha) = \mu_{\env\gk\gA,\alpha}$. Indeed, by item \emph {3}  of Fact~\ref{factOmAbsrait}, we~have 
\hbox{for $\gamma \in \env\gk\gA$,} $\gamma\alpha=\mu\iAk(\gamma)\cdot\alpha=\theta_{E}(\nu\te \alpha)(\gamma)$.
\end{proof}

\vspace{.01em} 
\begin{theorem}\label{thAlgStfSpbSte}\emph{(\Ste \algs and \spbs \algsz)}\\
Every \spb and \stfe \klgz~$\gA$ is \stez. More \prmtz, if $\vep\iAk = \sum x_i\te y_i \in \env\gk\gA $ is the \idst of~$\gA$, then $\big((x_i), (y_i) \big)$ is a \styc of~$\gA\sur\gk$.
\\
In brief, a \stfe \alg is \spb \ssi it is \stez.
\end{theorem}
NB: \Prmtz, the link between the two notions is obtained by the relation 
linking the \idst and the \sycsz, as is apparent in the direct \thrf{thAlgSteIdmSpb} and in the converse \thoz.

\begin {proof}
Let $x \in \gA$, then $(x\te1)\vep\iAk = \sum_i xx_i \te y_i$ is in $\Ann(\rJ\iAk)$, so by Lemma~\ref{lemTracAlgEnv}, we have $\sum_i xx_iy_i = \sum_i \Tr\iAk(x_ix)y_i$.\\
As $\sum_i x_iy_i = 1$, this gives $x = \sum_i \Tr\iAk(x_ix)y_i$. The result follows by the \carn of  \ases given in Fact~\ref{factCarAste}.
\end {proof}

The following \tho strengthens the previous \tho and shows that the existence of a \idst is a very strong condition of finiteness.

\pagebreak	

\begin{theorem}\label{thSepProjFi}
Let~$\gA$ be a \spb \klgz. \\
Suppose that~$\gA$ has a \ixe{\sycz}{coordinate system} in the following sense. We have a discrete set $I$, a family $(a_i)_{i\in I}$ in~$\gA$ and a family~$(\alpha_i)_{i\in I}$ in the dual \kmo $\Asta=\Lin_\gk(\gA,\gk)$, such that for all $x\in\gA$ we have

\snic{x=\sum_{i\in J_x}\alpha_i(x)a_i.}

Here $J_x$ is a finite subset of $I$, and every $\alpha_i(x)$ for $i\in I\setminus J_x$ is null.
\\
Then,~$\gA$ is \stfez, therefore \stez.%
\\
This is the case, for example, if $\gk$ is a \cdi and if $\gA$ is a \pf \klgz.
\end{theorem}
%
\begin{proof} 
Regarding the special case, the quotient \alg has a finite or countable basis of \momsz, by the theory of \bdgsz.
\\
Let $\vep=\sum_{k=1}^rb_k\te c_k$ be the \idstz.
We have $\vep\cdot x= x\cdot \vep$ for every $x\in\gA$,  and $\sum_{k=1}^rb_kc_k=1$.
\\
For $\alpha \in \Asta$ and $x \in \gA$, by applying $1\te \alpha$ to~$x \cdot \vep  = \vep \cdot x$ we obtain 

\snic{\som_k x b_k \alpha(c_k) = \som_k b_k \alpha(xc_k).
}

By denoting by $J$ the finite subset $J = \bigcup J_{c_k}$, we obtain for each~$k$
$$\preskip.4em \postskip.2em 
 c_k = \som_{i \in J} \alpha_i(c_k) a_i.
$$
We then write
$$\preskip.2em \postskip.4em  \mathrigid 2mu 
 x = \sum_{k \in \lrbr} \!\!xb_kc_k \,= \sum_{k \in \lrbr, i \in J} 
x b_k \alpha_i(c_k) a_i\, =\!\! \sum_{i\in J, k\in\lrbr} \!\!\alpha_i(c_kx)b_ka_i.
$$
%
This now gives a finite \syc for $\gA$, with the \elts~$b_ka_i$ and the forms $x\mapsto \alpha_i(c_kx)$ for $(i,k)\in J\times \lrbr$.
\end{proof}

\comm
Note that, when we have a \syc for a module, the module is \pro in the usual sense.
The \dfn of a \syc for a module $M$ amounts to saying that $M$ is \isoc to a direct summand 
of the module $\Ae{(I)}$. The latter module, freely generated by $I$, is \pro because $I$ is discrete.
\\
\emph{In \clamaz}, every \pro module has a \sycz, because all the sets are discrete, so the previous \tho applies: every \spb \klg which is a \pro \kmo is \stfez. By the same token \emph{every \spb \alg over a \cdi or over a \zedr \ri is \stfez}. 
\eoe

In the case of a \apf over a \cdiz, \thrfs{thSepIversen}{thSepProjFi} give the following result.

\begin{corollary}\label{corthSepProjFicdi}
For $f_1$, \dots, $f_s\in\kXn$ when $\gk$ is a \cdiz, \propeq
\begin{enumerate}
\item The quotient \alg $\gA=\kux$ is \stez.
\item The quotient \alg is \spbz.
\item The matrix $\rja(\ux)$, transposed of the Jacobian matrix of the \sypz, is surjective.
\end{enumerate} 
\end{corollary}

We will now show that a \spb \alg looks a lot like a diagonal \algz, 
including  
when the base \ri is arbitrary.
\\
Consider the diagonal \klgz~$\gk^n$. Let $(e_1, \ldots, e_n)$ be its canonical basis and~$p_i : \gk^n \to \gk$ be the \coo form in relation to~$e_i$.  Then we have 

\snic{e_i\in\BB(\gk^n)$, $p_i \in \Hom_\gk(\gk^n,
\gk)$,  $p_i(e_i) = 1$ and $xe_i = p_i(x)e_i \;\Tt x \in \gk^n.}

In a way, we are about to generalize the above result to \spb \algsz.


\begin {lemma}\label{lemIdmHomSpb}\emph{(Characters of a \spb \algz)}\\
Let~$\gA$ be a \spb \klg with~$\gk\subseteq\gA$.
\begin{enumerate}
\item Let $\imath:\gk\to\gA$ be the canonical injection.
If $\varphi \in \Hom_\gk(\gA,\gk)$, 
  $\imath\circ\varphi$ is a \prr with image~$\gk.1$,
so

\snic{\gA= \gk.1 \oplus\Ker\varphi \;\hbox{ and }\; \Im(\Id_\gA - \imath\circ \varphi)=\Ker\varphi.}

In fact the \id $\Ker\varphi$ is generated by an \idm of~$\gA$. We will denote by~$\vep_\varphi$ the \cop \idmz.
\item
For $\varphi$, $\varphi' \in \Hom_\gk(\gA,\gk)$, we have
$\varphi'(\vep_\varphi) = \varphi(\vep_{\varphi'})$.
\\
This \eltz, denoted by~$e_{\{\varphi,\varphi'\}}$, is an \idm of~$\gk$ and we have
\vspace{-1mm}
\[ 
\begin{array}{ccc} 
\vep_\varphi\vep_{\varphi'} =
e_{\{\varphi,\varphi'\}} \vep_\varphi =
e_{\{\varphi,\varphi'\}} \vep_{\varphi'}=\varphi(\vep_\varphi\vep_{\varphi'})=
\varphi'(\vep_\varphi\vep_{\varphi'}),     \\[2mm] 
\gen {\Im(\varphi - \varphi')}_\gk = \gen {1-e_{\{\varphi,\varphi'\}}}_\gk\; \;\hbox{and}
\;\;\Ann_\gk(\varphi - \varphi') = \gen {e_{\{\varphi,\varphi'\}}}_\gk.
 \end{array}
\]

\item
Consequently we have the \eqvcs

\snic {
e_{\{\varphi,\varphi'\}} = 1 \iff
\vep_\varphi = \vep_{\varphi'} \iff \varphi = \varphi',\,
$  and $\,
e_{\{\varphi,\varphi'\}} = 0 \iff \vep_\varphi \vep_{\varphi'} = 0.
}

\vspace{1mm}
\item
If~$\gk$ is connected, two \idms $\vep_\varphi$, (for $\varphi \in
\Hom_\gk(\gA, \gk)$), are equal or \ortsz.
\end{enumerate}
\end{lemma}
%
\begin {proof}
Let $\vep\iAk = \sum x_i \otimes y_i$. 
We know that $a\cdot\vep\iAk = \vep\iAk\cdot a$ for every $a \in \gA$, that $\sum x_i \otimes y_i = \sum y_i \otimes x_i$ and that $\sum_i x_iy_i = 1$.

 \emph{1.} The first assertion is valid for every \crc of every \algz~$\gA.$
It remains to see that $\Ker\varphi$ is generated by an \idmz.
We consider the \homo of \klgs $\nu=\mu\iAk\circ(\varphi\te\Id_\gA):\env\gk\gA\to\gA$, and the \elt $\vep=\nu(\vep\iAk)$. Thus $\vep=\som_i \varphi(x_i)y_i$ is an \idm and we obtain the \egts

\snic{\varphi(\vep)=\som_i \varphi(x_i)\varphi(y_i)=\varphi(\som_ix_iy_i)=\varphi(1)=1.}

Therefore $1-\vep\in\Ker\varphi$.
\\
By applying $\nu$ to the \egt $\som_i ax_i \otimes y_i = \som_i x_i \otimes ay_i$, we obtain $\varphi(a)\vep = a\vep$.
Therefore $a\in\Ker\varphi$ implies $a=(1-\vep )a$, and $\Ker\varphi=\gen{1-\vep}$.

 \emph{2.}
We have, for $a \in \gA$,
$$
\varphi'(a) \varphi'(\vep_\varphi) = \varphi'(a \vep_\varphi) = 
\varphi'(\varphi(a) \vep_\varphi) = \varphi(a)\varphi'(\vep_\varphi).
\eqno (\star)
$$
For $a = \vep_{\varphi'}$, we obtain $\varphi'(\vep_\varphi) = \varphi(\vep_{\varphi'})\varphi'(\vep_\varphi)$. 
By symmetry,
$\varphi(\vep_{\varphi'}) = \varphi'(\vep_\varphi)$. 
Denote by $e$ this \idm of~$\gk$. By \dfnz, we have $a\vep_\varphi = \varphi(a) \vep_\varphi$. By making $a = \vep_{\varphi'}$, we obtain $\vep_{\varphi'}\vep_\varphi = e\vep_\varphi$.
\\ 
Finally, let $\fa = \gen {\Im(\varphi - \varphi')}$. 
The relation  $(\star)$ 
shows that $\fa e = 0$. Moreover $1-e = (\varphi - \varphi')(\vep_\varphi) \in \fa$. Therefore $\fa = \gen {1-e}_\gk$ and $\Ann_\gk(\fa) = \gen {e}_\gk$.

 \emph{3} and \emph{4.}
Result from the previous item.
\end{proof}

\begin{lemma} \label{lemSspbdiag}
 \emph{(\Spb sub\alg of a diagonal extension)}
\\
Let~$\gk$ be a \emph{nontrivial connected} \riz,~$\gB=\gk^n$, $p_i:\gB\to\gk$ be the $i^{\rm th}$ canonical \prnz, $e_i$ be the \idm defined by $\Ker p_i=\gen{1-e_i}$ $(i\in\lrbn)$.
For a finite subset $I$ of $\lrbn$ we let $e_I = \sum_{i \in I} e_i$.\\
Let~$\gA$ be a \spb \klg with $\gk\subseteq\gA\subseteq\gk^n$ and $\pi_i$ be the restriction of $p_i$ to~$\gA$ for $i \in \lrbn$.  
\begin{enumerate}
\item We consider the \eqvc relation over $\lrbn$ defined by $\pi_i = \pi_j$.
The corresponding partition  $\cP$ 
is a finite set of finite subsets of~$\lrbn$. 
For $J\in\cP$ we denote by $\pi_J$ the common value of the $\pi_j$'s for $j \in J$.
\item
$\gA$ is a free \kmo with basis 
$\sotq{e_J}{J \in \cP}$.
\item $\Asta$ is a free \kmo with basis $\sotq{\pi_J}{J \in \cP}=\Hom_\gk(\gA, \gk)$.
\end{enumerate}
\end{lemma}

\smallskip 
\begin{proof}
\emph {1.} As~$\gk$ is nontrivial and connected, every \idm of~$\gB$ is of the form $e_I$ for a unique finite subset $I$ of $\lrbn$.
\\
 Let $i \in \lrbn$. 
By Lemma~\ref{lemIdmHomSpb} there exists one and only one \idmz~$\vep_i$ of~$\gA$ such that $\pi_i(\vep_i) = 1$ and $a \vep_i = \pi_i(a)\vep_i$ for every $a\in\gA$. 
This \idm is also an \idm of~$\gB$ so of the form $e_{J_i}$ for a finite subset~$J_i$ of $\lrbn$. 
Since $\pi_i(\vep_i)=p_i(e_{J_i}) = 1$, we have $i \in J_i$, and the union of the~$J_i$'s is~$\lrbn$.
Two distinct $J_i$ are disjoint by the last item of Lemma~\ref{lemIdmHomSpb}.
The~$J_i$'s therefore form a finite partition formed of finite subsets of~$\lrbn$.
\\
If $\pi_i=\pi_j$, then $\vep_i=\vep_j$ so $J_i=J_j$.
If $J_i=J_j$, then $\vep_i=\vep_j$ and $\pi_i(\vep_j)=1$.
Item \emph{2} of Lemma~\ref{lemIdmHomSpb} gives $1\in\Ann_\gA(\pi_i-\pi_j)$,
so $\pi_i=\pi_j$.

\emph{2.} Results from item \emph{1.}

\emph{3.} 
Let $\varphi \in \Hom_\gk(\gA,\gk)$. The $\varphi(e_J)$'s are \idms of~$\gk$. As~$\gk$ is connected, we have $\varphi(e_J) = 0$ or $1$. But the $(e_J)_{J \in \cP}$'s form a \sfioz, therefore there is only one $J \in \cP$ for which~$\varphi(e_J) = 1$ and consequently $\varphi = \pi_J$. The rest is \imdz.
\end{proof}

\section{\aGsz, \gnl theory}\label{secAGTG}

The theory developed by Artin
considers 
a finite group $G$ of \autos of a \cdi $\gL$, 
calls~$\gK$ the subfield of the fixed points of $G$ and 
proves that~$\gL$ is a Galois extension of~$\gK$, with $G$ as the
 Galois group.
 
In the current section we give the \gnn of Artin's theory for commutative \ris instead of \cdisz.
A good idea of \gui{how this can work} is already given by the following significant small example, which shows that the hypothesis \gui{\cdiz} is not required.

\mni
\textbf{A small example to start off}
\\
Let~$\gA$ be a commutative \riz, $\sigma \in \Aut(\gA)$ be an \auto 
of order~$3$, and $G$ be the group that it generates. Suppose that there exists an $x\in\gA$  such that $\sigma(x) - x \in \gA^{\times}$. Let~$\gk = \gA^{G}$ be the sub\ri of fixed points. Then, $(1, x, x^2)$ is a basis of~$\gA$ over~$\gk$. Indeed, let $V$ be the Vandermonde matrix
$$
V = \cmatrix {1 & x & x^2 \cr
1 & \sigma(x) & \sigma(x^2) \cr
1 & \sigma^2(x) & \sigma^2(x^2) \cr}
= \cmatrix {1 & x_0 & x_0^2 \cr
1 & x_1 & x_1^2 \cr
1 & x_2 & x_2^2 \cr}
\quad \hbox {with} \quad x_i = \sigma^i(x)
.$$
Let $\varepsilon = \sigma(x) - x$.
Then, $\det(V) = (x_1-x_0)(x_2-x_1)(x_2-x_0)$ is \ivz:

\snic{\det(V) =
\big(\sigma(x)-x\big) \cdot \sigma\big(\sigma(x)-x\big) \cdot \sigma^2\big(x-\sigma(x)\big)
= - \varepsilon \sigma(\varepsilon) \sigma^2(\varepsilon).}

For $y \in \gA$, we want to write $y=\lambda_0+\lambda_1x+\lambda_2x^2$ with each $\lambda_i\in \gk$. We then \ncrt have 
$$
\cmatrix {y\cr \sigma(y)\cr \sigma^2(y)\cr} =
\cmatrix {1 & x & x^2 \cr
1 & \sigma(x) & \sigma(x^2) \cr
1 & \sigma^2(x) & \sigma^2(x^2) \cr}
\cmatrix {\lambda_0\cr \lambda_1\cr \lambda_2\cr}
.$$
However, the above \sli has one and only one solution in~$\gA$. Since the solution is unique, $\sigma(\lambda_i) = \lambda_i$, \cad $\lambda_i \in \gk$ ($i=0,1,2$). \\
Finally, $(1, x, x^2)$ is indeed a~$\gk$-basis of~$\gA$.
\eoe

\subsec{Galois correspondence, obvious facts}

\vspace{3pt}
This can be considered as a resumption of Proposition~\ref{defiCorGal}.
\vspace{-1pt}

\begin{fact}
\label{lemGal1bis} \emph{(Galois correspondence, obvious facts)}\\
Consider a finite group $G$ of \autos of a \riz~$\gA$.
We use the notations defined in~\ref{NOTAStStp}. In particular, $\Ae H=\Fix_\gA(H)$ for a subgroup~$H$ of $G$. Let $\gk=\gA^{G}$.
\begin{enumerate}
\item  If $H\subseteq H'$ are two subgroups of $G$, then $\Ae H\supseteq\Ae {H'}$, and if $H$ is the subgroup generated by $H_1\cup H_2$, then $\Ae {H}= \Ae {H_1}\cap \Ae {H_2}$.
\item
$H\subseteq \Stp(\Ae H)$ for every subgroup $H$ of $G$.
\item If $\sigma\in G$ and $H$ is a subgroup of $G$
then

\snic{\sigma(\Ae H)=\Ae {\sigma H\sigma^{-1}}.}
\item  If~$\gC \subseteq \gC'$ are two \kslgs of~$\gA$, then $\Stp(\gC)\supseteq\Stp(\gC')$, and if~$\gC$ is the \kslg generated by~$\gC_1\cup \gC_2$, then 

\snic{\Stp(\gC)= \Stp(\gC_1)\cap \Stp(\gC_2).}

\item $\gC\subseteq \Ae {\Stp(\gC)}$ for every \kslgz~$\gC$ of~$\gA$.
\item After any \gui{go-come-go motion,} we end up with the resulting set of the first \gui{go}:

\snic{    \Ae H = \Ae {\Stp(\Ae H)} \; \et \; \Stp(\gC)= \Stp\big(\Ae {\Stp(\gC)}\big).}
\end{enumerate}
\end{fact}
\begin{proof}
The last item is a direct consequence of the previous ones, which are \imdsz. Likewise for all the \gui{dualities} of this type.
\end{proof}

\vspace{-.7em}
\pagebreak	

\subsec{A natural \dfnz}

Let $\cG=\cG_G$ be the set of finite (\cad detachable) subgroups of $G$, \hbox{and $\cA=\cA_G$} be the set of sub\ris of~$\gA$ which are of the form $\Fix(H)$ for some \hbox{$H\in\cG$}. Consider the restrictions of $\Fix$ and $\Stp$ to the sets $\cG$ and $\cA$. 
We are interested in determining under which conditions we thus obtain two inverse bijections between $\cG$ and $\cA$, and in giving a nice \carn of 
sub\algs belonging to $\cA$.
In the case where~$\gA$ is a \cdiz, Artin's theory shows that we find ourselves in the classical Galois situation:~$\gA$ is a Galois extension of the subfield~$\gk=\Ae G$, $G$ is the Galois group of this extension and $\cA$ is the set of all the \stfe subextensions of~$\gA$.

This \gui{Artin-Galois} theory has then been \gnee to an arbitrary commutative \riz~$\gA$, where certain conditions are imposed on the group $G$ and on the \kslgs of~$\gA$.

Actually, we want the corresponding notion of a \aG to be sufficiently stable.
In particular,   when we replace~$\gk$ by a nontrivial quotient~$\gk\sur\fa$ and $\gA$ by~$\gA\sur{\fa\gA}$, we wish to maintain the notion of a \aGz. Therefore two \autos of~$\gA$ present in $G$ must not be able to become a single \auto upon passage to the quotient.

This leads to the following \dfnz.

\begin{definition}\label{defaG}
\emph{(Well-separated maps, separating  \autosz, \aGsz)}
\index{algebra!Galois ---}
\index{well-separated!maps} 
\begin{enumerate}
\item
Two maps $\sigma$, $\sigma'$ from a set $E$ to a \riz~$\gA$ are said to be \emph{well-separated} 
if

\snic{\gen{\,\sigma(x)-\sigma'(x)\;;\;x\in E\,}_{\!\gA}=\gen{1}.}

\item An \auto $\tau$ of~$\gA$ is said to be \ixc{separating}{\autoz} if it is well-separated from~$\Id_\gA$.

\item A finite group $G$ that operates on~$\gA$ is said to be \ixc{separating}{group of \autosz} if the \elts $\sigma\neq 1_G$ of $G$ are separating (it amounts to the same to say that every pair of distinct \elts of $G$ gives two well-separated \autosz).
\\
We will also say that $G$ operates \emph{in a separating way} on~$\gA$.
\item A \emph{\aGz} is by \dfn a triple $(\gk,\gA,G)$, 
where~$\gA$ is a \riz, $G$ is a finite group operating on~$\gA$ in a separating way, and~$\gk=\Fix(G)$.
\end{enumerate}
\end{definition}

\comms 
\\
1) As for the \dfn of a \aGz, we did not want to forbid a finite group operating on the trivial \riz,
 and consequently we do not define $G$ as a group of \autos of~$\gA$, but as a finite group operating on~$\gA$.%
\footnote{The unique \auto of the trivial \ri is separating, and every finite group operates on the trivial \ri in order to make it a \aGz.}
In fact, the \dfn implies that $G$ always operates faithfully on~$\gA$ (and thus can be identified with a subgroup of $\Aut(\gA)$) except in the case where the \ri is trivial. This presents several advantages.
\\
On the one hand, a \aG remains Galoisian, \emph{with the same group~$G$}, for every \edsz; it is possible that we do not know if a \edsz~$\gk\to\gk'$, appearing in the middle of a \demz, is trivial or not. \\
On the other hand, the fact of not changing groups is more convenient for any \eds anyway.

2) We have imposed the condition $\gk\subseteq \gA$, which is not in the usual categorical style. The readers will be able to restore a more categorical definition, if they wish, by saying by saying that the morphism $\gk\to\gA$ establishes an \iso between $\gk$ and $\gA^G$. This will sometimes be \ncrz, for example in item~\emph{2} of Fact~\ref{factEdsAG}.
\eoe

\exls \\
1) If~$\gL/\gK$ is a Galois extension of \cdisz,
then the triple 
$$\preskip.3em \postskip.0em 
\big(\gK,\gL,\Gal(\gL/\gK)\big) 
$$
 is a \aGz.

2) We will show a little further (\thref{thAduAGB}) that for a \spl \polu $f\in\kT$, the triple $(\gk,\Adu_{\gk,f},\Sn)$ is a \aGz.

3) An \auto $\sigma$ of a \alo $\gA$ is separating \ssi there exists some $x\in\gA$ such that $x-\sigma(x)$ is \ivz. 
\eoe

\medskip
The notions of a separating \auto and of a \aG have been developed in order to satisfy the following fundamental facts.

\begin{fact}\label{factEdsAG}~
\begin{enumerate}
\item A separating \auto $\sigma$ of a \riz~$\gA$ provides by \eds $\rho:\gA\to\gB$ a separating \auto $\rho\ist(\sigma)$ of~$\gB$.
\item If $(\gk,\gA,G)$ is a \aG and if $\rho:\gk\to\gk'$ is a \ri \homoz, then $(\gk',\rho\ist(\gA),G)$ is a \aGz.
\end{enumerate}

\end{fact}
\begin{proof}
Item \emph{1}, as well as item \emph{2} in the case of a \eds by \lonz, are easy and left to the reader. 
\\
The \dem of the \gnl case for item \emph{2} will have to wait until \thref{corAGextsca}.
\end{proof}

\begin{plcc}\label{plcc.aGs}\relax \emph{(\aGsz)}
\\
Let $G$ be a finite group operating on a \klgz~$\gA$ with $\gk\subseteq \gA$.\\
Let $S_1$, $\dots$, $S_n$ be \moco of~$\gk$.\\
Then, $(\gk,\gA,G)$ is a \aG \ssi each triple $(\gk_{S_i},\gA_{S_i},G)$ is a \aGz.
\end{plcc}

\facile

\vspace{-.7em}
\pagebreak	

\subsec{Dedekind's lemma}
\label{subsecLDAC}\index{Dedekind!lemma}

Let~$\gA$ be a commutative \riz. Consider the $m^{\rm th}$ power \Algz~$\Ae m$. Its \elts will be ragarded as column vectors and the laws are the product laws

\snic{\cmatrix{a_1 \cr\vdots\cr  a_m}\star\cmatrix{b_1 \cr\vdots\cr  b_m} =
\cmatrix{a_1\star b_1 \cr\vdots\cr  a_m\star b_m},\quad a \cmatrix{a_1 \cr\vdots\cr  a_m}=\cmatrix{aa_1 \cr\vdots\cr  aa_m}.}


\begin{lemma}
\label{lemDA}
Let $C$ be a finite subset of $\Ae m$ which \gui{separates the rows};
\cad~$\gen{x_i- x_j \,;\,x\in C}_\gA=\gen{1}$ (for $i\neq j\in\lrbm$).
Then, the \Alg generated by $C$ is equal to $\Ae m$.
\end{lemma}

\vspace{.1em}
\begin{proof}
The fundamental remark is that in the \Amo generated by $1_{\Ae m}$ and  $x=\tra{[\,x_1 \,\cdots\, x_m\,]}$ there are the vectors

\snic{x-x_2 \,1_{\Ae m}=\tra{[\,x_1-x_2\;0\;* \,\cdots\, *\,]}$
and 
$-x+x_1 \,1_{\Ae m}=\tra{[\,0\;x_1-x_2\;* \,\cdots\, *\,]}.}

Therefore, when we suppose that the \id generated by the $x_1-x_2$'s 
contains~$1$, this implies that in the \Amo generated by $C$ there is a vector $g^{1,2}$ of the type $\tra{[\,1\;0\;g^{1,2}_3 \,\cdots\,  g^{1,2}_m\,]}$ and a vector $g^{2,1}$ of the type $\tra{[\,0\;1\;g^{2,1}_3 \,\cdots\,  g^{2,1}_m\,]}$. The general case is similar replacing $1$ and $2$ with two integers $i\neq j\in\lrbm$.
\\
We deduce that $\tra{[\,1\;0\;0 \,\cdots\, 0\,]}= g^{1,2}\cdot g^{1,3}\cdots g^{1,m}$ is in the \Alg generated by $C$.
Similarly, each vector of the canonical basis of $\Ae m$ will be in the \Alg  generated by $C$.
We actually get that $\Ae m$ is the image of a matrix whose columns are the products of at most $m$ columns in~$C$.
\end{proof}
\begin{notations}
\label{notas2.1}\relax \emph{(Context of Dedekind's lemma)}
\begin{enumerate}\itemsep0pt
\item [--] $\gA$ is a commutative \riz.
\item [--] $(M,\cdot,1)$ is a \moz.
\item [--] $\tau=(\tau_1,\tau_2,\ldots,\tau_m)$ is a list of $m$ \homosz, pairwise well-separated, 
           of $(M,\cdot,1)$ in $(\gA,\cdot,1)$.
\item [--] For $z\in M$ we denote by $\tau(z)$ the \elt of $\Ae m$ defined by
$$\preskip.2em \postskip.0em 
\tau(z)=\tra{[\,\tau_1(z) \,\cdots\, \tau_m(z)\,]}. 
$$
\end{enumerate}
\end{notations}

\begin{theorem}
\label{thDA}
\emph{(Dedekind's lemma)}\\
Using the notations in~\ref{notas2.1} there exist $y_1$, \ldots, $y_r\in M$ such that the matrix
$$\preskip.4em \postskip.0em \ndsp
[\,\tau(y_1) \mid \cdots \mid \tau(y_r)\,]= \big(\tau_i(y_j)\big)_{i\in\lrbm,j\in\lrbr} 
$$
 is surjective. 
 
\emph{Weak form.} In particular,
$\tau_1$, \ldots, $\tau_m$ are~$\gA$-\lint independent.
\end{theorem}
%
\begin{proof}
This is deduced from Lemma~\ref{lemDA} by noting that, since $\tau(xy) = \tau(x)\tau(y)$, the \Alg generated by the $\tau(x)$ coincide with the \Amo generated by the~$\tau(x)$.
\end{proof}

\rems\\
 1) Let $F=\big(\tau_i(y_j)\big)_{ij}\in\Ae{m\times r}$.
The \lin independence of the rows means that $\cD_m(F)$ is faithful, whereas the surjectivity of $F$ means that $\cD_m(F)$ contains $1$. Sometimes, Dedekind's lemma is called \gui{Artin's theorem} 
or the \gui{independence of characters lemma,} when one has in view the case 
where~$\gA$ is a \cdiz. In fact, it is only when $\gA$ is a \zed \ri that we can deduce  \gui{$\cD_m(F)=\gen{1}$} from  
          \gui{$\cD_m(F)$ is faithful.}

2) The integer $r$  can be controlled from the data in the problem.
\eoe

\subsec{Artin's \tho and first consequences}

\begin{definota}
\label{NotaAGAL}
Let~$\gA$ be a \klg with $\gk\subseteq \gA$.
%
\begin{enumerate}\itemsep0pt
\item We can equip the \kmo $\Lin_\gk(\gA,\gA)$ with an \Amo structure  by the external law
$$\preskip.2em \postskip.2em \ndsp
(y,\varphi)\mapsto \big(x\mapsto y\varphi(x)\big), \quad\gA\times
\Lin_\gk(\gA,\gA)\to\Lin_\gk(\gA,\gA)\,.
$$
We then denote this \Amo by $\LIN_\gk(\gA,\gA)$.
\end{enumerate}
Let $G=\so{\sigma_1=\Id,\sigma_2, \ldots ,\sigma_n}$ be a finite group operating (by~$\gk$-\autosz) on~$\gA$.
\begin{enumerate}\itemsep0pt\setcounter{enumi}{1}
\item
The \Ali $\iota_G:\prod_{\sigma\in G}\gA \to \LIN_\gk(\gA,\gA)$ is defined by 
$$\preskip.2em \postskip.2em \ndsp
\iota_G \big((a_\sigma)_{\sigma\in G}\big)=\sum_{\sigma\in G} a_\sigma
\sigma\,. 
$$
\item
The \kli $\psi_G:\env\gk\gA \to \prod_{\sigma\in G}\gA$ is defined by
$$\preskip.2em \postskip.2em 
\psi_G(a\otimes b)=\big(a\sigma(b)\big)_{\sigma\in G}\,.
$$
This is a \homo of \Algs (on the left-hand side).
\perso{Il y a aussi une action of $G$ qui is conserv\'ee in l'\isoz, \`a
pr\'eciser if vraiment \ncrz.}
\end{enumerate}
\end{definota}

\begin{fact}\label{factdefAGAL} With the above notations, 
and the left-structure
for the~\Amo $\env\gk\gA$, we have the following results.
\begin{enumerate}\itemsep1pt
\item Saying that $\iota_G$ is an \iso means that $\LIN_\gk(\gA,\gA)$ is a free \Amo whose $G$ is a basis.
\item If~$\gA$ is \ste of constant rank over~$\gk$, saying that $\env{\gk}{\gA}$ is a free \Amo of finite rank means that~$\gA$ diagonalizes itself.
\item Saying that $\psi_G$ is an \iso means \prmt the following.\\
The \Amo $\env{\gk}{\gA}$ is free of rank $\#{G}$, with a basis $\cB$ such that, after \eds from~$\gk$ to~$\gA$, the \ali $\mu_{\gA,a}$, which has become $\mu_{\env{\gk}{\gA},1\te a}$, is now diagonal over the basis $\cB$, 
with matrix
$$\preskip.3em \postskip.1em 
 \Diag\big(\sigma_1(a),\sigma_2(a),\ldots,\sigma_n(a)\big)
$$
 for any $a\in\gA$.
\end{enumerate}
\end{fact}

\vspace{-.5em}
\pagebreak


\begin{lemma}\label{lemArtin}~\\
Let $G=\so{\sigma_1=\Id,\sigma_2, \ldots ,\sigma_n}$ be a finite group operating on a \riz~$\gA$ and  let $\gk=\Ae G$. 
For $y\in\gA$, let $y\sta$ be the \elt of $\Asta$ defined by $x\mapsto \Tr_G(xy)$. \Propeq
\begin{enumerate}
\item $(\gk,\gA,G)$ is a \aGz.
\item  There exist $x_1$, \ldots, $x_r$, $y_1$, \ldots, $y_r$ in~$\gA$ such that for every $\sigma\in G$ we have
\begin{equation}\label{thAeq4}\preskip-.6em \postskip.2em
\som_{i=1}^r x_i \sigma(y_i)= \formule{1\;\;\hbox{if}\;\sigma=\Id\\[.2em]
0\;\;\hbox{otherwise}\,.}
\end{equation}
\end{enumerate}
In this case we have the following results. 
\begin{enumerate} \setcounter {enumi}{2}
\item 
For $z\in\gA$, we have $z=\sum_{i=1}^r \Tr_G(zy_i)\, x_i = \sum_{i=1}^r \Tr_G(zx_i)\, y_i$. 
\\
In other words, $\gA$ is a \ptf \kmo 
and 
$$\preskip.2em \postskip.2em 
\big((\xr),(y_1\sta,\ldots,y_r\sta)\big)\;\hbox{  and  }\;\big((\yr),(x_1\sta,\ldots,x_r\sta)\big) 
$$
are \sycsz.
\item 
The form
$\Tr_G: \gA \rightarrow \gk$ is dualizing and surjective.
\item 
For $\sigma \in G$, let $\varepsilon_\sigma = \sum_i \sigma(x_i) \otimes
y_i \in \env\gk\gA $. Then,
$(\varepsilon_\sigma)_{\sigma \in G}$ is an~$\gA$-basis \gui{on the left-hand side} of $\env\gk\gA$. In addition, for $a$, $b \in \gA$, We have 
$$\preskip.2em \postskip.3em\ndsp 
b \otimes a = \sum_\sigma  b\sigma(a) \varepsilon_\sigma
, 
$$
and the image of this basis $(\varepsilon_\sigma)_\sigma$ under $\psi_G : \env\gk\gA \to \prod_{\tau \in G}\gA$ is the canonical $\gA$-basis $(e_\sigma)_{\sigma\in G}$ of $\prod_{\tau \in G}\gA$.
Consequently, $\psi_G$ is an \iso of \Algsz.

\end{enumerate}
\end{lemma}


\begin {proof}
 \emph{1} $\Rightarrow$ \emph{2.} 
By Dedekind's lemma, there exist an integer $r$ and \elts \hbox{$x_1$, \ldots, $x_r$, $y_1$, \ldots, $y_r\in\gA$} 
such that 
$$ \sum_{i=1}^r \,x_i\,\cmatrix{\sigma_1(y_i) \cr\sigma_2(y_i) \cr\vdots 
\cr\sigma_n(y_i)}=
\cmatrix{1\cr0\cr\vdots\cr0},
$$
meaning, for $\sigma \in G$, precisely \Eqnsz~(\ref{thAeq4}).

 \emph{2} $\Rightarrow$ \emph{1.}
For $\sigma \ne \Id$, we have $\sum_{i=1}^r x_i \big(y_i - \sigma(y_i)\big) = 1$, which proves that $\sigma$ is separating.

 \emph{3.}
For $z \in \gA$, we have the \egts
\[ \preskip.3em \postskip.3em \arraycolsep2pt
\begin{array}{rccccc} 
\sum_{i=1}^r \Tr_G(zy_i) \, x_i&=&\sum_{i=1}^r \sum_{j=1}^n \sigma_j(zy_i) x_i 
&=& \\[.2em] 
 \sum_{j=1}^n \sigma_j(z) \big( \sum_{i=1}^r \sigma_j(y_i)x_i \big)&=&
\sigma_1(z)\cdot 1+\sum_{j=2}^n\sigma_j(z)\cdot 0&=&z. 
\end{array}
\]

\emph{3} $\Rightarrow$ \emph{4.} By item \emph{1} of \thref{factCarDua}.

\emph {5.}
We have $\psi_G(\varepsilon_\sigma) = \big(\sum_i
\sigma(x_i)\tau(y_i)\big)_\tau = e_\sigma$. Let us now show the \egt with respect to $b\otimes a$. Given the chosen  \Amo structure  on the left-hand side, we can assume that $b = 1$. Then
$$\preskip.3em \postskip.3em\arraycolsep2pt
\begin{array} {rcl}
\sum_\sigma  \sigma(a) \varepsilon_\sigma &=&
\sum_\sigma  \sigma(a) \sum_i \sigma(x_i) \otimes y_i =
\sum_i \Tr_G(ax_i) \otimes y_i
\\[.2em] 
&=& \sum_i 1 \otimes \Tr_G(ax_i)y_i = 
     1 \otimes \sum_i \Tr_G(ax_i)y_i = 1 \otimes a.
\end{array}
$$
This shows that $(\varepsilon_\sigma)_\sigma$ is a \sgr of the \Amo $\env\gk\gA$.
As its image under $\psi_G$ is the canonical $\gA$-basis of $\prod_{\tau\in G}\gA$, this \sys is free over~$\gA$. The rest follows from this.
\end{proof}

\rem Here is an alternative \dem of the surjectivity of the trace 
(item~\emph{4}). For $z = 1$, $1 = \sum_{i=1}^r t_i x_i$ with $t_i = \Tr_G(y_i) \in \Tr_G(\gA) \subseteq \gk$. Let us introduce the \gui{normic} \pol $N(T_1, \ldots, T_r)$:

\snic{N(T_1, \ldots, T_r) =\rN_G \big(\som_{i=1}^r T_i x_i\big) = 
\prod\nolimits_{\sigma \in G} \big(T_1\sigma(x_1) + \cdots + T_r\sigma(x_r)\big).
}

It is a \hmg \pol of degree $n \ge 1$, invariant under $G$, therefore with \coes in $\gk$: $N(\und{T}) = \sum_{|\alpha|=n} \lambda_\alpha \und{T}^\alpha$ with $\lambda_\alpha \in \gk$. Consequently, for $u_1$, \ldots, $u_r \in \gk$, we have $N(u_1, \ldots, u_r) \in \gk u_1 + \cdots + \gk u_r$.
In particular 

\smallskip \centerline{~~~~{\small $1 =\rN_G(1)= \rN_G\big(\som_{i=1}^r t_i x_i\big)= N(t_1, \ldots, t_r) 
\in \gk t_1 + \cdots + \gk t_r \subseteq \Tr_G(\gA).$} \eoe}


\begin{theorem}
\label{thA}\emph{(Artin's \thoz, \aGs version)}\index{Artin!theorem} \\
Let $(\gk,\gA,G)$ be a \aG (notations in~\ref{NotaAGAL}).
\begin{enumerate}
\item The \kmoz~$\gA$ is \prc $\#G$, 
and~$\gk$ is a direct summand in~$\gA$.
\item
There exist $x_1$, \ldots, $x_r$ and $y_1$, \ldots, $y_r$ such that for all $\sigma$, $\tau\in G$ we have
\begin{equation}\label{thAeq5} \preskip-.2em \postskip.2em
\forall\sigma,\tau \in G\qquad \som_{i=1}^r \tau(x_i) \sigma(y_i)= 
\formule{1\;\;\mathrm{if}\;\sigma=\tau\\
0\;\;\mathrm{otherwise}.}
\end{equation}
\item The form $\Tr_G$ is dualizing.
\item
The map $\psi_G : \env{\gk}{\gA} \to \prod_{\sigma\in G} \gA$ is an \iso of \Algsz. In particular,~$\gA$ diagonalizes itself.
\item 
\begin{enumerate}
\item $\rC{G}(x)(T)=\rC{\gA\sur\gk}(x)(T)$,
 $\Tr_G=\Tr\iAk $ and
 $\rN_G=\rN\iAk $, 
\item $\gA$ is \ste over~$\gk$. 
\end{enumerate}
\item If~$\gA$ is a \cdiz, it is a Galois extension of~$\gk$,
and we have~$G=\Gal(\gA\sur\gk)$.
\end{enumerate}
\end{theorem}
\begin{proof}
In this \demz, for $x\in\gA$, we let $\Tr(x)=\Tr_G(x)$, and $x\sta$ is the~$\gk$-\lin form $z\mapsto\Tr(zx)$.

Lemma~\ref{lemArtin} proves items~\emph{1} (besides the rank question), \emph{3} and \emph{4.} It also proves item~\emph {2}, because~(\ref{thAeq5}) clearly results from (\ref{thAeq4}).
\\
Note that~$\gk$ is a direct summand  
in~$\gA$ by item~\emph{3} of Lemma~\ref{lemIRAdu}.%
\footnote{Or more directly, by the surjectivity of the trace (which results from \thref{factCarAste}~\emph{1}.
Indeed, let $x_0 \in \gA$ such that $\Tr(x_0) = 1$. We have~$\gA = \gk\cdot 1 \oplus \Ker x_0\sta$, because every $y \in \gA$ can be written as $y = x_0\sta(y)\cdot 1 + (y - x_0\sta(y)\cdot 1)$ with $y - x_0\sta(y)\cdot 1\in\Ker x_0\sta$.}

Let us see that~$\gA$ is indeed of constant rank $n$.
Item \emph{4} shows that, after \eds from~$\gk$ to~$\gA$, the \kmoz~$\gA$ becomes free of rank $\#{G}$. Thus~$\gA$ is indeed of constant rank $n$ over~$\gk$; the \polmu of the \kmoz~$\gA$ \gui{does not change} under the \edsz%
~$\gk\to\gA$ (injective), it is therefore itself equal to 
$T^n$.\footnote{Actually, its \coes are transformed into \gui{themselves,} viewed in~$\gA$.}

\emph{5a} (and so \emph{5b}) Since $\psi_G$ is an \iso of \Algs (item \emph{4}),~$\gA$ diagonalizes itself.
We then deduce from Fact~\ref{factdefAGAL} item~\emph{3}, the \egt 
$$
\preskip.3em \postskip.3em 
\rC{G}(x)(T)=\rC{\gA/\gk}(x)(T). 
$$
This is true for the \pols  in $\AT$, therefore also in $\kT$.

\emph{6.} First of all, the \riz~$\gk$ is \zed by Lemma~\ref{lemZrZr2}. It is therefore a \cdiz, because it is connected and reduced. The extension is \'etale. It is normal, because every $x\in\gA$ annihilates 
$\rC G(x)(T)$, and this \pol can be decomposed into a product of \lin factors in $\AT$.
\end{proof}
%

\rem The computation that follows can clarify things, despite it not being \ncrz.
\\
Note that by item \emph{3} of Lemma~\ref{lemArtin}, the \kmoz~$\gA$ is the image of the \mprn
$$\preskip.2em \postskip.4em
\;\;P=(p_{ij})_{i,j\in \lrbr}=\big(y_i\sta(x_j)\big)_{i,j\in \lrbr}=
\big(\Tr(y_ix_j)\big)_{i,j\in \lrbr}\,.
$$
Also recall \Eqrf{thAeq5}:
${\sum_{i=1}^r \tau(x_i) \sigma(y_i)=
\formule{1\;\;\mathrm{if}\;\sigma=\tau\\[.2em]
0\;\;\mathrm{otherwise}}.
}$

Then let 
\[\preskip-.4em \postskip.4em 
\begin{array}{cccccc} 
X&=&\cmatrix{\sigma_1(x_1)&\sigma_1(x_2)&\cdots& \sigma_1(x_r) \cr
\sigma_2(x_1)&\sigma_2(x_2)&\cdots &\sigma_2(x_r) \cr\vdots&\vdots&  &\vdots 
\cr
\sigma_n(x_1)&\sigma_n(x_2)&\cdots &\sigma_n(x_r) }
&\,\hbox{ and }  \\[2.5em] 
Y&=&~\cmatrix{\sigma_1(y_1)&\sigma_1(y_2)&\,\cdots\, &\sigma_1(y_r) \cr
\sigma_2(y_1)&\sigma_2(y_2)&\cdots &\sigma_2(y_r) \cr\vdots&\vdots&  &\vdots 
\cr
\sigma_n(y_1)&\sigma_n(y_2)&\cdots &\sigma_n(y_r) }.& 
\end{array}
\]
By \Eqn (\ref{thAeq5}), we have $X\tra{Y}=\In$ and $P=\tra{Y}X$. 
\\
By Proposition~\ref{propImProjLib}, this means that the \kmoz~$\gA$, 
becomes free of rank~$n$, with for basis the $n$ rows of $Y$, after \eds from~$\gk$ to~$\gA$. In other words, the \Amo $\env\gk\gA$, seen as an image of the matrix~$P$ \gui{with \coes in~$\gA$} is a free \Asub of rank~$n$ of $\Ae r$,
and it is a direct summand.
\eoe


\begin{corollary}\label{corAGfree} \emph{(Free \aGz)}\\
Let $(\gk,\gA,G)$ be a free \aGz, and $n=\#{G}$.
If $\ub= (\bn)$ in~$\gA$, we define $M_\ub \in \Mn(\gA)$ by 
${M_\ub = \big(\sigma_i(b_j)\big)_{i,j\in\lrbn}.}$

Then, for two \syss $\ub$, $\ub'$ of $n$ \elts of~$\gA$ we obtain

\snic {
\tra M_\ub\, M_{\ub'} = \Tr_G(b_ib'_j)_{i,j\in\lrbn}.
}

Consequently, we obtain the following results.
\begin{itemize}
\item $\det(M_\ub)^2 = \disc(b_1, \ldots, b_n)$.
\item The \sys $(\bn)$ is a~$\gk$-basis of~$\gA$ \ssi the matrix $M_\ub$ is \ivz.
\item  In this case, if $\ub'$ is the dual basis of $\ub$ with respect to the trace-valued bi\lin form, then the matrices $M_\ub$ and $M_{\ub'}$ are inverses of one another.
\end{itemize}    
\end{corollary}

\vspace{.1em}
\rem In the situation where $\gA$ is a \cdiz, {Dedekind's lemma} in its original form asserts that the \gui{Dedekind matrix} $M_\ub$ is \iv when~$(\ub)$ is a basis of $\gA$ as a \kevz.\eoe


\begin{theorem}\label{corAGextsca}\emph{(\Eds for \aGsz)}\\
Let $(\gk,\gA,G)$ be a \aGz, $\rho:\gk\to\gk'$ be an \alg and~$\gA'=\rho\ist(\gA)$.
\begin{enumerate}
\item The group $G$ operates naturally over~$\gA'$ and $(\gk',\gA',G)$ is a \aGz.
\item The \gui{Galois theory} of $(\gk',\gA',G)$ is deduced by \eds of that of $(\gk,\gA,G)$, in the following sense: for each finite subgroup $H$ of~$G$, the natural \homo $\rho\ist(\Ae H)\to \gA'^H$ is an \isoz.
\end{enumerate}
\end{theorem}
\begin{proof} \emph{1.}
We easily see that $G$ acts on~$\gA'$ in a separating way. It remains to show that~$\gk'$ is the sub\ri of $G$-invariant \elts of~$\gA'$.\\
Let $\Tr = \Tr_G$. We see $\Tr$ as a $\gk$-\endo of~$\gA$, which by \eds gives the $\gk'$-\endo  $\Id_{\gk'} \otimes \Tr$ of $\gA'$.\\
Let~$y \in \gA'^{G}$. For~$z \in \gA'$, since $y$ is~$G$-invariant, we have the \egt 

\snic{(\Id_{\gk'} \otimes \Tr)(yz) =
y\,(\Id_{\gk'} \otimes \Tr)(z).}

By taking $z_0 = 1_{\gk'} \otimes x_0$, where~$x_0 \in \gA$ satisfies $\Tr(x_0) = 1$, we obtain the desired membership

\snic {
y = (\Id_{\gk'} \otimes \Tr) (yz_0) \in
\gk'\otimes_\gk \gk = \gk'
.}

\sni \emph{2.} Results from item \emph{1.} Indeed, consider the \aG $(\Ae H,\gA,H)$ and the \eds  $\varphi:\Ae H\to \gk'\te_\gk \Ae H=\rho\ist(\Ae H)$. We obtain the \egt  $\varphi\ist(\gA)=\gA'$. So  $\big(\rho\ist(\Ae H),\gA',H\big)$ is a \aG and
$\gA'^H=\rho\ist(\Ae H)$. 
\end{proof}
%


In the following \thoz, we could have expressed the hypothesis by saying that the finite group $G$ operates over the \ri $\gA$, and that $\gk$ is a sub\ri of~$\Ae H$. 

\pagebreak	

\begin{theorem}\label{thAGACar}\emph{(\Carns of  \aGsz)}\\
Let $G$ be a finite group operating over a \klgz~$\gA$ with $\gk\subseteq \gA$. 
\Propeq
\begin{enumerate}
\item $(\gk,\gA,G)$ is a \aG (in particular,~$\gk = \Ae G$).
\item $\gk=\Ae G$, and there exist $x_1$, \dots, $x_r$, $y_1$, \dots, $y_r$ in~$\gA$ such that we have for every $\sigma\in G$
$$\preskip-.2em \postskip.4em \ndsp \sum_{i=1}^r x_i \sigma(y_i)=
\formule{1\;\;\mathrm{if}\;\sigma=\Id\\[.2em]
0\;\;\mathrm{otherwise}.}
$$
\item $\gk=\Ae G$,~$\gA$ is finite over~$\gk$, and for every finite \sgr $(a_j)_{j\in J}$ of~$\gA$ as a \kmoz, there exists a family $(b_j)_{j \in J}$ in~$\gA$ such that we have for all $\sigma$, $\tau\in G$

\snic{\sum_{j \in J} \tau(a_j) \sigma(b_j)=
\formule{1\;\;\mathrm{if}\;\sigma=\tau\\[.2em]
0\;\;\mathrm{otherwise}.}
}
\item $\gk=\Ae G$, and $\psi_G : \env{\gk}{\gA} \to \prod_{\sigma\in G} \gA$ is an \iso of \Algsz.
\item $\gA$ is \stfe over $\gk$, and $G$ is a basis of~$\LIN_\gk(\gA,\gA)$.
%
\end{enumerate}
\end{theorem}

\vspace{.1em}
\begin{proof}
We have already seen \emph{1} $\Leftrightarrow$ \emph{2} and \emph{1} $\Rightarrow$ \emph{4} (Lemma~\ref{lemArtin}).
\\
The implication \emph{3} $\Rightarrow$ \emph{2} is clear.

\emph{2} $\Rightarrow$ \emph{3.}
We express $x_i$ in terms of $a_j$: $x_i = \sum_j u_{ij} a_j$ with $u_{ij} \in \gk$. Then,

\snic{\som_j \sigma\big(\som_i u_{ij} y_i\big) a_j =
\som_{j,i} u_{ij} \sigma(y_i) a_j =
\som_i \sigma(y_i) x_i = \delta_{\Id, \sigma},}

hence the result by taking $b_j = \sum_i u_{ij} y_i$.

 \emph{2} $\Rightarrow$ \emph{5.}
Let us first note that if $\varphi \in \LIN_\gk(\gA,\gA)$ is written as $\varphi = \sum_\sigma a_\sigma \sigma$, then by evaluating at $y_i$, by multiplying by $\tau(x_i)$ and by summing over the $i$'s, we get
$$
\preskip-.4em \postskip.3em \ndsp
\som_i \varphi(y_i) \tau(x_i) =
\som_{i,\sigma} a_\sigma \sigma(y_i) \tau(x_i) = a_\tau. 
$$
This shows on the one hand that $G$ is~$\gA$-free. On the other hand, this leads to believe that every $\varphi \in \LIN_\gk(\gA,\gA)$ is written as $\varphi = \sum_\sigma a_\sigma \sigma$ with $a_\sigma = \sum_i \varphi(y_i) \sigma(x_i)$. Let us verify this by evaluating $\varphi' := \sum_\sigma a_\sigma\sigma$ at $x \in \gA$,

\snac{\varphi'(x) = \som_{i,\sigma} \varphi(y_i) \sigma(x_i) \sigma(x) =
\som_i \Tr_G(x_ix) \varphi(y_i) =
\varphi(\sum_i \Tr_G(x_ix) y_i) = \varphi(x).}

\sni \emph{5} $\Rightarrow$ \emph{2.}
Since~$\gk \subseteq \gA$, we have an inclusion $\Asta \hookrightarrow \LIN_\gk(\gA,\gA)$.  Let us first show that $\Ae G \subseteq \gk$ (we will then have the \egtz). Each $\sigma \in G$ is $\Ae G$-\lin so, since $G$ generates $\LIN_\gk(\gA,\gA)$ as an \Amoz, each \eltz~$\varphi$ of~$\LIN_\gk(\gA,\gA)$ is $\Ae G$-\linz. In particular, each $\alpha \in \Asta$ is $\Ae G$-\linz.
Let $\big((x_i), (\alpha_i)\big)$ be a \syc of the \kmoz~$\gA$. As~$\gA$ is a faithful \kmoz, by Proposition~\ref{propAnnul}, there exists a family $(z_i)$ in~$\gA$ such that $1 = \sum_i
\alpha_i(z_i)$. Then, if $x \in \Ae G$, $x = \sum_i \alpha_i(z_i) x = \sum_i \alpha_i(z_ix)$ belonging to~$\gk$.

Let us then show that for each $\alpha \in \Asta$, there exists a unique~$a \in \gA$ such that $\alpha = \sum_{\sigma \in G} \sigma(a)\sigma$, \cad such that $\alpha$ is the $\gk$-\lin form~$x \mapsto \Tr_G(ax)$. Since $G$ is an~$\gA$-basis of $\LIN_\gk(\gA,\gA)$, we have $\alpha = \sum_\sigma a_\sigma \sigma$ with $a_\sigma \in \gA$. Let $a = a_\Id$. By writing, for $\tau \in G$, $\tau \circ \alpha = \alpha$, we obtain $\tau(a_\sigma) = a_{\sigma\tau}$, in particular $a_\tau = \tau(a)$, hence the desired \egt $\alpha = \sum_{\sigma \in G} \sigma(a)\sigma$. In passing, we have just proven that the \kli 

\snic{\gA \to
\Asta$, $a \mapsto \Tr_G(a\bullet)}

is an \iso of \kmosz. 
We can therefore define a \sys $(y_i)$ by the \egts $\alpha_i = \Tr_G(y_i\bullet)$. Then, for $x \in \gA$ we obtain

\snic{
x = \som_i \alpha_i(x) x_i =
\som_{i,\sigma} \sigma(y_ix) x_i =
\som_{\sigma} \big(\sum_i x_i\sigma(y_i)\big)\sigma(x),
}

\cad $\Id = \sum_{\sigma} \big(\sum_i x_i\sigma(y_i)\big)\sigma$.
But as $G$ is~$\gA$-free, the expression of $\Id \in G$ is reduced to $\Id$, so $\sum_i x_i\sigma(y_i) = 1$ if $\sigma = \Id$, $0$ otherwise.

NB: Since $\som_i x_iy_i = 1$, we have the \egts

\snic {
\Tr(x) = \sum_i\alpha_i(x_ix) =
\sum_{i,\sigma} \sigma(x_iy_i) \sigma(x) =
\sum_\sigma \sum_i \sigma(x_iy_i) \sigma(x) = \Tr_G(x).
}

\mni \emph{4} $\Rightarrow$ \emph{2.}
Let $z=\sum_i x_i \otimes y_i$ be the \elt of $\env\gk\gA$ defined by: $\psi_G(z)$ is the \elt of $\prod_{\sigma \in G} \gA$  every component of which 
is null, except that of index~$\Id$ which is equal to~$1$. This means precisely that $\sum_i x_i\sigma(y_i) = 1$ if $\sigma = \Id$ and $0$ otherwise.
\end{proof}
%

The case of free \aGs is described in the following corollary, which is an \imde consequence of the previous more \gnl results.

\begin{corollary}\label{corAGAfreeCar}\emph{(\Carns of  free \aGsz)}\\
Let $G$ be a finite group operating on a \klgz~$\gA$ with $\gk\subseteq \gA$. \\
Assume that $\gA$ is free over $\gk$, of rank $n=\abs G$, with $\ux=(\xn)$ as its basis.
\Propeq
\begin{enumerate}
\item $(\gk,\gA,G)$ is a \aG (in particular,~$\gk = \Ae G$).
\item The matrix $M_\ux=\big(\sigma_i(x_j)\big)_{i,j\in\lrbn}$ is \iv (we have indexed the group~$G$ by $\lrbn$). 
\item The form $\Tr_G$ is dualizing.
\item $\gk=\Ae G$, and there exist $y_1$, \dots, $y_n$ in~$\gA$ such that we have for every $\sigma\in G$
$$\preskip-.4em \postskip.2em\ndsp 
\sum_{i=1}^n x_i \sigma(y_i)=
\formule{1\;\;\hbox{if}\;\sigma=\Id\\[.2em]
0\;\;\hbox{otherwise}.} 
$$
\item The group $G$ is an~$\gA$-basis of~$\LIN_\gk(\gA,\gA)$.
\item $\gk=\Ae G$, and $\psi_G : \env{\gk}{\gA} \to \prod_{\sigma\in G} \gA$ is an \iso of \Algsz.
\end{enumerate}
In this case we have the following results. 
\begin{enumerate} \setcounter{enumi}{6}
\item 
In items {4} and {3},
\begin{itemize}
\item we obtain the $y_i$'s as the solution of $M_\ux . \tra {[\,y_1\;\cdots\; y_n\,]} = \tra {[\,1\;0\cdots\; 0\,]}$, where $M_\ux$ is defined as in item 2, with $\sigma_1 = Id$,
\item  $(y_1\sta,\ldots,y_n\sta)$ is the dual basis of $(x_1,\dots,x_n)$.
\end{itemize}
\item Item {6} can be specified as follows.\\
For $\sigma \in G$, we let $\varepsilon_\sigma = \sum_i \sigma(x_i) \otimes
y_i \in \env\gk\gA $. Then, $(\varepsilon_\sigma)_{\sigma \in G}$ is an~$\gA$-basis for the left-structure 
of $\env\gk\gA$. In addition, for $a$, $b \in \gA$, we have
$$\preskip.2em \postskip.2em\ndsp 
b \otimes a = \sum_\sigma  b\sigma(a) \varepsilon_\sigma 
$$
and the image of this basis $(\varepsilon_\sigma)_\sigma$ under $\psi_G : \env\gk\gA \to \prod_{\tau \in G}\gA$ is the canonical $\gA$-basis $(e_\sigma)_{\sigma\in G}$ of $\prod_{\tau \in G}\gA$.
\end{enumerate}
Finally, we underline the following items, in which we do not suppose that~$\gA$ is free over $\Ae G$. 
\begin{itemize}
\item When $\gA$ is a \cdi (historical background of Artin's \thoz), if a group $G$ operates faithfully over $\gA$, the \alg $(\Ae G,\gA,G)$ is always Galoisian, $\Ae G$ is a \cdi and $\gA$ is free of rank~$n$ over $\Ae G$.
%
%
\item When $\gA$ is a \alrdz, the \alg $(\Ae G,\gA,G)$ is Galoisian \ssi $G$ operates faithfully over the residual field~$\gA/\Rad\gA$.  
In this case, $\Ae G$ is a \alrd and $\gA$ is free of rank~$n$ over $\Ae G$.    
\end{itemize} 

\end{corollary}

Naturally, we strongly encourage the reader to give a more direct and shorter \dem of the previous corollary. It is \egmt possible to deduce the \gnl results of the particular results stated in the case where $\gA$ is a \alrdz, which could themselves be deduced from the \cdis case.


\begin{theorem} 
\label{corDAexplicite2} \emph{(The Galois correspondence for a \aGz)}
Let $(\gk,\gA,G)$ be a nontrivial \aGz, and $H$ be a finite subgroup of $G$.
\begin{enumerate}
\item 
The triple $(\Ae H,\gA,H)$ is a \aGz, $\Ae H$ is \ste over~$\gk$, of constant rank $[\Ae H:\gk]=\idg{G:H}$. 

\item If $ H'\supseteq H$ is a finite subgroup of $G$, $\Ae H$ is \stfe over~$\Ae {H'}$, of constant rank $[\Ae H:\Ae {H'}]=\idg{H':H}$.
\item  
We have $H= \Stp(\Ae H)$.
\item  
The map $\Fix_\gA$ restricted to the finite subgroups of $G$ is injective.
\item 
If $H$ is normal in $G$, $(\gk,\Ae H,G/H)$ is a \aGz. 
\end{enumerate}
\end{theorem}

\vspace{.1em}
\begin{proof}
\emph{1.}  
 Since $H$ is a separating group of \autos of~$\gA$, $(\Ae H,\gA,H)$ is a \aGz. So~$\gA$ is a \stfe $\Ae H$-\alg of constant rank $\#H$. 
Therefore $\Ae H$ is \stfe over~$\gk$, of constant rank equal to~$\idg{G:H}$ (\thref{propTraptf}). In addition, it is \ste by Fact~\ref{propTrEta}.

\emph{2.} We apply \thref{propTraptf}. 

\emph{3.} The inclusion $H\subseteq \Stp(\Ae H)$ is obvious. Let $\sigma\in \Stp(\Ae H)$ and~$H'$ be the subgroup generated by $H$ and $\sigma$. We have~$\idg{H':H}=[\Ae H:\Ae {H'}]$, but $\Ae H=\Ae {H'}$, therefore~$H'=H$ and~$\sigma\in H$.

\emph{4.} 
Results from \emph{3.}

\emph{5.} First of all, for $\sigma\in G$, we have $\sigma(\Ae H) = \Ae H$. If we let $\ov\sigma$ be the restriction of $\sigma$ to $\Ae H$, we obtain a morphism of groups $G \to \Aut_\gk(\Ae H)$, $\sigma \mapsto \ov\sigma$, whose kernel is $H$ by item~\emph {3.} The quotient group $G/H$ is therefore realized as a subgroup of $\Aut_\gk(\Ae H)$.
\\
Let  $x \in \gA$ satisfy $\Tr_H(x) = 1$,   $(a_1, \ldots, a_r)$ be a \sgr of~$\gA$ as a \kmoz, and  $b_1$, \ldots, $b_r$ be some \elts such that 
for all $\sigma$, $\tau\in G$ we have 
${\sum_{i=1}^r \tau(a_i) \sigma(b_i)= 
\formule{1\;\;\mathrm{if}\;\sigma=\tau\\
0\;\;\mathrm{otherwise}.}
}$.
We then define, for $i \in \lrb{1..r}$, the \elts of $\Ae H$,
${
a'_i = \Tr_H(xa_i), \hbox{ and } b'_i = \Tr_H(b_i).
}$
\\
We easily verify that for $\sigma\in G$ we have

\snic{\sum_{i=1}^r a'_i \sigma(b'_i)= 
\formule{1\;\;\mathrm{if}\;\sigma \in H\\
0\;\;\mathrm{otherwise}.}
}

Thus, when applying item~\emph{2} of \thref{thAGACar}, $(\gk, \Ae H, \G/H)$ is a \aGz. 
\end{proof}

\Thref{corDAexplicite2} above establishes the Galois correspondence between finite subgroups of $G$ on the one hand and \gui{certain} \stes \kslgs of~$\gA$ on the other.
An exact bijective correspondence will be established in the following subsection when~$\gA$ is connected.

However, beforehand we give a few \suls results.

\begin{proposition}\label{propAdiagAH}
Let $(\gk,\gA,G)$ be a \aG  and $H$ be a finite subgroup of $G$.
\begin{enumerate}
\item $\gA$ diagonalizes $\Ae H$.
\item For $b\in\Ae H$, the \polcar of $b$ (over~$\gk$, in $\Ae H$) is given by
${
\rC{\Ae H\!/\gk}(b)(T) = \prod_{\sigma \in G/H}\big(T - \sigma(b)\big)
}$
(the subscript $\sigma \in G/H$ means that we take exactly one $\sigma$ from each left coset of~$H$, and we note that $\sigma(b)$ does not depend on the chosen representative $\sigma$).
\end{enumerate} 
\end{proposition}
%
\begin{proof}
Recall that~$\gA$ diagonalizes itself, as the \iso $\psi_G:\env\gk\gA \to \prod_{\sigma\in G}\gA$ shows. We consider this product as the \alg of functions  $\cF(G,\gA)$. It is provided with a natural action of $G$ on the left-hand side as follows

\snic{\sigma\in G,\, w\in \cF(G,\gA): \sigma \cdot w\in \cF(G,\gA)
\hbox{ defined by } \tau \mapsto w(\tau\sigma).  
}

Similarly $G$ acts on the left-hand side over the \Alg  $\env\gk\gA = \gA\te_\gk\gA$ via
$\Id \te G$. We then verify that $\psi_G$ is a $G$-morphism, \cad that for $\tau \in G$, the following diagram commutes.
$$
\preskip.0em \postskip.4em 
\hspace{.03cm}\xymatrix {
\gA\te_\gk\gA \ar[d]^{\Id \te\tau} \ar[r]^(0.35){\psi_G} &
       \cF(G,\gA)=\prod_{\sigma\in G}\gA \ar[d]^{w \mapsto \tau\cdot w}\\
\gA\te_\gk\gA \ar[r]^(0.35){\psi_G} &\cF(G,\gA)=\prod_{\sigma\in G}\gA\\
} 
$$

\emph{1.}
Consider the commutative diagram
$$\preskip.4em \postskip.4em 
\hspace{.05cm}\xymatrix {
\gA\te_\gk\Ae H \ar[d] \ar[r]^(0.33){\varphi_H} &
       \cF(G/H,\gA)=\prod_{\sigma\in G/H}\gA \ar[d]\\
\gA\te_\gk\gA \ar[r]^(0.36){\psi_G}_(0.36){\sim} &\cF(G,\gA)=\prod_{\sigma\in G}\gA\\
}
$$

On the right-hand side, the vertical arrow is injective, and it identifies $\cF(G/H,\gA)$ with the subset~$\cF(G,\gA)^H$ of~$\cF(G,\gA)$ (constant functions over the left cosets of $H$ in $G$).  
\\
On the left, the vertical arrow (corresponding to the injection $\Ae H \hookrightarrow \gA$) is also an injection because $\Ae H$ is a direct summand in~$\gA$ viewed as an $\Ae H$-module. 
Finally, $\varphi_H$ is defined by $a\te b \mapsto \big(a\sigma(b)\big)_{\sigma \in G/H}$. 
\\
Then, $\varphi_H$ is an \iso of \Algsz. Indeed, $\varphi_H$ is injective, and for the surjectivity, it suffices to see that $(\gA\te_\gk\gA)^{\Id\te H} = \gA\te_\gk\Ae H$. 
This is given by \thref{corAGextsca} for the \aG $(\Ae H,\gA,H)$ and the \eds $\Ae H\hookrightarrow \gA$. 

\emph{2.} This results from item \emph{1} and from the following lemma.  
\end{proof}
%

\begin{lemma}\label{lemPolCarDiag}
Let~$\gA$ and~$\gB$ be two \klgs where $\gB$ is \stfe of constant rank $n$. Assume that $\gA$ diagonalizes~$\gB$ by means of an \iso  

\snic{\psi:
\gA\te_\gk \gB \lora \Ae n}

given by \gui{\coosz} denoted by $\psi_i : \gB \to \gA$.
\\
Then, for $b \in \gB$, we have an \egt 

\snic{\rC{\gB/\gk}(b)(T)=\prod_{i=1}^n \big(T - \psi_i(b)\big),}

if we transform the left-hand side (which is an \elt of $\kT$) into an \elt of~$\gA[T]$ via~$\gk \to \gA$.
\end{lemma}
%
\begin{proof}
Immediate by the computation of the \polcar of an \elt in a diagonal \algz.
\end{proof}
%

\subsect{The Galois correspondence 
when $\gA$ is connected}{The Galois correspondence in the connected case}

The reader is invited to revisit Lemma~\ref{lemSspbdiag}.

\begin{theorem}
\label{thCorGalGen}
If $(\gk,\gA,G)$ is a nontrivial \aG and if~$\gA$ is \emph{connected}, the Galois correspondence establishes a decreasing bijection between
\begin{itemize}
\item on the one hand, the set of detachable subgroups of $G$,
\item and on the other hand, the set of \kslgs of~$\gA$ which are \spbsz.
\end{itemize}
The latter set is \egmt that of the sub\algs of~$\gA$ which are \stes over~$\gk$.
\end{theorem}

\begin {proof}
Let~$\gk\subseteq\gA'\subseteq\gA$ with~$\gA'$ \spbz. By letting $H = \Stp(\gA')$, we must show that~$\gA' = \Ae H$. We of course have~$\gA'\subseteq \Ae H$. \\
Consider the product \Algz~$\gC = \prod_{\sigma \in G}\gA \simeq \gA^n$ with $n = \#G$. 
\\
Let $p_\sigma : \gC \to \gA$ be the \prn defined by $p_\sigma\big((a_\tau)_\tau\big) = a_\sigma$.  
Recall the \iso of \Algs $\psi_G : \gA\otimes_\gk\gA \to \gC,\;a \otimes b \mapsto \big(a\sigma(b)\big)_{\sigma \in G}$.
\\
Since~$\gA$ is a \ptf \kmoz, the canonical morphism~$\gA\otimes_\gk\gA' \to \gA\otimes_\gk\gA$ is injective. By composing it with $\psi_G$, we obtain an injective morphism of \Algs $\gA\otimes_\gk\gA' \to \gC$. In the 
above notation, we will identify~$\gA\otimes_\gk\gA'$ with its image~$\gB$ in~$\gC \simeq \gA^n$.
\\
Since~$\gA'$ is a \spb \klgz,~$\gB$ is a \spb \Algz. We can therefore apply Lemma~\ref{lemSspbdiag}. If we denote by $\pi_\sigma$ the restriction of~$p_\sigma$ to~$\gB$, we must identify the \eqv relation over $G$ defined by $\pi_\sigma = \pi_{\sigma'}$. 
For $a' \in \gA'$, $1 \otimes a'$ corresponds by $\psi_G$ to $\big(\tau(a')\big)_\tau$, so $\pi_\sigma(1 \otimes a') = \sigma(a')$. 
Consequently, $\pi_\sigma = \pi_{\sigma'}$ \ssi $\sigma $ and $\sigma'$ coincide over~$\gA'$ or, by \dfn of~$H$, \ssi $\sigma^{-1}\sigma' \in H$, \cad $\sigma H = \sigma' H$. 
We deduce that the \eqvc classes are the left cosets of $H$ in $G$. 
With the notations of Lemma~\ref{lemSspbdiag}, we therefore have~$\gB = \bigoplus_J \gA e_J$, where~$J$ describes $G/H$. 
By using the~$\gA$-basis $(e_J)_J$ of~$\gB$, we then see that~$\gB = \gC^H$.
\\
It remains to \gui{return} 
to~$\gA$. Via the inverse image under $\psi_G$, we have 
$$
(\gA\otimes_\gk\gA)^{\Id\otimes H} = \gA\otimes_\gk\gA'\,.
$$
In particular,~$\gA\otimes_\gk\Ae H \subseteq \gA\otimes_\gk\gA'$.
By applying $\Tr_G \otimes \Id_\gA$ to this inclusion and by using the fact that $\Tr_G : \gA \to \gk$ is surjective, we obtain the inclusion

\snic{\gk\otimes_\gk\Ae H \subseteq \gk\otimes_\gk\gA'$, \cad~$\;\Ae H \subseteq
\gA'\,.}

Thus~$\Ae H = \gA'$, as required.

Finally, since the \klgs~$\Ae H$ are \stes and the \stes \algs are \spbsz, it is clear that the \spbs \kslgs of~$\gA$ coincide with the \stes \kslgsz.
\end{proof}

\rem The theory of \aGs does not really require the use of \spb \algsz, even for the previous \tho that we can state with only \stes sub\algs of~$\gA$. For a \dem of the \tho without using \spb \algsz, see Exercises~\ref{lemPaquesFerrero} and~\ref{exoCorGalste}. 
Nevertheless the theory of \spb \algsz, noteworthy in itself, sheds an interesting light on the \aGsz. 
 \eoe

\subsec{Quotients of \aGsz}

\vspace{.1em}

\begin{proposition}\label{prop1GalQuo}
\emph{(Quotient of a \aG by an invariant \idz)} Let $(\gk,\gC,G)$ be a \aGz, $\fc$ be a $G$-invariant \id of~$\gC$ and~$\fa = \fc\cap\gk$. 
\begin{enumerate}
\item The triple $(\gk\sur\fa, \gC\sur\fc, G)$ is a \aGz. 
\item This \aG is naturally \isoc to that obtained from $(\gk,\gC,G)$ by means of the \edsz~$\gk\to\gk\sur\fa$. 

%
%
\end{enumerate}
\end{proposition}

\begin{proof} \emph{1.} The group $G$ operates on~$\gC\sur\fc$ because $\fc$ is (globally) invariant.
Let us show that the natural injective \homoz~$\gk\sur\fa \to (\gC\sur\fc)^G$ is surjective. If $x \in \gC$ is $G$-invariant modulo $\fc$, we must find an \elt of~$\gk$ equal to $x$ modulo $\fc$. Consider $x_0 \in \gC$ satisfying $\Tr_G(x_0) = 1$; then $\Tr_G(xx_0)$ satisfies:

\snic{x = \som_{\sigma \in G} x\sigma(x_0) \equiv 
\som_{\sigma \in G}  \sigma(x)\sigma(x_0) = \Tr_G(xx_0) \mod \fc\,.}

Thus $(\gC\sur\fc)^G =\gk\sur\fa$. Finally, it is clear that $G$ operates in a separating way over~$\gC\sur\fc$.

\emph{2.} The \edsz~$\gk\to\gk\sur\fa$ gives $(\gk\sur\fa\!,\gC\sur{\fa\gC}\!,G)$ (\aGz), with $\fa\gC\subseteq\fc$. We must verify that $\fc=\fa\gC$. \\
The \prn $\pi:\gC\sur{\fa\gC}\to\gC\sur\fc$ is a~$\gk\sur\fa$-\lin surjective map between \pro modules, so~$\gC\sur{\fa\gC}\simeq \gC\sur\fc\oplus\Ker\pi$. As the two modules have the same constant rank $\# G$, the \polmu of $\Ker\pi$ is equal to~$1$, therefore $\Ker\pi=0$ (\thref{th ptf sfio}).
\end{proof}

In the \dfn that follows, we do not need to suppose that $(\gk,\gC,G)$ is a \aGz.

\begin{definition}
\label{defQuoDeGal}\label{defIdmGalAdu}\index{ideal!Galoisian ---}
\index{Galois quotient!of an \algz} \hspace*{-1.3em}
Let $G$ be a finite group that operates on \hbox{a \klg $\gC$}.
\begin{enumerate}
\item An \idm  of~$\gC$ is said to be \ixc{Galoisian}{idempotent in an \alg provided with a finite group of \autosz} 
if its orbit under $G$ is a \sfio (this requires that this orbit is a finite set, or, equivalently, that the subgroup~$\St_G(e)$ is detachable).
\item An \id of~$\gC$ is said to be \ixc{Galoisian}{ideal} when it is generated by the \cop \idm of a Galoisian \idm $e$.
\item In this case, the group $\St_G(e)$ operates on the \algz~$\gC[1/e]\simeq\aqo{\gC}{1-e}$, and~$\big(\gk,\gC[1/e],\St_G(e)\big)$ is called a \emph{Galois quotient} of $(\gk,\gC,G)$.
\end{enumerate} 
\end{definition}

\smallskip 
\begin{fact}\label{factQuoDeGal}
With the hypotheses of \Dfnz~\ref{defQuoDeGal}, if $\so{e_1,\ldots,e_r}$ is the orbit of $e$, the natural \kliz~$\gC\to \prod_{i=1}^r \gC[1/ {e_i}]$ is an \iso of \klgsz. Moreover, the $\St_G(e_i)$'s are pairwise conjugated by \elts of $G$ that permute the \klgsz~$\gC[1/ {e_i}]$ (they are therefore pairwise \isocz). In particular~$\gC\simeq\gC[1/e]^r$.
\end{fact}
\smallskip 
\begin{theorem}\label{corAGQuo} \emph{(Galois quotients of \aGsz)}\\
Every Galois quotient of a \aG is a \aGz.
\end{theorem}
\begin{proof} See \thref{thADG1Idm} (Galois quotients of \apGsz).
\end{proof}

\pagebreak	
\Exercices

\begin{exercise}
\label{exochapAlgStricFiLecteur}
{\rm  We recommend that the proofs which are not given, or are sketched, or left to the reader, etc, be done. But in particular, we will cover the following cases. \perso{\`a compl\'eter \`a la fin}
\begin{itemize}\itemsep0pt
\item  Prove  \Tho \rref{factSDIRKlg}.
\item  Prove  Fact~\rref{factEdsAlg}.
\item  Prove  the \plgref{plcc.aGs} for \aGsz.
\item  Verify Fact~\rref{factdefAGAL}. 
\end{itemize}
}
\end{exercise}

\vspace{-1em}
\begin{exercise}
\label{exothEtalePrimitif}
{\rm  Give a detailed \dem of item \emph{2} of \thref{thEtalePrimitif}.
 
}
\end{exercise}

\vspace{-1em}
\begin{exercise}
\label{exoKnEltPrimitif}
{\rm
Consider the product \Algz~$\gB = \Ae n$.
 
\emph{1.} Under what condition does some $x \in \gB$ satisfy~$\gB = \gA[x]$? \\
In this case, prove that $(1, x, \ldots, x^{n-1})$ is an~$\gA$-basis of~$\gB$.
 
\emph{2.} If~$\gA$ is a \cdiz, under what condition does~$\gB$ admit a primitive \eltz?
}
\end{exercise}

\vspace{-1em}
\begin{exercise}
\label{exoIdmChangBase}
{\rm Let~$\gK$ be a nontrivial \cdiz,
$\gB$ be a reduced \stf \Klg and $v$ be an \idtrz.\\
Consider the \Llgz~$\gB(v)\eqdefi \gK(v)\otimes_\gK\gB$.
Prove the following results.

\emph{1.}~$\gB(v)$ is \stfe over~$\gK(v)$.
 
\emph{2.}   If~$\gB$ is \'etale over~$\gK$,~$\gB(v)$ is \'etale over~$\gK(v)$.
 
\emph{3.} Every \idm of~$\gB(v)$ is in fact in~$\gB$.
\entrenous{exo \ref{exoIdmChangBase}. Cela rappelle des choses in le cours:
\thrf{thIdmEtale}. Mais la preuve directe
pour le cas des \cdis peut \^etre un bon entra\^{\i}nement;
peut-\^{e}tre il faudrait demander
de \gnr lorsque~$\gK$ is un \ri arbitraire: la situation
est contr\^ol\'ee par le fait que~$\gK$ is \icl in~$\gL$,
et il semble que l'on  a besoin of~$\gB$ \'etale sur~$\gK$.
}
}
\end{exercise}

\vspace{-1em}
\begin{exercise}
\label{exo1SepFact}
{\rm  If~$\gK$ is a \sply factorial \cdiz, so is~$\gK(v)$, where $v$ is an \idtrz.
\perso{L'exo \ref{exo1SepFact} semble redoutable}
\\
NB: we do not assume that~$\gK$ is finite or infinite.
}
\end{exercise}

\vspace{-1em}
\begin{exercise}\label{exoBiquadratique}  {(The \ris of integers of the extension
$\QQ(\sqrt a) \subset \QQ(\sqrt a,\sqrt 2)$)}
\\
{\rm
Let~$\gK \subseteq \gL$ be two number fields and~$\gA \subseteq \gB$ be their \ris of integers; here we give an elementary example where~$\gB$ is not a free~$\gA$-module.

 \emph {1.}
Let $d \in \ZZ$ be squarefree. Determine the \ri of integers of $\QQ(\sqrt d)$.

Let $a \in \ZZ$ squarefree with $a \equiv 3 \bmod 4$. 
Let~$\gK = \QQ(\sqrt a)$,~$\gL = \gK(\sqrt 2)$, \hbox{and~$\beta = \sqrt 2 {1 + \sqrt a \over 2}$}. 
We define $\sigma \in \Aut(\gL/\gK)$ and $\tau \in \Aut\big(\gL/\QQ(\sqrt 2)\big)$, by\\
\centerline{$\sigma  (\sqrt 2) = -\sqrt 2\;$  and $\;\tau ( \sqrt a) = -\sqrt a$.}

 \emph {2.}
Verify that $\beta \in \gB$ and compute $(\sigma\tau)(\beta)$.

 \emph {3.}
We want to show that $(1, \sqrt 2, \sqrt a, \beta)$ (which is a $\QQ$-basis of $\gL$) is a $\ZZ$-basis of~$\gB$.  Let $z = r + s\sqrt 2 + t\sqrt a + u\beta \in \gB$ with $r$, $s$, $t$, $u \in \QQ$. \\
Considering $(\sigma\tau)(z)$, show that $u \in \ZZ$ then that $r$, $s$, $t \in \ZZ$.

 \emph {4.}
Express~$\gB$ as a \ptf \Amoz. 
Verify that it is \isoc to its dual.
}
\end{exercise}

\vspace{-1em}
\begin{exercise}\label{exoTensorielDiscriminant}
 {(Discriminant of the tensor product)}\\
{\rm
Let $\gA$, $\gA'$ be two free \klgs of ranks $n$, $n'$, $(\ux) = (x_i)$ be a family of $n$ \elts of~$\gA$, $(\ux') = (x'_j)$ be a family of $n'$ \elts of~$\gA'$. Let~$\gB = \gA\otimes_\gk\gA'$ and $(\ux\otimes\ux')$ be the family $(x_i \otimes x'_j)$ of $nn'$ \elts of~$\gB$. Prove  the \egt

\snic {
\Disc_{\gB\sur\gk}(\ux\otimes\ux') =
\Disc_{\gA\sur\gk}(\ux)^{n'} \, \Disc_{\gA'\sur\gk}(\ux')^n.
}

}
\end{exercise}

\vspace{-1em}
\begin{exercise}\label{CyclicNormalBasis}
 {(Normal basis of a cyclic extension)} 
\\
{\rm  
Let~$\gL$ be a \cdiz, $\sigma \in \Aut(\gL)$ of order $n$ and~$\gK = \gL^\sigma$ be the field of invariants under $\sigma$. Prove  that there exists an $x \in \gL$ such that $\big(x, \sigma(x), \cdots, \sigma^{n-1}(x)\big)$ is a~$\gK$-basis of~$\gL$; we then speak of a \emph{normal basis} of~$\gL\sur\gK$ (defined by $x$).
}

\end{exercise}

\vspace{-1em}
\begin{exercise}\label{exoHomographieOrdre3} \label{NOTAAn} {(Homography of order $3$ and \uvle \eqn with Galois group~$\rA_3$)}  
{\rm We denote by $\rA_n$ the subgroup of even permutations of~$\Sn$.
 Let~$\gL=\gk(t)$ where~$\gk$ is a \cdi and $t$ is \idtrz. 
\begin{enumerate}
\item
Check that \smashbot{$A = \crmatrix {0 &-1\cr 1 & 1\cr}$} is of order $3$ in $\PGL_2(\gk)$ and explain the origin of this matrix.
\end{enumerate}
We denote by $\sigma\in \Aut_\gk\big(\gk(t)\big)$ the \auto of order $3$ associated with $A$ (see \Pbmz~\ref{exoLuroth1}, we~have~$\sigma(f)=f({-1\over t+1})$), and $G = \gen
{\sigma}$.
\begin{enumerate}\setcounter{enumi}{1}
\item Compute $g=\Tr_G(t)$ and show that~$\gk(t)^G = \gk(g)$.

\item
Let $a$ be an \idtr over~$\gk$ and $f_a(T) = T^3 - aT^2 - (a+3)T - 1 \in \gk(a)[T]$.  Prove  that $f_a$ is \irdz, with Galois group $\rA_3$. 

\item
Prove  that the \pol $f_a(X)$ is a \gui{generic \pol with Galois group $\rA_3$} in the following sense: if~$\gL/\gK$ is a Galois extension with Galois group $\rA_3$ ($\gL$ being a \cdiz), there exists a primitive \elt of~$\gL/\gK$ whose \polmin is $f_\alpha(X)$ for some value of $\alpha \in \gK$.
\end{enumerate}
}
\end{exercise}

\vspace{-1em}
\begin{exercise}
\label{exoGroupAlgebra} {(Algebra of a finite commutative group)}\\
{\rm  
Let~$\gk$ be a commutative \riz, $G$ be a commutative group of order $n$ and~$\gA = \gk[G]$ be the \emph{\alg of the group $G$}, \cad $\gA$ admits $G$ as a~$\gk$-basis and the product in~$\gA$ of two \elts of $G$ is their product in $G$.\footnote{The \dfn works also for \emph{the \alg $\gk[M]$ of a \mo $M$}.}\index{algebra@\algz!of a monoid} 
\begin{enumerate}
\item Determine $\Ann(\rJ\iAk)$, its image under $\mu\iAk$ and the trace  form over~$\gA$.
\item Prove  that \propeq
\begin{itemize}
\item
$n$ is \iv in~$\gk$.

\item
$\gA$ is \stez.

\item
$\gA$ is \spbz.
\end{itemize}

\item Prove  that~$\gk[G]$ is a Frobenius \algz.
\end{enumerate}

}
\end{exercise}

\vspace{-1em}
\begin{exercise}\label{exo1Frobenius}
{(A finite monogenic \alg is a Frobenius \algz)}\\
{\rm  Let $f = X^n + a_{n-1} X^{n-1} + \cdots + a_0 \in \gk[X]$ and~$\gA = \aqo{\kX}{f} = \gk[x]$.  Consider the \lin form $\lambda : \gA \to \gk$ defined by $x^{n-1} \mapsto 1$ and $x^i \mapsto 0$ for $i < n-1$. 
We will show that $\lambda$ is dualizing and that $\Tr\iAk = f'(x)\centerdot\lambda$. 
\\
To that effect, we append an \idtr $Y$. 
The \sys $(1, x, \ldots,x^{n-1})$ is a basis of~$\gA[Y]\sur{\gk[Y]}$. Let $\wi\lambda : \gA[Y] \to \gk[Y]$ be the extension of $\lambda$ and define the $\gk[Y]$-\lin map  $\varphi : \gA[Y] \to \gk[Y]$, by $\varphi(x^i)=Y^i$ for $i \in\lrb{0.. n-1}$.
\begin{enumerate}
\item 
Prove  that
$\qquad\forall g \in \gA[Y],\quad f(Y) \wi\lambda(g) = \varphi\big((Y - x)g\big)\qquad(*)$

\item 
We define the (triangular Horner) basis $(b_0, \ldots, b_{n-1})$ of $\gA\sur\gk$ by

\snic {\arraycolsep2pt
\begin{array}{ccc} 
b_0  &  = & x^{n-1} + a_{n-1} x^{n-2} + \cdots + a_2x + a_1,  \\[1mm] 
b_1   & =   & x^{n-2} + a_{n-1} x^{n-3} + \cdots + a_3x + a_2  \,,
\end{array}
 }

and so on: $b_i = x^{n-i-1} + \cdots + a_{i+1}$ and $b_{n-1} = 1$.
We have
$$
f'(Y) = \frac{f(Y) - f(x)}{Y - x}
= \frac{f(Y)}{Y - x} = b_{n-1}Y^{n-1} + \cdots + b_1Y + b_0.
$$
Applying \Egt $(*)$ to $g_i = x^if'(Y)$, show that $(b_0\centerdot\lambda, \ldots, b_{n-1}\centerdot\lambda)$
is the dual basis of $(1, x, \ldots, x^{n-1})$.
Conclude the result.

\item 
Prove  that $\Tr\iAk = f'(x)\centerdot\lambda$.
\end{enumerate}

}

\end{exercise}

\vspace{-1em}
\begin{exercise}\label{exoFrobeniusAlgExemples}
 {(Frobenius \algsz: \elr examples and counterexamples)}\\
{\rm  
Throughout the exercise,~$\gk$ is a commutative \riz.

\emph {1.}
Let $f_1$, \ldots, $f_n \in \gk[T]$ be \polusz. Prove  that the quotient \klg $\gk[X_1, \ldots, X_n]\sur{\gen {f_1(X_1), \ldots, f_n(X_n)}}$ 
is Frobenius and free of finite rank.

\emph {2.}
Let~$\gA = \gk[X,Y]\sur{\gen{X,Y}^2} = \gk[x,y]$. Describe $\Asta$ as a \pf \Amoz. Deduce that~$\gA$ is not a Frobenius \algz.

\emph {3.}
Consider the analogous question to the previous one with~$\gA = \gk[X,Y]\sur{\gen {X,Y}^n}$ for $n \ge 2$ and~$\gB = \gk[X,Y]\sur{\gen {X^2,XY^{n+1},Y^{n+2}}}$ \hbox{for $n \ge 0$}.

}
\end{exercise}

\vspace{-1em}
\begin{exercise}\label{exoAlgMonogeneJAnnJ} {(The \id $\rJ\iAk$ for a monogenic \klgz~$\gA$)}\\
{\rm  
Let~$\gA = \gk[x]$ be a monogenic \klg and $\env{\gk}{\gA} = \gA\te_\gk\gA$ be its enveloping \algz. 
Let $y = x\te 1$, $z = 1\te x$, such that $\env\gk\gA = \gk[y,z]$.
We know that $\rJ\iAk=\gen{y-z}$.
Suppose $f(x) = 0$ for some $f \in \gk[X]$ (not \ncrt \monz) and consider  the \smq \pol  
 $f^\Delta(Y,Z) = \big(f(Y) - f(Z)\big) /(Y-Z)$. It satisfies the \egt $f^\Delta(X,X) = f'(X)$.

\emph{1.} 
Let $\delta = f^\Delta(y,z)$. Prove  that $\delta \in \Ann(\rJ\iAk)$
and that $\delta^2 = f'(y)\delta = f'(z)\delta$.

\emph{2.} 
Suppose that $1 \in \gen {f, f'}$.\\
\emph{2a.} 
Prove  that~$\gA$ is \spbz: make the \spt \idm explicit. \\
\emph{2b.} Prove  that  $\rJ\iAk=\gen{f^\Delta(y,z)}$
and that $f^\Delta(y,z) = f'(y)\vep\iAk = f'(z)\vep\iAk$.

\rem $\gA$ is not \ncrt \stfez.

}

\end{exercise}

\vspace{-1em}
\pagebreak	

\begin{exercise}\label{exoJacInversibleAlgSeparable}
 {(Complete intersection, Jacobian, Bézoutian and \sptz)}\\
{\rm  
In this exercise, the number of \idtrs is equal to the number of \polsz. 
We define the \emph{Bézoutian} of $(f_1, \ldots, f_n)$ where each $f_i\in \gk[\uX] = \gk[\Xn]$ by
$$\preskip-.6em \postskip.2em 
\beta_{\uY,\uZ}(\uf) = \det \BZ_{\uY,\uZ}(\uf), 
$$
so that $\beta_{\uX,\uX}(\uf) = \J_\uX(\uf)$.
\index{Bézoutian!determinant of a \sypz}

We denote by $\gA = \gk[\xn]$ a \tf \klg and $\env{\gk}{\gA} = \gk[\uy,\uz]$ its enveloping \algz. Suppose that $f_i(\ux) = 0$ for all $i$.

\emph {1.}
In the case where $\J_\ux(f_1, \ldots, f_n)\in\Ati$, provide a direct \dem of the fact that~$\gA$ is a \spb \algz.

\emph {2.}
We define in $\env{\gk}{\gA}$
$$\preskip-.0em \postskip.2em 
\;\;\;\vep = \J_\uy(\uf)^{-1} \beta_{\uy,\uz}(\uf) = \beta_{\uy,\uz}(\uf) \J_\uz(\uf)^{-1}
. 
$$
 Verify that $\beta_{\uy, \uz}(\uf)$ and $\vep$ are \gtrs of $\Ann(\rJ\iAk)$
and that $\vep$ is the \idst of $\gA$.

\emph {3.}
Give examples.
}

\end{exercise}

\vspace{-1em}
\begin{exercise}\label{exoSepAlgDedekindLemma} 
{(Separation of morphisms over a \spb \algz)}\\
{\rm  
Let~$\gk$ be a commutative \ri and~$\gA$,~$\gB$ be two \klgs with~$\gA$ \spbz.
For any arbitrary function $f : \gA \to \gB$, we define $\Ann_\gB(f) = \Ann_\gB \gen {f(\gA)}$.

\emph {1.}
Prove  that to every morphism $\varphi \in \Hom_\gk(\gA,\gB)$ is attached 
a pair of finite families $(a_i)_{i \in I}$, $(b_i)_{i \in I}$, with $a_i \in \gA$, $b_i \in \gB$, satisfying the following \prtsz:
\begin{itemize}
\item $\sum_i b_i\varphi(a_i) = 1$
\item $\sum_i \varphi(a)b_i \te a_i = \sum_i b_i \te aa_i$ for every $a \in \gA$.
\end{itemize}

\emph {2.}
If the pair of families $(a'_j)_j$, $(b'_j)_j$ is attached to the morphism $\varphi' \in \Hom_\gk(\gA,\gB)$, show that
$$
\preskip-.4em \postskip.4em\ndsp 
\sum_i b_i\varphi'(a_i) = \sum_j b'_j\varphi(a'_j)
, 
$$
and that the latter \eltz, denoted by $e$, is an \idm of~$\gB$ having the following \prt of \gui{separation of morphisms}

\snic {
\Ann_\gB(\varphi - \varphi') = \gen {e}_\gB, \qquad
\gen {\Im(\varphi - \varphi')}_\gB  = \gen {1-e}_\gB
.}

\emph {3.}
Let $\varphi_1$, \ldots, $\varphi_n \in \Hom_\gk(\gA,\gB)$ and, for $i$, $j\in\lrbn$, $e_{ij} = e_{ji}$ be the \idm defined by $\Ann_\gB(\varphi_i - \varphi_j) = \gen {e_{ij}}_\gB$; in particular, $e_{ii} = 1$. We say that a matrix~$A \in \MM_{n,m}(\gB)$ is a \emph {Dedekind \evn matrix} for the~$n$ morphisms $\varphi_1$, \ldots, $\varphi_n$ if each column of $A$ is of the form $\tra {[\, \varphi_1(a) \; \cdots \;\varphi_n(a)\,]}$ for some $a \in \gA$ (depending on the column).  Prove  the existence of a Dedekind \evn matrix whose image contains the vectors $\tra{ [\,e_{1i} \; \cdots \; e_{ni}\,]}$.  In particular, if $\Ann_\gB(\varphi_i - \varphi_j) = 0$ for $i \ne j$, such a matrix is surjective.

}
\end{exercise}

\vspace{-1em}
\begin{exercise}\label{exoArtinAnOtherProof}
{(Another \dem of Artin's \thoz, item 2)}\\
{\rm
The context is that of \thref{thA}:  $(\gk,\gA,G)$ is a \aG and  we want to show the existence of $a_1$, \ldots, $a_r$, $b_1$, \ldots, $b_r \in \gA$ such that for every $\sigma \in G$ we~have
$$\preskip-.6em \postskip.2em\ndsp 
\sum_{i=1}^r a_i \sigma(b_i)=
\formule{1\;\;\mathrm{if}\;\sigma=\Id\\
0\;\;\mathrm{otherwise}.} 
$$
For $\tau \in \G$, $\tau \ne \Id$, show that there exist $m_\tau$ and $x_{1,\tau}$, \ldots, $x_{m_\tau,\tau}$, $y_{1,\tau}$, \ldots, $y_{m_\tau,\tau}$ in~$ \gA$ such that
$$\preskip-.3em \postskip.3em\ndsp 
\;\;\;\sum_{j=1}^{m_\tau} x_{j,\tau} \tau(y_{j,\tau}) = 0, \qquad
\sum_{j=1}^{m_\tau} x_{j,\tau} y_{j,\tau} = 1
.
$$
Conclude the result.
}
\end{exercise}

\vspace{-1em}
\begin{exercise}\label{exoGalExemples} 
{(\aGsz: a few \elrs examples)}\\
{\rm  
Let $(e_1, \ldots, e_n)$ be the canonical basis of~$\gk^n$.  
We make $\rS_n$ act on~$\gk^n$ by permutation of the \coosz: $\sigma(e_i) = e_{\sigma(i)}$ for $\sigma \in \rS_n$.

\emph {1.}
Let $G \subset \rS_n$ be a transitive subgroup of cardinality $n$.
\vspace{-.3em}
\begin{enumerate}\itemsep0pt
\item [\emph{a.}]
Prove  that $(\gk, \gk^n, G)$ is a \aGz.
\item [\emph{b.}]
Give examples.
\end{enumerate}

\emph {2.}
Let~$\gB = \gk(e_1 + e_2) \oplus \gk(e_3 + e_4) \subset \gk^4$ and $G = \gen{(1,2,3,4)}$. \\
Determine  $\Stp_{\rS_4}(\gB)$ and $H = \Stp_G(\gB)$. Do we have~$\gB = (\gk^4)^H$?

\emph {3.}
Let $(\gk, \gA, G)$ be a \aGz. The group $G$ operates naturally on~$\gA[X]$. 
\vspace{-.3em}
\begin{enumerate}\itemsep0pt
\item [\emph{a.}]
Prove  that $(\gk[X], \gA[X], G)$ is a \aGz.
\item [\emph{b.}]
Let~$\gB = X\gA[X] + \gk$ ($\gB$ therefore consists of the \pols of~$\gA[X]$ whose constant \coe is in~$\gk$). Then,~$\gB$ is a \kslg of~$\gA[X]$ which is not of the form $\gA[X]^H$ except in a special case.
\end{enumerate}

}
\end{exercise}

\vspace{-1.5em}
\begin{exercise}
\label{lemPaquesFerrero}
{\rm  Let~$\gk\subseteq\gB\subseteq\gC$ with~$\gB$ \ste over~$\gk$ and~$\gC$ \stfe
over~$\gk$. Suppose that $\rg_\gk(\gB)=\rg_\gk(\gC)$ (\cad  $\gC$ and~$\gB$ have the same \polmu over~$\gk$). 
Then prove that~$\gB=\gC$.
}
\end{exercise}

\vspace{-1em}
\begin{exercise}
\label{exoCorGalste}
{\rm Base yourself on Exercise~\ref{lemPaquesFerrero} and prove the Galois correspondence (\thref{thCorGalGen})
between the finite subgroups of $G$ and the \stes \kslgs of~$\gA$
when $\gA$ is connected. 
}
\end{exercise}

\vspace{-1em}
\begin{exercise}\label{exoIdeauxGlobInvariants} 
{(\aGsz: globally invariant \idsz)}\\
{\rm  
Let $(\gA, \gB, G)$ be a \aGz. We say that an \id $\fc$ of~$\gB$ is \emph{globally invariant} if $\sigma(\fc) = \fc$ for every $\sigma \in G$.

\emph {1.}
Prove  that $\fc$ is generated by invariant \eltsz, \cad by \elts of~$\gA$.

\emph {2.}
More precisely, consider the two transformations between \ids of~$\gA$ and \ids of~$\gB$: $\fa \mapsto \fa\gB$ and $\fc \mapsto \fc \cap \gA$.  Prove  that they establish a non-decreasing bijective correspondence between \ids of~$\gA$ and globally invariant \ids of~$\gB$.
}
\end{exercise}


\begin{problem}
\label{exoLuroth1} \label{NOTAPGL} (L\"uroth's \thoz)\\
{\rm  Let~$\gL=\gk(t)$ where~$\gk$ is a \cdi and $t$ an \idtrz.
If $g=u/v\in\gL$ is a nonconstant \ird fraction  
 ($u,v\in\gk[t]$, coprime),
we define the \emph{height} of $g$ (with respect to $t$) by $\hauteur_t(g)\eqdefi \max\big(\deg_t(u),\deg_t(v)\big)$.
\begin{enumerate}
\item \emph{(Direct part of L\"uroth's \thoz)} 
Let~$\gK=\gk(g)\subseteq\gL$.
Prove  that~$\gL/\gK$ is an \agq extension of degree $d=\hauteur(g)$. More \prmtz, $t$ is \agq over~$\gK$ and its \polmin is, up to multiplicative factor 
in~$\gK\eti$, equal to $u(T)-g v(T)$.
Thus, every nonconstant \coe of $\Mip_{\gK,t}(T)$, $a\in\gK=\gk(g)$ is of the form $a=\fraC{\alpha g+\beta}{\gamma g+\delta}$ with $\alpha\delta-\beta\gamma\in\gk\eti$, and~$\gk(a)=\gk(g)$.  
\item Let $f\in\gL$ be an arbitrary \eltz.
Give an explicit formula using the resultants to express~$f$ as a~$\gK$-\coli of $(1,t,\ldots,t^{d-1})$.
\item If $h$ is another \elt of~$\gL\setminus\gk$ 
show that

\snic{\hauteur\big(g(h)\big)=\hauteur(g)\hauteur(h).}
%
Prove  that every~$\gk$-\alg \homoz~$\gL\to\gL$ is of the form $f\mapsto f(h)$ for some $h\in\gL\setminus\gk$. Deduce a precise description of $\Aut_\gk(\gL)$ by means of fractions of height $1$.

\item We denote by $\PGL_n(\gA)$ the quotient group $\GLn(\gA)/\Ati$ (where $\Ati$ is identified with the subgroup of \iv homotheties via $a\mapsto a\In$).  
To a matrix

\snic{A = \cmatrix {a &b\cr c&d\cr} \in \GL_2(\gA),}

we associate the~$\gA$-\autoz\footnote{We denote by $\gA(t)$ the Nagata \ri of $\gA$ which is obtained from $\gA[t]$ by inverting the primitive \polsz.}

\snic{\varphi_A : \gA(t) \to \gA(t),\quad t \mapsto {at + b\over ct + d}. }

We have  $\varphi_A \circ \varphi_B =
\varphi_{BA}$ and $\varphi_A = \Id \Leftrightarrow A =
\lambda\I_2 \;(\lambda \in \Ati)$. Thus $A \mapsto \varphi_A$ defines an injective \homo $\PGL_2(\gA)\eo \to \Aut_\gA\big(\gA(t)\big)$.
\\
 Prove  that in the \cdi case we obtain an \isoz.
\item \emph{(Converse part of L\"uroth's \thoz)}
Let $g_1$, \ldots, $g_r\in\gL\setminus \gk$. Prove  that $\gk(g_1,\ldots,g_r)=\gk(g)$ for a suitable $g$. It suffices to treat the~$n=2$ case. We show that $\gL$ is \stf over~$\gK_1=\gk(g_1,g_2)$. We must then have~$\gK_1=\gk(g)$
for any nonconstant \coe $g$ of~$\Mip_{\gK_1,t}(T)$.
\\
NB: Since~$\gL$ is a finite dimensional $\gk(g_1)$-\evcz,
every subfield of~$\gL$ strictly containing~$\gk$ is, in \clamaz, \tfz, therefore of the form~$\gk(g)$.
Our formulation of the converse part of L\"uroth's \tho give the \cov meaning of this assertion.%
\index{height!of a rational fraction}%
\index{Luroth@L\"uroth!\thoz} 

\end{enumerate}
}
\end{problem}

\vspace{-1em}
\begin{problem}\label{exoBuildingFrobAlgebra}
{(Differential operators and Frobenius \algsz)}\\
{\rm  
In the first questions,~$\gk$ is a commutative \riz. The \emph{Hasse derivative} of order $m$ of a \pol of $\kX$ is \fmt defined by $f^{[m]} = {1 \over m!}f^{(m)}$. Similarly, for $\alpha \in \NN^n$, we define $\partial^{[\alpha]}$ over $\kuX = \gk[X_1, \ldots, X_n]$ by
\index{Hasse derivative}
$$
\partial^{[\alpha]} f = {1\over\alpha!} 
{\partial^\alpha f \over \partial X^\alpha}
\quad \hbox {with} \quad \alpha! = \alpha_1!  \,\cdots\, \alpha_n!,\quad
f \in \kuX 
.$$
We then have $\partial^{[\alpha]} (fg) = \sum_{\beta+\gamma=\alpha}
\partial^{[\beta]}(f) \, \partial^{[\gamma]}(g)$. We denote by $\delta^{[\alpha]} : \kuX \to \gk$ the \lin form $f \mapsto \partial^{[\alpha]}(f)(0)$. Thus, $f = \sum_\alpha \delta^{[\alpha]}(f)X^\alpha$. We deduce, by letting $\alpha\le\beta$ for $X^\alpha \divi X^\beta$,
$$
\preskip.2em \postskip.4em 
X^\alpha \centerdot \delta^{[\beta]} = \formule{
\delta^{[\beta-\alpha]} \hbox{ if } \alpha\le\beta
\\
0 \hbox{ otherwise,}}
\qquad
\partial^{[\alpha]}(X^\beta) = \formule{
X^{\beta-\alpha} \hbox{ if } \alpha\le\beta
\\
0 \hbox{ otherwise}.}
 $$
Let $g = \sum_\beta b_\beta X^\beta$.  Evaluating the \emph{\dil \polz}~\hbox{$\sum_\beta b_\beta \partial^{[\beta]}$} at $(\uze)$ , we obtain a \lin form $\delta_g : \kuX \to \gk$, $\delta_g = \sum_\beta b_\beta\delta^{[\beta]}$, then an \id $\fa_g$ of~$\kuX$
$$
\preskip.1em \postskip.4em 
\fa_g = \sotq{f \in \kuX} {f\centerdot \delta_g = 0} \eqdf {\rm def}
\sotq{f \in \kuX} {\delta_g(fu) = 0\ \forall u \in \kuX}
. 
$$
We thus obtain a Frobenius \klg $\kuX\sur{\fa_g}$ (with $\delta_g$ dualizing).

\emph {1.}
Let $f=\sum_\alpha a_\alpha X^\alpha$, $g=\sum_\beta b_\beta X^\beta$. We let $\partial_f : \kuX \to \kuX$ be the \dil operator associated with~$f$, \cad $\partial_f = \sum_\alpha a_\alpha \partial^{[\alpha]}$. Check the following relation between the operator $\partial_f$ and the \lin form $\delta_g$
$$\preskip.4em \postskip.4em \ndsp 
\sum_\gamma (f \centerdot \delta_g)(X^\gamma) X^\gamma =
\partial_f(g) = \sum_{\alpha \le \beta} a_\alpha b_\beta X^{\beta-\alpha}
. 
$$
Deduce that $f\centerdot \delta_g = 0 \iff \partial_f(g) = 0$.

Now we must note that the law $f*g = \partial_f(g)$ provides the additive group $\kuX$ with a  $\kuX$-module structure (in particular because $\partial_{f_1f_2} = \partial_{f_1} \circ \partial_{f_2}$). 
But as~\hbox{$X^\alpha * X^\beta = X^{\beta-\alpha}$} or $0$, certain authors use~$X^{-\alpha}$ instead of $X^{\alpha}$; they provide~$\kuX$ with a 
$\gk[\uX^{-1}]$-module structure. Other authors permute~$\uX$ and~$\uX^{-1}$; they provide~$\gk[\uX^{-1}]$ with a $\kuX$-module structure such that the \idz~$\fa_g$
 (annihilator of $g \in \gk[\uX^{-1}]$) is an \id of a \pol \ri with \idtrs with exponents~$\ge 0$. In the latter formalism, a \polz~$f$ with \idtrs with exponents~$\ge 0$ therefore acts on a \pol $g$ having its \idtrs with exponents $\le 0$ to provide a \pol $f*g$ having \idtrs with exponents $\le 0$ (by deleting the \moms containing an exponent $> 0$).  Thus, if $g = X^{-2} + Y^{-2} + Z^{-2}$, the \id $\fa_g$ of~$\gk[X,Y,Z]$ contains for example $XY$, $X^2 - Y^2$ and every \hmg \pol of degree~$\ge 3$.

\emph {2.}
Let $d \ge 1$. Study the special case of the Newton sum $g = \sum_i X_i^{-d}$, \cad $\delta_g : f \mapsto \sum_i {1 \over d!} {\partial^d f \over \partial X_i^d}(0)$, the sum of the \coes over $X_1^d$, \ldots, $X_n^d$.

In the remainder, we fix $g = \sum_\beta b_\beta X^\beta$, or according to taste, $g = \sum_\beta b_\beta X^{-\beta}$.

\emph {3.}
Prove  that we have an inclusion $\fb \subseteq \fa_g$ for some \id $\fb = \gen {X_1^{e_1}, \cdots, X_n^{e_n}}$ with integers $e_i \ge 1$. In particular, $\kuX\sur{\fb}$ is a free \kmo of finite rank and~$\kuX\sur{\fa_g}$ is a \tf \kmoz.

\emph {4.}
Define a~$\gk$-\lin map $\varphi : \kuX\sur\fb \to \kuX$ such that $\Ker\varphi = \fa_g\sur\fb$. We can therefore compute $\fa_g$ if we know how to solve \lin \syss over~$\gk$.

\emph {5.}
Suppose that~$\gk$ is a \cdi and so~$\gA := \kuX\sur{\fa_g}$ is a finite dimensional \kevz. Prove  that $(\gA, \delta_g)$ is a Frobenius \klgz.

}

\end{problem}

\vspace{-1em}
\begin{problem}\label{exoTh90HilbertAdditif} {(Hilbert's \tho $90$, additive version)}\ihi \\
{\rm
Let $(\gk,\gA,G)$ be a \aG where $G = \gen {\sigma}$ is cyclic
of order $n$.

 \emph {1.}
Considering an \elt $z \in \gA$ of trace $1$, we will show that
$$
\preskip.4em \postskip.4em 
 \gA = \Im(\Id_\gA - \sigma) \oplus \gk z, \qquad
\Im(\Id_\gA - \sigma) = \Ker\Tr_G.
$$
Consequently $\Im(\Id_\gA - \sigma)$ is a \stl \kmo of rank $n-1$.
You can use the family of \endos $(c_i)_{i \in \lrb{0..n}}$,
$$
\preskip.4em \postskip.4em \ndsp
c_0 = 0,\  c_1(x)=x,\ c_2(x) = x+\sigma(x),\ \ldots,\
c_i(x) = \sum_{j=0}^{i-1} \sigma^j(x),\ \ldots
$$

\sni \emph {2.}
For $x\in\gA$ prove that to be of the form $y - \sigma(y)$, 
it is necessary and sufficient that $\Tr_G(x) = 0$.

 \emph {3.}
More \gnltz, let $(c_\tau)_{\tau\in G}$ be a family in~$\gA$. Prove  that there exists an \eltz~$y$ such that $c_\tau = y-\tau(y)$ \ssi the family satisfies the following additive cocycle condition: for all $\tau_1$, $\tau_2\in G$: $c_{\tau_1\tau_2} = \tau_1(c_{\tau_2}) + c_{\tau_1}$.

 \emph{4.}
Assume that $n$ is a prime number $p$ and that $p=0$ in~$\gk$. Prove  the existence of some $y \in \gA$ such that $\sigma(y) = y+1$. \\
Deduce that $(1, y, \dots, y^{p-1})$ is a~$\gk$-basis of~$\gA$ and that the \polcar of $y$ is of the form $Y^p - Y - \lambda$ with~$\lambda \in \gk$.  \\
We therefore have~$\gA = \gk[y] \simeq \aqo{\gk[Y]}{Y^p - Y - \lambda}$ (Artin-Schreier extension).

 \emph{5.}
Give a converse of the previous item.
}
\end{problem}

\vspace{-1em}
\begin{problem}\label{exoHeitmannGaloisExemple} 
 {(\aGsz: study of an example)}
{\rm  
Consider a \riz~$\gB$ in which $2$ is \ivz, with  $x$, $y \in \gB$ and $\sigma \in \Aut(\gB)$ of order $2$ satisfying $x^2 + y^2 = 1$, $\sigma(x) = -x$ \hbox{and $\sigma(y) = -y$}. We can take as an example the \riz~$\gB$ of continuous functions over the unit circle $x^2 + y^2 = 1$ and for~$\sigma$ the involution $f \mapsto \{(x,y) \mapsto f(-x,-y)\}$.  Let~$\gA = \gB^{\gen {\sigma}}$ (sub\ri of the \gui{even functions}).

\emph {1.}
Prove  that $(\gA, \gB, \gen\sigma)$ is a \aGz.
\\
Consequently,~$\gB$ is a \Amrc $2$.

\emph {2.}
Let $E = \gA x + \gA y$ (submodule of the \gui{odd functions}). \\
Check that~$\gB = \gA \oplus E$ and that $E$ is a \Amrcz~$1$.

\emph {3.}
Let $x_1 = 1$, $x_2 = x$, $x_3 = y$ such that $(x_1, x_2, x_3)$ is a \sgr of the \Amoz~$\gB$. Make $y_1$, $y_2$, $y_3 \in \gB$ explicit as in Lemma~\ref{lemArtin}, \cad $\big((x_i)_{i \in \lrb {1..3}}, (y_i)_{i \in \lrb {1..3}}\big)$ is a \stycz. \\
Deduce a \mprn $P \in \MM_3(\gA)$ of rank $2$ with~$\gB \simeq_\gA \Im P$.

\emph {4.}
Let $R = \Cmatrix{.4em} {x &-y\cr y & x\cr} \in \SL_2(\gB)$. Prove  that this \gui{rotation} $R$ induces an \iso of \Amos between $E^2$ and~$\Ae2$

\snic {
\cmatrix {f\cr g\cr} \mapsto R \cmatrix {f\cr g\cr} =
\cmatrix {xf-yg\cr yf+xg\cr}
.}

Consequently (next question), $E\te_\gA E \simeq \gA$; 
verify that $f \te g \mapsto fg$ realizes an \iso of \Amos of $E \te_\gA E$ over~$\gA$.

\emph {5.} For some \Amo $M$ ($\gA$ arbitrary), let 

\snic{M^{2\te}=M\te_\gA M,\,  M^{3\te}=M\te_\gA M\te_\gA M, \hbox{ etc}\dots}

Let $E$ be an \Amo satisfying $E^n \simeq \gA^n$ for some $n \ge 1$. Prove  that $E$ is a \Amrc $1$ and that $E^{n\te} \simeq \gA$.

\emph {6.}
Let $\fa$ be the \id of~$\gA$ defined by $\fa = \gen {xy, x^2}$.  Check that $\fa^2 = x^2\gA$ (so if $x$ is \ndzz, $\fa$ is an \iv \id of~$\gA$), that $\fa\gB$ is principal and finally, that $\fa$, regarded as an \Asub of~$\gB$, is equal to $xE$.

\emph {7.}
Let~$\gk$ be a nontrivial \ri with $2\in\gk\eti$ and $\gB = \aqo{\gk[X,Y]}{X^2 + Y^2 - 1}$. We write $\gB=\gk[x,y].$
We can apply the above by taking $\sigma$ as defined by $\sigma(x) = -x$ and $\sigma(y) = -y$. Suppose that $\alpha^2+\beta^2=0\Rightarrow \alpha=\beta=0$ in~$\gk$ (for example if~$\gk$ is a \cdi and $-1$ is not a square in~$\gk$).
\begin{enumerate}
\item [\emph {a.}]
Prove  that $\Bti = \gk^{\times}$; illustrate the importance of the hypothesis \gui{of reality} made about~$\gk$.

\item [\emph {b.}]
Prove  that $\fa$ is not principal and so $E$ is not a free \Amoz. Deduce that~$\gB$ is not a free \Amoz.
\end{enumerate}

\emph {8.}
Let~$\gB$ be the \ri of (real) continuous functions over the unit circle $x^2 + y^2 = 1$ and $\sigma$ the involution $f \mapsto \{(x,y) \mapsto f(-x,-y)\}$.
Prove  that $\fa$ is not principal and that~$\gB$ is not a free \Amoz.
}

\end{problem}

\sol

\exer{exothEtalePrimitif}
We have~$\gB=\Kxn$, with $\dex{\gB:\gK}=m$. We will perform a computation that shows that the \Klgz~$\gB$ is monogenic or contains an \idm $e\neq 0,1$. 
In the second case,~$\gB\simeq \gB_1\times \gB_2$, with $\dex{\gB_i:\gK}=m_i<m $, $m_1+m_2=m$, which allows us to conclude by \recu on $m$.
\\
If we are able to treat the $n=2$ case, we are done, because~$\gK[x_1,x_2]$ is \'etale over~$\gK$, so either we replace~$\gK[x_1,x_2]$ with~$\gK[y]$ for some $y$, or we find an \idm $e\neq 0,1$ within it. The \dem of item \emph{1} of \thref{thEtalePrimitif} shows that an \'etale \Klgz~$\gK[x,z]$ is monogenic if~$\gK$ contains an infinite sequence of distinct \eltsz.
It uses a \pol $d(a,b)$ which, evaluated in~$\gK$ must give an \iv \eltz.
If we do not have any information on the existence of an infinite sequence of distinct \elts of~$\gK$, we enumerate the integers of~$\gK$ until we obtain $\alpha,\beta$ in~$\gK$ with $d(\alpha,\beta)\in\gK\eti$, or until we conclude that the \cara is equal to a prime number $p$. We then enumerate the powers of the \coes of~$f$ and of $g$ (the \polmins of $x$ and $z$ over~$\gK$) until we obtain $\alpha,\beta$ in~$\gK$ with $d(\alpha,\beta)\in\gK\eti$, or until we conclude that the field~$\gK_0$ generated by the \coes of~$f$ and $g$ is a finite field. In this case,~$\gK_0[x,z]$ is a reduced finite $\gK_0$-\algz. It is a reduced finite \riz, so either it is a finite field, of the form~$\gK_0[\gamma]$, and~$\gK[x,z]=\gK[\gamma]$, or it contains an \idm $e\neq 0,1$. 

\smallskip 
\rem The reader will be able to verify that the \dem transformation that we put the \gui{$\gB$ is a \cdiz} case through is precisely the implementation of the \elgbm of \zedrs \risz. 
In fact the same machinery also applies to the \cdiz~$\gK$ and provides the following result: a \ste \alg over a \zedr \riz~$\gK$ (\Dfnz~\ref{defdualisante}) is a finite product of \stes \Klgsz.    
\eoe


\exer{exoKnEltPrimitif}
\emph{1.} We write $x=(\xn)=\sum_{i=1}^nx_ie_i$ and identify~$\gA$ with a sub\ri of~$\gB$ by $1\mapsto(1,\ldots,1)$.
By writing $e_i \in \gA[x]$, we obtain that the \elts $x_i - x_j$ are \ivs for $j \ne i$. Conversely, if $x_i - x_j$ is \iv for all~$i \ne j$, we have~$\gB = \gA[x]=\gA \oplus \gA x\oplus \cdots \oplus \gA x^{n-1} $ (Lagrange interpolation, Vandermonde \deterz).

 \emph{2.} If and only if $\#\gA \ge n$.


\exer{exoIdmChangBase} 
\emph{1} and \emph{2.} If $(a_1,\ldots,a_\ell)$ is a basis of~$\gB$ over~$\gK$, it is also a basis of~$\gB(v)$ over~$\gK(v)$.

\emph{3.}
Let $b/p$ be an \idm of~$\gB(v)$: we have $b^2=bp$. If $p(0)=0$, then $b(0)^2=0$, and since~$\gB$ is reduced, $b(0)=0$. We can then divide $b$ and $p$ by $v$. Thus, we can assume that $p(0)\in\gK\eti$. By dividing $b$ and $p$ by $p(0)$ we are reduced to the case where $p(0)=1$. We then see that $b(0)$ is \idmz. We denote it by $b_0$ and let~$e_0=1-b_0$. Let us write $e_0b=vc$. We multiply the \egt $b^2=bp$ by $e_0=e_0^2$ and we obtain $v^2c^2=vcp$. So $vc(p-vc)=0$, and since the \pol $p-vc$ has $1$ as its constant term, so is \ndzz, this gives us $c=0$. Therefore $b=b_0b$. Let us reason modulo $e_0$ for a moment: we have $b_0\equiv1$ so $b$ is primitive and the \egt $b^2=bp$ is simplified to $b\equiv p \mod e_0$. 
This gives the \egt $b=b_0b=b_0p$ in~$\gB(v)$ and so $b/p=b_0$.


\exer{exoBiquadratique} ~\\
\emph{1.} Classical: it is $\ZZ[\sqrt d]$ if $d \equiv 2$ or $3 \bmod 4$ and $\ZZ[{1 + \sqrt d \over 2}]$ if $d \equiv 1 \bmod 4$.

\emph{2.} We have~$\gA = \ZZ[\sqrt a]$.
We have $\beta^2 = {a+1 \over 2} + \sqrt a \in \gA$, therefore $\beta$ is integral over $\gA$, then over~$\ZZ$. Actually, $(\beta^2 - {a+1 \over 2})^2 = a$ and $\beta$ is a root of $X^4 - (a+1) X^2 + ({a-1 \over 2})^2$.
We thus find $(\sigma\tau)(\beta) = \beta - \sqrt 2$.

\emph{3.}
We find $(\sigma\tau)(z) = r - (s+u)\sqrt 2 + u\beta$ then $z + (\sigma\tau)(z) = 2r + u\sqrt {2a}$. This last \elt of $\QQ(\sqrt {2a})$ is integral over $\ZZ$, hence in $\ZZ[\sqrt {2a}]$ because $2a \equiv 2 \bmod 4$.
\\
Hence $u \in \ZZ$ (and $2r \in \ZZ$). We replace $z$ with $z - u\beta$ which is integral over $\ZZ$, \cad~$z = r + s\sqrt 2 + t \sqrt a$. We have $\sigma(z) = r-s\sqrt 2 + t\sqrt a$, $\tau(z) = r+s\sqrt 2 - t\sqrt a$; by using $z+\sigma(z)$ and $z+\tau(z)$, we see that $2r$, $2s$, $2t \in \ZZ$. Let us use

\snuc {
z\sigma(z) = x + 2rt\sqrt a,\ z\tau(z) = y + 2rs\sqrt 2,
\hbox { with } x = r^2 - 2s^2 + at^2,\  y = r^2 + 2s^2 - at^2
.}

We therefore have $x$, $y \in \ZZ$ then $x + y = 2r^2 \in \ZZ$, $x - y = 2at^2 - (2s)^2 \in \ZZ$, so $2at^2 \in \ZZ$. From $2r$, $2r^2 \in \ZZ$, we deduce $r \in \ZZ$. Similarly, from $2t$, $2at^2 \in \ZZ$ (using that $a$ is odd), we see that $t \in \ZZ$, and then finally $s \in \ZZ$. Phew!

Thanks to the $\ZZ$-basis of~$\gB$, we obtain $\Disc_{\gB/\ZZ} = 2^8 a^2$.

\emph{4.}
We have
%
$\gB = \ZZ \oplus \ZZ\sqrt a \oplus \ZZ\sqrt 2 \oplus \ZZ \beta =
\gA \oplus E$ with $E = \ZZ\sqrt 2 \oplus \ZZ \beta
.
$\\
%
We also have $2E = \sqrt2\,\fa$ with $\fa = 2\,\ZZ  \oplus \ZZ(\sqrt a -1) = \gen {2, \sqrt a-1}_\gA$. This proves on the one hand that $E$ is an~$\gA$-module, and on the other that it is \isoc to the \idz~$\fa$ of~$\gA$. Consequently, $E$ is a \Amrcz~$1$. The expression~$\gB =\gA\oplus E$ certifies that~$\gB$ is a \ptf \Amoz, written as a direct sum of a free \Amo of rank $1$ and of a \mrcz~$1$. 
In \gnlz, the \id $\fa$ is not principal, so $E$ is not a free \Amoz.
Here is a small sample of values of $a \equiv 3 \bmod 4$; we have underlined those values for which the ideal $\fa$ is principal:

\snic {
-33,\
-29,\
-21,\
-17,\
-13,\
-5,\
\und {-1},\
\und {3}, \
\und {7},\
\und {11},\
15,\
\und {19},\
\und {23},\
\und {31},\
35
.}

In the case where $\fa$ is not principal,~$\gB$ is not a free $\gA$-module. Otherwise, $E$ would be stably free of rank $1$, therefore free (see Exercise~\ref{exoStabLibRang1}).
Finally, we always have $\fa^2 = 2\gA$ (see below), so $\fa \simeq \fa^{-1} \simeq \fa\sta$. 
Consequently~$\gB \simeq_\gA \gB\sta$. 
\\
Justification of $\fa^2 = 2\gA$: always in the same context ($a \equiv 3 \bmod 4$ so~$\gA = \ZZ[\sqrt {a}\,]$), we have for $m \in \ZZ$

\snic {
\gen {m, 1+\sqrt a} \gen {m, 1-\sqrt a} = \pgcd(a-1,m)\,\gA
.}

Indeed, the left \id is generated by $\big(m^2, m(1 \pm \sqrt a), 1-a\big)$, 
each one a
multiple of the gcd. This \id contains $2m = m(1 + \sqrt a) + m(1 - \sqrt a)$, so it contains the \elt $\pgcd(m^2, 2m, 1-a) = \pgcd(m, 1-a)$, (the \egt is  due to $a \equiv 3 \bmod 4$). For $m = 2$, we have $\gen {2, 1+\sqrt a} = \gen {2, 1-\sqrt a} = \fa$ and $\pgcd(a-1, 2) = 2$.


\exer{exoTensorielDiscriminant}
We consider~$\gB = \gA\otimes_\gk\gA'$ as an \Algz, a \eds to $\gA$ of the \klgz~$\gA'$; it is free of rank $n'$. We therefore have at our disposal a stack of free \algsz~$\gk \to \gA \to \gB$ and the transitivity formula of the discriminant provides
$$\preskip.0em \postskip.3em 
 \Disc_{\gB\sur\gk}(\ux\otimes\ux') = \Disc_{\gA\sur\gk}(\ux)^{n'}
\cdot \rN_{\gA\sur\gk}\big(\Disc_{\gB\sur\gA}(1 \otimes \ux')\big)
.$$
But $\Disc_{\gB\sur\gA}(1 \otimes \ux') = \Disc_{\gA'\sur\gk}(\ux')$. As it is an \elt of~$\gk$, its norm~$\rN_{\gA\sur\gk}$ has the value~$\Disc_{\gA'\sur\gk}(\ux')^n$. Ultimately we obtain the \egt
$$\preskip.2em \postskip.5em 
\Disc_{\gB\sur\gk}(\ux\otimes\ux') =
\Disc_{\gA\sur\gk}(\ux)^{n'} \, \Disc_{\gA'\sur\gk}(\ux')^n.
$$

\exer{CyclicNormalBasis} 
We will use the following classical result on \lin \algsz. Let~$E$ be a
finite dimensional  \Kev  and $u\in\End_\gK(E)$. If $d$ is the degree of the \polmin of $u$, there exists an $x \in E$ such that the \elts  
$x$, $u(x)$, \ldots, $u^{d-1}(x)$ are~$\gK$-\lint independent.\\
Here $\dex{\gL : \gK} = n$, and $\Id_\gL$, $\sigma$, \dots, $\sigma^{n-1}$ are~$\gK$-\lint independent, so the \polmin of $\sigma$ is $X^n - 1$, of degree $n$. We apply the above result.


\exer{exoHomographieOrdre3}
\emph{1.} $A$ is the companion matrix of the \pol $X^2 - X + 1 = \Phi_6(X)$, so $A^3=-\I_2$ in $\GL_2(\gk)$ and $A^3 = 1$ in $\PGL_2(\gk)$. 

\emph{2.}
We know by Artin's \tho that~$\gk(t)/\gk(t)^G$ is a Galois extension with Galois group $\rA_3$.
The computation gives 

\snic{g = t+\sigma(t)+\sigma^2(t)= {t^3 - 3t - 1 \over t(t+1)}.}

We obviously have $g\in\gk(t)^G$ and $t^3 - gt^2 - (g+3)t - 1=0$. 
Therefore, (direct part of L\"uroth's \thoz, \Pbmz~\ref{exoLuroth1}) $\dex{\gk(t):\gk(g)}=3$, and~$\gk(t)^G=\gk(g)$.

 \emph{3.}
Since~$\gk(a)\simeq\gk(g)$ and $f_g(t)=0$, the extension~$\gk(a)\to\aqo{\gk[T]}{f_a}$ is a copy of the extension~$\gk(g)\to\gk(t)$.

 \emph{4.}
Let $\sigma$ be a \gtr of $\Aut(\gL/\gK)$. This question amounts to saying that we can find some $t\in\gL\setminus \gK$ such that $\sigma(t)={-1\over t+1}~(*)$. Since $t$ must be of norm $1$, we seek it of the form $t={\sigma(u)\over u}$.
The computation then shows that $(*)$ is satisfied under the condition that $u\in\Ker(\Tr_G)$. It remains to show that there exists some $u\in\Ker(\Tr_G)$ such that ${\sigma(u)\over u}\notin\gK$. This amounts to saying that the restriction of $\sigma$ to $E =\Ker(\Tr_G)$ is not a homothety. However, $E\subseteq\gL$ is a \Ksv of dimension $2$, stable under $\sigma$.
By Exercise~\ref{CyclicNormalBasis}, the \Kevz~$\gL$ admits a \gtr for the \endo $\sigma$. This \lin \alg \prt remains true for every stable subspace by $\sigma$.

\exer{exoGroupAlgebra} 
The \elts $g \te h$ form a~$\gk$-basis of $\env{\gk}{\gA}$.  \\
Let $z = \sum_{g,h} a_{g,h} g\te h$ with $a_{g,h} \in \gk$.  Then, $z \in \Ann(\rJ\iAk)$ \ssi we~have~$g' \cdot z = z \cdot g'$ for every $g' \in G$. We obtain $a_{g,h} = a_{1,gh}$, so $z$ is a~$\gk$-\coli of the $z_k \eqdf {\rm def} \sum_{gh = k} g \otimes h$.  \\
Conversely, we see that $z_k \in \Ann(\rJ\iAk)$ and we have $z_k = k \cdot z_1 = z_1 \cdot k$. \\
So
$\Ann(\rJ\iAk)$ is the \kmo generated by the $z_k$'s, and it is the \Amo (or the \id of $\env{\gk}{\gA}$) generated by $z_1 = \sum_g g\te g^{-1}$.
\\
The image under $\mu\iAk$ of $\Ann(\rJ\iAk)$ is the \id $n\gA$. 
Regarding the trace, we have $\Tr\iAk(g) = 0$ if $g \ne 1$.
Therefore $\Tr\iAk\!\big(\sum_g a_g g\big) = na_1$. 
\\
If $a = \sum_g a_g g$ and $b = \sum_g b_g g$, then $\Tr\iAk(ab) = n\sum_g a_g b_{g^{-1}}$.
\\
The \eqvcs of item \emph{2} are therefore clear, and in the case where $n\in\gk\eti$, the \idst is $n^{-1}\sum_g g\te g^{-1}$.

\emph {3.}
Let $\lambda : \gk[G] \to \gk$ be the \lin form \gui{\coo over $1$.}  For $g$, $h \in G$, we~have~$\lambda(gh) = 0$ if $h \ne g^{-1}$, and~$1$ otherwise. 
So, $\lambda$ is dualizing and  
$(g^{-1})_{g \in G}$ is the dual basis of $(g)_{g \in G}$ with respect to $\lambda$.  
We have $\Tr_{\gk[G]\sur\gk} = n\cdot\lambda$.

\exer{exo1Frobenius}
\emph{1.} 
It suffices to do it for $g \in \{1, x, \ldots, x^{n-1}\}$, which is a basis of $\gA[Y]$ over $\gk[Y]$.  
The right-hand side of $(*)$ with $g=x^{i}$ is 

\snic{h_i=\varphi\big((Y-x)x^i\big) =
\varphi(Yx^i -x^{i+1}) = Y^{i+1} - \varphi(x^{i+1}).}

If $i < n-1$, we have $\varphi(x^{i+1}) = Y^{i+1}$, so $h_i=0$. For $i = n-1$, we have

\snic {
\varphi(x^n) = -\varphi(a_0 + a_1x + \cdots + a_{n-1}x^{n-1}) = 
-(a_0 + a_1Y + \cdots + a_{n-1}Y^{n-1}),
}

and $h_n(Y)=f(Y)$, which gives the result.

 \emph{2.}
For $i < n$, we have 

\snic{ 
\begin{array}{ccc} 
f(Y)\wi\lambda\big(x^i f'(Y)\big) =
\varphi\big((Y-x)x^i f'(Y)\big) =  \varphi\big(x^i f(Y)\big) = Y^if(Y), \hbox{ \cad}  \\[1mm] 
\wi\lambda\big(x^i f'(Y)\big) = \sum_{j < n} \lambda(x^i b_j) Y^j = Y^i. 
\end{array}
}

Therefore $(b_j\centerdot\lambda)(x^i)=\lambda(x^i b_j) = \delta_{ij}$. Thus, $\lambda$ is dualizing.

 \emph{3.}
For two dual bases $(e_i)$, $(\alpha_i)$, we have $\Tr\iAk = \sum e_i\centerdot\alpha_i$. With the two dual bases $(1, x, \ldots, x^{n-1})$ and $(b_0\centerdot\lambda, b_1\centerdot\lambda, \ldots, b_{n-1}\centerdot\lambda)$ we obtain

\snic {
\Tr\iAk =  b_0\centerdot\lambda + x b_1\centerdot\lambda + \cdots + 
x^{n-1} b_{n-1}\centerdot \lambda = f'(x)\centerdot\lambda.
}


\exer{exoFrobeniusAlgExemples} 
\emph {1.}
The \klgz~$\gA := \gk[X_1, \ldots, X_n]\sur{\gen {f_1(X_1), \ldots, f_n(X_n)}}$ is the tensor product of the $\aqo{\gk[X_i]}{f_i(X_i)}$ which are Frobenius \algsz, so~$\gA$ is a Frobenius \algz.
Precision with $d_i = \deg(f_i)$. The \klg $\gA$ is free of rank $d_1 \cdots d_n$, the \moms $x^{\alpha}=x_1^{\alpha_1} \cdots x_n^{\alpha_n}$ with $\alpha_i < d_i$ for every $i$ form a~$\gk$-basis, and the
\lin form \gui{\coo over $x_1^{d_1-1} \cdots x_n^{d_n-1}$} 
is dualizing.

\emph {2.}
Let $\delta_0$, $\delta_x$, $\delta_y$ be three \lins forms over~$\gk[X,Y]$ defined by

\snic {
\delta_0(f) = f(0), \quad \delta_x(f) = f'_X(0), \quad \delta_y(f)=f'_Y(0).
}

Viewed as linear forms over $\gA$ 
they define a~$\gk$-basis of $\Asta$, a dual basis of the~$\gk$-basis $(1, x, y)$ of~$\gA$. We have 

\snic {
x \centerdot \delta_x = y \centerdot \delta_y = \delta_0, 
}


and so $\Asta = \gA\centerdot\delta_x + \gA\centerdot\delta_y$. Let us show that $G = \Cmatrix{.3em} {x\cr -y\cr}$ is a  \mpn of $\Asta$ for $(\delta_x, \delta_y)$. 
We must observe that for $u$, $v$ in~$\gA$ we have the implication
$$
\preskip.4em \postskip.4em 
u \centerdot \delta_x + v \centerdot \delta_y = 0 \;\;\Longrightarrow\;\;
\cmatrix {u\cr v\cr} \in \gA \Cmatrix{.3em} {x\cr -y\cr}. 
$$
By multiplying $u \centerdot \delta_x + v \centerdot \delta_y = 0$ by $x$, we obtain $u \centerdot \delta_0 + (xv) \centerdot \delta_y = 0$; we evaluate at $1$ and we put $x := 0$ to obtain $u(0,y) = 0$, \cad $u \in \gA x$. Similarly, $v \in \gA y$. If we write $u = xr$, $v = ys$, we obtain $r \centerdot \delta_0 + s \centerdot \delta_0 = 0$, \cad $r + s = 0$, as required.

The \idd $\cD_1(G) = \gen {x,y}$ is nonzero, with a null square, so it cannot be generated by an \idmz. Consequently, the \Amo $\Asta$ is not \proz. A fortiori, it is not free.

\exer{exoAlgMonogeneJAnnJ} 
\emph {1.}
We have $(y-z)f^\Delta(y,z) = 0$ so $\delta := f^\Delta(y,z) \in \Ann(\rJ)$.
We then know that for $\alpha \in \env{\gk}{\gA}$, we have $\alpha\delta = \mu\iAk(\alpha) \cdot \delta = \delta \cdot \mu\iAk(\alpha)$. We apply this result to $\alpha = \delta$ by noticing that $\mu\iAk(\delta) = f'(x)$.

 \emph {2.}
We write $f(Y) - f(Z) = (Y-Z)f'(Y) - (Y-Z)^2 g(Y,Z)$, which gives us in the \alg $\env{\gk}{\gA}$ the \egt $(y-z)f'(y) = (y-z)^2g(y,z)$. 
We write the \egt $1=uf+vf'$ in~$\AX$.
Then $f'(y)v(y)=1$, so $y-z = (y-z)^2v(y)g(y,z)$. \\
When $a = a^2b$, the \elt $ab$ is \idm and~$\gen {a} = \gen {ab}$. Therefore~$\rJ=\gen{e}$ with the \idm $e = (y-z)v(y)g(y,z)$.
\\
We have $f^\Delta(Y,Z) = f'(Y) - (Y-Z)g(Y,Z)$, so
$$\preskip.2em \postskip.3em 
f^\Delta(y,z) = f'(y) - (y-z)g(y,z) = f'(y)\big(1  - (y-z)v(y)g(y,z)\big) = 
f'(y) (1-e). 
$$


\exer{exoJacInversibleAlgSeparable} 
\emph {1.}
Let $f'_{ij}=\partial f_i/\partial X_j$ and we write, in~$\gk[\uY,\uZ]$,
$$\preskip.4em \postskip.4em\ndsp 
f_i(\uY) - f_i(\uZ) - \sum_j (Y_j - Z_j) f'_{ij}(\uY) =: - g_i(\uY,\uZ) \in \gen {Y_1-Z_1, \ldots,
Y_n-Z_n}^2. 
$$
 In $\env{\gk}{\gA}$, by letting $A=\JJ_\uy(f_1, \ldots, f_n)$ we obtain 

\snuc {
A \Cmatrix{.3em} {y_1-z_1\cr \vdots\cr y_n-z_n\cr} =
\Cmatrix{.3em} {g_1(\uy,\uz) \cr \vdots\cr g_n(\uy,\uz)\cr}
\quad \hbox {with} \quad
g_i(\uy,\uz) \in \gen {y_1-z_1, \ldots, y_n-z_n}^2 = \rJ\iAk^2.
}

By inverting $A$, we obtain $y_i - z_i \in \rJ\iAk^2$, \cad $\rJ\iAk = \rJ\iAk^2$.

\emph {2.}
As $\mu\iAk\big(\beta_{\uy,\uz}(\uf)\big) = \J_\ux(\uf)$, we have $\mu\iAk(\vep) = 1$. \\
As $\vep \in \Ann(\rJ\iAk)$, $\vep$ is the \idm \gtr of $\Ann(\rJ\iAk)$. \\
Finally,~$\beta_{\uy, \uz}(\uf)$, which is associated with $\vep$, is also a \gtr of $\Ann(\rJ\iAk)$.

\emph {3.}
Let $f_1$, \ldots, $f_n$ be in $\gk[\uX]$ and $\delta = \J_{\uX}(\uf)$. Let us invert $\delta$ with an in\-de\-ter\-mi\-nate~$T$. 
Then, in $\gk[\uX,T]$ we obtain
$$\preskip-.4em \postskip.4em
\JJ_{\uX,T}(\uf,\delta T-1) =
\bordercmatrix [\lbrack\rbrack]{
             &\partial_{X_1}  &\cdots         &\partial_{X_n} &\partial_T\cr
\quad f_1    &                &               &       &0 \cr
\ \quad\vdots&                &\JJ_{\uX}(\uf) &       &\vdots\cr
\quad f_n    &                &               &       &0\cr
\delta T-1   &\star           &\cdots         &\star  &\delta\cr
}
$$ 
and so $\J_{\uX,T}(\uf,\delta T-1) = \delta^2$.  Let 
$$\preskip.3em \postskip.4em 
\gA =
\aqo{\gk[\uX]}{\uf}\;\hbox{  and  }\;\gB = \gA[\delta^{-1}] =
\aqo{\gk[\uX,T]}{\uf,1-\delta T}. 
$$
Then the Jacobian of the \sys $(\uf,\delta T-1)$ which defines $\gB$ is \iv in $\gB$ and so $\gB$ is a \spb \algz.

\exer{exoSepAlgDedekindLemma} 
The~\Blgz~$\gB\te_\gk\gA$ is \spbz. We have a transformation (\uvle \prt of the \edsz)

\snic{\Hom_\gk(\gA,\gB) \to \Hom_\gB(\gB\te_\gk\gA,\gB)$, $\psi \mapsto
\ov\psi,}

defined by $\ov\psi(b\otimes a) = b\psi(a)$.

\emph {1.}
We then consider the \idm $\vep_{\ov\varphi} \in \gB\te_\gk\gA$ of 
Lemma~\ref{lemIdmHomSpb}, and we write it in the form $\vep_{\ov\varphi} = \sum_{i \in I} b_i\te a_i$.

\emph {2.}
Directly results from Lemma~\ref{lemIdmHomSpb}: the \idm $e$ is none other than $e_{\{\ov\varphi, \ov{\varphi'}\}}$.

\emph {3.}
Since the horizontal juxtaposition of Dedekind \evn matrices is a Dedekind \evn matrix, it suffices to show that there exists one, say $A_1$, whose image contains the vector $v := \tra{[ \,e_{11} \;\cdots \; e_{n1}\,]}$. \\
Let $\big((a_j)_{j \in \lrbm}, (b_j)_{j \in \lrbm}\big)$ be the pair attached to $\varphi_1$.
We put in column $j$ of~$A_1$ the vector $\tra {[\,\varphi_1(a_j) \; \cdots \; \varphi_n(a_j)\,]}$. We then have $A_1 \tra {[ \,b_1 \; \cdots \; b_m\,]} = v$.


\exer{exoArtinAnOtherProof}
By hypothesis, for each $\tau\in G \setminus \{\Id\}$ there exist $n_\tau\in\NN$ and $x_{1,\tau}$, \ldots, $x_{n_\tau,\tau}$, $y_{1,\tau}$, \ldots, $y_{n_\tau,\tau} \in \gA$ such that 
%
$1 = \sum_{j=1}^{n_\tau} x_{j,\tau}\big(y_{j,\tau} - \tau(y_{j,\tau})\big)$.
%
Let $s_\tau = \sum_{j=1}^{n_\tau} x_{j,\tau} \tau(y_{j,\tau})$ such that~$\sum_{j=1}^{n_\tau} x_{j,\tau}y_{j,\tau} = 1+s_\tau$, then we define $x_{n_\tau+1,\tau} = -s_\tau$ and $y_{n+1,\sigma} = 1$. Then, with $m_\tau=1+n_\tau$ 

\snic {
\sum_{j=1}^{m_\tau} x_{j,\tau} \tau(y_{j,\tau}) = s_\tau - s_\tau = 0, \qquad
\sum_{j=1}^{m_\tau} x_{j,\tau} y_{j,\tau} = 1+s_\tau - s_\tau = 1.
}

Fixing $\sigma\in G$, we obtain the product

\snic{
\prod_{\tau \in G \setminus \{\Id\}}
  \sum_{j=1}^{m_\tau} x_{j,\tau} \sigma(y_{j,\tau})=
\formule{1\;\;\mathrm{if}\;\sigma=\Id\\
0\;\;\mathrm{otherwise}.}
}

The development of the product provides two families $(a_i)$ and $(b_i)$ indexed by the same set (each $a_i$ is the product of some $x_{j,\tau}$'s and $b_i$ is the product of the corresponding $y_{j,\tau}$'s) satisfying
$$\preskip.0em \postskip.3em\ndsp 
\sum_{i=1}^r a_i \sigma(b_i)=
\formule{1\;\;\hbox{if}\;\sigma=\Id\\
         0\;\;\hbox{otherwise}.} 
$$


\exer{exoGalExemples} 
\emph {1.}
As $G$ acts transitively over $\lrb {1..n}$, we have $(\gk^n)^G = \gk$. In addition, $G$ being of cardinality $n$, a permutation $\sigma \in G\setminus\so{\Id}$ has no fixed point. We deduce that $\sum_{\sigma \in G} e_i\sigma(e_i) = 0$ if $\sigma \in G \setminus \{\Id\}$, and $1$ otherwise.\\
By taking $x_i = y_i = e_i$, the conditions of Lemma~\ref{lemArtin} are satisfied and $(\gk, \gk^n, G)$ is a \aGz.

The map $G \to \lrb{1..n}$, $\sigma \mapsto \sigma(1)$, is a bijection. The action of $G$ on $\lrb {1..n}$ is \ncrt \isoc to the action of $G$ on itself by translations. If $n$ is fixed, we can take for $G$ the group generated by an $n$-cycle.

\emph {2.}
We have $\Stp_{\rS_4}(\gB) = \gen {(1,2), (3,4)}$ and $H = \Stp_G(\gB) = \{\Id\}$;
so $(\gk^4)^H = \gk^4$.

\emph {3.}
The first item is immediate.
Suppose~$\gB = \gA[X]^H$ and let $a \in \gA$.\\
Then $aX \in \gB$, so $aX$ is invariant under $H$, \cad $a$ is invariant under $H$. \\
Recap:~$\gA = \Ae H$ so $H = \{\Id\}$ then~$\gA[X] = X\gA[X] + \gk$, \cad~$\gA = \gk$ and $G = \{\Id\}$. Besides this very special case,~$\gB$ is not of the form~$\gA[X]^H$.


\exer{lemPaquesFerrero} 
Assume \spdg that~$\gB$ and~$\gC$ are free of rank $n\in\NN$: it indeed suffices to check the conclusion after \lon at \eco and we have the local structure \tho for \mptfs readily available.
If $n=0$  then~$\gk$ is trivial, we can therefore assume that $1\leq n$. Consider a basis $\cC=(c_1,\ldots,c_n)$ of~$\gC$ and a basis $\cB=(b_1,\ldots,b_n)$ of~$\gB$ (over~$\gk$), and  write the matrix $B\in\Mn(\gk)$ of $\cB$ over $\cC$. The fact that the $b_i$'s form a basis implies that $B$ is injective, \cad $\delta=\det B$ is \ndz
(\thref{prop inj surj det}). Moreover, $\delta\gC\subseteq\gB$.
\\
Let us compare $\Tr\iBk(x)$ and $\Tr\iCk(x)$ for some $x\in\gB$. Consider~$\gk'=\gk[1/\delta]\supseteq\gk$. The two~$\gk'$-\algs obtained by \edsz,~$\gB[1/\delta]$ and~$\gC[1/\delta]$, are the same, and the trace is well-behaved under \edsz, so $\Tr\iBk(x)$ and $\Tr\iCk(x)$ are equal because they are equal in~$\gk'$.
But then

\snic{\disc \gB\sur\gk=\disc\iBk(\cB)=\disc\iCk(\cB)=
\delta^{2}\disc\iCk(c_1,\ldots,c_n).}

Finally, since $\disc \gB\sur\gk$ is \ivz, so is $\delta$ and~$\gB=\gC$.


\exer{exoCorGalste}
First of all, note that since~$\gk$ is connected, all the \mptfs over~$\gk$ are of constant rank. Also recall that the Galois correspondence is already established when~$\gk$ is a \cdiz.
\\
We must show that if~$\gk\subseteq\gB\subseteq\gA$ with~$\gB$ \stez, then

\snic{ \gB=\gC\eqdefi\Fix\big(\Stp(\gB)\big).}

By Lemma~\ref{lemPaquesFerrero}, it suffices to show that~$\gB$ and~$\gC$ have the same rank. In \clama we conclude by noting that after \eds to any field,~$\gB$ and~$\gC$ have the same rank since the Galois correspondence is established for fields.
\\
Via a dynamic rereading of this classical argument we obtain a constructive proof.
This is linked to the formal \nst (\thref{thNSTsurZ}).


\exer{exoIdeauxGlobInvariants} 
Let $(x_i)$, $(y_i)$ be two \syss of \elts of~$\gB$ as in Lemma~\ref{lemArtin}.

\emph {1.}
We know that for $x \in \gB$, $x = \sum_i \Tr_G(xy_i) x_i$. If $x \in \fb$, then $xy_i \in \fb$, and as~$\fb$ is globally invariant, $\Tr_G(xy_i) \in \fb$. \\
Recap~: $\fb$ is generated by the invariant \elts $\Tr_G(xy_i)$ for $x \in \fb$.

\emph {2.}
Let $\fa$ be an \id of~$\gA$; it is clear that $\fa\gB$ is globally invariant.\\
We must show that $\fa\gB \cap \gA = \fa$. This comes from the fact that $\gA$ is a direct summand 
in~$\gB$ (as an~\Amoz). Indeed, let~$\gB = \gA \oplus E$, so $\fa\gB = \fa \oplus \fa E$. If $x \in \fa\gB \cap \gA$, we write $x = y+z$ with $y \in \fa$ and $z \in \fa E \subseteq E$; we then have $x$, $y \in \gA$, so $z \in \gA$, and as $z \in E$, $z = 0$. Consequently, $x = y \in \fa$.

Conversely, if $\fb \subseteq \gB$ is globally invariant, we must show that $(\fb\cap \gA)\gB = \fb$; but this is what has been shown in the previous question.


\prob{exoBuildingFrobAlgebra} 
\Gnltz, the \lin form $\delta_g$ is passed on to the quotient modulo the \id $\fa_g$ that it defines. 
In addition, if $\delta_g(\ov u\,\ov v) = 0$ over~$\gA = \kuX\sur{\fa_g}$ for every $\ov v$, then $\delta_g(u v) = 0$ for every $v$, so $u \in \fa_g$, \cad $\ov u = 0$. Therefore $\Ann_\gA(\delta_g) = 0$.

For $i \in \lrbn$, let $\delta_i^m = \delta_{X_i^m}$ (\coo over $X_i^m$). \\
In particular, $\delta_i(f) = {\partial f \over \partial X_i}(0)$, and we define $\delta_0 : \kuX \to \gk$ by $\delta_0(f) = f(0)$.

\emph {1.}
Easy computation.

\emph {2.}
We verify that $f*g = 0$ \ssi $f_m*g = 0$ for every \hmg component 
 $f_m$ of~$f$. In other words the \id $\fa_g$ is \hmg (this is always the case if~$g$ is \hmgz).\\
It is also clear that for $i \ne j$, $X_iX_j * g = 0$, and for $|\alpha| > d$, $X^\alpha * g = 0$.  
If $f = \sum_i a_iX_i^m + \cdots$ is \hmg of degree $m \le d$, we have $f * g = \sum_i a_i X_i^{-(d-m)}$. 
 
If $m < d$, we therefore have $f*g = 0$ \ssi $a_i = 0,\;\Tt i$, \cad {if $f \in \gen {X_iX_j, i \ne j}$}.
 
If $m = d$, we have $f*g = 0$ \ssi $\sum_i a_i = 0$, \cad {if $f \in \gen {X_iX_j, i \ne j} + \gen {X_i^d - X_1^d, i \in \lrb {2..n}}$}, because $\sum_i a_iX_i^d = \sum_i a_i(X_i^d - X_1^d)$. 
\\
Recap: we have obtained a \sgr of $\fa_g$ consisting of $n(n-1) \over 2$ \pogs of degree $2$ and of $n-1$ \pogs of degree $d$
$$\preskip.4em \postskip.4em 
\fa_g = \gen {X_iX_j, i < j} + \gen {X_i^d - X_1^d, i \in \lrb{2..n}}
. 
$$
Let~$\gA = \kuX\sur\fa_g = \gk[x_1, \ldots, x_n]$. Then
$$\preskip.2em \postskip.4em 
1,\quad x_1,\ \ldots,\ x_n,\quad x_1^2,\ \ldots,\ x_n^2, 
\quad\ldots\quad
x_1^{d-1},\ \ldots,\ x_n^{d-1}, \quad x_1^d 
$$
is a~$\gk$-basis of~$\gA$ of cardinality $(d-1)n+2$. The~$\gk$-dual basis
of $\Asta$ is
$$\preskip.4em \postskip.4em 
\delta_0,\quad \delta_1,\ \ldots,\ \delta_n,\quad \delta_1^2,\ \ldots,\ \delta_n^2, 
\quad\ldots\quad
\delta_1^{d-1},\ \ldots,\ \delta_n^{d-1}, \quad \delta_g 
$$
and we have
$$\preskip-.2em \postskip.4em 
x_i^m \centerdot\delta_g = \delta_i^{d-m} \hbox { for } m \in \lrb{1..d-1},
\quad 
x_i^d \centerdot\delta_g = \delta_0
. 
$$
Therefore $\Asta = \gA\centerdot\delta_g$ and $\delta_g$ is dualizing.

\emph {3.}
If we take $e_i$ strictly greater than the exponent of $X_i$ in the set of \moms of $g$, we have $X_i^{e_i} * g = 0$.

\emph {4.}
Let $f \in \kuX$. We have seen that $f\centerdot \delta_g = 0$ \ssi $\partial_f(g) = 0$. 
\\
So the~$\gk$-\lin map~$\kuX \to \kuX$, $f \mapsto \partial_f(g)$, 
passes to the quotient modulo $\fb$ to define a~$\gk$-\lin map $\varphi : \kuX\sur\fb \to \kuX$.

\emph {5.}
The~$\gk$-\lin map~$\gA \to \Asta$, $f \mapsto f\centerdot \delta_g$,
is injective and as~$\gA$ and~$\Asta$ are \hbox{$\gk$-\evcsz} of the same finite dimension, it is an \isoz.

\prob{exoTh90HilbertAdditif}
\emph{1.}  As if by magic we let $\theta(x) =
\sum_{i=0}^{n-1} \sigma^i(z) c_i(x)$ (thanks to Hilbert\ihiz). We will check that

\snic {
\sigma\big(\theta(x)\big) = \theta(x) + \Tr_G(x)z - x
\quad \hbox {or} \quad
x = (\Id_\gA - \sigma)\big(\theta(x)\big) + \Tr_G(x)z.
}

So,
$(\Id_\gA - \sigma) \circ \theta$ and $x \mapsto
\Tr_G(x)z$ are two \orts \prrs with sum~$\Id_\gA$, \hbox{hence $\gA = \Im(\Id_\gA - \sigma) \oplus \gk z$}. For the \vfnz, write $c_i$ for $c_i(x)$ and~$y = \theta(x)$. We have $\sigma(c_i) = c_{i+1} - x$, $c_n = \tr_G(x)$ and

\snic {
\begin{array} {rcl}
\sigma(y) &=& \sum_{i=0}^{n-1} (c_{i+1} - x) \sigma^{i+1}(z) =
\sum_{i=0}^{n-1} c_{i+1} \sigma^{i+1}(z) -
\sum_{i=0}^{n-1} x\sigma^{i+1}(z)   \\[1mm]
&=& (y + \Tr_G(x)z) - x\Tr_G(z) = y + \Tr_G(x)z - x .
\end{array}
}

Since $\Tr_G(z)=1$, $z$ is a basis of~$\gk z$
(if $az=0$, \hbox{then $0=\Tr_G(az)=a$}), so $\Im(\Id_\gA - \sigma)$ is indeed \stl of rank $n-1$.

 \emph{2.} It is clear that $\Im(\Id_\gA - \sigma) \subseteq \Ker\Tr_G$. The other inclusion results from the previous item.

 \emph{3.}
The reader can verify this by letting $y = \sum_\tau c_\tau \tau(z)$. There is a link with  question \emph{1}: for fixed $x$ with $\Tr_G(x) = 0$, the family $\big(c_i(x)\big)$ is an additive $1$-cocycle under the condition that $\lrb{0..n-1}$ and $G$ are identified via $i \leftrightarrow \sigma^i$.

\emph {4.}
The \elt $-1$ has null trace, hence the existence of $y\in\gA$ such that
$-1=y-\sigma(y)$. We then have, for every $i \in \ZZ$, $\sigma^i(y) = y+i$%
, and $\sigma^j(y) - \sigma^i(y) = j-i$ is \iv for $i \not\equiv j
\bmod p$.
\\ 
Let $y_i = \sigma^i(y)$, $(i \in \lrb{0..p-1})$. The Vandermonde matrix of $(y_0, y_1, \ldots, y_{p-1})$ is \iv and consequently $(1, y, \ldots, y^{p-1})$ is \hbox{a~$\gk$-basis} of~$\gA$. Let~$\lambda = y^p
- y$. Then $\lambda\in\gk$ since

\snic {
\sigma(\lambda) = \sigma(y)^p - \sigma(y) = (y+1)^p - (y+1) = y^p  - y = \lambda.
}

The \polcar of $y$ is $(Y-y_0) (Y-y_1)\cdots (Y-y_{p-1})$ and this \pol is equal to $f(Y) = Y^p - Y - \lambda$ (because the $y_i$'s are roots of~$f$ and $y_i - y_j$ is \iv for $i \ne j$).

\emph {5.}
Let~$\gk$ be a \ri with $p=_\gk0$. Fix $\lambda \in \gk$ and let $\gA = \aqo{\gk[Y]}{f} = \gk[y]$, where $f(Y) = Y^p-Y-\lambda$. Then, $y+1$ is a root of~$f$, and we can define $\sigma \in \Aut(\gA/\gk)$ by $\sigma(y) = y+1$. The \elt $\sigma$ is of order $p$ and the reader will check that $(\gk,\gA, \gen {\sigma})$ is a \aGz.


\prob{exoHeitmannGaloisExemple} 
\emph {1.}
Consider the \idz~$\gen {x-\sigma(x), y-\sigma(y)} \eqdefi \gen {2x,2y}$.
Since $2$ is \ivz, it is the \idz~$\gen {x,y}$, and it contains $1$ because $x^2 + y^2=1$.
Thus,~$\gen {\sigma}$ is separating.

\emph {2.}
For all $f\in \gB$, we have $f = (xf)x + (yf)y$. If~$f$ is odd \cad if $\sigma(f) = -f$, we~have~$xf$,~$yf \in \gA$, so $f \in \gA x + \gA y$ and $E = \sotq {f\in\gB}{\sigma(f) = -f}$.  
The \egtz~$\gB = \gA \oplus E$ stems from the \egt $f = \fraC {f+\sigma(f)}2 + \fraC {f-\sigma(f)}2$ for $f \in \gB$.\\  
\emph{Other \demz.} We know that there exists a $b_0 \in \gB$ of trace $1$ and that the kernel of the \lin form $\gB \to \gA$ defined by $b \mapsto \Tr(b_0b)$ is a \supl subspace of~$\gA$
in~$\gB$. Here we can take $b_0 = 1/2$, again we find $E$ as a \supl subspace.

\emph {3.}
This is a matter of finding $y_1$, $y_2$, $y_3 \in \gB$ such that $\sum_{i=1}^3 x_i\tau(y_i)=1$ for $\tau = \Id$ and $0$ otherwise. We notice that
$$\preskip.2em \postskip.9em
1\cdot 1 + x\cdot x + y\cdot y = 2, \quad\hbox{and}\quad 
1\cdot \sigma(1) + x\cdot \sigma(x) + y\cdot \sigma(y) = 0,
$$
hence we obtain a solution when taking $y_i = x_i/2$. Letting \smash{$X = \Cmatrix{.3em} {1 & x & y\cr {1\over 2} &-{x\over 2} &-{y\over 2}\cr}$}, we have~$X \tra {X} = \I_2$ and $\tra{X} X = P$ with $P = 
\Cmatrix{.4em} {1 & 0 &0\cr 0 &x^2 &xy\cr 0 & xy &y^2 \cr}$. The matrix $P$ is a \prr of rank $2$ whose image is \isoc to the \Amoz~$\gB$.\\
Note: we deduce that $E$ is \isoc to the image of the projector
\smashbot{$\Cmatrix{.3em} {x^2 &xy\cr  xy &y^2}$} and that~$\gB\te_\gA E$ is \isoc to~$\gB$ as a \Bmoz.

\emph {4.}
Easy.

\emph {5.}
The \iso $E^n \simeq \Ae n$ proves that $E$ is a \mrcz~$1$.
Applying $\Al{n}$ we obtain $E^{n\te} \simeq \gA$.\\
Note: for more details see Section~\ref{sec ptf loc lib}, the \dem of Proposition~\ref{prop puissance ext}, \Egtz~(\iref{eqVik}) on \paref{eqVik} and \Egtz~(\iref{eqfactPicStab2}) on \paref{eqfactPicStab2}.

\emph {6.}
The \egt $1 = x^2 + y^2$ implies $\fa^2 = \gen {x^2 y^2, x^3y, x^4} = x^2 \gen {y^2, xy, x^2} = x^2\gA
$,
and $\fa\gB = xy\gB + x^2\gB = x(y\gB + x\gB) = x\gB$. In~$\gB$, $\fa = x(y\gA + x\gA) = xE$. Therefore if~$x$ is \ndzz, $\fa \simeq_\gA E$ via the multiplication by $x$.

\emph {7a.} We have~$\gk[x]\simeq\kX$ and~$\gA=\gk[x^2,xy,y^2]$.
We consider~$\gB$ as a free $\gk[x]$-module of rank $2$, with basis $(1,y)$, and we let $\rN : \gB \to \gk[x]$ be the norm. 
\\
For $a$, $b \in \gk[x]$ we obtain

\snic {
\rN(a+by) = (a+by)(a-by) = a^2 + (x^2-1)b^2
.}

As $\rN(x)=x^2$, $x$ is \ndz (Lemma~\ref{lemIRAdu} item~\emph{2}).
Moreover, $a + by\in\Bti$  \ssi $a^2 + (x^2-1)b^2 \in \gk\eti$.  
Suppose $b$ has formal degree $m\geq 0$ and~$a$ has formal degree $n\geq 0$. Then, $(x^2-1)b^2 = \beta^2 x^{2m+2} + \dots$ and $a^2 = \alpha^2 x^{2n} + \dots$ 
Since $a^2 +(x^2-1)b^2\in\gk\eti$, we obtain
\begin{itemize}
\item   if $n>m+1$, $\alpha^2=0$ so $\alpha=0$ and $a$ can be rewritten with formal degree $<n$,
\item if $n<m+1$, $\beta^2=0$ so $\beta=0$ and  
\begin{itemize}
\item if $m=0$, $b=0$ and $a=\alpha\in\gk\eti$, or
\item  if $m>0$, $b$ can be rewritten with formal degree $<m$,
\end{itemize}
\item if $n=m+1$ (which implies $n>0$),  $\alpha^2 + \beta^2= 0$ so $\alpha=\beta=0$ and $a$ can be rewritten with formal degree  $<n$.
\end{itemize}
We conclude by \recu on $m+n$ that if $a + by\in\Bti$, then $b = 0$ and $a \in \gk\eti$.

We notice that if $-1= i^2$ in~$\gk$, then $(x+iy)(x-iy) = 1$ and we obtain an \inv $x + iy$ which is not a constant.

\emph {7b.}
Let us show that $\fa$ is not principal. 
As $\fa \simeq_\gA E$, it will follow that $E$ is not a free \Amoz, and $\gB$ is not free either, because otherwise $E$ would be \stl of rank $1$, therefore free.\\
Suppose $\fa = a\gA$ with $a \in \gA$. By extending to~$\gB$, we obtain $\fa\gB = a\gB$. But we have seen that $\fa\gB = x\gB$, and since $x$ is \ndzz, $x = ua$ with $u \in \gB\eti = \gk\eti$. This would imply~$x \in \gA$, which is not the case because $\gk$ is nontrivial.

\emph {8.}
We reuse the preceding \dem   to show that $\fa$ is not principal, but here~$\gB\eti$ no longer consists of only constants, for example the (continuous) function $(x,y) \mapsto x^2 + 1$ is \ivz. From the point where $x = ua$ \hbox{and $u \in \gB\eti$}, we reason as follows. Since $u$ is an \iv \elt of $\gB$, its absolute value is bounded below by an \elt $>0$, and $u$ 
         is everywhere $>0$, or everywhere $<0$. 
As $x$ is odd and $a$ even,  $a$ and $x$ are identically zero; a contradiction.



\Biblio

A \cov study of \stfes (not \ncrt commutative) associative \algs over a \cdi 
can be found in \cite[Richman]{ri82} and in \cite[Chapitre~IX]{MRR}.

Proposition~\ref{propIdemMini} is found in \cite{MRR} which introduces the terminology of \emph{\sply factorial field}. 
See also \cite[Richman]{ri81}.

Lemma~\ref{lemSqfDec} for squarefree \fcn over a perfect \cdi admits a subtle \gnn in the form of an \gui{\algo for \spl \fcnz} over an arbitrary \cdiz; see \cite[th.~IV.6.3, p.~162]{MRR} and~\cite[Lecerf]{Lecerf-Factsep}.

The notions of a \aG and of a \spb \alg were introduced by Auslander \& Goldman in \cite[1960]{AG}. The core of the theory of \aGs is found in Chase, Harrison \& Rosenberg's paper \cite[1968]{CHR}. 
A book that presents this theory is \cite{DI}. Almost every argument in \cite{CHR} is already of an \elr and \cov nature.

The result given in Exercise~\ref{lemPaquesFerrero} is due to Ferrero and Paques in \cite{FP}.

\Pbmz~\ref{exoBuildingFrobAlgebra} is inspired by Chapter 21 (Duality, Canonical Modules, and Gorenstein Rings) of \cite{Eis} and in particular by Exercises 21.6 \hbox{and 21.7}.

\newpage \thispagestyle{CMcadreseul}
\incrementeexosetprob


\chapter[The dynamic method]{The dynamic method
}
\label{ChapGalois}
\vspace{-1.2cm}
{\LARGE \bf
\hspace*{1cm} \nst\\[.3cm] 
\hspace*{1cm} Splitting field\\[.3cm]
\hspace*{1cm} Galois theory}
\perso{compil\'e le \today}

\vspace{1.2cm}
\minitoc

\Intro

The first section of this chapter gives \gnl \cov versions of the \nst for a \syp over a \cdi (we will be able to compare \Thos \rref{thNstNoe}, \rref{thNstClassCof1} and \rref{thNstClassCof2},  to \Thos \rref{thNstfaibleClass}
and \rref{thNstClass}).
We \egmt give a simultaneous \iNoe positioning \tho (\thref{thNoetSimult}).

This is a significant example of a reformulation of a result from \clama \emph{in a more \gnl framework}: \clama admits that every field has an \agq closure. 
This means it does not have to deal  with the problem of the exact meaning of Hilbert's Nullstellensatz when such an algebraic closure is not available. But the question does get asked, and we can offer a perfectly reasonable answer: the algebraic closure is not really necessary. Rather than looking for the zeros of a polynomial system in an algebraic closure, we can look for them in finite algebras over the field given at the start.

\smallskip  
We then tackle another \pbz: that of \cot interpreting the classical discourse on the \agq closure of a field. 
The problem might seem to largely involve the use of Zorn's lemma, which is necessary for the construction of the global object.
 Actually, a more delicate \pb arises well beforehand, at the moment the splitting field of an individual \pol is constructed. 

The \tho from \clama stating that every \spl \pol of $\KT$ has a \stf \cdr over $\gK$ (in which case the Galois theory applies), is only valid from a \cof point of view under hypotheses regarding the possibility of factorizing the \spl \pols (cf. \cite{MRR} and in this work \thref{thResolUniv} on the one hand and Corollary~\ref{propIdemMini} on the other).
Our goal here is to give a \cov Galois theory for an arbitrary \spl \pol in the absence of such hypotheses.

The counterpart is that we must not consider the \cdr of a \pol as a usual \gui{static} object, but as a \gui{dynamic} object. This phenomenon is inevitable, because we must manage the ambiguity that results from the impossibility of knowing the Galois group of a \pol by an infallible method. 
Moreover, the disorientation produced by this shift to a dynamic perspective is but one example of the \gnle lazy \evn method: \emph{nothing comes of over-exhausting ourselves to know the whole truth when a partial truth is sufficient for the stakes of the ongoing computation.}

\smallskip 
In Section~\ref{subsecDyna}, we give a heuristic approach to the dynamic method, which forms a cornerstone of the new methods in \cov \algz.

Section~\ref{secBoole}, dedicated to  \agBsz, is a short introduction to the \pbs that will have to be dealt with in the context of a \adu over a \cdi when it is not connected.

Section~\ref{secadu} continues the theory of \adu already started in Section~\ref{sec0adu}. Without assuming that the \pol is \splz, the \adu has several interesting \prts that are preserved upon passage to a \gui{Galois quotient.} When summarizing these \prts we have been brought to introduce the notion of a \emph{\apGz}.

Section~\ref{subsecCDR} gives a \cov and dynamic approach to the \cdr of a \pol over a \cdiz, without a \spt hypothesis regarding the \polz.

The dynamic Galois theory of a \spl \pol over a \cdi is developed in Section~\ref{secThGB}.

\medskip
The current chapter can be read \imdt after Sections~\ref{secGaloisElr} and~\ref{sec2GaloisElr}, bypassing Chapters~\ref{chap mpf} and~\ref{chap ptf0}, if we restrict the \adu to the \cdis case (which would in fact simplify some of the \demsz).
However, it seemed natural to us to develop the material with respect to the \adu in a more \gnl framework, which requires the notion of a \mrc over an arbitrary commutative \riz.

\pagebreak
\pagestyle{CMheadings}

\section{The \nst without \agq closure}
\label{secNstSCA}

In this chapter, which is dedicated to the question \gui{how can we constructively recover the results from classical mathematics that are based on the existence of an algebraic closure, even when it is missing?,}
 it seemed logical to have a new look at 
the \nst and the \Noe position (\thref{thNstfaibleClass}) in this new framework.

\subsec{The case of an infinite basis field}

We claim that \thref{thNstfaibleClass} can be copied virtually word-for-word, by simply deleting the reference to an \cac that contains~$\gK$.

We no longer \ncrt see the zeros of the \syp considered in finite extensions of the \cdi $\gK$, but we construct \stfes nonzero \Klgs (\cad that are finite dimensional \Kevsz) which account for these zeros; in \clama the zeros are in the quotient fields of these \Klgsz, 
and such quotient fields are easily seen to exist by applying
 \TEM since it suffices to consider a strict \id that is of maximal dimension as a \Kevz. 

\begin{theorem}\label{thNstfaibleClassSCA}\emph{(Weak \nst and \iNoe position, 2)}\\
Let $\gK$ be an infinite \cdi and {\mathrigid 2mu $(\lfs)=(\uf)$ be a \syp in the \alg $\KuX=\KXn$ ($n\geq 1$). Let \smash{$\ff=\gen{\uf}_\KuX$} and $\gA=\KuX\sur\ff$.}

\emph{Weak \nstz.}\\  
Either $\gA=0$, or there exists a nonzero quotient of $\gA$ which is a \stfe \Klgz.

\emph{\iNoe postion.}\\
More \prmtz, we have a well-defined integer $r\in\lrb{-1..n}$ with the following \prtsz.
\begin{enumerate}

\item 
Either $r=-1$ and $\gA=0$  (\cad $\gen{\uf}=\gen{1}$). 
In this case, the \sys  $(\uf)$ does not admit any zero in any nontrivial \Klgz.

\item 
Or $r=0$, and $\gA$ is a nonzero \stfe \Klg (in particular, the natural \homo $\gK\to\gA$ is injective). 

\item 
Or $r\geq 1$, and there exists a $\gK$-\lin \cdv (the new variables are denoted $\Yn$) satisfying the following \prtsz.
\begin{itemize}
\item [$\bullet$] We have $\ff\,\cap\,\gK[\Yr]=0$. In other words, the \pol \ri $\gK[\Yr]$ can be identified with a sub\ri of $\gA$.
\item [$\bullet$] Each $Y_j$ for $j\in\lrb{r+1..n}$ is integral over $\gK[\Yr]$ modulo $\ff$ and the \ri $\gA$ is a $\gK[\Yr]$-\mpfz.
\item [$\bullet$] There exists an integer $N$ such that for each $(\alpha_1,\ldots,\alpha_r)\in\gK^r$, the quotient \alg $\aqo\gA{Y_1-\alpha_1,\ldots,Y_r-\alpha_r}$ is a nonzero \Kev of finite dimension $\leq N$.
\item [$\bullet$] We have \itfs $\ff_j\subseteq\gK[Y_1,\ldots,Y_j]$ $(j\in\lrb{r..n})$ with the following inclusions and \egtsz.
\[\preskip.2em \postskip.4em 
\begin{array}{ll} 
 \gen{0}=\ff_r\subseteq\ff_{r+1} \subseteq \ldots\subseteq\ff_{n-1}\subseteq\ff_n=\ff   \\[1mm] 
\ff_j\subseteq\ff_\ell\,\cap\,\gK[Y_1,\ldots,Y_j] &(j<\ell,\ j,\ell \in\lrb{r..n})   \\[1mm] 
\rD(\ff_j) = \rD\left(\ff_\ell\,\cap\,\gK[Y_1,\ldots,Y_j]\right) &(j<\ell,\ j,\ell  \in\lrb{r..n})   
 \end{array}
\]

\end{itemize}
\end{enumerate}
\end{theorem}

%
\begin{proof}
We essentially reason as in the \dem of \thref{thNstfaibleClass}.
To simplify we keep the same variable names at  each step of the construction.
Let $\ff_n=\ff$.
\begin{itemize}
\item Either $\ff=0$, and $r=n$ in item \emph{3.} 
\item Or there is a nonzero \pol among the $f_i$'s, we make a \lin \cdv that renders it \mon in the last variable,  and we compute the resultant \id $\fRes_{X_n}(\ff_n)=\ff_{n-1}\subseteq\gK[X_1,\ldots,X_{n-1}]\cap \ff_n$.
\\
Since $\ff_n\cap \gK[X_1,\ldots,X_{n-1}]$ and $\ff_{n-1}$ have the same nilradical, they are simultaneously null.
\item  If $\ff_{n-1}=0$, item~\emph{3} or \emph{2} is satisfied with $r=n-1$.  
\item  Otherwise, we iterate the process. 
\item  If the process halts with $\ff_r=0$, $r\geq0$,  item~\emph{3} or \emph{2} is satisfied with this value of~$r$.  
\item  Otherwise, $\ff_0=\gen{1}$ and the computation has allowed us to construct $1$ as an \elt of $\ff$.
\end{itemize}
There are two things left for us to verify. 
\\
First of all, that $\gA$ is a $\gK[\Yr]$-\mpfz. It is clear that it is a \tf module, the fact that it is \pf is therefore given by \thref{propAlgFinPresfin}.  
\\
Then, that when we specialize the $Y_i$'s ($i\in\lrbr$) in some $\alpha_i\in\gK$, \hbox{the \Kevz} obtained is \pf (so finite dimensional) and nonzero. 
\Thref{factSDIRKlg} on changing the base \ri gives us the fact that, after specialization, the \alg remains a \mpfz, so that the obtained \Kev is indeed finite dimensional.
We must show that it is nonzero. 
However, we notice that, by assuming the \cdvs already made at the start, all the computations done in $\gK[\Yn]$ specialize, \cad remain unchanged, if we replace the \idtrs $Y_1$, \dots, $Y_r$ by the scalars $\alpha_1$, \dots, $\alpha_r$. The conclusion $\ff\;\cap\;\gK[\Yr]=0$ is replaced by the same result specialized in the $\alpha_i$'s, \cad \prmt what we wanted. 
\\
We can obtain the same conclusion in the more scholarly form below. This specialization is a change of the base \ri $\gK[\Yr]\to\gK$. 
Apply item~\emph{1c} of the \gnl \eli lemma~\ref{LemElimAff} with 

\snic{\gk=\gK[\Yr]$, $\gC=\gA$ and $\gk'=\gK.}

The \eli \id and the resultant \id in $\gk$ are null, therefore after \eds the resultant \id remains null in~$\gK$. 
\\
Therefore, the same thing holds for the \eli \idz, and the natural \homo $\gK\to\aqo\gA{Y_1-\alpha_1,\ldots,Y_r-\alpha_r}$ is injective.

Let us end by explaining why the integer $r$ is well-defined. First of all the case $r=-1$ is the only case where $\gA=0$, then for $r\geq 0$, it is possible to show that $r$ is the maximum number of \elts \agqt independent over~$\gK$ in~$\gA$ (see the proofs \thrfs{propDimKXY}{thDKAG}). 
\end{proof}

\rems
1)
We have used resultant \ids $\fRes(\fb)$ (\thref{thElimAff}) instead of \ids $\fR(g_1,\ldots,g_s)$ (with $g_1$ \mon and $\gen{g_1,\ldots,g_s}=\fb$), introduced in Lemma~\ref{lemElimPlusieurs}. However, Lemma~\ref{lemElimParametre}
 shows that the latter \ids would do just as well.\\
2) 
For any arbitrary \homo $\gK[\Yr]\to\gB$, when $\gB$ is a reduced \Klgz, the last argument in the \dem of the \tho works, which tells us that $\gB\subseteq\gB\otimes_{\gK[\Yr]}\gA$.\\
3)
The last item of \emph{3} recalls the workings of the \dem by \recu which constructs the \itfs $\ff_j$ to reach the \iNoe position. This also gives a certain description of the \gui{zeros} of the \syp (more delicate than in the case where we have an \cac $\gL$ that contains $\gK$, and where we describe the zeros with \coos in $\gL$, as in \thref{thNstfaibleClass}).  
\eoe

\entrenous{L'argument final pour la sp\'ecialisation non nulle in la \dem of the \tho pourrait sembler un peu subtil if on ne le donnait que sous la forme savante.
Sous la forme concr\`ete: il is clair que tous les calculs se sp\'ecialisent,
il is en fait tr\`es simple. 
}

\smallskip  It remains for us to lift the restriction introduced by the consideration of an infinite \cdi $\gK$. For this we need a \cdv lemma that is a bit more \gnlz, using Nagata's trick.

\subsec{Changing variables}

\vspace{-.2em}

\begin{definition}\label{defiCDV}
We call a \ix{change of variables} in the \pol \ri $\kuX=\kXn$ an \auto $\theta$ of this \klgz. If each~$\theta(X_i)$ is denoted by $Y_i$, the $Y_i$'s are called \emph{the new variables}. Each~$Y_i$ is a \pol in the $X_j$'s, and each $X_i$ is a \pol in~the~$Y_j$'s. 
\end{definition}

The most frequently used are the \gui{\lin changes of variables,} in which we include, despite their name, the translations and all the affine transformations.

\smallskip 
\comm A non\lin \cdvz, like for instance 
$$\preskip.4em \postskip.4em
{(X,Y)\mapsto(X+Y^2,Y),}
$$
does not respect the \gmt in the intuitive sense. For example a line is transformed into a parabola; the \agq \gmt of the affine plane is not an extension of the affine \gmtz, it directly contradicts it!
It is only in the context of projective spaces that we find what we expect: the \autos of the projective plane, from the \agq \gmt point of view, are \ncrt \linz, and the notion of a \gui{(straight) line} reclaims its rights.
\eoe

\bni {\bf Pseudo\polusz}~
\rdb\label{polpseudunit}

Let $\gk$ be a connected \riz. A \pol in $\kT$ is said to be \ixc{pseudomonic}{polynomial} (in the variable $T$) if it is of the form $\som_{i=0}^pa_kT^{k}$ with $a_p$ \ivz.

In \gnlz, without assuming that $\gk$ is connected, a \pol in $\kT$ is said to be \emph{pseudomonic} (in the variable $T$) if there exists a \sfio $(e_0,\ldots ,e_r)$ such that, for each $j$, when taking $\gk[1/e_j]=\gk_j$, the \pol is expressible in the form $\som_{k=0}^{j}a_{k,j}T^{k}$ with $a_{j,j}$ \iv in~$\gk_j$.

A \pol in  $\kXn=\kuX$ is said to be
\emph{pseudomonic in the variable} $X_n$ if it is pseudomonic as an \elt of $\gk[X_1,\ldots ,X_{n-1}][X_n]$.\label{polpseudounitaire}%
\index{polynomial!pseudomonic ---}

NB: See also the notion of a \lot \polu in Exercise~\ref{exoPolLocUnitaire}.

Recall that a \pol of $\kXn$ is said to be \emph{primitive} when its \coes generate the \id $\gen{1}$.
Also recall that if $\gk$ is reduced, we have the \egt $\kXn\eti=\gk\eti$ (Lemma~\ref{lemGaussJoyal}).

\begin{fact}\label{factZedPrimFin}
Let $\gK$ be a reduced \zed \ri and $P\in \gK[T]$. \Propeq
\begin{enumerate}
  \item [--] The \pol $P$ is \ndzz.
  \item [--] The \pol $P$ is primitive.
  \item [--] The \pol $P$  is pseudomonic.
  \item [--] The quotient \alg $\aqo{\gK[T]}{P}$ is finite over $\gK$.
\end{enumerate}
\end{fact}
\begin{proof}
The \eqvcs are clear in the \cdis case. To obtain the \gnl result we can apply the \elgbm of  \zedr \ris (\paref{MethodeZedRed}).
\end{proof}
%

\subsubsection*{A simple and efficient lemma} 

\begin{lemma}
\label{lemNoether} \emph{(\Cdvs lemma \`a la Nagata)}\\
Let $\gK$ be a \zedr \ri and  $g\in\KuX=\KXn$ be a \ndz \eltz.
\begin{enumerate}
\item There exists a \cdv such that, by calling the new variables $Y_1$, $\ldots $, $Y_{n}$, the \pol $g$ becomes pseudomonic in $Y_n$.
Consequently the \Klg $\aqo{\KuX}{g}$ is finite over $\gK[Y_1,\ldots ,Y_{n-1}]$.
\item When $\gK$ is an infinite \cdiz, we can take a \lin \cdvsz.
\item The result also applies to a finite family of \ndz \pols of $\KuX$ (they can be made simultaneously pseudo\mons by the same \cdvz).
\end{enumerate}
\end{lemma}
\begin{proof}
For the case of an infinite \cdi see Lemma~\ref{lemCDV}.
\\
In the \gnl case we can assume that $\gK$ is a \cdi and we make a \cdv \gui{\`a la Nagata.} For example with three variables, if the \pol $g$ is of degree $< d$ in each of the variables $X$, $Y$,~$Z$, we make the \cdv $X\mapsto X,\,Y\mapsto Y+X^d,\,Z\mapsto Z+X^{d^2}.$
Then, seen as an \elt of $\gK[Y,Z][X]$, $g$ has become pseudo\mon in~$X$.

 Item \emph{3} is left to the reader.
\end{proof}

\subsec{The \gnl case}

By reasoning as we did for \thref{thNstfaibleClassSCA} and by using the \cdvs of the previous lemma we obtain the  \gnle form of the weak \nst and of the \Noe position in \comaz.

\begin{theorem}\label{thNstNoe}\emph{(Weak \nst and \iNoe position, 3)}\\
With the same hypotheses as in \thref{thNstfaibleClassSCA} but 
by only supposing that the \cdi $\gK$ is nontrivial, we get the same conclusions, except that the \cdv is not \ncrt \linz.
\end{theorem}

\begin{definition}\label{defiDimsyp} 
Consider the case $1\notin\gen{\lfs}$ of the previous \thoz. 
\begin{enumerate}
\item We say that the \cdv (which eventually changes nothing at all) has put the \id $\ff$ in \emph{\Noe position}.
\item The integer $r$ that intervenes in the \Noe positioning is called the \ixc{dimension}{of an affine \vrtz} \emph{of the \sypz}, or of the \vrt defined by the \sypz, or of the quotient \alg $\gA$.
By convention the null \alg is said to be of dimension~$-1$.%
\index{dimension!of a \syp over a \cdiz}\index{dimension!of an \apf over a \cdiz}
\end{enumerate}
\end{definition}
\rems 1) It is clear by the \tho 
that $r=0$ \ssi the quotient \alg is finite nonzero, which implies (Lemma~\ref{lemZrZr1}) that it is a nontrivial \zed \riz.  
\\
Conversely, if $\gA$ is \zed and $\gK$ nontrivial, Lemma~\ref{lemZrZr2} shows that the \ri $\gK[\Yr]$ is \zedz, which implies that $r\leq 0$ (\hbox{if $r>0$,} then an \egt $Y_r^m\big(1+Y_rQ(\Yr)\big)=0$ implies that $\gK$ is trivial).
Therefore there is no conflict with the notion of a \zed \riz.
Let us however note that the null \alg is still a \zed \riz.
 
 2) The link with the \ddk will be made in \thref{thDKAG}.
 
 3) A \gui{non-\noeez} version of the previous \tho for a \zedr \ri $\gK$ is given in Exercise~\ref{exothNst1-zed}.
\eoe

\vspace{-.2em}
\pagebreak

\begin{theorem}
\label{thNoetSimult} \emph{(Simultaneous \iNoe position)}
\\
Let $\ff_1$, \ldots, $\ff_k$ be \itfs of $\KuX=\KXn$. 
\begin{enumerate}\itemsep0pt
\item  There exist integers $r_1, $\ldots, $r_k\in\lrb{-1..n}$ and a \cdv such that, by calling $Y_1$, \ldots, $Y_{n}$ the new variables, we have for each $j\in\lrbk$ the following situation.
\\
If $r_j=-1$, then $\ff_j=\gen{1}$, otherwise
\begin{enumerate}
\item 
$\gK[Y_1,\ldots , Y_{r_j}]\,\cap\, \ff_j=\so{0}$,
\item 
for $\ell>r_j$, $Y_\ell$ is integral modulo $\ff_j$ over $\gK[Y_1,\ldots , Y_{r_j}]$.
\end{enumerate}
When $\gK$ is infinite, we can take a \lin \cdvsz.

\item If $\gen{1}\neq\rD(\ff_1)\supset \rD(\ff_2)\supset \cdots \supset \rD(\ff_k)$ with the strictly increasing dimensions $r_j$, we can insert radicals of \itfs such that the obtained sequence of dimensions is $0$, $1$, \ldots, $n$.  
\end{enumerate}
\end{theorem}
NB: In item \emph{1}, we say that the \cdv (which eventually changes nothing at all) has simultaneously put the \ids $\ff_1,\ldots,\ff_k$ in \Noe position.

\begin{proof}  \emph{1.}
The same \dem as for the previous \tho works considering the fact that a \cdv can simultaneously render a finite number of nonzero \pols \mons in the last variable.
\\
\emph{2.} Let $\gA_i=\gK[X_1,\ldots,X_i]$.
Suppose for example that $\ff_1$ is of dimension~$2$ and $\ff_2$ of dimension $5$. We have to insert \ids of dimensions~$3$ and $4$. Suppose \spdg that the $\ff_i$'s are in \Noe position with respect to $\Xn$.
\\
We have by hypothesis $\gA_2\,\cap\,\ff_1=0$, with \polus 

\snic{h_3\in\gA_2[X_3]\,\cap\,\ff_1$, $h_4\in\gA_2[X_4]\,\cap\,\ff_1$, \ldots, $h_n\in\gA_2[X_n]\,\cap\,\ff_1.}

We then have the following inclusions,

\snac{\fh_1=\ff_2+\gen{h_5,h_4} \supseteq \fh_2=\ff_2+\gen{h_5}
\supseteq \ff_2 \et \rD(\ff_1)\supseteq \rD(\fh_1)\supseteq \rD(\fh_2)\supseteq \rD(\ff_2),}

with $\fh_1$ of dimension $3$ and $\fh_2$ of dimension $4$, both in \Noe position with respect to $(\Xn)$. 
\end{proof}
%

\subsec{The actual \nstz}
\vspace{3pt}
In \Thosz~\ref{thNstfaibleClassSCA} (infinite \cdiz) and~\ref{thNstNoe} (arbitrary \cdiz) the \nst is in the weak form; \cad the proven \eqvc is between, on the one hand,
\\
$\bullet$ the \syp 
does not have any zero in any finite nonzero \Klgz,\\
 and on the other,
\\
$\bullet$ the corresponding quotient \alg is null.

The \gnl \nst states under what condition  a \pol is annihilated at the zeros of a \sypz. Here, since we do not have an \cac at our disposal, we will consider the zeros in the finite \Klgs and we obtain two \nsts depending on whether we only consider the reduced \Klgs or not.
 
These two \thos \gns from a \cof point of view (with explicit~\gui{either-or's}) the classical \nst stated in the form of \thref{thNstClass}.

\begin{theorem}\label{thNstClassCof1} \emph{(Classical \nstz, \gnl \cov version)}\\
Let $\gK$ be a \cdi and 
$f_1$, \dots, $f_s$, $g$ be in $\KXn$. 
Consider the quotient \alg $\gA=\aqo\KuX\lfs$. 
\begin{enumerate}
\item Either there exists a nonzero quotient $\gB$ of $\gA$ which is a reduced finite \Klg with $g\in\gB\eti$ (a fortiori $g\neq0$ in $\gB$).
\item Or $g$ is nilpotent in $\gA$ (in other words, there exists an integer $N$ such that $g^N\in\gen{\lfs}_\KuX$).
\end{enumerate}
\end{theorem}
%
\begin{proof}
We use Rabinovitch's trick. We introduce an additional \idtr $T$ and we notice that $g$ is nilpotent in~$\gA$ \ssi the quotient \algz~$\gA'$ for the \syp $(\lfs,1-gT)$ is null. We end with the weak \nstz: if $\gA'\neq 0$, we find a nonzero quotient $\gB'$ of $\gA'$ which is a finite dimensional \Kevz. As $g$ is \iv in $\gA'$, it is also \iv in $\gB'$ and in $\gB=\gB'\red$, and as $\gB\neq 0$, $g\neq 0$ in $\gB$.
\end{proof}
%

\begin{theorem}\label{thNstClassCof2} \emph{(\nst with multiplicities)}\\
Let $\gK$ be a \cdi and  
$f_1$, \dots, $f_s$, $g$ be in $\KXn$. 
Consider the quotient \alg $\gA=\aqo\KuX\lfs$. 
\begin{enumerate}
\item Either there exists a quotient $\gB$ of $\gA$ which is a finite dimensional \Kev  with $g\neq0$ in $\gB$.
\item Or $g=0$ in $\gA$ (in other words, $g\in\gen{\lfs}_\KuX$).
\end{enumerate}
\end{theorem}
%
\begin{Proof}{\Demo using \bdgsz. }
If when placing in the \Noe position we have $r=0$, the result is clear. The delicate point is when~$r\geq 1$. Suppose the \id is in \Noe position. We consider an  \eli  order for the variables $(\Yr)$ and the normal form of $g$ with respect to the corresponding \bdg of $\ff$. For \gui{everything to remain as is} after a specialization $Y_i\mapsto \alpha_i=\ov{Y_i}$ in a quotient \ri $\gL$ of $\gK[\Yr]$, it suffices that the leading \coes in the \bdg of $\ff$ and in the normal form of $g$ (those \coes are \elts of $\gK[\Yr]$) specialize in \iv \elts of $\gL$. If we have at our disposal enough distinct \elts in $\gK$ to find suitable $\alpha_i$'s in $\gK$ we can take $\gL=\gK$, otherwise we consider the product $h$ of all the leading \coes previously considered, and we replace $\gK[\Yr]$ with a nonzero quotient~$\gL$, \stf over $\gK$, in which $h$ is \iv (this is possible by \thref{thNstClassCof1}, applied to $h$ with no equation $f_i$). 
The~solution to our \pb is then given by the \alg
$$\preskip.4em \postskip.4em 
\gB=\gL\otimes_{\gK[\Yr]}\gA, 
$$
which is a quotient of $\gA$ \stf over $\gK$.
\end{Proof}
%

\subsec{Syzygies}
Another important consequence of the \cdv lemma (Lemma~\ref{lemNoether}) is the following \thoz.

\begin{theorem}
\label{thpolcohfd}
Let $\gK$ be a discrete \zedr \riz.
\begin{enumerate}
\item Every \pf \Klg is a \fdi \coriz.
\item Consequently every \mpf over such an \alg is \coh and strongly discrete.
\end{enumerate}
\end{theorem}
\begin{proof}
We prove the first item for $\KXn$ in the case where $\gK$ is a \cdiz. 
The \zed \ri case 
is deduced from it by the usual technique (\elgbm \num2). Then item~\emph{2} is a consequence 
of \thref{propCoh2}. 
\\
We give a \dem by \recu over $n$, the $n=0$ case being clear. Suppose $n\geq 1$ and let $\gB=\KXn$. We must show that an arbitrary \itf $\ff=\gen{f_1,\ldots ,f_s}$ is \pf and detachable.\\
If $\ff=0$ then it is clear, otherwise we can assume by applying Lemma~\ref{lemNoether} that $f_s$ is \mon in $X_n$ of degree $d$.
If $s=1$, the annihilator of~$f_1$ is null, and therefore also the module of syzygies 
for $(f_1)$.  The \id $\ff$ is detachable thanks to  Euclidean division with respect to~$X_n$. \\
If $s\geq 2$, let $\gA=\gK[X_1,\ldots ,X_{n-1}]$. The \ri $\gA$ is \fdi \coh by \hdrz. Let $R_i$ be the syzygy 
that corresponds to the \egt $f_if_s-f_sf_i=0$ ($i\in\lrb{1..s-1}$).  
Modulo the syzygies
$R_i$ we can rewrite each $X_n^kf_i=g_{k,i}$, for~$k\in\lrb{0..d-1}$ and~$i\in\lrb{1..s-1}$ as vectors in the free \Amo $L\subseteq\gB$ with basis $(1,X_n, \ldots ,X_n^{d-1})$. Modulo the syzygies 
$R_i$ every syzygy  
for $(f_1,\ldots ,f_s)$ with \coes in~$\gB$ can be rewritten as a syzygy  
for
\perso{Il semble que le contenu of the \tho \ref{thpolcohfd} soit partie int\'egrante
de la preuve originelle of the \nst by Hilbert, \`a moins que ce soit la preuve du
\tho "de la base". \`A v\'erifier and \`a signaler.
En tout cas, notre preuve of the \nst n'utilise pas ce \thoz.}

\snic{V=(g_{0,1},\ldots ,g_{d-1,1},\ldots ,g_{0,s-1},\ldots ,g_{d-1,s-1})\in
L^{d(s-1)}}

with \coes in $\gA$. As $L$ is a free \Amoz, it is \fdi \cohz. We have in particular a finite number of $\gA$-syzygies 
for $V$ that generate them all.
Let us call them $S_1$, \ldots, $S_\ell$.
Each $\gA$-syzygy 
 $S_j$ for~$V$ can be read as a $\gB$-syzygy
$S'_j$ for $(f_1,\ldots ,f_{s})$. Finally, the syzygies $R_i$ and $S'_j$ generate the \Bmo of the syzygies for $(f_1,\ldots ,f_{s})$.\\
Concerning the \fdi \crcz, we proceed in the same way. To test if an \elt of $\gB$ is in $\ff$ we start by dividing it by~$f_s$ with respect to $X_n$. We then obtain a vector in the \Amo $L$ for which we must test whether it belongs to the submodule generated by the~$g_{i,j}$'s.
\end{proof}

\vspace{-.8em}
\pagebreak

\section{The dynamic method}
\label{subsecDyna}

\begin{flushright}
{\em I do not believe in miracles.
}\\
A \cof mathematician.
\end{flushright}

In \clama \dems of existence are rarely explicit.
Two essential obstacles appear each time that we try to render such a \dem explicit.

The first obstacle is the application of \TEMz. For instance, if you consider the \dem that every univariate \pol over a field $\gK$ admits a \dcn into prime factors, you have a kind of \algo whose key ingredient is: if $P$ is irreducible all is well, if $P$ can be decomposed into a product of two factors of degree $\geq 1$, all is still well, by \hdrz.
 Unfortunately the disjunction used to make the \dem work \gui{$P$ is irreducible or $P$ can be decomposed into a product of two factors of degree $\geq 1$} is not explicit in \gnlz. In other words, even if a field is defined \cotz, we cannot be sure that this disjunction can be made explicit by an \algoz. Here we find ourselves in the presence of a typical case where \TEM \gui{is an issue,} because the existence of an irreducible factor cannot be the object of a \gnl \algoz.

The second obstacle is the application of Zorn's lemma, which allows us to \gnr to the uncountable case the usual \dems by \recu in the countable case.

For example in Modern Algebra by van der Waerden the second pitfall is avoided by limiting ourselves to the countable \agq structures.

\ss However, we have two facts that are now well established from experience:
\begin{itemize}
\item  The \emph{universal} concrete results proven by the dubious abstract methods above have never been contradicted. We have even very often successfully extracted unquestionable \prcos from them.
This would suggest that even if the abstract methods are in some way incorrect or contradictory, they have until now only been used with a sufficient amount of discernment.
\item  The key concrete results proven by the dubious abstract methods have not been invalidated either. On the contrary, they have often been validated by \algos proven \cotz.%
\footnote{On this second point, our assertion is less clear. If we return to the example of the \dcn of a \pol into prime factors, it is impossible to achieve the result \algqt over certain fields.}
\end{itemize}

\ss Faced with this slightly paradoxical situation: the abstract methods are a priori dubious, but they do not fundamentally deceive us when they give us a result of a concrete nature. There are two possible reactions.

Either we believe that the abstract methods are fundamentally correct because they reflect a \gui{truth,} some sort of \gui{ideal Cantor universe} in which exists the true semantic of \mathsz.
This is the stance taken by Platonic realism, defended for instance by G\"odel.

Or we think that the abstract methods truly are questionable. But then, unless we believe that \maths falls within the domain of magic or of miracles, it must be explained why \clama makes such few mistakes. If we believe in neither Cantor, nor miracles, we are led to believe that the abstract \dems of concrete results \ncrt contain sufficient \gui{hidden ingredients} to construct the corresponding concrete \demsz.

This possibility of \cot certifying concrete results obtained by dubious methods, if we manage to execute it systematically enough, is in line with Hilbert's program.\ihi

The dynamic method in \cov \alg is a \gnle method for decrypting
abstract \dems from \clama when they use \gui{ideal} objects whose existence relies on non-\cov principles: \TEM and the axiom of choice.
The ambition of this new method is to \gui{give a \cov semantic for the usually practiced \clamaz.}

We replace the abstract objects from \clama with incomplete but concrete specifications of these objects. This is the \cov counterpart of the abstract objects. For example a \emph{finite \ippz} (a notion that will be introduced in Section~\ref{subsecMoco}) is given by a finite number of \elts in the \id and a finite number of \elts in its complement.
This constitutes an incomplete but concrete specification of a \idepz.

More \prmtz, the dynamic method aims at giving a systematic interpretation of classical \dems that use abstract objects by rereading them as \prcos with respect to \cov counterparts of these abstract objects.

This is in keeping with the thought-process behind certain techniques developed in Computer Algebra. 
Here we are thinking about  \gui{lazy \evnz,} or  \gui{dynamic \evnz,} \cad  lazy \evn managed as a tree structure, as in the D5 system \cite{D5} which performs this tour de force very innocently: compute with certainty in the \agq closure of an arbitrary field, even though we know that this object (the \agq closure) cannot be constructed in all \gntz. 

In the current chapter an incomplete specification of the \cdr of a \spl \pol over a field $\gK$ will be given by a \Klg $\gA$ and a finite group of \autos $G$ of this \algz.
In $\gA$ the \pol can be decomposed into \lin factors such that a \cdr is a quotient of~$\gA$, and $G$ is an approximation of the Galois group in a suitable sense (in particular, it contains a copy of the Galois group). We will explain how to compute with such an approximation without ever making a mistake: when an oddity occurs, we know how to better the approximation during computation and to make the oddity disappear.

\subsect{Splitting fields and Galois theory in \clamaz}{Classical Galois theory}\label{GalPolCla}

In this subsection we will offer a possible presentation of the \cdr of an arbitrary \pol and of the Galois theory of a \spl \pol in \clamaz.
This allows us to understand the 
\gui{detours} 
that we will be obligated to take to have an entirely \cov theory.

\smallskip
If $f$ is a \poluz, we work with the \adu of $f$, $\gA=\Adu_{\gK,f}$ in which $f(T)=\prod_i(T-x_i)$, with $\Sn$ as a group of \autos 
(see Section~\ref{sec0adu}).

This \alg being a finite dimensional \Kevz, all the \ids are themselves finite dimensional \Kevs and we have the right to consider a strict \id $\fm$ of maximum dimension as a \Kev (all of this by applying \TEMz).
This \id is automatically a \idemaz.
The quotient \alg $\gL=\gA/\fm$ is then a \cdr for~$f$.
The group $G=\St(\fm)$ operates on $\gL$ and the fixed field of $G$, $\gL^G=\gK_1$, possesses the two following \prtsz:
\begin{itemize}
\item $\gL/\gK_1$ is a Galois extension with $\Gal(\gL/\gK_1)\simeq G$.
\item  $\gK_1/\gK$ is an extension obtained by successive additions of $p^{\rm th}$ roots, where $p=\car(\gK)$.
\end{itemize}

\smallskip
Moreover, if $\gL'$ is another \cdr for $f$ with $f=\prod_i(T-\xi_i)$ 
in~$\gL'[T]$, we have a unique \homo of \Klgs $\varphi:\gA\to\gL'$
satisfying the \egts $\varphi(x_i)=\xi_i$ for $i\in\lrbn$. We can then show that~$\Ker\varphi$, which is a \idema of $\gA$, is \ncrt 
a conjugate of $\fm$ under the action of~$\Sn$.
Thus the \cdr is unique, up to \iso (this \iso is not unique if $G\neq\so{\Id}$).

Finally, when $f$ is \splz, the situation is simplified because the \adu is \'etale, and $\gK_1=\gK$.

\smallskip
The previous approach is possible from a \cof point of view if the field $\gK$ is \sply factorial and if the \pol $f$ is \splz, because then, since the \adu $\gA$ is \'etale, it can be decomposed into a finite product of \'etale fields over $\gK$ (Corollary~\ref{propIdemMini}).

But when the field is not \sply factorial, we face an a priori insurmountable obstacle, and we cannot hope to systematically and \algqt  obtain a \cdr that is \stf over~$\gK$.

If the \cara is finite and if the \pol is not \splz, we need stronger \fcn \prts to construct a \cdr (the question is delicate, and very well presented in \cite{MRR}).

\subsect{Lazily bypassing the obstacle}{Bypassing the obstacle}\label{GalPolComa}

What is \gnlt proposed in Computer Algebra is, for instance in the case of a \spl \polz, at the very least to avoid computing a \uvl resolvent $R$ (as in \thref{thResolUniv}) whose degree,~$n!$, promptly renders computations impractical.

Here, we find ourselves in the most \gnl framework possible, and we avoid all recourse to the \fcn of the \pols which can turn out to be impossible, or which, when possible, has the risk of being too costly.

The idea is to use the \adu $\gA$, or a Galois quotient $\aqo{\gA}{1-e}$, with a Galoisian \idm $e$ (see \paref{defQuoDeGal}) as a substitute for $\gL$. This \gui{dynamic \cdrz} can be managed without too many \pbs because each time something strange happens, which indicates that the substitute of $\gL$ is not entirely satisfying, we are able to \gui{immediately repair the oddity} by computing a Galoisian \idm that refines the previous one, and in the new approximation of the \cdrz, the strange thing has disappeared.

To develop this point of view we will need to better know the \aduz, and Section~\ref{sec1adu} is dedicated to this objective.

Moreover, we will study in Section~\ref{subsecCDR} a dynamic and \cov version of the \cdr of a (not \ncrt \splz) \polz.

\section{Introduction to \agBs} \label{secBoole}

A \ix{lattice} is a set $\gT$ equipped with an order relation $\leq$ for which there exist a minimum \eltz, denoted by $0_\gT$, a maximum \eltz, denoted by $1_\gT$, and every pair of \elts $(a,b)$ admits a least upper bound, denoted by $a\vu b$, and a greatest lower bound, denoted by $a\vi b$.
A mapping from one lattice to another is called a \emph{lattice \homoz} if it respects the operations $\vu$ and $\vi$ as well as the constants $0$ and $1$.
The lattice is called a \ixy{distributive}{lattice} when each of the two operations $\vu$ and $\vi$ is distributive with respect to the other.

We will give a succinct study of the structure of \trdis and of structures that relate back to them in Chapter~\ref{chapTrdi}.

\begin{propdef}\label{defiBoole}
\emph{(\agBsz)}\index{algebra!Boolean ---}\index{Boolean!algebra}
\begin{enumerate}
\item By \dfn a \ri $\gB$ is a \emph{\agBz} \ssi every \elt is idempotent. Consequently $2=_\gB 0$ (because $2=_\gB 4$).
\item We can define over $\gB$ an order relation $x\preceq y$ by: $x$ is a multiple of~$y$, \cad $\gen{x}\subseteq\gen{y}$.
Then, two arbitrary \elts admit a lower bound, their lcm $x\vi y=xy$, and an upper bound, their gcd $x\vu y=x+y+xy$.
We thus obtain a \trdi with~$0$ as its minimum \elt and $1$ as its maximum \eltz.
\item For every $x\in\gB$, the \elt $x'=1+x$ is the unique \elt that satisfies the \egts $x\vi x'=0$ and $x\vu x'=1$, we call it \emph{the complement} of $x$.
\end{enumerate}
\index{complement!(in a \agBz)}
\end{propdef}

{\it Notation conflict.} Here we find ourselves with a conflict of notation. Indeed,  \dve in a \ri leads to a notion of the gcd, which is commonly denoted by $a\vi b$, because it is taken as a lower bound ($a$ divides $b$ being understood as \gui{$a$ smaller than $b$} in the sense of the \dvez). This conflicts with the gcd of the \elts in a \agBz, which is an upper bound. This is due to the fact that the order relation 
has been reversed, so that the elements $0$ and $1$ of the Boolean algebra are indeed the minimum and the maximum in the lattice. 
This inevitable conflict
 will appear in an even stronger sense when we will consider the \agB of the \idms of a \ri $\gA$.
\eoe

\smallskip 
Even though all the \elts of a \agB are \idms we will keep the terminology \gui{\sfioz\footnote{It would be more natural to say: fundamental \sys of \ort \eltsz.}} for a finite family $(x_i)$ of pairwise \ort \elts (\cadz~$x_ix_j=0$ for $i\neq j$) with sum  $1$.
This convention is all the more justified in that we will mainly preoccupy ourselves with the \agB that naturally appears in commutative \algz: that of the \idms of a \ri $\gA$.

\penalty-2500
\subsec{Discrete \agBsz}
\label{subsec1AGBDiscretes}

\begin{proposition}\label{propBoolFini} \emph{(Every discrete \agB behaves in computations as the \alg of the detachable subsets of a finite set)}\\
Let $(r_1,\ldots,r_m)$ be a finite family in a \agB $\gB$.\\
Let $s_i=1-r_i$ and, for a finite subset $I$ of $\{1,\ldots ,m\}$, let $r_I=\prod_{i\in I}r_i\prod_{j\notin I}s_j$. 
\begin{enumerate}
\item  The $r_I$'s form a \sfio and they generate the same \agB as the $r_i$'s.
\item Suppose that $\gB$ is discrete. Then, if there are exactly $N$ nonzero \eltsz~$r_I$, the \agsB generated by the $r_i$'s is \isoc to the \alg of finite subsets of a set with $N$ \eltsz.
\end{enumerate}
\end{proposition}

As a corollary we obtain the following fact and the fundamental structure \tho that summarizes it.
Recall that we denote by $\Pf(S)$ the set of finite subsets of a set $S$.

In a discrete \agB an \elt $e$ is called an \ix{atom} if it satisfies one of the following \eqve \prtsz.
\begin{itemize}
\item $e$ is minimal among the nonzero \eltsz.
\item  $e\neq 0$ and for every $f$, $f$ is \ort or greater than $e$.
\item  $e\neq 0$ and for every $f$, {\mathrigid 2mu $ef=$ $0$ or $e$,
or  $ef=0$ or $e(1-f)=0$}.
\item $e\neq 0$ and the \egt {\mathrigid 2mu $e=e_1+e_2$ with $e_1e_2=0$ implies $e_1=0$ or $e_2=0$}.
\end{itemize}
We also say that $e$ is \emph{indecomposable}. It is clear that an \auto
of a discrete \agB preserves the set of atoms and that for two atoms~$e$ and $f$, we have $e=f$ or $ef = 0$.\index{indecomposable!\elt in a \agBz}

\begin{theorem}\label{corpropBoolFini}\label{factagb} \emph{(Structure \thoz)}
\begin{enumerate}
\item Every finite \agB is \isoc to the \alg of the detachable subsets of a finite set. 

\item 
More \prmtz, for a \agB $C$ \propeq
\begin{enumerate}
\item $C$ is finite.
\item $C$ is discrete and \tfz.
\item The set $S$ of atoms is finite, and  $1_C$ is the sum of this set.
\end{enumerate}
In such a case $C$ is \isoc to the \agB $\Pf(S)$.
\end{enumerate}
\end{theorem}

\subsect{\agB of the \idms of a commutative \riz}{\agB of the \idms of a \riz}
\label{subsecAGBIDMSAC}
\begin{fact}\label{propB(A)} 
The \idms of a \ri $\gA$ form a \agBz, denoted by $\BB(\gA)$, with the operations $\vi$, $\vu$, $\lnot$ and $\oplus$ given by

\snic{r\vi s=rs$, $\;r\vu s =r+s-rs\;$ , $\;\lnot\ r=1-r\;$ and $\;r\oplus s=(r-s)^2.}


If $\gA$ is a \agBz, $\BB(\gA)=\gA$.
If $\varphi :\gA\to\gB$ is a morphism of \risz, its restriction to $\BB(\gA)$ gives a morphism $\BB(\varphi ):\BB(\gA)\to\BB(\gB)$.
\end{fact}
%
\begin{proof}
It suffices to show that if we equip the set $\BB(\gA)$ with the laws $\oplus$ and~$\vi$ we obtain a \agB with $0_\gA$ and $1_\gA$ as neutral \eltsz. The computations are left to the reader.
\end{proof}

\medskip \Thref{corpropBoolFini} has the following immediate consequence.

\begin{fact}\label{factZerloc2}
\Propeq
\begin{enumerate}
\item The \agB of the \idms $\BB(\gA)$ is finite.
\item The \ri $\gA$ is a finite product of nontrivial connected \risz.
\end{enumerate}
\end{fact}
\begin{proof}
It suffices to show that \emph{1} implies \emph{2.}
If $e$ is an atom of $\BB(\gA)$, the \riz~$\gA[1/e]$ is nontrivial and connected.
If $\BB(\gA)$ is finite, the finite set~$A$ of its atoms forms a \sfio of $\gA$, and we have a canonical \iso $\gA\to\prod_{e\in A}\gA[1/e]$.
\end{proof}

\rem If $\BB(\gA)$ has a single \eltz, $\gA$ is trivial and the finite product is an empty product. This also applies to the following corollary.
\eoe

\vspace{-.2em}
\pagebreak	

\begin{corollary}\label{lemZerloc2}
\Propeq
\begin{enumerate}
\item $\BB(\gA)$ is finite and $\gA$ is \zedz.
\item   $\gA$ is a finite product of nontrivial \zed local \risz.
\end{enumerate}
\end{corollary}

\subsec{Galoisian \elts in a \agB}

\begin{definition}
\label{defIdmGalAgb}~
\index{Boolean G-alg@Boolean $G$-algebra}
\index{Boolean G-alg@Boolean $G$-algebra!transitive ---}
\index{transitive!Boolean G-alg@Boolean $G$-algebra}
\begin{enumerate}
\item If $G$ is a group that operates over a \agB $C$, we say that the pair $(C,G)$ is a $G$-\agBz.
\item An \elt $e$ of a $G$-\agB $C$ is said to be \ixc{Galoisian}{\elt in a \agBz} if its orbit under $G$ is a \sfioz.
\item A $G$-\agB is said to be \emph{transitive} if $0$ and $1$
are the only \elts fixed by $G$.
\end{enumerate}
\end{definition}

\Dfn \ref{defQuoDeGal} of Galoisian \idms  agrees with the previous \dfn
when a finite group~$G$ acts on a \klg $\gC$ and when we consider the action
of $G$ over the \agB $\Bo(\gC)$.

Now we study the case where the group is finite and the \alg discrete.

\begin{fact}
\label{factAduIdmA}
Let $G$ be a finite group and $C$ be a transitive, discrete and nontrivial $G$-\agBz. Let $e\neq 0$ in $C$, and $\so{e_1,\ldots ,e_k}$ be the orbit of $e$ under~$G$.
\Propeq
\begin{enumerate}
\item \label{E1} The \elt $e$ is Galoisian.
\item \label{E2} For all $i>1$, $e_1e_i=0$.
\item \label{E21} For all $\sigma\in G$, $e\sigma(e)=e$  or $0$.
\item \label{E22} For all $i\neq j\in\so{1,\ldots ,k}$, $e_ie_j=0$.
\end{enumerate}
\end{fact}
\begin{proof}
Item \emph{\ref{E1}} clearly implies the others. Items \emph{\ref{E2}} and \emph{\ref{E22}} are easily \eqvs and imply item \emph{\ref{E21}}. Item \emph{\ref{E21}} means that for every $\sigma$, $\sigma(e)\geq e$ or $\sigma(e)e=0$. If we have $\sigma(e)\geq  e$ for some $\sigma$, then we obtain 

\snic{e\leq \sigma(e)\leq \sigma^2(e)\leq \sigma^3(e)\leq \ldots\,, }

which gives us $e=\sigma(e)$ when considering an $\ell$ such that $\sigma^\ell=1_G$.
Therefore, item~\emph{\ref{E21}} implies item~\emph{\ref{E2}}.
Finally, if item~\emph{\ref{E22}} is satisfied, the sum of the orbit is an \eltz~$>0$ fixed by~$G$ therefore equal to~$1$.
\end{proof}

\begin{lemma}\label{lemIdmGalAgb} \emph{(Meeting of two Galoisian \eltsz)}\\
Let $G$ be a finite group and $C$ be a nontrivial discrete $G$-\agBz. Given two Galoisian \elts $e$, $f$ in $(C,G)$, let

\snic{G.e=\so{e_1,\ldots,e_m}$, $E=\St_G(e)$, and $F=\St_G(f).}
\begin{enumerate}
\item There exists a $\tau\in G$ such that $f\tau(e)\neq0$.
\item
If $e\leq f$, then $E\subseteq F$ and $f=\sum_{i:e_i\leq f}e_i=\sum_{\sigma\in F/E}\sigma(e)$.
\end{enumerate}
Suppose $C$ is transitive and $ef\neq0$. We obtain the following results.
\begin{enumerate}\setcounter{enumi}{2}
\item
The \elt $ef$ is Galoisian, with stabilizer $E\,\cap\, F$, and the orbit $G.ef$ consists of nonzero \elts of $(G.e)(G.f)$.  In particular, $G.ef$ generates the same \agsB of $C$ as $G.e\,\cup\, G.f$.
\item If  $E\subseteq F$, then $e\leq f$.
\end{enumerate}
\end{lemma}
\begin{proof}
\emph{1.} Indeed, $f=\sum_ife_i$.

\emph{2.}
\Gnltz, for $x' = \sigma(x)$ where $x \ne 0$ satisfies $x \le f$, let us show 
$$x' \le f \;\;\buildrel {[a]} \over \Longrightarrow\;\;
fx' \ne 0 \;\;\buildrel {[b]} \over \Longrightarrow\;\;
\sigma(f) = f \;\;\buildrel {[c]} \over \Longrightarrow\;\;
x' \le f.\leqno(\star)
$$
We obtain $[a]$ by multiplying $x' \le f$ by $x'$, $[b]$ by multiplying $x' \le \sigma(f)$ (deduced from $x \le f$) by $f$ and by using the fact that $f$ is Galoisian and finally $[c]$ by applying~$\sigma$ to $x \le f$. The assertions of $(\star)$ are therefore \eqvcsz. 
We deduce~$\St_G(x) \subseteq \St_G(f)$. If in addition, $1 = \sum_{x' \in G.x} x'$, then
$$
\preskip.0em \postskip.4em \ndsp
f = \sum_{x' \in G.x} fx' = \sum_{x'\in G.x | x' \le f} x' =
\sum_{\sigma \in F/\St_G(x)} \sigma(x).
$$
This applies to $x = e$.

\emph{3.} Let $G.f=\so{f_1,\ldots,f_p}$.
For $\sigma\in G$ there exist $i$, $j$ such that 

\snic{e\,f\,\sigma(ef) = e\,f\,e_i\,f_j,}

which is equal to $ef$ if $\sigma\in E\cap F$ and to $0$ otherwise. By Fact~\ref{factAduIdmA}, $ef$ is therefore a Galoisian \elt with stabilizer $E \cap F$. Now assume $e_if_j\neq0$. Then, by item~\emph{1}, there exists a $\tau\in G$ such that $\tau(ef) e_if_j\neq0$. Since~$e$ and~$f$ are Galoisian, this implies $\tau(e)=e_i$ and $\tau(f)=f_j$, so~$e_if_j \in G.ef$.

\emph{4.} Immediately results from \emph{3.}
\end{proof}

The paradigmatic application  of the next \tho is the following.
We have a nontrivial connected \ri $\gk$, $(\gk,\gC,G)$ is a \pGa (cf. 
\Dfnz~\ref{defADG1}) or \aG and we take $C=\BB(\gC)$.

\begin{theorem}
\label{lemAduIdmA} \emph{(Galois structure \thoz, 1)} 
Let $G$ be a finite group and $C$ be a transitive, discrete and nontrivial $G$-\agBz.
\begin{enumerate}
\item \label{lemAduIdmA-1}  \emph{(Structure of the transitive finite $G$-\agBsz)}
\\
The \alg $C$ is finite \ssi there exists an atom $e$. In this case the structure of $(C,G)$ is entirely \caree by  $E=\St_G(e)$.
\\
More \prmtz, the \idm $e$ is Galoisian, $G.e$ is the set of atoms, {$C\simeq\Pf(G.e)$,} $G$ operates over $G.e$ as it does over $G/E$, and over $C$ as it does \hbox{over $\Pf(G/E)$}. In particular, $|C|=2^{\idg{G:E}}.$
\\
We will say that $e$ is a \emph{Galoisian \gtr of~$C$}.
\item \label{lemAduIdmA-2}  Every finite family of \elts of $C$ generates a finite $G$-sub\algz.
\item \label{lemAduIdmA-3}  The \agB $C$ cannot have more than $2^{|G|}$ \eltsz.
\item \label{lemAduIdmA-4}  Let $e$ and $f$ be Galoisian \eltsz, $E=\St_G(e)$ and $F=\St_G(f)$.
\begin{enumerate}
\item There exists a $\sigma\in G$ such that $f\sigma(e)\neq0$.
\item If $ef\neq0$, $ef$ is a Galoisian \gtr of the $G$-sub\agB of $G$ generated by $e$ and $f$, and  $\St_G(ef)=E\cap F$.
\item If $e\leq f$ (\cad $fe=e$), then $E\subseteq F$ and $f=\sum_{\sigma\in F/E}\sigma(e)$.
\end{enumerate}
\item \label{lemAduIdmA-5}  \emph{(\Carn of the Galoisian \elts in a finite $G$-sub\algz)}
\\
Let $e$ be a Galoisian \elt and $f$ be a sum of $r$ \elts of $G.e$, including~$e$. Let $E=\St_G(e)$  and $F=\St_G(f)$. Then \propeq
\begin{enumerate}
\item \label{lemAduIdmA-5a}$f$ is Galoisian.
\item \label{lemAduIdmA-5b}$E\subseteq F$ and $f=\sum_{\sigma\in F/E}\sigma(e)$.
\item \label{lemAduIdmA-5c}$\abs{F}=r\times \abs{E}$.
\item \label{lemAduIdmA-5d}$\abs{F}\geq r\times\abs{E} .$
\end{enumerate}
%
\end{enumerate}
\end{theorem}

\vspace{.1em}
\begin{proof}
\emph{\ref{lemAduIdmA-1}.} If $C$ is finite there exists an atom.
If $e$ is an atom, for every $\sigma\in G$, we have $e\,\sigma(e)=0$ or $e$, so $e$ is Galoisian (Fact~\ref{factAduIdmA}). The rest follows by taking into account \thref{factagb}.

\emph{\ref{lemAduIdmA-2}.} Consider the \agsB $C'\subseteq C$ generated by the orbits of the \elts of the given finite family. $C'$ is \tf and discrete therefore finite. Consequently its atoms form a finite set $S=\so{e_1,\ldots,e_k}$ and $C'$ is \isoc to the \agB of the finite subsets of $S$ 

\snic{C'=\sotq{\sum_{i\in
F}e_i}{F\in\cP_k}.}

Clearly, $G$ operates on $C'$. For $\sigma\in G$, $\sigma (e_1)$ is an atom, so~$e_1$ is Galoisian (Fact~\ref{factAduIdmA}~\emph{\ref{E21}}) and $(e_1,\ldots,e_k)$ is its orbit.

\emph{\ref{lemAduIdmA-3}.} Results from \emph{\ref{lemAduIdmA-1}} and \emph{\ref{lemAduIdmA-2}}.

\emph{\ref{lemAduIdmA-4}.} Already seen in Lemma~\ref{lemIdmGalAgb}.

\emph{\ref{lemAduIdmA-5}.} We write $\sigma_1=1_G$, $G.e=\so{\sigma_1.e,\ldots,\sigma_k.e}$ with $k=\idg{G:E}$, as well as~$f=\sigma_1.e+\cdots +\sigma_r.e$.
\\
\emph{\ref{lemAduIdmA-5a}} $\Rightarrow$ \emph{\ref{lemAduIdmA-5b}.} We apply item~\emph{\ref{lemAduIdmA-4}.}
\\
\emph{\ref{lemAduIdmA-5b}} $\Rightarrow$ \emph{\ref{lemAduIdmA-5a}.} For $\tau\in F$, $\tau.f=f$. 
 \\
For $\tau\notin F$, $F.e\cap (\tau F).e=\emptyset$, and so~$f\,\tau (f)=0$.
\\
\emph{\ref{lemAduIdmA-5b}} $\Rightarrow$ \emph{\ref{lemAduIdmA-5c}.} We have $F.e=\so{1_G.e,\sigma_2.e,\ldots,\sigma_r.e}$, and since $E$ is the stabilizer of $e$, we obtain $\abs{F}=r\times \abs{E}$.
\\
\emph{\ref{lemAduIdmA-5d}} $\Rightarrow$ \emph{\ref{lemAduIdmA-5b}.}
We have $F=\sotq{\tau}{\tau\so{\sigma_1.e,\ldots, \sigma_r.e}=\so{\sigma_1.e,\ldots,\sigma_r.e}}$. Hence the inclusion $F.e\subseteq \so{\sigma_1.e,\ldots,\sigma_r.e}$, and 
$F.e=\so{\sigma_1.e,\ldots,\sigma_s.e}$ with $s\leq r\leq k$. The stabilizer of $e$ for the action of $F$ on $F.e$ is equal to $E\cap F$.
Therefore
$$\abs{F}=\abs{F.e}\,\abs{E\cap F}= s\,\abs{E\cap F} \leq r\,\abs{E\cap
F}\leq r\,\abs{E}.$$
Therefore if $\abs{F}\geq r\,\abs{E}$, we have $\abs{F.e}=r$ and $\abs{E}=\abs{E\cap F}$, \cad $E\subseteq F$   
and $F.e=\so{\sigma_1,\ldots ,\sigma_r}$.
\end{proof}

Under the hypotheses of \thref{lemAduIdmA} we can compute a Galoisian \eltz~$e_1$ such that $G.e_1$ and $G.e$ generate the same \agBz, by means of \Algoz~\vref{algidmgal}.
\pagebreak

\begin{algor}[Computation of a Galoisian \elt and of its stabilizer.] \label{algidmgal}
\Entree $e$: nonzero \elt of a \agB $C$; $G$: finite group of \autos of $C$; $S=\St_G(e)$.
\\ \# \, \emph{Suppose that $0$ and $1$ are the only fixed points for the action of $G$ on $C$.}
\Sortie $e_1$: a Galoisian \elt of $C$ such that $G.e_1$ generates the same \agB as $G.e$; $H$: the stabilizer subgroup of $e_1$.
\Varloc $h$: in $C$; $\sigma$: in $G$;
$L$: list of \elts of $\ov{G/S}$;\\
$E$: corresponding set of \elts of $G/S$;
\\ \# \, $G/S$  is the set of left cosets of $S$. 
\\ \# \, $\ov{G/S}$ is a \sys of representatives of the left cosets of $S$ \Debut
\hsu $E \aff \emptyset$; $L \aff [\,]$; $e_1\aff e$; 
\hsu \pur{\sigma}{\ov{G/S}}
\hsd $h\aff e_1\sigma(e)$;
\hsd \sialors{h\neq 0} $e_1\aff h$; $L\aff L\bullet [\sigma]$;
$E \aff E \cup \so{\sigma S}$
\hsd \finsi;
\hsu \finpour;
\hsu $H\aff \St_G(E)$\quad  \# \, $H =\sotq{\alpha\in G} 
 {\forall\sigma\in L,    \alpha\sigma\in \bigcup_{\tau\in L}\tau S}$.
\Fin
\end{algor}

\begin{Proof}{Correctness proof of the \algoz. }
We denote by $G/S$ a \sys of representatives of the left cosets of $S$. Let us write $e_1=e \sigma_2(e)\cdots \sigma_r(e)$ where the~$\sigma_i$'s are all the $\sigma$'s which have successfully passed the test $h\neq 0$ in the \algo (and $\sigma_1=\Id$). We want to show that $e_1$ is an atom of $C'$ (the \agB generated by $G.e$), which is the same as saying that for all $\sigma\in G/S$ we have $e_1\sigma(e)=e_1$ or $0$ (since $C'$ is generated
by the $\tau(e)$'s). However, $\sigma$ has been tested by the \algoz, therefore either $\sigma$ is one of the $\sigma_i$'s, in which case $e_1\sigma(e)=e_1$, or $g\sigma(e)=0$ for some \idm $g$ which divides $e_1$, and a fortiori $e_1\sigma(e)=0$.\\
Let us show that the stabilizer $H$ of $e_1$ indeed satisfies the required condition. We have $e_1=\prod_{\tau\in L}\tau(e)$, and for $\sigma \in G$ we have the \eqvcs

\snic{
\sigma\in \bigcup_{\tau\in L}\tau S \;\Longleftrightarrow\;
e_1\sigma(e)=e_1\;\Longleftrightarrow\; e_1\leq
\sigma(e),\qquad\mathrm{and}}

\snic{
\sigma\notin \bigcup_{\tau\in L}\tau S \;\Longleftrightarrow\; e_1\sigma(e)=0.}

For $\alpha\in G$ we have $\alpha(e_1)=\prod_{\tau\in L}\alpha\big(\tau(e)\big)$. This is an \elt of the orbit of $e_1$, it is equal to $e_1$ \ssi $e_1\leq \alpha(e_1)$, \ssi $e_1\leq \alpha\big(\sigma(e)\big)$ for each $\sigma$ in $L$. Finally, for some arbitrary $\sigma$ in $G$, $e_1\leq \alpha\big(\sigma(e)\big)$  \ssi $\alpha \sigma$ is in $\bigcup_{\tau\in L}\tau S$.
\end{Proof}

Note that the \elt $e_1$ obtained as a result of this computation depends on the order in which the finite set  $G/S$ is enumerated and that there is no (intrinsic) natural order on this set.

\medskip \exl
We can ask ourselves if there exists a relation between the stabilizer~$S$ of $e$ and the stabilizer $H$ of a Galoisian \elt $e_1$ associated with $e$.
Here is an example that shows that there is no close relation, with $G=\rS_6$ operating on $\Adu_{\QQ,f}$ with the \pol $f(T)=T^6-4T^3+7$.
We consider the \idm $e=1/6(x_5^3x_6^3 - 2x_5^3 - 2x_6^3 + 7)$ that we compute from a \fcn of the \polmin of the \elt $x_5+x_6$  (cf. Proposition~\ref{propRvRel1}). \\
We find $\St(e)=S=\gen{(1432),(12),(56)}\simeq \rS_4\times \rS_2$ with $|S|=48$, and

\snac{\St(e_1)=H=\gen{(24),(123456)}=(\gen{(13),(135)}\times
\gen{(24),(246)})\rtimes \gen{(14)(25)(36)}
}

 with $H\simeq (\rS_3\times \rS_3)\rtimes \rS_2$, $|H|=72$, and $S\cap H=\gen{(24),(1234)(56)}$ dihedral of order $8$.
\perso{cela serait bien of give un autre example aboutissant a un subgroup $H$ compl\`etement diff\'erent (par example $H=S$).}\\
In short, $H$
(not even the conjugacy class of $H$ in $G$) 
 cannot be computed solely from $S$. Indeed, the list $L$ of left cosets selected by the \algo does not only depend on subgroup $S$ of~$G$ but also on the way in which $G$ operates on $C$.
\eoe

\section{The \adu  (2)}
\label{sec1adu}
\label{secadu}

Here is a small reading guide for the end of this chapter.

In Section~\ref{secGaloisElr}, we have seen that if $\gk$ is an infinite discrete field, if $f$ is \spl and if we are able to decompose a Galois resolvent into a product of \ird factors, then the \adu $\gA$ is \isoc to $\gL^r$, where $\gL$ is a \cdr for $f$  
 and $r=\idg{\Sn:G}$, where $G$ is a subgroup of $\Sn$ which is identified with $\Gal(\gL/\gk)$. Moreover, $\dex{\gL:\gk}=\abs G$.

We will see that this ideal situation can serve as a guideline for a lazy approach to the construction of a \cdrz.
What replaces the complete \fcn of a Galois resolvent is the discovery or the construction of a Galoisian \idmz.
Then, we have a situation analogous to the ideal situation previously described: $\gA\simeq\gB^r$, where
$\gB$ is a {Galois quotient} of~$\gA$, equipped with a group of \autos that can be identified with a subgroup $G$ of $\Sn$, with $\dex{\gB:\gk}=\abs G$ and $r=\idg{\Sn:G}$.

\smallskip

\Grandcadre{
Throughout Section~\ref{secadu},  $\gk$ is a commutative \riz, \\
$f=T^n+\sum_{k=1}^n(-1)^{k} s_kT^{n-k} \in \kT$ is \mon of degree $n$,\\
and $\gA =\Adu_{\gk,f}$ is the \adu of $f$ over $\gk$.
}

\vspace{-.7em}
\pagebreak

Recall that the \adu

\snic{\gA =\Adu_{\gk,f}=\aqo{\gk[\uX]}{S_1-s_1,\ldots,S_n-s_n}=\kuX/\cJ(f)}

(where the $S_i$'s are the \elr \smq \pols in the $X_i$'s) is the \alg which solves the \uvl \pb linked to the \dcn of the \pol $f$ into a product of factors $T-\xi_j$ (cf. Fact~\ref{factEvident}).
The \kmo $\gA=\Adu_{\gk,f}$ is free, and a basis is formed by the \gui{\momsz} $x_1^{d_1}\cdots x_{n-1}^{d_{n-1}}$ such that for $k\in\lrb{0..n-1}$, we have $d_k\leq n-k$ (see Fact~\ref{factBase}). We will denote this basis by~$\cB(f)$.
\label{notaBaseadu}

By a change of the base \riz, we obtain the following important fact (to be distinguished from Fact~\ref{factAduAdu}).
\begin{fact}
\label{factChangeBase}\emph{(Changing the base \ri for a \aduz)}
Let $\rho:\gk\to\gk_1$ be a \klgz. Let $f^{\rho}$ be the image of $f$ in $\gk_1[T]$. Then, the \alg  $\rho\ist(\Adu_{\gk,f})=\gk_1\otimes_\gk \Adu_{\gk,f}$ is naturally \isoc to  $\Adu_{\gk_1,f^{\rho}}$.
\end{fact}

\subsec{Galois quotients of \apGsz} \label{subsecIdGal}

\medskip If $\gC$ is a \klgz, we denote by $\Aut_\gk(\gC)$ its group of \autosz.

We now give a \dfn that allows us to place the \adu in a framework that is a little more \gnl and useful.

\begin{definition}
\label{defADG1} \emph{(\apGsz)}\\
A \emph{\apGz} is given by a triple $(\gk,\gC,G)$ where
\begin{enumerate}
\item $\gC$ is a \klg with $\gk\subseteq\gC$, $\gk$ a direct summand in $\gC$,
\item $G$ is a finite group of $\gk$-\autos of $\gC$,
\item $\gC$ is a \pro \kmo of constant rank $|G|$,
\item for every $z\in\gC$, we have $\rC {\gC/\gk}(z)(T)=\rC {G}(z)(T)$.
\end{enumerate}
\index{algebra!pre-Galois ---}
\end{definition}
\rem Recall that by Lemma~\ref{lemIRAdu}, if $\gB$ is a faithful \stfe \klg, then $\gk$ (identified with its image in $\gB$) is a direct summand in $\gB$. Consequently item~\emph{1} above results from item~\emph{3}.  \eoe

\smallskip 
\exls  1) From what we already know on the \adu (Section~\ref{sec0adu}) and by Lemma~\ref{lemPolCarAdu},
for every \polu $f$, the triple $(\gk,\Adu_{\gk,f},\Sn)$ is a \apGz.

2) Artin's \tho (Theorem~\ref{thA}) shows that every \aG is a \apGz.
\eoe

\medskip  The reader should refer to \paref{defQuoDeGal} for the \dfns of a Galoisian \idmz, of a Galoisian \id and of a Galois quotient.

\pagebreak	

\begin{theorem}
\label{thADG1Idm} \emph{(Galoisian structure \thoz, 2)}\\
Consider a \apG $(\gk,\gC,G)$. Let $e$ be a Galoisian \idm of~$\gC$, and~$\so{e_1,\ldots ,e_m}$ its orbits under $G$. Let~$H$ be the stabilizer of $e=e_1$ and $r=|H|$, so that $rm=|G|$. Let $\gC_i={\gC}[1/e_i]$ \hbox{for $(i\in\lrbm)$}. Finally, let $\pi:\gC\to\gC_1$ be the canonical projection.
\begin{enumerate}
\item \label{thADG1Idm1} $(\gk,\gC_1,H)$ is a \apG (in other words a Galois quotient of a \apG is a \apGz).
\item \label{thADG1Idm2} The $\gC_i$'s are pairwise
\isoc \klgsz, and $\gC\simeq \gC_1^m$.
\item \label{thADG1Idm3} The \alg $\gC_1$ is a \kmrc $r=|H|$. The restriction of $\pi$ to $\gk$, and even to $\gC^G$, is injective, and $\gk$ (identified with its image in $\gC_1$) is a direct summand in~$\gC_1$.
\item \label{thADG1Idm4} The group $H$ operates on $\gC_1$ and $\gC_1^H$ is canonically \isoc to $\gC^G$; more \prmtz, $\gC_1^H=\pi(\gC^H)=\pi(\gC^G)$.
\item \label{thADG1Idm5} For all $z\in\gC_1$, $\rC {\gC_1/\gk}(z)(T)=\rC {H}(z)(T)$.
\item \label{thADG1Idm6} Let $g_1$ be a Galoisian \idm of  $(\gk,\gC_1,H)$, $K$ its stabilizer in~$H$, and $g'\in e_1\gC$ be such that $\pi(g')=g_1$. Then, $g'$ is a Galoisian \idm of~$(\gk,\gC,G)$, its stabilizer is $K$, and we have a canonical \isoz~$\aqo{\gC_1}{1-g_1}\simeq\aqo{\gC}{1-g'}$.
\item \label{thADG1Idm7} If $(\gk,\gC,G)$ is a \aGz, then so is $(\gk,\gC_1,H)$.
\end{enumerate}
\end{theorem}

\begin{proof} Item~\emph{\ref{thADG1Idm1}} is a partial synthesis of items~\emph{\ref{thADG1Idm2}},
\emph{\ref{thADG1Idm3}}, \emph{\ref{thADG1Idm4}}, \emph{\ref{thADG1Idm5}.}

Item~\emph{\ref{thADG1Idm2}} is obvious. An immediate consequence is the first assertion of item~\emph{\ref{thADG1Idm3}.} Let $(\tau_1,\tau_2, \ldots, \tau_m)$ be a \sys of representatives for $G/H$, with $\tau_1=\Id$ and $\tau_i(e_1)=e_i$. Let us show that the restriction of $\pi$ to $\gC^G$ is injective. If $a\in\gC^G$ and
$e_1a=0$, then, by transforming by the $\tau_j$'s, all the~$e_ja$'s are null, and hence so is their sum, $a$. 
Finally, $\pi(\gk)$ is a direct summand in $\gC_1$ by Lemma~\ref{lemIRAdu}.

\emph{\ref{thADG1Idm4}.} 
Let us first show {\mathrigid 2mu $\gC_1^H=\pi(\gC^H)$.  Let $z\in\gC_1^H$ and $u\in\gC$ such that $\pi(u) = z$}.  Since $z\in\gC_1^H$, for $\sigma\in H$, $\sigma(u)\equiv u$ mod $\gen{1-e_1}$, \cad $e_1\sigma(u)=e_1 u$ or, since $\sigma(e_1) = e_1$, $\sigma(e_1u)=e_1 u$.  By letting $y = e_1u$, we see that~$y$ is $H$-invariant and $\pi(y) = z$.
\\
Let us now show that $z\in\pi(\gC^G)$. Let 

\snic{v=\sum_i\,\tau_i(y)
=\sum_i\,\tau_i(e_1y) =\sum_i\,e_i\tau_i(y).}

As $\pi(e_i)=\delta_{1i}$, we have $\pi(v)=\pi(y)$. The \elt $v$ constructed thus is independent of the \sys of representatives for $G/H$. Indeed, if $(\tau'_i)$ is another \sys of representatives, even if it means reordering the indices, we can assume that $\tau'_iH = \tau_iH$, and so, $y$ being $H$-invariant, $\tau'_i(y) = \tau_i(y)$.\\
For $\sigma\in G$, the $(\sigma \circ \tau_i)$'s form a \sys of representatives for~$G/H$,
so~$\sigma(v) = v$: the \elt $v$ is $G$-invariant.

\emph{\ref{thADG1Idm5}.} 
We have a \dcn $\gC=\gC'_1\oplus\cdots \oplus \gC'_m$, where $\gC'_j = e_j\gC$ is a~\ptf \kmo of rank $r$ and the restriction $\pi : \gC'_1 \to \gC_1$ is an \iso of \kmosz. For all $y \in \gC$, we~have
{
$$\preskip-.0em \postskip-.4em\mathrigid 2mu 
\rC{\gC/\gk}(y)(T) = \prod\limits_{j=1}^m \rC{\gC'_j/\gk}(e_jy)(T)  \; \hbox{ and }\;
\rC{G}(y)(T) = \prod\limits_{j=1}^m \prod\limits_{\tau\in H} \big(T - (\tau_j\circ \tau)(y)\big).
 $$
}

Let $y$ be the unique \elt of $\gC'_1$ such that $\pi(y) = z$. The \egt on the left-hand side gives 

\snic{\rC{\gC/\gk}(y)(T) = T^{(m-1)r} \rC{\gC_1/\gk}(z)(T).}

Next,
apply $\pi$ to the \egt on the right-hand side by letting $(\tau_j\circ\tau)(y) \in \gC'_j$ (use $y = e_1y$ and apply $\tau_j\circ\tau$). We then obtain 

\snic{\rC{G}(y)(T) = T^{(m-1)r} \rC{H}(z)(T).}

Hence $\rC{\gC_1/\gk}(z)(T) = \rC{H}(z)(T)$.

\emph{\ref{thADG1Idm6}.} 
Taking into account the fact that the restriction of $\pi$ to $e_1\gC$ is an \iso we have $g'^2=g'=g'e_1$.
Similarly for $\sigma\in H$ we~have $\sigma(g')=g'$ if~$\sigma\in K$, or $g'\sigma(g')=0$ if $\sigma\notin K$.
Finally, for $\sigma\in G\setminus H$, $e_1\sigma(e_1)=0$, and so~$g'\sigma(g')=0$.
This shows that $g'$ is a Galoisian \idm of $\gC$ with stabilizer~$K$.
The canonical \iso is immediate.

\emph{\ref{thADG1Idm7}.} Item \emph{4} implies that $\gk$ is the set of fixed points. It remains to see that $H$ is separating.
If $\sigma\in H=\St(e)$ is distinct of the \idtz, we have some \elts $a_i$ and $x_i\in\gC$ such that $\sum_ia_i(\sigma(x_i)-x_i)=1$.
This \egt remains true if we localize at $e$.
\end{proof}

\subsect{Case where the \agB of a \uvl \dcn \alg is discrete}{Case where the \agB of a \adu is discrete}
\label{subsec2AGBDiscretes}

It is desirable that one can test the \egt of two \idmsz\hbox{ $e_1$, $e_2$} in the \adu $\gA$, 
which is the same as knowing how to test~$e=0$ for an arbitrary \idm of $\gA$ (as in every additive group).
However, $e\gA$ is a \ptf \kmo and $e=0$ \ssiz $\rR {e\gA}(X)=1$
(\thref{th ptf sfio} item~\emph{\iref{remRang3}}).
Since the \polmu $\rR {e\gA}$ can be explicitly computed, we can test the \egt of two \idms in $\gA$ \ssi we can test the \egt of two \idms in $\gk$. The above argument works in a slightly more \gnl framework and we obtain the following result.

\begin{fact}
\label{factagbdisc}
If $\BB(\gk)$ is a discrete \agBz, so is~$\BB(\gA)$.
More \gnlt, if  $\gC$ is a \stfe \klgz, and if $\BB(\gk)$ is discrete, then $\BB(\gC)$ is discrete.
\end{fact}

\begin{fact}
\label{factIdemStable}
If $(\gk,\gC,G)$ is a \apGz, every \idm $e$ of~$\gC$ fixed by $G$ is an \elt of~$\gk$.
\end{fact}
\begin{proof}
The \polcar $\rC {G}(e)=(T-e)^{|G|}$ belongs to $\kT$, so its constant \coez, which is equal to $\pm e$, is in~$\gk$.
\end{proof}

\vspace{-.7em}
\pagebreak	

\begin{fact}
\label{factconnexe}
Let $(\gk,\gC,G)$ be a \apG with $\gk$ connected and nontrivial, then
\begin{enumerate}
\item \label{factconnexe.1}  $0$ and $1$ are the only \idms of $\gC$ fixed by $G$,
\item \label{factconnexe.2}$\BB(\gC)$ is discrete,
\item \label{factconnexe.3}every atom of $\BB(\gC)$ is a Galoisian \idmz,
\item \label{factconnexe.4}two atoms are conjugated under $G$,
\item \label{factconnexe.5}an \idm $e\neq0$ is Galoisian \ssi its orbit under $G$ is formed of pairwise \ort \eltsz,
\item \label{factconnexe.6}if $f$ is an \idm $\neq0$, the \id $\gen{1-f}$ is Galoisian \ssi its orbit under $G$ is formed of pairwise \com \idsz.
\end{enumerate}
\end{fact}
\begin{proof}
Items~\emph{\ref{factconnexe.1}} and \emph{\ref{factconnexe.2}} clearly result from Facts~\emph{\ref{factIdemStable}} and~\emph{\ref{factagbdisc}.}

\emph{\ref{factconnexe.3}.} If $e$ is an atom, so is $\sigma(e)$, therefore $\sigma(e) = e$ or $e\sigma(e) = 0$.
Thus two \elts of the orbit of $e$ are \ortz, so the sum of the orbit of $e$ is a nonzero \idm fixed by $G$; it is equal to~$1$.

\emph{\ref{factconnexe.4}.} If $e'$ is another atom, it is equal to the sum of the $e'e_i$'s, 
where $e_i$ ranges over the orbit of $e$, and since the $e_i$'s are atoms, each $e_ie'$ is zero or equal to $e_i$.

\emph{\ref{factconnexe.5}.} See Fact~\ref{factAduIdmA}.

\emph{\ref{factconnexe.6}.} Stems from \emph{\ref{factconnexe.5}} since $\gen{1-f,1-f'}=\gen{1-ff'}$ for \idms $f$ and~$f'$.
\end{proof}

Theorem~\ref{factagb} implies that the \agB $\BB(\gC)$ is finite \ssi the in\dcps \idms form a finite set (they are \ncrt pairwise \ortsz) and if they generate $\BB(\gC)$.

\medskip
\comm \label{ensborn}
A set $X$ is said to be \emph{bounded} if we know an integer $k$ which is an upper bound of the number of elements in $X$, \cad more \prmtz, if for every finite family $(b_i)_{i\in\lrb{0..k}}$ in $X$, we have $b_i=b_j$ for two distinct indices.
In \clama this implies that the set is finite, but from a \cof point of view many distinct situations can occur.%
\index{bounded!set}%
\index{set!bounded ---}
\\
A common situation is that of a bounded and discrete \agB $C$ for which we do not know of an atom with certainty.
The \itfs of $C$, all principal, are identified with \elts of $C$, so $C$ is identified with its own Zariski lattice $\Zar C$.{\footnote{For a commutative \ri $\gk$, $\Zar\gk$ is the set of radicals of \itfs of $\gk$ (Section~\ref{secZarAcom}). It is a distributive lattice. In \clamaz, $\Zar\gk$ is identified with the lattice of quasi-compact open sets of the spectral space  
$\Spec\gk$ (Section~\ref{secEspSpectraux}).}}
Moreover, in \clama the atoms are in bijection with the prime \ids (all maximal) of $C$ via $e\mapsto \gen{1-e}$. Thus the set of atoms of $C$ (supposed bounded) is identified with $\Spec C$. We therefore once again find in this special case the following \gnl fact: the Zariski lattice is the \covz, practical and \gui{point-free} version of the Zariski 
spectrum, a topological space whose points can turn out to be inaccessible from a \cof point of view. But this situation, although familiar, is perhaps more troubling in the case of a discrete and bounded topological space. This is typically a compact space for which we do not have a good description via a dense countable subset, therefore which is not included in the context of compact metric spaces \`a la Bishop (cf.~\cite{B67,BB85}).
\eoe

\medskip
Here is a corollary of the Galois structure \tho (\thref{lemAduIdmA}) in the context of \apGsz.

\begin{proposition}
\label{corlemAduIdmA}
Let $(\gk,\gC,G)$  be a \apG  with $\gk$ connected. For an \idm $h$ of $\gC$
\propeq
\begin{enumerate}
\item $h$ is a Galoisian \idmz.
\item $\gC[{1/h}]$ is a \pro \kmo of rank equal to $\St_G(h)$.
\item $\gC[{1/h}]$ is a \pro \kmo of rank less than or equal to $\St_G(h)$.
\end{enumerate}
\end{proposition}
\begin{proof}
We use \thref{lemAduIdmA}. By item~\emph{\ref{lemAduIdmA-2}} of this \tho we can assume that there exists a Galoisian \idm $e$ such that $h$ is equal to a sum $e_1+\cdots +e_r$ of \elts of the orbit $G.e$.  We have \isos of
\kmos $e\gC\simeq \gC[{1/e}]$ and $\gC\simeq (e\gC)^{\abs{G.e}}$, so $e\gC$ is \prc $\idg{G:G.e} = \abs{\St_G(e)}$. We deduce that the \kmo 
$$
\gC[{1/h}]\simeq h\gC = e_1\gC\oplus\cdots\oplus e_r\gC\simeq (e\gC)^r
$$
is \pro of rank $r\times \abs{\St_G(e)}$. We then apply item~\emph{\ref{lemAduIdmA-5}} of \thref{lemAduIdmA} with
$f = h$.
\\
Therefore, here item \emph{2} (resp.\ item \emph{3}) means the same thing as item~\emph{5c} (resp.\ item~\emph{5d}) in \thref{lemAduIdmA}.
\end{proof}

\subsec{Discriminant
}
\label{subsecDiscrim}

Recall that in $\gA=\Adu_{\gk,f}$ we have $\disc(f)= \prod\nolimits_{1\leq i<j\leq n}(x_i-x_j)^2$ and \smashtop{$\Disc\iAk =\disc(f)^{n!/2}$}.

In the following \thoz, we speak of the $\gA$-\mdi $\Om{\gk}{\gA}$ of the \klg $\gA$.
It actually suffices to know that the \mdi of a \apf is \isoc to the cokernel of the transpose of the Jacobian matrix of the \syp that defined the \algz. 
For more details on this subject see \thrfs{thDerivUniv}{thDerivUnivPF}.

\begin{theorem}\label{lemAdu1}
 Let $J$ be the Jacobian of the \sys of $n$ equations with $n$ unknowns defining the \adu $\gA=\Adu_{\gk,f}$.
\begin{enumerate}
\item
\begin{enumerate}
\item \label{i1lemAdu1}
We have $J=\prod_{1\leq i<j\leq n}(x_i-x_j)$ in $\gA$.
\item \label{i2lemAdu1}
 We have $J^2=\disc(f)\in\gk$.
\end{enumerate}
\item In particular, \propeq
\begin{enumerate}\itemsep0pt
\item \label{i3alemAdu1}
  $\Disc\iAk $ is \iv (resp.\,\ndzz) in  $\gk$.
\item  $\disc(f)$ is \iv (resp.\,\ndzz) in  $\gk$.
\item  $J$ is \iv (resp.\,\ndzz) in  $\gA$.
\item The $x_i-x_j$'s are \ivs (resp.\,\ndzsz) in $\gA$.
\item \label{i3elemAdu1}
 $x_1-x_2$ is \iv  (resp.\,\ndzz) in $\gA$.
\item \label{i3flemAdu1}
 $\Om{\gk}{\gA}=0$ (resp.\,$\Om{\gk}{\gA}$ is a \gui{torsion} \Amoz, \cad annihilated by a \ndz \eltz).
\item \label{i4lemAdu1}
 $\Sn$ is a separating group for $\gA$ (resp.\,for  $\Adu_{\Frac(\gk),f}$).
\end{enumerate}

%
\item \label{i5lemAdu1} The analogous \eqvcs are valid for every Galois quotient of the \aduz.

\end{enumerate}
 \end{theorem}
\begin{proof}
Item \emph{1a} is easy by \recu on $n$, with the exact sign if we consider the \sys which we used in the \dfn of the \aduz.
For example, here is the computation for $n=4$
{
$$\preskip.4em \postskip-.4em\ndsp\mathrigid 2mu 
\begin{array}{rclcl}
\arraycolsep .01em \mathrigid 2mu
J&=   &
\dmatrix{1&1&1&1\cr
\som_{i\neq 1}x_i&\som_{i\neq 2}x_i & \som_{i\neq 3}x_i&
\som_{i\neq 4}x_i  \cr
\som_{i,j\neq 1}x_ix_j&\som_{i,j\neq 2}x_ix_j & \som_{i,j\neq 3}x_ix_j&
\som_{i,j\neq 4}x_ix_j  \cr
x_2x_3x_4&x_1x_3x_4 &  x_1x_2x_4 &x_1x_2x_3  }
\\[8mm]
\phantom{J}&  = &
\Dmatrix{3pt}{1&0&0&0\\[.5mm]
{\som_{i\neq 1}x_i} & x_1-x_2 & x_1-x_3&x_1-x_4  \cr
{\sum\limits_{i,j\neq 1}x_ix_j} & (x_1-x_2){\sum\limits_{i\neq 1,2}x_i} &
(x_1-x_3){\sum\limits_{i\neq 1,3}x_i} & (x_1-x_4){\sum\limits_{i\neq 1,4}x_i}
\cr
x_2x_3x_4&(x_1-x_2) x_3x_4 &  (x_1-x_3)x_2x_4 &(x_1-x_4)x_2x_3  }
   \\[10mm]
\phantom{J}& = &(x_1-x_2)(x_1-x_3)(x_1-x_4)
\dmatrix{
 1 & 1&1  \cr
x_3+x_4 &
x_2+x_4 & x_2+x_3  \cr
 x_3x_4 &  x_2x_4 &x_2x_3  }
\end{array} 
$$
}

etc\ldots

We deduce from it item \emph{1b}, then the \eqvc of items \emph{2a} through \emph{2e}.

\emph{2f.} Since $\Om{\gk}{\gA}$ is an \Amo \isoc to the cokernel of the transpose of the Jacobian matrix, we obtain that $\Ann(\Om{\gk}{\gA})$ and $J\gA$ have the same nilradical (Lemma~\ref{fact.idf.ann}). Finally, the \elt $J$ is \ndz (resp.\,\ivz) \ssi the \id $\sqrt{J\gA}$ contains a \ndz \elt 
(resp.\,contains~$1$).

\emph{2g.} Suppose that $f$ is \spl (resp.\,\ndzz), if $\sigma\in \Sn$ is distinct from $\Id_\gA$, there is some $i\in\lrbn$ such that $x_{\sigma i}\neq x_i$.
Since $x_{\sigma i}- x_i$ is \iv (resp.\,\ndzz),~$\sigma$ is separating (resp.\,separating once we invert the \discriz).
For the converse, consider for example the transposition~$\sigma$ that swaps~$1$ and~$2$. We clearly have $\gen{g-\sigma(g)\vert g\in\gA}=\gen{x_1-x_2}$. The result follows.

\emph{3.} Clear since the \adu is always \isoc to a power of any of its Galois
quotients.
\end{proof}

\vspace{-.7em}
\pagebreak	

\subsec{Fixed points}
\label{subsecptsfix}

Let $\di(f)=\prod_{i<j\in\lrbn}(x_i+x_j)\in\gk$.

It is clear that $\di(f)$ is congruent to $\prod_{i<j\in\lrbn}(x_i-x_j)$ modulo $2$,
which gives \smashtop{$\gen{2,\di(f)^2}=\gen{2,\disc(f)}$}.

\begin{theorem}
\label{theoremAdu1} \emph{(\Uvl splitting \alg and fixed points)}\\
Let $\fa:= \Ann_{\gk}\big(\geN{2,\di(f)}\big)$.
Then 
${\Fix(\Sn)\subseteq\gk+\fa\gA.}$
In particular, if $\fa=0$ and a fortiori if $\Ann_{\gk}\big(\geN{2,\disc(f)}\big)=0$, we obtain
$\Fix(\Sn)=\gk$.
\end{theorem}
\begin{proof}
It suffices to prove the first assertion. \\
Let us consider the case where  $n=2$ with $f(T) = T^2 - s_1T + s_2$. 
\\
An \elt $z=c+dx_1\in\gA$ (with $c$, $d\in\gk$) is invariant under $\mathrm{S}_2$  \ssi $d(x_1-x_2)=d(s_1-2x_1)=0$, or yet again if
$ds_1=2d=0$, but we~have~$\di(f)=s_1$.
\\
We then proceed by \recu on $n$.
For the Cauchy modules we use the notations of Section~\ref{sec0adu}.
For $n>2$ consider the \riz~$\gk_1=\gk[x_1]\simeq\aqo{\gk[X_1]}{f(X_1)}$ and the \pol $g_2(T)=f_2(x_1,T) $ which is in $\gk_1[T]$. We identify  $\Adu_{\gk_1,g_2}$ with $\Adu_{\gk,f}$ (Fact~\ref{factAduAdu}). 
To switch from the expression of an \elt $y\in\gA$ over the basis $\cB(g_2)$ ($\gA$ seen as a $\gk_1$-module) to its expression over the basis $\cB(f)$ ($\gA$ seen as a \kmoz), it suffices to express each coordinate, which is an \elt of $\gk_1$, over the $\gk$-basis $(1,x_1,\ldots,x_1^{n-1})$ of $\gk_1$.
Let us also take note that $\di(f)=(-1)^{n-1}g_2(-x_1)\di(g_2)$ by a direct computation.
Therefore, if we let $\fa_1=\Ann_{\gk_1}(\gen{2,\di(g_2)})$, we obtain $\fa_1\gA\subseteq \fa\gA$ and~$\fa_1\subseteq \fa\gk_1$.
\\
Let us move on to the actual \recuz.
\\
Let $y\in \gA$ be a fixed point of $\Sn$, and consider it as being an \elt of the \adu  $\Adu_{\gk_1,g_2}$.
Since $y$ is invariant under~$\mathrm{S}_{n-1}$, we have 
$y\in\gk_1+\fa_1\gA$, and so $y\equiv h(x_1) \mod \fa_1\gA$ for some $h\in\kX$.
A fortiori $y\equiv h(x_1) \mod \fa\gA$.
It remains to see that $h(x_1)\in\gk+\fa\gA$.
Since~$y$ is invariant under~$\Sn$, by permuting~$x_1$ and $x_2$ we obtain the congruence

\snic{\qquad\qquad\qquad h(x_1)\equiv y \equiv h(x_2)
\quad\mod \fa\gA\qquad\qquad\qquad(*)}

Let $h=\sum_{i=0}^{n-1}+c_{i}X^{i}\in\kX$.
Note that $h(x_1)$ is a reduced expression over the canonical basis $\cB(f)$. Regarding $h(x_2)$, to obtain the reduced expression, we must, in the term $c_{n-1}x_2^{n-1}$, replace $x_2^{n-1}$ with its expression over the canonical basis, which results from $f_2(x_1,x_2)=0$.
\\
This rewriting sparks the apparition of the term $-c_{n-1}x_1^{n-2}x_2$, and this implies by~$(*)$ that~$c_{n-1}\in \fa$. But then, $h(x_2)-c_{n-1}x_2^{n-1}$ and $h(x_1)-c_{n-1}x_1^{n-1}$ are reduced expressions of two \elts equal modulo $\fa\gA$.
Therefore, the $c_i$'s for $i\in\lrb{1..n-2}$ are in $\fa$, and we saw that $c_{n-1}\in\fa$.
\end{proof}

\rem
In the $n=2$ case, the above study shows that as soon as~$\fa\neq 0$ the \riz~$\Fix(\mathrm{S}_2)= \gk\oplus
\fa\,x_1=\gk+\fa\gA$ strictly contains $\gk$.
\\
A computation in the~$n=3$ case gives the same converse: if $\fa\neq0$, the \riz~$\Fix(\mathrm{S}_3)$ strictly contains~$\gk$. We indeed find an \elt

\snic{v=x_1^2 x_2+  s_1 x_1^2 + (s_1^2+s_2) x_1 + s_2 x_2\neq 0}

(one of its coordinates over $\cB(f)$ is equal to $1$) such that~$\Fix(\rS_3)= \gk\oplus \fa\,v$.
However, for~$n\geq 4$, the situation becomes complicated.
\eoe

\smallskip 
We obtain as a corollary the following \thoz.

\begin{theorem}\label{thAduAGB} 
If $f$ is a \spl \pol of $\kT$, the \adu $\Adu_{\gk,f}$, as well as every Galois quotient, is a \aGz.
\end{theorem}

\begin{proof}
By the structure \tho \ref{thADG1Idm} (item \emph{7}) it suffices to show that $\Adu_{\gk,f}$ is Galoisian. However, we have just proven the fixed point condition, and the separating automorphisms condition is given in \thref{lemAdu1}.
\end{proof}

By Artin's \tho \ref{thA}, and in light of the previous \thoz, we know that every \adu for a \spl \polz, or every Galois quotient of such a \klgz, diagonalizes itself. We examine this question in further detail in the following subsection.  
Even with regard to the precise result that we have just mentioned, it is interesting to see things work \gui{concretely} for a \aduz.

\subsec{Separability}
\label{subsecsep}

When the \pol $f\in \gk[T]$ is \splz, its \adu $\gA = \Adu_{\gk,f} = \gk[\xn]$ is \stez, by Fact~\ref{factDiscriAdu}. The following \tho then simply recalls \thref{prop2EtaleReduit} regarding \ases in the current context.

\begin{theorem}
\label{theoremAdu2} Suppose $f$ is \splz.
\begin{enumerate}
\item
The nilradical $\rD_\gA(0)$ is the \id generated by $\rD_\gk(0)$. In particular, if $\gk$ is reduced, so is $\gA$.
\item 
For every reduced \alg $\gk\vers{\rho}\gk'$, the \alg $\rho\ist(\gA)\simeq\Adu_{\gk',\rho(f)}$ is reduced.
\end{enumerate}
\end{theorem}
\subsubsection*{Diagonalization of a \aduz}

\begin{theorem}
\label{theoremAdu3} \emph{(Diagonalization of a \aduz)} \\
Let $\varphi:\gk\to\gC$ be an \alg in which \emph{$f$ can be completely factorized}, \cad $\varphi(f)=\prod_{i=1}^n(T-u_i)$.
Also suppose that $f$ is \spl over~$\gC$, \cad the $u_i-u_j$'s are \iv for $i \ne j$.\\
Let $\gC\te_\gk\gA\simeq\Adu_{\gC,\varphi(f)}$,
and, for $\sigma\in\Sn$, let $\phi_\sigma: \gC\te_\gk\gA \to\gC$ be the unique
\homo of \Clgs which sends each $1_\gC\otimes x_i$ to $u_{\sigma i}$.
\\
Let $\Phi: \gC\te_\gk\gA\to\gC^{n!}$ be the $\gC$-\homo defined by $y\mapsto\big(\phi_\sigma(y)\big)_{\sigma\in\Sn}$.
\begin{enumerate}
\item 
$\Phi$ is an \isoz: $\gC$ diagonalizes $\gA.$
\item 
More \prmtz, in $\gC\te_\gk\gA$, write $x_i$ instead of $1_\gC\te x_i$, $u_i$ instead of $u_i\te 1_{\gA}$ (in accordance with the  \Clg structure of $\gC\te_\gk\gA$) and let $g_\sigma=\prod_{j\neq \sigma i} (x_i-u_j)$.  
Then, 

\snic{\phi_\sigma(g_\sigma)=\pm \varphi\big(\disc(f)\big)=\pm\disc\big(\varphi(f)\big),
}


and $\phi_\sigma(g_\tau)=0$ for $\tau\neq \sigma$, so that
if we let $e_\sigma=g_\sigma/\phi_\sigma(g_\sigma)$, the~$e_\sigma$'s form the \sfio corresponding to the \isoz~$\Phi$.
\item Moreover, $x_ie_\sigma=u_{\sigma i}e_\sigma$, so that the basis $(e_\sigma)$ of the \Cmo $\gC\te_\gk\gA$ is a common diagonal basis for the multiplications by the~$x_i$'s.
\end{enumerate}
In particular, when $f$ is \splz, the enveloping \alg

\snic{\env\gk\gA=\gA\otimes_\gk\gA\simeq\Adu_{\gA,f}}

is canonically \isoc to $\Ae {n!}$; $\gA$ diagonalizes itself.
\end{theorem}
NB: We will however be careful when letting $\Adu_{\gA,f}=\gA[u_1,\ldots,u_n]$ since the $x_i$'s are already taken as \elts of~$\gA$.
\begin{proof}
\emph{1.}
The two \algs are, as \Cmosz, \isoc to $\gC^{n!}$ and $\Phi$ is a \Cli whose surjectivity is all that we need to prove.
The surjectivity results by the Chinese remainder \tho from the $\Ker\phi_\sigma$'s being pairwise \comz: $\Ker\phi_\sigma$ contains $x_i-u_{\sigma i}$,  $\Ker\phi_\tau$ contains $x_i-u_{\tau i}$, 
therefore  $\Ker\phi_\sigma+\Ker\phi_\tau$ contains the $u_{\sigma i}-u_{\tau i}$'s, and there is at least one index $i$ for which $\sigma i\neq \tau i$, which shows that $u_{\sigma i}-u_{\tau i}$ is \ivz.

\emph{2.} The \sfio corresponding to the \iso $\Phi$ is the unique solution of the \sli
$\phi_\sigma(e_\tau)=\delta_{\sigma,\tau}$ (where $\delta$ is the Kronecker symbol). \\
However, the \egts $\phi_\sigma(g_\sigma)=\pm \varphi\big(\disc(f)\big)$ and $\phi_\sigma(g_\tau)=0$ are easy. 

\emph{3.} Fix $i$. The \egt $x_ig_\sigma=u_{\sigma i}g_\sigma$ results from the fact that in $g_\sigma$ there is already a product of the $x_i-u_j$'s for $j\neq \sigma i$, so $(x_i-u_{\sigma i})g_\sigma$ is a multiple of $\varphi(f)(x_i)$, which is null. 
\end{proof}

\rem 
Actually, \gnlt speaking, $\Phi$ is a \ali whose \deter can be computed with respect to the natural bases: the square of this \deter is a power of $\varphi\big(\disc(f)\big)$ and we therefore find that~$\Phi$ is an \iso \ssi $\varphi\big(\disc(f)\big)$ is \iv in~$\gC$. For this, and for a \gui{complete converse,} see Exercise~\ref{exoIdentiteDiscriminantale}.
\eoe

\perso{il ne semble pas que l'on  arrive \`a faire d\'ecouler le
\thrf{theoremAdu2} of the \ref{theoremAdu3}.}

\medskip The previous \tho implies the following result: if $\gA$ is a \adu for a \spl \polz, then every \Alg diagonalizes~$\gA$.
We now give a \gnn of this result for a Galois quotient of $\gA$.

\pagebreak	

\begin{theorem}
\label{theoremAdu4} \emph{(Diagonalization of a Galois quotient of a \aduz)}
Let $e$ be a Galoisian \idm of $\gA$, 

\snic{\gB=\aqo{\gA}{1-e}=\gk[\yn]\;\hbox{ and }\;G=\St_\Sn(e),}

(we have denoted by $y_i=\pi(x_i)$ the class of $x_i$ in $\gB$).
Let $\phi:\gB\to\gC$ be a \ri \homoz. Let $u_i=\phi(y_i)$.
Consider the \Clg  

\snic{\phi\ist(\gB)\simeq\gC\otimes_\gk\gB\simeq\aqo{\Adu_{\gC,f}}{1-\phi(e)}}

obtained from the \klg $\gB$ by \edsz.
For~$\sigma\in G$ let $\phi_\sigma:\gC\otimes_\gk\gB\to\gC$ be the unique \homo of \Clgs which sends each $1_\gC\otimes y_i$ to $u_{\sigma i}$.
Let $\Phi:\gC\otimes_\gk\gB \to\gC^{|G|}$ be the \homo of~\Clgs defined by $z\mapsto\big(\phi_\sigma(z)\big)_{\sigma\in G}$.
\begin{enumerate}
\item If $\phi\big(\disc(f)\big)\in\gC\eti$, $\Phi$ is an \isoz, so $\gC$ diagonalizes $\gB$.
\item In particular, if $f$ is \splz, $\gB\otimes_\gk\gB$ is  canonically \isoc to~$\gB^{|G|}$, 
\cad $\gB$ diagonalizes itself.
\end{enumerate}
\end{theorem}
\begin{proof}
The two $\gC$-\algs are \pro \Cmos of constant rank $|G|$ and $\Phi$ is a \Cli whose surjectivity is all that we need to prove.
In $\gC\otimes_\gk\gB$ we write $y_i$ instead of $1_\gC\otimes y_i$ and $u_i$ instead of $u_i\otimes 1_{\gB}$.
The surjectivity results by the Chinese remainder \tho from the fact that the $\Ker\phi_\sigma$'s are pairwise \comz: $\Ker\phi_\sigma$ contains $y_i-u_{\sigma i}$,  $\Ker\phi_\tau$ contains $y_i-u_{\tau i}$, so $\Ker\phi_\sigma+\Ker\phi_\tau$ contains the $u_{\sigma i}-u_{\tau i}$. However, there is at least one index $i$ for which $\sigma i\neq \tau i$ and $u_{\sigma i}-u_{\tau i}$ is \iv because $\phi\big(\disc(f)\big)$ is the product of the $(u_j-u_k)^2$'s for $1\leq j<k\leq n$.
\end{proof}

\subsec{Triangular structure of Galoisian \idsz}
\label{subsecidGTri}

In this subsection we prove \thref{thidGTri} which implies that the structure of the \id $\cJ(f)$, which is a \gui{triangular} structure (in the Lazard sense) when we consider the Cauchy modules as \gtrsz, remains a triangular structure for all the Galoisian \ids of the \adu in the case of a \spl \pol over a \cdiz.

\begin{lemma}\label{lemBaseDisc}
Let $\gk'$ be a \klg which is a \mptf of constant rank $m$, $x \in \gk'$ and $r(T) \in \gk[T]$ be the \polcar of~$x$ over~$\gk$. 
If $\disc(r) \in \gk\eti$, then $\gk' = \gk[x]$ and $(1, x, \ldots, x^{m-1})$ is \hbox{a $\gk$}-basis of $\gk'$. 
\end{lemma}

\begin{proof}
The case where $\gk'$ is free of rank $m$ has been proven in \ref{propdiscTra}. In the \gnl case, consider a \sys of \eco of $\gk$ such that each \lon makes of $\gk'$ a free \kmo of rank $m$.
\end{proof}

\vspace{-.5em}
\pagebreak

\begin{theorem}
\label{thidGTri}
Let  $(\gk,\gC,G)$ be a \aG with
\begin{itemize}
\item   $\gC=\kxn\simeq\kXn\sur\fa$,
\item  $G$ operates on $\{\xn\}$ and
\item  the $x_i-x_j$'s are \ivs for $i\neq j$.
\end{itemize}
A typical example of this situation: $\gC$ is a Galois quotient of the \adu of a \spl \polz.
\\
Let

\snic{
\begin{array}{l} 
G_i=\sotQ{\sigma\in G}{\sigma(x_k)=x_k,\, k\in \lrb{1..i}}
\hbox{ for }i \in \lrb{0..n} \hbox{ (so } G_0 = G \hbox{)}, \\[1mm]
 r_i(T)=\prod\nolimits_{\sigma\in G_{i-1}/G_i}\big(T-\sigma(x_i)\big)
\quad  \hbox{ for }  i \in\lrbn,
\end{array}
}

where $G_{i-1}/G_i$ designates a \sys of representatives of the left cosets. Let $d_i=\idg{G_{i-1}:G_i}$.
\\
We then have the following results.
\begin{enumerate}
\item  $\gk[x_1,\ldots ,x_i]=\Fix(G_i)$ and $G_i=\Stp(\gk[x_1,\ldots ,x_i]).$
\item  The \pol $r_i(T)$ is \mon with \coes in $\gk[x_1, \ldots, x_{i-1}]$, of degree $d_i$.
Let $R_i(X_1,\ldots ,X_i) \in \gk[X_1,\ldots ,X_i]$ be a \polu in $X_i$ of degree $d_i$ such that $R_i(x_1,\ldots ,x_{i-1},X_i)=r_i(X_i)$.
\item  
$\fa_i=\fa\cap\gk[X_1,\ldots ,X_i]$ is generated by $R_1(X_1)$, $\ldots$, $R_i(X_1,\ldots ,X_i)$.
\end{enumerate}
Consequently each \alg $\gk[x_1,\ldots ,x_i]$ is both a free $\gk[x_1,\ldots ,x_{i-1}]$-module of rank $d_i$ and a free \kmo of rank $\idg{G:G_i}$ at the same time, and each of the \ids $\fa_i$ is a triangular \id (in the Lazard sense) of $\gk[X_1,\ldots ,X_i]$.
\end{theorem}
\begin{proof}
The group $G_1$ is a separating group of \autos of \riz~$\gC$. Let $\gk_1=\gC^{G_1}$. We know that $\gC$ is a \pro
$\gk_1$-module of constant rank $|G_1|$ and that $\gk[x_1]\subseteq\gk_1$. Moreover, $\gk_1$ is a direct summand 
in~$\gC$, therefore it is a \kmrc $d_1=\deg_T(r_1).$
\\
The \id $\fa_1$ is formed by all the $R\in\gk[X_1]$'s such that $R(x_1) = 0$.\\
Therefore, $R\big(\sigma(x_1)\big)=0$ for every $\sigma\in G/G_1$. In other words $R$ is a multiple of each~$T- \sigma(x_1)$, and since the $x_i-x_j$'s are \ivsz, $R$ is a multiple of~$r_1$.  Thus $\fa_1=\gen{r_1(X_1)}$ and $\gk[x_1] \simeq \aqo{\gk[X_1]}{r_1(X_1)}$.
\\
Proposition~\ref{propAdiagAH} gives us the \egt

\snic{\rC{\gk_1/\gk}(x_1)(T) = \prod_{\sigma \in G/G_1}\big(T - \sigma(x_1)\big)=r_1(T).}

This implies that the \polcar $\rC{\gk_1/\gk}(x_1)(T)$ is \splz, and Lemma~\ref{lemBaseDisc} says that $(1,x_1,\dots,x_1^{d_1-1})$ is a basis of $\gk_1$.
\\
Thus $\gk[x_1]=\gk_1=\Fix(G_1)$ and $(\gk[x_1],\gC,G_1)$ is a \aGz. 

Then, $\gC=\gk_1[x_2,\ldots ,x_n]$ with $G_1$ operating on $\{x_2,\ldots ,x_n\}$ and the \hbox{$x_i-x_j$'s} being \ivsz.
The whole previous process works identically when replacing $\gk$ by $\gk_1$, $G$ by $G_1$, $x_1$ by $x_2$ and $G_1$ by $G_2$. 
The result then follows by \recuz.
\end{proof}

\vspace{-.7em}
\pagebreak

\section{Splitting field of a \pol over a \cdi}
\label{subsecCDR}

In this section we give a \cov and dynamic approach to the \cdr of a \polu over a \cdiz, in the absence of a factorization \algo of the \polsz.

\Grandcadre{
In Section~\ref{subsecCDR},  $\gK$ is a nontrivial \cdiz,\\
  $f$ is a \polu of degree $n$ and\\
  $\gA=\Adu_{\gK,f}=\KXn\sur{\cJ(f)}=\Kxn$.
}

\smallskip The quotients of the \adu $\gA$ are finite \Klgsz, therefore they are \zed \risz.

\subsect{\gui{Reduced} Galois quotients of the \aduz}{Good quotients of the \aduz}

We have placed quotation marks around \gui{reduced} because, a priori, one does not speak of a Galois quotient \emph{that is} reduced, but of a Galois quotient \emph{that one} reduces (by killing the nilpotents).

Given Fact~\ref{factDiscriAdu}, if the \pol $f$ is \spl the \adu is \'etale, therefore reduced, and every \id generated by an \idm is equal to its nilradical (since the quotient \ri is reduced).
We can then replace in the statements that follow each \id $\DA(1-e)=\sqrt{\gen{1-e}}$ with the \id $\gen{1-e}$.

\smallskip 
In the following lemma we know by hypothesis that $\gB$ is \stfe over~$\gK$, but we do not \ncrt know a basis of~$\gB\red$ as a~\Kevz.
The goal is then to give a \gui{satisfying enough} description of~$\gB\red$ as a quotient of the \aduz.

\begin{lemma}
\label{lemAquQuoRed}
Let $\gB$ be a \stfe \Klgz. Suppose that $f$ can be entirely decomposed in $\gB\red$ and that $\gB\red$ is generated by the corresponding zeros of $f$. Then, there exists an \idm $e$ of $\gA=\Adu_{\gK,f}$ such that $\gB\red\simeq\gA\sur{\DA(1-e)}$.
\end{lemma}
\begin{proof}
Let $y_1$, \dots, $y_n\in\gB$ such that $f(T)=\prod_i(T-\ov{y_i})$ in $\gB\red$.  There exists a unique \homoz~$\lambda :\KXn\to\gB$ which sends the~$X_i$'s to the~$y_i$'s. 
Let $\fb$ be the (\tfz) \id of $\gB$ generated by $\lambda\big(\cJ(f)\big)$. We then have~$\fb\subseteq \DB({0})$, and $\gB' :=\gB\sur{\fb}$ is a strictly finite \Klg satisfying~$\gB\red\simeq \gB'\red$. We thus obtain a diagram 

\snic {
\xymatrix @R = 0.4cm {
\KuX\ar[d]_{\lambda}\ar@{->>}[r] &\gA\ar[d]^{\varphi}\ar@{->>}[dr]^{\psi} 
\\
\gB \ar[r]         &\gB'\ar[r]   &\;\gB'\red = \gB\red 
\\
}}

in which $\varphi$ is the unique \homo which sends $x_i$ to the class of $y_i$. 
Since $\gB'$ is \stfez, $\Ker\varphi$ is a \itf of $\gA$, and there exists a $d \ge 0$ such that $(\Ker\varphi)^d = (\Ker\varphi)^{d+1}$ therefore $(\Ker\varphi)^d$ is generated by an \idm $1-e$. The result follows because on the one hand, $\psi$ is surjective, and on the other, $\Ker\psi = \DA(\Ker \varphi) = \DA\big((\Ker \varphi)^d\big) = \DA(1-e)$.
\end{proof}

\rem
Note that $\psi$ is surjective, but a priori $\Ker\psi$ is not a \itf of $\gA$. Symmetrically, a priori $\varphi$ is not surjective, but $\Ker\varphi$ is a \itf of $\gA$.
\eoe

\medskip 
In \clama a \cdr (\Dfref{defCorpsdesRacines}) for a \polu $f$ over a \cdi $\gK$ is obtained as a quotient of the \adu $\gA$ by a \idemaz. Such an \id exists: take a strict \id that is a~\Kev of maximal dimension, by \TEMz. 

In \coma we obtain the following more precise \tho (to be compared with \thref{propUnicCDR}).

\begin{theorem}
\label{thUnici} ~
\begin{enumerate}
\item \label{i1thUnici}
\Propeq
\begin{enumerate}
\item \label{i1athUnici} There exists in $\gA=\Adu_{\gK,f}$ an in\dcp \idmz~$e$.
\item \label{i1bthUnici} There exists an extension $\gL$ of $\gK$ that is a \cdr of $f$ and  denoted by $\gB\red$ where $\gB$ is a \stfe \Klgz.
\item \label{i1cthUnici} The \agB $\BB(\gA)$ is finite.
\end{enumerate}
\item \label{i2thUnici} In this case every \cdr of $f$ is \isoc to $\gA\sur{\DA(1-e)}$, and it is discrete.
\end{enumerate}
\end{theorem}
\begin{proof}
The \eqvc of \emph{1a} and \emph{1c} is valid in the \gnl context of \agBs (\thref{lemAduIdmA}). 
It is clear that \emph{1a} implies \emph{1b.} Conversely if we have a \cdrz~$\gL=\gB\sur{\DB(0)}$, where~$\gB$ is a~\stfe \Klgz, Lemma~\ref{lemAquQuoRed} provides an \idm $e$, and it is in\dcp because $\gL$ is connected.
\\
Let us look at item \emph{\ref{i2thUnici}}. Let $\gM$ be a \cdr for~$f$.
Write

\snic{f(T)=\prod_{i=1}^n(T-\xi_i)$ in $\gM.}

By the \uvl \prt of $\Adu_{\gK,f}$, there exists a unique \homo of~\Klgs $\varphi :\gA\to\gM$ such that $\varphi (x_i)=\xi_i$ for $i\in\lrbn$.  Let $(e_\ell)_{\ell=1,\ldots ,k}$ be the orbit of $e$. It is a \sfioz, so $\big(\varphi(e_\ell)\big)_{\ell=1,\ldots ,k}$ \egmtz, and since $\gM$ is a \cdi this implies that there is some $j$ for which $\varphi(e_\ell)=\delta_{j,\ell}$ (Kronecker symbol).
\\
Then, $\gen{1-e_j}\subseteq\Ker\varphi$, so $\gM$ is a quotient of $\gA\sur{\DA(1-e_j)}$, which is a \cdiz.  As $\gM$ is nontrivial, this implies $\gM\simeq \gA\sur{\DA(1-e_j)}$.  Finally, all the $\gA\sur{\DA(1-e_\ell)}$ are pairwise \isocz.
\end{proof}

\comm
In \cite{MRR}, it is shown that every enumerable \cdi possesses an \agq closure. However, a \cdr for $f$, which therefore exists, does not \ncrt possess a finite basis as a \Kevz, in the \coma sense, and we do not know of a \cof uniqueness \tho for such a \cdrz.
We can describe as follows an analogous procedure to that of \cite{MRR} to obtain a \cdr for $f$.
First of all we construct an enumeration~$(z_m)_{m\in\NN}$ of the \aduz.
Next we construct a sequence of \itfs $(\fa_m)$ of $\gA$ by letting $\fa_0=0$, and $\fa_{m+1}=\fa_m+\gen{z_m}$ if $\fa_m+\gen{z_m}\neq \gen{1}$, and $\fa_{m+1}=\fa_m$ otherwise (the test works because we can compute a basis for the \Kev $\fa_m+\gen{z_m}$). 
Then,
the \id $\bigcup_m\fa_m$ is a \idema of $\gA$, and the quotient is a \cdrz, which is discrete. Our point of view is slightly different. We do not a priori start from an enumerable field, and even in the case of an enumerable field, we do not favor one enumeration over another. We would rather answer questions about the splitting field as they arise, as we shall see in the following \thoz.
\eoe  
           
\medskip
The following \tho  
explains how to bypass the difficulty of the nonexistence of the \cdr in \comaz. The \cdr of $f$ is replaced by an \gui{approximation} given in the form of a reduced quotient $(\aqo{\gA}{1-e})\red$ of the \aduz, where $e$ is a Galoisian \idmz.
\\
We rely on the following fact which is already established in the \gnl context of \zed \ris (Lemma~\ref{lemme:idempotentDimension0}).
We recall a direct \demz.

\emph{For all $y\in\gA=\Adu_{\gK,f}$, there exists an \idm $e_y\in\gK[y]\subseteq\gA$ such that~$y$ is \iv
modulo $1-e_y$ and nilpotent modulo $e_y$.\label{factIdmAduPolmin}
}
\begin{proof}
Let $P(T)$ be the \polmin of $y$. There exists an \iv \eltz~$v$ of $\gK$
such that $vP(T)=T^k\big(1-TR(T)\big)$ with $k\geq0$. The \idm $e_y$ is~$\big(yR(y)\big)^k$.
\end{proof}
%

\begin{theorem}
\label{thdivzeridm} \emph{(Dynamic management of a \cdrz)}\\
Let $(z_i)_{i\in I}$ be a finite family of \elts of $\Adu_{\gK,f}=\gA$. There exists a Galoisian \idm $e$ of $\gA$ such that by letting $\gB=\aqo{\gA}{1-e}$ each $\pi(z_i)$ is null or \iv in the quotient \alg $\gB\red$ (here, $\pi:\gA\to\gB\red$ is the canonical projection).
\end{theorem}
\begin{proof}
For each $i\in I$ there is an \idm $g_i\in\gA$ such that $z_i$ is \iv modulo
$1-g_i$ and nilpotent modulo $g_i$.  Applied to the family of the $g_i$'s \thref{lemAduIdmA} gives a Galoisian \idm $e$, such that for each~$i$, $1-e$ divides $g_i$ or $1-g_i$. Therefore, in the quotient \alg $\gB=\aqo{\gA}{1-e}$
 each~$\pi(z_i)$ is nilpotent or \ivz.
\end{proof}

\rems 1) 
The reader may worry that we do not a priori have a finite \sgr of the \id $\DA(1-e)$ available. Consequently the finite \alg $\gB\red$ is not \ncrt a finite dimensional \Kev in the \cof sense. Actually the nilpotents can also be managed dynamically. We have in $\gB=\aqo{\gA}{1-e}$ a test of nilpotence and if a nilpotent \elt $x$ is revealed,
       we can replace $\gB$ with its quotient by the \idz~$\fa$
generated by the orbit of $x$ under the action of $G=\St_\Sn(e)$. Then, $\gB\sur{\fa}$ is finite dimensional and $G$ operates on~$\gB\sur{\fa}$.

 2) In \thref{thdivzeridm} it can be in our best interest to saturate the family $(z_i)_{i\in I}$ by the action of $\Sn$ in order to make manifest in $\gB$ all the possible \gui{scenarios.}
\eoe

\subsec{Uniqueness of the \cdr}

The uniqueness \tho of the \cdr admits an \gui{operative} \cov version 
        (which always works, 
even if we do not dispose of an in\dcp \idm in the \aduz) in the following form.

\begin{theorem}
\label{propUnicite} \emph{(Uniqueness of the \cdrz, dynamic version)}\\
Let $\gB_1$, $\gB_2$ be two nonzero \stfes \Klgs for which the \pol $f$ can be decomposed into a product of \lin factors in $(\gB_1)\red$ and $(\gB_2)\red$. Moreover suppose that $(\gB_i)\red$ is generated by the corresponding zeros of $f$. Then, there exists a Galoisian \idm $e$ of $\gA$ such that, with the \alg $\gB=\aqo{\gA}{1-e}$, we have two integers $r_i$ such that $(\gB_i)\red \simeq \gB\red^{r_i}$.
\end{theorem}
\begin{proof}
Lemma~\ref{lemAquQuoRed} gives \idms $e_1$, $e_2 \in \gA$ such that

\snic{(\gB_i)\red \simeq \gA\sur{\DA(1-e_i)}\quad (i=1,\,2)}

\Thref{lemAduIdmA} item~\emph{2} gives a Galoisian \idm $e$ \hbox{and $r_1$, $r_2\in\NN$}
such that $\aqo\gA{1-e_i}\simeq \gB^{r_i}$. Therefore $(\gB_i)\red\simeq \gB\red^{r_i}.$
\end{proof}

\incertain{

\subsubsection*{And le group des \autosz?}

\entrenous{Paragraphe incertain. Il serait bon of faire subir
\`a un nilpotent arbitraire in une \apG un traitement of the m\^eme style que
celui qui a \'et\'e fait aux \idms in le \thref{lemAduIdmA}.
En attendant voici quelques petites choses with l'\adu qui semblent
raisonnables.

Mais peut \^etre le fin mot of l'affaire is plut\^ot of the cot\'e des \gui{\fcns \splsz} qui sont by example tr\`es bien trait\'ees in \cite{MRR} chapter 6.}
Lorsque $f$ n'est pas \splz, $\disc\,f=0$ and $z=x_1-x_2$
est un \elt non \iv of~$\gA$.

On consid\`ere then l'\idm $e_z\neq 1$ of~$\gA$.

If $e_z=0$ tous les $x_i-x_j$ sont nilpotents. En passant au quotient
par l'\id qu'ils engendrent, on obtient une \Klg \isoc to
$\gK[x_1]\simeq\aqo{\gK[X_1]}{f(X_1}$.

Alors, le \cdr of $f$
est \isoc to $\gK[x_1]\red$. C'est une extension purement 
radicielle of $\gK$?????????

\smallskip If $e_z\neq 0$ on consid\`ere un \idm galoisien
$e$ correspondant (\thrf{lemAduIdmA} and \algo \ref{algidmgal}).
In le quotient of Galois $\aqo{\gA}{1-e}$, certains des $\pi(x_i-x_j)$
sont nilpotents, of autres ($\pi(z)$ by example) \ivsz.
On passe au quotient by l'\id generated by les $\pi(x_i-x_j)$ nilpotents.
Let $\gC$ la \Klg obtenue. On remplace  $G=\St(e)$ by le group de
permutation
des $x_i$ \gui{restants} and l'on obtient un group  $H$
de $\gK$-\autos of $\gC$. Let l'\alg $\gL=\gC^H$.

Then  $(\gL,\gC,H)$ is une \aGz, et
$\gL$ is une \gui{approximation} of une extension purement radicielle de
$\gK$, plus \prmt $\gL\red$ is un corps, c'est une extension
purement radicielle of $\gK$?????

}

\section{Galois theory of a \spl \pol over a \cdi} 
\label{secThGB}

\medskip 
\Grandcadre{
In Section~\ref{secThGB},  $\gK$ is a nontrivial \cdiz,
\\  $f$ is a \spl \polu of degree $n$ and $\gA=\Adu_{\gK,f}$.\\
We highlight the fact that $f$ \emph{is not assumed to be irreducible}.
}

\medskip 
Recall that for a \sply factorial field, every \spl \pol has a splitting field (Corollary~\ref{propIdemMini}), unique up to \auto (\thref{propUnicCDR}).
We are now interested in the case where the field \emph{is not} \sply factorial (or even in the case where the \fcn of the \spl \pols is too costly).

Here, as promised, we give the \cov and dynamic version of the Galois theory of a \spl \pol over a \cdiz.

\subsect{Existence and uniqueness of the dynamic and static \cdrsz}{Existence and uniqueness of the \cdrz}

Fact~\ref{factDiscriAdu} (or Corollary~\ref{corcorlemEtaleEtage}) assures us that $\gA$ is an \'etale \Klgz.
The same goes for its Galois quotients.
\Thref{thUnici} can be re-expressed as follows.

 
\THO{thUnici}
{\emph{(Separable \polz: when a \cdr exists and is a \stfe extension)}
\label{thExistence}
\begin{enumerate}
\item \label{thExistence1}
\Propeq
\begin{enumerate}
\item \label{thExistence1a} There exists in $\gA=\Adu_{\gK,f}$ an in\dcp \idm $e$.
\item \label{thExistence1a'}  There exists a \stfe extension  $\gL$ of $\gK$ that is a \cdr of $f$.
\item \label{thExistence1c} The \agB $\BB(\gA)$ is finite.
\end{enumerate}
\item \label{thExistence2}
In this case every \cdr of $f$ is a Galois extension of~$\gK$, \isoc to~$\gA[{1/e}]$.
\end{enumerate}

}

\medskip
Item~\emph{\ref{thExistence2}} also results from the fact that if
a \cdr exists and is \stf over $\gK$, two \cdrs are \isoc (\thref{propUnicCDR}).

The uniqueness \tho \ref{propUnicite} can be re-expressed as follows.

\THO{propUnicite} {\emph{(Uniqueness of the \cdr of a \spl \polz, dynamic version)}
Given two nonzero \stfes \Klgs $\gB_1$ and~$\gB_2$ in which $f$ can be decomposed into a product of \lin factors and which are generated by the corresponding zeros of $f$, there exists a Galois quotient $\gB=\gA[1/e ]$ of the \adu and two integers $r_i$ such that
$\gB_1\simeq\gB^{r_1}$ and~$\gB_2\simeq\gB^{r_2}$.
}

\subsect{Structure of the Galois quotients of the \adu}{Galois quotients of the \adu}

For the remainder of Section~\ref{secThGB} we fix the following notations.
\begin{notations}
\label{notas1.3}  \emph{(Context of a Galois quotient)}\\
Let $e$ be a Galoisian \idm of $\gA=\Adu_{\gK,f}$, $\fb=\gen{1-e}_\gA$.
Let

\snic{\gB=\gA\sur\fb=\gA[{1/e}] ,\;  \pi=\pi_{\gA,\fb}: \gA\to\gB,\;\hbox{ and }\;G=\St_\Sn(e).}

Let $(e_1,\ldots ,e_m)$ be the orbit of $e$ under $\Sn$. Each \Klg
$\gA[{1/e_i}]$ is \isoc to~$\gB$. The group $G$ operates on~$\gB$.  
\end{notations}

Note that for $y\in\gB$, the \pol $\Mip_y(T)$ is \spl (because $\gB$ is \'etale over $\gK$). In addition $y$ is \iv \ssi $\Mip_y(0)\neq 0$.  Also note that a \itf of $\gB$ (different from $\gen{1}$) is a Galoisian \id \ssi its orbit under $G$ is formed of pairwise \com \ids (every \itf is generated by an \idmz, and Fact~\ref{factconnexe}).

The structure \tho \ref{thADG1Idm} reads as follows, taking into account \thrfs{thdivzeridm}{thAduAGB}.

\begin{theorem}
\label{thStruc3} \emph{(Galois structure \thoz, 3)}\\
 In the context of~\ref{notas1.3} we obtain the following results.
\begin{enumerate}
\item \label{thADG3Idm3} \label{thADG3Idm2}
$(\gK,\gB,G)$ is a Galois quotient of $(\gK,\gA,\Sn)$. \\
In particular, $\gB$ is a finite dimensional \Kev  $\abs{G}$ and for every $y\in\gB$,
$\rC {\gB/\gK}(y)(T)=\rC G(y)(T)$. In addition, $\Fix(G)=\gK$.
\item \label{thADG3Idm1} 
We have an \iso of \Klgs $\gA\simeq \gB^m$.
\item \label{thADG3Idm4} If $\gB$ is connected, it is a \cdr for $f$ and a Galois extension of $\gK$ with $G$ as its Galois group.
\item \label{thADG3Idm8} 
Let $(y_i)$ be a finite family  of \elts of $\gB$.  There exists a Galoisian \idm $e_\gB$ of $(\gK,\gB,G)$ such that in $\gB[{1/e_\gB}]$, each $y_i$ is null or \ivz.
\item \label{thADG3Idm7} 
The restriction $\pi : e\gA \to \gB$ is a $\gK$-\lin \iso and establishes a biunivocal correspondence 
between the Galoisian \idms of $(\gK,\gA,\Sn)$ contained in $e\gA$ and those of $(\gK,\gB,G)$.
The stabilizers and residual quotients are preserved; 
\cad if $e_\gA \in e\gA$ and $e_\gB \in \gB$ are two Galoisian \idms that correspondent to each other, then $\St_\Sn(e_\gA) = \St_G(e_\gB)$ and 
$\gA[{1/e_\gA}] \simeq \gB[{1/e_\gB}]$. 
\end{enumerate}
\end{theorem}
NB: In what follows, we only give the statements for the relative situation, the absolute situation is indeed the special case where $e=1$.


\begin{lemma}
\label{Thidgal} \emph{(Resolvent and \polminz)} 
  Context~\ref{notas1.3},   $y\in\gB$.
\begin{enumerate}
\item \label{Thidgal5a} $\Rv_{G,y}(T)$ has \coes in $\gK$.
\item \label{Thidgal5b} $\Mip_y$ divides $\Rv_{G,y}$ which divides a power of $\Mip_y$.
\item \label{Thidgal5c} $\rC {\gB/\gK}(y)(T)= \rC G(y)(T) = \Rv_{G,y}(T)^{|\St_G(y)|}$.
\end{enumerate}
\end{lemma}
\begin{proof}
\emph{\ref{Thidgal5a}.} Consequence of item~\emph{\ref{thADG3Idm3}} in the structure \thoz. 
 
\emph{2.} We deduce that  $\Mip_{y}$ divides $\Rv_{G,y}$, because $\Rv_{G,y}(y)=0$, and since each~$y_i$ annihilates $\Mip_y$, the product of the $T-y_i$'s divides a power of~$\Mip_y$.
 
\emph{\ref{Thidgal5c}.} The second \egt is obvious, and the first is in item~\emph{\ref{thADG3Idm3}} of the structure \thoz.
\end{proof}

\subsec{Where the computations take place} 

Recall that $f$ is a \spl \polu of $\KT$ with $\gK$ a nontrivial \cdiz.

\Grandcadre{
We denote by $\gZ_0$ the sub\ri of $\gK$ generated by the \coes of $f$ \\
and by $1/\!\disc(f)$. We denote by $\gZ$ the \cli of $\gZ_0$ in $\gK$.
}

Here we highlight the fact that \gui{all the computations take place, and all the results are written, in the \ri $\gZ$,} since this follows from \thrfs{thIdmEtale}{thidGTri}.\footnote{It follows that if $\gK$ is a \gnl field (see Section~\ref{secAloc1}), the  questions of computability are actually discussed entirely in $\Frac(\gZ)=\Frac(\gZ_0)\te_{\gZ_0}\gZ={(\gZ_0\sta)^{-1}}\gZ$,
and $\Frac(\gZ)$ is discrete if $\gZ_0$ is itself a discrete \riz.
As $\gZ_0$ is a \tf \riz, it certainly is, in \clamaz, an effective (also called computable) \ri with an explicit \egt test, in the sense of recursion theory via Turing machines. 
\\
But this last result is not a truly satisfying approach to the reality of the computation. It is indeed akin to results in classical mathematics of the form \gui{every recursive real number admits a recursive development into a continued fraction}, a theorem that is evidently false from a practical point of view, since to implement it, one must first know whether the number is rational or not.}
\\
These \thos give us in the current framework items~\emph{1},~\emph{2} and~\emph{4} of the following \thoz.
As for item~\emph{3}, it is an \imde consequence of item~\emph{2.}

\begin{theorem}\label{thZsuffit} \emph{(The sub\ri $\gZ$ of $\gK$ is sufficient)} 
\begin{enumerate}
\item Let $\gZ_1$ be an intermediate \ri between $\gZ$ and $\gK$ (for example $\gZ_1=\gZ$).\\ 
Then the \adus 

\snic{\Adu_{\gZ_0,f}\subseteq\Adu_{\gZ,f}\subseteq\Adu_{\gZ_1,f}\subseteq\Adu_{\gK,f}}

are \aGs (with respect to their base \risz, and with the group $\Sn$).
\item Every \idm of $\gA=\Adu_{\gK,f}$ is in $\Adu_{\gZ,f}$: its coordinates over the basis $\cB(f)$ are in $\gZ$.
\item  The \emph{Galois theories of $f$} over $\gZ$, over $\gZ_1$, over $\Frac(\gZ)$ and over $\gK$ are \emph{identical}, in the following sense.
\begin{enumerate}
\item Every Galois quotient of $\Adu_{\gZ_1,f}$ is obtained by \eds to $\gZ_1$ from a Galois quotient of $\Adu_{\gZ,f}$.
\item Every Galois quotient of $\Adu_{\Frac(\gZ),f}$ is obtained by \eds to $\Frac(\gZ)$ from a Galois quotient of $\Adu_{\gZ,f}$.
This \eds is in fact the same thing as passing to the total \ri of fractions of the Galois quotient.
\item Every Galois quotient of $\Adu_{\gK,f}$ is obtained by \eds to $\gK$ from a Galois quotient of $\Adu_{\Frac(\gZ),f}$.
\end{enumerate} 
\item Let $e$ be a Galoisian \idm of $\gA$ and $\gZ_1$ be the sub\ri of $\gZ$ generated by $\gZ_0$ and the coordinates of $e$ over $\cB(f)$.
Then, the triangular structure of the \id of $\gZ_1[\Xn]$ generated by $1-e$ and the Cauchy modules is made explicit by means of \pols with \coes in $\gZ_1$.\\
For those that know \bdgsz: the \bdg (for a lexicographical monomial order) of the \id that defines the corresponding approximation of the \cdr of $f$ is formed of \polus with \coes in $\gZ_1$.
\end{enumerate}

\end{theorem}

Note the simplifications in the following special cases. If $\gK=\QQ$ 
and $f \in \ZZ[T]$ \monz, then $\gZ=\gZ_0=\ZZ[1/\!\disc(f)]$. Similarly, for $q$  a prime power and $\gK$ the field of rational fractions $\gK = \FF_q(u)$ we have, if~$f\in \FF_q[u][T]$ is \monz, $\gZ=\gZ_0=\FF_q[u][1/\!\disc(f)]$.

\medskip
\rems 1) Experiments  suggest not only that \gui{$\gZ$ is sufficient,} but that in fact all the results of computations (\coes of an \idm over $\cB(f)$, a \bdg of a Galoisian \idz) only use as \denos \elts whose square divides the \discri of $f$. \perso{Help! La th\'eorie des nombres
ne pourrait elle pas venir \`a la rescousse?}

 2) \emph{Absolute Galois theory of a \polz.}
Given a \spl \pol $f\in\KT$, rather than considering $\gK$ and the \cli $\gZ$ of~$\gZ_0$ in $\gK$,
we can consider $\gK'=\Frac(\gZ_0)$ and the \cli $\gZ'$ of~$\gZ_0$ in $\gK'$. \eoe

\subsec{Changing the base \riz, modular method}

Since everything takes place in $\gZ$, one can look at what happens after an arbitrary \eds $\varphi:\gZ\to \gk$.

It is possible for example that $\gk$ is a \gui{simple} \cdi and that we know how to compute $\Gal_\gk\big(\varphi(f)\big)$;
\cad identifying an in\dcp \idmz~$e'$ in $\Adu_{\gk,\varphi(f)}$.
This group will \ncrt be (\isoc to) a subgroup of the unknown Galois group $\Gal_\gZ(f)=\Gal_\gK(f)$.

Suppose that we have computed a Galois quotient $\gB$ of $\Adu_{\gZ,f}$ with a group $G\subseteq\Sn$. 

If $e$ is the Galoisian \idm of $\Adu_{\gZ,f}$ corresponding to $\gB$, we can reduce back  to the case 
where $\varphi(e)$ is a sum of conjugates of $e'$ and where 
$$\preskip.2em \postskip.4em 
\Gal_\gk\big(\varphi(f)\big)= H=\St_G(e')\subseteq G. 
$$
As this is true for every Galois quotient of $\Adu_{\gZ,f}$, we obtain a double inclusion
\begin{equation}\preskip-.4em \postskip.4em
\label{eqmetMod} H\subseteq \Gal_{\gZ,f}\subseteq G
\end{equation}
 except that the group $\Gal_{\gZ,f}$ is only defined up to 
conjugation,
and that it can a priori remain forever unknown.

This type of information, \gui{the Galois group of $f$ over $\gK$, up to conjugation, contains $H$} is outside of the dynamic method that we have presented, because this one takes a step in the other direction: giving information of the type \gui{the Galois group of $f$ over $\gK$, up to conjugation, is contained in~$G$.}

It is therefore a priori interesting to use the two methods in parallel, in the hope of completely determining $\Gal_\gK(f)$.

Replacing the field $\gK$ by a sub\ri is important from this point of view as we dispose of a lot more morphisms of \eds from the domain $\gZ$  than from the domain $\gK$.

\smallskip In particular, it is often useful to work modulo $\fp$, a \idema of $\gZ$.
This method is called a \emph{modular method}.

This method seems to have been invented by Dedekind for the determination of the Galois group of $f$ over $\QQ$ when $f\in\ZZ[T]$.
Note that in this case a \idema of $\gZ=\gZ_0=\ZZ[1/\disc(f)]$ is given by a prime number $p$ which does not divide $\disc(f)$.

\subsec{Lazy Galois theory}

The structure \tho \vref{thStruc3} and Lemma~\ref{Thidgal} (which gives a few details) are the theoretical \cof results that allow a lazy \evn of the \cdr and of the Galois group of a \spl \polz.

Please note that the term \gui{lazy} is absolutely not pejorative. It simply indicates that we can work with complete confidence in the splitting field of a \spl \pol over $\gK$, even in the absence of any \fcn \algo of the \pols over $\gK$. Indeed, any anomalies with the \alg $\gB$, the \gui{ongoing} approximation of the \cdr of $f$, for instance the presence of a nonzero zerodivisor, can be exploited to significantly improve our knowledge of the Galois group and of the \cdrz. A Galoisian \idm that is strictly a multiple of the \gui{ongoing} \idm $e$ can indeed be computed.
In the new \aGz, which is a quotient of the previous one, all the previously made computations remain valid, by passage to the quotient. Moreover, the number of significant improvements that may occur this way does not exceed the maximum length of a chain of subgroups of~$\Sn$.

We therefore develop a \gui{Galoisian} variant of the D5 system, which was the first computer algebra system to compute, both systematically and without risk of errors, in the \agq closure of a field in the absence of  a \fcn \algo for \pols (see \cite[Duval\&al.]{D5}).

Here, in contrast to what happens with the D5 system, the dynamic aspect of things does not consist in \gui{opening separate branches of computation} each time an anomaly occurs, but in improving the approximation of the \cdr (and of its Galois group) that constitutes the ongoing Galois quotient of the \adu each time.

\subsubsec{The basic \algoz}

We can rewrite \Algo \ref{algidmgal} in the current setting as follows, when we have an \elt $y$ neither null nor \iv in the Galois quotient $\gB=\gA\sur\fb$ at our disposal.

\begin{algoR}[Computation of a Galoisian \id and of its stabilizer]\label{alg2idmgal}
\Entree $\fb$: Galoisian \id of a \adu $\gA$ for a \spl \polz, $\fb$ is given by a finite \sgrz; $y$: nonzero \dvz  in 
$\gB=\gA\sur\fb$; $G=\St_\Sn(\fb)$;
$S_y=\St_G(y)$.
\Sortie $\fb'$: a Galoisian \id of $\gB$ containing $y$;
$H$:  $\St_G(\fb')$.
\Varloc $\fa$: \itf of $\gB$;
$\sigma$: \elt of $G$;
\\$L$: list of \elts of $\ov{G/S_y}$;
\hsu \# \, $\ov{G/S_y}$ is a \sys of representatives of the left cosets of $S_y$\\
$E$: corresponding set of \elts of $G/S_y$;
\hsu \# \, $G/S_y$  is the set of left cosets of $S_y$. 
\Debut
\hsu $E \aff \emptyset$; $L \aff [\,]$;  $\fb' \aff \gen{y}$;
\hsu \pur{\sigma}{\ov{G/S_y}}
\hsd $\fa \aff \fb' + \gen {\sigma(y)}$;
\hsd \sialors{1\notin\fa} $\fb' \aff \fa$;  
$L \aff L\bullet [\sigma]$; $E \aff E \cup \so{\sigma S_y}$
\hsd \finsi;
\hsu \finpour;
\hsu $H\aff \St_G(E)$
\quad \# \, $H = \sotq{\alpha\in G} 
  {\forall\sigma\in L, \alpha\sigma\in \bigcup_{\tau\in L}\tau S_y}$.
\Fin
\end{algoR}

The \id $\fb$ is given by a finite \sgrz, and $G=\St_\Sn(\fb)$.  Let~$e$ be the \idm of $\gB$ such that $\gen{1-e}_\gB=\gen{y}_\gB$, and $e'$ be a Galoisian \idm such that $G.e$ and $G.e'$ generate the same \agBz.  
\\
We are looking to compute the Galoisian \id $\fc$ of $\gA$ which gives the new Galois quotient $\gA\sur\fc \simeq \gB\sur{\fb'}$, where $\fb'=\gen{1-e'}_\gB$, \cad $\fc=\pi_{\gA,\fb}^{-1}(\fb')$.

In \Algo \ref{algidmgal} we find the product of $e$ and a maximum number of conjugates, avoiding obtaining a null product. 

Here we do not compute $e$, nor $\sigma(e)$, nor $e'$, because experimentation often shows that the computation of $e$ is very long (this \idm often occupies a lot of memory space, significantly more than $e'$). We then reason with the corresponding \ids $\gen{1-e}=\gen{y}$ and $\gen{1-\sigma(e)}=\gen{\sigma(y)}$. It follows that in the \algo the product of the \idms is replaced by the sum of the \idsz.

Moreover, as we do not compute $e$, we replace $\St_G(e)$ by $\St_G(y)$, which is contained in $\St_G(e)$, \gnlt strictly so. Nevertheless, experience shows that, even though $\ov{G/S_y}$ is larger, the whole computation 
is faster. We leave it to the reader to show that the last assignment in the \algo indeed provides the desired group $\St_G(\fb')$; \cad that the subgroup $H$ of~$G$, the stabilizer of $E$ in $G/S_y$, is indeed equal to $\St_G(\fb')$.

\subsubsec{When a relative resolvent factorizes}

Often an anomaly in a Galois quotient of the \adu corresponds to the observation that a relative resolvent can be factorized. We therefore treat this case in all \gnt to reduce it to a case where a nonzero zerodivisor is present.

\begin{proposition}
\label{propRvRel1} \emph{(When a relative resolvent factorizes)}\\
In the context of \ref{notas1.3} let $y\in\gB$ and $G.y=\so{y_1,\ldots,y_r}$.
\begin{enumerate}
\item \label{RvRenum5}  If $\Mip_y=R_1R_2$ with $R_1$ and $R_2$ of degrees $\geq 1$, $R_1(y)$ and $R_2(y)$ are nonzero zerodivisors, and there exists an \idm $e$ such that $\gen{e}=\gen{R_1(y)}$ and $\gen{1-e}=\gen{R_2(y)}$.
\item \label{RvRenum6}  If $\deg(\Mip_y)<\deg(\Rv_{G,y})$, then one of the $y_1-y_i$'s divides zero (we can therefore construct an \idm $\neq 0,1$ of $\gB$).
\item \label{RvRenum7}  If  $P$ is a strict divisor of $\Rv_{G,y}$ in $\gK[T]$, then at least one of the two following case occurs:
\begin{itemize}
\item  $P(y)$ is a nonzero zerodivisor, we are in item {\ref{RvRenum5}.}
\item an \elt ${y_1-y_i}$ is a nonzero zerodivisor, we are in item {\ref{RvRenum6}.}
\end{itemize}
\end{enumerate}
\end{proposition}
\begin{proof}
\emph{\ref{RvRenum5}.} Since $\Mip_y$ is \splz, $R_1$ and $R_2$ are \comz.
With a Bézout relation $U_1R_1+U_2R_2=1$, let $e=(U_1R_1)(y)$ and $e'=1-e$. We have~$ee'=0$, so~$e$ and~$e'$ are  \idmsz. We also \imdt have 

\snic{eR_2(y)=e'R_1(y)=0$,
$eR_1(y)=R_1(y)$ and $e'R_2(y)=R_2(y).
}

 Therefore  $\gen{e}=\gen{R_1(y)}$ and $\gen{1-e}=\gen{R_2(y)}$.

\emph{\ref{RvRenum6}.} The \dem that shows that over an integral \ri a \polu of degree $d$ cannot have more than $d$ distinct roots is reread as follows. 
\\
Over an arbitrary \riz, if a \polu $P$ of degree $d$ admits some $(a_1,\ldots ,a_d)$ as zeros with each
$a_i-a_j$ \ndz for $i \ne j$, then we~have~$P(T)=\prod(T-a_i)$. Therefore if $P(t)=0$ and $t$ is distinct from the $a_i$'s, 
at least two of the $t-a_i$'s are nonzero \dvzsz. 
We apply this to the \polmin $\Mip_y$ which has more zeros in $\gB$ than its degree (those are the $y_i$'s). This gives a nonzero \dvz $y_j-y_k$, and via some $\sigma\in G$ we transform $y_j-y_k$ into a $y_1-y_i$.

\emph{\ref{RvRenum7}.} If $P$ is a multiple of $\Mip_y$, we are in item \emph{\ref{RvRenum6}}.\\
 Otherwise, $\pgcd(\Mip_y,P)=R_1$ is a strict divisor of $\Mip_y$, and $R_1\neq1$ because we~have~$\pgcd\big((\Mip_y)^k,P\big)=P$ for large enough $k$. Therefore $\Mip_y=R_1R_2$, with~$\deg(R_1)$ and $\deg(R_2)\geq1$.
We are in item~\emph{\ref{RvRenum5}}.
\end{proof}

From this we deduce the following corollary which generalizes item \emph{\ref{thADG3Idm8}} of the structure \tho \vref{thStruc3}.

\pagebreak	

\begin{theorem}
\label{corRvRel1} In the context of \ref{notas1.3} 
let $(u_j)_{j\in J}$ be a finite family in $\gB$. There exists a Galoisian \id $\fc$ of $\gB$ such that, by letting $H=\St_G(\fc)$, $\gC=\gB\sur{\fc}$, and $\beta:\gB\to\gC$ be the canonical \prnz, we~have
\begin{enumerate}
\item Each $\beta(u_j)$ is null or \ivz.
\item In $\gC$, $\Mip_{\beta(u_j)}(T)=\Rv_{H,\beta(u_j)}(T)$.
\item The $\Mip_{\beta(u_j)}$'s are pairwise equal or \comz.
\end{enumerate}
\end{theorem}

\rem
In the previous \tho it is sometimes in our best interest to saturate the family $(u_j)_{j\in J}$ by the action of $G$ (or of $\Sn$ by lifting the $u_j$
to~$\gA$) in order to make manifest in $\gC$ all the possible \gui{scenarios.}
\eoe

\medskip
\exl We reuse the example of \paref{exemple1Galois}. We ask {\tt Magma} what it thinks about the \elt {\tt x5 + x6}. Finding that the resolvent is of degree~$15$ (without having to compute it) whilst the \polmin is of degree~$13$, it struggles to reduce the oddity and obtains a Galois quotient of the \adu of degree $48$ (the \cdr of degree $24$ is not yet reached) with the corresponding group. The computation is almost instantaneous.
Here is the result.

{\small\label{example2Galois}
\begin{verbatim}
y:=x5+x6;
MinimalPolynomial(y);
  T^13 - 13*T^12 + 87*T^11 - 385*T^10 +
     1245*T^9 - 3087*T^8 + 6017*T^7 - 9311*T^6 + 11342*T^5 
     - 10560*T^4 + 7156*T^3 - 3284*T^2 + 1052*T - 260
//new Galois algebra, computed from deg(Min)<deg(Rv) :
Affine Algebra of rank 6 over Rational Field
Variables: x1, x2, x3, x4, x5, x6
Quotient relations:
  x1 + x2 - 1,
  x2^2 - x2 + x4^2 - x4 + x6^2 - x6 + 3,
  x3 + x4 - 1,
  x4^4 - 2*x4^3 + x4^2*x6^2 - x4^2*x6 + 4*x4^2 - x4*x6^2 + x4*x6 -
        3*x4 + x6^4 - 2*x6^3 + 4*x6^2 - 3*x6 - 1,
  x5 + x6 - 1,
  x6^6 - 3*x6^5 + 6*x6^4 - 7*x6^3 + 2*x6^2 + x6 - 1
Permutation group acting we have set of cardinality 6
Order = 48 = 2^4 * 3
    (1, 2)
    (3, 5)(4, 6)
    (1, 3, 5)(2, 4, 6)
\end{verbatim}
}

Certain special cases of the situation examined in Proposition~\ref{propRvRel1} are used as exercises.
Each time the goal  is to obtain more precise information on what happens when we reduce the observed oddity.
See Exercises~\ref{propGaloisJordan}, \ref{propRvRelDec} and \ref{propRvRelDecMin}.

\subsubsec{When the triangular structure is missing}

Consider some \elts $\alpha_1,\ldots ,\alpha_\ell$ of $\gB$ and the nested \Klgs 
  $$\gK\subseteq\gK_1=\gK[\alpha_1]\subseteq\gK_2= \gK[\alpha_1,\alpha_2] \subseteq\cdots
\subseteq\gK_\ell= \gK[\alpha_1,\ldots ,\alpha_\ell]\subseteq \gB.$$

For $i=2,\ldots ,\ell$ the structure of $\gK_i $ as a $\gK_{i-1}$-module can be made explicit by different techniques.
If $\gB$ is a \cdr for $f$, all the $\gK_i$'s are fields and each of the modules is free.

If one of these modules is not free, then we can construct an \idmz~$\neq0,1$ in $\gB$
by using the same technique as for the \dem of the structure \tho \vref{th1Etale}, item \emph{2b}.

\smallskip
Using the Gr\"obner basis technique can turn out to be efficient, with the \id that defines $\gB$ as a quotient of $\KXn$.
We introduce some variable names $a_i$ for the $\alpha_i$'s and we choose a lexicographical order with $a_1< \cdots <a_\ell<X_1< \cdots <X_n$.

If $\gB$ is a field the \bdg must have a triangular structure. To each $\alpha_i$ must correspond one and only one \pol in the \bdgz, $P_i(a_1,\ldots, a_i)$ monic in $a_i$.

If this triangular structure is not respected for the variable $a_i$, we are certain that $\gK_{i-1}$ is not a field, and we can 
         explicitly construct a nonzero zerodivisor  
in this \Klgz.

Actually, let $P(a_1,\ldots ,a_i)$ be a \pol that appears in the \bdg and that is not monic in $a_i$. Its leading \coe as a \pol in $a_i$ is a \pol $Q(a_1,\ldots ,a_{i-1})$ which \ncrt gives a nonzero zerodivisor $Q(\alpha_1,\ldots ,\alpha_{i-1})$
in the \zede \alg $\gK_{i-1}\simeq \gK[a_1,\ldots, a_{i-1}]\sur{\fa}$,
where $\fa$ is the \id generated by the first \polsz, in the variables $a_1$, \ldots, $a_{i-1}$, that appear in the \bdgz.
Otherwise, we could multiply $P$ by the inverse of $Q$ modulo $\fa$,
and reduce the result modulo~$\fa$,
and we would obtain a \polu in~$a_i$ that precedes $P$ in the lexicographical ordering, and that would render the presence of~$P$ pointless.

 
\Exercices

\begin{exercise}
\label{exoChapGalLecteur}
{\rm  We recommend that the proofs which are not given, or are sketched, or left to the reader, etc, be done. But in particular, we will cover the following cases.
\begin{itemize}
\item \label{exopropBoolFini} Prove Propositions~\ref{defiBoole},
\ref{propBoolFini} and \thref{corpropBoolFini}.
\item Explain Fact~\ref{factChangeBase}.
\end{itemize}
}
\end{exercise}

\vspace{-1em}

\pagebreak	

\begin{exercise}
\label{exothSteCdiClass} (Structure of finite \algs over a field, classical version, dynamic \cov version)
\\
{\rm \emph{1.} Prove in \clama the following result.
\\
\emph{Every finite \alg over a field is a finite product of finite local \algsz. 
}
\\
\emph{2.} Explain why we cannot hope to obtain a \prco of this result, even if we assume that the field is discrete.
\\
\emph{3.} Propose a \cov version of the previous result.
 
}
\end{exercise}

\vspace{-1em}
\begin{exercise}
\label{exothNst1-zed}   
{\rm  Show that the \elgbm \num2  (\paref{MethodeZedRed})
applied to the \dem of \thref{thNstNoe} gives the following result, \eqv to \thref{thNstNoe} in the case of a nontrivial \cdiz.

\smallskip   \textbf{\Thref{thNstNoe} bis}
\label{thNst1-zed}  (Weak \nst and \iNoe position, \zedr \ris case) 
{\it Let $\gK$ be a \zedr \riz, $\ff$ be a \itf of $\KuX=\KXn$ and~$\gA=\KuX/\ff$ be the quotient \algz. Then,
there exists a \sfio $(e_{-1},e_{0},\ldots,e_{n})$ of $\gK$ and a \cdv such that, naming the new variables $Y_1,\ldots ,Y_{n}$, and letting 

\snic{\gK_r=\gK[1/e_r]\;\;\hbox{ and }\;\;\gA_r=\gA[1/e_r]=\gK_r\otimes_{\gK}\gA\simeq \gK_r[\uX]\sur{\ff\ \gK_r[\uX],}
}

 we have the following results.
\begin{enumerate}
\item  $\gA_{-1}=0$ and $\gK\cap \ff=e_{-1}\gK$.
\item  $\gA_0$ is a free $\gK_0$-module of finite rank $\geq 1$. 
\item  For $r=1$, $\ldots$, $n$ we have
\begin{itemize}
  \item [$\bullet$] $\gK_r[Y_1,\ldots , Y_{r}]\cap \ff=0$. In other words the \alg $\gK_r[Y_1,\ldots , Y_{r}]$ can be considered as a $\gK_r$-sub\alg of $\gA_r$.
  \item [$\bullet$] $\gA_r$ is a \mpf over $\gK_r[Y_1,\ldots ,Y_{r}]$.
  \item [$\bullet$] There exists an integer $N$ such that for each $(\alpha_1,\ldots,\alpha_r)\in\gK_r^r$, the $\gK_r$-\alg 
 
\snic{\gB_r=\aqo{\gA_r}{Y_1-\alpha_1,\ldots,Y_r-\alpha_r}
}

 is a quasi-free $\gK_r$-module of finite rank $\leq N$, and the natural \homoz~\hbox{$\gK_r\to \gB_r $} is injective.
\end{itemize}
\end{enumerate}
In particular, the \Klg $\gA$ is a \mpf over the \gui{\pollez} sub\alg $\gA=\bigoplus_{r=0}^n\gK_r[Y_1,\ldots ,Y_{r}]$.
We say that the \cdv (which eventually changes nothing at all) has put the \id in \emph{\Noe position}.\\
Finally, the \sfio that intervenes here does not depend on the \cdv that puts the \id in \Noe position.
}
}
\end{exercise}

\vspace{-1em}
\begin{exercise}\label{exoMatMag3}
{(Magic squares and commutative \algz)}
\\
{\rm  
In this exercise we provide 
 an application of commutative \alg to a combinatorial \pbz; the free \crc that intervenes in the \Noe positioning of question \emph {2} is an example of the Cohen-Macaulay \prt in a graded environment.   
 A \emph {magic square} of size $n$ is a matrix of $\Mn(\NN)$ for which the sum of each row and each column is the same. The set of these magic squares is an additive sub\mo of $\Mn(\NN)$; we will admit here that it is the \mo generated by the~$n\,!$ permutation matrices.
We are interested in counting the magic squares of size $3$ of fixed sum $d$. 
Here are the $6$ permutation matrices of $\MM_3(\NN)$

\snic {
\begin {array} {rcl}
P_1 = \cmatrix {1&0&0\cr 0&1&0\cr 0&0&1},\
P_2 = \cmatrix {0&0&1\cr 1&0&0\cr 0&1&0},\
P_3 = \cmatrix {0&1&0\cr 0&0&1\cr 1&0&0}
\\
\noalign {\smallskip}
P_4 = \cmatrix {0&1&0\cr 1&0&0\cr 0&0&1},\
P_5 = \cmatrix {0&0&1\cr 0&1&0\cr 1&0&0},\
P_6 = \cmatrix {1&0&0\cr 0&0&1\cr 0&1&0}
\\
\end {array}
}

They are linked by the relation $P_1+P_2+P_3 = P_4+P_5+P_6$. 
Let $\gk[(x_{ij})_{i,j \in \lrb{1..3}}]$ be the \pol \ri with nine \idtrs
where $\gk$ is an arbitrary \riz. 
We identify the matrix $M = (m_{ij}) \in \MM_3(\NN)$ with the \mom $\prod_{i,j} x_{ij}^{m_{ij}}$, denoted \smashtop{$\ux^M$};
for example~\smashtop{$\ux^{P_1} = x_{11}x_{22}x_{33}$}.

\emph {1.}
Let $U_1$, \ldots, $U_6$ be six \idtrs over $\gk$ and $\varphi : \gk[\uU] \twoheadrightarrow \gk[\ux^{P_1}, \ldots, \ux^{P_6}]$ defined by $U_i \mapsto \ux^{P_i}$. \\
We want to show that $\Ker\varphi$ is {the \id $\fa = \gen {U_1U_2U_3 - U_4U_5U_6}$}.
\begin {enumerate}\itemsep0pt
\item [\emph {a.}]
Show, for $a$, $b$, $c$, $d$, $e$, $f \in \NN$ and $m = \min(a,b,c)$, that

\snic {
U_1^a U_2^b U_3^c U_4^d U_5^e U_6^f \equiv
U_1^{a-m} U_2^{b-m} U_3^{c-m} U_4^{d+m} U_5^{e+m} U_6^{f+m} \bmod \fa
}

\item [\emph {b.}]
Let $\fa^\bullet$ be the \ksmo of $\gk[\uU]$ with the \moms not divisible by $U_1U_2U_3$ as its basis. 
Show that $\gk[\uU] = \fa\oplus\fa^\bullet$ and that $\Ker\varphi = \fa$.
\item [\emph {c.}]
Deduce that the number $M_d$ of magic squares of size $3$ and of sum~$d$ is equal to
${d+5 \choose 5} - {d+2 \choose 5}$ 
since the convention is that
${i \choose j} = 0$ for $i < j$.
\end{enumerate}

\emph {2.}
Let $\gB = \gk[\uU]\sur\fa = \gk[\uu]$. 
\vspace{-1pt}
\begin{enumerate}\itemsep0pt
\item [\emph {a.}]
Define a \Noe position $\gA = \gk[v_2,v_3,u_4,u_5,u_6]$ of $\gB$ \hbox{where $v_2$, $v_3$} are \lin forms in $\uu$, such that $(1, u_1, u_1^2)$ is an $\gA$-basis \hbox{of $\gB = \gA \oplus \gA u_1 \oplus \gA u_1^2$}.

\item [\emph {b.}] 
Deduce that the number $M_d$ is also equal to ${d+4 \choose 4} + {d+3 \choose 4} + {d+2 \choose 4}$ (MacMahon's formula, which in passing gives an \idt between binomial \coesz).
\end{enumerate}

\emph {3.}
Suppose that $\gk$ is a \cdiz. We want to show that the \ri $\gB$, regarded as the \ri $\gk[\ux^{P_1}, \ldots, \ux^{P_6}]$, is \icl (see also \Pbm \ref{exoFullAffineMonoid}).  Let $E \subset \MM_3(\ZZ)$ be the \ZZsmo of magic squares (analogous \dfnz) and the sub\ri $\gB_{11} \subset \gk[x_{ij}^{\pm 1}, i,j \in \lrb{1..3}]$

\snic {
\gB_{11} = \gk[\ux^{P_1}, \ux^{P_6}]
[\ux^{\pm P_2}, \ux^{\pm P_3}, \ux^{\pm P_4}, \ux^{\pm P_5}]
}

such that $\gB_{11} \subset \gk[\,\ux^M \,\vert\, M \in E,\ m_{11} \ge 0\,]$.
\vspace{-1pt}
\begin{enumerate}\itemsep0pt
\item [\emph {a.}]
Verify that $\gB$ and $\gB_{11}$ have the same quotient field, which is the quotient field $\gk(E)$, the field of rational fractions over $\gk$ with $5$ \idtrsz.
\item [\emph {b.}]
Show that $\gB_{11}$ is \iclz.
\item [\emph {c.}]
For $i$, $j \in \lrb{1..3}$, define a \ri $\gB_{ij}$ analogous to $\gB_{11}$ and deduce that $\gB$ is \iclz.
\end{enumerate}

}

\end {exercise}

\vspace{-1em}
\begin{exercise}
\label{exoNonCyclique}
{\rm  Give a direct \dem (not using a reductio ad absurdum) that if a \cdi has two \autos that generate a noncyclic finite group, the field contains some $x\neq 0$, all the powers of which are pairwise distinct, 
\cad it is not a root of unity.
}
\end{exercise}

\vspace{-1em}
\begin{exercise}\label{exoIdentiteDiscriminantale}
{(A \gui {discriminantal} \idtz)}\\ 
{\rm  
Let $n \ge 1$. Let $E$ be the set of $\alpha \in \NN^n$
such that $0 \le \alpha_i < i$ for $i \in \lrbn$; it is a set of cardinality $n!$ which we order by \gui{factorial numeration,} \cad $\alpha \preceq \beta$ if $\sum_i \alpha_i i! \le \sum_i \beta_i i!$. We order the \smq group $\rS_n$ by the lexicographic ordering, $\In$ being the smallest permutation. Consider $n$ \idtrs over $\ZZ$ and define a matrix $M \in \MM_{n!}(\ZZ[\ux])$, indexed by $\rS_n \times E$,

\snic {
M_{\sigma, \alpha} = \sigma(\ux^\alpha), \qquad 
\sigma \in \rS_n, \quad \alpha\in E.
}

Thus for $n = 3$:
$$\preskip-.20em \postskip.3em 
M = \cmatrix {
1& x_2& x_3& x_2x_3&  x_3^2& x_2x_3^2\cr
1& x_3& x_2& x_2x_3&  x_2^2& x_2^2x_3\cr
1& x_1& x_3& x_1x_3&  x_3^2& x_1x_3^2\cr
1& x_3& x_1& x_1x_3&  x_1^2& x_1^2x_3\cr
1& x_1& x_2& x_1x_2&  x_2^2& x_1x_2^2\cr
1& x_2& x_1& x_1x_2&  x_1^2& x_1^2x_2\cr
} 
$$
\emph {1.}
Show that $\det(M) = \delta^{n!/2}$ with $\delta = \prod_{i<j}(x_i-x_j)$.

\emph {2.}
Let $s_1$, \ldots, $s_n \in \ZZ[\ux]$ be the $n$ \elr \smq functions, $F(T)$ be the \uvl \pol $F(T) = T^n - s_1T^{n-1} + \cdots + (-1)^n s_n$, and $U \in \MM_{n!}(\ZZ[\ux])$ be the trace-valued matrix, indexed by $E \times E$, with term $\Tr_{\rS_n}(\ux^{\alpha+\beta})$, $\alpha, \beta\in E$.
\\  
Let $f \in \gk[T]$ be a \polu of degree $n$, $\gA=\Adu_{\gk,f}$.\\
Prove the \egt $\Disc_\gk \gA = \disc(f)^{n!/2}$ (also found in Fact~\ref{factDiscriAdu}). 
\\
And conversely?

\emph {3.}
Revisit \thref{theoremAdu3}.
}

\end {exercise}

\vspace{-1em}
\begin{exercise}\label{exoTnADU}
 {(The \adu of the \pol $f(T) = T^n$)}\\
{\rm  
Let $f(T) = T^n$ and $\gA =\Adu_{\gk,f} = \gk[x_1, \ldots, x_n]$.
 Describe the structure of~$\gA$.
}
\end {exercise}

\vspace{-1em}
\begin{exercise}\label{exoNilIndexInversiblePol}
{(Invertible \pols and nilpotency indices)} \\
{\rm  
Here we propose a quantitative version of item~\emph{4} of Lemma~\ref{lemGaussJoyal}. Let $\gk$ be a commutative \riz, $f$, $g \in \kX$ satisfying $fg = 1$ \hbox{and $f(0) = g(0) = 1$}. We write $f = \sum_{i=0}^n a_i X^i$, $g = \sum_{j=0}^m b_j X^j$.  Show that 

\snic {
a_1^{\alpha_1} a_2^{\alpha_2} \cdots a_n^{\alpha_n} 
b_1^{\beta_1} b_2^{\beta_2} \cdots b_m^{\beta_m}  = 0
\quad \hbox{if}\quad \sum_i i\alpha_i + \sum_j j\beta_j > nm.
}

In particular, for $i \ge 1$, we have $a_i^{\lceil (nm+1)/i\rceil}
= 0$ and therefore $a_1^{nm + 1} = 0$.

}

\end{exercise}

\vspace{-1em}
\begin{exercise}\label{exoADUpthRoot}
{(The \adu of the \pol $f(T) = T^p - a$ in \cara $p$)}\\
{\rm  
Let $p$ be a prime number, $\gk$ be a \ri in which $p\cdot 1_\gk=0$ and $a \in \gk$.  \\ 
Let $f(T) = T^p - a \in \gk[T]$, $\gA =\Adu_{\gk,f} = \gk[x_1, \ldots,
x_p]$, $\gk[\alpha] = \aqo{\gk[T]}{f}$, such that $T-a=(T-\alpha)^p$.  Let $\varphi : \gA \twoheadrightarrow
\gk[\alpha]$ be the $\gk$-morphism $x_i \mapsto \alpha$. 
Make the \id $\Ker\varphi$ explicit and describe the structure of the \klg $\gA$.
\\
NB: if $\gk$  is a \cdi and $a$ is not a $p^{\rm th}$ power in $\gk$, by Exercise~\ref{exoPrimePowerRoot}, the \pol
$f(T)$ is \ird and $\gk[\alpha]$ is a field of \dcn of $f$ over $\gk$.

}

\end{exercise}

\vspace{-1em}
\begin{exercise}\label{exoA5GaloisGroup} \\
 {(The trinomial $T^5+5bT\pm 4b$ where $b= 5a^2-1$,
with Galois group $\rA_5$)}\\
{\rm  
Consider a trinomial $T^5 + bT + c$. We will determine $b$, $c$ such that its discriminant is a square and obtain an \ird \pol with Galois group $\rA_5$ as an illustration of the modular method. 
\\
We use the \egt $\disc_T(T^5 + bT + c) = 4^4b^5 + 5^5c^4$ (see \Pbmz~\ref{exoDiscriminantsUtiles}).

\emph {1.} To force the \discri into being a square in $\ZZ$, explain why what follows is reasonable: $b \aff 5b$, $c \aff 4c$, then $f_a(T) = T^5 + 5(5a^2-1)T \pm 4(5a^2 - 1)$. The \discri is then the square $2^8 5^6 a^2(5a^2 - 1)^4$.

\emph {2.}
Taking $a=1$, we obtain $f_1(T) = T^5 + 20T \pm 16$ in $\ZZ[T]$. By examining the \fcns of $f_1$ modulo $3$ and $7$, show that $f_1$ is \ird with Galois group $\rA_5$. Deduce that for $a \equiv 1 \bmod 21$, $f_a$ is \ird with Galois group $\rA_5$. Show that the same thing holds for $f_a$ given as a \pol with \coes in the field of rational fractions $\QQ(a)$.

}

\end{exercise}

\vspace{-1em}
\begin{exercise}
\label{propGaloisJordan}
{\rm
\emph{(When a resolvent admits a zero in the base field)}\\
In the context of \ref{notas1.3} let $y\in\gB$,  $G.y=\so{\yr}$ and $g(T)=\Rv_{G,y}(T)$.
\begin{enumerate}
\item \label{i1propGaloisJordan}
Suppose that $a\in\gK$ is a simple zero of $g$. 
\begin{enumerate}
\item \label{i1apropGaloisJordan} $\fc=\gen{y-a}_\gB$ is a Galoisian \id of $(\gK,\gB,G)$
.
\item \label{i1bpropGaloisJordan}  If $\beta:\gB\to{\gC=\gB/\fc}$ is the canonical \prnz, and if $H=\St_G(\fc)$ is the new approximation of the Galois group, then $\beta(y_1)=a$ and for $j\neq 1$, $\Rv_{H,y_j}$ divides $g(T)/(T-a)$ (as usual we identify~$\gK$ with a subfield of $\gB$ and $\beta(\gK)$ with a subfield of $\gC$).
\end{enumerate}
\item \label{i2propGaloisJordan}
Suppose that $a\in\gK$ is a zero of $g$ with multiplicity $k$.
\begin{enumerate}
\item \label{i2apropGaloisJordan}
There exist $j_2$, \ldots, $j_k\in\lrb{2..r}$ such \hbox{that $\fc=\gen{y_1-a,y_{j_2}-a,\ldots,y_{j_k}-a}$} is a minimal Galoisian \id among those that contain $y-a$.\\
Let $j_1=1$. 
Show that, for $j\neq j_1$, \ldots, $j_k$, $y_j-a$ is \iv modulo $\fc.$
\item \label{i2bpropGaloisJordan}
Let $\beta:\gB\to{\gC=\gB/\fc}$ be the canonical \prnz, and  $H=\St_G(\fc)$. \\
Then $\beta(y_{j_1})=\cdots=\beta(y_{j_k})=a$, and for $j\neq j_1$, \ldots, $j_k$, the resolvent $\Rv_{H,y_j}$ divides $g(T)/(T-a)^k$.
\end{enumerate}
\item \label{i3propGaloisJordan}
Suppose that $\fc$ is a Galoisian \id of $\gB$ and that $\St_G(y)$ contains $\St_G(\fc)$, then $g(T)$ admits a zero in $\gK$.
\end{enumerate}
\rem Item~\emph{1} 
justifies the \gui{Jordan method} for the computation of the Galois group. See \paref{Jordan}.
\eoe

}
\end{exercise}

\vspace{-1em}
\begin{exercise}
\label{propRvRelDec}
{(When we know the decomposition into prime factors of a \spl resolvent)}
{\rm In the context of \ref{notas1.3} let $y\in\gB$ and $G.y=\so{y_1,\ldots,y_r}$.
\\
Suppose that $\Rv_{G,y}=\Mip_y=R_1\cdots R_\ell$, with the $R_i$ being \irds and $\ell>1$.
Compute a Galoisian \idm $e$ of $\gB$ with the following \prtsz, where we let $(\gK,\gC,H)$ be the corresponding Galois quotient and $\beta:\gB\to\gC$ be the canonical \prnz.
\begin{enumerate}
\item For each $i\in\lrbr$, the \pol $\Mip_{\beta(y_i)}$ is equal to one of the $R_j$'s.
\item The group $H$ operates over $\so{\beta(y_1),\ldots,\beta(y_r)}$. 
\item The orbits are of length $d_1=\deg(R_1),\ldots,d_\ell=\deg(R_\ell)$.
\item This situation recurs in every Galois quotient of $(\gK,\gC,H)$.
\end{enumerate}
\rem  Exercise~\ref{propRvRelDec}
is the basis of the \gui{McKay-Soicher method} for the computation of the Galois group. See \paref{soicher}.
\eoe

}
\end{exercise}

\vspace{-1em}
\begin{exercise}
\label{propRvRelDecMin}
{\rm
\emph{(When a \polmin strictly divides a resolvent)}\\
In the context of \ref{notas1.3} let $y\in\gB$ and $G.y=\so{y_1,\ldots,y_r}$.\\
Suppose that $g(T)=\Rv_{G,y}(T)\neq\Mip_y(T)$.
Let $(\gK,\gC,H)$ be a Galois quotient (with the canonical \prn $\beta:\gB\to\gC$) in which each $\beta(y_i)$ admits a \polmin equal to its resolvent.\\
Show that for the different zeros $\beta(y_j)$ of $g_1(T)=\Mip_{\beta(y_1)}(T)$ in $\gC$, the fibers $\beta^{-1}\big(\beta(y_j)\big)$ all have the same number of \eltsz, say $n_1$. \\
In addition,  $g_1^{n_1}$ divides $g$ and $g/g_1^{n_1}$ is \com with $g_1$.
}
\end{exercise}


\sol


\exer{exothSteCdiClass} \emph{1.} 
This results from the fact that a connected \zed \ri is local and from the fact that, by \TEMz, we know the in\dcp \idms of the \algz, which form a \sfioz.

 \emph{2.} In the case of an \alg $\aqo\KX f$ with \spl $f$, finding the \idms is the same as factoring the \polz. But there does not exist any \gnl \fcn \algo for \spl \polsz.

 \emph{3.} A \cov version consists in asserting that, concerning a computation, we can always \gui{act as though} the result (proven by means of \TEMz) were true.
This \emph{dynamic version} is expressed as follows. \\
\emph{Let $\gK$ be a \zed \ri (special case: a \cdiz).\\
Let $(x_i)_{i\in I}$ be a finite family of \elts in an integral \Klg $\gB$ (special case:  a finite \Klgz).
\\
There exists 
a \sfio $(e_1,\ldots,e_n)$ such that in each component 
$\aqo{\gB}{1-e_j}$, each $x_i$ is nilpotent or \ivz.} \\
We prove this result as follows: Lemma~\ref{lemZrZr1} tells us that $\gB$ is \zedz; we conclude by the \zed splitting lemma (Lemma~\ref{thScindageZed}).


\exer{exoMatMag3} 
We easily check that the $\ZZ$-syzygy module 
between the matrices $P_1$, \ldots, $P_6$ is generated by $(1,1,1, -1,-1,-1)$.
We will also use the fact that the number of \moms of degree $d$ in $n$ variables is ${d + n-1 \choose n-1}$.

\emph {1a.}
Let $S(Y,Z) = \sum_{i+j=m-1} Y^i Z^j$, so $Y^m - Z^m = (Y-Z)S(Y,Z)$.
In this \egtz, we make $Y = U_1U_2U_3$, $Z = U_4U_5U_6$.
We obtain the desired result by multiplying by $U_1^{a-m} U_2^{b-m} U_3^{c-m} U_4^d U_5^e U_6^f$.

\emph {1b.}
We clearly have $\fa\subseteq\Ker\varphi$. The \egt $\gk[\uU] = \fa + \fa^\bullet$ results from item~\emph{1a.} It therefore suffices to see that $\Ker\varphi \cap \fa^\bullet = \{0\}$, \cad that the restriction of $\varphi$ to $\fa^\bullet$ is injective. As $\varphi$ transforms a \mom into a \momz, it suffices to see that if two \moms $U_1^a \cdots U_6^f$ and $U_1^{a'} \cdots U_6^{f'} \in \fa^\bullet$ have the same image under $\varphi$, then they are equal. \\
We have $(a,b,c, \ldots, f) = (a',b',c', \ldots, f') + k(1,1,1,-1,-1,-1)$ with $k \in \ZZ$, and \hbox{as $\min(a,b,c) =
\min(a',b',c') = 0$}, we have $k = 0$, which gives the \egt of the two \momsz.

\emph {1c.}
The number $M_d$ that we are searching for is the dimension over $\gk$ of the \hmg component 
 of degree $d$ of $\gk[\ux^{P_1}, \ldots, \ux^{P_6}]$ or (via $\varphi$) that of $\fa^\bullet_d$. \\
But we also have $\gk[\uU] = \fb\oplus\fa^\bullet$ where $\fb$ is the (monomial) \id generated by the \moms divisible by $U_1U_2U_3$ (in some way, $\fb$ is an initial \id of $\fa$).  \\
We therefore have $\gk[\uU]_d = \fb_d \oplus \fa^\bullet_d$ and

\snic {
\dim_\gk \gk[\uU]_d = {d+5 \choose 5}, \quad
\dim_\gk \fb_d = {d+5-3 \choose 5}, \quad
M_d = \dim_\gk \fa^\bullet_d = {d+5 \choose 5} - {d+2 \choose 5}.
}

\emph {2a.}
We define $V_2, V_3$ by $U_2 = U_1 + V_2$, $U_3 = U_1 + V_3$. \\
The \pol $U_1U_2U_3 - U_4U_5U_6$ given in $\gk[U_1, V_2, V_3, U_4, U_5, U_6]$ becomes monic in $U_1$ of degree $3$.
We leave it up to the reader to check the other details.

\emph {2b.}
The number we are looking for is also $M_d = \dim_\gk \gB_d$. But we have

\snic{\gB_d =
\gA_d \oplus \gA_{d-1} u_1 \oplus \gA_{d-2} u_1^2 \simeq 
\gA_d \oplus \gA_{d-1} \oplus \gA_{d-2}.}

It suffices to use the fact that $\gA$ is a \pol \ri over $\gk$ with $5$ \idtrsz.
As an indication, for $d = 0$, $1$, $2$, $3$, $4$, $5$, $M_d= 1$, $6$, $21$, $55$, $120$, $231$.

\emph {3a.}
The $\ZZ$-module $E$ is free of rank $5$: $5$ arbitrary matrices among $\{P_1, \ldots, P_6\}$ form a $\ZZ$-basis of it.

\emph {3b.}
Since $P_1 + P_2 + P_3 = P_4 + P_5 + P_6$, we have

\snic {
\gB_{11} = 
\gk[\ux^{P_1}][\ux^{\pm P_2}, \ux^{\pm P_3}, \ux^{\pm P_4}, \ux^{\pm P_5}] =
\gk[\ux^{P_6}][\ux^{\pm P_2}, \ux^{\pm P_3}, \ux^{\pm P_4}, \ux^{\pm P_5}]
.}

We then see that $\gB_{11}$ is a localized \ri 
 of $\gk[\ux^{P_1}, \ux^{P_2}, \ux^{P_3}, \ux^{P_4}, \ux^{P_5}]$, which is a \pol \ri with $5$ \idtrs over $\gk$, so \iclz.

\emph {3c.}
We define $\gB_{ij}$ such that it is contained in $\gk[\,\ux^M \,\vert\, M \in E,\ m_{ij} \ge 0\,]$. For example, for $(i,j) = (3,1)$, the matrices $P_k$ with a null coefficient in position $(3,1)$ are those other than $P_3$, $P_5$, which leads to the \dfn of $\gB_{31}$:

\snic {
\gB_{31} = \gk[\ux^{P_3}, \ux^{P_5}]
[\ux^{\pm P_1}, \ux^{\pm P_2}, \ux^{\pm P_4}, \ux^{\pm P_6}].
}

We then have the \egt $\gB = \bigcap_{i,j} \gB_{ij}$, and as the $\gB_{ij}$'s are all \icl with the same quotient field $\Frac\gB$, the \ri $\gB$ is \iclz.


\exer{exoIdentiteDiscriminantale} 
\emph {2.}
We write $U = \tra{M} M$ and take the \deterz.\\
This gives $\Disc_\gk \gA = \disc(f)^{n!/2}$ from $\det(M) = \delta^{n!/2}$. 
Conversely, since this is a matter of \idas in $\ZZ[\ux]$, the \egt $(\det M)^2=(\delta^{n!/2})^2$
implies $\det M=\pm\delta^{n!/2}$. 

\emph {3.}
In \thref{theoremAdu3}, let us not assume that $f$ is \spl over $\gC$.
By hypothesis, {we have  $\varphi(f)(T) = \prod_{i=1}^n(T- u_i)$}. 
\\
With $\gA = \gk[\xn] = \Adu_\gk(f)$, we then have a morphism of $\gC$-\algs $\Phi : \gC\te_\gk\!\gA \to \gC^{n!}$ which performs $1 \te x_i \mapsto (u_{\sigma(i)})_{\sigma\in \rS_n}$. \\
The canonical $\gk$-basis $\cB(f)$ of $\gA$ is a $\gC$-basis of $\gC\te_\gk\gA$ and the matrix of $\Phi$ for this basis (at the start) and for the canonical basis of $\gC^{n!}$ (at the end) is the above matrix $M$ where $x_i$ is replaced by $u_i$. We deduce that $\Phi$ is an \iso \ssi $\varphi\big(\disc(f)\big)\in\gC\eti$, \cad if $f$ is \spl over $\gC$.

Finally, let us only suppose that an \alg $\varphi:\gk\to \gC$ diagonalizes~$\gA$. This means that we give $n!$ characters $\Adu_{\gC,\varphi(f)}\to\gC$ which, when put together, give an \iso of $\gC$-\algs of $\Adu_{\gC,\varphi(f)}$ over $\gC^{n!}$.
\\
Since there exists a character $\Adu_{\gC,\varphi(f)}\to\gC$, the \pol $\varphi(f)(T)$ completely factorizes in $\gC$.
\\
Finally, the \discri of the canonical basis of $\Adu_{\gC,\varphi(f)}$ is $\varphi\big(\disc(f)\big)^{n!/2}$ and the \discri of the canonical basis of $\gC^{n!}$ is $1$. 
Therefore, $f$ is \spl over~$\gC$.


\exer{exoTnADU} 
We have $\gA = \aqo {\gk[\Xn]}{S_1, \ldots, S_n}$ where $S_1$, \ldots, $S_n$ are the $n$ \elr \smq functions of $(\Xn)$; the \id $\gen {S_1, \ldots, S_n}$ being \hmgz, the \klg $\gA$ is graded (by the degree). Let $\gA_d$ be its \hmg component of degree $d$ and $\fm = \gen {\xn}$; we therefore have $\gA = \gA_0 \oplus \gA_1 \oplus \gA_2 \oplus \dots$ with $\gA_0 = \gk$
and

\snic {
\fm^d = \gA_d \oplus \gA_{d+1} \oplus \dots, \quad
\fm^d = \gA_d \oplus \fm^{d+1}.
}

Since $x_i^n = 0$, we have $\fm^{n(n-1)+1} = 0$, so $\gA_d = 0$ for $d \ge n(n-1)+1$. Recall the basis $\cB(f)$ of $\gA$, formed from the \elts $x_1^{\alpha_1} \ldots x_n^{\alpha_n}$ with $0 \le \alpha_i < n-i$. 
For all~$d$, the \hmg component $\gA_d$ of degree $d$ is a free \kmo whose basis is the set of the $x_1^{\alpha_1} \ldots x_n^{\alpha_n}$ with $0 \le \alpha_i < n-i$ and $|\alpha| = d$. The cardinality of this basis is the \coe of degree $d$ in the \pol $S(t) \in \ZZ[t]$

\snic {
S(t) = 1 (1+t) (1 + t + t^2) \cdots (1 + t + \cdots + t^{n-1}) =
\prod_{i=1}^n \frac {t^i-1} {t-1}.
}

Indeed, a multi-index $(\alpha_1, \ldots, \alpha_n)$ such that $0 \le \alpha_i < n-i$ and $|\alpha| = d$ is obtained by choosing a \mom $t^{\alpha_n}$ of the \pol $1 + t + \cdots + t^{n-1}$, a \mom $t^{\alpha_{n-1}}$ of the \pol $1 + t + \cdots + t^{n-2}$ and so on, the product of these \moms being~$t^d$. We thus obtain the Hilbert-Poincar\'e series $S_\gA(t)$ of $\gA$,

\snic{
S_\gA(t) \eqdefi \som_{i=0}^{\infty}\dim_\gk \gA_d\ t^d 
\eqdf {\rm here} \som_{0 \le \alpha_i < n-i} t^{|\alpha|} = S(t).
}

The \pol $S$ is a \polu of degree $e= 1 + \cdots + n-1= n(n-1)/2$.
We have $S(1) = n!$, in accordance with $S(1) = \dim_\gk\gA$.

\Deuxcol{.8}{.15}
{Variant. Let $\gB = \gk[S_1, \ldots, S_n] \subset \gC = \gk[\Xn]$. \\
Then $\gC$
is a free $\gB$-module with the $\uX^\alpha = X_1^{\alpha_1} \cdots X_n^{\alpha_n}$ as its basis, with $0 \le \alpha_i < n-i$. This basis is above the basis $\cB(f)$ of~$\gA$ over $\gk$ if we consider that we have a commutative diagram where each vertical arrow is a reduction modulo $\gen {S_1, \ldots, S_n}$.  
}
{~

$\xymatrix @R = 0.5cm @C = 0.8cm
{\gB\ar@{->>}[d]\ar@{->}[r] &\gC\ar@{->>}[d]\\
\gk\ar@{->}[r] &\gA \\}$ 
}

Writing $\gC = \bigoplus_{\alpha} \gB \uX^\alpha$, 
  with the shift
$S_{\gB\uX^\alpha}(t) = t^{|\alpha|}S_\gB(t)$, we have the following \egt between the Hilbert-Poincar\'e series

\snic {
S_\gC = S_\gA\ S_\gB  \qquad \hbox {with} \qquad
S_\gA = \som_{0 \le \alpha_i < n-i} t^{|\alpha|}.
}

However, it is easy to see that

\snic {
S_\gC(t) = \fraC {1} {(1-t)^n}, \;
S_\gB(t) = \prod_{d=1}^n \fraC {1} {1 - t^d}, 
\; \hbox { and so } \,
S_\gA(t) = \frac {S_\gC} {S_\gB} = \prod_{d=1}^n \frac {1-t^d} {1-t},
}

once again giving us the result for $S_\gA$.

Let us now move onto the powers of the \id $\fm$.  \\
Let $\varphi : \gA\twoheadrightarrow \gk$ be the character $x_i \mapsto 0$ with kernel $\fm = \gen {\xn}$. \\
We have $\gA=\kxn=\gk\oplus\fm$, $\fm\subseteq\DA(0)\subseteq \Rad(\gA)$
and for $z \in \gA$, 

\snic{z\in\Ati\iff\varphi(z)\in\gk\eti\iff z\in\gk\eti\oplus\fm.}

We have $\DA(0)=\rD_\gk(0)\oplus\fm$, $\Rad(\gA)=\Rad(\gk)\oplus\fm$.\\
Since $\fm \subseteq \Rad(\gA)$ is \tfz, we have 
$$
\preskip.2em \postskip.3em 
\fm^d = \fm^{d+1} \iff \fm^d
= 0 
$$
(Lemma~\ref{lemLocaliseFini}), which, since $\fm^d = \gA_d \oplus \fm^{d+1}$, is equivalent to $\gA_d = 0$. We deduce that $\fm^{e+1} = 0$.

\rem if $\gk$ is local, then so is $\gA$, and $\Rad\gA=\varphi^{-1}(\Rad\gk)$.\eoe


\exer{exoNilIndexInversiblePol} 
Consider the \pol \ri $\gC =\gk[\an,b_1,\dots,b_m]$, let $f(X) = 1 + \sum_{i=1}^n a_i X^i$, $g(X) = 1 + \sum_{j=1}^m b_j
X^j$, and $\fc=\rc_\gC(fg - 1)$. We assign to $a_i$ the weight $i$ and to $b_j$ the weight $j$. The \coe of degree $k$ of $fg-1$ is \hmg of degree $k$, so the \id $\fc$ is \hmgz. \\
Let $\gC' = \gC\sur\fc$. This \klg $\gC'$ is graded via the above weight and we must show that $\gC'_d = 0$ for $d > nm$. It is clear that $\gC'_d = 0$ for large enough $d$. 
We will determine the Hilbert-Poincar\'e series $S_{\gC'}$ of $\gC'$ (which here is a \polz)

\snic {
S_{\gC'}(t) \eqdefi \sum_{d \ge 0} \dim_\gk \gC'_d\, t^d =
\frac{ \prod_{d=1}^{n+m} (1-t^d) }{ \prod_{i=1}^{n} (1-t^i)\prod_{j=1}^{m} (1-t^j)}\;.
}

To prove this \egtz, we construct $\gC$ and $\gC'$ in a different way. \\
We consider $n+m$ \idtrs $(\Xn, \Ym)$, and let $(\an)$ be the \elr \smq functions of $(\Xn)$, and $(b_1,\dots,b_m)$ be the \elr \smq functions of $(\Ym)$. Since

\snuc {
\prod_{i=1}^n (T+X_i) \prod_{j=1}^m (T+Y_j) =
(T^n + a_1 T^{n-1} + \cdots + a_n) (T^m + b_1 T^{m-1} + \cdots + b_m), 
}

we see, by letting $a_0=b_0 = 1$, that $\sum_{i+j=d} a_ib_j$ is the $d^{\rm th}$ \elr \smq function of $(\Xn, \Ym)$. As $(\an, b_1,\dots,b_m)$ are \agqt independent over $\gk$, we can consider that $\gC$ is the following graded sub\alg

\snic {
\gC = \gk[\an, b_1,\dots,b_m]  \subset \gD = \gk[\Xn, \Ym],
}

and that the \id $\fc$ of $\gC$ is generated by the $n+m$ sums $\sum_{i+j=d} a_ib_j$, which are the \elr \smq functions of $(\Xn, \Ym)$.
\\
The \alg $\gD$ is free over $\gC$ of rank $n!m!$, as for a double \aduz. More \prmtz, here are some bases. \\
The $\uX^\alpha = X_1^{\alpha_1} \cdots X_n^{\alpha_n}$ for $0 \le \alpha_i < n-i$ form a basis of $\gk[\uX]$ over $\gk[\ua]$. \\
The $\uY^\beta = Y_1^{\beta_1} \cdots Y_m^{\beta_m}$ with $0 \le \beta_j < m-j$ form a basis of $\gk[\uY]$ over $\gk[\ub]$. 
\\
Thus, the $\uX^\alpha \uY^\beta$ form a basis of $\gD = \gk[\uX,\uY]$ over $\gC = \gk[\ua,\ub]$.\\
Finally, by the \eds $\gC\to\gC' = \gC\sur\fc$, the $\ux^\alpha \uy^\beta$ form a basis {of $\gD' = \gD\sur{\fc\gD} = \gk[\ux,\uy]$ over $\gC'$}.

\Deuxcol{.77}{.18}
{
We have a commutative diagram at our disposal where each vertical arrow is a reduction modulo $\fc$. Our aim is to determine the Hilbert-Poincar\'e series $S_{\gC'}$ of $\gC'$ given that we know those \hbox{of $\gD'$, $\gC$} and~$\gD$ (because $\gC$ and $\gD$ are \pol \risz, \hbox{and $\gD'$} is the \adu of $T^{n+m}$ over~$\gk$).
}
{\vspace{-3mm}
$$\xymatrix @R = 0.5cm @C = 0.8cm
{\gC\ar@{->>}[d]\ar@{->}[r] &\gD\ar@{->>}[d]\\
\gC'\ar@{->}[r] &\gD' \\}$$ 
}

We conclude the computations in the following simple manner. \\
We write $\gD = \bigoplus_{\alpha, \beta} \gC\, \uX^\alpha\uY^\beta$,
so 

\snic {  
S_{\gD}(t) = F(t) S_{\gC}(t)  \quad \hbox {with} \quad
F(t) = \sum_{\alpha, \beta} t^{|\alpha| + |\beta|} =
\sum_{\alpha} t^{|\alpha|} \sum_{\beta} t^{|\beta|},
}

and we also have $S_{\gD'}(t) = F(t) S_{\gC'}(t)$. We have seen in Exercise~\ref{exoTnADU} that

\snic {
F(t) = \prod_{i=1}^n \fraC {1-t^i} {1-t}\;\prod_{j=1}^m \fraC {1-t^j} {1-t},
\qquad
S_{\gD'}(t) = \prod_{d=1}^{n+m} \fraC {1-t^d} {1-t}.
}

Then let $S_d(t) = (1-t^d)/(1-t)$. It is a \pol of degree $d-1$ \hbox{and $S_d(1) = d$}. We have therefore obtained
$$\preskip.2em \postskip.0em
S_{\gC'}(t) = {S_1 S_2 \cdots S_{n+m} \over 
S_1 S_2 \cdots S_{n}\,  S_1 S_2 \cdots S_{m}}, 
$$
with
$$\preskip.0em \postskip.4em 
\deg S_{\gC'} = {(n+m-1)(n+m) - (n-1)n - m(m-1) \over 2} = nm. 
$$
Thus, as desired, $\gC'_k = 0$ for $k > nm$. \\
Please note that
$\dim_\gk \gC' = S_{\gC'}(1) = {n+m \choose n}$. 

\exer{exoADUpthRoot} 
For each $i\in\lrbp$ the restriction $\varphi : \gk[x_i] \to\gk[\alpha]$ is an \isoz.  
Consider the \id 

\snic{\fm = \gen {x_i - x_j, i,j \in\lrbp}= \gen {x_1 - x_i, i \in\lrb{2..p}}.}

Then $\gA=\gk[x_1]\oplus \fm$, hence $\fm=\Ker\varphi$.
\\
Actually we can regard $\gA$ as the \adu $\Adu_{\gk[x_1],g}$ for the \pol $g(T)=f(T)/(T-x_1)=(T-x_1)^{p-1}$ over the \ri $\gk[x_1]$ which brings us back to Exercise~\ref{exoTnADU}.
In particular 

\snic{\fm^{1+(p-1)(p-2)/2}=0$, 
$\DA(0)=\rD_{\gk[x_1]}(0)\oplus \fm$ and $\Rad(\gA)=\Rad(\gk[x_1])\oplus \fm.}


\exer{exoA5GaloisGroup} 
\emph {1.}
The goal of the operation $b \aff 5b$, $c \aff 4c$ is to replace $4^4b^5 + 5^5c^4$ by $4^4 5^5 (b^5 + c^4)$; by imposing $c = \pm b$, we obtain $4^4 5^5 b^4 (b + 1)$ which is easy to turn into a square by imposing $5(b + 1) = a^2$. To avoid the denominator $5$, we impose $5(b + 1) = (5a)^2$ instead, \cad $b = 5a^2-1$.

\emph {2.}
For $a \in \QQ\sta$, the \pol $f_a(T) \in \QQ[T]$ is \splz.  Modulo the small prime numbers we find the following \dcns of $f_1(T) = T^5 + 20 T + 16\vep$, \hbox{with $\vep \in \{\pm 1\}$}, into \irds factors
\[ 
\begin{array}{rccl} 
\mod &  2 &:& T^5 \\[1mm] 
\mod &  3 &:& f_1(T) \\[1mm] 
\mod &  5 &:& (T + \vep)^5\\[1mm] 
\mod &  7 &:& (T+2\vep)(T+3\vep)(T^3 + 2\vep T^2 + 5T + 5\vep)  
 \end{array}
\] 
The result modulo $3$ proves that $f_1(T)$ is \ird over $\ZZ$. Its Galois group~$G$ is a transitive subgroup of $\rA_5$ that contains a $3$-cycle (given the reduction modulo $7$). This implies $G = \rA_5$. Indeed, a transitive subgroup of $\rS_5$ containing a $3$-cycle is equal to $\rS_5$ or $\rA_5$. As for $\QQ(a)$ as a base field, the \pol $f_a(T)$ is \ird in $\QQ[a][T]$ since its reduction modulo $a = 1$ is \ird in $\QQ[T]$. Therefore it is \ird in $\QQ(a)[T]$.
Using the fact that its discriminant is a square and the reduction modulo $a=1$, we obtain that its Galois group is $\rA_5$.
\\
The readers might ask themselves the following question: For every $a \in \ZZ\setminus \{0\}$, is the \pol $f_a(T)$ \ird with Galois group $\rA_5$?

\emph{Possible experiment}. 
\\Here is the distribution of the types of permutation of the transitive subgroups of~$\rS_5$. 
\\For the $7$ types that appear in $\rS_5$, we use the following notation

\snuc {
t_1 = (1^5),\, t_2 =(2,1^3),\, t_{22} = (2^2,1),\, t_3 = (3,1^2),\, 
t_{3,2} = (3,2),\, t_4 = (4,1),\ t_5 = (5).
}

Thus $t_{22}$ is the type of the double-transpositions, $t_3$ that of the $3$-cycles,
etc{\dots} The announced table:
\vskip2pt
\snic {
\begin {array} {cccccc}
\hline\\[-2.5mm]
G &\rC{5} &\ASL_1(\FF_5) &\AGL_1(\FF_5) &\rA_5 &\rS_5
   \\[1mm]\hline\\[-2.5mm]
\#G &5 &10 &20  &60 &120
   \\[1mm]\hline\\[-2.5mm]
&t_1^1\ t_5^4 &t_1^1\ t_{22}^5\ t_5^4 &t_1^1\ t_{22}^5\ t_4^{10}\ t_5^4 &
    t_1^1\ t_{22}^{15}\ t_3^{20}\ t_5^{24} &
    t_1^1\ t_2^{10}\ t_{22}^{15}\ t_3^{20}\ t_{32}^{20}\ t_4^{30}\ t_5^{24}
\\[1mm]\hline
\end {array}
}
\vskip2pt

For example in the last row, under $\rA_5$, $t_1^1\ t_{22}^{15}\ t_3^{20}\ t_5^{24}$ means that $\rA_5$ contains the identity, $15$ double-transpositions, $20$ $3$-cycles and $24$ $5$-cycles ($1 + 15 + 20 + 24 = 60$).
\\
The reader will be able to experimentally test Cebotarev's density \tho with the help of a Computer Algebra \sysz. 
We must examine the \fcn of $f_1(T)$ modulo \gui{a lot} of primes $p$ and compare the distribution obtained from the types of \fcn with that of the types of permutation of $\rA_5$.  
\\
The author of the exercise has considered the first $120$
prime numbers --- other than $2$ and $5$ which divide $\disc(f_1)$ ---
and his program has found the following distribution
$$\preskip-.2em \postskip.4em 
t_{22}^{33}\ t_3^{38}\ t_5^{49}. 
$$
This means that we have found a \fcn of type $t_{22}$ ($2$ \ird factors of degree $2$, $1$ \ird factor of degree $1$) $33$ times, a \fcn of type $t_{3}$ $38$ times and a \fcn of type $t_{5}$ $49$ times (no \fcn of type~$t_1$). A distribution to be compared with that of $\rA_5$. As for the type $t_1$, the smallest prime $p$ for which $f_1(T) \bmod p$ is entirely decomposed is $p = 887$. Finally, when treating $1200$ primes instead of $120$, we find the distribution
$$\preskip.3em \postskip.1em 
 t_1^{16}\ t_{22}^{304}\ t_3^{428}\ t_5^{452}.
$$

\exer{propGaloisJordan}
\emph{1a.} We need to show that $\gen{y_1-a}+\gen{y_j-a}=\gen{1}$ for $j\in \lrb{2..r}$. For example in the quotient $\aqo{\gB}{y_1-a,y_2-a}$ the \pol $g(T)=\prod(T-y_j)$ has two factors equal to $T-a$ which implies that $g'(a)=0$. As $g'(a)$ is \iv by hypothesis (which remains true in a quotient), we indeed have $0=1$ in the quotient.

 \emph{1b.} We easily see that $H=\St(y_1)$.
Therefore $H$ operates over $\so{\beta(y_2),\ldots,\beta(y_r)}$. \\
However, $g(T)/(T-y_1)=\prod_{j=2}^r(T-y_j)$ in $\gB$, so $g(T)/(T-a)=\prod_{j=2}^r\big(T-\beta(y_j)\big)$ in $\gC$.
\perso{ce n'est pas clair if l'orbite of $\beta(y_2)$ sous $H$ contient tous les $\beta(y_j)$ with $j\neq1$}

 \emph{2a.} It is clear that $y_1-a$ is a nonzero  \dvz in $\gB$. A minimal Galoisian \id $\fc$ containing $\gen{y_1-a}$ is obtained by adding as many conjugates of $\gen{y_1-a}$ as possible under the condition of not reaching the \id $\gen{1}$.
The \id $\fc$ is therefore of the form $\gen{y_j-a\,|\,j\in J}$ for a subset $J$ of $\lrbr$. It remains to see if the number of $j$'s such that $y_j-a\in\fc$ is $k$. 
However, for every index $j$, the \elt $y_j-a$ is null or \iv modulo $\fc$. Since $g(T)=\prod_j\big(T-\beta(y_j)\big)$,
and since $a$ is a zero with multiplicity $k$ of $g$, the number of $j$'s such that $\beta(y_j)=a$ is equal to $k$  
(let $g(a)=g'(a)=\cdots=g^{(k-1)}(a)=0$ and $g^{(k)}(a)$ be \ivz).

 \emph{2b.} Reason as in \emph{1b.}

 \emph{3.} The Galois quotient $\gC=\gB/\fc$ is obtained with its group $H=\St_G(\fc)$.
By hypothesis $\ov {y_1}\in\Fix(H)$ so $\ov{y_1}\in\gK$. Let $a$ be the \elt of $\gK$ in question.
In~$\gC$ we have $g(T)=\prod_j(T-\ov{y_j})$, so ${g(a)}=0$.
Finally, $\gK $ is identified with its image in~$\gC$.

\medskip
\exl Here is an example with $\deg f=6$. We ask {\tt Magma} to compute the \polmin of $y=x_4+x_5x_6$, and then to factorize it. If $g$ is the first factor, $z=g(y)$ is a nonzero \dvzz. We launch \Algoz~\ref{alg2idmgal} with~$z$.
We therefore obtain the new approximations of the \cdr and of the Galois group by treating the oddity \gui{$z$ is a nonzero \dvzz,} but we can observe a posteriori that $z$ has multiplicity $6$ in its resolvent and that $\gen{z}$ is Galoisian.

\label{example3Galois}
{\footnotesize
\begin{verbatim}
f:= T^6 - 3*T^5 + 4*T^4 - 2*T^3 + T^2 - T + 1;
y:=x4+x5*x6; pm:=MinimalPolynomial(y);
  T^60 - 46*T^59 + 1035*T^58 - 15178*T^57 + 163080*T^56 + ... + 264613
Factorization(pm);
    <T^6 - 4*T^5 + 8*T^4 - 6*T^3 + T + 1, 1>,
    ...
z:=Evaluate(T^6 - 4*T^5 + 8*T^4 - 6*T^3 + T + 1,y);
20*x4^3*x5^3*x6^3 - 15*x4^3*x5^3*x6^2 - 15*x4^3*x5^2*x6^3 +
  11*x4^3*x5^2*x6^2 + 2*x4^3*x5^2*x6 + 2*x4^3*x5*x6^2 + x4^3*x5*x6 - ...
// z divides 0, we compute the new Galois quotient
Affine Algebra of rank 6 over Rational Field
Variables: x1, x2, x3, x4, x5, x6
Quotient relations:
  x1 + x2 + x3 - x6^5 + 2*x6^4 - x6^3 - x6^2 - 1,
  x2^2 + x2*x3 - x2*x6^5 + 2*x2*x6^4 - x2*x6^3 - x2*x6^2 - x2 + x3^2 -
      x3*x6^5 + 2*x3*x6^4 - x3*x6^3 - x3*x6^2 - x3 + x6^5 - 2*x6^4 +
      x6^3 + x6^2,
  x3^3 - x3^2*x6^5 + 2*x3^2*x6^4 - x3^2*x6^3 - x3^2*x6^2 - x3^2 +
      x3*x6^5 -  2*x3*x6^4 + x3*x6^3 + x3*x6^2 - x6^5 + 2*x6^4 - x6^3 -
      x6^2 + 1,
  x4 + x5 + x6^5 - 2*x6^4 + x6^3 + x6^2 + x6 - 2,
  x5^2 + x5*x6^5 - 2*x5*x6^4 + x5*x6^3 + x5*x6^2 + x5*x6 - 2*x5 -
      x6^4 + 2*x6^3 - x6^2 - x6,
  x6^6 - 3*x6^5 + 4*x6^4 - 2*x6^3 + x6^2 - x6 + 1
Permutation group G2 acting we have set of cardinality 6
Order = 72 = 2^3 * 3^2
    (1, 4)(2, 5)(3, 6)
    (1, 2)
    (2, 3)
Degree(MinimalPolynomial(z)); 55
#Orbit(z,G); 60
\end{verbatim}
}
{\small

\exer{propRvRelDec}
Notice that the $y_i-y_j$'s for $i\neq j$ are \ivsz, and that this remains true in every Galois quotient.



\Biblio
The dynamic method is clearly presented, for the first time it seems, in the article by Paul Lorenzen \cite{Lor1953} published in 1953, which uses the equivalent of the closed covering principle \ref{prcfgrl}. See on this subject articles \cite{CLN2019} and \cite{CLN2020}.

Theorem \ref{thpolcohfd} says that a polynomial ring over a zerodimensional reduced ring is stronbly discrete and coherent. It admits a remarkable generalization to strongly discrete \coh \adpsz: see \cite{Yengui} and \cite{DVY2015}. 

The versions that we have given of the \nst \gui{without \agq closure} can be found in a related form in \cite[VIII.2.4, VIII.3.3]{MRR}.

The intrinsic difficulty of the \pb of the \iso of two \agq closures of a field is illustrated in \cite[Sander, Theorem~26]{sand}, which shows that, in the presence of \TEM but in the absence of the axiom of dependent choice, it is impossible to prove in $\ZF$ that two \agq closures of $\QQ$ are \isocz.

The treatment of  Galois theory based on Galois quotients of the \adu dates back to Jules Drach \cite[1898]{Drach} and to Ernest Vessiot \cite[1904]{Vessiot}. Here is an extract of the introduction of the latter article, which speaks in the language of the time about Galois quotients of the \aduz:

\gui{\'Etant donn\'ee une \eqn alg\'ebrique, que l'on consid\`ere comme remplac\'ee par le syst\`eme $(S)$ des relations entre les racines $\xn$ et les \coesz,
on \'etudie d'abord le probl\`eme fondamental suivant: \emph{Quel parti peut-on tirer de la connaissance de certaines relations $(A)$ entre $\xn$, en n'employant que des op\'erations rationnelles?} Nous montrons que l'on peut d\'eduire du syst\`eme $(S,A)$ un syst\`eme analogue, dont le syst\`eme $(S,A)$ admet toutes les solutions, et qui est, comme nous le disons, \emph{automorphe}: ce qui veut dire que ses diverses solutions se d\'eduisent de l'une quelconque d'entre elles par les substitutions d'un groupe $G$, qui est dit \emph{le groupe associ\'e au syst\`eme}, ou simplement le \emph{groupe du syst\`eme}. On remarquera que $S$ est d\'ej\`a un syst\`eme automorphe, ayant le groupe g\'en\'eral pour groupe associ\'e. D\`es lors, si l'on se place du point de vue de Galois,
 \dots\ 
on voit que l'on peut se limiter \`a ne consid\'erer que des syst\`emes $(S,A)$ rationnels and automorphes.}%
\footnote{This quote translates as: \gui{Given an \agq \eqnz, that we consider as replaced  by the \sys $(S)$ of relations between the roots $\xn$ and the \coesz,
we first study the following fundamental \pbz: \emph{What subset can we extract from the knowledge of certain relations $(A)$ between $\xn$, by only employing rational operations?} We show that we can deduce from the \sys $(S,A)$ an analogous \sysz, for which the \sys $(S,A)$ admits all the solutions, and which is, as we say, \emph{automorphic}: which means that its diverse solutions are deduced from any one of them by the substitutions of a group $G$, which is said to be \emph{the group associated with the \sysz}, or simply the \emph{group of the \sysz}. We will notice that $S$ is already an automorphic \sysz, with the \gnl group being its associated group. From then on, if we take Galois' point of view,
\dots\ 
we see that we can limit ourselves to only considering rational and automorphic \syss $(S,A)$.}}

\medskip 
The \adu is dealt with in  considerable  detail
 in Chapter 2 of the book \cite[1989]{PZ}.

Among the good modern works that present all of  classical Galois theory, we cite \cite{Tignol} and~\cite{CoxGal}.

The \gui{dynamic Galois theory} presented in detail in this chapter is presented in \cite[D\'{\i}az-Toca]{DiJSC} and \cite[D\'{\i}az-Toca\&al.]{DiLQ,DiL10}. 

Regarding \thref{theoremAdu1} on the fixed points of $\Sn$ in the \aduz, the \gui{$f$ is \splz} case belongs to folklore. 
We find it, with a proof related to the one given here, in
 Lionel Ducos' thesis \cite{Du}. We have given another proof of it in \thref{thResolUniv} for the \cdis case.
The refinement that we give is found in \cite{DiLQ}, it is inspired by \cite{PZ} (see \Tho 2.18 page~46, Corollary 3.6 page 49
 and the following remark, page~50).

\Thref{thidGTri}, published in \cite{DiLQ} under a restrictive hypothesis, \gnss a result given separately in the \adu over a field case by L. Ducos \cite{Duc} and by P. Aubry and A. Valibouze~\cite{AV}. Our  method of proof is closer to that of L. Ducos, but it is different because the framework is more \gnlz: we start off with an arbitrary commutative \riz.

A related version of \thref{theoremAdu3} is found in \cite[lemme~II.4.1]{Du}.

\smallskip \rdb\label{Jordan}\label{soicher}
Regarding the explicit methods of computing Galois groups over $\QQ$ recently developed in Computer Algebra we  refer to \cite[Geissler\&Kl\"uners]{GeKl}.

The modular method, made popular by van der Waerden, is due to Dedekind (letter addressed to Frobenius on June 18, 1882, see \cite[Brandl]{Brandl}).

The Stauduhar \cite{Sta} and Soicher-McKay \cite{SoM} methods are based on computations of resolvents and on the knowledge of the transitive subgroups of the groups $\Sn$. These have been tabulated up to $n=31$ \cite[Hulpke]{Hulpke}. 
In most of the existing \algos the computation determines the Galois group of an \ird \polz, without computing the splitting field.\\
See however \cite[Kl\"uners\&Malle]{KlMa} and \cite[Valibouze\&al.]{AV,ORV,Val}.

Moreover, let us cite the remarkable \poll time computability result \cite[Landau\&Miller]{LaMi} regarding the solvability by radicals.

\smallskip
Alan Steel \cite{steel,steel2} was inspired by D5 to implement a very efficient \gui{dynamic} \agq closure of $\QQ$ in {\tt Magma}. The efficiency depends on him not using a  \fcn \algo for the \pols of $\ZZ[X]$, nor an \algo of representation of the finite extensions by means of primitive \eltsz.
Nevertheless he uses \fcn \algos modulo $p$ to control the process. The process is dynamic in the sense that the progressively constructed closure depends on the user's questions.
The author however does not give (and could not do so in his chosen framework) an implementation of the splitting field of a \pol (let us say \spl for the sake of simplification) over a \gui{\gnlz} field.

\smallskip 
For the Computer Algebra \sys {\tt Magma}, see \cite[Bosma\&al.]{Bosma,Cannon}.

\newpage \thispagestyle{CMcadreseul}
\incrementeexosetprob


\chapter{Flat modules} 
\label{chap mod plats} 
\perso{compil\'e le \today}
\minitoc
\newpage	
\Intro

\begin{flushright}
{\em Dear \eltsz, \\
if you aren't free, \\ 
it isn't my fault.}\\ 
A flat module. \\
\end{flushright}

Flatness is a fundamental notion of commutative \algz,
introduced by Serre in~\cite{Serregaga}.

In this chapter we introduce the notion of a flat module, of a flat \alg and of a \fpte \algz, and prove some of the essential \prts of these objects.

An integral \ri whose \itfs are flat is called a \ddpz.\index{Pru@Pr\"ufer!domain}
This is another fundamental notion of commutative algebra which will only be introduced here. It will be further developed
 in Chapter~\ref{ChapAdpc}.

\section{First \prtsz} \label{secPlatDebut}

\subsec{\Dfn and basic \prtsz}
We give an \elr \dfn and later develop the relationship with the 
exactness of the functor $M\otimes \bullet$. 
\begin{definition} 
\label{def.plat}  Consider an \Amo $M$.
\begin{enumerate}
\item A \emph{syzygy in $M$} is given by $L\in \Ae {1\times n}$ and $X\in M^{n\times 1}$ which satisfy $LX=0$.
\item We say that \emph{the syzygy $LX=0$ is explained in~$M$} 
if we find $Y\in M^{m\times 1}$ and a matrix $G\in \Ae {n\times m}$ that satisfy
\begin{equation}\label{eqdef.plat} 
LG=0\quad {\rm and } \quad GY = X\,.
\end{equation}
%
\item The \Amo $M$ is called a {\em  \mplz}\index{module!flat ---}\index{flat!module} if every syzygy  
in $M$ is explained in $M$.
(Intuitively speaking: if there is a syzygy between \elts of $M$, the module is not to blame.) 
\end{enumerate}
\end{definition}

\rems 1) In items \emph{1} and \emph{2} the symbol $0$ is specified  implicitly by the context. In \emph{1} it is $0_M$, whereas in \emph{2} it is $0_{\Ae {m\times 1}}$. 

2) In item \emph{2}, when we say that the syzygy $LX=0$ is explained in $M$, we mean that the explanation \gui{does not touch $L$.}
However, the \egts given by the matrix \eqn $LG=0$ take place in $\gA$ and not in $M$.
\eoe

\medskip 
\exls 1) If $M$ is free and \tfz,\footnote{Or more \gnlt if $M$ is freely generated by a discrete set, \cad \smash{$M\simeq \gA^{(I)}=\bigoplus_{i\in I}\gA$} with $I$ discrete. For another generalization see Exercise \ref{propfreeplat}.} it is flat: if $LX=0$,
we write~$X=GY$ with a column vector $Y$ which forms a basis, and $LX=0$ implies~$LG=0$.

 2)
If $M=\bigcup_{i\in I}M_i$ with $\forall i,j\in I,\;\exists k\in I,\;
M_k\supseteq M_i\cup M_j$ (we then say that~$M$ is a \emph{filtering union} of the $M_i$'s),
and if each $M_i$ is flat, then $M$ is~flat.
\index{union!filtering ---}
\index{filtering!union}

 3)
Let $a$ be a \ndz \elt in  $\gA$, $M$ be an \Amo and $u\in M$ such that~$au=0$.
If this syzygy is explained in $M$, we write $u=\sum_ia_iu_i$ ($a_i\in\gA$, $u_i\in M$) with each $aa_i=0$, so $u=0$. 
Thus in a flat module, every \elt annihilated by a \ndz \elt is null.  

 4) (Continued) 
The \emph{torsion submodule} of a module $M$ is the module 

\snic{N=\sotq{x\in M}{\exists a\in \Reg(\gA),\;ax=0},} 

where $\Reg(\gA)$ designates the filter of the \ndz \elts of $\gA$. This torsion submodule is the kernel of the morphism of \eds to $\Frac\gA$ for the module $M$.
The torsion submodule of a flat module is reduced to $0$.
\\
\emph{When the \ri $\gA$ is integral}, we say that a module is \emph{torsion-free} 
if its torsion submodule is reduced to $0$. Over a Bézout domain, or more \gnlt over a \ddpz, a module is flat \ssi it is torsion-free (Exercise \ref{exoPlatsLecteur} and \thref{thPruf} item \emph{2b}). 
\\
Later we give a \gnn of the notion of a torsion-free module for an arbitrary commutative \ri (\Dfref{def.locsdz}).%
\index{torsion!submodule}%
\index{torsion-free!module}%
\index{module!torsion-free ---}

 5) We will see (Proposition~\ref{prop.itfplat}) that a \tf flat \id $\fa$ is \lopz, which implies $\Al 2 \fa=0$ (\thrf{propmlm}). Thus, when~$\gA$ is a nontrivial integral \ri and $\gB=\gA[x,y]$, the \id $\fa=\gen{x,y}$ is an
example of a \Bmo that is torsion-free, but not flat (since $\Al 2_\gB \fa=\gA$ by Example on \paref{belexemple}).
 In fact, the relation \smashbot{$[\,y\;-x\,]\cmatrix{x\cr y}=0$} is not explained in $\fa$, but in~$\gB$.   
\eoe


\smallskip 
The following proposition says that the \gui{explanation} which is given for the syzygy  $LX=0$ in the \dfn of a \mpl extends to a finite number of syzygies.
\begin{proposition} 
\label{propPlat1} 
Let $M$ be a flat \Amoz. 
Consider a family of $k$ syzygies, written in the form $LX=0$,  where~$L\in \Ae {k\times n}$ and~$X\in M^{n\times 1}$. Then, we can find an integer~$m$, a vector $Y\in M^{m\times 1}$ and a matrix $G\in \Ae {n\times m}$ satisfying the \egts
$$\preskip.2em \postskip.4em
 GY = X\quad {\rm and } \quad LG=0.
$$
\end{proposition}
\begin{proof}
Let $L_1$, \ldots, $L_k$ be rows of $L$. 
The syzygy $L_1X=0$  is explained by two matrices $G_1$ and $Y_1$ and by two \egtsz~$X = G_1Y_1$ and~$L_1G_1=0$. 
The syzygy $L_2X=0$ is rewritten as $L_2G_1Y_1=0$  \cad $L'_2Y_1=0$. 
This syzygy is explained in the form $Y_1 = G_2Y_2$ and $L'_2G_2=0$.\\ 
Therefore $X = G_1Y_1 = G_1G_2Y_2$. With $L_1G_1G_2=0$  and $L_2G_1G_2=L_2'G_2=0$.
The column vector $Y_2$ and the matrix $H_2=G_1G_2$ therefore explain the two syzygies  $L_1X=0$  and  $L_2X=0$. \\
All that remains is to iterate the process.
\end{proof}

The following theorem reformulates Proposition~\ref{propPlat1} in the language of categories.
 The \dem is a translation exercise left to the reader.

\begin{theorem} 
\label{thPlat1} \emph{(\Carn of  flat modules, 1)}
\perso{des dessins?}\\
For some \Amo $M$ \propeq
\begin{enumerate}
\item The module $M$ is flat.
\item Every \ali from a \mpf $P$ to $M$ is factorized by a free module of finite rank.
\end{enumerate}
\end{theorem}
%
%

\begin{theorem} 
\label{cor pf plat ptf}  
An \Amo $M$ is \pf and flat \ssi it is \ptfz.

\end{theorem}
\begin{proof}
The condition is \ncr by the following remark.
It is sufficient, because the \idt of $M$ is factorized by a free \Amo $L$ of finite rank. Then, the composition $L\to M \to L$ is a  \prn whose image is \isoc to~$M$.
\end{proof}
%

It is \imd that the \Amo $M\oplus N$ is flat \ssi the modules~$M$ and~$N$ are flat.

The following proposition gives a slightly better result (see also \thref{propPlatQuotientdePlat} and Exercise~\ref{propfreeplat}).

\begin{proposition}\label{propSuExPlat}
Let $N\subseteq M$ be two \Amosz. If $N$ and $M/N$ are flat, then $M$ is flat.  
\end{proposition}
\begin{proof}
Write $\ov x$ the object $x$ (defined over $M$) considered modulo $N$.
Consider a syzygy $LX=0$ in $M$. Since $M/N$ is flat, we obtain $G$ over $\gA$ and $Y$ over $M$ such that $LG=0$  and $G\ov Y=\ov X$.
Consider the vector $X'=X-GY$ over $N$. We have $LX'=0$, and since $N$ is flat, we obtain $H$ over $\gA$ and $Z$ over $N$ such that $LH=0$
and $HZ=X-GY$.\\
 Thus the matrix $\blocs{.6}{.5}{.8}{0}{$G$}{$H$}{}{}$
and the vector  $\blocs{.4}{0}{.6}{.5}{$Y$}{}{$Z$}{}$ 
explain the relation~$LX=0$.
\end{proof}
%

\begin{fact} 
\label{fact.plat}  Let $S$ be a \mo of the \ri $\gA$.
\begin{enumerate}
\item The localized ring $\gA_S$ is flat as an \Amoz. 
\item If $M$ is an $\gA$-\mplz, then $M_S$ is flat as an \Amo and as an $\gA_S$-module.
\end{enumerate}

\end{fact}
%
\begin{proof} It suffices to prove item \emph{2.}
If we have a syzygy $LX=0$ in the \Amo $M_S$, we write $X=X'/s$ and we have a syzygy $uLX'=0$ in $M$ (with $u,s\in S$). We therefore find $Y'$ over $M$ and $G$ over $\gA$ such that $GY'=X'$ in $M$ and $uLG=0$ in $\gA$. This implies, for $Y=Y'/(su)$, the \egt $uGY=X$ in $M_S$, such that $uG$ and $Y$ explain the relation $LX=0$ in $M_S$. 
We can construct an analogous proof by starting with
 the syzygy in $M_S$ considered as an $\gA_S$-module. 
\end{proof}

\vspace{-.7em}
\pagebreak	

\subsec{Local-global principle}

Flatness is a local notion in the following sense.

\begin{plcc} 
\label{plcc.plat}  
{\em  (For flat modules)}\\
Let $S_1$, $\ldots$, $S_r$  be \moco of a \ri $\gA$, and let $M$ be \hbox{an \Amoz}.
%
\begin{enumerate}
\item A syzygy $LX=0$ in~$M$ is explained in~$M$ \ssi it is explained in each of the $M_{S_i}\!$'s. 
\item The module $M$ is flat over $\gA$ \ssi each of the $M_{S_i}\!$'s is flat over~$\gA_{S_i}$.
\end{enumerate}
\end{plcc}
\begin{proof} It suffices to prove the first item.
The \gui{only if} is given by Fact~\ref{fact.plat}. \\
Let us prove the other implication. Let $LX=0$ be a syzygy between \elts of~$M$ (where~$L\in \Ae {1\times n}$  and $X\in M^{n\times 1}$). 
We want to find $m\in\NN,\; Y\in M^{m\times 1}$ and a matrix~$G\in \Ae {n\times m}$ which satisfy \Eqrf{eqdef.plat}.
We have a solution~$(m_i,Y_i,G_i)$  for~$(\ref{eqdef.plat})$ in each 
localized \ri $\gA_{S_i}$. 
\\
We can write~$Y_i=Z_i/s_i$,~$G_i=H_i/s_i$ with~$Z_i\in M^{m_i\times 1}$, $H_i\in \Ae {n\times m_i}$ and some~$s_i\in S_i$ that are suitable. 
We then have~$u_iH_iZ_i=v_iX$ in $M$ and~$u_iLH_i=0$ in $\gA$ for some~$u_i$ and~$v_i\in S_i$. 
We write~$\sum_{i=1}^{r} b_i v_i =1$ in~$\gA$. 
For~$G$ we take the matrix obtained by juxtaposing the matrices~$b_i u_iH_i$ in a row, and for $Y$ we take the vector obtained by superposing the vectors~$Z_i$ in a column.  
We obtain $ GY = \sum_{i=1}^rb_i v_iX=X$ in $M$,  and $LG=0$  in~$\gA$. 
\end{proof}

The corresponding principle in \clama is the following.

\begin{plca} 
\label{plca.plat}  
{\em  (For flat modules)}
\begin{enumerate}

\item A syzygy $LX=0$ in $M$ is explained in~$M$ \ssi it is explained in $M_{\fm}$ for every \idemaz~$\fm$. 
\item An \Amo $M$ is flat \ssi for every \idemaz~$\fm$, the module $M_{\fm}$ is 
 flat over $\gA_{\fm}$.
\end{enumerate}
\end{plca}
\begin{proof}
It suffices to show the first item. However, the fact that a syzygy $LX=0$ can be explained is a \carf \prt (\Dfref{defiPropCarFini}). 
 We therefore apply Fact~\ref{factPropCarFin} which allows us to pass from the \plgc to the corresponding \plgaz.
\end{proof}

\subsec{Other \carns of flatness}

We will now consider \emph{syzygies over $M$ with \coes in another module $N$} and we will show that when $M$ is flat,
every syzygy with \coes in any module $N$ is explained in~$M$.

\pagebreak	

\begin{definition}\label{defiRDLGen}
Let $M$ and $N$ be two \Amosz.  \\
For $L=[\,a_1\,\cdots\,a_n\,]\in N^{1\times n}$ and $X=\tra[\,x_1\,\cdots\,x_n\,]\in M^{n\times 1}$,
let
$$\preskip.2em \postskip.2em \ndsp
L\odot X\;\eqdefi\;\sum_{i=1}^na_i\te x_i\;\in N\te M. 
$$

\begin{enumerate}
\item If $L\odot X=0$ we say that we have a syzygy between the $x_i$'s with \coes in~$N$.
\item We say that \emph{the syzygy $L\odot X=0$ is explained in~$M$} if we have $Y\in M^{m\times 1}$ and a matrix $G\in \Ae {n\times m}$ which satisfy
\begin{equation}\label{eqdef2.plat}\preskip.3em \postskip.4em 
LG =_{N^ {1\times m}}0 \quad {\rm and } \quad X=_{M^{n\times 1}}GY\,.
\end{equation}
%
\end{enumerate}
\end{definition}

\rem 1) When we say that the syzygy $LX=0$  is explained in $M$, we mean that the explanation \gui{does not touch $L$.}

2) We  note that in \gnl the \egt $L\odot G Y=LG\odot Y$ is assured for every matrix $G$ with \coes in $\gA$ because $a\otimes \alpha y=a\alpha\otimes y$ when~$a\in N$, $y\in M$ and~$\alpha\in\gA$.
\eoe

\begin{proposition}\label{propPlatRdl}
Let $M$ and $N$ be two \Amosz. 
\\
If $M$ is a flat \Amo every syzygy with \coes in $N$ is explained in~$M$.  
\end{proposition}
\begin{proof}
We assume that we are given a syzygy $L\odot X=0$ with~$L=[\,a_1\,\cdots\,a_n\,]\in N^{1\times n}$ and $X=\tra[\,x_1\,\cdots\,x_n\,]\in M^{n\times 1}$.

\emph{Case where $N$ is free of finite rank.} Proposition~\ref{propPlat1} gives the result.

\emph{Case where $N$ is \pfz.} \\
Write $N=P/R=\gA^k\sur{(\gA c_1+\cdots\gA c_r)}$.
The $a_i$'s are given by the~$b_i$'s of~$P$. The relation $L\odot X=0$ means that $\sum_i b_i\otimes x_i\in R\otimes M\subseteq P\otimes M$, \cad we have an \egt
$$\preskip-.2em \postskip.4em\ndsp 
\sum_i b_i\otimes x_i+\sum_\ell c_\ell\otimes z_{\ell}=0 
$$ 
in $P\otimes M$. We then observe that when we explain in $M$ this syzygy (regarding the $x_i$'s and the $z_\ell$'s)
 with \coes in the free module $P$, we explain at the same time the syzygy $L\odot X=0$ with \coes in~$N$.

\emph{Case of an arbitrary \Amo $N$.} 
\\
A relation $L\odot X=\sum_i a_i\otimes x_i=0$ comes from a finite computation, in which only a finite number of \elts of $N$ and of relations between these \elts intervene.
There exist therefore a \mpf $N'$, a \ali $\varphi: N'\to N$ and some $b_i\in N'$ such that on the one hand $\varphi(b_i)=a_i$ ($i\in\lrbn$), and on the other hand $\sum_i b_i\otimes x_i=0$ in $N'\otimes M$.
We then observe that when, in $M$, we explain this syzygy
with \coes in $N'$ (which is a \mpfz), we explain at the same time the syzygy $L\odot X=0$ with \coes in~$N$. 
\end{proof}

\vspace{-.7em}
\pagebreak

\begin{theorem} \emph{(\Carn of  flat modules, 2)}
\label{thplatTens} \\
For an \Amo $M$  \propeq
\begin{enumerate}
\item \label{i1thplatTens} 
The module $M$ is flat.
\item \label{i2thplatTens} 
  For all \Amos $N$, every syzygy between \elts of $M$ with \coes in $N$ is explained in~$M$.
\item \label{i3thplatTens} 
  For every \itf $\fb$ of $\gA$ the canonical map $\fb\otimes_\gA M\rightarrow M$ is injective (this therefore establishes an \iso from $\fb\otimes_\gA M $ to~$\fb M$).
\item \label{i4thplatTens} 
  For all \Amos  $N\subseteq N'$, the canonical \ali 
$$\preskip.3em \postskip.0em 
N\otimes_\gA M\rightarrow N'\otimes_\gA M 
$$  
is injective.
\item \label{i5thplatTens} 
  The functor $\bullet\otimes M$ preserves exact sequences.
\end{enumerate}
\end{theorem}
\begin{proof}
The implication \emph{\ref{i5thplatTens}}  $\Rightarrow$ \emph{\ref{i3thplatTens}} is trivial. 

\emph{\ref{i4thplatTens}}  $\Rightarrow$ \emph{\ref{i5thplatTens}.}
Short exact sequences are preserved by the functor~$\bullet\otimes M$. 
However every exact sequence decomposes into short exact sequences (see \paref{sexaseco}).

 \emph{\ref{i1thplatTens}} $\Leftrightarrow$ \emph{\ref{i3thplatTens}.} By the null tensor lemma \ref{lem-tenul}.

 \emph{\ref{i1thplatTens}} $\Rightarrow$ \emph{\ref{i2thplatTens}.} 
This is Proposition~\ref{propPlatRdl}.

 \emph{\ref{i2thplatTens}} $\Leftrightarrow$ \emph{\ref{i4thplatTens}.}
By the null tensor lemma \ref{lem-tenul}.
\end{proof}

The previous \tho admits some important corollaries. 
\begin{corollary}\label{corPlatTens} \emph{(Tensor product)}\\
 The tensor product of two \mpls is a \mplz.
\end{corollary}
\begin{proof}
Use item \emph{\ref{i4thplatTens}} of \thref{thplatTens}.
\end{proof}
%

\begin{corollary}\label{cor2PlatTens} \emph{(Other basic constructions)}\\ 
 The tensor, exterior and \smq powers of a \mpl are \mplsz.
\end{corollary}
\facile
%

\begin{corollary}\label{cor3PlatTens} \emph{(Intersection)}\\
Let $N_1,\ldots, N_r$ be submodules of a module $N$ and $M$ be a \mplz.
Since $M$ is flat, for $N'\subseteq N$, we identify $N'\otimes M$
with its image in $N\otimes M$.
Then we have the \egt 
$$\preskip.0em \postskip.2em\ndsp 
\;\;\;\;\left(\bigcap_{i=1}^r N_i\right) \otimes M= \bigcap_{i=1}^r (N_i\otimes M). 
$$ 
\end{corollary}
\begin{proof}
The exact sequence
$$\preskip.0em \postskip.2em\ndsp 
\;\;\;\;0\to \bigcap_{i=1}^r N_i\to N\to \bigoplus_{i=1}^r(N\sur {N_i}) 
$$
is preserved by the tensor product with $M$ and the module $(N\sur {N_i})\te M$ is identified with $(N\te M)\sur{(N_i\te M)}$.
\end{proof}

\vspace{-.7em}
\pagebreak

\begin{corollary}\label{corPlatEds} \emph{(\Edsz)}
Let $\rho:\gA\to\gB$ be an \algz. If~$M$ is a flat \Amoz, then $\rho\ist (M)$ is a flat \Bmoz. 
\end{corollary}
\begin{proof} 
Note that for a \Bmoz~$N$, we have 
$$\preskip.2em \postskip.3em 
N\te_\gB\rho\ist (M)\simeq N\te_\gB\gB\te_\gA M\simeq N\te_\gA M. 
$$
We then apply item \emph{4} of \thref{thplatTens}.
Note that the last tensor product is equipped with a  \Bmo structure via~$N$.
\end{proof}

\medskip 
\rem Without mentioning it, we have just used a \gnee form of associativity of the tensor product whose \dem we leave to the reader.
The form in question is the following.\\
First we say that an abelian group $P$ is an $(\gA,\gB)$-bimodule if it is equipped with two external laws which respectively make an \Amo and {a \Bmoz}, and if these two structures are compatible in the following sense: for all $a\in\gA$, $b\in\gB$ and $x\in P$, we have $a(bx)=b(ax)$.\\
In such a case, if $M$ is a \Bmoz, then the tensor product $M\te_\gB P$ can \hbox{itself} be equipped with a structure of an $(\gA,\gB)$-bimodule  by letting, for $a\in\gA$, $a(x\te y)=_{M\te_\gB P}x\te ay$.\\
Similarly, when $N$ is an \Amoz, the tensor product $P\te_\gA N$ can \hbox{itself} be equipped with a structure of an $(\gA,\gB)$-bimodule  by letting, \hbox{for $b\in\gB$}, $b(y\te z)=_{P\te_\gA N}by\te z$.
\\
Finally, under these hypotheses, there exists a unique \ali (for the  structure of an $(\gA,\gB)$-bimodule) $\varphi:(M\te_\gB P)\te_\gA N\to M\te_\gB (P\te_\gA N)$ which satisfies 
$$\preskip-.1em \postskip.4em 
\varphi\big((x\te y)\te z\big)=x\te (y\te z) 
$$
for all $x\in M$, $y\in P$, $z\in N$, and $\varphi$ is an \isoz.
\eoe

\subsec{Flat quotients}

\begin{theorem}\label{propPlatQuotientdePlat} \emph{(Flat quotients)}\\
Let $M$ be an \Amoz, $K$ be a submodule and $N=M/K$, with the exact sequence
$$\preskip-.2em \postskip.3em
0 \to K\vers\imath M\vers \pi N \to 0. 
$$ 
\begin{enumerate}
\item If $N$ is flat, for every module $P$, the sequence
$$\preskip.3em \postskip.3em
0 \to K\te P\vers{\imath_P} M\te P\vers {\pi_P} N\te P \to 0 
$$ 
is exact ($\imath_P=\imath\te\rI_P$, $\pi_P=\pi\te\rI_P$).
\item If $N$ and $M$ are flat, $K$ is flat.
\item If $N$ and $K$ are flat, $M$ is flat.
\item If $M$ is flat, \propeq
\begin{enumerate}
\item $N$ is flat.
\item For every \itf $\fa$, we have $\fa M \cap K = \fa K$.

\vspace{-.2em}
\pagebreak	

\item Every \itf $\fa$ gives an exact sequence
$$\preskip.2em \postskip.4em
0 \to K/\fa K\vers{\imath_\fa} M/\fa M\vers {\pi_\fa} N/\fa N\to 0. 
$$ 
\end{enumerate}
\end{enumerate}
\end{theorem}
%
NB: Item \emph{3} has already been the object of Proposition~\ref{propSuExPlat}. Here we give it another \demz, leaving it up to the reader to compare them.
\begin{proof} 
\emph{1.} \emph{Case where $P$ is \tfz}. We write $P$ as a quotient of a finite free module $Q$ with a short exact sequence
$$\preskip.2em \postskip.4em
0 \to R\vers a Q\vers p P \to 0. 
$$ 
We then consider the following commutative diagram in which all the horizontal and vertical sequences are exact because $N$ and $Q$ are flat

\smallskip 
\centerline{\small
\xymatrix@C=3em@R=1.6em {
&   &     & 0\ar[d]
\\
 &K\te R\ar[d]_{a_K} \ar[r]^{\imath_R} 
&M\te R\ar[d]_{a_M}\ar[r]^{\pi_R} &
N\te R\ar[d]_{a_N}\ar[r]&0
\\
0\ar[r]&K\te Q\ar[d]_{p_K} \ar[r]^{\imath_Q} 
&M\te Q\ar[d]_{p_M}\ar[r]^{\pi_Q} &
N\te Q\ar[r]&0
\\
 &K\te P\ar[d] \ar[r]^{\imath_P} &M\te P\ar[d]  & & 
\\
&0  &  0  &  
\\
}}

We must show that $\imath_P$ is injective. This is a special case of the snake lemma, which we can prove  by \gui{diagram chasing.} \\
Suppose $\imath_P(x)=0$. We write $x=p_K(y)$ and $v=\imath_Q(y)$.
We have $p_M(v)=0$, so we write $v=a_M(z)$. 
\\
As $\pi_Q(v)=0$, we have 
$a_N(\pi_R(z))=0$, so $\pi_R(z)=0$. \\
Therefore we write $z=\imath_R(u)$ and we have 
$$\preskip.3em \postskip.4em 
\imath_Q(a_K(u))=a_M(\imath_R(u))=a_M(z)=v=\imath_Q(y), 
$$
and since $\imath_Q$ is injective, $y=a_K(u)$, hence $x=p_K(y)=p_K(a_K(u))=0$.

\emph{\Gnl case.} One possibility is to describe $P$ as a quotient of a flat module $Q$ (see Exercise~\ref{propfreeplat} 
on the subject) 
in which case the previous \dem remains unchanged. We can also do without this slightly cumbersome construction as follows. Let us show that $\imath_P$ is injective. Let $x=\sum_ix_i\te y_i\in K\te P$ such that 
$\imath_P(x)=_{M\te P}0$, \cad $\sum_ix_i\te y_i=_{M\te P}0$.
\\
By \dfn of the tensor product, there exists a \smtf $P_1\subseteq P$ such that we also have $\sum_ix_i\te y_i=_{M\te P_1}0$. By the already examined case, we have $\sum_ix_i\te y_i=_{K\te P_1}0$, and this implies $\sum_ix_i\te y_i=_{K\te P}0$.

\smallskip 
\emph{2} and \emph{3.} Let $\fa$ be an arbitrary \itfz. Since $N$ is flat, we have by item~\emph{1} a commutative diagram with exact sequences
$$\preskip.4em \postskip.4em 
\xymatrix@C=3em@R=1.6em {
  &  & &0\ar[d] & 
\\
0\ar[r]&\fa\te K\ar[d]^{\varphi_K} \ar[r]^{\imath_\fa} &\fa\te M\ar[d]^{\varphi_M}\ar[r]^{\pi_\fa} &\fa\te N\ar[d]^{\varphi_N}\ar[r] & 0
\\
0\ar[r] &K\ar[r]^{\imath} &M\ar[r]^\pi &N\ar[r] & 0.
\\
}
$$

If $M$ is flat, $\varphi_M$ is injective, hence so is $\varphi_M\circ \imath_\fa$, then $\varphi_K$. By item \emph{3} of \thref{thplatTens} we conclude that $K$ is flat.

If $K$ is flat, $\varphi_K$ is injective and a short diagram chase shows that $\varphi_M$ is injective. Let $x\in\fa\te M$ with $\varphi_M(x)=0$. 
As $\varphi_N(\pi_\fa(x))=0$, we have $\pi_\fa(x)=0$
and we can write $x=\imath_\fa(y)$. 
Then $\imath(\varphi_K(y))=\varphi_M(x)=0$, \hbox{so $y=0$}, thus $x=0$.

\emph{4a} $\Rightarrow$ \emph{4b.} Since $M$ and $N$ are flat, so is $K$ and the top row 
of the previous diagram gives the exact sequence
$$\preskip.2em \postskip.4em
0 \to \fa K
\vvers{{\imath\,\frt{\fa\! K}}} 
\fa M
\vvers {{\pi\,\frt{\fa\! M}}} \fa N\to 0.\eqno(+) 
$$
 However, the kernel of  $\pi\frt{\fa\! M}$ is by \dfn $\fa M\cap K$.

\emph{4b} $\Leftrightarrow$ \emph{4c.} The sequence
$$\preskip.2em \postskip.4em
0 \to K/\fa K\vers{\imath_\fa} M/\fa M\vers {\pi_\fa} N/\fa N\to 0  
$$ 
is obtained from the exact sequence $0\to K \to M \to N$
by \eds to~$\gA/\fa$. Saying that it is exact is the same as saying that $\imath_\fa$ is injective. However, an \elt
$\ov x\in K/\fa K$ is sent to $0$ \ssi {we have $x\in \fa M\cap K$}.

\emph{4b} $\Rightarrow$ \emph{4a.}
Since $\fa K=\fa M\cap K$ the sequence $(+)$ is exact. Consider the following commutative diagram with exact sequences, for which we must show that $\varphi_N$ is injective.
$$
\preskip-.4em \postskip.6em 
\xymatrix@C=3em@R=1.6em {  &  & 0\ar[d]& & 
\\
&\fa\te K\ar[d]^{\varphi_K} \ar[r]^{\imath_\fa} &\fa\te M\ar[d]^{\varphi_M}\ar[r]^{\pi_\fa} &\fa\te N\ar[d]^{\varphi_N}\ar[r] & 0
\\
0\ar[r] &\fa K\ar[d]\ar[r]^{\imath\,\frt{\fa\! K}} &\fa M\ar[d]\ar[r]^{\pi\,\frt{\fa\! M}} &\fa N\ar[r]\ar[d] & 0
\\
& 0& 0&0
\\
} 
$$
This is obtained by a short diagram chase. 
If $\varphi_N(x)=0$, we write $x=\pi_\fa(y)$. As $\pi\,\frt{\fa\! M}(\varphi_M(y))=0$, we \hbox{have $z\in\fa K$} such that $\varphi_M(y)=\imath\,\frt{\fa\! K}(z)$, we write $z=\varphi_K(u)$, with $\varphi_M(\imath_\fa(u))=\varphi_M(y)$, and since $\varphi_M$ is injective, $y=\imath_\fa(u)$ and $x=\pi_\fa(y)=0$.
\end{proof}
\vspace{-.7em}
\pagebreak	

\begin{corollary}\label{corpropPlatQuotientdePlat} \emph{(A flat \algz)}
Let $f\in\AuX = \gA[\Xn]$ and~$\Aux=\aqo\AuX{f}$. The \Amo $\Aux$
is flat \ssiz $\rc(f)^2 = \rc(f)$, \cad \ssi the \id $\rc(f)$ is generated by an \idmz.
\end{corollary}
%
\begin{proof}
The \Amo
 $\Aux$ is flat \ssi for every \itfz~$\fa$ of $\gA$ we have
 {
\fbox{$\gen {f} \cap \fa[\uX] = f\fa[\uX]\;(*)$}.
}

If $\Aux$ is flat, we obtain, for $\fa = \rc(f)$, 
\hbox{that $\rc(f)^2 = \rc(f)$}, 
because~$f \in \gen {f} \cap \fa[\uX]$.

Conversely, let us suppose that $\rc(f)^2 = \rc(f)$ and show that $\Aux$ is flat. The \idm $e$ such that $\gen{e}=\gen{\rc(f)}$ splits the \ri into two components. 
 In the first we have $f=0$ and the result is clear. In the second, $f$ is primitive.  
Now suppose that $f$ is primitive. 
\\
By the
Dedekind-Mertens lemma,\footnote{Actually, this refers to a variant, with essentially the same \demz, which we leave to the reader.} for every \Amoz~$M$ the \Ali
$M[\uX] \vvers {\times f} M[\uX]$ is injective. Applied to
$M = \gA/\fa$, this gives $(*)$. Indeed, let $M[\uX] = \AuX/\fa[\uX]$
and suppose 
that $g \in \gen {f} \cap \fa[\uX]$. Then $g = fh$ for some $h\in \AuX$, and $\overline h$ is in the kernel of $\AuX/\fa[\uX] \vvers {\times f} \AuX/\fa[\uX]$, therefore $\overline h = 0$, \cad $h \in \fa[\uX]$, and $g \in f\fa[\uX]$. 
\end{proof}
%

%
\section{\Tf flat modules}

In the \mtf case, flatness is a more \elr \prtz.

\begin{lemma} 
\label{lem.plat-tf} 
Consider a \tf \Amo $M$, and let $X\in M^{n\times 1}$ be a column vector whose \coos $x_i$ generate $M.$ 
The module $M$ is flat \ssi for every syzygy $LX=0$  (where  $L\in \Ae {1\times n}$), we can find two matrices $G$, $H\in \Mn(\gA)$ which satisfy the \egts
$$\preskip.2em \postskip.3em 
H+G = \I_n,\;\; LG=0\;\; {\rm and } \;\; HX = 0. 
$$
In particular, a cyclic module $M=\gA y$ is flat \ssi 
$$\preskip.3em \postskip.4em
\forall a\in\gA, (\,ay=0\;\Longrightarrow\;\exists s\in\gA,\;\;as=0 \hbox{ and } sy=y\,)\,.
$$ 
\end{lemma}
\rem The symmetry between $L$ and $X$ in the statement is only apparent; the module $M$ is generated by the \coos of $X$, while the \ri $\gA$ is not generated (as a submodule) by the \coos of $L$.  
\eoe
\begin{proof}
We reduce an arbitrary syzygy $L'X'=0$ to a syzygy $LX=0$ by expressing $X'$ in terms of $X$. A priori we should write $X$ in the form $G_1Y$ with 
$LG_1=0$. 
\\
As $Y=G_2X$, we take $G=G_1G_2$ and $H=\I_n-G$.
\end{proof}

\rem For cyclic modules, by letting $t=1-s$, we obtain conditions on $t$ rather than on $s$
$$\preskip.0em \postskip.4em
a=at \qquad   \hbox{and} \qquad  ty=0,$$
which implies that the annihilator $\fa$ of $y$ satisfies $\fa^2=\fa$.
In fact, by \thref{propPlatQuotientdePlat}, $\gA\sur\fa$ is flat over $\gA$ \ssi for every \itf $\fb$ we have the \egt $\fa\cap\fb=\fa\fb$.
\eoe


\smallskip  Here is a \gnn of Lemma~\ref{lem.plat-tf} in the same style of Proposition~\ref{propPlat1}.
\begin{proposition} 
\label{propPlat2} 
Let $M$ be a \tf flat \Amoz, and  $X\in M^{n\times 1}$ be a column vector that generates~$M.$  
Let there be a family of $k$ syzygies expressed in the form $LX=0$ where~$L\in \Ae {k\times n}$ and $X\in M^{n\times 1}$. Then, we can find a matrix~$G\in \Mn(\gA)$ which satisfies the \egts
$$\preskip.2em \postskip.4em 
LG=0\;\; \hbox{ and } \;\; GX = X. 
$$
\end{proposition}
\begin{proof}
Identical to the \dem of Proposition~\ref{propPlat1}.
\end{proof}

A \cof substitute for the \prt according to which every \evc over a field admits a basis (only true in \clamaz) is the fact that every \evc over a discrete field is flat.
More \prmtz, we have the following result.

\begin{theorem}\label{propEvcPlat} \Propeq
\begin{enumerate}
  \item \label{i1propEvcPlat}  Every \Amo $\aqo\gA a$ is flat.
  \item \label{i2propEvcPlat}  Every \Amo is flat.
  \item \label{i3propEvcPlat}  The \ri $\gA$ is reduced \zedz.
\end{enumerate}
\end{theorem}
\begin{proof}
\emph{\ref{i1propEvcPlat}} $\Rightarrow$ \emph{\ref{i3propEvcPlat}}. 
If $\aqo{\gA}{a}$ is flat, then $\gen{a} = \gen{a}^2$ and if it is true for every $a$, then~$\gA$ is reduced \zedz.
\\
\emph{\ref{i3propEvcPlat}} $\Rightarrow$ \emph{\ref{i2propEvcPlat}}. Let us first treat the case of a \cdiz. 
\\
Consider a syzygy $LX=a_1x_1+\cdots+a_nx_n=0$ for some \elts $x_1$, \ldots, $x_n$ of an \Amo $M$. If all the~$a_i$'s are null the relation is explained with $Y=X$ and~$G=\I_n$: $LG=0$ and~$GY=X$.
If one of the~$a_i$'s is \ivz, for instance $a_1$, let $b_j=-a_1^{-1}a_j$ for $j\neq1$. We have $x_1=b_2x_2+\cdots+b_nx_n$ and $a_1b_j+a_j=0$ for $j>1$.
The syzygy is explained by $Y=\tra [\, x_2\; \cdots\; x_n\, ] $ and by the following matrix $G$ because $LG=0$ and~$GY=X$,

\snic{G=\cmatrix{
b_2&b_3&\cdots&b_n\cr
1& 0&\cdots&0\cr
0&\ddots& &\vdots\cr
\vdots& &\ddots&0\cr
0&\cdots &0&1
}.
}


For a reduced \zed \riz, we apply the \elgbm \num2 (\paref{MethodeZedRed}) which brings us back to the case of a \cdiz.
\end{proof}

\vspace{-.5em}
\pagebreak	

NB: This justifies the term \gui{absolutely flat} for reduced \zedz.


\begin{lemma} 
\label{lem.plat2-tf} 
Same context as in Lemma~\ref{lem.plat-tf}.
 If $\gA$ is a \alo and  $M$ is flat, we obtain under the hypothesis $LX=0$ the following alternative. The vector $L$ is null, or one of the $x_i$'s \lint depends on the others (it can therefore be deleted from the list of \gtrs of~$M$).
\end{lemma}
%
\begin{proof}
This is a \gui{\deter trick.} We note that $\det(G)=\det(\I_n-H)$  can be written as $1+\sum_{i,j}b_{i,j} h_{i,j}$. Therefore $\det(G)$ or 
one of the $h_{i,j}$'s is \ivz. In the first case $L=0$. In the second,  since $HX=0$, one of the vectors $x_i$ is a \coli  
of the others.
\end{proof}

The same \dem in the case of an arbitrary \ri gives the following result.


\begin{lemma} 
\label{lem.plat3-tf} 
Same context as in Lemma~\ref{lem.plat-tf}.
 If $M$ is flat and $LX=0$, there exist \eco $s_1$, \dots, $s_\ell$ such that over each of the \ris \smash{$\gA[1/{s_j}]$} we have $L=0$, or one of the $x_i$'s is a \coli of the others.
\end{lemma}

In \clamaz, Lemma~\ref{lem.plat2-tf} implies the following fact. 
\begin{factc}\label{factplatlocalclass}
A \tf \mpl over a \alo is free and a basis can be extracted from any \sgrz. 
\end{factc}

From Lemma~\ref{lem.plat3-tf}, we obtain the following.

\begin{factc}\label{factplattfclass}
A \tf \mpl over an integral \ri is \ptfz. 
\end{factc}

Here is a \cov version of Fact$\etl$~\ref{factplatlocalclass}.

\begin{proposition} 
\label{prop.mtf loc plat1} 
Let  $\gA$ be a \alo and $M$ be a  flat \Amo generated by $(\xn)$. 
Suppose that $M$ is \fdi or that the existence of nontrivial syzygies is explicit in $M$.
Then, $M$ is freely generated by a finite sequence $(x_{i_1},\ldots,x_{i_k})$ 
(with $k\geq 0$).
\end{proposition}
\begin{proof}
First suppose that $M$ is \fdiz, we can then find a finite sequence of integers $1\leq i_1<\cdots<i_{k}\leq n$ (where $k\geq 0$) such that none of the~$x_{i_\ell}$'s is a \coli of the others, and~$(x_{i_1},\ldots,x_{i_k})$ generates~$M$. 
To simplify the notation, suppose from now on that $k=n$, \cad none of the~$x_i$'s is a \coli of the others. Lemma~\ref{lem.plat2-tf}  then tells us that every syzygy between the~$x_i$'s is trivial.
\\
Now suppose that the existence of nontrivial syzygies is explicit in $M$, \cad for every family of \elts of $M$, we know how to tell whether there is a nontrivial syzygy between these \elts and how to provide one if necessary. Then, by using Lemma~\ref{lem.plat2-tf} we can delete the superfluous \elts one after the other in the $(x_i)$ family without changing the module $M$, until all that remains is a subfamily without a nontrivial syzygy (a limiting case is provided by the empty subset when the module is null).
\end{proof}
\comm Note that the \dem uniquely uses the hypothesis \gui{$M$ is \fdiz,} or \gui{the existence of nontrivial syzygies is explicit in $M$} with families extracted from the \sgr $(x_i)$.
Moreover, each of these hypotheses is trivially true in \clamaz. 
\eoe

\smallskip 
Now here is a \cov version of Fact$\etl$~\ref{factplattfclass}.

\begin{proposition} 
\label{prop.mtf integre plat} 
Let  $\gA$ be an integral \ri and $M$ be a  flat \Amo generated by $(\xn)$. 
Suppose that for every finite subset~$J$ of $\lrbn$ the existence of nontrivial syzygies between $(x_j)_{j\in J}$ is explicit in $M$ (in other words, by passing to the quotient field we obtain a finite dimensional \evcz).
Then, $M$ is \ptfz.
\end{proposition}
\begin{proof}
Suppose \spdg that $\gA$ is nontrivial. By using Lemma~\ref{lem.plat3-tf} we obtain the following alternative. 
Either $(\xn)$ is a basis, or after \lon at \eco the module is generated 
by $n-1$ of the $x_j$'s. We conclude by \recu on $n$: indeed, the syzygies after \lon at $s$
with $s\neq 0$ are the same as those over $\gA$. 
\\
Note that for $n=1$, either $(x_1)$ is a basis, or $x_1=0$.
\end{proof}

\section{Flat principal \ids}

A \ri $\gA$ is said to be {\em  \sdzz} if we have:\index{without zerodivisors!ring ---}
\index{ring!without zerodivisors}
\begin{equation}\label{eqSDZ}\preskip.4em \postskip.4em
\forall a,b\in \gA\quad \big(ab=0 \; \Rightarrow \; (a=0\;{\rm or}\;b=0)\big)
\end{equation}
An integral \ri (in particular a \cdiz) is \sdzz.
A discrete \ri \sdz is integral.
A nontrivial \ri is integral \ssi it is discrete and \sdzz.

\begin{lemma} \emph{(When a principal \id is flat)}
\label{lemIdps} ~
\begin{enumerate}
\item  A \idpz, or more \gnlt a cyclic \Amo $\gA  a$, is a \mpl \ssi 
$$\preskip.1em \postskip.2em 
\forall x\in\gA\quad \big(xa=0\;\Rightarrow\; \exists 
z\in\gA\;(za=0\;\;\hbox{and} 
\;\; xz=x)\big). 
$$%
\item  If $\gA$ is local, an \Amo  $\gA a$ 
is flat \ssi
$$\preskip.1em \postskip.2em
\forall x\in\gA\quad \big(xa=0\;\Rightarrow\; (x=0 \;\;\hbox{or}\;\; 
a=0)\big).
$$
\item   Let $\gA$ be a \aloz, if $\gA$ is discrete, or if we have a test to answer the question \gui{is $x$ \ndzz?,} 
then, an \id $\gen{a}$ is flat \ssi $a$ is null or \ndzz.
\item For a local \ri $\gA$ \propeq
\begin{enumerate}
\item    Every \idp is flat.
\item   The \ri is \sdzz.
\end{enumerate}
\end{enumerate}
\end{lemma}
\begin{proof}
Lemma~\ref{lem.plat-tf} gives item \emph{1}. The computation for item~\emph{2} results from it, because~$z$ or $1-z$ is \ivz. The rest is clear. 
\end{proof}
%
We similarly have the following \eqvcsz.
%
\begin{lemma} 
\label{corlemIdps} 
For a \ri $\gA$, \propeq
\begin{enumerate}
\item Every \idp of $\gA$ is flat. 
\item If $xy=0$,  we have
       $\Ann\,x+\Ann\,y=\gA.$
\item If $xy=0$, there exist \moco $S_{i}$ such that in each of the localized \ris $\gA_{S_{i}}$, $x$ or $y$ becomes null. 
\item If $xy=0$, there exists a $z\in\gA$ with $zy=0$ and $xz=x$. 
\item For all $x,\,y \in \gA, \quad \Ann\,xy =  \Ann\,x+\Ann\,y $.
\end{enumerate}
\end{lemma}


The \prt for a \ri to be without zerodivisors behaves badly under patching and that for a module to be flat is well-behaved under localization and patching.
 This justifies the following \dfnz.

\begin{definition} 
\label{def.locsdz}~
\begin{enumerate}
\item  A \ri $\gA$ is said to be {\em a \lsdz} (principal ideals are flat)  
when it satisfies the \eqves \prts of Lemma~\ref{corlemIdps}.
\item An \Amo $M$ is said to be \emph{torsion-free} when all of its cyclic submodules are flat (see Lemma~\ref{lemIdps}). 
\end{enumerate}
\index{pf-ring}%
\index{torsion-free!module}%
\index{module!torsion-free ---}%
\end{definition}

\rems
\\
1) The torsion submodule of a torsion-free module is reduced to $0$.
Our \dfn is therefore a little more constraining 
than the more usual definition, which
 says that a module is torsion-free when its torsion module is reduced to $0$.
We will note that the two \dfns coincide when the \ri $\gA$ is a \qiriz.
 
2) Every submodule of a torsion-free module is torsion-free, which is not the case in \gnl when we replace \gui{torsion-free} by \gui{flat.}
 
\rdb
3) In the French literature, the term \gui{locally without zerodivisors} is often used for a pf-ring.%
\index{pf-ring}%
\index{locally!ring --- without zerodivisors}%
\index{ring!locally without zerodivisors (or pf-\riz)}%

4) A \lsd is reduced. 
 
5) A local \ri is a \lsd \ssi it is \sdzz. 
 
6) The field of reals  {\em is not} \sdz (\emph{nor} a \lsdz): it is a local \ri for which we do not know how to explicitly perform the implication \pref{eqSDZ} on \paref{eqSDZ}.

7) In \clama a \ri is a \lsd \ssi it becomes integral after \lon at every \idep (Exercise~\ref{exoClamlsdz}). 
\eoe

\vspace{-.2em}
\pagebreak

\begin{lemma} 
\label{lem.platsdz}  Let $\gA$ be a  \lsd  
and $M$ be a flat \Amoz.
\begin{enumerate}
\item  The module $M$ is torsion-free.
\item  The annihilator $(0:y)$ of any $y\in M$ is \idmz.
\end{enumerate}
\end{lemma}
\begin{proof}
\emph{1.} Suppose $ay=0$, $a\in\gA$, $y\in M$.
Since $M$ is flat we have \elts $x_i$ of $M$, \elts $b_i$ of $\gA$, and an \egt $y=\sum_{i=1}^n b_i x_i$ in~$M$, with $ab_i=0$ ($i\in\lrbn$) in $\gA$.  
\\
For each $i$, since $ab_i=0$, there exists a $c_i$ such that $ac_i=a$ and $c_ib_i=0$.
 Let $c=c_1\cdots c_{n}$. Then, $a=ca$ and $cy=0$.
 
\emph{2.} Indeed, when $ay=0$, then $a=ca$ with $c\in(0:y)$.
\end{proof}


Using item \emph{2} of Lemma~\ref{lem.platsdz} and the fact that an \idm \itf is generated by an \idm (Lemma~\ref{lem.ide.idem}) we obtain the following result.

\begin{fact}\label{factLsdzCo}
Let $\gA$ be a \ri in which the annihilator of every \elt is \tfz.
\begin{enumerate}
\item $\gA$ is a \lsd \ssi it is a \qiriz.
\item $\gA$ is \sdz \ssi it is integral.
\end{enumerate}
In particular, a \coh \lsd is a \qiriz.
\end{fact}

Note that item \emph{2} is obvious in \clamaz, where the hypothesis \gui{the annihilator of every \elt is \tfz} is superfluous.

\section{\Tf flat \idsz} 
\label{secIplatTf}

We now study the flatness of \itfsz.
In \clamaz, the following proposition is an \imd corollary of Proposition~\ref{prop.mtf loc plat1}. In \comaz, it is \ncr to provide a new \demz, which gives \algqs information of a different nature from that given in the \dem of Proposition~\ref{prop.mtf loc plat1}. Indeed, we no longer make the same hypotheses regarding the discrete \crc of things.
\begin{proposition} 
\label{prop.itf plat local}{\em (\Tf flat \ids over a \aloz)}\\
Let $\gA$ be a \aloz, $x_1$, \ldots, $x_n\in\gA$ and $\fa=\gen{x_1,\ldots,x_n}$.
\begin{enumerate}
\item If $\fa$ is principal, it is generated by one of the $x_j$'s. \emph{(Bézout is always trivial over a local \riz.)}
\item  If $\fa$ is flat, it is principal, generated by one of the~$x_j$'s.
\item  Suppose that $\gA$ is discrete, or that we have a test to answer the question \gui{is $x$ \ndzz?} 
Then, a \itf is flat \ssi it is free of rank $0$ or~$1$.
\end{enumerate}
\end{proposition}
\begin{proof}
\emph{1.}  We have $\fa=\gen{x_1,\ldots,x_n}=\gen{z}$, $z=a_1x_1+\cdots+a_nx_n$, 
$zb_j=x_j$,
so $z(1-\sum_ja_jb_j)=0$. 
If $1-\sum_ja_jb_j$ is \ivz, $\fa=0=\gen{x_1}$.
If~$a_jb_j$ is \iv $\fa=\gen{x_j}$.
\\
\emph{2.}  Consider the syzygy $x_2x_1+(-x_1)x_2=0$. Let 
\halfsmashbot{$G=\cmatrix{ a_1 & \ldots & a_n \cr b_1 & \ldots & b_n}$} be a matrix such that $G\cmatrix{ x_1 \cr\vdots \cr x_n}=\cmatrix{ x_1 \cr  x_2}$ 
and 
$[\,x_2\;-x_1\,]\,G= [\, 0\,\,\, 0\,]$.  \\
If $a_1$ is \ivz, the \egt  $a_1x_2=b_1x_1$ shows that $\fa=\gen{x_1,x_3,\ldots,x_n }$.
\\ 
\hbox{If $1-a_1$} is \ivz, the \egt $a_1x_1+\cdots+a_nx_n=x_1$ shows that %
we have $\fa=\gen{x_2,x_3,\ldots,x_n }$. \\
We finish by \recu on~$n$.
\\
\emph{3.} Results from \emph{2} and from Lemma~\ref{lemIdps}, item~\emph{3}.
\end{proof}

Recall that a \itf $\fa$ of a \ri $\gA$ is said to be {\em \lopz} if there exist \moco $S_1$, \ldots, $S_n$ of $\gA$ such that each $\fa_{S_j}$ is principal in $\gA_{S_j}$.
The proposition that follows shows that \emph{every \tf flat \id is \lopz}.
Its \dem follows directly from that for the local case.

\begin{proposition} 
\label{prop.itfplat}  
{\em (\Tf flat \ids over an arbitrary \riz)}\\
Every \tf flat \id is \lopz. More \prmtz,
 if
 $\fa=\gen{x_1,\ldots,x_n}\subseteq\gA$, \propeq
\begin{enumerate}
\item  The \id  $\fa$ is a \mplz.
\item  After \lon at  suitable \mocoz, the \id $\fa$  is flat and principal.
\item  After \lon at suitable \ecoz, the \id  $\fa$  is flat and principal, generated by one of the~$x_i$'s.
\end{enumerate}
\end{proposition}
\begin{proof}
We obviously have \emph{3} $\Rightarrow$ \emph{2}. We have \emph{2} $\Rightarrow$ \emph{1} by the \plgref{plcc.plat}. To show \emph{1} $\Rightarrow$ \emph{3} we reuse the \dem of item \emph{2} of Proposition~\ref{prop.itf plat local}.
Consider the syzygy $x_2x_1+(-x_1)x_2=0$. Let \halfsmashbot{$G=\cmatrix{ a_1 & \ldots & a_n \cr b_1 & \ldots & b_n}$} be a matrix such that $G\cmatrix{ x_1 \cr\vdots \cr x_n}=\cmatrix{ x_1 \cr\vdots \cr x_n}$  
and $[\,x_2\;-x_1\,]\,G= [\, 0\,\,\, 0\,]$. With the localized \ri 
$\gA[1/a_1]$  the \egt $a_1x_2=b_1x_1$ shows that $\fa=_{\gA[1/a_1]}\gen{x_1,x_3,\ldots,x_n }$. With the localized \ri 
$\gA[1/(1-a_1)]$  the \egt $a_1x_1+\cdots+a_nx_n=x_1$ shows that $\fa =_{\gA[1/(1-a_1)]} \gen{x_2,x_3,\ldots,x_n }$. 
We finish by \recu on~$n$.
\end{proof}

\vspace{-.7em}
\pagebreak

\subsec{\Anars and \adpsz}

The following \dfn of \adpsz, based on flatness, is due to Hermida and S\'anchez-Giralda
\cite{HS}.

\begin{definition} 
\label{defArit} \emph{(\Anarsz)} 
A \ri $\gA$ is said to be \emph{arithmetic} if every \itf is \lopz.%
\index{ring!arithmetic ---}%
\end{definition}

\begin{propdef} 
\label{prop.itfplat 2}\emph{(\adpsz)}%
\index{Pru@Pr\"ufer!\ri}\index{ring!Pru@Pr\"ufer ---}\\
\Propeq
\begin{enumerate}
\item [1a.] Every \itf of $\gA$ is flat. 
\item [1b.] Every \id of $\gA$ is flat. 
\item [1c.] For all \itfs $\fa$ and $\fb$ of $\gA$, the canonical \ali \hbox{$\fa\te\fb\to\fa\fb$} is an \isoz. 
\item [2a.] The \ri $\gA$ is \losd and \ariz. 
\item [2b.] The \ri $\gA$ is reduced and \ariz. 
\end{enumerate}
A \ri satisfying these \prts is called a \emph{Pr\"{u}fer ring}.

\end{propdef}

\begin{proof}
The \eqvc between \emph{1a} and \emph{1c} is given by \thref{thplatTens} (item \emph{\ref{i3thplatTens}}). The \eqvc of \emph{1a} and \emph{1b} is \imdez. We already know that \hbox{\emph{1a} $\Rightarrow$ \emph{2a}}, and the implication \emph{2a} $\Rightarrow$ 
\emph{2b}  is clear. 

\emph{2b} $\Rightarrow$ \emph{2a.}
Let $x$, $y$ be such that $xy=0$. There exist $s$, $t$ with $s+t = 1$, $sx \in \gen {y}$ and~$ty \in \gen{x}$. Therefore $sx^2 = 0$ and $ty^2 = 0$ then ($\gA$ is reduced) $sx = ty = 0$.

\emph{2a} $\Rightarrow$ \emph{1a}. 
After suitable \lonsz, the \id becomes principal, and therefore flat, since the \ri is a \lsdz. We finish by the \plgref{plcc.plat} for flat modules. 
\end{proof}

\subsec{Local-global principle}

The different notions previously introduced are local in the sense of the following concrete local-global principle. The proofs are based on the basic local-global principle and are left to the reader.\iplg

\begin{plcc} 
\label{plcc.arith}  \emph{(\Anarsz)}\\
Let $S_1$, $\ldots$, $S_n$  be \moco of a \ri $\gA$ and $\fa$ be an \id of~$\gA$.   
We have the following \eqvcsz.
\begin{enumerate}
\item  The \id $\fa$ is \lop \ssi each of the $\fa_{S_i}\!$'s is \lopz.
\item  The \ri $\gA$ is a \lsd \ssi each of the $\gA_{S_i}\!$'s is a \lsdz.
\item  The \ri $\gA$ is \ari \ssi each of the $\gA_{S_i}\!$'s is \ariz.
\item  The \ri $\gA$ is a \adp \ssi each of the~$\gA_{S_i}\!$'s is a \adpz.
\end{enumerate}
\end{plcc}

\vspace{-.4em}
\pagebreak

\subsec{\Lgb machinery}

An ordered set $(E,\leq)$ is said to be \ixc{totally ordered}{set} if for all $x$, $y$ we have $x\leq y$ or $y\leq x$. 
A priori we do not assume it to be discrete and we therefore do not have a test for  strict in\egtz.

 For \alosz, Proposition~\ref{prop.itf plat local} gives the following result.

\begin{lemma}\label{lemme:anarloc}  \emph{(Local \anarsz)}
\begin{enumerate}
\item  A \ri $\gA$ is local and \ari \ssi for 
all~$a$,~$b\in\gA$,
we have $a \in b\gA$ or~$b \in a\gA$. \Eqvtz, every \itf is principal and the set of \itfs is totally ordered with respect to the inclusion.
\item Let $\gA$ be a local \anarz. For two arbitrary \idsz~$\fa$ and~$\fb$, if~$\fa$ is not contained in $\fb$, then $\fb$ is contained in $\fa$.
Therefore in \clamaz, the \gui{set} of all the \ids is totally ordered with respect to the inclusion.
\end{enumerate}
\end{lemma}

Thus, \ari \alos are the same thing as local Bézout \risz.
They have already been studied in Section~\ref{secBézout} (\paref{secpfval}).

The ease with which we prove \prts for \anars is mostly due to the following \lgbe machinery.

\mni{\bf \Lgb machinery of \anars}
\label{MetgenAnar}\imla
\\
\emph{When we have to prove a \prt regarding an \anar and that a finite family of \elts $(a_i)$ of the \ri intervenes in the computation, we start by proving the result in the local case.
We can therefore suppose that the \ids $\gen{a_i}$ are totally ordered by inclusion. 
In this case the \dem is in \gnl very simple.
Moreover, since the \ri is \ariz, we know that we can return to the previous situation after \lon at a finite number of \ecoz. We can therefore conclude if the \prt to be proven obeys a \plgcz.}

\medskip Here is an application of this machinery.

\begin{proposition}\label{propIddsAnar} \emph{(Determinantal \ids over an \anarz)}
\\
Let $\gA$ be a coherent \anarz, $M$ be a matrix $\in\Ae{n\times m}$ and $\fd_k=\cD_{\gA,k}(M)$ its \idds ($k\in\lrbp$ with $p=\inf(m,n)$). There exist \itfs $\fa_1,\ldots,\fa_p$
such that 

\snic{
\fd_1=\fa_1  ,\; \fd_2=\fd_1\fa_1\fa_2,\; \fd_3=\fd_2\fa_1\fa_2\fa_3,\; \ldots
}

\end{proposition}
%
\begin{proof}
Let $\mathfrak{b}_k=(\mathfrak{d}_k:\mathfrak{d}_{k-1})$ for all $k$, so $\mathfrak{b}_1=\mathfrak{d}_1$. We have $\mathfrak{b}_k\mathfrak{d}_{k-1}=\mathfrak{d}_k$ because the ring is arithmetic and coherent. 
Let  $\mathfrak{c}_k=\mathfrak{b}_1\cap\cdots\cap\mathfrak{b}_k$ for $k\geq 1$. This is a nonincreasing sequence of finitely generated ideals. Let  $\mathfrak{a}_1,\ldots,\mathfrak{a}_p$ be finitely generated ideals
satisfying $\mathfrak{a}_1=\mathfrak{d}_1$ and $\mathfrak{a}_k\mathfrak{c}_{k-1}=\mathfrak{c}_k$ for $k\geq 2$. It is sufficient to prove the equalities \fbox{$\mathfrak{c}_k\mathfrak{d}_{k-1}=\mathfrak{d}_k$}. This is clear for $k= 1$.
\\ 
If $\mathbf{A}$ is a local arithmetic ring, the matrix admits a reduced Smith form (Proposition~\ref{propPfVal}). Let $p=\inf(m,n)$
\\
The algorithm that produces the reduced Smith form in the local case and the previous local-global machinery of arithmetic rings provide us with a system of comaximal elements $(s_1,\ldots,s_r)$ such that, over each ring $\mathbf{A}[1/s_i]$, the matrix $M$ admits a reduced Smith form with the diagonal sub-matrix  $\mathrm{Diag}(c_1,c_2,\dots,c_p)$ and $c_1\mid c_2\mid \dots\mid c_{p}$. Moreover, for $k\geq 1$, $\mathfrak{d}_k=\gen{c_1\cdots c_k}$.
\\
It is sufficient to prove  $\mathfrak{c}_k\mathfrak{d}_{k-1}=\mathfrak{d}_k$ after \lon at these \ecoz.
Since $\mathfrak{b}_k\mathfrak{d}_{k-1}=\mathfrak{d}_k$, we get $\mathfrak{c}_k\mathfrak{d}_{k-1}\subseteq \mathfrak{d}_k$.
\\
In order to prove the other inclusion
 let us show that for all $k\geq 1$ we have $c_k\in\mathfrak{c}_k$ (this implies $\mathfrak{c}_k\mathfrak{d}_{k-1}\supseteq \mathfrak{d}_k$).
We have  $c_k\mathfrak{d}_{k-1}=\mathfrak{d}_k$, thus
$c_k\in\mathfrak{b}_k$. Moreover $c_k$ is multiple of $c_i\in\mathfrak{b}_i$ for  $i\leq {k-1}$, so $c_k\in\mathfrak{b}_1\cap\cdots\cap\mathfrak{b}_{k-1}$.
This finishes the proof.
\end{proof}

We will return to \anars and \adps in greater length in Chapter~\ref{ChapAdpc}.

\penalty-2500
\section{Flat algebras}
\label{secAlplates}

Intuitively speaking, an \Alg $\gB$ is flat when the homogeneous \slis over $\gA$ have \gui{no more} solutions in $\gB$ than in~$\gA$, and it is faithfully flat if this assertion is \egmt true for nonhomogeneous \slisz. More precisely, we adopt the following \dfnsz. 

\begin{definition}\label{defiAlgPlate} Let $\rho:\gA\to\gB$ be an \Algz.
\begin{enumerate}
\item $\gB$ is said to be \emph{flat (over $\gA$)} when every $\gB$-\lin \rde between \elts of $\gA$ is a $\gB$-\lin combination of $\gA$-\lin \rdes between these same \eltsz.
In other words, for every \lin form $\psi:\gA^n\to\gA$, we require that $\Ker \rho\ist(\psi) =\gen{\rho (\Ker \psi )}_\gB$.
\\
We will also say that \emph{the \ri \homo $\rho$ is flat}.%
\index{flat!algebra}\index{algebra!flat ---}%
\index{flat!ring homomorphism}
\item A \emph{flat} \Alg $\gB$ is said to be \emph{\fptez} if for every \lin form $\psi:\gA^n\to\gA$ and all $a\in\gA$, when the equation $\psi(X)=a$ 
admits a solution in $\gB$ (\cad $\exists X\in\gB^n,\,\big(\rho\ist(\psi)\big)(X)=\rho(a)$), 
then it admits a solution in $\gA$. 
\\
We will also say that \emph{the \ri \homo $\rho$ is \fptz}.
\end{enumerate}
\index{faithfully flat!algebra}
\index{faithfully flat!ring homomorphism}
\index{algebra!faithfully flat ---}
\end{definition}

For a \fpte \Algz, when considering the case where $n=1$ 
and $\psi=0$, we see that $\rho(a)=0$ implies $a=0$. Thus, $\rho$ is an injective \homoz.
We therefore say that $\gB$ is a \emph{\fpte extension} of~$\gA$.
We can then identify $\gA$ with a sub\ri of $\gB$ and the condition on the nonhomogeneous \lin equation is easier to formulate:
it is exactly the same equation that we seek to solve in $\gA$ or~$\gB$.

\begin{fact}\label{factplatplat} ~ \\
An \Alg $\gB$ is flat \ssi $\gB$ is a flat \Amoz. 
\end{fact}
\begin{proof}
Translation exercise left to the reader.
\end{proof}

{\bf Fundamental examples. }
The following lemma provides some examples.

\vspace{-.1em}
\pagebreak

\begin{lemma}\label{lemExAlgPlFiPl}~
\begin{enumerate}
\item 
A \lon morphism
$\gA\to S^{-1}\gA$ gives a flat \Algz.
\item 
If $S_1$, \ldots, $S_n$ are \moco of $\gA$ and if $\gB=\prod_i\gA_{S_i}$,
the canonical \gui{diagonal} \homo $\rho:\gA\to\gB$ gives a \fpte \algz.
\item If $\gk$  is \zedrz, every \klg $\gL$ is flat.
\end{enumerate}
\end{lemma}
\begin{proof}
\emph{1.} See Fact~\ref{fact.transporteur} or Facts \ref{factplatplat} and~\ref{fact.plat}.
\\
\emph{2.} This results from the basic \plg (we could even say that \und{it is} the basic \plgz).\iplg
\\
\emph{3.} Results from \ref{factplatplat} and from the fact that every \Kmo is flat (\thref{propEvcPlat}). 
\end{proof}

\rems Regarding item \emph{3} of the previous lemma.
 
1) 
 It seems difficult to replace $\gk$ in the hypothesis with a (Heyting) field that we do not assume to be \zedz.
 
2)
See \thref{thSurZedFidPlat} for the \fpte question. 
\eoe


\medskip 
In the following proposition, an analog of Propositions~\ref{propCoh1} (for \corisz) and \ref{propPlat1} (for \mplsz), we pass from an \eqn to a \sys of \eqnsz.
To lighten the text, we act as though we have an inclusion $\gA\subseteq\gB$ (even if $\gB$ is only assumed to be flat), in other words we do not specify that when we pass into $\gB$, everything must be transformed by means of the \homo $\rho:\gA\to\gB$.

\begin{proposition}\label{lemAlgPlate} 
Let $M\in\Ae{n\times m}$,  $C\in\Ae{n\times 1}$ and $\gB$ be a flat \Algz. 
\begin{enumerate}
\item Every solution in $\gB$ of the \hmg \sli $MX=0$ is a $\gB$-\lin combination of solutions in~$\gA$.
\item If in addition $\gB$ is \fptez, and if the \sys $MX=C$ admits a solution in $\gB$, it admits a solution in~$\gA$.
\end{enumerate}
\end{proposition}
\begin{proof} The \dfns of flat and \fptes \Algs concern the \slis with a single equation.
To solve a \gnl \sli we apply the usual technique: we start by solving the first equation, then substitute the \gnl solution of the first equation into the second, and so forth.
\end{proof}
%

\begin{proposition}\label{propIdsAlgPlate}~ \\
Let $\gA\vers{\rho}\gB$ be a flat \Alg and $\fa,\,\fb$ be two \ids of~$\gA$.
\begin{enumerate}
\item The natural \Bli $\rho\ist(\fa)\to \rho(\fa)\gB$ is an \isoz.
\end{enumerate}
In the remainder we identify $\rho\ist(\fc)$ with the \id $\rho(\fc)\gB$ for every \id $\fc$ of~$\gA$.
\begin{enumerate} \setcounter{enumi}{1}
\item We have $\rho\ist(\fa\cap \fb)=\rho\ist(\fa)\cap\rho\ist(\fb)$. 
\item If in addition $\fa$ is \tfz, we have 
$\rho\ist(\fb: \fa)=\big(\rho\ist(\fb):\rho\ist(\fa)\big)$. 
\end{enumerate}
\end{proposition}
\begin{proof}
The first two items result from analogous facts regarding \mpls (\thref{thplatTens} item \emph{\ref{i4thplatTens}} and Corollary~\ref{cor3PlatTens}).
\\
\emph{3.} If $\fa=\gen{\an}$, then $\fb:\fa=\bigcap_i(\fb:a_i)$, 
therefore given item~\emph{2} we are reduced to the case of a \idp $\gen{a}$. 
We then consider the exact sequence  
$$\preskip-.5em \postskip.3em 
0\to\,\fb: a\,\llra\gA\vvers a\gA\sur\fb,
$$
we tensor with $\gB$
and we obtain the exact sequence (use the flatness and Fact~\ref{factSuitExTens})
$$\preskip-.5em \postskip.3em 
0\to\rho\ist(\fb: a)\llra\gB\vvvers {\rho(a)}\gB\sur{\rho\ist(\fb)}, 
$$
which gives the desired result. 
\end{proof}
%

\begin{theorem} 
\label{thExtPlat} 
Let $\rho:\gA\rightarrow \gB$ be an \algz.  
\Propeq
\begin{enumerate}
\item \label{i1thExtPlat} $\gB$ is a flat \Algz.
\item \label{i2thExtPlat} $\gB$ is a flat \Amoz.
\item \label{i3thExtPlat} For every flat \Amo $M$, the \Amo  $\rho\ist (M)$  is flat.
\item \label{i4thExtPlat} For every \itf $\fa$ of $\gA$, the canonical \Ali
$$\preskip.2em \postskip.0em 
\gB\te_\gA\fa\simeq\rho\ist(\fa)\to\fa\gB 
$$
is an \isoz.
\item \label{i5thExtPlat} For all \Amos $N\subseteq M$, the \Bli $\rho\ist(N)\to\rho\ist(M)$  is injective.
\item  \label{i6thExtPlat} For every \Ali $\psi:M\to P$, the natural \Bli 
$$\preskip.2em \postskip.0em 
\rho\ist\big(\Ker(\psi)\big)\lora\Ker\big(\rho\ist (\psi)\big) 
$$
is an \isoz.
\item \label{i7thExtPlat} For every exact sequence of \Amos $M\vers{f} N\vers{g} P$ 
the sequence  
$$\preskip.1em \postskip.2em 
\rho\ist(M)\vvvers{\rho\ist(f)} \rho\ist(N)\vvvers{\rho\ist (g)} \rho\ist (P) 
$$
is an exact sequence of \Bmosz.
\end{enumerate}
\end{theorem}
Item \emph{\ref{i5thExtPlat}} allows us to identify  $\rho\ist(P)$ with a \Bsmo of $\rho\ist(Q)$ each time that we have two \Amos $P\subseteq Q$ and that $\gB$ is flat over~$\gA$.
\begin{proof} 
The reader will verify that the \eqvcs are clear by what we already know (Fact~\ref{factplatplat}, \thref{thplatTens}, Corollary~\ref{corPlatTens}). We  note that Proposition~\ref{lemAlgPlate} gives item \emph{\ref{i6thExtPlat}} in the case of a \ali between free modules of finite rank.
\end{proof}

The following proposition generalizes Propositions~\ref{fact.hom egaux} and \ref{fact.homom loc pf}.
 
\begin{proposition}\label{propPlateHom}
Let $\rho:\gA\to\gB$ be a flat \Alg and $M$, $N$ 
be \Amosz.
If $M$ is \tf (resp.\ \pfz), the natural \Bli 
$$\preskip-.20em \postskip.4em 
\rho\ist\big(\Lin_\gA(M,N)\big)\to \Lin_\gB\big(\rho\ist(M),\rho\ist(N)\big) 
$$
is injective (resp.\ is an \isoz). 
\end{proposition}
%
\begin{proof}
Consider an exact sequence 
$$\preskip.2em \postskip.3em 
\qquad K\lora \gA^k\lora M\to 0,\qquad(*) 
$$
corresponding to the fact that $M$ is \tf (if $M$ is \pf the module $K$ is also free of finite rank).  \\
Let $M_1=\rho\ist(M)$, $N_1=\rho\ist(N)$ and $K_1=\rho\ist(K)$.
First we have the exact sequence
$$\preskip-.4em \postskip.4em 
\qquad K_1\lora \gB^k\lora M_1\to 0.\qquad(**) 
$$
Next we obtain the exact sequences below. The first comes from $(*)$, the last one comes from $(**)$ and the second results from the first by \eds since $\gB$ is flat over~$\gA$. 
\[\preskip.6em \postskip.7em \arraycolsep2pt
\begin{array}{cccccccccc} 
0  &\to& \Lin_\gA(M,N) &\to& \Lin_\gA(\gA^k,N)\simeq N^k &\to&\Lin_\gA(K,N)
\\[3mm] 
0  &\to&\rho\ist\big(\Lin_\gA(M,N)\big)&\to&\rho\ist(\Lin_\gA(\gA^k,N)\simeq N_1^k   &\to&   \rho\ist\big(\Lin_\gA(K,N)\big)\\[1mm] 
 &&\dar{} &&\dar{} &&\dar{}  \\[1mm] 
0  &\to& \Lin_\gB(M_1,N_1)&\to&\Lin_\gB(\gB^k,N_1)\simeq N_1^k&\to&\Lin_\gB(K_1,N_1) 
 \end{array}
\]
 In addition, we have natural \gui{vertical} \Blis from the second to the third exact sequence, and the diagrams commute.
The second vertical arrow is an \iso (the identity of $N_1^k$ after the canonical identifications).
This implies that the first vertical arrow (the \Bli that we are interested in) is injective.
\\
If $M$ is \pf and if $K\simeq \gA^{\ell}$, the two \Bmos on the right are \isoc to $N_1^{\ell}$ and the corresponding vertical arrow is an \isoz. This implies that the first vertical arrow is an \isoz.
\end{proof}

Retrospectively the given \dem for Proposition~\ref{fact.homom loc pf} seems quite complicated.
The new \dem given here in a more \gnl framework is conceptually simpler.

\section{Faithfully flat \algsz}
\label{secAlfidplates}

We have already said that if $\gA\vers\rho\gB$ is a \fpte \algz, $\rho$ is injective. 
\perso{autres \prts analogues? Radical of Jacobson?}
It is \egmt clear that $\rho$  \emph{reflects the units}, \cad 
\index{reflects the units!\homo that ---}

\snic{\rho(a)\in\gB\eti\;\Longrightarrow\;a\in\Ati.}

\medskip 
 We now present a few \caras \prtsz. In what follows we will take note of the \eqvc of items \emph{1}, \emph{2a},  \emph{3a} and \emph{4.} 

\begin{theorem} 
\label{thExtFidPlat} \emph{(\Carns of  \fpt \algsz)}\\
Let $\rho:\gA\rightarrow \gB$ be a flat \algz.  
\Propeq
\begin{enumerate}
\item [1.] 
The \alg $\gB$ is \fptez.
\item [2a.] 
The \homo $\rho$ is injective, and when identifying $\gA$ with a sub\ri of $\gB$, for every \itf $\fa$ of $\gA$ we have

\snic{ \fa\gB\,\cap\,\gA\,=\,\fa.}

\item [2b.] 
Similarly with an arbitrary \id of $\gA$.
\item [3a.] 
For every \itf $\fa$ of $\gA$ we have the implication

\snic{ 1_\gB\in \rho\ist(\fa)\;\Longrightarrow\;1_\gA\in\fa.} 
\item [3b.]  
For every \itf $\fa$ of $\gA$, if  
$\rho\ist(\gA\sur\fa)=0$, then $\gA\sur\fa=0$.
\item [3c.] 
For all \Amos $N\subseteq M$, if $\rho\ist(N)=
\rho\ist(M)$, then~$N=M$.
\item [3d.] 
For every \Amo $M$, if  $\rho\ist(M)=0$, then~$M=0$.
\item [3e.] 
For every \Amo $M$ the natural \Ali $M\to\rho\ist(M)$ is injective.
%
%
\item [4.] 
The \eds from $\gA$ to $\gB$ \emph{reflects the exact sequences}.
\\
In other words, given an arbitrary sequence of \Amos 

\snic{N\vers{f} M\vers{g} P,}
 
 it is exact if the sequence of \Bmos  

\snic{\rho\ist (N)\vvvers{\rho\ist (f)} \rho\ist (M)\vvvers{\rho\ist (g)} \rho\ist (P)}
 
is exact.
\end{enumerate}
\end{theorem}
%
\begin{proof}
Item \emph{1} implies that $\rho$ is injective.
Once this has been shown,  \emph{2a} is a simple reformulation of \emph{1},
and it is easy to show that \emph{2a} is equivalent to \emph{2b}.

\emph{3a} $\Rightarrow$ \emph{1}.
We start by noticing that the implication is still valid if we replace the \itf $\fa$ by an arbitrary \id $\fc$. Indeed, \hbox{if $1\in\rho\ist(\fc)$}
we will \egmt have $1\in\rho\ist(\fc')$ for a \itf $\fc'$ contained in $\fc$.
\\
Now let $\fa=\gen{a_1,\ldots,a_n}$ and $c\in\gA$. 
The equation $\sum_ia_ix_i=c$ admits a solution \ssi $c\in\fa$, \cad $1\in(\fa:c)_\gA$. 
Since $\gB$ is flat, we have 
$\big(\rho\ist(\fa):\rho(c)\big)_\gB=\rho\ist(\fa:c)$ (Proposition~\ref{propIdsAlgPlate}).
If $\sum_i\rho(a_i)y_i=\rho(c)$ admits a solution in $\gB$,
then $1\in\big(\rho\ist(\fa):\rho(c)\big)_\gB$, so the hypothesis \emph{3a} implies that $1\in(\fa:c)$, \cad $\sum_ia_ix_i=c$ admits a solution in~$\gA$.

 The implications  \emph{3e} $\Rightarrow$ \emph{3d} $\Rightarrow$ \emph{3b} are trivial.

\emph{3d} $\Rightarrow$ \emph{3c}. Consider the module $M/N$.
The module $\rho\ist(N)$ is identified with a submodule of $\rho\ist(M)$ and $\rho\ist(M\sur N)$ is identified with $\rho\ist(M)\sur{\rho\ist(N)}$. The result follows.

\emph{3c} $\Rightarrow$ \emph{3d}. We take $N=0$.

\emph{3a} $\Leftrightarrow$ \emph{3b}. Same reasoning. 

\emph{1} $\Rightarrow$ \emph{3e}.
 We identify $\gA$ with a sub\ri of $\gB$. \\
 Let $x\in M$ such that $1\te x=0$ in $\rho\ist(M)$. 
 Since $\gB$ is a flat \Amoz, this syzygy is explained in the \Amo $\gB$: there exist $u_1$, \dots, $u_n\in\gB$ and $a_1$, \dots, $a_n\in\gA$ such that $\sum_ia_iu_i=1$ and $a_ix=0$ for $i\in\lrbn$. 
The equation in the $y_i$'s, $\sum_ia_iy_i=1$, admits a solution in~$\gB$, so it admits one in $\gA$. Hence~$x=0$.

 \emph{4} $\Rightarrow$ \emph{3d.}
We make $N=P=0$ in the sequence $N\to M\to P$.
It is exact after \eds to $\gB$, so it is exact.

\emph{1} $\Rightarrow$ \emph{4.}
Suppose that the sequence of \Bmos is exact. We must show that the sequence of \Amos is exact. First of all $g\circ f=0$, because the \Bli $P\to\rho\ist(P)$ is injective, and the diagrams commute. Next, since $\gB$ is flat, we can identify $\rho\ist(\Ker g)$ with $\Ker\rho\ist(g)$ and~$\rho\ist(\Im f)$ with $\Im\rho\ist(f)$. We are back in item~\emph{3c.}
\end{proof}
Given \thref{propEvcPlat}, we obtain as a consequence of the \carnz~\emph{2a} the following \thoz.

\begin{theorem}\label{thSurZedFidPlat}
Every extension of a \cdi or of a \zedr \ri is \fptez. 
\end{theorem}
%
\begin{proof}
We have $\gk\subseteq\gA$ with $\gk$ \zedrz. We know that the extension is flat by \thref{propEvcPlat}. We must show that if $\fa$ is a \itf of $\gk$, then $\fa\gA\cap\gk=\fa$. However, $\fa=\gen{e}$ for an \idmz~$e$; the membership of an \elt $x$ in an \id $\gen {e}$ ($e$ \idmz) being \caree by the \egtz~\hbox{$x = xe$}, it is independent of the \riz. In other words, for an \idm  $e$ of a \ri $\gB\subseteq \gB'$, we always have $e\gB' \cap \gB = e\gB$.
\end{proof}

As a special case of the \carn \emph{3a} we obtain the following corollary. 

\begin{corollary}\label{corpropFidPlatLoc}
Let $\rho$ be a flat \homo between \alosz. It is \fpt \ssi it \emph{reflects the units},
\cad $\rho^{-1}(\gB\eti)=\Ati$. 
\end{corollary}

A \homo between \alos that reflects the units is called a \emph{local \homoz}.
\index{local!homomorphism}
\index{homomorphism!local ---}

The \dems of the two following facts result from simple considerations about the preservation and about the \gui{reflection} of the exact sequences.
The details are left to the reader.

\begin{fact}\label{factAlgPlate} \emph{(Transitivity)}
Let $\gB$ be an \Alg and $\gC$ be a \Blgz.
\begin{enumerate}
\item If $\gB$ is flat over $\gA$ and $\gC$ flat over $\gB$, then
$\gC$ is flat over~$\gA$.
\item If $\gB$ is \fpte over $\gA$ and $\gC$ \fpte over $\gB$, then
$\gC$ is \fpte over~$\gA$.
\item If $\gC$ is \fpte over $\gB$ and flat over $\gA$, then
$\gB$ is flat over~$\gA$. 
\item If $\gC$ is \fpte over $\gB$ and  over $\gA$, then
$\gB$ is \fpte over~$\gA$. 
\end{enumerate}
\end{fact}

\begin{fact}\label{factAlgPlate2} \emph{(Changing the base \riz)}\\
Let $\gB$  and $\gC$ be two \Algsz, and $\gD=\gB\otimes_\gA\gC$.
\begin{enumerate}
\item If $\gC$ is flat over $\gA$, $\gD$ is flat over~$\gB$.
\item If $\gC$ is \fpte over $\gA$, $\gD$ is \fpte over~$\gB$.
\end{enumerate}
\end{fact} 

\vspace{-.6em}
\pagebreak

\begin{plcc}\label{plcc.algfptes}  \emph{(Localization at the source, flat \algsz)}\\
Let $\rho:\gA\to\gB$ be an \alg and $S_1$, $\ldots$, $S_r$  be \moco of~$\gA$.
\begin{enumerate}
\item The \alg $\gB$ is flat over $\gA$ \ssi for each~$i$, 
$\gB_{S_i}$ is flat over~$\gA_{S_i}$.
\item The \alg $\gB$ is \fpte over $\gA$ \ssi for each~$i$, 
the \alg $\gB_{S_i}$ is \fpte over~$\gA_{S_i}$. 
\end{enumerate}
\end{plcc}
\begin{proof}
We introduce the \fpte \Alg $\gC=\prod_i\gA_{S_i}$ that gives by \eds the \fpte
\Blg $\gD=\prod_i\gB_{S_i}$. It remains to apply Facts~\ref{factAlgPlate} and~\ref{factAlgPlate2}.
\end{proof}
%

The following \tho \gnss the \plgcs that assert the local \crc (in the \cof sense)
  of certain \prts of finiteness for modules.

\begin{theorem}\label{propFidPlatTf}
Let $\gA\vers \rho\gB$ be a \fpte \Algz.
\\
Let $M$ be an~\Amo \hbox{and $M_1=\rho\ist(M)\simeq\gB\otimes_\gA M$}.
\begin{enumerate}
\item The \Amo $M$ is flat \ssi the \Bmo $M_1$ is flat.
\item The \Amo $M$ is \tf \ssi the \Bmo  $M_1$ is \tfz.
\item If the \Bmo $M_1$ is \cohz, the \Amo $M$ is \cohz.
\item The \Amo $M$ is \pf \ssi the \Bmoz~\hbox{$M_1$} is \pfz.
\item The \Amo $M$ is \ptf \ssi the \Bmoz~\hbox{$M_1$} is \ptfz.
\item If the \Bmo $M_1$ is \noez, the \Amo $M$ is \noez.
\end{enumerate}
\end{theorem}
\begin{proof} In items \emph{1}, \emph{2}, \emph{4}, \emph{5}, we already know that any \eds preserves the concerned \prtz. Therefore all that remains to be proven are the converses.
\\
\emph{1.} Consider an exact sequence $N\vers f Q \vers g P$ of \Amosz.
We want to show that it is exact after tensorization by $M$.
We know that it is exact after tensorization by $\gB\te M$. However, $\gB\te \bullet$
reflects the exact sequences.

\emph{2.} Consider some \elts $y_i\in\rho\ist(M)$ ($i\in\lrbn$) that generate this module. These \elts are constructed as $\gB$-\lin combinations of a finite family of \elts  $1\te x_j$  ($x_j\in M$, $j\in\lrbm$).
This implies that the \Ali $\varphi:\Ae m\to M$ which sends the canonical basis to~$(x_j)_{j\in\lrbm}$ is surjective after tensorization by $\gB$. However, $\gB$ is \fptez,
therefore $\varphi$ is surjective.

\emph{3.} Let $N=\gA x_1+\cdots+\gA x_n$ be a \smtf of $M$.
Consider the corresponding surjective \Ali $\Ae n\to N$, let $K$ be its kernel.
The exact sequence $0\to K\to \Ae n\to N\to 0$ gives by \eds an exact sequence (because $\gB$ is flat).
Since $\rho\ist(M)$ is \cohz, $\rho\ist(K)$ is \tfz. 
It remains to apply item~\emph{2.}

\emph{4.} Same reasoning as for item~\emph{3.}

\emph{5.} A module is \ptf \ssi it is flat and \pfz.

\emph{6.} 
Consider an ascending sequence $(N_k)_{k\in\NN}$ of \smtfs of~$M$ and  extend the scalars to $\gB$. Two consecutive terms $\rho\ist(N_\ell)$ and~$\rho\ist(N_{\ell+1})$ are equal. Since $\gB$ is \fptez, we also have the \egts $N_\ell=N_{\ell+1}$.
\end{proof}
%

The following \tho  generalizes the \plgcs that asserts the local \crc  (in the \cof sense) of certain \prts of finiteness for \algsz.

\begin{theorem}\label{propFidPlatPrAlg}~
\\
\Deuxcol{.63}{.3}{
Let $\rho:\gA\to\gB$ be a \fpte \Alg and $\varphi:\gA\to\gC$ be an \Algz. 
\\
Let $\gD=\rho\ist(\gC)$ be the \fpte \Clg obtained by \edsz.}
{\vspace{-5mm}
\hspace{5mm}\xymatrix @R=18pt @C=12pt
{
\gA  \ar@{->}[d]_\rho\ar@{->}[rr]^{\varphi} 
   &\ar@{~>}[d]^{\rho\ist}
   &\gC\ar@{->}[d]
   \\                      
\gB  \ar@{->}[rr]_{\rho\ist(\varphi)} && \gD  \\
}
}

\vspace{-1mm}
In order for $\gC$ to have one of the \prts below as an \Alg it is necessary and sufficient that $\gD$ possesses the same \prt as a \Blgz:
\begin{itemize}
\item finite (as a module),
\item \pf as a module, 
\item \stfez, 
\item flat, 
\item \fptez, 
\item \stez,
\item \spbz,
%
%
\item \tf (as an \algz), 
\item \pf (as an \algz). 
%
%
\end{itemize}
\end{theorem} 
\begin{proof} 
The first three \prts are \prts of modules and thus falls within \thref{propFidPlatTf}.

\emph{Flat, \fptes \algsz}.
We apply Facts~\ref{factAlgPlate} and~\ref{factAlgPlate2}.

\emph{\Ste \algsz}. We already have the \eqvc for the \stf \crcz. If $\gB$ is free over $\gA$ we use the fact that the \discri is well-behaved under \edsz, and we conclude by using the fact that a \fpte extension reflects the units.
\\
In the \gnl case we return to the free case by \lon at \ecoz, or we invoke \thref{thAlgStfSpbSte}: 
a \stfe \alg is \spb \ssi it is \stez.

\emph{\Spb \algsz}. Consider at the commutative diagram in Fact \ref{factSpbEds} (beware, the names change). 
The vertical arrow on the right is obtained by \fpte \eds from the one on the left. They are therefore simultaneously surjective.  

\emph{\Tf \algsz}. The fact of being \tf or \pf is preserved by any \edsz. Let us take a look at the converse. 
\\
We identify $\gA$ with a sub\ri of $\gB$ and $\gC$ with a sub\ri of $\gD$.\\
Let $\gA_1=\varphi(\gA)$ and $\gB_1=\rho\ist(\varphi)(\gB)$.
Since $\gD=\gB\otimes _\gA\!\gC$ is \tf over $\gB$, and since every \elt of $\gD$ is expressible as a $\gB$-\lin combination of \elts of $\gC$, we can write $\gD=\gB_1[\xm]$ with some $x_i\in\gC\subseteq\gD$. This gives an exact sequence

\snic{\gB[\Xm]\vvvvers{
\rho\ist(\varphi) ,\, X_i\mapsto x_i
}
\gD\lora 0.
}

 We will show that $\gC=\gA_1[\xm]$. Indeed, the exact sequence above is obtained by \fpte \eds from the sequence
 
\snic{\gA[\Xm]\vvvvers{
\varphi ,\, X_i\mapsto x_i
}
\gC\lora 0.
}

\emph{\Fp \algsz}. 
\\
Let us begin with an \elr \gnl but useful remark about quotient \algs $\kuX\sur\fa$. We can regard $\kuX$ as the free \kmo having as its basis the family of \moms $(X^{\ual})_{\ual\in\NN^m}$. If $f\in\fa$, we then obtain the \egt

\snic{f\cdot\kuX=\sum_\ual (X^\ual f)\cdot\gk.}

Therefore the \id $\fa$ is the \ksmo of $\kuX$ generated by all the $X^\ual f$, where~$\ual$ ranges over $\NN^m$ and $f$ ranges over a \sgr of~$\fa$.  
\\
Let us then return to the \dem by continuing with the same notations as in the previous item.\\ 
Suppose that $\gD=\gB_1[\xm]\simeq\aqo{\gB[\uX]}{\lfs}$.
In the remainder we consider an \eqn $f_j=0$ as a syzygy between the \moms present in $f_j$.
Since the \Bmo $\gD$ is obtained by flat \eds of the \Amo $\gC$, the $\gB$-\lin \rdez~$f_j$  is a $\gB$-\lin combination of $\gA$-\lin \rdes $f_{j,k}$ (between the same \moms viewed in~$\gC$).
Each \egt $f_{j,k}(\ux)=0$ can also be read as an $\gA$-\agq \rde (a relator) between  $x_i$'s $\in\gC$.
Consider then the \Asub of $\gA[\uX]$ generated by all the $X^\ual f_{j,k}$'s.
By \eds from $\gA$ to $\gB$ the sequence of \Amos
$$\preskip.3em \postskip.4em\ndsp 
0\to \sum_{j,k,\ual}(X^\ual  f_{j,k})\cdot\gA\to\AuX\to\gC\to 0\quad(*) 
$$
gives the exact sequence of \Bmos
$$\preskip.4em \postskip.4em\ndsp 
0\to \sum_{j,k,\ual}(X^\ual  f_{j,k})\cdot\gB\to\BuX\to\gD\to 0. 
$$
Indeed, $\sum_{j,k,\ual}(X^\ual  f_{j,k})\cdot\gB=\sum_{j,k}f_{j,k}\cdot \BuX=\sum_{j}f_{j}\cdot \BuX=\fa$. Therefore,
since the extension is \fptez, the sequence $(*)$ is itself exact. Finally, since $\sum_{j,k,\ual}(X^\ual  f_{j,k})\cdot\gA=\sum_{j,k}  f_{j,k}\cdot\AuX$,
$\gC$ is a \pf \Algz.
\end{proof}
%

\Exercices{ 

\begin{exercise}
\label{exoPlatsLecteur}
{\rm  We recommend that the proofs which are not given, or are sketched, or left to the reader, etc, be done. But in particular, we will cover the following cases.
\begin{itemize}\itemsep0pt
\item \label{exooPlatsLecteur0} 
Over a Bézout domain, a module is flat \ssi it is torsion-free.
\item \label{exooothPlat1} 
Prove \thref{thPlat1}.

\item \label{exocorlemIdps} 
Prove Lemma~\ref{corlemIdps}.
\item \label{exothExtPlat} 
Prove Fact~\ref{factplatplat} and \thref{thExtPlat}.
\item Prove Facts~\ref{factAlgPlate} and~\ref{factAlgPlate2}.

\end{itemize}
}
\end{exercise}

\vspace{-1em}
\begin{exercise}
\label{exothPlat1}
{\rm Let  $\pi:N\to M$ be a surjective \aliz.
\\
 1. If $M$ is flat,  for every  \mpf $P$, the natural \ali

\snic{\Lin_\gA(P,\pi):\Lin_\gA(P,N)\to \Lin_\gA(P,M)}

is surjective. \emph{(Particular case of \thref{propPlatQuotientdePlat}.)}
\\
2. Suppose that  $N=\Ae{(I)}$, a free module over a discrete set $I$. If the previous \prt is satisfied, $M$ is flat. 

\comm In \comaz, an arbitrary module $M$ is not \ncrt a quotient of a module  $N=\Ae{(I)}$ as above, but this is true in the case where $M$ is discrete, by taking $I=M$. If we don't need $I$ discrete, we look at Exercise~\ref{propfreeplat}.
\eoe
}
\end{exercise}

\vspace{-1em}
\begin{exercise} 
\label{exoplatFitt} 
{\rm  Let $M$ be a \tf \Amoz. Prove that if $M$ is flat its \idfs are \idmsz. 
}
\end{exercise}

\vspace{-1em}
\begin{exercise} 
\label{exoClamlsdz} 
{\rm  Show using \clama 
         that a \ri is a \lsd \ssi it becomes integral after \lon at any \idepz.
} 
\end{exercise}

\vspace{-1em}
\begin{exercise} 
\label{exoClamAri} 
{\rm  Show using \clama that a \ri is \ari \ssi it becomes a Bézout \ri after \lon at any \idepz.
} 
\end{exercise}

\vspace{-1em}
\begin{exercise}
\label{exoIlops0}
{\rm  The image of a \lop \id under a \ri \homo is a \lop \idz.
Prove that the analogous result for \iv \ids is not always true. 
 }
\end{exercise}

\vspace{-1em}
\begin{exercise} 
\label{exoIlops2} 
{\rm 
If $\fa = \gen{x_1,\dots, x_k}$ is \lopz, then
   $\fa^n = \gen{x_1^n,\dots, x_k^n}$.
   \\
Compute a \mlp for $(x_1^n,\dots, x_k^n)$ from a \mlp for $(x_1,\dots, x_k)$. 
\\
Explicate the membership of $x_1^{n_1}\cdots x_k^{n_k}\in\gen{x_1^n,\dots, x_k^n}$ when $n=n_1+\cdots+n_k$.
} 
\end{exercise}

\vspace{-1em}
\begin{exercise}
\label{exoMLP}
{\rm Given $n$ \elts in an \anar give an \algo that constructs a \mlp for those \elts from \mlps for pairs of \elts only. 
 
}
\end{exercise}

\vspace{-1em}
\begin{exercise}
\label{exoIlops}
{\rm  Consider two \itfs $\fa$ and $\fb$ of a \ri $\gA$,
generated respectively by $m$ and $n$ \eltsz.
Let $f$, $g\in\AX$ of degrees $m-1$ and $n-1$ with $\rc(f)=\fa$ and $\rc(g)=\fb$.
 
\emph{1.}
Show that if $\fa$ is \lopz, we have $\fa\fb=\rc(fg)$ such that $\fa\fb$ is generated by
$n+m-1$ \elts (localize and use Corollary~\ref{corlemdArtin}~\emph{4}).
 
\emph{2.}
Show that if $\fa$ and $\fb$ are \lopsz, $\fa\fb$ is \lopz. Explain how to construct a \mlp for the \coes of $fg$ from two \mlpsz, respectively for the \gtrs of $\fa$ and for those of $\fb$.
}
\end{exercise}

\vspace{-1em}
\begin{exercise} 
\label{exogcdlcm} 
{\rm We are interested in the eventual \egt 
\begin{equation}\label{eqexogcdlcm}\preskip.4em \postskip.4em
\fa\,\fb =(\fa\cap\fb)(\fa+\fb)
\end{equation}
for two \itfs $\fa$ and $\fb$ of a \ri $\gA$.
  
\emph{1.}  Show that the \egt is satisfied if $\fa+\fb$ is \lopz.
If in addition $\fa$ and $\fb$ are \lopsz, then $\fa\cap \fb$ is \lopz.
 
\emph{2.}  Suppose $\gA$ is integral. Show that if the \egt is satisfied when $\fa$ and $\fb$ are \idps then the \ri is \ariz. 
 
\emph{3.}  Show that \propeq

\vspace{-.5em}
\begin{itemize}\itemsep0pt
\item $\gA$ is a \adpz.
\item $\gA$ is a \lsd and \Eqrf{eqexogcdlcm} is satisfied for \idpsz.
\item $\gA$ is a \lsd and \Eqrf{eqexogcdlcm} is satisfied for \itfsz. 
\end{itemize}
} 
\end{exercise}

\vspace{-1em}
\begin{exercise}
\label{exoSECSci} {\rm  (See also Exercise \ref{exoPgcdPpcm}) 
 Let $\fa$, $\fb$, $\fc$ be \itfsz. 
Prove the following statements.
 
\emph{1.} If $\fa+\fb$ is \lopz, then $(\fa:\fb)+(\fb:\fa)=\gen{1}$.
 
\emph{2.} If $(\fa:\fb)+(\fb:\fa)=\gen{1}$, then
\hsu 
\emph{a.}
$(\fa+\fb):\fc = (\fa:\fc)+(\fb:\fc)$.
\hsu 
\emph{b.} $\,\fc :(\fa~\cap~ \fb)=(\fc:\fa)+(\fc:\fb)$.
\hsu 
\emph{c.} $(\fa+\fb)(\fa~\cap~ \fb) = \fa~\fb$.
\hsu
\emph{d.} $\fc\,(\fa~\cap~ \fb)= \fc~ \fa~ \cap~ \fc~ \fb$.
\hsu 
\emph{e.} $\,\fc+(\fa~\cap~ \fb)=(\fc+\fa)\cap (\fc+\fb)$.
\hsu 
\emph{f.} $\,\fc~\cap\,(\fa+ \fb)=(\fc~\cap~ \fa)+(\fc~\cap~ \fb)$.
\hsu 
\emph{g.} The following short exact sequence (where $\delta(x)=(x, -x)$ and $\sigma(y, z)=y+z$) is split:
$$\preskip-.4em \postskip.0em 
0 \longrightarrow \fa~\cap~ \fb \vers{\delta} \fa \times \fb 
\vers{\sigma} \fa+\fb
\longrightarrow 0. 
$$
}
\end{exercise}

\vspace{-2.5em}
\begin{exercise}
\label{exoGaussien} (Gaussian \risz)
{\rm A \ri $\gA$ is said to be \emph{Gaussian} when for all \pols $f$, $g\in\AX$,
we have the \egt $\rc(fg)=\rc(f)\rc(g)$. 
Prove the following statements.
 
\emph{1.} Every \anar is Gaussian (see Exercise \ref{exoIlops}).

\emph{2.} A Gaussian integral \ri is a \adpz. 

\emph{3.} A Gaussian reduced \ri is a \adpz.  A Gaussian \qiri is a \adpc  (see \thref{th.adpcoh}).
}
\end{exercise}

\vspace{-1.1em}
\pagebreak

\begin{exercise} 
\label{exononanar} (A useful \ri for counterexamples)\\
{\rm Let $\gK$ be a nontrivial \cdi and $V$ be a \Kev of dimension $2$.
Consider the \Klg $\gA=\gK\oplus V$  defined by $x$, $y\in V \Rightarrow xy=0$. Show that every \elt of $\gA$ is \iv or nilpotent (\cad $\gA$ is local \zedz), and that the \ri is \coh but not \ariz. However, every \itf that contains a \ndz \elt is equal to $\gen{1}$, a fortiori it is \ivz. 
} 
\end{exercise}

\vspace{-1em}
\begin{exercise}
\label{exoVascon}
{\rm  Let $\gA$ be a \dcd \coh \aloz. \\
Let $\fm=\Rad\gA$ and suppose that $\fm$ is flat over $\gA$.

\emph{1.} Show that $\gA$ is integral.

\emph{2.} Show that $\gA$ is a \advz.

NB: We do not assume that $\gA$ is nontrivial. 
}
\end{exercise}

\vspace{-1em}
\begin{exercise}\label{exoQuotientPlat}
{(Flat quotient of a flat module: a direct \demz)} \\
{\rm
Provide a direct \dem of the following implication of \thref{propPlatQuotientdePlat}: Let~$M$ be a flat \Amo and $K$ be a submodule of $M$ satisfying $\fa M \cap K = \fa K$ for every finitely generated \id $\fa$; then $M/K$ is flat.
}

\end {exercise}

\vspace{-1em}
\begin{exercise}
\label{propfreeplat}
{\rm  This exercise starts with a long introductory text.
A single question is posed, at the very end.
We specify in the following \dfn the construction of the direct sum 
$\bigoplus_{i\in I}M=M^{(I)}$ for an arbitrary (not \ncrt discrete) set $I$ and a module $M$.\footnote{Concerning the general notion of a family of sets
indexed by an arbitrary set, see 
 \cite[page 18]{MRR}; the construction of the direct sum of an arbitrary family of \Amos is explained on pages 54 et 55.} This allows us to show that every module is a quotient of a flat module (actually a free module, not necessarily projective from a constructive point of view!).

\smallskip 
{\bf Definition.} \label{deffree}
Let $I$ be an arbitrary set and $M$ be an  \Amoz.
We define the \emph{direct sum $M^{(I)}$} as a quotient set of the set of finite formal sums $\oplus_{k\in \lrbn}(i_k,x_{k})$, where $i_k\in I$ \hbox{and $x_k \in M $} for each~\hbox{$k\in\lrbn$}: such a formal sum is defined as being \prmt the family $(i_k, x_{k})_{k\in \lrbn}$.
\\
The \eqvc relation  that defines the \egt over $M^{(I)}$ is the \eqvc relation generated by the following \gui{\egtsz}:
\begin{itemize}
\item associativity and commutativity of the formal sums: we can reorder the family as we wish,
\item if $i_k=_Ii_\ell$ then $(i_k,x_{k})$ and $(i_\ell,x_{\ell})$  can be replaced by $(i_k,x_k+x_\ell)$   (\gui{contraction} of the list); we can write this rewriting in the following form: if $i=_Ij$ then $(i,x_i)\oplus(j,x_j)=(i,x_i+x_j)$; 
\item every term $(i,0_{M})$ can be deleted. 
\end{itemize}
The addition over $M^{(I)}$ is defined by concatenation, and the external law is defined by $a\cdot \oplus_{k\in \lrbn}(i_k,x_{k})=\oplus_{k\in \lrbn}(i_k,ax_{k})$.\\
Finally, the \emph{\Amo freely generated by $I$} is the module  $\Ae{(I)}$.

\smallskip The direct sum solves the corresponding \uvl \pbz, which we can schematize by the following graph for a family $(\varphi_i)_{i\in I}$ of \alis from $M$ to an arbitrary module $N$.

\smallskip \centerline{\small
\xymatrix @C=1.5cm @R=.6cm 
          {
                &&& M_i=M&\jmath_i(x)=(i,x) 
                \\
N\ar@{<-}[urrr]^{\varphi_i} \ar@{<-}[drrr]_{\varphi_j}\ar@{<-}[ddrrr]_{\varphi_\ell}
\ar@{<--}[rr]^(.6){\varphi!} && P\ar@{<-}[ur]_{\jmath_i}\ar@{<-}[dr]^{\jmath_j} \ar@{<-}[ddr]_(.5){\jmath_\ell}&& P=M^{(I)} \\
  &&& M_j=M 
          \\
&&& M_\ell=M 
    \\
}}
 
\smallskip 
Let $I$ be an arbitrary set with at least one element.
The \Amo  $\Ae {(I)}$, is called the \emph{module freely generated by the set $I$.} It solves the corresponding \uvl \pbz, which we can schematize by the following graph for a family $x=(x_i)_{i\in I}$ in an arbitrary module $N$.

\smallskip \centerline{\small
\xymatrix @C=1.5cm @R=.6cm 
          {
                &&& I \ar[dl]^{\jmath}
                &x(i)=x_i\\
N\ar@{<-}[urrr]^{x} 
\ar@{<--}[rr]^(.6){\psi!} &&\Ae{(I)}   &&\jmath(i)=(i,1)\\
}}

\smallskip Note that as a consequence, \emph{if $(x_{i})_{i\in I}$ is an arbitrary \sgr of the module~$N$, the latter is \isoc to a quotient of $\Ae{(I)}$}.

Let $I$ be an arbitrary set and $M$ be an \Amoz.
Prove that the module $M^{(I)}$ is flat \ssi  $M$ is flat.
In particular this shows that the free module $\Ae{(I)}$ is flat.
}
\end{exercise}

}

\sol

\exer{exoPlatsLecteur}

\emph{Over a Bézout domain $\gZ$, a module $M$ is flat \ssi it is torsion-free}.

We know that the condition is \ncrz. Let us prove that it is sufficient.\\
Consider a syzygy $LX$ in $M$ with $L=\dex{a_1\;\cdots\;a_n}$ and $X=\tra{\dex{x_1\;\cdots\;x_n}}$.
If the $a_i$'s are all null, we have $L\,\In=0$ and $\In X= X$, which explains $LX=0$ in $M$.\\
Otherwise, we write $\sum_i a_i u_i=g$ and $gb_i=a_i$, where $g$ is the gcd of the $a_i$'s. 
\\
We have $g(\sum_ib_ix_i)=0$, and since $M$ is torsion-free $\sum_ib_ix_i=0$. 
\\
The matrix $C=\big((u_ib_j)_{i,j\in\lrbn}\big)=UB$ with $B=\fraC1 g L$, is a \mlp \hbox{for $(\an)$}. 
Let $G=\In-C$, we have $CX=0$ \hbox{and $LC=L$}, \hbox{so $LG=0$}
\hbox{and $GX=X$}, which explains $LX=0$ in $M$.

\exer{exothPlat1}

\emph{1.} Let $\mu:P\to M$ be a \aliz. We know (\thref{thPlat1}) that $\mu$ factorizes through a \tf free module $L$: $\mu=\lambda\circ \psi$.

\snic {
\xymatrix {
P\ar[d]_{?} \ar[r]^{\psi}\ar[dr]^{\mu} &L\ar[d]^{\lambda} &
\\
N \ar[r]^{\pi} &M\ar[r] &0 
\\
}}

Since $L$ is free, we can write $\lambda=\pi\circ \nu$ with a \ali $\nu:L\to N$, and so $\mu=\pi\circ \varphi$ for $\varphi= \nu\circ \psi$.

\emph{2.} If the \prt is satisfied with $N=\gA^{(I)}$, where $I$ is a discrete set, we consider an arbitrary \ali $\mu:P\to M$ with $P$ \pfz. We write $\mu=\pi\circ \varphi$ with a \ali $\varphi:P\to N$. There then exists a finite subset $I_0$ of $I$ such that for each \gtr $g_j$ of $P$, $\varphi(g_j)$ has null \coos outside of $I_0$. This shows that we can factorize $\mu$ via the free module of finite rank $\gA^{(I_0)}$. Therefore by \thref{thPlat1}, $M$ is flat.


\exer{exoplatFitt}\\
Consider a \mtf $M$ with a \sgr $(\xn)$. 
Let $X=\tra{\lst{x_{1}\;\cdots\; x_{n}}}$. 
For $k\in\lrbzn$ and $k+r=n$,
a typical \gtr of~$\cF_k(M)$ is $\delta=\det(L)$ where $L\in\MM_r(\gA)$ and $LY=0$, for a column vector extracted from $X$: $Y=\tra{\lst{x_{i_1}\;\cdots\;x_{i_r}}}$.
\\
We must show that $\delta\in\cF_k(M)^{2}$. Actually we will show that $\delta\in\delta\,\cF_k(M)$.
\\
Suppose \spdg that $(i_1,\dots,i_r)=(1,\dots,r)$. We apply Proposition~\ref{propPlat2}.
We therefore have a matrix $H\in\MM_{r,n}$ with $HX=Y$ and $LH=0$.\\
Let $H'=\I_{r,r,n}=\blocs{.9}{.7}{.9}{0}{$\I_{r}$}{$0$}{}{ }$, and $K=H'-H$.
We have 

\snic{KX=Y-Y=0\;\hbox{  and  }\;LK=LH'=\blocs{.9}{.7}{.9}{0}{$L$}{$0$}{}{ }\,.}

\snii
Let $K'$ be the matrix formed by the first $r$ columns of $K$. Then $L=LK'$ \hbox{and $\det(L)=\det(L)\det(K')$}, and since $KX=0$, we have $\det(K')\in\cF_k(M).$  


\exer{exoClamlsdz}
Suppose the \ri $\gA$ is a \lsdz. 
Let $\fp$ be a \idep and $xy=0$ in $\gA_\fp$.
There exists a $u\notin\fp$ such that $uxy=0$ in~$\gA$. Let $s$ and $t\in\gA$ such that $s+t=1$, $sux=0$ and $ty=0$ in $\gA$. The \elts $s$ and~$t$ 
        cannot both be in $\fp$  
(otherwise $1\in\fp$). If $s\notin\fp$, then since $sux=0$, we obtain $x=_{\gA_\fp}0$. If $t\notin\fp$, then since $ty=0$, we obtain $y=_{\gA_\fp}0$. 
Thus $\gA_\fp$ is an integral \riz.
\\
Now suppose that every localized \ri $\gA_\fp$ at every \idema $\fp$
is integral and suppose that $xy=_\gA0$. For some arbitrary \idema $\fp$ we have $x=_{\gA_\fp}0$ \hbox{or $y=_{\gA_\fp}0$}. 
In the first case let $s_\fp\notin\fp$ such that $s_\fp x=_\gA 0$. Otherwise let $t_\fp\notin\fp$ such that $t_\fp y=_\gA 0$. The family of the $s_\fp$'s or $t_\fp$'s generates the \id $\gen{1}$ (because otherwise all $s_\fp$'s or $t_\fp$'s would be in some \idemaz).
\\
There is therefore a finite number of $s_i$'s satisfying $s_ix=0$  (in $\gA$) and a finite number of $t_j$'s satisfying $t_jy=0$, with an \eqn $\sum_ic_is_i+\sum_jd_jt_j=1$. 
\\
We take $s=\sum_ic_is_i$, $t=\sum_jd_jt_j$ and we obtain
$sx=ty=0$ and $s+t=1$.


\exer{exoClamAri}
We begin by recalling the following: by item \emph{3} of \thref{propmlm}, an \id $\gen{a,b}$ of a \ri $\gA$ is \lop \ssi we can find $s$, $t$, $u$, $v\in\gA$ such that $s+t=1$, $sa=ub$ and $tb=va$.
\\
Suppose the \ri $\gA$ is \ariz. 
Let $\fp$ be a \idepz. For $a$, $b\in\gA_\fp$ we want to show that $a$ divides $b$ or $b$ divides~$a$ (see Lemma~\ref{lemBezloc}). \Spdg we can take $a$ and~$b$ in $\gA$. Then let $s$, $t$, $u$, $v$ be as above. The \elts $s$ and~$t$ cannot both be in $\fp$ (otherwise $1\in\fp$). If $s\notin\fp$, then $a=_{\gA_\fp} s^{-1}ub$ so $b$ divides $a$ in~$\gA_\fp$. If $t\notin\fp$, then $a$ divides $b$ in~$\gA_\fp$. 
\\
Now suppose that every localized \ri $\gA_\fp$ 
at every \idema $\fp$ is a local Bézout \ri and let $a,\,b\in\gA$. 
\\
For an arbitrary \idema $\fp$, we have that $b$ divides $a$ or $a$ divides $b$ in~$\gA_\fp$. 
In the first case let $s_\fp\notin\fp$ and $u_\fp\in\gA$ such that $s_\fp a=_\gA u_\fp b$. Otherwise let $t_\fp\notin\fp$ and $v_\fp$ such that $t_\fp b=_\gA v_\fp a$. The family of the $s_\fp$'s or $t_\fp$'s generates the \idz~$\gen{1}$ (because otherwise all~$s_\fp$'s or $t_\fp$'s would be in some \idemaz).
\\
Therefore there is a finite number of $s_i$'s, $u_i$'s satisfying 
$s_ia=u_ib$  (in $\gA$) and a finite number of $t_j$'s, $v_j$'s satisfying 
$t_jb=v_ja$, with an \eqn $\sum_ic_is_i+\sum_jd_jt_j=1$. 
\\
We take $s=\sum_ic_is_i$, $u=\sum_ic_iu_i$, $t=\sum_jd_jt_j$, $v=\sum_jd_jv_j$ and we obtain the \egts $s+t=1$, $sa=ub$ and $tb=va$.
\\
For an \id with a finite number of \gtrsz, we can reason analogously, or use the result of Exercise~\ref{exoMLP}.


\exer{exoIlops0} \\
The image of the \idp $\gen{60}$ of $\ZZ$ under the \homo $\ZZ\to\ZZ/27\ZZ$ is the \idz~$\gen{3}$ which does not contain any \ndz \eltz, and which is not \ivz. Actually, as a $\ZZ/27\ZZ$-module, the \id $\gen{3}$ is not even \pro (its annihilator~$\gen{9}$ is not \idmz).\\
When $\rho:\gA\to\gB$ is a flat \algz, the image of an \id $\fa\subseteq \gA$ is \isoc to $\rho\ist(\fa)\simeq\gB\otimes_\gA\fa$. Therefore if $\fa$ is \ivz, as it is \pro of rank $1$, its image is also a \pro module of rank $1$.


\exer{exoIlops2} \\
We first note that a product of \lops \ids is always \lot principal, because after suitable \come \lonsz, each \id becomes principal, and so does their product.
\\
We are then content with the $\fa=\gen{a,b}$ case and with the $\gen{a^{4},b^{4}}$ example. It will be clear that the computation technique is easily \gneez.\\
We start with $sa=ub$, $tb=va$ and $s+t=1$. Therefore  $s^{4}a^{4}=u^{4}b^{4}$ and $t^{4}b^{4}=v^{4}a^{4}$. Since $\gen{s^{4},t^{4}}=\gen{1}$ (which is obtained by writing $1=(s+t)^{7}$), we indeed obtain that the \id $\gen{a^{4},b^{4}}$ is \lopz. 
\\
Let us show, for example, that  $a^{2}b^{2}\in\gen{a^{4},b^{4}}$.\\
We write $s^{2}a^{2}=u^{2}b^{2}$ and $t^{2}b^{2}=v^{2}a^{2}$. Therefore $s^{2}a^{2}b^{2}=u^{2}b^{4}$ and $t^{2}a^{2}b^{2}=v^{2}a^{4}$.\\
Finally, $1=(s+t)^{3}=s^{2}(s+3t)+t^{2}(t+3s)$. Therefore 

\snic{a^{2}b^{2}=(t+3s)v^{2}a^{4}+(s+3t)u^{2}b^{4}.}


\exer{exoMLP} We do not need to assume that the \ri is \ariz.\\
We will show that if in a \ri $\gA$ each pair $(a_i,a_j)$ admits a \mlpz,
the same goes for the $n$-tuple $(\an)$.
\\
This can be compared to Dedekind's proof of \thref{th1IdZalpha}, which concerns only invertible ideals, because over an integral ring the invertible ideals are precisely the nonzero locally principal ideals.
\\
Also note that the result is a priori clear: by successive \come \lonsz, every leaf of each branch of an a priori very large computation tree will be a \idpz. This will show that the \id $\gen{\an}$ is always generated by one of the $a_i$'s after \lons at \ecoz.
What we are aiming for here is rather a practical computation of the \mlpz.
\\
We proceed by \recu on $n$. 

Let us show the \recu step for the passage of $n=3$ to $n+1=4$.\label{corexoPruf3}\\
Consider $a_1$, $a_2$, $a_3$, $a_4\in\gZ$. \\
By \hdr we have a matrix
$C=\cmatrix{x_1&x_2&x_3\cr y_1&y_2&y_3\cr z_1&z_2&z_3\cr }$ suitable for $(a_1,a_2,a_3)$, and  matrices $\cmatrix{c_{11}&c_{14}\cr d_{11}&d_{14}}$, $\cmatrix{c_{22}&c_{24}\cr d_{22}&d_{24}}$, $\cmatrix{c_{33}&c_{34}\cr d_{33}&d_{34}}$ respectively suitable for  $(a_1,a_4)$,  $(a_2,a_4)$ and  
$(a_3,a_4)$.
Then we will check that the \und{trans}p\und{ose} of the following matrix is suitable for $(a_1,a_2,a_3,a_4)$
$$
\cmatrix{c_{11}x_1&c_{22}y_1&c_{33}z_1&d_{11}x_1+d_{22}y_1+d_{33}z_1 
     \cr c_{11}x_2&c_{22}y_2&c_{33}z_2&d_{11}x_2+d_{22}y_2+d_{33}z_2  
     \cr c_{11}x_3&c_{22}y_3&c_{33}z_3& d_{11}x_3+d_{22}y_3+d_{33}z_3 
     \cr c_{14}x_1
        &c_{24}y_2 
        & c_{34}z_3
        & d_{14}x_1+d_{24}y_2+d_{34}z_3}
$$
First of all, we must check that the trace of the matrix is equal to $1$, \cad
 
\snic{t=c_{11}x_1+c_{22}y_2+c_{33}z_3+d_{14}x_1+d_{24}y_2+d_{34}z_3=1,}

\snii but $c_{11}+d_{14}=1=c_{22}+d_{24}=c_{33}+d_{34}$ so $t=x_1+y_2+z_3=1$.
\\
We must check that each row of the transposed matrix is proportional to $\lst{a_1\;a_2\;a_3\;a_4}$.
Two cases arise. First of all, we consider one of the first three rows,
for instance 
the row $\lst{c_{11}x_1\;c_{11}x_2\;c_{11}x_3\;c_{14}x_1}$. Both of the following types of equalities must be satisfied

\snic{a_1c_{11}x_2=a_2c_{11}x_1, \quad \hbox{and}\quad a_1c_{14}x_1=a_4c_{11}x_1.}

For the first \egt we use $a_2x_1=a_1x_2$ and for the second $a_1c_{14}=a_4c_{11}$.\\
Finally, we must verify that $[\,a_1\;a_2\;a_3\;a_4\,]$ is proportional to the transpose of

\snac{\cmatrix{d_{11}x_1+d_{22}y_1+d_{33}z_1\\
d_{11}x_2+d_{22}y_2+d_{33}z_2\\
d_{11}x_3+d_{22}y_3+d_{33}z_3\\
d_{14}x_1+d_{24}y_2+d_{34}z_3}.}

\snii This results on the one hand from the proportionality of $[\,a_1\;a_2\;a_3\,]$ to each of the rows~\hbox{$[\,x_i\;y_i\;z_i\,]$}, and on the other hand from the proportionality of the rows $[\,a_i\;a_4\,]$ to the rows $[\,d_{i1}\;d_{i4}\,]$.

To complete the \demz, note that the passage of $n-1$ to $n$ (for any $n>2$) is perfectly analogous.

\exer{exoIlops} We write $\fa=\gen{\am}$, $\fb=\gen{\bn}$.
We can assume \hbox{that $f=\sum_{k=1}^{m}a_kX^{k-1}$} and $g=\sum_{h=1}^{n}b_hX^{h-1}$.

\emph{1.} Let $F$ be a \mlp for $(\am)$. If $\rc(f)=\fa$, we have \eco $s_i$ (the diagonal of $F$) and \pols $f_i\in\AX$ (given by the rows of $F$) that satisfy the \egts $s_if=a_if_i$ in $\AX$. In addition, the \coe \hbox{of $X^{i-1}$} in $f_i$ is equal to $s_i$, so $\rc(f_i)\supseteq \gen{s_i}.$
\\ 
By letting $\gA_i=\gA\big[\fraC1{s_i}\big]$, we have $\rc(f_i)=_{\gA_i}\! \gen{1}$ and the \egts

\snic{s_i\rc(fg)=\rc(a_if_ig)=a_i\rc(f_ig)=_{\gA_i}\!a_i\rc(g)=_{\gA_i}\!\rc(a_if_i)\rc(g)=s_i\rc(f)\rc(g) }

\snii (the third \egt comes from Corollary~\ref{corlemdArtin}~\emph{4}
because $\rc(f_i)=_{\gA_i}\! \gen{1}$).\\
Hence the \egt $\rc(fg)=\rc(f)\rc(g)=\fa\fb$ because it is true in each $\gA_i$.
 
\emph{2.} If $g$ is also \lop we obtain  $t_jb=b_jg_j$ \hbox{in $\AX$}, with $\rc(g_j)\supseteq \gen{t_j}$ and \com $t_j$ in $\gA$.
We therefore have 

\snic{s_it_j\rc(fg)=_{\gA_{ij}}\!a_ib_j\rc(f_ig_j)=_{\gA_{ij}}\!\gen{a_ib_j}.}

This tells us that the \id $\rc(fg)=\fa\,\fb$ becomes principal after $mn$ \come \lonsz. As this \id admits $m+n-1$ \gtrs (the \coes of~$fg$) there is a \mlp for these \gtrsz.\\ 
To compute it, we can use the \dem of the implication \emph{1} $\Rightarrow$  \emph{3} in \thref{propmlm}. 
This \dem is quite simple, as well as the computation it implies. But if we examine in detail what is going to happen, we realise that in the \dem below we have used the Gauss-Joyal lemma: over the \ri $\gA_{ij}$, we have $1\in\rc(f_i)\rc(g_j)$ because $1\in\rc(f_i)$ and $1\in\rc(g_j)$. This lemma admits several \elr \dems (see \ref{lemGaussJoyal} and \ref{corlemdArtin}), but none of them gives a simple formula that allows us to provide the \coli of the \coes of $fg$ equal to $1$, from two \colis of the \coes of~$f$ and of those of $g$.
\\
We would be grateful to any reader who is able to indicate to us a
 short direct computation, for example in the case where the \ri is integral with explicit \dvez.\footnote{Please note that in the case of an integral \ri with explicit \dvez, a \mlp is known from its only diagonal \eltsz, which can simplify computations.} 

\exer{exogcdlcm} We write $\fa=\gen{\an}$, $\fb=\gen{\bbm}$.
\\
We will use the result of Exercise \ref{exoMLP} which shows that if every \id with two \gtrs is \lopz, then every \itf is \lopz. 

\emph{1.} In Exercise \ref{exoPgcdPpcm} item \emph{4} we have shown that $1 \in (\fa : \fb) + (\fb :\fa)$, $\fa \cap \fb$ is \tf and  $\fa\fb = (\fa \cap \fb)(\fa + \fb)$.\\
If $\fa+\fb$ is \lopz, there is a \sys of \eco such that by inverting any one of them, the \id is generated by some $a_k$ or some $b_\ell$. But if  $\fa+\fb=\gen{a_k}\subseteq \fa$, we have $\fb\subseteq \fa$, so $\fa\cap\fb=\fb$, \lop by hypothesis.
Thus $\fa\cap\fb$ is \lop because it is \lop after \lon at \ecoz.

\emph{2.} If the \ri is integral and if  $(\fa+\fb)(\fa\cap\fb)=\fa\fb$ for $\fa=\gen{a}$ and $\fb=\gen{b}$ (\hbox{where $a,b\neq 0$}), we get that $\gen{a,b}(\fa\cap\fb)=\gen{ab}$, so $\gen{a,b}$ is \iv (and also~\hbox{$\gen{a}\cap\gen{b}$} at the same time). When it is satisfied for all $a,b\neq 0$, the \ri is \ariz. 

\emph{3.} The only delicate implication consists in showing that if $\gA$ is a \lsd and if $(\fa+\fb)(\fa\cap\fb)=\fa\fb$ when $\fa=\gen{a}$ and $\fb=\gen{b}$ then the \ri is \ariz, in other words every \id $\gen{a,b}$ is \lopz. 

If $\gen{a,b}(\fa\cap\fb)=\gen{ab}$, we write $ab=au+bv$ with $u$ and $v\in\fa\cap\fb$:

\snic{u=ax=by,\;v=az=bt,\;\;\hbox{hence} \;\;au+bv=ab(y+z)=ab.}

Since the \ri is a \lsdz, from the \egtz~\hbox{$ab(y+z-1)=0$}, we deduce three \come \lons in which we obtain~\hbox{$a=0$}, $b=0$ and $1=y+z$ respectively. In the first two cases $\gen{a,b}$ is principal. In the last case $\gen{a,b}$ is \lop (localize at $y$ or at~$z$).


\exer{exoSECSci} We write $\fa=\gen{\an}$, $\fb=\gen{\bbm}$.

\emph{1.} Proven in item \emph{4} of Exercise \ref{exoPgcdPpcm}.

\emph{2.} Now suppose $(\fa:\fb)+(\fb:\fa)=\gen{1}$, \cad we have $s$, $t\in\gA$ with 

\snic{s+t=1$, $s\fa\subseteq \fb$, $t\fb\subseteq \fa.}

\emph{2a.}
$(\fa+\fb):\fc = (\fa:\fc)+(\fb:\fc)$. In this \egt as in the following (up to~\emph{2f}), an inclusion is not obvious (here it is $\subseteq $). Proving the non-obvious inclusion  comes down to solving a \sli (here, given some $x$ such \hbox{that $x\fc\subseteq \fa+\fb$}, we look for $y$ and $z$ such that $x=y+z$, $y\fc\subseteq \fa$ and $z\fc\subseteq \fb$). \\
We can therefore use the basic \plg with the \ecoz~$s$ and~$t$.\\
When we invert $s$, we get $\fa\subseteq \fb$, and if we invert $t$, we get $\fb\subseteq \fa$. In both cases the desired inclusion becomes trivial.

For the record: 
\emph{2b.} $\,\fc :(\fa~\cap~ \fb)=(\fc:\fa)+(\fc:\fb)$.
\,
\emph{2c.} $(\fa+\fb)(\fa~\cap~ \fb) = \fa~\fb$.
\\
\emph{2d.} $\fc\,(\fa~\cap~ \fb)= \fc~ \fa~ \cap~ \fc~ \fb$.
\,
\emph{2e.} $\,\fc+(\fa~\cap~ \fb)=(\fc+\fa)\cap (\fc+\fb)$.
\\
\emph{2f.} $\,\fc~\cap\,(\fa+ \fb)=(\fc~\cap~ \fa)+(\fc~\cap~ \fb)$.

\emph{2g.} The short exact sequence below (where $\delta(x)=(x, -x)$ and $\sigma(y, z)=y+z$) is split:
$$\preskip-.4em \postskip.4em 
0 \longrightarrow \fa~\cap~ \fb \vers{\delta} \fa \times \fb 
\vers{\sigma} \fa+\fb
\longrightarrow 0. 
$$
We want to define $\tau:\fa+\fb\to\fa\times \fb$ such that $\sigma\circ \tau=\Id_{\fa+\fb}$.\\
If $\fa\subseteq \fb$, we can take $\tau(b)=(0,b)$ for all $b\in\fb=\fa+\fb$.
If $\fb\subseteq \fa$, we can take $\tau(a)=(a,0)$ for all $a\in\fa=\fa+\fb$. \\
In the first case this implies $s\tau(a_i)=\big(0,\sum_jx_{ij}b_j\big)$ and $s\tau(b_j)=(0,sb_j)$.
\\
In the second case this implies $t\tau(b_j)=\big(\sum_iy_{ji}a_i,0\big)$ and $t\tau(a_i)=(ta_i,0)$.
\\
We therefore try to define $\tau$ by the following formula which coincides with the two previous ones in the two special cases.

\snic{\tau(a_i)=\big(ta_i,\sum_jx_{ij}b_j\big),\;\;\tau(b_j)=\big(\sum_iy_{ji}a_i,sb_j\big) .}

For this attempt to succeed, it is necessary and sufficient that when $\sum_i\alpha_ia_i=\sum_{j}\beta_jb_j$, we have the \egt

\snic{\sum_i\alpha_i\big(ta_i,\sum_jx_{ij}b_j\big)=\sum_{j}\beta_j\big(\sum_iy_{ji}a_i,sb_j\big) .}

\vspace{-.3em}
\pagebreak

For the first coordinate, this results from the following computation (and similarly for the second coordinate).
$$\preskip.4em \postskip.4em \ndsp 
\sum_i\alpha_ita_i=t\sum_i\alpha_ia_i=t\sum_{j}\beta_jb_j=\sum_{j}\beta_jtb_j=\sum_{j}\beta_j\sum_iy_{ji}a_i. 
$$
Finally, the \egt $\sigma\circ \tau=\Id_{\fa+\fb}$ is satisfied because it is satisfied when restricted to $\fa$ and $\fb$ (\imd computation).


\exer{exoGaussien}  
 
\emph{1.} Proven in Exercise \ref{exoIlops}.

\emph{2.} Let $a$, $b$, $c$, $d\in\gA$. Let $\fa=\gen{a,b}$\\
Consider $f=aX+b$ and $g=aX-b$. We get $\gen{a,b}^{2}=\geN{a^{2},b^{2}}$, \cad

\snic{ab=ua^{2}+vb^{2}.}

\snii
When considering $f=cX+d$ and $g=dX+c$, we obtain $\gen{c,d}^{2}=\geN{c^{2}+d^{2},cd}$. In other words $c^{2}$ and $d^{2}\in\geN{c^{2}+d^{2},cd}$.  
\\
Let $\fb=\gen{ua,vb}$. We have $ab\in\fa\fb$. It suffices to show that $\fa^{2}\fb^{2}=\geN{a^{2}b^{2}}$ because this implies that $\fa$ is \iv (we treat the case $a$, $b\in\Atl$).
However, we have

\snic{a^{2}b^{2}\in\fa^{2}\fb^{2}=\geN{a^{2},b^{2}}\geN{u^{2}a^{2},v^{2}b^{2}}.}

We therefore need to show that $u^{2}a^{4}$ and $v^{2}b^{4}\in\geN{a^{2}b^{2}}.$
Let $u_1=ua^{2}$ and  $v_1=vb^{2}$. We have $u_1+v_1=ab$ and $u_1v_1\in\geN{a^{2}b^{2}}$. Therefore  $u_1^{2}+v_1^{2}\in\geN{a^{2}b^{2}}$ also.\\
Since $u_1^{2}\in\geN{u_1^{2}+v_1^{2},u_1v_1}$, we indeed get $u_1^{2}\in\geN{a^{2}b^{2}}$ (likewise \hbox{for $v_1^{2}$}).

\emph{3.} The \egts of item \emph{2} are all satisfied.
\\
Let us first show that the \ri is a \lsdz. 
\\
Assume $cd=0$. Since $c^{2}\in\geN{c^{2}+d^{2},cd}$, we have $c^{2}=x(c^{2}+d^{2})$, \cad

\snic{xd^{2}=(1-x)c^{2}.}

\snii
We deduce that $xd^{3}=0$, and as $\gA$ is reduced, $xd=0$. Similarly $(1-x)c=0$.

Let us now see that the \ri is \ariz. We start from arbitrary $a$, $b$ and we want to show that $\gen{a,b}$ is \lopz.
By item \emph{2} we have an \id $\fc$ such that $\gen{a,b}\fc=\geN{a^{2}b^{2}}$. We therefore have $x$ and $y$ with

\snic{\gen{a,b}\gen{x,y}= \geN{a^{2}b^{2}} \;\hbox{ and } ax+by=a^{2}b^{2}.}

We write $ax=a^{2}b^{2}v$ and $by=a^{2}b^{2}u$.
From the \egt $a(ab^{2}v-x)=0$, we deduce two \come \lonsz, in the first $a=0$, in the second $x=ab^{2}v$.
We therefore suppose \spdg that $x=ab^{2}v$ and, symmetrically $y=ba^{2}u$.
This gives 

\snic{\gen{a,b}\gen{x,y}=ab\gen{a,b}\gen{au,bv}=\geN{a^{2}b^{2}}.}

\snii
We can also suppose \spdg that $\gen{a,b}\gen{au,bv}=\gen{ab}$.
\\
We also have $ax+by=a^{2}b^{2}(u+v)$ and since $ax+by=a^{2}b^{2}$, we suppose \spdg that $u+v=1$. 
\\
Since $a^{2}u=abu'$, we suppose \spdg that $au=bu'$. 
\\
Symmetrically $bv=av'$, and since $u+v=1$, $\gen{a,b}$ is \lopz.

\exer{exoVascon}
\emph{1.} Let $a\in\gA$ and $a_1$, \dots,  $a_n\in \gA$  generate $\fa=\Ann(a)$.
 If one of the $a_i$'s is in $\Ati$, we obtain $a=0$ and $\fa=\gen{1}$.
It remains to treat the case where all the~$a_i$'s are in $\fm$.
 Let $b$ be one of the $a_i$'s. Since $\fm$ is flat and $b\in\fm$, the equality $ab=0$ 
 gives us \eltsz~\hbox{$c_1$, \dots, $c_m\in\fa$} and $b_1$, \dots, $b_m\in\fm$ \hbox{with
  $b=\sum_{i\in\lrbm}c_ib_i$}. Therefore $b\in\fa\fm$, which gives $b=\sum_{i\in\lrbn}a_iz_i$
for some $z_i\in\fm$.
Hence a matrix \egt 
$$
[\,a_1\;\cdots\;a_n\,]=M\,[\,a_1\;\cdots\;a_n\,]\quad \hbox{ with } M\in\Mn(\fm).
$$
Thus $[\,a_1\;\cdots\;a_n\,] (\In-M)=[\,0\;\cdots\;0\,]$ with $\In-M$ \ivz, so $\fa=0$.

\emph{2.} Consider $a$, $b\in\gA$. We must prove that one divides the other. 
If one of the two is \ivz, the case is closed. 
It remains to examine the case where $a$ and $b\in\fm$. 
\\ 
We consider a matrix 
$$
P=\cmatrix{a_1&\cdots&a_n\cr b_1&\cdots&b_n}
$$ 
whose columns generate the module $K$, the kernel of $(x,y)\mapsto bx-ay$. In particular we have $a_ib=b_ia$ for each $i$. If one of the $a_i$'s or $b_i$'s is \ivz, the case is closed. 
It remains to examine the case where the $a_i$'s and $b_i$'s are in $\fm$.
\\
 Let $(c,d)$ be one of the $(a_i,b_i)$'s. Since $\fm$ is flat and $a$, $b\in\fm$, the equality 
$cb-da=0$ gives 
$$\preskip-.4em \postskip.4em
\cmatrix{c\cr d}=\cmatrix{c_1&\cdots&c_m\cr d_1&\cdots&d_m} \cmatrix{y_1\cr\vdots\cr y_m} \quad \hbox {with the }y_i's\in\fm\; \hbox{ and the }\cmatrix{c_j\cr d_j}'s\in K.
$$
By expressing the \smash{$\cmatrix{c_j\cr d_j}$} as \colis of the columns of $P$
we obtain
$$
\cmatrix{c\cr d}=P \cmatrix{z_1\cr\vdots\cr z_n} \quad \hbox {with every }z_i\in\fm.
$$
Hence $P=PN$ with a matrix $N\in\Mn(\fm)$, so $P=0$. This implies that $(a,b)=(0,0)$, and $a$ divides $b$ (actually, in this case, $\gA$ is trivial).

\hum{on peut se demander if l'hypoth\`ese  \gui{\dcdz} is vraiment \ncrz.}
  

\exer{exoQuotientPlat}
Let $a_i \in \gA$ and $x_i \in M$ satisfy $\sum_{i=1}^n a_i x_i \equiv 0 \bmod K$, a relation that we must explain. Let $\fa = \gen {\ua}$ such that $\fa K = \sum_i a_iK$; 
since $\sum_i a_ix_i \in \fa M \cap K = \fa K$, we have an equality 
$\sum_i a_ix_i = \sum a_iy_i$ where each
 $y_i \in K$.
We therefore have, with $z_i = x_i-y_i$, the relation $\sum_i a_iz_i = 0$ in $M$. Since $M$ is flat, this relation produces a certain number of vectors of $M$, say $3$ for simplicity, denoted by $u$, $v$, $w$ and $3$ sequences of scalars $\underline\alpha =(\alpha_1, \ldots,\alpha_n)$,
$\underline\beta =(\beta_1, \ldots,\beta_n)$ and
$\underline\gamma =(\gamma_1, \ldots,\gamma_n)$,  satisfying

\snic {
(z_1, \ldots, z_n) = (\alpha_1, \ldots,\alpha_n)\,u +
 (\beta_1, \ldots,\beta_n)\,v + (\gamma_1, \ldots,\gamma_n)\,w
}

\snii
and $\scp {\ua}{\underline\alpha} = \scp {\ua}{\underline\beta} = 
\scp {\ua}{\underline\gamma} = 0$.

\snii
Since $z_i \equiv x_i \bmod K$, we obtain our sought explanation in $M/K$:

\snic {
(x_1, \ldots, x_n) \equiv (\alpha_1, \ldots,\alpha_n)\,u +
 (\beta_1, \ldots,\beta_n)\,v + (\gamma_1, \ldots,\gamma_n)\,w
\;\mod K.
}


\exer{propfreeplat}
\Spdg suppose that $I$ is finitely enumerated. In other words $I=\so{i_1,\dots,i_n}$.
Let $P=M^{(I)}$. Any element of $P$ can be written in the form 
$x=\oplus_{k\in\lrbn}(i_k,x_k)$.

\emph{First suppose that the module $M$ is flat}, and consider a syzygy in $P$

\snic{0=\sum_{\ell\in\lrbm}a_\ell x_\ell=\sum_{\ell\in\lrbm}a_\ell\big(\oplus_{k\in\lrbn} (i_k,x_{k,\ell})\big)=\oplus_{k\in\lrbn}(i_k,y_k)}

with $y_k=\sum_{\ell\in\lrbm}a_\ell x_{k,\ell}$.\\
By \dfn of  \egt in $P$, since $\oplus_{k\in\lrbn}y_k=0$, we are in (at least) one of the two following cases:
\begin{itemize}
\item all the $y_k$'s are null,
\item two indices are equal in $I$: $i_k=_Ii_h$ for $h$ and $k$ distinct \hbox{in $\lrbn$}. 
\end{itemize} 
The first case is treated like that of a direct sum over a finite $I$.
The second case  reduces to the first by \recu on $n$.

\emph{Now suppose that $P$ is flat} and consider for instance the index~\hbox{$i_1\in I$} and a syzygy $\sum_{\ell\in\lrbm}a_\ell x_\ell=0$ in $M$.
We explain this syzygy in $P$ by writing 

\snic{(i_1,x_\ell)=_P\sum_{j\in \lrbp} g_{\ell,j}z_{j}$ with $\sum_{\ell\in\lrbm} a_\ell g_{\ell,j}=_\gA0$ for all $j.}

We re-express $z_{j}=\oplus_{k\in\lrbn}(i_k, y_{k,j})$, which gives

\snic{(i_1,x_\ell)=_P\oplus_{k\in\lrbp}\big(i_k,\sum_{j\in \lrbp}g_{\ell,j}y_{k,j}\big).}

By \dfn of the \egt in $P$, we are in (at least) one of the two following cases:
\begin{itemize}
\item for each $\ell$, we have  $x_\ell=\sum_{j\in \lrbn} g_{\ell,j}y_{1,\ell}$ in $M$,
\item we have in $I$: $i_1=_Ii_h$ for some $h\neq 1$ \hbox{in $\lrbn$}. 
\end{itemize} 
In the first case we have in $M$ the \egts that suit us.
The second case is reduced to the first by \recu on $n$. 
%

\Biblio

\ddps were introduced by H. Pr\"ufer in 1932 in \cite{Prufer}.
Their central place in multiplicative \id theory is showcased in the book of reference on the subject \cite{Gil}. 
Even though they were introduced in a very concrete way as the integral \ris in which every nonzero \itf is \ivz, 
this \dfn is often set aside in the modern literature for the following purely abstract alternative, 
        which only works  
in the presence of non-\cofs principles (\TEM and the axiom of choice): 
the \lon at any \idep gives a \advz. 

\Anars were introduced by L.\ Fuchs in 1949 in~\cite{Fuchs}. 
 
In the case of a non-integral \riz, the \dfn that we have adopted for \adps is due to Hermida and S\'anchez-Giralda \cite{HS}. 
It is the one that seemed the most natural to us, given the central importance of the concept of flatness in commutative \algz. Another name for these \ris in the literature is \emph{\ri of weak global dimension less than or equal to one}, which is rather inelegant. 
Moreover, we often find in the literature a Pr\"ufer \ri defined as a \ri in which every \id containing a \ndz \elt is \ivz. 
They are therefore almost \anarsz, but the behavior of the \ids that do not contain \ndz \elts seems utterly random (cf.\ Exercise~\ref{exononanar}).

A fairly complete presentation of arithmetic rings and Pr\"ufer rings, written in the style of constructive mathematics, can be found in
 \cite[Ducos\&al.]{dlqs} and~\cite[Lombardi]{lom99}.

A very comprehensive survey 
on variations of the notion of integral Pr\"ufer rings obtained by deleting the hypothesis of integrity
 is given in~\cite[Bazzoni\&Glaz]{BG}, including Gaussian \ris (Exercise~\ref{exoGaussien}).

\newpage \thispagestyle{CMcadreseul}

\incrementeexosetprob


\chapter{\Alosz, or just about}
\label{chap Anneaux locaux} 
\perso{compil\'e le \today}
\minitoc


\section{A few \cov \dfnsz}
\label{secAloc1}

In \clama a \alo is often defined as a \ri having a single \idemaz.
In other words the non-\iv \elts form an \idz.
This second \dfn has the advantage of being simpler (no quantification over the set of \idsz).
However, it lends itself fairly poorly to an algorithmic treatment
 because of the negation contained in \gui{non-\iv \eltsz.} 
This is the reason why we adopt the definition given on \paref{eqAloc} in constructive mathematics: if the sum of two \elts is \ivz, one of the two is \ivz.

We now find ourselves obligated to inflict a few unusual \dfns on the classic reader, in line with the \dfn of a \aloz.
Rest assured, on other planets, in other solar systems, no doubt the symmetric situation is taking place.
There, \maths has always been \covs and they have just barely discovered the interest of the abstract Cantorian point of view. 
An author in the new style is in the process of explaining that for them it is much simpler to regard a \alo as a \ri having a single \idemaz. Would the reader then put in the effort to follow what they are saying?

\subsec{The Jacobson radical, \alosz, fields}

Recall that for a \ri $\gA$ we denote by $\Ati$ the multiplicative group of \iv \eltsz, also called the group of units.

An \elt $x$ of a \ri $\gA$ is said to be \ix{non-invertible} if it satisfies{\footnote{Here we will use a slightly weakened version of negation. 
For a \prtz~$\sfP$ affecting \elts of the \ri $\gA$ or of an \Amo $M$, we consider the \prt $\sfP':=(\sfP\Rightarrow 1=_\gA0)$.
It is the negation of $\sfP$ when the \ri is not trivial.
Yet it often happens that a \ri constructed in a \dem can be trivial without one knowing.
To do an entirely \cof treatment of the usual classical \dem in such a situation (the classical \dem excludes the case of the trivial \ri by an ad hoc argument) our weakened version of negation then turns out to be \gnlt useful.
A \cdi  does not \ncrt satisfy the axiom of  discrete sets, $\forall x,y\;\big(x=y\;\hbox{or}\;\lnot(x=y)\big)$,
but it satisfies its weak version
$$\preskip.0em \postskip.2em 
\forall\, x,\,y,\;\big(x=y\;\hbox{ or }\;(x=y)'\big), 
$$
since if $0$ is \ivz, then $1=0$.
\label{footnoteNegation}}} the following implication

\snic{x\in\Ati\; \Rightarrow \; 1=_\gA0.
}

\smallskip
In the trivial \ri the \elt $0$ is both \iv and non-\iv at the same time.

For an arbitrary commutative \riz, the set of \elts  $a$ of $\gA$ which satisfy
\rdb
\begin{equation}\label{eqDefRadJac}
\forall x\in \gA~~~ 1+ax\in\Ati
\end{equation}
is called the \ixy{Jacobson}{radical} of $\gA$. It will be denoted by $\Rad(\gA)$.
It is an \id because if $a$, $b\in\Rad\gA$, we can write, for $x\in\gA$
\index{Jacobson radical!of a \riz}
$$\preskip.4em \postskip.4em 
1+(a+b)x=(1+ax)(1+(1+ax)^{-1}bx), 
$$
which is 
the product of two 
\iv \eltsz.

In a \alo the Jacobson radical is equal to the set of non-\iv \elts (the reader is invited to find a \prcoz). \index{Jacobson!radical}
In \clama the Jacobson radical is \care as follows.

\begin{lemmac}
\label{lemcRadJ}
The Jacobson radical is equal to the intersection of the \idemasz. 
\end{lemmac}
\begin{proof}
If $a\in\Rad\gA$ and $a\notin\fm$ with $\fm$ a \idemaz, we have $1\in\gen{a}+\fm$
which means that for some $x$, $1+xa\in\fm$, so $1\in\fm$: a contradiction.\\
If $a\notin\Rad\gA$, there exists an $x$ such that $1+xa$ is non-\ivz. Therefore there exists a strict \id containing $1+xa$. By Zorn's lemma there exists a \idemaz~$\fm$ containing $1+xa$, and $a$ cannot be in $\fm$ because otherwise we would have $1=(1+xa)-xa\in\fm$.\\
The reader will notice that the \dem actually says this: an \elt $x$ is in the intersection of the \idemas \ssi the following implication is satisfied:
$\gen{x,y}=\gen{1}
\Rightarrow  \gen{y}=\gen{1}$.
\end{proof}

\rem We have reasoned with a nontrivial \riz. If the \ri is trivial the intersection of the (empty) set of  \idemas is indeed equal to~$\gen{0}$.
\eoe

\medskip \rdb\label{corpsdeHeyting}
A \ixy{Heyting}{field}, or simply a \ix{field}, is by \dfn a \alo in which every non-\iv \elt is null, in other words a \alo whose Jacobson radical is reduced to~$0$.\index{Heyting!field}

In particular, a \cdiz, therefore also the trivial \riz, is a field.
The real numbers form a field that {\em is not} a \cdiz.{\footnote{We use the negation in
italics to indicate that the corresponding assertion, here it would be \gui{$\RR$ is a \cdiz,}
is not provable in \comaz.}} The same remark applies for the field $\QQ_p$ of $p$-adic numbers or that of the formal Laurent series $\gk(\!(T)\!)$ when $\gk$ is a \cdiz.

The reader will check that a field is a \cdi \ssi it is \zedz.

\smallskip
The quotient of a \alo by its Jacobson radical is a field, called the {\em residual field of the \aloz}.
\index{field!residual --- of a \aloz}

\begin{lemma}\label{lemZeDRaD}
If $\gA$ is \zedz, $\Rad\gA=\DA(0)$.
\end{lemma}
\begin{proof}
The inclusion $\Rad\gA\supseteq\DA(0)$ is always true. Now if $\gA$ is \zed and $x\in\Rad\gA$, since we have an \egt $x^\ell(1-ax)=0$, it is clear that $x^{\ell}=0$.
\end{proof}
%

\begin{lemma}\label{lemRadAX}
For all $\gA$, $\Rad(\AX)=\DA(0)[X]$.
\end{lemma}
\begin{proof}
If $f\in\Rad(\AX)$, then $1+Xf(X)\in\AX\eti$. We conclude with Lemma~\ref{lemGaussJoyal}~\emph{\iref{i4lemPrimitf}.} 
\end{proof}
%

\begin{fact}\label{fact1Rad} Let $\gA$ be a \ri and $\fa$ be an \id contained in $\Rad\gA$.
\begin{enumerate}
\item \label{i1fact1Rad}
$\Rad\gA=\pi_{\gA,\fa}^{-1}(\Rad\big(\gA\sur{\fa})\big)\supseteq\DA(\fa)$.
\item \label{i2fact1Rad}
$\gA$ is local \ssi $\gA\sur{\fa}$ is local.
\item \label{i3fact1Rad}
$\gA$ is local  and $\fa=\Rad\gA$ \ssi $\gA\sur{\fa}$ is a field.
\end{enumerate}
\end{fact}

The following fact describes a construction that forces a \mo to be inverted and an \id to be radicalized (for more details, see the subsection ``Duality in commutative rings'' on  \paref{secIDEFIL} and following, and Section~\ref{subsecMoco}).

\begin{fact}\label{fact2Rad}
Let $U$ be a \mo and $\fa$ be an \id of $\gA$.
Consider the \mo $S=U+\fa$. Let $\gB=S^{-1}\gA$ and $\fb=\fa\gB$.
\begin{enumerate}
\item The \id $\fb$ is contained in $\Rad\gB$.
\item The \ri $\gB\sur{\fb}$ is \isoc to $\gA_U\sur{\fa\gA_U}$.
\end{enumerate}
\end{fact}

 By \dfn a \ixy{\dcd local}{ring} is a \alo whose residual field is a \cdiz.
Such a \ri $\gA$ can be \care by the following axiom
\begin{equation}\label{eqDefAlrd}\preskip-.20em 
\forall x\in \gA \qquad x\in \Ati  \;\; {\rm  or} \;\;
    1+x\gA\,\subseteq\,  \Ati
\end{equation}
(the reader is invited to write its \prcoz).

For example the \ri of $p$-adic integers, although {\em non}-discrete, is \dcdz.

 We obtain a \emph{non}-\dcd \alo when taking $\gK[u]_{1+\gen{u}}$, where~$\gK$ is a \emph{non}-discrete field (for example the field of formal series~$\gk(\!(t)\!)$, where~$\gk$ is a \cdiz).

\medskip \comm
The slightly subtle difference that separates local rings from residually discrete local rings can also be found, by permuting addition and multiplication, in the difference that separates rings without zerodivisors from integral rings.
\\
In classical mathematics a ring without zerodivisors is integral; however the two notions do not have the same algorithmic content, and it is for this reason that they are distinguished in constructive mathematics.
\eoe

\vspace{-.3em}
\pagebreak

\begin{definition}\label{residzed}
\index{ring!residually \zed ---}%
\index{residually zero-dimensional!\riz} \hspace*{-1em}
 A \ri $\gA$ is said to be \emph{\plcz} when the \emph{residual ring} $\gA\sur{\Rad\gA}$ is \zedz. Likewise for \emph{\rdt connected rings}.
\end{definition}

Since a field is \zed \ssi it is a \cdiz, a \alo is \dcd \ssi it is \plcz.

\medskip
\comm In \clama a \ri $\gA$ is said to be semi-local if $\gA\sur{\Rad\gA}$ is \isoc to a finite product of \cdisz. This implies that it is a \plc \riz.
Actually the hypothesis of finiteness presented in the notion of a semi-local \ri is rarely decisive. Most of the \thos from the literature 
concerning semi-local rings applies to residually zero-dimensional rings, or even to local-global rings
(Section~\ref{secAlocglob}).
For a possible \dfn of semi-local \ris in \coma
see Exercises~\ref{exo1semilocal} and~\ref{exo2semilocal}.
\eoe

\subsec{Prime and maximal \idsz}
In \comaz, an \id of a \ri $\gA$ is called a \emph{maximal \idz} when the quotient \ri is a field.{\footnote{We have until now uniquely used the notion of a \idema in the context of proofs in \clamaz. 
A \cov \dfn would have been required sooner or later. Actually this notion is only rarely used in \comaz. As a \gnl rule, it is advantageously replaced by considering the Jacobson radical, for example in the \alos case.}} 
An \id is called a \emph{\idepz} when the quotient \ri is \sdzz.%
\index{ideal!maximal ---}%
\index{ideal!prime ---}%
\index{prime!\id of a commutative \riz}%
\index{maximal!ideal}

 These \dfns coincide with the usual \dfns in the context of \clamaz, except that we tolerate the trivial \ri as a field and hence the \id $\gen{1}$ as a \idema and as a \idepz.

In a nontrivial \riz, an \id is strict, maximal and detachable \ssi the quotient \ri is a nontrivial \cdiz, it is strict, prime and detachable \ssi the quotient \ri is a nontrivial integral \riz.

\medskip \comm \label{CommIdeps}
It is not without a certain apprehension that we declare the \id $\gen{1}$ both prime and maximal.
This will force us to say \gui{strict \idepz} or \gui{strict \idemaz} in order to speak of the \gui{usual} \ideps and \idemasz.
Fortunately it will be a very rare occurrence.
\perso{ce fut une d\'ecision difficile \`a prendre, fin 2006,
elle semble  convenablement motiv\'ee}

We actually think that there was \emph{a casting error right at the beginning}.
         To force a field or an integral \ri to be nontrivial, 
something that seemed eminently reasonable a priori, has unconsciously led mathematicians to transform numerous \cof arguments into reductio ad absurdum arguments. To prove that an \id constructed in the process of a computation is equal to $\gen{1}$, we have made it a habit to reason as follows: if it wasn't the case, it would be contained in a \idema and the quotient would be a field, which case we reach the contradiction $0=1$.
This argument happens to be a reductio ad absurdum simply because we have made the casting error: we have forbidden the trivial \ri from being a field. Without this prohibition, we would present the argument as a direct argument of the following form: let us show that every \idema of the quotient \ri contains $1$.
We will come back to this point in Section~\ref{subsecLGIdeMax}.

Moreover, as we will essentially use \ideps and \idemas heuristically, our transgression of the usual prohibition regarding the trivial \ri will have practically no consequence on reading this work. 
In addition, the reader will be able to see that this unusual convention does not force a modification of most of the results established specifically in \clamaz, like the \plgaz\etoz~\ref{plca.basic},
 Fact\eto \ref{factPropCarFin} or Lemma\eto \ref{lemcRadJ}:
it suffices for instance\footnote{Fact\eto \ref{factMoco} could \egmt be treated according to the same schema, by deleting the restriction to the nontrivial case.} for the \lon at a \idep $\fp$ to define it as the \lon at the filter 

\snic{S\eqdefi \sotq{x\in\gA}{x\in\fp \Rightarrow 1\in\fp}.}

Fundamentally we think that \maths is purer and more elegant when we avoid using negation (this radically forbids reductio ad absurdum arguments for example). It is for this reason that you will not find any \dfns that use negation in this book.\footnote{If such a definition could be found, it would be in a framework where the negation is equivalent to a positive assertion, because the considered \prt is decidable.}
\eoe

\subsect{The Jacobson radical and units in an integral extension}%
{The Jacobson radical and units in an integral extension}

\begin{theorem}\label{thJacUnitEntieres} Let $\gk\subseteq\gA$ with $\gA$ integral over $\gk$.
\begin{enumerate}
\item If $y\in\Ati$, then $y^{-1}\in\gk[y]$.
\item $\gk\eti=\gk\cap\Ati$.
\item $\Rad\gk=\gk\cap\Rad\gA$ and the \homo $\gA\to\gA\sur{\Rad(\gk)\gA}$ reflects the units.\footnote{Recall that we say that a \homo $\rho:\gA\to\gB$ reflects the units when~$\rho^{-1}(\gB\eti)=\Ati$.}
\end{enumerate}
\end{theorem}
\begin{proof}
\emph{1.} Let $y,z\in\gA$ such that $yz=1$. We have an \rdi for $z$: $z^n=a_{n-1}z^{n-1}+\cdots+a_0$ ($a_i\in\gk$).
By multiplying by $y^n$ we obtain $1=yQ(y)$ so $z=Q(y)\in\gk[y]$.

\emph{2.} In particular, if $y\in\gk$ is \iv in $\gA$, its inverse $z$ is in~$\gk$.

\emph{3.} Let $x\in \gk\cap\Rad\gA$, for all $y\in\gk$, $1+xy$ is invertible in $\gA$ therefore also in~$\gk$. This gives the inclusion $\Rad\gk\supseteq\gk\cap\Rad\gA$.\\
Let $x\in\Rad\gk$ and $b\in\gA$. We want to show that $y=-1+xb$ is \ivz. We write an \rdi for $b$
$$\preskip.4em \postskip.4em 
b^n+a_{n-1}b^{n-1}+\cdots+a_0=0, 
$$
we multiply by $x^n$ and replace $bx$ with $1+y$. We get a \pol in $y$ with \coes in $\gk$:
$y^n+\cdots+(1+a_{n-1}x+\cdots+a_0x^n)=0$. Therefore,~$yR(y)=1+xS(x)$ is \iv in $\gk$, and $y$ is \iv in $\gA$.\\
Now let $y\in\gA$ which is \iv modulo $\Rad(\gk)\gA$. A fortiori it is \iv modulo $\Rad\gA$, so it is \ivz.
\end{proof}

\begin{theorem}\label{thJacplc}
Let $\gk\subseteq\gA$ with $\gA$ integral over $\gk$.
\vspace{-3pt}
\begin{enumerate}
\item $\gA$ is \zed \ssi $\gk$ is \zedz.
\item $\gA$ is \plc \ssi $\gk$ is \plcz. In this case $\Rad\gA=\DA\big(\Rad(\gk)\gA\big)$.
\item If $\gA$ is local, so is $\gk$.
\end{enumerate}
\perso{point 2: Le lying over nous dit aussi que
$\Rad\gk=\gk\cap\Rad(\gk)\gA$. Aurait-on  $\Rad\gA=\sqrt{\Rad(\gk)\gA}$
d\`es que $\gA$ is enti\`ere sur $\gk$?}

\end{theorem}
\begin{proof} \emph{1.} Already known (Lemmas~\ref{lemZrZr1} and \ref{lemZrZr2}).

\emph{2.} By passage to the quotient, the integral morphism $\gk\to\gA$ gives an integral morphism $\gk\sur{\Rad\gk}  \to \gA\sur{\Rad\gA}$,  which is injective because $\Rad\gk=\gk\cap\Rad\gA$ (\thref{thJacUnitEntieres}). Therefore, the two \ris are simultaneously \zedsz.
In this case, let $\fa = \Rad(\gk)\,\gA \subseteq \Rad\gA$. We have an integral morphism 
$$\preskip-.4em \postskip.4em 
\gk\sur{\Rad\gk}  \to \gA\sur{\fa}, 
$$
so $\gA\sur{\fa}$ is \zedz, such that its Jacobson radical is equal to its nilpotent radical (Lemma~\ref{lemZeDRaD}), \cad $\Rad(\gA)\sur\fa = \DA(\fa)\sur\fa$, and so~$\Rad\gA = \DA(\fa).$

 \emph{3.} Results from \thref{thJacUnitEntieres}, item \emph{2.}
\end{proof}
%

\section{Four important lemmas}
\label{secAloc2}
First we give some variants of the \gui{\deter trick} often called \gui{Naka\-yama's lemma.}
In this lemma the important thing to underline is that the module $M$ is \tfz.
\index{Nakayama's lemma}
\index{determinant trick}

\CMnewtheorem{lemNak}{Nakayama's lemma}{\itshape}

\begin{lemNak}\label{lemNaka}  
\emph{(The \deter trick)}%
\index{Nakayama!lemma} \\
Let $M$ be a \tf \Amo and $\fa$ be an \id of $\gA$.
\begin{enumerate}
\item  If $\fa\,M=M$, there exists an $x\in\fa$ such that $(1-x)\,M=0$.
\item  If in addition $\fa\subseteq\Rad(\gA)$, then $M=0$.
\item  If $N\subseteq M$,  $\fa\,M+N=M$ and $\fa\subseteq\Rad(\gA)$, then $M=N$.
\item  If $\fa\subseteq\Rad(\gA)$ and $X\subseteq M$ generates $M/\fa M$ as an
$\gA/\fa$-module, then~$X$ generates $M$ as an \Amoz.
\end{enumerate}
\end{lemNak}

\begin{proof}
We prove item \emph{1}  and leave the others as an exercise, as easy consequences.
Let $V\in M^{n\times 1}$ be a column vector formed with \gtrs of $M$.
The hypothesis means that there exists a matrix $G\in\Mn(\fa)$ satisfying $GV=V$.
Therefore $(\I_n-G)V=0$, and by premultiplying by the cotransposed matrix of $\I_n-G$, we obtain $\det(\I_n-G)V=0$. However, $\det(\I_n-G)=1-x$ with $x\in\fa$.
\end{proof}

\Mptfs are locally free in the following (weak) sense: they become free when we localize at a \idepz. Proving this is the same as provinging the {\em local freeness lemma} (below) which states that a \mptf over a \alo is free.

\CMnewtheorem{lemlilo}{Local freeness lemma}{\itshape}
\begin{lemlilo}\label{lelilo}\index{Local freeness lemma}%
Let $\gA$ be a \aloz. Every \mptf over $\gA$ is free of finite rank.
Equivalently, every matrix $F\in\GAn(\gA)$ is similar to a standard \mprn
$$\I_{r,n}=\cmatrix{
\I_r&   &0_{r,n-r}  \cr
0_{n-r,r}&   &0_{n-r}
}.$$
\end{lemlilo}

\rem
The matrix formulation obviously implies the first, more abstract, formulation. 
Conversely if $M\oplus N=\Ae n$, saying that $M$ and $N$ are free (of ranks $r$ and $n-r$) is the same as saying that there is a basis of $\Ae n$ whose first $r$ \elts form a basis of $M$ and last $n-r$ a basis of $N$, consequently the \prn over $M$ parallel to $N$ is expressed over this basis by the matrix~$\I_{r,n}$.
\eoe

\begin{Proof}{First \demz, (usual classic \demz). }
We denote by $x\mapsto \ov{x}$ the passage to the residual field. If $M\subseteq \Ae n$ is the image of a \prn matrix~$F$ and if $\gk$ is the residual field we consider a basis of $\gk^n$ which begins with columns of $\ov{F}$ ($\Im \ov{F}$ is a \lin subspace of dimension $r$) and ends with columns of $\I_n-\ov{F}$ ($\Im  (\I_n-\ov{F})=\Ker \ov{F})$. When considering the corresponding columns of $\Im
{F}$ and  $\Im (\I_n- {F})=\Ker F$ we obtain a lift of the residual basis in $n$ vectors whose \deter is \rdt \ivz, therefore \ivz.
These vectors form a basis of~$\Ae n$ and over this basis it is clear that the \prn admits as a matrix~$\I_{r,n}$.
\\
Note that in this \dem we extract a maximal free \sys among the columns of a matrix with \coes in a field. This is usually done by the Gauss pivot method.
This therefore requires that the residual field be discrete.
\end{Proof}
\begin{Proof}{Second \demz, (\dem by Azumaya). }
 In contrast to the previous \demz, this one does not assume that the \alo is \dcdz.
We will diagonalize the matrix $F$.
The \dem works with a not \ncrt commutative \aloz.\\
Let us call $f_1$ the column vector $F_{1..n,1}$ of the matrix $F$,
$(e_1,\ldots ,e_n)$ the canonical basis of~$\Ae n$ and $\varphi$ the \ali represented by $F$.\\
-- First case, $f_{1,1}$ is invertible. Then,
$(f_1,e_2,\ldots ,e_n)$ is a basis of $\Ae n$. With respect to this basis, the \ali $\varphi$ has a matrix
$$
G=\cmatrix{
    1   &      L      \cr
    0_{n-1,1}&   F_1      }.
$$
By writing $G^2=G$, we obtain $F_1^2=F_1$ and $LF_1=0$.
We then define the matrix $P=\cmatrix{
    1   &      L      \cr
    0_{n-1,1}&   \I_{n-1}      }$ and we obtain the \egts
$$
\begin{array}{rcl}
PGP^{-1}& =  &
\cmatrix{
    1   &      L      \cr
    0_{n-1,1}&   \I_{n-1}      }
\cmatrix{
    1   &      L      \cr
    0_{n-1,1}&   F_1      }
\cmatrix{
    1   &      -L      \cr
    0_{n-1,1}&   \I_{n-1}      }   \\[4mm]
& =  &\cmatrix{
    1   &      0_{1,n-1}      \cr
    0_{n-1,1}&   F_1      }.
\end{array}
$$
-- Second case, $1-f_{1,1}$ is invertible.
We apply the previous computation to the matrix $\In-F$, which is therefore similar to a matrix
$$\preskip.4em \postskip.4em
A=\cmatrix{
    1   &      0_{1,n-1}      \cr
    0_{n-1,1}&   F_1      },
$$
with $F_1^2=F_1$, which means that $F$ is similar to a matrix
$$\preskip.4em \postskip.0em\In-A=\cmatrix{
    0   &      0_{1,n-1}      \cr
    0_{n-1,1}&   H_1      },
$$
with $H_1^2=H_1$.\\
We finish the proof by \recu on $n$.
\end{Proof}

\comm \label{comment lelilo} 
From the classical point of view, all the sets are discrete, and the corresponding hypothesis is superfluous in the first \demz.
The second \dem must be considered superior to the first as its \algq content is more \uvl than that  of the first (which can only be rendered completely explicit when the \alo is \dcdz). \eoe


\medskip \rdb
The following lemma can be considered as a variant     
of the local freeness lemma.

\CMnewtheorem{lemnllo}{Lemma of the \lnl map}{\itshape}
\begin{lemnllo}\label{lelnllo}%
\index{Lemma of the locally simple map}%
Let $\gA$ be a \alo and~$\psi $ be a \ali between free \Amos of finite rank.
\Propeq 
\begin{enumerate}
\item $\psi$ is \nlz.
\item $\psi$ is \lnlz.
\item $\psi$ has a finite rank $k$.
\end{enumerate}
\end{lemnllo}

\begin{proof}
\emph{2}  $\Rightarrow$ \emph{3.} The \egt $\psi\,\varphi\,\psi=\psi$ implies that the \idds of $\psi$ are \idmsz.
By Lemma~\ref{lem.ide.idem} these \ids are generated by \idmsz.
Since an \idm of a \alo is \ncrt equal to $0$ or $1$, and that $\cD_0(\psi
)=\gen{1}$ and $\cD_{r}(\psi )=\gen{0}$ for large enough $r$, there exists an integer $k\geq 0$ such that
$\cD_k(\psi )=\gen{1}$ and $\cD_{k+1}(\psi )=\gen{0}$.

\emph{3}  $\Rightarrow$ \emph{1.} By hypothesis $\cD_k(\psi )=\gen{1}$, so the minors of order $k$ are \com and since the \ri is local one of the minors of order $k$ is \ivz.
As $\cD_{k+1}(\psi )=\gen{0}$, the result is then a consequence of the freeness lemma \ref{lem pf libre}.
\end{proof}

Note that the term \lnl map is partly justified by the previous lemma. Also note that \thref{theoremIFD} can be considered as more \gnl than the previous lemma.

\CMnewtheorem{lemnbgl}{Local number of \gtrs lemma}{\itshape}
\begin{lemnbgl} 
\label{lemnbgtrlo}\index{Local number of \gtrs lemma}~\\
Let $M$ be a \tf \Amoz.
\begin{enumerate}
\item Suppose $\gA$ is local. 
\begin{enumerate}
\item [a.\phantom{*}] The module $M$ is generated by $k$ \elts \ssi its \idf $\cF_k(M)$ is equal to $\gA$.
\item [b.\phantom{*}] If in addition $\gA$ is \dcd and $M$ is \pfz, the module admits a \pn matrix whose every \coe is in the maximal \id $\Rad\gA$.
\end{enumerate}
\item \Gnltz, for any $k\in\NN$ \propeq
\begin{enumerate}
\item [a.\phantom{*}] $\cF_k(M)$ is equal to $\gA$.
\item [b.\phantom{*}] There exist \eco $s_j$ such that after \eds to each of the $\gA[1/s_j]$, $M$ is generated by $k$ \eltsz.
\item [c.\phantom{*}] There exist \moco $S_j$ such that each of the $M_{S_j}$ is generated by $k$ \eltsz.
\item [d*.] After \lon at any \idepz,  $M$ is generated by $k$ \eltsz.
\item [e*.] After \lon at any \idemaz,  $M$ is generated by $k$ \eltsz.
\end{enumerate}
\end{enumerate} 
\end{lemnbgl}
\begin{proof} It suffices to prove the \eqvcs for a \mpf
 due to Fact~\ref{facttfpf}.
\\
 Suppose $M$ is generated by $q$ \elts and let $k'=q-k$.
\\
\emph{1.}
The condition is always \ncrz, even if the \ri is not local.
Let $A\in\gA^{q\times m}$ be a \pn matrix for $M$.
If the \ri is local and if  $\cF_k(M)=\gA$, since the minors of order $k'$ are \comz, one of them is \ivz. 
By the \iv minor lemma~\ref{lem.min.inv}, the matrix $A$ is \eqve to a matrix
$$\preskip-.2em \postskip.4em 
 \cmatrix{
   \I_{k'}   &0_{k',m-k'}      \cr
    0_{k,k'}&     A_1},
$$
and so, the matrix $A_1\in\Ae {k\times (m-k')}$ is also a \mpn of~$M$.
Finally, if the \ri is \dcdz, we can reduce the number of \gtrs until the corresponding \pn matrix has all of its \coes in the radical.

 \emph{2.}
\emph{a}  $\Rightarrow$ \emph{b}.
The same \dem shows that we can take, for $s_j$, the minors of order $k'$ of~$A$.
\\
\emph{b}  $\Rightarrow$ \emph{c}.
Immediate.
\\
\emph{c}  $\Rightarrow$ \emph{a}.
Saying that $\cF_k(M)=\gA$ comes down to solving the \sli $\sum_\ell x_\ell s_\ell =1$, where the unknowns are the $x_\ell$'s and where the $s_\ell$'s are the minors of order $k'$ of the matrix  $A$. We can therefore apply the basic \lgb principle.
\\
\emph{a}  $\Rightarrow$ \emph{d}.
Results from  \emph{1}. 
\\
\emph{d}   $\Rightarrow$ \emph{e}.
Trivial.
\\
\emph{e}  $\Rightarrow$ \emph{a}.
This can only be proven in \clama (hence the star that we attached to \emph{d}  and \emph{e}).
We prove it by reductio ad absurdum, by proving the contrapositive. If $\cF_k(M) \neq \gA$ let $\fp$ be a strict \idema containing~$\cF_k(M)$. After \lon at $\fp,$ we obtain $\cF_k(M_\fp)\subseteq \fp\gA_\fp \neq \gA_\fp$, and so~$M_\fp$ is not generated by $k$ \eltsz.
\end{proof}

\comm  This lemma gives the \emph{true meaning} of the \egtz~$\cF_k(M)=\gA$; we can say that $\cF_k(M)$ \gui{measures} the possibility for the module to be locally generated by $k$ \eltsz. Hence the following \dfnz.
\\
See also Exercises~\ref{exoAutresIdF}, \ref{exoNbgenloc} and~\ref{exoVariationLocGenerated}.
\eoe

\begin{definition}
\label{deflocgenk}
A \mtf is said to be \ixd{locally generated by~$k$ \eltsz}{module} when it satisfies the \eqv \prts of Item \emph{2} in the local number of \gtrs lemma. \index{locally!module --- generated by $k$ \eltsz}
\end{definition}


\section{\Lon at \texorpdfstring{$1+\fa$}{1 + a}}
\label{secLoc1+fa} 

\medskip 
\Grandcadre{Let $\fa$ be an \id of $\gA$, $S := 1 + \fa$,  $\jmath:
\gA \to \gB:=\gA_{1+\fa}$ \\ be the canonical \homoz, and
$\fb := \jmath(\fa)\gB$.}

\medskip 
Note that $\fb$ is identified with $S^{-1}\fa$
(Fact \ref{fact.sexloc}) and that $1+\fb\subseteq\gB\eti$ 
(Fact~\ref{fact2Rad}).
\begin{lemma}\label{lemLoc1+a} 
 {(Quotient of powers of $\fa$ in the localized \ri  
$\gA_{1+\fa}$)}
\\
Under the previous hypotheses we have the following results. 
\begin{enumerate}
\item $\Ker\jmath\subseteq\fa$, $\gB = \jmath(\gA) + \fb$ and the canonical \homo $\gA\sur\fa\to\gB\sur\fb$ is an \isoz. 
\item The \lon at $1+\fa$ is the same as the \lon at $1+\fa^n$ ($n\geq1$), so $\Ker\jmath\subseteq\fa^n$, $\gB = \jmath(\gA) + \fb^n$ and $\gA\sur{\fa^n}\simeq\gB\sur{\fb^n}$.
\item  For all $p$, $q \in \NN$, $\jmath$ induces an \iso $\fa^p\sur{\fa^{p+q}} \simarrow \fb^p\sur{\fb^{p+q}}$ \hbox{of \Amosz}.
\end{enumerate}
\end{lemma}
%
\begin{proof}
\emph{1.} The inclusion $\Ker\jmath\subseteq\fa$ is \imdez.
\\
The fact that the \homo $\gA\sur\fa\to\gB\sur\fb$ is an \iso relies on 
two \eqv \uvl \pbs being solved: in the first we must annihilate the \elts of 
$\fa$, in the second, we also need to invert the \elts of $1+\fa$, but inverting~$1$ is costless. Finally, the surjectivity of this morphism means precisely that $\gB = \jmath(\gA) + \fb$.

\emph{2.} The \mos  $1+\fa$ and $1+\fa^n$ are \eqvs  because $1-a$ divides $1-a^n$.

\emph{3.} Let $\fb^q=S^{-1}\fa^q=\fa^q\gB$.
By multiplying $\gB = \jmath(\gA) + \fb^q$ by~$\fa^p$, we obtain $\fb^p = \jmath(\fa^p) + \fb^{p+q}$.
Therefore, the map $\jmath$ induces a surjection \hbox{of \Amosz}
$\fa^p \twoheadrightarrow \fb^p\sur{\fb^{p+q}}$. 
It remains to see that its kernel is~$\fa^{p+q}$.
\hbox{If $x \in \fa^p$} satisfies $\jmath(x) \in \fb^{p+q}$, 
there exists an~$s \in 1+\fa$ such that $sx \in \fa^{p+q}$,
and since~$s$ is \iv modulo $\fa$, it is also \iv modulo $\fa^{p+q}$,
and so $x \in \fa^{p+q}$.
\end{proof}
%

\CMnewtheorem{lemlofi}{Localized finite \ri lemma}{\itshape}
\begin{lemlofi}\label{lemLocaliseFini}%
If $\fa$ is a \itf and $n\in\NN\etl$, we have the \eqvcs%
\index{Localized finite ring lemma} 
$$\preskip.3em \postskip-.1em 
\fb^n = \fb^{n+1}
\iff \fb^n = 0 \iff \fa^n = \fa^{n+1}. 
$$
In this case, 
\begin{enumerate}
\item we have $\;\fa^n = \Ker \jmath=\gen{1-e}$ with $e$ \idmz, such that
$$\preskip.3em \postskip.3em 
\gB=\gA_{1 + \fa} = \gA[1/e]=\aqo{\gA}{1-e}, 
$$
\item
if in addition $\gA$ is a \klgz,
then $\gA\sur\fa$ is finite over $\gk$ \ssi $\gB$ is finite over~$\gk$.
\end{enumerate} 
\end{lemlofi}


%
\begin{proof}
If $\fb^n = \fb^{n+1}$, then $\fb^n$ is a \tf \idm ideal, 
so $\fb^n=\gen{\vep}$ 
with~$\vep$ being an \idmz.  
But since $\varepsilon \in \fb$, the \idm $1 - \varepsilon$
is \ivz, therefore equal to $1$, \cad $\varepsilon = 0$, so $\fb^n=0$.
 The third equivalence comes from $\fb^n\sur{\fb^{n+1}}\simeq\fa^n\sur{\fa^{n+1}}$
(Lemma~\ref{lemLoc1+a}).

 \emph {1.} 
Since $\fa^n$ is a \tf \idmz, $\fa^n=\gen{1-e}$ with $e$ an \idmz. The rest then stems from Fact~\ref{fact.loc.idm}.

\emph {2.} 
If $\gB$ is a \tf \kmoz, so is $\gA\sur\fa \simeq
\gB\sur\fb$. Conversely, suppose that $\gA\sur\fa$ is a \tf \kmo and let us consider the filtration of $\gB$ by the powers of $\fb$
$$\preskip.4em \postskip.4em 
0 = \fb^n \subseteq \fb^{n-1} \subseteq\cdots\subseteq \fb^2 
\subseteq \fb \subseteq \gB. 
$$
Then, each quotient $\fb^i\sur{\fb^{i+1}}$ is a $\gB\sur\fb$-\mtfz,
or an $\gA\sur\fa$-\mtfz, and consequently a \tf \kmoz. We deduce that $\gB$ is a \tf \kmoz.
\end{proof}
%

\CMnewtheorem{lemlozed}{Localized \zed \ri lemma}{\itshape} 
\begin{lemlozed} \index{Localized zero-dimensional ring lemma} 
\label{lemLocalisezeddim}~\\
Let $\fa$ be a \itf of $\gA$ such that the localized \ri $\gB=\gA_{1+\fa}$ is \zedz.  Then, there exist an integer $n$ and an \idm $e$ such that 
$$\preskip-.2em \postskip.3em 
\;\;\fa^n =\gen{1-e}\quad\hbox{  and  }\quad\gA_{1 + \fa} = \gA\big[\frac 1{e}\big]=\aqo{\gA}{1-e}. 
$$
If in addition $\gA$ is a \tf \klg with $\gk$ \zed (for example a \cdiz), then $\gB$ is finite over $\gk$. 
\end{lemlozed}
%
\begin{proof} We apply the localized finite \ri lemma; since $\gB$ is \zed and $\fb$ \tfz, there exists an integer $n$ such that $\fb^n = \fb^{n+1}$.  
\\
We end with the weak \nst \ref{thNst0} {because $\gB=\aqo\gA{1-e}$} is  \hbox{a \tfz} \klgz.   
\end{proof}

\rem
Let $\fa$ be a \itf of a \ri $\gA$ such that the localized \riz~$\gA_{1+\fa}$ is
\zedz.  
The natural map $\gA \to \gA_{1 + \fa}$ is therefore surjective with kernel $\bigcap_{k \ge 0} \fa^k=\fa^m$ with $m$ such that $\fa^m = \fa^{m+1}$. In addition, $\fa^m$ is generated by an \idm $1-e$ and $\gA_{1+\fa}=\gA[1/e]$. We then have
$$\preskip-.4em \postskip.4em \ndsp
\bigcap_{k \ge 0} \fa^k = \big(0 : (0 : \fa^\infty)\big). 
$$
This remark can be useful for computations. Suppose that $\gA = \kuX\sur\ff$ where $\kuX = \gk[\Xn]$ is a \pol \ri with $n$ \idtrs over a \cdi $\gk$ and $\ff = \gen {f_1, \ldots, f_s}$ is a \itfz.  Let $\fa$ be a \itf of $\kuX$ and $\ov\fa$ be its image in $\gA$. Then, if $\gA_{1 + \ov\fa}$ is \zedz, the composition $\kuX \to \gA_{1 + \ov\fa}$ is surjective and its kernel is expressed in two ways

\snic {
\bigcap_{k \ge 0} (\ff + \fa^k) = \big(\ff : (\ff : \fa^\infty)\big).
}

The right-hand side formula can turn out to be more efficient by computing $(\ff : \fa^\infty)$ as follows
$$\preskip-.0em \postskip.0em \ndsp
 \qquad
(\ff : \fa^\infty) = \bigcap_{j=1}^r (\ff : g_j^\infty)
\hbox { if } \fa = \gen {g_1, \ldots, g_r}.
$$
\vspace{-1em}
\eoe

\vspace{2em}
 \comm
In \clama a \idep $\fa$ of the \riz~$\gA$ is said to be \emph{isolated} if it is both minimal and maximal in the set of \ideps of~$\gA$. 
In other words if it does not compare to any other \idep for the inclusion relation. Saying that $\fa$ is maximal amounts to saying that~$\gA/\fa$ is \zedz. Saying that $\fa$ is minimal amounts to saying that~$\gA_S$ is \zedz, where $S=\gA\setminus \fa$. But if $\fa$ is assumed to be maximal, the \moz~$S$ is the 
saturated \mo 
of $1+\fa$.\\
Conversely suppose that $\gB=\gA_{1+\fa}$ is \zedz. Then $\gA/\fa$ is \egmt \zed since $\gA/\fa\simeq\gB/\fa\gB$.
Thus, when~$\fa$ is also \tfz, we find ourselves with a special case of the localized \zed \ri lemma \ref{lemLocalisezeddim}.  
It is worth noting that in the literature the isolated \ideps \gnlt intervene in the context of \noes \ris and that therefore in \clama they are automatically \tfz.
\eoe 


\newpage
\section{Examples of \alos in algebraic geometry}
\label{secExlocGeoAlg}

Here we propose to study in a few cases \gui{the local \alg in a zero of a \sypz.}
We fix the following context for all of Section~\ref{secExlocGeoAlg}.

\Grandcadre{$\gk$ is a \riz, $\uf = f_1,\ldots,f_s\in\kuX= \kXn$,\\[1mm]
$\gA=\aqo{\kuX}{\uf}=\kxn$,  \\[1mm]
$(\uxi)=(\xin)\in\gk^n$  is a zero of the \sysz,\\[1mm]
$\fm_\uxi= \gen{x_1-\xi_1,\ldots,x_n-\xi_n}_\gA$ is
the \id of the point $\uxi$,\\[1mm]
$J(\uX)=\JJ_\uX(\uf)$ is the Jacobian matrix of the \sysz.} 

Recall that $\gA=\gk\oplus\fm_\uxi$ (Proposition~\ref{prdfCaracAlg}).
More \prmtz, we have with the \evn at $\uxi$ a split exact sequence of \kmos
$$\preskip.0em \postskip.3em 
0\to\fm_\uxi\to\gA\vvvvvers{\ov g\,\mapsto\, g(\uxi)}\gk\to0, 
$$
and two \homos of \klgs 
$
\gk\to\gA\to\gk 
$ 
which when composed give $\Id_\gk$.

Also recall (\thref{thidptva}) that $\fm_\uxi$ is a \pf \Amo
(the \mpn is explicitly given).

\subsec{Local \alg at a zero}
\label{subsecAlgLocZer}

In the following \dfn the terminology \emph{local \alg at $\uxi$}
must not be ambiguous. We do not claim that it is a \aloz, we simply mimic the construction of the given local \alg in the case where~$\gk$ is a field.

\begin{definition}\label{notaAlgLocEnxi}
\emph{(Local \alg at a zero of a \sypz)}
\\
The \ri $\gA_{1+\fm_\uxi}$  is called \emph{the local \alg at $\uxi$ of the \syp  $\uf$}. We also use the shorthand notation $\gA_\uxi$ instead of $\gA_{1+\fm_\uxi}$.%
\index{local algebra!in a zero of a \sypz}
\end{definition}

We denote by $\xi:\gA\to\gk$ the  \evn at $\uxi$. It is factorized through the \lon \hbox{at $1+\fm_\uxi$} 
and we obtain a \crc $\gA_{\uxi}\to\gk$. We therefore have $\gA_{\uxi}=\gk\oplus\fm_\uxi\gA_{\uxi}$ and canonical \isos 
$$\preskip.2em \postskip.2em 
\gA_{\uxi}\big/{(\fm_\uxi\gA_{\uxi})}\simeq \gA\big/{\fm_\uxi} \simeq \gk 
$$

\begin{fact}\label{factAlgLocEnxi}
 \emph{(If $\gk$ is a \cdiz, the \alg $\gA_{\uxi}$ is a \aloz)}
\begin{enumerate}
\item Let $\gk$ be a \alo with $\Rad\gk=\fp$, 
 $\fM=\fp\gA+\fm_\uxi$ and $\gC=\gA_{1+\fM}$.
Then, $\gC$ is a \alo with $\Rad(\gC)=\fM \gC$
and $\gC\sur{\Rad\gC}\simeq \gk\sur{\fp}$.
 
\item If $\gk$ is a \cdiz, we have the following results.
\begin{enumerate}
\item The \ri $\gA_{\uxi}$ is a \alo with $\Rad\gA_{\uxi}=\fm_\uxi\gA_{\uxi}$  and its residual field is (canonically \isoc to)~$\gk$.
\item The \ris $\gA$ and $\gA_{\uxi}$ are \coh \noesz, and $\gA$ is \fdiz.
\item $\bigcap_{r\in\NN}\bigl(\fm_\uxi\gA_{\uxi}\bigr)^r=0.$
\end{enumerate}
\end{enumerate}
\end{fact}
\begin{proof}
\emph{1.} We have $\gC\sur{\fM\gC} \simeq \gA\sur{\fm_\xi} = \gk/\fp$ by item \emph{2}  of Fact~\ref{fact2Rad}, then use item~\emph{3}  of Fact~\ref{fact1Rad}.
 
\emph{2a.} Results from \emph{1.}
 
\emph{2b.} The \ri $\gA$ is \fdi and \coh by \thref{thpolcohfd}.
We deduce that $\gA_{\uxi}$ is \cohz.
\\
For the \noet we refer the reader to \cite[VIII.1.5]{MRR}.
 
\emph{2c.} Given items \emph{2a}  and \emph{2b}, this is a special case of Krull's intersection \tho (\cite[VIII.2.8]{MRR}).
\end{proof}
%

\subsubsection*{Tangent space at a zero}%
\index{tangent space}\index{space!tangent ---}

In what follows we write $\partial_jf$ for $\Dpp{f}{X_j}$. Thus the Jacobian matrix of the \sysz, which we have denoted by $J=J(\uX)$, is visualized as follows
$$
\bordercmatrix [\lbrack\rbrack]{
    & X_1                     & X_2                     &\cdots  & X_n \cr
f_1 & \partial_1 {f_1} &\partial_2 {f_1}  &\cdots  &\partial_n {f_1} \cr
f_2 & \partial_1 {f_2} &\partial_2 {f_2}  &\cdots  &\partial_n {f_2} \cr
\;\vdots & \vdots                  &                         &        & \vdots              \cr
f_i & \vdots                  &                         &        & \vdots              \cr
\;\vdots & \vdots                  &                         &        & \vdots              \cr
f_s & \partial_1 {f_s} &\partial_2 {f_s}  &\cdots  &\partial_n {f_s} \cr
}~=~J.
$$
The congruence below is \imdez,
for $f \in \gk[\uX]$,
\begin{equation}\label {eqModM2}
f(\uX) \equiv f(\uxi) +
\som_{j=1}^n (X_j - \xi_j)\, \partial_j f(\uxi)
\;\mod \gen {X_1-\xi_1, \ldots, X_n-\xi_n}^2
\end{equation}
By specializing $\uX$ in $\ux$ we obtain in $\gA$ the fundamental congruence
\begin{equation}\label {eq2ModM2}
f(\ux) \equiv f(\uxi)+ 
\som_{j=1}^n  \,(x_j - \xi_j) \,\partial_j f(\uxi)
 \,\;\mod {\fm_\uxi}^2
\end{equation}
We leave it up to the reader to verify that the kernel of
$J(\uxi)$  only depends on the \id $\gen{\lfs}$ and on the point $\uxi$. It is a \ksmo of $\gk^n$ which can be called \emph{the tangent space at $\uxi$ to the affine scheme  
over $\gk$ defined by~$\gA$}. We will denote it by $\rT_\uxi(\gA\sur\gk)$ or $\rT_\uxi$.

This terminology is reasonable in \gaq (\cad when $\gk$ is a \cdiz), 
at least in the case where $\gA$ is integral. In that case we have a \vrt defined as an intersection of hypersurfaces $f_i=0$, and the tangent space at $\uxi$ of the \vrt 
is the intersection of the tangent spaces at the hypersurfaces that define it.

In this same situation (\cdi as the basis), the zero $\uxi$ of the \syp is called a \emph{regular point} or a \emph{non-singular point} (of the affine scheme or yet again of the corresponding \vrtz) when the dimension of  the tangent space at $\uxi$ is equal to the dimension\footnote{If $\gA$ is integral, this dimension
does not depend on $\uxi$ and can be defined via a Noether position. In the \gnl case, the \ddk of the \riz~$\gA_{\uxi}$ must be considered.}
of the \vrt at the point $\uxi$.
A point that is not regular is called \emph{singular}.

\smallskip We now give a more abstract interpretation of the tangent space, in terms of \dvn spaces. This works with an arbitrary commutative \ri $\gk$.

For a \klg $\gB$ and a \crc $\xi:\gB\to\gk$ we define a \emph{$\gk$-\dvn at the point $\xi$ of $\gB$} as a $\gk$-\lin form $d:\gB\to\gk$ which satisfies Leibniz's rule, \cad by letting $f(\xi)$ for $\xi(f)$%
\index{derivation!at a point (a character) of an \algz}
$$
d(fg)=f(\xi)d(g)+g(\xi)d(f).
$$
This implies in particular $d(1)=0$ (writing $1=1\times 1$), and so $d(\alpha)=0$ for $\alpha\in\gk$.
We will denote by $\DBxk$ the \kmo of the $\gk$-\dvns of~$\gB$ at the point $\xi$.
\\
This notation is slightly abusive. 
Actually if we let $\gk'$ be the \ri $\gk$ provided with the \Bmo structure  given by $\xi$, the notation of \Dfnz~\ref{defiDeriv} would be $\Der \gk\gB{\gk'}$, and in fact equipped with the structure \hbox{of a \Bmoz}.

\smallskip  
We will see that the tangent space at $\uxi$ of $\gA$ and the \kmo of the $\gk$-\dvns of~$\gA$ at~$\xi$ are naturally \isocz. 

\begin{proposition}\label{propTangent}
\emph{ ($\rT_\uxi(\gA\sur\gk)$, $\DAxk$,
and $(\fm_\uxi/{\fm_\uxi}^2)\sta$)}\\
Let $\fm=\fm_{\uxi}$ and recall the notation $\rT_\uxi(\gA\sur\gk)=\Ker J(\uxi)$.
\begin{enumerate}
\item For $u=(\un)\in\gk^n$, let $D_u : \gk[\uX]\to\gk$ be the $\gk$-\lin form defined by
$$\preskip-.2em \postskip.2em\ndsp 
D_u(f) = \sum_{j=1}^n {\partial_j f}(\uxi)\; u_j. 
$$
It is a \dvn at the point $\uxi$, we have $u_j = D_u(X_j) = D_u(X_j - \xi_j)$,
and the map 
$$\preskip-.2em \postskip.2em 
 u \mapsto D_u, \; \gk^n\to \DkXxk
$$
 is a $\gk$-\lin \isoz.

\item If $u \in \Ker J(\uxi) \subseteq\gk^n$, then $D_u$ passes to the quotient modulo $\gen {f_1, \ldots, f_s}$ and provides a $\gk$-\dvn at the point $\uxi$, $\Delta_u : \gA\to\gk$. \\
We have 
$u_j = \Delta_u(x_j) = \Delta_u(x_j - \xi_j)$, and the map  
$$\preskip.2em \postskip.2em
u \mapsto \Delta_u, \; \Ker J(\uxi)\to \DAxk
$$
 is a $\gk$-\lin \isoz.

\item In addition, $\Delta_u(\fm^2) = 0$ and we obtain, by restriction to $\fm$ and passage to the quotient modulo $\fm^2$, a $\gk$-\lin form $\delta_u:\fm\sur{\fm^2} \to \gk$. 
We thus construct a \kli $u \mapsto \delta_u$ of $\Ker J(\uxi)$ in $(\fm\sur{\fm^2})\sta$.

\item Conversely, to $\delta \in (\fm\sur{\fm^2})\sta$, we associate $u \in \gk^n$ defined by 

\snic {u_j = \delta\big((x_j-\xi_j) \bmod\fm^2\big).}

Then, $u$ belongs to $\Ker J(\uxi)$.

\item The two maps defined in 3 and 4, 

\snic{\Ker J(\uxi) \to (\fm\sur{\fm^2})\sta\quad$ and 
$\quad(\fm\sur{\fm^2})\sta \to \Ker J(\uxi),}

are reciprocal $\gk$-\lin \isosz.

\end{enumerate}
\end{proposition}

\begin {proof}
\emph {1.}
Simple \vfn left to the reader.

\emph {2.}
For any $u\in \gk^n$, we easily verify that the set

\snic {
\sotq {f \in \gk[\uX]}{D_u(f) = 0 \hbox { and } f(\uxi) = 0}
}

is an \id of $\gk[\uX]$. If $u\in\Ker J(\uxi)$, we have $D_u(f_i) = 0$
by \dfn \hbox{(and $f_i(\uxi) = 0$)}; we deduce that $D_u$ is null over $\gen {f_1, \ldots, f_s}$.

\emph {3.}
To see that $\Delta_u(\fm^2) = 0$, we use $\Delta_u(fg) =
f(\uxi)\Delta_u(g) + g(\uxi)\Delta_u(f)$ and $f(\uxi) = g(\uxi) =
0$ for $f$, $g \in \fm$.

\emph {4.}
The congruence (\ref{eq2ModM2}) for $f = f_i$ is $\sum_{j=1}^n (x_j - \xi_j)
{\partial_j f_i}(\uxi) \in \fm^2$. Applying $\delta$, this gives   the \egt $\sum_{j=1}^n u_j {\partial_j f_i}(\uxi) = 0$, \cad $u \in \Ker J(\uxi)$.

\emph {5.}
Let $\delta \in (\fm\sur{\fm^2})\sta$ and $u \in \Ker J(\uxi)$ be the corresponding \eltz; it must be shown that $\delta_u = \delta$, which is the same as checking, for $f \in \fm$,

\snic {
\delta(f \bmod\fm^2) = \sum_{j=1}^n {\partial_j f}(\uxi)
\delta\big((x_j-\xi_j)\bmod\fm^2\big)
,}

but this stems from $(\ref{eq2ModM2})$.

Conversely, let $u \in \Ker J(\uxi)$ and $v \in \Ker J(\uxi)$ be the \elt
corresponding to $\delta_u$; it must be shown that $v = u$; which is the same as checking $\delta_u\big((x_j - \xi_j)\bmod\fm^2\big) = u_j$, an \egt which has already been observed.
\end{proof}

\rem Note that the \dfn which we have given for the tangent space $\rT_\uxi(\gA\sur\gk)$, natural and intuitive, portrays it as a submodule of $\gk^n$, where $n$ is the number of \gtrs of the \pf \klg $\gA$. Therefore, its more abstract \dfnz~$\Der\gk\gA\xi$, or~$\fm_\uxi/{\fm_\uxi}^2$,
which is more intrinsic, must be preferred since it only depends on the \klgz~$\gA$ and on the \crc $\xi:\gA\to\gk$, without taking into account the \pn chosen for~$\gA$ (actually only the structure of the localized \algz~$\gA_\uxi$ intervenes).  
\eoe

\subsubsection*{Cotangent space at a zero}

\Gnltz, we also have the dual notion of a \emph{cotangent space at $\uxi$}.
We will define it here 
as the cokernel of the transposed matrix $\tra{\,J(\uxi)}$. 
Actually, it is a \kmo which is intrinsically attached to the \alg $\gA$ and to the \crcz~$\xi$,
because it can also be defined formally as \gui{the space of  \diles at the point $\uxi$.} We will not be developing this notion here. 

The fundamental \tho that follows implies that the tangent space is canonically \isoc to the dual of the cotangent space (Fact~\ref{factDualReflexif}~\emph{2}  applied to $\tra J$ gives $(\Coker\tra J)\sta\simeq \Ker J$ since $\tra{(\tra J)}=J$). 
However, when we work with an arbitrary \ri $\gk$, the cotangent space is not \ncrt \isoc to the dual of the tangent space.

When a \Bmo $M$ admits a \mpn $W$ over a \sgr $(\yn)$, if $\fb$ is an \id of $\gB$, by the base \ri change $\pi_{\gB,\fb}:\gB\to\gB\sur\fb$, we obtain the~$\gB\sur\fb$-module $M\sur{\fb M}$ with the \mpn $W\mod\fb$ over the \sgr $(\ov{y_1},\ldots,\ov{y_n})$. 
\\
With the \Amo $M=\fm_\uxi$ and the \id $\fb=\fm_\uxi$, we obtain for the \mpn of the \kmo $\fm_\uxi/{\fm_\uxi}^2$ over 
$(\ov{x_1-\xi_1},\ldots,\ov{x_n-\xi_n})$, the matrix~$\ov W=W\mod\fm_\uxi$, with the matrix $W$ given in \thref{thidptva}. The latter matrix, up to null columns, is the matrix $\tra{J(\uxi)}$. The \tho that follows states the same thing in a precise manner.

\begin{theorem}\label{thCotangent}
 \emph{(Cotangent space at $\uxi$ and $\fm_\uxi \big/ {{\fm_\uxi}^2}$)}
Let  $(e_i)_{i\in\lrbn}$ be the canonical basis of $\gk^n$. Consider the \kli

\snic{\varphi : \gk^n \twoheadrightarrow {\fm_\uxi}/{{\fm_\uxi}^2},\quad
e_j \mapsto (x_j - \xi_j) \bmod {\fm_\uxi}^2.}

Then, $\varphi$ induces an \iso of \kmos $\;\Coker\tra{J(\uxi)} \simarrow 
{\fm_\uxi}/{{\fm_\uxi}^2}$.
\\
Thus, we have a canonical \iso $\Coker\tra{J(\uxi)} \simarrow 
{\fm_\uxi}\gA_\uxi/{({{\fm_\uxi}\gA_\uxi})^2}$. 
\end{theorem}

\begin{proof} Suppose \spdg that $\uxi=\uze$ and use the notations of \thref{thidptva}. The \mpn of $\fm_\uze$ for the \sgr $(\xn)$ is the matrix $W = [\,R_\ux\,|\,U\,]$ with $U(\uze)=\tra{J(\uze)}$. As the matrix $R_\ux\mod\fm_\uze$ is null, we obtain the stated result.
\\
 The last assertion is given by Lemma~\ref{lemLoc1+a}~\emph{3}. 
\end{proof}

\begin{definition}\label{defiEspCot}
We define \emph{the cotangent space at $\uxi$} as being the \kmo $\fm_\uxi\gA_\uxi/{({\fm_\uxi\gA_\uxi})^2}$, for which only the structure of the local \alg at $\uxi$ intervenes.
 \end{definition}

 In the remainder of Section~\ref{secExlocGeoAlg},
 we will study a few examples of local \algs at zeros of \sypsz, without assuming that we \ncrt have a \cdi to begin with; $\gk$ is only a commutative \riz.
Here we only seek to illustrate the \gmq  situation by freeing ourselves, if possible, of the hypothesis \gui{\cdiz,} but without aiming to give the most \gnl framework possible.

\subsec{Local \ri at an isolated point}
\label{subsecZedLocPtIsole}

The idea that drives this subsection comes from \gaq where the \alo at $\uxi$ is \zed \ssi the point $\uxi$ is an isolated zero, and where the isolated zero is simple \ssi the tangent space is reduced to $0$.

\pagebreak

\begin{theorem}
\label{thJZS} \emph{(A simple isolated zero)}\\
In the context described at the beginning of Section~\ref{secExlocGeoAlg}, \propeq
\begin{enumerate}
\item The natural morphism $\,\gk\to \gA_{\uxi}\,$  is an \iso (in other words, the \id $\fm_\uxi$ is null in $\gA_{\uxi}$). In short, we write $\gk=\gA_{\uxi}$.
\item The matrix $\tra J(\uxi)$  is surjective, \cad $1\in\cD_n\big(J(\uxi)\big)$.
\item The cotangent space at $\uxi$ is null, \cad $\fm_\uxi={\fm_\uxi}^2$.
\item The \id $\fm_\uxi$ is generated by an \idm  $1-e$  of  $\gA$.
In this case the natural morphisms $\gk\to\gA[1/e]\to\gA_{\uxi}$ are \isosz.
\item There exists a $g\in\gA$ such that $g(\uxi)=1$ and $\gA[1/g]=\gk$.

\end{enumerate}
If in addition $\gk$ is a \cdi (or a \zedr \riz), we also have the \eqvc with the following \prtz.
\begin{enumerate}\setcounter{enumi}{5}
\item The tangent space $\rT_\uxi$ is null.
\end{enumerate}
 
\end{theorem}
Here is how we can describe that previous situation more intuitively: the local \alg at $\uxi$ is a \gui{connected component 
of $\gA$} (\cad the \lon at $\uxi$ is the same as the \lon at an \idmz~$e$) \gui{reduced to a simple point} (\cad this \klg is \isoc to~$\gk$).
 In terms of \vgqsz, item \emph{5}  means that there is a Zariski open set
containing the point $\uxi$ in which the \vrt is reduced to this point.
  
\begin{proof}
\emph{1}  $\Leftrightarrow$ \emph{3.}  By the localized finite \ri lemma%
~\ref{lemLocaliseFini} with $n=1$.
 
\emph{2}  $\Leftrightarrow$ \emph{3.}  By \thref{thCotangent}.
 
\emph{3}  $\Leftrightarrow$ \emph{4.} By the lemma of the \tf \idm \idz~\ref{lem.ide.idem}. \\
We then obtain the desired \isos by Fact~\ref{fact.loc.idm}, and therefore item~\emph{5}  with $g=e$.
 
\emph{5}  $\Rightarrow$ \emph{1.} 
The \egt $g(\uxi)=1$ means \hbox{that $g\in1+\fm_\uxi$}. Thus the \ri $\gA_\uxi$ is a localized \ri 
of $\gA[1/g]=\gk$, and it is equal to $\gk$ since $\gA_{\uxi}=\gk\oplus\fm_\uxi\gA_{\uxi}$.

\emph{3}  $\Leftrightarrow$ \emph{6.} (Discrete field case.) Since the tangent space is the dual of the cotangent, \emph{3}  always implies \emph{6.} 
Over a \cdi a matrix is surjective \ssi its transposed matrix is injective, this gives the \eqvc of \emph{3}  and \emph{6}  (when considering the matrix ${J(\uxi)}$). 
\end{proof}
\rem The difference between the case $s$ (number of \eqnsz) $=n$ (number of \idtrsz) and the case $s>n$ is scarcely visible in the previous \thoz, but it is important. If we tweak 
a \sys with $s=n$ and if the base field is \acz, a simple zero continues to exist, slightly 
tweaked. 
In the $s>n$ case, a tweak \gnlt makes the zero disappear. 
But this is another story, because the notion of such a tweak needs to be defined in \algz.
\eoe

For the \cdi case, here is a result in the same style as \thref{thJZS}, but more \gnl and more precise. This can \egmt be seen as a local version of Stickelberger's \tho (\thrfs{thSPolZed}{thStickelberger}).
Please note however that, unlike what takes place for Stickelberger's \thoz, the \dem of \thref{thJZScdi} does not involve the \nst or the \Noe position. However, a \cdv \`a la Nagata intervenes in the call of \thref{thNst0} 
for the implication \emph{7}  $\Rightarrow$ \emph{8.}

\begin{theorem}\label{thJZScdi} \emph{(Isolated zero)}
Suppose that $\gk$ is a \cdiz.  \Propeq
\begin{enumerate}
\item  The \alg $\gA_\uxi$ is finite over $\gk$.
\item  The \alg  $\gA_\uxi$ is integral over $\gk$.
\item  The \alg  $\gA_\uxi$ is \zedez.
\item  The \id $\fm_\uxi$ is nilpotent in  $\gA_\uxi$.
\item  There exists an $ r\in\NN$ such that $\,{\fm_\uxi}^r={\fm_\uxi}^{r+1}$.
\item  There exists an  $ r\in\NN$ such that the \id ${\fm_\uxi}^r$ is generated by an \idm $1-e$,
the morphism $\gA \to \gA_\uxi$ is surjective, \hbox{and ${\aqo\gA {1-e}\simeq\gA_\uxi\simeq\gA[1/e]}$}.
\item  There exists a $g\in \gA$ such that $g(\uxi)=1$ and $\gA[1/g]=\gA_\uxi$.
\item  There exists a $g\in \gA$ such that $g(\uxi)=1$ and $\gA[1/g]$ is local and \zedz.
\item  There exists an $h\in \gA$ such that $h(\uxi)=1$ and $\gA[1/h]$ is finite over $\gk$.
\end{enumerate}
In this case,  $\gA_\uxi$ is \stfe over $\gk$, $(\gA_\uxi)\red=\gk$,
and if $m=[\gA_{\uxi}:\gk]$, for all $\ell \in\gA_\uxi$, we have $\rC{\gA_\uxi/\gk}(\ell )(T)=\big(T-\ell (\uxi)\big)^m$.
\end{theorem}
%
\begin{proof}
The localized finite \ri lemma%
~\ref{lemLocaliseFini}, applied with $\fa=\fm_\uxi$,
shows that \emph{4}  is equivalent to \emph{5}  and implies \emph{1.}
 
 \emph{3}  $\Rightarrow$ \emph{4.} By the localized \zed \ri lemma~\ref{lemLocalisezeddim}.
 
We have \emph{1}  $\Rightarrow$ \emph{2}, and since $\gk$ is a \cdiz, \emph{2}  $\Rightarrow$ \emph{3}.
 
Thus items \emph{1}  to \emph{5}  are \eqvsz.
 
Item \emph{5}  implies that $\fm_\uxi^{r}$ is \idmz. Therefore \emph{5}  $\Rightarrow$ \emph{6}  by the \tf \idm \id lemma  \ref{lem.ide.idem} and Fact~\ref{fact.loc.idm}. 

Note that $e\in1+\fm_\uxi^{r}\subseteq 1+\fm_\uxi$,
so $e(\uxi)=1$. Therefore \emph{6}  implies~\emph{7}  \hbox{with $g=e$}.
 
\emph{7}  $\Rightarrow$ \emph{8.} The \alg $\gA[1/g]=\gA_\uxi$ is local and \tfz, and \trf by \thref{thNst0}.
 
\emph{8}  $\Rightarrow$ \emph{9.} Take $h=g$.
 
\emph{9}  $\Rightarrow$ \emph{1.} 
Because $\gA_\uxi$ is a localized \ri of $\gA[1/h]$.

In this case  $\gA_\uxi$ is \stfe over $\gk$ because it is a finite and \pf \alg (\thref{propAlgFinPresfin}).
\\
Finally, the \egt  
$\rC{\gA_\uxi/\gk}(\ell )(T)= (T-\ell (\uxi) )^m$
comes from the fact that  $\ell -\ell (\uxi)$ is in~$\fm$, so is nilpotent in~$\gA_\uxi$,
therefore it admits~$T^m$ as a \polcarz.
\end{proof}

\vspace{-.7em}
\pagebreak

\begin{definition}\label{defiZerIsole} \emph{(Isolated zero of a \syp over a \riz)}
\begin{enumerate}
\item The zero $\uxi$ of the \sys is a \emph{simple isolated zero} 
(or \emph{simple zero}) if $\gA_{\uxi}=\gk$.%
\index{simple isolated zero!of a \sypz}%
\index{simple!isolated zero}
\item The zero $\uxi$ of the \sys is an \emph{isolated zero} if  $\gA_{\uxi}$ is finite over $\gk$.
\index{isolated zero!of a \sypz}
\item If in addition $\gk$ is a \cdiz, the dimension of  $\gA_{\uxi}$  as a \kev is called the \emph{multiplicity} of the isolated zero~$\uxi$.%
\index{multiplicity!of an isolated zero (field case)}
\end{enumerate}
\end{definition}
\rem
Item \emph{1}  is an abbreviation by which we mean \prmt that the canonical \homos $\gk\to\gA_{\uxi}\to\gk$ are \isosz.\\
In item \emph{3}  we see that over a \cdiz, an isolated zero is simple 
\ssi it is of multiplicity~$1$. 
\eoe
%



\subsec{Local \ri at a non-singular point of a complete intersection curve} 
\label{subsecLocPtLisse1}

We always consider the context defined at the beginning of Section~\ref{secExlocGeoAlg}, and we assume $s=n-1$.
In other words we now have a \sys of $n-1$ \polles \eqns with $n$ unknowns and we expect the corresponding \vrt to be \gui{a curve.}

We will see that if the zero $\uxi$ of the curve is non-singular in the intuitive sense that the cotangent space at the point $\uxi$ is a \pro \kmo of rank $1$, then the \gui{local} situation matches our expectation, \cad matches what the non-singular points of the curves in \dile \gmt have accustomed us to.

\begin{theorem} {\emph{(The \id of a non-singular point of a \lot complete intersection curve)}}
       \label{thPointLisseCourbeIC}
When $s=n-1$ \propeq
\begin{enumerate}

\item
The point $\uxi$ is non-singular in the sense that $J(\uxi)$ is a matrix of rank $n-1$ over~$\gk$. 
\item
The cotangent space at $\uxi$, $\fm_\uxi/{\fm_\uxi}^2$, is a  \pro \kmo of rank~$1$. 
\item
The \id $\fm_\uxi$ is a \pro \Amo of rank $1$. 
\item
The \id $\fm_\uxi\gA_\uxi$ is a \pro $\gA_\uxi$-module of rank $1$. 
\item
The \id \smash{$\fm_\uxi\gA_\uxi$ is a free $\gA_\uxi$-module} of rank~$1$. 
\item
The cotangent space at $\uxi$, \smash{$\fm_\uxi/{\fm_\uxi}^2$}, is a free \kmo of rank~$1$.%
\end{enumerate}
\end{theorem}
%
\begin{proof} Recall that for a \ri $\gB$, a \Bmo $M$ and an \idz~$\fb$ of~$\gB$ we obtain by \eds $\gB\sur\fb\otimes_\gB M\simeq  M\sur{\fb M}$.
In particular, if $\fc$ is an \id of $\gB$ we obtain 
$(\gB\sur\fb)\otimes_\gB \fc\simeq  \fc\sur{\fb \fc}$.
\\ 
But the natural surjective \Bli $\fb \,\otimes\, \fc\to\fb \fc$ is not always an \iso (it is the case if one of the two \ids is flat).

\emph{1}  $\Leftrightarrow$ \emph{2.} 
Indeed, $\tra J(\uxi)$ is a \mpn of the cotangent space.

 \emph{3}  $\Rightarrow$ \emph{4.}  
Indeed,  the $\gA_\uxi$-module $\fm_\uxi\gA_\uxi$ is obtained from the \Amo $\fm_\uxi$ by \eds from~$\gA$ to~\smash{$\gA_\uxi$}.

 \emph{4}  $\Rightarrow$ \emph{2}  and  \emph{5}  $\Rightarrow$ \emph{6.}
Indeed, the \kmo 
\smash{$\fm_\uxi/{\fm_\uxi}^2\simeq \fm_\uxi\gA_\uxi/({\fm_\uxi\gA_\uxi})^2$} is obtained from the $\gA_\uxi$-module $\fm_\uxi\gA_\uxi$ by \eds from~$\gA_\uxi$ to~$\gk\simeq \gA_\uxi/{\fm_\uxi\gA_\uxi}$ (see the first sentence of this proof).

\emph{2}  $\Leftrightarrow$ \emph{3.} This results from the consideration of the \mpn of~$\fm_\uxi$ as an \Amo given to \thref{thidptva} and to Lemma~\ref{lemDnRz}. 
\\
To simplify the presentation let us treat the case $n=4$ 
with $\uxi=\uze$. 
\\
We have four variables $X_i$ and three \pols 

\snic{
\begin{array}{rcl} 
 f_1(\uX) & =  &  X_1a_1(\uX)+X_2a_2(\uX)+X_3a_3(\uX)+X_4a_4(\uX), \\[1mm] 
 f_2(\uX) & =  & X_1b_1(\uX)+X_2b_2(\uX)+X_3b_3(\uX)+X_4b_4(\uX) , \\[1mm] 
 f_3(\uX) & =  & X_1c_1(\uX)+X_2c_2(\uX)+X_3c_3(\uX)+X_4c_4(\uX). 
\end{array}
}

A \mpn of~${\fm_\uze}$ over  $(x_1,x_2,x_3,x_4)$  is 

\snic{W(\ux)=
\cmatrix{
x_2&x_3&0& x_4&0&0&a_1(\ux)&b_1(\ux)&c_1(\ux)
\cr
-x_1&0&x_3& 0&x_4&0&a_2(\ux)&b_2(\ux)&c_2(\ux)
\cr
0&-x_1&-x_2&0&0&x_4&a_3(\ux)&b_3(\ux)&c_3(\ux)
\cr
0&0&0&-x_1&-x_2&x_3&a_4(\ux)&b_4(\ux)&c_4(\ux)
},
}

or yet again $W(\ux)=[\,R_\ux\mid U(\ux)\,]$ with

\snic{ U(\ux)=
\cmatrix{
a_1(\ux)&b_1(\ux)&c_1(\ux)
\cr
a_2(\ux)&b_2(\ux)&c_2(\ux)
\cr
a_3(\ux)&b_3(\ux)&c_3(\ux)
\cr
a_4(\ux)&b_4(\ux)&c_4(\ux)
}\;
$ and $\;\tra J(\uze)=U(\uze).
}

We want to show that 
$W(\ux)$
(\mpn of the \Amoz~${\fm_\uze}$) \hbox{and $W(\uze)$} (\mpn of the \kmo ${\fm_\uze}/{\fm_\uze}^2$)
are simultaneously of rank $n-1=3$.\\
Refer to Lemma~\ref{lemDnRz}.
 Item \emph{3}  gives the \egt $\cD_4\big(W(\ux)\big)=0$
(because~$\cD_4\big(U(\ux)\big)=0$), and since 
%
$U(\uze)=U(\ux)\mod {\fm_\uze}$,
%
item \emph{2}  gives the \eqvc

\snic{1\in\cD_{\gA,3}\big(W(\ux)\big) \iff 1\in \cD_{\gk,3}\big(U(\uze)\big)\iff 1\in\cD_{\gk,3}\big(W(\uze)\big).}

\emph{1}  $\Rightarrow$ \emph{5.}
We reuse the previous notations with $n=4$ and $\uxi=\uze$. Since the matrix $\tra J(\uze)=U(\uze)$ is of rank $n-1$, there exist $\lambda_1$, \ldots, $\lambda_4\in\gk$ such that 
$$\preskip.4em \postskip.4em 
\det \big(V(\uze)\big)= 1 ,\quad \hbox{where} \quad V(\ux)=\cmatrix{
a_1(\ux)&b_1(\ux)&c_1(\ux)&\lambda_1
\cr
a_2(\ux)&b_2(\ux)&c_2(\ux)&\lambda_2
\cr
a_3(\ux)&b_3(\ux)&c_3(\ux)&\lambda_3
\cr
a_4(\ux)&b_4(\ux)&c_4(\ux)&\lambda_4
}\;. 
$$
We deduce that $\det \big(V(\ux)\big)\in 1+\fm_\uxi$, and so $V(\ux)\in\GL_4(\gA_\uxi)$. However, 
$$\preskip.4em \postskip.4em \ndsp
[\,x_1\;x_2\;x_3\;x_4\,]\,V=[\,0\;\;0\;\;0\;\;y\,]\;\hbox{  with  }\;y=\sum_i\lambda_ix_i. 
$$
This shows that $\gen{x_1,x_2,x_3,x_4}=\gen{y}$ in $\gA_\uxi$.
Finally, $y$ is \ndz since the module~$\fm_\uxi$ is of rank $1$.
\end{proof}

%
\entrenous{In le \tho pr\'ec\'edent, on is vraiment surpris que cela remonte jusqu'\`a l'\ri global and pas seulement \`a l'\aloz.
}

We will denote by $M^{\otimes_\gB r}$ the $r^{\rm th}$ tensor power of the \Bmo $M$.

\vspace{-.1em}
\pagebreak

\begin{theorem}\label{th2PtlisseCourbeIC}
Suppose the \eqv \prts of \thref{thPointLisseCourbeIC} satisfied, denote by $\Omega$ the cotangent space $\fm_\uxi/{\fm_\uxi}^2$ and consider an \elt $p$  of  $\fm_\uxi$ that is a $\gk$-basis of  $\Omega$.
\begin{enumerate}
\item For each $r>0$, the natural \kli $\Omega^{\otimes_\gk r}\to{\fm_\uxi}^r/{\fm_\uxi}^{r+1}$ is an \isoz.\\
In other terms, the graded \klg $\bigoplus_{r\in\NN}{\fm_\uxi}^r/{\fm_\uxi}^{r+1}$ associated with the pair $(\gA,\fm_\uxi)$ is (naturally) \isoc to the \smq \algz~$\gS_\gk(\Omega)$ of the \kmo $\Omega$, itself \isoc to $\gk[X]$ because $\Omega$ is free of rank~$1$. 
\item If $\gk$ is a nontrivial \cdiz, $\gA_\uxi$ is a \emph{discrete \adv (DVR)} in the following sense: every nonzero \elt of $\gA_\uxi$ is uniquely expressed in the form $up^{\ell}$ for some $\ell\geq0$ and $u\in\Ati$.%
\index{discrete valuation ring (DVR)}\index{ring!discrete valuation ---}%
\index{valuation!discrete --- ring (DVR)}
\end{enumerate}
\end{theorem}
%
\begin{proof} 
Let $\fm=\fm_\uxi$. Also notice that for some \pro \kmo of rank $1$, the \smq \alg is equal to the tensor \algz.

\emph{1.} We have a natural \iso $\fm^{\otimes_\gA r}\simarrow \fm^r$ because $\fm$ is flat. By the \eds $\gA\to\gA\sur\fm=\gk$, the \Amos $\fm$ and  $\fm^r$ give the \kmos $\fm/\fm^2$ and  $\fm^r/\fm\fm^r=\fm^r/\fm^{r+1}$. 
\\
Since the \eds commutes with the tensor product, we deduce that the natural \homo $\left(\fm/\fm^2\right)^{\otimes_\gk r}\to \fm^r/\fm^{r+1}$ is an \iso of \kmosz.\\
Since the \kmo $\fm/\fm^2$ admits the $\gk$-basis $p  \mod \fm^2$, the \kmo  $\fm^r/\fm^{r+1}$ admits the basis $p^r\mod \fm^{r+1}$. Hence an \iso of \klgs
$$\preskip-.1em \postskip.3em \ndsp
\kX\simarrow\bigoplus_{r\in\NN}{\fm_\uxi}^r/{\fm_\uxi}^{r+1} =\gS_\gk(\Omega), 
$$
given by $X\mapsto p$. In practice, given the filtration
$$\preskip.3em \postskip.3em 
\fm^r\subset\cdots\subset\fm^2\subset\fm\subset\gA, 
$$
every quotient of which is a free \kmo of rank $1$, the quotient~$\gA\sur{\fm^r}$ admits as its $\gk$-basis $(1,p\dots,p^{r-1})$, with for $\ell<r$ 
the $\gk$-submodule~$\fm^{\ell}\sur{\fm^r}$ which admits the basis  $(p^{\ell},\dots,p^{r-1})$.

\emph{2.} By Fact~\ref{factAlgLocEnxi}~\emph{2}  we obtain the result thanks to the following computation: if $x\in\gA_\uxi$ is nonzero, it is nonzero in some $\gA_{\uxi}/{\fm^r}$. Given the previous filtration there exists a minimum $\ell$ such that $x\in \fm^{\ell}$. If $x\equiv ap^{\ell}\mod\fm^{\ell+1}$ with $a\in\gk\eti$, we write $x=p^\ell (a+vp)$ with $v\in\gA$ and $u=a+vp$ is \iv in~$\gA_\uxi$. 
\end{proof}
%
 
\EXL{The monomial curve $t \mapsto (x_1=t^4, x_2=t^5, x_3=t^6)$}
\label{exlcourbemonomiale}

For setwise coprime $n_1$, $n_2$, $n_3 \in \NN^*$ 
we define the monomial curve $(x_1 = t^{n_1}, x_2 = t^{n_2}, x_3 = t^{n_3})$, immersed in the affine space of dimension $3$.
\\
By \dfnz, the \id of this parameterized curve is, for a \riz~$\gk$, the kernel of the morphism $\gk[X_1, X_2, X_3] \to \gk[T]$ defined by $X_i \mapsto T^{n_i}$.
\\
We can show that this \id  is always defined over $\ZZ$ and generated by three \gtrsz. Here we have chosen (see the comment at the end) the special case where~$(n_1,n_2,n_3) = (4,5,6)$, a case for which two relators suffice:

\snic{x_1^3 = x_3^2\quad$ and $\quad x_2^2 = x_1x_3.}

(Left as an exercise for the reader.) Let

\snic{\gA = \gk[x_1,x_2,x_3]=\aqo{\gk[X_1,X_2,X_3]}{X_1^3- X_3^2,X_2^2 - X_1X_3}
}

be the \ri of the curve. For $t_0 \in \gk$, we consider the point 

\snic{(\uxi) = (\xi_1,
\xi_2, \xi_3) = (t_0^{4}, t_0^{5}, t_0^{6})  ,}

with its \id $\fm = \gen {x_1-\xi_1, x_2-\xi_2, x_3-\xi_3}_\gA$. 
The condition for the point $\uxi$ to be non-singular, in the sense that the Jacobian matrix $J$ evaluated at $\uxi$ is of rank $2$, is given by $t_0 \in \gk^\times$, because $\cD_2(J)=\gen{4t_0^{11}, 5t_0^{12},6t_0^{13}}$.
From now on suppose that $t_0\in\gk\eti$. 
A \mpn  of $\fm$  for the \sgr $(x_1-\xi_1, x_2-\xi_2, x_3-\xi_3)$ is given by

\snic{W = \cmatrix {
x_2 - \xi_2  &x_3 - \xi_3  &0            &x_1^2 + \xi_1x_1 + \xi_1^2 & -x_3\cr
-x_1 + \xi_1 &      0      &x_3 - \xi_3  &0                &x_2 + \xi_2 \cr
0            &-x_1 + \xi_1 &-x_2 + \xi_2 &-x_3 - \xi_3     &  -\xi_1\cr
}.}

We know that it is of rank $2$.
We observe that $W_2, W_3 \in \gen {W_1, W_5}$.
We therefore obtain a new, simpler \mpn $V$ with only the columns $W_1, W_4, W_5$.
Recall on the one hand that for $B \in \Ae {n \times m}$, we have $(\Ae n/\Im
B)\sta \simeq \Ker \tra B$ (Fact \ref{factDualReflexif}); and on the other hand (Exercise \ref{exoMatriceCorangUn}) that for a matrix $A \in \Mn(\gA)$ of rank $n-1$,  $\Ker A = \Im
\wi{A}$ is a direct summand in $\Ae n$.
 By applying this to $B = V$ and $A = \tra \,V$, we obtain
$$\preskip.2em \postskip.1em 
\fm\sta \simeq (\Ae 3/\Im V)\sta \simeq \Ker\tra V = \Im\tra {\,\wi V} 
$$
with $\Im\tra {\,\wi V}$ a direct summand in $\Ae 3$.
\\
We thus explicitly produce the \Amo $\fm\sta$ of constant rank $1$ as a direct summand in $\Ae 3$. 
\eoe

\smallskip 
\comm
\Gnlt a sub\mo  $M$ of $(\NN,+,0)$ has a finite complement $G$ \ssi it is generated by a setwise coprime list of integers (for example with the above monomial curve we define $ M= n_1\NN + n_2\NN + n_3\NN$ generated by $\so{n_1,n_2,n_3}$).
We say that the integers of $G$ are the \emph{holes} of the \mo $M$. \\
Their number $g := \#G$ is called the \emph{genus} 
of $M$.
\\
We always have $[\,2g, \infty\,[ \;\subseteq M$.  The \mos $M$ for which $2g-1 \in G$ are said to be \emph{\smqsz}. This terminology accounts for the fact that, in this case, the interval $\lrb{0 .. 2g-1}$ contains as many holes as non-holes, and that they are interchanged by the symmetry $x \mapsto (2g-1) - x$.
\\
For example, for coprime $a$ and $b$, the \mo ${a\NN + b\NN}$ is \smq of genus $g = {(a-1)(b-1) \over 2}$. We know how to \carar the \mos $n_1\NN + n_2\NN + n_3\NN$ that are
\smqs combinatorially. We prove that this is the case \ssi the \id of the curve $(x_1=t^{n_1}, x_2=t^{n_2}, x_3=t^{n_3})$ is generated by $2$ \eltsz. For example $4\NN + 5\NN + 6\NN$ is \smqz, of genus $4$, and its holes are $\so{1, 2, 3, 7}$.
\eoe

\section{\Dcp \risz}
\label{secRelIdm}

The \ris which are \isoc to finite products of \alos play an important role in the classical theory of Henselian \alos (for example in \cite{Ray} or \cite{Laf}).
Such \ris are called \emph{decomposed \risz} and a \alo is said to be Henselian (in \clamaz) if every finite extension is a decomposed \riz.

In this section we give an introductory fragment 
of the \cov approach for the notion of a decomposed \riz. In fact, since we would like to avoid the \fcn \pbsz, we will introduce the notion, \cot more pertinent, of a \dcp \riz.

Everything begins with this simple but important remark: in a commutative \ri the \idms are always \gui{isolated.}

\begin{lemma}\label{lemIdmIsoles}
In a commutative \ri $\gA$ two \idms equal modulo $\Rad\gA$ are equal.
\end{lemma}
\begin{proof}
We show that the \homo $\BB(\gA)\to\BB(\gA\sur{\Rad\gA}\!)$ is injective. If an \idm $e$ is in $\Rad\gA$, $1-e$ is \idm and invertible, therefore equal to $1$.
\end{proof}

\rem This does not hold at all in a noncommutative context; the \idms of a \ri of square matrices $\Mn(\gA)$ are the \mprnsz; over a field we obtain, for instance by fixing the rank to $1$, a connected \vrt of dimension $>0$ without any isolated points (if $n\geq2$).
\eoe

\subsec{\Dcp \eltsz}

\begin{definition}\label{defiRelIdm} Let $\gA$ be a \ri and $a\in\gA$.
 The \elt $a$ is said to be \emph{\dcpz\footnote{Some caution must be exercised here regarding this terminology as it comes into conflict with the notion of an in\dcp \idm insofar as every \idm is a \dcp \elt of the \riz.}} if there exists
an \idm $e$ such that
$$\preskip.4em \postskip.4em 
\formule{
a \mod \gen{1-e} 
\,\mathrm{is\;\iv\;in} \,\aqo{\gA}{1-e}
\et
\\[.5mm]
 a \mod \gen{e}\,\in \,\Rad(\aqo{\gA}{e}).}
$$
\end{definition}%
\index{decomposable!element in a ring}

Recall when underlining the analogies that an \elt $a$ has a quasi-inverse \ssi there exists an \idm $e$ such that
$$\preskip.4em \postskip.4em 
\formule{
a \mod \gen{1-e} \,\mathrm{is\;\iv\;in} \,\aqo{\gA}{1-e} \et
\\[.5mm]
 a \mod \gen{e} = 0\;\mathrm{in} \;\aqo{\gA}{e},} 
$$
and that an \elt $a$ has as its annihilator an \idm  \ssi there exists an \idm $e$ such that

\snic{\formule{
a \mod \gen{1-e} \,
\mathrm{is\;regular\;in} \,\aqo{\gA}{1-e} \et
\\[.5mm]
 a \mod \gen{e} = 0\;\mathrm{in} \;\aqo{\gA}{e}.}
 }

\begin{proposition}\label{prop1DecEltAnneau}
 An \elt $a$ of $\gA$ is \dcp \ssi there exists a $b$ such that
\begin{enumerate}
\item  $b(1-ab)=0$,
\item  $a(1-ab)\in\Rad\gA$.
\end{enumerate}
In addition, the \elt $b$ satisfying these conditions is unique, and $ab=e$ is the unique \idm of $\gA$ satisfying $\gen{a}=\gen{e}\mod \Rad\gA$.
\end{proposition}
%
\begin{proof} Suppose $a$ is \dcpz.
Then, in the product $\gA = \gA_1 \times \gA_2$, with~$\gA_1 = \aqo\gA{1-e}$
and~$\gA_2 = \aqo\gA{e}$, we have $e = (1,0)$, $a = (a_1, a_2)$, with~$a_1\in\gA_1\eti$ and $a_2 \in \Rad(\gA_2)$.
We let $b = (a_1^{-1}, 0)$, and we indeed get

\snic{b(1-ab)=(b,0)-(b,0)(1,0)=0_\gA\;$ and $\;a(1-ab)=(0,a_2)\in\Rad\gA.}


Suppose that an \elt $b$ satisfies

\snic{\formule{
 b(1-ab)=0 \et
\\[.5mm]
 a(1-ab)\in\Rad\gA.}
 }

Then, the \elt $ab=e$ is an \idm and $a$ is \iv modulo $1-e$.
Moreover, modulo $e$ we have $a=a(1-e)$ which is in $\Rad\gA$,
so $a\mod e$ is in  $\Rad(\aqo{\gA}{e})$.

Let us take a look at the uniqueness. If $b(1-ab)=0$ and $a(1-ab)\in\Rad\gA$, then $e=ab$ is an \idm  such that $\gen{a}=\gen{e}\mod \Rad\gA$. 
This characterizes it as an \idm of  $\gA\sur{\Rad\gA}$, so as an \idm of $\gA$.
The \egts $be=b$  and $ba=e$ imply that $\big(b+(1-e)\big)\big(ae+(1-e)\big)=1$.
The \elt $b+(1-e)$ is therefore uniquely determined as the 
inverse of $ae+(1-e)$.
Consequently, the \elt $b$ is itself uniquely determined.
\end{proof}
%

\begin{definition}\label{defiAdcp}
 We say that the \ri $\gA$ is \emph{\dcpz} if every \elt is
\dcpz.%
\index{decomposable!\riz}%
\index{ring!decomposable ---}
\end{definition}

\begin{fact}\label{factDCP}~
\begin{enumerate}
\item \label{i1factDCP} A product of \ris is \dcp \ssi each of the factors is \dcpz.
\item A \zed \ri is \dcpz. A \dcd \alo is \dcpz. A connected \dcp \ri is local and \dcdz.
\item \label{i4factDCP} The structure of a \dcp \ri is purely equational (it can be defined by means of  composition laws subjected to universal axioms).
\end{enumerate}
\end{fact}

\begin{proof}
\emph{\ref{i4factDCP}}. We add to the laws of commutative \ris two laws

\snic{a\mapsto b \;\hbox{ and  }\; (a,x)\mapsto y,}

with the axioms $b=b^2a$ and $\big(1+x(a^2b-a)\big)y=1$.
Hence $a^2b-a\in\Rad\gA$.

\emph{\ref{i1factDCP}.}  Results from item \emph{\ref{i4factDCP}}.
\end{proof}

\rem If we let $b=a\esh$, then $(a\esh)\esh=b\esh=a^2b$ and $\big((a\esh)\esh\big)\esh=a\esh$. In addition, $(a\esh)\esh$ and $a\esh$ are quasi-inverses of one another.
 \eoe


\subsec{Lifting \idmsz}

\begin{definition}\label{defi2RelIdm} Let $\gA$ be a \riz.
\begin{enumerate}
\item We say that \emph{the \ri $\gA$ lifts the \idmsz} if the natural \homo 
$$
\preskip-.35em \postskip.2em 
\BB(\gA)\to\BB(\gA\sur{\Rad\gA}\!) 
$$
is bijective, in other words if every \idm of the quotient $\gA\sur{\Rad\gA}$ is lifted at an \idm of $\gA$.
\item We say that the \ri $\gA$ is \emph{decomposed} if it is \dcp and if $\BB(\gA)$ is bounded.
\end{enumerate}
\index{ring!decomposed ---}
\index{decomposed!\riz}
\index{ring!that lifts the \idmsz}
\end{definition}
%

\begin{proposition}\label{propDecEltAnneau}\hspace*{-1em}
\Propeq
\begin{enumerate}
\item  $\gA$ is \plc and lifts the \idmsz.
\item $\gA$ is \dcpz.
\end{enumerate}
\end{proposition}

\begin{proof}
\emph{1}  $\Rightarrow$ \emph{2.}
Since $\gA\sur{\Rad\gA}$ is \zedrz, there exists an \idm $e$ of $\gA\sur{\Rad\gA}$ such that $\gen{a}=\gen{e}\mod \Rad\gA$.
This \idm is lifted at an \idm of $\gA$, that we continue to call $e$.\\
The \elt $a+(1-e)$ is \iv in  $\gA\sur{\Rad\gA}$, so in $\gA$. Therefore,~$a$ is \iv in $\aqo{\gA}{1-e}$.
Finally, since  $\gen{a}=\gen{e}\mod \Rad\gA$, we obtain $a\in\Rad(\aqo{\gA}{e}).$

 \emph{2}  $\Rightarrow$ \emph{1.} 
Let $\pi:\gA\to\gA\sur{\Rad\gA}$ be the canonical \prnz.
Every \elt $a$ of~$\gA$ satisfies $\gen{\pi(a)}=\gen{\pi(e)}$ for an \idm $e$ of $\gA$. 
The quotient is therefore \zedz. 
Let us show that $\gA$ lifts the \idmsz. \\
If~$\pi(a)$ is  \idm and if $e$ is the \idm such that $\gen{\pi(a)}=\gen{\pi(e)}$, 
then~$\pi(a)=\pi(e)$.
\end{proof}

\comm It is now easy to see that in \clama a \ri is decomposed
\ssi it is \isoc to a finite product of \alosz. \eoe


\vspace{-.3em}


\section{\Lgb \risz}
\label{secAlocglob}

In this section we introduce a notion which generalizes both that of a \alo and that of a \zed \riz.
This sheds light on a number of facts that are common to both these classes of \risz, such as, for instance, the fact that \mptfs are quasi-free.


\subsec{\Dfns and the \plgc}


\begin{definition}
\label{defalgb}\label{defipseudolocal}~
\index{ring!local-global ---}%
\index{polynomial!primitive --- by values}%
\index{primitive by values!\pol ---}%
\begin{enumerate}
\item
We say that a \pol $f\in\AXn$ \emph{represents (in $\gA$)
the \elt $a\in\gA$} if there exists an  $\ux\in\Ae n$ such that $f(\ux)=a$.

\item
We say that a \pol $f\in\AXn$ is \emph{\ppvz} if the values of $f$ generate the \id $\gen{1}$ (the variables being evaluated in~$\gA$).
\pagebreak

\item A \ri $\gA$ is said to be \ixc{local-global}{ring}
 if every \pppv represents an \invz.
\end{enumerate}
\end{definition}

\rem Every \pppv is primitive, therefore if a \ri has the \prt that every primitive \pol represents an \invz, it is a \lgb \riz. This corresponds to a \dfn in the literature (strongly  U-\ird \riz) which has preceded that of \algbz. 
\eoe

\begin{fact}
\label{factlgb1} ~
\begin{enumerate}
\item \label{i1factlgb1} A \ri $\gA$ is \lgb  \ssi $\gA\sur{\Rad(\gA)}$ is \lgbz.
\item A finite product of \ris is \lgb \ssi each of the \ris is \lgbz.
\item A \alo is \lgbz.
\item \label{i4factlgb1} A \plc \ri is \lgbz.
\item \label{i5factlgb1} A quotient of a \lgb \ri (resp.\ \plcz) is \lgb (resp.\ \plcz).
\item \label{i6factlgb1}
Let $\gA$ be a non-decreasing filtering union \ri of sub\ris $\gA_i$,
\cad for all $i, j$, there exists a $k$ such that $\gA_i \cup \gA_j \subseteq
\gA_k$. Then, if each $\gA_i$ is \lgbz, so is $\gA$.
\end{enumerate}
\end{fact}
\begin{proof}
We leave the first three items as an exercise.

\emph{\ref{i4factlgb1}.} Given item \emph{\ref{i1factlgb1}}, it suffices to treat the case of a \zedr \riz. This case reduces to the (obvious) case of a \cdi by the  \elgbm \num2.

\emph{\ref{i5factlgb1}.} Let us consider the \lgb case (the other case is obvious).
Let~$\gA$ be a \lgb \riz, $\fa$ be an \id and $f\in\gA[\uX]$ be a \pppv in $\gA\sur\fa$.
Therefore there are some values $p_1$, \ldots, $p_m$ of $f$  and some $a\in\fa$ such that $\gen{p_1,\ldots ,p_m,a}=\gen{1}$. The \pol $g(\uX,T)=Tf(\uX)+(1-T)a$
is therefore \ppvz. Since $\gA$ is \lgbz,
there is a value $tf(\ux)+(1-t)a$  of $g$ which is invertible. The value $f(\ux)$ is thus invertible modulo $\fa$.

\emph{\ref{i6factlgb1}.}
Let $P\in\AXn$ be \ppvz: $1 = uP(\ux) + vP(\uy) +\dots$. \\ 
By considering $u$, $\ux$, $v$, $\uy$, \dots\, and the \coes of $P$, we see that there is a sub\ri  $\gA_i$ such that $P \in \gA_i[\uX]$ and such that $P$ is \ppv over~$\gA_i$.  Thus, $P$ represents an \inv over $\gA_i$, a fortiori over $\gA$.
\end{proof}

For a \pol the \prts of representing an \inv or of being \ppv are of
 \carfz, 
 as indicated in the following lemma.
 
 \pagebreak

\begin{lemma}\label{lemAlgb0}
Let $S$ be a \mo of $\gA$ and $f\in\gA[X_1, \ldots, X_m]$ be a \polz.
\begin{enumerate}
\item  The \pol $f$ represents an \inv in $\gA_S$ \ssi there exists an $s\in S$
such that $f$ represents an \inv in~$\gA_s$.
\item  The \pol $f$ is \ppv in $\gA_S$ \ssi there exists an $s\in S$
such that $f$ is \ppv in~$\gA_s$.
\end{enumerate}
\end{lemma}
\begin{proof}
We only prove item \emph{1.} Let $F(\uX,T) \in \gA[\uX,T]$
be the homogenization of $f(\uX)$ at a large enough degree.
The hypothesis is equivalent to the existence of $\ux \in \Ae m$ and $t$, $u \in S$ such that $F(\ux,t)$ divides $u$ in $\gA$. Letting $s = tu$, the \elts $t$ and $F(\ux,t)$ are \ivs in
$\gA_s$ so $f$ represents an \inv in~$\gA_s$.
\end{proof}
%

\begin{lemma}\label{lemAlgb1}
Let $s \in \gA$ and $\fb$ be an \id of $\gA$ with $1 \in \gen {s} + \fb$.  
\begin{enumerate}
\item If $f$ represents an \inv in $\gA_s$ there exists a $\uz \in \Ae m$ such that $1 \in \gen {f(\uz)} + \fb$.
\item If $f$ is \ppv in $\gA_s$ there exists a finite number of \elts $\uz_j$, ($j\in\lrbk$), in~$\Ae m$
such that  $1 \in \gen {f(\uz_j)\mid j\in\lrbk} + \fb$.
\end{enumerate}
\end{lemma}
\begin{proof}
\emph{1.} Let $F(\uX,T) \in \gA[\uX,T]$ be the homogenization of~$f(\uX)$ at large enough degree $d$. The hypothesis is that $F(\ux, t)$ divides $u$ in~$\gA$ for some $\ux \in
\Ae m$ and $t, u \in s^\NN$. There exists an $a$ such that $ta \equiv 1 \bmod\fb$
so

\snic{
a^d F(\ux, t) = F(a\ux, at) \equiv F(a\ux, 1) = f(a\ux) \bmod \fb,
}

hence  $a^du \in \gen {f(\uz)} + \fb$ with $\uz = a\ux$.
But $1 \in \gen {a^du} + \fb$ therefore $1 \in \gen {f(\uz)} + \fb$.

 We can present the same argument \gui{without computation} as follows.
\\
We have $\gA_s\sur{(\fb\gA_s)} \simeq (\gA\sur\fb)_s$. Since $1 \in \gen {s} + \fb$, $s$ is \iv in $\gA\sur\fb$, and so $\gA_s\sur{(\fb\gA_s)} \simeq \gA\sur\fb$. Since $f$ represents an \inv in $\gA_s$, a fortiori it represents an \inv in $\gA_s\sur{(\fb\gA_s)} \simeq \gA\sur\fb$, \cad $f$
represents an \inv modulo $\fb$.

\emph{2.} Similar to Item \emph{1.} 
\end{proof}

We will use in the remainder a slightly more subtle \plgc that we state in the form of a lemma. See also Exercise \ref{exopropAlgb1}.

\begin{lemma}\label{lempropAlgb1}
Let $S_1$, $\dots$, $S_n$ be \moco of $\gA$ and $f\in\gA[X_1, \ldots, X_m]$ be a \polz. \Propeq
\begin{enumerate}
\item [1.] The \pol $f$ is \ppvz.
\item [2.] In each of the \ris $\gA_{S_i}$,  the \pol $f$ is \ppvz.
\item [3.]$\!\!\!\!$\eto For every \idema $\fm$  of $\gA$, $f$ represents an \inv in~$\gA\sur\fm$.
\end{enumerate}
In particular, if $f$ represents an \inv in each localized \ri $\gA_{S_i}$,  $f$ is \ppvz.
\end{lemma}

\begin{proof}
The implications \emph{1}  $\Rightarrow$ \emph{2}   $\Rightarrow$ \emph{3\eto}$\!$
are \imdesz. The implication  \emph{3\eto}  $\!\!\Rightarrow$ \emph{1}  is easy in \clamaz.

Here is  a direct and \cov \dem of \emph{2}  $\Rightarrow$ \emph{1.} It is a matter of decrypting the classical \dem of \emph{3.\eto}  $\!\!\Rightarrow$ \emph{1}, by using the method that will be explained in Section~\ref{subsecLGIdeMax}.
To simplify the notations but \spdgz, we will prove the special case where $f$ represents an \inv in each localized \ri $\gA_{S_i}$.
\\
We therefore dispose of \eco $(s_1, \ldots, s_n)$ such that in each localized \riz~$\gA_{s_i}$, 
the \polz~$f$ represents an \inv (Lemma~\ref{lemAlgb0}).
\\
By applying Lemma~\ref{lemAlgb1} we successively obtain, for $k = 0, \ldots, n$, 
\[\preskip.3em \postskip.2em
1 \in \gen {f(\und {z_1}), \ldots, f(\und {z_k}),\ s_{k+1}, \ldots, s_n}.
\]  After $n$ steps: $1 \in \gen {f(\und {z_1}), \ldots, f(\und {z_n})}$.
\end{proof}

\begin{proposition}\label{propAlgb1} \hspace*{-.8em}
 \Propeq
\begin{enumerate}
\item The \ri $\gA$ is \lgbz.
\item For every \pol $f\in\AXn$, if there exists a \sys of \eco $(s_1,\ldots,s_k)$
such that $f$ represents an \inv in each~$\gA_{s_i}$, then $f$ represents an \invz.
\item For every \pol $f\in\AXn$, if there exist \moco $S_i$ such that $f$  is \ppv in each $\gA_{S_i}$, then~$f$ represents an \invz.
\end{enumerate}
\end{proposition}

\begin{proof}
Given Lemmas~\ref{lemAlgb0} and \ref{lempropAlgb1}, 
it suffices to show that if $f$ is \ppv there exist \eco such that $f$ represents an \inv in each localized \riz.
To simplify the notation, we will write everything using a single variable.
We obtain $x_1$, \dots, $x_r \in\gA$ such that $1\in\gen{f(x_1), \ldots, f(x_r)}$.  Let~$s_i =
f(x_i)$, then the \pol $f$ represents an \inv in $\gA_{s_i}$.
\end{proof}

\rdb
By the Gauss-Joyal lemma (\ref{lemGaussJoyal}) the primitive \pols form a filter $U\subseteq\AX$.
We call the \ri 
 $\gA(X)=U^{-1}\gA[X]$ the \ix{Nagata ring}.

\begin{fact}\label{factLocNagata}
We use the above notation.
\begin{enumerate}
  \item $\gA(X)$ is faithfully flat over $\gA$.
  \item $\gA(X)$ is a \algbz.
\end{enumerate}
\end{fact}
\begin{proof} 
\emph{1.} It is clear that $\gA(X)$ is flat over $\gA$ 
(it is a \lon of $\AX$, which is free with a discrete basis). 
We then use the \carn \emph{3a}  in \thref{thExtFidPlat}. Let $\fa=\gen{\an}$ be a \itf of~$\gA$ such that $1\in\fa\gA(X)$. We must show that $1\in\fa$. The hypothesis gives  $f_1$, \ldots, $f_n\in\AX$ such that the \pol $f=\sum_ia_if_i$ is primitive, \cad $1\in\rc_\gA(f)$. However, the \id $\rc_\gA(f)$ is contained in $\fa$.

\emph{2.} We proceed in three steps.\\
a) Let us first show that every primitive \pol  $P(T)\in\gB[T]$ where $\gB:=\gA(X)$ represents an \iv \eltz.
Indeed, let $P(T) = \sum_i  Q_iT^i$ be such a \polz.
We can suppose \spdg that the~$Q_i$'s are in $\gA[X]$. We have \pols $B_i$ such that $\sum_i  B_i(X)Q_i(X)$ is primitive. A fortiori the \coes of the $Q_j$'s are \comz. 
\\
Then, for $k > \sup_i\!\big(\deg_X(Q_i)\big)$, since $P(X^k)$ has for \coes all the \coes of the $Q_j$'s (\KRAz's trick), it is a primitive \pol of $\gA[X]$, \cad an \iv \elt of $\gB$.\\
b) Let us show the same \prt for an arbitrary number of variables.
Consider a primitive \pol $Q(\Ym)\in\BuY$. By Kronecker's trick, by letting $Y_j=T^{n^j}$ with large enough $n$, we obtain a \pol $P(T)$ whose \coes are those of $Q$, which brings us back to the previous case.\\
c) Finally, consider a \pppv $Q$ with $m$ variables over $\gB$. Then,  $Q$ is primitive and we can apply item~b).
\end{proof}

\subsec{Remarkable \lgbes \prtsz}

\begin{plcc}
\label{thlgb1} Let $S_1$, $\ldots $, $S_r$ be \moco of a \algb $\gA$.
\begin{enumerate}
\item If two matrices of $\Ae {m\times n}$ are \eqves over each of the $\gA_{S_i}$, then they are \eqvesz.
\item If two matrices of $\Mn(\gA)$ are similar over each of the $\gA_{S_i}$, then they are similar.
\end{enumerate}
\end{plcc}
\begin{proof}
\emph{1.} Let $F$ and  $G$ be the matrices, then by hypothesis there exists some \sys of \eco  $(s_1, \ldots , s_r)$  and matrices $U_1$, \ldots, $U_r$, $V_1$, \ldots, $V_r $ such that for each~$i$ we have
     $U_i F = G V_i$ and 
     $\det(U_i) \det(V_i) = s_i$.
Let us introduce \idtrs  $(x_1, \ldots , x_r)=(\ux)$, and consider the matrices   

\snic{U=U(\ux) = x_1\,U_1+\cdots +x_r\,U_r$ and $V=V(\ux) = x_1\,V_1+\cdots
+x_r\,V_r.
}


We have  $U F = G V,$ and $\det(U) \det(V)$ is a \pol in the $x_i$'s that satisfy the hypotheses of \Dfnz~\ref{defalgb}; it suffices to evaluate $(x_1, \ldots , x_r)$ successively at $(1,0,\ldots ,0)$, \ldots, $(0,\ldots ,0,1)$.
Therefore there exists some $\ual\in\Ae r$ such that the \eltz~$\det\big(U(\ual)\big) \det\big(V(\ual)\big)$ is \ivz.

\emph{2.} The same \demz, with $U_i=V_i$ and $U=V$, works.
\end{proof}

We have the following corollary.

\pagebreak

\begin{plcc}
\label{thlgb2} Let $S_1$, $\ldots$, $S_r$ be \moco of a \algb $\gA$.
\begin{enumerate}
\item If two \mpfs are \isoc over each of the~$\gA_{S_i}$'s, then they are \isocz.
\item Every \mptf is quasi-free.
\end{enumerate}
\end{plcc}
\begin{proof}
\emph{1.} We consider \mpns and \carar the fact that the modules are \isoc by the \eqvc of associated matrices (Lemma~\ref{lem pres equiv}).
We then apply item~\emph{1}  of the \plgrf{thlgb1}.

\emph{2.} We apply item~\emph{1.}
Consider a quasi-free module that has the same \idfsz, we know that the two modules become free after \lon at \eco (and the rank is the same each time because they have the same \idfsz).
\end{proof}

Let us also mention the following principles.

\begin{plcc}
\label{thlgb3} Let $\gA$ be a  \algbz.
\begin{enumerate}
\item Let $S_1$, $\ldots$, $S_r$ be \mocoz, $M$ be a \mpf and $N$ a \mtfz. If $N$ is a quotient of $M$ over each of the~$\gA_{S_i}$'s, then $N$ is a quotient of $M$.
\item A module locally generated by $m$ \elts is generated by $m$ \eltsz.
\end{enumerate}
\end{plcc}

%
\begin{proof}
It suffices to prove item \emph{1}  because a module is generated by $m$ \elts \ssi it is a quotient of a free module of rank $m$.
\\
We will continue the \dem after the next two lemmas.
\end{proof}

\begin{lemma}\label{lem1thlgb3}
Let $M$ be a \pf \Amoz, $N$ be a \tf \Amoz, $S$ be a \mo of $\gA$ and $\varphi:M_S\to N_S$ be a surjective \Aliz. 
\begin{enumerate}
\item There exist  $s\in S$ and  $\psi\in\Lin_\gA(M,N)$ such that $s\varphi=_{\gA_S}\psi_S$.
\item There exists a $v\in S$ such that $vN\subseteq\psi(M)$.
\item There exists a matrix $Q$ of syzygies
satisfied by the \gtrs of~$N$ such that, when considering the module $N'$ admitting~$Q$ as a \mpnz, the map~$\psi$ is decomposed as follows

\snic{M\vers{\theta}N'\vers{\pi}N,}

($\pi$ is the canonical \prnz), with in particular $vN'\subseteq\theta(M)$
(a fortiori~$\theta_S$ is surjective). 
\end{enumerate}
\end{lemma}
%
\begin{proof}
Item \emph{1}  is a reformulation of Proposition~\ref{fact.homom loc pf} (which affirms only slightly more, in a more \gnl case). Item \emph{2}  easily stems from it.
\\
\emph{3.} We have $N=\gA y_1+\cdots+\gA
y_n$, and $M=\gA x_1+\cdots+\gA x_m$, with a \mpn $P$.  
\\
For the \fcn by $\theta$ to exist, it suffices that among the columns of the matrix $Q$ we find the sygygies 
which are \gui{images of the columns of $P$ by $\psi$} (they are syzygies
between the $y_k$'s once we have expressed the $\psi(x_j)$'s in terms of the $y_k$'s). \\
For $vN'\subseteq\theta(M)$ to hold, it suffices that among the columns of the matrix $Q$ we find the syzygies
 expressing that the $vy_k$'s are in $\gA\psi(x_1)+\cdots+\gA\psi(x_m)$ (once we have expressed the $\psi(x_j)$'s in terms of the $y_k$'s).
\end{proof}
%

\begin{lemma}\label{lempfthlgb3}
The \plgc \ref{thlgb3} is correct if $N$ is itself a \mpfz. 
\end{lemma}
%
\begin{proof} The hypothesis gives a surjective \ali $\varphi_i:M_{S_i}\to N_{S_i}$. By items~\emph{1}  and \emph{2}  of Lemma \ref{lem1thlgb3} we have $s_i,v_i\in S_i$ and a \ali $\psi_i:M\to N$ such that $s_i\varphi_i=(\psi_i)_{S_i}$ and $v_iN\subseteq\psi_i(M)$.
Each \ali $\psi_i$ is represented by two matrices $K_i$ and $G_i$ which make the suitable diagrams commute (see Section~\ref{secCatMpf}).
$$
\preskip-.8em \postskip.0em 
\xymatrix {
\gA^{p} \ar[r]^{P} \ar[d]_{K_i} & \Ae m \ar[d]^{G_{i}}
\ar@{->>}[r]^{\pi_M}
                   & M \ar[d]^{\psi_i} \\
\gA^{q} \ar[r]_{Q}                    & \Ae n \ar@{->>}[r]_{\pi_{N}}
                   & N \\} 
$$
Consider $r$ unknowns $a_i$ in $\gA$ and the map $\psi = \sum a_i \psi_i$ corresponding to the matrices  $K=\sum a_i  K_i$ and $G=\sum a_i  G_i$.
$$
\preskip.2em \postskip.2em 
\xymatrix {
\gA^{p} \ar[r]^{P} \ar[d]_{K} & \Ae m \ar[d]^{G}
\ar@{->>}[r]^{\pi_M}
                   & M \ar[d]^{\psi} \\
\gA^{q} \ar[r]_{Q}                    & \Ae n \ar@{->>}[r]_{\pi_{N}}
                   & N \\} 
$$
The fact that $\psi$ is surjective means that the matrix \smashtop{$H=\blocs{.7}{.8}{.6}{0}{$G$}{$Q$}{}{}$} is surjective, \cad $\cD_n(H)=\gen{1}$.  
We therefore introduce the \idtrs $c_\ell$ to construct an arbitrary \coli of the maximal minors $\delta_\ell$ of the matrix $H$. 
This \coli $\sum_\ell  c_\ell\delta_\ell$ is a \pol in the $a_i$'s and $c_\ell$'s. By hypothesis, this \pol represents $1$ over each of the $\gA\big[\frac 1{s_iv_i}\big]$, thus, since the \ri
is \lgbz, it represents an \inv (Proposition~\ref{propAlgb1}).
\end{proof}
%

\emph{End of the \dem of the \plgc \ref{thlgb3}.}
\\ 
 We have $N=\gA y_1+\cdots+\gA y_n$, and $M=\gA x_1+\cdots+\gA x_m$, with a \mpn $P$.  
For each $i\in\lrbr$ we apply Lemma~\ref{lem1thlgb3} with the \mo $S_i$ and the surjective \ali  $\varphi_i:M_{S_i}\to N_{S_i}$ given in the hypothesis. We obtain a \ali $\psi_i:M\to N$, a matrix~$Q_i$ of syzygies
satisfied by the $y_k$'s, a \aliz~$\theta_i:M\to N'_i$ (where~$N'_i$ is the \mpf corresponding to~$Q_i$), \elts $s_i$, $v_i\in S_i$ with $s_i \varphi_i=(\psi_i)_{S_i}$, and finally~$\psi_i$ factorizes through~$\theta_i:M\to N'_i$ with~$v_iN'_i\subseteq\theta_i(M)$.  
\\
We then consider the \pf module $N'$ corresponding to the matrix of syzygies 
$Q$ obtained by juxtaposing the matrices~$Q_i$, such that~$N'$ is a quotient of each~$N'_i$.
\\
As $N$ is a quotient of $N'$, we have brought the \pb back to the case where $N$ is itself \pfz, a case that has been treated in Lemma~\ref{lempfthlgb3}.
\eop


\subsec{Congruential \syssz}
\label{subsecSysCong}

An important stability \prt of \algbs is  stability by integral extension.

\begin{theorem}\label{thLgbExtEnt}
Let $\gA\subseteq\gB$ with $\gB$ integral over $\gA$.
If $\gA$ is \lgbz, then so is $\gB$. \perso{\ssi?}
\end{theorem}

The \dem is left until \paref{propSysCong}, after a detour via  congruential \risz.

\begin{definition}\label{defiSysCong}
A subset $C$ of a \ri $\gA$ is called a \emph{congruential \sysz} if it satisfies the following \prtz: if $s_1+s_2=1$ in $\gA$ and if~$c_1$, $c_2\in C$, then there exists a $c\in C$ such that $c\equiv c_1\mod s_1$ and $c\equiv c_2\mod s_2$.
\index{system!congruential ---}
\index{congruential!system}
\end{definition}

\rems 1) It amounts to the same thing to say: if $\fa_1$ and $\fa_2$ are two \com \ids of $\gA$ and if $c_1$, $c_2\in C$, then there exists a $c\in C$ such that $c\equiv c_1\mod \fa_1$ and $c\equiv c_2\mod \fa_2$.\\
2) The \elt $c'=c_2s_1+c_1s_2$ is the natural candidate for $c\in\gA$ satisfying the congruences $c\equiv c_1\mod s_1$ and $c\equiv c_2\mod s_2$. 
We therefore must have some \eltz~$c$ of~$C$ such that $c\equiv c'\mod s_1s_2$.
\eoe

\medskip \rdb 
\exl \label{defiSuslinSet}
Let $(\ub) = (\bn)$ be a sequence in a \ri $\gB$. The \emph{Suslin set of $(\bn)$} is the following subset of $\gB$:%
\index{Suslin!set of a finite sequence}%

\snic {
\Suslin(\ub) = \sotq{u_1b_1 + \cdots + u_nb_n} {
(\un)\hbox{ is $\EE_n(\gB)$-completable}},
} 

($(\un)$ is the first row of a matrix of $\EE_n(\gB)$).

If one of the $u_i$'s is \ivz, then $u_1b_1 + u_2 b_2 + \cdots + u_n b_n \in \Suslin(\ub)$ and we therefore have $\{\bn\} \subseteq \Suslin(\bn) \subseteq \gen {\bn}$.

Let us show that the set $\Suslin(\ub)$ is always congruential. 
\\
Indeed, for $E$, $F \in \EE_n(\gB)$ and two comaximal \elts $s$, $t$ of $\gB$, there exists a $G \in \EE_n(\gB)$ satisfying $G \equiv E \bmod s$ and $G \equiv F \bmod t$.

 Let $f$, $g_1$, \ldots, $g_n  \in \AX$ with $f$ \monz, and $\gB=\aqo\AX f$. 
Then the Suslin set of $(\ov {g_1},\dots,\ov {g_n})$ plays an important role in the study of the 
unimodular \poll vectors (cf. Lemma~\ref{lemSuslin1}).
\eoe

\begin{fact}\label{fact1SysCong}
For every \pol $P\in\AXn$ the set $V_P$ of values of $P$ is a congruential \sys  ($V_P=\sotq{P(\ux)}{\ux\in\Ae n}$).
\end{fact}
\begin{proof}
Let $s$, $t$ be two \eco and $\ux$, $\uy$ be in $\Ae n$. The Chinese remainder \tho
 gives us some $\uz\in\Ae n$ such that $\uz\equiv \ux\mod s$ and  $\uz\equiv \uy\mod t$. Then, we have $P(\uz)\equiv P(\ux)\mod s$ and  $P(\uz)\equiv P(\uy)\mod t$.
\end{proof}

\vspace{-.7em}
\pagebreak

\begin{fact}\label{fact2SysCong}
Let $C$ be a congruential \sysz. If $\fa_1$, \ldots, $\fa_\ell$ are pairwise \com \ids and if $c_1$, \ldots, $c_\ell\in C$, then there exists a $c\in C$ such that~$c\equiv c_j\mod \fa_j$ for $j\in\lrb{1..\ell}$.
\end{fact}
\begin{proof}
This is the usual \dem of the Chinese remainder \thoz, 
adapted to the current situation.
We proceed by \recu on $\ell\geq2$. The base case is by \dfnz.
If $\ell>2$ we consider the pairwise \com \ids $\fa_1,\ldots,\fa_{\ell-2}$ and $\fa_{\ell-1}\fa_{\ell}$.
Let $e\in C$ such that $e\equiv c_{\ell-1}\mod\fa_{\ell-1}$ and $e\equiv c_\ell\mod\fa_\ell$. \\
By \hdrz, we find $c$ in $C$ such that $c\equiv c_k\mod\fa_k$ for~$k\in\lrb{1..\ell-2}$ and~$c\equiv e\mod\fa_{\ell-1}\fa_{\ell}$. A fortiori, $c\equiv c_{\ell-1}\mod\fa_{\ell-1}$ and~$c\equiv c_\ell\mod\fa_\ell$.
\end{proof}
%

\begin{fact}\label{fact3SysCong}
Let $C$ be a congruential \sysz, $w_1$, \ldots, $w_n$ be \elts of~$C$ and~$(e_1,\ldots,e_n)$ be a \sfioz.
Then, the \elt $w=e_1w_1+\cdots+e_nw_n$ is in~$C$.
\end{fact}
\begin{proof}
We have $w\equiv w_i \mod 1-e_i$, and the $\gen{1-e_i}$'s are pairwise \comz, but $w$ is the unique \elt satisfying these congruences since $\bigcap_i\gen{1-e_i}=\gen{0}$. It remains to apply the previous fact.
\end{proof}
%
\begin{definition}\label{defiAnneauCongruentiel}
A \ri $\gA$ is said to be \ixc{congruential}{ring} if every congruential \sys that generates the \id $\gen{1}$ contains an \iv \eltz.%
\index{ring!congruential ---}
\end{definition}

\begin{lemma}\label{lemAnneauCongruentiel}~
\begin{enumerate}
\item \label{i3lemAnneauCongruentiel} Let $\fa\subseteq\Rad\gA$. Then, the \ri $\gA$ is congruential \ssi the \ri $\gA\sur{\fa}$ is congruential.
\item \label{i1lemAnneauCongruentiel} Every \plc \ri is congruential.
\item  \label{i2lemAnneauCongruentiel} Every congruential \ri is \lgbz.
\end{enumerate}
 \end{lemma}
\begin{proof}
\emph{\ref{i3lemAnneauCongruentiel}.}
We use the fact that \elts are \com (resp.\ \ivsz) in $\gA$ \ssi they are \com (resp.\ \ivsz) in~$\gA\sur{\fa}$.
 
\emph{\ref{i1lemAnneauCongruentiel}.}
Let us suppose that $\gA$ is \plcz. It suffices to show that~$\gA\sur{\Rad\gA}$ is congruential. Let $W$ be a congruential \sys of $\gA\sur{\Rad\gA}$ such that~$\gen{W}=\gen{1}$. Let $w_1$, \ldots, $w_n\in W$ with $\gen{w_1,\ldots,w_n}=\gen{1}$. There exists a \sfio  $(e_1,\ldots,e_n)$ such that we have~$\gen{e_1w_1+\cdots+e_nw_n}=\gen{1}$ (Lemma \ref{lemZerRed} item~\emph{5}).
We conclude with Fact~\ref{fact3SysCong} that~$W$ contains the \iv \elt $e_1w_1+\cdots+e_nw_n$.
 
\emph{\ref{i2lemAnneauCongruentiel}.}
Let us suppose that $\gA$ is congruential and let $P$ be a \pppvz. Since the values of $P$ form a congruential \sysz, a value of $P$ is \ivz.
\end{proof}

\vspace{-.4em}
\pagebreak

\subsec{Stability by integral extension}

As an \imd corollary of Lemma \ref{lemAnneauCongruentiel} we have the following result.
\begin{corollary}\label{CorSysCong}
Let $\gB$ be a \stfe \alg over a \cdiz~$\gA$ and $W$ be a congruential \sys in $\gB$ such that~$\gen{W}=\gen{1}_\gB$.
\\
Then, the set $\rN\iBA(W)$ contains an \iv \eltz.
\end{corollary}
\begin{proof}
We know that $\gB$ is \zedz, so it is congruential (Lemma~\ref{lemAnneauCongruentiel}). Since $W$ is congruential and generates the \id $\gen{1}$, it contains an \iv \eltz. Finally, the norm of an \iv \elt is \ivz.
\end{proof}
%

\begin{proposition}\label{propSysCong}
Let $\gB$ be a \stfe \alg over a \riz~$\gA$ and $W$ be a congruential \sys in $\gB$. If~$1\in\gen{W}$,
then,   \hbox{$1\in\gen{\rN\iBA(W)}$}.
\end{proposition}
\begin{proof} 
\emph{1.} A congruential \sys remains congruential by passage to a quotient \riz. 
If we read the conclusion of Corollary~\ref{CorSysCong} in the (weaker) form~\hbox{$1\in\gen{\rN\iBA(W)}$}, we observe that it is in an adequate form to be subjected to the \cov machinery with \idemas which will be explained on \paref{MethodeIdemax} in Section~\ref{subsecLGIdeMax}, and which is used to prove that an \id contains $1$. We therefore obtain the desired result.\imlma
\end{proof}

\rems \\
1) In \clama we would also say this: if $1\notin \gen{\rN\iBA(W)}_\gA$, this \id would be contained in a \idema $\fm$  of $\gA$. But Corollary~\ref{CorSysCong}, applied with the \cdi $\gA/\fm$ and the \stfe \alg $\gB/\fm\gB$, shows that it is impossible. \\
The \cov machinery with \idemas precisely aims at decrypting this type of abstract \dem and at transforming it into an \algo which constructs $1$ as an \elt of~$\gen{\rN\iBA(W)}_\gA$ from the hypotheses.\imlma

2) As an example, if $(\ub) = (b_1, \ldots, b_q)$  is a \sys of \eco in $\gB$, we have $1 \in \gen{\rN_{\gB/\gA}(w) \mid w \in \Suslin(\ub)}_\gA$,
since the set $\Suslin(\ub)$ is congruential. 
But we will refrain from believing \hbox{that $1 \in \gen {\rN_{\gB/\gA}(b_1), \ldots,
\rN_{\gB/\gA}(b_q)}_\gA$}. 
\\
A famous instance of this \prt is a result due to Suslin regarding \poll vectors, given in Lemma~\ref{lemSuslin1}. In this lemma, $\gB$ is of the form $\aqo{\gA[X]}{v}$ with $v \in \gA[X]$ a \poluz. A complete decrypting will be provided in the \dem of the lemma in question.
\eoe

\begin{Proof}{\Demo of \thrf{thLgbExtEnt}. }
Let us first treat the case where~$\gB$ is free of finite rank, say $\ell$, over $\gA$.  
Let $P\in\BXn$ be a \pppvz. We want some $\ub\in\gB^n$ with $P(\ub)$ \ivz. We consider the congruential \sys $W$ of the values of $P$. By hypothesis \hbox{we have $1\in\gen{W}$}. Proposition~\ref{propSysCong} then says that~$\gen{\rN\iBA(W)}_\gA=\gen{1}_\gA$.\\
But $\rN\iBA\big(P(b_1,\ldots,b_n)\big)$ is a \pol with $n\ell$ variables in $\gA$ if we express each $b_i\in \gB$ over an $\gA$-basis of $\gB$, and $\gA$ is \lgbz, so there exists a~$\ub\in\gB^n$ such that $\rN\iBA\big(P(\ub)\big)$ is \ivz, and this implies that $P(\ub)$ is \ivz.\\
In the \gnl case where $\gB$ is only assumed to be integral over $\gA$, let us consider in $\gB$ the \tf \Aslgs $\gB_i$; $\gB$ is its increasing filtering union. Since $\gB$ is integral over $\gA$, so is $\gB_i$, therefore it is a quotient of an \Alg which is a free \Amo of finite rank. By the first case, and in virtue of item \emph{\ref {i5factlgb1}}  of Fact~\ref{factlgb1}, each $\gB_i$ is \lgbz.
Finally, by the last item of Fact~\ref{factlgb1}, $\gB$ is \lgbz.
\end{Proof}
%


\Exercices

\begin{exercise}
\label{exoNilRad}
{\rm  Prove in \clama that the nilradical of a \ri is equal to the intersection of its \idepsz.
}
\end{exercise}

\vspace{-1em}
\begin{exercise}
\label{exoRadJacSat}
{\rm If $\fa$ is an \id of $\gA$ we let $\JA(\fa)$ be its \emph{Jacobson radical}, \cad the inverse image of $\Rad(\gA\sur\fa)$ under the canonical projection $\gA\to\gA\sur\fa$.
 Let~$\fa$ be an \id of $\gA$. Show that $\JA(\fa)$ is the greatest \id $\fb$
 such that the \moz~$1+\fb$ is contained in the saturated \mo of $1+\fa$.
} 
\end{exercise}

\vspace{-1em}
\begin{exercise}
\label{exoRadJacLoc}
{\rm  Prove in \coma that the Jacobson radical of a \alo coincides with the set of non\iv \eltsz, 
and that it is the unique \id $\fa$ satisfying
\begin{itemize}\itemsep0pt
\item $\fa$ is maximal
\item $1\in\fa$ implies $1=0$.
\end{itemize}
}
\end{exercise}

\vspace{-1em}
\begin{exercise}
\label{exoRadNonCom}
{\rm Let $\gA$ be a noncommutative \riz, $a,b\in\gA$. Prove the following statements.
 \emph{1.} If $a$ admits a left-\inv $c$, then $c$ is a right-\inv  of $a$ \ssi $c$ is unique as a left-\inv of $a$.

%
 \emph{2.} If $1-ab$ admits a left-\inv $u$, then $1-ba$ also admits a left-\inv $v$. Idea: if $ab$ and $ba$ are \gui{small,} $u$ must be equal to $1+ab+abab+\dots $, and~$v$ equal to $1+ba+baba+\cdots =1+b(1+ab+abab+\cdots)a$. 

 \emph{3.}  If for all $x$, $1-xa$ is left-\ivz, then for all $x$, $1-xa$ is right-\ivz. 

 \emph{4.}  \Propeq
\begin{itemize}\itemsep0pt
\item For all $x$, $1-xa$ is left-\ivz.
\item For all $x$, $1-xa$ is right-\ivz. 
\item For all $x$, $1-xa$ is \ivz.  
\item For all $x$, $1-ax$ is left-\ivz.
\item For all $x$, $1-ax$ is right-\ivz. 
\item For all $x$, $1-ax$ is \ivz.  
\item For all $x,y$, $1-xay$ is \ivz.  
\end{itemize}
The \elts $a$ that satisfy these \prts form a 
        two-sided \id 
called the Jacobson radical of $\gA$.

}
\end{exercise}

\vspace{-1em}

\pagebreak

\begin{exercise}\label{exoLocalFreenessLemma}
 {(A freeness lemma)} 
{\rm
Let $(\gA, \fm)$ be an integral \alo with residual field $\gk$, with quotient field $\gK$. Let $E$ be a \tf \Amoz; suppose that the \kev $E/\fm E = \gk \otimes_\gA E$ and the \Kev $\gK \otimes_\gA E$ have the same dimension $n$. Show that $E$ is a free $\gA$-module of rank $n$. \\
Better: if $(x_1, \ldots, x_n) \in E^n$ is a residual basis, it is an $\gA$-basis of $E$.
}
\end{exercise}

\vspace{-1em}
\begin{exercise}\label{exoNakayamaAppli} {(A consequence of Nakayama's lemma)}
\\
{\rm
Let  $E$ be a \pf \Amo and $a \in \Rad(\gA)$ be an $E$-regular \eltz.
Suppose that the $\gA/a\gA$-module $E/aE$ is free of rank $n$.
Show that $E$ is free of rank $n$. More \prmtz, let $e_1$, \ldots, $e_n \in E$, if $(\overline {e_1}, \ldots, \overline {e_n})$ is an~$\gA/a\gA$-basis of $E/aE$, then $(e_1, \ldots, e_n)$
is an $\gA$-basis of~$E$.
}
\end{exercise}

\vspace{-1em}
\begin{exercise}
\label{exoItfmonogene}
{\rm  Let $\gA$ be a \aloz. 
Prove the following statements.
If $\gen{b}=\gen{a}$,  there exists an \iv \eltz~$u$ such that $ua=b$. If $\fa=\gen{\xn}=\gen{a}$, there exists an index $i$ such that $\fa=\gen{x_i}$.
}
\end{exercise}

\vspace{-1em}
\begin{exercise}
\label{exoZeroSimpleIntersComp}
{\rm Give a detailed direct proof 
of \thref{thJZS} when $n=s$.  
}
\end{exercise}

\vspace{-1em}
\begin{exercise}
\label{exopropZerdimLib}
{\rm Here certain items of \thref{propZerdimLib} are revisited, now supposing that the \ri $\gA$ is  \plcz. The reader is invited to provide \dems which are independent from the results obtained for \algbsz.

 \emph{1.} Every  \ptf \Amo is quasi-free.

 \emph{2.} Every matrix $G\in\gA^{q\times m}$ of rank $\geq k$ is \eqve to a matrix

\snic{
\cmatrix{
    \I_{k}   &0_{k,m-k}      \cr
    0_{q-k,k}&  G_1      },}

with $\cD_r(G_1)=\cD_{k+r}(G)$ for all $r\geq 0$.
The matrices are \elrt \eqves if $k<\sup(q,m)$.

 \emph{3.} Every \mpf \lot generated by $k$ \elts is generated by $k$ \eltsz.
}
\end{exercise}

\vspace{-1em}
\begin{exercise}\label{exoSLnEn} 
(If $\gA$ is local, $\SLn(\gA)=\En(\gA)$)
\\
{\rm  Let $\gA$ be a \aloz. Show that every matrix $B\in\SLn(\gA)$ is produced from \elr matrices (in other words,~$B$ is \elrt \eqve to the matrix~$\In$). Inspiration may come from the \dem of the local freeness lemma. See also Exercise~\ref{exoLgb3}.
}
\end{exercise}

\vspace{-1em}
\begin{exercise}
\label{exoNbgenloc}
{\rm  \emph{1.}
Prove that a \tf \Amo $M$ is locally generated by $k$ \elts (\Dfnz~\ref{deflocgenk})
\ssi $\Vi_\Ae {k+1}M=0$. Inspiration may come from the case $k=1$  treated in \thref{propmlm}.
 
\emph{2.} Deduce that the annihilator $\Ann\big(\Vi_\Ae {k+1}M\big)$ and the \idfz~$\cF_k(M)$ have the same radical.
}
\end{exercise}

\vspace{-1em}
\begin{exercise}\label{exoVariationLocGenerated}
{(Variation on the locally generated theme)}\\
{\rm
Let $M$ be a \tf \Amoz, with two \sgrs $(\xn)$ \hbox{and $(\yr)$} with $r \le n$. We want to explicate a family $(s_I)$ of ${n \choose r}$ \ecoz, indexed by the $I \in \cP_{r,n}$, such that $s_IM \subseteq \gen{(x_i)_{i\in I}}$. Note that over each localized \ri
 $\gA[s_I^{-1}]$, the module $M$ is generated by the $(x_i)_{i\in I}$'s.

\emph {1.}
Let $A$ and $B\in \Mn(\gA)$. 
\begin{enumerate}\itemsep=0pt
\item [a.]
Explicate the membership

\snic {
\det(A+B) \in \cD_{n-r}(B) + \cD_{r+1}(A).
}
\item [b.]
Deduce that $1 \in \cD_{n-r}(\In-A) + \cD_{r+1}(A)$. 
\item [c.]
In particular, if $\rg(A) \le r$, then $\rg(\In-A) \ge n-r$.
\item [d.]
Let $a_1$, \dots,  $a_n \in \gA$, $\pi_I = \prod_I a_i$, $\pi'_J = \prod_J (1-a_j)$.\\
Show that the $(\pi_I)_{\#I=r+1}$'s and $(\pi'_J)_{\#J=n-r}$'s form a \sys of ${n+1 \choose r+1}$ \com \eltsz.
\end{enumerate}

\emph {2.}
Prove the result stated at the beginning of the exercise by making the family $(s_I)$ explicit.

\emph {3.}
Let $E$ be a \tf \Amo locally generated by~$r$ \eltsz. For any \sgr $(\xn)$, there exist  \ecoz~$t_j$ such that each of the localized modules 
$E_{t_j}$ is generated by $r$ \elts among the~$x_i$'s.

\emph {4.}
Let $E = \gen {\xn}$ be a \tf \Amo and $A \in \Mn(\gA)$ satisfying $\ux\,A = \ux$ with $\rg(A) \le r$. Show that $E$ is locally generated by~$r$ \eltsz.
Study a converse.
}

\end{exercise}

\vspace{-1em}
\begin{exercise}
\label{exoMorAdcp}
{\rm If $\gA$ and $\gB$ are two \dcp \ris we say that a \ri \homo $\varphi:\gA\to\gB$ is a \emph{\dcp \ri morphism} if, for all $a$, $b\in\gA$ satisfying $b(1-ab)=0$ and $a(1-ab)\in\Rad\gA$, we have in~$\gB$, \hbox{with $a'=\varphi(a)$} and $b'=\varphi(b)$,   $b'(1-a'b')=0$ and $a'(1-a'b')\in\Rad\gB$ (cf. Proposition~\ref{prop1DecEltAnneau}).
\begin{itemize}
\item [\emph{1.}] Show that $\varphi$ is a \dcp \ri morphism \ssiz
$\varphi(\Rad\gA)\subseteq\Rad\gB$.
\item [\emph{2.}] Study the injective and surjective \dcp \ri morphisms. In other terms, precise the notions of a \dcp sub\ri (considered as a single word) and of a \dcp quotient \riz. 
\end{itemize}
\index{morphism!\dcp \ri ---} 
}
\end{exercise}

\vspace{-1em}
\begin{exercise}
\label{exoMachDec} (\Elr local-global machinery of \dcp \risz)
\\
{\rm
The fact that one can systematically split a \dcp \ri into two components
leads to the following \gnl method.
\\
 {\it Most of the \algos that work with the \dcd \alos can be modified to work with the \dcp \risz, by splitting the \ri into two components 
each time the \algo written for the \dcd \alos
uses the test
\gui{is this \elt \iv or in the radical?}
In the first component
 the \elt in question is \ivz, in the second it is in the radical.}
\\
 Actually we rarely have  the occasion to use this \elr machinery, the main reason being that a more \gnl (but less \elrz) \lgb machinery applies with an arbitrary \riz, as it will be explained in Section~\ref{secMachLoGlo}.\imlb

}
\end{exercise}

\vspace{-1em}
\begin{exercise}
 \label{exopropAlgb1} (\Pol \lot representing an \invz, Lemma~\ref{lempropAlgb1})\\
{\rm Item \emph{3}  of this exercise gives a reinforced version of Lemma~\ref{lempropAlgb1}. The approach used here is due to Lionel Ducos.
\\
Let $\gA$ be a \riz, $d \in \NN$ and $e = d(d+1)/2$.

\emph{1.}
Here,  $s$ is an \idtr over $\ZZ$. Construct $d+1$ \pols $a_i(s) \in \ZZ[s]$ \hbox{for $i \in \lrb{0..d}$}, satisfying for every $P \in \AuX = \gA[\Xn]$ of degree $\le d$:
$$
s^e P(s^{-1}\uX) = a_0(s)P(s^0\uX) + a_1(s)P(s^1\uX) + \cdots + a_d(s)P(s^d\uX). 
\leqno(\star_d)
$$

\emph{2.}
For $s \in \gA$, $\ux \in \gA^n$ and $P\in\AuX$ of total degree $\leq d$, show that

\snic{
s^eP(\ux/s) \in\gen{P(\ux),P(s\ux), \ldots,P(s^{d}\ux)} \subseteq \gA.
}

\emph{3.} 
Let $S$ be a \mo and $P\in\AuX$. Suppose that~$P$ represents an \inv in $\gA_S$. Show that $S$ meets the \id generated by the values of $P$.

} 
\end{exercise}

\vspace{-1em}
\begin{exercise}
 \label{exoLgb2} 
 {\rm 
 (See also Exercise \ref{exoUAtoUB})
Let $\gA$ be a \algb and $M$ \hbox{an \Amoz}.
 
 \emph{1.} For every \id $\fa$,  the canonical \homo $\Ati\to(\gA\sur{\fa})\eti$ is surjective.
 
 \emph{2.} If $x$, $y\in M$ and $\gA x=\gA y$, there exists an \inv $u$ such that~$x=uy.$

} 
\end{exercise}

\vspace{-1em}
\begin{exercise}\label{exoLgb3} 
(If $\gA$ is \lgbz, $\SLn(\gA)=\En(\gA)$)\\
{\rm
Let $\gA$ be a \algbz, and $(a_1,\ldots,a_n)$ a \vmd ($n\geq2$).
 
 \emph{1.} Show that there exist $x_2$, \ldots, $x_n$ such that $a_1+\sum_{i\geq2}x_ia_i\in\Ati$.
\\ 
 \emph{2.} Deduce (for $n\geq2$) that every \vmd transforms into the vector $(1,0,\ldots,0)$ by \elr manipulations.
\\ 
 \emph{3.} Deduce that $\SLn\gA=\En\gA$.
} 
\end{exercise}

\vspace{-1em}
\begin{exercise}
\label{exo1semilocal} (\Slgbs \risz, 1)%
\index{ring!semi-local ---}%
\index{semi-local!\riz}%
\index{ring!strict semi-local ---}%
\index{strict semi-local!\riz}
\\
{\rm\emph{1.}
For a \ri $\gB$, prove that \propeq 

\vspace{-.5em}
\begin{itemize}\itemsep.0em
\item [~\emph{a.}]
If $(x_1, \ldots, x_k)$ is \umdz, there exists a \sys of \orts \idms $(e_1, \ldots, e_k)$  such that $e_1 x_1+\cdots + e_k x_k$ is \ivz.

\item [~\emph{b.}]
Under the same hypothesis, there exists a splitting $\gB \simeq \gB_1 \times \cdots \times \gB_k$ such that the component
 of $x_i$ in $\gB_i$ is \iv for $i \in \lrbk$.

\item [~\emph{c.}]
Same as in \emph{a}, but with $k = 2$.

\item [~\emph{d.}]
For all $x\in\gB$, there exists an \idm $e\in\gB$ such that $x+e$ is \ivz.\end{itemize}

\vspace{-.5em}
Note that at item \emph{a}, $(e_1,\ldots,e_k)$ is a \sfio since $1\in\gen{e_1,\ldots,e_k}$.\\
The \ris satisfying these \eqves \prts have been called \gui{clean rings} in \cite[Nicholson]{Nic}. 
\index{ring!clean ---}\index{clean!ring}

\emph{2.}
Clean \ris are stable under quotient and under finite product. Every \alo is clean.
\\
\emph{3.}
If $\gB\red$ is clean, the same goes for~$\gB$. Deduce that a \zed \ri is clean.
\\
\emph{4.}
If $\gB\red$ is clean,
 $\gB$ lifts the \idms of $\gB/\Rad\gB$.

We say that a \ri $\gA$ is \emph{\slgbz} if the \ri $\gB=\gA\sur{\Rad\gA}$ is clean. We say that it is \emph{\smlz} if it is \slgb and if $\BB(\gA\sur{\Rad\gA}\!)$ is a bounded \agBz.
}
\end{exercise}

\vspace{-1em}
\pagebreak

\begin{exercise}
\label{exo2semilocal} (\Slgbs \risz, 2) 
{\rm Prove the following statements.

 \emph{1.}  
A \alo is \smlz. 
\\ 
 \emph{2.}  
A \slgb and residually connected \ri is local.
\\ 
 \emph{3.} 
A \plc \ri is \slgbz.
\\ 
 \emph{4.}  
A \slgb \ri is \lgbz.
\\ 
 \emph{5.}  
The \slgbs \ris are stable under quotient and under finite product. 
\\ 
 \emph{6.} In \clamaz, a \ri is \sml \ssi it has a finite number of \idemasz.  
}
\end{exercise}

\vspace{-1em}
\begin{exercise}
\label{exoNagatalocal} 
(\Prts of the Nagata \riz)
{\rm See also Exercise \ref{exoPrufNagata}.\\ 
Let $\gA$ be a \ri and $U\subseteq\AX$ be the \mo of primitive \polsz.\\
Let $\gB=U^{-1}\AX=\gA(X)$ be the Nagata \ri of $\AX$.
 
 \emph{0.} 
Give a direct \dem of the fact that $\gB$ is \fpt over $\gA$.
 
 \emph{1.} 
$\gA\cap\gB\eti=\Ati$.
 
 \emph{2.} 
$\Rad\gA= \gA\cap\Rad\gB$ and $\Rad\gB=U^{-1}(\Rad\gA)[X]$.
 
 \emph{3.} 
$\gB\sur{\Rad\gB}\simeq  (\gA\sur{\Rad\gA}\!)(X)$.
 
 \emph{4.} 
If $\gA$ is local (resp.\ local and \dcdz), then $\gB$ is local (resp.\ local and \dcdz).
 
 \emph{5.} 
If $\gA$ is a field (resp.\ a \cdiz), then $\gB$ is a field (resp.\ a \cdiz).

} \end{exercise}

\vspace{-1em}
\begin{exercise}
 \label{exo2Nagata} (Nagata \ri with several \idtrsz)\\
 {\rm  Let $U$ be the set of primitive \pols of $\gA[X,Y]$.
 
 \emph{1.} Show that $U$ is a filter. 
\\  
Let $\gA(X,Y)=U^{-1}\gA[X,Y]$, we call it the \emph{Nagata \ri of $\gA[X,Y]$}.
 
 \emph{2.} Show that the canonical map $\gA[X,Y]\to\gA(X,Y)$ is injective
 and that we have a natural \iso $\gA(X,Y)\simarrow \gA(X) (Y)$.
 
 \emph{3.} Generalize the results of Exercise~\ref{exoNagatalocal}.
 } \end{exercise}

\vspace{-1em}
\begin{exercise}\label{exoAlgMon}
{(Algebra of a \mo and binomial \idsz)}\index{algebra!of a monoid}\\
{\rm Let $(\Gamma,\cdot,1_{\Gamma})$ be a commutative \mo denoted multiplicatively, and $\gk$ be a commutative \riz. 
\\
\hbox{The \emph{\alg of $(\Gamma,\cdot,1_{\Gamma})$ over $\gk$}}, denoted by $\gk[(\Gamma,\cdot,1_\Gamma)]$ or simply $\kGa$, is formed from the free \kmo over $\Gamma$ (if $\Gamma$ is not assumed to be discrete, see Exercise~\ref{propfreeplat}). If $\gk$ is nontrivial, we identify every \elt $\gamma$ of $\Gamma$ with its image in the free module. In case of doubt regarding $\gk$, we should denote by $1_{\gk}\gamma$ instead of $\gamma$ this \elt of~$\kGa$.
\\
The product law $\times$ of $\kGa$ is obtained by letting $\gamma\cdot\gamma'=\gamma\times \gamma'$ and by extending by $\gk$-bilinearity.   
Note that $1_{\gA}1_\Gamma=1_{\gk[\Gamma]}$. In practice, we identify $\gk$ with a sub\ri of $\kGa$, and we identify the three $1$'s above.
\begin{enumerate}
\item Prove that the \klg $\kGa$, considered with the map  

\snic{\iota_{\gk,\Gamma}:\Gamma\to\kGa, \, \gamma\mapsto 1_{\gk}\gamma,}

gives the solution to the \uvl \pb summarized in the picture below. 
\end{enumerate}

To sum up, we say that \emph{$\kGa$ is the \klg freely generated by the multiplicative \mo $\Gamma$}.
\vspace{-1.2em}
\PNV{\Gamma}{\iota_{\gk,\Gamma}}{\psi}{\kGa}{\theta}{\gL}{commutative \mos~~~~~}{\mos morphisms}{\klgsz}

\vspace{-2em}

When the law of $\Gamma$ is denoted additively, we denote by $X^{\gamma}$ the \elt of $\kGa$ image of $\gamma\in\Gamma$ such that we now have the natural expression $X^{\gamma_1}X^{\gamma_2}=X^{\gamma_1+\gamma_2}$.
\\
For example, when $\Gamma=\NN^{r}$ is the additive \mo freely generated by a set with $r$ \eltsz, we can see the \elts of $\NN^{r}$ as multiexponents and $\kGa=\gk[(\NN^{r},+,0)]\simeq \kXr$. Here  $X^{(m_1,\dots,m_r)}=X_1^{m_1}\cdots X_r^{m_r}$.
\\
When $\Gamma=(\ZZ^{r},+,0)$, we can again see the \elts of $\ZZ^{r}$ as multiexponents and $\gk[\ZZ^{r}]\simeq \gk[\Xr,\fraC1{X_1},\dots,\fraC1{X_r}]$ as the Laurent \pol \riz.
 
Now suppose that $(\Gamma,\cdot,1)$ is a \mo given by \gtrs and relations.
Let~$G$ be the set of the \gtrsz.  
\\
The relations 
are of the form $\prod_{i\in I}  g_i^{k_i}=\prod_{j\in J}  h_j^{\ell_j}$ for finite families 

\snic{(g_i)_{i\in I}$ and $(h_j)_{j\in J}$ in $G$,
and $(k_i)_{i\in I}$ and $(\ell_j)_{j\in J}$ in $\NN.}

Such a relation
 can be encoded by the pair $\big((k_i,g_i)_{i\in I},(\ell_j,h_j)_{j\in J}\big)$.
\\
If we hope to control things, $G$ and the set of relations 
better be enumerable and discrete. From the point of view of the computation, 
the central role
is taken up by the finite presentations.
\\
{\bf Notation.} To visualize a finite \pnz, for instance with $G=\so{x,y,z}$ and relations
 $xy^{2}=yz^{3}$, $xyz=y^{4}$ we write in multiplicative notation

\smallskip \centerline{\fbox{$\Gamma=_\mathrm{CM}\scp{x,y,z}{xy^{2}=yz^{3},xyz=y^{4}}\qquad$(*)}, }

 and in additive  notation

\smallskip \centerline{\fbox{$\Gamma=_\mathrm{CM}\scp{x,y,z}{x+2y=y+3z,x+y+z=4y}$}.}

The index $\mathrm{CM}$ is added for \gui{commutative \moz.}  

\vspace{-.4em}
\begin{enumerate} 
\setcounter{enumi}{1}
\item Show that $\kGa\simeq\gk[(g)_{g\in G}]/\fa$, where $\fa$ is the \id generated by the differences of \moms $\prod_{i\in I}{g_i}^{k_i}-\prod_{j\in J}
{h_j}^{\ell_j}$
(for the relations $\prod_{i\in I}  g_i^{k_i}=\prod_{j\in J}  h_j^{\ell_j}$ given in the \pn of $\Gamma$). Such an \id is called a \emph{binomial \idz}.
With the example $(*)$ above, we can therefore write

\smallskip \centerline{\fbox{$\kGa=_{\gk\mathrm{-\algsz}}\scp{x,y,z}{xy^{2}=yz^{3},xyz=y^{4}}$\qquad(**)}.}

In other words,  $\Gamma=_\mathrm{CM}\scp{thingy}{bob}$ implies $\kGa=_{\gk\mathrm{-\algsz}}\scp{thingy}{bob}$.
\end{enumerate}

}

\end {exercise}

\vspace{-1.1em}
\pagebreak

\sol

\exer{exoRadNonCom}
\emph{1.} If $c$ is right-\iv and left-\iv then it is the unique left-\inv because $c'a=1$ implies $c'=c'ac=c$. \\
Conversely, since $ca=1$, we have $(c+1-ac)a=ca+a-aca=1$. 
Therefore $c+1-ac$ is a left-\invz, and if there is uniqueness, $1-ac=0$.

 \emph{2.} We check that $v=1+bua$ suits.

\emph{3.} If $u(1-xa)=1$, then $u=1+uxa$, therefore it is left-\ivz. Thus $u$ is right- and left-\ivz, and so is $1-xa$.


\exer{exoLocalFreenessLemma}
Let $x_1$, \ldots, $x_n\in E$ such that $(\ov {x_1}, \ldots, \ov {x_n})$ is a $\gk$-basis of~$E/\fm E$. By Nakayama, the $x_i$'s generate $E$. Let $u : \Ae n \twoheadrightarrow E$ be the surjection 
$e_i \mapsto x_i$. By \eds to $\gK$, we obtain a surjection $U : \gK^n \twoheadrightarrow \gK\otimes_\gA E$ between two \evcs of same dimension $n$, thus an \isoz. 
\\
\Deuxcol{.7}{.15}{ Since $\Ae n\hookrightarrow \gK^n$, we deduce that $u$ is injective.
Indeed, if $y \in \Ae n$ satisfies $u(y) = 0$, then $1 \otimes u(y) = U(y) = 0$ in $\gK\otimes_\gA E$, therefore $y = 0$, cf.\, the diagram on the right.}
{$
\xymatrix @R = 7pt {
\Ae n \ar@{(->}[d] \ar@{->>}[r]_u & E\ar[d] \\
\gK^n \ar@{->>}[r]^{U~~}              & \gK\otimes_\gA E \\
}
$}

Recap: $u$ is an \iso and $(x_1, \ldots, x_n)$ is an $\gA$-basis of $E$.

\exer{exoNakayamaAppli}
By Nakayama, $(e_1, \ldots, e_n)$ generates the $\gA$-module $E$.\\
Let $L = \Ae n$ and $\varphi : L \twoheadrightarrow E$ be the (surjective) \ali that transforms the canonical basis of $L$ into $(e_1, \ldots, e_n)$. By hypothesis, $\overline\varphi : L/aL \to E/aE$ is an \isoz.
Let us show that $\Ker\varphi = a\Ker\varphi$.
Let $x \in L$ with $\varphi(x) = 0$; \hbox{we have $\overline\varphi(\overline x) = 0$}, so $\overline x = 0$, \cad $x \in aL$, say $x = ay$ with $y \in L$. But $0 = \varphi(x) = a\varphi(y)$ and $a$ being $E$-regular, $\varphi(y) = 0$. We indeed have $\Ker\varphi \subseteq a\Ker\varphi$.
Since $E$ is \pfz, $\Ker\varphi$ is \tfz, and we can apply Nakayama  to the \egt $\Ker\varphi = a\Ker\varphi$. We obtain $\Ker\varphi = 0$: $\varphi$ is an \isoz.

\exer{exoVariationLocGenerated} 
\emph {1a., b., c.}
The idea is to develop $\det(A+B)$ as a multi\lin function of the columns of $A+B$. The result is a sum of $2^n$ \deters of matrices obtained by mixing columns $A_j$, $B_k$ of $A$ and $B$.
We write

\snic {
\det(A_1+B_1, \dots, A_n+B_n) = \sum_{2^n} \det(C_1, \dots, C_n) 
\quad \hbox {with $C_j = A_j$ or $B_j$.}
}

For $J\in\cP_n$, let $\Delta^{\rm col}_J$ be the \deter when $C_j = B_j$ \hbox{for $j \in J$} \hbox{and $C_j = A_j$} otherwise. With this notation, we therefore have

\snic {
\det(A+B) = \sum_J \Delta^{\rm col}_J.
}

If $\#J \ge n-r$, then $\Delta^{\rm col}_J \in \cD_{n-r}(B)$; otherwise $\#\ov {J} \ge r+1$ and so $\Delta^{\rm col}_J \in \cD_{r+1}(A)$.

\emph {1d.}
Consider $A = \Diag(\an)$.

\emph {2.}
We write $\ux=\uy\,U$ with $U \in \gA^{r\times n}$, $\uy=\ux\,V$ with $V \in
\gA^{n\times r}$. \\
Let $A = VU$, $B = \In-A$. We have \framebox
[1.1\width][c]{$\ux\,B = 0$} and $\rg(B) \ge n-r$ since $\rg(A) \le r$.
The framed \egt shows, for $I \in \cP_{r,n}$ and $\nu$ minor of $B$ over the rows of $\ov I$, the inclusion $\nu M \subseteq \gen{(x_i)_{i\in I}}$, and we are done because $1 \in \cD_{n-r}(B)$.

\Prmtz, let $\Delta^{\rm row}_J$ be the \deter of the \gui{mixed} matrix whose {\em rows} of index $i \in J$ are the corresponding {\em rows} of $B$ and the {\em rows} of index $i \in \ov J$ are those of $A$. For $J \supseteq \ov I$, $\Delta^{\rm row}_J$ is a \coli of minors of $B$ over the rows of $\ov I$. 

Thus let
$$\preskip-.4em \postskip.4em\ndsp 
s_I = \sum_{J \mid J \supseteq \ov I} \Delta^{\rm row}_J. 
$$
Then on the one hand, $s_IM \subseteq \gen {(x_i)_{i\in I}}$, and on the other, since
$\rg(B) \ge n-r$,
$$\preskip.4em \postskip.4em\ndsp 
1 = \sum_{I \in \cP_{r,n}} s_I. 
$$

\emph {3.}
Clear by using the successive \lons lemma (Fact~\ref{factLocCas}).

\emph {4.}
If a matrix $A \in \Mn(\gA)$ exists as indicated, the \dem of item \emph {2}  applies \hbox{with $B=\In-A$}.

The converse is problematic because the constraint $\rg(A) \le r$ is not \lin in the \coes of $A$. However, we succeed in reaching it for $r=1$ by other means, (see \thref{propmlm}).


\exer{exoMorAdcp}
\emph{1.}
The condition $b'(1-a'b')=0$ is obtained by $\varphi\big(b(1-ab)\big)=0$. \\
Suppose that $\varphi$ is a \dcp \ri morphism and let us show {that $\varphi(\Rad\gA)\subseteq\Rad\gB$}: let $a \in \Rad\gA$, then $b = 0$ (by uniqueness of $b$), \hbox{thus $b' = 0$} and $a'=a'(1-a'b')\in\Rad\gB$.\\
Conversely, suppose $\varphi(\Rad\gA)\subseteq\Rad\gB$. If $a$, $b\in\gA$ satisfy $b(1-ab)=0$ and $a(1-ab)\in\Rad\gA$, then $\varphi\big(a(1-ab)\big)=a'(1-a'b')\in\Rad\gB$.

\exer{exopropAlgb1}
\emph{1.} 
It is sufficient and necessary that the $a_i$'s satisfy the \egt $(\star_d)$ for the \moms of total degree $\leq d$. Let $M=M(\uX)=\uX^{\alpha}$ be such a \mom  with $\abs\alpha=j\leq d$. Since $M(s^{r}\uX)=s^{rj}M$, we want to reach

\snic {
s^e s^{-j}\,M = a_0(s) \,M + a_1(s) s^{j}\,\uX + \cdots + a_d(s)s^{dj}\,M, 
}

\cad after simplification by $M$ and multiplication by $s^j$

\snic {
s^e = a_0(s)s^j + a_1(s)s^{2j} + \cdots + a_d(s)s^{(d+1)j} =
\sum_{i=0}^{d} a_i(s)(s^j)^{i+1}.
}

Let us introduce the \pol $F(T) \in \ZZ[s][T]$ defined by $F(T) = T\,\sum_{i=0}^{d} a_i(s)T^{i}$. \\
Then $\deg_T F \le d+1$ and $F$ performs the interpolation $F(0) = 0$ and $F(s^j) = s^e$ for $j \in
\lrb{1..d}$. However, a \pol $F \in \ZZ[s][T]$ which satisfies this interpolation is the following
$$
F(T) = s^e - (s^0-T)(s^1-T)(s^2-T) \cdots (s^d-T). \leqno (\#_d)
$$
Full astern. Consider the \pol defined by the \egt $(\#_d)$.
It is of degree $d+1$ in $T$, null in $T=0$, therefore it is of the form

\snic {
F(T) = T\, \sum_{i=0}^{d} a_i(s)T^{i}, \;\;\;\; \hbox {with}\;\;
a_0(s),\, \ldots,\, a_d(s) \in \ZZ[s].
}

These \pols $a_i(s)$ have the desired \prtz.

\emph{2.}
The required membership is deduced from the \egt $(\star_d)$ by evaluating $\uX$ at $\ux$.

\emph{3.} Suppose that $P$ is of total degree $\leq d$.
The fact that $P(\ux/s) \in (\gA_S)\eti$ means, in~$\gA$, that $y = s^e P(\ux/s)$ divides an \elt $t$ of $S$.  By item~\emph{2}, $y$ is in the \id generated by the values of~$P$; the same goes for $t$.

\exer{exoLgb2}{
\emph{1.} Let $b\in\gA$ \iv modulo $\fa$. There exists an $a\in\fa$ such that $1\in\gen{b,a}$. The \pol $aT+b$ takes the \come values $a$, $a+b$, thus it represents an \inv $b'=at+b$. Then, $b'\equiv b \mod \fa$ with $b'$ \ivz.
\\
\emph{2.} We write $x=ay$, $y=bx$, so $(1-ab)x=0$. \\
Since $b$ is \iv modulo $1-ab$, there exists a $u\in\Ati$  such that $u\equiv b \mod {1-ab}$. Then $ux=bx=y$.
}

\exer{exo1semilocal}~\\
For two \orts \idms $e$, $e'$, we have $\gen {ex,e'x'} = \gen {ex + e'x'}$.  
\\
Therefore for $(e_1, \ldots, e_k)$, we have $\gen{e_1x_1 + \cdots + e_kx_k} = \gen {e_1x_1,\dots,e_kx_k}$.
\\
Thus, $e_1x_1 + \cdots + e_kx_k$ is \iv
\ssi $e_1x_1$, \ldots, $e_kx_k$ are \comz. Consequently, in the context of \emph{1a}, let $y_i \in \gen {x_i}$ with \com $(y_1,\ldots, y_k)$ (a fortiori $(x_1, \ldots, x_k)$ is \comz); if \idms $(e_1, \ldots, e_k)$ work for $(y_1, \ldots, y_k)$, they also work for $(x_1, \ldots, x_k)$. Even if  $x_i$ is replaced by $u_ix_i$. We will therefore be able to assume $\sum x_i = 1$.
\\
For two \idms $e$, $e'$, we have $e \perp e'$ \ssi $1-e$, $1-e'$ are \comz.

 \emph{1.}
\emph {c}  $\Rightarrow$ \emph {d.} By taking $x_1 = x$, $x_2 = 1+x$, $e = e_2 = 1-e_1$, we have $e_1x_1 + e_2x_2 = x+e$.
\\
\emph {d}  $\Rightarrow$ \emph {c.} We can assume $1 = -x_1 + x_2$; we let $x = x_1$.
\\
Then, $e + x = (1-e)x + e(1+x) = (1-e)x_1 + ex_2$.
\\
\emph{a}  $\Leftrightarrow$ \emph {b.} Easily obtained by letting $\gB_i = \aqo{\gB}{1-e_i}$.
\\
\emph {c}  $\Rightarrow$ \emph {a}  (or \emph {d} $ \Rightarrow$ \emph {a}).
 By \recu on~$k$. 
We can assume $1 = \sum_i x_i$: there exists some \idm $e_1$ such that $e_1x_1 + (1-e_1)(1-x_1)$ is \ivz.  We have $1 \in \gen {x_2, \ldots, x_k}$ in the quotient $\aqo{\gB}{e_1}$ that also possesses the \prtz~\emph{d}; therefore, by \recuz, there exists $(e_2, \ldots, e_k)$ in $\gB$ forming a \sfio in the quotient $\aqo{\gB} {e_1}$ with $e_2x_2 + \cdots + e_kx_k$ \iv in $\aqo{\gB}{e_1}$. Then, $(e_1, (1-e_1)e_2, \ldots, (1-e_1)e_k)$ is a \sfio of $\gB$ and $e_1x_1 + (1-e_1)e_2x_2 + \cdots + (1-e_1)e_kx_k$ is \iv in~$\gB$.

 \emph{2.}
Easy.

 \emph{3.}
Let $x \in \gB$; there exists some \idm $e \in \gB\red$ such that $e + \ov x$ is \iv in~$\gB\red$. We lift $e$ at some \idm $e' \in \gB$. Then, $e' + x$ lifts $e + \ov x$ so is \ivz.  Let $\gB$ be a \zed \riz; even if we need to replace $\gB$ by~$\gB\red$, we can assume that $\gB$ is reduced; if
$x \in \gB$, there exists some \idm $e$ such that $\gen {x} = \gen {1-e}$; then $e + x$ is \ivz.

 \emph{4.}
Let $a\in\gB$ be an \idm \elt in $\gB\sur{\Rad\gB}$ and $b=1-a$.  \\
Since
$\gen{a,b}=1$, there exist two \orts \idms $e$ and $f$ in $\gB$ such \hbox{that $ae+bf$} is \ivz. Since $\gen{e,f}=1$, we have $f=1-e$. Now, we reason in the quotient. The \sys  $(ae, bf, af, be)$ is a \sfioz. As $ae+bf$ is \ivz, we have $ae+bf=1$, \hbox{hence $af=be=0$}.  Finally, (in the quotient) $a=e$ and $b=f$.

\exer{exo2semilocal}
A \ri $\gA$ is local \ssi $\gA\sur{\Rad\gA}$ is local;
a \ri $\gA$ is \slgb \ssi $\gA\sur{\Rad\gA}$ is \slgbz.

\emph {1.}
A \alo satisfies item~\emph {1d}  of the previous exercise with $e=0$ or $e=1$. 

\emph {2.}
$\gA\sur{\Rad\gA}$ is connected, semi-local thus local (use item \emph {1d}  of the previous exercise knowing that $e=0$ or $1$); so $\gA$ is local.

\emph {4.}
Is proven for the residual \ri and results from the following observation.\\
If $f$ is a \pol in $n$ \idtrs and $(e_1,\ldots, e_k)$ is a \sfioz, then for $(\xk)$ in $\gA^n$, since the \evn \homo commutes with the direct products, we have the \egt
$$\preskip.2em \postskip.2em 
f(e_1 x_1+\cdots +e_k
x_k) = e_1 f(x_1)+\cdots +e_k f(x_k). 
$$
 
\emph {6.}
A \ri $\gA$ has a finite number of \idemas \ssi 
          it is the case for $\gA\sur{\Rad\gA}$.  
In \clamaz, $\gA\sur{\Rad\gA}$ is a finite product of fields.

\exer{exoNagatalocal} ~\\
In the following $f=\sum_ib_iX^i\in\AX$ and $g=\sum_ic_iX^i\in U$, with $1=\sum_ic_iu_i$.

\sni
\emph{0}. Let $\uT$ be a set of 
\idtrs over $\gA$ and $\gA(\uT)$ be the Nagata \riz.
 We know \hbox{that $\gA(\uT)$} is flat over $\gA$ and we show that every \sli over $\gA$ that admits a solution over $\gA(\uT)$ admits a solution over $\gA$.
Thus let the \sli  $A x = b$ with $A\in\Ae{n\times m}$ and $b \in \gA^n$. Suppose the existence of a solution over $\gA(\uT)$; it is of the form $P/D$ with $P \in \gA[\uT]^m$ \hbox{and $D \in \gA[\uT]$} being a primitive \polz. We therefore have $A\,P = D\,b$ over~$\gA[\uT]$.
\\
Let us write $P = \sum_\alpha x_\alpha \uT^\alpha$ with $x_\alpha \in \gA^m$ and $D = \sum_\alpha a_\alpha \uT^\alpha$ where the $a_\alpha \in \gA$ are \comz. The \egt $A\,P = D\,b$ gives $A\, x_\alpha = a_\alpha b$ for each $\alpha$.
\\
If $\sum u_\alpha a_\alpha = 1$, the vector $x = \sum_\alpha u_\alpha x_\alpha$ is a solution of the \sys $Ax = b$.

\sni\emph{1}. Let $a\in\gA$ be \iv in $\gB$. There exist $f,\,g$ such that $af=g$, so $a$ and $f$ are primitive: $a\in\Ati$.

\sni\emph{2}. Let us show $\Rad\gA\subseteq\Rad\gB$. Let $a\in\Rad\gA$, we want to show that $1+a (f/g)$ is \iv in $\gB$, \cad $g+af\in U$. We want $1\in\gen{(c_i+ab_i)_i}$;
but this \id contains $\sum_iu_i(c_i+ab_i)=1+az\in \Ati.$ \\
We therefore know that $\Rad\gB\supseteq U^{-1}(\Rad\gA)[X]$.
Let $h=\sum_{i=0}^na_iX^{i}$. Let us show \hbox{that $h\in\Rad\gB$} implies $a_n\in\Rad\gA$. 
\\
We will deduce by \recu that $h\in(\Rad\gA)[X]$.
\\
Consider $a\in\gA$, take $f=a$ and $g=X^n-a(h-a_nX^n)$.
Clearly \hbox{$g\in U$}, so $g+fh=(1+aa_n)X^n$ must be \iv in $\gB$, \cad $1+aa_n$ must be in~$\Ati$.

\exer{exoAlgMon}
\emph{(Algebra of a \mo and binomial \idsz)}\\
\emph{1.} First of all we prove that $\kGa$ is indeed a \klg and that $\iota_{\gk,\Gamma}$ is a \mo morphism. Then, if $\alpha:\Gamma\to\gA$ is a \mo morphism, there is a priori a unique way to extend it to a morphism $\wi\alpha$ of \klgs from $\kGa$ to~$\gA$: let 
$\wi\alpha\big(\sum_{\gamma\in I}a_{\gamma}\gamma\big)=\sum_{\gamma\in I}a_{\gamma}\alpha(\gamma)$ (here, $I$ is a finitely enumerated subset of~$\Gamma$). 
\\
We then prove that $\wi\alpha$ is indeed a morphism of \klgsz. The readers are invited to prove all the details when $\Gamma$ is not assumed to be discrete, by basing themselves on Exercise~\ref{propfreeplat}.

\emph{2.} This is a \gnl result of \uvl \algz,
 because here we are in the framework of purely equational \agq structures. To obtain a \klg by means of \gtrs and relations 
 given by \egts of \momsz, we can first construct the similarly defined \moz, then the \alg freely generated by this \moz.
\\
If we do not want to invoke such a \gnl result, we can simply observe that the computation procedures in $\kGa$ with 
$\Gamma=_\mathrm{CM}\scp{thingy}{bob}$ are identical to those in $\gA=_{\gk\mathrm{-\algsz}}\scp{thingy}{bob}$.



\Biblio

The reader will certainly find our will to give to the trivial \ri every \prt under the sun a little arbitrary, especially through our use of a weakened version of negation (cf. footnote~\ref{footnoteNegation} \paref{footnoteNegation}).
We hope to convince them of the practical use of such a convention by way of the examples.
On the proper use of the trivial \riz, see~\cite[Richman]{Ri88}.

The \gui{\dem by Azumaya} of the local freeness lemma \ref{lelilo}
 is extracted from the \dem of the Azumaya \tho III.6.2 in \cite{MRR},
in the case that concerns us here.
In other words, we have given the \gui{matrix} content of the \dem of the local freeness lemma in \cite{MRR}.

Monomial curves (example on \paref{exlcourbemonomiale}) are treated in \cite{Kun}, Chapter V,
Example~3.13.f.

Decomposed \ris play an important role in the classical theory of Henselian \alos for example in the works \cite{Ray} or \cite{Laf}.

A \algb is sometimes called a \gui{ring with many units} in the literature. 
\Algbs have been particularly studied in \cite[Estes \& Guralnick]{EG}. 
Other \gui{\ris with many units} have appeared long beforehand, under the terminology \gui{unit-irreducible rings} (see for example \cite{vdK}). Those are the \ris $\gA$ for which the following \prt is satisfied: if two \pols of $\AX$ represent an \invz, then their product also represents an \invz. \Egmt introduced were the \gui{primitive} or \gui{strongly U-irreducible} \ris which are the \ris for which the following \prt is satisfied: every primitive \pol represents an \invz. They are special \algbsz. In the \dem of Fact \ref{factLocNagata} we have shown that a Nagata \ri is always \gui{primitive.}%
\index{primitive!\riz}\index{ring!primitive ---}

Concerning the Nagata \ri $\gA(X)$, given Fact~\ref{factLocNagata} and the good \prts of \algbsz, it is not surprising that this \ri plays a crucial role for the uniform solution of \slis with parameters over a \cdi and more \gnlt for the uniform computations \gui{in a reasonable amount of time} over arbitrary commutative \ris (see~\cite[D\'{\i}az-Toca\&al.]{DiGL,DiGLQ}).

\newpage \thispagestyle{CMcadreseul}
\incrementeexosetprob


\chapter{\Mptfsz, 2}
\label{chap ptf1}\label{ChapThStrBa}
\perso{compil\'e le \today}
\minitoc

\subsection*{Introduction}
\addcontentsline{toc}{section}{Introduction}

Here we continue the study of \mptfs started in Chapter~\ref{chap ptf0}.

In Section~\ref{sec ptf loc lib} we readdress the question regarding the \carn of \mptfs as \lot free modules, \cad regarding the local structure \thoz.

Section~\ref{subsecCalRang} is dedicated to the \ri of ranks over~$\gA$.
In the usual theory in \clama the rank of a \mptf is defined as a \lot constant function over the Zariski spectrum.
Here we give an \elr theory of the rank which does not require \idepsz.

In Section~\ref{secAppliLocPtf} we give some  simple applications of the local structure \thoz.

Section~\ref{secGrassman} is an introduction to Grassmannians.

In Section~\ref{subsecClassifMptfs} we introduce the \gnl \pb of completely classifying \mptfs over a fixed \ri $\gA$.
This classification is a fundamental and difficult \pbz, which does not admit a \gnl \algq solution.

Section~\ref{secAppliIdenti} presents a nontrivial example for which this classification can be obtained.

\newpage
\section{The \mptfs are \lot free}
\label{sec ptf loc lib}
We continue the theory of \mptfs after Section~\ref{subsec det ptf}. We ask however that the reader forgets what was learnt in Section~\ref{secMPTFlocLib}:
the \carn by the \idfsz, the local structure \tho \vref{prop Fitt ptf 2} and the considerations regarding the rank linked to \idfs as well as \thref{corth.ptf.sfio} whose \dem depends on the local structure \thoz.

Actually, all of the results of Sections~\ref{subsec det ptf} and~\ref{sec ptf loc lib} could be obtained by \lon arguments at \eco since we have already obtained the local structure \tho for \mptfs (\thrfs{theoremIFD}{prop Fitt ptf 2}) by exterior \alg methods.

We nevertheless think that the \gui{more global} point of view developed in this chapter is itself interesting, and, in a way, simpler, as highlighted by the \elr \dem of the matrix \tho
\ref{th matproj} which summarizes (and specifies) all the previous structure \thosz. There also the exterior \alg is an indispensable tool, but it seems better used, in a less invasive way.

\subsect{Complements on exterior powers of a \mptfz}{Complements on exterior powers}
\label{subsec comp ext}

The following lemma is \imdz.
\begin{lemma}
\label{lem ptf RM2}
Let $P$ be a free \Amo of rank $h$ and $\varphi\in\End(P)$ be a \dig \endoz, with a matrix similar to $\Diag(\lambda_1, \ldots, \lambda_h)$, then for the \polfon of $\varphi$ we get

\snic{\rF{\varphi}(X)\eqdefi\det(\Id_{P[X]}+X\varphi)=(1+\lambda_1X)\,\cdots\,(1+\lambda_hX).}
\end{lemma}

We now establish the crucial result.
\begin{proposition}
\label{prop puissance ext} \emph{(Exterior powers)}\\
Let $P$ be a \mptfz.
%
\begin{enumerate}
\item The $k^{\rm th}$ exterior power of $P$, denoted by $\Al{k}P$, is also a \mptfz. 
If $P=\Im(F)$ for $F\in\GA(\gA)$, the module $\Al{k}P$ is (\isoc to) the image of the
\mprn $\Al kF$.
\item If $\varphi$ is an \endo of $P$, the \polfon $\rF{{\Al{k}\!\varphi}}(X)$ only depends on $k$ and on the \pol $\rF{\varphi}(X)$.
In particular, the \polmu of $\Al{k} P$ only depends on $k$ and on the \polmu of~$P$.
\item
\begin{enumerate}
\item If $P$ is of constant rank  $h< k$, the module $\Al{k}P$ is null.
\item If $P$ is of constant rank $h\ge k$, the module $\Al{k}P$ is of constant rank~$h \choose k$.
\item In this case, if $\varphi$ is an \endo whose \polfon is $\rF{\varphi}=(1+\lambda_1X)\,\cdots\,(1+\lambda_hX)$, we have 
$$\preskip.2em \postskip.0em
\rF{\Al{k}\!\varphi}(X)=
\prod_{1\leq i_1<\cdots<i_k\leq h}(1+\lambda_{i_1}\cdots\lambda_{i_k}X) .$$
\end{enumerate}
\item If a \mprn  $F$ has as its image a \mrc $k$, then $\cD_{k+1}(F)=0$.
\end{enumerate}
\end{proposition}
\begin{proof}
\emph{1.} Let $M$ and $N$ be two \Amos and consider the first exterior powers of their direct sum $M\oplus N$. By examining the universal \pb that the $k^{\rm th}$ exterior power of a module solves, we obtain the canonical \isos
$$\arraycolsep2pt\preskip.4em \postskip.2em
\begin{array}{rcl}
\Al{2}(M\oplus N) &\simeq& \Al{2}M \oplus(M\otimes N) \oplus \Al{2}N \\[1mm]
\Al{3}(M\oplus N)&\simeq&\Al{3}M \oplus \big(\,(\Al{2}M)\otimes N\big)  \oplus
\big(M\otimes (\Al{2}N)\big)  \oplus \Al{3}N, \end{array}$$
and more \gnlt
\begin{equation}\label{eqVik}\preskip.2em \postskip.4em\ndsp
\Al{m}(M\oplus N)\,\,\simeq\quad
\bigoplus_{k=0}^m \big((\,\Al{k}M)\otimes (\,\Al{m-k}N )\big)
\end{equation}
(with $\Al0M=\gA$ and $\Al1M=M$). In particular, if $P\oplus Q\simeq \Ae m $, $\Al kP$ is a direct summand in $\Al k\Ae m \simeq\gA^{m \choose k}$.
We also see that if  $P=\Im(F)$ for some \mprn $F$, $\Al kP$ is (\isoc to) the image of the \mprn $\Al kF$, because this matrix represents the \idt over $\Al kP$ and $0$ over all the other summands of the direct sum.

 \emph{2.} We can assume $P = \Im(F)$, where $F \in \GA_n(\gA)$, and $n \ge k$.
\\
We therefore have $P \oplus Q = \gA^n$ with $Q = \Ker(F)$.
The \endo $\varphi$ extends into an \endo $\varphi_1:\Ae n \rightarrow \Ae n $, null over~$Q$, with matrix~$H$  satisfying $FHF=H$, and we have $\rF{\varphi}(X)=\rF{\varphi_1}(X)=\det(\I_n+XH)$.\label{Preuvepoint2prop puissance ext}
Then, we see that~$\Al k\!\varphi_1$ is an extension of $\Al k\!\varphi$,
null over the terms distinct from $\Al kP$
in the direct sum explicated in the \dem of item \emph{1}.  The matrix of~$\Al k\!\varphi_1$ is none other than~$\Al kH$.
\\
We therefore want to show that $\det\big(\I_{n \choose k}+X\;\Al kH\big)$ only depends on~$k$ and on~$\det(\I_n+XH)$. We are therefore brought back to the case of a free module, and this case has been treated in Proposition~\ref{propPolCarPuissExt}.

 \emph{3.}
This item results from the previous one, since the \gui{\pro of rank $k$} case can be deducted from the \gui{free of rank $k$} case.
Note that items \emph{3a} and \emph{3b} both say that when $P$ is of constant rank $h$,  $\Al kP$ is of constant rank~$h \choose k$ (which is equal to $0$ if $h<k$).
We have only separated them in order to give the result in a more visible form.

 \emph{4.}
This is equivalent to the fact that $\Al{k+1}P$ is null, which is item~\emph{3a.}
\end{proof}

\rems  (Consequences of Proposition~\ref{prop puissance ext}.)
\label{rem puissance ext}

 1) Let $\rR{P}(X)=r_0+r_1X+\cdots+r_nX^n$. Each $r_hP$
is a \mrc $h$ over $A[1/r_h]$, which gives, as a consequence of item \emph{3}, for $k>0$,
$$\preskip-.20em \postskip.4em\ndsp 
\rR{\Al k(r_hP)}(X)
=X^{h \choose k} \, \hbox{ over }\,\gA[1/r_h]. 
$$
By writing $P=\bigoplus_hr_hP$ and $\gA=\prod_h \gA[1/r_h]$ we obtain
\[\arraycolsep2pt\preskip.4em \postskip.4em
\begin{array}{rcl}
 \rR{\Al k\!P}(X) & =  &  r_0+\cdots+r_{k-1}+ r_kX+\cdots+
r_{k+j}X^{k+j \choose k} + \cdots+
r_{n}X^{n \choose k} \\[-.2em]
  & =  & \som_{h=0}^n r_{h}X^{h\choose k} .
\end{array}
\]
We also have by convention $\Al0P=\gA$ and thus also $\rR{\Al0\!P}(X)=X$ (so that the previous formula applies it must be agreed that ${n \choose 0}=1$ for all~$n\geq0$).

 2) If we let $\Vi P$ be the exterior \alg of $P$, the reader will show by an analogous computation that

\snic{\rR{\Vi \!P}(X)= r_0 X +r_{1}X^2 +\cdots+ r_k X^{2^k}  + \cdots+
r_{n}X^{2^n}
.}

3) We can compute  $\rF{\Al k(\varphi)}$ from $\rF{\varphi}$ as follows. \\
Since $\rF{\varphi}(0)=1$ and $\deg(\rF{\varphi})\leq n$,
if $\psi$ is the \endo of $\Ae n $ having as its matrix the companion matrix $C$ of $X^n\rF{\varphi}(-1/X)$, we obtain $\rF{\varphi }=\rF{\psi}$. Therefore
$$ \preskip-.4em \postskip.0em
\rF{\Al k\!\varphi}= \rF{\Al k\!\psi} =
\det \big(\I_{n \choose k}+X\,\Al kC\big)
\eqno\eoq  $$

\medskip From the previous remarks we deduce the following proposition.
\begin{proposition}
\label{prop rank ext carac}
Let $P$ be a \mptfz, and $k\leq h$ be two integers $>0$.
\Propeq
\vspace{-2pt}
\begin{enumerate}
\item  The module $P$ is of constant rank $h$.
\item  The module  $\Vi P$ is of constant rank $2^h$.
\item  The module  $\Al k P$ is of constant rank $h \choose k$.
\end{enumerate}
With $h=0$, the \prts 1 and 2 are \eqvesz.
\end{proposition}

\subsec{Case of the modules of constant rank}
\label{subsec rank constant}
\begin{theorem}
\label{th rg const loc free}
Let $P$ be a \Amrc $h$ with~$n$ \gtrsz,
(\isoc to the) image of a \prr $F\in \GAn(\gA)$. Then the~$n\choose h$ principal minors $(s_i)$ of order $h$ of $F$ satisfy
\begin{itemize}
\item [--] $\sum_i s_i = 1$,  and
\item [--] each $\gA_{s_i}$-module  $P_{s_{i}}$ is free of rank $h$, the matrix $F$ seen as a matrix with \coes in $\gA_{s_i}$ is similar to the standard \mprn $\I_{h,n}$.
\end{itemize}
\end{theorem}
\begin{proof}
The sum of the principal minors $s_i$ of order $h$ of $F$ is equal to $1$ since~\hbox{$\det(\I_n+XF)=(1+X)^h$}.\\
Moreover, since every minor of order $h+1$ is null (Proposition~\ref{prop puissance ext}), we can apply the freeness  
lemma \ref{lem pf libre} 
to each localized module~$P_{s_{i}}$, 
which is \isoc to the image of the matrix $F$ seen as a matrix with \coes in~$P_{s_{i}}$ (by
Proposition~\ref{propPtfExt}).
\end{proof}

\rem
\label{rem rg const loc free}
In the previous \thoz, it is possible that $s_i$ is nilpotent for certain values of $i$, therefore that $\gA_{s_i}$ is trivial. The fact of not excluding these zero \lons
is inevitable when we do not dispose of a test to know whether an \elt of $\gA$ is nilpotent or not. This justifies the natural convention given in the remark on \paref{conven rgc}.
\eoe

\subsec{\Gnl case}
\label{subsec cas general}
\begin{theorem}
\label{th ptf loc free}
Let $P$ be a \ptf \Amo with $n$ \gtrsz. Then for each \idm $\ide_h(P)$ there exist $n\choose h$ \elts $(s_{h,i})$ of $\gA$ 
with the following \prts
\begin{enumerate}
\item [--] $\sum_i s_{h,i} = \ide_h(P),$
\item [--] each $\gA_{s_{h,i}}$-module  $P_{s_{h,i}}$ is free of rank $h$.
\end{enumerate}
In particular, for every \mptf with $n$ \gtrsz, there exist $2^n$ \eco $v_\ell$ such that each $P_{v_\ell}$ is free.
\end{theorem}
\begin{proof}
We first localize by inverting $\ide_h(P)$ to be reduced to \thref{th rg const loc
free}. We then localize a little more in accordance with the latter \thoz. Fact~\ref{factLocCas} regarding the successive \lons applies.
\end{proof}

The following \tho summarizes \thrfs{th rg const loc free}{th ptf loc free}, and the converse given by the \plgrf{plcc.cor.pf.ptf}.

\begin{theorem}
\label{thptfloli}
An \Amo  $P$ is \ptf \ssi there exist \eco $s_1$, \ldots, $s_\ell$ such that each $P_{s_i}$ is free over $\gA_{s_i}$.
It is \pro of rank $k$ \ssi there exist \eco $s_1$, \ldots, $s_\ell$ such that each $P_{s_i}$ is free of rank $k$ over $\gA_{s_i}$.
\end{theorem}

A practical form of \thref{th ptf loc free} is its matrix form. 

\begin{theorem} 
\emph{(Explicit matrix form of \thrfs{th.ptf.loc}{th.ptf.idpt})}
\label{th matproj}\label{propPTFDec} \\
 Let  $\gA$ be a \riz, $F\in \Mn(\gA)$ with $F^2=F$ and $P$ be the \mptf image of $F$ in $\Ae n$. We define the \elts $r_h$ of $\gA$ for $h\in\lrb{0..n}$ by the \egtsz

\snic{\rR{P}(1+X):=\det(\In+XF),\quad \rR{P}(X)=:r_0+r_1X+\cdots+r_nX^n.
}

We have the following results.
\begin{enumerate}
\item  The family $(r_h)_{h=0,\ldots,n}$ is a \sfio of $\gA$.
\item  For $h\in\lrb{0..n-1}$ and for any minor $u$ of order $h+1$ of $F$, we have $ r_hu=0 $.
\item  If the $t_{h,i}$'s are principal minors of order $h$  of  $F$, by letting~$s_{h,i}=r_ht_{h,i}$ we obtain the following,
\begin{enumerate}\itemsep0pt
\item [--] the sum (for fixed $h$) of the $s_{h,i}$'s is equal to $r_h$,
\item [--] each $\gA_{s_{h,i}}$-module  $P_{s_{h,i}}$ is free of rank $h,$
\item [--] the matrix $F$ is similar to the matrix $\I_{h,n},$ over $\gA_{s_{h,i}}$
\item [--] the $s_{h,i}$'s are \comz, \prmt $\som_{h,i}s_{h,i}=1$.
\end{enumerate}
\end{enumerate}
\end{theorem}


\rem
\label{rem th matproj}
\Thref{th matproj} summarizes \Thos \ref{th decomp ptf}, \ref{th rg const
loc free} and \ref{th ptf loc free} which have preceded it. 
It is even slightly more precise. Thus it is not uninteresting to provide a purely matrix \dem of it that concentrates all of the previous \dems together, especially since it is particularly \elrz.

\begin{Proof}{Matrix \dem of the matrix \thoz. }
\\
\emph{1.}
This results from  $\rR{P}(1)=1$ (obvious) and  $\rR{P}(XY)=\rR{P}(X)\rR{P}(Y)$ which becomes apparent as follows
$$\preskip3pt\postskip3pt
\arraycolsep2pt\begin{array}{rclc}
\rR{P}(1+X)\rR{P}(1+Y)&  = &\det(\In+XF)\det(\In+YF)   & =     \\[1mm]
\det\big((\In+XF)(\In+YF)\big)&  = & \det(\In+(X+Y)F+XYF^2)  & =     \\[1mm]
\det(\In+(X+Y+XY)F)&  = & \rR{P}\big((1+X)(1+Y)\big).  &
\end{array}$$

\noi \emph{2.} The matrix $r_hF$ has as its \polfon $\det(\In+r_hXF)$. In the \ri $\gA_{r_h}$, we have $1=r_h$ and 

\snic{\det(\In+r_hXF)=\det(\In+XF)=\rR{P}(1+X)=(1+X)^h.}

Within the \ri $\gA_{r_h}$ we are therefore reduced to proving item~\emph{2}  for the case where $r_h=1$ and $\det(\In+XF)=(1+X)^h$, which we assume from now on.
We must show that the minors of order $h+1$ of $F$ are all null.
The minors of order~$h+1$ are the \coes of the matrix~$\Al{h+1}F=G$.
Since~$F^2=F$, we also have~$G^2=G$.
Moreover, for any square matrix~$H$,
the \polcar of~$\Al{k}H$ only depends on~$k$ and on the \polcar of~$H$ (Proposition~\ref{propPolCarPuissExt}).
By applying this to compute the \polcar of~$G$, we can replace~$F$ with the matrix $\I_{h,n}$ which has the same \polcar as $F$. Since the matrix $\Al{h+1}\I_{h,n}$ is null, its \polcar is $X^{h+1 \choose n}$, so, by Cayley-Hamilton, the matrix $G$ is nilpotent, and since it is \idmez, it is null.

 \emph{3.} Results from \emph{1}, \emph{2}  and from the freeness lemma \ref{lem pf libre}.
\end{Proof}


\subsec{Modules of constant rank: some precisions}
\label{secDecEqdim}

The following two results are now easy and we leave their \dem as an exercise.

\begin{proposition}
\label{prop rgconstant local 2}\label{prop rgconstant local}
\emph{(\Pro modules of constant rank)}\\
For some \Amo $P$  \propeq
\begin{enumerate}
\item $P$ is \prcz~$h$.
\item \label{i2prop rgconstant local}
There exist \eco $s_i$ of $\gA$ such that each $P_{s_i}$ is free of rank $h$ over~$\gA_{s_i}$.
\item $P$ is \ptf and for every \elt $s$ of $\gA$, if $P_s$ is free over $\gA_s$,
it is of rank~$h$.
\item $P$ is \pfz, $\cF_h(P)=\gen{1}$ and  $\cF_{h-1}(P)=0.$
\item $P$ is \isoc to the image of a \mprn of rank~$h$.
\end{enumerate}
\end{proposition}
In addition, if $P$ is generated by $n$ \eltsz, the number of  \eco in item \emph{\ref{i2prop rgconstant local}} is bounded above by ${n\choose h}.$

\begin{proposition}
\label{prop sfio unic}{\em  (Localized modules of constant rank and uniqueness of the \sfioz)}\\
Let $P$ be a \ptf \Amoz. Let $r_h=\ide_h(P)$. Let $s$ be an \elt of $\gA$.
\begin{enumerate}
\item The localized module $P_s$ is \pro of rank $h$ \ssiz 
$r_h/1=1$ in~$\gA_s$, \cad if $r_hs^m=s^m$ in $\gA$ for some exponent~$m$.
\item If $s$ is an \idmz, that means that $r_h$ divides $s$, or yet again that $1-r_h$ and~$s$ are two \ort \idmsz.
\item Finally, if $ (s_0,\ldots,s_n) $ is a \sfio such that each $P_{s_h}$ is of rank $h$ over $\gA_{s_h}$, then $ r_h=s_h$ for each $h\in\lrbzn$.
\end{enumerate}
\end{proposition}

In the following proposition we make the link between our \dfn and the usual \dfn (in \clamaz) of a \pro module of rank $k$. The \dem of this \eqvc is however not \cov (nor can it be).

\begin{proposition}\label{prop ptfrangconstant}
Let $k$ be a non-negative integer, $P$ be a \mptf over a nontrivial \ri $\gA$ and $\fa$ be an \id contained in $\Rad\gA$. Then \propeq
\begin{enumerate}
\item [1.\phantom{$^*$}] $P$ is of rank $k$, \cad $\rR{P}(X)=X^k$
\item [2.$^*$] For all \idema $\fm$ of $\gA$, the vector space obtained from $P$ by extending the scalars to the residual field $\gA/\fm$ has dimension~$k$.
\item [3.\phantom{$^*$}]  $\rR{P}(X)\equiv X^k$ modulo $\fa[X]$.
\end{enumerate}
\end{proposition}
\begin{proof}
From a classical point of view, the implication  \emph{2}  $\Rightarrow$ \emph{3}  is \imdez; it suffices to recall that the intersection of the \idemas is the Jacobson radical of $\gA$.
Note that from a \cov point of view, condition \emph{2}  is a priori too weak, for lack of \idemasz.
\\
Moreover, \emph{1}  trivially implies \emph{2}  and \emph{3}. 
\\
Conversely, if $\rR{P}(X)=X^k$ modulo $\fa[X]$, since the \idms are always isolated (Lemma~\ref{lemIdmIsoles}), the \egt takes place in $\AX$.
\end{proof}

\vspace{-.7em}
\pagebreak

\begin{theorem}
\label{propRgConstant2} \emph{(Modules of constant rank $k$ as submodules of~$\Ae k$)}\\
Suppose that over $\Frac \gA$ every \mrc $k$ is free.
Then every \pro \Amo of constant rank $k$ is \isoc to a submodule of $\Ae k$.
\end{theorem}
\begin{proof}
By the \dlg lemma \ref{propIsoIm}, we can assume that the module is an image of a \prr $F\in\GAn(\gA)$ of rank $k$ and that there exists a matrix~$P$ in $\GLn(\Frac \gA)$ such that $PFP^{-1}=\I_{k,n}$. We \hbox{have $P=Q/a$} 
\hbox{with $Q\in\Mn(\gA)$} and~$a\in\Reg\gA$; thus $\det Q=a^n\det P$ is also \ndz in~$\gA$.  
We define a matrix~$Q_1$ as
$$\preskip-.4em \postskip.4em
Q\cdot F\,=\,\I_{k,n}\cdot Q\,=
\blocs{.8}{.6}{.8}{.6}{$\I_k$}{$0$}{$0$}{$0$}
\cdot  Q\,=\,
\blocs{1.4}{0}{.8}{.6}{$Q_1$}{}{$0$}{}\;.
$$
However, the image of $Q\cdot F$ is \isoc to the image of $F$ because $Q$ is injective, and the image of $\I_{k,n}\cdot Q$ is clearly \isoc to the image of~$Q_{1}= Q_{1..k, 1..n}$.
\end{proof}

\rem 1) The previous \tho applies to the \qiris and more \gnlt to every \ri $\gA$ such that $\Frac\gA$ is \zedz, or even simply \lgbz.
This is the case, for example, for reduced \fdis \coh \noes \ris (see \Pbmz~\ref{exoAnneauNoetherienReduit}).  
In \clama one proves that for every \noe \riz~$\gA$, $\Frac\gA$ is \rdt \zedz, so we can apply the \tho to it. We do not know of a \cov analogue of this \thoz. 
\\
2) For further details regarding the $k=1$ case see \thref{propRgConstant3}.

\subsec{Generic case} \label{subsec cas generique}

What do we call the \gnq case, regarding a \pro module with $n$ \gtrsz?
We consider the \ri
$$\Gn =\ZZ[(f_{i,j})_{i,j\in\lrbn}]/\cGn,\label{NOTAGN}$$
where the $f_{i,j}$'s are \idtrsz,  $F$ is the matrix $(f_{i,j})_{i,j\in\lrbn}$ and $\cGn$ is the \id defined by the $n^2$ relations 
 obtained when writing $F^2=F$.
With \coes in this \ri $\Gn$, we have the matrix $F$ whose image in $\Gn^{n}$ is what deserves to be called {\em  the \gnq \pro module with $n$ \gtrsz}.

Let us reuse the notations of \thref{th matproj} in this particular case. 
\\
Saying that
$r_hr_k=0$ in $\Gn$ (for $0\leq h\not= k\leq n$) signifies that we have a membership  
$$\preskip-.4em \postskip.6em
r_h({F})r_k({F})
\in \cGn\qquad (*)
$$
 in the \ri $\ZZ[(f_{i,j})_{i,j\in\lrbn}]$.
This implies an \ida that allows us to express this membership. This \ida is naturally valid in all the commutative \risz.
It is therefore clear that if the membership $(*)$ is satisfied in the \gnq case, it implies $r_hr_k=0$ for any \mprn over an arbitrary commutative \riz.

The same holds for the \egts $r_hu=0$ when $u$ is a minor of order $h+1$.

In short: if \thref{th matproj} is satisfied in the \gnq case, it is satisfied in every case.
As is often the case, we therefore observe that important \thos of commutative \alg do nothing other than affirm the existence of certain particular types of \idasz.

\section[The \ri of generalized ranks $\HO(\gA)$]
{\texorpdfstring{The semi\ri $\HOp(\gA)$, and the \ri of generalized ranks $\HO(\gA)$}
 {H0+(A),  and the \ri of generalized ranks H0(A)}}
\label{subsecCalRang}

For a free module, by passing from the rank $k$ to the \polmu $X^k$, we pass from the additive notation to the  multiplicative notation.
For a \gnl \mptfz, we can conversely consider 
a \gui{\gne rank} of the module, which is the (purely formal) logarithm in base
$X$ of its \polmuz. Although this is just a simple play on notation, it so happens that computations with the ranks are facilitated by it.
Let us explain how this works.

\subsection*{The semi\ri of ranks}
\label{SemiAnneauRangs}

Recall that we say that a \pol $R(X)=r_0+r_1X+\cdots +r_nX^n$ is multiplicative when $R(1)=1$ and $R(XY)=R(X)R(Y)$.
It amounts to the same to say that the $r_i$'s form a \sfioz, or that $R(X)$ is the \polmu of a \mptfz.

\begin{notation}
\label{notaHO+}
		{\rm  We denote by $\HOp (\gA)$ the set of \iso classes of quasi-free modules over $\gA$, and $[P]_{\HOp (\gA)}$ (or $[P]_\gA$, or even $[P]$) the class of such a module $P$ 
in $\HOp (\gA)$. The set $\HOp (\gA)$ is equipped with a \emph{semi\riz\index{semi-ring} structure\footnote{This means that the structure is given by an addition, commutative and associative, a multiplication, commutative, associative and distributive with respect to addition, with a neutral $0$ for the addition and a neutral $1$ for the multiplication.
For example $\NN$ is a semi\riz.}
for the inherited laws of $\oplus$ and $\otimes$: $[P\oplus Q]=[P]+[Q]$ \hbox{and $[P\otimes Q]=[P]\cdot[Q]$}. For an \idm $e$ we will also write $[e]$ instead of~$[e\gA]$, when the context is clear.
The neutral \elt for the multiplication is~$[1]$.}
}
\end{notation}

Every quasi-free module $P$ is \isoc to a unique module\footnote{We also have (Exercise \ref{exoSfio}) $P\simeq e_1\gA \oplus \cdots \oplus e_n\gA$ with \smash{$e_k=\sum_{j=k}^nr_j$}, and $e_k$ divides~$e_{k+1}$ for $1 \le k < n$.}
$$\preskip.4em \postskip.4em
{(r_1\gA)\oplus (r_2\gA)^2 \oplus \cdots \oplus
(r_n\gA)^n,}
$$
where the $r_i$'s are \ort \idmsz, because then $\ide_{i}(P)=r_i$.
We therefore have= $[P]=\sum_{k=1}^n k\,[r_k]$ and its \polmu is  
$$\preskip.4em \postskip.2em
\rR{P}(X)=r_0+r_1X+\cdots +r_nX^n
$$
with $ r_0 =  1 - (r_1 +\cdots + r_n)$.

But while $\rR{P\oplus Q}=\rR{P}\rR{Q}$, we have $[P\oplus Q]=[P]+[Q]$: this assures the passage from the multiplicative notation to the additive notation.
Thus the \gui{logarithm in base $X$} of the multiplicative \pol $r_0+r_1X+\cdots +r_nX^n$is defined as the \elt $\sum_{k=1}^n k\,[r_k]$ of $\HOp (\gA)$.

\begin{definition}\label{defiRang}
If $M$ is a \ptf \Amo we call \emph{(\gnez) rank} and we denote by $\rg_\gA(M) $ or $\rg(M)$ the unique \elt of $\HOp(\gA)$ which has the same \polmuz.\index{rank!(generalized) --- of a \mptfz}
\end{definition}
\label{NOTAHO+}

Thus if $\rR{M}(X)=r_0+r_1X+\cdots +r_nX^n$, then $\rg(M)=\sum_{k=1}^n k\,[r_k]$.

The zero module is \care by $\rg(M)=0$ (\thref{th ptf sfio}).

If $\gA$ is nontrivial, then $[1]\neq[0]$ and $\NN$ is identified with the subsemi\ri of~$\HOp(\gA)$ generated by $[1]$ by means of the injection $n\mapsto n\,[1]$.
The above \dfn therefore does not conflict with the notion of rank for \mrcsz, previously defined.

Also note that when $\gA$ is trivial we have $\HOp (\gA)={0}$: this indeed conforms with the convention according to which the zero module over the trivial \ri has for rank any integer, since in $\HOp (\gA)$, $k=0$, or if we prefer, the two \polmus $1$ and $X^k$ are equal over the trivial \riz.

\medskip
\rem A rule of practical computation regarding ranks is the following

\snic{[r]+[r']=[r+r'] \;\;\;\mathrm{if}\;\;\; rr'=0,}

\cad more \gnlt
\begin{equation}\label{eqsomHO}\preskip.2em \postskip.2em
[r]+[r']=[r\vu r']+ [r\vi r']=[r\oplus r']+2\,[r\vi r']
\end{equation}
where the laws $\vu,$ $\vi$ and $\oplus$ are those of the \agB of the \idms of the \riz:
$r\oplus r'=r+r'-2rr',$ $r\vu r'=r+r'-rr'$ and $r\vi r'=rr'$. Note that the two \idms $r\oplus r'$ and $r\vi r'$ are \ortsz, of sum~$r\vu r'$, and that the meaning of the \egt (\ref{eqsomHO}) is given by the following \isos
$$\preskip.2em \postskip.4em
{r\gA\oplus r'\gA \;\simeq\; (r\vu r')\gA\oplus (r\vi r')\gA \;\simeq\;
(r\oplus r')\gA\oplus \big((r\vi r')\gA\big)^2.} \eqno \eoq
$$

\subsection*{Exponential notation}
\label{ExpoRangs}

Please note that $a^n$ is the result of evaluating the multiplicative \polz~$X^n$ at the point $a$: $X^n(a)=a^{\log_X(X^n)}$.

Thus, for a multiplicative \pol $R(X)=\sum_{k=0}^ne_kX^k$, whose logarithm in base $X$ is the \elt $r=\sum_{k=0}^nk\,[e_k]$, we adopt the following legitimate notations 
\begin{itemize}
\item $a^r=\sum_{k=0}^ne_ka^k=R(a)$,
\item and for an \Amo $M$,  $M^r=\bigoplus_{k=0}^ne_kM^k$. 
\end{itemize}
This is not a fancy: we indeed have 
\begin{itemize}
\item $a^{r+r'}=a^r a^{r'}$, $a^{r r'}=(a^r)^{r'}$,
\item  $M^{r+r'}\simeq M^r \times M^{r'}$ and $M^{r r'}\simeq M^r \otimes M^{r'} \simeq (M^r)^{r'}$,
\end{itemize}
for arbitrary $r,r'$ in $\HOp\gA$.

\subsection*{Symmetrization}
\label{SymetrisationRangs}\rdb

The additive \mo  $\HOp(\gA)$ is regular: either by McCoy's lemma (Corollary \ref{corlemdArtin}), or by one of the two uniqueness \thos \ref{prop unicyc} and \rref{prop quot non iso},
or finally by item \emph{\iref{th ptf sfio item reg}}  of \thref{th ptf sfio}.

The semi\ri $\HOp (\gA)$  can therefore be considered as 
a subsemi\ri of the \ri obtained by symmetrization.
This \ri is called the \emph{ring of (\gnesz) ranks  of \mptfsz} over $\gA$,
and we denote it by $\HO(\gA)$.\index{ring of (\gnesz) ranks!over $\gA$}
 \label{notaHO}
 
Every \elt of  $\HO(\gA)$ is expressed in the form $\som_{k\in J} k \,[r_k]$
where the $r_k$'s are pairwise \ort \idms and $J$
is a finite subset of $\ZZ\setminus\so{0}$.
\\
The expression is unique in the following sense: 
if  $\som_{k\in J} k \,[r_k]=\som_{k\in J'}
k \,[r'_k],$ \hbox{then  $r_k=r'_k$} if $k\in J\cap J',$ and the others are null.

\subsection*{Multiplication of the ranks}
\label{MultiplicationRangs}

We have defined a multiplication on $\HOp\gA$, as the law inherited from the tensor product. This implies that for two \idms $e$ and $e'$ \hbox{we have $[e]\cdot [e']=[ee']$}. The other product computations are deducted from it by distributivity. Hence the following fact.

\begin{fact}\label{fact0ProdRangs}
The \elt $1$ is the only 
\iv \elt of $\HOp(\gA)$. 
\end{fact}
\begin{proof}
If $r=\sum_k k[r_k]$ and $s=\sum_k k[s_k]$, then $rs=\sum_kk(\sum_{i,j,ij=k}[r_is_j])$. By uniqueness of the expression, if $rs=1=1[1]$, then
$r_1s_1=1$ therefore $r_1=s_1=1$.
\end{proof}

We can ask ourselves which is the corresponding law on multiplicative \polsz; the readers will convince themselves that is it the law
$$\preskip.4em \postskip.4em 
\big(R(X),R'(X)\big)\mapsto R\big(R'(X)\big)=R'\big(R(X)\big). 
$$
We also have the following fact which stems from the upcoming Proposition~\ref{prop Produit tensoriel}.

\begin{fact}\label{factProdRangs}
If $P$ and $Q$ are two \mptfsz, then $P\otimes Q$ is a \mptf and $\rg(P\otimes Q)=\rg( P)\cdot\rg(Q)$.
\end{fact}

\subsection*{Order relation over ranks}
\label{Relation-dordreRangs}

The natural order relation associated with the \mo structure of $\HOp\gA $ is described in the following proposition.

\begin{propdef}
\label{def rank inferieur}~
\begin{enumerate}
\item For $s$, $t\in\HO\gA$ we define $s\leq t$ by $\exists r\in\HOp\gA $, $s+r= t$.
\item This relation 
gives to $\HO\gA$ an \emph{ordered \riz}\index{ring!ordered ---} structure,%
\footnote{This means that $\geq$ is a partial order relation \emph{compatible}
with the laws $+$ and $\times$, more precisely\\ 
$\bullet$ $1\geq 0$,\\ 
$\bullet$ $x\geq 0$ and $y\geq 0$
imply $x+y$ and $xy\geq 0$,\\ 
$\bullet$ $x\geq y \Leftrightarrow x-y\geq 0$.
} and $\HOp\gA$ is the non-negative 
subset of $\HO\gA$. 
\item \label{i3def rank inferieur}
Let $P$  and $Q$ be \mptfsz, \propeq
\begin{enumerate}
\item  $\rg(P) \leq \rg(Q)$.
\item $\rR{P}$ divides $\rR{Q}$ in $\AX $.
\item $\rR{P}$ divides $\rR{Q}$ in $\BB(\gA)[X]$.
\item For all $s\in \gA$, if $P_s$ and $Q_s$ are free, then the rank of $P_s$ is less or equal to that of $Q_s$.
\item For all $k>i, \; \; \ide_k(P)\cdot \ide_i(Q)=0.$
\item For all $ k,\; \;  \ide_k(P)\cdot\sum_{i\geq k}\ide_i(Q)=\ide_k(P).$
\end{enumerate}
\end{enumerate}
\end{propdef}

\exl
Let us suppose that $P\oplus R=Q$ and that we know the ranks of $P$ and $Q$, we want to compute the rank of $R$.\\
We have $\rg{P}=\sum_{i=0}^n i \,[r_i]$ and $\rg{Q}=\sum_{j=0}^m j \,[s_j]$.
We write
\[\arraycolsep2pt
\begin{array}{ccccccc}
 \rg{P} & = \big(\sum_{i=0}^n i \,[r_i]\big)\big(\sum_{j=0}^m \,[s_j]\big) & = &
 \sum_{i,j} i \,[r_is_j] & \leq\\[1mm]
 \rg{Q} & = \big(\sum_{j=0}^m j \,[s_j]\big)\big(\sum_{i=0}^n \,[r_i]\big) & = &
 \sum_{i, j} j \,[r_is_j]. \\
\end{array}
\]

The $r_is_j$'s form a \sfio and we obtain by subtraction, without having to think, the \egts 

\snac{\quad \rg(Q)-\rg(P)=\rg(R)=\som_{i\leq j}
(j-i) \,[r_is_j] =\som_{k=0}^m k \big(\som_{j-i=k}\,[r_is_j]\big). \hfill\eoq}

\medskip
In the remainder of the text, we definitively drop the use of the word \gui{\gnez}
when we speak of the rank of a \mptfz.

\medskip
\rem In \clamaz,  $\HO(\gA)$ is often defined as the \ri of \lot constant (\cad continuous) functions from $\SpecA$ to~$\ZZ$. A more detailed comment on the subject can be found on \paref{comHOclassique}.\eoe

\goodbreak
\subsection*{Other uses for the rank}
\label{AutresRangs}

\begin{notations}\label{notaRangs} ~
\begin{enumerate}
\item If $\varphi\in\Lin_{\gA}(P, Q)$ with $P$, $Q$ \ptfsz, and if $\Im\varphi$ is a direct summand
in $Q$ we will denote $\rg(\Im\varphi)$ by $\rg\varphi$.
\item If $p(X)$ is a pseudomonic \pol of $\AX$, we can define its degree $\deg p$ as an \elt of $\HOp (\gA)$.
\end{enumerate}
\end{notations}

For item \emph{1}  we have $\Ker \varphi$ which is a direct summand in $P$ 
and we obtain the \gnns of well-known \egts in the case of \evcs over a \cdi 

\snic{\rg(\Ker \varphi)+\rg \varphi= \rg P\;\hbox{ and} \;\rg(\Ker \varphi)+\rg Q=\rg(\Coker \varphi)+\rg P\,.}

In addition, in case of free modules, and for some rank $r\in\NN$, we indeed find the notion of rank of a matrix as defined in~\ref{defRangk}. 

 As for item \emph{2}, note that for two pseudo\polusz~$p$ and~$q$ we have the \egt $\deg pq=\deg p+\deg q$.
\\
This notion of degree is naturally extended to \lot \polus
defined in Exercise~\ref{exoPolLocUnitaire}.

\section{Some applications of the local structure \thoz}
\label{secAppliLocPtf}

In this section we consider results regarding \mptfs and some \alis between those.

Given the local structure \tho for \mptfsz, and since the \deter and the related \pols
are well-behaved by change of base \ri (Fact \ref{fact.det loc}), we have almost systematically all the desired results by means of the \dem given in the following frame. 

\Grandcadre{\DebP\,\, In the free module case, the result is easy to establish. \quad \eop}

We will not always mention it in this section.

NB: if in the hypothesis there is a \lnl \ali between two different modules, by the local structure \tho we are reduced to the case of a \nl \aliz.

\smallskip The \dem works each time the result to be established is true \ssi it is true after \lon at \ecoz.

\medskip
\rem If we must prove a result which, in the case of free modules, comes down to \idas we can also suppose that the
endomorphisms are \digsz.
The argument here is different from the local structure \thoz. The fact is that to check an \ida it suffices to do so on a Zariski open set of the parameter space, and a \gnq matrix is \dig by Proposition~\ref{propMatGenDiag}.\eoe

\subsect{Trace of an \endo and new expression for the \polfon}{Fundamental \polz}
\label{subsec trace}
Recall that if $M$ and $N$ are two \Amosz, we denote by $\theta_{M,N}$ the canonical \ali

\snic{\theta_{M,N}:M\sta\otimes_\gA N\to\Lin_\gA(M,N), \; (\alpha \otimes y)\mapsto
(x\mapsto \alpha(x)y).}

Also recall the following results (Fact \ref{factDualPTF} and Proposition~\ref{propAliPtfs}).

\mni{\it Let $P$ be a \mptfz. \label{lem iso cano}
\begin{enumerate}
\item  $\theta_{P,N}$ is an \iso  of $P\sta\otimes_\gA N$ in $\Lin_\gA(P,N)$.
\item  $\theta_{N,P}$ is an \iso  of $N\sta\otimes_\gA P$ in $\Lin_\gA(N,P)$.
\item  The canonical \homo $P\rightarrow P^{\star\star}$  is an \isoz.
\item  The canonical \homo

\snic{\varphi \mapsto \tra{\varphi}~;~\Lin_\gA(N,P)\,\to\, \Lin_\gA(P\sta,N\sta),}

is an \isoz.
\end{enumerate}
}

\medskip If $P$ is a \mptfz, recall that the \emph{trace} of the \endo $\varphi$ of $P$ (denoted by $\Tr_P(\varphi)$) is the \coe in $X$ of the \polfon $\rF{\varphi}(X)$.
It can \egmt be defined from the natural \ali
$$\preskip-.4em \postskip.4em 
\tr_P:P\sta\otimes_\gA P\rightarrow \gA~:~\alpha \otimes y\mapsto \alpha(y)
, 
$$
and of the canonical \iso $\theta_{P}:P\sta\te_\gA P\to\End(P)$, as follows
$$\preskip.4em \postskip.4em
\Tr_P=\tr_P \circ\, {\theta_{P}}^{-1}.
$$
(The reader will be able to observe that the two \dfns coincide in the case of a free module, or to refer to Fact~\ref{factMatriceEndo}.)

When $P$ and $Q$ are \ptfsz, the trace also allows us to define a canonical duality between $\Lin_\gA(P,Q)$ and $\Lin_\gA(Q,P)$ by means of the bi\lin map $(\varphi,\psi)\mapsto \Tr(\varphi \circ\psi)=\Tr(\psi\circ\varphi)$. This duality can also be defined by the canonical \iso $(P\sta\otimes_\gA Q)\sta\simeq P\otimes_\gA Q\sta$.

\pagebreak

\begin{proposition}
\label{prop trace polfon}
Let $\varphi$ be an \endo of a \mptf $P$ with $n$ \gtrsz.
The \coes of the \polfon of $\varphi$ are given by
$$\preskip.2em \postskip.0em\ndsp
\rF{\varphi}(X) = 1+\som_{h\in\lrbn} \Tr\big(\,\Al h\varphi\,\big) \,X^h\,.
$$
\end{proposition}

\begin{proposition}
\label{prop trace surjective}
If $P$ is a faithful \mptfz, then the trace \Ali $\Tr_P:\End(P)\to \gA$ is surjective.
\end{proposition}

\subsec{Tensor product}
\label{subsec Produit tensoriel}
\begin{proposition}
\label{prop Produit tensoriel}
We consider two \ptf \Amos $P$ and $Q$. Let $\varphi$ and $\psi$ be
endomorphisms of $P$ and $Q$. The module $P\otimes_\gA Q$  is a \mptfz.
\begin{enumerate}
\item We have the \egt 

\snic{\det ({\varphi \otimes
\psi})=(\det {\varphi})^{\rg Q}\,(\det {\psi})^{\rg P}
\eqdefi \rR Q(\det \varphi) \rR P(\det \psi).}

\item The \polfon $\rF{\varphi\otimes \psi}(X)$ of  $\varphi\otimes_\gA \psi$ only depends on $\rg({P})$, on~$\rg({Q})$, on $\rF{\varphi},$ and on $\rF{\psi}$.
\item If
$\rF{\varphi}=(1+\lambda_1X)\,\cdots\,(1+\lambda_mX)$, and
$\rF{\psi}=(1+\mu_1X)\,\cdots\,(1+\mu_nX)$,
we have the \egt
$\rF{\varphi\otimes \psi}(X) =
\prod_{i,j}(1+\lambda_{i}\mu_{j}X)$.

\item
In particular,  $\rg(P\otimes Q)=\rg(P)\rg(Q).$
\end{enumerate}
\end{proposition}

Please note that the last \egt can be rewritten as

\snic{\ide_h(P\otimes Q)= \sum_{ jk=h}\,\ide_j(P)\ide_{k}(Q).
}

\smallskip
Also note that the previous proposition  could be \gui{directly} proven without making use of the local structure \thoz, with a \dem copied from that which was given for exterior powers (Proposition~\ref{prop puissance ext}).

\subsec{Ranks and \alisz}
\label{subsecSURJISO}

\begin{proposition}
\label{prop epi rank constant}\label{prop mono rang}
Let $\varphi :P\rightarrow Q$ be a \ali between \mptfsz.
\begin{enumerate}
\item If $\varphi$ is surjective, then $P\simeq \Ker\varphi\oplus Q$. 
If in addition $\rg(P)= \rg(Q)$, then~$\varphi$ is an \isoz.
\item If $\varphi$ is injective, then $\rg(P)\leq \rg(Q)$.
\end{enumerate}
\end{proposition}

\begin{proof}
In item \emph{2}, it suffices to prove the in\egt after \lon at an \elt $s$ which renders the two modules free.
 As \lon preserves injectivity, we can conclude by the free module case (see Corollary~\ref{corprop inj surj det} and the remark that follows).
\end{proof}

\vspace{-.7em}
\pagebreak

\begin{corollary}
\label{corEpiRC}
Let $P_1 \subseteq P_2 \subseteq P$ with $P_1$ a direct summand 
in $P$. 
    \\
Then $P_1$ is a direct summand in $P_2$.
Consequently, if the modules are \ptfsz, we have the \eqvc

\snic{
 \rg(P_1)=\rg(P_2)\;
 \Longleftrightarrow
\;P_1=P_2.}

 If in addition $P_1\oplus Q_1  =P_2\oplus Q_2=P$, we have the \eqvcs

\snic{
 \rg(P_1)=\rg(P_2)\;\Longleftrightarrow \;\rg(Q_1)=\rg(Q_2)\;
 \Longleftrightarrow
\;P_1=P_2.}
\end{corollary}

\subsec{Transitivity formulas}
\label{subsecTransPtf}

\begin{notation}
\label{notaTraceDetCarAlg}
{\rm
Let $\gB$ be an \Alg, \stfe over $\gA$. Then $[\gB:\gA]=\rg_\gA(\gB)$.
}
\end{notation}

Recall that by Fact~\ref{factTransptf},
if $\gB$ is \stfe over $\gA$, and if~$P$ is a \ptf \Bmoz,
then~$P$ is also a \ptf \Amoz.

When we take $P$ to be a quasi-free module over $\gB$, by considering its rank over $\gA$ it defines a \homo from the additive group $\HO\gB$ to the additive group $\HO\gA$.
This \homo is called the \emph{restriction \homoz} and it is denoted by~$\rRs{\gB/\!\gA}$.\index{restriction!homomorphism}
We thus obtain a contravariant functor from a subcategory of commutative \ris to that of Abelian groups.
This is the category whose morphisms are the $\rho:\gA\to\gB$ which make of~$\gB$ a \asf over $\gA$.

Moreover, $\HO$ defines a covariant functor from the category of commutative \ris to that of semi\risz, since by \edsz, 
a quasi-free module gives a quasi-free module.

As $\HO(\gC)$ is completely \care by $\BB(\gC)$ (for a categorical formulation, see Exercise~\ref{exoAnneaudesrangs}), items \emph{1}  and \emph{2}  of the following fact completely describe the two functors which we have just spoken of. 

\begin{fact}\label{factHORestriction}
Let $\rho:\gA\to\gB$ be an \algz. 
\begin{enumerate}
\item For $e\in\BB(\gA)$, we have $\HO(\rho)([e]_\gA)=\big[\rho(e)\big]_\gB$ in $\HO\gB$. \\
In particular $\HO(\rho)$ is injective (resp.\,surjective, bijective) \ssi
the restriction of $\rho$ to $\BB(\gA)$ and $\BB(\gB)$ is injective 
(resp.\,surjective, bijective). 
\end{enumerate}
Now suppose that $\gB$ is \stfe over $\gA$.
\begin{enumerate}\setcounter{enumi}{1}
\item For $e\in\BB(\gB)$, 
$\Rs\iBA([e]_\gB)=\rg_\gA(e\gB)$, and $\Rs\iBA(1)=[\gB:\gA]$.
\item If a \Bmo $P$ is quasi-free over both $\gA$ and $\gB$, we simply obtain
%
$\rRs{\gB/\!\gA}([P]_{\gB})=[P]_{\gA}$.
%
\end{enumerate} 
\end{fact}
\rem If $\gA$ is connected and contains $\ZZ$, we may pretend to consider $\HO(\gA)\simeq\ZZ$ as a sub\ri of $\gA.$ In item \emph{2}  above we then see that {$\Rs\iBA([e]_\gB)=\big[\Tr\iBA(e)\big]_\gA$}
(it suffices to consider the case where $e\gB$ is free and has a free direct
complement within $\gB$). 
\eoe

The following lemma generalizes \thref{Th.transitivity} 
(which handled the free case).

\begin{lemma}\label{lem1TransPtf} \emph{(Transitivity formulas for the trace and the \deterz)} \\
Let $\gB$  be a \stfe \Alg and $P$ be a \ptf \Bmoz. Let $u_\gB:P\to P$ be a \Bliz, that we denote as
$u_\gA$ when we consider it as an \Aliz.
Then we have the fundamental \egts 

\snic{\det_\gA(u_\gA)=\rN\iBA\big(\det_\gB(u_\gB)\big)\;\hbox{ and }\;\Tr(u_\gA)=\Tr\iBA\big(\Tr(u_\gB)\big).}
\end{lemma}
\begin{proof}
Localizing at \eco of $\gA$ we can assume that $\gB$ is a free \Amoz, of rank $k$.
We write

\snic{P\oplus N=L\simeq\gB^n\simeq\Ae{ nk},}

\sni(the last \iso is an \iso of \Amosz).
We consider $v=u\oplus\Id_N\in\End_\gB(L)$.
Then, by \dfn of the \deterz, we obtain the \egts $\det_\gB(u_\gB)=\det_\gB(v_\gB)$ and $\det_\gA(u_\gA)=\det_\gA(v_\gA)$. We can therefore apply the transitivity formula of \thref{Th.transitivity}.\\
The reasoning for the trace is similar.
\end{proof}

\begin{corollary}\label{lem2TransPtf} Let $\gA\vers{\rho}\gB$  be a \stfe \algz,
 $P$ a  \ptf  \Bmo and $u_\gB\in\End_\gB(P)$.
\begin{enumerate}
\item $\rC{u_\gA}(X) = \rN_{\gB [X]/\!{\AX }} \big(
\rC{u_\gB}(X)\big)$.
\item  $\rF{u_\gA}(X) = \rN_{\gB [X]/\!{\AX }} \big(
\rF{u_\gB}(X)\big)$.
\item  In particular, the \polmus of $P$ over $\gA$ and $\gB$ are linked by

\snic{\rR{P_\gA}(X)= \rN_{\gB[X]/\!{\AX }}\big(\rR{P_\gB}(X)\big) }
\item  The restriction \homo satisfies

\snic{\rRs{\gB/\!\gA}\big(\rg_\gB(P)\big)=\rg_\gA(P).}
\item  If $P$ is a \ptf \Amoz, then

\snic{\rg_\gB\big(\rho\ist(P)\big)=\HO(\rho)\big(\rg_\gA(P)\big), \hbox{ and }\rg_\gA\big(\rho\ist(P)\big)=[\gB:\gA]\,\rg_\gA(P).}

\end{enumerate}
\end{corollary}
\begin{proof}
Items \emph{1}, \emph{2}, \emph{3}  result from the previous lemma.
\\
\emph{4.}
Item \emph{3}  tells us that the \polmu of $P$ over $\gA$ only depends on the \polmu of $P$ over $\gB$. We can therefore assume that $P$ is quasi-free over $\gB$ and we apply the \dfn of the \homo $\Rs\iBA$.
\\
Item \emph{5}  is left to the reader.
\end{proof}

Another corollary is given by the following \thoz.

\begin{theorem}\label{corthTransPtf}
Let   $\gB$ be a \stfe \Alg  and  $\gC$ be a \stfe \Blgz.
Then  $\gC$ is a \stfe \Alg and

\snic{
[\gC:\gA]=\rRs{\gB/\!\gA}([\gC:\gB]).}

In particular, if $\HO(\gA)$ is identified with a sub\ri of $\HO(\gB)$, and if the rank of $\gC$ over $\gB$ is an \elt of $\HO(\gA)$, we have

\snic{
[\gC:\gA]= [\gB:\gA]\, [\gC:\gB].}
\end{theorem}

\vspace{-.5em}
\pagebreak

\subsec{\Pro modules of rank 1}
\label{subsec ptf rank 1}

\begin{fact}
\label{factRANG1}
A matrix $F\in\Mn(\gA)$ is \idme and of rank $1$ \ssi $\Tr(F)=1$ and $\Al2F=0$.
\end{fact}
\facile

\begin{proposition}
\label{propRANG1}
Let $P$ be a \pro \Amo of constant rank~$1$.
\begin{enumerate}
\item The canonical \homos

\snic{\gA\rightarrow \End(P)$, $\,a\mapsto \mu_{P,a}\,$  and
$\,\End(P)\rightarrow \gA$, $\varphi \mapsto \Tr(\varphi)}

are two reciprocal \isosz.
\item For all $\varphi \in \End(P)$, we have  $\det(\varphi)=\Tr(\varphi)$.
\item The canonical \homo $P\sta\otimes_\gA P\rightarrow \gA$ is an \isoz.

\end{enumerate}
\end{proposition}

\begin{proposition}
\label{th ptrg1}
Let $M$ and $N$ be two \Amosz. \\
If $N\te_\gA M$ is \isoc to $\gA$, then $M$ is a \pro module of rank $1$ and  $N$ is \isoc to~$M\sta$.
\end{proposition}
\begin{proof}
Let $\varphi$ be an \iso of $N\otimes_\gA M$  over $\gA$. Let $u=\sum_{i=1}^n c_i\otimes a_i$ be the \elt of $N\otimes M$ such that $\varphi(u)=1$.
We have two \isos from $N\otimes M\otimes M$ to $M$, constructed from $\varphi$.

\snic{c \otimes a  \otimes b\mapsto \varphi(c\otimes a)\, b \quad  
\hbox{and} \quad
c \otimes a\otimes b\mapsto \varphi(c\otimes b)\, a .}

This gives an \iso  $\sigma :M\rightarrow M$ satisfying

\snic{\sigma\big(\varphi(c\otimes a)\, b\big)=\varphi(c\otimes b)\,a \;\hbox{ for all }c \in N,\; a,b\in M, \hbox{ hence}}

\snic{\sigma(x)=
\sigma\big(\sum_i \varphi(c_i\te a_i) x\big) =
\sum_i \varphi(c_i\otimes x) a_i,\hbox{ and }
x=\sum_i \varphi(c_i\te x) \sigma^{-1}(a_i).}

This shows that $M$ is \ptfz, with the \syc $\big((\un),(\aln)\big)$, where $u_i=\sigma^{-1}(a_i)$ \hbox{and $\alpha_i(x)=\varphi(c_i\te x)$}.
Similarly, $N$ is \ptfz. 
\\
But $1=\rg({N\otimes M})=\rg({N})\rg({M})$,
so $M$ and $N$ are of rank~1 (Fact~\ref{fact0ProdRangs}).
Finally, $N\otimes M\sta\te M\simeq N\simeq M\sta$.
\end{proof}

\section{Grassmannians}
\label{secGrassman}

\subsec{The \gnq \ris 
\texorpdfstring{$\Gn$ and $\Gnk$}{Gn and Gn,k}}
\label{subsubsec AGBR}

We have defined the \ri $\Gn=\Gn(\ZZ) =\ZZ[(f_{ij})_{i,j\in\lrbn}]/\cGn$
on \paref{subsec cas generique}.

Actually the construction is functorial and we can define $\Gn(\gA)$ for every commutative \ri $\gA$:  $\Gn(\gA)=\gA[(f_{ij})_{i,j\in\lrbn}]/\cGn\simeq\gA\otimes_\ZZ\Gn.$
\\
Let~$r_k=\ide_k(\Im\,F)$ where $F$ is the matrix $(f_{i,j})$ in  $\Gn(\gA)$.

If we impose in addition that the rank must be equal to $k$, we introduce \hbox{the \id $\cGnk=\cGn+\gen{1-r_k}$} and we obtain the \ri

\snic{\Gnk=\ZZ[F]\sur{\cGnk}\simeq \Gn[1/r_k]\simeq\aqo{\Gn}{1-r_k}.}

We also have the version relativized to $\gA$:
$$
 \Gnk(\gA)=\gA[F]/\cGnk\simeq \Gn(\gA)[1/r_k] \simeq
\gA\otimes_\ZZ\Gnk.
$$
The \ri $\Gn(\gA)$ is \isoc to the product of the \ris $\Gnk(\gA).$

\Grandcadre{In the present subsection devoted to $\Gnk$ we let $h=n-k$.}

If $\gK$ is a field, the \ri $\Gnk(\gK)$ can be considered as the \ri of \coos of the affine \vrt $\GAnk(\gK)$ whose points are the pairs~$(E_1,E_2)$ of subspaces of $\gK^n$ satisfying the \egts $\dim(E_1)=k$ and~$\gK^n=E_1\oplus E_2$.
In \agq \gmtz, there are a few forceful arguments to affirm that the \ri  $\Gnk(\gK)$ has all the right \prts that we can imagine, with respect to the fact that the \vrt  $\GAnk(\gK)$ is a homogeneous space for an action of the \lin group.

We will find these results again \gui{by hand} and by freeing ourselves of the hypothesis \gui{$\gK$ is a field.}

By using the suitable \lons at the principal minors of order~$k$ of the matrix $F=(f_{ij})$ (the sum of these minors is equal to $1$  in~$\Gnk(\gA)$), we will establish a few essential \prts of the functor~$\Gnk$.

\begin{theorem} \emph{(The functor $\Gnk$)}
\label{thBnk}
\begin{enumerate}
\item  \label{thBnk1} There exist \eco $\mu_i$ of the \ri $\Gnk(\gA)$ such that each localized \ri
$ \Gnk(\gA)[1/\mu_i]$  is \isoc to the \ri  
$$\preskip.2em \postskip.2em 
\gA[(X_j)_{j\in\lrb{1..2hk}}][1/\delta] 
$$
for some $\delta$ which satisfies $\delta(\uze)=1$. 
\item \label{thBnka} The natural \homo $\gA\rightarrow \Gnk(\gA)$ is injective.
\item  \label{thBnkf} If $\varphi \,:\gA\rightarrow \gB$ is a \homoz, the kernel of $\Gnk(\varphi)$ is generated by $\Ker\,\varphi$. In particular, if $\varphi$ is injective, so is  $\Gnk(\varphi)$.
\end{enumerate}
\end{theorem}

\begin{corollary}\label{corthBnk}
Let $\gK$ be a \cdi and $\gA$ be a \riz.
\begin{enumerate}
\item  \label{corthBnkc} The \ri  $\Gnk(\gK)$ is integral, \iclz,
\cohz, \noe regular, of \ddk $2kh$.
\item  \label{thBnkd} If $\gA$ is an integral \ri (resp.\,reduced,  \qiriz,  \lsdz,  normal, 
 \coh \noez,  \coh \noe regular) the same goes for $\Gnk(\gA)$.
\item  \label{corthBnke} The \ddk of  $\Gnk(\gA)$ is equal to that of $\gA[X_1,\ldots ,X_{2hk}]$.
\item  \label{corthBnkb} The \ri $\Gnk=\Gnk(\ZZ)$ is integral, \iclz, \coh \noez, regular, of  (Krull) dimension $2k h +1$.
\end{enumerate}
\end{corollary}

\comm In the corollary we have used the notion of a \anor and that of a \ddk that we have not yet defined (see Sections~\ref{subsecEntiers} and \ref{secDefConsDimKrull}).
Finally, a \cori is called regular when every \mpf admits a finite \pro resolution (for this last notion see \Pbmz~\ref{exoFossumKumarNori}). \eoe

\vspace{-.5em}
\pagebreak

\begin{Proof}{\Demo sketch. }
If we make a principal minor of order $k$ of~$F$ \ivz, then the \ri $\Gnk(\gA)$ becomes \isoc to a localized \ri of a \pol \ri over $\gA$, therefore inherits all of the nice \prts  of~$\gA$.
For the fact that $\Gnk(\ZZ)$ is integral there is an added subtleness, because this is not a local \prtz.
\end{Proof}

We now develop the above sketch.
For the \cdi case we start with the following result.
\begin{lemma}
\label{lemPSES1}
Let $\gK$ be a \cdi and $(E_1,E_2)$ be a pair of \supls subspaces of dimensions $k$ and $h$ in $\gK^n$. \\
Suppose that in the matrix $\bloc{ \I_k}{L}{C}{\I_h}$, the first $k$ columns generate $E_1$ and the last $h$ columns generate $E_2$.
\begin{enumerate}
\item The matrices $L$ and $C$ are entirely determined by the pair $(E_1,E_2)$.
\item The matrix $\I_k-LC$  is \iv (we denote its \inv by $V$).
\item The \mprn over $E_1$ \paralm to $E_2$ is equal to 

\snic{F=\bloc{ V}{-V\,L}{C\,V}{-C\,V\,L}.}
\end{enumerate}
\end{lemma}
\begin{proof}
The uniqueness is clear. Let $F=\bloc{ V}{L'}{C'}{W}$ be the matrix of the considered projection.
It is \caree by the \egt

\snic{F\,\bloc{ \I_k}{L}{C}{\I_h}=\bloc{ \I_k}{0}{C}{0}\,,}

\cad  

\snic{V+L'C=\I_k$, $\,VL+L'=0$, $\,C'+WC=C\,$  and  $\,C'L+W=0,}

which is equivalent to
$$
L'=-VL, \;W=-C'L, \;C'(\I_k-LC)=C\,\hbox{ and }\,V\,(\I_k-LC)=\I_k,
$$
or  $(\I_k-LC)^{-1}=V$ , $C'=CV$,  $L'=-VL$, and~$W=-CVL$.
\end{proof}

This can be \gns to the case of a \mprn of rank $k$ over an arbitrary commutative \ri as follows, which is a variant common to both the freeness lemma and the local freeness lemma.

\CMnewtheorem{lelib2}{Second freeness lemma}{\itshape}
\begin{lelib2}\label{leli2} ~\\
Let $F$ be a \prr in $\GAn(\gA)$; recall that $k + h = n$.
\begin{enumerate}
\item If  $\rg(F)\leq k$ and  if a principal minor of order $k$ is \ivz, then the matrix $F$ is similar to a standard \mprn $\I_{k,n}.$
\item More \prmtz, suppose that $F=\bloc{V}{L'}{C'}{W}$ with
$V \in \GL_k(\gA)$.
Let

\snic{B=\bloc{V}{-L'}{C'}{\I_{h}-W}.}

Then, by letting $L=V^{-1}L'$ and $C=-C'V^{-1}$, the matrix $B$ is \ivz, with inverse
$\bloc{\I_k}{L}{C}{\I_{h}}$. In addition, we have the \egts

\snic{\begin{array}{c}
B^{-1}\,F\,B=\I_{k,n},\, W=C'V^{-1}L',\,  V=(\I_k-LC)^{-1},   \\[1mm]
\det V=\det(\I_h-W )\, \hbox{ and } \,
\I_{h}-W=(\I_{h}-CL)^{-1}.
\end{array}
}
\item Conversely, if  $L\in \Ae{ k{\times}h}$,
  $C\in \Ae{h{\times}k}$ and if $\I_k-LC$  is \iv
with \inv $V$, then the matrix

\snic{F=\bloc{V}{V\,L}{-C\,V}{-C\,V\,L}}

is a \prn of rank $k$; it is the \prn over the free submodule~$E_1$ generated by the first $k$ columns of
$\bloc{\I_k}{L}{C}{\I_h}$, \paralm to the free submodule $E_2$ generated by the last $h$ columns of this matrix.
\end{enumerate}

\end{lelib2}

\begin{proof} See Exercise~\ref{exoleli2} and its solution.
\perso{L'\egt $\det\,V=\det\,W_1$ is d\'emontr\'ee in l'exo sous
l'hypoth\`ese que $\det\,V$ is \ivz,
 mais il is probable qu'elle is vraie
pour toute \mprn of rank $k$.
}
\end{proof}

What we have gained relative to the first freeness lemma~\ref{lem pf libre} is that~$F$ is similar to $\I_{k,n}$ instead of simply being \eqve
(however, see Exercise~\ref{exo2.4.1}), and most importantly, the precisions obtained here will be useful to~us.

The previous lemma can be reformulated as follows, in a more abstract form, but essentially \eqve (albeit less precise).
\begin{lemma}
\label{lemthBnk} \emph{(The \ri $\Gnk(\gA)$ is almost a \pol \riz)}
\\
We consider the \gnq matrix $F=(f_{ij})_{i,j\in\lrbn}$ in the \ri $\Gnk(\gA)$. Let $\mu=\det\big((f_{ij})_{i,j\in\lrbk}\big)$ be its leading principal minor of order $k$.
Moreover let $\gA[L,C]$ be the \pol \ri in $2kh$ indeterminates, seen as \coes of two matrices $L$ and $C$ of respective types $k{\times}h$ and~$h{\times}k$. Finally, let $\delta=\det(\I_k-LC)\in\gA[L,C]$.
\\
Then the localized \ris $\Gnk(\gA)[1/\mu]$ and $\gA[L,C][1/\delta]$ are naturally \isocz.
\end{lemma}
\begin{proof}
Let us write \smashtop{$F$ in the form $\bloc{V}{L'}{C'}{W}$}  with $V\in\MM_k(\gA)$.  When we invert $\mu=\det(V)$ we get   $V\in\GL_k(\gA[1/\mu])$.
\\Let $L=V^{-1}L'$  and  $C=-C'V^{-1}$. 
By item~\emph{2} of Lemma~\ref{leli2},  we have $\delta=\det(\I_k-LC)\in\Ati$. This defines a morphism of \Algs from $\gA[L,C][1/\delta]$ to $\Gnk(\gA)[1/\mu]$.\\
In the other direction: to $L$ and $C$  with $\delta$ being \iv we associate the 
matrix  $F=\bloc{V}{V\,L}{-C\,V}{-C\,V\,L}$
(with $V=(\I_k-LC)^{-1}$). \\
The corresponding homomorphism goes from $\Gnk(\gA)[1/\mu]$ to $\gA[L,C][1/\delta]$.\\
By composing these morphisms we find the \idt in both cases.
\end{proof}

\begin{Proof}{\Demo of \thref{thBnk}. }
\\
\emph{\ref{thBnk1}.} This item is deducted from the previous lemma since the sum of the principal minors of order $k$ of $F$ is equal to $1$ in $\Gnk(\gA)$.

\emph{\ref{thBnka}.} Consider the $\gA$-\homo $\psi\,:\gA[(f_{i,j})]\rightarrow \gA$ with specialization at $\I_{k,n}$ defined by $\psi(f_{i,j})=1$ if $ i=j\in\lrbk$ and $=0$ otherwise. \\
It is clear that $\psi\big(\cGnk(\gA)\big)=0$. This proves that $\gA\cap\cGnk(\gA)=0$ because if $a$ is in this intersection, $a=\psi(a)=0$.

\emph{\ref{thBnkf}.} The kernel of $\varphi_{L,C}\,:\gA[L,C]\rightarrow\gB[L,C]$ (the natural extension of $\varphi$) is generated by the kernel of $\varphi$. The \prt remains true after \lonz. Then it remains true by gluing \lons at \mocoz. Therefore in our case we glue by saying that $\Ker\,\Gnk(\varphi)$ is generated by~$\Ker\,\varphi$.
\end{Proof}
\begin{Proof}{\Demo of Corollary \ref{corthBnk}. }
\\
\emph{\ref{thBnkd}.} Besides the integrality question this results from item~\emph{\ref{thBnk1}}  of \thref{thBnk},
because all the considered notions are stable by $\gA\leadsto \gA[X]$ and stem from the basic \plgz. As for the integrality, it is deducted from the result in the case of a \cdiz; if $\gA$ is integral and $S=\Reg(\gA)$, then $\gK=\Frac\gA=\gA_S$ is a \cdi and item \emph{\ref{thBnkf}}  of \thref{thBnk} allows us to conclude.\iplg

\emph{\ref{corthBnke}.} Given the \plgc for the \ddk (\paref{thDdkLoc}), it suffices to show that $\gA[L,C]$ and $\gA[L,C][1/\delta]$ have the same dimension, which result from Lemma~\ref{lemLocMemeKdim} below.

\emph{\ref{corthBnkc}.}  Given items \emph{\ref{thBnkd}}  and \emph{\ref{corthBnke}}  it remains to show that $\Gnk(\gK)$ is integral. In order to do so recall that $\SLn(\gK)$ operates transitively over $\GAnk(\gK)$, which means that every \mprn of rank $k$ and of order $n$ can be expressed in the form $S\cdot\I_{k,n}\cdot S^{-1}$ with $S\in\SLn(\gK)$.
Let us introduce the \ri of \coos of the \vrt  $\SLn(\gK)\subseteq\Mn(\gK)$:

\snic{
\Sln(\gK)=\aqo{\gK[(s_{i,j})_{i,j\in \lrbn}]}{1-\det\,S}.
}

To the surjective map

\snic{\theta_\gK:\SLn(\gK)\rightarrow\GAnk(\gK)\,:\,
S\mapsto S\cdot\I_{k,n}\cdot S^{-1},}

corresponds the $\gK$-\homo

\snic{\wi{\theta}_\gK:\Gnk(\gK)\to\Sln(\gK),}

which sends each  $f_{i,j}$ onto the \coe
$i,j$ of the matrix  $S\cdot\I_{k,n}\cdot S^{-1}$.
\\
It is well-known that $\Sln(\gK)$ is integral, and it therefore suffices to show that~$\widetilde{\theta}_\gK$ is injective.
As $\theta_\gL$ is surjective for every finite extension $\gL$ of $\gK$, every \elt of $\Ker\,\widetilde{\theta}_\gK$ is nilpotent (by the \nstz\footnote{Here we give a \cov \dem by assuming that $\gK$ is contained in an \cdacz. We could adapt it to the \gnl case.}).
However, $\Gnk(\gK)$ is reduced, so $\widetilde{\theta}_\gK$ is injective.

\emph{\ref{corthBnkb}.} Results from the other items (for the \ddkz, \thref{corthValDim} is also required).
\end{Proof}

\begin{lemma}\label{lemLocMemeKdim}
Using the previous notations the \ri $\gA[L,C][1/\delta]$ is a monogenic integral extension
 of a \pol \ri over $\gA$ with $2kh$ \idtrsz.
Consequently
$
\Kdim \gA[L,C][1/\delta] = \Kdim\gA[X_1,\ldots,X_{2kh}].
$
\end{lemma}
\begin{proof} Let $L = (l_{ij})_{i \in\lrbk,j \in\lrbh}$, $C = (c_{ij})_{ i\in\lrbh,j \in\lrbk}$.
The \pol $\delta$ is of degree $2m$ with $m =\min(h,k)$ and contains the \mom

\snic{(-1)^m l_{11} \ldots l_{mm} c_{11} \ldots c_{mm}.}

The localized \ri $\gA' = \gA[L,C][1/\delta]$ can be obtained by adjoining an \idtrz~$t$: $\gA'
= \aqo{\gA[L, C, t]}{t\delta - 1}$. We can put the \pol $g=t\delta - 1$ in \iNoe position.
Indeed, with the \cdv
$$
l'_{ii} = l_{ii} + t, \; c'_{ii} = c_{ii} + t, \  i \in\lrbm, \quad
l'_{ij} = l_{ij},\; c'_{ij} = c_{ij} \hbox { if } i\neq j,
$$
the \pol $g$ becomes, up to sign, \mon in $t$. Therefore $\gA'$ is a monogenic integral extension of  $\gA[L',C']$.
We conclude with \thref{cor2thKdimMor}.
\end{proof}

We will now study tangent spaces to grassmannians.
For this we need to define the concept itself.

We therefore begin with a heuristic introduction to abstract categorical and functorial notions. 
The readers unfamiliar with the language of categories will have to 
skip  
this introduction, in which we give practically no \demsz, and simply try to convince themselves 
from the given examples that the notion of a tangent space to a functor at a point is, all in all, quite reasonable, which will allow us them 
 to then see the beautiful application of this concept to grassmannians.

\subsec{Affine schemes, tangent spaces}

\subsubsec{\nst and equivalence of two categories}
\label{subsecNstMorphismes}

Let $(\uf)=(\lfs)$  be a \syp in $\kXn=\kuX$, and let $\gA=\kxn=\kux$ be the corresponding quotient \algz.

We have seen on \paref{ZerosCrucial} the crucial identification 
\Grandcadre{$\Hom_{\gk} (\gA,\gk)=\cZ(\uf,\gk)\subseteq \gk^n$} 
 between the zeros over $\gk$  of the \syp $(\uf)$ and the \crcs  of the \algz~$\gA$.
If $\gk$ is reduced, we obviously have $\Hom_{\gk} (\gA,\gk)=\Hom_{\gk} (\Ared,\gk)$.%
 \index{algebraic variety!over an \cacz}
 \Grandcadre{Now suppose that $\gk$ is a discrete \cacz.}
Such a set of zeros $\cZ(\uf,\gk)\subseteq \gk^n$ is then called an \emph{\vgq  over~$\gk$}. 
 
Let $\gA$ and $\gB$ be two quotient \klgs corresponding to two \syps $(\uf)$ and $\ug$ in $\kuX=\kXn$. The \nst (Corollary~\ref{corthNstClass}) tells us that the two reduced \algs $\gA\red$ and $\gB\red$ are equal \ssi they have the same \vrt of zeros in $\gk^n$:
\Grandcadre{$\cZ(\uf,\gk)=\cZ(\ug,\gk)\iff \rD_{\kuX}(\uf)=\rD_{\kuX}(\ug)\iff \gA\red=\gB\red$}
This observation is the first step in the development of the \eqvc between the category of \rpfs \klgs on the one hand, and that of \vgqs over $\gk$ on the other.
 
For the \eqvc to be complete, we must also treat the morphisms. Therefore we provide a preliminary study regarding the \algz~$\Ared$.

Notice that every \elt $p$ of $\kuX$ defines a \polle function  $\gk^n\to\gk,\;\uxi\mapsto p(\uxi)$, and that an \elt of $\Ared$ defines (by restriction) a function $\cZ(\uf,\gk)\to\gk$; indeed, if $p\equiv q\mod\rD_{\kuX}(\uf)$, a power \hbox{of $p-q$} is in the \id $\gen{\uf}$, so the restrictions of the \polles functions $p$ and $q$ to~$\cZ(\uf,\gk)$ are equal.
But in the case where $\gk$ is an \cacz, we have the converse; if the restrictions of $p$ and $q$ to $\cZ(\uf,\gk)$ are equal, $p-q$ vanishes over $\cZ(\uf,\gk)$, and by \nstz, a power \hbox{of $p-q$} is in \hbox{the \id $\gen{\uf}$}. 

\rdb
Thus, $\Ared$ can be interpreted as an \alg of functions over the \vgq that it defines, namely $A=\cZ(\uf,\gk)=\Hom_{\gk} (\gA,\gk)$. The \klg structure of $\Ared$ is indeed that of this \alg of functions.
These functions $\cZ(\uf,\gk)\to \gk$ are called the \emph{\frgsz}. 
\index{function!regular ---}\index{regular!function}%

Similarly, if $\gA=\kxn$ and $\gC=\gk[\ym]$ are the quotient \algs corresponding to two \syps 
$$ 
(\uf) \hbox{  in } \kXn \hbox{ and } (\uh) \hbox{ in } \kYm,
$$
\rdb
if $A=\cZ(\uf,\gk)\subseteq \gk^n$ and $C=\cZ(\uh,\gk)\subseteq \gk^m$ are the corresponding \vgqsz, we define a \emph{\aregz} from~$A$ to~$C$ as the restriction to $A$ and $C$ of a \poll map~\hbox{$\varphi:\gk^n\to\gk^m$} which sends~$A$ to~$C$.%
\index{regular!map}\index{map!regular ---}
\\
The \aregs are, by \dfnz, \emph{the morphisms from $A$ to~$C$ in  category of the \agq \vrts over $\gk$}.
We will denote by $\Mor_\gk(A,C)$ the set of these morphisms.
 
The above map $\varphi$  is given by a \sys $(F_1,\dots,F_m)$ in~$\kuX$, or by the \homo $F:\kuY\to\kuX,\;Y_j\mapsto F_j$.\\
Let $\varphi_1:A\to C$ be the restriction of $\varphi$; if $\gamma:C\to \gk$ is a \frgz, then the composite function $\gamma\circ \varphi_1:A\to\gk$ is a \frgz, and the map $\psi_1:\gamma\mapsto \gamma\circ \varphi_1$ can be seen as a map from~$\gC\red$ to~$\Ared$. Actually, this map is none other than the \homo which comes from $F$ by passage to the quotients. 

In the opposite direction, we can see that every \homo $\psi_1:\gC\red\to\Ared$ comes from a \homo $\psi:\gC\to\gA$, and that $\psi$ defines a \areg $\varphi:A\to C$, sometimes called the \emph{co-morphism} of $\psi$.
This takes place as follows: via the identifications  
$A=\Hom_{\gk} (\gA,\gk)$ \hbox{and $C=\Hom_{\gk} (\gC,\gk)$},
we simply have the \egt $\varphi(\uxi)=\uxi\circ \psi$ (which makes of~$\varphi$ the \gui{transpose} of $\psi$).\index{co-morphism}

Finally,  $\Mor_\gk(A,C)$, is naturally identified with  $\Hom_{\gk}(\gC\red,\gA\red)$, 
 identification that we express under the form of an \egtz:
\Grandcadre{$\Mor_\gk(A,C)=\Hom_{\gk}(\gC\red,\gA\red)$.}
However, note that the direction of the arrows is reversed.

Let us consider as a special case the case where $A$ is the \vgq  
        reduced to a point, 
associated with the \algz~$\gk$, 
corresponding to the empty \syp over \emph{the variable-free \polle \alg $\gk$}. If preferred, here we can see $\gk$ as the quotient $\aqo\kX X$, corresponding to the point~$\so 0$, sub\vrt of the \vgq $V=\gk$
associated with the \alg $\kX$.
\\
In these conditions, the framed \egt above admits as a special case \hbox{${C=\Mor_\gk(\so{0},C)=
\Hom_{\gk}(\gC\red,\gk).}$}
We have come full circle!

The summary of this study is the following: we can entirely reduce the  consideration of \vgqs over an \cac to the study of \rpf \klgsz. This is an interpretation in finite terms (finite \syps over~$\gk$ for the objects as well as for the morphisms) of objects a priori slightly more mysterious, and certainly more infinite. In categorical terms: we can advantageously replace the category of \agq \vrts over $\gk$ with the opposite category to that of \rpf \klgsz. There is a natural \eqvc between these two categories.

\subsubsec{Affine schemes}
\label{subsecSchAff}
 
Now we take a big leap into abstraction. 
First of all we admit that \vrts can have multiplicities. For example the intersection of a circle with a line must be a double point, and  not only a point, when the line is tangent to the circle. Consequently, it is sometimes counterproductive to limit ourselves to reduced \klgsz.

We also admit that our \ri is not \ncrt 
 an \cac but an arbitrary commutative \riz. In which case the points of the \vrt over $\gk$ are not sufficient to characterize what we want to consider as an abstract \agq \vrt defined over $\gk$ (by allowing multiplicities).
For example the abstract circle is certainly represented by the \ZZlg

\snic{\ZZ[x,y]=\aqo{\ZZ[X,Y]}{X^2+Y^2-1},}

but it is not its points over $\ZZ$ that will give us much information. On the contrary, it is its points over all the \ZZlgs (\cad over all the commutative \risz) that matter. Similarly an abstract \emph{double circle} is certainly represented by the \ZZlg

\snic{\ZZ[x',y']=\aqo{\ZZ[X,Y]}{(X^2+Y^2-1)^2},}

but we would not know how to distinguish a simple circle from a double circle if we only consider the points over the reduced \ris (the \ris without multiplicity).

We now can define 
the category of \emph{affine schemes over the commutative \ri $\gk$}. This could simply be the opposite category to the category of \klgsz; that whose objects are the \klgs and whose arrows are the \homos of \klgsz.

But there exists a strictly more meaningful (and elegant?) \eqve description:
\emph{an affine scheme over the commutative \ri $\gk$ is known when its zeros over all the \klgs are known}. 
In other words, the \klgz~$\gA$ defines an affine scheme which is nothing other that the functor $\Hom_{\gk}(\gA,\bullet)$ from the category of \klgs to the category of sets.

A \homo of \klgs $\gB\to\gA$ defines a natural transformation from the functor
$\Hom_{\gk}(\gA,\bullet)$ to the functor  $\Hom_{\gk}(\gB,\bullet)$: natural transformations of functors are \gui{in the right direction,}
\cad from the zeros of $\gA$ to the zeros of $\gB$.

If we do not want to abstract up too high, we can limit ourselves to \pf \klgsz, which is quite enough to make very beautiful abstract \agq \gmt (\cad not limited to \agq \gmt over \cdisz).

\subsubsec{Tangent space at a point of a functor}
\label{subsubsecTanFonct}%
\index{tangent space}\index{space!tangent ---}

First of all recall the notion of a tangent space to a \syp at a zero of the \sys introduced in Section~\ref{secExlocGeoAlg}.

Take the example of the sphere as an affine scheme defined over $\QQ$.
This scheme is associated with the \QQlg 

\snic{\gA=\QQ[x,y,z]=\aqo{\QQ[X,Y,Z]}{X^2+Y^2+Z^2-1}.}

If $\uxi=(\alpha,\beta,\gamma)\in\QQ^3$ is a zero of $\gA$ over $\QQ$, \cad a rational point of the sphere, we associate to it
\begin{itemize}
\item the \id $\fm_\uxi=\gen{x-\alpha,y-\beta,z-\gamma}_\gA$,
\item  the local \alg $\gA_\uxi=\gA_{1+\fm_\uxi}$, and 
\item the tangent space $\rT_\uxi(\gA/\QQ)
\simeq\Der\QQ\gA\uxi$, 
\end{itemize}
which is a \Qev canonically \isoc to $(\fm_\uxi/{{\fm_\uxi}^2})\sta$ or to $(\fm_\uxi\gA_\uxi/{{\fm_\uxi\gA_\uxi}^2})\sta$.

More intuitively but \eqvz ly
(Proposition~\ref{propTangent}), a tangent vector to the sphere at $\uxi$ is simply given by a $(u,v,w)\in\QQ^3$ which satisfies $u\alpha+v\beta+w\gamma=0$, \cad by letting $f=X^2+Y^2+Z^2-1$, 
$$\preskip.4em \postskip.4em 
u\Dpp f X(\uxi)+v\Dpp f Y(\uxi)+w\Dpp f Z(\uxi)=0.
$$
Here is now a new way to see this tangent space, which we express in terms of corresponding affine schemes, \cad of the functor $\Hom_\QQ(\gA,\bullet)=\cZ(f,\bullet)$. For this we must formally introduce an infinitesimal which we denote by $\vep$, \cad consider the \QQlg $\QQ[\vep]=\aqo{\QQ[T]}{T^2}$ ($\vep$ is the class of $T$ modulo~$T^2$).

The point $\uxi$ is seen as a \crc of $\gA$, \cad as the \elt $\wi\uxi:g\mapsto g(\uxi)$ 
of $\Hom_\QQ(\gA,\QQ)$. 
We then ask ourselves what are the $\lambda$  \elts of the set $\Hom_\QQ(\gA,\QQ[\vep])$
that \gui{lift $\wi\uxi$,} in the sense that when composed with the evaluation of $\vep$ at $0$, from $\QQ[\vep]$ to $\QQ$, we once again obtain $\wi\uxi$.

\centerline{\xymatrix {
             &\QQ[\vep]\ar[d]^{\vep := 0}  \\
\gA\ar@{-->}[ur]\ar[r]^{\wi\uxi} &\QQ                   \\
}} 

Such an \elt is a priori given by a zero of $f$ over $\QQ[\vep]$ which recovers~$\uxi$ when we  evaluate $\vep$ at $0$, \cad a triple 
$(\alpha + a\vep, \beta + b\vep, \gamma + c\vep)$,  
 with $f(\alpha + a\vep, \beta + b\vep, \gamma + c\vep) = 0$ in $\QQ[\vep]$. But this  means precisely that $(a,b,c)$ is a tangent vector to the sphere at $\uxi$. 

It is simply the substance of the mundane observation which states that
\gui{the \dile is 
the \lin component of the increase in the function}:

\snic{f(\uxi+\vep V)=f(\uxi)+\vep \,\rd f(\uxi)(V)  \mod \vep^2.}

This zero $\uxi+\vep(a,b,c)$ of $\gA$ in $\gk[\vep]$ defines a \homo
$\gA\to\gk[\vep]$ via  $x\mapsto \alpha+a\vep$, $y\mapsto \beta+b\vep$, $z\mapsto \gamma+c\vep$.
\\
This \homo sends $g$ to $g(\uxi)+a\Dpp g X(\uxi)+b\Dpp g Y(\uxi)+c\Dpp g Z(\uxi)$, since

\snic{g\big(\uxi+\vep (a,b,c)\big)=g(\uxi)+\vep \,\rd g(\uxi)(a,b,c)  \mod \vep^2.}

\smallskip 
The reader will be able to check that this little computation, which we have just performed on a small example, works for any zero of any \syp based on any commutative \riz.

However, we need to at least add how to interpret the  \kmo structure over the tangent space at a zero of a \syp over a \ri $\gk$ in terms of the functor $\Hom_\gk(\gA,\bullet)$. 

Here also the use of our small example shall be sufficient. 

In the category of \QQlgsz, the fiber product of the \gui{restriction arrow} 

\snic{\QQ[\vep]\to\QQ$, $\vep\mapsto0}

which itself is the \alg

\snic{\QQ[\vep]\times_\QQ\QQ[\vep]\simeq \QQ[\vep_1,\vep_2] 
\quad \hbox{with }  \vep_1^2=\vep_1\vep_2=
\vep_2^2=0,}

equipped with the two \gui{\prn} \homos 

\snic{\QQ[\vep_1,\vep_2]\vers{\pi_1}\QQ[\vep],\;\vep_1\mapsto \vep,\; \vep_2\mapsto 0\quad$ and $\quad\QQ[\vep_1,\vep_2]\vers{\pi_2}\QQ[\vep],\;\vep_2\mapsto \vep,\; \vep_1\mapsto 0,~~~~~~}

and with the \gui{restriction} arrow
$$\preskip.4em \postskip.4em 
\QQ[\vep_1,\vep_2]\to\QQ,\quad\vep_1\mapsto 0,\;\vep_2\mapsto 0. 
$$
There is also a natural \homo \gui{of addition}
$$
\preskip.4em \postskip.4em 
\QQ[\vep_1,\vep_2]\to\QQ[\vep],\quad\vep_1\mapsto\vep,\; \vep_2\mapsto\vep, 
$$
which commutes with the restrictions.

When we give two zeros $\uxi+\vep V_1$ and $\uxi+\vep V_2$ of $\gA$ in $\QQ[\vep]$, given the \cara \prt of the fiber product in the category of \QQlgsz, the two corresponding \homos $\gA\to\QQ[\vep]$ uniquely factorize to give a \homo from $\gA$ to $\QQ[\vep_1,\vep_2]$, the \gui{fiber product of the two,} which corresponds to the zero $\uxi+\vep_1 V_1+\vep_2 V_2$ of $\gA$ in $\QQ[\vep_1,\vep_2]$. 

Finally, by composing this {fiber product} \homo with the addition \homo $\QQ[\vep_1,\vep_2]\to\QQ[\vep]$, we obtain the \homo corresponding to the zero $\uxi+\vep (V_1+V_2)$. We have therefore come full circle, the addition of tangent vectors has been described in purely categorical terms.

\entrenous{un joli dessin semblerait utile pour faire comprendre ce qui se passe ici}

Let us recap. In the case of the functor that is an affine scheme defined by a \syp over a \ri $\gk$ with its quotient \alg $\gA$, there is a canonical identification between $\rT_\uxi(\gA\sur\gk)$ and the set of points of~$\gA$ over $\gk[\vep]$ that lift $\uxi$, when we identify $\uxi$ and $\uxi+\vep V$ with the corresponding \elts of $\Hom_\gk(\gA,\gk)$ and $\Hom_\gk(\gA,\gk[\vep])$.
In addition, the \kmo structure in the second interpretation is given by the \gui{addition} provided by the \homo 
$$\preskip.4em \postskip.2em 
\gk[\vep_1,\vep_2]\simeq\gk[\vep]\times_\gk \gk[\vep]\to\gk[\vep],\quad
\vep_1\mapsto\vep,\;\vep_2\mapsto\vep, 
$$
($\vep^2=\vep_1^2=\vep_2^2=\vep_1\vep_2
=0$). 

Note that the \gui{external law,} multiplication by the scalar $a$, comes from the \homo 
$$\preskip.0em \postskip.4em 
\gk[\vep]\vers{\lambda_a}\gk[\vep],\quad b+\vep c\mapsto b+\vep ac. 
$$
The formal mechanism of addition described thus will work with any other functor which will itself  be willing to transform the fiber products (in the category of \QQlgsz) into fiber products 
(in the category of sets).

Thus the notion of a tangent space at a point of a functor\footnote{Functor from category of \klgs to the category of sets.} \gnss to the other schemes over a \ri $\gk$, because they are \gui{good functors.} That is, the Grothendieck schemes (that we will not define here) are good functors. The Grassmannian functors (which have already been defined) are such schemes.

\pagebreak	
\subsec{Tangent spaces to the Grassmannians}

\subsubsec{Projectors and ranks}
\label{subsecPrRg}

Two easy facts before entering the heart of the matter.
Consider a module~$E$.
Two \prrs $\pi_1$, $\pi_2 : E \to E$  are said to be \ixc{orthogonal}{projectors} if they satisfy $\pi_1 \circ \pi_2 = \pi_2 \circ \pi_1 = 0$.

\begin{fact}\label{lemSomProjOrt}
If $\pi_1$, $\pi_2 : E \to E$ are \ort \prrs of images $E_{1}$  
and $E_{2}$, then $\pi_1 + \pi_2$ is a \prr and its image is $E_1 \oplus E_2$. Consequently, when $E$ is a \mptfz,  we obtain

\snic{\rg({\pi_{1}+\pi_{2}}) = \rg({E_{1}\oplus E_{2}}) = \rg {E_1} + \rg {E_2} .}
\end{fact}

\begin{fact} \label{lemProjProj}
Let $\pi_1$, $\pi_2 \in \End_\gA(E)$ be two \prrs of images $E_1$ and $E_{2}$.
Then the \Ali

\snic{
\Phi : \End(E)\to\End(E),\; \varphi \mapsto \pi_2\circ\varphi\circ\pi_1,}

is a \prr whose image is \isoc to $\Lin_{\gA}(E_1, E_2)$. Consequently, when $E$ is a \mptfz, we obtain the \egt

\snic{\rg\Phi  = \rg E_1 \cdot \rg E_2.}
\end{fact}

\subsubsec{Affine Grassmannian}
\label{subsecGrassAff}

This subsection is dedicated to the determination of the tangent space at a point to the functor $\gA \mapsto \GAn(\gA)$.
Recall that the acronym $\GA$ is used for ``Affine Grassmannian.''
The \gmq interpretation of a point $P$ of  $\GAn(\gA)$ is given by the ordered pair $(E,F)=(\Im P,\Ker P)$ of submodules as a direct sum in $\Ae n$.

More \gnltz, if $\gk $ is a \ri given as a reference (in usual geometry it would be a \cdiz) and if $M$ is a fixed \ptf \kmoz, we can consider the category of \klgs and the  
functor $\gA \mapsto \GA_M(\gA)$, where $\GA_M(\gA)$ designates the set of ordered pairs~$(E,F)$ of submodules as a direct sum in the extended module~$\gA\otimes _{\gk}M$, which we will denote by~$M_{\gA}$.
Such a pair can be represented by the projection~$\pi:M_{\gA}\to M_{\gA}$ over $E$ \paralm to $F$.
The affine Grassmannian~$\GA_M(\gA)$ can therefore be seen as the subset of \idm \elts in~$\End_{\gA}(M_{\gA})$.
It is this point of view that we adopt in the following.

To study the tangent space we must consider the \Alg $\gA[\varepsilon]$ where $\varepsilon$ 
is the \gnq \elt with null square.
First of all we give the statement of the usual Grassmannian $\GAn(\gA)$.

\vspace{.3em}
\pagebreak

\begin{theorem}\label{prop1TanGrassmann} \emph{(Tangent space to an affine Grassmannian)}\\
Let $P\in \GAn(\gA)$  be a \prr of image $E$ and of kernel~$F$.
For $ H \in \MM_n(\gA)$ we have the following \eqvcz.

\snic{P + \varepsilon H\in \GAn(\gA[\varepsilon]) \quad  \Longleftrightarrow
\quad  H = HP + PH .}
%

Let us associate to $P$ the \Ali $\wh P : \MM_n(\gA) \to \MM_n(\gA)$ defined by

\snic{\widehat P(G) = P\,G\,(\In-P) + (\In-P)\,G\,P.}

We have the following results.
\begin{enumerate}
\item [--] The \Alis 

\snic{\pi_1 : G\mapsto P\,G\,(\In-P)\;$ and
$\;\pi_2 : G\mapsto (\In-P)\,G\,P}

are \ort \prrsz.
In particular,  $\widehat P$ is a \prrz.
\item [--] For $H \in \MM_n(\gA)$, we have $H = PH + HP$ \ssi $H\in\Im\wh P$.
\item [--] The module $\Im\wh{P}$ is canonically  \isoc to $\Lin_{\gA}(E,F)\oplus\alb \Lin_{\gA}(F,E)$.
In particular, $\rg(\Im\wh P)= 2 \rg E  \cdot \rg  F $.
\end{enumerate}
In brief, the tangent space at the $\gA$-point $P$ to the functor $\GAn$ is canonically \isoc to the \mptf  $\Im\wh{P}$ (via $H \mapsto P+\varepsilon H$), itself canonically \isoc to $\Lin_{\gA}(E,F)\oplus\alb \Lin_{\gA}(F,E)$.
\end{theorem}
\begin{proof}
The first item is immediate. Let $V_{P}$ be the submodule of matrices~$H$ which satisfy $H=HP+PH$. This module is canonically \isoc to the tangent space that we are looking for. 
A simple computation shows that~$\pi_1$ and~$\pi_2$ are \ort \prrsz. Therefore $\widehat P$ is a \prrz. The following \egt is clear:  $P\widehat P(G) + \widehat P(G)P = \widehat P(G)$.
Therefore $\Im \wh P \subseteq V_{P}$.
Moreover, if~$H = PH+HP$, we have $PHP = 0$, so $\wh P(H) = PH + HP=H$.
Thus~$V_{P} \subseteq \Im \wh P$.
In brief $V_{P}=\Im \wh P=\Im\pi_{1} \oplus \Im\pi_{2}$: \trf by applying Fact~\ref{lemProjProj}.
\end{proof}

We now give the \gnl statement (the proof is identical).

\begin{proposition}\label{prop2TanGrassmann}
Let $\pi\in \GA_M(\gA)$ be a \prr of image $E$ and of kernel~$F$.
For $ \eta \in \End_{\gA}(M_{\gA})$ we have the \eqvc

\snic{\pi + \varepsilon \eta\in \GA_M(\gA[\varepsilon]) \quad  \Longleftrightarrow
\quad  \eta = \pi\eta + \eta\pi .}
%

We associate to $\pi$ the \Ali $\wh \pi : \End(M_{\gA}) \to \End(M_{\gA})$ defined by
$
\widehat \pi(\psi) = \pi\,\psi\,(\I-\pi) + (\I-\pi)\,\psi\,\pi
$. Then
\begin{enumerate}
\item [--] The \alis $\pi_1 : \psi\mapsto \pi\, \psi\,(\I-\pi )$ and $\pi_2 : \psi\mapsto (\I-\pi )\,\psi\, \pi $ are \ort \prrsz. 
In particular,  $\widehat \pi$ is a \prrz.
\item [--] An \Ali $\eta \in \End(M_{\gA})$ satisfies $\eta = \pi\eta + \eta\pi$ \ssi $\eta\in\Im\wh \pi$.
\item [--] The module $\Im\wh{\pi}$ is canonically \isoc to $\Lin_{\gA}(E,F)\oplus\alb \Lin_{\gA}(F,E)$.
In particular, $\rg(\Im\wh \pi )= 2 \rg E  \cdot \rg  F$.
\end{enumerate}
In short the tangent space at the $\gA$-point $\pi$ to the functor $\GA_M$ is canonically \isoc to the \mptf  $\Im\wh{\pi}$ (via $\eta \isosim \pi+\varepsilon \eta$), itself canonically \isoc to $\Lin_{\gA}(E,F)\oplus\alb \Lin_{\gA}(F,E)$.
\end{proposition}

\subsubsec{Projective Grassmannian}
\label{subsecGrassProj}

This subsection is dedicated to the determination of the tangent space at a point to the functor $\gA \mapsto \GGn(\gA)$, where $\GGn(\gA)$ designates the set of submodules which are direct summands in $\Ae{n}$.

\begin{fact}\label{fact1GrassProj} \emph{(The space of  \prrs that have the same image as a given \prrz)}
Let $P\in \GGn(\gA)$ be a \prr of image $E$. Let $\Pi_E$ be the set of \prrs that have $E$ as their image, and $V=\Ae n$.
Then $\Pi_E$ is an affine subspace of~$\Mn(\gA)$, having for \gui{direction} the \ptf \Amoz~$\Lin_\gA(V/E,E)$ (naturally identified with an \Asub of $\Mn(\gA)$).
We specify the result as follows.
\begin{enumerate}
\item \label{i1fact1GrassProj}
Let $Q \in \GGn(\gA)$ be another \prrz. \\
Then $Q\in \Pi_E$ \ssi $PQ = Q$ and $QP = P$. \\
In this case, the difference $N = Q-P$ satisfy the \egts $PN = N$ and~$NP = 0$, and so $N^2 = 0$.
\item \label{i2fact1GrassProj}
Conversely, if $N \in \Mn(\gA)$ satisfy $PN = N$ and $NP = 0$ (in which case $N^2 = 0$),
then $Q := P+N$ is in $\Pi_E$.
\item \label{i5fact1GrassProj}
In short, the set $\Pi_E$ is identified with the \Amo
$\Lin_\gA(V/E, E)$ via the affine map

\snic{\Lin_\gA(V/E, E)\to \Mn(\gA),\;\varphi\mapsto P+j\circ\varphi\circ \pi,}

where $j:E\to V$ is the canonical injection and $\pi:V\to V/E$ is the canonical \prnz.
\end{enumerate}
%
Additional information.
\begin{enumerate}\setcounter{enumi}3
\item \label{i3fact1GrassProj}
If $Q \in \Pi_E$, $P$ and $Q$ are conjugated in $\Mn(\gA)$.
More \prmtz, by letting $N = Q-P$, we have $(\In + N)(\In - N) = \In$ and $(\In - N) P (\In + N) = Q$.
\item \label{i4fact1GrassProj}
If $Q \in \Pi_E$, then
for all $t \in \gA$, we have $tP + (1-t)Q\in \Pi_E$.
\end{enumerate}

\end{fact}
\begin{proof}
\emph{\ref{i1fact1GrassProj}.}
$N^2 = 0$ as seen when multiplying $PN = N$ by $N$ on the left-hand side. 

\emph{\ref{i5fact1GrassProj}.}
The conditions $PN = N$ and $NP = 0$ over the matrix $N$ is equivalent to the inclusions $\,\Im N\subseteq E=\Im P$
and $E\subseteq \Ker N$. 
\\
The matrices $N$ of this type form an \Amo $\wi E$ which can be identified with the module $\Lin_\gA(\Ker P, \Im P)$
 \gui{by restriction of the domain and of the image.}
\\
More intrinsically, this module $\wi E$ is also identified with $\Lin_\gA(V/E, E)$ via the \ali
 $\Lin_\gA(V/E, E) \to \Mn(\gA),\;\varphi\mapsto j\circ\varphi\circ \pi$,
which is injective and admits $\wi E$ as its image.

\emph{\ref{i3fact1GrassProj}.} $(\In - N)\, P \,(\In + N) = P\,(\In + N) = P + PN = P + N = Q$.
\end{proof}

\vspace{-.7em}
\pagebreak

\begin{fact}\label{fact2GrassProj}
Let $E\in\GGn(\gA)$ and $E'\in\GGn(\gA[\varepsilon])$ which gives $E$ by the specialization $\varepsilon\mapsto 0$ (in other words $E'$ is a point of the tangent space at $E$ to the functor $\GGn$). Then
$E'$ is \isoc to the module obtained from $E$ by \edsz:
$E'\simeq \gA[\varepsilon]\otimes _\gA E$.
\end{fact}
\begin{proof}
By \thref{propComparRedRed}, a \mptf $M$ over a \riz~$\gB$ is characterized, up to \isoz, by its \gui{reduction} $M\red$ (\cad the module obtained by \eds 
to~$\gB\red$).
However, $E'$ and $\gA[\varepsilon]\otimes _\gA E$ have the same reduction $E\red$ to $(\gA[\varepsilon])\red\simeq\Ared$.
\perso{demander une preuve directe en exercice?}
\end{proof}
%

\begin{theorem}\label{prop3TanGrassmann} \emph{(Tangent space to a projective Grassmannian)}\\
Let $E\in \GGn(\gA)$ be an \Asub which is a direct summand in $\Ae{n}=V$.
Then the tangent space at the $\gA$-point $E$ to the functor $\GGn$ is canonically \isoc to $\Lin_\gA(E,V/E)$.
More \prmtz, if $\varphi\in\Lin_\gA(E,V/E)$ and if we let
$$
E_\varphi=\sotq{x+\varepsilon h}{x\in E,\,h\in V,\,h\equiv\varphi(x)\mod E}
,$$
then $\varphi\mapsto E_\varphi$ is a bijection from the module $\Lin_\gA(E,V/E)$ to the set of matrices $E'\in \GGn(\gA[\varepsilon])$ that give $E$ when we specialize $\varepsilon$ at $0$.
\end{theorem}
%
\begin{proof}
Let $E \in \GG_n(\gA)$ and $\varphi \in \Lin_\gA(E, V/E)$.
\\
\emph{Let us first show that $E_\varphi$ is in $\GG_n(\gA[\varepsilon])$ and above~$E$.}
Let us fix a matrix $P \in \GA_{n}(\gA)$ satisfying $E = \Im P$. We therefore have $V = E \oplus\Ker P$ and an isomorphism $V/E \simeq\Ker P \subseteq V$.
We can therefore lift the \aliz~$\varphi$ at a matrix $H \in \MM_{n}(\gA) = \End(V)$
in accordance with the diagram
$$\preskip-0.0ex\xymatrix @C=28pt @R = 16pt {
&V \ar@{->}[r]^{ H} \ar@{->>}[d]    & V  \\
&E \ar@{->}[r]^{ \varphi}          & V/E \ar@{(->}[u] &
}$$

\vspace{-.6em}The matrix $H$ satisfies $PH = 0$ and $H(\I_{n}-P) = 0$, \cad $HP = H$.
\\
It suffices to show that $P + \varepsilon H$ is a \prr of image $E_\varphi$. \\
For the inclusion $\Im(P + \varepsilon H) \subseteq
E_\varphi$, let $(P + \varepsilon H)(y + \varepsilon z)$ with $y, z \in V$:

\snic {
(P + \varepsilon H)(y + \varepsilon z) = Py + \varepsilon (Hy + Pz) =
Py + \varepsilon (HPy + Pz) = x + \varepsilon h
,}

with $x = Py \in E$, $h = Hx + Pz$. Since $x \in E$, we have $\varphi(x) = Hx$,
and  
so $h \equiv \varphi(x) \bmod E$.  
For the converse inclusion, let $x + \varepsilon h \in E_\varphi$ 
        and let us show that 
$(P + \varepsilon H)(x + \varepsilon h) = x + \varepsilon h$:

\snic {
(P + \varepsilon H)(x + \varepsilon h) = Px + \varepsilon(Hx + Ph)
.}

As $x \in E$, we have $Px = x$. We need to see that $Hx + Ph = h$,  but $h$ is of the form $h = Hx + y$ with $y \in E$, so $Ph = 0 + Py = y$ and we indeed have the \egt $h = Hx + Ph$.
\\
Finally, it is clear that $P + \varepsilon H$ is a \prr

\snic {
(P+\varepsilon H) (P + \varepsilon H) = P^2 + \varepsilon(HP + PH) =
P +\varepsilon H.
}


\emph{Let us show the surjectivity of $\varphi \mapsto E_\varphi$}. Let $E' \subseteq \gA[\varepsilon]^{n}$, a direct summand, above~$E$.
 Then $E'$ is the image of a \prr $P + \varepsilon H$ and we have
$$\preskip.4em \postskip.4em 
(P + \varepsilon H)(P + \varepsilon H) = P^2 + \varepsilon (HP + PH),
\;\;  \hbox{so} \; P^2 = P \;\hbox{and} \; HP + PH = H
, 
$$
which gives $PHP = 0$ (multiply $HP + PH = H$ by $P$ on the right-hand side, for instance). We therefore see that~$P$ is a \prr of image~$E$ (because $E'$, for $\varepsilon := 0$, it is~$E$). We replace $H$ with $K = HP$, which satisfies

\snic {
KP = (HP)P = K, \qquad PK = P(HP) = 0.
}

This does not change the image of $P + \varepsilon H$, \cad $\Im(P + \varepsilon H) = \Im(P + \varepsilon K)$. To see this, it suffices (and is necessary) to show that

\snic {
(P + \varepsilon H)(P + \varepsilon K) = P + \varepsilon K, \qquad
(P + \varepsilon K)(P + \varepsilon H) = P + \varepsilon H
.}

On the left-hand side, we obtain $P + \varepsilon(HP + PK) = P + \varepsilon(HP + 0) = P + \varepsilon K$; On the right-hand side, $P + \varepsilon(KP + PH) = P + \varepsilon(K + PH) = P + \varepsilon(HP + PH) = P + \varepsilon H$.
\\
The matrix $K$ satisfies $KP = K$, $PK = 0$ and represents a \ali $\varphi : E \to \Ae {n}/E$ with $E' = {\rm
Im}(P + \varepsilon K) = E_\varphi$.


\emph{Let us prove the injectivity of $\varphi \mapsto E_\varphi$}. Therefore suppose 
$E_\varphi = E_{\varphi'}$. We fix a \prr $P \in \GG_n(\gA)$ with image~$E$ and we encode $\varphi$ as~$H$, $\varphi'$ as~$H'$ with
$$\preskip-.4em \postskip.4em 
HP = H, \quad PH = 0, \qquad H'P = H', \quad PH' = 0
. 
$$
As $P + \varepsilon H$ and $P + \varepsilon H'$ have the same image, we have the \egts

\snic{
(P + \varepsilon H)(P + \varepsilon H') = P + \varepsilon H' \quad\hbox{and}\quad 
(P + \varepsilon H')(P + \varepsilon H) = P + \varepsilon H
.}

The \egt on the right gives $H = H'$, so $\varphi = \varphi'$.
\end{proof}


\rem
The \prn $\GA_{n}\to\GGn$ associates to~$P$ its image $E = \Im P$. Here is how the tangent spaces and the \prn (with $F = \Ker P$) are organised


\snac{\xymatrix @C = 0.4cm
{
\rT_P(\GA_{n},\gA)\ar[d] \ar@{-}[r]^(0.39){_\sim}
&
\Lin_\gA(E,F)\oplus\Lin_\gA(F,E)\ar@{-}[r]^(0.42){_\sim} &
\{H \in \MM_{n}(\gA) \mid H = HP + PH \}
     \ar[d]^{ H \mapsto K = HP}
\\
\rT_E(\GGn, \gA) \ar@{-}[r]^(0.45){_\sim} &
\Lin_\gA(E, \Ae {n}\!/E) \ar@{-}[r]^(0.33){_\sim} &
\{K \in \MM_{n}(\gA) \mid KP = K,\, PK = 0\}
\\
}}

\hfill\eoe

\section[Grothendieck and Picard groups]{Grothendieck and  Picard groups}
\label{subsecClassifMptfs}

Here we tackle the \gnl \pb of the complete classification of \mptfs over a fixed \ri $\gA$.

This classification is a fundamental but difficult \pbz, which does not admit a \gnl \algq solution.

We start by stating some waymarks for the case where all the \mrcs are free.

In the following subsections we give a very small introduction to classic tools that allow us to apprehend the \gnl \pbz.

\subsec{When the \mrcs are free}\label{subsecMrcLibre}

Let us begin with an \elr remark.
\begin{fact}
\label{factRgcstLib}
A \pro \Amo of rank $k$ is free \ssi it is generated by $k$ \eltsz.
\end{fact}
\begin{proof}
The condition is clearly \ncrz. Now suppose the module generated by $k$ \eltsz.
The module is therefore the image of a \prn matrix $F\in\Mk(\gA)$. By hypothesis $\det(\I_k+XF)=(1+X)^k$. In
particular, $\det F=1$, so $F$ is \ivz, and since $F^2=F$, this gives $F=\I_k$.
\end{proof}

Here is another easy remark.
\begin{fact}
\label{factRgcstLib2}
Every \Amrc is free \ssi every \pro \Amo is quasi-free.
\end{fact}
\begin{proof}
The condition is clearly sufficient. If every \Amrc is free and if $P$ is \proz, let $(r_0,\ldots ,r_n)$ be the corresponding \sfioz.
Then $P_k=r_kP\oplus (1-r_k)\Ae k$ is a \pro \Amo of rank $k$ therefore free.
Let $B_k$ be a base, the \gui{component} $r_kB_k$ is in $r_kP$, and $r_kP\simeq \left(r_k\gA\right)^k$. Since $P$ is the direct sum of the $r_kP$, it is indeed quasi-free.
\end{proof}

\begin{proposition}\label{thZerdimLib}
Every \mrc over a \algb is free. 
\end{proposition}

\begin{proof}
Already seen in \thref{thlgb2}.
\end{proof}

\begin{theorem}
\label{thBézoutLib}
Every \mptf over a Bézout domain is free.
Every \mptf of constant rank over a Bézout \qiri is free.
\end{theorem}
\begin{proof}
Let us consider the integral case. A \pn matrix of the module can be reduced to the form  $\cmatrix{T&0\cr0&0}$ where $T$ is triangular with \ndz \elts over the diagonal (see Exercise~\ref{exoBézoutstrict}).
As the \idds of this matrix are \idms the \deter $\delta$ of $T$ is a \ndz \elt which satisfies $\delta\gA=\delta^2\gA$. Thus $\delta$ is \iv and the \pn matrix is \eqve to $\cmatrix{\I_k&0\cr0&0}$.\\
For the \qiri case we apply the  \elgbmd explained on \paref{MethodeQI}.
\end{proof}

Let us take note of another important case:  $\gA=\BXn$ where $\gB$ is a Bézout domain. This is a remarkable extension of the Quillen-Suslin \thoz, due to Bass (for $n=1$), then Lequain and Simis~\cite{LS}.
The \tho will be proved in Section~\ref{sec.Etendus.Valuation}.

\subsec{\texorpdfstring{$\GKO(\gA)$, $\KO(\gA)$ and $\KTO(\gA)$} {GK0(A), K0(A) and \KTO(\gA)}}
\label{secGKO}

Let $\GKO \gA$ be the set of 
        \iso classes
of \mptfs over $\gA$.
It is a semi\ri for the inherited laws of $\oplus$ and $\otimes$.
The \textsf{G}  of $\GKO$ is in tribute to Grothendieck. \label{NOTAGKO}

Every \elt of $\GKO \gA$ can be represented by an \idme matrix with \coes in~$\gA$. Every \ri \homo $\varphi:\gA\to\gB$ induces a \homo $\GKO \varphi : \GKO \gA\to \GKO \gB$.
So $\GKO$ is a covariant functor from the category of commutative \ris to the category of semi\risz.
We have $\GKO (\gA_{1}\times\gA_{2})  \simeq  \GKO \gA_{1} \times \GKO \gA_{2}$.
The passage from a \pro module to its dual defines an involutive \auto of  $\GKO \gA$.

If $P$ is a \ptf \Amo we can denote by $[P]_{\GKO \gA}$ the \elt of $\GKO \gA$ that it defines.

The subsemi\ri of  $\GKO \gA$ generated by $1$ (the class of the \mptf $\gA$) is \isoc to $\NN$, except in the case where $\gA$ is the trivial \riz.
As a subsemi\ri of  $\GKO \gA$ we also have the one generated by the  \iso classes of the modules $r\gA$ where $r\in\BB(\gA)$,  \isoc to $\HOp (\gA)$.
We easily obtain the \iso $\HOp (\gA)\simeq\GKO\big(\BB(\gA)\big)$. Moreover, the rank defines a surjective \homo of semi\ris  $\GKO \gA\to \HOp (\gA)$, and the two \homos $\HOp (\gA)\to\GKO \gA\to \HOp (\gA)$ are composed according to the identity.

The \ix{Picard group} $\Pic \gA$ is the subset of $\GKO \gA$ formed by the \iso classes of the \mrcs $1$.
By Propositions~\ref{propRANG1} and~\ref{th ptrg1} this is the group of \iv \elts of the semi\ri  $\GKO \gA$
(the \gui{inverse} of $P$ is the dual of $P$). \label{NOTAPic}

\rdb
The (commutative) additive \mo of $\GKO \gA$ is not always regular.
To obtain a group, we symmetrize the additive \mo $\GKO\gA$ and we obtain the \ix{Grothendieck group} that we denote by $\KO \gA$.\label{NOTAK0}

The class of the \mptf $P$ in $\KO \gA$ is denoted by $[P]_{\KO (\gA)}$, or~$[P]_{\gA}$, or even $[P]$ if the context allows it. 
Every \elt of $\KO \gA$ is written in the form $[P]-[Q]$.
More \prmtz, it can be represented under the two forms
\begin{itemize}
\item [$\bullet$] [projective] - [free] on the one hand,
\item [$\bullet$] [free] - [projective] on the other hand.
\end{itemize}
Indeed

\snic{[P] - [Q] = [P \oplus P'] - [Q \oplus P']
          = [P \oplus Q'] - [Q \oplus Q'],}

with a choice of $P \oplus P'$ or $Q \oplus Q'$ being free.

The defined product in $\GKO \gA$ gives by passage to the quotient a product in $\KO \gA$, which therefore has a commutative \ri structure.\footnote{When the \ri $\gA$ is not commutative, there is no more multiplicative structure over $\GKO \gA$. This explains why the usual terminology is Grothendieck group and not Grothendieck \riz.}

The classes of two \mptfs $P$ and $P'$ are equal in~$\KO \gA$ \ssi
there exists an integer $k$ such that $P\oplus\Ae k\simeq P'\oplus\Ae k$.
We say in this case that $P$ and $P'$ are \ixc{stably \isoc}{modules}.

Two stably \isoc quasi-free modules are \isocz, so that $\HO \gA $  is identified with a sub\ri of $\KO\gA$, and when $P$ is quasi-free, there is no conflict between the two notations $[P]_\gA$ (above and \paref{notaHO+}).

\rdb\label{NOTAKTO}
Two stably \isoc \mptfs $P$ and $P'$ have the same rank since $\rg(P\oplus\Ae k)=k+\rg(P)$. 
Consequently, the (\gnez) rank of the \mptfs defines a surjective \ri \homo $\rg_\gA:\KO \gA\to\HO \gA$.
Let $\KTO \gA$ be its kernel. The two \homosz~$\HO \gA\to\KO \gA\to\HO \gA$ are composed in terms of the identity, in other words the map $\rg_\gA$ is a \crc of the $\HO(\gA)$-\alg $\KO\gA$ and we can write

\snic{\KO (\gA)=\HO (\gA)\oplus \KTO (\gA).}

\smallskip
The structure of the \id $\KTO \gA$ of $\KO \gA$ concentrates a good part of the mystery of classes of stable \iso of \mptfsz, since~$\HO \gA$ presents no mystery (it is completely decrypted by~$\BB(\gA)$). 
In this framework the following result can be useful (cf. \Pbmz~\ref{exoLambdaGammaK0}).

\begin{proposition}\label{propKTO}
The \id  $\KTO \gA$ is the nilradical of $\KO \gA$.
\end{proposition}
%

Finally, note that if $\rho:\gA\to\gB$ is a \ri \homoz, we obtain correlative \homos

\snic{
\KO \rho:\KO \gA\to\KO \gB, \quad \KTO \rho:\KTO \gA\to\KTO \gB\quad \hbox{and}
\quad \HO \rho:\HO \gA\to\HO \gB.}

And  $\KO$, $\KTO$ and $\HO$ are functors.

\subsec{The Picard group}
\label{subsecPicGp}

The Picard group is not affected by the passage to the classes of stable \isoz, because of the following fact.

\begin{fact}
\label{factPicStab}
Two stably \isoc \mrcs $1$ are \isocz.
In particular, a stably free module of rank $1$ is free.
More \prmtz, for a  \pro module  $P$ of consant rank $1$ we have
\begin{equation}
\label{eqfactPicStab0}
P\simeq\Al{k+1}(P\oplus\Ae k) .
\end{equation}
In particular, $\Pic\gA$ is identified with a subgroup of $(\KO\gA)\eti$.
\end{fact}
\begin{proof}
Let us prove the \isoz: it results from \gnl \isos given in the proof of Proposition~\ref{prop puissance ext} (\eqrf{eqVik}).
For arbitrary \Amos $P$, $Q$, $R$, \ldots, the consideration of the universal \prt that defined the exterior powers leads to
$$
\preskip.3em \postskip.3em
\arraycolsep2pt\begin{array}{rcl}
\Al2(P\oplus Q)& \simeq  & \Al2P\oplus(P\te Q)\oplus \Al2Q,
\\[1.4mm] \mathrigid 1mu
\Al3(P\oplus Q\oplus R)& \simeq  & \Al3P\oplus\Al3Q\oplus\Al3R\oplus\big(\Al2P\te Q\big)
\oplus \cdots \oplus (P\te Q\te R)
 ,
\end{array}
$$
with the following \gnl formula by agreeing that $\Al{0}(P_i) = \gA$
\begin{equation}\label{eqfactPicStab}
\Al k\big(\bigoplus\nolimits_{i=1}^mP_i\big)\simeq
\bigoplus_{\som_{i=1}^{m} k_i= k}
\Big(
\big(\, \Al{k_1}P_1 \big)\te\cdots \te
\big(\, \Al{k_m}P_m \big)
\Big).
\end{equation}
In particular, if $P_1$, \ldots, $P_r$ are \mrcsz~$1$ we obtain
\begin{equation}
\label{eqfactPicStab2}
\Al r (P_1\oplus\cdots\oplus P_r) \simeq P_1\otimes \cdots\otimes  P_r.
\end{equation}
It remains to apply this with the direct sum
 $P\oplus\Ae k=P\oplus\gA\oplus \cdots \oplus\gA$.
 The \iso of \Eqrf{eqfactPicStab0} is then obtained with the \Ali
 $P\to \Al{k+1}(P\oplus\Ae k)$,
 $x\mapsto x\vi 1_1 \vi 1_2 \vi \cdots \vi 1_k$,
 where the index represents the position in the direct sum
 $\gA\oplus \cdots \oplus\gA$.\\
 The last affirmation is then clear since we have just shown that the map $\GKO\gA\to\KO\gA$, restricted to  $\Pic\gA$, is injective.
\end{proof}

\rem The reader will be able to compare the previous result and its \dem with Exercise~\ref{exoStabLibRang1}.
\eoe

\smallskip 
We deduce the following \thoz.
\begin{theorem}
\label{thPicKTO} \emph{($\Pic \gA$ and $\KTO \gA$)} 
Suppose that every $\gA$-\mrc $k+1$ ($k\geq1$) is \isoc to a module $\Ae k\oplus Q$. Then the map from $(\Pic \gA,\times )$ to $(\KTO \gA,+)$ defined by
$$[P]_{\Pic \gA}\mapsto [P]_{\KO\! \gA} - 1_{\KO\! \gA} $$
is a group \isoz.
In addition, $\GKO\gA=\KO\gA$ and its structure is entirely known from that of $\Pic\gA$.\perso{Il semble qu'en \gnl
l'application ne soit pas un \homo of groups.}
\end{theorem}
\begin{proof}
The map is injective by Fact~\ref{factPicStab},
and surjective by hypothesis. It is a group \homo because $\gA\oplus(P\otimes Q) \simeq P\oplus Q$, 
\egmt in virtue of Fact~\ref{factPicStab}, since 
$$\preskip.0em \postskip.4em 
\;\;\Al2\big(\gA\oplus(P\otimes Q)\big)\simeq P\otimes Q
\simeq \Al2(P\oplus Q). 
$$

\vspace{-1.2em}
\end{proof}

Note that the law of $\Pic \gA$ is inherited from the tensor product whilst that of~$\KTO \gA$ is inherited from the direct form. We will see in Chapter~\ref{chapKrulldim} that the hypothesis of the \tho is satisfied for
\ris of \ddkz~$\leq1$.

\medskip \rdb
\comm \label{comHOclassique}
We have seen in Section~\ref{subsecCalRang} how the structure of $\HO(\gA)$ directly stems from that of the \agB $\BB(\gA)$.
\\
From the \clamaz' point of view the \agB $\BB(\gA)$ is the \alg of the open and closed sets in $\Spec\gA$ (the set of \ideps of $\gA$ equipped with a suitable topology, cf. Chapter~\ref{chapKrulldim}).
An \elt of  $\BB(\gA)$ can therefore be seen as the \cara function of an open-closed set of $\Spec\gA$.
Then the way in which we construct $\HO(\gA)$ from $\BB(\gA)$ shows that $\HO(\gA)$ can be seen as the \ri of functions with integer values, integral \colis of \elts in $\BB(\gA)$.
It follows that $\HO(\gA)$ is identified with the \alg of \lot constant functions, with integer values, over $\Spec\gA$.
Still from the point of view of \clama the (\gnez) rank of a \ptf \Amo $P$ can be seen as the function (with values in $\NN$) defined over $\Spec\gA$ as follows:
to a \idep $\fp$ we associate the rank of the free module  $P_\fp$ (over a \alo all the \mptfs are free). The \ri $\HO(\gA)$ is indeed obtained simply by symmetrizing the semi\ri $\HOp (\gA)$ of the ranks of \ptfs \Amosz.
\eoe

\subsubs{Picard group and class group of a \riz}
Let us consider the multiplicative  \mo of the \emph{\tf \ifrsz} of the \ri $\gA$, formed by the \tf \Asubs of the total \ri of fractions $\Frac\gA$.
We will denote this \mo by $\Ifr\gA$. \label{NOTAIfr}

More \gnlt a \emph{\ifrz} of $\gA$ is an \Asub $\fb$ of $\Frac\gA$
such that there exists some \ndz $b$ in $\gA$ satisfying $b\,\fb\subseteq\gA$.

In short we can see $\Ifr\gA$ as the \mo obtained from that of the \itfs of $\gA$ by forcing the invertibility of the \idps generated by \ndz \eltsz.%
\index{fractional!ideal}%
\index{ideal!fractional ---}

An \id $\fa\in\Ifr\gA$ is sometimes said to be \emph{integral} if it is contained in $\gA$,
in which case it is a \itf of $\gA$ in the usual sense.

An arbitrary \id $\fa$ of $\gA$  is \iv like an \id of $\gA$ (in the sense of \Dfnz~\ref{defiiv}) \ssi it is an \iv \elt in the \mo $\Ifr\gA$.
Conversely every \id of $\Ifr\gA$ \iv in this \mo is of the form $\fa/b$, where $b\in\gA$ is \ndz and $\fa$ is an \iv \id of $\gA$.
The \iv \elts of $\Ifr\gA$ form a group, the \emph{group of \iv \ifrs of~$\gA$}, which we will denote by $\Gfr\gA$.%
\index{group of invertible fractional ideals}

As an \Amoz, an \iv \ifr is \prc $1$.\label{pageclassgroup}
Two \iv \ids are \isoc as \Amos if they are equal modulo the subgroup of \iv \idps (i.e. generated by a \ndz \elt of $\Frac\gA$). We denote by $\Cl \gA $ the quotient group, that we call the \emph{group of classes of \iv \idsz}, or simply the \emph{class group} of the \ri $\gA$, and we
obtain a well-defined natural map $\Cl \gA\to\Pic\gA$.%
\index{class group!of a \riz}%
\index{group of classes of invertible ideals}%
\index{class of invertible ideals}
\label{NOTAGfr}

Moreover, let us consider an integral and \iv \id $\fa$.
Since $\fa$ is flat, the natural map $\fa\te_\gA\fb\to\fa\fb$ is an \isoz, for any \idz~$\fb$
(\thref{thplatTens}). Thus, the map~$\Cl \gA\to\Pic\gA$ is a group \homoz, and it is clearly an injective \homoz, so~$\Cl\gA$ is identified with a subgroup of $\Pic\gA$.

These two groups are often identical as the following \tho shows, which results from the previous considerations and from \thref{propRgConstant2}.

\begin{theorem}
\label{propRgConstant3} \emph{(Modules of constant rank $1$ as \ids of $\gA$)}\\
Suppose that over $\Frac \gA$ every \pro module of rank $1$ is free.
\begin{enumerate}
\item Every \pro \Amo of rank $1$  is \isoc to an \iv \id of~$\gA$.
\item Every \pro \id of rank $1$ is \ivz.
\item The group of classes of \iv \ids is naturally \isoc to the Picard group.
\end{enumerate}
 \end{theorem}
%
\begin{proof}
\Thref{propRgConstant2} shows that every \pro module of rank $1$ is \isoc to an \id  $\fa$. It therefore remains to see that such an \id is \ivz.
Since it is \lop it suffices to show that it contains a \ndz \eltz.
For this we consider an integral \id $\fb$ \isoc to the inverse of $\fa$ in $\Pic\gA$. The product of these two \ids is \isoc to their tensor product (because $\fa$ is flat) so it is a free module, thus it is a \idp generated by a \ndz \eltz.
\end{proof}

NB: With regard to the comparison of $\Pic\gA$ and $\Cl\gA$ we will find a more \gnl result in Exercise~\ref{exoPicAPicFracA}.

\subsect{The semi\ris \texorpdfstring{$\GKO(\gA)$, $\GKO(\Ared)$
and $\GKO(\gA\sur{\Rad \gA})$} {GK0(A), GK0(Ared) and GK0\big(A/Rad(A)\big)}}{Semi-\ris $\GKO(\gA)$, $\GKO(\Ared)$ and $\GKO(\gA\sur{\Rad \gA})$}
\label{subsecComparRed}

In this subsection we use $\Rad \gA$, the Jacobson radical of $\gA$, which is defined on \paref{eqDefRadJac}. We compare the \mptfs defined over $\gA$, those defined over $\gA'=\gA\sur{\Rad \gA}$ and those defined over $\Ared$.

The \eds from $\gA$ to $\gB$ transforms a \mptf defined over $\gA$ into a \mptf over $\gB$. From a \prn matrix point of view, this corresponds to considering the matrix transformed by the \homo $\gA\to\gB$.

\begin{proposition}
\label{propComparRedJac} ~\\
The natural \homo from $\GKO(\gA)$ to $\GKO(\gA\sur{\Rad \gA})$ is injective,
which means that if two {\mptfs $E$, $F$}  over~$\gA$ are \isoc over $\gA'=\gA\sur{\Rad \gA}$,
they also are over $\gA$. More precisely, if two \idmes matrices $P$, $Q$ of the same format are conjugated over $\gA'$, they also are over $\gA$, 
with an \auto which lifts the residual conjugation \autoz. 
\end{proposition}
\begin{proof}
Denote by $\ov{x}$ the object $x$ seen modulo $\Rad \gA$. Let $C\in\Mn(\gA)$
be a matrix such that $\ov{C}$ conjugates $\ov{P}$ with $\ov{Q}$. 
Since $\det C$ is \iv modulo $\Rad \gA$, $\det C$ is \iv in $\gA$ and $C\in\GL_n(\gA)$. So we have $\ov{Q}=\ov{C\,P\,C^{-1}}$. Even if it means replacing 
$P$ with $C\,P\,C^{-1}$ we can assume $\ov{Q}=\ov{P}$ and $\ov C=\In$. In this case we search an
invertible matrix $A$ such that $\ov A=\In$ and $APA^{-1}=Q$.
\\
We remark that $QP$ encodes an \Ali from $\Im P$ to $\Im Q$ that \rdt gives the \idtz.
Similarly $(\In-Q)(\In-P)$ encodes an \Ali from  $\Ker P$ to $\Ker Q$ that \rdt gives the \idtz. Taking inspiration from the \dlg lemma (Lemma~\ref{propIsoIm}), this leads us to the matrix $A=QP+(\In-Q)(\In-P)$ which realizes $AP=QP=QA$  and $\ov{A}=\In$, so $A$ is \iv and $APA^{-1}=Q$.\\
For two \rdt \isoc \mptfsz~$E$ and~$F$ we use the \dlg lemma which allows us to realize 
 $\ov{E}$ and~$\ov{F}$ as images of \idmes conjugated matrices of the same format.
\end{proof}

As for the reduction modulo the nilpotents, we obtain in addition the possibility to lift every \mptf on account of Corollary~\ref{corIdmNewton}. Hence the following \thoz.

\begin{theorem}
\label{propComparRedRed}
The natural \homo $\GKO(\gA)\to\GKO(\Ared)$ is an \isoz. More precisely, we have the following results.
\begin{enumerate}
\item 
\begin{enumerate}
\item Every \idme matrix over $\Ared$ is lifted to an \idme matrix over~$\gA$.

\item Every \mptf over $\Ared$ comes from a \mptf over~$\gA.$

\end{enumerate}
\item

\begin{enumerate}
\item If two \idme matrices of the same format are conjugated over $\Ared$, they also are over $\gA$, 
with an \auto which lifts the residual conjugation \autoz. 

\item   Two \mptfs over $\gA$ \isoc over $\Ared$, are also \isoc over $\gA.$

\end{enumerate}
\end{enumerate}
\end{theorem}

\subsec{The Milnor square}
A commutative square (in an arbitrary category) of the following style 
$$\preskip.0em \postskip.2em
\xymatrix @C=1.2cm{
A\,\ar[d]^{i_1}\ar[r]^{i_2}   & \,A_2\ar[d]^{j_2}   \\
A_1\,\ar[r]    ^{j_1}    & \,A'  \\
}
$$
is called a \ix{Cartesian square} if it defines $(A,i_1,i_2)$ as the  limit (or inverse limit, or \pro limit) of $(A_1,j_1,A')$, $(A_2,j_2,A')$.
In an equational category we can take 
\[\preskip.2em \postskip.4em
A=\sotq{(x_1,x_2)\in A_1\times A_2}{j_1(x_1)=j_2(x_2)}.
\]
The reader will verify for example that given $\gA\subseteq \gB$ and an \id $\ff$ of $\gA$ which is also an \id of $\gB$ (in other words, $\ff$ is contained in the conductor of $\gB$ into $\gA$), we have a Cartesian square of commutative \risz, defined below.
$$\preskip-.2em \postskip.2em
\xymatrix @C=1.2cm{
\gA\,\ar@{->>}[d]\ar[r]   & \,\gB\ar@{->>}[d]   \\
\gA\sur{\ff}\,\ar[r]    & \gB\sur{\ff} \\
}
 $$
Let $\rho:\gA\to\gB$ be a \homoz, $M$ be an \Amo and $N$ be a \Bmoz. Recall that an \Ali $\alpha
:M\to N$ is a \emph{\eds morphism} (cf.\ \Dfnz~\ref{defAliAliExtScal}) \ssi the natural \Bli  $\rho\ist(M)\to N$ is an \isoz.

In the entirety of this subsection we consider in the category of commutative \ris the \gui{Milnor square} below on the left, denoted by $\cA$, in which~$j_2$ is surjective,
$$\preskip-.24em \postskip.4em
\xymatrix @C=1.2cm{
\ar@{}[dr]|*+<10pt>[o][F]{\cA}
\gA\,\ar[d]^{i_1}\ar[r]^{i_2}   & \,\gA_2\ar@{->>}[d]^{j_2}   \\
\gA_1\,\ar[r]    ^{j_1}    & \,\gA'  \\
}
\quad
\xymatrix @C=1.2cm{
M\,\ar[d]^{\psi_1}\ar[r]^{\psi_2}   & \,M_2\ar@{->>}[d]^{\varphi_2}   \\
M_1\,\ar[r]    ^{\varphi_1}    & \,M'  \\
}
\quad
\xymatrix @C=1.2cm{
E\,\ar[d]\ar[r]   & \,E_2\ar@{->>}[d]^{{j_2}\ist}   \\
E_1\,\ar[r]^{h\circ {j_1}\ist}    & \,E'  \\
}
 $$
Given an \Amo $M$, an $\gA_1$-module $M_1$, an $\gA_2$-module $M_2$, an $\gA'$-module $M'$  and a Cartesian square of \Amos as the one illustrated in the center above, the latter is said to be \emph{adapted to $\cA$}, if the $\psi_i$'s and $\varphi_i$'s are \eds morphisms.

Given an $\gA_1$-module $E_1$, an $\gA_2$-module $E_2$, and an \iso \hbox{of $\gA'$-modules}
$$\preskip.2em \postskip.2em 
h:{j_1}\ist(E_1)\to {j_2}\ist(E_2)=E', 
$$
let $M(E_1,h,E_2)=E$ (above on the right-hand side) be the  \Amo limit of
the diagram
$$\preskip.2em \postskip.4em 
\big(E_1,h\circ {j_1}\ist,{j_2}\ist(E_2)\big),\big(E_2,{j_2}\ist,{j_2}\ist(E_2)\big) 
$$
Note that a priori the obtained Cartesian square is not \ncrt adapted to $\cA$.

\begin{theorem}
\label{thCarreMil1} \emph{(Milnor's \thoz)}
\begin{enumerate}
\item Suppose that $E_1$ and $E_2$ are \ptfsz, then
\begin{enumerate}
\item $E$ is \ptfz,
\item the Cartesian square is adapted to $\cA$: the natural \homos ${j_k}\ist(E)\to E_k$ ($k=1,2)$ are \isosz.
\end{enumerate}
\item Every \mptf over $\gA$ is obtained (up to \isoz) by this procedure.
\end{enumerate}
\end{theorem}
We will need the following lemma.
\begin{lemma}
\label{lemCarreMil1}
Let $A\in\Ae{m\times n}$, $A_k= i_k(A)$ ($k = 1, 2$),
$A'=j_1(A_1)=j_2(A_2)$, $K=\Ker A\subseteq\Ae n $,  $K_i=\Ker A_i$ ($i=1,2$),
 $K'=\Ker A'$. Then $K$ is the  limit (as an \Amoz) of $K_1\to K'$ and $K_2\to K'$.
\end{lemma}
\begin{proof}
Let $x\in\Ae n $, $x_1 = {j_1}\ist(x) \in \gA_1^n$, 
$x_2 = {j_2}\ist(x) \in \gA_2^n$. Since $x\in K$ \ssi $x_i\in K_i$ for $i=1,2$,
$K$ is indeed the desired  limit.
\end{proof}
The reader will notice that the lemma does not apply in \gnl to the submodules that are images of matrices.
\begin{Proof}{\Demo of \thref{thCarreMil1}. }
\emph{2.} If $V\oplus W=\Ae n $, let $P$ be the \pro matrix over $V$ \paralm to $W$. If $V_1$, $V_2$, $V'$ are the modules obtained by \eds to $\gA_1$, $\gA_2$ and $\gA'$, they are identified with kernels of the matrices $P_1=i_1(\In-P)$,  $P_2=i_2(\In-P)$,
$P'=j_2(\In-P_2)=j_1(\In-P_1)$, and the lemma applies: $V$ is the  limit of $V_1\to V'$ and  $V_2\to V'$. The \iso $h$ is then $\Id_{V'}$.
This \gui{miracle}  takes place thanks to the identification of ${j_i}\ist(V_i)$ and $\Ker P_i$.

\emph{1a.} Let $P_i\in\MM_{n_i}(\gA_i)$ be a \prr with image \isoc to $E_i$ ($i=1,2$). We dispose of an \iso of $\gA'$-modules 
 from $\Im\big(j_1(P_1)\big)\in\MM_{n_1}(\gA')$ 
 \hbox{to $\Im\big(j_2(P_2)\big)\in\MM_{n_2}(\gA')$}. 
 Let $n=n_1+n_2$. By the \dlg lemma~\ref{propIsoIm} there exists a matrix $C\in{\En(\gA')}$
 realizing the conjugation

\snic{
\Diag(j_1(P_1),0_{n_2}) = C\,\,\Diag\big(0_{n_1},j_2(P_2)\big)\,\,C^{-1}.
}

Since $j_2$ is surjective (ha ha!), $C$ is lifted to a matrix $C_2\in\En(\gA_2)$. Let

\snic{Q_1=\Diag(P_1,0_{n_2}),$ $\;Q_2=C_2\,\Diag(0_{n_1},P_2)\,C_2^{-1},}

such that $j_1(Q_1)=j_2(Q_2)$ (not bad, right?). There then exists a unique matrix $Q\in\Mn(\gA)$ such that $i_1(Q)=Q_1$ and $i_2(Q)=Q_2$. 
The uniqueness of~$Q$  assures $Q^2=Q$, and the previous lemma applies to show that~$\Im Q$ is \isoc to~$E$ (hats off to you, Mr Milnor!).

\emph{1b.} Results from the fact that $Q_k=i_k(Q)$ and $\Im Q_k\simeq \Im P_k\simeq E_k$
for $k=1,2$.
\end{Proof}

The following fact is purely categorical and left to the good will of the reader.
\begin{fact}
\label{factCarMil2}
Given two Cartesian squares adapted to $\cA$
as found below, it amounts to the same thing to take a \aliz~$\theta$ from~$E$ to~$F$
or to take three \alis (for the corresponding \risz) $\theta_1:E_1\to
F_1$,
 $\theta_2:E_2\to F_2$ and  $\theta':E'\to F'$ 
which make the adequate squares commutative.
$$\preskip.2em
\xymatrix@R=0.5cm@C=1.5cm{
E\ar[dd]\ar@{-->}[drr]_>>>>>>>>>>>>>>\theta\ar[r]   &
E_1\ar[dd]\ar[drr]^{\theta_1}             \\
            &            &  F\ar[r]\ar[dd]    & F_1\ar[dd] \\
E_2\ar[drr]_{\theta_2}\ar[r]        &  E'\ar[drr]^<<<<<<<<<<{\theta'}        &
\\
            &            &  F_2\ar[r]         & F'           \\
}
$$
\end{fact}

\begin{corollary}
\label{corthCarreMil1} ~\\
Consider two 
modules $E=M(E_1,h,E_2)$ and $F=M(F_1,k,F_2)$ like in \thref{thCarreMil1}.
Every \homo $\psi$ of $E$ in $F$ is obtained using two $\gA_i$-module \homos $\psi_i:E_i\to F_i$ compatibles with $h$ and $k$ in the sense that the diagram below is commutative. The \homo $\psi$ is an \iso \ssi
$\psi _1$ and $\psi _2$ are \isosz.
$$\preskip.2em 
\xymatrix @C=1.2cm{
{j_1}\ist(E_1)\,\ar[d]^{h}\ar[r]^{{j_1}\ist(\psi_1)}   & \, {j_1}\ist(F_1)
\ar[d]^{k}   \\
{j_2}\ist(E_2)\,\ar[r]^{{j_2}\ist(\psi_2)}        & \,{j_2}\ist(F_2)  \\
}
$$
\end{corollary}

\section[Identification of points
on the affine line]{A nontrivial example: identification of points on the affine line}
\label{secAppliIdenti}

\vspace{3pt}

\subsec{Preliminaries}

Consider a commutative \ri $\gk$, the affine line over $\gk$
\emph{corres\-ponds to} the \klg $\gk[t]=\gB$.
Given $s$ points $\alpha_1$, \ldots, $\alpha_s$ of $\gk$ and orders of multiplicity $e_1$, \ldots, $e_s\geq 1$,
we formally define a \klg $\gA$ which represents the result of the identification of these points with the given multiplicities.
$$
\gA=\sotQ{ f\in\gB}{ f(\alpha_1)=\cdots =f(\alpha_s),\,\,
f^{[\ell]}(\alpha_i)=0,\,\ell\in\lrb{1,e_i}, \,i\in\lrbs   }
$$
In this \dfn $f^{[\ell]}$ represents the \emph{Hasse derivative} of the \pol $f(t)$, \cad $f^{[\ell]}=f^{(\ell)}\sur{\ell!}$ (formally, because the \cara can be finite). The Hasse derivatives allow us to write a Taylor formula for any \ri $\gk$.
\index{Hasse derivative}

Let $e=\sum_i e_i$, $x_0=\prod_i(t-\alpha_i)^{e_i}$ and $x_{\ell}=t^\ell\,x_0$
for $\ell\in\lrb{0..e-1}$. Suppose $e>1$ without which $\gA=\gB$.
It is clear that the $x_\ell$'s are in $\gA$.

We also assume that the $\alpha_i-\alpha_j$'s are \ivs for $i\neq j$.
We then have 
by the Chinese remainder \tho a surjective \homo

\snic{\varphi:\gB\to\prod_{i}\big(\vphantom{2^{2^2}}\aqo{\gk[t]}{(t-
\alpha_i)^{e_i}}\big)}

whose kernel is the product of the \idps  $(t-\alpha_i)^{e_i}\gB$, \cad the \idz~$x_0\gB$.

\begin{lemma}
\label{lemIdentiMil}~
\begin{enumerate}
\item $\gA$ is a \tf \klgz, more \prmtz,
$\gA=\gk[x_0,\ldots ,x_{e-1}].$
\item $\gB=\gA\oplus \bigoplus_{1\leq \ell<e}\gk \,t^\ell$ as a \kmoz.
\item 
The conductor of $\gB$ into $\gA$, $\ff=(\gA:\gB)$ is given by
$$\preskip.2em \postskip.2em
\ff=\gen{x_0}_\gB=\gen{x_0,\ldots ,x_{e-1}}_\gA .
$$
\end{enumerate}
\end{lemma}
\begin{proof}
Let $f\in\gB$, we express it \gui{in base $x_0$,}
 $f=r_0+r_1 x_0 + r_2 x_0^2+\cdots $ 
with $\deg r_i<\deg x_0 =e$.
For $i\geq 1$, by writing $r_ix_0^i=(r_i\,x_0)x_0^{i-1}$ we see that
$r_ix_0^i\in\gk[x_0,\ldots,x_{e-1}]$.
This proves that
$$\preskip.2em \postskip.1em
\ndsp\gB=\gk[x_0,\ldots,x_{e-1}]+ \left(\bigoplus_{1\leq \ell<e}\gk
\,t^\ell\right).
$$
Let $f\in\gA$ which we write $g+h$ in the previous \dcnz.. We therefore have~$h$ in $\gA$, and if $\beta $ is the common value of the $h(\alpha_i)$'s, we obtain 
the \egt  $\varphi (h-\beta)=0$. Therefore $h-\beta \in x_0\gB$, and since $h\in \bigoplus_{1\leq \ell<e}\gk \,t^\ell$ (the \kmo of the \pols of degree $<e$ and without a constant term), we 
obtain
$h-\beta=0$ then $h=\beta=0$, so $f\in\gk[x_0,\ldots ,x_{e-1}]$.
\\
In conclusion $\gA=\gk[x_0,\ldots ,x_{e-1}]$, items \emph{1}  and \emph{2}  are proven.
\\
By multiplying the \egt of item \emph{2}  by $x_0$ we obtain

\snic{x_0\gB=x_0\,\gA\;\oplus\;\bigoplus_{\ell\in\lrb{1..e-1}}x_\ell\,\gk,}

then the \egt $x_0\,\gB=\gen{x_0,\ldots ,x_{e-1}}_\gA$,
which implies $x_0\,\gB\subseteq\ff$. 
Finally, let 
 $f\in\ff$,
and so $f\in\gA$, and $f=\lambda +g$ with $\lambda \in\gk$ and $g\in
\gen{x_0,\ldots ,x_{e-1}}_\gA$. We deduce that $\lambda \in\ff$, which implies $\lambda =0$; indeed,  $\lambda t\in\gA$, if $\beta$ is the common value of the $\lambda \alpha_i$'s, we have $\varphi(\lambda t-\beta )=0$,
so $\lambda t-\beta\in x_0\,\gB$, and since~$x_0$ is a \polu of degree $\geq 2$, $\lambda =0$.
\end{proof}

\subsubsection*{A Milnor square}

In the situation described in the previous subsection we have the following Milnor square
$$\preskip-.4em \postskip.2em
\xymatrix @C=1.2cm{
\gA\,\ar@{->>}[d]\ar[r]   & \,\gB\ar@{->>}[d]^\varphi&{\kern-125pt}=\gk[t]   \\
\gk=\gA\sur{\ff}\,\,\ar@{>->}[r]^\Delta    & \gB\sur{\ff}&
{\kern-40pt}
\simeq \prod_i\left(\aqo{\gk[t]}{(t-\alpha_i)^{e_i}}\right)
\\
}
 $$
In what follows we are interested in \Amrcsz~$r$ obtained by gluing the \Bmo  $\gB^r$ and the \kmo  $\gk^r$ together using a~$(\gB\sur{\ff})$-\iso
$$\preskip.2em \postskip.4em 
h:\Delta \ist(\gk^r)\to\varphi \ist(\gB^r), 
$$
as described before \thrf{thCarreMil1}.
\\
We have denoted by $M(\gk^r,h ,\gB^r)$ such an \Amoz.

Actually, $\Delta \ist(\gk^r)$ and $\varphi \ist(\gB^r)$ are both identified with $(\gB\sur{\ff})^r$, and the \iso $h$ is identified with an \elt of

\snic{\GL_r(\gB\sur{\ff})\simeq \prod_{i=1}^s\GL_r(\aqo{\gk[t]}{(t-
\alpha_i)^{e_i}}).}

We will use these identifications in the remainder of the text without mentioning them, and,
for the sake of convenience, we will encode $h^{-1}$ (and not~$h$) by the $s$ corresponding matrices $H_i$ (with $H_i\in \GL_r(\aqo{\gk[t]}{(t-\alpha_i)^{e_i}})$).   The module  $M(\gk^r,h ,\gB^r)$  will be denoted by
$M(H_1,\ldots ,H_s)$.

In the case where the \mrcs over $\gk$ and $\gB=\gk[t]$ are always free, Milnor's \tho affirms that we thus obtain (up to \isoz) all the \mrcs $r$ over~$\gA$.

In the following subsection we give a complete description of the category of  \mrcs over $\gA$ obtained by such gluings, in a special case. The one where all the multiplicities are equal to $1$.

\subsec{Identification of points without multiplicities}

We now apply the previous conventions by supposing that the multiplicities $e_i$ are all equal to~$1$.
\begin{theorem}
\label{thIdenPtsSimples}
With the previous conventions.
\begin{enumerate}
\item The module $M(H_1,\ldots ,H_s)$, (with $H_i\in \GL_r(\aqo{\gk[t]}{t-
\alpha_i})\simeq\alb \GL_r(\gk)$) is identified with the \Asub  $M'(H_1,\ldots ,H_s)$
of $\gB^r$ consisting of the \eltsz~$f$ of~$\gB^r$ such that

\snic{\forall 1\leq i< j\leq s,\;\;H_i\cdot f(\alpha_i)=H_j\cdot f(\alpha_j).}

In particular,   $M'(H_1,\ldots,H_s)=M'(HH_1,\ldots ,HH_s)$ if $H\in\GL_r(\gk)$.
\item Let, for $i\in\lrbs$,  
\[ 
\begin{array}{ccc} 
G_i\in\alb
\GL_{r_1}(\aqo{\gk[t]}{t-\alpha_i})\simeq \GL_{r_1}(\gk)&
\hbox{ and }
\\[1mm] 
H_i\in\alb \GL_{r_2}(\aqo{\gk[t]}{t-\alpha_i})\simeq
\GL_{r_2}(\gk).    
\end{array}
\]
An \Ali $\phi$  from  $M(G_1,\ldots ,G_s)$
to $M(H_1,\ldots ,H_s)$
can be encoded by a matrix $\Phi\in \gB^{r_2\times r_1}$
satisfying, for  $1\leq i< j\leq s$,
\begin{equation}\label{eqthIdenPtsSimples}
H_i\cdot \Phi(\alpha_i)\cdot G_i^{-1} = H_j\cdot \Phi(\alpha_j)\cdot G_j^{-1}.
\end{equation}
 Such a matrix sends  $M'(G_1,\ldots ,G_s)$ to $M'(H_1,\ldots ,H_s)$.
The \Ali $\phi$ is an \iso \ssi $r_1=r_2$ and the $\Phi(\alpha_i)$'s are \ivsz.
\end{enumerate}
\end{theorem}
\begin{proof}
The first item has no incidence on the results that follow, and it is left to the reader. The second item is an immediate consequence of Lemma~\ref{lemCarreMil1} and of Corollary~\ref{corthCarreMil1}.
\end{proof}

In the following \tho we suppose that
\begin{itemize}
\item  $\gk$ is reduced,
\item  the \mrcs over $\gk[t]$ are all free,
\item  the square matrices with determinant $1$ 
are products of \elr matrices, 
\cad $\SL_n(\gk)=\En(\gk)$ for every $n$.
\end{itemize}

For example $\gk$ can be a \cdiz, a reduced \zed \ri or an integral Euclidean \riz.
Also note that if the \mrcs over $\gk[t]$ are free, it is a fortiori true for the \mrcs over $\gk$.

\begin{theorem}
\label{thIdenPtsSimples2}
For $a\in \gk$ let $J_{r,a}\eqdefi\Diag(1,\alb\ldots ,1,a)
\in\MM_r(\gk)$. Under the previous hypotheses we obtain the complete classification of the \mrcs over the \ri $\gA$ (we use the previous notations and conventions).
\begin{enumerate}
\item The modules of constant rank $M(H_1,\ldots ,H_s)$ and $M(G_1,\ldots ,G_s)$ are \isoc
\ssi $\det(H_j^{-1}\cdot H_1)=\det(G_j^{-1}\cdot G_1)$ for all~$j$.
\item Every \Amrc $r$ is \isoc to a unique module
\[\preskip-.2em \postskip.4em
M_r(a_2,\ldots ,a_s)\eqdefi M(\I_r,\alb J_{r,a_2},\ldots ,J_{r,a_s}),
\] where the $a_i$'s are in $\gk^\times$. In addition
$$\preskip.4em \postskip.4em
M_r(a_2,\ldots ,a_s)   \simeq   \Ae{r-1}\oplus  M_1(a_2,\ldots ,a_s).
$$
\item Finally, the structure of $\GKO\gA$ is specified by
$$\preskip.4em \postskip.4em
\arraycolsep2pt\begin{array}{rcl}
M_1(a_2,\ldots ,a_s)\otimes M_{1}(b_2,\ldots ,b_s)&  \simeq &  M_1(a_2b_2,\ldots
,a_sb_s)    \\[1mm]
M_1(a_2,\ldots ,a_s)\oplus M_{1}(b_2,\ldots ,b_s)&  \simeq &  \gA\oplus
M_1(a_2b_2,\ldots ,a_sb_s )
\end{array}$$
In particular, $\Pic(\gA)\simeq (\gk\eti)^{s-1}$.
\end{enumerate}

\end{theorem}
\begin{proof}
\emph{1.} In case of an \iso all the matrices in \Eqnz s~(\ref{eqthIdenPtsSimples})
are \ivsz, and it amounts to the same thing to ask

\snic{  H_j^{-1}\cdot H_1\cdot \Phi(\alpha_1) \cdot G_1^{-1}\cdot G_j =
\Phi(\alpha_j) }

for $j\in\lrb{2..s}$. Since $\Phi=\Phi(t)$ is \ivz, its \deter is an \iv \elt of $\gk[t],$ so of $\gk,$ and all the
$\det \Phi(\alpha_i)$'s are equal to~$\det \Phi$.
Consequently the two modules can only be \isoc if

\snic{  \det(H_j^{-1}\cdot H_1)=\det(G_j^{-1}\cdot G_1)}

for all $j$
(this proves in particular the uniqueness of the sequence $a_2,\ldots ,a_s$ when $M_r(a_2,\ldots ,a_s)$ is \isoc to a given \mrcz).
Conversely if this condition is satisfied, we can find an \elr matrix $\Phi$ which realizes the above conditions. It indeed suffices to have

\snic{
\Phi(\alpha_1)=\I_r$ and $\Phi(\alpha_j)= H_j^{-1}\cdot H_1 \cdot G_1^{-1}\cdot
G_j,}

 which we obtain by applying the following lemma.
\\
The end of the \dem is left to the reader. Recall%
:  if $Q=P_1\oplus P_2\simeq \gA\oplus P$ (the $P_i$'s are \prcs $1$), we have
$P\simeq \Vi_\Ae 2Q\simeq P_1\te_\gA P_2$.
\end{proof}

\vspace{-.8em}
\pagebreak

\begin{lemma}
\label{lemthIdenPtsSimples2}
Let $\alpha_1$, \ldots, $\alpha _s$ in a commutative \ri $\gk$ with the \ivs \hbox{differences $\alpha_i-\alpha_j$}  for $i\neq j$. 
Given
 $A_1$, \ldots, $A_s\in \EE_r(\gk)$, there exists a {matrix $A\in\EE_r(\gk[t])$} such {that $A(\alpha_i)=A_i$} for each~$i$.
\end{lemma}
\begin{proof}
If a matrix $A\in\EE_r(\gk[t])$ is evaluated in $s$ matrices $A_1$, \ldots, $A_s$, and a matrix $B\in\EE_r(\gk[t])$  is evaluated in $s$ matrices $B_1$, \ldots, $B_s$, then  $AB$ is evaluated \hbox{in  
$A_1B_1$,} \ldots, $A_sB_s$. Consequently, it suffices to prove the lemma when the~$A_i$'s are all equal to $\I_r$ except for one which is an \elr matrix.
In this case we can make an interpolation \`a la Lagrange since the \hbox{\elts $\alpha_i-\alpha_j$} are \ivsz.
\end{proof}

\Exercices

\vspace{-1em}
\begin{exercise}
\label{exoChapPtf2Lecteur}
{\rm  We recommend that the \dems which are not given, or are sketched, or
left to the reader,
etc, be done.
But in particular, we will cover the following cases.
\begin{itemize}
\item \label{exopropdiverschap6}
Prove Propositions~\ref{prop rgconstant local}
and~\ref{prop sfio unic}.
\item \label{exodef rank inferieur} 
Prove the \eqvcs in Proposition~\ref{def rank inferieur}~\emph{\ref{i3def rank inferieur}}.

\item Prove Corollary~\ref{lem2TransPtf}.
\item \label{exolemSomProjOrt} Prove Facts~\ref{lemSomProjOrt}
and \ref{lemProjProj}.

\end{itemize}
}
\end{exercise}

\vspace{-1em}
\begin{exercise}\label{exoleli2}
{\rm  
  Check the computations in the second local freeness lemma~\ref{leli2}.  
}
\end{exercise}

\vspace{-1em}
\begin{exercise}
\label{exo7.1} (Magic formula to diagonalize a \mprnz)  \\
{\rm Let $n$ be a fixed integer.
If $\alpha\in\cP_{n}$ (set of finite subsets of $\lrbn$), we consider the canonical \prr obtained from $\In$ by annihilating the diagonal \elts whose index is not in $\alpha$. 
We denote it by $\I_{\alpha}$.
Let $F\in\GAn(\gA)$ be a \prrz, we will explicate a family $(F_{\alpha})$ indexed by $\cP_{n}$ with matrices  satisfying the \gui{conjugations}
$$\preskip.4em \postskip-.2em
FF_{\alpha}=F_{\alpha}\I_{\alpha}\eqno(\dag)
$$
as well as the \ida
$$\preskip.2em \postskip.4em
\som_{\alpha}\det F_{\alpha}=1\eqno(\ddag)
$$
This result provides a new uniform method to explicate the local freeness of a \mptfz:
we take the \lons at the \eco $\det(F_\alpha)$, since over the \ri $\gA[1/\det(F_\alpha)]$ we have $F_{\alpha}^{-1}FF_{\alpha}=\I_{\alpha}$.\\
We will see that this is realized by the family defined as follows
$$
F_{\alpha}=F\,\I_{\alpha}+(\In-F)(\In-\I_{\alpha}).
$$
 For example if $\alpha=\lrbk$,
we have the following block \dcns
$$
\I_{\alpha}=\I_{k,n}=\bloc{\I_{k}}{0}{0}{0}, \quad F=\bloc{F_{1}}{F_{2}}{F_{3}}{F_{4}}, \quad F_{\alpha}=\bloc{F_{1}}{-F_{2}}{F_{3}}{\I_{n-k}-F_{4}}.
$$
\begin{enumerate}
\item Show $(\ddag)$. Hint: for two square matrices $A$ and $B$ of order $n$, we develop the \deter $\det(A+B)$ as a multi\lin function of the columns of $A+B$. We obtain a sum of $2^n$ \deters of matrices obtained by shuffling columns of $A$ and columns of $B$. We apply this remark with $A=F$ and $B=\In-F$.
\item If $f$ and $e$ are two \idms in a not \ncrt commutative \riz,
and if we let $f*e=fe+(1-f)(1-e)$, show that $f(f*e)=fe=(f*e)e$.
With $f=F$ and $e=\I_{\alpha}$, we obtain $f*e=F_\alpha$ which gives \egt $(\dag)$ above.
\item We now study a few \egts which make $\det F_\alpha$ intervene.
Let $\beta$ be the \cop of $\alpha$
\begin{itemize}
  \item Show that $(1-2f)(1-e-f)=(1-e-f)(1-2e)=f*e$
  \item Show that $(1-2e)^2=(1-2f)^2=1$.
  \item With $f=F$ and $e=\I_{\alpha}$, we obtain
  $(\det F_\alpha)^2=\big(\det(\I_{\beta}-F)\big)^2$.
  \item Verify that $(1-e)f(1-e)+e(1-f)e=(e-f)^2$.
  \item If we let $\mu_\alpha$ be the principal minor extracted from $F$ on the indices belonging to $\alpha$, and $\mu'_\beta$ be the principal minor extracted from $\I-F$ on the indices belonging to $\beta$, show that $(\det F_\alpha)^2=\mu_\alpha\mu'_\beta$.\\
Hint:
  for the above example with  $f=F$ and $e=\I_{\beta}$ the \egt in the previous item gives
$$
\bloc{F_{1}}{0}{0}{\I_{n-k}-F_{4}}= (\I_{\beta}-F)^2
$$
\end{itemize}
\end{enumerate}

NB. This uniform method of \din of \mprns gives a shortcut for the local freeness lemma and for the structure \tho which affirms that a \mptf is \lot free in the strong sense.
We have taken the time to prove this structure \tho twice.
Once by the \idfs in Chapter~\ref{chap ptf0},
the other more structurally, in the previous chapter. 
We hope that the readers will not hold it against us for 
   subjecting them 
to substantially less \elr \dems in the course than that of Exercise~\ref{exo7.1}.
It is because magic formulas certainly are nice things, but they sometimes hide the profound meaning of more elaborate \demsz.
}
\end{exercise}

\vspace{-1em}
\begin{exercise}\label{exo2Diag} (Generalization of the previous exercise to the \din of matrices annihilating a split \spl \polz)
\\
{\rm  Let $a$, $b$, $c\in\gA$ such that $(a-b)(a-c)(b-c)\in\Ati$, \cad the \pol 

\snic{f(T)=(T-a)(T-b)(T-c)}

is \splz, and let $A\in\Mn(\gA)$ be a matrix such that $f(A)=0$. Consider the Lagrange \pols $f_a(T)=\fraC{(T-b)(T-c)}{(a-b)(a-c)}$, \ldots\, that satisfy $f_a+f_b+f_c=1$. Let $A_a=f_a(A), \,A_b=f_b(A), \,A_c=f_c(A)$.
\begin{enumerate}\itemsep0pt
\item Show that $AA_a=aA_a$, \cad every column vector $C$ of $A_a$ satisfies $AC=aC$.
\item Deduce that if a matrix $P$ has column vectors of $A_a$ or $A_b$ or $A_c$ as its column vectors, then $AP=PD$, where $D$ is a diagonal matrix with $a$, $b$ or $c$ as its diagonal \eltsz.
\item Using $1=\det(\In)=\det(A_a+A_b+A_c)$ and using the multilinearity of the \deter as a function of the column vectors, show that there exist~$3^n$ matrices $P_i$ that satisfy
\begin{itemize}
\item $\sum_i\det(P_i)=1$.
\item In $\gA[1/\det(P_i)]$, the matrix $A$ is similar to a diagonal matrix with $a$, $b$ or $c$ as its diagonal \eltsz.
\end{itemize}
\item If the \polcar of $A$ is equal to $(T-a)^m(T-b)^p(T-c)^q$,
show that several matrices $P_i$ are null and that the sum $\sum_i\det(P_i)=1$ can be restricted to a family of matrices indexed by a finite set with $\fraC{(m+p+q)!}{m!p!q!}$ \eltsz.
\end{enumerate}
}
\end{exercise}

\vspace{-1em}
\begin{exercise} \label {exoJacobienneP2=P}
       (Jacobian of the system $P^2 - P = 0$)\\
{\rm
Let $\Mn(\gA) \to \Mn(\gA)$ be the map defined by $P \mapsto P^2 - P$. Its \dile at a point$P\in\GAn(\gA)$ is

\snic{\varphi_P : \Mn(\gA) \to \Mn(\gA),\,
  H \mapsto HP + PH - H.}

If we identify  $\Mn(\gA)$ and $\Ae {n^2}$, $\varphi_P $ is given by the Jacobian matrix at the point~$P$ of the $n^2$ \eqns $P ^2 - P  = 0$.\\
By considering

\snic{\gA = \Gn(\ZZ) = \aqo{\ZZ[(X_{ij})_{i,j \in\lrbn}]} {P ^2 - P }$  with $P  = (X_{ij}),}

by \thref{prop1TanGrassmann}, the tangent space of the affine scheme $\GAn$ at the point~$P$ is canonically identified with

\snic{\Ker \varphi_P = \sotq{H \in \Mn(\gA)}{HP + PH = H} = \Im \pi_P,}

 where
$\pi_P \in \GAn(\gA)$ is the \prr defined by

\snic{\pi_P(H) = PH(\In-P) + (\In-P)HP= PH+HP-2PHP.}

This brings us to studying the relations between $\varphi_P$ and $\pi_P$.
Illustrate what is stated regarding the Jacobian matrix and the identification of $\Mn(\gA)$ and $\Ae {n^2}$ for $n = 2$.
In \gnl show the \egts

\snic{\varphi_P \circ \pi_P= \pi_P \circ \varphi_P=0$,
 $\;(\varphi_P)^2= \In - \pi_P$,  $\;(\varphi_P)^3 = \varphi_P,}

\snic{\Ker \varphi_P = \Ker (\varphi_P)^2 = \Im \pi_P\;$
 and
$\;\Im \varphi_P = \Im (\varphi_P)^2 = \Ker \pi_P.}
}
\end{exercise}

\vspace{-1em}
\begin{exercise}\label{exopropfideles}
{\rm Prove the following local \carn of  faithful \mptfsz.
For some \Amo $P$, \propeq
\begin{enumerate}
\item [$(a)$] $P$ is \ptf and faithful.
\item [$(b)$] There exist \eco $s_i$ of $\gA$ such that each $P_{s_i}$ is free of rank $h\geq 1$ over~$\gA_{s_i}=\gA[1/s_i]$.
\item [$(c)$] $P$ is \ptf and for every \elt $s$ of $\gA$, if $P_s$ is free over the \ri $\gA_s$, it is of rank~$h\geq 1$.
\end{enumerate}
}
\end{exercise}

\vspace{-1em}
\begin{exercise}
\label{exoRangphi}
{\rm Let $\varphi:P\to Q$ be an \Ali between \mptfs and $r\in\HOp\gA$.
Express $\rg(P)\leq r$ and  $\rg(P)\geq r$ in terms of the \idds of a \mprn with image $P$.
 }
\end{exercise}

\vspace{-1em}
\pagebreak

\begin{exercise}\label{exoP1FracRat}
{(Projective line and rational fractions)}\\
{\rm  
\emph{1.}
Let $\gk$ be a \riz, $P$, $Q \in \gk[u,v]$ be two \pogs of degrees
$p, q$. Define
$$\preskip0em \postskip-.4em 
g(t) = P(t,1), \quad \wi g(t) = P(1,t), \quad
h(t) = Q(t,1), \quad \wi h(t) = Q(1,t). 
$$
\begin {enumerate}\itemsep0pt
\item [\emph {a.}]
Show that $\Res(g,p,h,q) = (-1)^{pq}\Res(\wi g,p, \wi h,q)$,
value that we denote by $\Res(P,Q)$.
\item [\emph {b.}]
Show the inclusion

\snic {
\Res(P,Q) \gen {u,v}^{p+q-1} \subseteq \gen {P,Q}
}
\end {enumerate}

\emph {2.}
Recall that $\GA_{2,1}(\gk)$ is 
the subset of $\GA_{2}(\gk)$ formed by the \prrs of rank~$1$; we have a projection $F \mapsto \Im F$ from $\GA_{2,1}(\gk)$ to $\PP^1(\gk)$.
\\
  When $\gk$ is a \cdi and $f \in \gk(t)$ is a rational fraction, we associate to~$f$ the \gui{morphism,} denoted also by $f$, $\PP^1(\gk) \vers{f} \PP^1(\gk)$, which realizes $t \mapsto f(t)$
(for the usual inclusion $\gk\subseteq\PP^1(\gk)$).\\
How do we \gnr to an arbitrary \ri $\gk$?
\\
Explain how we can lift this morphism $f$ to a \pol map, illustrated below by a doted arrow.
$$\preskip.0em \postskip.4em 
\xymatrix {
\GA_{2,1}(\gk)\ar[d]\ar[dr]\ar@{-->}[r]  &\GA_{2,1}(\gk)\ar[d]\\
\PP^1(\gk)\ar[r]^{f}                     &\PP^1(\gk) \\
} 
$$
 
 \emph {3.}
Treat the examples $f(t) = t^2$, $f(t) = t^d$ and  $f(t) = (t^2 + 1)/t^2$.
How is a homography $f(t) = {at + b \over ct +d}$ lifted ($ad-bc\in\gk\eti$)?
}

\end{exercise}

\vspace{-1em}
\begin{exercise}\label{exoConiqueFondamentale}
{(The fundamental conic or Veronese embedding $\PP^1 \to \PP^2$)}\\
{\rm  
When $\gk$ is a \cdiz, the Veronese embedding $\PP^1(\gk) \to \PP^2(\gk)$
is defined by

\snic {
(u : v) \mapsto (X=u^2 : Y=uv : Z=v^2)
.}

Its image is the \gui{fundamental conic} of $\PP^2$ with \eqn 

\snic{\left| \matrix
{X & Y\cr Y &Z\cr}\right| = XZ - Y^2 = 0.}

Analogously to Exercise~\ref{exoP1FracRat} (see also \Pbm \ref{exoVeroneseMorphism}), show that we can lift the  Veronese morphism to a \pol map, illustrated below by a dotted arrow.

\snic {
\xymatrix {
\GA_{2,1}(\gk)\ar[d]_{F \mapsto \Im F}\ar[dr]\ar@{-->}[r]  &\GA_{3,1}(\gk)\ar[d]^{F \mapsto \Im F}\\
\PP^1(\gk)\ar[r]^{\rm Veronese}             &\PP^2(\gk) \\
}}

Your obtained lift must apply to an arbitrary \ri $\gk$.
}

\end{exercise}

\vspace{-1em}
\begin{exercise}
\label{exoProjecteurCorangUn} (\Mprns of corank $1$) \,
 {\rm Let $n\geq2$.
\begin{enumerate}
\item  [\emph{1.}] Let $P \in \GA_{n,n-1}(\gA)$. Show that
 $P + \wi{P} = \In$.

\item [\emph{2.}] If $P \in \GAn(\gA)$ satisfies $P + \wi {P} = \In$, then $P$ is of rank  $n-1$.

\item [\emph{3.}] If $P \in \Mn(\gA)$ satisfies $\det(P) = 0$ and $P + \wi {P} = \In$, then $P\in \GA_{n,n-1}(\gA)$.

\end{enumerate}

}
\end{exercise}

\vspace{-1em}
\begin{exercise}\label {exoMatriceCorangUn}
{\rm
In this exercise, $A \in \Mn(\gA)$ is a matrix of {\it corank $1$}, \cad of rank $n-1$. Using Exercise~\ref{exoProjecteurCorangUn}, show the following items.
\begin{enumerate}

\item [\emph{1.}]
$\Im A = \Ker \wi {A}$ (\pro module of rank $n-1$).

\item [\emph{2.}]
$\Im \wi {A} = \Ker A$ (\pro module of rank $1$).

\item [\emph{3.}]
$\Im \tra A = \Ker \tra {\wi {A}}$ (\pro module of rank $n-1$).

\item  [\emph{4.}]
$\Im \tra {\wi {A}} = \Ker \tra A$ (\pro module of rank $1$).

\item  [\emph{5.}]
The \pro modules of rank $1$, $\Ae n  / \Im A$ and $\Ae n  / \Im
\tra A$, are duals of one another. In short, from a matrix $A$ of corank $1$, we construct two  
dual \pro
modules of rank $1$

\snic{\begin{array}{c}
\Ae n /\Im A = \Ae n /\Ker \wi {A} \simeq \Im \wi {A} = \Ker A,
\\[1mm]
\Ae n /\Im \tra {A} = \Ae n  / \Ker \tra{\wi{A}} \simeq
\Im \tra{\wi{A}} = \Ker \tra{A}.
\end{array}
}
\end{enumerate}
}
\end{exercise}

\vspace{-1em}
\begin{exercise}\label{exoIntersectionSchemasAffines}
{(Intersection of two affine schemes over $\gk$)}\\
{\rm  
This exercise belongs to the informal setting of affine schemes over a \ri $\gk$ \gui{defined} on \paref{subsecSchAff}.
Let $\gA=\kxn$, $\gB=\gk[\yn]$ be two quotient \klgs corresponding to two \syps $(\uf)$, $(\ug)$ in $\kXn$. 
\\ 
Let $A$ and $B$ be the corresponding affine schemes.
The intersection scheme $A\cap B$ is defined as being associated with the \klg $\aqo\kuX{\uf,\ug} \simeq \gA\otimes_{\kuX}\gB$
(note that the tensor product is taken over $\kuX$).

\Deuxcol{.55}{.4}
{\small 
\gui{Justify} this definition by basing yourself on the picture opposite.
\\
In a \gui{Euclidean} coordinate system, the picture includes the ellipse $\big(\fraC x a\big)^2 + y^2 = 1$, \cad  $f(x,y) = 0$ with $f = x^2 + a^2y^2 - a^2$, and the circle $g(x,y) = 0$ with $g = (x-c)^2 + y^2 - (c-a)^2$.}
{~

\centerline{\includegraphics[width=4cm]{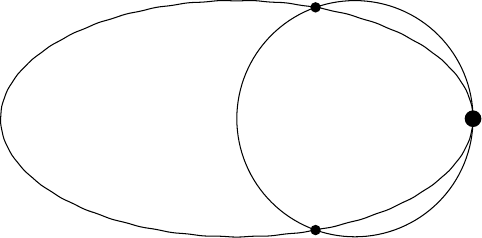}}
}
}
\end{exercise}

\vspace{-1em}
\begin{exercise}\label{exoPolPseudoUnitaire}
{(Pseudo\mon \polsz)}\\
{\rm  
Recall that a \pol $p(t) = \sum_{k\ge 0} a_kT^k\in \kT$ is said to be pseudo\mon if there exists a \sfio $(e_0,\ldots ,e_r)$ such that over each~$\gk[1/e_j]$,~$p$ is a \pol of degree $j$ with its \coe of degree $j$ being \iv (see  \paref{polpseudunit}). Such a \pol is primitive and this notion is stable under product and morphism.

\emph{1.}
Verify that $a_k = 0$ for $k > r$ and that $\geN {(1 - \som_{j>k} e_j)a_k}  = \gen {e_k}$ for $k \in \lrb{0..r}$.  In particular, $\gen {a_r} = \gen {e_r}$ and the $e_k$'s are unique or rather the \pol $\sum_k e_k X^k$ is unique (we can add null \idmsz).

\emph{2.}
Let $P = \aqo{\gA[T]}{p}$. Show that $P$ is a \ptf \Amo whose \pol rank is $\rR{P}(X) = \sum_{k=0}^r e_kX^k$; we \egmt have $\deg p = \sum_{k=1}^r k[e_k]$ (cf. item \emph {2}  of \ref{notaRangs}).
In a similar vein, see Exercise~\ref{exoPolLocUnitaire}.

}

\end{exercise}

\vspace{-1em}
\begin{exercise}\label{exoPolLocUnitaire}
 {(\Lot \mon \polsz)}\\
{\rm
\emph {1.}
Let $\fa \subseteq \gA[T]$ be an \id such that $\gA[t] = \gA[T]\sur\fa$
is a free \Amo of rank~$n$. Let $f \in \gA[T]$  be the \polcar of $t$ in $\gA[t]$. Show that $\fa = \gen {f}$. In particular, $1, t, \ldots, t^{n-1}$ is an $\gA$-basis of $\gA[t]$.

\emph {2.}
Analogous result by replacing the hypothesis \gui{$\gA[T]\sur\fa$ is a free \Amo of rank~$n$} with  \gui{$\gA[T]\sur\fa$ is a \mrc $n$.}

A \pol $f \in \gA[T]$ of degree $\leq r$ is said to be \emph{\lot \monz} if there exists a \sfio
$(e_0,\ldots,e_r)$ such that $f$ is \mon of degree $d$ in $\gA[1/e_d][T]$ for each $d\in\lrb{0..r}$.%
\index{polynomial!locally monic ---}%
\index{locally!monic polynomial}
Thus, for each $d\in\lrb{0..r}$, the \pol $f_d:= e_df$ is \mon of degree~$d$ modulo $\gen{1-e_d}$.
It is clear that this \dfn does not depend on the formal degree $r$ chosen for $f$, and that over a connected \riz, a \lot \mon \pol is \monz.

\emph {3.}
Characterize a \lot \mon \pol using its \coesz.

\emph {4.}
The \polcar of an \endo of a \mptf $M$ is \lot \mon and the corresponding \sfio is given by the $\ide_i(M)$'s.

\emph {5.}
Let $S_1$, \ldots, $S_m$ be \moco of $\gA$. Show that if $f$ is \lot \mon (for example \monz) over each $S_i^{-1}\gA$, it also is \lot \mon over $\gA$.

\emph {6.}
If $f \in \gA[T]$ is \lot \monz, show that the \ri $\gA[t] = \aqo{\gA[T]}{f}$ is a quasi-free \Amo and that $f$ is the \polcar of $t$.

\emph {7.}
Conversely, if $\fa \subseteq \gA[T]$ is an \id such that
the \ri $\gA[t] = \gA[T]\sur\fa$ is a \ptf \Amoz, then $\fa = \gen {f}$.
In particular, if a monogenic \Alg is a \ptf \Amoz, it is a quasi-free \Amoz.

\emph {8.}
For $g\in\AT$  \propeq
\begin{itemize}
\item $g$ can be written as $uf$ with $u\in\Ati$ and $f$ \lot \monz.
\item $g$ is pseudo\monz.
\item $\aqo{\AT}{g}$ is a \ptf \Amoz.
\end{itemize}

\emph {9}\eto$\!\!$.
Prove in \clama that a \pol is \lot \mon \ssi it becomes \mon after \lon at any \idepz.
}
\end{exercise}

\vspace{-1em}
\begin{exercise}\label{exoOneRankPtfIdeal}
{(Invertible modules and \mrcs $1$)}\\
{\rm  
We propose a slight variation with respect to \thref{propRgConstant3}.
\\
\emph {1.}
Let there be two commutative \ris $\gA \subseteq \gB$. 
The \Asubs of $\gB$ form a multiplicative \moz, with neutral \elt $\gA$.
Show that an \Asub $M$ of~$\gB$   \emph{invertible} in this \mo  is \tf and that for every \Asubz~$M'$ of~$\gB$ the canonical \homo $M\otimes_\gA M' \to M.M'$ is an \isoz.  Consequently, the \iv \Asubs of $\gB$ are \Amrcs $1$.
\\
\emph {2.}
Let $S\subseteq\Reg(\gA)$ be a \mo  and $\fa$ be a \lop \idz. Suppose that $S^{-1}\fa$  is an \iv \id of $S^{-1}\gA$; show that $\fa$ is an \iv \id of $\gA$. This is the case, for example, if $S^{-1}\fa$ is a free $S^{-1}\gA$-module.

}

\end{exercise}

\vspace{-1em}
\begin{exercise}\label{exoPicAPicFracA}
{(The exact sequence with $\Pic\gA$ and $\Pic\gK$, where $\gK = \Frac\gA$)}
\\
{\rm  
Let $\gA$ be a \ri and $\gK = \Frac\gA$. Define natural group morphisms

\snic {
1 \rightarrow \Ati \to \gK\eti \to \Gfr(\gA) \to \Pic\gA \to \Pic\gK
,}

and show that the obtained sequence is exact. Consequently, we have an exact sequence 
$$\preskip-.4em \postskip.4em 
1 \to \Cl(\gA) \to \Pic\gA \to \Pic\gK. 
$$
If $\Pic\gK$ is trivial, we obtain an \iso $\Cl(\gA) \simeq \Pic\gA$, and thus we once again find \thref{propRgConstant3}.

}

\end{exercise}

\vspace{-1em}
\begin{exercise}
 \label{exoAnneaudesrangs}
 {\rm
Show that $\HO\gA$ is the \ri \gui{generated by}
$\BB(\gA)$, the \agB of \idms of $\gA$, in the sense of adjoint functors.\\
More \prmtz, if $B$ is a \agBz, the \ri $\wi B$ freely generated by $B$ is given with a \homo of \agBs $\eta_B:B\to\BB(\wi B)$ such that for every \ri $\gC$ the map described below is a bijection:

\vspace{-.2em}
\Deuxcol{.6}{.38}
{
\[ 
\begin{array}{ccc} 
 \Hom_{\mathrm{Rings}}(\wi B,\gC)\lora \Hom_{\mathrm{Boolean\, alg.}}\big(B, \BB(\gC)\big)  \\[1mm] 
 \varphi\longmapsto \BB(\varphi)\circ \eta_B     
 \end{array}
\]
}
{\xymatrix @R = 0.2cm {
\wi {B}\ar[dd]_\varphi  &B\ar[dd]
  \ar[dr]^{\eta_B} \\
\ar@{~>}[r]
&  & \BB(\wi B)\ar[dl]^{\BB(\varphi)}  \\
\gC  & \BB(\gC) \\
}
}

\vspace{-.5em}
Then show that $\wi{\BB(\gA)}\simeq\HO\gA$.
}
\end{exercise}

\vspace{-1em}
\begin{exercise}
 \label{exoHOAclama}
 {\rm
Prove in \clama that $\HO(\gA)$ is canonically \isoc to the \ri of \lot constant (\cad continuous) functions from $\SpecA$ to~$\ZZ$.
}
\end{exercise}

\vspace{-1em}
\begin{exercise}\label{exoFonctDet1}  (The \deter as a functor)\\
{\rm
We have defined the \deter of an \endo of a \mptfz. We will see that more \gnlt we can define the \deter as a functor from the category of \pro modules to that of \pro modules of rank $1$.
No doubt, the simplest \dfn of the \deter of a \mptf is the following.

{\it Definition:
\vspace{-3pt}
\begin{enumerate}\itemsep0pt
\item [$(a)$] Let $M$ be a \ptf \Amo generated by $n$ \eltsz.
\\ Let $r_h=\ide_h(M)$ ($h\in\lrbzn$) and $M\ep{h}=r_hM$. Define $\det(M)$ by

\snic{ \det(M) := r_0\gA\oplus M\ep{1}\oplus \Al2M\ep{2}\oplus \cdots\oplus
 \Al nM\ep{n}
.}

We will also use the suggestive notation $\det(M) =\Al{\rg(M)}M$ by using the rank $\rg(M)=\sum_{k=1}^{n}k\,[\ide_k(M)]\in\HO\gA$.
\item [$(b)$] If $\varphi ~:M\rightarrow N$ is a \homo of \ptf \Amosz, with $s_h=\ide_h(N)$, we define $\det(\varphi)$ as a \homo of $\det(M)$ in $\det(N)$ sending $\Al hM\ep{h}$ to $\Al hN\ep{h}$ by $x\mapsto s_h (\Al h \!\varphi)(x)$.
\end{enumerate}
}

We will note that when $x\in\Al hM\ep{h}$ we have $x=r_hx$. 
\vspace{-3pt}
\begin{enumerate}\itemsep0pt
\item The module $\det(M)$ is a \mrc $1$, and we have the \egts $r_h\det(M)=\det(M)_{r_h}=\Al hM\ep{h}$.
More \gnltz, for every \idm $e$, we have $e\det(M)=\det(M_e)$.
\item  The previous \dfn provides a functor that commutes with the \lon and transforms the direct sums into tensor products.
Deduce that the functor $\det$ induces a surjective morphism from $(\KO\gA,+)$ to $\Pic\gA$.
\item  A \homo between \mptfs is an \iso \ssi its \deter is an \isoz.
\item  For an \endo of a \mptfz, the new \dfn of the \deter coincides with the previous one if we identify $\End(L)$ with~$\gA$ when $L$ is a \mrcz~$1$.
\end{enumerate}
}
\end{exercise}

\vspace{-1em}
\begin{exercise}
\label{exoFonctDet3}
{\rm
Show that, up to \isoz, the \deter functor is the only functor from the category of  \ptf \Amos in itself which possesses the following \prtsz:
\begin{itemize}\itemsep0pt
\item it transforms every arrow $\varphi:\gA\rightarrow \gA$ in itself,
\item  it transforms the direct sums into tensor products,
\item  it commutes to the \eds for every change of basis $\alpha:\gA\to\gB$.
\end{itemize}
\perso{Il semblerait que l'exo \ref{exoFonctDet3}
ait quelques rapports
with le \deter d'un complexe and les preuves t\'el\'egraphiques des
\prts of ce \deter donn\'ees in le \emph{Alg\`ebre locale,
multiplicit\'es} of J.-P. Serre.
}

}
\end{exercise}

\vspace{-1em}
\begin{exercise}
\label{exoIDDPTF1} (Determinantal \ids of a \ali between \mptfsz)\, 
{\rm
Let $\varphi : M\rightarrow N$ be a \homo between \mptfsz.
Let us write $M\oplus M'\simeq \Ae m $, $N\oplus N'\simeq \Ae n $, and \hbox{extend $\varphi$ to}  

\snic{\psi:M\oplus M'\to N\oplus N'\;$ 
with $\;\psi(x+x')=\varphi(x)\;  ( x\in M , \, x'\in M' ).}

Show that, for each integer $h$, the \idd $\cD_h(\psi)$ only depends on~$h$ and on $\varphi$. We denote it by $\cD_h(\varphi)$ and we call it \emph{the \idd of order $h$ of $\varphi$}.
}%
\index{determinantal ideals!of a \ali (\mptfsz)}%
\index{ideal!determinantal ---}
\end{exercise}

\begin{notation}\label{notaRgfi}
 {\rm  Let $r=\sum_{k=1}^{n}k\,[r_k]\in \HOp (\gA)$.
Applying the previous exercise, we call \emph{\idd of type $r$ for $\varphi$} and we denote by $\cD_r(\varphi)$ the \id

\snic{r_0\gA+r_1\cD_1(\varphi)+\cdots + r_n\cD_n(\varphi).
}

The notations $\rg(\varphi)\geq k$ and  $\rg(\varphi)\leq k$
for the \alis between free modules of finite rank are generalized as follows to the \alis between \mptfsz: let $\rg(\varphi)\geq r$  if $\cD_r(\varphi)=\gen{1}$,
$\rg(\varphi)\leq r$  if~$\cD_{1+r}(\varphi)=\gen{0}$, and
$\rg(\varphi)= r$ if $\rg(\varphi)\leq r$ and $\rg(\varphi)\geq r$.
\\
NB: see Exercise~\ref{exoLocSimpPtf}.%
\index{rank!of a linear map}
}
\end{notation}

\begin{exercise}
\label{exoIDDPTF2} 
{\rm  (Continuation of Exercise \ref{exoIDDPTF1}) \, Let $r\in\NN\etl$.
%
\begin{enumerate}
\item    If \smash{$M\vvers{\varphi}  N\vvers{\varphi'} L$} are \alis between \mptfsz, we have
$\cD_r(\varphi'\varphi)\subseteq \cD_r(\varphi')\cD_r(\varphi).$
\item   If $S$ is a \mo of $\gA$, then $\big(\cD_r(\varphi)\big)_S=\cD_{r}(\varphi_S)$.
\item   For every $s\in \gA$ such that $M_s$ and $N_s$  are free, we have $\big(\cD_r(\varphi)\big)_s=\cD_{r}(\varphi_s)$. 
In addition, this \prt characterizes the \id $\cD_r(\varphi)$.
\end{enumerate}
Let $r=\sum_{k\in\lrbn} k [r_k]\in\HOp\gA$. 
\begin{enumerate}\setcounter{enumi}{3}
\item [\emph{4.}] Redo the previous items of the exercise in this new context.

\end{enumerate}

}
\end{exercise}

\vspace{-1.1em}
\pagebreak

\begin{exercise}
 \label{exoLocSimpPtf}
 {\rm  (With notations \ref{notaRgfi}) \,
 Let $\varphi:M\to N$ be a \ali between \ptf \Amosz.
Prove that the following \prts are equivalent.
\\
 \emph{1.} $\varphi$ is \lnlz.
\\
 \emph{2.} $\varphi$ has a well-defined rank in $\HOp(\gA)$.
\\
 \emph{3.} After \lon at \ecoz, the modules are free and the \ali is \nlz.
 } 
\end{exercise}

\vspace{-1em}
\begin{exercise}
\label{exoCoMatCoRang1}
{\rm
Let $A \in \Ae {n \times m}$; if $A$ is of rank $m-1$, we can explicate a finite \sys of \gtrs of the submodule $\Ker A \subseteq \Ae n$ without using neither an \egt test nor a membership test.
In fact, under the only (weaker) hypothesis $n \ge m-1$, we uniformly define a matrix $A' \in \Ae {m \times N}$ with $N = {n \choose m-1}$ which is \gui{a kind of comatrix of $\gA$.} This matrix satisfies $\Im A' \subseteq \Ker A$ as soon as~$A$ is of rank $\leq m-1$, with
\egt when $A$ is of rank $m-1$.

We can define $A' \in \Ae {m \times N}$ via the exterior \algz: we see $A$ as a \ali $u : \Ae m \to \Ae n$
and we consider $u' = \Al {m-1}(\tra u) : \Al {m-1}(\Ae n) \to \Al {m-1}(\Ae m)$. In the canonical bases,
$\Al {m-1}(\Ae n) = \Ae N$ and $\Al {m-1}(\Ae m) = \Ae m$, so $u'$ is represented by a matrix $A' \in \Ae {m \times N}$. To explicate this matrix $A'$, we order the set of $N = {n \choose m-1}$ subsets $I$ of $\lrb{1..n}$ of cardinality $m-1$ such that their complements are in increasing lexicographic order; the columns of $A'$ are indexed by this set of subsets, as follows

\snic {
a'_{j,I} = (-1)^{k_I+j} \det(A_{I, \{1..m\} \setminus \{j\}}), \qquad
\hbox {$k_I$ being the number of $I$.}
}

For example, if $m = 2$, then $N = n$, and $A'=\cmatrix {
a_{n,2}  & -a_{n-1,2} & \cdots & \pm a_{1,2} \cr
-a_{n,1} & a_{n-1,1} & \cdots & \mp  a_{1,1}
}.$

\emph{1.}
For $i \in \lrb{1..n}$, we have $(AA')_{i,I} = (-1)^{k_I+1} \det(A_{\{i\} \cup
I, \{1..m\}})$.  In particular, if~$\cD_m(A) = 0$, then $AA' = 0$.

\emph{2.}
If $n = m$, then $A' = \wi{A}$ (the comatrix of $A$).

\emph{3.}
If $A$ is of rank $m-1$, then $\Im A' = \Ker A$; in particular, $A'$ is of rank $1$.

\emph{4.}
Every \stl module of rank $1$ is free. We will be able to compare with Fact~\ref{factPicStab} and with Exercise~\ref{exoStabLibRang1}.

\emph{5.} If $B$ is a matrix satisfying $ABA=A$, then $P=BA$ is a
\mprn satisfying $\Im(\In-P)=\Ker P=\Ker A$. This provides another way to answer the question: give a finite \sys of \gtrs of $\Ker A$.
Compare this other solution to that of the current exercise.
To compute the matrix $P$, we will be able to use the method explained in Section~\ref{secCramer} (\thref{propIGCram}). 
Another method, considerably more economical, can be found in \cite[D\'{\i}az-Toca\&al.]{DiGLQ} (based on \cite[Mulmuley]{Mul}).

}

\end{exercise}

\vspace{-1em}
\begin{exercise}
\label{exoVarProj} (\Hmg \pols and $\Pn(\gk)$)
\\ 
{\rm Let $(\lfs)=(\uf)$ in $\gk[\Xzn]
$ be a \hmg \sypz. We seek to define the zeros of $(\uf)$ in $\Pn(\gk)$. 
Let $P$ be a point of $\Pn(\gk)$, \cad 
a \pro \kmo of rank $1$ which is a direct summand in
$\gk^{n+1}$. Show that if a \sgr of $P$ annihilates $(\uf)$, then every \elt of $P$
annihilates~$(\uf)$.
 
}
\end{exercise}

\vspace{-1.3em}
\pagebreak

\begin{exercise}\label{exoGLnTangent}
 {(Tangent space to $\GLn$)}\\
{\rm Determine the tangent space at a point to the functor $\gk \mapsto \GLn(\gk)$.
}

\end{exercise}

\vspace{-1em}
\begin{exercise}\label{exoSLnTangent} {(Tangent space to $\SLn$)}
\\
{\rm  
Determine the tangent space at a point to the functor $\gk \mapsto \SLn(\gk)$.
}

\end{exercise}

\vspace{-1em}
\begin{exercise} \label{exoTangentJ0ConeNilpotent}
       (Tangent space at $J_0$ to the nilpotent cone)  \,  {\rm Let $\gk$ be a \riz.\\
Let $(e_{ij})_{i,j\in\lrbn}$ be the canonical basis of $\Mn(\gk)$ and $J_0 \in \MMn(\gk)$ be the standard Jordan matrix. For example, for $n = 3$, $J_0 =\cmatrix {0 &1 &0\cr 0 & 0& 1\cr 0& 0& 0\cr}$.

\emph{1.}
We define $\varphi : \Mn(\gk) \to \Mn(\gk)$ by 
$\varphi(H) = \sum_{i+j = n-1} J_0^i H J_0^j$.\\ 
Determine $\Im \varphi$.

\emph{2.}
Give a direct complement of $\Im \varphi$ in $\Mn(\gk)$,
         then give $\psi : \Mn(\gk) \to \Mn(\gk)$ satisfying $\varphi \circ
\psi \circ \varphi = \varphi$. 
Show that  $\Ker\varphi$ is free of rank $n^2-n$ and give a basis of this module. 

\emph{3.}
Consider the functor $\gk \mapsto \{N \in \Mn(\gk) \ | \
N^n = 0\}$. Determine the tangent space at $J_0$ to this functor.

}
\end{exercise}

\vspace{-1em}
\begin{exercise} \label{exoprop1TanGrassmann} (Complement of \thref{prop1TanGrassmann})
{\rm Let $\gA[\vep]=\AT/\geN{T^2}$. \\
Let $P$, $H \in \MM_n(\gA)$.
Show that the matrix $P + \varepsilon H$ is \idme \ssi

\snic{P^2 = P  \  $ and $ \  H = HP + PH.}

Generalize to an abstract noncommutative \ri with an \idm $\varepsilon$ in the center of the \riz.

\comm The example of the \ri $\Mn(\gA)$ shows that in the noncommutative case the situation for the \idms is quite different from the one in the commutative case where $\BB(\gA)=\BB(\Ared)$
(Corollary~\ref{corIdmNewton})
and where the \idms are \gui{isolated} (Lemma~\ref{lemIdmIsoles}).}\eoe
\end{exercise}


\vspace{-1em}
\begin{problem}\label{exoAnneauCercle1} {(The \ri of the circle)}
\\
{\rm
Let $\gk$ be a \cdi of \cara $\ne 2$, $f(X,Y) = X^2 + Y^2 - 1 \in \gk[X,Y]$. It is an \ird and smooth \polz, \cad $1 \in \gen {f, \Dpp{f}{X}, \Dpp{ f}{ Y}}$ (explicitly, we~have~$-2 = 2f -
X\Dpp{f}{X} - Y\Dpp{f}{Y}$). \\ 
It is therefore licit to think that the \ri $\gA = \gk[X,Y]/\gen {f} = \gk[x,y]$ is an integral \adpz. This will be proven in \Pbmz~\ref{exoArithInvariantRing} (item~\emph{4}). 
\\ 
Let $\gK$ be its field of fractions and let $t = {y \over x-1} \in \gK$.
\begin{enumerate}\itemsep0pt
\item
Show that $\gK = \gk(t)$; \gmqt justify how to find $t$ (parameterization of a conic having a $\gk$-rational point)
and make $x$, $y$ explicit in terms of~$t$.
\item
Let $u = (1+t^2)^{-1}$, $v = tu$.  Verify that the \cli of $\gk[u]$ in $\gK = \gk(t)$~is

\snic {
\gk[x,y] = \gk[u, v] =
\sotq {h(t)/(1 + t^2)^s} {h \in \gk[t],\  \deg(h) \le 2s}
.}

In particular, $\gA = \gk[x,y]$ is \iclz.  Explain how the $\gk$-circle
$x^2 + y^2 = 1$ is the \pro line $\PP^1(\gk)$ deprived \gui{of the $\gk$-point} $(x,y) = (1, \pm i)$.

\item
If $-1$ is a square in $\gk$, show that $\gk[x,y]$ is a localized \ri
 $\gk[w,w^{-1}]$ (for some $w$ to explicate)
of a \pol \ri over $\gk$, therefore a Bézout \riz.

\item
Let $P_0 = (x_0,y_0)$ be a $\gk$-point of the circle $x^2 + y^2 = 1$ and $\gen {x-x_0, y-y_0} \subseteq \gA$ be its \idemaz.  Verify that $\gen {x-x_0, y-y_0}^2$ is a \idp of \gtr $xx_0 + yy_0 - 1$.
\Gmq interpretation of $xx_0 + yy_0 - 1$?

\item
\mbox{\parbox[t]{.6\linewidth}{Here $(x_0, y_0) = (1,0)$. Describe the computations allowing to explicate the (projection) matrix\hsd
 $~P = {1\over 2} \cmatrix {1-x & -y\cr -y & 1+x}$ 
\\
as a \mlp for the pair $(x-1, y)$. The exact sequence

\snic {
\Ae 2 \vvvers{\I_2-P} {\Ae 2}\vvvvers{(x-1, y)} \gen {x-1, y}\to 0
}

allows us to realize the (\ivz) \id of the point~$(1,0)$ as the image of the \prr $P$ of rank~$1$.
\\
Comment on the opposite picture which is its \gmq counterpart (vector line bundle of the circle).}%
\hspace*{.05\linewidth}%
\parbox[t]{.35\linewidth}{\begin{center}
\includegraphics*[width=3.5cm]{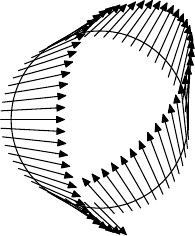} $\qquad\qquad$
\end{center}}}
\item
Suppose that $-1$ is not a square in $\gk$
and that we see $\gk[x,y]$ as  
a free $\gk[x]$-\alg of rank $2$, with basis $(1,y)$.
Explicate the norm and verify, for
$z = a(x) + b(x)y \ne 0$, the \egt

\snic {
\deg \rN_{\gk[x,y]/\gk[x]} (z) = 2\max(\deg a, 1 + \deg b).
}

 In particular, $\deg \rN_{\gk[x,y]/\gk[x]} (z)$ is even.  Deduce the group 
$\gk[x,y]^{\!\times}$, the fact that $y$ and $1 \pm x$ are \irds in $\gk[x,y]$, and that the \idz~$\gen {x-1, y}$ of the point~$(1,0)$ is not principal (\cad the line bundle above is not trivial).
\end{enumerate}
}
\end{problem}

\vspace{-1em}
\begin{problem}\label{exoLambdaGammaK0}
{(The operations $\lambda_t$ and $\gamma_t$ over $\KO(\gA)$)} \\
{\rm  
If $P$ is a \mptf over $\gA$, let, for $n\in\NN$, $\lambda^n (P)$ or $\lambda^n ([P])$ be the
class of~$\Al{n}P$ in~$\KO(\gA)$ and we have the fundamental \egt 

\snic{\qquad\lambda^n(P\oplus Q) =
\sum_{p+q = n}\lambda^p(P)\lambda^q(Q).\qquad(*)}

We \egmt define the \pol $\lambda_t(P) \in \KO(A)[t]$ by $\lambda_t(P) = \sum_{n\ge 0} \lambda^n(P) t^n$. It is a  \pol of constant term $1$ that we consider in the \ri of formal series $\KO(A)[[t]]$. Then

\snic{\lambda_t(P) \,\in\, 1 + t\KO(\gA)[[t]]\;\subseteq\;(\KO(\gA)[[t]])\eti.}

By $(*)$ we have $\lambda_t(P \oplus Q) = \lambda_t(P)\lambda_t(Q)$, which allows us to extend $\lambda_t$ to a morphism $(\KO(\gA), +) \to (1 + t\KO(\gA)[[t]], \times)$. Thus if $P$, $Q$ are two \mptfsz, for $x = [P] - [Q]$, we have by \dfn
$$
\lambda_t(x) = {\lambda_t(P) \over  \lambda_t(Q)} =
{1 + \lambda^1(P)t + \lambda^2(P)t^2 + \cdots \over
1 + \lambda^1(Q)t + \lambda^2(Q)t^2 + \cdots}
$$
sequence that we will denote by $\sum_{n\ge 0} \lambda^n(x)t^n$, with $\lambda^0(x) = 1$, $\lambda^1(x) = x$. 

Grothendieck has \egmt defined over $\KO(\gA)$ another operation $\gamma_t$ by the \egt
$$\preskip.4em \postskip.4em 
\gamma_t(x) = \lambda_{t/(1-t)}(x), 
$$
for $x\in\KO(\gA)$. This is licit because the multiplicative subgroup $1 + t\KO(\gA)[[t]]$ is stable by the substitution $t \leftarrow t/(1-t)$. 
This substitution $t \leftarrow t/(1-t)$ leaves the term invariant in $t$, let 

\snic{\gamma_t(x)\,=\,1+tx+t^2\big(x+\lambda^2(x)\big)+\cdots\,=\, \sum_{n\ge 0}
\gamma^n(x)t^n.}

\emph {1.} Give $\lambda_t(p)$ and $\gamma_t(p)$ for $p\in\NN\etl$.
Let  
$x \in\KTO\gA$. 
Show that $\gamma_t(x)$ is a \polz~$t$. 
By using $\gamma_{t}(-x)$, deduce that $x$ is nilpotent.

\emph {2.}
Show that $\KTO(\gA)$ is the nilradical of the \ri $\KO(\gA)$.

 We have $\rg\big(\lambda^n(x)\big)=\lambda^n(\rg x)$ and thus we dispose of a formal sequence $\rgl_t(x)$
with \coes in $\HO\gA$ defined by 

\snic{\rgl_t(x) = \lambda_t(\rg x)=\sum_{n\ge 0}
\rg\big(\lambda^n(x)\big)t^n.}

If $x\in\HO\gA$ this simply gives $\rgl_t(x)=\lambda_t(x)$.
 
\emph {3.}  If $x = [P]$, recall that $(1+t)^{\rg x} = \rR{P}(1+t) = \rF{\Id_P}(t)\in\BB(\gA)[t]$.
Show that, when we identify $\BB(\gA)$ with $\BB(\HO\gA)$ by letting $e=[e\gA]=[e]$ for~$e\in\BB(\gA)$, we obtain $\rgl_t(x) = 1+t\rg x=(1+t)^{\rg x}$ if $0\leq \rg x \leq 1$.\\
Then show that $\rgl_t(x) = (1+t)^{\rg x}$ for every $x \in \KO(\gA)$.

\emph {4.}
We define $\rgg_t(x) = \gamma_t(\rg x) = \sum_{n\ge 0} \rg\big(\gamma^n(x)\big)t^n$.
\\
Show that $\rgg_t(x) = (1-t)^{-\rg x}$ for every $x \in \KO(\gA)$, or for $x = [P]$ that $\rgg_t(x) = \rR{P}\big(1/(1-t)\big) = \rR{P}(1-t)^{-1}$. In addition, if $0\leq\rg x \le 1$, we obtain the \egt $\rgg_t(x) = 1+xt/(1-t) = 1 + xt + xt^2 + \dots$

\emph {5.}
For all $x$ of $\KO\gA$,  ${\gamma_t(x)}{(1-t)^{\rg(x)}}$ is a \polz.

\emph {6.}
Show the reciprocity formulas between $\lambda^n$ and $\gamma^n$ for $n \ge 1$

\snic {
\gamma^n(x) = \sum_{p=0}^{n-1} {n-1 \choose p} \lambda^{p+1}(x),
\qquad
\lambda^n(x) = \sum_{q=0}^{n-1} {n-1 \choose q} (-1)^{n-1-q}\gamma^{q+1}(x).
}

}

\end{problem}

\vspace{-1em}
\begin{problem}\label{exoApplicationProjectiveNoether}
 {(The \pro map of \Noe and the \mrcs$1$ direct summands in $\gk^2$)}\\
 {\rm  
Fix a \ri $\gk$, two \idtrs $X$, $Y$ over $\gk$ and an integer $n \ge1$. 
Given two $n$-sequences of \elts of $\gk$, $x = (\xn)$ and $y = (\yn)$, we associate with them an $(n+1)$-sequence $z = z(x,y) = (z_0, \ldots, z_n)$ as follows

\snic {
\prod_{i=1}^n (x_i X + y_i Y) = 
z_0X^n + z_1X^{n-1}Y + \cdots + z_{n-1}XY^{n-1} + z_nY^n
.}

Thus, we have $z_0 = x_1\cdots x_n$, $z_n = y_1\cdots y_n$, 
and for example, for $n = 3$,

\snic {
z_1 = x_1x_2y_3 + x_1x_3y_2 + x_2x_3y_1,\quad
z_2 = x_1y_2y_3 + x_2y_1y_3 + x_3y_1y_2
.}

For $d \in \lrb{0..n}$, we easily check that $z_d(y,x) = z_{n-d}(x,y)$ and that we have the following formal expression thanks to the \elr \smq functions with $n$ \idtrs $(S_0 = 1, S_1, \ldots, S_n)$:

\snic {
z_d = x_1\cdots x_n S_d(y_1/x_1, \ldots, y_n/x_n)
.}

In particular, $z_d$ is \hmg in $x$ of degree $n-d$, and \hmg in $y$ of degree~$d$. We can give a direct \dfn of $z_d$ as follows

\snic {
z_d = \sum_{\#I = n-d} \prod_{i \in I} x_i \prod_{j \in \lrbn\setminus I} y_j 
.}

If $\gk$ is a \cdiz, we have a map $\psi : (\PP^1)^n =
\PP^1 \times \cdots \times \PP^1 \to \PP^n$, said to be \Noez ian, defined by
$$
\psi : \big((x_1 : y_1), \dots, (x_n : y_n)\big) \mapsto (z_0 : \cdots : z_n)
\leqno (\star)
$$
We make the \smq group $\rS_n$ act on the product $(\PP^1)^n$ by permutation of the \coosz; then the map $(\star)$ above, which is $\rS_n$-invariant, intervenes in \agq \gmt to make $(\PP^1)^n /\rS_n$ and
$\PP^n$ isomorphic. \perso{voir Shafarevich vol 2 exo 2 page 232}

\emph {1.}
Show that for $P_1$, \ldots, $P_n$, $Q_1$, \ldots, $Q_n$ in $\PP^1$,
we have 

\snic{\psi(P_1, \ldots, P_n) = \psi(Q_1, \ldots, Q_n) \iff 
(Q_1, \ldots, Q_n)$ is a permutation \hspace*{3.2cm} of $(P_1, \ldots, P_n).}

We now want, $\gk$ being an arbitrary \riz, to formulate the map $(\star)$ in terms of \kmrcs$1$.  
\\
Precisely, let $L = \gk X\oplus \gk Y \simeq \gk^2$, and let

\snic{S_n(L)
= \gk X^n \oplus \gk X^{n-1}Y\oplus \cdots \oplus \gk XY^{n-1} \oplus \gk Y^n
\simeq \gk^{n+1}}

be the \hmg component 
of degree $n$ of $\gk[X,Y]$.  If $P_1$, \ldots, $P_n \subset L$ are $n$ \pro \ksmos of constant rank $1$ which are direct summands, 
we want to associate with them, functorially, a \ksmo $P = \psi(P_1, \ldots,
P_n)$ of $S_n(L)$, \prc$1$ and a direct summand. 
 Of course, we must have

\snic{\psi(P_1, \ldots, P_n) = \psi(P_{\sigma(1)}, \ldots, P_{\sigma(n)})}

for every permutation $\sigma \in \rS_n$.  In addition, if each $P_i$ is free with basis $x_iX + y_iY$, then $P$ must be free with basis $\sum_{i=0}^n z_i X^{n-i}Y^i$, in order to find $(\star)$.

\emph {2.}
Show that if each $(x_i, y_i)$ is \umdz, the same holds for $(z_0, \ldots, z_n)$.

\emph {3.}
Define $\psi(P_1, \ldots, P_n) \subset S_n(L)$ thanks to the module $P_1 \te_\gk \cdots \te_\gk P_n$ and to the $\gk$-\lin map $\pi : L^{n\te} \twoheadrightarrow S_n(L)$,

\snic {
\pi : \bigotimes_{i=1}^n (x_i X + y_iY) \longmapsto
\prod_{i=1}^n (x_i X + y_iY)
.}

\emph {4.}
Let $\gk[\uZ] = \gk[Z_0, \ldots, Z_n]$, $\gk[\uX, \uY] = \gk[X_1,Y_1, \ldots, X_n,Y_n]$. \\
What to say about the $\gk$-morphism $\varphi : \gk[\uZ] \to \gk[\uX, \uY]$ defined by
$$\preskip.4em \postskip.4em \ndsp
Z_d \longmapsto z_d = \sum_{\#I = n-d} \prod_{i \in I} X_i \prod_{j \in \lrbn\setminus I} Y_j 
~~~? 
$$
NB: $\varphi$ is the co-morphism of $\psi$.

}

\end{problem}


\vspace{-1em}
\begin{problem}\label{exoTh90HilbertMultiplicatif} 
(Hilbert's \tho 90, multiplicative form)
 \\
{\rm
Let $G$ be a finite group acting on a commutative \ri $\gB$; a \emph{$1$-cocycle} of $G$ over $\Bti$ is a family $(c_\sigma)_{\sigma \in G}$ such that $c_{\sigma\tau} = c_\sigma \sigma(c_\tau)$; consequently, $c_\Id = 1$. For every \elt $b \in \Bti$, $(\sigma(b)b^{-1})_{\sigma \in G}$ is a $1$-cocycle called a \emph{$1$-coboundary}. \\
Let $\zcoho$ be the set of $1$-cocycles of $G$ over $\Bti$; it is a subgroup of the (commutative) group of all the maps of $G$ in $\Bti$ equipped with the final product.
The map $\Bti \to \zcoho$, $b \mapsto (\sigma(b)b^{-1})_{\sigma \in G}$, is a morphism; let $\bcoho$ be its image and we define \emph{the first group of cohomology of~$G$ over $\Bti$} 

\snic {
\hcoho = \zcoho\sur\bcoho 
.}

Finally, we define the (\gnlt noncommutative) \ri  $\tgaBG$ as being \hbox{the \Bmoz} with basis $G$,
equipped with the product $(b\sigma) \cdot (b'\sigma') = b\sigma(b') \sigma\sigma'$.
Then $\gB$ becomes a $\tgaBG$-\alg via $(\sum_{\sigma} b_\sigma \sigma) \cdot
b = \sum_{\sigma} b_\sigma \sigma(b)$.
\\
We call $\tgaBG$  \emph{the twisted group \alg of the group $G$}.

Let $(\gA,\gB,G)$ be a \aGz. The aim of the \pb is to associate with every $1$-cocycle $c = (c_\sigma)_{\sigma \in G}$ a \Amrc $1$ denoted by $\gB_c^G$ and to show that $c \mapsto \gB_c^G$ defines an injective morphism of $\hcoho$ in $\Pic(\gA)$. In particular, if $\Pic(\gA)$ is trivial, then every $1$-cocycle of $G$ over $\Bti$ is a coboundary. 

\emph {1.}
Show that $\tgaBG \to \End_\gA(\gB)$, $\sigma \mapsto \sigma$ is an \iso of \Algsz.

\emph {2.}
Let $c \in \zcoho$. We define $\theta_c : \tgaBG \to \tgaBG$
by $\theta_c(b\sigma) = bc_\sigma \sigma$. 
\begin{itemize}
\item [\emph {a.}]
Verify $\theta_c \circ \theta_{d} = \theta_{cd}$; deduce that $\theta_c$ is an $\gA$-\auto of $\tgaBG$.
\item [\emph {b.}]
Show that if $c \in \bcoho$, then $\theta_c$ is an interior \autoz.
\end{itemize}

\emph {3.}
Let $c \in \zcoho$. Consider the action from $\tgaBG$ to $\gB$ \gui{twisted} by $\theta_c$, \cad $z \cdot b = \theta_c(z)\, b$; let $\gB_c$ be this $\tgaBG$-module, $\gB_c^G$ be the set of \elts of~$\gB$ invariant under $G$ (for this action twisted by $\theta_c$), and

\snic{
\pi_c = \sum_{\sigma \in G} c_\sigma\, \sigma \in \End_\gA(\gB)
.}

Verify that $\gB_c^G$ is an \Asub of $\gB$.  Show that $\pi_c$ is a surjection from $\gB$ to $\gB_c^G$ by explicating a section; deduce that $\gB_c^G$ is a direct summand in $\gB$ (as an \Amoz).

\emph {4.}
We will show that for every $c \in \zcoho$, $\gB_c^G$ is a \Amrcz~$1$.
\begin{enumerate}
\item [\emph {a.}]
Verify that $\gB_c^G \gB_{d}^G \simeq \gB_{cd}^G$ and $\gB_c^G \otimes_\gA
\gB_{d}^G \simeq \gB_{cd}^G$.
\item [\emph {b.}]
Show that if $c \in \bcoho$, then $\gB_c^G \simeq \gA$.
Conclude the result.
\item [\emph {c.}]
Show that $c \mapsto \gB_c^G$ induces an injective morphism from $\hcoho$ into~$\Pic(\gA)$. 
\end{enumerate}

\emph {5.}
In the case where $\gA$ is a \zed \ri (for example a \cdiz), show that every $1$-cocycle $(c_\sigma)_{\sigma \in G}$ is the coboundary of some $b \in \Bti$.

\emph {6.}
Suppose that $G$ is cyclic of order $n$, $G = \gen {\sigma}$, 
and that $\Pic(\gA) = 0$. \\
Let $x \in B$; show that $\rN\iBA(x) = 1$ \ssi there exists a $b \in \Bti$ such \hbox{that $x = \sigma(b)/b$}.

}
\end{problem}

\vspace{-1em}
\begin{problem}\label{exoSegreMorphism} {(The Segre morphism in a special case)}
\\
{\rm
Let $\gA[\uX,\uY] = \gA[X_1, \ldots, X_n, Y_1, \ldots, Y_n]$. Consider the \id $\fa = \gen {X_iY_j - X_jY_i}$, \cad the \id $\cD_2(A)$, where $A$ is the \gnq matrix \halfsmashbot{$\cmatrix {X_1 & X_2 & \cdots & X_n\cr Y_1 & Y_2 & \cdots & Y_n\cr}$}. We want to show that $\fa$ is the kernel of the morphism
$$
\varphi : \gA[\uX, \uY] \to \gA[T,U,\uZ] = \gA[T,U, Z_1, \ldots, Z_n],
\quad X_i\to TZ_i, \;\; Y_i\to UZ_i,
$$
where $T$, $U$, $Z_1$, \ldots, $Z_n$ are new \idtrsz. Let us agree to say that a \hbox{\mom $m \in \gA[\uX,\uY]$} is normalized if $m$ is equal to $X_{i_1}\cdots X_{i_r} Y_{j_1} \cdots Y_{j_s}$ 
{with $1 \le i_1 \le \cdots \le i_r \le j_1 \le \cdots \le j_s \le n$} (the indices of $\uX$ are smaller than that of $\uY$). Let $\fa_{\rm nor}$ be the \Asub of $\gA[\uX, \uY]$ generated by the normalized \momsz.
\begin{enumerate}
\item
If $m$, $m'$ are normalized, show that $\varphi(m) = \varphi(m') \Rightarrow m = m'$.
Deduce that $\Ker\varphi \cap \fa_{\rm nor} = \{0\}$.

\item
Show that we have a direct sum of \Amosz:
$\gA[\uX,\uY] = \fa \oplus \fa_{\rm nor}$

\item
Deduce that $\fa = \Ker\varphi$.
In particular, if $\gA$ is reduced (resp.\,\sdzz), then
$\fa$ is radical (resp.\,prime).

\end{enumerate}
}
\end{problem}

\vspace{-1em}
\comm
The morphism $\varphi$ induces, by co-morphism, a morphism between affine spaces
$$\preskip-.4em \postskip.4em
\psi : \AA^2(\gA) \times \An(\gA) \to \MM_{2,n}(\gA) \simeq \AA^{2n}(\gA),\; \big((t,u),z\big)
\mapsto \cmatrix {tz_1 & \cdots & tz_n\cr uz_1 & \cdots & uz_n\cr}.
$$
If $\gA$ is a field, the image of $\psi$ is 
the zero set 
$\cZ(\fa)$, and $\psi$ induces at the projective spaces level an inclusion $\PP^1(\gA) \times \PP^{n-1}(\gA) \to \PP^{2n-1}(\gA)$ (called \gui{embedding}).
\\
More \gnltz, by completely changing the notations, with \idtrs $X_1$, \ldots, $X_n$, $Y_1$, \ldots, $Y_m$,
$Z_{ij}$, $ i \in\lrbn$, $j \in\lrbm$, consider the morphism $\varphi : \gA[\uZ] \to \gA[\uX, \uY]$, $Z_{ij} \to X_iY_j$.
We show \hbox{that $\Ker\varphi = \cD_2(A)$} where~\hbox{$A \in \MM_{n,m}(\gA[\uZ])$}
is the \gnq matrix.
The morphism $\varphi$ induces, by co-morphism, a morphism between affine spaces 
$$\preskip.4em \postskip.4em 
\psi : \AA^n(\gA) \times \AA^m(\gA) \to \MM_{n,m}(\gA) \simeq \AA^{nm}(\gA), \;\big((x_i)_i,(y_j)_j\big)
\mapsto (x_iy_j)_{ij}, 
$$
whose image is the 
zero set $\cZ\big(\cD_2(A)\big)$. If $\gA$ is a \cdiz,  
$\psi$ induces an injection 
$\PP^{n-1}(\gA) \times \PP^{m-1}(\gA) \to \PP^{nm-1}(\gA)$: it is the Segre embedding. This allows us to realise $\PP^{n-1} \times \PP^{m-1}$ \emph{as a projective \agq sub\vrt of} $\PP^{nm-1}$  (in a precise sense that we do not specify here).
\\
If $\gA$ is arbitrary, let $E \in \PP^{n-1}(\gA)$,
$F \in \PP^{m-1}(\gA)$; $E$ is thus a direct summand  in $\Ae n$,
of rank~$1$; similarly for $F$. Then $E\otimes_\gA F$ is canonically identified with
a submodule of $\Ae n \otimes_\gA \Ae m \simeq \Ae {nm}$, a direct summand,
of rank~$1$. By letting $\psi(E, F) = E \otimes_\gA F$, we thus obtain a map from~\hbox{$\PP^{n-1}(\gA) \times \PP^{m-1}(\gA)$} to~\hbox{$\PP^{nm-1}(\gA)$} which \gui{extends} the map previously defined: if $x \in \Ae n$, $y \in \Ae m$ are unimodular, the same holds for~\hbox{$x \otimes y \in \Ae n \otimes_\gA \Ae m$}, and by letting $E = \gA x$, $F = \gA y$, we have $E \otimes_\gA F = \gA (x \otimes y)$.
\eoe

\vspace{-1em}
\begin{problem}\label{exoVeroneseMorphism} {(The Veronese morphism in a special case)}
\\
{\rm
Let $d \ge 1$, $\AuX = \gA[X_0, \ldots, X_d]$ and $\fa = \gen {X_iX_j -
X_kX_\ell, i+j = k+\ell}$. We will show that the \id $\fa$ is the kernel of the morphism

\snic {
\varphi : \AuX \to \gA[U,V], \qquad \varphi(X_i) = U^{d-i} V^i
.}

where $U$, $V$ are two new \idtrsz. We define another \id~$\fb$

\snic {
\fb = \gen {X_iX_j - X_{i-1}X_{j+1}, 1 \le i \le j \le d-1}
}

\vspace{-.8em}
\begin{enumerate}
\item
Show that

\snic {
\Ker\varphi \cap (\gA[X_0,X_d] + \gA[X_0,X_d]X_1 + \cdots +
\gA[X_0,X_d]X_{d-1}) = \{0\}
}

\item
Show that we have a direct sum of \Amosz

\snic {
\AuX = \fb \oplus \gA[X_0,X_d] \oplus \gA[X_0,X_d]X_1
\oplus \cdots \oplus \gA[X_0,X_d]X_{d-1}
}

\item
Deduce that $\fa = \fb = \Ker\varphi$.
In particular, if $\gA$ is reduced (resp.\,\sdzz), then
$\fa$ is radical (resp.\,prime).

\end{enumerate}
}
\end{problem}

\vspace{-1em}
\comm
More \gnltz, let $N = {n+d \choose d} = {n+d \choose n}$ and $n +
1 + N$ \idtrs $U_0$, \ldots, $U_n$, $(X_\alpha)_\alpha$,  where the indices $\alpha \in \NN^{n+1}$ are such that $|\alpha| = d$. We dispose of a morphism $\varphi : \gA[\uX] \to \gA[\uU]$, $X_\alpha \mapsto \uU^\alpha$ (the special case studied here is $n = 1 \mapsto N = d+1$); its kernel is the \id
$$
\preskip.4em \postskip.4em 
\fa = \gen {X_\alpha X_\beta - X_{\alpha'} X_{\beta'},
\alpha + \beta = \alpha' + \beta'}. 
$$
By co-morphism, $\varphi$ induces a morphism between affine spaces 

\snic{\psi :
\AA^{n+1}(\gA) \to \AA^N(\gA), \;u = (u_0, \ldots, u_n) \mapsto
(u^\alpha)_{|\alpha| = d}.}

If $\gA$ is a \cdiz, the image of $\psi$ is the zero set 
$\cZ(\fa)$ and we can show that $\psi$ induces an injection
$\PP^{n}(\gA) \to \PP^{N-1}(\gA)$: it is the Veronese embedding of degree~$d$. 
\\
Even more \gnltz, let $E$ be a direct summand 
in $\Ae {n+1}$, of rank~$1$.,
The \hmg component
 of degree $d$ of the \smq \algz~$\gS_\gA(E)$, which we denote by $\gS_\gA(E)_d$, is identified with a submodule
of $\gS_\gA(\Ae {n+1})_d \simeq \gA[U_0, \ldots, U_n]_{d}$ (\hmg component of degree $d$), a direct summand of rank $1$. \\
If we let $\psi(E) = \gS_\gA(E)_d$, we thus \gui{extend} the map $\psi$ previously defined.
\eoe

\vspace{-1em}
\begin{problem}\label{exoVeroneseMatrix}
 {(Veronese matrices)}\\
{\rm
Let two \pol \ris $\kuX = \kXn$ and $\kuY = \kYm$.  To every matrix $A
\in \gk^{m\times n}$, which represents a \ali $\gk^n \to \gk^m$, we can associate (watch the reversal), a $\gk$-morphism $\varphi_A : \kuY \to \kuX$ constructed as follows: let $X'_1$, \ldots, $X'_m$
be the $m$ \lin forms of~$\kuX$ defined as follows.

\snic {
\hbox{If}\quad \cmatrix {X'_1\cr \vdots \cr X'_m} = A \cmatrix {X_1\cr \vdots \cr X_n},
\quad \hbox {then }\quad
\varphi_A : f(\Ym) \mapsto f(X'_1,\ldots,X'_m).
} 

It is clear that $\varphi_A$ induces a \kli $A_d : \kuY_d \to \kuX_d$ between the \hmg components of degree $d \ge 0$, and that the restriction \hbox{$A_1 : \kuY_1 \to \kuX_1$} has as its matrix in the bases $(\Ym)$ and $(\Xn)$, the \emph{transpose} \hbox{of $A$}.  \hbox{The \kmoz} $\kuX_d$ is free of rank $n' = {n-1+d \choose d} $; it possesses a natural bases, that of the \moms of degree $d$, which we can choose to order lexicographically with $X_1 > \cdots > X_n$.  Similarly for $\kuY_d$ with its basis of $m' = {m-1+d \choose m-1}$ \momsz.  Let $V_d(A) \in \gk^{m' \times n'}$ be the \emph{transpose} of the matrix of the \endo $A_d$ in these bases (such \hbox{that $V_1(A) = A$}) and we say that $V_d(A)$ is the Veronese extension of $A$ in degree $d$. 
\\
 For example, let $n = 2$, $d = 2$, so $n' = 3$; if \smash{$A = \cmatrix {a &b\cr c & d\cr}$}, we obtain the matrix $V_2(A) \in \MM_3(\gk)$ as follows

\snic {
\cmatrix {x'\cr y'\cr} = A \cmatrix {x\cr y\cr} =
\cmatrix {ax + by\cr cx + dy\cr}, \;
\cmatrix {x'^2\cr x'y'\cr y'^2\cr} = 
\cmatrix {a^2 & 2ab & b^2\cr ac & ad+bc & bd\cr c^2 & 2cd & d^2\cr}
\cmatrix {x^2\cr xy\cr y^2\cr} 
.}

\emph {1.} 
If $A, B$ are two matrices for which the product $AB$ has a meaning, 
check the \egts $\varphi_{AB} = \varphi_B \circ \varphi_A$ and  $V_d(AB) =
V_d(A)V_d(B)$ for every $d \ge 0$.  
Also check that $V_d(\tra {A}) = \tra{V_d(A)}$.

\emph {2.}
If $E$ is a \kmoz, the \emph{$d$-Veronese transform of $E$}
 is the \kmo $\gS_\gk(E)_d$, \hmg component of degree $d$ of the \smq \alg $\gS_\gk(E)$.\\
If $E$ is a direct summand
 in $\gk^{n}$, then $\gS_\gk(E)_d$ is identified with a direct summand in 
 $\gS_\gk(\gk^n)_d \simeq\gk[X_1, \ldots, X_n]_d$  
(see also \Pbmz~\ref{exoVeroneseMorphism}).  Show that the image under $V_d$ of a \prr is a \prr and that we have a commutative diagram
$$
\preskip.4em \postskip.4em \ndsp
\xymatrix @C = 1.5cm{
\GA_n(\gk)\ar[d]_{\Im}\ar[r]^{V_d}         &\GA_{n'}(\gk)\ar[d]^{\Im}\\
\GG_n(\gk)\ar[r]^{d\rm -Veronese} &\GG_{n'}(\gk) \\
}
\qquad \hbox {with} \quad n' = {n-1+d\choose d} = {n-1+d \choose n-1} 
$$

\emph {3.}
Show that if $A$ is a \prr of rank $1$, the same holds for $V_d(A)$.  More \gnltz, if $A$ is a \prr of rank $r$, then $V_d(A)$ is a \prr of rank~${d+1-r \choose r-1}$.

}

\end{problem}

\vspace{-1em}
\begin{problem}\label{exoFossumKumarNori}
{(Some examples of finite projective resolutions)}\\
{\rm  
Given $2n+1$ \elts $z$, $x_1$, \dots, $x_n$,  $y_1$, \dots, $y_n$, of a \ri $\gA$, we define a sequence of matrices $F_k \in \MM_{2^k}(\gA)$, for $k \in \lrb{0..n}$, as follows

\snic {
F_0 = \cmatrix {z}, \qquad F_k = 
\cmatrix {F_{k-1} & x_k\I_{2^{k-1}}\cr y_k\I_{2^{k-1}} &\I_{2^{k-1}} - F_{k-1}\cr}
.}

Thus with $\ov z = 1-z$,
$$\preskip-.4em \postskip.2em 
F_1 = \cmatrix {z & x_1\cr y_1 & \ov z}, \qquad
F_2 = \cmatrix {
z & x_1 & x_2 & 0 \cr
y_1 & \ov z &0  &x_2 \cr 
y_2 &0 &\ov z & -x_1 \cr 
0 & y_2 & -y_1 & z \cr}
. 
$$

\emph {1.}
Check that $F_k^2 - F_k$ is the scalar matrix with the term $z(z-1) +
\sum_{i=1}^k x_iy_i$. Also show that $\tra{F_n}$ is similar to
$\I_{2^n} - F_n$ for $n \ge 1$.
Consequently, if $z(z-1) + \sum_{i=1}^n x_iy_i = 0$,
then $F_n$ is a \prr of rank $2^{n-1}$.

\smallskip \emph {2.}
We define three sequences of matrices 
$$
\preskip.3em \postskip.3em 
U_k, V_k \in \MM_{2^{k-1}}(\gA) \,
(k \in \lrb{1..n}),\quad G_k \in \MM_{2^k}(\gA)  \,(k \in \lrb{0..n}), 
$$
 as follows: $U_1 = \cmatrix {x_1}$, $V_1 = \cmatrix {y_1}$, $G_0 =
\cmatrix {z}$ and
$$\preskip.4em \postskip.2em 
U_k = \cmatrix {U_{k-1} & x_k\I \cr y_k\I & -V_{k-1}}, \;
V_k = \cmatrix {V_{k-1} & x_k\I \cr y_k\I & -U_{k-1}}, \;
G_k = \cmatrix {z\I & U_k \cr V_k & \ov z\I}
. 
$$
Thus,
$$\preskip-.6em \postskip-.4em 
U_2 = \cmatrix {x_1 & x_2\cr y_2 & -y_1}, \;
V_2 = \cmatrix {y_1 & x_2\cr y_2 & -x_1}, \;
G_2 = \cmatrix {
z & 0 & x_1 & x_2 \cr
0 & z & y_2  &-y_1 \cr 
y_1 &x_2 &\ov z & 0 \cr 
y_2 & -x_1 & 0 & \ov z \cr}
. 
$$
\begin{itemize}
\item [\emph {a.}]
Verify that $G_n$ and $F_n$ are conjugated by a permutation matrix.

\item [\emph {b.}]
Verify that $U_kV_k$ is the scalar $\sum_{i=1}^k x_iy_i$ and that $U_kV_k =
V_kU_k$.

\item [\emph {c.}]
For $n \ge 1$, if $z(z-1) + \sum_{i=1}^n x_iy_i = 0$, show that $G_n$ (therefore $F_n)$ is a \prr of rank $2^{n-1}$.
\end{itemize}

\emph {3.}
Let $M$ be an \Amoz. A \emph {finite projective resolution} of $M$ is an exact sequence of \mptfs $0 \rightarrow P_n \to \cdots \to P_1 \to P_0 \twoheadrightarrow M \to 0$; we say that $n$ is \emph{the length of the resolution}. In this case, $M$ is \pfz. 
\begin{itemize}
\item [\emph {a.}]
Consider two finite projective resolutions of $M$ that we can assume to be of the same length,
$$
\preskip.3em \postskip.3em 
\arraycolsep2pt\begin{array}{ccccccccccccccc}
0 &\rightarrow& P_n& \to & P_{n-1}& \to & \cdots & \to & P_1 & \to & P_0
& \to & M& \rightarrow& 0,
\\[.2em]
0 &\rightarrow& P'_n& \to & P'_{n-1}& \to & \cdots & \to & P'_1 & \to & P'_0
& \to & M&  \rightarrow& 0.
\end {array}
 $$
By using Exercise~\ref{exoSchanuelVariation}, show that we have in $\KO(\gA)$ the following \egt
$$\preskip.4em \postskip.4em \ndsp
\sum_{i=0}^n (-1)^i [P_i] = \sum_{i=0}^n (-1)^i [P'_i].\leqno (\star)
$$

\textit{Note.} Exercise~\ref{exoSchanuelVariation} provides a much more
precise result.
\eoe

 \textit{Definition and notation.} For a module $M$ which admits a finite projective resolution we let $[M] \in \KO(\gA)$ be the common value of $(\star)$ (even if $M$ is not \ptfz). We then define \emph{the rank of $M$} as that of $[M]$ and we have $\rg M = \sum_{i=0}^n (-1)^i \rg P_i\in\HO(\gA)$.\index{rank!of a module which admits finite projective resolution}

\item [\emph {b.}]
Let $M$ be an \Amo admitting a finite projective resolution; suppose \hbox{that $aM = 0$} with $a \in \Reg(\gA)$. Show that
$\rg(M) = 0$ \cad \hbox{that $[M] \in \KTO(\gA)$}.
\end{itemize}

If $\gk$ is an arbitrary \riz, we define the \ri 

\snic{\gB_n =
\gk[z, \ux, \uy] = \aqo{\gk[Z, \Xn, \Yn]}{Z(Z-1) + \sum_{i=1}^n X_iY_i}}

Thus $\gB_0 \simeq \gk\times\gk$. Let $\fb_n$ be the \id $\gen {z, \xn}$.

\emph {4.}
Show that the localized \ris
 $\gB_n[1/z]$ and $\gB_n[1/(1-z)]$ are \elr localized \ris
(\cad obtained by inverting a single \eltz) of a \pol \ri over $\gk$ with $2n$ \idtrsz.
Show that $\aqo {\gB_n}{x_n} \simeq \gB_{n-1}[y_n]\simeq\gB_{n-1}[Y]$.

\emph {5.}
For $n=1$, define a projective resolution of the $\gB_1$-module $\gB_1\sur{\fb_1}$ of length~$2$ and verify that
$[\gB_1\sur{\fb_1}] \in \KTO(\gB_1)$.

\emph {6.}
For $n=2$, define a projective resolution of the $\gB_2$-module $\gB_2\sur{\fb_2}$ of length~$3$

\snic {
0 \rightarrow \Im F_2 \to \gB_2^4 \to \gB_2^3 \vvvers{[z,x_1,x_2]} 
\gB_2 \twoheadrightarrow \gB_2\sur{\fb_2} \to 0
,}

and verify that $[\gB_2\sur{\fb_2}] \in \KTO(\gB_2)$.

\emph {7.}
Explicate a permutation $\sigma \in \rS_{2^n}$ such that the $n+1$ first \coes of the first row of the matrix $F'_n = P_\sigma F_n P_\sigma^{-1}$ are $z$, $x_1$, \dots, $x_n$ ($P_\sigma$ is the permutation matrix $\sigma$).

\emph {8.}
For $n=3$, define a projective resolution of the $\gB_3$-module
$\gB_3\sur{\fb_3}$ of length~$4$

\snic {
0 \rightarrow \Im (\I_8 -F'_3) \to \gB_3^8 \to \gB_3^7 \to 
\gB_3^4 \vvvvers{[z,x_1,x_2,x_3]}  \gB_3 \twoheadrightarrow \gB_3\sur{\fb_3} \to 0
,}

and verify that $[\gB_3\sur{\fb_3}] \in \KTO(\gB_3)$.

\emph {9.}
And in general?

}

\end{problem}


\vspace{-.3em}

\sol


\exer{exoleli2}{ 
We roughly 
rewrite the second proof of the local freeness lemma.
Let $\varphi$ be the \ali which has as its matrix $F$.
Let $f_j$ be the column $j$ of the matrix $F$, and $(e_1,\ldots ,e_n)$ be the canonical basis of $\Ae n$.
\pagebreak

By hypothesis, $(f_1,\ldots ,f_k,e_{k+1},\ldots ,e_n)$ is a basis of $\Ae n$.
The corresponding change of coordinate matrix is $B_1=\bloc{V}{0}{C'}{\I_{h}}.$
Since $\varphi(f_i)=\varphi\big(\varphi(e_i)\big)=\varphi(e_i)=f_i$, with respect to this basis,
$\varphi$ has a matrix of the type \halfsmashtop{$\bloc{\I_k}{X}{0}{Y}$}.
The computation gives
$$
B_1^{-1}=\bloc{V^{-1}}{0}{C}{\I_{h}},\qquad
G=B_1^{-1}\,F\,B_1=\bloc{\I_k}{L} {0}{-C'V^{-1}L'+W},
$$
where $L=V^{-1}L'$, and $C=-C'V^{-1}$.\\
Since $\cD_{k+1}(G)=0$, we have
$G=\bloc{\I_k}{L} {0}{0}$, therefore $W=C'V^{-1}L'$. 
\\
Let $B_2=\bloc{\I_k}{-L}{0}{\I_{h}}$, we have $B_2^{-1}=\bloc{\I_k}{L}{0}{\I_{h}}$, then $B_2^{-1}\,G\,B_2=\I_{k,n}$.\\
Finally, we obtain $B^{-1}\,F\,B=\I_{k,n}$
with
$$\preskip.4em \postskip.0em
B=B_1\,B_2=\bloc{V}{0}{C'}{\I_{h}}\cdot\bloc{\I_k}{-L}{0}{\I_{h}}=
\bloc{V}{-L'}{C'}{\I_h-W}
$$
and
$$\preskip.0em \postskip.4em
B^{-1}=B_2^{-1}\,B_1^{-1}=\bloc{\I_k}{L}{0}{\I_{h}}\cdot\bloc{V^{-
1}}{0}{C}{\I_{h}}=
\bloc{V^{-1}+LC}{L}{C}{\I_{h}}.
$$
The \egt $F^2=F$ gives in particular $V=V^2+L'C'$.\\
Therefore
$\I_k=V\,(\I_k+L'C'V^{-1})=V\,(\I_k-LC)$, and finally $V^{-1}=\I_k-LC$. Therefore as stated $B^{-1}=\bloc{\I_k}{L}{C}{\I_h}$.\\
Before proving the  statement regarding $\I_h-W$, let us prove the converse. \\
The double \egt
$$\preskip.3em \postskip.4em
\bloc{\I_k}{L}{C}{\I_h}
= \bloc{\I_k-LC}{L}{0}{\I_h} \,\bloc{\I_k}{0}{C}{\I_h}=
 \bloc{\I_k}{L}{0}{\I_h}\,\bloc{\I_k}{0}{C}{\I_h-CL}
$$
shows that $\I_k-LC$ is \iv \ssi $\I_h-CL$ is \iv \ssiz$\bloc{\I_k}{L}{C}{\I_h}$ is \ivz. This also gives 
$$
\preskip-.4em \postskip.2em
\det\bloc{\I_k}{L}{C}{\I_h}=\det(\I_k-LC)=\det(\I_h-
CL)\,.
$$
The computation then gives
$$\preskip.0em \postskip.1em
{\bloc{\I_k}{L}{C}{\I_h}}^{-1}=\bloc{V}{-VL}{-CV}{\I_h+CVL},
$$
hence
$$\preskip-.3em \postskip.4em
{\bloc{\I_k}{L}{C}{\I_h}}^{-1} \cdot \bloc{\I_k}{0}{0}{0} \cdot
\bloc{\I_k}{L}{C}{\I_h}
=  \bloc{V}{VL}{-CV}{-CVL}\,,
$$
which establishes the converse.
Finally, the \egt $B^{-1}\,F\,B=\I_{k,n}$ implies $B^{-1}\,(\I_n-F)\,B=\I_n-\I_{k,n}$, which gives 
$$\preskip.1em \postskip.4em
\bloc{\I_k-V}{-L'}{-C'}{\I_h-W}
={\bloc{\I_k}{L}{C}{\I_h}}^{-1} \cdot \bloc{0}{0}{0}{\I_h}\cdot
\bloc{\I_k}{L}{C}{\I_h}
$$
and we find ourselves in the \smq situation, therefore $(\I_h-W)^{-1}=\I_h-CL$ \hbox{and $\det\,V=\det(\I_h-W)$}.
}


\exer{exoJacobienneP2=P}{
Note $g$ (resp.\,$d$) as the left-multiplication (resp.\,right-multiplication) by $P$. We then have $g^2=g,\,d^2=d,\,gd=dg,\,\varphi=g+d-1$ and $\pi=g+d-2gd$.
}


\exer{exoP1FracRat} 
\emph {1a.}
The \gui{\hmg Sylvester matrix} $S$ is defined as that of the \ali  $(A,B)\mapsto PA+QB$ over the \hbox{bases $(u^{q-1},\ldots,v^{q-1})$} for $A$ (\pog of degree $q-1$), 
$(u^{p-1},\ldots,v^{p-1})$ for $B$ (\pog of degree $p-1$) and $(u^{p+q-1},\ldots,v^{p+q-1})$ for $PA+QB$ (\pog of degree $p+q-1$).\\
By making $v=1$, we see that $\tra S=\Syl(g,p,h,q)$, hence $\det(S)=\Res(g,p,h,q)$. \\
By making $u=1$, we see that $\tra S$ is almost the matrix $\Syl(\wi g,p,\wi h,q)$: the order of the rows, the order of the $q$ first columns and the order of the last $p$ must be reversed.
Hence the stated result because $(-1)^{\lfloor q/2\rfloor + \lfloor p/2\rfloor + \lfloor (p+q)/2\rfloor}
= (-1)^{pq}$.

\emph {1b.}
The \egt $S \wi S=\Res(P,Q)\,\I_{p+q}$ means that, if $k+\ell=p+q-1$, $u^kv^{\ell}\Res(P,Q)$ is a \coli of the column vectors of the matrix $S$.
That therefore exactly gives the required inclusion, which is after all just the \hmg version of the usual inclusion.

\emph {2.}
We write $f$ in the \ird form $f = a/b$ with $a$, $b \in \gk[t]$,
and we homogenize~$a$ \hbox{and $b$} in degree $d$ (maximum of the degrees of $a$ and $b$)
to obtain two \pogs $A$, $B \in \gk[u,v]$ of degree $d$. \\
If $\gk$ is an arbitrary \riz, we ask that $\Res(A,B)$ be \ivz.
That is \ncr for the fraction to remain well-defined after every \edsz.
Let us then see that the morphism $f$ is first defined at the \umd vector level

\snic {
(\xi : \zeta) \mapsto \big(A(\xi,\zeta) : B(\xi,\zeta)\big)
.}

This makes sense because if $1 \in \gen {\xi , \zeta}$, then $1 \in \gen {A(\xi,\zeta),B(\xi,\zeta)}$ after item~\emph{1b.}

To get back up to level $\GA_{2,1}(\gk)$, we take two new \idtrsz~\hbox{$x$, $y$} by thinking about the matrix $\cmatrix {xu & yu\cr xv & yv\cr}$.
As $\gen {u,v}^{2d-1} \subseteq \gen {A,B}$, we can write
$$(xu + yv)^{2d-1} = E(x,y,u,v)A(u,v) + F(x,y,u,v)B(u,v)
$$ 
with
$E$ and $F$ \hmgs in $(x,y,u,v)$. \\
Actually, $E$ and $F$ are bi\hmgs in $\big((x,y), (u,v)\big)$,
of degree $2d-1$ in $(x,y)$, of degree $d-1$ in $(u,v)$. 
As $EA$ is bi\hmgz, of bidegree $(2d-1, 2d-1)$, there exists 
(see the justification below)  
some \pog $\alpha'$ in $4$ variables,
$\alpha' = \alpha'(\alpha, \beta, \gamma, \delta)$, 
such that:

\snic {
EA = \alpha'(xu, yu, xv, yv), \quad \deg(\alpha')=2d-1
.}

Likewise with $FA$, $EB$, $FB$ to produce $\beta'$, $\gamma'$, $\delta'$. We then consider the matrices

\snic {
\xymatrix {\cmatrix {xu & yu\cr xv & yv\cr} \ar@{~>} [r] &
\cmatrix {\alpha & \beta\cr \gamma & \delta\cr}}, 
\qquad
\xymatrix {\cmatrix {EA & FA\cr EB & FB\cr} \ar@{~>} [r] &
\cmatrix {\alpha' & \beta'\cr \gamma' & \delta'\cr}
}.}

         The lifting we are looking for is then
$\cmatrix {\alpha & \beta\cr \gamma & \delta\cr} \mapsto
\cmatrix {\alpha' & \beta'\cr \gamma' & \delta'\cr}$.
\\
Note: $\alpha'$, $\beta'$, $\gamma'$, $\delta'$ are \pogs in $(\alpha, \beta, \gamma, \delta)$, of degree $2d-1$,
such that
$$\preskip.4em \postskip.4em \ndsp
\left| \matrix {\alpha & \beta\cr \gamma & \delta\cr}\right|
\hbox { divides }
\left| \matrix {\alpha' & \beta'\cr \gamma' & \delta'\cr}\right|,
\qquad
\alpha + \delta - 1 \hbox { divides } \alpha' + \delta' - 1 
. 
$$

\emph{Justification of the existence of $\alpha'$}. \\
This rests on the following simple fact: $u^iv^jx^ky^\ell$ is a monomial in $(xu, yu,
xv, yv)$ \ssi $i+j = k+\ell$; indeed, if this \egt is satisfied, there is a matrix $\cmatrix {m & n\cr r & s\cr} \in \MM_2(\NN)$ such
that the sums of rows are $(i,j)$ and the sums of columns are 
$(k,l)$.
A schema to help with the reading:

\snic {
\bordercmatrix[\lbrack\rbrack] {
   & k  & \ell \cr
i  & m  & n \cr
j  & r  & s \cr}
\qquad
\cmatrix {xu & yu\cr xv & yv\cr}
,}

and then

\snic {
u^iv^jx^ky^\ell = u^{m+n}v^{r+s}x^{m+r}y^{n+s}
= (xu)^m (yu)^n (xv)^r (yv)^s.
}

We deduce that a bi\hmg \pol in $\big((x,y), (u,v)\big)$, of bidegree $(d,d)$, is the \evn at $(xu, yu, xv, yv)$ of a \pog of degree $d$.

\emph {3.}
For $f(t) = t^2$, we obtain the lift
$$\preskip.4em \postskip.4em 
\cmatrix {\alpha & \beta\cr \gamma & \delta\cr} \mapsto
\cmatrix {\alpha^2(\alpha + 3\delta) & \beta^2(3\alpha + \delta)\cr 
\gamma^2(\alpha + 3\delta) & \delta^2(3\alpha + \delta)\cr}
. 
$$
More \gnltz, we develop $(\alpha + \delta)^{2d-1}$
in the form $\alpha^d S_d(\alpha,\delta) + \delta^d S_d(\delta,\alpha)$,
and we obtain the lift

\snic {
\cmatrix {\alpha & \beta\cr \gamma & \delta\cr} \mapsto
\cmatrix {\alpha^d S_d(\alpha,\delta) & \beta^d S_d(\delta, \alpha)\cr 
\gamma^d S_d(\alpha,\delta) & \delta^d S_d(\delta, \alpha)\cr}
.}

If $H = \cmatrix {a & b\cr c & d\cr}$, we obtain the following lift of $f(t) = {at + b \over ct +d}$:
$$\preskip-.5em \postskip.4em
\cmatrix {\alpha & \beta\cr \gamma & \delta\cr} \mapsto H \cmatrix {\alpha & \beta\cr \gamma & \delta\cr} H^{-1}.$$


\exer{exoConiqueFondamentale} \emph{(The fundamental conic or Veronese embedding $\PP^1 \to \PP^2$)}\\
We proceed as in Exercise~\ref{exoP1FracRat}, but it is simpler because, since $\gen {u,v}^2 = \gen {u^2, uv, v^2}$, the map $(u : v) \mapsto (u^2 : uv : v^2)$ is well-defined at the \vmd level.

We introduce $(x, y)$ with the matrix $\cmatrix {\alpha & \beta\cr
\gamma & \delta\cr} \leftrightarrow \cmatrix {xu & yu\cr xv & yv\cr}$ in mind.
We develop $(xu + yv)^2 = x^2 u^2 + 2xy uv + y^2 v^2$, sum of $3$ terms which will be the $3$ diagonal terms of a  matrix of $\GA_{3,1}(\gk)$, then we complete such that each column the ad-hoc multiple of the vector $\tra[\,u^2\;uv\;v^2\,]$. Which gives

\snic {
\cmatrix {
x^2 u^2 & 2xy u^2 & y^2 u^2 \cr
x^2 uv  & 2xy uv  & y^2 uv \cr
x^2 v^2 & 2xy v^2 & y^2 v^2 \cr
}
\qquad
F = \cmatrix {
\alpha^2     & 2\alpha\beta  & \beta^2 \cr
\alpha\gamma & 2\alpha\delta  & \beta\delta \cr
\gamma^2     & 2\gamma\delta & \delta^2 \cr
}.}

The lift $\GA_{2,1}(\gk) \to \GA_{3,1}(\gk)$ is 
$\cmatrix {\alpha & \beta\cr \gamma & \delta\cr} \mapsto F$.
\\
We of course have $\Tr(F) = (\alpha + \delta)^2 = 1$,
$\cD_2(F) \subseteq \gen {\alpha\delta - \beta\gamma} = 0$,
and $F$ is a \prr of rank~$1$.

\exer{exoProjecteurCorangUn} \emph{(\Mprns of corank $1$)}\\
\emph{1.}
We provide two solutions for this question. The first consists of using the expression of the adjoint in terms of the starting matrix; the second proof uses the \lonz.

For $A \in \Mn(\gA)$ we have the classical expression of $\wi{A}$ as a \pol in~$A$

\snic{\wi A = (-1)^{n-1}Q(A) \;$  with $\;XQ(X)=\rC{A}(X)-\rC{A}(0).}

Apply this to a \prr $P$ of rank $n-1$.
We get
$$\rC{P}(X) =  (X-1)^{n-1}X, \,Q(X)=(X-1)^{n-1}\hbox{ and }
 (P - \In)^{n-1}  = (-1)^{n-1} \wi P.$$
Since $(\In-P)^{n-1}=\In-P$, we obtain
$ P + \wi {P} = \In.$

Here is the proof by \lonz.
By the local structure \tho for \mptfs (\thref{prop Fitt ptf 2} or \thref{th ptf loc free}), there exist \come \lons such that over each localized \riz,
 $P$ is similar to $\I_{r,n} $, where the integer $r$ a priori depends on the \lonz. Here, since $P$ is of rank  $n-1$, we have $r=n-1$ or $1=0$. Therefore $P + \wi {P} = \In$ over each localized \riz, and the \egt is also globally true by the basic \plgz.\iplg

\emph{2.}
Let us see the \dem by \come \lonsz.  Over the localized \ri
 $\gA_s$, the \prr $P$ is similar to $Q_s=\I_{r,n}$, where $r$ depends on $s$.
We have $Q_s + \wi {Q_s} = \In$. 
\\
If $r < n-1$, then
$Q_s+\wi {Q_s}= \I_{r,n}$.  If $r = n$, then $Q_s+\wi {Q_s}  = 2\,\In$. 
\\
Recap: if $r\neq n-1$, then $1=0$ and the rank is also equal to $n-1$.
Consequently over all the localized \ris
 $\gA_s$, the \prr $P$ is of \hbox{rank $n-1$}, and therefore also globally.

\emph{3.}
It suffices to multiply $P + \wi {P} = \In$ by $P$ to obtain $P^2 = P$.

\exer {exoMatriceCorangUn}{
There exists a $B \in \Mn(\gA)$ such that $ABA = A$, in order for $AB$ to be a \prr with the same image as~$A$, so of rank $n-1$, and for $BA$ to be a \prr with the same kernel as $A$, therefore \egmt of rank $n-1$. We define $P$ and $Q \in \Mn(\gA)$ by $AB = \In - P$, and $BA = \In - Q.$
\\
Thus $P$, $Q\in\GA_{1,n}(\gA)$, with
$A = (\In-P)A = A(\In-Q)$.

\emph{1.}
We have $\det A = 0$, \cad $\wi {A} A = A \wi {A} = 0$, so $\Im A
\subseteq \Ker \wi {A}$.  \\
Next
$\wi{AB} = \wi{\In - P}
= P$ (because $P\in\GA_{1,n}(\gA)$), and the \egt  $\wi {B} \wi {A} = P$ proves that

\snic{\Ker \wi {A} \subseteq \Ker P = \Im(\In - P) = \Im A.}

Conclusion: $\Ker \wi {A}  = \Im A = \Im(\In - P)$.

\emph{2.}
By reasoning as in item~\emph{1}, we obtain $\Im \wi {A} \subseteq \Ker A = \Ker (BA)   = \Im Q$, then $\wi {A} \wi {B} = \wi {BA}   = \wi {\In - Q}
= Q$, and $\Ker A = \Im \wi {A} = \Im Q$.

\emph{3.}
We apply item~\emph{1}  to $\tra A$, so $\Im \tra A = \Ker \tra
{\wi{A}}$.  Then, we explicate the \gui{left} \prr (of rank~$1$)
associated with $\tra A$. We have
$$\preskip.4em \postskip.4em 
\tra A \tra B \tra A = \tra A,
\quad \hbox {which we write} \quad
\tra (BA) \tra A = \tra A
\quad \hbox {with} \quad
\tra (BA) = \In - \tra Q . 
$$
This left-\prr is therefore $\tra Q$.

\emph{4.}
Similarly, item~\emph{2}  gives $\Im \tra {\wi {A}} = \Ker \tra A$. We explicate the \gui{right-\prrz} (of rank~$1$) associated with $\tra A$, we obtain $\tra P$, whence the stated result.

\emph{5.}
Finally

\snic{
\Ae n  / \Im A = \Ae n  / \Im(\In - P) \simeq \Im P, \qquad
\Ker \tra A = \Im \tra P,}

so the two modules (\pros of rank $1$) are indeed duals of one another.
\\
Remark: we can \egmt use

\snic{\Ae n  / \Im \tra A = \Ae n  / \Im (\In - \tra Q) \simeq
\Im \tra Q, \qquad
\Ker A = \Im Q
,}

to see that the two modules (\pros of rank $1$) $\Ae n  / \Im \tra A$
and $\Ker A$ are indeed duals.
}


\exer{exoIntersectionSchemasAffines} \emph{(Intersection of two affine schemes over $\gk$)}\\
First of all we notice that surjective arrows $\kuX\vers{\pi_1}\gA$ and  $\kuX\vers{\pi_2}\gB$ in the category of \pf \klgs are seen, from the point of view of the schemes, as \gui{inclusions} $A\vers{\iota_1}\gk^n$ and $B\vers{\iota_2}\gk^n$, \hbox{where $\gk^n$} is interpreted as the affine scheme corresponding to $\kuX$. The \dfn of the intersection by tensor product is therefore in fact a \dfn as the push out of the two arrows $\pi_1$ and $\pi_2$ in the category \hbox{of \pf \klgsz}, or as the pull back of the two arrows $\iota_1$ and $\iota_2$ in the category of affine schemes over $\gk$.    

The center of the ellipse, the center of the circle and the double point of intersection have for respective \coos $(0,0)$, $(c,0)$ and $(a,0)$.
The computation of other points of intersection gives

\snic{x = a(2ac + 1 - a^2)/(a^2 - 1)$ and $y^2 = 4ac(a^2 - ac - 1) / (a^2 - 1)^2.}

From the point of view of quotient \algs we obtain

\snic {
\gA = \aqo{\gk[X,Y]}{f},\quad \gB = \aqo{\gk[X,Y]}{g}, \quad
\gC = \aqo{\gk[X,Y]}{f,g}.
}

Which gives the morphisms

\snic {
\xymatrix @R=10pt @C=25pt{
              & \gK=\gk[X,Y]\ar[dl]_(.6){\pi_1} \ar[dd]\ar[dr]^(.6){\pi_2} \\
\gA\ar[dr]    &              &  \gB\ar[dl] \\
              & \gC = \gA\otimes_\gK \gB \\
}
\quad
\xymatrix @R=10pt @C=8pt{
                   & \gk^2   \\
\{f = 0\}\ar[ur]^{\iota_1}  &              & \{g = 0\}\ar[ul]_{\iota_2} \\
            &\{f = 0\} \cap \{ g = 0\}\ar[ul]\ar[uu]\ar[ur] \\
}
}

If $\gk$ is a \cdi and if 
$4ac(a^2 - ac - 1) (a^2 - 1)\in \gk\eti$,
the \klgs $\gA$ and~$\gB$ are integral, but not $\gC$: we have an \iso

\snic{\gC\simarrow\gk[\zeta]\times \gk[\vep],\;\hbox{ where }\,\vep^2=0\;
\hbox{ and }\,\zeta^2=4ac(a^2 - ac - 1) /(a^2 - 1)^2.}

The \alg $\gC$ is a \kev of dimension $4$, corresponding to the affine scheme formed by two points of multiplicity $1$ (defined over $\gk$ or over a quadratic extension of $\gk$) and a point of multiplicity $2$ (defined over $\gk$).

\exer{exoPolPseudoUnitaire} \emph{(Pseudo\mon \polsz)}\\
Recall that for an \idm $e$, we have $\gen {a,e} = \gen {(1-e)a + e}$; if $e'$ is another \idm \ort to $e$, we have $\gen {\ov a} = \gen {\ov {e'}}$ in $\aqo{\gA}{e}$ \ssi $\gen {(1-e)a} = \gen {e'}$ in $\gA$.

\emph {1.}
For $k > r$, we have $a_k = 0$ in each component, so in $\gA$.  The \elt $a_r$ is null in $\aqo{\gA}{e_r}$, \iv in $\aqo{\gA}{1-e_r}$ therefore $\gen{a_r} = \gen {e_r}$.
\\
 Similarly in $\aqo{\gA}{e_r}$, we have $\gen {\ov
{a_{r-1}}} = \gen {\ov {e_{r-1}}}$ thus $\gen {(1-e_r)a_{r-1}} = \gen
{e_{r-1}}$, and so on.

\emph {2.}
Localize at each of the $e_i$'s.

\exer{exoPolLocUnitaire} \emph{(\Lot \mon \polsz)}\\
\emph {1.}
As $f(t) = 0$, we have $f \in \fa$, hence a surjective \Ali $\aqo{\gA[T]}{f} \twoheadrightarrow \gA[T]\sur\fa$ between two free \Amos of the same rank $n$: it is an \iso
(Proposition~\ref{propDimMod1}), so $\fa = \gen {f}$.

\emph {2.}
The \polcar $f$ of $t$ is \mon of degree $n$ because $\gA[t]$ is of constant rank  $n$.
As $f(t) = 0$, we have $f \in \fa$, hence a surjective \Ali $\aqo{\gA[T]}{f} \twoheadrightarrow \gA[T]\sur\fa$, of a free \Amo of rank $n$ over an $\gA$-\mrc $n$;
it is therefore an \iso (Proposition~\ref{prop epi rank constant}),
\hbox{so $\fa = \gen {f}$}.

\emph {3.}
Let $f = \sum_{i=0}^r a_i T^i=\sum_{i=0}^r f_r$ be a \lot\mon \pol of formal degree $r$, with the \sfio
$(e_0,\ldots,e_r)$, \hbox{and $fe_d=f_d$} \mon of degree $d$ modulo $\gen{1-e_d}$ for each $d\in\lrb{0..r}$.
\\
Then $a_r=e_r$ is \idmz. Then $f-f_r=(1-e_r)f$ is \lot \mon of formal degree $r-1$ and we can end by descending \recu on $r$ to compute the $e_d$'s from $f$.
If the \ri is discrete we obtain a test to decide if a given \pol is \lot \monz: each of the successively computed $e_d$'s must be \idm and the sum of the $e_d$'s must be equal to~$1$.


\exer{exoOneRankPtfIdeal} \emph{(Invertible modules and \mrcs $1$)}\\
\emph {1.} There exists an \Asub $N$ of $\gB$ such that $M.N=\gA$.\\
We have $(\xn)$ in $M$ and $(\yn)$ in $N$ such that $1=\sum_ix_iy_i$ and
$x_iy_j\in\gA$. We verify that $M = \sum_i \gA x_i$ and $N =
\sum_i \gA y_i$.  Let $\sum_kz_k\otimes z'_k$ in $M\otimes_\gA
M'$. We have, by noticing that $y_iz_k\in N.M = \gA$

\snic {\arraycolsep2pt
\begin{array}{rclcl} 
 \sum_kz_k\otimes z'_k &  = & \sum_{k,i}x_iy_iz_k\otimes z'_k  & =  & \sum_{k,i} x_i\left(y_iz_k\right)\otimes z'_k    \\[1mm] 
  & =  &  \sum_{k,i} x_i\otimes \left (y_iz_k\right)z'_k  &
=  & \sum_ix_i\otimes \left(y_i \sum_k  z_k z'_k\right),    
 \end{array}
}

therefore the canonical surjection $M\otimes_\gA M'\to M.M'$ is injective.

\emph {2.}
It must be shown that $\fa$ contains a \ndz \elt (Lemma~\ref{lemIdproj}~\emph{\iref{i6lemIdproj}}), which is \imdz.


\exer{exoPicAPicFracA} \emph{(The exact sequence with $\Pic\gA$ and $\Pic\gK$, where $\gK = \Frac\gA$)}\\
Defining the sequence is obvious; thus, the map $\gK\eti \to \Gfr(\gA)$ is that which to $x \in \gK\eti$ associates the principal fractional \id $\gA x$. No issues either to verify that the composition of two consecutive morphisms is trivial. 
\\
\emph{Exactness in $\gK\eti$:} if $x\in
\gK\eti$ is such that $\gA x = \gA$, then $x\in\Ati$.  
\\
\emph{Exactness in $\Gfr(\gA)$:} if $\fa \in \Gfr(\gA)$ is free, it means that it is principal \cad of the form $\gA x$ with $x\in\gK\eti$.

Only the exactness in $\Pic\gA$ is more delicate. \Gnltz, if $P$ is {a \ptf \Amoz}, then the canonical map $P \to \gK\te_\gA P$ is injective because $P$ is contained in a free \Amoz. Thus let $P$ be an $\gA$-\mrc $1$ such that $\gK\te_\gA P \simeq \gK$. Then $P$ is injected into $\gK$, then into $\gA$ (multiply by a denominator), \cad $P$ is \isoc to an integral \idz~$\fa$ of $\gA$. Similarly, the dual $P\sta$ is \isoc to an integral \id $\fb$ of $\gA$. \\
We have $\gA \simeq P\te_\gA P\sta \simeq \fa
\te_\gA \fb \simeq \fa\fb$, so $\fa\fb$ is generated by a \ndz \elt $x \in \gA$. We have $x \in \fa$ so $\fa$ is an \iv \idz: we have found an \iv \id $\fa$ of $\gA$ such that $\fa \simeq P$.


\exer{exoCoMatCoRang1}
\emph{1}  and \emph{2.}
Immediate.

\emph{3.}
Consider the short sequence $\Ae N \vers {A'}
\Ae m \vers {A}  \Ae n$; it is \lot exact, so it is globally exact.

\emph{4.}
Every \stl module of rank $1$ can be given in the form $\Ker A$ where $A \in \Ae {n \times (n+1)}$ is a surjective matrix $\Ae {n+1} \vvers {A}  \Ae n$. 
Since $1 \in \cD_n(A)$, we apply question \emph {3}  with $m = n+1$. We obtain $A' \in \Ae {(n+1) \times 1}$ of rank $1$ with $\Im A' = \Ker A$; so the column $A'$ is a basis of $\Ker A$.


\exer{exoVarProj} \emph{(\Hmg \pols and $\Pn(\gk)$)}\\
Let $f\in\gk[\Xzn]$ be a \pog of degree $m$ and (to simplify) 
$P=\gen{\ua,\ub,\uc}\subseteq\gk^{n+1}$
be a direct summand of rank $1$. 
Suppose that $f(\ua)=f(\ub)=f(\uc)=0$ and that we want to show that $f(\ux)=0$ if $\ux=\alpha\ua+\beta\ub+\gamma\uc$. 
The matrix of $(\ua, \ub, \uc)$ is of rank $1$, therefore the $a_i$'s, $b_j$'s, $c_k$'s are \comz. It therefore suffices to prove the \egt after \lon at one of these \coosz.
For example over $\gk[1/a_0]$ we have $\ux=(\alpha+\frac{b_0}{a_0}\beta+\frac{c_0}{a_0}\gamma)\ua=\lambda\ua$,
and \hbox{so $f(\ux)=\lambda^mf(\ua)=0$}.


\exer{exoGLnTangent} \emph{(Tangent space to $\GLn$)}\\
Consider the \klg $\gk[\vep]=\aqo{\kT}{T^2}$.\\ 
Let  $A\in\GL_n(\gk)$ and $H \in \Mn(\gk)$. We have 
$A + \vep H=A(\In+\vep A^{-1}H)$, and $\In+\vep M$ is invertible, with \inv
$\In-\vep M$, for every $M\in\Mn(\gk)$. Therefore 
 $A + \vep H \in \GLn(\gk)$ for any $H$.
Thus, the tangent space $\rT_A(\GL_n)$ is \isoc to $\Mn(\gk)$.
\\
NB: $(A + \vep H)^{-1}=A^{-1}-\vep A^{-1}H A^{-1}$.


\exer{exoSLnTangent} \emph{(Tangent space to $\SLn$)}\\
We use the \klg $\gk[\vep]$ of Exercise~\ref{exoGLnTangent}.  
 For $A$, $H \in \Mn(\gk)$, we have $\det(A + \vep H) =
\det(A) + \vep \Tr(\wi{A}H)$. We deduce

\snic {
\det(A + \vep H) = 1  \iff  (\det(A) = 1  \hbox { and } \Tr(\wi{A}H) = 0) 
.}

We therefore have, for $A \in \SLn(\gk)$, 
$\rT_A(\SLn) = \sotQ {H \in \Mn(\gk)} {\Tr(\wi{A}H) = 0}$.\\
Let us show that $\rT_A(\SLn)$ is a free $\gk$-module of rank $n^2 -1$. \\
Indeed, the $\gk$-\lin \auto $H \mapsto AH$ of $\Mn(\gk)$ transforms $\In$ into~$A$ and bijectively applies
$\rT_{\In}(\SLn)$ over $\rT_A(\SLn)$, since we can verify it by writing $\Tr(H) = \Tr(\wi{A}\, AH)$.  Finally, $\rT_{\In}(\SLn)$ is the \ksmo of $\Mn(\gk)$ made of the matrices of null trace (which is indeed free of rank $n^{2}-1$).
\\
NB: $H \mapsto HA$ was also possible, because $\Tr(AH\, \wi{A}) = \Tr(\wi{A}\, AH) = \Tr(H)$.


\exer{exoTangentJ0ConeNilpotent} \emph{(Tangent space at $J_0$ to the nilpotent cone)}\\
\emph{1.}
We easily see that $\varphi(H)J_0 = J_0\varphi(H)$. If $\gk$ was a field, we could deduce that $\varphi(H)$
is a \pol in $J_0$. The direct computation gives

\snic{\varphi(e_{ij}) = \cases {
0 & if $i < j$ \cr
J_0^{n-1-(i-j)} & otherwise.}}

In particular, $\,
\varphi(e_{i1}) = J_0^{n-i}$.
We therefore have $\Im \varphi = \bigoplus_{k=0}^{n-1} \gk J_0^{k}$.

%
\emph{2.}
For $k \in \lrb{0.. n-1}$, the matrix $J_0^k$ has null \coesz, except for those that are on the $k^{\rm th}$ up-diagonal,
 all equal to $1$.
We can therefore take as direct complement 
of $\Im\varphi$ the submodule generated by the $e_{ij}$'s, with $j < n$ (we therefore omit \hbox{the $e_{in}$'s} which corresponds to the last position of the up-diagonals of the $J_0^k$'s).
We then define $\psi$ by

\snic{\psi(e_{ij}) = \cases {
0 & if $j < n$ \cr
e_{i1} & if $j = n$ \cr}\quad \hbox {or} \quad
\psi(H) = H \tra {J_0^{n-1}}.}

We easily verify that $\psi(J_0^{n-i}) = e_{i1}$ for $i \in \lrbn$, then $(\varphi \circ \psi)(A) = A$ \hbox{if $A \in \Im\varphi$}, and finally $\varphi \circ \psi \circ \varphi = \varphi$. By miracle, we also have $\psi \circ \varphi \circ \psi = \psi$. \\
We have $e_{ij} - e_{i'j'} \in \Ker\varphi$ as soon as $i'-j' = i-j$ ($i' \ge j'$, $i \ge j$) and we obtain a basis of $\Ker\varphi$ by considering the $n(n-1) \over 2$ matrices $e_{ij}$ with $i < j$ and the $n(n-1) \over 2$
matrices $e_{i1} - e_{i+r,1+r}$, $r\in \lrb{1.. n-i}$, $i \in \lrb{1.. n-1}$.

%
\emph{3.}
We use the \klg $\gk[\varepsilon]\simeq\aqo\kT{T^2}$.  
For $A, H \in \Mn(\gk)$, we have 

\snic{(A + \varepsilon H)^n = A^n + \vep\som_{i+j=n-1} A^i H A^j.}

For $A = J_0$, we find that the tangent space \gui{to the nilpotent cone} is $\Ker \varphi$ which is a free module of rank $n^2 - n$ (it is the dimension of the nilpotent cone).


\prob{exoAnneauCercle1} \emph{(The \ri of the circle)}\\
\emph{1.}
Naively: let $f = f(x,y) \in \gk[x,y]$ be a conic, \cad a \pol of degree $2$, and $(x_0, y_0)$ be a $\gk$-point of $\so{f(x,y) = 0}$.

\addhabille1\Habillage{\includegraphics{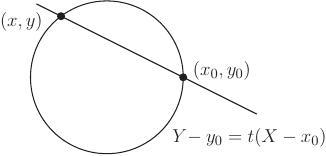}}{0}{0pt}
The classical trick of parameterization consists in defining $t$ by $y-y_0$ $=$ $t(x-x_0)$ and, in the \eqn
$$f(x,y) = f\big(x, t_0 + t(x-x_0)\big) = 0,$$
in looking for $x$ in terms of $t$. This \eqn admits $x = x_0$ as a solution, hence the other solution in the rational form.
\endHabillage

Algebraically speaking, we suppose that $f$ is \irdz. Let $\gk[x,y] = \gk[X,Y]/\gen {f}$. We obtain $\gk(x,y) = \gk(t)$ with $t = (y-y_0) / (x-x_0)$.
Here, the reader will compute the expressions of $x$, $y$ in terms of $t$: $x = {t^2 - 1 \over t^2 + 1}$, $y = {-2t \over t^2 + 1}$.
\\
\Gmqtz, the \elts of $\gk[x,y]$ are precisely the rational fractions defined everywhere on the projective line  $\PP^1(\gk)$ (parameterized by $t$) except maybe at the \gui{point} $t = \pm i$.

\emph{2.} 
We have $x = 1-2u$, $y = -2v$, so $\gk[x,y] = \gk[u,v]$. 
The \egt $\gk[x,y]=\gk[u,v]$ 
is not difficult and is left to the reader. What is more difficult, is to show \hbox{that $\gk[u,v]$} is the \cli of $\gk[u]$ in $\gk(t)$. We refer the reader to~Exercise~\ref{exoAnneauOuvertP1}.
\\
\Gmqtz, the poles of $x$ and $y$ are $t = \pm i$, which confirms that $x$, $y$ are integral over $\gk[(1+t^2)^{-1}] = \gk[u]$.
Algebraically, we have $x = 1 - u$, $y^2 = -1-x^2 \in \gk[u]$, and $x$, $y$ are indeed integral over $\gk[u]$.

\emph{3.}  
If $i^2 = -1$, we have $(x+iy)(x-iy) = 1$.\\
By letting $w = x+iy$, we have $\gk[x,y] = \gk[w, w^{-1}]$.

\emph{4.} 
We apply the standard method at a smooth point of a planar curve.
We write 
$$f(X,Y) - f(x_0,y_0) = (X-x_0)u(X,Y) + (Y-y_0)v(X,Y)
$$ 
with here $u = X+x_0$, $v = Y+y_0$; the matrix $A = \cmatrix {y-y_0 & x+x_0\cr
x_0-x & y+y_0\cr}$ is therefore a \mpn of $(x-x_0, y-y_0)$ with $1 \in \cD_1(A)$. Let us explicate the membership $1 \in \cD_1(A)$:

\snic {
(-y_0)(y-y_0) + x_0(x+x_0) + x_0(x_0-x) + y_0(y+y_0) = 2
.}

This leads to the matrix $B = \frac 1   2  \Cmatrix{2pt} {-y_0 & x_0\cr x_0 & y_0}$; this one satisfies $ABA = A$ and the desired matrix $P$ is $P = \I_2 - AB = \wi {AB}$
{
$$
AB = {1\over 2} \Cmatrix{2pt} {y-y_0 & x+x_0\cr x_0-x & y+y_0\cr}
\Cmatrix{2pt} {-y_0 & x_0\cr x_0 & y_0} = {1\over 2}
\Cmatrix{2pt} {x_0x - y_0y + 1 & y_0x + x_0y \cr y_0x + x_0y &
          -x_0x + y_0y + 1}.
$$
}

\vspace{-.5em} Hence the \gnl expression of $P$,
$
P = \frac 1   2 \cmatrix {-x_0x + y_0y + 1 & -(y_0x + x_0y) \cr
-(y_0x + x_0y) & x_0x - y_0y + 1}
$, for $x_0=1,\,y_0=0\,:
\frac 1   2 \cmatrix {1-x & -y\cr -y & 1+x}$.
Thus, $P$ is a \prr of rank $1$, \mpn of $(x-x_0, y-y_0)$.  
As $P$ is \smqz, \Egtz~(\iref{eqpmlm}) of Proposition~\ref{pmlm} has as consequence 
that $(x-x_0)^2 + (y-y_0)^2$ is a \gtr of $\gen {x-x_0, y-y_0}$ 
with $(x-x_0)^2 + (y-y_0)^2 = -2(x_0x + y_0y - 1)$.

\Gmqtz, $xx_0 + yy_0 - 1 = 0$ is the tangent line to the circle $x^2 + y^2 = 1$ at the \hbox{point $P_0 = (x_0, y_0)$}.
For those who know the divisors: the divisor of the zeros-poles of this tangent is the principal divisor $2P_0 - 2P_{t = \pm i}$, which corresponds to the fact that the square of the \id $\gen {x-x_0, y-y_0}$ is principal.

Variant I: we directly treat the case of the point $(x, y)=(1,0)$ (see the following question) then we use the fact that the circle is a group to pass from the point $(1,0)$ to an arbitrary point $P_0 = (x_0, y_0)$.  Thus, we dispose of the \gui{rotation} \auto 

\snic {
\cmatrix {x\cr y\cr} \mapsto
\crmatrix {x_0 & -y_0 \cr y_0 & x_0\cr} \cmatrix {x\cr y\cr}
\; \hbox {which realises} \;
\crmatrix {x_0 & -y_0 \cr y_0 & x_0\cr} \cmatrix {1\cr 0\cr} =
\cmatrix {x_0\cr y_0\cr}
.}

We consider its inverse $R$
$$\preskip.4em \postskip-.1em 
R = \crmatrix {x_0 & y_0 \cr -y_0 & x_0\cr}, \;
  R \cmatrix {x\cr y\cr} =\cmatrix {x'\cr y'\cr}, \;
 R \cmatrix {x_0\cr y_0\cr}= \cmatrix {1\cr 0\cr}
, 
$$
so that
$$\preskip.3em \postskip.4em 
  R \cmatrix {x-x_0\cr y-y_0\cr}=\cmatrix {x' - 1\cr y'\cr},
\quad \hbox{hence} \quad
\gen {x'-1, y'} = \gen {x-x_0, y-y_0}
. 
$$
As $\gen {x'-1, y'}^2 = \gen {x'-1}$, we obtain
$\gen {x-x_0, y-y_0}^2 = \gen {x_0x + y_0y - 1}$.

Variant II:
we provide another justification of the invertibility of $\gen{x-x_0, y-y_0}$ which does not directly use the fact that the circle is smooth. We consider $\gk[x,y]$ as an extension of degree $2$ of $\gk[x]$, by using $(1, y)$ as the basis.  We dispose of a $\gk[x]$-\auto $\sigma$ which transforms $y$ into $-y$. \\
We consider the norm $\rN$ of $\gk[x,y]$ over $\gk[x]$. For $z = a(x) + b(x)y$, we have

\snic {
\rN(z) = z\sigma(z) = (a + by)(a - by) = a^2 - (1-x^2)b^2 =
a^2 + (x^2-1)b^2
.}

The idea to invert $\gen {x-x_0, y-y_0}$ is to multiply it by its $\gk[x]$-conjugate. Let us show the following \egtz, certificate of the invertibility of the \id $\gen {x-x_0, y-y_0}$,

\snic {
\gen { x-x_0, y-y_0}\ \gen {x-x_0, y+y_0} = \gen {x-x_0}
.}

Indeed, the \gtrs of the left-product are

\snic {
(x-x_0)^2, \quad (x-x_0)(y+y_0), \quad
(x-x_0)(y-y_0), \quad y^2 - y_0^2 = x_0^2 - x^2
.}

Hence $\gen {x-x_0, y-y_0}\ \gen {x-x_0, y+y_0} =
(x-x_0) \gen { g_1, g_2, g_3, g_4}$ with

\snic {
g_1 = x-x_0, \quad g_2 = y+y_0, \quad g_3 = y-y_0, \quad g_4 = x+x_0
.}

But $\langle g_1, g_2, g_3, g_4\rangle$ contains ${g_4 - g_1\over 2} = x_0$
and ${g_2 - g_3\over 2} = y_0$ therefore it contains $1 = x_0^2 + y_0^2$.

\emph{5.} 
By brute force, by using only using $1 \in \gen {x-1, x+1}$ on the right-hand side,

\snuc {
\gen {x-1, y}\,\gen {x-1, y} = \gen {(x-1)^2, (x-1)y, y^2} =
(x-1)\,\gen {x-1, y, -(x+1)} = \gen {x-1}.
}

We divide this \egt by $x-1$: $\gen {x-1, y}\,\gen {1, {y\over x-1}}=\gen {1}$ and let

\snic {
x_1 = x-1, \quad x_2 = y, \; y_1 = 1, \; y_2 = {y \over x-1},
\; \hbox { such that} \; x_1y_1 + x_2y_2 = -2,
}

which leads to the projection matrix $P$ of rank $1$

\snic {
P = \frac{-1}2 \cmatrix {y_1\cr y_2\cr} [x_1, x_2] =
\frac{-1}2 \cmatrix {
x_1 y_1  &  x_2 y_1\cr
x_1 y_2  &  x_2 y_2\cr
} =
\frac{1}2 \cmatrix {
1 -x   &  -y\cr
-y  & 1+x\cr
}}

\emph{6.} 
Let $\rN = \rN_{\gk[x,y]/\gk}$.  For $a, b \in \gk[x]$, $\rN(a + by)
= a^2 + (x^2-1)b^2$.  The \egt to prove on the degrees is obvious if $a$ or $b$ is null. Otherwise, we write, with $n = \deg a$ \hbox{and $m = 1 + \deg b$}, $a^2 =
\alpha^2 x^{2n} + \dots$, $(x^2 - 1)b^2 = \beta^2 x^{2m} + \dots$ ($\alpha,
\beta \in \gk\sta$). The case \hbox{where $2n \ne 2m$} is easy. If $2n = 2m$, then
$\alpha^2 + \beta^2 \ne 0$ (because $-1$ is not a square in $\gk$), and so the \pol $a^2 + (x^2 - 1)b^2$ is of degree $2n = 2m$.
\\
If $a + by$ is \iv in $\gA$, $\rN(a+by) \in \gk[x]^{\!\times} =
\gk\sta$; hence $b = 0$ then $a$ is constant. 
Recap: $\gk[x,y]^{\!\times} = \gk\sta$. This is specific to the fact that $-1$ is not a square in $\gk$ because if $i^2 = -1$, the \egt $(x + iy)(x-iy) = 1$ shows the existence of \ivz s other than the constants.
\\
Let us show that $y$ is \irdz. \\
If $y = zz'$, then $\rN(y) = \rN(z)\rN(z')$, \cad $x^2 - 1 = (x-1)(x+1) = \rN(z)\rN(z')$. But in $\gk[x]$, $x \pm 1$ are not associated with a norm (a nonzero norm is of even degree). Therefore $\rN(z)$ or $\rN(z')$ is a constant, \cad $z$ or $z'$ is \ivz.
Similarly, $1 \pm x$ are \irdsz.
\\
We will use the \egt

\snic {
y^2 = (1-x)(1+x), \; \hbox { analogous to } \,
2\cdot 3 = (1 + \sqrt {-5}) (1 - \sqrt {-5}) \hbox { in }
\ZZ[\sqrt {-5}],
}

to see that $\gen {x-1, y}$ is not a \idp: an \egt $\gen {x-1, y} = \gen {d}$ would entail $d \divi x-1$, $d \divi
y$, \cad $d$ \ivz, and thus $1 \in \gen {x-1, y}$, which is not the case.


\vspace{-.4em}
\pagebreak	

\prob{exoLambdaGammaK0} \emph{(The operations $\lambda_t$ and $\gamma_t$ over $\KO(\gA)$)}\\
\emph {1.} We have $\lambda_t(\gA)=\lambda_t(1)=1+t$ and $\gamma_t(1)=1/(1-t)$,  
so $\lambda_t(p)=(1+t)^p$ and $\gamma_t(p)=1/(1-t)^p$ for $p\in\NN\etl$.
\\
We write $x$ in the form $[P]-[\gA^p]=P-p$ for a certain $p\in\NN\etl$, \hbox{with $P$}  of constant rank $p$.
By \dfn $\gamma_t([P]) = \sum_{n=0}^p \lambda^n(P) t^n/(1-t)^n$,
we have

\snic {
\gamma_t(x) =\frac{\gamma_t([P])}{\gamma_t(p)} =\sum_{n=0}^p \lambda^n(P) t^n(1-t)^{p-n}.
}

Thus $\gamma_t(x)$ is a \pol of degree $\le p$ in $t$. 
\\
Note:
$\gamma^p(x) = \sum_{n=0}^p \lambda^n(P) (-1)^{p-n}
= (-1)^p \sum_{n=0}^p \lambda^n(P) (-1)^{n}= (-1)^p \lambda_{-1}(P)$.
\\
We have $\gamma_{t}(x) \gamma_{t}(-x)  = 1$ and as they are \pols of $\KO(\gA)[t]$, their \coes of degree $> 0$ are nilpotent (Lemma~\ref{lemGaussJoyal} and Exercise~\ref{exoNilIndexInversiblePol}).
In particular the \elt $x$, which is the \coe of degree $1$ of $\gamma_{t}(x)$, is nilpotent.

\emph {2.}
Let $x \in \KO(\gA)$ be nilpotent, then $\rg x$ is a nilpotent \elt of $\HO(\gA)$. But this last \ri is reduced (actually, \qiriz); thus $\rg x = 0$.

\emph {3.} Suppose $\rg x=[e]$ for some \idm $e$.\\
We have $\Al n (e\gA)=0$ for $n\geq2$, therefore $\lambda_t([e])=1+[e]t$.
By \dfn of 
$a^r$ \hbox{for $a\in\gB$} and $r\in\HO\gB$, we obtain $(1+t)^{[e]}=(1-e)+e(1+t)=1+et$. \\
By direct computation we also obtain $\rR {e\gA}(t)=(1-e)+te$.\\
Finally, we have by convention
$\BB(\gA)\subseteq\HO\gA$
with the identification $e=[e]$.  
\\
We then obtain the \gnl \egt for $x = [P]$ by using the \sfio formed by the \coes of 
$\rR{P}$ and by noting that the two members are morphisms from $\KO(\gA)$ to $1 + t\KO(\gA)[[t]]$. \\
Let us also note that $\lambda_t(p)=(1+t)^p$ for $p\in\NN\etl$ is the desired \egt when $\rg x\in\NN\etl$.

\emph {4.}
Is obtain from item \emph{1}  by replacing $t$ by $t/(1-t)$.

\emph {5.}
Some $x \in \KO(\gA)$ is of the form $y+r$ with $r=\rg x \in\HO\gA$ and $y\in\KTO\gA$.\\
Then $\gamma_{t}(x)=\gamma_{t}(y)(1-t)^{-r}$.

\emph {6.}
Recall the two following formulas, for $d \ge 1$,

\snic {
{1 \over (1-t)^d} = \sum_{k\ge 0} {k+d-1 \choose d-1} t^k,
\qquad
(1-t)^{-d} = \sum_{k\ge 0} {-d \choose k} (-t)^k.
}

They are related by the \egt 

\snic {
{k+d-1 \choose d-1}  = {k+d-1 \choose k}  =
{-d \choose k} (-1)^k.
}

By \dfnz,

\snic {
\gamma_t(x) = 1 + \sum_{d \ge 1} {\lambda^d(x) t^d \over (1-t)^d} =
1 + \sum_{d \ge 1, k \ge 0} \lambda^d(x) t^d 
{k+d-1 \choose d-1} t^k.
}

For $n \ge 1$, the \coe $\gamma^n(x)$ of $t^n$ is

\snic {
\sum_{k+d=n} \lambda^d(x) {k+d-1 \choose d-1} 
\quad \hbox {i.e. with $p = d-1$} \quad
\sum_{p = 0}^{n-1} \lambda^{p+1}(x) {n-1 \choose p}. 
}

The other \egt is deduced via the \eqvc $\gamma_t = \lambda_{t/(1-t)} \iff \lambda_t = \gamma_{t/(1+t)}$.


\prob{exoApplicationProjectiveNoether} \emph{(The \pro map of \Noe and the \mrcs$1$ direct summands in $\gk^2$)}\\
\emph{1.}
Uniqueness of the factorization up to order of the factors and up to \iv \elts.

\emph{2.}
The product of primitive \pols is a primitive \polz, cf. Lemma~\ref{lemGaussJoyal} (Poor man's Gauss-Joyal). We have the more precise result which consists of the inclusion of \ids

\snic {
\gen {x_1, y_1} \cdots \gen {x_n, y_n} \subseteq 
\rD_\gk(\gen {z_0, \ldots, z_n}).
}


We can deduce it from the following fact: if $f$, $g$ are two \pols with an \idtrz, the product of a \coe of $f$ and of a \coe of $g$ is integral over the \id generated by the \coes of the product $fg$ (see Lemma~\ref{lemthKroicl}), 
and in particular it is in the radical of this \idz.  
\\
We can also use the following approach: for $I \subseteq \lrbn$, let $I'$ be its complement, $x_I = \prod_{i \in I} x_i$, $y_I = \prod_{i \in I} y_i$. For $d = \#I$ and  $N = {n \choose d}$, we will show an \egt
$$
\prod\nolimits_{\#I = d}(T - x_Iy_{I'}) = T^N +
\som_{j=1}^{N} a_j T^{N-j}, \quad a_j \in \gen {z_0, \ldots, z_n}.
\leqno (\star')
$$
By making $T=x_Iy_{I'}$, we will have $(x_Iy_{I'})^N \in \gen {z_0, \ldots, z_n}$,
therefore showing the stated inclusion of \idsz. To prove $(\star')$, we first examine the case where all the $y_i$'s are equal to $1$. We write, by letting $S_1(x), \ldots, S_n(x)$ be the  \elr \smq functions of $(\xn)$

\snic {
\prod_{\#I = d}(T - x_I) = T^N + \sum_{j=1}^{N} b_j T^{N-j}, \qquad
b_j = f_j\big(S_1(x), \ldots, S_n(x)\big).
}

A careful examination shows that $f_j$ is a \pol of degree $\le j$ in $(S_1, \ldots, S_n)$.  
Let us replace in this last \egt $x_i$ by $x_i/y_i$ and multiply by $(y_1\cdots y_n)^N$; we obtain, with $U =
y_1\cdots y_nT$ and $s_i = S_i(x_1/y_1, \ldots, x_n/y_n)$

\snic {
\prod_{\#I = d}(U - x_Iy_{I'}) = U^N + 
\sum_{j=1}^{N} (y_1\cdots y_n)^j f_j(s_1, \ldots, s_n) U^{N-j}.
}


Let $s_1^{\alpha_1} \cdots s_n^{\alpha_n}$ be a monomial of $f_j(s_1, \ldots,
s_n)$; since $\sum_i \alpha_i \le \deg f_j \le j$, we obtain, by remembering that $z_n = y_1\cdots y_n$, an \egt

\snic {
z_n^j s_1^{\alpha_1} \cdots s_n^{\alpha_n} =
z_n^{\alpha_0} (z_ns_1)^{\alpha_1} \cdots (z_ns_n)^{\alpha_n} =
z_n^{\alpha_0} z_{n-1}^{\alpha_1} \cdots z_0^{\alpha_n}
\hbox { with } \alpha_0 = j - \sum_i \alpha_i.
}

Since $j \ge 1$, one of the exponents $\alpha_i$ above is not null and we indeed have the membership to $\gen {z_0, \ldots, z_n}$, then the \egt $(\star')$.

\emph{3.}
Let $E = P_1 \te_\gk \cdots \te_\gk P_n \subset L^{n\te}$; it is a \mrcz~$1$. Let us show that the restriction of $\pi$ to $E$ is injective and that $\pi(E)$ is a direct summand in $S_n(L)$. 
This will indeed prove that $\pi(E)$ is a $\gk$-point of $\PP^n$.  Thanks to a finite number of \come
\lonsz, we are brought back to the case where each $P_i$ is free with basis $x_iX + y_iY$. Then each $(x_i, y_i)$ is \umd and $\sum_{i=0}^n z_i X^{n-i}Y^i$ is a \umd basis of $\pi(E)$. This proves on the one hand that $\pi\frt E$ is injective (since it transforms a basis of $E$ into a \umd vector of $S_n(L)$) and that $\pi(E)$ is a direct summand in $S_n(L)$.

\emph{4.}
It seems that $\varphi$ is injective, \cad $(z_0, \ldots, z_n)$ are \agqt independent over $\gk$. The image by $\varphi$ is the graded sub\ri $\gA = \gk[z_0, \ldots, z_n] \subset \gk[\uX,\uY]$ (the \hmg component of an \elt of $\gA$ is in $\gA$); if $f \in \gA$ is \hmg of degree $m$, we have $m \equiv 0 \bmod n$, and for arbitrary $t_1, \ldots, t_n$ 

\snic {
f(t_1X_1, t_1Y_1, \ldots, t_nX_n, t_nY_n) = (t_1\ldots t_n)^{m/n}
f(X_1,Y_1, \ldots, X_n,Y_n)
.}

Finally, $\gA$ is invariant under the action of the \smq group $\rS_n$ which acts on $\gk[\uX, \uY]$ by

\snic {
\sigma \cdot f(X_1,Y_1, \ldots, X_n,Y_n) =
f(X_{\sigma(1)},Y_{\sigma(1)}, \ldots, X_{\sigma(n)},Y_{\sigma(n)})
.}

These last two \prts probably characterize $\gA$.


\prob{exoTh90HilbertMultiplicatif} \emph{(Hilbert's \tho 90, multiplicative form)}\\
We fix once and for all an \elt $b_0 \in \gB$ of trace $1$.

\emph {1}  and \emph {2.}
No difficulty. The fact that $\theta_c$ is multiplicative exactly translates the fact that $c$ is a $1$-cocycle.

\emph {3.}
The action of $G$ over $\gB$ twisted by the $1$-cocycle $c$ is $\sigma\cdot_c b = c_\sigma\sigma(b)$; the fact that this is an action is exactly the condition of $1$-cocyclicity of $c$. Indeed

\snic {
\tau\cdot_c(\sigma\cdot_c b) = \tau\cdot_c c_\sigma \sigma(b) = 
c_\tau\tau\big(c_\sigma \sigma(b)\big) = 
c_\tau\tau(c_\sigma)\, (\tau\sigma)(b) = 
c_{\tau\sigma}\, (\tau\sigma)(b) = (\tau\sigma) \cdot_c b
.}

We will notice that $\pi_c = \sum_\sigma c_\sigma\, \sigma$ is some sort of $G$-trace relatively to the action of $G$ twisted by $c$.

We therefore have $\gB_c^G = \sotq {b \in \gB} {c_\sigma \sigma(b) = b}$.  By using the fact that $c$ is a $1$-cocycle, we find that $\tau \circ \pi_c = c_\tau^{-1}\, \pi_c$; we deduce that $c_\tau \tau(z) = z$ for every $z \in \Im\pi_c$, \cad $\Im\pi_c \subseteq \gB_c^G$. We define $s : \gB_c^G \to \gB$ by $s(b) = bb_0$. Then $\pi_c \circ s = \Id_{\gB_c^G}$; indeed, for $b \in \gB_c^G$,

\snic {
\pi_c(b_0b) = \sum_\sigma c_\sigma\sigma(bb_0) =
\sum_\sigma c_\sigma\sigma(b)\sigma(b_0) = \sum_\sigma b\sigma(b_0) 
= b\Tr_{\gB\sur\gA}(b_0) = b
.}

From the \egt $\pi_c \circ s = \Id_{\gB_c^G}$, we deduce that $\pi_c$ is a surjection from $\gB$ to $\gB_c^G$, that~$s$ is injective and that $\gB = s(\gB_c^G) \oplus \Ker\pi_c \simeq \gB_c^G \oplus \Ker\pi_c$.
In particular, $\gB_c^G$ is \hbox{a \ptf \Amoz}.

\rem Let us consider $s : b \mapsto b_0b$ in $\End_\gA(\gB)$, then $(\pi_c \circ s)\big(\pi_c(z)\big) = \pi_c(z)$
for all $z \in \gB$, \cad $\pi_c \circ s \circ \pi_c = \pi_c$.
Consequently $\pi'_c \eqdf {\rm def} \pi_c \circ s = \sum_\sigma c_\sigma \sigma(b_0\bullet)$ is a \prrz;
we could certainly compute its trace and find $1$, which would prove that $\pi'_c$ is a \prr of rank $1$.

\emph {4.}
Let $c$, $d$ be two $1$-cocycles, $x \in \gB_c^G$, $y \in \gB_d^G$, so $c_\sigma\sigma(x) = x$, $d_\sigma\sigma(y) = y$; we easily verify that $xy \in \gB_{cd}^G$. \\
Hence an \Ali $\gB_{c}^G \otimes_\gA \gB_{d}^G \to \gB_{cd}^G$, $x\otimes y \mapsto xy$, denoted by~$\mu_{c,d}$.
\\
Let $(x_i)$, $(y_i)$ be two \syss of \elts of $\gB$ like in Lemma~\ref{lemArtin} and let $\varepsilon = \sum_i x_i \otimes y_i = \sum_i y_i \otimes x_i$ (separability \idmz). Recall that $\varepsilon \in \Ann(\rJ)$, which translates to

\snic {
\forall\ b \in \gB\quad
\sum_i bx_i \otimes y_i = \sum x_i \otimes by_i  \quad\hbox{in}\quad
\env\gA\gB \eqdf {\rm def} \gB\otimes_\gA \gB
.}

We also have, for $b$, $b' \in \gB$

\snic {
\Tr_{\gB\sur\gA}(bb') =\sum_i \Tr_{\gB\sur\gA}(bx_i)\Tr_{\gB\sur\gA}(b'y_i)
.}

We will show that $z \mapsto (\pi_c \otimes \pi_d)(b_0z\varepsilon)$,
$\gB_{cd}^G \mapsto \gB_{c}^G \otimes_\gA \gB_{d}^G$ and $\mu_{c,d}$ are   
         reciprocals of one another. 
In one direction,

\snic {
(\pi_c \otimes \pi_d)(b_0z\varepsilon) = \sum_i a_i \otimes b_i, 
\quad \hbox {with} \quad
a_i = \sum_\sigma c_\sigma\sigma(b_0zx_i), \quad
b_i = \sum_\tau c_\tau\tau(y_i)
,}

and we have

\snic {
\sum_i a_ib_i = \sum_{\sigma,\tau} \sigma(b_0z) c_\sigma d_\tau
\sum_i \sigma(x_i) \tau(y_i)
,}

and since the internal sum (over $i$) evaluates to $1$ or $0$, there remains,
for $z \in \gB_{cd}^G$ 

\snic {
\sum\limits_i a_ib_i = \sum\limits_{\sigma} \sigma(b_0z) c_\sigma d_\sigma =
\sum\limits_{\sigma} \sigma(b_0) \sigma(z) (cd)_\sigma =
\sum\limits_{\sigma} \sigma(b_0) z = z \Tr_{\gB\sur\gA}(b_0) = z
.}

In the other direction, let $x \in \gB_c^G$ and $y \in \gB_d^G$. Then,
since $\varepsilon \in \Ann(\rJ)$, we can write

\snic {
(\pi_c \otimes \pi_d)(b_0 xy\varepsilon) = \sum_i a_i \otimes b_i, 
\; \hbox {with} \;
a_i = \sum_\sigma c_\sigma\sigma(b_0xx_i), \;
b_i = \sum_\tau d_\tau\tau(yy_i)
.}

By using 
$$\preskip.4em \postskip.4em
c_\sigma \sigma(b_0xx_i) = c_\sigma \sigma(x) \sigma(b_0x_i)
= x\sigma(b_0x_i) \hbox{ and } d_\tau \tau(yy_i) = d_\tau
\tau(y)\tau(y_i) = y\tau(y_i),
$$
we get
\[\preskip.0em \postskip.4em 
\begin{array}{ccc} 
  \sum_i a_i \otimes b_i = 
\sum_i x\Tr_{\gB\sur\gA}(b_0x_i) \otimes y\Tr_{\gB\sur\gA}(y_i)  =   \\[1mm] 
(x \otimes y) \cdot
\big(\sum_i \Tr_{\gB\sur\gA}(b_0x_i)\Tr_{\gB\sur\gA}(y_i) \otimes 1\big) =
(x \otimes y) \cdot \big(\Tr_{\gB\sur\gA}(b_0) \otimes 1\big) = x\otimes y
. 
\end{array}
\]
Item~\emph {a}  is proved. \\
For item~\emph {b}, let there be a $1$-cocycle, coboundary of $b_1 \in \Bti$, $c_\sigma = \sigma(b_1)b_1^{-1}$. \\
Then $b \in \gB_c^G$ \ssi for every $\sigma$, $c_\sigma\sigma(b) = b$,
\cad $\sigma(b_1b) = b_1b$ \cad $b_1b \in \gA$; so $\gB_c^G = b_1^{-1}\gA$.
We deduce that $\gB_c^G \otimes \gB_{c^{-1}}^G \simeq \gA$, 
\hbox{so $\gB_c^G$}
is an \Amrc $1$. \\
Moreover $c \mapsto \gB_c^G$ induces a morphism $\zcoho \to \Pic(\gA)$.

It remains to show that if $\gB_c^G$ is free, \cad $\gB_c^G = \gA b_1$
with $b_1 \in \gB$ and $\Ann_\gA(b_1) = 0$, then $c$ is a coboundary.
But $\gB_{c^{-1}}^G$, being the \inv of $\gB_c^G$ is also free, $\gB_{c^{-1}}^G = \gA b_2$, and 
$\gB_c^G \gB_{c^{-1}}^G = \gB_1^G = \gA$.
We therefore have $\gA b_1b_2 = \gA$, then $b_1$, $b_2$ are \ivs in~$\gB$ (and $\gA b_2 = \gA b_1^{-1}$). Then $c_\sigma^{-1} \sigma(b_2) = b_2$, \cad $c$ is the coboundary of $b_2$.

\emph {5.}
Since $\gA$ is a \zed \riz, $\Pic(\gA) = 0$ so $\hcoho = 0$.

\emph {6.}
Let $c_\tau = x\sigma(x) \cdots\sigma^{i-1}(x)$ with $i \in \lrb{1..n}$ and $\tau = \sigma^i$.\\
Thus, $c_\Id = \rN_{\gB\sur\gA}(x) = 1$, $c_\sigma = x$, $c_{\sigma^2} = x\sigma(x)$. \\
It is a $1$-cocycle: $c_\sigma \sigma(c_{\sigma^i}) = c_{\sigma^{i+1}}$, \cad $c_\sigma\sigma(c_\tau) = c_{\sigma\tau}$, then $c_{\sigma^j}\sigma^j(c_\tau) = c_{\sigma^j \tau}$.

\prob{exoSegreMorphism} \emph{(The Segre morphism in a special case)}\\
It is clear that $\fa \subseteq \Ker\varphi$.

\emph{1.}
Let $m = X_{i_1} \cdots X_{i_r} Y_{j_1} \cdots Y_{j_s}$, $m' = X_{i'_1}
\cdots X_{i'_{r'}} Y_{j'_1} \cdots Y_{j'_{s'}}$ with
$$\preskip-.0em \postskip.2em
1 \le i_1 \le \cdots \le i_r \le j_1 \le \cdots \le j_s \le n, \;
1 \le i'_1 \le \cdots \le i'_{r'} \le j'_1 \le \cdots \le j'_{s'} \le n.
$$
The \egt $\varphi(m) = \varphi(m')$ provides
$$\preskip.2em \postskip.3em
T^r U^s Z_{i_1} \ldots Z_{i_r} Z_{j_1} \ldots Z_{j_s} =
T^{r'} U^{s'} Z_{{i'}_{\!1}} \ldots Z_{{i'}_{\!r'}} Z_{j'_1} \ldots Z_{{j'}_{\!s'}}
.
$$
Therefore $r = r'$, $s = s'$ then $i_k = i'_k$ and $j_\ell = j'_\ell$.
Ultimately $m = m'$.
\\
Let $s = \sum_\alpha a_\alpha m_\alpha$ be an $\gA$-\lin combination of normalized \moms such that $\varphi(s) = 0$. As the \moms $\varphi(m_\alpha)$ are pairwise distinct, we have $a_\alpha = 0$, \cad $s = 0$.

\emph{2.}
Since $X_iY_j \equiv X_jY_i \mod {\fa}$, we see that every \mom is equivalent modulo~$\fa$ to a normalized \momz. We therefore get $\gA[\uX,\uY] = \fa + \fa_{\rm nor}$.
As $\fa \subseteq \Ker\varphi$, the sum is direct by the previous question.

\emph{3.}
Let $h \in \Ker\varphi$ which we decompose into $h = f + g$ with $f \in \fa$,
$g \in \fa_{\rm nor}$.\\
Since $\fa \subseteq \Ker\varphi$, we have $g \in
\Ker\varphi$, so $g = 0$.  Conclusion: $h = f \in \fa$, which proves $\Ker\varphi \subseteq \fa$, then $\Ker\varphi =\fa$.

\prob{exoVeroneseMorphism} \emph{(The Veronese morphism in a special case)}
\\
It is clear that $\fb \subseteq \fa \subseteq \Ker\varphi$.

\emph{1.}
         Let $f$ be in the intersection; $f$ is of the form $f = f_0 + \sum_{i=1}^{d-1} f_i X_i$ with $f_i \in \gA[X_0,X_d]$; We write that $\varphi(f) = 0$

\snic {
f_0(U^d,V^d) + f_1(U^d,V^d)U^{d-1}V + \cdots + f_{d-1}(U^d,V^d)UV^{d-1} = 0.
}

This is of the form, in $\gA[U][V]$, $h_0(V^d) + h_1(V^d)V + \cdots + h_{d-1}(V^d)V^{d-1} = 0$; by examining in this \egt the exponents of $V$ modulo $d$, we obtain $h_0 = h_1 = \cdots = h_{d-1} = 0$. Recap: $f_i = 0$
then $f = 0$.

\emph{2.}
We work modulo $\fb$ by letting

\snic{\Aux = \AuX/\fb$, $\gB = \gA[x_0,x_d]
+ \gA[x_0,x_d]x_1 + \cdots + \gA[x_0,x_d]x_{d-1} \subseteq \Aux.}

We will show that $\gB$ is an \Aslgz; as it contains the $x_i$'s, it is all \hbox{of the $\Aux$}. It suffices to prove that $x_ix_j \in \gB$ for $i \le j \in\lrb{1.. d-1}$, because the other products are in $\gB$ by \dfn of $\gB$. We use the syzygies
$x_ix_j  = x_{i-1}x_{j+1}$ for $i \le j \in\lrb{1.. d-1}$. We have $x_0x_k \in \gB$ for every $k$; we deduce $x_1x_j \in \gB$ for all $j \in\lrb{1.. d-1}$ and it is still true for $j = d$ and $0$ by \dfn of $\gB$. We then deduce that $x_2x_j \in \gB$ for $j \in\lrb{2.. d-1}$, and so on.
\\
The obtained \egt $\gB = \Aux$ is written as

\snic{\AuX = \fb + \big(\gA[X_0,X_d]
\oplus \gA[X_0,X_d]X_1 \oplus \cdots \oplus \gA[X_0,X_d]X_{d-1}\big),}

and the $+$ represents a direct sum by item~\emph{1}  (since $\fb \subseteq \Ker\varphi$).

\emph{3.}
Let $h \in \Ker\varphi$ which we decompose into $h = f + g$ as above. \\
Since $f\in\fb \subseteq \Ker\varphi$, we have $g \in \Ker\varphi$, so $g = 0$.
Conclusion: $h = f \in \fb$, which proves $\Ker\varphi \subseteq \fb$, then $\Ker\varphi = \fb = \fa$.


\prob{exoVeroneseMatrix} \emph{(Veronese matrices)}\\
\emph {2.} 
It is clear that $V_d(P)$ is a \prr if $P$ is also a \prrz, and the diagram is commutative for functorial reasons.  \\
We can add the following precision: if $P$, $Q \in \Mn(\gk)$ are two \prrs such that $\Im P \subseteq \Im
Q$, then $\Im V_d(P) \subseteq \Im V_d(Q)$. Indeed, \hbox{we have $\Im P \subseteq \Im Q$}  
\ssi $QP = P$, and we deduce that $V_d(Q)V_d(P) = V_d(P)$, 
\cad $\Im V_d(P) \subseteq \Im V_d(Q)$.

\emph {3.} 
It suffices to do this \lotz, \cad to compute $V_d(A)$ when $A$ is a standard \prr $\I_{r,n}$. If $A=\Diag(a_1, \ldots, a_n)$, then $V_d(A)$ is diagonal, with diagonal the $n'$ monomials $a^{\alpha}$ with $\abs{\alpha}= d$. In particular, for $A = \I_{r,n}$, we see that $V_d(A)$ is a standard projection, of rank the number of $\alpha$
such that $\alpha_1 + \cdots + \alpha_r = d$, \cad ${d+1-r \choose r-1}$, and $V_d(\I_{1,n}) = \I_{1,n'}$.


\prob{exoFossumKumarNori} \emph{(Some examples of finite projective resolutions)}\\
\emph {1.}
The computation of $F_k^2 - F_k$ is done by \recu and poses no \pbz.
For the conjugation ($n \ge 1$), we use

\snic {
\cmatrix {0& -\I\cr \I & 0} \cmatrix {A & B\cr C& D} \cmatrix {0& \I\cr -\I& 0} =
\cmatrix {D & -C\cr -B & A} 
.}

For $\cmatrix {A & B\cr C& D} = F_n$, this provides a conjugation between $F_n$ and $\I_{2^n} - \tra{F_n}$.\\
When $z(z-1) + \sum_{i=1}^n x_iy_i = 0$, the \prrs $F_n$ and $\I_{2^n} - {F_n}$ have for image \mptfs $P$ and $Q$ with $P\oplus Q\simeq \Ae{2^n}$ and $P\simeq Q\sta$. \\
 Therefore $2\rg(P)=2^{n}$, and since $mx=0\Rightarrow x=0$ for $m\in\NN\etl$ and $x\in\HO \gA$, we obtain $\rg(P)=2^{n-1}$.

\emph {2.}
The computation of $U_kV_k$ and $V_kU_k$ is done by \recuz.  The fact that $F_n$ and $G_n$ are conjugated by a permutation matrix is left to the sagacity of the reader. \\
For example, $G_2 = P_\tau F_2 P_\tau^{-1}$ for $\tau = (2,4,3) =(3,4)(2,3)$, and $G_3 = P_\tau F_3 P_\tau^{-1}$ {for $\tau = (2, 4, 7, 5)(3, 6)=(3,6)(2,4)(4,7)(5,7)$}.  \\
Regarding the constant rank $2^{n-1}$ we can invoke item \emph{1}, or make the direct computation after \lon at $z$ and at $\ov z = 1-z$. 





\emph {3a.}
Direct use of the referenced exercise.

\emph {3b.}
Let $S$ be the \mo $a^\NN$. We can localize a finite projective resolution of $M$ over $\gA$ to obtain one over $S^{-1}\gA$ 

\snic{0 \rightarrow S^{-1}P_n
\to \cdots \to S^{-1}P_1 \to S^{-1}P_0 \twoheadrightarrow S^{-1}M \to
0.}

As $aM = 0$, we have $S^{-1} M = 0$, therefore $\sum_{i=0}^n (-1)^i \rg (S^{-1}P_i) = 0$. But the natural morphism $\HO(\gA) \to \HO(S^{-1}\gA)$ is injective. \\
Therefore $\sum_{i=0}^n (-1)^i \rg P_i = 0$.

\emph {4.}
The localized \ri $(\gB_n)_z$ contains all the $y'_i= y_i/z$, and since $z(1-z) = \sum_i x_iy_i$, \hbox{we have $1-z = \sum_i x_iy'_i$}.
Therefore $z \in \gk[\xn, y'_1, \ldots, y'_n]$ and  $1 - \sum_i x_iy'_i\in(\gB_n)_z\eti$. We then verify that 

\snic {
(\gB_n)_z =  \gk[\xn, y'_1, \ldots, y'_n]_{s}
\quad \hbox {with} \quad s = 1 - \sum x_iy'_i
.}

Similarly, $(\gB_n)_{1-z} = \gk[x'_1, \ldots, x'_n, \yn]_{1 - \sum x'_iy_i}$
with $x'_i = x_i/(1-z)$.

\emph {5.}
For $n \ge 1$, every \elt $a \in \{z, \xn\}$ is \ndz and $a(\gB_n/\fb_n)
= 0$.  As $F_1=\cmatrix {z & x_1\cr y_1 & \ov z}$ is a \prrz, we have $[z,x_1]
F_1 = [z,x_1]$. The reader will check that $\Ker\,[z,x_1] = \Ker F_1 = \Im (\I_2 - F_1)$; hence the exact sequence

\snic {
0 \rightarrow \Im(\I_2-F_1) \to \gB_1^2
\vvers{[z,x_1]} \gB_1 \twoheadrightarrow \gB_1\sur{\fb_1} \to 0
.}

We indeed have $\rg (\gB_1\sur{\fb_1}) = 1-2+1 = 0$.

\emph {6.}
Let $A$ be the matrix constituted of the first 3 rows of $\I_4 - F_2$

\snic {
A = \cmatrix {1-z & -x_1 & -x_2 & 0\cr
-y_1 & z & 0 & -x_2\cr -y_2 & 0 & z & x_1\cr}
.}

It is clear that $A F_2 = 0$ and $[z,x_1,x_2]A = 0$. The reader will check that the sequence below is exact

\snic {
0 \rightarrow \Im F_2 \to \gB_2^4 \vvers {A} \gB_2^3 \vvvvvers{[z,x_1,x_2]} 
\gB_2 \twoheadrightarrow \gB_2\sur{\fb_2} \to 0
.}

We indeed have $\rg (\gB_2\sur{\fb_2}) = 1-3+4-2 = 0$.

\emph {7.}
Immediate given the \dfn of $F_n$.

\emph {8.}
Consider the upper half of the matrix $\I_8 - F'_3$ and delete its last (zero) column to obtain a matrix $A$ of format $4 \times 7$. Let $B$ be the matrix of format $7 \times 8$ obtained by deleting the last row of $F'_3$.  Then the brave reader will check the exactness of

\snic {
0 \rightarrow \Im (\I_8-F'_3) \to \gB_3^8 \vvers {B} \gB_3^7 \vvers {A}
\gB_3^4 \vvvvvers{[z,x_1,x_2,x_3]}  \gB_3 \twoheadrightarrow \gB_3\sur{\fb_3} \to 0
.}

We have $\rg(\gB_3\sur{\fb_3}) = 1-4+7-8+4 = 0$.

\emph {9.}
There is an exact sequence (let $\gB = \gB_n$, $\fb = \fb_n$):

\snic {
L_{n+1} \vvers {A_{n+1}} L_n \vvers {A_n} L_{n-1} \vvers {A_{n-1}} \cdots
\lora L_2 \vvers {A_2} L_1 \vvers{A_1} L_0=\gB 
\twoheadrightarrow \gB\sur{\fb}
.}

where $L_r$ is a free module of rank $\sum_{i \in I_r} {n+1 \choose i}$ with
$I_r = \sotq {i \in \lrb{0..r}}{i \equiv r \bmod 2}$.  In particular, $L_1 = \gB^{n+1}$ and $L_n = L_{n+1} = \gB^{2^n}$.  \\
As for the matrices $A_r$, we have $A_1 = [z,\xn]$, and the matrix $A_r$ is extracted from~$F_n$ if $r$ is odd,
and extracted from $\I - F_n$ otherwise. We have $A_{n+1} = F_n$ for even $n$, and $A_{n+1} = \I - F_n$ for odd $n$.
\\
By letting $P_{n+1} = \Im A_{n+1}$, the \Bmo $\gB\sur{\fb}$ admits a projective resolution of length $n+1$ of the following type

\snic {
0 \rightarrow P_{n+1} \to L_n = \gB^{2^n} \vers {A_n} L_{n-1} \vvers {A_{n-1}} \cdots
\to L_2 \vers {A_2} L_1 \vers{A_1} L_0=\gB 
\twoheadrightarrow \gB\sur{\fb}
.}

($P_{n+1}$ of constant rank  $2^{n-1}$).\\
	The explicit expression of the rank of $L_i$ confirms that $[\gB\sur{\fb}] \in
\KTO(\gB)$.  \\
We have $\rg L_{n-1} + \rg L_{0} = \rg L_{n-2} + \rg L_1 =
\cdots = 2^n$ (in particular, if $n = 2m+1$, then $\rg L_m =
2^{n-1}$).

Note:
If $\gk$ is a \cdiz, we can show that $\KTO(\gB_n) \simeq \ZZ$
with as a \gtr $[\gB_n\sur{\fb_n}]$.  We deduce that the \id $\KTO(\gB_n)$ has a null square; \gnltz, let $\gA$ be a \ri satisfying
$\KTO(\gA)=\ZZ x \simeq \ZZ$, then $x^2 = mx$ with $m \in \ZZ$, so $x^{k+1} = m^k x$ for $k \ge 1$, since $x$ is nilpotent (see \Pbmz~\ref{exoLambdaGammaK0}), there is some $k\ge 1$ such that $m^k x = 0$, so $m^k = 0$, then $m=0$ and $x^2 = 0$.


\Biblio

\Thref{th rg const loc free} specifies \Tho 2 in \cite{Bou} Chap.~II~{\S}5. \perso{d\'evelopper ce genre of r\'ef\'erences, pr\'eciser aussi le titre of the livre concern\'e}

Section~\ref{secAppliIdenti} is based on the articles \cite[Chervov\&Talalaev]{CT03,CT03bis} 
which examine the \gui{Hitchin \syssz} over the singular curves.

\Pbmz~\ref{exoLambdaGammaK0} is inspired from a non-published article of 
R.~G.~Swan: \emph{On a theorem of Mohan, Kumar and Nori}.

\Pbmz~\ref{exoTh90HilbertMultiplicatif} comes from an exercise of Chapter 4 of \cite{Jensen}.

In \Pbmz~\ref{exoFossumKumarNori}, the matrix $F_k$ occurs in the article: \emph{Vector bundles over
Spheres are Algebraic}, {\sc R. Fossum}, Inventiones Math. {\bf 8}, 222--225 (1969).
The \ri $\gB_n$ is a classic in algebraic K-theory.

\newpage \thispagestyle{CMcadreseul}
\incrementeexosetprob


\newcommand \dotdiv {\; {{\,.\,} \over {} } \;}

\chapter[Distributive lattices, lattice-groups]{Distributive lattices Lattice-groups}
\label{chapTrdi}
\perso{compil\'e le \today}
\minitoc

\Intro
\pagestyle{CMExercicesheadings}

This chapter begins with an introductory section which fixes the formal \agq framework of \trdis and of \agBsz.

\smallskip
The \trdis are important in commutative \alg for several reasons.

On the one hand the theory of \dve has as its \gui{\id model} the theory  of \dve of natural integers. If we take as the order relation   $a \preceq b$, the relation \gui{$a$ is a multiple of $b$,}
 we obtain that $\NN$ is a \trdi with: minimum \elt $0$, maximum \elt $1$,
the  supremum $a\vu b$ equal to the gcd and the infimum $a\vi b$ equal to the lcm. A few beautiful \prts of \dve in $\NN$ are expressed in modern terms by saying that the \ri $\ZZ$ is a Bézout \ri (see Sections~\ref{secApTDN} and 
\ref{secBézout}).
The \id numbers in number theory have been created by Kummer to fill the gap between the theory of \dve in the 
ring of integers of a  number field
and that in $\NN$.
The ring of integers of a number field is not a Bézout \ri in \gnlz, but its \itfsz%
\footnote{%
What for Kummer was \gui{the ideal gcd of several numbers} has been replaced in modern language by the corresponding finitely generated ideal. This tour de force, due to Dedekind, was one of the first intrusions of the \gui{actual} infinite in mathematics.%
} 
form a \trdiz, and their nonzero \itfs form the non-negative submonoid of an \grl (see Section~\ref{secGpReticules}) which re-establishes the well-ordering of things.
The \ris whose \itfs form a \trdi are called \anars (treated elsewhere in Sections~\ref{secIplatTf} and \ref{secAnars}).
Their \iv \ids also form the non-negative submonoid of an \grlz.
The theory of GCD-domains  (Section~\ref{secAnnPgcd})
also finds its natural framework in the context of \grlsz.

On the other hand the \trdis intervene as the \cov counterpart of diverse and various spectral spaces which are imposed as powerful tools of the abstract \algz.
The discussion on this subject is particularly enlightening when we consider the Zariski lattice of a commutative \riz,
relatively unknown, which serves as a \cov counterpart to the very famous Zariski spectrum:
spectral space that we could believe indispensable to the \ddk theory and to the Grothendieck scheme theory.
A systematic study of the Zariski lattice will be given in Chapter~\ref{chapKrulldim} regarding the Krull dimension, with a heuristic introduction in Section~\ref{secEspSpectraux}.
In Section~\ref{secZarAcom} we define the Zariski lattice of a commutative \ri $\gA$ essentially with respect to the construction of the reduced \zede closure $\Abul$ (\paref{secClotureZEDR}) of the \riz. This construction can be regarded as a construction \paral to that of the \agB
generated by a \trdi (see \thref{thZedGenEtBoolGen}).
The global object $\Abul$ constructed thus essentially contains the same information as the product of \risz~$\Frac(\gA\sur\fp)$ for all the \ideps $\fp$ of $\gA$.
We get this even though in the \gnl situation we do not \cot have access to the \ideps of a \ri individually.

Another reason to be interested in \trdis is the \cov (or intuitionist) logic in which the set of truth values of classical logic, that is $\so{\Vrai,\Faux}$,
which is a \agB with two \eltsz, is replaced with a more mysterious 
\trdiz.%
\footnote{Actually the truth values of \coma do not strictly speaking form a set, but a class.
Nevertheless the \cof logical connectives act on those truth values with the same \agq \prts as the $\vi$, the $\vu$ and the $\to$ of \agHsz. See the discussion on \paref{P(X)}.\label{NoteValVer}}
The \cov logic will be addressed in the Appendix (see \paref{chapPOM}), particularly in Sections~\iref{secAffirmerProuver} and \iref{secBHK}.
In Section~\ref{secEntRelAgH} of the previous chapter we implement the tools that serve as the framework for a formal \agq study of this logic: the \entrels and the \agHsz.
It is remarkable that Heyting defined those \algs in the first attempt to describe the intuitionist logic formally, and that there has not been a comma to add since.
Moreover, \entrels and \agHs are also useful 
in the \gnle study of \trdisz. For example it is sometimes important to be able to say that the Zariski lattice of a \ri is a \agHz.


\section{Distributive lattices and \agBsz}
\label{secTrdis}\label{secAGB}
\pagestyle{CMheadings}


In an ordered set $X$ we let, for some $a\in X$,
\begin{equation}\label{eqda}
\dar a=\sotq{x\in X}{x\leq a},\quad \uar a=\sotq{x\in X}{x\geq a}.
\end{equation}

We call a finite non-decreasingly ordered list $(a_0,\ldots,a_n)$ of elements of $X$ a \emph{non-decreasing chain}.
The number $n$ is called the \emph{length} of the chain. By convention the empty list is a chain of length $-1$.
\index{chain!in an ordered set}
\index{length of a chain!in an ordered set}

\begin{definition}
\label{deftrdi}~
\begin{enumerate}
\item A \ix{lattice} is a set $\gT$ equipped with an order relation $\leq$ for which every finite family admits an upper bound and a lower bound. Let $0_\gT$ be the minimum of $\gT$ (the upper bound of the empty family) and  $1_\gT$ be the maximum of $\gT$. Let $a\vu b$ be the upper bound of $(a,b)$ and~$a\vi b$ be its lower bound.
  \item A map from one lattice to another is called a \emph{lattice \homoz} if it respects the laws~$\vu$ and~$\vi$ and the constants $0$ and $1$.
  \item The lattice is called a \ixy{distributive}{lattice} when each of the two laws~$\vu$ and~$\vi$ is
distributive with respect to the other.
\end{enumerate}
\end{definition}

The axioms of lattices can be formulated with universal \egts uniquely regarding the two laws $\vi$ and $\vu$ and the two constants $0_\gT$ and $1_\gT$. The order relation is then defined by $a\leq_\gT b\equidef a\vi b=a$. Here are those axioms.
\[\arraycolsep2pt
\begin{array}{rclcrcl}
a \vu a   & = & a &\qquad\qquad& a \vi a   & = & a \\
a \vu b   & = &b\vu a &\quad\quad& a \vi b   & = &b\vi a \\
(a \vu b)\vu c   & =&a \vu(b\vu c)  &&  (a \vi b)\vi c   & =&a \vi(b\vi c) \\
(a \vu b)\vi a   & = &a && (a \vi b)\vu a   & = &a \\
a\vu 0_{\gT}   & =  &a & &a\vi 1_{\gT}   & =  &a
\end{array}
\]

We thus obtain a purely equational theory, with all the related facilities.
For example we can define a lattice by \gtrs and relations.
Similarly for the \trdisz.

\rdb
In a lattice, one distributivity implies the other.
Suppose for instance that $a\vi(b\vu c)=(a\vi b)\vu(a\vi c)$, for all $a,b,c$. Then the other \dit results from the following computation
\label{DistriTrdi}
\[
\begin{array}{c}
 (a\vu b)\vi(a\vu c)= \big((a\vu b)\vi a\big)\vu \big((a\vu b)\vi c\big)=a \vu \big((a\vu b)\vi c\big)=
 \\[1mm]
 a \vu\big(( a\vi c)\vu ( b\vi c) \big)= \big(a \vu( a\vi c)\big)\vu ( b\vi c)=a\vu  ( b\vi c).
\end{array}
\]

In a discrete lattice we have a test for $a\leq b$, since this relation is equivalent to $a\vi b=a$.

The subgroups of a group (or the \ids of a commutative \riz) form a lattice \wirt the inclusion, but it is not a \trdi in \gnlz.

A totally ordered set%
\footnote{Recall that this is a set $E$ equipped with an order relation $\leq$ for which we have, for all $x$ \hbox{and $y\in E$},
$x\leq y$ or $y\leq x$. This does not imply that the \egt is decidable.} 
is a \trdi if it possesses a maximum \elt and a minimum \eltz.
Let ${\bf n}$ be the totally ordered set with $n$ \eltsz.
A map between two totally ordered lattices $\gT$ and $\gS$ is a \homo \ssi it is non-decreasing and~$0_\gT$  and $1_\gT$ have as their images $0_\gS$ and $1_\gS$.

If $\gT$ and $\gT'$ are two \trdisz, the set $\Hom(\gT,\gT')$ of \homos from $\gT$ to $\gT'$ is equipped with a natural order structure given~by
$$ \varphi \leq \psi \equidef  \forall x\in \gT\;
\;
\varphi(x) \leq \psi(x)
$$

A cartesian product of \trdis is a \trdi (for the product laws $\vi$ and $\vu$, which gives the product partial order relation).

\rdb\label{NOTAtrdiopp}
For every \trdi $\gT$, if we replace the order relation $x\leq_\gT y$ by the symmetric
relation  $y\leq_\gT x$ we obtain the \ixy{opposite}{lattice} $\gT\eci$ with an exchange of $\vi$ and $\vu$ (we sometimes say \emph{dual lattice}).

If $A\in\Pfe(\gT)$  with a \trdi $\gT$ we will let
$$\preskip.4em 
 \Vu A:=\Vu_{x\in A}x\qquad {\rm and}\qquad \Vi A:=\Vi_{x\in A}x.
$$

\subsec{Quotient lattices, \idsz, filters}

\index{distributive lattice!quotient ---}
If an \agq structure is defined by laws of composition of different arities and by axioms that are universal \egts (such as groups, rings and distributive lattices), a quotient structure is obtained when we have an equivalence relation and when the laws of composition \gui{pass to the quotient.}
If we look at the structure as defined by \gtrs and relations (which is always possible), we obtain a quotient structure by adding relations.

A \emph{quotient lattice $\gT'$ of a lattice $\gT$} can \egmt be given by a binary relation $\preceq$ over $\gT$ satisfying the following \prts
\begin{equation}\label{eqPreceq}\preskip.3em \postskip.4em
\left.
\begin{array}{rcl}
a\leq b&  \Longrightarrow  & a\preceq b   \\
a\preceq b,\,b\preceq c&  \Longrightarrow  & a\preceq c   \\
a\preceq b,\,a\preceq c&  \Longrightarrow  & a\preceq b\vi c   \\
b\preceq a,\,c\preceq a&  \Longrightarrow  & b\vu c\preceq a
\end{array}
\right\}
\end{equation}
The relation $\preceq$ then induces a lattice structure over the quotient set~$\gT'$ obtained with the new \egtz%
\footnote{%
The fact that, when passing to the quotient, we change only the equality relation and not the objects is simpler than the classical approach, and is more consistent with the (Gaussian) tradition and with machine implementation. No doubt the popular success of equivalence classes as objects of the quotient set is largely due to the fortunate fact that in the case of a quotient group $G/H$, in additive notation for example, we have $(x+H)+(y+H) = (x+y)+H$ where the symbol $+$ has three different meanings. However, things are less fortunate in the case of quotient rings. For example, $(3+7\ZZ)(2+7\ZZ)$ is contained within, but is not equal to $6+7\ZZ$.
}
$$\preskip.2em \postskip.4em
(a,b\in\gT)\quad:\quad\quad a=_{\gT'}b \equidef (a\preceq b \;\mathrm{ and }\;
b\preceq a)
$$
Naturally if $\gT$ is distributive, the same goes for $\gT'$.

If $\varphi :\gT\rightarrow \gT'$ is a \trdi \homoz, $\varphi^{-1}(0)$ is called an \emph{\id of $\gT$}. An \id $\fb $ of $\gT$ is a subset of $\gT$ subjected to the following constraints
\begin{equation}\label{eqIdeal}\preskip-.4em 
\left.
\begin{array}{rcl}
   & &  0 \in \fb    \\
x,y\in \fb & \Longrightarrow   &  x\vu y \in \fb    \\
x\in \fb ,\; z\in \gT& \Longrightarrow   &  x\vi z \in \fb    \\
\end{array}
\right\}
\end{equation}
(the last is rewritten as $(x\in \fb ,\;y\leq x)\Rightarrow y\in
\fb$).
A \emph{\idpz} is an \id generated by a single \elt $a$,
it is equal to $\dar a$.%
\index{ideal!of a \trdiz}%
\index{ideal!principal --- (of a \trdiz)}\label{NOTAdara}%
\index{principal!ideal of a \trdiz}

The \id $\dar a$, equipped with the laws $\vi$ and $\vu$ of $\gT$, is a \trdi in which the maximum \elt is $a$.  The canonical injection $\dar a\rightarrow \gT$ is not a morphism of \trdis because the image of $a$ is not equal to $1_\gT$. However, the 
map $\gT\rightarrow \dar a,\;x\mapsto x\vi a$ is a surjective morphism, which therefore defines $\dar a$ as a quotient structure.

\rdb \label{NOTAuara}
The notion opposite to that of an \id is the notion of a {\em filter}. The principal filter generated by $a$ is equal to $\uar a$.%
\index{filter!of a distributive lattice}%
\index{filter!principal --- of a distributive lattice}%
\index{principal!filter of a distributive lattice}

The \emph{\id generated} by a subset $J$ of $\gT$ is equal to
$$
\cI_\gT(J)=\big\{x\in\gT \;\big \vert \; \Ex J_0\in \Pfe(J),\,x\leq \Vu J_0\big\}.
$$
Consequently \emph{every \itf is principal}.
\label{NOTAidtrdi}

If $A$ and $B$ are two subsets of $\gT$ let
\begin{equation}\label{eqvuvi}
A\vu B=\big\{ a\vu b \mid a\in A,\,b\in B\,\big\}  \; \hbox{ and } \; A\vi
B=\big\{ a\vi b \mid a\in A,\,b\in B\,\big\}.
\end{equation}
Then the \id generated by two \ids $\fa$ and $\fb$ is equal to
\begin{equation}\label{eqSupId}
\cI_\gT(\fa\cup \fb) = \fa\vu\fb.
\end{equation}
The set of \ids of $\gT$ itself forms a \trdiz%
\footnote{Actually we need to introduce a restriction to truly obtain a set, in order to have a well-defined procedure to construct concerned \idsz.
For example we can consider the set of \ids obtained from \idps via certain predefined operations, 
such as countable unions and intersections. This is the same problem
 as the one indicated in footnote~\vref{NoteValVer}.} 
\wirt the inclusion, with, for lower bound of $\fa$ and $\fb$, the \id 
\begin{equation}\label{eqInfId}\preskip.5em \postskip.4em
\fa\cap \fb=\fa\vi\fb.
\end{equation}

Thus the operations $\vu$ and $\vi$ defined in (\ref{eqvuvi}) correspond to the supremum and the infimum in the lattice of \idsz.

 \rdb\label{NOTAfitrdi}
We will denote by
$\cF_\gT(S)=\sotq{\,x\in\gT}{\Ex S_0\in \Pfe(S),\,x\geq \Vi S_0}$
the filter of~$\gT$ generated by the subset $S$.
\\
When we consider the lattice of filters, 
we must pay attention as to what the reversing of the order relation produces: $\ff\cap\ffg=\ff\vu\ffg$ is the infimum of the filters~$\ff$ and $\ffg$, whereas their supremum is equal to $\cF_\gT(\ff\cup \ffg)=\ff\vi \ffg$.

\rdb
 The \emph{quotient lattice of $\gT$ by the \id $\fa$}, denoted by $\gT/(\fa=0)$, is defined as the \trdi generated by the \elts of $\gT$ with as its relations the true relations in $\gT$ on the one hand, and the relations $x=0$ for the $x\in \fa$ on the other. It can also be defined by the following preorder relation
$$ 
\label{trquoideal}\preskip.4em \postskip-.1em
a\leq_{\gT/(\fa=0)}b\;\;
\equidef \;\;
\exists x\in \fa \;\;a\leq  x\vu b.
$$
This gives
$$\preskip.1em \postskip.4em
 a\equiv b\;\;\mod\;(\fa=0)\quad \Longleftrightarrow  \quad \exists
x\in \fa
\;\;\;a\vu x=b\vu x.
$$
In particular, the \homo of passage to the quotient
$$\varphi:\gT\to\gT'=\gT/(\fa=0)$$
satisfies $\varphi^{-1}(0_{\gT'})=\fa$. In the case of the quotient by a principal \id $\dar a$ we obtain $\gT/(\dar a=0)\simeq\uar a$ with the morphism $y\mapsto y\vu a$ from $\gT$ to $\uar a$.

\begin{proposition}
\label{propIdealFiltre} Let $\gT$ be a \trdi and $(J,U)$ be a pair of subsets of $\gT$.
Consider the quotient $\gT'$ of $\gT$ defined by the relations $x=0$ for each $x\in J$, and $y=1$ for each $y\in U$. Then the in\egt $a\leq_{\gT'}b$ is satisfied \ssi
there exist $J_0\in\Pfe( J)$  and  $U_0\in\Pfe( U)$  such that
\begin{equation}\label{eqpropIdealFiltre}\preskip.4em \postskip.4em
a \vi \Vi U_0 \; \leq_\gT\; b \vu \Vu J_0.
\end{equation}
We will denote by $\gT/(J=0,U=1)$ this quotient lattice $\gT'.$
\end{proposition}

We see in the example of totally ordered sets that a quotient structure of a \trdi is not \gnlt
\caree by the equivalence classes of $0$ and $1$.

\rdb
Let $\fa$ be an \id and $\ff$ be a filter of $\gT$.
We say that $\fa$ is \emph{$\ff$-saturated} if we have
$$\preskip-.4em \postskip.4em
(g\in \ff,\; x\vi g \in \fa) \Longrightarrow x\in \fa,
$$
we say that $\ff$ is \emph{$\fa$-saturated} if we have
$$\preskip.3em \postskip.4em
(a\in \fa,\; x\vu a \in \ff) \Longrightarrow x\in \ff.
$$
If $\fa$ is $\ff$-saturated and $\ff$ is $\fa$-saturated we say that $(\fa,\ff)$
is a \emph{\paz} in~$\gT$. %
\index{saturated pair}%
\index{saturated!pair}%
\index{saturated!$\ff$- --- ideal }%
\index{saturated!$\fa$- --- filter}%
 When $(\fa,\ff)$  is a \paz, we have the \eqvcs
$$\preskip.3em \postskip.0em 
1\in \fa\; \;\Longleftrightarrow\;\;  0\in \ff
\;\; \Longleftrightarrow\;\;  (\fa,\ff)=(\gT ,\gT ). 
$$

\begin{fact}\label{factIdFiAsTrdi}
Let $\varphi:\gT \rightarrow \gT_1$ be a \trdi \homoz.
The \id $\fa=\varphi^{-1}(0)$ and the filter $\ff=\varphi^{-1}(1)$ form a \paz.
Conversely, if $(\fa,\ff)$ is a \pa of $\gT$,
the \homo of passage to the quotient $\pi:\gT \rightarrow \gT/(\fa=0,\ff=1)$ satisfies $\pi^{-1}(0)=\fa$
and $\pi^{-1}(1)=\ff$.
\end{fact}

\vspace{-.5em}
\subsec{\agBsz}

In a \trdi an \elt $x'$ is called a \emph{complement}
of $x$ if we have $x\vi x'=0$ and $x\vu x'=1$. If it exists the complement of $x$ is unique. It is then often denoted by  $\lnot x$.
\index{complement!(in a \trdiz)}

Recall that by \dfn a \ri $\gB$ is a \agB \ssi every \elt is idempotent. We then define an order relation $x\preceq y$ by: $x$ is a multiple of $y$, \cad $\gen{x}\subseteq\gen{y}$.\\
We thus obtain a \trdi in which every \elt $x$ admits as its complement $x'=1+x$
(cf.~Proposition~\ref{defiBoole}).

We have the following converse.

\begin{proposition}\label{defiBooleTrdi}
\emph{(Boolean algebras)}
\begin{enumerate}
\item On a \trdi in which every \elt $x$ admits a complement, denoted by $\lnot x$, we can define  a \agB structure by letting

\snic{xy=x\vi y\;\hbox{  and  }\;x\oplus y=(x\vi \lnot y)\vu (y\vi \lnot x).}

We once again find $x\vu y=x\oplus y\oplus xy$ and $\lnot x =1\oplus x$.
\item Every \homo of \trdis between two \agBs is a \homo of \agBsz, and it respects the passage to the complement.
\end{enumerate}
\end{proposition}

\entrenous{
De la m\^eme mani\`ere que in la Proposition~\ref{propBoolFini}, \emph{tout \trdi discret se comporte in les calculs comme un
produit finite d'ensembles totalement ordonn\'es finis}. Mais if l'on veut en parler il faut un \tho qui pr\'ecise cela,
ou au moins un commentaire
}

\subsec{\agB generated by a \trdiz}

Let us begin with a few remarks on the \elts that have a complement in a \trdiz.
If $a$ admits a complement $a'$, since $b=(b\vi a)\vu(b\vi a')$ for every $b\in\gT$,
the canonical \homo
$$
\gT\to \gT\sur{(a=1)}\times \gT\sur{(a'=1)}
$$
is injective. Moreover this morphism is onto because for
 $x,y\in \mathbf{T}$, defining $z=(x\wedge a)\vee(y\wedge a')$,
we get $z\wedge a =x\wedge a$, i.e. $z\equiv x \mod (a=1)$, and in the same way $z\equiv y \mod (a'=1)$.  Conversely, we have the following result which shows the similarity between an \idm in a commutative \ri and an \elt having a complement in a \trdi
(see Fact~\ref{lemCompAnnComm}).

\begin{lemma}\label{lemCompTrdi} 
For every \iso $\lambda:\gT\to\gT_1\times \gT_2$, there exists a (unique)
\elt $a\in\gT$ such that
\begin{enumerate}
\item $a$ has a complement $\lnot a$,
\item the composed \homo $\gT\to\gT_1$ identifies $\gT_1$ with $\gT\sur{(a=0)}$ and with $\gT\sur{(\lnot a=1)}$,
\item the composed \homo $\gT\to\gT_2$ identifies $\gT_2$ with $\gT\sur{(a=1)}$ and with $\gT\sur{(\lnot a=0)}$.
\end{enumerate} 
\end{lemma}

\begin{proof}
The \elt $a$ is given by $\lambda(a)=(0_{\gT_1},1_{\gT_2})$.
\end{proof}

When two \elts $a$ and $b$ have complements $\lnot a$ and $\lnot b$, the \ix{De Morgan's laws} are satisfied
\begin{equation}\preskip.2em 
\label{eqMorgan}
\lnot(a\vi b)=\lnot a \vu \lnot b \;\hbox{  and }\; \lnot(a\vu b)=\lnot a \vi \lnot b.
\end{equation}

By \dfnz, the \emph{\agB freely generated by the \trdi $\gT$}
is given by a pair $(\Bo(\gT),\lambda)$, where $\Bo(\gT)$ is a \agBz, and where $\lambda:\gT\to\Bo(\gT)$ is a \trdi \homo satisfying the following universal \prtz.
\\
\emph{Every \trdi \homo $\psi$ from $\gT$ to a \agB $\gB$ is uniquely factored in the form $\varphi\circ \lambda$.}

\vspace{-1.6em}
\PNV{\gT}{\lambda}{\psi}{\Bo(\gT)}{\varphi}{\gB}{\trdisz}{}{\agBsz}

\vspace{-1.8em}
Since we are in the context of purely equational \agq structures, this \agB can be constructed from 
 $\gT$ by forcefully adding a unary law $a\mapsto \lnot a$ and by imposing the axioms $a\vi \lnot a=0$, $a\vu \lnot a=1$.

In other words $\Bo(\gT)$ can be defined as a \agB obtained by \gtrs and relations. 
The \gtrs are the \elts of~$\gT$ and the relations are those that are true in $\gT$: of the form $a\vi b=c$ or~$a\vu b= d$, not to mention $0_{\Bo(\gT)}=0_\gT$ and $1_{\Bo(\gT)}=1_\gT$.

This description is however somewhat vague so we will construct $\Bo(\gT)$ at turtle speed to see things more clearly.

\begin{lemma}\label{lemRajouCompl}
Let $\gT$ be a \trdi and $a\in\gT$. Consider the \trdi
$$\preskip.0em \postskip.4em
\gT[a\bul]\eqdefi\gT\sur{(a=0)}\times \gT\sur{(a=1)}$$
and $\lambda_a:\gT\to\gT[a\bul]$ be the canonical \homoz. 
\begin{enumerate}
\item  The \homo $\lambda_a$ is injective and  $\lambda_a(a)=(0,1)$ admits $(1,0)$ as its complement.
\item
For every \homo $\psi :\gT\to\gT'$ such that $\psi(a)$ admits a complement,
there exists a unique \homo $\varphi:\gT[a\bul]\to \gT'$
such that~$\varphi\circ \lambda_a=\psi$.

\vspace{-1.2em}
\Pnv{\gT}{\lambda_a}{\psi}{\gT[a\bul]}{\varphi}{\gT'}{ }{}{$\psi(a)$ admits a complement}

\vspace{-1.4em}
\end{enumerate}
\end{lemma}

\begin{proof}
Lemma~\ref{lemCompTrdi} gives $\gT'\simeq \gT'\sur{(\psi(a)=0)}\times \gT'\sur{(\psi(a)=1)}$, hence
the \homo $\varphi$ and the uniqueness. The injectivity of $\lambda_a$ is not obvious but it is a grand classic:
if $x\vi a=y\vi a$ and $x\vu a=y\vu a$, then
$$\preskip.3em \postskip.3em 
x=(x\vu a)\vi x=(y\vu a)\vi x=(y\vi x)\vu(a\vi x). 
$$
Symmetrically
$y=(y\vi x)\vu(a\vi y)$, so $x=y$ since $a\vi x=a\vi y$.
\end{proof}

\begin{corollary}\label{corlemRajouCompl}
Let $a_1$, \ldots, $a_n\in\gT$.
\begin{enumerate}\itemsep1pt
\item The lattice $\gT[a_1\bul][a_2\bul]\cdots[a_n\bul]$ is independent, up to \isoz, in the order of the $a_i$'s.
It will be denoted by $\gT[a_1\bul,a_2\bul,\ldots,a_n\bul]$.
\item A possible description is the following
$$\preskip.2em \postskip.1em
\gT[a_1\bul,a_2\bul,\ldots,a_n\bul]
\;\simeq\;
\prod\nolimits_{I\in\cP_n}{\gT\sur{\bigl( (a_i=0)_{i\in I}, (a_j=1)_{j\in \lrbn\setminus I}\bigr)}}.
$$
\item The natural \homo $\gT\to\gT[a_1\bul,a_2\bul,\ldots,a_n\bul]$ is injective.
It uniquely factors every \homo $\psi$ from $\gT$ to a \trdi $\gT'$ such that the $\psi(a_i)$'s admit a complement.
\end{enumerate}
\end{corollary}

\begin{theorem}\label{thBoolGen}
\emph{(\agB freely generated by a \trdiz)}
For every \trdi $\gT$ there exists a \agBz, denoted by $\Bo(\gT)$, with a \homo $\lambda:\gT\to\Bo(\gT)$, which
uniquely factorizes every \homo $\psi :\gT\to\gB$ to a \agBz.
This pair $(\Bo(\gT),\lambda)$ is unique up to \isoz.
We have in addition the following \prtsz.
\vspace{-.1em}
\begin{enumerate}\itemsep0pt
\item [--] The \homo  $\lambda$ is injective.
\item [--] We have $\Bo(\gT)=\gT[(a\bul)_{a\in\gT}]$.
\end{enumerate}
\end{theorem}
\begin{proof}
It remains to see that the (filtering) colimit of $\gT[a_1\bul,a_2\bul,\ldots,a_n\bul]$ is indeed a \agBz. This results from De Morgan's laws.
\end{proof}

\exl Suppose that $\gT$ is a lattice of detachable subsets of a set $E$, in the sense that if $A$ and $B$ are \elts of $\gT$, then so are $A \cup B$ and $A\cap B$ (with in addition $\emptyset$ and $E$ as \elts of $\gT$). Then $\Bo(\gT)$ identifies with the set of finite Boolean combinations of \elts of $\gT$ and it is a \agB of subsets of $E$. \eoe

\smallskip 
\comm In \clamaz, every \trdi is \isoc to a lattice of subsets of a set. This provides an alternative \gui{construction} of the \agB $\Bo(\gT)$.\eoe

\section{Lattice-groups}
\label{secGpReticules}

\subsec{First steps}

In this book we limit ourselves, for the ordered groups, to the case of commutative groups.

\begin{definition}\label{defiGpReticule}
We call an \emph{ordered group} an Abelian group $G$ equipped with a partial order relation \emph{compatible} with the group law, \cad in additive notation,
$$\preskip-.2em \postskip.4em
\forall a,x,y\in G\qquad x\leq y\; \Longrightarrow\; a+x\leq a+y.
$$
An ordered group is called a \emph{lattice-group} when two arbitrary \elts admit a lower bound, which we will denote by $x\vi y$.
If \ncrz, we specify the structure by writing $(G,0,+,-,\vi)$.
A \emph{morphism of \grlsz} is a group \homo which respects the law $\vi$.%
\index{ordered group}\index{lattice-group}
\end{definition}

An Abelian group equipped with a compatible total order (we say {a \emph{totally ordered group}) is an \grlz.
The totally ordered group morphisms are then the non-decreasing group \homosz.\index{totally ordered!group}

An \emph{\sgrlz} of an \grl $G$ is by \dfn a stable subgroup for the lattice law $\vi$. For that it is not sufficient for the induced order relation on the subgroup to make a lattice of it.%
\index{lattice-subgroup} 

A guiding idea in the theory of \grls is that \emph{an \grl behaves in computations as a product of totally ordered groups}. This will be constructively translated by the \prf \vref{prcfgrl}.

\medskip
\exls
1) (Careful, multiplicative notation!) The set $\QQ^{>0}$ of strictly positive rationals is an \grl with as its  positive subset the \mo $(\NN^{>0},1,\times )$.
The example of this multiplicative structure is paradigmatic.
We have an \iso of \grls $\QQ^{>0}\simeq \ZZ^{(P)}$,\label{NotaZZP} where $P$ is the set of prime numbers, $\ZZ^{(P)}=\bigoplus_{p\in P}\ZZ$ and the order is induced by the product order. This is just another way to express the fundamental \tho of arithmetic \gui{every natural number is uniquely expressible as a product of powers of prime numbers.}
It is by wanting to make multiplication for integers of number fields look like multiplication in $\NN^{>0}$ at all costs that mathematicians have been brought to invent the \gui{\id gcd numbers.}

\rdb
2)
If $(G_i)_{i\in I}$ is a family of \grls with a discrete indexing set $I$, we define the \emph{\sdoz} of the family, denoted by $\boxplus_{i\in I}G_i$%
\label{NotaSDirOr}, which is an \grl with as subjacent group the group $\bigoplus_{i\in I}G_i$, the law $\vi$ being defined \cooz wise. If $I=\lrb{1..3}$ we will let $G_1\boxplus G_2\boxplus G_3 $.%
\index{orthogonal direct sum}
\\
For example $\ZZ^{(P)}=\boxplus_{p\in P}\ZZ$. 
\\
We \egmt define the product $\prod_{i\in I}G_i$ in the usual way, and it is the product in the category of \grlsz.
When $I$ is a finite set, the \grls  $\boxplus_{i\in I}G_i$
and $\prod_{i\in I}G_i$ are naturally \isocz.

3) 
If $(G_i)_{i\in I}$ is a family of totally ordered discrete groups with for $I$ a totally ordered discrete set we define the {lexicographic sum} of this family, it is the totally ordered discrete group~$G$ whose subjacent group is $\bigoplus_{i\in I}G_i$ and the order relation is the lexicographical order: $(x_i)_{i\in I}< (y_i)_{i\in I}$ \ssi $x_{i_0}< y_{i_0}$ for the smallest index $i_0$ such that $x_{i_0}\neq y_{i_0}.$
\eoe

\medskip
In an \grl the translations are \autos of the order structure, 
hence the \dit rule
\begin{equation}\preskip.4em \postskip.4em
\label{eq1grl}
{x+(a\vi b)\,=\,(x+a)\vi(x+b).}
\end{equation}
We also see that the bijection $x\mapsto-x$ reverses the order, and thus that two arbitrary \elts $x$, $y$ also admits an upper bound
$$
\preskip.4em \postskip.4em 
x\vu y=-\big((-x)\vi(-y)\big), 
$$

\vspace{-.7em}
\pagebreak	

with $x+y-(x\vu y)=(x+y)+\big((-x)\vi(-y)\big)=(x+y-x)\vi (x+y-y)$, so
\begin{eqnarray}
\label{eq2grl}
x+y&=&(x\vi y)+(x\vu y),\\
x+(a\vu b)&=&(x+a)\vu(x+b).\label{eq3grl}
\end{eqnarray}

However, a minimum \elt and a maximum \elt are missing to obtain a lattice.

\subsec{Remarkable identities in the \grlsz}

\Grandcadre{In this subsection $G$ is an \grl
and $G^+$ is the sub\mo  \\ of $G$ formed from the 
                        non-negative \eltsz.}

Let
$
x^+=x\vu 0, \, x^-=(-x)\vu 0  \hbox{ and }  \abs{x}=x\vu(-x). 
$
We respectively call them the \ix{positive part}, the \ix{negative part} 
and the \ix{absolute value} of $x$.

\begin{theorem}\label{th0GpRtcl} \emph{(Distributivity in the \grlsz)}\\
In an \grl the laws $\vi$ and $\vu$ are distributive with respect to one another.
\end{theorem}
\begin{proof}
It suffices to show $x\vu(y_1\vi y_2)= (x\vu y_1)\vi (x\vu y_2)$. By translating \hbox{by $-x$}, we are reduced to $x = 0$, \cad to $(y_1 \vi y_2)^+ = y_1^+ \vi y_2^+$.  
\\
The in\egt $(y_1 \vi y_2)^+ \le y_1^+ \vi y_2^+$ is \imdez.  
\\
Let $y=y_1\vi y_2$.
The \elt $y_i + y^+ - y$ is  $\ge y_i$ and $\geq 0$, so $\geq y_i^+$.
\\
Hence
$y_i^+ + y \le y_i + y^+$. Then
$(y_1^+ + y) \vi (y_2^+ + y) \le (y_1 + y^+) \vi (y_2 + y^+)$, \cad $(y_1^+ \vi
y_2^+) + y \le (y_1 \vi y_2) + y^+$, \cad $y_1^+ \vi y_2^+ \le y^+$.
\end{proof}

Two \elts $x$, $y$ are said to be \emph{disjoint} or \emph{\ortz} if $\abs{x}\vi\abs{y}=0$.

\begin{lemma}\label{lemx+x-}\relax Let $x$, $y\in G$.
\vspace{-1mm}
\begin{eqnarray}
\label{i1lemx+x-}\relax
&x=x^+ - x^-,\quad x^+ \perp x^-,\quad \abs{x}=x^++x^-=x^+\vu x^-\in G^+& \\
\label{i2lemx+x-}\relax
&x\leq y \; \iff\; x^+\leq y^+ \et y^-\leq x^-,\quad x=0\;\iff\;\abs{x}=0&%
\end{eqnarray}
\end{lemma}
\begin{proof} ({\ref{i1lemx+x-}}).
First of all $x^+-x=(x\vu 0)-x=(x-x)\vu\big(0+ (-x)\big)=x^-$. \\
Still by \dit we obtain

\snac{x^+ + x^-= (x\vu 0)+((-x)\vu 0)= (x-x)\vu (x+ 0)\vu\big(0+(-x)\big)\vu(0+ 0)=x^+\vu x^- .}

 Finally, since $x^+ + x^-=(x^+\vu x^-)+(x^+\vi x^-)$, this gives $x^+\vi x^-=0$.

 ({\ref{i2lemx+x-}}).
Left to the reader.
\end{proof}
%

\begin{lemma}\label{lemGlem} \emph{(Gauss' lemma)}
Let $x$, $y$, $z\in G^+.$
\vspace{-1mm}
\begin{eqnarray}
\label{i1lemGlem}
(x\perp y \et x\leq y+z) &\;\Longrightarrow\;& x\leq z \\
\label{i2lemGlem}
x\perp y &\Longrightarrow& x\vi (y+z) = x\vi z\\
\label{i3lemGlem}
(x\perp y \et x\perp z) &\Longrightarrow& x\perp (y+z) \\
\label{i4lemGlem}
(x\perp y \et x\leq z \et y\leq z) &\Longrightarrow& x+y\leq z
\end{eqnarray}
\end{lemma}
\begin{proof}
{(\ref{i1lemGlem}).} We have $x\leq z+x$ because $z\geq0$ and $x\leq z+y$ by hypothesis, 
therefore $x\leq (z+x)\vi(z+y)=z+(x\vi y)=z$.

{(\ref{i2lemGlem}).}
Let $x' = x\vi (y+z)$. It suffices to see that $x' \le x\vi z$.
We have $x' \geq 0$, $x' \le x$ so $x' \perp y$. We can apply the previous item to the in\egt $x' \le y+z$: it provides $x' \le z$, as desired.

{(\ref{i3lemGlem}).}
Direct consequence of the previous item.

{(\ref{i4lemGlem}).}
 Because $x+y=x\vu y$ and $x\vu y\leq z$.
\end{proof}
%

\begin{corollary}\label{corlemGlem} Let $x$, $y$, $z\in G,\, n\in\NN\etl$.
\vspace{-1mm}
\begin{eqnarray}
\label{i1corlemGlem}
&(x=y - z\hbox{, }y\geq 0\hbox{, }z\geq 0 \hbox{, and } y\perp z) \,\Longleftrightarrow\, (y=x^+ \hbox{ and } z=x^-) &\;\;\;\;\;\;\;\\
\label{i11corlemGlem}
&(x\geq 0,\;y\geq 0, \hbox{ and }x\perp y)\;\Longrightarrow\;x\perp ny & \\
\label{i12corlemGlem}
& (nx)^+=nx^+,\; (nx)^-=nx^- ,\; \abs{nx}=n\abs{x}& \\
\label{i2corlemGlem}
&\;nx=0\;\Longrightarrow\;x=0 & \\
\label{i3corlemGlem}
&n(x\vi y)=nx\vi ny,\; \; \;n(x\vu y)=nx\vu ny &
\end{eqnarray}
\end{corollary}
\begin{proof}
{(\ref{i1corlemGlem}).} 
It remains to show $\Longrightarrow$. We have $x^++z=x^-+y$.
By applying Gauss' lemma, we obtain $y\leq x^+$ (because $y\perp z$)
and $x^+\leq y$ (because $x^+\perp x^-$).

{(\ref{i11corlemGlem}).} 
 Results from {(\ref{i2corlemGlem}).} 

{(\ref{i12corlemGlem}).} 
 By (\ref{i1corlemGlem}) and (\ref{i11corlemGlem})
since $nx=nx^+-nx^-$ and $nx^+\perp nx^-$.

{(\ref{i2corlemGlem}).} 
 By (\ref{i12corlemGlem}) since the implication is true for $x\geq0$.

{(\ref{i3corlemGlem}).} 
 The \elts $b=x\vu y$, $a=x\vi y$, $x_1=x-a$ and $y_1=y-a$
are \cares by the following relations
$$\preskip.4em \postskip.4em 
x_1\geq 0,\;y_1\geq 0,\;x=x_1+a,\;y=y_1+a,\;x_1\perp y_1,\;a+b=x+y. 
$$
We multiply everything by $n$.
\end{proof}

\subsec{Simultaneous congruences, covering principle by quotients}

\begin{definition}\label{defiCongru}
If $a\in G$, we define  \ixc{congruence modulo $a$}{in an \grlz} as follows

\snic{x\equiv y\; \mod\; a \equidef \exists n\in\NN^*,\;\abs{x-y}\leq n\,\abs{a}.}

We denote by $\cC(a)$ the set of $x$'s congruent to $0$ modulo $a$.
\end{definition}

\begin{fact}\label{factCongru}
The set $\cC(a)$ is an \sgrl of $G$ and the lattice laws pass to the quotient in $G/\cC(a)$. 
\\
Thus, the canonical map $\pi_a:G\to G/\cC(a)$ is a morphism of \grlsz,
and every \grl morphism $G\to G'$ which annihilates~$a$ is factorized by~$\pi_a$.
\end{fact}

The meaning of the congruence $x\equiv 0 \mod a$ is therefore that every \grl morphism $G\vers{\varphi} G'$ that annihilates $a$ annihilates~$x$.%
\footnote{In fact, by direct computation, if $\varphi(a)=0$, then $\varphi(\abs{a})=\abs{\varphi(a)}=0$, and $\abs{\varphi(x)}=\varphi(\abs{x})\leq\varphi(n\abs{a})=n \varphi(\abs{a})=0$, so~$\varphi(x)=0$.}

\medskip The following lemma has an \ari Chinese remainder \tho flavor (see \thref{thAnar} item \emph{5}) for the \grlsz, but only a flavor. It is distinctly simpler.

\begin{lemma}\label{lemChinoisGRL} \emph{(Lemma of simultaneous congruences)}\\
Let $(\xn)$ in $G^+$ and $(\an)$ in $G$.
\begin{enumerate}
\item If the in\egts
$ \abs{a_i - a_j} \leq x_i+x_j,\;i,j\in\lrbn,
$
 are satisfied there exists some $a\in G$ such that $ \abs{a - a_i} \leq x_i,\;i\in\lrbn
$. Moreover
\begin{itemize}
\item If the $a_i$'s are in $G^+$ we have a solution $a$ in $G^+$.
\item If $\Vi_ix_i=0$, the solution $a$ is unique.
\end{itemize}
\item Similarly, if $a_i\equiv a_j \mod x_i+x_j$ for $i$, $j\in\lrbn$,
there \hbox{exists an $a\in G$} such that $a \equiv a_i \mod x_i,\;i\in\lrbn$. Moreover
\begin{itemize}
\item If the $a_i$'s are in $G^+$ we have a solution $a$ in $G^+$.
\item If $\Vi_ix_i=0$, the solution $a$ is unique.
\end{itemize}
\end{enumerate}
\end{lemma}
\begin{proof} It suffices to prove item \emph{1}. Let us first take a look at uniqueness.
If we have two solutions $a$ and $a'$ we will have $\abs{a-a'}\leq2x_i$ for each $i$,
so $\abs{a-a'}\leq2\Vi_ix_i$.\\
Let us move on to existence. We treat the case where the $a_i$'s are in $G^+$.
This is actually a matter of showing that the hypotheses imply the in\egt $\Vu_i(a_i-x_i)^+\leq \Vi_i(a_i+x_i)$.
It suffices  to verify that for each $i$, $j$, we have $(a_i-x_i)\vu 0 \leq a_j+x_j$. However, $0\leq a_j+x_j$, and $a_i-x_i \leq a_j+x_j$ by hypothesis.
\end{proof}
%

\begin{lemma}\label{lemRecfermebasique}
Given a finite family $(a_j)_{j\in J}$ in an \grl $G$ and a finite subset $P$ of $J\times J$, there exists a finite family $(x_i)_{i\in I}$ in~$G$ such that
\begin{enumerate}
\item $\Vi_{i\in I}x_i=0$.
\item Modulo each of the $x_i$'s, for each $(j,k)\in P$, we have $a_j\leq a_k$ or~$a_k\leq a_j$.
\end{enumerate}
\end{lemma}
%
\begin{proof}
Let $y_{j,k}=a_j-(a_j\vi a_k)$ and $z_{j,k}=a_k-(a_j\vi a_k)$.
We have $y_{j,k}\vi z_{j,k}=0$. Modulo $y_{j,k}$, we have $a_j=a_j\vi a_k$,
\cad $a_j\leq a_k$, and modulo $z_{j,k}$, we~have~$a_k\leq a_j$. 
\\
By expanding by \dit the sum $0=\sum_{(j,k)\in P}(y_{j,k}\vi z_{j,k})$ we obtain some $\,\Vi_{i\in I}x_i$, where each $x_i$ is a sum $\sum_{j,k}t_{j,k}$, with one of the two \elts $y_{j,k}$ or $z_{j,k}$ as~$t_{j,k}$. Modulo such a $x_i$ each of the $t_{j,k}$'s is null (because they are $\geq 0$ and their sum is null). We are therefore indeed in the stated situation.
\end{proof}
%


The next principle is a kind of analogue, for \grlsz, of the basic \plg for  commutative \risz.\iplg

Actually this is a simple special case of item \emph{2} of Lemma~\ref{lemChinoisGRL} when the $a_i$'s are all zeros: we apply uniqueness.

\pagebreak

\begin{prvq}\label{prcfgrl} \emph{(For  \grlsz)}\\
Let $a$, $b\in G$, $x_1$, \ldots, $x_n\in G^+$
with $\Vi_ix_i=0$. Then  $a\equiv b \mod x_i$ for each~$i$
\ssi $a=b$.
\\
Consequently, given Lemma~\ref{lemRecfermebasique}, to demonstrate an \egt $a=b$ we can always suppose that the (finite number of) \elts which occur in a computation for a \dem of the \egt are comparable, if we need it to do the \demz.
The principle applies just as well for  in\egts as for  \egts since $a\leq b$ is equivalent to $a\vi b=a$.
\end{prvq}

\rem In slightly more abstract terms, we could have said that the canonical \grl morphism $G\to\prod_i G\sur{\cC(x_i)} $ is injective,
and the comment that concludes the covering principle by quotients can be paraphrased as follows: in  computations, an \grl always behaves like a product of totally ordered groups.
\eoe

\medskip
In the Riesz \tho that follows we will note that the \gui{there exists} are abbreviations for explicit formulas which result from the \demz.
Thus the \tho can be seen as a family of \idas in $G$, under certain sign conditions (which are in the hypothesis).
It is also possible to regard this \tho as a family of \gui{pure} \idas in $G^+$, \cad without any sign condition.
In this case $G^+$ must be considered as an \agq structure for which we add the operation $x \dotdiv y\eqdefi x-(x\vi y)$ (well-defined over $G^+$ despite the fact that it calls upon the $-$ operation of $G$).

\begin{theorem}\label{th1GpRtcl} \emph{(Riesz \thoz)}\\
Let $G$ be an \grl and $u$, $x_1$, \dots, $x_n$, $y_1$, \dots, $y_m$ in $G^+$.
\begin{enumerate}
\item \label{i1th1GpRtcl}
If $u\leq \sum_jy_j$, there exist $u_1$, \ldots, $u_m\in G^+$ such that $u_j\leq y_j$
for $j\in\lrbm$ and $u=\sum_ju_j$.
\item \label{i2th1GpRtcl}
If $\sum_ix_i=\sum_jy_j$, there exists $(z_{i,j})_{i\in\lrbn,j\in \lrbm} $ in $G^+$ such that for all $i$, $j$ we have $ \sum_{k=1}^m z_{i,k}=x_i$ and $\sum_{\ell=1}^n z_{\ell,j}=y_j$.
\end{enumerate}
\end{theorem}
\begin{Proof}{Direct \demz, but clever. }\\
\emph{\ref{i1th1GpRtcl}.} It suffices to prove it for $m=2$ (easy \recu on $m$). If $u\leq y_1+y_2$, 
let 
$u_1=u\vi y_1$ and $u_2=u-u_1$. We need to prove $0\leq u_2\leq y_2$.  
However, $u_2=
u-(u\vi y_1)=u+\big((-u)\vu (-y_1)\big)=(u-u)\vu (u-y_1)\leq y_2$.

 \emph{\ref{i2th1GpRtcl}.}
For $n=1$ or $m=1$ there is nothing to do.
For $n=2$, it is given by item~\emph{\ref{i1th1GpRtcl}.}
Therefore let us suppose $n\geq3$. Let $z_{1,1}=x_1\vi y_1$, $x'_1=x_1-z_{1,1}$
and~$y'_1=y_1-z_{1,1}$. We have $x'_1+x_2+\cdots+x_n=y'_1+y_2+\cdots+y_m $.
\\
Since $x'_1\vi y'_1=0$, Gauss' lemma gives $x'_1\leq y_2+\cdots+y_m$.
\\
By item~\emph{\ref{i1th1GpRtcl}} we can write $x'_1=z_{1,2}+\cdots+z_{1,m}$
with each $z_{1,j}\leq y_j$, \cad $y_j=z_{1,j}+y'_j$ and $y'_j\in G^+$.
Therefore $x_2+\cdots+x_n=y'_1+y'_2+\cdots+y'_m$. 
\\
This therefore allows us to perform an \recu on~$n$.

\emph{\Demo by the covering principle by quotients.}\\
It suffices to prove item \emph{2.}
By applying the principle~\ref{prcfgrl}, we can assume that the group is totally ordered.
Suppose for example $x_1\leq y_1$. Let $z_{1,1}=x_1$, $z_{1,k}=0$ for $k\geq2$.
We replace $y_1$ with $y_1-x_1=y'_1$. We are reduced to solving the \pb for $x_2$, \ldots, $x_n$ and~$y'_1$, $y_2$, \ldots, $y_m$. Gradually, we thus decrease $n+m$ until  $n=1$ \hbox{or $m=1$}, in which case everything is clear.
\end{Proof}
%

\begin{fact}\label{factGpRtcl} \emph{(Other identities in the \grlsz)}\\
Let $x$, $y$, $x'$, $y'$, $z$, $t\in G$, $n\in\NN$, $x_1$, \dots, $x_n\in G$.

\vspace{-.10em} 
\begin{enumerate}
\item \label{i1factGpRtcl} $x+y =\abs{x-y} +2(x\vi y)$.
\item $(x\vi y)^+=x^+\vi y^+$, $(x\vi y)^-=x^-\vu y^-$, \\ $(x\vu y)^+=x^+\vu y^+$, $(x\vu y)^-=x^-\vi y^-$.
\item $2(x \vi y)^+ \leq (x+y)^+ \leq x^++y^+$.
\item $\abs{x+y} \leq \abs{x}+\abs{y}\;:\;$
$\abs{x}+\abs{y}=\abs{x+y}+2(x^+\vi y^-) +2( x^-\vi y^+)$.
\item $\abs{x-y} \leq \abs{x}+\abs{y}\;:\;$
$\abs{x}+\abs{y}=\abs{x-y}+2(x^+\vi y^+) +2( x^-\vi y^-)$.
\item $\abs{x+y}\vu\abs{x-y}=\abs{x}+\abs{y}$.
\item $\abs{x+y}\vi\abs{x-y}=\abS{\abs{x}-\abs{y}}$.

\item $\abs{x-y}=(x\vu y)-(x\vi y)$.
\item $\abs{(x\vu z)-(y\vu z)}+\abs{(x\vi z)-(y\vi z)}= \abs{x-y}.$
\item $\abs{x^+ - y^+} + \abs{x^- - y^-} = \abs{x-y}$.
\item \label{i11factGpRtcl} $x\leq z \;\Longrightarrow\; (x\vi y)\vu z= x\vi (y\vu z)$.
\item $x+y=z+t \;\Longrightarrow\; x+y=(x\vu z)+(y\vi t)$.
\item \label{i13factGpRtcl} $n\, x\geq \Vi_{k=1}^n (k y+(n-k)x)  \,\Longrightarrow\,x\geq y$.
\item $\Vu_{i=1}^nx_i =  \sum_{k=1}^n(-1)^{k-1}
        \big(\sum_{I\in \cP_{k,n}}\Vi_{i\in I}x_i\big)$.
\item  $x\perp y\,\Longleftrightarrow\, \abs{x+y}=\abs{x-y}
\,\Longleftrightarrow\, \abs{x+y}=\abs{x}\vu \abs{y}$.
\item  $x\perp y\,\Longrightarrow\, \abs{x+y}=\abs{x}
+\abs{y}=\abs{x}\vu \abs{y}$.
\item \label{i15bisfactGpRtcl} $(x\perp y,\,x'\perp y,\,x\perp y',\,x'\perp y',\,x+y=x'+y') \,\Longrightarrow\, (x=x', \, y=y')$.
\item  \label{i16bisfactGpRtcl} We define $\Tri(\ux)=[\Tri_1(\ux), \Tri_2(\ux),\ldots,\Tri_n(\ux)]$, where 

\snic{\Tri_k(\xn) =\Vi_{I\in \cP_{k,n}}\left(\Vu_{i\in I}x_i\right) \quad(k\in\lrbn).}

We have the following results.
\begin{enumerate}
\item $\Tri_k(\xn) =\Vu_{J\in \cP_{n-k+1,n}}\big(\Vi_{j\in J}x_j\big)$, $(k\in\lrbn)$.
\item $\Tri_1(\ux)\leq \Tri_2(\ux)\leq \cdots \leq \Tri_n(\ux).$
\item If the $x_i$'s are pairwise comparable, the list $\Tr(\ux)$ is the list of the~$x_i$'s non-decreasingly ordered (it is not \ncr that the group be discrete).
\end{enumerate}
\end{enumerate}
Suppose $u$, $v$, $w\in G^+$.
\vspace{-.25em} 
\begin{enumerate}\setcounter{enumi}{18}\itemsep=1pt
\item \label{i17factGpRtcl} $u\perp v\,\Longleftrightarrow\, u+v=\abs{u-v}$.
\item $(u+v)\vi w \leq (u\vi w)+(v\vi w)$.
\item $(x+y)\vu w \leq (x\vu w)+(y\vu w)$.
\item $v\perp w  \,\Longrightarrow\,(u+v)\vi w = u\vi w$.
\item $u\perp v  \,\Longrightarrow\,(u+v)\vi w = (u\vi w)+(v\vi w)$.
\end{enumerate}
\end{fact}
\begin{proof}
All of this is just about \imd in a totally ordered group, by reasoning case-by-case. \Trf by the principle~\ref{prcfgrl}.
\end{proof}

\rems
\\
1) An implication like, for instance,

\snic{( u\vi v=0,\ u\geq0,\ v\geq0)\,\Longrightarrow\, u+v=\abs{u-v}}

(see item~\emph{\ref{i17factGpRtcl}})
can be seen as the result of an \idt which expresses, for a certain integer $n$, that $n\abs{u+v-\abs{u-v}}$ is equal to an expression which combines $u^-$, $v^-$ and $\abs{u\vi v}$ by means of the laws $\vu$, $\vi$ and $+$.
Actually, the \egt given in item~\emph{\ref{i1factGpRtcl}} directly settles the question without a sign hypothesis on $u$ and $v$: $\abs{u+v-\abs{u-v}}=2\abs{u\vi v}$.

2) There is an important difference between the usual \idasz, which are ultimately \egts between \pols in a free commutative \ri over \idtrsz, $\ZZ[\Xn]$, and the \idas in the \grlsz. The latter are certainly \egts between expressions that we can write in an \grl freely generated by a finite number of \idtrsz, but the structure of such a free \grl is distinctly more difficult to decrypt than that of a \pol \riz, in which the objects have a normalized expression. The comparison of two expressions in $\ZZ[\Xn]$ is \gui{easy} in so far as we bring each of them to normal form. The task is much more difficult in the free \grlsz, for which there is no unique normal form (we can reduce every expression to a supremum of infima of \colis of \idtrsz, but there is no uniqueness).
\eoe

\subsec{Partial \dcnz, \dcncz}

\begin{definition}\label{defiDecPar}
Let $(a_i)_{i\in I}$ be a finite family of non-negative \elts in a discrete \grlz~$G$.
\begin{enumerate}
\item We say that this family admits a \emph{\dcnpz} if we can find a finite family $(p_j)_{j\in J}$ of pairwise \ort non-negative \elts such that each $a_i$ is of the form $\sum_{j\in J}r_{i,j}p_j$ with all the $r_{i,j}\in\NN$.
The family $(p_j)_{j\in J}$ is then called a \emph{\bdpz} for the family~$(a_i)_{i\in I}$.
\item Such a \dcnp is called a \emph{\dcncz} if the $p_j$'s are \emph{\irdsz} (an \elt $q>0$ is said to be \ird if an \egt $q=c+d$ in $G^+$ implies $c=0$ or~$d=0$).
\item We say that an \grl  \emph{admits \dcnpsz} if it is discrete and if every finite family of non-negative \elts admits a \dcnpz.
\item We say that an \grl  \emph{admits \dcncsz} if it is discrete and if every non-negative \elt admits a \dcncz.
\item We say that an \grl  \emph{admits \dcnbsz} when for  all $x\geq 0$ there exists an integer $n$ such that, when $x=\sum_{j=1}^ny_j$ with each $y_j\geq 0$, 
at least one of the $y_j$'s is zero.
\item  An \grl is said to be \emph{\noez} if every non-increasing sequence of non-negative \elts admits two equal consecutive terms.
\end{enumerate}%
\index{decomposition!partial ---}%
\index{decomposition!complete ---}%
\index{decomposition!bounded ---}%
\index{irreducible!\elt in an \grlz}%
\index{noetherian!lattice-group}%
\index{lattice-group!admitting partial decompositions}%
\index{lattice-group!admitting bounded decompositions}%
\index{lattice-group!admitting complete decompositions}%
\end{definition}

\exls ~\\
An empty family, or a 
          family of null \eltsz,
 admits the empty family as a \bdpz.
 \\
The \grl $\ZZ^{(\NN)}$  admits \dcncsz.
\\
The \grls $\QQ^n$ ($n\geq 1$)  admit partial but not complete \dcnsz.
 \\ 
The \grl $\QQ[\sqrt 2]$ does not admit \dcnps (consider the finite family $(1,\sqrt2)$). \\
The lexicographical product $\ZZ\times \ZZ$ does not admit \dcnpsz. 
 \\
More \gnlt a totally ordered group admitting \dcnps is \isoc to a subgroup of $\QQ$.
\eoe

\medskip
It is clear that an \grl admitting \dcncs admits \dcnbs and that an \grl admitting \dcnbs is \noez.

In an \grl admitting \dcnpsz, two \dcnps for two finite families of $G^+$ admit a common refinement for the union of two families: here we mean that a \bdp $(q_1,\ldots,q_s)$ refines another if it is a \bdp for this other.

\begin{proposition}\label{propRaffinementCommun}
In an \grlz, if an \elt $>0$ admits a \dcncz, it is unique up to the order of the factors.
\end{proposition}
\begin{proof}
It suffices to show that if an \ird \elt $q$ is bounded above by a sum $\sum_ip_i$ of \ird \elts it is equal to one of them.
\\
However, we then have $q=q\vi \sum_ip_i$, and since $q\vi p_j=0$ or $p_j$, we can conclude with Gauss' lemma (\egrf{i2lemGlem}).
\\  
Note that we do not need to assume that the group is discrete.
\end{proof}
%

\begin{proposition}\label{factgrldcntot}
Let $G$ be an \grl admitting \dcncsz.
\begin{enumerate}
\item The \ird \elts of $G^+$ form a detachable subset $P$, and~$G$ is \isoc to the \ort direct sum $\ZZ^{(P)}$. 
\item The group $G$ admits \dcnbs (and a fortiori is \noez).
%
\end{enumerate}
\end{proposition}

\begin{proof} \emph{1.}
The irreducibility test is given by the \dcnc of the \elt to be tested.  
The isomorphism is obtained from the uniqueness of the \dcnc (up to the order of the factors).

\emph{2.} Let $x\in G^+$. Let us write $x=\sum_{j\in J}n_jp_j$ with the \ird $p_j$'s and $n_j\in\NN$, and let $n=\sum_jn_j$.
Then if $x=\sum_{k=1}^{n+1}x_k$ with non-negative $x_k$'s, one $x_k$ is \ncrt zero (consider the \dcn of each $x_k$ as a sum of \ird \eltsz).
\end{proof}

In \clamaz, a discrete \noe \grl admits \dcncsz. This result cannot be obtained \cotz. Nevertheless we obtain a \dcnpz.

\begin{theorem}\label{th2GpRtcl}
\emph{(Partial \dcn under \noe condition)}\\
A discrete and \noe \grl $G$ admits \dcnpsz.
\end{theorem}

For the \demz, we will use the following lemma.

\begin{lemma}\label{decomp2} (under the hypotheses of \thref{th2GpRtcl})\\
For $a \in G^+$ and $p_1$, \ldots, $p_m>0$  pairwise \ortsz, we can find pairwise \ort \elts $a_0$, $a_1$, \ldots, $a_m$ in $ G^+ $ satisfying the following \prtsz.
\begin{itemize}
\item [{1.}] $ a = \sum_{i=0}^{m}a_{i}$.
\item [{2.}] For all $i \in \lrbm $, there exists an integer $n_i \geq 0$ such that $a_{i} \leq n_{i}p_i$.
\item [{3.}] For all $ i \in \lrbm$,  $ a_0 \vi  p_i = 0$.
\end{itemize}
\end{lemma}
\begin{proof}
For each $i$, we consider the non-decreasing sequence $(a\vi np_i)_{n\in\NN}$ bounded above by~$a$.
There exists an $n_i$ such that $a\vi n_ip_i =a\vi (n_{i}+1)p_i.$
We then take $a_i= a\vi n_ip_i$. If~$a=a_i+b_i$, we have $b_i\vi p_i=0$ because $a_i\leq a_i+(b_i\vi p_i)\leq a\vi (n_{i}+1)p_i=a_i$.
The~$a_i$'s are $\leq a$, pairwise \orts and $\geq0$ so $a\geq\Vu_ia_i=\sum_ia_i$.
Thus, we write in $G^+$ $a = a_1 + \cdots + a_n + a_0$, with $a_{i} \leq n_{i}p_i$
for $i \in \lrbm$.
Finally, we have $b_i = a_0 + \sum_{j \ne i} a_j$, so $a_0 \le b_i$, then $a_0 \vi p_i \le b_i \vi p_i = 0$. As $a_i \le n_ip_i$, we a fortiori have $a_0 \vi a_i = 0$.
\end{proof}
\begin{Proof}{\Demo of \thref{th2GpRtcl}. }\\
By \recu on the number $m$ of \elts of the family. \\
$\bullet $ Suppose $m = 2$, consider the \elts $x_{1}$, $x_{2}$.
For ease of notation, let us call them $a$ and $b$.
Let $L_1 = \left [a,b \right]$, $m_1=1$, $E_{1,a}= [1,0]$, and~$E_{1,b}= [0,1]$.
The \algo proceeds in steps, at the beginning of step $k$ we have a natural integer $m_k$ and three lists of equal length: $L_k$, a list of non-negative \elts of $G$, $E_{k,a}$ and $E_{k,b}$, two lists of natural integers.
At the end of the step the integer $m_k$ and the three lists are replaced with a new integer and new lists, which are used at the next step (unless the \algo terminates).
The \gnl idea is the following: if $x$, $y$ are two consecutive non-\ort terms of  $L_k$, we replace in $L_k$ the 
segment $(x,y)$ with the segment  $(x-(x \vi y), x \vi y, y-(x \vi y))$ (by omitting the first and/or the last term if it is null).
We will denote this procedure as follows:

\snic{R:(x,y)\mapsto\;$
the new segment (of length 1, 2 or $3).}

Note that $x+y>\big(x-(x \vi y)\big)+ x \vi y+ \big(y-(x \vi y)\big)$.
\\
We have to define a loop-invariant. More precisely the conditions satisfied by the integer $m_k$ and the three lists are the following:
\begin{itemize}
\item $a$ is equal to the \coli of \elts of $L_k$
with \coes given by  $E_{k,a}$,
\item $b$ is equal to the \coli of \elts of $L_k$ with \coes given by $E_{k,b}$,
\item if $L_k=[x_{k,1},\ldots ,x_{k,r_k}]$ the \elts $x_{k,j}$ and $x_{k,\ell}$  are \orts as soon as
\begin{itemize}
\item  $j<m_k$ and $\ell\neq j$  or
\item  $j\geq m_k$ and $\ell\geq j+2$
\end{itemize}
\end{itemize}
In short, the $x_{k,j}$'s are pairwise \ortsz, except perhaps for certain 
pairs $(x_{k,j},x_{k,j+1})$ with $j\geq m_k$.
These conditions constitute \emph{the loop-invariant}.
It is clear that they are (trivially) satisfied at the start. 
\\
The \algo terminates at step $k$ if the \elts of $L_k$ are pairwise \ortsz.
In addition, if the \algo does not terminate at step $k$, we have the in\egt $\sum_{x\in L_k}x \,>\, \sum_{z\in L_{k+1}}z$, therefore the decreasing chain condition assures the termination of the \algoz.\\
It remains to explain the development of a step and to verify 
the loop-invariant.
In order to not manipulate too many indices, we make a slight abuse of notation and  write
$L_k = \left [p_1,\ldots, p_n \right]$, 
$E_{k,a} =\left [\alpha_1,\ldots, \alpha_n \right]$ 
and~$E_{k,b} = \left [\beta_1,\ldots, \beta_n \right]$.\\
The segment $(x,y)$  of $L_k$ which is treated by the procedure $R(x,y)$ is the following:
we consider the smallest index $j$ (\ncrt $\geq m_k$) such that $p_j\vi p_{j+1}\neq 0$ and we take $(x,y)=(p_{j},p_{j+1})$. If such an index does not exist, the \elts of~$L_k$ are pairwise \orts and the \algo is terminated.
Otherwise we apply the procedure  $R(x,y)$ and we update the integer (we can take $m_{k+1}=j$) and the three lists. 
\\
For example by letting $q_j=p_j\vi p_{j+1}$,
$p'_j=p_j-q_j$ and $p'_{j+1}=p_{j+1}-q_j$, if~$p'_j\neq 0\neq
p'_{j+1}$, we will have
$$\preskip.2em \postskip.4em\arraycolsep2pt
\begin{array}{rcl}
L_{k+1}& =  &  \left [p_1,\ldots,p_{j-1},p'_j,q_j,p'_{j+1},p_{j+2},\ldots ,
p_n \right]  , \\[.2em]
E_{k+1,a}& =  &  \left [\alpha_1,\ldots, \alpha_{j-1},\alpha_j,\alpha _j+\alpha
_{j+1},\alpha _{j+1},\alpha_{j+2},\ldots \alpha_n \right]  , \\[.2em]
E_{k+1,b}& =  &  \left [\beta_1,\ldots, \beta_{j-1},\beta_j,\beta _j+
\beta _{j+1},\beta _{j+1},\beta_{j+2},\ldots \beta_n \right] .
\end{array}
$$
We verify without difficulty in each of the four possible cases that the loop-invariant is preserved.

$\bullet $  If $m > 2$,  by \hdrz, we have for $(x_1, \ldots, x_{m-1})$ a \bdp  $(p_1,\ldots ,p_n)$.
By applying Lemma~\ref{decomp2} with $x_m$ and $(p_1, \ldots ,p_n) $ we write $x_m = \sum_{i=0}^{n}a_i$.\\
The case of two \elts gives us for each $(a_i , p_i), \; i \in \lrbn$, a \bdp $S_i$. 
Finally, a \bdp for $(x_1, \ldots ,x_m)$ is the concatenation of $S_i$'s and of $a_0$.
\end{Proof}

\rem It is easy to convince ourselves that the \bdp computed by the \algo is minimal: every other \bdp for $(x_1, \ldots ,x_m)$ would be obtained by decomposing certain \elts of the previous basis.

\section{GCD-monoids, GCD-domains}
\label{secAnnPgcd}

Let $G$ be an \grlz.
Since $a\leq b$ \ssi $b\in a + G^+$, the order relation is \caree by the sub\mo $G^+$.
The \egt $x=x^+ - x^-$ shows that the group $G$ can be obtained by symmetrization of the \mo $G^+$,
and it amounts to the same thing to speak of an \grl or of a \mo satisfying certain particular \prts (see \thref{lem1MonGcd}).

We would therefore have had good reason to begin with the theory of  \gui{non-negative submonoids of an \grlz}  rather than with \grlsz.
We would therefore have had good reasons to start by the theory of objects of the type \gui{non-negative sub\mo of an \grlz} rather than by that of \grlsz.
Indeed, in an \grl the order relation must be given at once in the structure, whereas in its non-negative subset, only the law of the \mo intervenes, exactly as in the multiplicative theory of non-negative integers.

It is therefore solely for reasons of comfort in \dems that we have chosen to start with \grlsz.

\subsec{Non-negative sub\mo of an \grlz}

\begin{theorem}\label{lem1MonGcd}
For a commutative \mo $(M,0,+)$ to be the non-negative sub\mo of an \grlz,
it is sufficient and necessary that  conditions~\ref{i1lem1MonGcd},~\ref{i2lem1MonGcd} and~\ref{i3lem1MonGcd} below are satisfied. In addition, we can replace  condition~\ref{i3lem1MonGcd} with condition~\ref{i4lem1MonGcd}.
\begin{enumerate}
\item \label{i1lem1MonGcd}
The \mo is \emph{regular}, \cad $x+y=x+z\Rightarrow y=z$.
\item \label{i2lem1MonGcd}
The preorder relation $x\in y+M$ is an order relation. In other words, we have 
$x+y=0\Rightarrow x=y=0$.\\
We denote it by $y\leq_M x$, or if the context is clear, by $y\leq x$.
\item \label{i3lem1MonGcd}
Two arbitrary \elts admit an upper bound,
\cad

\snic{\forall a,b\;\exists c\;\;\uar c=(\uar a)\cap (\uar b).}
\item \label{i4lem1MonGcd}
Two arbitrary \elts admit a lower bound,
\cad

\snic{\forall a,b\;\exists c\;\;\dar c= (\dar a)\cap (\dar b).}
\end{enumerate}
\index{regular!monoid}
\end{theorem}
\begin{proof}
A priori condition~\emph{\ref{i3lem1MonGcd}} for a particular pair $(a,b)$ is stronger than condition~\emph{\ref{i4lem1MonGcd}} for the following reason: if $a$, $b\in M$, the set of \elts of $M$ less than $a$ and $b$ is contained in $X=\dar(a+b)$.
On this set $X$, the map $x\mapsto a+b-x$ is a bijection that reverses the order and therefore exchanges  supremum and infimum when they exist.
However, in the other direction, the infimum in $X$ (which is the absolute infimum) can a priori 
only be transformed into a supremum for the order relation restricted to the subset $X$, which need not be a global upper bound.
\\
Nevertheless, when condition \emph{\ref{i4lem1MonGcd}} is satisfied for all $a$, $b\in M$, it implies condition~\emph{\ref{i3lem1MonGcd}}.
Indeed, let us show that $m = a+b - (a \vi b)$ is the supremum of $(a, b)$ in $M$ by considering some $x \in M$ such that $x \ge a$ and $x \ge b$. We want to show that $x \ge m$, \cad by letting $y = x \vi m$, that $y \ge m$.
However, $y$ is an upper bound of $a$ and $b$, and $y\in X$.
Since $m$ is the supremum of $a$ and $b$ in~$X$, we indeed have $m\leq y$.\\
The rest of the \dem is left to the reader.
\end{proof}

The previous \tho leads to the notion of a GCD-\moz.
As this notion is always used for the multiplicative \mo of the \ndz \elts of a commutative \riz, we pass to the  multiplicative notation, and we accept that the divisibility relation defined by the \mo is only a preorder relation, in order to take into account the group of units.

\subsec{GCD-monoids}

In multiplicative notation, a commutative \mo $M$  is regular when, for all $a$, $x$, $y\in M$, the \egt $ax=ay$ implies $x=y$.

\begin{definition}\label{defiMoGcd}
We consider a commutative \moz, multiplicatively denoted by $(M,1,\cdot)$.
We say that \emph{$a$ divides $b$} when $b\in a\cdot M$, we also say that~$b$ \emph{is a multiple of $a$}, and we write $a\divi b$. The \mo $M$ is called a \emph{GCD-monoid} when the two following \prts are satisfied%
\index{GCD-monoid}
\begin{enumerate}
\item $M$ is regular.
\item Two arbitrary \elts admits a gcd, \cad

\snic{\forall a,b,\;\exists g,
\;\forall x,
\quad (x \divi a \et x \divi b) \iff x\divi g.}
\end{enumerate}
\end{definition}

Let $U$ be the group of \iv \elts (it is a sub\moz), also called \emph{group of units}. Two \elts $a$ and $b$ of $M$ are said to be \ixc{associated}{elements in a \moz} if there exists an \iv \elt $u$ such that $ua=b$. This is an equivalence relation (we say \gui{the \emph{association} relation}) and the \mo structure passes to the quotient. Let $M/U$ be the quotient \moz. It is still a regular \moz, and the \dve relation, which was a preorder relation on $M$, becomes an order relation on $M/U$. \index{association}

By \thref{lem1MonGcd}, we obtain the following result.

\begin{theorem} \label{thmoGCD}
With the previous notations, a regular commutative \mo $M$ is a GCD-monoid \ssi $M/U$
is the non-negative submonoid of an \grlz.
\end{theorem}

In multiplicative notation, the \dcnsz, partial or complete, are called \emph{\fcnsz}. We then speak of \emph{\bdfz} instead of \bdpz.

Similarly we use the following terminology: a GCD-\mo $M$ 
\emph{satisfies the divisor chain condition} if the \grl $M/U$ is \noez, \cad if in every sequence of \elts $(a_n)_{n\in\NN}$ of $M$ such that $a_{k+1}$ divides $a_k$ for every $k$, there are two associated consecutive terms.%
\index{divisor chain condition}

A GCD-\mo $M$ is said 
\emph{to admit \fabsz}
if $M/U$ admits \dcnbsz, \cad if for each $a$ in $M$ there exists an integer $n$ such that for every \fcn $a=a_1\cdots a_n$ of $a$ in $M$, one of the $a_i$'s is a unit.
It is clear that such a \mo satisfies the divisor chain condition.
\index{factorization!bounded ---}

\subsec{GCD-\risz}
\label{subsecAnnPgcd}\index{GCD-ring}\index{GCD-domain}

We call a \emph{GCD-\riz} a commutative \ri for which the multiplicative \mo of \ndz \elts is a GCD-\moz.
We define in the same way \emph{a bounded \fcn \riz}
or \emph{a \ri which satisfies the divisor chain condition}.%
\index{factorization!bounded --- ring}\index{ring!bounded factorization ---}%

A GCD-domain for which $\Reg(\gA)/\Ati$ admits \faps is called a \emph{GCD-domain admitting \fapsz}. Recall that in particular, the corresponding \grl must be discrete, which here means that $\Ati$ must be a detachable subset of 
$\Reg(\gA)$.%
\index{factorization!partial ---}\index{domain!GCD- --- admitting partial factorizations}

An GCD-domain  for which $\Reg(\gA)/\Ati$ admits  complete \fcns is called a \emph{unique factorization domain}, or a \emph{UFD}. In this case we speak rather of \emph{total \fcnz}.%
\index{factorization!complete ---}%
\index{factorization!total ---}%
\index{unique factorization domain}%
\index{UFD}%
\index{domain!unique factorization ---}

\smallskip
Other than the \gnl results on the GCD-\mos  (which are the translation in multiplicative language of the corresponding results in the \grlsz), we establish some specific facts about GCD-\risz, because the addition intervenes in the statements. They could have been extended to \qiris without difficulty.

\begin{fact}\label{factBezGCD}~
\vspace{-.25em} 
\begin{enumerate}\itemsep=1pt
\item An GCD-domain whose group of units is detachable and which satisfies the divisor chain condition admits \faps (\thref{th2GpRtcl}).
\item A Bézout \ri is a  GCD-\riz.
\item A PID is an GCD-domain  which satisfies the divisor chain condition. 
If the group of units is detachable, the \ri admits \fapsz.
\item If $\gK$ is a nontrivial \cdiz, $\KX$ is a Bézout domain, admits \fabsz, and the group of units is detachable. In particular, the \ri $\KX$ admits \fapsz.
\item The \ris $\ZZ$, $\ZZ[X]$ and $\QQ[X]$ are UFD (Proposition~\ref{propZXfactor}).
%
%
\end{enumerate}
\end{fact}
\facile

\begin{theorem}\label{thAgcdNormal}
Every GCD-domain is integrally closed.
\end{theorem}
\begin{proof}
The \dem of Lemma~\ref{lemZintClos} can be reused word for word.
\end{proof}

We leave to the reader the \dem of the following facts (for \ref{factAXiclgcd}, \KROz's \tho must be used).

\pagebreak	
\begin{fact}\label{factLocaliseGCD}
Let $\gA$  be a GCD-domain and $S$ be a \moz. Then~$\gA_S$ is a GCD-domain, and for $a$, $b\in\gA$ a gcd in~$\gA$ is a gcd in~$\gA_S$.
\end{fact}

We will say that a sub\mo $V$ of a \mo $S$ is \emph{saturated} (in $S$) if $xy\in V$ and  $x$, $y\in S$ imply $x\in V$.
In the literature, we also find \emph{$V$ is factorially closed in $S$}.
A \mo $V$ of a commutative \ri $\gA$ is therefore saturated \ssi
it is saturated in the multiplicative \moz~$\gA$.%
\index{monoid!saturated --- in another}%
\index{saturated!submonoid}%
\index{factorially closed!submonoid}

\begin{fact}\label{factSouspgcdsat}
A saturated sub\mo $V$ of a  GCD-\mo (resp.\,admitting \fabsz) $S$
is a  GCD-\mo (resp.\,admitting \fabsz) with the same gcd and lcm as in $S$.
\end{fact}

\begin{fact}\label{factAXiclgcd}~
\\
 Let $\gA$ be a nontrivial \icl \ri and $\gK$ be its quotient field.
\\
The multiplicative \mo of the \polus of $\AuX=\AXn$ is naturally identified with a saturated sub\mo of $\KuX\etl\!\sur\gK\!\eti$.
\\
In particular, the multiplicative \mo of the \polus of $\AuX$ is a  GCD-\mo admitting \fabsz.
\end{fact}

\penalty-2500
\subsubsec{GCD-domains \ddi$1$}

\begin{definition}\label{defiDim1Integre}
A \qiri $\gA$ is said to be \emph{\ddi$1$} if for every \ndz \elt $a$ the quotient $\aqo{\gA}{a}$ is \zedz.
\index{dimension@(Krull) dimension $\leq1$!\qiri of ---}
\end{definition}

\rem Under the hypothesis that $a$ is \ndzz, we therefore obtain that for all $b$, there exist $x$, $y\in\gA$ and $n\in\NN$ such that

\snic{\;\;\qquad\qquad b^n(1+bx)+ay=0.\qquad\qquad\;\;(*)}

If we no longer make anymore hypotheses about $a$, we can consider the \idm $e$ that generates $\Ann(a)$, and we  then have an \egt of the type $(*)$, but by replacing $a$ by $a+e$, which is \ndzz. This \egt gives, after a  multiplication by $a$ that makes $e$ disappear,

\snic{\qquad\qquad a(b^n(1+bx)+ay)=0\qquad\qquad(+).}

We thus obtain an \egt in accordance with that given in Chapter~\ref{chapKrulldim} where a \cov \dfn \cov of the sentence \gui{$\gA$ is a \ri of \ddk at most $r$} appears, for an arbitrary \ri $\gA$ (see item~\emph{\iref{i3corKrull}} of Proposition~\ref{corKrull}).
\eoe

\begin{lemma}\label{lemDim1-1}\relax
\emph{(A \fcn in dimension $1$)}
\begin{enumerate}
\item Let $\fa$ and $\fb$ be two \ids in a \ri $\gA$ with $\gA\sur\fa$ \zed and
$\fb$ \tfz.
 Then we can write

\snic{\fa=\fa_1\fa_2\;$ with $\;\fa_1+\fb=\gen{1}$ and $\,\fb^n\subseteq \fa_2}

for a suitable integer $n$. This writing is unique and we have 

\snic{\fa_1+\fa_2=\gen{1},$
   $\;\fa_2=\fa+\fb^n=\fa+\fb^{m}$ for every $m\geq n .}  

\item The result applies if $\gA$ is a \qiri \ddi$1$, $\fa$ is \ivz, and $\fb$ is \tfz.
In this case $\fa_1$ and $\fa_2$ are \ivsz.
In particular, $\fa+\fb^n$ is \iv for large enough $n$.
\end{enumerate}

\end{lemma}
\begin{proof}
 It suffices to prove item \emph{1.}
\\ 
\emph{Existence and uniqueness of the \fcnz}.
Consider a triple $(\fa_1,\fa_2,n)$ susceptible of satisfying the hypotheses. Since $\fa_1$ and $\fa_2$ must contain $\fa$, we can reason modulo $\fa$, and therefore suppose $\gA$ is \zed with the \egt $\fa_1\fa_2=\gen{0}$.
\\
 Let $\fa_1+\fb=\gen{1}$ imply $\fa_1+\fb^{\ell}=\gen{1}$ for every exponent $\ell\geq1$. In particular, $\gA=\fa_1\oplus\fa_2=\fa_1\oplus \fb^m$ for every $m\geq n$. This forces, with $e$ being \idmz, $\fa_1=\gen{1-e}$ and $\fa_2=\fb^m=\gen{e}$ for $m$ such that $\fb^m=\fb^{m+1}$ (see Lemma~\ref{lemfacile} and item~\emph{\iref{LID003}} of Lemma~\ref{lemme:idempotentDimension0}).
\end{proof}

\rem Item~\emph{2} is valid without assuming that $\gA$ is a \qiriz. This will become clear after the \gnl \cov \dfn of the Krull dimension, since for every \ndz \elt $a$, if $\gA$ is \ddi$ 1$, the \ri $\aqo\gA a$ is \zedz. 
\eoe

\begin{proposition}\label{lemGCDLop}
Let $\gA$ be a GCD-domain; then every \lop \id is principal.
\end{proposition}
\begin{proof}
Let $\fa = \gen {\an}$ be \lop and $d = \pgcd(\an)$. Let us show that $\fa = \gen{d}$. There exists a \sys of \eco  $(s_1, \ldots, s_n)$  with $\gen {\an} = \gen{a_i}$ in $\gA_{s_i}$. It suffices to see that $\gen {\an} = \gen {d}$ in each $\gA_{s_i}$ because this \egtz, \lot true, will be true globally.
But~$\gA_{s_i}$ remains a  GCD-domain, and the gcds do not change. Therefore, in~$\gA_{s_i}$, we obtain $\gen {\an} = \gen {a_i} = \gen {\pgcd(\an)} = \gen {d}$.
\end{proof}

\vspace{-.4em}
\begin{theorem}\label{propGCDDim1}
A GCD-domain \ddi$ 1$ is a Bézout domain.
\end{theorem}
\begin{proof}
Since $\gen{a,b} = g \gen{a_1,b_1}$ with $\pgcd(a_1,b_1) = 1$, it suffices to show that $\pgcd(a,b) = 1$ implies $\gen{a,b} = \gen{1}$.  However, $\pgcd(a,b) = 1$ implies $\pgcd(a,b^n) = 1$ for every $n \ge 0$.  Finally, after item~\emph{2} of Lemma~\ref{lemDim1-1}, for large enough $n$, $\gen {a, b^n}$ is \iv therefore \lopz,  and \trf by Proposition~\ref{lemGCDLop}.
\end{proof}
%


\subsubsec{Gcd in a \pol \riz}

\smallskip If $\gA$ is a  GCD-domain and $f\in\AX$ we let $\G_X(f)$ or $\G(f)$ be a gcd of the \coes of $f$ (it is defined up to unit \elts multiplicatively) and we call it the \ix{G-content} of $f$.
A \pol whose G-content is equal to $1$ is said to be \ix{G-primitive}.

\pagebreak	
\begin{lemma}\label{lemGcont}
Let $\gA$ be a GCD-domain, $\gK$ be its quotient field  and~$f$ be a nonzero \elt of $\KX$.
\vspace{-.25em} 
\begin{itemize}\itemsep=0pt
\item We can write $f=af_1$ with $a\in\gK$ and $f_1$ as G-primitive in $\AX$.
\item This expression is unique in the following sense:
for another expression of the same type $f=a'f_1'$, there exists a $u\in\Ati$ such that $a'=ua$ and $f_1=uf_1'$.
\item  $f\in\AX$ \ssi $a\in\gA$, in this case $a=\G(f)$.
\end{itemize}
\end{lemma}
\facile

\begin{proposition}
\label{propLG} \emph{(Gauss' lemma, another)} 
Let $\gA$ be a GCD-domain and $f$, $g\in\AX$. Then $\G(fg)=\G(f)\G(g)$. In particular, the product of two G-primitive \pols is a~\hbox{G-primitive} \polz.
\end{proposition}
%
\begin{proof} Let $f_i$ and $g_j$ be the \coes of $f$ and $g$.
It is clear that $\G(f)\G(g)$ divides $\G(fg)$.
By \dit the gcd of the $f_ig_j$'s is equal to $\G(f)\G(g)$, but Proposition~\ref{propArm} implies that $\G(fg)$  divides the $f_ig_j$'s therefore their gcd.
\end{proof}
%

\begin{corollary}\label{corpropLG}
Let $\gA$ be a GCD-domain, $\gK$ be its quotient field and $f$, $g\in\AX$. Then $f$ divides $g$ in $\AX$ \ssi $f$ divides~$g$ in $\KX$ and $\G(f)$ divides $\G(g)$ in $\gA$.
\end{corollary}
\begin{proof}
The \gui{only if} results from Gauss' lemma.
For the \gui{if} we can suppose that $f$ is G-primitive.
If $g=hf$ in $\KX$, we can write $h=ah_1$ where $h_1\in\AX$ is G-primitive and $a\in\gK$. By Gauss' lemma, we have $fh_1$ G-primitive.
By applying Lemma~\ref{lemGcont} to the \egt $g = a(h_1f)$, we obtain $a \in \gA$, then $h\in\AX$.
\end{proof}
%


Recall that if $\gA$ is a reduces \riz, $\AX\eti=\Ati$
(Lemma~\ref{lemGaussJoyal} \emph{\iref{i4lemPrimitf}}).
In particular, if $\gA$ is a nontrivial domain and if the group of units of $\gA$ is detachable, the same goes for~$\AX$.

\begin{theorem}\label{thAXgcd}
Let $\gA$ be a GCD-domain and $\gK$ be its quotient field.
\begin{enumerate}
\item $\AXn$ is a GCD-domain.
\item If $\gA$ admits \fapsz, the same goes for~$\AX$.
\item If $\gA$ satisfies the divisor chain condition, the same goes for~$\AX$.
\item If $\gA$ admits \fabsz, the same goes for~$\AX$.
\item If $\AX$ is a UFD, the same goes for~$\AXn$ \emph{(Kronecker)}.
\end{enumerate}
\end{theorem}
\begin{proof}
\emph{1.} It suffices to treat the case $n=1$. 
Let $f$, $g\in\AX$. \\
Let us express $f=af_1$, $g=bg_1$, with G-primitive $f_1$ and $g_1$. Let $c=\pgcd_\gA(a,b)$ \hbox{and $h=\pgcd_\KX(f_1,g_1)$}. We can assume \spdg \hbox{that $h$} is in $\AX$ and that it is G-primitive. Then,
by using Corollary~\ref{corpropLG}, we verify that $ch$ is a gcd of $f$ and $g$ in $\AX.$
 
Items~\emph{2}, \emph{3} and \emph{4} are left to the reader.
 
\emph{5.} It suffices to treat the case $n=2$ and to know how to detect if a \pol admits a strict factor. 
We use the \KRA trick.
To test the \polz~$f(X,Y)\in\gA[X,Y]$, assumed of degree $< d$ in $X$, 
we consider the \polz~$g(X)=f(X,X^d)$. A \dcnc of $g(X)$ allows us to know if there exists a strict factor of $g$ of the form $h(X,X^d)$ (by considering all the strict factors of $g$, up to association), which corresponds to a  strict factor of $f$. For some precisions see Exercise~\ref{exoKroneckerTrick}.
\end{proof}
%

\begin{corollary}\label{corthAXgcd}
If $\gK$ is a nontrivial \cdiz, $\KXn$ is a GCD-domain, admitting \fabs and  \fapsz.  The group of units is $\gK\eti$. Finally, $\KXn$ ($n\geq 2$) is a UFD \ssi $\KX$ is a UFD.
\end{corollary}

\section{Zariski lattice of a commutative \riz}
\label{secZarAcom}

\vspace{4pt}
\subsec{Generalities}

Recall the notation $\DA(\fa)$ with some precisions.
\begin{notation}
\label{notaZA}\label{defZar}
{\rm  If $\fa$ is an \id of $\gA$, let $\DA(\fa)=\sqrt{\fa}$ be the nilradical of~$\fa$. If $\fa=\gen{x_1,\ldots ,x_n}$ let $\DA(x_1,\ldots ,x_n)$ for $\DA(\fa)$. \hbox{Let $\ZarA$} be the set of $\DA(x_1,\ldots ,x_n)$ (for $n\in\NN$ and $x_1$, \ldots, $x_n\in\gA$).
}
\end{notation}

We therefore have $x\in\DA(x_1,\ldots ,x_n)$ \ssi a power of $x$ is a member of $\gen{x_1,\ldots ,x_n}$.

The set $\ZarA$ is ordered by the inclusion relation.

\begin{fact}
\label{factZar}
$\ZarA$ is a \trdi with
$$\preskip.4em \postskip.3em
\begin{array}{rcl}
\DA(0)=0_{\ZarA}, && \DA(\fa_1)\vu\DA(\fa_2)=\DA(\fa_1+\fa_2), \\[.2em]
\DA(1)=1_{\ZarA},  &\quad   &
\DA(\fa_1)\vi\DA(\fa_2)=\DA(\fa_1\,\fa_2).
\end{array}$$
We call it the \emph{Zariski lattice of the \ri $\gA$}.%
\index{lattice!Zariski ---}
\end{fact}

In \clama $\DA(x_1,\ldots ,x_n)$ can be seen as a \oqc of $\SpecA$: the set of \ideps $\fp$ of $\gA$ such that at least one of the $x_i$'s does not belong to $\fp$, and $\ZarA$ is identified with the lattice of the \oqcs of $\SpecA$.
For more details on the subject see Section~\ref{secEspSpectraux}.

\begin{fact}\label{fact1Zar}~
\begin{enumerate}
\item For every morphism $\varphi:\gA\to\gB$, we have a natural morphism $\Zar\varphi$ \hbox{from $\ZarA$} \hbox{to $\Zar\gB$}, and we thus obtain a functor from the category of commutative \ris to that of \trdisz.
\item For every \ri $\gA$ the natural \homo $\ZarA\to\Zar\Ared$ is an \isoz, so that we can identify the two lattices.
\item The natural \homo $\Zar(\gA_1\times \gA_2)\to\Zar\gA_1\times \Zar\gA_2$ is an \isoz.
\item For a \agB $\gB$, the map $x\mapsto \DB(x)$ is an \iso from $\gB$ to $\Zar\,\gB$.
\end{enumerate}
\end{fact}

\begin{fact}\label{factZarABol}
\Propeq
\begin{enumerate}
\item $\ZarA$ is a \agBz.
\item $\gA$ is \zedz.
\end{enumerate}
\end{fact}

\begin{proof} Recall that a \trdi \gui{is} a \agB \ssi every \elt admits a complement
(Proposition~\ref{defiBooleTrdi}).
\\
Suppose \emph{2.} Then for every \itf $\fa$, there exist an \idm $e$ and an integer $n$ such that $\fa^n=\gen{e}$. Therefore $\DA(\fa)=\DA(e)$.
Moreover, it is clear that $\DA(e)$ and $\DA(1-e)$ are complements in $\ZarA$.
\\
Suppose \emph{1.} Let $x\in\gA$ and $\fa$ be a \itf of $\gA$ such that $\DA(\fa)$ is the complement of $\DA(x)$ in $\ZarA$. Then there exist $b\in\gA$ and $a\in\fa$ such \hbox{that $bx+a=1$}. As $xa=x(1-bx)$ is nilpotent we obtain an \hbox{\egt $x^n(1+cx)=0$}.
\end{proof}

\begin{fact}\label{fact2Zar} \label{corthfactorisation2}
Let $a\in\gA$ and $\fa\in\ZarA$.
\begin{enumerate}
\item The \homo $\Zar\pi:\ZarA\to\Zar(\aqo{\gA}{a})$,
where $\pi:\gA\to \aqo{\gA}{a}$ is the canonical projection, is surjective, and it allows us to identify $\Zar(\aqo{\gA}{a})$ with the quotient lattice $\Zar(\gA)\sur{(\DA(a)=0)}$.
%
More \gnltz, $\Zar(\gA\sur{\fa})$ is identified with $\Zar(\gA)\sur{(\fa=0)}$.
\item  The \homo $\Zar j:\ZarA\to\Zar(\gA[1/a])$, where $j:\gA\to \gA[1/a]$ is the canonical \homoz, is surjective and it allows us to identify~$\Zar(\gA[1/a])$ with the quotient lattice $\Zar(\gA)\sur{(\DA(a)=1)}$.
\item  For some \id $\fc$   and some \mo $S$ of  $\gA$ we have a natural \iso

\snic{\Zar(\gA_S\sur{\fc\gA_S})\;\simeq\; \Zar(\gA)\sur{(\fb=0,\ff=1)},}

where $\fb$ is the \id of $\ZarA$ generated by the $\DA(c)$'s for $c\in\fc$, and $\ff$ is the filter of $\ZarA$ generated by the $\DA(s)$'s for $s\in S$.
\end{enumerate}
\end{fact}


\subsec{Duality in the commutative \risz}\label{secIDEFIL} \perso{Je ne sais vraiment
pas where doit se trouver cette section in le livre.}

\subsubsec{Annihilating and inverting simultaneously}
\label{secAnnEtInv}

In the \trdis we exchange the roles of $\vi$ and $\vu$ by passing to the opposite lattice, \cad by reversing the order relation.

In the commutative \risz, a fecund duality also exists between the addition and the multiplication,
more mysterious when we try to exchange their roles.


Recall that a saturated \mo is called a \emph{filter}.
The notion of filter is a dual notion to that of \idz, just as important.

The ideals are the inverse images of $0$ under the homomorphisms. They serve to pass to the quotient, \cad to annihilate elements by force. The filters are the inverse images of the group of units under the homomorphisms. They serve to localize, i.e. to render elements invertible by force.

Given an \id $\fa$ and a \mo $S$ of the \ri $\gA$ we may want to annhilate the \elts of $\fa$ and invert the \elts of $S$.
The solution of this \pb is given by consideration of the following \riz.

\begin{definota}\label{defiASa}
Let (by abuse) $\gA_{S}\sur{\fa}$ or $S^{-1}\gA\sur{\fa}$ be the \ri whose \elts are given by the pairs $(a,s)\in\gA\times S$, with the \egt $(a,s)=(a',s')$ in $\gA_{S}\sur{\fa}$ \ssi there exists an $s''\in S$  such that $s''(as'-a's)\in \fa$ (we will write $a/s$ for the pair $(a,s)$).
\end{definota}

The fact that $\gA_{S}\sur{\fa}$ defined thus answers the posed \pb signifies that the following  \fcn  \tho is true (see the analogous Facts~\ref{factUnivQuot} and~\ref{factUnivLoc}).

\begin{fact}
\label{thfactorisation2} \emph{(Factorization \thoz)}\\
With the above notations, let $\psi:\gA\to\gB$ be a \homoz. Then $\psi$ is factorized by $\gA_S\sur{\fa}$
\ssi $\psi(\fa)\subseteq \so{0}$ and $\psi(S)\subseteq \gB^\times$. In this case, the \fcn is unique.
\end{fact}

\vspace{-1em}
\pun{\gA}{\lambda}{\psi}{\gA_S\sur{\fa}}{\theta}{\gB}{$\psi(\fa)\subseteq \so{0}$ and $\psi(S)\subseteq \gB^\times$}

\vspace{-1em}
Naturally we can also solve the \pb by first annihilating $\fa$ then by inverting (the image of) $S$, or by  first inverting $S$ then by annihilating (the image of) $\fa$. We thus obtain canonical \isos
$$\preskip.4em 
\gA_{S}\sur{\fa} \;\simeq \;\big(\pi_{\gA,\fa}(S)\big)^{-1}(\gA\sur{\fa})
\;\simeq\; (\gA_{S})\sur {\left(j_{\gA,S}(\fa)\gA_{S}\right)}.
$$

\subsubsec{Dual \dfnsz}

The duality between \ids and filters is a form of duality between  addition and  multiplication.

This is easily seen from the respective axioms that are used to define the \ids (resp.\,\idepsz) and the filters (resp.\,prime filters)
$$\arraycolsep3pt\begin{array}{rclcrcl}
     \mathrm{ideal}   &\fa   &   & \qquad\qquad   & \mathrm{filter}
& \ff & \\[1mm]
       &\vdash   &  0\in\fa & \quad   &   & \vdash & 1\in\ff\\
   x\in\fa, \, y\in\fa     &\vdash   & x+y\in\fa  & \quad
  &x\in\ff, \, y\in\ff   & \vdash &  xy\in\ff\\
   x\in\fa      &\vdash   & xy\in\fa  & \quad
& xy\in\ff  & \vdash & x\in\ff\\[1mm]
     \hbox{prime ---}   &  &   & \quad   & \hbox{prime ---}  &   & \\   xy\in\fa
&\vdash   & x\in\fa \lor y\in\fa  & \quad   &   x+y\in\ff     &\vdash   & x\in\ff
\lor y\in\ff
\end{array}$$

Note that according to the above \dfnz,  
$\fa$ is both a \idep and a prime filter of $\gA$.
This convention can seem strange, but it happens to be the most practical one: an \id is prime \ssi the quotient \ri is \sdzz, a filter is prime \ssi the localized \ri is a \aloz.
With regard to ideals we have already commented on this
on \paref{CommIdeps}.%
\index{filter!prime ---}%
\index{prime!filter}

We will adopt the following \dfn for a \emph{maximal filter}: the localized \ri is a \zed \alo (when the \ri is reduced: a \cdiz). 
In particular, every maximal filter is prime.
We will essentially make use of this \dfn as a heuristic.
\index{filter!maximal ---}\label{labfima}
\index{maximal!filter}

Now suppose the \ri $\gA$ is nontrivial.
Then a detachable strict \id (resp.\,a detachable strict filter) is prime \ssi its complement is a filter (resp.\,an \idz).
We once again find in this case the familiar ground in \clamaz.

Generally in \clama the complement of a strict \idep is a strict prime filter and vice versa, therefore the complement of a strict \idema is a minimal prime filter, and the complement of a strict \fima is a \idemiz. The prime filters therefore seem more or less useless and have a tendency to disappear from the scene in \clamaz.


\subsubsec{Saturated pairs}

A good way to understand the duality is to simultaneously treat \ids and filters. For this we introduce the notion of a \emph{\paz}, analogous to that which we have given for  \trdisz.

\begin{definition}
\label{defpa}
Let $\fa$ be an \id and $\ff$ be a filter of $\gA$.
We say that $\fa$ is \emph{$\ff$-saturated} if we have
$$\preskip.0em \postskip.4em
(as\in\fa,\,s\in\ff)\Longrightarrow a\in\fa,
$$
we say that $\ff$ is \emph{$\fa$-saturated} if we have
$$\preskip.3em \postskip.4em
(a+s\in\ff,\,a\in\fa)\Longrightarrow s\in\ff.
$$
If $\fa$ is $\ff$-saturated and $\ff$ is $\fa$-saturated we say that $(\fa,\ff)$
is a \emph{\paz} in~$\gA$.%
\index{saturated!pair}%
\index{saturated!$\ff$- --- ideal}%
\index{saturated!$\fa$- --- filter}%
\end{definition}

To recap the axioms for the \pas (note that the last condition can be rewritten as $\fa+\ff=\ff$).
$$\arraycolsep3pt\begin{array}{rclcrcl}\preskip.1em \postskip.3em
       &\vdash   &  0\in\fa & \phantom{aa}   &   & \vdash & 1\in\ff\\[.2em]
   x\in\fa, \, y\in\fa     &\vdash   & x+y\in\fa  &
  &x\in\ff, \, y\in\ff   & \vdash &  xy\in\ff\\[.2em]
   x\in\fa      &\vdash   & xy\in\fa  &
& xy\in\ff  & \vdash & x\in\ff\\[.2em]
   xy\in\fa,\,y\in\ff      &\vdash   & x\in\fa  &
& x+y\in\ff,\,y\in\fa  & \vdash & x\in\ff
\end{array}$$

\begin{fact}
\label{factPaire}~
\begin{enumerate}\itemsep0pt

\item For every \homo $\varphi:\gA\to\gB$, the pair $\big(\Ker\varphi,\varphi^{-1}(\gB^{\times})\big)$ is a \paz.

\item Conversely if $(\fa,\ff)$ is a \pa and if $\psi:\gA\to\gA_{\ff}\sur{\fa}=\gC$ designates the canonical \homoz, we have {\mathrigid 2mu $\Ker\psi=\fa$} and {\mathrigid 2mu $\psi^{-1}(\gC^{\times})=\ff$}.

\item Let $\varphi:\gA\to\gC$ be a \homo and $(\fb,\ffg)$ be a \pa of $\gC$, then
$\big(\varphi^{-1}(\fb),\varphi^{-1}(\ffg)\big)$ is a \pa of $\gA$.

\end{enumerate}
\end{fact}

\begin{fact}
\label{factPaire2}
Let $(\fa,\ff)$  be a \paz.
\begin{enumerate}\itemsep0pt
\item  $\gA_{\ff}\sur{\fa}$ is local \ssi $\ff$ is a prime filter
(\cad \ssi $\gA_{\ff}$ is local).
\item  $\gA_{\ff}\sur{\fa}$ is \sdz \ssi $\fa$ is a prime \id 
(\cad \ssi $\gA\sur{\fa}$ is \sdzz).
\end{enumerate}
\end{fact}

\begin{definition}
\label{defRaffine}
If $(\fa,\ff)$ and $(\fb,\ffg)$ are two \pas of $\gA$ we say that
\emph{$(\fb,\ffg)$ refines  $(\fa,\ff)$}
 and we write it $(\fa,\ff)\leq (\fb,\ffg)$ when $\fa\subseteq\fb$ and $\ff\subseteq\ffg$.
\end{definition}\index{refine}

The following lemma describes the \pa \gui{generated} (in the sense of the refinement relation) by a pair of subsets of $\gA$.
Actually it suffices to treat the case of a pair formed by an \id and a \moz.

\begin{lemma}
\label{lempaireengendree}
Let $\fa$ be an \id and $\ff$ of $\gA$ be a \moz.
\begin{enumerate}\itemsep0pt

\item The \pa $(\fb ,\ffg )$ generated by $(\fa ,\ff )$ is obtained as follows
\[\preskip.2em \postskip.2em
\fb =\sotq{x\in\gA}{\exists s\in \ff,\, xs\in \fa },
\hbox{ and }\,
\ffg =\sotq{y\in\gA}{\exists u\in\gA,\, uy\in \fa +\ff}.
\]

\item If $\ff \subseteq\Ati$, then
$\fb=\fa$ and $\ffg $ is the filter obtained by saturating the \moz~\hbox{$1+\fa $}. In this case, $\gA_{\ffg }\sur{\fa }=\gA\sur{\fa }$.

\item If $\fa =0$, then
$\fb =\sotq{x\in\gA}{\exists s\in \ff,\,xs=0}=\sum_{s\in\ff}(0:s)$, and  $\ffg$ is the saturation of $\ff$.
In this case, $\gA_{\ffg }\sur{\fb }=\gA_{\ff }$.
If in addition $\ff=s^{\NN} $, $\fb =(0:s^\infty)$.
\end{enumerate}
\end{lemma}

\rdb \label{NOTASatu}

An important case is that of the filter obtained by saturation of a \moz~$S$.
We introduce the notation $\sat{S}$, or, if \ncrz, $\satu S \gA$
for this filter.

\subsubsection*{Incompatible ideal and filter}

For any \pa $(\fa,\ff)$ we have the following \eqvcsz.
\begin{equation}\label{eqIncompatibles}\preskip.3em \postskip.3em
\fa=\gA\;\Longleftrightarrow\; 1\in\fa\;\Longleftrightarrow\;
0\in\ff\;\Longleftrightarrow\; \ff=\gA\;\Longleftrightarrow\;
\gA_\ff\sur{\fa}=\so{0}.
\end{equation}

An \id $\fa$ and a filter $\ff$ are said to be \emph{incompatible} when they generate the pair~\hbox{$(\gA,\gA)$}, \cad when $0\in\fa+\ff$.

An \id $\fa$ and a filter $\ff$ are said to be \emph{compatible} if they 
{satsify $(0\in\fa+\ff\Rightarrow1=0)$}.
If the \ri is nontrivial this also means $\fa\cap\ff= \emptyset$.
In this case we can both annihilate the \elts of $\fa$ and render the \elts of $\ff$ \ivs without reducing the \ri to $0$.%
\index{compatible!ideal and filter}%
\index{compatible!\paz}%
\index{incompatible!\paz}%
\index{incompatible!ideal and filter}

\begin{fact}\label{factPremComp} \hspace*{-.6em}
Let $\fa$ be an \id and $\ff$ be a compatible filter.
\\
If $\fa$ is prime, it is $\ff$-saturated, if $\ff$ is prime,
it is $\fa$-saturated.
\end{fact}


\begin{fact}
\label{corlempaireengendree} \emph{(The lattice of \pasz)} 
The \pas of $\gA$  have a lattice structure for the refinement relation, 
such that
\begin{enumerate}
\item [--] The minimum \elt is $(\so{0},\Ati)$ and the maximum \elt  $(\gA,\gA)$.
\item [--] $(\fa,\ff)\vu(\fb,\ffg)$ is the \pa generated by $(\fa+\fb,\ff\,\ffg)$.
\item [--] $(\fa,\ff)\vi(\fb,\ffg)=(\fa\cap\fb,\ff\cap\ffg)$.
\end{enumerate}
\end{fact}


\addcontentsline{toc}{subsubsection}{Ideals and filters in a localized quotient \riz}
\begin{fact}
\label{propPaires} \emph{(Ideals and filters in a localized quotient \riz)}
 Let  $(\fa,\ff)$ be a \pa of $\gA$ and $\pi:\gA\to\gB=\gA_\ff\sur{\fa}$ be the canonical map.
Then
\begin{enumerate}
\item
The map $(\fb,\ffg)\mapsto \big(\pi^{-1}(\fb),\pi^{-1}(\ffg)\big)$
is a non-decreasing bijection (for the refinement relations) between on the one hand, the \pas of $\gB$, and on the other, the \pas of $\gA$ which refine $(\fa,\ff)$.
\item If $(\fb,\ffg)$ is a \pa of $\gB$ the canonical map
$$\preskip.2em \postskip-.10em
\gA_{\pi^{-1}(\ffg)}\sur{\pi^{-1}(\fb)} \;\longrightarrow\; \gB_{\ffg}\sur{\fb}$$
is an \isoz.
\item In this bijection
\begin{enumerate}
\item [--] the \id $\fb$ is prime \ssi $\pi^{-1}(\fb)$ is prime,
\item [--] every prime \id of $\gA$ compatible with $\ff$ and containing $\fa$ is obtained,
\item [--] the filter $\ffg$ is prime \ssi $\pi^{-1}(\ffg)$ is prime,
\item [--] every prime filter of $\gA$ compatible with $\fa$ and containing $\ff$ is obtained.
\end{enumerate}
%
\end{enumerate}
\end{fact}

We deduce the following instructive comparison on the duality between \ids and filters.

\vspace{-.5em}
\DeuxCol
{
\begin{fact}
\label{factQuoIDFI}
Let $\fa$ be a strict \id of $\gA$ and $\pi:\gA\to
\gA/\fa$ be the corresponding \homoz.
\begin{enumerate}
\item The map $\fb\mapsto\pi^{-1}(\fb)$ is a non-decreasing bijection between \ids of $\gA/\fa$ and \ids of $\gA$ containing $\fa$.
In this bijection the \ideps correspond to the \idepsz.
\item The map $\ffg\mapsto\pi^{-1}(\ffg)$ is a non-decreasing bijection between filters of $\gA/\fa$ and $\fa$-saturated filters of $\gA$.
\item  In this bijection the strict prime filters of $\gA/\fa$ correspond exactly to the prime filters of $\gA$
compatible with~$\fa$.
\end{enumerate}
\end{fact}

}
{
\begin{fact}
\label{factQuoFIID}
Let $\ff$ be a strict filter of $\gA$ and $\pi:\gA\to
\gA_\ff$ be the corresponding \homoz.
\begin{enumerate}
\item The map $\ffg\mapsto\pi^{-1}(\ffg)$ is a non-decreasing bijection between filters of $\gA_\ff$ and filters of $\gA$ containing $\ff$.
In this bijection the prime filters correspond to the prime filters.
\item The map $\fb\mapsto\pi^{-1}(\fb)$ is a non-decreasing bijection between \ids of $\gA_\ff$ and $\ff$-saturated \ids of $\gA$.
\item In this bijection the strict \ideps of $\gA_\ff$ correspond exactly to the \ideps of $\gA$
compatible with~$\ff$.
\end{enumerate}
\end{fact}

}

\subsec{Closed covering principles}

The duality between \ids and filters suggests that a dual principle of the \plg must be able to function in commutative \algz.
First of all note that the \ids of $\ZarA$ bijectively correspond to the radical \ids (\cad equal to their nilradical) of $\gA$ via 
$$\preskip.4em \postskip.4em 
\fa \;(\hbox{\id of }\ZarA) \mapsto \sotq{x\in\gA}{\DA(x)\in\fa }. 
$$
In addition, the \ideps correspond to the \idepsz.

For filters, things are not quite so perfect, but for a filter $\ff$ of~$\gA$, the set $\sotq{\DA(x)}{x\in\ff}$ generates a filter of $\ZarA$, and this gives an injective map which is  bijective for the first filters.

Let us return to the \plg and look at what it means in the lattice $\ZarA$. When we have \moco $S_1$, \ldots, $S_n$ of $\gA$, it corresponds to filters $\ff_i $ of $\ZarA$
(each generated by the $\DA(s)$'s for $s\in S_i$)
which are \gui{\comz}  in the sense that $\bigcap_i\ff_i=\so{1_\ZarA}$.
In this case the natural \homos
$$\preskip.4em \postskip.0em\ndsp 
\gA\to\prod_i\gA_{S_i}\quad \hbox{and}
\quad\ZarA\to\prod_i\ZarA\sur{(\ff_i=1)} 
$$
are injective.

By duality, we will say that a \sys of \ids $(\fa_1,\ldots,\fa_n)$ constitutes a \ix{closed covering} of $\gA$ when $\bigcap_i\DA(\fa_i)=\so{0_\ZarA}$,
\cad  when $\prod_i\fa_i\subseteq\DA(0)$.
In this case the natural \homos
$$\preskip.4em \postskip.2em\ndsp 
\gA/\DA(0)\to\prod_i\gA\sur{\DA(\fa_i)}\quad\hbox{and}
\quad \ZarA\to\prod_i\ZarA\sur{(\DA(\fa_i)=0)} 
$$
 are injective.

We will say that a \prt $\sfP$ (regarding objects related to a \riz~$\gA$) satisfies the \gui{closed covering principle} when:\\
\emph{each time that \ids $\fa_i$ form a closed covering of $\gA$, the \prt $\sfP$ is true for $\gA$ \ssi it is true after passage to the quotient by each of the~$\fa_i$'s.
}

For example we easily obtain (see also Lemma~\ref{lemNilpotProd}).

\begin{prcf}\label{prcf1}\hspace*{-.6em} \emph{(Nilpotent, \com \eltsz)}\\
Consider a closed covering $(\fa_1,\ldots,\fa_r)$ of the \ri $\gA$. Let $x_1$, \dots, $x_n\in\gA$, $\fb$, $\fc$ be two \ids and $S$ be a \moz.
\begin{enumerate}
\item The \mo $S$ contains $0$  \ssi it contains $0$ modulo each~$\fa_i$.
\item We have $\fb\subseteq\sqrt{\fc}$ \ssi  $\fb\subseteq\sqrt{\fc}$
modulo each~$\fa_i$.
\item The \elts $x_1$, \dots, $x_n$ are \com \ssi they are \com modulo each~$\fa_i$.
\end{enumerate}
\end{prcf}
\begin{proof}
It suffices to prove item \emph{2.} Suppose that $\DA(\fb)\leq \DA(\fc)\vu \DA(\fa_i)$,
so $\DA(\fb)\leq \Vi_i \left(\DA(\fc)\vu \DA(\fa_i)\right) = \DA(\fc)\vu \left(\Vi_i\DA(\fa_i)\right)=\DA(\fc)$.
\end{proof}
\rem However, there is no \prf for the solutions of \slisz.
Indeed, consider $u,\,v\in\gA$ such \hbox{that $uv=0$}.
The \sli (with $x$ as the unknown)  

\snic{ux=u,$  $vx=-v,}

admits a solution modulo $u$ (namely $x=-1$) and a solution
modulo $v$ (namely $x=1$). But in the case of the \ri
$\gA=\ZZ[u,v]=\aqo{\ZZ[U,V]}{UV}$ the \sli has no solution in $\gA$.
\eoe

\begin{prcf}\label{prcf2} \emph{(\Tf modules)}\\ 
Consider a closed covering $(\fa_1,\ldots,\fa_r)$ of the \ri $\gA$.
Suppose \hbox{that $\prod_i\fa_i=0$} (this is the case if $\gA$ is reduced).
An \Amo $M$ is \tf  \ssi it is \tf modulo each~$\fa_i$.
\end{prcf}
\begin{proof}
Suppose \spdg that $r=2$. Let $g_1$, \ldots, $g_k$ be \gtrs modulo $\fa_1$,
and  $g_{k+1}$, \ldots, $g_\ell$ be \gtrs modulo $\fa_2$. Let $x\in M$. 
We write
$x=\sum_{i=1}^k\alpha_ig_i+\sum_{j=1}^p\beta_jx_j$ with $\alpha_i\in\gA$,
$\beta_j\in\fa_1$, $x_j\in M$.
\\ 
Each $x_j$ is written as a \coli of  $g_{k+1}$, \ldots, $g_\ell$
modulo $\fa_2$.  \hbox{Since $\fa_1\fa_2=0$}, we obtain $x$ as a \coli \hbox{of
$g_1$, \ldots, $g_\ell$}.
\end{proof}
%

\begin{prcf}\label{prcf3} \emph{(\Mptfsz)}
Consider a closed covering $(\fa_1,\ldots,\fa_r)$ of the \ri
$\gA$, a matrix  $F\in\Ae{m\times n}$, a \itf $\fa$ and a \mpfz~$M$.
\begin{enumerate}
\item The matrix $F$ is of rank $\geq k$  \ssi it  is of rank $\geq k$
modulo each~$\fa_i$.
\end{enumerate}
Suppose $\bigcap_i\fa_i=0$ (it is the case if $\gA$ is reduced). Then
\begin{enumerate} \setcounter{enumi}{1}
\item The matrix $F$ is of rank $\leq k$  \ssi it is of rank $\leq k$ modulo each~$\fa_i$.
\item The \itf $\fa$ is generated by an \idm  \ssi it is generated by an \idm modulo each~$\fa_i$.
\item The matrix $F$ is \lnl  \ssi it is \lnl modulo each~$\fa_i$.
\item The module $M$ is \ptf  \ssi it is \ptf  modulo each~$\fa_i$.
\end{enumerate}
\end{prcf}
\begin{proof}
Item \emph{1} results from the \prf \ref{prcf1} by considering the \idd of order $k$. Item~\emph{2} comes from the fact that if a \idd is null modulo each $\fa_i$, it is null modulo their intersection.
Item~\emph{5} is a reformulation of item~\emph{4} 
which is a consequence of item~\emph{3.}\\
Let us prove item~\emph{3.} Suppose \spdg that $r=2$.
We use the  the lemma of the \id generated by an \idm (Lemma~\ref{lem2ide.idem}).
We have

\snic{\fa+(0:\fa)_{\gA\sur{\fa_i}}=\gA\sur{\fa_i}\;\; (i=1,2). }

This means that $\fa+\fa_i+(\fa_i:\fa)=\gA$, and since $\fa_i\subseteq(\fa_i:\fa)$,
we have $1\in\fa+(\fa_i:\fa)$. By taking the product we get $1\in\fa+(\fa_1:\fa)(\fa_2:\fa)$ and since
$$
\preskip.2em \postskip.4em 
(\fa_1:\fa)(\fa_2:\fa)\subseteq  (\fa_1:\fa)\cap(\fa_2:\fa) =((\fa_1\cap\fa_2):\fa)=(0:\fa), 
$$
 we obtain $1\in\fa+(0:\fa)$.
\end{proof}
%

\subsec{Reduced \zede closure of a commutative \riz}
\label{secClotureZEDR}

Let us begin with some results regarding a sub\ri $\gA$ of a \zedr \riz. The reader can refer to the study of \zedr \ris on \paref{subsecAzedred} and revisit \Egtsz~(\iref{eqQuasiInv}) for the \carn of a quasi-inverse.

If in a \ri an \elt $c$ admits a quasi-inverse, we denote it \hbox{by $c\bul$},
and we denote by $e_c=c\bul c$ the \idm associated with~$c$ which satisfies the \egts $\Ann(c)=\Ann(e_c)= \gen {1-e_c}$.

\begin{lemma}\label{lem1SousZedRed}
\emph{(\Ri generated by a quasi-inverse)}
\begin{enumerate}
\item Let $a\in\gA\subseteq\gB$. Suppose that $\gA$ and $\gB$ are reduced and that $a$ admits a quasi-inverse in $\gB$.
Then

\snic{\gB\supseteq\gA[a\bul]\simeq\aqo{\gA[a\bul]}{1-e_a}\times \aqo{\gA[a\bul]}{e_a}=
\gA_1\times \gA_2.}

In addition
\begin{enumerate}
\item We have a well-defined natural \homo $\gA[1/a]\to\gA_1$, and it is an \isoz. In particular, the natural \homoz~\hbox{$\gA\to\gA_1$} has as its kernel $\Ann_\gA(a)$.
\item
The natural \homo $\gA\to\gA_2$ is surjective, its kernel is the intersection \hbox{$\fa=\gA\cap e_a\gA[a\bul]$} and satisfies the double inclusion
$$\preskip.3em \postskip.2em 
 \Ann_\gA\big(\Ann_\gA(a)\big)\supseteq \fa \supseteq \DA(a). \eqno(*)
$$
\end{enumerate}
In short $\gA[a\bul]\simeq\gA[1/a]\times \gA\sur{\fa}$.
\item Conversely for every \id $\fa$ of $\gA$ satisfying $(*)$, the \elt $(1/a,0)$ is a quasi-inverse of (the image of) $a$ in the \ri $\gC=\gA[1/a]\times \gA\sur{\fa}$ and the canonical \homo of $\gA$ in $\gC$ is injective.
\end{enumerate}
\end{lemma}
\begin{proof}
The \iso $\gA[a\bul]\simeq \gA_1\times \gA_2$ only means that $e_a$ is an \idm in $\gA[a\bul]$. Let
$\pi_i:\gA[a\bul]\to\gA_i$ be the canonical \homosz.

\emph{1b.} Let $\mu$ be the composed \homo $\gA\lora \gA[a\bul] \lora \gA_2$.
In $\gA_2$,  
we have \hbox{$a\bul=e_aa\bul=0$}, so
$\gA_2=\gA\sur{\left(\gA\cap {e_a \gA[a\bul]}\right)}$. Thus~\hbox{$\fa=\gA\cap e_a\gA[a\bul]$}.
In~$\gA[a\bul]$, we have $a=e_aa$, so $\mu(a)=\pi_2(a)=\pi_2(e_aa)=0$,  
and $a\in\fa$.
\\
As $\gB$ is reduced, the three \ris $\gA[a\bul]$, $\gA_1$ and $\gA_2$ are also reduced. 
Therefore $\gen{a}\subseteq \fa$ implies $\DA(a)\subseteq \fa$.\\
Finally,
$\fa\,\Ann_\gA(a) \subseteq \gen{e_a}\Ann_\gA(a) = 0$, so
     $\fa\subseteq\Ann_\gA\big(\Ann_\gA(a)\big)$.


\emph{1a. }
Since $aa\bul=_{\gA_1}1$, we have a unique \homo $\lambda:\gA[1/a]\to\gA_1$
obtained from the composed \homo $\gA\to \gA[a\bul] \to \gA_1$, and $\lambda$ is clearly surjective. Consider an \elt $x/a^n$ of $\Ker\lambda$. Then $\lambda(ax)=0$, so $\pi_1(ax)=0$. As we also have $\pi_2(ax)=0$, we deduce that $ax=0$,  
so $x=_{\gA[1/a]}0$. Thus $\lambda$ is injective.

\emph{2.} The image of $a$ in $\gC$ is $(a/1,0)$, so $(1/a,0)$ is indeed its quasi-inverse. Now let $x\in\gA$ whose image in $\gC$ is $0$.
On the one hand $x=_{\gA[1/a]}0$, so $ax=_\gA0$. On the other hand $x\,\Ann_\gA(a)=0$ so $x^2=_\gA0$, and $x=_\gA0$.
\end{proof}

\comm We see
that the notation $\gA[a\bul]$ presents a priori a possible ambiguity, at least when $\DA(a)\neq\Ann_\gA\big(\Ann_\gA(a)\big)$. \eoe

\begin{lemma}\label{lem2SousZedRed}
If $\gA\subseteq\gC$ with $\gC$ \zedrz, the smallest \zed sub\ri of $\gC$ containing $\gA$ is equal to $\gA[(a\bul)_{a\in\gA}]$. More \gnlt if $\gA\subseteq\gB$ with $\gB$ reduced, and if every \elt of $\gA$ admits a quasi-inverse in $\gB$, then the sub\ri $\gA[(a\bul)_{a\in\gA}]$ of $\gB$ is \zedz.
In addition, every \elt of $\gA[(a\bul)_{a\in\gA}]$ is of the form 
$$\preskip.3em \postskip-.2em\ndsp 
\som_j{  a_jb_j\bul e_j}, \;\hbox{ with} 
$$
\begin{itemize}
\item the $e_j$'s are pairwise \orts \idms of $\gA[(a\bul)_{a\in\gA}]$,
\item  $a_j,b_j\in\gA$ and $b_j{b_j\bul} e_j=e_j$ for every $j$, 
\end{itemize}
such that \smash{$\big(\som_j{  a_jb_j\bul e_j}\big)\bul=\som_j{  a_j\bul b_j e_j}$}.
\end{lemma}
NB: Care will be taken, however, that we do not always have {\mathrigid 2mu $a_j{a_j\bul} e_j=e_j$}. We must therefore a priori replace $e_j$ with $e'_j=a_j{a_j\bul} e_j$ to obtain an expression of the same type as the previous one. We will also be able to note that every \idm of $\gA[(a\bul)_{a\in\gA}]$ is expressible in the form $e_c\prod_i(1-e_{d_i})$ for a~$c$ and some $d_i\in\gA$.
\begin{proof}
Among the \elts of $\gB$, those that are expressed as a sum of products~$ab\bul$ with $a$, $b\in\gA$ clearly form a sub\ri of $\gB$, which is therefore equal to $\gA[(a\bul)_{a\in\gA}]$.
Moreover, $ab\bul=ab\bul e_b$.
By considering the \agB generated by the $e_b$'s which intervene in a finite sum of the previous type, we deduce that every \elt of $\gA[(a\bul)_{a\in\gA}]$ can be expressed in the form
$$\preskip.3em \postskip-.30em\ndsp 
\som_j{ \, \left(\som_i a_{i,j}b_{i,j}\bul\right)\,e_j},\;\;\hbox{such that} 
$$
\begin{itemize}\itemsep1pt
\item the $e_j$'s are pairwise \ort \idms in $\gA[(a\bul)_{a\in\gA}]$,
\item  $a_{i,j},b_{i,j}\in\gA$, and $b_{i,j}b_{i,j}\bul e_j=e_j$, for all $i,j$.
\end{itemize}
Note that $b_{i,j}\bul$ is the inverse of  $b_{i,j}$ in $\gA[(a\bul)_{a\in\gA}][1/e_j]$, and 
 we can perform the computation
as for a usual sum of fractions $\sum_i a_{i,j}/b_{i,j}$.
For example to simplify a term with a sum of three \elts let us take 

\snic{  (a_1b_1\bul+a_2b_2\bul+a_3b_3\bul) e.}

Since $b_2b_2\bul e=b_3b_3\bul e=e$, we have $a_1b_1\bul e=a_1b_2b_3(b_1b_2b_3)\bul e$,
and

\snic{(a_1b_1\bul+a_2b_2\bul+a_3b_3\bul) e = (a_1b_2b_3+a_2b_1b_3+a_3b_1b_2)(b_1b_2b_3)\bul e = d c\bul e,}

which admits for quasi-inverse $c d\bul e$.
\end{proof}

Recall that $\gB\red$ designates the quotient of a \ri $\gB$ by its nilradical.

In the following lemma we observe what happens when we forcefully add a quasi-inverse to an \elt of a \riz. It is an operation neighboring localization, when we forcefully add an inverse of an \eltz, but slightly more delicate.

\begin{lemma}
\label{lemZedClotBasicStep} Let $\gA$ be a \ri and $a\in\gA$.
\begin{enumerate}
\item Consider the \ri $\aqo{\gA[T]}{aT^2-T,a^2T-a}=\gA[a\efl  ]$
and the canonical \homo $\lambda_a:\gA\to\gA[a\efl  ]$ ($a\efl$ designates the image of $T$). Then
for every \homo $\psi : \gA\to\gB$ such that $\psi(a)$ admits a quasi-inverse
there exists a unique \homo $\varphi:\gA[a\efl  ]\to \gB$
such that $\varphi\circ \lambda_a=\psi$.

\vspace{-1.6em}
\Pnv{\gA}{\lambda_a}{\psi}{\gA[a\efl  ]}{\varphi}{\gB}{ }{}{$\psi(a)$ admits a quasi-inverse}

\vspace{-1.2em}
\item  In addition, $aa\efl$ is an \idm and $\gA[a\efl  ]\simeq \gA[1/a] \times \aqo{\gA}{a}$.
\item  If $\gB$ is reduced we have a unique factorization via $(\gA[a\efl  ])\red$.
\end{enumerate}
In the rest of the text we denote by $\gA[a\bul]$ the \ri $(\gA[a\efl  ])\red$.
\label{NOTAAbul}
\begin{enumerate}\setcounter{enumi}{3}
\item
We have $\gA[a\bul] \simeq \Ared[1/a]\times \gA\sur{\DA(a)}$.
If $\gA$ is reduced the canonical \homo $\gA\to\gA[a\bul]$ is injective.
\item \label{ite5lemZedClotBasicStep}
$\Zar(\gA[a\bul]) =\Zar(\gA[a\efl])$ is identified with $(\ZarA)[\DA(a)\bul].$
\end{enumerate}
\end{lemma}

\begin{proof}
Left to the reader. The last item results from Lemma~\ref{lemRajouCompl}
and from Fact~\ref{fact2Zar}.
\end{proof}

\begin{corollary}\label{corlemZedClotBasicStep}
Let $a_1$, \ldots, $a_n\in\gA$.
\begin{enumerate}
\item The \ri $\gA[a_1\bul][a_2\bul]\cdots[a_n\bul]$ is independent, up to unique \isoz, of the order of the $a_i$'s.
It will be denoted by $\gA[a_1\bul,a_2\bul,\ldots,a_n\bul]$.
\item A possible description is the following
$$\preskip.4em \postskip.3em
\gA[a_1\bul,a_2\bul,\ldots,a_n\bul]
\simeq
\big(\gA[T_1,T_2,\ldots,T_n]\sur{\fa}\big)\red$$
with $\fa=\gen{(a_iT_i^2-T_i)_{i=1}^n,(T_ia_i^2-a_i)_{i=1}^n}$.
\item Another possible description is
$$\preskip.2em \postskip.0em
\gA[a_1\bul,a_2\bul,\ldots,a_n\bul]
\simeq
\prod\nolimits_{I\in\cP_n}{\left(\aqo{\gA}{(a_i)_{i\in I}}\right)\red[1/\alpha_I]}
$$
with $\alpha_I=\prod_{j\in \lrbn\setminus I}a_j$.
\end{enumerate}
\end{corollary}

\vspace{-.7em}
\pagebreak	
\begin{theorem}\label{thZedGen}
\emph{(Reduced \zede closure of a commutative \riz)}
For every \ri $\gA$ there exists a \zedr \ri $\Abul$ with a \homo $\lambda:\gA\to\Abul$, which uniquely factorizes every \homo $\psi :\gA\to\gB$ to a \zedr \riz.
This pair $(\Abul,\lambda)$ is unique up to unique \isoz.

\vspace{-1.7em}
\PNV{\gA}{\lambda}{\psi}{\Abul}{\varphi}{\gB}{commutative \risz}{}{\zedr \risz}

\vspace{-2.2em}
In addition\index{closure!reduced zero@\zedr ---}
\begin{enumerate}
\item [--] The natural \homo $\Ared\to \Abul$ is injective.
\item [--] We have $\Abul=\gA\red[(a\bul)_{a\in\Ared}]$.
\end{enumerate}
\end{theorem}
\begin{proof}
This is a corollary of the previous lemmas. We can suppose that $\gA$ is reduced. 
The uniqueness result (Corollary~\ref{corlemZedClotBasicStep}) allows for a construction of a colimit (which mimics a filtering union) based on the extensions of the type~\hbox{$\gA[a_1\bul,a_2\bul,\ldots,a_n\bul]$}, and \trf by Lemma~\ref{lem2SousZedRed}.
\end{proof}

\comms 1) A priori, since we are dealing with purely equational structures, the \uvle \zedr closure of a \ri exists and we could construct it as follows: 
we \fmt add the unary operation $a\mapsto a\bul$ and we force $a\bul$ to be a quasi-inverse of $a$. Our \dem has also allowed us to give a simplified precise description of the constructed object and to show the injectivity in the reduced case.

2) In \clamaz, the \zedr closure $\Abul$ of a \ri $\gA$ can be obtained as follows.
First of all we consider the product $\gB=\prod_\fp \Frac(\gA\sur{\fp})$, where $\fp$
ranges over all the \ideps of $\gA$. As $\gB$ is a product of fields, it is \zedrz.
Next we consider the smallest \zed sub\ri of $\gB$ containing the image of $\gA$ in $\gB$ under the natural diagonal \homoz.
\\
We then understand the importance of the earlier construction  of $\Abul$.
It allows us to have explicit access to something which looks like \gui{the set of all the} \ideps
of~$\gA$ (those of \clamaz) without needing to construct any one of them individually.
The assumption is that the classical mathematical reasoning that manipulates unspecified arbitrary prime ideals of the ring
 $\gA$ (\gnlt inaccessible objects) can be reread as arguments about the \ri $\Abul$: a mystery-free object!
\eoe

\medskip\rdb
\exls \label{exlCLZED}~
\\
1) Here is a description of the \zedr closure of $\ZZ$.\\
First of all, for $n\in\NN\etl$ the \ri $\ZZ[n\bul]$ is \isoc to $ \ZZ[1/n]\times \prod_{p|n}\FF_p$,
where $p$ indicates \gui{$p$ prime,}
and $\FF_p=\ZZ\sur{p\ZZ}$. Next, $\ZZ\bul$ is the colimit (that we can regard as a non-decreasing union) of the~$\ZZ[(n!)\bul]$'s.

2) Here is a description of the \zedr closure of $\ZZ[X]$.\\
First of all,  if $Q$ is a square-free \poluz, and if $n\in\NN\etl$ is a multiple of $\disc(Q)$,
the \ri $\ZZ[X][n\bul,Q\bul]$ is \isoc to
$$\preskip.4em \postskip.2em
\ZZ[X,1/n,1/Q]\times \prod_{p\divi n}\FF_p[X,1/Q]
\times \prod_{P\divi Q}\aqo{\ZZ[X,1/n]}{P} \times\!\!\!
\prod_{p\divi n,R\divi Q}\!\!\!\aqo{\FF_p[X]}{R}~~
$$
with $p$ standing for \gui{$p$ prime,} $P$ standing for \gui{\ird $P$ in $\ZZ[X]$,} and~$R\divi Q$ standing for \gui{\ird $R$ in $\FF_p[X]$ divides $Q$ in $\FF_p[X]$.} 
\\
Next, we pass to the colimit of the \ris 
$\ZZ[X][u_n\bul,Q_n\bul]$ (here, it is a non-decreasing union), where $Q_n$ 
is the squarefree part of the product of the
first $n$ \elts in an enumeration of squarefree \polus of $\ZZ[X]$, and where~\hbox{$u_n=n!\disc(Q_n)$}. 
\\
Note that we thus obtain a \ri by which all the natural \homos $\ZZ[X]\to\Frac(\ZZ[X]\sur\fp)$ are factorized for all the \ideps $\fp$ of $\ZZ[X]$: such a $\Frac(\ZZ[X]\sur\fp)$ is indeed either $\QQ(X)$, or some $\aqo{\QQ[X]}{P}$, or some $\FF_p(X)$, or some $\aqo{\FF_p[X]}{R}$.

3) The (\cot well-defined) \ri $\RR\bul$ is certainly one of the more intriguing 
objects in the world \gui{without LEM} that constitutes constructive mathematics.
Naturally, in \clamaz, $\RR$ is \zed and $\RR\bul=\RR$. 
\eoe

\begin{theorem}\label{thZedGenEtBoolGen}
For every \ri $\gA$ we have natural \isos 

\snic{\Bo(\ZarA)\simeq \BB(\Abul)\simeq\Zar(\Abul).}

\end{theorem}
\begin{proof}
This results from the last item of Lemma~\ref{lemZedClotBasicStep}, and from the fact that the two constructions can be seen as colimits of \gui{constructions at a step}
$\gE\leadsto\gE[a\bul]$ ($\gE$ is a \ri or a \trdiz).
\end{proof}

Note that if we adopted the notation $\gT\bul$ for $\Bo(\gT)$ we would have the pretty formula $(\ZarA)\bul\simeq\Zar(\Abul)$.

\begin{proposition}\label{propClZdrLoc}
Let $\gA$ be a \riz, $\fa$ be an \id and $S$ be a \moz.  
\begin{itemize}
\item The two \ris $(\gA\sur\fa)\bul$
and $\Abul\sur{\rD(\fa\Abul)}$ are canonically \isocz.
\item  The two \ris $(\gA_S)\bul$ and $(\Abul)_S$ are canonically \isocz.
\end{itemize}
 
\end{proposition}
%
\begin{proof} 
Note that $(\Abul)_S$ is \zedr as a \lon of a \zedr \riz. Similarly,
$\Abul\sur{\rD(\fa\Abul)}$ is \zedrz. Let us write the \dem for the \lonsz.
Consider the natural \homos 
$$
\preskip-.1em \postskip.3em 
\gA \to \gA_S \to (\gA_S)\bul  \quad\hbox {and}\quad
\gA  \to \Abul \to (\Abul)_S
. 
$$
The \homo $\gA \to \Abul$ uniquely \gui{extends} to a \homo $\gA_S \to (\Abul)_S$, and by the \uvl \prt of the \zedr closure, provides a unique morphism $(\gA_S)\bul \to (\Abul)_S$ which renders the ad hoc commutative diagram. Similarly, the morphism $\gA \to \gA_S$ gives birth to a unique morphism $\Abul \to (\gA_S)\bul$ which  extends to a morphism $(\Abul)_S\to (\gA_S)\bul$. By composing these two morphisms, by uniqueness, we obtain the identity twice.
\end{proof}
%

\section[Entailment relations]{
Entailment relations and \agHsz}
\label{secEntRelAgH}
\vspace{4pt}
\subsec{A new look at \trdisz}

A particularly important rule for distributive lattices, known as the
 \emph{cut}, is the following
\begin{equation}\label{coupure1}
 \bigl(\,x\vi a\; \leq\;  b\,\bigr)\quad\&\quad  \bigl(\,a\; \leq\; x\vu  b\,\bigr)
\quad \Longrightarrow \quad a \leq\;  b.
\end{equation}
To prove it we write $ x\vi a\vi b=x\vi a$  and $a= a\vi(x\vu b)$ so
$$ a=(a\vi x)\vu(a\vi b)=(a\vi x\vi b)\vu(a\vi b)=a\vi b
$$

\begin{notation}
\label{notaVupVda}
{\rm
For a \trdi $\gT$ we denote by $A \vda B$ or $A \vdash_\gT B$ the relation defined as follows over the set $\Pfe(\gT)$

\snic{A \vda B \; \; \equidef\; \; \Vi A\;\leq \;
\Vu B.}
}
\end{notation}
Note that the relation  $A \vdash B$ is well-defined over  $\Pfe(\gT)$ because the laws $\vi$  and $\vu$  are associative, commutative and idempotent.
Note that $\; \emptyset  \vda \{x\}$ implies $x=1 $ and that $\; \{y\} \vda \emptyset$ implies $y=0$.
This relation satisfies the following axioms, in which we write $x$ for $\{x\}$ and $A, B$  for $A\cup B$.
$$\arraycolsep3pt\begin{array}{rcrclll}
&    & a  &\vda& a    &\; &(R)     \\[1mm]
A \vda B &   \; \Longrightarrow \;  & A,A' &\vda& B,B'   &\; &(M)     \\[1mm]
(A,x \vda B)\;\;
\&
\;\;(A \vda B,x)  &  \; \Longrightarrow \; & A &\vda& B &\;
&(T).
\end{array}$$
\rdb
We say that the relation is \emph{reflexive}, \label{remotr} \emph{monotone} and
\emph{transitive}.
The third rule (transitivity) can be seen as a \gnn of the rule (\ref{coupure1}) and is \egmt called the \emph{cut} rule.
\index{cut}

Let us also quote the following so-called \emph{distributivity} rules:
$$\preskip.4em \postskip.4em
\arraycolsep3pt\begin{array}{rcl}
(A,\;x \vda B)\;\& \;(A,\;y \vda B)  &  \;  \Longleftrightarrow  \; &
A,\;x\vu y \vda B  \\[1mm]
(A\vda B,\;x )\;\&\;(A \vda B,\;y)  &  \; \Longleftrightarrow \; &
A\vda B,\;x\vi y
\end{array}$$

An interesting way to approach the question of \trdis defined by \gtrs and relations is to consider the relation $A \vda B$ defined over the set $\Pfe(\gT)$ of finitely enumerated subsets of a \trdiz~$\gT$.
Indeed, if $S\subseteq \gT$ generates $\gT$ as a lattice, then the knowledge of the relation~$\vda$ over $\Pfe(S)$ suffices to characterize without ambiguity the lattice~$\gT$, because every formula over $S$ can be rewritten, either in \gui{conjunctive normal form} (infimum of supremums in~$S$) or in \gui{disjunctive normal form} (supremum of infimums in~$S$). Therefore if we want to compare two \elts of the lattice generated by~$S$ we write the first in disjunctive normal form, the second in conjunctive normal form, and we observe that
$$ \Vu_{i\in I}\big(\Vi A_i \big)\; \leq \; \Vi_{j\in J}\big(\Vu B_j
\big)
\quad \Longleftrightarrow\quad  \Tt i\in I,\ \Tt j \in J,\;\;
A_i \vda  B_j
$$

\rdb
\begin{definition}
\label{defEntrel}
For an arbitrary set $S$, a relation over $\Pfe(S)$ which is reflexive, monotone and transitive is called an {\em \entrelz}.
\end{definition}

The following \tho is fundamental. It says that the three \prts of \entrels are exactly what is needed for the interpretation in the form of a \trdi to be adequate.

\begin{theorem}
\label{thEntRel1} {\rm  (Fundamental \tho of the \entrelsz)} \\
Let $S$ be a set with an \entrel $\vdash_S$ over $\Pfe(S)$. We consider the \trdi $\gT$ defined by \gtrs and relations as follows: the \gtrs are the \elts of $S$ and the relations are the
$$\preskip.4em \postskip.4em A\; \vdash_\gT \;  B
$$
each time that $A\; \vdash_S \; B$.  Then, for all $A$,
 $B$ in $\Pfe(S)$, we have
$$\preskip.4em \postskip.4em  A\; \vdash_\gT \;  B
\; \Longrightarrow \; A\; \vdash_S \;  B.
$$
\end{theorem}
\begin{proof}
We give an explicit description of the \trdi $\gT$. The \elts of $\gT$ are represented by those of $\Pfe\big(\Pfe(S)\big)$, \cad $X$'s of the form

\snic{X=\{A_1,\dots,A_n\}}

(intuitively $X$ represents $\Vi_{i\in\lrbn}\Vu {A_i}$).
We then inductively define the relation $A\preceq Y$ for $A\in \Pfe(S)$ and $Y\in \Pfe\big(\Pfe(S)\big)$ as follows
\begin{itemize}
\item If $B\in Y$ and $B\subseteq A$ then $A\preceq Y$.
\item If we have $A\vdash_S y_1,\dots,y_m$ and $A,y_j\preceq Y$ for $j=1$, \ldots, $m$ then $A\preceq Y$.
\end{itemize}
We easily show that if $A\preceq Y$ and $A\subseteq A'$ then we have $A'\preceq Y.$ We deduce that $A\preceq Z$ if $A\preceq Y$ and $B\preceq Z$ for all $B\in Y$. We can then define $X\leq Y$ by
\gui{$A\preceq Y$ for every $A\in X$.} We finally verify that $\gT$ is a distributive lattice{\footnote{More \prmtz, as $\preceq$ is only a preorder, we take for $\gT$ the quotient of $\Pfe\big(\Pfe(S)\big)$ with respect to the \eqvc relation: $X\preceq Y$ and $Y\preceq X$.}}
with respect to the operations ($0$-aries and binaries)
\begin{equation}\label{eqentrel}
\left.\arraycolsep2pt
\begin{array}{rclcrcl}
1  & =  &  \emptyset  &\qquad\quad   &  0  &  = & \so{\emptyset}    \\[1mm]
X\vi Y  &  = & X\cup Y && X\vu Y  &  = & \sotq{ A \cup B} {A\in X,~B\in Y}
\end{array}
\right|
\end{equation}
For this we show that if $C\preceq X$ and $C\preceq Y$, then we have $C\preceq X\vi Y$ by induction on the proofs of $C\preceq X$ and $C\preceq Y$.
\\
We notice that if $A\vdash_S y_1,\dots,y_m$ and $A,y_j\vdash_S B$ for all $j$, then we obtain $A\vdash_S B$ by using the cut rule $m$ times. From this, it results that if we have $A\vdash_\gT B$,
that is that $A\preceq \{\{b\}~|~b\in B\}$, then we have $A\vdash_S B$.
\end{proof}

\begin{corollary}\label{corthEntRel1} \emph{(\Fp distributive lattice)}
\begin{enumerate}
\item A \trdi freely generated by a finite set $\,E\,$ is finite.
\item A \pf \trdi is finite.
\end{enumerate}
\end{corollary}
\begin{proof}
\emph{1.} Consider the minimal \entrel on $E$. It is defined by 
$$\preskip.2em \postskip.2em
A\vdash_E B \;\equidef  \;\exists x\in A\cap B.
$$
We then consider the \trdi corresponding to this \entrel via \thref{thEntRel1}.
It is \isoc to a subset 
of $\Pfe\big(\Pfe(E)\big)$, the one which is represented by the lists $(A_1,\ldots,A_k)$ in $\Pfe(E)$ such that two $A_i$ with distinct indices are incomparable with respect to the inclusion. The laws are obtained from~(\ref{eqentrel}), by simplifying the obtained lists when they do not satisfy the criteria of incompatibility.

\emph{2.} If we impose a finite number of relations between the \elts of $E$, we have to pass to a quotient lattice of the free \trdi over $E$. The \eqvc relation generated by these relations and compatible with the lattice laws is decidable because the structure is defined by only using a finite number of axioms.
\end{proof}\rdb

\rems \label{remagBfree}
1) Another proof of item \emph{1} could be the following.
The \agB freely generated by the \trdi $\gT$ freely generated by $E$ is the \agB $\gB$ freely generated by $E$. The latter can easily be described by the \elts of $\Pfe\big(\Pfe(E)\big)$, without any passage to the quotient:
the subset $\{A_1,\dots,A_n\}$ intuitively represents 
$\Vu_{i\in\lrbn}\left(\Vi {A_i}\vi \Vi {A'_i}\right)$, by designating by $A'_i$ the subset of $E$ formed by the $\lnot x$'s for the $x\notin A_i$. Therefore $\gB$ has $2^{2^{\#E}}$ \eltsz. Finally, we have seen that $\gT$ is
identified with a  \strdi of $\gB$ (\thref{thBoolGen}).
\\
2) The given proof of item~\emph{2} uses an altogether \gnl argument. In the case of \trdis we can more \prmt refer to the description of the quotients given on \paref{eqPreceq}.
\eoe

\subsect{Duality between finite \trdis and finite ordered sets}{Duality}\label{SpecTrdiFi}

If $\gT$ is a \trdi let
$\SpecT\eqdefi\Hom(\gT,\Deux)$. It is an ordered set called the \emph{(Zariski) spectrum of $\gT$}.
An \elt $\alpha$ of $\SpecT$ is \care by its kernel. In \clama such a kernel is called a \idepz.
From the \cof point of view it must be detachable.
Here we are interested in the case where $\gT$ is finite, which implies that~$\SpecT$ is \egmt finite (in the \cof sense).%
\index{spectrum!of a \trdiz}%
\index{ideal!prime --- (\trdiz)}%
\index{prime!ideal of a distributive lattice}

If $\varphi:\gT\to\gT'$ is a \homo of \trdis  and if $\alpha\in\Spec\gT'$,
then  $\alpha\circ \varphi\in\SpecT$.
This defines a non-decreasing map from~$\Spec\gT'$ to~$\SpecT$, 
denoted by $\Spec\alpha$, called the \gui{dual} of~$\varphi$.

Conversely, let $E$ be a finite ordered set. Let $E\sta$ be the set of \emph{initial sections} of $E$, \cad the set of finite subsets of $E$ that are stable under the operation $x\mapsto \dar x$. This set, ordered by the relation $\supseteq$, is a finite \trdiz,  a sublattice of the lattice $\Pf(E)\eci$ (the opposite lattice of $\Pf(E)$).

\begin{fact}\label{factTrdiFiniDual}
The number of \elts of a finite ordered set $E$ is equal to the maximum length of a strictly increasing chain of \elts of~$E\sta$.
\end{fact}
\begin{proof} It is clear that a strictly monotone chain of \elts of $E\sta$ (therefore of finite subsets of $E$) cannot have more than $1+\# E$ \eltsz. Its \gui{length} is therefore $\leq \#E$. Regarding the 
reverse in\egtz, we verify it for $E=\emptyset$ (or for a singleton), then we perform an \recu on $\#E$, by regarding an ordered set with $n$ \elts ($n\geq1$) as an ordered set \hbox{with $n-1$} \elts that we extend by adding a maximal \eltz.
\end{proof}

If $\psi:E\to E_1$ is a non-decreasing map between finite ordered sets, then for every $X\in E_1\sta$, $\psi^{-1}(X)$ is an \elt of $E\sta$. This defines a \homo $E_1\sta\to E\sta$ denoted by $\psi\sta$, called the \gui{dual} of $\psi$.

\begin{theorem}\label{thDualiteFinie} \emph{(Duality between finite \trdis and finite ordered sets)}

\vspace{-.5em}
\begin{enumerate}
\item For every finite ordered set $E$ let us define $\nu_E:E\to\Spec(E\sta)$ by
$$\preskip.4em \postskip.4em 
\nu_E(x)(S)=0 \;\;\hbox{if}\;\; x\in S,\;\;1\;\;\hbox{otherwise}. 
$$
Then, $\nu_E$ is an \iso of ordered sets. In addition, for every non-decreasing map $\psi:E\to E_1$, we have $\nu_{E_1}\circ \psi=\Spec(\psi\sta)\circ \nu_E.$
\item For every finite \trdi $\gT$ let us define $\iota_\gT:\gT\to(\SpecT)\sta$ by
$$\preskip-.6em \postskip.4em 
\iota_\gT(x)=\sotq{\alpha \in \SpecT}{\alpha(x)=0}.
$$
Then, $\iota_\gT$ is an \iso of \trdisz. In addition, for every morphism $\varphi:\gT\to \gT'$, we have ${\iota_{\gT'}\circ \varphi=(\Spec\varphi)\sta\circ \iota_\gT.}$

\end{enumerate}
\end{theorem}
\begin{proof}
See Exercise~\ref{exoTreillisDistributifFini}.
\end{proof}

In other terms, the categories of finite \trdis and of finite ordered sets are anti\eqvesz. The anti\eqvc is given by the contravariant functors $\Spec\bullet$ and $\bullet\sta$, and by the natural transformations~$\nu$ and~$\iota$ defined above.
The \gnn of this anti\eqvc of categories to the case of not necessarily finite distributive lattices will briefly be addressed on \paref{SpecTrdi}.

\subsec{Heyting \algsz}
\label{secagH}

A \trdi $\gT$ is called an \ix{implicative lattice} or a {\em \agHz}  when there exists a binary operation $\im$ satisfying for all $a,\,b,\,c$,
\index{algebra!Heyting ---}\index{Heyting!algebra}
\begin{equation}\label{eqAgHey}\preskip.4em \postskip.4em
a\vi b \leq c \;\;\Longleftrightarrow \;\; a \leq  (b\im c)\,.
\end{equation}
This means that for all $b$, $c\in\gT$, the \emph{conductor \idz}%
\index{ideal!conductor --- (in a \trdiz)}\label{NOTATransp2}%
\index{conductor!of one ideal into another (\trdiz)}
$$\preskip.2em \postskip.4em
(c:b)_\gT \eqdefi \sotq{x\in\gT}{x\vi b\leq c}
$$
is principal, its \gtr being denoted by $b\im c$.
Therefore if it exists, the operation~$\im$ is uniquely determined by the structure of the lattice.
We then define the unary law $\neg x = x\im 0$.
The structure of a \agH can be defined as purely equational by giving suitable axioms, described in the following fact.

\begin{fact}\label{fact1AGH}
A lattice $\gT$ (not assumed distributive) equipped with a law $\im$ is a \agH \ssi the following axioms are satisfied
$$
\preskip.4em \postskip.4em
\arraycolsep2pt\begin{array}{rcl}
a\im a&=   &1,    \\[.2em]
a\vi(a\im b)&=   &a\vi b,    \\[.2em]
b\vi(a\im b)&=   & b,   \\[.2em]
a\im(b\vi c)&=   &(a\im b)\vi(a\im c).
\end{array}$$
\end{fact}

Let us also note the following important facts.

\begin{fact}\label{fact2AGH} In a \agH we have
$$\preskip.4em \postskip.4em
\arraycolsep2pt\begin{array}{rclcrcl}
a\leq b &\Leftrightarrow&  a\im b =1,\\[.2em]
a\im (b\im c) &=& (a\vi b) \im c, &\phantom{aa}&
         a\im b &\leq& \neg b\im \neg a,
\\[.2em]
(a\vu b)\im c &=& (a\im c)\vi(b\im c), &&
 a&\leq    & \neg\neg a,
         \\[.2em]
\neg\neg\neg a&=    &\neg a, &&
a\im b &\leq& (b\im c) \im (a\im c),~~
             \\[.2em]
\neg(a\vu b)&=   & \neg a\vi \neg b,&&
\neg a\vu b&\leq    & a\im b.
\end{array}$$
\end{fact}

Every finite \trdi is a \agHz, because every \itf is principal. A special important case of \agH is a \agBz.

A \emph{\homo of \agHsz} is a \homo of \trdis $\varphi :\gT\to\gT'$
 such that $\varphi(a\im b)=\varphi(a)\im\varphi(b)$ for all $a$, $b\in\gT$.

The following fact is immediate.

\begin{fact}
\label{factQuoAgH}
Let $\varphi:\gT\to\gT_1$ be a \homo of \trdisz, with~$\gT$ and $\gT_1$ being \agHsz. Let $a\preceq b$ for $\varphi(a)\leq_{\gT_1}\varphi(b)$. Then $\varphi$ is a \homo of \agHs \ssi we have for all $a$, $a'$, $b$, $b'\in\gT$

\snic{
a\preceq a'\Longrightarrow (a'\im b)\preceq(a\im b) , \quad
\hbox{and}\quad
b\preceq b'\Longrightarrow (a\im b)\preceq(a\im b').
}
\end{fact}

\begin{fact}
\label{factQuoAgH2}
If $\gT$ is a \agH every quotient $\gT/(y=0)$ (\cad every quotient by a principal \idz) is also a \agHz.
\end{fact}
\begin{proof}
Let $\pi:\gT\to\gT'=\gT/(y=0)$ be the canonical projection. We have
\[ \arraycolsep4pt
\begin{array}{ccccccccccccc} 
\pi(x)\vi\pi(a)\,\leq_{\gT'}\, \pi(b)  &  \Longleftrightarrow 
&  \pi(x \vi a)\,\leq_{\gT'}\, \pi(b)  &  \Longleftrightarrow  \\[1mm]
  x\vi a \,\leq\, b\vu
y  & \Longleftrightarrow   &  
x\,\leq\, a\im(b\vu y).  
 \end{array}
\]
 However, $y\,\leq\, b\vu y\,\leq\,a\im(b\vu y)$, 
 therefore 
 
 \snic{\pi(x)\vi\pi(a)\,\leq_{\gT'}\,\pi(b)\;\Longleftrightarrow\;
x\,\leq\, \big(a\im(b\vu y)\big)\vu y,}
 
\cad $\pi(x)\leq_{\gT'}\pi\big(a\im(b\vu y)\big)$, which shows that $\pi\big(a\im(b\vu y)\big)$ holds for $\pi(a)\im\pi(b)$ in $\gT'$.
\end{proof}

\rems
1) The notion of a \agH is reminiscent of the notion of a \cori in commutative \algz. Indeed, a \cori can be \care as follows: the intersection of two \itfs is a \itf and the conductor of a \itf into a \itf is a \itfz.  If we \gui{reread} this for a \trdi by recalling that every \itf is principal we obtain a \agHz. 

 2) Every \trdi $\gT$ generates a \agH naturally. In other words we can formally add a \gtr for every \id $(b:c)$.  But if we start from a \trdi which happens to be a \agHz, the \agH which it generates is strictly greater. Let us take for example the lattice $\Trois$ which is the free \trdi with a single \gtrz. The \agH that it generates is therefore the free \agH with one \gtrz.  But it is infinite (cf.  \cite{Johnstone}).  A~contrario the Boolean lattice generated by $\gT$ (cf.\  \thref{thBoolGen}) remains equal to~$\gT$ when it is Boolean. \eoe

\Exercices{

\begin{exercise}
\label{exoTrdiLecteur}
{\rm  We recommend that the \dems which are not given, or are sketched, or
left to the reader,
etc, be done.
But in particular, we will cover the following cases. \perso{\`a compl\'eter}
\begin{itemize}\itemsep0pt
\item \label{exotreillisquotient}
Show that the relations \pref{eqPreceq} on \paref{eqPreceq} are exactly what is needed to define a quotient lattice.
\item \label{exopropIdealFiltre} Prove Proposition~\ref{propIdealFiltre}.
\item \label{exocorlemRajouCompl}
Prove Corollary~\ref{corlemRajouCompl}.
\item
Prove Facts~\ref{factBezGCD}, \ref{factLocaliseGCD},
\ref{factSouspgcdsat} and \ref{factAXiclgcd}.
\item Prove Fact~\ref{fact1Zar} and all the numbered facts between \ref{fact2Zar} and \ref{factQuoFIID} (for Fact~\ref{factZar} see Exercise~\ref{exoTreillisZariski}).
\item \label{exoexlCLZED} Prove what is affirmed in the examples on \paref{exlCLZED}.
\item \label{exofactAGH} Prove Facts~\ref{fact1AGH} and~\ref{fact2AGH}.
\end{itemize}
}
\end{exercise}

\vspace{-1em}
\begin{exercise}
\label{exoRegleCoupure}
{\rm
Let $\gT$ be a \trdi and $x\in\gT$. We have seen (Lemma~\ref{lemRajouCompl}) that
$$\preskip-.4em \postskip.4em 
\lambda_x:\gT\to\gT[x\bul] \eqdefi\gT\sur{(x=0)}\times \gT\sur{(x=1)} 
$$
is injective, which means: if $y\vi x=z\vi x$ and $y\vu x=z\vu x$, 
then $y =z$.
\\  
Show that we can deduce the cut rule (\ref{coupure1}).
}
\end{exercise}

\vspace{-1em}
\begin{exercise}
\label{exo0GpRtcl}
{\rm  Let $\gA$ be an integral \ri and $p$, $a$, $b \in \Reg(\gA)$, with $p$ irreducible. Suppose that $p \divi ab$, but $p \nedivi a$, $p \nedivi b$. Show that $(pa, ab)$ does not have a gcd. Show that in $\ZZ[X^2,X^3]$ the \elts $X^2$ and $X^3$ admit a gcd, but no lcm, and that the \elts $X^5$ and $X^6$ do not have a gcd.
}
\end{exercise}

\vspace{-1em}
\begin{exercise}
\label{exo2GpRtcl} (Another \dfn of \grlsz)\\
{\rm   Show that the axioms that must satisfy a subset $G^+$ of a group $(G,0,+,-)$ to define a compatible lattice-order are
\begin{itemize}
\item $G=G^+-G^+$,
\item $G^+\cap -G^+=\so{0}$,
\item $G^+ + G^+\subseteq G^+$,
\item $\forall a,b\;\exists c,\;\;c+G^+=(a+G^+)\cap (b+G^+)$.
\end{itemize}
}
\end{exercise}

\vspace{-1em}
\begin{exercise}
\label{exoLemmeGaussAX}
{(Another proof of Gauss' lemma)}\\
{\rm
In the context of Proposition~\ref{propLG}, show that $\G(fg)=\G(f)\G(g)$ with the help of a \dem based on the \DKM lemma  \ref{lemdArtin}.
}
\end{exercise}

\vspace{-1em}
\begin{exercise}\label{exoKroneckerTrick}
{(\KRAz's trick)}
{\rm  
Let $d$ be a fixed integer $\ge 2$.
\\
\emph {1.}
Let $\AuX_{<d} \subset \AuX = \gA[\Xn]$ be the \Asub constituted from \pols $P$ such that $\deg_{X_i} P < d$ for every $i \in \lrbn$, and $\gA[T]_{<d^n}\subset \gA[T]$ be the one formed from the \pols $f \in \gA[T]$ of degree $< d^n$.
\\
Show that $\varphi : P(\Xn) \mapsto P(T, T^d, \ldots, T^{d^{n-1}})$ induces an \iso of \Amos between the \Amos $\Aux_{<d}$ and $\gA[T]_{<d^n}$.

\emph {2.} 
We assume that $\AX$ is a UFD.
Let $P \in \AuX_{<d}$ and $f = \varphi(P) \in \AT_{d^n}$. Show that any factorization of $P$ in $\AuX$ can be found by a finite procedure from those of $f$ in $\AT$.
}
\end{exercise}

\vspace{-1em}
\begin{exercise}
\label{exoTreillisZariski}
{\rm
Verify Fact~\ref{factZar}, \cad  $\ZarA$ is a \trdiz.
Show that this \trdi can be defined by \gtrs and relations as follows. The \gtrs are the symbols $\rD(a)$, $a \in \gA$, with the \sys of relations:

\snic {
\rD(0) = 0,\quad  \rD(1) = 1,\quad
\rD(a + b) \le \rD(a) \vu \rD(b),  \quad \rD(ab) = \rD(a) \vi \rD(b).
}
}
\end{exercise}

\vspace{-1em}
\begin{exercise}\label{exoClosedCover1}
{\rm
The context is that of the closed covering principle~\ref{prcf2}. We consider a closed covering of the \ri $\gA$ by \ids \hbox{$\fa_1$, \ldots, $\fa_r$}.  We do not suppose that $\prod_i\fa_i=0$, but we suppose that each $\fa_i$ is \tfz.  Show that an \Amo $M$ is \tf \ssi it is \tf modulo each~$\fa_i$.
}
\end{exercise}

\vspace{-1em}
\begin{exercise}\label{exoAbul}
 (The \ri $\Abul$) {\rm We are in the context of \clamaz.
\\
Let $\gA$ be a \ri and $\varphi:\gA\to\Abul$ be the natural \homoz.

\emph{1.} Show that the map $\Spec\varphi:\Spec\Abul\to\Spec\gA$ is a bijection and that for $\fq\in\Spec\Abul$, the natural \homo $\Frac(\gA\sur{\varphi^{-1}(\fq)})\to\Abul\sur\fq$ is an \isoz.

\emph{2.} The \ri $\Abul$ is identified with the \zedr sub\ri of 

\snic{\wi\gA\eqdefi\prod_{\fp\in\SpecA}\Frac(\gA\sur\fp)}

generated by (the image of) $\gA$.
}

\end{exercise}

\vspace{-1em}
\begin{exercise}\label{exoMinA} (Minimal prime ideals)
\\
{\rm We are in the context of \clamaz. A \idep is said to be minimal if it is minimal among the \idepsz. Let $\Min\gA$ be the subspace of $\SpecA$ formed by the \idemisz. 
Recall that we have defined a \emph{maximal filter} as a filter whose localized \ri is a \zedr \aloz. 
In item \emph{1} of this exercise we make the link with the most usual \dfnz.

\emph{1.}
Show that a strict filter $\ff$ is maximal among the strict filters \ssi for every $x\notin \ff$ there exists an $a\in\ff$ such that $ax$ is nilpotent.
Another possible \carn is that the localized \ri $\ff^{-1}\gA$ is local, \zed and nontrivial. 
In particular, every strict maximal filter among the strict filters is prime. 

NB: reformulation of the first \cara \prt for the \cop \idepz: a \idep $\fp$ is minimal \ssi for all~\hbox{$x\in \fp$}, there exists an $a\notin\fp$ such that $ax$ is nilpotent.

\emph{2.} The dual notion of the Jacobson radical is the intersection filter of the maximal filters (\cad the complement of the union of the minimal \idepsz). It can be \care as follows in \clama (compare with Lemma~\ref{lemcRadJ} and its proof): it is the set of~\hbox{$a\in\gA$}  \gui{nilregular} in the following sense
\begin{equation}\label{Eqnilreg}
\forall y \in\gA\quad\;ay  \;\mathrm{nilpotent}\;\Rightarrow \; y
\;\mathrm{nilpotent}.
\end{equation}
In particular, in a reduced \riz, it is the set of \ndz \eltsz.
 
}
\end{exercise}

\vspace{-1em}
\begin{exercise}\label{exoFreeBooleAlgebra} 
{(Boolean \alg freely generated by a finite set)}
\\
{\rm
Let $E = \{x_1, \ldots, x_n\}$ be a finite set.

\emph{1.}
Show that the \agB $\gB$ freely generated by $E$ identifies with the \alg 
\[\FF_2[X_1, \ldots, X_n]\sur\fa = \FF_2[x_1, \ldots, x_n]\] with $\fa = \gen {(X_i^2 - X_i)_{i=1}^n}$.

\emph {2.}
Define two \gui{natural} $\FF_2$-bases of $\gB$, indexed by $\Pf(E)$, one being monomial and the other being a \sfioz.
Express one in terms of the other.
}
\end{exercise}

\vspace{-1em}
\begin{exercise}
\label{exotrdifree}
{\rm Give a precise description of \trdis freely generated by sets with $0$, $1$, $2$ and $3$ \eltsz. In particular, specify the number of their \eltsz.
}
\end{exercise}

\vspace{-1em}
\begin{exercise}
\label{exoTreillisDistributifFini}
{\rm
We detail the \dem of \thref{thDualiteFinie}.

\emph{1.} We use (as in the course) the order relation $\supseteq$ over $E\sta$ 
(the set of initial sections of the finite ordered set $E$).\\
If $S_1, S_2 \in E\sta$, what are $S_1 \vi S_2$, $S_1 \vu S_2$ equal to?

\emph{2.} What is the order relation over the set of \ideps of $\gT$ corresponding to the order which has been defined for $\SpecT$?

\emph{3.}
Prove item~\emph{1} of the \thoz.  \\
We will start by verifying that for $S \in E\sta$, $S$ generates a \idep \ssi $S$ is of the form $\dar x$ with $x \in E$; then that $\Ker \nu_E(x) = \cI_{E\sta}(\dar x)$.

\emph{4.}  How to construct $E\sta$ from $E$?
Treat the following example

\snic {
E = \vcenter {\xymatrix @R = 5pt @C = 5pt{
   & &&d \\
c  &     & b\ar@{->}[ur] \\
   & a\ar@{->}[lu]\ar@{->}[ru] \\
}}
\qquad
a < b<d,\  a < c.
}

Study the case where $E$ is totally ordered, and the case where $E$ is ordered by the \egt relation.

\emph{5.} Prove item \emph{2} of the \thoz.

\emph{6.}
Consider the same questions for the opposite order over $E\sta$ and  adapt the order over $\Spec(E\sta)$.
}
\end{exercise}

\vspace{-1em}
\begin{exercise}
\label{exoIVpgcd} {\rm
Let $a$, $b$ be nonzero in an integral \riz. 
Suppose that the \idz~\hbox{$\gen{a,b}$} is \iv and that $a$ and $b$ admit a lcm $m$.\\
Show that $\gen{a,b}$ is a \idpz.}
\end{exercise}

\vspace{-1em}
\begin{exercise}
\label{exoFactFini} {(A UFD with only a finite number of \ird \eltsz)}\\ {\rm
Show that a UFD with only a finite number of \ird \elts is a PID. }
\end{exercise}

\vspace{-1em}
\begin{exercise}
\label{exoPrincipalIntersecSoucorps} {(An interesting intersection)}\\
 {\rm Let $\gk$ be a \cdiz. We consider the intersection 
$$\gA = \gk(x,y)[z] \cap \gk(z,x+yz).$$
They are two sub\ris of $\gk(x,y,z)$.
The first is a PID, the second is a \cdiz.
Show that $\gA = \gk[z,x+yz]$, \isoc to $\gk[z,u]$.
Thus the intersection is not a PID, not even a Bézout \riz.}
\end{exercise}


\vspace{-1em}
\begin{problem}
\label{exoGpRtclQuotient} (Quotient lattice-groups, solid subgroups)\\
{\rm
In an ordered set $E$, if $a\leq b$, we call 
\emph{the segment with endpoints $a$ and $b$} the subset $\sotq{x\in E}{a\leq x\leq b}$. We denote it by $[a,b]_E$ or $[a,b]$. 
A subset $F$ of $E$ is said to be \ixc{convex}{subset of an ordered set} when the implication $a$, $b\in F\Rightarrow[a,b]\subseteq F$ is satisfied.
\\
A subgroup $H$ of an \grl is said to be \ixc{solid}{subgroup} if it is a convex \sgrlz.\index{solid subgroup!of an \grlz} 
We will see that this notion is the analogue for \grls of that of an ideal for rings. 

 \emph{1.} A subgroup $H$ of an ordered group $G$ is convex \ssi the order relation over $G$ passes to the quotient in $G/H$, \cad more \prmt $G/H$ is equipped with an ordered group structure for which $(G/H)^{+}=G^{+}+H$.
  We also say \emph{isolated subgroup} for \gui{convex subgroup of an ordered group.}\index{isolated subgroup!of an ordered group}\index{isolated!subgroup}

 \emph{2.}  
The kernel $H$ of a morphism of \grls $G\to G'$ is a solid subgroup of $G$.

 \emph{3.} 
Conversely, if $H$ is a solid subgroup of an \grl $G$, the law $\vi$ passes to the quotient, it defines an \grl structure over $G/H$, and the canonical surjection from~$G$ to~$G/H$ is a morphism of \grls which factorizes every morphism of source $G$ which vanishes over $H$. 

 \emph{4.}  
We have defined in \ref{defiCongru} the \sgrl $\cC(x)$.\\  
Show that $\cC(x)\cap\cC(y)=\cC(\abs{x}\vi\abs{y})$, 
and that the solid subgroup generated by $x_1$, $\ldots$, $x_n\in G$ is equal to $\cC(\abs{x_1}+\cdots+\abs{x_n})$. 
In particular, the set of solid \emph{principal} subgroups, 
\cad of the form $\cC(a)$, is \gui{almost} a \trdi (in \gnl a maximum \elt is missing).

}
\end{problem}

\vspace{-1em}
\begin{problem}
\label{exoSgpPolaire} (Polar subgroups, \ort direct summands)\\
{\rm
 \emph{1.} 
 If $A$ is an arbitrary subset of an \grl $G$ let

\snic{A\epr:=\sotq{x\in G\,}{\,\forall a\in A,\;\abs{x}\perp \abs{a}}.}

Show that $A\epr$ is always a solid subgroup.
\\
Show that, as usual in this type of situation, we have

\snic{ A\subseteq(A\epr)\epr, \; (A\cup B)\epr=A\epr\cap B\epr, \; A\subseteq B\Rightarrow B\epr\subseteq
A\epr\;\hbox{ and }\;{{{A\epr}\epr}\epr}=A\epr.}

\emph{2.} 
A solid subgroup $H$ of an \grl is called a \ix{polar subgroup} when ${{H\epr}\epr}=H$.
We also say \emph{a polar} instead of \gui{a polar subgroup.}
\\
   A subgroup $H$ is said to be an \emph{\ort direct summand} 
   when  $G=H\oplus H\epr$ (direct sum of subgroups in an Abelian group), in which case $G$ is naturally \isoc to $H\boxplus H\epr$. 
   We also say that $G$ is the \ix{internal \ort direct sum} of $H$ and $H\epr$ and let (by abuse) $G=H  \boxplus H\epr$.
\\
Show that an \ort direct summand is always a polar subgroup.
\\
Show that if $G=H\boxplus K$ (with $H$ and $K$ identified with subgroups of $G$) and if~$L$ is a solid subgroup, then $L=(L\cap H)\boxplus(L\cap K)$.

\emph{3.} 
\Gnltz, we say that $G$ is the \ixx{internal \ort direct sum} {of a family of \sgrlsz}  $(H_i)_{i\in I}$, indexed by a discrete set $I$, 
when we have $G=\sum_{i\in I}H_i$ and when the $H_i$'s are pairwise \ortsz.
In this case, each~$H_i$ is a polar subgroup of $G$ and we have a natural \iso of \grls $\boxplus_{i\in I}H_i\simeq G$.
We write (by abuse)  $G=\boxplus_{i\in I}H_i$.
\\
Suppose that an \grl is an \ort direct sum of a family of polar subgroups 
$(H_i)_{i\in I}$, as well as of another family $(K_j)_{j\in J}$. 
Show that these two \dcns admit a common refinement.
\\
Deduce that if the components of a \dcn as an \ort direct sum are nontrivial \emph{in\dcpsz} subgroups, that is, 
which do not admit a strict \ort direct summand,
then the \dcn is unique, up to bijection of the set of indices.
}
\end{problem}

\vspace{-1em}
\begin{problem}\label{exoAutourGaussJoyal} 
        {(Revisiting Gauss-Joyal)\iJG}\\
{\rm
Let $u : \gA \to \gT$ ($\gA$ is a commutative \riz, $\gT$ a \trdiz) satisfying

\snic {
u(ab) = u(a) \vi u(b),\quad u(1) = 1_\gT, \quad
u(0) = 0_\gT, \quad u(a+b) \le u(a) \vu u(b).
}

For $f = \sum_i a_i X^i \in \gA[X]$, we let

\centerline {$u(f) = u\big(\rc(f)\big) \eqdefi \Vu_{i} u(a_i).$}

\emph{1.} 
 Prove that \gui{it is well-defined,} \cad that $u(f)$ 
 only depends on $\rc(f)$.

\smallskip  We want to \gui{directly} prove (in particular, without using Lemma~\ref{lemGaussJoyal}), the following version of the Gauss-Joyal lemma

\smallskip \centerline{{\sf LGJ:} $\quad u(fg) = u(f)
\vi u(g)$.}

\smallskip 
\emph{2.} 
Verify that if $g = \sum b_jX^j\in \gA[X]$ the result is equivalent to $u(a_ib_j) \le u(fg)$.

\emph{3.} 
What does {\sf LGJ} say if $\gT = \so{\Vrai,\Faux}$ and $u(a) = (a \ne 0)$?

\emph{4.} 
Taking inspiration from the classical proof of the result of the previous question, prove~{\sf LGJ}.

\emph{5.} 
What does {\sf LGJ} say if $\gT = \Zar\gA$ and $u(a) = \rD_\gA(a)$?
}
\end{problem}

\vspace{-1em}
\begin{problem}
\label{exoQiClot} (\qiri closure of a commutative \riz)%
\index{closure!pp-ring ---}
 \\
{\rm Taking inspiration from the \zedr closure, give a construction of the \qiri closure $\Aqi$ of an arbitrary commutative \ri $\gA$.
\\
The following universal \pb needs to be solved: 
\Pnv{\gA}{\lambda}{\psi}{\Aqi }{\varphi}{\gB}{~}{\homos of \risz}{morphisms of \qiris}

\vspace{-3mm}where \emph{the \qiri morphisms are the \ri \homos which respect the law $a\mapsto e_a$} ($e_a$ is the \idm satisfying $\gen{1-e_a}=\Ann(a)$). 
Hereinafter, we will speak of \emph{\qiri morphism}.
\\
A \qiri closure of a \ri $\gA$ \gui{a priori} exists, from the simple fact that the theory of \qiris is purely equational. Indeed, for any \sys of \gtrs and of relations (a relation is an \egt between two terms constructed from \gtrsz, of $0$ and of $1$, by using the laws $+,-,\times,  a\mapsto e_a$), there exists some \qiri \gui{the most \gnl as possible} corresponding to this \pnz:
we take over the set of terms the smallest \eqvc relation which respects the axioms and which places in the same \eqvc class two terms related by a given relation at the start.
Under these conditions the \ri $\Aqi$ is simply the \qiri generated by the \elts of $\gA$ with for relations all the \egts $a+b=c$, $a\times b=d$, $a=-a'$ true in $\gA$. 
\\
But we want a precise description, as for the \zedr closure.
\\
We will then prove the following results.
\index{morphism!of \qirisz} 

\emph{1.} \emph{(\qiri morphisms)}
\begin{itemize}
\item [\emph{a.}] A morphism $\varphi: \gA\to\gB$ 
                         is a \qiri morphism  
\ssi it transforms every \ndz \elt into
a \ndz \eltz.  
In this case, it uniquely extends to a morphism $\Frac(\varphi) : \Frac(\gA) \to \Frac(\gB)$.

\item [\emph{b.}] A \qiri morphism is injective \ssi its restriction to $\BB(\gA)$ is injective.
\item [\emph{c.}] There exist injective \homos between \qiris that are not
\qiri morphisms.
\item [\emph{d.}] Every \homo between \zedr \ris is a
\qiri morphism.
\item [\emph{e.}]  If $\gA$ is a \qiriz, $\BB (\Frac\gA)$ is identified with $\BB(\gA)$ and the injection $\gA\to\Frac(\gA)$ is a \qiri morphism.
\end{itemize}

\emph{2.} 
 We have natural \ri \homos $\Ared\to\Aqi\to\Frac(\Aqi)\to\Abul$.
 \\
They are all injective and the natural \homo $\Frac(\Aqi)\to\Abul$ is an \isoz.

\emph{3.} If $\gA\subseteq\gC$ with $\gC$ a \qiriz, the smallest \sqiri of $\gC$ containing $\gA$ is equal to $\gA[(e_a)_{a\in\gA}]$.

\emph{4.} 
 If we identify $\Ared$ with its image in $\Abul$, we can identify $\Aqi$ with the sub\ri of $\Abul$ generated by $\Ared$ and by the \idms $e_x$ for $x\in\Ared$. 
 
\Grandcadre{In what follows we suppose \spdg that $\gA$ is reduced.}

\emph{5.} We refer to  Corollary~\ref{corlemZedClotBasicStep} for the description of the finite steps of the construction of $\Abul$. 
Given item \emph{4}, we get a description of the finite steps of a possible construction of $\Aqi$.
\\
For $a_1$, \dots, $a_n\in \gA$, we have an injection $\gA\to \gA[a_1\bul, \cdots, a_n\bul]=\gC$.\\ 
Let $e_i$ be the \idm $a_ia_i\bul$, $\gB=\gA[e_1, \ldots, e_n]\subseteq \gC$, and $e_I = \prod_{i\in I}(1-e_i) \prod_{j\notin I}e_j$ for $I \in \cP_n$.
Prove the following results.
\begin{itemize}
\item [\emph {a.}] The family $(e_I)_{I\in \cP_n}$ is a \sfio of~$\gB$ and $\gen {1-e_I}_\gB = \gen{(e_i)_{i\in I},\ (1-e_j)_{j\notin I}}_\gB$.

\item [\emph {b.}]
$\Ann_\gB(a_i) = \gen {1-e_i}_\gB$.

\item [\emph {c.}]
$\gA \cap \gen {e_i, \in I}_\gB = \DA(a_i, i \in I)$.

\item [\emph {d.}]
By letting $\fa'_I = (\DA(a_i, i \in I) : \prod_{j\notin I} a_j)$, we have $\gA\cap \gen {1-e_I}_\gB = \fa'_I$ and an \iso $\gB \simeq \prod_{I\in\cP_n}\gA\sur{\fa'_I}$.

\item [\emph {e.}]
The \ri $\gC$ is a \lon of the \ri $\gB$:
$\gC=\gB_s$ with \ndz $s \in \gB$.

\end{itemize}

In particular, let $a\in \gA$ and $\gA[e_a] \subseteq \gA[a\bul]$ with $e_a = aa\bul$. \\
Then, $\Ann_{\gA[e_a]}(a) = \gen {1-e_a}$, $\gA[e_a] \simeq \gA\sur{\Ann_\gA(a)} \times \gA\sur{\DA(a)}$, with $e_a \leftrightarrow (1,0)$, 
and $\gA[a\bul]$ is the localized \ri $\gA[e_a]_s$ with \ndz $s = 1-e_a + a$.
\\
In what follows let $\gA_{\bra{\an}}$ for $\gA[a_1a_1\bul,\ldots,a_na_n\bul]$

\emph{6.}
Let $\varphi : \gA \to \gD$ be a morphism with $\gD$ reduced, $a \in \gA$ and $b = \varphi(a)$. Suppose that $\Ann_\gD(b) = \gen {1-e_b}_\gD$ with \idm $e_b$. Show that we can extend $\varphi$  to a morphism of $\gA_{\bra a}  \to \gD$ realizing $e_a \mapsto e_b$.\\
However, in \gnlz, for $a$, $b \in \gA$, the \ris $\gA_{\bra{a,b}}$ and $(\gA_{\bra a})_{\bra b}$ are not \isocz.

\emph{7.} 
 Give a precise description of $\ZZ_\mathrm{pp}$.\\ 
 Explain why the \homo $\ZZ_\mathrm{pp}\to(\ZZ_\mathrm{pp})_\mathrm{pp}$ is not an \isoz.

\emph{8.} 
 (In \clamaz) If $\gA$ is \qiriz, and $\imath:\gA\to\Frac\gA$ is the canonical injection, then $\Spec \imath$ establishes a bijection between $\Spec (\Frac\gA)$ and~$\Min \gA$.

\emph{9.} 
 (In \clamaz) For every \ri $\gA$, there is a natural bijection between $\Min(\Aqi )$ and $\Spec \gA$.

\smallskip 
\comm
Despite $\ZZ$ being a \qiriz, $\ZZ_\mathrm{pp}$ is not \isoc to $\ZZ$. This is understood by observing that the natural \prn $\ZZ\to\ZZ\sur{15\ZZ}$ is not a \qiri morphism. This situation is different from that of the \zedr closure: 
 this comes from the fact that the quasi-inverse $b$ of an \eltz~$a$, when it exists, is unique and simply defined by two \eqnsz~\hbox{$ab^2=b$} and~\hbox{$a^2b=a$}, which implies that every \ri \homo respects quasi-inverses.
\eoe

}
\end{problem}

}

\sol{

\exer{exoRegleCoupure}
Indeed, $(a \vi b) \vi x = a \vi x$ since $x \vi a \le b$.\\
And $(a \vi b) \vu x = (a \vu x) \vi (b \vu x) = a \vu x$
because $a \vu x \le  b \vu x$ (since $a \le x \vu b$).\\
Therefore $a \vi b = a$, i.e. $a \le b$.

\exer{exo0GpRtcl} We write $a \sim b$ to indicate that $a$ and $b$ are associated.
Let us prove the following form (which is actually stronger if the \dve in $\gA$ is not explicit): if~$p$ is \irdz, $p \divi ab$ and $(pa, ab)$ has a gcd $d$, then $p \divi a$ or $p \divi b$.  
We have $p \divi pa$ and~$p \divi ab$, so $p \divi d$. Furthermore $a \divi pa$, $a \divi ab$, so $a \divi d$. 
Let $a' = d/a \in \gA$. As~$d \divi pa$, we have $a' \divi p$. 
But $p$ being \irdz, we either have $a' \sim 1$, or $a' \sim p$. \\
In the first case, \hbox{$d \sim a$}, and as $p \divi d$, we have $p \divi a$. 
In the second case, we have \hbox{$d \sim ap$}, thus~\hbox{$ap \divi ab$}, \cad $p \divi b$.
\\
In $\ZZ[X^2, X^3]$, $X^2$ is \irdz, $X^2 \divi X^3 \cdot X^3$ but $X^2 \nedivi X^3$, so~\hbox{$X^2 \cdot X^3$} and $X^3 \cdot X^3$ do not have a gcd. A fortiori they do not have a lcm.
\\
Finally, the gcd of $X^2$ and $X^3$ in $\ZZ[X^2,X^3]$ is $1$, if they had a lcm it would therefore be $X^5$, but $X^5$ does not divide $X^6$.


\exer{exoLemmeGaussAX}
Let $\G(\fa)$ be the gcd of the \gtrs of a \itf $\fa$. We easily observe that it is well-defined. Next, the \dit \hbox{$a(b\vi c)=ab\vi ac$} is generalized in the form \hbox{$\G(\fa)\G(\fb)=\G(\fa\fb)$} for two \itfsz~$\fa$ and~$\fb$. Finally, for two \pols \hbox{$f$, $g\in\AX$}, \DKM states that

\snic{\rc(f)^{p+1}\rc(g)=\rc(f)^{p}\rc(fg)$ for $p\geq\deg g.}

As~\hbox{$\G(f)=\G\big(\rc(f)\big)$} we obtain $\G(f)^{p+1}\G(g)=\G(f)^{p}\G(fg)$, and since they are \elts of the \riz, we can simplify to obtain $\G(fg)=\G(f)\G(g)$.


\exer{exoKroneckerTrick} 
\emph {1.}
Let $\uX^\alpha = X_1^{\alpha_1} \cdots X_n^{\alpha_n}\in\AuX_{<d}$, then

\snic { 
\varphi(\uX^\alpha) = T^a  \quad \hbox {with} \quad
a = \alpha_1 + \alpha_2 d + \cdots + \alpha_n d^{n-1}.
}

We thus see that $a < d^n$. The numbering in base $d$ proves that $\varphi$ transforms {the $\gA$-basis} of $\AuX_{<d}$ constituting of $\uX^\alpha$'s with $\alpha_i < d$ into the $\gA$-basis $(1, T, \ldots, T^{d^n-1})$ of $\gA[T]_{<d^n}$.

\emph {2.}
 Let us recall that $\AX\eti=\Ati=\AuX\eti$. Here we assume that $\AT$ is
 a UFD. 
If $P = QR\in \AuX_{<d}$ then $Q$ and $R \in \AuX_{<d}$, and $\varphi(Q)$ and $\varphi(R) \in \gA[T]_{<d^n}$.  Since $\varphi(QR) = \varphi(Q)\varphi(R)$, and $f=\varphi(P)$ has only finitely many factors (in $\AX\etl/\Ati$), it is sufficient to test for each factor $g(T)$ of $f(T)$ 
if $\varphi^{-1}(g)$ is a factor of $P$. This is possible since $\gA$ is supposed to be \dveez.


\exer{exoClosedCover1}{
We reduce to $r = 2$. The hypothesis  \gui{$M$ is \tf modulo $\fa_i$,} provides a \tf submodule $M_i$ of $M$  such that $M = M_i + \fa_i M$. By substituting the value of $M$ in the right-hand side, we obtain
$$
M = M_i + \fa_i M_i + \fa_i^2 M = M_i + \fa_i^2 M.
$$
By iterating, we obtain for $k \ge 1$, $M = M_i + \fa_i^k M$.
By substituting $M = M_2 + \fa_2^k M$ in $M = M_1 + \fa_1^k M$, we obtain $M = M_1 + M_2 + (\fa_1\fa_2)^k M$. But $\fa_1$, $\fa_2$ are \tf and $\fa_1\fa_2 \subseteq \DA(0)$, so there exists some $k$ such that $(\fa_1\fa_2)^k = \{0\}$, and consequently~\hbox{$M = M_1+M_2$} is \tfz.
}


\exer{exoAbul} 
We can assume that $\gA$ is a reduced sub\ri of $\Abul$. 
\\ \emph{1.}
Let $\fp$ be a \idep of $\gA$; the canonical morphism $\gA\to \gK=\Frac(\gA\sur\fp)$ can be factorized 
 through $\Abul$:

\centerline{\xymatrix @R=10pt @C=25pt{
\gA\ar[d]\ar[dr] \\
\Abul\ar@{-->}[r]_(.3){\pi_\fp}  & \Frac(\gA\sur\fp) \\
}}

The morphism $\pi_\fp$ is surjective because for $a \in \gA\setminus\fp$, we have $1/a = \pi_b(a\bul)$ in $\gK$. Its kernel $\fq = \Ker\pi_p$ is a \idema of $\Abul$; we then have $\gA\sur\fp \subseteq \gK \simeq \Abul\sur\fq$, so the natural arrow $\gA\sur\fp \to \Abul\sur\fq$ being injective, $\fp = \fq\cap\gA$.
We thus dispose of two transformations 

\snic{\Spec\Abul \to \Spec\gA$, $\fq
\mapsto \fq\cap\gA$, and $\Spec\gA \to \Spec\Abul$, $\fp \mapsto \Ker\pi_p,}

which are inverses of one another.
Indeed, if  $\fq\in\Spec\Abul$ and $\fp=\fq\cap\gA$, then~\hbox{$\gK = \Abul\sur\fq$} (because $a\bul = 1/a$ for $a \in \gA\setminus\fp$) so $\Ker\pi_\fp = \fq$.
\\
\emph{2.} 
By item \emph{1} the \homo $\Abul \to \wi\gA$ that factorizes the natural \homo $\gA\to  \wi\gA$ is injective, because its kernel is the intersection of all the \ideps of $\Abul$. We identify $\gA\subseteq\Abul\subseteq\wi\gA$. Lemma~\ref{lem2SousZedRed} describes the smallest \zedr sub\ri of $\wi\gA$ containing $\gA$. We see that this is indeed of $\Abul$ (by the construction of $\Abul$).
\\
\emph{Another \demz, left to the reader.} Let $\gA_1$ be the smallest \zedr sub\ri of $\wi\gA$ containing $\gA$. We then prove that this object satisfies the desired \uvl \prtz.


\exer{exoMinA} \emph{1.} The first \carn of strict maximal filters among the strict filters is \imdez: it is the same as saying that every attempt to make $\ff$ grow 
                         by adding an exterior \elt $x$ to it 
fails, because the filter generated by~$\ff$ and $x$ contains $0$.
\\
Let us then prove that a maximal strict filter among the strict filters is prime. 
\\
Let $x$, $y\in\gA$ with $x+y\in\ff$. We want to show that $x\in\ff$ or $y\in\ff$. If~\hbox{$x\notin\ff$}, there exist $a\in\ff$ and $n\in\NN$ such that $a^nx=0$,
so $a^n(x+y)=a^ny\in\ff$ therefore $y\in\ff$. 
\\
Let us now show that the localized \ri is \zedz, \cad (since the \ri is local) that every noninvertible \elt is nilpotent. A noninvertible \elt in 
the localized \ri is a multiple of $x/1$ with $x\notin \ff$. It suffices to see that~$x/1$ is nilpotent in $\ff^{-1}\gA$, but there exists an $a\in\ff$ such that $ax$ is nilpotent in~$\gA$, and $a$ is \iv in the localized \riz.
\\
Let us finally show that if  $\ff^{-1}\gA$ is local \zed and nontrivial, then $\ff$ is strict, maximal among the strict filters.
Indeed, some $x\notin\ff$ is not invertible, therefore is nilpotent in the localized \riz, which means that there exists an $a\in\ff$ such that $ax$ is nilpotent in $\gA$.

\emph{2.} 
Let $S$ be the subset defined by \Eqvrf{Eqnilreg}.
If $a\in S$ and $a\notin\ff$ with $\ff$ a maximal filter, we have $0\in a^\NN\ff$ which means that for some $x\in\ff$ and $n\in\NN$, $xa^n=0$, so, since $a\in S$, $x$ is nilpotent; a contradiction.\\
If $a\notin S$, there exists some non-nilpotent $x$ such that $xa$ is nilpotent. Therefore there
exists a strict filter containing $x$. By Zorn's lemma there exists a maximal filter $\ff$ containing $x$, and $a$ cannot be in $\ff$ because otherwise $xa$ and therefore $0$ would be in $\ff$.

\exer{exoFreeBooleAlgebra} \emph {1.}
Clearly results from the definition of a \agB as a \ri where all the \elts are \idmz, provided that we verify that the constructed object is indeed a \agBz, which offers no difficulty. Notice that $\gB$ is isomorphic~to

\snic {
\FF_2[x_1] \otimes_{\FF_2} \cdots  \otimes_{\FF_2} \FF_2[x_n],
}

which is the direct sum of $n$ \agBs freely generated by a single \gtr in the category of \agBsz.
Indeed, the direct sum of two \agBs $\gB$, $\gB'$ is the \agB $\gB\otimes_{\FF_2}\gB'$.

\emph {2.}
The monomial $\FF_2$-basis of $\gB$ is $(m_I)$ with $m_I = \prod_{i \in I} x_i$.  It is of cardinality~$2^n$, so $\gB$ is of cardinality~$2^{2^n}$. We define $e_I$ by $e_I = m_I \prod_{j \notin I} (1 + x_j)$; we easily verify that $(e_I)$ is a \sfioz, that~\hbox{$m_I e_J = e_J$} if~\hbox{$I \subseteq J$}, and $0$ otherwise. \\
We have the same expression $e_I = \sum_{J \,|\, J \supseteq I} m_J$ and $m_I = \sum_{J \,|\, J \supseteq I} e_J$ (which confirms that $(e_I)$ is an $\FF_2$-basis of $\gB$).
\\
With respect to the description given in this course, $x_1^{\varepsilon_1} \cdots x_n^{\varepsilon_n}$ corresponds to the following \elt of $\Pf\big(\Pf(E)\big)$: $\so{\sotq{x_i}{\varepsilon_i = 1}}$.


\exer{exotrdifree} The \trdi freely generated by $\emptyset$  is the lattice $\Deux$.\\
The \trdi freely generated by $\so{a}$ is $\so{0,a,1}$.\\
The \trdi freely generated by $\so{a,b}$ ($a\neq b$) is: {\mathrigid 2mu $\so{0,\,a\vi b,\, a,\, b, a\vu b,\, 1}$}.\\
The \trdi freely generated by $\so{a,b,c}$ $(a\neq b\neq c\neq a)$ contains: 
\[\preskip.2em \postskip.4em 
\begin{array}{ccc} 
 0,\, 1,\, a,\,b,\,c,\,a\vu b,\,a\vu c,\,b\vu c,\,a\vu b\vu c,\,a\vi b,\,a\vi c,\,b\vi c,\,a\vi b\vi c,   \\[.3em] 
a\vi(b\vu c),\,b\vi(a\vu c),\,c\vi(a\vu b),\,(a\vu b)\vi(a\vu c),\,(a\vu b)\vi(b\vu c),\,\\[.3em] 
(a\vu c)\vi(b\vu c),\,(a\vu b)\vi(a\vu c)\vi(b\vu c). \end{array}
\]

\exer{exoTreillisDistributifFini}
\emph{1.} By \dfn of an initial section the intersection and the union of two initial sections is another initial section.   \\
Therefore in $E\sta$: $S_1 \vi S_2 = S_1 \cup S_2$, $S_1 \vu S_2 = S_1 \cap S_2$, $\emptyset =
1_{E\sta}$ and $E = 0_{E\sta}$.

\emph{2.} It amounts to the same to give $\alpha\in\SpecT$ or the \idep $\Ker\alpha$.
This leads us to order the set of prime \ids of $\gT$ by the relation $\supseteq$. \\
Indeed, if~$\alpha$, $\beta : \gT \to \{0,1\}$ are two morphisms, we have the \eqvc
$$\preskip.4em \postskip.0em 
 \alpha \le \beta \iff \Ker\alpha \supseteq
\Ker\beta.
$$
\emph{3.} We have
$$\preskip.2em \postskip.4em 
~~~\Ker\nu_E(x) = \sotq {S \in E\sta} {x \in S} =
\sotq {S \in E\sta} {\dar x \subseteq S} =
\sotq {S \in E\sta} {S \le \dar x}
, 
$$
\cad $\Ker\nu_E(x) = \cI_{E\sta}(\dar x)$. We indeed have the \eqvcs
$$\preskip.4em \postskip.4em \mathrigid 4mu
 x \le y \iff \dar x \subseteq \dar y \iff \dar y \le \dar x
\iff \cI(\dar y) \subseteq \cI(\dar x) \iff
\cI(\dar x) \le \cI(\dar y).
$$
Moreover, in $E\sta$: $S_1 \vi S_2 \le \dar x \Rightarrow (S_1 \le \dar x) \hbox { or } (S_2 \le \dar x)$ (because the first in\egt means $\dar x \subseteq S_1 \cup S_2$, \cad $x \in S_1 \cup S_2$), and since $\dar x \ne 1_{E^*} = \emptyset$,  $\dar x$ generates a \idepz.  
Conversely, let $\fp$ be a \idep of $E\sta$.
Being finite, it is principal: $\fp = \cI_{E\sta}(S)$ with $S \ne 1_{E\sta}$, \cad $S$ is nonempty. It must be shown that $S$ is of the form $\dar x$.
If $S =\{x_1, \ldots, x_n\}$, we have $S = (\dar x_1)\cup\cdots\cup(\dar x_n)$, \cad \hbox{$(\dar x_1) \vi \cdots \vi (\dar x_n) = S$}.
As $S$ generates a \idepz, there exists some~$i$ such that $\dar x_i \le S$, \cad $S \subseteq \dar x_i$, then $S = \dar x_i$.

 \emph{4.}
We determine $E\sta$ by noticing that every initial section is a union of subsets $\dar x$. The picture of the lattice $E\sta$ is the following

\snic {
\xymatrix @R = 5pt @C = 4pt{
       & \emptyset\\
       & \dar a = \{a\}\ar@{-}[u] \\
\dar c =\{a,c\}\ar@{-}[ur] &&\dar b =\{a,b\}\ar@{-}[ul] \\
       &\{a,b,c\}\ar@{-}[ul]\ar@{-}[ur] &   & \dar d = \{a,b,d\}\ar@{-}[ul] \\
       &                              &\{a,b,c,d\}\ar@{-}[ul]\ar@{-}[ur] \\
}
}

If $E$ is totally ordered, then $E\sta = \{\, \dar x \;\vert\; x \in E\,\}
\cup \{\emptyset\}$ is also totally ordered and $\#E\sta = 1 + \#E$.
If $\gT$ is a finite totally ordered set, then 

\snic{\Spec\gT = \big\{\cI_{\gT}(a) \;\vert\;  a \in \gT \setminus \{1_\gT\}\big\}$, and $\#\Spec\gT = \#\gT - 1.}

If $E$ is ordered by the \egt relation, $E\sta = \cP(E)$ ordered by $\supseteq$. As for $\Spec(E\sta)$, it is the set $\cI_{\cP(E)}(\{x\})$ with $x \in E$ (which is indeed isomorphic to $E$).

 \emph{5.}
The reader will verify that by letting, for $a \in \gT$, $\wh {a} = \sotq {\fp \in \Spec\gT}{a \in \fp}$, we obtain an initial section, that every initial section of $\Spec\gT$ is of this form, and finally that $a \le b \iff \wh {a} \le \wh {b}$.

\emph{6.}
We now consider $E\sta$ and $\Spec \gT$ with the order relation $\subseteq$.
\\
Then $S_1 \vi S_2 = S_1 \cap S_2$, $S_1 \vu S_2 = S_1 \cup S_2$, $\emptyset = 0_{E\sta}$, $E =
1_{E\sta}$.  For~\hbox{$x \in E$}, we let~\hbox{$\wi x = E \setminus \uar x = \sotq {y \in E} {y \not\geq x}$}: this \elt of $E\sta$ satisfies, for~\hbox{$S \in E\sta$}, the \eqvc $x \notin S \iff S \subseteq \wi x$.  
We have $\wi x \ne 1_{E\sta} = E$, and $\wi x$ generates a prime \id of the lattice $E\sta$: $S_1
\vi S_2 \le \wi x \Rightarrow S_1 \le \wi x \hbox { or } S_2 \le \wi x$ 
(indeed, the hypothesis is~\hbox{$\uar x \subseteq (E \setminus S_1) \cup (E \setminus S_2)$}, therefore for example $x \notin S_1$, \cad~\hbox{$S_1 \subseteq \wi {x}$)}.  We have the \eqvc $x \le y \iff \wi x \subseteq \wi y$.  We prove that every prime \id of $E\sta$ is of the form $\wi x$, so the ordered set $E$ is isomorphic, via~\hbox{$x \mapsto\cI_{E\sta} (\wi x)$}, to the set of prime \ids of $E\sta$, ordered by inclusion.
$$\preskip.4em \postskip.4em 
\xymatrix @R = 5pt @C = 4pt{
       & \{a,b,c,d\} \\
\wi c =\{a,b,d\}\ar@{-}[ur] &&\wi d =\{a,b,c\}\ar@{-}[ul] \\
       &\{a,b\}\ar@{-}[ul]\ar@{-}[ur] &   & \wi b = \{a,c\}\ar@{-}[ul] \\
       &                              &\{a\}\ar@{-}[ul]\ar@{-}[ur] \\
       &                              &\wi a = \emptyset \ar@{-}[u] \\
} 
$$


\exer{exoIVpgcd}
Since  $ \gen{a,b} $ is invertible we have $s,$ $t$, $u$, $v$  with  $sa = ub$,  $tb = va$ and $s+t = 1$.
\\
Since $m$ is the lcm of $a$ and $b$ we can write

\snic{m= ab' = ba'  \;\hbox{ and }\;   ab/m = g = b/b' = a/a'.}            
           
Thus $sa = mx = ab'x$ and $tb = m=ba'y$, which give $s=b'x$ and $t=a'y$.
\\
Therefore $b'x+a'y = 1$,
  $bx + ay = gb'x + ga'y = g$
and consequently $ \gen{a,b}=\gen{g} $.             


\exer{exoFactFini} \emph{(A UFD with only a finite number of \ird \eltsz)}\\ 
Let $(p_i)_{i\in I}$ be the finite family of distinct \ird \elts (up to association).
\\
We must show that $\gA$ is a Bézout \riz. In order to do so, it suffices to show that if $a$ and $b\in\Atl$ have as their gcd $1$, then $\gen{a,b}=\gen{1}=\gA$.
We write 
$$
\ndsp a=\prod_{i\in A}p_i^{\alpha_i}, \;b=\prod_{j\in B}p_i^{\beta_j}, \hbox{with}  \;\alpha_i\hbox{'s},  \,\beta_j\hbox{'s}>0 \,\hbox{ and }\,A\cap B =\emptyset.
$$ 
Let $C=I\setminus(A\cup B)$ and $c=\prod_{k\in C}p_k$. We show that $a+bc\in\Ati$. 
\\
Indeed, for $i\in A$, $p_i$ divides $a$, therefore it cannot divide~\hbox{$a+bc$}, otherwise it would divide~\hbox{$bc=(a+bc)-a$}. Similarly, for $j\in B\cup C$, $p_j$ cannot divide~\hbox{$a+bc$}, otherwise it would divide $a=(a+bc)-bc$.
Thus $a+bc$  is not divisible by any \ird \eltz.


\exer{exoPrincipalIntersecSoucorps} \emph{(An interesting intersection)}\\
Consider the \evn \homo $$\varphi:\gk[z,u]\to \gk[z,x+yz] , \;z\mapsto z, \,u\mapsto x+yz.$$ 
It is surjective by construction. It is injective because, for $f=f(z,u)$, by evaluating $\varphi(f)$  in $\gk[x,y,z]$ we obtain $\varphi(f)(x,0,z)=f(z,x)$. It is therefore indeed an \isoz.
\\
In what follows we can therefore let~\hbox{$u=x+yz$}, with $\gk[z,x+yz]=\gk[z,u]$ where $z$ and $u$
play the role of distinct \idtrsz. 
\\
Moreover we notice that $\gk[z,u][y]=\gk[x,y,z]$. As $\gk[z,u]$ is a  GCD-domain, this implies that two \elts of $\gk[z,u]$ have gcd $1$ in $\gk[z,u]$ \ssi they have gcd $1$ in $\gk[x,y,z]$.

Now let $h\in\gA$ be an arbitrary element that we write in the form of an \ird fraction $f(z,u)/g(z,u)$ in $\gk(z,u)$, and in the form of a fraction~\hbox{$a/b$} ($a\in\gk[x,y,z]$, $b\in\gk[x,y]$) as an \elt of $\gk(x,y)[z]$. This last fraction can itself be written
in irreducible form, that is so that the gcd 
of~$a$ and $b$ in $\gk[x,y,z]$ is equal to $1$. By uniqueness of the expression of a fraction in reduced form, we therefore have a constant $\gamma\in\gk\etl$ such that $f(z,u)=\gamma a(x,y,z)$ and~\hbox{$g(z,u)=\gamma b(x,y)$}. 
\\
It remains to show that the denominator $g(z,x+yz)$ is a constant. By making $z=0$ in the \egt $g(z,x+yz)=\gamma b(x,y)$ we obtain $$g(0,x)=\gamma  b(x,y)=c(x).$$
Finally, by making $(z,y)=(1,-x)$ in the \egt $g(z,x+yz)=c(x)$, we obtain~\hbox{$c(x)=g(1,0)$}.


\prob{exoAutourGaussJoyal}
The first item is left to the reader. Let $fg = \sum_k c_k X^k$.

 \emph{2.}
We easily have $u(fg) \le u(f) \vi u(g)$. \\
Indeed, $c_k = \sum_{i+j = k}
a_ib_j$, so $u(c_k) \le \Vu_{i+j = k} u(a_ib_j) \le \Vu_{i} u(a_i)
 = u(f)$ (we have used $u(ab) \le u(a)$).
\\
If we dispose of the Gauss-Joyal lemma, 
then $u(a_ib_j) \le u(a_i)
\vi u(b_j) \le u(f) \vi u(g) = u(fg)$.  Conversely, if we know how to prove $u(a_ib_j) \le u(fg)$ for all $i,j$, then

\snic {
\Vu_{i,j} u(a_ib_j) \le u(fg),
\; \hbox {i.e. by \ditz} \;
\bigr(\Vu_{i} u(a_i)\bigl) \vi \bigr(\Vu_{j} u(b_j)\bigr)
\le u(fg),
}

\cad $u(f) \vi u(g) \le u(fg)$.

 \emph{3.}
If $\gA$ is integral, the same goes for $\gA[X]$.

 \emph{4.}
Let us show by decreasing \recu on $i_0+j_0$ that $u(a_{i_0}b_{j_0}) \le u(fg)$.
It is true if $i_0$ or $j_0$ is large because then $a_{i_0}b_{j_0} = 0$.  
We write the \dfn of the product of two \polsz

\snic {
a_{i_0}b_{j_0} = c_{i_0+j_0}
- \sum\limits_{i+j = i_0+j_0 \atop i > i_0} a_ib_j
- \sum\limits_{i+j = i_0+j_0 \atop j > j_0} a_ib_j.
}

We apply $u$ by using on the one hand $u(\alpha + \beta + \cdots) \le u(\alpha) \vee u(\beta) \vee \dots$ and on the other hand $u(\alpha\beta) \le u(\alpha)$ to obtain

\snic {(\star)\;:\;
u(a_{i_0}b_{j_0}) \le u(c_{i_0+j_0}) \vee
\Vu_{i > i_0} u(a_i) \vee \Vu_{j > j_0} u(b_j).}

We thus dispose of an in\egt $x \le y$ which we write as $x \le x \wedge y$. In other words, in $(\star)$, we reinsert $u(a_{i_0}b_{j_0})$ in the right-hand side, which gives, by  \ditz

\snic {
u(a_{i_0}b_{j_0}) \le u(c_{i_0+j_0}) \vee
\Vu_{i > i_0} \big(u(a_i) \wedge u(a_{i_0}b_{j_0})\big) \vee
\Vu_{j > j_0} \big(u(b_j) \wedge u(a_{i_0}b_{j_0})\big).
}

By using $u(a_i) \wedge u(a_{i_0}b_{j_0}) \le u(a_i) \wedge u(b_{j_0})$ and $u(b_j) \wedge u(a_{i_0}b_{j_0}) \le u(b_j) \wedge u(a_{i_0})$, and (by definition) $u(c_{i_0+j_0}) \le u(fg)$, we bound $u(a_{i_0}b_{j_0})$ above by

\snic {
u(fg) \vee \Vu_{i > i_0} u(a_{i}b_{j_0}) \vee
\Vu_{j > j_0} u(a_{i_0}b_{j}).
}

By \recu on $i_0, j_0$, $u(a_{i}b_{j_0}) \le u(fg)$,
$u(a_{i_0}b_{j}) \le u(fg)$.\\ Hence $u(a_{i_0}b_{j_0}) \le u(fg)$.

 \emph{5.}
In this case $a_ib_j \in \rD_\gA(c_k, k = 0, \ldots)$, which is the usual Gauss-Joyal lemma.


\prob{exoQiClot} \emph{(\qiri closure of a commutative \riz)}

Preliminary remark: if in a \ri $\gA$, $\Ann(a) = \gen {e'_a}$ with~$e'_a$ \idmz, then $e'_a$ is the unique $e'$ such that

\snic {
e' a = 0,\quad e'+a\hbox { is \ndzz} \quad\hbox{and}\quad  e' \hbox { is \idmz.}
}

Indeed, $e' = e'e'_a$ (because $e'a=0$) and $(e'+a)e' = (e'+a)e'_a$ ($=e'$) hence $e'=e'_a$.

\emph{1.} Let $\gA$, $\gB$ be \qiris and a \qiri morphism $\varphi:\gA\to\gB$.
\\
\emph{1a.} 
If $a\in\gA$ is \ndzz, $e_a=1$ so $e_{\varphi(a)}=1$ therefore $\varphi(a)$ is \ndzz. Conversely, let $\psi:\gA\to\gB$ be a \ri \homo which transforms every \ndz \elt into a \ndz \eltz.
Let $a\in\gA$, $b=\psi(a)$ and $f=\psi(1-e_a)$.  \\
Then $fb = \psi\big((1-e_a)a\big) = 0$, $f+b = \psi(1-e_a + a)$ is \ndz and $f^2 = f$, and \hbox{so $f = 1-e_b$}.
\\
\emph{1b.} Suppose $\varphi(x)=0$, then $e_{\varphi(x)}=0$, \cad $\varphi(e_x)=0$. Therefore if $\varphi|_{\BB(\gA)}$ is injective, we obtain $e_x=0$, \cad $x=0$.
\\
\emph{1c.} We consider the unique \homo  $\rho:\ZZ\to\prod_{n>0}\aqo\ZZ{2^n}$. Then $\rho$ is injective but $\rho(2)$ is not \ndzz.
\\
\emph{1d.} The \homo preserves the quasi-inverses, therefore also the associated \idms because $e_a=aa\bul$ if $a\bul$ is the quasi-inverse of $a$.  
\\
\emph{1e.} 
Results \imdt from Fact~\ref{factQoQiZed}.

 \emph{2.} Since $\Aqi$ is reduced, there is a unique \ri \homo $\Ared\to\Aqi$ which factorizes the two canonical \homos $\gA\to\Ared$ and $\gA\to\Aqi$.
Since $\Abul$ is a \qiriz, there is a unique \qiri morphism $\Aqi\to\Abul$ that factorizes the two canonical \homos $\gA\to\Aqi$ and $\gA\to\Abul$.
 Since the morphism $\Aqi\to\Abul$ transforms a \ndz \elt into a \ndz \eltz, and since a \ndz \elt in a (reduced or not) \zed \ri is \ivz, there exists a unique \homo $\Frac(\Aqi)\to\Abul$ which factorizes the two canonical \homos $\Aqi\to\Frac(\Aqi)$ and $\Aqi\to\Abul$.
\\
Similarly, for every \homo $\gA\to\gB$ with $\gB$ being \zedrz, we first obtain a unique \qiri \hbox{morphism $\Aqi\to\gB$} (which factorizes what is needed), then a unique morphism $\Frac(\Aqi)\to\gB$ which factorizes the two \homos $\gA\to\Frac(\Aqi)$ and $\gA\to\gB$. \\
 In other words, since $\Frac(\Aqi)$ is \zedrz, it solves the \uvl \pb of the \zedr closure for $\gA$. Consequently the \homo $\Frac(\Aqi)\to\Abul$ that we have constructed is an \isoz.  

 \emph{3.} This item is copied from Lemma~\ref{lem2SousZedRed} which concerns the \zedr \risz: the reader could also just about copy the \demz.

 \emph{4.} First of all note that the natural \homo $\Ared\to\Aqi$ is injective because the \homo $\Ared\to\Abul$ is injective and there is \fcnz.
We can therefore identify $\Ared$ with a sub\ri of $\Aqi$, 
which is itself identified with a sub\ri of $\Frac(\Aqi)$ that we identify with $\Abul$.
In this framework $\Aqi$ \ncrt contains $\Ared$ while the \elts $e_x=xx\bul$ for all $x\in\Ared$ since the morphism $\Aqi\to\Abul$ is a \qiri and is injective.\\ 
Let $\gB$ be the sub\ri of $\Abul$ generated by $\Ared$ and the \idms $(e_x)_{x\in\Ared}$.
It remains to see that the inclusion $\gB\subseteq\Aqi$ is an \egtz.
\\
It is clear that $\Frac(\gB) = \Frac(\Aqi)$. On the one hand, as $\gB$ is a \qiriz, the injection $\Ared \to \gB$ provides a (unique) \qiri morphism $\varphi : \Aqi \to\gB$ such that $\varphi(a) = a$ for every $a \in \Ared$.
Since the morphism is a \qiri morphism, we deduce that $\varphi(e_a) = e_a$ for every $a \in \Ared$, then  $\varphi(b)=b$ for all~\hbox{$b\in\gB$}. Let $x \in \Aqi$; we want to show that $x \in \gB$; as $x \in \Frac(\gB)$, there exists a \ndz $b \in \gB$ such that $bx \in \gB$ therefore $\varphi(bx) = bx$ \cad $b\varphi(x) = bx$; as $b$ is \ndz in $\gB$, it is \ndz in $\Frac(\gB)$, a fortiori in $\Aqi$, so $x = \varphi(x) \in \gB$.

\emph{5a} and \emph{5b.} Easy.

\emph{5c.}
Since $a_j = a_je_j$, we have, for $j \in I$, $a_j \in \gen{e_i, i \in I}_\gB = \gen {e}_\gB$ with $e$ the \idm $1 - \prod_{i \in I}(1-e_i)$.  But in a reduced \riz, every \idm generates a radical \id

\snic {
b^m \in \gen{e} \Rightarrow b^m(1-e) = 0 \Rightarrow b(1-e) = 0 
\Rightarrow b = be \in \gen {e}.
}

Therefore $\DA(a_i, i \in I) \subseteq \gen{e_i, i \in I}_\gB$.  \\
Let us now show that $\gA \cap \gen {e_i, \in I}_\gC \subseteq \DA(a_i, i \in I)$. Let $x \in \gA \cap \gen {e_i, \in I}_\gC$; by returning to the initial \dfn of $\gC$, we have $x \in \gen{a_iT_i, i \in I}_{\gA[\uT]} + \fc$.
Let us work on the reduced \ri $\gA' = \gA\sur{\DA(a_i, i \in I)}$; we then have 

\snic{\ov x \in \rD_{\gA'[\uT]}(a_kT_k^2 - T_k, a_k^2T_k - a_k, k \in
\lrbn).}

Since $\gA' \to \gA'[\ov a_1\bul, \ldots, \ov a_n\bul]$ is injective, we have $\ov x = 0$ \cad $x \in \DA(a_i, i \in I)$.

\emph{5d.} 
Let $\pi$ be the product $\prod_{j\notin I} a_j$. Let $x \in \gA\cap \gen{1-e_I}_\gB$; since $\pi (1-e_j) = 0$ for $j \notin I$, we have $\pi x \in \gen{e_i, i \in I}_\gB$, so, by \emph{5c)}, $\pi x \in \DA(a_i, i \in I)$, {\cad  $x \in \fa'_I = (\DA(a_i, i \in I) : \pi)$}.

Conversely, let $x \in \fa'_I$; we write $x = \pi' x + (1-\pi')x$ with $\pi' = \prod_{j \notin I} e_j$.  
 We have $1-\pi' \in \gen{1-e_j, j \notin I}$. As for $\pi' x$, we notice that in $\gC$, $\gen{e_j}_\gC = \gen{a_j}_\gC$, so $\pi'x \in \gen{\pi x}_\gC \subseteq \rD_\gC(a_i, i \in I) \subseteq \gen {e_i, i \in I}_\gC$. \\
Recap: $x \in \gen {(e_i)_{i\in I}, (1-e_j)_{j \notin I}}_\gC = \gen {1-e_I}_\gC$. \\
But $\gA \cap \gen{1-e_I}_\gC = \gA \cap \gen {1-e_I}_\gB$, so $x \in \gen {1-e_I}_\gB$.

Finally, $\gB$ is \isoc to the product of $\gB\sur{\gen{1-e_I}_\gB}$
and $\gB\sur{\gen{1-e_I}_\gB} \simeq \gA\sur{\fa'_I}$.

\emph{5e.} 
Take $s = \sum_I e_I \prod_{j\notin I}a_j = \sum_I \prod_{i\in I} (1-e_i) \prod_{j\notin I}a_j$: $s$ is the unique \elt of $\gB$ which is equal to $\prod_{j\notin I}a_j$ over the component $e_I = 1$.

\emph{6.} In the \iso $\gA[e_a] \simeq \gA\sur{\Ann_\gA(a)} \times \gA\sur{\DA(a)}$, we have $e_a = (1,0)$ and so $(\ov x, \ov y) = xe_a + y(1-e_a)$. We then consider the map 

\snic{\gA\times\gA \to \gD$, $(x,y)
\mapsto \varphi(x)e_b + \varphi(y) (1-e_b).}

It is a \ri morphism and since $\gD$ is reduced, it passes to the quotient modulo $\Ann_\gA(a) \times \DA(a)$.
\\ 
Let us now compare $\gA_{\bra{a,b}}$ and $(\gA_{\bra a})_{\bra b}$.
We find

\snic {\arraycolsep2pt
\begin {array} {rcl}
\gA_{\bra{a,b}} &\simeq& 
\gA\sur{(0 : ab)} \times \gA\sur{(\rD(b) : a)} \times
\gA\sur{(\rD(a) : b)} \times \gA\sur{\rD(a,b)},
\\[1mm]
(\gA_{\bra a})_{\bra b} &\simeq&
\gA\sur{(0 : ab)} \times  \gA\sur{\rD\big((0 : a) + \gen {b}\big)} \times
\gA\sur{(\rD(a) : b)} \times \gA\sur{\rD(a,b)}.
\end {array}
}
Finally, note that $\rD\big((0 : a) + \gen {b}\big)$ is contained in $(\rD(b) : a)$ but that the inclusion can be strict. Take for example $\gA=\ZZ$, $a = 2p$, $b = 2q$ where~$p$ and~$q$ are two distinct odd primes.
We use $(x : y) = x/\pgcd(x,y)$ for  $x$, $y \in \ZZ$.\\
Then $\ZZ_{\bra{a,b}} \simeq \ZZ \times \ZZ\sur{q\ZZ} \times \ZZ\sur{p\ZZ} \times \ZZ\sur{2\ZZ}$,
but
$(\ZZ_{\bra a})_{\bra b} \simeq \ZZ \times\ZZ\sur{2q\ZZ} \times \ZZ\sur{p\ZZ} \times \ZZ\sur{2\ZZ}$. In the first \riz, $\Ann(a)$ is generated by $(0,0,1,1)$. In the second (the first \ri is a quotient), $\Ann(a)$ generated by the \idm $(0,q,1,1)$. 

\emph{8.} 
Recall (Exercise~\ref{exoMinA}) that a \idep $\fp$ of a \ri $\gA$ is minimal \ssi for all $x\in\fp$, there exists an $s \in \gA\setminus\fp$ such that $sx^n = 0$ for a certain $n$ (if $\gA$ is reduced, we can take $n = 1$).
\\
First, a \idemi of $\gA$ remains a strict \idep in~$\Frac(\gA)$ 
(this does not use the fact that $\gA$ is a \qiriz), 
\cad $\fp \cap \Reg(\gA) = \emptyset$: if $x \in \fp$, there exist $s \notin \fp$ and~\hbox{$n \in \NN$} such that $sx^n = 0$, which proves that $x \notin \Reg(\gA)$.
\\
Conversely, for $\fq$ a \idep of $\Frac(\gA)$, let us prove that $\fp = \fq\cap\gA$ is a \idemi of $\gA$. Let $x \in \fp$; then $x + 1-e_x$ is \ndz in $\gA$, so \iv in $\Frac(\gA)$, therefore $1-e_x \notin \fp$. Then $x(1-e_x)=xe_x(1-e_x)=0$: we have found $s=1-e_x \notin \fp$ such that $sx = 0$.

\emph{9.}
By Exercise~\ref{exoAbul}, the injection $\gA\to\Abul$ induces a bijection $\Spec\Abul\to\Spec\gA$; but $\Abul = \Frac(\Aqi)$ and $\Aqi$ is a \qiriz.\\ 
Therefore, by item \emph{8} applied to $\Aqi$, $\Spec\Abul$ is identified with $\Min(\Aqi)$, hence the natural bijection between $\Spec\gA$ and $\Min(\Aqi)$.


}

\Biblio

Some reference books on the study of lattices are \cite{Birkhoff}, \cite{Grae} and~\cite{Johnstone}. In \cite{Johnstone} the focus is essentially on \trdisz, which are the objects that primarily interest us. This book presents the theory of locales. The notion of a \ix{locale} is a \gnn of that of a topological space. The structure of a locale is given by the \trdi of its open sets, but the open sets are no longer \ncrt sets of points.
This is the reason why a locale is sometimes called a \emph{pointless topological space} \cite[Johnstone]{Joh}. The author \gnlt tries to give \prcos and explicitly signals the \thos whose proof uses the axiom of choice.

In abstract \algz, spectral spaces are omnipresent, foremost among which we include the Zariski spectrum of a commutative \riz.
From the \cof point of view they are very peculiar locales which \gui{lack points.}
We will quickly present this notion in Section~\iref{secEspSpectraux} of Chapter~\ref{chapKrulldim} devoted to the Krull dimension.

An elegant \dem of \thref{thAXgcd} (if $\gA$ is a GCD-domain the same goes for $\AX$) is found in \cite[th.\,IV.4.7]{MRR}.

The origin of \entrels is found in the Gentzen sequent calculus, which is the first to place a focus on the cut (the rule $(T)$).
The link with \trdis has been highlighted in \cite[Coquand\&al.]{cc,cp}.
The fundamental \tho of the entailment relations \rref{thEntRel1} is found in~\cite{cc}.
In fact, its first appearance seems to date back to the article by Paul Lorenzen \cite{Lor1951}, which studies the relations between formal logic and distributive lattices.

The dynamic method is clearly presented, for the first time it seems, in the article by Lorenzen \cite{Lor1953}, which uses the equivalent of the closed covering principle \ref{prcfgrl}. See on this subject articles \cite{CLN2019} and \cite{CLN2020}.

We find the terminology of an \emph{implicative lattice} in \cite{Curry} and that of a \emph{\agHz} in \cite{Johnstone}.

A basic book for the theory of \grls and of (not \ncrt commutative) lattice-group \ris is \cite{BKW}.
We have said that a guiding idea in the theory of \grls is that an \grl behaves in  computations like a product of totally ordered groups. This is translated in \clama by the representation \tho which affirms that every (Abelian) \grl is \isoc to an \sgrl of a product of totally ordered groups (\Tho 4.1.8 in the cited book).

The \grls that are \Qevs somewhat constitute the purely \agq version of the theory of Riesz spaces. Every good book on  Riesz spaces starts by developing the purely \agq \prts of these spaces, which are copied (with very similar, if not identical proofs) from the theory of (Abelian) \grlsz. See for example \cite{Zaanen}.

In the exercises of Bourbaki (Commutative \algz, Diviseurs) an integral Bézout \ri is called a \emph{anneau bezoutien}, a GCD-domain is called a \emph{anneau pseudo-bezoutien}, and a \ddp is called a \emph{anneau pr\"uferien}.

\newpage \thispagestyle{CMcadreseul}
\incrementeexosetprob


\chapter{Pr\"ufer and Dedekind \risz} 
\label{ChapAdpc}
\perso{compil\'e le \today}
\vskip-1em

\minitoc

\subsection*{Introduction} 
\addcontentsline{toc}{section}{Introduction}

The usual \dfns of \adk do not lend themselves to an \algq treatment.

First, the notion of \noet is delicate (from the \algq point of view). Secondly, the questions of \fcn \gnlt demand extremely strong hypotheses. 
For example, even if $\gK$ is a quite explicit \cdiz, there is no \gnl method (valid over all \cdisz) for factorizing the \pols of~$\KX$.

Thus, an essential aspect of the theory of \adksz, namely that the integral closure of a \adk in a finite extension of its quotient field remains a \adkz, no longer works in full \gnt (from an \algq point of view) if we require the complete \fcn of the \ids (see for example the treatment of this question in~\cite{MRR}).

Moreover, even if a \fac is theoretically feasible (in the \ris of integers of number fields for example), we very quickly encounter \pbs of a prohibitive complexity 
such as that of factorizing the discriminant (an impossible task in practice if it has several hundred digits).
Also Lenstra and Buchmann, \cite{LB}, 
proposed to work in the \ris of integers without having a $\ZZ$-basis at our disposal.
An important \algq fact is that it is always easy to obtain a \emph{\fapz} for a family of natural numbers, that is a \dcn of each of these numbers into a product of factors taken in a family of pairwise coprime numbers (see \cite[Bernstein]{Ber1}, and \cite[Bernstein]{Ber2} for an implementation with the \ids of number fields,
see also \Pbmz~\rref{exoPlgb2}).

A goal of this chapter is to show the \gnl validity of such a point of view and to propose tools in this framework.

\smallskip A crucial and simplifying role in the theory is played by the \anars (in accordance with an intuition of
Gian~Carlo~Rota~\cite{Rota}), that are the \ris in which the lattice of \ids is distributive, and by the \emph{\mlpsz}, which are the matrices that explicate the computational machinery of the \lop \itfsz, in an essentially \eqv way to what Dedekind \cite{Ddk2} estimated to be a fundamental \prt of \ris of integers in the number fields (see \cite[Avigad]{Avi} and item \emph{3$\,'.$} of our Proposition~\ref{propItfLocprinc}).

\smallskip The willingness to put off implementing, for as long as possible, \noees hypotheses has \egmt prompted us to develop a \cof treatment of several important points of the theory in a simpler and less rigid framework than that of \adksz. This is the context of \ris that have the two following \prts
\begin{itemize}
\item  the \itfs are \pros (this characterizes what we call a \emph{\adpcz}),
\item  the \ddk is at most $1$.
\end{itemize}
As the reader will observe, the proofs do not become more complicated, on the contrary, by this weakening of the hypotheses.

\smallskip Similarly, we have been brought to study the \fap\adps 
(in the local case, they are the \ddvs whose group of valuation is \isoc to a subgroup of~$\QQ$). We think that these \ris constitute the natural framework suggested by Buchman and Lenstra~\cite{LB}.

\smallskip Finally, for what concerns the \adksz, we have freed ourselves of the usual hypothesis of integrity (because it is hardly preserved from an \algq point of view by \agq extension) and we have left the \fac of the
\itfs in the background (for the same reason) in favor of the only \noe character. The \noet implies the \fap of families of \itfsz, which itself implies the dimension $\leq 1$ in the \cov form.

\section{\Anarsz} 
\label{secAnars}

Recall that an \anar is a \ri whose \itfs are \lops (see Section~\ref{secIplatTf}).
We begin with a few results regarding the \lop \ids in an arbitrary \riz.

\subsec{\Lop \idsz, principal \lon matrix}

We take up \thref{propmlm} again (stated for the \lmo \mtfsz) in the framework of  \lop \idsz.

\begin{proposition} 
\label{propItfLocprinc} {\em (\Lop  \idsz)}\\
Let $x_1$, \ldots, $x_n\in\gA$. 
\Propeq 
\begin{enumerate}
\item  The \id $\fa=\gen{\xn}$ is \lopz.
\item  There exist $n$ \eco $s_i$ of $\gA$ such that for each $i$, after \lon at $s_i$, $\fa$ becomes principal, generated by $x_i$. 
\item  There exists a \mlp for $(\xn)$, that is a matrix $A = (a_{ij}) \in \Mn(\gA)$ that satisfies
\begin{equation}\label{eqmlp}
\left\{\arraycolsep2pt
 \begin{array}{rcl}
   \;\;\sum a_{ii}&=&1\\[1mm]
   \;\;a_{\ell j}x_{i}& =& a_{\ell i}x_{j} \qquad \forall i,j,\ell \in
 \lrbn
 \end{array}
\right.
\end{equation}
 Note: The last line is read as follows:
 for each row~$\ell$, the minors of order~$2$ of the matrix $\cmatrix{ a_{\ell 1}&\cdots &a_{\ell n}\cr x_1&\cdots &x_n}$ are null.

\item  $\Vi_\Ae 2(\fa)=0$.
\item  $\cF_1(\fa)=\gen{1}$.
\end{enumerate}
In the case where one of the $x_k$'s is \ndz the existence of the matrix $A$ in item 3 has the same meaning as the following item.
\begin{enumerate}
\item [3$\,'$.] 
There exist $\gamma_1$, \ldots, $\gamma_n$ in $\Frac\gA$ such that $\sum_i\gamma_ix_i=1$ and each of the $\gamma_ix_j$'s is in $\gA$
 \emph{(Dedekind formulation)}.
\end{enumerate}
\end{proposition}
%
\begin{proof}
The only new thing is the formulation \emph{3$\,'$.}
If for example $x_1\in\Reg(\gA)$ and if we dispose of $A$, let $\gamma_i=a_{i1}/x_1$. 
Conversely, if we dispose of the $\gamma_i$'s, let $a_{ij}=\gamma_ix_j$.
\end{proof}

\medskip 
The following proposition takes up and adds details to Proposition~\ref{pmlm}.
The results could be obtained more directly, by using the Dedekind formulation, when one of the $x_k$'s is \ndzz.

\begin{proposition}\label{pilps1}
Let $\fa = \gen{\xn} $ be a \lop \id of~$\gA$ and $A= (a_{ij})$ be a \mlp for $(\xn)$.
We have the following results.
\begin{enumerate}
\item 
$[\,x_1\; \cdots\;x_n\,]\; A = 
[\,x_1\; \cdots\;x_n\,]$.
\item 
   Each $x_i$ annihilates $\cD_2(A)$ and $A^2 - A$.
\item 
  Let $\gA_i = \gA[1/a_{ii}]$, we have $\fa =_{\gA_i} \gen{x_i}$.
\item 
 $\gen{\xn} \gen{a_{1j},\dots,a_{nj}} = \gen{x_j}$.
\item 
 More \gnltz, if $a = \sum\alpha_{i}x_i$ and $\tra{[\,y_1\; \cdots\;y_n\,]} = A~ \tra{[\,\alpha_1\; \cdots\;\alpha_n\,]}$, then

\snic{\displaystyle
\gen{\xn}\gen{\yn} = \gen{a}
.}

 In addition, if $\Ann(\fa) = 0$, the matrix $\tra{\!A}$ is a \mlp for $(\yn)$.

\item 
 In particular, if $ \sum\alpha_{i}x_i = 0$ and $\tra{[\,y_1\; \cdots\;y_n\,]} = A~ \tra{[\,\alpha_1\; \cdots\;\alpha_n\,]}$, then

\snic{\displaystyle
\gen{\xn}\gen{\yn} = 0.
}
\item 
 Consider the \lin form ${\ux}: (\alpha_i)\mapsto \sum_i\alpha_ix_i$ associated with $(\xn)$, let $\fN=\Ann\, \gen{\xn}$ and $\fN^{(n)}$ be the cartesian product

\snic{\sotq{(\nu_1,\dots,\nu_n)} {\nu_i\in \fN,\,i\in\lrbn} \subseteq \Ae n.}

Then $\Ker{\ux}=\Im(\I_n-A)+\fN^{(n)}$.
\item 
For $i\in\lrb{1..n-1}$ the intersection $\gen{x_1,\dots,x_i}\cap\gen{x_{i+1},\dots,x_n}$ is the \id generated by the $n$ \coes of the row vector 
   
\snic{
[\,x_1\;\cdots\;x_i\;0\;\cdots\; 0\,](\I_n-A) =
-[\,0\;\cdots\;0\;x_{i+1}\;\cdots\; x_n\,](\I_n-A).
}
\end{enumerate}
\end{proposition}

\begin{proof}
Item \emph{3} is clear, items \emph{4} and \emph{6} are special cases of the first part of item \emph{5}.\\
Items \emph{1}, \emph{2} and the first part of item \emph{5} have been shown for the \mlmosz.

 \emph{5.} It remains to show that, when $\Ann(\fa)=0$, $\tra{A}$ is a \mlp for~$(\yn)$. Indeed, on the one hand $\Tr(\tra{A})=1$, and on the other hand, since $\fa\cD_2(A) = 0$, we have $\cD_2(A) = 0$, or $A_i \wedge A_j = 0$, $A_i$ being the column $i$ of $A$. As the vector $y := \tra{[\,y_1\; \cdots\;y_n\,]}$ is in~$\Im A$, we also have~$y \wedge A_j = 0$, which translates that~$\tra{A}$ is a \mlp for~$(\yn)$.

\emph{7.} The inclusion $\Ker\ux\subseteq\Im(\I_n-A)+\fN^{(n)}$ results from item \emph{6} and the reverse inclusion of item~\emph{1.}

\emph{8.} Results from \emph{7} by noticing that taking an \elt $a$ of
{the \id   $\fb=\gen{x_1,\dots,x_i}\cap\gen{x_{i+1},\dots,x_n}$} is the same as taking an \elt
$$
\preskip.4em \postskip.4em 
 (\alpha_1,\dots,\alpha_n)\in\Ker\ux~:~a=\alpha_1x_1+\cdots +\alpha_ix_i=-\alpha_{i+1}x_{i+1}-\cdots
-\alpha_nx_n.
$$
Thus, $\fb$ is generated by the \coes of  $[\,x_1\;\cdots\;x_i\;0\;\cdots\; 0\,] (\I_n-A)$.
\end{proof}

\begin{corollary}\label{corpilps1}
Let $\fa = \gen{\xn}$ be a \itf of $\gA$.  
\begin{enumerate}
\item If $\fa$ is \lopz, for every \itf $\fc$ contained in~$\fa$, there exists a \itf $\fb$ such that $\fa\fb=\fc$.
\item Conversely, if $n=2$ and if there exists some \id $\fb$ such that $\gen{x_1}=\fa\fb$, then $\fa$ is \lopz.
\item The \id $\fa$  is a \mrc $1$ \ssi it is \lop and faithful.
In this case, if $A$ is a \mlp for $(\xn)$, it is a \mprn of rank $1$ and $\fa\simeq \Im A$.
\item The \id $\fa$ is \iv \ssi it is \lop and contains a \ndz \eltz.
\end{enumerate}
\end{corollary}
\begin{proof}
\emph{1}, \emph{3},  \emph{4.} See Lemma~\ref{lemIdproj}, which gives slightly more \gnl results.
These items also result from the previous proposition, items \emph{5} and \emph{7.} 
\\
\emph{2.} In $\fb$ we must have $u_1$ and $u_2$ such that on the one hand $u_1x_1+u_2x_2=x_1$, %
so $(1-u_1)x_1=u_2x_2$, and on the other hand $u_1x_2 \in \gen{x_1}$. 
 When we invert the \eltz~$u_1$,~$x_1$ generates $\fa$, and when we invert $1-u_1$, it is $x_2$ that generates $\fa$.
\end{proof}
%

\subsec{First \prtsz}

Recall that a \ri is \coh \ssi on the one hand the intersection of two \itfs is a \itfz, and on the other hand the annihilator of every \elt is \tf (\thref{propCoh4}).  Consequently, by using item \emph{8} of Proposition~\ref{pilps1}, we obtain

\begin{fact} 
\label{factAnarCoh} 
In an \anar the intersection of two \itfs is a \itfz. An \anar is \coh \ssi 
the annihilator of every \elt is \tfz.
\end{fact}

Every quotient and every localized \ri of an \anar is an \anarz.

In a \fdi \riz, the \dve relation is explicit.
We have the (remarkable) converse for  \anarsz.
\pagebreak	

\begin{proposition}\label{aritfdi}
An \anar is \fdi \ssi the \dve relation is explicit.
More precisely, in an arbitrary \riz, if an \id $\gen {b_1, \ldots, b_n}$ is \lop and if $A = (a_{ij})$ is a \mlp for $(b_1,\dots,b_n)$, we have the \eqvc

\snic {
c \in \gen {b_1, \dots, b_n}  \iff  a_{jj}c \in \gen {b_j}
\hbox { for every } j.
}

In particular, we have $1\in\gen{b_1,\dots,b_n}$ \ssi for all
$j$, $b_j$ divides $a_{jj}$.
\end{proposition}

\begin{proof}
If $a_{jj}c = u_jb_j$, then $c = \sum_j u_jb_j$. Conversely, if $c \in \gen {b_1, \dots, b_n}$, then for each $j$ we get
$$
\preskip-.2em \postskip.4em 
a_{jj}c \in \gen {a_{1j},\dots,a_{nj}} \gen {b_1, \dots, b_n} = \gen {b_j}. 
$$

\vspace{-1em}
\end{proof}

In the following \tho we give a few possible \carns of \anarsz.
The simplest \carn of  \anars is no doubt the one given in item \emph{1b}. Since an \id $\gen{x,y}$ is \lop \ssi there is a \mlp for $(x,y)$, condition \emph{1b} means
\Grandcadre{$\forall x,y\in\gA\;\;\exists u,a,b\in\gA,\quad\quad ux=ay,\;(1-u)y=bx$,}
which is also exactly what item \emph{2c} says.

\begin{theorem} 
\label{thAnar} \emph{(\Carns of  \anarsz)}\\
For a \ri $\gA$ \propeq
\begin{enumerate}
\item [1a.] $\gA$  is \ari (every \itf is \lopz).
\item [1b.] Every \id $\fa=\gen{x_1,x_2}$ is \lopz.
\item [2a.] For all \itfs $\fb\subseteq \fa$, there exists some \itf $\fc$ such that $\fa\fc=\fb.$
\item [2b.] For every \id $\fa=\gen{x_1,x_2}$, there exists some \itf $\fc$ such  
that~$\fa\fc=\gen{x_1}$.
\item [2c.] $\forall x_1,x_2\in \gA$ the following \sli $BX=C$ admits a solution
\begin{equation}\label{eqSLI}
[ \,B\mid C \,] =\cmatrix{ ~x_1 & x_2  &  0    &\vert&x_1 \cr
                                ~x_2 & 0    &  x_1  &\vert&~0~~}   
\end{equation}
\item [2d.] $\forall x_1,x_2\in \gA$ there exists a $u\in \gA$ such that
$$ \gen{x_1}\; \cap \; \gen{x_2}\; =\; \gen{(1-u)x_1,ux_2}.
$$
%
\item [3.] For all \itfs $\fa$ and $\fb$, the following short exact sequence is split
$$ 0\longrightarrow \gA/(\fa\cap \fb) \vers{\delta}\gA/\fa\times 
\gA/\fb\vers{\sigma}
\gA/(\fa+\fb)\longrightarrow 0 
$$
where $\delta : \ov x_{\fa\cap\fb} \mapsto (\ov x_\fa, \ov x_\fb)$ and 
$\sigma : (\ov y_\fa, \ov z_\fb) \mapsto \ov {(y-z)}_{\fa+\fb}.$
\item [4.] For all \itfs $\fa$ and $\fb$, $(\fa:\fb)+(\fb:\fa)=\gen{1}$. 
\item [5.] {\em (Chinese remainder \thoz, \ari form)}\\ 
If $(\fb_k)_{k=1,\ldots,n}$ is a finite family of \ids of $\gA$ and $(x_k)_{k=1,\ldots,n}$ is a family of \elts of $\gA$ satisfying $x_k\equiv x_\ell\; \mod\,\fb_k+\fb_\ell$ for %
all $k,\,\ell,$ then there exists some $x\in\gA$ such that  $x\equiv x_k\; \mod\,\fb_k$ for all~$k.$   
\item [6.] The lattice of \ids of $\gA$ is a \trdiz. 
\end{enumerate}
\end{theorem}
\begin{proof}
\emph{1b} $\Rightarrow$ \emph{1a}.
If we have a \itf with $n$ \gtrsz, successive \lons (each time at \ecoz) make it principal.\\
Consider item \emph{2a}. Let $\fa=\gen{x_1,\ldots ,x_n}$ and $\fb=\gen{y_1,\ldots ,y_m}$.
If $\fc$ exists, for each $j=1,\ldots,m$ there exist \elts $a_{i,j}\in \fc$ such that 

\snic{\som_i a_{i,j} x_i =y_j.}


Moreover, for each $i,i',j$ we must have $a_{i,j}x_{i'}\in \fb$, which is expressed by the existence of \elts $b_{i,i',j,j'}\in \gA$ satisfying 

\snic{\som_{j'}b_{i,i',j,j'}y_{j'}=a_{i,j}x_{i'}.
}


Conversely, if we can find some \elts $a_{i,j}$ and $b_{i,i',j,j'}\in \gA$ satisfying the \lin \eqns above (in which the $x_i$'s and $y_j$'s are \coesz), then the \id
$\fc$ generated by the $a_{i,j}$'s indeed satisfies $\fa\fc=\fb$.
Thus, finding $\fc$ comes down to solving a \sliz.\\
It follows that to prove \emph{1a} $\Rightarrow$ \emph{2a}
we can use suitable \lonsz:  the two \ids $\fa$  and $\fb$ become principal, one being included in the other, in which case $\fc$ is obvious.

We easily verify that the \prts \emph{1b}, \emph{2b}, \emph{2c} and  \emph{2d} are \eqves (taking into account the previous remark for \emph{1b}).

To show that \emph{1a} implies \emph{3}, \emph{4}, \emph{5} and  \emph{6}, note that each of the \prts considered can be interpreted as the existence of a solution of a certain \sliz, and that this solution is obvious when the \ids that intervene are principal and totally ordered for the inclusion.

It remains to show the converses.

\emph{3} $\Rightarrow$ \emph{2c}  and \emph{4} $\Rightarrow$ \emph{2c}. \\
Consider in \emph{3} or \emph{4} the case where $\fa=\gen{x_1}$ 
and $\fb=\gen{x_2}$.
 
\emph{5}  $\Rightarrow$ \emph{1b}. 
Let $a$, $b\in \gA$. Let 

\snic{c=a+b, \; 
\fb_1=\gen{a},\;\fb_2=\gen{b}, \; \fb_3=\gen{c}, \; x_1=c, \;x_2=a\hbox{ and }x_3=b.}

We have $\fb_1+\fb_2=\fb_1+\fb_3=\fb_3+\fb_2=\gen{a,b}$.\\
The congruences $x_i\equiv x_k\; \mod\,\fb_i+\fb_k$ are satisfied, so there exist $u$, $v$, $w$ in~$\gA$ such that 
$$\preskip.0em \postskip-.2em 
c+ua = a+vb = b+wc, 
$$
hence
$$\preskip.0em \postskip.4em 
wb=(1+u-w)a,\; (1-w)a=(1+w-v)b. 
$$
Therefore the \id $\gen{a,b}$ is \lopz.

\emph{6}  $\Rightarrow$ \emph{1b}. 
Take the \prt of \dit $\fa+(\fb\cap \fc)=(\fa+\fb)\cap (\fa+\fc)$, with  $\fa=\gen{x}$, $\fb=\gen{y}$ and $\fc=\gen{x+y}$. We therefore have $y\in \gen{x}+(\gen{y}\cap\gen{x+y})$, that is,  there exist $a$, $b$, $c$ such that $y=ax+by$, $by=c(x+y)$. 
Hence $cx=(b-c)y$ and $(1-c)y=(a+c)x$. Thus, $\gen{x,y}$ is \lopz.
\end{proof}

The \iso $\gA/\fa\oplus \gA/\fb\simeq \gA/(\fa+\fb)\oplus \gA/(\fa\cap\fb)$
which results from item~\emph{3} of the previous \tho admits the following \gnnz.
\begin{corollary}\label{corthAnar}
Let $(\fa_i)_{i\in\lrbn}$ be a family of \itfs of an \anar $\gA$.
Let
\[\preskip.4em \postskip.4em 
\begin{array}{ccc} 
  \fb_1=\sum_{k=1}^n\fa_k,\; \fb_2=\sum_{1\leq j<k\leq n}(\fa_j\cap \fa_k),\;\dots  \\[.3em] 
\fb_r=\sum_{1\leq j_1<\cdots<j_r\leq n}(\fa_{j_1}\cap \cdots\cap \fa_{j_r}),\;\dots,\;
  \fb_n=\bigcap_{k=1}^n\fa_k.  
  \end{array}
\] 
Then we have $\fb_n\subseteq\cdots\subseteq \fb_1$ with an \iso
$$\preskip.3em 
\bigoplus\nolimits_{k=1}^n\gA/\fa_k \;\simeq \;\bigoplus\nolimits_{k=1}^n\gA/\fb_k.
$$
\end{corollary}

By bringing this result closer to \thref{prop unicyc} we obtain a complete classification of \Amos of this type.
We can also compare with Fact~\ref{factGpRtcl}~\emph{\iref{i16bisfactGpRtcl}.}

\begin{corollary}\label{cor2thAnar}
Let $\gB$ be a \fpte \Algz. If $\gB$ is an \anar (resp.\,a \adpz, a \adpcz), then so is $\gA$. 
\end{corollary}
%
\begin{proof} Since $\gA\subseteq\gB$, if $\gB$ is reduced, so is $\gA$.
\Thref{propFidPlatTf}~\emph{3} implies that if $\gB$ is \cohz, so is $\gA$. It remains to show the result for an \gui{\anarz.} 
Consider $x$, $y\in\gA$. We must show that there exist $u$, $a$, $b\in\gA$ such %
that $ux=ay$ and $(1-u)y=bx$.
This is actually a \sli with \coes in $\gA$, with the unknowns $(u,a,b)$. 
However, this \sys admits a solution in $\gB$ and~$\gB$ is \fpte over $\gA$, so it admits a solution in $\gA$. 
\end{proof}
%

\subsec{Multiplicative structure of \itfsz}

Recall that we denote by $\Ifr\gA$ the multiplicative \mo of \tf \ifrs of an arbitrary \ri $\gA$ (see \paref{NOTAIfr}). 

A priori an inclusion $\fa\subseteq\fb$ in $\Ifr\gA$ does not imply the existence of a \ifr $\fc\in \Ifr\gA$ such that $\fb\fc=\fa$. But this is satisfied in the case of \anarsz.

\rdb
For $\fa$ and $\fb$ in $\Ifr\gA$, let \fbox{$\fa\div\fb=\sotq{x\in\Frac\gA}{x\fb\subseteq\fa}$}.\label{NOTAfadivfb} 

\begin{lemma}\label{lemIfrCoh} 
Let $\gA$ be a \emph{\coriz}. 
\begin{enumerate}
\item $\Ifr\gA$ is a lattice \wir inclusion relation, the supremum is given by the sum and the infimum by the intersection. 
\item  $\Ifr\gA$ is a \trdi \ssi the \ri is \ariz.
\item  Concerning \iv \elts of $\Ifr\gA$.
\begin{enumerate}
\item If $\fa\,\fa'=\gA$ in $\Ifr\gA$,  we have $\fa'\fc=\fc\div \fa$ and $\fa(\fc\div \fa)=\fc$ for \hbox{all $\fc\in\Ifr\gA$}. In particular $\gA\div\fa$ is the inverse of $\fa$.
\item A \ifr $\fraC \fa a$ (where $\fa$ is a \itf of $\gA$) is \iv in $\Ifr\gA$ \ssi  $\fa$ is an \iv \idz.
\item If $\fa(\gA\div \fa)=\gA$, $\fa$ is \iv in $\Ifr\gA$.
\end{enumerate}
\end{enumerate}
Let $\fa$, $\fb\in\Ifr\gA$ with $b\in\fb\cap\Reg\gA$. 
Suppose that $\gA$ is \icl in~$\Frac\gA$.
\begin{enumerate}\setcounter{enumi}{3}
\item We have $\fa\div\fb\in\Ifr\gA$.  
\item If in addition $\fa\subseteq\fb\subseteq\gA$, then we have $\fa\div\fb=\fa:\fb$. 
\end{enumerate}
\end{lemma}

\begin{proof}
Every \elt of $\Ifr\gA$ is written in the form $\fraC \fa a$ for some \itfz~$\fa$ of~$\gA$ and some $a\in\Reg\gA$. In addition $\fraC \fa a\, \fraC \fb b =\fraC{\fa\,\fb}{ab}$. Finally, the neutral \elt of the \mo is $\gA=\gen{1}$.
This shows items \emph{1}, \emph{2} and \emph{3b.}

\emph{3a.} We have $\fa\fa'\fc=\fc $ so $\fa'\fc\subseteq \fc\div\fa$ and $\fc=\fa\fa'\fc\subseteq\fa(\fc\div\fa)=\fc$. 
\\
If $x\in\fc\div\fa$, \cad $x\fa\subseteq \fc$, then $x\gA=x\fa\fa'\subseteq \fa'\fc$, so $x\in\fa'\fc$.

\emph{3c.} 
With $\fa=\gen{a_1,\dots,a_k}\subseteq \gA$, suppose that $\fa(\gA\div \fa)=\gA$. \\
There exist $x_1$, \dots, $x_k\in (\gA\div \fa)$ such that $\sum_ix_ia_i=1$ and $x_ia_j\in\fa$ for \hbox{all $i,j$}. 
We can write the $x_i$'s in the form $\fraC{b_i}c$ with the same denominator~$c$.
We obtain $\sum_ia_ib_i=c$ and $a_ib_j\in\gen{c}$ for all $i,j$.
\\
Thus by letting $\fb=\gen{b_1,\dots,b_k}$ we obtain $\fa \,\fb=\gen{c}$.

\emph{5.} The inclusion $\fa:\fb\subseteq \fa\div \fb$ is \imdez. 
Conversely, if some $x\in\gK$ satisfies~\hbox{$x\fb\subseteq \fa$}, we need to show that $x\in\gA$.\\
As $\gA$ is \icl in $\Frac\gA$, we apply item \emph{3} of Fact~\ref{factEntiersAnn}, \hbox{with $M=\fb$} \hbox{and $\gB=\Frac\gA$}, \hbox{because $x\fb\subseteq \fa\subseteq \fb$}.

\emph{4.} Results from item \emph{5} because we are brought back to the case treated in item~\emph{5}, and in a \coriz, the conductor $\fa:\fb$ is \tf if $\fa$ and $\fb$ is \tfz.
\end{proof}

The following \tho says that the multiplicative structure of the \mo of \iv \ids of an \anar has all the desired \prtsz. 

Recall that by Lemma~\ref{lemIdproj}, a \itf is \prc $1$ \ssi it is \lop and faithful.  

\begin{theorem}\label{thiivanar}
In an \anar the faithful \itfs form  the non-negative sub\mo of an \grlz.
The lattice laws are $\fa\vi\fb=\fa+\fb$ and $\fa\vu\fb=\fa\cap\fb$. 
\\
The \iv \ids (\cad the \itfs that contain a \ndz \eltz) form the non-negative sub\mo of an \sgrl of the previous \grlz.  
\end{theorem}
\begin{proof}
This results from Corollary~\ref{corpilps1}, from \thref{thAnar} and from \thref{lem1MonGcd}.
\end{proof}

Actually the two groups  coincide as soon as $\gA$ is a \qiriz, or more \gnlt when the \mrcs $1$ over $\Frac\gA$
are free (\thref{propRgConstant3}, item \emph{2}).

\section{Integral \elts and \lonz} 
\label{subsecEntiers}

The following \dfn generalizes \Dfnz~\ref{defEntierAnn0} in two directions.

\begin{definition} 
\label{defPropACO} Let $\varphi:\gA\to\gC$ be a \homo between commutative \ris and $\fa$ be an \id of $\gA$. 
\begin{enumerate}
\item  An \elt $x\in\gC$  is said to be \ixc{integral}{element over an ideal} over $\fa$ if there exists an integer $k\geq 1$  such that 
$$\preskip-.3em \postskip.3em 
x^k=\varphi(a_1)x^{k-
1}+\varphi(a_2) x^{k-2}+\cdots +\varphi(a_k) \eqno(*)
$$
with each~$a_h\in\fa^h$. 
In the case where $\gC=\gA$, this is equivalent to $\big(\fa+\gen{x}\!\big)^k=\fa\big(\fa+\gen{x}\!\big)^{k-1}$.
We also say that the \egt $(*)$  is \emph{an \rdiz} of $x$ over $\fa$.%
\index{dependence relation!integral ---}

\item  An \id $\fa$ of $\gA$ is said to be \ix{integrally closed} in~$\gC$ if every \elt of~$\gC$ integral over $\fa$  is in~$\varphi(\fa)$.%
\index{ideal!integrally closed ---}
\item  The \ri $\gA$ is said to be \ixc{normal}{ring} if every \idp of $\gA$ is \icl in $\gA$.%
\index{ring!normal}
\end{enumerate}
\end{definition}

In all  cases, a normal \ri is \icl in its total \ri of fractions. We have the following partial converse.

\begin{fact}\label{lemNormalIcl}
A \qiri is normal \ssi it is \icl in its total \ri of fractions.
\end{fact}
%
\facile

It is clear that every normal \ri is reduced (because a nilpotent is integral over $\gen{0}$). We even have a little better.

\begin{lemma} 
\label{lemiclplat} 
Every normal \ri is a \lsdz.
More \prmtz, we have for every \ri  $\gA$  the implications  
{1} $\Rightarrow$ {2} $\Rightarrow$ {3.}
\begin{enumerate}
\item  Every \idp is \icl (\cad $\gA$ is normal).
\item  For all $x$, $y\in \gA,\; \; x^2\in \gen{xy} $ implies $x\in \gen{y}.$
\item  Every \idp is flat (\cad $\gA$ is a \lsdz).
\end{enumerate}
\end{lemma}
\begin{proof}
Note that the \id $0$ is \icl \ssi the \ri is reduced. 
We obviously have \emph{1} $\Rightarrow$ \emph{2}, and \emph{2} implies that the \ri is reduced.
Suppose~\emph{2} and let $x$, $y\in \gA$ such that $xy=0$. We have $x^2 =x(x+y)$ 
therefore $x\in \gen{x+y}$, e.g. $x=a(x+y)$. 
Then $(1-a)x=ay$, $ay^2=(1-a)xy=0$, 
and since the \ri is reduced, $ay=0$,
then $(1-a)x=0$. 
\end{proof}

\vspace{-.5em}
\begin{fact}\label{factEntiers}
Let $x$ be an \elt and $\fa$ be an \id of $ \gA$.
For the \prts that follow we have {2} $\Rightarrow$ {1}, and {1} $\Rightarrow$ {2} if $\fa$ is faithful and \tfz. 
\begin{enumerate}\itemsep0pt
\item
The \elt $x$ is integral over the \idz~$\fa$. 
\item 
There exists a \tf faithful \Amo $M$ such that \linebreak {$xM \subseteq \fa M$}. 
\end{enumerate}
\end{fact}
\begin{proof} (Compare with the \dem of Fact~\ref{factEntiersAnn}.)\\
\emph{2} $\Rightarrow$ \emph{1}. 
Consider a matrix $A$ with \coes in $\fa$ that represents~$\mu_{M,x}$ (multiplication by $x$ in $M$) over a finite \sgr of~$M$. If~$f$ is the \polcar of $A$, we have by the Cayley-Hamilton \tho $0=f(\mu_{M,x})=\mu_{M,f(x)}$, and since the module is faithful, $f(x)=0$.
 
\emph{1} $\Rightarrow$ \emph{2}. If we have an \rdi of degree $k$ of $x$ over $\fa$ we take $M = (\fa+\gen{x})^{k-1}$.
\end{proof}

\rdb

Let $\fa$ be an \id of $\gA$ and $t$ be an \idtrz. Then the sub\alg $\gA[\fa t]$  of $\gA[t]$, \cad \prmt
$$\preskip.1em \postskip.4em
\gA[\fa t] = \gA \oplus \fa t \oplus \fa^2 t^2 \oplus \cdots
$$
is called the \emph{Rees \alg of the \id $\fa$}.%
\index{Rees algebra!of the ideal $\fa$}%
\label{NOTARees}

The \dem of the two following facts is left to the reader.
\begin{fact}\label{fact2Entiers}
Let $\fa$ be an \id of $ \gA$. 
\begin{enumerate}\itemsep0pt
\item For $x \in \gA$, \propeq
\begin{enumerate}\itemsep0pt
\item
The \elt $x $ is integral over the \id $\fa$ of $\gA$. 
\item 
The \pol $xt$ is integral over the sub\alg $\gA[\fa t]$ of $\gA[t]$. 
\end{enumerate}\itemsep0pt
\item More precisely
\begin{enumerate}\itemsep0pt
\item  If $\ov \fa$ is the set of \elts of $\gA$ integral over $\fa$, then the \cli of $\gA[\fa t]$ in $\gA[t]$ is the sub\ri $\gA[\ov\fa t]$.%
\index{integral closure!of the ideal $\fa$ in $\gA$}
\item In particular,  $\overline \fa$ is an \id of $\gA$, called the \emph{\cli of the \id $\fa$ in $\gA$}.
We denote it by $\Icl_\gA(\fa)$ or $\Icl(\fa)$.%
\end{enumerate}
\end{enumerate}
\end{fact}

\vspace{-.5em}
\begin{fact}\label{fact3Entiers}
Let $\fa$ and $\fb$ be two \ids of $\gA$.
\begin{enumerate}\itemsep0pt
\item $\Icl\big(\Icl(\fa)\big)=\Icl(\fa)$.
\item $\fa\Icl(\fb)\subseteq\Icl(\fa)\Icl(\fb)\subseteq\Icl(\fa\fb).$
\end{enumerate}
\end{fact}

We now revisit two important results which have been already established.
Item \emph{2c} of \KROz's \tho \ref{thKro} gives \prmt the following result.

\begin{lemma}\label{lemthKroicl}\emph{(\KROz's \thoz, reformulated)}\\
Suppose that we have in $\AT$ an \egt

\snic{f=\som_{i=0}^nf_iT^{n-i},\;\;  g=\som_{j=0}^m g_j T^{m-j} \;\; \hbox{and} \;\; h=fg=\som_{r=0}^{m+n} h_rT^{m+n-r}.
}

Let $\gk$ be the sub\ri of $\gA$ generated by the~$f_ig_j$'s.
Then, each $f_ig_j$ is integral over the \id $\rc_{\gk}(h)$ of $\gk$.
\end{lemma}

Note that item \emph{2c} of \KROz's \tho \ref{thKro} tells us \prmt this:
\emph{there exists some \pog $R_{i,j}\in \ZZ[Y,H_0,\ldots,H_p]$
 (all the variables have the same weight $1$), monic in $Y$, such that}
 
 \snic{R_{i,j}(f_ig_j,h_0,\dots,h_p)=0.}

\smallskip 
Here is a new version of the Lying Over (Lemma~\ref{lemLingOver}).%
\index{lying over}

\pagebreak

\begin{lemma} \emph{(Lying Over, more precise form)}\label{lemLingOver2}\\
 Let  $\gA\subseteq\gB$ with $\gB$ integral over $\gA$ and $\fa$ be an \id of $\gA$, then $\fa\gB\cap\gA\subseteq\DA(\fa)$.
More \prmtz, every \elt of $\fa\gB$ is integral over $\fa$.
\end{lemma}
\begin{proof}
We textually rework the proof of Lemma~\ref{lemLingOver}. If $x\in\fa\gB$, we~have~$x=\sum a_ib_i$, with $a_i\in \fa,\;b_i\in \gB$.  The $b_i$'s generate an \Aslgz~$\gB'$ which is finite.  
Let $G$ be a finite \sgr (with $\ell$ \eltsz) of the~\Amo $\gB'$.  
Let $B_i\in\MM_\ell(\gA)$ be a matrix that expresses the multiplication by~$b_i$ over $G$.  
The multiplication by~$x$ is expressed by the matrix $\sum a_iB_i$, which is with \coes in $\fa$. The \polcar of this matrix, which annihilates~$x$ (because $\gB'$ is a faithful \Amoz), therefore has its \coe of degree $\ell-d$ in~$\fa^d$.

We could \egmt apply Fact~\ref{factEntiers} by taking $M = \gB'$. Indeed, as $x \in \fa\gB'$, we have $x\gB' \subseteq \fa\gB'$ and so $x$ is integral over $\fa$.
\end{proof}

We now examine the relationships between \prts of the type \gui{integral over} and \lonsz.

\begin{fact} 
\label{fact.loc.normal} 
Let $\fa$ be an \id of $\gA$,  
$S$ be a \mo of $\gA$ and $x\in\gA$.
\begin{enumerate}
\item   The \elt $x/1\in\gA_S$ is integral over $\fa_S$ \ssi there exists a $u\in S$ such that $xu$ is integral over $\fa$ in $\gA$.

\item  If $\gA$ is normal, then so is $\gA_S$.

\end{enumerate}
Let $\gB\supseteq\gA$ be a \fpte \algz.
\begin{enumerate}\setcounter{enumi}{2}
\item  If $\gA'$ is the \cli of $\gA$ in $\gB$,
then $\gA'_{S}$ is the \cli of $\gA_S$ in~$\gB_{S}$.

\item  If $\gB$ is normal, then $\gA$ is normal.
\end{enumerate}
\end{fact}
\begin{proof} We only prove item \emph{1.}
In the \dem we confuse an \elt of $\gA$ and its image in $\gA_S$ to alleviate the notation.
If an \egt $x^k=a_1x^{k-1}+ a_2 x^{k-2}+ \cdots +a_k$ is performed in $\gA_S$ with each $a_j\in(\fa\gA_S)^j$,
we obtain \gui{by reducing all the fractions to the same \deno and by getting rid of the \denoz} an \egt 

\snic{sx^k=b_1x^{k-1}+ b_2 x^{k-2}+ \cdots +b_k}

in $\gA_S$ with $s\in S$ and each $b_j\in\fa^j$. This means an \egt in $\gA$ after multiplication by another \eltz~$s'$ of~$S$. We can also multiply by $s'^ks^{k-1}$ and we obtain with~$u=ss'$ an \egt  

\snic{(xu)^k=c_1(xu)^{k-1}+ c_2 (xu)^{k-2}+ \cdots +c_k}

in~$\gA$ with each $c_j\in\fa^j$.
\end{proof}

The fact that a \ri is normal  is a local notion, in the following sense.

\pagebreak	

\begin{plcc} 
\label{plcc.normal} \emph{(Normal \risz)}\\
Let $S_1$, $\ldots$, $S_n$  be \moco of a \ri $\gA$, $x\in \gA$ and $\fa$ be an \id of $\gA$.   
\begin{enumerate}
\item  The \elt $x$ is integral over $\fa$  \ssi it is integral over each of the~$\fa_{S_i}$'s.
\item  The \id $\fa$ is \icl in $\gA$ \ssi each of the~$\fa_{S_i}$'s is \icl in~$\gA_{S_i}$.
\item  The \ri $\gA$ is normal \ssi each of the $\gA_{S_i}$'s is normal.
\end{enumerate}
\end{plcc}
\begin{proof}
It suffices to prove item \emph{1}, the passage from the local to the global.
\perso{on devrait le laisser to the reader?}
We obtain by applying Fact~\ref{fact.loc.normal} for each $i\in\lrbn$ some $s_i\in S_i$ such that $s_ix$ is integral over the \id $\fa$ in $\gA$.
We can suppose that all the \rdis have the same degree $k$.
Let us write these \rdis
$$
\preskip.0em \postskip.4em
(s_ix)^k \in \som_{h=1}^k \fa^h (s_ix)^{k-h},
\qquad i\in\lrbn.
$$  
A \coli of these relations based on an \egt $\som_{i=1}^nb_is_i^k=1$ gives us an \rdi of $x$ over $\fa$ in~$\gA$.
\end{proof}

Note that since the \prt in item \emph{1} is of \carfz, Lemma \ref{factPropCarFin} says that
the previous \plgc is equivalent in \clama to the corresponding \plga (in which the \lon intervenes at any \idema of $\gA$).

\section{Pr\"ufer \risz} 
\label{secAdP}

Recall that a \ri is said to be Pr\"ufer when its \ids are flat, or if it is \ari and reduced, or if it is \ari and a \lsd (Proposition~\ref{prop.itfplat 2}).

\begin{propdef}\label{defiAdv}
We call a \ri $\gA$ satisfying one of the following \eqv \prts a \ixx{valuation}{ring}.%
\index{ring!valuation ---}
\begin{enumerate}
\item $\gA$ is a reduced local Bézout \riz.
\item $\gA$ is a local \adpz.
\item $\gA$ is reduced and satisfies: for all $a$, $b\in\gA,$ $a\divi b$ or $b \divi a$.
\end{enumerate}
If $\gK=\Frac\gA$, the quotient group $\gK\eti/\Ati$ is equipped with the total order relation $\ov x \divi \ov y$ defined by $\exists a\in\Reg(\gA)$, $y=ax$. This totally ordered group is called the \ixx{valuation}{group} of $\gA$.%
\index{group!valuation ---}    
\end{propdef}

In addition, $\gA$ is then \sdzz.

\medskip
\exl \label{exempleadvgroupe}
Let $\gk$ be a nontrivial \cdi and $(\Gamma,\cdot,1_{\Gamma})$ be a totally ordered discrete group.
We construct a \klg which is a \ddv with $\Gamma$ as its valuation group as follows. First of all consider the \klg $\gA=\gk[\Gamma^{+}]$ described in Exercise \ref{exoAlgMon}.
\\
For an \elt $a=\sum_{i}a_i\gamma_i$ of $\Atl$ we define $v(a)$ as the smallest $\gamma_i$ that intervenes in the expression of $a$ (we have taken the pairwise distinct $\gamma_i$'s, and $a_i\neq 0$). We then prove that $v(ab)=v(a)v(b)$, which implies that $\gA$ is integral. We also let $v(0)=+\infty$.
Finally, our \adv is the sub\ri $\gV=\sotq{\fraC a b}{a\in\gA, b\in\Atl, v(a)\geq v(b)}$ of $\Frac\gA$.
\eoe

\medskip 
We now give a few other \cara \prts of \adpsz, which add to those that we can obtain from \thref{thAnar} for \anarsz.

\begin{theorem} 
\label{thPruf} \emph{(\Carns of  \adpsz)}\\
For some \ri $\gA$ \propeq
\begin{enumerate}
\item [1a.] $\gA$ is an \ari \lsd (\cad a \adpz).
\item [1b.] $\gA$ is a \lsd and for all $x$, $y$ there exist $n\in\NN\etl$ and an \id $\fb$ such that $\gen{x,y}\fb=\gen{x^n}$.
\item [2a.] Every submodule of a flat \Amo is flat.
\item [2b.] $\gA$ is  a \lsd and every torsion-free module is flat.
\item [3a.] An arbitrary \sli $BX=C$, as soon as the \idds of $[\,B \mid C\,]$ are equal to those of $B$, admits a solution. 
\item [3b.] Likewise if we  limit ourselves to $B\in \Ae {2\times 3}$ and $C\in \Ae {2\times 1}$.
\item [4a.] Every \id is \iclz.
\item [4b.] Every \itf is \iclz. 
\item [4c.] Every \id $\gen{x,y}$ is \iclz.
\item [4d.] $\gA$ is normal and for all $ x,y\in\gA$, we have $xy\in\gen{x^2,y^2}$.
\item [5a.] If $\fa$, $\fa'$ and $\fc$ are \itfsz, we have the implication

\snic{ \fa+\fa'\subseteq\fc$, 
 $\;\fa\fc\subseteq\fa'\fc$ $\;\Longrightarrow\;$  $\fa\subseteq\fa'.}
\item [5b.] If $\fa$, $\fa'$ and $\fc$ are \itfsz, we have the implication

\snic{ \Ann(\fa+\fa')\supseteq\Ann(\fc)$, 
  $\;\fa\fc\subseteq\fa'\fc$ $\;\Longrightarrow\;$  $\fa\subseteq\fa'.}

\end{enumerate}
\end{theorem}
\begin{proof} We first take care of \eqvcs between \emph{1}, \emph{2} and \emph{3}.
\\
The implications \emph{1a} $\Rightarrow$ \emph{1b},  \emph{2a} $\Rightarrow$ \emph{1a} and  
\emph{3a} $\Rightarrow$ \emph{3b} are obvious.

\emph{1b} $\Rightarrow$ \emph{1a.} Results from Lemma~\ref{lemleszlop} below.

\emph{3b} $\Rightarrow$ \emph{1a.} The \ri is \ari because the \sli (\ref{eqSLI}) in \thref{thAnar} admits a solution. 
In addition, the \ri is reduced: if~$a^2=0$, the \sli $\so{\,ax=0,\,0x=a\,}$ admits a solution because it corresponds~to 
$$  
   B=\cmatrix{a&0&0\cr0&0&0},\,\,C = \cmatrix{0\cr a}\,\,\,\,
\mathrm{with}\,\,\cD_2([\,B\,\vert\, C\,])=\cD_2(B)=0\,\,!
$$

\emph{1a} $\Rightarrow$ \emph{3b}. First of all suppose that the \ri is local.
Therefore the \ri is \sdz and every \itf is principal. 
Then, \trf by Lemma~\ref{lem.solsli} below. 
In the \gnl case, the proof of the lemma can be reproduced after \lons at suitable \mocoz, 
and since this is a matter of solving a \sli the basic \plg applies.\iplg

\emph{2b} $\Rightarrow$ \emph{2a.} A flat module is torsion-free (Lemma \ref{lem.platsdz}). 
Every submodule of a torsion-free module is torsion-free, therefore flat.

\emph{1a} $\Rightarrow$ \emph{2b.} 
Let $M$ be a torsion-free module over a \adpz.
We want to show that it is flat.
Suppose first of all that the \ri is local. \\
Let $LX=0$ be a syzygy with $L=[\,a_1\;\cdots\; a_m\,]$ in $\gA$ and $X\in M^{m\times 1}$. \Spdgz, suppose \hbox{that $a_i=b_ia_1$} \hbox{for $i>1$}. 
The syzygy is rewritten as $a_1y=0$ %
with $y=x_1+b_2x_2+\cdots+b_mx_m$. 
The cyclic submodule $\gA y$ is flat and the \ri is local therefore $a_1=0$ or $y=0$. In the first case $L=0$. In the second case $X=HX$ and $LH=0$ with the following triangular matrix $H$
$$ \preskip.4em \postskip.4em H=
\cmatrix{ 
   0 &   -b_2   &   -b_3   &  \ldots   &   -b_m   \cr 
   0 &    1     &    0     &  \ldots   &   0   \cr 
\vdots&    \ddots      &  \ddots    &    \ddots  &   \vdots   \cr 
\vdots &       &  \ddots    & \ddots    &   0   \cr 
   0 & \ldots &  \ldots    &  0   &    1  \cr 
}.
$$
In the case of an arbitrary \adpz, we repeat the previous reasoning by using the \lons (at \ecoz) which make the \id  
$\gen{a_1,\ldots,a_m }$ generated by one of the $a_i$'s.

We now pass to the \eqvcs between \emph{1}, \emph{4} and \emph{5}.
\\
The implications \emph{4a} $\Leftrightarrow$ \emph{4b} $\Rightarrow$ \emph{4c} $\Rightarrow$ \emph{4d} and \emph{5b} $\Rightarrow$ \emph{5a} are \imdesz.

 \emph{4d} $\Rightarrow$ \emph{1a.} 
The \ri $\gA$ is a \lsd (Lemma \ref{lemiclplat}). 
It therefore suffices to show that every \id $\fa=\gen{x,y}$ is \lopz.
We have $xy=ax^2+by^2$, and $z=ax$ satisfies $z^2=zy-aby^2$. 
Therefore, since the \ri is normal, $ax=a'y$ for a certain $a'$. 
Similarly, $by=b'x$ for a certain~$b'$.
Therefore, $xy(1-a'-b')=0$.
The \elts $1-a'-b'$, $a'$ and $b'$ are \comz. When we invert $1-a'-b'$, we obtain $xy=0$, and after two new \lonsz, $x=0$ or $y=0$, so $\fa$ is principal. When we invert~$a'$, we obtain $\fa=\gen{x}$ because $a'y=ax$. Likewise when we invert $b'$.
  
 \emph{1a} $\Rightarrow$ \emph{4b}.   
Let $x\in \gA$ be integral over a \itf $\fa$. We have for a certain~\hbox{$n\in\NN,$} $\fa (\fa +\gen{x})^n=(\fa +\gen{x})^{n+1}$. Since the \ri is \ariz, we have a \itf $\fb$ such that  $(\fa +\gen{x})\fb=\gen{x}$. Therefore by multiplying \hbox{by $\fb^n$} we obtain $x^n\fa =x^n(\fa +\gen{x})$ which means that there exists some $y\in \fa$ such that $x^{n+1}=x^ny$, \cadz~\hbox{$x^{n}(y-x)=0$}. 
Since the \ri is a \lsdz, this implies that after \come \lons we have $x=0$ or $y-x=0$, and in each case $x\in \fa$.

\emph{5a} $\Rightarrow$ \emph{4b.} 
Let $x\in \gA$ be integral over some \itf $\fa$. We have for a certain~$n\in\NN,$ $\fa (\fa +\gen{x})^n=(\fa +\gen{x})^{n+1}$. We apply several times the simplification \prt with the \id $\fc=\fa+\gen{x}$ and we obtain at the end of the process~\hbox{$\fa+\gen{x}\subseteq\fa$}.

\emph{4b} $\Rightarrow$ \emph{5b}.
Let $\fc ,\,\fa ,\,\fa' $ be three \itfs satisfying the hypothesis in~\emph{5b}.
Let~$x$ be an \elt of $\fa$ and $X$ be a column vector column formed by a \sgr of $\fc$. Since $x\fc \subseteq \fa'\fc $, there exists a matrix $G \in \Mn(\fa')$ such \hbox{that $xX=GX$}, \cad $(x\In - G) X=0$. If $P$ is the \polcar of $G$, we have on the one hand $P(x)X=0$, and on the other $P(x) \in x^n + \fa'$.
\\ 
Therefore $P(x)\in\Ann(\fc)\subseteq \Ann(\fa +\fa')$ and $P(x)\in \fa +\fa'$. Hence $P(x)^2=0$, \hbox{then
$P(x)=0$}. This is an \rdi of $x$ over $\fa'$. Therefore $x\in \fa'$.
\end{proof}

\begin{lemma}\label{lemleszlop}
In a \lsdz, if %
{we have $\gen{x,y}\fb=\gen{x^n}$} with $n\geq1$, 
then $\gen{x,y}$ is \lopz.
\end{lemma}
\begin{proof}
It suffices to solve this \pb after \come \lonsz. The \lsd \crc of the \ri will be used to manufacture these \lonsz.

We have an \egt $\gen{u,v}\gen{x,y}=\gen{x^n}$ with $x^n=ux+vy$, 
$ux=u_1x^n$, $vx=ax^n$ and $uy=bx^n$.  We get

\snic {
(u_1y - bx)x^n  = 0, \quad 
(u_1x + ay - x)x^n = (ux + vy - x^n)x = 0.
}

We therefore have \come \lons in which $x=0$ and the result is clear. In the latter, $u_1y = bx$ and $u_1x+ay = x$ \cad $(1-u_1)x=ay$.
Thus, $\gen{x,y}$ is \lopz.
\end{proof}
%
\begin{lemma}\label{lem.solsli}
Let $\gA$ be an arbitrary \riz, $B\in \Ae {m\times n}$ and $C \in \Ae {m\times 1}$. The \sli $BX=C$ admits a solution in $\Ae {n\times 1}$ when the following conditions are realized for all $k\in\lrb{1..\inf(m,n)}$
\begin{enumerate}
\item  
The \idd $\cD_k(B)$ is of the form $\delta_k\gA$, where $\delta_k$ is a minor of order $k$.
\item  $\delta_k$ satisfies the condition: 
$ \forall y\in \gA\; \; (y\delta_k=0\; \Rightarrow \; \big(\delta_k=0\; \lor\; y=0)\big)$.
\item  $\cD_k([\,B \mid C\,]) = \cD_k(B)$.
\end{enumerate}
\end{lemma}
\begin{proof}
We begin with $k=\inf(m,n)$. We write the \idt \`a la Cramer 

\snic{\delta_k\times C=\delta_k\,\times $ (a \coli of the columns of $B),}


which results from the nullity of \idds of index $k+1$ and from the fact that $\cD_k([\,B \mid C\,])$ is generated by $\delta_k$. Given \emph{2}, we are in one of the following two cases
\begin{enumerate}[--]
	\item we can simplify by dividing everything by $\delta_k$, so $C \in \Im B$.
	\item $\delta_k=0$, or $k=1$ in which case $C \in \Im B$ (because $B = C = 0$),
or~$k \ge 2$, and we can perform an \recu by replacing $k$ by $k-1$.
\end{enumerate}

\vspace{-1em}
\end{proof}

\vspace{-.4em}
\pagebreak	

\subsec{Extensions of \adpsz}

The fact that a normal \ri is a \lsd means that \lot it behaves like a \ri \sdzz. Actually the machinery of \come \lons at work in the \dfn of a \lsd often allows us to return to the integral case, as we have already seen in the \dem of Lemma~\ref{lemleszlop}.

We have the following important \thoz, which is a \gnn of the analogous result obtained in number theory (\thref{th1IdZalpha}).

\begin{theorem}\label{thExtEntPruf} \emph{(Normal integral extension of a \adpz)}.\\
Let $\gA\subseteq\gB$ with $\gB$ being normal and integral over $\gA$ and $\gA$ being Pr\"ufer.
Then $\gB$ is a \adpz. 
\end{theorem}
\begin{proof}
We will show that every \id  $\gen{\alpha,\beta}$ is \lopz. 

Let us first consider the case of an \id $\gen{a,\beta}$ with $(a,\beta)\in\gA\times\gB$.
We can then reuse almost word for word the \dem \gui{\`a la Dedekind\footnote{This would also work with the \dem \`a la Kronecker.}} of  \thref{th1IdZalpha}.
\\
Let $f \in \gA[X]$ be \mon and vanishing at $\beta$. 
We write $f(X)=(X-\beta)h(X)$ where~$h \in \gB[X]$. We therefore have $f(a X)=(a X-\beta )h(a X)$, which we write as $f_1=g_1h_1$.  
Let $\fc=\rc_\gA(f_1)$,  $\fb=\rc_\gB(h_1)$ and $\fa=\rc_\gB(g_1)=\gen{a,\beta}$.\\
If $\deg(f)=n$, we have $a^n\in\fc$. 
Let $\fc'$ be a \itf of $\gA$ with $\fc\fc'=a^n\gA$.
 \\
By using \KROz's \tho (reformulated in Lemma~\ref{lemthKroicl}), we obtain
$\fc\gB\subseteq\fa \fb \subseteq\Icl_\gB(\fc)$ and so
$$\preskip.4em \postskip.4em \mathrigid 2mu
 a^n\gB=(\fc\gB) (\fc'\gB)\subseteq \fa \fb (\fc'\gB)\subseteq\Icl_\gB(\fc)(\fc'\gB)\subseteq\Icl_\gB(\fc\fc')
=\Icl_\gB(a^n)=a^n\gB.
$$
Therefore $\fa \fb(\fc'\gB)=a^n\gB$ and $\fa$ is \lop by Lemma~\ref{lemleszlop}.

Let us pass to the \gnl case, with $\alpha,\beta\in\gB$. 
If $\gB$ is integral, we can suppose that~$\alpha\neq0$ and we find $\gamma\neq0$ in $\gB$ such that $\alpha\gamma=a\in\gA$, which brings us to the \pb already treated. 
\perso{se pose la question of d\'epenser of l'espace pour la version \gnlez,
renvoyer en exo la fin?}


\smallskip It remains to see the more delicate case where we do not suppose that $\gB$ is integral.
In this case we have $f_i = h_{i-1} - \alpha h_i$ for $i \in \lrbn$
(by convention, $h_{-1}=0$ and~$h_n=0$). We let $\beta_i = - \beta h_i$ for
$i \in \lrb{0.. n-1}$. Then using the first case, we know how to compute
 \mlps $M_i$ for the pairs
 $(f_i, \beta_i) \in \gA \times \gB$~:
$$\preskip0 pt \postskip0pt 
\forall i \in \lrb{0..n-1} \quad
\left\{
\begin{array}{l} \Tr(M_i) = 1 \\[2pt]
 M_i\cdot  {}^t(-\beta_i,f_i) = 0
\end{array}
\right.
$$ 

\noindent
But $(-\beta_i, f_i) \equiv h_i(\beta, -\alpha) \mod h_{i-1}\gB^2$, so
we have
%
\begin{equation}\preskip3pt \label{equation:extensionModulo}
\forall i \in \lrb{0..n-1} \quad
h_i M_i\cdot {}^t(\beta,-\alpha) \equiv 0 \mod h_{i-1} \gB^2.
\postskip2pt
\end{equation}
\noindent
Let us prove by induction that there exist \elts $\zeta_0,\dots,\zeta_{n-1} \in \gB$ and matrices 
$\wi{M_0}, \dots, \wi{M_{n-1}} \in M_2(\gB)$ such that

\hfil
$\displaystyle
\forall i \in \lrb{0..n-1} \quad
\left\{
\begin{array}{l}
\Tr(\wi{M_i}) = 1-\zeta_0 \cdots \zeta_i \\[2pt]
\wi{M_i}\cdot {}^t(\beta,-\alpha) = 0 \\[1pt]
\zeta_i h_i = 0
\end{array}
\right.
$

\noindent
The result is true for $i=0$~: we take $\wi{M_0} = M_0$, and we get $\zeta_0$
from \eqnz~(\ref{equation:extensionModulo})
because $\gB$ is a \lsdz.
In order to go from rank $i$ to rank $i+1$, we multiply
$h_{i+1} M_{i+1}\cdot {}^t(\beta,-\alpha) \equiv 0 \mod h_i \gB^2$ par~$\zeta_i$.  We get
$$\preskip.2em \postskip.2em 
h_{i+1} \zeta_i M_{i+1}. {}^t(\beta,-\alpha) = 0
 \; \hbox{ et a fortiori }\;
h_{i+1} \zeta_0 \cdots \zeta_i M_{i+1}. {}^t(\beta,-\alpha) = 0  
$$
in $\gB$.

So there exists 
  $\zeta_{i+1} \in \gB$ such that

\hfil
$0 = (1-\zeta_{i+1}) \zeta_0 \cdots \zeta_i M_{i+1}\cdot {}^t(\beta,-\alpha)$
\quad and \quad
$0 = \zeta_{i+1} h_{i+1}$

\noindent
We let $\wi{M_{i+1}} = \wi{M_i} + (1-\zeta_{i+1}) \zeta_0 \cdots \zeta_i
M_{i+1}$. One verifies equalities  $\Tr(\wi{M_{i+1}}) = 1-\zeta_0
\cdots \zeta_{i+1}$ et $\wi{M_{i+1}}\cdot {}^t(\beta,-\alpha) = 0$.
So the \hdr is satisfied for $i+1$.

\noindent Finally, in rank  $n-1$, we have $h_{n-1} = f_n = 1$,
so $\zeta_{n-1} = \zeta_{n-1} h_{n-1} = 0$, which gives  $\Tr(\wi{M_{n-1}}) =
1$,
and the matrix $\wi{M_{n-1}}$ is a \mlp
for the pair $(\alpha,\beta)$.  
\end{proof}

\rem This \demz, like  that of Lemma~\ref{lemleszlop}, is more formidable than it seems. It manages to treat in a single way the case where $\alpha=0$, the case where $\alpha$ is \ndzz, and \gui{all the intermediary cases.}
\eoe

\medskip 
We also have the following easy result.

\begin{theorem}\label{thSurAdp}
Let $\gA\subseteq\gB\subseteq\Frac\gA$. 
\begin{enumerate}
\item If $\gA$ is a \lsdz, the same goes for $\gB$.
\item If $\gA$ is \ariz, the same goes for $\gB$.
\item If $\gA$ is a Pr\"ufer \riz, the same goes for $\gB$.
\end{enumerate}
\end{theorem}
\begin{proof}
Item \emph{1} is left to the reader. 
\\
\emph{2.} 
Let $x$, $y \in \gB$. There exists a $d \in \Reg(\gA)$  such that  $x_1 = dx$, and $y_1 = dy$ are in $\gA$. Then $d (x,y) = (x_1,y_1)$, and a \mlp in $\gA$ for $(x_1,y_1)$ is also a \mlp for~$(x, y)$.
\end{proof}

The two previous \thos are linked to two classic results in the \noe framework (cf.~\cite[page~17]{FJ}):

\mni
{\bf Krull-Akizuki \thoz.} 
\emph{If $\gA$ is a \adk and $\gL$ is a finite extension of the quotient field of $\gA$, then the \cli of~$\gA$ in $\gL$ is a \adkz.}

\mni\rdb
{\bf Grell-Noether \thoz.}\label{thGrellNoether} 
\emph{If $\gA$ is a \adkz, then every \ri contained between $\gA$ and its quotient field is a Dedekind \riz.}

\medskip 
Given the \carn of  \adks (in \clamaz) as integral \noe \adpsz, we see that we have established the non-\noees and non-integral versions of these two \thosz.

We will later prove that in the analogous circumstances, the \ddk of 
$\gB$ is always less than or equal to that of $\gA$, which this time is linked to the \carn of  \adks as \icl \noe \ris \ddi$1$.   

\newpage
\section{\Coh Pr\"ufer \risz} 
\label{secAdPcoh}

\vspace{3pt}
\subsec{First \prtsz}

Recall that over a \qiri a \itf is faithful \ssi it contains a \ndz \elt (see Corollary~\ref{corlemQI}). Actually every \itf contains an \elt that has the same annihilator.
In particular, over a \qiri a \pro \itf is \iv \ssi it is faithful.

\smallskip 
After having provided \carns of  \adps (see Proposition and \Dfnz~\ref{prop.itfplat 2} and  \thref{thPruf}), here are some for  \adpcsz; the reader will find others in Exercise~\ref{exocaracPruCoh}.

\begin{theorem}
\label{th.adpcoh} \emph{(\Carns of  \adpcsz)}\\
For any \ri $\gA$, \propeq
\begin{enumerate}
\item \label{i1th.adpcoh} $\gA$ is a \adpcz.
\item  $\gA$ is an \ari \qiriz.
\item  Every \itf  is \proz.
\item  Every \id with two \gtrs is \proz.
\item  $\gA$ is a \qiri and every \id $\gen{a,b}$ with $a\in\Reg\gA$ is \ivz.
\item \label{i6th.adpcoh}  $\gA$ is a \qiri and every faithful \itf is a \mrcz~$1$.
\end{enumerate}
\end{theorem}
\begin{proof} \emph {1} $\Leftrightarrow$ \emph {2.}  Use Fact~\ref{factLsdzCo}.\\
\emph{3} $\Rightarrow$ \emph{4.} Trivial.
\\
\emph{4} $\Rightarrow$ \emph{2.}
\Thref{thAnar} gives the implication for the \lop character of \idsz.
Moreover a \ri is a \qiri \ssi the \idps are projective.
\\
The implications \emph {1} $\Rightarrow$ \emph {3},  \emph {5},  \emph {6} come from the \carn of  \pro \ids as \lop \ids whose annihilator is an \idm and that of the \iv \ids as \lop \ids containing a \ndz \elt (Lemma~\ref{lemIdproj}, items \emph{2} and \emph{6}).
\\
For the converses, recall that a \idp is \pro \ssi its annihilator is generated by an \idm (Lemma~\ref{lemIdpPtf}), and we can look at the solution of Exercise~\ref{exocaracPruCoh}.
We can also examine these converses in the integral case, where they are clear, and use the \elgbmd of  \qirisz.
\end{proof}

In the local case we obtain the following result (trivial in \clamaz, but meaningful from a \cof point of view).

\begin{fact}\label{factValCoh}
A \adv is \coh \ssi it is integral. We call it a \emph{\ddvz}.%
\index{valuation!domain}\index{domain!valuation ---} 
\end{fact}
\begin{proof}
A \adp is \coh \ssi it is a \qiriz. A \alo is connected. A connected \ri is integral \ssi it is a \qiriz. 
\end{proof}

In this case $\gK=\Frac\gA$ is a \cdi and for all $x\in\gK\eti$, $x$ or~$1/x$ is in~$\gA$. Generally, we call a sub\ri satisfying the preceding \prt a \emph{\adv of the \cdi $\gK$}. It is clear that it is a \ddvz.%
\index{valuation!ring of a \cdiz}\index{ring!valuation ---}

The following stability \prts are easy.

\begin{fact}\label{factAdpc}  ~
\begin{enumerate}
\item 
A \zedr \ri is a \adpcz.

\item 
  A \ri obtained by localization  of a \adpc is a \adpcz.
 A reduced quotient \ri of a \adpc by a \itf is a \adpcz.

\item 
A \ri is a \adpc \ssi it has the same \prt after \lon at \mocoz.

\end{enumerate}
\end{fact}

Recall: item \emph{1} is valid for  \qiris and item \emph{2} for  \anarsz. 

\begin{fact}\label{factAdpcdisc}  Let $\gA$ be a \adpcz.
\begin{enumerate}
\item 
$\gA$ is discrete \ssi $\BB(\gA)$ is discrete.
\item 
$\gA$ is \fdi \ssi it is \dveez. 
\end{enumerate}
\end{fact}

\subsec{Kernel, image and cokernel of a matrix}


\begin{theorem}
\label{ThImMat} Let $\gA $ be a \adpcz. 
\begin{enumerate}
\item 
The image of a matrix $F\in\gA^{n\times m}$ is \isoc to a direct sum of $n$ \itfsz.
\item 
Every \smtf of a \mptf is a \mptfz.
\item 
The kernel of a \ali between \mptfs is a direct summand (therefore \ptfz).
\item 
Every \mpf is a direct sum of its torsion submodule (which is \pfz) and of a \smptfz.
\item 
Every \pro module of rank $k\geq0$ is  \isoc to a direct sum of $k$ \iv \idsz.
\item 
Every \pro module of rank  $\leq k$ is \isoc to a direct sum of $k$ \itfsz.
\end{enumerate}
\end{theorem}
Note: we do not require that $\gA$  be discrete.
\begin{proof}
Consider an arbitrary \ali $\varphi:\Ae m\rightarrow \Ae n$.

 \emph{1.} 
 We treat the case of the module $M=\Im\,\varphi \subseteq\Ae n$. Let $\pi_n:\Ae n\rightarrow \gA$ be the last \coo form. The \id $\pi_n(M)=\fa_n$ is \tf therefore \proz, and the surjective induced map $\pi'_n:M\rightarrow \fa_n$  is split, and 
$$\preskip-.4em \postskip.3em
M\simeq \Ker \pi'_n\oplus \Im \pi'_n=(M\cap \Ae {n-1})\oplus \fa_n.
$$
We end the proof by \recu on $n$:
$M\cap \Ae {n-1}$ is \tf since it is \isoc to a quotient of $M$.
We therefore obtain $M\simeq \fa_1\oplus \cdots \oplus \fa_n$.

 \emph{2.} Results \imdt from \emph{1.} 

 \emph{3.} This results from the image of the \ali being a \mptfz.

 \emph{4.} We treat the case of the module $N=\Coker \varphi$.

\emph{Let us first consider the case where $\gA$ is local}, \cad it is a \ddvz. \\
The matrix of $\varphi$ is put in Smith form (Proposition~\ref{propPfVal}). Since the \ri is integral, $N$ is a direct sum of a free module (corresponding to the zero diagonal \elts in the reduced Smith form) and of a torsion submodule, itself a direct sum of submodules $\aqo{\gA}{d_i}$ corresponding to the \ndzs diagonal \eltsz.

\emph{Next let us consider the case where $\gA$ is integral}.\\
By means of a finite number of \lons at \ecoz, say~$s_1$, $\ldots$, $s_r$, we are brought back to the situation of the local case (Smith reduction of the matrix).
Since $\Ann_\gA(s_i)= \gen{0}$ or $\gen{1}$, and since the \lons at $0$ are useless, we can suppose that the $s_i$'s are in~$\Reg(\gA)$.\\ 
Denote by $T$ the torsion submodule of $N$ and take a look at what happens after \lon at $S_i=s_i^{\NN}$. We easily observe that the torsion submodule of $N_{S_i}$ is equal to $T_{S_i}$. 
Thus, $T$ is \pf because it is \pf after \lon at the $S_i$'s.
It is a direct summand in $N$ because $T_{S_i}$ is a direct summand in $N_{S_i}$ for each $i$: the canonical injection $T\to N$ admits a left-inverse by the \plgrf{plcc.scinde}. 
Finally, the module $N/T$, which is \ptf after \lon at the $S_i$'s, is indeed \ptfz.
\\  
We therefore obtain what we wanted, with a little bonus: the module~$T$ becomes, after \lon at each of the \elts $s_j$ of a comaximal \sys $(s_1,\ldots, s_r)$, a direct sum of cyclic torsion modules, 
\cad \isoc to $\aqo{\gA[1/s_j]}{u_{k,j}}$, with $u_{k,j}\in\Reg(\gA)$.

\emph{Finally, let us consider the \gnl case, where $\gA$ is a \qiriz}.
\\
Starting from the \dem of the integral case, the \elgbmd of  \qiris produces a \sfio $(e_1,\ldots,e_r)$ such that the result is attained in each of the components $e_iN$ (regarded as $\gA[1/e_i]$-module). This \imdt gives the global result.

 \emph{5.} 
In the case where $\gA$ is integral, this results from item~\emph{1} since each \id in the \dcn into a direct sum is of rank $0$ or $1$.
\\
We can deduce the \gnl case by the \elgbmd of \qirisz.
Here is another \demz,%
\footnote{More scholarly or less scholarly, it is difficult to say. This is a matter of taste.} independent of the \dem of item~\emph{1.}
If $M$ is of constant rank $k \geq 1$, then its dual $M\sta$ is also of constant rank, their annihilators are null, and there exists a $\mu \in M\sta$ such that $\Ann(\mu) = \gen{0}$
(see Lemma~\ref{lemQI}). Then $\mu(M)$ is an \iv \id of $\gA$ because its annihilator is \egmt null. Moreover, $M \simeq \Ker \mu \oplus \Im \mu$, which proves that $\Ker\mu$ is \ptf of constant rank $k-1$. We finish by \recuz.

 \emph{6.} Consider $M$ as a direct sum of its components of constant rank,  and apply item \emph{4} to each of them.
\end{proof}

\subsec{Extensions of \adpcsz}%
\index{algebraic!primitively --- element over a \riz}
\label{subsecExtAdpC}

An \elt $x$ of an \Alg $\gB$ is said to be \emph{primitively \agq over $\gA$} if it annihilates a primitive \pol of $\AX$. 
        After a change of base \riz,
a primitively \agq \elt remains primitively
\agqz. The \prt for an \elt to be primitively \agq is local in the following sense.

\begin{plcc} 
\label{plcc.agq} 
\emph{(Primitively \agq \eltsz)} 
Let $S_1$, $\ldots$, $S_n$  be \moco of a \ri $\gA$, $\gB$ be an \Alg and $x\in \gB$.   
Then $x$ is primitively \agq over $\gA$  \ssi it is primitively \agq over each of the $\gA_{S_i}$'s.
\end{plcc}
%
\begin{proof}
We need to show that the condition is sufficient. We have \eco $s_1$, $\ldots$, $s_n$
($s_i\in S_i$) and \pols $f_i\in\AX$ such that $s_i\in\rc(f_i)$ and~$f_i(x)=0$.
If $d_i\geq\deg_X(f_i)+1$, we consider the \pol 

\snic{f=f_1+X^{d_1}f_2+X^{d_1+d_2}f_3+\cdots.}

We then have $f(x)=0$ and $\rc(f)=\sum_{i=1}^n\rc(f_i)=\gen{1}$.
\end{proof}
%

\begin{lemma}\label{lemEmmanuel} \emph{(Emmanuel's trick)} 
Let $\gB$ be a \ri and $\gA$ be a sub\riz. Let $\gA'$ be the \cli of $\gA$ in $\gB$ and $s$ be an \elt of $\gB$ which annihilates a \pol $f(X)=\sum_{k=0}^na_kX^k\in\AX$.\\
Let $g(X)=\sum_{k=1}^nb_kX^{k-1}$ be the \pol $f(X)/(X-s)$. 
\begin{enumerate}
\item The \elts $b_i$ and $b_is$ are in $\gA'$.
\item In $\gA'$ we obtain
\vspace{-1mm}

\snic{
\gen{a_0,\ldots,a_n}=\rc(f)\subseteq\rc(g)+\rc(sg)=\gen{b_1,\ldots,b_n,b_1s,\ldots,b_ns}.}
\item In $\gA'[s]$ the two \ids are equal. 
\end{enumerate}
\end{lemma}
\begin{proof}
Since $f(X)=(X-s)g(X)$, \KROz's \tho tells us that the~$b_i$'s and $b_is$'s are integral over $\gA$. We have 

\snic{b_n=a_n$, $b_{n-1}=b_ns+a_{n-1}$, \ldots,
$b_{1}=b_2s+a_{1}$,
$0=b_1s+a_{0}.}

Therefore each $a_i\in\rc(g)+\rc(sg)$ and, in $\gA'[s]$, step by step, we obtain $b_n\in\rc(f)$, $b_{n-1}\in\rc(f)$, \ldots, $b_{1}\in\rc(f)$.
\end{proof}
%

\begin{theorem}\label{th.2adpcoh}\label{i7th.adpcoh} 
\emph{(Another \carn of  \adpcsz, see also Exercises~\ref{exocaracPruferC} and \ref{exocaracPruCoh})}
A \ri $\gA$ is a \adpc \ssi it is a \qiriz, \icl in $\Frac\gA$, and if every \elt of $\Frac\gA$ is primitively \agq over~$\gA$. 
\end{theorem}
\begin{proof} 
Suppose that $\gA$ is a \adpcz. It remains to show that every \elt of $\Frac\gA$ is primitively \agq over $\gA$.
Let $x=a/b\in\Frac\gA$. 
There is a  principal \lon matrix for $(b,a)$: \smashbot{$\cmatrix {s & u\cr v & t\cr} \in \MM_2(\gA)$},  with $s + t = 1$, $sa = ub$ and $va = tb$.
\\
This gives $sx-u=0$ and $t=vx$. Thus, $x$ annihilates the primitive \polz~$-u+sX+X^2(t-vX)$, or if we prefer $t-(u+v)X+sX^2$. 

Let us prove the converse. It suffices to consider only the integral case. 
We need to show that every \id $\gen{a,b}$ is \lopz.
Suppose \spdg that $a$, $b\in\Reg(\gA)$. The \elt $s=a/b$ annihilates a primitive \pol $f(X)$. Since $\rc(f)=\gen{1}$ in $\gA$, by Lemma~\ref{lemEmmanuel} (items~\emph{1} and~\emph{2}), we have \com \elts $b_1$, $\ldots$, $b_n$, $b_1s$, $\ldots$, $b_ns$ in $\gA$. \\
We then have $s\in\gA[1/b_i]$ and $1/s\in\gA[1/(b_is)]$: in each of the \come \lonsz, $a$ divides $b$ or $b$ divides~$a$.
\end{proof}

The \tho that follows contains a new \dem of the stability of the integral \adps by integral  and \icl extension (see \thref{thExtEntPruf}).
It seems disconcertingly easy when compared to that given without the \cohc hypothesis.

\begin{theorem}\label{factAdpIntExt}
If $\gB$ is a normal \qiriz, and an integral extension of a \adpc $\gA$, then $\gB$ is a \adpcz.
\end{theorem}
%
\begin{proof} Let us first consider the case where \emph{$\gB$ is integral and nontrivial.}
Let $s\in\Frac\gB$. It suffices to show that $s$ is primitively \agq over $\gB$.
We have a nonzero \pol $f(X)\in\AX$ such that $f(s)=0$.
\\
\emph{Case where $\gA$ is a Bézout domain.} We divide $f$ by $\rc(f)$ and we obtain a primitive \pol which annihilates $s$.
\\
\emph{Case of a \ddpz.} After \lon at \ecoz, the \id $\rc(f)$ is generated by one of the \coes of $f$, the first case applies.

 In the \gnl case, the \elgbmd of  \qiris brings us back to the integral case. 
\end{proof}

Now here is the analogue of Proposition~\ref{propAECDN}, which described the \ri of integers of a  number field.  In the case where $\gA$ is a Bézout domain, we could have repeated almost word for word the same \demsz.
Also note that \thref{Thextent} studies a similar situation with a slightly weaker hypothesis.
See also item \emph{1} of \Pbmz~\ref{exoLemmeFourchette}.

\begin{theorem}\label{thAESTE}
 \emph{(\Ri of integers in an \agq extension)}\\
Let $\gA$ be a \adpcz, $\gK=\Frac(\gA)$, $\gL\supseteq\gK$ be a reduced \ri integral over $\gK$ and $\gB$ be the \cli of $\gA$ in $\gL$.
\begin{enumerate}
\item \label{i1thAESTE} $\Frac\gB=\gL=(\Reg\gA)^{-1}\gB$ and  $\gB$ is a \adpcz.
\item  \label{i2thAESTE} If $\gL$ is \stfe over $\gK$ and if  $\gA$ is \fdiz, $\gB$ is \fdiz.
\end{enumerate}
If in addition $\gL$ is \'etale over $\gK$, we obtain
\begin{enumerate} \setcounter{enumi}{2}
\item  \label{i3thAESTE} If  $\gA$ is \noez, the same goes for $\gB$.
\item  \label{i4thAESTE} If  $\gA$ is a \adk (\Dfn \ref{defDDK}), so is $\gB$.
\item  \label{i5thAESTE} If $\gL=\Kx=\aqo\KX f$ with $f\in\AX$ \mon and $\disc_X(f) \in\Reg\gA$,  %
then $\,\fraC1\Delta\gA[x]\subseteq \gB\subseteq\gA[x]$ ($\Delta=\disc_X(f)$). \\
In particular $\gA[x][\fraC1\Delta]=\gB[\fraC1\Delta].$
\item  \label{i6thAESTE} If in addition $\disc_X(f)\in\Ati$, we have $\gB=\gA[x]$ \ste over $\gA$.
\end{enumerate}
 
\end{theorem}
%
\begin{proof}
\emph{\ref{i1thAESTE}.} Direct consequence of Fact~\ref{factReduitEntierQi} and of \thref{factAdpIntExt}.

 \emph{\ref{i2thAESTE}.}
Since $\gB$ is a \adpz, it suffices to know how to test the \dve in $\gB$, \cad testing that an \elt of $\gL$ is a member of $\gB$.  Let~\hbox{$y \in \gL$} and~$Q \in \KY$ be its (\monz) \polmin  over $\gK$. Then $y$ is integral over~$\gA$ \ssi $Q \in \AY$: in the non-\imd sense, %
let $P \in \AY$ be \mon such that $P(y) = 0$, then $Q$ divides $P$ in $\KY$ and Lemma~\ref{lem0IntClos} implies that $Q\in\AY$. 
\\
Note: we might as well have used the \polcarz, but the \dem that uses the \polmin works in a more \gnl framework 
(it suffices for $\gL$ to be \agq over $\gK$ and for us to know how to compute the \polminsz).

\emph{\ref{i5thAESTE}.} 
In the case where $\gA$  is a  Bézout domain and $\gL$ is a field, we apply \thref{Thextent}. 
The result in the \gnl case is then obtained from this \dem by using the \lgbes machineries of \qiris and of \anarsz.\imlgdz\imla

\emph{\ref{i3thAESTE}.} We carry out the \dem under the hypotheses of item \emph{5}. This is not restrictive because by \thref{thEtalePrimitif}, $\gL$ is a product of monogenic \'etale \Klgsz.
\\
Let $\fb_1\subseteq \fb_2 \subseteq\cdots\subseteq \fb_n\subseteq \dots$ be a sequence of \itfs of $\gB$ that we write as $\fb_n=\gen{G_{n}}_\gB$ with $G_n \subseteq G_{n+1}$; we define

\snic {
L_n = \disc_X(f) \cdot \left(\sum_{g \in G_n} \gA g\right) \subseteq \gA[x]
.}

Then $L_1\subseteq L_2\subseteq \cdots\subseteq L_n\subseteq \dots$ is a sequence of \tf \Asubs of $\gA[x]$. 
However, $\gA[x]$ is a free \Amo of finite rank (equal to $\deg(f)$), so \noez. 
We finish by noting that if $L_m=L_{m+1}$, then $\fb_m=\fb_{m+1}$.

\emph{\ref{i4thAESTE}.} Results from \emph{\ref{i2thAESTE}} and \emph{\ref{i3thAESTE}.}

\emph{\ref{i6thAESTE}.}  It is clear that $\gB=\gA[x]$.
\end{proof}

\rem
 The previous \tho applies in two important cases in the history of commutative \algz. 
 \\
 The first case is that of the \ris of integers of number fields, %
with $\gA=\ZZ$ and $\gB$ being the \ri of integers of a number field (a case already examined in Section~\ref{secApTDN}).
 \\
The second case is that of \agq curves. Consider a \cdi $\gk$, the PID $\gA=\gk[x]$ and a \pol $f(x,Y)\in\gk[x,Y]$ \mon in $Y$, \irdz, with $\disc_Y(f)\neq0$. Let $\gK=\gk(x)$.
\\
The \ri $\gA[y]=\gk[x,y]=\aqo{\gk[x,Y]}f$ is integral. 
The planar curve $\cC$ of equation $f(x,Y)=0$ can have singular points, in which case $\gA[y]$ is not \ariz.
But the \cli $\gB$ of $\gA$ in $\gK[y]=\aqo{\gK[Y]}f$ is indeed a \ddp (\thref{thcohdim1}). The field $\gK[y]$ is called the field of functions of $\cC$.  The \ri $\gB$ corresponds to a curve (which is no longer \ncrt plane) without a singular point, with the same field of functions as $\cC$.   
\eoe

\section{\qiris \ddi1}
\label{subsecQiDim1}

Most \gui{classical} \thos regarding \dDks are already valid for \adpcs \ddi $1$, or even for \anarsz. 
We prove a certain number of them in this and the following section.

In this section the results relate to the \qiris \ddi $1$.

The following \tho is a special case of Bass' \gui{stable range} for which we will give \gnl versions (\Thosz~\ref{Bass0} and~\ref{Bass}).

\begin{theorem}
\label{thK1-SLE}
Let $n\geq 3$ and $\tra[\,x_1 \;\cdots\;x_n \,]$ be a \vmd over a \qiri $\gA$ \ddi $1$. This vector is the first column of a matrix of $\En(\gA)$. In particular, $\SLn(\gA)$ is generated by $\En(\gA)$ and $\SL_2(\gA)$ for $n\geq 3$.
For $n\geq 2$ every \vmd is the first column of a matrix of~$\SLn(\gA)$.
\end{theorem}
\begin{proof}
The annihilator of $\gen{\xn}$ is null, so we can by \mlrs transform the vector $v= {\tra[\,x_1 \;\cdots\;x_n \,]}$ into a \vmd ${\tra[\,y_1\;x_2 \;\cdots\;x_n \,]}$, with~$y_1\in\Reg(\gA)$ (cf. Lemma~\ref{lemQI}).
\\
Consider the \riz~$\gB=\aqo{\gA}{y_1}$. This \ri is \zed and the vector $v$ 
becomes equal to ${\tra[\,0\;x_2 \;\cdots\;x_n \,]}$ still \umdz.
{\mathrigid 2mu Since~\hbox{$n\geq 3$}}, we can transform ${\tra[\,x_2 \;\cdots\;x_n \,]}$ into ${\tra[\,1\;0 \;\cdots\;0 \,]}$ by \mlrs in $\gB$ (Exercise \ref{exoSLnEn}).
This gives in $\gA$: ${\tra[\,y_1\;1+ay_1 \;z_3\;\cdots\;z_n \,]}$,
hence, still by \mlrsz, 
${\tra[\,y_1\;1 \;z_3\;\cdots\;z_n \,]}$, then ${\tra[\,1\;0 \;\cdots\;0 \,]}$.
\end{proof}

The following \tho \gnss the analogous result already obtained in number theory (Corollary~\ref{corpropZerdimLib}). Item~\emph{1} concerns the \iv \idsz. Item~\emph{2} applies to all the \itfs of a \adpc \ddi$1$. A \gnn is proposed in \thref{th1.5}.

\begin{theorem}\label{th1-5} \emph{(One and a half \thoz)}%
\index{one and a half!\thoz}
\\
Let $\gA$ be a \qiri \ddi1 and $\fa$ be a \lop \id (thus \ptfz).
\begin{enumerate}
\item 
If  $a\in\fa\cap\Reg(\gA)$, there exists a $b\in\fa$ such that $\fa=\gen{a^n,b}$ for every $n\geq1$. 
\item 
There exists an $a\in\fa$ such that $\Ann(a)=\Ann(\fa)$. For such an $a$ there exists a $b\in\fa$ such that $\fa=\gen{a^n,b}$ for every $n\geq1$. 
\end{enumerate}
\end{theorem}
\begin{proof}
The \dem of item~\emph{1} is identical to that of Corollary~\ref{corpropZerdimLib} which gave the result in number theory.\\
\emph{2.} Every \itf $\fa$ contains an \elt $a$ such that $\Ann(a)=\Ann(\fa)$ (Corollary~\ref{corlemQI}). We pass to the quotient $\aqo{\gA}{1-e}$ where $e$ is the \idm such that $\Ann(a)=\Ann(e)$ and we apply item~\emph{1.} 
\end{proof}
%

\begin{proposition} \label{avant.dekinbe}
Let $\gA$ be a \qiri \ddi1, 
whose Jacobson radical contains a \ndz \eltz, and $\fa$ be an \iv \idz. Then $\fa$ is principal.
\end{proposition}

\begin{proof}
Let $y \in \Rad(\gA)$ and $x \in \fa$ both be \ndzz.
Then $\fa\cap \Rad(\gA)$ contains $a=xy$ which is \ndzz.
By the one and a half \thoz, there exists a $z\in\fa$ such that $\fa=\gen{a^2,z}$.
Therefore $a=ua^2+vz$ which gives $a(1-ua)=vz$ and since $a\in\Rad(\gA)$, $a\in\gen{z}$ so $\fa=\gen{z}$.
\end{proof}

We now revisit the following classical result, in which we will get rid of the \noee hypothesis:
{\em if $\gA$ is an integral \noe \ri \ddi $1$ and $\fa$, $\fb$ are two \ids with $\fa$ \iv and $\fb\neq 0$, then there exists a $u\in\Frac(\gA)$ such that $u\,\fa\subseteq\gA$ %
and $u\fa+\fb=\gen{1}$.}

\begin{lemma}
\label{lemEvitCo}
Let $\gA$ be a \qiri (for example a \adpcz)  \ddi $1$. Let $\fa$ be an \iv \id of $\gA$ and $\fb$ be an \id containing a \ndz \eltz. 
Then there exists an \iv \elt $u$ in $\Frac(\gA)$ such that $u\fa\subseteq\gA$ and $u\fa+\fb=\gen{1}$.
\end{lemma}

\begin{proof}
We carry out the proof in the integral case, leaving it to the readers to apply the \elgbmd of pp-rings.
To facilitate this task, we do not assume $\gA$ to be nontrivial and we put \gui{\ndzz} when in the nontrivial case we would have put \gui{nonzero.}\\
We look for $a$ and $b$ \ndz such that $\fraC{b}{a}\,\fa\subseteq\gA$, \cad $b\,\fa\subseteq a\gA$, and $\gA=\fraC{b}{a}\,\fa+\fb$. 
If $c$ is a \ndz \elt of $\fb$, as the condition should also be realized when $\fb$ is the \id $c\gA$, we must find $a$ and $b$ \ndz such that $b\,\fa\subseteq a\gA$ and $\gA=\fraC{b}{a}\,\fa+c\gA$. 
If steps are taken so that~$a\in\fa$, we will have~$b\in \fraC{b}{a}\fa$, 
and it therefore suffices to realize the conditions $b\,\fa\subseteq a\gA$ 
and~$\gA=\gen{b,c}$. This is what we will do. \\
Let $c \in \fa\cap\fb$ be a \ndz \elt (for example the product of two \ndz \eltsz, one in $\fa$ and the other in $\fb$).
By the one and a half \thoz, there exists an~$a\in\fa$ such that $\fa=\gen{a,c^2}=\gen{a,c}$. 
If $a=0$ the \id $\fa=\gen{c}$ is \idm therefore equal to $\gen{1}$ and there was therefore no need to overexert ourselves:\footnote{Note however that we are not supposed to know in advance if an \iv \id of $\gA$ contains $1$, therefore we have not tired ourselves entirely for nothing, the computation has told us that $1\in\fa.$} we could have chosen $b=a=1$. 
\\
We therefore suppose that $a$ is \ndzz. 
Since $c\in\fa$, we have an \egt $c=\alpha a+\beta c^2$ which gives $c(1-\beta c)=\alpha a$. Let $b=1-\beta c$ such that $\gA=\gen{b,c}$. We obtain $b\,\fa=b\gen{a,c}=\gen{ba,bc} =a\gen{b,\alpha}\subseteq a\gA$.
If $b$ is \ndz we therefore have won, and if $b=0$, then $1\in\gen{c}$ and there was no need to tire ourselves. 
\end{proof}

\begin {proposition}
\label{prop-a/ab} 
Let $\fa$ be an \iv \id of an integral \ri $\gA$ \ddi1. For every nonzero \id $\fb$
of $\gA$, we have an \iso \hbox{of \Amosz} 
$\fa/\fa\fb \simeq \gA/\fb$.
\end{proposition}
\begin{proof}
By Lemma~\ref{lemEvitCo}, there exists an integral \id~$\fa'$ in the class\footnote{See \paref{pageclassgroup}.} of~$\gA\div\fa$ such \hbox{that $\fa' + \fb = \gA$}; \hbox{we have $\fa\fa' = x\gA$} with $x \in \Reg\gA$. The multiplication by~$x$, $\mu_x : \gA \to \gA$, induces an \iso 

\snic{\gA/\fb \simarrow x\gA/x\fb = \fa'\fa/\fa'\fa\fb.}

Let us now consider the canonical map 

\snic{f : \fa'\fa \to \fa/\fa\fb}

\snii
which associates to $y \in \fa'\fa\subseteq \fa$ the class of $y$ modulo $\fa\fb$.
Let us show that $f$ is surjective:
\hbox{indeed, $\fa' + \fb = \gA \Rightarrow \fa'\fa + \fa\fb = \fa$},
 so every \elt of $\fa$ is congruent to an \elt of $\fa'\fa$ modulo $\fa\fb$. Let us finally examine $\Ker f=\fa'\fa \cap \fa\fb$. Since $\fa$ is \ivz, $\fa'\fa \cap \fa\fb = \fa(\fa' \cap
\fb)$, and finally $\fa' + \fb = \gA$ entails \hbox{that $\fa' \cap \fb = \fa'\fb$},
so $\Ker f=\fa'\fa\fb$.  We thus have \isos of
\Amos
$$\preskip-.0em \postskip-.2em 
\;\;\;\gA/\fb \simeq x\gA/x\fb= \fa'\fa/\fa'\fa\fb \simeq \fa/\fa\fb, 
$$
hence the result.
\end{proof}

\begin{corollary} \label{corprop-a/ab}
Let $\gA$ be an integral \ri with $\Kdim\gA \le 1$, $\fa$ be an \iv \id and $\fb$ be a nonzero \idz. We then have an exact sequence of \Amos
$$\preskip-.4em \postskip.4em
0 \to \gA/\fb \to \gA/\fa\fb \to \gA/\fa \to 0.
$$
\end{corollary}

\begin{lemma} \label{lemRadJDIM1}
\emph{(Jacobson radical of a domain \ddi$1$)}
\\ 
Let $\gA$ be an integral \ri \ddi$1$.  
\begin{enumerate}
\item For every nonzero $a$ in $\gA$ we have  $\Rad(\gA)\subseteq \sqrt[\gA]{a\gA}$.
\item For \tf $\fb$ containing $\Rad(\gA)$, we have 
$$\preskip.3em \postskip.3em 
\Rad(\gA)=\fb \big(\Rad(\gA):\fb\big). 
$$
\item If $\Rad(\gA)$ is an \iv \idz, $\gA$ is a Bézout domain.
\end{enumerate} 
\end{lemma}
%
\begin{proof} Let $\fa=\Rad(\gA)$.\\
\emph{1.}
Let $x\in\fa$. The \ri $\aqo\gA a$ is \zedz, so there exist $y$, $z\in\gA$ \hbox{and $m\in\NN$} such that $x^{m}(1+xz)=ay$. 
As $x\in\Rad(\gA)$, we have $1+xz\in\Ati$, therefore $x^{m}\in a\gA$ and $x\in \sqrt[\gA]{a\gA}$.

\emph{2.} If $\fa=0$ it is clear, otherwise the \ri $\gA/\fa$ is \zedrz, so the \itfz~$\fb$ is equal to an \id $\gen{e}$ modulo $\fa$, with $e$ \idm modulo $\fa$. Therefore $\fb=\fb+\fa=\fa+\gen{e}$, then $(\fa:\fb)=\fa+\gen{1-e}$, and finally
$$\preskip-.4em \postskip.4em 
\fb(\fa:\fb)=(\fa+\gen{e})(\fa+\gen{1-e})=\fa. 
$$

\emph{3.} Let $\fc_1$ be a nonzero \itfz.
We define $\fb_1=\fc_1+\fa$  and $\fc_2=(\fc_1:\fb_1)$.  By item~\emph{2} since $\fa$ is \ivz, so is $\fb_1$. 
If $\fb_1\fb'=\gen{b}$ ($b$ is \ndzz), all the \elts of $\fc_1\fb'$ are divisible by $b$. We then consider $\fd=\fraC 1 b \fc_1\fb'$, therefore~\hbox{$\fd\fb_1=\fc_1$} and~$\fd$ is \tfz. Clearly we have $\fd\subseteq \fc_2$. Conversely, if $x\fb_1\subseteq \fc_1$ then $bx=x\fb_1\fb'\subseteq b\fd$, \hbox{so $x\in\fd$}. In short $\fc_2=\fd$ and we have established the \egt $\fb_1\fc_2=\fc_1$, with $\fc_2$ \tfz.
By iterating the procedure we obtain an ascending sequence of \itfsz~$(\fc_k)_{k\in\NN}$ with $\fc_{k+1}=(\fc_k:\fb_k)$ and $\fb_k=\fc_k+\fa$.
\\
Actually $\fc_2=\big(\fc_1:(\fc_1+\fa)\big)=(\fc_1:\fa)$, then $\fc_3=(\fc_2:\fa)=(\fc_1:\fa^{2})$ and more \gnltz~\hbox{$\fc_{k+1}=(\fc_1:\fa^{k})$}. \\
Let $a\neq 0$ in $\fc_1$. By item \emph{1}, $\fa\subseteq \sqrt{a\gA}$. However, $\fa$ is \tfz, therefore the inclusion~\hbox{$\fa\subseteq \sqrt{a\gA}$} implies that for a certain $k$, $\fa^{k}\subseteq a\gA\subseteq \fc_1$, \hbox{therefore $\fc_{k+1}=\gen{1}$}.
\\
When $\fc_{k+1}=\gen{1}$, we have $\fc_1=\prod_{i=1}^{k}\fb_i$, which is \iv as a product of \iv \idsz.
\\
We have shown that every nonzero \itf is \ivz, so the \ri is a \ddpz, and by Proposition~\ref{avant.dekinbe} it is a Bézout \riz.
\end{proof}

\section[\Coh Pr\"ufer \ris \ddi1]{\Coh Pr\"ufer \ris of dimension \texorpdfstring{$\leq1$}{<=1}}
\label{secAdPcohDim1}

\vspace{3pt}
\subsec{When a \adp is a Bézout \riz}
We now \gns a classical result often formulated as follows:\footnote{See the \cov \dfn of a Dedekind \ri on page \pageref{defDDK}.}
{\em an integral Dedekind \ri having a finite number of \idemas is a PID.}

\begin{theorem}\label{dekinbe}
Let $\gA$ be a \adpc \ddi $1$ and whose Jacobson radical contains a \ndz \eltz. Then $\gA$ is a Bézout \riz.
\end{theorem}

\begin{proof}
Let $\fb$ be a \itfz. There exists a $b\in\fb$ such that $\Ann\,\fb=\Ann\, b=\gen{e}$ with $e$ \idmz. Then $\fa=\fb+\gen{e}$ contains the \ndz \elt $b+e$: it is \iv and $\fb=(1-e)\fa$.
It suffices to show that $\fa$ is  principal. 
This results from Proposition~\ref{avant.dekinbe}.
\end{proof}

The previous \tho and the following are to be compared with \thref{propGCDDim1} which affirms that a GCD-domain \ddi$1$ is a Bézout \riz.

\subsec{An important \carnz}

 
The result given in \thref{thcohdim1} below is important: 
the three computational machineries of normality, coherence and dimension at most $1$ combine to provide the machinery of principal localization of finitely generated ideals.

\begin{theorem}\label{thcohdim1}
A normal, \coh \ri $\gA$ \ddi$1$ is a \adpz. 
\end{theorem}
\begin{proof}
Let us start by noticing that $(\gA\div\fa\fb) = (\gA\div\fa)\div\fb$.
\\
Since $\gA$ is a \qiriz, it suffices to treat the integral case and to finish with the \elgbmd of \qirisz.
We therefore suppose that $\gA$ is a domain and we show that every \itf  $\fa$ containing a \ndz \elt is \ivz.\\
Let us consider $(\gA\div \fa) \in \Ifr\gA$ and $\fb=\fa (\gA\div \fa)$, which is a  \tf (integral) \id of $\gA$; we want to show that $\fb=\gA$.
Let us first show that %
 $\gA\div \fb=\gA$.
Let $y \in \gA\div \fb$, hence $y (\gA\div \fa) \subseteq (\gA\div \fa)$.
Since $\gA\div \fa$ is a faithful module (it contains $1$) and is \tfz, $y$ is integral over $\gA$ (see Fact~\ref{factEntiersAnn}) so $y\in\gA$ because $\gA$ is normal.
\\
By \recuz, by using $\gA\div \fb^{k+1} = (\gA\div \fb)\div \fb^k$, we obtain $\gA\div \fb^k = \gA$ for every $k \geq 1$.

Let us fix  a \ndz \elt $x \in \fb$.  By Lemma~\ref{lemDim1-1}, there exists a $k \in\NN\sta$ such that $\fb' := \gen{x}+\fb^k$ is \ivz. Consequently $\fb' (\gA\div\fb') = \gA$. Finally, as $\fb^k\subseteq\fb'\subseteq\fb$, we have $\gA\div \fb'=\gA$, hence $\fb' =\gA$ then $\fb=\gA$.
\end{proof}

\medskip
\exl Other than the example of the \advs given on \paref{exempleadvgroupe}, 
which can have an arbitrary \ddkz, there are other natural examples of \ddps which are not of dimension~$\leq 1$.
\\
The \emph{\ri of  integer-valued \polsz}
is the sub\ri of $\QQ[X]$ formed by the \pols $f(X)$ such that $f(x)\in\ZZ$ for all $x\in\ZZ$. We easily show that it is a free \ZZmo admitting as its basis the combinatorial \pols $ {x} \choose{n}$ for $n\in\NN$. The ideal generated by the $ {x} \choose{n}$'s for $n\geq 1$ is not \tfz. One can show that an integer-valued \pol can be evaluated at an arbitrary $p$-adic integer, which provides an uncountable set of \idepsz. This \ri is a \ddp of dimension two, but the \dem of this result is not simple, especially if we ask that it be \covz.
On this subject see \cite[Ducos]{DucPVE} and \cite[Lombardi]{LomPVE}.
\eoe

\subsec{The structure of \mpfsz}

\begin{theorem}
\label{thPTFDed}
Let $\gA$ be a \adpc \ddiz~$1$. Every \pro module $M$ of constant rank $k \geq 1$ is \isoc to \hbox{$\Ae {k-1} \oplus \fa$}, where $\fa$ is an \iv \idz. In particular, it is generated by $k+1$ \eltsz.
Finally, since $\fa\simeq \Al k M$, the \iso class of $M$  as an \Amo determines that of $\fa$.
\end{theorem}
\begin{proof}
By \thref{ThImMat}, $M$ is a direct sum of $k$ \iv \idsz. It therefore suffices to treat the case $M\simeq \fa\oplus \fb$, with \iv \ids $\fa$ and~$\fb$.
By Lemma~\ref{lemEvitCo}, we can find an \id $\fa_1$ such that $\fa_1\simeq \fa$ (as~$\gA$-modules) and $\fa_1+\fb=\gen{1}$ (as \idsz). We then have the short exact sequence
$$\preskip.8pt
\gen{0} \longrightarrow \fa_1 \fb = \fa_1\cap \fb\vers{\delta} \fa_1\oplus \fb
\vers{\sigma}
\fa_1+\fb=\gA\longrightarrow \gen{0},
\postskip2.5pt
$$
where $\delta(x)=( x, -x)$ and $\sigma(x, y)=x+y$.
Finally, since this sequence is split, we obtain
$M\simeq \fa\oplus \fb \simeq \fa_1\oplus \fb \simeq \gA\oplus (\fa_1 \cap \fb)
=\gA\oplus (\fa_1 \fb) $.
\end{proof}

An \imde consequence is the following structure \thoz.
\begin{corollary}\label{corthPTFDed}
Let $\gA$ be a \adpc \ddi$1$. Every \mptf is \isoc to a direct sum

\snic{r_{1}\gA \oplus r_{2}\Ae 2  \oplus \cdots   \oplus r_{n}\Ae n\oplus \fa,}

where the $r_{i}$'s are \ort \idms (some can be null) and $\fa$ is a \itfz.
\end{corollary}

\begin{proposition}\label{propAriCohZed}
Let $\gA$ be a \zed \anarz. Every matrix admits a reduced Smith form. Consequently every \pf \Amo   is \isoc to a direct sum of cyclic modules~$\aqo\gA{a_k}$. 
\end{proposition}

%
\begin{proof} 
If $\gA$ is local, it is a local Bézout \ri and the matrix admits a reduced Smith form (Proposition~\ref{propPfVal}), which gives the result. 
By following the \dem of the local case, and by applying the \lgbe machinery of  \anars (\paref{MetgenAnar}), we produce a family of \eco $(s_1,\ldots,s_r)$ such that the result is guaranteed over each \ri $\gA[1/s_i]$.\imla
\\
Since $\gA$ is \zedz, every principal filter is generated by an \idm (Lemma~\ref{lemme:idempotentDimension0}~\emph{\iref{LID002}}). Consider the \idms $e_i$ corresponding to the $s_i$'s, then a \sfio $(r_j)$ such that each $e_i$ is a sum of certain $r_j$'s.
\\
The \ri is written as a finite product $\prod\gA_j$ with the Smith reduction over each $\gA_j$.
The result is therefore guaranteed. 
\end{proof}
\rems
For the uniqueness of the \dcnz, see \thref{prop unicyc}. Moreover, the \dem shows that the  reduction can be done with products of \elr matrices. 
Finally, a \gnn is proposed in Exercise~\ref{exoAnarlgb}.
\eoe

\begin{corollary}\label{corpropAriCohZed}
Let $\gA$ be an \anar \ddi $1$. Every \pf torsion \Amo is \isoc to  a direct sum of cyclic modules~$\aqo\gA{b,a_k}$ with $b\in\Reg(\gA)$.
\end{corollary}
%
\begin{proof}
The module is annihilated by an \elt $b\in\Reg(\gA)$. We consider it as an $\aqo\gA b$-module and we apply Proposition~\ref{propAriCohZed}.
\end{proof}

We can now synthesize \thrfs{thPTFDed}{ThImMat}, and Corollary~\ref{corpropAriCohZed} as follows.
We leave it up to the reader to give the statement for the \qiri case (\cad for \adpcsz).

\begin{theorem}\label{thMpfPruCohDim}\emph{(\Tho of the invariant factors)}\\
Over a \ddp $\gA$ \ddi$1$, every \mpf is a direct sum
\begin{itemize}
\item of a \ptf \Amoz, null or of the form $\Ae r\oplus\fa$ ($r\geq0$, $\fa$ an \iv \idz),
\item   and of its torsion submodule, which is \isoc to a direct sum of cyclic modules~$\aqo\gA{b,a_k}$ with $b\in\Reg(\gA)$.
\end{itemize}
In addition
 \begin{itemize}
\item the \id $\fa$ is uniquely determined by the module, 
\item we can assume that the \ids $\gen{b,a_k}$ are totally ordered with respect to the inclusion relation, and the \dcn of the torsion submodule is then unique in the precise sense given in \thref{prop unicyc}.
\end{itemize}
\end{theorem}

\rems
1) In particular, the structure \tho for  \mpfs over a PID (Proposition~\ref{propPfPID}) is valid for every  Bézout domain \ddi$1$.
\\
2) For a torsion module $M$, the \ids $\gen{b,a_k}$ of the previous \tho are the \ix{invariant factors} of $M$, in accordance with the \dfn given in \thref{prop unicyc}. 
\eoe

\subsubsection*{Reduction of matrices}

The following \tho gives a reduced form for a column matrix, \`a la Bézout. It would be interesting to \gnr it to an arbitrary matrix.

\pagebreak	
\begin{theorem}\label{mat33}
Let $\gA$ be a \adpc \ddi$1$ and $x_1$, \dots, $x_n\in \gA$.  
There exists a matrix $M\in\GL_n(\gA)$ such that 

\snic{M\,\tra[\,x_1 \;\cdots\; x_n\,]=\tra[\,y_1 \;y_2 \;0 \;\cdots\; 0\,].}
\end{theorem}
\begin{proof}
It suffices to treat the case where $n=3$.
\\
If $e$ is an \idmz, then $\GL_n(\gA) \simeq \GL_n(\gA_e) \times \GL_n(\gA_{1-e})$:
even if it entails localizing by inverting both the annihilating \idm of $\gen{x_1,x_2,x_3}$ and its complement,  we can therefore assume that {\mathrigid 2mu $\Ann(\gen{x_1,x_2,x_3})=\gen{0}$}.
\\
Let $A$ be a \mlp for $(x_1,x_2,x_3)$.\\
The module $K=\Im(\I_3-A)$ is the kernel of \lin form associated with the row vector $X=[\,x_1 \;x_2\; x_3\,]$ and it is a \pro module of rank $2$ as a direct summand in $\Ae 3$. \Thref{thPTFDed} tells us that $K$ contains a free submodule of rank $1$ as a direct summand in $\Ae 3$, that is a module~$\gA v$ where $v$ is a \vmd of  $\Ae 3$.
By \thref{thK1-SLE}, this vector is the last column of an \iv matrix~$U$; the last \coe of~$XU$ is null and the matrix $M = \tra {U}$ is the matrix we were looking for.
\end{proof}



\section{Factorization of \itfs} 
\label{secAdPfactpar}

\vspace{3pt}

\subsec{\Gnl factorizations} 
\label{subsecfactogene}

In a \gnl \anar it seems that we do not have any \fcn results that go beyond what stems from the fact that the \iv \ids (\cad the \itfs containing a \ndz \eltz) form a GCD-\moz, and more \prmt the non-negative sub\mo of an \grlz.

For example the Riesz \tho can be reread as follows.

\begin{theorem}\label{thRieszAnar} \emph{(Riesz \tho for  \anarsz)}\\
Let $\gA$ be an \anarz, $(\fa_i)_{i\in\lrbn}$ and $(\fb_j)_{j\in\lrbm}$ be \iv \ids such that $\prod_{i=1}^n\fa_i=\prod_{j=1}^m\fb_j$. 
Then there exist \iv \ids $(\fc_{i,j})_{i\in\lrbn,j\in\lrbm}$ such that we have for all $i$ and all $j$,
$$\ndsp \preskip.4em 
\fa_i=\prod_{j=1}^m\fc_{i,j}\,\hbox{ and }\,\fb_j=\prod_{i=1}^n\fc_{i,j}. 
$$
\end{theorem}

\penalty-2500
\subsec{Factorizations in dimension 1} 
\label{subsecfactodimun}

\begin{theorem}
\label{thFactor}
In a \adpc \ddi $1$, we consider two \itfs $\fa$ and $\fb$ with $\fa$ \ivz. 
\\
Then we can write

\snic{\fa=\fa_1\fa_2\;$ with $\;\fa_1+\fb=\gen{1}$ and $\,\fb^n\subseteq \fa_2,}

for some suitable integer $n$. This expression is unique and we have 

\snic{\fa_1+\fa_2=\gen{1},$
   $\;\fa_2=\fa+\fb^n=\fa+\fb^{n+1}.}
\end{theorem}

\begin{proof}
This is a special case of Lemma~\ref{lemDim1-1}.
\end{proof}

\rem We do not need to assume that the \ids are detachable.

\begin{theorem}
\label{thFactor2}
We consider in a \adpc \ddi$1$ some pairwise \com \itfs $\,\fp_1$, $\dots$, $\fp_n$, and an \iv \id $\fa$.\\
We can write $\fa=\fa_0\cdot \fa_1\cdots\, \fa_n$ with the pairwise \com \itfs $\,\fa_0$, $\dots$, $\fa_n$ and, for $j \geq 1$, $\fp_j^{m_j}\subseteq \fa_j$ with $m_j$ as the suitable integer.
\\
This expression is unique and we have $\,\fa_j=\fa+\fp_j^{m_j}=\fa+\fp_j^{1+m_j}$.
\end{theorem}
\begin{proof} 
By \recu by using \thref{thFactor} with $\fb \in \{\fp_1,\dots,\fp_n\}$.
\end{proof}

\subsec{Pr\"ufer \ris admitting \fapsz} 
\label{subsecDedfactpar}
Let us re-express the \dfn of \dcnps  (given for \grlsz) in the framework of the \mo of \iv \ids of a \adpc $\gA$  (this is the non-negative sub\mo of the \grl formed by the \iv \elts of $\Ifr(\gA)$).

\begin{definition}
Let $F=(\fa_1,\dots,\fa_n) $ be a finite family of \iv \ids in a \ri $\gA$. We say that $F$ admits a \emph{\fapz} if there exists a family $P=(\fp_1,\dots,\fp_k)$ of pairwise \com \iv \ids such that every \id $\fa_j$  can be written in the form: $\fa_j = \fp_1^{m_{1j}}\cdots \fp_k^{m_{kj}} $ (certain $m_{ij}$'s can be null).
We then say that $P$ is a \emph{\bdfz} for the family $F$.%
\index{partial factorization!basis}
\end{definition}

For the \mo $\Ifr(\gA)$ to be discrete we need to assume that $\gA$ is \fdiz.
This leads to the following \dfnz.

\begin{definition}
\label{defPrufFP}
A \ri is called a \emph{\fap\adpz} if it is a \fdi \adpcz\footnote{By Proposition~\ref{aritfdi} an \anar is \fdi \ssi the \dve relation is explicit.}
and if every finite family of \iv \ids admits a \fapz.%
\index{ring!partial factorization Pr\"{u}fer ---}%
\index{factorization!partial ---}\index{Pru@Pr\"ufer!partial factorization --- ring}
\end{definition}

\begin{lemma}
\label{lemK1}
A \fap\adp  is \ddiz~$1$.
\end{lemma}

\begin{proof}
We consider a \ndz \elt $y$. We want to show that $\aqo{\gA}{y}$ is \zedz. \\
For this we take a \ndz $x$ and we want to find $a\in\gA$ and $n\in\NN$ such that $x^n(1-ax)\equiv 0 \mod y$. 
The \fap of $(x,y)$ gives 

\snic{\gen{x} = \fp_{1}^{\alpha_1} \cdots \fp_{i}^{\alpha_i}\fq_{1}^{\beta_1} \cdots
\fq_{j}^{\beta_j} = \fa\fb, \hbox{  and }
\gen{y} = \fp_{1}^{\gamma_1} \cdots \fp_{i}^{\gamma_i} \fh_{1}^{\delta_1}
\cdots \fh_{k}^{\delta_k} = \fc\fd}

with all the  exponents $>0$.
There exists an $n\geq0$ such that $\fa^n$ is a multiple of $\fc$ which gives $\gen{x^n} = \fc\ffg$.
As $\gen{x} + \fd = 1$, there exists some~$a \in \gA$ such that~$1 - ax \in \fd$. We therefore have $\gen{y}=\fc\fd\supseteq  \fc\ffg\fd = \gen{x^n}\fd \supseteq \gen{x^n(1 - ax)}$,
\cadz~$x^n(1-ax)\equiv 0 \mod y$.
\end{proof}

\subsec{Dedekind \risz} 
\label{subsecDedfactot}

\begin{definition}\label{defDDK}
We call a \fdi and \noe \coh \adp a \emph{\adkz}. 
A \emph{\dDkz} is an integral \adk (or yet again a connected \adkz).%
\index{ring!Dedekind ---}%
\index{Dedekind!\riz}  
\end{definition}

\begin{theorem}\label{prop2DDK}
A \adk is a \fap\adpz, so \ddi $1$. 
\end{theorem}
\begin{proof}
\Thref{th2GpRtcl} gives the \fap result in the framework of \trdis and we finish with Lemma~\ref{lemK1}.
\end{proof}
%

\begin{theorem}\label{propDDK} \emph{(\Carns of  \adksz)}\\
For some \ri $\gA$ \propeq
\begin{enumerate}
\item $\gA$ is a \adkz.
\item $\gA$ is an \ariz, \noez, \qiri \dveez.
\item $\gA$ is a normal \qiri \ddi $1$,  \dveez, which is \coh and \noez.
\end{enumerate}
\end{theorem}
\begin{proof}
Since $\gA$ is a \coh Pr\"ufer \ri \ssi it is \ari and a \qiriz, and since an \anar is \fdi \ssi it is \dveez, items \emph{1} and \emph{2} are \eqvsz.
\\
The implication \emph{1} $\Rightarrow$  \emph{3} results from \thref{prop2DDK}, and \thref{thcohdim1} gives the converse (it simply suffices to add \fdi and \noe in the hypothesis and the conclusion).
\end{proof}
%

\begin{definition}
Let $\fa $ be an \id of a \ri $\gA$. We say that $\fa$ admits a \emph{\facz} if it is of the form $\fa = \fp_1^{m_{1}}\cdots \fp_k^{m_{k}} $ ($m_{i}>0$, $k>0$) where  \ids $\fp_i$ are detachable, strict and maximal (in other words, each 
\ri $\gA\sur{\fp_i}$ is a nontrivial \cdiz).%
\index{total factorization!of an \id in a \riz}%
\index{factorization!total ---}%
\end{definition}

\begin{thdef}\label{thdefDDKTOT}%
\index{ring!total \fcn Dedekind ---}%
\index{Dedekind!total \fcn --- \riz}%
\index{total factorization!Dedekind \riz}%
\index{factorization!total ---}%
\emph{(Total \fcn Dedekind ring)} 
\\
For a nontrivial \fdi \qiri $\gA$, \propeq
\begin{enumerate}
\item Every \idp $\gen{a}\neq\gen{1}$ with $a\in\Reg\gA$ admits a \facz.
\item The \ri $\gA$ is a \adkz, and every \iv \id $\neq\gen{1}$ admits a \facz.
\end{enumerate}
Such a \ri is called a \emph{\fac\adkz}.
\end{thdef}
%
\begin{proof} We need to show that \emph{1} implies \emph{2.} 
We treat the integral case (the \qiri case is then easily deduced). 
\\
We refer to Exercise~\ref{exoDecompIdeal} and its solution. We see that every \itf containing a \ndz \elt is \ivz, and that it admits a \facz. \Thref{th.adpcoh} then tells us that $\gA$ is a \adpcz. It remains to see that it is \noez. 
Consider a \itf and its \fac  $\fa = \fp_1^{m_{1}}\cdots \fp_k^{m_{k}} $. Every \itf  $\fb\supseteq\fa$ is of the form $\fp_1^{n_{1}}\cdots \fp_k^{n_{k}}$ with all the $n_i\in\lrb{0..m_i}$. 
Every ascending sequence of \itfs starting with $\fa$ therefore admits two consecutive equal terms. 
\end{proof}

\rem Exercise~\ref{exoDecompIdeal} uses no complex theoretical paraphernalia. 
So it is possible to expose the theory of \adks by starting with the previous \thoz, which promptly leads to the essential results. The main drawback of this approach is that it is based on a \fac \prt which is not \gnlt satisfied from a \cov point of view, even by the PIDs, and which does not \gnlt extend to the integral extensions.  
\eoe

\medskip 
Recall that we have already established \thref{thAESTE} regarding the finite extensions of \adksz. 

We can add the following more precise result.

\begin{theorem} \label{lemthAESTE} \emph{(A computation of \cliz)}\\
Let $\gA$ be a \adkz, $\gK=\Frac(\gA)$, $\gL\supseteq\gK$ be an \'etale \Klg and $\gB$ be the \cli of $\gA$ in $\gL$.\\
Suppose that $\gL=\aqo\KX f$ with \mon $f\in\AX$ and $\disc_X(f) \in\Reg\gA$ (which is not really restrictive). If $\gen{\disc_X(f)}$ admits a \facz, and if for each \idema $\fm$ of this \fcnz, the residual field $\gA/\fm$ is perfect, then $\gB$ is a \ptf \Amoz.
\end{theorem}
%
\begin{proof} 
As $\gA$ is a \qiriz, it suffices to treat the case where $\gA$ is integral (\elgbmd of  \qirisz), so~$\gK$ is a \cdiz.
The hypothesis $\gL=\aqo\KX f$ with \mon $f\in\AX$ and $\disc_X(f) \in\Reg\gA$ is not really restrictive because by \thref{thEtalePrimitif}, $\gL$ is a product of monogenic \'etale \Klgsz.
We can even suppose that $\gL$ is an \'etale field over $\gK$  (\elgbm of \zedr \risz).\\
Let $\Delta=\disc_X(f)$.
By item~\emph{\ref{i5thAESTE}} of \thref{thAESTE} we have the inclusions

\snic{\gA[x]\subseteq\gB\subseteq\fraC 1 {\Delta}\,\gA[x].}

Thus $\gB$ is a submodule of the \tf \Amo $\fraC 1 {\Delta}\gA[x]$. By \thref{ThImMat}, if $\gB$ is \tfz, it is \ptfz.   
\\
We have  $\gA[x,\fraC 1 {\Delta}]=\gB[\fraC 1 {\Delta}]$, so $\gB$ is \tf after \lon at~$\Delta^{\NN}$. It remains to show that $\gB$ is \tf after \lon at $S=1+\Delta\gA$. The \ri $\gA_S$ is a Bézout \ri (\thref{dekinbe}). If $\fp_1$, \dots, $\fp_r$ are the \idemas that intervene in the \fac of $\Delta$, the \mos $1+\fp_i$ are \com in $\gA_S$, and it suffices to show that $\gB$ is \tf after \lon at each of the $1+\fp_i$.
We are thus brought back to the case treated in Lemma~\ref{lemlemthAESTE} that follows. 
\end{proof}

Note that a local \dDk $\gV$ is just as much a \fdi \noe \ddvz, 
as a local PID with detachable~$\gV\eti$. 
The following lemma in addition requires that the radical be principal (which is  automatic in \clamaz). 
In this case we will also say that $\gV$ is a \emph{discrete \advz} or a \emph{DVR}, according to the classical terminology;
and any generator of $\Rad\gV$ is called a \emph{regular parameter}.%
\index{ring!discrete valuation ---}%
\index{valuation!discrete ---}\index{DVR}%
\index{regular!parameter}%
\index{parameter!regular ---}%
\index{discrete valuation ring (DVR)}

\begin{lemma} \label{lemlemthAESTE}
Let $\gV$ be a local \dDk with $\Rad\gV=p\gV$ and with perfect residual field $\gk=\aqo\gV p$. Let  $f\in\VX$ be an \ird \monz, therefore $\Delta=\disc_X(f) \in\Reg\gV$.
Let $\gK=\Frac(\gV)$, $\gL=\Kx=\aqo\KX f$, and~$\gW$ be the \cli of~$\gV$ in~$\gL$. Then $\gW$ is \tf over~$\gV$.
\end{lemma}
%
\begin{proof} Since $\gk$ is perfect, by Lemma~\ref{lemSqfDec}, for every \polu $f_i$ of $\VX$ we know how to compute the \gui{squarefree subset} of $\ov {f_i}$ ($f_i$ taken modulo $p$), \cad a \spl \pol $\ov {g_i}$ in~$\gk[X]$ which divides $\ov {f_i}$, and whose power is a multiple of $\ov {f_i}$.
\\
The strategy is to add some \elts $x_i\in\gW$ to $\gV[x]$ until we obtain a \ri $\gW'$ whose radical is an \iv \idz.
When this is realized, we know by Lemma~\ref{lemRadJDIM1} that $\gW'$ is a \ddpz, therefore that it is \iclz, thus equal to $\gW$.\\
To \gui{construct} $\gW'$ (\tf over $\gV$) we will use in an \recu the following fact, initialized with $\gW_1=\gV[x]$ ($x_1=x$, $r_1=1$). 

{\it Fact. Let $\gW_k=\gV[x_1,\dots,x_{r_k}]\subseteq \gW$, then 
$$\preskip.3em \postskip.2em 
\Rad(\gW_k)=\geN{p,g_1(x_1),\dots,g_k(x_{r_k})}, 
$$
where $\ov{g_i}$ is the squarefree subset of $\ov{f_i}$ and $f_i$ is the \polmin over~$\gK$ of the integer~$x_i$.}
\begin{Proof}{}
\Thref{thJacplc} states that $\Rad(\gW_k)=\rD_{\gW_k}(p\gW_k)$.
This \id is the inverse image of $\rD_{\gW_k/p\gW_k}(0)$ and we have $\gW_k/p\gW_k=\gk[\ov{x_1},\dots,\ov{x_{r_k}}]$. As the $g_i(x_i)$'s are nilpotent modulo $p$ by construction, they are in the nilradical $\rD_{\gW_k}(p\gW_k)$. Now it suffices to verify that the \klg 

\snic{\Aqo{\gk[\ov{x_1},\dots,\ov{x_{r_k}}]}{\ov {g_1}(\ov{x_1}), \dots, 
\ov {g_{r_k}}(\ov{x_{r_k}})}}

is reduced. Actually $\gW_k$ is a \tf $\gV$-submodule of $\fraC1\Delta\gV[x]$,
therefore is free finite over $\gV$. Consequently $\gW_k/p\gW_k$ is strictly finite over $\gk$, and it is \'etale because it is generated by some \elts that annihilate \spl \pols over $\gk$ (\thref{corlemEtaleEtage}). 
\end{Proof}

Given this, since $\gW$ is a \ddpz, we know how to invert the \itf $\Rad(\gW_k)$ in $\gW$. \\
This means computing some \elts $x_{r_{k}+1}$, \dots, $x_{r_{k+1}}$ of $\gW$ and a \itf $\ffg_k$ in the new \ri $\gW_{k+1}$ such that the  \id $\ffg_k\Rad(\gW_k)$ is  principal (and nonzero). 
However, it is possible that the \gtrs of $\Rad(\gW_k)$ do not generate the \id $\Rad(\gW_{k+1})$ of $\gW_{k+1}$, which forces an iteration of the process.\\
The ascending sequence of $\gW_k$'s is an ascending sequence of \tf $\gV$-modules contained in $\fraC1\Delta \gV[x]$, therefore it admits two equal consecutive terms. In this case we have reached the required goal.
\hfill\eop
\end{proof}

\begin{plcc} 
\label{plcc.ddk} \emph{(Dedekind \risz)}\\
Let $s_1$, $\ldots$, $s_n$  be \eco of a \ri $\gA$.   
Then
\begin{enumerate}
\item The \ri $\gA$ is \fdi \noe \coh \ssi each of the~$\gA_{s_i}$'s is \fdi \noe \cohz.
\item The \ri $\gA$ is a \adk \ssi each of the~$\gA_{s_i}$'s is a \adkz.
\end{enumerate}
 
\end{plcc}
%
\begin{proof} We already know that the \plgc works for the \adps  and for the \coris with \mocoz. The same goes for the \ris or \noes modules (a \dem is given with the \plgref{plcc.ptf}).
\\
It remains to examine the \gui{\fdiz} \prt in the case of \ecoz. 
Let $\fa$ be a \itf and~$x\in\gA$.
It is clear that if we have a test for~$x\in\fa\gA_{s_i}$ for each of the~$s_i$'s, this provides a test for $x\in\fa\gA$. The difficulty is in the other direction: if $\gA$ is \fdi and if $s\in\gA$, then $\gA[1/s]$ is \fdiz. It is not true in \gnlz, but it is true for the \noe \corisz. Indeed, the membership $x\in\fa\gA[1/s]$ is equivalent to $x\in(\fa:s^{\infty})_\gA$. However, the \id  $(\fa:s^{\infty})_\gA$ is the union of the ascending sequence of \itfs $(\fa:s^{n})_\gA$, and as soon as $(\fa:s^{n})_\gA=(\fa:s^{n+1})_\gA$, the sequence becomes constant.
\end{proof}
%

%
%
%
%
%

\Exercices

\begin{exercise} \label {exoDetTrick1}
       (Another \gui{determinant trick})\\
{\rm 
Let $E$ be a faithful $\gA$-module generated by $n$ \elts and $\fa \subseteq \fb$ be two \ids of $\gA$ satisfying $\fa E = \fb E$.
Show that $\fa\fb^{n-1} = \fb^n$.
}
\end{exercise}

\vspace{-1em}
\pagebreak	

\begin{exercise}\label{exoMatLoc22}
{(Principal \lon matrices in $\MM_2(\gA)$)}   
{\rm Let $x$, $y \in\gA$.

\emph {1.} 
Show that the \id $\gen {x,y}$ is \lop \ssi there exists a matrix~$B \in \MM_2(\gA)$ of
trace $1$ satisfying $[\,x\; y\,]B = 0$; in this case, $A = \wi B$ is a \mlp for $(x, y)$.

\emph {2.}
Let $z \in \gA$; suppose that there exists an \id $\fb$ such that $\gen{x,y}\fb = \gen {z}$.  Show that there exists a $B \in \MM_2(\gA)$ such that $z[x\; y\,]B = 0$ and $z\big(1 - \Tr(B)\big) = 0$.

\emph {3.}
Deduce from the previous questions another \dem of Lemma~\ref{lemleszlop}.

}

\end{exercise}

\vspace{-1em}
\begin{exercise}
\label{exoPrufNagata}\label{exoWhenNagataRingIsArith}
($\gA$ is \ari $\Leftrightarrow$  $\gA(X)$ is a Bézout \ri)  {\rm See also Exercise~\ref{exoBézoutKdim1TransfertArX}.
\\
Let $\gA$ be a \ri and $\gA(X)$ be the Nagata \riz. 

 \emph{1.}
Show that for $a,b\in\gA$, $a\mid b$ in $\gA$ \ssi $a \mid b$ in $\gA(X)$.

 \emph{2.} If $\gA$ is an \anar and $f\in\AX$, we have in $\gA(X)$

\snic{\gen{f}=\rc_\gA(f)\gA(X).}

Also show that $\gA(X)$ is a Bézout \riz.

 \emph {3.}
Let $x$, $y \in \gA$. Show that if $\gen{x,y}$ is \lop in $\gA(X)$, it is \lop in $\gA$ (use Exercise~\ref{exoMatLoc22}).
In particular, if $\gA(X)$ is \ariz, the same goes for $\gA$.
A fortiori, if $\ArX$ is \ariz, the same goes for $\gA$.

\emph {4.}
Conclude the result. 

Note. Concerning the \ri $\gA(X)$ see Fact~\ref{factLocNagata} and Exercise~\ref{exoNagatalocal}. 
}
\end{exercise}

\vspace{-1em}
\begin{exercise} 
\label{exoanars} (A few other \caras \prts of \anarsz)\\
{\rm  
For any \ri $\gA$, \propeq
\begin{description}\itemsep0pt
\item [\phantom{.2}$(1)$] The ring $\gA$ is an \anarz. 
\item [$(2.1)$] For all \ids  $\fa$, $\fb$ and $\fc$ we have 
$\fa\cap(\fb+\fc)=(\fa\cap \fb)+(\fa\cap \fc).$ 
\item [$(2.2)$] As above but limiting ourselves to the \idpsz. 
\item [$(2.3)$] As above but limiting ourselves to the case $\fb=\gen{x}$, $\fc=\gen{y}$ and $\fa=\gen{x+y}.$
\item [$(3.1)$] For all \ids  $\fa$, $\fb$ and $\fc$ we have $\fa+(\fb\cap \fc)=(\fa+\fb)\cap (\fa+\fc).$
\item [$(3.2)$] As above but limiting ourselves to the \idpsz. 
\item [$(3.3)$] As above but limiting ourselves to the case $\fa=\gen{x}$, $\fb=\gen{y}$ and $\fc=\gen{x+y}.$  
\item [$(4.1)$] For all \itfs $\fa$, $\fb$ and $\fc$ we have $(\fb+\fc):\fa=(\fb:\fa)+(\fc:\fa).$
\item [$(4.2)$] As above with \idps $\fb$ and $\fc$, and $\fa=\fb+\fc$.
\item [$(5.1)$] For every \id $\fa$ and all the \itfs $\fb$ and $\fc$ we have the \egt
$$\preskip-.6em \postskip.2em 
 \fa:(\fb\cap \fc)=(\fa:\fb)+(\fa:\fc).
$$
\item [$(5.2)$] As above with \idps $\fb$ and $\fc$, and $\fa=\fb\cap \fc$. 
\end{description}
} 
Hint: to prove that the conditions are \ncrs we use the \gnl method explained on \paref{MetgenAnar}.
\end{exercise}

\vspace{-1em}
\begin{exercise} 
\label{exoNormalClass} 
{\rm Prove in \clama that a \ri is normal \ssi it becomes normal when we localize at an arbitrary \idep (recall that in the integral case, normal means \icl in its quotient field).
}
\end{exercise}

\vspace{-1.2em}
\pagebreak	

\begin{exercise}\label{FermetureAlgZariski}
 {(\Agq closure: a \tho due to Zariski)}\\
{\rm  
Let $\gK \subseteq \gL$ be two \cdisz, $\gK'$ be the \agq closure of $\gK$ in $\gL$. Then the \agq closure of $\gK(\Xn)$ in $\gL(\Xn)$ is $\gK'(\Xn)$; an analogous result holds if we replace \agq closure by \spb \agq closure.
}

\end{exercise}

\vspace{-1em}
\begin{exercise}\label{exoNotAbsIntClos}
{(A lack of integrality by \edsz)}\\
{\rm  
Let $\gk$ be a \cdi of \cara $p \ge 3$, $a \in \gk$ and $f = Y^2 - f(X)
\in \gk[X,Y]$ with $f(X) = X^p - a$.

\emph {1.}
Show that $Y^2 - f(X)$ is \ix{absolutely \irdz}, that is that for every overfield $\gk' \supseteq \gk$, the \pol $Y^2 - f(X)$ is \ird in $\gk'[X,Y]$.

Let $\gk[x,y] = \aqo {\gk[X,Y]}{Y^2 - f(X)}$ and $\gk(x,y) = \Frac(\gk[x,y])$.

\emph {2.}
Show that $\gk$ is \agqt closed in $\gk(x,y)$ and that for every \agq extension $\gk'$ of $\gk$, we have $\gk' \te_\gk \gk(x,y) = \gk'(x,y)$.

\emph {3.}
Suppose that $a \notin \gk^p$. Show that $\gk[x,y]$ is \icl and that~$\gk(x,y)$ is not a field of rational fractions with one \idtr over $\gk$.

\emph {4.}
Suppose $a \in \gk^p$ (for example $a = 0$). Show that $\gk[x,y]$ is not \icl and explicate $t \in \gk(x,y)$ such that $\gk(x,y) = \gk(t)$.

}

\end{exercise}

\vspace{-1em}
\begin{exercise}
\label{exoAnneauOuvertP1} (The \ri of functions over the projective line minus a finite number  of points)\\
{\rm 
We informally use in this exercise the notions of an affine scheme and of a projective line which have already been discussed in Sections~\ref{sec1Apf} and~\ref{secGrassman} (see  pages \pageref{subsecNstMorphismes} to \pageref{subsubsecTanFonct}).\\
If $\gk$ is a \cdiz, the \klg  {of \pol functions defined over the affine line $\AA^1(\gk)$} is $\gk[t]$. 
If we think of $\AA^1(\gk) \cup \{\infty\}  = \PP^1(\gk)$, the \elts of~$\gk[t]$ are then the rational fractions over $\PP^1(\gk)$ which are \gui{defined everywhere, except maybe at $\infty$.}
\\
Let $t_1$, \ldots, $t_r$ be points of this affine line (we can have $r=0$). We equip $\AA^1(\gk) \setminus \{t_1, \ldots, t_r\}$ (affine line minus $r$ points) with a structure of an affine \vrt by forcing the invertibility of the $t-t_i$'s, \cad by defining

\snic {
\gB = \gk\big[t, (t-t_1)^{-1}, \ldots, (t-t_r)^{-1}\big] 
 \simeq 
\aqo{\gk[t,x]} {F(t,x)},
}

with $F(t,x) = (t-t_1) \cdots (t-t_r)\cdot x - 1$. This \klg $\gB$ then appears as the \alg of rational fractions over $\PP^1(\gk)$ defined everywhere except at the points~$\infty$ and~$t_i$.  It is an \acl and even a Bézout \ri (indeed, it is a localized \ri of $\gk[t]$).
\\
Analogously, for $n$ points  $t_1$, \ldots, $t_n$ of the affine line (with $n \ge 1$ this time), we can consider the \klg 

\snic {
\gA = \gk\big[(t-t_1)^{-1}, \ldots, (t-t_n)^{-1}\big]\subseteq \gk(t).
}

This \ri $\gA$ is a localized \ri
 of $\gk[(t-t_1)^{-1}]$ (which is \isoc to $\kX$) since by 
letting $v = (t-t_1)^{-1}$, we have $t-t_i = \big((t_1-t_i)v + 1\big)/v$.
So, 
$$\preskip.8pt
{\gA=\gk\big[v,\big((t_1-t_2)v+1\big)^{-1},\dots,\big((t_1-t_n)v+1\big)^{-1}\big]\subseteq \gk(v)=\gk(t)}.
\postskip2.5pt
$$
The \klg $\gA$  is therefore an \acl (and even a Bézout domain). 
By letting $p(t) = (t-t_1) \cdots (t-t_n)$, we easily have the \egt
$$
\preskip.4em \postskip.4em 
\gA = \gk[1/p, t/p, \ldots, t^{n-1}/p]. 
$$
The \klg $\gA$, constituted of rational fractions $u/p^s$ with $\deg(u) \le ns$, appears as that of rational fractions defined everywhere over $\PP_1(\gk)$ (including at the point $t = \infty$) except eventually at the points $t_i$.  In short, we can agree that~$\gA$ is the \klg of the \gui{functions} defined over the projective line minus the points $t_1$, \dots, $t_n$.

We study in this exercise a more \gnl case where $p$ is a \polu of degree $n \ge 1$.

Let $\gk$ be a \cdi and $p(t) = t^n + a_{n-1} t^{n-1} + \cdots + a_1 t + a_0\in\gk[t]$ $(n\geq1)$,
where~$t$ is an \idtrz. Let $x_i = t^i /p$.
 \\
Show that the \cli of $\gk[x_0]$ in $\gk(t)$ is the \klg

\snic{\gA = \gk[x_0, \ldots, x_{n-1}] = 
\sotq {u/p^s} {s\in\NN,\,u \in \gk[t],\  \deg(u) \le ns}.
}

In addition, $\Frac(\gA) = \gk(t)$. 
}
\end{exercise}

\vspace{-1em}
\begin{exercise}\label{exoPresentationAlgOuvertP1}
 {(A presentation of the \alg of functions over the projective line minus a finite number of points)}\\
{\rm  
The context is that of Exercise~\ref{exoAnneauOuvertP1}, but this time $\gk$ is an arbitrary \riz. Let $p = a_nt^n + \cdots + a_1t + a_0 \in \gk[t]$ be a \polu ($a_n = 1$) and

\snic {
\gA = \gk[1/p, t/p, \ldots, t^{n-1}/p].
}


Let $x_i = t^i/p$ for $i \in \lrb{0..n-1}$. We can write $\gA = \gk[\uX]/\fa$ where $(\uX) = (X_0, \ldots, X_{n-1})$ and $\fa$ is the \id of the relators between $(x_0, \ldots, x_{n-1})$.  It will be convenient to define $x_n$ by $x_n = t^n/p$; we therefore have $x_j = x_0t^j$ and

\snic {
\sum_{i=0}^n a_i x_i = 1  \quad \hbox {or yet} \quad
x_n = 1 - \sum_{i=0}^{n-1} a_i x_i.
}

The \egt on the right-hand side proves that $x_n \in \gA$.

\emph {1.}
Prove that the following family $R$ gives relators between the $x_j$'s.

\snic {
R \ : \quad
x_ix_j = x_kx_\ell  \qquad \hbox {for $i+j = k+\ell$}, \qquad
0 \le i, j, k, \ell \le n.
}

We define the subfamily $R_{\rm min}$ with ${n(n-1) \over 2}$ relators.

\snic {
R_{\rm min}\ : \quad x_i x_j = x_{i-1}x_{j+1}, \qquad 1 \le i \le j \le n-1.
}

\emph {2.}
Show that the family $R_{\rm min}$ (so $R$ also) generates the \id of the relators between the $x_i$'s for $i \in \lrb{0..n-1}$. 
In other words, if we let $\varphi : \gk[\uX] \to \gk[t,1/p]$ be the morphism defined by $X_i \mapsto x_i$ for $i \in \lrb {0..n-1}$, this means that $\Ker\varphi$ is generated by

\snic {
X_iX_j - X_{i-1}X_{j+1}, \qquad 1 \le i \le j \le n-1 \quad (\hbox{with }  X_n := 1 -
\sum_{i=0}^{n-1} a_iX_i)
.}

You may use the \kmo $\gk[X_0] \oplus \gk[X_0]X_1 \oplus \cdots \oplus \gk[X_0]X_{n-1}$.

}

\end{exercise}

\vspace{-1em}
\begin{exercise}\label{exoEntiersEmmanuel} {(Emmanuel's trick)}
\\
{\rm
Give a direct proof of item \emph{1} of Lemma~\ref{lemEmmanuel} without using \KROz's \thoz.}
\end{exercise}

\vspace{-1em}
\begin{exercise}\label{exoKroneckerTheorem} 
{(Another \dem of \KROz's \thoz)}\\
{\rm 
Consider the \pols
\[\preskip.4em \postskip.4em 
\begin{array}{ccc} 
  f(T) = a_0T^n + \cdots + a_n, \quad g(T) = b_0T^m + \cdots + b_m \;\;\hbox{and}  \\[.3em] 
h(T)
= f(T)g(T) = c_0T^{n+m} + \cdots + c_{n+m}.
 \end{array}
\] 
Kronecker's \thoz~\ref{thKro} affirms that each product $a_ib_j$ is integral over the \ri $\gA = \ZZ[c_0, \ldots, c_{n+m}]$.  
\\
It suffices to treat the case where the $a_i$'s and $b_j$'s are \idtrsz. Then in a \ri containing $\ZZ[a_0,\dots,a_n,b_0,\dots,b_m]$ we have

\snic {
f(T) = a_0(T-x_1) \cdots (T-x_n), \qquad
g(T) = b_0(T-y_1) \cdots (T-y_m)
.}

\emph {1.}
By using Emmanuel's trick (Lemma~\ref{lemEmmanuel}, with the \dem given in Exercise~\ref{exoEntiersEmmanuel}, independent of Kronecker's \thoz), show that for all~$I \subseteq \lrb{1..n}$, $J \subseteq \lrb{1..m}$, the product $a_0b_0 \prod_{i\in I}x_i \prod_{j\in J}y_j$ is integral over $\gA$. 

\emph {2.}
Conclude the result.

}
\end{exercise}

\vspace{-1em}
\begin{exercise}\label{exoAnneauCoinceBézout} 
{(Intermediary \ri $\gA \subseteq \gB \subseteq\Frac(\gA)$, Bézout case)}
 \\
{\rm  
Let $\gA$ be a Bézout domain, $\gK$ be its quotient field and $\gB$ be an intermediary \ri $\gA \subseteq \gB \subseteq \gK$.  Show that $\gB$ is a localized \ri of $\gA$ (therefore a Bézout \riz).
}
\end{exercise}

\vspace{-1em}
\begin{exercise}\label{exoGrellNoether}
{(Intermediary \riz, Pr\"ufer case)}\\
{\rm  
In this exercise we generalize the result of Exercise~\ref{exoAnneauCoinceBézout}
in the case where $\gA$ is a \ddp and we detail \thref{thSurAdp}. This is therefore a variation around the Grell-\Noe \tho (\paref{thGrellNoether}).

\emph {1.}
Let $x \in \gK=\Frac\gA$.

\vspace{-.6em} 
\begin{itemize}\itemsep0pt
\item [\emph {a.}]
Show that there exists an $s \in \Reg(\gA)$  such that $sx \in \gA$ and $1-s \in \gA x$.

\item [\emph {b.}]
Let $t \in \gA$ such that $tx = 1-s$. For every intermediary \ri  $\gA'$ between~$\gA$ and~$\gK$, show that $\gA'[x] = \gA'_s \cap \gA'_t$.
In particular, $\gA[x] = \gA_s \cap \gA_t$. Consequently, $\gA[x]$ is \iclz, and it is a \ddpz.
\end{itemize}

\vspace{-.5em} 
\emph {2.}
Show that every \tf \Aslg $\gB$ of $\gK$  is the intersection of a finite number of localized \ris of $\gA$ of the form $\gA_s$ with $s \in \gA$. Consequently, $\gB$ is \iclz, and it is a \ddpz.

\emph {3.}
Deduce that every intermediary \ri between $\gA$ and $\gK$ is a \ddpz.

\emph {4.}
Give an example of an \icl \ri $\gA$,  with an intermediary \ri $\gB$ between~$\gA$ and~$\Frac(\gA)$ which is not \icl (in particular, $\gB$ is not a localized \ri of $\gA$).

}

\end{exercise}

\vspace{-1.1em}
\begin{exercise}\label{exoPrimitivementAlg}
 {(To be primitively \agqz)}\\
{\rm  
Let $\gA = \aqo{\ZZ[A,B,U,V]}{AU+BV-1} = \ZZ[a,b,u,v]$ and $\gB = \gA[1/b]$. \\
Let $x = a/b$.
Show that $x$ is primitively \agq over $\gA$, but that~$y=2x$ is not.

}

\end{exercise}

\vspace{-1em}
\begin{exercise}\label{exocaracPruferC}  {(\Carns of the \adpcsz, 1)}
\\
{\rm  
Let $\gA$ be a \qiri and $\gK=\Frac\gA$.
\Propeq

\vspace{-1pt}
\begin{itemize}\itemsep0pt
\item [\emph {1.}]
$\gA$ is a Pr\"ufer \riz.

\item [\emph {2.}]
$\gA$ is normal and $x \in \gA[x^2]$ for every $x \in \gK$.

\item [\emph {3.}]
Every \ri $\gA[y]$ where $y\in\gK$ is normal.

\item [\emph {4.}]
Every intermediary \ri between $\gA$ and $\gK$ is normal.

\item [\emph {5.}]
$\gA$ is normal and $x \in \gA+ x^2 \gA$ for every $x \in \gK$.
\end{itemize}

}
\end{exercise}

\vspace{-1em}
\begin{exercise}
\label{exocaracPruCoh} (\Carns of the \adpcsz, 2)\\
{\rm  For any \qiri $\gA$,
\propeq
\vspace{-1pt}
\begin{itemize}
\item [\emph {1.}] $\gA$ is a \adpz.
\item [\emph {2.}]  Every \itf containing a \ndz \elt is \ivz.
\item [\emph {3.}]  Every \id $\fa=\gen{x_1,x_2}$ with $x_1$, $x_2\in\Reg(\gA)$ is \ivz.
\item [\emph {4.}]  For all $a$, $b\in\gA,$ we have $\gen{a,b}^2=\gen{a^2,b^2}=\gen{a^2+b^2,ab}.$
\item [\emph {5.}]  For all $f$, $g\in\gA[X]$, we have $\rc(f)\rc(g)=\rc(fg)$.
\end{itemize}
 
}
\end{exercise}

\vspace{-1em}
\begin{exercise}\label{exoAnarlgb} (A \gnn of Proposition~\ref{propAriCohZed})\\
{\rm  Let $\gA$ be a \lgb \adpc (e.g.\ \rdt \zedz). 
\begin{itemize}
\item [\emph {1.}]
 Every matrix is \eqve to a matrix in Smith form (\cad $\gA$ is a Smith \riz).
\item  [\emph {2.}] Every \pf \Amo is characterized, up to \isoz, by its \idfsz. Actually it is \isoc to a direct sum of cyclic modules~$\gA\sur{\fa_k}$ with \idps $\fa_1\subseteq\cdots\subseteq\fa_n$ ($n\geq0$).
\\ 
Note: We can naturally deduce an analogous \gnn of Corollary~\ref{corpropAriCohZed}.
\end{itemize} 
}
\end{exercise}

\vspace{-1em}
\begin{exercise}\label {exoReductionIdeal} 
  {(Reduction \id of another \idz)}
{\rm  
\vspace{-1pt}
\begin{itemize}
\item [\emph {1.}] Let $E$ be an $\gA$-module generated by $n$ \eltsz, $b \in \gA$
and $\fa$ be an \id such that~$bE \subseteq \fa E$.
Show that there exists a $d= b^n + a_{1} b^{n-1} + \cdots + a_{n-1}b + a_n$, 
with the~$a_i \in \fa^i$, that annihilates $E$.
\end{itemize}
We say that an \id $\fa$ is a \emph{reduction} of an \id $\fb$ if $\fa \subseteq \fb$ and if $\fb^{r+1} = \fa\fb^r$ for a certain exponent $r$ (it is then true for all the larger exponents).
\begin{itemize} 
\item [\emph {2.}]
Let $f$, $g \in \gA[\uX]$. Prove that $\rc(fg)$ is a reduction of $\rc(f)\rc(g)$.
\item [\emph {3.}]
In $\gA[X,Y]$, show that $\fa = \gen {X^2, Y^2}$ is a reduction of $\fb = \gen {X,Y}^2$.\\
Show that $\fa_1 = \gen {X^7, Y^7}$ and $\fa_2 = \gen {X^7, X^6Y + Y^7}$ are reductions of %
the \idz~$\fb' = \gen {X^7, X^6Y, X^2Y^5, Y^7}$. Give the smallest possible exponents.
\item [\emph {4.}]
Let $\fa \subseteq \fb$ be two \ids with $\fb$ \tfz. Show that $\fa$ is a reduction of~$\fb$ \ssi $\Icl(\fa) = \Icl(\fb)$.
\end{itemize}
}
\end{exercise}

\vspace{-1em}
\begin{exercise}  \label{exolemNormalIcl} {(Normal \qiriz)}\\
{\rm  
Here is a light \gnn of Fact~\ref{lemNormalIcl}.
By \Pbmz~\ref{exoAnneauNoetherienReduit} the hypothesis is satisfied for the \fdi reduced \coh \noe \ris (in \clama they are the reduced \noe \risz).

Consider a reduced \ri $\gA$. Suppose that its total \ri of fractions is \zedz. 

\emph{1.} If $\gA$ is normal, it is a \qiriz.

\emph{2.} The \ri $\gA$ is normal \ssi it is \icl in $\Frac\gA$.

}

\end{exercise}

\vspace{-1em}
\begin{exercise}  \label{exoEntierSurIX} {(Integral \pol over $\fa[X]$)}

{\rm  
Let $\gA \subseteq\gB$ be two \risz, 
$\fa$ be an \id of $\gA$ and $\fa[X]$ be the \id of $\gA[X]$ constituting of \pols with \coes in $\fa$. 
For $F \in \gB[X]$, show that $F$ is integral \hbox{over $\fa[X]$} \ssi each \coe of $F$ is integral over $\fa$.
}

\end{exercise}

\vspace{-1.2em}
\pagebreak	

\begin{exercise}  \label{exosdirindec} {(Indecomposable modules)}

{\rm  
We say that a module $M$ is \emph{indecomposable} if the only 
direct summand submodules of~$M$
 are $0$ and $M$.
The goal of the exercise is to prove that over a \fac\dDkz, every \mpf is a direct sum of a  finite number of indecomposable modules, this \dcn being unique up to the order of terms when the  module is a torsion module.

\emph{1.} Let $\gA$ be a \ri and $\fa$ be an \idz. If the \Amo $M=\gA/\fa$ is a direct sum of two submodules $N$ and $P$ we have $N=\fb/\fa$, $P=\fc/\fa$ with  \com $\fb\supseteq \fa$ and~$\fc\supseteq \fa$. More \prmtz, $\fb=\gen{b}+\fa$, $\fc=\gen{c}+\fa$, where $b$ and $c$ are \cop \idms modulo~$\fa$. 

\emph{2.} Let $\gZ$ be a \dDkz. 
 
\emph{2a.} Show that a \mrc $1$ is indecomposable.
 
\emph{2b.} Show that a cyclic module $\gZ/\fa$ with $\fa$ \tfz, $\neq \gen{0},\,\gen{1}$ is indecomposable \ssi $\fa=\fp^{m}$ for some \idema $\fp$ and some~$m\geq 1$.
 
\emph{2c.} Deduce that if $\gZ$ admits \facsz, every \mpf is a direct sum of a finite number of indecomposable modules.

\emph{3.} When the module is a torsion module, show the uniqueness of the \dcn with a meaning to be specified. 
}

\end{exercise}

\vspace{-1em}
\begin{problem}\label{exoArithInvariantRing} 
{(Sub\ri of invariants under a finite group action and \ari \crcz)}
{\rm Note: See also \Pbmz~\ref{exoGaloisNormIdeal}.\\  
\emph {1.}
If $\gA$ is a \anorz, every \lop \id 
is \iclz. Consequently, if $f$, $g \in \AX$ with $\rc(fg)$
\lopz, then $\rc(f)\rc(g) = \rc(fg)$.

\emph {2.}
Suppose that $\gA$ is normal and that $\gB\supseteq\gA$ is integral over $\gA$.
If $\fa$ is a \lop \id of~$\gA$, then $\fa\gB\cap\gA = \fa$.

\emph {3.}
Let $(\gB, \gA, G)$ where $G \subseteq \Aut(\gB)$ is a finite group and $\gA =\Fix_\gB(G)= \gB^G$. If $\fb$ is an \id of $\gB$, let $\rN'_G(\fb) = \prod_{\sigma \in G} \sigma(\fb)$ (it is an \id of $\gB$)  
 and $\rN_G(\fb) = \gA \cap \rN'_G(\fb)$. 

Suppose that $\gB$ is normal and that $\gA$ is a Pr\"ufer \ri (therefore normal). 
\begin{itemize}
\item [{\it a.}]
For $b \in \gB$, prove that $\rN'_G(b\gB)=\rN_G(b)\gB$ and $\rN_G(b\gB)=\rN_G(b)\gA$.

\item [{\it b.}]
If $\fb$ is a \itf of $\gB$, show that $\rN_G(\fb)$ is a \itf of $\gA$ and that $\rN'_G(\fb) = \rN_G(\fb)\gB$. 
You can write $\fb = \gen {b_1, \ldots, b_n}$ and introduce $n$ \idtrs $\uX = (X_1, \ldots, X_n)$ and consider the normic \pol $h(\uX)$

\snic {
h(\uX) = \prod
_{\sigma \in G} h_\sigma(\uX) 
\quad \hbox {with} \quad
h_\sigma(\ux) = \sigma(b_1) X_1 + \cdots +  \sigma(b_n) X_n.
}

\item [{\it c.}]
For \itfs $\fb_1$, $\fb_2$ of $\gB$, we obtain $\rN_G(\fb_1\fb_2) = \rN_G(\fb_1) \rN_G(\fb_2)$.

\item [{\it d.}]
A \itf $\fb$ of $\gB$ is \iv \ssi $\rN_G(\fb)$ is \iv in $\gA$.
\end{itemize}

Note: We know that $\gB$ is a \adp (\thref{thExtEntPruf}); in the case where $\gB$ is integral, question \emph {3d} provides a new proof for it.

\emph {4.}
Let $\gk$ be a \cdi with $2\in\gk\eti$ and $f(X) \in \gk[X]$ be a \spb \poluz. The \pol $Y^2 - f(X) \in \gk[X,Y]$ is absolutely \ird (see Exercise~\ref{exoNotAbsIntClos}); let $\gk[x,y] =\aqo {\gk[X,Y]}{Y^2 - f(X)}$. 
Show that $\gk[x,y]$ is a \adpz.

}

\end{problem}

\vspace{-1em}
\begin{problem}\label{exoFullAffineMonoid}
 {(Full sub\mos of $\NN^n$)}
\\
{\rm  
Let $M\subseteq\NN^n$ be a sub\moz; for a \ri $\gk$, let $\gk[M]$ be the \klg of the \mo $M$. It is the \kslg of $\gk[\NN^n] \simeq \gk[\ux] = \gk[\xn]$ generated by the monomials $\ux^m = x_1^{m_1} \cdots x_n^{m_n}$ for $m \in M$.
We say that $M$ is a \emph{full sub\moz} of $\NN^n$ if for $m \in M$, $m' \in \NN^n$, we have $m+m' \in M \Rightarrow m' \in M$.

\emph {1.}
The subgroup of $\ZZ^n$ generated by $M$ is equal to $M-M$, and if $M$ is full, then $M = (M-M)\cap\NN^n$. Conversely, if $L \subseteq \ZZ^n$ is a subgroup, then the \mo $M = L\cap\NN^n$ is a full sub\mo of $\NN^n$.

\emph {2.}
Let $M \subseteq \NN^n$ be a full sub\mo and $\gk$ be a \cdiz.
\begin{enumerate}
\item [\emph {a)}]
Let $\gA=\gk[M] \subseteq \gB=\kux$. Show that if $a \in
\gA\setminus\{0\}$,  $b \in \gB$, and~$ab \in\gA $,  then~$b \in\gA$.

\item [\emph {b)}]
Let $\gA\subseteq\gB$ be two domains satisfying: if $a \in \gA\setminus\{0\}$,  $b \in \gB$, and~$ab \in\gA $,  then~$b \in\gA$.  

\begin{enumerate}
\item [\emph {i.}]
Show that $\gA = \gB\cap\Frac(\gA)$; deduce that if $\gB$ is \iclz, the same goes for $\gA$.
\item [\emph {ii.}]
In particular, if $M \subseteq \NN^n$ is a full sub\moz, then $\gk[M]$ is \icl for every \cdi $\gk$.
\item [\emph {iii.}]
More \gnltz, if $\gB\subseteq\gC$ is \icl in $\gC$, then $\gA$ is \icl in $\gC\cap\Frac(\gA)$.
\end{enumerate}
\end{enumerate}

\emph {3.}
Let $M \subseteq \NN^n$ be the sub\mo of magic squares (see Exercise~\ref{exoMatMag3}); then $\gk[M]$ is \icl for every \cdi $\gk$.

}

\end{problem}

\vspace{-1em}
\begin{problem}\label{exoBaseNormaleAlInfini} 
 {(Normal basis at infinity)}\\
{\rm 
A \ddv $\gB$ with quotient field $\gK$ is a DVR if $\gK\eti\!\sur{\gB\eti}\simeq \ZZ$ (\iso of ordered groups).
A regular parameter is  every \elt $b\in\gB$ such that $v(b)=1$, where
$v:\gK\eti\to\ZZ$ is the map defined via the previous \iso (this map $v$ is also called a \emph{valuation}).
Every \elt $z$ of $\gK\eti$ is then of the form $ub^{v(z)}$ with $u\in\gB\eti$.%

Let $\gk$ be a \cdiz, $t$ be an \idtr over $\gk$, $\gA = \gk[t]$, $\gA_\infty = \gk[t^{-1}]_{\langle t^{-1} \rangle}$, and~$\gK = \Frac(\gA) = \gk(t) = \Frac(\gA_\infty) = \gk(t^{-1})$.  If $L$ is a finite dimensional \Kevz, we study in this \pb the intersection of an $\gA$-lattice of~$L$ 
and of an $\gA_\infty$-lattice of~$L$ (see the \dfns question \emph {2}),
an intersection which is always a finite dimensional \kevz. 

In the theory of algebraic function fields this study is at the basis of the determination of Riemann-Roch spaces, however, when certain integral closures are known by bases; as a subproduct, we determine the algebraic closure of $\gk$ in a finite extension of $\gk(t)$.
\\
The \ri $\gA_\infty$ is a DVR; let $v : \gK \to \ZZ \cup \{\infty\}$ be the corresponding valuation, defined by $v = -\deg_t$, and we fix $\pi = t^{-1}$ as regular parameter.  If $x = \tra {[x_1, \ldots, x_n]}$, let $v(x) = \min_i v(x_i)$.  This allows us to define a \emph{modular reduction}
$$\preskip.4em \postskip.3em 
\gK^n
\setminus \so{0} \to \gk^n \setminus \so{0},\; \quad x \mapsto \xi = \ov {x}, 
$$
with
$\xi_i = (x_i/\pi^{v(x)}) \mod \pi \in \gk$.  
\\
\Gnltz, if $\gV$ is a \adv of a \cdi $\gK$, of residual field $\gk$, we have a \emph{reduction}

\snic {
\PP^{m}(\gK) \to \PP^{m}(\gk),  \qquad
(x_0 : \ldots : x_m) \mapsto (\xi_0 : \ldots : \xi_m)
\quad \hbox {with $\xi_i = \ov {x_i/x_{i_0}}$}
,}

where $x_{i_0} \,\vert\, x_i$ for all $i$; the \elt $(\xi_0: \ldots: \xi_n)\in\PP^{m}(\gk)$ is well-defined: it corresponds to a \vmd of $\gV^{m+1}$. In short we have an \gui{\isoz} $\PP^{m}(\gV)\simeq\PP^{m}(\gK)$ and a {reduction} $\PP^{m}(\gV)\to\PP^{m}(\gk)$.
\\
		Here the choice of the regular parameter $\pi = t^{-1}$ gives a direct \dfn of the reduction  $\gK^n \setminus\{0\} \to \gk^n \setminus \{0\}$, without having to change the coordinates on the projective line to understand what is happening at infinity.

We will say that a matrix $A \in \GL_n(\gK)$ with columns $(A_1, \ldots, A_n)$ is \emph{$\gA_\infty$-reduced} if the matrix $\ov A \in \MM_n(\gk)$  is in $\GL_n(\gk)$.

\emph {1.}
Let $A \in \GL_n(\gK)$ of columns $A_1, \ldots, A_n$. Show that $\sum_{j=1}^n v(A_j) \le v(\det A)$.

\emph {2.}
Let $A \in \GL_n(\gK)$;  compute $Q \in \GL_n(\gA)$ such that $AQ$ is $\gA_\infty$-reduced. Or yet again, let $E \subset \gK^n$ be an $\gA$-lattice, \cad a free $\gA$-module of rank $n$; then $E$ admits an $\gA_\infty$-reduced $\gA$-basis (a basis $(A_1, \ldots, A_n)$ such that $(\ov {A_1}, \ldots, \ov {A_n})$ is a $\gk$-basis of $\gk^n$). 
You can start with the example $A = \cmatrix {\pi^2 &\pi\cr 1 & 1\cr}$.

\emph {3.}
For $P \in \GL_n(\gA_\infty)$, prove the following points.
\begin{itemize}
\item [\emph {a.}]
$P$ is a $v$-isometry, \cad $v(Px) = v(x)$ for every $x \in \gK^n$.

\item [\emph {b.}]
For $x \in \gK^n \setminus \{0\}$, $\ov {Px} = \ov{P}\,\ov{x}$.

\item [\emph {c.}]
If $A \in \GL_n(\gK)$ is $\gA_\infty$-reduced, the same goes for $PA$.
\end{itemize}

\emph {4.}
Let $A \in \GL_n(\gK)$ be triangular. What is the meaning of \gui{$A$ is $\gA_\infty$-reduced}?

\emph {5.}
Let $A \in \GL_n(\gK)$. Show that there exists a $Q \in \GL_n(\gA)$, $P \in \GL_n(\gA_\infty)$ and integers $d_i \in \ZZ$ such that $PAQ = \Diag(t^{d_1}, \ldots, t^{d_n})$; moreover, if we order the $d_i$'s by increasing order, they are unique.

\emph {6.}
Let $L$ be a \Kev of dimension $n$, $E \subset L$ be an $\gA$-lattice,  and $E'\subset L$ be an $\gA_\infty$-lattice. 
\begin{itemize}
\item [\emph {a.}]
Show that there exist an $\gA$-basis $(e_1,\ldots, e_n)$ of $E$, an $\gA_\infty$-basis $(e'_1, \ldots, e'_n)$ of $E'$ and integers $d_1$, \ldots, $d_n \in \ZZ$ satisfying $e'_i = t^{d_i} e_i$ for $i\in \lrb{1..n}$. Moreover, the $d_i$'s ordered in increasing order only depend on $(E, E')$.
\item [\emph {b.}]
Deduce that $E \cap E'$ is a finite dimensional \kevz.
More \prmtz,
$$\preskip.2em \postskip-.1em\ndsp 
E \cap E' = \bigoplus_{d_i \ge 0} \bigoplus_{j=0}^{d_i} \gk t^j e_i
, 
$$
and in particular,
$$\preskip.1em  
\ndsp
\dim_\gk(E\cap E') = \sum_{d_i \ge 0} (1+d_i) = \sum_{d_i \ge -1} (1+d_i). 
$$
\end{itemize}

\vspace{-1.2em}
\pagebreak

\emph {7.}
Suppose that $\gL$ is a finite $\gK$-extension of degree $n$. We define integral closures in $\gL$: $\gB$ that of $\gA$, $\gB_\infty$ that of $\gA_\infty$ and $\gk'$ that of $\gk$. We say that a basis $(\ue) = (e_1, \ldots, e_n)$ of $\gB$ over $\gA$ is \emph {normal at infinity} if there exist $r_1$, \ldots, $r_n \in \gK^*$ such that $(r_1e_1, \ldots,r_ne_n)$ is an $\gA_\infty$-basis of $\gB_\infty$. Show that the \elts of the basis $(\ue)$ \gui{integral at infinity,} that is which are  members of $\gB_\infty$, form a $\gk$-basis of the extension~$\gk'$.

\emph {8.}
Let $\gk=\QQ$, $\gL = \aqo{\gk[X,Y]}{X^2 + Y^2} = \gk[x,y]$, $\gA =
\gk[x]$. \\
Show that $(y+1, y/x)$ is an $\gA$-basis of $\gB$ but that it is not normal at infinity.  Explicate an $\gA$-basis of $\gB$ normal at infinity.

}
\end{problem}

\vspace{-1em}
\begin{problem}\label{exoHyperEllipticFunctionRing}
{(\Ri of functions of an affine hyper-elliptic curve having a single point at infinity)}\\
{\rm  
Here we will use a notion of a \emph{norm of an \idz} in the following context: $\gB$ being a free \Alg of finite rank $n$ and $\fb$ being a \itf of $\gB$, the norm of $\fb$ is the \id 

\snic{\rN_{\gB\sur\gA}(\fb)=\rN(\fb) \eqdefi \cF_{\gA,0}(\gB\sur\fb)\subseteq\gA.}

 It is clear that for $b \in \gB$, $\rN(b\gB) = \rN_{\gB\sur\gA}(b)\gA$, that $\rN(\fa\gB) = \fa^n$ for $\fa$ a \itf of $\gA$ and that $\fb_1 \subseteq \fb_2 \Rightarrow \rN(\fb_1) \subseteq \rN(\fb_2)$.%
\index{norm!of an ideal}

Let $\gk$ be a field of \cara $\ne 2$ and $f = f(X)\in \gk[X]$ be a \spl \polu of odd degree $2g+1$. The \pol $Y^2 - f(X) \in \gk[X,Y]$ is absolutely \irdz; let $\gB = \aqo {\gk[X,Y]}{Y^2 - f(X)} = \gk[x,y]$ and $\gA = \gk[x] \simeq
\gk[X]$. 
The \ri $\gB$ is integral, it is a free \Amo of basis $(1,y)$. 
For $z=a + by$ with $a$, $b \in \gA$, let $\ov z=a-yb$, and $\rN = \rN_{\gB\sur\gA}$: $\rN(z) = z\ov z = a^2 - fb^2$. 
\\
The goal of the \pb is to parameterize the nonzero \itfs of $\gB$, to show that $\gB$ is a \adp and to study the group $\Cl(\gB)$ of classes of \iv \ids of $\gB$.
\\
If $\fb$ is a \itf of $\gB$, its \emph{content} is the \idf $\cF_{\gA,1}(\gB\sur\fb)$. 
\\
To two \elts $u$, $v \in \gA$ satisfying $v^2 \equiv f \bmod u$, we associate the \Asub \hbox{of $\gB$}: $\fb_{u,v} = \gA u + \gA (y-v)$. We have $u \neq  0$ because $f$ is \splz. We will sometimes make the \pol $w \in \gA$ intervene so that $v^2 - uw = f$ and we will write $\fb_{u,v,w}$ instead of $\fb_{u,v}$ (even if $w$ is completely determined by $u,v$). 

\emph {1.}
Show that $\fb_{u,v}$ is an \id of $\gB$ and that $\fb_{u,v} = \gA u \oplus \gA (y-v)$. Conversely, for $u, v \in \gA$, if $\gA u + \gA (y-v)$ is an \id of $\gB$, then $v^2 \equiv f \bmod u$.

\emph {2.}
Show that $\gA \to \gB\sur{\fb_{u,v}}$ induces an \iso $\gA\sur{\gA u} \simeq \gB\sur{\fb_{u,v}}$; consequently, $\Ann_\gA(\gB\sur{\fb_{u,v}})=\gA u$.
Deduce \gui{the uniqueness of $u$}

\snic {u_1,\,u_2
\hbox{ \mons and } \fb_{u_1,v_1} = \fb_{u_2,v_2} \;\; \Longrightarrow  \;\;u_1 = u_2.
}


Also prove that $\rN(\fb_{u,v}) = u\gA$ and that $v$ is unique modulo $u$

\snic {
\fb_{u,v_1} = \fb_{u,v_2}  \iff  v_1 \equiv v_2 \mod u
.}

\emph {3.}
Show that
$$\preskip.0em \postskip.3em 
 \fb_{u,v,w}\fb_{w,v,u} = \gen {y-v}_\gB, \qquad
\fb_{u,v}\fb_{u,-v} = \gen {u}_\gB.
$$
Consequently, the \id $\fb_{u,v}$ is \ivz.  \\
In addition, for $u = u_1u_2$ satisfying $v^2 \equiv f \bmod u$, we have $\fb_{u,v} = \fb_{u_1,v}\,\fb_{u_2,v}$.

\emph {4.}
Let $\fb$ be a nonzero \itf of $\gB$. 
\begin{itemize}
\item [\emph {a.}]
Show that there exist two unique \polus $d$, $u \in \gA$ and $v \in \gA$ 
with $v^2 \equiv f \bmod u$, so that $\fb = d\,\fb_{u,v}$. Consequently, $\fb$ is an \iv \id (so $\gB$ is a \adpz).  In addition, $v$ is unique modulo $u$, therefore unique if we impose $\deg v < \deg u$.

\item [\emph {b.}]
Deduce that $\fb\ov\fb = \rN(\fb)\gB$ then that the norm is multiplicative over the \idsz.

\item [\emph {c.}]
Show that $\gB\sur\fb$ is a finite dimensional \kevz.\\
Show that $\dim_\gk (\gB\sur\fb) = \dim_\gk (\gA\sur\fa)$ with $\fa = \rN(\fb)$. This integer will be  
denoted by $\deg(\fb)$. Prove that $\deg(\fb_{u,v}) = \deg u$, that $\deg (\fb) = \deg \rN(\fb)$, and finally that $\deg$ is additive, \cad $\deg(\fb_1\fb_2) = \deg(\fb_1) + \deg(\fb_2)$.
\end{itemize}

Let $u$, $v \in \gA$ with $v^2 \equiv f \bmod u$. We say that \emph{the pair $(u,v)$ is reduced} if~$u$ is \mon and \framebox [1.1\width]{$\deg v < \deg u \le g$}.
By abuse of language, we also say that $\fb_{u,v}$ is reduced.
For example, if $(x_0, y_0)$ is a point of the hyper-elliptic curve $y^2 = f(x)$, its \id $\gen {x-x_0, y-y_0}$ is a reduced \id (take $u(x) = x-x_0$, $v = y_0$).

\emph {5.}
Show that every nonzero \itf of $\gB$ is associated with a reduced ideal of $\gB$ (two \ids $\fa$ and $\fa'$ are said to be \emph{associated} if there exist two \ndz \elts $a$ and $a'$ such that $a\fa'=a'\fa$, we then let $\fa\sim\fa'$).

\emph {6.}
In this question, for a nonzero \itf $\fb$ of $\gB$, we designate by~$\rN(\fb)$ the \polu \gtr of the \id $\rN_{\gB/\gA}(\fb)$.  Let $\fb_{u,v}$ be a reduced \idz.
\begin{itemize}\itemsep0pt
\item [\emph {a.}]
Let $z \in \fb_{u,v} \setminus \{0\}$ such that $u = \rN(\fb_{u,v}) \divi \rN(z)$, \cad $\rN(z)/\rN(\fb_{u,v})$ is a \polz. Show that
$$\preskip.0em \postskip.4em
\deg \left(\rN(z)/\rN(\fb_{u,v})\right) \ge \deg u,
$$
with \egt \ssi $z \in \gk\eti u$. 

\item [\emph {b.}]
Let $\fb'$ be a \itf of $\gB$ satisfying $\fb' \sim \fb_{u,v}$. Show that $\deg(\fb') \ge \deg(\fb_{u,v})$ with  \egt \ssi $\fb' = \fb_{u,v}$. In summary, in a class of invertible \ids of $\gB$, there is therefore one and only one \id of minimum degree: it is the unique reduced \id of the class.
\end{itemize}
\emph {7a.}
Show that the affine curve $y^2 = f(x)$ is smooth; more precisely, by letting $F(X,Y) = Y^2 - f(X) \in \gk[X,Y]$, show that $1 \in \gen {F, F'_X, F'_Y}$; this uniquely uses the fact that $f$ is \spb and that the \cara of $\gk$ is not $2$, not the fact that $f$ is of odd degree.
\\
If $\gk$ is \acz, we thus obtain a biunivocal correspondence between the points $p_0 = (x_0,y_0)$ of the affine curve $y^2 = f(x)$ and the DVRs $\gW$ of $\gk(x,y)$ containing $\gB = \gk[x,y]$: to $p_0$, we associate its local \ri $\gW$ and in the other direction, to $\gW$ we associate the point $p_0 = (x_0,y_0)$ such  
that $\gen {x-x_0, y-y_0}_\gB = \gB \cap \fm(\gW)$.

\emph {b.}
We now study \gui{the points at infinity of the smoothed projective curve,} at infinity relative to the model $y^2 = f(x)$. \Agqtz, these are \advs for $\gk(x,y)$ not containing $\gB$ (but containing $\gk$ of course).  
Let \smash{$\gA_\infty = \gk[x^{-1}]_{\gen {x^{-1}}}$} be the DVR. Show that there exists one and only one \ri $\gB_\infty$, $\gA_\infty \subseteq \gB_\infty \subseteq \Frac(\gB) = \gk(x,y)$, having $\Frac(\gB)$ as quotient field. Show that $\gB_\infty$ is a DVR, that $\gB_\infty\sur{\fm(\gB_\infty)} \simeq \gA_\infty\sur{\fm(\gA_\infty)} \simeq \gk$ and that it is the only point at infinity.
}
\end{problem}

\vspace{-1.2em}
\pagebreak	

\begin{problem}\label{exoTrifolium}
{(Trifolium: \cli and parameterization)}
\\
{\rm  Let $\gk$ be a \cdi and 

\snic{F(X,Y) = (X^2+Y^2)^2 + \alpha X^2Y + \beta Y^3,}

with $\alpha \ne \beta$ in $\gk$.

\Deuxcol{.50}{.45}
{
We study the curve $F(x,y) = 0$, its singular points, its field of functions 

\snic{\gL = \gk(x,y)}

(we will show that $F$ is \irdz), its \ri of functions $\gk[x,y]$, the \cli $\gB$ of $\gk[x,y]$ in $\gL$~\ldots\ {etc}\ \dots \\
Note that $F(-X,Y) = F(X,Y)$ and therefore that the involution $(x,y) \mapsto (-x,y)$ leaves the curve $F(x,y) = 0$ invariant.
\\
Opposite is  an example of such a curve.}
{~\vspace{-5mm}

\includegraphics*[width=5cm]{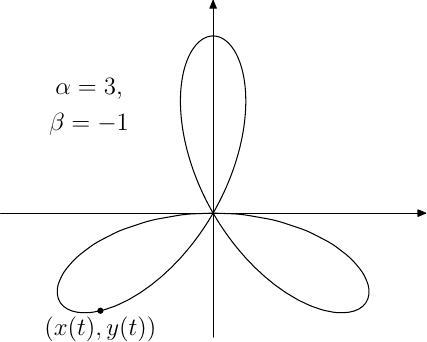}}

\vspace{2mm}
\emph {1.}
Show that $F$ is an absolutely \ird \polz. More \gnltz: let $\gk$ be an integral \riz, $\gk[\uT]$ be a \pol \ri with several \idtrs and $F \in \gk[\uT]$, $F = F_N + F_{N+1}$ with nonzero \hmg $F_N$, $F_{N+1}$, of degrees $N$, $N+1$, respectively. Then, in every \fcn $F = GH$, one of the two \pols $G$ or $H$ is \hmgz; finally, if $\gk$ is a field, then $F$ is \ird \ssi $F_N$, $F_{N+1}$ are coprime.

\emph {2.}
Determine the singular points of the curve $F = 0$.

Let $\gL = \gk(x,y)$ and $\gB$ be the \cli of $\gk[x,y]$ in $\gL$.

\emph {3.}
Let $t = y/x$ be such that $\gL = \gk(x,t)$.  
\begin{itemize}
\item [\emph {a.}]
Determine a primitive \agq \eqn of $t$ over $\gk[x]$. \\
Let $G(X,T) = a_4 T^4 + \cdots + a_1 T + a_0 \in \gk[X][T]$, with $a_i=a_i(X) \in \gk[X]$,
such a primitive \polz, therefore satisfying $G(x,t) = 0$. Prove that $(x, t)=(0,0)$ is a nonsingular point of the curve $G = 0$.

\item [\emph {b.}]
Using Emmanuel's trick (Lemma~\ref{lemEmmanuel}), determine the integral \elts  $b_4$, \ldots, $b_1$ associated with $(G, t)$ with $\gA = \gk[x]$ as a base \riz.
Deduce a \mlp for $(x,y)$ and describe the \id $\fq$ of $\gB$ such that $\gen{x}_\gB = \fq \gen {x,y}_\gB$.
\end{itemize}

\emph {4.}
Show that $\gL = \gk(t)$ and express $x$, $y$ as \elts of $\gk(t)$.  

\emph {5.}
Determine the \cli $\gB$ of $\gk[x,y]$ in $\gL$.
\begin{itemize} 
\item [\emph {a.}]
Show that $\gB = \gk[g_0, g_1]$ with $g_0 = 1/(1+t^2)$ and $g_1 = tg_0$.
Express $x$, $y$ in $\gk[g_0, g_1]$. What is \gui{the \eqnz} relating $g_0$ and $g_1$?
\item [\emph {b.}]
Show that $(1, y, b_3t, b_2t)$ is an $\gA$-basis of $\gB$.
\item [\emph {c.}]
Prove that $\dim_\gk  \gB\sur{\gen{x,y}_\gB} = 3$.
\end{itemize}

\emph {6.}
Let $\gV$ be the \advz\footnote{A sub\ri $\gV$ of a \cdi $\gL$ is called a \emph{\adv of $\gL$} if for all $x\in\gL\eti$ we have $x\in\gV$ or $x^{-1}\in\gV$.} of $\gL$ defined by the nonsingular point $(0, 0)$ of the curve $G = 0$. It is the only \adv of $\gL$ containing $\gk$ and such \hbox{that $x$, $t \in \Rad\gV$} (and so $y \in \Rad\gV$ also). \\
Consider the \idep $\fp_1 = (\Rad\gV) \cap \gB$. Show that

\snic {
\fp_1 = \gen {x, y, b_4t, b_3t, b_2t, b_1t} = \gen {g_0-1, g_1} 
\quad \hbox{and} \quad \gB\sur{\fp_1} = \gk
,}

and prove that $\fp_1^2 = \gen {g_0-1, g_1^2}$.

\emph {7.}
Determine the \fcn in $\gB$ of the \id $\gen {x,y}_\gB$ as a product of \idepsz. The response is not uniform at $(\alpha, \beta)$, unlike the determination of the \cli $\gB$ of $\gA$.

\emph {8.}
Repeat the questions by only assuming that $\gk$ is an \icl \ri and that $\beta-\alpha\in\gk\eti$.

}

\end{problem}

\sol

\exer{exoDetTrick1}{We need to show the inclusion $\fb^n\subseteq\fa\fb^{n-1}$.
Let $(x_1, \ldots, x_n)$ be a \sgr of $E$, $X=\tra{ [\,{x_1\;\cdots\; x_n}\,]}$,  $b_1$, \ldots, $b_n \in \fb$ and  $B=\Diag(b_1, \ldots, b_n)$.
Since $b_ix_i \in \fa E$ ($i\in\lrbn$), there exists an  $A \in \Mn(\fa)$ such that $B\,X = A \,X$. \\
Let $C = B - A$. We have $C \,X = 0$, and since $E$ is faithful, $\det C  = 0$.
Expanding this \deterz, we obtain $b_1 \cdots b_n + a = 0$ 
with $a \in \fa\fb^{n-1}$ (since $\fa\subseteq\fb$).
}

\exer{exoMatLoc22}\emph{(Principal \lon matrices in $\MM_2(\gA)$)}\\
\emph {1.}
Immediate, because if $B = \cmatrix {b_{11}&b_{12}\cr b_{21}&b_{22}\cr}$, then $\wi B = \crmatrix {b_{22}&-b_{12}\cr -b_{21}&b_{11}\cr}$ and $[x \; y\,]B = [x' \; y'\,]$ with

\snic {
x' = -\left|\matrix {-b_{21} & b_{11}\cr x & y\cr}\right|, \quad
y' = \left|\matrix {b_{22} & -b_{12}\cr x & y\cr}\right|
.}

\emph {2.}
We have $u$, $v \in \fb$ with $z = ux + vy$ and $ux$, $uy$, $vx$, $vy$ are multiples of $z$, which we write as $\crmatrix {y\cr -x\cr} [v \; -u\,] = zB$.
As $[x \; y\,] \crmatrix {y\cr -x\cr} = 0$, we have $[x \; y\,]zB = 0$; in \hbox{addition $\Tr(zB) = yv + xu = z$}.

\emph {3.}
In the lemma in question, $z = x^n$ and the \ri is a \lsdz.  The \egts $x^n[x\; y\,]B = 0$ and $x^n\big(1 - \Tr(B)\big) = 0$ provide two \come \lons of $\gA$: one in which $x^n = 0$, in which case $x = 0$ because the \ri $\gA$ and its localized \ris are reduced, and the other in which $[x\; y\,]B = 0$ and $\Tr(B) = 1$. In each one of them, $\gen {x,y}$ is \lop therefore it is \lop in $\gA$.

\exer{exoPrufNagata}  \emph{1.} Indeed, $\gA(X)$ is \fpt over $\gA$.

 \emph{2.}
Let $f=\sum_{k=0}^n a_kX^k\in\AX$
. 
For each $k$, we have, in $\gA$, an \egt 

\snic{\gen{a_0,\ldots,a_n}\gen{b_{0,k},\ldots,b_{n,k}}=\gen{a_k}
$  with $ a_0b_{0,k}+\cdots +a_nb_{n,k}=a_k.}

Consider then the \pol $g_k=\sum_{j=0}^n b_{j,k}X^{n-j}$. 
All of the \coes of $fg_k$ are in $\gen{a_k}$. We can therefore write $fg_k=a_kh_k$ with the \coe of degree~$k$ in~$h_k$ equal to $1$.
This implies that in $\gA(X)$, $a_k\in\gen{f}$. However, we have $f\in\gen{a_0,\ldots,a_n}$ in $\AX$. Thus, in  $\gA(X)$, $\gen{f}=\gen{a_0,\ldots,a_n}$.
\\
We deduce that $\gA(X)$ is a Bézout \riz, because for $f_0$, \ldots, $f_m\in\AX$ of degrees~$<d$, a consequence of the previous result is that in $\gA(X)$

\snic{\gen{f_0,\ldots,f_m}=\gen{f_0+X^df_1+\cdots+X^{dm}f_m}.}

\emph{3.}
By Exercise~\ref{exoMatLoc22}, $(x,y)$ admits a \mlp over~$\gB$ \ssi there exists a $B \in \MM_2(\gB)$ of trace $1$ satisfying $[\,x\; y\,]B = [\,0\; 0\,]$.
\\
So let $B \in \MM_2\big(\gA(X)\big)$ satisfy $[\,x\; y\,] B = [\,0\; 0\,]$ and $\Tr(B) = 1$.  \\
By multiplying the \coes of $B$ by a common \denoz, we obtain some \elts $p$, $q$, $r$, $s$ of $\AX$ such that \halfsmashbot{$[\,x\; y\,] \cmatrix {p &q\cr r &s\cr} = [\,0\; 0\,]$} and $p + s$ is primitive.  
We therefore have (with $p=\sum_k p_kX^{k}$, \dots):  $[\,x\; y\,] \cmatrix {p_i &q_i\cr r_i &s_i\cr} = [\,0\; 0\,]$.
As $p+s$ is primitive, we have $u_i \in \gA$ such that $\sum u_i(p_i+s_i) = 1$. Let $B' = \sum_i u_i\cmatrix{p_i&q_i\cr r_i &s_i}\in\MM_2(\gA)$: we obtain $[\,x\; y\,] B' = [\,0\; 0\,]$ with $\Tr(B')=1$.

\emph{4.} $\ArX$ is \ari $\Rightarrow$ $\gA$ is \ari and 

\snic{ \gA $ is \ari $ \,\iff\,\gA(X) $ is \ari  $\,\iff\,\gA(X)$ is a Bézout \riz$.}

The last \eqvc also results from the \plgref{thlgb3}.
In addition, the \mo of  \dve in $\gA(X)$, \cad $\gA(X)/\gA(X)\eti$, is \isoc to the \mo of  \itfs of $\gA$.

\exer{FermetureAlgZariski} 
We only show the first item. It is clear that $\gK'(\uX)$ is \agq over $\gK(\uX)$. \\
Conversely, let $z \in \gL(\uX)$ be \agq over $\gK(\uX)$, then there exists some nonzero $a \in \gK[\uX]$ such that $az$ is integral over $\gK[\uX]$, a fortiori over $\gL[\uX]$. As $\gL[\uX]$ is a GCD-domain, we have $az \in \gL[\uX]$. Moreover, we know that the \cli of $\gK[\uX]$ in~$\gL[\uX]$ is $\gK'[\uX]$ (Lemma~\ref{lemPolEnt}); so $az \in \gK'[\uX]$ then $z = (az)/a \in \gK'(\uX)$.


\exer{exoNotAbsIntClos} 
\emph {1.}
Immediate.

In what follows, we will use the fact that $(1,y)$ is a $\gk[x]$-basis of $\gk[x,y]$;
it is also a $\gk(x)$-basis of $\gk(x,y)$ and the extension $\gk(x,y)\sur{\gk(x)}$ is a Galois extension of the group $\gen {\sigma}$ where $\sigma : \gk(x,y) \to \gk(x,y)$ is the involutive $\gk(x)$-\auto which realizes $y \mapsto -y$.

\emph {2.}
Let $z = u(x) + yv(x) \in\gk(x,y)$ be \agq over $\gk$. \\
Then $z+\sigma(z) = 2u$ and $z\sigma(z) = u^2 - fv^2$ are \agq over $\gk$ and in $\gk(x)$ so in $\gk$. Hence $u \in \gk$, $v = 0$ and $z = u \in \gk$.

\emph {3.}
As $a \notin \gk^p$, we easily see that $f(X)$ is \ird in $\gk[X]$.
Let us show that $\gk[x,y]$ is the \cli of $\gk[x]$ in~$\gk(x,y)$. \\
Let $z = u(x) + yv(x) \in\gk(x,y)$ be integral over $\gk[x]$. \\
Then $z+\sigma(z) = 2u$ and $z\sigma(z) = u^2 - fv^2$ are in $\gk(x)$ and integral over $\gk[x]$, so in $\gk[x]$. 
Thus $u$ and $fv^2 \in \gk[x]$. By using the fact that $f$ is \irdz, we see that $v \in \gk[x]$. 
Recap: $z \in \gk[x,y]$.

\emph {4.}
Let $\alpha=a^{1/p}\in \gk$, hence $f(X) = (X-\alpha)^p$.
Let $t = y/(x-\alpha)^{p-1\over 2}$. \\
Then $t^2 = x-\alpha$,
therefore $x \in \gk[t]$, and $y = t(x-\alpha)^{p-1\over 2} = t^p \in \gk[t]$.\\
Therefore $\gk [x,y]  \subseteq  \gk [t] $ and $\gk(x,y) = \gk(t)$. We see that $t$ is integral over $\gk[x]$,
but %
that $t\notin\gk[x,y] = \gk[x] \oplus \gk[x]y$. The \cli of $\gk[x]$ (or that of $\gk[x,y]$) in~$\gk(x,y)$ is $\gk[t]$ (which indeed contains $x$ and $y$).


\exer{exoAnneauOuvertP1} Recall that $x_0=\fraC1 p$. The \egt

\snic{\gk[x_0, \ldots, x_{n-1}] = \sotq {u/p^s} {u \in \gk[t],\  \deg(u) \le ns}}

is easy by noticing that $t^n x_0 \in \gk[x_0, \ldots, x_{n-1}]$ since

\snic{{t^n \over p} = 1 + {t^n - p \over p} \;\in \;\sum_{i=0}^{n-1}\gk\,{ {t^i}  \over p}.}

Let us write that $t$ is \agq over $\gk(x_0)$ as a root in $T$ of the \pol

\snic{p(T)x_0 - 1=x_0T^n + x_0a_{n-1}T^{n-1} + \cdots+x_0 a_1T + (x_0a_0 - 1).}

The \elts determined by \gui{Emmanuel's trick} (see Lemma~\ref{lemEmmanuel} or Exercise~\ref{exoEntiersEmmanuel}) are 
\[\preskip-.2em \postskip.4em 
\begin{array}{rcl} 
  x_0t, &x_0t^2 + x_0a_{n-1}t,   & x_0t^3 + x_0a_{n-1}t^2 + x_0a_{n-2}t,  \\[1mm] 
&\ldots,&x_0t^{n-1} + \cdots + x_0a_2t . 
\end{array}
\]
Thus, $t^k x_0$ is integral over $\gk[x_0]$ for $k \in \lrb {0..n-1}$ and $\gk[x_0, \ldots, x_{n-1}] \subseteq \gA$.

It remains to show that $\gA \subseteq \gk[x_0, \ldots, x_{n-1}]$. We use the inclusion

\snic{\gk[x_0]\subseteq\gV_\infty := \gk[1/t]_{1+\gen{1/t}}.}

This last \ri is comprised of rational fractions of degree $\le 0$, \cad defined at $t=\infty$. It is \isoc to $\gk[y]_{1+\gen{y}}$ so it is \iclz, and $\gA\subseteq\gV_\infty$. 
The \ri $\gV_\infty$ is called \gui{the local \ri of the point $t = \infty$.}
\\
Let $z \in \gk(t)$ be an integral rational fraction over $\gk[x_0]$. By multiplying an \rdi of $z$ over $\gk[x_0]$ by $p^N$ with large enough $N$, we obtain

\snic {
p^N z^m + b_{m-1} z^{m-1} + \cdots + b_1 + b_0 = 0, \qquad b_i \in \gk[t]
.}

This entails that $p^N z$ is integral over $\gk[t]$ and therefore belongs to $\gk[t]$ ($\gk[t]$ is \iclz). Moreover, $z \in \gV_\infty$, \cad $\deg z \le 0$.  Ultimately, $z$ is a rational fraction of degree $\le 0$ whose \deno divides a power of $p$, %
so $z \in \gk[x_0, \ldots, x_{n-1}]$.

Finally, we have $x_1 = tx_0$ so $t = x_1/x_0 \in \Frac(\gA)$ then $\gk(t) = \Frac(\gA)$.


\exer{exoPresentationAlgOuvertP1} 
\emph {1.}
Immediate.

\emph {2.}
Let $\fb$ be the \id generated by each of the $X_iX_j - X_{i-1}X_{j+1}$ for $1 \le i \le j \le n-1$ and $E$ be the \kmo

\snic {
E = \gk[X_0] \oplus \gk[X_0]X_1 \oplus \cdots \oplus \gk[X_0]X_{n-1}
.}

We will prove that $E \cap \Ker\varphi = 0$ and that $\gk[\uX] = E + \fb$. As $\fb \subseteq \Ker\varphi$, we will obtain $\gk[\uX] = E \oplus \fb$. Let $y \in \Ker\varphi$ which we write as $y = y_1 + y_2$ with $y_1 \in E$ %
and $y_2 \in \fb$.
By applying $\varphi$, we obtain $\varphi(y_1) = 0$, so $y_1 = 0$, then $y = y_2 \in \fb$. We have obtained $\Ker\varphi \subseteq \fb$, hence the \egt $\Ker\varphi = \fb$.

$\bullet$ \emph{Justification of $E \cap \Ker\varphi = 0$}. Let $f \in E$

\snic {
f = f_0(X_0) + f_1(X_0)X_1 + \cdots + f_{n-1}(X_0)X_{n-1}
.}

Let $\varphi(f) = 0$,

\snic {
\varphi(f) = f_0(1/p) + f_1(1/p)t/p + \cdots + f_{n-1}(1/p)t^{n-1}/p = 0
.}

By multiplying each $f_i(1/p)$ by $p^N$, with large enough $N$, we obtain $g_i(p)\in \gk[p]$,

\snic {
pg_0(p) + g_1(p)t + \cdots + g_{n-1}(p)t^{n-1} = 0
.}


But $(1, t, \ldots, t^{n-1})$ is a basis of $\gk[t]$ over $\gk[p]$, so the $g_{k}$'s $ = 0$, then $f = 0$.

$\bullet$
\emph{Justification of $\gk[\uX] = E + \fb$}. \\
Letting $\gk[x_0, \ldots, x_{n-1}] = \gk[\uX]/\fb$ and $E' = \gk[x_0] + \gk[x_0]x_1 + \cdots + \gk[x_0]x_{n-1}$,
this amounts to showing that $\gk[\ux] = E'$. Moreover $E'$ contains $x_n = 1 - \sum_{i=0}^{n-1} a_ix_i$. It therefore suffices to prove that $E'$ is a sub\riz, or  \hbox{that $x_ix_j \in E'$} \hbox{for $i$, $j \in \lrb {0..n-1}$}.  By \dfnz, it contains $x_0^2$, $x_0x_1$, \dots, $x_0x_{n-1}$ and therefore  also contains~$x_0x_n$. But $x_1x_j = x_0x_{j+1}$ for $j \in \lrb {1..n-1}$, so $E'$ contains these~$x_1x_j$'s and therefore also contains~$x_1x_n$,  and by using $x_2x_j = x_1x_{j+1}$ \hbox{for $j \in \lrb {2..n-1}$}, we see that~$E'$ contains all the~$x_2x_j$'s. And so forth.

\rem
The author of the exercise proceeded as follows, for some \cdiz~$\gk$: he used an additional \idtr 
 $X_n$ and chose, for $\gk[X_0, X_1, \ldots, X_n]$, the graded monomial reversed lexicographical order by ordering the \idtrs as follows: $X_0 < X_1 < \cdots < X_n$.  We then observe that
the trivial \id of the \id $\gen {R_{\rm min}} + \gen {1 - \sum_{i=0}^n a_iX^i}$ is the monomial \id generated by the monomials
$$
\preskip.4em \postskip.4em 
X_n  \hbox { and } X_iX_j  \hbox { for }  1 \le i \le j \le n-1
.\leqno (\star)
$$
The \kev generated by the monomials wich are  nondivisible by a monomial of~$(\star)$ is the \kev $E = \gk[X_0] \oplus \gk[X_0]X_1 \oplus \cdots \oplus \gk[X_0]X_{n-1}$.
It is the space that appears in the above solution (in which $\gk$ is an arbitrary \riz, not \ncrt a \cdiz).
\eoe


\exer{exoEntiersEmmanuel}
By multiplying the initial \eqn by $a_n^{n-1}$, we obtain $a_n s$ integral over~$\gA$. Let us then express the initial \eqn as follows

\snic {
(a_n s + a_{n-1}) s^{n-1} + a_{n-2} s^{n-2} + \cdots + a_1 s + a_0 = 0,
\, \hbox { with } b = b_{n-1}= a_n s + a_{n-1}
,}

and let us consider the \ri $\gA[b]$. Thus,~$s$ annihilates a \pol of~$\gA[b][X]$ whose leading \coe is~$b$; by what precedes, $bs$ is integral over~$\gA[b]$. But~$b$ is integral over~$\gA$ so $bs = a_n s^2 + a_{n-1} s$ is integral over~$\gA$.  
\\
The following step consists in writing the initial \eqn in the form

\snic {\mathrigid 2.5mu
cs^{n-2} + a_{n-3} s^{n-3} + \cdots + a_1 s + a_0 = 0, \; \hbox { with } c=  b_{n-2}= a_n s^2 + a_{n-1} s + a_{n-2}.
}


\exer{exoKroneckerTheorem}
\emph {1.}
We write $\lrb{1..n} \setminus I = \{i_1, i_2, \ldots\}$. By using Lemma~\ref{lemEmmanuel}, we see that the \coes of $h_1(T) = h(T)/(T - x_{i_1})$ are integral over $\gA$, that those of $h_2(T) = h_1(T)/(T- x_{i_2})$ are integral over $\gA[{\rm coeffs.\ of\ } h_1]$, therefore integral over~$\gA$ and so on and so forth. Therefore by letting $q(T) = \prod_{i'\notin I}(T - x_{i'}) \prod_{j'\notin J}(T - y_{j'})$, the \coes of the \pol $h(T)/q(T)$ are integral over $\gA$. The constant \coe of this last \pol is $\pm a_0b_0 \prod_{i\in I}x_i \prod_{j\in J}y_j$.

\emph {2.}
\Elr \smq functions: we have $a_i = \pm a_0 S_i(\ux)$, $b_j = \pm b_0 S_j(\uy)$, so~$a_ib_j$ is integral over $\gA$.

\exer{exoAnneauCoinceBézout}
Let $S \subseteq \gA\setminus \{0\}$ be the set of  \denos $b$ of the \elts of~$\gB$ written in the form $a/b$ with $a, b \in \gA$, $b\ne 0$ and $1 \in \gen {a,b}$. It is clearly a \moz. To show that $\gB = \gA_S$, it suffices to prove that $S^{-1} \subseteq \gB$.\\
Let $a/b \in \gB$ be expressed irreducibly; there exist $u, v \in \gA$ such that $1 = ua + vb$ which implies $1/b = u(a/b) + v \in \gA\gB + \gA \subseteq \gB$.

\exer{exoGrellNoether}
We want to show that an intermediary \ri  between~$\gA$ and~$\gK$ is a \adpz.
Every \elt of $\gK$ is primitively \agq over any arbitrary intermediary \ri between $\gA$ and $\gK$. It remains to
prove that the intermediary \ri is \icl in order to apply \thref{th.2adpcoh}.

\emph {1.}
If $x = a/b$, with $a$, $b \in \gA$, there exists a  \mlp 
for~$(b,a)$, $\cmatrix {s & c\cr t & 1-s\cr} 
\in \MM_2(\gA)$,  with 
$sa = cb$ and $ta = (1-s)b$.\\
Therefore $x = c/s = (1-s)/t$ and $x \in \gA'_s \cap \gA'_t$. 
Conversely, if $x' \in \gA'_s \cap \gA'_t$, there \hbox{is $a'$, $b' \in \gA'$} \hbox{and $n$, $m \in \NN$} such that $x' = a'/s^n = b'/t^m$. Therefore, for $u$, $v \in \gA$, since $1/t = x/(1-s)$ we have
$$\preskip.2em \postskip.4em
x' = \frac {a'}{s^n} = \frac{b' x^m}{(1-s)^m} =
\frac {ua' + vb' x^m}{us^n + v(1-s)^m}\;.
$$
It suffices to take $us^n + v(1-s)^m = 1$ to observe that $x' \in \gA'[x]$.

\emph {2.}
Let $\gB \subseteq \gK$ be an \Alg generated by $n$ \elts $(n\geq1)$. \\
We write $\gB = \gA'[x]$, where $\gA'$ is an \Alg generated by $n-1$ \eltsz. By item~\emph{1}, there exist $s$, $t \in \gA$ such that $\gA'[x] = \gA'_s \cap \gA'_t$.\\
By \recuz, there exist $u_1$, \ldots, $u_k \in \gA$ such that $\gA' = \gA_{u_1} \cap \cdots \cap \gA_{u_k}$. \\
Then,
$\gA'_s = \gA_{su_1} \cap \cdots \cap \gA_{su_k}$ and 
$\gA'_t = \gA_{tu_1} \cap \cdots \cap \gA_{tu_k}$, so

\snic {
\gB = \gA_{su_1} \cap \cdots \cap \gA_{su_k} \cap
\gA_{tu_1} \cap \cdots \cap \gA_{tu_k}.
}

\emph {3.}
Let $\gB$ be an intermediary \ri and $x \in \gK$ be integral over $\gB$.
Then $x$ is integral over a \tf \Aslgz, therefore it belongs to this \tf \Aslgz, therefore to $\gB$, \cad $\gB$ is \iclz.

\emph {4.}
Let $x,y$ be two \idtrs over a \cdi $\gk$ and $\gA = \gk[x,y]$.\\
Let $\gB =\gk[x,y,(x^2+y^2)/xy]$. Then $\gA$ is \icl but not $\gB$: indeed, $x/y$ and $y/x$ are integral over $\gB$ (their sum and their product are members of~$\gB$) but $x/y$ and $y/x \notin\gB$ as we can easily prove, thanks to a homogeneity argument.


\exer{exoPrimitivementAlg} 
We have $bx - a = 0$ with $1 = ua + vb$.
The reader will check that \hbox{if $f(Y) \in \gA[Y]$} satisfies $f(y) = 0$, then $f$ is a multiple, in $\gA[Y]$, 
of $bY - 2a$. Therefore $\rc(f) \subseteq \gen {2a, b}$ and as $1 \notin \gen {2a, b}$, $y$ is not primitively \agq over $\gA$.


\exer{exocaracPruferC} 
The implications \emph {4} $\Rightarrow $ \emph {3} $\Rightarrow $ \emph {2}
and \emph {5} $\Rightarrow $ \emph {2} are trivial. \Thref{thSurAdp} gives \emph {1} $\Rightarrow $ \emph {4} 
and \thref{thPruf}~\emph{4d} (\paref{thPruf}) gives \emph {1}~$\Rightarrow$~\emph{5.}

  \emph
{2} $\Rightarrow $ \emph {1.} 
$x$ is  primitively \agq over $\gA$, we apply \thref{th.2adpcoh}.

\exer{exocaracPruCoh}
We already know that \emph{1} $\Rightarrow$ \emph{2} $\Rightarrow$ \emph{3}
and
\emph{1} $\Rightarrow$ \emph{5.} 

Let us show that \emph{3} implies that the \ri is \ariz.
Consider an \id with two arbitrary \gtrs $\fa=\gen{y_1,y_2}$ and let $r_i$ be the \idm annihilator of $y_i$.
Consider the \ort \idmsz:  $e=(1-r_1)(1-r_2)$,  $f=r_1(1-r_2)$, \hbox{and $g=r_2$}. We have $e+f+g=1$.
If we invert $f$ or $g$, one of the $y_i$'s is null and \hbox{the \id $\fa$} becomes  principal.
To see what happens if we invert $e$, consider the \ndz \elts $x_1=(1-e)+ey_1$ \hbox{and $x_2=(1-e)+ey_2$}. The \id $\fb=\gen{x_1,x_2}$ is \iv in $\gA$.
Then let $u$, $v$, $w$ be such \hbox{that $ux_1=vx_2$} \hbox{and $(1-u)x_2=wx_1$}. We multiply by~$e$ and we obtain $uey_1=vey_2$ \hbox{and $(1-u)ey_2=wey_1$}, which implies that the \id $\fa\gA_e=\gen{ey_1,ey_2}\gA_e$ is \lopz.

\emph{5} $\Rightarrow$ \emph{4.} First consider $f=aX+b$, $g=aX-b$, then $f=aX+b$, $g=bX+a$.

\emph{4} $\Rightarrow$ \emph{3.} Let $\fa=\gen{a,b}$, with \ndz $a$ and $b$.  Let $\alpha$, $\beta$ such that $ab=\alpha a^2+\beta b^2,$ and let $\fb=\gen{\alpha a,\beta b}$.  We have $ab \in \fa\fb$, therefore

\snic{
a^2b^2 \in\fa^2\fb^2= \gen{a^2,b^2}\gen{\alpha^2a^2,\beta^2b^2}.
}

Let us show the \egt $\gen{a^2b^2} = \fa^2\fb^2$, which will imply that $\fa$ is \ivz.  \\
Letting $u=\alpha a^2,$ $v=\beta b^2$, it suffices to show \hbox{that $u^2=\alpha^2a^4$} \hbox{and $v^2=\beta^2b^4$} are in $\gen{a^2b^2}$.  By \dfnz, $u+v = ab \in \fa\fb$ and $uv \in \gen {a^2b^2}$.
\\
Therefore $u^2 + v^2 = (u+v)^2 - 2uv \in \gen {a^2b^2}$. As $u^2$, $v^2 \in\gen{u^2+v^2,uv}$, we indeed have $u^2$, $v^2 \in \gen {a^2b^2}$.


\exer{exoAnarlgb} We give the \dem for the integral case. The \qiri case is deduced from it by applying the usual \elgbmdz.
\\
\emph{1.} Let $M\in \Ae{n\times m}$, $p=\inf(m,n)$. 
Proposition~\ref{propIddsAnar} gives us \lop \ids $\fa_i$ such that

\snic{
\cD_{\gA,1}(M)=\fa_1  ,\; \cD_{\gA,2}(M)=\fa_1^2\fa_2,\; \cD_{\gA,3}(M)=\fa_1^3\fa_2^2\fa_3,\; \cD_{\gA,4}(M)=\fa_1^4\fa_2^3\fa_3^2\fa_4,\;\ldots }

Since the \ri is \lgbz, the \lop \ids $\fa_j$ are principal 
(\plgref{thlgb3}). 
\\
Let $\fa_j=\gen{a_j}$ and consider the matrix $M'\in \Ae{n\times m}$ in Smith form, whose diagonal \elts are
$a_1$, $a_1a_2$, \ldots, $a_1a_2\cdots a_p$.
\\
As in the \dem of Proposition~\ref{propIddsAnar} the \algo that produces the reduced Smith form in the local case and the \lgbe machinery of \anars provides us with a comaximal \sys $(s_1,\ldots,s_r)$ such that, over each $\gA[1/s_i]$, the matrix $M$ admits a reduced Smith form.\imla 
By comparing the \idds we see that this reduced form can always be taken equal to $M'$ (here is where the fact that over an integral \riz, two \gtrs of a \idp are always associated intervenes).
\\
Thus, $M$ and $M'$ are \eqves over each $\gA[1/s_i]$.
\Trf by the \plgref{thlgb1} that they are \eqvesz.

\emph{2.} Immediate consequence of \emph{1.}

\entrenous{Cela semble probable que cela marche for every \anar \lgbz. 
}

\exer{exoReductionIdeal}
\emph{1.} We write $E = \gA x_1 + \cdots + \gA x_n$ therefore $\fa E = \fa x_1 + \cdots
+ \fa x_n$. By using $bE \subseteq \fa E$, we obtain a matrix $A \in \Mn(\fa)$ such that 

\snic{b \tra { [\,x_1\; \cdots\; x_n\,] } = A\tra {[\,x_1\; \cdots\; x_n\,]}.}

\snii
It then suffices to let $d = \det(b\In - A)$.

 \emph{2.}
If $\deg(g) \le m$, we know that $\rc(f)^{m+1}\rc(g) = \rc(f)^m \rc(fg)$
(Lemma~\ref{lemdArtin}).
By multiplying by $\rc(g)^m$, we obtain $\big(\rc(f)\rc(g)\big)^{m+1} = \rc(fg)\big(\rc(f)\rc(g)\big)^m$.

 \emph{3.}
We have $\fb^2 = \fa\fb$, $\fb'^5 = \fa_1\fb'^4$ and $\fb'^4 = \fa_2\fb'^3$.

 \emph{4.}
Suppose $\fb^{r+1} = \fa\fb^r$. We apply the first question with $E = \fb^r$ and $b \in \fb$. We obtain $d = b^n + a_{1} b^{n-1} + \cdots + a_{n-1} b + a_n \in \Ann(\fb^r)$ with $a_i \in \fa^i$.
\\
As $d \in \fb$ and $d \in \Ann(\fb^r)$, we have $d^{r+1} = 0$, which is an \rdi of $b$ over $\fa$.\\
For the converse, let $\fb$ be integral over $\fa$. For $b \in \fb$, by writing an \rdi of $b$ over $\fa$. We obtain $n$ such that $b^{n+1} \in \fa \fb^n$. However, if we have two \ids $\fb_1, \fb_2 \subseteq \fb$ with $\fb_i^{n_i+1} \subseteq \fa\fb^{n_i}$, we have $(\fb_1 + \fb_2)^{n_1+n_2+1} \subseteq \fa\fb^{n_1 + n_2}$.
\\
By using a finite \sgr of $\fb$, we obtain an exponent $r$ with the inclusion~\hbox{$\fb^{r+1} \subseteq \fa\fb^r$}%
.


\exer{exolemNormalIcl}
Let $\gK=\Frac\gA$.

\emph{1.} Let $a\in\gA$ and $e_a$ be the \idm of $\gK$ such that $\Ann_\gK(a)=\Ann_\gK(e_a)$. The \eltz~$e_a$ is integral over~$\gA$, so $e_a\in\gA$,
and $\Ann_\gA(b)=\Ann_\gK(b)\cap\gA$ for every $b\in\gA$.

\emph{2. Direct implication.}  The computation is \imdz.

\emph{2. Converse implication.} Let $a$ be integral over the \idp $\gen{b}$ in $\gA$.
Let us express the \rdi of $a$ over $\gen{b}$.
$$\preskip.4em \postskip.4em
 a^n=b(u_{n-1}a^{n-1}+u_{n-2}ba^{n-2}+\cdots+u_0b^{n-1}).\eqno(*)
$$ 
We have $(1-e_b)a^n=0$, therefore since $\gA$ is reduced $(1-e_b)a=0$.
We introduce the \ndz \elt $b_1=b+(1-e_b)$. Then the \elt $c=a/b_1\in\gK$ is integral over~$\gA$. Indeed, the \egt $(*)$ remains true when replacing $b$ by~$b_1$ and the~$u_i$'s \hbox{by $e_bu_i$'s}, because in the component $e_b=1$ we obtain $(*)$ and over the component~\hbox{$e_b=0$} we obtain $0=0$.
\\
Therefore $c$ is in $\gA$, and $a=e_ba=e_bb_1c=bc$. 

\exer{exoEntierSurIX}
Let $T$ be a new \idtr over $\gB$. For $b \in\gB$, we will use the result (similar to Fact~\ref{fact2Entiers}): $b$ is integral over the \id $\fa$ \ssi $bT$ is integral over the sub\riz
%
$\gA[\fa T] \eqdefi \gA\oplus \fa T \oplus \fa^2 T^2 \oplus \dots$
of $\gB[T]$.

\snii
Let us take a look at the difficult case. Let $F\in\gB[X]$ be integral over $\fa[X]$, we must show that each \coe of $F$ is integral over $\fa$. We write an \rdi

\snic {
F^n + G_1F^{n-1} + \cdots + G_{n-1}F + G_n = 0, \quad
G_k=G_k(X) \in  (\fa[X])^k = \fa^k[X].
}

\snii
We therefore have an \egt in $\gB[X][T]$ with some $Q_i$'s in $\gB[X]$

\snic {\mathrigid 2.5mu
T^n + G_1T^{n-1} + \cdots + G_{n-1}T + G_n = 
(T - F) (T^{n-1} + Q_1T^{n-2} + \cdots + Q_{n-1}).
}

\snii
We replace $T$ by $1/(TX)$ and we multiply by $(TX)^n = TX \times (TX)^{n-1}$, which gives

\snic {\mathrigid 2.5mu
1 + XTG_1 + \cdots + X^nT^nG_n = (1 - XTF) 
(1 + XTQ_1 + \cdots + X^{n-1}T^{n-1}Q_{n-1}).
}

\snii
We now look at this \egt in $\gB[T][X]$. \\
If $b$ is a \coe of $F$,  $bT$ is a \coe in $X$ of $1 - XTF$ and $1$ is a \coe in $X$ of $1 + XTQ_1 + \cdots + X^{n-1}T^{n-1}Q_{n-1}$. By Kronecker's \thoz, the product $bT = bT \times 1$ is integral over the \ri generated by the \coes (in $X$) of the \pol $1 + XTG_1 + \cdots + X^nT^nG_n$. But the \coe in $X^k$ of this last \pol is in $\gA[\fa T] = \gA \oplus \fa T \oplus \fa^2 T^2 \oplus \dots$ and therefore $bT$ is integral over $\gA[\fa T]$ and consequently $b$ is integral over $\fa$.


\exer{exosdirindec} \emph{(Indecomposable modules)}\\
\emph{1.} Everything takes place modulo $\fa$. We therefore consider the quotient \ri $\gB=\gA/\fa$. Then the result is obvious (Lemma~\ref{lemfacile}).

\emph{2a.} If $M=N\oplus P$, $N$ and $P$ are \prcs and the sum of the ranks is equal to $1$, therefore one of the two is null.

\emph{2b.} We refer to item~\emph{1.} If the module is decomposable, we have $\fa\subseteq \fb$ and $\fc$ with~$\fb$ and~$\fc$ \tf \comz. These \ids are therefore obtained from the \fac of $\fa$ as two partial products of this \fcnz.
\\
Thus, we cannot have $\fb$ and $\fc$ \com if the \fac of $\fa$ makes only one \idema intervene. 
\\
Otherwise the \fac of $\fa$ provides two \com \idsz~$\fb$ and $\fc$ such \hbox{that $\fb\fc=\fa$}. Therefore $\fb+\fc=\gZ$ and $\fb\cap\fc=\fa$ which gives $\gZ/\fa=\fb/\fa\oplus \fc/\fa$.
\\
Actually, if $\fa=\prod_{i=1}^{k}\fq_i=\prod_{i=1}^{k}\fp_i^{m_i}$ is the \fac of $\fa$, we obtain by \recu on $k$ that
  $\gZ/\fa=\bigoplus_{i=1}^{k}\fq_i/\fa$.

\emph{2c.} Results from the previous considerations and from the structure \tho for \mpfs over a \dDkz.

\emph{3.} The uniqueness can be stated as follows: if $M$ can be expressed in two ways as a sum of indecomposable modules, there is an \auto of $M$ which sends the modules of the first \dcn to those of the second. 
\\
If a torsion \mpf $M$ is decomposed into direct sums of indecomposable modules, each term of the sum is itself \pf and with torsion. It is therefore of the form $\gZ/\fp^{m}$ by item \emph{1.} 
\\
By the Chinese remainder \tho we return to the case where only one \idema intervenes in the direct sum, and the uniqueness then results from \thref{prop unicyc}.
\\
Note also that in the case of a \fac PID, the uniqueness is valid for the \dcn of every \mpfz. 



\prob{exoArithInvariantRing}  
Hereinafter the word \gui{\lotz} means \gui{after \lon at \ecoz.}

\emph {1.}
The \id $\fa$ is \lopz, therefore since $\gA$ is normal, \lot \iclz, so it is \icl (\plgref{plcc.normal}). We end with Lemma~\ref{lemthKroicl} (variant of \KROz's \thoz).

\emph {2.}
If $x \in \fa\gB \cap \gA$, then $x$ is integral over the \id $\fa$ (Lying Over, Lemma~\ref{lemLingOver2}) therefore in~$\fa$ by the previous question.

\emph {3a.}
If $a = \rN_G(b)$, we have $\rN_G(b\gB) = a\gB \cap \gA = a\gA$.

\emph {3b} and \emph{3c.}
The \itf $\fa = \rc_\gA(h)$ is \lopz, so $\rc_\gB(h) = \fa\gB$ is a \lop \id of $\gB$. By the first question, we have

\snic {
\prod_\sigma \rc_\gB(h_\sigma) = \rc_\gB(h),
\quad \hbox {i.e.} \quad \rN'_G(\fb) = \fa\gB.
}

By question \emph {2}, $\fa =  \rN_G(\fb)$.
Next we note that

\snic {
\rN_G(\fb_1\fb_2)\gB = \rN'_G(\fb_1\fb_2) = \rN'_G(\fb_1)\rN'_G(\fb_2)
= \rN_G(\fb_1)\rN_G(\fb_2)\gB,
}

hence the result when taking the intersection with $\gA$.

\emph{3d.} This results from item \emph{2} and from the two following facts.
\\
$\bullet$ If $b\in\gB$ is \ndz then $a = \rN_G(b)\in\gA$ is \ndz in $\gA$:
indeed, it is a product of \ndz \elts in $\gB$ therefore it is \ndz in $\gB$.
\\
$\bullet$  If $a \in \gA$ is \ndz in $\gA$ then it is \ndz in $\gB$. Indeed, let $x \in \gB$ such that $ax = 0$. We want to show that $x=0$. We consider the \pol

\snic{\rC{G}(x)(T)  =    \prod_{\sigma\in G}\big(T-\sigma(x)\big).}

As $a\sigma(x)=0$ for each $\sigma$, the \coes of $\rC{G}(x)(T)$ are annihilated by $a$ therefore null, except for the leading \coez. 
Thus $x^{\abs{G}}=0$, but $\gB$ is normal therefore reduced.

\emph {4.}
Let $\gk(x,y) = \Frac \gk[x,y]$.  We will use the fact that $(1,y)$ is a $\gk[x]$-basis of $\gk[x,y]$; it is also a $\gk(x)$-basis of $\gk(x,y)$ and the extension $\gk(x,y)\sur{\gk(x)}$ is a Galois extension of the group $\gen {\sigma}$ where $\sigma : \gk(x,y) \to \gk(x,y)$ is the involutive $\gk(x)$-\auto which realizes $y \mapsto -y$.
Let us show that $\gk[x,y]$ is the \cli of $\gk[x]$ in $\gk(x,y)$. Let $z = u(x) + yv(x) \in\gk(x,y)$ be integral over $\gk[x]$. Then $z+\sigma(z) = 2u$ and $z\sigma(z) = u^2 - fv^2$ are in $\gk(x)$ and integral over $\gk[x]$ therefore in $\gk[x]$. We therefore have \hbox{$fv^2 \in \gk[x]$}. By using the fact that $f$ is \spbz, we see that $v \in \gk[x]$. Recap: $z \in \gk[x,y]$. Therefore $\gk[x,y]$ is \iclz.
We apply the preceding with $\gA = \gk[x]$, $\gB = \gk[x,y]$, $G = \gen {\sigma}$.

\prob{exoFullAffineMonoid}~\\ 
\emph {2a)}
Let $a = \sum_\alpha a_\alpha \ux^\alpha$, $b = \sum_\beta b_\beta
\ux^\beta$. 
\\
We must show that $\beta \in M$ for each $\beta$ such that $b_\beta \ne 0$. We can assume~$b$ nonzero. Let $a_\alpha\ux^\alpha$ be the leading monomial of $a$ for the lexicographical order \hbox{and $b_\beta\ux^\beta$} be that of $b$. The leading monomial of $ab$ is $a_\alpha b_\beta\ux^{\alpha+\beta}$, therefore $\alpha+\beta \in M$.
\\
As $\alpha\in M$ and as $M$ is full, we have $\beta\in M$. We then start again by replacing $b$ by $b' = b-b_\beta\ux^\beta$ which satisfies $ab' \in \gk[\ux]$.
We obtain $b' \in \gk[\ux]$ and finally $b \in \gk[\ux]$.


\prob{exoBaseNormaleAlInfini}~\\
\emph {1.}
If $A = (a_{ij})$, then $\det A = \sum_{\sigma \in \rS_n} a_{\sigma(1)1} \cdots
a_{\sigma(n)n}$ and

\snic {
v(a_{\sigma(1)1} \cdots a_{\sigma(n)n}) \ge v(A_1) + \cdots + v(A_n)
.}

We deduce that $v(\det A) \ge v(A_1) + \cdots + v(A_n)$.

\emph {2.}
For the matrix given as an example: we have $\det(A) = \pi^2 - \pi \ne 0$.
\\
But
$\ov A = \cmatrix {0 & 0\cr 1 & 1\cr}$ is not invertible. By realizing $A_1 \leftarrow A_1 - A_2$, we obtain the \egt $A' = \cmatrix {\pi^2 - \pi &\pi \cr 0 & 1\cr}$ and this time $\ov {A'} = \cmatrix {1 & 0\cr 0 & 1\cr}$ is invertible.

\label{allnonzero}
Here is the \gnl method: if $\det \ov A \ne 0$, $A$ is $\gA_\infty$-reduced and that is all. Otherwise, there are some $\lambda_1, \ldots, \lambda_n \in \gk$, not all zero, such that $\lambda_1 \ov {A_1} + \cdots + \lambda_n \ov {A_n} = 0$. We consider a column $A_j$ with $\lambda_j \ne 0$ and $v(A_j)$ minimum; to simplify, we can suppose that it is $A_1$ and that $\lambda_1 = 1$ (even if it entails dividing the relation by $\lambda_1$); we then perform the \elr \op

\snic {
A_1 \leftarrow A'_1 = A_1 + \sum_{j=2}^n \lambda_j \pi^{v(A_1)-v(A_j)} A_j
.}

In this sum, by only making the $A_j$'s for which $\lambda_j \ne 0$ intervene, each exponent of $\pi$ is $\ge
0$. This is therefore a $\gk[\pi^{-1}]$-\elr \op on the columns, \cad $\gk[t]$-\elrz, and we do not change the $\gk[t]$-module generated by the columns. Moreover, $v(A'_1) > v(A_1)$; indeed, (by remembering that $\lambda_1 = 1$):

\snic {
A'_1 / \pi^{v(A_1)} = s \eqdf {\rm def} 
\som_{\lambda_j \ne 0} \lambda_j A_j / \pi^{v(A_j)}
,}

and $v(s) > 0$ since by hypothesis $\sum_{\lambda_j \ne 0} \lambda_j \ov {A_j}= 0$. 
We iterate this process which eventually stops because at each step, the  sum
 $\sum_j v(A_j)$ strictly increases while being bounded above by $v(\det A)$, invariant under the above operations.

\emph {3.}
Let $y = Px$, \cad $y_i = \sum_j p_{ij} x_j$; we have $v(p_{ij}) \ge 0$,
$v(x_j) \ge v(x)$ so $v(y_i) \ge v(x)$ then $v(y) \ge v(x)$.
By symmetry, $v(y) = v(x)$. The remainder poses no more difficulties.

\emph {4.}
$A$ is $\gA_\infty$-reduced \ssi every (\ncrt nonzero) diagonal coefficient divides (in the $\gA_\infty$ sense) all the coefficients of its column.

\emph {5.}
Even if it entails replacing $A$ by $AQ$ with suitable $Q \in \GL_n(\gA)$, we can suppose that $A$ is $\gA_\infty$-reduced. We will realize some \ops $A \leftarrow PA$ with $P \in \GL_n(\gA_\infty)$ (\cad consider the $\gA_\infty$-lattice generated by the rows of $A$), which does not modify the $\gA_\infty$-reduced \crc of $A$.  There exists a $P \in \GL_n(\gA_\infty)$ such that $PA$ is upper triangular and we replace $A$ by $PA$.  Let $L_1, \ldots, L_n$ be the rows of $A$; we then realize the $\gA_\infty$-\elr \op 

\snic {
L_1 \leftarrow L_1 - {a_{12} \over a_{22}} L_2
\qquad  \hbox {recall : }  a_{22} \divi_{\gA_\infty} a_{12}
,}

which brings a $0$ in position $a_{12}$ (and the new matrix is always triangular \hbox{and $\gA_\infty$-reduced)}. We continue in order to annihilate all the coefficients of the first row (except $a_{11}$); we can then pass to the second row and so on and so forth in order to obtain a diagonal matrix (by constantly using the fact that in an $\gA_\infty$-reduced triangular matrix, each diagonal coefficient $\gA_\infty$-divides all the coefficients of its column). As $\gA_\infty$
is a DVR, we can make sure that the final obtained diagonal matrix is $\Diag(\pi^{d_1}, \ldots, \pi^{d_n})$ with $d_i \in \ZZ$.

\emph {6a.}
Let $\und\varepsilon$ be an $\gA$-basis of $E$, $\und{\varepsilon'}$ be an $\gA_\infty$-basis of $E'$  and $A = \Mat_{\und\varepsilon, \und{\varepsilon'}}(\Id_L)$. Then there exist $P \in \GL_n(\gA_\infty)$ and $Q \in \GL_n(\gA)$ such that {\mathrigid 2mu $PAQ = \Diag(t^{-d_1}, \ldots, t^{-d_n})$}.  Let
$\ue$ and $\und{e'}$ be defined by  $\Mat_{\ue,\und\varepsilon}(\Id_L) = Q$, $\Mat_{\und{\varepsilon'},\und{e'}}(\Id_L) = P$. 
\\
Then $\ue$ is an $\gA$-basis of $E$, $\und{e'}$  an $\gA_\infty$-basis of $E'$ and $e_i = t^{-d_i} e'_i$.

\emph {6b.}
Since $t^je_i = t^{j-d_i} e'_i$, it is clear that $t^j e_i \in E \cap E'$ for $0 \le j \le d_i$. 
Conversely, let $y \in E \cap E'$ which we express as 

\snic {
y = \sum_i a_ie_i = \sum_i a'_i t^{d_i} e_i, \quad\hbox{with }
a_i \in \gA \hbox{ and } a'_i \in \gA_\infty
,}

and therefore $a_i = a'_i t^{d_i}$. \\
If $d_i < 0$, we obtain $a_i = a'_i = 0$, and if $a_i \ne 0$,  $0 \le \deg a_i \le d_i$. Hence the stated $\gk$-basis.

\emph {7.}
First of all $\gk' = \gB \cap \gB_\infty$, so $\gB$ and $\gB_\infty$ are $\gk'$-\evcsz. Let us show that each $r_i \in \gA_\infty$ and that in addition, if $e_i \notin \gB_\infty$, then  $v(r_i) > 0$, \hbox{\cad $\deg(r_i) < 0$}. 
If $e_i \in \gB_\infty$, we have $e_i \in \gB \cap\gB_\infty = \gk'$, so also $e_i^{-1} \in \gk'$; consequently $r_i = e_i^{-1} (r_ie_i) \in \gB_\infty$ therefore $r_i \in \gB_\infty \cap \gK = \gA_\infty$.
\\
If $e_i \notin \gB_\infty$, we write $e_i = r_i^{-1}(r_ie_i)$, an \egt which proves that $r_i^{-1} \notin \gA_\infty$ (let us not forget that $r_ie_i \in \gB_\infty$) so $v(r_i^{-1}) < 0$, \cad $v(r_i) > 0$.

Now let $c \in \gk'$ which we express in the $\gA$-basis $(e_i)$ and the $\gA_\infty$-basis $(r_ie_i)$

\snic {
c = \sum_i a_ie_i = \sum_i a'_i r_ie_i, \quad a_i \in \gA, \quad
a'_i \in \gA_\infty, \quad a_i = a'_i r_i
.}

For the $i$'s such that $e_i \in \gk'$, as $r_i \in \gA_\infty$, we have $a_i = a'_ir_i \in \gA \cap \gA_\infty = \gk$. It remains to see that for $e_i \notin \gk'$, $a_i = 0$; the \egt $a_i = a'_ir_i$ and the fact that $a_i \in \gA$, $a'_i \in \gA_\infty$ and $\deg(r_i) < 0$ then entail $a_i = a'_i = 0$.
Recap: the $e_i$'s which are in $\gk'$ form a $\gk$-basis of $\gk'$.

\emph {8.}
By letting $i = y/x$, we have $i^2 = -1$ and

\snic {
\cmatrix {1 & x\cr 0 & 1\cr} \cmatrix {1\cr i} = \cmatrix {y+1\cr i}
.}

The matrix on the left-hand side has determinant $1$, therefore $(1,i)$ and $(y+1,i)$ are two bases of the same $\gA$-module. But $y+1$ is not integral over $\gA_\infty$ (because $x$ is integral \hbox{over $\gk[y] = \gk[y+1]$} and is not integral over $\gA_\infty$).
The basis $(1, i)$ is normal at infinity but not the basis $(y+1,i)$.

\prob{exoHyperEllipticFunctionRing} 
\emph {1.}
Let $z = y-v$, $(1,z)$ is an $\gA$-basis of $\gB$ and $\gA u \cap \gA z = \{0\}$. To show that $\fb_{u,v} = \gA u \oplus \gA z$ is an \idz, it suffices to see that $z^{2}\in\fb_{u,v}$. \\
However, $y^2 = (z+v)^2 = z^2 + 2vz + v^2$, \cad $z^2 + 2vz + uw = 0$.

\emph {2.}
As $(1,z)$ is an $\gA$-basis of $\gB$ and $(u, z)$ is an $\gA$-basis of $\fb_{u,v}$, we obtain the \egt $\gA \cap \fb_{u,v} = u\gA$. Moreover, every \elt of $\gB$ is congruent modulo $z$ to an \elt of $\gA$, therefore $\gA \to \gB\sur{\fb_{u,v}}$ is surjective with kernel $u\gA$. \\
The matrix $M$ of $(u, y-v)$ over $(1, y)$ is $M = \cmatrix {u & -v\cr 0 & 1}$ with $\det(M) = u$,
which gives $\rN(\fb_{u,v}) = u\gA$. We \egmt see that the content of $\fb_{u,v}$ is $1$. The other points are easy.

\emph {3.}
We have $\fb_{u,v,w} = \gA u \oplus \gA z$, $\fb_{w,v,u} = \gA w \oplus \gA z$. The product of these two \ids is generated (as an \id or \Amoz) by the four \elts $uw$, $uz$, $wz$, $z^2$, all multiples of $z$ (because $z^2 + 2vz + uw = 0$). It therefore suffices to see that

\snic {
z \in \gen {uw, uz, wz, z^2}_\gB = \gen {uw, uz, wz, 2vz}_\gB =
\gen {uw, uz, wz, vz}_\gB
.}

However, $v^2 - uw = f$ is separable, therefore $1 \in \gen {u,w,v}_\gA$, and $z \in \gen {uz,wz,vz}_\gB$.

As for $\fb_{u,-v}$ it is $\gA u \oplus \gA \ov z$ with $z\ov z = uw$ and $z + \ov z = -2v$. The  product $\pi$ of the two \ids $\fb_{u,v}$ and $\fb_{u,-v}$ is equal to $\gen{u^2, u\ov z, uz, z\ov z}$, with $z\ov z = uw$, so $\pi\subseteq \gen{u}$. Finally, $-2uv = uz + u\ov z \in \pi$ and therefore $\pi\supseteq \gen{uv, u^2, uw}=u\gen {v, u, w}=\gen{u}$.

Finally, with $u = u_1u_2$, we have $\fb_{u_1,v}\fb_{u_2,v} = \gA u + \gA u_1z +
\gA u_2z + \gA z^2$ clearly included in $\gA u + \gA z = \fb_{u,v}$. As $z^2 + 2vz + uw = 0$ we obtain

\snic {\mathrigid1mu
\gA u + \gA u_1z + \gA u_2z + \gA z^2 =
\gA u + \gA u_1z + \gA u_2z + \gA vz = 
\gA u + (\gA u_1 + \gA u_2 + \gA v)z
.}

Hence $\fb_{u_1,v}\fb_{u_2,v} = \fb_{u,v}$; indeed, $1\in\gen{u_1,u_2, v}_\gA$ because $v^2 - u_1u_2w = f$ is separable, therefore $\gen{u_1,u_2, v}_\gA z = \gA z$.

\emph {4a.}
Let $\fb$ be a nonzero \itf of $\gB$. As a free \Amo of rank~$2$, it admits an $\gA$-basis $(e_1, e_2)$ and we let $M = \cmatrix {a & b\cr 0 &d\cr}$ be the matrix of $(e_1, e_2)$ over $(1,y)$.
We write that $\fb$ is an \idz, \cad $y\fb \subseteq \fb$: 
the membership $ye_1 \in \gA e_1 \oplus \gA e_2$ gives a multiple $a$ of $d$ and the membership $ye_2 \in \gA e_1 \oplus \gA e_2$ gives a multiple $b$ of $d$.
Ultimately, $M$ is of the form $M = d \cmatrix {u & -v\cr 0 &1\cr}$ and we obtain $\fb = d\fb_{u,v}$. We see that $\gen {d}_\gA$ is the content of $\fb$, and $d$ is unique if we impose $d$ as unitary.


\emph {4b.}
We have seen that $\fb = d\fb_{u,v}$ therefore $\ov\fb = d\fb_{u,-v}$ then $\fb\ov\fb = d^2u\gB$. \\
But we also have $\rN(\fb) = d^2u\gA$ because $d\cmatrix {u & -v\cr 0 & 1\cr}$ is the matrix of an $\gA$-basis of~$\fb$ over an $\gA$-basis of $\gB$.
We deduce that $\fb\ov\fb = \rN(\fb)\gA$. Then, for two nonzero \ids $\fb_1, \fb_2$ of $\gB$

\snic {
\rN(\fb_1\fb_2)\gB = \fb_1\fb_2\ov {\fb_1\fb_2} =
\fb_1\ov{\fb_1} \fb_2\ov{\fb_2} = \rN(\fb_1)\rN(\fb_2)\gB
,}

hence $\rN(\fb_1\fb_2) = \rN(\fb_1)\rN(\fb_2)$ since the three \ids are \idps of $\gA$.

\emph {4c.}
First of all, if $\fb$ is a nonzero \itf of $\gB$, it contains a \ndz \elt $b$ and $a = \rN(b) = b\wi b$ is a \ndz \elt of $\fb$ contained in $\gA$. We then have a surjection $\gB\sur{a\gB} \twoheadrightarrow \gB\sur{\fb}$ and as $\gB\sur{a\gB}$ is a finite dimensional $\gk$-\evcz, the same goes for $\gB\sur{\fb}$.

If $d \in \gA \setminus \{0\}$, we have an exact sequence

\snic {
0 \to \gB\sur{\fb'} \simeq d\gB\sur{d\fb'} \to \gB\sur{d\fb'} \to
\gB\sur{d\gB} \to 0
.}

We deduce that $\deg(d\fb') = \deg(\fb') + \deg(d\gB) = \deg(\fb') + \deg(d^2)$.
In particular, for $\fb' = \fb_{u,v}$ and $\fb = d\fb_{u,v}$, we obtain 

\snic {
\deg(\fb) = \deg(u) + \deg(d^2) = \deg \rN(\fb)
.}

This shows that $\deg$ is additive.

\emph {5.}
We first provide a reduction algorithm of $(u,v)$ satisfying $v^2 \equiv f \bmod u$. Even if it entails replacing $v$ by $v \bmod u$, we can assume that $\deg v < \deg u$. If $\deg u \le g$, then, by rendering $u$ \monz, $(u,v)$ is reduced. Otherwise, with $v^2 - uw = f$ let us show that $\deg w < \deg u$; this will allow us to consider $\wi u := w$, $\wi v := (-v) \bmod \wi u$, having the \prt $\fb_{u,v} \sim \fb_{\tilde u, \tilde v}$ and to iterate the process $(u,v) \leftarrow (\wi u, \wi v)$ until we obtain the inegality $\deg u \le g$. To show $\deg u > g \Rightarrow \deg w < \deg u$, we consider the two following cases; either $\deg(uw) > 2g+1 = \deg f$, in which case the \egt $f + uw = v^2$ provides $\deg(uw) = 2\deg v < 2\deg u$ so $\deg w < \deg u$; or $\deg(uw) \le 2g+1$, in which case $\deg w \le 2g+1 - \deg u < 2g+1 - g$ so $\deg w \le g < \deg u$.

Every \id $\fb_{u,v}$ is therefore associated with a reduced \id and as every nonzero \itf $\fb$ of $\gB$ is associated with an \id $\fb_{u,v}$, $\fb$ is therefore associated with a reduced \idz.

\emph {6a.}
Let $w$ satisfy $v^2 - uw = f = y^2$; as $(u,v)$ is reduced, we have

\snic {
\deg v < \deg u \le g < g+1 \le \deg w
\quad \hbox {and} \quad
\deg u + \deg w = 2g + 1
.}

Let $y' = y-v$ and $z = au + by'$ with $a, b \in \gA$. 
\\
We have $y' + \ov {y'} = -2v$, $y \ov {y'} = -(y^2 - v^2) = uw$, so

\snic {
\rN(z)= z\ov z = a^2u^2 + aub(y' + \ov {y'}) +
b^2 y'\ov {y'} = u(a^2u - 2vab + b^2w)
,}

hence 
$\rN(z)/\rN(\fb_{u,v}) = \rN(z)/u = a^2u - 2vab + b^2w$, a
\pol whose degree we need to bound from below. Consider the special case $b = 0$
(therefore $a \ne 0$) in which case $\rN(z)/u = a^2u$, of degree $2\deg a + \deg u \ge \deg u$. Here we see that the \egt $\deg (\rN(z)/u) = \deg u$ is reached \ssi $\deg a = 0$, \cad \ssi $z \in \gk\eti u$. 
\\
There is also the special case $a = 0$ (therefore $b \ne 0$) in which case $\rN(z)/u = b^2 w$,
which is of degree $2\deg b + \deg w > \deg u$.

Therefore it remains to show that for $a \ne 0$, $b \ne 0$, we have $\deg(\rN(z)/u)
> \deg u$. We introduce  $\alpha = \deg a \ge 0$, $\beta = \deg b \ge 0$ and

\snic {\mathrigid1mu
d_1 = \deg(a^2u) = 2\alpha + \deg u, \
d_2 = \deg(vab) = \alpha + \beta + \deg v, \
d_3 = \deg(b^2w) = 2\beta + \deg w
.}


We have $d_1 + d_3 \equiv \deg u + \deg w = 2g+1 \bmod 2$ so $d_1 \ne d_3$. Also, $\alpha \ge \beta \Rightarrow d_1 > d_2$ and finally $\beta \ge \alpha \Rightarrow d_3 > \max(d_1,d_2)$.  
\\
If $d_3 > \max(d_1,d_2)$, then
$\deg (\rN(z)/u) = d_3 \ge \deg w > \deg u$. If $d_3 \le \max(d_1,d_2)$, then
$\alpha > \beta$, so $d_1 > d_2$, then $d_1 > \max(d_2,d_3)$. 
\\
We therefore have $\deg (\rN(z)/u) = d_1 = 2\alpha + \deg u \ge 2 + \deg u > \deg u$.

\emph {6b.}
We have $\fb' = d\fb_{u_1,v_1}$ and $\deg(\fb') = 2\deg(d) + \deg(\fb_{u_1,v_1})$.
We can therefore assume \hbox{that $d = 1$}.  We have $c$, $c_1 \in \gB \setminus \{0\}$
with $c \fb_{u,v} = c_1\fb_{u_1,v_1}$, which we denote by $\fb$. We have $\rN(\fb) = u\rN(c) = u_1\rN(c_1)$. The minimum degree of the $\rN(z)/\rN(\fb)$'s for $z \in \fb \setminus \{0\}$ is $\deg u$ and it is uniquely reached for $z \in \gk\eti cu$.  \\
For $z = c_1u_1 \in \fb$, we have $\rN(z) = u_1^2 \rN(c_1)$
therefore $\rN(z)/\rN(\fb) = 
\fraC{u_1^2 \rN(c_1)}{u_1\rN(c_1)}= u_1$.
We therefore have $\deg u_1 \ge \deg u$, \cad $\deg(\fb_{u_1,v_1}) \ge \deg(\fb_{u,v})$.
The \egt is only possible \hbox{if $c_1u_1 \in \gk\eti cu$}. In this case, $u\fb_{u_1,v_1} = u_1\fb_{u,v}$. Since the content of $u\fb_{u_1,v_1}$ is~$u$, and since that of $u_1\fb_{u,v}$ is $u_1$, the previous \egt entails $u = u_1$ then $v = v_1$.

\emph {7a.}
We have $F'_X(X,Y) = -f'(X)$, $F'_Y(X,Y) = 2Y$. \\
As $\car(\gk)\neq 2$, we obtain $f(X)\in\gen {F, F'_X, F'_Y}$, then $1\in\gen {F, F'_X, F'_Y}$.

\emph {7b.}
We realize the change of variable $\gx = 1/x$ in 

\snic {y^2 = f(x) = x^{2g+1} + a_{2g} x^{2g}
+ \cdots + a_1 x + a_0,}

and we multiply by $\gx^{2g+2}$ to obtain

\snic {
\gy^2 = \gx + a_{2g}\gx^2 + \cdots + a_0\gx^{2g+2} = \gx \big(1 + \gx h(\gx)\big)
\quad \hbox {with} \quad \gy = y\gx^{g+1}
.}


Recap: the change of variable $\gx = 1/x$, $\gy = y/x^{g+1}$ gives $\gk(x) = \gk(\gx)$ and $\gk(x,y) = \gk(\gx, \gy)$, and $\gy$ is integral over $\gk[\gx]$, a fortiori over $\gA_\infty$. 
\\
Let $\gB_\infty = \gk[\gx,\gy]_{\gen {\gx,\gy}}$; in this localized \riz,
 we have $\gen {\gx, \gy} = \gen {\gy}$ since $\gx = \frac{\gy^2 } {1 + \gx h(\gx)}$.  Conclusion: $\gB_\infty$ is a  DVR with regular parameter $\gy$.

Finally, let $\gW$ be a \adv for $\gk(x,y)$ containing $\gk$.  
\\
If $x \in \gW$, then
$\gk[x] \subset \gW$. Then  $y$, integral over $\gk[x]$, is in $\gW$, therefore $\gB  \subset \gW$.\\
If $x \notin \gW$,
we have $x^{-1} \in \fm(\gW)$, so $\gA_\infty = \gk[x^{-1}]_{\gen {x^{-1}}}
\subset \gW$, and $\gW = \gB_\infty$.


\prob{exoTrifolium} 
Let $\vep$ be the unit defined by \smashtop{\framebox[1.1\width][c]{$\vep = \beta-\alpha$}}.

\emph {1.}
We decompose $G$ and $H$ into \hmgs components $G_i, H_j$,

\snic {
G = G_a + \cdots + G_b, \  a \le b, \quad
H = H_c + \cdots + H_d, \  c \le d
.}

The  lower  \hmg component
 of $GH$, of degree $a+c$, is $G_aH_c$ while 
the  upper  \hmg component of $GH$, of degree $b+d$, is $G_bH_d$.
We deduce that $a+c=N$, $b+d=N+1$; we cannot have $a < b$ and $c < d$  at the same time (because we would then have $a+c+2 \le b+d$, \cad $N+2 \le N+1$). If $a=b$, then $G$ is \hmgz, if $c=d$ it is $H$.
Suppose that $F_N$, $F_{N+1}$ are coprime and let $F = GH$ be a \fcnz; for example, $G$ is \hmg of degree $g$; we deduce that $H = H_{N-g} + H_{N+1-g}$ and that $F_N = GH_{N-g}$,
$F_{N+1} = GH_{N+1-g}$: $G$ is a common factor of $F_N, F_{N+1}$, so $G$ is \ivz. The converse is easy.

The \pols $(X^2 + Y^2)^2$ and $\alpha X^2Y + \beta Y^3 = Y(\alpha X^2 + \beta
Y^2)$ are coprime \ssi the \pols $X^2 + Y^2$ and $\alpha X^2 + \beta Y^2$ are coprime; \cad \ssi $\alpha \ne \beta$.

\emph {2.}
The reader will verify that $(0,0)$ is the only singular point; we have the more precise result
$$\preskip-.3em \postskip.0em 
\vep^2 X^5, \vep^2 Y^5 \in \gen {F, F'_X, F'_Y}. 
$$

\emph {3.}
Let $Y = TX$ in $F(X,Y)$. We obtain $F(X, TX) = X^3 G(X,T)$ with

\snic {
G(X, T) = XT^4 + \beta T^3 + 2X T^2 + \alpha T + X
.}

The \pol $G$ is primitive (in $T$) and $(x=0,t=0)$ is a simple point of the curve $G = 0$. With $a_4 = x$, $a_3 = \beta$, $a_2 = 2x$, $a_1 = \alpha$, $a_0 = x$, we consider the integral elements (by Emmanuel's trick)
$$\preskip.4em \postskip.4em 
b_4 = a_4,\quad  b_3 = a_3 + tb_4,\quad b_2 = a_2 + tb_3,\quad
b_1 = a_1 + tb_2 
. 
$$
Thus, $b_4 = x$, $b_3 = \beta+y$ and $b_2 = 2x + (\beta+y)y/x$.

It is clear that $a_4$, $a_3$, \dots, $a_0 \in \sum_i \gA b_i + \sum_i \gA tb_i$.
As $a_3-a_1 = \vep$ is \ivz,  there are $u_i, v_i \in \gA$ such that $1 = \sum_i u_i b_i + \sum_i v_i tb_i$. We formally write (without worrying about the nullity of a $b_i$)
$$\preskip.0em\postskip.2em \ndsp
t = {b_1t \over b_1} = \cdots = {b_4t \over b_4} =
{\sum_i v_ib_it \over \sum_i v_i b_i} = 
{\sum_i u_ib_it \over \sum_i u_i b_i}.
$$
Thus, $t= y/x = a/b = c/d$ with $a$, $b$, $c$, $d \in \gB$ and $a+d = 1$. 
\\
The \egts $by = ax$, $dy = cx$,  $a+d = 1$ are those coveted. 
\\
Thus we obtain $\fq\gen{x,y}_\gB = \gen {x}_\gB$ with $\fq =\gen {d,b}_\gB$. Here by letting

\snic {
a = b_2t - b_4t, \quad  b = b_2 - b_4,\quad
c = b_3t - b_1t, \quad  d = b_3 - b_1 
,}


we have $\vep = a+d$. By letting $g_0 = 1/(1+t^2)$, $g_1 = tg_0$, we find $b = \vep g_1$, $d = \vep g_0$, so $\fq = \gen {g_0,g_1}_\gB$. We will show (question \emph {5}) that $\gB = \gk[g_0,g_1]$,
so $\gB\sur\fq = \gk$.

\emph {4.}
A \gmq idea leads to the \egt $\gk(t) = \gk(x,y)$. It is the parameterization of the trifolium. The \pol defining the curve is of degree $4$ and the origin is a singular point of multiplicity $3$. Therefore a rational line passing through the origin intersects the curve at a rational point.
Algebraically, this corresponds to the fact that the \pol $G(X,T)$ is of degree $1$ in $X$

\snic {
G(T,X) = (T^4 + 2T^2 + 1)X + \beta T^3 + \alpha T =
(T^2 + 1)^2 X + T(\beta T^2 + \alpha)
,}

hence
$$\preskip-.4em \postskip.4em
x = -{t(\beta t^2 + \alpha) \over (t^2 + 1)^2}, \quad
y = tx = -{t^2(\beta t^2 + \alpha) \over (t^2 + 1)^2}.
$$
At $t = 0$, we have $(x,y) = (0,0)$. What are the other values of the parameter $t$ for which $\big(x(t), y(t)\big) = (0,0)$? \\
We have to first find the zeros of $x(t)$, a rational fraction of height~$4$. There is the value $t = \infty$, for which $y(t) = -\beta$.  \\
If $\alpha = 0$, we only have two zeros of $x$: $t=0$ (of multiplicity $3$) and $t = \infty$ (of multiplicity $1$).  
\\
If
$\beta = 0$, we only have two zeros of $x$: $t=0$ (of multiplicity $1$)
and $t = \infty$ (of multiplicity $3$).  
\\
If $\beta \ne 0$, we have two other zeros of $x$ (eventually coinciding): $t = \pm \sqrt {-\alpha/\beta}$.
We can render this more uniform by making the quadratic character of $-\alpha\beta$ intervene, see question \emph {7.}

Remark: in all the cases, at $t = \infty$, we have $(x,y) = (0,-\beta)$.

\emph {5.}
We know by Exercise~\ref{exoAnneauOuvertP1} that $\gk[g_0,g_1]$ is an \aclz, the \cli of $\gk[g_0]$ in $\gk(t)$. To obtain a $\gk$-relator between $g_0$ and $g_1$, we substitute $t = g_1/g_0$ in
the expression $g_0 = 1/(1+t^2)$, which gives $g_0^2 - g_0 + g_1^2 = 0$ and confirms that $g_1$ is integral over $\gk[g_0]$. At $t=0$, we have $(g_0,g_1) = (1,0)$; this point is a nonsingular point of the curve $g_0^2 - g_0 + g_1^2 = 0$. Actually the conic $C(g_0, g_1) = g_0^2 - g_0 + g_1^2$ is smooth over every \ri since

\snic {
1 = -4C + (2g_0-1)\Dpp{C}{g_0} + 2g_1 \Dpp{C}{g_1}.
}


The same goes for the homogenized conic denoted by $C$, $C =
g_0^2 - g_0g_2 + g_1^2$, which satisfies $\gen {g_0,g_1,g_2}^2 \subseteq \gen {C, \Dpp{C}{g_0}, \Dpp{C}{g_1}, \Dpp{C}{g_2}}$
$$
g_0 = -\Dpp{C}{g_2}, \quad g_1^2 = C + (g_0-g_2)\Dpp{C}{g_2}, \quad
g_2 = -\Dpp{C}{g_0} - 2\Dpp{C}{g_2}.
$$
We dispose of $\PP^1 \to \PP^2$ defined by $(u : v) \mapsto (g_0 : g_1 : g_2) = (u^2 : uv : u^2 + v^2)$ whose image is the \hmg conic $C = 0$; more or less, this is a Veronese embedding $\PP^1 \to \PP^2$ of degree $2$.

Moreover, the \dcn into simple \elts provides the following expressions of $x$, $y$, $b_3t$, $b_2t$ in $\gk[g_0,g_1]$

\snic {
\begin {array} {c}
x = \vep g_0g_1 - \beta g_1, \quad
y = \vep g_0^2 + (2\beta-\alpha)g_0 - \beta =  (g_0-1)(\beta - \vep g_0)
\\[2pt]
b_2t = 2y + (\beta+y)t^2  = -\alpha + \beta g_0 - \vep g_0^2, \quad
b_3t = (2\beta-\alpha)g_1 - \vep g_0g_1
.
\end {array}
}

We see that $g_0$ is integral over $\gk[y]$, therefore integral over $\gk[x]$;
as $g_1$ is integral over~$\gk[g_0]$, it also is integral over $\gk[x]$. We have just obtained the \egt $\gB = \gk[g_0,g_1]$.

First consider $\gk[y] \subset \gk[g_0] \subset \gk[g_0,g_1]$; it is clear that $(1,g_0)$ is a basis of $\gk[g_0]$ over $\gk[y]$ and $(1,g_1)$ is a basis of $\gk[g_0,g_1]$ over $\gk[g_0]$, therefore $(1, g_0, g_1, g_0g_1)$ is a basis of $\gk[g_0,g_1]$ over $\gk[y]$ (but not over $\gA = \gk[x]$).  
\\
Let us show that $(1, y, b_3 t, b_2t)$ is an $\gA$-basis, let $E$ be the generated \Amoz. By using $y-b_2t = \vep(g_0-1)$ and $x+b_3t = \vep g_1$, we see that $g_0$, $g_1 \in E$. Finally, $E$ contains $x + \beta g_1 = \vep g_0g_1$, so $g_0g_1 \in E$ and $E = \gk[g_0,g_1] = \gB$.

An \iv \id $\fb$ of $\gB$ contains a \ndz \elt therefore $\gB\sur\fb$ is a finite dimensional \kevz, which allows us to define $\deg \fb$ by $\deg \fb = \dim_\gk \gB\sur\fb$; we then have (see Proposition~\ref{prop-a/ab} and its Corollary~\ref{corprop-a/ab}) $\deg(\fb\fb') = \deg(\fb) + \deg(\fb')$.  We deduce that $\deg \gen {x,y}_\gB = 4-1 = 3$.

\emph {6.}
We have $\fp_1 = \gen {g_0-1, g_1}$, therefore to show the \egt $\fp_1^2 = \gen{g_0-1, g_1^2}$, it suffices to see that $g_0-1 \in \gen {(g_0-1)^2, g_1^2}$. This results from the \egt $1-g_0 = (1-g_0)^2 + g_1^2$ which stems from $g_0^2 - g_0 + g_1^2 = 0$.

\emph {7.}
Let $X = UY$ in $F(X,Y)$. We obtain $F(UY,Y) = U^3H(U,Y)$
with 
$$\preskip.4em \postskip.4em 
H(U,Y) = YU^4 + (2Y + \alpha) U^2 + Y + \beta, \qquad
H(U,0) = \alpha U^2 + \beta
. 
$$
This \pol $H = a'_4U^4 + a'_2U^2 + a'_0$ is primitive in $U$ (we have $a'_2 = 2a'_4
+ \alpha$ \hbox{and $a'_0 = a'_4 + \beta$} therefore $\epsilon = a'_0 - a'_2 +
a'_4$). It satisfies $H(u,y) = 0$ with $u = x/y$; 
the associated \eltz~$b'_3$ determined by Emmanuel's trick
 is $x$ and we therefore have $b'_3u \in \gB$ with
$$\preskip-.4em \postskip.4em 
b'_3u = x^2/y = \vep g_0 - \beta - y.
$$
In root $t$ of $\beta t^2 + \alpha = 0$, we have $g_0 = \beta/\vep$ and $g_1^2 =
-\alpha\beta / \vep^2$, which renders the introduction of the \id $\fa =
\gen {\vep g_0 - \beta, \vep^2 g_1^2 + \alpha\beta}$ natural. We verify the \egt

\snic {
\gen {y, x^2/y}_\gB = 
\gen {\vep g_0 - \beta,  \vep^2 g_1^2 + \alpha\beta}
.}

We then have $\gen {x,y}_\gB = \fp_1\fa$ and $\deg \fa = 2$.  If $-\alpha\beta$
is not a square, then $\fa$ is prime. Otherwise, we have $\fa = \fp_2\fp_3$ with $\fp_2, \fp_3$ expressed with the two square roots of $-\alpha\beta$. We have $\fp_2 = \fp_3$ \ssi the two square roots are confused; this happens when $\alpha\beta = 0$ for example or in \cara $2$.
Finally, for $\alpha = 0$, we have $\fp_1 = \fp_2 = \fp_3$.



\penalty-2500    
\Biblio

\vspace{2pt}

Regarding the genesis of the theory of \ids of number fields developed by Dedekind, we can read the articles of H.\ Edwards \cite{Edw81} and of J.\ Avigad~\cite{Avi}.

The  \ddps were introduced by H.\ Pr\"ufer in 1932 in \cite{Prufer}.
Their central role in the multiplicative theory of \ids is highlighted in the reference book on the subject \cite{Gil}. 
See also the bibliographic comments at the end of Chapter~\ref{chap mod plats}.

In the classical literature a \adpc is often called a 
\ixy{semihereditary}{ring} (according to item \emph{3} in \thref{th.adpcoh}), which is not very pretty. 
These \ris are signaled as important in \cite{CE}.
The \cov proof of item \emph{1} of \thref{ThImMat} is given in Chapter 1, Proposition~6.1.%
\index{semihereditary!\riz}%
\index{hereditary!\riz}

A \ixy{hereditary}{ring} is a \ri in which every \id is \proz. This notion is badly defined in \coma because of the non-legitimate quantification \gui{every \idz.}
An example of such a non-\noe \ri is the sub\ri of a countable product of
the  field $\FF_2$, formed by the sequences which are either
 null almost everywhere,
or equal to $1$ almost everywhere.  
The most interesting case is that of \noe \adpcsz, which we described in \clama as the \ris in which every \id is \ptfz.
Our \dfn of a \adk (freed from the integrity constraint) corresponds exactly (in \clamaz) to the notion of a \noe hereditary \riz.

Fairly comprehensive presentations on  \anars and  \adps written in the style of \coma can be  found in the articles \cite[Ducos\&al.]{dlqs} and~\cite[Lombardi]{lom99}.

\gui{Emmanuel's trick} of Lemma~\ref{lemEmmanuel} appears in Emmanuel Hallouin's PhD thesis \cite{TheseHallouin}. 

\Thref{th.2adpcoh} is due to Gilmer and Hoffmann \cite{GHo}.
\Thref{dekinbe} for the case of a  \ddp is given by Heitman and Levy in~\cite{HeLe}.
\Thref{thcohdim1} has been proven in \clama by Quentel in~\cite{Que}. The \cov \dem is due to I. Yengui.

\Thref{thPTFDed} is classical (Steinitz's \thoz)
for  \adksz.
It has been \gne for the \ddps having the one and a half \prt in \cite[Kaplansky]{Kap52} and 
\index{one and a half!\tho}
\cite[Heitmann\&Levy]{HeLe}. A detailed inspection of our \dem  actually shows that the hypothesis \gui{\ri \ddi $1$} could be weakened to \gui{\ri having the one and a half \prtz.}

We find \thref{thMpfPruCohDim} (see also Exercise~\ref{exoAnarlgb}) in   \cite[Brewer\&Klinger]{BrKl} for the integral case. It has been \gne to the case of a \qiri in \cite[Couchot]{Couc1}.

Lemma~\ref{lemRadJDIM1} and \thref{lemthAESTE} are due to Claire T\^ete and Lionel Ducos. 

\Pbmz~\ref{exoBaseNormaleAlInfini} is based on the article~\cite[Hess]{Hess02}.

\thref{thIntClosStab} says that
if $\gA$ is \iclz, the same goes for~$\AX$. A constructive proof of the same result for normal rings is given in \cite{cl2016}.

\newpage \thispagestyle{CMcadreseul}

\incrementeexosetprob


\chapter{Krull dimension}
\label{chapKrulldim}
\perso{compil\'e le \today}
\minitoc

\subsection*{Introduction}
\addcontentsline{toc}{section}{Introduction}

In this chapter we introduce the \ddk in its \elr \cov version and we compare it to the corresponding classical notion.

Next we establish the first \prts of this dimension. The ease with which we obtain the \ddk of a \pol \ri over a \cdi shows that the \cov version of the \ddk can be seen as a conceptual simplification of the usual classical version. 

We then apply the same type of ideas to define the \ddk of a \trdiz, that of a morphism of commutative \risz, then the valuative dimension of commutative \risz.

We establish a few basic important \thos regarding these notions.

We finish by indicating the \cov versions of the usual classical notions of Lying Over, Going Up, Going Down and Incomparability, with some applications.

\section{Spectral spaces}
\label{secEspSpectraux}

In this section, we describe the classical approach to Krull dimension.

For us, this is above all a matter of heuristics. It is for this reason that we give no proofs.
This will have no incidence in the rest of the book.
Indeed, the \cof aspect of the spectral spaces is entirely concentrated in the \trdis obtained by duality. In particular, the \cov aspect of the \ddk is entirely concentrated in the \ddk of the \trdis and it can be defined completely independently from the spectral spaces.

Nevertheless the heuristic given by the spectral spaces is essential to the understanding of the small miracle that will happen with the introduction of the dual \cov notions. 
This small miracle will only fully be realized in the following chapters, where we will see the transformation of many beautiful abstract theorems into algorithms.

\subsec{The Zariski lattice and the Zariski spectrum}
\label{subsecZarZar}

Recall that we denote by $\DA(\fa)$ the nilradical of the \id $\fa$ in the \ri $\gA$ and that the Zariski lattice $\ZarA$ is the set of the $\DA(x_1,\ldots ,x_n)$'s (for~$n\in\NN$ and $x_1$, \ldots, $x_n\in\gA$).
We therefore have $x\in\DA(x_1,\ldots ,x_n)$ \ssi a power of $x$ belongs to $\gen{x_1,\ldots ,x_n}$.
The set $\ZarA$, ordered by the inclusion relation,
is a \trdi with 

\snic{\DA(\fa_1)\vu\DA(\fa_2)=\DA(\fa_1+\fa_2)\;
\hbox{ and } \;\DA(\fa_1)\vi\DA(\fa_2)=\DA(\fa_1\,\fa_2).}

\begin{definition}
\label{nota Spec(A)}\relax
We  call the set of strict prime \ids of the \riz~$\gA$ the \ix{Zariski spectrum} of  $\gA$  and we denote it by~$\SpecA$. We equip it with the topology that has as its basis of open sets the $\fD_\gA(a)=\sotq{\fp\in\SpecA}{a\notin\fp}$.\\
We denote by $\fD_\gA(\xn)$  the set  $\fD_\gA(x_1)\cup\cdots \cup\fD_\gA(x_n)$.

\noi For $\fp\in\SpecA$ and $S=\gA\setminus\fp$, we denote $\gA_S$ by $\Ap$  
(the ambiguity between the two contradictory notations $\Ap$ and $\gA_S$ is removed in practice by the context).
\end{definition}

In \clamaz, we then obtain the following result.

\begin{theoremc}\label{thZarZar} ~
\begin{enumerate}
\item The \oqcs of $\SpecA$ are the open sets $\fD_\gA(x_1,\ldots ,x_n)$.
\item The map $\DA(x_1,\ldots ,x_n)\mapsto \fD_\gA(x_1,\ldots ,x_n)$ is well-defined.
\item It is an \iso of \trdisz.
\end{enumerate}
\end{theoremc}

\subsec{Spectrum of a \trdiz}
\label{subsecSpecT}

The Zariski spectrum is the paradigmatic example of a \ix{spectral space}.
Spectral spaces were introduced by Stone \cite{Sto} in 1937.%
\index{spectral!space}%
\index{space!spectral ---}

They can be \cares as the topological spaces satisfying the following \prtsz
\begin{itemize}
\item the space is quasi-compact,
\item every open set is a union of \oqcsz,
\item the intersection of two \oqcs is a \oqcz,
\item for two distinct points, there is an open set containing one of them but not the other,
\item every \ird closed set is the adherence of a point.
\end{itemize}
The \oqcs then form a \trdiz, the supremum and the infimum being the union and the intersection, respectively. A continuous map between spectral spaces is said to be \ixc{spectral}{map} if the inverse image of every \oqc is a \oqcz.
Stone's fundamental result can be stated as follows.

\emph{In \clama the category of spectral spaces and spectral maps is anti-\eqve to the category of \trdisz. }

\smallskip  Here is how this works.

First of all if $\gT$ is a \trdiz,
a \emph{\idepz} is an \id $\fp$ which satisfies
$$\begin{array}{llllll}\preskip.0em \postskip.4em
x\vi y \in \fp \;\Rightarrow\; (x\in\fp\;\mathrm{or}\; y\,\in\fp) ,&   &
1_\gT\notin\fp.
\end{array}$$
\rdb\label{SpecTrdi}%
The \emph{spectrum} of $\gT$, denoted by $\SpecT$, is then defined as the space whose points are the \ideps of $\gT$ and which has a basis of open sets given by the subsets $\fD_{\gT}(a):=\sotq{\fp\in\SpecT}{a\notin\fp}$ for $a\in\gT$.%
\index{spectrum!of a \trdiz}
\\
If $\varphi:\gT\to\gV$ is a morphism of \trdisz,
we define the map 

\snic{\Spec\varphi:\Spec\gV\to\SpecT,\quad \fp\mapsto\varphi^{-1}(\fp).}

It is a spectral map and all of this defines $\Spec$ as a contravariant functor.
\\
We show that the $\fD_{\gT}(a)$'s are all the \oqcs of $\SpecT$.
Actually \Thoz\eto \ref{thZarZar} applies to every \trdi $\gT$:
{\it \begin{enumerate}
\item The \oqcs of $\SpecT$ are exactly the $\fD_\gT(u)$'s.
\item The map $u\mapsto \fD_\gT(u)$ is well-defined and it is an \iso of \trdisz.
\end{enumerate}
}\label{NOTASpecT}
In the other direction, if $X$ is a spectral space we let $\OQC(X)$ be the \trdi formed by its \oqcsz.
If $\xi:X\to Y$ is a spectral map, the map 

\snic{\OQC(\xi):\OQC(Y)\to\OQC(X),\quad U\mapsto\xi^{-1}(U)}

is a \homo of \trdisz.
This defines $\OQC$ as a contravariant functor.

The stated anti-\eqvc of categories is defined by the functors~$\Spec$ and $\OQC$.
It generalizes  the anti-\eqvc given in the finite case by  \thref{thDualiteFinie}.

Note that the empty spectral space corresponds to the lattice $\Un$, and that a reduced spectral space at a point corresponds to the lattice $\Deux$.

\subsec{Spectral subspaces}
\label{subsecSesSpec}

By \dfnz,
a subset $Y$ of a spectral space $X$ is a \ix{spectral subspace} 
	if $Y$ is a spectral space by the induced topology 
and if the canonical injection $Y\to X$ is spectral.
\\
This notion is actually exactly the dual notion of the notion of a quotient \trdiz.
In other words a spectral map $\alpha:Y\to X$ identifies $Y$ with a spectral subspace of $X$ \ssi the \homo of \trdis $\OQC(\alpha)$ identifies $\OQC(Y)$ to a quotient \trdi of~$\OQC(X)$.

The closed subspaces of $X$ are spectral and correspond to the quotients by the \idsz. More \prmt an \id $\fa$ of $\OQC(X)=\gT$ defines the closed set~$\fV_\gT(\fa)=\sotq{\fp \in X}{\fa \subseteq \fp}$,
(provided we identify the points of~$X$ with the \ideps of $\OQC(X)$)
and we then have a canonical \iso
$$\OQC(\fV_\gT(\fa)\big)\simeq\OQC(X)\sur{(\fa=0)}.$$
The \ird closed sets correspond to the \ideps of $\OQC(X)$.

Finally, the \oqcs correspond to the quotients by principal filters
$$\OQC(\fD_\gT(u)\big)\simeq\OQC(X)\sur{(\uar u=1)}.$$

\subsec{A heuristic approach to the \ddkz}
\label{subsecHKdim}

Note moreover that the Zariski spectrum of a commutative \ri is naturally identified with the spectrum of its Zariski lattice.

In \clamaz, the notion of \ddk can be defined, for an arbitrary spectral space~$X$, as the maximal length of the strictly increasing chains of \ird closed sets.

An intuitive way to apprehend this notion of dimension is the following. The dimension can be \caree by \recu by saying that on the one hand, the dimension $-1$ corresponds to the empty space, and on the other hand, for $k\geq 0$, a space $X$ is \ddi$k$ \ssi for every \oqc $Y$, the boundary of $Y$ in $X$ is \ddi$k-1$ (this boundary is closed therefore it is a spectral subspace of $X$).

Let us see, for example, for a commutative \ri $\gA$, how we can define the boundary of the open set $\fD_\gA(a)$ in $\Spec\gA$. The boundary is the intersection of the adherence of $\fD_\gA(a)$ and of the \cop closed set of $\fD_\gA(a)$, which we denote by $\fV_\gA(a)$.
The adherence of  $\fD(a)$ is the intersection of all the $\fV(x)$'s that contain $\fD(a)$, \cad such that  $\fD(x)\cap\fD(a)=\emptyset$.
\\
As $\fD(x)\cap\fD(a)=\fD(xa)$, and as we have $\fD(y)=\emptyset$ \ssi $y$ is nilpotent, we obtain a heuristic approach to the \id \gui{Krull boundary of~$a$,}  which is the \id generated by $a$ on the one hand
(which corresponds to $\fV(a)$),
and by all the $x$'s such that $xa$ is nilpotent on the other hand (which corresponds to the adherence of $\fD(a)$).

\section[A \cov \dfnz]{\Cov definition and first consequences}
\label{secDefConsDimKrull}

In \clamaz, the \ddk of a commutative \ri is defined as the maximum (eventually infinite) of the lengths of the strictly increasing chains of strict \ideps (beware, a chain $\fp_0\subsetneq \cdots\subsetneq \fp_\ell$ is said to be of length $\ell$).
Since the complement of a \idep is a prime filter, the \ddk is also the maximum of the lengths of the strictly increasing chains of prime filters.

As this \dfn is impossible to manipulate from an \algq point of view, we replace it in \coma by an \eqv \dfn (in \clamaz) but of a more \elr nature.

The quantification over the set of \ideps of the \ri is then replaced by a quantification over the \elts of the \ri and the  non-negative integers. Since this discovery (surprisingly it is very recent) the \thos that make the \ddk intervene have been able to become fully integrated into \coma and into Computer Algebra.

\begin{definition}
\label{defZar2} Let  $\gA$ be a commutative \riz, $x\in\gA$ and $\fa$ be a \itfz.
\begin{enumerate}
\item [$(1)$] The \ix{Krull upper boundary} of $\fa$ in $\gA$ is the quotient \ri
\begin{equation}\label{eqBKAC}\preskip.4em \postskip.5em
\gA_\rK^{\fa}:=\gA/\JK_\gA(\fa)  \quad \hbox{where} \quad
 \JK_\gA(\fa):=\fa+(\sqrt{0}:\fa).
\end{equation}
Write $\JK_\gA(x)$ for $\JK_\gA(x\gA)$ and $\gA_\rK^{x}$ for $\gA_\rK^{x\gA}$. This \ri is called the \emph{upper boundary of $x$ in $\gA$}. \\
We will say that $\JK_\gA(\fa)$ is \emph{the Krull boundary \id of $\fa$ in $\gA$.}%
\index{ideal!Krull boundary ---}
\item [$(2)$] The \ix{Krull lower boundary} of $x$ in $\gA$ is the localized \riz
\begin{equation}\label{eqBKAS}\preskip.4em \postskip.5em
\gA^\rK_{x}:=\SK_\gA(x)^{-1}\!
\gA \quad\hbox{where}\quad \SK_\gA(x)=x^\NN(1+x\gA).
\end{equation}
We will say that $\SK_\gA(x)$ is the \emph{Krull boundary monoid of $x$} in~$\gA$.%
\index{monoid!Krull boundary ---}%
\index{Krull boundary!ideal}\index{Krull boundary!monoid}
\end{enumerate}
\end{definition}

Recall that in \clama the \ddk of a \ri is $-1$ \ssi the \ri does not admit any \idepsz, which means that it is trivial.

The following \tho then gives  
in \clama an \elr inductive \carn of the \ddk of a commutative \riz.

\begin{theoremc}
\label{thDKA} For a commutative \ri  $\gA$ and an integer $k\geq 0$ \propeq
\begin{enumerate}
\item  The \ddk of $\gA$ is $\leq k$.
\item  For all $x\in \gA$ the \ddk of $\gA_\rK^{x}$ is $\leq k-1$.
\item  For all $x\in \gA$ the \ddk of $\gA^\rK_{x}$ is $\leq k-1$.
\end{enumerate}
\end{theoremc}
Note: 
this is a \tho of \clama which cannot admit a \prcoz. \eoe

In the \dem that follows all the prime or maximal \ids and filters are taken in the usual sense in \clamaz: they are strict. 
\begin{proof}
Let us first show the \eqvc of items \emph{1} and \emph{3.}
Recall that the \ideps of $S^{-1}\gA$ are of the form $S^{-1}\fp$ where $\fp$ is a \idep of $\gA$ which does not intersect $S$ (Fact~\ref{factQuoFIID}).
The \eqvc then clearly results from the two following statements.\\
(a) Let $x\in\gA$, if $\fm$ is a maximal \id of $\gA$ it always intersects $\SK_\gA(x)$. Indeed, if $x\in\fm$ it is clear and otherwise, $x$ is \iv modulo $\fm$ which means that $1+x\gA$ intersects $\fm$.\\
(b) 
Let $\fa$ be an \idz, $\fp$ be a \idep with $\fp\subset\fa$ and $x\in\fa\setminus\fp$;
if $\fp\cap \SK_\gA(x)$ is nonempty, then $1 \in \fa$.
Indeed, let $x^n(1+xy)\in\fp$; since $x\notin\fp$, we~have~$1+xy\in\fp\subset\fa$, which gives, with $x \in \fa$, $1\in\fa$.\\
Thus, if $\fp_0\subsetneq \cdots \subsetneq \fp_\ell$ is a chain with $\fp_\ell$ maximal, it is shortened by at least its last term when we localize at $\SK_\gA(x)$, and it is only shortened by its last term if $x\in\fp_\ell\setminus\fp_{\ell-1}$.

The \eqvc of items \emph{1} and \emph{2} is proven dually, by replacing the \ideps by the prime filters.
Let $\pi:\gA\to\gA/\fa$ be the canonical \prnz. We notice that the prime filters of $\gA/\fa$ are exactly the $\pi(S)$'s, where $S$ is a prime filter of $\gA$ that does not intersect $\fa$  (Fact~\ref{factQuoIDFI}).
It then suffices to prove the two dual statements of (a) and (b) which are the following.\\
(a') Let $x\in\gA$, if $S$ is a maximal filter of $\gA$ it always intersects $\JK_\gA(x)$. Indeed, if $x\in S$ it is clear and otherwise, since $S$ is maximal, $Sx^\NN$ contains~$0$, which means that there is an integer $n$ and an \elt $s$ of $S$ such that $sx^n=0$.
Then $(sx)^n=0$ and $s\in (\sqrt{0}:x)\subseteq \JK_\gA(x)$.\\
(b') 
Let $S'$ be a prime filter contained in a filter $S$ and $x\in S\setminus S'$. If~$S'\cap \JK_\gA(x)$ is nonempty, then $S = \gA$. Indeed, let $ax+b\in S'$ with~$(bx)^n=0$. Then, since $x\notin S'$, we have $ax\notin S'$ and, given that $S'$ is prime,~$b\in S'\subseteq S$. As $x\in S$, we obtain $(bx)^n=0\in S$.
\end{proof}

In \coma we replace the usual \dfn given in \clama by the following more \elr \dfnz.
\begin{definition}
\label{defDiKrull}
The \ixc{Krull dimension}{of a commutative \riz}
(denoted by $\Kdim$) \label{NOTAKdim}
of a commutative \ri $\gA$ is defined by \recu as follows
\begin{enumerate}
\item  $\Kdim\gA=-1$ \ssi $\gA$ is trivial.
\item  For $k\geq 0$,  $\Kdim\gA\leq k$ means $\forall  x\in \gA,\,\Kdim(\gA^\rK_{x})\leq k-1$.
\end{enumerate}
\end{definition}

Naturally, we will say that $\gA$ is infinite dimensional
 \ssi for every integer $k\geq 0$ we have the implication  $\Kdim\gA\leq k\Rightarrow 1=_\gA0$.

\medskip The following lemma immediately results from the \dfnsz.

\begin{lemma}
\label{lemZEdDef}
A \ri is \zed \ssi it is \ddi $0$.
\end{lemma}

Note that the terminology \gui{\zed \riz} therefore constitutes a slight abuse of language because  affirming that the dimension is less than or equal to~$0$ leaves the possibility of a dimension equal to $-1$ open, which means that the \ri is trivial.

\medskip \exls \label{exlKdim} ~
\\
1) If $x$ is nilpotent or \iv in $\gA$, 
         the boundary \id and the boundary \mo of $x$
in $\gA$ are both equal to $\gA$. The two boundary \ris are trivial.

 2)
For $x\neq 0$, $1$, $-1$ in $\ZZ$, the boundary \ris $\ZZ_\rK^x=\ZZ/x\ZZ$ and $\ZZ_x^\rK=\QQ$ are \zedsz.
We therefore find that $\Kdim\ZZ\leq1$ again.

 3)
Let $\gK$ be a field contained in a discrete \cac $\gL$. Let $\fa$ be a \itf of $\gK[X_1,\ldots ,X_n]$ and $\gA=\gK[X_1,\ldots ,X_n]/\fa$.
Let $V$ be the affine \vrt corresponding to $\fa$ in  $\gL^n$ and $W$ be the sub\vrt of~$V$ defined by $f$. Then the \gui{boundary of~$W$ in~$V$,} defined as the intersection of~$W$ with the Zariski closure of its complement in~$V$, is the affine \vrt corresponding to the \riz~$\gA_\rK^f$.
We abbreviate this as

\snic {
\hbox {boundary}_V\, \cZ(f) = \cZ_V(\hbox {boundary of $f$}).
}

\sni 4) Let $\gA$ be integral and $k\geq 0$: $\Kdim\gA\leq k$ is equivalent to $\Kdim(\gA\sur{a\gA})\leq k-1$ for every \ndz $a$ (use the Krull boundary \idsz).

 5) Let $\gA$ be a \dcd local \ri and $k\geq 0$: $\Kdim\gA\leq k$ is equivalent to $\Kdim \gA[1/a]\leq k-1$ for all $a\in\Rad\gA$ (use the Krull boundary \mosz).
\eoe

\medskip
\comms
1) The advantage of the \cov \dfn of the \ddk with respect to the usual \dfn is that it is simpler (no quantification over the set of \idepsz) and more \gnl (no need to assume the axiom of choice).
However, we have only defined the sentence \gui{$\gA$ is \ddi$k$.}

2)
In the context of \clamaz. The \ddk of $\gA$ can be defined as an \elt of $\so{-1}\cup\NN\cup\so{+\infty}$ by letting

\snic{\Kdim\gA=\inf\sotq{k\in\ZZ,\,k\geq -1}{\Kdim\gA\leq k},}

(with $\inf\,\emptyset_\ZZ=+\infty$). This \dfn based on the \cov 
\dfnz~\ref{defDiKrull} is \eqve to the usually given \dfn via  chains of \ideps (see \Thoc\ref{thDKA}).

3)
From the \cov point of view, the previous method does not define the \ddk of $\gA$ as an \elt of $\so{-1}\cup\NN\cup\so{+\infty}$.
Actually it so happens that the concept in question is \gnlt not \ncr (but the reader must take our word for it).
\perso{remarque  qui semble in\'evitable pour
r\'epondre \`a l'angoisse bien naturelle of la lectrice}
\\
The most similar point of view to \clama would be to look at $\Kdim\gA$ as a subset of $\NN\cup\so{-1}$, defined by

\snic{\sotq{k\in\ZZ,\,k\geq -1}{\Kdim\gA\leq k}.}

We then reason with (eventually empty) final subsets of $\NN\cup\so{-1}$, the order relation is given by the reversed inclusion, the upper bound by the intersection and the lower bound by the union.
\\
This approach finds its limit with the \gui{counterexample} of the real number field (see the comment on \paref{remDKRR}).
\eoe

\medskip
We use in \coma the following \emph{notations}, to be closer to the classical language

\begin{notation}
\label{notaKdiminf}
{\rm  Let $\gA$, $\gB$, $(\gA_i)_{i\in I}$,
$(\gB_j)_{j\in J}$ be commutative \ris (with $I$, $J$ finite).
\begin{enumerate}
\item [--]  $\Kdim\gB\leq \Kdim\gA$ means $\forall \ell\geq -1\;
(\Kdim\gA\leq \ell\;\Rightarrow \Kdim\gB\leq \ell)$.
\item [--]  $\Kdim\gB= \Kdim\gA$ means  $\Kdim\gB\leq  \Kdim\gA$ and
$\Kdim\gB\geq  \Kdim\gA$.
\item  [--] $\sup_{j\in J}\Kdim\gB_j\leq  \sup_{i\in I}\Kdim\gA_i$ means
$$\preskip.1em \postskip.3em
\forall\ell\geq -1\quad  \big(\,\&_{i\in I} \, \Kdim\gA_i\leq
\ell \;\Rightarrow\;\&_{j\in J}\, \Kdim\gB_j\leq
\ell  \,\big).
$$
\item  [--] $\sup_{j\in J}\Kdim\gB_j=  \sup_{i\in I}\Kdim\gA_i$ means
$$\preskip.1em \postskip.0em
\forall\ell\geq -1\quad \big(\, \&_{i\in I} \,\Kdim\gA_i\leq
\ell\;\Leftrightarrow\;\&_{j\in J}\,\Kdim\gB_j\leq
\ell  \,\big).
$$
\end{enumerate}
}
\end{notation}

\subsec{Iterated boundaries, \susisz, \cop sequences}

\Dfn \ref{defDiKrull} can be rewritten in terms of \idasz.
For this, we introduce the notion of a \emph{\susiz}.

\begin{definition}\label{notaBordsIteres}
For a sequence  $(\ux)=(\xzk)$ in $\gA$ we define the iterated Krull boundaries as follows.

 1.  An \gui{iterated} version of the \mo $\SK_\gA(x)$: the set
\begin{equation}\label{eqMonBordKrullItere}
\SK_\gA(\xzk):=
x_0^\NN(x_1^\NN\cdots(x_k^\NN (1+x_k\gA) +\cdots)+x_1\gA) + x_0\gA)
\end{equation}
is a \moz.
For an empty sequence, we define $\SK_\gA()=\so{1}$.%
\index{monoid!iterated Krull boundary ---}%
\index{Krull boundary!iterated --- monoid}%

 2. We define two variants for the iterated Krull boundary \idz.
\\
 --- 2a) The \id $\JK_\gA(\xzk)=\JK_\gA(\ux)$ is defined as follows
\begin{equation}\label{eq1IdBordKrullItere}
\JK_\gA()=\so{0}, \;
\JK_\gA(x_0, \ldots, x_k) =
\big(\DA\big(\JK_\gA(x_0, \ldots, x_{k-1})\big) : x_k\big) + \gA x_k.%
\index{ideal!iterated Krull boundary ---}%
\index{Krull boundary!iterated --- ideal}%
\end{equation}
 --- 2b)  The \id $\IK_\gA(\xzk)=\IK_\gA(\ux)$ is defined as follows
\begin{equation}\label{eqIdBordKrullItere}
\IK_\gA(\ux):=\sotq{y\in\gA}{\,0\in
x_0^\NN\big(\cdots\big(x_k^\NN (y+x_k\gA) +\cdots\big) + x_0\gA\big)}
\end{equation}
For an empty sequence, we define $\IK_\gA()=\so{0}$.

\end{definition}

We will show (Lemma~\ref{fact2BordKrullItere}) that the two \gui{iterated boundary} \ids defined above have the same nilradical.

\begin{definition}
\label{defSeqSing}
A sequence $(x_0,\ldots ,x_k)$ in  $\gA$ is said to be \emph{singular} 
if $0\in\SK_\gA(\xzk)$, in other words if $1\in\IK_\gA(\xzk)$, \cad if there exist $a_0$, \ldots, $a_k\in \gA$ and $m_0$, \ldots, $m_k\in \NN$ such that
\begin{equation}\label{eqsing}\preskip.4em \postskip.0em
x_0^{m_0}(x_1^{m_1}(\cdots(x_k^{m_k} (1+a_k x_k) +
\cdots)+a_1x_1) + a_0x_0) =  0
\index{singular!sequence}\index{sequence!singular ---}
\end{equation}%
\end{definition}

\vspace{-.6em}
\pagebreak	        

\begin{proposition}
\label{corKrull}
For a commutative \ri $\gA$  and an integer $k\geq 0$, \propeq
\begin{enumerate}
\item The \ddk of $\gA$ is $\leq k.$
\item For all $x\in \gA$ the \ddk of $\gA_\rK^{x}$ is $\leq k-1$.
\item \label{i3corKrull} Every sequence $(\xzk)$ in $\gA$  is \singz.
\item \label{i4corKrull} For all $x_0$, \ldots, $x_k \in \gA$ there exist $b_0$, \ldots, $b_k\in \gA$ such that
\begin{equation}\label{eqCG}
\left.\arraycolsep2pt
\begin{array}{rcl}
\DA(b_0x_0)& =  &\DA(0),    \\
\DA(b_1x_1)& \leq  & \DA(b_0,x_0),  \\
\vdots~~~~& \vdots  &~~~~  \vdots \\
\DA(b_k x_k )& \leq  & \DA(b_{k -1},x_{k -1}),  \\
\DA(1)& =  &  \DA(b_k,x_k ).
\end{array}
\right\}
\end{equation}
%
\item  \label{i5corKrull} For all $x_0$, \ldots, $x_k \in \gA$, by letting $\pi_i = \prod_{j < i} x_j$ for $i \in \lrb {0..k+1}$ (so $\pi_0=1$), there exists an $n\in\NN$ such that

\snic{\pi_{k+1}^n \in \gen {
\pi_{k}^n x_k^{n+1},\
\pi_{k-1}^n x_{k-1}^{n+1},\ \ldots,\
\pi_1^n x_1^{n+1},\
\pi_0^n x_0^{n+1}}.\qquad\qquad\quad}%
\end{enumerate}
\end{proposition}
For example, for $k=2$ item~\emph{\ref{i4corKrull}} corresponds to the following graph in $\ZarA$.
$$
\SCO{\DA(x_0)}{\DA(x_1)}{\DA(x_2)}{\DA(b_0)}{\DA(b_1)}{\DA(b_2)}
$$

\begin{proof}
The \eqvcs for  dimension $0$ are \imdes by application of the \dfnsz.
 
\emph{1} $\Leftrightarrow$ \emph{3.} Suppose the \eqvc established for  dimension at most $k$ and for every commutative \riz. We then see that $S^{-1}\gA$ is \ddi$k$ \ssi we have
\\
\emph{for all $x_0$, \ldots, $x_k\in \gA$ there exist $a_0$, \ldots, $a_k\in \gA$, $s\in S$ and $m_0$, \ldots, $m_k\in \NN$ such that}
\begin{equation}\label{eqsing3}
x_0^{m_0}(x_1^{m_1}\cdots(x_k^{m_k} (s+a_k x_k) +
\cdots + a_1x_1) + a_0x_0) =  0.
\end{equation}
Note that with respect to \Eqrf{eqsing}, some $s\in S$ has replaced the $1$ in the center of the expression on the left-hand side.
\\
It therefore remains to replace $s$ by an arbitrary \elt of $\SK_\gA(x_{k+1})$, \cad an \elt of the
form $x_{k+1}^{m_{k+1}} (1+a_{k+1} x_{k+1})$.
\\
The \eqvc between \emph{2} and \emph{3} is proven analogously.
 
\emph{3} $\Rightarrow$ \emph{4.} We take $b_k=1+a_k x_k$, then $b_{\ell -1} = x_\ell^{m_\ell} b_\ell+ a_{\ell -1}x_{\ell -1}$, successively for $\ell=k,$ \ldots, $1$.
 
\emph{4} $\Rightarrow$ \emph{2.} Proof by \recuz.
The implication for  dimension $\leq 0$ is clear. Suppose it established for  dimension $< k $. Assume  \prtz~\emph{4} and let us show that for all $x_0$ the dimension of $\gB= \gA\ul{x_0}$  is $<k$.
\\
By \hdr it suffices to find, for all $x_1$, \ldots, $x_k$, some \elts $b_1$, \ldots, $b_k$  such that
$$\left.\arraycolsep2pt
\begin{array}{rcl}
\rD_\gB(b_1x_1)& =  & \rD_\gB(0)  \\
\vdots~~~~& \vdots  &~~~~  \vdots \\
\rD_\gB(b_k x_k )& \leq  & \rD_\gB(b_{k -1},x_{k -1})  \\
\rD_\gB(1)& =  &  \rD_\gB(b_k,x_k ).
\end{array}
\right\}
$$
However, by hypothesis we have some \elts $b_0$, \ldots, $b_k$  such that
$$
\left.\arraycolsep2pt
\begin{array}{rcl}
\DA(b_0x_0)& =  &\DA(0)    \\
\DA(b_1x_1)& \leq  & \DA(b_0,x_0)  \\
\vdots~~~~& \vdots  &~~~~  \vdots \\
\DA(b_k x_k )& \leq  & \DA(b_{k -1},x_{k -1})  \\
\DA(1)& =  &  \DA(b_k,x_k ).
\end{array}
\right\}
$$
and the in\egts with $\DA$ imply the same in\egts with $\rD_\gB$.
 The second in\egt means that $(b_1x_1)^m \in \gen {b_0, x_0}$ (for a certain $m$); the first tells us that $b_0x_0$ is nilpotent therefore $\gen {b_0,x_0} \subseteq \JK_\gA(x_0)$. Recap: $b_1x_1$ is nilpotent in $\gB$.
\\
We could also prove \emph{4} $\Rightarrow$ \emph{3} by a direct, slightly more tedious, computation.
 
\emph{3} $\Leftrightarrow$ \emph{5.} In the \dfn of a \susiz, we can replace all the exponents $m_i$ by their maximum $n$. Once this is acquired, item \emph{5} is a simple reformulation of item \emph{3.}
\end{proof}

We could therefore have given a \dfn by \recu of the \ddk based on the upper boundaries $\gA_\rK^{x}$ rather than on the lower boundaries $\gA^\rK_{x}$:
we have just obtained a direct \prco (without using \Thoc\vref{thDKA}) of the \eqvc between the two possible inductive \dfnsz.

\medskip \rem
The \sys of in\egts (\ref{eqCG}) in item~\emph{\ref{i4corKrull}} of Proposition~\ref{corKrull} establishes an interesting and \smq relation between the two sequences~$(b_0,\ldots ,b_k)$ and $(x_0,\ldots ,x_k)$. \\
When $k=0$,
this means $\DA(b_0)\vi\DA(x_0)=0$ and $\DA(b_0)\vu\DA(x_0)=1$,
that is that the two \elts $\DA(b_0)$ and $\DA(x_0)$ are complements of one another in the lattice $\ZarA$.
In $\Spec\gA$ this means that the corresponding basic open sets are \copsz. \\
We therefore introduce the following terminology: when the sequences $(b_0,\ldots ,b_k)$ and $(x_0,\ldots ,x_k)$ satisfy the in\egts (\ref{eqCG}) we say that they are two \emph{\cop sequences}.%
\index{complementary sequences!in a commutative \riz}%
\eoe

\begin{fact}\label{fact0BordKrullItere}
Let $(\ux) = (\xn)$ and $(\uy) = (\ym)$ be two sequences of \elts of $\gA$,
$\gA \to \gA'$ be a morphism and $(\uxp)$ be the image of $(\ux)$ in~$\gA'$.
\begin{enumerate}
\item
We have the \eqvcs
$$\preskip.4em \postskip.4em
\exists z \in\IK_\gA(\ux) \cap \SK_{\gA}(\uy)   \iff
1 \in \IK_\gA(\ux, \uy) \iff 0 \in \SK_\gA(\ux,\uy).
$$

\item
If $\gA \to \gA'$ is surjective, the image of $\SK_\gA(\ux)$ is $\SK_{\gA'}(\uxp)$.

\item
If $\gA' = S^{-1}\gA$, with $S$ being a \mo of $\gA$, then $S^{-1} \IK_{\gA}(\ux) = \IK_{\gA'}(\uxp)$.

\end{enumerate}
\end{fact}

\begin{fact}\label{fact1BordKrullItere}
Let $\fa$ be an \id of $\gA$, $Z \subseteq \gA$ be an arbitrary subset and~$x \in \gA$. 
$$\preskip.0em \postskip.4em 
x^\NN(Z + \gA x) \hbox { meets } \fa
\;  \iff \;
Z  \hbox { meets } (\fa : x^\infty) + \gA x. 
$$
\end{fact}


\begin{lemma}\label{lem1BordKrullItere} \emph{(Krull boundary \ids \`a la Richman)}\\
For a sequence $(\ux) = (\xn)$ of \elts of $\gA$, the iterated boundary \id $\IK_\gA(\ux)$ can be defined recursively as follows
$$\preskip.4em \postskip.2em
\IK_\gA()=\so{0}, \qquad
\IK_\gA(x_1, \ldots, x_n) =
\big(\IK_\gA(x_1, \ldots, x_{n-1}) : x_n^\infty\big) + \gA x_n.
$$
For example,
$$\preskip.3em \postskip.4em
\IK_\gA(x_1) = (0 : x_1^\infty) + \gA x_1, \;
\IK_\gA(x_1, x_2) =
\big(\big((0 : x_1^\infty) + \gA x_1\big) : x_2^\infty\big) + \gA x_2.
$$
\end{lemma}

\begin{proof}
We temporarily define
$$\preskip.4em \postskip.4em
N() = \{0\}, \quad
N(x_1, \ldots, x_n) = \big(N(x_1, \ldots, x_{n-1}) : x_n^\infty\big) + \gA x_n
$$
Take $n = 3$ to fix the ideas. Then, for $y \in \gA$, we have the \eqvcs
$$
\begin{array} {ll}
0 \in x_1^\NN\big(x_2^\NN\big(x_3^\NN(y + \gA x_3) + \gA x_2\big) + \gA x_1\big)
&\quad  \Longleftrightarrow \\[1mm]
x_2^\NN\big(x_3^\NN(y + \gA x_3) + \gA x_2\big) \hbox { meets } N(x_1)
&\quad  \Longleftrightarrow \\[1mm]
x_3^\NN(y + \gA x_3) \hbox { meets }
\big(N(x_1) : x_2^\infty\big) + \gA x_2 \eqdefi  N(x_1, x_2)
&\quad  \Longleftrightarrow \\[1mm]
y \in \big(N(x_1,x_2) : x_3^\infty\big) + \gA x_3 \eqdefi
  N(x_1, x_2, x_3),
\end{array}
$$
which proves that $\IK_\gA(x_1, x_2, x_3) = N(x_1, x_2, x_3)$.
\end{proof}

\begin{lemma}\label{lem2BordKrullItere} 
\emph{Iterating boundary \idsz}\\
Let $(\ux) = (\xn)$ and $(\uy) = (\ym)$ be two sequences of \elts of~$\gA$. 
Let $\gA' = \gA/\IK_\gA(\ux)$ and let $(\uyp) = (\ypm)$ be the image of $(\uy)$ in~$\gA'$.
\begin{enumerate}
\item
The kernel of the (surjective) canonical morphism $\gA \to \gA'\big/\IK_{\gA'}(\uyp)$ is the \id $\IK_\gA(\ux, \uy)$.
\item
We define $\gA_0 = \gA$ and $\gA_i = \gA_{i-1}/\IK_{\gA_{i-1}}(x_i)$ for $ i\in \lrbn$. Then the kernel of the (surjective) canonical morphism $\gA \to \gA_n$ is the \id $\IK_\gA(\ux)$.
\end{enumerate}
\end{lemma}

\begin{proof}
It suffices to prove the first item for $n = 1$. Let $x = x_1$.
\\
Let $z \in \gA$  and $z'$ be its image in $\gA' = \gA/\IK_{\gA}(x)$.
We have the \eqvcs

\snic{\begin{array} {ll}
z=0 \;\;\mathrm{ in }\; \gA'/\IK_{\gA'}(\uyp)
&\quad  \Longleftrightarrow \\[1mm]
0 \in {y'_1}^\NN\big(\cdots\big({y'_m}^\NN (z' + y'_m\gA') +\cdots\big) + y'_1\gA'\big)
&\quad  \Longleftrightarrow \\[1mm]
{y_1}^\NN\big(\cdots\big({y_m}^\NN (z + y_m\gA) +\cdots\big) + y_1\gA\big)
\hbox { rencontre } \IK_\gA(x)
&\quad  \Longleftrightarrow \\[1mm]
0 \in x^\NN({y_1}^\NN\big(\cdots\big({y_m}^\NN (z + y_m\gA) +\cdots\big) + y_1\gA) + x\gA\big)
&\quad  \Longleftrightarrow \\[1mm]
z \in \IK_\gA(x, \uy).
\end{array}
}
\vspace{-1.5em}
\end{proof}

\begin{fact}\label{fact2BordKrullItere}
For every sequence $(\ux)$ of \elts of $\gA$, the \ids $\IK_\gA(\ux)$ and~$\JK_\gA(\ux)$ have the same nilradical.
\end{fact}

\begin{proof}
For every \id $\fa$ and all $x \in \gA$, we easily prove that the nilradical of the \id $(\fa : x^\infty)$ is $(\DA(\fa) : x)$. By using $\DA(\fb + \fc) = \DA(\DA(\fb) + \fc)$, we deduce that the \ids $(\fa : x^\infty) + \gA x$ and $(\DA(\fa) : x) + \gA x$ have the same nilradical.
The stated result is then deduced by \recu on the length of the sequence $(\ux)$ by using the recursive \dfn of the two iterated boundary \idsz.
\end{proof}

\begin{lemma}\label{lem3BordKrullItere}
Let $S$ be a \mo of $\gA$, $\gA' = S^{-1} \gA$, $x \in \gA$, $x'$ be its image in $\gA'$ and $V=\SK_{\gA'}(x')$.  Then the canonical morphism $\gA \to V^{-1}\gA'$ is a \lon morphism\footnote{See  \Dfnz~\ref{defHomloc}.} and induces an \iso \hbox{of $T^{-1}\gA$} \hbox{over $ V^{-1}\gA'$}, where $T$ is the \mo $x^\NN(S + \gA x)$.
\end{lemma}

\begin{proof}
The image in $V^{-1}\gA'$ of the \elt $s + ax \in S + \gA x$ is \iv since we can write $s+ ax = s(1 + ax/s)$ (with a few notation abuses). Hence a (canonical) morphism $\varphi : T^{-1}\gA \to V^{-1}\gA'$. \\ 
Moreover, since $S \subseteq T$, we have a morphism $\gA' \to T^{-1}\gA$. The image under this morphism of $1 + xa/s
\in 1 + x\gA'$ is \iv because $1 + xa/s = (s + xa)/s$, hence a morphism $\varphi' : V^{-1}\gA' \to T^{-1}\gA$.\\
We prove without difficulty that $\varphi$ and $\varphi'$ are inverses of one another.
\end{proof}

\begin{corollary}\label{cor1BordKrullItere} 
\emph{Iterating boundary \mosz}\\
Let $(\ux) = (\xn)$ and $(\uy) = (\ym)$  
in $\gA$, $\gA' = \SK_\gA(\uy)^{-1}\gA$, \hbox{and $(\uxp) = (\xpn)$} be the image of $(\ux)$ in~$\gA'$.
\\
Then, the morphism $\gA \to \SK_{\gA'}(\uxp)^{-1}\gA'$ gives by \fcn an \iso $\SK_{\gA}(\ux,\uy)^{-1}\gA \simarrow \SK_{\gA'}(\uxp)^{-1}\gA'$.
\end{corollary}

\subsec{A regular sequence \gui{is not} singular}

Item~\emph{4} of the following proposition implies that a regular sequence that does not generate the \id $\gen{1}$ is nonsingular, which explains the title of this subsection.

An advantage of the iterated Krull boundaries \`a la Richman is that over a \coh \noe \ri they are \itfsz. Another advantage is given by item \emph{1} in the following proposition.

\begin{proposition}
\label{lemRegsing} \emph{(Regular sequences and \ddkz)}\\
Let $(\xn)$ be a regular sequence in $\gA$ and $(\yr)$ be another sequence. 
\begin{enumerate}
\item \label{i1lemRegsing} We have $\IK_\gA(\xn)=\gen{\xn}$.
\item \label{i2lemRegsing} The sequence $(\xn,\yr)$ is singular in $\gA$ \ssi the sequence $(\yr)$ is singular in $\aqo\gA\xn$.
\item \label{i3lemRegsing} The following implication is satisfied for every $k\geq -1$,

\snic{\Kdim \gA \leq n+k\;\Longrightarrow\;\Kdim \aqo\gA\xn \leq k.}

If $1\notin\gen{\xn}$, we have ${n+\Kdim \aqo\gA\xn \leq \Kdim\gA}$. 
\item \label{i4lemRegsing} 
If the sequence $(\xn)$ is \egmt singular, we have $\gen{\xn}=\gen{1}$. \\
Consequently if $\Kdim\gA\leq n-1$ every regular sequence of length $n$ generates the \id $\gen{1}$.
\end{enumerate}
\end{proposition}
%
\begin{proof}
\emph{1.} Immediate computation taking into account the recursive \dfn given in Lemma~\ref{lem2BordKrullItere} (item~\emph{2}).\\
\emph{2.} We apply item \emph{1} of Lemma~\ref{lem2BordKrullItere}.\\
\emph{3.} Results from item~\emph{2.}\\
\emph{4.} Special case of item~\emph{2}, with the empty sequence $(\yr)$.
\end{proof}
%

\subsec{Lower bounds of the \ddk}

It can be comfortable, sometimes even useful, to define the statement \gui{$\gA$ is of \ddk $\geq k$.}

First of all $\Kdim\gA\geq 0$ must mean $1\neq0$. For $k\geq1$, a possibility would be to ask: \gui{there exists a sequence $(x_1,\ldots ,x_k)$ which is not singular.}
Note that from the \cov point of view this affirmation is stronger than the negation of \gui{every sequence  $(x_1,\ldots ,x_k)$ is singular.}

A \ri then has a well-defined \ddk if there exists an integer~$k$ such that the \ri is both of \ddk
$\geq k$ and of \ddkz~$\leq k$.

The annoying thing is the negative \crc of the assertion

\vspace{.2em}
\centerline{\gui{the sequence $(x_1,\ldots ,x_k)$ is not singular.}}

Anyway here it seems impossible to avoid the use of the negation, because we do not see how we could define $\Kdim\gA\geq 0$ other than by the negation $1\neq0$. Naturally, in the case where $\gA$ is a discrete set, \gui{$x\neq0$} loses its negative \crcz, and the statement \gui{there exists a sequence $(\ux)=(x_1,\ldots ,x_k)$ such that $0\notin\IK_\gA(\ux)$} does not strictly speaking contain any negation.

However, note that the \dfns of $\Kdim\gA\leq k$ and $\Kdim\gA\geq k$ use an alternation of quantifiers that introduces an infinity (for an infinite \riz). Consequently the \dfn cannot \gnlt be certified by a simple computation: a proof is needed.

Note that for the \ri $\RR$, if we use the strong negation (of positive \crcz), 
 for which $x\neq0$ means \gui{$x$ is \ivz,} to define the sentence $\Kdim\RR\geq k$,
then it is absurd that $\Kdim\RR\geq1$. But we cannot \cot prove $\Kdim\RR\leq0$ (see the comment on \paref{remDKRR}).

\section{A few \elr \prts of the \ddkz}
\label{secKrullElem}

The stated facts in the following proposition are easy (note that we use the notation introduced in \ref{notaKdiminf}).

\begin{proposition}
\label{propDdk0} Let $\gA$ be a \riz, $\fa$ be an \id and $S$ be a \moz.
\begin{enumerate}
\item A \susi remains singular in $\gA\sur{\fa}$ and $\gA_S$.
\item $\Kdim \gA\sur{\fa} \leq \Kdim \gA$, $\Kdim \gA_S \leq \Kdim \gA$.
\item $\Kdim (\gA\times \gB) = \sup( \Kdim \gA, \Kdim\gB)$.
\item $\Kdim\gA=\Kdim\Ared$.
\item If $a$ is \ndz in $\Ared$ (a fortiori if it is \ndz in $\gA$), %
then $\Kdim \gA\sur{a\gA}\leq \sup(\Kdim\gA,0)-1$.
\item If $a\in\Rad\gA$, then $\Kdim \gA[1/a] \leq \sup(\Kdim\gA,0) -1$.
\end{enumerate}
 \end{proposition}

\exl We give a \ri $\gB$ for which $\Frac(\gB)$ is of \ddk $n>0$,
but $\gB\red=\Frac(\gB\red)$ is \zedz.
\\
Consider $\gB = \gA\sur{x\fm}$, where $\gA$ is local \dcdz,
$\fm=\Rad\gA$ and~$x \in \fm$.  The \ri $\gB$ is local, $\Rad\gB=\fm' = \fm/x\fm$ and $\gB/\fm'=\gA/\fm$.\\ 
If $\overline x = 0$, then $x \in x\fm$,
\cad $x(1-m) = 0$ with $m \in \fm$, which implies $x = 0$.
\\
For $y\in \fm$ we have $\ov y\,\ov x=0$. 
Therefore if $\ov y \in\Reg \gB\cap\fm'$, we obtain $x=0$. 
\\
However, we have $\ov y\in\fm'$ or $\ov y\in\gB\eti$, 
therefore if $x\neq 0$ and~$\ov y \in\Reg \gB$, we obtain $\ov y\in\gB\eti$. In other words, if $x\neq 0$, $\gB=\Frac(\gB)$.
\\
Take $\gA = \gk[\xzn]_{\gen {\xzn}}$ where $\gk$ is a nontrivial \cdiz, 
\hbox{and $x = x_0$}.
We then have 

\snic{\aqo\gA{x_0} \simeq \gk[x_1, \ldots, x_n]_{\gen {x_1, \ldots, x_n}}\hbox{   and  }\Kdim\aqo\gA{x_0} = n.}

As $\ov {x_0}^{2}=0$ in $\gB$, we have $\gB\red \simeq \aqo\gA{x_0}$ and therefore $\Kdim\gB=n$.\\
Finally, $\Frac(\gB\red)=\gk(\xn)$ is a \zed \cdiz.
\\
Geometrically, we have considered the \ri of a \vrt \gui{with multiplicities} consisting at a point immersed in a hyperplane of dimension $n$, and we have localized at this immersed point. 
\\
Note: in \clamaz, if $\gC$ is \noe and reduced, $\Frac(\gC)$ is a finite product of fields, therefore \zedz. For a \cov version we  refer to \Pbmz~\ref{exoAnneauNoetherienReduit} and to \cite[Coquand\&al.]{cls}.
\eoe

\begin{plcc}\label{thDdkLoc}
\emph{(For the \ddkz)}\\
Let $S_1$, $\ldots$, $S_n$ be \moco of a \ri $\gA$ and $k\in\NN$.
\begin{enumerate}
\item A sequence is singular in $\gA$ \ssi it is singular in each of the $\gA_{S_i}$'s.
\item The \ri $\gA$ is \ddi$k$ \ssi the $\gA_{S_i}$'s are \ddi$k$.
\end{enumerate}
\end{plcc}
We could have written $\Kdim\gA=\sup_{i}\Kdim\gA_{S_i}=\Kdim\prod_{i}\gA_{S_i}$.
\begin{proof}
It suffices to prove the first item.
Consider a sequence $(\xzk)$ in~$\gA$.
We look for $a_0$, \ldots, $a_k\in \gA,$ and $m_0$, \ldots, $m_k\in \NN$ such that
$$\preskip.4em \postskip.4em
x_0^{m_0}\big(x_1^{m_1}\cdots\big(x_k^{m_k} (1+a_k x_k) +
\cdots + a_1x_1\big) + a_0x_0\big) =  0.
$$
An equation of this type at the $a_j$'s can be solved in each of the $\gA_{S_i}$'s. We notice that if in a \ri  $\gA_{S_i}$ we have a solution for certain exponents $m_0$, \ldots, $m_k$ then we also have a solution for any \sys of larger exponents.
Therefore by taking a \sys of exponents that bound from above each of those obtained separately for each  $\gA_{S_i}$, we obtain a unique \lin equation in the $a_j$'s which has a solution in each  $\gA_{S_i}$.
We can therefore apply the basic \plg \ref{plcc.basic}.\iplg
\end{proof}

As the \prt for a sequence to be singular is of \carfz, item \emph{1} in the previous \plgc actually always applies with a family of \ecoz, which corresponds to a finite covering of the Zariski spectrum by basic open sets.

In the case of a finite covering by closed sets, the result still holds.

\begin{prcf}\label{thDdkRecFer} \emph{(Krull dimension)}\\
Let $\gA$ be a \riz, $k$ be an integer $\geq 0$, and  $\fa_1$, \ldots, $\fa_r$ be \ids of $\gA$.\\
First we assume that the $\fa_i$'s form a closed covering of $\gA$.
\begin{enumerate}
\item A sequence $(\xzk)$ is singular in $\gA$ \ssi it is singular in each of the $\gA\sur{\fa_i}$'s.
\item The \ri $\gA$ is \ddi$k$ \ssi each of the $\gA\sur{\fa_i}$'s is \ddi$k$.
\end{enumerate}
More \gnltz, without a hypothesis on the $\fa_i$'s we have
\begin{enumerate}\setcounter{enumi}{2}
\item The \ri $\gA\big/\!\bigcap_i\fa_i$ is \ddi$k$ \ssi each of the $\gA\sur{\fa_i}$'s is \ddi$k$.\\
This can be abbreviated to 

\snic{\Kdim\gA\big/\!\prod_i\fa_i=\Kdim\gA\big/\!\bigcap_i\fa_i=\sup_{i}\Kdim\gA\sur{\fa_i}=\Kdim\prod_{i}\gA\sur{\fa_i}.}
\end{enumerate}

\end{prcf}
\begin{proof}
It suffices to prove item \emph{1}. The sequence $(\xzk)$ is singular \ssi the \mo $\SK(\xzk)$ contains $0$. \\
In addition, $\SK_{\gA\sur{\fa_i}\!}(\xzk)$ is none other than $\SK_\gA(\xzk)$ considered modulo~$\fa_i$. \Trf by the \prf \vref{prcf1}.
\end{proof}

\goodbreak
\begin{theorem} \emph{(One and a half \thoz)}%
\index{one and a half!\tho}
\label{th1.5}\relax
\begin{enumerate}
\item
\begin{enumerate}
\item If $\gA$ is \zedz, or more \gnlt if $\gA$ is \lgbz, every \lmo module is cyclic.
\item  If $\gA$ is \zedz, every \ptf \id is generated by an \idmz.
\end{enumerate}
\item Let $\gA$ be \ddi$k$, let $(x_1,\ldots ,x_k)$ be a regular sequence and $\fb$ be a \lop \id containing $\fa=\gen{x_1,\ldots ,x_k}$. Then there exists a $y\in\fb$ such that
$$
\preskip.3em \postskip.3em 
 \fb=\gen{y,x_1,\ldots ,x_k}=\gen{y}+\,\fb\fa=\gen{y}+\,\fa^m
$$
for any exponent $m\geq 1$.
\item Let $\gA$ be such that  $\gA\sur{\Rad\gA}$ is \ddi$k$, %
let $(x_1, \ldots, x_k)$ be a regular sequence in $\gA\sur{\Rad\gA}$ and $\fb$ be a \ptf \id of $\gA$ containing $\fa=\gen{x_1,\ldots ,x_k}$ then there exists a $y\in\fb$ such that
$$
\preskip-.4em \postskip.2em 
\fb=\gen{y,x_1,\ldots ,x_k}=\gen{y}+\,\fb\fa=\gen{y}+\,\fa^m 
$$
for any exponent $m\geq 1$.
\end{enumerate}
\end{theorem}
\begin{proof}
Item \emph{1a} is a reminder (see item \emph{\iref{ite4propZerdimLib}} of \thref{propZerdimLib}
for the \zed \ris and item \emph{2} of \thref{thlgb3} for the \algbsz).
 
 \emph{1b.} Recall that in an arbitrary \ri a \ptf \idz~$\fa$ has as its annihilator an \idmz~$h$. In $\aqo\gA h$, $\fa$ is faithful, therefore so is $\fa^k$, for all $k\geq1$. In $\gA[1/h]$, $\fa=0$. Therefore $\Ann(\fa^k)=\Ann(\fa)=\gen{h}$ 
for $k\geq1$.
\\
In the \zed case, since a \ptf \id is \lopz, it is principal by item \emph{1a}, let us denote it by $\gen{x}$. We know that for large enough $k$, $\gen{x}^k=\gen{e}$ with $e$ \idmz. By the preliminary remark $\Ann (x)=\Ann (e)= \gen{1-e}$.
In $\aqo\gA {1-e}$, $x$ is \ivz, so $\gen{x}=\gen{1}$;  
in $\aqo\gA {e}$, $x$ is null; thus in $\gA$, $\gen{x}=\gen{e}$.
 
\emph{3.} Results from \emph{2} by Nakayama's lemma.
 
\emph{2.}
The \id $\fb$ seen as an \Amoz, after \eds to $\gA\sur{\fa}$, becomes the module $\fb\sur{\fb\fa}$ and it remains \lmoz. Since the quotient \ri $\gA\sur{\fa}$ is \zedz, item~\emph{1a} tells us that $\fb\sur{\fb\fa}$ is generated by an \elt $y$. This means that $\fb=\gen{y}+\,\fb\fa$ and the other \egts \imdt follow.
\end{proof}

\rem
In the case of dimension $1$ and of an \iv \idz, item \emph{2} of the previous \tho is often called the \gui{One and a half \thoz.} See Corollary~\ref{corpropZerdimLib} and \thref{th1-5}.
\eoe

\section{Integral extensions}\label{secDDKExtEn}

\begin{proposition}
\label{propDKEXENT}
Let $\gA\subseteq\gB$ be \ris with $\gB$ integral over $\gA$.
Every finite sequence of \elts of $\gA$ that is singular in $\gB$ is singular in $\gA$.
In particular, $\Kdim\gA\leq \Kdim\gB$. 
\end{proposition}
Note: the opposite in\egt is proven a little later (\thref{cor2thKdimMor}).
\begin{proof} 
Suppose for example that the sequence $(x,y)\in\gA$ is singular in $\gB$, \cad 

\snic{\exists a,b\in\gB, \;\exists m,\ell\in\NN,\;  x^\ell\big(y^m(1+ay)+bx\big)=0.}

We want to realize the same type of \egtz, with some \elts $a'$, $b'$ of $\gA$ instead of  \elts $a$, $b$ in $\gB$. The intuitive idea is to transform the previous \egt by the \gui{norm} operation.
Consider some \polusz~$f$,~$g\in\AT$ that annihilate $a$ and $b$.
Let $\gB_1=\aqo{\gA[T,T']}{f(T),g(T')}$.
Let $\alpha$ and $\beta$ be the classes of $T$ and $T'$ in $\gB_1.$
The sub\ri $\gA[a,b]$ of~$\gB$ is a quotient of $\gB_1=\gA[\alpha,\beta]$, via an $\gA$-\homo which sends~$\alpha$ and~$\beta$ to~$a$ and~$b$.
In addition, $\gB_1$ is a free module of finite rank over $\gA$ which allows for a \dfn of the norm and the cotransposed \elt of an arbitrary \elt of~$\gB_1[X,Y]$.
Then let

\snic{U  (\alpha,\beta,X,Y)=X^\ell\big(Y^m(1+\alpha Y)+\beta X\big)\;$ and $\;V(X,Y)=\rN_{\gB_1[X,Y]/\gA[X,Y]}(U).}

By Lemma~\ref{lemNormeSuiteSing}, $V(X,Y)$ is of the form

\snic{X^{p}\big(Y^q\big(1+A(Y) Y\big)+B(X,Y) X\big),}

with $A\in\gA[Y]$, $B\in\gA[X,Y]$.
Moreover let $W(\alpha,\beta,X,Y)\in\gB_1[X,Y]$ be the cotransposed \elt of $U(\alpha,\beta,X,Y)$. 
By specializing $X,Y,\alpha,\beta$ at $x,y$ in $\gA$ and $a,b$ in $\gB$, we obtain an \egt in $\gB$

\snic{V(x,y)=x^{p}\big(y^q\big(1+A(y) y\big)+B(x,y) x\big)=U(a,b,x,y)W(a,b,x,y),}

which ends the \dem since $V(x,y)=0$ is an \egt in~$\gA$:
note that we have $U(a,b,x,y)=0$ in $\gB$ but that~$U(\alpha,\beta,x,y)$ is perhaps nonzero in $\gB_1$.
\end{proof}

\begin{lemma}\label{lemNormeSuiteSing}
Let $\gC$ be a free \Alg of finite rank over $\gA$, $(\czn)$ in $\gC$ and $(\Xzn)=(\uX)$ be a list of \idtrsz. 
Let
$$\preskip.4em \postskip.4em
U(\uX) =X_0^{k_0}\big(X_1^{k_1}\big(\cdots(X_n^{k_n} (1+c_nX_n) +\cdots)+c_1X_1\big) + c_0X_0\big)\in\gC[\uX].
$$
Then $V(\uX)\eqdefi\rN_{\gC[\uX]/\gA[\uX]}(U(\uX)\big)$ is of the form
$$\preskip.3em \postskip.3em
V(\uX)=
X_0^{\ell_0}\big(X_1^{\ell_1}\big(\cdots(X_n^{\ell_n}
(1+a_nX_n) +\cdots)+a_1X_1\big) + a_0X_0\big)\in\gA[\uX],
$$
with $a_n\in\gA[X_n]$,  $a_{n-1}\in\gA[X_n,X_{n-1}]$, \ldots,
$a_0\in\gA[\uX]$.
\end{lemma}
\begin{proof}
First of all the norm $\rN(1+c_nX_n)$ is a \pol  $h(X_n)\in\gA[X_n]$ which satisfies $h(0)=1$, therefore which can be expressed in the form $1+a_n(X_n)X_n$. Next we use the multiplicativity of the norm, and an \evn at $X_{n-1}=0$
to show that $\rN\big(X_n^{k_n}(1+c_nX_n)+c_{n-1}X_{n-1}\big)$ is of the form
$$\preskip.0em \postskip.3em 
X_n^{\ell_n}\big(1+a_n(X_n)X_n\big)+a_{n-1}(X_n,X_{n-1})X_{n-1}. 
$$
And so on and so forth. The skeptical or meticulous reader can formulate a proof by \recu in good and due form.
\end{proof}

\section{Dimension of geometric \risz}
\label{secDimGeom}

\vspace{.3em}
\subsec{\Pol \ris over a \cdiz}

A first important result in the theory of \ddk is the dimension of \pol \ris over a \cdiz.
\begin{theorem}
\label{propDimKXY}
If $\gK$ is a nontrivial \cdiz, the \ddk of the \pol \ri $\gK[X_1,\ldots,X_\ell]$ is equal to $\ell$.
\end{theorem}

We first establish the following result which needs a precise \dfnz.
Some \elts $x_1$, \ldots, $x_\ell$ of a \Klg with \zed $\gK$ are said to be \emph{\agqt dependent over $\gK$} if they annihilate a primitive \polz%
\footnote{The notion introduced here generalizes the notion of a primitively algebraic element introduced on \paref{subsecExtAdpC}. If $\gK$ were not zero-dimensional, it would be reasonable to use a more restrictive terminology such as
\gui{primitively \agq \rdez.} It is also clear that the \plgref{plcc.agq} can be \gne in the case of~$\ell$ \eltsz.}
 $f\in\gK[X_1,\ldots, X_\ell]$.

\begin{proposition}
\label{propAdepSing}
Let $\gK$ be a \cdiz, or more \gnlt a \zed \riz, $\gA$ be a \Klgz, and $x_1$, \ldots, $x_\ell \in \gA$ be \agqt dependent over $\gK$.
Then the sequence $(x_1,\dots,x_\ell)$ is singular.
\end{proposition}
\begin{proof}
We treat the case of a \cdiz. The \gnl case is then deduced by applying the \elgbm \num2 (\paref{MethodeZedRed}).
\\
Let $Q(x_1,\ldots,x_\ell)= 0$ be an \agq \rde over $\gK$.
Let us put a lexicographical order on the nonzero \moms $\alpha_{p_1,\ldots, p_\ell}\allowbreak x_1^{p_1}x_2^{p_2}\cdots x_\ell^{p_\ell}$ of $Q$, in accordance with the \gui{words} $p_1\,p_2\,\ldots \,p_\ell$.
We can assume that the \coe of the smallest nonzero \mom equal to $1$ (here we use the hypothesis that the field is discrete, because we assume that we can determine for each $\alpha_{p_1,\ldots ,p_\ell}$ wether it is null or \ivz).
Let $x_1^{m_1}x_2^{m_2}\cdots x_\ell^{m_\ell}$ be this \momz.
By following the lexicographical order, we see that we can express $Q$  in the form
$$\preskip-.0em \postskip.4em
\arraycolsep2pt\begin{array}{rcl}
Q& =   &  x_1^{m_1}\cdots x_\ell^{m_\ell}+
x_1^{m_1}\cdots x_\ell^{1+m_\ell}R_\ell  +x_1^{m_1}\cdots x_{\ell-1}^{1+m_{\ell-
1}}R_{\ell-1} \\[.8mm]
&   & +\cdots+ x_1^{m_1}x_2^{1+m_2}R_2+ x_1^{1+m_1}R_1                 
\end{array}$$
where $R_j\in \gK[x_k\,;\,k\geq j]$. Then $Q=0$ is the desired \egtz.
\end{proof}

\begin{Proof}{Proof of \thref{propDimKXY}. }
We first note that the sequence $(X_1,\ldots,X_\ell)$ is regular, which shows that the \ddk of $\gK[X_1,\ldots,X_\ell]$ is $\geq \ell$. We can also directly see that it is nonsingular: in \Egref{eqsing}   with $x_i=X_i$ the left-hand side is nonzero (consider the \coe of $X_1^{m_1}X_2^{m_2}\cdots \allowbreak X_\ell^{m_\ell}$).
To prove that the dimension of $\gK[X_1,\ldots,X_\ell]$  is $\leq \ell$,
it suffices, given Proposition~\ref{propAdepSing}, to show that $\ell+1$ \elts of $\gK[X_1,\ldots,X_\ell]$ are always \agqt dependent over $\gK$. Here is an \elr \dem of this classical result. \label{prvalgdep}
Let $y_1,\ldots ,y_{\ell+1}$ be these \eltsz, and $d$ be a bound on their degrees. For some integer $k\geq 0$ consider the list ${L}_k$ of all \hbox{the $y_1^{\delta_{1}}\cdots y_{\ell+1}^{\delta_{\ell+1}}$} such that $\sum_{i=1}^{\ell+1} \delta_i\leq k$. The number of \elts of the list ${L}_k$ is equal to ${k +\ell+1}\choose{k }$: this is a \pol of degree $\ell+1$ in $k$. The \elts of ${L}_k$ live in the \evc $E_{\ell,kd}$ of the \elts of $\gK[X_1,\ldots,X_\ell]$ of degree $\leq k \,d$, which is of dimension ${kd +\ell}\choose{kd}$: this is a \pol of degree~$\ell$ in~$k$. Thus for large enough $k$, the cardinal of ${L}_k$ is greater than the dimension of the \evc $E_{\ell,kd}$ containing $L_k$, 
therefore there is a \rdl between the \elts of  ${L}_k$. This provides an \agq \rde between the~$y_i$'s.
\end{Proof}

\comm   \label{remDKRR}
The proof of Proposition~\ref{propAdepSing} cannot \cot provide the same result for the field of reals $\RR$ (which \emph{is not} discrete). Actually it is impossible to realize the test of zero-dimensionality for $\RR$: \perso{l\'eg\`ere modif nov 06}
$$
\preskip.3em \postskip.3em 
\forall x\in \RR \;\;\exists a\in \RR\;\;\exists n\in \NN,\,
x^n\, (1-ax) = 0. 
$$
This would indeed mean that for every real number $x$, we know how to find a real $a$ such that $x(1-ax)=0$.
If we have found such an $a$, we obtain
\vspace{-.1em}
\begin{enumerate}\itemsep0pt
\item [--]  if $ax$ is \iv then $x$ is \ivz,
\item  [--] if  $1-ax$ is \iv then $x=0$.
\end{enumerate}
However, the alternative \gui{$ax$ or $1-ax$ is \ivz} is explicit over $\RR$.
Thus providing the test of zero-dimensionality amounts to the same as providing the test for \gui{is $x$ null or \ivz~?}
But this is not possible from the \cov point of view. 
Moreover, we can show that it is impossible to have a nonsingular sequence of length $1$, if we take $y\neq 0$ in the strong sense of \gui{$y$ is \ivz} (in the \dfn of \gui{nonsingular}).
Indeed, if we have some $x$ such that for all
$ a\in \RR $ and $ n\in \NN,$  $x^n\,(1-ax)\, \in \RR\eti $,
we get a contradiction: if $a=0$ then $x\in\RR\eti$, therefore there exists a  $b$ such that $1-bx=0$. 
\eoe

\vspace{-.3em}
\subsec{An interesting corollary}

\vspace{-.3em}
\begin{lemma}\label{lemahbonvraiment}
A \ri generated by $k$ \elts is of finite \ddkz.
\end{lemma}
\begin{proof}
Since the dimension can only decrease by passage to a quotient, it suffices to show that $\ZZXk$ is of \ddkz~$\leq2k+1$ (actually this \ri is of \ddk $k+1$ by \thref{corthValDim}). \\
Let $(h_1,\ldots,h_{2k+2})$ be a sequence of $2k+2$ \elts in $\ZZXk=\ZZuX$. We need to show that it is singular.\\
The sequence $(h_1,\ldots,h_{k+1})$ is singular in $\QQXk=\QQuX$. This means that the iterated boundary \id $\IK_\QQuX(h_1,\ldots,h_{k+1})$ contains $1$. 
\\
By {getting rid of the \denosz} we obtain that $\IK_\ZZuX(h_1,\ldots,h_{k+1})$ contains an integer~$d>0$.
Therefore the \ri $\gB=\ZZuX\big/{\IK_\ZZuX(h_1,\ldots,h_{k+1})}$ is a quotient of the \ri $\gC=(\aqo{\ZZ}{d})[\uX]$. As $\aqo{\ZZ}{d}$ is \zedz, the %
sequence $(h_{k+2},\ldots,h_{2k+2})$ is singular in $\gC$ (Proposition~\ref{propAdepSing}), in other words the \id $\IK_\gC(h_{k+2},\ldots,h_{2k+2})$ contains $1$. A fortiori $\IK_\gB(h_{k+2},\ldots,h_{2k+2})$ contains $1$. Finally, the \riz

\snic{\ZZuX\sur{\IK_\ZZuX(h_1,\ldots,h_{2k+2})}= \gB\sur{\IK_\gB(h_{k+2},\ldots,h_{2k+2})}}

is trivial.
\end{proof}
%

\subsec{Geometric \risz}\rdb

We call a \emph{\gmq \riz} a \ri $\gA$ that is a \pf \Klg with $\gK$ as a nontrivial \cdiz.\index{ring!geometric ---}

\Thref{thNstNoe} of \iNoe position affirms that such a quotient \ri is a finite integral extension of a \ri $\gB=\KYr$ contained in $\gA$ (here, $Y_1$, \dots, $Y_r$ are \elts of $\gA$ \agqt independent over $\gK$).

\begin{theorem}
\label{thDKAG}
Under the previous hypotheses, the \ddk of the \ri $\gA$ is equal to $r$.
\end{theorem}
\begin{proof}
\thref{propDimKXY} shows that $\Kdim\gB\leq r$. 
We can get the fact that $r+1$ \elts of $\gA$ are \agqt dependent over $\gK$ in the same style as described on \paref{prvalgdep} for a \poll \algz.
This will give $\Kdim\gA\leq r$.
For the \ddk to \hbox{be $\geq r$} results from Proposition~\ref{propDKEXENT}.\\
NB: \thref{cor2thKdimMor}, which implies $\Kdim\gA=\Kdim\gB$, gives another \demz.
\end{proof}

\section{Krull dimension of \trdisz}
\label{secDDKTRDIS}

As previously mentioned, the \ddk of a commutative \ri $\gA$ is none other than that \ddk of the spectral space $\SpecA$, at least in \clamaz.

In \coma we introduce the \ddk of a \trdi $\gT$ so that it is equal, in \clamaz, to the \ddk of its spectrum $\SpecT$.
The \dem of this \egt is very nearly identical to the one which we gave for  commutative \risz.
We will not repeat it, since in any case, we will always use the \ddk of a \trdi via the \cov \dfn that follows.

\pagebreak	        

\begin{definition}\label{defiDDKTRDI}~
\index{Krull dimension!of a distributive lattice}%
\index{complementary sequences!in a distributive lattice}
\begin{enumerate}
\item Two sequences $(\xzn)$ and $(\bzn)$ in a \trdi $\gT$ are said to be \emph{complementary} if
\begin{equation}\label{eqC2G}
\left.\arraycolsep3pt
\begin{array}{rcl}
 b_0\vi x_0& =  & 0    \\
 b_1\vi x_1& \leq  &  b_0\vu x_0  \\
\vdots~~~~& \vdots  &~~~~  \vdots \\
 b_n\vi  x_n & \leq  &  b_{n -1}\vu x_{n -1}  \\
 1& =  &   b_n\vu x_n
\end{array}
\right\}
\end{equation}
A sequence that has a \cop sequence will be said to be \emph{singular}.%
\index{singular!sequence}\index{sequence!singular ---}
\item For $n\geq0$ we will say that the \trdi $\gT$ is \emph{of \ddk at most $n$} if every sequence $(\xzn)$ in $\gT$ is singular.
Moreover, we will say that the \trdi $\gT$ is of \ddk $-1$ if it is trivial, \cad if $1_\gT=0_\gT$.
\end{enumerate}
\end{definition}
\hmodeHabillage{\hbox{\xyrowsp=3pt
$\SCO{x_0}{x_1}{x_2}{b_0}{b_1}{b_2}$}}{0}{-15pt}%
For example, for $k=2$ item \emph{1} corresponds to the following graph in $\gT$.
\\
We will write $\Kdim\gT\leq n$ when the \ddk is at most $n$.
\label{NOTAKdimTrdi}\\
It is obvious that a lattice has the same \ddk as the opposite lattice.
We also immediately see that a lattice is \zed \ssi it is a \agBz.
\\
Also, a totally ordered set of $n$ \elts has for \ddkz~$n-2$.
\endHabillage

\begin{fact}\label{factDDKTRDI}
Let $S$ be a subset of $\gT$ that generates $\gT$ as a \trdiz. Then
$\gT$ is \ddi$n$ \ssi every sequence $(\xzn)$ in $S$ admits a \cop sequence in $\gT$.
\end{fact}
\begin{proof}
Let us illustrate the computations on a sufficiently \gnl example in the case~$n=4$. 
\\
We verify that if $(x_0,x_1,x_2,x_3,x_4)$ admits $(a_0,a_1,a_2,a_3,a_4)$ as a \cop sequence, and if  $(x_0,x_1,y_2,x_3,x_4)$ admits $(b_0,b_1,b_2,b_3,b_4)$ as a \cop sequence, then the sequence $(x_0,x_1,x_2\vu y_2,x_3,x_4)$ admits the \cop sequence~$(a_0\vu b_0,a_1\vu b_1,a_2\vi b_2,a_3\vi b_3,a_4\vi b_4)$. 
\\
Dually, the sequence $(x_0,x_1,x_2\vi y_2,x_3,x_4)$ admits the \cop sequence~$(a_0\vu b_0,a_1\vu b_1,a_2\vu b_2,a_3\vi b_3,a_4\vi b_4)$.
\\ 
The same computation would work for an arbitrary $x_i$ (instead of $x_2$ above) in an arbitrary finite sequence. Thus if each sequence $(z_0,\ldots,z_n)$ in $S$ admits a \cop sequence in $\gT$, the same will hold for every sequence of $n+1$ terms in the lattice generated by $S$.
\end{proof}
%

\begin{fact}\label{corfactDDKTRDI}
A commutative \ri has the same \ddk as its Zariski lattice.
\end{fact}
\begin{proof}
The \demz, based on Fact~\ref{factDDKTRDI}, is left to the reader. 
Another \dem will be given later in the form of Lemma~\ref{lemsutfidele}.
\end{proof}

We can also access the \ddk via the Krull boundary \ids as for the commutative \risz.

\goodbreak
\begin{definition}\label{defiBordKrullTrdi}~
\begin{enumerate}
\item The lattice $\gT\ul x=\gT/(\JK_\gT(x)=0)$, where
\begin{equation}\label{eqbordsupTrdi}
\JK_\gT(x)\,=\,\dar x \,\vu\, (0:x)_\gT
\end{equation}
is called \emph{the (Krull) upper boundary of $x$ in $\gT$}. We also say that the \id $\JK_\gT(x)$ is \emph{the Krull boundary \id of $x$ in $\gT$.} \index{ideal!Krull boundary --- (\trdiz)}
\index{Krull upper boundary!(\trdiz)}
\item More \gnltz, for a sequence $(\ux)$ in $\gT$, the iterated Krull boundary \id $\JK_\gT(\ux)$ is defined by \recu as follows: $\JK_\gT()=\so{0},$ and
\begin{equation}\label{eq1IdBordKrullItereTrdi}
\JK_\gT(x_0, \ldots, x_k) =
\big(\JK_\gT(x_0, \ldots, x_{k-1}) : x_k\big)_\gT \vu \dar x_k.%
\index{ideal!iterated Krull boundary --- (\trdiz)}
\end{equation}
\end{enumerate}
\end{definition}

\begin{fact}\label{factDdkTrdiBord}
Let   $n\in\NN$  and $\gT$ be a \trdiz. 
\begin{enumerate}
\item A sequence $(\xzn)$ in $\gT$ is singular \ssi the iterated boundary \id $\JK_\gT(\xzn)$ contains $1$.
\item We have $\Kdim\gT\leq n$ \ssi for every $x$, 
$\Kdim \gT\ul x\leq n-1$. 
\end{enumerate}
\end{fact}

\begin{fact}\label{propDDKagH}
In a \agHz, every iterated Krull boundary \id is principal: $\JK_\gT(x)=\dar{\big(x\vu\lnot x\big)}$ and more \gnltz,
\begin{equation}
   \!\!\!  \JK_\gT(\xzn) =\dar{\big(x_{n}\vu \big(x_{n}\im( \cdots (x_1 \vu (x_1 \im (x_0\vu \neg x_0)))\cdots)\big)\big)}
\end{equation}
 \end{fact}
%
%

\begin{lemma}\label{lemZarAHeyt}
Let $\fa$ and $\fb$ be two \itfs of a \ri $\gA$. \\
In the lattice $\ZarA$, the \elt $\DA(\fa)\im\DA(\fb)$ exists \ssi the \idz~$(\fb:\fa^{\infty})$  has the same radical as a \itfz.
\end{lemma}
%
\begin{proof}
In a \trdiz, the \elt $u\im v$  exists if the \id $(v:u)$ is principal (its \gtr is then denoted by $u\im v$). 
However, for some \tf \id $\fa$,  $\big(\DA(\fb):\DA(\fa)\big)=\DA(\fb:\fa^{\infty})$. 
Hence the stated result.
\end{proof}
%
\begin{lemma}\label{lem2ZarAHeyt}
Suppose that $\ZarA$ is a \agHz. \\
For $(\xzn)$ in $\gA$, we have the \egt 

\snic{\DA\big(\JK_\gA(\xzn)\big)=\JK_\ZarA\big(\DA(x_0),\ldots,\DA(x_n)\big) .}
\end{lemma}
\facile 

\begin{proposition}\label{propNoetAgH}
Let $\gA$ be a \noe \coriz.
\begin{enumerate}
\item If $\fa$ and $\fb$ are two \itfsz,  the \id $(\fb:\fa^{\infty})$ is \tfz.
\item $\ZarA$ is a \agHz, with $\DA(\fa)\im\DA(\fb)=\DA(\fb:\fa^{\infty})$.  
\item The iterated Krull boundary \ids defined on \paref{eq1IdBordKrullItere} have the same radical as the \itfsz.
\item If in addition $\gA$ is  \fdiz, $\ZarA$ is discrete and we dispose of a test to decide if a sequence in
$\gA$ admits a \cop sequence.
\end{enumerate}
\end{proposition}

\begin{proof} 
\emph{1.} Let  $\fa$, $\fb\in \Zar\gA$.  
Let $\fJ_k=(\fb:\fa^k)$.
Since $\gA$ is \cohz, each \id $\fJ_k$ is \tfz. Since $\gA$ is \noez, the sequence admits two consecutive equal terms, for example of indices $p$ and $p+1$, from which it is clear that it becomes stationary.
We then have $\fJ_p=(\fb:\fa^{\infty})$.
\\
\emph{2.} Consequence of \emph{1} given Lemma~\ref{lemZarAHeyt}.
\\
\emph{3.} Results by \recu of \emph{2} given Fact~\ref{propDDKagH} and Lemma~\ref{lem2ZarAHeyt}.
\\
\emph{4.}
Results from \emph{2}, from Fact~\ref{propDDKagH}  and from Lemma~\ref{lem2ZarAHeyt}.
\end{proof}

\section{Dimension of morphisms}
\label{secKdimMor}

\subsec{Definition and first \prtsz}

\begin{definition}\label{defiKdimMor}
Let $\rho:\gA\to\gB$ is a \ri \homoz. \emph{The \ddk of the morphism} $\rho$
is  the \ddk of the \ri $\Abul\otimes_\gA\!\gB$ obtained by \eds (\thref{factSDIRKlg}) from $\gA$ to its \zedr closure $\Abul$ (\thref{thZedGen}).
\end{definition}

\exls ~\\
1) If $\gk$ is \zedz,  we have seen that $\Kdim\kXn\leq n$. We deduce that the \ddk of the morphism $\gA\to\AXn$ is~$\leq n$, with \egt if $\gA$ is nontrivial.
\\
2) If $\gB$ is an integral \Algz, after \eds the \alg is integral over $\Abul$, therefore \zedez. Thus, the morphism $\gA\to\gB$ is \zedz.
\eoe

\begin{lemma}\label{lem1KdimMor}
Let $\gB$  and $\gC$ be two \Algsz. 
Then by \eds we obtain $\Kdim(\gC\to\gC\otimes_\gA\!\gB)\leq \Kdim(\gA\to\gB)$ in the following cases.
\begin{enumerate}
\item $\gC$ is a quotient of $\gA$, or a localized \ri of $\gA$, or the quotient of a localized \ri of $\gA$.
\item  $\gC$ is a finite product of \ris of the previous type.
\item  $\gC$ is a filtering colimit of \ris of the previous type.  
\end{enumerate}

\end{lemma}
\begin{proof} We use the observation $\gC\bul\otimes_\gC(\gC\otimes_\gA\gB)\simeq\gC\bul\otimes_\gA\gB \simeq \gC\bul\otimes_{\Abul}(\Abul\otimes_\gA\gB)$.
We then prove that the functor $\gB\mapsto\gB\bul$ transforms a quotient into a quotient, a localized \ri into a localized \ri (Proposition~\ref{propClZdrLoc}),
a finite product into a finite product, and a filtering colimit into a filtering colimit. Moreover, the \eds also commutes with all these constructions. Finally, the Krull dimension can only decrease by these constructions.
\end{proof}

\rem{It is not true that $\gC\otimes_{\gA}\gB$ is \zed as soon as the three \ris are \zedz.
For example we can take $\gA$ to be a \cdi and~$\gB=\gC=\gA(X)$.
Then $\Kdim(\gC\otimes_\gA\gB)=1$ (see Exercise \ref{exoKdimSomTr}). 
It follows that the \edsz, even in the case of a \fpte extension, can strictly increase the Krull dimension of morphisms.
A contrario we have the following \plgcz. 
\eoe}

\begin{plcc}\label{plcc.KdimMor}\relax
Let $S_1$, $\ldots$, $S_n$  be \moco of a \ri $\gA$, 
$k\geq-1$ be an integer and $\gB$ be an \Algz.
The \ddk of the morphism $\gA\to\gB$ is $\leq k$ \ssi 
the \ddk of each of the morphisms $\gA_{S_i}\to\gB_{S_i}$ is $\leq k$.  
\end{plcc}
%
\begin{proof}
As $(\gA_{S_i})\bul \simeq (\gA\bul)_{S_i}$
(Proposition~\ref{propClZdrLoc}), we obtain

\snic{(\gA_{S_i})\bul\otimes _{\gA_{S_i}}\gB_{S_i}
\simeq(\gA\bul\otimes_{\gA}\gB)_{S_i},}


and we are brought back to the \plgref{thDdkLoc}.
\end{proof}

\medskip 
The goal of this section is to show, for a morphism $\rho:\gA\to\gB$, the in\egt
$$
\fbox{$1+\Kdim\gB\,\leq\, (1+\Kdim\gA)(1+\Kdim\rho)$}.
$$

Note that for $\Kdim\gA\leq0$ we trivially have $\Kdim\gB=\Kdim\rho$. 
We then treat a simple but nontrivial case to get a clear picture.
The truly simple case would be the one where $\gA$ is integral and $\Kdim\gA\leq1$.
As the \dem is unchanged, we will only suppose that $\gA$ is a \qiriz, which will make the rest easier.

\begin{proposition}\label{prop1KdimMor} 
Let  $\rho:\gA\to\gB$ be a morphism, with $\gA$ a \qiriz.\\
If $\Kdim \rho \leq n$ and $\Kdim\gA\leq 1$, then $\Kdim\gB\leq 2n+1$.
\end{proposition}
\begin{proof}
Let $\uh=(h_0,\ldots,h_{2n+1})$ be a sequence of $2n+2$ \elts in $\gB$. We need to show that it is singular.\\
By hypothesis the \ri $\Abul\otimes_\gA\gB$ is \ddi$n$.\\
Let $\gK=\Frac\gA$ be the total \ri of fractions. It is \zedr and generated by $\gA$ as a \zedr \riz, therefore it is a quotient of $\Abul$.
We conclude that the sequence $(h_0,\ldots,h_{n})$ is singular in~$\wi\gB=\gK\otimes_\gA\!\gB$. \\
This means that the iterated boundary \id $\IK_{\wi\gB}(h_0,\ldots,h_{n})$ contains $1$, and {by getting rid of the \denosz}, that $\IK_{\gB}(h_0,\ldots,h_{n})$ contains some $a\in\Reg(\gA)$.\\
Therefore  $\gB_0=\gB\sur{\IK_\gB(h_0,\ldots,h_{n})}$ is a quotient of $\gB\sur{a\gB}=\gA\sur{a\gA}\otimes_\gA\!\gB$.
Since $a$ is \ndz and $\Kdim\gA\leq1$, the quotient $\gA\sur{a\gA}$ is \zedz, so $(\gA\sur{a\gA})\red$ is a quotient of $\Abul$ and the \ri $(\gB_0)\red$ is a quotient of  $\Abul\otimes_\gA\!\gB$.
We deduce that the sequence $(h_{n+1},\ldots,h_{2n+1})$ is singular in $(\gB_0)\red$, therefore also in $\gB_0$.
\\
Therefore the \ri $\gB\sur{\IK_{\gB}(\uh)}=\gB_0\sur{\IK_{\gB_0}(h_{n+1},\ldots,h_{2n+1})}$ is  trivial.
\end{proof}

To pass from the \qiri case to the \gnl case we want to say that every reduced \ri can behave in the computations like an integral \ri provided we replace $\gA$ with
$$
\gA\sur{\Ann_\gA (a) }\times \gA\sur{\Ann_\gA(\Ann_\gA( a) )}
$$
when an \algo asks to know if the annihilator of $a$ is equal to $0$ or~$1$.
The important thing in this construction is that the closed covering principle for the \susis applies since the product of the two \ids
$\Ann_\gA (a)$
and $\Ann_\gA(\Ann_\gA (a) )$ is null.

This type of  \dem will probably be easier to grasp when we will 
familiarize ourselves with the basic local-global machinery explained on
\paref{MethodeIdeps}.\imlb
Here we do not proceed by successive \come \lons but by successive \gui{closed coverings.}

\smallskip
Actually we will not introduce a dynamic computation tree as such, we will instead construct a universal object. This universal object is a \gui{\cov finitary approximation} of the product of all the quotients of $\gA$ by its \idemisz, a slightly too ideal object of \clama to be considered \cotz, at least in the form that we just defined: actually, if~$\gB$ is this product and if $\gA_1$ is the natural image of $\gA$ in $\gB$, then the \uvl \ri that we are constructing should be equal to the \qiri closure of~$\gA_1$ in $\gB$, at least in \clamaz.

\subsec{The minimal \qiri closure of a reduced \riz}

In what follows we denote by $a\epr$ the annihilator \id of the \elt $a$ when the context is clear (here the context is simply the \ri in which we must consider~$a$). We will also use the notation $\fa\epr$ for the annihilator of an \idz~$\fa$.

The facts stated below are \imdsz.
\begin{eqnarray}
\fa & \subseteq & (\fa\epr)\epr \label{ann1}\\
\fa\subseteq\fb & \Longrightarrow & \fb\epr \subseteq \fa\epr\label{ann2}\\
\fa\epr & = & \big((\fa\epr)\epr\big)\epr\label{ann3}\\
(\fa+\fb)\epr & = & \fa\epr \cap \fb\epr\label{ann4}\\
\fa\epr\subseteq\fb\epr & \Longleftrightarrow &(\fa+\fb)\epr=\fa\epr
\label{ann5}\\
\fa\epr\subseteq\fb\epr & \Longleftrightarrow & (\fb\epr)\epr \subseteq (\fa\epr)\epr\label{ann6}\\
(\fa\epr:\fb)&=& (\fa\fb)\epr \label{ann7}
\\
(\gA\sur{\fa\epr})\big/{{~\ov\fb~}\epr}& = &  \gA\sur{(\fa\fb)\epr} \label{ann8}
\end{eqnarray}

\rems 1) An \id $\fa$ is an annihilator (of another \idz) \ssi $\fa=(\fa\epr)\epr$.
\\
2) The inclusion $\fa\epr+\fb\epr \subseteq (\fa\cap\fb)\epr$ can be strict, even if $\fa=\fa_1\epr$ and~$\fb=\fb_1\epr$.
Take for example $\gA=\ZZ[x,y]=\aqo{\ZZ[X,Y]}{XY}$, $\fa_1=\gen{x}$ and $\fb_1=\gen{y}$. Then, $\fa=\fa_1\epr=\gen{y}$,  $\fb=\fb_1\epr=\gen{x}$, $\fa\epr+\fb\epr=\gen{x,y}$, and~$(\fa\cap\fb)\epr=\gen{0}\epr=\gen{1}$.\eoe

\medskip
If we assume that $\gA$ is reduced, we also have the following results.
\begin{eqnarray}
\sqrt{\fa\epr}\quad=\quad\fa\epr & = & (\sqrt{\fa})\epr\quad=\quad(\fa^2)\epr \label{ann9}\\
 (\fa\fb)\epr & = &  (\fa\cap \fb)\epr \label{ann10}
\\
\fa\epr\subseteq\fb\epr & \Longleftrightarrow & (\fa\fb)\epr = \fb\epr
\label{ann11}
\end{eqnarray}

\begin{lemma}\label{lem20MorRc}
Let $\gA$ be a reduced \ri and $a\in\gA$.
We define
$$\gA_{\so{a}}\eqdefi\gA\sur{a\epr}\times \gA\sur{({a\epr})\epr}$$
and we let $\psi_a:\gA\to\gA_{\so{a}}$ be the canonical \homoz. 
\begin{enumerate}
\item $\psi_a(a)\epr$ is generated by the \idm $(\ov 0,\wi 1)$, so $\psi_a(a)\epr=(\ov 1,\wi 0)\epr$.
\item $\psi_a$ is injective (we can identify $\gA$ with a sub\ri of $\gA_{\so{a}}$).
\item Let $\fb$ be an \id  in $\gA_{\so{a}}$, then the \id $\psi_a^{-1}(\fb\epr)=\fb\epr\cap\gA$ is an annihilator \id in $\gA$.
\item The \ri $\gA_{\so{a}}$ is reduced.
\end{enumerate}
\end{lemma}
\begin{proof} \emph{1.}  We have $\psi_a(a)=(\ov a,\wi 0)$, where $\ov x$ is the class modulo $a\epr$ and $\wi x$ is the class modulo $(a\epr)\epr$. If $c=(\ov y, \wi z)$, the \egt $\psi_a(a)c=0$ means $\ov{ya}=0$,
\cad $ya^2=0$, or yet $ya=0$, \cad $\ov y=\ov 0$.
\\
\emph{2.} If $xa=0$ and $xy=0$ for every $y\in a\epr$ then $x^2=0$ so $x=0$.
\\
\emph{3.} Let $\psi_1:\gA \to\gA\sur{a\epr}$ and $\psi_2:\gA \to\gA\sur{(a\epr)\epr}$ be the two \prnsz.
We have~$\fb=\fb_1\times \fb_2$. If $x\in\gA$ we have

\snic{\psi_a(x)\in\fb\epr \;\Longleftrightarrow\;\psi_1(x)\fb_1=0 \et \psi_2(x)\fb_2=0,}

\cad  $x\in\psi_1^{-1}(\fb_1\epr)
\cap \psi_2^{-1}(\fb_2\epr)$. \Egrf{ann8} tells us that each $\psi_i^{-1}(\fb_i\epr)$ is an annihilator \idz.
\Trf by \Egrf{ann4}.
\\
\emph{4.} In a reduced \riz, every annihilator \id $\fb\epr$ is radical: indeed, if~$x^2\fb=0$, then $x\fb=0$. 
\end{proof}

\vspace{-.6em}
\pagebreak	        

\begin{lemma}\label{lem3MorRc}
Let $\gA$ be reduced and $a$, $b\in\gA$. Then with the notations of Lemma~\ref{lem20MorRc} the two \ris $(\gA_{\so{a}})_{\so{b}}$ and $(\gA_{\so{b}})_{\so{a}}$ are canonically \isocz.
\end{lemma}
\begin{proof}
The \ri $(\gA_{\so{a}})_{\so{b}}$ can be described symmetrically as follows

\snic{\gA_{\so{a,b}}=\gA\sur{(ab)\epr}\times \gA\sur{(ab\epr)\epr}\times \gA\sur{(a\epr b)\epr}\times \gA\sur{(a\epr b\epr)\epr},}

and if $\psi:\gA\to\gA_{\so{a,b}}$ is the canonical \homoz, it is clear that we have~$\psi(a)\epr=(1,1,0,0)\epr$ and $\psi(b)\epr=(1,0,1,0)\epr$.
\end{proof}

\rem The case where $ab=0$ is typical: when we meet it, we would like to split the \ri into components where things are \gui{clear.}
The previous construction then gives the three components

\snic{\gA\sur{(ab\epr)\epr}, \; \gA\sur{(a\epr b)\epr}\, \hbox{ and } \,\gA\sur{(a\epr b\epr)\epr}.}

In the first one, $a$ is \ndz and $b=0$, in the second one $b$ is \ndz and $a=0$, and in the third one $a=b=0$.
\eoe

\medskip
The following lemma regarding  \qiris is copied from Lemma~\ref{lem2SousZedRed} which concerned  \zedr \ris (the reader will also be able to just about copy the \demz).

\begin{lemma}\label{lem2qi}
If $\gA\subseteq\gC$ with $\gC$ a \qiriz, the smallest \sqiri of $\gC$ containing $\gA$ is equal to $\gA[(e_a)_{a\in\gA}]$,
where $e_a$ is the \idm of $\gC$ such that $\Ann_\gC(a)=\gen{1-e_a}_\gC$. More \gnlt if $\gA\subseteq\gB$
with reduced $\gB$, and if every \elt $a$ of $\gA$ admits an annihilator in $\gB$ generated by an \idm $1-e_a$,
then the sub\ri $\gA[(e_a)_{a\in\gA}]$ of $\gB$ is a \qiriz.
\end{lemma}

\begin{thdef}\label{thAmin} \emph{(Minimal \qiri closure)}
\\
Let $\gA$ be a reduced \riz.
We can define a \ri $\Amin$ as a filtering colimit by iterating the basic construction which consists in replacing~$\gE$ (the \gui{current} \riz, which contains $\gA$) by

\snic{\gE_{\so{a}}\eqdefi\gE\sur{a\epr}\times \gE\sur{({a\epr})\epr}=\gE\sur{\Ann_\gE (a) }\times \gE\sur{\Ann_\gE(\Ann_\gE (a) )},}

 when $a$ ranges over $\gA$.
\begin{enumerate}
\item  This \ri $\Amin$ is a \qiriz, contains $\gA$ and is integral over $\gA$.
\item For all $x\in\Amin$,
$
x\epr\cap\gA$ is an annihilator \id in $\gA$.
\end{enumerate}
This \ri $\Amin$ is called the \emph{minimal \qiri closure of $\gA$}.
 \\
When $\gA$ is not \ncrt reduced, we will take $\gA\qim\eqdefi (\Ared)\qim$.
\end{thdef}
\index{closure!minimal \qiri ---}
\begin{proof} \emph{1.} By Lemma~\ref{lem2qi}, it suffices to add an \idm $e_a$ for each $a\in\gA$ to obtain a \qiriz. The colimit is well-defined due to the relation of commutation given by Lemma~\ref{lem3MorRc}.
\\
For item~\emph{2} note that $x$ is obtained at a finite stage of the construction, and that $x\epr\cap\gA$ stops changing from the moment where $x$ is reached because the successive \homos are injections. We can therefore call upon item \emph{3} of Lemma~\ref{lem20MorRc}.
\end{proof}
\rem
We can ask ourselves if $\Amin$ could not be \care by a \uvl \prt related to item \emph{2.}
\eoe

\medskip 
By \gui{iterating} the description of  $(\gA_{\so{a}})_{\so{b}}$ given in the \dem of Lemma~\ref{lem3MorRc} we obtain the following description of each \ri obtained at a finite stage of the construction of $\Amin$
(see Exercise~\ref{exoAminEtagesFinis}).

\begin{lemma}\label{lem4MorRc}
Let $\gA$ be a reduced \ri and $(\ua) = (\an)$ be a sequence of~$n$ \elts of $\gA$.  For $I\in\cP_n$, let $\fa_I$ be the \id

\snic {
\fa_I = \big(\prod_{i\in I} \gen{a_i}\epr \prod_{j\notin I} a_j\big)\epr
= \big(\gen{a_i, i \in I}\epr \prod_{j\notin I} a_j\big)\epr
.}

Then $\Amin$ contains the following \riz, a product of $2^n$ quotient \ris of~$\gA$ (some eventually null)

\snic {
\gA_{\so\ua} = \prod_{I\in\cP_n} \gA\sur{\fa_I}.
} 
\end{lemma}

\begin{fact}\label{fact2Amin}~
\begin{enumerate}
\item Let $\gA$ be a \qiriz. 
\begin{enumerate}
\item $\Amin=\gA$.
\item $\AX$ is a \qiriz, and $\BB(\gA)=\BB(\AX)$.
%
%
\end{enumerate}

\item For every \ri $\gA$ we have a canonical \iso

\snic{\Amin[\Xn]\simeq(\AXn)\qim.}
\end{enumerate}
\end{fact}
\begin{proof} \emph{1a.} Results from the construction of $\Amin.$
\\
\emph{1b.} The result is clear for integral \risz.
We can apply the \elgbmd \num1 (\paref{MethodeQI}).
We could also use McCoy's lemma, Corollary~\ref{corlemdArtin}~\emph{2.}
\\
\emph{2.} We suppose \spdg that $\gA$ is a reduced \riz.
It also suffices to treat the case of a single variable.
Given Lemma~\ref{lem3MorRc} we can \gui{start} the construction of $\AX\qim$ with the constructions $\gE\leadsto \gE_{\so{a}}$ for some $a\in\gA$.
But if $\gE=\gB[X]$ and $a\in\gA\subseteq\gB$ then $\gE_{\so{a}}=\gB_{\so{a}}[X]$.
Thus $\Amin[X]$ can be seen as a first step of the construction of $\AX\qim$. But since by item \emph{1} $\Amin[X]$ is a \qiri and that for a \qiri $\gC$ we have $\gC=\gC\qim$, the construction of $\AX\qim$ is completed with $\Amin[X]$.
\end{proof}

\comm
In practice, to use the \ri $\Amin$, we only need the finite stages of the construction.
We can however note that even a single stage of the construction is  a little mysterious, insofar as the \ids $a\epr$ and $(a\epr)\epr$ are difficult to handle. It is only in the case of \coris that we know how to describe them with finite \sgrsz.
Actually if the \ri is \noez, the construction must end in a finite number of steps (at least from the point of view of \clamaz), and needs to replace the \ri with the product of its quotients with the \idemisz. Here we are in a situation where the construction of $\Amin$ meeting the standards of \coma seems more complicated than the result in \clama (at least if the \ri is \noez). 
Nevertheless, since we do not need to know the \idemisz, our method is more \gnl (it does not need \TEMz).
In addition, its complication is mostly apparent. When we use $\gA\sur{a\epr}$ for example, we actually make computations in $\gA$ by forcing $a$ to be \ndzz, \cad by forcefully annihilating every $x$ that presents itself and that annihilates~$a$.
When we use $\gA\sur{(a\epr)\epr}$, things are less easy, because a priori, we need a proof (and not simply the result of a computation) to certify that an \elt $x$ is in $(a\epr)\epr$.
It is a fact that the use of \idemis in a proof of \clama can in \gnl be made innocuous (\cad \covz) by using $\Amin$ (or another \uvl \ri of the same type\footnote{$\Amin$ corresponds to using all the quotients by the \idemisz, $\Frac(\Amin)$ corresponds to using all the quotient fields of these quotients.}),
even if we do not dispose of other means to \gui{describe an \id $(a\epr)\epr$} than the one of applying the \dfnz.
\eoe

\subsec{Application}

\begin{corollary}\label{corthAmin}
Let  $\rho:\gA\to\gB$ be a morphism of finite \ddkz. We \gui{extend the scalars} from $\gA$ to $\Aqim$: we obtain
$\gB'=\Aqim\otimes_\gA\!\gB$  and let $\rho':\Aqim\to\gB'$ be the natural morphism.
\\
Then $\Kdim\Aqim=\Kdim\gA$, $\Kdim\gB=\Kdim\gB'$ and $\Kdim\rho'\leq\Kdim\rho$.
\end{corollary}
\begin{proof}
The first two items result from the fact that in the construction of the \ri $\Aqim$,
at each \elr step
$$
\gE \quad \rightsquigarrow
\quad\gE\sur{\Ann_\gE (a) }\times \gE\sur{\Ann_\gE(\Ann_\gE a )},
$$
the product of the two \ids is null, which is found again after tensorization by~$\gB$. Therefore, the closed covering principle for the \ddkz~\vref{thDdkRecFer} applies. 
Finally, the in\egt $\Kdim\rho'\leq\Kdim\rho$ results from Lemma~\ref{lem1KdimMor}.
\end{proof}

\rem \Gnlt the \ri $\Frac\Aqim$ seems a better concept than~$\Abul$ to replace the quotient field in the case of a non-integral reduced \riz.
In the case where $\gA$ is a \qiriz, we indeed have $\Aqim=\gA$, so~$\Frac\Aqim=\Frac\gA$, while $\Abul$ is in \gnl
significantly more cumbersome (as the example $\gA=\ZZ$ shows).
\eoe

\begin{corollary}\label{prop2KdimMor} 
Let  $\rho:\gA\to\gB$ be a morphism. \\
If $\Kdim \rho \leq n$ and~$\Kdim\gA\leq 1$, then
$\Kdim\gB\leq 2n+1$.
\end{corollary}
\begin{proof}
This clearly results from Proposition~\ref{prop1KdimMor} and from Corollary \ref{corthAmin}.
\end{proof}
%
\entrenous{\'Etrange que seule la dimension of $(\Aqim)\bul\otimes_{\gA}\gB$
soit en fin of compte pertinente ici? Non, ce que dit cette preuve \gui{raffin\'ee}
du corollaire \ref{prop2KdimMor}, c'est que l'on  pourrait consid\'erer d'une part la longueur maximale
des cha\^{\i}nes of premiers in $\gB$ au dessus d'un premier minimal of $\gA$, d'autre part
 la longueur maximale
des cha\^{\i}nes of premiers in $\gB$ au dessus d'un premier maximal of $\gA$
et faire la somme des deux.}

\begin{theorem}\label{thKdimMor}
Let $\rho:\gA\to\gB$ be a morphism.\\ 
If $\Kdim \rho \leq n$ and $\Kdim\gA\leq m$, then
$\Kdim\gB\leq mn+m+n$.
 \end{theorem}
%
\begin{proof}
We perform a \dem by \recu on $m$. The case $m=0$ is trivial.
The proof given for $m=1$ in the case where $\gA$ is a \qiri 
(Proposition~\ref{prop1KdimMor}), which 
relied on the dimension $0$ case to prove the result in dimension $m = 1$, can easily be adapted to pass from dimension $m$ to dimension $m + 1$.
 We copy the \dem in the case where~$\gA$ is a \qiriz.\\
To pass to the case of an arbitrary \ri we use Corollary~\ref{corthAmin}.\\
We therefore suppose that $\gA$ is a \qiri and we consider a sequence $(\uh)=(h_0,\ldots,h_{p})$ in $\gB$ with  $p=(m+1)(n+1)-1$. We need to show that it is singular.\\
By hypothesis the \ri $\Abul\otimes_\gA\!\gB$ is \ddi$n$.
The total \ri of fractions $\gK=\Frac\gA$  is \zedrz, and it is generated by~$\gA$ as a \zedr \riz, so it is a quotient of~$\Abul$.
We conclude that the sequence $(h_0,\ldots,h_{n})$ is singular in the \riz~$\wi\gB=\gK\otimes_\gA\!\gB$. \\
This means that the iterated boundary \id $\IK_{\wi\gB}(h_0,\ldots,h_{n})$ contains $1$, and {by getting rid of the \denosz} that $\IK_\gB(h_0,\ldots,h_{n})$ contains some $a\in\Reg(\gA)$.
Therefore the \ri $\gB_0=\gB\sur{\IK_\gB(h_0,\ldots,h_{n})}$ is a quotient of~$\gB\sur{a\gB}=\gA\sur{a\gA}\otimes_\gA\!\gB$.
Since $a$ is \ndz and $\Kdim\gA\leq m$, the quotient $\gA\sur{a\gA}$ is \ddi$m-1$. The natural \homo $\gA\sur{a\gA}\to\gB\sur{a\gB}$ remains \ddi$n$ (Lemma~\ref{lem1KdimMor}). Therefore, by \hdrz, the sequence $(h_{n+1},\ldots,h_{p})$ is singular in $\gB\sur{a\gB}$.
Therefore the sequence~$(h_{n+1},\ldots,h_{p})$ is singular in $\gB_0$.
\\
In conclusion, the \ri $\gB\sur{\IK_{\gB}(\uh)}= \gB_0\sur{\IK_{\gB_0}(h_{n+1},\ldots,h_{p})}$ is  trivial.
\end{proof}

\begin{corollary}\label{corthKdimMor}
Suppose $\Kdim\gA\leq m$.  Then

\snic{\Kdim\AXn\leq mn+m+n.}
\end{corollary}
\begin{proof}
We know that if $\gK$ is \zedrz, $\Kdim\KXn\leq n$.
Thus $\Kdim(\gA\to\AXn)\eqdefi\Kdim\Abul[\Xn]\leq n$. We apply \thref{thKdimMor}.
\end{proof}
We dispose equally of a lower bound of $\Kdim \AXn$.
\begin{lemma}\label{lemKdimAxn}
For every nontrivial \ri $\gA$ and all $n>0$ we have 

\snic{n+\Kdim\gA\leq \Kdim\AXn.}

More \prmtz, the following implication is satisfied for every $k\geq -1$ and for every \riz

\snic{\Kdim\AXn\leq n+k\;\Longrightarrow\;\Kdim\gA\leq k}
\end{lemma}
%
\begin{proof} Immediate consequence of Proposition~\ref{lemRegsing}.
\end{proof}

\vspace{-.7em}
\pagebreak	        

\begin{theorem}\label{cor2thKdimMor}
Consider an \alg $\rho:\gA\to\gB$.
\begin{enumerate}
\item Suppose that $\gB$ is generated by \elts which are primitively  \agq over $\gA$,
then $\Kdim \rho\leq0$ and so $\Kdim\gB\leq\Kdim\gA$.
\item If $\rho$ is injective and $\gB$ is integral over $\gA$, then $\Kdim\gB=\Kdim\gA$.
\end{enumerate}
 \end{theorem}
\begin{proof}
\emph{1.} Given Fact~\ref{factZedPrimFin}, the \ri $\Abul\otimes_\gA\gB$ is \zedz,
in other words $\Kdim\rho\leq0$. \Trf by \thref{thKdimMor}.
\\
\emph{2.} By item \emph{1} and Proposition~\ref{propDKEXENT}.
\end{proof}
For a more direct \dem of the in\egt $\Kdim\gB\leq\Kdim\gA$, see Exercise~\ref{exoInclusionBordLionel}.


\section{Valuative dimension}
\label{secValdim}

\vspace{3pt}
\subsec{Dimension of \advsz}

Recall that a \adv is a reduced \ri in which we have, for all $a,\,b$: $a$ divides $b$ or $b$ divides $a$. In other words it is a Bézout and reduced \aloz.
A \adv is a normal, local \ri  \sdzz.
It is integral \ssi it is \cohz.

It is clear that the Zariski lattice of a \adv is a totally ordered set.
\begin{fact}\label{fact1SeqSingTD}
In a \trdi if a subsequence of $(\ux)=(\xn)$ is singular,
the sequence $(\ux)$ is singular.
\end{fact}
\begin{proof}
We consider a \susi $(\yr)$, with a \cop sequence $(b_1,\ldots,b_r)$.
Let us add a term $z$ to $(\yr)$. To obtain a \cop sequence from it, we proceed as follows. If $z$ is placed at the end, we add $1$ at the end of $(b_1,\ldots,b_r)$. If $z$ is placed at the start, we add $0$ at the start of~$(b_1,\ldots,b_r)$. If $z$ is intercalated between $y_i$ and $y_{i+1}$
we intercalate~$b_i$ between~$b_i$ and~$b_{i+1}$.
\end{proof}
%
\begin{fact}\label{fact2SeqSingTD}
Let  $(\ux)=(\xn)$ in a \trdiz. If we have $x_1=0$, or~$x_n=1$, or~$x_{i+1}\leq x_i$ for some $i\in\lrb{1..n-1}$,
then the sequence $(\ux)$ is singular.
\end{fact}
\begin{proof}
We apply the previous fact by noting that $(0)$ and $(1)$ are two \cop sequences and that the sequence  $(x_i,x_{i+1})$ with $x_{i+1}\leq x_i$ admits $(0,1)$
as the \cop sequence.
\end{proof}

To recap: the \cov meaning of the phrase \gui{the number of \elts of $E$ is bounded by $k$} (which we denote by $\#E\leq k$) is that for every finite list of $k+1$ \elts in $E$, two of them are equal.

\pagebreak	        

\begin{lemma}\label{lemSeqSingTD}
For a \emph{non-decreasing} sequence $(\ua)=(a_1,\ldots, a_{n})$ in a totally ordered lattice \propeq
\begin{enumerate}
\item The sequence is singular.
\item $a_1=0$, or $a_{n}=1$, or there exists an $i\in\lrb{1..n-1}$
such that $a_{i}= a_{i+1}$.
\item The number of \elts in $(0,\an,1)$ is bounded by $n+1$.
\end{enumerate}
\end{lemma}
\begin{proof}
\emph{1} $\Rightarrow$ \emph{2.}
Let us do the computation for the case $n=3$ by leaving the \recu to the skeptical reader.
Consider a \cop sequence $(b_1,b_2, b_3)$.
We have
\[\preskip-.2em \postskip.1em
\arraycolsep3pt
\begin{array}{ccc}
1  & \leq  &  a_3\vu b_3 \\[.2em]
a_3\vi b_3  & \leq  &  a_2\vu b_2 \\[.2em]
a_2\vi b_2  & \leq  &  a_1\vu b_1 \\[.2em]
a_1\vi b_1  & \leq & 0
 \end{array}
\]
Thus, $a_1=0$ or $b_1=0$. \\
If $b_1=0$, then $a_1\vu b_1=a_1\geq a_2\vi b_2$.
Therefore $a_2\leq a_1$ or $b_2\leq a_1$. In the first case, $a_1=a_2$.
In the second case, $b_2\leq a_1\leq a_2$ therefore $a_2\vu b_2=a_2$.
This implies $a_3\leq a_2$ or $b_3\leq a_2$. In the first case,~$a_2=a_3$.
In the second case, $b_3\leq a_2\leq a_3$, therefore $a_3\vu b_3=a_3=1$.

 \emph{2} $\Rightarrow$ \emph{1.} By Fact~\ref{fact2SeqSingTD}.

 \emph{3} $\Rightarrow$ \emph{2.} If we have two equal \elts in a non-decreasing sequence, then there are also two consecutively equal \eltsz.
\end{proof}

\smallskip
The following \tho gives a precise and \elr \cov interpretation of the \ddk of a totally ordered set. It directly results from Fact~\ref{fact2SeqSingTD} and from Lemma~\ref{lemSeqSingTD}.

\begin{theorem}\label{thKdimTDTO}
For a totally ordered \trdi $\gT$, \propeq
\begin{enumerate}
\item $\gT$ is \ddi$n$.
\item The number of \elts of $\gT$ is bounded by $n+2$ ($\#\gT\leq n+2$).
\item For every non-decreasing sequence $(\xzn)$ in $\gT$, we have $x_0=0$, or $x_n=1$, or $x_{i+1}= x_i$ for some $i\in\lrb{0,n-1}$.
\end{enumerate}
\end{theorem}

Note that the previous \tho applies to the Zariski lattice of a valuation \riz.
We now present two very simple and useful facts regarding  \advsz.

\begin{fact}\label{fact1ValRing}
In a \adv let $u_1$, \ldots, $u_m$ be \elts with $\sum_iu_i=0$ (and~$m\geq2$).
Then there exists a $j\neq k$ and an \iv \elt $v$ such that $\gen{u_1,\ldots,u_m}=\gen{u_j}=\gen{u_k}$ and $vu_j=u_k$.
\end{fact}
\begin{proof}
First of all there exists a $j$ such that $\gen{u_1,\ldots,u_m}=\gen{u_j}$. Then for each~$k$ let $v_k$ be an \elt such that $u_k=v_ku_j$, with $v_j=1$. \\
We obtain the \egt $u_j(1+\sum_{k\neq j}v_k)=0$.
Therefore $u_j=0$ or $1+\sum_{k\neq j}v_k=0$.
If $u_j=0$, we can take all the $v_k$'s equal to $1$.\\
If $1+\sum_{k\neq j}v_k=0$, one of the $v_k$'s is \iv since $\gV$ is local.
\end{proof}
%

\begin{fact}\label{fact2ValRing}
Let $\gV$ be a \adv and a sequence $(\an)$ in~$\gV\etl$.
For exponents $p_i$ all $>0$,  let $a=\prod_{i=1}^n a_i^{p_i}$.
Then there exists some~$j\in\lrbn$ such that $\rD_\gV(a)=\rD_\gV(a_j)$.
\end{fact}
\begin{proof}
Consider a $j$ such that $a_i$ divides $a_j$ for all $i\in\lrbn$.
Then $a_j$ divides $a$ which divides $a_j^p$, where $p=\sum_{i=1}^n p_i$.
\end{proof}
We will need the following combinatorial lemma.

\begin{lemma}\label{lemBornes}
Let $E\subseteq F$ be two sets. We suppose that for every sequence $(x_0,\ldots,x_m)$ in $F$, one of the following two alternatives takes place
\begin{itemize}
\item there exist $i<j\in\lrb{0..m}$ such that $x_i=x_j$,
\item there exists an $i\in\lrb{0..m}$ such that $x_i\in E$.
\end{itemize}
Then  $\#E\leq\ell$ implies $\#F\leq\ell+m$.
\end{lemma}
\begin{proof}
We consider a sequence $(y_0,\ldots,y_{\ell+m})$ in $F$. We need to show that there are two equal terms. We consider the first $m+1$ terms. Either two of them are equal, and the case is closed, or one of the terms is in $E$.
In this case, we delete the term which is in $E$ from the sequence $(y_0,\ldots,y_{\ell+m})$ and we consider the first $m+1$ terms of this new sequence. 
Either two of them are equal, and the case is closed, or one of the terms is in $E$ \ldots\, In the worst case, we follow the procedure till the end and we finally obtain $\ell+1$ terms in $E$ and two of them are equal.
\end{proof}
%

\begin{theorem}\label{th1Valdim}
Let $\gV$ be a  \ddvz, $\gK$ be its quotient field, $\gL\supseteq \gK$ be a \cdi of transcendence  degree $\leq m$ over~$\gK$,
and~$\gW\supseteq \gV$ be a \adv of $\gL$. Then $\Kdim\gW\leq \Kdim\gV+m$.
\end{theorem}
\begin{proof} 
We need to show that if $\,\Kdim\gV\leq n$, then $\,\Kdim\gW\leq n+m$.
\\
Since these are \advsz, we must  simply show that

\snic{\#\Zar\gV\leq n+2\;$ implies $\;\#\Zar\gW\leq n+m+2.}

(See \thref{thKdimTDTO}.)
It therefore suffices to show that the hypotheses of Lemma \ref{lemBornes} are satisfied for the integers $\ell=n+2$ and $m$, and for the sets~$E=\Zar\gV$ and $F=\Zar\gW$.
\\
Let $\gV'=\gW\cap\gK$. Since $\gV'$ is a localized \ri
 of $\gV$, we have $\Kdim\gV'\leq\Kdim\gV$. We are thus brought back to the case where $\gV=\gW\cap\gK$,
which implies $\Zar\gV\subseteq\Zar\gW$.
\\
Now let $x_0$, \ldots, $x_{m}\in\Reg \gW$, denoted by $\gW\etl$. \\
Consider an \agq \rde over $\gK$ for $(x_0,\ldots,x_{m})$. We can suppose that the \coes of the \pol $P\in\gK[X_0,\ldots,X_m]$ which gives this \agq \rde are in $\gV\cap\gK\eti=\gV\etl$. By letting, for $p \in \NN^{m+1}$, $x^p = x_0^{p_0}\cdots x_{m}^{p_m}$, Fact~\ref{fact1ValRing} gives us~$p$ and~$q$ distinct in $\NN^{m+1}$ such that~$ax^{p}$ and~$bx^{q}$ are associated in $\gW$ with~$a$,~$b\in\gV\etl$. 
By simplifying by~$x^{p \wedge q}$, we can assume $p \wedge q = 0$. Since~$a$ divides~$b$ or~$b$ divides~$a$, we can assume that~$b = 1$. We therefore have $ax^p$ associated with $x^q$ in~$\gW$.  If $q = 0$, then each $x_j$ contained in $x^p$ (there is at least one) is \iv in $\gW$, \cad
$\rD_\gW(x_j) = \rD_\gW(1)$.  Otherwise, Fact~\ref{fact2ValRing} applied to~$x^q$ gives us some~$x_j$ present in~$x^q$ such that~$\rD_\gW(x^q) = \rD_\gW(x_j)$; applied to $ax^p$, this tells us that $\rD_\gW(ax^p) = \rD_\gW(a)$ or $\rD_\gW(x_k)$ with~$x_k$ present in~$x^p$;
we therefore have~$\rD_\gW(x_j) = \rD_\gW(a)$, or $\rD_\gW(x_j) = \rD_\gW(x_k)$. The \dem is complete.
\end{proof}
%

\subsec{Valuative dimension of a commutative \riz}

\begin{definition}\label{defiValdim}~
\begin{enumerate}
\item If $\gA$ is a \qiriz, the \ix{valuative dimension} is defined as follows. Let $d\in\NN\cup\so{-1}$ and $\gK=\Frac\gA$. We say that the valuative dimension of $\gA$ is less than or equal to $d$ and we write $\Vdim\gA\leq d$ if for every finite sequence $(\ux)$ in $\gK$ we have $\Kdim \Aux \leq d$.

\item In the \gnl case we define \gui{$\Vdim\gA\leq d$} by \gui{$\Vdim\Amin\leq d$.}
\end{enumerate}
\end{definition}
We \imdt have
\begin{itemize}
  \item $\Kdim\gA\leq\Vdim\gA$,
  \item $\Vdim\gA=-1$ \ssi $\gA$ is trivial,
  \item $\Vdim\gA\leq0$ \ssi $\Kdim\gA\leq0$,
  \item if $\gA$ is a \qiri then 
\begin{itemize}
\item $\Kdim\gA=\Vdim\gA$ \ssi $\Kdim \gB  \leq \Kdim\gA$ for every intermediary \ri $\gB$  between $\gA$ and~$\Frac\gA$,
\item if $\gB$ is intermediary between $\gA$ and~$\Frac\gA$, we have $\Vdim\gB\leq\Vdim\gA$. 
\end{itemize}
\end{itemize}

\smallskip The following fact results directly from the construction of $\Amin$.
\begin{fact}\label{fact1Amin}
If $\gA$ is an \anarz, then so is $\Amin$.
\end{fact}

\begin{lemma}\label{lem1Valdim}
If $\gA$ is an \anarz, we have $\Kdim\gA = \Vdim\gA$.
\end{lemma}
\begin{proof}
Since $\Kdim\gA=\Kdim\gA\qim$, and since $\gA\qim$ is an \anar 
 if $\gA$ is arithmetic, it suffices to treat the case where $\gA$ is a \qiriz.
We then apply \thref{th.2adpcoh} which says that every \elt of $\Frac\gA$ is primitively \agq over $\gA$, and \thref{cor2thKdimMor} which says that in such a case $\Kdim \gB  \leq \Kdim\gA$ for every intermediary \ri $\gB$ between $\gA$ and~$\Frac\gA$.
\end{proof}

\rem Here is the end of a less scholarly proof (for the case where $\gA$ is an \ari \qiriz). We first suppose that $\gA$ is local, \cad it is an integral \advz.
For every $x=a/b\in\Frac\gA$, we have the alternative: $b$ divides $a$, in which case $x\in\gA$, or~$a$ divides $b$, that is, $ac=b$ in which case~$c$ is \ndz and $x=1/c$ such that~$\gA[x]$ is a localized \adv of $\gA$, so $\Kdim\gA[x]\leq\Kdim\gA$. We finish by \recu on the number of \elts of $\Frac\gA$ which we add to~$\gA$. 
Finally, in the \gnl case, we re-express the previous \demz. We replace the alternative \gui{$b$ divides $a$ or $a$ divides $b$} by the creation of two \come \lons of $\gA$. In the first $b$ divides $a$, in the second $a$ divides $b$.
\eoe

\begin{lemma}\label{lemVdimKdim}
Let $\gA$ be an integral \riz, $n\geq 1$ and $k\geq -1$.\\
If $\Kdim\AXn\leq n+k$,
then for all $x_1$, \dots, $x_n$ in $\Frac\gA$, we have $\Kdim\Axn\leq k$.
\end{lemma}
%
\begin{proof} We introduce the intermediary \ris

\snic{\gB_0=\AXn, \;\gB_1=\gA[x_1,X_2,\dots,X_n],\;\dots,
\;\gB_n=\Axn.}

For $i\in\lrbn$, let $\varphi_i$ be the \homo of \evn $\gB_{i-1}\to\gB_i$ defined by $X_i\mapsto x_i$. If $x_i=a_i/b_i$, the kernel $\Ker\varphi_i$ contains $f_i=b_iX_i-a_i$.\\ 
Let $i\in\lrb{0..n-1}$. 
Since $b_{i+1}\in\Reg\gA[(x_j)_{1\leq j\leq i}]$, we have $f_{i+1}\in\Reg\gB_{i}$ (McCoy's lemma, Corollary~\ref{corlemdArtin}).
Therefore, by item~\emph{5} of Proposition~\ref{propDdk0}, we have $\Kdim\aqo{\gB_{i}}{f_{i+1}}\leq \Kdim\gB_i-1$. 
Finally, since $\gB_{i+1}$ is a quotient of~$\aqo{\gB_{i}}{f_{i+1}}$, we obtain $\Kdim\gB_{i+1} \le \Kdim\gB_i-1$.
\end{proof}

In the following proposition, as we will see a little later, the three \prts are actually \eqves (\thref{thValDim} item~\emph{2}).

\begin{proposition}\label{propVdimKdim}
Let $\gA$ be an integral \ri and $n\geq 1$, then we have for the following items the implications {1} $\Rightarrow$  {2} $\Rightarrow$  {3.}  
\begin{enumerate}
\item We have $\Kdim\AXn\leq 2n$.  
\item For all $x_1$, \dots, $x_n$ in $\Frac\gA$, we have $\Kdim\Axn\leq n$.
\item We have $\Vdim\gA\leq n$.
\end{enumerate}
 
\end{proposition}
%
\begin{proof}
\emph{1} $\Rightarrow$ \emph{2.} Special case of Lemma~\ref{lemVdimKdim}.

\emph{2} $\Rightarrow$ \emph{3.}
We consider an arbitrary sequence $(\yr)$ in $\Frac\gA$, then an arbitrary sequence $(\xzn)$ in $\gB=\Ayr$. We need to prove that the sequence $(\xzn)$ is singular in $\gB$. It suffices to show that it is singular in $\gC=\gA[\xzn]$, or that the sequence $(\xn)$ is singular in $\gC/\IK_\gC(x_0)$. \\
We write $x_0=a_0/b_0$ with $b_0\in\Reg\gA$.
If $a_0=0$, we are done. \\
If $a_0$ is \ndzz, then $\IK_\gC(x_0)= x_0 \gC\supseteq  a_0\gC$. Therefore $\gC/\IK_\gC(x_0)$ is a quotient of $\aqo\gC{a_0}$ which is equal to $\aqo\Axn{a_0}$, which is \ddi $ n-1$. Thus $\gC/\IK_\gC(x_0)$ is \ddi $n-1$, and the sequence $(\xn)$ is singular in $\gC/\IK_\gC(x_0)$.
\end{proof}
\vspace{-.7em}
\pagebreak	        

\subsec{Valuative dimension of a \pol \riz}

The aim of this subsection is to prove the \egt

\snic{\fbox{$\Vdim\AXn = n + \Vdim\gA $}, }

for all $n\geq 1$. We deduce the same \egt for the Krull dimensions in the case of an \anarz.

\underline{By \dfnz}, this \egt of dimensions means   the following  \eqvc
\begin{equation}
\label{eqVdimAXn}
\fbox{$\forall k\geq -1,\;\;\Vdim\gA\leq k\iff \Vdim\AXn \leq  n +   k$}.
\end{equation}

Thus the first framed \egt does not quite stick for the trivial \ri
(we should say that the dimension of the trivial \ri is $-\infty$ rather \hbox{than $-1$}).

\emph{Preliminary remark.} 
Given that $\Vdim\gA=\Vdim\Amin$ by \dfnz, and that $\Amin[\Xn]\simeq(\AXn)\qim$ (Fact~\ref{fact2Amin}), 
it suffices to treat the case where $\gA$ is a \qiriz, and by the \elgbmd of \qirisz, it suffices to treat the integral case.
In the rest of the subsection, we will therefore sometimes use the saving phrase \gui{we can \spdg suppose that the \ri is integral,} or sometimes, if we want to explain the functioning of the \elr \lgbe machinery, \gui{we can \spdg suppose that the \ri is a \qiriz.} 
\eoe

\begin{fact}\label{factVdimAXnfacile}
In \pref{eqVdimAXn}, the converse implication (from right to left) is correct. 
\end{fact}
%
\begin{proof}
Suppose \spdg that $\gA$ is integral. \\ Let $[\uX]=[\Xn]$. 
Suppose $\Vdim\AuX\leq n+k$. Let $\gB=\gA[\yr]$,
with $y_i\in\Frac\gA$ for $i\in\lrbn$. 
We want to prove that $\Kdim \gB \leq k$.
\\
However, $\BuX=\AuX[\yr]$ with the $y_i$'s in $\Frac(\AuX)$.\\
Therefore $\Kdim\BuX\leq n+k$, and by Lemma~\ref{lemKdimAxn}, $\Kdim\gB\leq k$.
\end{proof}

We now study the difficult direct implication in (\ref{eqVdimAXn}).
In \clama we have the following result: \\
$(*)$ \emph{the valuative dimension of an integral \ri $\gA$ is also the maximum of the dimensions of valuation \ris containing $\gA$ and contained in its quotient field.}

This affirmation $(*)$ is no longer true in \gnl from a \cov point of view (by lack of \advsz), but it is a direct consequence (in \clamaz) of Corollary~\ref{prop4ValDim}, which is therefore a \cov version of~$(*)$.

\pagebreak	        

\begin{lemma}\label{lem2ValDim}
Let $x_0$, $x_1$, \ldots, $x_\ell$, $u$, $v$, $\alpha$ be \idtrs over a \ri $\gA$, $P_0(\alpha)$, \ldots, $P_\ell(\alpha) \in \gA[\alpha]$ and $Q_0(\alpha^{-1})$, \ldots, $Q_\ell(\alpha^{-1}) \in \gA[\alpha^{-1}]$. For some $m_i$, $n_i \in \NN$, we define $P = P(\alpha)$ and $Q = Q(\alpha^{-1})$ as follows

\snac{\arraycolsep2pt \begin{array}{rcl}
P &=&
x_0^{m_0}(x_1^{m_1}(\cdots (x_\ell^{m_\ell}
  (u + P_\ell(\alpha)x_\ell) + \cdots) + P_1(\alpha)x_1) + P_0(\alpha)x_0),
\\[1mm]
Q &=&
x_0^{n_0}(x_1^{n_1}(\cdots (x_\ell^{n_\ell}
  (v + Q_\ell(\alpha^{-1})x_\ell) + \cdots) + Q_1(\alpha^{-1})x_1) + Q_0(\alpha^{-1})x_0).
\end{array}
}

If $P$ is of formal degree $p$ (in $\alpha$), $Q$ of formal degree $q$ (in $\alpha^{-1}$), we consider the resultant
$$\preskip-.1em \postskip.4em
R = \Res_{\alpha}(\alpha^q Q,q, P,p )  \in \gA[x_0, \ldots, x_\ell, u, v].
$$
Then, by letting $r_i = qm_i + pn_i$ and $w = u^q v^p$, $R$ is of the form

\snac{ 
R = x_0^{r_0}(x_1^{r_1}(\cdots (x_\ell^{r_\ell} (w + a_\ell x_\ell) + \cdots) + a_1x_1) + a_0x_0)
\quad \hbox {with $a_i \in \gA[\ux,u,v]$}.
}
\end{lemma}
\begin{proof}
Writing $\Res_{\alpha,q,p}(U,V)$ in place of $\Res_{\alpha}(U,q, V,p )$, we suppose $n = 1$ and we let $x = x_0$, $y = x_1$, such that $P = x^{m_0} S$, $\alpha^q Q = x^{n_0} T$, with

\snic{S = y^{m_1} (u + P_1(\alpha) y) + P_0(\alpha)x, \quad
T = y^{n_1} (v\alpha^q + T_1(\alpha) y) + T_0(\alpha)x.}

We obtain $R = x^{r_0} \Res_{\alpha,q,p}(T,S)$,  $r_0 = qm_0 + pn_0$. By letting $x := 0$ we have

\snic{\arraycolsep2pt
\begin{array} {rl}
\Res_{\alpha,q,p}(T,S)_{x := 0} &= \Res_{\alpha,q,p}(T_{x:=0} ,S_{x := 0})
\\
&=
\Res_{\alpha,q,p}(y^{n_1} (v\alpha^q + T_1(\alpha)y), y^{m_1} (u + P_1(\alpha)y)\big)
\\
&=
y^{r_1} \Res_{\alpha,q,p}(v\alpha^q + T_1(\alpha)y, u + P_1(\alpha) y),
\\
\end{array}
}

with $r_1 = qm_1 + pn_1$. By letting $y := 0$ we have

\snic{
\Res_{\alpha,q,p}(v\alpha^q + T_1(\alpha)y, u + P_1(\alpha) y)_{y:=0} =
\Res_{\alpha,q,p}(v\alpha^q, u)  = u^q v^p,}

which gives the stated result.
\end{proof}

\begin{proposition}\label{prop3ValDim}
Let $\gA\subseteq\gB$, $(\ux)=(x_0,\ldots,x_n)$ be a sequence in $\gA$ and~$\alpha_0$,~$\beta_0$ in $\gB$ such that $\alpha_0\beta_0=1$. Suppose that the sequence is singular in~$\gA[\alpha_0]$ and~$\gA[\beta_0]$, then it is singular in~$\gA$.
\end{proposition}
\begin{proof}
We apply the previous lemma by specializing $u$ and $v$ in $1$. Since the \polsz~$P(\alpha)$ and $\alpha^q Q(\alpha^{-1})$ have a common root $\alpha_0$ in $\gB$, their resultant is null (Lemma~\ref{lem0Resultant}).
\end{proof}
%
\begin{corollary}\label{prop4ValDim}
Let $a$ and $b$ be \ndz \elts of a \qiri $\gA$.
Then $\Vdim\gA=\sup\big(\Vdim\gA[\fraC a b],\Vdim\gA[\fraC b a]\big)$.
\end{corollary}
\begin{proof}
The in\egts $\Vdim\gA[\fraC a b]\leq \Vdim\gA$ and $\Vdim\gA[\fraC b a]\leq \Vdim\gA$ result from the \dfn of the valuative dimension.
\\
Finally, suppose that $\Vdim\gA[\fraC a b]\leq n$ and $\Vdim\gA[\fraC b a]\leq n$ for some $n\in\NN$. Let $(\xzn)$ be a sequence in $\gA$. It is singular in $\Vdim\gA[\fraC a b]$ and~$\Vdim\gA[\fraC b a]$, therefore it is singular in $\gA$ by Proposition~\ref{prop3ValDim}.
\end{proof}

\vspace{-.7em}
\pagebreak

\begin{proposition}\label{th0ValDim}
For every \ri $\gA$  and all $n\geq 1$, we have
$$\preskip.3em \postskip.4em
\Vdim\gA[\Xn] \leq n + \Vdim\gA.$$
\end{proposition}
\begin{proof}
We need to show that if $\Vdim\gA\leq k$ then $\Vdim\gA[\Xm] \leq k+m$.
By Fact~\ref{fact2Amin}, it suffices to treat the case where 
$\gA$ is a \qiriz. 
\\
We first suppose that $\gA$ is integral. We re-express the proof of \thref{th1Valdim} and we use the dynamic method. Each time that we have a disjunction of the type \gui{$a$ divides $b$ or $b$ divides $a$} we introduce the \ris $\gC[ \fraC a b]$ and $\gC[\fraC b a]$, where~$\gC$ is the \gui{current} \riz.
At each leaf of the tree constructed thus we have a \ri $\gA[u_1,\ldots,u_\ell]\subseteq\Frac\gA$ in which the considered sequence is singular. We conclude by Proposition~\ref{prop3ValDim} that the sequence is singular in~$\gA$. 
 \\
In the case where $\gA$ is a \qiri we can call upon the \elgbmd of \qirisz. We can also reason more directly: $a$ and~$b$ produce the decomposition of \gui{the current \riz} $\gC$ in a product of four components. In three of them, $a$ or $b$ is null and everything is easy. In the fourth one, $a$ and $b$ are \ndzs and we are brought back to the integral case.
\end{proof}

As corollaries we obtain the following \thosz.

\begin{theorem}\label{thValDim} 
For a \ri $\gA$, we have the following \eqvcsz.
\begin{enumerate}
\item If $n\geq 1$ and $k\geq -1$, then
$$\preskip.2em \postskip.3em
\Vdim\gA\leq k\iff\Vdim\AXn\leq n+k. 
$$
In other words, $\Vdim\AXn=n+\Vdim\gA$.
\item If $n\geq 0$, then
$$\preskip.2em \postskip.3em
\Vdim\gA\leq n\iff \Kdim \AXn \leq 2n.
$$
In the case where $\gA$ is a \qiriz, it is also \eqv to:\\
for all $x_1$, \dots, $x_n$ in $\Frac\gA$, we have $\Kdim\Axn\leq n$.
\end{enumerate}
\end{theorem}
\begin{proof}
\emph{1.} Proved in Fact~\ref{factVdimAXnfacile} and Proposition~\ref{th0ValDim}.
\\
\emph{2.} The case $n=0$ has already been done. Let us look at the case $n\geq 1$.
The direct implication results from item~\emph{1} because $\Kdim\AXn\leq\Vdim\AXn$. The converse implication is given (in the integral case, but it is not restrictive) in Proposition~\ref{propVdimKdim}. 
\end{proof}
%

\begin{theorem}\label{corthValDim}~
\begin{enumerate}
\item If $\gA$ is an \anar of finite \ddk we have
\[
\Vdim \gA[\Xn]= \Kdim \gA[\Xn]\leq n+\Kdim \gA.
\] 
with \egt if $\gA$ is nontrivial.
\item $\Vdim \ZZ[\Xn]= \Kdim \ZZ[\Xn]=1+n$.
\item Every \ri generated by $n$ \elts is of valuative dimension (therefore of \ddkz) $\leq1+n$.
\item \label{i5corthValDim} Let $\gA$ be a \qiri generated by $n$ \elts and $\gB$ be an intermediary \ri between $\gA$ and $\Frac\gA$. Then $\Vdim\gB\leq1+n$.
\end{enumerate}
\end{theorem}
\begin{proof}
Item~\emph{1} results from the most \gnl \tho (\thref{cor0thValDim})
 and item~\emph{2} is a special case.

\emph{3.} The \ri $\gA$ is a quotient of $\ZZXn$, so $\gA[Y_1,\dots,Y_{n+1}]$ is a quotient of $\ZZXn[Y_1,\dots,Y_{n+1}]$ which is of \ddkz~$2n+2$ by item~\emph{2} Therefore $\Vdim\gA\leq n+1$
by item~\emph{2} of \thref{thValDim}.

\emph{4.} Consequence of item~\emph{3} since $\Vdim\gA\leq n+1$.
\end{proof}
%

\begin{theorem} \label{cor0thValDim} 
For a \ri $\gA$ \ddi $n$ ($n\geq 1$) \propeq
\begin{enumerate}
\item $\Vdim\gA=\Kdim\gA$.
\item For all $k\geq 1$, $\Kdim(\AXk)\leq k+\Kdim\gA$
.
\item $\Kdim(\AXn)\leq n+\Kdim\gA$
.
\end{enumerate}
Morevover if $\gA$ is nontrivial we can replace $\leq $ by $=$ in items 2 and 3.\\
When $\Vdim\gA=\Kdim\gA$, for all $ k\geq 1$, we have the \egt 

\snic{\Kdim(\AXk)=\Vdim(\AXk).}

\end{theorem}
\begin{proof}
Note that we do not assume that the \ddk of $\gA$ is  exactly known.
\\
\emph{1} $\Rightarrow$ \emph{2.} We fix some $k\geq 1$ and we need to show that for every $m\geq -1$, we have~$\Kdim\gA\leq m\Rightarrow\Kdim(\AXk)\leq m+k$.
\\
We have $\Vdim(\gA)\leq m$, so $\Vdim(\AXk)\leq m+k$ by Proposition~\ref{th0ValDim},
\hbox{therefore $\Kdim(\AXk)\leq m+k$} because we still have $\Kdim\gB\leq \Vdim\gB$.

\emph{2} $\Rightarrow$ \emph{3.} This is the special case where $k=n$.
 
\emph{3} $\Rightarrow$ \emph{1.} Suppose $\Kdim\gA\leq m$ and we need to show $\Vdim\gA\leq m$. \Spdg $0\leq m\leq n$. If $m=n$ \trf by item \emph{2} of \thref{thValDim}.
If $n= m+r$, we have $\Kdim(\AXn)\leq n+m$ by hypothesis. As $(X_{m+1},\dots,X_n)$ is singular of length $r$, item \emph{3} of Proposition~\ref{lemRegsing} gives us $\Kdim(\AXm)\leq n+m-r= 2m$ and \trf by item \emph{2} of \thref{thValDim}.

The last statement is left to the reader.
\end{proof}
%

\section{Lying Over, Going Up and Going Down}
\label{secGoingUp}

In this section we are interested in understanding in \cof terms certain \prts of commutative \ris and of their morphisms which are introduced in \clama  via the notions of Zariski spectrum or of spectral morphism (corresponding to a \ri \homoz).

As the goal of the current book is to develop the \cof framework, we will not prove that the \elr \dfns that we propose are \eqves to the \dfns usually given in \clamaz.

By making our \cov \dfns work we hope to obtain \cov versions of several \thos of \clamaz, truly usable in practice.
Actually, it is what will happen systematically in the following chapters.

\subsec{Lifting \ideps (Lying Over)}

In \clama we say that a \homo $\alpha:\gT\to\gV$ of \trdis \gui{has the lifting \prt of \idepsz} when the dual \homo $\Spec\alpha:\Spec\gV\to\Spec\gT$ is surjective, in other words when every \idep of $\Spec\gT$ is the inverse image of a \idep of $\Spec\gV$. To abbreviate we also say that the morphism is \gui{Lying Over.}
We will give a pertinent \cof \dfn without using the dual \homoz. For the \eqvc in \clama with the \dfn via the spectra, see Exercise~\ref{exoLYOV}.

\begin{definition}\label{defiLYO}~
\begin{enumerate}
\item A \homo $\alpha:\gT\to\gV$ of \trdis is said to be \emph{Lying Over} when it is injective. It amounts to the same to say that $\alpha$  reflects ineqalities: for all $a$, $b\in\gT$, $\alpha(a)\leq\alpha(b)$ implies $a\leq b$.
\item A commutative \ri \homo $\varphi:\gA\to\gB$ is said to be \emph{Lying Over} when the \homo  $\Zar\varphi:\Zar\gA\to\Zar\gB$ is injective.
\end{enumerate}
\end{definition}
\index{Lying Over!morphism}

\rdb\label{remLY}
\rem We also have the following \eqv formulations for the Lying Over morphisms.
\begin{itemize}
\item For the \trdisz:
\begin{itemize}
\item For all $b\in\gT$, $\alpha^{-1}(\dar \alpha(b)\big)=\dar b$.
\item For every \id $\fa$ of $\gT$, $\alpha^{-1}\big(\cI_{\gV}(\alpha(\fa))\big)=\fa$.
\end{itemize}
\item For the commutative \risz:
\begin{itemize}
\item For all the \itfs $\fa,\,\fb$ of $\gA$ we have the implication

\centerline{$\varphi(\fa)\subseteq\varphi(\fb)\gB\
\Longrightarrow\ \fa\subseteq\DA(\fb).$}
\item  For every \itf $\fa$ of $\gA$ we have $\varphi^{-1}(\gen {\varphi(\fa)}) \subseteq \DA({\fa})$.
\item  For every \id $\fa$ of $\gA$ we have $\varphi^{-1}\big(\DB(\gen{\varphi(\fa)})\big)=\DA({\fa})$. \eoe
\end{itemize}
\end{itemize}

\begin{fact}\label{factLOInt}
Let $\gB\supseteq\gA$ be an extension. If $\gB$ is integral or \fpte  (over $\gA$), the  inclusion morphism 
$\gA\to\gB$ is Lying Over.
\end{fact}
%
\begin{proof}
The first case is a simple reformulation of Lemma~\ref{lemLingOver} (Lying Over). In the second case, for every \itf $\fa$ of $\gA$, we have \hbox{$\fa\gB\cap\gA=\fa$}.
\end{proof}
%

\subsec{Going Up}

In \clama we say that a \homo $\alpha:\gT\to\gV$ of \trdis 
\gui{has the going up 
\prt for  chains of \idepsz} 
when we have the following \prtz.

\emph{If $\fq\in\Spec\gV$ and $\alpha^{-1}(\fq)=\fp$, every chain $\fp_1\subseteq\cdots\subseteq\fp_n$ of \ideps of $\Spec\gT$ with $\fp_1=\fp$ is the inverse image of a chain $\fq_1\subseteq\cdots\subseteq\fq_n$ of \ideps of $\Spec\gV$ with $\fq_1=\fq$.
}

\smallskip Naturally we could limit ourselves to the case $n=2$.
In this case the definition can be reread in the following way.

\emph{If $\fq\in\Spec\gV$ and $(\Spec\alpha)(\fq)=\fp$, and if we note
 $$\alpha':(\gV/(\fq=0)\to \gT/(\fp=0)$$ 
 the induced morphism,  
then the dual morphism $$\Spec\alpha':\Spec(\gV/(\fq=0))\to\Spec(\gT/(\fp=0))$$ is onto.} 

So we come back to the Lying Over.

 Here are the \cov \dfns in terms of \trdis and of commutative \risz.

\begin{definition}\label{defiGoingup}~
\begin{enumerate}
\item A \homo $\alpha:\gT\to\gV$ of \trdis is said to be \emph{Going Up} when for all $a,c\in\gT$ and $y\in\gV$ we have
$$\preskip.4em \postskip.4em 
\alpha(a)\leq\alpha(c)\vu y \quad\Longrightarrow\quad\exists x\in\gT\; (a\leq c \vu x\et \alpha(x)\leq y). 
$$
\item A \homo $\varphi:\gA\to\gB$ of commutative \ris is said to be \emph{Going Up} when the \homo  $\Zar\varphi:\Zar\gA\to\Zar\gB$ is Going~Up.
\end{enumerate}
\end{definition}
\index{Going Up!morphism}

\rems
1)
For item \emph{1}, if $\fa=\alpha^{-1}(0_\gV)$ and $\gT_1=\gT\sur{(\fa=0)}$,
then $\alpha$ is Going Up \ssi  $\alpha_1:\gT_1\to\gV$ is going~up.
\\
For item \emph{2}, if $\gT=\ZarA$, then $\gT_1\simeq\Zar(\varphi(\gA)\big)$. We deduce, by letting $\gA_1=\varphi(\gA)$, that $\varphi$ is Going Up \ssi  $\varphi_1:\gA_1\to\gB$ is going~up.

 2)
For the \trdisz, if $\alpha^{-1}(0)=0$ and if $\alpha$ is Going Up, then it is Lying Over. For the commutative \risz, if $\Ker\varphi\subseteq\DA(0)$ and if $\varphi$ is Going Up,  then it is Lying Over.
\eoe

\begin{proposition}\label{prop1Gup}
If $\gB$ is an integral \Algz, the morphism $\gA\to\gB$ is Going Up.
\end{proposition}
\begin{proof}
By the previous remark we can assume $\gA\subseteq\gB$.
We then know that the \homo is Lying Over, that is we know that $\ZarA\to\Zar\gB$ is injective, so we can identify $\ZarA$ with a sublattice of $\Zar\gB$.
We need to show that given $a_1$, \dots, $a_n$, $c_1$, \dots, $c_q$ in $\gA$ and $y_1$, \dots, $y_p$ in $\gB$ satisfying

\snic{\DB(\ua)\leq\DB(\uc)\vu\DB(\uy),}

we can find a sequence $(\ux) $ in $\gA$ such that
$$\preskip.4em \postskip.4em 
\DA(\ua)\leq\DA(\uc)\vu\DA(\ux)\;\hbox{  and  }\;
\DB(\ux)\leq \DB(\uy). 
$$
Let $\fb=\DB(\uy)$,  $\fa=\fb\cap\gA$, $\gB_1=\gB\sur{\fb}$ and $\gA_1=\gA\sur{\fa}$. We consider the integral extension $\gB_1\supseteq\gA_1$.
The hypothesis is now that ${\rD_{\gB_1\!}(\ua)\leq\rD_{\gB_1\!}(\uc).}$
\\
By Lying Over we know that this implies that ${\rD_{\gA_1\!}(\ua)\leq\rD_{\gA_1\!}(\uc).}$
This means that for each $i\in\lrbn$ we have some $x_i\in\fa$ such that $\DA(a_i)\leq\DA(\uc)\vu\DA(x_i)$. We have therefore attained the sought goal with $(\ux)=(\xn)$.
\end{proof}
%

\subsec{Going Down}

In \clama we say that a \homo $\alpha:\gT\to\gV$ of \trdis \gui{has the going down \prt for chains of \idepsz} when  
the opposite morphism $\alpha\eci:\gT\eci\to\gV\eci$ is Going Up.
In other words
we have the following \prtz.

\emph{If $\fq\in\Spec\gV$ and $\alpha^{-1}(\fq)=\fp$, every chain $\fp_1\subseteq\cdots\subseteq\fp_n$ of \ideps of $\Spec\gT$ with $\fp_n=\fp$ is the inverse image of a chain $\fq_1\subseteq\cdots\subseteq\fq_n$ of \ideps of $\Spec\gV$ with $\fq_n=\fq$.
}

\smallskip Naturally we could limit ourselves to the case $n=2$.
And our \cov \dfn is the notion  opposite to Going Up: we reverse the order relation.

\begin{definition}\label{defiGoingdown}~
\begin{enumerate}
\item A \homo $\alpha:\gT\to\gV$ of \trdis is said to be \emph{Going Down} when the same \homo for the opposite lattices $\gT\eci$ and $\gV\eci$ is Going Up.
In other words for all $a,c\in\gT$ and $y\in\gV$ we have
$$\preskip.4em \postskip.4em 
\alpha(a)\geq\alpha(c)\vi y \quad\Longrightarrow\quad\exists x\in\gT\; (a\geq c \vi x \et \alpha(x)\geq y). 
$$
\item A \homo $\varphi:\gA\to\gB$ of commutative \ris is said to be \emph{Going Down} when the \homo  $\Zar\varphi:\Zar\gA\to\Zar\gB$ is Going Down.
\end{enumerate}
\end{definition}
\index{Going Down!morphism}

\rems 1) The \dfn in item {1} comes down to saying that the image by $\alpha$ of the conductor \id $(a:c)_\gT$
generates the \id $(\alpha(a):\alpha(c)\big)_{\gV}$. So if the \trdis are \agHsz, it means that the lattice \homo is also a \homo of \agHsz.

 2) Same remarks as for Going Up. \\
 If $\ff=\alpha^{-1}(1_\gV)$ and $\gT_2=\gT\sur{(\ff=1)}$, then $\alpha$ is Going Down \ssi  $\alpha_2:\gT_2\to\gV$ is Going Down.
\\
This gives for the commutative \risz: if $S=\varphi^{-1}(\gB\eti)$ and $\gA_2=\gA_S$, then $\varphi$ is  Going Down \ssi  $\varphi_2:\gA_2\to\gB$ is Going Down.
\\
For the \trdisz, if $\alpha^{-1}(1)=1$ and $\alpha$ is Going Down, then it is Lying Over. For the commutative \risz, if $\varphi^{-1}(\gB\eti)\subseteq\Ati$ and $\varphi$ is Going Down, then it is Lying Over.
\eoe

\begin{theorem}\label{propGupLY}
If a \homo $\alpha:X\to Y$ (of \trdis or of commutative \risz) is Lying Over and Going Up,
or if it is Lying Over and Going Down, we have $\Kdim X\leq\Kdim Y$.
\end{theorem}
\rem This is the case, for example, when the \ri $\gB$ is an integral extension of $\gA$. We thus find
Proposition~\ref{propDKEXENT} again. For the flat extensions, see Proposition~\ref{prop1Gdown}. 
\eoe
\begin{proof}
It suffices to treat the Going Up case with lattices.
\\
Suppose $\Kdim Y\leq n$ and consider a sequence $(\azn)$ in $X$. %
We have in $Y$ a \cop sequence $(y_0,\ldots,y_n)$ of $\alpha(\ua)$

\snic{\alpha(a_0)\vi y_ 0\leq0,\;\ldots,\;\alpha(a_{n})\vi y_n\leq\alpha(a_{n-1})\vu y_{n-1}, \; 1\leq\alpha(a_n)\vu y_n.}

We will construct a \cop sequence $(\xzn)$ of $(\ua)$ in $X$. 
At step $n$, 
by Going Up, there exists an $x_ n\in X$ such that

\snic{1\leq a_n\vu x_n$ and $\alpha(x_n)\leq y_n.}

This gives at the stage $n-1$ the in\egtz:
$\alpha(a_{n}\vi x_n)\leq \alpha(a_{n-1})\vu y_{n-1}.$
\\
By Going Up there exists an $x_ {n-1}\in X$ such that

\snic{a_{n}\vi x_n\leq a_{n-1}\vu x_{n-1}\hbox{  and  }\alpha(x_{n-1})\leq y_{n-1}.}

We continue in the same way until stage $0$, where this time we need to use the Lying Over.
\end{proof}
%

\begin{lemma}\label{lem1Gdown}
For a \ri \homo $\varphi:\gA\to\gB$ to be Going Down it is necessary and sufficient that for all $c$, $a_1$, \ldots, $a_q\in\gA$ and $y\in\gB$ such that~$\varphi(c)y\in\DB(\varphi(\ua)\big)$, there exist some \elts $x_1$, \dots, $x_m\in\gA$ such that
$$\preskip-.2em \postskip.4em 
\DA(c)\vi\DA(\ux)\leq\DA(\ua)\;\et\;\DB(y)\leq\DB(\varphi(\ux)\big). 
$$
\end{lemma}
\begin{proof}
In the \dfn we have replaced an arbitrary \elt $\DA(\uc)$ of~$\ZarA$ and an arbitrary \elt $\DB(\uy)$ of $\Zar\gB$ by \gtrsz~$\DA(c)$ and~$\DB(y)$.
As the \gtrs $\DA(c)$ (resp.\ $\DB(y)$) generate $\ZarA$ (resp.\ $\Zar\gB$) by finite suprema, the rules of \dit imply that the restriction to these \gtrs is sufficient (computations left to the reader).
\end{proof}
%

\begin{proposition}\label{prop1Gdown}
A \homo $\varphi:\gA\to\gB$ of commutative \ris is Going Down in the following two cases.
\begin{enumerate}
\item $\gB$ is a flat \Algz.
\item $\gB\supseteq \gA$ is a domain integral over $\gA$, 
and $\gA$ is \iclz.
\end{enumerate}
\end{proposition}
\begin{proof} We assume the hypotheses of Lemma~\ref{lem1Gdown}, with an \egt in~$\gB$,
$$\ndsp
\varphi(c)^\ell y^\ell +\som_{i=1}^qb_i\varphi(a_i)=0\eqno(*)
$$
\emph{1.}
We consider $(*)$ as a 
          $\gB$-syzygy 
between the \elts $c^\ell$, $a_1$, \ldots, $a_q$.
We express that it is a $\gB$-\lin combination of $\gA$-syzygies.
These relations are written as $x_jc^\ell+\sum_{i=1}^q u_{j,i}a_i=0$ for $j\in\lrbm$, with the~$x_j$'s and the~$u_{j,i}$'s in $\gA$.
Hence $\DA(cx_j)\leq\DA(\ua)$, and 
$\DA(c)\vi\DA(\ux)\leq\DA(\ua)$. 
Finally,~$y^\ell$
is a $\gB$-\lin combination of the $\varphi(x_j)$'s, hence $\DB(y)\leq\DB(\varphi(\ux)\big)$.

 \emph{2.}
By $(*)$, $(cy)^\ell \in \gen{\ua} \,\gB$. By the Lying Over~\ref{lemLingOver2}, $(cy)^\ell$, and a fortiori~$cy$, is integral over $\gen{\ua}_\gA$. We write an \rdi for~$cy$ over the \id $\gen{\ua}_\gA$ in the form $f(cy)=0$ with

\snic{f(X)=X^k+\sum_{j=1}^k\mu_jX^{k-j}
\qquad \hbox {where }\mu_j\in\gen{\ua}_{\!\gA}^{\,j}.}

Moreover, $y$ annihilates a \polu $g(X)\in\AX$. Consider in~$(\Frac\gA)[X]$ the \mon gcd $h(X)=X^{m}+x_1X^{m-1}+\cdots+x_m$ of the two \pols $f(cX)$ and $g(X)$. Since $\gA$ is \iclz, \KROz's \tho says that $x_j\in\gA$, and the \egt $h(y)=0$ gives $y\in\DB(\ux)$. \\
It remains to see that $cx_j\in\DA(\ua)$ for $j\in\lrbm$. By formally replacing~$X$ with~$Y/c$, 
we get that the \pol 

\snic{h_c(Y)=Y^{m}+cx_1Y^{m-1}+\cdots+c^mx_m}

divides~$f(Y)$ in $(\Frac\gA)[Y]$. \KROz's \tho (under the form of Lemma~\ref{lemthKroicl}) tells us that $cx_j\in\DA(\mu_1,\ldots,\mu_k)$.\\
Finally, as~$\DA(\mu_1,\ldots,\mu_k)\leq\DA(\ua)$, we indeed have $cx_j\in\DA(\ua)$.
\end{proof}
%

\subsection*{Incomparability}
In \clama we say that a \homo $\alpha:\gT\to\gT'$ of \trdis \gui{has the incomparability \prtz} when the fibers of the dual \homo $\Spec\alpha:\Spec\gT'\to\Spec\gT$ are constituted of pairwise incomparable \eltsz. 
In other words, for $\fq_1$ and $\fq_2$ in~$\Spec\gT'$, if $\alpha^{-1}(\fq_1)=\alpha^{-1}(\fq_2)$ and $\fq_1\subseteq\fq_2$, then $\fq_1=\fq_2$.

The corresponding \cov \dfn is that the morphism $\gT\to\gT'$ is \zedz.

We have already given the \dfn of the dimension of a morphism in the case of commutative \risz. An analogous \dfn can be provided for the \trdisz, but we will not be using it.

The principal consequence of the incomparability situation for a \homo $\varphi:\gA\to\gB$ is the fact that $\Kdim\gB\leq\Kdim\gA$. This is a special case of \thref{thKdimMor} with the important \thref{cor2thKdimMor}.

\Exercices{

\begin{exercise}
\label{exoKdimLecteur}
{\rm  We recommend that the \dems which are not given, or are sketched, or
left to the reader,
etc, be done.
But in particular, we will cover the following cases.
\begin{itemize}\itemsep0pt
\item \label{exoDdk1} Prove Proposition~\ref{propDdk0}.
\item \label{exobord}
Prove what is stated in Examples on \paref{exlKdim}.
\item \label{exocorfactDDKTRDI} Prove Fact \ref{corfactDDKTRDI}.
\item \label{exopropDDKagH}
Prove Facts \ref{factDdkTrdiBord} and \ref{propDDKagH}.
\item \label{exopropNoetAgH} 
Check the details in the \dem of Proposition \ref{propNoetAgH}.
\item \label{exolem2qi} 
Prove Lemma~\ref{lem2qi} using the \dem of Lemma~\ref{lem2SousZedRed} as inspiration.

%
%
%
\item \label{exolem1Gdown} Check the details in the \dem of Lemma~\ref{lem1Gdown}.
\end{itemize}
}
\end{exercise}

\vspace{-1em}
\begin{exercise}
 \label{exoIdFilPrem}
 {\rm  If $\ff$ is a filter of the \ri $\gA$, let us define its \emph{complement}~$\wi{\ff}$
 as being $\sotq{x\in\gA}{x\in\ff\Rightarrow0\in\ff}$.
 In particular, we still have $0\in\wi\ff$, even if $0\in\ff$.
 Similarly, if $\fa$ is an \id of the \ri $\gA$, let us define its \emph{complement} $\bar{\fa}$
 as being $\sotq{x\in\gA}{x\in\fa\Rightarrow1\in\fa}$. Show that if $\ff$ is a prime filter its complement is an \idz. If in addition $\ff$ is detachable,
 then $\fa$ is a detachable \idepz. Also show the dual affirmations.

} \end{exercise}

\vspace{-1em}
\pagebreak	        
\begin{exercise}
\label{exoDimKX}
{\rm  \emph{1.} If the sequence $(\Xn)$ is singular in the \ri $\AXn$,
then $\gA$ is trivial.
\\
\emph{2.} Let $k\in\NN$. Prove that if $\AX$ is a \ri  \ddi$k$ then
$\gA$ is \ddi$k-1$.
Thus obtain once again item~\emph{1}.
}
\end{exercise}

\vspace{-1em}
\begin{exercise}
\label{exoDimKX1n}
{\rm  Prove that if $\gK$ is a \ri of \ddk exactly equal to $0$ then  $\gK[X_1,\ldots ,X_n]$ is of \ddk exactly equal to $n$.
}
\end{exercise}

\vspace{-1em}
\begin{exercise}\label{exoPartitionUnite}
 {(Partition of unity associated with an open covering of the spectrum)}\\
{\rm  
Let $\gA$ be a \ri and $(U_i)_i$ be an open covering of $\Spec(\gA)$. 
Show in \clama that there exists a family $(f_i)_i$ of \elts of $\gA$ with $f_i = 0$ except for a finite number of indices $i$ and

\snic {
(\star) \qquad\qquad
\DA(f_i) \subseteq U_i, \qquad   \sum_i f_i = 1
.}

Remark: thus, we replace every open covering of $\Spec(\gA)$ by a finite \sys of \elts of $\gA$ which \gui {cover} $\gA$
(since their sum is equal to $1$), without \gui{losing information} since $(\star)$ confirms once again that $(U_i)_i$ is a covering.
}
\end{exercise}

\vspace{-1em}
\begin{exercise}
\label{exoKdimGeom}
{\rm  For a \apf $\gA$ over a nontrivial \cdiz, let us call the \gui{\Noe dimension of $\gA$} the number of \agqt independent variables after a \iNoe position.
\\
\emph{1.} Let $f\in\gA\supseteq\KYr=\KuY$ ($\gA$ integral over $\KuY$). 
\\
\emph{1a.} Show that the boundary \id of $f$ contains a $g\in\KuY\setminus\so{0}$. 
\\
\emph{1b.} 
Deduce that the Krull boundary \ri $\gA\sur{\JK_\gA(f)}$ is a quotient of a \apf whose \Noe dimension is $\leq r-1$.
\\
\emph{2.}  
Deduce a direct \dem of the \egt of Krull and \Noe dimensions of the \apfs over a nontrivial \cdiz.
}
\end{exercise}

\vspace{-1em}
\begin{exercise}
\label{exolemLocMemeKdim}
{\rm  \emph{1.} Let $\gK$ be a nontrivial \cdiz, $\KuX=\KXn$ and~$f\in \KuX\setminus\so{0}$,
then $\Kdim \KuX[1/f] =n$.
\\
\emph{2.} 
More \gnltz, give a sufficient condition on the \pol $\delta\in\AuX$ for us to have $\Kdim(\AuX[1/\delta])=\Kdim\AuX$ (see the \dem of Lemma~\ref{lemLocMemeKdim}).
}
\end{exercise}

\vspace{-1em}
\begin{exercise}\label{exoPruferKdim} {(\Carn of  integral \adps \ddi$1$)}
\\
{\rm
Let $\gA$ be an \icl \riz.
\\
\emph{1.} Show that if $\Kdim\AX\leq2$, then $\gA$ is a \adpz, by showing that every \elt of~$\Frac\gA$ is primitively \agq over $\gA$.
\\
\emph{2.} Show that $\gA$ is a \adp \ddi1 \ssi \hbox{$\Kdim\AX\leq2$}.
\\
\emph{3.} Can we \gnr to a normal \riz?
}
\end{exercise}

\vspace{-1em}
\begin{exercise}\label{exoMultiplicativiteIdeauxBords}
{(A multiplicative \prt of boundary \idsz)}\\
{\rm  
\emph{1.}
For $a, b \in \gA$ and two sequences $(\ux)$, $(\uy)$ of \elts of $\gA$,
show that

\snic {
\IK_\gA(\ux, a, \uy)\,\IK_\gA(\ux, b, \uy) \subseteq \IK_\gA(\ux, ab, \uy)
.}

\emph{2.}
Deduce that $\IK_\gA(a_1b_1, \ldots, a_nb_n)$ contains the product  $\prod_{\uc} \IK_\gA(\uc)$,
 in which the sequence $(\uc) = (c_1, \ldots, c_n)$ ranges over the set of $2^n$ sequences such that $c_i = a_i$ %
or $c_i = b_i$ for each $i$.
}

\end{exercise}

\vspace{-1em}
\begin{exercise}\label{exoInclusionBordLionel}
{(Boundary \ids and \agq relations)}
\\
{\rm  
\emph {1.}
We consider the lexicographical order over $\NN^n$. Let $\alpha = (\alpha_1, \ldots, \alpha_n) \in\NN^n$. Prove, for $\beta > \alpha$, that

\snic {
\uX^\beta \in \gen {X_1^{1+\alpha_1},\ X_1^{\alpha_1} X_2^{1+\alpha_2},\
X_1^{\alpha_1} X_2^{\alpha_2} X_3^{1+\alpha_3},\ \cdots,\ 
X_1^{\alpha_1} X_2^{\alpha_2} \cdots X_{n-1}^{\alpha_{n-1}}X_n^{1+\alpha_n}}
.}

\emph {2.}
Let $\gA$ be a reduced \riz, $(\ux) = (\xn)$ be a sequence in $\gA$ and $P = \sum_\beta a_\beta \uX^\beta$ in~$\gA[\uX]$, which annihilates $\ux$.
\begin{itemize}
\item [\emph {a.}]
Show, for $\alpha\in\NN^n$, that
$a_\alpha \prod_{\beta < \alpha} \Ann(a_\beta) \subseteq \IK(\ux)$.

\item [\emph {b.}]
Deduce
$$\preskip-.4em \postskip.4em\ndsp 
\prod_\beta \IK(a_\beta) \subseteq \IK(\ux) + 
\prod_\beta \Ann(a_\beta). 
$$
\end{itemize}
\emph {3.}
Let $\gA \to \gB$ be an \alg with reduced $\gB$ 
        and let $x\in\gB$ be primitively \agq over $\gA$: 
$\sum_{i=0}^d a_i x^i = 0$ with $a_i \in \gA$ and $1 \in \gen{a_i, i \in \lrb{0..d}}$. Deduce from the previous question that $\IK_\gB(x)$ contains the image of $\prod_{i=0}^d \IK_\gA(a_i)$.

 \emph {4.}
Deduce a new \dem of \thref{cor2thKdimMor}:
if every \elt of $\gB$ is primitively \agq over $\gA$, then
$\Kdim\gB\leq\Kdim\gA$.

}

\end{exercise}

\vspace{-1em}
\begin{exercise}\label{exoExtEntiereIdealBord}
 {(Integral extension of the boundary \id $\IK$)}\\
{\rm  
Let $\gA \subseteq \gB$ be an integral extension of \risz.

\emph {1.}
If $\fa$ is an \id of $\gA$, $\fb$ is an \id of $\gB$, show that

\snic {
\gA \cap (\fb + \fa\gB) \subseteq \rD_\gA(\fa + \gA\cap\fb)
.}

\emph {2.}
Deduce, for $a_0$, \ldots, $a_d \in \gA$,
$$\preskip.2em \postskip.4em 
\gA\cap\IK_\gB(a_0, \ldots, a_d) \subseteq 
\rD_\gA\big(\IK_\gA(a_0, \ldots, a_d)\big)
. 
$$
\emph {3.}
Give a new \dem of the fact that $\Kdim\gA \le \Kdim\gB$, see Proposition~\ref{propDKEXENT} and \thref{propGupLY}. Compare with Exercise~\ref{exoExtEntiereMonoideBord}.
}

\end {exercise}

\vspace{-1em}
\begin{exercise}\label{exoExtEntiereMonoideBord}
 {(Integral extension of the boundary \mo $\SK$)}\\
{\rm  
Let $\gA \subseteq \gB$ be an integral extension of \risz.

\emph {1.}
Let $\fa$ be an \id of $\gA$ and $S \subseteq \gA$ be a \moz.
Show that
$$\preskip.3em \postskip.2em 
S + \fa\gB \subseteq \satu {(S + \fa)}{\gB}
. 
$$
\emph {2.}
Deduce, for $a_0$, \ldots, $a_d \in \gA$,
$$
\preskip.3em \postskip.4em 
\SK_\gB(a_0, \ldots, a_d) \subseteq 
\satu{\SK_\gA(a_0, \ldots, a_d)}{\gB}
. 
$$
\emph {3.}
Give a new \dem of the fact that $\Kdim\gA \le \Kdim\gB$.
}

\end{exercise}

\vspace{-1em}
\begin{exercise}
\label{exoKdimSomTr}
{\rm  Let $\gK$ be a nontrivial \cdiz. Denote   $(\Xn)$  by $(\uX)$ and $(\Ym)$ and $(\uY)$
Let $\gA=\gK(\uX)\otimes _\gK\gK(\uY)$.
We intend to determine the \ddk of $\gA$.
\vspace{-.5em}
\begin{itemize}\itemsep0pt
\item [\emph{1.}] $\gA$ is the \lon of $\gK[\uX,\uY]$ at ${S=(\KuX)\etl(\KuY)\etl}$.
It is also a \lon of $\gK(\uX)[\uY]$ and of $\gK(\uY)[\uX]$. 
Consequently $\Kdim\gA\leq\inf(m,n)$.
\item [\emph{2.}] Suppose $n\leq m$. Show that the sequence $(X_1-Y_1,\ldots,X_n-Y_n)$ is a regular sequence in $\gA$.
\end{itemize}
 Conclude that $\Kdim\gA = \inf(n,m)$.
 
}
\end{exercise}


\vspace{-1em}
\begin{exercise}\label{exoDualiteBords}
{(Prime \idsz, boundaries and duality)}\\
{\rm  
Let $\fp_0\subsetneq\fp_1\subsetneq \cdots \subsetneq\fp_{d-1} \subsetneq\fp_d\subsetneq\gA$ be a chain of detachable \ideps with $x_1\in\fp_1\!\setminus\fp_0$, $x_2\in\fp_2\!\setminus\fp_1$, \dots, $x_{d}\in\fp_{d}\!\setminus\fp_{d-1}$, according to the following diagram.

\snic {
\def \foo {\ar@{.}[ld]|\notin\ar@{.}[rd]|\in}
\def \subnot {\subsetneq}
\xymatrix @R=15pt @C=8pt {
     &x_1\foo    &      &x_2\foo &&&
         &x_{d-1}\foo    &      &x_d\foo   \\
0\in\fp_0&\subnot &\fp_1 & \subnot & \fp_2 
&\subnot\cdots\subnot &\fp_{d-2}&\subnot& \fp_{d-1} &\subnot& \fp_d&\!\!\!\!\!\!\not\ni1 }}

\emph {1.}
Show that $\IK(x_1, \ldots, x_i) \subseteq \fp_i$ for $i\in\lrb{0..d}$.  Therefore $\IK(x_1, \ldots, x_d) \subseteq \fp_d$. In addition, if $x_{d+1} \notin \fp_d$, then $\IK(x_1, \ldots, x_d, x_{d+1}) \subseteq \fp_d + \gA x_{d+1}$. Consequently, %
if $x_{d+1}\notin \fp_d$ and $1 \in \IK(x_1, \ldots, x_d, x_{d+1}) $, then $1 \in \fp_d + \gA x_{d+1}$.

\emph {2.}
Consider the \cop prime filters $\ff_i = \gA\setminus\fp_i$ for $i\in\lrb{0..d}$. We have the dual diagram of the  previous one.

\snic {
\def \foo {\ar@{.}[ld]|\notin\ar@{.}[rd]|\in}
\def \subnot {\subsetneq}
\xymatrix @R=15pt @C=8pt {
     &x_d\foo    &      &x_{d-1}\foo &&&
         &x_2\foo    &      &x_1\foo   \\
1\in\ff_d&\subnot &\ff_{d-1} & \subnot & \ff_{d-2} 
&\subnot\cdots\subnot &\ff_2&\subnot& \ff_1 &\subnot& \ff_0&\!\!\!\!\!\!\not\ni0 \\
}}

Show that $\SK(x_{i+1}, \ldots, x_d) \subseteq \ff_i$ for $i \in \lrb {0..d}$. Therefore $\SK(x_{1}, \ldots, x_d) \subseteq \ff_0$. In addition, if $x_{0} \notin \ff_0$, \cad if $x_0 \in \fp_0$, then $\SK(x_0, x_1, \ldots,
x_d) \subseteq x_0^\NN \ff_0$. Consequently, if $x_0\notin \ff_0$ and $0 \in \SK(x_0,x_1,\ldots, x_d)$, then $0 \in x_0^\NN \ff_0$.

Note:  $\fp_d + \gA x_{d+1}$ is the \id generated by $\fp_d$ and $x_{d+1}$, dually $x_0^\NN\ff_0$ is the \mo generated by $\ff_0$ and $x_0$.

}

\end{exercise}

\vspace{-1em}
\begin{exercise}\label{exoEliminationEtBord}
{(Elimination and boundary \ids in \pol \risz)}\\
{\rm  
 Here is a detailed \dem of the in\egt $\Kdim\gA[T] \le 1 + 2\Kdim\gA$ (Section~\ref{secKdimMor}), 
with a few further results. 
\Spdg  $\gA$ is assumed to be reduced.

\emph {1.}
Let $f \in \gA[T]$ be a \pol such that the annihilator of each \coe is generated by an \idmz. For $g \in \gA[T]$, define $R \in \gA[X,Y]$ such that $\Ann(R) = 0$ %
and $R(f,g) = 0$: note that the \pol $\Res_T(f(T)-X, Y-g(T)\big)$ solves the question when $f$ is \mon of degree
$\geq1$ (why?), and use Lemma~\ref{lemQI}.

\emph {2.}
By using Exercise~\ref{exoInclusionBordLionel}, show that if $R = \sum_{i,j} r_{ij}X^i Y^j$, we have

\snic {
\prod_{i,j} \IK_{\AT}(r_{ij}) \subseteq \IK_{\AT}(f,g).
}

\emph {3.}
By using a \ri of type $\gA_{\so\ua}$ (Lemma~\ref{lem4MorRc} and Exercise~\ref{exoAminEtagesFinis}),
find the in\egt $\Kdim\gA[T] \le 1 + 2\Kdim\gA$.

\emph {4.}
Show the following more precise result: for a reduced \ri $\gA$ \hbox{and $f$, $g \in \gA[T]$}, the \id  $\rD_\AT\left(\IK_{\gA[T]}(f,g)\right)$ contains a finite product of boundary \ids $\IK_\gA(a)$, $a \in \gA$.

\emph {5.}
More \gnltz: if $\gA[\uT] = \gA[T_1, \ldots, T_r]$ and $f_0$, \ldots, $f_r \in \gA[\uT]$, then the nilradical of the boundary \id $\IK_{\gA[\uT]}(f_0, \ldots, f_r)$ contains a finite product of boundary \ids $\IK_\gA(a_i)$, with~$a_i \in \gA$. We once again deduce that $1 + \Kdim\gA[\uT] \le (1+r)(1 + \Kdim\gA)$.

}

\end{exercise}

\vspace{-1.2em}
\pagebreak	        
\begin{exercise}\label{exoIdealBordPolynomes}
{(Boundary \ids of \polsz)}
{\rm  
Continued from Exercise~\ref{exoEliminationEtBord}.
 \\
\emph {1.}
Let $x$, $y \in \gB$ and $(z_j)$ be a finite family in $\gB$ satisfying $\prod_j \IK(z_j) \subseteq \IK(x,y)$.
Show that for $(\bn)$ in $\gB$, $\prod_j \IK(z_j,\bn) \subseteq \IK(x,y,\bn)$.

\emph {2.}
Let $T$ be an \idtr over a \ri $\gA$.
\begin{itemize}
\item [\emph {a.}]
For $(\an)$ in $\gA$, prove that $\IK_\gA(\an)\gA[T] = \IK_{\gA[T]}(\an)$. 

\item [\emph {b.}]
Show that the boundary \id of $2d$ \pols of $\gA[T]$ contains, up to radical, a product of boundary \ids of $d$ \elts of $\gA$.
\\
Consequently $\Kdim\gA < d \Rightarrow \Kdim\gA[T] < 2d$; this is another form of the in\egt $\Kdim\gA[T] \le 1 + 2\Kdim\gA$.
\end{itemize}

\emph {3.}
How can we \gnr the first and second item?

}
\end{exercise}

\vspace{-1em}
\begin{exercise}
\label{exoKdimEspanol} (Another \dfn of the \ddk of  \trdisz, see \cite[Espa\~nol]{Espa08})
{\rm  In an ordered set, a sequence $(\xzn)$ is called a \emph{chain of length $n$} if we have $x_0\leq x_1\leq\cdots\leq x_{n}$. 
In a \trdiz, two chains $(\xzn)$ and $(b_0,\dots,b_n)$ are said to be \emph{linked}, if there exists a chain $(c_1,\dots,c_n)$ with 
\begin{equation}\label{eqEspanol}
\left.\arraycolsep2pt
\begin{array}{rcccl}
 x_0\vi b_0& =  & 0    \\
 x_1\vi b_1& = & c_1 &=  &  x_0\vu b_0  \\
\vdots~~~~& \vdots &\vdots& \vdots & ~~~~  \vdots \\
 x_n\vi  b_n & = & c_n & =& x_{n -1}\vu b_{n -1}  \\
&& 1& =  &   x_n\vu b_n
\end{array}
\right\}
\end{equation}
Please compare 
with \Dfn \ref{defiDDKTRDI} for the \cop sequences.
Also note that if the sequences $(\xzn)$, $(b_0,\dots,b_n)$ and $(c_1,\dots,c_n)$ are linked by \eqns (\ref{eqEspanol}), then they are chains.

 \emph{1.} If in a \trdi we have $x \le y$ and $x\vu  a\geq  y \vi b$, then we can explicate~$a'$ and $b'$ such that

\snic {
x\vi a' =  x\vi a, \qquad  y \vu b' = y\vu b, \qquad
x\vu  a' =  y \vi b'
.}

Therefore from a left-configuration (by still assuming that $x \le y$), 
we can construct a right-configuration
$$
\left\{\,\arraycolsep2pt
\begin{array}{rcl}
x\vi a   & =  &  p  \\
x\vu  a & \geq  & y \vi b  \\
q & =  &   y\vu b
\end{array}
\right.
\qquad\qquad
\left\{\,\arraycolsep2pt
\begin{array}{rcl}
x\vi a'   & =  &  p  \\
 x\vu  a' & =  & y \vi b'  \\
 q& =  &   y\vu b'
\end{array}
\right.
$$

 \emph{2.} In a \trdiz, a chain $(\xzn)$ has a \cop sequence \ssi there exists a chain which is linked to it.

\sni \emph{3.} For a \trdi $\gT$ \propeq
\begin{enumerate}
\item [\emph{a.}] $\gT$ has \ddk $\leq n$.
\item [\emph{b.}] Any chain of length $n$ has a \cop sequence.  
\item [\emph{c.}] Any chain of length $n$ has a linked une chain.
\end{enumerate}
}
\end{exercise}

\vspace{-1em}
\begin{exercise}\label{exoAminEtagesFinis}
 {(A few results on the finite stages of $\Amin$)}
\\
{\rm  Let $\gA$ be a reduced \riz. For \ids $\fa$, $\fb$ of $\gA$ let $\fa\diamond\fb = (\fa\epr\fb)\epr =
(\fa^{\perp\perp} : \fb)$.

\emph {1.}
Prove that $\gA\sur{\fa\diamond\fb}$ is a reduced \ri in which $\fa$ is null and $\fb$ faithful.

\emph{2.}
Prove that $(\gA\sur{\fa_1\diamond\fb_1})\sur{(\ov{\fa_2}\diamond\ov{\fb_2})}
\simeq \gA\sur{\fa_3\diamond\fb_3}$ with $\fa_3 = \fa_1+\fa_2$,
$\fb_3 = \fb_1\fb_2$.

\emph {3.}
Let $(\ua) = (a_1, \ldots, a_n)$ in $\gA$.
In Lemma~\ref{lem4MorRc} we have defined (for $I\in\cP_n$)
$$\preskip.4em \postskip.4em \ndsp
\fa_I = \gen {a_i, i \in I} \diamond \prod_{j\notin I}a_j
\qquad
\gA_{\so\ua} = \prod_{I\in\cP_n} \gA\sur{\fa_I}. 
$$
Thus, modulo $\fa_I$, $a_i$ is null for $i \in I$ and \ndz for $i \notin I$.
Finally, let $\vep_i$ be the \idm of $\gA_{\so\ua}$ whose \coo in  $\gA\sur{\fa_I}$ is $1$ if $i \in I$, $0$ if $i \notin I$.
\begin{itemize}
\item [\emph {a.}]
Prove that the intersection (and a fortiori the product) of the \ids $\fa_I$ is null; consequently, the morphism $\gA \to \gA_{\so\ua}$ is injective and $\Kdim\gA = \Kdim\gA_{\so\ua}$.

\item [\emph {b.}]
Prove that $\Ann_{\gA_{\so\ua}}(a_i) = \gen{\vep_i}_{\gA_{\so\ua}}$.
\end{itemize}

}
\end{exercise}

\vspace{-1em}
\begin{exercise}
\label{exoAmin} (A few results on $\Amin$)%
\index{regular!morphism}%
\index{morphism!regular ---}
~{\rm  See \Pbm \ref{exoQiClot} for $\Aqi$.
\\
 A \ri \homo $\gA\to\gB$ is said to be \emph{\regz} when the image of every \ndz \elt is  a \ndz \eltz.
\fbox{Let $\gA$ be a reduced \riz.}

\vspace{-.5em}
\begin{enumerate}\itemsep0pt
\item Let $\theta:\gA\to\gB$ be a \reg morphism and $a\in\gA$. If  $a\epr$ is generated by an \idm $e$, then $\theta(a)\epr$ is generated by the \idm $\theta(e)$. 
\\
In particular, as already mentioned in \Pbm \ref{exoQiClot}, a morphism between \qiris is a \qiri morphism \ssi it is \regz.
\item The natural morphism $\Aqi \to \Amin$  is \reg and surjective.
\item For $a\in\gA$, the natural morphism $\psi_a:\gA\to\gA_{\so{a}}$ is \regz. 
\item The natural morphism $\psi:\gA\to\Amin$ is \reg and  
the natural morphism $\ZZ\to\ZZ_\mathrm{qi}$ is not \regz.
\end{enumerate}
}
\end{exercise}

\vspace{-1em}
\begin{exercise}
\label{exolemVdimKdim}
{\rm  Explicate the \dem of Lemma~\ref{lemVdimKdim} in terms of \susisz. 
}
\end{exercise}

\vspace{-1em}
\begin{exercise}
\label{exothValDim} (A \gnn of \thref{thValDim})\\
{\rm  
For $\gA\subseteq\gB$ and $\ell\in \NN$, if for every sequence $(\ux)=(x_0,\ldots,x_\ell)$ in $\gB$, we have a primitive \pol of $\gA[\uX]$ which annihilates $(\ux)$, then  $\Vdim\gB\leq\ell + \Vdim\gA$. 
 
}
\end{exercise}

\vspace{-1em}
\begin{exercise}
\label{exoLyingOverClassique} 
{(Lying Over morphism)}
\\
{\rm Prove what is affirmed in the remark following the \dfn of the Lying Over on \paref{remLY}.

}
\end{exercise}

\vspace{-1em}
\begin{exercise}\label{exoLYOV} {(Lying Over morphism, 2)}
\\
{\rm  In the category of finite ordered sets, it is clear that a morphism is surjective \ssi it is an epimorphism.
This therefore corresponds, for the dual \trdisz, to a monomorphism, which here means an injective homomorphism, \cad a Lying Over morphism.
\\
Give a \dem in \clama of the \eqvcz, for some \homo $\alpha:\gT\to\gT'$ of \trdisz, between: $\alpha$ is  Lying Over on the one hand, and $\Spec \alpha:\Spec \gT'\to\Spec\gT$ is surjective, on the other hand. 
\\
Idea: use \emph{Krull's lemma}, which can be easily proven \`a la Zorn:
\emph{If in a \trdi we have an \id $\fa$ and a filter $\ff$ which do not intersect, there exists a \idep containing $\fa$ whose complement is a filter containing $\ff$.
}}
\end{exercise}

\vspace{-1.2em}
\pagebreak	        
\begin{exercise}
\label{exodefiGoingup} {(Going Up, Going Down morphisms)}\\
{\rm  Prove what is stated in the remark following the \dfn of Going Up on \paref{defiGoingup} (use the description of the quotient lattice $\gT\sur{(\fa=0)}$ given on \paref{trquoideal}).
Do the same thing for Going Down.
 }
\end{exercise}



\vspace{-1em}
\begin{problem}
\label{exoAnneauNoetherienReduit}
{(Annihilator of an \id in a reduced \noe \riz)}\\
{\rm
We consider a reduced \ri $\gA$ such that every ascending sequence of \ids of the form $\DA(x)$ has two equal consecutive terms.

\emph{1.}
Let $\fa$ be an \id of $\gA$ such that we know how to test for $y\in\gA$ if $\Ann(y)\fa=0$
(and in case of a negative answer provide the corresponding certificate).
\\
\emph{1a.} If some $x \in \fa$ satisfies $\Ann(x)\fa \ne 0$,
 determine some $x' \in \fa$ such that $\DA(x) \subsetneq \DA(x')$.
\\
\emph{1b.}
Deduce the existence of some $x \in \fa$ such that $\Ann(x) = \Ann(\fa)$.

\emph{2.}
Suppose moreover that every \ndz \elt of $\gA$ is \ivz, and that for all $y$, $z$ we know how to test if $\Ann(y)\Ann(z)=0$.
Show that $\Kdim\gA \le 0$.

\emph{3.}
Let $\gB$ ne a \fdi \coh \noe \riz. Show that 
 $\Frac(\gB\red)$ is a \zed \riz.
\\
Note: in \clama $\gB$ admits a finite number of \idemis $\fp_1$, \dots, $\fp_k$ and $\Frac(\gB\red)$ is \isoc to the finite product of corresponding fields: $\Frac(\gA/\fp_1)\times \cdots\times \Frac(\gA/\fp_k)$. However, in \gnlz, we have no \algq  access to the $\fp_i$'s.
}
\end{problem}

\vspace{-1em}
\begin{problem}
\label{exoContractedInclusion} (Lying Over, Going Up, Going Down, examples)
\\
{\rm
\emph {1.}
Let $\gA\subseteq\gB$ be an inclusion of \ris such that, as an \Amoz, $\gA$ is a direct factor in $\gB$. Show that $\fa\gB \cap \gA = \fa$ for every \id $\fa$ of~$\gA$. In particular, $\gA \hookrightarrow \gB$ is Lying Over.

\emph {2.}
Let $G$ be a finite group acting on a \ri $\gB$ with $\abs G 1_\gB$ invertible in~$\gB$. Let $\gA = \gB^G$ be the sub\ri of fixed points. We define the \emph{Reynolds operator} $R_G : \gB \to \gA$:%
\index{Reynolds operator}
$$
\preskip-.4em \postskip.4em \ndsp
R_G(b) = {1 \over \abs G} \sum_{g \in G} g(b). 
$$
Prove that $R_G$ is an $\gA$-\prr of image $\gA$; in particular, $\gA$ is a direct summand (as an \Amoz) in $\gB$.

\emph {3.}
Let $\gA \hookrightarrow\gB$ with $\gA$ as a direct summand (as an \Amoz) in $\gB$.  Provide a direct \dem of $\Kdim\gA \le \Kdim\gB$.

\emph {4.}
Let $\gk$ be a nontrivial \cdiz, $\gA = \gk[XZ,YZ] \subset\gB = \gk[X,Y,Z]$.
Then $\gA$ is a direct summand in $\gB$, therefore $\gA \hookrightarrow \gB$ is Lying Over. But $\gA \hookrightarrow \gB$ is neither Going Up nor Going~Down.

}
\end{problem}

\vspace{-1em}
\begin{problem}
\label{exoChainesPotPrem} (Potential chains of \idepsz)\index{chain!potential --- of \idepsz}
\\
{\rm Over a \ri $\gA$ we call a \ix{potential chain of \idepsz}, or \emph{potential chain} a list  $[(I_0,U_0),\ldots,(I_n,U_n)]$, where the $I_j$'s and $U_j$'s are subsets of $\gA$ (\cad each $(I_j,U_j)$ is a potential \idep of $\gA$). A potential chain is said to be \emph{finite} if the $I_j$'s and $U_j$'s are finitely enumerated subsets.
\\
A potential chain is said to be \emph{complete} if the following conditions are satisfied

\vspace{-.5em}
\begin{itemize}\itemsep0pt
\item the $I_j$'s are \ids and the $U_j$'s are \mosz,
\item $I_0\subseteq I_1\subseteq \cdots\subseteq I_n$ and  $U_0\supseteq U_1\supseteq \cdots\supseteq U_n$,
\item $I_j+U_j=U_j$ for each $j$.
\end{itemize}

\vspace{-.5em}
We say that the potential chain $[(I_0,U_0),\ldots,(I_n,U_n)]$ \emph{refines} the chain $[(J_0,V_0),\alb\ldots,\alb(J_n,V_n)]$ if we have the inclusions $J_k\subseteq I_k$ and $V_k\subseteq U_k$ for each $k$.

 \emph{1.} Show that every potential chain generates a complete potential chain (in the sense of the refinement relation). More \prmtz, from $[(I_0,U_0),\ldots,(I_n,U_n)]$, we successively construct 

\vspace{-.5em}
\begin{itemize}\itemsep0pt
\item $\fa_j=\gen{I_j}$, $\fb_j=\sum_{i\leq j}\fa_i$ ($j\in\lrb{0..n}$),
\item $\ff_n=\cM(U_n)+\fb_n$,\footnote{Recall that $\cM(A)$ is the \mo generated by the subset $A$.} $\ff_{n-1}=\cM(U_{n-1}\cup \ff_n)+\fb_{n-1}$, \ldots, $\ff_{0}=\cM(U_{0}\cup \ff_1)+\fb_{0}$.
\end{itemize}

\vspace{-.5em}
And we consider $[(\fb_0,\ff_0),\ldots,(\fb_n,\ff_n)]$.  

 \emph{2.} We say that a potential chain $\cC$ \emph{collapses} if in the complete chain that it generates $[(\fb_0,\ff_0),\ldots]$ we have $0\in\ff_0$.
Show that a sequence $(x_1,\ldots,x_n)$ is \sing \ssi the potential chain $[(0,x_1),(x_1,x_2),\ldots,(x_{n-1},x_n),(x_n,1)]$ collapses.

 \emph{3.} Show in \clama that a potential chain $\cC$ of $\gA$ collapses \ssi it is impossible to find \ideps    $\fp_0\subseteq\fp_1\subseteq\cdots\subseteq\fp_n$, such that 
the chain $[(\fp_0,\gA\setminus \fp_0),\ldots,\alb(\fp_n,\gA\setminus \fp_n)]$ refines the chain $\cC$.

 \emph{4.} Given a potential chain $\cC=[(I_0,U_0),\ldots,(I_n,U_n)]$, we \emph{saturate} it by adding, in $I_k$, (resp.\ in $U_k$) every $x\in\gA$ which, added to $U_k$ (resp.\ to $I_k$) would lead to a collapse.
Thus a potential chain collapses \ssi its saturated chain is $[(\gA,\gA),\ldots,(\gA,\gA)]$.
\\
Show that we obtain thus a potential chain $[(J_0,V_0),\ldots,(J_n,V_n)]$ which refines the complete chain generated by~$\cC$.
\\
Show in \clama that $J_k$ is the intersection of the \ideps that appear in position $k$ in a chain of \ideps which refines $\cC$
(as in the previous question). Also prove the dual statement for $V_k$.   
}
\end{problem}

}

\sol{

\exer{exoDimKX}
\emph{2.} Consider a sequence of length $k$ in $\gA$, to it we add $X$ at the start, and it becomes \sing in $\gA[X]$. We then get rid of $X$ in the corresponding \Egref{eqsing}.
\\
Note: we can also invoke item~\emph{3} of Proposition~\ref{lemRegsing}.

\exer{exoDimKX1n}
We can assume that $\gK$ is reduced ($\gK\red[\uX]=\KuX\red$ has the same dimension as $\KXn$).
Two possibilities are then offered. The first is to rewrite the \dem given in the case of a \cdi by using Exercise~\ref{exoZDpiv} and \plgrf{thDdkLoc}.
The second is to apply the \elgbm \num2.


\exer{exoPartitionUnite} 
We write each $U_i$ in the form $U_i = \bigcup_{j \in J_i} \DA(g_{ij})$. Saying that the~$\DA(g_{ij})$'s cover $\Spec(\gA)$ means that $1 \in \gen {\,\DA(g_{ij})\mid j \in J_i,\ i \in I}$,
hence an \egt $1 = \sum_{j,i} u_{ji} g_{ij}$, the $u_{ji}$'s being null except for a finite number of them (\hbox{\cad $i\in I_0$}, $j\in J_i$, where $I_0$ and the $J_i$'s are finite).
Let $f_i = \sum_{j\in J_i} u_{ji}g_{ij}$. We obtain $\DA(f_i) \subseteq U_i$  because for $\fp \in \DA(f_i)$, we \hbox{have $f_i \notin\fp$}, therefore an index $j$ such that $g_{ij} \notin \fp$, \cad $\fp \in \DA(g_{ij})
\subseteq U_i$. And  $\sum_{i\in I_0} f_i = 1$.

\exer{exoKdimGeom} ~\\
\emph{1a.} We write an \rdi of $f$ over $\KYr$ 

\snic{f^n+a_{n-1}f^{n-1}+\cdots+a_kf^k=0,}

with $n\geq1$, the $a_i$'s  in $\KYr$ and $a_k\neq0$. 
The \egt $(a_k+bf)f^{k}=0$ shows that $a_k+bf\in(\DA(0):f)$ (even if $k=0$). Therefore $a_k\in \JK_\gA(f)$.

\exer{exolemLocMemeKdim} \\
\emph{1.} We write $\KuX[1/f]=\aqo{\gK[\uX,T]}{1-fT}$.
A \iNoe position of the nonconstant \pol $1-fT$ brings us to an integral extension of $\gK[\Yn]$.
\\
\emph{2.} We write $\AuX[1/\delta]=\aqo{\gA[\uX,T]}{1-\delta T}$. We seek to apply \thref{cor2thKdimMor} 
 	to integral extensions.
On the one hand we want  $\delta $ to be \ndzz, for the \homo  $\AuX\to\AuX[1/\delta ]$ to be injective, and on the other hand we want to be able to put the \pol $1-\delta T$ into \Noe position, 
             for the \ri $\AuX[1/\delta ]$ to be 
integral over a \ri $\gA[\Yn]$. 
\\
The first condition means that the \id $\rc(\delta )$ is faithful (McCoy, Corollary \ref{corlemdArtin}). 
\\
The second condition is satisfied if we are in the same situation as for Lemma~\ref{lemLocMemeKdim}  

\vspace{-.5em}
\begin{itemize}\itemsep0pt
\item $\delta$ is of formal degree $d$,
\item  one of the \moms of degree $d$, relating to a subset of variables $(X_i)_{i\in I}$, has as its \coe an \elt of $\Ati$,
\item  and it is the only \mom of degree $d$ in the variables $(X_i)_{i\in I}$ present in $\delta$.
\end{itemize}   
Indeed, the \cdv \gui{$X'_i = X_i + T$ if $i \in I$, $X'_i = X_i$ otherwise,} then renders the \pol $1 - \delta T$ \mon in $T$ (up to \invz).
Note that in this case the \pol $\delta$ is primitive and the first condition is \egmt satisfied.


\exer{exoPruferKdim} \emph {1.} Consider $s=a/b\in\Frac\gA$ with $b$ \ndzz.\\
The sequence $(bX-a,b,X)$ is \sing in $\AX$. This gives an \egt in~$\AX$ of the following type
$$
\preskip.3em \postskip.3em 
 (bX-a)^{k_1}\big(b^{k_2}\big(X^{k_3}\big(1+Xp_3(X)\big)+bp_2(X)\big)+(bX-a)p_1(X)\big)=0.
$$
Since $\AX$ is integral, we can delete the factor $(bX-a)^{k_1}$, after which we specialize $X$ in $s$. We get
$$\preskip.3em \postskip.0em 
b^{k_2}\big(s^{k_3}\big(1+sp_3(s)\big)+bp_2(s)\big)=0, 
$$
and since $b$ is \ndzz,
$$\preskip-.20em \postskip.2em 
s^{k_3}\big(1+sp_3(s)\big)+bp_2(s)=0.
$$
Thus $s$ annihilates 
$g(X)=X^{k_3}\big(1+Xp_3(X)\big)+bp_2(X)$
\hbox{and $f(X)=bX-a$}.\\
Finally, since the \coe of $X^{k_3}$ in $g$ is of the form $1+bc$, we obtain \hbox{that $1\in\rc(f)+\rc(g)=\rc(f+X^2g)$}.

\emph{2.} Results from \emph{1} and from the \gnl results on the dimension of $\AX$, for an arbitrary \ri and for a \adpz.

\emph{3.} The answer seems to be yes.


\exer{exoMultiplicativiteIdeauxBords} 
\emph {1.}
It suffices to show, for two \ids $\fa$, $\fb$ and two \eltsz~\hbox{$u$, $v \in \gA$}, that
$$\preskip-.6em \postskip.2em
\begin{array} {c}
\big((\fa : u^\infty) + \gA u\big)\,\big((\fb : u^\infty) + \gA u\big) \subseteq
(\fa\fb : u^\infty) + \gA u \et
\\[.2em]
\big((\fa : u^\infty) + \gA u\big) \, \big((\fa : v^\infty) + \gA v\big) \subseteq
(\fa : (uv)^\infty) + \gA uv.
\end{array}
$$
The first inclusion stems from $(\fa : u^\infty)\,(\fb : u^\infty) \subseteq (\fa\fb : u^\infty)$ and the second from $(\fa : u^\infty) + (\fa : v^\infty) \subseteq (\fa : (uv)^\infty)$.


\exer{exoInclusionBordLionel} 
\emph {1.}
As $\beta > \alpha$, $\uX^\beta$ is a multiple of one of the following \moms

\snic {
X_1^{\alpha_1} X_2^{\alpha_2} \cdots X_{n-1}^{\alpha_{n-1}}X_n^{1+\alpha_n},
\;
X_1^{\alpha_1} X_2^{\alpha_2} \cdots X_{n-1}^{1+\alpha_{n-1}},
\; \dots,\; 
X_1^{\alpha_1}  X_2^{1+\alpha_2},\;
X_1^{1+\alpha_1}.
}

\emph {2a.}
Let $y\in \prod_{\beta < \alpha} \Ann(a_\beta)$; by letting $Q(\uX) = yP(\uX)$, we have

\snic{Q(\uX) = ya_\alpha\uX^\alpha + \sum_{\beta>\alpha} ya_\beta
\uX^\beta\quad$ and $\quad Q(\ux) = 0.}

To show that $ya_\alpha\in\IK(\ux)$, we can therefore suppose that we have $y=1$ and~\hbox{$P(\uX) = a_\alpha\uX^\alpha + \sum_{\beta>\alpha} a_\beta\uX^\beta$}.
By using the \egt $P(\ux) = 0$ and the first question, we obtain $a_\alpha \in \IK(\ux)$.

\emph {2b.}
First, since $\gA$ is reduced, we have $\IK(a) = \Ann(a) + \gA a,\;\forall a \in \gA$. Next, we use the following remark: let $\fc$ be an \id and $2m$ \ids $\fa_1$, $\fb_1$, \ldots, $\fa_m$, $\fb_m$ such that $\fa_1\cdots\fa_{k-1}\fb_k\subseteq \fc$ for every $k \in \lrbm$. Then we obtain the inclusion
$$\preskip.2em \postskip.2em 
(\fa_1+\fb_1) \cdots (\fa_m+\fb_m) \subseteq \fc + \fa_1\cdots\fa_m. 
$$

Indeed, by \recu on $m$, if
$\mathrigid 2mu
(\fa_1+\fb_1) \cdots (\fa_{m-1}+\fb_{m-1}) \subseteq \fc + \fa_1\cdots\fa_{m-1}$, we deduce
$$
\preskip-.1em \postskip.2em\mathrigid 2mu 
\fa_1+\fb_1) \cdots (\fa_m+\fb_m) \subseteq \fc + 
\fa_1\cdots\fa_{m-1}\fa_m + \fa_1\cdots\fa_{m-1}\fb_m
\subseteq \fc  + 
\fa_1\cdots\fa_{m-1}\fa_m+\fc, 
$$
hence the stated inclusion.
\\
Let us apply this to $\fc = \IK(\ux)$ and to the \ids $\fa_\beta = \Ann(a_\beta)$,
$\fb_\beta = \gA a_\beta$.  \\
As $\Ann(a_\beta) + \gA a_\beta = \IK(a_\beta)$, we obtain the desired inclusion.

\emph {3.}
Direct application with $n = 1$.

\emph {4.}
We can suppose that $\gA$, $\gB$ are reduced even if it entails replacing $\gA \to \gB$ with $\gA\red \to \gB\red$  (every $z\in\gB$ remains primitively \agqz). We can also suppose that $\gA \subseteq \gB$ even if it entails replacing~$\gA$ with its image in $\gB$. Let us show that $\Kdim\gA \le m \Rightarrow \Kdim\gB \le m$ by \recu on $m$. It suffices to show, for $x \in \gB$, that $\Kdim(\gB\sur{\IK_\gB(x)}) \le m-1$; but~$\IK_\gB(x)$ contains an \id $\fa$ of $\gA$, finite products of boundary \ids $\IK_\gA(a)$, $a\in \gA$.
\\
We therefore have an \alg $\gA\sur\fa \to \gB\sur{\IK_\gB(x)}$ to which we can apply the \hdr since $\Kdim \gA\sur\fa \le m-1$.

\exer{exoExtEntiereIdealBord} 
\emph {1.}
We use the integral extension $\ov\gA = \gA\sur{\gA\cap\fb} \hookrightarrow \ov\gB = \gB\sur{\fb}$. 
\\
Let $a \in \gA \cap (\fb + \fa\gB)$;
the Lying Over (\ref{lemLingOver}) with $\ov\gA \subseteq \ov\gB$, gives $\ov a^n \in \ov\fa$, \hbox{\cad $a^n \in \fa + \fb$} and as $a \in \gA$, $a^n \in \fa + \gA\cap\fb$.

\emph {2.}
By \recu on $d$. Let

\snuc {
\fa = \IK_\gA(a_0, \ldots, a_{d-1}), \
\fa' = \IK_\gA(a_0, \ldots, a_d), \
\fb = \IK_\gB(a_0, \ldots, a_{d-1}), \
\fb' = \IK_\gB(a_0, \ldots, a_d)
.}

We therefore have by \dfn $\fa' = (\fa : a_d^\infty)_\gA + \gA a_d$ and $\fb' = (\fb : a_d^\infty)_\gB + \gB a_d$. We want to show that $\gA\cap\fb' \subseteq \rD_\gA(\fa')$. Item~\emph{1} gives $\gA\cap\fb' \subseteq \rD_\gA(\fc)$ \hbox{with $\fc = \gA a_d + \gA \cap (\fb : a_d^\infty)_\gB = \gA a_d + (\gA\cap\fb : a_d^\infty)_\gA$}. 
\\
By \recuz, $\gA\cap\fb
\subseteq \rD_\gA(\fa)$, therefore
$$\preskip.0em \postskip.4em
\fc \subseteq \gA a_d + (\rD_\gA(\fa) : a_d^\infty)_\gA 
\subseteq  \rD_\gA(\gA a_d + (\fa : a_d^\infty)_\gA)
\eqdf {\rm def} \rD_\gA(\fa')
,
$$
hence $\gA\cap\fb' \subseteq \rD_\gA(\fa')$. 

\exer{exoExtEntiereMonoideBord} 
\emph {1.}
Let $t \in S + \fa\gB$; \cad   $t+s \in \fa\gB$ with $s \in S$. So, $t+s$ is integral over $\fa$, \cad is a zero of a \polu 
${P(X) \in X^n + \fa X^{n-1} + \cdots +
\fa^{n-1} X + \fa^n .}$
We write $P(T+s) = TQ(T) + P(s)$. 
Thus $P(s) \in s^n + \fa$ and $tQ(t) \in S + \fa$. 

\emph {2.}
By \recu on $d$.  Let $V = \SK_\gB(a_0, \ldots, a_d) = a_0^\NN(\SK_\gB(a_1, \ldots, a_d) + a_0\gB)$; the \recu provides $\SK_\gB(a_1, \ldots, a_d) \subseteq \satu{\SK_\gA(a_1, \ldots, a_d)}{\gB}$ so

\snic {
V \subseteq 
a_0^\NN( \satu{\SK_\gA(a_1, \ldots, a_d)}{\gB} + a_0\gB) 
\subseteq
a_0^\NN \satu{(\SK_\gA(a_1, \ldots, a_d) + a_0\gB)}{\gB}
.}

The first question provides

\snic {
V \subseteq a_0^\NN \satu{(\SK_\gA(a_1, \ldots, a_d) + a_0\gA)}{\gB}
\subseteq \satu{(a_0^\NN(\SK_\gA(a_1, \ldots, a_d) + a_0\gA)\big)}{\gB}
,}

\cad $V \subseteq \satu{\SK_\gA(a_0, \ldots, a_d)}{\gB}$.


\exer{exoKdimSomTr} 
\emph{2.} The quotient \ri $\aqo\gA{X_1-Y_1}$ can be seen as the \lon of~$\gK[\Xn,Y_2,\ldots,Y_m]$ at the \mo 

\snic{S_1=(\KXn)\etl(\gK[X_1,Y_2,\ldots,Y_m])\etl.}

It is therefore integral. In the same way, we describe the successive quotients.


\exer{exoDualiteBords} 
\emph {1.}
Let $\fa_i := \IK(x_1, \ldots, x_i)$, with $\fa_{i+1} = (\fa_i : x_{i+1}^\infty) + \gA x_{i+1}$. By \recuz, $\fa_i \subseteq \fp_i$:  $x_{i+1}\notin\fp_i$ gives $(\fp_i : x_{i+1}^\infty) \subseteq \fp_i$, then $\fa_{i+1} \subseteq \fp_i + \gA x_{i+1}$, so $\fa_{i+1}\subseteq \fp_{i+1}$. The rest poses no difficulties.

\emph {2.}
By letting $S_i = \SK(x_{i+1}, \ldots, x_d)$, we have $S_d = 1$ and $S_{i-1} = x_i^\NN(S_i + \gA x_i)$. Step by step, we prove 
$S_i \subseteq \ff_i$  by using $x_i \in \fp_i$ and $x_i \in f_{i-1}$:

\snic {
S_{i-1} = x_i^\NN(S_i + \gA x_i) \subseteq x_i^\NN(\ff_i + \fp_i) = x_i^\NN \ff_i 
\subseteq x_i^\NN\ff_{i-1} \subseteq \ff_{i-1}
.}

The rest poses 
          no difficulties.


\exer{exoEliminationEtBord}  
If $f$ is \mon of degree $n \ge 1$, the \pol $R(X,Y)$ defined in the statement of the question is $Y$-\mon of degree $n$,
therefore $\Ann(R) = 0$, and $R(f,g) = 0$ because 
$R \in \gen {f(T)-X, Y-g(T)}_{\gA[T,X,Y]}$.

\emph {1.} 
Let $f = \sum_{k=0}^n a_kT^k$. By Lemma~\ref{lemQI} there exists a \sfio $(t_n, t_{n-1}, \ldots, t_0, t_{-1})$ such that:
\\
-- in the component $t_k = 1$
for $k \in\lrbzn$, we have $a_i = 0$ for $i > k$ and 
                   $a_k$  \ndzz;\\
-- in the component $t_{-1} = 1$, we have $f = 0$, \cad $t_{-1}f = 0$ and even $\Ann(f) = \gen {t_{-1}}$.  
\\
Let $m$ be the formal degree of $g$. For $1 \le k \le n$, we let

\snic {
R_k(X,Y) = t_k\Res_T(t_kf(T) - X, k, Y-g(T), m).
}

We define $R_0(X,Y) = t_0(t_0f(T)-X)$ and $R_{-1}(X,Y) = t_{-1}X$.  For $k \in \lrb{-1..n}$, we have $\Ann(R_k) = \gen {1-t_k}$ and $R_k(f,g) = 0$. Thus by letting $R = \sum_{k=-1}^n R_k(X,Y)$, we have $\Ann(R) = 0$ and $R(f,g) = 0$.  

\emph {2.}
Direct application of the referenced exercise.

\emph {3.}
By \recu on the Krull dimension of $\gA$. We can replace $\gA$ by a \ri $\gA' := \gA_{\so\ua}$ such that the annihilator (in $\gA'$) of each \coe of $f$ is generated by an \idm (recall that $\Kdim\gA = \Kdim\gA'$). \\
Then, if $\fa=\prod_{i,j}\IK_\gA(r_{ij})$,  the \ri $\gA[T]\sur{\IK_\AT(f,g)}$ is a quotient of $(\gA\sur\fa)[T]$. As $\Kdim(\gA\sur\fa) < \Kdim\gA$, we obtain by \hdr

\snic {
\Kdim (\gA\sur\fa)[T] \le 1 + 2\Kdim(\gA\sur\fa) \le
1 + 2(\Kdim\gA - 1)$,\ \  then $ \
}

\snic {
\Kdim \gA[T] \le 2 + \Kdim \gA[T]\big/{\IK_\AT(f,g)} \le 
2 + 1 + 2(\Kdim\gA - 1) = 1 +  2\Kdim\gA. 
}

\emph {4.}
We preserve the notations of the previous questions. Each $\IK_{\gA'}(r_{ij})$ contains a finite product of boundary \ids of $\gA$ (Exercise~\ref{exoInclusionBordLionel}) therefore the product of the $\IK_{\gA'}(r_{ij})$'s contains an \id $\fa$ of $\gA$, a finite product of boundary \ids of $\gA$. \\
Thus $\fa \subset \gA[T] \cap \IK_{\gA'[T]}(f,g) \subseteq \rD_{\gA[T]}(\IK_{\gA[T]}(f,g)\big)$ (Exercise~\ref{exoExtEntiereIdealBord}).

\exer{exoIdealBordPolynomes}
\emph {1.}
By \recu on $n$, the case $n=0$ being the hypothesis. Let us add an \elt $b$ to $\bn$ and let $\fb'_j = \IK(z_j, \bn, b)$.
\\
By \dfn $\fb'_j = \gB b + (\fb_j : b^\infty)$ with $\fb_j = \IK(z_j, \bn)$; the product of the $\fb_j$'s is contained in $\IK(x,y,\bn)$ (by \recuz). By using inclusions of the type $(\fb : b^\infty)(\fb' : b^\infty) \subseteq (\fb\fb' : b^\infty)$, we obtain

\snic {\arraycolsep2pt
\begin {array}{rcl}
\prod_j \fb'_j &\subseteq& \gB b + \prod_j (\fb_j : b^\infty)
\subseteq \gB b + \big( \prod_j \fb_j : b^\infty\big)
\\[.4em]
&\subseteq& \gB b + (\IK(x,y,\bn) : b^\infty) = \IK(x,y,\bn,b).
\\
\end{array}
}

\emph {2a.} Results from the fact that for two \ids $\fa, \fb$ of $\gA$, we have 

\snic{(\fa : \fb)_\gA\,\AT = (\fa : \fb)_{\AT}.}

\emph {2b.}
For two \ids $\fa, \fb$, let $\fa \Subset \fb$ for $\fa \subseteq\rD(\fb)$. We reason by \recu on $d$,
the case $d = 1$ appearing in Exercise~\ref{exoEliminationEtBord}. 
\\
Consider $2(d+1)$ \pols $p$, $q$,  $g_1$, \ldots, $g_{2d} \in \gA[T]$. There exist  $a_j$'s $\in \gA$
such that $\prod_j \IK_\gA(a_j) \Subset \IK_{\gA[T]}(p,q)$ (the case $d = 1$).
By the first question,

\snic {
\prod_j \IK_{\gA[T]}(a_j, g_1, \ldots, g_{2d}) \Subset 
\IK_{\gA[T]}(p,q, g_1, \ldots, g_{2d})
.}

It suffices therefore to show, for $a \in \gA$, that a boundary \id $\IK_{\gA[T]}(a, g_1, \ldots, g_{2d})$ contains, up to radical, a product of boundary \ids of $d+1$ \elts of $\gA$.  
Let $\ov\gA= \gA\sur{\IK(a)}$ and $\varphi : \gA[T] \to \ov\gA[T] \simeq \gA[T]\sur{(\IK_{\gA[T]}(a)\big)}$ be the \homo of passage to the quotient. By \recuz, the boundary \id $\IK_{\ov\gA[T]}(\ov {g_1}, \ldots, \ov{g_{2d}})$ contains, up to radical, a product $\prod_j \fa_j$ where each $\fa_j$ is a boundary \id of $d$ \elts of~$\ov\gA$. By taking the inverse image under $\varphi$, we obtain

\snic {
\prod_i \varphi^{-1}(\fa_i) \subseteq
\varphi^{-1}\big(\prod_i \fa_i\big) \Subset
\varphi^{-1}\big(\IK_{\ov\gA[T]}(\ov {g_1}, \ldots, \ov{g_{2d}})\big)
.}

By using Lemma~\ref{lem2BordKrullItere}, we have on the one hand

\snic {
\varphi^{-1}\big(\IK_{\ov\gA[T]}(\ov {g_1}, \ldots, \ov{g_{2d}})\big) =
\IK_{\gA[T]}(a, g_1, \ldots, g_{2d})
,}

and on the other hand $\varphi^{-1}(\fa_i)$ is a boundary \id of $d+1$ \elts of $\gA$ (the first \elt being $a$).  This shows that $\IK_{\gA[T]}(a, g_1, \ldots, g_{2d})$ contains up to radical, a product of boundary \ids of $d+1$ \elts of $\gA$.

\emph {3.}
If $\gA[\uT] = \gA[T_1, \ldots, T_r]$, the boundary \id of $(r+1)d$ \pols of $\gA[\uT]$ contains, up to radical, a product of boundary \ids of $d$ \elts of $\gA$.  
\\
Consequently, $\Kdim\gA < d \Longrightarrow \Kdim\gA[\uT] < (r+1)d$,
\cad

\snic{\Kdim\gA[\uT]+1\, \le \,(r+1)(\Kdim\gA +1) .}

\exer{exoKdimEspanol}
\emph{1.} 
We take $a'=y\vi a$ and $b'=x\vu b\vu a'$.  Then $x \vi a' = x \vi y \vi a = x \vi a$ (because $x \le y$).  Then $y \vu b' = y \vu x \vu b \vu a' = (x \vu a') \vu (y \vu b) = y \vu b$ (the last \egt uses $x \vu a' \le y$ which stems from $x \le y$ and $a' \le y$, a fortiori $x \vu a' \le y \vu b$).
\\
Ir remains to see that $y \vi b' = x \vu a'$; we have the \idt for all $y, b, z$, $y \vi (b \vu z) = y \vi z'$ with $z' = (y \vi b) \vu z$ that we use with $z = x \vu a'$.  But we have $y \vi b \le x \vu a'$ because the hypothesis is $y \vi b \le x \vu a$, so $y \vi b \le (x \vu a) \vi y = (x \vi y) \vu (y \vi a) \le x \vu a'$. Therefore $z' = x \vu a'$ and $y \vi b' = y \vi (x \vu a')$.
Finally, $y \vi (x \vu a')=x \vu a'$ because $x \vu a' \le y$
(by using $x \le y$ and $a' \le y$).

\emph{2.} By \emph{1} by \recu on $n$.

\emph{3.} Item \emph{3a} implies item \emph{3c} by item \emph{2}. Item \emph{3c} implies item~\emph{3b} because a linked chain is a particular case of \cop sequence. 
In order to see that \emph{3b} implies  \emph{3a}, 
 \hbox{let $y_0,\dots,y_{n}$} be arbitrary. We define
$x_0=y_0$, $x_i=y_i\vu x_{i-1}$  ($i\in\lrbn$).
Let $(a_0,\dots,a_{n})$ be a complementary sequence of $(x_0,\dots,x_{n})$.
We define $b_0=a_0$  and $b_{i} = a_{i}\vu  x_{i-1}$  \hbox{for $i\in\lrbn$}. We have $x_{i}\vu  a_i=y_i\vu b_i$ for $i\in\lrb{0..n}$. Thus $0 = x_0\vi a_0=y_0\vi b_0$ and $1=x_{n}\vu a_{n}=y_n\vu b_n$. Let us see the intermediary inequalities. For $i\in\lrb{1..n}$ 
we have $x_{i}\vi a_{i}\leq  x_{i-1}\vu a_{i-1}$, so
\[ 
y_{i} \vi a_{i} \leq  x_{i} \vi a_{i}  \leq  x_{i-1}\vu a_{i-1}= y_{i-1}\vu b_{i-1}
\]
Then we have
\[ y_{i}\vi b_{i}=  y_{i}\vi (a_{i}\vu  x_{i-1})=
(y_{i}\vi a_{i})\vu (y_{i}\vi x_{i-1})\leq (y_{i}\vi a_{i})\vu x_{i-1}
\]
As the two terms after $\leq $ are bounded by $x_{i-1}\vu a_{i-1}= y_{i-1}\vu b_{i-1}$ we get the inequality $y_{i}\vi b_{i}\leq y_{i-1}\vu b_{i-1}$.


\exer{exoAminEtagesFinis} 
First of all, for every \id $\fc$, the \ri $\gA\sur{\fc\epr}$ is reduced.  \\
Let us show that $(\fa_1\epr\fa_2\epr)\epr = (\fa_1+\fa_2)^{\perp\perp}$: the \egt $\fa_1\epr\cap\fa_2\epr = (\fa_1+\fa_2)\epr$ implies that the \ids $\fa_1\epr\cap \fa_2\epr$, $\fa_1\epr\fa_2\epr$ and $(\fa_1+\fa_2)\epr$ have the same nilradical therefore the same annihilator. 
\\
We deduce that

\snic {
(\fa_1\epr\fa_2\epr\fb)\epr = (\fa_1+\fa_2) \diamond \fb.
}

\vspace{-.3em}
Indeed

\snic {
(\fa_1\epr\fa_2\epr\fb)\epr = \big((\fa_1\epr\fa_2\epr)\epr : \fb\big) =
\big((\fa_1+\fa_2)^{\perp\perp} : \fb\big) = (\fa_1+\fa_2) \diamond \fb.
}

\emph {1.}
As $\fa\epr\fb \subseteq \fa\epr$, we have $\fa\diamond\fb \supseteq \fa^{\perp\perp} \supseteq \fa$. Let $x\in \gA$ such that in the quotient we have $\ov x\,\ov\fb = 0$, that is $x\fb \subseteq \fa\diamond\fb$, \cad  $x\fb\fa\epr\fb = 0$. We therefore have $x\fb\fa\epr = 0$, that is $x\in \fa\diamond\fb$, \cad $\ov x =
0$.

\emph {2.}
We have

\snic {
(\gA\sur{\fa_1\diamond\fb_1})\sur{(\ov{\fa_2}\diamond\ov{\fb_2})} \simeq
\gA\sur{(\fa_1\epr\fa_2\epr\fb_1\fb_2)\epr} =
\gA\sur{\big((\fa_1+\fa_2)\diamond(\fb_1\fb_2)\big)}. 
}


\exer{exoAmin}
\emph{1.} Let $y\in\gB$ and suppose $y\theta(a)=0$. Let $e'=\theta(e)$.
We must show that $y=ye'$. Since $e+a$ is \ndzz, $e'+\theta(a)$ is \ndzz.
\hbox{However, $y(e'+\theta(a)\big)=ye'=ye'(e'+\theta(a)\big)$} because $e'$ is \idmz.

\emph{2.} The \homo  $\Aqi \to \Amin$ comes from the \uvl \prt of $\Aqi$.
It is surjective because $\Amin=\gA[(e_x)_{x\in\gA}]$ and because the morphism
$\Aqi \to \Amin$ is a \qiriz.

\emph{3.} Let $x$ be \ndz in $\gA$ and $u=(\ov y, \wi z)\in\gA_{\so{a}}=\gA\sur{a\epr}\times \gA\sur{({a\epr})\epr}$,
with $ux=0$. We must show that $u=0$, \cad $\ov y=\ov 0$ and $\wi z=\wi 0$. 
\\
We have $xy\in a\epr$, \hbox{\cad $xay=0$}, so $ay=0$, then $\ov y=\ov 0$. 
\\
To see that $\wi z=\wi0$ we consider an arbitrary \eltz~$t$ of $a\epr$ and we must show that $zt=0$.
However, $\wi{xz}=\wi0$, so $xzt=0$, then $zt=0$.

\emph{4.} If $a\in\gA$ is \ndzz, it remains \ndz at the finite stages of the construction of~$\Amin$ by item \emph{3} 
and this is sufficient for it to be \ndz in $\Amin$.
If the natural \homo $\ZZ\to\ZZ_\mathrm{qi}$ were \reg all the \homos from $\ZZ$ to \qiris would be \regs given the \uvl \prt of $\ZZ_\mathrm{qi}$. However, the surjection $\ZZ\to\aqo\ZZ n$ is not a \reg \homo for $n\geq2$.
Note that the argument applies to every \ri $\gA$ for which there exists a \ndz \elt $x$ such that $\aqo\gA x$ is a \qiri and is nontrivial.

\exer{exolemVdimKdim}
Let us write the computation for $n=k=2$.
\\
Let $x_1={a_1\over b_1}$, $x_2={a_2\over b_2}\in\Frac\gA$ and $s=\big(P(x_1,x_2),(Q(x_1,x_2),(R(x_1,x_2)\big)$ be a sequence in $\gA[x_1,x_2]$, with $P$, $Q$, $R\in\gA[X_1,X_2]$.
We must show that the sequence~$s$ is singular. 
Let $\gA_1=\gA[x_1]$. We know that the sequence

\snic{(b_1X_1-a_1,b_2X_2-a_2,P,Q,R)=(f_1,f_2,P,Q,R)}

is \sing in $\gA[X_1,X_2]$, which gives an \egt

\snic{f_1^{m}(f_2^{m}(P^{m}(Q^{m}(R^{m}(1+AR)+BQ)+CP)+Df_2)+Ef_1)=0}

in $\gA[X_1,X_2]$. Since $b_1\in\Reg\gA$, we have $f_1\in\Reg\gA[X_1,X_2]$ (McCoy's lemma, Corollary~\ref{corlemdArtin}). We therefore simplify the \egt by $f_1^{m}$, then we evaluate it in $\gA_1[X_2]$ by the morphism $X_1\mapsto x_1$. We obtain the following \egt in  $\gA_1[X_2]$

\snic{f_2^{m}(p^{m}(q^{m}(r^{m}(1+ar)+bq)+cp)+df_2)=0,}

with $p=P(x_1,X_2)$,  $q=Q(x_1,X_2)$, \dots, $d=D(x_1,X_2)$.\\
Since $b_2\in\Reg\gA_1$, we have $f_2\in\Reg\gA_1[X_2]$. We can therefore simplify the \egt by~$f_2^{m}$, then evaluate it in $\gA[x_1,x_2]$ by the morphism $X_2\mapsto x_2$.
We obtain an \egt which says that the sequence $s$ is singular.

\exer{exoLyingOverClassique}
Let \emph{a}, \emph{b} and \emph{c} be the three \prts for the commutative \risz. The \eqvc of \emph{a} and \emph{b} is easy.
The implication \emph{a} $\Rightarrow$ \emph{c} has been given as a remark after the Lying Over (Lemma~\ref{lemLingOver}).
\\
 \emph {c} $\Rightarrow$ \emph {a.}
In \clama $\DA({\fa})$ is the intersection of the \ideps that contain $\fa$. We therefore want to show that for every prime \id $\fp$ such \hbox{that $\fa\subseteq\fp$}, we have $\varphi^{-1}(\gen {\varphi(\fa)}) \subseteq\fp$. Let $\fq$ be a prime \id of $\gB$ above $\fp$, \hbox{\cad $\varphi^{-1}(\fq)=\fp$}. Then, $\gen {\varphi(\fa)} \subseteq \gen {\varphi(\fp)} \subseteq \fq$, hence $\varphi^{-1}(\gen {\varphi(\fa)}) \subseteq\fp$.



\prob{exoAnneauNoetherienReduit}
\emph {1a.}
We have some nonzero $a \in \Ann(x)\fa$, in particular $ax = 0$. 
\\
Let us show that $a \notin \rD(x)$: if $a^n \in \gen {x}$, then $a^{n+1} \in \gen {ax} = 0$, and so $a = 0$. 
\\
Therefore $\rD(x) \subsetneq \rD(x,a) = \rD(ax, a+x) = \rD(a+x)$: we take $x' = a+x$ (which is indeed in $\fa$).

\emph{1b.}
Let $x_0 = 0$. If $\Ann(x_0)\fa = 0$, that is $\fa=0$, then $\Ann(x_0) \subseteq \Ann(\fa)$,
\hbox{so $\Ann(x_0) = \Ann(\fa)$}. In this case we let $x_i=x_0$ for every $i\geq0$. Otherwise, there is some $x_1 \in \fa$ with $\rD(x_0) \subsetneq \rD(x_1)$.  If $\Ann(x_1)\fa = 0$, then $\Ann(x_1) \subseteq \Ann(\fa)$, \hbox{so $\Ann(x_1) = \Ann(\fa)$}. In this case we let $x_i=x_1$ for every $i\geq1$. Otherwise, there is some $x_2 \in \fa$ with $\rD(x_1) \subsetneq \rD(x_2)$ \ldots\,\ldots
\\
This way we construct a non-decreasing infinite sequence of \ids $\rD(x_i)$, which is stationary as soon as two consecutive terms are equal, in which case the initial \pb is solved.\footnote{The proof that the \algo terminates under the \cov \noe hypothesis which has just been given is a little confusing. Spontaneously we would have preferred to say: 
the \algo needs to end some day because otherwise, we would have a strictly increasing infinite sequence. The problem with this last argument is that it is an argument by contradiction.
Here we have used the \noe hypothesis in \cov form and this provided us with the means to know a priori when the \algo will terminate.
This delicate point sends us back to the discussion about the Markov principle (Annex~\paref{principeMarkov}).}

\emph{2.}
Let $y\in\gA$.
By the hypothesis, we apply item \emph{1} with the \id $\fa=\Ann(y)$
and we know how to determine some $x \in \Ann(y)$ such that $\Ann(x) = \Ann(\Ann(y)\big)$, \cad $xy =
0$ and $\Ann(x)\Ann(y) = 0$.
We then have $(\Ann(y) \cap \Ann(x)\big)^2 \subseteq \Ann(x)\Ann(y) = 0$, \hbox{so $\Ann(x) \cap \Ann(y) = 0$} (the \ri is reduced). Let us show that $x+y$ is \ndzz; suppose $z(x+y) = 0$. By multiplying by $y$, $zy^2 = 0$, so $zy = 0$, then $zx = 0$, so $z \in \Ann(x) \cap \Ann(y) = 0$. Consequently, $x + y$ is \iv and this \elt is
in the boundary \id of $y$ since $x\in\Ann(y)$.

\emph{3.}
For every \ri $\gC$, every \ndz \elt of $\Frac(\gC)$ is \ivz.
We can apply the result of item \emph{2} to the \ri $\gC=\Frac(\gB\red)$.
\\
Indeed, the first hypothesis that needs to be checked is that every ascending sequence of \ids of the form $\rD_\gC(x_n/y_n)$ ($x_n\in\gB\red$, $y_n\in\Reg(\gB\red)$) admits two equal consecutive terms.
However, in $\gC$ we have the \egt $\rD_\gC(x_n/y_n)=\rD_\gC(x_n)$, and \trf by the fact that in $\gB$, the ascending sequence $\gen{\xzn}_\gB$ admits two equal consecutive terms.\\
The second hypothesis is that we know how to test, for $\fraC x u,\fraC y v\in\gC$, 

\snic{\Ann\big(\fraC x u\big)\Ann\big(\fraC y v\big)=0?}

which is the same thing as $\Ann(x)\Ann(y)=0$ in $\gB\red$. 
However,  in $\Zar\gB$ we have the \egt $\Ann_{\gB\red}(x)=\rD_\gB(x)\im \rD_\gB(0)$,
 and we know that $\Zar\gB$ is a discrete \agH (Proposition~\ref{propNoetAgH}). 



\prob{exoContractedInclusion}
\emph{1.} Let $\pi : \gB \to \gA$ be an $\gA$-\prr of image $\gA$. \\
Let
$a \in \fa\gB \cap\gA$, $a = \sum_i a_ib_i$ with $a_i \in \fa$,
$b_i \in \gB$; so $a = \pi(a) = \sum_i a_i\pi(b_i) \in \fa$.

 \emph{2.}
It is clear that $R_G$ is $\gA$-\lin and that $R_G(a) = a$ for all $a \in \gA$. The rest stems from this.

 \emph{3.}
Let us suppose that $\Kdim\gB \le d$ and show that $\Kdim\gA \le d$. 
\\
Let $a_0$, \ldots, $a_d \in \gA$; as $\Kdim\le d$, there exists an $n \ge 0$ such that
$$\preskip.2em \postskip.2em
(a_0 \ldots a_d)^n \in \gen {c_{d}, c_{d-1}, \cdots, c_0}_\gB
\quad \hbox {with} \quad
c_i = (a_0 \ldots a_{i-1})^n a_{i}^{n+1}.
$$
But $\gen {c_{d}, c_{d-1}, \cdots, c_0}_\gB \cap \gA = \gen {c_{d}, c_{d-1},
\cdots, c_0}_\gA$. Therefore $\Kdim\gA \le d$.

 \emph{4.} (\Demo in \clamaz)
\\
We graduate $\gB$ by $\deg X=\deg Y=1$ and $\deg Z = -1$. Then $\gA$ is the \hmg component of degree $0$, so is a direct summand in $\gB$.
\\
Let $\fq' =  \gen{Z}_\gB$  (it is a \idepz) and $\fp' := \gA \cap \fq' = \gen {XZ,YZ}$.  \\
Let $\fp =  \gen{XZ}_\gA$; it is a \idep with $\fp \subset \fp'$ but there does not exist a \idep $\fq$ of $\gB$ contained in $\fq'$ and above $\fp$. Thus $\gA\subseteq\gB$ is not Going~Down.
\\
Let $\fq = \gen {X,Y^2Z-1}_\gB$ (it is a \idepz) and $\fp := \gA \cap \fq = \gen {XY}$.  
\\
Let $\fp' = \gen {XZ, YZ}_\gA$; it is a \idep with $\fp \subset \fp'$ but there does not exist a \idep $\fq'$
of $\gB$ containing $\fq$ and lying over $\fp'$ (a \idep lying over $\fp'$ must contain $Z$, or $X$ and $Y$). Thus $\gA\subseteq\gB$ is not Going~Up.

}

\Biblio

A very good presentation of the \ddk from the point of view of \clama is found in \cite{Eis}.

The spectral spaces were introduced by Stone \cite{Sto} in 1937.
The theory of spectral spaces is at the heart of this book \cite{Johnstone}.

An important Hochster \tho \cite{Hoc} states that every spectral space is homeomorphic to the spectrum of a commutative \riz.
A pointless version of the Hochster \tho is: every \trdi is \isoc to the Zariski lattice of a commutative \ri (for a non\cov \dem see~\cite[Banaschewski]{Ba}). The delicate point is to know how to construct a \ri whose Zariski lattice is a given finite ordered set.

The \cov \dfn of the \ddk of \trdis and commutative \ris dates back to the works of Andr\'e Joyal \cite{Joyal71,Joyal} and Luis Espa\~nol~\cite{espThesis,esp,esp83,esp86,Es,Espa08}.
Joyal's idea was to construct for each integer $\ell\geq1$, from the \trdi $\gT$, a \trdi $\gT_\ell$,
that satisfies an adequate \uvl \prt such that, in \clamaz, the \ideps of $\gT_\ell$ can be identified with the chains $\fp_0\subseteq\cdots\subseteq\fp_\ell$ of \ideps of~$\gT$ (the inclusions being not \ncrt strict).
This then allows for a \dfn of $\Kdim\gT\leq\ell$ by means of a \prt relating $\gT_{\ell+1}$ and~$\gT_\ell$.
Finally, the \ddk of a commutative \ri can be defined in the same way as the dimension of its Zariski lattice. For further details on the subject, see the \hbox{articles \cite[Coquand\&al.]{cl,cl1,clq2}}.

\Thref{thDKA}, which gives an \elr inductive \carn of the \ddk of a commutative \riz, is found in \cite[Coquand\&al.]{clr03}.
The \carn in terms of \idas given in Proposition~\ref{corKrull} are found in \cite[Coquand\&Lombardi]{cl} and~\cite[Lombardi]{lom}.
\\
Even though the \carn in terms of \cop sequences is already present for the \trdis in \cite{cl}, it only appears for the commutative \ris in \cite[Coquand\&al.]{clq2}.

Additional results on the treatment of the \ddk in integral extensions are found in \cite[Coquand\&al.]{CDLQ06}

The classical theory of valuative dimension can be found in \cite{Jaffard2},
\cite[Chap. 5, \S 30]{Gil} and \cite[Cahen]{Cah90}.
Regarding the valuative dimension and \thref{corthValDim}, a very elegant \cof treatment is given in the integral case by T. Coquand in~\cite{coqval}. 

The result of Exercise~\ref{exoInclusionBordLionel} is due to Lionel Ducos. \Pbmz~\ref{exoAnneauNoetherienReduit} is directly related to the article \cite[Coquand\&al.]{cls}.
 \Pbmz~\ref{exoChainesPotPrem} is drawn from the articles \cite[Brenner]{Brenner}, \cite[Coquand\&Lombardi]{cl1} and \cite{lom}.
 A variant for \trdis is found in~\cite{cl}. 

\newpage \thispagestyle{CMcadreseul}
\incrementeexosetprob


\chapter{The number of \gtrs of a module}
\label{chapNbGtrs}
\perso{compil\'e le \today}

\vskip-1em

\minitoc

\subsection*{Introduction}
\addcontentsline{toc}{section}{Introduction}

In this chapter we establish the \elrz, non-\noee and \cov version of some \gui{grand} \thos of commutative \algz.

These \thosz, due in their original version to \KRNz, Bass, Serre, Forster and Swan, regard the number of radical \gtrs of a \itfz, the number of \gtrs of a module, the possibility of producing a free submodule as a direct summand in a module, and the possibility of simplifying \isosz, in the following style: if $M\oplus N\simeq M'\oplus N$, then $M\simeq M'$.

Decisive progress was made by Heitmann \cite[(1984)]{Hei84} who proved how to get rid of \noees hypotheses.

Further progress was made by T.\ Coquand who proved in several articles how to obtain all the classical results (sometimes in a stronger form) by means of \dems that are both \covs and \elrsz.

The \dems given here are essentially those of \cite[Coquand]{Coq3,coq07} and of \cite[Coquand\&al.]{clq,clq2}.

\section[Kronecker's \tho and Bass' stable range]{\KROz's \tho and Bass' stable range (non-\noees versions of Heitmann)}
\label{secKroBass}

\subsec{\KROz's \thoz}
\label{subsecKro}

\KRNz's \thoz\footnote{This theorem of Kronecker is different from the one given in Chapter~\ref{chapGenerique}.} is usually stated in the following form (\cite{Kronecker}): an \agq \vrt in $\CC^n$ can always be defined by $n+1$ equations.

For Kronecker it was more about replacing a \sys of arbitrary equations in $\QQXn$ with an \gui{\eqvz} \sys having at most~$n+1$ equations. The \eqvc of two \syss as seen by Kronecker is translated in the current language by the fact that the two \ids have the same nilradical.
It is by using the \nst that we obtain the above formulation in the language of \gui{\agq \vrtsz.}

In the version proven in this section, we give the formulation \`a la Kronecker by replacing the \ri $\QQXn$
with an arbitrary \ri of \ddkz~$\leq n$.

The following lemma, although terribly trivial, is an essential key.
\begin{lemma}
\label{gcd}
For $u,v\in\gA$ we have
$$\begin{array}{c}
\DA(u,v) =   \DA(u+v,uv)   =   \DA(u+v)\vu\DA(uv)
\end{array}.$$
In particular, if $uv\in\DA(0)$, 
then $\DA(u,v)=\DA(u+v)$.
\end{lemma}
\begin{proof}
We obviously have $\gen{u+v,uv}\subseteq\gen{u,v}$,
therefore $\DA(u+v,uv) \subseteq\DA(u,v)$. Moreover,
$u^2=(u+v)u-uv\in\gen{u+v,uv}$, so $u\in\DA(u+v,uv)$.
\end{proof}

Recall that two sequences that satisfy the inequalities (\iref{eqCG}) in Proposition~\ref{corKrull} are said to be \copsz.

\begin{lemma}
\label{lemKroH}
Let $\ell\geq 1$.
If $(b_1,\ldots ,b_\ell)$ and  $(x_1\ldots ,x_\ell)$ are two \cop sequences in $\gA$
then for every $a\in\gA$ we have
$$\DA(a,b_1,\dots,b_\ell) = \DA(b_1+ax_1,\dots,b_\ell+ax_\ell),$$
\cad $a\in \DA(b_1+ax_1,\dots,b_\ell+ax_\ell)$.
\end{lemma}
\begin{proof}
We have the in\egts
$$\arraycolsep2pt
\begin{array}{rcl}
\DA(b_1x_1)& =  &\DA(0)    \\
\DA(b_2x_2)& \leq  & \DA(b_1,x_1)  \\
\vdots~~~~& \vdots  &~~~~  \vdots \\
\DA(b_\ell x_\ell )& \leq  & \DA(b_{\ell -1},x_{\ell -1})  \\
\DA(1)& =  &  \DA(b_\ell,x_\ell ).
\end{array}
$$
We deduce these
$$\arraycolsep2pt
\begin{array}{rcl}
\DA(ax_1b_1)& =  &\DA(0)    \\
\DA(ax_2b_2)& \leq  & \DA(ax_1,b_1)  \\
\vdots~~~~& \vdots  &~~~~  \vdots \\
\DA(ax_\ell b_\ell )& \leq  & \DA(ax_{\ell -1},b_{\ell -1})  \\
\DA(a)& \leq   &  \DA(ax_\ell,b_\ell ).
\end{array}
$$
We therefore have by Lemma~\ref{gcd}
$$\arraycolsep2pt
\begin{array}{rcl}
\DA(a)& \leq   & \DA(ax_\ell+b_\ell)\vu \DA(ax_\ell b_\ell)\\
\DA(ax_\ell b_\ell)& \leq &\DA(ax_{\ell-1}+b_{\ell-1})\vu \DA(ax_{\ell-1} b_{\ell-
1}) \\
\vdots~~~~& \vdots  &~~~~  \vdots  \\
\DA(ax_3b_3)& \leq  & \DA(ax_2+b_2)\vu \DA(ax_2 b_2)\\
\DA(ax_2b_2)& \leq  & \DA(ax_1+b_1)\vu \DA(ax_1 b_1) = \DA(ax_1+b_1).
\end{array}
$$
Therefore finally
$$\arraycolsep2pt
\begin{array}{rcl}
\DA(a)&\leq &\DA(ax_1+b_1)\vu \DA(ax_2+b_2)\vu \cdots  \vu \DA(ax_\ell+b_\ell)\\
&= &\DA(ax_1+b_1,ax_2+b_2,\ldots ,ax_\ell+b_\ell).
\end{array}
$$

\vspace{-1.5em}
\end{proof}

\vspace{-.7em}
\pagebreak	        

\begin{theorem}
\label{thKroH} \emph{(Non-\noe \KRNz-Heitmann \thoz,
with the Krull dimension)}
\begin{enumerate}
\item Let $n\geq 0$.
If $\Kdim \gA <n$ and $b_1$, \dots, $b_n\in\gA$,  there exist $x_1$, \dots, $x_n$ such that for every $a\in\gA$, $\DA(a,b_1,\dots,b_n) = \DA(b_1+ax_1,\dots,b_n+ax_n)$.
\item Consequently, in a \ri \ddi$n$, every \itf has the same nilradical as an \id generated by at most $n+1$ \eltsz.
\end{enumerate}
\end{theorem}

\begin{proof}
\emph{1.} Clear by Lemma~\ref{lemKroH} and Proposition~\ref{corKrull}
(if $n=0$, the \ri is trivial and $\DA(a)=\DA(\emptyset)$).\\
\emph{2.} Stems from \emph{1} because it suffices to iterate the procedure. Actually, if $\Kdim \gA\leq n$  and~$\fa=\DA(b_1,\ldots,b_{n+r})$  ($r\geq 2$), we finally obtain

\snic{\fa=\DA(b_1+c_1,\ldots,b_{n+1}+c_{n+1})}

with $c_i\in\gen{b_{n+2},\ldots ,b_{n+r} }$.
\end{proof}

\subsec{Bass' \texorpdfstring{\gui{stable range}}{"stable range"} \,\thoz, 1}
\label{subsecBass}

\begin{theorem}
\label{Bass0} \emph{(Bass' \thoz, with the Krull dimension, without \noetz)} 
Let $n\geq 0$. If $\Kdim \gA <n$, then $\Bdim\gA<n$. \\
Abbreviated to: $\Bdim\gA\leq \Kdim \gA$. In particular, if $\Kdim\gA<n$, every stably free \Amo of rank $\geq n$ is free (see \thref{corBass2}).
\end{theorem}

\begin{proof} Recall that $\Bdim\gA<n$ means that for all $b_1$, \dots, $b_n\in\gA$,
there exist some~$x_i$'s such that the following implication is satisfied

\snic{\forall
a\in\gA\quad (
1 \in\gen{a,b_1,\dots,b_n } \,\Rightarrow\, 1
\in\gen{b_1+ax_1,\dots,b_n+ax_n}).}

This results directly from the first item in \thref{thKroH}.
\end{proof}

\vspace{.7em}

%
%

\subsec{The local \KRO \thoz}
\label{Kroloc}

\begin{propdef}\label{propdefdisjointes}
In a \ri we consider two sequences $(a_0, \dots, a_d)$ and $(x_0, \dots, x_d)$ such that
$$\arraycolsep2pt\left\{\begin{array}{rcl}
a_0 x_0 & \in & \rD(0) \\
a_1 x_1 & \in & \rD(a_0,x_0) \\
a_2 x_2 & \in & \rD(a_1,x_1) \\
a_3 x_3 & \in & \rD(a_2,x_2) \\
& \vdots & \\
a_d x_d & \in & \rD(a_{d-1},x_{d-1}) 
\end{array}\right.
$$
We will say that these two sequences are \ixc{disjoint}{sequences}.
Then for $0 \leq i < d$, we have
\index{sequences!disjoint ---}

\snic{
\rD(a_0, \dots, a_i, x_0, \dots, x_i, a_{i+1} x_{i+1})
= \rD(a_0 + x_0, \dots, a_i + x_i).
}
\end{propdef}
\begin{proof}
%
%
An inclusion is obvious. To prove the converse inclusion, we use the \egts $\rD(a_i,x_i) = \rD(a_i x_i,a_i+x_i)$.

We then successively get
{
\vspace{1mm}\mathrigid2mu
$$ \hspace{0mm}
{\arraycolsep1pt 
\begin{array}{rcccccl}
a_0 x_0 & \in & \rD(0) & = & \rD( ) & =  & \rD()\\
& & & & & & ~\VRT{~\supseteq~~} \\
a_0, x_0, a_1 x_1 & \in & \rD(a_0,x_0) & = & \rD(a_0 x_0, a_0 + x_0)
         & =  & \rD(a_0 + x_0)\\
& & & & & & ~~~~~~\VRT{~\supseteq~~} \\
a_1, x_1, a_2 x_2 & \in & \rD(a_1,x_1) & = & \rD(a_1 x_1, a_1+x_1)
        & \subseteq & \rD(a_0 + x_0, a_1 + x_1)\\
& & & & & & ~~~~~~~~~~~\VRT{~\supseteq~~} \\
a_2, x_2, a_3 x_3 & \in & \rD(a_2,x_2) & = & \rD(a_2 x_2, a_2+x_2)
        & \subseteq & \rD(a_0+x_0, a_1 + x_1, a_2 + x_2)\\
\vdots~~~~~~& \vdots & \vdots& \vdots &\vdots & \vdots & ~~~~~~~~~~~\vdots \\
a_i, x_i, a_{i+1} x_{i+1} & \in
& \rD(a_i,x_i) & = & \rD(a_i x_i, a_i+x_i)
       & \subseteq & \rD(a_0 + x_0, \dots, a_i + x_i).
\end{array}
}
$$
}
\end{proof}

Note that the sequences $(a_0, \dots, a_d)$ and $(x_0, \dots, x_d)$
are \cops \ssi they are disjoint and $1\in\gen{a_d,x_d}$.

\begin{theorem}\label{thKroLoc}
Let $\gA$ be a \dcd \alo \ddi$d$, with Jacobson radical $\fm$.
We suppose that~$\fm$ is \emph{radically \tfz}, \cad there exist $z_1$, \dots, $z_n\in\gA$ such %
that $\fm=\DA(z_1, \dots, z_n)$.
Then $\fm$ is radically generated by $d$~\eltsz.%
\index{radically \tfz!ideal}\index{ideal!radically \tf ---}
\end{theorem}
\begin{proof}
Since $\Kdim \gA\leq d$ and $\fm$ is radically \tfz, \KRNz's \thoz~\ref{thKroH} tells us that~$\fm= \rD(x_0,\dots,x_d)$.
In addition, there exists a %
\cop sequence $(\ua) = (a_0,\dots,a_d)$ of $(\ux) = (x_0,\dots,x_d)$.
In particular (disjoint sequences), for every $i \leq d$, we have

\snic{
\rD(a_0, \dots, a_{i-1}, x_0, \dots, x_{i-1}, a_i x_i)
= \rD(a_0 + x_0, \dots, a_{i-1} + x_{i-1})
,}

but also (\cop sequences)
$
1 \in \gen{a_d, x_d}
$.
This shows that $a_d$ is \iv since $x_d \in \fm$. Let $i$ be the smallest index such that $a_i$ is \iv
(here we use the hypothesis that $\fm$ is detachable).\\
We then get $a_0$, \dots, $a_{i-1} \in \fm$, then

\snic{
\rD(x_0, \dots, x_{i-1}, x_i)
\subseteq \rD(a_0 + x_0, \dots, a_{i-1} + x_{i-1})
\subseteq \fm,}

and finally

\snic{
\begin{array}{c}
 \fm = \rD(x_0, \dots, x_{i-1}, x_i, x_{i+1},\dots, x_d)
\subseteq \hspace*{5cm}    \\[1mm]
\hspace*{4cm}
\rD(\underbrace{a_0 + x_0, \dots, a_{i-1} + x_{i-1}, x_{i+1},
\dots, x_d}_{d \hbox{\scriptsize ~\eltsz}})
\subseteq \fm.
\end{array}
}
\end{proof}

\rem For a \gnn see Exercises~\ref{exoKroLocvar} and \ref{exoKroLocvarbis}.
\eoe

\section{Heitmann dimension and Bass' \thoz}
\label{subsecDimHeit}

We will introduce a new dimension, which we will call the Heitmann dimension of a commutative \riz. Its \dfn will be copied from the inductive \dfn of the Krull dimension, and we will denote it by $\Hdim$.
Beforehand, we introduce the dimension $\Jdim$ defined by Heitmann. 

\begin{definota}
\label{notaJA}\label{defHeit}~
\begin{enumerate}

\item [--] If $\fa$ is an \id of $\gA$ we let $\JA(\fa)$ be its \ixc{Jacobson radical}{of an ideal}, \cad the inverse image of $\Rad(\gA\sur\fa)$ by the canonical projection $\gA\to\gA\sur\fa$.%
\index{radical!Jacobson ---}\index{Jacobson!radical}

\item [--] If $\fa=\gen{\xn}$ we denote $\JA(\fa)$ by $\JA(\xn)$. In particular,~$\JA(0)=\Rad \gA.$

\item [--] Let $\HeA$ be the set of \ids  $\JA(\xn)$.
We call it the \ixx{Heitmann}{lattice} of the \riz~$\gA$.

\item [--] We define $\Jdim\gA$ as equal to $\Kdim(\Heit\gA)$.
\end{enumerate}
\end{definota}

We therefore have $x\in\JA(\fa)$ \ssi for every $y\in\gA$, $1+xy$ is \iv
modulo~$\fa$. In other words

\snic{x\in\JA(\fa)\iff 1+x\gA\subseteq \sat{(1+\fa)},}

and $\JA(\fa)$ is the largest \id $\fb$ such that  $1+\fb\subseteq \sat{(1+\fa)}$. \\
We therefore have $\sat{\big(1+\JA(\fa)\big)}=\sat{(1+\fa)}$ and $\JA\big(\JA(\fa)\big)=\JA(\fa)$. \\
In particular $\JA\big(\JA(0)\big)=\JA(0)$ and the \ri $\gA/\Rad\gA$ has its Jacobson radical reduced to $0$.

\begin{lemma}
\label{gcd2} ~
\begin{enumerate}
\item For an arbitrary \id $\fa$ we have  $\JA(\fa)=\JA\big(\DA(\fa)\big)=\JA\big(\JA(\fa)\big)$.\\
Consequently,  $\HeA$ is a quotient \trdi of $\ZarA$.
\item For $u$, $v\in\gA$ we have

\snic{\JA(u,v)\; =  \; \JA(u+v,uv)  \;=\;  \JA(u+v)\vu\JA(uv).}

In particular, if $uv\in\JA(0)$, then $\JA(u,v)=\JA(u+v)$.
\end{enumerate}
\end{lemma}
\begin{proof}
We have $\fa\subseteq\DA(\fa)\subseteq\JA(\fa)$, therefore $\JA(\fa)=\JA\big(\DA(\fa)\big)=\JA\big(\JA(\fa)\big)$. 
\\
The \egt
$\DA(u,v)= \DA(u+v,uv)$ therefore implies  $\JA(u,v) =\JA(u+v,uv)$.
\end{proof}

\rdb
\comm \label{NOTAJdim}
The $\Jdim$ introduced by Heitmann in \cite{Hei84} corresponds to the following spectral space $\Jspec \gA$: it is the smallest spectral subspace of~$\Spec \gA$ containing the set $\Max \gA$ of the \idemas of $\gA$.
This space can be described as the adherence of $\Max \gA$ in $\Spec \gA$ for the constructible topology. A topology having as its \sgr of open sets the $\fD_\gA(a)$'s and their complements $\fV_\gA(a)$.
Heitmann notices that the dimension used in Swan's \tho or in \SSO theorem, namely the dimension of $\Max \gA$, only works well in the case where this space is \noez. In addition, in this case, the dimension of~$\Max \gA$  is that of a spectral space, the space $\jspec \gA$ formed by the \ideps which are intersections of \idemasz. However, in the \gnl case, the subspace $\jspec \gA$ of $\Spec\gA$ is no longer spectral, and so, according to a remark which he qualifies as philosophical,
 $\jspec \gA$ must be replaced by the spectral space that naturally offers itself as a spare solution, namely~$\Jspec \gA$. Actually, 
$\Jspec \gA$ is identified with the spectrum of the \trdiz~$\HeA$
(see \cite[\Tho 2.11]{clq2}),
and the \oqcs of~$\Jspec \gA$ form a quotient lattice of~$\ZarA$,
canonically \isoc to~$\HeA$.
In \comaz, we therefore define $\Jdim\gA$ as equal \hbox{to $\Kdim(\Heit\gA)$}. 
\eoe

\medskip The \dfn of the Heitmann dimension given below is quite natural, insofar as it mimics the \cov \dfn of the Krull dimension by replacing $\DA$ with $\JA$. 

\begin{definition}
\label{defHei2} Let  $\gA$  be a commutative \riz, $x\in\gA$ and $\fa$ be a \tf \idz. 
The \emph{Heitmann boundary of $\fa$ in $\gA$} is the quotient \ri $\gA_\rH^\fa:=\gA\sur{\JH_\gA(\fa)}$ with
         $$\preskip.0em \postskip.4em\JH_\gA(\fa):=\fa+(\JA(0):\fa).$$
This \id is called the \emph{Heitmann boundary \id of $\fa$ in $\gA$}.
\\
We also let  
$\JH_\gA(x):=\JH_\gA(x\gA)$ and $\gA_\rH^x:=\gA\sur{\JH_\gA(x)}$.%
\index{Heitmann boundary!quotient \riz, \idz}\index{ideal!Heitmann boundary ---}
\end{definition}

\begin{definition}
\label{defDHA} The \ix{Heitmann dimension} of $\gA$ is defined by \recu as follows
\begin{enumerate}
\item [--] $\Hdim \gA=-1$ \ssi $1_\gA=0_\gA$.
\item [--] Let $\ell\geq 0$,  we have the \eqvc

\snic{\Hdim \gA\leq \ell$  $\;\Longleftrightarrow\;$ for every $x\in\gA$,
$\Hdim(\gA_\rH^x)\leq \ell-1.}

\end{enumerate}
\end{definition}

This dimension  is less than or equal to the $\Jdim$ defined by Heitmann in~\cite{Hei84}, \cad the \ddk of the \trdi $\Heit(\gA)$.

\begin{fact}
\label{factKdimHdim} ~
\begin{enumerate}
\item The Heitmann dimension can only decrease by passage to a quotient \riz.
\item The Heitmann dimension is always less than or equal to the Krull dimension.
\item More precisely $\Hdim \gA\leq \Kdim\big(\gA/\JA(0)\big)\leq \Kdim\gA.$
\item Finally, $\Hdim \gA\leq0$ \ssi $\Kdim\big(\gA/\JA(0)\big)\leq 0$ (\cad  $\gA$ is \plcz).
\end{enumerate}
\end{fact}
\begin{proof}
\emph{1.} By \recu on $\Hdim\gA$(%
\footnote{Actually by \recu on $n$ such that $\Hdim\gA\leq n$.}) by noticing that for every $x\in\gA$, the \ri $(\gA/\fa)_\rH^x$ is a quotient of~$\gA_\rH^x$.
 
\emph{2.} By \recu on $\Kdim\gA$ (by using \emph{1}) by noticing that $\gA_\rH^x$ is a quotient of~$\gA_\rK^x$.
  
\emph{3} and \emph{4.} Let $\gB=\gA/\JA(0)$. Then $\rJ_\gB(0)=\gen{0}$, and we have $\gA_\rH^x\simeq \gB_\rH^{\ov x}=\gB_\rK^{\ov x}$ for all~$x\in\gA$.
\end{proof}

\subsec{Bass' \texorpdfstring{\gui{stable range}}{"stable range"} \,\thoz, 2}

\begin{theorem}
\label{Bass} \emph{(Bass' \thoz, with the Heitmann dimension, without \noetz)} 
Let $n\geq 0$. If $\Hdim \gA <n$, then $\Bdim\gA<n$. 
\\
Abbreviated to: $\Bdim\gA\leq \Hdim \gA$. In particular if $\Hdim\gA<n$, every stably free \Amo of rank $\geq n$ is free.
\end{theorem}

\begin{proof}
The same \dem would give \thref{Bass0} by replacing the Heitmann boundary with the Krull boundary.
Recall that $\Bdim\gA<n$ means that for all $b_1$, \dots, $b_n\in\gA$, there exist some~$x_i$'s such that the following implication is satisfied:

\snic{\forall
a\in\gA\quad (
1 \in\gen{a,b_1,\dots,b_n } \,\Rightarrow\, 1
\in\gen{b_1+ax_1,\dots,b_n+ax_n}).}
 
Recall that $1\in\gen{L}$ is equivalent to $1\in\JA(L)$ for every list $L$.
We reason by \recu on $n$.\\
When $n=0$ the \ri is trivial and $\JA(1)=\JA(\emptyset)$.\\
Suppose $n\geq 1$. Let $\fj=\JH_\gA(b_n)$.
The \hdr gives $x_1$, \ldots, $x_{n-1}\in\gA$ such that
\begin{equation}\preskip.3em \postskip.5em
\label{eqBass1}
1\in \gen{b_1+x_1a,\ldots , b_{n-1}+ x_{n-1} a} \quad \mathrm{in} \quad
\gA/\fj .
\end{equation}
Let $L=(b_1+x_1a,\ldots ,b_{n-1}+ x_{n-1} a)$. An arbitrary \elt of $\fj$ is written in the form $b_{n}y+x$ with $xb_n\in\JA(0)$. The membership (\ref{eqBass1}) therefore means that there exists an $x_n$ such that
$$
\preskip.2em \postskip.0em
x_nb_n\in\JA(0)\quad  \mathrm{and}\quad
1\in \gen{L,b_n,x_n}.
$$
A fortiori
\begin{equation}\preskip.0em 
\label{eqBass3}
 1\in \JA(L,b_n,x_n)=\JA(L,b_n)\vu\JA(x_n).
\end{equation}
Note that by hypothesis $1\in\gen{a,b_1,\dots,b_n}=\gen{L,b_n,a}$.
Therefore
\begin{equation}
\label{eqBass4}
1\in\JA(L,b_n,a)=\JA(L,b_n)\vu\JA(a).
\end{equation}
As $\JA(x_{n}a)=\JA(a)\vi\JA(x_{n})$,  (\ref{eqBass3}) and (\ref{eqBass4}) give by distributivity
$$
1\in \JA(L,b_n)\vu\JA(x_{n}a)=\JA(L,b_n,x_na)
.$$
Since $b_n\,x_n\,a\in\JA(0)$, Lemma~\ref{gcd2} gives $\JA(b_n,x_na)=\JA(b_n+x_na)$, and so 
$$\preskip.0em \postskip.2em
1\in\JA(L,b_n+x_na)=\JA(L,b_n,x_na),
$$
as required.
\end{proof}

\vspace{-.7em}
\pagebreak	

\subsec{\texorpdfstring{\gui{Improved}}{"Improved"}\, variant of Kronecker's \thoz}

\begin{lemma}\label{HenriLemma} 
Let $a$, $c_1$, \ldots, $c_m$, $u$, $v$, $w \in \gA$ and $Z=(c_1,\ldots,c_m)$.
\begin{enumerate}
\item
 $a \in \DA(Z)  \iff  1 \in \gen {Z}_{\gA[a^{-1}]}$.
\item 
$\big(w\in \Rad(\gA[a^{-1}])$ and $a\in \DA(Z,w)\big)$
$\Longrightarrow$ $a\in \DA(Z)$.
\item 
$\big(uv\in \Rad(\gA[a^{-1}])$ and $a\in \DA(Z,u,v)\big)$ $\Longrightarrow$ $a\in \DA(Z,u+v).$ 
\end{enumerate}
\end{lemma}

\begin{proof}
\emph{1.} Immediate.

\emph{2.} Suppose $a\in \DA(Z,w)$ and work in the \ri $\gA[a^{-1}]$. \\
We have $1 \in \gen {Z}_{\gA[a^{-1}]} +
\gen{w}_{\gA[a^{-1}]}$, and as $w$ is in $\Rad(\gA[a^{-1}])$, this implies that $1 \in \gen {Z}_{\gA[a^{-1}]}$, \cad $a\in \DA(Z)$.

\emph{3.}
Results from item \emph{2} because $\DA(Z,u,v)=\DA(Z,u+v,uv)$ (Lemma~\ref{gcd}).
\end{proof}

\rem We can ask ourselves if the \id $\Rad\gA[a^{-1}]$ is the best possible. The answer is yes. The implication of item \emph {2} is satisfied (for every $Z$) by replacing $\Rad\gA[a^{-1}]$ with $\fJ$ \ssi $\fJ\subseteq \Rad\gA[a^{-1}]$. \eoe

\begin{lemma}
\label{thCor2.2Heit}\relax ~\\
Suppose that $\Hdim\gA [1/a]< n$, $L\in \Ae n$ and $\DA(b)\leq\DA(a)\leq \DA(b, L)$.
Then there exists an~$X\in \Ae n$ such that $\DA(L+bX)=\DA(b, L)$, which is equivalent to~$b \in\DA(L+bX)$, or to $\,a \in\DA(L+bX)$.
In addition, we can take $X=aY$ with~$Y\in \Ae n$.
\end{lemma}
\begin{proof}
\emph{Preliminary remark.}
If $\DA(L+bX)=\DA(b, L)$, we have 
$a\in\DA(L+bX)$ because $a\in \DA(b,L)$. 
Conversely, if $a\in\DA(L+bX)$, we have $b \in
\DA(L+bX)$ (since $b\in\DA(a)$),
so $\DA(L+bX)=\DA(b, L)$.

We reason by \recu on $n$. The case $n=0$ is trivial.  
\\
Let $L=(b_1,\ldots,b_n)$ and we start by looking for $X\in \Ae n$.  \\
Let $\fj=\JH_{\gA[1/a]}(b_n)$ and $\gA'=\gA/(\fj\cap\gA)$, where $\fj\,\cap\,\gA$ stands for \gui{the inverse image of $\fj$ in $\gA$.}
We have an identification $\gA[1/a]/\fj=\gA'[1/a]$.
\\
As $\Hdim\gA'[1/a]< n-1$, we can apply the \hdr to $\gA'$ and $(a,b,b_1,\ldots ,b_{n-1})$, by noticing that $b_n = 0$ in $\gA'$. 
We then obtain $x_1$, \dots, $x_{n-1}\in \gA$ such that, by letting~$Z=(b_1+bx_1,\ldots, b_{n-1}+ bx_{n-1})$, we have $\rD(Z)=\rD(b, b_1,\dots,b_{n-1})$ in $\gA'$. By the preliminary remark, this last \egt is equivalent to $a \in \rD_{\gA'}(Z)$, which, by Lemma~\ref{HenriLemma}~\emph{1}, means that $1 \in \gen{Z}$ in $\gA'[1/a]$,
\cad $1 \in \gen{Z} + \fj$ in $\gA[1/a]$. By \dfn of the Heitmann boundary, this means that there exists an $x_n$, which we can choose in $\gA$, such that $x_nb_n\in\Rad\gA[1/a]$ and $1 \in \gen{Z, b_n,x_n}_{\gA[1/a]}$. \\
 We therefore have $a\in\DA(Z,b_n,x_n)$. But we also have $a\in\DA(Z,b_n,b)$, since
$$\preskip-.2em \postskip.4em
\gen {Z, b_n, b} = \gen {b_1,, \ldots, b_{n-1}, b_n, b} 
\eqdefi {\gen{L, b}},
$$
and since $a \in \DA(L,b)$ by hypothesis.
Recap: $a \in\DA(Z,b_n,x_n)$, $a \in\DA(Z,b_n,b)$ so $a \in\DA(Z,b_n,bx_n)$.
The application of Lemma~\ref{HenriLemma}~\emph{3} with $u = b_n$, $v = bx_n$ provides $a \in\DA(Z,b_n+bx_n)$, \cad $a \in \DA(L + bX)$ where $X = (x_1, \ldots, x_n)$.
\\
Finally, if $b^p\in\gen{a}_\gA$, we can apply the result with $b^{p+1}$ instead of $b$ since $\DA(b)=\DA(b^{p+1})$. Then $L+b^{p+1}X$ is re-expressed as $L+baY$.
\end{proof}

For $a\in\gA$, we always have $\Hdim\gA[1/a]\leq \Kdim\gA[1/a]\leq \Kdim\gA$. Consequently the following \tho improves \KRNz's \thoz.

\begin{theorem}
\label{KroH2} \emph{(\KRNz's \thoz, Heitmann dimension)}
\begin{enumerate}
\item Let $n\geq 0$. If  $a$, $b_1$, $\dots$, $b_n\in\gA$ and $\Hdim\gA[a^{-1}] <n$,
then there exist~$x_1$,~$\dots$,~$x_n\in \gA$ such that
$$\preskip.4em \postskip.4em 
\DA(a,b_1,\dots,b_n) =\DA(b_1+ax_1,\dots,b_n+ax_n). 
$$
\item Consequently, if $a_1$, \ldots, $a_r$, $b_1$, \dots, $b_n\in\gA$ and $\Hdim\gA[1/a_i] <n$ for $i\in\lrbr$, then there exist $y_1$, \dots, $y_n\in \gen{a_1,\ldots,a_r}$ such that
$$\preskip.4em \postskip.4em 
\DA(a_1,\ldots ,a_r,b_1,\dots,b_n) =\DA(b_1+y_1,\dots,b_n+y_n). 
$$
\end{enumerate}
\end{theorem}
\begin{proof}
\emph{1.} Direct consequence of Lemma \ref{thCor2.2Heit} by making $a=b$. 

 \emph{2.} Deduced from \emph{1} by \recu on $r$: 

\snac{
\begin{array}{cclcl} 
  \fa&=&\DA(a_1,\ldots,a_r,b_1,\dots,b_n)& =  & \DA(a_1,\ldots ,a_{r-1},b_1,\dots,b_n) \vu
\DA(a_r)  \\[1mm] 
  & =  &  \fb\vu\DA(a_r), \quad \hbox{with}&& \\[1mm] 
\fb  & =  & \DA(b_1+z_1,\dots,b_n+z_n) 
\end{array}
}

where $z_1$, \dots, $z_n\in \gen{a_1,\ldots ,a_{r-1}}$, so $\fa=\DA(a_r,b_1+z_1,\dots,b_n+z_n)$, and we once again apply the result.
\end{proof}

\section[\SSO and Forster-Swan \thosz]{Serre's Splitting Off theorem, the Forster-Swan theorem, and Bass' cancellation theorem}
\label{secSOSFSa}

In this section, we describe the matrix \prts of a \ri that allow us to make \SSO theorem  and the Forster-Swan \tho 
(control of the number of \gtrs of a \mtf in terms of the number of local \gtrsz) work.

The following sections consist in developing results that show that certain \ris satisfy the matrix \prts in question.
The first \ris that appeared (thanks to Serre and Forster) were the \noe \ris with certain dimension \prts (the \ddk for Forster and the dimension of the maximal spectrum for Serre and Swan).
Later Heitmann showed how to get rid of the \noet regarding the \ddkz, and gave the guiding ideas to do the same for the dimension of the maximal spectrum. 
In addition Bass also introduced a \gnn in which he would replace the \ddk by the maximum of the Krull dimensions for the \ris associated with a partition of the Zariski spectrum in constructible subsets.
Finally, Coquand brought a \gui{definitive} light to these questions by generalizing the results and by treating them \cot thanks to two subjacent notions to the previous \demsz:  $n$-stability on the one hand and Heitmann dimension on the other.
The purely matrix aspect of the \pbs to be solved has clearly been highlighted in a review paper by Eisenbud-Evans \cite[(1973)]{Eisenbud}. The present section can be considered as a non-Noetherian and constructive approach to these works.

\begin{definition}\label{defiSdimGdim}
Let $\gA$ be a \ri and $n\geq 0$ be an integer.
\begin{enumerate}
\item We write $\Sdim\gA< n$ if, for every matrix $F$ of rank $\geq n$, there is a \umd \coli of the columns.\\
In other words $1=\cD_n(F)\Rightarrow \exists X,\, 1=\cD_1(FX)$.
\item We write $\Gdim\gA< n$
when the following \prt is satisfied. For every matrix $F=[\,C_0\,|\,C_1\,|\,\dots\,|\,C_p\,]$  (the~$C_i$'s are the columns, and let $G=[\,C_1\,|\,\dots\,|\,C_p\,]$)
such that $1=\cD_1(C_0)+\cD_n(G)$, there is a \coli  $C_0+\sum_{i=1}^{p}\alpha_i C_i$, which is \umdz. 
\end{enumerate}
\end{definition}

In the acronym $\Sdim$, $\mathsf{S}$ refers to \gui{splitting} or to \gui{Serre} and is justified by \thref{thSerre}.
Similarly, in $\Gdim$, $\mathsf{G}$ refers to \gui{generators} and is justified by \thref{thSwan}.

The notations $\Sdim\gA<n$ and  $\Gdim\gA<n$ are justified by the
following obvious implications, for every $n\geq 0$, 

\snic{\Sdim\gA<n\Rightarrow\Sdim\gA<n+1$ and $\Gdim\gA<n\Rightarrow\Gdim\gA<n+1.}

Note that $\cD_n(F) \subseteq \cD_1(C_0) + \cD_n(G)$, and consequently the hypothesis for~$F$ in $\Sdim\gA< n$ implies the hypothesis for~$F$ in $\Gdim\gA< n$.
Moreover the conclusion in $\Gdim\gA< n$ is stronger.
This gives the following item~\emph{2}.

\begin{fact}\label{factSDimGdimHdim}~
\begin{enumerate}
\item  $\Sdim\gA< 0\iff\Gdim\gA< 0\iff$ the \ri $\gA$ is trivial.
\item For all $n\geq 0$, we have $\;\Gdim\gA< n\,\Longrightarrow\,\Sdim\gA< n$.
\\
Abbreviated to $\Sdim\gA\leq \Gdim\gA$.
\item If $\gB=\gA/\fa$, we have 
$\,\Sdim\gB\leq \Sdim\gA$ and $\,\Gdim\gB\leq \Gdim\gA$.
\item We have  $\Sdim\gA= \Sdim\gA/\!\Rad\gA$ and  $\Gdim\gA= \Gdim\gA/\!\Rad\gA$.
\item If $\gA$ is $n$-stable (Section~\ref{secSUPPORTS}), then $\Gdim\gA<n$ (\thref{matrix}).\\
Abbreviated to $\Gdim\gA\leq \Cdim\gA$. 
\item If $\Hdim\gA<n$, then $\Gdim\gA<n$ (\thref{MAINCOR}).\\
Abbreviated to $\Gdim\gA\leq \Hdim\gA$. 
\end{enumerate}
 
\end{fact}
%
\begin{proof}
It suffices to prove items \emph{3} and \emph{4.}
Item \emph{4} is clear because an \elt of $\gA$ is \iv in $\gA$ \ssi it is \iv in $\gA/(\Rad\gA)$.

\emph{3 for $\Sdim$.} Let $F\in\Ae{m\times r}$ with $\cD_n(F)=1$ modulo $\fa$. If $n>\inf(m,r)$ we obtain $1\in\fa$ and all is well.
Otherwise, let $a\in\fa$ such that $1-a\in \cD_n(F)$. \\
Consider the matrix $H\in\Ae{(m+n)\times r}$ obtained by superposing $F$ and the matrix $a\In$ followed by $r-n$ null columns. \\
We have $1-a^{n}\in \cD_n(F)$, so $1\in\cD_n(H)$. A \coli of the columns of $H$ is \umdz. The same \coli of the columns of $F$ is \umd modulo $\fa$. 

\emph{3 for $\Gdim$.} The same technique works, but here it suffices to consider the matrix $H\in\Ae{(m+1)\times r}$ obtained by inserting the row $[\,a\;0\;\cdots\;0\,]$ underneath~$F$.
\end{proof}

The \dem of the following fact helps to justify the slightly surprising \dfn chosen for $\Gdim\gA<n$.
 
\begin{fact}\label{factGdimBdim}
For all $n\geq 0$, we have $\Gdim\gA<n\Rightarrow\Bdim\gA < n$. \\
Abbreviated to $\Bdim\gA\leq \Gdim\gA$.
\end{fact}
%
\begin{proof}
For example with $n=3$.
Consider $(a,b_1,b_2,b_3)$ with $1=\gen{a,b_1,b_2,b_3}$.
We want some $x_i$'s such that $1=\gen{b_1 + ax_1 ,  b_2 + ax_2,  b_3 + ax_3}$.
Consider the matrix 
$F=\cmatrix{b_1 & a & 0 & 0\cr
b_2 & 0 & a & 0\cr
b_3&  0&  0 & a
}= [\,C_0\mid G\,]$ with $G=a\I_3$.
We have  

\snic{1=\cD_1(C_0)+\cD_3(G)$,\quad  \cad $1=\gen{b_1,b_2,b_3}+\gen{a^{3}},}

because $1=\gen{b_1,b_2,b_3}+\gen{a}$.
By applying the \dfn of $\Gdim\gA<3$ to $F$, we obtain a \vmd $\tra[\,b_1 + ax_1\; b_2 + ax_2\; b_3 + ax_3\,]$.
\end{proof}
%

\subsec{\SSO \thoz}

The following version of Serre's \tho is relatively easy, the delicate part being to establish that $\Sdim  \gA<k$ for a \ri $\gA$. Modulo  \thrfs{matrix}{MAINCOR} 
we obtain the truly strong versions of the \thoz.
\begin{theorem}
\label{thSerre} \emph{(\SSO \thoz, with $\Sdim$)}\\
Let $k\geq 1$ and $M$ be a \pro \Amo of rank $\geq k$, or more \gnlt \isoc to the image of a matrix of rank $\geq k$.\\
Suppose that $\Sdim  \gA< k$.
Then $M\simeq N\oplus \gA$ for a certain module~$N$ \isoc to the image of a matrix of rank $\geq k-1$. 
\end{theorem}

\begin{proof}
Let $F\in\Ae{n\times m}$ be a matrix with $\cD_{k}(F)=1$.
By \dfnz, we have a vector $u=\tra{[\,u_1\,\cdots \,u_n\,]}\in\Im F$ which is \umd in~$\Ae n$. Therefore $\gA u$ is a free submodule of rank $1$ and a direct summand in~$\Ae n$, and a fortiori in~$M$.
More \prmtz, if $P\in\GAn(\gA)$ is a \prr of image~$\gA u$, we obtain $M=\gA u\oplus N$ with 

\snic{N=\Ker(P)\cap M=(\In-P)(M)=\Im\big((\In-P)\,F\big).}

It remains to see that $(\In-P)\,F$ is of rank $\geq k-1$.
Even if it entails localizing and making a change of basis, we can suppose that $P$ is the standard projection $\I_{1,n}$.
Then the matrix $G=(\In-P)\,F$ is the matrix $F$ in which we have replaced its first row by $0$, and it is clear that $\cD_k(F)\subseteq \cD_{k-1}(G)$.

\end{proof}

Thus,  if $M$ is the image of $F\in\Ae {n\times m}$ of rank $\geq k$,
we obtain a \dcn $M= N\oplus L$ where $L$ is free of rank $1$ as a direct summand in~$\Ae n$ and $N$ \isoc to the image of a matrix of rank $\geq k-1$. Now if~$F$ is of greater rank, we can iterate the procedure and we have the following corollary (with the correspondence $h \leftrightarrow k-1$).
	
\begin{corollary}\label{corthSerre}
Let $\gA$ be a \ri such that $\Sdim \gA\le h$, and~$M$ be a module \isoc to the image of a matrix of rank $\geq h+s$.  Then~$M$ contains as a direct summand a free submodule of rank~$s$. More \prmtz, if~$M$ is the image \hbox{of $F\in\Ae {n\times m}$} of rank $\geq h+s$, we have $M= N\oplus L$ where $L$ is free of rank $s$ and a direct summand in $\Ae n$, and $N$ is the image of a matrix of \hbox{rank $\geq h$}.
\end{corollary}

\subsec{The Forster-Swan \thoz} \label{subsecFoSw}

Recall that a \mtf $M$ is said to be \lot generated by~$r$ \elts if $\cF_r(M)=1$.
On this subject see the local number of \gtrs lemma (Lemma~\ref{lemnbgtrlo}).

The Forster-Swan \tho below was first established for the \ddk ($\Kdim$ instead of $\Gdim$). 
The version presented here is relatively easy, and the delicate part is to establish that $\Gdim  \gA\leq \Kdim  \gA$ for every \riz~ $\gA$. 
Modulo \thrfs{matrix}{MAINCOR} we obtain the known better versions of the \thoz, under an entirely \cov form. 

\begin{theorem}
\label{thSwan} \emph{(Forster-Swan \thoz, with $\Gdim$)} Let $k\geq 0$ and $r\geq 1$.
\\
If $\Gdim  \gA\leq   k$, or even only if $\Sdim  \gA\leq   k$ and $\Bdim  \gA\leq   k+r$, and if a \tf \Amoz~$M$ is   \lot generated by $r$ \eltsz, then it is generated by $k+r$ \eltsz.\\ 
In the first case, more \prmtz, if $M$ is generated by $y_1$, \dots, $y_{k+r+s}$, we can compute $z_1$, \ldots, $z_{k+r}$ in $\gen{y_{k+r+1},\ldots,y_{k+r+s}}$ such that $M$ is generated by $y_1+z_1$, $\ldots$, $y_{k+r}+z_{k+r}$.
\end{theorem}
\begin{proof}
Since $M$ is \tf and $\cF_r(M)=1$, $M$ is the quotient of a \mpf $M'$ satisfying $\cF_r(M')=1$. We can therefore suppose that~$M$ is \pfz.
\\
Starting from a \sgr with more than $k+r$ \eltsz, we are going to replace it with a \sgr of the stated form minus an \eltz.
Therefore let $(y_0,y_1,\dots,y_p)$ be a \sgr of $M$ \hbox{with $p\geq k+r$}, and $F$ be a \mpn of $M$ for this \sysz.
Then by hypothesis $1 = \cF_{r}(M)= \cD_{p+1-r}(F)$, and since $p+1-r\geq k+1$ we have $1 = \cD_{k+1}(F)$.

\emph{First case.}  
Let $L_0$, \ldots, $L_p$ be the rows of $F$. We apply the \dfn of $\Gdim\gA<k+1$ with  
the transposed matrix of $F$ (which is of \hbox{rank $\geq k+1$}).
We obtain some $t_i$'s such that the row $L_0+t_1L_1+\cdots +t_pL_p$ is \umdz.
Replacing the row~$L_0$ with the row $L_0+t_1L_1+\cdots +t_pL_p$ amounts to the same as replacing the \sgr $(y_0,y_1,\dots,y_p)$ with  

\snic{(y_0,y_1-t_1y_0,\dots,y_p-
t_py_0)=(y_0,y'_1,\dots,y'_p).}

Since the new row $L_0$ is \umdz, a suitable \coli of the columns is of the form $\tra{[\,1\;y_1\,\cdots \,y_p\,]}$. This means that we have $y_0+y_1y'_1+\cdots +y_py'_p=0$ in $M$, 
and thus that $(y'_1,\dots,y'_p)$ generates~$M$.

\emph{Second case.} We apply the \dfn of $\Sdim\gA<k+1$ with the matrix $F$. We obtain  a \umd \coli of columns, and we add this column in the first position in front of~$F$. 
Then, by applying Fact~\ref{corBass} with $\Bdim\gA<k+r+1\leq p+1$, by \mlrrsz, we obtain a new \mpn of~$M$ (for another \sgrz) with the first column equal to $\tra[\,1\,0\,\cdots\,0\,]$. This means that the first \elt of the new \sgr is null.
\end{proof}

\Thref{thSwan} is obviously valid by replacing the \ri $\gA$ with the \ri $\gA/\Ann(M)$ or  $\gA/\cF_0(M)$.
We propose in \thref{thSwan2} a slightly more subtle refinement.

\hum{J'ai pr\'ef\'er\'e mettre $1\in\cD_1(F)$ plutt que d'affirmer les choses \`a partir of $k=0$. J'ai introduit $q$.}

\begin{proposition}
\label{basic2} Let $F=[\,C_0\,|\,C_1\,|\,\dots\,|\,C_p\,]\in\Ae{n\times (p+1)}$ (the~$C_i$'s are the columns) and $G=[\,C_1\,|\,\dots\,|\,C_p\,]$, such that~$F=[\,C_0\,\vert\,G\,]$.
\\
If $1=\cD_1(F)$ and if we have $\Gdim(\gA/\cD_{k+1}(F))<k$ for $k\in\lrbq$, then there exist $t_1$, \dots, $t_p$ such that the vector $C_0+t_1C_1+\cdots+t_pC_p$ is \umd modulo $\cD_{q+1}(F)$.
\end{proposition}
%
\begin{proof}
First consider the \ri $\gA_2=\gA/\cD_2(F)$. Since $\cD_1(F)=1$  and  $\Gdim(\gA_2)<1$,  we obtain \hbox{some $t_{1,i}$'s} and $C_{1,0}=C_0+t_{1,1}C_1+\cdots+t_{1,p}C_p$ such that $\cD_1(C_{1,0})=1$ modulo~$\cD_2(F)$, \cad $\cD_1(C_{1,0}) +\cD_2(G)=1$. We change $F$ into $F_1$ by replacing $C_0$ with $C_{1,0}$ without changing $G$. Note that we have $\cD_i(F_1)=\cD_i(F)$ for every $i$.
\\
We then consider the \ri $\gA_3=\gA/\cD_3(F_1)$ with $\Gdim(\gA_3)<2$. 
\\
Since $\cD_1(C_{1,0}) +\cD_2(G)=1$, we obtain $C_{2,0}=C_{1,0}+t_{2,1}C_1+\cdots+t_{2,p}C_p$ such that $\cD_1(C_{2,0})=1$ modulo $\cD_3(F)$, \cad $\cD_1(C_{2,0}) +\cD_3(G)=1$. We change $F_1$ into $F_2$ by replacing $C_{1,0}$ with $C_{2,0}$ without changing $G$. We once again have $\cD_i(F_2)=\cD_i(F)$ for every $i$.
\\
We continue as above until we obtain a vector $C_{q,0}$ of the form $C_0+t_1C_1+\cdots+t_pC_p$ \umd modulo  $\cD_{q+1}(F)$.
\end{proof}

\vspace{-.7em}
\pagebreak

\begin{theorem}
\label{thSwan2} \emph{(Forster-Swan \thoz, more \gnlz, with $\Gdim$)}\\
Let  $M$  be a \tf module over $\gA$. Let $\ff_k=\cF_k(M)$ be its \idfsz.
Suppose that $1\in\ff_m$ (\cad $M$ is \lot generated by $m$ \eltsz) and that for $k\in\lrb{0..m-1}$,  we have $\Gdim(\gA/\ff_k)<  m-k$. 
Then $M$ is generated by~$m$ \eltsz. 
More \prmtz, if $M=\gen {y_1,\dots,y_{m+s}}$, we can compute some~$z_i$'s \hbox{in $\gen {y_{m+1},\dots,y_{m+s}}$} such that $M=\gen {y_1+z_1,\dots,y_{m}+z_{m}}$.
\end{theorem}
\begin{proof}
Since $\ff_0$ annihilates $M$, we can replace $\gA$ with $\gA/\ff_0$, or, which amounts to the same thing, suppose that $\ff_0=\cF_0(M)=0$, which we do from now on.\\
Starting with a \sgr of $M$ with more than $m$ \elts  
         we are going to replace it by
a \sgr of the stated form minus an \eltz.
Therefore let $(y_0,y_1,\dots,y_p)$ be a \sgr of~$M$ \hbox{with $p\geq m$}.

When the module is \pf we reason as for \thref{thSwan}.
\\
Let $F$ be a \mpn of $M$ for the considered \sgrz.
We have $\ff_{k+1}=\cD_{p-k}(F)$, and in particular $1\in\ff_p=\cD_1(F)$. The hypotheses of Proposition~\ref{basic2} are then satisfied with $q=p$ for the transposed matrix of $F$.
If $L_0$, \ldots, $L_p$ are the rows of $F$, we obtain some $t_i$'s with  $L_0+t_1L_1+\cdots +t_pL_p$ \umd modulo $\cD_{p+1}(F)=\ff_0=0$. 
The remainder of the argument is as in \thref{thSwan}.

The reasoning in the case where $M$ is only assumed to be \tf consists in showing that $M$ is the quotient of a \pf module which has a \mpn supporting with success the \dem of Proposition~\ref{basic2}.
Let $\uy=[\,y_0\;\cdots \;y_p\,]$.   
Every syzygy between the $y_i$'s is of the form $\uy\,C=0$ for some $C\in\gA^{p+1}$.
\\
The \idf $\ff_{p+1-i}$ of~$M$ is the \id $\Delta_i$, the sum of the \idds $\cD_i(F)$, for \hbox{$F\in \gA^{(p+1)\times n}$} that satisfy $\uy\,F=0$, \cad for the matrices that are \gui{syzygy matrices for $(y_0, \ldots, y_p)$.}\\
By the hypotheses, we have $\Delta_1=1$ and  $\Gdim(\gA/\Delta_{k+1})<k$ \hbox{for $k\in\lrbp$}.
\\
The fact that $\Delta_1=1$ is observed over a syzygy matrix $F_1$.  
\\
Consider the matrix $\tra F_1$  and the \riz~$\gA_2=\gA/\Delta_2$.  
\hbox{As $\Gdim(\gA_2)<1$}, we obtain a \coli $C_{1,0}$ of the columns of $\tra F_1$  \umd modulo~$\Delta_2$, 
\cad such that  $1= \cD_1(C_{1,0})+ \Delta_2$.
More \prmtz, we obtain $C_{1,0}=\tra F_1 X_1$ with $X_{1}=\tra[\,1\;x_{1,1}\;\cdots\;x_{1,p}\,]$.
\\
The \egt $1= \cD_1(C_{1,0})+ \Delta_2$ provides an \elt $a\in\Delta_2$ obtained as a \coli of a finite number of minors of order $2$ of syzygy matrices, 
and so $a\in\cD_2(F_2)$ for a syzygy matrix  $F_2$. 
Then consider the  matrix $F'_2=[\,F_1\,|\,F_2\,]$. 
For the transposed matrix of $F'_2$ we first obtain that the column $C_2=\tra {F'_2} \,X_1$ is \umdz. 
We replace the first column of $\tra {F'_2}$ by $C_2$, which gives a matrix $\tra {F''_2}$ suitable for the hypotheses of $\Gdim\gA_3<2$ (where $\gA_3=\gA/\Delta_3$), \cad $1=\cD_1(C_{2})+\cD_2(F''_2)$. We ultimately obtain 
a \coli $C_{2,0}$ of the columns of $\tra {F'_2}$  \umd  modulo~$\Delta_3$, \cad such that  $1= \cD_1(C_{2,0})+ \Delta_3$.
More \prmtz, $C_{2,0}=\tra {F'_2}\, X_2$ with $X_{2}=\tra[\,1\;x_{2,1}\;\cdots\;x_{2,p}\,]$.
And so forth.\\
We ultimately obtain a syzygy matrix  for $\uy$, 
$$\preskip.4em \postskip.4em
F=\lst{F_1\,|\,\cdots\,|\,F_p}
$$
and a vector $X_{p}=\tra[\,1\;x_{p,1}\;\cdots\;x_{p,p}\,]$ with the \coli $\tra F\, X_p$ \umd (since it is \umd modulo $\Delta_{p+1}=\ff_0=0$).\\
The remainder of the argument is as in \thref{thSwan}.
\end{proof}

\comm \Thref{thSwan2} with $\Hdim$ or $\Kdim$ instead of $\Gdim$ has as an easy consequence in \clama some much more abstract statements, which seem much more scholarly. For example the usual statement of the Forster-Swan \thoz\footnote{Corollary 2.14 (page 108) in \cite{Kun} or \Tho 5.8 (page 36) in \cite{Mat}.
In addition, the authors replace $\gA$ with $\gA\sur{\Ann(M)}$, which costs nothing.}
 (stated in the case where $\Max\gA$ is \noez)
uses the maximum,\footnote{Recall that $\jspec\gA$ designates the subspace of $\SpecA$ formed by the \ideps which are intersections of \idemasz.} 
for $\fp\in \jspec\gA$ of $\mu_\fp(M)+\Kdim(\gA/\fp)$:  here $\mu_\fp(M)$ is the minimum number of \gtrs of $M_\fp$.
This type of statement suggests that the \ideps that are intersections of \idemas play an essential role in the \thoz. In reality, it is not \ncr to scare children with $\jspec\gA$,
because this abstract \tho is exactly \eqv (in the considered case, and in \clamaz) to \thref{thSwan2} for the~$\Jdim$, which in this envisaged case is equal to the~$\Hdim$.
In addition, from a strictly  practical point of view it is unclear how to access  the quite mysterious maximum of the $\mu_\fp(M)+\Kdim(\gA/\fp)$.
By contrast, the hypotheses of \thref{thSwan2} are susceptible to a constructive proof, which in this case will lead to an algorithm making it possible to explicate the conclusion.

\subsec{Bass' cancellation \thoz}
\label{secSFS}

\begin{definition}\label{defiSimplifiMod}
Given two modules $M$ and $L$ we say that $M$ \emph{is cancellative for $L$} if $M\oplus L\simeq N\oplus L$ implies $M \simeq N$.%
\index{cancellative!module}%
\index{module!cancellative ---}
\end{definition}

\begin{lemma}\label{lemSimplifiMod}
Let $M$ and $L$ be two \Amosz.
In the following statements we have
  1 $\Leftrightarrow$ 2 and 3 $\Rightarrow$ 2.
\begin{enumerate}
\item $M$ is cancellative for $L$.
\item For every decomposition $M\oplus L=M'\oplus L'$ with $L'\simeq L$, there exists an \auto $\sigma$ of $M\oplus L$ such that $\sigma(L')=L$.
\item For every decomposition $M\oplus L=M'\oplus L'$ with $L'\simeq L$, there exists an \auto $\theta$ of $M\oplus L$  such that $\theta(L')\subseteq M$.
\end{enumerate}
\end{lemma}
\begin{proof} The \eqvc of \emph{1} and \emph{2} 
is a game of photocopies.
\\
 \emph{1 $\Rightarrow$ 2.} Suppose $M\oplus L=M'\oplus L'$.
Since $L\simarrow L'$, we obtain an \iso $M\oplus L\simarrow M'\oplus L$,
so $M\simarrow M'$, and by performing the sum we obtain an \iso $M\oplus L\simarrow M'\oplus L'$,
\cad an \auto of $M\oplus L$ which sends $L$ to~$L'$.
\\
 \emph{2 $\Rightarrow$ 1.} Suppose $N\oplus L\simarrow M\oplus L$.
This \iso sends $N$ to~$M'$ and~$L$ to~$L'$, such that $M\oplus L=M'\oplus L'$.
Therefore there is an \autoz~$\sigma$ of $M\oplus L$ which sends $L$ to $L'$, and say $M$ to $M_1$. Then,

\snic{N\simeq M'\simeq (M'\oplus L')/L'= (M\oplus L)/L' = (M_1\oplus L')/L'
\simeq M_1\simeq M.}

 \emph{3 $\Rightarrow$ 2.} 
 Since $\theta(L')$ is a direct summand in $M\oplus L$, it is a direct summand in $M$, which we write as $M_1\oplus \theta(L')$.
Let $\lambda$ be the \auto of $M\oplus L$ which swaps $L$ and $\theta(L')$ by fixing $M_1$.
Then $\sigma=\lambda\circ \theta$ sends $L'$ to $L$.
\end{proof}

Recall that an \elt $x$ of an \emph{arbitrary} module $M$ is said to be \umd when there exists a \lin form $\lambda\in M\sta$ such that $\lambda(x)=1$. It amounts to the same as saying that $\gA x$ is free (of basis $x$) and a direct summand in $M$ (Proposition~\ref{propSplittingOffAlgExt}).

\begin{theorem}
\label{thBassCancel2} \emph{(Bass' cancellation \thoz, with $\Gdim$)}\\
Let $M$ be a \ptf \Amo of rank $\geq k$.
If $\Gdim \gA< k$,
then $M$ is cancellative for every \ptf \Amoz:
 if $Q$ is \ptf and $M\oplus Q\simeq N\oplus Q$, then~$M\simeq N$.
\perso{on ne pourra sans doute pas affaiblir l'hypoth\`ese
au cas where $M$ is \isoc to l'image of a matrice $F$ of rank $\geq k$}
\end{theorem}

\begin{proof}
Suppose that we have shown that $M$ is cancellative for $\gA$. \\
Then,
since $M\oplus\Ae\ell$ also satisfies the hypothesis, we show by \recu on $\ell$ that $M$ is cancellative for $\Ae{\ell+1}$.
As a result $M$ is cancellative for every direct summand in $\Ae{\ell+1}$.\\
Finally,  $M$ is cancellative for $\gA$ because it satisfies item \emph{3} of Lemma~\ref{lemSimplifiMod} for~$L=\gA$.
Indeed, suppose that $M=\Im F\subseteq\Ae n$, where $F$ is a \mprn (of rank $\geq k$),
and let $L'$ be a direct summand \hbox{in $M\oplus \gA$}, \isoc to $\gA$:
$L'=\gA(x,a)$ with $(x,a)$ \umd in  $M\oplus \gA$.
Since every \lin form over $M$ extends to $\gA^{n}$, there exists a form $\nu\in (\Ae n)\sta$ such that $1\in\gen{\nu(x),a}$.
By Lemma~\ref{forbass} below, \hbox{with $x=C_0$}, there exists \hbox{a $y\in M$} such that $x'=x+a y$
is \umd in $M$. Consider a \hbox{form $\mu\in M\sta$} such that $\mu(x')=1$. We then define an \auto $\theta$ \hbox{of $M\oplus \gA$} as follows
$$\theta=\cmatrix{1&0\cr -a\mu&1}\cmatrix{1&y\cr0&1}\;
\cad\; \cmatrix{m\cr b}\mapsto \cmatrix{m+by\cr \mu(x)b-a\mu(m)}.$$
Then
$\theta(x,a)=(x',0)$, so $\theta(L')\subseteq M$.
\Trf by Lemma~\ref{lemSimplifiMod}.
\end{proof}

In the following lemma, which ends the \dem of \thref{thBassCancel2}, 
         we use the notations
of Proposition~\ref{basic2}, the matrix $F=[\,C_0\,|\,C_1\,|\,\dots\,|\,C_p\,]$ being that of the previous \thoz.
\begin{lemma}\label{forbass}
If  $\Gdim \gA <k$ and $\cD_k(F) = 1 =\DA(C_0 )\vu \DA(a)$, then there \hbox{exist $t_1$, \dots, $t_p$} such that 

\snic{1 = \DA(C_0 +at_1C_1+\cdots+at_pC_p).}
\end{lemma}

\begin{proof}
Consider the matrix $[\,C_0\,|\,aC_1\,|\,\dots\,|\,aC_p\,] $, obtained by replacing $G$ by $aG$ in $F$.
As $\DA(C_0)\vu \cD_k(G)=1=\DA(C_0)\vu\DA(a)$, we indeed have by \dit $\DA(C_0)\vu \cD_k(aG)=1$. 
\Trf since $\Gdim\gA<k$.
\end{proof}


\subsec{A simple \cara \prt for $\Gdim\gA<n$}

In order to prove $\Gdim\gA<n$ for a \ri $\gA$ it suffices to verify the conclusion (in the \dfn of $\Gdim\gA<n$) for particularly simple matrices.
This is the subject of the following proposition.

\begin{proposition}\label{propGdimGdim}
For a \ri $\gA$ we have $\,\Gdim\gA<n$ \ssi  
for every matrix $V\in\MM_{n+1}(\gA)$ of the form
$$
V=\cmatrix{
b&c_1&\cdots&\cdots&c_n\cr
b_1&a&0&\cdots&0\cr
\vdots&0&\ddots&\ddots&\vdots\cr
\vdots&\vdots&\ddots& \ddots&0\cr
b_n&0&\cdots&0&a\cr
} = [\,V_0\mid V_1 \mid \dots \mid V_n\,]
,
$$
and for every $d \in \gA$ such that $1 = \gen{b,a,d}$,
there exist  $x_i$'s $\in \gA$ such that
$$\preskip.3em \postskip.4em
1 = \cD_1(V_0+x_1V_1+\cdots+x_nV_n)+\gen{d}.
$$
\end{proposition}
\rem
Instead of using an \elt $d$ subjected to the constraint $1=\gen{a,b,d}$, we could have used a pair $(u,v)$ not subjected to any constraint and replaced $d$ \hbox{by $1+au+bv$} in the conclusion. In this form, it is particularly obvious that if the condition above is satisfied for the \ri $\gA$, it is satisfied for every quotient of $\gA$. \eoe
\begin{proof}
To show that the condition is \ncrz, we reason with the quotient \ri $\gB=\aqo\gA d$ and we consider the matrix 

\snic{F=V=[\,V_0\mid V_1 \mid \dots \mid V_n\,].}

With the notations of \Dfnz~\ref{defiSdimGdim} we have $p=n$, $F=[\,C_0\mid G\,]$,  \hbox{and $C_i=V_i$} for $i\in\lrbzn$. 
\\
Since $1 = \gen{b,a,d}$ in $\gA$, we have $1 = \gen{b,a^{n}}\subseteq \cD_1(C_0)+\cD_n(G)$ in~$\gB$, 
and the hypothesis of the \dfn is satisfied. Since $\Gdim\gB<n$, we obtain~$x_i$'s in $\gA$ such that 
$$ 
\preskip.0em \postskip.3em 
 1=\cD_1(C_0+x_1C_1+\cdots+x_nC_n)\hbox { in } \gB .
$$
Hence the desired conclusion in $\gA$.

To prove the converse we proceed in two steps.
First of all recall that if the condition is satisfied for the \ri $\gA$, it is satisfied for every quotient of $\gA$.
We will actually use this condition \hbox{with $d=0$} (the hypothesis over $V$ then becomes of the same type as that which serves to define~\hbox{$\Gdim<n$}), with the \ri $\gA$ and certain of its quotients.

\emph{First step: the case where the matrix $F$ has $n+1$ columns, \cad $p=n$.} 
With $F\in\Ae{m\times (n+1)}$, we have by hypothesis a \lin form $\varphi_0:\Ae m \to\gA$ and an  $n$-multi\lin alternating form $\psi : (\Ae m)^n \to\gA$ such that 
$$
\preskip.3em \postskip.3em 
1=\varphi_0(C_0)+\psi(C_1, \dots, C_n). 
$$
For $j\in\lrbn$ let $\varphi_j:\Ae m \to\gA$ be the \lin form 
$$\preskip.3em \postskip.3em 
 X\mapsto \psi(C_1, \dots,
C_{j-1},X,C_j, \dots, C_n).
$$
By letting $a = \psi(C_1, \ldots, C_n)$, we then have
\begin{itemize}
\item $\varphi_1(C_1)=\cdots=\varphi_n (C_n) = a$,
\item  $\varphi_i(C_j)=0$ if $1\leq i\neq j\leq n$
\end{itemize}
Considering the matrix of the $\varphi_i(C_j)$'s, we obtain
$$
V = [\,V_0\mid \dots \mid V_n\,] :=
\Cmatrix{2pt}{
\varphi_0(C_0)&\varphi_0(C_1)&\cdots&\cdots&\varphi_0(C_n)\\[.3em]
\varphi_1(C_0)&a&\phantom{m}0\phantom{m}&\cdots&0\\[.3em]
\vdots&0&\ddots& \ddots&\vdots\\[.1em]
\vdots&\vdots&\ddots& \ddots&0\\[.3em]
\varphi_n(C_0)&0&\cdots&\phantom{m}0\phantom{m}&a
}
,$$
that is $V = [\,\varphi(C_0)\mid \dots \mid \varphi(C_n)\,]$ by letting $\varphi(Z)=\Cmatrix{2pt}{\varphi_0(Z)\cr\vdots\cr\varphi_n(Z)}$.\\
We can apply the hypothesis with $d=0$. We find $x_1$, \dots, $x_n\in\gA$
such that the vector $V_0+x_1V_1+\cdots+x_nV_n$ is \umdz.
This vector is equal \hbox{to $\varphi(C_0+x_1C_1+\cdots+x_nC_n)=\varphi(C)$}.
Since this vector is \umd and since $\varphi$ is \linz, the vector $C$ is itself \umdz.

\emph{Second step: the \gnl case.}\\
As $1=\cD_1(C_0 )+\cD_n(G)$, we have a family $(\alpha_i)_{i\in \lrbq}$ of subsets with~$n$ \elts of $\lrbp$ such that $1 = \cD_1(C_0 )+ \som_i \cD_n(G_{\alpha_i})$, where $G_{\alpha_i}$ is the extracted matrix of $G$ by uniquely considering the columns whose index is in~$\alpha_i$.
Let $C_{0,0}=C_0$ and $J_\ell = \som _{i>\ell} \cD_n(G_{\alpha_i})$. We then apply the case of the first step successively with $\ell=1$, \ldots, $q$ to obtain 

\snic{1= \cD_1(C_{0,\ell})=\cD_1(C_{0,\ell-1}) + \cD_n(G_{\alpha_\ell})
\hbox{ in }\gA/J_\ell}

and therefore $\cD_1(C_{0,q}) = 1$ in~$\gA$.
\\
Note that in this second step, we use the result of the first step with quotient \ris of $\gA$.
\end{proof}
%

\section{Supports and \texorpdfstring{$n$}{n}-stability}
\label{secSUPPORTS}

In Section~\ref{secManipElemCol} we will establish \thos regarding the \mlrs on matrices. They will have as corollaries some grand \thos due to Serre, Forster, Bass and Swan. 
We will give them in two similar but nevertheless different versions.
We do not think that they can be reduced to a unique form.

The first version is based on the notion of $n$-stability.
This version leads inter alia to a sophisticated result due to Bass in which a partition of the Zariski spectrum intervenes in a finite number of subsets which are all of small dimension (smaller than the \ddk of the \riz). This result will be used in Chapter~\ref{ChapMPEtendus} to prove Bass' theorem (\thref{thBassValu}) regarding the extended modules.

The second version uses the Heitmann dimension, introduced in Section~\ref{subsecDimHeit}, less than or equal to the \ddkz, but for which we do not know of an analogue of Bass' sophisticated version.

Section~\ref{secSUPPORTS} gives a few \ncr preliminaries for the first version based on the $n$-stability. 

\subsec{Supports, dimension, stability}

\begin{definition}\label{defiSupportT}
A \emph{\sutz} over a \ri $\gA$ in a \trdiz~$\gT$ is a map $D:\gA\to\gT$, which satisfies the following axioms%
\index{support!over a commutative \riz}
\vspace{-2pt}
\begin{itemize}
\item [$\qquad\bullet$] $D(0_\gA)=0_\gT,\;\;\; D(1_\gA)=1_\gT$,
\item [$\qquad\bullet$]  $D(ab)=D(a)\vi D(b),$
\item [$\qquad\bullet$]  $D(a+b)\leq D(a)\vu D(b)$.
\end{itemize}
Let $D(\xn)=D(x_1)\vu\cdots\vu D(x_n)$.
\end{definition}

It is clear that $\DA:\gA\to\ZarA$ is a \sutz, called the \emph{Zariski \sutz}.
The following lemma shows that the Zariski \sut is the \gui{free} \sutz.%
\index{support!Zariski ---}

\begin{lemma}\label{lemAnodin}
For every \sut $D$ we have
\begin{enumerate}
\item \label{i1lemAnodin} $D(a^m)=D(a)$ for $m\geq1$, $D(ax)\leq D(x),$ $D(a,b)=D(a+b,ab).$
\item \label{i2lemAnodin} $\gen{\xn}= \gen{\yr} $ implies $D(\xn)=D(\yr).$
\item \label{i3lemAnodin} $\DA(y)\leq\DA(\xn)$ implies $D(y)\leq D(\xn).$
\item \label{i4lemAnodin} There exists a unique \homo $\theta$ of \trdis which makes the following diagram commute:

\vspace{-1.2em}
\Pnv{\gA}{\DA}{ D }{\ZarA}{\theta}{\gT}{ }{\suts}{\homos of \trdisz}
\end{enumerate}
\end{lemma}
\facile

Thus every \sut $D:\gA\to\gT$  such that $D(\gA)$ generates $\gT$ as a \trdi is obtained by composing the Zariski \sut with a passage to the quotient $\ZarA\to\ZarA\sur{\sim}$ by an \eqvc relation compatible with the lattice structure.

Denote $D(\fa) $ by $D(\xn)$ if $\fa=\gen{\xn}$.
We say that a vector $X\in\Ae n$ is \emph{$D$-\umdz} if $D(X)=1$.
\index{Dunimo@$D$-unimodular!vector}

\subsubsection*{Dimension of a \sutz, \KRNz's \thoz}

\begin{definition}\label{defi2SuitComp} 
Given two sequences $(\xzn)$ and~$(\bzn)$ in~$\gA$ and a \sutz~$D$ over $\gA$, we say that the two sequences are \emph{$D$-\copsz} if we have the following in\egts
\begin{equation}\label{eqC1G}
\left.\arraycolsep2pt
\begin{array}{rcl}
 D(b_0x_0)& =  & D(0)    \\
 D(b_1x_1)& \leq  &  D(b_0,x_0)  \\
\vdots ~~~~& \vdots  &~~~~  \vdots \\
 D(b_n x_n )& \leq  &  D(b_{n -1},x_{n -1})  \\
 D(1)& =  &   D(b_n,x_n )
\end{array}
\right\}
\end{equation}
The \sut $D$ is said to be \emph{of \ddk $\leq n$} if every sequence $(\xzn)$ in $\gA$
admits a $D$-\cop sequence. Let $\Kdim(D)\leq n$.%
\index{Dcomple@$D$-complementary sequences}%
\index{Krull dimension!of a support}%
\index{complementary sequences!for a support}
 \end{definition}

For example for $n=2$ the \cop sequences correspond to the following picture in $\gT$.
$$\preskip-.6em \postskip.1em
\SCO{D(x_0)}{D(x_1)}{D(x_2)}{D(b_0)}{D(b_1)}{D(b_2)}
$$

\rem Note that $\Kdim\gA=\Kdim(\DA)$.
\eoe

\smallskip 
The \dem of the following lemma can be copied from that of Lemma~\ref{lemKroH} by replacing $\DA$ by $D$. \KROz's \tho is then a direct consequence.

\begin{lemma}
\label{lem2KroH}
Let $\ell\geq 1$.
If $(b_1,\ldots ,b_\ell)$ and  $(\xl)$ are two $D$-\cop sequences in $\gA$,
then for every $a\in\gA$ we have
$$\preskip.3em \postskip.2em
D(a,b_1,\dots,b_\ell) =  D(b_1+ax_1,\dots,b_\ell+ax_\ell),
$$
\cad $D(a)\leq  D(b_1+ax_1,\dots,b_\ell+ax_\ell)$.
\end{lemma}

\vspace{-.5em}
\pagebreak

\begin{theorem}
\label{th2KroH} \emph{(\KRNz's \thoz, for the \sutsz)}\\
If $D$  is a \sut of \ddk $\leq n$, for every \itfz~$\fa$
there exists an \id $\fb$ generated by $n+1$ \elts such that $D(\fa)=D(\fb)$.
Actually, for all $b_1$, \ldots, $b_{n+r}$ ($r\geq 2$), there exist $c_j\in\gen{b_{n+2},\ldots ,b_{n+r} }$
such that $D(b_1+c_1,\ldots,b_{n+1}+c_{n+1})=D(b_1,\ldots,b_{n+r})$.
\end{theorem}

\subsubsection*{Faithful supports}

In this subsection we prove in particular that the \ddk of a \ri (which we already know is equal to the dimension of its Zariski support) is equal to that of its Zariski lattice: here 
 we keep the promise 
made in \ref{corfactDDKTRDI}.

\begin{definition}\label{defisutfidele}
A \sut $D:\gA\to\gT$ is said to be \emph{faithful} if $\gT$ is generated by the image of $D$ and if, for every $a\in\gA$ and $L\in \Ae m$, the in\egt $D(a)\leq D(L)$ implies the existence of a $b\in\gen{L}$ such that $D(a)\leq  D(b)$.
\index{support!faithful ---}
\index{faithful!support}
\end{definition}

For example the Zariski \sut $\DA$ is always faithful.

\medskip Let $D:\gA\to\gT$ be a support. If the image of $\gA$ generates $\gT$, since we have the \egt $D(a_1)\vi\cdots \vi D(a_n)=D(a_1\cdots a_n)$, every \elt of $\gT$ can be written in the form $D(L)$ for a list $L$ of \elts of $\gA$.

\begin{lemma}\label{lemsutfidele}
If $D$ is faithful and $\Kdim\gT<k$ then $\Kdim(D)<k$.
In particular the \ddk of a \ri is equal to that of its Zariski lattice.
\end{lemma}
\begin{proof}
Let $(a_1,\ldots,a_k)$ be a sequence in $\gA$. We must show that it admits a $D$-\cop sequence. \\
Since $\Kdim\gT<k$, the sequence $\big(D(a_1),\ldots,D(a_k)\big)$ has a \cop sequence $\big(D(L_1),\ldots,D(L_k)\big)$ in $\gT$ with lists in $\gA$ for $L_i$,
$$\arraycolsep2pt
\begin{array}{rcl}
 D(a_1)\vi D(L_1)& =  & D(0)    \\
 D(a_2)\vi D(L_2)& \leq  &  D(a_1,L_1)  \\
\vdots~~~~~~~~& \vdots  &~~~~  \vdots \\
 D(a_k)\vi D(L_k )& \leq  &  D(a_{k -1},L_{k -1})  \\
 D(1)& =  &   D(a_k,L_k).
\end{array}
$$
Since $D$ is faithful, there exists a $c_k$ in $\gen{a_k,L_k}$
such that $D(1)\leq D(c_k)$, which gives $b_k\in\gen{L_k}$
such that $D(1)\leq D(a_k,b_k)$.\\
Note that we have
$$
D(a_kb_k)=D(a_k)\vi D(b_k )\leq D(a_k)\vi D(L_k )\leq D(a_{k -1},L_{k -1}).
$$
Since $D$ is faithful, we have $c_{k-1}\in\gen{a_{k-1},L_{k-1}}$ with $D(a_kb_k)\leq D(c_{k-1})$, which gives $b_{k-1}\in\gen{L_{k-1}}$ such that $D(a_kb_k)\leq D(a_{k-1},b_{k-1})$.\\
And so forth. Ultimately, we have constructed a sequence $(b_1,\ldots,b_k)$ which is $D$-\cop to $(a_1,\ldots,a_k)$.
\end{proof}

\vspace{-.7em}
\pagebreak

\subsubsection*{$n$-stable \sutsz}

We now abstract the \prt described in Lemma~\ref{lem2KroH} for the \cop sequences in the following form.


\begin{definition}\label{defisutnstab}~
\begin{enumerate}
\item Let $n\geq 1$. A \sut $D:\gA\to\gT$ is said to be \emph{$n$-stable} when, for 
all $a\in\gA$ and $L\in \Ae n$, there exists an $X\in \Ae n$
such that $D(L,a)=  D(L+aX)$, \cad $D(a)\leq  D(L+aX)$.%
\index{support!nstab@$n$-stable ---}%
\index{nstab@$n$-stable!support}
\item \label{defiAnneaunStable}
The \ri $\gA$ is said to be \emph{$n$-stable} if its Zariski \sut $\DA$ is $n$-stable.
We will write $\Cdim\gA<n$ to say that $\gA$ is $n$-stable.%
\index{ring!nstab@$n$-stable ---}%
\index{nstab@$n$-stable!\riz}
\item The \ri $\gA$ is said to be $0$-stable if it is trivial.
\end{enumerate}

\end{definition}

In the acronym $\Cdim$, $\mathsf{C}$ alludes to \gui{Coquand.}

Naturally, if  $\Kdim(D)< n$ then $D$ is $n$-stable.
In particular, with the free support $\DA$, we obtain $\Cdim\gA\leq \Kdim\gA$.
Moreover, Kronecker's \tho applies (almost by \dfnz) to every $n$-stable \sutz.

The notation $\Cdim\gA<n$ is justified by the fact that if $D$ is $n$-stable, it \hbox{is $(n+1)$-stable}.  
Finally, item \emph{3} in the \dfn was given for the sake of clarity, but it is not really \ncrz: by reading item~\emph{1} \hbox{for $n=0$}, we obtain that for every $a\in\gA$, $D(a)\leq D(0)$.

\medskip 
\exls ~\\ 
1) A \advz, or more \gnlt a \ri $\gV$ which satisfies \gui{$a\mid b$ or $b\mid a$ for all $a$, $b$,} is $1$-stable, even in infinite \ddkz.
For all $(a, b)$ it suffices to find some $x$ such that $\gen{a, b} = \gen{b + x a}$.
\hbox{If $a = q b$}, we have $\gen{a, b} = \gen{b}$ and we take $x=0$.
If $b = q a$, we have $\gen{a, b} = \gen{a}$ and we take $x=1-q$.

2) A Bézout domain is $2$-stable. More \gnltz, a strict Bézout \ri (see Section~\ref{secBézout} on \paref{secpfval} and Exercise~\ref{exoAnneauBézoutStrict})  is $2$-stable. 
More \prmtz, for $a$, $b_1$, $b_2 \in \gA$, there exist $x_1$, $x_2$ such \hbox{that $a\in \gen {b_1 + x_1a, b_2 + x_2a}$}, \cad $\gen {a, b_1, b_2} = \gen {b_1 + x_1a, b_2 + x_2a}$.
\\
Indeed, by question \emph {1.c} of the exercise, there exist \com $u_1$ and $u_2$  such that $u_1{b_1} + u_2{b_2} = 0$.
We take $x_1$, $x_2$ such that $u_1x_1 + u_2x_2 = 1$ and we obtain the \egt

$\qquad\qquad\quad 
a = u_1b_1 + a + u_2b_2 =  
u_1(b_1 + x_1a) + u_2(b_2 + x_2a)
$.
\eoe

\hum{apparemment pas besoin of la pr\'ecision with le localis\'e
vu que $\Cdim\gA[1/u]\leq \Cdim\gA$.}
\begin{fact}\label{factStBdim}
 We always have  $\,\Bdim\gA\leq \Cdim\gA$.
\end{fact}
%
\begin{proof}
If $\gA$ is $n$-stable, then $\Bdim \gA  < n$: indeed, we apply the \dfn with $(a,\an)$ in $\gA$ satisfying $1 \in \gen{a, \an}$.
\end{proof}
%

\begin{fact}\label{lemsutnstab}
If $D$ is $n$-stable, for every $a\in\gA$
 and $L\in \Ae n$, there exists an $X\in \Ae n$
such that $D(L,a)=  D(L+a^2X)$, \cad $D(a)\leq  D(L+a^2X)$.
\end{fact}

Indeed, $D(a)=D(a^2)$ and $D(L,a)=D(L,a^2)$.

\subsec{Constructions and patchings of \sutsz}

\begin{definition}\label{lemsutHeit} ~\\
The map $\JA:\gA\to\Heit\gA$ defines the \emph{Heitmann \sutz}.
\index{support!Heitmann ---}
\end{definition}

\rem
A priori $\Kdim\DA=\Kdim\gA\geq \Kdim\JA \geq \Jdim\gA$.
We lack examples that would show that the two in\egts can be strict.
\eoe

\begin{lemma}\label{lem2GaussJoyal} \emph{(Variant of the Gauss-Joyal lemma~\ref{lemGaussJoyal})\iJG}\\
If $D$ is a \sut over $\gA$, we obtain a \sut $D[X]$ over $\gA[X]$ by letting $$\preskip.2em \postskip.3em D[X](f)=D\big(\rc(f)\big).$$
\end{lemma}
\begin{proof}
Lemma~\ref{lemGaussJoyal}  gives $\DA\big(\rc(fg)\big)=\DA\big(\rc(f)\big)\vi \DA\big(\rc(g)\big)$.
\end{proof}
%
\begin{lemma}\label{lemsutquo}\emph{(\Sut and quotient)}
Let $D:\gA\to\gT$  be a \sut and $\fa$ be a \itf of $\gA$.
We obtain a \sut
$$\preskip.2em \postskip.4em
D/\fa:\gA\to\gT\sur{\fa}\eqdefi \gT\sur{(D(\fa)=0)}
$$
by composing $D$ with the projection $\Pi_{D(\fa)}:\gT\to\gT\sur{(D(\fa)=0)}$. 
\begin{enumerate}
\item $\DA/\fa$ is canonically \isoc to $\rD_{\gA\sur{\fa}}\circ \Zar (\pi_\fa)$,
where $\pi_\fa$ is the canonical map $\gA\to\gA\sur{\fa}$.
\item If $D$ is faithful, then so is $D/\fa$.
\item If $D$ is $n$-stable, then so is $D/\fa$.\\
In particular $\Cdim\gA/\fa\leq \Cdim\gA$.
\end{enumerate}
\end{lemma}
\begin{proof} 
Recall that $\Pi_{D(\fa)}(x)\leq\Pi_{D(\fa)}(y)\iff x\vu {D(\fa)}\leq y\vu {D(\fa)}$.
\\
\emph{1.} Results from Fact~\ref{fact2Zar}. \\
\emph{2.} Let $D'=D\sur{\fa}$. Let $a\in\gA$ and $L$ be a vector such that~\hbox{$D'(a)\leq D'(L)$}.
We seek some $b\in\gen{L}$  such that~$D'(a)\leq D'(b)$. By \dfn of~$D'$ we have~\hbox{$D(a)\leq D(L,\fa)$}, and since $D$ is faithful, there exists a $c\in\gen{L}+\fa$ such that~\hbox{$D(a)\leq D(c)$}, which gives some $b\in L$ such that $D(a)\leq D(b,\fa)$, in other words~$D'(a)\leq D'(b)$.\\
\emph{3.} Let $a \in \gA$ and $L \in \Ae n$. We seek $X \in \Ae n$ such that~\hbox{$D'(a) \leq  D'(L + aX)$}, \cad $D(a) \vu D(\fa)  \leq  D(L + aX) \vu D(\fa)$.
However, we have some $X$ which is suitable for $D$, that is $D(a) \leq  D(L + aX)$, therefore it is suitable for $D'$.
\end{proof}

Dually we have the following lemma.
\begin{lemma}\label{lemsutloc}\emph{(\Sut and \lonz)}
Let $D:\gA\to\gT$ be a \sut and $u$ be an \elt of $\gA$.
We obtain a \sut
$$\preskip.1em \postskip.4em 
D[1/u]:\gA\to\gT[1/u]\eqdefi\gT\sur{(D(u)=1)} 
$$
by composing $D$ with $j_{D(u)}:\gT\to\gT\sur{(D(u)=1)}$.
\begin{enumerate}
\item $\DA[1/u]$ is canonically \isoc to $\rD_{\gA[1/u]}\circ \Zar (\iota_u)$,
where $\iota_u$ is the canonical map $\gA\to\gA[1/u]$.
\item If $D$ is faithful, then so is $D[1/u]$.
\item If $D$ is $n$-stable, then so is $D[1/u]$. \\
In particular $\Cdim\gA[1/u]\leq \Cdim\gA$.
\end{enumerate}
\end{lemma}
\begin{proof}
Recall that $j_{D(u)}(x)\leq j_{D(u)}(y)\iff x\vi {D(u)}\leq y\vi {D(u)}.$
\\
\emph{1.} Results from Fact~\ref{fact2Zar}. \\
\emph{2.} Let $D'=D[1/u]$. Let $a\in\gA$ and $L$ be a vector such that $D'(a)\leq D'(L)$.
 By \dfn of $D'$ we have $D(au)=D(a)\vi D(u)\leq D(L)$. Since $D$ is faithful, there exists a $b\in\gen{L}$ such that $D(au)\leq D(b)$, \cadz~$D'(a)\leq D'(b)$.\\
\emph{3.} As for Lemma~\ref{lemsutquo} by replacing $D/\fa$ and $\vu$ by $D[1/u]$  and~$\vi$.
\end{proof}
%

\begin{lemma}\label{lemPartitionSpec}\label{corlemPartitionSpec}~
\begin{enumerate}
\item Let $D:\gA\to\gT$ be a  \sut and $b\in \gA$. 
\begin{enumerate}
\item $D/b$ and $D[1/b]$ are $n$-stable \ssi $D$ is $n$-stable.
\item If $D$ is faithful and if $\gT/b$ and $\gT[1/b]$ are of \ddk $< n$, then $D$ is $n$-stable.
\end{enumerate}
\item Let $\gA$ be a \ri and $b\in\gA$. Then $\aqo\gA b$ and $\gA[1/b]$ are $n$-stable \ssi $\gA$ is $n$-stable. \\
Abbreviated to: $\Cdim\gA=\sup\big(\Cdim\aqo\gA b,\Cdim\gA[1/b]\big)$.
%
%
\end{enumerate}

\end{lemma}
\begin{proof} It suffices to show the direct implication in item~\emph{1a.}\\
Let $a\in \gA$ and~$L\in\Ae n$. 
Since $D/b$ is $n$-stable, we have some $Y\in\Ae n$ such that~\hbox{$D(a)\leq  D(L+aY)$} in ${\gT/\big(D(b)=0\big)}$, \cad in~$\gT$,
$$
D(a)\leq D(b)\vu D(L+aY).\eqno(*)
$$
Next we apply the $n$-stability of $D[1/b]$ with $ab$ and $L+aY$ which provides some $Z\in\Ae n$ such that $D(ab)\leq  D(L+aY+abZ)$ in ${\gT/\big(D(b)=1\big)}$. 
\\
In~$\gT$, by letting $X = Y+bZ$, this is expressed as
$$D
(ab)\vi D(b) \leq D(L+aX), \quad \hbox {i.e.}\quad D(ab)\leq D(L+aX).
\eqno (\#)
$$
But we have $\gen{b,L+aX}=\gen{b,L+aY}$, therefore $D(b,L+aX)= D(b,L+aY)$.
The in\egts $(*)$ and $(\#)$ are then expressed as

\snic{D(a)\leq D(b)\vu D(L+aX) \et D(a)\vi D(b)\leq D(L+aX).}

This implies (by \gui{cut,} see \paref{coupure1}) that $D(a)\leq  D(L+aX)$.
\end{proof}
%

\subsubsection*{Constructible partitions of the Zariski spectrum}

A \ix{constructible} subset of  $\SpecA$ is a Boolean combination of open sets of basis $\fD(a)$. In \clamaz, if we equip the set~$\SpecA$ with the \gui{constructible topology} having as its basis of open sets the constructible subsets, we obtain a spectral space, the \ix{constructible spectrum} \emph{of the \riz~$\gA$}, which we can identify with $\Spec\Abul$. 

From a \cof point of view, we have seen that we can replace $\SpecA$ (an object a little too ideal) by the lattice~$\ZarA$ (a concrete object), \isoc in \clama to the lattice of \oqcs of $\Spec\gA$.
When we pass from the Zariski topology to the constructible topology in \clamaz, we pass from $\ZarA$ to $\Bo(\ZarA)\simeq\Zar(\Abul)$ in \coma (for this last \isoz, see \thref{thZedGenEtBoolGen}).

Hyman Bass took interest in the partitions of the constructible spectrum.
An \elr step of the construction of such a partition consists in the replacement of a \ri $\gB$ by the two \ris $\aqo{\gB}{b}$ and $\gB[1/b]$, for an \elt $b$ of $\gB$.
An important remark made by Bass is that these two \ris can each have a strictly smaller \ddk than that of~$\gB$,
whereas certain \prts of the \riz, to be satisfied in $\gB$, only need to be satisfied in each of its two children.  
This is the case for the $n$-stability of the free \sutz. In any case, this is the analysis that T.~Coquand made from a few pages of \cite{BASS}.
\perso{citer les pages en question}

In \clamaz, from any covering of $\Spec\gA$ by open sets of the constructible topology, we can extract a finite covering, which we can refine into a finite partition by some \oqcs (\cad some finite Boolean combinations of open sets with basis~$\fD(a)$).
These are a lot of high caliber abstractions, but the result is extremely concrete, and this is the result that interests us in practice. 

We define in \coma a \emph{constructible partition of the Zariski spectrum} by its dual version, which is a \sfio in the \agBz~\hbox{$\Zar\Abul=\Bo(\ZarA)$}.\\
In practice, an \elt of $\Zar\Abul$ is given by a double list in the \ri $\gA$

\snic{(a_1,\dots,a_\ell;u_1,\dots,u_m)=(I;U)}

that defines the following \elt of $\Zar\Abul$ 

\snic{\Vi_i{\lnot\DAbul(a_i)}\vi\Vi_j\DAbul(u_j)=\lnot\DAbul(a_1,\dots,a_\ell) \vi \DAbul(u), \hbox{ where }u=\prod_j u_j.}

To this \elt $(I;U)$ is associated the \ri $(\aqo\gA I) [1/u]$.\footnote{In \clama $\Spec(\aqo\gA I) [1/u]=\bigcap_{a\in I}\fV(a)\cap\bigcap_{v\in U}\fD(v)$, where~$\fV(a)$ designates the complement of $\fD(a)$.} A \sfio of $\Bo(\ZarA)$ can then be obtained as a result of a tree construction which starts with the double list $(0;1)$ and which authorizes the replacement of a list $(I;U)$ by two double lists $(I,a;U)$ and $(I;a,U)$ for some $a\in\gA$.

The following crucial \tho is a corollary of item \emph{2} of Lemma~\ref{corlemPartitionSpec}.
\begin{theorem}\label{thPartitionSpec}
Consider a constructible partition of~$\,\SpecA$,
described as above by a family $(I_k;U_k)_{k\in\lrbm}$. 
Let $\fa_k$ be the \idz~$\gen{I_k}$ and $u_k$ be the product of the \elts of $U_k$.
\begin{enumerate}
\item If $D:\gA\to\gT$ is a \sutz, and if all the $(D/\fa_k)[1/u_k]$'s are $n$-stable, then $D$ is $n$-stable.
\item In particular, if each \ri $\gA[1/{u_k}]/{\fa_k}$ is $n$-stable (for example if its \ddk is $<n$), then $\gA$ is $n$-stable.
\end{enumerate}
\end{theorem}
%

\rems ~

1) The paradigmatic case of  an $n$-stable \ri is given in the previous \tho when each \ri  $\gA[1/u_i]/{\fa_i}$ is of \ddkz~$<n$.

2) Every constructible partition of $\SpecA$ can be refined in the partition described by the $2^{n}$ \cop pairs formed from a finite list~\hbox{$(\an)$} in $\gA$.

3) Analogous tree constructions appear in Chapter~\ref{chap gen loc} in the framework of the basic \plgcz, but there are other \risz, localized \ris denoted by $\gA_{\cS(I;U)}$, that intervene then.
\eoe


\penalty-2500
\section{\Elr column operations}
\label{secManipElemCol}

\vspace{3pt}
In this section we establish  analogous \thos in two different contexts.
The first uses the stability of a support, the second uses the Heitmann dimension.

The reader can visualize most of the results of the chapter in the following picture, keeping in mind \Thos \ref{thSerre}, \ref{thSwan},  \ref{thSwan2} and~\ref{thBassCancel2}.

An arrow such that $\Sdim\lora\Gdim$ is added for $\Sdim\gA\leq \Gdim\gA$.

\newcommand\URRU[1]%
  {\ar@{}[urr]|{\rotatebox{20}{\hbox{$\lllra$}}}^{\rotatebox{20}{#1}}}

\newcommand\DRRD[1]%
  {\ar@{}[drr]|{\rotatebox{-25}{\hbox{$\lllra$}}}_{\rotatebox{-25}{#1}}}

\newcommand\DRRU[1]%
  {\ar@{}[drr]|{\rotatebox{-20}{\hbox{$\lllra$}}}^{\rotatebox{-20}{#1}}}

\newcommand\URRD[1]%
  {\ar@{}[urr]|{\rotatebox{15}{\hbox{$\lllra$}}}_{\rotatebox{15}{#1}}}

\vspace{5pt}
\Grandcadre{
\xymatrix @R=0.2cm @C=1cm{
\Sdim\DRRU{Fact \ref{factSDimGdimHdim} \emph{2}}&&&&\Hdim\DRRU{Fact \ref{factKdimHdim}}
\\
                                       &&\Gdim\URRU{Thm. \ref{MAINCOR}}
                                         \DRRD{Thm. \ref{matrix}}
                                       &&&&\Kdim
\\
\Bdim\URRD{Fact  \ref{factGdimBdim}}   &&&&\Cdim\URRD{Def. \ref{defisutnstab}}
\\
}
}


\subsec{With the stability of a \sutz}


\Grandcadre{In this subsection, $D:\gA\to\gT$ is a fixed \sutz}

We fix the following notations, analogous to those used to define $\Gdim\gA<n$ in \Dfnz~\ref{defiSdimGdim}.
\begin{notation}
\label{nota0Matrix}
{\rm Let $F=[\,C_0\,|\,C_1\,|\,\dots\,|\,C_p\,]$ be a matrix in $\Ae{m\times (p+1)}$ (the~$C_i$'s are the columns) and $G=[\,C_1\,|\,\dots\,|\,C_p\,]$, such that~$F=[\,C_0\,\vert\,G\,]$.
}
\end{notation}


Notice that for every $n$ we have $\DA\big(C_0,\cD_n(F)\big)=\DA\big(C_0,\cD_n(G)\big)$, and a fortiori $D\big(C_0,\cD_n(F)\big)=D\big(C_0,\cD_n(G)\big)$.

\pagebreak	

\begin{lemma}
\label{cor0main}
Suppose that $D$ is $n$-stable and take the notation \ref{nota0Matrix} with   $m=p=n$. Let $\delta=\det(G)$.
There exist $x_1$, \dots, $x_n$ such that

\snic{D(C_0,\delta)\leq  D\big(C_0+\delta(x_1C_1+\cdots+x_nC_n)\big).}

\end{lemma}

\begin{proof}
It suffices to realize $D(\delta)\leq  D\big(C_0+\delta(x_1C_1+\cdots+x_nC_n)\big)$, \cad 

\snic{D(\delta)\leq D(C_0+\delta G X)$ for some $X\in\Ae n.}

Let $\wi{G}$ be the adjoint matrix of $G$ and $L = \wi{G}C_0$.
For any $X\in \Ae n$, we have $\widetilde {G}(C_0+\delta GX) = L+\delta^{2} X$, so $\DA(L+\delta^{2} X)\leq \DA(C_0+\delta GX)$, and a fortiori~\hbox{$D( L+\delta^{2} X) \leq  D(C_0+\delta GX)$}.
Since  $D$ is $n$-stable, by Fact~\ref{lemsutnstab}, we have some $X\in \Ae n$ such that $D(\delta)\leq  D(L+\delta^2X)$.\\
Therefore  $D(\delta) \leq D(C_0+\delta GX)$, as required.
\end{proof}

\begin{theorem}
\label{matrix}  
\emph{(Coquand's \thoz, 1: Forster-Swan and others with the $n$-stability)} 
We have  $\Gdim\gA\leq \Cdim\gA$.
Consequently, \SSOz, Forster-Swan' and Bass' cancellation \thos (\ref{thSerre}, \ref{thSwan},  \ref{thSwan2},   \ref{thBassCancel2}) apply with the~$\Cdim$.
\end{theorem}
\begin{proof} We assume $\Cdim\gA<n$ and we prove $\Gdim\gA< n$.
We use the \carn of $\Gdim\gA< n$ given in Proposition~\ref{propGdimGdim}.
Lemma \ref{cor0main} with the support $D=\rD_{\gA/\! \gen{d}}$ tells us that the equivalent property described in~\ref{propGdimGdim} is satisfied if $\Cdim\aqo\gA d<n$.
We conclude by observing that
$\Cdim\aqo\gA d\leq \Cdim\gA$.
\end{proof}
%

\begin{theorem}
\label{matrixC}  \emph{(Coquand's \thoz, 2: \elrs column operations, support and $n$-stability)}
With the notations \ref{nota0Matrix}, let $n\in\lrbp$. If $D$ is $n$-stable
 there exist $t_1$, \dots, $t_p\in \cD_n(G)$ such that

\snic{
 D\big(C_0,\cD_n(G)\big)\leq  D(C_0+t_1C_1+\cdots+t_pC_p).}
 
\end{theorem}
The \dem of this \tho as a consequence of Lemma~\ref{cor0main} is analogous to the \dem of the difficult implication in Proposition~\ref{propGdimGdim}, in a slightly different context. The result is stronger because Proposition~\ref{propGdimGdim} is only interested in the special case given in  Corollary \ref{t0basic}, with in \hbox{addition $D=\DA$}.
\begin{proof} 
We need to find $t_1$, \ldots, $t_p$ in $\cD_n(G)$ such that,
 for every minor $\nu$ of order $n$ of $G$, we have $D(C_0,\nu)\leq D(C_0+t_1C_1+\cdots+t_pC_p)$.  
\\ 
Actually it suffices to know how to realize

\snic{D(C_0,\delta)\leq  D\big(C_0+\delta (x_1C_1+\cdots+x_pC_p)\big)}

for \emph{one} minor $\delta$ of order $n$ of $G$, and as previously mentioned, for this $D(\delta)\leq  D\big(C_0+\delta (x_1C_1+\cdots+x_pC_p)\big)$ is sufficient.\\
Indeed in this case, we replace $C_0$ by $C'_0=C_0+\delta (x_1C_1+\cdots+x_pC_p)$ in~$F$ (without changing $G$), and we can pass to another minor~$\delta'$ of~$G$ for which we will obtain $x'_1$, \dots, $x'_p$, satisfying

\snic{D(C_0,\delta,\delta')\leq D(C'_0,\delta')\leq
 D\big(C'_0+\delta'(x'_1C_1+\cdots+x'_pC_p)\big)=  D(C''_0),}

with $C''_0=C_0+t''_1C_1+\cdots+t''_pC_p$ and so forth.
\\
To realize the in\egt

\snic{D(\delta)\leq  D\big(C_0+\delta (x_1C_1+\cdots+x_pC_p)\big)}

  for some minor $\delta$ of order $n$ of $G$, we use Lemma~\ref{cor0main} with the extracted matrix~$\Gamma$ corresponding to the minor~$\delta$, and for $C_0$ we limit ourselves to the rows of~$\Gamma$, 
which gives us a vector $\Gamma_0$.
We obtain \hbox{some $X\in\Ae n$} such that

\snic{
D(\delta)\leq  D(\Gamma_0+\delta \Gamma X)\leq D(C_0+\delta G Z).
}

where $Z\in \gA^{p}$ is obtained by completing $X$ with $0$'s.
\end{proof}

Still with the notations \ref{nota0Matrix}, we obtain as a corollary the following result, which implies, when $D=\DA$,
that $\Gdim\gA\leq \Cdim\gA$.
\begin{corollary}
\label{t0basic}
With the notations \ref{nota0Matrix}, let $n\in\lrbp$. \\
If $D$ is $n$-stable and $1=D\big(C_0,\cD_n(G)\big)$, 
there exist $t_1$, \dots, $t_p$ such that the vector $C_0+t_1C_1+\cdots+t_pC_p$ is $D$-\umdz. 
\end{corollary}

\subsec{With the Heitmann dimension}
\begin{lemma}
\label{mainlemma}
We consider a matrix of the form
$$
\cmatrix{
b_0&c_1&\cdots&\cdots&c_n\cr
b_1&a&0&\cdots&0\cr
\vdots&0&\ddots&\ddots&\vdots\cr
\vdots&\vdots&\ddots& \ddots&0\cr
b_n&0&\cdots&0&a\cr
}
,$$
for which we denote the columns   by $V_0$, $V_1$, \ldots, $V_n$.
\\ If $\Hdim\gA < n$ and
$1 = \DA(b_0,a)$,
then there exist $x_1$, \dots, $x_n\in a\gA$ such that
$ 1 = \DA(V_0+x_1V_1+\cdots+x_nV_n).$
\end{lemma}

\vspace{.1em}
\begin{proof}
The \dem is by \recu on $n$.  For $n=0$, it is clear.  \\
If $n>0$, let $\fj=\IH_\gA(b_n)$.  We have $b_n\in \fj$ and $\Hdim\gA/\fj< n-1$, therefore by \hdrz, we can find $y_1$, \dots, $y_{n-1}\in \gA$ such that
$$\preskip.4em \postskip.4em
1 = \rD(U_0+ay_1U_1+\cdots+ay_{n-1}U_{n-1})\quad \mathrm{in}\;\gA/\fj,\eqno (\alpha)
$$
where $U_i$ designates the vector $V_i$ minus its last \cooz.\\
Let $U'_0=U_0+ay_1U_1+\cdots+ay_{n-1}U_{n-1}$, we have $\DA(U'_0,a)=\DA(U_0,a)$.  
The \egt $(\alpha)$  means that there exists a $y_n$ such that $b_ny_n\in \JA(0)$ and
$$\preskip.4em \postskip.4em
1 = \DA(U'_0) \vu \DA(b_n,y_n).\eqno (\beta)
$$
Let $V'_0 = V_0+ay_1V_1+\cdots+ay_{n-1}V_{n-1}+ay_nV_n$. The lemma is proven if $1 \in \DA(V'_0)$. 
Notice that $V'_0$ minus its last \coo is the vector $U'_0 + a_ny_nU_n$ and that its last \coo is $b_n + a^2 y_n$, hence the tight game that comes with $b_n$, $a$, $y_n$.
We have
$$\preskip.4em \postskip.4em
\DA(U_0'+ay_nU_n) \vu \DA(a) = \DA(U'_0,a) = \DA(U_0,a)
\supseteq\DA(b_0,a) = 1,\eqno (\gamma)
$$
and, by $(\beta)$,
$$\preskip.4em \postskip.4em
\DA(U_0'+ay_nU_n)\vu \DA(b_n,y_n) = \DA(
U_0') \vu \DA(b_n,y_n) = 1.\eqno (\delta)
$$
Next $(\gamma)$ and $(\delta)$ imply
$$\preskip.4em \postskip.4em
 \DA(U_0'+ay_nU_n) \vu \DA(b_n,a^2y_n) = 1 =
\JA(U_0'+ay_nU_n,b_n,a^2y_n),\eqno (\eta)
$$
and by Lemma~\ref{gcd2}, since $b_na^2y_n\in \JA(0)$,
$$\preskip.4em \postskip.4em
1 = \JA(U_0'+ay_nU_n,b_n+a^2y_n),$$
\cad $1=\DA(V'_0)$.
\end{proof}

\begin{theorem}
\label{MAINCOR} \emph{(Coquand's \thoz, 3: Forster-Swan and others with the Heitmann dimension)} 
We have  $\Gdim\gA\leq \Hdim\gA$. Consequently,  \SSOz, Forster-Swan' and Bass' cancellation \thos apply with the~$\Hdim$ (\Thos \ref{thSerre}, \ref{thSwan},  \ref{thSwan2},  
\ref{thBassCancel2}).
\end{theorem}
\begin{proof} We use the \carn of $\Gdim\gA< n$ given in Proposition~\ref{propGdimGdim}.
Lemma \ref{mainlemma} tells us that the equivalent property described in~\ref{propGdimGdim} is satisfied if $\Hdim\aqo\gA d<n$.
We conclude by noticing that 
$\Hdim\aqo\gA d\leq \Hdim\gA$.
\end{proof}

\REM{Final}
 All the \thos of commutative \alg which we have proven in this chapter are ultimately brought back to  \thos regarding  matrices and their \mlrsz.
\eoe

\Exercices

\begin{exercise}
\label{exo16Lecteur}
{\rm   
Explicate the computation that gives the \dem of \thref{thKroH} in the %
case $n=1$.
}
\end{exercise}

\vspace{-1em}
\begin{exercise}\label{exoRegularSequence1} {(A \prt of \ndzes sequences)}
\\
 {\rm  
Let $(\an)$ be a \ndze sequence of $\gA$  and $\fa = \gen {\an}$ ($n\geq1$).
\\
\emph {1.}
Show that  $(\ov {a_1}, \ldots, \ov{a_n})$ is an $(\gA\sur\fa)$-basis of $\fa\sur{\fa^2}$.
\\
\emph {2.}
Deduce, when $1\notin\fa$, that $n$ is the minimum number of generators of the \id $\fa$.
For example, if $\gk$ is a nontrivial \ri and $\gA = \gk[\Xm]$,
then for $n \le m$, the minimum number of generators of the \id $\gen {\Xn}$ is $n$.
}
\end{exercise}

\vspace{-1em}
\begin{exercise}\label{exoNbGensIdealBis}
\label{exoNbGensIdeal} (Number of \gtrs  of $\fa/{\fa^2}$ and of $\fa$)
 \\ 
{\rm
Let $\fa$ be a \itf of $\gA$ with 
$\fa/\fa^2 = \gen {\ov {a_1}, \cdots,\ov{a_n}}$. 
\\
 \emph {1.}
Show that  $\fa$ is generated by $n+1$ \eltsz.
\\
 \emph {2.}
Show that  $\fa$ is \lot generated by $n$ \elts in the following precise sense: there exists an $s \in \gA$ such that over the two localized \ris $\gA_s$ and $\gA_{1-s}$, 
$\fa$ is generated by $n$ \eltsz.
\\
 \emph {3.}
Deduce that if $\gA$ is \lgb (for example if $\gA$ is \rdt \zedz), then $\fa$ is generated by $n$ \eltsz.
}
\end{exercise}

\vspace{-1.2em}
\pagebreak	
\begin{exercise}
\label {exoGensPolIdeal}
{\rm
\emph{1.}
Let $E$ be an \Amo and $F$ be a \Bmoz.
If $E$ and $F$ are generated by $m$ \eltsz, the same goes for the $(\gA\times\gB)$-module $E\times F$.
\\
 \emph{2.}
Let $\fa \subseteq \gA[X]$ be an \id containing a \spl \pol $P = \prod_{i = 1}^s (X - a_i)$. Let $\ev_{a_i} :
\gA[X] \twoheadrightarrow \gA$ be the evaluation morphism 
that specializes $X$ in $a_i$. Suppose that each $\fa_i := \ev_{a_i}(\fa)$ is generated by $m$ \eltsz. Show that  $\fa$ is generated by $m+1$ \eltsz.
\\
 \emph{3.}
Let $\gK$ be a \cdi and $V \subset \gK^n$ be a finite set. 
Show that  the \id

\snic{\fa(V) = \sotq{ f \in \KXn}{ \forall\ w \in V,\ f(w) = 0}}

is generated by $n$ \elts (note that this bound does not depend on $\#V$ and that the result is clear for $n=1$).

}
\end{exercise}

\vspace{-1em}
\begin{exercise}\label{exoCubiqueGaucheP3}
 {(The left cubic of $\PP^3$, image of $\PP^1$
under the Veronese embedding of degree $3$)} 
{\rm  
The base \ri $\gk$ is arbitrary, except in the first question where it is a \cdiz. We define the Veronese morphism $\psi : \PP^1 \to \PP^3$ by

\snic {
\psi : (u : v) \mapsto (x_0 : x_1 : x_2 : x_3) \quad\hbox {with}\quad 
x_0=u^3,\ x_1=u^2v,\  x_2=uv^2,\  x_3=v^3
.}

\emph{1.}
Show that  $\Im\psi = \cZ(\fa)$ where $\fa = \gen {D_1, D_2, D_3}=\cD_2(M)$ with the matrix 

\snic{M=\cmatrix {X_0 & X_1 & X_2\cr X_1 & X_2 & X_3\cr},}

\snic{D_1 = X_1X_3-X_2^2,\; D_2 = -X_0X_3+X_1X_2,\; 
D_3 = X_0X_2-X_1^2.}

\emph{2.}
Show that $\fa$ is the kernel of $\varphi : \gk[X_0,X_1,X_2,X_3] \to \gk[U,V]$, $X_i \mapsto U^{3-i}V^i$.  In particular, if $\gk$ is integral, $\fa$ is prime and if $\gk$ is reduced, $\fa$ is radical.
We will show that by letting

\snic {
\fa^\bullet = \gA \oplus \gA X_1 \oplus \gA X_2
\quad \hbox {with} \quad \gA = \gk[X_0,X_3],
}

we get
$\gk[X_0,X_1,X_2,X_3] = \fa + \fa^\bullet \; \hbox { and } \;
\ker\varphi \cap \fa^\bullet = 0$.

\emph{3.}
Show that $\fa$ cannot be generated by two \gtrsz.

\emph{4.}
Explicate a \hmg \pol $F_3$ of degree $3$ such that 
$\DA  (\fa) = \DA  ({D_1,F_3})$. In particular,
if $\gk$ is reduced, $\fa = \DA  ({D_1,F_3})$.
}

\end{exercise}

\vspace{-1em}
\begin{exercise}
\label{exo0KroLocvar}
{\rm
Show that if two sequences are disjoint (see \paref{propdefdisjointes}) they remain disjoint when we multiply one of the sequences by an \elt of the \riz.
}
\end{exercise}

\vspace{-1em}
\begin{exercise}\label{exoJPFurter}
{(Transitivity of the action of $\GL_2(\gk[x,y])$ on the \syss of two \gtrs of $\gen {x,y}$)} 
{\rm  The result of question \emph{1} is due to Jean-Philippe Furter, of the Universit\'e de La Rochelle.
\\
Let $\gk$ be a \riz, $\gA = \gk[x,y]$ and $p$, $q \in \gA$ satisfying 
$\gen {p,q} = \gen {x,y}$.

\emph {1.}
Construct a matrix $A \in \GL_2(\gA)$ such that $A \cmatrix {x\cr y} = \cmatrix {p\cr q}$ and $\det(A) \in \gk^\times$.

\emph {2.}
We write $p = \alpha x + \beta y + \dots$, $q = \gamma x + \delta y + \dots$
with $\alpha$, $\beta$, $\gamma$, $\delta \in \gk$. 
\begin{itemize}\itemsep0pt
\item [\emph {a.}]
Show that $\cmatrix {\alpha & \beta\cr \gamma &\delta} \in \GL_2(\gk)$.

\item [\emph {b.}]
Let $G \subset \GL_2(\gA)$ be the intersection of $\SL_2(\gA)$ and of the kernel of the \homo \gui{reduction modulo $\gen {x,y}$} $\GL_2(\gA) \to \GL_2(\gk)$. The subgroup $G$ is distinguished in $\GL_2(\gA)$.
The subgroup $G\,\GL_2(\gk) = \GL_2(\gk)\,G$ of $\GL_2(\gA)$ operates transitively on the \syss of two \gtrs of $\gen {x,y}$.
\end{itemize}

\emph {3.}
Let $p = x + \sum_{i+j=2} p_{ij}x^iy^j$, $q = y + \sum_{i+j=2}
q_{ij}x^iy^j$. We have $\gen {x,y} = \gen {p,q}$
\ssi the following \eqns are satisfied

\snic { \arraycolsep2pt
\begin{array}{rclc} 
p_{20}p_{02} + p_{02}q_{11} + q_{02}^2  &  = &  
p_{20}p_{11} + p_{02}q_{20} + p_{11}q_{11} - p_{20}q_{02} + q_{11}q_{02} & = 
\\[.5mm] 
p_{20}^2 + p_{11}q_{20} + q_{20}q_{02}  & =  &  0 
\end{array}
}

\emph {4.}
Generalize the result of the previous question.

}

\end{exercise}

\vspace{-1em}
\begin{exercise}\label{exoSdimSmithRing}
{(About Smith \ris and  $\Sdim$)}
\\
{\rm For the notions of a strict Bézout \ri and a Smith \riz, see Section~\ref{secBézout} on \paref{secpfval} and Exercises~\ref{exoAnneauBézoutStrict} and~\ref{exoSmith}. 
Exercise \ref{exoSmith} gives a direct solution of item \emph{5.}

\emph {1.}
If $\gA$ is a Smith \riz, we have $\Sdim\gA \le 0$.\\
Deduce $\Sdim\ZZ$, $\Bdim\ZZ$, $\Gdim\ZZ$ and $\Cdim\ZZ$.

In questions \emph {2} and \emph {3}, the \ri $\gA$ is arbitrary.

\emph {2.}
Let $A \in \MM_2(\gA)$ and $u \in \Ae2$ be a \umd vector.
Show that $u \in \Im A$ \ssi there exists a $Q \in \GL_2(\gA)$
such that $u$ is the first column of $AQ$.

\emph {3.}
Let $A \in \MM_2(\gA)$ of rank $\ge 1$. Then $A$
is \eqve to a diagonal matrix \ssi $\Im A$ contains a \umd vector.

\emph {4.}
Let $\gA$ be a strict Bézout \riz. Show that $\Sdim\gA \le 0$ \ssiz$\gA$ is a Smith \riz.

\emph {5.}
Deduce that a \ri $\gA$ is a Smith \ri \ssi
it is a strict Bézout \ri and if for \com $a$, $b$, $c$ the matrix $\cmatrix {a & b\cr 0 &c\cr}$
has a \umd vector in its image. This last condition can be expressed by the condition known as the Kaplansky condition:
$
1 \in \gen {a,b,c}  \Rightarrow
\hbox { there exist }p,\,q\hbox{ such that }1 \in \gen {pa, pb + qc}. 
$

Remark: we dispose of the \elr \carnz: $\gA$ is a strict Bézout \ri \ssi for every $a$, $b\in \gA$, there exist $d$ and \com $a'$, $b'$ such that $a=da'$ \hbox{and $b=db'$}. If we add the Kaplansky condition above, we obtain an \elr \carn of  Smith \risz.

}

\end{exercise}


\penalty-2500
\sol

\exer{exo16Lecteur} 
The given \dem says this.
Since $\Kdim \gA\leq 0$, there exists some $x_1$ such that $b_1x_1\in\DA(0)$ and $1\in\DA(b_1,x_1)$. A fortiori  $b_1ax_1\in\DA(0)$ and $a\in\DA(b_1,ax_1)$. Lemma~\ref{gcd} tells us that $\DA(b_1,ax_1)=\DA(b_1+ax_1)$, so $a\in\DA(b_1+ax_1)$.


\exer{exoRegularSequence1}
\emph {1.}
Let $b_1$, \dots, $b_n \in \gA$ such that $\sum_i \ov{b_i} \ov{a_i} = 0$ in $\fa\sur{\fa^2}$. In other words $\sum_i b_ia_i = \sum_i c_ia_i$ with $c_i \in \fa$. By Lemma~\ref{PetitLemmeAlterne}, there exists an alternating matrix $M \in \Mn(\gA)$ such that $[\,b_1-c_1 \; \cdots \; b_n-c_n\,] =  [\,a_1 \; \cdots \; a_n\,]M$.\\
Hence $b_i - c_i \in \fa$, and so $b_i \in \fa$. 

 \emph{Same thing, presented more abstractly.} We know that a \mpn of the \Amo $\fa$ for the \sgr $(\an)$ is $R_{\ua}$.  By changing the base \ri $\gA\to\gA\sur\fa$, this gives a null \mpn ($R_\ua\mod\fa $)  of the $\gA\sur\fa$-module $\fa/{\fa^2}$ for $(\ov{a_1},\ldots,\ov{a_n})$, which means that this \sys is a basis.

\emph {2.}
If $(\yp)$ is a \sgr of the \id $\fa$, $(\ov {y_1}, \ldots, \ov {y_p})$ is a \sgr of the free $(\gA\sur\fa)$-module $\fa/{\fa^2}$ of rank $n$. Therefore, if $p<n$,  $\gA\sur\fa$ is trivial.

\exer{exoNbGensIdeal}
\emph{1.} By letting $\fb = \gen {a_1, \cdots, a_n}$, the \egt $\fa/\fa^2 = \gen {\ov {a_1}, \cdots,\ov{a_n}}$ means that $\fa = \fb + \fa^2$. We then have $(\fa/\fb)^2 = (\fa^2 + \fb)/\fb = \fa/\fb$, and the \itf $\fa/\fb$ of~$\gA/\fb$ is \idmz, therefore generated by an \idmz. Therefore there exists some $e \in \fa$, \idm modulo $\fb$, such that $\fa = \fb + \gen {e}$: $\fa=\gen{a_1, \ldots, a_n, e}$.

\emph {2.}
With the same notations we see 
that $(1-e)\fa \subseteq \fb + \gen {e^2 - e} \subseteq \fb$.\\
Therefore in $\gA_{1-e}$, $(\an)$ generates $\fa$ whereas in $\gA_e$,
$1 \in \fa$.

\emph{Variant.}
We introduce $S = 1 + \fa$ and work on $\gA_S$:  $\fa\gA_S\subseteq\Rad(\gA_S)$ and so, by Nakayama, a \sgr of $\fa_S\sur{\fa_S^2}$ is also a \sgr of $\fa_S$. We therefore have $\fa_S = \fb_S$, hence the existence of some $s \in S$ such \hbox{that $s\fa \subseteq\fb$}. \hbox{In
$\gA_s$, $(\an)$} generates $\fa$, whereas in $\gA_{1-s}$, $1 \in \fa$
($s\in1+\fa$, \hbox{so $1-s \in \fa$}).

\exer{exoGensPolIdeal}
Item \emph{1} is obvious, and we deduce \emph{2} since $\aqo{\fa}{P} \simeq \fa_1 \times \cdots \times \fa_s$.
We deduce \emph{3} by \recu on $n$. We observe that the Chinese remainder \tho used in item \emph{2} is concretely realized by the interpolation \`a la Lagrange.
\\
 Note: see also Exercise~\ref{exoGensIdealEnsFini}.



\exer{exoCubiqueGaucheP3} 
\emph {1.}
$\psi$ is \hmg of degree $3$. Let $p = (x_0:x_1:x_2:x_3)$ \hbox{in $
\cZ(\fa)$}. If $x_0 \ne 0$, we are brought back to $x_0 = 1$, so $(x_0, x_1, x_2, x_3) = (1, x_1, x_1^2, x_1^3) = \psi(1 : x_1)$. If $x_0 = 0$, then
$x_1 = 0$, then $x_2 = 0$, so $p = \psi(0 : 1)$.

\emph {2.}
Let $\kux = \kuX\sur\fa$ and $\ov\gA = \gk[x_0,x_3]$. Showing the \egt $\kuX = \fa + \fa^\bullet$ amounts to showing that $\kux = \ov\gA + \ov\gA x_1 + \ov\gA x_2$.  We have the relations $x_1^3 = x_0^2x_3 \in \ov\gA$, \hbox{and $x_2^3 = x_0x_3^2 \in \ov\gA$}, therefore $\ov\gA[x_1,x_2]$ is the
$\ov\gA$-module generated by the $x_1^ix_2^j$'s for \hbox{$i$, $j \in \lrb{0..2}$}. But we also have $x_1x_2 = x_0x_3$, $x_1^2 = x_0x_2$, $x_2^2 = x_1x_3$, which completes the proof of $\ov\gA[x_1,x_2] = \ov\gA + \ov\gA x_1 + \ov\gA x_2$.

Let $h = a + bX_1 + cX_2 \in \fa^\bullet$ satisfy $\varphi(h) = 0$
($a, b, c \in \gA = \gk[X_0,X_3]$). We therefore have
$$\preskip.0em \postskip.4em
a(U^3,V^3) + b(U^3,V^3)U^2V + c(U^3,V^3)UV^2 = 0
.$$
By letting $p(T) = a(U^3,T)$, $q(T) = b(U^3,T)U^2$, $r(T) = c(U^3,T)U$, we obtain the \egt $p(V^3) + q(V^3)V + r(V^3)V^2 = 0$, and an examination modulo $3$ of the exponents in~$V$ of $p$, $q$, $r$ provides $p=q=r=0$. Hence $a=b=c=0$, \cad $h = 0$.
Now, \hbox{if $f \in \ker\varphi$}, by writing $f = g + h$ with
$g \in \fa$, $h \in \fa^\bullet$, we obtain $h \in \ker\varphi \cap \fa^\bullet = 0$, so $f = g \in \fa$.

\emph {3.}
Let $E = \fa/\!\gen{\uX}\fa$. It is a $\aqo{\kuX}{\uX}$-module generated by $d_i = \ov {D_i}$.
 In other words  $E = \gk d_1 + \gk d_2 + \gk d_3$.  Moreover, $d_1$, $d_2$, $d_3$ are $\gk$-\lint 
independent. Indeed, if $ad_1 + bd_2 +cd_3 = 0$, then $aD_1 + bD_2 + cD_3 \in \gen{\uX}\fa$,
which for homogeneity reasons gives $aD_1 + bD_2 + cD_3 = 0$, then $a = b = c = 0$. Therefore $E$ is free of rank $3$ over $\gk$.  If $G$ is a \sgr of $\fa$, \hbox{then $\ov G$} is a \sgr of the \kmo $E$, therefore $\#\ov G \ge 3$, a fortiori $\# G \ge 3$.

\emph {4.}
Let $F_3 = X_0D_2 + X_1D_3 = -X_0^2X_3 + 2X_0X_1X_2 - X_1^3  \in 
\gen {D_2, D_3}$.
 We have
$$\preskip.4em \postskip.4em
\begin {array} {c}
D_2^2 = -(X_3F_3 + X_1^2D_1) \in \gen {D_1, F_3}, \quad 
D_3^2 = -(X_1F_3 + X_0^2D_1) \in \gen {D_1, F_3}, \\[1mm]
D_2D_3 = X_0X_1D_1 + X_2F_3  \in \gen {D_1, F_3} \;\;
\hbox{then}\\[1mm]
\gen {D_1, D_2, D_3}^2 \subseteq \gen {D_1, F_3} \subseteq
\gen {D_1, D_2, D_3}, \; \hbox {hence} \;
\sqrt {\gen {D_1, D_2, D_3}} = \sqrt {\gen {D_1, F_3}}.
\end {array}
$$


\exer{exoJPFurter}
 \emph {1.}
Let us first notice that for $m_{ij} \in \gA=\gk[x,y]$, an \egt
$$
\cmatrix {m_{11} & m_{12}\cr m_{21} & m_{22}} \cmatrix {x\cr y} =
\cmatrix {0\cr 0}
$$ 
entails $m_{ij} \in \gen {x,y}$. Moreover, we will use the following \idts for $2\times 2$ matrices:
$\det(A + B) = \det(A) + \det(B) + \Tr(\wi AB)$ and

\snic {
\hbox {for } H = \crmatrix {v\cr -u} [\,y\;-x\,], \quad
\Tr(\wi AH) = [\,u\;v\,]\,A\cmatrix {x\cr y} 
.}

By hypothesis, we have $A$, $B \in \MM_2(\gA)$ such that 
$$A \cmatrix
{x\cr y} = \cmatrix {p\cr q}\hbox{  and  }B \cmatrix {p\cr q} = \cmatrix {x\cr y}$$
therefore $(BA - \I_2) \cmatrix {x\cr y} = \cmatrix {0\cr 0}$.  Thus, modulo $\gen
{x,y} = \gen {p,q}$, we have $BA \equiv \I_2$. Therefore $a = \det(A)(0,0) \in \gk\eti$
and we can express, with $u$, $v \in \gA$, $\det(A) = a + up + vq$. Let $H = \crmatrix {v\cr -u} [\,y\;-x\,]$. We have $H\cmatrix {x\cr y} = \cmatrix {0\cr 0}$, $\det(H) = 0$, and we change $A$ to $A' = A-H$.  Then $A' \cmatrix {x\cr y} = \cmatrix {p\cr q}$ and 
$$
\det(A') = \det(A) + \det(H) -
\Tr(\wi AH) = a + up + vq - [\,u\;v\,] \cmatrix {p\cr q} = a.
$$

\emph {2.}
We decompose $A$ into \hmgs components: $A = A_0 + A_1 + \dots$, and we examine the \egt $A\cmatrix {x\cr y} = \cmatrix {p\cr q}$. \\
The examination of the \hmg component of degree $1$ gives 
\smash{$A_0 = \cmatrix {\alpha & \beta\cr \gamma &\delta}$}, and we know that $\det(A) = \det(A_0) \in \gk^\times$. 
\\
We then can write $A_0(A_0^{-1}A)\cmatrix {x\cr y} = \cmatrix {p\cr q}$ with $A_0 \in \GL_2(\gk)$ and $A_0^{-1}A \in G$.

\emph {3.}
We write $A\cmatrix {x\cr y} = \cmatrix {p\cr q}$ with $A \in G$.
For degree reasons, we obtain an \egt $A = \I_2 + xB + yC$ with $B$, $C \in \MM_2(\gk)$.
We then have
$$ 
\begin{array}{ccc} 
 \cmatrix {p\cr q} = A\cmatrix {x\cr y} = 
\cmatrix {x\cr y} + B\cmatrix {x^2\cr xy} + C\cmatrix {xy\cr y^2}
= \\[1.2em] 
 \cmatrix {x + b_{11}x^2 + (c_{11}  + b_{12})xy + c_{12}y^2 \\[.3em]
y + b_{21}x^2 + (c_{21} + b_{22})xy + c_{22}y^2} 
\end{array}
\eqno(\star)
$$
Moreover, we notice that the \coe of $\det(A)-1$ in $x^iy^j$ is a \pog of degree $i+j$ in the \coes of $B$ and $C$

\snic {
\det(A)-1 = \Tr(B)x + \Tr(C)y + \det(B)x^2 + \Tr(\wi BC)xy + \det(C)y^2
.}


If $\gk$ was an \ac field, we could give the following argument.
The \egt $\det(A) = 1$ defines a \pro sub\vrt $V \subset \PP^{8-1}$ ($2 \times 4$ \coes for $(B, C)$); on the other hand $(\star)$ defines a morphism $V \to \PP^{6-1}$ ($6$ for the \coes of $p-x$, $q-y$). The image of this morphism is the set $W$ defined by the \eqns of the statement.

The \ri $\gk$ being arbitrary, we carefully examine the \eqns $(\star)$; by using $\Tr(B) = \Tr(C) = 0$, we can express $B$ and $C$ in terms of the \coes of $p$ and~$q$

\snic {
B = \cmatrix {p_{20} & p_{11} + q_{02}\cr q_{20} & -p_{20}}, \qquad
C = \cmatrix {-q_{02} & p_{02}\cr p_{20}+q_{11} & q_{02}}
.}

We thus construct a section $s : W \to G$ of the map $(\star)$, and in fact the three \eqns of $W$ appearing in the statement are, up to sign, $\det(C)$, $\Tr(\wi BC)$ \hbox{and $\det(B)$}.

\exer{exoSdimSmithRing} 
\emph {1.}
A \gui{diagonal} rectangular matrix of rank $\ge 1$ has in its image a \vmd (this for every \riz).  Let $A$ be a matrix of rank $\geq 1$, if $\gA$ is a Smith \riz, $A$ is \eqve to a \gui{diagonal} matrix~$D$, therefore $\Im D$ contains a \vmdz, and \egmt $\Im A$.

We therefore have $\Sdim\ZZ = 0$. Moreover, $\Cdim\ZZ \le 1$
($\ZZ$ is $2$-stable because $\ZZ$ is a Bézout domain). Finally, $\Bdim\ZZ > 0$ because $1 \in \gen {2,5}$ without finding some $x \in \ZZ$ such that $1 \in \gen {2 + 5x}$. \\
Recap: $\Bdim\ZZ = \Gdim\ZZ = \Cdim\ZZ = 1$ but $\Sdim\ZZ = 0$.

\emph {2.}
If $u = Av$, then $v$ is \umdz. \hbox{Therefore $v = Q\cdot e_1$} \hbox{with $Q \in \SL_2(\gA)$} and $u$ is the first column of $AQ$. The other direction is immediate.

\emph {3.}
Suppose that $\Im A$ contains a \vmdz. By item \emph{2}, we have~$A \sim B$ with $B\cdot e_1$ \umdz. Therefore the space of rows of $B$ contains a vector of the form $[\,1\ *\,]$. 
Item \emph{2} for $\tra B$ gives

\snic {
\tra {B} \sim \cmatrix {1 &*\cr * &*\cr} \sim \cmatrix {1 & 0\cr 0 &*\cr}, \hbox{ diagonal.}
}

Recap: $A$ is \eqve to a diagonal matrix. The other direction is \imdz.

\emph {4.}
Let $\gA$ be a strict Bézout \ri with $\Sdim\gA \le 0$. We show that every triangular matrix~$M \in \MM_2(\gA)$ is \eqve to a diagonal matrix. 
\\
We can write $M = dA$ with $A$ of rank $\ge 1$ (because $\gA$ is a strict Bézout \riz).
\\
Since $\Sdim\gA \le 0$, $\Im A$ contains a \vmd therefore is \eqve to a diagonal matrix $D$. Ultimately $M \sim dD$.

\emph {5.}
Now easy.


\vspace{-.5em}


\pagebreak

\Biblio

If we stick to the \cof aspect of the results, the whole chapter is essentially due to T. Coquand, with at times the help of the authors of the book that you are holding in your hands.
This is a remarkable success for the \cov approach of the theory of the \ddkz. Without this approach, it was simply unthinkable to obtain in a \gnl \cov form the \gui{grand} classical theorems proven here.
In addition, this approach has guided the development of a new dimension, which we call Heitmann dimension, thanks to which the remarkable Heitmann non-\noe results were able to be improved further, namely the \gnl non-\noe version of \SSO \tho and of the Forster-Swan \thoz.

\smallskip
\KRNz's \tho is usually stated in the following form: an algebraic \vrt in $\CC^n$ can always be defined by $n+1$ \eqnsz.
\\
It was extended to the case of \noe \ris by van der Waerden \cite{vW}
in the following form: in a \noe \ri of \ddkz~$n$, every \id has the same nilradical as an \id generated by at most $n+1$ \eltsz.
\\
Kronecker's version was improved by various authors in the articles \cite[Storch]{Stor} and \cite[Eisenbud\&Evans]{E2} which showed that $n$ \eqns \gnlt suffice.
A \prco of this last \tho is in \cite[Coquand\&al.]{cls}.
Moreover, we still do not know whether every curve in the complex space of dimension $3$ is an intersection of two surfaces or not (see \cite{Kun}, Chapter~5).

Lemma~\ref{thCor2.2Heit} is Heitmann's Corollary 2.2 \cite{Hei84}, 
(for us, the $\Hdim$ replaces the $\Jdim$) which leads to \thref{KroH2} (Heitmann's improvement of Kronecker's \thoz).

The local Kronecker's \tho \ref{thKroLoc} is due to Lionel Ducos \cite{Duc08}.

\smallskip
\emph{Note regarding the \gui{stable range.}} \Thref{Bass} is due to Bass in the \noe case (with the dimension of the maximal spectrum, which in this case coincides with the $\Jdim$ and the $\Hdim$) and to Heitmann in the non-\noe case with the $\Jdim$. \Thref{Bass0} is a non-\noee version, but with the \ddkz, of~\thref{Bass}.

\smallskip
\emph{Note regarding the $\Jdim$.} In \cite{Hei84} Heitmann introduces the $\Jdim$ for a not \ncrt \noe \ri as the correct substitute for the dimension of the maximal spectrum $\Max\gA$. It is the dimension of $\Jspec\gA$, the smallest spectral subspace of $\SpecA$ containing $\Max\gA$.
He establishes Bass' \gui{stable range}  \tho for this dimension.
However, for the \thos by Serre and Forster-Swan,
he has to use an ad hoc dimension, the upper bound of the $\Jdim(\gA[1/x])$'s for $x\in\gA$. 
As this ad hoc dimension is bounded above by the Krull dimension anyway, he is then able to obtain, in particular, a non-Noetherian version of the cited grand theorems for the Krull dimension. 

\smallskip \emph{Note regarding Serre's \tho and the Forster-Swan \thoz.} Serre's \tho is in \cite[Serre]{Serre}.
The Forster-Swan \tho (\noe version) is in \cite[Forster]{Forster} for the \ddk and in \cite[Swan]{Swan} for the dimension of the maximal spectrum. 
Non-\noe versions for the Krull dimension are due to Heitmann \cite{Hei76,Hei84}.
Finally, the article by Eisenbud-Evans \cite{Eisenbud} has greatly helped to clarify matters regarding Forster, Swan and Serre's theorems.
\\
Sections~\ref{secSOSFSa} and \ref{secManipElemCol} (second part: Heitmann dimension) are inspired by the outline of~\cite[Eisenbud\&Evans]{Eisenbud} and \cite[Heitmann]{Hei84}. These sections give \cov versions of Serre's (Splitting Off), Forster-Swan's, and Bass' (cancellation) \thosz.
This improves (even without taking into account the \cof aspect of the \demz) all the known \thos on the subject, by answering positively for the Heitmann dimension (and a fortiori for the $\Jdim$) a question left open by Heitmann.

\smallskip
\emph{Note regarding the $\Hdim$}. 
The Heitmann dimension, denoted by $\Hdim$, was introduced in \cite{clq}
(see also \cite{clq2}). Fundamentally it is the dimension that makes the \dems work in the article by Heitmann \cite{Hei84}.
The fact that it is better a priori than the $\Jdim$ is not the core point.
It is rather the fact that Serre's and Forster-Swan's \thos work with the $\Hdim$, and so a fortiori  with the $\Jdim$, which gives the complete non-\noe version of these \thosz, which had been conjectured by Heitmann.
\\
In the case of a \noe \riz, the $\Hdim$, Heitmann's $\Jdim$ and the dimension of the maximal spectrum $\Max \gA$ which intervenes in Serre's and Swan's \cite{Swan} \thos are the same
(refer to \cite{clq2,Hei84}).

\smallskip
\emph{Note regarding  $n$-stability}. 
The notion of a \sut dates back to Joyal \cite{Joyal} and Espa\~ nol \cite{espThesis}, who use it to give a constructive \carn of the \ddk of  commutative \risz.
It is used systematically in the recent articles by T.~Coquand.
In Section~\ref{secSUPPORTS} and in the first part of Section~\ref{secManipElemCol} the notion of an  $n$-stable \sut is decisive.
It was invented by T.~Coquand \cite{coq07} to update the \cof content of Bass' rhetoric on the finite partitions of $\SpecA$ in \cite{BASS}.

\smallskip
The version of Bass' cancellation \tho for the $\Hdim$ was first proved by Lionel Ducos \cite{DQ}. The \dem that we give is based on \cite{clq2} instead.

\smallskip Regarding Exercise~\ref{exoJPFurter}, Murthy, in \cite{Murthy03}, proved the following \gnl result. Let $\gA = \gk[\xm]$ be a \pol \ri ($\gk$ being a commutative \riz) and $r \ge 1$ be fixed.
Suppose, for every $n \in \lrbr$, that every \vmd of $\Ae n$ is completable and consider, \hbox{for $n \le \inf(r,m)$}, the set of \syss of $r$ \gtrs of the \id $\gen {\xn}$ of $\gA$, such as for example
$(\xn, 0, \ldots, 0)$ where there are $r-n$~zeros. Then the \hbox{group $\GL_r(\gA)$} operates transitively on this set (Murthy's result is actually much more precise).

\newpage \thispagestyle{CMcadreseul}
\incrementeexosetprob



\chapter{The \plgz}
\label{chapPlg} \label{chap gen loc}
\perso{compil\'e le \today}
\minitoc

\subsection*{Introduction}
\addcontentsline{toc}{section}{Introduction}

In this chapter, we discuss a few important methods directly related to what is commonly called the \plg
in commutative \algz.

In Section~\ref{subsec loc glob conc} we develop it in the form of \plgcsz. This is to say that certain \prts are globally true as soon as they are locally true. Here the term locally is taken in the \cof sense: after \lon at a finite number of \mocoz.

In Section~\ref{subsec loc glob abs}, we establish the corresponding \plgasz, by using, inevitably, non-\cov \demsz: here locally is taken in the abstract sense, \cad after \lon at any \idepz.

In Section~\ref{secColleCiseaux}, we explain the construction of \gui{global} objects from objects of the same nature only defined locally.

\smallskip
Sections \ref{secMachLoGlo}, \ref{subsecLGIdeMax} and \ref{subsecLGIdepMin} are devoted to the \gui{dynamic and \cov decryption} of methods used in abstract \algz.
Recall that in Section~\ref{subsecDyna} we presented the \gnl philosophy of this dynamic method.

In Section~\ref{secMachLoGlo}, we discuss the \cov decryption of abstract methods that fall within a \gnl framework of the type \gui{\plgz.} We give a \gnl statement (but inevitably a little informal) for this, and we give simple examples, which could be treated more directly. 
The truly pertinent examples will come in Chapter~\ref{ChapMPEtendus}.

This dynamic method is a fundamental tool of \cov \algz.
We could have written this work by starting with this preliminary explanation and by systematically using this decryption.
We preferred to start by developing everything that could be directly developed,
by establishing the \plgcs that usually allow us to avoid using the dynamic decryption as such.
In short, rather than highlighting the magic at work in classical \alg we preferred to first show a different kind of magic at work in \cov \alg under the general slogan:
\gui{why make things complicated when you can make them simple?}

In Section~\ref{subsecLGIdeMax}, we analyze the method of abstract \algz, which consists in \gui{seeing what happens when we quotient by an arbitrary \idemaz.}

In Section~\ref{subsecLGIdepMin}, we analyze the method
which consists in \gui{seeing what happens when we localize at an arbitrary \idemiz.}

In Sections \ref{secPlgcor} and \ref{secPlgprof2}, we examine to what extent certain \plgs remain valid when we replace in the statements the lists of \eco by lists of depth $\geq 1$ or of depth $\geq 2$.

\section{Comaximal \mosz, coverings}
\label{subsecMoco}

We treat in Section~\ref{subsec loc glob conc} concrete versions of principles of the local-global type.
For these concrete versions, the \lon have to be done in a finite number of \eco (or of \mocoz) of $\gA$:
\emph{if the considered \prt is true after \lon at a finite number of \ecoz, then it is true.}

We introduce a \gnnz.
\begin{definition}
\label{def.moco} 
We say that
{\em the \mos $S_1$, \ldots, $S_n$ of the \ri $\gA$ cover the \moz~$S$}
if $S$ is contained in the saturated \mo of each $S_i$ and if an \id of $\gA$ that intersects each of the $S_i$'s always intersects $S$, in other words if we have

\snic{\forall s_1\in S_1 \;\dots\; \forall s_n\in S_n \;\;
\exists$ $a_1$, \ldots, $a_n\in \gA\quad \sum_{i=1}^{n} a_i s_i \in S.
}
\end{definition}

Monoids are \com if they cover the \moz~$\so{1}$.

In \clama (with the axiom of the \idepz)%
{\footnote{The axiom of the \idep affirms that every strict \id of a \ri is contained in a \idepz. This is a weakened version of the axiom of choice.
In the  classical set theory $\ZF$, the axiom of choice is equivalent to the axiom of the \idemaz, which states that every strict \id of a \ri is contained in a \idemaz. This is a little stronger that the axiom of the \idepz.
The latter is equivalent to the fact that every consistent formal theory admits a model (this is the compactness \tho in classical logic).
In classical set theory with the axiom of choice, the axiom of the prime ideal becomes a theorem and is called \gui{Krull's lemma.}%
\index{Krull's lemma}\index{axiom of the \idepz}}} 
we have the \carn given in the following lemma.
For some \mo $S$, we denote by $U_S$ the subset of $\Spec \gA$ defined~by

\snic{U_S=\sotq{\fp\in \Spec \gA}{\fp\cap S=\emptyset}.}

If $S$ is the \mo generated by the \elt $s$, we denote $U_S$ by $U_s$.
From a \cof point of view, $\Spec \gA$ is a topological space known via      
       its basis of open sets 
$U_s=\fD_\gA(s)$ but whose points are often difficult to access.

Recall that we denote by $\sat{S}$ the saturated \mo of the \mo $S$.

\goodbreak
\begin{lemmac}
\label{lemdefmoco}~
\begin{enumerate}
\item For every \mo $S$ we have $\sat S = \bigcap_{\fp\in U_S}(\gA\setminus \fp)$.
Consequently for two \mos $S$ and $T$, $\sat S\subseteq\sat{T}\Leftrightarrow U_{T}\subseteq U_S$.
\item $S_1$, \ldots, $S_n$ are \com \ssi  $\Spec \gA=\bigcup_iU_{S_i}$.
\item $S_1$, \ldots, $S_n$ cover the \mo $S$ \ssi $U_S=\bigcup_iU_{S_i}$.
\end{enumerate}
\end{lemmac}
\begin{proof} \emph{1.} Results from the Krull lemma:  if an \id $\fa$ does not intersect a \mo $S$, there exists a \idep $\fp$ such that $\fa\subseteq\fp$ and $\fp\cap S=\emptyset$.\\
\emph{2.} We can assume that $\gA$ is not trivial.
If the \mos are \com and if $\fp$ is a \idep not belonging to any of the $U_{S_i}$'s, there is in each~$S_i$ an \elt $s_i$ of $\fp$, therefore by the \dfn of the \mocoz,~$1\in\fp$, a contradiction.
Conversely assume \hbox{that $\Spec \gA=\bigcup_iU_{S_i}$} and let $s_1\in S_1$, \ldots, $s_n\in S_n$. If $\gen{s_1,\ldots ,s_n}$ does not contain~$1$, it is contained in a \idep $\fp$. Therefore $\fp$ is in none of \hbox{the $U_{S_i}$'s}, a contradiction.
\end{proof}

The following lemma is a variation on the theme: \index{covering}%
\emph{a covering of coverings is a covering}.
It is also a \gnn of Fact~\ref{factLocCas}. The corresponding computations are immediate.
In \clama it would be even faster via Lemma\eto\ref{lemdefmoco}.
\goodbreak
\begin{lemma} \index{Successive localizations lemma, 2}
\label{lemAssoc} \emph{(Successive \lons lemma, 2)}
\begin{enumerate}
\item \emph{(Associativity)} If the \mos $S_1$, \ldots, $S_n$ of the \ri $\gA$ cover the \mo $S$  and if each $S_\ell$ is covered by \mos $S_{\ell,1},\ldots ,S_{\ell,m_\ell}$, then the $S_{\ell,j}$'s cover $S$.
\item \emph{(Transitivity)}
\begin{enumerate}
\item Let  $S$ be a \mo of the \ri $\gA$ and $S_1,\ldots,S_n$ be \mos of the \ri $\gA_S$.
For $\ell\in\lrbn$ let $V_\ell$ be the \mo of $\gA$ formed by the numerators of the \elts of  $S_\ell$. Then
 the \mos $V_1,\ldots ,V_n$ cover $S$ \ssi the \mos  $S_1$, \ldots, $S_n$ are \comz.
\item More \gnlt let $S_0,\ldots ,S_n$ be \mos of the \ri  $\gA_S$ and for $\ell=0,\ldots,n$ let $V_\ell$ be the \mo of  $\gA$ formed by the numerators of \elts of  $S_\ell$.
Then the \mos $V_1,\ldots ,V_n$ cover  $V_0$ \ssi  $S_1$, \ldots, $S_n$ cover $S_0$ in $\gA_S$.
\end{enumerate}
\end{enumerate}
\end{lemma}

\begin{definota}
\label{nota mopf}
Let $U$ and $I$ be subsets of the \ri $\gA$. Let $\cM(U)$ be the \mo generated by $U$, and  $\,\cS(I,U)$ be the \mo

\snic{\cS(I,U)= \gen {I}_\gA + \cM(U).}

The pair $\fq=(I,U)$ is also called a \emph{\ippz}, and we write (by abuse) $\gA_{\fq}$ for $\gA_{\cS(I,U)}$.
Similarly we let

\snic{
\cS(a_1,\ldots,a_k;u_1,\ldots,u_\ell) = \gen{a_1,\ldots ,a_k}_\gA + \cM(u_1,\ldots,u_\ell). 
}

We say that such a \mo {\em admits a finite description}.
The pair

\snic{(\so{a_1,\ldots,a_k},\so{u_1,\ldots,u_\ell})}

is called a \emph{finite \ippz}.\index{ideal!potential prime ---}%
\index{ideal!finite potential prime ---}
\end{definota}

It is clear that for $u=u_1\cdots u_\ell$, the \mos $\cS(a_1,\alb\ldots,\alb a_k;u_1,\ldots,u_\ell)$ \hbox{and $\cS(a_1,\ldots,a_k;u)$} are \eqvsz, \cad have the same saturated \moz.

\medskip \rem
The \ipp  $\fq=(I,U)$ is constructed for the following goal: \emph{when we localize at $\cS(I,U)$, we obtain $U\subseteq \gA_\fq\eti$  and $I\subseteq \Rad(\gA_\fq)$.} \\
Similarly, for every \idep $\fp$ such that $I\subseteq\fp$ and $U\subseteq \gA\setminus \fp$, we have $U\subseteq \gA_\fp\eti$  and $I\subseteq \Rad(\gA_\fp)$. The pair $\fq=(I,U)$ therefore represents partial information on such a \idepz. It can be considered as an approximation of $\fp$. This explains the terminology of a \ipp and the notation $\gA_{\fq}$.\\ 
We can compare the approximations of $\fp$ by finite \ipps with approximations of a real number by rational intervals.
\eoe

\pagebreak	

\begin{lemma} \index{Successive localizations lemma, 3}
\label{lemRecouvre} \emph{(Successive \lons lemma, 3)}\\
Let $U$ and $I$ be subsets of the \ri $\gA$ and $a\in \gA$,
then the \mos

\snic{\cS(I;U,a)\eqdefi\cS(I,U\cup\so a)\;$ and $\;\cS(I,a;U)\eqdefi\cS(I\cup\so a,U)}

cover the \mo $\cS(I,U)$.\\
In particular, the \mos $S=\cM(a)=\cS(0;a)$ and $S'=\cS(a;1)= 1+a\gA$  are  \comz.
\end{lemma}
\begin{proof}
Let $x\in \cS(I;U,a)$, $y\in \cS(I,a;U)$. 
We need to see that $\gen {x,y}$ meets $\gen {I} + \cM(U)$, or that $\gen {x,y} + \gen {I}$
meets $\cM(U)$.  \\
We have $k \ge 0$, $u, v \in \cM(U)$ and $z \in \gA$ such that
$x \in ua^k + \gen {I}$ and $y \in v - az + \gen {I}.$
Modulo $\gen {x,y} + \gen {I}$, $ua^k \equiv 0$, $v \equiv az$ so $uv^k \equiv 0$, \cad $uv^k \in \gen
{x,y} + \gen {I}$ with $uv^k \in \cM(U)$.
\end{proof}

\comm
The previous lemma is fundamental. It is the \cov counterpart of the following banal observation in \clamaz:
after localizing at a \idep every \elt is found to be either \iv or in the radical.
When dealing with this type of argument in a classical \demz, most of the time it can be interpreted \cot by means of this lemma. Its \dem is very simple, in the image of the banality of the observation made in the classical \demz. But here there is a true computation. We can in fact ask whether the classical \dem avoids this computation. A detailed analysis shows that no: it is found in the \dem of Lemma\eto\ref{lemdefmoco}.
\eoe

\medskip The examples given in the following lemma are frequent.
\begin{lemma}
\label{exaMoco}
Let $\gA$ be a \riz,  $U$ and $I$ be subsets of $\gA$, and $S=\cS(I,U)$.
\begin{enumerate}
\item  If $s_1$, \ldots, $s_n\in \gA$ are \ecoz, the \mosz~$\cM(s_i)$ are \comz.
 More \gnltz, if $s_1$, \ldots, $s_n\in \gA$ are \eco in $\gA_S$, the \mos $\cS(I;U,s_i)$ cover the \mo $S$.

\item  Let $s_1$, \ldots, $s_n\in \gA$. The \mos
$$
\preskip.2em \postskip.2em
\begin{array}{l}
S_1=\cS(0;s_1), \,\,S_2=\cS(s_1;s_2),\,\,S_3=\cS(s_1,s_2;s_3),\, \ldots ,  \\[1mm]
S_n=\cS(s_1,\ldots,s_{n-1};s_n)\, \, \, \mathrm{and}\,  \,\,
S_{n+1}=\cS(s_1,\ldots,s_{n};1)
\end{array}
$$
are \comz.\\
More \gnltz, the \mos
$$
\preskip.2em \postskip.2em
\begin{array}{l}
V_1=\cS(I;U,s_1),\,\, V_2=\cS(I,s_1;U,s_2),\,\, V_3=\cS(I,s_1,s_2;U,s_3),\,\,
\ldots ,\\[1mm]
 V_n=\cS(I,s_1,\ldots,s_{n-1};U,s_n)\, \,\, \mathrm{and}\,  \, \,
V_{n+1}=\cS(I,s_1,\ldots,s_{n};U)
\end{array}
$$
cover the \mo $S=\cS(I,U)$.

\item 
If $S$, $S_1$, $\ldots$, $S_n\subseteq \gA$ are \moco and if  $a\in\gA$, then the \mos
$\cS(I;U,a),$ $\cS(I,a;U)$, $S_1$, $\dots$, $S_n$ are \comz.
\end{enumerate}
\end{lemma}
\begin{proof}
Items \emph{2} and \emph{3} result \imdt from Lemmas~\ref{lemAssoc} and~\ref{lemRecouvre}. 
\\
\emph{1.} The first case results from the fact that for $k_1$, \ldots, $k_n \ge 1$, we have,
for large enough $k$, $\gen {s_1, \ldots, s_n}^k \subseteq \geN {s_1^{k_1},
\ldots, s_n^{k_n}}$ (e.g. $k = \sum_i (k_i-1) + 1$).
\\
For the \gnl case, let $t_1$, \ldots, $t_n$ with $t_i \in \cS(I;U, s_i)$; we want to show that $\gen {t_1, \ldots, t_n}$ meets $S = \cS(I,U)$. By \dfnz, there is some $u_i \in \cM(U)$ and $k_i \ge 0$ such that $t_i \in u_is_i^{k_i} + \gen {I}$; by letting $u = u_1 \cdots u_n \in \cM(u)$, we obtain $us_i^{k_i} \in \gen {t_i} + \gen {I} \subseteq \gen {t_1, \ldots, t_n} + \gen {I}$. Therefore for large enough $k$,

\snic {
u\gen {s_1, \ldots, s_n}^k \subseteq u\geN {s_1^{k_1}, \ldots, s_n^{k_n}}
\subseteq \gen {t_1, \ldots, t_n} + \gen {I}.
}

But as $s_1$, $\ldots$, $s_n$ are \eco in $\gA_S$, there is some~\hbox{$s \in S$} such that $s \in \gen {s_1, \ldots, s_n}$; therefore $us^k \in \gen {t_1, \ldots, t_n} + \gen {I}$, \cad  $\gen {t_1, \ldots, t_n}$ meets $us^k + \gen {I} \subseteq S$.
\end{proof}

\section{A few \plgcsz}
\label{subsec loc glob conc}
\vspace{4pt}
\subsec{Linear \syssz}
\label{subsecLGSLI}

The following \plgc is a slight \gnn of the \plgrf{plcc.basic}
(basic \plgcz), which only concerned item \emph{4} below in the case of free modules of finite rank.\iplg 
 Actually the essential result has already been given in the \plgrf{plcc.basic.modules} (\plgc for modules). We give 
 the \dems 
again 
to emphasize their great simplicity.

Let $M_1$, $\dots$, $M_\ell$, $P$ be \Amosz. We say that a map $\Phi:M_1\times \cdots \times M_\ell\to P$ is
\ixc{homogeneous}{map}
if there exist integers $r_1$, $\ldots$, $r_\ell$ such that we identically have \hbox{$\Phi(a_1x_1,\ldots ,a_\ell x_\ell)=a_1^{r_1}\cdots a_\ell^{r_\ell}\Phi(x_1,\ldots
,x_\ell)$}. In such a case, the map $\Phi$ \gui{passes to the \lonsz}:
it can be naturally extended to a map

\snic{\Phi_S:S^{-1}M_1\times
\cdots \times S^{-1}M_\ell\to S^{-1}P}

for any \mo $S$.
The prototype of a \hmg map is a map given by \hmg \pols in the \coos when the modules are free of finite rank.

\begin{plcc}
\label{plcc.sli}
Let $S_1$, $\dots$, $S_n$ be \moco of $\gA$, $M$, $N$, $P$ be \Amosz, $\varphi$, $\psi$
be \alis  from~$M$ to~$N$, $\theta:N\to P$ be \aliz, and $x$, $y$ be \elts of~$N$. We write
$\gA_i$ for $\gA_{S_i}$, $M_i$ for $M_{S_i}$, etc.
Then we have the following \eqvcsz.
\begin{enumerate}
\item  Concrete patching of the \egts
$$\preskip.3em \postskip.4em 
x=y\quad   {\rm in} \quad   N  \quad  \Longleftrightarrow
\quad
\forall i\in\lrbn\;\; x/1= y/1  \quad   {\rm in} \quad   N_{i}. 
$$

\item  Concrete patching of the \egts of \alis
\[\preskip.3em \postskip.4em 
\begin{array}{ccc} 
 \varphi =\psi \quad   {\rm in} \quad   \Lin_\gA(M,N)  \qquad
\Longleftrightarrow    \\[.3em] 
\quad
\forall i\in \lrbn\;\; \varphi/1= \psi/1  \quad   {\rm in} \quad
\Lin_{\gA_{i}}(M_{i},N_{i}). \end{array}
\]

\item  Concrete patching of the \ndzs \elts
\[\preskip.3em \postskip.4em 
\begin{array}{ccc} 
x \hbox{ is \ndz in }  N \qquad \Longleftrightarrow   \\[.3em] 
\forall i\in \lrbn
\;\;x/1\; \hbox{ is \ndz in } N_{i}.
\end{array}
\]

\vspace{-.7em}
\pagebreak	

\item  Concrete patching of the solutions of \slis
$$\preskip.4em \postskip.4em 
x\in\Im \varphi  \quad  \Longleftrightarrow \quad
\forall i\in \lrbn\;\; x/1\in\Im \varphi_{i}. 
$$

\item  Concrete patching of the solutions of \slis under \hmg conditions. 
 Let  $(\Phi_\ell)$ be a finite family of \hmgs maps
$$\preskip.2em \postskip.2em \mathrigid 2mu
\Phi_\ell:\Lin_\gA(M,N)\times N\to Q_\ell,\hbox{ or }\Phi_\ell:\Lin_\gA(M,N)\to
Q_\ell ,\hbox{
or  }\Phi_\ell: N\to Q_\ell. 
$$
Then
\[\preskip-.2em \postskip.4em 
\begin{array}{ccc} 
\left(\big(\&_\ell \; \Phi_\ell(\varphi ,y)=0\big)\;\Rightarrow \;
 y\in\Im \varphi\right) \quad \quad \Longleftrightarrow     
                        \\[.3em] 
\forall i\in \lrbn\;
\left(\big(\&_\ell \; \Phi_\ell(\varphi ,y)
=_{Q_{\ell,i}}
0\big)\;\Rightarrow \;
 y/1\in\Im \varphi_{i} \right).
 \end{array}
\]
where we have written $Q_{\ell,i}$ for $(Q_\ell)_{S_i}$.

\item   \label{plcc.sex} Concrete patching of the exact sequences.
The sequence

\snic{
M\vvers{\varphi}N\vvers{\theta}P}

is exact \ssi the sequences

\snic{
M_{i}\vvers{\varphi_{S_i}}N_{i}\vvers{\theta_{S_i}}P_{i}}

are exact for $i\in \lrbn$.

\item  Concrete patching of direct summands in the \mpfsz.
Here $M$ is a \smtf of a \mpfz~$N$.
\[\preskip-.4em \postskip.4em 
\begin{array}{ccc} 
M \hbox{ is a direct summand in  } N \Longleftrightarrow \\[.2em] 
\forall i\in \lrbn\;M_i  \hbox{ is a direct summand in } N_i. 
 \end{array}
\]
\end{enumerate}
\end{plcc}
\begin{proof}
The conditions are \ncrs because of Fact~\ref{fact.sexloc}.
A direct verification is \imde anyway.
Let us prove that the local conditions are sufficient.

\emph{1.} Suppose that $x/1=0$ in each $N_{i}$. 
For suitable $s_i\in S_i$ we~therefore have~$s_ix=0$ in $N$.
As $\sum_{i=1}^{n} a_i s_i =1$, we obtain $x=0$ in~$N$.

\emph{2.} Immediate consequence of \emph{1.}

\emph{3.} Suppose that $x/1$ is \ndz in each $N_{i}$. 
Let $a\in \gA$ with $ax=0$ in $\gA$, therefore also $ax/1=0$ in each~$N_{i}$. We therefore have~$a/1=0$ in each $\gA_{i}$, so also in~$\gA$.

\emph{4.}
Suppose that the equation $\varphi (z)=x$ admits a solution $z_i$  in each~$M_{i}$. We can write $z_i=y_i/s_i$ with $y_i\in M$ and $s_i\in S_i$. We therefore have~$u_i\varphi(y_i)=s_iu_ix$ in $N$ with $u_i\in S_i$. 
As $\sum_{i=1}^{n} a_i s_iu_i =1$, let~$z=\sum_{i=1}^{n} a_i u_iy_i$. We  obtain $\varphi(z)=x$ in~$N$.

\emph{5.} This is a simple variant of \emph{4}. The homogeneity of the $\Phi_\ell$'s intervenes so that the local \prt is well-defined, and so that it results from the global \prtz.

\emph{6.} This is a special case of the previous item.

\emph{7.} Let $\rho:N\rightarrow N/M$ be the canonical \prnz. The module $N/M$ is \egmt a \mpfz. The module $M$ is a direct summand  in~$N$ \ssi $\rho$ is right-\ivz. We can therefore conclude by the \plgrf{plcc.scinde}.
 \end{proof}
\rdb

\rem  \label{remplgcsli}
We can see that item \emph{5}, a simple variant of item \emph{4}, implies all the others as special cases. Moreover, item \emph{4} results from item \emph{1} with $y=0$ by considering the module
$\big(N/\varphi(M)\big)_{S_i}\simeq N_{S_i}/\varphi_{S_i}(M_{S_i})$.
We could therefore have stated item \emph{1} as the only basic principle and, from it, deduce items \emph{2} to \emph{6} as corollaries.
Finally, item \emph{7} also directly results from item \emph{4}
(see the \dem of the \plgz~\ref{plcc.scinde}).
\eoe

\subsec{Finiteness \prts for  modules}
\label{subsecLGFiMo}

The usual finiteness \prts of modules have a local \crcz.
Most have already been proven, we summarize.
\begin{plcc}
\label{plcc.ptf}
{\em  (Concrete patching of finiteness \prts for  modules)}
Let $S_1$, $\dots$, $S_n$  be \moco of $\gA$  and  $M$ be an \Amoz.   Then we have the following \eqvcsz.
\begin{enumerate}\itemsep0pt\parsep0pt
\item  $M$ is  \tf \ssi each of the $M_{S_i}$'s is an $\gA_{S_i}$-\mtfz.
\item $M$ is \pf \ssi each of the $M_{S_i}$'s is an~$\gA_{S_i}$-\mpfz.
\item  $M$ is flat \ssi each of the $M_{S_i}$'s is an~$\gA_{S_i}$-\mplz.
\item \label{iptfplcc.ptf} $M$ is \ptf \ssi each of the $M_{S_i}$'s is an~$\gA_{S_i}$-\mptfz.
\item $M$ is projective of rank $k$ \ssi each of the $M_{S_i}$'s is a \pro $\gA_{S_i}$-module of rank~$k$.
\item $M$ is \coh \ssi each of the $M_{S_i}$'s is an $\gA_{S_i}$-\comoz.
\item  $M$ is \noe \ssi each of the $M_{S_i}$'s is a \noe $\gA_{S_i}$-module.
\end{enumerate}
\end{plcc}
\begin{proof}
\emph{1.} See the \plgrf{plcc.tf}.

 \emph{2.} See the \plgrf{plcc.pf}.

 \emph{3.} See the \plgrf{plcc.plat}.

 \emph{4.} See the \plgrf{plcc.cor.pf.ptf}. We can also use the fact that a \mpf is \pro \ssi it is flat (and apply items \emph{2} and \emph{3}).

 \emph{5.} Results from item \emph{4} and from the fact that the \pol rank can be locally computed (it is equal to $X^k$ \ssi it is equal to $X^k$ after \lon at \mocoz).

 \emph{6.} See the \plgrf{plcc.coh}.

 \emph{7.} We exhibit the \dem for the \noet \cot defined \`a la Richman-Seidenberg. Let us limit ourselves to the case of two \come \lons at $S_1$ and $S_2$.
Consider a non-decreasing sequence $(M_k)_{k\in\NN}$ of \smtfs of $M$. 
 It admits an infinite subsequence $\big(M_{\sigma(k)}\big)_{k\in\NN}$, where~$\sigma(k)<\sigma(k+1)\, \forall \, k,$
with $M_{\sigma(k)}=M_{\sigma(k)+1}$ after \lon at $S_1$ for all~$k$.
Consider the infinite sequence $M_{\sigma(k)}$ seen in $M_{S_2}$. It admits two equal consecutive terms $M_{\sigma(k)}$ and~$M_{\sigma(k+1)}$.
So~$M_{\sigma(k)}$ and  $M_{\sigma(k)+1}$ are equal both in~$M_{S_1}$ and~$M_{S_2}$. Therefore they are equal in $M$.
\end{proof}

\subsec{\Prts of  commutative \risz}
\label{subsecLGA1}

We recall a few results already established regarding the local \crc of a few interesting \prts for  commutative \risz, in the sense of the \lon at \mocoz.
\begin{plcc}
\label{plcc.propaco}
{\em  (Concrete patching of \prts of commutative \risz)} 
Let $S_1$, $\dots$, $S_n$ be \moco and~$\fa$ be an \id of~$\gA$.
Then we have the following \eqvcsz.
\begin{enumerate}\itemsep0pt\parsep0pt
\item  $\gA$ is  \coh \ssi each  $\gA_{S_i}$ is  \cohz.
\item  $\gA$ is a \lsd \ssi each  $\gA_{S_i}$ is a \lsdz.
\item  $\gA$ is a \qiri \ssi each  $\gA_{S_i}$ is a \qiriz.
\item  $\gA$ is reduced \ssi each  $\gA_{S_i}$ is reduced.
\item  The \id $\fa$ is \lop \ssi each  $\fa_{S_i}$ is \lopz.
\item  $\gA$ is  \ari \ssi each  $\gA_{S_i}$ is  \ariz.
\item  $\gA$ is  a Pr\"ufer \ri \ssi each  $\gA_{S_i}$ is  a Pr\"ufer \riz.
\item  The \id $\fa$ is \icl \ssi each  $\fa_{S_i}$ is \iclz.
\item  $\gA$ is normal \ssi each  $\gA_{S_i}$ is normal.
\item  $\gA$ is of \ddk $\leq k$ \ssi  each  $\gA_{S_i}$ is of \ddkz~$\leq k$.
\item  $\gA$ is  \noe \ssi each  $\gA_{S_i}$ is  \noez.
\end{enumerate}
\end{plcc}

Moreover recall that for \lons at \ecoz, the \plgc also applies for the notions of a \adk and of a \fdi \noe \cori (\plgref{plcc.ddk}).

\penalty-2500
\subsec{Concrete \lgb principles for \algsz}
\vspace{3pt}
\subsubsection*{Localization at the source}

\begin{plcc}\label{plcc.etale}\label{plcc.apf}
Let $S_1$, $\dots$, $S_n$ be \moco of a \ri $\gk$ and $\gA$ be a \klgz. Then
\propeq
\begin{enumerate}
\item $\gA$ is \tf (resp.\ flat, \fptez, \pfz, finite, integral, \stfez, \spbz, \stez) over~$\gk$.
\item Each of the \algs $\gA_{S_i}$ is \tf (resp.\ flat, \fptez, \pfz, finite, integral, \stfez, \spbz, \stez) over~$\gk_{S_i}$.
\end{enumerate}
Similarly if $\gA$ is \stfe and if $\lambda\in\Asta$, 
then $\lambda$ is dualizing \ssi each of the forms $\lambda_{S_i}$ is dualizing. 
\end{plcc}
\begin{proof} 
\emph{1} $\Leftrightarrow$ \emph{2.} 
We introduce the \fpte \klg $\prod_i\gk_{S_i}$. It then suffices to apply \thref{propFidPlatPrAlg}.

The question of the dualizing form (when $\gA$ is \stfez) is a question of \iso of modules and stems from the \plgcs for  modules
(by taking into account Fact~\ref{factSpbEds}). 
\end{proof}

\paragraph{Localization at the sink}~

There are also the \plgs that correspond to \prts said to be \gui{local in~$\gA$.}
Here we need \lons at \eco (\moco are not sufficient).

\begin{plcc}\label{plcc2.apf}~\\
Let $\gA$ be a \klg and $s_1$, \ldots, $s_m$ be \eco of $\gA$. Then \propeq
\begin{enumerate}
\item $\gA$ is \tf (resp.\ \pfz,  flat) over~$\gk$.
\item Each of the \algs $\gA_{s_i}$ is \tf (resp.\ \pfz,  flat) over~$\gk$.
\end{enumerate}
\end{plcc}
\begin{proof}
First of all if $\gA=\kxn=\kXn\sur\fa$ and $s=S(\ux)$ (where~\hbox{$S\in\kuX$}), then~$\gA_s=\gk[x_1,\ldots,x_n,t]$ with $t=1/s$ in $\gA_s$, which also gives 
$$\preskip-.40em \postskip.3em 
\gA_s=\gk[X_1,\ldots,X_n,T]\sur{(\fa+\gen{TS(\uX)-1})}. 
$$
Thus the \prt of being \tf or \pf is stable by \lon at an \elt (but it is not stable for a \lon at an arbitrary \moz).
\\
Regarding the flatness, as $\gA_s$ is flat over $\gA$, if $\gA$ is flat over~$\gk$,~$\gA_s$ is flat over $\gk$ (Fact~\ref{factAlgPlate}).

Now suppose that $\sum_is_iu_i=1$ in $\gA$.
\\
First of all let us see what we obtain if each of the \klgs $\gA_{s_i}$ is \tfz. We can suppose that the \gtrs are derived from \elts of~$\gA$ (by considering the corresponding fraction of \denoz~$1$).
Let us make a single list $(\xn)$ with all these \elts of $\gA$.
The reader will then observe by a small computation that $\gA$ is generated by
$$\preskip.3em \postskip.3em 
{(\xn,s_1,\ldots,s_m,u_1,\ldots,u_m)=(y_1,\ldots,y_p),\hbox{  with  }p=n+2m.}
$$
Now let us consider the case where all the \algs $\gA_{s_i}$ are \pfz.
We consider some \idtrs $Y_i$ corresponding to the list $(y_1,\ldots,y_p)$ defined above. We write $s_i=S_i(\ux)$, $u_i=U_i(\ux)$  (\pols  in~$\gk[\ux]$).
\\
For the common \sgr $(\xn)$ that we have just considered, and for each $i\in\lrbm$, we have a corresponding \sypz, say $F_i$, in  $\gk[\uX,Y_{n+i},T_i]$, which allows us to define the \iso
$$\preskip.3em \postskip.3em 
{\gk[\uX,Y_{n+i},T_i]\sur{\fa_i}\to \gA_{s_i},}
$$
with $\fa_i=\gen{F_i,Y_{n+i}-S_i(\uX),Y_{n+i}T_i-1}$. For each $f\in F_i$ there is an exponent $k_f$ such that $s_i^{k_f}f(\ux)=0$ in $\gA$. 
We can take all the $k_f$'s equal, say, to~$k$.
\\
We then consider the following \syp in $\gk[Y_1,\ldots,Y_p]$, %
with~$Y_j=X_j$ for $j\in\lrbn$. First of all we take all the $Y_{n+i}^{k}f(\uX)$'s for $f\in F_i$ and $i\in\lrbm$.\\
Next we write the relations $Y_{n+i}-S_i(\uX)$'s and $Y_{n+m+i}-U_i(\uX)$'s for the indices $i\in\lrbm$. Finally, we take the relation that corresponds to~$\sum_iu_is_i=1$, \cad $\sum_{i=1}^mY_{n+i}Y_{n+m+i}-1$.
\\
 The readers will do the computation to convince themselves that we indeed have a faultless description of the \klgz~$\gA$. The contrary would have been surprising, even immoral, since we have transcribed all that we could have known about the situation. The key was that this could have been expressed by a finite \sys of relations over a finite \sys of \idtrsz. Actually we proceeded exactly as in the \dem of the \plgrf{plcc.pf} for the \mpfsz.

Regarding the flatness, consider $(\an)$ in $\gk$ and $(\xn)$ in $\gA$ such that $\sum_ix_ia_i=0$. We want to show that $(\xn)$ is an $\gA$-\lin combination of \rdls  in $\gk$. We know that this is true after \lon at each of the $s_k$'s. We therefore have an exponent $N$ such that for each $k$ we have an \egt 
$$\preskip.4em \postskip.4em\ndsp
s_k^N(\xn)=\som_{j=1}^{p_j}b_{k,j}(x_{1,k,j},\dots,x_{n,k,j}),
$$
 ($x_{i,k,j}\in\gk$, $b_{k,j}\in\gA$)
with $\sum_ix_{i,k,j}a_i=0$.
We finish by taking an $\gA$-\lin combination of the $s_k^N$'s equal to $1$.    
\end{proof}

\section{A few abstract \plgsz}
\label{subsec loc glob abs}
An essential tool in classical \alg is the \lon at (the complement of) a prime \idz. This tool is a priori difficult to use \cot because we do not know how to construct the \ideps which intervene in the classical \demsz, and whose existence relies on the axiom of choice. However, we 
 observe that those prime ideals are generally used in proofs by contradiction, and this gives an explanation of the fact that the use of these \gui{ideal} objects can be avoided and even interpreted \cot (see Section~\ref{secMachLoGlo}).


The \plga in commutative \alg is an informal principle according to which certain \prts regarding  modules over  commutative \ris are true if and only if they are true after \lon at any \idepz.

We now recall a few cases where the \plga applies in \clamaz, by explaining the link with the corresponding concrete principles.

An abstract version of the \plgc \ref{plcc.sli} is the following.
\begin{plca}
\label{plca.ring}
Let  $\varphi$, $\psi$ be \alis $M\to N$, $\theta$ be a \ali $N\to P$, and $x$, $y$ be \elts of~$N$.
Then we have the following \eqvcsz.
\begin{enumerate}
\item  Abstract patching of the \egts
$$\preskip.3em \postskip.4em 
 x=y\quad   {\rm in} \quad   N  \quad  \Longleftrightarrow
\quad
\forall \fp\in\Spec\gA\;\; x/1= y/1  \quad   {\rm in} \quad   N_{\fp}. 
$$

\item  Abstract patching of the \egts of \alis
\[\preskip.3em \postskip.4em 
\begin{array}{ccc} 
 \varphi =\psi \quad   {\rm in} \quad   \Lin_\gA(M,N)  \qquad
\Longleftrightarrow    \\[.3em] 
\quad
\forall \fp\in\Spec\gA\; \; \varphi/1= \psi/1  \quad   {\rm in} \quad
\Lin_{\Ap}(M_{\fp},N_{\fp}). \end{array}
\]

\item  Abstract patching of the \ndz \elts
\[\preskip.3em \postskip.4em 
\begin{array}{ccc} 
x \hbox{ is \ndz in }  N \qquad \Longleftrightarrow   \\[.3em] 
\forall \fp\in\Spec\gA\;
\;\;x/1\;$ is \ndz in $N_{\fp}.
\end{array}
\]

\item  Abstract patching of the solutions of \slis
$$\preskip.4em \postskip.4em 
x\in\Im \varphi  \quad  \Longleftrightarrow \quad
\forall \fp\in\Spec\gA\; \; x/1\in\Im \varphi_{\fp}. 
$$

\item  Abstract patching of the solutions of \slis under \hmg conditions. Let  $(\Phi_\ell)$ be a finite family of \hmg maps
$$\preskip.2em \postskip.2em \mathrigid 2mu
\Phi_\ell:\Lin_\gA(M,N)\times N\to Q_\ell,\hbox{ or }\Phi_\ell:\Lin_\gA(M,N)\to
Q_\ell ,\hbox{
or  }\Phi_\ell: N\to Q_\ell. 
$$
Then
\[\preskip-.4em \postskip.4em 
\begin{array}{ccc} 
\left(\big(\&_\ell \; \Phi_\ell(\varphi ,y)=0\big)\;\Rightarrow \;
 y\in\Im \varphi\right) \quad \quad \Longleftrightarrow \\[.3em] 
\forall \fp\in\Spec\gA\;
\left(\big(\&_\ell \; \Phi_\ell(\varphi ,y)
=_{Q_{\ell,\fp}}
0\big)\;\Rightarrow \;
 y/1\in\Im \varphi_{\fp} \right), \end{array}
\]
where we have written $Q_{\ell,\fp}$ for $(Q_\ell)_{\fp}$.

\item Abstract patching of the exact sequences.
The sequence

\snic{
M\vers{\varphi}N\vers{\theta}P}

is exact \ssi the sequence

\snic{
M_{\fp}\vers{\varphi_{\fp}}N_{\fp}\vers{\theta_{\fp}}P_{\fp}}

is exact for every $\fp\in\Spec\gA\; $.

\item  Abstract patching of direct summands in  \mpfsz.
Here $M$ is a \smtf of a \mpfz~$N$.
\[\preskip.2em \postskip.4em 
\begin{array}{ccc} 
M \hbox{ is a direct summand in  } N \Longleftrightarrow \\[.2em] 
\forall \fp\in\Spec\gA\;M_\fp  \hbox{ is a direct summand in } N_\fp. 
 \end{array}
\]
\end{enumerate}

\end{plca}

\vspace{-.7em}
\pagebreak	
\begin{Proof} {Proofs (nonconstructive). }
The conditions are \ncrs because of Fact \ref{fact.sexloc}.
A direct verification is actually \imdez.
For the converses, we assume \spdg that the \ri $\gA$ is nontrivial. It suffices to treat item \emph{4}
(see the remark on \paref{remplgcsli}).
Actually we have already established item \emph{6}, which implies item \emph{4}, in the \plga \rref{plca.basic.modules}, but we think that it is  usefull to give two distinct classical \dems (the second is the one given in Chapter~\ref{chapSli})
and to compare their degree of effectivity.

 {\em First \demz}.
\\
Suppose  $x\notin\Im\varphi$, it amounts to the same as saying that $x\neq 0$ in $N/\varphi(M)$. Since for a \idep $\fp$ we have $\big(N/\varphi(M)\big)_\fp\simeq N_\fp/\varphi_\fp(M_\fp)$, it suffices to prove item \emph{1} with $y=0$. We reason by contradiction by assuming \hbox{$x\neq 0$} in~$N$. In other words $\Ann_\gA(x)\neq \gen{1}$, and there exists a $\fp\in\Spec \gA$ which contains~$\Ann_\gA(x)$.
Then, 
since~$\big(\Ann_\gA(x)\big)_\fp=\Ann_{\gA_\fp}(x/1)$, we obtain~$x\neq_{N_\fp}0$.

 {\em Second \demz}.
\\
The \prt $x\in\Im\varphi$ is of \carfz. We can therefore apply Fact\etoz~\ref{factPropCarFin} which says (in \clamaz) that for a \carf \prtz, the \plgc (\lon at \mocoz) is \eqv to the \plga (\lon at all the \idemasz).
\end{Proof}

\comms
\label{comment plga1}~
\\
 1) It seems impossible that the second proof, which is too general, can ever be made into a constructive proof. 
 The first \dem is not \gui{\gnltz} \cov either, but there exist some cases where it is. For this it suffices to satisfy the following conditions, in the case of item \emph{4.}
\begin{enumerate}
\item [--] The module $N$ is \pf and the module $M$ is \tfz.
\item [--] The \ri $\gA$ is \coh and \fdiz.
\item [--]
For every strict \itf $\fa$ of $\gA$ we know how to construct a \idep $\fp$ containing~$\fa$.
\end{enumerate}
The last two conditions are satisfied when $\gA$ is a \pf \alg over $\ZZ$ or over a \gui{fully factorial} field (see~\cite{MRR}).

 2)
This allows us, for example, to give another \prco of the explicit matrix form \tho (\thref{th matproj}).
As mentioned on \paref{subsec cas generique}, it suffices to treat the \gnq case and to show certain \egtsz~\hbox{$r_ir_j=0$} and $r_hu=0$. As the \ri $\Gn$ is a \pf
\alg over $\ZZ$, we can show these \egts by applying the \rca of the \egtsz. We are therefore brought back to the case of a local  \ri 
obtained as a localization of $\Gn$, 
and in this case the \egts are true since the module is free by applying the local freeness lemma.

 3) In practice, we can understand the \plga \ref{plca.ring} in the following intuitive form: to prove a \tho of commutative \alg whose meaning is that a certain \sli over a commutative \riz~$\gA$ admits a solution, it suffices to treat the case where the \ri is local. It is a principle of the same type as the Lefschetz principle: to prove a \tho of commutative \alg whose meaning is that a certain \ida takes place, it suffices to treat the case where the \ri is the complex number field (or any sub\ri that suits us best, in fact).
This remark is developed in Section~\ref{secMachLoGlo}.

 4)
In the article \cite{Bass}, Hyman Bass makes the following comment regarding a \noe
version of the \plga \ref{plca.ring}, item~\emph{7.}\\
{\em The latter result, elementary as it is, seems to defy any proof which does not either use, or essentially reconstruct, the functor 
 $\mathrm{Ext}^1.$}\\
This comment is surprising, in view of the perfectly trivial \crc of our \dem of the corresponding concrete principle, which computes nothing that resembles an $\mathrm{Ext}^1.$
Actually, when the goal is to show that a short exact sequence splits,
it seems that the efficient
 computational machinery of the $\mathrm{Ext}$'s is often useless, and that it can be short-circuited by a more \elr argument.

 5) The \plga above also works by uniquely using the \lon at any \idemaz, as seen in the \plga \rref{plca.basic.modules}. But this is not really useful because the \lons at the \idemas are 
the least extensive 
(among the \lons at the \idepsz).
However, there are cases where the classical reasoning uniquely uses \lons at \idemisz. They are more subtle \dems that are more difficult to decrypt \cotz. 
We will elaborate on this in Section~\ref{subsecLGIdepMin}.

 6) 
As mentioned on \paref{plcc.tf}, the \plga for  \mtfs does not work: just because a module is \tf after \lon at every \idep does not mean it is \ncrt \tfz.
The same would hold for the \prcc of  \mpfs or for that of  \comosz.
This denotes a certain superiority of the \plgcs over the \plgasz.
\eoe

\newpage
\section{Concrete patching of objects}
\label{secColleCiseaux}

\vspace{4pt}
\subsec{Glue and  scissors} \label{subsecColleCiseaux}

Here we give a brief discussion regarding  patching methods in \dile \gmt and their translations in commutative \algz.

\smallskip   
First of all we examine the possibility of constructing a smooth manifold from local charts, 
\cad by a patching of open sets $U_i$ of $\RR^n$ by means of diffeomorphisms (or \isosz)
$\varphi_{ij}:U_{ij}\to U_{ji}$:   $U_{ij}$ is an open set of $U_i$ and $\varphi_{ji}=\varphi_{ij}^{-1}$.

{~

\centerline{\includegraphics[width=11cm]{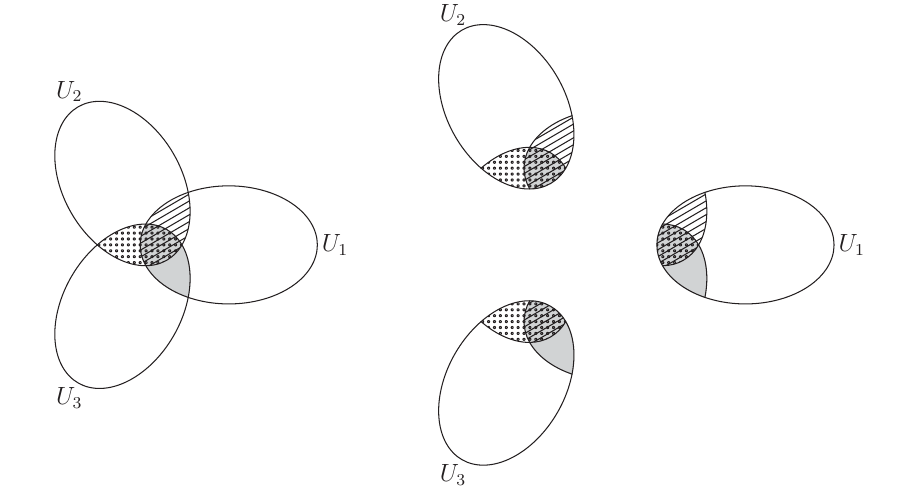}}
}

We will consider the simple case where the \vrt is obtained by only patching a finite number of open sets of $\RR^n$.

In this case the condition to fulfil is that the morphisms of \rcms must be \emph{compatible between them three by three}.
This \prmt means the following.
For each triple of distinct indices $(i,j,k)$ we consider the open set $U_{ijk}= U_{ij}\cap U_{ik}$ (therefore with $U_{ijk}=U_{ikj}$). The \cpbt means on the one hand that, for each $(i,j,k)$, the restriction $\Frt{\varphi_{ij}} {U_{ijk}}$ establishes an \iso from $U_{ijk}$ to $U_{jik}$, and on the other hand that if we compose the \isos 

\snic{U_{ijk}\vvvvers{ \Frt{\varphi_{ij}} {U_{ijk}} }U_{jik}\quad$ and  $\quad
U_{jki}\vvvvers{ \Frt{\varphi_{jk}} {U_{jki}} }U_{kji} }

we obtain the \iso 
$U_{ijk}\vvvvers{ \Frt{\varphi_{ik}} {U_{ijk}} }U_{kij}$: 
$\Frt{\varphi_{ik}}{\bullet}=\Frt{\varphi_{jk}}{\bullet}\circ \Frt{\varphi_{ij}}{\bullet}\,$.

\smallskip 
If we try to do the same thing in commutative \algz, we will consider some \ris $\gA_i$ (corresponding to the \ris $\Cin(U_i)$) and some \eltsz~\hbox{$f_{ij}\in\gA_i$}.
The \ri $\Cin(U_{ij})$ would correspond to $\gA_i[1/f_{ij}]$ and the patching morphism $\varphi_{ij}$ to an \iso $\omega_{ij}:\gA_i[1/f_{ij}]\to\gA_j[1/f_{ji}]$.
We will \egmt have to formulate some three-by-three compatibility conditions.
We then hope to construct a \ri $\gA$ and some \elts $f_i\in\gA$, such that $\gA_i$ could be identified with $\gA[1/f_i]$, $f_{ij}$ with \gui{$f_j$ seen in $\gA[1/f_i]$,} and~$\omega_{ij}$ with the identity between $\gA[1/f_i][1/f_j]$ and  $\gA[1/f_j][1/f_i]$.

Unfortunately, this does not always work well. The \ri $\gA$ that is supposed to patch the $\gA_i$'s does not always exist (however, if it exists it is well-determined, up to unique \isoz).

The first example of this obvious failure of the patching is in \pro space. The complex \pro space $\Pn(\CC)$ is obtained by patching affine charts $\CC^n$, but the corresponding \ris of functions, \isoc to~$\CC[\Xn]$, do not patch
together: 
there are no \pol functions defined over $\Pn(\CC)$, besides the constants,
and by localizing the \ri $\CC$ there is no chance of obtaining the \ri 
$\CC[\Xn]$.

This illustrates the fact that  \agq \gmt is much more rigid than 
          $\Cin$ \gmtz.

This unpleasant phenomenon is at the origin of the creation of Grothendieck's schemes, which are the  abstract objects \fmt obtained by patching \ris along 
patching morphisms 
when the three-by-three compatibility conditions 
are satisfied, but whose patching no \ri wants to perform.

\smallskip  
Let us now consider the question of the patching of 
        vector bundles when they are locally defined 
over a fixed smooth \vrt $U$, covered by a finite number of open sets $U_i$. Let $U_{ij}=U_i\cap U_j$. 
The vector bundle $\pi:W\to U$ 
that we want to construct, whose every fiber is \isoc to a given \evcz~$F$, is known a priori only by its restrictions~\hbox{$\pi_i:W_i\to U_i$}. 
In order to patch, we need
 patching diffeomorphisms~\hbox{$\psi_{ij}: W_{ij}\to W_{ji}$}

 where $W_{ij}=\pi_i^{-1}(U_{ij})$. 
These morphisms must first of all respect the structure of the \evc fiber by fiber. In addition, again, we need three-by-three compatibility conditions,
analogous to those which we have defined in the first case.

Now if we pass to the analogous case in commutative \algz, we must start from a \ri $\gA$ with a \sys of \eco $(f_1,\ldots,f_\ell)$. Let $\gA_i=\gA[1/f_i]$ and $\gA_{ij}=\gA[1/f_if_j]$.
For each index $i$, we give the \gui{module of the sections of the fiber $\pi_i:W_i\to U_i$,} \cad an \hbox{$\gA_i$-module $M_i$}. The  $\psi_{ij}$'s are now represented by \isos of $\gA_{ij}$-modules 

\snic{\gA_{ij}\otimes_{\gA_i}M_i\vvers{\theta_{ij}} \gA_{ji}\otimes_{\gA_j}M_j
\simarrow M_{ij}=M_{ji}.}

We will see in the following subsections that this time everything goes well: 
if the three-by-three compatibility conditions 
are satisfied, 
we indeed have an \Amo $M$ that \gui{patches} the $\gA_i$-modules $M_i$.

\subsec{A simple case}\label{subsecLGFAC}

\begin{theorem}\label{thRecolSousMod}
Let $\gA$ be an integral \ri with quotient field $\gK$, 
$N$ be a \torf \Amoz, $S_1$, \ldots, $S_n$ be \moco of $\gA$ and for each $i\in\lrbn$ let $M_i$ be $\gA_{S_i}\!$-submodule of $S_i^{-1}N\subseteq \gK\te_\gA N$. Suppose that for each $i,j\in\lrbn$ we have $S_j^{-1}M_i=S_i^{-1}M_j$ (seen as \Asubs \hbox{of $\gK\te_\gA N$}). Then we have the following results. 
\begin{enumerate}
\item There exists a unique \Asub $M$ of $N$ such that we have $S_i^{-1}M=M_i$ for each $i\in\lrbn$.
\item This submodule $M$ is equal to the intersection of the $M_i$'s.
\item If the $M_i$'s are \tf (resp.\,\pfz, \cohsz, \ptfsz), the same goes for $M$.
\end{enumerate}
\end{theorem}

\begin{proof}
\emph{1} and \emph{2.} Let $P=\bigcap_iM_i$. First of all $P\subseteq N$ because one \elt of the intersection is of the form

\snic{\frac{x_1}{s_1}=\cdots=\frac{x_n}{s_n}=\frac{\som_ia_ix_i}{\som_ia_is_i}=\som_ia_ix_i\quad\mathrm{if}\;\;\som_ia_is_i=1\;\;\hbox{in}\;\gA}

(with  $x_i\in N$, $s_i\in S_i$ for $i\in\lrbn$).
\\
Let us show that
the module $P$ satisfies the required conditions.
 \\ 
First of all $P \subseteq M_i$ so $S_i^{-1} P \subseteq M_i$ for each $i$. 
Conversely, let $x_1\in M_1$ for example, we want to see that $x_1$ is in $S_1^{-1}P$.\\
Since $S_j^{-1}M_1=S_1^{-1}M_j$, there exists a $u_{1,j}\in S_1$ such that $u_{1,j}x_1\in M_j$. By letting $s_1=\prod_{j\neq1}u_{1,j}$, we indeed obtain $s_1x_1\in\bigcap_iM_i$.
\\
Now let us prove the uniqueness. \\
 Let $Q$ be a module satisfying the required conditions. We have $Q \subseteq S_i^{-1}Q=M_i$ and thus $Q\subseteq P$. Then consider the sequence $Q \to P\to 0$.
Since it is exact after \lon at \mocoz, it is exact (\plgref{plcc.basic.modules}), \cad the inclusion \homo is surjective, so $Q=P$.

Finally, item \emph{3} results from already established \plgcsz. 
\end{proof}

If we do not assume that the \ri is integral and the module is torsion-free, the previous \tho is a little more delicate. This will be the object of the \plgref{plcc.modules 2}. 

\subsec{Patching of objects in  modules}


Let $\gA$ be a commutative \riz, $(S_i)_{i\in\lrbn}$ be \moco of~$\gA$.
Let $ \gA_i:= \gA_{S_i}$ and  $ \gA_{ij}:= \gA_{S_iS_j}$ ($i\neq j$) such that $ \gA_{ij}= \gA_{ji}$. Let~\hbox{$\alpha_i: \gA\to  \gA_i$} and $\alpha_{ij}: \gA_i\to  \gA_{ij}$ be natural
\homosz. 

In the remainder, notations like $(M_{ij})_{i<j\in\lrbn}$ and $(\varphi_{ij})_{i\neq j\in\lrbn})$ mean that we have $M_{ij}=M_{ji}$ but (a priori) not $\varphi_{ij}=\varphi_{ji}$.

\pagebreak	

\begin{plcc}
\label{plcc.modules 1} {\em (Concrete patching of \elts in a module, and of \homos between modules)}  
\begin{enumerate}
\item \label{i1plcc.modules 1} Let $(x_i)_{i\in\lrbn}$ be an \elt of  $\prod_{i\in\lrbn}  \gA_i$.
So that there exists some~\hbox{$x\in  \gA$} satisfying $\alpha_i(x)=x_i$ in each $ \gA_i$, it is sufficient and necessary for each~\hbox{$i<j$} we have $\alpha_{ij}(x_i)=\alpha_{ji}(x_j)$ in $ \gA_{ij}$. In addition, this $x$ is then uniquely determined.
In other terms the \ri $ \gA$ (with the \homos $\alpha_{i}$) is the limit of the diagram
$$\preskip.2em  
\big(( \gA_i)_{i\in\lrbn},( \gA_{ij})_{i<j\in\lrbn};(\alpha_{ij})_{i\neq j\in\lrbn}\big) 
$$
$$\preskip-.2em \postskip.4em\hspace*{.1em}
\xymatrix @C=2em @R=1.5em
          {  &&& \gA_i \ar[r]^{\alpha_{ij}}\ar[ddr]^(.7){\alpha_{ik}}
                & \gA_{ij}\\
\gC\ar[urrr]^{\psi_i} \ar[drrr]_{\psi_j}\ar[ddrrr]_{\psi_k}
\ar@{-->}[rr]_(.6){\psi!} &&\gA \ar[ur]_(.4){\alpha_i}\ar[dr]^(.4){\alpha_j} \ar[ddr]_(.5){\alpha_k} &&\\
  &&& \gA_j \ar[uur]_(.7){\alpha_{ji}}\ar[dr]
          &\gA_{ik}\\
&&& \gA_k\ar[ur] \ar[r]_{\alpha_{kj}} 
   & \gA_{jk}\\
}
$$

\item \label{i2plcc.modules 1}
Let $M$ be an \Amoz.
Let $M_i:=M_{S_i}$ and  $M_{ij}:=M_{S_iS_j}$ ($i\neq j$) such that $M_{ij}=M_{ji}$. Let $\varphi_i:M\to M_i$ and $\varphi_{ij}:M_i\to M_{ij}$ be the natural \alisz.
Then the \Amo $M$ (with the \alis $\varphi_i:M\to M_i$) is the limit of the diagram
$$\preskip.3em \postskip.4em 
\big((M_i)_{i\in\lrbn},(M_{ij})_{i<j\in\lrbn};(\varphi_{ij})_{i\neq j\in\lrbn}\big). 
$$

\item \label{i3plcc.modules 1}
Let $N$ be another module, let $N_i:=N_{S_i}$, $N_{ij}:=N_{S_iS_j}$. 
 For each $i\in\lrbn$ let $\psi_i:M_i\rightarrow N_i$ be an $\gA_i$-\lin map.
So that there exists an \Ali $\psi: M\rightarrow N$ satisfying $\psi_{S_i}=\psi_i$ for each $i$, it is sufficient and necessary, for each $i<j$, for the two \alis $(S_j)^{-1}\psi_i$ and $(S_i)^{-1}\psi_j$  from $M_{ij}$ to $N_{ij}$ to be equal. In addition, the \ali $\psi$ is then uniquely determined.
$$\preskip-.2em \postskip.4em\hspace*{3em}
\xymatrix @C=2em @R=1em
          {
    & M_i \ar[ddr]^(.3){\varphi_{ij}}\ar[rrr]^{\psi_{i}}&&&   
                  N_i \ar[ddr]^{\phi_{ij}}&\\
M\phantom{_{ij}} \ar[ur]^{\varphi_i}\ar[ddr]_{\varphi_j}\ar@{-->}[rrr]^(.7){\psi!}
&&& ~\phantom{_i} N\phantom{_j}\ar[ur]^{\phi_i}\ar[ddr]^(.3){\phi_j} \\
&& M_{ij}\ar@{-->}[rrr] &&&N_{ij}&\\
                 & M_j \ar[ur]^{\varphi_{ji}}\ar[rrr]^{\psi_{j}}&&&   
                  N_j \ar[ur]^{\phi_{ji}} 
}
$$
In other terms the \Amo $\Lin_\gA(M,N)$ is the  limit of the diagram formed by the $\Lin_{\gA_i}(M_i,N_i)$'s, the $\Lin_{\gA_{ij}}(M_{ij},N_{ij})$'s and the natural \alisz.

\end{enumerate}
\end{plcc}

\begin{proof}
 \emph{\ref{i1plcc.modules 1}.} Special case of \emph{\ref{i2plcc.modules 1}.}
 
\emph{\ref{i2plcc.modules 1}.} Let $(x_i)_{i\in\lrbn}$ be an \elt of $\prod_{i\in\lrbn}  M_i$. We need to show that for some 
 $x\in M$ satisfying $\varphi_i(x)=x_i$ in each $ M_i$ to exist, it is sufficient and necessary that for each $i<j$ we have $\varphi_{ij}(x_i)=\varphi_{ji}(x_j)$ in $ M_{ij}$. In addition, this $x$ must be unique.
\\
 The condition is clearly \ncrz. Let us show that it is sufficient.
\\
Let us show the existence of $x$. 
There exist $s_i$'s in $S_i$ and $y_i$'s in $M$ 
 such that we have $x_i=y_i/s_i$ in each $M_i$.
If $\gA$ is integral, 
$M$ torsion-free and \hbox{each $s_i\neq 0$},
we have in the module obtained by \eds to the quotient field

\snic{\frac{y_1}{s_1}=\frac{y_2}{s_2}=\cdots
=\frac{y_n}{s_n}=\frac{\som_ia_iy_i}{\som_ia_is_i}= \som_ia_iy_i=x\in M,}

with $\sum_ia_is_i=1$. In the \gnl case we do just about the same thing.
For each pair $(i,j)$  with $i\neq j$, the fact that $x_i/1=x_j/1$ in $M_{ij}$ means that for certain $u_{ij}\in S_i$ and  $u_{ji}\in S_j$ we have $s_j u_{ij} u_{ji} y_i = s_i u_{ij} u_{ji} y_j $.
Let~\hbox{$u_i=\prod_{k\neq i}u_{ik}\in S_i$}. We have $s_j u_{i} u_{j} y_i = s_i u_{i} u_{j} y_j $.
Let $(a_i)$ be \elts of $\gA$ such that $\sum_i a_i s_i u_i =1$.
Let~\hbox{$x=\sum a_i  u_i y_i$}.
We need to show that $x/1=x_i$ in $M_i$ for each $i$. For example for $i=1$, we write the following \egts in $M$
$$
\preskip.0em \postskip.4em
\begin{array}{c}
s_1 u_1 x\,= \,s_1 u_1 \som_i a_i u_i y_i\, = \,
\som_i a_i s_1 u_1 u_i y_i    \\[1mm]
 =\, \som_i a_i s_i u_1 u_i y_1\,= \,
\big(\som_i a_i s_i u_i\big)u_1y_1 \,=\, u_1y_1.
\end{array}
$$
Thus $s_1 u_1 x=u_1y_1$ in $M$ and $x=y_1/s_1$ in $M_{S_1}$.
\\ 
Finally, the uniqueness of $x$  results from the \prcc of \egtsz.

\emph{\ref{i3plcc.modules 1}.} The composites of the \alis $M\to M_i\to N_i$ are compatible with the natural \alis $N_i\to N_{ij}$. We conclude with the fact that $N$ is the  limit of the diagram of the $N_i$'s and $N_{ij}$'s 
(item \emph{2}).
\end{proof}

{\bf A delicate point} (regarding item \emph{\ref{i3plcc.modules 1}}). If $M$ is a \pf \Amo or if $\gA$ is integral and $M$ \tfz, the natural $\gA_i$-\lin maps $\Lin_\gA(M,N)_{s_i}\to \Lin_{\gA_i}(M_i,N_i)$ are \isos (see Propositions \ref{fact.homom loc pf} and \ref{propPlateHom}).
 \\
 In the \gnl case, the notation $\psi_{s_i}$ is made ambiguous because it can either represent an \elt of $\Lin_{\gA_i}(M_i,N_i)$ or an \elt of $\Lin_\gA(M,N)_{s_i}$, and the natural \ali $\Lin_\gA(M,N)_{s_i}\to \Lin_{\gA_i}(M_i,N_i)$ is a priori only injective if $M$ is \tfz. This ambiguity can be a source or error. Especially as $\Lin_\gA(M,N)$ then appears as a limit of two essentially distinct diagrams: 
the one based on the $\Lin_{\gA_i}(M_i,N_i)$'s (the most interesting of the two) and the one based on the $\Lin_\gA(M,N)_{s_i}$'s.
\eoe

{\bf An example of a \rcm of \eltsz.}
  Given that the \deters of \endos of free modules are well-behaved under \lonz, given the \tho that affirms that the \mptfs are \lot free (in the strong sense) and given the previous \plgcz, we obtain the possibility of \emph{defining} the \deter of an \endo of a \mptf by only using \deters of \endos between free modules, after suitable \come \lonsz.
In other words the following fact can be established independently of the theory of \deters developed in Chapters~\ref{chap ptf0} and~\ref{chap ptf1}.

\sni
{\bf Fact.} \label{prop carc loc det} 
\emph{For an \endo $\varphi$ of a \ptf \Amoz~$M\!$, there exists a unique \elt $\det\varphi$ satisfying the following \prtz: if~$s\in \gA$ is such that the module $M_s$ is free, then $(\det\varphi)_s=\det(\varphi_s)$ in~$\gA_s$.}
\eoe

\subsec{Patching of modules}

The \rcm principle \ref{plcc.modules 2} that follows specifies which conditions are needed in order for the limit of an  analogous \sys of modules to fall within the framework indicated in the \plgrf{plcc.modules 1}.

\begin{definition}
\label{defAliloc}
Let $S$ be a \mo of $\gA$, $M$ be an \Amo and $N$ be an~\hbox{$\gA_S$-module}. An \Ali $\alpha
:M\to N$ is called a \emph{\lon morphism at $S$} if it is a morphism of \eds from $\gA$ to $\gA_S$ for $M$ (see \paref{pageChgtBase}).%
\index{morphism!localization --- (modules)}%
\index{localization!morphism (modules)}%
\end{definition}

In other words, if $\alpha:M\to N$ is a \lon morphism at $S$, and if~\hbox{$\beta_{M,S} :M\to M_S$} is the natural \aliz, the unique \Ali \hbox{$\varphi : M_S\to N$}
such that $\varphi \circ\beta_{M,S} =\alpha $ is an \isoz.
A \lon morphism at $S$ can be \care by the following conditions:
\begin{enumerate}
\item [--] $\forall x,x'\in M,\;\;\big(\alpha(x)=\alpha(x')\;\Longleftrightarrow\; \exists s\in
S,\;sx=sx'\big)$,
\item [--] $\forall y\in N,\;\exists x\in M,\;\exists s\in S,\;\;sy=\alpha(x)$.
\end{enumerate}

\begin{plcc}
\label{plcc.modules 2} {\em (Concrete patching of modules) }\\
 Let $S_1$, $\dots$, $S_n$ be \moco of $\gA$. 
 \\
Let $\gA_i=\gA_{S_i}$, $\gA_{ij}=\gA_{S_iS_j}$ and $\gA_{ijk}=\gA_{S_iS_jS_k}$.
We give in the category of \Amos a commutative diagram~$\fD$

\snac{\big((M_i)_{i\in  I}),(M_{ij})_{i<j\in  I},(M_{ijk})_{i<j<k\in  I};(\varphi_{ij})_{i\neq j},(\varphi_{ijk})_{i< j,i\neq k,j\neq k}\big)}

\snii
as in the following figure.

\smallskip {\small\hspace*{8em}{
$
\xymatrix @R=2em @C=7em{
 M _i\ar[d]_{\varphi _{ij}}\ar@/-0.75cm/[dr]^{\varphi _{ik}} &
     M _j\ar@/-1cm/[dl]^{\varphi _{ji}}\ar@/-1cm/[dr]_{\varphi _{jk}} &
        M _k\ar@/-0.75cm/[dl]_{\varphi _{ki}}\ar[d]^{\varphi _{kj}} &
\\
 M _{ij} \ar[rd]_{\varphi _{ijk}} & 
    M _{ik}  \ar[d]^{\varphi _{ikj}} & 
      M _{jk}  \ar[ld]^{\varphi _{jki}} 
\\
   &  M _{ijk} 
}
$
}}

\vspace{-1em}
Suppose that %
\begin{itemize}
\item For all $i$, $j$, $k$ (with $i<j<k$), $M_i$ is an $\gA_i$-module,  $M_{ij}$ is an~$\gA_{ij}$-module and $M_{ijk}$ is an~$\gA_{ijk}$-module.
Recall that according to our notation conventions we let $M_{ji}=M_{ij}$, $M_{ijk}=M_{ikj}=\dots$

\item For $i\neq j$,  $\varphi_{ij}:M_i\to M_{ij}$ is a \molo at $S_j$ (seen in $\gA_i$).
\item For $i\neq k$, $j\neq k$ and $i<j$, $\varphi_{ijk}:M_{ij}\to M_{ijk}$ is a  \molo at $S_k$ (seen in $\gA_{ij}$).
\end{itemize}
 Then, by letting $\big(M,(\varphi_i)_{i\in\lrbn}\big)$ be the limit of the diagram, each morphism $\varphi_i:M\to M_i$ is a \molo at~$S_i$.
In addition $\big(M,(\varphi_{i})_{i\in\lrbn}\big)$ is,
up to unique \isoz, the unique \sys that makes the  diagram commutative and that makes each $\varphi_i$ a \molo at $S_i$.
\end{plcc}

\begin{proof} 
The first item does not depend on the fact that the $S_i$'s are \comz. Indeed the construction of a limit of \Amos for an arbitrary diagram is stable by flat \eds (because this is the kernel of a \ali between two products). 
\\
However, if we take as a \eds the \moloz~\hbox{$\gA\to \gA_i$}, the diagram can be simplified as follows 
$$\preskip.4em \postskip-.20em 
\xymatrix @R=.5cm @C=.8cm{
 M _i\ar[d]_{\varphi _{ij}}\ar[dr]^{\varphi _{ik}} &               
\\
 M _{ij} \ar[rd]_{\varphi _{ijk}} & 
    M _{ik}  \ar[d]^{\varphi _{ikj}}  
\\
   &  M _{ijk} 
} 
$$
and it trivially admits the limit $M_i$.

To prove the uniqueness, we reason \spdg with a \sys of \eco $(s_1,\dots,s_n)$.
Let $\big(N,(\psi_i)\big)$ be a competitor. 
Since~$M$ is the limit of the diagram, there is a unique  \Ali  $\lambda:N\to M$ such that
$\psi_i=\varphi_i\circ \lambda$ for every $i$. Actually \hbox{we have $\lambda(v)=\big(\psi_1(v),\dots,\psi_n(v)\big)$}. Let us show that $\lambda$ is injective. If $\lambda(v)=0$ all the $\psi_i(v)$'s are null, and since $\psi_i$ is a  \molo at~$s_i$, there exist  exponents $m_i$ such that $s_i^{m_i}v=0$. 
Since the $s_i$'s are \comz, we have $v=0$. As $\lambda$ is injective we can assume $N\subseteq M$ and $\psi_i=\varphi_i\frt N$. 
Let us show that $N=M$. Let $x\in M$. As~$\psi_i$ and~$\varphi_i$ are two  \molos at $s_i$, there is an exponent $m_i$ such that $xs_i^{m_i}\in N$. Since the $s_i$'s are \comz, $x\in N$.
\end{proof}

\rem To understand why the comaximality condition is really \ncr for the uniqueness, let us examine the following \gui{overly simple} example.  
With the \ri $\ZZ$, and the unique \elt $s=2$, let us take for $M$ \hbox{a free $\ZZ[1/2]$-module} with basis~$(a)$ (where~$a$ is an arbitrary individual object). 
For clarity, let $M'$ be the \ZZmo $M$.\\
Also consider the free \ZZmo $N$ with basis $(a)$. Consider two \lon morphisms at $2^{\NN}$, $\varphi:M'\to M$ and $\psi:N\to M$. They both send $a$ to $a$. Thus $M'$ and $N$ are not \isoc as \ZZmos and 
the uniqueness 
does not hold.
If we had taken $s=1$ we could have defined two distinct \lon morphisms at $1$, namely $\phi_1:N\to N,\;a\mapsto a$, and $\phi_2:N\to N,\;a\mapsto -a$, and the uniqueness would be guaranteed in the sense required in the statement. \eoe

In practice, we often construct a module by taking some \hbox{$\gA_i$-modules~$M_i$} and by patching them via their \lons $M_{ij}=M_i[1/s_j]$. In this case the modules $M_{ij}$ and $M_{ji}$ are distinct, and we must give for each $(i,j)$ an \iso of $\gA_{ij}$-modules $\theta_{ij}:M_{ij}\to M_{ji}$.
This gives the following variant, in which the modules $M_{ijk}$ are not given in the hypothesis, but where we indicate the compatibility conditions that the $\theta_{ij}$'s need to satisfy.

\PLCC{plcc.modules 2}{{\em (Concrete patching of modules) }
 Let $S_1$, $\dots$, $S_n$ be \moco of $\gA$. 
 \\
Let $\gA_i=\gA_{S_i}$, $\gA_{ij}=\gA_{S_iS_j}$ and $\gA_{ijk}=\gA_{S_iS_jS_k}$.
\\
Assume we are given some $\gA_i$-modules $M_i$ and let

\snic{M_{j\ell}=M_j[1/s_\ell] \hbox{ and } M_{jk\ell}=M_j[1/s_ks_\ell] \hbox{ for all distinct } j,k,\ell
\in\lrbn,}

such that $M_{jk\ell}=M_{j\ell k}$, 
with the \lon morphisms  

\snic{\varphi_{j\ell}: M_{j}\to M_{j\ell}\;$
and $\;\varphi_{j\ell k}: M_{j\ell}\to M_{j\ell k}.}

Also assume we are given some morphisms of $\gA_{ij}$-modules 
$\theta_{ij}:M_{ij}\to M_{ji}$.
\\
Let $\theta_{ij}^{k}:M_{ijk}\to M_{jik}$ be the morphism of $\gA_{ijk}$-modules obtained by \lon at $s_k$ from $\theta_{ij}$.
Finally, we suppose that the following compatibility relations are satisfied 
\begin{itemize}
\item $\theta_{ji}\circ \theta_{ij}=\Id_{M_{ij}}$ for $i\neq j\in\lrbn$,
\item for distinct $i$, $j$, $k$ in $\lrbn$, by circularly composing
$$\preskip.1em \postskip.4em 
M_{ijk} \vvvers{\theta_{ij}^{k}} M_{jik}=M_{jki}
 \vvvers{\theta_{jk}^{i}} M_{kji}=M_{kij} \vvvers{\theta_{ki}^{j}}M_{ikj} 
$$
we must obtain the identity.  

%
\end{itemize}
Then, if $\big(M,(\varphi_i)_{i\in\lrbn}\big)$ is the limit of the diagram
$$\preskip.3em \postskip.4em 
\big((M_i)_{i\in \lrbn}),(M_{ij})_{i\neq j\in \lrbn};(\varphi_{ij})_{i\neq j},(\theta_{ij})_{i\neq j}\big), 
$$
each morphism $\varphi_i:M\to M_i$ is a \molo at $S_i$. \\
In addition, $\big(M,(\varphi_{i})_{i\in\lrbn}\big)$ is, up to unique \isoz, the unique \sys that makes the diagram commutative and that makes each $\varphi_i$ a \molo at $S_i$.
{\small
$$\preskip.2em \postskip-.4em\hspace*{-.2em}
\xymatrix @R=1em @C=.15em{
&M_i\ar[ddl]_(.4){\varphi_{ij}}\ar[ddrrr]^(.4){\varphi_{ik}} 
&&&&
M_j\ar[ddlll]_(.4){\varphi_{ji}}\ar[ddrrr]^(.4){\varphi_{jk}} 
&&&&
M_k\ar[ddlll]_(.4){\varphi_{ki}}\ar[ddr]^{\varphi_{kj}} 
\\
\\
M_{ij}\ar@<0.4ex>[rr]^{\theta_{ij}} \ar[dddrr]_{\varphi_{ijk}}    
&&
M_{ji}\ar@<0.4ex>[ll]^{\theta_{ji}} \ar[ddrrr]^(.7){\varphi_{jik}}
&\hspace{3.5em}&
M_{ik}\ar@<0.4ex>[rr]^{\theta_{ik}}  \ar[dddll]_(.6){\varphi_{ikj}}
&&
M_{ki}\ar@<0.4ex>[ll]^{\theta_{ki}} \ar[dddrr]^(.6){\varphi_{kij}}
&\hspace{3.5em}&
M_{jk}\ar@<0.4ex>[rr]^{\theta_{jk}} \ar[ddlll]_(.7){\varphi_{jki}}
&&
M_{kj}\ar@<0.4ex>[ll]^{\theta_{kj}} \ar[dddll]^{\varphi_{kji}}
\\
\\
&\hspace{2em}&     
&&&
M_{jik}\ar@{-->}[drrr]^{\theta_{jk}^{i}}  
&&&  
&\hspace{2em}
\\
&&M_{ijk}\ar@{-->}[urrr]^{\theta_{ij}^{k}}     
&&&  
&&&
M_{kij}\ar@{-->}[llllll]^{\theta_{ki}^{j}}  
&
}
$$
}
}
\begin{proof}
Note that the diagram above is commutative by construction, except eventually the bottom triangle in dotted lines, each time that it is possible to join two modules using two different paths: for example~\hbox{$\varphi_{ij}\circ \varphi_{ijk}= \varphi_{ik}\circ \varphi_{ikj}$} and $\theta_{ij}^{k}\circ \varphi_{ijk}=\varphi_{jik}\circ \theta_{ij}$. \\  
Here we need to convince ourselves that the indicated compatibility conditions are exactly what is \ncr and sufficient to be brought back to the situation described in the \plgref{plcc.modules 2}.
\\
For this, when $i<j<k$ we only 
keep 
$M_{ij}$, $M_{ik}$, $M_{jk}$ \hbox{and $M_{ijk}=M_{ikj}$}. 
\\
This forces us to replace
\[ 
\begin{array}{rclcrcl} 
 \varphi_{ji} & :  & M_j\to M_{ji} & \;\;\hbox{with}\;\;&\gamma_{ji}=\theta_{ji}\circ \varphi_{ji} &: & M_j\to M_{ij}, \\[1mm] 
 \varphi_{ki}  & :  & M_k\to M_{ki}  & \;\hbox{with}\;&\gamma_{ki}=\theta_{ki}\circ \varphi_{ki}  &: & M_k\to M_{ik},  \\[1mm] 
\varphi_{kj}   & :  & M_k\to M_{kj}  & \;\hbox{with}\;& \gamma_{kj}=\theta_{kj}\circ \varphi_{kj} &: &  M_k\to M_{jk}, \\[1mm] 
 \varphi_{jki}   & :  & M_{jk}\to M_{jik}  & \;\hbox{with}\;&\gamma_{jki}=\theta_{ji}^{k}\circ \varphi_{jki}  &: & M_{jk}\to M_{ijk}.   
 \end{array}
\] 
So far everything is taking place unhindered (in relation to the modules with two and three indices that we chose to preserve): the squares 
$(M_i,M_{ij},M_{ijk},M_{ik})$ and~$(M_j,M_{ij},M_{ijk},M_{jk})$ are commutative and the arrows are \molosz.
\\
It is only with the two \molos $M_k\to M_{ijk}$ that we will see the \pbz.
{\small
$$\preskip-.2em \postskip-.2em \hspace*{0mm}
\xymatrix @R=1.2em @C=.15em{
&M_i\ar[ddl]_(.4){\varphi_{ij}}\ar[ddrrr]^(.3){\varphi_{ik}} 
&&&&
M_j\ar[ddlllll]_(.3){\gamma_{ji}}\ar[ddrrr]^(.3){\varphi_{jk}} 
&&&&
M_k\ar[ddlllll]_(.3){\gamma_{ki}}\ar[ddl]^{\gamma_{kj}} 
\\
\\
M_{ij} \ar[ddrrr]_{\varphi_{ijk}}    
&&
&\hspace{2.5em}&
M_{ik}  \ar[ddl]_{\varphi_{ikj}}
&&
&\hspace{2.5em}&
M_{jk} \ar[ddlllll]^{ \gamma_{jki}}
&&
\\
\\
&&&M_{ijk}     
&&&  
&&&
}
$$
}

These two \molos are now imposed, namely the one that passes through $M_{ik}$, which must be 
$$\preskip.3em \postskip.3em 
{\varphi_{ikj}\circ \gamma_{ki}=\varphi_{ikj}\circ\theta_{ki}\circ \varphi_{ki}=\theta_{ki}^{j}\circ \varphi_{kij}\circ \varphi_{ki},}
$$
and the one that passes through $M_{jk}$, which must be 

\snic{\gamma_{jki}\circ \gamma_{kj}=\theta_{ji}^{k}\circ \varphi_{jki}\circ\theta_{kj}\circ \varphi_{kj} 
=\theta_{ji}^{k}\circ\theta_{kj}^{i}\circ\varphi_{kji}\circ \varphi_{kj}.}

As $\varphi_{kij}\circ \varphi_{ki}=\varphi_{kji}\circ \varphi_{kj}$, the fusion is successful if $\theta_{ki}^{j}=\theta_{ji}^{k}\circ\theta_{kj}^{i}$.
 \\
Actually the condition is also \ncr because \gui{every \molo is an epimorphism}: if $\psi_1\circ \varphi=\psi_2\circ \varphi$ with $\varphi$ being a \moloz, then $\psi_1=\psi_2$. 
\end{proof}
%

\subsec{Patching of \homos between \risz}
\label{subsecRecolAnn}

\begin{definition}
\label{defHomloc}
Let $S$ be a \mo of $\gA$.
A morphism $\alpha  :\gA\to \gB$ is called a \ix{localization morphism at $S$} if every morphism~\hbox{$\psi:\gA\to\gC$} such that $\psi(S)\subseteq \gC\eti$ can be uniquely factored by $\alpha$.%
\index{morphism!localization --- (\risz)}%
\index{localization!morphism (\risz)}%
\end{definition}

\rem 
If $\alpha :\gA\to\gB$ is a \lon \homoz, and if~\hbox{$S=\alpha^{-1}(\gB^\times)$}, then
$\gB$ is canonically \isoc to $\gA_S$.
Moreover, a \molo can also be \care as follows
\begin{enumerate}
\item [--] $\forall x,x'\in \gA\;(\alpha(x)=\alpha(x')\;\iff\; \exists s\in
S\;sx=sx')$
\item [--] $\forall y\in \gB,\;\exists x\in \gA,\;\exists s\in S\;\;sy=\alpha(x)$.
\eoe
\end{enumerate}

\medskip In the theory of schemes developed by Grothendieck, the \molos $\gA\to\gA[1/s]$ play a preponderant role.

We have already discussed at the beginning of this section (Section~\ref{secColleCiseaux}) the impossibility of patching \ris in \gnlz, with the example of $\Pn(\CC)$, which leads to the \dfn of schemes.

The possibility of defining a category of  schemes as \gui{patchings of \risz} ultimately relies on the following \prcc for  \homos between \risz.
The \dem of the principle is very simple. The important thing is that the morphism is uniquely defined using \lons and that the \cpbt conditions are themselves described via more advanced \lonsz.

\begin{plcc}
\label{plcc.RecolHomAnn}  {\em (Patching of morphisms of \risz) } 
Let $\gA$ and $\gB$ be two \risz, $s_1$, $\ldots$, $s_n$ be \eco of $\gA$ and $t_1$, $\ldots$, $t_n$ be \eco of $\gB$. Let 
$$\preskip.3em \postskip.3em
\gA_i=\gA[1/s_i],\; \gA_{ij}=\gA[1/s_is_j],\;  \gB_i=\gB[1/t_i]\; \mathit{and} \;\gB_{ij}=\gB[1/t_it_j].
$$
For each $i\in\lrbn$, let $\varphi_i:\gA_{i}\to\gB_{i}$ 
be a \homoz. Suppose that the following \cpbt conditions are satisfied: for $i\neq j$
the two \homos  $\beta_{ij}\circ \varphi_i:\gA_{i}\to\gB_{ij}$ and $\beta_{ji}\circ\varphi_j:\gA_{j}\to\gB_{ij}$ 
can be factorized via  $\gA_{ij}$ and give the same \homo $\varphi_{ij}:\gA_{ij}\to\gB_{ij}$ 
(see the diagram).
$$\preskip-.4em \postskip.4em\hspace*{3em}
\xymatrix @C=2em @R=1em 
          {
    & \gA_i \ar[ddr]^(.3){\alpha_{ij}}\ar[rrr]^{\varphi_{i}}&&&   
                  \gB_i \ar[ddr]^{\beta_{ij}}&\\
\gA\phantom{_{ij}} \ar@{..>}[ur]^{\alpha_i}\ar@{..>}[ddr]_{\alpha_j}\ar@{-->}[rrr]^(.7){\varphi!}
&&& ~\phantom{_i} \gB\phantom{_j}\ar@{..>}[ur]^{\beta_i}\ar@{..>}[ddr]^(.3){\beta_j} \\
&& \gA_{ij}\ar@{.>}[rrr]^(.3){\varphi_{ij}}
&&&\gB_{ij}&\\
                 & \gA_j \ar[ur]^{\alpha_{ji}}\ar[rrr]^{\varphi_{j}}&&&   
                  \gB_j \ar[ur]^{\beta_{ji}} 
}
$$
Then there exists a unique \homo $\varphi:\gA\to\gB$ such that for each $i$, we have $\varphi_i\circ \alpha_i=\beta_i\circ \varphi$.
\end{plcc}

\begin{proof}
The \cpbt conditions are clearly \ncrsz.
Let us show that they are sufficient. 
By the \plgrf{plcc.modules 1},  $\gB$ is the limit of the diagram of the $\gB_i$'s, $\gB_{ij}$'s and $\beta_{ij}$'s. The \cpbt conditions imply that we also have the \egts 

\snic{\beta_{ij}\circ (\varphi_i\circ \alpha_i)=
\beta_{ji}\circ (\varphi_j\circ \alpha_j)}

which are the conditions  guaranteeing the existence and the uniqueness of~$\varphi$.
\end{proof}

\section{The basic \cov local-global machinery} 
\label{secMachLoGlo}\imlb

\begin{flushright}
{\em Therefore localize at any \idepz.
}\\
A classical mathematician
\end{flushright}

Recall that we presented in Section~\ref{subsecDyna}
the \gnl philosophy of the dynamic method in \cov \algz.
\\
We now indicate how several proofs using the \plg in abstract \alg can be decrypted into \prcos leading to the same results in an explicit form.

In Section~\ref{subsecLGIdeMax} we will focus on the decryption of abstract \dems that use the quotients by all the \idemas and in Section~\ref{subsecLGIdepMin} we will focus on the decryption of abstract \dems that use \lons at all the \idemisz.

\subsec{Decryption of classical proofs using \lon at all primes}
A typical argument of \lon works as follows in \clamaz. 
When the \ri is local a certain \prt $\sfP$ is satisfied in virtue of a sufficiently concrete proof. 
When the \ri is not local, the same \prt is still true (from a non\cof classical point of view) because it suffices to satisfy it \lotz.
This in virtue of an \plgaz.

We examine with some attention the first \demz. We then see certain computations appear that are feasible in virtue of the following principle

\sni \centerline{$\forall x\in\gA\;\;\; x\in\Ati\; \vee\; x\in \Rad(\gA)$,}

\sni a principle which is applied to \elts $x$ derived from the \dem itself.
In other words, the given classical \dem in the local case provides us with a \prco under the  hypothesis of a \dcd \aloz.
Now here is our \cof dynamic decryption. In the case of an arbitrary \riz, we repeat the same \demz, by replacing each disjunction  \gui{$x$ is \iv or $x$ is in the radical} with the consideration of the two \ris $\gA_{\cS(I,x;U)}$ and $\gA_{\cS(I;x,U)}$, where $\gA_{\cS(I,U)}$ is the \gui{current} \lon of the starting \ri $\gA$, at this point in the \demz.
When the initial \dem is deployed thus, we will have constructed in the end a certain, finite because the \dem is finite, number of localized \risz~$\gA_{S_i}$, 
for which the \prt is true.
From a \cof point of view,  we obtain at least the \gui{quasi-global} result, \cad after \lon at \mocoz, in virtue of Lemma \ref{lemRecouvre}.
We then call upon a \plgc to conclude the result.

Our decryption of the classical proof is made possible by the fact that the \prt $\sfP$  is of \carf (see Section \ref{secPLGCBasic} from page \pageref{defiPropCarFini},  and Section~\ref{secSalutFini}):
it is preserved by \lonz, and if it is true after \lon at a \mo $S_i$, it is \egmt true after \lon at some $s_i\in S_i$.

The complete decryption therefore contains two essential ingredients. The first is the decryption of the given \dem in the local case which allows us to obtain a quasi-global result. The second is the \prco of the \plgc corresponding to the \plga used in \clamaz.
In all the examples that we have encountered, this \prco offers no difficulty because the \dem found in the classical literature already gives the concrete argument, at least in a telegraphic form
(except sometimes in Bourbaki, where the concrete arguments are skilfully hidden).

The \gnl conclusion is that the classical \dems \gui{by \plgaz} are already \covsz, if we bother to read them in detail. This is good news, other than the fact that this confirms that no supernatural miracles take place in \mathsz.

\smallskip The method indicated above therefore gives, as a corollary of Lemma~\ref{lemRecouvre}, the following \gnl decryption principle, which \emph{allows us to automatically obtain a global \cov version (or at least quasi-global) of a \tho from its local version}.

\mni {\bf Local-global machinery with \idepsz.}\label{MethodeIdeps}\\
{\it  When we reread a \prcoz, given for the case of a \dcd \aloz, with an arbitrary \ri $\gA$, such that at the start we consider it as $\gA=\gA_{\cS(0;1)}$ and at each disjunction (for an \elt $a$ that occurs during the computation in the local case)
$$\preskip.3em \postskip.3em 
a\in\Ati\; \vee\; a\in \Rad(\gA),
$$
we replace the \gui{current} \ri $\gA_{\cS(I,U)}$ with the two \ris
 $\gA_{\cS(I;U,a)}$ and $\gA_{\cS(I,a;U)}$ (in each of which the computation can be continued), at the end of the rereading we obtain a finite family of \ris $\gA_{\cS(I_j,U_j)}$ with the \com \mosz~${\cS(I_j,U_j)}$ and finite $I_j,\;U_j$. In each of these \risz, the computation has been successfully continued and has produced the desired result.
}

\medskip Please take note that if \gui{the current \riz} is $\gA'=\gA_{\cS(I;U)}$
 and if the disjunction relates to
$$\preskip-.2em \postskip.2em 	 
{b\in{\gA'}^{\times}\; \vee\; b\in \Rad(\gA'),} 
$$
with $b=a/(u+i)$, $a\in\gA$, $u\in\cM(U)$ and $i\in \gen{I}_\gA$, then the localized \ris $\gA_{\cS(I;U,a)}$ and $\gA_{\cS(I,a;U)}$ must be considered.

 In what follows we will speak of the \lgbe machinery with \ideps as we do of the \gui{basic \lgbe machinery.}\imlb

\subsect{Examples of the basic \lgbe machinery}{Examples}
\label{chapApPlg}\label{chapComplements}
\vspace{3pt}
\subsubsec{First example}
We want to prove the following result.

\begin{lemma}\label{lemPrimUnitaire}
Let $f\in\gA[X]$ be a primitive \pol and $r\in\gA$ be a \ndz \elt with $\Kdim\gA\leq1$.
Then the \id $\gen{f,r}$ contains a \poluz.
\end{lemma}
\begin{proof}
We begin by proving the lemma in the case where $\gA$ is a \dcd \aloz.
We can write $f=f_{1}+f_{2}$ with $f_{1}\in(\Rad\gA)[X]$ and $f_{2}$ pseudo\monz.
Moreover, for every $e\in\Rad \gA$ we have an \egtz~\hbox{$r^m(e^m(1+ye)+zr)=0$},
so $r$ divides $e^m$. Consequently $r$ divides a power of $f_{1}$, say with exponent $N$.
We have ${f_{2}}^N=(f-f_{1})^N\in\gen{f,f_{1}^N}\subseteq\gen{f,r}$.
Then,~$f_2^N$ provides the \polu required.

For an arbitrary \ri we re-express the previous \dem dynamically. For example if $f=aX^2+bX+c$, we explicate the previous reasoning in the following form.

\vspace{2pt}
\snac{\xymatrix {
        &\cS(0;1)\ar[ld]\ar[rd] \\
\dessus {\cS(0;a)} {a^\NN} &    &\dessus{\cS(a;1)} {1 + \langle a\rangle} \ar[ld]\ar[rd] \\
        &\dessus {\cS(a;b)} {b^\NN  + \langle a\rangle}&  &
                                    \dessus {\cS(a,b ; 1)} {1 + \langle a, b\rangle}
                                           \ar[ld]\ar[rd] \\
        &    &\dessus {\cS(a,b;c)} {c^\NN + \langle a,b\rangle}&&
                 \dessus {\cS(a,b,c;1)} {1 + \langle a,b,c\rangle} \\
}
}

\smallskip Either $a$ is \ivz, or it is in the radical. 
If $a$ is \ivz, then we take $f_{2}=f,\,f_{1}=0$.

Otherwise, either $b$ is \ivz, or it is in the radical.
If $b$ is \ivz, then we take $f_{2}=bX+c, \,f_{1}=aX^2$.

Otherwise, either $c$ is \ivz, or it is in the radical.
If  $c$ is \ivz, then we take $f_{2}=c, \,f_{1}=aX^2+bX$.
Otherwise $\gen{1}=\gen{a,b,c}\in\Rad \gA$ so the \ri is trivial.

See above the graph of the tree of the successive \lonsz. The \moco are found at the leaves of the tree, the last one contains $0$ and does not intervene in the computation.

Let us complete the \dem by indicating how we construct a \polu in the \id $\gen{f,r}$ of ${\gA_{\cS(I,U)}[X]}$ from two \polus $g$ and $h$ in the \ids $\gen{f,r}$ of ${\gA_{\cS(I,y;U)}[X]}$ and ${\gA_{\cS(I;y,U)}[X]}$.
On the one hand we have
$$\preskip.4em \postskip.2em 
sg=sX^m+g_{1}\hbox{  with  }\deg g_1 < m,\, 
s\in\cS(I,y;U)\hbox{ and }sg\in\gen{f,r}_{\AX}, 
$$
and on the other hand
$$\preskip.2em \postskip.3em 
th=tX^n+h_{1}\hbox{  with  }\deg h_1 < n,\, 
t\in\cS(I;y,U)\hbox{ and }th\in\gen{f,r}_{\AX}. 
$$
The \pols $sX^ng$ and $tX^mh$ of formal degree $n+m$ have for \fmt leading \coes $s$ and $t$. By taking $us+vt\in \cS(I,U)$, the work ends with $usX^ng+vtX^mh$.
\end{proof}

\subsubsect{Second example: a quasi-global result obtained from a given \dem for a \alo}{Second example: a quasi-global result}
\label{quasiglobaldynamique}

{\it Dynamic reread of the local freeness lemma.}
The dynamic rereading of  \gui{Azumaya's \demz} (\paref{lelilo}) of the local freeness lemma gives a new \dem of the \tho which states that the \mptfs are \lot free, with the following precise formulation.

\textit{If $F\in\Mn(\gA)$ is a \mprnz, there exist $2^n$ \eco $s_{i}$ such that over each $\gA_{s_{i}}$, the matrix is similar to a standard \mprnz. More \prmtz, for each $k=0,\ldots,n$ there are $n \choose k$ \lons at which the matrix is similar to $\I_{k,n}$.}

\rdb
First recall (see the graph below) how the computation tree for a \alo with a matrix $F$ in $\MM_3(\gA)$ is presented.

\def \AAA {\ar@{-}[dll]\ar@{-}[drr]}
\def \BBB {\ar@{-}[dllll]\ar@{-}[drrrr]}
\snuc{\hspace*{-3mm}
\xymatrix  @C = 1pt @R = 1.2cm{
&&&&&&&&&&&&&&& 1\ar@{-}[dllllllll]\ar@{-}[drrrrrrrr]
\\
&&&&&&&\BBB2 &&&&&&&&&&&&&&&&\BBB3
\\
&&&\AAA4 &&&&&&&&\AAA5 &&&&&&&&\AAA6 &&&&&&&&\AAA7
\\
&8 &&&&9 &&&&10 &&&&11 &&&&12 &&&&13 &&&&14 &&&&15
} \label{arbrebinaire}
}

\medskip
At the point $1$ the computation starts with the test \gui{$f_{11}$ or $1-f_{11}$ is \ivz} (note that the disjunction is \gnlt not exclusive, and the test must only certify that one of the two possibilities takes place).
If the test certifies that~$f_{11}$ is \ivz, we follow the left branch, we go to $2$ where we make a base change  that allows us to reduce the matrix to the form 
{\blocs{.5}{1}{.5}{1}{$1$}{$0$}{$0$}{$G$}}  with $G\in\MM_2(\gA)$ and $G^2=G$.
If the test certifies that $1-f_{11}$ is \ivz, we follow the right branch, we go to $3$ where we make a base change that allows us to reduce the matrix to the form \smashbot{\blocs{.5}{1}{.5}{1}{$0$}{$0$}{$0$}{$H$}} with $H^2=H$.

If we reach $2$, we test the \elt $g$ in position $(1,1)$ in $G$.
According to the result, we head towards $4$ or $5$ to make a base change that reduces us to one of the two forms $\Cmatrix{4pt}{1&0&0\cr0&1&0\cr0&0&a}$  with $a^2=a$, or $\Cmatrix{4pt}{1&0&0\cr0&0&0\cr0&0&b}$ with $b^2=b$.
If we reach $3$, we test the \elt $h$ in position $(1,1)$ in~$H$.
According to the result, we head towards $6$ or $7$ to make a base change that reduces us to one of the two forms $\Cmatrix{4pt}{0&0&0\cr0&1&0\cr0&0&c}$ with $c^2=c$, or $\Cmatrix{4pt}{0&0&0\cr0&0&0\cr0&0&d}$ with $d^2=d$.

In all cases, we finish with an invertibility test that certifies that the \idm is equal to $1$ or to $0$, which gives one of the $8$ possible diagonal \mprns (with only $0$'s and $1$'s on the diagonal).

If we dynamically reread this computation with an arbitrary \riz, we obtain the following \come \lonsz. \\
At the start at $1$, we have the \riz~\hbox{$\gA_{1}=\gA$}. At $2$ and $3$ we have the \come \lons $\gA_{2}=\gA_{1}[1/f_{11}]$ and $\gA_{3}=\gA_{1}[1/(1-f_{11})]$.
At $4$ and $5$ we have the following \come \lons of $\gA_{2}$: $\gA_{4}=\gA_{2}[1/g]$ and $\gA_{5}=\gA_{2}[1/(1-g)]$. At $6$ and $7$ we have the following \come \lons of $\gA_{3}$: $\gA_{6}=\gA_{3}[1/h]$ and $\gA_{7}=\gA_{3}[1/(1-h)]$. \\
We move on to the final level. At $8$ and $9$ we create the following \come \lons of $\gA_{4}$: $\gA_{8}=\gA_{4}[1/a]$ and $\gA_{9}=\gA_{4}[1/(1-a)]$.  At $10$ and $11$ we create the following \come \lons of $\gA_{5}$: $\gA_{10}=\gA_{5}[1/b]$ and~\hbox{$\gA_{11}=\gA_{5}[1/(1-b)]$} etc. 

Ultimately, by considering the \denos  $d_{i}$ ($i=8,\ldots,15$)
of the fractions created in the different branches (for example $d_{11}$ is the \deno in $\gA$ of the fraction $1/f_{11}(1-g)(1-b)$, where $g\in\gA_{2}$ and~\hbox{$b\in\gA_{5}$}), we obtain eight comaximal \elts of $\gA$, and for each of the \lons we obtain the corresponding reduced diagonal form of the starting matrix~$F$.

In other words, the dynamic rereading of the \dem given in the case of a \alo created eight \ecoz, where the abstract classical \dem would instruct us to localize at all the \idemasz, which could take quite some time.

\penalty-2500
\section{Quotienting by all the \idemas}
\label{subsecLGIdeMax}
\begin{flushright}
{\em A \ri that has no \idemas is reduced to $0$.
}\\
A classical mathematician
\end{flushright}

In the literature we find a certain number of \dems in which the author proves a result by considering \gui{the passage to the quotient by an arbitrary \idemaz.} The analysis of these \dems shows that the result can be understood as the fact that a \ri obtained from more or less complicated constructions is actually reduced to $0$.
For example, if we want to prove that an \id $\fa$ of $\gA$ contains $1$,
we reason by contradiction, we consider a \idema $\fm$ that would contain $\fa$, and we find a contradiction by making a computation in the residual field $\gA\sur{\fm}$.

This comes down to applying the quoted principle:
a \ri that has no \idemas is reduced to $0$.

The idea of presenting the reasoning as a \dem by contradiction is the result of an occupational bias.
Proving that a ring is reduced to $0$ is a fact of a concrete nature (we must prove that $1 = 0$ in the considered ring), and not a contradiction, and the computation performed in the field~$\gA\sur{\fm}$ only leads to a contradiction because one day we decided that, in a field, $1=0$ is prohibited. But the computation has nothing to do with such a prohibition.
The computation in a field uses the fact that every \elt is null or \ivz, but not the fact that this disjunction would be exclusive.

Consequently, the dynamic reread of the \dem by contradiction in a \prco is possible according to the following method. 
Let us follow the computation that we are required to do as if the \ri $\gA\sur{\fa}$ were truly a field. Each time that the computation demands to know if an \elt $x_{i}$ is null or \iv modulo $\fa$, let us bet on $x_{i}=0$ and add it to $\fa$. After a while, we find that $1=0$ modulo the constructed \idz. Instead of losing courage in the face of such a contradiction,
let us take a look at the good side of things. For example we have just observed that $1\in\fa+\gen{x_{1},x_{2},x_{3}}$. This is a positive fact and not a contradiction. We have actually just computed an inverse~\hbox{$y_{3}$} of $x_{3}$ in $\gA$ modulo $\fa+\gen{x_{1},x_{2}}$.
We can therefore examine the computation that the classical \dem requires us to do when $x_{1},x_{2}\in\fm$ and $x_{3}$ is \iv modulo $\fm$. Except that we do not need~$\fm$ since we have just established that $x_{3}$ is \iv
modulo $\fa+\gen{x_{1},x_{2}}$.

Contrary to the strategy that corresponds to the \lon at any \idepz, we do not try to deploy all of the computation tree 
that seems to reveal itself to us.
We only use quotients, and for this we systematically follow the \gui{to be  null} branch 
(modulo $\fm$)
rather than the \gui{to be \ivz} branch. 
 This creates more and more advanced successive quotients. When a so-called contradiction appears, \cad when a computation reaches a certain result of positive nature, we backtrack by taking advantage of the information that we have just collected: an \elt has been certified \iv in the previous quotient.

For example, with a deployed tree of the type of that of \paref{arbrebinaire} and by taking as its \gnl context the \ri $\gA\sur{\fa}$, if every time the right branch corresponds to $x=0$ and the left one to an \iv $x$, it is necessary to start by following the path $1\to3\to7\to15$ and to consider the successive quotients. At $15$ the computation gives us a positive result which allows us to backtrack to $7$ to follow the branch $7\to14$. At~$14$ a positive result allows us to backtrack to point $3$ (by the path $14\to 7\to 3$) by knowing that the \eltz~$a_3$ that produces the disjunction at this point is actually \ivz.
We can then follow the proposed computation for the branch $3\to6\to13$.
At~$13$ the classical \dem gives us a so-called contradiction, actually a positive result in the considered quotient at $6$.

We will ultimately follow the path

\snic{1\to 3\to 7\to 15\to 7\to 14 \to 3\to 6\to 13\to 6\to 12\to 
 1\to 2\to 5\to 11\to 5\to 10\to 2\to 4\to 9\to 4 \to 8\to 1.}

We will uniquely compute in quotients of $\gA\sur\fa$ and the final result is that $1=0$ in $\gA\sur\fa$,
\cad $\fa=\gA$, which was the pursued objective.

Note that during the first passage to the point $7$, we work with the \ri 
\smash{$\gA_{1,3,7}=\gA\sur{(\fa+\gen{a_1,a_3,a_7})}$}. Arriving at $15$, we learn that this \ri is trivial therefore that $a_7$ is \iv in $\gA_{1,3}=\gA\sur{(\fa+\gen{a_1,a_3})}$. At~$14$, we learn that $\gA_{1,3}$ is trivial, \cad $a_3$ is \iv in the \ri \hbox{$\gA_{1}=\gA\sur{(\fa+\gen{a_1})}$}.
We therefore head for the point $6$ with both the \ri $\gA_{1}$ and an inverse of~$a_3$ in hand \ldots\, 
Thus in repeated passages to the same point we are not working with the same ring, because new information accumulates as we progress through the computations.

The argument of passage to the quotient by all the \idemas of $\gA\sur{\fa}$ (assumed by contradiction non-reduced to $0$), which seemed a little magical, is thus replaced by a very concrete computation, implicitly given by the classical \demz.
Let us summarize the previous discussion.

\goodbreak
\mni {\bf Local-global machinery with \idemasz.}\label{MethodeIdemax}
\imlma\\
{\it To reread a classical \dem that proves by contradiction that a \ri $\gA$ is trivial by assuming the contrary, then by considering a \idemaz~$\fm$ of this \riz, by making a computation in the residual field and by finding the contradiction $1=0$, proceed as follows. First ensure that the \dem becomes a \prco that $1=0$ under the additional hypothesis that $\gA$ is a discrete field. Secondly, delete the additional hypothesis and follow step-by-step the previous \dem by favoring the branch $x=0$ each time that the disjunction \gui{$x=0$ or $x$ is \ivz} is required for the rest of the computation.
Each time that we prove $1=0$ 
we have actually showed that in the previously constructed quotient \riz, 
the last \elt to have undergone the test was \ivz, 
which allows us to backtrack to this point to follow the branch \gui{$x$ is \ivz} 
according to the proposed \dem for the \iv case (which is now certified).
If the considered \dem is sufficiently uniform (experience shows that it is always the case), the computation obtained as a whole is finite and ends at the desired conclusion.
}

\medskip \exl \\
The following crucial lemma was the only truly non\cof ingredient in the solution by \Sus of Serre's \pbz. We will expose this solution beginning on \paref{subsecSuslin} (see namely the \dem of \thref{th4SusQS}).
Here, we give the \dem of the crucial lemma by Suslin in \clamaz, then its \cof decryption.

\pagebreak	
\begin{lemma}\label{lemSuslin1}
Let $\gA$ be a \riz, $n$ be an integer $\geq2$ and $U=\tra[\,v_1\;\cdots\;v_{n}\,]$ be a \umd vector in $\gA[X]^{n\times 1}$ with $v_{1}$ \monz. \\
Let $V=\tra[\,v_2\;\cdots\;v_{n}\,]$.
There exist matrices $E_1$, \ldots, $E_\ell \in \EE_{n-1}(\gA[X]) $, such that, by letting $w_i$ be the first \coo of the vector $E_i V$, the \id $\fa$ below contains $1$

\snic{\fa=\gen{\Res_X(v_1,w_1), \Res_X(v_1,w_2), \ldots,\Res_X(v_1,w_\ell) }_\gA.}
\end{lemma}
\begin{proof}
If $n=2$, we have $u_1v_1+u_2v_2=1$ and since $v_1$ is \monz, $\Res(v_1,v_2)\in\Ati$:

\snac{\Res(v_1,v_2)\Res(v_1,u_2)=\Res(v_1,u_2v_2)=\Res(v_1,u_2v_2+u_1v_1)=\Res(v_1,1)=1.}

 If $n\geq3$, let $d_1=\deg v_1$. We suppose \spdg that the $v_i$'s are formal \pols of degrees $d_i<d_1$ ($i\geq2$). At the start we have some \pols $u_i$ such that $u_1v_1+\cdots+u_nv_n=1$.

\sni{\it Suslin's classical \demz.} We show that for every \idemaz~$\fm$, we can find a matrix $E_\fm\in \EE_{n-1}(\gA[X])$ such that, by letting~$w_\fm$ be the first \coo of $E_\fm V$, we have $1\in\gen{\Res_X(v_1,w_\fm)}$ modulo~$\fm$.
For this we work over the field $\gk=\gA\sur{\fm}$.
By using the Euclidean \algoz, the gcd $w_\fm$ of the $v_i$'s ($i\geq2$) is the first \coo of a vector obtained by \elr  manipulations over $V$. We lift the \elr matrix that was computed in $\EE_{n-1}(\gk[X])$ at a matrix $E_\fm\in\EE_{n-1}(\AX)$. Then, since $v_1$ and $w_\fm$ are coprime, the resultant $\Res_X(v_1,w_\fm)$ is nonzero in the field $\gA/\fm$.

{\it \Cov \dem (by decryption).}\\
We perform a \dem by \recu on the smallest of the formal degrees~$d_i$, which we denote by $m$ (recall that $i\geq2$).
To fix the ideas suppose that it is $d_2$. \\
Basic step: if $m=-1$, $v_2=0$ and by an \elr transformation we put $u_3v_3+\cdots+u_nv_n$ in position $2$, which brings us to the case $n=2$.
\\
Inductive step: from $m-1$ to $m$.
Let $a$ be the \coe of $v_2$ of degree $m$ and $\gB$ be the \ri  $\aqo{\gA}{a}$. In this \ri the \hdr is satisfied.
Thus, we have matrices  $E_1$, \ldots, $E_\ell\in \EE_{n-1}(\gB[X]) $, such that, by letting $\wi{w_i}$ be the first \coo of $E_iV$, we have the \egt
$$\preskip.2em \postskip.3em 
\gen{\Res_X(v_1,\wi{w_1}), \Res_X(v_1,\wi{w_2}), \ldots,\Res_X(v_1,\wi{w_\ell}) }_\gB = \gen{1}. 
$$
This means, by lifting the matrices in $\EE_{n-1}(\AX)$ without changing their name, and by letting $w_i$ be the first \coo of $E_iV$ that we have
$$\preskip.2em \postskip.3em
\gen{a,\Res_X(v_1,w_1), \Res_X(v_1,w_2), \ldots,\Res_X(v_1,w_\ell)}_\gA=\gen{1}.
$$
Then consider $\fb=\gen{\Res_X(v_1,w_1), \Res_X(v_1,w_2), \ldots,\Res_X(v_1,w_\ell) }_\gA$,
and \hbox{$\gC=\gA\sur{\fb}$}. Since $a$ is \iv in $\gC$, we can by an \elr manipulation replace $v_3$ with a \pol $v'_3=v_3-qv_2$ with $\deg v_3'\leq m-1$.
We apply the \hdr with the \ri $\gC$, we have \elr matrices $E'_1$, \ldots, $E'_q\in \EE_{n-1}(\gC[X])$ that we lift in $\EE_{n-1}(\AX)$ without changing their name. If $w'_1$, \ldots, $w'_q$ are the corresponding \pols (for  each $j$, $w'_j$ is the first \coo of $E'_jV$), we obtain
$$\preskip.3em \postskip.3em
1\in\gen{\Res_X(v_1,w_1), \ldots,\Res_X(v_1,w_\ell),\Res_X(v_1,w'_1), \ldots,\Res_X(v_1,w'_q)  }_\gA.
$$
\end{proof}

\comm Now let us see why this elegant \dem is indeed a decryption of that of Suslin according to the method indicated beforehand. Let $a_2=u_2v_2+\cdots+u_nv_n$.

When we want to treat the vector $V$ over a discrete field by the Euclidean \algoz, we have to do divisions. One division depends on the degree of the dividend (the \pol by which we divide). In the dynamic decryption, we therefore have tests to do on the \coes of the dividend to determine its degree. If we choose to start with the division of $v_3$ by $v_2$, the indicated method therefore requires us to first consider the case where $v_2$ is identically null. Note that this corresponds to the basic step of the \recuz. 

Let $\fa_1=\gen{(v_{2,i})_{i\in\lrb{0..d_2}}}$ be the \id generated by the \coes of $v_2$.
\\
If $v_2$ is identically null, we have the resultant $r_1=\Res(v_1,a_2)=\Res(v_1,w_1)$ (\ivz) with $w_1$ which is of the first \coo type of $E_1V$ for an explicit matrix $E_1\in\EE_{n-1}$.
\\
Naturally, this is only true modulo $\fa_1$, which gives $\fa_1+\gen{r_1}=\gen{1}$.
Let $\fa_2=\gen{(v_{2,i})_{i\in\lrb{1..d_2}}}$. We have established that $\fa_2+\gen{r_1}+\gen{v_{2,0}}=\gen{1}$.

We now reason modulo $\fb_2=\fa_2+\gen{r_1}$. Since $v_{2,0}$ is \iv and $v_2=v_{2,0}$, we can reduce to $0$ the vector $v_3$ by \elr manipulations then put in position $3$ an \elt equal to $a_2$ modulo $\fb_2$,
then bring it back to position $2$.
We therefore have a matrix $E_2\in\EE_{n-1}$ with~$w_2$ being the first \coo of $E_2V$ and $\Res(v_1,w_2)=r_2$ is \iv in $\gA/\fb_2$, \cad $\fa_2+\gen{r_1}+\gen{r_2}=\gen{1}$.
Let $\fa_3=\gen{(v_{2,i})_{i\in\lrb{2..d_2}}}$. We have just established 
that $\fa_3+\gen{r_1,r_2}+\gen{v_{2,1}}=\gen{1}$.

We now reason modulo $\fb_3=\fa_3+\gen{r_1,r_2}$.
Since $v_{2,1}$ is \iv and $\fa_3=0$, we can reduce the vector $v_3$ to a constant by \elr manipulations (corresponding to the division of $v_3$ by $v_2$), then bring it in position $2$.
We find ourselves in the situation studied previously (where $v_2$ was reduced to a constant). We therefore know how to compute two new \elr matrices $E_3$ and $E_4$ such that, by letting $w_3$ and $w_4$ be their first \coosz, and $r_i=\Res(v_1,w_i)$, we obtain $\fa_3+\gen{r_1,r_2,r_3,r_4}=\gen{1}$.
\\
Let $\fa_4=\gen{(v_{2,i})_{i\in\lrb{3..d_2}}}$. We have established that $\fa_4+\gen{r_1,r_2,r_3,r_4}+\gen{v_{2,2}}=\gen{1}$.

We now reason modulo $\fb_4=\fa_4+\gen{r_1,r_2,r_3,r_4}$.
Since $v_{2,2}$ is \iv and $\fa_4=0$, we can reduce the vector $v_3$ to the degree $1$ by \elr manipulations (corresponding to the division of $v_3$ by $v_2$),
then bring it in position $2$. We find ourselves in the situation studied previously (where $v_2$ was of degree $1$).
\,\ldots\,\ldots\, We obtain $\fa_4+\gen{r_1,r_2,\ldots,r_8}=\gen{1}$.
\\
Let $\fa_5=\gen{(v_{2,i})_{i\in\lrb{4..d_2}}}$.
We have established that $\fa_5+\gen{r_1,r_2,\ldots,r_8}+\gen{v_{2,3}}=\gen{1}$.

And so on and so forth \, \ldots \,\ldots

The important part of this is that the inverses of leading \coes of successive $v_2$ that appear in the \algo are always computed as \elts of the \ri and not by a \lon procedure.
Each time they are only \iv modulo a certain specified \idz, but it does not matter, the \id grows by incorporating the authorized resultants but decreases by expelling intruders that are the \coes of $v_2$.
\eoe

\penalty-2500
\section{Localizing at all the minimal \idepsz}
\label{subsecLGIdepMin}
\begin{flushright}
{\em A \ri that has no minimal \ideps is reduced to $0$.
}\\
A classical mathematician
\end{flushright}

The readers are now called upon to convince themselves of the correctness of the following method, by replacing in the previous section addition by multiplication and passage to the quotient by localization.

\mni {\bf Local-global machinery with \idemisz.}\label{MethodeIdemin}\\
{\it
To reread a classical \dem that proves by contradiction that a \riz~$\gA$ is trivial by assuming the contrary, then by considering a \idemi of this \riz, by making a computation in the localized \ri (which is local and \zedz, therefore a field in the reduced case) and by finding the contradiction $1=0$, proceed as follows.
\\
First ensure that the \dem becomes a \prco of the \egt $1=0$ under the additional hypothesis that $\gA$ is local and \zedz. Secondly, delete the additional hypothesis and follow step-by-step the previous \dem by favoring the \gui{$x$ is \ivz} branch each time that the disjunction \gui{$x$ is nilpotent or $x$ is \ivz} is required for the rest of the computation.
Each time that we prove $1=0$ we actually have shown that in the previously constructed localized \riz, the last \elt to be subjected to the test was nilpotent, which allows us to backtrack to this point to follow the \gui{$x$ is nilpotent} branch according to the proposed \dem for the nilpotent case (which is now certified).
If the considered \dem is sufficiently uniform (experience shows that this is always the case),  the computation obtained as a whole is finite and ends at the desired conclusion.

}

\medskip \exl
A quite spectacular example is given in the next chapter with the \cof decryption of an abstract \dem of Traverso's \tho regarding seminormal \risz.

\section{Local-global principles in depth $1$}
\label{secPlgcor}

Until now the different variants of the \plg were based on the families of \ecoz, that is on the finite families that generate the \id $\gen{1}$.
A weaker notion is sufficient for  questions of regularity: these are the finite families that generate a faithful \idz, or more \gnlt an \Erg \idz.

We say that they are families of \prof $\geq 1$.
In  Section \ref{secPlgprof2}, we will examine what we call the families of \prof $\geq 2$.

\begin{definition} \label{defiCoreg1} 
\begin{enumerate}
\item A finite family $(\an)$ of a \ri $\gA$ is called a \ix{system of coregular elements} if the \id $\gen{\an}$ is faithful.\footnote{Not to be confused with the notion of a coregular sequence introduced by Bourbaki, as a dual notion of that of a \srgz.}
\\
We also say that \emph{the \id $\fa$, or the list $(\an)$, is of \prof $\geq 1$}, and we express this in the form $\Gr_\gA(\an)\geq 1$.%
\index{coregular!\eltsz}%
\index{depth!finite family of --- $\geq 1$}
\item Let $E$ be an \Amoz.
\begin{itemize}
\item We say that an \elt $a\in\gA$ is  \emph{\Ergz} (or \emph{\ndz for~$E$}) if  

\snic{\forall x\in E,\; 
(ax=0\;\Longrightarrow\; x=0).}

\item A finite family $(\an)$ is said to be  \emph{once \Ergez} if

\snic{\forall x\in E,\; 
\big((a_1x=0,\,\dots,\,a_nx=0)\;\Longrightarrow\; x=0\big).}

We also say that the $a_i$'s are \emph{\cor for $E$}. 
\\
We express this in the form $\Gr_\gA(\an,E)\geq 1$. 
\item A \itf  $\fa\subseteq \gA$ is said to be \emph{\Ergz} if some (every) \sgr of $\fa$ is once \Ergz.
We also say that \emph{the depth of~$E$ relative to $\fa$ is greater than or equal to $1$}, and we express this in the form $\Gr_\gA(\fa,E)\geq 1$.%
\index{E-reg@\Ergz!ideal}\index{E-reg@\Ergz!\eltz}
\end{itemize}
\end{enumerate}
\end{definition}

Thus $\Gr_\gA(\ua)\geq 1$ means $\Gr_\gA(\ua,\gA)\geq 1$.
In what follows, we will often only give the statement with $\Gr_\gA(\ua,E)\geq 1$.

\medskip 
\rem The notation $\Gr(\fa,E)$ comes from \cite{Nor}. In this wonderful book,
Northcott defines the \gui{true grade} \`a la Hochster as the better 
non-\noe substitute for the usual depth.
\eoe

\begin{fact} \label{factdefiCoreg1}~
\begin{itemize}
\item The product of two \Ergs \itfs is \Ergz.
\item If $\fa\subseteq \fa'$ with $\fa$ \Ergz, then $\fa'$ is \Ergz.
\end{itemize}
\end{fact}

\begin{lemma} \label{lemCoreg1} \emph{($(a,b,ab)$ trick for  depth $1$)}\\
Suppose that the \ids $\gen{a,c_2, \dots, c_n}$ and $\gen{b,c_2, \dots, c_n}$ are \Ergsz. Then the \idz~$\gen{ab,c_2, \dots, c_n}$ is \Ergz. 
\end{lemma}
%
\begin{proof}
Let $x\in E$ such that $abx=c_2x=\cdots=c_nx=0$.\\
Then~$abx=c_2bx=\cdots=c_nbx=0$,
so~$bx=0$. So $bx=c_2x=\cdots=c_nx=0$, therefore~$x=0$.
\end{proof}

We have the following \imd corollary.\footnote{We also could have noticed that for large enough $q$, the \id $\gen{\an}^{q}$, which \hbox{is \Ergz}, is contained in the \id $\gen{a_1^p,\dots,a_n^p}$.} 
%
\begin{lemma} \label{lemCoreg2}
Let $\gen{\an}$ be an \Erg \id and let $p_i\in\NN$. \\
Then the \id $\geN{a_1^{p_1},\dots,a_n^{p_n}}$ is \Ergz. 
\end{lemma}

We can compare the following \plg to items \emph{1} and \emph{3} of the \plgref{plcc.sli}.

Note that the statement \gui{$\fb$ is \Ergz} is stable under \lon when~$\fb$ is \tfz. This gives the implication in the direct sense for 
item~\emph{\ref{i3plcc.regularite}} in the following \plgz.

\begin{plcc}
\label{plcc.regularite}
\emph{(Localizations in \prof $\geq 1$)}\\
Let  $b$,  $a_1$, \dots, $a_n\in\gA$,   and~$\fb$ be a \itfz. Let $\gA_i=\gA[1/a_i]$.
\begin{enumerate}
\item Suppose that the~$a_i$'s are \corz.  
\begin{enumerate}
\item  We have $x=0$ in $\gA$   \ssi $x=0$ in each $\gA_i$.
\item  The \elt $b$ is \ndz \ssi it is  \ndz in $\gA_i$ for each~$i$.
\item   The \id $\fb$ is faithful \ssi it is faithful in $\gA_i$ for each~$i$.
\end{enumerate}%

\item Let $E$ be an \Amo and let $E_i=E[1/a_i]$. \\
Suppose that the \id $\gen{\an}$ is \Ergz.
\begin{enumerate}
\item \label{i1plcc.regularite} We have $x=0$ in $E$   \ssi $x=0$ in each $E_i$.
\item \label{i2plcc.regularite} The \elt $b$ is \Erg \ssi it is  $E_i$-\ndz for each~$i$.
\item  \label{i3plcc.regularite} The \id $\fb$ is \Erg \ssi it is  $E_i$-\ndz for each~$i$.
\end{enumerate}%
\end{enumerate}
\end{plcc}
\begin{proof} It suffices to treat item \emph{2.}
\\
\emph{\ref{i1plcc.regularite}.} If $x=0$ in  $E_i$ there is an exponent $k_i$ such that $a_i^{k_i}x=0$ in~$E$. We conclude by Lemma~\ref{lemCoreg2}  (with the module~$\gA x$) that~$x=0$.

\noindent \emph{\ref{i3plcc.regularite}.} Suppose that $\fb$ is $E_i$-\ndz for each~$i$, and~$\fb\,x=0$.
Then~$x=0$ in each~$E_i$, so~$x=0$ by item~\emph{\ref{i1plcc.regularite}.}
\end{proof}

We often implicitly use the following lemma, which is a variant of Lemma \ref{factLocCas} stated for the \syss of \ecoz.

\begin{fact} \emph{(Lemma of successive \core \lonsz)}
\label{factLocCasreg}\index{Successive \lons lemma, depth 1} \\
If $\Gr_\gA(s_1, \ldots, s_n,E)\geq 1$ and if for each $i$,
we have \elts
$
s_{i,1},\, \ldots ,\,
s_{i,k_i},
$
\cor for $E[1/s_i]$,
then the $s_is_{i,j}$'s are \cor for~$E$.
\end{fact}
%
\begin{proof}
Let $\fb$ be the \id generated by the $s_is_{i,j}$'s. By item \emph{2c} of the \plgref{plcc.regularite}, it suffices to prove that it is \Erg after \lon at \cor \elts for $E$. The $s_i$'s are suitable. 
\end{proof}
%

\subsec{McCoy's \thoz}

As an application of the \plgref{plcc.regularite}, we give a new \dem of McCoy's \tho (\ref{prop inj surj det} item~\emph{2}).

\CMnewtheorem{thoMccoy}{McCoy's \thoz}{\itshape}
\begin{thoMccoy} \label{propInjIdd} 
A matrix $M\in\Ae{m\times n}$ 
is injective \ssi
the \idd $\cD_n(M)$ is faithful.  
\end{thoMccoy}
%
\begin{proof}
The implication \gui{if} is simple. Let us show that if the matrix~$M$ is injective, the \idz~$\cD_n(M)$ is faithful. We perform an \recu on the number of columns.
Since~$M$ is injective, the \coes of the first column (which represents the image of the first basis vector),
 generate a faithful \idz.
By the \plgrf{plcc.regularite}, it therefore suffices to prove that~$\cD_n(M)$ is faithful over the \riz~$\gA_a=\gA[1/a]$, where~$a$ is a \coe of the first column. 
\\
Over this \ri it is clear that the matrix $M$ is \eqve to a matrix of the form $\blocs{.4}{.6}{.4}{.9}{$1$}{$0$}{$0$}{$N$}$. 
In addition~$N$ is injective
therefore by \hdr the \id $\cD_{n-1}(N)$ is faithful over~$\gA_a$. 
Finally~$\cD_{\gA_a,n-1}(N)=\cD_{\gA_a,n}(M)$.    
\end{proof}

\rems ~\\ 
1) The \dem also gives that if~$m<n$ and~$M$ is injective, then the \ri is trivial. Indeed at each step of the \recuz, when we replace~$M$ with~$N$ the difference~$m-n$ remains constant. Therefore if~$m<n$ we obtain \gui{at the base step} an injective map from~$\Ae0$ in~$\Ae{n-m}$ which implies~$1=0$ in~$\gA$. This is in accordance with the \gnl statement of \thref{propInjIdd}, because for $m<n$, $\cD_n(M)=0$, and if $0$ is a \ndz \eltz, the \ri is trivial.
 
2) We often find in the literature McCoy's \tho stated as follows, in a contrapositive form (in appearance).
\\
 \emph{If the \id is not faithful, the map is not injective}. 
 \\
Or more precisely. 
\\
\emph{If a nonzero \elt $x\in\gA$ annihilates $\cD_n(M)$, there exists a nonzero column vector $C\in\Ae{m\times 1}$ such that $MC=0$}.
\\
Unfortunately, this statement can only be proven with classical logic, and the existence of the vector $C$ cannot result from a \gnl \algoz. Here is a counterexample, well-known by numerical analysts.
If $M$ is a matrix with real \coes with $m<n$, we do not know how to produce a nonzero vector in its kernel so long as we do not know the rank of the matrix. For example for $m=1$ and $n=2$, we give two reals $(a,b)$, and we look for a pair $(c,d)\neq (0,0)$ such that $ac+bd=0$. If the pair~$(a,b)$ is a priori indistinguishable from the pair $(0,0)$, it is impossible to provide a suitable pair~$(c,d)$ so long as we have not established whether $\abs a+\abs b$ is null or not.
\Cov variants of the contraposition are proposed in Exercises~\ref{exoMcCoyContr1} and~\ref{exoMcCoyContr2}.
\eoe

\section{Local-global principles in depth 2}\label{secPlgprof2}

\begin{definition} \label{defiProf2} Let $a_1$, \dots, $a_n\in\gA$ and $E$ be an \Amoz.
\begin{itemize}
\item The list $(\ua)=(\an)$ is said to be \emph{of \prof $\geq 2$} 
if it is of \profz~\hbox{$\geq 1$} and if, for every list $(\ux)=(\xn)$ in $\gA$ proportional%
\footnote{Recall that this means that the \deters $\Dmatrix{.2em}{a_i&a_j\\ x_i&x_j}$ are all null.} 
to $(\ua)$, there exists an $x\in\gA$ such that $(\ux)=x (\ua)$. 
\\
We express this in the form $\Gr_\gA(\ua)\geq 2$ or $\Gr(\ua)\geq 2$.
\item   The list $(\ua)=(\an)$ is said to be \emph{twice \Ergez}
 if $\Gr_\gA(\ua,E)\geq 1$ and if, for every list $(\ux)=(\xn)$ in $E$ proportional to $(\ua)$ there exists an $x\in E$ such that $(\ux)=(\ua)x$. \\
We express this in the form $\Gr_\gA(\an,E)\geq 2$ or $\Gr(\ua,E)\geq 2$. 
\\
We also say%
\footnote{Eisenbud speaks of the \prof of $\fa$ over $E$, and Matsumura of the $\fa$-\prof of~$E$. The terminology adopted here is that of Bourbaki.} that \emph{the depth of~$E$ relative to $(\an)$ is greater than or equal to 2}.
\end{itemize}
\index{depth!finite family of --- $\geq 2$}
\end{definition}

\exls 1) In an integral \ri a list $(a,b)$ with $a,b\in\Reg(\gA)$ is of \profz~\hbox{$\geq 2$}
\ssi $\gen{a}\cap\gen{b}=\gen{ab}$, \cad $ab$ is the lcm of $a$ and $b$ in the sense of  \dvez.
\\
2) In a  GCD-domain a list $(\an)$ is of \profz~$\geq 2$ \ssi $1$ is the gcd of the list.
\\
3) If $n=1$ and the list is reduced to the single term $a$, $\Gr(a,E)\geq 2$ means that each $y\in E$ is of the form $y=ax$, \cad $aE=E$. 
In particular $\Gr_\gA(a)\geq 2$ means $a\in\Ati$.
\\
4) Every list of \eco is of \prof $\geq 2$ (by the basic \plgz). 
\eoe

\medskip It is clear that $\Gr(\ua)\geq 2$ means $\Gr(\ua,\gA)\geq 2$.
In the remainder this exempts us from duplicating the statements: we represent them with $\Gr(\ua,E)\geq 2$ for an arbitrary module~$E$ whenever  possible.

\begin{propdef} \label{lem1prof2}~\\
Let $(\ua)=(\an)$ and $(\ub)=(\br)$ in $\gA$ and $E$ be an \Amoz. If $\Gr_\gA(\ua,E)\geq 2$ and $\gen{\ua}\subseteq \gen{\ub}$, then $\Gr_\gA(\ub,E)\geq 2$. \\
Consequently, we say that \emph{a \itf $\fa$ is twice \Ergz}
 if every finite \sgr of $\fa$ is twice  \Erg
 (it suffices to verify it for a single one).
We express this in the form $\Gr_\gA(\fa,E)\geq 2$.
\end{propdef}
\begin{proof} It suffices to prove the two following facts.
\begin{itemize}
\item  If  $\Gr(\an,E)\geq 2$, then $\Gr(\an,b,E)\geq 2$. 
\item If {\mathrigid 2mu $a\in\gen{\an}$ and $\Gr(\an,a,E)\geq 2$, then $\Gr(\an,E)\geq 2$}. 
\end{itemize}
This indeed first shows that we can replace a \sgr of a \itf by another without changing \gui{the \prof $\geq 2$} and then that when we replace $\fa$ by a larger \itfz, the \prof $\geq 2$ is preserved.
\\
Let us consider the first item. We have a list $(\xn,y)$ in $E$ proportional to $(\an,b)$. We find some $x$ (unique, in fact) such that $(\ux)=(\ua)x$.
We must show that $bx=y$. However, $a_iy=bx_i$ and $bx_i=ba_ix$ for $i\in\lrbn$. 
\\
Therefore $a_i(y-bx)=0$ and we conclude that $y=bx$ because $\Gr(\ua,E)\geq 1$.
\\
The second item is left to the reader.
\end{proof}
%

\begin{lemma} \label{lemtrickprof2} \emph{($(a,b,ab)$ trick for  depth 2)}\\
Suppose that the lists $({\an,a})$ and $({\an,b})$ are twice \Ergesz.
 Then the list~$({\an,ab})$ is twice \Ergez. 
\end{lemma}
%
\begin{proof}
We already know that $({\an,ab})$ is once \Ergez.
\\
Let $(\xn,y)$ be a list in $E$ proportional to $({\an,ab})$.
The list $(x_1b,\dots,x_nb,y)$ is proportional to $({\an,a})$.
So there exists a~$z\in E$ such that 
$$
x_1b=a_1z,\,\dots,\,x_nb=a_nz,\,y=az
$$
This implies the list $(\xn,z)$ is proportional to $({\an,b})$. So there exists an  $x\in E$ such that 
$$
x_1=a_1x,\,\dots,\,x_n=a_nx,\,z=bx\;\hbox{ and a fortiori }y=abx
$$

\vspace{-1.5em}
\end{proof}
%


\begin{plcc} \label{plcc0Profondeur2} \label{plcc1Profondeur2} \emph{(For  \dve and  \icl \risz, \lons in depth 2)}
Consider a family $(\us)=(s_1,\dots,s_n)$ in $\gA$  with $\Gr_\gA(\us,E)\geq 2$. Let $\gA_{i}=\gA[\fraC1{s_i}]$ and $E_i=E[\fraC1{s_i}]$.
\begin{enumerate}
\item Let $a\in\gA$ be a \Erg \elt and $y\in E$. Then $a$ \gui{divides} $y$ in $E$ \ssi  $a$ divides $y$ after \lon at each $s_i$.
\item Let $(\bbm)$ in $\gA$. Then $\Gr_{\gA}(\bbm,E)\geq 2$ \ssi $\Gr_{\gA_i}(\bbm,E_i)\geq 2$ for each $i$.
\item \label{I2plccAnnDivl} Suppose that $\gA$ is integral and $\Gr_\gA(\us)\geq 2$. The \ri $\gA$ is \icl \ssi each \riz~$\gA_i$ is \iclz.
\end{enumerate}
\end{plcc}
\begin{proof} \emph{1.} Suppose that $a$ divides $y$ after \lon at $s_i$. We have $ax_i=u_iy$ in~$E$ for some $u_i=s_i^{n_i}$ and some $x_i\in E$. The list of the $u_i$'s is twice  \Erge  (Lemma \ref{lemtrickprof2}).
We have $au_jx_i=u_iu_jy=au_ix_j$ and as $a$ is \Ergz, $u_jx_i=u_ix_j$. Therefore we have some $x\in E$ such that $x_i=u_ix$ for each $i$.
This gives $u_iax=u_iy$ and as $\Gr(u_1,\dots,u_n,E)\geq 1$, we obtain $ax=y$.

\emph{2.} 
Consider in $\gA$ a sequence $(\xm)$ proportional to $(\bbm)$.
We seek some $x\in E$ such that $x_\ell=xc_\ell$ for every $\ell \in\lrbm$.
In each $E_i$ we find some $y_i$ such that $x_\ell=y_i c_\ell$ for every $\ell \in\lrbm$. 
This means that we have some $u_i\in s_i^{\NN}$ and some $z_i\in E$ such that $u_ix_\ell=z_i c_\ell$ in $E$ for every $\ell \in\lrbm$. It suffices to show that there exists some $z\in E$ such that  $z_i=u_i z$ for each $i$,
because then $u_i(x_\ell-z c_\ell)=0$ for each~$i$ (and the $u_i$'s are \cor for~$E$). It therefore suffices to show that the $z_i$'s form a family proportional to the~$u_i$'s, \cad $u_iz_j=u_jz_i$ for all $i,j\in\lrbn$.
However, we know that the $c_\ell$'s are \cor for $E$ (by the \plgref{plcc.regularite}). Therefore it suffices to show that we have the \egts $u_iz_jc_\ell=u_jz_ic_\ell$, but the two members are equal to $u_iu_jx_\ell$.
 
\emph{3.} Let $x$ and $y$ in $\gA$ with $y$ integral over the \id $x\gA$.
This remains true for each localized \ri $\gA_i$, which is \iclz.
Therefore $x$ divides $y$ in  each~$\gA_i$. Therefore by item \emph{1} with $E=\gA$, $x$ divides $y$ in $\gA$. 
\end{proof}
%
\begin{fact} \emph{(Successive \lons lemma, with depth 2)}
\label{lelosuccprof}\index{Successive \lons lemma, depth 2} \\
If $\Gr_\gA(s_1, \ldots,s_n,E)\geq 2$ and if for each $i$ we have a list $(s_{i,1}, \ldots, s_{i,k_i})$ in~$\gA$ which is  twice   $E[1/s_i]$-\ndzez,
then the \sys of the $s_{i}s_{i,j}$'s is  twice  \Ergz.
\end{fact}
\begin{proof}
Applying \ref{plcc0Profondeur2} \emph{2.}, it suffices to verify that the~$s_is_{ij}$'s are~twice  \Ergs
after \lon at \elts that form a list twice  \Ergez. 
This works with the list of the $s_i$'s. 
\end{proof}
%

\begin{lemma} \label{lem2prof2}
Let $(\ua)=(\an)$ and  $(\ub)=(\br)$ in $\gA$ and $E$ be an \Amoz. Let $(\ua \star \ub)$ be the finite family of the $a_ib_j$'s.
\\
If $\Gr_\gA(\ua,E)\geq 2$ and $\Gr_\gA(\ub,E)\geq 2$ then $\Gr_\gA(\ua\star\ub,E)\geq 2$.
\\
In terms of \itfsz: \\
$\bullet$ if $\Gr_\gA(\fa,E)\geq 2$ and $\Gr_\gA(\fb,E)\geq 2$ then $\Gr_\gA(\fa \fb,E)\geq 2$. 
\end{lemma}
\begin{proof}
Applying \ref{plcc0Profondeur2} \emph{2.}, it suffices to show that the family of~$a_ib_j$'s is twice \Erge after \lon in each $a_i$.
E.g., when localizing in $a_1$, the list of $a_1b_j$'s
generate the same ideal  \id as the $b_j$'s, and this \id is twice \Ergz. 
\end{proof}

\subsec{Patchings in depth 2}

The following \dfn allows us to simplify the writing of certain \dems a little.


\begin{definition} \label{defiProfMon2}  
\emph{(\Sys of \mos twice~\hbox{\Ergz})}\\
 A \sysz~$(S_1, \dots, S_n)=(\uS)$ of \mos of $\gA$ is said to be  \emph{twice \Ergz} 
if for all~$s_1\in S_1,\,\dots,\,s_n\in S_n$, we have~$\Gr_\gA(s_1,\dots,s_n,E)\geq 2$.
\end{definition}

The most important case is the \sys of \mos $(s_1^{\NN},\dots,s_n^{\NN})$ when $\Gr_\gA(s_1,\dots,s_n,E)\geq 2$.

We now re-express the \plgref{plcc.modules 1} by replacing the hypothesis according to which the \mos are \com by a weaker hypothesis (\sys of \mos twice \ndzz).

The context is the following.
Consider a \sys of \mos $(\uS)=(S_i)_{i\in\lrbn}$.
\\
Let $\gA_i:= \gA_{S_i}$ and  $\gA_{ij}:= \gA_{S_iS_j}$ ($i\neq j$)
such that $ \gA_{ij}= \gA_{ji}$. 
\\
Let~\hbox{$\varphi_i: \gA\to  \gA_i$} and $\varphi_{ij}: \gA_i\to  \gA_{ij}$ be the natural \homosz. 
\\
In what follows notations like $(E_{ij})_{i<j\in\lrbn}$ and $(\varphi_{ij})_{i\neq j\in\lrbn})$ mean that we have $E_{ij}=E_{ji}$ but (a priori) not $\varphi_{ij}=\varphi_{ji}$.

\begin{plcc}
\label{plcc.modules1bis}  {\em (Concrete patching of elements of a module in \prof 2) }  
Consider the context described above.
\begin{enumerate}
\item \label{i1plcc.modules1bis} 
Suppose that $(\uS)$ is twice \ndzz. 
Consider an \elt $(x_i)_{i\in\lrbn}$ of  $\prod_{i\in\lrbn}  \gA_i$.
So that there exists some~\hbox{$x\in  \gA$} satisfying \hbox{$\varphi_i(x)=x_i$} in each $ \gA_i$, it is sufficient and necessary that for each~\hbox{$i<j$} we have $\varphi_{ij}(x_i)=\varphi_{ji}(x_j)$ in $ \gA_{ij}$. In addition, this $x$ is then uniquely determined.
In other terms the \ri $ \gA$ (with the \homos $\varphi_{i}$) is the limit of the diagram
$$\preskip.2em \postskip.4em 
\big(( \gA_i)_{i\in\lrbn},( \gA_{ij})_{i<j\in\lrbn};(\varphi_{ij})_{i\neq j\in\lrbn}\big). 
$$

\item \label{i2plcc.modules1bis} 
Let $E$ be an \Amoz. Suppose that $(\uS)$ is twice \Ergz.
\\
Let $E_i:=E_{S_i}$ and  $E_{ij}:=E_{S_iS_j}$ ($i\neq j$)
such that $E_{ij}=E_{ji}$. 
\\
Let $\varphi_i:E\to E_i$ and $\varphi_{ij}:E_i\to E_{ij}$ be the natural \alisz.
Then the pair $\big(E,(\varphi_{i})_{i\in\lrbn}\big)$ gives the limit of the following diagram in the category of  \Amos
$$\preskip.3em \postskip-.4em 
\big((E_i)_{i\in\lrbn},(E_{ij})_{i<j\in\lrbn};(\varphi_{ij})_{i\neq j\in\lrbn}\big). 
$$
 
$$\preskip-.2em \postskip-.8em
\!\!\!\!\!\!\!\!
\xymatrix @C=3.5em @R=1.5em 
          {
                &&& E_i \ar[r]^{\varphi_{ij}}\ar[ddr]_(.3){\varphi_{ik}}
                & E_{ij}\\
F\ar[urrr]^{\psi_i} \ar[drrr]^{\psi_j}\ar[ddrrr]_{\psi_k}
\ar@{-->}[rr]^(.6){\psi!} &&E \ar[ur]_{\varphi_i}\ar[dr]^{\varphi_j} \ar[ddr]_(.5){\varphi_k} &&\\
  &&& E_j \ar[uur]_(.7){\varphi_{ji}}\ar[dr]
          &E_{ik}\\
&&& E_k \ar[ur]\ar[r]_{\varphi_{kj}} 
   & E_{jk}\\
}
$$
\end{enumerate}
\end{plcc}

\begin{proof}
 \emph{\ref{i1plcc.modules1bis}.} Special case of \emph{\ref{i2plcc.modules1bis}.}

\emph{\ref{i2plcc.modules1bis}.} Let $(x_i)_{i\in\lrbn}$ be an \elt of  $\prod_{i\in\lrbn}  E_i$. We must show the following \eqvcz: there exists an  $x\in  E$ satisfying $\varphi_i(x)=x_i$ in each $ E_i$ \ssi for each $i<j$ we have $\varphi_{ij}(x_i)=\varphi_{ji}(x_j)$ in $E_{ij}$. In addition, this $x$ must be unique.
\\
The condition is clearly \ncrz. Let us prove that it is sufficient.
\\
Let us show the existence of $x$. Let us first note that there exist some $s_i$'s in $S_i$ and some $y_i$'s in $E$ 
 such that we have $x_i=y_i/s_i$ in each $E_i$.
\\
If $\gA$ is integral, $E$ is \torf and \hbox{each $s_i\neq 0$}, we have in the \evc obtained by \eds to the quotient field the \egtsz

\snic{\frac{y_1}{s_1}=\frac{y_2}{s_2} = \cdots = \frac{y_n}{s_n},}

and given the hypothesis regarding the $s_i$'s there exists some $x\in E$ such that $xs_i=y_i$ for each $i$. 
\\
In the \gnl case we do just about the same thing.
\\
For each pair $(i,j)$  with $i\neq j$, the fact that $x_i/1=x_j/1$ in $E_{ij}$
means that for certain $u_{ij}\in S_i$ and  $u_{ji}\in S_j$ we have $s_j u_{ij} u_{ji} y_i = s_i u_{ij} u_{ji} y_j $.
For each~$i$, let $u_i\in S_i$ be a common multiple of the $u_{ik}$'s (for $k\neq i$). 
\\
We then have $(s_j u_{j}) (u_{i}  y_i) = (s_i u_{i}) (u_{j} y_j) $.
Thus the vector of the $u_{i}  y_i$'s is proportional to the vector of the $s_i u_{i}$'s.
Since the \sys $(\uS)$ is twice \Ergz, 
there exists \hbox{some $x\in E$} such that $u_{i}  y_i= s_i u_{i} x$ for every $i$, which gives the \egtsz~\hbox{$\varphi_i(x)=\fraC{u_{i}  y_i}{s_i u_{i}}=\fraC{  y_i}{s_i} =x_i$}.
\\
Finally, this $x$ is unique because the $S_i$'s are $E$-\corz.
\end{proof}

Now here is a variant of the \plgref{plcc.modules 2}.
This variant appears this time as a converse of the previous \plgz.

\begin{plcc}
\label{plcc.modules2bis} {\em (Concrete patching of modules in \prof 2) }
Let $(\uS)=(S_1, \dots, S_n)$ be a \sys of \mos of $\gA$. 
Let {\mathrigid 1mu$\gA_i=\gA_{S_i}$, $\gA_{ij}=\gA_{S_iS_j}$ and $\gA_{ijk}=\gA_{S_iS_jS_k}$}.
Suppose that a commutative diagram 

\snic{\big((E_i)_{i\in \lrbn}),(E_{ij})_{i<j\in \lrbn},(E_{ijk})_{i<j<k\in \lrbn};(\varphi_{ij})_{i\neq j},(\varphi_{ijk})_{i< j,i\neq k,j\neq k}\big)}

(as in the figure below) is given in the category of  \Amos with the following \prtsz. 
\begin{itemize}
\item For all $i$, $j$, $k$ (with $i<j<k$), $E_i$ is an $\gA_i$-module,  $E_{ij}$ is an~$\gA_{ij}$-module and $E_{ijk}$ is an~$\gA_{ijk}$-module.
Recall that according to our conventions of notation we let $E_{ji}=E_{ij}$, $E_{ijk}=E_{ikj}=\dots$

\item For $i\neq j$,  $\varphi_{ij}:E_i\to E_{ij}$ is a \molo at $S_j$ (see in $\gA_i$).
%
\item For $i\neq k$, $j\neq k$ and $i<j$, $\varphi_{ijk}:E_{ij}\to E_{ijk}$ is a \molo at $S_k$ (seen in $\gA_{ij}$).
\end{itemize}

\smallskip {\small\hspace*{6em}{
$
\xymatrix @R=2em @C=7em{
 E _i\ar[d]_{\varphi _{ij}}\ar@/-0.75cm/[dr]^{\varphi _{ik}} &
     E _j\ar@/-1cm/[dl]^{\varphi _{ji}}\ar@/-1cm/[dr]_{\varphi _{jk}} &
        E _k\ar@/-0.75cm/[dl]_{\varphi _{ki}}\ar[d]^{\varphi _{kj}} &
\\
 E _{ij} \ar[rd]_{\varphi _{ijk}} & 
    E _{ik}  \ar[d]^{\varphi _{ikj}} & 
      E _{jk}  \ar[ld]^{\varphi _{jki}} 
\\
   &  E _{ijk} 
}
$
}}

\vspace{-.3em}

 Then, if $\big(E,(\varphi_i)_{i\in\lrbn}\big)$ is the limit of the diagram, we have the following results.
\begin{enumerate}
\item 
Each morphism $\varphi_i:E\to E_i$ is a \molo at~$S_i$.  
\item 
The \sys $(\uS)$ is twice \Ergz.  
\item  The \sys $\big(E,(\varphi_{i})_{i\in\lrbn}\big)$ is, up to unique \isoz, the unique \sys  $\big(F,(\psi_{i})_{i\in\lrbn}\big)$ with the $\psi_i\in\Lin_\gA(F,E_i)$ satisfying the following items:
\begin{itemize}
\item the diagram is commutative,
\item each $\psi_i$ is a \molo at $S_i$,
\item the \sys $(\uS)$ is twice \Frgz.
\end{itemize}
\end{enumerate}
\end{plcc}
%
\begin{proof} \emph{1.} This \prt is valid with no hypothesis on the considered \sys of \mos (see the \dem of the \plgref{plcc.modules 2}).

\emph{2.}  Consider some $s_i\in S_i$ and a sequence $(\bmx_i)_{i\in\lrbn}$ in $E$ proportional to $(s_i)_{i\in\lrbn}$. Let $\bmx_i=(x_{i1},\dots,x_{in})$. The proportionality of the two sequences means that $s_ix_{jk}=s_jx_{ik}$ in $E_k$ for all $i$, $j$, $k$. Let $\bmx=(\fraC{x_{ii}}{s_i})_{i\in\lrbn}$. Next we prove that $s_i\bmx=\bmx_{i}$, \cad $s_i\fraC{x_{jj}}{s_j}=x_{ij}$ in each $E_j$. Indeed, this results from the \egt of proportionality  $s_ix_{jk}=s_jx_{ik}$ for $k=j$.

\emph{3.} 
Since~$E$ is the limit of the diagram, there is a unique \Ali  $\psi:F\to E$ such that
$\psi_i=\varphi_i\circ \psi$ for all~$i$.
\\ 
Actually we have $\psi(y)=\big(\psi_1(y),\dots,\psi_n(y)\big)$.
\\ 
Let us first show that $\psi$ is injective. \hbox{If $\psi(y)=0$}, all the $\psi_i(y)$'s are null, and since $\psi_i$ is a \molo at~$S_i$, there exist some $s_i\in S_i$ such that $s_iy=0$. 
Since $(\uS)$ is an \Frg \sysz, we have $y=0$. 
\\
As $\psi$ is injective we can suppose that $F\subseteq E$ and $\psi_i=\varphi_i\frt F$. \\
In this case showing that $\psi$ is bijective comes down to showing that $F=E$. 
Let $\bmx\in E$. As~$\psi_i$ and~$\varphi_i$
are two  \molos at $S_i$, there are $u_i\in S_i$ such that $u_i\bmx\in F$. Since $(\uS)$ is twice \Frgz,
and since the sequence of the $u_i\bmx$'s is proportional to the sequence of the $u_i$'s, there exists some $y\in F$ such \hbox{that $u_i\bmx=u_iy$} for every $i$, so $y=\bmx\in F$.
\end{proof}

\vspace{-.7em}
\pagebreak

\Exercices

\begin{exercise}
\label{exoLocMon1}
{\rm  
Let $S_1$, \ldots, $S_n$, $S$ be \mos of $\gA$ such that $S$ is contained in the saturated \mo of each $S_i$. 
\Propeq

 \emph{1.} 
The $S_i$'s cover $S$.

 \emph{2.} 
The $S_i$'s are \com in $\gA_S$.
}
\end{exercise}

\vspace{-1em}
\begin{exercise}
\label{exoLocMon2}
{\rm Let $I$ be an \id and $U$ be a \mo of $\gA$.
Let $S = \cS(I , U)$. 

 \emph{1.} 
 In $\gA_S$, the \mo $\cS(I; U,a)$ is \eqv to $\cS(I ; a) = I + a^\NN$.

 \emph{2.} 
 In $\gA_S$, the \mo $\cS(I,a ; U)$  is \eqv to $\cS(a ; 1) = 1 +\lra a$.
 }
\end{exercise}

\vspace{-1em}
\begin{exercise}
\label{exoLocMon3}
{\rm  Give a \dem of Lemma~\ref{lemRecouvre} based on the previous two exercises.}
\end{exercise}

\vspace{-1em}
\begin{exercise}
 \label{exoComCom}  {\rm Let $\gA$ be a \riz.
 
\emph{1.} For  $a_1$, \ldots, $a_n$ in
$\gA$, if $a_1 \cdots a_n \in \Rad\gA$, the \mos $1 + \gen {a_i}$ are \comz.

 \emph{2.}
 If $\fa_1$, \ldots, $\fa_\ell$ are \ids of $\gA$, the \mosz~\hbox{$1+\fa_i$} cover the \moz~\hbox{$1+\prod_i\fa_i$}.
} 
 \end{exercise}

\vspace{-1em}
\begin{exercise}
\label{exoIdepDyna} (In accordance with the \dfn of the \idepsz, if a product of factors is in a potential \idepz, 
we can open branches of computation
in each of which 
at least one of the factors is in the new potential \idepz)

\noindent 
{\rm We reuse the notations of \Dfn \ref{nota mopf}. 
Consider two subsets $I$ and~$U$ of~$\gA$ and the corresponding \mo $\cS(I,U)$. Let $a_1$, \dots, $a_k\in\gA$ 
for which we have

\snic{\prod_{i=1}^ka_i\in \gen{I}_{\gA_{\cS(I,U)}}.}

\emph{1.} Show that the \mos $\cS(I\cup\so{a_i},U)$ cover the \mo $\cS(I,U)$.

\emph{2.}  If we have $a_i-a_{j}\in \cS(I,U)$, then  $a_{j}$  is \iv in $\gA_{\cS(I\cup\so{a_i},U)}$.
 
\emph{3.}  Suppose that for each $j\in\lrbk$, we have an \auto of $\gA$ that fixes the \mo $\sat{\cS(I,U)}$ and that sends $a_1$ to $a_j$, and that each of the $\gA_{\cS(I\cup\so{a_i},U)}$'s is trivial, then
$\gA_{\cS(I,U)}$ is trivial. 
}
\end{exercise}

\vspace{-1em}
\begin{exercise}\label{exoMonoidesComax1}
{\rm
Let $S = (S_1, \ldots, S_n)$ be a family of \mos of $\gA$.

  \emph{1.}
Consider the family $S'$ obtained from $S$ by repeating each $S_i$ a certain number of times (at least once)
$$\preskip.2em \postskip.4em
\;\;\;S' = (S_1, S_1, \ldots, S_2, S_2, \ldots, S_n, S_n, \ldots)
$$
Show that $S$ is a family of \moco of $\gA$ if and only if the same goes for $S'$.

  \emph{2.}
Consider a second family $U = (U_1, \ldots, U_m)$ of \mos of $\gA$. Suppose that for each $i \in\lrb{1..n}$, there exists a $j \in\lrb{1..m}$ such that $S_i \subseteq U_j$ and for each $j \in \lrb{1..m}$ there exists an $i \in\lrb{1..n}$ such that $U_j \supseteq S_i$. Show that if $U$ is a family of \moco of $\gA$, the same goes for $S$.
}
\end{exercise}

\vspace{-1.1em}
\pagebreak

\begin{exercise}
\label{exoKroLocvar} (Variation on the local \KRNz's \thoz, \paref{thKroLoc})\\
{\rm To solve the exercise, we observe that the desired result is a \gui{quasi-global} statement that we can obtain by rereading the \dem of the local \KRNz's \thoz.
\\
 Let $x_0$, \dots, $x_d\in\gA$ and $\fa = \DA(x_0,\dots,x_d)$. If $\Kdim \gA \leq d$ and $\Kdim \gA/\fa \leq 0$,
there exist some \elts $s_0$, \dots, $s_d \in \gA$ and some \ids $\fb_0$, \dots, $\fb_d$, each generated by~$d$ \eltsz, such that\footnote{f.s.o.i.: \sfioz.}

\snic{
(s_0, \dots, s_d) \, \hbox{ is a f.s.o.i. of } \gA/\fa
\;\;\;  \hbox{and} \;\;\;
\forall i,~~s_i \fa \subseteq \sqrt{\fb_i} \subseteq \fa}

(locally, $\fa$ is maximal and radically generated by $d$~\eltsz).
}
\end{exercise}

\vspace{-1em}
\begin{exercise}
\label{exoKroLocvarbis}
(Second variation on the local \KRNz's \thoz)\\
{\rm
Let $\gA$ be a \ri and $\fa$ be a \tf  \idz.
If $\Kdim \gA/\fa \leq 0$ and $\Kdim \gA_{1+\fa}\leq d$, %
there exist some \elts $s_0,$ \dots, $s_d \in \gA$ and some \ids $\fb_0$, \dots, $\fb_d \subseteq \fa$, each generated by~$d$ \eltsz, such that

\snic{
(s_0, \dots, s_d)\, \hbox{ is a s.f.i.o. of } \gA/\fa
\;\;\;  \hbox{and} \;\;\;
\forall i,~~s_i \fa \subseteq \sqrt{\fb_i}. }
 }
\end{exercise}

\vspace{-1em}
\begin{exercise}
\label{exoFonctDet2}
{\rm
Given the \prcc of the modules (\plgc \ref{plcc.modules 2}), and given the canonical \iso 

\snic{\big(\Lin_\gA(M,N)\big)_S\rightarrow \Lin_{\gA_S}(M_S,N_S)}

in the case of \mpfs (Proposition~\ref{fact.homom loc pf}), we have local \carns for the \deter of a \mptf and that of a \homo between \mptfs (see Exercise~\ref{exoFonctDet1}).

 \emph{1.} 
The module $\det(M)$ is \care up to unique \iso by the following \prtz: if $s\in \gA$ is such that $M_s$ is free, then $\det(M)_s\simeq\det(M_s)$, with compatible \isos when we make a more advanced \lonz.{\footnote{This  \prmt means: if $s''=ss'$, then the \iso $(\det(M))_{s''}\simeq\det(M_{s''})$ is given by the \lon of the \iso $(\det(M))_s\simeq\det(M_s)$.}}

 \emph{2.} 
 If $\varphi ~:M\rightarrow N$ is a \homo of  $\gA$-\mptfsz, the \homo $\det(\varphi)$  is \care by the following \prtz: if $s\in \gA$ is such that $M_s$ and $N_s$ are free, then $\det(\varphi)_s=\det(\varphi_s)$ (modulo the canonical \isosz).
}
\end{exercise}

\vspace{-1em}
\begin{exercise}
\label{exoRecolle2PTF} {\rm Let $n\geq3$,   $s_1$, $s_2$ be two \eco of $\gA$.
We propose to concretely patch two \mptfs $P_1$ and $P_2$ respectively defined over $\gA_{s_1}$ and $\gA_{s_2}$ which have \isoc extensions to $\gA_{s_1s_2}$. By using the \dlg lemma, we can suppose that they are images of  \mprns $F_1$ and $F_2$ over $\gA_{s_1s_2}$ conjugated by means of a product of \elr matrices.

 \emph{1.} Let $E\in\En(\gA_{s_1s_2})$.
Show that there exists an $E_1\in\En(\gA_{s_1})$ and $E_2\in\En(\gA_{s_2})$ such that $E=E_1E_2$ over $\gA_{s_1s_2}$.

 \emph{2.} Let $F_1\in\Mn(\gA_{s_1})$ and $F_2\in\Mn(\gA_{s_2})$ be two  \mprns over $\gA_{s_1s_2}$ conjugated by means of a matrix $E\in\En(\gA_{s_1s_2})$. What to do?
}
\end{exercise}

\vspace{-1em}
\begin{exercise}
\label{exoMcCoyContr1} {(Contrapositive McCoy's \thoz, distressing version)}\\
{\rm 
Let $\gA$ be a nontrivial discrete \ri and $M\in\Ae{m\times n}$ be a matrix.
\begin{enumerate}
\item If $\cD_n(M)$ is faithful, $M$ is injective.
\item If we know an integer $k<n$ and some nonzero $x\in\gA$, such that 

\snic{x\cD_{k+1}(M)=0\hbox{  and  }\cD_k(M)\hbox{  is faithful},}

\snii then we can construct a nonzero vector in the kernel of $M$.  
\end{enumerate}
 
}
\end{exercise}

\vspace{-1em}
\begin{exercise}
\label{exoMcCoyContr2} (Contrapositive McCoy's \thoz, digestible version)\\
{\rm Let $\gA$ be a nontrivial discrete \coh \ri and $M\in\Ae{m\times n}$ be a matrix.
\begin{enumerate}
\item Either $\cD_n(M)$ is faithful, and $M$ is injective.
\item Or we can construct in the kernel of $M$ a vector with at least a \coo in~$\Atl$.  
\end{enumerate} 
}
\end{exercise}

%

\vspace{-1em}
\begin{exercise}
\label{exolemtrickprof2}
{\rm We notice that some \dfns of  depth $1$ and of  depth~2 are given in terms that do not make the additive structure of the considered \ri intervene, but only its
multiplicative structure, \cad the \mo $(\gA,\times,1)$.  

As the statement of Lemma~\ref{lemtrickprof2} does not make use of the additive structure either, we can hope for a purely multiplicative \dem of this lemma. The inspection of the \dem given in the course shows that this is not the case. We therefore propose to the reader to find a \dem of Lemma~\ref{lemtrickprof2} which works for any \moz.
 
}
\end{exercise}


\vspace{-1em}
\begin{problem}\label{exoChasserIdeauxPremiers1}
{(Avoiding the prime \idsz)}\\
{\rm  
 In this \pbz, we examine how to \cot decrypt a classical \dem that uses as a basis tool \gui{go see what happens in the fields $\Frac(\gA\sur\fp\!)$ for all the \ideps $\fp$ of $\gA$.}

  \emph{1.}
Let $t$ be an \idtr and $a$, $b$, $c_1$, \ldots, $c_n \in \gA$ such that $(at+b, c_1, \ldots, c_n)$ is a \vmd over $\gA[t,t^{-1}]$. We want to show that $ab \in \DA(c_1, \ldots, c_n)$. The following \demz, typical in \clamaz, uses \TEM and the axiom of choice.  
If $ab \notin \DA(c_1, \ldots, c_n)$, there exists a prime \id $\fp$ with $c_i \in \fp$ for $i \in \lrbn$ and $ab \notin \fp$. Over the field $\gK = \Frac(\gA\sur\fp\!)$, since $\ov a$ is nonzero, the \eqn $at + b = 0$ has a unique solution $t = -\ov b/\ov a$, which is nonzero because $\ov b$ is nonzero; we can then define a morphism $\varphi : \gA[t,t^{-1}] \to \gK$ by $t \mapsto -\ov b/\ov a$; $\varphi$ transforms the \vmd $(at+b, c_1t, \ldots, c_nt)$ into the null vector of $\gK^{n+1}$. A contradiction.
\\
What do you think?

 \emph{2.} If $\gB$ is a reduced \ri describe the units of $\gB[t,1/t]$. 
\\
We will be able to show that if $p$, $q\in\gB[t]$ satisfies $pq=t^m$ (with $p=\sum_k p_kt^k$ and~\hbox{$q=\sum_k q_kt^k$}), then $1 \in \rc(p)$,  $1 \in \rc(q)$ and

\snic {
p_kp_\ell=q_kq_\ell=0 \hbox { if } k\neq\ell,\quad
p_kq_\ell=0 \hbox { if } k+\ell\neq m, \quad
p_i = q_i=0 \hbox { if } i>m.
}

Consequently, for $k\in\lrb{0..m}$, $\gen {p_k}$ is generated by some \idm $e_k$. We then have a \sfio available $(e_0,\ldots,e_m)$ in $\gB$ such that $\gen{e_k}=\gen{p_k}=\gen{q_{m-k}}$ for $k\in\lrb{0..m}$ and

\snic{e_kp=e_kp_kt^k, \; e_kq=e_kq_{m-k}t^{m-k}\;$ and $\;e_k=e_kp_kq_{m-k}.}

The result is clear when the \ri is integral, so the reflex in \clama is to use some \idepsz.
A possible solution to \cot decrypt this reasoning is to use the formal \nst (\thref{thNSTsurZ}).

}
\end{problem}


\sol

\exer{exoLocMon1}~\\
$\emph {2} \Rightarrow \emph {1.}$ Let $s_1$, \ldots, $s_n$ with
$s_i \in S_i$. We want $b_i$'s $\in \gA$ such that~\hbox{$b_1s_1 + \cdots + b_ns_n
\in S$}. The fact that the $S_i$'s are \com in $\gA_S$ provides an $s \in S$ and some $a_i$'s $\in \gA$ such that $a_1s_1 + \cdots + a_ns_n = s$ in $\gA_S$; therefore there exists a  
 $t \in S$ such that, in $\gA$, $(ta_1)s_1 + \cdots + (ta_n)s_n = ts \in S$.

\exer{exoLocMon2}~\\
\emph{1.} An \elt $s$ of $\cS(I; U,a)$
is of the form $x + ua^k$. 
We see that $s$ divides the \elt $xu^{-1} + a^k \in \cS(I ; a)$ in~$\gA_S$.

\emph{2.} An \elt $s$ of $\cS(I,a ; U)$ is of the form $x + ya + u$. Let $x' =u^{-1}x$ %
and $y' = u^{-1}$.
Then $s$ divides in $\gA_S$ the \elt  $x' + y'a + 1$, which divides $1 + y''a$,  %
where $y'' = (1 + x')^{-1} y'$.

\exer{exoLocMon3}
Let $S = \cS(I,U)$, $S_1 = \cS(I; U,a)$ and $S_2 = \cS(I,a ; U)$.  
We must prove that $S_1$ and $S_2$ are \com in $\gA_S$.  In $\gA_S$, $S_1$ is \eqv to $I + a^\NN$, and~$S_2$ is \eqv  to $1 + \langle a\rangle$.  Let us use the following \idt

\snic{y^k(x + a^k) + \bigl(\som_{j < k} y^j a^j\bigr)(1 - ya) =
y^k(x + a^k) + 1 - y^ka^k = 1 + y^kx.}

Applied to $x \in I$, it proves that $x + a^k$ and $1 - ya$ are \com in $\gA_S$ (since $1 + y^kx \in 1 + I$
and since $I$ is contained in the radical of $\gA_S$).


\exer{exoComCom} \emph{1.} For $j\in\lrbn$ let $b_j=1-a_jx_j$ in the \moz~\hbox{$1+a_j\gA$}. Let~\hbox{$a=\prod_ia_i$}.
We must show that the \id $\fm=\gen{b_1,\dots,b_n}$ contains $1$.
\\
However, modulo  $\fm$  we have $1=a_jx_j$, therefore $1=a\prod_ix_i=ax$. Thus $1-ax\in\fm$, \hbox{but $1-ax\in\Ati$} because $a\in\Rad\gA$.

\emph{2.}
It is clear that $S=1+\prod_i\fa_i\subseteq1+\fa_j=S_j $ for each $j$. We must therefore prove (Exercise \ref{exoLocMon1}) that the $1+\fa_j$ given in $\gA_S$ are \comz. However, the product $\prod_i\fa_i$, seen in $\gA_S$, is in $\Rad\gA_S$.
Therefore it suffices to apply item \emph{1.}



\exer{exoIdepDyna} The hypothesis means that we have 
          a $u\in\cM(U)$ 
and a $j\in\gen{I}_\gA$ such that 
%
$(u+j)\,\prod_{i=1}^ka_i\in\gen{I}_\gA$, 
%
 or $u\,\prod_{i=1}^ka_i\in\gen{I}_\gA$.

\sni
\emph{1.} 
\emph{First solution, by direct computation.}

\noindent Consider $x_i\in S_i=\cS(I\cup\so{a_i},U)$ and look for a \coli which is in  $S=\cS(I,U)$.
For each $i$ we write
$$\preskip.2em \postskip.1em
{x_i=u_i+j_i+a_iz_i\hbox{  with  }u_i\in \cM(U),\;j_i\in\gen{I}_\gA\hbox{ and }z_i\in\gA.}
$$
In the product
$$\preskip-.2em \postskip.2em\ndsp
{u\,\prod_{i=1}^k(x_i-(u_i+j_i))=u\;\prod_{i=1}^ka_iz_i\in\gen{I}_\gA ,}
$$
we re-express the left-hand side in the form
$$\preskip.3em \postskip.2em\ndsp
{\som_{i=1}^kc_ix_i\pm u\,\prod_{i=1}^k(u_i+j_i)}
$$
and we obtain, by moving $\pm u\,\prod_{i=1}^k(u_i+j_i)$ to the right-hand side, the desired membership $\som_{i=1}^kc_ix_i\in \cS(I,U)$.

\emph{Second solution, conceptual.}
\\ 
It is clear that $S \subseteq S_i$. It therefore suffices (Exercise~\ref{exoLocMon1}) to show that the  $S_i$'s are \com in $\gA_S$.
In $\gA_S$, (Exercise~\ref{exoLocMon2}) the \mos $S_i$ and $1 + \gen {a_i}$ have the same saturated \moz.
 It therefore suffices to see that the \mos $1 + \gen {a_i}$ are comaximal in $\gA_S$.
Moreover, we know that, seen in $\gA_S$, $I$ is contained in $\Rad(\gA_S)$. We therefore apply item \emph{1} of Exercise~\ref{exoComCom}.

\emph{2.} Clear since $a_{j}\in \cS(I,U)+\gen{a_i}\subseteq \cS(I\cup\so{a_i},U)$.

\emph{3.} If one of the $\gA_{S_i}$'s is trivial, all of them are also trivial because they are pairwise \isocz. Since the $S_i$'s cover $S$, $\gA_S$ is itself trivial.

\exer{exoMonoidesComax1}
 \emph{1.}
It suffices to show it for $S' = (S_1, S_1, S_2, \ldots, S_n)$.  \\
Suppose that the family $S$ is comaximal. Let $s'_1$, $s''_1 \in S_1$, and $s_i \in S_i$ for $i \in \lrb{2..n}$. The \elts $s'_1s''_1$, $s_2$, \ldots, $s_n$ are \comz, and since $s'_1s''_1 \in \gen {s'_1, s''_1}$, the same goes for $s'_1$, $s''_1$, $s_2$, $\ldots$, $s_n$. In the other direction, suppose that $S'$ is comaximal and let $s_i \in S_i$ for $i \in \lrb{1..n}$; then $s_1$, $s_1$, $s_2$, $\ldots$, $s_n$ are comaximal therefore the same goes for $s_1$, $s_2$, $\ldots$, $s_n$.

  \emph{2.}
By repeating some of the $S_i$'s and the $U_j$'s, we obtain two families $S'$, $U'$ of \mos of $\gA$, indexed by the same interval $\lrb{1..p}$ and satisfying $S'_k \subseteq U'_k$ %
for $k \in \lrb {1..p}$. Since $U$ is comaximal, the same goes for $U'$ so for $S'$ then for~$S$.


\exer{exoKroLocvar}
As $\Kdim \gA\leq d$, there exists a sequence $(\underline{a}) = (a_0,\dots,a_d)$ \cop to $(\underline{x}) = (x_0,\dots,x_d)$.
Therefore (disjoint sequences), for every $i \leq d$, we have

\vspace{-.3em}\snic{\rD(a_0, \dots, a_{i-1}, x_0, \dots, x_{i-1}, a_i x_i)
= \rD(a_0 + x_0, \dots, a_{i-1} + x_{i-1}).}

As $\Kdim \gA/\fa \leq 0$ and $\fa=\DA(\fa)$, we also have $\gA = \gA a_i + (\fa : a_i)$ for all~$i$. We then construct the triangle

\snuc{\arraycolsep2pt\begin{array}{rl}
& \gA \\
  = & \gA a_0 + (\fa : a_0) \\
  = & \gA a_0 + (\fa : a_0)a_1 + (\fa : a_0)(\fa : a_1) \\
  & \qquad \vdots \\
  = & \gA a_0 + (\fa : a_0)a_1 + \cdots + (\fa : a_0)\cdots(\fa : a_{d-1}) a_d
 + (\fa : a_0)\cdots(\fa : a_d). 
\end{array}
}

Now, we write
$$\preskip.0em \postskip.4em
1 = b_0 a_0 + b_1 a_1 + \cdots + b_d a_d + t
$$
with $b_i \in (\fa : a_0)\cdots(\fa : a_{i-1})$ and $t \in (\fa : a_0)\cdots(\fa : a_d)$. For $i \leq d$, on the one hand we have

\snic{
\begin{array}{c}
b_i \gen{x_0, \dots, x_{i-1}}
\subseteq \rD\big( b_i(a_0+x_0), \dots, b_i(a_{i-1}+x_{i-1}) \big)
  \hspace*{3cm}    \\[1mm] \hspace*{3cm}
\subseteq \rD(  b_ia_0,x_0, \dots, b_ia_{i-1},x_{i-1} ) \subseteq \rD( \fa ) =\fa,
\end{array}
}

and on the other hand

\snic{
\begin{array}{c}
b_i a_i x_i \in b_i \rD(a_0+x_0,\dots, a_{i-1}+x_{i-1})
  \hspace*{5cm}    \\[1mm] \hspace*{4cm}
\subseteq \rD\big( b_i(a_0+x_0), \dots, b_i(a_{i-1}+x_{i-1}) \big)
\subseteq \fa.
\end{array}
}

Now let $s_i = b_i a_i$. Thus we reach
$$\preskip.2em \postskip.2em
{
s_i \gen{x_0, \dots, x_{i-1}, x_i}
\subseteq \rD\big(b_i(a_0+x_0),\dots, b_i(a_{i-1}+x_{i-1}) \big) \subseteq \fa
,}
$$
then
$$\preskip.2em \postskip.2em\mathrigid2mu
\begin{array}{c}
s_i \gen{x_0, \dots, x_{i-1}, x_i, x_{i+1}, \dots, x_d}
\subseteq
  \hspace*{6cm}    \\[1mm] \hspace*{3cm}
\rD\big( \underbrace{b_i(a_0+x_0), \dots, b_i(a_{i-1}+x_{i-1}),
     s_i x_{i+1}, \dots, s_i x_d}_{\hbox{\scriptsize generate }\fb_i
\hbox{\scriptsize ~(def.) }} \big)
\subseteq \fa.
\end{array}
$$
Therefore there exists an \id $\fb_i$ generated by $d$ \elts satisfying

\snic{
s_i \fa \subseteq \rD(s_i \fa) \subseteq \rD( \fb_i ) \subseteq \fa.}

We end the \dem by using $1 \in \gen{a_d,x_d}$. We get

\snic{
t \in t \gen{a_d, x_d} \subseteq \gen{ta_d, x_d} \subseteq \fa,
}

so much so that the sum of the $s_i$'s is equal to $1 \mod \fa$.
Moreover, for $i>j$,

\snic{
s_i s_j \in b_i a_j \gA \subseteq \fa,
}

which allows us to conclude that $(s_0,\dots,s_d)$ is a \sfio of~$\gA/\fa$.

\exer{exoKroLocvarbis}
To begin with, \KRNz's \tho gives that $\sqrt{\fa}$ is radically generated by $d+1$ \eltsz.
Next we apply the result of Exercise \ref{exoKroLocvar} to $\sqrt{\fa}$ in the localized \ri $(1+\fa)^{-1} \gA$.
This provides some~$s_i$'s forming a \sfio modulo~$\sqrt{\fa}$ and some~$\fb_i \subseteq \sqrt{\fa}$.
Even if it entails taking multiples of powers of the~$s_i$'s, we can impose that $(s_0, \dots, s_d)$ is a \sfio of~$\gA/\fa$.
Even if it entails taking powers of the \gtrs of the~$\fb_i$'s, we can impose $\fb_i \subseteq \fa$.

\exer{exoRecolle2PTF} \emph{1.} \cite[page 208 Proposition 1.14]{Lam06}.

\emph{2.} We have $F_1E=EF_2$. We write $E=E_1E_2$, so $F_1E_1E_2=E_1E_2F_2$
 and 

\snic{\wi{E_1}F_1E_1=_{\gA_{s_1s_2}}E_2F_2\wi{E_2}.}

The matrix $\wi{E_1}F_1E_1$ (resp.\ $E_2F_2\wi{E_2}$) is a \mprn over $\gA_{s_1}$ (resp.\ over~$\gA_{s_2}$) because $\wi{E_1}E_1=\In$ (resp.\ $\wi{E_2}E_2=\In$). By the \plg of \rcm of the \elts in a module (here $\Mn(\gA)$),
there exists a unique matrix~\hbox{$F\in\Mn(\gA)$} which is equal to $\wi{E_1}F_1E_1$ over $\gA_{s_1}$ and to $E_2F_2\wi{E_2}$ over $\gA_{s_2}$.
To prove that $F^2=F$, it suffices to prove it over $\gA_{s_1}$ and $\gA_{s_2}$. Let $P=\Im F\subseteq\Ae n$.
\\
By construction, for $i=1,2$

\snic{P_{s_i}= _{\gA_{s_i}^n}\Im F_{s_i}\simeq _{\gA_{s_i}^n}\Im F_i\simeq _{\gA_{s_i}^n}P_i.}

\exer{exoMcCoyContr1} \emph{(Contrapositive McCoy's \thoz, distressing version)}\\
\emph{1.} Already seen.\\
\emph{2.}  We have $x\neq 0$,  $\cD_k(M)$ is faithful and the \ri is discrete, therefore there exists a minor~$\mu$ of order $k$ of $M$ such that $x\mu\neq 0$.
Suppose for example that $\mu$ is the north-west minor and let $C_1$, \dots, $C_{k+1}$ be the first columns of $M$, let~$\mu_i$ ($i\in\lrbk$) be the suitably signed \deters of the matrices extracted on the rows $\lrbk$ and the previous columns, except the column of index $i+1$. Then the Cramer formulas give the \egt $\sum_{i=1}^{k}x \mu_iC_i+x\mu C_{k+1}=0$. As $x\mu\neq 0$, this gives a nonzero vector in the kernel of $M$.

\comm In \clamaz, if $\cD_n(M)$ is not faithful, as $\cD_0(M)=\gen{1}$ is faithful, there exists some $k<n$ such that $\cD_k(M)$ is faithful and~$\cD_{k+1}(M)$ is not faithful. Still in \clamaz, if $\cD_{k+1}(M)$ is not faithful, there exists some $x\neq 0$ such that $x\cD_{k+1}(M)=0$.
For these things to become explicit, we need for example to dispose of a test for the faithulness of the \itfsz, in a very strong sense. 
\eoe


\exer{exoMcCoyContr2}\emph{(Contrapositive McCoy's \thoz, digestible version)}\\
Since the \idds are \itfsz, and the \ri is \cohz, their annihilators are \egmt \itfsz, and we can test the nullity of a \itf because the \ri is discrete. The hypotheses of Exercise~\ref{exoMcCoyContr1} are therefore satisfied.\\
Note: The alternative \gui{\emph{1} or \emph{2}} is exclusive because the \ri is nonzero, this justifies the \gui{either, \dots, or} of the statement.







\prob{exoChasserIdeauxPremiers1}
\emph{1.}
There is no miracle: a certificate for $ab \in \DA(c_1, \ldots, c_n)$ can be obtained from a certificate of unimodularity of $(at+b, c_1t, \ldots, c_nt)$ in $\gA[t,t^{-1}]$.
\\ 
By replacing $\gA$ by $\gA_1=\gA\sur{\DA(c_1, \ldots, c_n)}$, we are brought back to $\gA$ being reduced \hbox{and $c_i = 0$}. The hypothesis is then $at+b$ is \iv in $\gA[t,t^{-1}]$, and the result to be shown is~\hbox{$ab = 0$} (\smq result in $a,b$ just like the hypothesis). \\
We have $(at+b)g(t) = t^e$ for some $g \in \gA[t]$ and some $e \in \NN$, therefore the \pol $at + b$ is primitive. To show $ab = 0$, it suffices to localize at $a$ then at~$b$. Over the localized \ri at $a$, 
we take $t = -b/a$ in $(at+b)g(t) = t^e$, we obtain $(-b/a)^e = 0$.
Thus, we~have~$b = 0$, then $ab = 0$. By symmetry, we obtain in $\gA_b$, $a = 0$ so $ab=0$. \\
Actually, if $ua+vb=1$, the \ri $\gA_1$ is split in two by the \idm $ua$.
 In the first component, $at+b=a$ with $a$ \ivz, in the second, $at+b=b$ with $b$ \ivz.
 
 \emph{2.} In \clamaz: if we pass to the quotient by a \idep the result is clear.
By continuity, the spectrum is partitioned into a finite number of open sets corresponding to the coveted \sfioz. 
\\
 A \cov \dem is given in \cite[Yengui]{Ye0}. The reader will also be able to draw from the \dem of item \emph{1}.
\\
 A method that we can systematically use consists in calling upon the formal \nst (\thref{thNSTsurZ}).
\\ 
In the current case, we note that the \pb comes down to proving that the $p_kp_\ell$'s are null for $k\neq\ell$ and that the $p_{m+r}$'s null for $r>0$. Once this is observed, since the $p_k$'s are \comz, we obtain a \sfio $(e_0,\ldots,e_m)$ such that $e_kp=e_kp_kt^k$ for all $k$, from which \trfz.
\\ 
The philosophy is the following: if we take all the \coes of the \pb as \idtrs over $\ZZ$, the hypothesis comes down to passing to the quotient by the radical $\fa$ of a \itfz, which represents the hypotheses. The goal is then to prove that the conclusions are \egmt in~$\fa$. For this it suffices to prove that it is indeed the case when we evaluate the \pb in an arbitrary finite field.
\\
 Here the \idtrs are $p_0$, \ldots, $p_n$, $q_0$, \ldots, $q_n$.\\
 For $p = \sum_{k=0}^n p_kt^k$ and $q = \sum_{k=0}^n q_kt^k$; we define the \pol
$$\ndsp
\preskip.4em \postskip.4em 
\sum r_jt^j\eqdefi pq-t^m
$$
(with all $r_j\in\ZZ[p_0,\ldots,p_n,q_0,\ldots,q_n]$)  and the \id $\fa$ is $\rD(r_0,\ldots,r_{2n})$. 
We will show that the $p_kp_\ell$'s and $q_kq_\ell$'s are in $\fa$ if $k\neq\ell$, 
that the $p_kq_\ell$'s are in $\fa$ for $k+\ell\neq m$ 
and that the $p_{m+r}$'s and $q_{m+r}$'s are in $\fa$ for $r>0$.
\\
 However, this directly results from item \emph{2} in the formal \nst (or then from item \emph{4} in Corollary~\ref{corthNSTsurZ}).
\\
 In \gmq terms: if $n\geq m$, the \vrt of the zeros of $pq-t^m$ over a field $\gK$ is a space formed of 
$m+1$ copies of $\gK\eti$ isolated from one another; over a reduced \ri the response is fundamentally the same, 
but the isolated components in the case of the fields here make a \sfio appear.

 Thus, the formal \nst (\thref{thNSTsurZ}) provides a \cov method to decrypt the hidden \algos in certain reasonings from \clamaz,  when the argument consists in seeing what happens in all the $\Frac(\gA\sur\fp\!)$'s for all the \ideps of $\gA$.


\Biblio

The dynamic method as it is explained in the \lgbe machinery with \ideps (Section~\ref{secMachLoGlo}) consists for the most part in flattening the computations implied by the \emph{method of dynamic evaluation} given in \cite[Lombardi]{Lom97}, a successor of the dynamic method implemented in \cite[Coste\&al.]{clr} for proofs of the \nst type, itself a successor of the dynamic evaluation \`a la D5 in Computer Algebra \cite[Duval\&al.]{D5}.\imlb
With respect to what is proposed in \cite{clr,Lom97}, the difference in the previous chapter is mainly that we have avoided the reference to formal logic.

In \clamaz, we find the \prcc of the \mptfs (item \emph{\ref{iptfplcc.ptf}} of the \plgref{plcc.ptf}) for example in \cite[Proposition~2.3.5 and Lemma~3.2.3]{Kni} (with an almost entirely \cov \demz) and in \cite[rule~1.14 of Chapter~IV]{Kun}.

The \cov treatment of \Susz's lemma \ref{lemSuslin1} is due to Ihsen Yengui \cite{Y1}, who gives the key to the \lgbe machinery with \idemasz.

The \lgbe machinery with \idemis is due to Thierry Coquand \cite[On seminormality]{coq}.

Exercises \ref{exoKroLocvar} and \ref{exoKroLocvarbis} are due to Lionel Ducos~\cite{Duc08}.\perso{penser \`a demander la r\'ef\'erence \`a Lionel}

The dynamic method was applied for the computation of \gui{dynamic Gr\"obner bases} by Yengui in \cite{Y3}.

\newpage \thispagestyle{CMcadreseul}
\incrementeexosetprob


\chapter[Extended projective modules]{Extended projective modules}
\label{ChapMPEtendus}
\perso{compil\'e le \today}
\vskip-1em

\minitoc

\subsection*{Introduction}
\addcontentsline{toc}{section}{Introduction}
In this chapter we \cot establish a few important results regarding the situations where the \mptfs over a \pol \ri are extended from the base \riz.

We especially treat Traverso-Swan's \tho (Section~\ref{sec.Traverso}), the patching \`a la Vaserstein-Quillen (Section~\ref{subsecQPatch}), Horrocks' \thos (Section~\ref{sec.Horrocks}), Quillen-\Susz's \tho (Section~\ref{sec.QS}), and in Section~\ref{sec.Etendus.Valuation}, Bass' \tho (\thref{thBass.Valuation}) and the Lequain-Simis \tho (\thref{thLSValu}). 


\section{Extended modules} \label{sec etendus}

Given an \alg $\gA\vers{\rho}\gB$, the \eds from $\gA$ to $\gB$ transforms a module $M$  over $\gA$ into a module
$\rho\ist(M)\simeq \gB\otimes_{\gA}M$ over $\gB$.
Recall that a \Bmo \isoc to such a module $\rho\ist(M)$ is said to be extended from~$\gA$. We also say that it comes from the \Amo $M$ by \edsz.

In the case of a \mpfz, from the point of view of the \mpnsz, this corresponds to considering the matrix transformed by the \homo $\rho$.

A \ncr condition for a \mpf to be extended is that its Fitting \ids are of the form $\rho(\fa_{i}) \gB$ for \itfs $\fa_{i}$ of~$\gA$.
This condition is realized for the \mptfs \ssi the \idms of $\gB$ are all images of \idms of~$\gA$.

\subsec{The \pb of the extension}

For the \mptfsz, the following \pb arises naturally given the morphism  $\GKO \rho:\GKO \gA\to \GKO \gB$.

\smallskip  \PB{\num1} Does every \mptf over $\gB$ come from a \mptf over $\gA$?
Or yet again: is $\GKO\rho$ surjective?

\medskip Recall that $\GKO \Ared=\GKO \gA$ and  $\GKO \gB\red=\GKO \gB$,
such that the \pb of the extension of the \mptfs can be narrowed down to the case of reduced \risz.
Moreover, if $\HO\rho:\HO\gA\to\HO\gB$ is not surjective, the answer to  \pb \num1 is negative \gui{for the wrong reason}
and the following \pb is then more natural.

\smallskip \PB{\num2}  Does every \mrc over $\gB$ come from a \mptf over $\gA$?

\medskip For the \mpfs the natural \gnn of the previous \pb is then the following.

\smallskip \PB{\num3} Does every \mpf  over $\gB$ whose \idfs are extensions of \itfs of $\gA$ come from a \mpf over $\gA$?

\subsec{The case of the \pol \risz}

Let $\gB=\AXr=\AuX$.
If $(\ua)\in\Ae r$ we denote by $\ev_\ua$ the evaluation \homo at $\ua$
$$\preskip-.2em \postskip.2em 
\ev_\ua\,:\,\gB\to\gA,\; p\mapsto p(\ua). 
$$
The two \homos $\gA\vers{j}\gB\vvers{\ev_\ua~}\gA$ are composed according to the identity.

Most of what follows in this subsection could be written in the more \gnl context of an \Alg $\gB$ possessing a character (cf. Proposition~\ref{prdfCaracAlg}).
For  \pol \ris we obtain the following results (with an obvious intuitive notation for $M(\uX)$).

\begin{fact}\label{factEtPol} With $\gB=\AuX$.
\begin{enumerate}
  \item A \Bmo  $M=M(\uX)$ is extended from $\gA$ \ssi it is \isoc to $M(\uze)$.
  \item In particular, if $M$ is \pf with a \mpn $G(\uX)\in\gB^{q\times m}$, Lemma \ref{lem pres equiv} implies that $M$ is extended from $\gA$ \ssi the matrices $H(\uX)$ and $H(\uze)$, where $H$ is illustrated below, are \eqves over the \ri $\gB$
$$\preskip.4em \postskip.6em
H(\uX)\;=\kern-7pt\raise7pt\hbox{$\begin{array}{c|p{35pt}|p{23pt}|p{23pt}|p{35pt}|}
\multicolumn{1}{c}{} & \multicolumn{1}{c}{$m$} & \multicolumn{1}{c}{$q$} &
\multicolumn{1}{c}{$q$} & \multicolumn{1}{c}{$m$} \\[2pt]
\cline{2-5}
\vrule height20pt depth13pt width0pt \rlap{\hskip6cm$q$}\; & \hfil $G(\underline{X})$  \hfil &\hfil
$0$\hfil &\hfil $0$\hfil &\hfil $0$\hfil \\
\cline{2-5}
\vrule height20pt depth13pt width0pt \rlap{\hskip6cm$q$}\; & \hfil $0$  &\hfil
${\rm I}_{q}$\hfil & \hfil $0$\hfil &\hfil$0$\hfil \\
\cline{2-5}
\end{array}$}
\quad \quad 
$$
\vskip0pt
\end{enumerate}
\end{fact}

\rem By Lemma~\ref{lem pres equiv} when the matrices $H(\uX)$ and $H(\uze)$ are \eqvesz, they are \elrt \eqvesz.
\eoe

\smallskip 
Regarding the \mptfs we obtain \homos of semi\ris which are composed 
according to the identity
$$\preskip.2em \postskip.4em
\GKO\gA\vvvers{\GKO j}\GKO\AuX\vvvvers{\GKO\ev_\ua}\GKO\gA.$$
Consequently $\GKO j$ is injective, and the phrase \gui{every \mptf over $\AuX$ is extended from $\gA$} means that $\GKO j$ is an \isoz, which we abbreviate to \gui{$\GKO\gA=\GKO\AuX$.}

Similarly, for the Grothendieck \ris
$$\preskip.2em \postskip.4em 
\KO\gA\vvers{\KO j}\KO\AuX\vvvvers{\KO\ev_\ua}\KO\gA,\quad \hbox{with}
\quad \KO(\ev_\ua)\circ \KO( j)=\Id_{\KO\gA}. 
$$

 Moreover we have the following \elr results, in which each \egt has the meaning that a natural morphism is an \isoz.

\begin{fact}\label{fact A->A[X]} With $\gB=\AuX$.
\begin{enumerate}
\item \label{i1fact A->A[X]}
  $\rD_{\gB}(0)=\DA(0)\gB$ (a \pol is nilpotent \ssi all its \coes are nilpotent). In particular, $\gB\red=\Ared[\uX]$.
\item \label{i2fact A->A[X]}
\label{item3fact A->A[X]} If $\gA$ is reduced, $\Bst=\Ati$. More \gnltz, $\Bst=\Ati+\DA(0)\gen{\uX}$.
\item \label{i3fact A->A[X]}
 $\BB(\gA)=\BB(\AuX)$ and $\HO \gA =\HO \AuX $.
\item \label{i4fact A->A[X]}  $\GKO \gA =\GKO \Ared $.
\item \label{i5fact A->A[X]}  $\GKO \gA =\GKO \AuX \iff \GKO \Ared =\GKO \Ared[\uX]$.%
\item \label{i6fact A->A[X]}   $\Pic \gA =\Pic \gB \iff \Pic \Ared =\Pic \Ared[\uX]$.
\end{enumerate}
\end{fact}
\begin{proof}
\emph{\ref{i1fact A->A[X]}} and \emph{\ref{i2fact A->A[X]}.} See Lemma \ref{lemGaussJoyal}.\\
\emph{\ref{i3fact A->A[X]}.} We must show that every \idm \pol is constant. This is done (in a single variable) by \recu on the formal degree of the \polz.\\
\emph{\ref{i4fact A->A[X]}.}  This is \thref{propComparRedRed}.\\
Items \emph{\ref{i5fact A->A[X]}}  and  \emph{\ref{i6fact A->A[X]}} result from items \emph{\ref{i1fact A->A[X]}} and \emph{\ref{i4fact A->A[X]}.}    
\end{proof}
%


\section[The Traverso-Swan's \thoz]{The Traverso-Swan's \thoz, seminormal \risz}
\label{sec.Traverso} 

This section is devoted to the study of the \ris $\gA$  for which the natural \homo from $\Pic\gA$ to $\Pic\AXr$ is an \iso (\cad the \mrcs $1$ over $\AXr$ are all extended from $\gA$).
The answer is given by the Traverso-Swan-Coquand \tho (\cite{Tra,Swan80,coq}):

\THo{(Traverso-Swan-Coquand)\\}
{
\Propeq
\begin{enumerate}
\item The \ri $\Ared$ is seminormal (\dfn \ref{defiseminormal}).
\item The natural \homo $\Pic\gA\vvers{\Pic j}\Pic\AX$ is an \isoz.
\item 
$\Tt r\geq1$, the natural \homo $\Pic\gA\to\Pic\AXr$ is an \isoz.
\item 
$\Ex r\geq1$, the natural \homo $\Pic\gA\to\Pic\AXr$ is an \isoz.
\end{enumerate}
}

\medskip
We will show \emph{1  $\Rightarrow$ 3} and \emph{2  $\Rightarrow$ 1}.
As a corollary, $\gA$ is seminormal \ssi $\AX$ is seminormal.

\subsec{Preliminaries}

First of all recall the following result (see Proposition~\ref{propImProjLib}).

\begin{lemma}
\label{lempropImProjLib}
 A \prn matrix of rank $1$, $P$, has a free image \ssi there exists a column vector $C$ and a row vector $L$ such that $LC=1$ and $CL=P$. In addition, $C$ and $L$  are unique, up to the product by a unit, 
under the only condition that $CL=P$.
\end{lemma}

%

Moreover recall that the natural morphism $\Pic\gA\to \Pic\AuX$ is an \iso \ssi the natural morphism  $\Pic\Ared\to \Pic\Ared[\uX]$ is an \iso (Fact~\ref{fact A->A[X]}~\emph{\ref{i6fact A->A[X]}}).

The two group \homos

\vspace{-.2em}\snic{\Pic\gA\vvvers{\Pic j}\Pic\AuX\vvvers{\Pic\ev_\uze}\Pic\gA} 
are composed 
according to the identity. The first is injective, the second surjective.
They are \isos \ssi the first is surjective, \ssi the second is injective. \perso{la structure of group n'intervient
pas in le raisonnement, cependant  l'injectivit\'e of the second
est plus facile \`a caract\'eriser in le cas des groups.}
\\
This last \prt means: every \idme square matrix $P(\uX)$ of rank $1$ over $\AuX$ which satisfies \gui{$\Im\big(P(\uze)\big)$ is free,} satisfies \gui{$\Im(P\big(\uX)\big)$ is free} in itself.
\\
Actually, if  $\Im\big(P(\uze)\big)$ is free,  the matrix 
$\Diag(P(\uze),0_1)$ is similar to a standard \prn matrix $\I_{1,n}=\Diag(1,0_{n-1,n-1})$ (\dlg lemma \ref{propIsoIm}). Hence the following lemma.

\begin{lemma}
\label{lemPicPic1} \Propeq
\begin{enumerate}
\item The natural \homo $\,\Pic\gA\to\Pic\AuX$ is an \isoz.
\item For every matrix  $M(\uX)=(m_{i,j})
\in\GAn(\AuX)$
such that $M(\uze)=\I_{1,n}$,
there exist
$f_1$, \ldots, $f_n$, $g_1$, \ldots, $g_n\in\AuX$ such that $m_{i,j}=f_ig_j$
for \hbox{all $i$, $j$}.
\end{enumerate}
\end{lemma}

Note that the hypothesis $M(\uze)=\I_{1,n}$ implies $\rg(M)=1$ because the \homo $\HO(\AuX)\to\HO(\gA)$ is an \isoz.

\mni{\bf Convention.}
\label{convPicPic}
We abbreviate the statement \gui{the natural morphism from~$\Pic\gA$ to~$\Pic\AuX$ is an \isoz} by writing:
 \gui{$\Pic\gA=\Pic\AuX$.}

\begin{lemma}\label{lemPicPic2}
Let $\gA\subseteq\gB$ be reduced \ris and  $f_1$, \ldots, $f_n$,  $g_1$, \ldots, $g_n$ be
\pols in $\BuX$ that satisfy the following \prts
$$\leqno{(*)} \quad\quad
\formule{f_1(\uze)=g_1(\uze)=1,  f_i(\uze)=g_i(\uze)=0 \;\;(i=2,\ldots,n),
\\[1mm]
m_{ij}\eqdefi f_ig_j \in\AuX  \;\;  (i,j=1,\ldots,n),
\\[1mm]
\sum_if_ig_i=1.}
$$
Under these hypotheses, the matrix $M:=(m_{ij})$ is a \mprn of rank $1$, $M(\uze)=\I_{1,n}$,
and \propeq
\begin{enumerate}
\item The module  $\Im M$ is free over $\AuX$, \cad extended from $\gA$.
\item The $f_i$'s and the $g_i$'s are in $\AuX$.
\item $f_1\in\AuX$.
\end{enumerate}
\end{lemma}

\begingroup\def\mou{}
\begin{proof}
\emph{3 $\Rightarrow$ 2.} The $g_j$'s are obtained from $f_1$ and from the $m_{1j}$'s by making divisions by non-decreasing powers, because the constant \coe of $f_1$ is equal to~$1$. Similarly, we then obtain the $f_i$'s from $g_1$ and from the $m_{i1}$'s.
The converse implication is trivial.\\
\emph{2 $\Leftrightarrow$ 1.} By Lemma~\ref{lempropImProjLib}, the \pb is to find suitable $f_i$'s and $g_j$'s from the matrix $(m_{ij})$.
However, these $f_i$'s and $g_j$'s exist in $\BuX$, and the condition $f_1(\uze)=1$ forces their uniqueness because the \ris are reduced (so the \ivs \elts in \pol \ri are constants).
\end{proof}

Lemmas \ref{lemPicPic1} and \ref{lemPicPic2} imply the following result.

\begin{corollary}\label{corlemPicPic2}
Let $\gA\subseteq\gB$ be two reduced \ris with  $\Pic\gB=\Pic\BuX$.  \Propeq
\begin{enumerate}
\item $\Pic\gA=\Pic\AuX$.
\item If \pols $f_1$, \ldots, $f_n$,  $g_1$, \ldots, $g_n$
in $\BuX$  satisfy the conditions~$(*)$ of Lemma~\ref{lemPicPic2},
then the $f_i$'s and the $g_i$'s are in $\AuX$.
\item If \pols $f_1$, \ldots, $f_n$,  $g_1$, \ldots, $g_n$
in $\BuX$ satisfy the conditions~$(*)$, then $f_1\in\AuX$.
\end{enumerate}

\end{corollary}
%
\subsec{Seminormal \risz}

An integral \ri $\gA$ is said to be \emph{seminormal} if, each time \hbox{that $b^2=c^3\neq 0$}, 
the \elt $a=b/c$ of $\Frac(\gA)$ is actually in $\gA$. 
In this case, $a^3=b$ \hbox{and $a^2=c$}.

\begin{definition}\label{defiseminormal}
An arbitrary \ri $\gA$ is said to be \ixc{seminormal}{\riz} if each time that $b^2=c^3$, there exists an $a\in\gA$ such that $a^3=b$ and $a^2=c$.%
\index{ring!seminormal ---}
\end{definition}

\begin{fact}
\label{factRed1}~
1. A seminormal \ri is reduced.\\
2. In a reduced \riz,  $x^2=y^2$ and  $x^3=y^3$ imply \hbox{$x=y$}.
\end{fact}
%
\begin{proof}
\emph{1.} If $b^2=0$,  then $b^2=0^3$,
hence $a\in\gA$ with $a^3=b$ and $a^2=0$, so $b=0$.\\
\emph{2.} In every \riz,  $(x-y)^3=4(x^3-y^3)+3(y^2-x^2)(x+y)$.
\end{proof}
Consequently the \elt $a$ in \Dfnz~\ref{defiseminormal}
is always unique. 
In addition, $\Ann(b)=\Ann(c)=\Ann(a)$.

\begin{fact}
\label{fact1seminormal}
Every normal \ri is seminormal.
\end{fact}
\begin{proof}
A \ri is normal when every \idp is \iclz.
Such a \ri is a \lsdz: if $uv=0$, there exists an $s$ such that $su=(1-s)v=0$
 (Lemma~\ref{lemiclplat}). Let $b$ and $c$ such that $b^3=c^2$, then~$c$ is integral over the \id $\gen{b}$, hence some $x$ such that $c=xb$,  
hence $b^3=c^2=x^2b^2$ \hbox{and $b^2(x^2-b)=0$}. Therefore there exists an $s$ such that $s(x^2-b)=0$ and $b^2(1-s)=0$. This gives  $b(1-s)=0$, then $(sx)^2=s^2b=sb=b$.
By letting $a=sx$, we get $a^2=b$, $a^3=bsx=bx=c$.
\end{proof}
\endgroup

\subsubsection*{The condition is necessary: Schanuel's example}

\begin{lemma}\label{lemSchaSeminor}
If $\gA$ is reduced and $\Pic\gA=\Pic\AX$, then $\gA$ is seminormal.
\end{lemma}

\begin{proof}
Let $b$, $c\in\gA$ with $b^2=c^3$. Let $\gB=\gA[a]=\gA+a\gA$ be a reduced \ri containing $\gA$, with $a^3=b, \,a^2=c$.
\\
Consider the \pols $f_i$ and $g_j$ ($i,j=1,2$) defined as follows

\snic{
f_1=1+aX,\; f_2=g_2=cX^2 \; \hbox{and} \;g_1=(1-aX)(1+cX^2).}

We have $f_1g_1+f_2g_2=1$, $f_1(0)=g_1(0)=1$, $f_2(0)=g_2(0)=0$, and each product $m_{ij}=f_ig_j$ is in $\AX$.
We apply Lemma~\ref{lemPicPic2}: the image of the matrix $(m_{ij})$ is free \ssi $f_1\in\AX$, \cad $a\in\gA$.
\end{proof}

\noi Note: For $\gB$ we can take 
$\left(\aqo{\gA[T]}{T^2-c, T^3-b} \right){}\!\red$. If a suitable \elt $a$ is already present in $\gA$, we obtain by uniquness $\gB=\gA$.

\subsec{The case of integral \risz}
We first treat the GCD-domains, then the normal \ris and finally the seminormal \risz.
\subsubsection*{The case of a GCD-domain}

Recall that a  GCD-domain is an integral \ri in which two arbitrary \elts admit a greatest common divisor, \cad an upper bound for the divisibility relation.
Also recall that if~$\gA$ is a GCD-domain, the same goes for the \pol \riz~$\AuX$.

\begin{lemma}
\label{lemPicGcd}
If $\gA$ is a GCD-domain, then $\Pic\gA=\so{1}$.
\end{lemma}

\begin{proof}
We use the \carn given in Lemma~\ref{lempropImProjLib}.\\
Let $P=(m_{ij})$ be an \idme matrix of rank $1$.
Since $\sum_i m_{ii}=1$, we can assume that $m_{1,1}$ is regular.
Let $f$ be the gcd of the \elts of the first row. We write $m_{1j}=fg_j$ with the gcd of the $g_j$'s equal to $1$.  The \egtz~\hbox{$m_{1,1}m_{ij}=m_{1j}m_{i1}$}  gives, by simplifying by $f$, $g_1m_{ij}=m_{i1}g_j$.
Thus,~$g_1$ divides all the $m_{i1}g_j$'s, and so also divides their gcd $m_{i1}$.
We write $m_{i1}=g_1f_i$. Since~\hbox{$g_1f_1=m_{1,1}=fg_1$}, this gives $f_1=f$. Finally,  $m_{1,1}m_{ij}=m_{1j}m_{i1}$ gives the \egtz~\hbox{$f_1g_1m_{ij}=f_1g_jg_1f_i$}, then $m_{ij}=f_ig_j.$
\end{proof}
 We then have the following corollary.

\pagebreak	

\begin{proposition}\label{thcorlemPicGcd}
If $\gA$ is a discrete field or a reduced \zed \riz, then $\Pic\gA=\Pic\AuX=\so{1}$.
\end{proposition}
\begin{proof}
Lemma~\ref{lemPicGcd} gives the result for the discrete fields.
It then suffices to apply the \elgbm \num2 (\paref{MethodeZedRed}).
\end{proof}

\subsubsection*{The case of a normal domain}

\begin{lemma}
\label{lemIntegclos}
If $\gA$ is a normal domain, then $\Pic\gA=\Pic\AuX$.
\end{lemma}
\begin{proof}
We use the \carn given in Corollary~\ref{corlemPicPic2}~\emph{3},
with here $\gA\subseteq\gK$, the quotient field of $\gA$. Let $f_i$ and $g_j$, $(i,j\in\lrbn)$ be the suitable \pols of $\KuX$. Then, since $f_1g_1=m_{1,1}\in\AuX$ and  $g_1(\uze)=1$, given \KROz's \tho \ref{thKro}, the \coes of $f_1$ are integral over the \ri generated by the \coes of~$m_{1,1}$. Thus $f_1 \in \gA[X]$.
\end{proof}

\rem As for Proposition~\ref{thcorlemPicGcd}, we can extend the result of Lemma~\ref{lemIntegclos} to the case of a reduced \ri $\gA$ integrally closed in a reduced \zed \ri $\gK\supseteq\gA$.
\eoe

\subsubsection*{The case of a seminormal integral \riz}

\begin{proposition}
\label{propIntSemin}
If $\gA$ is integral and seminormal, then $\Pic\gA=\Pic\AuX$.
\end{proposition}
\begin{Proof}{Start of the \demz. }
As in the \dem of Lemma~\ref{lemIntegclos},
we start with \pols $f_1(\uX)$, \ldots, $f_n(\uX)$, $g_1(\uX)$, \ldots, $g_n(\uX)$ in~$\KuX$ 
that satisfy the conditions $(*)$ of Lemma \ref{lemPicPic2}.
 We call $\gB$ the sub\ri of $\gK$ generated by $\gA$ and by the \coes of the $f_i$'s and the $g_j$'s, or,
what amounts to the same thing, generated by $\gA$ and the \coes of $f_1$. Then, given \KROz's \thoz, $\gB$ is a finite extension of $\gA$. Our goal is to show that $\gA=\gB$.
Let $\fa$ be the conductor of $\gB$ into $\gA$, 
\cad the set~\hbox{$\sotq{x\in\gB}{x\gB\subseteq\gA}$}. 
It is both an \id of $\gA$ and $\gB$. Our goal is now to show $\fa=\gen{1}$, \cad $\gC=\gA\sur{\fa}$ is trivial.
\end{Proof}

We start with two lemmas.

\begin{lemma}
\label{lemIntSemin1}
If $\gA\subseteq\gB$, $\gA$ is seminormal and $\gB$ is reduced, then the
conductor~$\fa$ of $\gB$ into $\gA$ is a radical \id of $\gB$.
\end{lemma}
\begin{proof}
We must show that if $u\in\gB$ and $u^2\in\fa$, then $u\in\fa$. So let $c\in\gB$. We must show that $uc\in\gA$. We know that $u^2c^2$ and $u^3c^3=u^2(uc^3)$ are in $\gA$ since $u^2\in\fa$.
Since $(u^3c^3)^2=(u^2c^2)^3$, we have some $a\in\gA$ such that~\hbox{$a^2=(uc)^2$} and $a^3=(uc)^3$. As $\gB$ is reduced, we obtain $a=uc$, and so $uc\in\gA$.
\end{proof}

\rem The \emph{seminormal closure} of a \ri $\gA$ in a reduced \riz~$\gB\supseteq \gA$ is obtained by starting from $\gA$ and adding the \elts $x$ of $\gB$  such that $x^2$ and $x^3$ are in the previously contructed \riz. Note that by Fact~\ref{factRed1}, $x$ is uniquely determined by the given~$x^2$ and~$x^3$. 
The \dem of the previous lemma can then be interpreted as a \dem of the following variant.%
\index{seminormal!closure in a reduced over\riz}
\eoe

\begin{lemma}
\label{lemIntSemin1bis}
Let $\gA\subseteq\gB$ be reduced, $\gA_{1}$ be the seminormal closure of  $\gA$ in~$\gB$, and $\fa$ be the conductor of $\gB$ into $\gA_{1}$.
Then, $\fa$ is a radical \id of~$\gB$.
\end{lemma}

\begin{lemma}
\label{lemIntSemin2}
Let $\gA\subseteq\gB$,  $\gB=\gA[\cq]$ be reduced and finite over $\gA$ and $\fa$ be the conductor of $\gB$ into $\gA$.
 Suppose that $\fa$ is a radical \idz, then it is equal to $\sotq{x\in\gA}{xc_1,\ldots ,xc_q\in\gA}$.
\end{lemma}
\begin{proof}
Indeed, if $xc_i\in\gA$, then  $x^\ell c_i^\ell\in\gA$ for all $\ell$, and so for some large enough~$N$,
$x^N y\in\gA$ for all $y\in\gB$, so $x$ is in the nilradical of $\fa$ (if $d$ is the upper bound of the degrees of the equations of integral dependence of the~$c_i$'s over~$\gA$, we can take $N=(d-1)q$).
\end{proof}

\begin{Proof}{End of the \dem of Proposition~\ref{propIntSemin}. }\\
We first give it in \clamaz. 
The natural classical reasoning would proceed by contradiction:  the \ri $\gC$ is trivial because otherwise, it would have a \idemi and the \lon at this \idemi would lead to a contradiction.\\
To avoid the non\cof \crc  of the argument by contradiction, we localize at a \fimaz, recalling our \dfn \gui{without negation} according to which a filter is maximal \ssi the localized \ri is a \zed \aloz.
In other words we tolerate for the \fimas of a \ri not only the complements of the \idemis but also the filter generated by $0$ which gives by \lon the trivial \riz.
In \clama a \ri  is then trivial \ssi its only \fima is the whole \ri (in other words, the filter generated by $0$).
\\ 
Let us insist on the fact that it is only in the previous affirmation that the \gui{classical} \crc of the argument is located.
Because the \dem of what follows is perfectly \covz: if $S$ is a \fima of $\gC$, then~\hbox{$0\in S$} (so $S=\gC$).
\\ 
Consider the inclusion ${{\gC = {\gA}\sur{\fa}} \subseteq {{\gB}\sur{\fa} = \gC'}}$.
Let $S$ be a \fima of~$\gC$, and $S_1$ be the corresponding \fima of $\gA$ (the inverse image of $S$  by the canonical \prnz).
Since $S$ is a \fimaz, and since~$\gC$ is reduced, $S^{-1}\gC=\gL$ is a reduced \zed \aloz, that is a discrete field, contained in the reduced \ri $S^{-1}\gC'=\gL'$.\\
If  $x$ is an object defined over $\gB$, let us denote by $\ov{x}$ what it becomes after the base change  $\gB\to\gL'$.
 Since $\gL$ is a \cdiz, $\gL[\uX]$ is a  GCD-domain, and the $\ov{f_i}$'s and $\ov{g_j}$'s are in $\gL[\uX]$.
 This means that there exists an $s\in S_1$  such that  $sf_1\in\AX$.
 By Lemma~\ref{lemIntSemin2}, this implies that $s\in\fa$. Thus $\ov s=0$  and $\ov s\in S$.
\end{Proof}

The \dem given above for Proposition~\ref{propIntSemin} is in fact quite simple.
It is however not entirely \cov and it seems to only treat the integral case.

\begin{Proof}{\Cov \dem of Proposition~\ref{propIntSemin}. }\\
We rewrite the \dem given in \clama by considering that the \fima $S$ of $\gC$ is a purely generic object which guides us in the \prcoz.
\\ 
Imagine that the \ri $\gC$ is a \cdiz, \cad that we have already done the \lon at a \fimaz. 
\\ 
Then, \pols ${F_i}$ and ${G_j}$ of $\gC[\uX]$ satisfying ${F_i}{G_j}=\ov{m_{ij}}$ and $F_1(\uze)=1$ are computed from the $\ov{m_{ij}}$'s according to an \algo that we deduce from the previously given \prcos for the case of  \cdis (Lemma~\ref{lemPicGcd}).
The uniqueness of the solution then forces the \egt $F_1=\ov{f_1}$, which shows that~$\ov{f_1}\in\gC[\uX]$, and therefore that $\gC$ is trivial.
\\ 
This \algo uses the disjunction \gui{$a$ is null or $a$ is \ivz,} for the \elts $a\in\gC$ which are produced by the \algo from the \coes of the \pols $\ov{m_{i,j}}$.  As $\gC$ is only a reduced \riz,
with neither a test for equality to 0 nor an invertibility test, the algorithm for discrete fields, if we execute it with $\gC$, must be replaced by a tree in which we open two branches each time a question 
 \gui{is $a$ null or invertible?} is asked by the \algoz.
\\ 
Here we are, facing a gigantic, but finite, tree.
Say that we have systematically placed the \gui{$a$ is \ivz}  branch on the left-hand side, and the \gui{$a=0$} branch on the right.
Let us look at what happens in the extreme left branch. 
\\ 
We have successively inverted $a_1$, \ldots , $a_p$ and we have obtained an $s$  that shows that the \riz~$\gC[1/(a_1\cdots a_p)]$ is trivial.
\\
 \emph{Conclusion: in the \ri $\gC,$ we have the \egt $a_1\cdots a_p=0$.}
\\
 Let us take one step back up the tree. 
\\ 
In the \ri $\gC[1/(a_1\cdots a_{p-1})]$, we know that $a_p=0$.
\\ 
The left branch should not have been opened.
Let us take a look at the computation in the branch $a_p=0$.
\\ 
Let us follow from here the extreme left branch.
\\ 
We have inverted  $a_1$, \ldots , $a_{p-1}$, then, say $b_1,\ldots
,b_k$ (eventually, $k=0$).
We obtain an $s$ that shows that the \riz~$\gC[1/(a_1\cdots a_{p-1} b_1\cdots  b_k)]$ is trivial.
\\
 \emph{Conclusion: in the \ri $\gC,$ we have the \egt $a_1\cdots a_{p-1} b_1\cdots  b_k=0$.}
\\ 
Let us take one step back up the tree. 
We know that $b_k=0$ (or, if $k=0$, $a_{p-1}=0$) in the \ri that was there just before the last branching; namely the \ri  $\gC[1/(a_1\cdots a_{p-1} b_1\cdots  b_{k-1})]$  (or, if $k=0$, $\gC[1/(a_1\cdots a_{p-2})]$).
  The left branch should not have been opened.
Let us look at the computation in the branch~\hbox{$b_k=0$} (or, if $k=0$, the branch $a_{p-1}=0$)\,\ldots
\\
 \emph{And so forth.} 
         When we follow the process all the way through, 
 we find ourselves at the root of the tree with the \riz~$\gC=\gC[1/1]$, which is trivial.
\end{Proof}

By using Lemma~\ref{lemIntSemin1bis} instead of Lemma~\ref{lemIntSemin1} we will obtain the following result, which is more precise than Proposition~\ref{propIntSemin}.

\pagebreak	

\begin{proposition}
\label{propIntSeminBis}
If $\gA$ is an integral \ri and $P$ is a \pro module of rank $1$ over $\AuX$ such that $P(\uze)$ is free,
there exist $c_{1}$, \ldots, $c_{m}$ in the quotient field of $\gA$ such that
\begin{enumerate}
\item $c_{i}^2$ and $c_{i}^3$ are in $\gA[(c_{j})_{j<i}]$ for $i=1,\ldots,m$,
\item $P$ is free over $\gA[(c_{j})_{j\leq m}][X]$.
\end{enumerate}
\end{proposition}
\rem Actually, only the quotient field of the sub\ri generated by the \coe present in a \prn matrix, whose image is \isoc to~$P$, intervenes.
\eoe

\subsec{General case}
\label{subsecTravSwanGeneral}

\begin{proposition}
\label{thseminormalCoq} $\!$\emph{(Coquand)}
Let $\gA\subseteq\gK$ with $\gK$ reduced.
\begin{enumerate}
\item
Given $f$ and $g\in\KuX^{n}$ that satisfy the conditions $(*)$ of Lemma~\ref{lemPicPic2},
we can construct $c_{1}$, \ldots, $c_{m}$ in  $\gK$ such that
\begin{enumerate}
\item [--]  $c_{i}^2$ and $c_{i}^3$ are in $\gA[(c_{j})_{j<i}]$ for $i\in\lrbm$,
\item [--]  $f$ and $g$ have their \coos in $\gA[(c_{k})_{k\in\lrbm}][\uX]$.
\end{enumerate}
\item
If   $\Pic\gK=\Pic\KuX$ and if $P$ is a \pro module of rank $1$ over $\AuX$, there exist $c_{1}$, \ldots, $c_{m}$ in  $\gK$ such that
\begin{enumerate}
\item [--]  $c_{i}^2$ and $c_{i}^3$ are in $\gA[(c_{j})_{j<i}]$ for $i\in\lrbm$,
\item [--]  $P\simeq P(\uze)$  over $\gA[(c_{k})_{k\in\lrbm}][\uX]$.
\end{enumerate}
\end{enumerate}
\end{proposition}
\begin{proof}
The \dem of Proposition~\ref{propIntSemin}, or of its more precise variant~\ref{propIntSeminBis},
is in fact a \dem of item \emph{1} above. Item \emph{2} is easily deduced.
\end{proof}

\begin{theorem}\label{thTSC} $\!$\emph{(Traverso-Swan-Coquand)}\\
If $\gA$ is a seminormal \riz, then $\Pic\gA=\Pic\AuX$.
\end{theorem}
\begin{proof}
We deduce it from the previous proposition by using the fact that there exists an over\ri $\gK$ of $\gA$ such that $\Pic\gK=\Pic\KuX$.
Indeed, every reduced \ri is contained in a reduced \zed \ri (\thref{thZedGen} or~\ref{thAmin})
$\gK$, which satisfies $\Pic\gK=\Pic\KuX=\so1$ (Proposition~\ref{thcorlemPicGcd}).
\end{proof}
%

\subsubsection*{A direct computation leading to the result}

As  is often the case when trying to implement on a machine a \cov \tho that has an elegant \demz, we are led to finding certain shortcuts in the computations that give a definitively simpler solution. But this solution partially hides at least the thought process that developed the \demz, if not the deep mechanism of the initial \demz.
See for example how Exercise~\ref{exo7.1} trivializes the \dem of the local structure \tho for \mptfsz.

This is what happened with Proposition~\ref{thseminormalCoq} which was finally realized by a quite \elr \algo in \cite[Barhoumi\&Lombardi]{BL07}, based on the theory of the resultant \id (see Section~\ref{subsecIdealResultant}) and of the subresultant modules.


\section{Patching \`a la Quillen-Vaserstein}
\label{subsecQPatch}

In this section we present the so called Quillen patching.
It is a deep result that could a priori seem a little too abstract (abusive usage of \idemasz) 
but which happens to make a lot of \cov sense.

The \dems that we give are (for the most part) copied from~\cite{Kun}.
We have replaced the \lon at any maximal \id with the \lon at \mocoz.

\begin{lemma}
\label{lem0PrepVaser}
Let $S$ be a \mo of the \ri $\gA$ and $P\in\AX$ be a \pol  such that $P=_{\gA_S[X]}0$ and $P(0)=0$. Then, there exists an $s\in S$  such that $P(sX)=0$.
\end{lemma}
\facile

Here is a slight variant.

\begin{fact}
\label{lem01PrepVaser}
Let $S$ be a \mo of the \ri $\gA$ and $P\in\gA_S[X]$ be a \pol  such that  $P(0)=0$. Then, there exist $s\in S$  and $Q\in\AX$ such \hbox{that $P(sX)=_{\gA_S[X]}Q$}.
\end{fact}
%
%

\begin{lemma}
\label{lem1PrepVaser}
Let $S$ be a \mo of the \ri $\gA$. Consider three matrices with \coes in $\AX$, $A_1,\,A_2,\,A_3$ such that the  product $A_1\,A_2$ has the same format as $A_3$. If $A_1\,A_2=_{\gA_S[X]}A_3$ and $A_1(0)\,A_2(0)=A_3(0)$, there exists an $s\in S$  such that $A_1(sX)\,A_2(sX)=A_3(sX)$.
\end{lemma}
\begin{proof}
Apply Lemma~\ref{lem0PrepVaser} to the coefficients of the matrix  $A_1\,A_2-A_3$.
\end{proof}

\begin{lemma}
\label{lem2PrepVaser}
Let $S$ be a \mo of the \ri $\gA$ and $C(X)\in\GL_p(\gA_S[X])$.  There exist $s\in S$  and $U(X,Y) \in
\GL_p(\gA[X,Y])$ such that $U(X,0)=\I_p$, and, over the \ri $\gA_S[X,Y]$,
$U(X,Y)=C(X+sY)C(X)^{-1}$.
\end{lemma}
\begin{proof}
Let $E(X,Y)=C(X+Y)C(X)^{-1}$. Let $F(X,Y)=E(X,Y)^{-1}$. We have $E(X,0)=\I_p$, so $E(X,Y)=\I_p+E_1(X)Y+\cdots+E_k(X)Y^k$. For some $s_1\in S$, the ${s_1}^jE_j$'s can be rewritten \gui{without denominator.}
We thus obtain a matrix $E'(X,Y)\in \MM_p(\gA[X,Y])$ such that $E'(X,0)=\I_p$ and, over $\gA_S[X,Y]$,
$E'(X,Y)=E(X,s_1Y)$. We proceed similarly with $F$ (and we can choose some common $s_1$).
We then have $E'(X,Y)F'(X,Y)=\I_p$ %
in $\MM_p(\gA_S[X,Y])$ and $E'(X,0)F'(X,0)=\I_p$. \\
By applying Lemma~\ref{lem1PrepVaser} in which we replace $X$ with $Y$ and $\gA$ with~$\AX$, we obtain $s_2\in S$ such that $E'(X,s_2Y)F'(X,s_2Y)=\I_p$. \\
Hence the desired result with $U=E'(X,s_2Y)$ and $s=s_1s_2$.
\end{proof}

\vspace{-.7em}
\pagebreak	

\begin{lemma}
\label{lem3PrepVaser}
Let $S$ be a \mo of  $\gA$ and $G\in \AX^{q\times m}$. If $G(X)$ and $G(0)$ are \eqves over $\gA_S[X]$, there exists an $s\in S$ such that $G(X+sY)$ and $G(X)$ are \eqves over $\gA[X,Y]$.
\end{lemma}
\begin{proof}
Let $G=C\,G(0)\,D$ with $C\in\GL_q(\gA_S[X])$ and $D\in\GL_m(\gA_S[X])$. We therefore have 
\[\preskip.6em \postskip.2em \arraycolsep2pt
\begin{array}{rcl} 
 G(X+Y) &  = &  C(X+Y)G(0)D(X+Y) \\[1mm] 
  &  = &   C(X+Y)C(X)^{-1}G(X)D(X)^{-1}D(X+Y).
 \end{array}
\]
By applying Lemma~\ref{lem2PrepVaser}, we obtain $s_1\in S$, $U(X,Y)\in\GL_q(\gA[X,Y])$ and $V(X,Y)\in\GL_m(\gA[X,Y])$,
such that 
$$\preskip.2em \postskip.2em 
U(X,0)=\I_q\,,\;\;   V(X,0)=\I_m \,, 
$$
and, over the \ri   $\gA_S[X,Y]$, 
$$\preskip.2em \postskip.2em 
 U(X,Y)=C(X+s_1Y)C(X)^{-1}   \hbox{ and }     V(X,Y)=D(X)^{-1}D(X+s_1Y).
$$
Therefore
$$\preskip-.4em \postskip.2em 
  G(X)    =    U(X,0)G(X)V(X,0),   
$$
and over the \ri  $\gA_S[X,Y]$  :  
$$\preskip.2em \postskip.2em 
  G(X+s_1Y)   =   U(X,Y)G(X)V(X,Y).
$$  
By applying Lemma~\ref{lem1PrepVaser}
(as in Lemma~\ref{lem2PrepVaser}), we obtain $s_2\in S$ such that $G(X+s_1s_2Y)=U(X,s_2Y)G(X)V(X,s_2Y)$.\\
Hence the result with $s=s_1s_2$.
\end{proof}

\begin{plcc}
\label{thPatchV}  $\!$\emph{(Vaserstein patching)}\\
Let $G$ be a matrix over $\AX$ and $S_1$, $\ldots$, $S_n$ be \moco of~$\gA$.
\begin{enumerate}
\item The matrices $G(X)$ and $G(0)$ are \eqves over $\AX$ \ssi  they are \eqves over $\gA_{S_i}[X]$ for each $i$.
\item Same result for \gui{the left-\eqvcz}: two matrices $M$ and $N$ with the same format over a commutative \ri are said to be \emph{left-\eqvesz} if there exists an \iv square matrix $H$ such that $H\,M=N$.
\end{enumerate}
\index{matrices!left-equivalent ---}
\index{left-equivalent!matrices}
\end{plcc}
\begin{proof}
\emph{1.} One sees easily that the set of $s\in\gA$ such that the matrix $G(X+sY)$ is \eqve to $G(X)$ over $\gA[X,Y]$ forms an \id of $\gA$. By applying Lemma \ref{lem3PrepVaser}, this \id contains an \elt $s_i$ in $S_i$ for each $i$, so it contains $1$, and $G(X+Y)$ is \eqve to $G(X)$. It remains to 
make~\hbox{$X=0$}.\\
\emph{2.} In all the previous \demsz, we can replace  \eqvc with  left-\eqvcz.
\end{proof}

\begin{plcc}
\label{thPatchQ}  $\!$\emph{(Quillen patching)}\index{Quillen}\index{Quillen patching}\\
Let $M$ be a \mpf over $\AX$ and $S_1$, $\ldots$, $S_n$ be \moco of $\gA$.
Then, $M$ is a module extended from $\gA$ \ssi each  $M_{S_i}$ is extended from~$\gA_{S_i}$.
\end{plcc}
\begin{proof} This is a corollary of the previous \tho because the \iso between the modules $M(X)$ and $M(0)$ can be expressed by the \eqvc of two matrices $H(X)$ and  $H(0)$ constructed from a \mpn $G$ of $M$ (see Fact~\ref{factEtPol}).
\end{proof}

\comm
The original formulation by Quillen\index{Quillen patching}, \eqve to the \plgrf{thPatchQ} in \clamaz, is the following: \emph{if $M_\fm$ is extended from $\gA_\fm$ after \lon at every \idema $\fm$, 
then $M$ is extended from $\gA$.}
To \cot  decipher 
a classical \dem based on the Quillen \rcm in the original formulation, we will have to call upon the basic \lgbe machinery explained in Section~\ref{secMachLoGlo}.\imlbz
\eoe


\subsec{A Roitman \thoz}
\label{secRoitman}
This subsection is devoted to the \dem of the following \thoz, which consists in a kind of converse of the Quillen \rcm \thoz.

\begin{theorem}\label{thRoitman}  $\!$\emph{(Roitman's \thoz)}\\
Let $r$ be an integer $\geq1$ and $\AuX=\AXr$. If every \ptf $\AuX$-module is extended from $\gA$, then every \lon $\gA_S$ of $\gA$ satisfies the same \prtz. 
\end{theorem}

\subsubsection*{The univariate case}
\begin{lemma}\label{lemRoitman}
If every \ptf $\AX$-module is extended from~$\gA$, then every \lon $\gA_S$ of $\gA$ satisfies the same \prtz.
\end{lemma}
\begin{proof}
\emph{Special case: 
$\gA_S$ is a \dcd \aloz.} 
\\
Let $\rho:\AX\to\gA_S[X]$ be the natural morphism.
Let $M \in\GAn(\gA_S[X])$.
Since $\gA_S$ is local, $M(0)$ is similar to a standard \prr $\I_{k,n}$.
We can therefore suppose \spdg that $M(0)=\I_{k,n}$, \cad $M(X)=\I_{k,n}+M'(X)$ with
$M'(X) \in\Mn(\gA_S[X])$ and $M'(0) = 0$. 
\\
Let~$v$ be the \gui{product of the \denosz} in the \coes of the entries of~$M'(X)$.
Since $M'(0) = 0$, we have a matrix $N' = N'(X)\in\Mn(\AX)$ such that $M'(vX)=_{\gA_S[X]}N'(X)^{\rho}$ and $N'(0)=0$.
\\
With $N(X)=\I_{k,n}+N'(X)$ we obtain~\hbox{$N(0)=\I_{k,n}$} and $M(vX)=N(X)^{\rho}$.
Since $M^2=_{\gA_S[X]}M$, we have \hbox{some $s\in S$} such that $s(N^2-N)=0$.
\\
As $(N^2-N)(0)=0$, we write $N^2-N=XQ(X)$.\\
Now $sXQ(X)=0$ implies $sQ(X)=0$.
A fortiori $sQ(sX)=0$,  
so $N(sX)^2=N(sX).$
However, the \mptfs over  $\AX$ are extended from $\gA$, therefore the \mprn $N(sX)$ has a kernel and an image \isoc to the kernel and to the image of $N(0)=\I_{k,n}$. Therefore $N(sX)$ is similar to $\I_{k,n}$: there exists a $G=G(X)\in\GL_n(\AX)$ such that
$$
G^{-1}(X)N(sX)G(X)=\I_{k,n}.
$$
By letting $H(X)=G(X)^{\rho}\in\GL_n(\gA_S[X])$, we obtain over $\gA_S[X]$ the \egt  
\[\preskip-.2em \postskip.2em
H^{-1}(X)M (svX)H(X)=\I_{k,n}
\] 
and therefore
$$\preskip.0em \postskip.4em
H^{-1}(X/sv)M(X)H(X/sv)=_{\gA_S[X]}\I_{k,n}
$$
with $H(X/sv)\in\GL_n(\gA_S[X])$.

 \emph{\Gnl case.} Let $P$ be an arbitrary \ptf $\gA_S[X]$-module. 
Let $\gB=\gA_S$.
As usual, $P(0)$ denotes the \Bmo obtained by \eds via the morphism $\ev_0:\BX\to\gB$.
We apply the basic \lgbe machinery (\paref{MethodeIdeps}) to the \cov \dem that we have just given in the special case.
We obtain \moco $V_1$, \ldots, $V_m$ of $\gB$ with $P\simeq_{\gB_{V_i}}P(0)$. We conclude with the Quillen \rcmz: $P\simeq_{\gB}P(0)$.\imlbz
\end{proof}

\medskip 
\rem \emph{To implement the \algo corresponding to this \demz.} 
Actually the only particular \prt that we have used in the \dem of the special case, it is that the \ptf $\gA_S$-modules are free. Therefore the implementation of the basic \lgbe machinery here is very \elrz. It consists in constructing \come \lons for which  
the matrix $M(0)$ becomes similar to a standard \mprnz, and to get the \algo given by the \dem of the special case running in each of these \lonsz.
Naturally, we end with the \algo corresponding to the \cov \dem of the Quillen \rcmz.\imlbz
\eoe
\subsubsection*{The multivariate case}
\begin{Proof}{\Demo of Roitman's \thoz~\ref{thRoitman}. }
We reason by \recu on $r$. The case $r=1$ has already been treated. Let us pass from $r\geq1$ to $r+1$. 
Consider a \moz~$S$ of a \ri $\gA$. Let $(\Xr)=(\uX)$.\\  
We have $\gA_S[\uX,Y]=\gA[\uX,Y]_S=(\gA[Y]_S)[\uX]$.
Let $P$ be a \ptf $\gA_S[\uX,Y]$-module. 
By the \hdr applied with the \ri $\gA[Y]$, $P$ is extended from $\gA[Y]_S=\gA_S[Y]$, that is, $P(\uX,Y)$ is \isoc to $P(\uze,Y)$ as an $\gA_S[\uX,Y]$-module, 
and by the case $r=1$ applied with the \ri $\gA$, $P(\uze,Y)$ is extended from~$\gA_S$.
\end{Proof}
%
\subsubsection*{A long-open question solved negatively}

That question is the following.

\emph{If every \ptf $\AX$-module is extended from $\gA$, is it always true that for any $r$ every \ptf $\AXr$-module is extended from $\gA$?}

A negative response is given in \cite[Corti\~nas \& al., (2011)]{CHWW}. 

\subsubsection*{A \lgb principle \`a la Roitman}

\begin{plcc} \label{plgcetendus}
Let $n$ and $r>0$. Consider the following \prt for a \ri $\gA$. 
{\sf P}$_{n,r}(\gA)\,:$ every \mrc $r$ over $\AXn$ is extended from $\gA$.\\
Let $S_1$, \ldots, $S_k$ be \moco of a \ri  $\gA$.
Then $\gA$ satisfies the \prt {\sf P}$_{n,r}$ \ssi each of the $\gA_{S_i}$'s satisfies it.\\
In particular $\gA$ is seminormal \ssi each of the $\gA_{S_i}$'s is seminormal.
\end{plcc}
%
\begin{proof}
The condition is \ncr by Roitman's \tho \ref{thRoitman},  
whose \dem remains valid if we limit ourselves to the \mrcs  $r$.
\\
The condition is sufficient by the Quillen \rcmz~\ref{thPatchQ}. 
\end{proof}


\section{Horrocks' \thoz}
\label{sec.Horrocks} 

The following lemma is a special case of Proposition~\ref{propFiniBon1}~\emph{4}.
\begin{lemma}\label{lemLocptfLib}
Let $S$ be a \mo of $\gA$ and $P$, $Q$ be \ptfs \Amos  such that
$P_{S}\simeq Q_S$. Then, there exists an $s\in S$ such that
$P_{s}\simeq Q_s$.
\end{lemma}

\begin{notation}\label{notaA<X>} 
Let $\ArX$ be the \ri $S^{-1}\AX$, where $S$ is the \mo of the \polus of $\AX$.\label{NOTAAlraX}
\end{notation}

\begin{theorem}\label{thHor0} $\!$\emph{(Local Horrocks' \thoz)}\\
Let $\gA$ be a \dcd \alo and $P$ be a \mptf over $\AX$.
If $P_S$ is free over $\ArX$, then $P$ is free over $\AX$ (so extended from~$\gA$).
\end{theorem}

We use the \dem by \cite[Nashier \& Nichols]{NaNi}
 which is almost \covz, as presented in \cite{Lam06} or \cite{IRa}.

We need a few preliminary results.
\perso{introduire qq part la terminologie \gui{sp\'ecialisation}}

\begin{lemma}\label{lemLocLocCom}
Let $\gA$ be a \riz, $\fm=\Rad\gA$  and $S\subseteq\AX$ be the \mo of the \polusz. The \mos $S$ and $1+\fm[X]$ are \comz.
\end{lemma}
\begin{proof}
Let $f(X)\in S$ and $g(X)\in 1+\fm[X]$.
The resultant $\Res_X(f,g)$ belongs to the \id $\gen{f,g} $ of $\AX$.
Since $f$ is \monz, the resultant is successfully subjected to the specialization $\gA\to\gA\sur\fm$. Therefore $\Res_X(f,g)\equiv\Res_X(f,1)=1  \mod \fm$.
\end{proof}
%
\begin{lemma}\label{lem0thHor0}
Let $\gA\subseteq\gB$,  $s\in\Reg(\gB)$, and $P$, $Q$ be two \ptf $\gB$-modules with $sQ\subseteq P\subseteq Q$. If $\gB$ and $\aqo{\gB}{s}$ are (not \ncrt \tfz) \pro $\gA$-modules, then the same goes for \hbox{the \Amoz~$Q/P$}.
\end{lemma}
\begin{proof}
Since $s$ is \ndz and since $Q$ and $P$ are submodules of a free module, the multiplication by $s$ (denoted $\mu_s$) is injective in $P$ and in $Q$.
We have the following exact sequences of \Amosz.
\[\arraycolsep2pt
\begin{array}{ccccccccccccccccc}
0 &\rightarrow& Q& \vvers{\mu_s} &P & \lora & P/sQ& \rightarrow& 0   \\[1mm]
0 &\rightarrow& sQ/sP& \llra &  P/sP& \lora & P/sQ& \rightarrow& 0   
\end{array}
\]
The \Amo $P$ is \pro by transitivity, the $\gB/s\gB$-module $P/sP$ is \proz, therefore by transitivity $P/sP$ is a \pro \Amoz.
We can then apply Schanuel's lemma (Lemma~\ref{lemScha}): 
$(P/sP)\oplus Q\simeq (sQ/sP) \oplus P$ as  \Amosz. Since $Q$ is a \pro \Amoz, the same goes for $sQ/sP$. But since $\mu_s$ is injective, $sQ/sP$ is \isoc to $P/Q$.
\end{proof}
%

\begin{lemma}\label{lem1thHor0} $\!$\emph{(Murthy \& Pedrini, \cite{MuPe})}
\\ 
Let $\gA$ be a \riz, $\gB=\AX$, $S$ be the \mo of the \polus of~$\AX$, $P$, $Q$ be two \mptfs over $\gB$, and $f\in S$.
\begin{enumerate}
\item If $fQ \subseteq P\subseteq Q$, then $P$ and $Q$ are stably \isocz.
\item If $P_S\simeq Q_S$,  then $P$ and $Q$ are stably \isocz.
\item If in addition $P$ and $Q$ are of rank $1$, then $P\simeq Q$.
\end{enumerate}
\end{lemma}
\begin{proof}
\emph{1.} Since $f$ is \monz, the \Alg $\aqo{\gB}{f}$ is a free \Amo of rank $\deg f$. The $\aqo{\gB}{f}$-module $Q/fQ$ is \ptf over $\gA$.  By the previous lemma, the \Amo $M=Q/P$ is \proz, and it is \tf over $\aqo{\gB}{f}$, so over $\gA$. Therefore $M[X]$ is a \ptf \Bmoz.
We have two exact sequences ($\mu_X$ is  multiplication by $X$)
\[\preskip.2em \postskip.4em
\arraycolsep2pt
\begin{array}{ccccccccccccccccc}
0 &\rightarrow& P& \llra &Q& \lora& M& \rightarrow& 0 ,  \\[1mm]
0 &\rightarrow& M[X]& \vvers{\mu_X} &M[X]& \lora & M& \rightarrow& 0.   
\end{array}
\]
By Schanuel's lemma (Lemma~\ref{lemScha}) we have $P\oplus M[X]\simeq Q\oplus M[X]$ as  \Bmosz. Since $M[X]$ is \ptf over $\gB$, $P$ and $Q$ are stably \isocz.

\emph{2.} We know that $P_f\simeq Q_f$ for some $f\in S$.\\
By hypothesis,  we have $F\in\GAn(\gB)$  and $G\in\GAn(\gB)$ with $P\simeq \Im F$, and~\hbox{$Q\simeq \Im G$}. 
We know that $F'=\Diag(F,0_m)$ and $G'=\Diag(G,0_n)$ are conjugated over $\gB_f$ (\dlg lemma~\ref{propIsoIm}). This means that there exists a matrix $H\in\MM_{m+n}(\gB)$ such that $HF'=G'H$ and $\det(H)=\delta$ divides a power of $f$. We then have $P_1=\Im(HF')\subseteq\Im G'$.
Then, by postmultiplying by $\wi H$, $(HF')\wi H=\delta G'$, which implies $\delta \Im G'\subseteq \Im(HF')$.
Since $H$ is injective, we have $P_1\simeq P$, and moreover $\Im G'=Q_1\simeq Q$. We can conclude with item \emph{1} since $\delta Q_1\subseteq P_1\subseteq Q_1$.

\emph{3.} The modules $P$ and $Q$ are of rank $1$ and stably \isocz, therefore \isoc (Fact~\ref{factPicStab}).
\end{proof}

\begin{Proof}{\Demo of \thref{thHor0}. } ~\\
Notations:  $\fm=\Rad\gA$, $\gk=\gA/\fm$ (\cdiz)
$\gB=\AX $, $n=\rg(P)$  ($n\in\NN$ since $\gB$ is connected), $U=1+\fm [X]$, and $\ov{E}$ be the object $E$ reduced modulo $\fm$.

\emph{1.} We show by \recu on $n$ that we have an \iso $P\simeq P_1\oplus\gB^{n-1}$. For $n=1$ it is trivial.

\textbf{Little lemma} {(see the \dem below)}\\
\emph{There exist $z$, $y_2$, \ldots, $y_n$,  $z_2$, \ldots, $z_n$ in $P$ such that $(z,y_2, \ldots,y_n)$ is a basis of $P_S$ over $\gB_S$ and $(\ov{z}, \ov{z_2}, \ldots, \ov{z_n})$ is a basis of $\ov{P}$ over $\ov{\gB}=\gk[X]$.}

The $\gB_U$-module $P_U$ is free with basis $(z,  z_2, \ldots, z_n)$: indeed, $\fm\subseteq\Rad\gB_U$, and modulo $\fm$, $(\ov{z}, \ov{z_2}, \ldots, \ov{z_n})$ generates  $\ov{P}=\ov{P_U}$, so $(z,  z_2, \ldots, z_n)$ generates~$P_U$ by Nakayama's lemma. Finally, a \mptf of rank~$n$ generated by $n$ \elts is free.\\
Let $P'=P/\gB z$. The two modules $P'_U$ and $P'_S$ are free. The \mosz~$U$ and $S$ are \com (Lemma~\ref{lemLocLocCom}), so $P'$ is \ptf over~$\gB$, hence $P\simeq P'\oplus \gB z$.
By \hdrz, $P'\simeq P_1\oplus \gB^{n-2}$, which gives $P\simeq P_1\oplus \gB^{n-1}$

\emph{2.} The \iso $P\simeq P_1\oplus \gB^{n-1}$ with $P_1$ of rank $1$ gives by \lon that $(P_1)_S$ is stably free.
We apply item \emph{3} of Lemma~\ref{lem1thHor0}: we obtain that $P_1$ is free.
\end{Proof}

\begin{Proof}{\Demo of the little lemma. }
Let $(\yn)$ in $P$ which is a $\gB_S$-basis of~$P_S$.
There exists a basis  $(\ov{z_1}, \ov{z_2}, \ldots, \ov{z_n})$ of $\ov P$, with the $z_i$'s in~$P$ such that $\ov{y_1}\in\kX\,\ov{z_2}$ 
(by dividing $\ov{y_1}$ by the gcd of its \coesz, we obtain a \vmdz,
and over a Bézout \riz, every \vmd is completable).
We look for $z$  in the form~\hbox{$z_1+X^ry_1$}.
It is clear that, for any~$r$,   $(\ov{z}, \ov{z_2}, \ldots, \ov{z_n})$ is a basis of $\ov P$.
Since $(\yn)$  is a $\gB_S$-basis of~$P_S$, there exists an $s\in S$  such that $sz_1=\sum_{i=1}^nb_iy_i$, with the $b_i$'s in~$\gB$. \\
Then, $sz=(b_1+sX^r)y_1+\sum_{i=2}^nb_iy_i$, and for large enough $r$, $b_1+sX^r$ is a \poluz: $(z, y_2, \ldots, y_n)$ is a $\gB_S$-basis of $P_S$.
\end{Proof}
We now give the global version.

\begin{theorem}\label{thHor}  $\!$\emph{(Affine Horrocks'  \thoz)}\\ 
Let $S$ be the \mo of the \polus of $\AX$ and $P$ be a \mptf over $\AX$. If $P_S$ is extended from $\gA$, then $P$  is extended from~$\gA$.
\end{theorem}
\begin{proof}
We apply the basic \lgbe machinery (\paref{MethodeIdeps}) with the \cov \dem of \thref{thHor0}.\imlb We obtain a finite family of \moco of $\gA$, $(U_i)_{i\in J}$,   with each localized module $P_{U_i}$ extended from $\gA_{U_i}$.
We conclude with the Quillen \rcmz\index{Quillen patching} (\plgc \ref{thPatchQ}).
\end{proof}

This important \tho can be completed by the following subtle result, which does not seem possible to extend 
to the \mpfsz.

\pagebreak

\begin{theorem}\label{th2Bass}  $\!$\emph{(Bass)}\\
Let $P$ and $Q$ be two \ptf \Amosz. If they are \isoc after \eds to $\ArX$, they are \isocz.
\end{theorem}
\begin{proof}
We reason with \mprns and similarities 
between these matrices that correspond to \isos between the image modules. Therefore implicitly, we systematically use the \dlg lemma \ref{propIsoIm}, without mentioning it.
\\ 
We start with $F$ and $G$ in $\GAn(\gA)$, conjugated over the \ri $\ArX$. The \mptfs are $P \simeq \Im F$ and $Q \simeq \Im G$. We therefore have a matrix~\hbox{$H\in\Mn(\AX)$}, with $\det(H)\in S$ (\mo of the \polusz), and $HF=GH$.\\
By letting $Y=1/X$, for large enough $N$, the matrix $Y^N H=H'$ is in~\hbox{$\Mn(\AY)$},
with $\det(H')=Y^r\big(1+Yg(Y)\big)=Y^rh(Y)$ where $h(0)=1$, and obviously $H'F=GH'$. In other words,~$F\sim G$~over the \ri $\gA[Y]_{Yh}$. \\
The \elts $Y$ and $h$ are \comz, so, by applying the \rcm \tho of the modules (\plgc \ref{plcc.modules 2}), there exists some $\AY$-module $M$ such that $M_Y$ is \isoc 
to~\gui{$P$ extended to $\AY_Y$,} and~$M_h$ is \isoc to \gui{$Q$ extended to $\AY_h$.} And~$M$ is  \ptf
since there are two \come \lons which are \mptfsz.
This provides a \mprn $E$ with \coes in $\AY$ such that $E\sim F$ over $\AY_Y$ and $E\sim G$ over $\AY_h$. Since $Y$ is a \poluz, Horrocks' \tho tells us that $\Im E$ comes by \eds from a \ptf \Amo $M'$.
Consequently, for all $a$, $b\in \gA$ the \gui{evaluated} matrices $E(a)$ and $E(b)$ are conjugated over $\gA$ (their images are both \isoc to $M'$).\\
Finally, $F\sim E(1)$ and $G\sim E(0)$ over $\gA$, therefore $F \sim G$ over $\gA$.
\end{proof}

\rem For the mathematician who 
wishes to implement the algorithm subjacent to the previous proof, we suggest
 suggest using \mpns (whose cokernels are the modules) rather than \mprns (whose images are the modules). This in particular avoids having to repetitively use an implementation of the \dlg lemma. 
\eoe

\medskip
We finish this section with a corollary of Lemma \ref{lem1thHor0}. This \tho is to be compared with \thref{th2QUILIND}.
\perso{Je n'ai pas of r\'ef\'erence pour ce \thoz}

\begin{theorem}\label{thStabLibPol}  $\!$\emph{(Concrete Quillen induction\index{Quillen induction!concrete ---, stably free case}, stably free case)}\\
Let $\cF$ be a class of \ris that satisfy the following \prtsz. 
\perso{donner un nom and rajouter une r\'{e}f\'{e}rence in la biblio}
\hsu\emph{1.} If $\gA\in\cF$, then $\gA\lra{X} \in\cF$.
\hsu
\emph{2.} If $\gA\in\cF$, every \pro $\gA$-module of constant rank is stably free.\\
Then, for $\gA\in\cF$ and  $r\in\NN$, every \pro $\AXr$-module of constant rank is stably free.
\end{theorem}
\begin{proof}
We proceed by \recu on $r$, the case $r=0$ being clear.\\
We pass from $r-1$ to $r$ ($r\geq 1$). Let $\gA$ be a \ri in the class $\cF$, and $P$ be a \mrc over $\AXr$.\\
Let $\gB=\gA[(X_i)_{i<r}]$, $\gC=\gA[X_r]$, and $V$ be the \mo of the \polus of $\gA[X_r]$.
Thus $\AXr\simeq \gB[X_r]\simeq \gC[(X_i)_{i<r}].$\\
The \ri $\gA\lra{X_r}=V^{-1}\gC$  is in the class $\cF$.\\
The $\gA\lra{X_r}[(X_i)_{i<r}]$-module $P_V$, which is \prcz, is \stl by \hdrz. \\
If $S$ is the \mo of the \polus of $\gB[X_r]$, we have $V\subseteq S$, and so~$P_S$ is \stl over the \riz~\hbox{$S^{-1}\gB[X_r]$}. By item \emph{2} of Lemma \ref{lem1thHor0}, $P$ is \stlz.
\end{proof}

\begin{corollary}\label{corthStabLibPol}
If $\gK$ is a \cdiz, every \mptf over $\KXr$ is stably free.
\end{corollary}
\begin{proof}
We apply the previous result with the class $\cF$ of the \cdisz:
 if~$\gK$ is a \cdiz, then $\gK\lra{X}=\gK(X)$ is also a \cdiz.
\end{proof}
%


\section{Solution to Serre's problem}
\label{sec.QS} 

In this section we present several \cov solutions to Serre's \pbz, in which $\gK$ is a \cdiz.

\Grandcadre{The \mptfs over $\KXr$ are free}


\subsec{\`A la Quillen\index{Quillen}}
\label{subsecQuillen}

The solution by Quillen of Serre's \pb is based on the Local Horrocks'  \tho and on the following \ix{Quillen induction} (see~\cite{Lam06}). 

\CMnewtheorem{IQa}{Abstract Quillen induction\index{Quillen induction!abstract ---}}{\itshape}
\begin{IQa}\label{InductionQabs}~\\
Let $\cF$ be a class of \ris that satisfies the following \prtsz.
\begin{description}
\item [{\rm (Q1)}]
If $\gA\in\cF$, then $\ArX \in\cF$.
\item [{\rm (Q2)}]
If $\gA\in\cF$, then $\gA_{{\fm}}\in\cF$ for every \idema ${\fm}$ of $\gA$.
\item [{\rm (Q3)}]
If $\gA\in\cF$ is local, every \ptf $\AX$-module is extended from $\gA$ (\cad free).
\end{description}
Then, for all $\gA\in\cF$ and all $r\geq1$, every \mptf over~$\AXr$ is extended from $\gA$. 
\end{IQa}

Actually, the \prts (Q1), (Q2) and (Q3) are first used by Quillen to obtain the case $r=1$, by using the Local Horrocks' \tho and the Quillen patching.
The \gui{\dem by \recuz} part is based over the case $r=1$, over (Q1) and over Horrocks' \tho (local or affine).

In what follows we isolate this \dem by \recuz, which we qualify as a \gui{concrete} Quillen induction.
We replace (Q3) with a stronger version (q3) which is case $r=1$.

In a posterior comment, we explain how we can actually somehow replace (q3) with (Q3) without losing the \cov \crc of the \demz.  

\subsubsection*{The \dem by \recu itself}

\begin{theorem} $\!$\emph{(Concrete Quillen induction\index{Quillen induction!concrete ---})} \label{thQUILIND}~\\
Let $\cF$ be a class of \ris that satisfy the following \prtsz.
\begin{description}
  \item[{\rm (q1)}] If $\gA\in\cF$, then $\ArX \in\cF$.
  \item[{\rm (q3)}] If $\gA\in\cF$, every \ptf $\AX$-module  is extended from~$\gA$.
\end{description}
Then, for all $\gA\in\cF$ and all $r\geq1$, every \mptf over~$\AXr$ is extended from $\gA$.
\end{theorem}
\begin{proof}
Let us pass from $r\geq1$ to $r+1$. Consider a \ptf $\gA[\Xr,Y]$-module $P=P(\Xr,Y)=P(\uX,Y)$. 
 Let
\begin{enumerate}
  \item [--] $P(\uX,0)$ be the $\AuX$-module obtained by the \homo $Y\mapsto 0$,
  \item [--] $P(\uze,Y)$ be the $\gA[Y]$-module obtained by the \homo $\uX\mapsto \uze$,
  \item [--] $P(\uze,0)$ be the \Amo obtained by the \homo $\uX,Y\mapsto \uze,0$.
\end{enumerate}
We must show that $P(\uX,Y)\simeq P(\uze,0)$ over $\gA[\uX,Y]$.
\\ 
We call $S$ the \mo of the \polus of $\gA[Y]$, 
that is contained in the \mo $S'$ of the \polus of $\AuX[Y]$. 
We then have
\begin{enumerate}
  \item  $P(\uX,Y)\simeq P(\uze,Y)$ over $\gA\lra{Y}[\uX]=\gA[\uX,Y]_{S}$ by \hdr since  $\gA\lra{Y}\in \cF$,
  \item  a fortiori $P(\uX,Y)\simeq P(\uze,Y)$ over $\gA[\uX]\lra{Y}=\gA[\uX,Y]_{S'}$,
  \item  $P(\uze,Y)\simeq P(\uze,0)$ over $\gA[Y]$ by the case $r=1$,
  \item  $P(\uze,0)\simeq P(\uX,0)$ over $\AuX$ by \hdrz,
  \item  by combining 2, 3 and 4, we have $P(\uX,Y)\simeq P(\uX,0)$ over $\gA[\uX]\lra{Y}$,
  \item  so, by the Affine Horrocks'  \thoz,  {\mathrigid 2mu {$P(\uX,Y)\simeq P(\uX,0)$} over {$\gA[\uX,Y]$}},
  \item  we combine this last \iso with the \iso between $P(\uX,0)$ and $P(\uze,0)$ over the \ri $\gA[\uX]$ obtained by \hdrz.
\end{enumerate}

\vspace{-1em}\end{proof}

\vspace{-.7em}
\pagebreak	

\begin{corollary}\label{corthQUILIND}  $\!$\emph{(Quillen-Suslin \thoz,  Quillen's \dem\index{Quillen})}\\
If $\gK$ is a \cdi (resp.\,a \zed \riz),
every \mptf over $\KXr$ is free (resp.\,quasi-free).
\end{corollary}

\begin{proof}
The concrete Quillen induction 
applies with the class $\cF$ of \cdisz: note that $\KX$ is a Bézout domain, therefore the \mptfs over $\KX$ are free, and a fortiori extended. 
We pass to the \zedr \ris by the \elgbm \num2. 
Finally, for the \zed \risz, we use the \egt $\GKO(\gA)=\GKO(\Ared)$.
\end{proof}
\rems
1) Recall that a \zed \ri is connected \ssi it is local.
If $\gK$ is such a \riz, every \mptf over $\KXr$ is free.
\\
 2) The concrete Quillen induction 
applies to the Bézout domain of \ddk $\leq 1$ 
(see Exercise \ref{exoBézoutKdim1TransfertArX}) and more \gnlt to the \adps of dimension $\leq1$ (see \thref{thMaBrCo}). This generalizes the case of the \dDks obtained by Quillen. For the case of  regular \noe \ris of \ddk $\leq2$ (which we do not treat in this book), see \cite{Lam06}. 
\eoe

\subsubsection*{(Q3) versus (q3)}

The abstract Quillen induction 
(which does not provide any \cov results) presents the advantage of using a hypothesis (Q3) weaker than the hypothesis (q3) used in the concrete induction. We now explain how we can \cot 
recoup the situation,
even for the hypothesis~(Q3).

\sni{\it The free case.}

 In the case where the class $\cF$ is such that the \mptfs are free, we observe that the hypothesis (q3) is actually useless. Indeed, let~$P$ be a \ptf $\AX$-module and $S$ be the \mo of the \polus of $\AX$. 
Then, by (q1) the $\ArX$-module $P_S$ is free, so extended from~$\gA$. But then, by the Affine Horrocks'  \thoz, the module $P$ is extended from~$\gA$. In other words we have proved the following particularly simple version, for the free case.

\begin{theorem}\label{th2QUILIND}  $\!$\emph{(Concrete Quillen induction, free case)}\index{Quillen induction!concrete ---, free case}\\
Let $\cF$ be a class of \ris that satisfies the following \prtsz.
\begin{description}
  \item[{\rm (q0)}] If $\gA\in\cF$, every \ptf $\gA$-module  is free.
  \item[{\rm (q1)}] If $\gA\in\cF$, then $\ArX \in\cF$.
\end{description}
Then, for all $\gA\in\cF$ and all $r\geq1$, every \mptf over~$\AXr$  is free.
\end{theorem}

\vspace{-.5em}
\pagebreak	
\sni{\it The \gnl case.}

We would have noticed that the \prt  (Q2) does not intervene in the concrete Quillen induction: 
this hypothesis is rendered useless by the hypothesis (q3).

The \prt  (Q2) however intervenes when we want to replace (q3) with (Q3), which is a hypothesis a priori weaker than (q3).

We think that this weakening of the hypothesis is always possible in practice, without losing the \cov \crc of the result. However, this is based on the basic \lgbe machinery (\lgbe machinery with \idepsz), and as the latter is a \dem method and not strictly speaking a \thoz, we were not able to formulate our concrete induction directly with (Q3), because we wanted a \tho in due form.

Let us move on to the explanation of the replacement of the strong hypothesis (q3) with the weak hypothesis (Q3). 

We re-express the hypothesis (Q2) in the following more \gnl form.

{\rm (q2)} If $\gA\in\cF$ and $S$ is a \mo of $\gA$, then $\gA_S \in\cF$.

We suppose that (Q3) is satisfied in the following form: 
under the hypothesis that $\gA$ is a \alrd in the class $\cF$ 
we have a \cov \dem of the fact that every \mptf $P$ over $\AX$ is extended, 
which is translated into a computation \algo 
(for the \iso between $P$ and $P(0)$) based on the \prts of the class $\cF$ and on the disjunction
$$\preskip.2em \postskip.4em
a\in\Ati  \quad\hbox{or}\quad a\in\Rad(\gA)
$$
for the \elts $a$ that occur during the \algoz.
Under these conditions the basic \lgb machinery applies.\imlb
Consequently for a \mptf $P$ over $\AX$ for an arbitrary \ri $\gA\in\cF$, 
the \dem given in the local \dcd case, 
followed step by step, provides us with \moco $S_1$, \ldots, $S_\ell$ such that for each of them, 
the module $P_{S_i}$ (over $\gA_{S_i}[X]$) is extended from~$\gA_{S_i}$. 
Note that for this method to work, the considered class of \ris must satisfy (q2), and that we can limit ourselves to the \lons at \mos $\cS(\an;b)$. 
It then only remains to apply the Quillen \rcmz\index{Quillen patching} 
(\plgcz~\ref{thPatchQ}) 
to obtain the desired result: the module $P$ is extended from $\gA$.


\subsec{\`A la Suslin, Vaserstein or Rao}
\label{subsecSuslin}

The solution by \Sus to Serre's conjecture consists in showing that every stably free module over $\KXr$ is free (Serre had already proven that every \mptf over $\KXr$ is stably free), in other words that the kernel of every surjective matrix is free, or that every \vmd is the first column of an \iv matrix (see Fact~\ref{factStablib} and Proposition~\ref{corpropStabliblib}).

 If \rdb\label{NOTAfGg} $\cG$ is a subgroup of $\GL_n(\gA)$ and $A$, $B\in \Ae{n\times 1}$, we will write~\hbox{$A\sims{\cG }B$} to say that there exists a matrix $H\in\cG $ such that $HA=B$.
It is clear that this is an \eqvc relation.

Recall that a \vmd $f\in\Ae {n\times 1}$ is said to be completable if it is the first column vector of a matrix $G\in\GL_n(\gA)$. This amounts to saying that we have  
$$\preskip-.2em \postskip.4em 
f\;\sims{\GLn(\gA)}\;\tra{[\,1\;0\;\cdots\;0\,]}. 
$$
The goal in this subsection is therefore to obtain a \cov \dem of the following \thoz.

\begin{theorem}\label{thSuslinQS}  $\!$\emph{(\Susz)}\\
Every \vmd $f$ with \coos in $\KXr=\KuX$ (where $\gK$ is a \cdiz) is completable.
\end{theorem}

We will give three distinct \demsz, in chronological order.


\subsubsection*{First \dem}

Here we follow very closely Suslin's original \demz. 
We only have to get rid of a non\cov usage of a generic \idemaz, and have already done this work when we gave a \cov \dem of Suslin's lemma (Lemma~\ref{lemSuslin1}) in Chapter~\ref{chapPlg}.


\begin{fact}\label{factM2}
Let $M$, $N\in\MM_2(\gA)$. We have $\Tr(M)\,\I_2=M+\wi M$
and

\snic{\det(M+N)=\det (M) + \Tr(\wi{M}\,N) + \det( N).}
\end{fact}
\begin{proof}
For the matrices in $\MM_2(\gA)$, the map $M\mapsto\wi M$ is \linz, so

\snac{\arraycolsep2pt
\begin{array}{rcl}
 \det (M+N)\, \I_2\;=\; (\wi M +  \wi N)(M+N)&=&
\wi M M + (\wi M N + \wi N M)+ \wi N N\\
  & = & \big(\det (M) + \Tr(\wi{M}\,N) + \det( N)\big)\,\I_2.
\end{array}}
\vspace{-1mm}
\end{proof}
%
\begin{lemma}\label{lem02SusQS}
Let $B\in\MM_2(\gA),\,H=H(X)\in\MM_2(\AX)$, $\gB$ be an \Alg and $x\in\gB$.
Let $C(X)=B+XH$.
Suppose $\det C=\det B=a$.
By letting $S=\I_2+ x\wi H(ax)\,B $, we then have  $S\in\SL_2(\gA)$ and $S\wi B=\wi C(ax)$.
\end{lemma}
\begin{proof}
Fact~\ref{factM2} gives $\det (C)=\det (B) + X\,(\Tr(\wi{H}\,B) + X \det H)$,
and therefore
$$\preskip.3em \postskip.3em
E(X)=\Tr(\wi{H}\,B) + X \det H=0.
$$
Let $H_1=H(ax)$ and $C_1=C(ax)$.\\
We have then $S\wi B=\wi B+ x\wi{H_1} B\wi B=\wi B+ax\wi{H_1}=
\wi{C_1}$ and
$$
\preskip.3em \postskip.0em 
\begin{array}{rcl}
\det(S)&=&1+x\Tr(\wi{H_1} B)+\det(x\wi{H_1}B)\\
&=&1+x\Tr(\wi{H_1} B)+x^2a\det(H_1) = 1+x E(ax)=1.
\end{array} 
$$

\vspace{-1.4em}
\end{proof}

\vspace{-.5em}
\pagebreak	
\begin{lemma}\label{lem2SusQS}  $\!$\emph{(\Susz's lemma)}\index{Suslin's lemma}\\
Let $u,\,v\in\AX$, $a\in\gA\,\cap\,\gen{u,v}$, $\gB$ be an \Alg and $b,\,b'\in\gB$. \\
If $b\equiv b' \mod a\gB$,  then $
\Cmatrix{2pt}{u(b)\cr v(b)}  \sims{\SL_2(\gB)}
\cmatrix{u(b')\cr v(b')}$.
\end{lemma}
\begin{proof}
Let $p, q\in\AX$ such that $up+vq=a$ and $x\in\gB$ such that $b'=b+ax$.
Consider the matrix $M=\Cmatrix{2pt}{p&q\cr-v&u}\in\MM_2(\AX)$.
We apply Lemma~\ref{lem02SusQS} with the matrices $B=M(b)$ and $C(X)=M(b+X)$.\\
Note that the first column of $\wi B$ is \smashbot{$\Cmatrix{2pt}{u(b)\cr v(b)}$}
and that the first column of  $\wi C(ax)$ is 
$\Cmatrix{2pt}{u(b')\cr v(b')}$. $\phantom{\cmatrix{b\cr b}^{2}}$
\end{proof}
%

\begin{lemma}\label{lem3SusQS}
Let $f\in \AX^{n\times 1}$, $\gB$ be an \Alg and $\cG$ be a subgroup of~$\GLn(\gB)$, then the set
$$\preskip.0em \postskip.3em
\fa=\sotQ{a\in\gA}{\Tt b,\,b'\in\gB, \big((b\equiv b' \mod a\gB) \;\Rightarrow\;f(b)\sims{\cG}f(b')\big)}
$$
is an \id of $\gA$.
\end{lemma}
\facile
%

\begin{theorem}\label{th4SusQS}
Let $n\geq2$, $f$ be a \vmd of $\AX^{n\times 1}$ with~$f_1$ \monz, $\gB$ be an \Algz, and $\cG\subseteq \GLn(\gB)$ be the subgroup generated by $\En(\gB)$ and $\,\SL_2(\gB)$.\footnote{$\SL_2(\gB)$ is embedded in $\GLn(\gB)$ by the injection  $A\mapsto \Diag(A,\I_{n-2}$.}
 Then, for all $b,\,b'\in\gB$, we have $f(b)\sims{\cG}f(b')$.
\end{theorem}
\begin{proof}
It suffices to show that the \id $\fa$ defined in Lemma~\ref{lem3SusQS} contains $1$.
For an \elr matrix $E=E(X)\in\EE_{n-1}(\AX)$, we consider the vector
$$
\cmatrix{g_2\cr\vdots\cr g_n}\;=\;E\,\cmatrix{f_2\cr\vdots\cr f_n}.
$$
We will show that the resultant $a=\Res_X(f_1,g_2)$, which is well-defined since $f_1$ is \monz, is an \elt of $\fa$. We will therefore finish by invoking Suslin's lemma \ref{lemSuslin1}.\\
Let us therefore show that $a\in\fa$. We just use the fact that $a\in\gen{f_1,g_2}\cap\gA$.
We take $b,\,b'\in\gB$ with $b\equiv b'\mod a\gB$.
We want to reach $f(b)\sims{\cG}f(b')$.
Note that for $i\geq2$ we have
\begin{equation}\label{eqth4SusQS}
\begin{array}{ccc}
g_i(b')-g_i(b)\in \gen{b'-b}\subseteq \gen{a}  \subseteq \gen{f_1(b),g_2(b)},
\\[1mm]
\hbox{\cad } ~~~~
g_i(b')\in g_i(b) + \gen{f_1(b),g_2(b)}.
\end{array}
\end{equation}
We then have a sequence of \eqvcs
$$\hss \mathrigid1.5mu
\Cmatrix{.18em}{f_1(b)\cr f_2(b)\cr f_3(b)\cr\vdots\cr f_n(b)} \sims{E(b)}
\cmatrix{f_1(b)\cr g_2(b)\cr g_3(b)\cr\vdots\cr g_n(b)} \sims{\En(\gB)}
\cmatrix{f_1(b)\cr g_2(b)\cr g_3(b')\cr\vdots\cr g_n(b')} \sims{\SL_2(\gB)}
\cmatrix{f_1(b')\cr g_2(b')\cr g_3(b')\cr\vdots\cr g_n(b')} \sims{E(b')^{-1}}
\cmatrix{f_1(b')\cr f_2(b')\cr f_3(b')\cr\vdots\cr f_n(b')}
\hss.$$
The second is given by \Eqrf{eqth4SusQS}, the third by Lemma~\ref{lem2SusQS} applied to $u=f_1$ and $v=g_2$.
\end{proof}
%

\begin{corollary}\label{corth4SusQS}
Let $n\geq2$, $f$ be a \vmd of $\AX^{n\times 1}$ with~$f_1$ \mon
and $\cG$ be the subgroup of~$\GLn(\AX)$ generated by $\En(\AX)$ and $\,\SL_2(\AX)$. Then  $f\sims{\cG}f(0)$.
\end{corollary}
\begin{proof}
In \thref{th4SusQS}, we take $\gB=\AX$, $b=X$ and~$b'=0$.
\end{proof}
%

\begin{corollary}\label{corth5SusQS}
Let $\gK$ be a \cdiz,  $n\geq2$, $f$ be a \vmd of~$\KuX^{n\times 1}$, where $\KuX=\KXr$, and~$\cG\subseteq \GLn(\KuX)$ be the subgroup generated by $\En(\KuX)$ and $\,\SL_2(\KuX)$. Then  $f\sims{\cG} \tra{[\,1\;0\;\cdots\;0\,] }$.
\end{corollary}
\begin{proof}
If $f_1=0$, we easily transform the vector $f$ into $\tra{[\,1\;0\;\cdots\;0\,]}$ by \mlrsz.
Otherwise, a \cdv allows us to transform $f_1$ into a pseudo\polu in $X_r$ (Lemma~\ref{lemNoether}).
We can therefore assume that $f_1$ is \mon in $X_r$, we apply Corollary~\ref{corth4SusQS} with the \riz~$\gA=\gK[X_1,\ldots,X_{r-1}]$, and we obtain $f\sims{\cG} f(X_1,\ldots,X_{r-1},0)$.
We conclude by \recu on~$r$.
\end{proof}

We have indeed obtained \thref{thSuslinQS}, actually with an interesting precision over the group~$\cG$.


\subsubsection*{Second \demz}

We now closely follow a \dem by Vaserstein~\cite{Vaserstein2} such as it is presented in~\cite{Lam06} but by using constructive arguments.

More \gnlt we are interested in the possibility of finding in the \eqvc class of a vector defined over $\AX$ a vector defined over $\gA$, in a suitable sense.

We will use the following lemma.
\begin{lemma}
\label{lemfij}
Let $\gA $ be a \ri  and $f(X)=\tra{[\,f_1(X)\;\cdots\;f_n(X)\,]}$ be a \vmd in
$\AX^{n{\times}1}$, with $f_1$ \mon of degree $\geq1$. \\
Then, the \id $\fa=\rc(f_2)+\cdots +\rc(f_n)$  contains~$1$.
\end{lemma}
\begin{proof} 
We have $1= u_1f_1$ in $\gA\sur\fa$. This \egt in the \ri $(\gA\sur\fa)[X]$,
with~$f_1$ \mon of degree $\geq1$ implies that $\gA\sur\fa$
is trivial (by \recu on the formal degree of~$u_1$).
\end{proof}

\vspace{-.7em}
\pagebreak	
\begin{theorem}\label{th2HorrocksLocal} $\!$\emph{(Little Horrocks' local \thoz)}\\
Let $n\ge 3$ be an integer, 
$\gA$ be a \dcd \alo and 
$f(X)=\tra{[\,f_1(X)\;\cdots\;f_n(X)\,]}$ be a \vmd in $\AX^{n{\times}1}$, with $f_1$ \monz. 
Then
$$
f(X)=\cmatrix{ f_1   \cr \vdots\cr \vdots \cr f_n}
\sims{\En(\AX)}
\cmatrix{ 1 \cr 0\cr \vdots \cr 0}\sims{\En(\gA)}
\Cmatrix{2pt}{ f_1(0)   \cr \vdots\cr \vdots \cr f_n(0)}
.
$$
 \end{theorem}

\begin{proof}
Let $d$ be the degree of $f_1$. 
By \mlrrsz, we bring the \pols $f_2$, \ldots, $f_n$ to being of degrees $< d$. 
Let $f_{i,j}$ be the \coe of~$X^j$ in $f_i$. 
The vector $\tra{[\,f_1(X)\;\cdots\;f_n(X)\,]}$ remains \umdz. 
If $d=0$, we are done. 
Otherwise given Lemma~\ref{lemfij} and since the \ri is local, one of the $f_{i,j}$'s for $i\in\lrb{2..n}$ is a unit. 
Suppose for example that $f_{2,k}$ is \ivz. 
We will see that we can find two \pols $v_1$ and $v_2$ such that the \pol $g_2=v_1f_1+v_2f_2$ is \mon of degre $d-1$.
If $k=d-1$, this works with $v_1=0$ and $v_2$ constant.
If $k<d-1$, consider the following disjunction
$$\preskip.4em \postskip.4em 
f_{2,d-1}\in\Ati\; \hbox{ or }\; f_{2,d-1}\in
\Rad(\gA ). 
$$
In the first case, we are reduced to $k=d-1$. In the second case the \pol $q_2=Xf_2-f_{2,d-1}f_1$ is of degree $\le d-1$ and satisfies: $q_{2,k+1}$ is a unit. We have gained some ground: it now suffices to iterate the process.
\\ 
Now we therefore have $g_2=v_1f_1+v_2f_2$  of degree $d-1$ and \monz. So we can divide $f_3$ by $g_2$ and we obtain $g_3=f_3-g_2q$ of degree $<d-1$ ($q\in\gA $), so the \pol
$$\preskip.3em \postskip.4em 
h_1=g_2+g_3=f_3+g_2(1-q)=f_3+(1-q)v_1f_1+(1-q)v_2f_2 
$$
 is \mon of degree $d-1$.
Thus, by an \mlrr we were able to replace
$\tra{[\,f_1\;f_2\;f_3\,]}$ with $\tra{[\,f_1\;f_2\;h_1\,]}$, with $h_1$ \mon of degree~$d-1$.
We can therefore by a sequence of  \mlrrs bring 
the vector $\tra{[\,f_1(X)\;\dots\;f_n(X)\,]}$, with $ f_1 $ \mon of degree $d$, to
$$\preskip.3em \postskip.4em 
\tra{[\,h_1(X)\;\dots\;h_n(X)\,]} \hbox{  with }  h_1  \hbox{ \mon of degree }  d-1. 
$$
We obtain the desired result by \recu on~$d$.
\end{proof}

 {\bf Terminology.} We consider a \sys of formal \pols $(f_i)$ with $\deg f_i=d_i$.
We then call the \gui{head \id of the \sys $(f_i)$} the \id of the formally leading \coes of the~$f_i$'s.

\penalty-2500
\begin{theorem}\label{th2HorrocksGlobal}
 $\!$\emph{(Little Horrocks' global \thoz)}\\
Let $n\ge 2$ be an integer, $\gA $ be a \ri and
$f\in \AX^{n{\times}1}$ be a \vmdz.
Suppose that the head \id of the $f_i$'s contains $1$.
Then
$$
f(X)=\cmatrix{ f_1   \cr \vdots \cr f_n}
\sims{\GLn(\AX)}
\cmatrix{ f_1(0) \cr \vdots \cr f_n(0)}=f(0).
$$
\end{theorem}
\begin{proof}
The case $n=2$ is an exception: if $u_1f_1+u_2f_2=1$, the \egt 

\snic{\cmatrix{u_1&u_2\cr -f_2&f_1}
\cmatrix{f_1\cr f_2}=\cmatrix{1\cr 0}}

 gives the required matrix, in $\SL_2(\AX)$.
\\ 
For $n\geq3$, we apply the basic \lgbe machinery (\paref{MethodeIdeps}) with the \cov \dem of \thref{th2HorrocksLocal}.\imlb 
We obtain a finite family of \mocoz, $(S_i)_{i\in J}$  in $\gA$, such that for each $i$ we have
$
f(X)\sims{\En(\gA_{S_i}[X])} f(0).
$
We conclude with the Vaserstein \rcm for the \eqvcs of matrices on the left-hand side (item~\emph{2} of the \plgrf{thPatchV}).
\end{proof}

 {\bf Conclusion.}
We have just obtained a (slightly weaker) variant of Corollary~\ref{corth4SusQS}, and this gives the \dem of Suslin's \thoz~\ref{thSuslinQS} in the same way as in the first solution.

\mni\comm The little Horrocks' global \tho can also be obtained as a consequence of the \gui{grand} Horrocks' global \tho \vref{thHor}.
Let $P=\Ker\tra{f(X)}$. By localizing at $f_1$, $P$ becomes free.
The Affine Horrocks' \tho tells us that $P$ is free, which means that {$f(X)\sim \tra{[\,1\;0\;\cdots\;0\,]}$} over~$\GLn(\AX)$.
\eoe


\subsubsection*{Third \demz}

We now closely follow a \dem by Rao.
This time we will not need any \recu on the number of variables to reach Suslin's \thoz.

\begin{lemma}\label{lemUMD}
We consider a vector $x=(x_1,\ldots ,x_n)\in\Ae n$ and $s\in\gA$. If $x$ is \umd over $\aqo{\gA}{s}\!$ and over $\gA[1/s]$, it is \umdz.
 \end{lemma}
\begin{proof}
Let $\fa=\gen{x_1,\ldots,x_n}$. We have $s^r\in\fa$ (for a certain $r$) and~\hbox{$1-as\in\fa$} (for a certain $a$). We write $1=a^rs^r+(1-as)(1+as+\cdots) \in \fa$.
\end{proof}
%

\begin{lemma}\label{lem1Rao}
Let $n\ge 2$ be an integer, $\gA $ be a \riz, and $f$ be a \vmd in $\AX^{n{\times}1}$: $f=\tra{[\,f_1(X)\;\cdots\;f_n(X)\,]}$.
For each $f_i$ of formal degree $d_i$, let $f_i\sta$ be the \emph{formal reciprocal \pol} $X^{d_i}f_i(1/X)$. Let $f\sta(X)=\tra{[\,f_1\sta(X),\;\cdots\;f_n\sta(X)\,]}$.
If $f\sta(0)$ is \umdz, the same goes for~$f\sta$.
\end{lemma}
\begin{proof}
By Lemma~\ref{lemUMD}, it suffices to prove that $f\sta(0)$ is \umd (it is true by hypothesis) and that $f\sta$ is \umd over $\gA[X,1/X]$, which comes from the \egt $\sum_iu_i(1/X)X^{-d_i}f_i\sta=1$ (where $\sum_iu_if_i=1$ in~$\AX$).
\end{proof}
%

\begin{theorem}\label{th1Rao} $\!$\emph{(Rao's \thoz, \cite{Rao85b})}\\
Let $n\ge 2$ be an integer, $\gA $ be a \riz, and $f=\tra{[\,f_1(X)\;\cdots\;f_n(X)\,]}$ be a \vmd in $\AX^{n{\times}1}$,
with $1$ in the head \id of the $f_i$'s. Then
$$\preskip.2em \postskip.2em
f\sims{\GLn(\AX)} f(0)\sims{\GLn(\gA)}f\sta(0)\sims{\GLn(\AX)} f\sta.
$$
If in addition one of the $f_i$'s is \monz, we have $f\sims{\GLn(\AX)} \tra[\,1\;0\;\cdots\;0\,]$.
\end{theorem} 
\begin{proof}
Since $f\sim f(0)$ by the little Horrocks' global \thoz, we deduce that $f\sim f(1)$. In addition,  $f\sta(0)$ is \umd therefore $f\sta$ is \umd (by Lemma~\ref{lem1Rao}).
Moreover, $1$ is in the head \id of the $f_i\sta$'s, which allows us to apply to $f\sta$ the little Horrocks' global \thoz.\\
We conclude:   $f\sim f(0)\sim f(1)= f\sta(1)\sim f\sta$.
\end{proof}

\comm
The same result is valid by replacing $\GLn$ by $\En$, but the \dem is strictly more delicate (see
\thref{th2Rao}).
\eoe

\mni {\bf Conclusion.} We then obtain Suslin's \tho \rref{thSuslinQS}
as follows.
We take for $\gA $ the \ri $\gK[X_1,\ldots,X_{r-1}]$ and we make a \cdv that renders pseudomonic one of the \polsz.
\\ 
Thus, 
\begin{itemize}
\item  on the one hand, the solution is much more \gui{efficient} than in the first two \dems since there is no longer an \recu on the number of variables,
\item  and on the other hand, the \tho is much more \gnlz.
\end{itemize}


\section[\Pro modules extended from  \anarsz]
{\Pro modules extended from  valuation or \ari \risz}  \label{sec.Etendus.Valuation} 

Recall that a \adv is a reduced \ri in which we have, for all $a,\,b$, $a$ divides $b$ or $b$ divides $a$. It is a normal, local \ri \sdzz.
\perso{Il is int\`egre \ssi il is \cohz, mais cela semble sans importance ici.
}

We begin with a useful result regarding  \advs and the \ddk
(we can also refer back to Exercise~\ref{exoPrufNagata}).

\begin{lemma}\label{lemV-V(X)} 
If  $\gA$ is a \advz, then so is $\gA(X)$.
If~$\gA$ is a \adv of finite \ddkz, then $\gA(X)$ has the same \ddkz.
\end{lemma}
\begin{proof}
If $\gA$ is a \advz, every $f\in\AX$ is expressible in the form $f=ag$ with $a\in\gA$ and $g\in\AX$ which admits a \coe equal to $1$.
In particular,~$g$ is \iv in $\gA(X)$.
If  $F_1=a_1g_1/u_1$ and  $F_2=a_2g_2/u_2$ are two arbitrary \elts of $\gA(X)$ (with $a_i\in\gA$ and $g_i,u_i$ primitive in $\AX$), then $F_1$ divides $F_2$ in $\gA(X)$ \ssi $a_1$ divides $a_2$ in $\gA$. Therefore \gui{the \dve is identical in $\gA$ and $\gA(X)$} and $\gA(X)$ is a \advz.
In addition, since the \itfs are principal, the canonical \homo $\ZarA\to\Zar\gA(X)$ is an \iso of \trdis (note: these are totally ordered sets),
which implies that the \ddk is the same. 
\end{proof}
%

\subsec{The univariate case} \label{subsec.Bass.Valuation} 

This subsection is devoted in the most part to a \cov \dem of the following Bass' \thoz.

\begin{theorem}\label{thBass.Valuation} 
 If $\gV$ is a \adv of finite \ddkz, every \ptf $\VX$-module is free.
\end{theorem}

We will actually prove slightly stronger variants: we can get rid off the hypothesis on the \ddkz,
and we have a version for \anarsz.

We start with a simple example.
\subsubsection*{A simple example}
\label{subsecZXfree}

\begin{proposition}\label{propZXfree}
Every \mptf over $\ZZX$ is free.
\end{proposition}
\begin{proof}
Let $M$ be a \ptf $\ZZX$-module.
First of all note that if $M$ is of rank $1$, it is free because $\ZZX$ is a GCD-domain (Lemma~\ref{lemPicGcd}). \\
Now suppose that $M$ is of rank $r>1$.
If we extend the scalars to $\QQX$, the module becomes free.
Therefore there exists an integer $d>0$ such that $M$ becomes free over $\ZZ[1/d][X]$. 
If $d=1,$ nothing needs to be done. Otherwise, let~\hbox{$p_1$, \ldots, $p_k$} be the prime factors of~$d$. \\
The \mos $d^{\NN}$, $1+p_1\ZZ$, \ldots,  $1+p_k\ZZ$ are \com (see the fundamental example on \paref{explfonda}). It therefore suffices to show that the modules $M_{1+p_i\ZZ}$ are free (therefore extended),
because then the Quillen \rcmz\index{Quillen patching}  implies that $M$ is extended from $\ZZ$, therefore free. \\
Let $p$ be any of the $p_i$'s.  Since $\ZZ_{1+p\ZZ}[X]$ is $2$-stable (lemma below), by applying \SSO \tho (\thref{thSerre}), 
we obtain an \iso $M_{1+p\ZZ}\simeq \ZZ_{1+p\ZZ}[X]^{r-1}\oplus N$, with $N$ being a \pro $\ZZ_{1+p\ZZ}[X]$-module of constant rank~$1$.
By the initial remark (which applies by replacing~$\ZZ$ by $\ZZ_{1+p\ZZ}$), $N$ is free, so $M$ is free.
\end{proof}

\vspace{-.8em}
\pagebreak	

\begin{lemma}\label{lemZ1+pZ} 
The \ri $\ZZ_{1+p\ZZ}[X]$ is $2$-stable.
\end{lemma}
\begin{proof}
Consider the partition of $\Spec(\ZZ_{1+p\ZZ}[X])$ attached to $\so{p}$:
more \prmtz,
the \ri $\ZZ_{1+p\ZZ}[X]$ is replaced by the two \risz

\snic{\ZZ_{1+p\ZZ}[X][1/p]\simeq\QQX$ and
$\aqo{(\ZZ_{1+p\ZZ}[X])}{p}\simeq \FF_{p}[X],}

which are of \ddk $1$.
\Thref{thPartitionSpec} then tells us that $\ZZ_{1+p\ZZ}[X]$ is $2$-stable.
\end{proof}

\rem
Actually the use of prime factors of $d$, although intuitively natural, introduce an unnecessary complication.
Indeed, the \mosz~$d^{\NN}$ and $1+d\ZZ$ being \comz, it suffices to prove that $M_{1+d\ZZ}$ is free. As $\ZZ_{1+d\ZZ}[X]$ is a GCD-domain, the previous reasoning applies if we know how to show that $\ZZ_{1+d\ZZ}[X]$ is $2$-stable. The \dem of Lemma~\ref{lemZ1+pZ} works by replacing $p$ by $d$, 
because $\ZZ_{1+d\ZZ}[X][1/d]\simeq\QQX$ and~$\aqo{\ZZ_{1+d\ZZ}[X]}{d}\simeq (\aqo{\ZZ}{d}\!)[X]$, which are of \ddk $1$ ($\aqo{\ZZ}{d}$ is \zedz).
\eoe

\subsubsection*{A more elaborate example}
\label{subsec2ZXfree}

Given the previous remark we leave the \dem of the following \gnn  to the reader.

\begin{proposition}\label{propBezIntXfree}
Let $\gA$ be an integral \ri of \ddkz~$\leq 1$, $d$ be an \elt of $\Reg(\gA)$, and $M$ be a
\ptf $\AX$-module.
\begin{enumerate}
\item $\gA_{1+d\gA}[1/d]=\Frac\gA$ is \zedz.
\item $\aqo{\gA_{1+d\gA}}{d}\simeq\aqo{\gA}{d}$ is \zedz.
\item $\gA_{1+d\gA}[X]$ is $2$-stable.
\item
\begin{enumerate}
\item  If $\gA$ is a Bézout \riz,  $M$ is free.
\item  If $\gA$ is seminormal,  $M$ is extended from $\gA$.
\end{enumerate}
\end{enumerate}
\end{proposition}

\subsubsection*{An example in finite \ddk $>0$}
\label{subsec1VXfree}

Let $\gV$ be an integral \ri with some \elts $a_1$, \ldots, $a_k$. Suppose

\snic{\rD_\gV(a_1)\leq \rD_\gV(a_2)\leq \cdots \leq \rD_\gV(a_k).}

The partition in constructible subsets of $\Spec\gV$ 
associated with this family contains only $k+1$ \elts

\snic{\rD_\gV(a_1),\; \rD_\gV(a_2)\setminus \rD_\gV(a_1),\; \dots,\;
\rD_\gV(a_k)\setminus \rD_\gV(a_{k-1}),\;
\rD_\gV(1)\setminus \rD_\gV(a_{k}),}

that correspond to the \risz

\snic{\gV[1/a_1]$, $(\aqo{\gV}{a_1}\!)[1/a_2]$, \ldots, $(\aqo{\gV}{a_{k-1}})[1/a_k]$
and $\aqo{\gV}{a_k}.}

Now suppose that these \ris are all \emph{\zedsz}. Then, we similarly have a partition into $k+1$ constructible subsets
 of $\Spec\VX$ and the corresponding \risz

\snic{\gV[1/a_1][X]$, $(\aqo{\gV}{a_1})[1/a_2][X]$, \ldots,
$(\aqo{\gV}{a_{k-1}})[1/a_k][X]$ and $(\aqo{\gV}{a_k})[X]}

 are all of \ddkz~$\leq 1$.
\Thref{thPartitionSpec} then tells us that~$\gV[X]$ is $2$-stable.
Therefore if $M$ is a \pro $\VX$-module of constant rank $r$, by \SSO \thoz, we obtain $M\simeq \VX^{r-1}\oplus N$, with~$N$ of constant rank $1$.

If $\gV$ is in addition a seminormal \ri (resp.\,a GCD-domain), then  $N$ is extended from $\gV$ (resp.\,then $N$ is free), so $M$ is extended from $\gV$ (resp.\,$M$ is free).

\smallskip The result \gui{$\VX$ is $2$-stable} is satisfied when $\gV$ is a valuation domain of \ddkz~$k$ for which we have a precise sufficient knowledge of the valuation group.
In \clama (with \TEM but without using the \ideps or the axiom of choice) we therefore already obtain the desired Bass' \tho for valuation domains of finite Krull dimension.

However, the result is not of an \algq nature if we do not know how to compute some suitable \elts $a_i$.

This difficulty will be bypassed dynamically.

\subsubsection*{\Cov \dem of Bass' \thoz}
\label{subsec2VXfree}

We need to establish the following theorem.

\begin{theorem}\label{thVX2stab}
If $\gV$ is a \advz, 
$\VX$ is $2$-stable.
\end{theorem}

We start with the following lemma (the \dem of the \tho is postponed until \paref{DemothVX2stab}).

\begin{lemma}\label{lemVX3stab}
Let $\gV$ be a \adv and $\gV'$ be the \sdv of $\gV$
generated by a finite  family of \elts of $\gV$.
Then, $\gV'[X]$ is $2$-stable.
\end{lemma}
\begin{proof}
Let $\gV_1$ be the sub\ri of $\gV$ generated by the finite family.
We define
$$\preskip.2em \postskip.4em
\gV'=\sotq{c/b}{c,b\in\gV_1,\;\mathrm{regular}\;b\;\mathrm{divides}\; c\; \mathrm{in}\; \gV}\subseteq\Frac(\gV_1).
$$
It is easily seen that $\gV'$ is a valuation domain. Moreover, since $\Kdim\gV_1\leq m$ for some $m$ (see Lemma~\ref{lemahbonvraiment}), we have also $\Kdim\gV'\leq m$.
Indeed, consider a sequence $(z_1,\dots,z_{m+1})$ in $\gV'$, write  $z_i=x_i/b$ with a common denominator $b$, and introduce a \cop sequence of  $(x_1,\dots,x_{m+1})$ in $\gV_1$.
A fortiori it is \cop in  $\gV'$.\\
Let $\ell_1$, $\ell_2$ and $a$ in $\gW=\gV'[X]$.
Let $L=(\ell_1,\ell_2)$ and $Q=(q_1,q_2)$.
We are searching for $q_1$, $q_2\in\gW$ that satisfy $\rD_{\gW}(a,L)=\rD_{\gW}(L+aQ)$.
If $\gV'$ were a \cdiz, we would have an \algo to compute $Q$ from $L$. By executing this \algoz, we would use the test \gui{$y=0$ or $y$ \ivz?} for some \elts $y\in\gV_1$ that occur during the computation (indeed, in the case where $\gV'$ is a \cdiz, some $y/z$ in $\gV'$ is null if $y$ is null, \iv if $y$ is \ivz, $z$ having been already certified as \ivz).\\
We can transform the \algo dynamically by replacing each test \gui{$y=0$ or $y$ \ivz?} by the splitting of \gui{the current \ri $\gA$,} which gives the two \ris $\gA[1/y]$ and ${\gA}\sur{\DA(y)}$.\\
At the beginning $\gA=\gV'$.
As in $\gV'$ the \elts are comparable with respect to \dvez, all the introduced \ris can be brought back to the standard form ${\gV'}\sur{\rD_{\gV'}( y_i)}\![1/y_{i-1}]$ ($i\in\lrb{2..k}$) for a finite family $(y_i)_{i\in\lrbk}$ of $\gV_1$,
with $y_{i-1}$ dividing $y_i$ in $\gV'$ for $i\in\lrb{2..k}$.
\\
Here we might have the impression of having succeeded insofar as we could say that: we now apply Lemma~\ref{lemPartitionSpec}.
\\
However, by reading the \dem of this lemma, we see that during a splitting $\gB\mapsto(\gB[1/b],\aqo{\gB}{b}\!)$, first the given $L$ and $a$ produce some~$Q$ for  $\aqo{\gB}{b}$, then $L+aQ$ and $ab$ produce some $R$ for $\gB[1/b]$, the final result being that $Q+bR$ suits for $L$ and $a$ in $\gB$.\\
Thus the dynamic of our transformed \algo must be more carefully controlled.%
\footnote{Otherwise,  the lemma could actually be proven without any hypothesis on $\gV$.}
What saves us is that in our dynamic use of Lemma~\ref{lemPartitionSpec}, the computations that start with $L$ and $a$ remain entirely in $\gV'\subseteq\Frac(\gV_1)$. Consequently, we can be certain not to fall into an infinite loop where the number of \ris ${\gV'}\sur{\rD_{\gV'}( y_i)}\![1/y_{i-1}]$ would increase indefinitely, which would prevent the \algo from halting.
Indeed, if $k> m$ (where $\Kdim\gV'\leq m$), the sequence $(y_1,\dots,y_k)$ is singular, and since $\rD_{\gV'}(y_{i-1})\leq\rD_{\gV'}(y_{i})$ and $\Zar\gV'$ is totally ordered, Lemma \ref{lemSeqSingTD} tells us that one of the following three situations occurs
 \begin{itemize}
\item $\rD_{\gV'}(y_{1})=0$, in which case the \riz~${\gV'}\sur{\rD_{\gV'}( y_2)}\![1/y_{1}]$ is trivial and the list is shortened by deleting $y_1$,
\item $\rD_{\gV'}(y_{m+1})=1$, in which case the \riz~${\gV'}\sur{\rD_{\gV'}( y_{m+1})}\![1/y_{m}]$ is trivial and the list is shortened by deleting $y_{m+1}$,
\item for some $i\in\lrb{2,m+1}$, we have the \egt $\rD_{\gV'}(y_{i-1})=\rD_{\gV'}(y_{i})$, in which case the \riz~${\gV'}\sur{\rD_{\gV'}( y_{i})}\![1/y_{i-1}]$ is trivial and the list is shortened by deleting $y_{i}$.
\end{itemize}
\vspace{-1em}\end{proof}

\rem  Thus, once $\gV_1$ is fixed, the \ri $\gV'$ behaves, with regard to the $2$-stability of $\gV'[X]$ as the \ri of finite \ddk \gui{$>0$ but entirely controlled} which was given in the previous subsection: the sequence of the $y_i$'s, limited to $m$ terms, behaves like the sequence of the $a_i$'s of the previous subsection, except that the $y_i$'s are produced dynamically as the \algo executes whereas the $a_i$'s were given at the start.
\eoe

\begin{Proof}{\Demo of \thref{thVX2stab}. } \label{DemothVX2stab}
Let $\ell_1$, $\ell_2$ and $a$ in $\VX$.
We search for $q_1,q_2\in\gV[X]$ satisfying $\rD_{\VX}(a,L)=\rD_{\VX}(L+aQ)$
(with $L=(\ell_1,\ell_2)$ and $Q=(q_1,q_2)$).
We apply Lemma~\ref{lemVX3stab} with the finite family constituted by the \coes of $\ell_1$, $\ell_2$ and $a$. We find $q_1,q_2$ in $\gV'[X]\subseteq\VX$.
\end{Proof}
%

\begin{theorem}\label{thBassValu} $\!$\emph{(Bass-Simis-Vasconcelos)}
If $\gV$ is a \advz,  every \ptf $\VX$-module is free.
\end{theorem}
\begin{proof}
Let  $M$ be a \mptf over $\VX$. Since $\VX$ is connected, $M$ has a constant rank $r\in\NN$.
Since $\VX$ is $2$-stable, \SSO \tho gives us that $M\simeq \VX^{r-1}\oplus N$, where $N$
is a \pro $\VX$-module of constant rank $1$.
It remains to show that $N\simeq\VX$.\\
If $\gV$ is integral we finish like this: since $\VX$ is a GCD-domain, $N\simeq \VX$.
In \gnl we can say: $\gV$ is normal, therefore every \mrc $1$ over $\VX$ is extended from $\gV$.
However, $\gV$ is local, in conclusion $N$ is free over $\VX$.
\end{proof}
%

\subsubsection*{The case of  \anarsz}
\label{subsec4VXfree}
\begin{theorem}\label{thBassAri} $\!$\emph{(Bass-Simis-Vasconcelos)}
If $\gA$ is an \anarz,  every \ptf $\AX$-module is extended from $\gA$.
\end{theorem}
\begin{proof}
First of all, since $\GKO(\gA)=\GKO(\Ared)$ and $\AX\red=\Ared[X]$, it suffices to prove the \tho in the reduced case, \cad for the \adpsz.
\\
Consider a \ptf $\AX$-module $M$.\\
In \clama we would apply the Quillen \rca \thoz: a \mptf over $\AX$ is extended because it is extended if we localize at an arbitrary \idep of $\gA$ (the \ri becomes a \advz).
\\
In \comaz, we rewrite the \prco given in the local case (for \thref{thBassValu}) by applying the basic \lgbe machinery.\imlb
\\
More \prmtz, suppose that in the local case (\cad for a \advz) we use the disjunction \gui{$a$ divides $b$ or $b$ divides $a$.} Since we are dealing with a \adpz, we know $u$, $v$, $s$, $t$ such that \hbox{$s+t=1$, $sa=ub$} and~$tb=va$.
If $\gB$ is the \gui{current} \riz, we consider the two \come \lons $\gB[1/s]$ and $\gB[1/t]$. In the first, $a$ divides $b$, and in the second,~$b$ divides $a$.
\\
Ultimately we obtain a finite family $(S_i)$ of \moco of $\gA$ such that after \lon at any of the $S_i$'s, the module $M$ becomes free, therefore extended. We conclude with the Quillen \rcm (\plgc \ref{thPatchQ}).
\end{proof}

\rems 1) We did not need to assume that the \adv was \dcd to make the \cov \dem of \Thosz~\ref{thVX2stab}, \ref{thBassValu} and \ref{thBassAri} work. This is especially translated by the fact that in the last \demz, the \moco are based on the disjunction (in a \aloz) \gui{$s$ or $1-s$ is \ivz} and are directly given by \ecoz.

 2)
In this type of passage from the local to the global, 
to make sure that the \algo halts, we have to make sure that the version given in the local case is \gui{uniform,} 
meaning that its execution is done in a number of steps that is bounded by a function of the discrete parameters of the input: 
the size of the matrix and the degrees of its \coesz. 
This is indeed the case here, modulo the \dem of Lemma~\ref{lemahbonvraiment}. 
Note that the fact that the \algo in the local case does not use any tests of \egt to $0$ greatly simplifies life and helps us to appreciate the validity of its dynamic implementation in the passage from the local to the global.
\eoe

\subsec{The multivariate case} 
\label{sec.LS} 
This subsection is devoted to the \cov \dem of the following Lequain-Simis \thoz.

\THo{(Lequain-Simis)}
{If $\gA$ is an \anarz, every \mptf over $\AXr$ is extended from $\gA$.
}

\medskip 
{\bf A dynamic comparison between the \ris $\gA(X)$ and $\ArX$} \label{dyn}

\smallskip 
In the following \thoz, we prove that for a \ri $\gA$ \ddi$d$, the \ri $\ArX$ dynamically behaves like the \ri $\gA(X)$ or like a \lon of a \ri  $\gA_S[X]$ for a \mo $S$   of $\gA$ with $\Kdim\gA_S\leq d-1$.

\begin{theorem}  $\!$\emph{(Dynamic comparison of $\gA(X)$ with $\ArX$)}
\label{compa} \\
Let $\gA$ be a \riz, $f=\sum_{j=0}^ma_jX^j \in \AX$ be a primitive \polz, and,   for $j\in\lrbm$,  $S_j=\cS_\gA^\rK(a_j)=a_j^{\NN}(1+a_j\gA)$ be the Krull boundary \mo of~$a_j$ in $\gA$.
Then, the \mos $f^{\NN}$, $S_1$, \ldots, $S_m$ are \com in $\ArX$.\\
In particular,  if $\Kdim \gA$ and $d\geq 0$, each \ri $\ArX_{S_j}$ is a \lon of a \ri $\gA_{S_j}[X]$ with
$\Kdim\gA_{S_j}\leq d-1$.
\end{theorem}

\begin{proof}
For $x_1$, \ldots, $x_m \in \gA$ and $n$, $d_1$, \ldots, $d_m \in \NN$, we must show that the following \elts of $\gA[X]$

\snic{
f^n,\; a_m^{d_m}(1 - a_mx_m),\; \dots, \;
a_1^{d_1}(1 - a_1x_1),}

generate an \id of $\gA[X]$ that contains a \poluz. We reason by \recu on $m$; it is obvious for $m = 0$ because $a_m = a_0$ is \ivz. \\
For $m \ge 1$ and $j\in \lrb{1..m-1}$, let 

\snic{a = a_m, \;\;x = x_m, \;\;d = d_m\;\;\hbox{and}\;\;
a'_j = a_j^{d_j}(1 - a_jx_j).}

Consider the quotient $\gB=\aqo{\gA}{a^d(1 - ax)}$; we must show that the family

\snic{\cF = (f^n, \; a'_{m-1},\; \dots,\; a'_1 )
}

generates an \id of $\BX$ which contains a \poluz.  
\\
Since $a^d(1 - ax) = 0$, $e = a^d x^d$ is an \idm and $\gen{e} = \gen
{a^d}$. \\
Let $\gB_e \simeq \aqo{\gB}{1-e}$ and $\gB_{1-e} \simeq \aqo{\gB}{e}$.
It suffices to show that $\gen {\cF}_{\gB_e[X]}$ and $\gen {\cF}_{\gB_{1-e}[X]}$ contain a \poluz.
\\
 In $\gB_e[X]$, it is \imd because $a$ is \ivz.  In $\gB_{1-e}[X]$, we have $a^d = 0$. Let $f = aX^m +
r$ with $r = \sum_{j=0}^{m-1} a_j X^j$.  
In $\gB$, for every exponent~$\delta$, the \elts of  $(a^\delta, a_{m-1}, \ldots, a_1, a_0)$ are comaximal. For $\delta = d$, we deduce that in $\gB_{1-e}[X]$, the \pol $r$ is primitive. Since $r = f - aX^m$ and $a^d = 0$, we have $r^d \in \gen {f}$ so $r^{dn} \in \gen {f^n}$. \\
We apply the \hdr to the \pol $r \in \gB_{1-e}[X]$ of (formal) degree $m-1$: the \id $\gen{r^{dn}, a'_{m-1}, \cdots, a'_1}$ of $\gB_{1-e}[X]$ contains a \poluz; therefore the same holds for the \id $\gen{f^n, a'_{m-1}, \cdots, a'_1}$.
\end{proof}

\rem The previous \tho seems to have fallen from the sky as if by magic.
Actually it is the result of a slightly complicated story. 
In the article~\cite{ELY07}, the following \tho was proved by starting with the special case of a \dcd \aloz, then by generalizing to an arbitrary \ri by means of the basic \lgbe machinery.\imlb

\smallskip 
{\bf Theorem.} \emph{Let $\gA$ be a \ri such that $\Kdim \gA\leq d\in\NN$.
Let $f \in \AX$  be a primitive \polz. There exist \moco $V_1$, \ldots, $V_\ell $ of $\ArX$ such that for each $i\in \lrbl$, either $f$ is \iv in  $\ArX_{V_i}$,  or~$\ArX_{V_i}$ is a \lon of an $\gA_{S_i}[X]$ with $\Kdim\gA_{S_i}<d$.}

By explicating the \algo contained in the \dem of this \thoz, we have obtained \thref{compa}.
\eoe

\rdb \label{MachDynAxAx}
\medskip 
{\bf Dynamic machinery with $\ArX$ and $\gA(X)$}

\smallskip 
The previous \tho allows us to implement a dynamic machinery of a new type.

Suppose that we have established a \tho for the \advs of \ddk $\leq n$.
We want the same \tho for the \risz~$\ArX$ when $\gA$ is a \adv of \ddk $\leq n$.

Suppose also that the \prt to be proven is stable by \lon and that it comes from a \plgcz.

We perform a \dem by \recu on the \ddkz. When the \ddk is null, $\gA$ is a \cdi and we have $\ArX=\gA(X)$, which is also a \cdiz, therefore the \tho applies.

Let us look at the passage from the dimension $k$ to the dimension $k+1$ ($k<n$).
Notice that $\gA(X)$ is a \adv with the same \ddk as~$\gA$ (Lemma~\ref{lemV-V(X)}). We assume that $\Kdim\gA\leq k+1$.
We have a \prco of the \tho for the \advs of \ddk $\leq n$, in particular it works for $\gA(X)$.
We try to make this \dem (\cad this \algoz) work with $\ArX$ instead of~$\gA(X)$.
This \dem uses the fact that in $\gA(X)$ the primitive \pols of~$\AX$ are \ivz. Each time that the initial \dem uses the \inv of such a \pol $f$, we call upon \thref{compa}, which replaces the \gui{current} \ri by \come \lonsz. In the first \lon the \pol $f$ has been inverted, and the \dem can be continued as if $\ArX$ were $\gA(X)$. In each of the other \lons we have replaced $\ArX$ by a localization
 of a \ri $\gA_{S_i}[X]$ with $\Kdim \gA_{S_i}\leq k$, and, \emph{if we are lucky}, the \hdr allows us to conclude.

Ultimately we have proven the \tho for  \come \lons of~$\ArX$. Since the conclusion stems from a \plgcz, we have proven the \tho for $\ArX$.

\subsubsection*{Application to the \tho of Maroscia and Brewer{ \& }Costa}

The dynamic machinery explained in the previous subsection applies for the first of the following results.
\begin{enumerate}
  \item [(i)] \emph{If $\gA$ is a \adv with $\Kdim\gA \leq 1$, then $\ArX$ is a \adp with $\Kdim\ArX \leq 1$.} \\
Indeed, it suffices to prove the conclusion \lot (here, after \lon of $\ArX$ at \mocoz).
However, \thref{compa} allows us to split the \ri $\ArX$ into components that behave (for the computation to be done) either like $\gA(X)$, or like a localized \ri of a~$\KX$ 
where $\gK$ is \zedrz. In the two cases we obtain a \adp of \ddk $\leq1$.
  \item [(ii)] \emph{If $\gA$ is a \adp with $\Kdim\gA \leq 1$, then so is $\ArX$.}\\
Indeed, it suffices to prove the conclusion \lot (here, after \lon of $\gA$ at \mocoz). We apply the \lgbe machinery of  \anars to the \dem of item~(i): the \ri $\gA$ is subjected to \come \lonsz, in each of which it behaves like a \advz.\imla
\end{enumerate}

As a consequence we obtain a special version
of the Lequain-Simis \tho by using the concrete Quillen induction\index{Quillen induction!concrete ---} (\thref{thQUILIND}).

\begin{theorem}\label{thMaBrCo}  $\!$\emph{(Maroscia, Brewer{ \& }Costa)}\\
 If $\gA$ is an \anar with $\Kdim\gA\leq1$,
every \mptf over $\AXr$ is extended from~$\gA$.
\end{theorem}
\begin{proof}
Since $\Ared[\uX]=\AuX\red$ and $\GKO(\gB)=\GKO(\gB\red)$, it suffices to treat the reduced case, \cad the case of the \adpsz.\\
Let us verify that the class of \adps of \ddkz~\hbox{$\leq 1$}
satisfies the hypotheses of \thref{thQUILIND}.
The first condition is item~(ii) above that we have just proven. \\
The second condition is that the \mptfs over $\AX$ are extended from $\gA$. This is the Bass-Simis-Vasconcelos \thoz.
\end{proof}

\rdb
\subsubsection*{The Lequain-Simis induction} \label{LSinduction}

For the purpose of generalizing the Quillen-Suslin \tho to  Pr\"ufer domains, 
and observing that this class of \ris is not stable under the passage from $\gA$ to $\ArX$, 
Lequain and Simis \cite{LS} 
have found a skillful way to bypass the difficulty by proving a new induction \tho \gui{\`a la Quillen,} 
suitably modified.

\pagebreak

\CMnewtheorem{ILSa}{Abstract Lequain-Simis induction}{\itshape}
\begin{ILSa}\label{LSindabs} ~\\
Let $\cF$ be a class of \ris that satisfy the following \prtsz.
\begin{description}
\item [{\rm (LS1)}]
If $\gA \in \cF$, every nonmaximal \idep $\fp$ of $\gA$
has a finite height.\footnote{I.e., $\Kdim(\gA_\fp)<\infty$.}
\item [{\rm (LS2)}]
If $\gA \in \cF $, then $ \AX_{{\fp}[X]}  \in \cF $ for every \idep
${\fp}$  of $\gA$.
\item [{\rm (LS3)}]
If $\gA \in \cF $, then $\gA_{\fp}\in\cF$ for every \idep  ${\fp}$ of $\gA$.
\item [{\rm (LS4)}]
If $\gA \in \cF $ is local,  every \mptf over $\AX$ is~free.
\end{description}
Then, for all $\gA\in\cF$ and all $r\geq1$, 
every \mptf  over~$\AXr$ is extended from $\gA$. 
\end{ILSa}

\smallskip 
Here note that if $\gA$ is local with $\Rad\gA={\fm}$,
then $\gA(X)=\AX_{{\fm}[X]}$.

\smallskip 
We propose a \gui{\cov variation} on the theme of the Lequain-Simis induction. This is an important application of our dynamic comparison between $\gA(X)$ and $\ArX$. This \cov induction \gui{\`a la Lequain-Simis} is due to I.~Yengui.

\smallskip 
\begin{theorem} $\!$\emph{(Yengui induction)}\label{induYengui} 
\\ Let $\cF$ be a class of commutative \ris of finite \ddk (not \ncrt bounded) which satisfies the following \prtsz.
\begin{description}
\item [{\rm (ls1)}] If $\gA \in \cF$, then $\gA(X) \in \cF$.
\item [{\rm (ls2)}] If $\gA \in \cF$, then $\gA_S \in \cF$ for every \mo $S$ of $\gA$.
\item [{\rm
(ls3)}] If $\gA \in \cF$, then every \ptf $\AX$-module is extended from~$\gA$.
\end{description}
Then, for all $\gA\in\cF$ and all $r\geq1$, 
every \mptf  over~$\AXr$ is extended from $\gA$.
\end{theorem}
Note: (ls1) replaces (LS2), (ls2) replaces (LS3) and (ls3) replaces (LS4).
\begin{proof}
Due to Fact~\ref{fact A->A[X]}~\emph{\ref{i5fact A->A[X]}},
we limit ourselves to the case of reduced \risz.
We reason by double \recu on the number $r$ of variables and over the \ddk $d$ of~$\gA$.
\\
The basic step for $r=1$ (arbitrary $d$) is given by (ls3), and \hbox{for $d=0$} (with arbitrary $r$) it is the Quillen-Suslin \thoz.
\\
We suppose that the result is  proven in $r$ variables
 for the \ris in~$\cF$. We consider the case of $r+1$ variables and we perform an \recu \dem on (an upper bound $d$ of) the \ddk of a \riz~\hbox{$\gA\in\cF$}.
\\
Therefore let $\gA$ be a \ri  of \ddk $\leq d+1$.
 Let $P$ be a \mptf over $\gA[\Xr,Y]=\gA[\uX,Y]$.
Let~$G=G(\uX,Y)$ be a \mpn of $P$ with \coes  
in~\hbox{$\gA[\uX,Y]$}.
Let $H(\uX,Y)$ be the matrix constructed from $G$ as in Fact~\ref{factEtPol}.
\\
By using the \hdr for $r$ and (ls1), we obtain that the matrices $H(\uX,Y)$ and $H(\uze,Y)$ are elementarily \eqves over~$\gA(Y)[\uX]$.
This means that there exist matrices $Q_{1}$, $R_{1}$ over $\gA[\uX,Y]$ such that
$$\preskip.4em \postskip.4em
\begin{array}{c}
Q_{1}H(\uX,Y)  =  H(\uze,Y)R_{1} \label{eqeq}\\[1mm]
\hbox{with}\quad
\det(Q_{1})   \hbox{ and }   \det(R_{1}) \;\hbox{ primitive  in } \gA[Y].
\end{array}
$$
We now show that $H(\uX,Y)$ and $H(\uze,Y)$ are \eqves %
over~$\ArY[\uX]$.
By the Vaserstein \rcm it suffices to show that they are \eqves over $\ArY_{S_{i}}[\uX]$ for \mocoz~$S_{i}$ of $\ArY$.
\\
We consider the primitive \pol $f=\det(Q_{1})\det(R_{1})\in\gA[Y]$, and we apply \thref{compa}.
If $f$ is of formal degree $m$, 
we obtain \mos $(S_i)_{i \in \lrbm}$ of $\gA$ such that 
the \mos $V=f^{\NN}$ and 
$(S_i)_{i\in \lrbm}$ are \com  in $\ArY$. 
In addition, $\Kdim\gA_{S_i}\leq d$ for $i\in \lrbm$.
\\
For the \ri localized at $V$, 
$\det(Q_{1})$ and $\det(R_{1})$ are \ivs in $\ArY_V$. This implies that $H(\uX,Y)$ and $H(\uze,Y)$ are \eqves over $\ArY_V[\uX]$.
\\
For a localized \ri at $S_i$ ($i\in \lrbm$), 
by \hdr over $d$ and by using (ls2),  $H(\uX,Y)$ and $H(\uze,0)$ are \eqves over
$\gA_{S_i}[\uX,Y]$.
A~fortiori $H(\uX,Y)$ and $H(\uze,Y)$ are \eqves over $\gA_{S_i}[\uX,Y]$, therefore also 
over $\ArY_{S_i}[\uX]$, which is a \lon of $\gA_{S_i}[Y][\uX]=\gA_{S_i}[\uX,Y]$.

 Thus, we have fulfilled the contract and we obtain \iv matrices~$Q$ and $R$ over $\ArY[\uX] \subseteq \AXY$ such that 
$$
Q\,H(\uX,Y)=H(\uze,Y)\,R.
$$
  Moreover, we know by (ls3) that $H(\uze,0)$ and $H(\uze,Y)$ are \eqves over $\gA[Y]\subseteq \AXY$, and by \hdr over $r$ that $H(\uze,0)$ and $H(\uX,0)$ are \eqves over $\AuX\subseteq \AXY$.
In conclusion $H(\uX,0)$ and~$H(\uX,Y)$ are \eqves over $\AXY$.
Therefore by the Affine Horrocks' \thoz, $P$ is extended from $\AuX$.\\
 Finally, by \hdr over $r$, $P(\uX,0)$ is extended from~$\gA$.
\end{proof}

\rem We have asked in (ls2) that the class $\cF$ is stable under \lon for any \moz. Actually in the \dem only \lons at Krull boundary \mos intervene, or by  inversion of a unique \elt (all this in an iterative way).
\eoe

\subsubsection*{Lequain-Simis in finite dimension}

\begin{corollary}~\label{corLSValu}\\
If $\gA$ is an \anar of finite \ddkz, every \mptf over $\AXr$ is extended from $\gA$.
\end{corollary}
\begin{proof}
We show that the class of \anars of finite Krull dimension satisfies the concrete Lequain-Simis induction.
 The condition (ls1) is given by Exercise~\ref{exoPrufNagata}, (ls3) by the Bass-Simis-Vasconcelos \thoz, and~(ls2) is clear.
\end{proof}

\subsubsection*{Local Lequain-Simis without the dimension hypothesis}

\begin{corollary}\label{corLSValu2} If $\gV$ is a \advz, every \mptf over $\VXr$ is extended from $\gV$ (\cad free).
\end{corollary}
%
\begin{proof}
 Let  $M$ be a \mptf over $\VXr$.
We must show that $M$ is free.
Let $F=(f_{ij})\in\GAq(\VXr)$ be a matrix whose image is \isoc to the module $M$. Let $\gV_1$ be the sub\ri of $\gV$
generated by the \coes of the \pols $f_{ij}$ and $\gV'$ be the \sdv of $\gV$ generated by $\gV_1$. 
Item~\emph{\iref{i5corthValDim}} of \thref{corthValDim} tells us that every \ri between $\gV_1$ and $\Frac\gV_1$,
in particular $\gV'$, is of finite \ddkz. We apply Corollary~\ref{corLSValu}.
\end{proof}

\subsubsection*{\Gnl Lequain-Simis \thoz}

\begin{theorem}\label{thLSValu}  $\!$\emph{(Lequain-Simis)}
If $\gA$ is an \anarz, every \mptf over $\AXr$ is extended from $\gA$.
\end{theorem}
\begin{proof}
This results from Corollary~\ref{corLSValu} (the local case) with the same \dem as as the proof which deduces \thref{thBassAri} from \thref{thBassValu}.
\end{proof}

\entrenous{ tr\`es incertain \ldots

{\bf Am\'elioration of Brewer \& Costa}
\label{subsecBreCos}

Il serait question of d\'emontrer in le cas int\`egre que if la cl\^oture
int\'egrale of $\gA$ is a \ri seminormal, then les \mptfs over $\AXr$
sont tous \'etendus depuis $\gA$ \ssi $\gA$ is seminormal.

In le \tho of Brewer\&Costa proprement dit, $\gA$ doit en plus v\'erifier
une hypoth\`ese myst\'erieuse difficile \`a d\'{e}crypter.
}


\vspace{5pt}
\section*{Conclusion: a few conjectures}
\addcontentsline{toc}{section}{Conclusion: a few conjectures}
  \label{sec.BQHConj} 

The solution to Serre's \pb has naturally led to a few conjectures about possibles \gnnsz.

We will cite the two most famous ones and refer to \cite[chap.V,VIII]{Lam06} for detailed information on the subject.

\smallskip  The first, and the strongest, 
is the \emph{Hermite \ris conjecture}, 
that can be stated in two \eqves forms, 
one local and another global, 
given the Quillen patching principle.
Recall that a \ri is called a \gui{Hermite \riz} when the \stl \tf modules are free, which amounts to saying that the \vmds are completable. 

 {\bf (H)} If $\gA$ is a Hermite \riz, then so is $\AX$.

 {\bf (H')} If $\gA$ is a \dcd \aloz, then $\AX$ is a Hermite \riz.

Bass' \gui{stable-range} gives a first approach of the \pb (see Proposition~\ref{propStabliblib}, Corollary~\ref{corBass} and \thref{corBass2}). Special cases are treated for example in \cite[Roitman]{Roi86} and \cite[Yengui]{Ye3,Ye4}, which treats the $n=1$ case of the following conjecture: over a \riz~$\gA$ of Krull dimension $\leq 1$, the stably free $\AXn$-modules are free.

\smallskip  The second is the \emph{Bass-Quillen conjecture\index{Quillen}}.

A \cori is called a \emph{regular \riz} if every \mpf admits a finite \pro resolution (for the \dfn and an example of a finite \pro resolution, see \Pbmz~\ref{exoFossumKumarNori}). 
For the Bass-Quillen conjecture there are also two \eqv versions, a local one and  a global one.

 {\bf (BQ)} If $\gA$ is a regular \coh \noe \riz,\footnote{Naturally, in \clama the hypothesis \gui{\cohz} is superfluous.} then the \mptfs over $\AXn$ are extended from $\gA$.

 {\bf (BQ')} If $\gA$ is a regular \coh \noe \dcd \aloz,\footnote{Naturally, in \clama the hypotheses \gui{\cohz} and \gui{\dcdz} are superfluous.} then the \mptfs over $\AXn$ are free.

Actually, since $\gA$ regular \noe implies $\AX$ regular \noez, it would suffice to prove the $n=1$ case. 
Partial results have been obtained.
For example, the conjecture is proven in Krull dimension $\leq2$, for arbitrary~$n$ 
(but at the moment we do not dispose of a \cov \demz).
We can a priori \egmt consider a non-\noe version for the regular \coris of fixed \ddk $\leq k$.


\Exercices

\begin{exercise}
 \label{exolem1SusQS}
 {\rm
Let $\fA$ be an \id of $\AX$ containing a \polu and $\fa$ be an \id of $\gA$.
Then $\gA\,\cap\, (\fA+\fa[X])$ is contained in $\DA\big((\gA\,\cap\,\fA)+\fa\big)$.
In particular, if $1\in \fA+\fa[X]$, then $1\in(\gA\,\cap\,\fA)+\fa$.
} 
\end{exercise}


\vspace{-1em}
\begin{exercise}
 \label{exolemLocLocCom} (Top-Bottom lemma)
 {\rm Let $\gA$ be a \ri and  $\fm=\Rad\gA$.
 \begin{enumerate}\itemsep-1pt
\item Let  $S\subseteq\AX$ be the \mo of the \polusz. The \mos $S$ and $1+\fm[X]$ are \comz.
\item Let $U\subseteq \AX$ be the \mo $\sotq{X^n+\sum_{k<n}a_kX^k}{n\in\NN,a_k\in\fm\; (k<n)}$.
The \mos $U$ and $1+\fm+X\AX$ are \comz.
\end{enumerate}
} 
\end{exercise}

\vspace{-1em}
\begin{exercise}
\label{exoSLequiv}
{\rm  The goal of the exercise is to show a result similar to the Vaserstein \rcm (\plgrf{thPatchV}) in which we replace $\GLn$ by~$\SLn$.
 
 \vspace{-.3em}
 \begin{enumerate}\itemsep-2pt
\item Let $\gB$ be a \ri and $S$ be \mo of $\gB$.
\begin{enumerate}\itemsep-1pt
\item Let $P\in\BY$ such that $P(0)=0$ and $P=0$ in $\gB_S[Y]$. Show that there exists an $s\in S$ such that $P(sY)=0$.
\item Let $H\in\Mn(\BY)$ such that $H(0)\in\SLn(\gB)$ and $H\in\SLn(\gB_S[Y])$. Show that there exists an $s\in S$ such that $H(sY)\in\SLn(\BY)$.
\end{enumerate}
\item Prove Lemma~\ref{lem2PrepVaser} by replacing
$\GL$ by $\SL$.
\item Prove the \plgrf{thPatchV} by replacing $\GL$ by $\SL$.
\end{enumerate}
} 
\end{exercise}


\vspace{-1em}
\begin{exercise}
 \label{exoCasParticHorrocks}
{\rm Let $\gA$ be a \dcd \alo and $\fb \subseteq \gA[X]$ be an \iv \id containing a \poluz. 
We want to show that $\fb$ is a \idpz.}

This constitutes a special case of the Local Horrocks' \tho (\thref{thHor0}): 
indeed, on the one hand $\fb$ is a \pro $\gA[X]$-module, 
and on the other hand, 
if $f \in\fb$ is a \poluz, 
then by localizing at $f$, 
$\fb_f = \gA[X]_f$, and so, 
by the Local Horrocks' \thoz, $\,\fb$ is a free $\gA[X]$-module. 
This exercise gives a \dem independent from the current one.
In the special case studied here, we add the assumption that~$\fb$ is generated by a \poluz.

{\rm  Let $\gA$ be a \riz, let $\fm=\Rad\gA$ and $\gk = \gA/\fm$.
Let $\fb \subseteq \AX$ be an \id containing a \poluz.
Let $\ov a$ be the reduction of $a$ modulo $\fm$.
 
\emph {1.}
Prove that every \polu of $\ov\fb\subseteq\kX$ can be lifted to a \polu of~$\fb$.

 Now suppose that $\gA$ is \dcd and local.

\emph {2.}
 Show the existence of a \polu $f \in \fb$ such that $\ov \fb = \langle \ov f\rangle$ in $\gk[X]$
and so $\fb = \langle f\rangle + \fb \cap \fm[X]$.

Now suppose that the \id $\fb$ is \ivz.

\emph {3.}
Show that $\fb \cap \fm[X] = \fb\fm[X]$.

\emph {4.}
 Consider the \ri $\aqo{\AX}{f}$.
Show that $\fm (\aqo{\fb}{f}) = \aqo{\fb}{f}$.
\\ Deduce that $\fb = \gen{f}$.

We propose a \gnnz.

\emph {5.}
 Does the \dem work with a \plc \ri $\gA$?

} 
\end{exercise}


\vspace{-1em}
\begin{exercise}\label{exoBézoutKdim1TransfertArX}
{(Brewer {\&} Costa \thoz: the case of  Bézout domains of dimension $\le 1$)}
{\rm  See also Exercise~\ref{exoPrufNagata} and \thref{thMaBrCo}.
 
     Let  $\cF$  be the  class of  Bézout domains of dimension $\le
1$, and $\gA\in\cF$. 
 
\emph {1.}
Show that $\Kdim\ArX \le 1$ (use Exercise~\ref{exoMultiplicativiteIdeauxBords}).

\emph {2.}
Deduce that $\ArX$ is a Bézout \riz.
 
\emph {3.}
The class $\cF$  satisfies the hypotheses of \thref{th2QUILIND} (concrete Quillen induction, free case).
 Thus,  every \ptf $\AXr$-module is free.

}

\end{exercise}

\vspace{-1em}
\begin{exercise}
\label{exoseminorlgb} (\Lgb principle for  seminormal \risz)\\
{\rm  We give a direct \dem of  principle \ref{plgcetendus} in the special case of  seminormal \lsdsz.
\emph{1.} In a \lsdz, if $xc=b$ and $b^2=c^3$, then
there exists a $z$ such that $zc=b$ and $z^2=c$, so $z^3=b$.

\emph{2.} Let $S_1$, \ldots, $S_n$ be \moco of a \ri  $\gA$.
Suppose that each of the $\gA_{S_i}$'s is a seminormal \lsdz. 
Show that~$\gA$ is a seminormal \lsdz.
 
}
\end{exercise}

\vspace{-1.2em}
\pagebreak	
\begin{exercise}\label{exoKdimBounded}
{(Rings satisfying some of the conditions of the section \gui{An example in finite \ddk $>0$} \paref{subsec1VXfree})}\\
{\rm
Let $a_1$, \ldots, $a_k \in \gA$, $a_0 = 0$,
$a_{k+1} = 1$. Let $\gA_1$, \ldots, $\gA_{k+1}$ be  the following \ris 

\snic{
\gA_i = \big(\aqo{\gA}{a_{i-1}}\!\big)[1/a_i] \quad \hbox {for $i \in \lrb{1..k+1}$}.
}

\snii
We will show that if each $\gA_i$ is \zedz, then $\Kdim\gA \le k$.
The same result holds with $\gA_i = \big(\gA/\rD_\gA(a_{i-1})\!\big)[1/a_i]$.

\snii
\emph {1.}
Let $a \in \gA$. If $\Kdim \gA[1/a] \le n$ and $\Kdim \aqo{\gA}{a} \le m$, then $\Kdim\gA \le n+m+1$.

\snii
\emph {2.}
Deduce the stated result.

}

\end{exercise}



\sol


\exer{exolem1SusQS}
{Let $\gB=\gA\sur{\gA\,\cap\,\fA}$, $\gB'=\AX\sur{\fA}$, $\fb=\ov{\fa}$,
$\fb'=\fb\,\gB'$. \\
The \ri $\gB'$ is an integral extension of $\gB$. We apply the Lying Over (\ref{lemLingOver}).\\
\emph{Another solution}.
Let $f\in\fA$ be \monz. Let $a\in\gA\,\cap\, (\fA+\fa[X])$,
there exists a $g\in\fA$ such that $g\equiv a \mod \fa$. Then $\Res(f,g)\equiv\Res(f,a)\mod\fa$. But $\Res(f,a)=a^{\deg f}$ and $\Res(f,g)\in\fA\,\cap\,\gA$.
}


\exer{exolemLocLocCom}Use the resultant.


\exer{exoSLequiv}
\emph{1b.} Let $P(Y)=1-\det\big(H(Y)\big)$. We apply item \emph{1a}.

\emph{2.} Lemma~\ref{lem2PrepVaser} provides us with a matrix $U(X,Y)\in
\GL_r(\gA[X,Y])$ such that 

\snic{U(X,0)=\I_r$ and, over $\gA_S[X,Y]$,
$U(X,Y)=C(X+sY)C(X)^{-1}.}

By item \emph{1}, there exists a  $t\in S$
such that $U(X,tY)\in \SL_r(\gA[X,Y])$.\\
Let $V(X,Y)=U(X,tY)$
and we replace $s$ by $st$. 

\emph{3.} Lemma~\ref{lem3PrepVaser} is successfully subjected to the replacement of $\GL$ (implicit in the word \gui{\eqvez}) by $\SL$.
Likewise for the Vaserstein \rcmz.


\exer{exoCasParticHorrocks}
\emph{1.}
We first show the following result: if we have $g$, $f \in \fb$ with
$\ov g$ \mon of degree $r$ and $f$ \mon of degree $r+1$, then
$\ov g$ can be lifted to a \polu of $\fb$
(of degree $r$). 
We write $g = aX^{r + \delta} + \dots$, with $\delta \in \NN$ and we show by \recu on $\delta$ that $\ov g$ can be lifted to a \polu in $\fb$. 
If $\delta = 0$, we have $a \equiv 1 \bmod \fm$ (because $\ov g$ is \monz), so $a$ is \iv and the \polu $a^{-1}g \in \fb$ lifts $\ov g$.  If $\delta \ge 1$, we have $a \in \fm$ (because $\ov g$ is \monz), and we consider $h = g - aX^{\delta-1}f \in \fb$. It is of the form $bX^{r + \delta-1} + \dots$, and it satisfies~\hbox{$\ov h = \ov g$}. We apply the \hdrz.
\\ 
It then suffices to show that for all $g \in \fb$ such that $\ov g$ is \mon of degree~$r$, the \idz~$\fb$ contains a \polu of degree $r+1$. By hypothesis,~$\fb$ contains a \polu $f$. If $\deg(f) \le r+1$, then the result is clear. If~\hbox{$n = \deg(f) > r+1$}, then the \pol $X^{n-(r+1)}\ov g$ is \mon of degree $n-1$, and by the first step,
$\fb$ contains a \polu of degree $n-1$. We~conclude by \recu on $n-r$.

 \emph{2.}
The \id $\ov \fb$ is a \itf of $\gk[X]$,
therefore $\ov\fb$ is principal.  As $\fb$ contains a \polu we can take the \mon \gtr $\ov h$ and we lift it to a \polu of $\fb$ by the previous question.

 \emph{3.} Let $f$ be \mon in $\fb$, and $\fb_1$ be the \id that satisfies $\fb\fb_1=\gen{f}$.\\
We consider  $\fb'=\fb_1(\fb\cap \fm[X])/f$ (it is an \id of $\AX$).
Then $\fb\fb' = \fb\cap \fm[X]$.
We have $f\fb'\subseteq\fm[X]$ and $f$ is \mon so $\ov{\fb'}=0$, \cad $\fb' \subseteq \fm[X]$. By multiplying by~$\fb$, we obtain $\fb \cap \fm[X] \subseteq \fb\fm[X]$, so $\fb \cap \fm[X] = \fb\fm[X]$.

 \emph{4.} We have

\snic{
\fm (\aqo{\fb}{f}) = \aqo{\fc}{f}
\; \hbox {with} \;
\fc = \fm\fb + \gen{f} = \fm[X]\fb +  \gen{f} =
\fm[X] \cap \fb + \gen{f} = \fb.}

The $\aqo{\AX}{f}$-module $\aqo{\fb}{f}$ is \tf and as $f$ is \monz, $\aqo{\AX}{f}$ is a \tf \Amoz.
We deduce that $\aqo{\fb}{f}$ is a \tf \Amoz.
By Nakayama's lemma we obtain $\aqo{\fb}{f} = 0$, \cad $\fb = \gen{f}$.


\exer{exoBézoutKdim1TransfertArX} 
\emph {1.}
We must show that for $f$, $g \in \gA[X]$, we have $1 \in \IK_{\ArX}(f, g)$.  Since $\gA$ is a Bézout domain,  every \pol of $\gA[X]$ is the product of an \elt of $\gA$ by a primitive \polz. By Exercise \ref{exoMultiplicativiteIdeauxBords}, it suffices to show that $1 \in \IK_{\ArX}(f, g)$, either when $f$ or $g$ is primitive, or when $f$ and $g$ are constants $a$, $b$. In the latter case, since $\Kdim\gA \le 1$, this stems from $1\in \IK_{\gA}(a,b) \subseteq \IK_{\ArX}(a,b)$.
\\ 
Therefore suppose that $f$ or $g$ is primitive, for example $f$. It suffices to show  
that $1 \in \IK_{\ArX}(f, g)$ after \lon at \mocoz. However, \thref{compa} provides boundary \mos $S_j$ in $\gA$ such that $f^\NN$ and the $S_j$'s are \com in $\ArX$. For the \lon at $f^\NN$, it is clear that $1 \in \IK(f,g)$.  
\\
 As for $S_j^{-1}\ArX$, it is a \lon of $\gA_{S_j}[X]$ with $\gA_{S_j}$ \zedz, which gives  $\Kdim \gA_{S_j}[X] \le 1$. Therefore $1 \in \IK(f,g)$ in $\gA_{S_j}[X]$,  and a fortiori in the localized \ri $S_j^{-1}\ArX$.

\emph {2.}
The \ri $\gA[X]$ is a  GCD-domain, so the same holds for its localized \riz~$\ArX$. 
As $\Kdim \ArX \le 1$, \thref{propGCDDim1} tells us that $\ArX$ is a Bézout \riz.

\emph {3.} We have proven the \prt (q1) and we already know that the \prt (q0) is satisfied (\thref{thBézoutLib}).


\exer{exoseminorlgb}
\emph{1.}  We have $x^2c^2=b^{2}=c^3$, so $c^2(x^2-c)=0$,
therefore $c(x^2-c)=0$. 
\\
Then let $s$, $t$ such that $s+t=1$, $sc=0$ and $t(x^2-c)=0$. Let $z=tx$. \\
We have $tc=c$, $z^2=t^2c=c$ and $zc=xtc=xc=b$.

\emph{2.} 
Suppose that each of the $\gA_{S_i}$'s is a seminormal \lsdz. Therefore $\gA$ is a \lsdz. 
 \\
 Let $b$, $c\in\gA$  with $b^2=c^3$. If the $\gA_{S_i}$'s are seminormal, there exist $x_i\in\gA_{S_i}$ such that $x_i^2=c$ and $x_i^3=b$, and so $x_ic=b$. This implies that there exists an $x\in\gA$ such \hbox{that $xc=b$}.
We conclude by item \emph{1}.

Note: There are seminormal \ris that are not \lsdsz: for example $\gk[x,y]$ with $xy=0$
where $\gk$ is a \cdiz.


\exer{exoKdimBounded} 
\\
\emph {1.}
Let $(\ux) = (x_0, \ldots, x_n)$, $(\uy) = (y_0, \ldots,y_m)$ be $n+m+2$ \elts of $\gA$.  
By considering the iterated boundary \mo of $(\uy)$ in $\aqo{\gA}{a}$, we obtain that $\SK_\gA(\uy)$ contains a multiple of $a$, say $ba$. By considering the iterated boundary \id of $(\ux)$ in $\gA[1/a]$, we obtain that $\IK_\gA(\ux)$ contains a power of $a$, say $a^e$.
\\
 Then $(ba)^e \in \IK_\gA(\ux) \cap
\SK_\gA(\uy)$, so $1 \in \IK_\gA(\ux,\uy)$ by Fact~\ref{fact0BordKrullItere}, item~\emph{1}.

\snii
\emph {2.}
By using the previous question, we show by \recu on $i \in \lrb{0..k+1}$ that we have $\Kdim\gA[1/a_i] \le i-1$; for $i = k+1$, we obtain $\Kdim\gA \le k$.




\Biblio

Carlo Traverso proved the theorem that bears his name in
  \cite{Tra},  for a reduced \noe \ri $\gA$ (with an additional restriction). For the integral case without the \noee hypothesis we can refer back to \cite[Querr\'e]{Querre}, \cite[Brewer\&Costa]{BC2} and \cite[Gilmer\&Heitmann]{GH}. 
The most \gnl case is given by \cite[Swan]{Swan80}.

Traverso-Swan's \tho over the seminormal \ris has been decrypted from the \cov point of view by Coquand in \cite{coq}.
The decryption began with the \elr \dem of Proposition~\ref{propIntSemin} as it is given here. 
This \dem is a (quite spectacular) simplification of the existing \dems in the literature. It was then necessary to bypass the argument of the consideration of a \idemi to obtain a complete \prco of the result. 
It is remarkable that, at the same time, the case of a non-integral ring could have been treated effortlessly, in contrast to what happens in Swan's proof in~\cite{Swan80}. For a detailed explanation of \cite{coq} see \cite[Lombardi\&Quitt\'e]{LQ06}.
For a \gui{simple} \algo that realizes the \tho in the univariate case, see~\cite[Barhoumi\&Lombardi]{BL07}. 
A direct \demz, in the same spirit, for the implication \gui{seminormal ring $\gA$ implies seminormal ring $\AX$} is found in \cite[Barhoumi]{Barh09}.

Roitman's \tho \ref{thRoitman} is found in \cite{Roi}.

As for the history of the resolution of Serre's problem over  \pol \risz, the reader can refer to Chapter~III of  Lam's book~\cite{Lam06} as well as the presentation by Ferrand to Bourbaki~\cite{Ferrand}.

The original \dems of Quillen-Suslin's \tho (solution to Serre's \pbz) are found in\index{Quillen}\index{Suslin} \cite[Quillen]{Qu} and \cite[Suslin]{Sus}. Horrocks' \thos have their source in \cite[Horrocks]{Hor}.

The \gui{Quillen patching} that appears in \cite{Qu} is often called the Quillen \plgz. A remarkable overview of the applications of this principle and of its extensions is found in \cite[Basu\&al.]{BaRaKha}.
Also read~\cite[Rao\&Selby]{RaSe}.

The \ri $\ArX$ played a great role in the solution of Serre's \pb by Quillen and in its successive \gnns (the \thos of Maroscia and Brewer{\&}Costa, and of Lequain\&Simis).
The \ri $\gA(X)$ proved to be an efficient tool for several results of commutative \algz. Refer to the article \cite[Glaz]{Glaz1} for a considerably comprehensive bibliography regarding these two \risz.

Lam's book \cite{Lam06} (which follows \cite{Lam})  is a gold mine regarding the extended \pro modules. It contains especially several \dems of Horrocks' \thos (local and affine), with all the details and all the \ncr references, at least from a \clamaz' point of view.

The Affine Horrocks' \tho (\thref{thHor}) was \cot proven, first (for a slightly weaker variant) in the article \cite[Lombardi\&Quitt\'e]{LQ02}, then in \cite[Lombardi,Quitt\'e\&Yengui]{LQY05}.
The version given on \paref{thHor} reuses the latter article by specifying all the details. It is based on the books by Kunz and Lam.

\Thref{thBassValu} by Bass-Simis-Vasconcelos (\cite{BASS,SV}) was decrypted from a \cov point of view by Coquand in~\cite{coq07}.

Regarding the \tho of Maroscia and Brewer{\&}Costa (\thref{thMaBrCo}), see the original articles \cite{BC,Ma}. 
 A \cov \dem can be found in \cite{LQY05}.
This \tho is a slight antecessor of the \tho of Lequain\&Simis.
The latter was mostly decrypted from a \cov point of view by I.~Yengui  \cite{BLY06,ELY07}.

Many \algos for Quillen-Suslin's \tho (the field case) have been proposed in Computer Algebra, \gnlt based on the \dem by \Susz.
 
Quillen-Suslin's \tho has been studied from the point of view of its \algq complexity 
in \cite[Fitchas\&Galligo]{FiGa} and \cite[Caniglia\&al.]{CCDHKS} 
(for efficient \algosz, but that seem yet to be implemented).

A new, simple and efficient algorithm for SuslinÍs theorem (complete a \vmd containing a \poluz) is given in \cite[Lombardi\&Yengui]{YL04} and improved in \cite[Mnif\&Yengui]{MnY}.

\newpage \thispagestyle{CMcadreseul}
\incrementeexosetprob


\chapter[Suslin's stability \thoz]{Suslin's stability \thoz, the field case} 
\label{ChapSuslinStab}
\perso{compil\'e le \today}
\minitoc

\Intro


In this chapter, we give an entirely \cov treatment of Suslin's stability \tho for the case of \cdisz. 

\newpage	
\section{The \elr group} \label{secGpEn}

\vspace{4pt}
\subsec{Transvections}
\label{subsecTransvections} 

Regarding the \elr group $\En(\gA)$, recall that it is generated by the \elr matrices  $\rE^{(n)}_{i,j}(a)=\rE_{i,j}(a)$. 

If we let $(e_{ij})_{1 \leq i, j \leq n}$ be the canonical basis of $\Mn(\gA)$, 
we have 

\snic{
\rE_{i,j}(a) = \I_n + a e_{ij},
\quad 
\dsp\rE_{i,j}(a)\,e_k =
\cases {e_k               &if $k \ne j$    \cr  
        e_j + a e_i &if $k = j$\cr}
        \quad (i\neq j),}

with for example
$$\rE_{2,3}(a) = \cmatrix {
1 & 0 & 0 & 0\cr
0 & 1 & a & 0\cr
0 & 0 & 1 & 0\cr
0 & 0 & 0 & 1\cr
}.
$$

For fixed $i$ (resp.\ for fixed $j$) the matrices $\rE_{i,j}(\bullet)$ commute, and form a subgroup of $\EE_n(\gA)$ \isoc to $(\Ae {n-1},+)$. For example
$$\preskip.4em \postskip.0em
\rE_{2,1}(a)\cdot \rE_{2,3}(b)\cdot \rE_{2,4}(c) = \cmatrix {
1 & 0 & 0 & 0\cr
a & 1 & b & c\cr
0 & 0 & 1 & 0\cr
0 & 0 & 0 & 1\cr
}
$$
and
$$\preskip.2em \postskip.4em 
\;\;\Cmatrix{4pt} {
1 & 0 & 0 & 0\cr
a & 1 & b & c\cr
0 & 0 & 1 & 0\cr
0 & 0 & 0 & 1\cr
} \cdot\Cmatrix{4pt} {
1 & 0 & 0 & 0\cr
a' & 1 & b' & c'\cr
0 & 0 & 1 & 0\cr
0 & 0 & 0 & 1\cr
}=
 \Cmatrix{4pt} {
1 & 0 & 0 & 0\cr
a+a' & 1 & b+b' & c+c'\cr
0 & 0 & 1 & 0\cr
0 & 0 & 0 & 1\cr
}.
$$

\medskip 
More \gnlt let $P$ be a \ptf \Amoz. We will say that a pair $(\lambda,w)\in P\sta\times P$ is \emph{\umdz} if $\lambda(w)=1$. 
In this case $w$ is a \umd \elt of $P$, $\lambda$  is a \umd \elt of $P\sta$ and the \Ali $\theta_P(\lambda\otimes w):P\to P$ defined by $x\mapsto\lambda(x)w$ is the projection over $L=\gA w$ \paralm
to $K=\Ker \lambda$, represented over $K\times L$ by the matrix%
\index{unimodular!pair}%
\index{pair!unimodular ---}

\snic{\bloc{0_{K\to K}}{0_{L\to K}}{0_{K\to L}}{1_{L\to L}}=
\bloc{0_{K\to K}}{0}{0}{\Id_L}.
}

\smallskip  If $u\in K$, the \Ali $\tau_{\lambda,u}:=\Id_P+\theta_P(\lambda\otimes u)$,
$x\mapsto x+\lambda(x)u$ is called a \ix{transvection}, 
it is represented over $K\times L$ by the matrix

\snic{ \bloc{1_{K\to K}}{(\lambda\otimes u)\frt{L}}{0_{K\to L}}{1_{L\to L}}
= \bloc{\Id_{K}}{(\lambda\otimes u)\frt{L}}{0}{\Id_L}.
}

For example, if $P=\Ae n$, an \elr matrix defines a transvection.

\rdb\smallskip\label{NOTAtransvec}  
Let $\GL(P)$ be the group of \lin \autos of $P$ and $\SL(P)$ be the subgroup of \endos of \deter $1$.
The subgroup of $\SL(P)$ generated by the transvections will be denoted by $\wi{\EE}(P)$.
The affine map
$$\preskip.4em \postskip.4em 
u\mapsto \tau_{\lambda,u},\;\Ker\lambda\to \End_\gA(P) 
$$
provides a \homo of the group $(\Ker\lambda,+)$ in the group $\wi{\EE}(P)$. 

In the case where $P=\Ae {n}$, if $\lambda$ is a \coo form, we find that the matrix of the transvection is a product of \elr matrices. For example,
with the vector $u=\tra{\lst{u_1\,u_2\,u_3\,0}}$:
$$\preskip.4em \postskip.4em 
\bloc{\I_{3}}{u'}{0}{1}\;=\;
\cmatrix{1&0&0&u_1\cr
0&1& 0&u_2\cr
0& 0 &1&u_{3}\cr
0&0 & 0 &1}\; =\;\dsp\prod_{i=1}^{3}\;\rE_{i,4}(u_i). 
$$

  However, note that a priori $\En(\gA)$ is only \emph{contained in} $\wi{\EE}(\Ae n)$. This shows that the \elr group is a priori deprived of clear \gmq meaning. As a crucial point, $\En(\gA)$ is a priori not stable under $\GLn(\gA)$-conjugation.

\subsec{Special matrices}
\label{subsecMatSpec} 

We now only speak of the groups $\En(\gA)$.

Let $u = \cmatrix {u_1\cr \vdots\cr u_n} \in \Ae {n \times 1}$ and $v = \cmatrix {v_1 & \cdots & v_n} \in \Ae {1 \times n}$ to which we associate the matrix $\I_n + uv \in \Mn(\gA)$. We will provide results specifying the membership of this matrix to the \elr group $\En(\gA)$.  
\\
Since {\mathrigid 2mu $\det(\I_n + uv) = 1 + \tr(uv) = 1 + vu$}, it is imperative to demand the \egt {\mathrigid 2mu $vu \eqdefi v_1u_1 + \cdots + v_nu_n = 0$}. 
In this case, we have {\mathrigid 2mu $(\I_n + uv)(\I_n - uv) = \I_n$}.

The transvections admit for matrices the matrices of this type, with $v$ being \umdz.
In addition, the set of these matrices $\I_n + uv$ (with $vu = 0$) is a stable set under $\GLn(\gA)$-conjugation. 
\\
For example, for $A \in \GLn(\gA)$, we obtain $A\, \rE_{ij}(a)\, A^{-1} = \I_n + a uv$, where $u$ is the column $i$ of $A$ and $v$ is the row $j$ of~$A^{-1}$.

Take care, however, that if we do not assume that $v$ is \umd these matrices do not in \gnl represent transvections. If neither $u$ nor $v$ is \umd the matrix does not even a priori represent an \elt of~$\wi\EE(\Ae n)$.

\begin{lemma}\label{lemE01Rao}
Suppose $u \in \Ae {n \times 1}$, $v \in \Ae {1 \times n}$ and $vu
= 0$. \\
Then 
$\cmatrix {\I_n + uv & 0\cr 0 & 1\cr} \in \EE_{n+1}(\gA)$.
\end{lemma}

\begin {proof}
We have a sequence of \elr operations on the right-hand side (the first uses the \egt $vu=0$)
$$\preskip.4em \postskip.4em 
\arraycolsep2pt
\begin{array}{cccccccccccccc}
\cmatrix {\I_n + uv & 0\cr 0 & 1\cr}  & \vers{\alpha}  &   
\cmatrix {\I_n + uv & -u\cr 0 & 1\cr} & \vers{\beta}\\[3mm]
   \cmatrix {\I_n  & -u\cr v & 1\cr}  & \vers{\gamma} & 
\cmatrix {\I_n & 0\cr v & 1\cr}  & \vers{\delta}  &  \cmatrix {\I_n & 0\cr 0 & 1}. 
\end{array} 
$$
This implies
 $\cmatrix {\I_n + uv & 0\cr 0 & 1\cr} = 
   \delta^{-1} \gamma^{-1} \beta^{-1} \alpha^{-1}$,
\cad

\snic{\cmatrix {\I_n + uv & 0\cr 0 & 1\cr} =\cmatrix {\I_n & 0\cr v & 1\cr} \cdot
\cmatrix {\I_n & -u\cr 0 & 1\cr} \cdot
\cmatrix {\I_n & 0 \cr -v & 1\cr} \cdot
\cmatrix {\I_n & u\cr 0 & 1\cr}. 
}

\vspace{-1em}
\end {proof}

A column vector $u$ is said to be \emph{special} if at least one of its \coos is null.
If $vu=0$ and if  $u$ is special we say that  $\In+uv$ is a \emph{special matrix}.\index{matrix!special ---}

\begin{corollary}\label{corlemE01Rao}
Let $u \in \Ae {n \times 1}$ and $v \in \Ae {1 \times n}$ satisfy $vu = 0$. If $u$ is special, then $\I_n + uv \in \En(\gA)$.
In other words every special matrix is in~$\En(\gA)$.
\end{corollary}

\begin{proof}
We can assume that $n\geq2$ and $u_n = 0$. 
Let \smash{$u = \cmatrix {\mathring {u} \cr 0\cr}$, $v = \cmatrix {\mathring {v} & v_n}$}, with $\mathring {u} \in
\Ae {(n-1) \times 1}$ and $\mathring v \in \Ae {1 \times (n-1)}$. Then
$\phantom{A^{A^A}}$
$$
\I_n + uv = 
\cmatrix {\I_{n-1} + \mathring{u}\,\mathring{v} & v_n\mathring{u}\cr 0 & 1\cr} =
\cmatrix {\I_{n-1} & v_n\mathring {u}\cr 0 & 1\cr} 
\cmatrix {\I_{n-1} + \mathring {u}\,\mathring{v} & 0\cr 0 & 1\cr}. 
$$
Since $\mathring {v}\,\mathring{u} = vu = 0$, Lemma~\ref{lemE01Rao} applies and $\I_n + uv \in \En(\gA)$.
\end{proof}

The special matrices are easily \gui{lifted} from a localized \ri $\gA_S$ to $\gA$ itself. 
More \prmtz, we obtain the following.

\begin{fact}\label{factE02Rao}
Let $S \subseteq \gA$ be a \moz, $u \in \gA_S^{n \times 1}$, $v \in \gA_S^{1 \times n}$ with $vu = 0$ and $u$ be special.
Then there exist $s \in S$, $\wi{u} \in \Ae {n \times 1}$, $\wi{v} \in \Ae {1
\times n}$ with $\wi{v}\wi{u} = 0$, $\wi{u}$ special and~\hbox{$u = \wi{u}/s$}, $v = \wi{v}/s$ over $\gA_S$.
\end{fact}

\begin{proof}
By \dfnz, $u = u'/s_1$, $v = v'/s_1$ with $s_1 \in S$, $u' \in \Ae {n \times 1}$ and $v' \in \Ae {1 \times n}$. The \egt $vu = 0$ provides some $s_2 \in S$ such that $s_2v'u' = 0$, and $u_i = 0$ provides some $s_3 \in S$ such that $s_3u'_i = 0$. Then $s = s_1s_2s_3$, $\wi{u} = s_2s_3u'$ and $\wi{v} = s_2s_3v'$ fulfill the required conditions.
\end{proof}

\begin{theorem}\label{lemE03Rao}
If $n \ge 3$, then $\wi\EE(\Ae n)=\En(\gA)$.
In particular, $\En(\gA)$ is stable under $\GLn(\gA)$-conjugation.
\\
 {\emph{Precisions:} Let $u \in \Ae {n \times 1}$, $v \in \Ae {1 \times n}$ with $vu
= 0$ and $v$  \umdz. Then, we can write $u$ in the form $u = u'_1 + u'_2 + \cdots + u'_N, $  
with $vu'_k = 0$ and each $u'_k$ has at most two nonzero components.
The matrix $\I_n + uv$ is then expressible as a product of special matrices
$$\preskip.4em \postskip.4em
\I_n + uv =  (\I_n + u'_1v)\  (\I_n + u'_2v)\ \cdots\ (\I_n + u'_Nv) 
$$
and consequently, it belongs to $\En(\gA)$.
}
\end{theorem}

\vspace{-.8em}
\pagebreak	

\begin{proof} The canonical basis of $\Ae n$ is denoted $(e_1, \ldots, e_n)$.
We have $a_1$, \ldots, $a_n$ in~$\gA$ such that $a_1v_1 + \cdots + a_nv_n = 1$.
\\
For $i \leq j$, let us define $a_{ij} \in \gA$ by $a_{ij} = u_i a_j - u_j a_i$.
Then

\snic{u = \som_{i < j} a_{ij} (v_j e_i - v_i e_j) =
\som_{i \leq j} a_{ij} (v_j e_i - v_i e_j).}

Indeed, for fixed $k$, the coefficient of $e_k$ in the right-hand sum is
$$\preskip.4em \postskip.4em 
{\arraycolsep2pt\begin{array}{lllllll}
\dsp 
\sum_{j \ge k} a_{kj} v_j - \sum_{i < k} a_{ik} v_i  &  = &  
\dsp
\sum_{j \ge k}\, (u_ka_j - u_ja_k) v_j - \sum_{i < k}\, (u_ia_k - u_ka_i) v_i
 \\[5mm]
  & =  & 
  \dsp 
  u_k \som_{j = 1}^n a_jv_j - a_k \som_{j = 1}^n u_jv_j\;=\;
   u_k. 
  \end{array}} 
$$
For $i < j$, we then define $u'_{ij} \in \Ae {n \times 1}$ by $u'_{ij} = a_{ij}(v_j e_i - v_i e_j)$. It is clear that $u'_{ij}$ has at most two nonzero components and that $vu'_{ij} = 0$. 
\end{proof}


\section{The Mennicke symbol} \label{secMennicke} 

\begin{lemma} \label{lemMennicke1}
Let $a, b$ be \com \elts in $\gA$. Then the \eqvc class in $\SL_3(\gA)/\EE_3(\gA)$ of the matrix $\cmatrix {a & b & 0\cr c & d & 0\cr 0 & 0 & 1\cr}$ does not depend on the choice of $c$ and $d$ satisfying $1 = ad - bc$.\\
We denote by $\meck {a}{b}$ the \elt of $\;\SL_3(\gA)/\EE_3(\gA)$ obtained thus. We call it the \ix{Mennicke symbol} of~$(a, b)$.
\end{lemma}

\begin{proof}
Let $A = \cmatrix {a&b\cr c&d\cr}$, 
$A' = \cmatrix {a&b\cr c'&d'\cr}$ with $ad - bc = ad' - bc' = 1$. Then

\snic{
AA'^{-1} = \cmatrix {a&b\cr c&d\cr}\cmatrix {d'&-b\cr -c'&a\cr} =
\cmatrix {1 & 0\cr cd' - c'd& 1\cr} ,  
}

and $\bloc{A}{0_{2,1}}{0_{1,2}}{1}{\bloc{A'}{0_{2,1}}{0_{1,2}}{1}}^{-1}=\bloc{AA'^{-1}}{0_{2,1}}{0_{1,2}}{1}$ is in $ \EE_3(\gA)$.
\end{proof}

\begin{proposition} \label{propMennicke1}
The Mennicke symbol satisfies the following \prtsz.
\begin{enumerate}
\item
If $a \in \Ati$, then $\meck {a}{b} = 1$ for all~$b \in \gA$.
\item
If $\gen{1} = \gen {a,b}= \gen {a',b}$ then
$1 \in \gen {aa',b}$ and $\meck {aa'}{b} = \meck {a}{b}\meck{a'}{b}$.
\item
If $1 \in \gen {a,b}$, then $\meck {a}{b} = \meck {b}{a} = \meck{a+tb}{b}$ for all~$t \in \gA$.

\end{enumerate}
\end{proposition}

\begin{proof}
\emph{1.} 
The matrix $\cmatrix {a & b\cr 0 & a^{-1}\cr}$ is a member of $\EE_2(\gA)$.

 \emph{2.} 
We have
$$\preskip.0em \postskip.0em
\cmatrix {a & b & 0\cr c & d & 0\cr 0 &0 & 1\cr} \sims{\EE_3(\gA)} 
\cmatrix {a & 0 & b\cr 0 & 1 & 0\cr c &0 & d\cr},
$$
and
$$\preskip.0em \postskip.4em
\cmatrix {a' & b & 0\cr c' & d' & 0\cr 0 &0 & 1\cr} \sims{\EE_3(\gA)} 
\cmatrix {a' & 0 & -b\cr c' & 0 & -d'\cr 0 &1 & 0\cr} \sims{\EE_3(\gA)} 
\cmatrix {a' & 0 & -b\cr c' & 0 & -d'\cr 0 &1 & a\cr}.
$$
The product $\meck {a}{b}\meck{a'}{b}$ is represented by the product of the matrices on the right-hand side, \cad by
$$
\cmatrix {aa' & b & 0\cr c' & 0 & -d'\cr ca' &d & 1\cr} \sims{\EE_3(\gA)} 
\cmatrix {aa' & b & 0\cr * & * & 0\cr ca' &d & 1\cr} \sims{\EE_3(\gA)} 
\cmatrix {aa' & b & 0\cr * & * & 0\cr 0 &0 & 1\cr},
$$
and therefore $\meck {a}{b}\meck{a'}{b} = \meck {aa'}{b}$.

 \emph{3.} 
If $ad - bc = 1$, then $\cmatrix {a &b\cr c&d\cr} \sims{\EE_2(\gA)} 
\cmatrix {-b &a\cr -d&c\cr},$ and so
$$
\meck {a}{b} = \meck{-b}{a} = \meck{-1}{a}\meck{b}{a} = \meck{b}{a}.
$$
Finally, $\cmatrix {a &b\cr c&d\cr} \sims{\EE_2(\gA)} 
\cmatrix {a + tb & b\cr c+td & d\cr},$ so $\meck {a}{b} = \meck{a+tb}{b}$.
\end{proof}

\begin{lemma}\label{lemII3.6}\relax \emph{(Local version)}\\ Let $\gA$ be a \dcd \alo and $f$, $g\in \AX $ be \com with $f$ \monz. Then we have 
$$\preskip-.2em \postskip.2em
\;\;\meck{f}{g}=\meck{f(0)}{g(0)}=1.$$
\end{lemma}

\begin{proof}
Let $af+bg=1$. First note that we can divide $b$ by~$f$ and that we then obtain an \egt $a_1f+b_1g=1$ with $\deg(b_1)<\deg(f)$, and so, since $f$ is \monz, $\deg(a_1)<\deg(g)$. 
Therefore assume \spdg that $\deg(b)<\deg(f)$ and $\deg(a)<\deg(g)$.\\
Let $r$ be the remainder of the Euclidean division of $g$ by $f$. 
Then $\meck{f}{g}=\meck{f}{r}$. 
In particular, if $\deg(f)=0$ we are done. Otherwise, we can assume $\deg(g)< \deg(f)$ and we reason by \recu on $\deg(f)$. 
Since~$\gA$ is local and \dcdz,  $g(0)\in\Ati$  or $g(0)\in\fm=\Rad\gA$.\\
First of all suppose that $g(0)$ is \ivz. Then
$$
 \meck{f}{g}=\meck{f-g(0)^{-1}f(0)g}{g},
$$
 so that we can assume that $f(0)=0$ and $f=Xf_1$.
Then 
$$
 \meck{Xf_1}{g}=\meck{X}{g}\meck{f_1}{g}=\meck{X}{g(0)}\meck{f_1}{g}
=\meck{f_1}{g}
$$
and the \dem ends by \recu since $f_1$ is \monz.
\\
Now suppose that $g(0)$ is in  $\fm$.  As $a(0)f(0)+b(0)g(0)=1$, we have  $a(0)f(0)\in 1+\fm\subseteq \Ati$, and 
so $a(0)\in\Ati$. However,
$$\preskip.4em \postskip.0em \left[\matrix{ 
  \phantom-f  &  g   &   0   \cr 
  -b  &  a   &   0   \cr
  \phantom-0  &  0   &   1
}\right]\equiv 
\left[\matrix{ 
  f-b  &  g+a   &   0   \cr 
  -b  &  a   &   0   \cr
  0  &  0   &   1
}\right]\;\;\mod\; \EE_3(\AX ),
$$
so
$$\preskip.0em \postskip.4em 
\meck{f}{g}=\meck{f-b}{g+a},
$$
with $f-b$ monic, 
 $\deg(f-b)=\deg (f)$, $\deg(g+a)<\deg(f)$ and 
$(g+a)(0)$ in~\hbox{$\fm+\Ati=\Ati$}. We are therefore brought back to the previous case.
\end{proof}
Our basic \lgbe machinery (\paref{MethodeIdeps}) applied to the previous local \demz, gives the following quasi-global lemma.\imlb 

\begin{lemma} 
\label{lemII3.6bis}\relax \emph{(Quasi-global version)}\\
Let $\gA$ be a \ri  and $f$, $g$ be \elts of $\AX $ \com with $f$ \monz. 
Then, there exists in $\gA$ a \sys of \eco  $(s_i)$  such that  in each localized \ri $\gA[1/s_i]$, 
we have the following \egt of the Mennicke symbols
 $$\preskip.2em \postskip-.2em
 \meck{f}{g}=\meck{f(0)}{g(0)}=1.
 $$
\end{lemma}

  
\section{Unimodular \pol vectors} 
\label{secCompVmdsPols}

\begin{notation}\label{notaGLIdeal} 
{\rm If $\fb$ is an \id of $\gB$, let $\GLn(\gB,\fb)$ be the subgroup of $\GLn(\gB)$ that is the kernel of the natural morphism $\GLn(\gB)\to\GLn(\gB\sur{\fb})$.
 We adopt an analogous notation for $\SLn$. \\
 Be careful for the group $\En$!
 Let~\hbox{$\En(\gB, \fb)$} be the normal subgroup generated by the $\rE_{ij}(b)$'s with
 $b \in \fb$. 
 }
\end{notation}

The group~\hbox{$\En(\gB, \fb)$} is a subgroup of the kernel of $\En(\gB)\to\En(\gB\sur{\fb})$, and in \gnlz, it is a strict subgroup.  However, in the case where $\gB = \gA[X]$ \hbox{and $\fb = \gen {X}$}, the two groups coincide. This result is given by the following lemma.

\begin{lemma}\label{lemE04Rao}
The group $\En(\AX, \gen{X})$ is the kernel of the canonical \homo $\En(\gA[X])\to \En(\gA[X]\sur{\gen {X}}\!) =  \En(\gA)$.
It is generated by the matrices of the type $\gamma\, \rE_{ij}(Xg)\, \gamma^{-1}$ with $\gamma \in \En(\gA)$ and $g \in \gA[X]$.
\end{lemma}

\begin{proof}
Let $H$ be this kernel. We will use the following \dcnz, valid in every group, of a product $\alpha_1\beta_1 \alpha_2\beta_2  \cdots \alpha_m\beta_m$, for example with $m=3$
$$\preskip.3em \postskip.4em
\bigl( \alpha_1\beta_1\alpha_1^{-1} \bigr)\,
\bigl( (\alpha_1\alpha_2)\beta_2(\alpha_1\alpha_2)^{-1} \bigr)\,
\bigl( (\alpha_1\alpha_2\alpha_3)\beta_3(\alpha_1\alpha_2\alpha_3)^{-1} \bigr)
\, (\alpha_1\alpha_2\alpha_3).
$$
Therefore let $E = E(X) \in H$, $E = \prod_{i=1}^m \rE_{i_k,j_k}(f_k)$ with $f_k \in \AX$. \\
We write $f_k = c_k + Xg_k$ with $c_k = f_k(0) \in \gA$ and
$$\preskip.3em \postskip.4em
\rE_{i_k,j_k}(f_k) = \alpha_k \beta_k, \quad \hbox {with} \quad
\alpha_k = \rE_{i_k,j_k}(c_k), \quad
\beta_k = \rE_{i_k,j_k}(Xg_k). 
$$
We finish by applying the \dcn given above and by using the \egtz~\hbox{$\alpha_1\cdots\alpha_m = E(0) = \I_n$}.
\end{proof}

\vspace{-.9em}
\pagebreak	

\begin{proposition}\label{propE0Rao}
Let $n \ge 3$, $s \in \gA$ and $E = E(X) \in \En(\gA_s[X],
\gen{X})$. There exist $k \in \NN$ and $E' = E'(X) \in \En(\AX,
\gen{X})$ satisfying $E'(X) = E(s^k X)$ over~$\gA_s[X]$.
\end{proposition}
\begin{proof}
We can suppose $E = \gamma \rE_{ij}(Xg) \gamma^{-1}$ with $\gamma \in \En(\gA_s)$ and $g \in \gA_s[X]$. Letting $u \in \gA_s^{n \times 1}$ be the column $i$ of $\gamma$ and $v \in \gA_s^{1 \times n}$ be the row $j$ of $\gamma^{-1}$, we have
$$\preskip.2em \postskip.4em
E(X) = \gamma \rE_{ij}(Xg) \gamma^{-1} = \I_n + (Xg) uv,
\quad vu = 0, \quad \hbox{$v$ unimodular}.
$$
\Thref{lemE03Rao} allows us to write 
$
u = u'_1 + u'_2 + \cdots + u'_N $ 
 {with $vu'_k = 0$ and~\hbox{$u'_k \in \gA_s^{n \times 1}$} 
has at most two nonzero components}.
We therefore have
$$\preskip.4em \postskip.4em
E(X) = (\I_n + (Xg) u'_1v)\  (\I_n + (Xg) u'_2v)\ \cdots\ (\I_n + (Xg) u'_Nv). 
$$
By using an analogous method to Fact~\ref{factE02Rao}, 
we easily prove that there exist $k \in \NN$, $\wi {g} \in \AX$, $\wi u_k \in \Ae {n \times 1}$ and $\wi {v} \in \Ae {1 \times n}$ such that we have over $\gA_s$ the \egts
$
g = \wi{g}/s^k $, $   u'_k = \wi u_k /s^k$,  $ 
v = \wi{v}/s^k$,  $\wi{v} \wi u_k = 0$,  
and $\wi u_k$ has at most two nonzero components.
Then let
$$\preskip.4em \postskip.4em
E'(X) = (\I_n + (X\wi{g}) \wi u_1\wi{v})\  
(\I_n + (X\wi{g}) \wi u_2 \wi{v})\ \cdots\ (\I_n + (X\wi{g}) \wi u_N \wi{v}) .
$$
By Corollary~\ref{corlemE01Rao}, each $\I_n + (X\wi{g}) \wi u_k \wi{v}$ belongs to $\En(\AX)$.  We therefore have $E'(X) \in \En(\AX)$, $E'(0) = \I_n$ and $E'(s^{3k}X) = E(X)$ over $\gA_s[X]$.
\end{proof}

%
\begin{lemma}\label{lemE1Rao}
Let {\mathrigid 2mu $n\geq3$ be an integer, $s\in\gA$ and $E=E(X)\in\En(\gA_s[X])$}. 
There exists an integer $k\geq0$ such that for all $a , \, b\in\gA$ congruent modulo $s^k$, the matrix $E^{-1}(aX)E(bX)$ is in  the image of the natural \homo  

\snic{\En(\AX,\gen{X})\longrightarrow\En(\gA_s[X],\gen{X}).}
\end{lemma}
Note: in short, but less precisely, if $a$ and $b$ are sufficiently \gui{close}, the \coes of the matrix $E^{-1}(aX)\,E(bX)$ have no more \denoz.
\begin{proof}
We introduce two new \idtrs $T$, $U$ and let
$$\preskip.4em \postskip.4em
E'(X,T,U)=E^{-1}\big((T+U)X\big)\,E(TX).
$$ 
We have $E'(X,T,0)=\In$. 
We apply Proposition~\ref{propE0Rao} with $F=E'$ by taking $\gA[X,T]$ instead of $\gA$ and $U$ instead of $X$: there exist a matrix~$G$ in $\En(\gA[X,T,U],\gen{U})$ and an integer $k\geq0$ such that 
$$\preskip.3em \postskip.4em 
E'(X,T,s^kU)=G(X,T,U) \;\hbox{ in }\;  \En(\gA_s[X,T,U],\gen{U}). 
$$
Therefore $G(X,T,U)=E^{-1}\big((T+s^kU)X\big) \, E(TX)$ over $\gA_s$, and if $b=a+s^kc$, then
$$\preskip-.2em \postskip.4em 
E^{-1}(aX)\,E(bX)=G(X,a,c) \; \hbox{ over } \; \gA_s . 
$$
We have $G(0,T,U)=\In$ over $\gA_s$, but not \ncrt over~$\gA$.
Let
$$\preskip.4em \postskip.4em 
H(X,T,U)=G^{-1}(0,T,U)\,G(X,T,U). 
$$
We then have  $H(0,T,U)=\In$ over $\gA$ and $H(X,T,U)=G(X,T,U)$ over~$\gA_s$.
We therefore obtain
$$\preskip.3em \postskip.4em 
E^{-1}(aX)\,E(bX)=H(X,a,c) \;\;\; \mathrm{in} \;\;
\; \En(\gA_s[X],\gen{X}), 
$$
with $H(X,a,c)\in \En(\AX,\gen{X})$.
\end{proof}
%

\begin{lemma}\label{lemE2Rao}
Let $n\geq3$ be an integer, $s\in\gA$ and 
$$\preskip.4em \postskip.3em
E=E(X)\in\GLn(\gA[X])\cap \En(\gA_s[X]).
$$ 
There exists an integer $k\geq0$ such that for all $a ,  b\in\gA$ congruent modulo~$s^k$, the matrix $E^{-1}(aX)E(bX)$ is in $\En(\gA[X],\lra X)$. 
\end{lemma}
\facile
%

\begin{lemma}\label{lemE3Rao}
Let $n\geq3$, $s$, $t$ be \com in $\gA$ and 
$$\preskip.4em \postskip.4em
E\in\GLn(\gA[X],\lra X)\cap \En(\gA_s[X]) \cap \En(\gA_t[X])
.$$
Then~$E\in \En(\gA[X])$.
\end{lemma}
\begin{proof}
By Lemma~\ref{lemE2Rao}, there exists some $k$ such that for all $a , \, b\in\gA$ congruent modulo $s^k$, or modulo $t^k$, the matrix $E^{-1}(aX)E(bX)$ is in  $\En(\gA[X],\lra X)$. Let $c\in\gA$ such that $c\equiv 0 \mod s^k$ and $c\equiv 1 \mod t^k$. \\
Then we write $E=E^{-1}(0\cdot X)\,E(c\cdot X)\,E^{-1}(c\cdot X)\,E(1\cdot X)$. 
\end{proof}
%

%


\section{Suslin's and Rao's \lgbs principles} 
\label{secPlgbRao}

Now we prove \cite[Lemma I 5.9 (page 26)]{GM}.
\begin{theorem}\label{lem159}
Let $n\geq 3$ and $A = A(X)\in \GL_n(\AX)$.

\smallskip 
\begin{enumerate}
\item 
If $A(0) = \In$, then
$
\fa=\sotq{s \in \gA}{  A\in \En(\gA_s[X]) }
$ 
is an \id of $\gA$. 

\item 
The set
$
\fa = \sotQ {s\in\gA}{A(X)\sims{\En(\gA_s[X])} A(0)}
$
is an \id of $\gA$.
\end{enumerate}
\end{theorem}

\begin{proof}
The two formulations are \eqvesz; we prove the second from the first by considering $A(X)A(0)^{-1}$.
\\ 
\emph{1.}  It is clear that $s \in \fa \Rightarrow as \in \fa$ for all $a \in \gA$.  Now let $s$, $t$ in $\fa$. We must show that $s+t\in\fa$, or that  $1\in\fa\gA_{s+t}$. In short, we suppose that $s$ and $t =1-s$ are in $\fa$, and we must show that $1\in\fa$.
\\ 
By \dfnz, we have $A \in\En(\gA_{s}[X])$ and $A \in\En(\gA_{t}[X])$; by Lemma~\ref{lemE3Rao}, we have $A \in\En(\gA[X])$, \cad $1 \in \fa$.
\end{proof}

This lemma could have been written in the form of the following \plgc  (very nearly \cite[Lemma I 5.8]{GM}).

\begin{plcc} 
\label{plcc.159} \emph{(For the elementary group)}\\
Let {$n\geq 3,$  $S_1$, $\ldots$, $S_k$  be comaximal monoids of $\gA$ and $A\in \GL_n(\AX)$, with $A(0)=\In$}. Then  
$$\preskip.3em \postskip.2em
A\in \En(\AX )\;\;\;\;\Longleftrightarrow\;\;\;\; 
\hbox{ for }\;i\in\lrbk,\quad A\in \En(\gA_{S_i}[X]).
$$
\end{plcc}

The following \tho re-expresses \cite[corollary II 3.8]{GM}.

\pagebreak	

\begin{theorem}\label{II3.8}\relax 
\emph{(Global version of Lemma~\ref{lemII3.6})}
\\ 
Let $n\geq 3,$ and $f$, $g\in \AX $ be \comz, with $f$ \monz. 
Then, we have the following \egt of  Mennicke symbols: $\meck{f}{g}=\meck{f(0)} {g(0)}.$ 
\end{theorem}
\begin{proof}
Let $af-bg=1$. Let $B=\cmatrix {f & g & 0\cr b & a & 0\cr 0 & 0 & 1\cr}$. \\
The \egt $\meck{f}{g}=\meck{f(0)}{g(0)}$ means $A=BB(0)^{-1}\in \EE_3(\AX )$. We obviously have $A(0)=\I_3$. The \plgc \ref{plcc.159} tells us that it suffices to prove the assertion after \lon at \ecoz~$(s_i)$, and Lemma~\ref{lemII3.6bis} has constructed such a family.
\end{proof}

\begin{corollary} 
\label{KuXMennickeTrivial}
\emph{(Triviality of the Mennicke symbol over $\KuX$)}\\
Let $\gK$ be a \cdiz, and  $f$, $g\in \KuX$ be \comz.
Then $\meck{f}{g} = 1$.
\end{corollary}

\begin{proof}
We reason by \recu on the number $r$ of variables in $\uX$. \\
The case $r=0$, \cad  $\KuX = \gK$ stems from $\EE_3(\gK) = \SL_3(\gK)$ ($\gK$ is a \cdiz).  For $r \ge 1$, we
suppose \spdg that $f$ is nonzero. A \cdv allows us to transform $f$ into a pseudo\polu in $X_r$ (Lemma~\ref{lemNoether}), say $f = ah$ with $a\in\gK^*$ and $h$ \mon in~$X_r$. Then, by letting $h_0 = h(X_1, \ldots, X_{r-1}, 0)$ and $g_0 = g(X_1, \ldots, X_{r-1}, 0)$, which are in $\gK[X_1, \ldots, X_{r-1}]$, we have $\meck{f}{g} = \meck{h}{g} = \meck{h_0}{g_0}$.
\end{proof}

At the end of this section the results are proved in the case of an integral \riz. They are actually true for an arbitrary \riz. 
For the \gnl case, we must refer back to \cite{Rao85a,Rao85b,Rao85c}.

\begin{theorem}\label{th3Rao}
Let $n\geq3$,  $\gA $ be an integral \ri and $f(X)$ be a \vmd in $\AX^n$.\\
Then the set
\smash{$
\fa=\sotQ{s\in\gA}{f(X)\sims{\En(\gA_s[X])} f(0)}
$}
is an \idz.
\end{theorem}

We express the same thing in the following \plgcz.

\CMnewtheorem{plgcRao}{Concrete Rao \lgb principle}{\itshape}
\begin{plgcRao}\label{plgc-Rao}\relax
Let $n\geq3$,  $\gA $ be an integral \riz,~$f(X)$ be \vmd in $\AX^n$, and  $S_1$, $\ldots$, $S_k$ be \moco of $\gA$. 
\Propeq 
\begin{enumerate}
\item ${f(X)}\sims{\En(\AX)} f(0)$.
\item ${f(X)}\sims{\En(\gA_{S_i}[X])} f(0)$ for each $i$.
\end{enumerate}
\end{plgcRao}

%
\begin{Proof}{\Demo of \thref{th3Rao}. }\\
We must show that the set
$$\preskip.0em \postskip.1em 
 \fa=\sotQ{s\in\gA}{ f(X)\sims{\En(\gA_s[X])} f(0)}
$$
is an \idz. Since all the computations in $\gA_s$ are valid in $\gA_{sa}$, we have:~$s\in\fa$ implies
$ as \in \fa$. Now let $s_1$ and $s_2$ be in $\fa$. We must show  $s_1+s_2\in\fa$, or   $1\in\fa\gA_{s_1+s_2}$. In short, we suppose that $s_1$ and $s_2=1-s_1$ are in $\fa$, and we must show~$1\in\fa$.
\\ 
By \dfnz, for $i=1$, $2$, we have a matrix $E_i=E_i(X)\in\En(\gA_{s_i}[X])$
such that $E_if(X)=f(0).$ We have $E_i(0)f(0)=f(0)$. 
Therefore, even if it entails replacing~$E_i$ by $E_i^{-1}(0)E_i$, we can assume that $E_i(0)=\In$.
\\ 
We introduce $E=E_1E_2^{-1}\in \En(\gA_{s_1s_2}[X],\gen{X})$, which gives an integer $k\geq0$ satisfying the conclusion of Lemma~\ref{lemE1Rao} for the matrix $E$ and for the two \lons $\gA_{s_1}\to \gA_{s_1s_2}$ and $\gA_{s_2}\to \gA_{s_1s_2}$.
\\
Let $c\in \gA$ with $c\equiv 1\mod s_1^k$ and $c\equiv 0\mod s_2^k$.  
We therefore have two matrices $E'_1\in \En(\gA_{s_1}[X],\lra X)$,
 $E'_2\in \En(\gA_{s_2}[X],\lra X)$ that satisfy
\begin{enumerate}
\item [--]  $E^{-1}(cX)E(X)=E'_2$  over $\gA_{s_1s_2}$ (since $c\equiv 1 \mod s_1^k$),
\item [--]  $E(cX)=E(cX)E(0\cdot X)=E'_1$ over $\gA_{s_1s_2}$ (since $c\equiv 0 \mod s_2^k$).
\end{enumerate}
We obtain  $E=E'_1E'_2=E_1E_2^{-1}$ over $\gA_{s_1s_2}$, and ${E'_1}^{-1}E_1=E'_2E_2$ over $\gA_{s_1s_2}$.
Since ${E'_1}^{-1}E_1=F_1$ is defined over $\gA_{s_1}$, that $E'_2E_2=F_2$ is defined over~$\gA_{s_2}$, that they are equal over $\gA_{s_1s_2}$, and that~$s_1$ and $s_2$ are \comz, there exists a unique matrix $F\in\Mn(\AX)$ which gives $F_1$ over $\gA_{s_1}$ and $F_2$ over~$\gA_{s_2}$.
Once again we must prove that $F\in\En(\AX)$ and~$Ff=f(0)$.
\\ 
The first item results from Lemma~\ref{lemE3Rao}.
To satisfy $Ff=f(0)$, we will assume that $\gA$ is integral, which legitimizes the following \egts over $\gA$ 
\[ \preskip.3em \postskip-.3em
\arraycolsep2pt
\begin{array}{cccccccc} 
Ff  & =  & {E'_1}^{-1}E_1f &= & {E'_1}^{-1}f(0) & =&\\[.25em] 
E^{-1}(cX)f(0)  & =  & E_2(cX)E_1^{-1}(cX)f(0) &= &E_2(cX)f(cX) &=&f(0).  
 \end{array}
\]
\end{Proof}
Note: In this \dem the last verification is the only place where we need to assume that the \ri is integral.

\begin{theorem}\label{th2Rao}
Let $n\geq3$, $\gA $ be a \ri and 
$f=\tra{\big(f_1(X),\ldots,f_n(X)\big)}$ be a \vmd in $\AX^n$,
with $1$ in the head \id of the $f_i$'s. Then
$$\preskip.2em \postskip.4em
f\;\sims{\En(\AX)} \;f(0)\;\sims{\En(\gA)}\;f\sta(0)\;\sims{\En(\AX)}\; f\sta.
$$
If one of the $f_i$'s is \monz, we have $f\sims{\En(\AX)}  
 \tra{[\,1\;0\;\cdots\;0\,]}$. 
\end{theorem}
\begin{proof}
The little Horrocks' local \tho (\thref{th2HorrocksLocal})
and Rao's \plg gives the first \eqvcz. Next we copy the \dem of Rao's \tho (\thref{th1Rao}) 
by replacing~$\GLn$ by $\En$.
\end{proof}
%

\begin{corollary} 
\label{TransitiviteEn}
\emph{(Transitivity of $\En$ for $n \ge 3$)}
\\
If $\gK$ is a \cdi and $\KuX = \gK[\Xr]$, 
then $\En(\KuX)$ acts transitively over the set of \vmds of $\KuX^n$ 
for~$n\geq 3$.
\end{corollary}

\begin{proof}
We reason by \recu on $r$. The case $r=0$ stems from the fact that $\gK$ is a \cdiz. \\
Let $r \ge 1$ and  $f = \tra{[\,f_1(\uX)\;\cdots\;f_n(\uX)\,]}$ be a \vmd of~$\KuX^n$.  
Let $\gA = \gK[X_1, \ldots, X_{r-1}]$.
One of the $f_i$'s is nonzero and a \cdv 
          allows
for the transformation into a pseudo\polu in~$X_r$ (Lemma~\ref{lemNoether}).  With $f_i$ \mon in~$X_r$, we apply \thref{th2Rao} to obtain
$$
\preskip-.2em \postskip.4em 
f\;\sims{\En(\KuX)}\; f(X_1, \ldots, X_{r-1},0). 
$$
This last vector is a \vmd of $\Ae n$. We apply the \hdrz. 
\end{proof}
%


Finally, the \dem that \thref{II3.8} implies Suslin's stability \tho is simple and \covz, as in \cite{GM}.

\begin{theorem}\label{thStabSuslCorps}
 \emph{(\Susz's stability \thoz, case of \cdisz)}\\
Let $\gK$ be a \cdiz.  For $n\geq 3$, we have $\SL_n(\KuX)= \En(\KuX)$.
\end{theorem}

\begin{proof}
Let us prove the following preliminary result.\\
\emph{For $A \in \GL_n(\KuX)$, there exist $P$, $Q\in \En(\KuX)$ such that} 

\snic{P\,A\,Q \in \GL_2(\KuX) \subseteq  \GL_n(\KuX) .\footnote{The inclusion $\GL_r \hookrightarrow \GL_n$
is defined as usual by $B \mapsto \Diag(B,\I_{n-r})$.}}

Indeed, let us consider the last row of $A$. It is a \vmdz, therefore (Corollary~\ref{TransitiviteEn}),
there exists a $Q_n \in \En(\KuX)$ such that the last row of~$A\,Q_n$ is $[\,0\;\cdots\;0\;1\,]$.
Hence $P_n \in \En(\KuX)$ such that the last column of $P_n(A\,Q_n)$ is $\tra{[\,0\;\cdots\;0\;1\,]}$, \cad  $P_n\,A\,Q_n \in \GL_{n-1}(\KuX)$. 
\\
By iterating, we find matrices $P$, $Q \in \En(\KuX)$ of the form 

\snic{P = P_3\cdots P_n$, $\;Q =
Q_n\cdots Q_3,}

such that $P\,A\,Q \in \GL_2(\KuX)$.
\\ 
If in addition $A \in \SLn(\KuX)$, we obtain $P\,A\,Q \in \SL_2(\KuX) \hookrightarrow \SL_3(\KuX)$. 
\\
We can~then consider its image in $\SL_3(\KuX)\sur{\EE_3(\KuX)}$.
\\
As the corresponding Mennicke symbol equals $1$ (Corollary~\ref{KuXMennickeTrivial}), we obtain $P\,A\,Q \in \EE_3(\KuX)$, and ultimately $A \in \EE_n(\KuX)$.
\end{proof}
%



\Exercices

\begin{exercise}
\label{exoChap18-1} 
{\rm  
Let $U \in \Ae {n \times m}$ and $V \in \Ae {m \times n}$.
\begin{itemize}
\item [\emph{1.}]
Prove, for $N \in \MM_n(\gA)$, that

\snic{(\I_m - VNU) (\I_m + VU) = \I_m + V\big( \I_n - N(\I_n + UV)\big)U.
}

Deduce that if $\I_n + UV$ is \iv with \inv $N$,   
then $\I_m + VU$ is \iv with \inv $\I_m - VNU$.

\item [\emph{2.}]
Deduce that $\I_n + UV$ is \iv if and only if $\I_m + VU$ is \ivz, and establish symmetrical formulas for their \invsz.

\item [\emph{3.}]
Show that $\det(\I_n + VU) = \det(\I_m + UV)$ in all cases.

\item [\emph{4.}]
Suppose that $\I_m + VU$ is \ivz. Show the following membership, due to Vaserstein.

\snic{\cmatrix {\I_n + UV & 0\cr 0 & (\I_m + VU)^{-1} \cr} \in \EE_{n+m}(\gA).
}

\smallskip  What happens when $VU = 0$? 

\end{itemize}
}
\end{exercise}

\vspace{-1em}
\begin{exercise}\label{exoMennicke}
{\rm 
With the notations of Lemma~\ref{lemMennicke1}, prove that the matrix $A'^{-1} A$ 
is of the form $\I_2 + uv$ with 
$u, v \in \Ae {2 \times 1}$, 
$vu = 0$ and $v$ \umdz.
}
\end{exercise}

\vspace{-1em}
\begin{exercise}\label{exoMennicke2}
{\rm 
Let $a$, $b$, $u$, $v \in \gA$ satisfy $1 = au + bv$.  Show, only using the \prts of the Mennicke symbol appearing in Proposition~\ref{propMennicke1}, that $\meck {a}{b} = \meck {u}{v} = \meck
{a-v}{b+u}$.
}
\end{exercise}

\vspace{-1em}
\begin{exercise}\label{exoStabFreeCoRang1}
{\rm 
A \stl \Amo $E$ of rank $r$ is said to be \emph{of type $t$} 
if~$E \oplus \Ae t \simeq \Ae {r+t}$. Here we are interested in the relations 
between on the one hand the \iso classes of the \stl modules of rank $n-1$, of type $1$, and on the other hand the $\GLn(\gA)$-set $\Um_n(\gA)$ consisting of the \vmds of $\Ae n$.
\begin{enumerate}\itemsep0pt

\item
Let $x \in \Um_n(\gA)$. Prove that the module $\Ae n / \gA x$ is \stl of rank $n-1$, of type $1$, and that for $x' \in \Um_n(\gA)$, we have $\Ae n / \gA x \simeq \Ae n / \gA x'$ \ssi $x \sims{\GL_n(\gA)}  x'$. 
Show that we thus obtain a (first) bijective correspondence: $x \longleftrightarrow \Ae n / \gA x$

\vbox {\advance \vsize by -4cm
\Grandcadre{$
\displaystyle{\Um_n(\gA) \over \GL_n(\gA)} \buildrel {(1)} \over \simeq
{\hbox{\stl modules of rank $n-1$, of type $1$} 
\over \hbox {\iso}}
$}}

What are the \vmds that correspond to a free module?

\item
Let $x \in \Um_n(\gA)$. Show that $x^\perp \eqdefi  \Ker \tra {x}$ 
is a \stl module of rank $n-1$, 
of type $1$, 
and that for $x' \in \Um_n(\gA)$, 
we have $x^\perp \simeq x'^\perp$ if and only if $x \sims{\GL_n(\gA)} x'$.  
Prove that we thus obtain a (second) bijective correspondence: 
$x \longleftrightarrow x^\perp$

\vbox {\advance \vsize by -4cm
\Grandcadre{$
\displaystyle{\Um_n(\gA) \over \GL_n(\gA)} \buildrel {(2)} \over \simeq
{\hbox{\stl modules of rank $n-1$, of type $1$} 
\over \hbox {\iso}}
$}}

\item
If $E$ is \stl of rank $r$ and of type $t$, the same goes for its dual $E\sta$. For $t = 1$, describe the involution of $\Um_n(\gA)/\GLn(\gA)$ induced by the involution $E \leftrightarrow E\sta$.

\item
Let $x$, $x'$, $y \in \Ae n$ such that $\tra {x}y = \tra {x'}y = 1$. Why do we have $x \sims{\GL_n(\gA)}  x'$? Explicate $g \in \GL_n(\gA)$ such that $gx = x'$, $g$ of the form $\I_n + uv$ with  
$vu = 0$ and $v$ \umdz. Deduce that for $n \ge 3$, $g \in \En(\gA)$, and so $x \sims{\En(\gA)}  x'$.

\end{enumerate}

}
\end{exercise}

\vspace{-1.2em}
\pagebreak	

\begin{exercise}\label{exoStabFree2}
       (Autodual \stl modules of type $1$)
{\rm 
\begin{enumerate}\itemsep0pt
\item
Let $a$, $b \in \gA$, $x = (x_1, \ldots, x_n) \in \Ae n$ with $n \ge 3$ and $ax_1 + bx_2$ be \iv modulo $\gen {x_3, \ldots, x_n}$ (in particular, $x$ is \umdz). \\
Let $x' = (-b,a, x_3, \ldots, x_n)$.
Explicate  $z \in \Ae n$ such that $\scp {x}{z} = \scp {x'}{z} = 1$.
Deduce, for $G = \GLn(\gA)$ (or better yet for $G = \En(\gA)$), that

\snic{
x \,\sims G \, x'\, \sims G \, (a,b, x_3, \ldots, x_n).
}

\item
Let $x, y \in \Ae 4$ such that $\scp {x}{y} = 1$.  Show that $x \sims {\EE_4(\gA)}  y$. 
In particular, the \stl module $x^\perp = \Ker \tra {x}$ is \isoc to its dual.

\item
Analogous question to the previous one by replacing $4$ with any even number $n \ge 4$. 

\end{enumerate}
}

\end{exercise}



\sol

\exer{exoChap18-1}{
\emph{2}. We establish the formulas 

\snic{N = (\I_n + UV)^{-1} = \I_n - UMV, $ $  M = (\I_m + VU)^{-1} = \I_m - VNU}

\emph{4}.
We know that $\I_n + UV$ is \ivz. Let $N = (\I_n + UV)^{-1}$,
$M  = (\I_m + VU)^{-1}$. We have therefore $N + UVN = \I_n = N + NUV$ and $M + VUM = \I_m = M + MVU$. \\
We realize the following \elr operations

\snic{\cmatrix {\I_n + UV & 0\cr 0 & M\cr} \cmatrix {\I_n  & -NU\cr 0 & \I_m\cr} =
\cmatrix {\I_n + UV & -U\cr 0 & M\cr},}

\snic{\cmatrix {\I_n + UV & -U\cr 0 & M\cr} \cmatrix {\I_n & 0\cr V & \I_m\cr} =
\cmatrix {\I_n & -U\cr MV & M\cr},
}

then

\snic{
\cmatrix {\I_n & -U\cr MV & M\cr} \cmatrix {\I_n & U\cr 0 & \I_m\cr} =
\cmatrix {\I_n & 0\cr MV & MVU + M\cr} = \cmatrix {\I_n & 0\cr MV & \I_m\cr},
}

and finally

\snic{
\cmatrix {\I_n & 0\cr MV & \I_m\cr} \cmatrix {\I_n & 0\cr -MV & \I_m\cr} =
\cmatrix {\I_n & 0\cr 0 & \I_m\cr}.
}

We therefore have explicated matrices $\alpha$, $\beta$, $\gamma$, $\delta \in
\EE_{n+m}(\gA)$ such that

\snic{\cmatrix {\I_n + UV & 0\cr 0 & (\I_m + VU)^{-1}\cr}
\,\alpha\,\beta\,\gamma\,\delta = \I_{n+m},}

hence
\[\preskip.0em \postskip.4em 
\begin{array}{l} 
\cmatrix {\I_n + UV & 0\cr 0 & (\I_m + VU)^{-1}\cr} = 
\delta^{-1}\, \gamma^{-1}\, \beta^{-1}\, \alpha^{-1} = \\[1.3em]
\cmatrix {\I_n & 0\cr MV & \I_m\cr} 
\cmatrix {\I_n & -U\cr 0 & \I_m\cr} 
\cmatrix {\I_n & 0 \cr -V & \I_m\cr} 
\cmatrix {\I_n & NU\cr 0 & \I_m\cr} . 
 \end{array}
\] 
In the special case where $VU = 0$, we have shown that

\snic{
\cmatrix {\I_n + UV
& 0\cr 0 & \I_m \cr} \in \EE_{n+m}(\gA).}

}

\exer{exoMennicke} 
By using $ad' = 1+bc'$, $ad = 1+bc$, we obtain for $A'^{-1}A$

\snic{
\Cmatrix{3pt} {d'&-b\cr -c'&a\cr} \cmatrix {a&b\cr c&d\cr} =
\cmatrix {ad'-bc & bd'-bd\cr ac-ac'&ad-bc'\cr} =
\cmatrix {1 + b(c'-c) & b(d'-d)\cr a(c-c') &1+b(c-c')\cr}.}

By replacing $b(c'-c)$ with $a(d'-d)$, we see that $A'^{-1}A = \I_2 + uv$
with

\snic{
u=\cmatrix {d'-d\cr c-c'\cr}, \,v=\cmatrix {a & b}, \, vu = 0, 
\hbox{ and }v  \hbox{ \umdz}.
}

\exer{exoMennicke2}{
We have $\meck {au}{b} = \meck {a}{b}\meck {u}{b}$.\\
But $au = 1-bv$ so $\meck {au}{b} = \meck {1-bv}{b} = \meck {1}{b} = 1$. 
Recap: $\meck {a}{b}\meck {u}{b} = 1$. \\
Similarly, $\meck {u}{b}\meck {u}{v} = 1$, so
$\meck {a}{b} = \meck {u}{v}$. \\
Finally, $(a-v)u + (b+u)v = 1$, 
so $\meck {a-v} {b+u} = \meck {u} {v}$.
}

\exer{exoStabFreeCoRang1}{
 \emph{1.} 
Let $y \in \Ae n$ such that $\tra {y} x = 1$.\\
We have $\Ae n = \gA x \oplus \Ker
\tra {y}$ and so $\Ae n/\gA x \simeq \Ker \tra {y}$ is \stlz.  \\
{If $x \sims{\GL_n(\gA)} x'$}, 
it is clear that $\Ae n / \gA x \simeq \Ae n / \gA x'$. \\
Conversely, let $\varphi : M = \Ae n/\gA x \to M' = \Ae n/\gA x'$ be an \isoz.   
We have $\Ae n \simeq M \oplus \gA x 
\simeq M' \oplus \gA x'$. 
We define $\psi : \gA x \to \gA x',\; 
ax \mapsto ax'$.  
Then $\varphi \oplus \psi$ seen in $\GLn(\gA)$ transforms $x$ into $x'$, so $x \sims{\GL_n(\gA)} x'$.
\\
 A \vmd $x \in \Ae n$ provides a free module $\Ae n / \gA x$ \ssi $x$ is part of a basis of $\Ae n$.

 \emph{2.}
Let $M = x^\perp$, $M' = x'^\perp$ and assume $M \simeq M'$.  By denoting by $\mathring M \subseteq (\Ae n)\sta$ the \ort of $M \subseteq \Ae n$, we have $\mathring M = \gA\tra {x}$ and $\mathring M' = \gA\tra {x'}$. If $\scp {x}{y} = 1$, $\scp {x'}{y'} = 1$, we have $\Ae n = \gA y \oplus M = \gA y' \oplus M'$, hence an \auto of $\Ae n$ transforming $M$ into $M'$ (send $y$ to $y'$), then by duality, an \auto $u$ of $(\Ae n)\sta \simeq \Ae n$ transforming $\gA\tra{x}$ into $\gA\tra{x'}$. We deduce $u(\tra{x}) = \varepsilon\tra{x'}$
with $\varepsilon \in \Ati$. Then, $\varepsilon^{-1}\tra{u}$ transforms $x$ into $x'$.

 \emph{3.}
If $G = E \oplus F$, then $G\sta \simeq E\sta \oplus F\sta$; with $G = \Ae {r+t} \simeq G\sta$, $F = \Ae {r} \simeq F\sta$, we obtain the result.
 The involution induced over $\Um_n(\gA)/\GLn(\gA)$ is the following: to the class modulo $\GLn(\gA)$ of $x \in \Um_n(\gA)$, we associate the class modulo $\GLn(\gA)$ of an \elt $y \in \Um_n(\gA)$ that satisfies $\scp {x}{y} = 1$. Naturally, there are several~$y$ that are suitable but their class modulo $\GLn(\gA)$ is well-defined.

 \emph{4.}
We have $\Ae n = \gA y \oplus x^\perp = \gA y \oplus x'^\perp$ hence $x^\perp \simeq x'^\perp \simeq \Ae n / \gA y$ so $x \sims{\GLn(\gA)} x'$.  To determine $g \in \GLn(\gA)$ realizing $gx = x'$, we use $\Ae n = \gA x \oplus y^\perp = \gA x' \oplus y^\perp$. \Gnltz, let $G = E \oplus F = E' \oplus F$; to explicate some \auto of $G$ mapping $E$ to $E'$, we proceed as follows. Let $\pi$ be the \prn over $E$, $\pi'$ be the \prn over $E'$ and $p = \I_G - \pi$, $p' = \I_G - \pi'$. \\
The \prrs $p$ and $p'$ have the same image $F$. Let $h = p'-p = \pi-\pi'$. \\
We obtain $h^2 = 0$ and $(\I_G - h) p (\I_G + h) = p'$, or $(\I_G - h) \pi (\I_G + h) = \pi'$. \hbox{Therefore $\I_G - h$} is an \auto of $G$ transforming $\Im\pi = E$ into $\Im\pi' = E'$. Here $E = \gA x$, $E' = \gA
x'$, $F = y^\perp$, so
$$\preskip.2em \postskip.4em
\pi(z) = \scp {z}{y}x, \quad \pi'(z) = \scp {z}{y}x', \quad
h(z) = \scp {z}{y} (x - x').
$$
The desired \auto of $\Ae n$ that transforms $x$ into $x'$ is therefore
$$\preskip.3em \postskip.4em
\I_n - h : z \mapsto z + \scp {z}{y} (x' - x) 
\quad \hbox {i.e.} \quad
\I_n - h = \I_n + uv 
$$
with $u = x'-x \in \Ae {n \times 1}$, $v = \tra {y} \in \Ae {1 \times n}$;
we indeed have $vu = 0$ and $v$ \umdz.

}

\vspace{-.5em}
\pagebreak	

\exer{exoStabFree2}{
\emph{1.}
The key to the \pb is found in the double \egtz, for some $u$ in~$\gA$,
$z_1 = u(a+x_2)$, $z_2 = u(b-x_1)$, which implies 

\snic{
z_1x_1 + z_2x_2 = u (ax_1 + bx_2) = z_1 b + z_2 (-a).} 

Let $u$ such that $u(ax_1 + bx_2) + z_3x_3 + \cdots + z_nx_n = 1$ and $z = (z_1, z_2, z_3, \ldots, z_n)$. We then have $\scp {z}{x} = \scp {z}{x'} = 1$. By Exercise~\ref{exoStabFreeCoRang1}, $x \sims G  x'$.  
As $(b, -a) \sims{\EE_2(\gA)}  (a,b)$, we have $x \sims G (a,b, x_3, \ldots, x_n)$.

 \emph{2.}
As $x_1y_1 + x_2y_2 + x_3y_3 + x_4y_4 = 1$, we have

\snic{(x_1,x_2,x_3,x_4)\, \sims G \,
(y_1,y_2,x_3,x_4)\, \sims G  \,
(y_1,y_2,y_3,y_4)} 

The rest of the question immediately stems from this.

 \emph{3.}
Analogous method to the previous question.

}




\Biblio
Section~\ref{secMennicke} and the \dem of \thref{thStabSuslCorps} follow very closely the presentation in \cite{GM}.
For the most part we have only transformed a few abstract \lgb arguments into concrete arguments via the use of the \lgbe machinery with \ideps explained in Section~\ref{secMachLoGlo}.\imlb

Section~\ref{secCompVmdsPols} is directly inspired by \cite[chapter VI, section 2]{Lam06}.

\newpage \thispagestyle{CMcadreseul}

\incrementeexosetprob

\appendix
\let\showchapter\oldshowchapter
\let\showsection\oldshowsection


\def\thechapter{A}


\def\chaptername{Annex}
\refstepcounter{chapter}

\chapter*{Annex. Constructive logic}
\addcontentsline{toc}{chapter}{\hspace{-0.8em}Annex. Constructive logic}
\mtcaddchapter
\markboth{Annex}{Constructive logic}
\label{chapPOM}
\perso{compil\'e le \today}

\minitoc

\newpage	
\Intro

This annex is devoted to presenting a few basic concepts of \coma in Bishop's style, illustrated by the three founding works \cite{B67,BB85,MRR}.

By constructive logic, we mean the logic of \comaz.

\setcounter{section}{0}

\section{Basic objects, Sets, Functions}
 Non-negative integers
 and constructions are two primitive notions.
They cannot be defined.

Other primitive notions are closely tied to common language and are difficult to place.
 For example the \egt of the number $2$ in two distinct occurrences.

The formalization of a piece of \maths can be used to better understand 
what we are doing 
          to it.  
However, to speak about a formalism it is necessary to understand a lot of things that are of the same type of complexity as the non-negative integers.
Thus, the formalism is only a tool and it cannot replace  intuition and  basic experience
 (for example the  non-negative integers,
the constructions): as powerful as a computer may be, it will never understand \gui{what it does,} or, as Ren\'e Thom used to say, \gui{All that is rigorous is insignificant.}

\subsubsec{Sets}

\noi A \emph{set} $(X,=_X,\neq_X)$ is defined by saying:

\noi --- how to construct an \elt of the set (we say that we have defined a \emph{preset} $X$)
\index{preset}

\noi --- what is the meaning of the \emph{\egtz} of two \elts of the set (we have to prove that it is indeed an \eqvc relation)

\noi --- what is the meaning of the \ix{distinction}{\footnote{This terminology \emph{is not} a
homage to Pierre Bourdieu. All in all, we prefer \emph{distinction} to \emph{non-\egtz}, which  presents the disadvantage of a negative connotation, and to \emph{in\egtz} which is rather used in the context of order relations. For the real numbers for example, it is the \egt and not the distinction that is a negative assertion.}} of two \elts of the set (we then say that the \elts are \emph{discernible} or \emph{distinct}). We need to show the following \prtsz:

\hspace*{5mm}-- $ \; (x\neq_X y \; \land \; x=_Xx' \; \land \; y=_Xy')\; \Rightarrow \;
x' \neq_X y'$,

\hspace*{5mm}-- $ \; x\neq_X x$  is impossible,

\hspace*{5mm}-- $ \; x\neq_X y\; \Rightarrow \; y \neq_X x$.

\ss Ordinarily, we drop the index $X$ for the symbols $=$ and $\neq $. If the distinction is not specified, it is implicitly defined as meaning the absurdity of the \egtz.

\ss A distinction relation is called a \emph{separation} relation if it satisfies the following \ix{cotransitivity} \prt (for three arbitrary \elts $x,y,z$ of $X$):\index{separation}

\hspace*{5mm}-- $\; \; x\neq_X y\; \Rightarrow \; (x \neq_X z\; \lor\;  y\neq_X z).$

\noi A separation relation $\neq_X$ is said to be \emph{narrow} if $x=_X y$ is equivalent to the absurdity of  $x\neq_X y$. In a set with a narrow separation, 
 distinction is often more important than \egtz. \index{separation!narrow ---}

\ss A set  $(X,=_X,\neq_X)$ is said to be \emph{discrete} \label{discret} if we have  $$ \forall x,y\in X\; (x=_X y \lor x\neq _X y).$$  In this case the distinction is a narrow separation and it is equivalent to the absurdity of the \egtz.%
\index{set!discrete ---}

\subsubsec{The  non-negative integers}

\noi The set $\NN=\so{0,1,2,\ldots}$ of  non-negative integers
is considered as a priori well-defined. However, note that \cot this is a \ix{potential infinity} and not an \ix{actual infinity}.%
\index{infinite!potiential ---}\index{infinite!actual ---} 
By the idea of a potential infinite we mean that the infiniteness of $\NN$ is apprehended as an essentially negative notion; we never stop exhausting the non-negative integers.
  On the contrary, the semantic of $\NN$ in \clama is that of a completed infinite, which exists \gui{somewhere,} at least in a purely ideal way.

A non-negative integer
 can be encoded in the usual way. The comparison of two integers given in a coded form can be made reliably.
In short, the set of non-negative integers
  is a discrete set and the order relation is \emph{decidable}
$$\preskip.4em \postskip.1em
\forall n,m\in\NN\; \; \; (n<m\; \; \lor\; \;  n=m\; \; \lor\; \;  n>m)$$

\subsubsec{Sets of pairs}

\noi When two sets are defined their \emph{Cartesian product} is also naturally defined: the fabrication of the pairs of objects is an \elr construction.  Equality and distinction over a Cartesian product are naturally defined.

\subsubsec{Functions}

\noi The set $\NN^\NN$ of sequences of non-negative integers   depends on the primitive notion of construction. An \elt of $\NN^\NN$ is a construction that takes as input an \elt of $\NN$ and gives as output an \elt of $\NN$.  The \egt of two \elts in $\NN^\NN$ is the \emph{extensional \egtz}
$$(u_n)=_{\NN^\NN}(v_n) \quad \mathrm{ signifies} \quad  \forall n\in\NN\;\;
u_n=v_n.
$$
Thus, the \egt between two \elts of $\NN^\NN$ a priori asks for an infinity of \gui{\elr computations,} actually the \egt demands a \demz.\index{extensional!equality}

\noi The distinction of two \elts of $\NN^\NN$ is the \emph{extensional distinction} relation
$$(u_n)\neq_{\NN^\NN}(v_n) \quad {\eqdef} \quad  \exists n\in\NN\;\;  u_n\neq
v_n.
$$
Thus, the distinction of two \elts of $\NN^\NN$ can be observed by a simple computation.\index{extensional!distinction}

\begin{example}
\label{exo neq in N^N}\relax
\emph{
The distinction of $\NN^\NN$ is a narrow separation relation. }
\end{example}

Cantor's diagonalization argument is \cofz. It shows that $\NN^\NN$ is \emph{much more complicated} than $\NN$. From a \cof point of view, $\NN$ and $\NN^\NN$ are only potential infinities: it holds no meaning to say that a potential infinity is \emph{greater} than another.

\mni \emph{Digression.}
When you say \gui{I give you a sequence of non-negative integers,}
  you must prove that the construction $\; n\mapsto u_n\; $ that you propose works for any input $n$. Moreover, when you say \gui{Let us consider an arbitrary sequence of non-negative integers
 $(u_n)_{n\in\NN}$,}
the only thing that you know for certain is that for all $n\in\NN$, you have $u_n\in\NN$, and that this $u_n$ is nonambiguous: you can for example conceive the sequence as given by an oracle.
Actually, you could a priori ask, symmetrically, what exactly is the construction $\; n\mapsto u_n$,
and a \dem that this construction works for every input $n$.

However, in the constructivism \`a la Bishop, we make no specific assumptions regarding \gui{what the legitimate constructions from $\NN$ to $\NN$ are,} nor 
on \gui{what \prmt is a \dem that a  construction works.} Thus we are in a dissymmetrical situation.

This dissymmetry has the following consequence. Everything you prove has a computational content, but everything you prove is also valid from a classical point of view. 
Classical mathematics could regard constructive mathematics as only speaking of constructive objects, and Bishop's constructive mathematics is certainly primarily interested in constructive objects (see \cite{bi}). But in fact, the \cov  \dems \`a la Bishop    work for any type of mathematical object.\footnote{\ldots~if there exist non\cov mathematical objects.} The \thos that we find in \cite{BB85} and \cite{MRR} are valid in \clamaz, but they also support the Russian \cov interpretation (in which all the mathematical objects are words from a formal language that we could fix once and for all) or yet again Brouwer's intuitionist philosophy,
 which has a significantly idealistic component. 
\eoe

\ss After this digression let us get back on topic: functions. \Gnltz,  a \emph{function} $f:X\rightarrow Y$ is a construction that takes as input some $x\in X$ and a \dem that $x\in X$, and gives as output some $y\in Y$ and a \dem that $y\in Y$. In addition, this construction must be \emph{extensional}
$$\preskip.4em \postskip.4em
x=_Xx'\Rightarrow f(x)=_Yf(x')\qquad  \mathrm{and} \qquad  f(x)\neq_Yf(x')
\Rightarrow x\neq_Xx'.
$$
When $X$ and $Y$  are well-defined sets, we consider (in \coma \`a la Bishop) that the set $\cF(X,Y)$ of functions~\hbox{$f:X\rightarrow Y$} is also well-defined. For the \egt and the distinction we take the usual extensional \dfnsz.

\ss A function $f:X\rightarrow Y$ is \emph{injective} if it satisfies
$$\preskip.3em \postskip.0em
f(x)=_Yf(x')\Rightarrow x=_Xx'\quad \mathrm{and} \quad x\neq_Xx'
\Rightarrow f(x)\neq_Yf(x'). 
$$

\subsubsec{Finite, bounded, \enums and \denbs sets}

\noi We now give a certain number of pertinent \cov \dfns related to the concepts of finite, infinite and  \denbs sets in \clamaz.\rdb
\begin{itemize}
\item  A set is said to be \emph{finite} if there is a bijection between this set and the set of integers $<n$ for a certain integer $n$ (this is the \dfn given \paref{Deux mots}).
\item  A  set $X$  is said to be \emph{finitely \enumz} if there is a surjective map $[0,n[\;\to X$ for some integer $n$ (this is the \dfn given \paref{Deux mots}).%
\item  A preset $X$ is said to be \emph{\enumz} if we have given a means to enumerate it that allows it to possibly be empty, which happens in practice as follows.\footnote{The \dfn given on \paref{Deux mots} is only for nonempty sets.} We give some $\alpha\in\so{0,1}^\NN$ and some operation $\varphi$ that satisfy the following two assertions:\\
-- if $\alpha(n)=1$ then $\varphi$ constructs from the input $n$ an \elt of $X$,\\
-- every \elt of $X$ is constructed as such.%
\index{set!enumerable ---}\index{enumerable!set}%
\index{set!finitely enumerable ---}\index{finitely enumerable!set}%
\item  A set is said to be \emph{\denbz} if it is \enum (as a preset) and discrete.%
\index{set!countable ---}\index{countable!set}
\item  If $n$ is a nonzero integer, we say that a set \emph{has at most $n$ \eltsz} if for every family $(a_i)_{i=0,\ldots,n}$ in the set there exist integers $h$ and $k$ ($0\leq h<k\leq n$) such that $a_h=a_k$.
\item  A set $X$ is \emph{bounded in number}, or \emph{bounded}, 
 if there exists some nonzero integer $n$ such that $X$ has at most $n$ \elts (this is the \dfn given \paref{ensborn}).
\index{bounded (in number)!set}
\index{set!bounded ---}
\item  A set $X$  is \emph{weakly finite} if for every sequence $(u_n)_{n\in\NN}$ in $X$ there exist $m$ and $p>m$ such that $u_m=u_p$.
\index{set!weakly finite ---}\index{weakly finite!set}
\item A set $X$ is \emph{infinite} if there exists an injective map $\NN\rightarrow X$.
\index{set!infinite ---}\index{infinite!set}
\end{itemize}

\begin{example}
\label{exo inf deno1}\relax
{\rm An infinite and \denb set can be put in bijection with~$\NN$.
}
\end{example}

\rdb
\subsubsec{Subsets of a set}\label{P(X)}\relax

\noi A subset of a set $(X,=_X,\neq_X)$ is defined by a \prt $P(x)$ \emph{regarding the \elts of $X$}, \cad satisfying 
$$
\forall x,y\in X\;\big( \; (\; x=y\; \land \; P(x)\;)\; \;
\Longrightarrow\; \;  P(y)\; \big).
$$
An \elt of the subset $\sotq{x\in X}{P(x)}$ is given by a pair $(x,p)$ where~$x$ is an \elt of $X$ and $p$ is a \dem that $P(x)$.{\footnote{For example, a nonnegative real number is \emph{slightly more than} a real number.}} Two \prts concerning the \elts of $X$ define the same subset when they are \eqvesz.

We can also present this as follows, which, although amounting to the same thing, causes a slightly milder headache to the newcomer. 
A subset of $X$ is given by a pair $(Y,\varphi)$  where  $Y$ is a set and $\varphi$ is an injective function of $Y$ into $X$.{\footnote{For example we can define the real numbers $\geq 0$ as those that are given by the Cauchy sequences of non-negative rational numbers.}}
Two pairs $(Y,\varphi)$ and $(Y',\varphi')$ define the same subset of $X$ if we have
$$
\forall y\in Y\; \exists y'\in Y'\; \varphi(y)=\varphi'(y')\;\;
\hbox{ and }\;\;
\forall y'\in Y'\; \exists y\in Y\; \varphi(y)=\varphi'(y').
$$

In \coma the subsets of a set $X$ are not considered to form a set, but a \ix{class}. This class is clearly not a set
(in the sense given earlier).
The intuition is the following: the sets are sufficiently well-defined classes so that we can universally or existentially quantify over their \eltsz. For this, it is necessary for the procedure of construction of \elts to be clear.

\ss Recall that a subset $Y$ of $X$ is said to be \emph{detachable} when we have a test for \gui{$x\in Y$?} when $x\in X$. The detachable subsets of $X$ form a set that can be identified with $\{0,1 \}^X$.

Constructively, we do not know of any detachable subsets of $\RR$, besides $\emptyset$ and $\RR$: \emph{there are no holes in the continuum without the logic of the excluded middle.}

\mni\rem  An interesting \cov variant for \gui{a subset $Y_1$ of $X$} is obtained by considering a pair $(Y_1,Y_2)$ of subsets of $X$ that satisfy the following two \prts
$$\forall x_1\in Y_1\; \forall  x_2\in Y_2\;  \; x_1\neq_X x_2\quad
\mathrm{and}\quad
\forall x\in X\; \lnot( x\notin Y_1\; \land\;  x\notin Y_2).
$$
The \emph{complement} is then given by the pair  $(Y_2,Y_1)$, which re-establishes a certain symmetry.\eoe

\sni\emph{The class of subsets of a set}~
\rdb
\\
Let  $\rP(X)$ be the class of subsets of the set $X$.
If we admitted $\rP(\{0\})$ as a set, then $\rP(X)$ would \egmt be a set and there would be a natural bijection between  $\rP(X)$  and $\cF\big(X,\rP(\{0\})\big)=\rP(\{0\})^X$.

This shows that all the difficulties with the set of subsets are focused on
 the class $\rP(\{0\})$, \cad the class of \emph{truth values}. In \clamaz, we admit that this class is a set with two \eltsz. This is the \emph{Law of  Excluded Middle} \TEMz:
$$
\rP(\{0\})=\{\{0\},\emptyset  \}
$$
(the class of truth values reduces to the set $\{\Vrai,\Faux \}$) and we obviously no longer have any issues with $\rP(X)$.


\section{Asserting means proving} \label{secAffirmerProuver}
In  \coma  truth is also the result of a construction. 
If $P$ is a mathematical assertion, we write \gui{$\; \vda P \; $} for \gui{we have a proof of $P$.}

\ss The \elr assertions can be tested by simple computations. 
For example, the comparison of two non-negative integers. 
When an assertion means an infinity of \elr assertions (e.g. the Goldbach conjecture\footnote{Every even number $\geq 4$ is the sum of two prime numbers.}), 
constructive mathematics consider it not to be a priori
\gui{true or false.} A fortiori, the assertions having an even greater logical complexity are not considered (from a \cof point of view) as having a priori the truth value $\Vrai$ or $\Faux$.

This must not be necessarily considered as a philosophical position 
concerning truth, but it is surely a mathematical position concerning mathematical assertions.
Actually, it is necessary to assume this position; in order to be of computational significance, all theorems must be proven constructively.

\mni \emph{Downright philosophical digression.} This position is also to be distinguished from the position that consists in saying that there certainly are different possible mathematical universes, for instance one in which the
 continuum hypothesis{\footnote{The continuum hypothesis  is, in classical set theory, the assertion that there is no cardinal strictly between that of $\NN$ and that of $\RR$, in other words, that every infinite subset of $\RR$ is equipotent to $\NN$ or to $\RR$.}}  is true, another in which it is false. This position is naturally perfectly defendable (Cantor, and no doubt G\"odel, would have rejected it in the name of a Platonic realism of  Ideas), but it is of little interest to \coma \`a la Bishop which have 
as its object of study an abstraction of the concrete universe of finite computations, with the idea that this abstraction must correspond as closely as possible to the reality that it wants to describe.
Thus, the continuum hypothesis is in this framework rather considered as empty of meaning, because it is vain to want to compare potential infinites according to their size. 
If we desire  to compare them according to their complexity, we quickly realize that there is no hope of 
defining a true total order relation on this mess.
Consequently, the continuum hypothesis today seems to be nothing other than a game of experts in the formal theory of $\ZF$. But each one of us is free to believe Plato, or even Cantor, or Zermelo-Frankel, or yet again Ð why not Ð to believe in the multiplicity of worlds. No one will ever be able to prove the latter wrong. In fact nothing says that the $\ZF$ game will not one day prove to be really useful, for instance in understanding certain subtle points of mathematics that have a concrete meaning.
\eoe

\section{Connectives and quantifiers} 
\label{secBHK}

Here we give the \gui{Brouwer-Heyting-Kolmogorov} explanation for the constructive meaning of the usual logical symbols. They are only informal explanations, not \dfnsz.\footnote{For Kolmogorov's point of view, more \prmt on \gui{the logic of \pbsz}, see \cite[Kolmogorov]{Kol} and \cite[Coquand]{CoqK}.}

These are \gui{detailed} explanations, as for the logical connectives 
and the \qtfsz,
regarding what we mean by the slogan \gui{asserting means proving.} When we write $\vda P$ we imply that we have a \cov \dem of $P$. We will make this explicit by giving a name, for example
$p$, to this \mathe object that is the \dem of $P$.
Then the explanations regard these particular objects $p$, but all of this remains informal.

\mni \textbf{Conjunction:} $\vda P \; \land\;  Q$  means: \gui{$\vda P$  and $\vda Q$}
(as for classical logic). In other terms: a \dem of $P \;\land\;  Q$ is a pair $(p,q)$ where $p$ is a \dem of $P$ and $q$ a \dem of~$Q$.

\mni \textbf{Disjunction:} $\vda P \; \lor\;  Q$  means: \gui{$\vda P$  or $\vda Q$} (which does not work with classical logic).
In other terms: a \dem of $P \, \lor\,  Q$ is a pair $(n,r)$ with $n\in\{0,1 \}$.
If $n=0$, $r$ must be a \dem of $P$, and if~$n=1$, $r$ must be a \dem of~$Q$.

\mni \textbf{Implication:} $\vda P \; \Rightarrow \;  Q$  has the following meaning: \\
a \dem of $\; P \; \Rightarrow \;  Q\; $ is a cons\-truction $p\mapsto q$ that transforms every \dem $p$ of $P$ into a \dem $q$ of~$Q$.

\mni \textbf{Negation:} $\lnot P $  is an abbreviation of $P \; \Rightarrow \; 0=_\NN1$.

\mni \textbf{Universal \qtfz:} (similar to  implication). \emph{A
\qtn is always a \qtn on the objects of a previously defined set.}
Let $P(x)$ be a \prt regarding the objects $x$ of a set $X$.
\\
Then   $\vda \forall x\in X\; \; P(x)$ has the following meaning:
we have a construction $(x,q)\mapsto p(x,q)$ that takes as input any pair $(x,q)$, where $x$ is an object and  $q$ is a \dem that $x\in X$, and gives as output a \demz~\hbox{$p(x,q)$} of the assertion~$P(x)$.

 For a \qtn on $\NN$,  giving a non-negative integer 
$x$ (in the standard form) suffices to prove that $x\in\NN$: the proof $q$ in the pair $(x,q)$ above can be omitted.

\begin{example}\label{exo log0}\relax
\emph{Suppose that the \prts $P$ and $Q$ depend on a variable $x\in \NN$.
Then a \dem of $\,  \forall x\in \NN \;\big(P(x) \, \lor\,  Q(x)\big)\, $ is a construction 
$\NN\ni x\mapsto \big(n(x),r(x)\big)$, where $n(x)\in\{0,1 \}$:
 if $n(x)=0$, $r(x)$ is a \dem of $P(x)$, and if $n(x)=1$, $r(x)$ is a \dem of~$Q(x)$.
}\eoe
\end{example}

\mni \textbf{Existential \qtfz:} (similar to disjunction) \emph{A \qtn is always a \qtn on the objects of a previously defined set.} Let $P(x)$ be a \prt regarding the objects $x$ of a set~$X$. Then    $\vda \exists x\in X\;  P(x)$
 has the following meaning: a \dem of $\exists x\in X\;  P(x)$ is a triple $(x,p,q)$ where $x$ is an object, $p$ is a \dem of $x\in X$,  and  $q$ a \dem of~$P(x)$.

\begin{example}
\label{exo log1}\relax
\emph{Let $P(x,y)$ be a \prt regarding the non-negative integers~$x$  and $y$.
 Then the assertion
$$\vda \forall x\in\NN\; \; \exists y\in\NN\;\;  P(x,y)$$
means: here is a pair $(u,p)$ where $u$ is a construction $u:x\mapsto y=u(x)$ from $\NN$ to $\NN$ and
$p$ is a \dem of $\vda \forall x\in\NN\; P\big(x,u(x)\big)$.
}\eoe
\end{example}

\begin{example}
\label{exo log2}\relax (Propositional logics)\\
{\rm The class of truth values in \coma is a \agHz.\\
NB: By $\rP(\so{0})$ being a class and not a set we simply mean that the connectives
 $\vi$, $\vu$ and $\im$
and the constants $\Vrai$ and $\Faux$ satisfy the axioms of the \agHsz.

In particular, let $A,\; B,\; C$ be \mathe \prtsz. We have the following  \eqvcsz.
\begin{itemize}
\item [$\vda$]  $\big((A\Rightarrow C)\; \land\; (B\Rightarrow C)\big)\; \;
\Longleftrightarrow \; \; \big((A\; \lor \; B)\; \Rightarrow \; C\big) $
\item [$\vda$]  $\big(A\Rightarrow (B\Rightarrow C)\big)\; \; \Longleftrightarrow \; \;
\big((A\; \land \; B)\; \Rightarrow \; C\big)$
\item [$\vda$]  $\lnot (A\; \lor B\; ) \; \; \Longleftrightarrow \; \;
        (\lnot A \land \lnot B)$
\item [$\vda$]  $(A\Rightarrow B)\; \; \Longrightarrow \; \;
(\lnot B\; \Rightarrow  \; \lnot A)$
\item [$\vda$]  $\lnot\lnot\lnot A\; \; \Longleftrightarrow \; \;\lnot A$
\end{itemize}
In addition, if we have $ \vda A\; \lor \lnot A$ and $\vda B\; \lor \lnot B$, then we have
\begin{itemize}
\item [$\vda$] $\lnot \lnot A\; \; \Longleftrightarrow \; \; A$
\item [$\vda$] $\lnot (A\; \land B\; ) \; \; \Longleftrightarrow \; \;
        (\lnot A \lor \lnot B)$
\item [$\vda$] $(A\Rightarrow B) \; \; \Longleftrightarrow \; \;(\lnot A \lor B)$
\eoe\end{itemize}
}
\end{example}

\rem
\label{rem propneg}\relax
Since  $\lnot\lnot\lnot A\,\Leftrightarrow \,\lnot A$, a \prt $C$ is \eqv to a \prt $\lnot B$ (for a  certain \prt $B$ not yet specified) \ssi
$\lnot\lnot C\,\Rightarrow \,C$. Thus, in \coma we can define the concept of \emph{negative \prtz}. In \clamaz, the concept is pointless since every \prt is negative. In \comaz, care must be taken because $\Vrai$ is also a negative \prtz, since $\Faux \Rightarrow \Faux$, $\lnot \Faux$  is equal to $\Vrai$. 
\eoe

\section{Mechanical computations}\label{AnnexeCalculsMec}

Here we discuss a point that classical mathematicians often fail to appreciate. A function from $\NN$ to $\NN$ is given by a construction. The usual constructions correspond to algorithmic programs that can run on an \gui{ideal} computer.{\footnote{A computer having all the space and time necessary for the considered computation.}} This leads to the notion of \emph{mechanical computations}. A function  $f\in\NN^\NN$ obtained by such a mechanical computation is called a \emph{recursive function}.
\\
The subset $\Rec\subset \NN^\NN$ formed by the recursive functions can then be described more formally as we will now explain.

Recall that a \emph{primitive recursive function} is a function $\NN^k\rightarrow \NN$ that can be defined by composition or by simple recurrence from primitive recursive functions already defined (we start with the constant functions and addition~$+$).
Let us denote by $\Prim_2$ the set of primitive recursive functions $\NN^2\rightarrow \NN$. We easily prove that $\Prim_2$ is an enumerable set.
\\
A function $\beta\in\Prim_2$ can be thought of as simulating the execution of a program as follows. For an input $n$ we compute~\hbox{$\beta(n,m)$} for $m=0$, $1$, $\ldots$ until $\beta(n,m)\neq 0$ (intuitively: until the program reaches the instruction {\sf Halt}).
 Then, the function $\alpha\in\Rec$ computed by the \gui{program} $\beta\in\Prim_2$ is: $f:n\mapsto \beta(n,m_n)-1$ where $m_n$ is the first value of $m$ such that $\beta(n,m)\neq 0$.

Thus, we obtain a surjective map from a subset $Rec$ of $\Prim_2$ onto $\Rec$, and   $\Rec$ can be identified with the preset $Rec$ equipped with the suitable \egt and distinction. This means that  $\Rec$ is defined as a \gui{quotient}({\footnote{Since $\Rec$ is the image of  $Rec$ under a surjective map.}}) of a subset of an enumerable set. The \elts  of the subset $Rec$ of $\Prim_2$ are defined by the following condition:
$$\beta\in Rec \; \eqdef \;(*)\;:\;  \forall n\in \NN \;\; \exists m\in
\NN\;\;\; \beta(n,m)\neq 0. 
$$
From a classical point of view, for any $\beta\in\Prim_2$, the above assertion~$(*)$ is true or false in the absolute, in reference to the logic of the excluded middle (or, if you prefer, to the actual infinity of $\NN$): the notion of a mechanical computation can thus be defined without any reference to a primitive notion of construction.

However, from a \cof point of view, the assertion $(*)$ must be proven, and such a \dem is itself a construction. Thus \emph{the notion of a mechanical computation depends on the notion of construction, which cannot be defined}.

To finish this section, let us note that the Russian constructivism \`a la Markov admits as a fundamental principle the \egt  $\Rec=\NN^\NN$, a principle sometimes called the \textbf{false Church's thesis}\index{Church's thesis!False ---}. See \cite{Be,BR} and \cite[Richman]{ri2}. The true  \textbf{Church's thesis}\index{Church's thesis} is that no 
automated system of computation
 will ever be able to compute other functions than the recursive functions: we will be able to improve the performances of computers, but no 
automated system of computation
 will be able to surpass what they know how to compute \gui{in principle} (\cad if they dispose of the necessary time and space).
The true Church's thesis is extremely likely, but obviously it is unlikely to have a \demz.

\penalty-2500
\section{Principles of omniscience}
A \emph{principle of omniscience} is a principle that, although true in \clamaz, clearly poses a \pb in \comaz, because it a priori assumes knowledge of what happens with a potential infinity. The word omniscience here is therefore valid for \gui{prescience of the potential infinite.} The
principles of omniscience in \gnl have strong counterexamples 
in Russian \comaz. They however cannot be disproven
 in \coma \`a la Bishop, because they are compatible with \clamaz.
\subsubsec{The Little Principle of Omniscience} 

\noi Let $\alpha=(\alpha_n)\in\{0,1\}^\NN$ be a {\it binary sequence}, \cad
a construction that gives for each non-negative integer
 (as input) an \elt of $\{0,1\}$ (as output). Consider the following assertions
$$\arraycolsep2pt
\begin{array}{rcl}
P(\alpha) & :~ & \alpha_n=1 \mathrm{\;  for \; some\;  } n,\\[1mm]
\lnot P(\alpha) & : & \alpha_n=0 \mathrm{\; for\;  all\;  }n,\\[1mm]
P(\alpha)\vee \lnot P(\alpha) & : & P(\alpha) \mbox{ or }\lnot P(\alpha),\\[1mm]
\forall \alpha \;\; \big(P(\alpha)\vee \lnot P(\alpha)\big) & : & \mathrm{for\;
every \;binary \;sequence\; } \alpha, \; P(\alpha)\mathrm{\;or\; }\lnot P(\alpha).
\end{array}
$$

A \cov \dem of $P(\alpha)\vee \lnot P(\alpha)$ should provide an algorithm that either shows that $\alpha_n=0$ for all $n$, or computes a non-negative integer $n$ 
such that $\alpha_n=1$.

Such an \algo is much too efficient, because it would allow us to automatically solve a great number of important conjectures. In fact we know that if such an \algo exists, it is certainly not \gui{mechanically computable}: a program that runs on a machine can surely not accomplish such a thing even when we impose the limitation on the input $\alpha$ that it be an explicit primitive recursive binary sequence. This impossibility is a grand \tho of computability theory, often indicated under the name \gui{undecidability of the Halting Problem.}


\mni
{\bf Undecidability of the Halting problem} (We cannot know everything)\\
{\it In three immediately \eqves forms:
\begin{itemize}
\item  We cannot automatically assure the halting of  programs: there exists no program \,$T$\, that can test if an arbitrary program \,$P$\, will eventually reach its Halt instruction.
\item  There exists no  program that can test if an arbitrary primitive recur\-sive sequence is identically null.
\item There exists no program \,$U$\, that takes as input two integers, gives as output a Boolean, and that enumerates all the programmable binary sequences (the sequence \,$n\mapsto U(m,n)$\, is the
\,$m^{\rm th}$ sequence enumerated by \,$U$).
\end{itemize}
 }

\medskip
Not only does this \thoz, in its last formulation, resemble Cantor's \tho which asserts that we cannot enumerate the set of binary sequences, but the (very simple) \dem is essentially the same.

\smallskip
Although the previous \tho does not a priori forbid the existence of an  effective 
but not mechanizable procedure to systematically solve this type of \pbz, it confirms the intuitive idea according to which new ingenuity 
 will always have to be shown to progress in our knowledge of the \mathe world.

\smallskip Thus, from a \cof point of view, we reject the {\it Limited
Principle of Omniscience}.\index{omniscience!LPO@\LPO}\index{LPO@\LPO}
\begin{description}
\item \LPOz:  If $(\alpha_n)$ is a binary sequence, then either there exists some $n$
such that $\alpha_n=1$, or $\alpha_n=0$ for every $n$.
\end{description}
%
Here it is in a more concentrated form.
\begin{description}
\item \LPOz: \hfill $\forall \alpha \in\NN^\NN, \; \;
(\alpha\not= 0 \;\lor \;\alpha= 0)$\hfill~
\end{description}

We will call an \emph{\elr \prtz} a \prt \eqve to

\snic{\exists n\; \alpha(n)\not=0}

for a certain $\alpha\in\NN^\NN$.

\medskip The principle \LPO has several \eqv forms. Here are a few of them.
\begin{enumerate}
\item  If $A$ is an \emph{\elrz} \prtz, we have $ \; A \lor \lnot A $.
\item Every sequence in $\NN$ is either bounded, or unbounded.
\item Every decreasing sequence in $\NN$ is constant from a certain rank.
\item From a bounded sequence in $\NN$ we can extract a constant infinite subsequence.
\item Every \enum subset of $\NN$ is detachable.
\item Every \enum subset of $\NN$ is either finite, or infinite.
\item For every double sequence of integers $\beta :\NN^2\rightarrow \NN$ we have
$$\preskip.4em \postskip.0em \forall n\; \exists m\;\;  \beta(n,m)=0\quad \lor\quad
\exists n\; \forall m\; \; \beta(n,m)\neq 0$$
\item Every detachable subgroup of $\ZZ$ is generated by a single \eltz.
\item Every subgroup of $\ZZ^p$ generated by an infinite sequence is \tfz.
\item $\forall x\in\RR$, $(\; x\not= 0 \; \lor\;  x= 0\; )$.
\item $\forall x\in\RR$, $(\; x> 0 \; \lor\;  x= 0 \; \lor\;  x<0\; )$.
\item Every monotone bounded sequence in $\RR$  converges.
\item From a bounded sequence in $\RR$ we can extract a convergent subsequence.
\item Every real number is either rational or irrational.
\item Every \tf \sevc of \,$\RR^n$\, admits a basis.
\item Every separable Hilbert space admits \\
-- either a finite Hilbert basis \\
-- or a \denb Hilbert basis.
\end{enumerate}

\subsubsec{The Lesser Limited Principle of Omniscience}\rdb

\noi  Another, weaker, principle of omniscience \LLPO  (Lesser Limited Principle of Omniscience) is the following.\index{omniscience!LLPO@\LLPO}\index{LLPO@\LLPO}
\begin{description}\label{LLPO}
\item \LLPOz:  If $A$ and $B$  are two \elrs \prtsz, we have

\snic{\lnot (A\; \land\;  B) \quad \Longrightarrow \quad(\lnot
A\;  \lor\;  \lnot B)}
\end{description}
\rdb
This principle \LLPO  has several \eqv forms.
\begin{enumerate}
\item $\forall\, \alpha, \,\beta$ non-decreasing sequences $\in\NN^\NN$, if $\, \forall
n\,
     \alpha(n)\beta(n)=0$, then $\, \alpha=0$ or~\hbox{$\beta=0$}.
\item \label{recinsep}$\forall \alpha,\beta\in\NN^\NN$, if $\; \forall
n,m\in\NN\;  \alpha(n)\neq
      \beta(m)\;  $ then $ \; \exists \gamma\in\NN^\NN\;  $ such that
$$\preskip.4em \postskip.0em 
\forall n,m\in\NN\quad
    \big(\gamma(\alpha(n))=0\; \land\;  \gamma(\beta(m))=1\big). 
$$
\item $\forall \alpha\in\NN^\NN$, $\exists k\in \{0,1\}$, ($\; \exists n\;
     \alpha(n)=0\; \Rightarrow \; \exists m\;
     \alpha(2m+k)=0)$.
\item  $ \forall x \in \RR \quad (\; x\leq 0\;  \lor\;  x\geq 0\; )$
 (this allows us to make many \dems by dichotomy with the real numbers.)
\item $\forall x,y \in \RR \quad (\; xy=0\; \Rightarrow\;  (\; x= 0 \; \lor\; y=0\; )\; )$.
\item The image of an interval $[a,b]\subset \RR$ under a uniformly continuous real function is an interval~$[c,d]$.
\item A uniformly continuous real function over a compact metric space attains its least upper bound and its greatest lower bound.
\item \kl (one of the versions of K\"onig's lemma) 
Every explicit, infinite, finitely branching tree 
has an infinite path. \label{K1LLPO}
\end{enumerate}

It is known that if an \algo exists for the third item it cannot be \gui{mechanically computable} (\cad recursive): we can construct mechanically computable $\alpha$ and $\beta$ satisfying the 
hypothesis, but for which no mechanically computable $\gamma$ satisfies the conclusion.
Similarly, Kleene's singular tree is an infinite \denb recursive finitely branching tree that has no infinite recursive path.
 This gives a \gui{recursive counterexample} for \klz.

\smallskip We will now prove the \eqvc
\kl $\Leftrightarrow$ \LLPOz.~\footnote{As for all the \dems in this annex, it is informal and we do not specify in which formal framework it could be written. The readers will notice in this \dem a use of a construction by \recu which actually stems from the Axiom of Dependent Choice, \gnlt considered as non-problematic in \comaz.}

An explicit infinite finitely branching tree can be defined by a set $A\subset \Lst(\NN)$ of lists of integers satisfying the following \prts (the first four corresponding to the notion of an explicit finitely branching tree).
\begin{itemize}
\item  The empty list $[\,]$ represents the root of the tree, it belongs to $A$,
\item  an $a=[a_1,\ldots,a_n ]\in A$ represents both a node of the tree 
and the path that goes from the root to the node,
\item  if $[a_1,\ldots,a_n ]\in A$ and $n\geq 1$, then $[a_1,\ldots,a_{n- 1}]\in A,$
\item  if $a=[a_1,\ldots,a_n ]\in A$ the $x$'s $\in\NN$ such that $[a_1,\ldots,a_n,x ]\in A$ form a segment $\sotq{x\in\NN}{x<\mu(a)}$
where $\mu(a)$ is explicitly given in terms of $a$: the branches stemming from $a$ are numbered $0,\ldots,\mu(a)-1$.
\item  For all $n\in\NN$ there is at least one $[a_1,\ldots,a_n]\in A$
(the tree is explicitly infinite).
\end{itemize}
Thus the subset $A$ of $\Lst(\NN)$ is  detachable (this is ultimately what the word \gui{explicit} means here), and $A$ is \denbz.

%
\begin{Proof}{\Demo of \kl $\Leftrightarrow$ \LLPOz.}\\
We use the variant of \LLPO given in item 1.\\
Assume \klz. Let $\alpha,\beta\in\NN^\NN$ as in item 1.
Consider the following tree. 
	The root has two children. They form two distinct paths that grow indefinitely without ever branching out, 
 until~\hbox{$\alpha(n)\neq 0$} or $\beta(n)\neq 0$ (if this ever occurs). If this occurs with~\hbox{$\alpha(n)\neq 0$}, we stop the left branch and we continue the one on the right. If it occurs with~$\beta(n)=0$, we do the opposite. Explicitly giving an infinite branch in this tree amounts to certifying in advance that $\alpha=0$ or $\beta=0$.

Conversely, assume \LLPOz. Consider an explicit infinite finitely branching tree. Suppose \spdg that the tree is binary: beyond a node there are at most two branches. We prove by \recu that we can select up to depth $n$ a path that reaches a node $K_n$ underneath which the tree is infinite.
This is true for $n=0$ by hypothesis. If this is true for $n$, there is at least one branch underneath the selected node $K_n$. If there are two, consider the sequences $\alpha_n$ and $\beta_n\in\NN^\NN$ defined as follows\\
--- $\alpha_n(m)=0$ if there is at least one branch of length $m$ below~$K_n$ going to the right-hand side, otherwise $\alpha_n(m)=1$ \\
--- $\beta_n(m)=0$ if there is at least a branch of length $m$ below~$K_n$ going to the left-hand side, otherwise $\beta_n(m)=1$.\\
By \hdr the sequences $(\alpha_n)_{n\in\NN}$ and $(\beta_n)_{n\in\NN}$ are non-decreasing and their product is null.  We apply item~1 of \LLPOz: one of the two sequences is null and this gives us the means to select the path on the right or the left.
\end{Proof}

\subsubsec{The Law of Excluded Middle}

\noi  The Law of Excluded Middle (\TEMz)\index{omniscience!LEM@\TEM}\index{LEM@\TEM}
 states that $P\vee \lnot P$ is true for every proposition $P$. This extremely strong principle of omniscience implies \LPOz. It implicitly assumes that sets such as $\NN$ or $\NN^\NN$ or even significantly more complicated, 
are \emph{actual infinities}. It also implies that every set $X$ is
discrete  if we define $x\neq_Xy$ as meaning~\hbox{$\lnot (x=_Xy)$}.

\section[Problematic principles \ldots]{Problematic principles in \coma}
\label{para pripro}\index{omniscience!principles}\relax
By a \emph{problematic principle} we mean a principle that, although satisfied in practice if we do \coma in Bishop's style, is \cov unprovable. In \clamaz, these principles are known as true or known as false. 
\\
For example, in practice, if some $\alpha \in\NN^\NN$ is \cot well-defined, it can be computed by a program. \\
In other words, in practice, the \textbf{false Church's thesis},
which we can write in the form \fbox{$\Rec=\NN^\NN$}, is satisfied in \comaz. But it cannot be proven in the minimalist framework of \coma \`a la Bishop, which is compatible with \clamaz,
because the false Church's thesis is a principle that is false in \clamaz, in virtue of a cardinality argument. However, Russian \coma takes it as a fundamental axiom.

\ss Here we will (briefly) only examine two problematic principles, both true in \clamaz. \rdb

\subsubsec{Markov's Principle} \label{principeMarkov}

\noi \emph{Markov's Principle}, \MPz, is the following
$$ \forall x \in \RR \quad (\lnot x= 0 \Rightarrow  x\neq 0).
$$

Asserting  \MP  amounts to saying: for every binary sequence $\alpha$, if it is impossible for all its terms to be null, then it must have a nonzero term.

Or even: if $A$ is an \elr \prt then $\lnot\lnot A \Rightarrow A$.

The Russian  \cov school admits \MPz. Actually, for some $\alpha\in\NN^\NN$, it seems impossible to give  a constructive \dem  of $\lnot (\alpha=0)$ without finding some $n$ such that $\alpha(n)\neq 0$. Thus \MP is valid from a practical point of view in the constructivism \`a la Bishop.  Note that \LPO  clearly implies~\MPz.

\subsubsec{Principles of uniform continuity}

\noi  The principle of uniform continuity asserts that every pointwise continuous function
over a compact metric space is uniformly continuous. It is equivalent to the same assertion in a special case, which is itself very close to one of the classical forms of K\"onig's lemma. 
It is of particular interest to study the mutual relations between the following problematic principles, especially as they frequently appear in classical analysis.

\begin{description}
\item [\UCp]  Every pointwise continuous function $f:X\rightarrow Y$,
with
$X$ a compact metric space and $Y$ a metric space, is uniformly continuous.
\item [\UC]  Every pointwise continuous function $f:\{0,1\}^\NN\rightarrow \NN$ is uniformly continuous.
\item [\Mini]  Every uniformly continuous real function $>0$ over a compact metric space is bounded below by a real $>0$.
\item [\Minim]  Every uniformly continuous real function $>0$ over a compact interval $[a,b]$ is bounded below by a real  $>0$.
\item [\Minip]  Every continuous real function $>0$  over a compact metric space is bounded below by a real $>0$.
\item  [\FAN]  An explicit binary tree $A$ that has no infinite path
(\cad $\forall \alpha\in \{0,1\}^\NN\; \exists m\in \NN\;
\alpha |^m\notin A$) is finite.
\end{description}

In the formulation of \FANz, we see that this  principle is seemingly 
related to
\LLPO (see the last \eqv form \kl cited on \paref{K1LLPO}). Actually, we can show that it is a consequence of \LPOz. But this is not a principle of omniscience.
Besides, it does not imply \LLPOz. In \comaz, \LLPO is obviously false in practice, whereas \KL is satisfied in practice, because each time that we know how to \cot prove that a finitely branching tree has no infinite path, 
we also know how to prove that it is finite.

\penalty-2500
\Exercices{

\Oui{\setcounter{exercise}{0}}
\vspace{-1em}
\begin{exercise}
\label{exoAnnexe1}
{\rm  Give \dems for  examples
 \iref{exo neq in N^N}, \iref{exo inf deno1}, \iref{exo log0}, \iref{exo log1}  
 and~\iref{exo log2}.
 }
\end{exercise}

\vspace{-1em}
\begin{exercise}
\label{exo inf deno2} 
{\rm
Explain why the notions of a finite set, a finitely enumerable set,  a bounded set,  a weakly finite set, and an enumerable bounded set cannot be identified in \comaz. 
Explain why these notions coincide if we admit \TEMz.
}
\end{exercise}

\vspace{-1em}
\begin{exercise}
\label{exoA22}
\emph{Prove a few of the \eqvcs mentioned for \LPOz.
}
\end{exercise}

\vspace{-1em}
\begin{exercise}
\label{exoA23}
\emph{Prove a few of the \eqvcs mentioned for \LLPOz.
}
\end{exercise}

}


%

\vspace{5pt}
\Biblio
\vspace{3pt}

The controversy on the nature and the use of the infinite in \maths was very strong at the beginning of the $20^{\rm th}$ century:
see for example Hilbert \cite[1926]{Hi2}\ihiz, Poincar\'e \cite[1909]{Po}, H. Weyl \cite[1918]{W1}, 
 \cite[1951]{Bro} and \cite[1987]{INS}).
The debate seems at first to have ended in favor of the point of view represented by classical logic.
Actually, since the 60s and especially since the publication of Bishop's book, the two points of view are considerably less contradictory than when they first appeared.
\\
 A few interesting references on this theme: \cite[1962]{Loren}, \cite[Fred Richman, 1990]{Ri90}, \cite[2007]{Dow2} and \cite[Per Martin-L\"of, 2008]{ML2008}.

Constructive logic is often called \gui{intuitionistic logic.}
It was developed as a formal \sys by A.\ Heyting.

The article \cite[Lorenzen, 1951]{Lor1951}  which informally uses the constructive point of view in the study of the formal system Principia Mathematica by Whitehead and Russell and makes the connection with the purely algebraic theory of distributive lattices deserves a thorough study.

There are pleasant presentations of such formal systems in the books \cite[1962]{Loren} and \cite[2001]{DNR}. 

The small book \cite[1995]{Dow} also gives an interesting informal presentation.

Concerning the discussion on the links between effectiveness 
and recursiveness, see \cite[Coquand]{CoqRec}, \cite[Heyting]{HeyRec} and \cite[Skolem]{SkoRec}.

The book \cite[1985]{Be} carries out a systematic study of several problematic principles in \comaz.
For Kleene's singular tree, see \cite[page 68]{Be} and \cite[1965]{KV}.

The development and the comparison of formal \syss able to serve as frameworks for the \coma employed in \cite{B67} or \cite{MRR} has been a very active research subject for a long time. 
We make sure to note the preponderant influence of the \cov theory of the types \tsbf{CTT} of Per Martin-L\"of,   
\cite{ML1973,ML} and \cite[1984]{PML}, 
and of the theory \tsbf{CZF} of Peter Aczel and Michael Rathjen (\cite[Aczel]{Acz} and \cite{AcRa}).
See also the recent developments in \cite[2014]{HoTT} and Thierry Coquand's webpage: \url{http://www.cse.chalmers.se/~coquand/}.

Let us also cite the beautiful book \cite[1998]{fef} which is inline with the propositions of Hermann Weyl.

For a discussion of the \gui{Fan Theorem}, see~\cite[Coquand]{CoqFT}.

The systematic study of the comparison  (in \comaz) of principles  of omniscience (such as \LPO or \LLPOz), as well as that of problematic principles (such as \MP or \FANz),  has recently been the subject of a major boom. On this subject, we can refer to \cite[Berger\&al.]{Berg,BeBr,BeIs} and \cite[Ishihara]{Ishi,Ishi2,Ishi3}.

\newpage \thispagestyle{CMcadreseul}
\incrementeexosetprob


\pagestyle{CMExercicesheadings}


\cleardoublepage \thispagestyle{CMcadreseul}
\refstepcounter{chapter}
\addcontentsline{toc}{chapterbis}{Tables of theorems}
\addtocontents{toc}{\vskip-0.8em}

\catcode`\@=11
\@makesindexhead{Tables of theorems}
\perso{compilé le \today}
\pagestyle{CMExercicesheadings}

\setlength\parskip{1pt}

\markboth{Tables of \thosz}{Tables of \thosz}


\newcommand \TTT[3]{%
\noindent 
\mbox{%
\parbox[b]{.80\textwidth}{\leftskip10pt\parindent-\leftskip\strut#3\dotfill\par}%
\hspace{.01\textwidth}%
\parbox[b]{.12\textwidth}{\hfill#1}%
\hspace{.01\textwidth}%
\parbox[b]{.06\textwidth}{\hfill#2}%
}%
}

\newcommand \TTref[2]{\TTT{\ref{#1}}{\pageref{#1}}{#2}\par}

\newcommand \TTpa[3]{\TTT{\ref{#3}}{\pageref{#1}}{#2}\par}

\newcommand \chap[1]{

\goodbreak
\mni\textbf{#1}}

\small
\begin{center}
    {\large\bf Dynamic methods}\label{machineries}
\end{center}

\vspace{2mm}
\TTT{~}{page}{Name} \par\vspace{1mm}

\TTpa{MethodeQI}{\Elr \lgb machinery of \qirisz}{subsecAnneauxqi}
\TTpa{MethodeZedRed}{\Elr \lgb machinery of reduced \zed \risz}{secKrull0dim}
\TTpa{subsecDyna}{The dynamic method}{subsecDyna}
\TTpa{MetgenAnar}{\Lgb machinery of \anars}{secIplatTf}
\TTpa{MethodeIdeps}{Basic \lgb machinery (with \idepsz)}{secMachLoGlo}
\TTpa{MethodeIdemax}{Dynamic machinery with \idemasz}{subsecLGIdeMax}
\TTpa{MethodeIdemin}{Dynamic machinery with \idemisz}{subsecLGIdepMin}
\TTpa{MachDynAxAx}{Dynamic machinery with $\ArX$ and $\gA(X)$}{sec.Etendus.Valuation}

\vspace{5mm}

\begin{center}
    {\large\bf Concrete local-global principles}
\end{center}
\vspace{2mm}
\TTref{plcc.basic}{Basic concrete local-global principle}
\TTref{PrTransfertBasic}{Basic transfer principle}
\TTref{plcc.coh}{\Coh modules}
\TTref{plcc.tf}{\Tf modules}
\TTref{plccRangMat}{Rank of a matrix}
\TTref{fact.lnl.loc}{\Lnls \lin maps} 
\TTref{plcc.basic.modules}{Exact sequences of modules} 
\TTref{plcc.basic.monoides}{For \mosz} 
\TTref{plcc.entier}{Integral \elts over a \riz}
\TTref{plcc.scinde}{\Prts of \alis between \mpfsz}
\TTref{plcc.pf}{\Fp modules}
\TTref{plcc.aqi}{\Qiri \ris (\qirisz)}
\TTref{plcc.cor.pf.ptf}{\Mptfsz}
\TTref{plcc.aGs}{Galois \algsz}
\TTref{plcc.plat}{Flat modules}
\TTref{plcc.arith}{\lsdsz, \ari \risz, Pr\"ufer \risz, \lops \idsz} 
\TTref{plcc.algfptes}{Flat or \fptes \algsz, \lon at the source} 
\TTref{thlgb1}{Similar or equivalent matrices (\lgb \riz)}
\TTref{thlgb2}{Isomorphic \pf modules (\lgb \riz)}
\TTref{thlgb3}{Quotient modules (\lgb \riz)}
\TTref{plcc.normal}{Normal \ris and \icl \idsz}
\TTref{plcc.agq}{Primitively \agq \eltsz}
\TTref{plcc.ddk}{Dedekind \risz}
\TTref{thDdkLoc}{Krull dimension of \risz}
\TTref{plcc.KdimMor}{Krull dimension of morphisms}\rdb
\TTref{plcc.sli}{Exact sequences and generalizations}
\TTref{plcc.ptf}{Finiteness \prts for modules}
\TTref{plcc.propaco}{\Prts of commutative \risz}
\TTref{plcc.apf}{Finiteness \prts of \algsz, \lon at the source} 
\TTref{plcc2.apf}{Finiteness \prts of \algsz, \lon at the sink}
\TTref{plcc.modules 1}{Concrete patching of \elts of a module, or of \homos between modules} \TTref{plcc.modules 2}{Concrete patching of modules}
\TTref{plcc.RecolHomAnn}{Concrete patching of \ri \homosz}
\TTref{plcc.regularite}{Concrete local-global principle for equality in depth 1}
\TTref{plcc0Profondeur2}{Concrete local-global principle for divisibility in depth 2}
\TTref{plcc.modules1bis}{Concrete patching of \elts in a module in \prof 2}
\TTref{plcc.modules2bis}{Concrete patching of modules in \prof 2}
\TTref{thPatchV}{Vaserstein \rcmz: equivalent matrices over $\gA[X]$}
\TTref{thPatchQ}{Quillen \rcmz: extended modules}
\TTref{plgcetendus}{Local-global principle \`a la Roitman}
\TTref{plcc.159}{Concrete \rcm in the \elr group}
\TTref{plgc-Rao}{Rao's concrete \lgb principle}
\label{plgvarios}

\vspace{5mm}
%


\begin{center}
    {\large\bf Closed covering principles}
\end{center}

\vspace{2mm}

\TTref{prcfgrl}{For \grlsz}
\TTref{prcf1}{Nilpotent, \com \eltsz}
\TTref{prcf2}{\Tf modules}
\TTref{prcf3}{Rank of a matrix,  \mptfsz}
\TTref{thDdkRecFer}{Krull dimension}
%

\vspace{5mm}
%

\begin{center}
    {\large\bf Stability under \eds}
\end{center}

\vskip2mm

\TTref{propPfExt}{\Tf and \pf modules, tensor products, \smqs and exterior powers, exterior \algz}
\TTref{fact.idf.change}{Fitting \idsz}
\TTref{propPtfExt}{\Mptfsz}
\TTref{fact.det loc}{Determinant, \polcarz, \polfonz, \polmuz,
cotransposed \endoz} 
\TTref{factEdsAlg}{\Tfz, \pfz, \stfes \algsz}
\TTref{factEdsDualisante}{Dualizing forms, Frobenius \algsz}
\TTref{factEdsEtale}{\Ste \algsz}
\TTref{factSpbEds}{\Spb \algsz}
\TTref{factEdsAG}{Separating automorphisms}
\TTref{corAGextsca}{Galois \algsz}
\TTref{factChangeBase}{\Uvl \dcn \algz}
\TTref{corPlatEds}{Flat modules}
\TTref{propFidPlatPrAlg}{Converses in the \fptes extensions case}

\vspace{5mm}
%
\begin{center}
\bigskip    {\large\bf Theorems}
\end{center}


\CHAP{\nameref{chapSli}}

\TTref{plcc.basic}{Basic concrete \lgb principle}
\TTref{lemGaussJoyal}{Gauss-Joyal lemma}
\TTref{propCoh4}{Characterization of \comosz}
\TTref{fact.sfio}{Fundamental \sys of \ort \idmsz}
\TTref{lem.ide.idem}{Lemma of the \tf \idm \idz}
\TTref{restes chinois}{Chinese remainder \thoz, \gnl form (for the \ari form see \thref{thAnar})}
\TTref{lem.min.inv}{Lemma of the \iv minor}
\TTref{lem pf libre}{Freeness lemma}
\TTref{lemCram}{Generalized Cramer formula}
\TTref{propIGCram}{Magic formula \`a la Cramer} 
\TTref{propFactDirRangk}{\Smtfs as direct summands of a free module} 
\TTref{prop inj surj det}{Injectivity and surjectivity criteria} 
\TTref{theoremIFD}{\Lnls matrices}
\TTref{Th.transitivity}{Transitivity formulas for the trace and the \deterz}
\TTref{thTransDisc}{Transitivity formula for the \discrisz}

\CHAP{\nameref{chapGenerique}}

\TTref{thSymEl}{\Elr \smq \polsz}
\TTref{lemdArtin}{\DKM lemma}
\TTref{thKro}{\KROz's \tho (1)}
\TTref{propUnicCDR}{Uniqueness of the \cdr  (\stf case)}
\TTref{thPrIs}{Isomorphism extension \thoz}
\TTref{thGaloiselr}{Galois correspondence}
\TTref{thResolUniv}{Construction of a splitting field}
\TTref{LemElimAffBasic}{Basic \eli lemma}
\TTref{thIntClosStab}{Integrally closed \pol \riz}
\TTref{thClotAlgQ}{Splitting field, primitive \elt \thoz}
\TTref{th1IdZalpha}{Every nonzero \itf of a number field is \ivz}
\TTref{th2IdZalpha}{Multiplicative structure of \itfs of a number field}
\TTref{propEvitementConducteur}{Dedekind's \thoz, \ids that avoid the conductor}
\TTref{thNstfaibleClass}{Weak \nst and \Noe position, see also \tho \ref{thNstNoe}}
\TTref{thNstClass}{Classical \nstz}
\TTref{thNSTsurZ}{\nst over $\ZZ$, formal \nstz}
\TTref{corthNSTsurZ}{\nst over $\ZZ$, formal \nstz, 2}
\TTref{thNewtonQuad}{Newton's method}
\TTref{corIdmNewton}{Residual \idms lemma}

\CHAP{\nameref{chap mpf}}

\TTref{lem pres equiv}{Matrices that present the same module}
\TTref{propsregpf}{An \id generated by a regular sequence is \pfz}
\TTref{thidptva}{The \id of a point is a \mpfz}
\TTref{propCoh2}{Coherence and finite \pn (see also Proposition \ref{propAliCoh})}
\TTref{prop unicyc}{Direct sum of cyclic modules (uniqueness)}
\TTref{prop quot non iso}{An isomorphic quotient is a quotient by $0$}
\TTref{thScindageQi}{Quasi integral splitting lemma.}
\TTref{thScindageZed}{\Zed splitting lemma}
\TTref{thZerDimRedLib}{The reduced zero-dimensional ring paradise}
\TTref{thSPolZed}  {\Zed system over a \cdiz}
\TTref{thStickelberger}  {Stickelberger's \tho (\zed \sysz)}
\TTref{fact.idf.ann}{Annihilator and first \idfz}
\TTref{LemElimAff}{General elimination lemma}
\TTref{thElimAff} {\Agq \eli \thoz: the resultant \idz}

\CHAP{\nameref{chap ptf0}}

\TTref{propdef ptf}{\Mptfsz}
\TTref{prop pf ptf}{\Mpn of a \mptfz}
\TTref{corlemScha}{Schanuel's lemma}
\TTref{propIsoIm}{\Dlg lemma}
\TTref{propZerdimLib}{\Zed freeness lemma: item \emph{\iref{ite2propZerdimLib}} of the \thoz}
\TTref{propZerdimLib}{Incomplete basis \thoz: item \emph{\iref{ite5propZerdimLib}} of the \thoz}
\TTref{corBass2}{Bass' \thoz, stably free modules, with $\Bdim$}
\TTref{prop Fitt ptf 2}{Local structure \tho for \mptfsz. 
See also \Thos \ref{theoremIFD},
\ref{corth.ptf.sfio}, \ref{th ptf loc free} and~\ref{th matproj}}
\TTref{factLocCas}{Successive \lons lemma, 1}
\TTref{propmlm}{\Lcy \tf modules, see also \ref{pmlm}}
\TTref{propdef det ptf}{Determinant of an \endo of a \mptfz}
\vspace{.15em}
\TTref{th ptf sfio} {The \sfio associated with a \mptfz}
\TTref{lem calculs}{Explicit computations: \deterz, \polcarz, etc}
\TTref{th decomp ptf}{Decomposition of a \mptf into a direct sum of modules of constant rank}

\CHAP{\nameref{chap AlgStricFi}}

\TTref{th1Etale}{Structure \tho for étale \Klgsz, 1}
\TTref{cor1lemEtaleEtage}{Separable \elts in a \Klgz}
\TTref{corlemEtaleEtage}{\Carn of  étale \Klgsz}
\TTref{thEtalePrimitif}{Primitive \elt \thoz}
\TTref{th2Etale}{Structure \tho for étale \Klgsz, 2}
\TTref{thClsep}{Separable closure}
\TTref{thNormalEtaleGalois}{\Carn of  Galois extensions}
\TTref{thCGSynthese}{Galois correspondence, synthesis}
\TTref{factSDIRKlg}{Direct sum in the category of  \klgsz}
\TTref{lemLingOver}{Lying Over: see also Lemma~\vref{lemLingOver2}}
\TTref{thNst0}{A weak \nstz}
\TTref{propAlgFinPresfin} {When a \klg is a \pf \kmoz}
\TTref{propTraptf}{Transitivity for  \stfes \algsz}
\TTref{factCarDua}{\Carn of the dualizing forms in the \stf case}
\TTref{factCarAste}{\Carn of  \asesz}
\TTref{prop2EtaleReduit}{If $\gk$ is reduced, so is any \ste \klgz}
\TTref{thIdmEtale}{Idempotents and \edsz,  \asesz}
\TTref{thAlgSteIdmSpb}{\Spt idempotent of a \asez}
\TTref{thSepIversen}{\Cara \prts of  \spb \klgsz}
\TTref{thAlgStfSpbSte}{A \stfe \spb \alg is \stez}
\TTref{thSepProjFi}{Finiteness \prt of  \spb \algsz}
\TTref{corthSepProjFicdi}{Over a \cdiz, a \pf \spb \alg is \stez.}
\TTref{thDA}{Dedekind's lemma}
\TTref{thA}{Artin's \thoz, \aGs version}
\TTref{corAGextsca}{\Eds for  \aGsz}
\TTref{thAGACar}{\Carns of  \aGsz}
\TTref{corAGAfreeCar}{\Carns of  free \aGsz}
\TTref{corDAexplicite2}{Galois correspondence, \aGs version}
\TTref{thCorGalGen}{Galois correspondence, connected \aGsz}
\TTref{corAGQuo}{Galois quotient of a \aGz}
\TTref{exoLuroth1}{L\"uroth's \tho (exercise)}

\CHAP{\nameref{ChapGalois}}

\TTref{thNstfaibleClassSCA}{Weak \nst and \Noe position, 2}
\TTref{thNstNoe}{Weak \nst and \Noe position, 3}
\TTref{thNoetSimult}{Simultaneous \Noe position}
\TTref{thNstClassCof1}{Classical \nstz, \gnl \cov version}
\TTref{thNstClassCof2}{\nst with multiplicities}
\TTref{thpolcohfd}{A \apf over a \cdi is a \fdi \coriz}
\TTref{corpropBoolFini}{Structure \tho for finite \agBsz}
\TTref{lemAduIdmA}{Galois structure \tho (1), $G$-Boolean \algsz}
\TTref{thADG1Idm}{Galois structure \tho (2), Galois quotients of \apGsz}
\TTref{lemAdu1}{\Uvl splitting \alg and \sptz. See also \thref{theoremAdu2}}
\TTref{theoremAdu1}{\Uvl splitting \alg and fixed points}
\TTref{thAduAGB}{The \adu as a \aGz}
\TTref{theoremAdu3}{Diagonalization of a \aduz, see also \thref{theoremAdu4}}
\TTref{thidGTri}{Triangular basis of the \id defining a \aGz}
\TTref{thUnici}{Eventual uniqueness of the \cdr of a \spl \polz}
\TTref{thdivzeridm}{Dynamic management of a \cdrz, see also \thref{corRvRel1}}
\TTref{propUnicite}{Uniqueness of the \cdrz, dynamic version}
\TTref{thStruc3}{Galois structure \tho (3), Galois quotients of the \adu of a \spl \pol over a \cdiz}
\TTref{thZsuffit}{Where the computations take place: the sub\ri $\gZ$ of $\gK$ is sufficient}
\TTref{exothNst1-zed}{\nst and Noether position, reduced \zed \ris case (exercise)}

\CHAP{\nameref{chap mod plats}}

\TTref{thPlat1}{\Carn of  flat modules, 1}
\TTref{cor pf plat ptf}{\Carn of  \mptfs by the flatness}
\TTref{thplatTens}{\Carn of  flat modules, 2}
\TTref{propPlatQuotientdePlat}{Flat quotients}
\TTref{thExtPlat}{\Carn of  flat \algsz}
\TTref{thExtFidPlat}{\Carn of  \fptes \algsz}
\TTref{thSurZedFidPlat}{Every extension of a \cdi is \fptez}
\TTref{propFidPlatTf}{Faithfully flat extensions and finiteness \prts of  modules}
\TTref{propFidPlatPrAlg}{Faithfully flat extensions and finiteness \prts of  \algsz}

\CHAP{\nameref{chap Anneaux locaux}}

\TTref{thJacUnitEntieres}{Jacobson radical and units of an integral extension}
\TTref{thJacplc}{Local \prts of integral extensions}
\TTref{lemNaka}{Nakayama's lemma (the determinant trick)}
\TTref{lelilo}{Local freeness lemma}
\TTref{lelnllo}{Lemma of the \lnl map}  
\TTref{lemnbgtrlo}{Local number of \gtrs lemma}
\TTref{lemLocaliseFini}{Localized finite \ri lemma}
\TTref{lemLocalisezeddim}{Localized zero-dimensional ring lemma}
\TTref{thCotangent}{Cotangent space at $\uxi$ and $\fm_\uxi / {{\fm_\uxi}^2}$}
\TTref{thJZS}{A simple isolated zero}
\TTref{thJZScdi}{Isolated zero}
\TTref{thPointLisseCourbeIC}{The \id of a non-singular point of a locally complete intersection curve. See also \thref{th2PtlisseCourbeIC}}
\TTref{thLgbExtEnt}{Integral extension of a \algbz}

\CHAP{\nameref{chap ptf1}}

\TTref{th rg const loc free}{The modules of constant rank are \lot free}
\TTref{th ptf loc free}{The \ptf modules are \lot free;
see also \Thos \ref{thptfloli} and \vref{th matproj}}
\TTref{propRgConstant2}{Modules of constant rank $k$ as submodules of $\Ae k$}
\TTref{corthTransPtf}{\Stf \Algsz: transitivity formula for the ranks}
\TTref{thBnk}{The functor $\Gnk$}
\TTref{leli2}{Second freeness lemma}  
\TTref{prop1TanGrassmann}{Tangent space to a Grassmannian, see also \thref{prop3TanGrassmann}}
\TTref{thBézoutLib}{Every \mrc over a Bézout \qiri is free}
\TTref{thPicKTO}{$\Pic \gA$ and $\KTO \gA$}
\TTref{propRgConstant3}{Picard group and group of classes}
\TTref{propComparRedRed}{$\GKO(\gA)\simeq\GKO(\Ared)$}
\TTref{thCarreMil1}{Milnor square}
\TTref{thIdenPtsSimples2}{A complete classification of $\GKO(\gA)$;
see also \thref{thIdenPtsSimples}}

\CHAP{\nameref{chapTrdi}}
\nobreak

\TTref{thBoolGen}{Boolean \alg freely generated by a \trdiz}
\TTref{th0GpRtcl}{Distributivity in the \grlsz}
\TTref{th1GpRtcl}{Riesz \tho (\grlsz)}
\TTref{th2GpRtcl}{Partial \dcn \tho under Noetherian condition}
\TTref{thAgcdNormal}{A GCD-domain is \iclz.}
\TTref{propGCDDim1}{A GCD-domain of dimension $\leq 1$ is a Bézout \riz}
\TTref{thAXgcd}{GCD-domains: $\gA$ and $\AX$}
\TTref{thZedGen}{\Zedr closure of a commutative \riz}
\TTref{thEntRel1}{Fundamental \tho of the \entrelsz}
\TTref{thDualiteFinie}{Duality \tho between \trdis and finite ordered sets}

\CHAP{\nameref{ChapAdpc}}\label{tdtChapAdpc}

\TTref{thAnar}{\Carns of  \anarsz, Chinese remainder \thoz}
\TTref{thiivanar}{Multiplicative structure of the \iv \ids in an \anarz}
\TTref{thPruf}{\Carns of  \adpsz}
\TTref{thExtEntPruf}{Normal integral extension of a \adpz}
\TTref{thSurAdp}{Overring of a \adpz}
\TTref{th.adpcoh}{\Carns of \adpcsz}
\TTref{ThImMat}{\Fp modules over a \adpcz}
\TTref{th.2adpcoh}{Another \carn of  \adpcsz}
\TTref{factAdpIntExt}{Finite extension of a \adpc (see also \thref{thAESTE})}
\TTref{thK1-SLE}{$\SL_3=\EE_3$ for a \qiri of dimension $\leq 1$}
\TTref{th1-5}{One and a half \thoz: \qiris of dimension~$\leq1$}
\TTref{dekinbe}{A Bézout \adpcz}
\TTref{thcohdim1}{A normal, \coh \ri of dimension $\leq1$ is a \adpz}
\TTref{thPTFDed}{\Pro modules of rank $k$ over a \ddp of dimension~\hbox{$\leq1$}}
\TTref{thMpfPruCohDim}{\Tho of the invariant factors:  \mpfs  over a \ddp of dimen\-sion~$\leq1$}
\TTref{mat33}{Reduction of a row matrix}
\TTref{thRieszAnar}{Riesz \tho for  \anarsz}
\TTref{thFactor}{Factorization of \itfs over a \adpc of dimension~$\leq1$
(see also \thref{thFactor2})}
\TTref{prop2DDK}{A \adk admits  \fapsz}
\TTref{propDDK}{\Carns of  \adksz}
\TTref{thdefDDKTOT}{Total factorization Dedekind rings}
\TTref{lemthAESTE}{A computation of integral closure (\adksz)}
%

\CHAP{\nameref{chapKrulldim}}

\TTref{thZarZar}{The duality between Zariski spectrum and Zariski lattice}
\TTref{thDKA}{\Elr \carn of the \ddkz}
\TTref{th1.5}{One and a half \tho (variant)}
\TTref{propDimKXY}{Krull dimension of a \pol \ri over a field}
\TTref{thDKAG}{Krull dimension and \Noe position}
\TTref{thAmin}{Minimal \qiri closure of a \riz}
\TTref{thKdimMor}{Krull dimension of a morphism}
\TTref{corthKdimMor}{Krull dimension of a \pol \riz}
\TTref{cor2thKdimMor}{Krull dimension of an integral extension}
\TTref{thKdimTDTO}{Krull dimension of a totally ordered set}
\TTref{th1Valdim}{Krull dimension of an extension of \advsz}
\TTref{thValDim}{Valuative dimension of a \pol \riz}
\TTref{corthValDim}{Krull dimension and \anarsz}
\TTref{propGupLY}{Going Up, Going Down and \ddkz}

\CHAP{\nameref{chapNbGtrs}}

\TTref{thKroH}{Non-Noetherian Kronecker-Heitmann theorem, with the Krull dimension}
\TTref{Bass0}{Bass' \gui{stable range} \thoz, with the Krull dimension, without \noetz}
\TTref{thKroLoc}{Kronecker's \thoz, local version}
\TTref{Bass}{Bass' \thoz, with the Heitmann dimension, without \noetz}
\TTref{KroH2}{Kronecker's \thoz, Heitmann dimension}
\TTref{thSerre}{\SSO \thoz, with $\Sdim$}
\TTref{thSwan}{Forster-Swan \thoz, with $\Gdim$}
\TTref{thSwan2}{\Gnl Forster-Swan \thoz, with $\Gdim$}
\TTref{thBassCancel2}{Bass' cancellation \thoz, with $\Gdim$}
\TTref{th2KroH}{Kronecker's \thoz, for the \sutsz}
\TTref{thPartitionSpec}{Constructible partition of the Zariski spectrum and $k$-stability}
\TTref{matrix}{Coquand's \thoz, 1: Forster-Swan and others with the $n$-stability}
\TTref{matrixC}{Coquand's \thoz, 2: elementary column operations and $n$-sta\-bi\-lity}
\TTref{MAINCOR}{Coquand's \thoz, 3: Forster-Swan and others
with the Heitmann dimension}

\CHAP{\nameref{chapPlg}}

Dynamic machineries and various \plgs are indicated pages \pageref{machineries}
and \pageref{plgvarios}.


\CHAP{\nameref{ChapMPEtendus}}

\TTref{thTSC}{Traverso-Swan-Coquand \thoz}
\TTref{thRoitman}{Roitman's \thoz}
\TTref{thHor0}{Local Horrocks'  \thoz}
\TTref{thHor}{Affine Horrocks'  \thoz}
\TTref{th2Bass}{Bass' \thoz}
\TTref{thStabLibPol}{Concrete Quillen induction, stably free case}
\TTref{InductionQabs}{Abstract Quillen induction}
\TTref{thQUILIND}{Concrete Quillen induction}
\TTref{corthQUILIND}{Quillen-Suslin \thoz, Quillen's \demz}
\TTref{th2QUILIND}{Concrete Quillen induction, free case}
\TTref{thSuslinQS}{Suslin's \thoz}
\TTref{th4SusQS}{(Another) Suslin's \thoz}
\TTref{th2HorrocksLocal}{Little Horrocks' theorem \`a la Vaserstein
 (and \thref{th2HorrocksGlobal})}
\TTref{th1Rao}{Rao's \thoz}
\TTref{thBass.Valuation}{(Another) Bass' \thoz}
\TTref{thVX2stab}{The $2$-stability of $\VX$, for \thref{thBass.Valuation}}
\TTref{thBassValu}{Bass-Simis-Vasconcelos \thoz (and \thref{thBassAri})}
\TTref{compa}{Dynamic comparison of $\gA(X)$ with $\ArX$}
\TTref{thMaBrCo}{Maroscia \& Brewer-Costa \thoz}
\TTref{LSindabs}{Abstract Lequain-Simis induction}
\TTref{induYengui}{Yengui induction}
\TTref{thLSValu}{Lequain-Simis \thoz}
\newpage \thispagestyle{CMcadreseul} 

\refstepcounter{chapter}

\newpage \thispagestyle{CMcadreseul}


\small

\cleardoublepage \thispagestyle{CMcadreseul}

\catcode`\@=11
\@makesindexhead{Index of notations}
\refstepcounter{chapter}
\addtocontents{toc}{\vskip-0.8em}
\addcontentsline{toc}{chapterbis}{Index of notations}
\markboth{Index of notations}{Index of notations}
\perso{compilé le \today}
\pagestyle{CMExercicesheadings}

\newbox\toto
\newbox\tata
\newlength\largeurtoto
\newlength\largeurtata

\newcommand \ttt[3]{\setbox\tata=\hbox{#1}%
\ifdim\wd\tata<.15\textwidth\relax%
\setlength{\largeurtoto}{.78\textwidth}%
\setlength{\largeurtata}{.15\textwidth}%
\else%
\setlength{\largeurtoto}{.93\textwidth}%
\addtolength{\largeurtoto}{-\wd\tata}%
\setlength{\largeurtata}{\wd\tata}%
\fi%
\setbox\toto=\hbox{\parbox[b]{\largeurtoto}{\leftskip10pt\parindent-\leftskip\strut#2\dotfill\par}}%
\smallskip\noindent\mbox{%
\parbox[b][\ht\toto][t]{\largeurtata}{#1}%
\parbox[b]{\largeurtoto}{\leftskip10pt\parindent-\leftskip\strut#2\dotfill\par}%
\hspace{.02\textwidth}%
\parbox[b]{.05\textwidth}{\hfill#3}%
}%
\par}

\newcommand
\NOTA[3]{\ttt{$#1\phantom{_{S}^{1}}$}{#2}{\pageref{#3}}

\vspace{-3pt}}

\begin{flushright}{page}\end{flushright}%

\CHAP{\nameref{chapMotivation}}
\NOTA{\Der{\RR}{\gB}{M}}{the \Bmo of the \dvns of $\gB$ in $M$}{Leibniz2}
\NOTA{\mathrm{Der}(\gB)}{the \Bmo of the \dvns of $\gB$}{Leibniz2}
\NOTA{\Om{\RR}{\gB}}{the \Bmo of the  \diles (of K\"ahler) of~$\gB$, see also \paref{thDerivUniv}}{dilesK}

\CHAP{\nameref{chapSli}}
\NOTA{\pi_{\gA,\fa}}{the canonical \homo $\gA\to\gA\sur{\fa}$}{NOTAgAfa}
\NOTA{\Ati}{the multiplicative group of  \iv \elts of $\gA$}{NOTAAst}
\NOTA{\gA_{S}}{(or $S^{-1}\gA$) the localized \ri of $\gA$ at $S$}{NOTAAst}
\NOTA{j_{\gA,S}}{the canonical \homo $\gA\to\gA_{S}$}{NOTAAst}
\NOTA{\sat{S}}{the saturated \mo of the \mo $S$}{NotaSatmon}
\NOTA{\gA[1/s]}{(or $\gA_{s}$) the localized \ri of $\gA$ at $s^\NN$}{NOTAA[1/s]}
\NOTA{(\fb:\fa)_\gA}{the conductor of the \id $\fa$ into the \id $\fb$}%
{NOTATransp}
\NOTA{(P:N)_\gA}{the conductor of the module $N$ into the module $P$}%
{NOTATransp}
\NOTA{\Ann_\gA(x)}{the annihilator of the \elt $x$}{NOTAAnn}
\NOTA{\Ann_\gA(M)}{the annihilator of the module $M$}{NOTAAnn}
\NOTA{(N:\fa)_M}{$\sotq{x\in M}{\fa x\subseteq N}$}{NOTAAnn2}
\NOTA{{(N:{\fa^\infty})}_M}{$\sotq{x\in M}
{\exists n\,\, \fa^n x\subseteq N}$}{NOTAAnn2}
\NOTA{\Reg \gA }{\mo of the \ndzs \elts of $\gA$}{NOTATotFrac}
\NOTA{\Frac\gA}{total \ri of fractions of $\gA$}{NOTATotFrac}
\NOTA{\Ae{m\times p}}{(or $\MM_{m,p}(\gA)$) matrices with $m$ rows and $p$ columns}{NOTAmatrices}
\NOTA{\Mn(\gA)}{$\MM_{n,n}(\gA)$}{NOTAmatrices}
\NOTA{\GLn(\gA)}{group of \ivs matrices}{NOTAmatrices}
\NOTA{\SLn(\gA)}{group of matrices with \deter $1$}{NOTAmatrices}
\NOTA{\GAn(\gA)}{\mprnsz}{NOTAmatrices}
\NOTA{\DA(\fa)}{(or $\sqrt{\fa}$) nilradical of the \id $\fa$ of $\gA$}{NOTADA}
\NOTA{\gA\red}{$\gA\sur{\DA({0})}$: reduced \ri associated with $\gA$}{NOTADA}
\NOTA{\rc_{\gA,\uX}(f)}{(or $\rc(f)$) \id of $\gA$, content of the \pol $f\in\AuX$}{NOTADA}
\NOTA{\rg_\gA(M)}{rank of a free module, see also the \gnns to the \mptfs pages \pageref{ModStabLibre}, \pageref{Nota2rang} and \pageref{defiRang}}{Nota1rang}
\NOTA{\Adj\,B\phantom{\wi{B}}}{(or $\wi{B}$) cotransposed matrix of $B$}{NOTACotrans}
\NOTA{\cD_k(G)}{\idd of order $k$ of the matrix $G$}{defIdDet}
\NOTA{\cD_k(\varphi)}{\idd of order $k$ of the \ali $\varphi$,
 see also \paref{exoIDDPTF1}}{fact.idd prod}
\NOTA{{\rg(\varphi)\geq k}}{notation
understood with \Dfn \ref{defRangk},
see also notation \ref{notaRgfi}}{defRangk}
\NOTA{{\rg(\varphi)\leq k}}{same thing}{defRangk}
\NOTA{\rE^{(n)}_{i,j}(\lambda)}{(or $\rE_{i,j}(\lambda)$) \elr matrix}{NOTAEn}
\NOTA{\En(\gA)}{\elr group}{NOTAEn}
\NOTA{\I_{k}}{identity matrix of order $k$}{lem.min.inv}
\NOTA{0_{k}}{square matrix of order $k$}{lem.min.inv}
\NOTA{0_{k,\ell}}{null matrix of type $k\times \ell$}{lem.min.inv}
\NOTA{\I_{k,q,m}}{standard simple matrix}{NOTAIkqm}
\NOTA{\I_{k,n}}{standard \mprnz}{NOTAIkqm}
\NOTA{A_{\alpha,\beta}}{extracted matrix}{NOTAextraite}
\NOTA{\Adj_{\alpha,\beta}(A)}{see notation \ref{notaAdjalbe}}{notaAdjalbe}
\NOTA{\cP_{\ell}}{set of finite subsets of $\{1,\ldots ,\ell\}$}{notaAdjalbe}
\NOTA{\cP_{k,\ell}}{subsets with $k$ \eltsz}{notaAdjalbe}
\NOTA{\GAnk(\gA)}{subsets of $\GAn(\gA)$: \mprns of rank $k$}{defiGrassmanniennes}
\NOTA{\GGnk(\gA)}{projective Grassmannian over $\gA$}{defiGrassmanniennes}
\NOTA{\GGn(\gA)}{projective Grassmannian over $\gA$}{defiGrassmanniennes}
\NOTA{\Pn(\gA)}{projective space of dimension $n$ over $\gA$}{defiGrassmanniennes}
\NOTA{\Diag(a_1,\ldots,a_n)}{
diagonal square matrix}{NOTADiag}
\NOTA{\Tr(\varphi)}{trace of $\varphi$ (\endo of $\Ae n$),
see also \paref{NOTAPolcar}}{NOTA1Polcar}
\NOTA{\rC\varphi(X)}{\polcar of $\varphi$ (idem), see also \paref{NOTAPolcar}}{NOTA1Polcar}
\NOTA{\dex{\gB:\gA}}{$\rg_\gA(\gB)$, see also  \paref{subsecTransRang} and \ref{notaTraceDetCarAlg}}{notaCTrN}
\NOTA{\Tr_{\gB/\!\gA}(a)}{trace of (the multiplication by) $a$, see also \vref{defi2STF}}{notaCTrN}
\NOTA{\rN_{\gB/\!\gA}(a)}{norm of $a$, see also \vref{defi2STF}}{notaCTrN}
\NOTA{\rC{\gB/\!\gA}(a)}{\polcar of (the multiplication by) $a$, see also \vref{defi2STF}}{notaCTrN}
\NOTA{\Gram_\gA(\varphi,\ux)}{
Gram matrix of  $(\ux)$ for $\varphi$}{defiGram}
\NOTA{\gram_\gA(\varphi,\ux)}{
Gram \deter of $(\ux)$ for $\varphi$}{defiGram}
\NOTA{\disc_{\gB/\!\gA}(\ux)}{
\discri of the family $(\ux)$}{defiDiscTra}
\NOTA{\Disc_{\gB/\!\gA}}{discriminant of a free extension}{defiDiscTra}
\NOTA{\Lin_\gA(M,N)}{\Amo of \alisz}{NOTAAlis}
\NOTA{\End_\gA(M)}{$\Lin_\gA(M,M)$}{NOTAAlis}
\NOTA{M\sta}{dual module of $M$}{NOTAAlis}
\NOTA{\AuX_d}{\Asub of $\AuX$ of the \pogs of degree~$d$}{NOTAAXd}
\CHAP{\nameref{chapGenerique}}
\NOTA{\Pf(E)}{set of  finite subsets of $E$}{NOTAPfPfe}
\NOTA{\Pfe(E)}{set of  finitely enumerated subsets of $E$}{NOTAPfPfe}
\NOTA{\Hom_\gA(\gB,\gB')}{set of  \homos of \Algs from $\gB$ to $\gB'$ }{def0Alg}
\NOTA{\mu_{M,b}}{(or $\mu_b$)  $y\mapsto by$, $\in\End_\gB(M)$ ($b\in \gB$, $M$ a \Bmoz)}{NOTAmux}
\NOTA{\cJ(f)}{\id of the symmetric relators}{definotaAdu}
\NOTA{\Adu_{\gA,f}}{\adu of $f$ over $\gA$}{definotaAdu}
\NOTA{\disc_{X}(f)}{discriminant of the \polu  $f$ of $\AX$}{eqDiscri}
\NOTA{\Tsc_g(f)}{Tschirnhaus transform of $f$ by $g$}{defiTschir}
\NOTA{\Mip_{\gK,x}(T)}{or $\Mip_{x}(T)$, \mon \polmin of $x$ (over the field $\gK$)}{defiSTF}
\NOTA{G.x}{orbit of $x$ under $G$}{NOTAStStp}
\NOTA{{G.x=\so{x_1,\ldots,x_k}}}{
orbit enumerated without repetition with $x_1=x$ }{NOTAStStp}
\NOTA{\St_G(x)}{(or $\St(x)$) stabilizer subgroup of the point $x$}{NOTAStStp}
\NOTA{\Stp_G(F)}{(or $\Stp(F)$) point by point stabilizer of the subset $F$}{NOTAStStp}
\NOTA{\idg{G:H}}{index of the subgroup $H$ in the group $G$: $\#(G/H)$}{NOTAStStp}
\NOTA{\Fix_E(H)}{(or $E^H$) subset of $E$ formed from the fixed points of $H$}{NOTAStStp}
\NOTA{\sigma\in G/H}{we take a $\sigma$ in each left coset of $H$ in $G$}{NOTAStStp}
\NOTA{\rC{G}(x)(T)}{$=\prod_{\sigma\in G}(T-\sigma(x))$}{NOTAStStp}
\NOTA{\rN_G(x)}{$=\prod_{\sigma\in G}\sigma(x)$}{NOTAStStp}
\NOTA{\Tr_G(x)}{$=\som_{\sigma\in G}\sigma(x)$}{NOTAStStp}
\NOTA{\Rv_{G,x}(T)}{resolvent of $x$ (relative to $G$)}{NOTAStStp}
\NOTA{\Aut_\gA(\gB)}{group of  $\gA$-\autos of $\gB$}{NOTAAutAB}
\NOTA{\Gal(\gL/\gK)}{idem, for a Galois extension}{defiCorGal}
\NOTA{\cG_{\gL/\gK}}{finite subgroups of $\Aut_\gK(\gL)$}{defiCorGal}
\NOTA{\cK_{\gL/\gK}}{\stfes $\gK$-subextensions of $\gL$}{defiCorGal}
\NOTA{\Gal_\gK(f)}{Galois group of the \spl \pol $f$}{NOTAGalKf}
\NOTA{\Syl_X(f,p,g,q)}{
Sylvester matrix of $f$ and $g$
in degrees $p$ and $q$}{secRes}
\NOTA{\Res_X(f,p,g,q)}{
resultant of the \pols $f$ and $g$
in degrees $p$ and $q$}{secRes}
\NOTA{\car(\gK)}{\cara of a field}{NOTACarK}
\NOTA{\Adj\iBA (x)}{or $\wi x$: cotransposed \eltz, see also \paref{eqeltcotransp}}{eqelt0cotransp}
\NOTA{(\gA:\gB)}{conductor of $\gB$ into $\gA$}{defiConducteur}
\NOTA{\fR_X(f,g_1,\ldots,g_r)}{}{lemElimPlusieurs}
\NOTA{\JJ_\uX (\uf)}{Jacobian matrix of a \sypz}{secNewton}
\NOTA{\J_{\uX}(\uf)}{Jacobian of a \sypz}{secNewton}
\NOTA{\idg{L:E}_\gA}{index of a \tf submodule in a free module}{exolemSousLibre}

\CHAP{\nameref{chap mpf}}
\NOTA{R_\ua}{matrix of the trivial relators}{secRelTrivSeqReg}
\NOTA{\scp{\ux}{\uz}}{$\sum_{i=1}^n x_iz_i$}{NOTAProdScal}
\NOTA{\fm_\uxi}{$\gen{x_1-\xi_1,\ldots,x_n-\xi_n}_\gA$: \id of the zero $\uxi$}{secExempleGeo}
\NOTA{M\otimes_\gA N}{tensor product of two \Amosz}{propPftens}
\NOTA{\Al k_\gA M}{$k^{\rm th}$ exterior power of $M$}{propPfPex}
\NOTA{\gS_\Ae k M}{$k^{\rm th}$ symmetric power of $M$}{propPfPex}
\NOTA{\rho\ist(M)}{\Bmo obtained from the \Amo $M$  by the \eds $\rho:\gA\rightarrow \gB$}{NOTArhosta}
\NOTA{\cF_n(M)}{or $\cF_{\gA,n}(M)$:  $n^{\rm th}$ \idf of the \tf \Amo $M$}{def ide fit}
\NOTA{\fRes_{X}(\ff)}{resultant \id of $\ff$
(with a \polu in $\ff$)}{thElimAff}
\NOTA{\cK_n(M)}{$n^{\rm th}$ Kaplansky \id of the \Amo  $M$}{exoAutresIdF}

\CHAP{\nameref{chap ptf0}}
\NOTA{\theta_{M,N}}{natural \Ali $M\sta\te_\gA N\to \Lin_\gA(M,N)$}{NOTAthetaMN}
\NOTA{\theta_{M}}{natural \Ali $M\sta\te_\gA M\to \End_\gA(M)$}{NOTAthetaM}
\NOTA{\Diag(M_1,\ldots,M_n)}{%
block diagonal matrix}{Notadiagblocs} 
\NOTA{\Bdim\gA<n}{stable range (of Bass) less than or equal to $n$}{defiStableRange}
\NOTA{\det\varphi}{\deter of the \endo $\varphi$ of a \mptfz }{NOTAdet}
\NOTA{\rC\varphi(X)}{\polcar of $\varphi$ \ldots\, (idem) }{NOTAPolcar}
\NOTA{\wi{\varphi}}{cotransposed \endo of $\varphi$ \ldots\, (idem) }{NOTACotransp}
\NOTA{\rF{\varphi}(X)}{\polfon of  $\varphi$, i.e.,
$\det(\Id_{P}+X\varphi)$}{NOTAPolfon}
\NOTA{\Tr_P(\varphi)}{trace of the \endo $\varphi$}{NOTAPolfon}
\NOTA{\rR{P}(X)}{\polmu of the \mptf $P$}{NOTAPolmu}
\NOTA{\ide_{h}(P)}{the \idm  associated with the integer $h$ and with the \pro module $P$}{NOTAide}
\NOTA{P\ep{h}}{component of the module $P$ in rank~$h$}{NOTAep}

\CHAP{\nameref{chap AlgStricFi}}
\NOTA{\rC{\gB/\!\gA}(x)(T)}{\polcar of (the multiplication by) $x$}{defi2STF}
\NOTA{\rF{\gB/\!\gA}(x)(T)}{\polfon of (the multiplication by) $x$}{defi2STF}
\NOTA{\rN\iBA(x)}{norm of $x$: \deter of the multiplication by $x$}{defi2STF}
\NOTA{\Tr\iBA(x)}{trace of (the multiplication by) $x$}{defi2STF}
\NOTA{a\centerdot\alpha}{$\alpha\circ \mu_a:x\mapsto \alpha(ax)$}{definotaAsta}
\NOTA{\Adj\iBA(x)}{or $\wi x$: cotransposed \eltz}{eqeltcotransp}
\NOTA{[\gB:\gA]}{$\rg_\gA(\gB)$, see also pages \pageref{notaCTrN} and \pageref{notaTraceDetCarAlg}}{subsecTransRang}
\NOTA{\Phi_{\gA/\gk,\lambda}}{$\Phi_\lambda(x,y) = \lambda(xy)$}{defdualisante}
\NOTA{\phi\otimes\phi'}{tensor product of bi\lin forms}{NOTAptfb}
\NOTA{\env\gk\gA}{$\gA\otimes_\gk\gA$, enveloping \alg of $\gA/\gk$}{NotaAGenv}
\NOTA{\rJ\iAk}{\id of $\env\gk\gA$}{eqDiAk}
\NOTA{\Delta\iAk}{$\Delta(x)=x\te1-1\te x$}{eqDiAk}
\NOTA{\mu\iAk}{$\mu\iAk \left(\sum_ia_i\te b_i\right) =\sum_ia_ib_i$}{eqmuAk}
\NOTA{\Der{\gk}{\gA}{M}}{the \Amo of the \dvns of $\gA$ in $M$}{defiDeriv}
\NOTA{\mathrm{Der}(\gA)}{the \Amo of the \dvns of $\gA$}{defiDeriv}
\NOTA{\Om{\gk}{\gA}}{the \Amo of the  \diles (of K\"ahler) of~$\gA$}{thDerivUniv}
\NOTA{\vep\iAk}{\idm that generates $\Ann(\rJ\iAk)$, if it exists}{thSepIversen}
\NOTA{\LIN_\gk\!(\gA,\gA)}{\Amo of  \klis from $\gA$ to $\gA$}{NotaAGAL}
\NOTA{\PGL_n(\gA)}{quotient group $\GLn(\gA)/\Ati$}{NOTAPGL}
\NOTA{\rA_n}{subgroup of  even permutations of~$\Sn$}{NOTAAn}
\NOTA{\gk[G]}{\alg of a group, or of a \moz}{exoGroupAlgebra}

\CHAP{\nameref{ChapGalois}}
\NOTA{\BB(\gA)}{\agB of the \idms of $\gA$}{propB(A)}
\NOTA{\cB(f)}{\gui{canonical} basis of the \aduz}{notaBaseadu}


\CHAP {\nameref{chap Anneaux locaux}}
\NOTA{\Rad(\gA)}{radical of Jacobson  of $\gA$}{eqDefRadJac}
\NOTA{\gA(X)}{Nagata \ri of $\gA[X]$}{factLocNagata}
\NOTA{\Suslin(\bn)}{Suslin set of $(\bn)$}{defiSuslinSet}

\CHAP{\nameref{chap ptf1}}
\NOTA{\Gn}{$\Gn =\ZZ[(f_{i,j})_{i,j\in \lrbn}]/\cGn$}{NOTAGN}
\NOTA{\cGn}{relations obtained when writing $F^2=F$}{NOTAGN}
\NOTA{\HOp(\gA)}{semi-\ri of ranks of quasi-free \Amosz}{notaHO+} 
\NOTA{[P]_{\HOp (\gA)}}{or $[P]_\gA$, or $[P]$: class of a quasi-free \Amo in $\HOp(\gA)$}{notaHO+}
\NOTA{\rg_\gA(M)}{(\gnez) rank of the \ptf \Amo $M$}{NOTAHO+}
\NOTA{\HO\gA}{\ri of the ranks over $\gA$}{notaHO}
\NOTA{\dex{\gB:\gA}}{$\rg_\gA(\gB)$, see also pages \pageref{notaCTrN} and \pageref{subsecTransRang}}{notaTraceDetCarAlg}
\NOTA{\Gn(\gA)}{$\Gn\otimes_\ZZ\gA$}{subsubsec AGBR}
\NOTA{\cGnk}{$\cGn+\gen{1-r_k}$, with (in $\Gn$) $r_k=\ide_k(\Im\,F)$}{subsubsec AGBR}
\NOTA{\Gnk}{$\Gnk =\ZZ[(f_{i,j})_{i,j\in \lrbn}]/\cGnk$
or $\Gn[1/r_k]$}{subsubsec AGBR}
\NOTA{\GAnk(\gA)}{\gui{sub\vrtz} of $\GAn(\gA)$: \prrs of rank $k$}{subsubsec AGBR}
\NOTA{\GKO\gA}{semi-\ri of the \iso classes of \mptfs over $\gA$}{NOTAGKO}
\NOTA{\Pic\gA}{group of \iso classes  of  \mrcs 1 over $\gA$}{NOTAPic}
\NOTA{\KO\gA}{Grothendieck \ri of $\gA$}{NOTAK0}
\NOTA{[P]_{\KO (\gA)}}{or $[P]_\gA$, or $[P]$: class of a \ptf \Amo in $\KO (\gA)$}{NOTAK0}
\NOTA{\KTO\gA}{kernel of the rank \homo $\rg:\KO \gA\to\HO \gA$}{NOTAKTO}
\NOTA{\Ifr\gA}{monoid of the \tf \ifrs of the \ri $\gA$}{NOTAIfr}
\NOTA{\Gfr\gA}{group of \iv \elts of $\Ifr\gA$}{NOTAIfr}
\NOTA{\Cl \gA}{group of classes of \iv \ids (quotient of $\Gfr\gA$
by the subgroup of  \iv \idpsz)}{NOTAIfr}

\CHAP{\nameref{chapTrdi}}
\NOTA{\dar a}{$\sotq{x\in X}{x\leq a}$, see also \vpageref{NOTAdara}}{eqda}
\NOTA{\uar a}{$\sotq{x\in X}{x\geq a}$, see also \vpageref{NOTAuara}}{eqda}
\NOTA{\gT\eci}{opposite lattice of the lattice $\gT$}{NOTAtrdiopp}
\NOTA{\cI_\gT(J)}{\id generated by $J$ in the \trdi $\gT$}{NOTAidtrdi}
\NOTA{\cF_\gT(S)}{filter generated by $S$ in the \trdi $\gT$}{NOTAfitrdi}
\NOTA{{\gT/(J=0,U=1)}}{particular quotient lattice}{propIdealFiltre} 
\NOTA{\Bo(\gT)}{\agB generated by the \trdi $\gT$}{thBoolGen}
\NOTA{\ZZ^{(P)}}{\ort direct sum of copies of $\ZZ$, indexed by $P$}{NotaZZP}
\NOTA{\boxplus_{i\in I}G_i}{\ort direct sum of ordered groups}{NotaSDirOr}
\NOTA{\cC(a)}{solid subgroup generated by $a$ (in an \grlz)}{defiCongru}\NOTA{\DA(x_1,\ldots ,x_n)}{$\DA(\gen{x_1,\ldots ,x_n})$: an \elt of $\ZarA$}{notaZA}
\NOTA{\ZarA}{Zariski lattice of $\gA$}{notaZA}
\NOTA{\gA_{S}\sur{\fa}}{(or $S^{-1}\gA\sur{\fa}$) we invert the \elts of $S$ and we annihilate the \elts of $\fa$}{defiASa}\NOTA{\satu S \gA}{or $\sat{S}$: the filter obtained by saturating the \mo $S$ in $\gA$}{NOTASatu}
\NOTA{\Abul}{reduced \zed \ri generated by $\gA$}{thZedGen}
\NOTA{A \vda B}{$\Vi A\;\leq \;\Vu B$: \entrel}{notaVupVda}
\NOTA{\SpecT}{spectrum of the finite \trdi $\gT$, see also \paref{SpecTrdi}}{SpecTrdiFi}
\NOTA{(b:a)_\gT}{the conductor ideal of $a$ in $b$ (\trdisz)}%
{NOTATransp2}
\NOTA{\gA_\mathrm{pp}}{\qiri closure of $\gA$}{exoQiClot}
\NOTA{\Min\gA}{subspace of $\SpecA$ formed by the \idemisz}{exoMinA}

\CHAP{\nameref{ChapAdpc}}

\NOTA{\fa\div\fb}{$\sotq{x\in\Frac\gA}{x\fb\subseteq\fa}$}{NOTAfadivfb}
\NOTA{\gA[\fa t]}{Rees \alg of the \id $\fa$ of $\gA$}{NOTARees}
\NOTA{\Icl_\gA(\fa)}{\cli of the \id $\fa$ in $\gA$}{fact2Entiers}

\CHAP{\nameref{chapKrulldim}}
\NOTA{\SpecA}{Zariski spectrum of the \ri $\gA$}{nota Spec(A)}
\NOTA{\fD_\gA(\xn)}{compact-open set of $\SpecA$}{nota Spec(A)}
\NOTA{\SpecT}{spectrum of the \trdi $\gT$}{SpecTrdi}
\NOTA{\fD_\gT(u)}{compact-open set of $\SpecT$}{SpecTrdi}
\NOTA{\OQC(\gT)}{\trdi of the compact-open sets of $\SpecT$}{SpecTrdi}
\NOTA{\JK_\gA(x)}{$\gen{x}+(\DA(0):x)$: Krull boundary \id of $x$ in $\gA$}{defZar2}
\NOTA{\JK_\gA(\fa)}{$\fa+(\DA(0):\fa)$: Krull boundary \id of $\fa$ in $\gA$}{defZar2}
\NOTA{\gA_\rK^x}{$\gA\sur{\JK_\gA(x)}$: upper boundary \ri of $x$ in $\gA$}{defZar2}
\NOTA{\SK_\gA(x)}{$x^\NN(1+x\gA)$: Krull boundary \mo of $x$ in $\gA$}{defZar2}
\NOTA{\gA^\rK_{x}}{$(\SK_\gA(x))^{-1}\gA$: lower boundary \ri of $x$ in $\gA$}{defZar2}
\NOTA{\Kdim\gA\leq r}{the Krull dimension of the \ri $\gA$ is $\leq r$}{NOTAKdim}
\NOTA{\Kdim\gA\leq\Kdim\gB}{}{notaKdiminf}
\NOTA{\SK_\gA(\xzk)}{iterated Krull boundary \moz}{notaBordsIteres}
\NOTA{\JK_\gA(\xzk)}{iterated Krull boundary \idz}{notaBordsIteres}
\NOTA{\IK_\gA(\xzk)}{iterated Krull boundary \idz, variant}{notaBordsIteres}
\NOTA{\Kdim\gT\leq r}{the \ddk of the \trdi $\gT$ is $\leq r$}{defiDDKTRDI}
\NOTA{\JK_\gT(x)}{$\dar x \,\vu\, (0:x)_\gT$: Krull boundary \id of $x$ in the \trdiz~$\gT$}{defiBordKrullTrdi}
\NOTA{\gT_\rK^x}{$\gT\sur{\JK_\gT(x)}$: upper boundary lattice of $x$}{defiBordKrullTrdi}
\NOTA{\JK_\gT(\xzk)}{iterated Krull boundary \id in a \trdiz}{eq1IdBordKrullItereTrdi}
\NOTA{\Kdim \rho}{\ddk of the morphism $\rho$}{defiKdimMor}
\NOTA{\gA_{\so{a}}}{$\gA\sur{a\epr}\times \gA\sur{({a\epr})\epr}$}{lem20MorRc}
\NOTA{\Amin}{minimal \qiri closure of $\gA$}{thAmin}
\NOTA{\Vdim\gA}{valuative dimension}{defiValdim}
\CHAP{\nameref{chapNbGtrs}}
\NOTA{\JA(\fa)}{Jacobson radical of the \id $\fa$  of $\gA$}{notaJA}
\NOTA{\JA(x_1,\ldots ,x_n)}{
$\JA(\gen{x_1,\ldots ,x_n})$: an \elt of $\HeA$}{notaJA}
\NOTA{\HeA}{Heitmann lattice of $\gA$}{defHeit}
\NOTA{\Jdim}{dimension of the Heitmann J-spectrum}{NOTAJdim}
\NOTA{\Max\gA}{subspace of $\SpecA$ formed by the \idemasz}{NOTAJdim}
\NOTA{\Jspec\gA}{$\Spec(\HeA)$: Heitmann J-spectrum}{NOTAJdim} 
\NOTA{\JH_\gA(x)}{$\gen{x}+(\JA(0):x)$: Heitmann boundary \id (of $x$ in $\gA$)}{defHei2}
\NOTA{\gA_\rH^x}{$\gA\sur{\JH_\gA(x)}$: the Heitmann boundary \ri of $x$}{defHei2}
\NOTA{\Hdim}{Heitmann dimension}{defDHA}
\NOTA{\Sdim\gA<n}{}{defiSdimGdim}
\NOTA{\Gdim\gA<n}{}{defiSdimGdim}
\NOTA{\Cdim\gA<n}{the \ri $\gA$ is $n$-stable}{defiAnneaunStable}

\CHAP{\nameref{chapPlg}}
\NOTA{\cM(U)}{the \mo generated by the \elt or the subset $U$ of $\gA$}{nota mopf}
\NOTA{\cS(I,U)}{$\sotq{v\in\gA}{\exists u\in\cM(U)\;\exists a\in \gen{I}_\gA,
\;  v=u+a }$}{nota mopf}
\NOTA{\cS(a_1,\ldots,a_k;u_1,\ldots,u_\ell)}
{$\cS(\so{a_1,\ldots,a_k},\so{u_1,\ldots,u_\ell})$}{nota mopf}

\CHAP{\nameref{ChapMPEtendus}}
\NOTA{\gA\lra{X}}{localized \ri of $\gA[X]$ at the \polusz}{NOTAAlraX}
\NOTA{A\sims{\cG }B}{there exists a matrix $H\in\cG $ such that $HA=B$}{NOTAfGg}

\CHAP{\nameref{ChapSuslinStab}}
\NOTA{\GL(P)}{group of  \lin \autos of $P$}{NOTAtransvec}
\NOTA{\wi{\EE}(P)}{subgroup  of $\GL(P)$ generated by the transvections}{NOTAtransvec}
\NOTA{\meck {a}{b}}{Mennicke symbol}{lemMennicke1}
\NOTA{\GLn(\gB,\fb)}{kernel of $\GLn(\gB)\to\GLn(\gB\sur{\fb})$}{notaGLIdeal}
\NOTA{\En(\gB, \fb)}{normal subgroup of $\En(\gB)$ generated by the $\rE_{ij}(b)$ with $b \in \fb$.}{notaGLIdeal}
\CHAP {Annex: \cov logic}
\NOTA{\rP(X)}{the class of subsets of $X$}{P(X)}

\newpage
\thispagestyle{CMcadreseul}
\cleardoublepage 
\rdb


\addtocontents{toc}{\vskip-0.8em}
\addcontentsline{toc}{chapterbis}{Index}
\markboth{Index}{Index}

\printindex

\finincrementeexosetprob
\end{document}